\documentclass[11pt,openany]{book}
\usepackage{geometry}                		% See geometry.pdf to learn the layout options. There are lots.
\usepackage{graphicx}				% Use pdf, png, jpg, or eps� with pdflatex; use eps in DVI mode						
\usepackage{amssymb}\vspace{0.25cm}
\usepackage{mathtools}
\usepackage{setspace}
\usepackage{tikz-cd}
\usepackage{amsmath, amsthm}
\usepackage[mathscr]{euscript}
\newcommand{\abs}[1]{\left\vert#1\right\vert}
\usepackage[normalem]{ulem}
\usepackage{manfnt}
\usepackage{pgfplots}
\usepackage{pifont}
\usepackage[safe]{tipa}
\usepackage{marvosym}
\usepackage{scalerel,stackengine}
\usepackage{titlesec}
\usepackage{makeidx}
\usepackage{centernot}
\usepackage{xspace}
\usepackage[noauto]{chappg}
\usepackage[OT2,T1]{fontenc}
\usepackage{hyperref}
\stackMath
\newcommand\reallywidehat[1]{%
\savestack{\tmpbox}{\stretchto{%
  \scaleto{%
    \scalerel*[\widthof{\ensuremath{#1}}]{\kern-.6pt\bigwedge\kern-.6pt}%
    {\rule[-\textheight/2]{1ex}{\textheight}}%WIDTH-LIMITED BIG WEDGE
  }{\textheight}% 
}{0.5ex}}%
\stackon[1pt]{#1}{\tmpbox}%
}
\newcommand\R{%
\mathbb{R}
}

\newcommand\C{%
\mathbb{C}
}
\newcommand\N{%
\mathbb{N}
}
\newcommand\Z{%
\mathbb{Z}
}
\newcommand\PP{%
\mathbb{P}
}
\newcommand\Q{%
\mathbb{Q}
}
\newcommand\HH{%
\mathbb{H}
}
\newcommand\F{%
\mathbb{F}
}

\newcommand\Y{%
\mathbb{Y}
}
\renewcommand\L{%
\mathcal{L}^1
}

\newcommand\sB{%
\mathcal{B}
}
\newcommand\sC{%
\mathcal{C}
}
\newcommand\sD{%
\mathcal{D}
}
\newcommand\sE{%
\mathcal{E}
}
\newcommand\sA{%
\mathcal{A}
}
\newcommand\sF{%
\mathcal{F}
}
\newcommand\sG{%
\mathcal{G}
}
\newcommand\sH{%
\mathcal{H}
}
\newcommand\sI{%
\mathcal{I}
}

\newcommand\sL{%
\mathcal{L}
}
\newcommand\sK{%
\mathcal{K}
}
\newcommand\sM{%
\mathcal{M}
}
\newcommand\sO{%
\mathcal{O}
}
\newcommand\sP{%
\mathcal{P}
}
\newcommand\sQ{%
\mathcal{Q}
}
\newcommand\sS{%
\mathcal{S}
}
\newcommand\sT{%
\mathcal{T}
}
\newcommand\sU{%
\mathcal{U}
}
\newcommand\sV{%
\mathcal{V}
}
\newcommand\sW{%
\mathcal{W}
}
\newcommand\sX{%
\mathcal{X}
}
\newcommand\sY{%
\mathcal{Y}
}
\newcommand\sZ{%
\mathcal{Z}
}
\newcommand\A{%
\mathcal{A}
}

\newcommand\B{%
\mathcal{B}
}

\newcommand\sR{%
\mathcal{R}
}
\newcommand\Ob{%
\text{Ob$_{}$}\hspace{0.03cm}
}
\newcommand\ob{%
\text{ob}
}
\newcommand\id{%
\text{id}
}
\newcommand\gro{%
\text{gro}
}
\newcommand\tran{%
\text{tran}\hspace{0.03cm}
}
\newcommand\holim{%
\text{holim}
}
\newcommand\holimx{%
\text{holim}\hspace{0.05cm}
}
\newcommand\hocolim{%
\text{hocolim}
}
\newcommand\hocolimx{%
\text{hocolim}\hspace{0.05cm}
}

\newcommand\ran{%
\text{ran$_{}$}\xspace
}
\newcommand\lan{%
\text{lan}
}
\newcommand\FGL{%
\text{FGL}
}
\newcommand\IPS{%
\text{IPS}
}
\newcommand\tot{%
\text{tot}\hspace{0.03cm}
}
\newcommand\res{%
\text{res}%\hspace{0.03cm}
}
\newcommand\cosin{%
\text{cosin}
}
\newcommand\sk{%
\text{sk}
}
\newcommand\cosk{%
\text{cosk}
}
\newcommand\tr{%
\text{tr}
}
\newcommand\Ex{%
\text{Ex }
}
\newcommand\Exx{%
\text{Ex}
}
\newcommand\UW{%
\text{UW}
}
\newcommand\HR{%
\text{HR}
}
\newcommand\eq{%
\text{eq}
}
\newcommand\coeq{%
\text{coeq}
}
\newcommand\Sd{%
\text{Sd}
}
\newcommand\sd{%
\text{sd}
}
\newcommand\con{%
\text{con}
}
\newcommand\cpt{%
\text{cpt}
}
\newcommand\cptx{%
\text{cpt\hspace{0.05cm}}
}
\newcommand\im{%
\text{im$_{}$}\hspace{0.05cm}
}
\newcommand\coim{%
\text{coim$_{}$}\xspace
}

\newcommand\barr{%
\text{bar}
}
\newcommand\iso{%
\text{iso}\hspace{0.05cm}
}
\newcommand\End{%
\text{End}
}
\newcommand\Tor{%
\text{Tor}
}
\newcommand\tor{%
\text{tor}
}
\newcommand\xdiv{%
\text{div}\hspace{0.05 cm}
}
\newcommand\divx{%
\text{div}
}
\newcommand\Ext{%
\text{Ext}
}
\newcommand\Extx{%
\text{Ext }
}
\newcommand\dux{%
\text{du\hspace{0.05cm}}
}
\newcommand\dimx{%
\text{dim\hspace{0.03cm}}
}
\newcommand\tel{%
\text{tel}\hspace{0.03cm}
}
\newcommand\telsub{%
\text{tel}
}
\newcommand\telx{%
\text{tel}
}
\newcommand\coker{%
\text{coker}\hspace{0.05cm}
}
\newcommand\cokersub{%
\text{coker}
}
\newcommand\Ph{%
\text{Ph}
}
\newcommand\Ann{%
\text{Ann}
}

\newcommand\circx{%
\hspace{0.03cm} \circ \hspace{0.03cm}
}
\newcommand\bA{%
\textbf{A}
}
\newcommand\bB{%
\textbf{B}
}
\newcommand\bb{%
\textbf{b}
}
\newcommand\tb{%
\text{b}
}

\newcommand\bC{%
\textbf{C}\xspace
}
\newcommand\bD{%
\textbf{D}\xspace
}
\newcommand\td{%
\text{d}\xspace
}
\newcommand\tD{%
\text{D}
}
\newcommand\bE{%
\textbf{E}\xspace
}
\newcommand\be{%
\textbf{e}\xspace
}
\newcommand\tE{%
\text{E}
}
\newcommand\bc{%
\textbf{c}\xspace
}
\newcommand\bcf{%
\textbf{cf}\xspace
}
\newcommand\bF{%
\textbf{F}\xspace
}
\newcommand\bff{%
\textbf{f}\xspace
}
\newcommand\bG{%
\textbf{G}\xspace
}
\newcommand\bg{%
\textbf{g}\xspace
}
\newcommand\bh{%
\textbf{h}\xspace
}
\newcommand\bH{%
\textbf{H}\xspace
}
\newcommand\tH{%
\text{H}
}
\newcommand\bI{%
\textbf{I}\xspace
}
\newcommand\bJ{%
\textbf{J}\xspace
}
\newcommand\bK{%
\textbf{K}\xspace
}
\newcommand\bk{%
\textbf{k}\xspace
}
\newcommand\bL{%
\textbf{L}\xspace
}
\newcommand\tL{%
\text{L}
}
\newcommand\bloc{%
\textbf{loc}
}
\newcommand\bM{%
\textbf{M}\xspace
}
\newcommand\tM{%
\text{M}\xspace
}
\newcommand\bm{%
\textbf{m}\xspace
}
\newcommand\bN{%
\textbf{N}\xspace
}
\newcommand\bn{%
\textbf{n}\xspace
}
\newcommand\bO{%
\textbf{O}\xspace
}
\newcommand\bP{%
\textbf{P}\xspace
}
\newcommand\bp{%
\textbf{p}\xspace
}
\newcommand\bR{%
\textbf{R}
}
\newcommand\tR{%
\text{R}
}
\newcommand\bQ{%
\textbf{Q}\xspace
}
\newcommand\bS{%
\textbf{S}\xspace
}
\newcommand\bss{%
\textbf{s}\xspace
}
\newcommand\bsin{%
\textbf{sin}\hspace{0.03cm}
}
\newcommand\bT{%
\textbf{T}\xspace
}
\newcommand\bt{%
\textbf{t}\xspace
}
\newcommand\tT{%
\text{T}
}
\newcommand\beS{%
\textbf{eS}
}
\newcommand\bU{%
\textbf{U}\xspace
}
\newcommand\bu{%
\textbf{u}\xspace
}
\newcommand\bV{%
\textbf{V}\xspace
}
\newcommand\tV{%
\text{V}
}
\newcommand\bv{%
\textbf{v}\xspace
}
\newcommand\bW{%
\textbf{W}\xspace
}
\newcommand\bw{%
\textbf{w}\xspace
}
\newcommand\bX{%
\textbf{X}\xspace
}
\newcommand\bx{%
\textbf{x}\xspace
}
\newcommand\bY{%
\textbf{Y}\xspace
}
\newcommand\bZ{%
\textbf{Z}\xspace
}
\newcommand\bco{%
\textbf{co}
}
\newcommand\bwco{%
\textbf{wco}
}
\newcommand\bzero{%
\textbf{0}\xspace
}
\newcommand\bone{%
\textbf{1}\xspace
}
\newcommand\btwo{%
\textbf{2}\xspace
}
\newcommand\bLEX{%
\textbf{LEX}\xspace
}
\newcommand\bEX{%
\textbf{EX}
}
\newcommand\bIND{%
\textbf{IND}\xspace
}
\newcommand\bPRO{%
\textbf{PRO}
}

\newcommand\bSET{%
\textbf{SET}\xspace
}
\newcommand\bSISET{%
\textbf{SISET}\xspace
}
\newcommand\bSSISET{%
\textbf{SSISET}\xspace
}
\newcommand\bBISISET{%
\textbf{BISISET}\xspace
}
\newcommand\bHSISET{%
\textbf{HSISET}\xspace
}
\newcommand\bHZEROSISET{%
\textbf{H$_0$SISET}\xspace
}
\newcommand\bHBISISET{%
\textbf{HBISISET}\xspace
}
\newcommand\bHZEROC{%
\textbf{H$_0$C}
}
\newcommand\bSIC{%
\textbf{SIC}\xspace
}
\newcommand\bBISIC{%
\textbf{BISIC}\xspace
}
\newcommand\bSITOP{%
\textbf{SITOP}\xspace
}
\newcommand\bSICG{%
\textbf{SICG}\xspace
}
\newcommand\bSICGH{%
\textbf{SICGH}\xspace
}
\newcommand\bFIL{%
\textbf{FIL}\xspace
}
\newcommand\bFILSP{%
\textbf{FILSP}
}
\newcommand\bCAT{%
\textbf{CAT}\xspace
}
\newcommand\bHCAT{%
\textbf{HCAT}\xspace
}
\newcommand\bCG{%
\textbf{CG}\xspace
}
\newcommand\bCW{%
\textbf{CW}\xspace
}
\newcommand\bHCW{%
\textbf{HCW}\xspace
}
\newcommand\bHC{%
\textbf{HC}\xspace
}
\newcommand\bHD{%
\textbf{HD}\xspace
}
\newcommand\bCWSP{%
\textbf{CWSP}\xspace
}
\newcommand\bHCWSP{%
\textbf{HCWSP}\xspace
}
\newcommand\bCONCW{%
\textbf{CONCW}\xspace
}
\newcommand\bHCONCW{%
\textbf{HCONCW}\xspace
}
\newcommand\bFCONCW{%
\textbf{FCONCW}\xspace
}
\newcommand\bHFCONCW{%
\textbf{HFCONCW}%\xspace
}
\newcommand\bPREHCONCW{%
\textbf{PREHCONCW}%\xspace
}
\newcommand\bCPTHAUS{%
\textbf{CPTHAUS}\xspace
}
\newcommand\bCONCWSP{%
\textbf{CONCWSP}\xspace
}
\newcommand\bHCONCWSP{%
\textbf{HCONCWSP}\xspace
}
\newcommand\bHSCONCWSP{%
\textbf{HSCONCWSP}\xspace
}
\newcommand\bNILCWSP{%
\textbf{NILCWSP}\xspace
}
\newcommand\bHNILCWSP{%
\textbf{HNILCWSP}\xspace
}
\newcommand\bCGH{%
\textbf{CGH}\xspace
}
\newcommand\bCOSIC{%
\textbf{COSIC}\xspace
}
\newcommand\bTOP{%
\textbf{TOP}\xspace
}
\newcommand\bHAUS{%
\textbf{HAUS}\xspace
}
\newcommand\CPTHAUS{%
\textbf{CPTHAUS}\xspace
}
\newcommand\bHTOP{%
\textbf{HTOP}\xspace
}
\newcommand\bMON{%
\textbf{MON}\xspace
}
\newcommand\bALG{%
\textbf{ALG}\xspace
}
\newcommand\bSALG{%
\textbf{S-ALG}\xspace
}
\newcommand\bCOALG{%
\textbf{COALG}\xspace
}
\newcommand\bGR{%
\textbf{GR}\xspace
}
\newcommand\bTOPGR{%
\textbf{TOPGR}\xspace
}
\newcommand\bGRD{%
\textbf{GRD}\xspace
}
\newcommand\bPRESPEC{%
\textbf{PRESPEC}\xspace
}
\newcommand\bBIPRESPEC{%
\textbf{BIPRESPEC}\xspace
}
\newcommand\bTOW{%
\textbf{TOW}\xspace
}
\newcommand\bAB{%
\textbf{AB}\xspace
}
\newcommand\bNIL{%
\textbf{NIL}\xspace
}
\newcommand\bRG{%
\textbf{RG}\xspace
}
\newcommand\bMOD{%
\textbf{MOD}\xspace
}
\newcommand\bAMOD{%
\text{$A$-}\textbf{MOD}\xspace
}

\newcommand\bMODA{%
\textbf{MOD}\text{-$A$}\xspace
}
\newcommand\AbMOD{%
\textbf{$A$-MOD}\xspace
}
\newcommand\GbMOD{%
\textbf{$G$-MOD}\xspace
}
\newcommand\bACT{%
\textbf{ACT}\xspace
}
\newcommand\bRMOD{%
\textbf{R-MOD}\xspace
}
\newcommand\RMOD{%
\text{$R$-\textbf{MOD}}\xspace
}
\newcommand\bFRC{%
\textbf{FRC}\xspace
}
\newcommand\bUFRC{%
\textbf{UFRC}\xspace
}
\newcommand\bSIMC{%
\textbf{SIMC}\xspace
}
\newcommand\bSIAB{%
\textbf{SIAB}\xspace
}
\newcommand\bSIGR{%
\textbf{SIGR}\xspace
}
\newcommand\bSIMD{%
\textbf{SIMD}\xspace
}
\newcommand\bTRI{%
\textbf{TRI}
}
\newcommand\bIDTRI{%
\textbf{IDTRI}
}
\newcommand\tPER{%
\text{PER}
}
\newcommand\bPER{%
\textbf{PER}
}
\newcommand\bSPEC{%
\textbf{SPEC}\xspace
}

\newcommand\bSEPPRESPEC{%
\textbf{SEPPRESPEC}\xspace
}
\newcommand\bINJPRESPEC{%
\textbf{INJPRESPEC}\xspace
}
\newcommand\bHWSPEC{%
\textbf{HWSPEC}\xspace
}
\newcommand\bHSPEC{%
\textbf{HSPEC}\xspace
}
\newcommand\bWPRESPEC{%
\textbf{WPREPSEC}\xspace
}
\newcommand\bHWPRESPEC{%
\textbf{HWPRESPEC}\xspace
}
\newcommand\bHWOMEGAPRESPEC{%
\textbf{HW$\Omega$-PRESPEC}\xspace
}
\newcommand\bUN{%
\textbf{UN}
}
\newcommand\bHT{%
\textbf{HT}
}
\newcommand\bCT{%
\textbf{CT}
}

\newcommand\bKU{%
\textbf{KU}\xspace
}
\newcommand\bku{%
\textbf{ku}
}
\newcommand\bi{%
\textbf{i}
}
\newcommand\bj{%
\textbf{j}
}
\newcommand\bl{%
\textbf{l}
}
\newcommand\br{%
\textbf{r}
}
\newcommand\bKO{%
\textbf{KO}\xspace
}
\newcommand\bBP{%
\textbf{BP}\xspace
}
\newcommand\bFGL{%
\textbf{FGL}\xspace
}
\newcommand\bMU{%
\textbf{MU}\xspace
}
\newcommand\bEM{%
\textbf{EM}\xspace
}
\newcommand\bER{%
\textbf{ER}\xspace
}
\newcommand\bACY{%
\textbf{ACY}\xspace
}
\newcommand\bWALD{%
\textbf{WALD}\xspace
}
\newcommand\bSIWALD{%
\textbf{SIWALD}\xspace
}
\newcommand\bBUN{%
\textbf{BUN}\xspace
}
\newcommand\bFIB{%
\textbf{FIB}\xspace
}
\newcommand\bPRIN{%
\textbf{PRIN}\xspace
}
\newcommand\bLCCS{%
\textbf{LCCS}\xspace
}
\newcommand\bLACT{%
\textbf{LACT}\xspace
}

\newcommand\bGM{%
\textbf{GM}\xspace
}
\newcommand\bGL{%
\textbf{GL}\xspace
}
\newcommand\bSL{%
\textbf{SL}\xspace
}
\newcommand\bSO{%
\textbf{SO}\xspace
}
\newcommand\bST{%
\textbf{ST}\xspace
}
\newcommand\bKL{%
\textbf{KL}\xspace
}
\newcommand\bPh{%
\textbf{Ph}\xspace
}
\newcommand\bBA{%
\textbf{BA}\xspace
}
\newcommand\bDL{%
\textbf{DL}\xspace
}
\newcommand\bOP{%
\textbf{OP}
}

\newcommand\bOPERc{%
\textbf{OPER$_\textbf{C}$}
}
\newcommand\bisox{%
\textbf{$[(\textbf{iso}\boldsymbol\Gamma)^\textbf{OP},\textbf{C}]$}
}
\newcommand\bMONx{%
\textbf{MON$_{[(\textbf{iso}\boldsymbol\Gamma)^\textbf{OP},\textbf{C}]}$}\xspace
}
\newcommand\pspsp{%
\textbf{ps}\hspace{0.03cm}\boldsymbol\Pi\text{-}\textbf{SP}
}
\newcommand\psosp{%
\textbf{ps}\hspace{0.03cm}\widehat{\sO}\text{-}\textbf{SP}
}

\newcommand\bSP{%
\textbf{SP}
}

\newcommand\bBV{%
\textbf{BV}
}
\newcommand\xBV{%
\text{BV}
}
\newcommand\bDelta{%
\boldsymbol\Delta
}
\newcommand\bdn{%
{\boldsymbol\Delta}_n
}
\newcommand\dcg{%
\boldsymbol\Delta\text{-}\textbf{CG}
}
\newcommand\bGamma{%
\boldsymbol\Gamma
}
\newcommand\bgamma{%
\boldsymbol\gamma
}
\newcommand\balpha{%
\boldsymbol\alpha
}
\newcommand\bLambda{%
\boldsymbol\Lambda
}
\newcommand\bPi{%
\boldsymbol\Pi
}
\newcommand\bPhi{%
\boldsymbol\Phi
}
\newcommand\bXi{%
\boldsymbol\Xi
}
\newcommand\bpi{%
\boldsymbol\pi
}
\newcommand\brho{%
\boldsymbol\rho
}
\newcommand\btheta{%
\boldsymbol\theta
}
\newcommand\bxi{%
\boldsymbol\xi
}
\newcommand\bbeta{%
\boldsymbol\eta
}
\newcommand\bphi{%
\boldsymbol\phi
}
\newcommand\bsigma{%
\boldsymbol\sigma
}
\newcommand\Gal{%
\text{Gal}
}
\newcommand\Cen{%
\text{Cen}\hspace{0.05cm}
}
\newcommand\Out{%
\text{Out}\hspace{0.03cm}
}
\newcommand\di{%
\text{di}
}

\newcommand\ini{%
\text{in}
}
\newcommand\mo{%
\text{mo}
}
\newcommand\ord{%
\text{ord}
}
\newcommand\spx{%
\text{sp}
}
\newcommand\st{%
\text{st}
}
\newcommand\ab{%
\text{ab}
}
\newcommand\cy{%
\text{cy}
}
\newcommand\car{%
\text{car}
}
\newcommand\ev{%
\text{ev}
}
\newcommand\si{%
\text{si}
}
\newcommand\dis{%
\text{dis}
}
\newcommand\dev{%
\text{dev}
}
\newcommand\pr{%
\text{pr}
}
\newcommand\pro{%
\text{pro}\hspace{0.05 cm}
}
\newcommand\prox{%
\text{pro}
}
\newcommand\proc{%
\text{pro$_c$}\hspace{0.05 cm}
}
\newcommand\prop{%
\text{pro$_p$ }
}
\newcommand\Pur{%
\text{Pur}
}
\newcommand\rel{%
\hspace{0.03cm} \text{rel}\hspace{0.05cm}
}
\newcommand\fr{%
\text{fr}\hspace{0.05cm}
}
\newcommand\ext{%
\text{ext}\hspace{0.03cm}
}
\newcommand\ST{%
\text{ST}%\hspace{0.03cm}
}
\newcommand\ad{%
\text{ad}
}
\newcommand\gen{%
\text{gen }
}

\newcommand\norx{%
\text{nor}
}
\newcommand\diam{%
\text{diam}\hspace{0.03 cm}
}
\newcommand\Acy{%
\text{Acy}
}
\newcommand\crg{%
\text{cr}
}

\newcommand\nil{%
\text{nil}\hspace{0.05 cm}
}
\newcommand\cl{%
\text{cl}
}
\newcommand\itr{%
\text{int}
}
\newcommand\itrx{%
\text{int }
}
\newcommand\itry{%
\text{int}\hspace{0.05cm}
}
\newcommand\thx{%
\text{th}
}
\newcommand\For{%
\text{For}
}
\newcommand\wt{%
\text{wt}\hspace{0.05cm}
}
\newcommand\pow{%
\text{pow}\hspace{0.05cm}
}
\newcommand\rank{%
\text{rank}\hspace{0.05cm}
}
\newcommand\gr{%
\text{gr}
}
\newcommand\grd{%
\text{grd}
}
\newcommand\mic{%
\text{mic}
}
\newcommand\osc{%
\text{osc}
}
\newcommand\pa{%
\text{pa}
}
\newcommand\tf{%
\text{tf}
}
\newcommand\we{%
\text{weak equivalence}
}
\newcommand\wes{%
\text{weak equivalences}
}
\newcommand\whe{%
\text{weak homotopy equivalence}
}

\newcommand\smc{%
\text{simplicial model category}
}
\newcommand\cd{%
\text{commutative diagram}
}
\newcommand\mc{%
\text{model category}
}

\newcommand\OP{%
\text{OP}
}
\newcommand\ner{%
\text{ner}\hspace{0.05cm}
}
\newcommand\nersub{%
\text{ner}
}
\newcommand\colim{%
\text{colim}
}
\newcommand\colimx{%
\text{colim}\hspace{0.05 cm}
}
\newcommand\Aut{%
\text{Aut}\hspace{0.05 cm}
}
\newcommand\Mor{%
\text{Mor$_{}^{^{}}$}\xspace
}
\newcommand\Nat{%
\text{Nat}
}
\newcommand\map{%
\text{map}
}
\newcommand\HOM{%
\text{HOM}
}
\newcommand\bHOM{%
\textbf{HOM}
}
\newcommand\sHOM{%
\sH\sO\sM
}
\newcommand\Hom{%
\text{Hom}
}

\newcommand\shom{%
\texttt{HOM}
}
\newcommand\BV{%
\text{BV}
}

\newcommand\spt{%
\text{spt} \hspace{0.03cm}
}
\newcommand\sptx{%
\text{spt }
}
\newcommand\quadx{%
\hspace*{1em}
}
\newcommand\bbox{%
\ \Box \ 
}
\newcommand\bboxsub{%
\hspace{0.05cm} \Box\xspace 
}
\newcommand\xbox{%
\hspace{0.05cm} \Box \hspace{0.05cm} 
}
\newcommand\ra{%
\rightarrow
}
\newcommand\lra{%
\longrightarrow
}
\newcommand\lla{%
\longleftarrow
}
\newcommand\rat{%
\rightarrowtail
} 
\newcommand\thra{%
\twoheadrightarrow
}
\newcommand\la{%
\leftarrow
}
\newcommand\ds{%
\displaystyle
}
\newcommand\bs{%
\backslash
}
\newcommand\un[1]{%
\underline{#1}\xspace
}
\newcommand\ov[1]{%
\overline{#1}
}
\newcommand\ohc{%
\overline{\text{hocolim}}\hspace{0.05cm}
}
\newcommand\Tee{%
\mathsf{T}
}
\newcommand\dn{%
\Delta[n]
}
\newcommand\dm{%
\Delta[m]
}
\newcommand\restr[2]{%
{#1}|{#2}
}
\newcommand\absx[1]{%
|{#1}|
}
\newcommand\ddn{%
\dot\Delta[n]
}
\newcommand\ddpn{%
{\dot\Delta}^n
}

\newcommand\dpn{%
\Delta^n
}
\newcommand\dpm{%
\Delta^m
}
\newcommand\dpz{%
\Delta^0
}
\newcommand\dpo{%
\Delta^1
}

\newcommand\ovdm{%
\overset{   \raisebox{0.05cm}{{\un{\hspace{0.3cm}}}}   \hspace{0.2cm}     }{\Delta^m}
}
\newcommand\ovdn{%
\overset{   \raisebox{0.05cm}{{\un{\hspace{0.3cm}}}}   \hspace{0.2cm}     }{\Delta^n}
}
\newcommand\ovdalpha{%
\overset{   \raisebox{0.05cm}{{\un{\hspace{0.3cm}}}}   \hspace{0.2cm}     }{\Delta^\alpha}
}

\newcommand\mdpn{%
\overset{\circ \hspace{0.5em}}{\Delta^n}
}
\newcommand\TO{%
\textbf{T}_{\mathcal{O}}
}
\newcommand\TOh{%
\textbf{T}_{\widehat{\mathcal{O}}}
}

\newcommand\pOhp{%
\textbf{ps}\hspace{0.05cm} \widehat{\sO}\text{-}\textbf{SP}
}
\newcommand\psg{%
\textbf{ps}\hspace{0.05cm} \boldsymbol\bGamma\text{-}\textbf{SP}
}

\newcommand\dk{%
\Delta[k]
}
\newcommand\ddk{%
\dot\Delta[k]
}
\newcommand\dz{%
\Delta[0]
}
\newcommand\ddz{%
\dot\Delta[0]
}
\newcommand\Dq{%
\Delta^?
}

\newcommand\dw{%
\Delta[1]
}
\newcommand\ddw{%
\dot\Delta[1]
}

\newcommand\aq{%
\abs{?}
}

\newcommand\aX{%
\abs{X}
}
\newcommand\aY{%
\abs{Y}
}

\newcommand\dsp{%
$\Delta$-separated\xspace
}
\newcommand\dcf{%
$\Delta$-cofibered\xspace
}

\newcommand\ps{%
proper special
}
\newcommand\Oinf{%
\Omega^\infty
}
\newcommand\Sinf{%
\Sigma^\infty
}
\newcommand\Ohs{%
$\widehat{\sO}$-space
}
\newcommand\Oh{%
\widehat{\sO}
}

\newcommand\dso{%
$\Delta$-separated $\sO$-spaces
}

\newcommand\dsep{%
$\Delta$-separated\xspace
}
\newcommand\mA{%
$A$\xspace
}
\newcommand\mB{%
$B$\xspace
}
\newcommand\mC{%
$C$\xspace
}
\newcommand\mD{%
$D$\xspace
}
\newcommand\mE{%
$E$\xspace
}
\newcommand\mF{%
$F$\xspace
}
\newcommand\mG{%
$G$\xspace
}
\newcommand\mH{%
$H$\xspace
}
\newcommand\mI{%
$I$\xspace
}

\newcommand\mK{%
$K$\xspace
}
\newcommand\mL{%
$L$\xspace
}
\newcommand\mM{%
$M$\xspace
}
\newcommand\mN{%
$N$\xspace
}
\newcommand\mO{%
$O$\xspace
}
\newcommand\mP{%
$P$\xspace
}
\newcommand\mQ{%
$Q$\xspace
}
\newcommand\mR{%
$R$\xspace
}
\newcommand\mS{%
$S$\xspace
}
\newcommand\mT{%
$T$\xspace
}
\newcommand\mU{%
$U$\xspace
}
\newcommand\mV{%
$V$\xspace
}
\newcommand\mW{%
$W$\xspace
}
\newcommand\mX{%
$X$\xspace
}
\newcommand\mY{%
$Y$\xspace
}
\newcommand\mZ{%
$Z$\xspace
}
\newcommand\pp{%
^{\prime\prime}
}
\newcommand\cH{%
\widecheck{H}
}
\newcommand\cX{%
\widecheck{X}
}
\newcommand\cY{%
\widecheck{Y}
}
\newcommand\cZ{%
\widecheck{Z}
}
\newcommand\as[1]{%
\abs{\sin{#1}}
}
\newcommand\hsx{%
\hspace{0.05cm}
}
\newcommand\hsy{%
\hspace{0.03cm}
}
\newcommand\hthree{%
\hspace{0.03cm}
}

\newcommand\dnd{%
\hspace{.1cm} | \hspace{-.375cm} \not \hspace{0.1cm}
}
\newcommand\hspnx{%
\hspace{-.25cm}
}
\newcommand\vspi{%
\\ \indent
}
\makeatletter
\DeclareRobustCommand\widecheck[1]{{\mathpalette\@widecheck{#1}}}
\def\@widecheck#1#2{%
    \setbox\z@\hbox{\m@th$#1#2$}%
    \setbox\tw@\hbox{\m@th$#1%
       \widehat{%
          \vrule\@width\z@\@height\ht\z@
          \vrule\@height\z@\@width\wd\z@}$}%
    \dp\tw@-\ht\z@
    \@tempdima\ht\z@ \advance\@tempdima2\ht\tw@ \divide\@tempdima\thr@@
    \setbox\tw@\hbox{%
       \raise\@tempdima\hbox{\scalebox{1}[-1]{\lower\@tempdima\box
\tw@}}}%
    {\ooalign{\box\tw@ \cr \box\z@}}}
\makeatother

\newtheoremstyle{propx}% name of the style to be used
  {4pt}% measure of space to leave above the theorem. E.g.: 3pt
  {0pt}% measure of space to leave below the theorem. E.g.: 3pt
  {\upshape}% name of font to use in the body of the theorem
  {20pt}% measure of space to indent
  {\bfseries}% name of head font
  {}% punctuation between head and body
  { }% space after theorem head; " " = normal interword space \footnotesize
  {}

\theoremstyle{propx}

\newtheorem{proposition}{\small PROPOSITION}

\theoremstyle{definition}
\newtheorem{x}{}[subsubsection]

\newcommand{\norm}[1]{\left\lVert #1 \right\rVert}

\def\Therefore{\boldsymbol{\text{ }
\leavevmode
\lower0.4ex\hbox{$\cdot$}
\kern-.5em\raise0.7ex\hbox{$\cdot$}
\kern-0.55em\lower0.4ex\hbox{$\cdot$}
\thinspace\text{ }}}

\title{\textbf{Topics in Topology and Homotopy Theory}}
\author{Garth Warner\\
Department of Mathematics\\
University of Washington}
\date{}									% Activate to display a given date or no date
\titleformat{\chapter}[display]
{\normalfont\filcenter\huge\bfseries}{}{0pt}{\large}
\titleformat{\chapter}[display]
{\normalfont\filcenter\huge\bfseries}{}{0pt}{\large}
\usepackage[OT2,OT1]{fontenc} 
\newcommand\cyr
{
\renewcommand\rmdefault{wncyr} 
\renewcommand\sfdefault{wncyss} 
\renewcommand\encodingdefault{OT2} 
\normalfont
\selectfont
}

\DeclareTextFontCommand{\textcyr}{\cyr}

\begin{document}

\maketitle                              % Print title page.
%\tableofcontents                        % Print table of contents

%\include{_Preface}
%\addcontentsline{toc}{chapter}{\protect\numberline{}PREFACE}

%\include{_Notation}
%\addcontentsline{toc}{chapter}{\protect\numberline{}NOTATION}

%\mainmatter                             % only in book class (arabic page #s)

%\maketitle
%\tableofcontents
%\chapter{Preimages}
%\label{chapter-preimages}

%In this and the next chapter we will turn our attention to continuum
%spaces, primarily spaces such as $\Re^n$.  The same ideas that

%\section{Preimage Planning}

%\begin{document}
%\maketitle
%\tableofcontents
%\setlength{\parindent}{0pt}
%%%dmc\flushleft
\titlespacing*{\chapter}{0pt}{-50pt}{40pt}
\setlength{\parskip}{0.1em}

\begingroup%%----------------------------------->>
\fontsize{11pt}{11pt}\selectfont

\qquad\qquad TOPICS IN TOPOLOGY AND HOMOTOPY THEORY\\
 %%dmc00A_\vspace{2.cm}

\qquad $\S0.\ $ \qquad CATEGORIES AND FUNCTORS\\
%%dmc00A_\vspace{0.3cm}

\qquad $\S1.\ $ \qquad COMPLETELY REGULAR HAUSDORFF SPACES\\
%%dmc00A_\vspace{0.3cm}

\qquad $\S2.\ $ \qquad CONTINUOUS FUNCTIONS\\
%%dmc00A_\vspace{0.3cm}

\qquad $\S3.\ $ \qquad COFIBRATIONS\\
%%dmc00A_\vspace{0.3cm}

\qquad $\S4.\ $ \qquad FIBRATIONS\\
%%dmc00A_\vspace{0.3cm}

\qquad $\S5.\ $ \qquad VERTEX SCHEMES AND CW COMPLEXES\\
%%dmc00A_\vspace{0.3cm}

\qquad $\S6.\ $ \qquad ABSOLUTE NEIGHBORHOOD RETRACTS\\
%%dmc00A_\vspace{0.3cm}

\qquad $\S7.\ $ \qquad $\sC$-THEORY\\
%%dmc00A_\vspace{0.3cm}

\qquad $\S8.\ $ \qquad LOCALIZATION OF GROUPS\\
%%dmc00A_\vspace{0.3cm}

\qquad $\S9.\ $ \qquad HOMOTOPICAL LOCALIZATION\\
%%dmc00A_\vspace{0.3cm}

\qquad $\S10.$ \qquad COMPLETION OF GROUPS\\
%%dmc00A_\vspace{0.3cm}

\qquad $\S11.$ \qquad HOMOTOPICAL COMPLETION\\
%%dmc00A_\vspace{0.3cm}

\qquad $\S12.$ \qquad MODEL CATEGORIES\\
%%dmc00A_\vspace{0.3cm}

\qquad $\S13.$ \qquad SIMPLICIAL SETS\\
%%dmc00A_\vspace{0.3cm}

\qquad $\S14.$ \qquad SIMPLICIAL SPACES\\
%%dmc00A_\vspace{0.3cm}

\qquad $\S15.$ \qquad TRIANGULATED CATEGORIES\\
%%dmc00A_\vspace{0.3cm}

\qquad $\S16.$ \qquad SPECTRA\\
%%dmc00A_\vspace{0.3cm}

\qquad $\S17.$ \qquad STABLE HOMOTOPY THEORY\\
%%dmc00A_\vspace{0.3cm}

\qquad $\S18.$ \qquad ALGEBRAIC K-THEORY\\
%%dmc00A_\vspace{0.3cm}

\qquad $\S19.$ \qquad DIMENSION THEORY\\
%%dmc00A_\vspace{0.3cm}

\qquad $\S20.$ \qquad COHOMOLOGICAL DIMENSION THEORY\\
%%dmc00A_\vspace{0.3cm}

%%dmc00A_\vspace{0.3cm}
\endgroup %%------------------------------------<<

\chapter*{PREFACE}

\setlength\parindent{2em}
%-------------------------
$\text{ }$\\[-1.25cm]

This book is addressed to those readers who have been through 
Rotman\footnote[2]{\textit{An Introduction to Algebraic Topology}, Springer Verlag (1988).} 
(or its equivalent), possess a wellthumbed copy of 
Spanier\footnote[3]{\textit{Algebraic Topology}, Springer Verlag (1989).}, 
and have a good background in algebra and general topology.

Granted these prerequisites, my intention is to provide at the core a state of the art treatment of the homotopical foundations of algebraic topology.  
The depth of coverage is substantial and I have made a point to include material which is ordinarily not included, for instance, an account of algebraic K$-$theory in the sense of Waldhausen.  
There is also a systematic treatment of ANR theory (but, reluctantly, the connections with modern geometric topology have been omitted).  
However, truly advanced topics are not considered (e.g., equivariant stable homotopy theory, surgery, infinite dimensional topology, etale K$-$theory, $\ldots$).  
Still, one should not get the impression that what remains is easy:  There are numerous difficult technical results that have to be brought to heel.

Instead of laying out a synopsis of each chapter, here is a sample of some of what is taken up.\\

\begin{tabular}{cll}
\indent (1)& \  &Nilpotency and its role in homotopy theory.\\
\indent (2)& \  &Bousfield's theory of the localization of spaces and spectra.\\
\indent (3)& \  &Homotopy limits and colimits and their applications.\\
\indent (4)& \  &The James construction, symmetric products, and the Dold$-$Thom theorem.\\
\indent (5)& \  &Brown and Adams representability in the setting of triangulated categories.\\
\indent (6)& \  &Operads and the May$-$Thomason theorem on the uniqueness of infinite\\
&&loop space machines.\\
\indent (7)& \  &The plus construction and theorems A and B of Quillen.\\
\indent (8)& \  &Hopkins' global picture of stable homotopy theory.\\
\indent (9)& \  &Model categories, cofibration categories, and Waldhausen categories.\\
\indent (10)& \  &The Dugundji extension theorem and its consequences.\\
\end{tabular}
\vspace{0.25cm}

A book of this type is not meant to be read linearly.  
For example, a reader wishing to study stable homotopy theory could start by perusing $\S12$ and $\S15$ and then proceed to $\S16$ and $\S17$ or a reader who wants to learn the theory of dimension could immediately turn to $\S19$ and $\S20$.  One could also base a second year course in algebraic topology on $\S3 - \S11$.  
Many other combinations are possible.

Structurally, each $\S$ has its own set of references (both books and articles).  
No attempt has been made to append remarks of a historical nature but for this, the reader can do no better than turn to 
Dieudonn\'e\footnote[2]{\textit{A History of Algebraic and Differential Topology 1900-1960}, Birkh\"auser (1989); see also, Adams, Proc. Sympos. Pure Math. 22 (1971), 1-22 and Whitehead, \textit{Bull. Amer. Math. Soc.} 8 (1963), 1-29.}.  
Finally, numerous exercises and problems (in the form of ``examples'' and ``facts'') are scattered throughout the text, most with partial or complete solutions.\\

\chapter*{PREFACE (bis)}
\setlength\parindent{2em}

\ \indent This project which started almost thirty years ago has for various reasons remained dormant now for almost twenty-five years.  At the time that this book was finished, it was very much up to date but, of course, since then there have been a number of developments which are not included.  
Still, there is a lot of material to be covered and the numerous detailed examples are a feature which sets it apart from other accounts.\\

\indent \un{N.B.} \ As regards model category theory, the author has written a greatly expanded exposition, Categorical Homotopy Theory, which does include more recent material and can be found at 
\url{https://sites.math.washington.edu/~warner/CHT_Warner.pdf}.
\\[1cm]

\[
\textbf{ACKNOWLEDGMENTS}
\]
\\

\indent\indent\textbullet \ My thanks to Mary Sheetz, who typed my original hand written manuscript.
\\

\indent\indent\textbullet \ My thanks to David Clark, who undertook the heroic task of converting the original manuscript, 
which was formatted in a now obsolete ``language'', to AMS-TeX.
\\

\indent\indent\textbullet \ My thanks to Judith Clare for her meticulous job of proofreading. 
%\indent ACKNOWLEDEMENT  David Clark undertook the heroic task of converting the original manuscript, which was formatted in a now obsolete ``language'', to AMS-TeX.

\chapter*{NOTATION}

\begin{tabular}{cll}
\qquad (1)& \quad &$\N$, the positive integers;\\
&&$\Z$, the integers;\\
&&$\Q$, the rational numbers;\\
%&&$\textbf{P}$, the irrational numbers;  
&&$\PP$, the irrational numbers; \\
&&$\R$, the real numbers;\\
&&$\C$, the complex numbers;\\
&&\textbf{H}, the quaternions;\\
&&$\bPi$, the prime numbers.\\
\qquad (2)& \quad &$\R^n = \R \times \cdots \times \R$ (n factors);\\
&&\textbf{D}$^n = \{x \in \R^n: \norm{x} \leq 1\}$;\\
&&\textbf{B}$^n = \{x \in \R^n: \norm{x} < 1\}$;\\
&&\textbf{S}$^{n-1} = \{x \in \R^n: \norm{x} = 1\}$;\\
&&\textbf{T}$^n = S^1 \times \cdots \times S^1$ ($n$ factors).\\
\qquad (3)& \quad &$\Delta^n = \{x \in \R^{n+1}: \ds\sum\limits_i x_i = 1 \ \& \  \forall \ i, x_i \geq 0\}$;\\
%&&$\mathring{\Delta}^n = \{x \in \R^{n+1}: \sum\limits_i x_i = 1 \ \& \ \forall \ i, x_i > 0\}$;\\
&&$\mdpn = \{x \in \R^{n+1}: \ds\sum\limits_i x_i = 1 \ \& \ \forall \ i, x_i > 0\}$;\\
&&$\dot{\Delta}^n = \{x \in \R^{n+1}: \ds\sum\limits_i x_i = 1 \ \forall \ i, x_i = 0\}$;\\
\qquad (4)& \quad &$\omega = $ first infinite ordinal; $\Omega =$ first uncountable ordinal.\\
\qquad (5)& \quad &cl = closure, fr = frontier, wt = weight, int = interior, osc = oscillation.\\
\qquad (6)& \quad &Given a set $\mathcal{S}$, $\chi_{_\mathcal{S}}$ is the characteristic function of $\mathcal{S}$ and $\#(\mathcal{S})$ is the\\
&&cardinality of $\mathcal{S}$.\\
\qquad (7)& \quad &Given a topological space $X$, $C(X)$ is the set of real valued continuous\\
&&functions on $X$ and $BC(X)$ is the set of real valued bounded continuous functions\\
&&on $X$.\\
\qquad (8)& \quad &Given a topological space $X$, $X_\infty$ is the one point compactification of $X$.\\
\qquad (9)& \quad &Given a completely regular Hausdorff space $X$, $\beta X$ is the Stone-Cech\\
&&compactification of $X$.\\
\qquad (10)& \quad &Given a completely regular Hausdorff space $X$, $\nu X$ is the $\R$-compactification\\
&&of $X$.\\
\end{tabular}

\newpage 
.

\pagenumbering{bychapter}
%\chapter{}
\newpage
\setcounter{page}{1}
\renewcommand{\thepage}{0-\arabic{page}}
\chapter{
$\boldsymbol{\S}$\textbf{0}.\quadx  CATEGORIES AND FUNCTORS}
\setlength\parindent{2em}
\setcounter{proposition}{0}
%%----------------------------------------------------------------------------------------------01
$\text{ }$\\[-1.25cm]

In addition to establishing notation and fixing terminology, background material from the theory relevant to the work as a whole is collected below and will be referred to as the need arises.

Given a category \bC, denote by $\Ob\bC$ its class of objects and by $\Mor\bC$ its class of morphisms.  
If $X, Y \in \Ob\bC$ is an ordered pair of objects, then $\Mor(X,Y)$ is the set of morphisms (or arrows) from $X$ to $Y$.  
An element $f \in \Mor(X,Y)$ is said to have 
\un{domain}
\index{domain} 
$X$ and 
\un{codomain}
\index{codomain}  
$Y$.  
One writes $f:X \ra Y$ or $X \overset{f}{\ra} Y$.  
Functors preserve the arrows, while cofunctors reverse the arrows, i.e., a cofunctor is a functor on $\bC^\OP$, the category opposite to \bC.

Here is a list of frequently occurring categories.

\indent\indent (1) \bSET, 
\index{\bSET} 
\index{category: (\bSET)}
the category of sets and 
$\bSET_*$, 
\index{$\bSET_*$} 
\index{category: ($\bSET_*$)}
the category of pointed sets.  
If $X$, $Y$ $\in \Ob\bSET$, then 
$\Mor(X,Y) = F(X,Y)$, the functions from $X$ to $Y$, and if 
$(X,x_0)$, $(Y,y_0)$ $\in \Ob\bSET_*$, then 
$\Mor((X,x_0),(Y,y_0)) = F(X,x_0;Y,y_0)$, the base point preserving functions from $X$ to $Y$.

\indent\indent (2) \bTOP, 
\index{\bTOP} 
\index{category: (\bTOP)}
the category of topological spaces, and 
$\bTOP_*$, 
\index{$\bTOP_*$} 
\index{category: ($\bTOP_*$)}
the category of pointed topological spaces.  
If $X$, $Y$ $\in \Ob\bTOP$, then 
$\Mor(X,Y) = C(X,Y)$, the continuous functions from $X$ to $Y$, and if 
$(X,x_0)$, $(Y,y_0)$ $\in \Ob\bTOP_*$, then 
$\Mor((X,x_0),(Y,y_0)) = C(X,x_0;Y,y_0)$, 
the base point preserving continuous functions from $X$ to $Y$.

\indent\indent (3) $\bSET^2$, 
\index{$\bSET^2$} 
\index{category: ($\bSET^2$)}
the category of pairs of sets, and 
$\bSET_*^2$, 
\index{$\bSET_*^2$} 
\index{category: ($\bSET_*^2$)}
the category of pointed pairs of sets.
If $(X,A)$ $(Y,B) \in \Ob \bSET^2$, then $\Mor((X,A),(Y,B)) = F(X,A;Y,B)$, 
the functions from $X$ to $Y$ that take $A$ to $B$, and if 
$(X,A,x_0)$ $(Y,B,y_0) \in \Ob \bSET_*^2$,  then 
$\Mor((X,A,x_0),(Y,B,y_0)) = F(X,A,x_0;Y,B,y_0)$, the base point preserving functions from $X$ to $Y$ that take $A$ to $B$.

\indent\indent (4) $\bTOP^2$, 
\index{$\bTOP^2$,} 
\index{category: ($\bTOP^2$,)}
the category of pairs of topological spaces, and 
$\bTOP_*^2$, 
\index{$\bTOP_*^2$,} 
\index{category: ($\bTOP_*^2$,)}
the category of pointed pairs of topological spaces.  If $(X,A)$ $(Y,B) \in \Ob \bTOP^2$, then 
$\Mor((X,A),(Y,B)) =$ $C(X,A;Y,B)$, the continuous functions from $X$ to $Y$ that take $A$ to $B$, and if 
$(X,A,x_0)$ $(Y,B,y_0) \in \Ob \bTOP_*^2$,  then 
$\Mor((X,A,x_0),(Y,B,y_0)) =$ $C(X,A,x_0;Y,B,y_0)$, 
the base point preserving continuous functions from $X$ to $Y$ that take $A$ to $B$.

\indent\indent (5) \bHTOP, 
\index{\bHTOP} 
\index{category: (\bHTOP)}
the homotopy category of topological spaces, and 
$\bHTOP_*$, 
\index{$\bHTOP_*$} 
\index{category: ($\bHTOP_*$)} 
the homotopy category of pointed topological spaces.   
If \ $X$, $Y$ $\in \Ob\bHTOP$,  then 
$\Mor(X,Y) = [X,Y]$, the homotopy classes in \ $C(X,Y)$ 
and if \ $(X,x_0)$, $(Y,y_0)$ $\in \Ob\bHTOP_*$, then \ 
$\Mor((X,x_0),Y,y_0)) = [(X,x_0);(Y,y_0)]$, 
the homotopy classes in $C(X,x_0;Y,y_0)$.

\indent\indent (6) $\bHTOP^2$, 
\index{$\bHTOP^2$} 
\index{category: ($\bHTOP^2$)}
the homotopy category of pairs of topological spaces, and
$\bHTOP_*^2$, 
\index{$\bHTOP_*^2$} 
\index{category: ($\bHTOP_*^2$)}
the homotopy category of pointed pairs of topological spaces.  
If $(X,A)$, $(Y,B)$ $\in \Ob \bHTOP^2$, 
%%----------------------------------------------------------------------------------------------02
then $\Mor((X,A),(Y,B)) =$ $[X,A;Y,B]$, the homotopy classes in $C(X,A;Y,B)$ and if 
$(X,A,x_0)$, $(Y,B,y_0) \in \Ob \bHTOP_*^2$, then 
$\Mor((X,A,x_0),(Y,B,y_0)) =$ $[X,A,x_0;Y,B,y_0]$, the homotopy classes in $C(X,A,x_0;Y,B,y_0)$.

\indent\indent (7) \bHAUS, 
\index{\bHAUS} 
\index{category: (\bHAUS)}
the full subcategory of \bTOP whose objects are the Hausdorff spaces and 
\CPTHAUS, 
\index{\CPTHAUS} 
\index{category: (\CPTHAUS)} 
the full subcategory of \bHAUS whose objects are the compact spaces.

\indent\indent (8) $\Pi X$, 
\index{$\bPi X$} 
\index{category: ($\bPi X$)}
the fundamental groupoid of a topological space $X$.

\indent\indent (9) \bGR, \bAB, \bRG (\bAMOD) or (\bMODA), 
\index{category: \bGR} 
\index{category: \bAB} 
\index{category: \bRG} 
\index{category: \bAMOD} 
\index{category: \bMODA}
the category of groups, abelian groups, rings with unit (left or right $A$-modules, $A \in \Ob\bRG$).

\indent\indent (10) \bzero,  
\index{category: (\bzero)}
the category with no objects and no arrows.  
\bone, 
\index{category: (\bone)}
the category with one object and one arrow.  
\textbf{2}, 
\index{category: (\textbf{2})}
the category with two objects and one arrow not the identity.

A category is said to be 
\un{discrete}
\index{discrete category} 
if all its morphisms are identities.  
Every class is the class of objects of a discrete category.

[Note: \ A category is 
\un{small}
\index{small category}
\index{category: small}
if its class of objects is a set; otherwise it is 
\un{large}.
\index{large category}
\index{category: large}   
A category is 
\un{finite (countable)}
\index{finite category}
\index{category: finite}
\index{countable category}
\index{category: countable}
if its class of morphisms is a finite (countable) set.]\\

\begingroup%%------------------------------------>>
\fontsize{9pt}{11pt}\selectfont
In this book, the foundation for category theory is the ``one universe'' approach taken by Herrlich-Strecker and Osborne referenced at the end of the $\S$).  The key words are ``set'', ``class'', and ``conglomerate''.  
Thus the issue is not only one of size but also membership (every set is a class and every class is a conglomerate).  
Example: $\{\Ob\bSET\}$ is a conglomerate, not a class (the members of a class are sets).
\\ \indent
[Note: \ A functor $F:\bC \ra \textbf{D}$ is a function from $\Mor \bC$ to $\Mor \bD$ that preserves identities and respects composition.  
In particular: $F$ is a class, hence $\{F\}$ is a conglomerate.]
\\ \indent
A 
\un{metacategory}
\index{metacategory} 
is defined in the same way as a category except that the objects and the morphisms are allowed to be conglomerates and the requirement that the conglomerate of morphisms between two objects be a set is dropped.  
While there are exceptions, most categorical concepts have metacategorical analogs or interpretations.  
Example:  The ``category of categories'' is a metacategory.
\\ \indent
[Note: \ Every category is a metacategory.  
On the other hand, it can happen that a metacategory is isomorphic to a category but is not itself a category.  
Still, the convention is to overlook this technical nicety and treat such a metacategory as a category.]\\
\endgroup%%------------------------------------<<

\vspace{0.25cm}

Given categories \bA, \bB, \bC and functors
$
\begin{cases}
\ T:\textbf{A} \ra \bC\\
\ S:\textbf{B} \ra \bC
\end{cases}
, \ 
$
the 
\un{comma category}
\index{comma category}  
$\abs{T,S}$ is the category whose objects are triples 
$
(X,f,Y): 
\begin{cases}
\ X \in \Ob \bA\\
\ Y \in \Ob \bB
\end{cases}
\& f \in \Mor(TX,SY)
$
and whose morphisms $(X,f,Y) \ra (X^\prime,f^\prime,Y^\prime)$ are the pairs
$
(\phi,\psi):
\begin{cases}
\ \phi \in \Mor(X,X^\prime)\\
\ \psi \in \Mor(Y,Y^\prime)
\end{cases}
$
for
%%----------------------------------------------------------------------------------------------03
which the square 
\begin{tikzcd}[ sep=large]
{TX} \ar{d}[swap]{T\phi} \ar{r}{f} &{SY} \ar{d}{S\psi}\\
{TX^\prime} \ar{r}[swap]{f^\prime} &{SY^\prime}
\end{tikzcd}
commutes.  
Composition is defined componentwise and the identity attached to $(X,f,Y)$ is $(\id_X,\id_Y)$.

\indent\indent $(A\backslash\bC)$  
\index{$A\backslash\bC$}
Let $A \in \Ob\bC$ and write $K_A$ for the constant functor 
$\bone \ra \bC$ with value $A$ $-$then 
$A\backslash\bC \equiv \abs{K_A,\id_{\bC}}$ is the category of 
\un{objects under $A$}
\index{category of objects under $A$}.

\indent\indent $(\bC/B)$   
\index{$\bC/B$}
Let $B \in \Ob\bC$ and write $K_B$ for the constant functor 
$\bone \ra \bC$ with value $B$ $-$then 
$\bC/B \equiv \abs{\id_{\bC}, K_B}$ is the category of
\un{objects over $B$}
\index{category of objects over $B$}.

Putting together $A\backslash\bC$ $\&$ $\bC/B$ leads to the category of 
\un{objects under $A$ and over \mB}:
\index{objects under $A$ and over $B$}
$A\backslash\bC/B$.
\index{$A\backslash\bC/B$}  
The notation is incomplete since it fails to reflect the choice of the structural morphism $A \ra B$.  
Examples: 
(1) \ $\emptyset\backslash\bTOP/* = \bTOP$; 
(2) \ $*\backslash \bTOP/* = \bTOP_*$;
(3) \ $A\backslash\bTOP/* = A\backslash\bTOP$;
(4) \ $\emptyset \backslash \bTOP/B = \bTOP/B$; 
(5) \ $B\backslash\bTOP/B = \bTOP(B)$, the ``exspaces'' of James (with structural morphism $\id_B$).

\label{5.0al}
\label{18.3}
The 
\un{arrow category}
\index{} $\bC(\ra)$ 
of \bC is the comma category $\abs{\id_{\bC},\id_{\bC}}$.  
Examples: 
(1) \ $\bTOP^2$ is a subcategory of  $\bTOP(\ra)$; 
(2) \  $\bTOP_*^2$ is a subcategory of $\bTOP_*(\ra)$.

[Note: \ There are obvious notions of homotopy in $\bTOP(\ra)$ or $\bTOP_*(\ra)$,  from which  $\bHTOP(\ra)$ or  $\bHTOP_*(\ra)$.] \\

\begingroup%%------------------------------------>>%%----------------------------------->>
\fontsize{9pt}{11pt}\selectfont
The comma category $\abs{K_A,K_B}$ is $\Mor(A,B)$ viewed as a discrete category.\\
\endgroup%%------------------------------------<< %%------------------------------------<<

A morphism $f:X \ra Y$ in a category \bC is said to be an 
\un{isomorphism}
\index{isomorphism} 
if there exists a morphism $g:Y \ra X$ such that $g \circ f = \id_X$ and $f \circ g = \id_Y$.  
If $g$ exists, then $g$ is unique.  
It is called the
\un{inverse}
\index{inverse (morphism)} 
of $f$ and is denoted $f^{-1}$.  
Objects $X, \ Y \in \Ob\bC$ are said to be 
\un{isomorphic}, 
\index{isomorphic (objects)}
written $X \approx Y$, provided that there is an isomorphism $f:X \ra Y$.  The relation 
``isomorphic to'' is an equivalence relation on $\Ob\bC$.\\

\begingroup%%------------------------------------>>%%----------------------------------->>
\fontsize{9pt}{11pt}\selectfont

The isomorphisms in \bSET are the bijective maps, in \bTOP the homeomorphisms, in \bHTOP the homotopy equivalences.  The isomorphisms in any full subcategory of \bTOP are the homeomorphisms.\\

\endgroup%%------------------------------------<< %%------------------------------------<<

Let 
$
\begin{cases}
\ F:\bC \ra \bD\\
\ G:\bC \ra \bD
\end{cases}
$
be functors $-$then a 
\un{natural transformation}
\index{natural transformation} 
$\Xi$ from \mF to \mG is a function that assigns to each $X \in \Ob\bC$ an element $\Xi_X \in \Mor(FX,GX)$ such that for every $f \in \Mor(X,Y)$ the square 
\begin{tikzcd}[ sep=large]
{FX} \ar{d}[swap]{F f} \ar{r}{\Xi_X} &{GX} \ar{d}{G f}\\
{FY} \ar{r}[swap]{\Xi_Y} &{GY}
\end{tikzcd}
commutes, $\Xi$ being termed a 
\un{natural isomorphism}
\index{natural isomorphism} 
if all the $\Xi_X$ are isomorphisms, in which case \mF and \mG are said to be 
\un{naturally isomorphic}
\index{naturally isomorphic (functors)} 
written, $F \approx G$.\\

%%----------------------------------------------------------------------------------------------04
Given categories 
$
\begin{cases}
\ \bC\\
\ \bD
\end{cases}
, \ 
$
the 
\un{functor category}
\index{functor category} 
$[\bC,\bD]$ is the metacategory whose objects are the functors 
$F:\bC \ra \bD$ and whose morphisms are the natural transformations $\Nat(F,G)$ from \mF to \mG.  
In general, $[\bC,\bD]$ need not be isomorphic to a category, although this will be true if \bC is small.

[Note: \ The isomorphisms in $[\bC,\bD]$ are the natural isomorphisms.]\\

Given categories 
$
\begin{cases}
\ \bC\\
\ \bD
\end{cases}
$
and functors 
$
\begin{cases}
\ K:\bA \ra \bC\\
\ L:\bD \ra \bB
\end{cases}
, \ 
$
there are functors \\
$
\begin{cases}
\ [K,\bD]:[\bC,\bD]  \ra [\bA,\bD]\\
\ [\bC,L]:[\bC,\bD]  \ra [\bC,\bB]
\end{cases}
$
defined by 
$
\begin{cases}
\ \text{precomposition}\\
\ \text{postcomposition}
\end{cases}
. \ 
$ 
If $\Xi \in \Mor([\bC,\bD])$, then we shall write
$
\begin{cases}
\ \Xi K\\
\ L \Xi
\end{cases}
$
in place of 
$
\begin{cases}
\ [K,\bD] \Xi\\
\ [\bC,L] \Xi
\end{cases}
$
, so $L(\Xi K) = (L\Xi)K$.\\

\vspace{0.5cm}

\begingroup%%------------------------------------>>%%----------------------------------->>
\fontsize{9pt}{11pt}\selectfont
There is a simple calculus that governs these operations:\\
\[
\begin{cases}
\ \Xi (K \circ K^\prime) = (\Xi K)K^\prime\\
\ (\Xi^\prime \circ \Xi)K = (\Xi^\prime K) \circ (\Xi K)
\end{cases}
\text{and} \quadx 
\begin{cases}
\ (L^\prime \circ L) \Xi = L^\prime (L \Xi)\\
\ L(\Xi^\prime \circ \Xi) = (L \Xi^\prime) \circ (L \Xi)
\end{cases}
.
\]

\endgroup%%------------------------------------<< %%------------------------------------<<
\vspace{0.25cm}

A functor $F:\bC \ra \bD$ is said to be 
\un{faithful}
\index{faithful (functor)}
(\un{full})
\index{full (functor)}
if for any ordered pair $X, \ Y \in \Ob\bC$, the map $\Mor(X,Y) \ra \Mor(FX,FY)$ is injective (surjective).  
If \mF is full and faithful, then $F$ is 
\un{conservative}
\index{conservative (functor)}, 
i.e., $f$ is an isomorphism iff $F f $ is an isomorphism.\\

\begingroup%%------------------------------------>>%%----------------------------------->>
\fontsize{9pt}{11pt}\selectfont
A category \bC is said to be 
\un{concrete}
\index{concrete (category)} 
if there exists a faithful functor 
$U:\bC \ra \bSET$.  
Example: \bTOP is concrete but \bHTOP is not.

[Note: \ A category is concrete iff it is isomorphic to a subcategory of \bSET.]\\

\endgroup%%------------------------------------<< %%------------------------------------<<

Associated with any object \mX in a category \bC is the 
functor $\Mor(X,-) \in \Ob[\bC,\bSET]$ and the 
cofunctor $\Mor(-,X) \in$ $\Ob[\bC^\OP,\bSET]$.  
If $F \in \Ob[\bC,\bSET]$ is a functor or if $F \in \Ob[\bC^\OP,\bSET]$ is a cofunctor, then the Yoneda lemma establishes a 
bijection $\iota_X$ between $\Nat(\Mor(X,-),F)$ or $\Nat(\Mor(-,X),F)$ and $FX$, viz. $\iota_X(\Xi) = \Xi_X(\id_X)$.  
Therefore the assignments 
$
\begin{cases}
\ X \ra \Mor(X,-)\\
\ X \ra \Mor(-,X)
\end{cases}
$
lead to functors 
$
\begin{cases}
\ \bC^\OP \ra [\bC,\bSET]\\
\ \bC \ra [\bC^\OP,\bSET]
\end{cases}
$
that are full, faithful, and injective on objects, the 
\un{Yoneda embeddings}
\index{Yoneda embeddings}.  
One says that \mF is 
\un{representable}
\index{representable (functor (cofunctor))} 
(by \mX) if \mF is naturally isomorphic to 
$\Mor(X,-)$ or $\Mor(-,X)$.  Representing objects are isomorphic.\\

\begingroup%%------------------------------------>>%%----------------------------------->>
\fontsize{9pt}{11pt}\selectfont
The forgetful functors 
$\bTOP \ra \bSET$,
$\bGR \ra \bSET$, 
$\bRG \ra \bSET$ 
are representable.  
The power set cofunctor $\bSET \ra \bSET$ is representable.\\

\endgroup%%------------------------------------<< %%------------------------------------<<

%%----------------------------------------------------------------------------------------------05
A functor $F:\bC \ra \bD$ is said to be an 
\un{isomorphism}
\index{isomorphism (functor)} 
if there exists a functor 
$G:\bD \ra \bC$ such that $G \circ F = \id_{\bC}$ and $F \circ G = \id_{\bD}$.  
A functor is an isomorphism iff it is full, faithful, and bijective on objects.  Categories \bC and \bD are said to be 
\un{isomorphic}
\index{isomorphic (categories)} 
provided there is an isomorphism $F:\bC \ra \bD$.

[Note: \ An isomorphism between categories is the same as an isomorphism in the ``category of categories''.]\\

\begingroup%%------------------------------------>>%%----------------------------------->>
\fontsize{9pt}{11pt}\selectfont
A full subcategory of \bTOP whose objects are the \mA spaces is isomorphic to the category of ordered sets and order preserving maps 
(reflexive + transitive = \un{order}).

[Note: \ An 
\un{\mA space}
\index{A space} %\index{\mA space}
is a topological space \mX in which the intersection of every collection of open sets is open.  
Each $x \in X$ is contained in a minimal open set $U_x$ and the relation $x \leq y$ iff $x \in U_y$ is an order on \mX.  
On the other hand, if $\leq$ is an order on a set \mX, then \mX becomes an \mA space by taking as a basis the sets 
$U_x = \{y: y \leq x\}$ $(x \in X)$.]\\

\endgroup%%------------------------------------<< %%------------------------------------<<

A functor $F:\bC \ra \bD$ is said to be an 
\un{equivalence}
\index{equivalence (functor)} 
if there exists a functor 
$G:\bD \ra \bC$ such that $G \circ F \approx \id_{\bC}$ and $F \circ G \approx \id_{\bD}$.  
A functor is an equivalence iff it is full, faithful, and  has a 
\un{representative image}, 
\index{representative image (functor)} 
i.e., for any $Y \in \Ob\bD$ there exists an $X \in \Ob\bC$ such that $FX$ is isomorphic to \mY.  
Categories \bC and \bD are said to be 
\un{equivalent}
\index{equivalent (categories)} 
provided there is an equivalence 
$F:\bC \ra \bD$.  The object isomorphism types of equivalent categories are in a one-to-one correspondence.

[Note: \ If \mF and \mG are injective on objects, then \bC and \bD are isomorphic (categorical 
``Schroeder-Bernstein'').]\\

\begingroup%%------------------------------------>>%%----------------------------------->>
\fontsize{9pt}{11pt}\selectfont
The functor from the category of metric spaces and continuous functions to the category of metrizable spaces and continuous functions which assigns to a pair $(X,d)$ the pair $(X,\tau_d)$, $\tau_d$ the topology on \mX determined by $d$, is an equivalence but not an isomorphism.

[Note: \ The category of metric spaces and continuous functions is not a subcategory of \bTOP.]\\

\endgroup%%------------------------------------<< %%------------------------------------<<

A category is 
\un{skeletal} 
\index{skeletal (category)} 
if isomorphic objects are equal.  Given a category \bC, a 
\un{skeleton} 
\index{skeleton (category)} 
of \bC is a full, skeletal subcategory $\ov{\bC}$ for which the inclusion  
$\ov{\bC} \ra \bC$ had a representative image (hence is an equivalence).  
Every category has a skeleton and any two skeletons of a category are isomorphic.  
A category is 
\un{skeletally small}
\index{skeletally small (category)} 
is it has a small skeleton.\\

\begingroup%%------------------------------------>>%%----------------------------------->>
\fontsize{9pt}{11pt}\selectfont
The full subcategory of \bSET whose objects are the cardinal numbers is a skeleton of \bSET.\\

\endgroup%%------------------------------------<< %%------------------------------------<<

A morphism $f:X \ra Y$ in a category \bC is said to be a 
\un{monomorphism}
\index{monomorphism (category)} 
if it is left cancellable with respect to composition, i.e., for any pair of morphisms $u,\  v:Z \ra X$ such that $f \circ u = f \circ v$, there follows $u = v$.

%%----------------------------------------------------------------------------------------------06
A morphism $f:X \ra Y$ in a category \bC is said to be an  
\un{epimorphism}
\index{epimorphism (category)} 
if it is right cancellable with respect to composition, 
i.e., for any pair of morphisms $u,\  v:Y \ra Z$ such that $u \circ f = v \circ f$, there follows $u = v$.

A morphism is said to be a 
\un{bimorphism}
\index{bimorphism (category)} 
if it is both a monomorphism and an epimorphism.  Every isomorphism is a bimorphism.  A category is said to be 
\un{balanced}
\index{balanced (category)} 
if every bimorphism is an isomorphism.  The categories \bSET, \bGR, and \bAB are balanced but the category \bTOP is not.\\

\begingroup%%------------------------------------>>%%----------------------------------->>
\fontsize{9pt}{11pt}\selectfont
In \bSET, \bGR, and \bAB , a morphism is a monomorphism (epimorphism) iff it is injective (surjective).  
In any full subcategory of \bTOP, a morphism is a monomorphism iff it is injective.  
In the full subcategory of $\bTOP_*$ whose objects are the connected spaces, there are monomorphisms that are not injective on the underlying sets (covering projections in this category are monomorphisms).  
In \bTOP, a morphism is an epimorphism iff it is surjective but in \bHAUS, a morphism is an epimorphism iff it has a dense range.  
The homotopy class of a monomorphism (epimorphism) in \bTOP need not be a monomorphism (epimorphism) in \bHTOP.\\
\endgroup%%------------------------------------<< %%------------------------------------<<

Given a category \bC and an object \mX in \bC, let $M(X)$ be the class of all pairs $(Y,f)$, where $f:X \ra Y$ is a monomorphism.  Two elements $(Y,f)$ and $(Z,g)$ of $M(X)$ are deemed equivalent if there exists an isomorphism 
$\phi:Y \ra Z$ such that $f = g \circ \phi$.  
A 
\un{representative class of monomorphisms}
\index{representative class of monomorphisms} 
in $M(X)$ 
is a subclass of $M(X)$ that is a system of representatives for this equivalence relation.  
\bC is said to be 
\un{wellpowered}
\index{wellpowered} 
provided that each of its objects has a representative class of monomorphisms which is a set.

Given a category \bC and an object \mX in \bC, let $E(X)$ be the class of all pairs $(Y,f)$, where $f:X \ra Y$ is an epimorphism.  Two elements $(Y,f)$ and $(Z,g)$ of $E(X)$ are deemed equivalent if there exists an isomorphism 
$\phi:Y \ra Z$ such that $g = \phi \circ f$.  
A 
\un{representative} 
\un{class of epimorphisms}
\index{representative class of monomorphisms} 
in $E(X)$ 
is a subclass of $E(X)$ that is a system of representatives for this equivalence relation.  
\bC is said to be 
\un{cowellpowered}
\index{cowellpowered} 
provided that each of its objects has a representative class of epimorphisms which is a set.\\

\begingroup%%------------------------------------>>%%----------------------------------->>
\fontsize{9pt}{11pt}\selectfont
\bSET, \bGR, \bAB, \bTOP (or \bHAUS) are wellpowered and cowellpowered.  
The category of ordinal numbers is wellpowered but not cowellpowered.\\

\endgroup%%------------------------------------<< %%------------------------------------<<

A monomorphism $f:X \ra Y$ in a category \bC is said to be 
\un{extremal}
\index{extremal monomorphism}
provided that in any factorization $f = h \circ g$, if $g$ is an epimorphism, then $g$ is an isomorphism.

An epimorphism $f:X \ra Y$ in a category \bC is said to be 
\un{extremal}
\index{extremal epimorphism}
provided that in any factorization $f = h \circ g$, if $h$ is a monomorphism, then $h$ is an isomorphism.\\

%%----------------------------------------------------------------------------------------------07
In a balanced category, every monomorphism (epimorphism) is extremal.  In any category, a morphism is an isomorphism iff it is both a monomorphism and an extremal epimorphism iff it is both an extremal monomorphism and an epimorphism.\\

\begingroup%%------------------------------------>>%%----------------------------------->>
\fontsize{9pt}{11pt}\selectfont
In \bTOP, a monomorphism is extremal iff it is an embedding but in \bHAUS, a monomorphism is extremal iff it is a closed embedding.  In \bTOP or \bHAUS, an epimorphism is extremal iff it is a quotient map.\\

\endgroup%%------------------------------------<< %%------------------------------------<<

A 
\un{source}
\index{source (category)} 
in a category \bC is a collection of morphisms $f_i:X \ra X_i$ indexed by a set $I$ and having a common domain.  An 
\un{$n$-source}
\index{n-source (category)} 
is a source for which $\#(I) = n$.

A 
\un{sink}
\index{sink (category)} 
in a category \bC is a collection of morphisms $f_i:X_i \ra X$ indexed by a set $I$ and having a common codomain.  An 
\un{$n$-sink}
\index{n-sink (category)} 
is a sink for which $\#(I) = n$.\\

A 
\un{diagram}
\index{diagram (category)} 
in a category \bC is a functor $\Delta:\bI \ra \bC$, where \bI is a small category, the 
\un{indexing category}.
\index{indexing category (for diagram in a category)}    
To facilitate the introduction of sources and sinks associated with $\Delta$, we shall write 
$\Delta_i$ for the image in $\Ob\bC$ of $i \in \Ob\bI$.

\indent\indent (lim) \ Let $\Delta:\bI \ra \bC$ be a diagram $-$then a source $\{f_i:X \ra \Delta_i\}$ is said to be 
\un{natural}
\index{natural (source)} 
if for each $\delta \in \Mor \ \bI$, say $i \overset{\delta}{\ra} j$, $\Delta\delta \circ f_i = f_j$.  
A \un{limit}
\index{limit (category)} 
of $\Delta$ is a natural source $\{\ell_i: L \ra \Delta_i\}$ with the property that if 
$\{f_i:X \ra \Delta_i\}$ is a natural source, then there exists a unique morphism 
$\phi:X \ra L$ such that $f_i = \ell_i \circ \phi$ for all $i \in \Ob\bI$.  
Limits are essentially unique.  
Notation: \ $L = \lim_{\bI}\Delta$ (or $\lim \Delta$).

\indent\indent (colim) \ Let $\Delta:\bI \ra \bC$ be a diagram $-$then a sink $\{f_i:\Delta_i \ra X\}$ is said to be 
\un{natural}
\index{natural (sink)} 
if for each $\delta \in \Mor \bI$, say $i \overset{\delta}{\ra} j$, $f_i = f_j \circ \Delta\delta$.  
A 
\un{colimit}
\index{colimit (category)} 
of $\Delta$ is a natural sink $\{\ell_i: \Delta_i \ra L\}$ with the property that if 
$\{f_i:\Delta_i \ra X\}$ is a natural sink, then there exists a unique morphism 
$\phi:L \ra X$ such that $f_i = \phi \circ \ell_i$ for all $i \in \Ob\bI$.  
Colimits are essentially unique.  
Notation: \ $L = \colim_{\bI}\Delta$ (or $\colimx \Delta$).

There are a number of basic constructions that can be viewed as a limit or a colimit of a suitable diagram.\\

Let \mI be a set; let \bI be the discrete category with $\Ob\bI = I$.  Given a collection $\{X_i:i \in I\}$ of objects in \bC, 
define a diagram $\Delta:\bI \ra \bC$ by $\Delta_i = X_i$ $(i \in I)$.

\indent\indent (Products) \ \  A limit  \ $\{\ell_i:L \ra \Delta_i\}$ \ of $\Delta$ is said to be a 
\un{product}
\index{product (category)} 
of the $X_i$.  
Notation: \ $L = \prod\limits_i X_i$ (or $X^I$ if $X_i = X$ for all i), $\ell_i = \pr_i$, the 
\un{projection}
\index{projection (product)} 
from $\prod\limits_i X_i$ to $X_i$.  
Briefly put: Products are limits of diagrams with discrete indexing categories.  In particular, the limit of a diagram having \bzero 
for its indexing category is a final object in \bC.

[Note: \ An object \mX in a category \bC is said to be 
\un{final}
\index{final (object in a category)} 
if for each object \mY there is exactly one morphism from \mY to \mX.]

%%----------------------------------------------------------------------------------------------08
\label{13.65}
\indent\indent (Coproducts) \ A colimit $\{\ell_i:\Delta_i \ra L\}$ of $\Delta$ is said to be a 
\un{coproduct} 
\index{coproduct (category)} 
of the $X_i$.  
Notation: \ $L = \coprod\limits_i X_i$ (or $I \cdot X$ if $X_i = X$ for all $i$), $\ell_i = \ini_i$, the 
\un{injection}
\index{injection (coproduct)} 
from $X_i$ to $\coprod\limits_i X_i$.  
Briefly put: Coproducts are colimits of diagrams with discrete indexing categories.  
In particular, the colimit of a diagram having \bzero 
for its indexing category is an initial object in \bC.

[Note: \ An object \mX in a category \bC is said to be 
\un{initial}
\index{initial (object in a category)} 
if for each object \mY there is exactly one morphism from \mX to \mY.]\\

\begingroup%%------------------------------------>>%%----------------------------------->>
\fontsize{9pt}{11pt}\selectfont
In the full subcategory of \bTOP whose objects are the locally connected spaces, the product is the product in \bSET equipped with the coarsest locally connected topology that is finer than the product topology.  In the full subcategory of \bTOP whose objects are the compact Hausdorff spaces, the coproduct is the Stone-\u Cech compactification of the coproduct in \bTOP.\\ 

\endgroup%%------------------------------------<< %%------------------------------------<<

Let \bI be the category 
$1 \  \bullet \overset{a}{\underset{b}{\rightrightarrows} } \bullet \ 2$.  
Given a pair of morphisms $u, \ v:X \ra Y$ in \bC, define a diagram $\Delta:\bI \ra \bC$ by 
$
\begin{cases}
\ \Delta_1 = X\\
\ \Delta_2 = Y
\end{cases}
$
$\&$ 
$
\begin{cases}
\ \Delta a = u\\
\ \Delta b = v
\end{cases}
.
$

\indent\indent (Equalizers) \ An 
\un{equalizer}
\index{equalizer} 
in a category \bC of a pair of morphisms $u, \ v:X \ra Y$ 
is a morphism $f:Z \ra X$ with $u \circ f = v \circ f$ such that for any morphism $f^\prime:Z^\prime \ra X$
with $u \circ f^\prime = v \circ f^\prime$ there exists a unique morphism 
$\phi:Z^\prime \ra Z$ such that $f^\prime = f \circ \phi$.  
The 2-source $X \overset{f}{\lla} Z \overset{u \circ f}{\lra} Y$ is a limit of $\Delta$ iff $Z \overset{f}{\ra} X$ is an equalizer of 
$u, \ v:X \ra Y$.  
 Notation: $Z = \eq(u,v)$.
 
[Note: \ Every equalizer is a monomorphism.  A monomorphism is 
\un{regular}
\index{regular monomorphism}
if it is an equalizer.  A regular monomorphism is extremal.  
In \bSET, \bGR, \bAB, \bTOP (or \bHAUS), an extremal monomorphism is regular.]

\indent\indent (Coequalizers) \ A 
\un{coequalizer}
\index{coequalizer} 
in a category \bC of a pair of morphisms $u, \ v:X \ra Y$ 
is a morphism $f:Y \ra Z$ with $f \circ u = f \circ v$ such that for any morphism $f^\prime:Y \ra Z^\prime$ 
with $f^\prime \circ u = f^\prime \circ v$ there exists a unique morphism 
$\phi:Z \ra Z^\prime$ such that $f^\prime = \phi \circ f$.  
The 2-sink $Y \overset{f}{\lra} Z \overset{f \circ u}{\lla} X$ is a colimit of $\Delta$ iff $Y \overset{f}{\ra} Z$ is a coequalizer of 
 $u, \ v:X \ra Y$.  
 Notation: $Z = \coeq(u,v)$.
 
[Note: \ Every coequalizer is an epimorphism.  
An epimorphism is 
\un{regular}
\index{regular epimorphism}
if it is a coequalizer.  A regular epimorphism is extremal.  
In \bSET, \bGR, \bAB, \bTOP, (or \bHAUS), an extremal epimorphism is regular.]\\

\begingroup%%------------------------------------>>%%----------------------------------->>
\fontsize{9pt}{11pt}\selectfont
There are two aspects to the notion of equalizer and coequalizer, namely: 
(1) \ Existence of $f$ and
(2) \ Uniqueness of $\phi$.  
Given (1), (2) is equivalent to requiring that $f$ be a monomorphism or an epimorphism.  
If (1) is retained and (2) is abandoned, 
then the terminology is
\un{weak equalizer}
\index{weak equalizer} 
or 
\un{weak coequalizer}
\index{weak coequalizer}.
%%----------------------------------------------------------------------------------------------09
For example, $\bHTOP_*$ has neither equalizers nor coequalizers but does have weak equalizers and weak coequalizers.\\
\endgroup%%------------------------------------<< %%------------------------------------<<

\label{12.12}
\label{12.34}
Let \bI be the category 
$1 \bullet \overset{a}{\ra} \underset{3}{\bullet} \overset{b}{\la} 2$.  
Given morphisms 
$
\begin{cases}
\ f:X \ra Z\\
\ g:Y \ra Z
\end{cases}
$
in \bC, define a diagram $\Delta:\bI \ra \bC$ by 
$
\begin{cases}
\ \Delta_1 = X\\
\ \Delta_2 = Y\\
\ \Delta_3 = Z
\end{cases}
$
$\&$
$
\begin{cases}
\ \Delta a = f\\
\ \Delta b = g
\end{cases}
.
$

\indent\indent (Pullbacks) \ 
Given a 2-sink 
$X \overset{f}{\ra} Z \overset{g}{\la} Y$,
a commutative diagram
\begin{tikzcd}[ sep=large]
{P} \ar{d}[swap]{\xi}  \ar{r}{\eta}  &{Y} \ar{d}{g}\\
{X} \ar{r}[swap]{f}  &{Z}
\end{tikzcd}
is said to be a 
\un{pullback square}
\index{pullback square} 
if for any 2-source
$X \overset{\xi^\prime}{\la} P^\prime \overset{\eta^\prime}{\ra} Y$
with $f \circ \xi^\prime = g \circ \eta^\prime$ 
there exists a unique morphism $\phi:P^\prime \ra P$ such that 
$\xi^\prime = \xi \circ \phi$ and $\eta^\prime = \eta \circ \phi$.  
The 2-source 
$X \overset{\xi}{\la} P \overset{\eta}{\ra} Y$ is called a 
\un{pullback}
\index{pullback} of the 2-sink 
$X \overset{f}{\ra} Z \overset{g}{\la} Y$.  
Notation: $P = X \times_Z Y$.  Limits of $\Delta$ are pullback squares and conversely.

Let \bI be the category 
$1 \bullet \overset{a}{\la} \underset{3}{\bullet} \overset{b}{\ra} 2$.  
Given morphisms 
$
\begin{cases}
\ f:Z \ra X\\
\ g:Z \ra Y
\end{cases}
$
in \bC, define a diagram $\Delta:\bI \ra \bC$ by 
$
\begin{cases}
\ \Delta_1 = X\\
\ \Delta_2 = Y\\
\ \Delta_3 = Z
\end{cases}
$
$\&$
$
\begin{cases}
\ \Delta a = f\\
\ \Delta b = g
\end{cases}
.
$

\indent\indent (Pushouts) \ 
Given a 2-source 
$X \overset{f}{\la} Z \overset{g}{\ra} Y$,
a commutative diagram
\begin{tikzcd}[ sep=large]
{Z} \ar{d}[swap]{f}  \ar{r}{g}  &{Y} \ar{d}{\eta}\\
{X} \ar{r}[swap]{\xi}  &{P}
\end{tikzcd}
is said to be a 
\un{pushout square}
\index{pushout square} if for any 2-sink
$X \overset{\xi^\prime}{\ra} P^\prime \overset{\eta^\prime}{\la} Y$ 
with 
$\xi^\prime \circ f = \eta^\prime \circ g$ 
there exists a unique morphism $\phi:P\ra P^\prime $ such that 
$\xi^\prime = \phi \circ \xi$ 
and 
$\eta^\prime = \phi \circ \eta$.  
The 2-sink 
$X \overset{\xi}{\ra} P \overset{\eta}{\la} Y$ 
is called a 
\un{pushout}
\index{pushout} 
of the 2-source 
$X \overset{f}{\la} Z \overset{g}{\ra} Y$.  
Notation: $P = X \underset{Z}{\sqcup} Y$.  Colimits of $\Delta$ are pushout squares and conversely.\\

\begingroup%%------------------------------------>>%%----------------------------------->>
\fontsize{9pt}{11pt}\selectfont
The result of dropping uniqueness in $\phi$ is 
\un{weak pullback}
\index{weak pullback} 
or 
\un{weak pushout}
\index{weak pushout}.  
Examples are the commutative squares that define fibration and cofibration in \bTOP.\\

\endgroup%%------------------------------------<< %%------------------------------------<<

Let \bI be a small category, $\Delta:\bI^\OP \times \bI \ra \bC$ a diagram.

\indent\indent (Ends) \ A source $\{f_i:X \ra \Delta_{i,i}\}$ is said to be 
\un{dinatural}
\index{dinatural (source)}
if for each 
$\delta \in \Mor \bI$, say 
$i \overset{\delta}{\ra} j$, 
$\Delta(\id,\delta) \circ f_i = \Delta(\delta,\id) \circ f_j$.
An 
\un{end}
\index{end} 
of $\Delta$ is a dinatural source 
$\{e_i:E \ra \Delta_{i,i}\}$ 
with the property that if 
$\{f_i:X \ra \Delta_{i,i}\}$ 
is a dinatural source, then there exists a unique morphism 
$\phi:X \ra E$ such that $f_i = e_i \circ \phi$ for all $i \in \Ob\bI$.  
Every end is a limit (and every limit is an end.)  
Notation: $E = \ds\int_i \Delta_{i,i}$ (or $\ds\int_{\bI} \Delta$).

\indent\indent (Coends)  \ A sink $\{f_i:\Delta_{i,i} \ra X\}$ is said to be 
\un{dinatural} \index{dinatural (sink)}
if for each 
$\delta \in \Mor \bI$, say 
$i \overset{\delta}{\ra} j$, 
$f_i \circ \Delta(\delta,\id) = f_j \circ \Delta(\id,\delta)$.
A 
\un{coend}
\index{coend} 
of $\Delta$ is a dinatural sink 
$\{e_i:\Delta_{i,i} \ra E\}$ 
with the property that if 
$\{f_i:\Delta_{i,i} \ra X\}$ 
is a dinatural sink, then there exists a unique
%%----------------------------------------------------------------------------------------------10
morphism 
$\phi:E \ra X$ such that $f_i = \phi \circ e_i $ for all $i \in \Ob\bI$.  
Every coend is a colimit (and every colimit is a coend.)  
Notation: $E = \ds\int^i \Delta_{i,i}$ (or $\ds\int^{\bI} \Delta$).\\

\begingroup%%------------------------------------>>%%----------------------------------->>
\fontsize{9pt}{11pt}\selectfont
Let 
$
\begin{cases}
\ F:\bI \ra \bC\\
\ G:\bI \ra \bC
\end{cases}
$
be functors $-$then the assignment $(i,j) \ra \Mor(Fi,Gj)$ defines a diagram 
$\bI^\OP \times \bI \ra \bSET$ and $\Nat(F,G)$ is the end $\ds\int_i \Mor(Fi,Gi)$.\\

\endgroup%%------------------------------------<< %%------------------------------------<<

\index{Integral Yoneda Lemma}
\textbf{\small INTEGRAL YONEDA LEMMA} \quadx
Let \bI be a small category, \bC a complete and cocomplete category $-$then for every $F$ in $[\bI^\OP,\bC]$, 
$\ds\int^i \Mor(-,i) \cdot F i \approx F \approx \ds\int_i F i^{\Mor(i,-)}$.\\

Let $\bI \neq \bzero$ be a small category $-$then \bI is said to be 
\un{filtered}
\index{filtered (small category)} 
if\\ 
%^
\indent\indent (F$_1$) \ Given any pair of objects i, j in \bI, there exists an object $k$ and morphisms 
$
\begin{cases}
\ i \ra k\\
\ j \ra k
\end{cases}
;
$

\indent\indent (F$_2$) \ Given any pair of morphisms $a, b:i \ra j$ in \bI, there exists an object $k$ and a morphism
$c:j \ra k$ such that $c \circ a = c \circ b$.

Every nonempty directed set $(I,\leq)$ can be viewed as a filtered category \bI, where $\Ob\ \bI = I$ and 
$\Mor(i,j)$ is a one element set when $i \leq j$ but is empty otherwise.\\

\label{16.2}
\label{16.7}
Example: Let $[\N]$ be the filtered category associated with the directed set of nonnegative integers.  
Given a category \bC, denote by $\bFIL(\bC)$
\index{$\bFIL(\bC)$} 
the functor category $[[\N],\bC]$ 
$-$then an object $(\bX,\bff)$ in $\bFIL(\bC)$ is a 
sequence $\{X_n,f_n\}$, where $X_n \in \Ob\bC$ $\&$ $f_n \in \Mor(X_n,X_{n+1})$, and a morphism 
$\phi:(\bX,\bff) \ra (\bY,\bg)$ in $\bFIL(\bC)$ is a sequence $\{\phi_n\}$, where $\phi_n \in \Mor(X_n,Y_n)$ $\&$ 
$g_n \circ \phi_n = \phi_{n+1} \circ f_n$.\\

\indent\indent (Filtered Colimits) \ A 
\un{filtered colimit}
\index{filtered colimit} 
in \bC is the colimit of a diagram 
$\Delta:\bI \ra$ $\bC$, where \bI is filtered.

\indent\indent (Cofiltered Limits) \ A 
\un{cofiltered limit}
\index{filtered limit} 
in \bC is the limit of a diagram 
$\Delta:\bI \ra$ $\bC$, where \bI is cofiltered.

[Note: \ A small category $\bI \neq \bzero$ is said to be 
\un{cofiltered}
\index{cofiltered (small category)} 
provided that $\bI^\OP$ is filtered.]\\

\begingroup%%------------------------------------>>%%----------------------------------->>
\fontsize{9pt}{11pt}\selectfont
A Hausdorff space is compactly generated iff it is the filtered colimit in \bTOP of its compact subspaces.  Every compact Hausdorff space is the cofiltered limit in \bTOP of compact metrizable spaces.\\

\endgroup%%------------------------------------<< %%------------------------------------<<

Given a small category \bC, a 
\un{path}
\index{path (small category)}  
in \bC is a diagram $\sigma$ of the form 
$X_0 \ra X_1 \la$ $\cdots \ra$ 
$X_{2n-1} \la X_{2n}$ $(n \geq 0)$.  
One says that $\sigma$ 
\un{begins}
\index{begins (path in a small category)} 
at $X_0$ 
and 
\un{ends}
\index{ends (path in a small category)} 
at $X_{2n}$.  The quotient of $\Ob\bC$ with respect to the equivalence 
relation  obtained by declaring that $X^\prime \sim X\pp$ iff there
%%----------------------------------------------------------------------------------------------11
exists a path in \bC which begins at $X^\prime$ and ends at $X\pp$ is the set $\pi_0(\bC)$ of 
\un{components}
\index{components (category)} of \bC, \bC being called 
\un{connected}
\index{connected (category)} 
when the cardinality of $\pi_0(\bC)$ is one.  The full subcategory of \bC determined by a component is connected and is maximal with respect to this property.  If \bC has an initial object or a final object, then \bC is connected.

[Note: \ The concept of ``path'' makes sense in any category.]\\

\begingroup%%------------------------------------>>%%----------------------------------->>
\fontsize{9pt}{11pt}\selectfont
Let $\bI \neq \bzero$ be a small category $-$then \bI is said to be 
\un{pseudofiltered}
\index{pseudofiltered (small category)} 
if \\
\indent\indent (PF$_1$) \ Given any pair of morphisms 
$
\begin{cases}
\ a:i \ra j\\
\ b:i \ra k
\end{cases}
$
in \bI, there exists an object $\ell$ and morphisms 
$
\begin{cases}
\ c:j \ra \ell\\
\ d:k \ra \ell
\end{cases}
$
such that $c \circ a = d \circ b$;

\indent\indent (PF$_2$) \ Given any pair of morphisms $a, b:i \ra j$ in \bI, there exists  a morphism 
$c:j \ra k$ such that $c \circ a = c \circ b$.

\bI is filtered iff \bI is connected and pseudofiltered.  \bI is pseudofiltered iff its components are filtered.\\

\endgroup%%------------------------------------<< %%------------------------------------<<

\label{15.31}

Given small categories 
$
\begin{cases}
\ \bI\\
\ \bJ
\end{cases}
\hspace{-.25cm},
$
a functor $\nabla:\bJ \ra \bI$ is said to be 
\un{final}
\index{final (functor)} 
provided that for every $i \in \Ob\bI$, the comma 
category $\abs{K_i,\nabla}$ is nonempty and connected.  
If \bJ is filtered and $\nabla:\bJ \ra \bI$ is final, then \bI is filtered.
\label{13.122}

[Note: \ A subcategory of a small category is 
\un{final}
\index{final (subcategory of a small category)} 
if the inclusion 
is a final functor.]

Let $\nabla:\bJ \ra \bI$ be final.  
Suppose that $\Delta:\bI \ra \bC$ is a diagram for which 
$\colimx \Delta \circ \nabla$ exists $-$then \ $\colimx \Delta$ exists and the arrow \ 
$\colimx \Delta \circ \nabla \ra \colimx \Delta$ \ is an isomorphism.  
Corollary:  If $i$ is a final object in \bI, then $\colimx \Delta \approx \Delta_i$.

[Note: \  Analogous considerations apply to limits so long as ``final'' is replaced throughout by ``initial''.]\\

\begingroup%%------------------------------------>>%%----------------------------------->>
\fontsize{9pt}{11pt}\selectfont
Let \bI be a filtered category $-$then there exists a directed set $(J,\leq)$ and a final functor $\nabla:\bJ \ra \bI$.\\

\endgroup%%------------------------------------<< %%------------------------------------<<

Limits commute with limits.  In other words, if $\Delta:\bI \times \bJ \ra \bC$ is a diagram, then under obvious assumptions %dmc
\[
\text{lim}_\bI \text{lim}_{\bJ} \Delta 
\hsx \approx \hsx 
\text{lim}_{\bI \times \bJ} \Delta 
\hsx \approx \hsx  
\text{lim}_{\bJ \times \bI} \Delta 
\hsx \approx \hsx  
\text{lim}_{\bJ}\text{lim}_{\bI} \Delta.
%\lim_{\bI} \lim_{\bJ} \Delta 
%\approx 
%\lim_{\bI \times \bJ} \Delta 
%\approx 
%\lim_{\bJ \times \bI} \Delta 
%\approx 
%\lim_{\bJ}\lim_{\bI} \Delta.
\]

Likewise, colimits commute with colimits.  In general, limits do not commute with colimits.  
However, if $\Delta:\bI \times \bJ \ra \bSET$ and if \bI is finite and \bJ is filtered, then the arrow 
$\colim_{\bJ} \lim_{\bI} \Delta \ra \lim_{\bI}\colim_{\bJ} \Delta$ is a bijection, so that in \bSET filtered colimits commute 
with finite limits.

[Note: \ In \bGR, \bAB, or \bRG filtered colimits commute with finite limits.  But, e.g., filtered colimits commute 
do not commute with finite limits in $\bSET^\OP$.]\\

\begingroup%%------------------------------------>>%%----------------------------------->>
\fontsize{9pt}{11pt}\selectfont
In \bAB (or any Grothendieck category), pseudofiltered colimits commute with finite limits.\\

\endgroup%%------------------------------------<< %%------------------------------------<<

%%----------------------------------------------------------------------------------------------12
A category \bC is said to be 
\un{complete}
\index{complete}
(\un{cocomplete})
\index{cocomplete}
if for each small category \bI, every $\Delta \in \Ob[\bI,\bC]$ has a limit (colimit).  
The following are equivalent.

\indent\indent (1) \ \bC is complete (cocomplete).\\
\indent\indent (2) \ \bC has products and equalizers (coproducts and coequalizers).\\
\indent\indent (3) \ \bC has products and pullbacks (coproducts and pushouts).\\
\indent\indent (4) \ \bC has a final object and multiple pullbacks  (initial object and multiple \\pushouts).

[Note: \ A source $\{\xi_i:P \ra X_i\}$ (sink $\{\xi_i:X_i \ra P\}$) is said be be a 
\un{multiple pullback}
\index{multiple pullback}
(\un{multiple pushout})
\index{multiple pushout} 
of a sink $\{f_i:X_i \ra X\}$ (source  $\{f_i:X \ra X_i\}$) 
provided that 
$f_i \circ \xi_i = f_j \circ \xi_j$ ($\xi_i \circ f_i = \xi_j \circ f_j$) $\forall$ 
$
\begin{cases}
\ i\\
\ j
\end{cases}
$
and if for any source $\{\xi_i^\prime:P^\prime \ra X_i\}$ (sink $\{\xi_i^\prime:X_i \ra P^\prime\}$) 
with 
$f_i \circ \xi_i^\prime = f_j \circ \xi_j^\prime$ ($\xi_i^\prime \circ f_i = \xi_j^\prime \circ f_j$) $\forall$
$
\begin{cases}
\ i\\
\ j
\end{cases}
$
, there exists a unique morphism $\phi:P^\prime \ra P$ ($\phi:P \ra P^\prime$) such that $\forall \ i$, 
$\xi_i^\prime = \xi_i \circ \phi$ ($\xi_i^\prime = \phi \circ \xi_i$).  
Every multiple pullback (multiple pushout) is a limit (colimit).]\\

\begingroup%%------------------------------------>>%%----------------------------------->>
\fontsize{9pt}{11pt}\selectfont
The categories \bSET, \bGR, and \bAB are both complete and cocomplete.  
The same is true of \bTOP and $\bTOP_*$ but not of \bHTOP and $\bHTOP_*$.
\\ \indent
[Note: \ \bHAUS is complete; it is also cocomplete, being epireflective in \bTOP.]\\
\endgroup%%------------------------------------<< %%------------------------------------<<

A category \bC is said to be 
\un{finitely complete}
\index{finitely complete} 
(
\un{finitely cocomplete})
\index{finitely cocomplete} 
if for each finite category \bI, every $\Delta \in \Ob[\bI,\bC]$ has a limit (colimit).  The following are equivalent.\\
%^
\indent\indent (1) \ \bC is finitely complete (finitely cocomplete).\\
\indent\indent (2) \ \bC has finite products and equalizers (finite coproducts and coequalizers).\\
\indent\indent (3) \ \bC has finite products and pullbacks (finite coproducts and pushouts).\\
\indent\indent (4) \ \bC has a final object and pullbacks (initial object and pushouts).\\

\begingroup%%------------------------------------>>%%----------------------------------->>
\fontsize{9pt}{11pt}\selectfont
The full subcategory of \bTOP whose objects are the finite topological spaces is finitely complete and finitely cocomplete but neither complete nor cocomplete.  A nontrivial group, considered as a category, has multiple pullbacks but fails to have finite products.\\
\endgroup%%------------------------------------<< %%------------------------------------<<

If \bC is small and \bD is finitely complete and wellpowered (finitely cocomplete and cowellpowered), 
then $[\bC,\bD]$ is wellpowered (cowellpowered).\\

\begingroup%%------------------------------------>>%%----------------------------------->>
\fontsize{9pt}{11pt}\selectfont
$\bSET(\ra)$, $\bGR(\ra)$, $\bAB(\ra)$, $\bTOP(\ra)$ (or $\bHAUS(\ra)$) are wellpowered and cowellpowered.

[Note: \ The arrow category $\bC(\ra)$ of any category \bC is isomorphic to $[\btwo,\bC]$.]\\
\endgroup%%------------------------------------<< %%------------------------------------<<

Let $F:\bC \ra \bD$ be a functor.

%%----------------------------------------------------------------------------------------------13
\indent\indent (a) \ \mF is said to preserve a limit 
$\{\ell_i:L \ra \Delta_i\}$ 
$(\text{colimit} \{\ell_i:\Delta_i \ra L\})$ of a diagram 
$\Delta:\bI \ra \bC$ if 
$\{F\ell_i:FL \ra F\Delta_i\}$ 
$(\{F\ell_i:F\Delta_i \ra FL\})$ 
is a limit (colimit) of the diagram 
$F \circ \Delta:\bI \ra \bD$.

\indent\indent (b) \ \mF is said to preserve limits (colimits) over an indexing category \bI if \mF preserves all limits (colimits) of diagrams $\Delta:\bI \ra \bC$.

\indent\indent (c) \ \mF is said to preserve limits (colimits) if \mF preserves limits (colimits) over all indexing categories \bI.\\

\begingroup%%------------------------------------>>%%----------------------------------->>
\fontsize{9pt}{11pt}\selectfont
The forgetful functor $\bTOP \ra \bSET$ preserves limits and colimits.  The forgetful functor $\bGR \ra \bSET$ preserves limits and filtered colimits but not coproducts.  The inclusion $\bHAUS \ra \bTOP$ preserves limits and coproducts but not coequalizers.  The inclusion $\bAB \ra \bGR$ preserves limits but not colimits.\\
\endgroup%%------------------------------------<< %%------------------------------------<<

There are two rules that determine the behavior of 
$
\begin{cases}
\ \Mor(X,-)\\
\ \Mor(-,X)
\end{cases}
$
with respect to limits and colimits.

\indent\indent (1) \ The functor $\Mor(X,-):\bC \ra \bSET$ preserves limits.  Symbolically, therefore, 
$\Mor(X,\lim \ \Delta) \approx$ $\lim(\Mor(X,-) \circ \Delta)$.

\indent\indent (2) \ The cofunctor $\Mor(-,X):\bC \ra \bSET$ converts colimits into limits.  Symbolically, therefore, 
$\Mor(\colimx \Delta,X) \approx$ $\lim(\Mor(-,X) \circ \Delta)$.\\

\index{Theorem: Representable Functor Theorem}
\index{Representable Functor Theorem}
\textbf{\small REPRESENTABLE FUNCTOR THEOREM} \quadx
Given a complete category \bC, a functor 
$F:\bC \ra \bSET$ is representable iff \mF preserves limits and satisfies the 
\un{solution set} \un{condition}:
\index{solution set condition} 
There exists a set $\{X_i\}$ of objects in \bC such that for each $X \in \Ob\bC$ and each $y \in FX$, 
there is an $i$, a $y_i \in FX_i$, and an $f:X_i \ra X$ such that $y = (Ff)y_i$.\\

\begingroup%%------------------------------------>>%%----------------------------------->>
\fontsize{9pt}{11pt}\selectfont
Take for \bC the category opposite to the category of ordinal numbers $-$then the functor $\bC \ra \bSET$ defined by 
$\alpha \ra *$ has a complete domain and preserves limits but is not representable.\\
\endgroup%%------------------------------------<< %%------------------------------------<<

Limits and colimits in functor categories are computed ``object by object''.  
So, if \bC is a small category, then \bD (finitely) complete $\implies$ $[\bC,\bD]$ (finitely) complete and 
\bD (finitely) cocomplete $\implies$ $[\bC,\bD]$ (finitely) cocomplete.\\

Given a small category \bC, put $\widehat{\bC} = [\bC^\OP,\bSET]$ $-$then 
$\widehat{\bC}$ is complete and cocomplete.  
The Yoneda embedding 
$Y_{\bC}:\bC \ra \widehat{\bC}$ preserves limits; it need not, however, preserve finite colimits.  
The image of \bC is ``colimit dense'' in $\widehat{\bC}$, i.e., every cofunctor 
$\bC \ra \bSET$ is a colimit of representable cofunctors.\\

An 
\un{indobject}
\index{indobject} 
in a small category \bC is a diagram $\Delta:\bI \ra \bC$, where \bI is filtered.  
Corresponding to an indobject $\Delta$, is the object $L_\Delta$ in $\widehat{\bC}$ defined by 
$L_\Delta = \colim(Y_{\bC} \circ \Delta)$.  
%%----------------------------------------------------------------------------------------------14
The 
\un{indcategory}
\index{indcategory} $\bIND(\bC)$ 
\index{$\bIND(\bC) $}
of \bC is the category whose objects are the indobjects and whose morphisms are the sets $\Mor(\Delta^\prime,\Delta\pp) = \Nat(L_{\Delta^\prime},L_{\Delta\pp})$.  The functor 
$L:\bIND(\bC) \ra$ $\widehat{\bC}$ 
that sends $\Delta$ to $L_\Delta$ is full and faithful (although in general not injective on objects), hence establishes an equivalence between $\bIND(\bC)$ and the full subcategory of $\widehat{\bC}$ 
whose objects are the cofunctors $\bC \ra \bSET$ which are filtered colimits of representable cofunctors.  
The category $\bIND(\bC)$ has filtered colimits; they are preserved by $L$, as are all limits.  
Moreover, in $\bIND(\bC)$, filtered colimits commute with finite limits.  If \bC is finitely cocomplete, then $\bIND(\bC)$ is complete and cocomplete.  
The functor 
$K:\bC \ra$ $\bIND(\bC)$ that sends \mX to $K_X$, where $K_X:\bone \ra \bC$ is the constant functor with value \mX, is full, faithful, and injective on objects.  
In addition, $K$ preserves limits and finite colimits.  
The composition 
$\bC \overset{K}{\lra} \bIND(\bC) \overset{L}{\lra} \widehat{\bC}$  is the Yoneda embedding $Y_\bC$.  
A cofunctor $F \in \Ob\widehat{\bC}$ is said to be 
\un{indrepresentable}
\index{indrepresentable} 
if it is 
naturally isomorphic to a functor of the form $L_\Delta$, $\Delta \in \Ob\bIND(\bC)$.  
An indrepresentable cofunctor converts finite colimits into finite limits and conversely, provided that \bC is finitely cocomplete.

[Note: \ The 
\un{procategory}
\index{procategory} 
$\bPRO(\bC)$
\index{$\bPRO(\bC)$} 
is by definition 
$\bIND(\bC^\OP)^\OP$.  Its objects are the 
\un{proobjects}
\index{proobjects} 
in \bC, i.e., the diagrams defined on cofiltering categories.]\\

\begingroup%%------------------------------------>>%%----------------------------------->>
\fontsize{9pt}{11pt}\selectfont
The full subcategory of \bSET whose objects are the finite sets is equivalent to a small category.  
Its indcategory is equivalent to \bSET and its procategory is equivalent to the full subcategory of \bTOP whose objects are the totally disconnected compact Hausdorff spaces.
\\ \indent
[Note: \ There is no small category \bC for which $\bPRO(\bC)$ is equivalent to \bSET.  This is because in \bSET, cofiltered limits do not commute with finite colimits.]\\
\endgroup%%------------------------------------<< %%------------------------------------<<

Given categories 
$
\begin{cases}
\ \bC\\
\ \bD
\end{cases}
, \ 
$
functors 
$
\begin{cases}
\ F:\bC \ra \bD\\
\ G:\bD \ra \bC
\end{cases}
$
are said to be an 
\un{adjoint pair}
\index{adjoint pair} if the functors 
$
\begin{cases}
\ \Mor \circ (F^\OP \times \id_{\bD})\\
\ \Mor \circ (\id_{\bC^\OP} \times G)
\end{cases}
$
from $\bC^\OP \times \bD$ to \bSET are naturally isomorphic, i.e., if it is possible to assign to each ordered pair
$
\begin{cases}
\ X \in \Ob\bC\\
\ Y \in \Ob\bD
\end{cases}
$
\label{12.33}
\label{13.13}
a bijective map $\Xi_{X,Y}:\Mor(FX,Y) \ra \Mor(X,GY)$ which is functorial in \mX and \mY.  When this is so, $F$ is a 
\un{left adjoint}
\index{left adjoint} 
for $G$ and \mG is a 
\un{right adjoint}
\index{right adjoint} 
for $F$.  
Any two left (right) adjoints for \mG (\mF) are naturally isomorphic.  Left adjoints preserve colimits; right adjoints preserve limits.  In order that $(F,G)$ be an adjoint pair, it is necessary and sufficient that there exist natural transformations 
$
\begin{cases}
\ \mu \in \Nat(\id_{\bC},G \circ F)\\
\ \nu \in \Nat(F \circ G,\id_{\bD})
\end{cases}
$
subject to 
$
\begin{cases}
\ (G \nu) \circ(\mu G) = \id_G\\
\ (\nu F) \circ (F \mu) = \id_F
\end{cases}
. \ 
$
The data $(F,G,\mu,\nu)$ is referred to as an 
\un{adjoint situation},
\index{adjoint situation} 
the natural transformations 
$
\begin{cases}
\ \mu:\id_{\bC} \ra G \circ F\\
\ \nu:F \circ G \ra \id_{\bD}
\end{cases}
$
being the 
\un{arrows of adjunction}.
\index{arrows of adjunction}

(UN) \ Suppose that \mG has a left adjoint \mF $-$then for each $X \in \Ob\bC$, each
%%----------------------------------------------------------------------------------------------15
$Y \in \Ob\bD$, and each $f:X \ra GY$, there exists a unique $g:FX \ra Y$ such that $f =$ $G g \circ \mu_X$.

[Note: \ When reformulated, this property is characteristic.]\\

\begingroup%%------------------------------------>>%%----------------------------------->>
\fontsize{9pt}{11pt}\selectfont
The forgetful functor $\bTOP \ra \bSET$ has a left adjoint that sends a set \mX to the pair $(X,\tau)$, where $\tau$ is the discrete topology, and a right adjoint that sends a set \mX to the pair $(X,\tau)$, where $\tau$ is the indiscrete topology.\\

\label{14.2}
Let \bI be a small category, \bC a complete and cocomplete category.  
Examples: 
(1) \ The constant diagram functor $K:\bC \ra [\bI,\bC]$ has a left adjoint, viz.,
$\colim:[\bI,\bC] \ra \bC$, and a right adjoint viz. $\lim:[\bI,\bC] \ra \bC$; 
(2) \ The functor $\bC \ra [\bI^\OP \times \bI,\bC]$ that sends \mX to 
$(i,j) \ra \Mor(i,j) \cdot X$ is a left adjoint for end and the functor that sends 
\mX to $(i,j) \ra X^{\Mor(j,i)}$ is a right adjoint for coend.\\
\endgroup%%------------------------------------<< %%------------------------------------<<

\index{Theorem: General Adjoint Functor Theorem}
\index{General Adjoint Functor Theorem}
\textbf{\small GENERAL ADJOINT FUNCTOR THEOREM} \quadx
Given a complete category \bD, a functor $G:\bD \ra \bC$ has a left adjoint iff \mG preserves limits and satisfies the 
\un{solution set} \un{condition}: 
\index{solution set condition (adjoint functor)}
For each $X \in \Ob\bC$, there exists a source $\{f_i:X \ra GY_i\}$ such that for every $f:X \ra GY$, there is an $i$ and a 
$g:Y_i \ra Y$ such that $f = G g \circ f_i$.\\

\begingroup%%------------------------------------>>%%----------------------------------->>
\fontsize{9pt}{11pt}\selectfont
The general adjoint functor theorem implies that a small category is complete iff it is cocomplete.\\
\endgroup%%------------------------------------<< %%------------------------------------<<

\index{Theorem: Kan Extension Theorem}
\index{Kan Extension Theorem}
\textbf{\small KAN EXTENSION THEOREM} \quadx
Given small categories 
$
\begin{cases}
\ \bC\\
\ \bD
\end{cases}
$
, a complete (cocomplete) category \bS, and a functor $K:\bC \ra \bD$, the functor 
$[K,\bS]:[\bD,\bS] \ra [\bC,\bS]$ has a right (left) adjoint ran (lan) and preserves limits (colimits).

[Note: \ If \mK is full and faithful, then ran (lan) is full and faithful.]
\\

Suppose that \bS is complete.  Let $T \in \Ob[\bC,\bS]$ $-$then  $\ran T$ is called the 
\un{right Kan} \un{extension}
\index{right Kan extension} 
of \mT along \mK.  In terms of ends, 
$(\ran T)Y = \ds\int_X TX^{\Mor(Y,KX)}$.  
There is a ``universal'' arrow $(\ran T) \circ K \ra T$.  It is a natural isomorphism if \mK is full and faithful.

Suppose that \bS is cocomplete.  Let $T \in \Ob[\bC,\bS]$ $-$then lan\mT is called the 
\un{left Kan} \un{extension}
\index{left Kan extension} 
of \mT along \mK.  In terms of coends, 
$(\lan T)Y = \ds\int^X \Mor(KX,Y) \cdot TX$.  
There is a ``universal'' arrow $T \ra (\lan T) \circ K$.  It is a natural isomorphism if \mK is full and faithful.
\\

Application: If \bC and \bD are small categories and if $F:\bC \ra \bD$ is a functor, then the precomposition functor 
$\widehat{\bD} \ra \widehat{\bC}$ has a left adjoint 
$\widehat{F}:\widehat{\bC} \ra \widehat{\bD}$ and 
$\widehat{F} \circ Y_{\bC} \approx Y_{\bD} \circ F$.

[Note: \ One can always arrange that $\widehat{F} \circ Y_{\bC} = Y_{\bD} \circ F$.]
\\

The construction of the right (left) adjoint of $[K,\bS]$ does not use the assumption that \bD is small, its role being to ensure that $[\bD,\bS]$ is a category. 
For example, if \bC is small
%%----------------------------------------------------------------------------------------------16
and \bS is cocomplete, then taking $K = Y_{\bC}$, the functor 
$[Y_{\bC},\bS]:[\widehat{\bC},\bS] \ra [\bC,\bS]$ has a left adjoint that sends $T \in \Ob[\bC,\bS]$ to 
$\Gamma_T \in \Ob[\widehat{\bC},\bS]$, where $\Gamma_T \circ Y_{\bC} = T$.  
On an object $F \in \widehat{\bC}$, 
$\Gamma_T F = \ds\int^X \Nat(Y_{\bC} X,F)\cdot TX = \ds\int^X FX \cdot TX$.  
$\Gamma_T$ is the 
\label{13.10}
\label{13.115}
\un{realization functor};
\index{realization functor} 
it is a left adjoint for the 
\un{singular functor}
\index{singular functor} 
$S_T$, the composite of the Yoneda embedding $\bS \ra [\bS^\OP,\bSET]$ 
and the precomposition functor $[\bS^\OP,\bSET] \ra [\bC^\OP,\bSET]$, thus 
$(S_T Y) X$ $=$ $\Mor(TX,Y)$.

[Note: \ The arrow of adjunction $\Gamma_T \circ S_T \ra \id_{\bS}$ is a natural isomorphism iff $S_T$ is full and faithful.]\\

\bCAT is the category whose objects are the small categories and whose morphisms are the functors between them: 
$\bC, \bD \in \Ob\bCAT$ $\implies$ $\Mor(\bC,\bD) = \Ob[\bC,\bD]$. \bCAT is concrete and complete and cocomplete.  
\bzero is an initial object in \bCAT and \bone is a final object in \bCAT.\\

\begingroup%%------------------------------------>>%%----------------------------------->>
\fontsize{9pt}{11pt}\selectfont
Let $\pi_0:\bCAT \ra \bSET$ be the functor that sends \bC to $\pi_0(\bC)$, the set of components of \bC; 
let $\dis:\bSET \ra$ $\bCAT$ be the functor that sends \mX to $\dis X$, the discrete category on \mX; 
let $\ob:\bCAT \ra$ $\bSET$ be the functor that sends \bC to \Ob\bC, the set of objects in \bC; 
let $\grd:\bSET \ra$ $\bCAT$ be the functor that sends \mX to $\grd X$, the category whose objects are the elements of \mX and whose morphisms are the elements of $X \times X$ $-$then $\pi_0$ is a left adjoint for dis, dis is a left adjoint for ob, and ob is a left adjoint for grd.
\\ \indent
[Note: \ $\pi_0$ preserves finite products; it need not preserve arbitrary products.]\\
\endgroup%%------------------------------------<< %%------------------------------------<<

\bGRD is the full subcategory of \bCAT whose objects are the groupoids, i.e., the small categories in which every morphism is invertible.  \ \ 
Example: \ The assignment \ 
$
\Pi:
\begin{cases}
\ \bTOP \ra \bGRD\\
\ X \ra \Pi X
\end{cases}
$
is a functor.\\

\label{13.9}
\label{14.68}
\begingroup%%------------------------------------>>%%----------------------------------->>
\fontsize{9pt}{11pt}\selectfont
Let $\iso:\bCAT \ra \bGRD$ be the functor that sends \bC to $\iso \bC$, the groupoid whose objects are those of \bC and whose morphisms are the invertible morphisms in \bC $-$then iso is a right adjoint for the inclusion 
$\bGRD \ra \bCAT$.  Let $\pi_1:\bCAT \ra \bGRD$ be the functor that sends \bC to $\pi_1(\bC)$, the 
\un{fundamental groupoid}
\index{fundamental groupoid} 
of \bC, 
i.e., the localization of \bC at $\Mor \bC$ $-$then $\pi_1$ is a left adjoint for the inclusion $\bGRD \ra \bCAT$.\\

\endgroup%%------------------------------------<< %%------------------------------------<<

\label{14.171}
$\bDelta$ is the category whose objects are the ordered sets 
$[n] \equiv \{0,1, \ldots, n\}$ $(n \geq 0)$ and whose morphisms are the order preserving maps.  
In $\bDelta$, every morphism can be written as an epimorphism followed by a monomorphism and a morphism is a monomorphism (epimorphism) iff it is injective (surjective).  
The 
\un{face operators}
\index{face operators} 
are the monomorphisms 
$\delta_i^n:[n-1] \ra [n]$ $(n > 0, 0 \leq i \leq n)$ defined by omitting the value $i$.  
The 
\un{degeneracy} 
\un{operators}
\index{degeneracy operators} 
are epimorphisms  
$\sigma_i^n:[n+1] \ra [n]$ $(n \geq 0, 0 \leq i \leq n)$ defined by repeating the value $i$.  
Suppressing superscripts, if $\alpha \in \Mor([m],[n])$ is not
%%----------------------------------------------------------------------------------------------17
the identity, then $\alpha$ has a unique factorization 
$\alpha = (\delta_{i_1} \circ \cdots \circ \delta_{i_p}) \circ (\sigma_{j_1} \circ \cdots \circ \sigma_{j_q})$, 
where 
$n \geq i_1 > \cdots > i_p \geq 0$, $0 \leq j_1 < \dots < j_q < m$, and $m + p = n + q$.  
Each $\alpha \in \Mor([m],[n])$ determines a linear transformation $\R^{m+1} \ra \R^{n+1}$ which restricts to a map 
$\Delta^\alpha: \dpm \ra$ $\dpn$.  Thus there is a functor 
$\Delta^?:\bDelta \ra$ $\bTOP$ that sends $[n]$ to $\dpn$ and $\alpha$ to $\Delta^\alpha$.  
Since the objects of $\bDelta$ are themselves small categories, there is also an inclusion 
$\iota:\bDelta \ra$ $\bCAT$.

\label{13.93}
Given a category \bC, write 
\bSIC 
\index{\bSIC} 
for the functor category $[\bDelta^\OP,\bC]$
and 
\bCOSIC 
\index{\bCOSIC} 
for the functor category $[\bDelta,\bC]$
$-$then by definition, a 
\un{simplicial object}
\index{simplicial object} 
in \bC is an object in \bSIC and a 
\un{cosimplicial object}
\index{cosimplicial object} 
in \bC is an object in \bCOSIC.  
Example: $Y_{\bDelta} \equiv \Delta$ is a cosimplicial object in $\widehat{\bDelta}$.

Specialize to \bC = \bSET $-$then an object in \bSISET is called a 
\un{simplicial set}
\index{simplicial set} 
and a morphism in \bSISET is called a 
\un{simplicial map}.
\index{simplicial map}  Given a simplicial set $X$, put $X_n =$ $X([n])$, so for 
$\alpha:[m] \ra [n]$, $X \alpha : X_n \ra X_m$.  If 
$
\begin{cases}
\ d_i = X \delta_i\\
\ s_i = X \sigma_i
\end{cases}
, \ 
$
then $d_i$ and $s_i$ are connected by the 
\un{simplicial identities}: 
\index{simplicial identities} 
\[
\begin{cases}
\ d_i  \circ d_j = d_{j-1} \circ d_i \hspace{0.5cm} (i < j)\\
\ s_i  \circ s_j = s_{j+1} \circ s_i \hspace{0.6cm} (i \leq j)
\end{cases}
\hspace{-.2cm}, \hspace{.75cm}
d_i \circ s_j = 
\begin{cases}
\ s_{j-1} \circ d_i  \hspace{.4cm} (i < j)\\
\ \id \hspace{1.5cm} (i = j \text{ or } i = j+1)\\
s_j \circ d_{i-1}\hspace{0.53cm}  (i > j+1)
\end{cases}
\hspace{-.2cm}.
\]
The \un{simplicial standard $n$-simplex}
\index{simplicial standard $n$-simplex} 
is the simplicial set $\dn = \Mor(-,[n])$, i.e., $\dn$ is the result of applying $\Delta$ to $[n]$, so for 
$\alpha:[m] \ra [n]$, $\Delta[\alpha]:\Delta[m] \ra \Delta[n]$.  
Owing to the Yoneda lemma, if $X$ is a simplicial set 
and if $x \in X_n$, then there exists one and only one simplicial map $\Delta_x:\dn \ra X$ that takes $\id_{[n]}$ to x.  
\bSISET is complete and cocomplete, wellpowered and cowellpowered.

Let $X$ be a simplicial set $-$then one writes $x \in X$ when one means $x \in \bigcup\limits_n X_n$.  
With this understanding, an $x \in X$ is said to be 
\un{degenerate}
\index{degenerate (element of a simplicial set)} 
if there exists an epimorphism $\alpha \neq \id$ and a $y \in X$ such that $x = (X\alpha)y$; otherwise, $x \in X$ is said to be
\un{nondegenerate}.
\index{nondegenerate (element of a simplicial set)}  
The elements of $X_0$ (= 
\un{vertexes}
\index{vertexes (of a simplicial set)} 
of X) are nondegenerate.  
Every $x \in X$ admits a unique representation $x = (X\alpha)y$ where $\alpha$ is an epimorphism and y is nondegenerate.  
The nondegenerate elements in $\dn$ are the monomorphisms $\alpha:[m] \ra [n]$ $(m \leq n)$.

A 
\un{simplicial subset}
\index{simplicial subset} 
of a simplicial set $X$ is a simplicial set $Y$ such that $Y$ is a subfunctor of $X$, i.e., $Y_n \subset X_n$ 
for all $n$ and the inclusion $Y \ra X$ is a simplicial map.  
Notation: $Y \subset X$.  The  
\un{$n$-skeleton}
\index{n-skeleton (simplicial set)} 
of a simplicial set $X$ is the simplicial subset 
$X^{(n)}$ $(n \geq 0)$ of $X$ defined by stipulating that $X_p^{(n)}$ is the set of all $x \in X_p$ for which there exists an 
epimorphism $\alpha: [p] \ra [q]$ $(q \leq n)$ and a $y \in X_q$ such that $x = (X\alpha)y$.  
Therefore
$X_p^{(n)} = X_p$ $(p \leq n)$; furthermore, $X^{(0)} \subset X^{(1)} \subset \cdots$ and 
$X = \colimx X^{(n)}$.  
A proper simplicial subset of $\dn$ is contained in $\dn^{(n-1)}$, the 
\un{frontier}
\index{frontier (simplicial subset)} 
$\ddn$ of $\dn$.  
Of course, 
%%----------------------------------------------------------------------------------------------18
\label{13.3}
$\ddz = \emptyset$.  $X^{(0)}$ is isomorphic to $X_0 \cdot \dz$.  
In general, let $X_n^{\#}$ be the set of nondegenerate elements of $X_n$.  
Fix a collection $\{\dn_x:x \in X_n^{\#}\}$ of simplicial standard $n$-simplexes indexed by $X_n^{\#}$ 
$-$then the simplicial maps 
$\Delta_x:\dn \ra X$ $(x \in X_n^{\#})$ determine an arrow 
$X_n^{\#}\cdot \dn \ra X^{(n)}$ and the commutative diagram
\begin{tikzcd}[ sep=large]
{X_n^{\#}\cdot \ddn} \ar{d} \ar{r} &{X^{(n-1)}} \ar{d}\\
{X_n^{\#}\cdot \dn} \ar{r} &{X^{(n)}}
\end{tikzcd}
is a pushout square.  
\label{13.28}
Note too that $\ddn$ is a coequalizer:  Consider the diagram
\[
\coprod\limits_{0 \leq i  < j \leq n} \Delta[n-2]_{i,j} \overset{u}{\underset{v}{\rightrightarrows}}
\coprod\limits_{0 \leq i \leq n} \Delta[n-1]_i,
\]
where $u$ is defined by the $\Delta[\delta_{j-1}^{n-1}]$ and $v$ is defined by the $\Delta[\delta_i^{n-1}]$ $-$then the 
$\Delta[\delta_i^n]$ define a simplicial map 
$f: \ds\coprod\limits_{0 \leq i \leq n} \Delta[n-1]_i \ra \dn$ that induces an isomorphism 
$\coeq(u,v) \ra$ $\ddn$.\\

\label{13.84}
\label{13.84a}
\label{13.89}
\begingroup%%------------------------------------>>%%----------------------------------->>
\fontsize{9pt}{11pt}\selectfont
Call $\bdn$ the full subcategory of $\bDelta$ whose objects are the $[m]$ $(m \leq n)$.  
Given a category \bC, 
denote by $\bSIC_n$ the functor category $[\bdn^\OP,\bC]$.  The objects of $\bSIC_n$ are the 
``$n$-truncated simplicial objects'' in \bC.  
Employing the notion of the Kan extension theorem, take for $K$ the inclusion 
$\bdn^\OP \ra \bDelta^\OP$ and write $\tr^{(n)}$
\index{$\tr^{(n)}$} 
in place of $[K,\bC]$, so 
$\tr^{(n)}:\bSIC \ra \bSIC_n$.  
If \bC is complete (cocomplete), then $\tr^{(n)}$ has a left (right) adjoint $\sk^{(n)}$ ($\cosk^{(n)}$). 
\index{$\sk^{(n)}$} 
\index{$\cosk^{(n)}$}
Put 
$sk^{(n)} = \sk^{(n)} \circ \tr^{(n)}$ (the 
\un{$n$-skeleton}),
\index{n-skeleton} 
$cosk^{(n)} = \cosk^{(n)} \circ \tr^{(n)}$ (the 
\un{$n$-coskeleton}).
\index{n-skeleton}  
Example: Let \bC = \bSET $-$then for any simplicial set $X$, $sk^{(n)} X \approx X^{(n)}$.\\
\endgroup%%------------------------------------<< %%------------------------------------<<

\label{5.0c}
\indent\indent (Geometric Realizations) \ 
The realization functor $\bGamma_{\Dq}$ is a functor 
$\bSISET \ra \bTOP$ such that $\bGamma_{\Dq} \circ \Delta = \Dq$.  
It assigns to a simplicial set $X$ a topological space 
$\aX = \ds\int^{[n]} X_n \cdot \dpn$, the 
\un{geometric realization}
\index{geometric realization (of a simplicial set)} 
of $X$, and to a simplicial map 
$f:X \ra Y$ a continuous function $\abs{f}:\aX \ra \aY$, the 
\un{geometric realization}
\index{geometric realization (of a simplicial map)} 
of $f$.  In particular, 
$\abs{\dn} = \dpn$ and $\abs{\Delta[\alpha]} = \Delta^\alpha$.  
There is an explicit description of $\aX$: 
Equip $X_n$ with the discrete topology and $X_n \times \dpn$ with the product topology $-$then $\aX$ can be identified with the quotient 
$\coprod\limits_n X_n \times \dpn / \sim$, the equivalence relation being generated by writing 
$((X\alpha)x,t) \sim (x,\Delta^\alpha t)$.  
These relations are respected by every simplicial map 
$f:X \ra$ $Y$.  
Denote by $[x,t]$ the equivalence class corresponding to $(x,t)$.  The projection $(x,t) \ra [x,t]$ of 
$\coprod\limits_n X_n \times \dpn$ onto $\aX$ restricts to a map 
$\coprod\limits_n X_n^{\#} \times \mdpn \ra \aX$ that is in fact a set theoretic bijection.  Consequently, if we attach to each 
$x \in X_n^{\#}$ the subset $e_x$ of $\aX$ consisting of all $[x,t]$ $(t \in \mdpn)$, then the collection 
$\{e_x: x \in X_n^{\#} (n \geq 0)\}$ partitions $\aX$.  It follows from this that a simplicial map $f:X \ra Y$ is 
injective (surjective) iff its geometric realization $\abs{f}:\aX \ra \aY$ is injective (surjective).  Being a left adjoint, the functor 
$\aq:\bSISET \ra \bTOP$ preserves colimits.  So, e.g., by taking the geometric realization of 
%%----------------------------------------------------------------------------------------------19
the diagram
\[
\coprod\limits_{0 \leq i  < j \leq n} \Delta[n-2]_{i,j} \overset{u}{\underset{v}{\rightrightarrows}}
\coprod\limits_{0 \leq i \leq n} \Delta[n-1]_i,
\]
and unraveling the definitions, one finds that $\abs{\ddn}$ can be identified with $\ddpn$.

\label{13.5}
[Note: \ It is also true that the arrow 
$\abs{\dm \times \dn} \ra \abs{\dm} \times \abs{\dn}$ associated with the geometric realization of the projections 
$
\begin{cases}
\ p_m:\dm \times \dn \ra \dm\\
\ p_n:\dm \times \dn \ra \dn
\end{cases}
$
is a homeomorphism but this is not an a priori property of $\aq$.]

\indent\indent (Singular Sets) The singular functor $S_{\Dq}$ is a functor $\bTOP \ra \bSISET$ 
that assigns to a topological space $X$ a simplicial set $\sin X$, the 
\un{singular set}
\index{singular set} 
of $X$: 
$\sin X([n]) =$ $\sin_n X$ $=$ $C(\dpn,X)$.  $\aq$ is a left adjoint for $\sin$.  
The arrow of adjunction $X \ra \sin \aX$ 
sends $x \in X_n$ to $\abs{\Delta_x} \in C(\dpn,\aX)$, where $\abs{\Delta_x}(t) = [x,t]$; it is a monomorphism.  
The arrow of adjunction $\abs{\sin X} \ra X$ sends $[x,t]$ to $x(t)$; it is an epimorphism.\\

\begingroup%%------------------------------------>>%%----------------------------------->>
\fontsize{9pt}{11pt}\selectfont
There is a functor $T$ from \bSIAB to the category of chain complexes of abelian groups: Take an $X$ and let 
$TX$ be 
$X_0 \overset{\partial}{\la}  X_1 \overset{\partial}{\la} X_2 \overset{\partial}{\la} \cdots$, 
where 
$\partial = \ds\sum\limits_0^n (-1)^id_i$ $(d_i:X_n \ra X_{n-1})$.  
That $\partial \circ \partial = 0$ is implied by the simplicial identities.  
One can then apply the homology functor $H_*$ and end up in the category of graded abelian groups.  
On the other hand, the forgetful functor $\bAB \ra \bSET$ has a left adjoint $F_{ab}$ that sends a set $X$ 
to the free abelian group $F_{ab}X$ on $X$.  Extend it to a functor $F_{ab}:\bSISET \ra \bSIAB$.  
In this terminology, the singular homology $H_*(X)$ of a topological space $X$ is $H_*(TF_{ab}(\sin X))$.\\
\endgroup%%------------------------------------<< %%------------------------------------<<

\indent\indent (Categorical Realizations) \ The realization functor $\Gamma_\iota$ 
\index{realization functor}
\index{realization functor! $\Gamma_\iota$ - realization functor} 
is a functor 
$\bSISET \ra \bCAT$ such that $\Gamma_\iota \circ \Delta = \iota$.  
It assigns to a simplicial set $X$ a small category 
$cX = \ds\int^{[n]} X_n \cdot [n]$ called the 
\un{categorical realization}
\index{categorical realization} 
of $X$.  
In particular, $c \Delta[n] = [n]$.   
In general, $cX$ can be represented as a quotient category $CX / \sim$.  
Here, $CX$ is the category whose objects are the elements of $X_0$ and whose morphisms are the finite sequences $(x_1, \ldots, x_n)$ of elements of $X_1$ such that 
$d_0x_i = d_1x_{i+1}$.  
Composition is concatenation and the empty sequences are the identities.  
The relations are 
$s_0x = \id_x$ $(x \in X_0)$ and $(d_0x) \circ (d_2x) = d_1x$ $(x \in X_2)$.

\indent\indent (Nerves) \ The singular functor $S_\iota$ is a functor $\bCAT \ra \bSISET$ that assigns to a small category \bC a simplicial set $\ner \bC$, the 
\un{nerve}
\index{nerve (of \bC in \bCAT)} 
of \bC: $\ner \bC([n]) = \nersub_n \bC$, the set of all diagrams in \bC of the form 
$X_0 \overset{f_0}{\lra} X_1 \ra \cdots \ra X_{n-1} \overset{f_{n-1}}{\lra} X_n$.  
Therefore, 
$\nersub_0 \bC = \Ob\bC$ 
and 
$\nersub_1 \bC = \Mor\bC$.  
$c$ is a left adjoint for $\ner$.  
Since $\ner$ is full and faithful, the arrow of adjunction 
$c \circ \ner \ra \id_{\bCAT}$ is a natural isomorphism.  
\label{13.107}
The 
\un{classifying space}
\index{classifying space (of \bC in \bCAT)} 
of \bC is the geometric realization of its nerve: 
$B\bC \equiv \abs{\ner \bC}$.  
Example: $B\bC \approx B\bC^\OP$.\\

\begingroup%%------------------------------------>>%%----------------------------------->>
\fontsize{9pt}{11pt}\selectfont
The composite $\Pi = \pi_1 \circ c$ is a functor $\bSISET \ra \bGRD$ that sends a simplicial set $X$ to its 
\un{fundamental groupoid}
\index{fundamental groupoid} 
$\Pi X$.  
Example: If $X$ is a topological space, then $\Pi X \approx \Pi (\sin X)$.\\
\endgroup%%------------------------------------<< %%------------------------------------<<

%%----------------------------------------------------------------------------------------------20
Let \bC be a small category.  
Given a cofunctor $F:\bC \ra \bSET$, the 
\un{Grothendieck} 
\un{construction}
\index{Grothendieck construction} 
on $F$ is the category $\gro_{\bC}F$ 
\index{$\gro_{\bC}F$} 
whose objects are the pairs $(X,x)$, where $X$ is an object in \bC 
with $x \in FX$, and whose morphisms are the arrows $f:(X,x) \ra (Y,y)$, where $f:X \ra Y$ is a morphism in \bC, with 
$(Ff)y = x$.  
Denoting by $\pi_F$ the projection $\gro_{\bC}F \ra \bC$, if \bS is cocomplete, then for any 
$T \in \Ob[\bC,\bS]$, 
$\Gamma_TF \approx$ 
$\colim(\gro_{\bC}F \overset{\pi_F}{\lra}$ 
$\bC \overset{T}{\lra} \bS)$.  
In particular, 
$F \approx \colim(\gro_{\bC}F \overset{\pi_F}{\lra}$ 
$\bC \overset{Y_{\bC}}{\lra}$ $\widehat{\bC})$.

[Note: \ The Grothendieck construction on a functor $F:\bC \ra \bSET$ is the category $\gro_{\bC}F$ 
whose objects are the pairs $(X,x)$, where $X$ is an object in \bC with $x \in FX$ and 
whose morphisms are the arrows $f:(X,x) \ra (Y,y)$, where $f:X \ra Y$ is a morphism in \bC with $(F f)x = y$.  
Example: $\gro_{\bC}\Mor(X,-) \approx X\backslash \bC$.]\\

\begingroup%%------------------------------------>>%%----------------------------------->>
\fontsize{9pt}{11pt}\selectfont
Let $\gamma:\bC \ra \bCAT$ be the functor that sends $X$ to $\bC/X$ 
$-$then the realization functor $\Gamma_\gamma$ assigns to each $F$ in $\widehat{\bC}$ its Grothendieck construction, i.e., 
$\Gamma_\gamma F \approx \gro_{\bC} F$.\\
\endgroup%%------------------------------------<< %%------------------------------------<<

A full, isomorphism closed subcategory \bD of a category \bC is said to be a 
\un{reflective}
\index{reflective (full, isomorphism closed subcategory)}
(\un{coreflective})
\index{coreflective (full, isomorphism closed subcategory)}
subategory of \bC if the inclusion $\bD \ra \bC$ has a left (right) adjoint \mR, a 
\un{reflector}
\index{reflector (of a full, isomorphism closed subcategory)}
(\un{coreflector})
\index{coreflector (of a full, isomorphism closed subcategory)}
for \bD.

[Note: \ A full subcategory \bD of a category \bC is 
\un{isomorphism closed}
\index{isomorphism closed (full subcategory)} 
provided that every object in \bC which is isomorphic to an object in \bD is itself an object in \bD.]\\

\begingroup%%------------------------------------>>%%----------------------------------->>
\fontsize{9pt}{11pt}\selectfont
\bSET has precisely three (two) reflective (coreflective) subcategories.  
\bTOP has precisely two reflective subcategories whose intersection is not reflective.  
The full subcategory of \bGR whose objects are the finite groups is not a reflective subcategory of \bGR.\\
\endgroup%%------------------------------------<< %%------------------------------------<<

\label{8.25} %dmc mnft
Let \bD be a reflective subcategory of \bC, $R$ a reflector for \bD $-$then one may attach to each $X \in \Ob\bC$ 
a morphism $r_X:X \ra RX$ in \bC with the following property: Given any $Y \in \Ob\bD$ and any morphism 
$f:X \ra Y$ in \bC, there exists a unique morphism $g:RX \ra Y$ in \bD such that $f = g \circ r_X$.  
If the $r_X$ are epimorphisms, then \bD is said to be an 
\un{epireflective}
\index{epireflective (subcategory)} 
subcategory of \bC.

[Note: \ If the $r_X$ are monomorphisms, then the $r_X$  are epimorphisms, so 
``monocoreflective'' $\implies$ ``epireflective''.]

A reflective subcategory \bD of a complete (cocomplete) category \bC is complete (cocomplete).

\label{13.4}
\label{13.101}
[Note: \ Let $\Delta:\bI \ra \bD$ be a diagram in \bD.

\indent\indent (1) To calculate a limit of $\Delta$, postcompose $\Delta$ with the inclusion $\bD \ra \bC$ and let 
$\{\ell_i:L \ra \Delta_i\}$ be its limit in \bC $-$then $L \in \Ob\bD$ and $\{\ell_i:L \ra \Delta_i\}$ is a limit of $\Delta$.

%%----------------------------------------------------------------------------------------------21
\indent\indent (2) To calculate a colimit of $\Delta$, postcompose $\Delta$ with the inclusion $\bD \ra \bC$ and let 
$\{\ell_i:\Delta_i \ra L\}$ be its colimit in \bC $-$then $\{r_L\circ \ell_i:\Delta_i \ra RL\}$ is a colimit of $\Delta$.]\\

\index{Theorem: Epireflective Characterization Theorem}
\index{Epireflective Characterization Theorem}
\label{1.19}
\label{8.10}
\textbf{\small EPIREFLECTIVE CHARACTERIZATION THEOREM} \quadx
If a category \bC is complete, wellpowered, and cowellpowered, then a full, isomorphism closed subcategory \bD of \bC is an epireflective subcategory of \bC iff \bD is closed under the formation in \bC of products and extremal monomorphisms.

[Note: \ Under the same assumptions on \bC, the intersection of any conglomerate of epireflective subcategories is epireflective.]\\

\label{1.20}
\label{1.22}
A full, isomorphism closed subcategory of \bTOP (\bHAUS) is an epireflective subcategory iff it is closed under the formation in \bTOP (\bHAUS) of products and embeddings (products and closed embeddings).

\indent\indent ($hX$) \bHAUS is an epireflective subcategory of \bTOP.  \ The reflector sends $X$ to its 
\un{maximal Hausdorff quotient}
\index{maximal Hausdorff quotient}
\index{$hX$} 
$hX$.

\label{19.4}
\indent\indent ($crX$) The full subcategory of \bTOP whose objects are the completely regular Hausdorff spaces is an epireflective subcategory of \bTOP.  The reflector sends $X$ to its 
\un{complete regularization}
\index{complete regularization} 
$crX$. 
\index{$crX$.}

\indent\indent ($\beta X$) The full subcategory of \bHAUS whose objects are the compact spaces is an epireflective subcategory of \bHAUS. 
Therefore the category of compact Hausdorff spaces is an epireflective subcategory of the category of completely regular Hausdorff spaces and the reflector sends \mX to $\beta X$,  
\index{$\beta X$} the 
\un{Stone-\u Cech compactification}
\index{Stone-\u Cech compactification} 
of \mX.

[Note: \ If \mX is Hausdorff, then $\beta(cr X)$ is its compact reflection.]

\indent\indent ($\nu X$) The full subcategory of \bHAUS whose objects are the $\R$-compact spaces is an epireflective subcategory of \bHAUS. 
Therefore the category of $\R$-compact spaces is an epireflective subcategory of the category of completely regular Hausdorff spaces and the reflector sends \mX to $\nu X$  
\index{$\nu X$} , the 
\un{$\R$-compactification}
\index{Real ($\R$)-compactification} 
of \mX.

[Note: \ If \mX is Hausdorff, then $\nu(cr X)$ is its $\R$-compact reflection.]\\

\begingroup%%------------------------------------>>%%----------------------------------->>
\fontsize{9pt}{11pt}\selectfont
A full, isomorphism closed subcategory of \bGR or \bAB is an epireflective subcategory iff it is closed under the formation of products and subgroups.  
Example: \bAB is an epireflective subcategory of \bGR, the reflector sending $X$ to its abelianization $X/[X,X]$.\\
\endgroup%%------------------------------------<< %%------------------------------------<<

\label{1.18}
If \bC is a full subcategory of \bTOP (\bHAUS), then there is a smallest epireflective subcategory of \bTOP (\bHAUS) containing \bC, the \un{epireflective hull}
\index{epireflective hull} of \bC.  
If $X$ is a topological space (Hausdorff topological space), then $X$ is an object in the epireflective hull of 
%%----------------------------------------------------------------------------------------------22
\bC in \bTOP (\bHAUS) iff there exists a set $\{X_i\} \subset \Ob\bC$ and an extremal monomorphism 
$f:X \ra \prod\limits_i X_i$.\\

\begingroup%%------------------------------------>>%%----------------------------------->>
\fontsize{9pt}{11pt}\selectfont
The epireflective hull in \bTOP (\bHAUS) of $[0,1]$ is the category of completely regular Hausdorff spaces (compact Hausdorff spaces).  
The epireflective hull in \bTOP of $[0,1]/[0,1[$ is the full subcategory of \bTOP whose objects satisfy the $\tT_0$ separation axiom.  
The epireflective hull in \bTOP (\bHAUS) of $\{0,1\}$ (discrete topology) is the full subcategory of \bTOP (\bHAUS)  whose objects are the zero dimensional Hausdorff spaces (zero dimensional compact Hausdorff spaces).  
The epireflective hull in \bTOP of $\{0,1\}$ (indiscrete topology) is the full subcategory in \bTOP whose objects are the indiscrete spaces.
\\ \indent
[Note: \ Let $E$ be a nonempty Hausdorff space $-$then a Hausdorff space  $X$ is said to be 
\un{$E$-compact}
\index{E-compact (E a nonempty Hausdorff space)} %\index{$E$-compact (E a nonempty Hausdorff space)}
provided that $X$ is in the epireflective hull of $E$ in 
\bHAUS.  
Example: A Hausdorff space is $\N$-compact iff it is $\Q$-compact iff it is $\PP$-compact.  There is no $E$ such that every Hausdorff space is $E$-compact.  In fact, given $E$, there exists a Hausdorff space $X_E$ with $\#(X_E) > 1$ 
such that every element of $C(X_E,E)$ is a constant.]\\
\endgroup%%------------------------------------<< %%------------------------------------<<

A morphism $f:A \ra B$ and an object $X$ in a category \bC are said to be 
\un{orthogonal}
\index{orthogonal (morphism and an arrow in a category)} 
$(f \perp X)$ 
\index{$(f \perp X)$} 
if the precomposition arrow $f^*:\Mor(B,X) \ra \Mor(A,X)$ is bijective.  Given a class $S \subset \Mor \bC$, $S^\perp$ 
is the class of objects orthogonal to each $f \in S$ and given a class $D \subset \Ob\bC$, $D^\perp$ is the class of morphisms orthogonal to each $X \in D$.  
\label{15.33}
One then says that a pair $(S,D)$ is an 
\un{orthogonal pair}
\index{orthogonal pair (class of morphisms, class of objects)} 
provided that $S = D^\perp$, and $D = S^\perp$.  
Example: Since $\perp\perp\perp = \perp$, for any $S$, $(S^{\perp\perp},S^\perp)$ is an orthogonal pair, and for any 
$D$, $(D^\perp,D^{\perp\perp})$ is an orthogonal pair.

\label{9.50}
[Note: \ Suppose that $(S,D)$ is an orthogonal pair $-$then 
(1) \ $S$ contains the isomorphisms of \bC;
(2) \ $S$ is closed under composition; 
(3) \ $S$ is 
\un{cancellable}
\index{cancellable (class of morphisms per an orthogonal pair)}, 
i.e., 
$g \circ f \in S$ $\&$ $f \in S$ $\implies$ $g \in S$ and 
$g \circ f \in S$ $\&$ $g \in S$ $\implies$ $f \in S$.  
In addition, if 
\begin{tikzcd}[ sep=large]
{A} \ar{d}[swap]{f} \ar{r} &{A^\prime} \ar{d}{f^\prime}\\
{B} \ar{r} &{B^\prime}
\end{tikzcd}
is a pushout square, then $f \in S$ $\implies$ $f^\prime \in S$, and if $\Xi \in \Nat(\Delta,\Delta^\prime)$, where 
$\Delta,\Delta^\prime:\bI \ra \bC$, then $\Xi_i \in S$ $(\forall \ i)$ $\implies$  $\colimx \Xi_i \in S$ 
(if $\colimx \Delta, \colimx \Delta^\prime$ exist).]

\label{9.1}
Every reflective subcategory \bD of \bC generates an orthogonal pair.  Thus, with $R:\bC \ra \bD$ the reflector, put 
$T = \iota \circ R$, where $\iota:\bD \ra \bC$ is the inclusion, and denote by $\epsilon:\id_{\bC} \ra T$ the associated natural transformation.  Take for $S \subset \Mor \bC$ the class consisting of those $f$ such that $Tf$ is an isomorphism and take for 
$D \subset \Ob\bC$ the object class of \bD, i.e., the class consisting of those $X$ such that $\epsilon_X$ is an isomorphism $-$then $(S,D)$ is an orthogonal pair.\\

\begingroup%%------------------------------------>>%%----------------------------------->>
\fontsize{9pt}{11pt}\selectfont
A full, isomorphism closed subcategory \bD of a category \bC is said to be an 
\un{orthogonal}
\index{orthogonal (full, isomorphism closed subcategory)} 
subcategory
%%----------------------------------------------------------------------------------------------23
of \bC if $\Ob\bD = S^\perp$ for some class $S \subset \Mor \bC$.  If \bD is reflective, then \bD is orthogonal but the converse is false (even in \bTOP).
\\ \indent
[Note: \ Let $(S,D)$ be an orthogonal pair.  
Suppose that for each $ X \in \Ob\bC$ there exists a morphism 
$\epsilon_X:X \ra TX$ in $S$, where $TX \in D$ $-$then for every $f:A \ra B$ in $S$ 
and for every $g:A \ra X$ there exists a unique $t:B \ra TX$ such that $\epsilon_X \circ g = t \circ f$.  
So, for any arrow $X \ra Y$, there is a commutative diagram 
\begin{tikzcd}[ sep=large]
{X} \ar{d} \ar{r}{\epsilon_X} &{TX} \ar{d}\\
{Y} \ar{r}[swap]{\epsilon_Y} &{TY}
\end{tikzcd}
, thus $T$ defines a functor $\bC \ra \bC$ and $\epsilon:\id_{\bC} \ra T$ is a natural transformation.  
Since 
$\epsilon T = T \epsilon$ is a natural isomorphism, it follows that $S^\perp = D$ 
is the object class of a reflective subcategory of \bC.]\\
\endgroup%%------------------------------------<< %%------------------------------------<<

\indent\indent ($\kappa-$DEF) \ Fix a regular cardinal $\kappa$ $-$then an object \mX in a cocomplete category \bC is said to be 
\un{$\kappa$-definite}
\index{kappa-definite, $\kappa$-definite} 
provided that $\forall$ \ regular cardinal 
$\kappa^\prime \geq \kappa$, $\Mor(X,-)$ preserves colimits over $[0,\kappa^\prime[$, so for every diagram 
$\Delta:[0,\kappa^\prime[ \ra \bC$, the arrow 
$\colim\  \Mor(X,\Delta_\alpha) \ra$ $\Mor(X,\colim\ \Delta_\alpha)$ is bijective.\\

\begingroup%%------------------------------------>>%%----------------------------------->>
\fontsize{9pt}{11pt}\selectfont
Given a group \mG, there is a $\kappa$ for which \mG is $\kappa-$definite and all finitely presented groups are 
$\omega-$definite.\\

\endgroup%%------------------------------------<< %%------------------------------------<<

\label{8.6}
\index{Theorem: Reflective Subcategory Theorem}
\index{Reflective Subcategory Theorem}
\textbf{\small REFLECTIVE SUBCATEGORY THEOREM} \quadx
Let \bC be a cocomplete category.  Suppose that $S_0 \subset \Mor \bC$ is a set with the property that for some $\kappa$, 
the domain and codomain of each $f \in S_0$ are $\kappa$-definite $-$then $S_0^\perp$ is the object class of a reflective subcategory of \bC.\\

\indent\indent ($P$-Localization) \ Let $P$ be a set of primes.  Let 
$S_P = \{1\} \cup \{n > 1:p \in P \implies$ $p  \not|  n\}$ $-$then 
a group $G$ is said to be 
\un{$P$-local}
\index{P-local (group)} 
if the map
$
\begin{cases}
\ G \ra G\\
\ g \ra g^n
\end{cases}
$
is bijective $\forall \ n \in S_P$.  $\bGR_P$, the full subcategory of \bGR whose objects are the $P$-local groups, is a reflective subcategory of \bGR.  
In fact, $\Ob\bGR_P = S_P^\perp$, where now $S_P$ stands for the homomorphisms 
$
\begin{cases}
\ \Z \ra \Z\\
\ 1 \ra n
\end{cases}
(n \in S_P).
$
The reflector 
$
L_P:
\begin{cases}
\ \bGR \ra \bGR_P\\
\ G \ra G_P
\end{cases}
$
is called \un{$P$-localization}. 
\index{P-localization (reflector from \bGR)}\\

\vspace{0.25cm}

\begingroup%%------------------------------------>>%%----------------------------------->>
\fontsize{9pt}{11pt}\selectfont
$P$-localization need not preserve short exact sequences.  For example 
$1 \ra A_3 \ra S_3 \ra S_3/A_3 \ra 1$, when localized at $P = \{3\}$, gives 
$1 \ra A_3 \ra 1 \ra 1 \ra 1$.\\
\endgroup%%------------------------------------<< %%------------------------------------<<

\label{13.23}
A category \bC with finite products is said to be 
\un{cartesian closed}
\index{cartesian closed}
provided that each of the functors 
$-\times Y:\bC \ra$ $\bC$ has a right adjoint $Z \ra Z^Y$, so 
$\Mor(X \times Y,Z) \approx$ $\Mor(X,Z^Y)$.  
The object $Z^Y$ is called an 
\un{exponential object}.
\index{exponential object}  
The 
\un{evaluation morphism}
\index{evaluation morphism} 
$\ev_{Y,Z}$ is the morphism 
$Z^Y \times Y \ra Z$ such that for every $f:X \times Y \ra Z$ there is a unique $g:X \ra Z^Y$ such that 
$f = \ev_{Y,Z} \circ (g \times \id_Y)$.\\

%%----------------------------------------------------------------------------------------------24
In a cartesian closed category:\\
%something is not aligning correctly - come back and fix dmcXXX
%\begin{align*}
\[
\begin{array}{lll}
&{(1)\  \ X^{Y \times Z} \approx (X^Y)^Z\text{;}} \hspace{1cm}
&{(3) \ \ X^{\coprod\limits_i Y_i} \approx \prod\limits_i (X^{Y_i})\text{;}}\\
&{(2) \ \ \left(\prod\limits_i X_i\right)^Y \approx \prod\limits_i (X_i^Y)\text{;}} \hspace{1cm}
&{(4) \ \ X \times \left(\coprod\limits_i Y_i\right) \approx \coprod\limits_i (X \times Y_i).}
\end{array}
\]
%\end{align*}

\vspace{0.05cm}

\begingroup%%------------------------------------>>%%----------------------------------->>
\fontsize{9pt}{11pt}\selectfont
\bSET is cartesian closed but $\bSET^\OP$ is not cartesian closed.  \bTOP is not cartesian closed but does have full, cartesian closed subcategories, e.g., the category of compactly generated Hausdorff spaces.
\\ \indent
[Note: \ If \bC is cartesian closed and has a zero object, then \bC is equivalent to $\bone$.  Therefore neither 
$\bSET_*$ nor $\bTOP_*$ is cartesian closed.]
\\ \indent
\bCAT is cartesian closed: $\Mor(\bC \times \bD,\bE) \approx \Mor(\bC,\bE^{\bD})$, where
$\bE^{\bD} = [\bD,\bE]$.  
\label{13.1}
\bSISET is cartesian closed: 
$\Nat(X \times Y,Z) \approx$ 
$\Nat(X,Z^Y)$, where $Z^Y([n]) = \Nat(Y \times \dn,Z)$.
\\ \indent
[Note: \ The functor $\ner: \bCAT \ra \bSISET$ preserves exponential objects.]\\
\endgroup%%------------------------------------<< %%------------------------------------<<

\quadx A 
\un{monoidal category} 
\index{monoidal category} 
is a category \bC equipped with a functor $\otimes: \bC \times \bC \ra \bC$ (the \un{multiplication} and an object $e \in \text{Ob}\bC$ (the \un{unit}), together with natural isomorphisms $R$, $L$, and $A$, where 
$
\begin{cases}
\ R_X: X \otimes e \ra X\\
\ L_X: e \otimes X \ra X
\end{cases}
$
and $A_{X,Y,Z}:X \otimes (Y\otimes Z) \ra (X \otimes Y) \otimes Z$, subject to the following assumptions.
\label{16.52}

\indent\indent  (MC$_1$) The diagram
\[
\begin{tikzcd}[ sep=large]
X \otimes (Y \otimes (Z \otimes W)) \arrow[r, "A"]   \arrow{d}[swap]{\id \otimes A} 
&(X \otimes Y) \otimes (Z \otimes W) \arrow[r, "A"] & ((X \otimes Y) \otimes Z) \otimes W \\
X \otimes ((Y \otimes Z) \otimes W) \ar{rr}[swap]{A}												
&& (X \otimes (Y \otimes Z)) \otimes W \ar{u}[swap]{A \otimes \id}
\end{tikzcd}
\]
commutes.

\indent\indent (MC$_2$) The diagram
\[
\begin{tikzcd}[sep=large]
{X \otimes ( e \otimes Y)}  \arrow[rr, "A"]   \arrow{d}[swap]{\id \otimes L} 
&& (X \otimes e) \otimes Y\arrow[d,"R \otimes \id"]\\
{X \otimes Y} \arrow[rr,shift right=0.5,dash] \arrow[rr,shift right=-0.5,dash] 										
&& {X \otimes Y}
\end{tikzcd}
\]
commutes.

[Note: \ The ``coherency'' principle then asserts that ``all'' diagrams built up from instances of $R$, $L$, $A$ (or their inverses), 
and is by repeated application of $\otimes$ necessarily commute.  In particular, the diagrams
\[
\begin{tikzcd}[sep=large]
{e \otimes ( X \otimes Y)}  \arrow[rr, "A"]   \arrow{d}[swap]{L} 
&& {(e \otimes X) \otimes Y} \arrow[d,"L \otimes \id"]\\
{X \otimes Y} \arrow[rr,shift right=0.5,dash] \arrow[rr,shift right=-0.5,dash] 										
&& {X \otimes Y}
\end{tikzcd}
\quadx
\begin{tikzcd}[sep=large]
{X \otimes ( Y \otimes e)}  \arrow[rr, "A"]   \arrow{d}[swap]{\id \otimes R} 
&&{ (X \otimes Y) \otimes e} \arrow[d,"R"]\\
{X \otimes Y} \arrow[rr,shift right=0.5,dash] \arrow[rr,shift right=-0.5,dash] 										
&& {X \otimes Y}
\end{tikzcd}
\]
%%----------------------------------------------------------------------------------------------25
commute and $L_\epsilon = R_\epsilon:e \otimes e \ra e$.]\\

\begingroup%%------------------------------------>>%%----------------------------------->>
\fontsize{9pt}{11pt}\selectfont
Any category with finite products (coproducts) is monoidal: 
Take $X \otimes Y$ to be $X \prod Y$ ($X \coprod Y$) and 
let $e$ be a final (initial) object.  
The category \bAB is monoidal: 
Take $X \otimes Y$ to be the tensor product and let $e$ be $\Z$.  
The category $\bSET_*$ is monoidal: 
Take $X \otimes Y$ to be the smash product $X \# Y$ and let $e$ be the two point set.\\
\endgroup%%------------------------------------<< %%------------------------------------<<

A \un{symmetry} 
\index{symmetry (monoidal category)} 
for a monoidal category \bC is a natural isomorphism  
$\Tee$, where $\Tee_{X,Y}:X \otimes Y \ra Y \otimes X$, such that 
$\Tee_{Y,X} \circ \Tee_{X,Y}:X \otimes Y \ra X \otimes Y$ is the identitiy, 
$R_X = L_X \circ \Tee_{X,e}$, and the diagram
\[
\begin{tikzcd}[ sep=large]
{X \otimes (Y \otimes Z)}  \ar{d}[swap]{\id \otimes \Tee} \ar{r}{A}
&{(X \otimes Y) \otimes Z} \ar{r}{\Tee}
&{Z \otimes (X \otimes Y)} \ar{d}{A}\\
{X \otimes (Z \otimes Y)}   \ar{r}[swap]{A}
&{(X \otimes Z) \otimes Y} \ar{r}[swap]{\Tee \otimes \id}
&{(Z \otimes X) \otimes Y}
\end{tikzcd}
\]
commutes.  A 
\label{16.53}
\un{symmetric monoidal category}
\index{symmetric monoidal category} 
is a monoidal category \bC endowed 
with a symmetry $\Tee$.  A monoidal category can have more than one symmetry (or none at all).  

[Note: The ``coherency'' principle then asserts that ``all'' diagrams built up from instances of $R$, $L$, $A$, $\Tee$ 
(or their inverses), and id by repeated application of $\otimes$ necessarily commute.]\\

\begingroup%%------------------------------------>>%%----------------------------------->>
\fontsize{9pt}{11pt}\selectfont

Let \bC be the category of chain complexes of abelian groups; let \bD be the full subcategory of \bC whose objects are the graded abelian groups.  \bC and \bD are both monoidal: Take $X \otimes Y$ to be the tensor product and let $e = \{e_n\}$ 
be the chain complex defined by $e_0 = \Z$ and $e_n = 0$ ($n \ne 0$).  If 
$
\begin{cases}
\ X = \{X_p\}\\
\ Y = \{Y_q\}
\end{cases}
$
and if 
$
\begin{cases}
\ x \in X_p\\
\ y \in Y_q
\end{cases}
, \ 
$
then the assignment 
$
\begin{cases}
\ X \otimes Y \ra Y \otimes X\\
\ x \otimes y \ra (-1)^{pq} (y \otimes x)
\end{cases}
$
is a symmetry for \bC and there are no others.  By contrast, \bD admits a second symmetry, namely the assignment 
$
\begin{cases}
\ X \otimes Y \ra Y \otimes X\\
\ x \otimes y \ra y \otimes x
\end{cases}
.
$
\\

\endgroup%%------------------------------------<< %%------------------------------------<<

A \un{closed category}
\index{closed category} 
is a symmetric monoidal category \bC 
with the property that each of the functors 
$-\otimes Y: \bC \ra$ $\bC$ has a right adjoint 
$Z \ra$ $\hom(Y,Z)$, so 
$\Mor(X \otimes Y,Z) \approx$ 
$\Mor(X,\hom(Y,Z))$.  The functor 
$\bC^\OP \times \bC \ra$ $\bC$ is called an 
\un{internal hom functor}.
\index{internal hom functor}  
The 
\un{evaluation morphism}
\index{evaluation morphism} 
$\ev_{Y,Z}$ is the morphism 
$\hom(Y,Z) \otimes Y \ra$ $Z$ such that for every 
$f:X \otimes Y \ra$ $Z$ there is a unique  
$g:X \ra$ $\hom(Y,Z)$ such that 
$f = \ev_{Y,Z} \circ (g \otimes \id_Y)$.  Agreeing to write $U_e$ for the functor $\Mor(e,-)$ (which need not be faithful), 
one has 
$U_e \circ \hom \approx$ $\Mor$.  Consequently, $X \approx \hom(e,X)$ and 
$\hom(X \otimes Y,Z) \approx$ $\hom(X,\hom(Y,Z))$.\\

%%----------------------------------------------------------------------------------------------26
A cartesian closed category is a closed category.  \bAB is a closed category but is not cartesian closed.\\

\begingroup%%------------------------------------>>%%----------------------------------->>
\fontsize{9pt}{11pt}\selectfont

\bTOP admits, to within isomorphism, exactly one structure of a closed category.  For let \mX and \mY be topological 
spaces $-$then their product $X \otimes Y$ is the cartesian product $X \times Y$ supplied with the final topology determined 
by the inclusions 
$
\begin{cases}
\ \{x\} \times Y \ra X \times Y\\
\ X \times \{y\} \ra X \times Y
\end{cases}
$
$(x \in X$, $y \in Y)$, the unit being the one point space.  
The associated internal hom functor $\hom(X,Y)$  sends $(X,Y)$ to 
$C(X,Y)$, where $C(X,Y)$ carries the topology of pointwise convergence.\\

\endgroup%%------------------------------------<< %%------------------------------------<<

Given a monoidal category \bC, a 
\un{monoid}
\index{monoid (in a monoidal category)} 
in \bC 
is an object $X \in \Ob\bC$ together with morphisms $m:X \otimes X \ra X$ and $\epsilon:e \ra X$ subject to the following assumptions.

\indent\indent (MO$_1$) \ The diagram
\[
\begin{tikzcd}[ sep=large]
{X \otimes (X \otimes X)} \ar{d}[swap]{\id\otimes m} \ar{r}{A} 
&{(X \otimes X) \otimes X} \ar{r}{m \otimes \id} 
&{X \otimes X} \ar{d}{m}\\
{X \otimes X} \ar{rr}[swap]{m} &&{X}
\end{tikzcd}
\]
commutes.

\indent\indent (MO$_2$) \ The diagrams
\[
\begin{tikzcd}[sep=large]
{e \otimes X} \ar{d}[swap]{L}  \ar{r}{\epsilon \otimes \id}  &{X \otimes X} \ar{d}{m}\\
{X} \arrow[r,shift right=0.5,dash] \arrow[r,shift right=-0.5,dash]  &{X}
\end{tikzcd}
\hspace{2cm} 
\begin{tikzcd}[sep=large]
{X \otimes X} \ar{d}[swap]{m}    &{X \otimes e} \ar{l}[swap]{\id \otimes \epsilon} \ar{d}{R}\\
{X} \arrow[r,shift right=0.5,dash] \arrow[r,shift right=-0.5,dash]  &{X}
\end{tikzcd}
\]
commute.

$\bMON_{\bC}$ 
\index{$\bMON_{\bC}$} 
is the category whose objects are the monoids in \bC and whose 
morphisms $(X,m,\epsilon) \ra (X^\prime,m^\prime,\epsilon^\prime)$ are the arrows $f:X \ra X^\prime$ such that 
$f \circ m = m^\prime \circ (f \otimes f)$ and $f \circ \epsilon = \epsilon^\prime$.\\

\begingroup%%------------------------------------>>%%----------------------------------->>
\fontsize{9pt}{11pt}\selectfont
$\bMON_{\bSET}$ 
\index{$\bMON_{\bSET}$} 
is the category of semigroups with unit.  
$\bMON_{\bAB}$ 
\index{$\bMON_{\bAB}$}  
is the category of rings with unit.\\

\endgroup%%------------------------------------<< %%------------------------------------<<

Given a monoidal category \bC, a 
\un{left action}
\index{left action (of a monoid)}
of a monoid \mX in \bC on an object $Y \in \Ob\bC$ is a morphism 
$l:X \otimes Y \ra Y$ such that the diagram
\[
\begin{tikzcd}[sep=large]
{X \otimes (X \otimes Y)} \ar{d}[swap]{\id \otimes l} \ar{r}{A}
&{(X \otimes X) \otimes Y} \ar{r}{m \otimes \id}
&{X \otimes Y} \ar{d}{l}
&{e \otimes Y} \ar{l}[swap]{e\otimes\id} \ar{d}{L}\\
{X \otimes Y} \ar{rr}[swap]{l}
&&{Y} \arrow[r,shift right=0.5,dash] \arrow[r,shift right=-0.5,dash] 
&{Y}
\end{tikzcd}
\]
%^
%%----------------------------------------------------------------------------------------------27
commutes.

[Note: \ The definition of a 
\un{right action}
\index{right action (of a monoid)} 
is analogous.]

$\bLACT_X$ 
\index{$\bLACT_X$} 
is the category whose objects are the left actions of \mX and whose morphisms 
$(Y,l) \ra (Y^\prime,l^\prime)$ are the arrows $f:Y \ra Y^\prime$ such that 
$f \circ l = l^\prime \circ (\id \otimes f)$.\\

\begingroup%%------------------------------------>>%%----------------------------------->>
\fontsize{9pt}{11pt}\selectfont
If \mX is a monoid in \bSET, then $\bLACT_X$ is isomorphic to the functor category $[\bX,\bSET]$, \bX the 
category having a single object $*$ with $\Mor(*,*) = X$.\\

\endgroup%%------------------------------------<< %%------------------------------------<<

A 
\un{triple}
\index{triple (in a category )} 
$\bT = (T,m,\epsilon)$ 
in a category \bC consists of a functor $T:\bC \ra \bC$ 
and natural transformations 
$
\begin{cases}
\ m \in \Nat(T \circ T,T)\\
\ \epsilon \in \Nat(\id_{\bC},T)
\end{cases}
$
subject to the following assumptions.

\indent\indent (T$_1$) \ The diagram
\[
\begin{tikzcd}[sep=large]
{T \circ T \circ T} \ar{d}[swap]{Tm} \ar{r}{mT} &{T \circ T} \ar{d}{m}\\
{T \circ T} \ar{r}[swap]{m} &{T}
\end{tikzcd}
\]
commutes.

\indent\indent (T$_2$) \ The diagrams
\[
\begin{tikzcd}[sep=large]
{T} \ar{d}[swap]{\id} \ar{r}{\epsilon T} &{T \circ T} \ar{d}{m}\\
{T} \arrow[r,shift right=0.5,dash] \arrow[r,shift right=-0.5,dash] &{T}
\end{tikzcd}
\hspace{2cm}  
\begin{tikzcd}[sep=large]
{T \circ T} \ar{d}[swap]{m}  &{T} \ar{l}[swap]{T\epsilon} \ar{d}{\id}\\
{T} \arrow[r,shift right=0.5,dash] \arrow[r,shift right=-0.5,dash] &{T}
\end{tikzcd}
\]
commute.

[Note: \ Formally, the functor category $[\bC,\bC]$ is a monoidal category: Take $F \otimes G$ to be $F \circ G$ 
and let $e$ be $\id_{\bC}$.  Therefore a triple in \bC is a monoid in $[\bC,\bC]$ 
\label{10.3}
(and a \un{cotriple} 
\index{cotriple (in a category )} 
in \bC is a monoid in $[\bC,\bC]^\OP$), a morphism of triples being a morphism in 
the metacategory $\bMON_{[\bC,\bC]}$.]

Given a triple $\bT = (T,m,\epsilon)$ in \bC, a \un{\bT-algebra} 
\index{algebra! \bT-algebra} 
is an object \mX in \bC 
and a morphism $\xi:TX \ra X$ subject to the following assumptions.\\
%^
\label{14.124}
\label{16.55}
\indent\indent (TA$_1$) \ The diagram
\[
\begin{tikzcd}[sep=large]
{T(TX)} \ar{d}[swap]{m_X} \ar{r}{T\xi}  &{TX}  \ar{d}{\xi}\\
{TX} \ar{r}[swap]{\xi} &{X}
\end{tikzcd}
\]
%%----------------------------------------------------------------------------------------------28
commutes.

\indent\indent (TA$_2$) \ The diagram
\[
\begin{tikzcd}[sep=large]
{X} \ar{d}[swap]{\id} \ar{r}{\epsilon_X}  &{TX}  \ar{d}{\xi}\\
{X}  \arrow[r,shift right=0.5,dash] \arrow[r,shift right=-0.5,dash]  &{X}
\end{tikzcd}
\]
commutes.

$\bT\text{-}\bALG$ 
\index{$\bT\text{-}\bALG$} 
is the category whose objects are the \bT-algebras and whose morphisms 
$(X,\xi) \ra (Y,\eta)$ are the arrows $f:X \ra Y$ such that $f \circ \xi = \eta \circ T f$.

[Note: \ If $\bT = (T,m,\epsilon)$ is a cotriple in \bC, then the relevant notion is 
\un{\bT-coalgebra}
\index{coalgebra! \bT-coalgebra} 
and the relevant category is $\bT\text{-}\bCOALG$ 
\index{$\bT\text{-}\bCOALG$}.]\\

\begingroup%%------------------------------------>>%%----------------------------------->>
\fontsize{9pt}{11pt}\selectfont
Take $\bC = \bAB$.  Let $A \in \Ob\bRG$.  
Define $T:\bAB \ra \bAB$ by $TX = A \otimes X$, $m \in \Nat(T \circ T,T)$ by 
$
m_X:
\begin{cases}
\ A \otimes (A \otimes X) \ra A \otimes X\\
\ a \otimes (b \otimes x) \ra ab \otimes x
\end{cases}
$
, $\epsilon \in \Nat(\id_{\bAB},T)$ by 
$
\epsilon_X:
\begin{cases}
\ X \ra A \otimes X\\
\ x \ra 1 \otimes x
\end{cases}
$
%%dmc00A_\vspace{0.1cm}
$-$then $\bT\text{-}\bALG$ is isomorphic to \bAMOD.\\

\endgroup%%------------------------------------<< %%------------------------------------<<

\label{13.19}
Every adjoint situation $(F,G,\mu,\nu)$ determines a triple in \bC, viz. $(G \circ F,G\nu F,\mu)$ 
(and a cotriple in \bD, viz. $(F \circ G, F\mu G, \nu)$).  
Different adjoint situations can determine the same triple.  
Conversely, every triple is determined by at least one adjoint situation, in general by many.  
One realization is the construction of Eilenberg-Moore:  
Given a triple $\bT = (T,m,\epsilon)$ in \bC, call $F_{\bT}$ the functor 
$\bC \ra$ $\bT\text{-}\bALG$ that sends 
$X \overset{f}{\ra}$ $Y$ to 
$(TX,m_X)$ $\overset{Tf}{\lra}$ $(TY,m_Y)$, call $G_{\bT}$ the 
functor 
$\bT\text{-}\bALG \ra$ $\bC$ that sends 
$(X,\xi) \overset{f}{\ra}$ $(Y,\eta)$ to  
$X \overset{f}{\ra} Y$, put $\mu_X = \epsilon_X$, 
and $\nu_{(X,\xi)} = \xi$ $-$then $F_{\bT}$ is a left adjoint for $G_{\bT}$ and this adjoint situation determines \bT.\\

\begingroup%%------------------------------------>>%%----------------------------------->>
\fontsize{9pt}{11pt}\selectfont
Suppose that \bC = \bSET, $\bD = \bMON_{\bSET}$.  Let $F:\bC \ra \bD$ be the functor that sends \mX to the free semigroup with unit on \mX $-$then $F$ is a left adjoint for the forgetful functor $G:\bD \ra \bC$.  The triple determined by 
this adjoint situation is $\bT = (T,m,\epsilon)$, where $T:\bSET \ra \bSET$ assigns to each \mX the set 
$TX = \ds\bigcup\limits_0^\infty X^n$, $m_X:T(TX) \ra TX$ is defined by concatenation and $\epsilon_X:X \ra TX$ by inclusion.  
The corresponding category of \bT-algebras is isomorphic to $\bMON_{\bSET}$.\\

\endgroup%%------------------------------------<< %%------------------------------------<<

Let $(F,G,\mu,\nu)$ be an adjoint situation.  
If \ $\bT = (G \circ F,G\nu F,\mu)$ is the associated triple in \bC, then the 
\un{comparison functor}
\index{comparison functor} 
$\Phi$ is the functor 
$\bD \ra \bT\text{-}\bALG$ that sends \mY to $(GY,G\nu_Y)$ and $g$ to $G g$.  
It is the only functor $\bD \ra \bT$-$\bALG$ 
for which $\Phi \circ F = F_{\bT}$ and $G_{\bT} \circ \Phi = G$.\\

\begingroup%%------------------------------------>>%%----------------------------------->>
\fontsize{9pt}{11pt}\selectfont

Consider the adjoint situation produced by the forgetful functor $\bTOP \ra$ $\bSET$ $-$then
$\bT\text{-}\bALG =$ $\bSET$ and the comparison functor $\bTOP \ra$ $\bSET$ is the forgetful functor.\\

\endgroup%%------------------------------------<< %%------------------------------------<<

%%----------------------------------------------------------------------------------------------29
Given categories
$
\begin{cases}
\ \bC\\
\ \bD
\end{cases}
\hspace{-.25cm}, \ 
$
a functor $G:\bD \ra \bC$ is said to be 
\un{monadic}
\index{monadic (functor)}
(\un{strictly monadic})
\index{strictly monadic (functor)}
provided that \mG has a left adjoint $F:\bC \ra \bD$ and the comparison functor $\Phi:\bD \ra \bT\text{-}\bALG$ 
is an equivalence (isomorphism) of categories.\\

\begingroup%%------------------------------------>>%%----------------------------------->>
\fontsize{9pt}{11pt}\selectfont
In order that \mG be monadic, it is necessary that \mG be conservative.  
So, e.g., the forgetful functor 
$\bTOP \ra \bSET$ is not monadic.  
If \bD is the category of Banach spaces and linear contractions and if 
$G:\bD \ra \bSET$ is the ``unit ball'' functor, then \mG has a left adjoint and is conservative, but not monadic.  
Theorems due to Beck, Duskin, and others lay down conditions that are necessary and sufficient for a functor to be monadic or strictly monadic.  
In particular, these results imply that if \bD is a ``finitary category of algebraic structures'', then the forgetful functor $\bD \ra \bSET$ is strictly monadic.  
Therefore the  forgetful functor from \bGR, \bRG, $\ldots$, to \bSET is strictly monadic.
\\ \indent
[Note: \ No functor from \bCAT to \bSET can be monadic.]\\
\endgroup%%------------------------------------<< %%------------------------------------<<

Among the possibilities of determining a triple $\bT = (T, m, \epsilon)$ in \bC by an adjoint situation, 
the construction of Eilenberg-Moore is ``maximal''.  The ``minimal'' construction is that of Kleisli: 
$\bKL(\bT)$ 
\index{Kleisli category}
\index{Kleisli construction}
\index{$\bKL(\bT)$} 
is the category whose objects are those of \bC, the morphisms from \mX to \mY being $\Mor(X,TY)$ with 
$\epsilon_X \in \Mor(X,TX)$ serving as the identity.  
Here, the composition of 
$
\begin{cases}
\ X \overset{f}{\ra} TY\\
\ Y \underset{g}{\ra} TZ
\end{cases}
$
in $\bKL(\bT)$ is $m_Z \circ Tg \circ f$ (calculated in \bC).  
If $K_{\bT}:\bC \ra \bKL(\bT)$ is the functor that sends 
$X \overset{f}{\ra} Y$ to 
\begin{tikzcd}[sep=large]
{X}  \ar{r}{\epsilon_Y \circ f} &{TY}
\end{tikzcd}
%$X \overset{\epsilon_Y \circ f}{\longrightarrow} TY$
and if $L_{\bT}:\bKL(\bT) \ra \bC$ is the functor that sends 
$X \overset{f}{\ra} TY$ to 
\begin{tikzcd}[sep=large]
{TX}  \ar{r}{m_Y \circ Tf} &{TY,}
\end{tikzcd}
%$TX \overset{m_Y \circ Tf}{\longrightarrow} TY$, 
then $K_{\bT}$ is a left adjoint for $L_{\bT}$ with arrows of adjunction 
$\epsilon_X$, $\id_{TX}$ and this adjoint situation determines $\bT$.

[Note: Let $G:\bD \ra \bC$ be a functor $-$then the 
\un{shape}
\index{shape (of a functor)}
of \mG is the metacategory $\bS_G$ whose objects are those of \bC, the morphisms from \mX to \mY being the conglomerate 
$\Nat(\Mor(Y,G-),\Mor(X,G-))$.  While ad hoc arguments can sometimes be used to show that $\bS_G$ is isomorphic to a 
category, the situation is optimal when \mG has a left adjoint $F:\bC \ra \bD$ since in this case $\bS_G$ is isomorphic to 
$\bKL(\bT)$, \bT the triple in \bC determined by \mF and \mG.]\\

\begingroup%%------------------------------------>>%%----------------------------------->>
\fontsize{9pt}{11pt}\selectfont
Consider the adjoint situation produced by the forgetful functor $\bGR \ra \bSET$ $-$then 
$\bKL(\bT)$ is isomorphic to the full subcategory of \bGR whose objects are the free groups.\\
\endgroup%%------------------------------------<< %%------------------------------------<<

A triple $\bT = (T, m, \epsilon)$ in \bC is said to be 
\un{idempotent}
\index{idempotent (triple)} 
provided  that $m$ is a natural isomorphism (hence $\epsilon T = m^{-1} = T \epsilon$).  If \bT is idempotent, then the 
comparison functor $\bKL(\bT) \ra \bT$-$\bALG$ is an equivalence of categories.  Moreover, 
$G_{\bT}:\bT$-$\bALG \ra \bC$ is full, faithful, and injective on objects.  
Its image is a reflective subcategory of \bC, the objects 
%%----------------------------------------------------------------------------------------------30
being those \mX such that $\epsilon_X:X \ra TX$ is an isomorphism.  
On the other hand, every reflective subcategory of \bC generates an idempotent triple.  
Agreeing that two idempotent triples 
$T$ and $T^\prime$ are equivalent if there exits a natural isomorphism $\tau:T \ra T^\prime$ such that 
$\epsilon^\prime = \tau \circ \epsilon$ (thus also 
$\tau \circ m = m^\prime \circ \tau T^\prime \circ T \tau$), 
the conclusion is that the conglomerate of reflective subcategories of \bC is in a one-to-one correspondence with the conglomerate of 
idempotent triples in \bC modulo equivalence.

\label{8.3}
\label{8.19}
\label{9.1}
\label{9.13}
\label{9.16}
\label{9.55a}
\label{9.61}
[Note: \ An idempotent triple $\bT = (T, m, \epsilon)$ determines an orthogonal pair $(S,D)$.  Let $f:X \ra Y$ be a 
morphism $-$then $f$ is said to be 
\un{\mT-localizing}
\index{localizing! \bT-localizing (idempotent triple)}
if there is an isomorphism $\phi:TX \ra Y$ such that $f = \phi \circ \epsilon_X$.  
For this to be the case, it is necessary and 
sufficient that $f \in S$ and $Y \in D$.  
If $\bC^\prime$ is a full subcategory of \bC and if 
$\bT^\prime = (T^\prime, m^\prime, \epsilon^\prime)$ is an idempotent triple in 
$\bC^\prime$, then \bT (or $T$) is said to 
\un{extend}
\index{extend (an idempotent triple)} 
$\bT^\prime$ (or $T^\prime$) provided that $S^\prime \subset S$ and 
$D^\prime \subset D$ (in general, $(S^\prime )^\perp \supset D \supset (D^\prime)^{\perp\perp}$, 
where orthogonality is 
meant in \bC).]\\

\begingroup%%------------------------------------>>%%----------------------------------->>
\fontsize{9pt}{11pt}\selectfont
Let $(F,G,\mu,\nu)$ be an adjoint situation $-$then the following conditions are equivalent: 
(1) $(G \circ F,G\nu F,\mu)$ is an idempotent triple;  
(2) $\mu G$ is a natural isomorphism; 
(3) $(F \circ G, F\mu G, \nu)$ is an idempotent cotriple;
(4) $\nu F$ is a natural isomorphism.  
And: $(1), \ldots, (4)$ imply that the full subcategory $\bC_\mu$ of \bC whose objects are the \mX such that $\mu_X$ is an 
isomorphism is a reflective subcategory of \bC and the full subcategory $\bD_\nu$ of \bD whose objects are the \mY such that $\nu_Y$ is an isomorphism is a coreflective subcategory of \bD.
\\ \indent
[Note: \ $\bC_\mu$ and $\bD_\nu$  are equivalent categories.]\\
\endgroup%%------------------------------------<< %%------------------------------------<<

Given a category \bC and a class $S \subset \Mor \bC$, a 
\un{localization of \bC at \mS}
\index{localization of \bC at \mS} 
is a pair $(S^{-1}\bC,L_S)$, where $S^{-1}\bC$ is a metacategory and $L_S:\bC \ra S^{-1}\bC$ is a functor such that 
$\forall \ s \in S$, $L_Ss$ is an isomorphism, $(S^{-1}\bC,L_S)$ being initial among all pairs having this property, 
i.e., for any 
metacategory \bD and for any functor $F:\bC \ra \bD$ such that $\forall \ s \in S$, $Fs$ is an isomorphism, 
there exists a unique functor $F^\prime:S^{-1}\bC \ra \bD$ such that $F = F^\prime \circ L_S$.  
$S^{-1}\bC$ exists, is unique up to isomorphism, 
and there is a representative that has the same objects as \bC itself.  
Example: Take $\bC = \bTOP$ and let $S \subset \Mor \bC$ be the class of homotopy equivalences $-$then 
$S^{-1}\bC = \bHTOP$.

[Note: \ If $\ov{S}$ is the class of all morphisms rendered invertible by $L_S$ (the 
\un{saturation}
\index{saturation (of a class of all morphisms)} 
of \mS), then the arrow $S^{-1}\bC \ra \ov{S}^{-1}\bC$ is an isomorphism.]\\

\begingroup%%------------------------------------>>%%----------------------------------->>
\fontsize{9pt}{11pt}\selectfont
Fix a class $I$ which is not a set.  Let \bC be the category whose objects are \mX, \mY, and $\{Z_i: i \in I\}$ and whose 
morphisms, apart from the identities, are $f_i:X \ra Z_i$ and $g_i:Y \ra Z_i$.  Take $S = \{g_i:i \in I\}$ $-$then 
$S^{-1}\bC$ is a metacategory that is not isomorphic to a category.

\label{15.14}
[Note: \ The localization of a small category at a set of morphisms is again small.]\\
\endgroup%%------------------------------------<< %%------------------------------------<<

%%----------------------------------------------------------------------------------------------31
\label{9.60}
Let \bC be a category and let $S \subset \Mor \bC$ be a class containing the identities of \bC and closed with respect to 
composition $-$then \mS is said to admit a 
\un{calculus of left fractions}
\index{calculus of left fractions} if

\indent\indent (LF$_1$) \ Given a 2-source $X^\prime \overset{s}{\la} X \overset{f}{\ra} Y$ $(s \in S)$, there exists a commutative square
\begin{tikzcd}[sep=large]
{X} \ar{d}[swap]{s} \ar{r}{f} &{Y} \ar{d}{t}\\
{X^\prime} \ar{r}[swap]{f^\prime} &{Y^\prime}
\end{tikzcd}
, where $t \in S$;

\indent\indent (LF$_2$) \ Given $f, g:X \ra Y$ and $s:X^\prime \ra X$ $(s \in S)$ such that $f \circ s = g \circ s$, there exists 
$t:Y \ra Y^\prime$ $(t \in S)$ such that $t \circ f = t \circ g$.

[Note: \ Reverse the arrows to define ``calculus of right fractions'' 
\index{calculus of right fractions}.]\\

\begingroup%%------------------------------------>>%%----------------------------------->>
\fontsize{9pt}{11pt}\selectfont
Let $S \subset \Mor \bC$ be a class containing the identities of \bC and closed with respect to composition such that 
$\forall \ (s,t)$: $t \circ s \in S$ $\&$ $s \in S$ $\implies$ $t \in S$ $-$then \mS admits a calculus of left fractions if every 
2-source 
$X^\prime \overset{s}{\la} X \overset{f}{\ra} Y$ $(s \in S)$ can be completed to a weak pushout square
\begin{tikzcd}[sep=large]
{X} \ar{d}[swap]{s} \ar{r}{f} &{Y} \ar{d}{t}\\
{X^\prime} \ar{r}[swap]{f^\prime} &{Y^\prime}
\end{tikzcd}
, where $t \in S$.  For an illustration, take $\bC = \bHTOP$ and consider the class of homotopy classes of homology equivalences.\\

\endgroup%%------------------------------------<< %%------------------------------------<<

\label{15.41}
Let \bC be a category and let $S \subset \Mor \bC$ be a class admitting a calculus of left fractions.  
Given $X, Y \in \Ob S^{-1}\bC$, $\Mor(X,Y)$ is the conglomerate of equivalence classes of pairs $(s,f)$: 
 $X \overset{f}{\ra} Y^\prime \overset{s}{\la} Y$, two pairs
 $
\begin{cases}
\ (s,f)\\
\ (t,g)
\end{cases}
$
 being equivalent iff there exists $u, v \in \Mor \bC$: 
 $
\begin{cases}
\ u \circ s\\
\ v \circ t
\end{cases}
\in S,
$
with $u \circ s = v \circ t$ and $u \circ f = v \circ g$ .  
Every morphism in $S^{-1}\bC$ can be represented in the form 
 $(L_S s)^{-1}L_Sf$ and if $L_S f = L_S g$, then there is an $s \in S$ such that $s \circ f = s \circ g$.
 
[Note: \ $S^{-1}\bC$ is a metacategory.  
To guarantee that $S^{-1}\bC$ is isomorphic to a category, it suffices to 
impose a \un{solution set condition}: 
\index{solution set condition ($S^{-1}\bC$ is isomorphic to a category)}
For each $X \in \Ob\bC$, there exists a source $\{s_i:X \ra X_i^\prime\}$ $(s_i \in S)$ such that for every $s:X \ra X^\prime$ 
$(s \in S)$, there is an $i$ and a $u:X^\prime \ra X_i^\prime$ such that $u \circ s = s_i$.  
This, of course, is automatic provided that 
$X\backslash S$, the full subcategory of $X\backslash \bC$ whose objects are the $s:X \ra X^\prime$ $(s \in S)$, has a 
final object.]\\

\begingroup%%------------------------------------>>%%----------------------------------->>
\fontsize{9pt}{11pt}\selectfont
If \bC is the full subcategory of $\bHTOP_*$ whose objects are the pointed connected CW complexes and if \mS is the class of pointed homotopy classes of pointed $n$-equivalences, then \mS admits a calculus of left fractions and satisfies the solution set condition.\\
\endgroup%%------------------------------------<< %%------------------------------------<<

Let $(F,G,\mu,\nu)$ be an adjoint situation.  
Assume: $G$ is full and faithful or, equivalently, that $\nu$ is a natural isomorphism.  
Take for $S \subset \Mor \bC$ the class consisting of those $s$ such that $Fs$ is an isomorphism 
(so, $F = F^\prime \circ L_S$) $-$then $\{\mu_X\} \subset S$ and $S$ admits a calculus 
%%----------------------------------------------------------------------------------------------32
of left fractions.  Moreover, \mS is saturated and satisifies the solution set condition (in fact, $\forall \ X \in \Ob\bC$, 
$X\backslash S$ has a final object, viz. $\mu_X:X \ra GFX$).  
Therefore $S^{-1}\bC$ is isomorphic to a category and 
$L_S:\bC \ra S^{-1}\bC$ has a right adjoint that is full and faithful, while $F^\prime:S^{-1}\bC \ra \bD$ is an equivalence.

[Note: \ Suppose that $\bT = (T,m,\epsilon)$ is an idempotent triple in \bC.  Let \bD be the corresponding reflective subcategory of \bC with reflector $R:\bC \ra \bD$, so $T = \iota \circ R$, where $\iota:\bD \ra \bC$ is the inclusion.  
Take for $S \subset \Mor \bC$ the class consisting of those $f$ such that $Tf$ is an isomorphism $-$then \mS is the class consisting of those $f$ such that $Rf$ is an isomorphism, hence $S$ admits a calculus of left fractions, is saturated, and satisfies the solution set condition.  
The Kleisli category of \bT is isomorphic to $S^{-1}\bC$ and $T$ factors as 
$\bC \ra S^{-1}\bC \ra \bD \ra \bC$, the arrow $S^{-1}\bC \ra \bD$ being an equivalence.]\\

\vspace{0.25cm}

\label{13.22}
\begingroup%%------------------------------------>>%%----------------------------------->>
\fontsize{9pt}{11pt}\selectfont
Let $(F,G,\mu,\nu)$ be an adjoint situation.  Put 
$
\begin{cases}
\ S = \{\mu_X\} \subset \Mor \bC\\
\ T = \{\nu_Y\} \subset \Mor \bD
\end{cases}
$
$-$then 
$
\begin{cases}
\ S^{-1}\bC\\
\ T^{-1}\bD
\end{cases}
$
are isomorphic to categories and
$
\begin{cases}
\ F\\
\ G
\end{cases}
$
induce functors 
$
\begin{cases}
\ F^\prime:S^{-1}\bC \ra T^{-1}\bD\\
\ G^\prime:T^{-1}\bD \ra S^{-1}\bC
\end{cases}
$
such that 
$
\begin{cases}
\ G^\prime \circ F^\prime \approx \id_{S^{-1}\bC}\\
\ F^\prime \circ G^\prime \approx \id_{T^{-1}\bD}
\end{cases}
, \ 
$
thus 
$
\begin{cases}
\ S^{-1}\bC\\
\ T^{-1}\bD
\end{cases}
$
are equivalent.  
In particular, when $G$ is full and faithful, $S^{-1}\bC$ is equivalent to \bD (the saturation of \mS 
being the class consisting of those $s$ such that $Fs$ is an isomorphism, 
i.e., $\ov{S}$ is the ``$S$'' considered above.).\\
\endgroup%%------------------------------------<< %%------------------------------------<<

Given a category \bC, a set $\sU$ of objects in \bC is said to be a 
\un{separating set}
\index{separating set (in a category)} 
if for every pair 
$X \overset{f}{\underset{g}{\rightrightarrows}} Y$ of distinct morphisms, there exists a $U \in \sU$ and a morphism 
$\sigma:U \ra X$ such that $f \circ \sigma \neq g \circ \sigma$.  
An object $U$ in \bC is said to be a 
\un{separator}
\index{separator (in a category)}
if $\{U\}$ is a separating set, i.e., if the functor $\Mor(U,-):\bC \ra \bSET$ is faithful.  
If \bC is balanced, finitely complete, and has a separating set, then \bC is wellpowered.  Every cocomplete, cowellpowered category with a separator is wellpowered and complete.  If \bC has coproducts, then a $U \in \Ob\bC$ is a separator iff 
each $X \in \Ob\bC$ admits an epimorphism $\coprod U \ra X$.

[Note: \ Suppose that \bC is small $-$then the representable functors are a separating set for $[\bC,\bSET]$.]\\

\begingroup%%------------------------------------>>%%----------------------------------->>
\fontsize{9pt}{11pt}\selectfont
Every nonempty set is a separator for \bSET.  $\bSET \times \bSET$ has no separators but the set 
$\{(\emptyset,\{0\}),(\{0\},\emptyset)\}$ is a separating set.  Every nonempty discrete topological space is a separator for 
\bTOP (or \bHAUS).  $\Z$ is s separator for \bGR and \bAB, while $\Z[t]$ is a separator for \bRG.  
In \bAMOD, \mA (as a left \mA-module) is a separator and in \bMODA, \mA (as a right \mA-module) is a separator.\\
\endgroup%%------------------------------------<< %%------------------------------------<<

Given a category \bC, a set $\sU$ of objects in \bC is said to be a 
\un{coseparating set}
\index{coseparating set (in a category)} 
if for 
%%----------------------------------------------------------------------------------------------33
 every pair 
$X \overset{f}{\underset{g}{\rightrightarrows}} Y$ of distinct morphisms, there exists a $U \in \sU$ and a morphism 
$\sigma:Y \ra U $ such that $\sigma \circ f \neq \sigma \circ g$.  
An object $U$ in \bC is said to be a 
\un{coseparator}
\index{coseparator (in a category)}
if $\{U\}$ is a coseparating set, i.e., if the cofunctor $\Mor(-,U):\bC \ra \bSET$ is faithful.  
If \bC is balanced, finitely cocomplete, and has a coseparating set, then \bC is cowellpowered.  
Every complete, wellpowered category with a coseparator is cowellpowered and cocomplete.  
If \bC has products, then a $U \in \Ob\bC$ is a coseparator iff 
each $X \in \Ob\bC$ admits a monomorphism $X \ra \prod U$.\\

\begingroup%%------------------------------------>>%%----------------------------------->>
\fontsize{9pt}{11pt}\selectfont
Every set with at least two elements is a coseparator for \bSET.  
Every indiscrete topological space with at least two elements is a coseparator for \bTOP.  $\Q/\Z$ is a coseparator for \bAB.  
None of the categories \bGR, \bRG, \bHAUS has a coseparating set.\\
\endgroup%%------------------------------------<< %%------------------------------------<<

\index{Theorem: Special Adjoint Functor Theorem}
\index{Special Adjoint Functor Theorem}
\textbf{\small SPECIAL ADJOINT FUNCTOR THEOREM} \quadx
Given a complete wellpowered category \bD which has a coseparating set, a functor $G:\bD \ra \bC$ has a left adjoint iff \mG preserves limits.\\

\begingroup%%------------------------------------>>%%----------------------------------->>
\fontsize{9pt}{11pt}\selectfont
\label{2.1}
A functor from \bSET, \bAB, or \bTOP to a category \bC has a left adjoint 
iff it preserves limits and a right adjoint iff it preserve colimits.\\
\endgroup%%------------------------------------<< %%------------------------------------<<

Given a category \bC, an object \mP in \bC is said to be 
\un{projective}
\index{projective (object in \bC)} 
if the functor $\Mor(P,-):\bC \ra \bSET$ preserves epimorphisms.  In other words: \mP is projective iff for each epimorphism 
$f:X \ra Y$ and each morphism $\phi:P \ra Y$, there exists a morphism $g:P \ra X$ such that $f \circ g = \phi$.  
A coproduct of projective objects is projective.

A category \bC is said to have 
\un{enough projectives}
\index{enough projectives} 
provided that for any $X \in \Ob\bC$ there is an epimorphism $P \ra X$, with \mP projective.  
If a category has enough projectives and a separator, then it has a projective separator.  If a category has coproducts and a projective separator, then it has enough projectives.\\

\begingroup%%------------------------------------>>%%----------------------------------->>
\fontsize{9pt}{11pt}\selectfont
The projective objects in the category of compact Hausdorff spaces are the extremally disconnected spaces.  
The projective objects in \bAB or \bGR are the free groups.  The full subcategory of \bAB whose objects are the torsion groups has no projective objects other than the initial objects.  In \bAMOD or \bMODA, an object is projective iff it is a direct summand of a free module (and every free module is a projective separator).\\
\endgroup%%------------------------------------<< %%------------------------------------<<

Given a category \bC, an object \mQ in \bC is said to be 
\un{injective}
\index{injective (object in \bC)} 
if the cofunctor $\Mor(-,Q):\bC \ra \bSET$ converts monomorphisms into epimorphisms.  In other words: \mQ is injective 
%%----------------------------------------------------------------------------------------------34
iff for each monomorphism $f:X \ra Y$ and each morphism $\phi:X \ra Q$, there exists a morphism 
$g:Y \ra Q$ such that $g \circ f = \phi$.  A product of injective objects is injective.

A category \bC is said to have 
\un{enough injectives}
\index{enough injectives} 
provided that for any $X \in \Ob\bC$, there is a monomorphism $X \ra Q$, with \mQ injective.  
If a category has enough injectives and a coseparator, then it has an injective coseparator.  
If a category has products and a injective coseparator, then it has enough injectives.\\

\begingroup%%------------------------------------>>%%----------------------------------->>
\fontsize{9pt}{11pt}\selectfont
The injective objects in the category of compact Hausdorff spaces are the retracts of products $\prod [0,1]$.  
The injective objects in the category of Banach spaces and linear contractions are, up to isomorphism, the $C(X)$, where 
\mX is an extremally disconnected compact Hausdorff space.  
In \bAB, the injective objects are the divisible abelian groups (and $\Q/\Z$ is an injective coseparator) but the only injective objects in 
\bGR or \bRG are the final objects.  The module $\Hom_{\Z}(A,\Q/\Z)$ is an injective coseparator in \bAMOD or \bMODA.\\
\endgroup%%------------------------------------<< %%------------------------------------<<

A 
\un{zero object}
\index{zero object} 
in a category \bC is an object which is both initial and final.  
The categories $\bTOP_*$, \bGR, and \bAB have zero objects.  
If \bC has a zero object $0_{\bC}$ (or 0), then for any ordered pair $X, \ Y \in \Ob\bC$ there exists a unique morphism $X \ra 0_{\bC} \ra Y$, the 
\un{zero morphism}
\index{zero morphism} 
$0_{XY}$ (or 0) in $\Mor(X,Y)$.  
It does not depend on the choice of a zero object in \bC.  
An equalizer (coequalizer) of an $f \in \Mor(X,Y)$ and $0_{XY}$ 
is said be be a 
\un{kernel}
\index{kernel (in a category)} 
(\un{cokernel}
\index{cokernel (in a category)}) 
of $f$.  
Notation: $\ker f$ ($\coker f$).

[Note: \ Suppose that \bC has a zero object.  Let $\{X_i:i \in I\}$ be a collection of objects in \bC for which 
$\prod\limits_i X_i$ and $\coprod\limits_i X_i$ exist.  
The morphisms $\delta_{ij}:X_i \ra X_j$ defined by
$
\begin{cases}
\ \id_{X_i} \ \ \ (i = j)\\
\ 0_{X_iX_j} \ (i \neq j)
\end{cases}
$
then determine a morphism 
\label{15.18}
\label{15.23}
$t:\coprod\limits_i X_i \ra$ $\prod\limits_i X_i$ such that 
$\pr_j \circ t \circ \ini_i$ $=$ $\delta_{ij}$.  
Example: Take $\#(I) = 2$ $-$then this morphism can be a monomorphism (in $\bTOP_*$), an epimorphism (in \bGR), or an isomorphism (in \bAB).]\\

\begingroup%%------------------------------------>>%%----------------------------------->>
\fontsize{9pt}{11pt}\selectfont
A 
\un{pointed category}
\index{pointed category} 
is a category with a zero object.\\
\endgroup%%------------------------------------<< %%------------------------------------<<

Let \bC be a category with a zero object.  
Assume that \bC has kernels and cokernels.  
Given a morphism 
$f:X \ra Y$, an 
\un{image}
\index{image (of a morphism)} \ 
(\un{coimage})
\index{coimage (of a morphism)} 
of $f$ is a kernel  of a cokernel (cokernel of a kernel) for $f$.  
Notation: $\im f$ ($\coim f$).  
There is a commutative diagram
\[
\begin{tikzcd}[sep=large]
{\ker f} \ar{r} &{X} \ar{d} \ar{r}{f} &{Y} \ar{r} &{\coker f}\\
&{\coim f} \ar{r}[swap]{\ov{f}} &{\im f,} \ar{u}
\end{tikzcd}
\]
where $\ov{f}$ is the morphism 
\un{parallel}
\index{parallel (morphism)} 
to $f$.
If parallel morphisms are isomorphisms, then \bC is said to be an 
\un{exact category}
\index{exact category}.

%%----------------------------------------------------------------------------------------------35
[Note: \ In general, $\ov{f}$ need be neither a monomorphism nor an epimorphism and $\ov{f}$ can be a bimorphism without being an isomorphism.]

A category \bC that has a zero object is exact iff every monomorphism is the kernel of a morphism, every epimorphism is the cokernel of a morphism, and every morphism admits a factorization: $f = g \circ h$ ($g$ a monomorphism, $h$ an epimorphism).  
Such a factorization is essentially unique.  
An exact category is balanced; it is wellpowered iff it is cowellpowered.  
Every exact category with a separator or a coseparator is wellpowered and cowellpowered.  
If an exact category has finite products (finite coproducts), then it has equalizers (coequalizers), hence is finitely complete (finitely cocomplete).\\

\begingroup%%------------------------------------>>%%----------------------------------->>
\fontsize{9pt}{11pt}\selectfont
\bAB is an exact category but the full subcategory of \bAB whose objects are the torsion free abelian groups is not exact.  
Neither \bGR nor $\bTOP_*$ is exact.\\
\endgroup%%------------------------------------<< %%------------------------------------<<

Let \bC be an exact category.

\indent\indent (EX) \ A sequence 
$\cdots \ra X_{n-1} \overset{d_{n-1}}{\lra} X_n \overset{d_{n}}{\lra} X_{n+1} \ra \cdots$ 
is said to be 
\un{exact}
\index{exact sequence} 
provided that $\im d_{n-1} \approx \ker d_n$ for all $n$.

[Note: \ A 
\un{short exact sequence}
\index{short exact sequence} 
is an exact sequence of the form 
$0 \ra$ 
$X^\prime \ra$ 
$X \ra$ 
$X\pp \ra$ $0$.]\\

\indent\indent (Ker-Coker Lemma) 
\index{Ker-Coker Lemma} \ 
Suppose that the diagram
\[
\begin{tikzcd}[sep=large]
&{X_1} \ar{d}{f_1} \ar{r} &{X_2} \ar{d}{f_2} \ar{r} &{X_3} \ar{d}{f_3} \ar{r} &{0}\\
{0} \ar{r} &{Y_1} \ar{r} &{Y_2} \ar{r} &{Y_3}
\end{tikzcd}
\]
is commutative and has exact rows $-$then there is a morphism $\delta: \ker f_3 \ra \coker f_1$, the 
\un{connecting morphism}
\index{connecting morphism}, 
such that the sequence
\[
\ker f_1 \ra \ker f_2 \ra \ker f_3 \overset{\delta}{\ra} \coker f_1 \ra \coker f_2 \ra \coker f_3
\]
is exact.  Moreover, if $X_1 \ra X_2$ ($Y_2 \ra Y_3$) is a monomorphism (epimorphism), then 
$\ker f_1 \ra \ker f_2$ $(\coker f_2 \ra \coker f_3)$ is a monomorphism (epimorphism).\\

\indent\indent (Five Lemma) 
\index{Five Lemma} \ 
Suppose that the diagram
\[
\begin{tikzcd}[sep=large]
{X_1} \ar{d}{f_1} \ar{r} &{X_2} \ar{d}{f_2} \ar{r} &{X_3} \ar{d}{f_3} \ar{r} &{X_4} \ar{d}{f_4} \ar{r} &{X_5} \ar{d}{f_5}\\
{Y_1} \ar{r} &{Y_2} \ar{r} &{Y_3} \ar{r} &{Y_4} \ar{r} &{Y_5}
\end{tikzcd}
\]
is commutative and has exact rows.

%%----------------------------------------------------------------------------------------------36
\indent\indent (1) \ If $f_2$ and $f_4$ are epimorphisms and $f_5$ is a monomorphism, then $f_3$ is an epimorphism.

\indent\indent (2)  If $f_2$ and $f_4$ are monomorphisms and $f_1$ is an epimorphism, then $f_3$ is a monomorphism.

\indent\indent(Nine Lemma) 
\index{Nine Lemma} 
\ Suppose that the diagram
\[
\begin{tikzcd}[sep=large]
&{0} \ar{d} &{0} \ar{d} &{0} \ar{d}\\
{0} \ar{r} &{X^\prime} \ar{d} \ar{r} &{X} \ar{r} \ar{d} &{X\pp} \ar{d} \ar{r} &{0}\\
{0} \ar{r} &{Y^\prime} \ar{d} \ar{r} &{Y} \ar{r} \ar{d} &{Y\pp} \ar{d} \ar{r} &{0}\\
{0} \ar{r} &{Z^\prime} \ar{d} \ar{r} &{Z} \ar{r} \ar{d} &{Z\pp} \ar{d} \ar{r} &{0}\\
&{0} &{0} &{0}
\end{tikzcd}
\]
is commutative, has exact columns, and an exact middle row $-$then the bottom row is exact iff the top row is exact.\\

\begingroup%%------------------------------------>>%%----------------------------------->>
\fontsize{9pt}{11pt}\selectfont
In an exact category \bC, there are two short exact sequences associated with each morphism $f:X \ra Y$, viz.
$
\begin{cases}
\ 0 \ra \ker f \ra X \ra \coim f \ra 0\\
\ 0 \ra \im f \ra Y \ra \coker f \ra 0
\end{cases}
$
.\\
\endgroup%%------------------------------------<< %%------------------------------------<<

\label{16.21}
An 
\un{additive category}
\index{additive category} 
is a category \bC that has a zero object and which is equipped with a 
function + that assigns to each ordered pair $f, g \in \Mor \bC$ having common domain and codomain, a morphism 
$f + g$ with the same domain and codomain satisfying the following conditions.

\indent\indent (ADD$_1$) \ On each morphism set $\Mor(X,Y)$, + induces the structure of an abelian group.

\indent\indent (ADD$_2$) \ Composition is distributive over +:
$
\begin{cases}
\ f \circ (g + h) = (f \circ g) + (f \circ h)\\
\ (g + h) \circ k = (g \circ k) + (h \circ k)
\end{cases}
.
$

\indent\indent (ADD$_3$) \ The zero morphisms are identities with respect to $+$: $0+f = f+0 = f$.\\
An additive category has finite products iff it has finite coproducts and when this is so, finite coproducts are finite products.

[Note: If \bC is small and \bD is additive, then $[\bC,\bD]$ is additive.]\\

%%----------------------------------------------------------------------------------------------37

\begingroup%%------------------------------------>>%%----------------------------------->>
\fontsize{9pt}{11pt}\selectfont
\bAB is an additive category but \bGR is not.  
Any ring with unit can be viewed as an additive category having exactly one object (and conversely).  
The category of Banach spaces and continuous linear transformations is additive but not exact.\\
\endgroup%%------------------------------------<< %%------------------------------------<<

An 
\un{abelian category}
\index{abelian category} 
is an exact category \bC that has finite products and finite coproducts.  
Every abelian category is additive, finitely complete, and finitely cocomplete.  A category \bC that has a zero object is abelian iff it has pullbacks, pushouts, and every monomorphism (epimorphism) is the kernel (cokernel) of a morphism.  
In an abelian category, 
$t:\coprod\limits_{i=1}^n X_i \ra \prod\limits_{i=1}^n X_i$ is an isomorphism.

[Note: \ If \bC is small and \bD is abelian, then $[\bC,\bD]$ is abelian.]\\

\begingroup%%------------------------------------>>%%----------------------------------->>
\fontsize{9pt}{11pt}\selectfont
\bAB is an abelian category, as is its full subcategory whose objects are the finite abelian groups but there are full subcategories of \bAB which are exact and additive, yet not abelian.\\
\endgroup%%------------------------------------<< %%------------------------------------<<

A 
\un{Grothendieck category}
\index{Grothendieck category} 
is a cocomplete abelian category \bC in which filtered colimits commute with finite limits or, equivalently, in which filtered colimits of exact sequences are exact.  
Every Grothendieck category with a separator is complete and has an injective coseparator, hence has enough injectives 
(however there exist wellpowered Grothendieck categories that do not have enough injectives).  
In a Grothendieck category, every filtered colimit of monomorphisms is a monomorphism, coproducts of monomorphisms are monomorphisms, and 
$t:\coprod\limits_i X_i \ra \prod\limits_i X_i$
is a monomorphism.

[Note: \ If \bC is small and \bD is Grothendieck, then $[\bC,\bD]$ is Grothendieck.]\\

\begingroup%%------------------------------------>>%%----------------------------------->>
\fontsize{9pt}{11pt}\selectfont
\bAB is a Grothendieck category but its full subcategory whose objects are the finitely generated abelian groups, while abelian, is not Grothendieck.  
If \mA is a ring with unit, then \bAMOD and \bMODA are  Grothendieck categories.\\
\endgroup%%------------------------------------<< %%------------------------------------<<

Given exact categories 
$
\begin{cases}
\ \bC\\
\ \bD
\end{cases}
\hspace{-.25cm}, \ 
$
a functor $F:\bC \ra \bD$ is said to be 
\un{left exact}
\index{left exact (functor)} 
(\un{right} 
\un{exact})
\index{right exact (functor)} 
if it preserves kernels (cokernels) and 
\un{exact}
\index{exact (functor)} 
if it is both right and left exact.  
$F$ is left exact (right exact) iff for every short exact sequence 
$0 \ra X^\prime \ra$ 
$X \ra$ 
$X\pp \ra$ $0$ in \bC, the sequence 
$0 \ra FX^\prime \ra$ 
$FX \ra$ 
$FX\pp$ 
$(FX^\prime \ra$ 
$FX \ra$ 
$FX\pp \ra 0$) is exact in \bD.  
Therefore $F$ is exact iff $F$ preserves short exact sequences or still, iff $F$ preserves arbitrary exact sequences.

[Note: \ $F$ is said to be 
\un{half exact}
\index{half exact (functor)} 
if for every short exact sequence 
$0 \ra X^\prime \ra X \ra X\pp \ra 0$ in \bC, the sequence 
$FX^\prime \ra FX \ra FX\pp$  is exact in \bD.]\\

\begingroup%%------------------------------------>>%%----------------------------------->>
%\linespread{1.5} %dmc
\fontsize{9pt}{11pt}\selectfont
The projective (injective) objects in an abelian category are those for which 
$\Mor(X,-)$ $(\Mor(-,X))$ is exact.  
In \bAB, $X \otimes -$ is exact iff \mX is flat or here, torsion free.  
If \bI is small and filtered and if \bC is Grothendieck, then $\colim:[\bI,\bC] \ra \bC$ is exact.\\
\endgroup%%------------------------------------<< %%------------------------------------<<

%%----------------------------------------------------------------------------------------------38
Given additive categories 
$
\begin{cases}
\ \bC\\
\ \bD
\end{cases}
\hspace{-.25cm}, \ 
$
a functor $F:\bC \ra \bD$ is said to be 
\un{additive}
\index{additive (functor)}
if for all $X, \ Y \in \Ob\bC$, the map $\Mor(X,Y) \ra \Mor(FX,FY)$ is a homomorphism of abelian groups.  
Every half exact functor between abelian categories is additive. An additive functor between abelian categories is 
left exact (right exact) iff it preserves finite limits (finite colimits).  
The 
\un{additive functor category}
\index{additive functor category} 
$[\bC,\bD]^+$ 
\index{$[\bC,\bD]^+$} 
is the full submetacategory of $[\bC,\bD]$ whose objects are the additive functors.  
There are Yoneda embeddings
$
\begin{cases}
\ \bC^\OP \ra [\bC,\bAB]^+\\
\ \bC \ra [\bC^\OP,\bAB]^+
\end{cases}
. \ 
$
If \bC and \bD are abelian categories with \bC small, if $K:\bC \ra \bD$ is additive, and if \bS is a complete (cocomplete) 
abelian category, then there is an additive version of Kan extension applicable to 
$
\begin{cases}
\ [\bC,\bS]^+\\
\ [\bD,\bS]^+
\end{cases}
. \ 
$
The functors produced need not agree with those obtained by forgetting the additive structure.\\

\begingroup%%------------------------------------>>%%----------------------------------->>
\fontsize{9pt}{11pt}\selectfont
Let \mA be a ring with unit viewed as an additive category having exactly one object $-$then 
\bAMOD is isomorphic to 
$[A,\bAB]^+$ and \bMODA  is isomorphic to $[A^\OP,\bAB]^+$.

[Note: \ A right $A$-module \mX and a left $A$-module \mY define a diagram 
$A^\OP \times A \ra \bAB$ (tensor product over $\Z$) and the coend $\ds\int^A X \otimes Y$ is $X \otimes_A Y$, the tensor product over $A$.]\\
\endgroup%%------------------------------------<< %%------------------------------------<<

\label{15.22}
If \bC is small and additive and if \bD is additive, then

\indent\indent (1) \ \bD finitely complete and wellpowered (finitely cocomplete and cowellpowered) 
$\implies$ $[\bC,\bD]^+$ wellpowered (cowellpowered); 

\indent\indent (2) \ \bD (finitely) complete $\implies$ $[\bC,\bD]^+$ (finitely) complete and 
\bD (finitely) cocomplete $\implies$ $[\bC,\bD]^+$ (finitely) cocomplete;

\indent\indent (3) \ \bD abelian (Grothendieck) $\implies$ $[\bC,\bD]^+$ $\implies$ abelian (Grothendieck).

[Note: \ Suppose that \bC is small.  If \bC is additive, then $[\bC,\bAB]^+$  is a complete Grothendieck category and if 
\bC is exact and additive, then $[\bC,\bAB]^+$ has a separator which as a functor $\bC \ra \bAB$ is left exact.]\\

Given a small abelian category \bC and an abelian category \bD, write 
$\bLEX(\bC,\bD)$
\index{$\bLEX(\bC,\bD)$} 
for the full, isomorphism closed subcategory of $[\bC,\bD]^+$ whose objects are the left exact functors.\\

\index{Theorem: Dervied Functor Theorem}
\index{Dervied Functor Theorem}
\textbf{\small DERIVED FUNCTOR THEOREM} \ 
If \bC is a small abelian category and if \bD is a wellpowered Grothendieck category, then $\bLEX(\bC,\bD)$ is a reflective 
subcategory of $[\bC,\bD]^+$.  As such, it is Grothendieck.  
Moreover, the reflector is an exact functor.

[Note: \ The reflector sends $F$ to its 
\un{zeroth right derived functor}
\index{zeroth right derived functor} 
$R^0F$.]\\

If \bC is a small abelian category, then $\bLEX(\bC,\bAB)$ is a Grothendieck category with a separator.  
Therefore $\bLEX(\bC,\bAB)$ has enough injectives.  Every injective object in 
%%----------------------------------------------------------------------------------------------39
$\bLEX(\bC,\bAB)$ is an exact functor.  The Yoneda embedding $\bC^\OP \ra [\bC,\bAB]^+$ is left exact.  It factors through 
$\bLEX(\bC,\bAB)$ and is then exact.

[Note: \ Since \bC is abelian, every object in $[\bC,\bAB]^+$ is a colimit of representable functors and every object in 
$\bLEX(\bC,\bAB)$ is a filtered colimit of representable functors.  
Thus $\bLEX(\bC,\bAB)$ is equivalent to $\bIND(\bC^\OP)$ and so 
$\bLEX(\bC,\bAB)^\OP$ is equivalent to $\bPRO(\bC)$.]\\

\begingroup%%------------------------------------>>%%----------------------------------->>
\fontsize{9pt}{11pt}\selectfont
The full subcategory of \bAB whose objects are the finite abelian groups is equivalent to a small category.  Its procategory is equivalent to the opposite of the full subcategory of \bAB whose objects are the torsion abelian groups.\\
\endgroup%%------------------------------------<< %%------------------------------------<<

Given an abelian category \bC, a nonempty class $\sC \subset \Ob\bC$ is said to be a 
\un{Serre class}
\index{Serre class (abelian category)} 
provided that for any short exact sequence 
$0 \ra X^\prime \ra X \ra X\pp \ra 0$ in \bC, $X \in \sC$ iff 
$
\begin{cases}
\ X^\prime\\
\ X\pp 
\end{cases}
\in \sC
$
or equivalently, for any exact sequence $X^\prime \ra X \ra X\pp$ in \bC, 
$
\begin{cases}
\ X^\prime\\
\ X\pp 
\end{cases}
\in \sC
$
$\implies X \in \sC$.

[Note: \ Since $\sC$ is nonempty, $\sC$ contains the zero objects of \bC.]

\label{15.5}
Given an abelian category \bC with a separator and a Serre class $\sC$, let 
$S_{\sC} \subset \Mor \bC$ be the class consisting of those $s$ such that 
$\ker s \in \sC$ and $\coker s \in \sC$ $-$then $S_{\sC}$ admits a 
calculus of left and right fractions and 
$S_{\sC}$ = $\ov{S}_{\sC}$, i.e., $S_{\sC}$ is saturated.  
The metacategory $S_{\sC}^{-1}\bC$ is isomorphic to a category.  
As such, it is abelian and 
$L_{S_{\sC}}:\bC \ra S_{\sC}^{-1}\bC$ is exact and additive.  
An object \mX in \bC belongs to $\sC$ iff $L_{S_{\sC}}X$ is a zero object.  
Moreover, if \bD is an abelian category and $F:\bC \ra \bD$ is an exact functor, then $F$ can be factored through 
$L_{S_{\sC}}$ iff all the objects of $\sC$ are sent to zero objects by $F$.

[Note: \ Suppose that \bC is a Grothendieck category with a separator $U$ $-$then for any Serre class $\sC$,
$L_{S_{\sC}}:\bC \ra S_{\sC}^{-1}\bC$ has a right adjoint iff $\sC$ is closed under coproducts, in which case 
$S_{\sC}^{-1}\bC$ is again Grothendieck and has $L_{S_{\sC}}U$ as a separator.]\\

\begingroup%%------------------------------------>>%%----------------------------------->>
\fontsize{9pt}{11pt}\selectfont
Take $\bC = \bAB$ and let $\sC$ be the class of torsion abelian groups $-$then $\sC$ is a Serre class and 
$S_{\sC}^{-1}\bC$ is equivalenct to the category of torsion free divisible abelian groups or still, 
to the category of vector spaces over $\Q$.\\
\endgroup%%------------------------------------<< %%------------------------------------<<

Given a Grothendieck category \bC with a separator, a reflective subcategory \bD of \bC  is said to be a 
\un{Giraud subcategory}
\index{Giraud subcategory} 
provided that the reflector $R:\bC \ra \bD$ is exact.  
Every Giraud subcategory of \bC is Grothendieck and has a separator.  
There is a one-to-one correspondence between the 
Serre classes in \bC which are closed under coproducts and the Giraud subcategories of \bC.

%%----------------------------------------------------------------------------------------------40
[Note: \ The Gabriel-Popescu theorem says that every Grothendieck category with a separator is equivalent to a 
Giraud subcategory of \bAMOD for some $A$.]\\

\begingroup%%------------------------------------>>%%----------------------------------->>
\fontsize{9pt}{11pt}\selectfont
Attached to a topological space \mX is the category $\bOP(X)$ whose objects are the open subsets of \mX and whose morphisms are the inclusions.  The functor category $[\bOP(X)^\OP,\bAB]$ is the category of abelian presheaves on \mX.  
It is Grothendieck and has a separator.  The full subcategory of $[\bOP(X)^\OP,\bAB]$ whose objects are the abelian sheaves on \mX is a Giraud subcategory.\\
\endgroup%%------------------------------------<< %%------------------------------------<<

Fix a symmetric monoidal category \bV $-$then a 
\un{\bV-category}
\index{category! \bV-category} 
\bM consists of a class 
$\sO$ (the \un{objects}) and a function that assigns to each ordered pair $X, Y \in O$ an object $\HOM(X,Y)$ in \bV 
plus morphisms 
$C_{X,Y,Z}:\HOM(X,Y) \otimes \HOM(Y,Z) \ra \HOM(X,Z)$, 
$I_X:e \ra \HOM(X,X)$ satisfying the following conditions.

\indent\indent (\bV-cat$_1$) \ The diagram
\[
\begin{tikzcd}[sep=large]
{\HOM(X,Y) \otimes (\HOM(Y,Z) \otimes \HOM(Z,W))} \ar{d}[swap]{A}  \ar{rr}{\id \otimes C}
&&{\HOM(X,Y) \otimes \HOM(Y,W)}\ar{dd}{C}\\
{(\HOM(X,Y) \otimes \HOM(Y,Z)) \otimes \HOM(Z,W)} \ar{d}[swap]{C \otimes \id}\\
{\HOM(X,Z) \otimes \HOM(Z,W)} \ar{rr}[swap]{C} 
&&{\HOM(X,W)}
\end{tikzcd}
\]
commutes.

\label{13.61}
\indent\indent (\bV-cat$_2$) \ The diagram
\[
\begin{tikzcd}[sep=large]
{e \otimes \HOM(X,Y)} \ar{d}[swap]{I \otimes \id} \ar{r}{L}
&{\HOM(X,Y)} \arrow[d,shift right=0.5,dash] \arrow[d,shift right=-0.5,dash]
&{\HOM(X,Y) \otimes e} \ar{l}[swap]{R} \ar{d}{id \otimes I}\\
{\HOM(X,X) \otimes \HOM(X,Y)} \ar{r} [swap]{C}
&{\HOM(X,Y)}
&{\HOM(X,Y) \otimes \HOM(Y,Y) } \ar{l}{C}
\end{tikzcd}
\]
commutes.

[Note: \ The opposite of a \bV-category is a \bV-category and the product of two \bV-categories is a \bV-category.]

The 
\un{underlying category}
\index{underlying category} 
$\bU\bM$ of a \bV-category \bM has for its class of objects the class $O$, 
$\Mor(X,Y)$ being the set $\Mor(e,\HOM(X,Y))$.  Composition 
$\Mor(X,Y) \times \Mor(Y,Z$) $\ra$ $\Mor(X,Z)$ is calculated from 
$e \approx$ $e \otimes e$ 
$\overset{f \otimes g}{\lra}$ $\HOM(X,Y) \otimes \HOM(Y,Z)$ 
$\ra \HOM(X,Z)$, while $I_X$ serves as the identity in $\Mor(X,X)$.

[Note: \ A closed category \bV can be regarded as a \bV-category (take $\HOM(X,Y) = \hom(X,Y))$ and 
$\bU\bV$ is isomorphic to \bV.]\\

%%----------------------------------------------------------------------------------------------41
Every category is a \bSET-category and every additive category is an \bAB-category.\\

\label{13.63}
\begingroup%%------------------------------------>>%%----------------------------------->>
\fontsize{9pt}{11pt}\selectfont
A morphism $F:\bV \ra \bW$ of symmetric monoidal categories is a functor $F:\bV \ra \bW$, a morphism 
$\epsilon:e \ra Fe$, and morphisms $T_{X,Y}:FX \otimes FY \ra F(X \otimes Y)$ natural in \mX, \mY such that the diagrams 
\[
\begin{tikzcd}[sep=large]
{F_e \otimes FX} \ar{rr}{T} &&{F(e \otimes X)} \ar{d}{FL}\\
{e \otimes FX} \ar{u}{\epsilon \otimes \id} \ar{rr}[swap]{L} &&{FX}\\
\end{tikzcd}
\indent\indent
\begin{tikzcd}[sep=large]
{FX \otimes F_e} \ar{rr}{T} &&{F(X \otimes e)} \ar{d}{FR}\\
{FX \otimes e} \ar{u}{\id \otimes \epsilon} \ar{rr}[swap]{R} &&{FX}\\
\end{tikzcd}
\]
\[
\begin{tikzcd}[sep=large]
{FX \otimes (FY \otimes FZ)} \ar{d}[swap]{\id \otimes T} \ar{rr}{A} &&{(FX \otimes FY) \otimes FZ} \ar{d}{T \otimes \id}\\
{FX \otimes F(Y \otimes Z)} \ar{d}[swap]{T} &&{F(X \otimes Y) \otimes FZ} \ar{d}{T}\\
{F(X \otimes (Y \otimes Z))} \ar{rr}[swap]{FA} &&{F((X \otimes Y) \otimes Z)}
\end{tikzcd}
\]
commute with $F\Tee_{X,Y} \circ T_{X,Y} = T_{Y,X} \circ \Tee_{FX,FY}$.
\\ \indent
Example: Given a symmetric monoidal category \bV, the representable functor $\Mor(e,-)$ determines a morphism 
$\bV \ra \bSET$ of symmetric monoidal categories.
\\ \indent
Let $F:\bV \ra \bW$ be a morphism of symmetric monoidal categories.  Suppose that \bM is a \bV-category.  
Definition: $F_*\bM$ is the \bW-category whose object class is $O$, the rest of the data being $F\HOM(X,Y)$, 
$F\HOM(X,Y) \otimes F\HOM(Y,Z) \overset{T}{\lra}$ 
$F(\HOM(X,Y) \otimes \HOM(Y,Z)) \overset{FC}{\lra}$ 
$F\HOM(X,Z)$, 
$e \overset{\epsilon}{\ra}$ 
$Fe \overset{FI}{\lra}$ 
$F\HOM(X,X)$.
\\ \indent
[Note: \ Take \bW = \bSET and $F = \Mor(e,-)$ to recover $\bU\bM$.]\\
\endgroup%%------------------------------------<< %%------------------------------------<<

Fix a symmetric monoidal category \bV.  Suppose given \bV-cateogories \bM, \bN $-$then a 
\un{\bV-functor}
\index{functor! \bV-functor} 
$F:\bM \ra \bN$ is the specification of a rule that assigns to each object \mX in \bM an object $FX$ in \bN and the specification of a rule that assigns to each ordered pair $X,\ Y \in O$ a morphism 
$F_{X,Y}:\HOM(X,Y) \ra \HOM(FX,FY)$ in \bV such that the diagram
\[
\begin{tikzcd}[sep=large]
{\HOM(X,Y) \otimes \HOM(Y,Z)} \ar{d}[swap]{F_{X,Y} \otimes F_{Y,Z}} \ar{r}{C}
&{\HOM(X,Z)} \ar{d}{F_{X,Z}}\\
{\HOM(FX,FY) \otimes \HOM(FY,FZ)} \ar{r} [swap]{C}
&{\HOM(FX,FZ)}
\end{tikzcd}
\]
commutes with $F_{X,X} \circ I_X = I_{FX}$.

[Note: \ The 
\un{underlying functor}
\index{underlying functor} 
$UF:\bU\bM \ra \bU\bN$ sends \mX to $FX$ and $f:e \ra \HOM(X,Y)$ to $F_{X,Y} \circ f$.]

Example: $\HOM:\bM^\OP \times \bM \ra \bV$ is a \bV-functor if \bV is closed.\\

%%----------------------------------------------------------------------------------------------42
\label{16.25}
A \bV-category is 
\un{small}
\index{small (\bV-category)} 
if its class of objects is a set;
otherwise it is 
\un{large}
\index{large (\bV-category)}.  
\bV-\bCAT 
\index{\bV-\bCAT}, 
the category of small \bV-categories and \bV-functors, is a symmetric monoidal category.\\

\begingroup%%------------------------------------>>%%----------------------------------->>
\fontsize{9pt}{11pt}\selectfont
Take \bV = \bAB $-$then an additive functor between additive categories ``is'' a \bV-functor.\\
\endgroup%%------------------------------------<< %%------------------------------------<<

Fix a symmetric monoidal category \bV.  
Suppose given \bV-categories \bM, \bN and \bV-functors 
$F, G:\bM \ra \bN$ $-$then a 
\un{\bV-natural transformation}
\index{natural transformation! \bV-natural transformation} 
$\Xi$ from $F$ to $G$ is a class of morphisms $\Xi_X:e \ra \HOM(FX,GX)$ for which the diagram
%\[
%\begin{tikzcd}[sep=large]
%{e \otimes \HOM(X,Y)}\ar{rr}{\Xi_X \otimes G_{X,Y}} &&{\HOM(FX,GX) \otimes \HOM(GX,GY)}\ar{dd}{C}\\
%\\
%{\HOM(X,Y)}\ar{dd}[swap]{R^{-1}} \ar{uu}{L^{-1}}   &&{\HOM(FX,GY)}\\
%\\
%{\HOM(X,Y) \otimes e} \ar{rr}[swap]{F_{X,Y} \otimes \Xi_Y} 
%&&{\HOM(FX,FY) \otimes \HOM(FY,GY)} \ar{uu}[swap]{C}
%\end{tikzcd}
%\]
\[
\begin{tikzcd}[sep=large]
{e \otimes \HOM(X,Y)}\ar{rr}{\Xi_X \otimes G_{X,Y}} &&{\HOM(FX,GX) \otimes \HOM(GX,GY)}\ar{d}{C}\\
{\HOM(X,Y)}\ar{d}[swap]{R^{-1}} \ar{u}{L^{-1}}   &&{\HOM(FX,GY)}\\
{\HOM(X,Y) \otimes e} \ar{rr}[swap]{F_{X,Y} \otimes \Xi_Y} 
&&{\HOM(FX,FY) \otimes \HOM(FY,GY)} \ar{u}[swap]{C}
\end{tikzcd}
\]
commutes.\\

\begingroup%%------------------------------------>>%%----------------------------------->>
\fontsize{9pt}{11pt}\selectfont
Assume that \bV is complete and closed.  Let \bM, \bN be \bV-categories with \bM small $-$then the category 
$\bV[\bM,\bN]$ whose objects are the \bV-functors $\bM \ra \bN$ and whose morphisms are the \bV-natural transformations is a \bV-category if 
$\HOM(F,G) = \ds\int_X \HOM(FX,GX)$, the equalizer of 
$\ds\prod\limits_{X \in O} \HOM(FX,GX) \rightrightarrows$ 
$\ds\prod\limits_{X^\prime,X\pp \in O} \hom(\HOM(X^\prime,X\pp),\HOM(FX^\prime,GX\pp))$.\\
\endgroup%%------------------------------------<< %%------------------------------------<<

Let \bC be a category with pullbacks $-$then an 
\un{internal category}
\index{internal category} 
(or a 
\un{category object})
\index{category object} 
in \bC consists of an object $M$, an object $O$, and morphisms 
$s:M \ra O$, 
$t:M \ra O$, 
$e:O \ra M$, 
$c:M \times_O M \ra M$ satisfying the usual category theoretic relations (here, 
\begin{tikzcd}[sep=large]
{M \times_O M} \ar{d} \ar{r} &{M} \ar{d}{t}\\
{M} \ar{r}[swap]{s} &{O}
\end{tikzcd}
).  Notation: $\bM = (M,O,s,t,e,c)$. 
\index{$(M,O,s,t,e,c)$}
\index{internal category,category object! (M,O,s,t,e,c)}

[Note: \ There are obvious notions of
\un{internal functor}
\index{internal functor} and 
\un{internal natural}
\un{transformation}.]
\index{internal natrual transformation}\\

\label{17.17}
\begingroup%%------------------------------------>>%%----------------------------------->>
\fontsize{9pt}{11pt}\selectfont
An internal category in \bSET is a small category.  
An internal category in \bSISET is a simplicial object in \bCAT.
\\ \indent
An internal category in \bCAT is a (small) 
\un{double category}.
\index{double category}
\\ \indent
[Note: \ Spelled out, such an entity consists of objects $X,Y, \ldots$, horizontal morphisms $f,g, \ldots$, 
vertical morphisms $\phi,\psi, \ldots$, and bimorphisms (represented diagramatically by squares).  
The objects and 
%%----------------------------------------------------------------------------------------------43
the horizontal morphisms form a category with identities $X \overset{h_X}{\lra} X$.  The objects and the vertical morphisms form a category with identities 
\begin{tikzcd}[sep=large]
{X} \ar{d}[swap]{v_X}\\
{X}
\end{tikzcd}
.  The bimorphisms have horizontal and vertical laws of composition
\begin{tikzcd}[sep=large]
{\bullet}\ar{d}\ar{r} &{\bullet}\ar{d}\ar{r} &{\bullet}\ar{d}\\
{\bullet}\ar{r} &{\bullet}\ar{r} &{\bullet}
\end{tikzcd}
,
\begin{tikzcd}[sep=large]
{\bullet}\ar{d}\ar{r} &{\bullet}\ar{d}\\
{\bullet}\ar{d}\ar{r} &{\bullet}\ar{d}\\
{\bullet}\ar{r} &{\bullet}
\end{tikzcd}
under which they form a category with identities
\begin{tikzcd}[ sep=small]
{X}\ar{dd}[swap]{\phi}\ar{rr}{h_X} &&{X}\ar{dd}{\phi}\\
&{\id_\phi}\\
{Y} \ar{rr}[swap]{h_Y} &&{Y}\\
\end{tikzcd}
,
\begin{tikzcd}[ sep=small]
{X}\ar{dd}[swap]{v_X}\ar{rr}{f} &&{Y}\ar{dd}{v_Y}\\
&{\id_f}\\
{X} \ar{rr}[swap]{f} &&{Y}\\
\end{tikzcd}
.  In the situation
\begin{tikzcd}[sep=large]
{\bullet}\ar{d}\ar{r} &{\bullet}\ar{d}\ar{r} &{\bullet}\ar{d}\\
{\bullet}\ar{d}\ar{r} &{\bullet}\ar{d}\ar{r} &{\bullet}\ar{d}\\
{\bullet}\ar{r} &{\bullet}\ar{r} &{\bullet}
\end{tikzcd}
, 
the result of composing horizontally and then vertically is the same as the result of composing vertically and then horizontally.  
Furthermore, horizontal composition of vertical identities gives a vertical identity and vertical composition of horizontal identities gives a horizontal identity.  Finally, the horizontal and the vertical identities
\begin{tikzcd}[ sep=small]
{X}\ar{dd}[swap]{v_X}\ar{rr}{h_X} &&{X}\ar{dd}{v_X}\\
&{\id_{v_X}}\\
{X} \ar{rr}[swap]{h_X} &&{X}\\
\end{tikzcd}
,
\begin{tikzcd}[ sep=small]
{X}\ar{dd}[swap]{v_X}\ar{rr}{h_X} &&{X}\ar{dd}{v_X}\\
&{\id_{h_X}}\\
{X} \ar{rr}[swap]{h_X} &&{X}\\
\end{tikzcd}
coincide.]
\\ \indent
Example: Let \bC be a small category $-$then db\bC is the double category whose objects are those of \bC, 
whose horizontal and vertical morphisms are those of \bC, and whose bimorphisms are the commutative squares in \bC.  
All sources, targets, identities, and compositions come from \bC.\\
\endgroup%%------------------------------------<< %%------------------------------------<<

Let \bC be a category with pullbacks.  Given  an object $O$ in \bC, an 
\un{$O$-graph}
\index{O-graph} %\index{$O$-graph} 
is an object \mA and a pair of morphisms 
$s, t:A \ra O$.  $O$-\bGR 
\index{$O$-\bGR } 
is the category whose objects are the $O$-graphs and whose morphisms 
$(A,s,t) \ra (A^\prime,s^\prime,t^\prime)$ are the arrows $f:A \ra A^\prime$ such that 
$s = s^\prime \circ f$, 
$t = t^\prime \circ f$.  If $A \times_O A^\prime$ is defined by the pullback square 
\begin{tikzcd}[sep=large]
{A \times_O A^\prime} \ar{d} \ar{r} &{A^\prime} \ar{d}{t^\prime}\\
{A} \ar{r}[swap]{s} &{O}
\end{tikzcd}
and if the structural morphisms are 
$A \times_O A^\prime \ra$ $A^\prime \overset{s^\prime}{\ra} O$, 
$A \times_O A^\prime \ra$ $A \overset{t}{\ra} O$, then $A \times_O A^\prime$ is an $O$-graph.  
Therefore  $O$-\bGR is a monoidal category: 
Take $A \otimes A^\prime$ to be $A \times_O A^\prime$ and let $e$ be 
$(O, \id_O, \id_O)$.  A monoid \bM in $O$-\bGR is an internal category in \bC with object element $O$.

Let \bC be a category with pullbacks.  Given an internal category \bM in \bC the 
\un{nerve} 
\index{nerve (internal category)}
%%----------------------------------------------------------------------------------------------44
$\ner \bM$ of \bM is the simplicial object in \bC defined by 
$\ner_0 \bM = O$, 
$\ner_1 \bM = M$, 
$\nersub_n \bM = M \times_O  \cdots \times_O M$ ($n$ factors).  At the bottom, 
$
\begin{cases}
\ d_0\\
\ d_1
\end{cases}
:\ner_1 \bM \ra \ner_0 \bM
$
is 
$
\begin{cases}
\ t\\
\ s
\end{cases}
, \ 
$
while higher up, in terms of the underlying projections, 
$d_0 = (\pi_1, \dots, \pi_{n-1})$, 
$d_n = (\pi_2, \dots, \pi_{n})$, 
%dmc - note the discrepancy in the next line with the online version which probably has two errors - revert back to orig
%$d_i = (\pi_1, \dots, c \circ (\pi_{n-i},\pi_{n-i+1}),\ldots, \pi_n)$ $(0 < i < n)$, %orig orig and imo much better
$d_i = (\pi_1, \dots, c \circ \pi_{n-i, n-i+1},\ldots, \pi_n)$ $(0 < i < n)$, 
and at the bottom, 
$s_0:\ner_0 \bM \ra \ner_1 \bM$ is $e$, while higher up, 
$s_i = e_i \circ \sigma_i$, where $\sigma_i$ inserts $O$ at the $n - i +1$ spot and $e_i$ is 
$\id \times_O \cdots \times_O e \times_O \cdots \times_O \id$ placed accordingly $(0 \leq i \leq n)$.

[Note: \ An internal functor $\bM \ra \bM^\prime$ induces a morphism 
$\ner \bM \ra \ner \bM^\prime$ of simplicial objects.]\\

\begingroup%%------------------------------------>>%%----------------------------------->>
\fontsize{9pt}{11pt}\selectfont
Suppose that \bC is a small category.  Consider $\ner \bC$ $-$then an element $f$ of $\nersub_n \bC$ is a diagram of the form 
$X_0 \overset{f_0}{\ra} X_1 \ra \cdots \ra X_{n-1} \overset{f_{n-1}}{\ra} X_n$ and 
\[
d_i f= 
\begin{cases}
\ X_1 \ra \cdots \ra X_n \hspace{5.48cm} \  (i =0)\\
\ X_0 \ra \cdots \ra X_{i-1}   
\begin{tikzcd}
{} \ar{rr}{f_i \circ f_{i-1}} &&X_{i+1}
\end{tikzcd}
\ra \cdots \ra X_n \qquad (0 < i < n)\\
\ X_0 \ra \cdots \ra X_{n-1} \hspace{5.25cm}  (i= n) 
\end{cases}
\hspace{-.3cm},
\]
$s_i f = X_0 \ra \cdots \ra X_i \overset{\id_{X_i}}{\ra}  X_i \ra \cdots \ra X_n$.  
The abstract definition thus reduces to these formulas since $f$ corresponds to the $n$-tuple $(f_{n-1}, \ldots, f_0)$.\\
\endgroup%%------------------------------------<< %%------------------------------------<<

Let \bC be a category with pullbacks.  Given an internal category \bM in \bC, a 
\un{left \bM-object}
\index{left \bM-object} 
is an object $T:Y \ra O$ in $\bC/O$ and a morphism 
$\lambda:M \times_O Y \ra Y$ such that
\[
\begin{tikzcd}[sep=large]
{M \times_O M \times_O Y} \ar{d}[swap]{\id \times_O \lambda} \ar{r}{c \times_O \id}
&{M \times_O Y} \ar{d}{\lambda}
&{O \times_O Y} \ar{l}[swap]{e \times_O \id} \ar{d}{L}\\
{M \times_O Y} \ar{r}[swap]{\lambda}
&{Y} \arrow[r,shift right=0.5,dash] \arrow[r,shift right=-0.5,dash] 
&{Y}
\end{tikzcd}
\]
and
\begin{tikzcd}[sep=large]
{M \times_O Y} \ar{d} \ar{r}{\lambda} &{Y} \ar{d}{T}\\
{M} \ar{r}[swap]{t} &{O}
\end{tikzcd}
commute, where $M \times_O Y$ is defined by the pullback square
\begin{tikzcd}[sep=large]
{M \times_O Y} \ar{d} \ar{r} &{Y} \ar{d}{T}\\
{M} \ar{r}[swap]{s} &{O}
\end{tikzcd}
.  
Example: Take \bC = \bSET $-$then \bM is a small category and the category of left \bM-objects is equivalent to the functor category $[\bM,\bSET]$.

[Note: A 
\un{right \bM-object}
\index{right \bM-object} 
is an object $S:X \ra$ $O$ in $\bC/O$ and a morphism 
$\rho:X \times_O M$ $\ra$ $X$ such that the analogous diagrams commute, where $X \times_O M$ is defined 
%%----------------------------------------------------------------------------------------------45
by the pullback square 
\begin{tikzcd}[sep=large]
{X \times_O M} \ar{d} \ar{r} &{M} \ar{d}{t}\\
{X} \ar{r}[swap]{S} &{O}
\end{tikzcd}
Example: Take \bC = \bSET $-$then \bM is a small category and the category of right \bM-objects 
is equivalent to the functor category $[\bM^\OP,\bSET]$.

\label{17.21}
Let \bC be a category with pullbacks.  
Given an internal category \bM in \bC and a left \bM-object \mY, the 
\un{translation category}
\index{translation category ($\tran Y$)} 
$\tran Y$ of \mY  is the category object 
$\bM_Y = (M_Y, O_Y, s_Y, t_Y, e_Y, c_Y)$ in \bC, where 
$M_Y = M \times_O Y$, 
$O_Y = Y$, $s_Y$ is the projection
$M \times_O Y \ra Y$, 
$t_Y$ is the action $\lambda:M \times_O Y \ra Y$, and 
$e_Y$, $c_Y$ are derived from $e:O \ra M$, $c:M \times_O M \ra M$.  
Example: Take \bC = \bSET, let \bM be a small category, and suppose that $G:\bM \ra \bSET$ is a functor 
$-$then $G$ determines a left \bM-object $Y_G$ and the translation category of $Y_G$ 
can be identified with the Grothendieck construction on $G$.\\

\begingroup%%------------------------------------>>%%----------------------------------->>
\fontsize{9pt}{11pt}\selectfont
Let $G$ be a semigroup with unit, \bG the category having a single object $*$ with $\Mor(*,*) = G$.  
Suppose that $Y$ is a left $G$-set, i.e., an object in $\bLACT_G$ or still, a left \bG-object.  
The translation category of \mY is $(G \times Y,Y,s_Y,t_Y,e_Y,c_Y)$, where 
$s_Y(g,y) = y$, 
$t_Y(g,y) = g \cdot y$,  
$e_Y(y) = (e,y)$, 
$c_Y((g_2,y_2),(g_1,y_1))$ $=$ $(g_2g_1,y_1)$.  
Specialize and let $Y = G$ 
$-$then the objects of the translation category of $G$ are the elements of $G$ and 
$\Mor(g_1,g_2) \approx \{g:gg_1 = g_2\}$.\\
\endgroup%%------------------------------------<< %%------------------------------------<<

\label{13.56}
\label{14.65}
Let \bC be a category with pullbacks.  
Given an internal category \bM in \bC, and a right \bM-object \mX and a left 
\bM-object \mY, the 
\un{bar construction}
\index{bar construction (internal category)} 
$\barr(X;\bM;Y)$ 
\index{$\barr(X;\bM;Y)$} 
on $(X,Y)$ is the simplicial object in \bC defined by 
$\barr_n(X;\bM;Y) = X \times_O \nersub_n\bM \times_O Y$.  
Note that $\rho$ appears only in $d_n$ and $\lambda$ appears only in $d_0$.  
The 
\un{translation category}
\index{translation category ($\tran (X,Y)$ )} 
$\tran (X,Y)$ 
\index{$\tran (X,Y)$} 
of $(X,Y)$ is the category object 
$\bM_{X,Y} = (M_{X,Y},O_{X,Y},s_{X,Y},t_{X,Y},e_{X,Y},c_{X,Y})$ in \bC, where
$M_{X,Y} = X \times_O M \times_O Y$,
$O_{X,Y} = X \times_O Y$, 
$s_{X,Y} = \rho \times_O \id_Y$, 
$t_{X,Y} = \id_X \times_O \lambda$, 
$e_{X,Y}$ $\&$ 
$c_{X,Y}$ being definable in terms of $e$ $\&$ $c$.  
Therefore 
$\barr(X;\bM;Y) \approx \ner \bM_{X,Y}$.  
Example: O can be viewed as a right \bM-object via 
$O \times_O M \overset{L}{\ra} M \overset{s}{\ra} O$ and as a left \bM-object via 
$M \times_O O \overset{R}{\ra} M \overset{t}{\ra} O$, and \mM can be viewed as a right \bM-object via 
$M \times_O M \overset{c}{\ra} M \overset{s}{\ra} O$ and as a left \bM-object via 
$M \times_O M \overset{c}{\ra} M \overset{t}{\ra} O$, so 
$\barr(O;\bM;O)$, 
$\barr(O;\bM;M)$,
$\barr(M;\bM;O)$, 
$\barr(M;\bM;M)$ are meaningful.\\

\label{14.122}
\begingroup%%------------------------------------>>%%----------------------------------->>
\fontsize{9pt}{11pt}\selectfont
Let $G$ be a group, \bG the groupoid having a single object $*$ with $\Mor(*,*) = G$.  
View $G$ as a left $G$-set $-$then $\barr(*;\bG;G)$ is isomorphic to the nerve of $\grd G$.  
In fact, the objects of $\grd G$ are the elements of $G$ and the morphisms of $\grd G$ are the elements of 
$G \times G$ 
$(s(g,h) = g$, 
$t(g,h) = h$, 
$\id_g = (g,g)$, 
$(h,k) \circ (g,h) $= $(g,k)$), thus 
$\nersub_n \grd G $= $G \times \cdots \times G$ ($n+1$ factors) and 
$d_i(g_0, \ldots, g_n) $= $(g_0, \ldots, \widehat{g_i}, \ldots, g_n)$, 
$s_i(g_0, \ldots, g_n) $=$ (g_0, \ldots, g_i, g_i, \ldots, g_n)$.  On the other hand, 
$\barr(*;\bG;G)$ is the nerve of the translation category of $G$.  The functor 
$\tran G \ra \grd G$ which is the identity on objects and sends a morphism $(g,h)$ 
%%----------------------------------------------------------------------------------------------46
in $\tran G$ to the morphism $(h,g\cdot h)$ in $\grd G$ induces an isomorphism 
$\ner \tran G \ra \ner \grd G$ of simplicial sets.  
For 
$(g_0, \ldots, g_n) \ra (g_n, g_{n-1}g_n, \ldots, g_0 \cdots g_n)$ is the arrow 
$\nersub_n \  \tran G \ra \nersub_n \  \grd G$, its inverse being 
$(g_0, \ldots, g_n) \ra (g_n g_{n-1}^{-1}, g_{n-1}g_{n-2}^{-1},\ldots, g_0)$.    
Both $\ner \tran G$ and $\ner \grd G$ are simplicial right $G$-sets, viz. 
$(g_0, \ldots, g_n) \cdot g = (g_0, \ldots, g_n g)$ and 
$(g_0, \ldots, g_n) \cdot g = (g_0g, \ldots, g_n g)$, and the isomorphism 
$\ner \tran G \ra \ner \grd G$ is equivariant.\\
\endgroup%%------------------------------------<< %%------------------------------------<<

\label{14.27}
\label{14.126}
\label{14.127a}
\label{14.140}
\label{14.162}
Let $\bT = (T,m,\epsilon)$ be a triple in a category \bC $-$then a 
\un{right \bT-functor}
\index{right \bT-functor} 
in a category \bV is a functor $F:\bC \ra \bV$
plus a natural transformation $\rho:F \circ T \ra F$ such that the diagrams 
\begin{tikzcd}[sep=large]
{F \circ T \circ T} \ar{d}[swap]{Fm} \ar{r}{\rho T} &{F \circ T} \ar{d}{\rho}\\
{F \circ T} \ar{r}[swap]{\rho} &{F}
\end{tikzcd}
,
\begin{tikzcd}[sep=large]
{F} \ar[equals]{rd} \ar{r}{F \epsilon} &{F \circ T} \ar{d}{\rho}\\
&{F}
\end{tikzcd}
commute and a 
\un{left \bT-functor}
\index{left \bT-functor} 
in a category \bU is a functor $G:\bU \ra \bC$ 
plus a natural transformation $\lambda:T \circ G \ra G$ such that the diagrams 
\begin{tikzcd}[sep=large]
{T \circ T \circ G} \ar{d}[swap]{mG} \ar{r}{T \lambda} &{T \circ G} \ar{d}{\lambda}\\
{T \circ G} \ar{r}[swap]{\lambda} &{G}
\end{tikzcd}
,
\begin{tikzcd}[sep=large]
{G} \ar[equals]{rd} \ar{r}{\epsilon G} &{T \circ G} \ar{d}{\lambda}\\
&{G}
\end{tikzcd}
commute.  The 
\un{bar construction}
\index{bar construction} 
$\barr(F;\bT,G)$ on $(F,G)$ is the simplicial object in 
$[\bU,\bV]$ defined by 
$\barr_n(F;\bT;G) = F \circ T^n \circ G$, where 
$d_0 = \rho T^{n-1}G$, 
$d_i = F T^{i-1}mT^{n-i-1}G$ $(0 < i < n)$,  
$d_n = FT^{n-1}\lambda$, and 
$s_i = FT^i\epsilon T^{n-i}G$.  
In particular: 
$\barr_1(F;\bT;G) = F \circ T \circ G$, 
$\barr_0(F;\bT;G) = F \circ G$, and 
$d_0, d_1:F \circ T \circ G \ra F \circ G$ are $\rho G$, $F \lambda$, while 
$s_0:F \circ G \ra F \circ T \circ G$ is $F \epsilon G$.

Example: If \mX is a \bT-algebra in \bC with structural morphism $\xi:TX \ra X$, then \mX determines a left \bT-functor $G:\bone \ra \bC$ and one writes $\barr(F;\bT;X)$ for the associated bar construction.\\

\begingroup%%------------------------------------>>%%----------------------------------->>
\fontsize{9pt}{11pt}\selectfont
Take \bV = \bC, $F = T$, $\rho = m$, and put $\tau = \epsilon TG$ (thus 
$\tau:T \circ G \ra T \circ T \circ G$).  There is a commutative diagram 
\[
\begin{tikzcd}[sep=large]
{T \circ G} \arrow[dd,shift right=0.5,dash] \arrow[dd,shift right=-0.5,dash]  \ar{rd}{\tau} \ar{rrr}{\lambda}
%{T \circ G}  \ar[equals]{dd} \ar{rd}{r} \ar{rrr}{\lambda}
&&&{G} \ar{ld}[swap]{\epsilon G} \arrow[dd,shift right=0.5,dash] \arrow[dd,shift right=-0.5,dash] \\
%&&&{G} \ar{ld}[swap]{\epsilon G} \ar[equals]{dd}\\
&{T \circ T \circ G} \ar{ld}{mG} \ar{r}{T\lambda} &{T \circ G} \ar{rd}[swap]{\lambda}\\
{T \circ G}  \ar{rrr}[swap]{\lambda} &&&{G} 
\end{tikzcd}
\]
from which it follows that $\lambda:T \circ G \ra G$ is a coequalizer of $(d_0,d_1) = (mG,T\lambda)$.  
Consider the string of arrows 
$T \circ T^n \circ G \overset{d_0}{\lra}$ 
$T \circ T^{n-1} \circ G \lra$ 
$\cdots \lra$
$T \circ T \circ G \overset{d_0}{\lra}$ 
$T \circ G \overset{\lambda}{\lra}$ 
$G \overset{\epsilon G}{\lra}$ 
$T \circ G \overset{s_0}{\lra}$ 
$T \circ T \circ G \ra \cdots \lra$ 
$T \circ T^{n-1} \circ G \overset{s_0}{\lra}$ 
$T \circ T^n \circ G$.  
Viewing $G$ as a constant simplicial object in $[\bDelta^{\OP},[\bC,\bV]]$, there are simplicial morphisms 
$G \ra \barr(T;\bT;G)$, 
$\barr(T;\bT;G) \ra G$ viz. 
$s_0^n \circ \epsilon G: G \ra T \circ T^n \circ G$, 
$\lambda \circ d_0^n:T \circ T^n \circ G \ra G$, and the composition 
$G \ra \barr(T;\bT;G) \ra G$ is the identity.  On the other hand, if 
$h_i:T \circ T^n \circ G \ra T \circ T^{n+1} \circ G$ 
%%----------------------------------------------------------------------------------------------47
is defined by 
$h_i = s_0^i(\epsilon T^{n-i+1}G)d_0^i$ $(0 \leq i \leq n)$, then 
$d_0 \circ h_0 = \id$, 
$d_{n+1} \circ h_n = s_0^n \circ \epsilon G \circ \lambda \circ d_0^n$, and
\[
d_i \circ h_j = 
\begin{cases}
\ h_{j-1} \circ d_i \quadx (i < j)\\
\ d_i \circ h_{i-1} \quadx (i = j > 0)\\
\ h_j \circ d_{i-1} \quadx (i > j + 1) 
\end{cases}
%$
, \ 
%$
s_i \circ h_j = 
\begin{cases}
\ h_{j+1} \circ s_i \quadx (i \leq j)\\
\ h_j \circ s_{i-1} \quadx (i > j)
\end{cases}
.
\]
\\ \indent
[Note: \ Take instead $\bU = \bC$, $G = T$, $\lambda = m$ $-$then with 
$\tau = FT\epsilon$, $\rho:F \circ T \ra F$ is a coequalizer of $(d_1,d_0) = (Fm, \rho T)$ 
and the preceding observations dualize.]\\
\endgroup%%------------------------------------<< %%------------------------------------<<

%%%%%%%%%%%%%%%%%%%%%%%%%%%%%%%%%%%%%%
%%%%%%%%%%%%%%%%%%%%%%%%%%%%%%%%%%%%%%
%%%%%%%%%%%%%%%%%%%%%%%%%%%%%%%%%%%%%%

\begin{center}
$\S \ 0$
\\[0.5cm]
$\mathcal{REFERENCES}$\\[-.2cm]
\end{center}

\[
\textbf{BOOKS}
\]

\begingroup
\fontsize{9pt}{11pt}\selectfont
\setlength\parindent{0 cm}

[1] \quad Ad\'amek, J., Herrlich, H., and Strecker, G., \textit{Abstract and Concrete Categories}, John Wiley (1990).
\\[-.2cm]

[2] \quad Borceux, F., \textit{Handbook of Categorical Algebra 1, 2, 3}, Cambridge University Press (1994).
\\[-.2cm]

[3] \quad Dubuc, E., \textit{Kan Extensions in Enriched Category Theory}, Springer Verlag (1970).
\\[-.2cm]

[4] \quad Gabriel, P. and Ulmer, F., \textit{Lokal Pr\"asentierbare Kategorien}, Springer Verlag  (1971).
\\[-.2cm]

[5] \quad Herrlich, H., and Strecker, G., \textit{Category Theory}, Heldermann Verlag (1979).
\\[-.2cm]

[6] \quad Kelly, G., \textit{Basic Concepts of Enriched Category Theory}, Cambridge University Press (1982).
\\[-.2cm]

[7] \quad MacLane, S., \textit{Categories for the Working Mathematician}, Springer Verlag (1971).
\\[-.2cm]

[8] \quad Osborne, M., \textit{Basic Homological Algebra}, Springer Verlag (2000).
\\[-.2cm]

[9] \quad Popescu, N. and Popescu, L., \textit{Theory of Categories}, Sijthoff and Noordhoff (1979).
\\[-.2cm]

[10] \quad Porter, J. and Woods, R., \textit{Extensions of Absolutes of Hausdorff Spaces}, Springer Verlag (1988).
\\[-.2cm]

[11] \quad Preuss, G., \textit{Theory of Topological Structures}, Reidel (1988).
\\[-.2cm]

[12] \quad Schubert, H., \textit{Categories}, Springer Verlag (1972)
.\\[-.2cm]
\endgroup

\[
\textbf{ARTICLES}
\]

\begingroup
\fontsize{9pt}{11pt}\selectfont
\setlength\parindent{0 cm}

[1] \quad Almkvist, G., Fractional Categories, \textit{Ark. Mat.} \textbf{7} (1968), 449-476.
\\[-.2cm]

[2] \quad Eilenberg, S. and Kelly, G., Closed Categories In: \textit{Proceedings of the Conference on Categorical}

\hspace{0.8cm}\textit{Algebra}, S. Eilenberg, et al. (ed.), Springer Verlag  (1966), 421-562.
\\[-.2cm]

[3] \quad Gabriel, P., Des Cat\'egories Ab\'eliennes, \textit{Bull. Soc. Math. France} \textbf{90} (1962), 323-448.
\\[-.2cm]

[4] \quad Grothendieck, A., Sur Quelques Points d'Alg\`ebre Homologique, \textit{Tohoku Math. J.} \textbf{9} (1957), 119-221.
\\[-.2cm]

[5] \quad Herrlich, H., Categorical Topology, 1971-1981, In: \textit{General Topology and its Relations to Modern}

\hspace{0.8cm}\textit{Analysis and Algebra V}, J. Novak (ed.), Heldermann Verlag (1982), 279-383.
\\[-.2cm]

[6] \quad Isbell, J., Structure of Categories, \textit{Bull Amer. Math. Soc.} \textbf{72} (1966), 619-655.
\\[-.2cm]

[7] \quad Kan, D., Adjoint Functors, \textit{Trans. Amer. Math. Soc.} \textbf{87} (1958), 294-329.
\\[-.2cm]

[8] \quad Kelly, G., A Unified Treatment of Transfinite Constructions for Free Algebras, Free Monoids, Colimits, 

\hspace{0.8cm}Associated Sheaves, \ldots, \textit{Bull. Austral. Math. Soc.} \textbf{22} (1980), 1-83.
\\[-.2cm]

[9] \quad Kelly, G., A Survey of Totality for Enriched and Ordinary Categories, \textit{Cahiers Topol. Geom. Diff.}

\hspace{0.8cm}\textit{Cat.} \textbf{27} (1986), 109-132.
\\[-.2cm]

[10] \quad Kelly, G. and Street, R., Review of the Elements of 2-Categories, In: \textit{Category Seminar}, G. Kelly 

\hspace{0.95cm}(ed.), Springer Verlag (1974), 75-103.
\\[-.2cm]

[11] \quad MacLane, S., Categorical Algebra, \textit{Bull. Amer. Math. Soc.} \textbf{71} (1965), 40-106.
\\[-.2cm]

[12] \quad MacLane, S., Topology and Logic as a Source of Algebra, \textit{Bull. Amer. Math. Soc.} \textbf{82} (1976), 1-40.
\\[-.2cm]

[13] \quad Mitchell, B., Rings with Several Objects, \textit{Adv. Math.} \textbf{8} (1972), 1 - 161.
\\[-.2cm]

[14] \quad Nakagawa, R., Categorical Topology, In: \textit{Topics in General Topology}, K. Morita and J. Nagata (ed.), 

\hspace{0.95cm}North Holland (1989), 563-623.
\\[-.2cm]

[15] \quad Tholen, W., Factorizations, Localizations, and the Orthogonal Subcategory Problem, \textit{Math. Nachr.} 

\hspace{0.95cm}\textbf{114} (1983), 63-85.
\\[-.2cm]

[16] \quad Tsalenko, M. and Shul'geifer, E.,  Categories, \textit{J. Soviet Math.} \textbf{7} (1977), 532-586.

\setlength\parindent{2em}

\endgroup

\pagenumbering{bychapter}
\chapter{
$\boldsymbol{\S}$\textbf{1}.\quadx  COMPLETETLY REGULAR HAUSDORFF SPACES}
\setlength\parindent{2em}
\setcounter{proposition}{0}
%%----------------------------------------------------------------------------------------------01
$\text{ }$\\[-1.25cm]

The reader is assumed to be familiar with the elements of general topology.  
Even so, I think it best to provide a summary of what will be needed in the sequel.  Not all terms will be defined; most proofs will be omitted.

Let $X$ be a locally compact Hausdorff space(LCH space).\\

\begin{proposition} \ 
A subspace of $X$ is locally compact iff it is locally closed, i.e., has the form $A \cap U$, 
where $A$ is closed and $U$ is open in $X$.\\
\end{proposition}

\begingroup%%----------------------------------->>
\fontsize{9pt}{11pt}\selectfont
The class of nonempty LCH spaces is closed under the formation in \bTOP 
of finite products and arbitrary coproducts.

[Note: \ An arbitrary product of nonempty LCH spaces is a LCH space iff 
all but finitely many of the factors are compact.]\\
\endgroup%%------------------------------------<<

In practice, various additional conditions are often imposed on a LCH space $X$.  
The connections among the most common of these can be summarized as follows:
\\[-.2cm]

\begin{tikzpicture}
\node at (-5,0) {compact metrizable};
\node at (-1.5,1.25) {metrizable};
\node at (-1.5,-1.25) {compact};
\node at (1.25,-2.5) {Lindel\"of};

\node at (3.75,1.25) {paracompact};
\node at (3.75,-1.25) {$\sigma$-compact};
\node at (7.25,1.25) {normal};

\draw[->] (-3.25,0) -- (-2.4,1.05) node[below] {$$};
\draw[->] (-3.25,0) -- (-2.2,-1.1) node[above] {$$};

\draw[->] (-.5,1.25) -- (2.6, 1.25) node[above] {$$};

\draw[->] (-.7,-1.25) -- (2.75, -1.25) node[above] {$$};
\draw[->] (-.7,-1.25) -- (0.55, -2.35) node[above] {$$};

\draw[->] (2.,-2.5) -- (2.8, -1.35) node[above] {$$};

\draw[->] (3.75,-1.05) -- (3.75, 1.05) node[above] {$$};

\draw[->] (4.9,1.25) -- (6.6, 1.25) node[above] {$$};
\end{tikzpicture}

\begingroup%%----------------------------------->>
\fontsize{9pt}{11pt}\selectfont
\textbf{EXAMPLE} \ Let $\Omega$ be the first uncountable ordinal and consider $[0,\Omega]$ (in the order topology) 
$-$then $[0,\Omega]$ is Hausdorff.  
And: 
(i) $[0,\Omega]$ is compact but not metrizable; 
(ii) $[0,\Omega[$ is locally compact and normal but not paracompact; 
(iii) $[0,\Omega] \times [0,\Omega[$ is locally compact but not normal.\\

\endgroup%%------------------------------------<<

Here are some important points to keep in mind.
\\[-.5cm]

\indent\indent (LCH$_1$) $X$ is completely regular,
 i.e., $X$ has enough real valued continuous functions to separate points and closed sets in the sense that for every point 
 $x \in X$ and every closed subset $A \subset X$ not containing $x$, there exists a continuous function 
 $\phi: X\ra [0,1]$ such that $\phi(x) = 1$, $\restr{\phi}{A} = 0$.

\indent\indent (LCH$_2$) $X$ is $\sigma$-compact iff $X$ possesses a \un{sequence of exhaustion}, i.e., an increasing sequence $\{U_n\}$ of relatively compact open sets $U_n \subset X$ such that $\overline{U}_n \subset U_{n+1}$ and $X = \bigcup\limits_n U_n$.
\index{sequence of exhaustion}
%%----------------------------------------------------------------------------------------------02
\label{2.12}
\label{19.6}

\indent\indent (LCH$_3$) $X$ is paracompact iff $X$ admits a representation $X = \coprod\limits_i X_i$, where the $X_i$ are pairwise disjoint nonempty open $\sigma$-compact subspaces of $X$.
\label{1.14}
\label{19.19}

\indent\indent (LCH$_4$) $X$ is second countable iff $X$ is $\sigma$-compact and metrizable.

\indent\indent\indent(a) If $X$ is metrizable, then $X$ is completely metrizable.

\indent\indent\indent(b) If $X$ is metrizable and connected, then $X$ is second countable.
\\

Let $X$ be a topological space $-$then a collection $\sS = \{S\}$ of subsets of $X$ is said to be:

\un{point finite}
\index{point finite} 
if each $x \in X$ belongs to at most finitely many $S \in \sS$;  

\un{neighborhood finite}
\index{neighborhood finite} 
if each $x \in X$ has a neighborhood meeting at most finitely many $S \in \sS$;

\un{discrete}\index{discrete}
 if each $x \in X$ has a neighborhood meeting at most one $S \in \sS$.

A collection which is the union of a countable number of 
$
\begin{cases}
\ \text{point finite}\\[-.1cm]
\ \text{neighborhood finite}\\[-.1cm]
\ \text{discrete}
\end{cases}
$
subcollections is said to be
\[
\begin{cases}
\ \text{\un{$\sigma$-point finite}} \\[-.1cm]
\ \text{\un{$\sigma$-neighborhood finite}} \\[-.1cm]
\ \text{\un{$\sigma$-discrete}}
\end{cases}
.
\]
\index{sigma, $\sigma$-point finite}
\index{sigma, $\sigma$-neighborhood finite}
\index{sigma, $\sigma$-discrete}

\begingroup%%----------------------------------->>
\fontsize{9pt}{11pt}\selectfont
A collection $S = \{\sS\}$ of subsets of $X$ is said to be 
\un{closure preserving}
\index{closure preserving} 
if for every subcollection $\sS_0 \subset \sS$, $\bigcup \overline{S}_0 = \overline{\bigcup S_0}$,  $\overline{S_0}$ the collection $\{\overline{S}: S \in \sS_0\}$.\\
\quadx A collection which is the union of a countable number of closure preserving subcollections is said to be 
\un{$\sigma$-closure preserving}.
\index{sigma, $\sigma$-closure preserving}

Every neighborhood finite collection of subsets of $X$ is closure preserving but the converse is certainly false since any collection of subsets of a discrete space is closure preserving.  A point finite closure preserving closed collection is neighborhood finite.  However, this is not necessarily true if ``closed'' is replaced by ``open'' as can be seen by taking $X = [0,1]$, $\sS = \{]0,1/n[\ : n \in \N\}$.\\
\endgroup%%------------------------------------<<

Let $\sS = \{\sS\}$ be a collection of subsets of $X$.  
The \un{order}
\index{order of a point wrt a set} 
of a point $x \in X$ with respect to $\sS$, written ord($x, \sS)$, is the cardinality of $\{S \in \sS: x \in S\}$.  $\sS$ is of 
\un{finite order}
\index{point of finite order } 
if ord($\sS) = \sup\limits_{x \in X} \text{ord($x, \sS) < \omega$}$.  
The 
\un{star}
\index{star (of a subset)}  
of a subset $Y \subset X$ with respect to $\sS$, written st($Y,\sS$), is the set $\bigcup \{S \in \sS: S \cap Y \neq \emptyset\}$.  $\sS$ is 
\un{star finite}
\index{star finite} 
if $\forall \ S_0 \in \sS: \#\{S \in \sS: S \cap S_0 \neq \emptyset\} < \omega$.

Suppose that $\sU = \{U_i: i \in I\}$ is a covering of $X$ $-$then a covering 
$\sV = \{V_j: j \in J\}$ of $X$ is a 
\un{refinement}
\index{refinement (of a covering)}  
(\un{star refinement})
\index{star refinement} 
of $\sU$ if each $V_j$ (st($V_j,\sV$)) is contained in some $U_i$ and is a 
\un{precise refinement}
\index{precise refinement}  
of $\sU$ if $I = J$ and $V_i \subset U_i$ for every $i$.  
If $\sU$ admits a point finite (open) or a neighborhood finite (open, closed) refinement, 
then $\sU$ admits a precise point finite (open) or neighborhood finite (open, closed) refinement.\\

%%----------------------------------------------------------------------------------------------03
To illustrate the terminology, recall that if $X$ is metrizable, then every open covering of $X$ has an open refinement that is both neighborhood finite and $\sigma$-discrete.

Let $X$ be a completely regular Hausdorff space (CRH space).

\indent\indent $\text{(C) \  $X$}$ is compact
\index{compact} 
iff every open covering of $X$ has a finite (neighborhood finite, point finite) subcovering.

\indent\indent $\text{(P) \  $X$}$ is paracompact
\index{paracompact}  
iff every open covering of $X$ \ has a neighborhood finite open\  (closed) refinement.

\indent\indent $\text{(M) \  $X$}$ is metacompact
\index{metacompact} 
iff every open covering of $X$ has a point finite open refinement.
\\

\label{6.33}
\begingroup%%----------------------------------->>
\fontsize{9pt}{11pt}\selectfont
The following conditions are equivalent to paracompactness.

\indent\indent (P$_1$) \  Every open covering of $X$ has a closure preserving open refinement.

\indent\indent (P$_2$) \  Every open covering of $X$ has a $\sigma$-closure preserving open refinement.

\indent\indent (P$_3$) \  Every open covering of $X$ has a closure preserving closed refinement.

\indent\indent (P$_4$) \  Every open covering of $X$ has a closure preserving refinement.
\\
\endgroup%%------------------------------------<<

\begin{proposition} \ %02
A LCH space $X$ is paracompact iff every open covering of $X$ has a star finite open refinement.
\end{proposition}
[Suppose that $X$ is paracompact.  
Given an open covering $\sU = \{U_i\}$ of $X$, choose a relatively compact open refinement 
$\sV = \{V_j\}$ of $\sU$ such that each $\overline{V}_j$ is contained in some $U_i$ 
$-$then every neighborhood finite open refinement of $\sV$ is necessarily star finite.]\\

\begingroup%%----------------------------------->>
\fontsize{9pt}{11pt}\selectfont
A collection $\sS = \{S\}$ of subsets of a CRH space $X$ is said to be 
\un{directed}
\index{directed 
(collection of sets)} 
if for all $S_1, \ S_2 \in \sS$, there exists $S_3 \in \sS$ such that $S_1 \cup S_2 \subset S_3$.

The following condition is equivalent to metacompactness.

\indent\indent $(\tM)_\text{D}$ \  Every directed open covering of $X$ has a closure preserving closed refinement.

Given an open covering $\sU$ of $X$, denote by $\sU_F$ the collection whose elements are the unions of the finite subcollections of $\sU$ $-$then $\sU_F$ is directed and refines $\sU$ if $\sU$ itself is directed.  
So the above characterization of metacompactness can be recast:

\indent\indent $(\tM)_\text{F}$ \  For every open covering $\sU$ of $X$, $\sU_F$ has a closure preserving closed refinement.

It is therefore clear that a LCH space $X$ is metacompact iff $X$ admits a representation $X = \bigcup\limits_i K_i$, where $\{K_i\}$ is a closure preserving collection of compact subsets of $X$.\\
\endgroup%%------------------------------------<<

A CRH space $X$ is said to be 
\un{subparacompact}
\index{subparacompact} 
if every open covering of $X$ has a $\sigma$-discrete closed refinement.

[Note: \  This definition is partially suggested by the fact that $X$ is paracompact iff every open covering of $X$ has a $\sigma$-discrete open refinement.]\\
%^

%%%%%%%%%%%%%%%%%%%%%%%%%%
%%----------------------------------------------------------------------------------------------04
%%%%%%%%%%%%%%%%%%%%%%%%%%
%^
Suppose that $X$ is subparacompact.  Let $\sU = \{U\}$ be an open covering of $X$ $-$then $\sU$ has a closed refinement $\sA = \bigcup\limits_n \sA_n$, where each $\sA_n$ is discrete.  Every $A \in \sA_n$ is contained in some $U_A \in \sU$.  The collection
\[
\sV_n \ = \ \{U_A - (\cup \sA_n - A): A \in \sA_n\} \  \cup \ \{U - \cup \sA_n: U \in \sU\}
\]
is an open refinement of $\sU$ and $\forall \ x \in X$ $\exists \ n_x$: ord($x,\sV_{n_x}) = 1$.\\

\begingroup%%----------------------------------->>
\fontsize{9pt}{11pt}\selectfont
\label{1.12}
\textbf{\small FACT} 
$X$ is subparacompact iff every open covering of $X$ has a $\sigma$-closure preserving closed refinement.\\
\endgroup%%------------------------------------<<

A CRH space $X$ is said to be 
\un{submetacompact}
\index{submetacompact} 
if for every open covering $\sU$ of $X$ there exists a sequence $\{\sV_n\}$ of open refinements of $\sU$ such that $\forall \ x \in X$ $\exists \ n_x$: ord$(x,\sV_{n_x}) < \omega$.\\

\begingroup%%----------------------------------->>
\fontsize{9pt}{11pt}\selectfont
\textbf{\small FACT} 
$X$ is submetacompact iff every directed open covering of $X$ has a $\sigma$-closure preserving closed refinement.\\
\endgroup%%------------------------------------<<

These properties are connected by the implications:
\\[-.2cm]
%\[
%\begin{tikzcd}%[ sep=small]
%&&&{\text{metacompact}} \ar{rr} &&{} &{\hspace{-1.3cm}\text{submetacompact}}\\
%{\text{compact}} \ar{r} &{\text{paracompact}} \ar{rru}\ar{rrd}\\
%&&&{\text{subparacompact}} \ar{rruu}
%\end{tikzcd}
%\]

\begin{tikzpicture}
\node at (-5,0) {compact};
\node at (-1.75,0) {paracompact};
\node at (2,1.25) {metacompact};
\node at (2,-1.25) {subparacompact};
\node at (6.5,0) {submetacompact};

\draw[->] (-4.15,0) -- (-2.95,0) node[below] {$$};
\draw[->] (-.55,.05) -- (0.8, 1.2) node[above] {$$};
\draw[->] (-.55,-.05) -- (0.6, -1.15) node[above] {$$};
\draw[->] (3.2, 1.25) -- (4.95, .05) node[above] {$$};
\draw[->] (3.4,-1.25) -- (4.95, -.05) node[above] {$$};
\end{tikzpicture}

%^
Each is hereditary with respect to closed subspaces and, apart from compactness, each is hereditary with respect to $F_\sigma$-subpsaces (and all subspaces if this is so of open subspaces).\\

\begingroup%%----------------------------------->>  TTTTHHHHOOOOMMMMAAAASSS PLANK
\fontsize{9pt}{11pt}\selectfont
\index{Thomas Plank}
\textbf{EXAMPLE} \ (\un{The Thomas Plank}) \quadx Let $L_0 = \{(x,0): 0 < x < 1\}$ and for $n \geq 1$, let $L_n = \{(x, 1/n):0 \leq x < 1\}$.  Put 
$X = \ds\bigcup\limits_0^\infty L_n$.  
Topologize $X$ as follows:  For $n \geq 1$, each point of $L_n$ except for $(0,1/n)$ is isolated, 
basic neighborhoods of $(0,1/n)$ being subsets of $L_n$ containing $(0,1/n)$ and having finite complements, 
while for $n = 0$, basic neighborhoods of $(x,0)$ are sets of the form 
$\{(x,0)\} \cup \{(x,1/m): m \geq n\}$ $(n = 1, 2, \ldots)$.  $X$ is a LCH space.  
Moreover, $X$ is metacompact:  
Every open covering of $X$ has an open refinement consisting of one basic neighborhood for each $x \in X$ and any such refinement is point finite since the order of each $x \in X$ with respect to it is at most three.  
But $X$ is not paracompact.  
In fact, $X$ is not even normal: 
$A = \{(0,1/n): n = 1, 2, \ldots\}$ and $B = L_0$ are disjoint closed subsets of $X$ and every neighborhood of $A$ contains all but countably many points of $\ds\bigcup\limits_1^\infty L_n$, while every neighborhood of $B$ contains uncountably many points of $\ds\bigcup\limits_1^\infty L_n$.  
Finally, $X$ is subparacompact.  
This is because $X$ is a countable union of closed paracompact subspaces.\\

\begin{tikzpicture}[scale=6,shift={(10,2)}]
%\begin{tikzpicture}[scale=6]%[scale=0.5,shift={(-5,-3)}]
\node[label={{}}] at (0,0) {$\circ$}; \draw[] (0,0) node[below] {$(0,0)$};
\node[label={{}}] at (1,0) {$\circ$}; \draw[] (1,0) node[below] {$(1,0)$};
\draw[green](0,0) -- (1,0);

\node[label={{}}] at (0,1) {\textbullet}; \draw[] (0,1) node[left] {$(0,1)$};
\node[label={{}}] at (1,1) {\textbullet}; \draw[] (1,1) node[right] {$(1,1)$};
\draw[violet](0, 1) -- (1, 1);
\node[label={{}}] at (0,0.5) {\textbullet}; \draw[] (0,0.5) node[left] {$(0,0.5)$};
\node[label={{}}] at (1,0.5) {\textbullet}; \draw[] (1,0.5) node[right] {$(1,0.5)$};
\draw[violet](0, 0.5) -- (1, 0.5);

\node[label={{}}] at (0,0.33) {\textbullet}; \draw[] (0,0.33) node[left] {$(0,0.33)$};
\node[label={{}}] at (1,0.33) {\textbullet}; \draw[] (1,0.33) node[right] {$(1,0.33)$};
\draw[violet](0, 0.33) -- (1, 0.33);

\node[label={{}}] at (0,0.2455) {\textbullet}; \draw[] (0,0.25) node[left] {$(0,0.25)$};
\node[label={{}}] at (1,0.2455) {\textbullet}; \draw[] (1,0.25) node[right] {$(1,0.25)$};
\draw[violet](0, 0.25) -- (1, 0.25);

%\node[label={{}}] at (0,0.2) {\textbullet}; \draw[] (0,0.2) node[left] {$(0,0.2)$};
%\node[label={{}}] at (1,0.2) {\textbullet}; \draw[] (1,0.2) node[right] {$(1,0.2)$};
\draw[violet](0, 0.2) -- (1, 0.2);

\draw[violet](0, 0.1667) -- (1, 0.1667);
\node[label={{}}] at (0,0.1667) {\textbullet};
\node[label={{}}] at (1,0.1667) {\textbullet};
\draw[violet](0, 0.1429) -- (1, 0.1429);
\node[label={{}}] at (0,0.1429) {\textbullet};
\node[label={{}}] at (1,0.1429) {\textbullet};
\draw[violet](0, 0.1250) -- (1, 0.1250);
\node[label={{}}] at (0,0.1250) {\textbullet};
\node[label={{}}] at (1,0.1250) {\textbullet};
\draw[violet](0, 0.1111) -- (1, 0.1111);
\node[label={{}}] at (0,0.1111) {\textbullet};
\node[label={{}}] at (1,0.1111) {\textbullet};
\draw[violet](0, 0.100) -- (1, 0.1000);
\node[label={{}}] at (0,0.1000) {\textbullet};
\node[label={{}}] at (1,0.1000) {\textbullet};
\draw[violet](0, 0.0909) -- (1, 0.0909);
\node[label={{}}] at (0,0.0909) {\textbullet};
\node[label={{}}] at (1,0.0909) {\textbullet};

\node[label={{}}] at (0.5,0.0636) {\textbf{.}};
\node[label={{}}] at (0.5,0.0455) {\textbf{.}};
\node[label={{}}] at (0.5,0.0273) {\textbf{.}};
\end{tikzpicture}
\endgroup%%------------------------------------<<
\\

%%----------------------------------------------------------------------------------------------05
%%%%%%%%%%%%%%%%%%%%%%%%%%%%%%%%%%%%%%%%%%%%%%%%%
\begingroup%%----------------------------------->> Burke PlankBurke PlankBurke PlankBurke PlankBurke PlankBurke Plank
\fontsize{9pt}{11pt}\selectfont
\index{Burke Plank} 
\textbf{EXAMPLE} \  (\un{The Burke Plank}) \   
Take $X = [0,\Omega^+[ \times [0,\Omega^+[ - \{(0,0)\}$, $\Omega^+$ the cardinal successor of $\Omega$.  For $0 < \alpha < \Omega^+$, put
\[
\begin{cases}
\ H_\alpha = [0,\Omega^+[ \times  \{\alpha\} \\
\ V_\alpha = \{\alpha\} \times  [0,\Omega^+[\ 
\end{cases}
\\[-.2cm].
\]
Topologize $X$ as follows: Isolate all points except those on the vertical or horizontal axis, the basic neighborhoods of
$
\begin{cases}
\ (0,\alpha) \\
\ (\alpha,0)
\end{cases}
$
being the subsets of
$
\begin{cases}
\ H_\alpha\\
\ V_\alpha\
\end{cases}
$
containing
$
\begin{cases}
\ (0,\alpha) \\
\ (\alpha,0)
\end{cases}
$
and having finite complements.  $X$ is a metacompact LCH space.  
But $X$ is not subparacompact.  
To see this, first observe that if $S$ and $T$ are subsets of $X$ such that 
$S \cap H_\alpha$ and 
$T \cap V_\alpha$ are countable for every $\alpha < \Omega^+$, then $X \neq S \cup T$.  
Let $\sU = \{H_\alpha: 0 < \alpha < \Omega^+\} \cup \{V_\alpha: 0 < \alpha < \Omega^+\}$.  
$\sU$ is an open covering of $X$ and the claim is: 
$\sU$ does not have a $\sigma$-discrete closed refinement 
$\sV = \bigcup\limits_n \sV_n$.
To get a contradiction suppose that such a $\sV$ does exist.  
Let $\sS_n$ and $\sT_n$ be the elements of $\sV_n$ which are contained in 
$\{H_\alpha: 0 < \alpha < \Omega^+\}$ and 
$\{V_\alpha: 0 < \alpha < \Omega^+\}$, 
respectively $-$then 
$\sV_n = \sS_n \cup \sT_n$.  
Write
$
\begin{cases}
\ S = \underset{n}{\bigcup} S_n\\
\ T = \underset{n}{\bigcup} T_n
\end{cases}
, \ 
$
where \ 
$
\begin{cases}
\ S_n = \bigcup \sS_n \\
\ T_n = \bigcup \sT_n
\end{cases}
. \ \ 
$
Since the $\sV_n$ are discrete, $S \cap H_\alpha$ and $T \cap V_\alpha$ are countable for every 
$\alpha < \Omega^+$, thus $X \neq S \cup T = \cup \sV$ and so $\sV$ does not cover $X$.

[Note: \  Why does one work with $\Omega^+$ rather than $\Omega$?  Reason: In general, if the weight of 
$X$ is $\leq \Omega$, then $X$ is subparacompact iff $X$ is submetacompact.]\\
\endgroup%%------------------------------------<<

\begingroup%%----------------------------------->>Isbell-Mr\'owka SpaceIsbell-Mr\'owka SpaceIsbell-Mr\'owka SpaceIsbell-Mr\'owka Space
\fontsize{9pt}{11pt}\selectfont
\textbf{EXAMPLE}  \ (\un{Isbell-Mr\'owka Space}) \   
\index{Isbell-Mr\'owka space}
Let $D$ be an infinite set.  
Choose a maximal infinite collection $\sS$ of almost disjoint countably infinite subsets of $D$, 
almost disjoint meaning that 
$\forall \ S_1 \neq S_2 \in \sS$, $\#(S_1 \cap S_2) < \omega$.  
Observe that $\sS$ is uncountable.  Put $\Psi(D) = \sS \cup D$.  
Topologize $\Psi(D)$ as follows: 
Isolate the points of $D$ and take for the basic neighborhoods of a point $S \in \sS$ 
all sets of the form $\{S\} \cup (S - F)$, \mF a finite subset of \mS.  $\Psi(D)$ is a LCH space.  
In addition: $\sS$ is closed and discrete, while $D$ is open and dense.  
Specialize and let $D = \N$ $-$then $X = \Psi(\N)$ is subparacompact, being a Moore space 
(cf. p. \pageref{1.1}),
 but is not metacompact.  In fact, since $\sS$ is uncountable, the open covering 
 $\{\N\} \cup \{\{\sS\} \cup S : S \in \sS\}$ cannot have a point finite open refinement.
\vspi
[Note: \  The Isbell-Mr\'owka space $\Psi(\N)$ depends on $\sS$.  Question: Up to homeomorphism how many distinct $\Psi(\N)$ are there?  Answer: $2^{2^\omega}$.]\\
\vspi
The coproduct of the Burke plank and the Isbell-Mr\'owka space provides an example of a submetacompact $X$ that is neither metacompact nor subparacompact.\\
\vspace{0.25cm}
\endgroup%%------------------------------------<<

\begingroup%%----------------------------------->>van Douwen Linevan Douwen Linevan Douwen Linevan Douwen Linevan Douwen Line
\fontsize{9pt}{11pt}\selectfont
\textbf{EXAMPLE} (\un{The van Douwen Line}) \ 
\index{van Douwen Line}
The object is to equip $X = \R$ with a first countable, separable topology that is finer than the usual topology (hence Hausdorff) 
and under which $X = \R$ is locally compact but not submetacompact.  
Given $x \in \R$, choose a sequence $\{q_n(x)\} \subset \Q$ such that $\abs{x - q_n(x)} < 1/n$.  
Next, let $\{C_\alpha: \alpha < 2^\omega\}$ be an enumeration of the countable subsets $C_\alpha$ of $\R$ with $\#(\overline{C}_\alpha) = 2^\omega$.  
For $\alpha < 2^\omega$, $N = 0, 1, 2, \ldots$, pick inductively a point
\endgroup

\begingroup
\fontsize{9pt}{11pt}\selectfont
\[
x_{\alpha_N} \in  \overline{C}_\alpha - (\Q \cup \{x_{\beta_M}: \beta < \alpha \text{ or } \beta = \alpha \text{ and } M < N\}).
\]
%%----------------------------------------------------------------------------------------------06
Put
\[
%\begin{equation*}
\begin{cases}
\ S_0 \ = \ \{x_{\alpha_0}: \alpha < 2^\omega \}\\
\ S_N \ = \ \{x_{\alpha_N}: \alpha < 2^\omega  \text{ and } C_\alpha \subset S_0\} \quadx (N = 1, 2, \ldots)\\
\end{cases}
%\end{equation*}
\]
and write $S$ in place of $\R - \ds\bigcup\limits_1^\infty S_N$.  
Observe that $\Q \cup S_0 \subset S$ and that the $S_N$ are pairwise disjoint.  
Given $x = x_{\alpha_N} \in \R - S$, choose a sequence 
$\{c_m(x)\} \subset C_\alpha (\subset S_0 \subset S)$ such that $\abs{x - c_m(x)} < 1/m$.  
Topologize $X = \R$ as follows:  Isolate the points of $\Q$ and take for the basic neighborhoods of
$ \begin{cases}
\ x \in S - \Q\\
\ x \in \R - S
\end{cases}
$
the sets
\begin{equation*}
\begin{cases}
\ K_k(x) = \{x\} \cup \{q_n(x): n \geq k\}\\
\ K_k(x) = \{x\} \cup \{c_m(x): m \geq k\} \cup \{q_n(c_m(x)): m \geq k, n \geq m\}\\
\end{cases}
(k = 1, 2, \ldots)
\end{equation*}
This prescription defines a first countable, separable topology on the line that is finer than the usual topology.  
And, since the $K_k$ are compact, it is a locally compact topology.  
However, it is not a submetacompact topology.  
Thus let $U_N = S \cup S_N$ $-$then $U_N$ is open and $\sU = \{U_N\}$ is an open covering of $X$.  
Consider any sequence $\{\sV_M\}$ of open refinements of 
$\sU$.  For $M = 1, 2, \ldots$ and $N = 1, 2, \ldots$, let 
$W_{MN} = \ds\bigcup \{V \in \sV_M: V \cap S_N \neq \emptyset\}$ and form 
$W_0 = S_0 \cap \ds\bigcap\limits_{M,N} W_{MN} =$ 
$S_0 - \ds\bigcup\limits_{M,N} (S_0 - W_{MN})$.  
Since $\#(S_0) = 2^\omega$ and since the 
$S_0 - W_{MN}$ are countable, $W_0$ is nonempty.  
But any $x_0$ in $W_0$ necessarily belongs to infinitely many distinct elements of 
$\sV_M$ $(M = 1, 2, \ldots)$.  
Consequently, the topology is not submetacompact.\\
\endgroup%%------------------------------------<<

\begingroup%%----------------------------------->>
\fontsize{9pt}{11pt}\selectfont
\textbf{\small JONES' LEMMA} \quadx
If a Hausdorff space $X$ contains a dense set $D$ and a closed discrete subspace $S$ with 
$\#(S) \geq 2^{\#(D)}$, then $X$ is not normal.\\
\endgroup%%------------------------------------<<

\begingroup%%----------------------------------->>
\fontsize{9pt}{11pt}\selectfont
Application:  The van Douwen line is not normal.

[In fact, each $S_N$ is closed and discrete with $\#(S_N) = 2^\omega$.]\\
\endgroup%%------------------------------------<<

Let $X$ be a LCH space.  Under what conditions is it true that $X$ metacompact $\implies$ $X$ paracompact?  For example, is it true that if $X$ is normal and metacompact , then $X$ is paracompact?  
This is an open question.  
There are no known counterexamples in ZFC or under any additional set theoretic assumptions.  
Two positive results have been obtained.

\indent\indent (1) 
(Daniels\footnote[2]{\vspace{.11 cm} \textit{Canad. J. Math.} \textbf{35} (1983), 807-823; 
see also \textit{Topology Appl.} 28 (1988), 113-125.}) 
A normal LCH space $X$ is paracompact provided that it is 
\un{boundedly} \un{metacompact},
\index{boundedly metacompact} 
i.e., every open covering of 
$X$ has an open refinement of finite order.

\label{1.11}
\indent\indent (2) 
(Gruenhage\footnote[3]{\vspace{.11 cm}\textit{Topology Proc.} \textbf{4} (1979), 393-405.}) 
A normal LCH space $X$ is paracompact provided that it is locally connected and submetacompact.\\
%%----------------------------------------------------------------------------------------------07

\begingroup%%----------------------------------->>
\fontsize{9pt}{11pt}\selectfont
Suppose that $X$ is normal and metacompact $-$then on general grounds all that one can say is this.  Consider any open covering $\sU$ of $X$: By metacompactness, $\sU$ has a point finite open refinement $\sV$ which, by normality, has a precise open refinement $\sW$ with the property that $\overline{\sW}$ is a precise closed refinement of $\sV$.\\

\textbf{\small FACT} \ 
Let $X$ be a CRH space.  Suppose that $X$ is submetacompact $-$then $X$ is normal iff every open covering of $X$ has a  precise closed refinement.\\
\endgroup%%------------------------------------<<

A Hausdorff space $X$ is said to be 
\un{perfect}
\index{perfect (space)} 
if every closed subset of $X$ is a $G_\delta$.  The Isbell-Mr\'owka space $\Psi(\N)$ is perfect; however, it is not normal (cf. \pageref{1.2}).

A Hausdorff space $X$ is said to be 
\un{perfectly normal}
\index{perfectly normal (space)} 
if it is perfect and normal.  
The ordinal space $[0,\Omega]$, while normal, is not perfectly normal since the point $\{\Omega\}$ is not a $G_\delta$.  
On the other hand, $X$ metrizable $\implies$ $X$ perfectly normal.  Every perfectly normal LCH space $X$ is first countable.

[Note: \  The assumption of perfect normality can be used to upgrade the strength of a covering property.\\
\indent\indent (1) \ 
(Arhangel'ski\u i\footnote[6]{\vspace{.11 cm}\textit{Soviet Math. Dokl.} \textbf{13} (1972), 517-520.}) 
Let $X$ be a LCH space.  If $X$ is perfectly normal and metacompact, then $X$ is paracompact.\\
\indent\indent (2) \ 
(Bennett-Lutzer\footnote[1]{\vspace{.11 cm}\textit{General Topology Appl.} \textbf{2} (1972), 49-54.}) 
Let $X$ be a LCH space.  If $X$ is perfectly normal and submetacompact, then $X$ is subparacompact.]

A CRH space $X$ is said to be 
\un{countably paracompact}
\index{countably paracompact} 
if every countable open covering of $X$ has a neighborhood finite open refinement.  
The ordinal space $[0,\Omega[$ is countably paracompact (being countably compact) and normal, whereas the ordinal space 
$[0,\Omega] \times [0,\Omega[$ is countably paracompact (being compact $\times$ countably compact $\equiv$ countably compact) but not normal.  
On the other hand, $X$ perfectly normal $\implies $ $X$ countably paracompact.

To recapitulate:\\
\[
\begin{tikzpicture}
\node at (2,0) {metrizable};
\node at (6,0) {normal};
\node at (11,0) {countably paracompact};

\node at (6,1.25) {paracompact};
\node at (6,-1.25) {perfectly normal};

\draw[->] (6,1.05) -- (6,.2) node[below] {$$};
\draw[->] (6,-1.05) -- (6,-.2) node[above] {$$};

\draw[->] (3,.1) -- (4.8, 1.25) node[above] {$$};
\draw[->] (3,-.1) -- (4.5, -1.25) node[above] {$$};

\draw[->] (7.2,1.25) -- (8.9, .1) node[above] {$$};
\draw[->] (7.5,-1.25) -- (8.9, -.1) node[above] {$$};
\end{tikzpicture}
\]
\\[-.5cm]

\begingroup%%----------------------------------->>
\fontsize{9pt}{11pt}\selectfont
\textbf{\small FACT} \ 
Suppose that $X$ is normal $-$then $X$ is countably paracompact iff every countable open covering of $X$ has a 
$\sigma$-discrete closed refinement.\\

%%----------------------------------------------------------------------------------------------08
So: In the presence of normality, $X$ subparacompact $\implies$ $X$ countably paracompact.  This implication is strict since the ordinal space $[0,\Omega[$ is normal and countably paracompact; however, it is not even submetacompact 
(cf. p \pageref{1.3}).  On the other hand:
(i) The ordinal space $[0,\Omega] \times [0,\Omega[$ is nonnormal and countably paracompact but not subparacompact; 
(ii) The Isbell-Mr\'owka space $\Psi(\N)$ is nonnormal and subparacompact but not countably paracompact 
(cf. p. \pageref{1.4}).
\vspi
[Note: \  To verify that $ X = [0,\Omega] \times [0,\Omega[$ is not subparacompact, let 
$A \ =  \ \{\Omega,\alpha): \alpha < \Omega\}$ and 
$B \ =  \ \{(\alpha, \alpha): \alpha < \Omega\}$ 
$-$then $A$ and $B$ are disjoint closed subsets of $X$.  
Therefore $X = U \cup V$, where $U = X - A$ and $V = X - B$.  
Since the open covering $\{U,V\}$ has no $\sigma$-discrete closed refinement, $X$ is not subparacompact.]\\
\endgroup%%------------------------------------<<

Is every normal LCH space countably paracompact?  This question is a reinforcement of the ``Dowker problem''.  Dropping the supposition of local compactness, a 
\un{Dowker space}
\index{Dowker space} 
is by definition a normal Hausdorff space which fails to be countably paracompact or, equivalently, whose product with $[0,1]$ is not normal.  
Do such spaces exist?  
The answer is ``yes'', the first such example within ZFC being a construction due to 
M.E. Rudin\footnote[2]{\textit{Fund. Math.} \textbf{73} (1971) 179-186; 
see also Balogh, \textit{Proc. Amer. Math. Soc.} \textbf{124} (1996), 2555-2560.}.  
Her example is not locally compact and only by imposing assumptions beyond ZFC has it been possible to produced locally compact examples.\\

\begingroup%%----------------------------------->>
\fontsize{9pt}{11pt}\selectfont
The ordinal space $[0,\Omega] \times [0,\Omega[$ is neither first countable nor separable.  Can one construct an example of a nonnormal countably paracompact LCH space with both of these properties''  The answer is ``yes''.
\vspi
Let $S$ and $T$ be subsets of $\N$.  Write $S \leq T$ if $\#(S -T) < \omega$; write $S < T$ if $S \leq T$ and $\#(T - S) = \omega$.\\
\endgroup%%------------------------------------<<

\begingroup%%----------------------------------->>
\fontsize{9pt}{11pt}\selectfont
\textbf{\small LEMMA \ (Hausdorff)} \ 
There exist collections
$ \begin{cases}
\ \sS^+ = \{S_\alpha^+ : \alpha < \Omega\} \\
\ \sS^- = \{S_\alpha^- : \alpha < \Omega\}
\end{cases}
$
of subsets of $\N$ with the following properties:

\indent\indent 
(1) $\forall \ \alpha: \#(\N - (S_\alpha^+ \cup S_\alpha^-)) = \omega$.\\
\indent\indent 
(2) $\forall \ \alpha$, $\forall \ \beta$: $\beta < \alpha \implies S_\beta^+ < S_\alpha^+$ and $S_\beta^- < S_\alpha^-$.\\
\indent\indent 
(3) $\forall \ \alpha$: $\#(S_\alpha^+ \cap S_\alpha^-) <  \omega$.\\
\indent\indent 
(4) $\forall \ \alpha$, $\forall \ n \in \N$: $\#\{\beta: \beta < \alpha \ \& \ S_\alpha^+ \cap S_\beta^- \subset F_n\} < \omega$ $(F_n = \{1, \ldots, n\})$.
\vspi
There is then no $H \subset \N$ such that $\forall \ \alpha$: $S_\alpha^+ \leq H$ and $S_\alpha^- \leq \N - H$.
\vspi
[We shall establish the existence of $\sS^+$ and $\sS^-$ by constructing their elements via induction on $\alpha$.  
Start by setting $S_0^+ = \emptyset$ and $S_0^- = \emptyset$.  Given $S_\alpha^+$ and $S_\alpha^-$, 
decompose $\N - (S_\alpha^+ \cup S_\alpha^-)$ into three infinite pairwise disjoint sets $N_\alpha^+$, $N_\alpha^-$, and $N_\alpha$.  
Put\\
\[
\begin{cases}
\ S_{\alpha + 1}^+ = S_\alpha^+ \cup N_\alpha^+\\
\ S_{\alpha + 1}^- = S_\alpha^- \cup N_\alpha^-
\end{cases}
\quadx (\implies \N - (S_{\alpha + 1}^+ \cup S_{\alpha + 1}^-) \supset N_\alpha).
\]
%%----------------------------------------------------------------------------------------------09
Then this definition handles successor ordinals $< \Omega$.  
Suppose now that $0 < \Lambda < \Omega$ is a limit ordinal.  
Choose a strictly increasing sequence 
$\{\alpha_i\} \subset [0,\Omega[$: $\alpha_1 = 0$, 
$\sup \alpha_i = \Lambda$.  
Fix $n_i \in \N$ such that 
$S_{\alpha_i}^+ \cap \ds\bigcup\limits_{j \leq i} S_{\alpha_j}^- \subset F_{n_i}$ 
and write $T_\Lambda^+$ for 
$\ds\bigcup\limits_i (S_{\alpha_i}^+ - F_{n_i})$.  
Note that $\forall \ \alpha < \Lambda$: 
$S_{\alpha}^+ < T_\Lambda^+$ and $\forall \ i$: 
$\#(T_\Lambda^+ \cap S_{\alpha_i}^-) < \omega$.  
If 
$I_i = \{\alpha: \alpha_i \leq \alpha < \alpha_{i+1} \ \& \ T_\Lambda^+ \ds\bigcap S_\alpha^- \subset F_i\}$ and if 
$I = \ds\bigcup\limits_iI_i$, then each $I_i$ is finite and so 
$I \cap [0,\alpha[$ is finite for every $\alpha < \Lambda$.  
Assign to each nonzero $\alpha \in I_i$ the infinite set 
$S_\alpha^ -  - \ds\bigcup \{S_{\alpha_j}^-: \alpha_j < \alpha\}$ 
and denote by $n(\alpha)$ its minimum element in $\N - F_i$.  
Relative to this data, define 
$S_\Lambda^+ = T_\Lambda^+ \cup \{n(\alpha): \alpha \in I (\alpha \neq 0)\}$.  
Then it is not difficult to verify that
\[
\begin{cases}
\ \forall \ \alpha < \Lambda: S_\alpha^+ < S_\Lambda^+ \ \text{and } \forall \ i: \ \#(S_\Lambda^+ \cap S_{\alpha_i}^-) < \omega\\
\ \forall \ n \in \N: \#\{\alpha: \alpha< \Lambda \ \& \ S_\Lambda^+ \cap S_{\alpha}^- \subset F_n\} < \omega
\end{cases}
\\[-.2cm].
\]
As for $S_\Lambda^-$, observe that $(\N - S_\Lambda^+) - \bigcup\limits_{j \leq i} S_{\alpha_j}^-$ is infinite, thus there exists an infinite set $L_\Lambda \subset (\N - S_\Lambda^+)$ such that $L_\Lambda \cap S_{\alpha_i}^-$ is finite for every $i$.  
Defining $S_\Lambda^- = \N - (S_\Lambda^+ \cup L_\Lambda)$, we have\\
\[
\begin{cases}
\ \forall \ \alpha < \Lambda: S_\alpha^- < S_\Lambda^-\\
\ S_\Lambda^+ \cap S_{\Lambda}^- = \emptyset, \ \#(\N - (S_\Lambda^+ \cup S_{\Lambda}^-)) = \omega
\end{cases}
\\[-.2cm],
\]
which completes the induction.  
There remains the assertion of nonseparation.  
To deal with it, assume that there exists an $H \subset \N$ such that $S_{\alpha}^+ - H$ and 
$S_{\alpha}^- \cap H$ are both finite for every $\alpha < \Omega$.  
Choose an $n \in \N$: $W = \{\alpha: S_{\alpha}^- \cap H \subset F_n\}$ is uncountable.  
Fix an $\alpha \in W$ with the property that $W \cap [0,\alpha[$ is infinite.  
If $S_{\alpha}^+ - H \subset F_m$, then 
$\{\beta: \beta < \alpha \ \& \ S_{\alpha}^+ \cap S_{\beta}^- \subset F_{\max (m,n)}\}$ contains 
$W \cap [0,\alpha[$.  Contradiction.]\\
\endgroup%%------------------------------------<<

%van Douwen Spacevan Douwen Spacevan Douwen Spacevan Douwen Spacevan Douwen Space
\index{van Douwen Space}
\begingroup%%----------------------------------->>
\fontsize{9pt}{11pt}\selectfont
\textbf{\small EXAMPLE \ \ (\un{van Douwen Space})} \ 
Let
\[
\begin{cases}
\ X^+ = \{+1\} \times ]0,\Omega[\\
\ X^- = \{-1\} \times ]0,\Omega[
\end{cases}
\]
and put $X = X^+ \cup X^- \cup \N$.  Topologize $X$ as follows: Isolate the points of $\N$ and take for the basic neighborhoods of a point
$ \begin{cases}
\ (+1,\alpha) \in X^+\\
\ (-1,\alpha) \in X^-
\end{cases}
$
all sets of the form\\[.1cm]
\[
\begin{cases}
\ K(+1,\alpha: \beta,F) = \{(+1,\gamma): \beta < \gamma \leq \alpha\} \cup ((S_\alpha^+ - S_\beta^+) - F)\\
\ K(-1,\alpha: \beta,F) = \{(-1,\gamma): \beta < \gamma \leq \alpha\} \cup ((S_\alpha^- - S_\beta^-) - F)
\end{cases}
\hspace{-.25cm},\\
\]
where $\beta < \alpha$ and $F \subset \N$ is finite.  
Since the $K(\pm 1,\alpha: \beta,F)$ are compact, $X$ is a LCH space. 
Obviously, $X$ is first countable and separable; in addition, $X$ is countably paracompact, $X^\pm$ being a copy of $]0,\Omega[$.  Still, $X$ is not normal.
\vspi
[Suppose that the disjoint closed subsets $X^+$ and $X^-$ can be separated by disjoint open sets 
$U^+$ and $U^-$.  Given $\alpha \in ]0,\Omega[$, select an ordinal $f(\alpha) < \alpha$ and a finite subset 
$F(\alpha) \subset \N$ such that $K(\pm,\alpha:f(\alpha),F(\alpha)) \subset U^\pm$.  
Choose $\kappa < \Omega$ and a cofinal 
$\sK \subset [0,\Omega[$ such that $\restr{f}{\sK} = \kappa$ (by ``pressing down'', i.e., Fodor's lemma).  
Put
\[ 
\begin{cases}
\ H^+ = \{S_\kappa^+ \cup (\N \cap U^+)) - S_\kappa^-\\
\ H^- = \{S_\kappa^- \cup (\N \cap U^-)) - S_\kappa^+
\end{cases}
\\[-.2cm].
\]
%%----------------------------------------------------------------------------------------------10
Then $H^+ \cap H^- = \emptyset$.  Let $\alpha < \Omega$ be arbitrary.  
Using the cofinality of $\sK$ and the relation 
$\restr{f}{\sK} = \kappa$, one finds that $S_\alpha^\pm \leq H^\pm$.  Contradiction.]\\
\endgroup%%------------------------------------<<

A CRH space $X$ is said to be 
\un{countably compact}
\index{countably compact} 
if every countable open covering of $X$ has a finite subcovering or, equivalently, if every neighborhood finite collection of nonempty subsets of $X$ is finite.  The ordinal space $[0,\Omega[$ is countably compact but not compact.  The van Douwen space is not countably compact but is countably paracompact.\\
\label{1.10}

\begingroup%%----------------------------------->>
\fontsize{9pt}{11pt}\selectfont
Associated with this ostensibly simple concept are some difficult unsolved problems.  
Sample: Within ZFC, does there exist a first countable, separable, countably compact LCH space $X$ that is not compact?  
This is an open question.  
\label{1.10}
But under CH, e.g., such an $X$ does exist (cf. p. \pageref{1.5}).  
Consider the assertion: Every perfectly normal, countably compact LCH space $X$ is compact.  
While innocent enough, this statement is undecidable in ZFC 
(Ostaszewski\footnote[2]{\textit{J. London Math. Soc.} \textbf{14} (1976), 505-516.}, 
Weiss\footnote[3]{\textit{Canad J. Math.} \textbf{30} (1978), 243-249}).\\
\endgroup%%------------------------------------<<%%------------------------------------<<

\begin{proposition} \ %03
$X$ is countably compact iff every point finite open covering of $X$ has a finite subcovering.
\end{proposition}

[Suppose that $X$ is countably compact.  
Let $\sU$ be a point finite open covering of $X$ $-$then, on general grounds, $\sU$ admits an irreducible subcovering $\sV$.  
This minimal covering must be finite: 
For otherwise there would exist an infinite subset $S \subset X$ such that each $x \in X$ has a neighborhood containing exactly one point of \mS, an impossibility.

Suppose that $X$ is not countably compact $-$then there exists a countably infinite discrete closed subset $D \subset X$, say $D = \{x_n\}$.  
Choose a sequence $\{U_n\}$ of nonempty open sets whose closures are pairwise disjoint such that 
$\forall \ n$: $x_n \in U_n$.  The collection $\{X - D, U_1, U_2, \ldots\}$ is a point finite open covering of $X$ which has no finite subcovering.]\\

A CRH space $X$ is said to be 
\un{pseudocompact}
\index{pseudocompact} 
if every countable open covering of $X$ has a finite subcollection whose closures cover $X$ or, equivalently, 
if every neighborhood finite collection of nonempty open subsets of $X$ is finite.  
The Isbell-Mr\'owka space $\Psi (\N)$ is pseudocompact but not countably compact 
(cf. p. \pageref{1.6}).\\

\begin{proposition} \ 
$X$ is a pseudocompact iff every real valued continuous function on $X$ is bounded.
\end{proposition}

[Suppose that $X$ is not pseudocompact $-$then there exists a countably infinite neighborhood finite collection $\{U_n\}$ of nonempty open subsets of $X$. 
Choose a point $x_n \in U_n$.  
%%----------------------------------------------------------------------------------------------11
Since $X$ is completely regular, there exists a continuous function $f_n:X \ra [0,n]$ such that 
$f_n(x_n) = n$, $\restr{f_n}{X - U_n} = 0$.  
Put $f = \sum\limits_n f_n$: $f$ is continuous and unbounded.]\\

A CRH space $X$ is said to be 
\un{countably metacompact}
\index{countably metacompact} 
if every countable open covering of $X$ has a point finite open refinement.  The ordinal space $[0,\Omega[$ is countably metacompact but not metacompact 
(cf. p. \pageref{1.7}).  
Every perfect $X$ is countably metacompact.

The relative positions of these conditions is shown by:\\
\[
\begin{tikzcd}%[ sep=small]
&{\text{compact}} \ar{d} \ar{r} &{\text{paracompact}} \ar{d} \ar{r}  &{\text{metacompact}} \ar{d}\\
&{\text{countably compact}} \ar{d} \ar{r} &{\text{countably paracompact}} \ar{r}  &{\text{countably metacompact}}\\
&{\text{pseudocompact}}\\
\end{tikzcd}
\]

\begingroup%%----------------------------------->>
\fontsize{9pt}{11pt}\selectfont
\textbf{\small FACT} \ 
$X$ is countably metacompact iff for every countable open covering 
$\sU$ of $X$ there exists a sequence $\{\sV_n\}$ of open refinements of 
$\sU$ such that $\forall \ x \in X$ $\exists \ n_x$: ord($x,\sV_{n_x}) < \omega$.
\vspi
[The point here is to show that the stated condition forces $X$ to be countably metacompact.  
Enumerate the elements of $\sU$: $U_n \ (n = 1, 2, \ldots$).  
Write $W_n$ for the set of all $x \in U_n$ such that 
$\forall \ m \leq n$ $\exists \ V \in \sV_m$: $x \in V$ and 
$V \not\subset \bigcup\limits_{i < n} U_i$.  
Then $\sW = \{W_n\}$ is a point finite open refinement of $\sU = \{U_n\}$.]\\

So: $X$ submetacompact $\implies$ $X$ countably metacompact.  The van Douwen line is not countably metacompact (inspect the argument used to establish nonsubmetacompactness).  The Tychonoff plank is countably metacompact but is neither submetacompact nor countably paracompact (cf. p. \pageref{1.8}).\\
\endgroup%%------------------------------------<<

\begin{proposition} \ 
If $X$ is a pseudocompact and either normal or countably paracompact, then $X$ is countably compact.
\end{proposition}

[Suppose that $X$ is normal.  
If $X$ is not countably compact, then there exists a countably infinite discrete closed subset 
$D \subset X$, say $D = \{x_n\}$. \  
By the Tietze extension theorem, there exists a continuous function $f:X \ra \R$ such that $f(x_n) = n$ $(n = 1, 2, \ldots)$. \  Contradiction.

Suppose that $X$ is countably paracompact.  
If $X$ is not countably compact, then there exists a countable open covering $\{U_n\}$ of $X$ that cannot be reduced to a finite covering.  
Let $\{V_n\}$ be a precise neighborhood finite open refinement of $\{U_n\}$ $-$then there exists a finite subset 
$F \subset \N$ such that $V_n = \emptyset$ iff $n \in F$.  
But $\bigcup\limits_n V_n = X$.  Contradiction.]\\
%%----------------------------------------------------------------------------------------------12

\begingroup%%----------------------------------->>
\fontsize{9pt}{11pt}\selectfont
\label{1.2}
\label{1.4}
\label{1.6}
\label{2.9}
\textbf{\small EXAMPLE}  \ 
The Isbell-Mr\'owka space $\Psi (\N)$ is not countably compact.  However, $\Psi (\N)$ is pseudocompact so, by the above, it is neither normal nor countably paracompact.

[Put $X = \Psi (\N)$ and suppose that $f:X \ra \R$ is continuous but unbounded.  
Since $\forall \ S \in \sS$, $\{\sS\} \cup S$ is compact, $\restr{f}{S}$ is bounded.  
This means that there exists a sequence $\{x_n\}$ of distinct points in $X$ such that
(i) \ $\abs{f(x_n)} \geq n$ and 
(ii) \ $\forall \ S \in \sS$, $\#(\{x_n\} \cap S) < \omega$.  
The maximality of $\sS$ then implies that $\{x_n\} \in \sS$.  Contradiction.]\\

\label{1.8}
\index{The Tychonoff Plank}
\textbf{\small EXAMPLE \  \ (\un{The Tychonoff Plank})} \ 
Let $X = [0,\Omega] \times [0,\omega] - \{(\Omega,\omega)\}$.  $X$ is not countably compact (consider $\{(\Omega,n): 0 \leq n < \omega\})$.  However, $X$ is pseudocompact so, by the above, it is neither normal nor countably paracompact.

[Suppose that $f:X \ra \R$ is continuous $-$then it suffices to show that $f$ extends continuously to 
$\{(\Omega,\omega)\}$.  
Because every real valued continuous function on $[0,\Omega[$ is constant on some tail $[\alpha,\Omega[$, $\forall \ n \leq \omega$, there exists $\alpha_n < \Omega$ and a constant $r_n$ such that $f(\alpha,n) = r_n$ $\forall \ \alpha \geq \alpha_n$.  Put $\alpha_0 = \sup \alpha_n$ $-$then $\alpha_0 < \Omega$.  
One can therefore let $f(\Omega,\omega) = r_\omega$.]\\
\endgroup%%------------------------------------<<

\begin{proposition} \ %06
If $X$ is countably compact and submetacompact, then $X$ is compact.
\end{proposition}

[Let $\sU$ be an open covering of $X$.  
Let $\{\sV_n\}$ be a sequence of open refinements of $\sU$ such that $\forall \ x \in X$ $\exists \ n_x$: $\text{ord}(x,\sV_{n_x}) < \omega)$.  
Write $A_{mn}$ for $\{x:  \text{ord}(x,\sV_n) \leq m\}$ $-$then $A_{mn}$ is a closed subspace of $X$, 
hence is countably compact, and $\sV_n$ is point finite on the $A_{mn}$.  
Proposition 3 therefore implies that $A_{mn}$ can be covered by finitely many elements of $\sV_n$.  
Every $x \in X$ is in some $A_{mn}$, so there is a countable covering of $X$ made up of elements from the sequence 
$\{\sV_n\}$.  This covering has a finite subcovering, thus so does $\sU$.]\\

\begingroup%%----------------------------------->>
\fontsize{9pt}{11pt}\selectfont
\label{1.3}
\label{1.7}
Consequently, the ordinal space $[0,\Omega[$ is not submetacompact.  
It then follows from this that the Tychonoff plank is not submetacompact (since $[0,\Omega[$  sits inside it as a closed subspace).\\
\endgroup%%------------------------------------<<

Let $X$ be a CRH space.  A 
\un{$\pi$-basis}
\index{pi, $\pi$-basis} 
for $X$ is a collection $\sP$ of nonempty open subsets of $X$ such that if $O$ is a nonempty open subset of $X$, then for some $P \in \sP$, $P \subset O$.\\

\textbf{\small LEMMA} \ \ 
Suppose that $X$ is Baire.  
Let $\sU$ be a point finite open covering of $X$ $-$then there exists a $\pi$-basis $\sP$ for $X$ such that 
$\forall \ P \in \sP$ and $\forall \ U \in \sU$, either $P \subset U$ or $P \cap U = \emptyset$.

[For $n = 1, 2, \ldots$, denote by $X_n$ the subset of $X$ consisting of those points that are in at most $n$ elements of $\sU$.  Each $X_n$ is closed and 
$X = \ds\bigcup\limits_n X_n$.  
Let $O$ be a nonempty open subset of $X$.  
Since $O = \ds\bigcup\limits_n  O \cap X_n$, there will be an $n$ such that $O \cap X_n$ has a nonempty interior.  
Let $n(O)$ be the smallest such $n$.  
Let $U_O \subset O \cap X_{n(O)}$ be a nonempty open subset 
%%----------------------------------------------------------------------------------------------13
of $X$.  
Choose $x_O \in U_O$ that belongs to exactly $n(O)$ elements of $\sU$ and write $P$ 
for their intersection with $U_O$ $-$then $\sP = \{P\}$ is a $\pi$-basis for $X$ with the stated properties.]\\

Suppose that $X$ is pseudocompact $-$then $X$ is Baire.  To see this, let $\{O_n\}$ be a decreasing sequence of dense open subsets of $X$.  
Let $U$ be a nonempty open subset of $X$.  
Inductively choose nonempty open sets 
$V_n$: $V_1 = U$ $\&$ $\overline{V}_{n+1} \subset U \cap O_n \cap V_n$.  
By pseudocompactness, 
$\ds\bigcap\limits_n \overline{V}_{n} \neq \emptyset$, hence 
$U \cap \bigl(\ds\bigcap\limits_n O_n \bigr) \neq \emptyset$.\\

\begin{proposition} \ %07
If $X$ is pseudocompact and metacompact, then $X$ is compact.
\end{proposition}

[Let $\sO$ be an open covering of $X$.  
Let $\sU = \{U\}$ be a point finite open refinement of $\sO$ with the property that
 $\overline{\sU} = \{\overline{U}\}$ refines $\sO$.  
 Use the lemma to determine a $\pi$-basis $\sP$ for $X$ per $\sU$.  
Fix $P_1 \in \sP$.  
Consider $\{U \in \sU: U \cap P_1 \neq \emptyset\}$.  
Since $U \cap P_1 \neq \emptyset$ $\implies$ 
$P_1 \subset U$ and since $\sU$ is point finite, it is clear that this is a finite set.  
If $X = \overline{\text{st}(P_1,\sU)}$, then finitely many elements of $\sO$ cover $X$ and we are done.  Otherwise, proceed inductively and, using the fact that $\sP$ is a $\pi$-basis for $X$, given 
$n \in \N$ choose $P_{n+1} \in \sP$ such that
\[
P_{n+1} \subset X - \bigcup\limits_{m \leq n} \overline{\text{st}(P_m,\sU)}.
\]
We claim that the process terminates, from which the result.  
Suppose the opposite $-$then, due to the pseudocompactness of $X$, $\{P_n\}$ cannot be neighborhood finite.  
Therefore there exists $x \in U_x \in \sU$ with $U_x \cap P_n \neq \emptyset$ for infinitely many $n$, contrary to construction.]\\

\begingroup%%----------------------------------->>
\fontsize{9pt}{11pt}\selectfont
One cannot replace ``metacompact'' by ``submetacompact'' in the preceding result: 
The Isbell-Mr\'owka space $\Psi(\N)$ is pseudocompact and submetacompact but not compact.  
However, the argument does go through under the weaker condition: Every open covering of $X$ has a $\sigma$-point finite open refinement.\\
\endgroup%%------------------------------------<<

\begin{proposition} \ 
If $X$ is normal and countably metacompact, then $X$ is countably paracompact.\\
\end{proposition}

One can check:\\
\indent\indent (CP)
\label{2.6}
\ $X$ is countably paracompact iff for every decreasing sequence $\{A_n\}$ of closed sets such that 
$\ds\bigcap\limits_n A_n = \emptyset$, there exists a decreasing sequence $\{U_n\}$ of open sets with $A_n \subset U_n$ for every $n$ and such that 
$\ds\bigcap\limits_n \overline{U}_n = \emptyset$.\\
\indent\indent (CM) \ $X$ is countably metacompact iff for every decreasing sequence $\{A_n\}$ of closed sets such that 
$\ds\bigcap\limits_n A_n = \emptyset$, there exists a decreasing sequence $\{U_n\}$ of open sets with 
$A_n \subset U_n$ for every $n$ and such that 
$\ds\bigcap\limits_n U_n = \emptyset$.

%%----------------------------------------------------------------------------------------------14
It remains only to note that for normal $X$, CP $\iff$ CM.\\

\begingroup%%----------------------------------->>
\fontsize{9pt}{11pt}\selectfont
If $X$ is the Tychonoff plank, then $X = Y \cup Z$, where 
$Y = \ds\bigcup\limits_{n < \omega} [0,\Omega] \times \{n\}$ and $Z = [0,\Omega[ \times \{\omega\}$.  
Since $Y$ is an open paracompact subspace of $X$ and $Z$ is a closed countably compact subspace of $X$, 
it is clear that $X$ is countably metacompact.  
Because $X$ is not countably paracompact, 
Proposition 8 allows one to infer once again that $X$ is not normal (cf. Proposition 5).\\
\endgroup%%------------------------------------<<

\label{6.36} %dmc mnft lower
A Hausdorff space $X$ is said to be 
\un{collectionwise normal}
\index{collectionwise normal} 
if for every discrete collection $\{A_i: i \in I\}$ of closed subsets of $X$ there exists a pairwise disjoint collection 
$\{U_i: i \in I\}$ of open subsets of $X$ such that 
$\forall \ i \in I$: $A_i \subset U_i$.

Of course, $X$ collectionwise normal $\implies$ $X$ normal.  
On the other hand, $X$ normal and countably compact $\implies$ $X$ collectionwise normal.  
So, the ordinal space $[0,\Omega[$ is collectionwise normal.  
However, it is not perfectly normal since the set of all limit ordinals 
$\alpha < \Omega$, while closed, is not a $G_\delta$.  
Rudin's Dowker space is collectionwise normal.\\

\textbf{\small LEMMA} \ \ 
Suppose that $X$ is collectionwise normal.  
Let $\{A_i: i \in I\}$ be a discrete collection of closed subsets of $X$ $-$then there exists a discrete collection 
$\{O_i: i \in I\}$ of opens subsets of $X$ such that $\forall \ i \in I$: $A_i \subset O_i$.

[Let $\{U_i: i \in I\}$ be a pairwise disjoint collection of open subsets of $X$ such that 
$\forall \ i \in I$: $A_i \subset U_i$.  
Choose an open set $U$ subject to 
$\ds\bigcup\limits_i A_i \subset U \subset \overline{U} \subset \bigcup\limits_i U_i$ and then put $O_i = U_i \cap U$.]\\

Suppose that $X$ is normal.  Let $\{A_n\}$ be a countable discrete collection of closed subsets of $X$ 
$-$then there exists a countable pairwise disjoint collection $\{U_n\}$ of open subsets of $X$ such that 
$\forall \ n$: $A_n \subset U_n$.  In fact, given $n \in \N$, choose a pair $(O_n,P_n)$ of disjoint open subsets of $X$ such that 
$O_n \supset A_n$, $P_n \supset \bigcup\limits_{m \ne n} A_m$ and then put 
$U_n = O_n \cap \ds\bigcap\limits_{m < n} P_m$.\\

\begin{proposition} \ %09
If $X$ is paracompact, then $X$ is collectionwise normal.
\end{proposition}

[Let \ $\{A_i: i \in I\}$ be a discrete collection of closed subsets of $X$.  \ 
Put \ $O_i = X - \ds\bigcup\limits_{j \neq i} A_j$ $-$then the collection $\{O_i: i \in I\}$ is an open covering of $X$, hence in view of the paracompactness of $X$, has a precise neighborhood finite closed refinement $\{C_i: i \in I\}$.  
If $U_i = X -  \ds\bigcup\limits_{j \neq i} C_j$, then 
$\{U_i: i \in I\}$ is a pairwise disjoint collection of open subsets of $X$ such that 
$\forall \ i \in I$: $A_i \subset U_i$.  Therefore $X$ is collectionwise normal.]\\

\begin{proposition} \ %10
If $X$ is collectionwise normal and metacompact, then $X$ is paracompact.
\end{proposition}

%%----------------------------------------------------------------------------------------------15
[It is enough to prove that a given point finite open covering $\sO = \{O\}$ of $X$ has a $\sigma$-discrete open refinement 
$\sU = \ds\bigcup\limits_n \sU_n$.  
Put $A_n = \{x: \text{ord}(x,\sO) \leq n\}$ $-$then $A_n$ is a closed subspace of $X$ and 
$X = \ds\bigcup\limits_n A_n$. 
Assign to each $x \in X$ the open set 
$O_x = \bigcap \{O \in \sO: x \in O\}$.  
Using the $O_x$, we shall construct the $\sU_n$ by induction.  
To start off, observe that $\{O_x \cap A_1: x \in A_1\}$ is a discrete collection of closed subsets of $X$ covering $A_1$.  
So, by collectionwise normality, there exists a discrete collection $\sU_1$ of open subsets of $X$ covering $A_1$ such that each element of $\sU_1$ is contained in some element of $\sO$.  
Proceeding, suppose that 
$\ds\bigcup\limits_{m = 1}^n \sU_m$ is a covering of $A_n$ by open subsets of $X$, 
each of which is contained in some element of $\sO$, with $\sU_m$ discrete.  Let
$U_n = \ds\bigcup  \{U: U \in \sU_m, 1 \leq m \leq n\}$ $-$then 
$U_n \supset A_n$ and $\{O_x \cap (A_{n+1} - U_n):x \in A_{n+1} - U_n\}$ is a discrete collection of closed subsets of $X$ covering $A_{n+1} - U_n$.
Once again,  by collectionwise normality, there exists a discrete collection $\sU_{n+1}$ of open subsets of $X$ covering $A_{n+1} - U_n$ such that each element of $\sU_{n+1}$ is contained in some element of $\sO$.  
And $A_{n+1} \subset \ds\bigcup\limits_{m = 1}^{n+1} \sU_m$ .]\\

\label{6.37} %dmc mnft
\begingroup%%----------------------------------->>
\fontsize{9pt}{11pt}\selectfont
Trifling modifications in the preceding argument allow one to replace ``metacompact'' by ``submetacompact'' and still arrive at the same conclusion.

Kemoto \footnote[2]{\textit{Fund. Math.} \textbf{132} (1989), 163-169.}
 has shown by very different methods that if a normal LCH space $X$ is submetacompact, then $X$ is subparacompact.  Example: The Burke plank is not normal.\\
\endgroup%%------------------------------------<<

Let $X$ be a LCH  space.  Does the chart
\[
\begin{tikzpicture}
\node at (0,0) {paracompact};
\node at (4,0) {collectionwise normal};
\node at (8,0) {normal};
\node at (4,-1.25) {pefectly normal};

\draw[->] (1.2,0) -- (2.05,0) node[below] {$$};
\draw[->] (6, 0) -- (7.15, 0) node[above] {$$};
\draw[->] (5.45, -1.25) -- (7.15, -.20) node[above] {$$};

\end{tikzpicture}
\]
admit any additional arrows?  We do know that there exists a paracompact $X$ that is not perfectly normal and a collectionwise normal $X$ that is not paracompact.
\label{1.16} %dmcxxx why arent these lining up perfectly? ie without the space

\indent\indent (Q$_a$) \ Is every normal LCH space $X$ collectionwise normal?

[There are counterexamples under MA + $\neg$ CH (cf. \pageref{1.9}).  Consistency has been established modulo the consistency of the existence of a supercompact cardinal.]

\indent\indent (Q$_b$) \ Is every perfectly normal LCH space $X$ collectionwise normal?

[This is undecidable in ZFC.]\\
\indent\indent (Q$_c$) \ Is every perfectly normal LCH space $X$ paracompact?

%%----------------------------------------------------------------------------------------------16
[The Kunen line under CH and the rational sequence topology over a CUE-set under MA +$\neg$ CH are counterexamples.  
However, under ZFC alone, the issue has not been resolved.]

These questions (and many others) are discussed by 
Watson\footnote[3]{In: \textit{Open Problems in Topology}, J. van Mill and G. Reed (ed.), North Holland (1990), 37-76.}.
\\[-.2cm]

\begingroup%%----------------------------------->>
\fontsize{9pt}{11pt}\selectfont
The construction of topologies by transfinite recursion is an important technique that can be used to produce a variety of illuminating examples.\\
\endgroup%%------------------------------------<<

\label{1.15}
\label{19.8}
\label{19.34}
\index{Kunen line} %Kunen lineKunen lineKunen lineKunen lineKunen lineKunen lineKunen lineKunen line
\begingroup%%----------------------------------->>
\fontsize{9pt}{11pt}\selectfont
\textbf{\small EXAMPLE \  [Assume CH] (\un{The Kunen Line})} \quadx
The object is to equip $X = \R$ with a first countable, separable topology that is finer than the usual topology (hence Hausdorff) 
and under which $X = \R$  is locally compact and perfectly normal but not Lindel\"of, hence not paracompact (since paracompact + separable = Lindel\"of).  
It will then turn out that the resulting topology is even hereditarily separable and collectionwise normal.
\vspi
Let $\{x_\alpha: \alpha < \Omega\}$ be an enumeration of $\R$ and put $X_\alpha = \{x_\beta: \beta < \alpha\}$, so $X_\Omega = \R$.  Let $\{C_\alpha: \alpha < \Omega\}$ be an enumeration of the countable subsets of $\R$ such that $\forall \ \alpha$: $C_\alpha \subset X_\alpha$.  We shall now construct by induction on $\alpha \leq \Omega$ a collection $\{\tau_\alpha: \alpha \leq \Omega\}$, where $\tau_\alpha$ is a topology on $X_\alpha$ (with closure operator cl$_\alpha$) subject to:
\\
\indent\indent (a) \ $\forall \ \alpha$: $\tau_\alpha$ is a first countable, zero dimensional, locally compact topology on $X_\alpha$ that is finer than the usual topology on $X_\alpha$ (as a subspace of $\R$) and, if $\alpha < \Omega$, is metrizable.
\\
\indent\indent (b) \ $\forall \ \beta < \alpha$: $(X_\beta,\tau_\beta)$ is an open subspace of $(X_\alpha,\tau_\alpha)$. 
\\
\indent\indent (c) \ $\forall \ \gamma \leq \beta < \alpha$:  
If $x_\beta \in \text{cl}_\R(C_\gamma)$, then $x_\beta \in \text{cl}_\alpha(C_\gamma)$.
\vspi
First, take $\tau_\alpha$ discrete if $\alpha \leq \omega$.  
Assume next that $\omega < \alpha \leq \Omega$.  If $\alpha$ is a limit ordinal, take for $\tau_\alpha$ the topology on $X_\alpha$ generated by 
$\ds\bigcup\limits_{\beta < \alpha} \tau_\beta$.  
If $\alpha$ is a successor ordinal, say $\alpha = \beta + 1$, then the problem is to define 
$\tau_\alpha$ on $X_\alpha = X_\beta \cup \{x_\beta\}$ and for that we distinguish two cases.
\\
\indent\indent ($*$) \ If there is no $\gamma \leq \beta$ such that $x_\beta \in \text{cl}_\R(C_\gamma)$, isolate $x_\beta$ and take for $\tau_\alpha$ the topology generated by $\tau_\beta$ and $\{x_\beta\}$.
\\
\indent\indent $\neg(*)$ \ Let $\{\gamma_n\}$ enumerate 
$\{\gamma \leq \beta: x_\beta \in \text{cl}_\R(C_\gamma)\}$, each 
$\gamma$ being listed $\omega$ times.  
Put $I_n = ]x_\beta - 1/n,x_\beta + 1/n[$ and pick a sequence $\{y_n\}$ of distinct points 
$y_n \in C_{\gamma_n} \cap I_n$.  
Choose a discrete collection 
$\{K_{n,\beta}\}$ of $\tau_\beta$-clopen compact sets 
$K_{n,\beta}$: $y_n \in  K_{n,\beta} \subset I_n$.  
To complete the induction, take for $\tau_\alpha$ the topology generated by 
$\tau_\beta$ and the sets $\{x_\beta\} \cup \ds\bigcup\limits_{m \geq n} K_{m,\beta}$ $(n = 1, 2, \ldots)$.
\vspi
It follows that $\R$ or still, 
$X_\Omega = \ds\bigcup\limits_{\alpha < \Omega} X_\alpha$ is a first countable, LCH space under $\tau_\Omega$.  
Because each $X_\alpha$ is $\tau_\Omega$-open, $X_\Omega$ is not Lindel\"of.  
Every $x \in X_\Omega$ has a countable clopen neighborhood.
\vspi
Claim: Let $S \subset \R$ $-$then $\#(\text{cl}_\R(S) - \text{cl}_\Omega(S)) \leq \omega$.
\vspi
%%----------------------------------------------------------------------------------------------17
[Fix a countable subset $C \subset S$ such that $\text{cl}_\R(C) =\text{cl}_\R(S)$.  \ 
Write 
$C = C_{\alpha_{0}}$ (some $\alpha_0 < \Omega$).  If $\alpha > \alpha_0$ and if $x_\alpha \in \text{cl}_\R(C)$, 
then $x_\alpha \in \text{cl}_\Omega(C)$.  
Therefore $\text{cl}_\R(S) - \text{cl}_\Omega(S) \subset \{x_\alpha: \alpha \leq \alpha_0\}$.]
\vspi
The fact that $X_\Omega$ is hereditarily separable is thus immediate.  
To establish perfect normality, suppose that 
$A \subset X_\Omega$ is closed $-$then it is a question of finding a sequence $\{U_n\} \subset \tau_\Omega$ such that 
$A = \ds\bigcap\limits_n U_n $ 
$= \ds\bigcap\limits_n \text{cl}_\Omega(U_n)$.  
Since $\R$ is perfectly normal, there exists a sequence $\{O_n\}$ of $\R$-open sets such that 
$\text{cl}_\R(A) = \ds\bigcap\limits_n O_n$ 
$= \ds\bigcap\limits_n \text{cl}_\R(O_n)$.
From the claim, $\text{cl}_\R(A) - A$ can be enumerated: $\{a_n\}$.  
Each $a_n \in X_\Omega - A$, so 
$\exists \ K_n \in \tau_\Omega$: 
$a_n \in K_n \subset X_\Omega - A$, $K_n$ clopen.  
Bearing in mind that $\tau_\Omega$ is finer than the usual topology on $\R$, we then have
\[
A  \ = \  \bigcap\limits_n O_n \ \cap \  \bigcap\limits_n (X_\Omega - K_n) \ = \ 
\bigcap\limits_n \text{cl}_\Omega(O_n)  \ \cap \  \bigcap\limits_n  (X_\Omega - K_n).
\]
The final point is collectionwise normality.  
But as CH is in force, Jones' lemma implies that $X_\Omega$, being separable and normal, has no uncountable closed discrete subspaces.
\vspi
\label{1.5}
[Note: \  $X_\Omega$ is not metacompact (cf. Proposition 10).  
However, $X_\Omega$ is countably paracompact (being perfectly normal).]\\
\vspi
Retaining the assumption CH and working with\\
\[
\begin{cases}
\ X_\Omega = \N \cup (\{0\} \times [0,\Omega[) \\
\ X_\alpha = \N \cup \{(0,\beta): \beta < \alpha\}
\end{cases}
\\[-.2cm],
\]
one can employ the foregoing methods and construct an example of a first countable, separable, countably compact, noncompact LCH space 
(cf. p. \pageref{1.10}).  
Recursive techniques can also be used in conjunction with set theoretic hypotheses other than CH to manufacture the same type of example.\\
\endgroup%%------------------------------------<<

A CRH space $X$ is said to be a 
\un{Moore Space}
\index{Moore space} 
if it admits a development.

[Note: \  A 
\un{development}
\index{development} 
for $X$ is a sequence $\{\sU_n\}$ of open coverings of $X$ such that $\forall \ x \in X$: $\{\st(x,\sU_n)\}$ is a neighborhood basis at $x$.]

Every Moore space is first countable and perfect.  
Any first countable $X$ that is expressible as a countable union of closed discrete subspaces $X_n$ is Moore, so, 
\label{1.1}e.g., the Isbell-Mr\'owka space $\Psi(\N)$ is Moore.\\

\begingroup%%----------------------------------->>
\fontsize{9pt}{11pt}\selectfont
\textbf{\small FACT} \ 
Suppose that $X$ is a Moore space $-$then $X$ is subparacompact.

[Let $\sO = \{O_i: i \in I\}$ be an open covering of $X$ 
$-$then the claim is that $\sO$ has a $\sigma$-discrete closed refinement.  
Fix a development $\{\sU_n\}$ for $X$.  Equip $I$ with a well ordering $<$ and put
\[
A_{i,n} \ = \ X - \left( \text{st}(X - O_i, \sU_n) \ \cup \  \bigcup\limits_{j < i} O_j\right) \ \subset \  O_i.
\]
Each $A_{i,n}$ is closed and their totality $\sA$ covers $X$.  
Denote by $\sA_n$ the collection $\{A_{i,n}:i \in I\}$ $-$then 
$\sA_n$ is discrete, so 
$\sA = \ds\bigcup\limits_n \sA_n$ is a $\sigma$-discrete closed refinement of $\sO$.]\\
\endgroup%%------------------------------------<<

%%----------------------------------------------------------------------------------------------18
\label{6.34}
\label{6.35}
The metrization theorem of Bing says:  $X$ is metrizable iff $X$ is a collectionwise normal Moore space.  
Equivalently:  $X$ is metrizable iff $X$ is a paracompact Moore space (cf. Proposition 9).\\ 

\begingroup%%----------------------------------->>
\fontsize{9pt}{11pt}\selectfont

The Kunen line is not a Moore space.  For if it were, then, being collectionwise normal, it would be metrizable, hence paracompact, which it is not.  Variant: The Kunen line is not submetacompact, therefore is not subparacompact (cf. the remark following the proof of Proposition 10), proving once again that it is not a Moore space.\\
\endgroup%%------------------------------------<<

Let $X$ be a LCH space.  If $X$ is locally connected, normal, and Moore, then $X$ is metrizable (Reed-Zenor).  
Proof: 
(1) \ $X$ Moore $\implies$ $X$ subparacompact;
(2) \ $X$ locally connected, normal, and subparacompact (hence submetacompact) $\implies$ $X$ paracompact 
(via the result of Gruenhage mentioned on p. \pageref{1.11}).  
Now cite Bing.\\

\begingroup%%----------------------------------->>
\fontsize{9pt}{11pt}\selectfont
Question: \quadx Is every locally compact normal Moore space metrizable?  
It turns out this question is undecidable in ZFC.

\indent\indent (1) \ Under V = L, every locally compact normal Moore space is metrizable.
\vspi
[Watson\footnote[2]{\textit{Canad. J. Math.} \textbf{34} (1982), 1091-1096.}
proved that under V = L, every normal submetacompact LCH space $X$ is paracompact.  
This leads at once to the result.]

\label{1.9}
\indent\indent (2) \ Under MA $+ \neg$ CH, there exist locally compact normal Moore spaces that are not metrizable.
\vspi
[Many examples are known that illustrate this phenomenon.  
A particularly simple case in point is that of the rational sequence topology over a CUE-set.  
By definition, a 
\un{CUE-set}
\index{CUE-set} 
\mS is an uncountable subset of $\R$ with the property that 
$\forall \ T \ \subset S$, there exists a sequence $\{U_n\}$ of open subsets of $\R$ such that 
$T = S \ \cap \  \bigl(\ds\bigcap\limits_n U_n\bigr)$, i.e., \mT is a relative $G_\delta$.
Assuming MA $+ \neg$ CH, it can be shown that every uncountable subset of $\R$ having cardinality 
$< 2^\omega$ is a CUE-set.  
This said, let $S$ be any uncountable subset of the irrationals of cardinality 
$< 2^\omega$.  
Put $X = (\Q \ \times \Q) \cup (S \times \{0\})$.  
Topologize $X$ as follows: Isolate the points of 
$\Q \ \times \Q$ and take for the basic neighborhoods of $(s,0)$ $(s \in S)$ the sets 
$\{(s,0)\} \cup \{(s_m,1/m): m \geq n\}$ $(n = 1, 2, \ldots)$, 
where $\{s_n\}$ is a fixed sequence of rationals converging to s in the usual sense.
$X$ is a separable LCH space.  
It is clear that $X$ is Moore but not metrizable, hence
(i) $X$ is perfect but not collectionwise normal and 
(ii) $X$ is subparacompact but not metacompact (since separable + metacompact $\implies$ Lindel\"of $\implies$ paracompact).  
Nevertheless, $X$ is normal.  
Indeed, given $T \subset S$, it suffices to produce disjoint open sets 
$U,V \subset X$: $U \supset T$ and $V \supset S - T$.  
Using the fact that \mS is a CUE-set, write 
$T = S \cap \bigl(\ds\bigcap\limits_n U_n \bigr)$ and 
$S - T = S \cap \bigl(\ds\bigcap\limits_n V_n \bigr)$, where 
$\{U_n\}$ and $\{V_n\}$ are sequences of open subsets of 
$\R$ : $\forall \ n$, $U_n \supset U_{n+1}$ $\&$ $V_n \supset V_{n+1}$.  
Choose open sets
%%----------------------------------------------------------------------------------------------19
$O_n$, $P_n \subset X$:\\
\[
\begin{cases}
\ T - V_n \subset O_n\\
\ (S - T) \cap \overline{O}_n = \emptyset
\end{cases}
,
\quad
\begin{cases}
\ (S - T) - U_n \subset P_n\\
\ T \cap \overline{P}_n = \emptyset
\end{cases}
.
\]
Then put\\
\[
\begin{cases}
\ U = \ds\bigcup\limits_n \big(O_n - \ds\bigcup\limits_{m \leq n} \overline{P}_m \big)\\
\ V = \ds\bigcup\limits_n \big(P_n - \ds\bigcup\limits_{m \leq n} \overline{O}_m \big)
\end{cases}
\\[-.3cm].]
\]
\\
\endgroup%%------------------------------------<<

\label{4.68}
\label{6.29}
A topological space $X$ is said to be 
\un{locally metrizable}
\index{locally metrizable} 
if every point in $X$ has a metrizable neighborhood.  
If $X$ is paracompact and locally metrizable, then $X$ is metrizable.  
Proof: Fix a neighborhood finite open covering $\sU = \{U_i: i \in I\}$ of $X$ consisting of metrizable $U_i$ and choose a development 
$\{\sU_i(n)\}$ for $U_i$ such that $\forall \ n$: $\sU_i(n+1)$ refines $\sU_i(n)$ $-$then the sequence 
$\big\{\ds\bigcup\limits_i\sU_i(1), \ds\bigcup\limits_i\sU_i(2), \ldots \big\}$ is a development for $X$.\\

\begingroup%%----------------------------------->>
\fontsize{9pt}{11pt}\selectfont
\textbf{\small FACT} \ 
Suppose that $X$ is submetacompact and locally metrizable $-$then $X$ is a Moore space.

[Under the stated conditions, every open covering of $X$ has a closed refinement that is neighborhood countable (obvious definition).  
Construct a $\sigma$-closure preserving closed refinement for the latter and thus conclude that $X$ is subparacompact 
(by the characterization mentioned on p. \pageref{1.12}).  
Suppose, then, that $X$ is subparacompact and locally metrizable or, more generally, locally developable in the sense that every 
$x \in X$ has a neighborhood $U_x$ with a development 
$\{\sU_n(x)\}$.  
Let  $\sV = \ds\bigcup\limits_n \sV_n$ be a 
$\sigma$-discrete closed refinement of $\{U_x: x \in X\}$.  
Assign to each $V \in \sV_n$ an element $x_V \in X$ for which 
$V \subset U_{x_V}$, put 
$U_V = \ X - (\cup \sV_n - V)$, and let 
$\sU_{m,n}(V) = U_V \cap \sU_m(x_V)$.  
The collection 
$\sU_{m,n} = \{U: U \in \sU_{m,n}(V) (V \in \sV_n)\} \cup \{X - \cup \sV_n\}$ is an open covering of $X$ and the sequence $\{\sU_{m,n}\}$ is a development for $X$.]\\
\endgroup%%------------------------------------<<

A 
\un{topological manifold}
\index{topological manifold} 
(or an 
\un{$n$-manifold})
\index{n-manifold} 
is a Hausdorff space $X$ for which there exists a nonnegative integer $n$ such that each point of $X$ has a neighborhood that is homeomorphic to an open subset of $\R^n$.

[Note: \  We shall refer to $n$ as the 
\un{euclidean dimension}
\index{euclidean dimension}
 of $X$.  
Homeomorphic topological manifolds have the same euclidean dimension 
(cf. p. \pageref{1.13}).]

Let $X$ be a topological manifold $-$then $X$ is a LCH space.  
As such, $X$ is locally connected.  The components of $X$ are therefore clopen.  Note too that $X$ is locally metrizable.\\

\begingroup%%----------------------------------->>
\fontsize{9pt}{11pt}\selectfont
\textbf{\small FACT} \ 
Let $X$ be a second countable topological manifold of euclidean dimension $n$.  Assume: $X$ is connected $-$then there exists a surjective local homeomorphism $\R^n \ra X$.\\
\endgroup%%------------------------------------<<

\begin{proposition} \ 
Let $X$ be a topological manifold $-$then $X$ is metrizable iff $X$ is paracompact.
\end{proposition}

%%----------------------------------------------------------------------------------------------20
[Note: \  Taking into account the results mentioned on 
p. \pageref{1.14},
 it is also clear that $X$ is metrizable iff each component of $X$ is $\sigma$-compact or, equivalently, iff each component of $X$ is second countable.]\\

\begingroup%%----------------------------------->>
\fontsize{9pt}{11pt}\selectfont
A topological manifold is a Moore space iff it is submetacompact.\\
\endgroup%%------------------------------------<<

\index{long line}
\begingroup%%----------------------------------->>
\fontsize{9pt}{11pt}\selectfont
\textbf{\small EXAMPLE  \ (\un{The Long Line})} \ 
Put $X = [0,\Omega[ \times [0,1[$ and order $X$ by stipulating that 
$(\alpha,x) < (\beta,y)$ if $\alpha < \beta$ or $\alpha = \beta$ and $x < y$.  
Give $X$ the associated order topology $-$then the 
\un{long ray}
\index{long ray} 
$L^+$ is $X - \{(0,0)\}$ and the 
\un{long line} 
$L$ is $X \coprod $X$ / \sim$, $\sim$ meaning that the two origins are identified.
Both $L$ and $L^+$ are normal connected 1-manifolds.  
Neither $L$ nor $L^+$ is $\sigma$-compact, so neither $L$ nor $L^+$ is metrizable.  
Therefore neither $L$ nor $L^+$ is Moore: Otherwise, Reed-Zenor would imply that they are metrizable.
Variant: Moore $\implies$ perfect, which they are not.  
So, neither $L$ nor $L^+$ is submetacompact.  
Finally, observe that $L$ is not homeomorphic to $L^+$.  
Reason: $L$ is countably compact but $L^+$ is not.\\
\endgroup%%------------------------------------<<

\index{Pr\"ufer Manifold}
\begingroup%%----------------------------------->>
\fontsize{9pt}{11pt}\selectfont
\textbf{\small EXAMPLE  \ (\un{The Pr\"ufer Manifold})} \ 
Assign to each $r \in \R$ a copy of the plane 
$\R_r^2 = \R^2 \times \{r\}$ $=$ $\{(a,b,r) \equiv (a,b)_r\}$. \ 
Denote by $\overline{L}_r$ the closed lower half plane in $\R_r^2$, $L_r$ the lower half plane in $\R_r^2$, and 
$\partial \overline{L}_r$ the horizontal axis in $\R_r^2$.
Let $H$ stand for the open upper half plane in $\R^2$. \ 
Put $X = H \cup \ds\bigcup\limits_r \overline{L}_r$.
Topologize $X$ as follows: 
Equip $H$ and each $L_r$ with the usual topology and take for the basic neighborhoods of a typical point 
$(a,0)_r \in \partial \overline{L}_r$ the sets $N(a:r:\epsilon)$, 
a given such being the union of the open rectangle in \  $\overline{L}_r$ 
with corners at $(a \pm \epsilon,0)_r$ and $(a \pm \epsilon, -\epsilon)_r$ and the open wedge consisting of all points within 
$\epsilon$ of $(r,0)$ in the open sector of $H$ bounded by the lines of slope \ 
$1/(a - \epsilon)$ and \ $1/(a + \epsilon)$ emanating from $(r,0)$.
So, e.g., the sequence $(r + 1/n,1/n(a + \epsilon))$ converges to $(a + \epsilon,0)_r$ in the topology of $X$ (although it converges to $(r,0)$ in the usual topology).
The subspace 
$H \cup \{(0,0)_r: r \in \R\}$ (which is not locally compact) is homeomorphic to the Niemytzki plane:
$
\begin{cases}
\ (x,y) \mapsto (x,y^2)\\
\ (0,0)_r \mapsto (r,0)
\end{cases}
\hspace{-.2cm}. \ 
$
$X$ is a connected 2-manifold.
Reason: A closed wedge with its apex removed is homeomorphic to a closed rectangle with one side removed. \ 
It is clear that $X$ is not separable.  
Moreover, $X$ is not second countable, hence is not metrizable (and therefore not paracompact).  
But $X$ is a Moore space: Let $\sU_n$ be the collection comprised of all open disks of radius $1/n$ in $H$ 
and the $L_r$ together with all the $N(a:r:1/n)$ $-$then $\{\sU_n\}$ is a development for $X$.  
This remark allows one to infer that $X$ is not normal: Otherwise, Reed-Zenor would imply that $X$ is metrizable.  
Explicitly, if
$
\begin{cases}
\ A = \{(0,0)_r: r \ \text{rational}\}\\
\ B = \{(0,0)_r: r \ \text{irrational}\} 
\end{cases}
\hspace{-.3cm},
$
then $A$ and $B$ are disjoint closed subsets of $X$ that fail to have disjoint neighborhoods.  
Since $A$ is countable, this means that $X$ cannot be countably paracompact.
However, $X$ is Moore, thus is subparacompact.  
Still, $X$ is not metacompact.  
For $X$ is locally separable (being locally euclidean) and locally separable + metacompact $\implies$ paracompact.  
Apart from all this, $X$ is contractible and so is simply connected.
\vspi
[Note: \  There are two other nonmetrizable, nonnormal connected $2$-manifolds associated with this construction.
\\
%%----------------------------------------------------------------------------------------------21
\indent\indent (1) \ Take two disjoint copies of 
$H \cup \ds\bigcup\limits_r \partial \overline{L}_r$ and identify the corresponding points on the various 
$\partial \overline{L}_r$.  The result is Moore and separable but has an uncountable fundamental group.
\\
\indent\indent (2) \  Take $H \cup \bigcup\limits_r \partial \overline{L}_r$ and $\forall \ r$ identify $(a,0)_r$ and $(-a,0)_r$.   
The result is Moore and separable but has a trivial fundamental group.]\\
\endgroup%%------------------------------------<<

\begingroup%%----------------------------------->>
\fontsize{9pt}{11pt}\selectfont
According to Reed-Zenor, every normal topological Moore manifold is metrizable.  
What happens if we drop ``Moore'' but retain perfection?  
In other words: Is every perfectly normal topological manifold metrizable?  
It turns out this question is undecidable in ZFC.
\\
\indent\indent (1) \ Under MA $+\neg$ CH, every perfectly normal topological manifold is metrizable.

[Lane\footnote[2]{\textit{Proc. Amer. Math. Soc.} \textbf{80} (1980), 693-696; 
see also Balogh-Bennett, \textit{Houston J. Math.} \textbf{15} (1989), 153-162.} 
proved that under MA $+\neg$ CH, every perfectly normal, locally connected LCH space $X$ is paracompact.  
This leads at once to the result.]
\\
\indent\indent (2) \ Under CH, there exist perfectly normal topological manifolds that are not metrizable.
\vspi
[Let $D = \{(x,y) \in \R^2: -1 < x < 1 \ \& \ 0 < y < 1\}$ $-$then the idea here is to coherently paste  
$\Omega$ copies of $[0,1[$ to $D$ via a modification of the Kunen technique 
(cf. p. \pageref{1.15}).  
So let $\{I_\alpha: \alpha < \Omega\}$ be a collection of copies of $[0,1[$ that are unrelated to $\overline{D}$ or to each other.  
Let $\{x_\alpha: \alpha < \Omega\}$ be an enumeration of $\overline{D} - D$.  
Put $X_\alpha = D \cup \bigl ( \ds\bigcup\limits_{\beta < \alpha} I_\beta \bigr)$ and 
$X = \ds\bigcup\limits_{\alpha < \Omega} X_\alpha$.  Let $\{C_\alpha: \alpha < \Omega\}$ 
be an enumeration of the countable subsets of $X$ such that 
$\forall \ \alpha$: $C_\alpha \subset X_\alpha$.  
Define a function $\phi:X \ra \overline{D}$ : $\restr{\phi}{D} = \text{id}_D$ $\&$ $\restr{\phi}{I_\alpha} = x_\alpha$.  
We shall now construct by induction on $\alpha < \Omega$ a topology $\tau_\alpha$ on $X_\alpha$ subject to:
\\
\indent\indent (a) \ $\forall \ \alpha: (X_\alpha,\tau_\alpha)$ is homeomorphic to $D$ and 
$\phi_\alpha = \restr{\phi}{X_\alpha}$ is continuous.\\
\indent\indent (b) \ $\forall \ \beta < \alpha: (X_\beta,\tau_\beta)$ is an open dense subspace of  $(X_\alpha,\tau_\alpha)$.\\
\indent\indent (c) \ $\forall \ \gamma \leq \beta < \alpha:$ If $x_\beta$ is a limit point of $\phi(C_\gamma)$ in $\overline{D}$, then every element of $I_\beta$ is a limit point of $C_\gamma$ in $(X_\alpha,\tau_\alpha)$ .
\vspi
Assign to $D = X_0$ the usual topology.  
If $\alpha$ is a limit ordinal, take for $\tau_\alpha$ the topology on $X_\alpha$ generated by 
$\ds\bigcup\limits_{\beta < \alpha} \tau_\beta$.  
Only condition (a) of the induction hypothesis requires verification.  
This can be dealt with by appealing to a generality: 
Any topological space expressible as the union of an increasing sequence of open subsets, each of which is homeomorphic to 
$\R^n$, is itself homeomorphic to $\R^n$ 
(Brown\footnote[3]{\textit{Proc. Amer. Math. Soc.} \textbf{12} (1961), 812-814.}).  
If $\alpha$ is a successor ordinal, say $\alpha = \beta + 1$, then 
$X_\alpha = X_\beta \cup I_\beta$ and the problem is to define $\tau_\alpha$ knowing $\tau_\beta$.
\vspi
Write $\N = \ds\coprod\limits_1^\infty N_k$: $\forall \ k, \#(N_k) = \omega$ and fix a bijection 
$\iota_k: N_k \ra \Q \ \cap \ ]-1,1[$.
\vspi
Claim: Let $\{U_n\}$ be a sequence of connected open subsets of $D$ and let $\{p_n\}$ be a sequence of distinct points of $D$ : $\forall \ n$,
\\[-1.1cm]
\endgroup

\begingroup
\fontsize{9pt}{11pt}\selectfont
\[
U_n \supset U_{n+1} \ \ \& \ \   D \ \cap \ \bigcap\limits_n \overline{U}_n = \emptyset, \ \  p_n \in U_n.
\]
%%----------------------------------------------------------------------------------------------22
Then there exists an embedding $\mu:D \ra D$ such that $D - \mu(D)$ is homeomorphic to $[0,1[$ and
\\
\indent\indent (i) \ $\forall \ k:$ Each point of $D - \mu(D)$ is a limit point of $\{\mu(p_n): n \in N_k\}$;
\\
\indent\indent (ii) \ $\forall \ n:$ $D - \mu(D)$ is contained in the interior of the closure of $\mu(U_n)$.
\vspi
[To begin with, there exists a homeomorphism $h:D \ra D$ such that 
$\forall \ n$: $h(U_n) \supset D_n$ $\&$ $h(p_n) \in D_n - D_{n+1}$, where $D_n = \{(x,y) \in D: 0 < y < 1/2n\}$.
Choose next a homeomorphism 
$g:D \ra D$ for which the second coordinate of $g(x,y)$ is again $y$ but for which the first coordinate of 
$g(h(p_n))$ is $\iota_k(n)$ $(n \in N_k, k = 1, 2, \ldots)$.
Each point of $\{(x,0): -1 < x < 1\}$ is therefore a limit point of $\{g(h(p_n))$: $n \in N_k\}$.  
Finally, if $F$ is the map with domain $D \cup \{(x,0): -1 < x < 1\}$ defined by
$
\begin{cases}
\ \restr{F}{D} = \text{id}_D\\
\ F(x,0) = (\abs{x},0)
\end{cases}
$
, then the image $D \cup \{(x,0): 0 \leq x < 1\}$, when given the quotient topology, is homeomorphic to $D$ via $f$, say.  
The embedding $\mu = f \circ g \circ h$ satisfies all the assertions of the claim.]
\vspi
To apply the claim, we must specify the $U_n$ and the $p_n$ in terms of $X_\beta$.  
Start by letting $U_n = \phi_\beta^{-1}(O_n(x_\beta))$, where $O_n(x_\beta)$ is the intersection of $\overline{D}$ 
with the open disk of radius $1/n$ centered at $x_\beta$.
Fix a bijection 
$\iota:[0,\beta] \ra \N$ and choose the $p_n \in U_n$ so that if 
$\gamma \leq \beta$ and if $x_\beta$ is a limit point of 
$\phi(C_\gamma)$ in $\overline{D}$, then $p_n \in C_\gamma \cap U_n$ for all $n \in N_{\iota(\gamma)}$.
By assumption, there is a homeomophism $\eta_\beta: X_\beta \ra D$. 
Use this to transfer the data from $X_\beta$ to $D$ and determine an embedding $\mu:D \ra D$.
Put $\mu_\beta = \mu \circ \eta_\beta$, write $D$ as 
$\mu_\beta(X_\beta) \cup (D - \mu_\beta(X_\beta))$ and let 
$\nu_\beta:I_\beta: \ra D - \mu_\beta(X_\beta)$ be a homeomorphism.
The pair $(\mu_\beta,\nu_\beta)$ defines a bijection 
$X_\alpha = X_\beta \cup I_\beta \ra D$.  
Take then for $\tau_\alpha$ the topology on $X_\alpha$ that renders this bijection a homeomorphism and thereby complete the induction.
\vspi
Given $X = \ds\bigcup\limits_{\alpha < \Omega} X_\alpha$ the topology generated by 
$\ds\bigcup\limits_{\alpha < \Omega} \tau_\alpha$ $-$then $X$ is a connected $2$-manifold.  
It is clear that $X$ is not Lindel\"of.  \ 
Because $X$ is separable (in fact is hereditarily separable), it follows that $X$ is not paracompact, thus is not metrizable.  
There remains the verification of perfect normality.  
Let \mA be a closed subset of $X$.  
Fix an $\alpha < \Omega$: $\overline{C}_\alpha = A$.   \ 
Choose a sequence $\{O_n\}$ of open subsets of $\overline{D}$ such that 
$\overline{\phi(C_\alpha)} = \ds\bigcap\limits_n O_n$ $=$ $\ds\bigcap\limits_n \overline{O}_n$.  
Obviously, $A \subset \phi^{-1}(\overline{\phi(C_\alpha)}) =$ 
$\ds\bigcap\limits_n \phi^{-1}(O_n) =$ $\bigcap\limits_n \overline{\phi^{-1}(O_n)}$. \ 
But thanks to condition (c) of the induction hypothesis, 
$\phi^{-1}(\overline{\phi(C_\alpha)}) - A$ is contained in $X_\alpha$.  
So write $X_\alpha - A = \ds\bigcup\limits_n K_n$, $K_n$ compact, and let $P_n$ be a relatively compact open subset of 
$X$ : $K_n \subset P_n \subset \overline{P}_n \subset X - A$.  
To finish, simply note that 
$A = \ds\bigcap\limits_n \ldots_n = \ds\bigcap\limits_n \overline{\overset{ }{\ldots}}_n$, 
$\ldots_n$ being $\phi^{-1}(O_n) - \overline{P}_n$.  Corollary: $X$ is not submetacompact.]
\vspi
The preceding construction is due to 
Rudin-Zenor\footnote[2]{\textit{Houston J. Math.} \textbf{2} (1976), 129-134.}.  
Rudin\footnote[3]{\textit{Topology Appl.} \textbf{35} (1990), 137-152.} 
employed similar methods to produce within ZFC an example of a topological manifold that is both normal and separable, yet is not metrizable.
\vspi
Is every normal topological manifold collectionwise normal?  Recall that this question was asked of an arbitrary LCH space $X$ on p. \pageref{1.16}.  Using the combinatorial principal $\diamond^+$, Rudin (ibid.) established the existence of a normal topological manifold that is not collectionwise normal.  On the other hand, since
%%----------------------------------------------------------------------------------------------23
the cardinality of a connected topological manifold is $2^\omega$, there are axioms that imply a positive answer but I shall not discuss them here.\\
\endgroup%%------------------------------------<<

Let $X$ be a topological space.  
A collection $\{\kappa_i: i \in I\}$ of continuous functions 
$\kappa_i:X \ra [0,1]$ is said to be a 
\un{partition of unity}
\index{partition of unity} 
on $X$ if the supports of the $\kappa_i$ form a neighborhood finite closed covering of $X$ and for every 
$x \in X$, $\ds\sum\limits_i \kappa_i(x) = 1$.  
If  \ $\sU = \{U_i: i \in I\}$ is a covering of $X$, then a partition of unity 
$\{\kappa_i: i \in I\}$ on $X$ is said to be  
\un{subordinate}
\index{subordinate (parition of unity)} 
to $\sU$ if $\forall \ i$: spt $\kappa_i \subset U_i$.

[Note: \   Given a map $f:X \ra \R$, the 
\un{support}
\index{support} 
of $f$, written spt$f$, is the closure of $\{x: f(x) \neq 0\}$.]

A 
\un{numerable}
\index{numerable (covering)} 
covering of $X$ is a covering that has a subordinated partition of unity. 
Examples:  Suppose that $X$ is Hausdorff $-$then
(1) Every neighborhood finite open covering of a normal $X$ is numerable;
(2) Every $\sigma$-neighborhood finite open covering of a countably paracompact normal $X$ is numerable;
(3) Every point  finite open covering of a collectionwise normal $X$ is numerable;
(4) Every open covering of a paracompact $X$ is numerable.

[Note: \   Numerable coverings and their associated partitions of unity allow one to pass from the ``local'' to the ``global'' without the necessity of imposing a paracompactness assumption, a point of some importance in, e.g., fibration theory.]

The requirement on the functions determining a numeration can be substantially weakened.

\label{6.39}
\indent\indent (NU) Suppose given a collection $\{\sigma_i: i \in I\}$ of continuous functions 
$\sigma_i:X \ra [0,1]$ such that 
$\ds\sum\limits_i \sigma_i(x) = 1$ 
$(\forall \ x \in X)$ $-$then there exists a collection 
$\{\rho_i: i \in I\}$ of continuous functions 
$\rho_i:X \ra [0,1]$ such that 
$\forall \ i \in I$ : cl$(\rho_i^{-1}(]0,1])) \subset \sigma_i^{-1}(]0,1])$ and 
(a) \ $\{\rho_i^{-1}(]0,1]):i \in I\}$ is neighborhood finite and
(b) \ $\ds\sum\limits_i \rho_i(x) = 1$ $(\forall \ x \in X)$.

[Of course, at any particular $x \in X$, the cardinality of the set of $i \in I$ such that $\sigma_i(x) \neq 0$ is $\leq \omega$.  Put $\mu = \sup\limits_i \sigma_i$ $-$then $\mu$ is strictly positive.  
Claim: $\mu$ is continuous.  
In fact, $\forall \ \epsilon > 0$, every $x \in X$ has a neighborhood 
$U$ : $\restr{\sigma_i}{U} < \epsilon$ for all but a finite number of $i$, 
thus $\mu$ agrees locally with the maximum of finitely many of the 
$\sigma_i$ and so $\mu$ is continuous.  
Let $\sigma = \ds\sum\limits_i \max\{0, \sigma_i - \mu / 2\}$ and take for 
$\rho_i$ the normalization $\max\{0, \sigma_i - \mu / 2\} / \sigma$.]\\

\begingroup%%----------------------------------->>
\fontsize{9pt}{11pt}\selectfont
Suppose that $H$ is a Hilbert space with orthonormal basis $\{e_i: i \in I\}$.  Let $X$ be the unit sphere in $H$ and set $\sigma_i(x) = \abs{\langle x,e_i \rangle}^2$ $(x \in X)$ $-$then the $\sigma_i$ satisfy the above assumptions.\\
\endgroup%%------------------------------------<<

\begin{proposition} \ 
Every numerable open covering $\sU = \{U_i: i \in I\}$ of $X$ has a numerable open refinement that is both  neighborhood finite and $\sigma$-discrete.
\end{proposition}

%%----------------------------------------------------------------------------------------------24
[Let $\{\kappa_i: i \in I\}$ be a partition of unity on $X$ subordinate to $\sU$.  Denote by $\sF$ the collection of all nonempty finite subsets of $I$.  Assign to each $F \in \sF$ the functions
$
\begin{cases}
\ m_F = \min \kappa_i \quadx (i \in F)\\
\ M_F = \max \kappa_i \quadx (i \notin F)
\end{cases}
$
and put $\mu = \max\limits_F (m_F - M_F)$, which is strictly positive.  
Write $\mu_F$ in place of $m_F - M_F - \mu / 2$, $\sigma_F$ in place of $\max\{0,\mu_F\}$ and set 
$V_F = \{x: \sigma_F(x) > 0\}$ $-$then 
$\overline{V}_F \subset \{x: m_F(x) > M_F(x)\} \subset \ds\bigcap\limits_{i \in F} U_i$.  
The collection 
$\sV = \{V_F: F \in \sF\}$ is a neighborhood finite open refinement of $\sU$ which is in fact 
$\sigma$-discrete as may be seen by defining $\sV_n = \{V_F: \#(F) = n\}$.   \ 
In this connection, note that 
$F^\prime \neq F^{\prime\prime}$ $\&$ 
$\#(F^\prime) = \#(F^{\prime\prime})$ $\implies$ 
$\{x: m_{F^\prime}(x) > M_{F^\prime}(x)\} \cap \{x: m_{F\pp}(x) > M_{F\pp}(x)\}  = \emptyset$.  
The numerability of 
$\sV$ follows upon considering the $\sigma_F / \sigma$ $(\sigma = \sum\limits_F \sigma_F)$.]\\

\begingroup%%----------------------------------->>
\fontsize{9pt}{11pt}\selectfont
Implicit in the proof of Proposition 12 is the fact that if $\sU$ is a numerable open covering of $X$, 
then there exists a countable numerable open covering 
$\sO = \{O_n\}$ of $X$ such that $\forall \ n$, $O_n$ is the disjoint union of open sets each of which is contained in some member of $\sU$.\\
\endgroup%%------------------------------------<<

\label{6.26}
\index{domino principle}
\begingroup%%----------------------------------->>
\fontsize{9pt}{11pt}\selectfont
\textbf{\small FACT  \ (\un{Domino Principle})} \ 
Let $\sU$ be a numerable open covering of $X$.  Assume:\\
\indent\indent (D$_1$) \  Every open subset of a member of $\sU$ is a member of $\sU$.\\
\indent\indent (D$_2$) \ The union of each disjoint collection of members of $\sU$ is a member of $\sU$.\\
\indent\indent (D$_3$) \ The union of each finite collection of members of $\sU$ is a member of $\sU$.
\vspi
Conclusion: \mX is a member of $\sU$.
\vspi
[Work with the \ $O_n$ introduced above, noting that there is no loss of generality in assuming that 
$O_n \subset O_{n+1}$.   \ Choose a precise open refinement 
$\sP = \{P_n\}$ of $\sO$ : $\forall \ n$, $\ov{P}_n \subset P_{n+1}$.  
Put  
$
Q_n =
\begin{cases}
\ P_n \qquad \qquad (n = 1, 2)\\
\ P_n - \ov{P}_{n-2} \ (n \geq 3) 
\end{cases}
$
and write 
$X = \ds\bigcup\limits_1^\infty Q_n =$ 
$\big( \ds\bigcup\limits_1^\infty Q_{2n-1}\big) \cup \big( \ds\bigcup\limits_1^\infty Q_{2n}\big) =$ 
$X_1 \cup X_2$.]\\
\endgroup%%------------------------------------<<

Let $X$ be a topological space $-$then by 
$
\begin{cases}
\ C(X)\\
\ C(X,[0,1])
\end{cases}
$
we shall understand the set of all continuous functions
$
\begin{cases}
\ f:X \ra \R\\
\ f:X \ra [0,1]
\end{cases}
.
$
Bear in mind that $C(X)$ can consist of constants alone, even if $X$ is regular Hausdorff.

\label{19.48}
A 
\un{zero set}
\index{zero set} 
in $X$ is a set of the form \  $Z(f) = \{x: f(x) = 0\}$, where $f \in C(X)$.  \ 
The complement of a zero set is a 
\un{cozero set}
\index{cozero set}.
Since $Z(f) = Z(\min\{1,\abs{f}\})$, $C(X)$ and $C(X,[0,1])$ determine the same collection of zero sets.  All sets of the form
$
\begin{cases}
\ \{x: f(x) \geq 0\}\\
\ \{x: f(x) \leq 0\}
\end{cases}
$
$(f \in C(X))$ are zero sets and all sets of the form
$
\begin{cases}
\ \{x: f(x) > 0\}\\
\ \{x: f(x) < 0\}
\end{cases}
$
$(f \in C(X))$ are cozero sets.  
The collection of zero sets in $X$ is closed under the formation of finite unions and countable intersections and the collection of cozero sets in $X$ is closed under the formation of countable unions and finite intersections.  
The union of a neighborhood finite collection
%%----------------------------------------------------------------------------------------------25
of cozero sets is a cozero set.  
On the other hand, the union of a neighborhood finite collection of zero sets need not be a zero set.  
But this will be the case if each zero set in the collection is contained in a cozero set, the totality of which is neighborhood finite.

[Note: \  Suppose that $X$ is Hausdorff $-$then $X$ is completely regular iff the collection of cozero sets in $X$ is a basis for $X$.  
Every compact $G_\delta$ in a CRH space is a zero set.  
If $X$ is normal, then
$
\begin{cases}
\ \text{closed $G_\delta = $ zero set}\\[-2pt]
\ \text{open $F_\sigma = $ cozero set}
\end{cases}
\hspace{-.2cm}, 
$
so if $X$ is perfectly normal, then 
$
\begin{cases}
\ \text{closed set = zero set}\\[-2pt]
\ \text{open set = cozero set}
\end{cases}
\hspace{-.4cm}.]
$

A
$
\begin{cases}
\ \text{zero set}\\[-2pt]
\ \text{cozero set}
\end{cases}
\hspace{-.4cm}
$
covering of $X$ is a covering consisting of 
$
\begin{cases}
\ \text{zero sets}\\[-2pt]
\ \text{cozero sets}
\end{cases}
\hspace{-.4cm}. 
$
The numerable coverings of $X$ are those coverings that have a neighborhood finite cozero set refinement.  

\label{4.67}
\label{6.2}
\noindent Example:  Every countable cozero set covering $\sU = \{U_n\}$ of $X$ is numerable.  
Proof:  
Choose $f_n \in C(X,[0,1])$: $U_n = f_n^{-1}(]0,1])$, put $\phi_n = 1/2^n \bullet f_n/1 + f_n$ 
$\&$ 
$\phi = \ds\sum\limits_n \phi_n$, let $\sigma_n = \phi_n / \phi$, and apply NU.
\label{5.0w}
\label{19.9}

[Note: \  Every countable cozero set covering \ $\sU = \{U_n\}$ of $X$ has a countable star finite cozero set refinement.  
Proof:  
Choose $f_n \in C(X,[0,1])$ : $U_n = f_n^{-1}(]0,1])$, put $f = \ds\sum\limits_n 2^{-n} f_n$ and define
\[
V_{m.n} =  f_n^{-1}(]0,1]) \cap 
\bigl( f^{-1}\bigl(\bigl]\frac{1}{m+1},1\bigr]\bigr) -  f^{-1}\bigl(\bigl[\frac{1}{m-1},1\bigr]\bigr)\bigr) 
\quadx (1 \leq n \leq m),
\]
with the obvious understanding that if $m = 1$ $-$then the collection $\{V_{m,n}\}$ has the properties in question.]\\

\label{6.3}
\label{19.28}
\label{19.47}
\textbf{\small LEMMA} \ \ 
Let $\sU =\{U_i: i \in I\}$ be a neighborhood finite cozero set covering of $X$ 
$-$then there exists a zero set covering $\sZ = \{Z_i: i \in I\}$ and a cozero set covering 
$\sV = \{V_i: i \in I\}$ such that 
$\forall \ i$: $Z_i \subset V_i \subset \overline{V}_i \subset U_i$.

[Choose a partition of unity $\{\kappa_i: i \in I\}$ on $X$ subordinate to $\sU$.  
Put $V_i = \kappa_i^{-1}(]0,1])$ and take for $Z_i$ the zero set of the function $\max\limits_i \kappa_i - \kappa_i$.]\\

Let $\sU = \{U_i: i \in I\}$ be a neighborhood finite cozero set covering of \mX; 
let $\sZ = \{Z_i: i \in I\}$ and $\sV = \{V_i: i \in I\}$ be as in the lemma.  
Denote by $\sF$ the collection of all nonempty finite subsets of $I$.  
Assign to each $F \in \sF$: $W_F = \ds\bigcap\limits_{i \in F} V_i \cap \big(X - \bigcup\limits_{i \notin F} Z_i\big)$.  
The collection $\sW = \{W_F: F \in \sF\}$ is a neighborhood finite cozero set covering of $X$ such that $\forall \ i$: $\text{st}(Z_i,\sW) \subset V_i$.  
Therefore $\{\text{st}(x,\sW) : x \in X\}$ refines $\sV$, hence $\sU$.  
Now repeat the entire procedure with $\sW$ playing the role of $\sU$.  
The upshot is the following conclusion.\\
%%----------------------------------------------------------------------------------------------26

\begin{proposition} \ %13
Every numerable open covering of $X$ has a numerable open star refinement that is neighborhood finite.\\
\end{proposition}

\begingroup%%----------------------------------->>
\fontsize{9pt}{11pt}\selectfont
\textbf{\small FACT} \ 
Let $\sU = \{U_i: i \in I\}$ be an open covering of $X$ $-$then $\sU$ is numerable iff there exists a metric space $Y$, an an open covering $\sV$ of $Y$, 
and a continuous function $f:X \ra Y$ such that $f^{-1}(\sV)$ refines $\sU$.
\vspi
[The condition is clearly sufficient.  
As for the necessity, let $\{\kappa_i: i \in I\}$ be a partition of unity on $X$ subordinate to $\sU$.  
Let $Y$ be the subset of $[0,1]^I$ comprised of those $y = \{y_i: i \in I\}$ : $\ds\sum\limits_i y_i = 1$.  
The prescription 
$d(y^\prime,y^{\prime\prime}) = \ds\sum\limits_i \abs{y_i^\prime - y_i^{\prime\prime}}$ is a metric on $Y$.
Define a continuous function $f:X \ra Y$ by sending $x$ to $\{\kappa_i(x): i \in I\}$.  
Consider the collection $\sV = \{V_i:i \in I\}$, where $V_i = \{y:y_i > 0\}$.]\\
\endgroup%%------------------------------------<<

\begingroup%%----------------------------------->>
\fontsize{9pt}{11pt}\selectfont
Application: Let $\sU = \{U_i: i \in I\}$ be an open covering of $X$ $-$then $\sU$ is numerable iff there exists a numerable open covering $\sO = \{O_i: i \in I\}$ of $\crg X$ such that $\forall \ i: \text{cr}^{-1}(O_i) \subset U_i$.\\
\endgroup%%------------------------------------<<

\begingroup%%----------------------------------->>
\fontsize{9pt}{11pt}\selectfont
\textbf{\small EXAMPLE}  \ 
Let $G$ be a topological group; let $U$ be a neighborhood of the identity in $G$ $-$then the open covering 
$\{xU: x \in G\}$ is numerable.\\
\endgroup%%------------------------------------<<

Suppose given a set $X$ and a collection $\{X_i: i \in I\}$ of topological spaces $X_i$.\\
\indent\indent (FT) \ Let $\{f_i: i \in I\}$ be a collection of functions $f_i:X_i \ra X$ $-$then the 
\un{final topology}
\index{final topology} 
on $X$ determined by the $f_i$ is the largest topology for which each $f_i$ is continuous.  \ 
The final topology is characterized by the property that if $Y$ is a topological space and if $f:X \ra Y$ is a function, 
then $f$ is continuous iff $\forall \ i$ the composition $f \circ f_i:X_i \ra Y$ is continuous.\\
\indent\indent (IT) \ Let $\{f_i: i \in I\}$ be a collection of functions $f_i:X \ra X_i$ $-$then the  
\un{initial topology}
\index{initial topology} 
on $X$ determined by the $f_i$ is the smallest topology for which each $f_i$ is continuous.  
The initial topology is characterized by the property that if $Y$ is a topological space and if $f:Y \ra X$ is a function, 
then $f$ is continuous iff $\forall \ i$ the composition $f_i \circ f:Y \ra X_i$ is continuous.\\

For example, in the category of topological spaces, coproducts carry the final topology and products carry the initial topology.  The discrete topology on a set $X$ is the final topology determined by the function $\emptyset \ra X$ and the indiscrete topology on a set $X$ is the initial topology determined by the function $X \ra *$.
If $X$ is a topological space and if $f:X \ra Y$ is a surjection, then the final topology on $Y$ determined by $f$ is the quotient topology, while if $Y$ is a topological space and if $f:X \ra Y$ is an injection, then the initial topology on $X$ determined by $f$ is the induced topology.\\

\begingroup%%----------------------------------->>
\fontsize{9pt}{11pt}\selectfont
\textbf{\small EXAMPLE}  \ 
Let $E$ be a vector space over $\R$ $-$then the 
\un{finite topology}
\index{finite topology}
on $E$ is the final topology determined by the inclusions $F \ra E$, 
where $F$ is a finite dimensional linear subspace of $E$ endowed with
%%----------------------------------------------------------------------------------------------27
its natural euclidean topology.  $E$, in the finite topology, is a perfectly normal paracompact Hausdorff space.  
Scalar multiplication $\R \times E \ra E$ is jointly continuous; vector addition $E \times E \ra E$ is separately continuous but jointly continuous iff $\dim E \leq \omega$.  
For a concrete illustration, put $\R^\infty = \bigcup\limits_0^\infty \R^n$, where $\{0\} = \R^0 \subset \R^1 \subset \cdots$.  The elements of $\R^\infty$ are therefore the real valued sequences having a finite number of nonzero values.  
Besides the finite topology, one can also give $\R^\infty$ the inherited product topology 
$\tau_P$ or any of the topologies $\tau_p$ $(1 \leq p \leq \infty)$ derived from the usual $\ell^p$ norm.   
It is clear that $\tau_P \subset \tau_{p^\prime} \subset \tau_{p^{\prime\prime}}$ 
$(1 \leq p^{\prime\prime} < p^{\prime} \leq \infty)$, each inclusion being proper.  
Moreover, $\tau_1$ is strictly smaller than the finite topology.  
To see this, let $U = \{x \in \R^\infty: \forall \ i, \abs{x_i} < 2^{-i}\}$ $-$then $U$ is a neighborhood of the origin in the finite topology but $U$ is not open in $\tau_1$.  
These considerations exhibit uncountably many distinct topologies on $\R^\infty$.  
Nevertheless, under each of them, $\R^\infty$ is contractible, so they all lead to the same homotopy type.
\vspi
[Note: \ The finite topology on $\R^\infty$ is not first countable, thus is not metrizable.]\\
\endgroup%%------------------------------------<<

\begin{proposition} \ %14
Suppose that $X$ is Hausdorff $-$then $X$ is completely regular iff $X$ has the initial topology determined by the elements of 
$C(X)$ (or, equivalently, $C(X,[0,1])$).
\end{proposition}

[Note: \  Therefore, if $\tau^{\prime}$ and $\tau^{\prime\prime}$ are two completely regular topologies on $X$, then $\tau^{\prime} = \tau^{\prime\prime}$ iff, in the obvious notation, $C^{\prime}(X) = C^{\prime\prime}(X)$.]\\

\begingroup%%----------------------------------->>
\fontsize{9pt}{11pt}\selectfont
When constructing the initial topology, it is not necessary to work with functions whose domain is all of $X$.
\vspi
Suppose given a set $X$, a collection $\{U_i: i \in I\}$ of subsets $U_i \subset X$, 
and a collection $\{X_i: i \in I\}$ of topological spaces $X_i$.  Let $\{f_i: i \in I\}$ be a collection of functions $f_i:U_i \ra X_i$ $-$then the 
\un{initial topology}
\index{initial topology} 
on $X$ determined by the $f_i$ is the smallest topology for which each $U_i$ is open and each $f_i$ is continuous.  
The initial topology is characterized by the property that if $Y$ is a topological space and if 
$f:Y \ra X$ is a function, then $f$ is continuous iff $\forall \ i$ the composition 
$f^{-1}(U_i) \overset{f}{\lra} U_i \overset {f_i}{\lra} X_i$ is continuous.\\
\endgroup%%------------------------------------<<

\begingroup%%----------------------------------->>
\fontsize{9pt}{11pt}\selectfont
\textbf{\small EXAMPLE}  \ 
Let $X$ and \mY be nonempty topological spaces $-$then the 
\un{join}
\index{join} 
$X * Y$ is the quotient of 
$X \times Y \times [0,1]$ with respect to the relations
$
\begin{cases}
\ (x,y^\prime,0) \sim (x,y^{\prime\prime},0)\\
\ (x^\prime,y,1) \sim (x^{\prime\prime},y,1)
\end{cases}
. \ 
$
Conventionally 
$
\begin{cases}
\ X * \emptyset = X\\
\ \emptyset * Y = Y
\end{cases}
, \ 
$
so $*$ is a functor $\bTOP \times \bTOP \ra \bTOP$.  
The projection
$
p:
\begin{cases}
\ X \times Y \times [0,1] \ra X * Y\\
\ (x,y,t) \mapsto [x,y,t]
\end{cases}
$
sends \ $X \times Y \times \{0\}$ (or $X \times Y \times \{1\}$) \ 
onto a closed subspace homeomorphic to $X$ (or \mY).  \ 
Consider $X * Y$ as merely a set.  
Let $t : X * Y \lra [0,1]$ be the function $[x,y,t] \mapsto t$; \  let
$
\begin{cases}
\ x: t^{-1}([0,1[) \ra X\\
\ y: t^{-1}(]0,1]) \ra Y
\end{cases}
$
be the functions
$
\begin{cases}
\ [x,y,t] \ra x\\
\ [x,y,t] \ra y
\end{cases}
$
$-$then the 
\un{coarse join}
\index{coarse join}
\label{3.2} 
$X *_c Y$ is $X * Y$ equipped with the initial topology determined by $t$, $x$, and $y$.  \ 
The identity map $X * Y \ra $X$ *_c Y$ is continuous; \ 
it is a homeomorphism if $X$ and $Y$ are compact Hausdorff  but not in general.  \ 
The coarse join $X *_c Y$ of Hausdorff $X$ and $Y$ is Hausdorff, thus so is $X * Y$.  
The join $X * Y$ of path connected $X$ and $Y$ is path connected, thus so is $X *_c Y$.  
Examples: (1)
%%----------------------------------------------------------------------------------------------28
The 
\un{cone}
\index{cone} 
$\Gamma X$ of $X$ is the join of $X$ and a single point; 
(2) The \un{suspension}
\index{suspension} 
$\Sigma X$ of $X$ is the join of $X$ and a pair of points.  
There are also coarse versions of both the cone and the suspension, say
$
\begin{cases}
\ \Gamma_c X\\
\ \Sigma_c X
\end{cases}
\hspace{-.25cm}. \ 
$
Complete the picture by setting 
$
\begin{cases}
\ X *_c \emptyset = X\\
\ \emptyset *_c Y= Y
\end{cases}
\hspace{-.25cm}.
$
\vspi
[Note: \  Analogous definitions can be made in the pointed category $\bTOP_*$.]\\
\endgroup%%------------------------------------<<

\begingroup%%----------------------------------->>
\fontsize{9pt}{11pt}\selectfont
\textbf{\small FACT} \ 
Let $X$ and $Y$ be topological spaces $-$then the identity map $X * Y \ra X *_c Y$  is a homotopy equivalence.
\vspi
[A homotopy inverse $X *_c Y \ra X * Y$ is given by 
$
[x,y,t] \ra
\begin{cases}
\ [x,y,0] \ \hspace{.86cm} \ (0 \leq t \leq 1/3)\\
\ [x,y,3t-1] \hspace{0.4cm} (1/3 \leq t \leq 2/3)\\
\ [x,y,1] \ \hspace{.87cm} \ (2/3 \leq t \leq 1)
\end{cases}
\hspace{-.3cm}.
$
Since the homotopy type of $X * Y$ depends only on the homotopy type of $X$ and \mY and since the coarse join is associative, 
it follows that the join is associative up to homotopy equivalence.]\\
\endgroup%%------------------------------------<<

\label{14.3}
\index{star construction}
\begingroup%%----------------------------------->>
\fontsize{9pt}{11pt}\selectfont
\textbf{\small EXAMPLE  \ (\un{Star Construction})} \ 
The cone $\Gamma X$ of a topological space $X$ is contractible and there is an embedding 
$X \ra \Gamma X$.  However, one drawback to the functor 
$\Gamma: \bTOP \ra \bTOP$ is that it does not preserve embeddings or finite products.  
Another drawback is that while $\Gamma$ does preserve \textbf{HAUS}, within \textbf{HAUS} it need not preserve complete regularity 
(consider $\Gamma X$, where $X$ is the Tychonoff plank).  
The star construction eliminates these difficulties.  
Thus put $\emptyset^* = \emptyset$ and for $X \neq \emptyset$, 
denote by $X^*$ the set of all right continuous step functions $f:[0,1[ \ra X$.  
So, $f \in X^*$ iff there is a partition $a_0 = 0 < a_1 < \cdots < a_n < 1 = a_{n+1}$ of $[0,1[$ such that $f$ is constant on 
$[a_i,a_{i+1}[$ $(i = 0, 1, \ldots, n)$.
There is an injection $i:X \ra X^*$ that sends $x \in X$ to $i(x) \in X^*$, 
the constant step function with value $x$.   
Given $a, b: 0 \leq a < b < 1$, $U$ an open subset of $X$, and $\epsilon > 0$, 
let $O(a,b,U,\epsilon)$ be the set of $f \in X^*$ such that $f$ is constant on 
$[a,b[$, $U$ is a neighborhood of $f(a)$, and the Lebesgue measure of $\{t \in[a,b[: f(t) \notin U\}$ is $< \epsilon$.  
Topologize $X^*$ by taking the $O(a,b,U,\epsilon)$ as a subbasis $-$then $i:X \ra X^*$ is an embedding, which is closed if $X$ is Hausdorff.  
The assignment $X \ra X^*$ defines a functor $\bTOP \ra \bTOP$ that preserves embeddings and finite products.  
It restricts to a functor $\textbf{HAUS} \ra \textbf{HAUS}$ that respects complete regularity.
\vspi
\label{3.5}
Claim: Suppose that $X$ is not empty $-$then $X^*$ is contractible and has a basis of contractible open sets.
\vspi
[Fix $f_0 \in X^*$ and define $H:X^* \times [0,1] \ra X^*$ by $H(f,T)(t) = $
$
\begin{cases}
\ f_0(t) \hspace{.5cm} (0 \leq t < T)\\
\ f(t) \hspace{.625cm} (T \leq t < 1)
\end{cases}
\hspace{-.2cm}.]
$
\\
\endgroup%%------------------------------------<<

An \un{expanding sequence} 
\index{expanding sequence} 
of topological spaces is a system consisting of a sequence of topological spaces $X^n$ linked by embeddings \ 
$f^{n,n+1}:X^n \ra X^{n+1}$.  \  
Denote by $X^\infty$ the colimit in \bTOP associated with this data $-$then for every $n$ there is an arrow \ 
$f^{n,\infty}:X^n \ra X^{\infty}$ and the topology on $X^\infty$ is the final topology determined by the $f^{n,\infty}$.  
Each $f^{n,\infty}$ is an embedding and 
$X^\infty =\ds \bigcup\limits_n f^{n,\infty}(X^n)$.  
One can therefore identify $X^n$ with $f^{n,\infty}(X^n)$ and regard the $f^{n,n+1}$ as inclusions.

%%----------------------------------------------------------------------------------------------29
\label{9.54}
[Note: \  If all the $f^{n,n+1}$ are open (closed) embeddings, then the same holds for all the $f^{n,\infty}$.]

If all the $X^n$ are $\tT_1$, then $X^\infty$ is $\tT_1$.  
If all the $X^n$ are Hausdorff, then $X^\infty$ need not be Hausdorff but there are conditions that lead to this conclusion.

\indent\indent (A) \ If all the $X^n$ are LCH spaces, then $X^\infty$ is a Hausdorff space.

[Let $x, y \in X^\infty$: $x \neq y$.  
Fix an index $n_0$ such that $x,y \in X^{n_0}$.  
Choose open relatively compact subsets $U_{n_0}, V_{n_0} \in X^{n_0}$: $x \in U_{n_0}$, $\&$ $y \in V_{n_0}$ with $\overline{U}_{n_0} \cap \overline{V}_{n_0} = \emptyset$.  
Since $\overline{U}_{n_0}$ and $\overline{V}_{n_0}$ 
are compact disjoint subsets of $X^{n_0 + 1}$, there exist open relatively compact subsets 
$U_{n_0 + 1}, V_{n_0 + 1} \subset X^{n_0 + 1}$ : 
$U_{n_0} \subset U_{n_0 + 1}$ $\&$ $V_{n_0} \subset V_{n_0 + 1}$, with 
$\overline{U}_{n_0 + 1} \cap \overline{V}_{n_0 + 1}  = \emptyset$.  
Iterate the procedure to build disjoint neighborhoods 
$U = \ds\bigcup\limits_{n \geq n_0} U_n$ and 
$V = \ds\bigcup\limits_{n \geq n_0} V_n$ of $x$ and $y$ in $X^\infty$.]

\label{14.7}
\label{14.9}
\label{14.71}
\indent\indent (B) \ Suppose that all the $X^n$ are Hausdorff.  
Assume: $\forall \ n$, $X^n$ is a neighborhood retract of $X^{n+1}$ $-$then $X^\infty$ is Hausdorff.

\label{14.108}
\indent\indent (C) \ If all the $X^n$ are normal (normal and countably paracompact, perfectly normal, collectionwise normal, paracompact) Hausdorff spaces and if $\forall \ n$, $X^n$ is a closed subspace of $X^{n+1}$, then $X^\infty$ is a normal (normal and countably paracompact, perfectly normal, collectionwise normal, paracompact) Hausdorff space.

[The closure preserving closed covering $\{X^n\}$ is absolute, so the generalities on 
p. \pageref{1.17} 
can be applied.]\\

\label{1.24}
\textbf{\small LEMMA} \ \ 
Given an expanding sequence of $\tT_1$ spaces, let 
$\phi:K \ra X^\infty$ be a continuous function such that 
$\phi(K)$ is a compact subset of $X^\infty$ $-$then there exists an index $n$ and a continuous function 
$\phi^n:K \ra X^n$ such that $\phi = f^{n, \infty} \circ \phi^n$.\\

\begingroup%%----------------------------------->>
\fontsize{9pt}{11pt}\selectfont
\textbf{\small EXAMPLE}  \ 
Working in the plane, fix a countable dense subset $S = \{s_n\}$ of $\{(x,y): x = 0\}$.  Put $X^n = \{(x,y): x > 0\} \cup \{s_0, \ldots, s_n\}$ and let $f^{n,n+1}:X^n \ra X^{n+1}$ be the inclusion $-$then $X^\infty$ is Hausdorff but not regular.\\
\endgroup%%------------------------------------<<

\index{Marciszewski Space}
\begingroup%%----------------------------------->>
\fontsize{9pt}{11pt}\selectfont
\textbf{\small EXAMPLE \ (\un{Marciszewski Space})} \ 
Topologize the set $[0,2]$ by isolating the points in $]0,2[$, basic neighborhoods of 0 or 2 being the usual ones.  
Call the resulting space $X_0$.  \ 
Given $n > 0$, topologize the set $]0,2[ \times [0,1]$ by isolating the points of 
$]0,2[ \times ]0,1]$ along with the point $(1,0)$, basic neighborhoods of 
$(t,0)$ $(0 < t < 1$ or $1 < t < 2)$ being the subsets of $L_n$ that contain $(t,0)$ and have a finite complement, 
where $L_n$ is the line segment joining $(t,0)$ and $(t + 1 - 1/n,1)$ $(0 < t < 1)$ or 
$(t,0)$ and $(t - 1 + 1/n,1)$ $(1 < t < 2)$.  Call the resulting space $X_n$.  \ 
Form $X_0 \coprod X_1 \coprod \ldots \coprod X_n$ and let $X^n$ be the quotient obtained by identifying points in $]0,2[$. \  Each $X^n$ is Hausdorff and there is an embedding \ $f^{n,n+1}:X^n \ra X^{n+1}$.  \ 
But $X^\infty$ is not Hausdorff.\\
\endgroup

%%============================================================================
%%============================================================================
\ifnum 1 > 0
\vspace{-.5cm}
\begingroup
\fontsize{9pt}{11pt}\selectfont
\indent
[The underlying point set of $X^\infty$ is \\
\[
X^\infty \ = \ 
\{0,0\} \cup \  
\bigg( (\{0\} \times ]0,2[) \ \cup \ \coprod\limits_{n = 0}^\infty \hsx ]0,2[ \ \times \ [0,1]\bigg) 
\  \cup \{2,0\}.
\]
The topology on $X^\infty$ is the final topology with respect to the embeddings $X^n \hookrightarrow X^\infty$  
i.e., the largest topology on $X^\infty$ such that the maps $X^n \hookrightarrow X^\infty$ are all continuous.
And\\
\[
X^n \ = \ X_0 \coprod X_1 \coprod \ldots \coprod X_n / \sim
\]  
carries  the quotient topology, where $]0,2[ \times \{0\}$ is identified $\forall \ n$
as stated supra.
\endgroup

\begingroup
\fontsize{9pt}{11pt}\selectfont
The space $X^\infty$ may be visualized as a book with a countably infinite number of pages 
each page, $]0,2[ \ \times \ [0,1]$,  identified along the spine, $\{0\} \times ]0,2[)$,  
along with points $\{0,0\}$ and $\{2,0\}$ adjoined at the top and bottom of the spine.
\vspi
%------------------------------------------------------------------------------------------------------------------
The diagram below `depicts' one of the pages, $X_0 \cup \ X_n$ of $X^n$, hence of $X^\infty$. 
The violet segments show the $L_n(x)$ for $x \in \ ]0,\delta[$.  
The green segments show the $L_n(x)$ for some $x \in \ ]2 - \epsilon,2[$.  
\vspi
Claim: \ 
Any neigbhorhood $\sO_{(0,0)}$ of $(0,0)$ has nonempty intersction with any neighborhood $\sO_{(2,0)}$ of $(2,0)$.  
Consequently $X^\infty$ is not Hausdorff.
\vspi
[Consider an open neighborhood $\sO_{(0,0)}$ of $\{0,0)\}$ in $X^\infty$.  
By definition of the topology, and it can assumed that it has the form:\\
\[
\sO_{(0,0)} \ = \ 
[0, \delta[ \ \cup \ \ \bigcup\limits_{n = 0}^\infty \bigcup\limits_{x \in \ ]0,\delta[} L_n^\prime(x) 
\quad \text{where}  \ \exists \ \delta > 0 \ \text{and} \  L_n^\prime(x) \subset L_n(x) \ \text{ is cofinite}.
\]
Similarly, let $\sO_{(2,0)}$ be a neighborhood of $(2,0)$:\\
\[
\sO_{(2,0)} \ = \ 
]2 - \epsilon,2] \ \cup \ \ \bigcup\limits_{n = 0}^\infty \bigcup\limits_{x \in \ ]2 - \epsilon,2[} L_n^\prime(x) 
\quad \text{where} \ \exists \ \epsilon > 0  \ \text{and} \   L_n^\prime(x) \subset L_n(x) \ \text{is cofinite}.
\]
\endgroup

\begingroup
\fontsize{9pt}{11pt}\selectfont
%================================================================
%------------------------------------------------------------------------------------------------------------------
\begin{tikzpicture}[scale=6]%[scale=0.5,shift={(-5,-3)}]
\node[label={{}}] at (0,0) {\textbullet}; \draw[] (0,0) node[below] {$0$};
\draw[violet](0,0) -- (2,0);
\node[label={{}}] at (0,1) {$\circ$}; \draw[] (0,1) node[above] {$1$};
\node[label={{}}] at (1,0) {\textbullet}; \draw[] (1,0) node[below] {$1$};
\node[label={{}}] at (2,0) {\textbullet}; \draw[] (2,0) node[below] {$2$};
\node[label={{}}] at (2,1) {$\circ$}; \draw[] (2,1) node[above] {$(2,1)$};
\node[label={{}}] at (1,-0.2) {$X_0 \coprod X_N / \sim$};
%----------------------------------
\draw[] (0.15,-.045) node[below] {$\delta$};
\node[label={{}}] at (0.15,0) {$\textbf{[}$};
\node[label={{}}] at (0.67,0.35) {$\small\textbf{\textcolor{black}{$\underset{(0 < x < \delta)}{L_{N}(x)}$}}$};
\draw[] (1.8,-.045) node[below] {$2 - \epsilon$};
\draw[] (1.8,0) node[] {{$]$}};
\draw[green](1.8166666667,0) -- (0.82666666667,1);
\draw[green](1.8333333333,0) -- (0.84333333333,1);
\draw[green](1.85,0) -- (0.86,1);
\draw[green](1.8666666667,0) -- (0.87666666667,1);
\draw[green](1.8833333333,0) -- (0.89333333333,1);
\draw[green](1.9,0) -- (0.91,1);
\draw[green](1.9166666667,0) -- (0.92666666667,1);
\draw[green](1.9333333333,0) -- (0.94333333333,1);
\draw[green](1.95,0) -- (0.96,1);
\draw[green](1.9666666667,0) -- (0.97666666667,1);
\draw[green](1.9833333333,0) -- (0.99333333333,1);
\node[label={{}}] at (1.75,0.45)  {$\small\textbf{\textcolor{black}{$\underset{(2 - \epsilon < x < 2)}{L_{N}(x)}$}}$};
%----------------------------------
\draw[black](1,0) -- (1,1);
%----------------------------------
\draw[violet](0.01,0) -- (1,1);
\draw[violet](0.011386138614,0) -- (1.0013861386,1);
\draw[violet](0.012772277228,0) -- (1.0027722772,1);
\draw[violet](0.014158415842,0) -- (1.0041584158,1);
\draw[violet](0.015544554455,0) -- (1.0055445545,1);
\draw[violet](0.016930693069,0) -- (1.0069306931,1);
\draw[violet](0.018316831683,0) -- (1.0083168317,1);
\draw[violet](0.019702970297,0) -- (1.0097029703,1);
\draw[violet](0.021089108911,0) -- (1.0110891089,1);
\draw[violet](0.022475247525,0) -- (1.0124752475,1);
\draw[violet](0.023861386139,0) -- (1.0138613861,1);
\draw[violet](0.025247524752,0) -- (1.0152475248,1);
\draw[violet](0.026633663366,0) -- (1.0166336634,1);
\draw[violet](0.02801980198,0) -- (1.018019802,1);
\draw[violet](0.029405940594,0) -- (1.0194059406,1);
\draw[violet](0.030792079208,0) -- (1.0207920792,1);
\draw[violet](0.032178217822,0) -- (1.0221782178,1);
\draw[violet](0.033564356436,0) -- (1.0235643564,1);
\draw[violet](0.03495049505,0) -- (1.024950495,1);
\draw[violet](0.036336633663,0) -- (1.0263366337,1);
\draw[violet](0.037722772277,0) -- (1.0277227723,1);
\draw[violet](0.039108910891,0) -- (1.0291089109,1);
\draw[violet](0.040495049505,0) -- (1.0304950495,1);
\draw[violet](0.041881188119,0) -- (1.0318811881,1);
\draw[violet](0.043267326733,0) -- (1.0332673267,1);
\draw[violet](0.044653465347,0) -- (1.0346534653,1);
\draw[violet](0.04603960396,0) -- (1.036039604,1);
\draw[violet](0.047425742574,0) -- (1.0374257426,1);
\draw[violet](0.048811881188,0) -- (1.0388118812,1);
\draw[violet](0.050198019802,0) -- (1.0401980198,1);
\draw[violet](0.051584158416,0) -- (1.0415841584,1);
\draw[violet](0.05297029703,0) -- (1.042970297,1);
\draw[violet](0.054356435644,0) -- (1.0443564356,1);
\draw[violet](0.055742574257,0) -- (1.0457425743,1);
\draw[violet](0.057128712871,0) -- (1.0471287129,1);
\draw[violet](0.058514851485,0) -- (1.0485148515,1);
\draw[violet](0.059900990099,0) -- (1.0499009901,1);
\draw[violet](0.061287128713,0) -- (1.0512871287,1);
\draw[violet](0.062673267327,0) -- (1.0526732673,1);
\draw[violet](0.064059405941,0) -- (1.0540594059,1);
\draw[violet](0.065445544554,0) -- (1.0554455446,1);
\draw[violet](0.066831683168,0) -- (1.0568316832,1);
\draw[violet](0.068217821782,0) -- (1.0582178218,1);
\draw[violet](0.069603960396,0) -- (1.0596039604,1);
\draw[violet](0.07099009901,0) -- (1.060990099,1);
\draw[violet](0.072376237624,0) -- (1.0623762376,1);
\draw[violet](0.073762376238,0) -- (1.0637623762,1);
\draw[violet](0.075148514851,0) -- (1.0651485149,1);
\draw[violet](0.076534653465,0) -- (1.0665346535,1);
\draw[violet](0.077920792079,0) -- (1.0679207921,1);
\draw[violet](0.079306930693,0) -- (1.0693069307,1);
\draw[violet](0.080693069307,0) -- (1.0706930693,1);
\draw[violet](0.082079207921,0) -- (1.0720792079,1);
\draw[violet](0.083465346535,0) -- (1.0734653465,1);
\draw[violet](0.084851485149,0) -- (1.0748514851,1);
\draw[violet](0.086237623762,0) -- (1.0762376238,1);
\draw[violet](0.087623762376,0) -- (1.0776237624,1);
\draw[violet](0.08900990099,0) -- (1.079009901,1);
\draw[violet](0.090396039604,0) -- (1.0803960396,1);
\draw[violet](0.091782178218,0) -- (1.0817821782,1);
\draw[violet](0.093168316832,0) -- (1.0831683168,1);
\draw[violet](0.094554455446,0) -- (1.0845544554,1);
\draw[violet](0.095940594059,0) -- (1.0859405941,1);
\draw[violet](0.097326732673,0) -- (1.0873267327,1);
\draw[violet](0.098712871287,0) -- (1.0887128713,1);
\draw[violet](0.1000990099,0) -- (1.0900990099,1);
\draw[violet](0.10148514851,0) -- (1.0914851485,1);
\draw[violet](0.10287128713,0) -- (1.0928712871,1);
\draw[violet](0.10425742574,0) -- (1.0942574257,1);
\draw[violet](0.10564356436,0) -- (1.0956435644,1);
\draw[violet](0.10702970297,0) -- (1.097029703,1);
\draw[violet](0.10841584158,0) -- (1.0984158416,1);
\draw[violet](0.1098019802,0) -- (1.0998019802,1);
\draw[violet](0.11118811881,0) -- (1.1011881188,1);
\draw[violet](0.11257425743,0) -- (1.1025742574,1);
\draw[violet](0.11396039604,0) -- (1.103960396,1);
\draw[violet](0.11534653465,0) -- (1.1053465347,1);
\draw[violet](0.11673267327,0) -- (1.1067326733,1);
\draw[violet](0.11811881188,0) -- (1.1081188119,1);
\draw[violet](0.1195049505,0) -- (1.1095049505,1);
\draw[violet](0.12089108911,0) -- (1.1108910891,1);
\draw[violet](0.12227722772,0) -- (1.1122772277,1);
\draw[violet](0.12366336634,0) -- (1.1136633663,1);
\draw[violet](0.12504950495,0) -- (1.115049505,1);
\draw[violet](0.12643564356,0) -- (1.1164356436,1);
\draw[violet](0.12782178218,0) -- (1.1178217822,1);
\draw[violet](0.12920792079,0) -- (1.1192079208,1);
\draw[violet](0.13059405941,0) -- (1.1205940594,1);
\draw[violet](0.13198019802,0) -- (1.121980198,1);
\draw[violet](0.13336633663,0) -- (1.1233663366,1);
\draw[violet](0.13475247525,0) -- (1.1247524752,1);
\draw[violet](0.13613861386,0) -- (1.1261386139,1);
\draw[violet](0.13752475248,0) -- (1.1275247525,1);
\draw[violet](0.13891089109,0) -- (1.1289108911,1);
\draw[violet](0.1402970297,0) -- (1.1302970297,1);
\draw[violet](0.14168316832,0) -- (1.1316831683,1);
\draw[violet](0.14306930693,0) -- (1.1330693069,1);
\draw[violet](0.14445544554,0) -- (1.1344554455,1);
\draw[violet](0.14584158416,0) -- (1.1358415842,1);
\draw[violet](0.14722772277,0) -- (1.1372277228,1);
\draw[violet](0.14861386139,0) -- (1.1386138614,1);
\draw[violet](0.15,0) -- (1.14,1);
%----------------------------------
%----------------------------------
%%%
\draw[densely dotted](0,0.05) -- (2.0,0.05);
\draw[densely dotted](0,0.10) -- (2.0,0.10);
\draw[densely dotted](0,0.15) -- (2.0,0.15);
\draw[densely dotted](0,0.20) -- (2.0,0.20);
\draw[densely dotted](0,0.25) -- (2.0,0.25);
\draw[densely dotted](0,0.30) -- (2.0,0.30);
\draw[densely dotted](0,0.35) -- (2.0,0.35);
\draw[densely dotted](0,0.40) -- (2.0,0.40);
\draw[densely dotted](0,0.45) -- (2.0,0.45);
\draw[densely dotted](0,0.50) -- (2.0,0.5);
\draw[densely dotted](0,0.55) -- (2.0,0.55);
\draw[densely dotted](0,0.60) -- (2.0,0.60);
\draw[densely dotted](0,0.65) -- (2.0,0.65);
\draw[densely dotted](0,0.70) -- (2.0,0.70);
\draw[densely dotted](0,0.75) -- (2.0,0.75);
\draw[densely dotted](0,0.80) -- (2.0,0.80);
\draw[densely dotted](0,0.85) -- (2.0,0.85);
\draw[densely dotted](0,0.90) -- (2.0,0.90);
\draw[densely dotted](0,0.95) -- (2.0,0.95);
\draw[densely dotted](0,1.0) -- (2.0,1.0);
\end{tikzpicture}
%%%}
%%%}
\vspi
Fix $\ds N > \max\bigg\{\frac{1}{\delta},\frac{1}{\epsilon}\bigg\}$ $-$then
\[
\begin{cases}
\  \text{right endpoint } L_N(x) = x + 1 - 1/N > 1 \qquad \forall \ x \in \ ]0,\delta[ \\[11pt]
\  \text{left  endpoint  } \hspace{0.25cm}  L_N(x) =  2 - 1 + 1/N < 1  \qquad \forall \ x \in \ ]2 - \epsilon,2[ 
\end{cases}
.
\] 
%----------\%----------%----------%----------%----------%----------%----------
%Set 
%\[
%\sU 
%\ = \ 
%\bigcup\limits_{x \in \ ]0,\delta[} L_N(x)  \qquad \text{(i.e., the violet region diagram)}.
%\]
Let  $\{y_k\}$ be a sequence of distinct numbers such that 
\[
2 - \epsilon  
\ < \ 
2 - \frac{1}{N}
\ < \  
y_k 
\ < \ 
2 \qquad \forall \ k = 0, 1, \ldots.
\]
E.g.: Take all rationals in $]2 - \frac{1}{N},2[$.
\vspi
So  $\{y_k\} \subset \ \sO_{(2,0)}$ and for each $k = 0, 1, \ldots$, there are but finitely many 
$x \in \ ]0,\delta[$ such that $L_N(x) \cap L^\prime_N(y_k) \neq \emptyset$.
Set
\[
F_k 
\ = \ 
\{x: x \in \ ]0,\delta[  \ \& \ L_N(x) \ \cap \ L_N^\prime(y_k) \ = \ \emptyset\}. 
\]
So
\[
\#(F_k) \ <  \ \infty \qquad \forall \ k = 0, 1, \ldots.
\]
And
\[
F \ = \ \bigcup\limits_{k=0}^\infty \hsx F_k \qquad \text{is countable}.
\]
Let $x_0 \in \ ]0,\delta[ - F$ $-$then
\[
L_N(x_0) \ \cap \ L_N^\prime(y_k) \ \neq \ \emptyset \ \quad \forall \ k.  
\]
\qquad\qquad $\implies$
\[
L_N^\prime(x_0) \ \cap \ L_N^\prime(y_{k_0}) \ \neq \ \emptyset \  \qquad \text{$\exists \ k_0$ (cofinality)}
\]
\qquad\qquad $\implies$
\[
\sO_{(0,0)}  \ \cap \ \sO_{(2,0)} \neq \ \emptyset.
\]
I.e., $X^\infty$ is not Hausdorff.]
\fi
\\
\endgroup%%------------------------------------<<

%%============================================================================
%%============================================================================

%%----------------------------------------------------------------------------------------------30
\begingroup%%----------------------------------->>
\fontsize{9pt}{11pt}\selectfont
\textbf{\small FACT} \ 
Suppose that
$
\begin{cases}
\ X^0 \subset X^1 \subset \cdots\\
\ Y^0 \subset Y^1 \subset \cdots
\end{cases}
$
are expanding sequences of LCH spaces $-$then $X^\infty \times Y^\infty = \text{colim }(X^n \times Y^n)$.\\
\endgroup%%------------------------------------<<

Let $X$ be a topological space $-$then a 
\un{filtration}
\index{filtration} 
on $X$ is a sequence $X^0, X^1, \ldots$ of subspaces of $X$ such that $\forall \ n$: $X^n \subset X^{n+1}$.  \ 
Here, one does not require that $\ds\bigcup\limits_n X^n = X$.  \ 
A 
\un{filtered space}
\index{filtered space} 
\bX 
is a topological space $X$ equipped with a filtraton $\{X^n\}$.  
A \un{filtered map}
\index{filtered map} 
$\bff:\bX \ra \bY$ of filtered spaces is a continuous function 
$f:X \ra Y$ such that $\forall \ n$ : $f(X^n) \subset Y^n$.
Notation: $\bff \in \bC(\bX,\bY)$.  
\textbf{FILSP}
\index{\textbf{FILSP}} 
is the category whose objects are the filtered spaces and whose morphisms are the filtered maps.  
\textbf{FILSP} is a symmetric monoidal category:  
Take $\bX \otimes \bY$ to be $X \times Y$ supplied with the filtration
$n \ra \ds\bigcup\limits_{p+q=n} X^p \times Y^q$, 
let $e$ be the one point space filtered by specifying that the initial term is $\neq \emptyset$, 
and make the obvious choice for $\top$.  There is a notion of homotopy in \textbf{FILSP}.  \ 
Write \textbf{I} for $I = [0,1]$ endowed with its skeletal filtration, i.e., $I^0 = \{0,1\}$, $I^n = [0,1]$ $(n \geq 1)$ $-$then filtered maps $\bff, \ \textbf{g}:\bX \ra \bY$ are said to be 
\un{filter homotopic}
\index{filter homotopic} 
if there exists a filtered map $\textbf{H}:\bX \otimes \textbf{I} \ra \bY$ such that 
$
\begin{cases}
\ H(x,0) = f(x)\\
\ H(x,1) = g(x)
\end{cases}
(x \in X).
$
\\

\begingroup%%----------------------------------->>
\fontsize{9pt}{11pt}\selectfont
Geometric realization may be viewed as a functor \ $\abs{?} : \textbf{SISET} \ra \textbf{FILSP}$ \ via consideration of skeletons.  
To go the other way, equip $\Delta^n$ with its skeletal filtration and let \textbf{$\bDelta^n$} be the associated filtered space.  
Given a filtered space \bX, write $\bsin \bX$ for the simplicial set defined by 
$\bsin \bX ([n])$ = \textbf{sin$_n$X} = \textbf{C$(\bDelta^n,\bX)$} $-$then the assignment 
$\bX \ra \bsin \bX$ is a functor $\textbf{FILSP} \ra \textbf{SISET}$ and $(\abs{?},\bsin)$ is an adjoint pair.\\
\endgroup%%------------------------------------<<

If \textbf{C} is a full subcategory of \bTOP (\textbf{HAUS}) and if $X$ is a topological space (Hausdorff topological space), then $X$ is an object in the monocoreflective hull of \textbf{C} in \bTOP (\textbf{HAUS}) iff there exists a set $\{X_i\} \subset \text{Ob}\textbf{C}$ and an extremal epimorphism $f:\coprod\limits_i X_i \ra X$ 
(cf. p. \pageref{1.18} ff.).
Example: The monocoreflective hull  in \bTOP of the full subcategory of \bTOP whose objects are the locally connected, connected spaces is the category of locally connected spaces.

[Note: \  The categorical opposite of ``epireflective'' is ``monocoreflective''.]\\

\begingroup%%----------------------------------->>
\fontsize{9pt}{11pt}\selectfont
\index{A spaces}
\textbf{\small EXAMPLE \  (\un{A Spaces})} \ 
The monocoreflective hull in \bTOP of $[0,1]/[0,1[$ is the category of A spaces.\\
\endgroup%%------------------------------------<<

\begingroup%%----------------------------------->>
\fontsize{9pt}{11pt}\selectfont
\index{sequential spaces}
\textbf{\small EXAMPLE \ (\un{Sequential Spaces})} \ 
A topological space $X$ is said to be \un{sequential} provided that a subset $U$ of $X$ is open iff every sequence converging to a point of $U$ is eventually in $U$.  Every first 
%%----------------------------------------------------------------------------------------------31
countable space is sequential.  On the other hand, a compact Hausdorff space need not be sequential (consider ($[0,\Omega]$).
Example: The one point compactification of the Isbell$-$Mr\'owka space $\Psi(\N)$ is sequential but there is no sequence in $\N$ converging to $\infty \in \overline{\N}$.  
If \textbf{SEQ} is the full, isomorphism closed subcategory of \bTOP whose objects are the sequential spaces, 
then \textbf{SEQ} is closed under the formation in \bTOP of coproducts and quotients.  
Therefore \textbf{SEQ} is a monocoreflective subcategory of \bTOP
(cf. p. \pageref{1.19}), 
hence is complete and cocomplete.  
The coreflector sends $X$ to its 
\un{sequential modification}
\index{sequential modification} 
\index{sX}
$sX$.  
Topologically, $sX$ is $X$ equipped with the final topology determined by the $\phi \in C(\N_\infty,X)$, 
where $\N_\infty$ is the one point compactification of $\N$ (discrete topology).  
The monocoreflective hull in \bTOP of $\N_\infty$ is \textbf{SEQ}, so a topological space is sequential iff it is a quotient of a first countable space.  
\textbf{SEQ} is cartesian closed: $C(s(X \times Y), Z) \approx C(X,Z^Y)$.  
Here, $s(X \times Y)$ is the product in \textbf{SEQ} (calculate the product in \bTOP and apply $s$).  
As for the exponential object $Z^Y$, given any open subset $P \subset Z$ and any continuous function $\phi:\N_\infty \ra Y$, put $O(\phi,P) = \{g \in C(Y,Z): g(\phi(\N_\infty)) \subset P\}$ and call $C_s(Y,Z)$ the result of topologizing $C(Y,Z)$ by letting the $O(\phi,P)$ be a subbasis $-$then $Z^Y = sC_s(Y,Z)$.
\vspi
[Note: \  Every CW complex is sequential.]\\
\endgroup%%------------------------------------<<

A Hausdorff space $X$ is said to be 
\un{compactly generated}
\index{compactly generated} 
provided that a subset $U$ of $X$ is open iff $U \cap K$ is open in $K$ for every compact subset $K$ of $X$. Examples
(1) \ Every LCH space is compactly generated; 
(2) \ Every first countable Hausdorff space is compactly generated; 
(3) \ The product $\R^\kappa$, $\kappa > \omega$, is not compactly generated.  A Hausdorff space is compactly generated iff it can be represented as the quotient of a LCH space.  Open subspaces and closed subspaces of compactly generated Hausdorff spaces are compactly generated, although this is not the case for arbitrary subspaces (consider $\N \cup \{p\} \subset \beta \N$, where $p \in \beta\N - \N$).  However, 
Arhangel'ski\u i\footnote[2]{\textit{Czech. Math. J.} \textbf{18} (1968), 392-395.} 
has shown that if $X$ is a Hausdorff space, then $X$ and all its subspaces are compactly generated iff for every $A \subset X$ and each $x \in \overline{A}$ there exists a sequence $\{x_n\} \subset A$: $\lim x_n = x$.  The product $X \times Y$ of two compactly generated Hausdorff spaces may fail to be compactly generated (consider $X = \R - \{1/2,1/3,\ldots\}$ and $Y = \R/\N$) but this will be true if one of the factors is a LCH space or if both factors are first countable.\\

\begingroup%%----------------------------------->>
\fontsize{9pt}{11pt}\selectfont
\index{sequential spaces}
\textbf{\small EXAMPLE \  (\un{Sequential Spaces})} \ 
A Hausdorff sequential space is compactly generated.  
In fact, a Hausdorff space is sequential provided that a subset $U$ of $X$ is open iff $U \cap K$ is open in $K$ for every second countable compact subset $K$ of $X$.\\
\endgroup%%------------------------------------<<

\begingroup%%----------------------------------->>
\fontsize{9pt}{11pt}\selectfont
\textbf{\small EXAMPLE}  \ 
Equip $\R^\infty$ with the finite topology and let $H(\R^\infty)$ be its homeomorphism group.  
%%----------------------------------------------------------------------------------------------32
Give $H(\R^\infty)$ the compact open topology $-$then $H(\R^\infty)$ is a perfectly normal paracompact Hausdorff space.  But $H(\R^\infty)$ is not compactly generated.
\vspi
[The set of all linear homeomorphisms $\R^\infty \ra \R^\infty$ is a closed subspace of $H(\R^\infty)$.  Show that it is not compactly generated.  Incidentally, $H(\R^\infty)$ is contractible.]\\
\endgroup%%------------------------------------<<

\label{9.91}
\label{9.93}
For certain purposes in algebraic topology, it is desirable to single out a full, isomorphism closed subcategory of \bTOP, small enough to be ``convenient'' but large enough to be stable for the ``standard'' constructions.
A popular candidate is the category \bCGH of compactly generated Hausdorff spaces 
(Steenrod\footnote[3]{\textit{Michigan Math. J.} \textbf{14} (1967), 133-152.}).  
Since \bCGH is closed under the formation in \textbf{HAUS} of coproducts and quotients, \bCGH is a monocoreflective subcategory of \textbf{HAUS} 
(cf. p. \pageref{1.20}).  
As such, it is complete and cocomplete.  
The coreflector sends $X$ to its 
\un{compactly generated modification}
\index{compactly generated modification} 
$kX$.
Topologically, $kX$ is $X$ equipped with the final topology determined by the inclusions $K \ra X$, $K$ running through the compact subsets of $X$.  
The identity map $kX \ra X$ is continuous and induces isomorphisms of homotopy and singular homology and cohomology groups.  
If $X$ and $Y$ are compactly generated, then their product in \bCGH is $X \times_k Y \equiv k(X \times Y)$.  
Each of the functors $- \times_k Y:\bCGH \ra \bCGH$ has a right adjoint 
$Z \ra Z^Y$, the exponential object $Z^Y$ being $kC(Y,Z)$, where $C(Y,Z)$ carries the compact open topology. 
\label{13.6} 
So one of the advantages of \bCGH is that it is cartesian closed.  
Another advantage is that if
$
\begin{cases}
\ X, X^\prime\\
\ Y, Y^\prime
\end{cases}
$
are in \bCGH and if 
$
\begin{cases}
\ f:X \ra X^\prime\\
\ g:Y \ra Y^\prime
\end{cases}
$
are quotient, then 
$f \times_k g:X \times_k Y \ra X^\prime \times_k Y^\prime$ is quotient.  
But there are shortcomings as well.  Item: The forgetful functor $\bCGH \ra \bTOP$ does not preserve colimits.  
For let $A$ be a compactly generated subspace of $X$ and consider the pushout square
\begin{tikzcd}%[ sep=small]
A \ar{d} \ar{r} &{*} \ar{d}\\
X \ar{r} &P
\end{tikzcd}
in \bCGH $-$then $P = h(X/A)$, the maximal Hausdorff quotient of the ordinary quotient computed in \bTOP.  
To appreciate the point, let $X = [0,1]$, $A = [0,1[$ $-$then $[0,1]/[0,1[$ is not Hausdorff and $h([0,1]/[0,1[)$ is a singleton.  Finally, it is clear that \bCGH is the monocoreflective hull in \textbf{HAUS} of the category of compact Hausdorff spaces.\\

\begingroup%%----------------------------------->>
\fontsize{9pt}{11pt}\selectfont

$\bCGH_*$, the category of pointed compactly generated Hausdorff spaces, is a closed category: 
Take $X \otimes Y$ to be the smash product $X \#_k Y$ 
(cf. p. \pageref{1.21}) 
and let $e$ be $\textbf{S}^0$.  Here, the internal hom functor sends $(X,Y)$ to the closed subspace of $kC(X,Y)$ consisting of the base point preserving continuous functions.\\

\textbf{\small FACT} \ 
Let $X$ be a CRH space.  
Suppose that there exists a sequence $\{\sU_n\}$ of open coverings of $X$ such that 
$\forall \ x \in X$: 
$K_x \equiv \ds\bigcap\limits_n \text{st}(x,\sU_n)$ is compact and $\{\text{st}(x,\sU_n)\}$ is a neighborhood basis of $K_x$ (i.e., any
%%----------------------------------------------------------------------------------------------33
open $U$ containing $K_x$ contains some $\text{st}(x,\sU_n))$ $-$then $X$ is compactly generated.
Example: Every Moore space is compactly generated.

[Note: \  
Jiang\footnote[2]{\textit{Topology Proc.} \textbf{11} (1986), 309-316.} 
has shown that any CRH space $X$ realizing this assumption is necessarily submetacompact.]\\
\endgroup%%------------------------------------<<

\label{14.120a}
In practice, it can be troublesome to prove that a given space is Hausdorff and while this is something which is nice to know, there are situations when it is irrelevant.  
We shall therefore englarge \bCGH to its counterpart in \bTOP, the category \bCG
\index{\bCG} %\symbol
of 
\un{compactly generated}
\index{compactly generated spaces}  spaces 
(Vogt\footnote[3]{\textit{Arch. Math.} \textbf{22} (1971), 545-555; 
see also Wyler, \textit{General Topology Appl.} \textbf{3} (1973), 225-242.}), 
by passing to the monocoreflective hull in \bTOP of the category of compact Hausdorff spaces.
It is thus immediate that a topological space is compactly generated iff it can be represented as the quotient of a LCH space.  Consequently, if $X$ is a topological space, then $X$ is compactly generated provided that a subset $U$ of $X$ is open iff $\phi^{-1}(U)$ is open in $K$ for every 
$\phi \in C(K,X)$, $K$ any compact Hausdorff space.
What has been said above in the Hausdorff case is now applicable in general, the main difference being that the forgetful functor $\bCG \ra \bTOP$ preserves co\-limits.  
Also, like \bCGH, \bCG is cartesian closed: $C(X \times_k Y,Z) \approx C(X,Z^Y)$.  Of course, $X \times_k Y \equiv k(X \times Y)$ and the exponential object $Z^Y$ is defined as follows.  
Given any open subset $P \subset Z$ and any continuous function $\phi:K \ra Y$, where $K$ is a compact Hausdorff space, 
put $O(\phi,P) = \{g \in C(Y,Z): g(\phi(K)) \subset P\}$ and call $C_k(Y,Z)$ the result of topologizing $C(Y,Z)$ 
by letting the $O(\phi,P)$ be a subbasis $-$then $Z^Y = kC_k(Y,Z)$.  
Example: A sequential space is compactly generated.

[Note: If $X$ and $Y$ are compactly generated and if $f:X \ra Y$ is a continuous injection, then $f$ is an extremal monomorphism iff the arrow $X \ra kf(X)$ is a homeomorphism, where $f(X)$ has the induced topology.  
Therefore an extremal monomorphism in \bCG need not be an embedding (= extremal monomorphism in \bTOP).  
Extremal monomorphisms in \bCG are regular.  
Call them 
\un{\bCG embeddings}.]
\index{compactly generated embeddings, \bCG embeddings}\\

\begingroup%%----------------------------------->>
\fontsize{9pt}{11pt}\selectfont
\textbf{\small EXAMPLE}  \ 
Partition $[-1,1]$ by writing 
$[-1,1] = \{-1\} \cup \ds\bigcup\limits_{0 \leq x < 1} \{x,-x\} \cup \{1\}$.  
Let $X$ be the associated quotient space $-$then $X$ is compactly generated (in fact, first countable).  
Moreover, $X$ is compact and $\mT_1$ but not Hausdorff; $X$ is also path connected.\\

\textbf{\small FACT} \ 
Let $X$ and Y be compactly generated $-$then the projections
$
\begin{cases}
\ $X$ \times_k Y \ra X\\
\ $X$ \times_k Y \ra Y
\end{cases}
$
are open maps.\\
\endgroup%%------------------------------------<<

%%----------------------------------------------------------------------------------------------34
Given any class $\sK$ of compact spaces containing at least one nonempty space, denote by \textbf{M} the monocoreflective hull of $\sK$ in \bTOP and let $R:\bTOP \ra \textbf{M}$ be the associated coreflector.  
If $X$ is a topological space, then a subset of $U$ of $RX$ is open provided that $\phi^{-1}(U)$ is open in 
$K$ for every $\phi \in C(K,X)$, $K$ any element of $\sK$.  
Write $\bDelta$-\bK
\index{$\bDelta$-\bK} %\symbol
for the full, isomorphism closed subcategory of \bTOP whose objects are those 
$X$ which are 
\un{$\Delta$-separated}
\index{Delta-separated, $\Delta$-separated} 
by $\sK$, i.e., such that 
$\Delta_X \equiv \{(x,x): x \in X\}$ is closed in $R(X \times X)$ 
$-$then $\bDelta$-\bK is closed under the formation in \bTOP of products and embeddings.  
Therefore $\bDelta$-\bK is an epireflective subcategory of \bTOP 
(cf. p. \pageref{1.22}).
Examples:
(1) \ Take for $\sK$ the class of all finite indiscrete spaces $-$then an $X$ in \bTOP is $\Delta$-separated by $\sK$ iff it is $\mT_0$;
(2) \ Take for $\sK$ the class of all finite spaces $-$then an $X$ in \bTOP is $\Delta$-separated by $\sK$ iff it is $\mT_1.$

[Note: \  Recall that a topological space $X$ is Hausdorff iff its diagonal is closed in $X \times X$ (product topology).]\\

\begingroup%%----------------------------------->>
\fontsize{9pt}{11pt}\selectfont
\index{sequential spaces}
\textbf{\small EXAMPLE \  (\un{Sequential Spaces})} \ 
Let $X$ be a topological space $-$then every sequence in $X$ has at most one limit iff $\Delta_X$ is sequentially closed in $X \times X$, i.e., iff $X$ is $\Delta_X$-separated by $\sK = \{\N_\infty\}$.  
When this is so, $X$ must be $\mT_1$ and if $X$ is first countable, then $X$ must be Hausdorff.

[Note: \  Recall that a topological space $X$ is Hausdorff iff every net in $X$ has at most one limit.]\\
\endgroup%%------------------------------------<<

If $K$ is a compact space, then for any $\phi \in C(K,X)$, $\phi(K)$ is a compact subset of $X$. 
In general, $\phi(K)$ is neither closed nor Hausdorff.\\
\indent\indent ($\sK_1$) \quadx A topological space $X$ is said to be $\sK_1$ provided that $\forall \ \phi \in C(K,X)$ $(K \in \sK)$, $\phi(K)$ is a closed subspace of $X$.

\indent\indent ($\sK_2$) \quadx A topological space $X$ is said to be $\sK_2$ provided that $\forall \ \phi \in C(K,X)$ $(K \in \sK)$, $\phi(K)$ is a Hausdorff subspace of $X$.

A topological space $X$ which is simultaneously $\sK_1$ and $\sK_2$ is necessarily 
$\Delta$-separated by $\sK$.

Specialize the setup and take for $\sK$ the class of compact Hausdorff spaces 
(McCord\footnote[2]{\textit{Trans. Amer. Math. Soc.} \textbf{146} (1969), 273-298; 
see also Hoffmann, \textit{Arch. Math.} \textbf{32} (1979), 487-504.}), 
so \bM = \bCG.  
Suppose that $X$ is $\sK_1$ (hence $\mT_1$) $-$then $X$ is $\sK_2$.  Proof: Let
$
\begin{cases}
\ x\\
\ y
\end{cases}
\hspace{-.3cm} \in \phi(K)
$
$ (\phi \in C(K,X))$: $x \neq y$, choose disjoint open sets 
$
\begin{cases}
\ U\\
\ V
\end{cases}
\hspace{-.3cm} \subset K \ : \ 
\begin{cases}
\ \phi^{-1}(x) \subset U\\
\ \phi^{-1}(y) \subset V
\end{cases}
$
and consider 
$
\begin{cases}
\ \phi(K) - \phi(K - U)\\
\ \phi(K) - \phi(K - V)
\end{cases}
\hspace{-.3cm} . \ 
$
Denote by $\bDelta\text{-}\bCG$
\index{$\bDelta$-\bCG} %symbol
the full subcategory of \bCG whose objects are $\Delta$-separated by $\sK$.  \ 
There are strict inclusions 
$\bCGH \subset$ $\bDelta\text{-}\bCG$ $\subset \bCG$.  
Example: Every first countable 
$X$ in $\bDelta$-\bCG is Hausdorff.\\

%%----------------------------------------------------------------------------------------------35
\textbf{\small LEMMA} \ \ 
Let $X$ be a $\Delta$-separated compactly generated space $-$then $X$ is $\sK_1$.

[Let $K,\  L \in \sK$; let $\phi \in C(K,X)$, $\psi \in C(L,X)$.  Since $\phi \times \psi: K \times L \ra $X$ \times_k X$ is continuous, $(\phi \times \psi)^{-1}(\Delta_X)$ is closed in $K \times L$.  Therefore $\psi^{-1}(\phi(K)) = \text{pr}_L((\phi \times \psi)^{-1} (\Delta_X))$ is closed in $L$.]\\

\label{14.35}
It follows from the lemma that every $\Delta$-separated compactly generated space $X$ is $\mT_1$.  
More is true: Every compact subspace $A$ of $X$ is closed in $X$. 
Proof: For any $\phi \in C(K,X)$ $(K \in \sK)$, 
$A \cap \phi(K)$ is a closed subspace of $A$, thus is compact, so 
$A \cap \phi(K)$ is a closed subspace of $\phi(K)$, implying that $\phi^{-1}(A) = \phi^{-1}(A \cap \phi(K))$ is closed in $K$.  
Corollary: The intersection of two compact subsets of $X$ is compact.

Equalizers in \bCGH and $\dcg$ are closed (e.g. retracts) but $\dcg$ is better behaved than \bCGH when it comes to quotients.  
\label{14.97a}
Indeed, if $X$ is in $\dcg$ and if $E$ is an equivalence relation on $X$, then $X/E$ is in $\bDelta$-\bCG iff 
$E \subset X \times_k X$ is closed.
To see this, let $p:X \ra X/E$ be the projection.  
Because $p \times_k p:X \times_k $X$ \ra X/E \times_k X/E$ is quotient, 
$\Delta_{X/E}$ is closed in $X/E \times_k X/E$ iff 
$(p \times_k p)^{-1}(\Delta_{X/E}) = E$ is closed in $X \times_k X$.  
Consequently, if $A \subset X$ is closed, then $X/A$ is in $\dcg$.

[Note: \  Recall that if $X$ is a topological space, then for any equivalence relation $E$ on $X$, 
$X/E$ Hausdorff $\implies E \subset X \times X$ closed and $E \subset X \times X$ closed plus 
$p:X \ra X/E$ open $\implies X/E$ Hausdorff.]

$\dcg$, like \bCG and \bCGH, is cartesian closed.  
For $\dcg$ has finite products and if $X$ is in \bCG and if $Y$ is in $\dcg$, 
then $kC_k(X,Y)$ is in $\dcg$.

[Note: \ Suppose that $B$ is $\Delta$-separated $-$then $\bCG/B$ is cartesian closed 
(Booth-Brown\footnote[2]{\textit{General Topology Appl.} \textbf{8} (1978), 181-195.}).]\\

\begingroup%%----------------------------------->>
\fontsize{9pt}{11pt}\selectfont
$\bCG_*$ 
\index{$\bCG_*$} %Symbol
and $\dcg_*$  
\index{$\bDelta$-\bCG$_*$} %Symbol
are the pointed versions of \textbf{CG} and $\dcg$.  Both are closed categories.

[Note: \  The 
\un{pointed exponential object}
\index{pointed exponential object} 
$Z^Y$ is $\text{hom}(Y,Z)$.]\\ 
\endgroup%%------------------------------------<<

\begingroup%%----------------------------------->>
\fontsize{9pt}{11pt}\selectfont
\textbf{\small EXAMPLE}  \ 
Let $X$ be a nonnormal LCH space.  Fix nonempty disjoint closed subsets $A$ and $B$ of $X$ that do not have disjoint neighborhoods $-$then $X/A$ and $X/B$ are compactly generated Hausdorff spaces but neither $X/A$ nor $X/B$ is regular.  Put $E = A \times A \cup B \times B \cup \Delta_X$.  The quotient $X/E$ is a $\Delta$-separated compactly generated space which is not Hausdorff.  
Moreover, $X/E$ is not the continuous image of any compact Hausdorff space.
\vspi
[Note: \  Take for $X$ the Tychonoff plank.  Let $A = \{(\Omega,n): 0 \leq n < \omega\}$ and $B = \{(\alpha,\omega): 0 \leq \alpha < \Omega\}$ $-$then $X/E$ is compact and all its compact subspaces are closed.  By comparison, the product $X/E \times X/E$, while compact, has compact subspaces that are not closed.]\\
\endgroup%%------------------------------------<<

%%----------------------------------------------------------------------------------------------36
\index{k-Spaces}
\begingroup%%----------------------------------->>
\fontsize{9pt}{11pt}\selectfont
\textbf{\small EXAMPLE \  (\un{$k$-Spaces})} \ 
The monocoreflective hull in \bTOP of the category of compact spaces is the category of $k$-spaces.  In other words, a topological space $X$ is a \un{$k$-space} provided that a subset $U$ of $X$ is open iff $U \cap K $ is open in $K$ for every compact subset $K$ of $X$.
Every compactly generated space is a $k$-space.  
The converse is false: Let $X$ be the subspace of $[0,\Omega]$ obtained by deleting all limit ordinals except $\Omega$ 
$-$then $X$ is not discrete.  Still, the only compact subsets of $X$ are the finite sets, thus $kX$ is discrete. 
The one point compactification $X_\infty$ of $X$ is compact and contains $X$ as an open subspace.  
Therefore $X_\infty$ is not compactly generated but is a $k$-space (being compact).  
The category of $k$-spaces is similar in many respects to the category of compactly generated spaces.  
However, there is one major difference: It is not cartesisan closed 
(\u Cin\u cura\footnote[3]{\textit{Topology Appl.} \textbf{41} (1991), 205-212.}).
\vspi
[Note: \  If $\sK$ is the class of compact spaces, then $\textbf{HAUS} \subset$ $\bDelta$-\bK and the inclusion is strict.  Reason: A topological space $X$ is in $\bDelta$-\bK iff every compact subspace of $X$ is Hausdorff.]\\
\endgroup%%------------------------------------<<

\label{14.5}
\label{14.69c}
\label{14.157}
\begingroup%%----------------------------------->>
\fontsize{9pt}{11pt}\selectfont
\textbf{\small FACT} \ 
Let $X^0 \subset X^1 \subset \cdots$ be an expanding sequence of topological spaces.  
Assume: $\forall \ n$, $X^n$ is in $\bDelta$-\bCG and is a closed subspace of $X^{n+1}$ $-$then $X^\infty$ is in $\bDelta$-\bCG.
\vspi
[That $X^\infty$ is in \bCG is automatic.  Let $K$ be a compact Hausdorff space; let $\phi \in C(K,X^\infty)$ $-$then, from the lemma on 
p. \pageref{1.24}, 
$\phi(K) \subset X^n$ $(\exists \ n)$ $\implies$ $\phi(K)$ is closed in $X^n$ $\implies$ $\phi(K)$ is closed in $X^\infty$.]\\
\endgroup%%------------------------------------<<

\index{weak products (example)}
\begingroup%%----------------------------------->>
\fontsize{9pt}{11pt}\selectfont
\textbf{\small EXAMPLE \  (\un{Weak Products})} \ 
Let $(X_0,x_0), (X_1,x_1), \ldots $ be a sequence of pointed spaces in $\dcg_*$.  
Put $X^n = X_0 \times_k \cdots \times_k X_n$ $-$then $X^n$ is in $\dcg_*$ with base point $(x_0, \ldots, x_n)$.  
The pointed map $X^n \ra X^{n+1}$ is a closed embedding.  
One writes $(\omega) \ds\prod\limits_1^\infty X_n$ in place of $X^\infty$ and calls it the 
\un{weak product}
\index{weak product} 
of the $X_n$.  
By the above, $(\omega) \ds\prod\limits_1^\infty X_n$ is in $\dcg_*$ 
(the base point is the infinite string made up of the $x_n$).
\vspi
\label{4.58}
[Note: \  The same construction can be carried out in \bTOP, the only difference being that $X^n$ is the ordinary product of $X_0, \ldots, X_n$.]\\
\endgroup%%------------------------------------<<

Every Hausdorff topological group is completely regular.  In particular, every Hausdorff topological vector space is completely regular.  Every Hausdorff locally compact topological group is paracompact.

[Note: \  Every topological group which satisfies the $\mT_0$ separation axiom is necessarily a CRH space.]\\

\begingroup%%----------------------------------->>
\fontsize{9pt}{11pt}\selectfont
\label{5.0ak}
\textbf{\small EXAMPLE}\label{2.10} \quadx
Take $G = \R^\kappa$ $(\kappa > \omega)$ $-$then $G$ is a Hausdorff topological group but $G$ is not compactly generated.  Consider $kG$: Inversion $kG \ra kG$ is continuous, as is multiplication $kG \times_k kG \ra kG$.  But $kG$ is not a topological group, i.e., multiplication $kG \times kG \ra kG$ is not continuous.  
In fact, $kG$, while Hausdorff, is not regular.
\\[-.2cm]
\endgroup%%------------------------------------<<

%%----------------------------------------------------------------------------------------------37
\begingroup%%----------------------------------->>
\fontsize{9pt}{11pt}\selectfont
Let $E$ be a normed linear space; 
let $E^*$ be its dual, i.e., the space of continuous linear functionals on $E$ $-$then $E^*$ is also a normed linear space.  
The elements of $E$ can be regarded as scalar valued functions on $E^*$.  
The initial topology on $E^*$ determined by them is called the 
\un{weak$^*$ topology}.
\index{weak$^*$ topology}  
It is the topology of pointwise convergence.  
In the weak$^*$ topology, $E^*$ is a Hausdorff topological vector space, thus is completely regular.  
If $\dim E \geq \omega$, then every nonempty weak$^*$ open set in $E^*$ is unbounded in norm.  
By contrast, Alaoglu's theorem says that the closed unit ball in $E^*$ is compact in the weak$^*$ topology (and second countable if $E$ is separable).  
However, the weak$^*$ topology is metrizable iff $\dim E \leq \omega$.
\vspi
[Note: \  Let $E$ be a vector space over $\R$ $-$then 
Kruse\footnote[2]{\textit{Math. Zeit.} \textbf{83} (1964), 314-320.} 
has shown that $E$ admits a complete norm (so that $E$ is a Banach space) iff 
$\dim E < \omega$ or $(\dim E)^\omega = \dim E$.  
Therefore, the weak$^*$ topology on the dual of an infinite dimensional Banach space is not metrizable.]\\
\endgroup%%------------------------------------<<

\label{14.69}
The forgetful functor from the category of topological groups to the category of topological spaces (pointed topological spaces) has a left adjoint 
$X \ra F_{gr}X$ $((X,x_0) \ra F_{gr}(X,x_0))$, where 
$F_{gr}X$ $(F_{gr}(X,x_0))$ is the 
\un{free topological group}
\index{free topological group} 
on $X((X,x_0))$.
Algebraically, $F_{gr}X$ $(F_{gr}(X,x_0))$ is the free group on $X$ $(X - \{x_0\})$.  
Topologically, $F_{gr}X$ $(F_{gr}(X,x_0))$ carries the finest topology compatible with the group structure for which the canonical injection $X \ra F_{gr}X$ 
$((X,x_0) \ra (F_{gr}(X,x_0))$ is continuous.  There is a commutative triangle
\begin{tikzcd}%[ sep=small]
X \ar{dr} \ar{r} &{F_{gr}X} \ar{d}\\
&{F_{gr}(X,x_0)}
\end{tikzcd}
and $F_{gr}(X,x_0) \approx F_{gr}X/\langle x_0\rangle$ ($\langle x_0\rangle$ 
the normal subgroup generated by the word $x_0$).  
On the other hand, 
$F_{gr}X \approx F_{gr}(X,x_0) \coprod \Z$ ($\coprod$ the coproduct in the category of topological groups) and, of course, 
$F_{gr}X \approx (F_{gr}(X \coprod *,*)$.

[Note: \  The arrow of adjunction $X \ra F_{gr}X$  $((X,x_0) \ra F_{gr}(X,x_0))$ is an embedding iff 
$X$ is completely regular and is a closed embedding iff $X$ is completely regular + Hausdorff 
(Thomas\footnote[2]{\textit{General Topology Appl.} \textbf{4} (1974), 51-72; 
see also \textit{Quaestiones Math.} 2 (1977), 355-377.}).]\\

\textbf{\small LEMMA} \ \ 
If $X$ is a compact Hausdorff space, then $F_{gr}X$ $(F_{gr}(X,x_0))$ is a Hausdorff topological group.\\

Application:  If $X$ is a CRH space, then $F_{gr}(X)$ $(F_{gr}(X,x_0))$ is a Hausdorff topological group.

[Consider 
$X \ra F_{gr}(\beta X)$  
$((X,x_0) \ra F_{gr}(\beta X,\beta x_0))$.]\\

%%----------------------------------------------------------------------------------------------38

\begingroup%%----------------------------------->>
\fontsize{9pt}{11pt}\selectfont
\textbf{\small EXAMPLE}  \ 
It is easy to construct nonnormal Hausdorff topological groups.  
Thus, given a topological space $X$, let $F_{gr}X$ be the free topological group on $X$ $-$then, for $X$ a CRH space, the arrow $X \ra F_{gr}X$ is a closed embedding and $F_{gr}X$ is a Hausdorff topological group, so $X$ not normal $\implies$ $F_{gr}X$  not normal.\\
\endgroup%%------------------------------------<<

\begingroup%%----------------------------------->>
\fontsize{9pt}{11pt}\selectfont
\textbf{\small FACT} \ 
Given a topological space $X$, $F_{gr}(X,x_0^\prime) \approx F_{gr}(X,x_0^{\prime\prime})$ $\forall \ x_0^\prime, x_0^{\prime\prime} \in X$.
\vspi
[Let 
$\mu^{\prime}:(X,x_0^\prime) \ra F_{gr}(X,x_0^{\prime})$, 
$\mu^{\prime\prime}:(X,x_0^{\prime\prime}) \ra F_{gr}(X,x_0^{\prime\prime})$
be the arrows of adjunction and consider the pointed continuous functions
$f^\prime :(X,x_0^\prime) \ra F_{gr}(X,x_0^{\prime\prime})$,
$f^{\prime\prime} :(X,x_0^{\prime\prime}) \ra F_{gr}(X,x_0^{\prime})$,
defined by 
$f^\prime(x) = \mu^{\prime\prime}(x) \mu^{\prime\prime}(x_0^\prime)^{-1}$,
$f^{\prime\prime}(x) = \mu^{\prime}(x) \mu^{\prime}(x_0^{\prime\prime})^{-1}$.]\\
\endgroup%%------------------------------------<<

The forgetful functor from the category of abelian topological groups to the category of topological spaces (pointed topological spaces) has a left adjoint 
$X \ra F_{ab}X$ $((X,x_0) \ra F_{ab}(X,x_0))$ and when given the quotient topology,
$F_{gr}X/[F_{gr}X,F_{gr}X] \approx F_{ab}X$ 
$(F_{gr}(X,x_0)/[F_{gr}(X,x_0),F_{gr}(X,x_0)] \approx F_{ab}(X,x_0))$.\\

%%%%%%%%%%%%%%%%%%%%%%%%%%%%%%%%%%%%%%
%%%%%%%%%%%%%%%%%%%%%%%%%%%%%%%%%%%%%%
%%%%%%%%%%%%%%%%%%%%%%%%%%%%%%%%%%%%%%

\begin{center}
$\S \ 1$
\\[0.5cm]
$\mathcal{REFERENCES}$\\
\end{center}

\[
\textbf{BOOKS}
\]

\begingroup
\fontsize{9pt}{11pt}\selectfont
\setlength\parindent{0 cm}

[1] \quad Brown, R., \textit{Topology}, Ellis Horwood (1988).
\\[-.2cm]

[2] \quad Dugundji, J., \textit{Topology}, Allyn and Bacon (1966).
\\[-.2cm]

[3] \quad Engelking, R., \textit{General Topology}, Heldermann Verlag (1989).
\\[-.2cm]

[4] \quad Fremlin, D., \textit{Consequences of Martin's Axiom}, Cambridge University Press (1984).
\\[-.2cm]

[5] \quad Hocking, J. and Young, G., \textit{Topology}, Dover (1988).
\\[-.2cm]

[6] \quad James, I., \textit{General Topology and Homotopy Theory}, Springer Verlag (1984).
\\[-.2cm]

[7] \quad Jech, T., \textit{Set Theory}, Academic Press (1978).
\\[-.2cm]

[8] \quad Kunen, K., \textit{Set Theory}, North Holland (1980).
\\[-.2cm]

[9] \quad Nabor, G., \textit{Set Theoretic Topology}, University Microfilms International (1977).
\\[-.2cm]

[10] \quad Nagata, J., \textit{Modern General Topology}, North Holland (1985).
\\[-.2cm]

[11] \quad Rudin, M., \textit{Lectures on Set Theoretic Topology}, American Mathematical Society (1975).
\\[-.2cm]

[12] \quad Sieradski, A., \textit{An Introduction to Topology and Homotopy}, PWS-Kent (1992).
\\[-.2cm]

[13] \quad Steen, L. and Seeback, J., \textit{Counterexamples in Topology}, Springer Verlag (1978).
\\[-.2cm]

[14] \quad Tall, F. (ed.), \textit{The Work of Mary Ellen Rudin}, New York Academy of Sciences (1993).
\\[-.2cm]

[15] \quad Willard, S., \textit{General Topology}, Addison Wesley (1970).
\\
\endgroup

\[
\textbf{ARTICLES}
\]

\begingroup
\fontsize{9pt}{11pt}\selectfont
\setlength\parindent{0 cm}

[1] \quad Arhangel'ski\u i, A. and Fedorchuk, V., The Basic Concepts and Constructions of General Topology, In: 

\hspace{.65cm} \textit{General Topology}, EMS \textbf{17}, Springer Verlag (1990), 1-90.
\\[-.2cm]

[2] \quad Arhangel'ski\u i, A., Compactness, In: \textit{General Topology}, EMS \textbf{50}, Springer Verlag (1996), 1-117.
\\[-.2cm]

[3] \quad Arhangel'ski\u i, A., Paracompactness and Metrization, In: \textit{General Topology}, EMS \textbf{51}, Springer Verlag 

\hspace{.65cm} (1995), 1-70.
\\[-.2cm]

[4] \quad Balogh, Z., On Collectionwise Normality of Locally Compact, Normal Spaces, \textit{Trans Amer. Math.} 

\hspace{.65cm} \textit{Soc.} \textbf{323} (1991), 389-411.
\\[-.2cm]

[5] \quad Burke, D., Covering Properties, In: \textit{Handbook of Set Theoretic Topology}, K. Kunen and J. Vaughan 

\hspace{.65cm} (ed.), North Holland (1984), 347-422.
\\[-.2cm]

[6] \quad Comfort, W., Deciding Some Undecidable Topological Statements, \textit{Ann. New York Acad. Sci.} \textbf{321} 

\hspace{.65cm} (1979), 9-26.
\\[-.2cm]

[7] \quad van Douwen, E., The Integers and Topology, In: op. cit. [5], 111-167.
\\[-.2cm]

[8] \quad Dow, A., Set Theory in Topology, In: \textit{Recent Progress in General Topology}, M. Hu\u sek and J. van Mill 

\hspace{.65cm} (ed.), North Holland (1992), 167-197.
\\[-.2cm]

[9] \quad Fleissner, W., The Normal Moore Space Conjecture and Large Cardinals, In: op. cit. [5], 733-760.
\\[-.2cm]

[10] \quad Nyikos, P., The Theory of Nonmetrizable Manifolds, In: op. cit. [5], 633-684.
\\[-.2cm]

[11] \quad Rudin, M., Dowker Spaces, In: op. cit. [5], 761-780.
\\[-.2cm]

[12] \quad Tall, F., Normality versus Collectionwise Normality, In: op. cit. [5], 685-732.
\\[-.2cm]

[13] \quad Watson, S., The Construction of Topological Spaces, In: op. cit. [8], 673-757.
\\[-.2cm]

[14] \quad Yasui, Y., Generalized Paracompactness, In: \textit{Topics in General Topology}, K. Morita and J. 
Nagata 

\hspace{.8cm} (ed.), North Holland (1989), 161-202.

\setlength\parindent{2em}

\endgroup

\chapter{
$\boldsymbol{\S}$\textbf{2}.\quadx  CONTINUOUS FUNCTIONS}
\setlength\parindent{2em}
\setcounter{proposition}{0}

%%----------------------------------------------------------------------------------------------01
$\text{ }$\\[-1.25cm]

Apart from an important preliminary, namely a characterization of the exponential objects in \bTOP, the emphasis in this $\S$ is on the properties possessed by $C(X)$, where $X$ is a CRH space.

A topological space $Y$ is said to be \underline{cartesian} if the functor \ $-\times Y:\bTOP \ra \bTOP$ has a right adjoint $Z \ra Z^Y$.  Example:  A LCH space is cartesian.\\

\begin{proposition} \ 
A topological space $Y$ is cartesian iff $- \times Y$ preserves colimits 
(cf. p. \pageref{2.1})
 or equivalently, iff $- \times Y$ preserves coproducts and coequalizers.
\end{proposition}

[Note: \  The preservation of coproducts is automatic and the preservation of coequalizers reduces to whether $- \times Y$ takes quotient maps to quotient maps.]\\

Notation:  Given topological spaces $X, Y, Z, \ \Lambda: F(X \times Y, Z) \ra F(X,F(Y,Z))$ is the bijection defined by the rule 
$\Lambda(f)(x)(y) = f(x,y)$.

Let $\tau$ be a topology on $C(Y,Z)$ $-$then $\tau$ is said to be 
\un{splitting}
\index{splitting} 
if 
$\forall \ X$, $f \in C(X \times Y,Z)$ $\implies$ $\Lambda (f) \in C(X,C(Y,Z))$ and $\tau$ is said to be 
\un{cosplitting}
\index{cosplitting} 
if $\forall \ X$, $g \in C(X,C(Y,Z))$ $\implies$ $\Lambda^{-1}(g) \in C(X \times Y,Z)$.\\

\textbf{\small LEMMA} \ 
If $\tau^\prime$ is a splitting topology on $C(Y,Z)$ and $\tau^{\prime\prime}$ is a cosplitting topology on $C(Y,Z)$, the $\tau^\prime \subset \tau^{\prime\prime}$.\\

Application: $C(Y,Z)$  admits at most one topology which is simultaneously splitting and cosplitting, the 
\un{exponential topology}.
\index{exponential topology}\\

\begingroup
\fontsize{9pt}{11pt}\selectfont
\textbf{\small EXAMPLE} \ 
$\forall \ Y \& \ \forall \ Z$, the compact open topology on $C(Y,Z)$ is splitting.\\
\endgroup

\begingroup
\fontsize{9pt}{11pt}\selectfont
\textbf{\small EXAMPLE} \ 
If $Y$is locally compact, then $\forall \ Z$ the exponential topology on $C(Y,Z)$ exists and is the compact open topology.

[Note: \  A topological space $Y$ is said to be 
\un{locally compact}
\index{locally compact} 
if $\forall$ open set $P$ and $\forall \  y \in P$, there exists a compact set $K \subset P$ with $y \in \itrx K$.  
Example:  The one point compactification $\Q_\infty$ of $\Q$ is compact but not locally compact.]\\
\endgroup

\begingroup
\fontsize{9pt}{11pt}\selectfont
\textbf{\small FACT} \ 
Let $Y$ be a locally compact space $-$then for all $X$ and $Z$, the operation of composition $C(X,Y) \times C(Y,Z) \ra C(X,Z)$ is continuous if the function spaces carry the compact open topology.\\
\endgroup

\begin{proposition} \ %02
A topological space $Y$ is cartesian iff the exponential topology on $C(Y,Z)$  exists for all Z.\\
\end{proposition}

%%----------------------------------------------------------------------------------------------02
\begingroup
\fontsize{9pt}{11pt}\selectfont
\textbf{\small EXAMPLE} \ 
A locally compact space is cartesian.\\
\endgroup

\begingroup
\fontsize{9pt}{11pt}\selectfont
\textbf{\small FACT} \ 
Suppose that $Y$ is cartesian.  
Assume: $\forall \ Z$, the exponential topology on $C(Y,Z)$ is the compact open topology $-$then  $Y$ is locally compact.\\
\endgroup

Let $Y$ be a topological space, $\tau_Y$ its topology $-$then the open sets in the 
\un{continuous} \un{topology}
\index{continuous topology} 
on $\tau_Y$ 
are those collections $\sV \subset \tau_Y$ such that 
(1) $V \in \sV$, $V^\prime \in \tau_Y$ $\implies$ 
$V^\prime \in \sV$ if $V \subset V^\prime$ and 
(2) $V_i \in \tau_Y$ $(i \in I), \  \ds\bigcup\limits_i V_i \in \sV$ $\implies$ 
$\exists \ i_1, \ldots, i_n:$ $V_{i_1} \cup \ldots \cup V_{i_n} \in \sV$.\\

\label{2.2}
\textbf{\small LEMMA} \  Let $f \in F(X,\tau_Y)$, where $X$ is a topological space and $\tau_Y$ has the continuous topology $-$then $f$ is continuous if $\{(x,y): y \in f(x)\}$ is open in $X \times Y$.\\

\begingroup
\fontsize{9pt}{11pt}\selectfont
Let T $= \{(P,y): y \in P\} \subset \tau_Y \times Y$ $-$then a topology on $\tau_Y$ is said to have 
\un{property T}
\index{property T} 
if T is open in $\tau_Y \times Y$.  
Example:  The discrete topology on $\tau_Y$ has property T.\\

\textbf{\small FACT} \ 
The continuous topology $\tau_Y$ is the largest topology in the collection of all topologies on $\tau_Y$ that are smaller than every topology on $\tau_Y$ which has property T.

[If $\tau_Y(\tT)$ is $\tau_Y$ in a topology having property T, then by the lemma, the identity function $\tau_Y(\tT) \ra \tau_Y$ is continuous if $\tau_Y$ has the continuous topology.]\\
\endgroup

\label{3.19}
Let $Y$ be a topological space $-$then $Y$ is said to be 
\underline{core compact}
\index{core compact} 
if $\forall$ open set $P$ and $\forall \ y \in P$, there exists an open set $V \subset P$ with $y \in V$ such that every open covering of $P$ contains a finite covering of $V$.  
Example:  A locally compact space is core compact.\\

\begingroup
\fontsize{9pt}{11pt}\selectfont
There exists a core compact space with the property that every compact subset has an empty interior 
(Hofman$-$Lawson\footnote[2]{\textit{Trans. Amer. Math. Soc.} \textbf{246} (1978), 285-310 (cf. 304-306).}).\\

\textbf{\small FACT} \ 
Equip $\tau_Y$ with the continuous topology $-$then $Y$ is core compact iff $\forall$ open set $P$ and $\forall \ y \in P$, 
there exists an open $\sV \subset \tau_Y$ such that $P \in \sV$ and $y \in \itr \cap \sV$.\\
\endgroup

\begingroup
\fontsize{9pt}{11pt}\selectfont
\label{2.3}
\textbf{\small EXAMPLE} \ 
A topological space $Y$ is core compact iff the continuous topology on $\tau_Y$ has property T.
\endgroup

Let $Y,Z$ be topological spaces $-$then the 
\un{Isbell topology}
\index{Isbell topology}
 on $C(Y,Z)$ is the initial topology on $C(Y,Z)$ determined by the 
$
e_Q : 
\begin{cases}
\ C(Y,Z) \ra \tau_Y \\
\ f \ra f^{-1}(Q)
\end{cases}
(Q \in \tau_Z), \ 
$
where $\tau_Y$ has the 
%%----------------------------------------------------------------------------------------------03
continuous topology.  
Notation: is$C(Y,Z)$.  
%dmc\symbindex{is$C(Y,Z)$}
Examples:  (1) is$C(Y,[0,1]/[0,1[) \approx \tau_Y$; 
(2) is$C(*,Z) \approx Z$.\\

\textbf{\small LEMMA} \  
The compact open topology on $C(Y,Z)$ is smaller than the Isbell topology.\\

\begingroup
\fontsize{9pt}{11pt}\selectfont
\textbf{\small EXAMPLE} \ 
$\forall \ Y$ $\&$ $\forall \ Z$, the Isbell topology on $C(Y,Z)$ is splitting.

[Fix an $f \in C(X \times Y,Z)$ and let $g \in \Lambda(f)$ 
$-$then the claim is that $g \in C(X,\text{is}C(Y,Z))$.  
From the definitions, this amounts to showing that $\forall \ Q \in \tau_Z$, $e_Q \circ g$ is continuous.  
Write 
$f^{-1}(Q)$ as a union of rectangles $R_i = U_i \times V_i \subset X \times Y$.  
Take an $x \in X$ and consider any 
$\sV$: $e_Q(g(x)) \in \sV$.  
Since $e_Q(g(x)) = \ds\bigcup\limits_i \{y:(x,y) \in R_i\}$, $\exists$ $i_k$ $(k = 1, \ldots, n)$ : 
$\ds\bigcup\limits_{k = 1}^n \{y:(x,y) \in R_{i_k}\} \in \sV$, so 
$\forall \ u \in \ds\bigcap\limits_{k = 1}^n U_{i_k}$, 
$e_Q(g(u)) \in \sV$.]\\

\textbf{\small FACT} \ 
Let $Y$ be a core compact space $-$then for all $X$ and $Z$, the operation of composition $C(X,Y) \times C(Y,Z) \ra C(X,Z)$ is continuous if the function spaces carry the Isbell topology.\\
\endgroup

\begin{proposition} \ %03
Let $Y$ be a topological space $-$then $Y$ is cartesian iff $Y$ is core compact.
\end{proposition}

[Necessity: \  Let $\tau_i$ run through the topologies on $\tau_Y$ which have property T and put $X_i = (\tau_Y,\tau_i)$.  
Form the coproduct $X = \ds\coprod\limits_i X_i$ and let $f:X \ra \tau_Y$ be the function whose restriction to each 
$X_i$ is the identity, where $\tau_Y$ carries the continuous topology $-$then $f$ is a quotient map 
(cf. p. \pageref{2.1}).  
Since $Y$ is cartesian, it follows from Proposition 1 that $f \times id_Y:X \times Y \ra \tau_Y \times Y$ is also quotient.  
But $X \times Y \approx \ds\coprod\limits_i X_i \times Y$ and, by hypothesis, T is open in $X_i \times Y$ $\forall \ i$.  Therefore T must be open in $\tau_Y \times Y$ as well, i.e., the continuous topology on $\tau_Y$ has property T, thus $Y$ is core compact (cf. p. 
\pageref{2.3}
).

Sufficiency: \  As has been noted above, the Isbell topology on $C(Y,Z)$ is splitting, so to prove that $Y$ is cartesian it suffices to prove that the Isbell topology on $C(Y,Z)$ is cosplitting when $Y$ is core compact (cf. Proposition 2).  
Fix $g \in C(X,\text{is}C(Y,Z))$ and put $f = \Lambda^{-1}(g)$.  
Given a point $(x,y) \in X \times Y$, let $Q$ be an open subset of $Z$ such that $f(x,y) \in Q$.  
Choose an open $P \subset Y : y \in P$ $\&$ $f(\{x\} \times P) \subset Q$.  
Because $Y$ is core compact, there exists an open $\sV \subset \tau_Y: P \in \sV$ and $y \in \itr \cap \sV$.  
But $e_Q(g(x)) \supset P$ $\implies$ $e_Q(g(x)) \in \sV$ and, from the continuity of $e_Q \circ g$, $\exists$ 
a neighborhood $O$ of $x: e_P(g(O)) \subset \sV$, hence $f(O \times \itr \cap  \sV) \subset Q$.]\\

Remark:  Suppose that $Y$ is core compact $-$then $\forall \ Z$, ``the'' exponential object $Z^Y$ is $\text{is}C(Y,Z)$, the exponential topology on $C(Y,Z)$ being the Isbell topology.

[Note: \  The Isbell topology and the compact open topology on $C(Y,Z)$ are one and the same if $Y$ is locally compact.\\

%%----------------------------------------------------------------------------------------------04
\begingroup
\fontsize{9pt}{11pt}\selectfont
\textbf{\small FACT} \
Let $f,g \in C(Y,Z)$.  Assume: $f,g$ are homotopic $-$then $f,g$ belong to the same path component of is$C(Y,Z)$.\\

\textbf{\small FACT} \ 
Let $f,g \in C(Y,Z)$.  Assume: $f,g$  belong to the same path component of is$C(Y,Z)$ $-$then $f,g$ are homotopic if $Y$ is core compact.\\
\endgroup

What follows is a review of the elementary properties possessed by $C(X,Y)$ when equipped with the compact open topology (omitted proofs can be found in 
Engelking\footnote[2]{\textit{General Topology}, Heldermann Verlag (1989).}).

Notation:  Given Hausdorff spaces $X$ and $Y$, let co$C(X,Y)$ stand for $C(X,Y)$ in the compact open topology.

[Note: \  The point open topology on $C(X,Y)$ is smaller than the compact open topology.  Therefore co$C(X,Y)$ is necessarily Hausdorff.  Of course, if $X$ is discrete, then ``point open'' = ``compact open''.]\\

\begin{proposition} \ %04
Suppose that $Y$ is regular $-$then  co$C(X,Y)$ is regular.\\
\end{proposition}

\begin{proposition} \ %05
Suppose that $Y$ is completely regular $-$then co$C(X,Y)$ is completely regular.\\
\end{proposition}

\begingroup
\fontsize{9pt}{11pt}\selectfont
\textbf{\small EXAMPLE} \ 
It is false that $Y$ normal $\implies $co$C(X,Y)$ normal.  Thus take $X = \{0,1\}$ (discrete topology) $-$then co$C(\{0,1\},Y) \approx Y \times Y$ and there exists a normal Hausdorff space $Y$ whose square is not normal (e.g., the Sorgenfrey line 
(cf. p. 
\pageref{2.4}
)).\\

O'Meara\footnote[3]{\textit{Proc. Amer. Math. Soc.} \textbf{29} (1971), 183-189.}
has shown that if $X$ is a second countable metrizable space and $Y$ is a metrizable space, then co$C(X,Y)$ is perfectly normal and hereditarily paracompact.\\

\textbf{\small EXAMPLE} \ 
The loop space $\Omega Y$ of a pointed metrizable space $(Y,y_0)$ is paracompact.\\
\endgroup

A Hausdorff space $X$ is said to be 
\un{countable at infinity}
\index{countable at infinity} 
if there is a sequence $\{K_n\}$ of compact subsets of $X$ such that if $K$ is any compact subset of $X$, then $K \subset K_n$ for some $n$.  
Example:  A LCH space is countable at infinity iff it is $\sigma$-compact.

[Note: \ $X$ countable at infinity $\implies X$ $\sigma$-compact.  
Example: $\PP$ is not $\sigma$-compact, hence it is not countable at infinity.]\\

\begingroup
\fontsize{9pt}{11pt}\selectfont
\textbf{\small FACT} \ 
Suppose that $X$ is countable at infinity.  Assume: $X$ is first countable $-$then $X$ is locally compact.\\

%%----------------------------------------------------------------------------------------------05
\textbf{\small EXAMPLE} \ 
$\Q$ is $\sigma -$compact but $\Q$ is not countable at infinity.\\

\textbf{\small EXAMPLE} \ 
Fix a point $x \in \beta\N - \N$ $-$then $X = \N \cup \{x\}$, viewed as a subspace of $\beta\N$, is countable at infinity but it is not first countable.

[Note: \  The compact subsets of $X$ are finite.  However $X$ is not compactly generated.]\\

\textbf{\small EXAMPLE} \ 
Let $E$ be an infinite dimensional Banach space $-$then $E^*$ in the weak$^*$ topology is countable at infinity.\\
\endgroup

\begin{proposition} \ %06
Suppose that $X$ is countable at infinity $-$then for every metrizable $Y$, co$C(X,Y)$ is metrizable.\\
\end{proposition}

\begin{proposition} \ %07
Suppose that $X$ is countable at infinity and compactly generated $-$then for every completely metrizable $Y$, co$C(X,Y)$ is completely metrizable.\\
\end{proposition}

Notation: Given a topological space $X$, write $H(X)$ for its set of homeomorphisms $-$then $H(X)$ is a group under composition.

Let us assume that $X$ is a LCH space.  Endow $H(X)$ with the compact open topology. \\
Question:  Is $H(X)$ thus topologized a topological group?  In general, the answer is ``no'' (cf. infra) but there are situations in which the answer is ``yes''.

[Note: \  The composition
$
\begin{cases}
\ H(X) \times H(X) \ra H(X) \\
\ (f,g) \mapsto g\circ f
\end{cases}
$
is continuous, so the problem is whether the inversion $f \ra f^{-1}$ is continuous.]

Remark: The evaluation 
$
\begin{cases}
\ H(X) \times X \ra X  \\ %dmcerrxxx note \ra H(X) missing in orig text
\ (f,x) \mapsto f(x)
\end{cases}
$
is continuous.

Given subsets $A$ and $B$ of $X$, put $\langle A,B\rangle = \{f \in H(X): f(A) \subset f(B)\}$ 
$-$then by definition, the collection $\{\langle K,U\rangle\}$ ($K$ compact and $U$ open) is a subbasis for the compact open topology on $H(X)$.\\

\begin{proposition} \ %08
If $X$ is a compact Hausdorff space, then $H(X)$ is a topological group in the compact open topology.
\end{proposition}

[For $f \in \langle K,U\rangle \Leftrightarrow f^{-1} \in \langle X - U,X - K\rangle$.]\\

\textbf{\small FACT} \ 
If $X$ is a compact metric space, then $H(X)$ is completely metrizable.\\

\textbf{\small LEMMA} \ 
Let $X$ be a locally connected LCH space $-$then the collection $\{\langle L,V\rangle\}$, where $L$ is compact $\&$ connected with $L \neq \emptyset$ and $V$ is open, constitute a subbasis for the compact open topology on $H(X)$.\\

%%----------------------------------------------------------------------------------------------06
\label{6.12}
\begin{proposition} \ %09
If $X$ is a locally connected LCH space $-$then $H(X)$ is a topological group in the compact open topology.
\end{proposition}

[Fix an $f \in H(X)$ and choose $\langle L,V\rangle$ per the lemma: $f^{-1} \in \langle L,V\rangle$.  
Determine relatively compact $O \ \& \ P$: 
$f^{-1}(L) \subset O \subset \ov{O}$ $\subset P \subset \ov{P}$ $\subset V$ 
$(\implies$ 
$f((X - O) \cap \ov{P})$ $\subset (X - L) \cap f(V))$.  
Let $x$ be any point such that $f(x) \in \itrx L$ $-$then 
$\langle \{x\}, \itrx L\rangle \ \cap$ 
$\langle(X - O) \cap \  \ov{P},(X - L) \ \cap \ f(V)\rangle$ 
is a neighborhood of $f$ in $H(X)$, call it $H_f$.  
Claim: $g \in H_f$ $\implies$ $g^{-1} \in \langle L,V\rangle$.  
To check this, note that 
$g((X - O) \cap \ \ov{P}) \subset$ $(X - L) \cap f(V)$ $\implies$ 
$L \cup (X - f(V)) \subset$ $g(O) \cup g(X - \ov{P})$.  
But $g(O)$, $g(X - \ov{P})$ are nonempty disjoint open sets, so \mL is contained in either $g(O)$ or $g(X - \ov{P})$ 
(\mL being connected).  
Since the containment $L \subset g(X - \ov{P})$ is impossible ($g(x) \in \itrx L$ and $x \notin X - \ov{P})$, 
it follows that $L \subset g(O)$ or still, $g^{-1}(L) \subset O \subset V$, i.e., $g^{-1} \in \langle L,V\rangle$.  
Therefore inversion is a continuous function.]\\

Application: The homeomorphism group of a topological manifold is a topological group in the compact open topology.\\

\begingroup%%----------------------------------->>
\fontsize{9pt}{11pt}\selectfont
\textbf{\small EXAMPLE} \  
Let $X = \{0,2^n (n \in \Z)\}$ $-$then in the induced topology from $\R$, \mX is a LCH space but $H(X)$ in the compact open topology is not a topological group.\\
\endgroup %%------------------------------------<<

Suppose that \mX is a LCH space, $X_\infty$ its one point compactification 
$-$then $H(X)$ can be identified with the subgroup of $H(X_\infty)$ consisting of those homeomorphisms $X_\infty \ra X_\infty$ which leave $\infty$ fixed.  
In the compact open topology, $H(X_\infty)$ is a topological group (cf. Proposition 8).  
Therefore $H(X)$ is a topological group in the induced topology.  As such, $H(X)$ is a closed subgroup of $H(X_\infty)$.

[Note: This topology on $H(X)$ is the complemented compact open topology.  
It has for a subbasis all sets of the form 
$\langle K,U \rangle$, where \mK is compact and \mU is open, as well as all sets of the form 
$\langle X -V,X - L \rangle$, where \mV is open and $L$ is compact.]\\

\begingroup%%----------------------------------->>
\fontsize{9pt}{11pt}\selectfont
An 
\un{isotopy} 
\index{isotopy} 
of a topological space \mX is a collection $\{h_t: 0 \leq t \leq 1\}$ of homeomorphisms of \mX such that 
$
\begin{cases}
\ h:X \times [0,1] \ra X \\
\ h(x,t) = h_t(x)
\end{cases}
$
is continuous.
\\ \indent
[Note: \ When \mX is a LCH space, isotopies correspond to paths in $H(X)$ (compact open topology).]\\
\endgroup %%------------------------------------<<

\begingroup%%----------------------------------->>
\fontsize{9pt}{11pt}\selectfont
\textbf{\small EXAMPLE} \  
A homeomorphism $h:\R^n \ra \R^n$ is said to be 
\un{stable}
\index{stable (homeomorphism)} 
if $\exists$ homeomorphisms $h_1, \ldots, h_k:\R^n \ra \R^n$ such that $h = h_1 \circ \cdots \circ h_k$, where each $h_i$ has the property that for some nonempty open $U_i \subset \R^n$, $\restr{h_i}{U_i} = \id_{U_i}$.  
Every stable homeomorphism of $\R^n$ is isotopic to the identity.
\\ \indent
[Take $k = 1$ and consider a homeomorphism $h:\R^n \ra \R^n$ for which $\restr{h}{U} = \id_U$.  
Define an isotopy $\{h_t: 0 \leq t \leq 1\}$ of $\R^n$ as follows.  Fix $u \in U$ and put
$
h_t(x) = 
\begin{cases}
\ h(x + 2tu) - 2tu  \hspace{1.0cm}  (0 \leq t \leq 1/2)\\
\ \ds\frac{1}{2 - 2t} h_{1/2}((2 - 2t)x) \hspace{0.3cm} (1/2 \leq t < 1)
\end{cases}
$
$\&$ $h_1(x) = x$.]\\
\endgroup %%------------------------------------<<

%%----------------------------------------------------------------------------------------------07
\begingroup%%----------------------------------->>
\fontsize{9pt}{11pt}\selectfont
\textbf{\small FACT} \  \ 
Equip $H(\R^n)$ with the compact open topology and write $H_\ST(\R^n)$ for the subspace of $H(\R^n)$ consisting of the stable homeomorphisms 
$-$then $H_\ST(\R^n)$ is an open subgroup of $H(\R^n)$.
\\ \indent
[Note: \ Therefore $H_\ST(\R^n)$ is also a closed subgroup of $H(\R^n)$ (since $H(\R^n)$ is a topological group in the compact open topology).]\\
\endgroup %%------------------------------------<<

\begingroup%%----------------------------------->>
\fontsize{9pt}{11pt}\selectfont
Application: The path component of $\id_{\R^n}$ in $H(\R^n)$ is $H_\ST(\R^n)$.
\\ \indent
[In view of the example, there is a path from every element of $H_\ST(\R^n)$ to $\id_{\R^n}$.  On the other hand, if 
$\tau:[0,1] \ra H(\R^n)$ is a path with $\tau(1) = \id_{\R^n}$ but $\tau(0) \notin H_\ST(\R^n)$, then 
$\tau^{-1}(H_\ST(\R^n))$ would be a nontrivial clopen subset of $[0,1]$.]
\\ \indent
[Note: \ It can be shown that $H(\R^n)$  is locally path connected (indeed, locally contractible 
(cf. p. \pageref{2.5}
)).]\\
\endgroup %%------------------------------------<<

\begingroup%%----------------------------------->>
\fontsize{9pt}{11pt}\selectfont
An isotopy $\{h_t: 0 \leq t \leq 1\}$ is said to be 
\un{invertible}
\index{invertible (isotopy)} 
if the collection 
$\{h_t^{-1}:0 \leq t \leq 1\}$ is an isotopy.\\
\endgroup %%------------------------------------<<

\begingroup%%----------------------------------->>
\fontsize{9pt}{11pt}\selectfont
\textbf{\small LEMMA} \  
An isotopy $\{h_t:0 \leq t \leq 1\}$ is invertible iff the function $H:X \times [0,1] \ra X \times [0,1]$ defined by the rule 
$(x,t) \ra (h_t(x),t)$ is a homeomorphism.
\\ \indent
[Note: \mH is necessarily one-to-one, onto, and continuous.]\\
\endgroup %%------------------------------------<<

\begingroup%%----------------------------------->>
\fontsize{9pt}{11pt}\selectfont
\textbf{\small FACT} \  
Let \mX be a LCH space $-$then every isotopy $\{h_t: 0 \leq t \leq 1\}$ of \mX is invertible.
\\ \indent
[Show first that $\forall \ x \in X$, $h_t^{-1}(x)$ is a continuous function of $t$.]\\
\endgroup %%------------------------------------<<

\begingroup%%----------------------------------->>
\fontsize{9pt}{11pt}\selectfont
\textbf{\small FACT} \  
Let \mX be a LCH space $-$then every isotopy $\{h_t: 0 \leq t \leq 1\}$ of \mX extends to an isotopy of $X_\infty$.
\\ \indent
[Define $\ov{h}_t:X_\infty \ra X_\infty$ by $\restr{\ov{h}_t}{X} = h_t$ $\&$ $\ov{h}_t(\infty) = \infty$.  
To verify that $\ov{h}$ is continuous, extend $H$ to $X_\infty \times [0,1]$ via the prescription 
$\ov{H}(\infty,t) = (\ov{h}_t(\infty),t)$ so $\ov{h} = \pi_\infty \circ \ov{H}$, where $\pi_\infty$ is the projection of 
$X_\infty \times [0,1]$ onto $X_\infty$.  Establish the continuity of $\ov{H}$ by utilizing the continuity of $H^{-1}$ (the substance of the previous result.]\\
\endgroup %%------------------------------------<<

\begingroup%%----------------------------------->>
\fontsize{9pt}{11pt}\selectfont
\textbf{\small EXAMPLE} \  
Every isotopy $\{h_t: 0 \leq t \leq 1\}$ of $\R^n$ extends to an isotopy of $\bS^n$.\\
\endgroup %%------------------------------------<<

Let \mX be a CRH space, $(Y,d)$ a metric space.  
Given $f \in C(X,Y)$ and 
$\phi \in C(X,\R_{>0})$, put $N_\phi(f) = \{g: d(f(x),g(x)) < \phi(x) \ \forall \ x\}$.

Observations: 
(1) If $\phi_1,  \phi_2, \in C(X,\R_{>0})$, then $N_\phi(f) \subset N_{\phi_1}(f) \cap N_{\phi_2}(f)$, where 
$\phi(x) = \min\{\phi_1(x),\phi_2(x)\}$;
(2) If $g \in N_\phi(f)$, then $N_\psi(g) \subset N_\phi(f)$, where $\psi(x) = \phi(x) - d(f(x),g(x))$.

Therefore the collection $\{N_\phi(f)\}$ is a basic system of neighborhoods at $f$.  
Accordingly, varying $f$ leads to a 
topology on $C(X,Y)$, the 
\un{majorant topology}
\index{majorant topology (on $C(X,Y)$)}.

%%----------------------------------------------------------------------------------------------08
[Note: \ Each $\phi \in C(X,\R_{>0})$ determines a metric $d_\phi$ on $C(X,Y)$, viz. 
$\ds d_\phi(f,g) = \min\left\{1, \sup\limits_{x \in X} \ds\frac{d(f(x),g(x))}{\phi(x)}\right\}$, and their totality defines the majorant topology on $C(X,Y)$, which is thus completely regular.  
However, in general, the majorant topology on $C(X,Y)$ need not be normal 
(Wegenkitt\footnote[2]{\textit{Ann. Global Anal. Geom.} \textbf{7} (1989), 171-178; 
see also van Douwen, \textit{Topology Appl.} \textbf{39} (1991), 3-32.}).]\\

\begingroup%%----------------------------------->>
\fontsize{9pt}{11pt}\selectfont
Here is a proof that $C(X,Y)$ (majorant topology) is completely regular.\label{2-8a}  
Fix a closed subset $A \subset C(X,Y)$ and an $f \in C(X,Y) - A$.  
Choose $\phi \in C(X,\R_{>0})$: 
$N_\phi(f) \subset C(X,Y) - A$.  
Define a function $\Phi:C(X,Y) \ra [0,1]$ by 
$\ds \Phi(g) = \sup\limits_{x \in X} \ds\frac{d(f(x),g(x))}{\phi(x)}$ if $g \in N_\phi(f)$ 
and let it be 1 otherwise $-$then $\Phi$ is continuous and $\Phi(f) = 0$, $\restr{\Phi}{A} = 1$.]
\\ \indent
[Note: \ The verification of the continuity of $\Phi$ hinges on the observation that 
$g \in \ov{N_\phi(f)}$ $\implies$ $d(f(x),g(x)) \leq \phi(x)$ $\forall \ x$, hence $\forall \ g \in \ov{N_\phi(f)} - N_\phi(f)$, 
$\sup\limits_{x \in X} \ds\frac{d(f(x),g(x))}{\phi(x)} = 1$.]\\
\endgroup %%------------------------------------<<

Example: Suppose that the sequence $\{f_k\}$ converges to $f$ in $C(\R^n,\R^n)$ (majorant topology) $-$then $\exists$ a compact $K \subset \R^n$ and an index $k_0$ such that $f_k(x) = f(x)$ $\forall \ k > k_0$ $\&$ 
$\forall  \ x \in \R^n - K$.\\

\begingroup%%----------------------------------->>
\fontsize{9pt}{11pt}\selectfont
\textbf{\small EXAMPLE} \  
Suppose that $f:\R^n \ra \R^n$ is a homeomorphism $-$then $f$ has a neighborhood of surjective maps in $C(\R^n,\R^n)$ (majorant topology).\\
\endgroup %%------------------------------------<<

\begingroup%%----------------------------------->>
\fontsize{9pt}{11pt}\selectfont
\textbf{\small EXAMPLE} \  
Equip $H(\R^n)$ with the majorant topology $-$then the path component of $\id_{\R^n}$ in $H(\R^n)$ consists of those homeomorphisms that are the identity outside some compact set.\\
\endgroup %%------------------------------------<<

\begingroup%%----------------------------------->>
\fontsize{9pt}{11pt}\selectfont
\textbf{\small FACT} \  
The majorant topology on $C(\R^n,\R^n)$ is not first countable.\\
\endgroup %%------------------------------------<<

\textbf{\small LEMMA} \  
The compact open topology on $C(X,Y)$ is smaller than the majorant topology.

[Fix a compact $K \subset X$, an open $V \subset Y$, and a continuous $f:X \ra Y$ such that $f(K) \subset V$.  Choose 
$\epsilon > 0$ such that $\forall \ y \in f(K)$, $d(y,y^\prime) < \epsilon$ $\implies$ $y^\prime \in V$.  
Let $\phi \in C(X,\R_{> 0})$ be the constant function $x \ra \epsilon$ $-$then $\forall \ g \in N_\phi(f)$, $g(K) \subset V$.]\\

Remark: The 
\un{uniform topology}
\index{uniform topology (on $C(X,Y)$ )} 
on $C(X,Y)$ is the topology induced by 
the metric 
$d(f,g) = \min\left\{1, \sup\limits_{x \in X} d(f(x),g(x)) \right\}$.  The proof of the lemma shows that the compact open topology on $C(X,Y)$ is smaller than the uniform topology (which in turn is smaller than the majorant topology).\\

%%----------------------------------------------------------------------------------------------09
\begingroup%%----------------------------------->>
\fontsize{9pt}{11pt}\selectfont
\textbf{\small FACT} \  
The compact open topology on $C(X,Y)$ equals the uniform topology if \mX is compact.\\

\textbf{\small FACT} \  
The uniform topology on $C(X,Y)$ equals the majorant topology if \mX is pseudocompact.\\
\endgroup %%------------------------------------<<

Let $M(Y)$ be the set of all metrics on \mY which are compatible with the topology of \mY $-$then the 
\un{limitation topology}
\index{limitation topology} 
on $C(X,Y)$ has for a neighborhood basis at $f$ the $N_m(f)$, 
$(m \in M(Y))$, where $N_m(f) = \left\{g: \sup\limits_{x \in X} m(f(x),g(x)) < 1\right\}$.

[Note: \ If $m_1$, $m_2 \in M(Y)$, then $N_{m_1 + m_2}(f) \subset N_{m_1}(f) \cap N_{m_2}(f)$, and if 
$g \in N_m(f)$, then 
$N_{\left(\ds\frac{2}{\epsilon}\right)m}(g) \subset N_m(f)$, where 
$m(f(x),g(x)) \leq 1 - \epsilon$ 
$\forall \ x$.]\\

\begingroup%%----------------------------------->>
\fontsize{9pt}{11pt}\selectfont
The limitation topology is defined by the metrics 
$(f,g) \ra \min\left\{1, \sup\limits_{x \in X} m(f(x),g(x))\right\}$ $(m \in M(Y))$, 
thus the uniform topology on $C(X,Y)$ is smaller than the limitation topology.\\
\endgroup %%------------------------------------<<

\textbf{\small LEMMA} \  
Suppose that \mX is paracompact $-$then the limitation topology on $C(X,Y)$ is smaller than the majorant topology.

[Fix $m \in M(Y)$ and let $f \in C(X,Y)$.  
By compatibility, $\forall \ x \in X$, $\exists$ $\epsilon(x) > 0$: 
$d(f(x),y) < \epsilon(x)$ $\implies$ $\ds m(f(x),y) < \frac{1}{4}$. 
Put 
$O_x = \left\{x^\prime:d(f(x),f(x^\prime)) < \ds\frac{\epsilon(x)}{2}\right\}$ $-$then $\{O_x\}$ is an open covering of \mX.  
Let $\{U_x\}$ be a precise neighborhood finite open refinement and choose a subordinated partition of unity $\{\kappa_x\}$.  
Definition: $\ds\phi = \sum\limits_x \ds\frac{\epsilon(x)}{2} \kappa_x$.  Consider now any $x_0 \in X$ and assume that 
$d(f(x_0),y) < \phi(x_0)$.  
Let $\kappa_{x_1} \ldots, \kappa_{x_n}$ be an enumeration of those $\kappa_x$ whose support contains $x_0$ and fix $i$ between 1 and $n$: 
$\ds\frac{\epsilon(x_j)}{2} \leq \ds\frac{\epsilon(x_i)}{2}$ $(j = 1, \ldots, n)$ to get 
$\ds\phi(x_0) \leq \frac{\epsilon(x_i)}{2}$.    But $x_0 \in U_{x_i} \subset O_{x_i}$.  Therefore 
$\ds d(f(x_i),f(x_0)) < \frac{\epsilon(x_i)}{2}$ $\ds \left (\implies m(f(x_i),f(x_0)) < \frac{1}{4}\right)$ 
$\implies$ $d(f(x_i,y) < \epsilon(x_i)$ 
$\implies$ $\ds m(f(x_i),y) < \frac{1}{4}$ 
$\implies$ $\ds m(f(x_0),y) < \frac{1}{2}$.  And this shows that $N_\phi(f) \subset N_m(f)$.]

[Note: \ In general, the limitation topology is strictly smaller than the majorant topology.  To see this, observe that 
$C(\R,\R)$ is a topological group under addition in the majorant topology.  On the other hand, there is a countable basis at a given $f \in C(\R,\R)$ (limitation topology) iff $f$ is bounded, thus $C(\R,\R)$ is not a topological group under addition in the limitation topology.]\\

\begingroup%%----------------------------------->>
\fontsize{9pt}{11pt}\selectfont
\textbf{\small FACT} \  
Take $X = Y$ $-$then in the limitation topology, $H(X)$ is a topological group.\\
\endgroup %%------------------------------------<<

\index{Refinement Principle}
\textbf{\small REFINEMENT PRINCIPLE} \  \ 
Let $(Y,d)$ be a metric space $-$then for any open covering $\sV = \{V\}$ of \mY, $\exists \ m \in M(Y)$ such that the collection $\{V_y\}$ is a refinement of $\sV$, where $V_y = \{y^\prime: m(y,y^\prime) < 1\}$.

%%----------------------------------------------------------------------------------------------10
[A proof can be found in 
Dugundji\footnote[2]{\textit{Topology}, Allyn and Bacon (1966), 196; 
see also Bessaga-Pelczy\'nski, \textit{Selected Topics in Infinite Dimensional Topology}, PWN (1975), 63.\vspace{0.11cm}}.]\\

\textbf{\small LEMMA} \  
Let $(Y,d)$ be a metric space $-$then for any $\delta \in C(Y,\R_{>0})$, $\exists \ m \in M(Y)$: $d(y,y^\prime) < \delta(y)$ 
whenever $m(y,y^\prime) < 1$.

[Choose an open covering $\sV = \{V\}$ of \mY such that the diameter of a given \mV is 
$\ds \leq \frac{1}{2} \inf \delta(V)$.  
Using the refinement principle, fix an $m \in M(Y)$ such that the collection $\{V_y\}$ refines $\sV$.  
If $(y,y^\prime)$ is a pair with $m(y,y^\prime) < 1$, then $V_y \subset V$ for some \mV, hence 
$y,y^\prime \in V$ $\implies$ $\ds d(y,y^\prime)  \leq \frac{1}{2} \delta(y) < \delta(y)$.]\\

\begin{proposition} \ %10
Take $X = Y$ $-$then the limitation topology on $H(X)$ is equal to the majorant topology.
\end{proposition}

[Fix $f \in H(X)$ and let $\phi \in C(X,\R_{>0})$.  
Thanks to the lemma, $\exists$ $m \in M(X)$: 
$d(x,x^\prime) < \phi \circ f^{-1}(x)$ whenever $m(x,x^\prime) < 1$.  
If $g \in H(X)$ and 
$\sup\limits_{x \in X} m(f(x),g(x)) < 1$, then $d(f(x),g(x)) < \phi \circ f^{-1}(f(x)) = \phi(x)$ $\forall \ x$, i.e., 
$N_\phi(f) \cap H(X)$ is open in $H(X)$ (limitation topology).]\\

Application: The homeomorphism group of a metric space is a topological group in the majorant topology.\\

\begingroup%%----------------------------------->>
\fontsize{9pt}{11pt}\selectfont
\textbf{\small EXAMPLE} \  
Let \mX be a second countable topological manifold of euclidean dimension $n$ $-$then in the majorant topology, $H(X)$ is a topological group.  Moreover, 
\u Cernavski\u i\footnote[3]{\textit{Math. Sbornik} \textbf{8} (1969), 287-333.}
 has shown that $H(X)$ is locally contractible.
 
[Note: \ \mX is metrizable (cf. $\S 1$, Proposition 11), so $\exists$ $d$: $(X,d)$ is a metric space.]\\
\endgroup %%------------------------------------<<

Notation: $\forall \ f \in C(X,Y)$, $\gr_f \subset X \times Y$ is its graph.

Given an open subset $O \subset X \times Y$, let $\Gamma_O = \{f:\gr_f \subset O\}$ $-$then the collection 
$\{\Gamma_O\}$ is a basis for a topology on $C(X,Y)$, the 
\un{graph topology}.
\index{graph topology}

[Note: \ In this connection, observe that $\Gamma_O \cap \Gamma_P = \Gamma_{O \cap P}$.]\\

\textbf{\small LEMMA}\label{2.7} \  
The majorant topology on $C(X,Y)$ is smaller than the graph topology.

[The function $(x,y) \ra \phi(x) - d(f(x),y)$ from $X \times Y$ to $\R$ is continuous, thus 
$O = \{(x,y): d(f(x),y) < \phi(x)\}$ is an open subset of $X \times Y$.  But $\Gamma_O = N_\phi(f)$.]\\

%%----------------------------------------------------------------------------------------------11
Rappel: \ A function $f:X \ra \R$ is 
\un{lower semicontinuous}
\index{lower semicontinuous}
(\un{upper semicontinuous}
\index{upper semicontinuous})
if for each real number c, 
$\{x:f(x) > c\}$ 
$(\{x:f(x) < c\})$ 
is open.  
Example: The characteristic function of a subset \mS of \mX is lower semicontinuous (upper semicontinuous) iff \mS is open (closed).\\

\index{Hahn's Einschiebungsatz}
\textbf{\small HAHN'S EINSCHIEBUNGSATZ} \  
Suppose that \mX is paracompact.  Let $g:X \ra \R$ be lower semicontinuous and $G:X \ra \R$ upper semicontinuous.  
Assume: $G(x) < g(x)$ $\forall \ x \in X$ $-$then $\exists$ a continuous function $f:X \ra \R$ such that 
$G(x) < f(x) < g(x)$ $\forall \ x \in X$.

[Put $U_r = \{x: G(x) < r\} \cap \{x:g(x) > r\}$ ($r$ rational).  
Each $U_r$ is open and $X = \ds\bigcup\limits_r U_r$.  
Let $\{\kappa_r\}$ be a partition of unity subordinate to $\{U_r\}$ and take $\ds f = \ds\sum\limits_r r \kappa_r$.]\\

\begingroup%%----------------------------------->>
\fontsize{9pt}{11pt}\selectfont
The following result characterizes the class of \mX satisfying the conditions of Hahn's einschiebungsatz.\\

\label{2.8}
\textbf{\small FACT} \  
Let \mX be a CRH space $-$then \mX is normal and countably paracompact iff for every lower semicontinuous 
$g:X \ra \R$ and upper semicontinuous $G:X \ra \R$ such that $G(x) < g(x)$ $\forall \ x \in X$, $\exists$ $f \in C(X,\R)$: 
$G(x) < f(x) < g(x)$ $\forall \ x \in X$.

[Necessity:  With $r$ running through the rationals, there exists a neighborhood finite open covering $\{O_r\}$ of \mX: 
$O_r \subset \{x:G(x) < r < g(x)\}$ $\forall \ r$ and a neighborhood finite open covering $\{P_r\}$ of \mX: 
$\ov{P}_r \subset O_r$, $\forall \ r$.  Fix a continuous function $f_r:X \ra [-\infty,r]$ such that 
$
f_r(x) = 
\begin{cases}
\ -\infty \hspace{0.3cm} \ (x \notin O_r)\\
\ r \hspace{0.75cm} \  (x \in \ov{P}_r)
\end{cases}
\hspace{-.2cm}. \ 
$
Put $f(x) = \sup\limits_r f_r(x)$ $-$then $f$ has the required properties.

Sufficiency: There are two parts.

\indent\indent \mX is normal.  Thus let \mA, \mB be disjoint closed subsets of \mX.  With \mG the characteristic function of \mA, let $g$ be defined by
$
\begin{cases}
\ g(x) = 1 \quadx (x \in B)\\
\ g(x) = 2 \quadx (x \notin B)
\end{cases}
: \ 
$
$g$ is lower semicontinuous, $G$ is upper semicontinuous, and $G(x) < g(x)$ $\forall \ x \in X$.  Choose $f \in C(X,\R)$ per the assumption and let $U = \{x: f(x) > 1\}$, $V = \{x: f(x) < 1\}$ $-$then 
$
\begin{cases}
\ U\\
\ V
\end{cases}
$
are disjoint open subsets of \mX and 
$
\begin{cases}
\ A \subset U\\
\ B \subset V
\end{cases}
, \ 
$
hence \mX is normal.

\indent\indent \mX is countably paracompact.  
Thus consider any decreasing  sequence $\{A_n\}$ of closed sets such that 
$\ds\bigcap\limits_n A_n = \emptyset$.  
Put $\ds g(x) = \ds\frac{1}{n + 1}$ $(x \in A_n - A_{n+1}, n = 0, 1, \ldots)$ 
$(A_0 = X)$: $g$ is lower semicontinuous.  
Take $f \in C(X,\R)$: $0 < f(x) < g(x)$ and let 
$\ds U_n = \{x:f(x) < \ds\frac{1}{n + 1}\}$ $-$then $\{U_n\}$ is a decreasing sequence of open sets with 
$A_n \subset U_n$ for every $n$ and $\ds\bigcap\limits_n U_n = \emptyset$.  
Since \mX is normal, this guarantees that \mX is also countably paracompact (via CP 
(cf. p. \pageref{2.6})
).]\\
\endgroup %%------------------------------------<<

\textbf{\small LEMMA} \  
Assume that \mX is paracompact and suppose given a neighborhood finite closed covering $\{A_j:j \in J\}$ of \mX and 
$\forall \ j$, a positive real number $a_j$ $-$then $\exists$ a continuous function $\phi:X \ra \R_{> 0}$ such that 
$\phi(x) < a_j$ if $x \in A_j$.

[The function from \mX to $\R$ defined by the rule $x \ra \min\{a_j: x \in A_j\}$ is lower semicontinuous and strictly positive.]\\

%%----------------------------------------------------------------------------------------------12
\begin{proposition} \ %11
The majorant topology on $C(X,Y)$ is independent of the choice of $d$ provided that \mX is paracompact.
\end{proposition}

[It suffices to show that the graph topology on $C(X,Y)$ is smaller than the majorant topology 
(cf. p. \pageref{2.7}).  
So fix an $f \in \Gamma_O$ and consider any $x_0 \in X$.  
Choose a neighborhood $U_0$ of $x_0$ and a positive real number $a_0$ such that $x \in U_0$ $\&$ $d(f(x_0),y) < 2a_0$ 
$\implies$ $(x,y) \in O$.
Choose further a neighborhood $V_0$ of $x_0$ such that $V_0 \subset U_0$ $\&$ $d(f(x_0),f(x)) < a_0$ $\forall \ x \in V_0$ 
$-$then $\{(x,y):x \in V_0 \ \& \ d(f(x),y) < a_0\} \subset O$.  
From this, it follows that one can find a neighborhood finite closed covering $\{A_j: j \in J\}$ of \mX and a set $\{a_j: j \in J\}$ 
of positive real numbers for which $\{(x,y): x \in A_j \ \& \ d(f(x),y) < a_j\} \subset O$.  
In view of the lemma, $\exists$ a continuous function, $\phi:X \ra \R_{>0}$ with $\phi(x) < a_j$ whenever $x \in A_j$, 
hence $N_\phi(f) \subset \Gamma_O$, i.e., every point of $\Gamma_O$ is an interior point in the majorant topology.]\\

To reiterate: If \mX is paracompact, then the majorant topology on $C(X,Y)$ equals the graph topology.

[Note: \ The assumption of paracompactness can be relaxed (see below).]\\

\begingroup%%----------------------------------->>
\fontsize{9pt}{11pt}\selectfont
Let \mX be a CRH space, $(Y,d)$ a metric space.  
Given $f \in C(X,Y)$ and a lower semicontinuous 
$\sigma:X \ra \R_{> 0}$, put $N_\sigma(f) = \{g:d(f(x),g(x)) < \sigma(x) \  \forall \ x\}$.
\\ \indent
Observations: 
(1) \ If $\sigma_1, \sigma_2:X \ra \R_{>0}$ are lower semicontinuous, then 
$N_\sigma(f) \subset N_{\sigma_1}(f) \cap N_{\sigma_2}(f)$, where $\sigma(x) = \min\{\sigma_1(x),\sigma_2(x)\}$;
(2) \ If $g \in N_\sigma(f)$, then $N_\tau(g) \subset N_\sigma(f)$, where $\tau(x) = \sigma(x) - d(f(x),g(x))$.
\\ \indent
[Note: \ The minimum of two lower semicontinuous functions is lower semicontinuous, so $\sigma$ is lower semicontinuous.  
On the other hand, the sum of two lower semicontinuous functions is lower semicontinuous.  
But $x \ra d(f(x),g(x))$ is continuous, thus $x \ra -d(f(x),g(x))$ is lower semicontinuous, so $\tau$ is lower semicontinuous.]
\\ \indent
Therefore, the collection $\{N_\sigma(f)\}$ is a basic system of neighborhoods at $f$.  
Accordingly, varying $f$ leads to a topology on $C(X,Y)$, the 
\un{semimajorant topology}.
\index{semimajorant topology}\\
\endgroup %%------------------------------------<<

\begingroup%%----------------------------------->>
\fontsize{9pt}{11pt}\selectfont
\textbf{\small LEMMA} \  
The semimajorant topology on $C(X,Y)$ is smaller than the graph topology.
\\ \indent
[Let $O = \{(x,y): d(f(x),y) < \sigma(x)\}$ $-$then $\Gamma_O$ is open in $C(X,Y)$.  
Proof: Fix $(x_0,y_0) \in O$, 
put $\ds \epsilon = \frac{1}{3}\big(\sigma(x_0) - d(f(x_0),y_0)\big)$, 
and note that the subset of \mO consisting of those $(x,y)$ such that 
$\sigma(x) > \sigma(x_0) - \epsilon$, 
$d(f(x),f(x_0)) < \epsilon$, and 
$d(y,y_0) < \epsilon$ is open.  
And: $N_\sigma(f) = \Gamma_O$.]\\
\endgroup %%------------------------------------<<

\begingroup%%----------------------------------->>
\fontsize{9pt}{11pt}\selectfont
\textbf{\small LEMMA} \  
The graph topology on $C(X,Y)$ is smaller than the semimajorant topology.
\\ \indent
[Fix an $f \in \Gamma_O$.  
Define a strictly positive function $\sigma:X \ra \R$ by letting $\sigma(x_0)$ be the supremum of those $a_0 \in ]0,1]$ 
for which $x_0$ has a neighborhood $U_0$ such that $x \in U_0$ $\&$ 
$d(f(x_0),y) < a_0$ $\implies$ $(x,y) \in O$.  
Since $N_\sigma(f) \subset \Gamma_O$, the point is to prove that $\sigma$ is lower semicontinuous, 
i.e., that $\forall \ c \in \R$, $\{x:c < \sigma(x)\}$ is open.  
This is trivial if $c \leq 0$ or $c \geq 1$, 
so take $c \in \ ]0,1[$ and fix $x_0$: $c < \sigma(x_0)$.  
Put $\epsilon = (\sigma(x_0) - c)/3$ $-$then 
%%----------------------------------------------------------------------------------------------13
$c + 2\epsilon < \sigma(x_0)$, thus $\exists$ a neighborhood $U_0$ of $x_0$ such that $x \in U_0$ $\&$ 
$d(f(x_0),y) < c + 2\epsilon$ $\implies$ $(x,y) \in O$.  
Supposing further that $x \in U_0$ $\implies$ 
$d(f(x_0),f(x)) < \epsilon$, one has $x \in U_0$ $\&$ $d(f(x),y) < c + \epsilon$ $\implies$ $(x,y) \in O$ $\implies$ 
$c < c + \epsilon \leq \sigma(x)$.]\\
\endgroup %%------------------------------------<<

\begingroup%%----------------------------------->>
\fontsize{9pt}{11pt}\selectfont
\textbf{\small FACT} \  
The semimajorant topology on $C(X,Y)$ equals the graph topology.\\
\endgroup %%------------------------------------<<

\begingroup%%----------------------------------->>
\fontsize{9pt}{11pt}\selectfont
A CRH space \mX is said to be a 
\un{CB space}
\index{CB space} 
if for every strictly positive lower semicontinuous $\sigma:X \ra \R$ there exists a strictly positive continuous function $\phi:X \ra \R$ such that $0 < \phi(x) \leq \sigma(x)$ $\forall \ x \in X$.
\\ \indent
Example: If \mX is normal and countably paracompact, then \mX is a CB space (cf. p. 
\pageref{2.8}
).
\\ \indent
Examples: 
(Mack\footnote[2]{\textit{Proc. Amer. Math. Soc.} \textbf{16} (1965), 467-472.}): 
(1) \ Every countably compact space is a CB space; 
(2) \ Every CB space is countably paracompact.\\
\endgroup %%------------------------------------<<

\begingroup%%----------------------------------->>
\fontsize{9pt}{11pt}\selectfont
\textbf{\small EXAMPLE} \  
The Isbell-Mr\'owka space $\Psi(\N)$ is a pseudocompact LCH space which is not countably paracompact 
(cf. p. \pageref{2.9} ), 
hence is not a CB space.\\
\endgroup %%------------------------------------<<

\begingroup%%----------------------------------->>
\fontsize{9pt}{11pt}\selectfont
\textbf{\small FACT} \  
The majorant topology on $C(X,Y)$ equals the graph topology $\forall$ pair $(Y,d)$ iff \mX is a CB space.
\\ \indent
[Necessity: Fix a strictly positive lower semicontinuous $\sigma:X \ra \R$.  
Specialized to the case $Y = \R$, the assumption is that the majorant topology on $C(X)$ equals the semimajorant topology, so working with 
$N_\sigma(0)$, $\exists$ $\phi$: $N_\phi(0) \subset N_\sigma(0)$ $\implies$ 
$(1 - \epsilon)\phi \in N_\phi(0) \subset N_\sigma(0)$ $(0 < \epsilon < 1)$ $\implies$ $0 < \phi(x) \leq \sigma(x)$ 
$\forall \ x \in X$, thus \mX is a CB space.
\\ \indent
Sufficiency: Since $N_\phi(f) \subset N_\sigma(f)$, the semimajorant topology on $C(X,Y)$ is smaller than the majorant topology.]\\
\endgroup %%------------------------------------<<

If $(Y,d)$ is a complete metric space, then co$C(X,Y)$ need not be Baire.  
Examples: 
(1) \ co$C([0,\Omega[,\R)$ is not Baire;
(2) \ co$C(\Q,\R)$ is not Baire.

[Note: \ Recall, however, that if \mX is countable at infinity and compactly generated, then co$C(X,Y)$ is completely metrizable (cf. Proposition 7), hence is Baire.]\\

\begin{proposition} \ %12
Assume: $(Y,d)$ is a complete metric space $-$then $C(X,Y)$ (majorant topology) is Baire.
\end{proposition}

[Let $\{O_n\}$ be a sequence of dense open subsets of $C(X,Y)$.  Let \mU be a nonempty open subset of $C(X,Y)$.  
Since $U \cap O_1$ is nonempty and open and since $C(X,Y)$ is completely regular 
(cf. p. \pageref{2-8a}),
$\exists$ $f_1 \in U \cap O_1$ $\&$ $\phi_1 \in C(X,\R_{>0})$: 
$\{g:d(f_1(x),g(x)) \leq \phi_1(x) \ \forall \ x\}$ $\subset$ $U \cap O_1$, where $\phi_1 < 1$.  
Next, $\exists$ $f_2 \in N_{\phi_1}(f_1) \cap O_2$ $\&$ $\phi_2 \in C(X,\R_{> 0})$ : $\{g:d(f_2(x),g(x) \leq
%%----------------------------------------------------------------------------------------------14
\phi_2(x) \  \forall x\}$ $\subset N_{\phi_1}(f_1) \cap O_2$, where $\phi_2 < \phi_1/2$.  
Proceeding, $\exists$ 
$f_{n+1} \in N_{\phi_n}(f_n) \cap O_{n+1}$ $\&$ $\phi_{n+1} \in C(X,\R_{>0})$ : 
$\{g: d(f_{n+1}(x),g(x)) \leq \phi_{n+1}(x) \ \forall \ x\}$ $\subset N_{\phi_n}(f_n) \cap O_{n+1}$, where
$\phi_{n+1} < \phi_n/2$.  
So $\forall \ x$, 
$\ds d(f_{n+1}(x),f_n(x)) \leq \ds\frac{1}{2^{n-1}}$, thus $\{f_n(x)\}$ is a Cauchy sequence in \mY.  
Definition: $f(x) = \lim f_n(x)$.  
Because the convergence is uniform, $f \in C(X,Y)$.  Moreover, 
$d(f_n(x),f(x)) \leq \phi_n(x)$ $\forall \ n$ $\&$ $\forall \ x$, which implies that 
$f \in U \cap \bigl(\ds\bigcap\limits_n O_n\bigr)$.]\\

\begingroup%%----------------------------------->>
\fontsize{9pt}{11pt}\selectfont
\textbf{\small FACT} \  
Assume: $(Y,d)$ is a complete metric space $-$then $C(X,Y)$ (limitation topology) is Baire.\\
\endgroup %%------------------------------------<<

Convention: 
Maintaining the assumption that \mX is a CRH space, $C(X)$ henceforth carries the compact open topology.

Let \mK be a compact subset of \mX.  
Put 
$p_K(f) = \sup\limits_K \abs{f}$ $(f \in C(X))$ $-$then $p_K:C(X) \ra \R$ is a seminorm on $C(X)$, 
i.e., $p_K(f) \geq 0$, 
$p_K(f + g) \leq p_K(f) + p_K(g)$,  
$p_K(cf) = \abs{c}p_K(f)$.

[Note: \ More is true, viz.  $p_K$ is multiplicative in the sense that $p_K(fg) \leq p_K(f) p_K(g)$.]

Remark: \ The initial topology on $C(X)$ determined by the $p_K$ as \mK runs through the compact subsets of \mX is the compact open topology.

[Note: \ In the compact open topology, $C(X)$ is a Hausdorff locally convex topological vector space.]\\

\begingroup%%----------------------------------->>
\fontsize{9pt}{11pt}\selectfont
\label{2.18}
\label{6.1}
Observation: If $K \subset X$ is compact and if $f \in C(K)$, then $\exists$ $F \in BC(X)$: $\restr{F}{K} = f$.  
Proof: Apply the Tietze extension theorem to \mK regarded as a compact subset of $\beta X$.\\
\endgroup %%------------------------------------<<

A CRH space \mX is said to be a 
\un{$k_\R$-space}
\index{k$_\R$-space} 
provided that a real valued function 
$f:X \ra \R$ is continuous whenever its restriction to each compact subset of \mX is continuous.  
Example: A compactly generated space \mX is a $k_\R$-space (but not conversely (cf. infra)).\\

\begingroup%%----------------------------------->>
\fontsize{9pt}{11pt}\selectfont
\textbf{\small EXAMPLE}\label{2.11} \  
Let \mX be a $k_\R$-space.  Assume \mX is countable at infinity $-$then \mX is compactly generated.
\\ \indent
[Fix a ``defining'' sequence $\{K_n\}$ of compact subsets of \mX with $K_n \subset K_{n+1}$ $\forall \ n$.  
Claim: A subset \mA of \mX is closed if $A \cap K_n$ is closed in $K_n$ for each $n$.  
For if not, then \mA has an accumulation point $a_0$: $a_0 \notin A$, which can be taken in $K_1$ (adjust notation).  
Choose a continuous function $f_1:K_1 \ra \R$ such that $f_1(A \cap K_1) = \{0\}$ and $f_1(a_0) = 1$.  
Extend $f_1$ to  a continuous function 
$f_2:K_2 \ra \R$ such that $f_2(A \cap K_2) = \{0\}$.  
Repeat the process to get a function $f:X \ra \R$ such that $f(x) = f_n(x)$ $(x \in K_n)$.  Since \mX is a $k_\R$-space, $f$ is continuous.  
This, however, is a contradiction: $f(A) = \{0\}$, $f(a_0) = 1$.]\\
\endgroup %%------------------------------------<<

\begingroup%%----------------------------------->>
\fontsize{9pt}{11pt}\selectfont
\textbf{\small FACT} \  
A $k_\R$-space \mX is compactly generated iff $kX$ is completely regular.
\\ \indent
[If \mX is a $k_\R$-space, then $C(X) = C(kX)$.  So, the supposition that $kX$ is completely regular forces 
$X = kX$ (cf. $\S 1$, Proposition 14).]
\\ \indent
%%----------------------------------------------------------------------------------------------15
[Note: \ Recall that in general, \mX completely regular $\centernot\implies$ $kX$ completely regular 
(cf. p. \pageref{2.10}).]\\
\endgroup %%------------------------------------<<

\begin{proposition} \ %13
$C(X)$ is complete as a topological vector space iff \mX is a $k_{\R}$-space.
\end{proposition}

[Necessity: Suppose that $f:X \ra \R$ is a real valued function such that $\restr{f}{K}$ is continuous $\forall$ compact 
$K \subset X$.  Let $f_K \in C(X)$ be an extension of $\restr{f}{K}$ $-$then $\{f_K\}$ is a Cauchy net in $C(X)$, thus is convergent, 
say $\lim f_K = F$.  But $f = F$.

Sufficiency: Let $\{f_i\}$ be a Cauchy net in $C(X)$ $-$then $\forall$ compact $K \subset X$, the net 
$\restr{\{f_i}{K}\}$ is Cauchy in $C(K)$, hence has a limit, call it $f_K$.  If $K_1 \subset K_2$, then $\restr{f_{K_2}}{K_1} = f_{K_1}$, so the prescription 
$f(x) = f_K(x)$ $(x \in K)$ defines a function $f:X \ra \R$.  Since \mX is a $k_\R$-space, $f$ is continuous.  
And: $\lim f_i = f$.]\\

\begingroup%%----------------------------------->>
\fontsize{9pt}{11pt}\selectfont
\textbf{\small EXAMPLE} \  
Let $\kappa$ be a cardinal $> \omega$ $-$then $\N^\kappa$ is a $k_\R$-space but $\N^\kappa$ is not compactly generated.

[Note: \ $\N^\omega$ is homeomorphic to $\PP$, thus is compactly generated.]\\

\label{2.13}
\textbf{\small FACT} \  
Suppose that the closed bounded subsets of $C(X)$ are complete $-$then \mX is a $k_\R$-space.\\
\endgroup %%------------------------------------<<

\begin{proposition} \ %14
$C(X)$ is metrizable iff \mX is countable at infinity (cf. Proposition 6).
\end{proposition}

[Let $d$ be a compatible metric on $C(X)$.  Put $U_n = \{f: d(f,0) < 1/n\}$.  Choose a compact $K_n \subset X$ and a positive 
$\epsilon_n:f(K_n) \subset \ ]-\epsilon_n,\epsilon_n[$ $\implies$ 
$f \in U_n$ $-$then for any compact subset \mK of \mX, $\exists$ $n: K \subset K_n$.  Therefore \mX is countable at infinity.]\\

\begin{proposition} \ %15
$C(X)$ is completely metrizable iff \mX is countable at infinity and compactly generated (cf. Proposition 7).
\end{proposition}

[If $C(X)$ is completely metrizable, then $C(X)$ is complete as a topological vector space, so \mX is a $k_{\R}$-space 
(cf. Proposition 13), thus \mX, being  countable at infinity, is compactly generated 
(cf. p. \pageref{2.11}).]\\

\begingroup%%----------------------------------->>
\fontsize{9pt}{11pt}\selectfont
A CRH space \mX is said to be 
\un{topologically complete}
\index{topologically complete} 
if \mX is a $G_\delta$ in 
$\beta X$ or still, if \mX is a $G_\delta$ in any Hausdorff space containing it as a dense subspace.  
Example: $\PP$ is topologically complete but $\Q$ is not.
\\ \indent
Examples: 
(1) Every completely metrizable space is topologically complete and every topologically complete metrizable space is completely metrizable; 
(2) Every LCH space is topologically complete.
\\ \indent
[Note: A topologically complete space is necessarily compactly generated and Baire 
(Engleling\footnote[2]{\textit{General Topology}, Heldermann Verlag (1989), 197-198.}).]
\\ \indent
%%----------------------------------------------------------------------------------------------16
Remark: \ It can be shown that Proposition 15 goes through if the hypothesis ``completely metrizable'' is weakened to 
``topologically complete'' 
(McCoy-Ntantu\footnote[3]{\textit{SLN} \textbf{1315} (1988), 75.\vspace{0.11cm}}).\\
\endgroup %%------------------------------------<<

\label{2.14}
\begingroup%%----------------------------------->>
\fontsize{9pt}{11pt}\selectfont
\textbf{\small EXAMPLE} \  
Let \mX be a LCH space.  Assume: \mX is paracompact $-$then $C(X)$ is Baire.

[Using LCH$_3$ 
(cf. p. \pageref{2.12}), 
write $X = \ds\coprod\limits_i X_i$, where the $X_i$ are pairwise disjoint nonempty open $\sigma$-compact subspaces of \mX.  
Each $X_i$ is countable at infinity and there is a homeomorphism 
$C(X) \approx \ds\prod\limits_i C(X_i)$.  
But the $C(X_i)$ are completely metrizable (cf. Proposition 15), hence are topologically complete, and it is a fact that a product of topologically complete spaces is Baire 
(Oxtoby\footnote[6]{\textit{Fund. Math.} \textbf{49} (1961), 157-166.}).]
\\ \indent
[Note: \ The paracompactness assumption on \mX cannot be dropped.  
Example: Take $X = [0,\Omega[$ $-$then $C(X)$ is not Baire.  
Proof: Since \mX is pseudocompact, 
$O_n = \ds\bigcup\limits_x\{f: n < f(x) < n + 1\}$ is a dense open subset of $C(X)$ and 
$\ds\bigcap\limits_n O_n = \emptyset$.]\\
\endgroup %%------------------------------------<<

\begingroup%%----------------------------------->>
\fontsize{9pt}{11pt}\selectfont
\textbf{\small FACT} \  
Suppose that \mX is first countable and $C(X)$ is Baire $-$then \mX is locally compact.\\
\endgroup %%------------------------------------<<

\index{Theorem: Stone-Weierstrass Theorem}
\index{Stone-Weierstrass Theorem}
\textbf{\small STONE-WEIERSTRASS THEOREM} \  
Let \mX be a compact Hausdorff space.  Suppose that $\sA$ is a subalgebra of $C(X)$ which contains the constants and separates the points of \mX $-$then $\sA$ is uniformly dense in  $C(X)$.\\

\begingroup%%----------------------------------->>
\fontsize{9pt}{11pt}\selectfont
\textbf{\small EXAMPLE} \  
Let $0 < a < b < 1$ $-$then every $f \in C([a,b])$ can be uniformly approximated by polynomials 
$\ds \sum\limits_1^d n_kx^k$, $n_k$ integral.
\\ \indent
[It is enough to show that $\ds f = \ds\frac{1}{2}$ can be so approximated.  
Given an odd prime $p$, put 
$\phi_p(x) =$ $\ds\frac{1}{p}(1 - x^p - (1 - x)^p)$ : $\phi_p$ is a polynomial with integral coefficients, no constant term, and 
$p\phi_p \ra 1$ uniformly on $[a,b]$ as $p \ra \infty$.  
Now write $p = 2q + 1$, note that 
$\ds \abs{\ds\frac{1}{2} - \frac{q}{p}} < \ds\frac{1}{p}$, and consider $q\phi_p$.]\\
\endgroup %%------------------------------------<<

\begin{proposition} \ %16
Suppose that \mX is a compact Hausdorff space $-$then $C(X)$ is separable iff \mX is metrizable.
\end{proposition}

[Necessity: If $\{f_n\}$ is a uniformly dense sequence in $C(X)$, then the
$\ds \big\{x: \abs{f_n(x)} > \ds\frac{1}{2}\big\}$ constitute a basis for the topology on \mX, therefore \mX is second countable, hence metrizable.

Sufficiency: Let $d$ be a compatible metric on \mX.  Choose a countable basis $\{U_n\}$ for its topology and put 
$f_n(x) = d(x,X - U_n)$ $(x \in X)$ $-$then the $f_n$ separate the points of \mX, thus the subalgebra of $C(X)$ generated by 1 and the $f_n$ is uniformly dense in $C(X)$, so
%%----------------------------------------------------------------------------------------------17
the same is true of the rational subalgebra of $C(X)$ generated by 1 and the $f_n$.  But the latter is a countable set.]\\

\begingroup%%----------------------------------->>
\fontsize{9pt}{11pt}\selectfont
\textbf{\small EXAMPLE} \  
Assume that \mX is not compact and consider $BC(X)$, viewed as a Banach space in the supremum norm:
$\norm{f} = \sup\limits_X \abs{f}$ $-$then $BC(X)$ can be identified with $C(\beta X)$ 
($f \ra \beta f$ : $\norm{f} = \norm{\beta f}$).  Since $\beta X$ is not metrizable, it follows that $BC(X)$ is not separable.
\\ \indent
[Note: \ To see that $\beta X$ is not metrizable, fix a point $x_0 \in \beta X - X$ and, arguing by contradiction, choose a sequence $\{x_n\} \subset X$ of distinct $x_n$ having $x_0$ for their limit.  
Put 
$A = \{x_{2n}\}$, $B = \{x_{2n+1}\}$ $-$then \mA and \mB are disjoint closed subsets of \mX, so, by Urysohn, $\exists$ 
$\phi \in BC(X)$ such that $0 \leq \phi \leq 1$ with $\phi = 1$ on \mA and $\phi = 0$ on \mB.  Therefore 
$1 = \phi(x_{2n}) \ra \beta \phi(x_0)$ $\&$ $0 = \phi(x_{2n + 1}) \ra \beta \phi(x_0)$, an absurdity.]\\
\endgroup %%------------------------------------<<

\begin{proposition} \ %17
$C(X)$ is separable iff \mX admits a smaller separable metrizable topology.
\end{proposition}

[Necessity: Fix a countable dense set $\{f_n\}$ in $C(X)$ $-$then $\{f_n\}$ separates points of \mX and the initial topology on \mX determined by the $f_n$ is a separable metrizable topology.  Reason: 
The arrow $X \ra \R^\omega$ defined by the rule $x \ra \{f_n(x)\}$ is an embedding.

Sufficiency: Let $X_0$ stand for \mX equipped with a smaller separable metrizable topology.  
Embed $X_0$ in 
$[0,1]^\omega$.  
Fix a countable dense set $\{\phi_n\}$ in $C([0,1]^\omega)$ (cf. Proposition 16) and put 
$f_n = \restr{\phi_n}{X_0}$ $-$then the sequence $\{f_n\}$ is dense in $C(X_0)$, thus $C(X_0)$ is separable.  
Indeed, given a compact subset $K_0$ of $X_0$ and $f_0 \in C(X_0)$, $\exists$ $\phi_0 \in C([0,1]^\omega)$ : 
$\restr{\phi_0}{K_0} = \restr{f_0}{K_0}$ $\&$ $\forall \ \epsilon > 0$, $\exists$ $\phi_n$ : $p_{K_0}(\phi_n - \phi_0) < \epsilon$ 
$\implies$ $p_{K_0}(f_n - f_0) < \epsilon$.  
Finally, the separability of $C(X_0)$ forces the separability of $C(X)$.  This is because a compact subset \mK of \mX is a compact subset of $X_0$ and the two topologies induce the same topology on \mK.]\\

Example: Take $X = \R$ (discrete topology) $-$then $C(X)$ is separable.\\

\begingroup%%----------------------------------->>
\fontsize{9pt}{11pt}\selectfont
\textbf{\small EXAMPLE} \  
If $X = \ds\bigcup\limits_n K_n$, where each $K_n$ is compact and metrizable, then $C(X)$ is separable.
\\ \indent
[There is no loss in generality in supposing that $K_n \subset K_{n+1}$ $\forall \ n$.  Choose a countable dense subset 
$\{f_{n,m}\}$ in $C(K_n)$ (cf. Proposition 16) and let $F_{n,m}$ be a continuous extension of $f_{n,m}$ to \mX $-$then the initial topology on \mX determined by the $F_{n,m}$ is a separable metrizable topology which is smaller than the given topology on \mX, so $C(X)$ is separable (cf. Proposition 17).]\\
\endgroup %%------------------------------------<<

\begingroup%%----------------------------------->>
\fontsize{9pt}{11pt}\selectfont
\textbf{\small FACT} \  
Let \mX be a LCH space $-$then $C(X)$ is separable and metrizable iff \mX is separable and metrizable.\\
\endgroup %%------------------------------------<<

\begingroup%%----------------------------------->>
\fontsize{9pt}{11pt}\selectfont
\textbf{\small FACT} \  
Let \mX be a LCH space $-$then $C(X)$ is separable and completely metrizable iff \mX is separable and completely metrizable.\\
\endgroup %%------------------------------------<<

%^
%^
%%----------------------------------------------------------------------------------------------18
\begin{proposition} \ %18
$C(X)$ is first countable iff \mX is countable at infinity.\\
\end{proposition}
%^

%^
\begin{proposition} \ %19
$C(X)$ is second countable iff \mX is countable at infinity and all the compact subsets of \mX are metrizable.
\end{proposition}

[Necessity: $C(X)$ second countable $\implies$ $C(X)$ first countable $\implies$ \mX countable at infinity (cf. Proposition 18).  
In addition, $C(X)$ second countable $\implies$ $C(X)$ separable.  
So, by Proposition 17, \mX admits a smaller separable metrizable topology which, however, induces the same topology on each compact subset of \mX.

Sufficiency: The hypotheses on \mX guarantee that $C(X)$ is separable (via the example above) and metrizable (cf. Proposition 14).]\\

\begingroup%%----------------------------------->>
\fontsize{9pt}{11pt}\selectfont
\textbf{\small EXAMPLE} \  
Let \mE be an infinite dimensional locally convex topological vector space.  
Assume: \mE is second countable and completely metrizable $-$then the Anderson-Kadec theorem says that \mE is homeomorphic to $\R^\omega$ 
(for a proof, see 
Bessaga-Pelczy\'nski\footnote[2]{\textit{Selected Topics in Infinite Dimensional Topology}, PWN (1975), \textbf{189}.\vspace{0.11cm}}).  
Consequently, if \mX is countable at infinity and compactly generated and if all the compact subsets of \mX are metrizable, then $C(X)$ is homeomorphic to $\R^\omega$.\\
\endgroup %%------------------------------------<<

\begingroup%%----------------------------------->>
\fontsize{9pt}{11pt}\selectfont
\textbf{\small FACT} \  
Suppose that \mX is second countable $-$then $C(X)$ is Lindel\"of.\\
\endgroup %%------------------------------------<<

Up until this point, the playoff between \mX and $C(X)$ has been primarily ``topological'', little use having been made of the fact that $C(X)$ 
is also a locally convex topological vector space.  
It is thus only natural to ask: Can one characterize those \mX for which $C(X)$ has a certain additional property (e.g., barrelled or bornological)?  
While this theme has generated an extensive literature, I shall present just two results, namely Propositions 20 and 21, these being due independently to 
Nachbin\footnote[3]{\textit{Proc. Nat. Acad. Sci. U.S.A.} \textbf{40} (1954), 471-474.\vspace{0.11cm}}
and 
Shirota\footnote[6]{\textit{Proc. Japan Acad. Sci.} \textbf{30} (1954), 294-298.}.\\

\begingroup%%----------------------------------->>
\fontsize{9pt}{11pt}\selectfont
\textbf{\small FACT} \  
$C(X)$ is reflexive iff \mX is discrete.

[Assuming that $C(X)$ is reflexive, its bounded weakly closed subsets are weakly compact.  Therefore the compact subsets of \mX are finite which means that $C(X)$ is a dense subspace of $\R^X$ (product topology).  
But the reflexiveness of $C(X)$ also implies that its closed bounded subsets are complete, hence \mX is a $k_{\R}$-space 
(cf. p. 
\pageref{2.13}).  
Thus $C(X)$ is complete (cf. Proposition 13), so $C(X) = \R^X$ and \mX is discrete.]\\
\endgroup %%------------------------------------<<

A subset \mA of \mX is said to be 
\un{bounding}
\index{bounding (subset)} 
if every $f \in C(X)$ is bounded on \mA.  
Example: \mX is pseudocompact iff \mX is bounding.\\

%%----------------------------------------------------------------------------------------------19
Given a subset \mW of $C(X)$, let $K(W)$ be the subset of \mX consisting of those $x$ with the property that for every neighborhood $O_x$ 
of $x$ there exists an $f \in C(X): f(X - O_x) = \{0\}$ $\&$ $f \notin W$.\\

\index{Bounding Lemma}
\textbf{\small BOUNDING LEMMA} \  
If \mW is a barrel in $C(X)$, then $K(W)$ is bounding.

[Suppose that $K(W)$ is not bounding and fix an infinite discrete collection $\sO = \{O\}$ of open subsets of \mX such that 
$O \ \cap \  K(W) \neq \emptyset$ \ $\forall \ O \in \sO$.  
Choose an element $O_1 \in \sO$.  
Since $O_1 \ \cap \  K(W) \neq \emptyset$,  
$\exists$ $f_1 \in C(X)$ : $f_1(X - O_1) = \{0\}$ $\&$ $f_1 \notin W$.  
On the other hand, \mW, being a barrel, is closed, so $\exists$ a compact $K_1 \subset X$ and a positive $\epsilon_1$ : 
$\{g:p_{K_1}(f_1 - g) < \epsilon_1\} \cap W = \emptyset$.  
Choose next an element $O_2 \in \sO$ : 
$O_2 \cap K_1 = \emptyset$ and continue.  
The upshot is that there exist sequences 
$\{O_n\}$, $\{f_n\}$, $\{K_n\}$, $\{\epsilon_n\}$ with the following properties: 
(1) $O_{n+1} \ \cap \  \bigl(\ds\bigcup\limits_{i = 1}^n K_i\bigr) = \emptyset$; 
(2) $f_n(X - O_n) = \{0\}$ $\&$ $f_n \notin W$; 
(3) $\{g:p_{K_n}(f_n - g) < \epsilon_n\} \cap W = \emptyset$.  
Take $c_1 = 1$ and determine 
$c_{n+1}$ : $0 < c_{n+1} < \ds\frac{1}{n+1}$, subject to the requirement that 
$c_{n+1} p_{K_{n+1}} \ds\big( \sum\limits_{i = 1}^n \frac{1}{c_i} f_i \big) < \epsilon_{n+1}$ $\forall$ $n$.  \ 
Put $f = \ds\sum\limits_{i = 1}^\infty \ds\frac{1}{c_i} f_i$ $-$then by 
(2) and the discreteness of $\{O_n\}$, $f$ is continuous, 
and (1) - (3) combine to imply that $c_{n+1} f \notin W$ $\forall \ n$, 
thus \mW does not absorb the function $f$, a contradiction.]\\

\index{Lemma of Determination}
\textbf{\small LEMMA OF DETERMINATION} \  
If \mW is a barrel in $C(X)$ and if $f$ is an element of $C(X)$ such that $f(x) = 0$ $\forall \ x \in U$, where $U$ is an open set containing $K(W)$, then $f \in W$.

[Suppose false.  Choose a compact $K \subset X$ and a positive $\epsilon$ : 
$\{g:p_K(f - g) < \epsilon\} \ \cap \  W = \emptyset$, and for each $x \in K - U$, choose a neighborhood $O_x$ of $x$ : 
$g(X - O_x) = \{0\}$ $\implies$ 
$g \in W$.  Fix $f_x \in C(X,[0,1])$ : $f_x(x) = 1$ $\&$ $\restr{f_x}{X - O_x} = 0$, and let 
$U_x = \{y: f_x(y) > 1/2\}$.  
The $U_x$ comprise an open covering of $K - U$, thus one can extract a finite subcovering $U_{x_1}, \ldots, U_{x_n}$.  
Put 
$\ds \kappa_{x_i} = \ds\frac{f_{x_i}}{\max\{1/2, f_{x_1} + \cdots + f_{x_n}\}}$ 
$(i = 1, \ldots, n)$ $-$then 
$\ds\sum\limits_{i = 1}^n \restr{\kappa_{x_i}}{K - U} = 1$.  
Since $\kappa_{x_i}(X - O_{x_i}) = \{0\}$, 
$c\kappa_{x_i} f \in W$ $(c \in \R)$, therefore 
$F = \kappa_{x_1}f + \cdots + \kappa_{x_n}f =$ 
$\ds\frac{1}{n}(n\kappa_{x_1}f + \cdots + n\kappa_{x_n}f) \in W$.  But by its very construction, 
$\restr{F}{K} = \restr{f}{K} \implies F \notin W$.]\\

\begin{proposition} \ %20
$C(X)$ is barrelled iff every bounding subset of \mX is relatively compact.
\end{proposition}

[Necessity: Rephrased, the assertion is that for any closed noncompact subset \mS of \mX, $\exists$ $f \in C(X)$ : $f$ is unbounded on \mS.  
Thus let $B_S = \{f:\sup\limits_S \abs{f} \leq 1\}$ $-$then $B_S$ is balanced and convex.  
Since $B_S$ is also closed and since the requirement that there be some $f \in C(X)$ which is unbounded on \mS amounts to the failure of $B_S$ to be absorbing, it need only be shown that $B_S$ does not contain a neighborhood of 0.  
Assuming the opposite, choose a compact \mK and a positive $\epsilon$
: $\{f:p_K(f) < \epsilon\} \subset B_S$.  
Claim: $S \subset K$.  Proof: 
%%----------------------------------------------------------------------------------------------20
If $x \in S - K$, $\exists$ $f \in C(X)$ : $f(K) = \{0\}$ $\&$ $f(x) = 2$, an impossibility.  Therefore \mS is compact (being closed), contrary to hypothesis.

Sufficiency: Fix a barrel \mW in $C(X)$ $-$then the contention is that \mW contains a neighborhood of 0.  
Owing to the bounding lemma, $K(W)$ is compact (inspect the definitions to see that $K(W)$ is closed).  
Accordingly, it suffices to produce a positive $\epsilon$ : $\{f:p_{K(W)}(f) < \epsilon\} \subset W$.  
To this end, consider $BC(X)$ viewed as a Banach space in the supremum norm.  
Because $BC(X)$ barrelled and $W \ \cap \  BC(X)$ is a barrel in $BC(X)$, $\exists$ $\epsilon > 0$ : 
$\norm{\phi} \leq 2\epsilon$ $\implies$ $\phi \in W$ ($\phi \in BC(X)$).  
Assuming that $p_{K(W)}(f) < \epsilon$, fix an open set \mU containing $K(W)$ such that 
$\abs{f(x)} < \epsilon$ $\forall \ x \in U$.  
Let 
$F(x) = \max\{\epsilon,f(x)\} + \min\{-\epsilon,f(x)\}$ $-$then $2F(x) = 0$ $(x \in U)$, thus the lemma of determination implies that $2F \in W$.  
But $\forall \ x \in X$, $\abs{2(f(x) - F(x))} < 2\epsilon$ $\implies$ 
$\norm{2(f - F)} \leq 2\epsilon$ $\implies$ $2(f - F) \in W$, so 
$\ds\frac{1}{2}(2F) + \ds\frac{1}{2}(2(f - F)) \in W$, 
i.e., $f \in W$.]\\

Example: $C([0,\Omega[)$ is not barrelled.\\

\begingroup%%----------------------------------->>
\fontsize{9pt}{11pt}\selectfont
\textbf{\small EXAMPLE} \  
If \mX is a paracompact LCH space, then $C(X)$ is Baire 
(cf. p. 
\pageref{2.14}
).  Since Baire $\implies$ barrelled, it follows from Proposition 20 that the bounding subsets of \mX are relatively compact.\\
\endgroup %%------------------------------------<<

Notation: Every $f \in C(X)$ can be regarded as an element of $C(X,\R_\infty)$, hence admits a unique continuous extension $f_\infty:\beta X \ra \R_\infty$.

[Note: \ Put $v_f X = \{x \in \beta X: f_\infty(x) \in \R\}$ $-$then the intersection 
$\ds\bigcap\limits_{f \in C(X)} v_f X$ is $vX$.]\\

\begingroup%%----------------------------------->>
\fontsize{9pt}{11pt}\selectfont
\label{2.16}
\label{2.15a}
\textbf{\small FACT} \  
The elements of $\beta X - vX$ are those $x$ with the property that there exists a $G_\delta$ in $\beta X$ containing $x$ which does not meet \mX.\\
\endgroup %%------------------------------------<<

Let \mW be a balanced, convex subset of $C(X)$ $-$then \mW is said to 
\un{contain a ball}
\index{contain a ball (convex subset)} 
if $\exists$ $r > 0$: 
$\{f:\underset{X}{\sup} \abs{f} \leq r\} \subset W$.

Example: Every balanced, convex bornivore \mW in $C(X)$ contains a ball.

[Given $f, \ g \in C(X)$ with $f \leq g$, let $[f,g] = \{\phi: f \leq \phi \leq g\}$.  
Since $\forall$ compact $K \subset X$, 
$p_K(\phi) \leq \max\{p_K(f),p_K(g)\}$, $[f,g]$ is bounded, thus is absorbable by \mW.  
In particular: $\exists$ $r > 0$ such that $[-r,r] \subset W$.]\\
%dmc both originals have $[-r1,r1] \subset W$

\begingroup%%----------------------------------->>
\fontsize{9pt}{11pt}\selectfont
\label{2.17}
\textbf{\small FACT} \  
Suppose that \mW contains a ball.  Let \mK be a compact subset of \mX.  
Assume: $f(K) = \{0\}$ $\implies$ $f \in W$ $-$then 
$\exists$ $\epsilon > 0$: $\{f:p_{K}(f) < \epsilon\} \subset W$.\\
\endgroup %%------------------------------------<<

Let \mW be a balanced, convex subset of $C(X)$ $-$then a compact subset \mK of $\beta X$ is said to be a 
\un{hold}
\index{hold (balanced, convex subset)} 
of \mW if $f \in W$ whenever $f_\infty(K) = \{0\}$.  
Example: $\beta X$ is a hold of \mW.\\

%%----------------------------------------------------------------------------------------------21
\textbf{\small LEMMA} \  
Suppose that \mW contains a ball $-$then a compact subset \mK of $\beta X$ is a hold of \mW provided that $f \in W$ whenever $f_\infty$ vanishes on some open subset \mO of $\beta X$ containing \mK.\\

Application: Under the assumption that \mW contains a ball, if \mK and $L$ are holds of \mW, then so is $K \cap L$.

[Consider any $f:f_\infty(O) = \{0\}$, where $O$ is some open subset of $\beta X$ containing $K \cap L$.  
Choose disjoint open subsets $U, \ V$ of $\beta X$ : $K \subset U$, $L - O \subset V$ and let $U^\prime$, $V^\prime$ be open subsets of $\beta X$: $K \subset U^\prime \subset \ov{U}^\prime$ $\subset U$, 
$L - O \subset V^\prime$ $\subset \ov{V}^\prime$ $\subset V$.  
Fix $\phi \in C(X,[0,1])$ : 
$\beta \phi(\ov{U}^\prime) = \{1\}$, 
$\beta \phi(\ov{V}^\prime) = \{0\}$.  
Note that $2f \phi$ vanishes on $(O \cup V^\prime) \cap X$.  
But $O \cup V^\prime \subset \ov{(O \cup V^\prime) \cap X}$ 
$\implies$ $(2 f \phi)_\infty (O \cup V^\prime) = \{0\}$.  
On the other hand, $L \subset O \cup V^\prime$, thus by the lemma, $2 f \phi \in W$.  
Similarly, $2 f (1 - \phi) \in W$.  Therefore 
$\ds f = \frac{1}{2}(2 f \phi) + \frac{1}{2}(2 f (1 - \phi)) \in W$.]\\

Let \mW be a balanced, convex subset of $C(X)$ $-$then the 
\un{support}
\index{support (balanced, convex subset of $C(X)$)} 
of \mW, 
written $\spt W$, is the intersection of all the holds of \mW.\\

\textbf{\small LEMMA} \  
Suppose that \mW contains a ball $-$then $\spt W$ is a hold of \mW.

[Since $\beta X$ is a compact Hausdorff space, for any open $O \subset \beta X$ containing  $\spt W$, 
$\exists$ holds $K_1, \ldots K_n$ of $W$ such that $\bigcap\limits_{i = 1}^n K_i \subset O$.]\\

\begin{proposition} \ %21
$C(X)$ is bornological iff \mX is $\R$-compact.
\end{proposition}

[Necessity: Assuming that \mX is not $\R$-compact, fix a point $x_0 \in vX - X$ $-$then 
the assignment 
$f \ra f_\infty(x_0)$ defines a nontrivial homomorphism 
$\widehat{x}_0:C(X) \ra \R$, which is necessarily discontinuous 
(cf. p. \pageref{2.15}).  
So, to conclude that $C(X)$ is not bornological, it suffices to show that $\widehat{x}_0$ takes bounded sets to bounded sets.  
If this were untrue, then there would be a bounded set $B \subset C(X)$ and a sequence $\{f_n\} \subset B$ such that 
$\widehat{x}_0(f_n) \ra \infty$.  
The intersection $\ds\bigcap\limits_n \{x \ \in \beta X: (f_n)_\infty(x) > (f_n)_\infty(x_0) - 1\}$ is a 
$G_\delta$ in $\beta X$ containing $x_0$, thus it must meet \mX 
(cf. p. \pageref{2.15a}) 
say at $x_{00}$ hence $f_n(x_{00}) \ra \infty$.  
But then, as \mB is bounded, $\ds\frac{f_n}{f_n(x_{00})} \ra 0$ in $C(X)$, which is nonsense.

Sufficiency: 
It is a question of proving that every balanced, convex bornivore \mW in $C(X)$ contains a neighborhood of 0.  
Because \mW contains a ball, the lemma implies that $\spt W$ is a hold of \mW, thus the key is to establish the containment 
$\spt W \subset X$ since this will allow one to say that $\exists$ $\epsilon > 0$ : $\{f: p_{\spt W} (f) < \epsilon\} \subset W$ 
(cf. p. \pageref{2.16}).  
So take a point $x_0 \in \beta X - X$ and choose closed subsets $A_1 \supset A_2 \supset \cdots$ of 
$\beta X$ : $\forall \ n$, $x_0 \in \itr A_n$ $\&$ 
$\bigl(\bigcap\limits_n A_n \bigr) \cap X = \emptyset$ (possible, \mX being 
$\R$-compact 
(cf. p. \pageref{2.17})).  
Claim: At least one of the $\beta X - \itr A_n$ is a hold of \mW 
($\implies$ $x_0 \notin \spt W$ $\implies$ $\spt W \subset X$).  
If not, then $\forall \ n$, 
%%----------------------------------------------------------------------------------------------22
$\exists$ $f_n$ : $(f_n)_\infty(\beta X - \itrx A_n) = 0$ $\&$ $f_n \notin W$.  
The sequence $\{X - A_n\}$ is an increasing sequence of open subsets of \mX whose union is \mX.  
Therefore $f = \sup\limits_n n\abs{f_n}$ is in $C(X)$.  
Fix 
$d > 0$ : $[-f,f] \subset dW$ $-$then 
$nf_n \in dW$ $\forall \ n$ $\implies$ $f_n \in W$ $\forall \ n \geq d$, a contradiction.]\\

\begingroup%%----------------------------------->>
\fontsize{9pt}{11pt}\selectfont
\textbf{\small LEMMA} \  
A subset \mA of \mX is bounding iff its closure $\beta X$  is contained in $v X$.\\
\endgroup %%------------------------------------<<

\begingroup%%----------------------------------->>
\fontsize{9pt}{11pt}\selectfont
\textbf{\small FACT} \  
If $C(X)$ is bornological, then $C(X)$ is barrelled.

[Note: \ Recall that in general, bornological $\centernot\implies$ barrelled and barrelled $\centernot\implies$ bornological.]\\
\endgroup %%------------------------------------<<

\begingroup%%----------------------------------->>
\fontsize{9pt}{11pt}\selectfont
Remark: There are completely regular Hausdorff spaces \mX whose bounding subsets are relatively compact but that are not $\R$-compact 
(Gillman-Henriksen\footnote[2]{\textit{Trans. Amer. Math. Soc.} \textbf{77} (1954), 340-362 (cf. 360-362).}).  
For such \mX, $C(X)$ is therefore barrelled but not bornological.\\
\endgroup %%------------------------------------<<

Given a closed subset \mA of \mX, let $I_A = \{f:\restr{f}{A} = 0\}$ $-$then $I_A$ is a closed ideal in $C(X)$.  
Examples: 
(1) \ $I_\emptyset = C(X)$; 
(2) \ $I_X = \{0\}$.\\

\textbf{\small SUBLEMMA} \  
Suppose that \mX is compact.  Let $I \subset C(X)$ be an ideal.  
Assume: $\forall \ x \in X$, $\exists$ $f_x \in I$: 
$f_x(x) \neq 0$ $-$then $I = C(X)$.

[$\forall \ x \in X$, $\exists$ a neighborhood $U_x$ of $x$ : $\restr{f_x}{U_x} \neq 0$.  
Choose points $x_1, \ldots, x_n$ : 
$X = \bigcup\limits_{i = 1}^n U_{x_i}$ and let $f = \sum\limits_{i = 1}^n f_{x_i}^2$ : $f \in I$ $\implies$ 
$1 = f \cdot \ds\frac{1}{f} \in I$ $\implies$ $I = C(X)$.]\\

\textbf{\small LEMMA} \ 
Suppose that \mX is compact.  Let $I \subset C(X)$ be an ideal and put $A = \bigcap\limits_{f \in I} Z(f)$.  
Assume: $A \subset U \subset Z(\phi)$, where $U$ is open and $\phi \in C(X)$ $-$then $\phi \in I$.

[The restriction $\restr{I}{X - U}$ is an ideal in $C(X - U)$ (Tietze), hence by the sublemma, equals $C(X - U)$.  Choose an 
$f \in I$ : $\restr{f}{X - U} = 1$ to get $\phi = f\phi \in I$.]\\

\begin{proposition} \ %22
Suppose that \mX is compact.  Let $I \subset C(X)$ be an ideal $-$then $\ov{I} = I_A$, where 
$A = \bigcap\limits_{f \in I} Z(f)$.
\end{proposition}

[Since $I \subset I_A$, it need only be shown that $I_A \subset \ov{I}$.  So let $f$ be a nonzero element of $I_A$ and fix $\epsilon > 0$.  
Choose $\phi \in C(X,[0,1])$ : $\{x:\abs{f(x)} \leq \epsilon/2\} \subset  Z(\phi)$ $\&$ 
$\{x:\abs{f(x)} \geq 3\epsilon/4\} \subset Z(1 - \phi)$.  
Because $A \subset \{x:\abs{f(x)} < \epsilon/4\} \subset Z(f\phi)$, the lemma gives $f \phi \in I$.  
And: $\norm{f - f\phi} = \sup\limits_X \abs{f - f\phi} < \epsilon$ $\implies$ $f \in \ov{I}$.]\\

\begin{proposition} \ %23
The closed subsets of \mX are in a one-to-one correspondence with the closed ideals of $C(X)$ via $A \ra I_A$.
\end{proposition}

%%----------------------------------------------------------------------------------------------23
[Due to the complete regularity of \mX, the map $A \ra I_A$ is injective.  
To see that it is surjective, it suffices to prove that for any closed ideal \mI in $C(X)$ : $I = I_A$, where 
$A = \bigcap\limits_{f \in I} Z(f)$.  
Obviously, $I \subset I_A$.  
On the other hand, $\forall$ compact $K \subset X$, the restriction $\restr{I}{K}$ is an ideal in $C(K)$ 
(cf. p. 
\pageref{2.18}), thus 
$\ov{I/K} = I_{A \cap K}$ (cf. Proposition 22), and from this it follows that $I_A \subset \ov{I} = I$.]\\

Application: \ The points of \mX are in a one-to-one correspondence with the closed maximal ideals of $C(X)$ via 
$x \ra I_{\{x\}}$.\\

\begingroup%%----------------------------------->>
\fontsize{9pt}{11pt}\selectfont
 By comparison, recall that the points of $\beta X$ are in a one-to-one correspondence with the maximal ideals of 
$C(X)$.
\\ \indent
[Note: \ Assign to each $x \in \beta X$ the subset $m_x$ of $C(X)$ consisting of those $f$ such that 
$x \in \cl_{\beta X}(Z(f))$ $-$then $m_x$ is a maximal ideal and all such have this form.  For details, see 
Walker\footnote[2]{\textit{The Stone-\u Cech Compactification}, Springer Verlag (1974), 18.}.]\\
\endgroup %%------------------------------------<<

A 
\un{character}
\index{character (of $C(X)$} 
of $C(X)$ is a nonzero multiplicative linear functional on $C(X)$, i.e., a homomorphism $C(X) \ra \R$ of algebras.\\

\textbf{\small LEMMA} \  
If $\chi:\R \ra \R$ is a nonzero ring homomorphism, then $\chi = \id_{\R}$.

[In fact, $\chi$ is order preserving and the identity on $\Q$.]\\

Application: Every ring homomorphism $C(X) \ra \R$ is $\R$-linear, thus is a character.\\

\textbf{\small LEMMA} \  
If $\chi:C(X) \ra \R$ is a character of $C(X)$, then $\forall \ f$, $\abs{\chi(f)} = \chi(\abs{f})$.

[For $\abs{\chi(f)}^2 = \chi(f)^2 = \chi(f^2) = \chi(\abs{f}^2) = \chi(\abs{f})^2$ and $\chi(\abs{f})$ is $\geq 0$.]\\

\begingroup%%----------------------------------->>
\fontsize{9pt}{11pt}\selectfont
By way of a corollary, if $\chi:C(X) \ra \R$ is a character of $C(X)$ and if $\chi(f) = 0$, then 
$\chi(\min\{1,\abs{f}\}) = 0$.  Proof: $2 \chi(\min\{1,\abs{f}\}) = \chi(1) + \chi(f) - \chi(\abs{1 - f}) =$ 
$1 - \abs{\chi(1 - f)} =$ 
$1 - 1 = 0$.\\
\endgroup %%------------------------------------<<

\begingroup%%----------------------------------->>
\fontsize{9pt}{11pt}\selectfont
\textbf{\small FACT} \  
Write $vf$ for the unique extension of $f \in C(X)$ to $C(vX)$ $-$then $C(X)$ ``is'' $C(vX)$ and the characters of $C(X)$ 
are parameterized by the points of $vX$: $f \ra vf(x)$ $(x \in vX)$.
\\ \indent
[If \mX is $\R$-compact and if $\chi:C(X) \ra \R$ is a character, then in the terminology of 
p. \pageref{2.19} 
$\&$ p. 
\pageref{2.20}, 
$\sF_\chi = \{Z(f):\chi(f) = 0\}$ is a zero set ultrafilter on \mX.  
Claim $\sF_\chi$ has the countable intersection property.  
Thus let $\{Z(f_n)\} \subset \sF_\chi$ be a sequence and put 
$\ds f = \sum\limits_1^\infty \frac{\min\{1,\abs{f_n}\}}{2^n}$ $-$then 
$\ds \bigcap\limits_1^\infty Z(f_n) = Z(f)$.  
To prove that $\chi(f) = 0$, write 
$\ds f = \sum\limits_{i=1}^n \frac{\min\{1,\abs{f_i}\}}{2^i} + g_n$, where 
$0 \leq g_n \leq 2^{-n}$, apply $\chi$ to get 
$\chi(f) = 
%%----------------------------------------------------------------------------------------------24
\chi(g_n) \leq 2^{-n}$, and let $n \ra \infty$.  
It therefore follows that $\cap \ \sF_\chi$ is nonempty, say $x \in \cap \ \sF_\chi$ 
(cf. p. \pageref{2.21}).  
And: $\chi(f - \chi(f)) = 0$ $\implies$ $x \in Z(f - \chi(f))$ $\implies$ $\chi(f) = f(x)$.]\\
\endgroup %%------------------------------------<<

Notation: $\widehat{C(X)}$ is the set of continuous characters of $C(X)$.

From the above, there is a one-to-one correspondence $X \ra \widehat{C(X)}$, viz. $x \ra \chi_x$, where 
$\chi_x(f) = f(x)$.\\

\begingroup%%----------------------------------->>
\fontsize{9pt}{11pt}\selectfont
\label{2.15}
If \mX is not $\R$-compact, then the elements of $vX - X$ correspond to the discontinuous characters of $C(X)$.\\
\endgroup %%------------------------------------<<

Topologize $\widehat{C(X)}$  by giving it the initial topology determined by the functions $\chi \ra \chi(f)$ $(f \in C(X))$ $-$then the correspondence $X \ra \widehat{C(X)}$ is a homeomophism (cf. $\S 1$, Proposition 14).\\

\begin{proposition} \ %24
Let 
$
\begin{cases}
\ X\\[-.1cm]
\ Y
\end{cases}
$
be CRH spaces.  Assume: 
$
\begin{cases}
\ C(X)\\[-.1cm]
\ C(Y)
\end{cases}
$
are isomorphic as topological algebras $-$then 
$
\begin{cases}
\ X\\[-.1cm]
\ Y
\end{cases}
$
are homeomorphic.
\end{proposition}

[Schematically, 
\begin{tikzcd}[sep=large]
{X}   \arrow[d,shift right=0.5,dash] \arrow[d,shift right=-0.5,dash] 
&{Y} \arrow[d,shift right=0.5,dash] \arrow[d,shift right=-0.5,dash] \\
{\widehat{C(X)}} \ar{r}  &{\widehat{C(Y)}} \ar{l}
\end{tikzcd}
and 
%$\leftrightarrow$ 
\begin{tikzcd}[ sep=small]
{} \ar{rr}  &&{} \ar{ll}
\end{tikzcd}
is a homeomorphism.]\\

\vspace{0.75cm}
\begingroup%%----------------------------------->>
\fontsize{9pt}{11pt}\selectfont
\textbf{\small FACT} \  
Let 
$
\begin{cases}
\ X\\
\ Y
\end{cases}
$
be CRH spaces.  Assume: 
$
\begin{cases}
\ C(X)\\
\ C(Y)
\end{cases}
$
are isomorphic as algebras $-$then 
$
\begin{cases}
\ vX\\
\ vY
\end{cases}
$
are homeomorphic.\\
\endgroup %%------------------------------------<<
%%%%%%%%%%%%%%%%%%%%%%%%%%%%%%%%%%%%%%

\begin{center}
$\S \ 2$
\\[0.5cm]
$\mathcal{REFERENCES}$\\
\end{center}

\[
\textbf{BOOKS}
\]

\begingroup
\fontsize{9pt}{11pt}\selectfont
\setlength\parindent{0 cm}

[1] \quad Arhangel'ski\u i, A., \textit{Topological Function Spaces}, Kluwer (1992).
\\[-.2cm]

[2] \quad Beckenstein, E., Narici, L., and Suffel, C., \textit{Topological Algebras}, North Holland (1977).
\\[-.2cm]

[3] \quad Burckel, R., \textit{Characterizations of C(X) among its Subalgebras}, Marcel Dekker (1972).
\\[-.2cm]

[4] \quad Carreras, P. and Bonet, J., \textit{Barrelled Locally Convex Spaces}, North Holland (1987).
\\[-.2cm]

[5] \quad Gierz, G. et al., \textit{A Compendium of Continuous Lattices}, Springer Verlag (1980).
\\[-.2cm]

[6] \quad Gillman, L. and Jerison, M., \textit{Rings of Continuous Functions}, Van Nostrand (1960).
\\[-.2cm]

[7] \quad Hirsch, M., \textit{Differential Topology}, Springer Verlag (1976).
\\[-.2cm]

[8] \quad Howes, N., \textit{Modern Analysis and Topology}, Springer Verlag (1995).
\\[-.2cm]

[9] \quad Kirby, R. and Siebenmann, L., \textit{Foundational Essays on Topological Manifolds, Smoothings,}

\hspace{.65cm} \textit{and Triangulations}, 
Princeton University Press (1977).
\\[-.2cm]

[10] \quad Margalef-Roig, J. and Outerelo Dominguez, E., \textit{Differential Topology}, North Holland (1992).
\\[-.2cm]

[11] \quad McCoy, R. and Ntantu, I., \textit{Topological Properties of Spaces of Continuous Functions}, Springer Verlag 

\hspace{.9cm}
(1988).
\\[-.2cm]

[12] \quad Schmets, J., \textit{Espaces de Fonctions Continues}, Springer Verlag (1976).
\\[-.2cm]

[13] \quad Semadeni, Z., \textit{Banach Spaces of Continuous Functions}, PWN (1971).
\\[-.2cm]

[14] \quad Todorcevic, S., \textit{Topics in Topology}, Springer Verlag (1997).
\\[-.2cm]

[15] \quad Walker, R., \textit{The Stone-\u Cech Compactification}, Springer Verlag (1974).
\\[-.2cm]

[16] \quad Weir, M., \textit{Hewitt-Nachbin Spaces}, North Holland (1975).
\\[-.2cm]

[17] \quad Wilansky, A., \textit{Topology for Analysis}, Ginn and Company (1970).
\\[-.2cm]

[18] \quad Wilansky, A., \textit{Modern Methods in Topological Vector Spaces}, McGraw-Hill (1978).
\\
\endgroup

\[
\textbf{ARTICLES}
\]

\begingroup
\fontsize{9pt}{11pt}\selectfont
\setlength\parindent{0 cm}

[1] \quad Arens, R., Topologies for Homeomorphism Groups, \textit{Amer. J. Math.} \textbf{68} (1946), 593-610.
\\[-.2cm]

[2] \quad Arens, R., A Topology for Spaces of Transformations, \textit{Ann. of Math.} \textbf{47} (1946), 480-495.
\\[-.2cm]

[3] \quad Arhangel'ski\u i, A., A Survey of $C_p$-Theory, \textit{Questions Answers Gen. Topology} \textbf{5} (1987), 1-109.
\\[-.2cm]

[4] \quad Arhangel'ski\u i, A., Spaces of Mappings and Rings of Continuous Functions, In:  \textit{General Topology}, EMS 

\hspace{.65cm}
\textbf{51}, Springer Verlag (1995), 71-156.
\\[-.2cm]

[5] \quad Blasco, J., On $\mu$-Spaces and $k_{\R}$-Spaces, \textit{Proc. Amer. Math. Soc.} \textbf{67} (1977), 179-186.
\\[-.2cm]

[6] \quad Bowers, P., Limitation Topologies on Function Spaces, \textit{Trans. Amer. Math. Soc.} \textbf{314} (1989), 421-431.
\\[-.2cm]

[7] \quad Collins, H., Strict, Weighted, and Mixed Topologies and Applications, \textit{Adv. Math.} \textbf{19} (1976), 207-237.
\\[-.2cm]

[8] \quad Day, B. and Kelly, G., On Topological Quotient Maps Preserved by Pullbacks or Products, 
\textit{Math.}

\hspace{.65cm}
\textit{ Proc. Cambridge Philos. Soc.} \textbf{67} (1970), 553-558.
\\[-.2cm]

[9] \quad Fitzpatrick, B. and Zhou, H-X., A Survey of Some Homogeneity Properties in Topology, \textit{Ann. New}

\hspace{.65cm}
\textit{York Acad. Sci.} 
\textbf{552} (1989), 28-35.
\\[-.2cm]

[10] \quad Hamstrom, M-E., Homotopy in Homeomorphism Spaces, TOP and PL, \textit{Bull. Amer. Math. Soc.} \textbf{80} 

\hspace{.9cm}
(1974), 207-230.
\\[-.2cm]

[11] \quad Isbell, J., General Function Spaces, Products and Continuous Lattices, \textit{Math. Proc. Cambridge} 

\hspace{.9cm}
\textit{Philos. Soc.} 
\textbf{100} (1986), 193-205.
\\[-.2cm]

[12] \quad Krikorian, N., A Note Concerning the Fine Topology on Function Spaces, \textit{Compositio Math.} \textbf{21} 

\hspace{.9cm}
(1969), 343-348.
\\[-.2cm]

[13] \quad Luukkainen, J., Approximating Continuous Maps of Metric Spaces into Manifolds by Embeddings, 

\hspace{.9cm}
\textit{Math. Scand.} \textbf{49} (1981), 61-85.
\\[-.2cm]

[14] \quad Michael, E., $\aleph_0$-Spaces, \textit{J. Math. Mech.} \textbf{15} (1966), 983-1002.
\\[-.2cm]

[15] \quad Myers, S., Spaces of Continuous Functions, \textit{Bull. Amer. Math. Soc.} \textbf{55} (1949), 402-407.
\\[-.2cm]

[16] \quad Vechtomov, E., Specification of Topological Spaces by Algebraic Systems of Continuous Functions, 

\hspace{.9cm}
\textit{J. Soviet Math.} \textbf{53} (1991), 123-147.
\\[-.2cm]

[17] \quad Vechtomov, E., Rings of Continuous Functions-Algebraic Aspects, \textit{J. Math. Sci.} \textbf{71} (1994), 2364-

\hspace{.9cm}
2408.
\\[-.2cm]

[18] \quad Vechtomov, E., Rings and Sheaves, \textit{J. Math. Sci.} \textbf{74} (1995), 749-798.
\\[-.2cm]

[19] \quad Warner, S., The Topology of Compact Convergence on Continuous Function Spaces, \textit{Duke Math. J.} 

\hspace{.9cm}
\textbf{25} (1958), 265-282.
\\[-.2cm]

[20] \quad Whittaker, J., On Isomorphic Groups and Homeomorphic Spaces, \textit{Ann. of Math.} \textbf{78} (1963), 74-91.

\setlength\parindent{2em}

\endgroup

\chapter{
$\boldsymbol{\S}$\textbf{3}.\quadx  COFIBRATIONS}
\setlength\parindent{2em}
\setcounter{proposition}{0}
%%----------------------------------------------------------------------------------------------01
$\text{ }$\\[-1.25cm]

The machinery assembled here is the indispensable technical prerequisite for the study of homotopy theory in \bTOP or $\bTOP_*$.

Let $X$ and $Y$ be topological spaces.  Let $A \ra X$ be a closed embedding and let $f:A \ra Y$ be a continuous function $-$then the 
\un{adjunction space}
\index{adjunction space} 
$X \sqcup_f Y$ corresponding to the 2-source 
\begin{tikzcd}[ sep=small]
X &A \arrow{l} \arrow{r}{f} &Y
\end{tikzcd}
 is defined by the pushout square 
\begin{tikzcd}
A \arrow{r}{f} \arrow{d} &Y \arrow{d}\\
X \arrow{r} &X \sqcup_f Y\\
\end{tikzcd}
, $f$ being the 
\un{attaching map}.
\index{attaching map}   
Agreeing to identify $A$  with its image in $X$, the restriction of the projection $p:X \coprod Y \ra X \sqcup_f Y$ to
$
\begin{cases}
\ X - A\\[-.2cm]
\ Y
\end{cases}
$
 is a homeomorphism of 
$
\begin{cases}
\ X - A\\[-.2cm]
\ Y
\end{cases}
$
 onto an 
$
\begin{cases}
\ \text{open}\\[-.2cm]
\ \text{closed}
\end{cases}
$
subset of $X \sqcup_f Y$ and the images 
$
\begin{cases}
\ p(X - A)\\[-.2cm]
\ p(Y)
\end{cases}
$ 
partition $X \sqcup_f Y$.

\label{19.17}
[Note: \ The adjunction space $X \sqcup_f Y$ is unique only up to isomorphism.  
For example, if $\phi:X \ra X$ is a homeomorphism such that $\restr{\phi}{A} = \id_A$, 
then there arises another pushout square equivalent to the original one.]

\label{6.31}
\label{18.12}
\indent\indent(AD$_1$)  If $A$ is not empty and if $X$ and $Y$ are connected (path connected), then $X \sqcup_f Y$ is connected (path connected).\\
\indent\indent(AD$_2$)  If $X$ and $Y$ are $\tT_1$, then $X \sqcup_f Y$ is $\tT_1$ but if $X$ and $Y$ are Hausdorff, then $X \sqcup_f Y$ need not be Hausdorff.\\
\indent\indent(AD$_3$)  If $X$ and $Y$ are Hausdorff and if $A$ is compact, then $X \sqcup_f Y$ is Hausdorff.\\
\indent\indent(AD$_4$)  If $X$ and $Y$ are Hausdorff and if $A$ is a neighborhood retract of $X$ such that each $x \in X - A$ has a neighborhood $U$ with $A \cap \overline{U} = \emptyset$, then $X \sqcup_f Y$ is Hausdorff.\\
\indent\indent(AD$_5$)  If \mX and $Y$ are normal (normal and countably paracompact, perfectly normal, collectionwise normal, paracompact) Hausdorff spaces, 
then $X \sqcup_f Y$ is a normal (normal and countably paracompact, perfectly normal, collectionwise normal, paracompact) Hausdorff space.\\
\label{14.4}
\label{14.69b}
\indent\indent(AD$_6$)  If $X$ and $Y$ are in \bCG, ($\dcg$), then $X \sqcup_f Y$ is in \bCG, ($\dcg$).\\

\begingroup
\fontsize{9pt}{11pt}\selectfont
\textbf{\small EXAMPLE} \quadx
Working with the Isbell-Mr\'owka space $\Psi(\N) = \sS \cup \N$, consider the pushout square
\begin{tikzcd}
\sS \arrow{r}{f} \arrow{d} &\beta\sS \arrow{d}\\
\Psi(\N) \arrow{r} &\Psi(\N) \sqcup_f \beta\sS\\
\end{tikzcd}
.   
Due to the maximality of $\sS$, every open covering of $\Psi(\N) \sqcup_f \beta \sS$ has a finite subcovering.  
Still, $\Psi(\N) \sqcup_f \beta \sS$ is not Hausdorff.\\
\endgroup
%dmc??? maximality of $\sS$ ???
%%----------------------------------------------------------------------------------------------02
 
The 
\un{cylinder functor}
\index{cylinder functor} 
$I$ is the functor 
$
I \ : \ 
\begin{cases}
\ \bTOP \ra \bTOP\\[-.2cm]
\ X \ra X \times [0,1]
\end{cases}
\hspace{-.2cm}, \ 
$
where $X \times [0,1]$ carries the product topology.  
There are embeddings 
$
i_t:
\begin{cases}
\ X \ra IX\\[-.2cm]
\ x \ra (x,t)
\end{cases}
(0 \leq t \leq 1)
$
and a projection 
$
p:
\begin{cases}
\ IX \ra X\\[-.2cm]
\ (x,t) \ra x
\end{cases}
\hspace{-.2cm}.
$
The 
\un{path space functor}
\index{path space functor} 
\mP is the functor 
$
P:
\begin{cases}
\ \bTOP \ra \bTOP\\[-.2cm]
\ X \ra C([0,1]),X)
\end{cases}
\hspace{-.2cm},
$
where $C([0,1]),X)$ carries the compact open topology.  
There is an embedding 
$
j:
\begin{cases}
\ X \ra PX\\[-.2cm]
\ x \ra j(x)
\end{cases}
\hspace{-.4cm},
$
with $j(x)(t) = x$, and projections 
$
p_t:
\begin{cases}
\ PX \ra X\\[-.2cm]
\ \sigma \ra p_t(\sigma)
\end{cases}
(0 \leq t \leq 1),
$
with $p_t(\sigma) = \sigma(t)$.  
$(I,P)$ is an adjoint pair: $C(IX,Y) \approx C(X,PY)$.  
Accordingly, two continuous functions 
$
\begin{cases}
\ f:X \ra Y\\[-.2cm]
\ g:X \ra Y
\end{cases}
$
determine the same morphism in \bHTOP, i.e., are homotopic $(f \simeq  g)$, iff $\exists$ $H \in C(IX,Y)$ such that 
$
\begin{cases}
\ H \circ i_0 = f\\[-.2cm]
\ H \circ i_1 = g
\end{cases}
$
or equivalently, iff $\exists$ $G \in C(X,PY)$ such that 
$
\begin{cases}
\ p_0 \circ G = f\\
\ p_1 \circ G = g
\end{cases}
\hspace{-.2cm}.
$
\\

Let $A$ and $X$ be topological spaces $-$then a continuous function $i: A \ra X$ is said to be a 
\un{cofibration}
\index{cofibration} 
if it has the following property:  
Given any topological space $Y$ and any pair $(F,h)$ of continuous functions 
$
\begin{cases}
\ F: X \ra Y \\[-.2cm]
\ h: IA \ra Y
\end{cases}
$
such that $F \circ i = h \circ i_0$, there is a continuous function $H:IX \ra Y$ such that 
$F = H \circ i_0$ and $H \circ I i = h$.  
Thus \mH is a filler for the diagram
\[
\begin{tikzcd}
A \arrow{rr}{i} \arrow{dd}[swap]{i_0} &&X \arrow{dl} \arrow{dd}{i_0}\\
&Y\\
IA \arrow{rr}[swap]{I i} \arrow{ru} &&IX \arrow[lu,dashed]
\end{tikzcd}
.
\]

[Note: \ One can also formulate the definition in terms of the path space functor, viz.\\
\[
\begin{tikzcd}
A \arrow{rr} \arrow{dd}[swap]{i} &&PY  \arrow{dd}{p_0}\\
\\
X \arrow{rr} \arrow[rruu,dashed] &&Y
\end{tikzcd}
.]
\]

%%----------------------------------------------------------------------------------------------03
A continuous function $i :A \ra X$ is a cofibration iff the commutative diagram
\begin{tikzcd}
{A} \arrow{r}{i} \arrow{d}[swap]{i_0} &{X}  \arrow{d}{i_0}\\
{IA} \arrow{r}[swap]{I i} &{IX}
\end{tikzcd}
is a weak pushout square.  Homeomorphisms are cofibrations.  
Maps with an empty domain are cofibrations.  
The composite of two cofibrations is a cofibration.\\

\label{4.65}
\begingroup%%----------------------------------->>
\fontsize{9pt}{11pt}\selectfont
\textbf{\small EXAMPLE} \quadx 
Let $p:X \ra B$ be a surjective continuous function.  
Consider $C_p = IX \amalg B/\sim$, where $(x^\prime,0) \sim (x\pp,0)$ 
$\&$ $(x,1) \sim p(x)$ (no topology).  Let $t:C_p \ra [0,1]$ be the function $[x,t] \ra t$; let $x:t^{-1}(]0,1[) \ra X$ be the 
function $[x,t] \ra x$; let $p:t^{-1}(]0,1]) \ra B$ be the function $[x,t] \ra p(x)$.  
Definition: The 
\un{coordinate topology}
\index{coordinate topology} 
on $C_p$ is the initial topology determined by $t, x, p$.  
There is a closed embedding $j:B \ra C_p$ which is a cofibration.  For suppose that 
$
\begin{cases}
\ F:C_p \ra Y\\
\ h:IB \ra Y
\end{cases}
$
are continuous functions such that $F \circ j = h \circ i_0$ $-$then the formulas $H(j(b),T) = h(b,T)$, 
\[
H([x,t],T) = 
\begin{cases}
\ F\left[x,t + \frac{T}{2}\right] \hspace{2.4cm} (t \geq 1/2, T \leq 2 - 2t) \\
\ h(p(x),2t + T - 2) \hspace{1.5cm} (t \geq 1/2, T \geq 2 - 2t)\\
\ F[x,t + tT] \hspace{2.52cm}  (t \leq 1/2)
\end{cases}
\]
specify a continuous function $H:IC_p \ra Y$ such that $F = H \circ i_0$ and $H \circ Ij = h$.
\\ \indent
[Note: \ $C_p$ also carries another (finer) topology 
(cf. p. \pageref{3.1}).  
When $X = B$ $\&$ $p = \id_X$, 
$C_p$ is $\Gamma_cX$, and when 
$B = *$ $\&$ $p(X) = *$, $C_p$ is $\Sigma_c X$ i.e., the coordinate topology is the coarse topology 
(cf. p. \pageref{3.2} ff.).]\\
\endgroup %%------------------------------------<<

\label{4.11}
\textbf{\small LEMMA}  \  
Suppose that $i:A \ra X$ is a cofibration $-$then $i$ is an embedding.

[Form the pushout square \ 
\begin{tikzcd}
{A} \arrow{r}{i} \arrow{d}[swap]{i_0} &{X}  \arrow{d}{F}\\
{IA} \arrow{r}[swap]{h} &{Y}
\end{tikzcd}
\ 
corresponding to the 2-source $IA \overset{i_0}{\lla} A \overset{i}{\lra} X$.  
The definitions imply that there is a continuous function $G:Y \ra IX$ such that 
$
\begin{cases}
\ G \circ F = i_0\\[-.2cm]
\ G \circ h = Ii
\end{cases}
$
and a continuous function $H:IX \ra Y$ such that
$
\begin{cases}
\ H \circ i_0 = F\\[-.2cm]
\ H \circ Ii = h
\end{cases}
\hspace{-.3cm}.
$
Because $H \circ G = \id_Y$, \mG is an embedding.  On the other hand, $h \circ i_1:A \ra Y$ is an embedding, hence 
$G \circ h \circ i_1: A \ra i(A) \times \{1\}$ is a homeomorphism.]\\

\begingroup%%----------------------------------->>
\fontsize{9pt}{11pt}\selectfont
 For a subspace \mA of \mX, the cofibration condition is local in the sense that if there exists a numerable covering 
$\sU = \{U\}$ of \mX such that $\forall \ U \in \sU$, the inclusion $A \cap U \ra U$ is a cofibration, then the inclusion 
$A \ra X$ is a cofibration (cf. p. \pageref{3.3}).\\
\endgroup %%------------------------------------<<

 When \mA is a subspace of \mX and the inclusion $A \ra X$ is a cofibration, the commutative diagram 
\begin{tikzcd}
{i_0A} \arrow{r} \arrow{d} &{IA}  \arrow{d}\\
{i_0 X} \arrow{r} &{i_0 X \cup IA}
\end{tikzcd}
is a pushout square and there is a retraction $r:IX \ra 
%%----------------------------------------------------------------------------------------------04
i_0X \cup IA$.  If $\rho:i_0X \cup IA \ra IX$ is the inclusion and if 
$
\begin{cases}
\ u:X \ra IX\\[-.2cm]
\ v:X \ra IX
\end{cases}
$
are defined by 
$
\begin{cases}
\ u = i_1\\[-.2cm]
\ v = \rho \circ r \circ i_1
\end{cases}
\hspace{-.2cm}, 
$
then \mA is the equalizer of $(u,v)$.  Therefore the inclusion $A \ra X$ is a closed cofibration provided that \mX is Hausdorff or in $\dcg$.\\

\label{4.70}

\begin{proposition} \ %01
Let \mA be a subspace of \mX $-$then the inclusion $A \ra X$ is a cofibration iff $i_0 X \cup IA$ is a retract of $IX$.\\
\end{proposition}

Why should the inclusion $A \ra X$ be a cofibration if $i_0X \cup IA$ is a retract of $IX$?  
Here is the problem.  
Suppose that $\phi:i_0X \cup IA \ra Y$ is a function such that $\restr{\phi}{i_0X}$ $\&$ $\restr{\phi}{IA}$ are continuous.  
Is $\phi$ continuous?  
That the answer is ``yes'' is a consequence of a generality (which is obvious if \mA is closed).\\

\textbf{\small LEMMA}  \  If $i_0 X \cup IA$ is a retract of $IX$, then a subset \mO of $i_0 X \cup IA$ is open in 
$i_0 X \cup IA$ iff its intersection with 
$
\begin{cases}
\ i_0 X\\[-.2cm]
\ IA
\end{cases}
$
is open in 
$
\begin{cases}
\ i_0 X\\[-.2cm]
\ IA
\end{cases}
\hspace{-.2cm}.
$

[Let $r$ be the retraction in question and assume that \mO has the stated property.  
Put $X_O = \{x:(x,0) \in O\}$.  
Write $U_n$ for the union of all open $U \subset X: A \cap U \times [0,1/n[ \ \subset O$.  Note that 
$A \cap X_O = A \cap \ds\bigcup\limits_1^\infty U_n$ and 
$X - \ds\bigcup\limits_1^\infty U_n \subset \ov{A}$.  
Claim: 
$X_O \subset \ds\bigcup\limits_1^\infty U_n$.  Turn it around and take an 
$x \in X - \bigcup\limits_1^\infty U_n$ $-$then 
for any $t \in ]0,1]$, $r(\ov{A} \times \{t\}) = A \times \{t\}$, so 
$r(x,t) \in (A - \ds\bigcup\limits_1^\infty U_n) \times [0,1] =$ 
$(A - X_O) \times [0,1] \subset$ 
$(X - X_O) \times [0,1]$ $\implies$ 
$(x,0) = r(x,0) \in (X - X_O) \times [0,1]$ $\implies$ $x \in X - X_O$, from which the claim.  
Thus $O = O^\prime \cup O\pp$, where 
$O^\prime = O \cap (A \times ]0,1])$ and 
$O\pp = (i_0X \cup IA) \cap \ds\bigcup\limits_1^\infty (X_O \cap U_n \times [0,1/n[)$ are open in $i_0X \cup IA$.]\\

\begingroup%%----------------------------------->>
\fontsize{9pt}{11pt}\selectfont
\textbf{\small EXAMPLE} \quadx 
Not every closed embedding is a cofibration: Take $X = \{0\} \cup \{1/n: n \geq 1\}$ and let $A = \{0\}$.  
Not every cofibration is a closed embedding: Take $X = [0,1]/[0,1[ \ = \{[0],[1]\}$ and let $A = \{[0]\}$.\\
\endgroup %%------------------------------------<<

\begingroup%%----------------------------------->>
\fontsize{9pt}{11pt}\selectfont
\textbf{\small EXAMPLE} \quadx 
Given nonempty topological spaces 
$
\begin{cases}
\ X\\[-.2cm]
\ Y
\end{cases}
, \ 
$
form their coarse join $X *_c Y$ $-$then the closed embeddings 
$
\begin{cases}
\ X\\[-.2cm]
\ Y
\end{cases}
\ra X *_c Y
$
are cofibrations.

[It suffices to exhibit a retraction $r:I(X *_c Y) \ra i_0(X *_c Y) \cup IY$.  To this end, consider 
$r([x,y,1],T) = ([x,y,1],T)$, 
\[
r([x,y,t],T)
\begin{cases}
\ \ds\bigg(\left[x,y,\frac{2t}{2 - T}\right],0\bigg) \hspace{1.25cm} \ \bigg(0 \leq t \leq \frac{2 - T}{2}\bigg) \\[11pt]
\ \ds\bigg([x,y,1],\frac{T + 2t - 2}{t}\bigg) \hspace{.93cm} \bigg(\frac{2 - T}{2} \leq t \leq 1\bigg)
\end{cases}
.]
\]
\\
\endgroup %%------------------------------------<<

%%----------------------------------------------------------------------------------------------05
\label{3.4}
\label{3.16}
\label{5.0a}
\label{5.0f}
\begingroup%%----------------------------------->>
\fontsize{9pt}{11pt}\selectfont
\textbf{\small FACT} \   
Let $X^0 \subset X^1 \subset \cdots$ be an expanding sequence of topological spaces.  
Assume: $\forall \ n$, the inclusion 
$X^n \ra X^{n+1}$ is a cofibration $-$then $\forall \ n$, the inclusion $X^n \ra X^\infty$ is a cofibration.

[Fix retractions $r_k:IX^{k+1} \ra i_0X^{k+1} \cup IX^k$.  Noting that $IX^\infty = \colim IX^n$, work with the 
$r_k$ to exhibit $i_0X^\infty \cup IX^n$ as a retract of $IX^\infty$.]\\
\endgroup %%------------------------------------<<

\textbf{\small LEMMA} \    
Let \mX and \mY be topological spaces; let $A \subset X$ and $B \subset Y$ be subspaces.  Suppose that the inclusions 
$
\begin{cases}
\ A \ra X\\[-.2cm]
\ B \ra Y
\end{cases}
$
are cofibrations $-$then the inclusion $A \times B \ra X \times Y$ is a cofibration.

[Consider the inclusions figuring in the factorization $A \times B \ra X \times B \ra X \times Y$.]\\

Given $t: 0 \leq t \leq 1$, the inclusion $\{t\} \ra [0,1]$ is a closed cofibration and therefore, for any topological space \mX, the embedding $i_t:X \ra IX$ is a closed cofibration.  Analogously, the inclusion $\{0,1\} \ra [0,1]$ is a closed cofibration and it too can be multiplied.\\

\begin{proposition} \ %02
Let 
\begin{tikzcd}%[ sep=small]
{Z} \ar{d}[swap]{f} \ar{r}{g} &{Y}\ar{d}{\eta}\\
{X}\ar{r}[swap]{\xi}  &{P}
\end{tikzcd}
be a pushout square and assume that $f$ is a cofibration $-$then $\eta$ is a cofibration.
\end{proposition}

[The cylinder function preserves pushouts.]\\

Application: Let $A \ra X$ be a closed cofibration and let $f:A \ra Y$ be a continuous function $-$then the embedding 
$Y \ra X \sqcup_f Y$ is a closed cofibration.\\

\label{12.24} %dmc mnft
The inclusion $\bS^{n-1} \ra \bD^n$ is a closed cofibration.  Proof: Define a retraction 
$r: I\bD^n \ra i_0\bD^n \cup I\bS^{n-1}$ by letting $r(x,t)$ be the point where the line joining 
$(0,2) \in \R^n \times \R$ and 
$(x,t)$ meets $i_0\bD^n \sqcup I\bS^{n-1}$.  Consequently, if $f:\bS^{n-1} \ra A$ is a continuous function, then the embedding $A \ra \bD^n \sqcup_f A$ is a closed cofibration.  
Examples: 
(1) The embedding $\bD^n \ra \bS^n$ of $\bD^n$ as the northern or southern hemisphere of $\bS^n$ is a closed cofibration; 
(2) The embedding $\bS^{n-1} \ra \bS^n$ of $\bS^{n-1}$ as the equator of $\bS^n$ is a closed cofibration, so 
$\forall \ m \leq n$, the embedding $\bS^m \ra \bS^n$ is a closed cofibration.\\

\begingroup%%----------------------------------->>
\fontsize{9pt}{11pt}\selectfont
\textbf{\small FACT} \   
Let $f:\bS^{n-1} \ra A$ be a continuous function.  Suppose that \mA is path connected $-$then $\bD^n \sqcup_f A$ is path connected and the homomorphism $\pi_q(A) \ra \pi_q(\bD^n \sqcup_f A)$ is an isomorphism if $q < n - 1$ and an 
epimorphism if $q = n - 1$.\\
\endgroup %%------------------------------------<<

\index{Theorem: Van Kampen Theorem}
\index{Van Kampen Theorem}
\begingroup%%----------------------------------->>
\fontsize{9pt}{11pt}\selectfont
\textbf{\small VAN KAMPEN THEOREM} \quadx 
Suppose that the inclusion $A \ra X$ is a closed cofibration.  Let
%%----------------------------------------------------------------------------------------------06
$f:A \ra Y$ be a continuous function $-$then the commutative diagram
\begin{tikzcd}[ sep=large]
{\Pi A} \ar{d} \ar{r}{\Pi f} &{\Pi Y} \ar{d}\\
{\Pi X}\ar{r} &{\Pi (X \sqcup_f Y)}
\end{tikzcd}
is a pushout square in \bGRD.

[Note: If in addition \mA, \mX, and \mY are path connected, then for every $x \in A$, the commutative diagram 
\begin{tikzcd}[ sep=large]
{\pi_1(A,x)} \ar{d} \ar{r}{f_*} &{\pi_1(Y,f(x))} \ar{d}\\
{\pi_1(X,x)}  \ar{r}  &{\pi_1(X \sqcup_f Y,f(x))}
\end{tikzcd}
is a pushout square in \bGR.]\\[.25cm]
\endgroup %%------------------------------------<<

Let \mA be a subspace of \mX, $i:A \ra X$ the inclusion.

\indent\indent (DR) \ \mA is said to be a 
\un{deformation retract}
\index{deformation retract}
of \mX if there is a continuous function $r:X \ra A$ such that $r \circ i = \id_A$ and $i \circ r \simeq \id_X$.

\label{4.13}
\indent\indent (SDR) \ \mA is said to be a 
\un{strong deformation retract}
\index{strong deformation retract}
of \mX if there is a continuous function $r:X \ra A$ such that $r \circ i = \id_A$ and $i \circ r \simeq \id_X \rel A$.

If $i_0X \cup IA$ is a retract of $IX$, then $i_0X \cup IA$ is a strong deformation retract of $IX$.  
Proof: Fix a retraction $r:IX \ra i_0X \cup IA$, say $r(x,t) = (p(x,t),q(x,t))$, and consider the homotopy 
$H:I^2X \ra IX$ defined by $H((x,t),T) = (p(x,tT),(1 - T)t + Tq(x,t))$.\\

\label{3.32}
\begin{proposition} \ %03
Let \mA be a closed subspace of \mX and let $f:A \ra Y$ be a continuous function.  Suppose that \mA is a strong deformation retract of \mX $-$then the image of \mY in $X \sqcup_f Y$ is a strong deformation retract of $X \sqcup_f Y$.\\
\end{proposition}

\begingroup%%----------------------------------->>
\fontsize{9pt}{11pt}\selectfont
\textbf{\small EXAMPLE} \quadx 
The house with two rooms\\
\vspace{0.25cm}
%$

\parbox{4cm}{
\begin{tikzpicture}[scale=0.25]%[scale=0.5,shift={(-5,-3)}]
\draw[]
(0,0) -- (10.8757,-1.31258);
\draw[]
(10.8757,-1.31258) -- (10.8757,10.2785);
\draw[]
(10.8757,10.2785) -- (0,11.5911);
\draw[]
(0,11.5911) -- (0,0);
\draw[]
(5.07142,2.81484) -- (15.9471,1.50226);
\draw[]
(15.9471,1.50226) -- (15.9471,13.0934);
\draw[]
(15.9471,13.0934) -- (5.07142,14.4059);
\draw[]
(5.07142,14.4059) -- (5.07142,2.81484);
\draw[]
(0,0) -- (5.07142,2.81484);
\draw[]
(10.8757,-1.31258) -- (15.9471,1.50226);
\draw[]
(0,11.5911) -- (5.07142,14.4059);
\draw[]
(10.8757,10.2785) -- (15.9471,13.0934);
\draw[]
(0,5.79555) -- (10.8757,4.48298);
\draw[]
(10.8757,4.48298) -- (15.9471,7.29781);
\draw[]
(15.9471,7.29781) -- (5.07142,8.61039);
\draw[]
(5.07142,8.61039) -- (0,5.79555);
\draw[blue]
(3.92571,6.74964) -- (5.73832,6.53088);
\draw[blue]
(5.73832,6.53088) -- (6.58356,7.00002);
\draw[blue]
(6.58356,7.00002) -- (4.77094,7.21878);
\draw[blue]
(4.77094,7.21878) -- (3.92571,6.74964);
\draw[blue]
(3.92571,12.5452) -- (5.73832,12.3264);
\draw[blue]
(5.73832,12.3264) -- (6.58356,12.7956);
\draw[blue]
(6.58356,12.7956) -- (4.77094,13.0143);
\draw[blue]
(4.77094,13.0143) -- (3.92571,12.5452);
\draw[blue]
(3.92571,6.74964) -- (3.92571,12.5452);
\draw[blue]
(5.73832,6.53088) -- (5.73832,12.3264);
\draw[blue]
(4.77094,7.21878) -- (4.77094,13.0143);
\draw[blue]
(6.58356,7.00002) -- (6.58356,12.7956);
\draw[purple]
(9.36355,6.09335) -- (11.1762,5.87459);
\draw[purple]
(11.1762,5.87459) -- (12.0214,6.34373);
\draw[purple]
(12.0214,6.34373) -- (10.2088,6.56249);
\draw[purple]
(10.2088,6.56249) -- (9.36355,6.09335);
\draw[purple]
(9.36355,0.297795) -- (11.1762,0.079032);
\draw[purple]
(11.1762,0.079032) -- (12.0214,0.548171);
\draw[purple]
(12.0214,0.548171) -- (10.2088,0.766935);
\draw[purple]
(10.2088,0.766935) -- (9.36355,0.297795);
\draw[purple]
(9.36355,6.09335) -- (9.36355,0.297795);
\draw[purple]
(11.1762,5.87459) -- (11.1762,0.079032);
\draw[purple]
(10.2088,6.56249) -- (10.2088,0.766935);
\draw[purple]
(12.0214,6.34373) -- (12.0214,0.548171);
\draw[]
(2.53571,7.20297) -- (4.34833,6.98421);
\draw[]
(4.34833,6.98421) -- (4.34833,12.7798);
\draw[]
(4.34833,12.7798) -- (2.53571,12.9985);
\draw[]
(2.53571,12.9985) -- (2.53571,7.20297);
\draw[]
(13.4114,0.0948384) -- (11.5988,0.313602);
\draw[]
(11.5988,0.313602) -- (11.5988,6.10916);
\draw[]
(11.5988,6.10916) -- (13.4114,5.89039);
\draw[]
(13.4114,5.89039) -- (13.4114,0.0948384);
\end{tikzpicture}
}
%\end{cases}
%\text{is a strong deformation retract of $[0,1]^3$.}
%$
%\vspace{0.15cm}
\hspace{2cm} is a strong deformation retract of $[0,1]^3$.\\
\endgroup %%------------------------------------<<
\vspace{0.5cm}

\label{3.21}
\label{4.15}
\label{12.26}
\textbf{\small LEMMA}  \  
Suppose that the inclusion $A \ra X$ is a cofibration $-$then the inclusion $i_0 X \cup IA \cup i_1 X \ra IX$ is a cofibration.

[Fix a homeomorphism $\Phi:I[0,1] \ra I[0,1]$ that sends 
$I\{0\} \cup i_0[0,1] \cup I\{1\}$ to $i_0[0,1]$ $-$then the homeomorphism $\id_X \times \Phi:I^2X \ra I^2X$ sends 
$i_0IX \cup I(i_0X \cup IA \cup i_1X)$ to $i_0IX \cup I^2A$.  
Since the inclusion $IA \ra IX$ is a cofibration, 
$i_0IX \cup I^2A$ is a retract of $I^2X$ and Proposition 1 is applicable.]

[Note: \ A similar but simpler argument proves that the inclusion $i_0X \cup IA \ra IX$ is a cofibration.]\\

%%----------------------------------------------------------------------------------------------07
\begin{proposition} \ %04
If \mA is a deformation retract of \mX and if $i:A \ra X$ is a cofibration, then \mA is a strong deformation retract of \mX.
\end{proposition}

[Choose a homotopy $H:IX \ra X$ such that $H \circ i_0 = \id_X$ and $H \circ i_1 = i \circ r$, where $r:X \ra A$ 
is a retraction.  Define a function $h:I(i_0X \cup IA \cup i_1X) \ra X$ by 
\[
\begin{cases}
\ h((x,0),T) = x \hspace{3.0cm}   (x \in X) \\
\ h((a,t),T) = H(a,(1-T)t) \hspace{.9cm}  (a \in A)\\
\ h((x,1),T) = H(r(x), 1- T) \hspace{.75cm}(x \in X)
\end{cases}
.
\]
Observing that $i_0X \cup IA \cup i_1X$ can be written as the union of $i_0X \cup A \times [0,1/2]$ and 
$A \times [1/2,1] \cup i_1 X$, the lemma used in the proof of Proposition 1 implies that $h$ is continuous.  
But the restriction of \mH to $i_0X \cup IA \cup i_1X$ is $h \circ i_0$, so there exists a continuous function $G:IX \ra X$ 
which extends $h \circ i_1$.  Obviously, $G \circ i_0 = \id_X$, $G \circ i_1 = i \circ r$, and $\forall \ a \in A$, 
$\forall \ t \in [0,1]$: $G(a,t) = a$.  Therefore \mA is a strong deformation retract of \mX.]\\

\begin{proposition} \ %05
If $i:A \ra X$ is both a homotopy equivalence and a cofibration, then \mA is a strong deformation retract of \mX.
\end{proposition}

[To say that $i:A \ra X$ is a homotopy equivalence means that there exists a continuous function $r:X \ra A$ such that 
$r \circ i \simeq \id_A$ an $i \circ r \simeq \id_X$.  
However, due to the cofibration assumption, the homotopy class of $r$ contains an honest retraction, thus \mA is a deformation retract of \mX or still, a strong deformation retract of \mX 
(cf. Proposition 4).]\\

\begingroup%%----------------------------------->>
\fontsize{9pt}{11pt}\selectfont
\index{The Comb (example)}
\textbf{\small EXAMPLE  \ (\un{The Comb})} \  
Consider the subspace \mX of $\R^2$ consisting of the union 
$([0,1] \times \{0\}) \cup$ $(\{0\} \times [0,1])$ and the line segments joining $(1/n,0)$ and $(1/n,1)$ $(n = 1, 2, \ldots)$ 
$-$then \mX is contractible.  
\\[.25cm]
\begin{tikzpicture} [scale=2.0,shift={(0,0)}]
\draw[black][-] (0,0) -- (2,0) node[below] {$$};
\draw[black][-] (0,0) -- (0,1) node[below] {$$};
\draw[green][-] (1,0) -- (1,1) node[below] {$$};
\draw[green][-] (.5,0) -- (.5,1) node[below] {$$};
\draw[green][-] (.25,0) -- (.25,1) node[below] {$$};
\draw[green][-] (.125,0) -- (.125,1) node[below] {$$};
\draw[green][-] (.0625,0) -- (.0625,1) node[below] {$$};
\draw[green][-] (.03125,0) -- (.03125,1) node[below] {$$};
\draw[green][-] (.015625,0) -- (.015625,1) node[below] {$$};
\node at (4.75,0.8) {Moreover, $\{0\} \times [0,1]$ is a deformation retract of \mX.  But it is not a};
\node at (4.8,0.5) {strong deformation retract.  Therefore the inclusion $\{0\} \times [0,1] \ra X$,};
\node at (4.15,0.2) {while a homotopy equivalence,  is not a cofibration. };
\end{tikzpicture}
%Moreover, $\{0\} \times [0,1]$ is a deformation retract of \mX.  
%But it is not a strong deformation retract.  
%Therefore the inclusion  
%$\{0\} \times [0,1] \ra X$, while a homotopy equivalence, is not a cofibration.
\\
\endgroup %%------------------------------------<<

Let \mA be a subspace of \mX $-$then a 
\un{Str{\o}m structure}
\index{Str{\o}m structure} 
on $(X,A)$ consists of a continuous function $\phi:X \ra [0,1]$ such that $A \subset \phi^{-1}(0)$ and a homotopy $\Phi:IX \ra X$ of $\id_X \rel A$ 
such that $\Phi(x,t) \in A$ whenever $t > \phi(x)$.

[Note: \ If the pair $(X,A)$ admits a Str{\o}m structure $(\phi,\Phi)$ and if \mA is closed in \mX, then 
$A = \phi^{-1}(0)$.  Proof: $\phi(x) = 0 \implies x = \Phi(x,0) = \lim\Phi(x,1/n) \in A$.]

If the pair $(X,A)$ admits a Str{\o}m structure $(\phi_0,\Phi_0)$ for which $\phi_0 < 1$ throughout \mX, then \mA is a strong deformation retract of \mX.  Conversely, if \mA is a strong deformation retract of \mX and if the pair $(X,A)$ admits a Str{\o}m structure $(\phi,\Phi)$, then the pair $(X,A)$ admits a Str{\o}m structure $(\phi_0,\Phi_0)$ for which $\phi_0 < 1$ 
throughout \mX.  
Proof: Choose a homotopy $H:IX \ra X$ of $\id_X \rel A$ such that $H \circ i_1(X) \subset A$ and put 
$\phi_0(x) = \min\{\phi(x),1/2\}$, $\Phi_0(x,t) = H(\Phi(x,t),\min\{2t,1\})$.\\

%%----------------------------------------------------------------------------------------------08
\label{14.16}
\index{Theorem: Cofibration Characerization Theorem}
\index{Cofibration Characerization Theorem}
\textbf{\small COFIBRATION CHARACTERIZATION THEOREM}  \ \  
The inclusion $A \ra X$ is a cofibration iff the pair $(X,A)$ admits a Str{\o}m structure $(\phi,\Phi)$.

[Necessity: Fix a retraction $r:IX \ra i_0X \cup IA$ and let 
$X \overset{p}{\lla} IX \overset{q}{\lra} [0,1]$ be the projections.  
Consider 
$\phi(x) = \sup\limits_{0 \leq t \leq 1} \abs{t - qr(x,t)}$, $\Phi(x,t) = p r(x,t)$.

Sufficiency: Given a Str{\o}m structure $(\phi,\Phi)$ on $(X,A)$, define a retraction $r:IX \ra i_0X \cup IA$ by
\[
r(x,t) = 
\begin{cases}
\ (\Phi(x,t),0) \hspace{1.85cm} (t \leq \phi(x))\\
\ (\Phi(x,t),t - \phi(x)) \hspace{.75cm} (t \geq \phi(x))
\end{cases}
.]
\]
\\

\label{9.107}
\label{12.10}
One application of this criterion is the fact that if the inclusion $A \ra X$ is a cofibration, then the inclusion 
$\ov{A} \ra X$ is a closed cofibration.  
For let $(\phi,\Phi)$ be a Str{\o}m structure on $(X,A)$ $-$then $(\phi,\ov{\Phi})$, where 
$\ov{\Phi}(x,t) = \Phi(x,\min\{t,\phi(x)\})$, is a Str{\o}m structure on $(X,\ov{A})$.  
Another application is that if the inclusion $A \ra X$ is a closed cofibration, then the inclusion $kA \ra kX$ is a closed cofibration.  
Indeed, a Str{\o}m structure on $(X,A)$ is also a Str{\o}m structure on $(kX,kA)$.\\

\begingroup%%----------------------------------->>
\fontsize{9pt}{11pt}\selectfont
\textbf{\small EXAMPLE} \quadx 
Let $A \subset [0,1]^n$ be a compact neighborhood retract of $\R^n$ $-$then the inclusion $A \ra [0,1]^n$ is a cofibration.\\
\endgroup %%------------------------------------<<

\label{3.15}
\label{3.27}
\begingroup%%----------------------------------->>
\fontsize{9pt}{11pt}\selectfont
\textbf{\small EXAMPLE} \quadx 
Take $X = [0,1]^\kappa$ ($\kappa > \omega$) and let $A = \{0_\kappa\}$, $0_\kappa$ the ``origin'' in \mX $-$then \mA is a strong 
deformation retract of \mX but the inclusion $A \ra X$ is not a cofibration (\mA is not a zero set in \mX).\\
\endgroup %%------------------------------------<<

\label{4.43}
\label{5.19}
\begingroup%%----------------------------------->>
\fontsize{9pt}{11pt}\selectfont
\textbf{\small FACT} \   
Let \mA be a nonempty closed subspace of \mX.  Suppose that the inclusion $A \ra X$ is a cofibration $-$then $\forall \ q$, 
the projection $(X,A) \ra (X/A,*_A)$ induces an isomorphism $H_q(X,A) \ra H_q(X/A,*_A)$, $*_A$ the image of \mA in 
$X/A$.
\\ \indent
[Note: \ With \mU running over the neighborhoods of \mA in \mX, show that 
$H_q(X,A) \approx \lim H_q(X,U)$ and then 
use excision.]\\
\endgroup %%------------------------------------<<

\label{3.13}
\label{14.6}
\label{14.8}
\label{14.70}
\begingroup%%----------------------------------->>
\fontsize{9pt}{11pt}\selectfont
\textbf{\small LEMMA}  \  
Let \mX and \mY be Hausdorff topological spaces.  Let \mA be a closed subspace of \mX and let $f:A \ra Y$ be a continuous function.  Assume: The inclusion $A \ra X$ is a cofibration $-$then $X \sqcup_f Y$ is Hausdorff.\\
\fontsize{9pt}{11pt}\selectfont
\endgroup %%------------------------------------<<

As we shall now see, the deeper results in cofibration theory are best approached by implementation of the cofibration characterization theorem.\\

\begin{proposition} \ %06
Let \mK be a compact Hausdorff space.  Suppose that the inclusion $A \ra X$ is a cofibration $-$then the inclusion 
$C(K,A) \ra C(K,X)$ is a cofibration (compact open topology).
\end{proposition}

%%----------------------------------------------------------------------------------------------09
[Let $(\phi,\Phi)$ be a Str{\o}m structure on $(X,A)$.  Define $\phi_K:C(K,X) \ra [0,1]$ by 
$\phi_K(f) = \sup\limits_K \phi \circ f$ and $\Phi_K:IC(K,X) \ra C(K,X)$ by $\Phi_K(f,t)(k) = \Phi(f(k),t)$ $-$then 
$(\phi_K,\Phi_K)$ is a Str{\o}m structure on $(C(K,X),C(K,A))$
.]\\

\begingroup%%----------------------------------->>
\fontsize{9pt}{11pt}\selectfont
\textbf{\small EXAMPLE} \quadx 
If \mA is a subspace of \mX, then the inclusion $PA \ra PX$ is a cofibration provided that the inclusion $A \ra X$ is a cofibration.\\
\endgroup %%------------------------------------<<

\begingroup%%----------------------------------->>
\fontsize{9pt}{11pt}\selectfont
\textbf{\small EXAMPLE} \quadx 
Take $A = \{0,1\}$, $X = [0,1]$ $-$then the inclusion $A \ra X$ is a cofibration but the inclusion $C(\N,A) \ra C(\N,X)$ is not a cofibration (compact open topology).

[The Hilbert cube is an AR but the Cantor set is not an ANR.]\\
\endgroup %%------------------------------------<<

\begin{proposition} \ %07
Let 
$
\begin{cases}
\ A \subset X\\[-.1cm]
\ B \subset Y
\end{cases}
\hspace{-.2cm}, \ 
$
with \mA closed, and assume that the corresponding inclusions are cofibrations $-$then the inclusion 
$A \times Y \cup X \times B \ra X \times Y$ is a cofibration.
\end{proposition}

[Let $(\phi,\Phi)$ and $(\psi,\Psi)$ be Str{\o}m structures on $(X,A)$ and $(Y,B)$.  Define $\omega:X \times Y \ra [0,1]$ 
by $\omega(x,y) = \min\{\phi(x),\psi(y)\}$ and define $\Omega:I(X \times Y) \ra X \times Y$ by 
\[
\Omega((x,y),t) = (\Phi(x,\min\{t,\psi(y)\}), \Psi(y,\min\{t,\phi(x)\})).
\]
Since \mA is closed in \mX, $\phi(x) < 1$ $\implies$ $\Phi(x,\phi(x)) \in A$, so $(\omega,\Omega)$ is a Str{\o}m structure on 
$(X \times Y,A \times Y \cup X \times B)$
.]

[Note: \ If in addition, \mA (\mB) is a strong deformation retract of \mX (\mY), then $A \times Y \cup X \times B$ is a strong deformation retract of $X \times Y$.  Reason: $\phi < 1$ $(\psi < 1)$ throughout \mX (\mY) $\implies$ 
$\omega < 1$ throughout $X \times Y$.]\\

\begingroup%%----------------------------------->>
\fontsize{9pt}{11pt}\selectfont
\textbf{\small EXAMPLE} \quadx 
If the inclusion $A \ra X$ is a cofibration, then the inclusion $A \times X \cup X \times A \ra X \times X$ need not be a cofibration.  To see this, let $X = [0,1]/[0,1[ \ = \{[0],[1]\}$, $A = \{[0]\}$ and, to get a contradiction, assume that the pair 
$(X \times X,A \times X \cup X \times A)$ admits a Str{\o}m structure $(\phi,\Phi)$.  
Obviously, 
$\phi^{-1}([0,1[) \supset$ $\ov{A \times X} \cup \ov{X \times A} =$ $X \times X$ (since $\ov{A} = X$), so there exists a 
retraction 
$r:X \times X \ra A \times X \cup X \times A$.  But 
$([1],[1]) \in \ov{\{([0],[1])\}}$ $\implies$ 
$r([1],[1]) \in \ov{\{r([0],[1])\}} =$ 
$\ov{\{([0],[1])\}} =$ 
$\ov{\{[0]\}} \times \{[1]\}$ 
$\implies$ 
$r([1],[1]) = $ 
$([0],[1])$ and 
$([1],[1]) \in \ov{\{([1],[0])\}}$ $\implies$ $\cdots$ $\implies$ $r([1],[1]) = ([1],[0])$.\\
\endgroup %%------------------------------------<<

\textbf{\small LEMMA}  \  
Let \mA be a subspace of \mX and assume that the inclusion $A \ra X$ is a cofibration.  Suppose that $K,L:IX \ra Y$ are continuous functions that agree on $i_0X \cup IA$ $-$then $K \simeq L \  \rel i_0X \cup IA$.

[The inclusion $i_0X \cup IA \cup i_1X \ra IX$ is a cofibration (cf. the lemma preceding the proof of Proposition 4).  
With this in mind, define a continuous function $F:IX \ra Y$ by $F(x,t) = K(x,0)$ and a continuous function 
$h:I(i_0X \cup IA \cup i_1X) \ra Y$ by 
$
\begin{cases}
\ h((x,0),T) = K(x,T)\\
\ h((x,1),T) = L(x,T)
\end{cases}
$
%%----------------------------------------------------------------------------------------------10
$\&$ $h((a,t),T) = K(a,T) = L(a,T)$.  
Since the restriction of \mF to 
$i_0X \ \cup \  IA \ \cup \  i_1X$ is equal to $h \circ i_0$, 
there exists a continuous function $H:I^2X \ra Y$ such that $F = H \circ i_0$ and 
$\restr{H}{I(i_0X \cup IA \cup i_1X)} = h$.  
Let $\iota:[0,1] \times [0,1] \ra [0,1] \times [0,1]$ be the involution $(t,T) \ra (T,t)$ $-$then 
$H \circ (\id_X \times \iota):I^2X \ra Y$ is a homotopy between \mK and $L \  \rel i_0X \cup IA$.]\\

\begin{proposition} \ %08
Let \mA and \mB be closed subspaces of \mX.  Suppose that the inclusions 
$
\begin{cases}
\ A \ra X\\[-.1cm]
\ B \ra X
\end{cases}
\hspace{-.2cm}, \ 
$
$A \ \cap B \  \ra X$ are cofibrations $-$then the inclusion $A \cup B \ra X$ is a cofibration.
\end{proposition}

[In \ $IX$, \  write $(x,t) \sim (x,0)$ $(x \in A \ \cap \  B)$, call $\widetilde{X}$ the quotient $IX/\sim$, \ and let 
$p:IX \ra \widetilde{X}$ be the projection.  
Choose continuous functions $\phi,\psi:X \ra [0,1]$ such that $A = \phi^{-1}(0)$, 
$B = \psi^{-1}(0)$.  Define $\lambda:X \ra \widetilde{X}$ by $\ds \lambda(x) = \left[x,\frac{\phi(x)}{\phi(x) + \psi(x)}\right]$ 
if $x \notin A \cap B$, $\lambda(x) = [x,0]$ if $x \in A \cap B$ $-$then $\lambda$ is continuous and 
$
\begin{cases}
\ \lambda(x) = [x,0] \text{ on \mA}\\[-.1cm]
\ \lambda(x) = [x,1] \text{ on \mB}
\end{cases}
.
$
Consider now a pair $(F,h)$ of continuous functions 
$
\begin{cases}
\ F:X \ra Y\\[-.1cm]
\ h:I(A \cup B) \ra Y
\end{cases}
$
for which $\restr{F}{A \cup B}= h \circ i_0$.  Fix homotopies 
$
\begin{cases}
\ H_A:IX \ra Y\\[-.1cm]
\ H_B:IX \ra Y
\end{cases}
$
such that 
$
\begin{cases}
\ \restr{H_A}{IA} = \restr{h}{IA}\\[-.1cm]
\ \restr{H_B}{IB} = \restr{h}{IB}
\end{cases}
$
$\&$ $F = H_A \circ i_0 = H_B \circ i_0$ and, using the lemma, fix a homotopy $H:I^2X \ra Y$ between $H_A$ and 
$H_B \ \rel i_0X \cup I(A \cap B)$.  With $\iota$ as in the proof above, the composite 
$H \circ (\id_X \times \iota)$ factors through $I^2X \overset{p \times \id}{\lra} I\widetilde{X}$, thus there is a continuous function 
$\widetilde{H}:I\widetilde{X} \ra Y$ that renders the diagram 
\begin{tikzcd}%[ sep=small]
{I^2X}\ar{d}[swap]{p \times \id} \ar{r}{\id_X \times \iota} &{I^2 X}\ar{d}{H}\\
{I\widetilde{X}} \ar{r}[swap]{\widetilde{H}} &{Y}
\end{tikzcd}
commutative.  An extension of $(F,h)$ is then given by the composite 
$\widetilde{H} \circ(\lambda \times \id):IX \ra I\widetilde{X} \ra Y$.]\\

\begingroup%%----------------------------------->>
\fontsize{9pt}{11pt}\selectfont
\label{3.11}
\textbf{\small FACT} \   
Let \mA and \mB be closed subspaces of a metrizable space \mX.  Suppose that the inclusions 
$A \cap B \ra A$, $A \cap B \ra B$, $B \ra X$, $A - B \ra X - B$ are cofibrations $-$then the inclusion $A \ra X$ is a cofibration.\\
\endgroup %%------------------------------------<<

Let \mA be a subspace of \mX.  Suppose given a continuous function $\psi:X \ra [0,\infty]$ such that $A \subset \psi^{-1}(0)$ and a homotopy $\Psi:I\psi^{-1}([0,1]) \ra X$ of the inclusion $\psi^{-1}([0,1]) \ra X \ \rel A$ such that $\Psi(x,t) \in A$ whenever 
$t > \psi(x)$ $-$then the inclusion $A \ra X$ is a cofibration.  Proof: Define a Str{\o}m structure $(\phi,\Phi)$ on $(X,A)$ by
$\phi(x) = \min\{2\psi(x),1\}$,
\[
\Phi(x,t) = 
\begin{cases}
\ \Psi(x,t) \hspace{2.7cm} (2 \psi(x) \leq 1)\\
\ \Psi(x,t(2 - 2\psi(x))) \hspace{.75cm} (1 \leq 2\psi(x) \leq 2)\\
x  \hspace{3.75cm} (\psi(x) \geq 1)
\end{cases}
.
\]
\\

%%----------------------------------------------------------------------------------------------11
\label{4.3}
\textbf{\small LEMMA}  \  
Let \mA be a subspace of \mX and asume that the inclusion $A \ra X$ is a cofibration.  Suppose that \mU is a subspace of \mX with the property that there exists a continuous function $\pi:X \ra [0,1]$ for which 
$\ov{A} \cap U \subset \pi^{-1}(]0,1]) \subset U$ $-$then the inclusion $A \cap U \ra U$ is a cofibration.

[Fix a Str{\o}m structure $(\phi,\Phi)$ on $(X,A)$.  
Set 
$\pi_0(x) = \inf\limits_{0 \leq t \leq 1} \pi(\Phi(x,t))$ $(x \in X)$.  
Define a continuous function $\psi:U \ra [0,\infty]$ by 
$\psi(x) = \phi(x)/\pi_0(x)$.  
This makes sense since $\phi(x) = 0$ $\implies$ $\pi_0(x) > 0$ $(x \in U)$.  
Next, $\psi(x) \leq 1$ 
$\implies$ $\pi_0(x) > 0$ $\implies$ $\pi(\Phi(x,t)) > 0$ $\implies$ $\Phi(x,t) \in U$ $(\forall \ t)$.  
One can therefore let 
$\Psi:I\psi^{-1}([0,1]) \ra U$ be the restriction of $\Phi$ and apply the foregoing remark to the pair $(U,A \cap U)$.]\\

Let \mA, \mU be subspaces of a topological space \mX $-$then \mU is said to be a 
\un{halo}
\index{halo} 
of \mA 
in \mX if there exists a continuous function $\pi:X \ra [0,1]$ (the 
\un{haloing function}
\index{haloing function}) such that 
$A \subset \pi^{-1}(1)$ and $\pi^{-1}(]0,1]) \subset U$.  For example, if \mX is normal (but not necessarily Hausdorff), then every neighborhood of a closed subspace \mA of \mX is a halo of \mA in \mX but in a nonnormal \mX, a closed subspace \mA of \mX may have neighborhoods that are not halos.

\label{3.26}
\indent\indent (HA$_1$) If \mU is a halo of \mA in \mX, then \mU is a halo of $\ov{A}$ in \mX.
\label{4.2}
\label{4.4}

\indent\indent (HA$_2$) If \mU is a halo of \mA in \mX, then there exists a closed subspace \mB of \mX: 
$A \subset B \subset X$, such that \mB is a halo of \mA in \mX and \mU is a halo of \mB in \mX.

[A haloing function for $\pi^{-1}([1/2,1])$ is $\max\{2\pi(x) - 1,0\}$.]

Observation: If the inclusion $A \ra X$ is a cofibration and if \mU is a halo of \mA in \mX, then the inclusion 
$A \ra U$ is a cofibration.

[This is a special case of the lemma.]\\

\begin{proposition} \ %09
If $j:B \ra A$ and $i:A \ra X$ are continuous functions such that $i$ and $i \circ j$ are cofibrations, then $j$ is a cofibration.
\end{proposition}

[Take $i$ and $j$ to be inclusions.  Using the cofibration characterization theorem, fix a halo \mU of \mA in \mX and a retraction $r:U \ra A$.  Since \mU is also a halo of \mB in \mX, the inclusion $B \ra U$ is a cofibration.  
Consider a commutative diagram
\begin{tikzcd}%[ sep=small]
{B} \ar{d}[swap]{j} \ar{r}{g} &{PY} \ar{d}{p_0}\\
{A} \ar{r}[swap]{F} &{Y}
\end{tikzcd}
.  To construct a filler for this, pass to its counterpart
\begin{tikzcd}%[ sep=small]
{B} \ar{d} \ar{r}{g} &{PY} \ar{d}{p_0}\\
{U} \ar{r}[swap]{F\circ r} &{Y}
\end{tikzcd}
over U, which thus admits a filler $G:U \ra PY$.  
The restriction $\restr{G}{A}:A \ra PY$ will then do the trick.]\\

\begingroup%%----------------------------------->>
\fontsize{9pt}{11pt}\selectfont
\index{Telescope Construction}
\textbf{\small EXAMPLE \  (\un{Telescope Construction})} \  
Let $X^0 \subset X^1 \subset \cdots $ be an expanding sequence of topological
%%----------------------------------------------------------------------------------------------12
spaces.  Assume: $\forall \ n$, the inclusion $X^n \ra X^{n+1}$ is a closed cofibration $-$then $\forall \ n$, the inclusion 
$X^n \ra X^\infty$ is a closed cofibration (cf. p. \pageref{3.4}).  
Write $\tel X^\infty$ for the quotient 
$\ds\coprod\limits_0^\infty X^n \times [n,n+1]/\sim$.  
Here, $\sim$ means that the pair $(x,n+1) \in X^n \times \{n+1\}$ is identified with the pair 
$(x,n+1) \in X^{n+1} \times \{n+1\}$.  
\label{5.0h}
One calls $\tel X^\infty$ the 
\un{telescope}
\index{telescope} 
of $X^\infty$.  
It can be viewed as a closed subspace of $X^\infty \times [0,\infty[$.  
The inclusion 
$\telsub_nX^\infty \equiv \ds\bigcup\limits_{k=0}^n X^k \times [k,k+1] \ra X^\infty \times [0,\infty[$ 
is a closed cofibration (cf. Proposition 8), so the same is true of the inclusion $\telsub_n X^\infty \ra \telsub_{n+1} X^\infty$ (cf. Proposition 9) and 
$\tel X^\infty = \colim \ \telsub_n X^\infty$.  
Denote by $p^\infty$ the composite 
$\tel X^\infty \ra X^\infty \times [0,\infty[ \ra X^\infty$.
\label{5.55c}
\\ \indent
Claim: $p^\infty$ is a homotopy equivalence.
\\ \indent
[It suffices to establish that $\tel X^\infty$ is a strong deformation retract of $X^\infty \times [0,\infty[$.  
One approach is to piece together strong deformation retractions 
$X^{n+1} \times [0,n+1] \ra$
$X^{n+1} \times\{n+1\} \cup X^n \times [0,n+1]$.]\\
\\
\endgroup %%------------------------------------<<
\label{4.29}

\begingroup%%----------------------------------->>
\fontsize{9pt}{11pt}\selectfont
Let 
$
\begin{cases}
\ X^0 \subset X^1 \subset \cdots\\
\ Y^0 \subset Y^1 \subset \cdots
\end{cases}
$
be expanding sequences of topological spaces.  Assume: $\forall \ n$, the inclusions 
$
\begin{cases}
\ X^n \ra X^{n+1}\\
\ Y^n \ra Y^{n+1}
\end{cases}
$
are closed cofibrations.  Suppose given a sequence of continuous functions $\phi^n:X^n \ra Y^n$ 
such that $\forall \ n$, the diagram 
\begin{tikzcd}[ sep=large]
{X^n} \ar{d}[swap]{\phi^n} \ar{r} &{X^{n+1}} \ar{d}{\phi^{n+1}}\\
{Y^n} \ar{r} &{Y^{n+1}}
\end{tikzcd}
commutes.  Associated with the $\phi^n$ is a continuous function 
$\phi^\infty:X^\infty \ra Y^\infty$ and a continuous function 
$\tel \phi:\tel X^\infty \ra \tel Y^\infty$, the latter being defined by 
\[
\tel \phi(x,n+t) = 
\begin{cases}
\ (\phi^n(x),n+2t) \in Y^n \times [n,n+1] \quadx (0 \leq t \leq 1/2)\\
\ (\phi^n(x),n+1) \in Y^{n+1} \times \{n+1\} \quadx (1/2 \leq t \leq 1)
\end{cases}
.
\]
\label{4.57}
There is then a commutative diagram
\begin{tikzcd}[ sep=large]
{\tel X^\infty} \ar{d}[swap]{\tel \phi} \ar{r} &{X^\infty} \ar{d}{\phi^\infty}\\
{\tel Y^\infty} \ar{r} &{Y^\infty}
\end{tikzcd}
.  The horizontal arrows are homotopy equivalences.  Moreover, $\tel \phi$ is a homotopy equivalence if this is the case of the 
$\phi^n$, thus, under these circumstances, $\phi^\infty:X^\infty \ra Y^\infty$ itself is a homotopy equivalence.

[Note: \ One can also make the deduction from first principles (cf. Proposition 15).]\\
\endgroup %%------------------------------------<<

\begin{proposition} \ %10
Let \mA be a closed subspace of a topological space \mX.  Suppose that \mA admits a halo \mU with $A = \pi^{-1}(1)$ for which there exists a homotopy $\Pi:IU \ra X$ of the inclusion $U \ra X$ $\rel A$ such that $\Pi \circ i_1(U) \subset A$ $-$then 
the inclusion $A \ra X$ is a closed cofibration.
\end{proposition}

[Define a retraction $r:IX \ra i_0X \cup IA$ as follows: 
(i) $r(x,t) = (x,0)$ $(\pi(x) = 0)$; 
(ii) $r(x,t) = (\Pi(x,2\pi(x)t),0)$ $(0 < \pi(x) \leq 1/2)$; 
(iii) $r(x,t) = (\Pi(x,t/2(1 - \pi(x))),0)$ $(1/2 \leq \pi(x) < 1$ $\&$ $0 \leq t \leq 2(1 - \pi(x)))$ and 
$r(x,t) = (\Pi(x,1),t - 2(1 - \pi(x)))$ $(1/2 \leq \pi(x) < 1$ $\&$ $2(1 - \pi(x)) \leq t \leq 1)$; 
(iv) $r(x,t) = (x,t)$ $(\pi(x) = 1)$.]\\

%%----------------------------------------------------------------------------------------------13
\begingroup%%----------------------------------->>
\fontsize{9pt}{11pt}\selectfont
\label{3.35}
\textbf{\small EXAMPLE} \quadx 
If \mA is a subcomplex of a CW complex \mX, then the inclusion $A \ra X$ is a closed cofibration.\\
\endgroup %%------------------------------------<<

A topological space \mX is said to be 
\un{locally contractible}
\index{locally contractible} 
provided that for any $x \in X$ and any neighborhood \mU of $x$ there exists a neighborhood $V \subset U$ of $x$ such that the inclusion $V \ra U$ is inessential.  
If \mX is locally contractible, then \mX is locally path connected.  
Example: $\forall \ X, X^*$ is locally contractible (cf. p. \pageref{3.5}).

[Note: \ The empty set is locally contractible but not contractible.]\\

\label{4.25}
\label{4.28}
\label{14.10}
\label{14.11}
\label{14.18}
\label{14.21}
\label{14.91}
\begingroup%%----------------------------------->>
\fontsize{9pt}{11pt}\selectfont
A topological space \mX is said to be 
\un{numerably contractible}
\index{numerably contractible} 
if it has a numerable  
covering $\{U\}$ for which each inclusion $U \ra X$ is inessential.  
Example: Every locally contractible paracompact Hausdorff space is numerably contractible.
\\ \indent
[Note: \ The product of two numerably contractible spaces is numerably contractible.]\\
\endgroup %%------------------------------------<<

\label{3.17}
\label{3.18}
\begingroup%%----------------------------------->>
\fontsize{9pt}{11pt}\selectfont
\textbf{\small FACT} \   
Numerable contractibility is a homotopy type invariant.  Proof: If \mX is dominated in homotopy by \mY and if \mY is numerably contractible, then \mX is numerably contractible.\\
\endgroup %%------------------------------------<<

\label{5.0p}
\begingroup%%----------------------------------->>
\fontsize{9pt}{11pt}\selectfont
Examples: 
(1) Every topological space having the homotopy type of a CW complex is numerably contractible; 
(2) If the $X^n$ of the telescope construction are numerably contractible, then $X^\infty$ is numerably contractible (consider 
$\tel X^\infty$).\\
\endgroup %%------------------------------------<<

A topological space \mX is said to be 
\un{uniformly locally contractible}
\index{uniformly locally contractible} 
provided that there exists a neighborhood \mU of the diagonal $\Delta_X \subset X \times X$ and a homotopy 
$H:IU \ra X$ between $\restr{p_1}{U}$ and $\restr{p_2}{U}$ rel$\Delta_X$, 
where $p_1$ and $p_2$ are the projections onto the first and second factors.  
Examples:
(1) $\R^n$, $\bD^n$, and $\bS^{n-1}$ are uniformly locally contractible; 
(2) The long ray $L^+$ is not uniformly locally contractible.\\

\begingroup%%----------------------------------->>
\fontsize{9pt}{11pt}\selectfont
\index{Stratifiable Spaces (example)}
\textbf{\small EXAMPLE  \ (\un{Stratifiable Spaces})} \  
Suppose that \mX is stratifiable and in NES(stratifiable) $-$then \mX is uniformly locally contractible.  Thus put 
$A = X \times i_0X \cup (I \Delta_X) \cup X \times i_1X$, a closed subspace of the stratifiable space $I(X \times X)$.  
Define a continuous function $\phi:A \ra X$ by 
$
\begin{cases}
\ (x,y,0) \ra x\\
\ (x,y,1) \ra y
\end{cases}
$
$\&$ $(x,x,t) \ra x$ $-$then $\phi$ extends to a continuous function $\Phi:O \ra X$, where \mO is a neighborhood of \mA in 
$I(X \times X)$.  
Fix a nieghborhood \mU of $\Delta_X$ in $X \times X$ : $IU \subset O$ and consider $H = \restr{\Phi}{IU}$.
\\ \indent
[Note: \ Every CW complex is stratifiable 
(cf. p. \pageref{3.6}) 
and in NES(stratifiable) 
(cf. p. \pageref{3.7}).  
Every metrizable topological manifold is stratifiable 
(cf. p. \pageref{3.8} ff.: metrizable  $\implies$  stratifiable) 
and, being an ANR 
(cf. p. \pageref{3.9}),  
is in NES(stratifiable) 
(cf. p. \pageref{3.10}: 
stratifiable $\implies$ perfectly normal + paracompact).]\\
\endgroup %%------------------------------------<<

\begingroup%%----------------------------------->>
\fontsize{9pt}{11pt}\selectfont
\textbf{\small FACT} \   
Let \mK be a compact Hausdorff space.  
Suppose that \mX is uniformly locally contractible $-$then $C(K,X)$ is uniformly locally contractible (compact open topology).\\
\endgroup %%------------------------------------<<

%%----------------------------------------------------------------------------------------------14

\textbf{\small LEMMA}  \  
A uniformly locally contractible topological space \mX is locally contractible.
\\ \indent
[Take a point $x_0 \in X$ and let $U_0$ be a neighborhood of $x_0$ $-$then $I\{(x_0,x_0)\} \subset$ $H^{-1}(U_0)$.  
Since $H^{-1}(U_0)$ is open in $IU$, hence open in $I(X \times X)$, there exists a neighborhood $V_0 \subset U_0$ of $x_0$ : 
$I(V_0 \times V_0) \subset H^{-1}(U_0)$.  
To see that the inclusion $V_0 \ra U_0$ is inessential, define 
$H_0:IV_0 \ra U_0$ by $H_0(x,t) = H((x,x_0),t)$.]

[Note: \ The homotopy $H_0$ keeps $x_0$ fixed throughout the entire deformation.  In addition, the argument shows that an open subspace of a uniformly locally contractible space is uniformly locally contractible.]\\

\begingroup%%----------------------------------->>
\fontsize{9pt}{11pt}\selectfont
\index{A Spaces (example)}
\textbf{\small EXAMPLE  \ (\un{A Spaces})} \  
Every A space is locally contractible.  In fact, if \mX is a nonempty A space, then $\forall \ x \in X$, $U_x$ is contractible, thus \mX has a basis of contractible open sets, so \mX is locally contractible.  But an A space need not be uniformly locally contractible.  Consider, e.g., $X = \{a,b,c,d\}$, where 
$
\begin{cases}
\ c \leq a\\
\ d \leq a
\end{cases}
$
\hspnx
, \ 
$
\begin{cases}
\ c \leq b\\
\ d \leq b
\end{cases}
\hspnx
.
$
\\[.25cm]
\endgroup %%------------------------------------<<

\begingroup%%----------------------------------->>
\fontsize{9pt}{11pt}\selectfont
\textbf{\small FACT} \   
Let \mX be a perfectly normal paracompact Hausdorff space. Suppose that \mX admits a covering by open sets \mU, each of which is  uniformly locally contractible $-$then \mX is  uniformly locally contractible.

[Use the domino principle.]\\
\endgroup %%------------------------------------<<

\label{14.117}
When is \mX uniformly locally contractible?  A sufficient condition is that the inclusion 
$\Delta_X \ra X \times X$ be a cofibration.  Proof: Fix a Str{\o}m structure $(\phi,\Phi)$ on the pair 
$(X \times X,\Delta_X)$, put 
$U = \phi^{-1}([0,1[)$ and define $H:IU \ra X$ by 
\[
H((x,y),t) = 
\begin{cases}
\ p_1(\Phi((x,y),2t)) \hspace{1.55cm} (0 \leq t \leq 1/2)\\
\ p_2(\Phi((x,y),2 - 2t)) \hspace{.75cm} \ (1/2 \leq t \leq 1)
\end{cases}
.
\]
\\

\label{6.8}
\begingroup%%----------------------------------->>
\fontsize{9pt}{11pt}\selectfont
\textbf{\small FACT} \   
Suppose that \mX is a perfectly normal Hausdorff space with a perfectly normal square $-$then \mX is uniformly locally contractible iff the diagonal embedding $X \ra X \times X$ is a cofibration.
\\ \indent
[Use Proposition 10, noting that $\Delta_X$ is a zero set.]\\

\label{3.12}
\label{4.9}
Application: If \mX is a CW complex or a metrizable topological manifold, then the diagonal embedding 
$X \ra X \times X$ is a cofibration.\\

\textbf{\small FACT} \   
Let \mA be a closed subspace of a metrizable space \mX such that the inclusion $A \ra X$ is a cofibration. Suppose that \mA and 
$X - A$ are uniformly locally contractible $-$then \mX is uniformly locally contractible.
\\ \indent
[Show that the inclusion $\Delta_X \ra X \times X$ is a cofibration by applying the result on p. \pageref{3.11} to the triple $(X \times X,\Delta_X,A \times A)$.]\\
\endgroup %%------------------------------------<<

%%----------------------------------------------------------------------------------------------15
\label{18.2} %dmc mnft
\begin{proposition} \ %11
Supppose that $A \subset X$ admits a halo \mU such that the inclusion $\Delta_U \ra U \times U$ is a cofibration.  
Assume that the inclusion $A \ra X$ is a cofibration $-$then the inclusion $\Delta_A \ra A \times A$ is a cofibration.
\end{proposition}

[Consider the commutative diagram \ 
\begin{tikzcd}%[ sep=small]
{A}\ar{d} \ar{r}{\Delta_A} &{A \times A}\ar{d}\\
{U} \ar{r}[swap]{\Delta_U} &{U \times U}
\end{tikzcd}
.  \ The vertical arrows are cofibrations, as is $\Delta_U$.  
That $\Delta_A$ is a cofibration is therefore implied by Proposition 9.]\\

\begin{proposition} \ %12
Let \mX be a Hausdorff space and suppose that the inclusion $\Delta_X \ra X \times X$ is a cofibration.  Let $f:X \ra [0,1]$ be a continuous function such that $A = f^{-1}(0)$ is a retract of $f^{-1}([0,1[)$ $-$then the inclusion $A \ra X$ is a closed cofibration.
\end{proposition}

[Write $r$ for the retraction $f^{-1}([0,1[) \ra A$, fix a Str{\o}m structure $(\phi,\Phi)$ on the pair 
$(X \times X,\Delta_X)$, and let $H:IU \ra X$ be as above.  
Define $\phi_f:X \ra [0,1]$ by $\phi_f(x) = \max\{f(x),\phi(x,r(x))\}$ $(f(x) < 1)$ $\&$ 
$\phi_f(x) = 1$ $(f(x) = 1)$ $-$then $\phi_f^{-1}(0) = A$.  
Put $H_f(x,t) = H((x,r(x)),t)$ to obtain a homotopy 
$H_f:I\phi_f^{-1}([0,1[) \ra X$ of the inclusion $\phi_f^{-1}([0,1[) \ra X \ \rel A$ such that 
$H_f\circ i_1(\phi_f^{-1}([0,1[)) \subset A$.  
Finish by citing Proposition 10.]\\

\label{14.112}
\label{14.132}
Application: Let \mX be a Hausdorff space and suppose that the inclusion $\Delta_X \ra X \times X$ is a cofibration.  
Let $e \in C(X,X)$ be idemptotent: $e \circ e = e$ $-$then the inclusion $e(X) \ra X$ is a closed cofibration.

[Define $f:X \ra [0,1]$ by $f(x) = \phi(x,e(x))$.]\\

\label{6.9}
\label{14.13}
\label{16.46}
So, if \mX is a Hausdorff space and if the inclusion $\Delta_X \ra X \times X$ is a cofibration, then for any retract \mA of \mX, the inclusion $A \ra X$ is a closed cofibration.  
\label{3.14}
In particular: $\forall \ x_0 \in X$, the inclusion $\{x_0\} \ra X$ is a closed cofibration, which as seen above, is a condition realized by every CW complex or metrizable topological manifold.

[Note: \ Let \mX be the Cantor set $-$then $\forall \ x_0 \in X$, the inclusion $\{x_0\} \ra X$ is closed but is not a cofibration.]\\

\label{14.12}
\begingroup%%----------------------------------->>
\fontsize{9pt}{11pt}\selectfont
\textbf{\small FACT} \   
Let \mX be in $\dcg$ and suppose that the inclusion 
$\Delta_X \ra X \times_k X$ is a cofibration $-$then for any retract \mA of \mX, the inclusion $A \ra X$ is a closed cofibration.
\\ \indent
[Rework Proposition 12, noting that for any continuous function $f:X \ra X$, the function $X \ra X \times_k X$ defined by 
$x \ra (x,f(x))$ is continuous.]\\
\endgroup %%------------------------------------<<

\label{12.17a}
\label{14.15}
\label{14.98}
\label{14.100}
\begingroup%%----------------------------------->>
\fontsize{9pt}{11pt}\selectfont
\textbf{\small LEMMA}  \  
Suppose that the inclusions 
$
\begin{cases}
\ A \ra X\\
\ A^\prime \ra X^\prime
\end{cases}
$ 
are closed cofibrations and that \mX is a closed subspace of $X^\prime$ with $A = X \cap A^\prime$.  
Let
$
\begin{cases}
\ f:A \ra Y\\
\ f^\prime:A^\prime \ra Y^\prime
\end{cases}
$ 
be continuous functions.  Assume that the diagram
%%----------------------------------------------------------------------------------------------16
\begin{tikzcd}[ sep=large]
{X}\ar{d} &{A}\ar{d} \ar{l} \ar{r}{f} &{Y}\ar{d}\\
{X^\prime} &{A^\prime} \ar{l} \ar{r}[swap]{f^\prime} &{Y^\prime}
\end{tikzcd}
commutes and that the vertical arrows are cofibrations $-$then the induced map 
$X \sqcup_f Y \ra X^\prime \sqcup_{f^\prime} Y^\prime$ is a cofibration and $(X \sqcup_f Y) \cap Y^\prime = Y$.
\\ \indent
[Consider a commutative diagram 
\begin{tikzcd}[ sep=large]
{X \sqcup_f Y} \ar{d} \ar{r}{g} &{PZ} \ar{d}{p_0}\\
{X^\prime \sqcup_{f^\prime} Y^\prime} \ar{r}[swap]{F^\prime}&{Z}
\end{tikzcd}
.   \ 
To construct a filler $H^\prime$ for this, work first with, 
\begin{tikzcd}[ sep=large]
{Y} \ar{d}\ar{r}
&{X \sqcup_f Y} \ar{r}{g} \ar{d}
&{PZ} \ar{d}{p_0}\\
{Y^\prime} \ar{r}
&{X^\prime \sqcup_{f^\prime} Y^\prime} \ar{r}[swap]{F^\prime}
&{Z}
\end{tikzcd}
\ to get an arrow $G:Y^\prime \ra PZ$.  Next, look at 
$
\begin{cases}
\ A^\prime \overset{f^\prime}{\lra} Y^\prime \overset{G}{\lra} PZ\\
\ X \lra X \sqcup_f Y \overset{g}{\lra} PZ
\end{cases}
. \ 
$
Since equality obtains on $A = X \cap A^\prime$, $\exists$ $G^\prime \in C(X \cup A^\prime,PZ)$ : 
$\restr{G^\prime}{A^\prime} = G \circ f^\prime$.  
But the inclusion $X \cup A^\prime \ra X^\prime$ is a cofibration (cf. Proposition 8), so the commutative diagram 
\begin{tikzcd}[ sep=large]
{X \cup A^\prime} \ar{d}\ar{rr}{G^\prime}
&&{PZ} \ar{d}{p_0}\\
{X^\prime} \ar{r}
&{X^\prime \sqcup_{f^\prime} Y^\prime} \ar{r}[swap]{F^\prime}
&{Z}
\end{tikzcd}
\ admits a filler $H:X^\prime \ra PZ$ which agrees with $G \circ f^\prime$ on $A^\prime$ and therefore determines 
$H^\prime:X^\prime: \sqcup_{f^\prime:} Y^\prime: \ra PZ$.]\\
\endgroup %%------------------------------------<<

\label{14.14}
\label{14.113}
\label{14.146}
\begingroup%%----------------------------------->>
\fontsize{9pt}{11pt}\selectfont
\textbf{\small FACT} \   
Let $A \ra X$ be a closed cofibration and let $f:A \ra Y$ be a continuous function.  
Suppose that 
$
\begin{cases}
\ X\\
\ Y
\end{cases}
$
are in $\dcg$ and that the inclusions 
$
\begin{cases}
\ \Delta_X \ra X \times_k X\\
\ \Delta_Y \ra Y \times_k Y
\end{cases}
$
are cofibrations $-$then the inclusion $\Delta_Z \ra Z \times_k Z$ is a cofibration, \mZ the adjunction space $X \sqcup_f Y$.
\\ \indent
[There are closed cofibrations 
$
\begin{cases}
\ A \times_k A \ra X \times_k A \cup A \times_k X\\
\ Y \times_k Y \ra Z \times_k Y \cup Y \times_k Z
\end{cases}
. \ 
$
Precompose these arrows with the diagonal embeddings, form the commutative diagram
\[
\begin{tikzcd}[ sep=large]
{X} \ar{d}
&{A} \ar{l} \ar{r}{f} \ar{d}
&{Y} \ar{d}\\
%\\
{X \times_k X} 
&{X \times_k A \cup A \times_k X} \ar{l} \ar{r}
&{Z \times_k Y \cup Y \times_k Z}
\end{tikzcd}
,
\]
and apply the lemma.]
\\ \indent
[Note: \ Proposition 7 remains in force if the product in \bTOP is replaced by the product in $\dcg$.  
Take $U = X$ in Proposition 11 to see that the inclusion $\Delta_A \ra A \times_k A$ is a cofibration.]\\
\endgroup %%------------------------------------<<

\begingroup%%----------------------------------->>
\fontsize{9pt}{11pt}\selectfont
Application: Let \mX and \mY be CW complexes.  Let \mA be a subcomplex of \mX and let $f:A \ra Y$ 
be a continuous function $-$then the inclusion $\Delta_Z \ra Z \times_k Z$ is a cofibration, \mZ the adjunction space $X \sqcup_f Y$.
\\ \indent
%%----------------------------------------------------------------------------------------------17
[The inclusions 
$
\begin{cases}
\ \Delta_X \ra X \times X\\
\ \Delta_Y \ra Y \times Y
\end{cases}
$
are cofibrations (cf. p. \pageref{3.12}), thus the same is true of the inclusions 
$
\begin{cases}
\ \Delta_X \ra X \times_k X\\
\ \Delta_Y \ra Y \times_k Y
\end{cases}
$
(cf. p. \pageref{3.13}).  \mZ itself need not be a CW complex but, in view of the skeletal approximation theorem, \mZ at least has the homotopy type of a CW complex.]\\
\endgroup %%------------------------------------<<

\begingroup%%----------------------------------->>
\fontsize{9pt}{11pt}\selectfont
\textbf{\small FACT} \   
Let $A \ra X$ be a closed cofibration and let $f:A \ra Y$ be a continuous function.  Suppose that 
$
\begin{cases}
\ X\\
\ Y
\end{cases}
$
are uniformly locally contractible perfectly normal Hausdorff spaces with perfectly normal squares $-$then $X \sqcup_f Y$ is uniformly locally contractible provided that its square is perfectly normal.

[Note: \ A priori, $X \sqcup_f Y$ is a perfectly normal Hausdorff space (cf. AD$_5$).]\\
\endgroup %%------------------------------------<<

A pointed space $(X,x_0)$ is said to be 
\un{wellpointed}
\index{wellpointed} 
if the inclusion $\{x_0\} \ra X$ is a cofibration.  $\ov{\Pi}X$ is the full subgroupoid of $\Pi X$ whose objects are the $x_0 \in X$ such that $(X,x_0)$ is wellpointed.  
\label{5.0e}
Example: Let \mX be a CW complex or a metrizable topological manifold $-$then $\forall \ x_0 \in X$, $(X,x_0)$ is wellpointed 
(cf. p. \pageref{3.14}).

[Note: \ Take $X = [0,\Omega]$, $x_0 = \Omega$ $-$then $(X,x_0)$ is not wellpointed.]\\

\begingroup%%----------------------------------->>
\fontsize{9pt}{11pt}\selectfont
\label{3.28}
The full subcategory of $\bHTOP_*$ whose objects are the wellpointed spaces is not isomorphism closed, i.e., if 
$(X,x_0) \approx (Y,y_0)$ in $\bHTOP_*$, then it can happen that the inclusion $\{x_0\} \ra X$ is a cofibration but the inclusion $\{y_0\} \ra Y$ is not a cofibration (cf. p. \pageref{3.15}).\\
\endgroup %%------------------------------------<<

\label{4.76}
\begingroup%%----------------------------------->>
\fontsize{9pt}{11pt}\selectfont
\textbf{\small EXAMPLE} \quadx 
Let \mX be a topological manifold $-$then $\forall \ x_0 \in X$, $(X,x_0)$ is wellpointed.\\
\endgroup %%------------------------------------<<

\begingroup%%----------------------------------->>
\fontsize{9pt}{11pt}\selectfont
\textbf{\small FACT} \   
Let \mK be a compact Hausdorff space.  
Suppose that$(X,x_0)$  is wellpointed $-$then $\forall \ k_0 \in K$, 
$C(K,k_0;X,x_0)$ is wellpointed (compact open topology).
\\ \indent
[Note: \ The base point in $C(K,k_0;X,x_0)$ is the constant map $K \ra x_0$.]\\
\endgroup %%------------------------------------<<

Given topological spaces
$
\begin{cases}
\ X\\[-.15cm]
\ Y
\end{cases}
, \ 
$
the base point functor $\ov{\Pi}X \times \Pi Y \ra \bSET$ sends an object $(x_0,y_0)$ to the set $[X,x_0;Y,y_0]$.  
To describe its behavior on morphisms, let 
$
\begin{cases}
\ x_0, x_1 \in X\\[-.15cm]
\ y_0,y_1 \in Y
\end{cases}
$
and suppose that both $(X,x_0)$ and $(X,x_1)$ are wellpointed.  Let $\sigma \in PX$:
$
\begin{cases}
\ \sigma(0) = x_0\\[-.15cm]
\ \sigma(1) = x_1
\end{cases}
$
$\&$ let $\tau \in PY$ :
$
\begin{cases}
\ \tau(0) = y_0\\[-.15cm]
\ \tau(1) = y_1
\end{cases}
$
$-$then the pair $(\sigma,\tau)$ determines a bijection 
$[\sigma,\tau]_\#:[X,x_0;Y,y_0] \ra [X,x_1;Y,y_1]$ that depends only on the path classes of
$
\begin{cases}
\ \sigma\\[-.15cm]
\ \tau
\end{cases}
$
in 
$
\begin{cases}
\ \Pi X\\[-.15cm]
\ \Pi Y
\end{cases}
$
\hspace{-.2cm}.  
Here is the procedure.  
Fix a homotopy $H:IX \ra X$ such that $H \circ i_0 = \id_X$, $H(x_1,t) = \sigma(1-t)$, and put 
$e = H \circ i_1$.  
Take an $f \in C(X,x_0;Y,y_0)$ and define a continuous function 
$F:i_0X \cup I\{x_1\} \ra X \times Y$ by 
$
\begin{cases}
\ (x,0) \ra (e(x),f(e(x)))\\[-.15cm]
\ (x_1,t) \ra (\sigma(t),\tau(t))
\end{cases}
$
$-$then the
%%----------------------------------------------------------------------------------------------18
diagram
\begin{tikzcd}%[ sep=small]
{i_0X \cup I\{x_1\}} \ar{d} \ar{r}{F} &{X \times Y} \ar{d}{p}\\
{IX} \ar{r}[swap]{G} &{X}
\end{tikzcd}
commutes, where $G(x,t) = H(x,1-t)$.  To construct a filler $H_f:IX \ra X \times Y$, let $q:X \times Y \ra Y$ be the projection, 
choose a retraction $r:IX \ra i_0X \cup I\{x_1\}$ and set $H_f(x,t) = (G(x,t),qF(r(x,t)))$.  
Write 
$f_\# = q \circ H_f \circ i_1 \in C(X,x_1;Y,y_1)$.  
Definition: $[\sigma,\tau]_\#[f] = [f_\#]$.  
The fundamental group $\pi_1(Y,y_0)$ thus operates to the left on 
$[X,x_0;Y,y_0]$ : $([\tau],[f]) \ra [\sigma_0,\tau]_\#[f]$, $\sigma_0$ the constant path in \mX at $x_0$.  
If $f$, $g \in C(X,x_0;Y,y_0)$ then $f \simeq g$ in \bTOP iff 
$\exists$ $[\tau] \in \pi_1(Y,y_0) \ :$ $[\sigma_0,\tau]_\#[f] = [g]$.  
\label{5.0ae}
Therefore the forgetful function $[X,x_0;Y,y_0] \ra [X,Y]$ passes to the quotient to define an injection 
$\pi_1(Y,y_0)\backslash [X,x_0;Y,y_0] \ra [X,Y]$ which, when \mY is path connected, is a bijection.  
The  forgetful function $[X,x_0;Y,y_0] \ra [X,Y]$ is one-to-one iff the action of $\pi_1(Y,y_0)$ on $[X,x_0;Y,y_0]$ is trivial.  
Changing \mY to \mZ by a homotopy equivalence in 
$\bTOP \ : \ 
\begin{cases}
\ Y \ra Z\\[-.15cm]
\ y_0 \ra z_0
\end{cases}
$
leads to an arrow $[X,x_0;Y,y_0] \ra$ $[X,x_0;Z,z_0]$.  It is a bijection.\\
\vspace{0.5cm}

\label{4.37}
\begingroup%%----------------------------------->>
\fontsize{9pt}{11pt}\selectfont
\label{4.38}
\label{5.1}
\label{5.0j}
\label{5.0af}
\label{5.8a}
\textbf{\small FACT} \   
Suppose that \mX and \mY are path connected.  Let $f \in C(X,Y)$ and assume that $\forall \ x \in X$, 
$f_*:\pi_1(X,x) \ra \pi_1(Y,f(x))$ is surjective $-$then $\forall \ x \in X$, $f_*:\pi_n(X,x) \ra \pi_n(Y,f(x))$ is injective 
(surjective) iff $f_*:[\bS^n,X] \ra [\bS^n,Y]$ is injective (surjective).\\
\endgroup %%------------------------------------<<

\textbf{\small LEMMA}  \  
Suppose that the inclusion $i:A \ra X$ is a cofibration.  Let $f \in C(X,X)$: $f \circ i = i$ $\&$ $f \simeq \id_X$ $-$then 
$\exists$ $g \in C(X,X)$: $g \circ i = i$ $\&$ $g \circ f \simeq \id_X \ \rel A$.

[Let $H:IX \ra X$ be a homotopy with $H \circ i_0 = f$ and $H \circ i_1 = \id_X$; let $G:IX \ra X$ be a homotopy with 
$G \circ i_0 = \id_X$ and $G \circ Ii = H \circ Ii$.  Define $F:IX \ra X$ by 
$
F(x,t) =
\begin{cases}
\ G(f(x),1 - 2t) \hspace{.75cm} (0 \leq t \leq 1/2)\\[-.15cm]
\ H(x,2t - 1) \hspace{1.25cm}  (1/2 \leq t \leq 1)
\end{cases}
$
and put
\[
k((a,t),T) = 
\begin{cases}
\ G(a,1 - 2t(1 - T)) \hspace{1.45cm}  (0 \leq t \leq 1/2)\\[-.15cm]
\ H(a,1 - 2(1 - t)(1 - T)) \hspace{.5cm} (1/2 \leq t \leq 1)
\end{cases}
\]
to get a homotopy $k:I^2A \ra X$ with $F \circ Ii = k \circ i_0$.  Choose a homotopy $K:I^2X \ra X$ such that $F = K \circ i_0$ 
and $K \circ I^2 i = k$.  Write $K_{(t,T)}:X \ra X$ for the function $x \ra K((x,t),T)$.  Obviously, 
$K_{(0,0)} \approx$ 
$K_{(0,1)} \approx$ 
$K_{(1,1)} \approx$ 
$K_{(1,0)}$, 
all homotopies being $\rel A$.  Set $g = G \circ i_1$ $-$then $g \circ f = F \circ i_0 = K_{(0,0)}$ is homotopic $\rel A$ to 
$K_{(1,0)} = F \circ i_1 = \id_X$.]\\

\begin{proposition} \ %13
 Suppose that 
$
\begin{cases}
\ i:A \ra X\\[-.15cm]
\ j:A \ra Y
\end{cases}
$
are cofibrations.  Let $\phi \in C(X,Y)$: $\phi \circ i = j$.  Assume that $\phi$ is a homotopy equivalence $-$then $\phi$ is a 
homotopy equivalence in $A\backslash \bTOP$.
\end{proposition}

%%----------------------------------------------------------------------------------------------19
[Since $j$ is a cofibration, there exists a homotopy inverse $\psi:Y \ra X$ for $\phi$ with $\psi \circ j = i$, thus, 
from the lemma, $\exists$ $\psi^\prime \in C(X,X)$ : $\psi^\prime \circ i = i$ 
$\&$ $\psi^\prime \circ \psi \circ \phi \simeq \id_X \ \rel i(A)$.  
This says that $\phi^\prime = \psi^\prime \circ \psi$ is a homotopy left inverse for $\phi$ under \mA.  
Repeat the argument with 
$\phi$ replaced by $\phi^\prime$ to conclude that $\phi^\prime$ has a homotopy left inverse $\phi\pp$ under \mA, hence that 
$\phi^\prime$ is a homotopy equivalence in $A\backslash\bTOP$ or still, that $\phi$ is a homotopy equivalence in 
$A\backslash\bTOP$.]\\

\label{3.22}
\label{4.77}
\label{5.0n}
Application: Suppose that 
$
\begin{cases}
\ (X,x_0)\\
\ (Y,y_0)
\end{cases}
$
are wellpointed.  Let $f \in C(X,x_0;Y,y_0)$ $-$then $f$ is a homotopy equivalence in \bTOP iff $f$ is a homotopy equivalence in 
$\bTOP_*$.\\

\begingroup%%----------------------------------->>
\fontsize{9pt}{11pt}\selectfont
\textbf{\small FACT} \   
Suppose that $(X,x_0)$ is wellpointed.  Let $f \in C(X,Y)$ be inessential $-$then $f$ is homotopic in $\bTOP_*$ to the function 
$x \ra f(x_0)$.\\
\endgroup %%------------------------------------<<

\textbf{\small LEMMA}  \  
Suppose given a \cd 
\begin{tikzcd}%[ sep=small]
{A} \ar{d}[swap]{\phi} \ar{r}{i} &{X} \ar{d}{\psi}\\
{B} \ar{r}[swap]{j} &{Y}
\end{tikzcd}
in which 
$
\begin{cases}
\ i\\[-.15cm]
\ j
\end{cases}
$
are cofibrations and 
$
\begin{cases}
\ \phi\\[-.15cm]
\ \psi
\end{cases}
$
are homotopy equivalences.  Fix a homotopy inverse $\phi^\prime$ for $\phi$ and a homotopy $h_A:IA \ra A$ between 
$\phi^\prime \circ \phi$ and $\id_A$ $-$then there exists a homotopy inverse $\psi^\prime$ for $\psi$ with 
$i \circ \phi^\prime = \psi^\prime \circ j$ and a homotopy $H_X:IX \ra X$ between $\psi^\prime \circ \psi$ and $\id_X$ such that 
$
H_X(i(a),t) = 
\begin{cases}
\ i(h_A(a,2t)) \hspace{.5cm} (0 \leq t \leq 1/2)\\[-.15cm]
\ i(a) \hspace{1.65cm}  \ (1/2 \leq t \leq 1)
\end{cases}
.
$

[Fix some $\psi^\prime$ with $i \circ \phi^\prime = \psi^\prime \circ j$ (possible, $j$ being a cofibration).  
Put $h = i \circ h_A$: 
$h \circ i_0 = $ 
$i \circ h_A \circ i_0 = $ 
$i \circ \phi^\prime \circ \phi = $
$\psi^\prime \circ j \circ \phi =$ 
$\psi^\prime \circ \psi \circ i$ $\implies$ $\exists$ $H:IX \ra X$ such that 
$\psi^\prime \circ \psi = H \circ i_0$ and 
$H \circ Ii = h$.  Put $f = H \circ i_1$: $f \circ i = i \circ h_A \circ i_1 = i$ $\&$ 
$f \simeq H \circ i_0 =$ $\psi^\prime \circ \psi \simeq \id_X$ $\implies$ 
$\exists$ $g \in C(X,X)$: $g \circ i = i$ $\&$
$g \circ f \simeq \id_X \ \rel i(A)$.  
Let $G:IX \ra X$ be a homotopy between $g \circ f$ and $\id_X \ \rel i(A)$.  
Define $H_X:IX \ra X$ by 
$
H(X,t) = 
\begin{cases}
\ g(H(x,2t))  \hspace{.62cm} (0 \leq t \leq 1/2)\\[-.15cm]
\ G(x,2t - 1) \hspace{.5cm} (1/2 \leq t \leq 1)
\end{cases}
\hspace{-.25cm}: \ 
$
$H_X$ is a homotopy between $g \circ \psi^\prime \circ \psi$ and $\id_X$ and $H_X \circ Ii = i \circ h_A^\prime$, where 
$h_A^\prime(a,t) = h_A(a,\min\{2t,1\})$ is a homotopy between $\phi^\prime \circ \phi$ and $\id_A$.  Make the substitution 
$\psi^\prime \ra g \circ \psi^\prime$ to complete the proof.]\\

\begin{proposition} \ %14
Suppose given a \cd 
\begin{tikzcd}%[ sep=small]
{A} \ar{d}[swap]{\phi} \ar{r}{i} &{X} \ar{d}{\psi}\\
{B} \ar{r}[swap]{j} &{Y}
\end{tikzcd}
in which 
$
\begin{cases}
\ i\\[-.15cm]
\ j
\end{cases}
$
are cofibrations and 
$
\begin{cases}
\ \phi\\[-.15cm]
\ \psi
\end{cases}
$
are homotopy equivalences $-$then $(\phi,\psi)$ is a homotopy equivalence in $\bTOP(\ra)$.
\end{proposition}

%%------------------------------------------------------------------------------20
[The lemma implies that $(\phi^\prime,\psi^\prime)$ is a homotopy left inverse for $(\phi,\psi)$ in 
$\bTOP(\ra)$.]\\

\label{13.29} %dmc mnft

\begingroup%%----------------------------------->>
\fontsize{9pt}{11pt}\selectfont
\textbf{\small EXAMPLE} \quadx 
Let 
$
\begin{cases}
\ f:X \ra Y\\
\ f^\prime:X^\prime \ra Y^\prime
\end{cases}
$
be objects in $\bTOP(\ra)$.  
Write $[f,f^\prime]$ for the set of homotopy classes of maps in $\bTOP(\ra)$ from $f$ to $f^\prime$.  
Question: Is it true that if 
$
\begin{cases}
\ f \simeq g\\
\ f^\prime \simeq g^\prime
\end{cases}
$
(in \bTOP), then $[f,f^\prime] = [g,g^\prime]$?  The answer is ``no''.  Let $f = g$ be the constant map $\bS^1 \ra (1,0)$; 
let $f^\prime:\bS^1 \ra \bD^2$ be the inclusion and let $g^\prime:\bS^1 \ra \bD^2$ be the constant map at 
$(1,0)$ $-$then $[f,f^\prime] \neq [g,g^\prime]$.\\
\endgroup %%------------------------------------<<

\begin{proposition} \ %15
Let 
\begin{tikzcd}%[ sep=small]
{X^0} \ar{d}\ar{r} &{X^1} \ar{d}\ar{r}&{\cdots}\\
{Y^0} \ar{r} &{Y^1} \ar{r}&{\cdots}
\end{tikzcd}
be a commutative ladder connecting two expanding sequences of topological spaces.  
Assume: $\forall \ n$, the inclusions 
$
\begin{cases}
\ X^n \ra X^{n+1}\\[-.15cm]
\ Y^n \ra Y^{n+1}
\end{cases}
$
are cofibrations and the vertical arrows $\phi^n:X^n \ra Y^n$ are homotopy equivalences $-$then the induced map 
$\phi^\infty:X^\infty \ra Y^\infty$ is a homotopy equivalence.
\end{proposition}

[Using the lemma, inductively construct a homotopy left inverse for $\phi^\infty$.]\\

\label{13.15}
\label{14.70b}
\label{14.109}
\begingroup%%----------------------------------->>
\fontsize{9pt}{11pt}\selectfont
\textbf{\small FACT} \   
Let $X^0 \subset X^1 \subset \cdots $ be an expanding sequence of topological spaces.  Assume: $\forall \ n$, the inclusion 
$X^n \ra X^{n+1}$ is a cofibration and that $X^n$ is a strong deformation retract of $X^{n+1}$ $-$then $X^0$ is a strong deformation retract of $X^\infty$.
\\ \indent
[Bearing in mind Proposition 5, recall first that the inclusion $X^0 \ra X^\infty$ is a cofibration (cf. p. \pageref{3.16}).  
Consider the commutative ladder 
\begin{tikzcd}%[ sep=small]
{X^0} \ar{d}\ar{r} &{X^0} \ar{d}\ar{r}&{\cdots}\\
{X^0} \ar{r} &{X^1} \ar{r}&{\cdots}
\end{tikzcd}
to see that the inclusion $X^0 \ra X^\infty$ is also a homotopy equivalence.]\\
\endgroup %%------------------------------------<<

\label{4.71}
\begingroup%%----------------------------------->>
\fontsize{9pt}{11pt}\selectfont
\textbf{\small FACT} \   
Let $X^0 \subset X^1 \subset \cdots $ be an expanding sequence of topological spaces.  Assume: $\forall \ n$, the inclusion 
$X^n \ra X^{n+1}$ is a cofibration and inessential $-$then $X^\infty$ is contractible.\\
\endgroup %%------------------------------------<<

\label{5.0ad}
\begingroup%%----------------------------------->>
\fontsize{9pt}{11pt}\selectfont
\textbf{\small EXAMPLE} \quadx 
Take $X^n = \bS^n$ $-$then $X^\infty = \bS^\infty$ is contractible.\\
\endgroup %%------------------------------------<<

\label{12.28}
Let $f:X \ra Y$ be a continuous function $-$then the 
\un{mapping cylinder}
\index{mapping cylinder} 
$M_f$ of $f$ is defined by the pushout square 
\begin{tikzcd}%[ sep=small]
{X} \ar{d}[swap]{i_0} \ar{r}{f} &{Y} \ar{d}\\
{IX} \ar{r} &{M_f}
\end{tikzcd}
.  Special case: The mapping cylinder of $X \ra *$ is $\Gamma X$, the 
\un{cone}
\index{cone} 
of \mX (in particular, 
$\Gamma \bS^{n-1} = \bD^n$, so $\Gamma \emptyset = *$).  
There is a closed embedding $j:Y \ra M_f$, a homotopy 
$H:IX \ra M_f$, and a unique continuous function $r:M_f \ra Y$ such that $r \circ j = \id_Y$ and 
$r \circ H = f \circ p$ $(p:IX \ra X)$.  One has 
%%----------------------------------------------------------------------------------------------21
\label{5.27d}
$j \circ r \simeq \id_{M_f} \ \rel j(Y)$.  
The composition $H \circ i_1$ is a closed embedding $i:X \ra M_f$ and $f = r \circ i$.\\

\begingroup%%----------------------------------->>
\fontsize{9pt}{11pt}\selectfont
Suppose that \mX is a subspace of \mY and that $f:X \ra Y$ is the inclusion $-$then there is a continuous bijection 
$M_f \ra i_0Y \cup IX$.  In general, this bijection is not a homeomorphism (consider $X = \ ]0,1]$, $Y = [0,1]$) but will be if \mX is closed or $f$ is a cofibration.\\
\endgroup %%------------------------------------<<

\textbf{\small LEMMA}  \  
$j$ is a closed cofibration and $j(Y)$ is a strong deformation retract of $M_f$.\\

\label{12.15}
\textbf{\small LEMMA}  \  
$i$ is a closed cofibration.

[Define  \ $F:X \coprod X \ra Y \coprod X$ \  by \  $F = f \coprod \id_X$ \  and form the pushout square
\begin{tikzcd}%[ sep=small]
{X \coprod X} \ar{d}[swap]{i_0}\ar{d}{i_1} \ar{r}{F} &{Y \coprod X} \ar{d}\\
{IX} \ar{r} &{IX \sqcup_F (Y \coprod X)}
\end{tikzcd}
$-$then $IX \sqcup_F (Y \coprod X)$ can be identified with $M_f$, $i$ becoming the composite of the closed cofibrations 
$X \ra Y \coprod X \ra IX \sqcup_F (Y \coprod X)$.]\\

It is a corollary that the embedding $i$ of \mX into its cone $\Gamma X$ is a closed cofibration.\\

\begingroup%%----------------------------------->>
\fontsize{9pt}{11pt}\selectfont
\label{9.53}
\label{12.40}
\index{mapping telescope}
\textbf{\small EXAMPLE} \quadx 
The \un{mapping telescope} is the functor $\tel:\bFIL(\bTOP) \ra \bFILSP$ defined on an object $(\bX,\bff)$ by 
$\tel(\bX,\bff) = \ds\coprod\limits_n IX_n / \sim$, where $(x_n,1) \sim (f_n(x_n),0)$, and on a morphism 
$\phi:(\bX,\bff) \ra (\bY,\bg)$ by $\tel\phi([x_n,t]) = [\phi_n(x_n),t]$.  
Let $\telsub_n(\bX,\bff)$ be the image of 
$\left(\ds\coprod\limits_{k \leq n-1} IX_k \right) \ds\coprod i_9 X_n$, 
so $\telsub_n(\bX,\bff)$ is obtained from $X_n$ via iterated application of the mapping cylinder construction.  
The embedding $\telsub_n(\bX,\bff) \ra \telsub_{n+1}(\bX,\bff)$ is a closed cofibration and 
$\tel(\bX,\bff) = \colim \ \telsub_n(\bX,\bff)$.  
There is a homotopy equivalence $\telsub_n(\bX,\bff) \ra X_n$, viz. the assignment 
$[x_k,t] \ra (f_{n-1} \circ \cdots \circ f_k)(x_k)$ $(0 \leq k \leq n - 1)$, $[x_n,0] \ra x_n$ and the diagram 
\begin{tikzcd}[ sep=large]
{\telsub_n(\bX,\bff)} \ar{d}\ar{r} &{\telsub_{n+1}(\bX,\bff)} \ar{d}\\
{X_n} \ar{r} &{X_{n+1}}
\end{tikzcd}
commutes.  Consequently, if all the $f_n$ are cofibrations, then it follows from Proposition 15 that the induced map 
$\tel(\bX,\bff) \ra \colim X_n$ is a homotopy equivalence.
\\ \indent
[Note: \ Up to homeomorphism, the telescope construction is an instance of the above procedure.]\\
\endgroup %%------------------------------------<<

\begin{proposition} \ %16
Every morphism in \bTOP can be written as the composite of a closed cofibration and a homotopy equivalence.\\
\end{proposition}

\begin{proposition} \ %17
Let $f:X \ra Y$ be a continuous function $-$then $f$ is a homotopy equivalence iff $i(X)$ is a strong deformation retract of 
$M_f$.
\end{proposition}

[Note that $f$ is a homotopy equivalence iff $i$ is a homotopy equivalence and quote Proposition 5.]\\

%%----------------------------------------------------------------------------------------------22
\label{5.34}
Let $f:X \ra Y$ be a continuous function $-$then the 
\un{mapping cone}
\index{mapping cone} 
$C_f$ of $f$ is defined by the pushout square
\begin{tikzcd}%[ sep=small]
{X} \ar{d}[swap]{i} \ar{r}{f} &{Y} \ar{d}\\
{\Gamma X} \ar{r} &{C_f}
\end{tikzcd}
.   \ Special case: The mapping cone of $X \ra *$ is $\Sigma X$, the 
\un{suspension}
\index{suspension} 
of \mX 
(in particular, $\Sigma \bS^{n-1} = \bS^n$, so $\Sigma \emptyset = \bS^0$).  There is a closed cofibration 
$j:Y \ra C_f$ and an arrow $C_f \ra \Sigma X$.  By construction, $j \circ f$ is inessential and for a any $g:Y \ra Z$ with 
$g \circ f$ inessential, there exists a $\phi:C_f \ra Z$ such that $g = \phi \circ j$.

[Note: \ The 
\un{mapping cone sequence}
\index{mapping cone sequence} 
associated with $f$  is given by 
$X \overset{f}{\ra} $ 
$Y \ra $ 
$C_f \ra $ 
$\Sigma X \ra $ 
$\Sigma Y \ra $ 
$\Sigma C_f \ra $ 
$\Sigma^2 X \ra $ 
$\cdots$.  
Taking into account the suspension isomorphism $\widetilde{H}_q(X) \approx \widetilde{H}_{q+1}(\Sigma X)$, there is an exact sequence 
\[
\cdots \ra \widetilde{H}_q(X) \ra \widetilde{H}_q(Y) \ra \widetilde{H}_q(C_f ) \ra \widetilde{H}_{q-1}(X) 
\ra \widetilde{H}_{q-1}(Y) \ra \cdots .]
\]
\vspace{0.1cm}

\label{3.1}%%dmc???
\begingroup%%----------------------------------->>
\fontsize{9pt}{11pt}\selectfont
The mapping cylinder and the mapping cone can be viewed as functors $\bTOP(\ra) \ra \bTOP$.  With this interpretation, $i$, $j$ and $r$ are natural transformations.
\\ \indent
[Note: \ Owing to AD$_4$, these functors restrict to functors $\bHAUS(\ra) \ra \bHAUS$.  Consequently, if \mX and \mY are in \bCGH, then for any continuous function $f:X \ra Y$, both $M_f$ and $C_f$ remain in \bCGH.  On the other hand, 
stability relative to \bCG or $\dcg$ is automatic.]\\
\endgroup %%------------------------------------<<

\begingroup%%----------------------------------->>
\fontsize{9pt}{11pt}\selectfont
\textbf{\small FACT} \   
Suppose that 
$
\begin{cases}
\ f:X \ra Y\\
\ g:X \ra Y
\end{cases}
$
are homotopic $-$then in $\bHTOP^2$, $(M_f,i(X)) \approx (M_g,i(X))$, and in \bHTOP, $C_f \approx C_g$.\\
\endgroup %%------------------------------------<<

\begingroup%%----------------------------------->>
\fontsize{9pt}{11pt}\selectfont
\textbf{\small FACT} \   
Let $f \in C(X,Y)$.  Suppose that $\phi:X^\prime \ra X$ $(\psi:Y \ra Y^\prime)$ is a homotopy equivalence $-$then the arrow 
$(M_{f \circ \phi},i(X^\prime)) \ra (M_f,i(X))$ $((M_f,i(X))  \ra (M_{\psi \circ f},i(X)))$ is a homotopy equivalence  
(in $\bTOP^2$) and the arrow $C_{f \circ \phi} \ra C_f$ ($C_f \ra C_{\psi \circ f}$) is a homotopy equivalence (in \bTOP).\\
\endgroup %%------------------------------------<<

\label{3.33}

\begingroup%%----------------------------------->>
\fontsize{9pt}{11pt}\selectfont
\textbf{\small EXAMPLE} \quadx 
The suspension $\Sigma X$ of \mX is the union of two closed subspaces 
$\Gamma^- X$ and $\Gamma^+ X$, each homeomorphic to the cone 
$\Gamma X$ of \mX, with 
$\Gamma^- X \cap \Gamma^+ X = X$ (identify the section $i_{1/2}X$ with \mX).  Therefore 
$\Sigma X$ is numerably contractible.  The commutative diagram 
\begin{tikzcd}[ sep=large]
{X} \ar{d}\ar{r} &{\Gamma^+ X} \ar{d}\\
{\Gamma^- X} \ar{r} &{\Sigma X}
\end{tikzcd}
is a pushout square and the inclusions 
$
\begin{cases}
\ \Gamma^- X \ra \Sigma X\\
\ \Gamma^+ X \ra \Sigma X
\end{cases}
$
are closed cofibrations.\\
\vspace{0.25cm}
\endgroup %%------------------------------------<<

\label{3.25}

\begingroup%%----------------------------------->>
\fontsize{9pt}{11pt}\selectfont
\textbf{\small FACT} \   
Let $f:X \ra Y$ be a continuous function.  Suppose that \mY is numerably contractible $-$then $C_f$ is numerably contractible.
\\ \indent
%%----------------------------------------------------------------------------------------------23
[The image of $X \times [0,1[$ in $C_f$ is contractible.  On the other hand, the image of $X \times ]0,1] \amalg Y$ 
in $C_f$ has the same homotopy type as \mY, hence is numerably contractible (cf. p. \pageref{3.17}).]
\\ \indent
[Note: \ \mY and $M_f$ have the same homotopy type, so \mY numerably contractible $\implies$ $M_f$ numerably contractible (cf. p. \pageref{3.18}).]\\ 
\endgroup %%------------------------------------<<

\label{9.52}
\label{12.37}
Let $X \overset{f}{\lla} Z \overset{g}{\lra} Y$ be a 2-source $-$then the 
\un{double mapping cylinder}
\index{double mapping cylinder} 
$M_{f,g}$ of $f$, $g$ is defined by the pushout square
\begin{tikzcd}%[ sep=small]
{Z \coprod Z} \ar{d}[swap]{i_0}\ar{d}{i_1} \ar{r}{f \coprod g} &{X \coprod Y} \ar{d}\\
{IZ} \ar{r} &{M_{f,g}}
\end{tikzcd}
.   \ The homotopy type of $M_{f,g}$ depends only on the homotopy classes of $f$ and $g$ and $M_{f,g}$ is homeomorphic to 
$M_{g,f}$.  
There are closed cofibrations 
$
\begin{cases}
\ i:X \ra M_{f,g}\\
\ j:Y \ra M_{f,g}
\end{cases}
$
and an arrow $M_{f,g} \ra \Sigma Z$.  
The diagram 
\begin{tikzcd}[ sep=large]
{Z} \ar{d}[swap]{f} \ar{r}{g} &{Y} \ar{d}{j}\\
{X} \ar{r}[swap]{i} &{M_{f,g}}
\end{tikzcd}
is homotopy commutative and if the diagram 
\begin{tikzcd}[ sep=large]
{Z} \ar{d}[swap]{f} \ar{r}{g} &{Y} \ar{d}{\eta}\\
{X} \ar{r}[swap]{\xi} &{W}
\end{tikzcd}
is homotopy commutative, then there exists a $\phi:M_{f,g} \ra W$ such that 
$
\begin{cases}
\ \xi = \phi \circ i\\
\ \eta = \phi \circ j
\end{cases}
\hspace{-.25cm}. \ 
$
Example: The double mapping cylinder of 
$X \la X \times Y \ra Y$ is $X * Y$, the 
\un{join}
\index{join} 
of \mX and \mY.

\label{4.16}
[Note: \ The mapping cylinder and the mapping cone are instances of the double mapping cylinder (homeomorphic models arise from the parameter reversal $t \ra 1 - t$).  Consideration of 
$
\begin{cases}
\ Z \times [0,1/2]\\
\ Z \times [1/2,1]
\end{cases}
$
leads to a pushout square
\begin{tikzcd}[ sep=large]
{Z} \ar{d} \ar{r} &{M_g} \ar{d}\\
{M_f} \ar{r} &{M_{f,g}}
\end{tikzcd}
.]\\
\vspace{0.25cm}

\begingroup%%----------------------------------->>
\fontsize{9pt}{11pt}\selectfont
\label{4.30}
\label{5.0z}
\index{Mapping Telescope (example)}
\textbf{\small EXAMPLE  \ (\un{The Mapping Telescope})} \  
$\tel(\bX,\bff)$ can be identified with the double mapping
cylinder of the 2-source 
$\ds\coprod\limits_{n \geq 0} X_{2n} \la \ds\coprod\limits_{n \geq 0} X_{n} \ra \ds\coprod\limits_{n \geq 0} X_{2n+1}$.  
Here, the left hand arrow is defined by $x_{2n} \ra x_{2n}$ $\&$ $x_{2n+1} \ra f_{2n+1}(x_{2n+1})$ and the right hand arrow is defined by 
$x_{2n+1} \ra x_{2n+1}$ $\&$ $x_{2n} \ra f_{2n}(x_{2n})$.\\
\endgroup %%------------------------------------<<

Every 2-source $X \overset{f}{\lla} Z \overset {g}{\lra} Y$ determines a pushout square
\begin{tikzcd}[ sep=large]
{Z} \ar{d}[swap]{f} \ar{r}{g} &{Y} \ar{d}{\eta}\\
{X} \ar{r}[swap]{\xi} &{P}
\end{tikzcd}
and there is an arrow $\phi:M_{f,g} \ra P$ characterized by the conditions 
$
\begin{cases}
\ \xi = \phi \circ i\\
\ \eta = \phi \circ j
\end{cases}
$
$\&$ $IZ \ra M_{f,g} \overset{\phi}{\ra} P \ = \ $ 
$
\begin{cases}
\ \xi \circ f \circ p\\[-.2cm]
\hspace{.65cm} ||\\
\ \eta \circ g \circ p
\end{cases}
\hspace{-.25cm}.
$
\\
\vspace{0.25cm}

\begin{proposition} \ %18
If $f$ is a cofibration, then $\phi:M_{f,g} \ra P$ is a homotopy equivalence in $Y\backslash \bTOP$.
\end{proposition}

%%----------------------------------------------------------------------------------------------24
[The arrow $M_f \ra IX$ admits a left inverse $IX \ra M_f$.]\\

\label{4.42}
\label{5.18}
\label{5.55d}
\label{9.34}
\label{13.11} %dmc mnft
\label{13.18} %dmc mnft
\label{14.70a} %dmc mnft
Application: Suppose that $f:X \ra Y$ is a cofibration $-$then the projection $C_f \ra Y/f(X)$ is a 
homotopy equivalence.

\label{5.6}
[Note: \ If in addition \mX is contractible, then the embedding $Y \ra C_f$ is a homotopy equivalence.  Therefore in this case the projection $Y \ra Y/f(X)$ is a homotopy equivalence.]\\

\begingroup%%----------------------------------->>
\fontsize{9pt}{11pt}\selectfont
\textbf{\small EXAMPLE} \quadx 
Let \mA be a nonempty finite subset of $\bS^n$ $(n \geq 1)$ $-$then $\bS^n/A$ has the homotopy type of the wedge of $\bS^n$ with $(\#(A) - 1)$ circles.
\\ \indent
[The inclusion $A \ra \bS^n$ is a cofibration (cf. Proposition 8).]\\
\endgroup %%------------------------------------<<

\label{5.0i}
\label{12.16}
\label{12.17}
Consider the 2-sources
$
\begin{cases}
\ X \la A \overset{f}{\ra} Y\\[-.2cm]
\ X \la A \underset{g}{\ra} Y
\end{cases}
\hspace{-.25cm} , \ 
$
where the arrow $A \ra X$ is a closed cofibration.   
Assume that $f \simeq g$ $-$then Proposition 18 implies that 
$X \sqcup_f Y$ and $X \sqcup_g Y$ have the same homotopy type $\rel Y$.  
Corollary: If $f^\prime : A \ra Y^\prime$ is a continuous function and if $\phi:Y \ra Y^\prime$ 
is a homotopy equivalence such that $\phi \circ f \simeq f^\prime$, then there is a homotopy equivalence 
$\Phi:X \sqcup_f Y \ra X \sqcup_{f^\prime} Y^\prime$ with $\restr{\Phi}{Y} = \phi$.\\

\begingroup%%----------------------------------->>
\fontsize{9pt}{11pt}\selectfont
\textbf{\small FACT} \   
Suppose that $A \ra X$ is a closed cofibration.  Let $f:A \ra Y$ be a homotopy equivalence $-$then the arrow 
$X \ra X \sqcup_f Y$ is a homotopy equivalence.\\
\endgroup %%------------------------------------<<

\begingroup%%----------------------------------->>
\fontsize{9pt}{11pt}\selectfont
Denote by $\abs{\Delta,\id}_{\bTOP}$ the comma category corresponding to the diagonal functor 
$\Delta:\bTOP \ra \bTOP \times \bTOP$ and the identity functor id on $\bTOP \times \bTOP$.  
So, an object in 
$\abs{\Delta,\id}_{\bTOP}$ is a 2-source $X \overset{f}{\la} Z \overset{g}{\ra} Y$ 
and a morphism of 2-sources is a commutative diagram 
\begin{tikzcd}[sep=large]
{X} \ar{d} &{Z}\ar{l}[swap]{f} \ar{d} \ar{r}{g} &{Y}\ar{d}\\
{X^\prime}  &{Z^\prime}\ar{l}{f^\prime}  \ar{r}[swap]{g^\prime} &{Y^\prime}
\end{tikzcd}
.  The double mapping cylinder is a functor $\abs{\Delta,\id}_{\bTOP} \ra \bTOP$.  It has a right adjoint 
$\bTOP \ra \abs{\Delta,\id}_{\bTOP}$, viz. the functor that sends \mX to the 2-source 
$X \overset{p_0}{\lla} PX \overset{p_1}{\lra} X$.\\
\endgroup %%------------------------------------<<

\label{3.20}
\label{5.0aa}
\label{14.76} %dmc mnft
\begingroup%%----------------------------------->>
\fontsize{9pt}{11pt}\selectfont
\textbf{\small FACT} \   
Let 
\begin{tikzcd}[sep=large]
{X} \ar{d} &{Z}\ar{l}[swap]{f} \ar{d} \ar{r}{g} &{Y}\ar{d}\\
{X^\prime}  &{Z^\prime}\ar{l}{f^\prime}  \ar{r}[swap]{g^\prime} &{Y^\prime}
\end{tikzcd}
be a \cd in which the vertical arrows are homotopy equivalences $-$then the arrow 
$M_{f,g} \ra M_{f^\prime,g^\prime}$ is a homotopy equivalence.\\
\endgroup %%------------------------------------<<

\label{4.18}
\label{4.26}
\label{5.4}
\begingroup%%----------------------------------->>
\fontsize{9pt}{11pt}\selectfont
Application: Suppose that 
$
\begin{cases}
\ A \ra X\\
\ A^\prime \ra X^\prime
\end{cases}
$
are closed cofibrations.  Let 
$
\begin{cases}
\ f:A \ra Y\\
\ f^\prime:A^\prime \ra Y^\prime
\end{cases}
$
be continuous
%%----------------------------------------------------------------------------------------------25
functions.  
Assume that the diagram
\begin{tikzcd}[ sep=large]
{X} \ar{d} &{A}\ar{l} \ar{d} \ar{r}{f} &{Y} \ar{d}\\
{X^\prime}  &{A^\prime} \ar{l} \ar{r}[swap]{f^\prime} &{Y^\prime}
\end{tikzcd}
commutes and that the vertical arrows are homotopy equivalences $-$then the induced map 
\label{5.0x}
$X \sqcup_f Y \ra X^\prime \sqcup_{f^\prime} Y^\prime$ is a homotopy equivalence.\\
\endgroup %%------------------------------------<<

\begingroup%%----------------------------------->>
\fontsize{9pt}{11pt}\selectfont
\textbf{\small EXAMPLE} \  \ 
Suppose that $X = A \ \cup \  B$, where 
$
\begin{cases}
\ A\\
\ B
\end{cases}
$
are closed and the inclusions 
$
\begin{cases}
\ A \cap B \ra A\\
\ A \cap B \ra B
\end{cases}
$
are cofibrations.  
Assume: \mA and \mB are contractible $-$then the arrow $\Sigma (A \cap B) \ra X$ is a homotopy equivalence.\\
\vspace{0.25cm}
\endgroup %%------------------------------------<<

\index{Segal-Stasheff Construction}
\begingroup%%----------------------------------->>
\fontsize{9pt}{11pt}\selectfont
\textbf{\small SEGAL-STASHEFF CONSTRUCTION}  \ \ 
Let \mX be a topological space.  
Fix a covering $\sU = \{U_i:i \in I\}$ of \mX.  Equip \mI with a well ordering $<$ and put 
$I[n] = \{[i] \equiv (i_0,\ldots,i_n):i_0 < \cdots < i_n\}$.  
Every strictly increasing $\alpha \in \Mor([m],[n])$ defines a map 
$I[n] \ra I[m]$.  
Set $U_{[i]} = U_{i_0} \cap \cdots \cap U_{i_n}$ and form 
$\sU([n]) = \ds\coprod\limits_{I[n]} U_{[i]}$, a coproduct in \bTOP.  
Give $\sU([n]) \times \Delta^n$ the product topology and call $B\sU$ the quotient 
$\ds\coprod\limits_n \sU([n]) \times \Delta^n / \sim$, the equivalence relation being generated by writing 
$((x,[i]),\Delta^\alpha t) \sim ((x,\alpha[i]),t)$.  Let $B\sU^{(n)}$ be the image of 
$\ds\coprod\limits_{m \leq n} \sU([m]) \times \Delta^m$ in $B\sU$, so $B\sU = \colimx B\sU^{(n)}$.  
The commutative diagram
\endgroup

\begingroup
\fontsize{9pt}{11pt}\selectfont
\[
\begin{tikzcd}[ sep=large]
{\ds\coprod\limits_{I[n]} U_{[i]} \times \ddpn} \ar{d} \ar{r} &{B\sU^{(n-1)}} \ar{d}\\
{\ds\coprod\limits_{I[n]} U_{[i]} \times \dpn} \ar{r} &{B\sU^{(n)}}
\end{tikzcd}
\]
is a pushout square in \bTOP and the vertical arrows are closed cofibrations.  
There is a projection $p_{\sU}:B\sU \ra X$ induced by the arrows 
$U_{[i]} \times \Delta^n \ra U_{[i]}$, i.e., $((x,[i]),t) \ra x$.  
Moreover, $p_{\sU}$ is a homotopy equivalence provided that $\sU$ is numerable.  
Indeed, any partition of unity $\{\kappa_i:i \in I\}$ on \mX subordinate to $\sU$ determines a continuous function 
$s_{\sU}:X \ra B\sU$ (since $\forall \ x$, $\#\{i \in I: x \in \sptx \kappa_i\} < \omega$).  
Obviously, $p_{\sU} \circ s_{\sU} = \id_X$ and $s_{\sU}\circ p_{\sU}$  can be connected to the identity on $B\sU$ via a linear homotopy.\\
\endgroup %%------------------------------------<<

\label{4.67a}
\begingroup%%----------------------------------->>
\fontsize{9pt}{11pt}\selectfont
\textbf{\small FACT} \   
Let 
$
\begin{cases}
\ X\\
\ Y
\end{cases}
$
be topological spaces and let $f:X \ra Y$ be a continuous function.  
Suppose that 
$
\begin{cases}
\ \sU = \{U_i: i \in I\}\\
\ \sV = \{V_i: i \in I\}
\end{cases}
$
are numerable coverings of 
$
\begin{cases}
\ X\\
\ Y
\end{cases}
$
such that $\forall \ i$: $f(U_i) \subset V_i$.  Assume: $\forall \ [i]$, the induced map 
$f_{[i]}:U_{[i]} \ra V_{[i]}$ is a homotopy equivalence $-$then $f$ is a homotopy equivalence.

\label{4.32}
[There is an arrow $F:B\sU \ra B\sV$ and a \cd 
\begin{tikzcd}[ sep=large]
{B\sU} \ar{d}[swap]{p_{\sU}} \ar{r}{F} &{B\sV} \ar{d}{p_{\sV}}\\
{X} \ar{r}[swap]{f} &{Y}
\end{tikzcd}
.  
Due to the numerability of $\sU$ and $\sV$, $p_{\sU}$ and $p_{\sV}$ are homotopy equivalences.  
Claim: $\forall \ n$, the restriction $F^{(n)}:B\sU^{(n)} \ra B\sV^{(n)}$ is a homotopy equivalence.  
This is clear if $n = 0$.  
For $n > 0$, consider the commutative 
%%----------------------------------------------------------------------------------------------26
diagram
\[
\begin{tikzcd}[ sep=large]
{\ds\coprod\limits_{I[n]} U_{[i]} \times \dpn} \ar{d}  
&{\ds\coprod\limits_{I[n]} U_{[i]} \times \ddpn} \ar{l} \ar{d} \ar{r}
&{B\sU^{(n-1)}} \ar{d}\\
{\ds\coprod\limits_{I[n]} V_{[i]} \times \dpn} 
&{\ds\coprod\limits_{I[n]} V_{[i]} \times \ddpn} \ar{l} \ar{r} 
&{B\sV^{(n-1)}}
\end{tikzcd}
\]
By induction, $F^{(n-1)}$ is a homotopy equivalence, thus $F^{(n)}$ is too.  Proposition 15 then implies that 
$F:B\sU \ra B\sV$ is a homotopy equivalence, so the same is true of $f$.]\\
\endgroup %%------------------------------------<<

Let $u,v:X \ra Y$ be a pair of continuous functions $-$then the 
\un{mapping torus}
\index{mapping torus} 
$T_{u,v}$ of $u,v$ is defined by the pushout square
\begin{tikzcd}[ sep=large]
{X \amalg X} \arrow[r,"v",swap,shift right=1] \arrow[r,"u", shift left=1]{v} \ar{d}[swap]{i_0}\ar{d}{i_1}
&{Y} \ar{d}\\
{IX} \ar{r} &{T_{u,v}}
\end{tikzcd}
.  There is a closed cofibration $j:Y \ra T_{u,v}$.  
From the definitions, $j \circ u \simeq j \circ v$ and for any $g:Y \ra Z$ 
with $g \circ u \simeq g \circ v$, there exists a $\phi:T_{u,v} \ra Z$ such that $g = \phi \circ j$.

[Note: \ If $u = v = \id_X$, then $T_{u,v}$ is the product $X \times \bS^1$.]\\

\begingroup%%----------------------------------->>
\fontsize{9pt}{11pt}\selectfont
\index{The Scorpion (example)}
\textbf{\small EXAMPLE \ (\un{The Scorpion})} \  
Let $\pi:\bS^n \ra \bD^n$ be the restriction of the canonical map $\R^{n+1} \ra \R^n$; let 
$p:\bD^n \ra \bD^n/\bS^{n-1} = \bS^n$ be the projection.  
Put $f = p \circ \pi$ $-$then $f:\bS^n \ra \bS^n$ is inessential.  
The \un{scorpion} $\sS^{n+1}$ is the quotient of $I\bS^n$ with respect to the relations $(x,0) \sim (f(x),1)$, 
i.e., $\sS^{n+1}$ is the mapping torus of $x \ra f(x)$ $\&$ $x \ra x$ $(x \in \bS^n)$.  
One may also describe $\sS^{n+1}$  as the quotient 
$\bD^{n+1}/\sim$, where $x \sim p(2x)$ $(x \in (1/2)\bD^n)$.  
Fix a point $x_0 \in (1/2)\bS^{n-1}$, let $L_0$ be the line segment from $x_0$ to $p(2x_0)$, and let $C_0$ be the circle $L_0/\sim$ $-$then the inclusion $C_0 \ra S^{n+1}$ is a homotopy equivalence, thus $\sS^{n+1}$ is a homotopy circle.  
The 
\un{dunce hat}
\index{dunce hat} 
$\sD^{n+1}$ is the quotient 
$\sS^{n+1}/C_0$.  It is contractible.\\
\endgroup %%------------------------------------<<

The formalities in $\bTOP_*$ run parallel to those in \bTOP, thus a detailed account of the pointed theory is unnecessary.  
Of course, there is an important difference between \bTOP and $\bTOP_*$: $\bTOP_*$ has a zero object but \bTOP does not.  
Consequently, if 
$
\begin{cases}
\ (X,x_0)\\[-.15cm]
\ (Y,y_0)
\end{cases}
$
are in $\bTOP_*$, then $[X,x_0;Y,y_0]$ is a pointed set with distinguished element $[0]$, the pointed homotopy class of the zero morphism, i.e., of the constant map $X \ra y_0$.  Functions $f \in [0]$ are said to be 
\un{nullhomotopic}:
\index{nullhomotopic} 
$f \simeq 0$.

\label{4.6}
[Note: \ The forgetful functor $\bTOP_* \ra \bTOP$ has a left adjoint $\bTOP \ra \bTOP_*$ that sends the space \mX 
to the pointed space $X_+ = X\coprod *$.]

The computation of pushouts in $\bTOP_*$ is expedited by noting that a pushout in \bTOP of a 2-source in $\bTOP_*$ 
is a pushout in $\bTOP_*$.  \ 
Examples: 
(1)  \ The pushout square \ \ 
\begin{tikzcd}[ sep=large]
{*} \ar{d} \ar{r} &{(Y,y_0)} \ar{d}\\
{(X,x_0)} \ar{r} &{X \vee Y}
\end{tikzcd}
defines the 
\un{wedge}
\index{wedge} 
$X \vee Y$; 
(2) The pushout square
%%----------------------------------------------------------------------------------------------27
\begin{tikzcd}[ sep=large]
{X \vee Y} \ar{d} \ar{r} &{*} \ar{d}\\
{X \times Y} \ar{r} &{X \# Y}
\end{tikzcd}
defines the 
\un{smash product}
\index{smash product} 
$X \# Y$.

[Note: \ Base points are suppressed if there is no need to display them.]\\

\label{14.138}
\begingroup%%----------------------------------->>
\fontsize{9pt}{11pt}\selectfont
The wedge is the coproduct in $\bTOP_*$.  If both of the inclusions 
$
\begin{cases}
\ \{x_0\} \ra X\\
\ \{y_0\} \ra Y
\end{cases}
$
are cofibrations and if at least one is closed, then the embedding 
$X \vee Y \ra X \times Y$ is a cofibration (cf. Proposition 7) and $X \vee Y$ is wellpointed (cf. Proposition 9).\\
\endgroup %%------------------------------------<<

\begingroup%%----------------------------------->>
\fontsize{9pt}{11pt}\selectfont
\textbf{\small FACT} \   
Suppose that 
$
\begin{cases}
\ (X,x_0)\\
\ (Y,y_0)
\end{cases}
$
are in $\bTOP_*$ $-$then $\forall \ n > 1$, there is a split short exact sequence
\[
0 \ra \pi_{n+1}(X \times Y,X \vee Y) \ra \pi_n(X \vee Y) \ra \pi_n(X \times Y) \ra 0.
\]
\endgroup %%------------------------------------<<

\begingroup%%----------------------------------->>
\fontsize{9pt}{11pt}\selectfont
Griffiths\footnote[2]{\textit{Quart. J. Math.} \textbf{5} (1954), 175-190.}
proved that if $(X,x_0)$ is a path connected pointed Hausdorff space which is both first countable and locally simply connected at 
$x_0$, then for any path connected pointed Hausdorff space $(Y,y_0)$, the arrow 
$\pi_1(X,x_0)*\pi_1(Yy_0) \ra \pi_1((X,x_0) \vee (Y,y_0))$ is an isomorphism.

[Note: \ \mX is 
\un{locally simply connected}
\index{locally simply connected} 
at $x_0$ provided that for any neighborhood \mU of $x_0$ there exists a neighborhood $V \subset U$ of $x_0$ such that the induced homomorphism 
$\pi_1(V,x_0) \ra \pi_1(U,x_0)$ is trivial.]

Eda\footnote[3]{\textit{Proc. Amer. Math. Soc.} \textbf{109} (1990), 237-241; 
see also Morgan-Morrison, \textit{Proc. London Math. Soc.} \textbf{53} (1986), 562-576.} 
has constructed an example of a path connected CRH space \mX which is locally simply connected at $x_0$ with the property that $\pi_1(X,x_0) = 1$ but $\pi_1((X,x_0) \vee (X,x_0)) \neq 1$.  Moral: The hypothesis of first countability cannot be dropped.\\
\endgroup %%------------------------------------<<

\index{The Hawaiian Earring (example)}
\begingroup%%----------------------------------->>
\fontsize{9pt}{11pt}\selectfont
\textbf{\small EXAMPLE \  (\un{The Hawaiian Earring})} \  
Let \mX be the subspace of $\R^2$ consisting of the union of the circles $X_n$, where $X_n$ has center $(1/n,0)$ and radius 
$1/n$ $(n \geq 1)$.  \ \ 
\parbox{4cm}{
\begin{tikzpicture}[scale=1.75]%[scale=0.5,shift={(-5,-3)}]
\draw[violet, thick]
(1,0) circle (1);
\draw[teal, thick]
(0.5,0) circle (0.5);
\draw[violet, thick]
(0.333333,0) circle (0.333333);
\draw[teal, thick]
(0.25,0) circle (0.25);
\draw[violet, thick]
(0.2,0) circle (0.2);
\draw[teal, thick]
(0.166667,0) circle (0.166667);
\draw[violet, thick]
(0.142857,0) circle (0.142857);
\draw[teal, thick]
(0.125,0) circle (0.125);
\draw[violet, thick]
(0.111111,0) circle (0.111111);
\draw[teal, thick]
(0.1,0) circle (0.1);
\draw[violet, thick]
(0.0909091,0) circle (0.0909091);
\draw[teal, thick]
(0.0833333,0) circle (0.0833333);
\draw[violet, thick]
(0.0769231,0) circle (0.0769231);
\draw[teal, thick]
(0.0714286,0) circle (0.0714286);
\draw[violet, thick]
(0.0666667,0) circle (0.0666667);
\draw[teal, thick]
(0.0625,0) circle (0.0625);
\draw[violet, thick]
(0.0588235,0) circle (0.0588235);
\end{tikzpicture}
}
\\
Take $x_0 = (0,0)$ $-$then \mX is first countable at $x_0$, \mX is not locally simply connected at $x_0$, 
the inclusion $\{x_0\} \ra X$ is not a cofibration, and the arrow 
$\pi_1(X,x_0)*\pi_1(X,x_0) \ra \pi_1((X,x_0) \vee (X,x_0))$ is injective but not surjective.  
Denote now by $X_0$ the result of assigning to \mX the final topology determined by the inclusions
$X_n \ra X$.  $X_0$ is a CW complex.  
Take $x_0 = (0,0)$ $-$then $X_0$ is not first countable at $x_0$, $X_0$ is locally simply connected at $x_0$, 
the inclusion $\{x_0\} \ra X$ is a cofibration, and the arrow 
$\pi_1(X,x_0)*\pi_1(X,x_0) \ra \pi_1((X,x_0)\vee (X,x_0))$ is an isomorphism (Van Kampen).\\
\endgroup %%------------------------------------<<

\begingroup%%----------------------------------->>
\fontsize{9pt}{11pt}\selectfont
\textbf{\small FACT} \   
Given a wellpointed space $(X,x_0)$, suppose that $X = A \cup B$, where $x_0 \in A \cap B$ and $A \cap B$ is contractible.  
Assume: The inclusions 
$
\begin{cases}
\ A \cap B \ra A\\
\ A \cap B \ra B
\end{cases}
$
$\&$ 
$
\begin{cases}
\ A \ra X\\
\ B \ra X
\end{cases}
$
are cofibrations.  Take
$
\begin{cases}
\ a_0 = x_0\\
\ b_0 = x_0
\end{cases}
$
$-$then the arrow $A \vee B \ra X$ is a pointed homotopy equivalence.\\

%%----------------------------------------------------------------------------------------------28
The smash product $\#$ is a functor $\bTOP_* \times \bTOP_* \ra \bTOP_*$.  It respects homotopies, 
thus the pointed homotopy type of $X \# Y$ depends only on the pointed homotopy types of \mX and \mY.  \label{3.29}If both of the inclusions 
$
\begin{cases}
\ \{x_0\} \ra X\\
\ \{y_0\} \ra Y
\end{cases}
$
are cofibrations and if at least one is closed, then $X \# Y$ is wellpointed.

[Note: \ Suppose that \mY is a pointed LCH space $-$then it is clear that the functor 
$-\#Y:\bTOP_* \ra \bTOP_*$ has a right adjoint $Z \ra Z^Y$ which passes to $\bHTOP_*$: $[X\# Y,Z] \approx [X,Z^Y]$, $Z^Y$ 
the set of pointed continuous functions from \mY to \mZ equipped with the compact open topology.  One can say more: In fact, 
Cagliari\footnote[2]{\textit{Proc. Amer. Math. Soc.} \textbf{124} (1996), 1265-1269.}
has shown that for any pointed \mY, the functor $-\#Y$ has a right adjoint in $\bTOP_*$ iff the functor $-\times Y$ has a right adjoint in \bTOP, i.e., iff \mY is core compact (cf. p. \pageref{3.19}).]
\\
\indent\indent ($\#_1$) \ $X \# Y$ is homeomorphic to $Y \# X$.
\\
\indent\indent ($\#_2$) \ $(X \# Y) \# Z$ is homeomorphic to $X \# (Y \# Z)$ if both \mX and \mZ are LCH spaces or if two of \mX, \mY, \mZ are compact Hausdorff.

[Note: \ The smash product need not be associative (consider $(\Q \# \Q) \# \Z$ and $\Q \# (\Q \# \Z)$).]
\\
\indent\indent ($\#_3$) \ $(X \vee Y) \# Z$ is homeomorphic to $(X \# Z) \vee(Y \# Z)$.
\\
\indent\indent $(\#_4$) \ $\Sigma (X * Y)$ is homeomorphic to $\Sigma X \# \Sigma Y$ if \mX and \mY are compact Hausdorff.
\\ \indent
[Note: \ The suspension can be viewed as a functor $\bTOP \ra \bTOP_*$.  This is because the suspension is the result of collapsing to a point the embedded image of a space in its cone.  
Example: $\bS^{m-1}*\bS^{n-1} = \bS^{m+n-1}$ $\implies$ $\bS^m \# \bS^n = \bS^{m+n}$.]
\endgroup %%------------------------------------<<

\begingroup%%----------------------------------->>
\fontsize{9pt}{11pt}\selectfont
All the homeomorphisms figuring in the foregoing are natural and preserve the base points.\\
\endgroup %%------------------------------------<<

\label{14.133a}
\textbf{\small LEMMA}  \  
\begingroup%%----------------------------------->>
\fontsize{9pt}{11pt}\selectfont
The smash product of two pointed Hausdorff spaces is Hausdorff.

The pushout square 
\begin{tikzcd}[ sep=large]
{X \vee Y} \ar{d} \ar{r} &{*} \ar{d}\\
{X \times_k Y} \ar{r} &{X \#_k Y}
\end{tikzcd}
\label{1.21}
defines the 
\un{smash product}
\index{smash product ($\bCG$, $\dcg$, or $\bCGH$)} 
$X \#_k Y$ in $\bCG$, $\dcg$, or $\bCGH$.  It is associative and distributes over the wedge.
\\ \indent
[Note: \ With $\#_k$ as the multiplication and $\bS^0$ as the unit, $\bCG_*$, $\bDelta$-$\bCG_*$, and $\bCGH_*$ 
are closed categories.]\\
\endgroup %%------------------------------------<<

\label{14.182}
\label{16.4}
\label{16.5}
The 
\un{pointed cylinder functor}
\index{pointed cylinder functor} 
$I:\bTOP_* \ra \bTOP_*$ is the functor that sends 
$(X,x_0)$ to the quotient $X \times [0,1]/\{x_0\} \times [0,1]$, i.e., $I(X,x_0) = IX/I\{x_0\}$.  
Variant: Let $I_+ = [0,1] \amalg *$ $-$then $I(X,x_0)$ is the smash product $X \# I+$.   
The 
\un{pointed path space functor}
\index{pointed path space functor} 
$P:\bTOP_* \ra \bTOP_*$ is the functor that sends 
$(X,x_0)$ to $C([0,1],X)$ (compact open topology), the base point for the latter being the constant path 
$[0,1] \ra x_0$.  As in the unpointed situation, $(I,P)$ is an adjoint pair.

%%----------------------------------------------------------------------------------------------29
Using \mI and \mP, one can define the notion of a pointed cofibration.  
Since all maps and homotopies must respect the base points, an arrow $A \ra X$ in $\bTOP_*$ may be a pointed cofibration without being a cofibration.  For example, 
$\forall \ x_0 \in X$, the arrow $(\{x_0\},x_0) \ra (X,x_0)$ is a pointed cofibration but in general the inclusion 
$\{x_0\} \ra X$ is not a cofibration.  On the other hand, an arrow $A \ra X$ in $\bTOP_*$ which is a cofibration, when considered as an arrow in \bTOP, is necessarily a pointed cofibration.  Pointed cofibrations are embeddings.  If 
$x_0 \in A \subset X$ and if $\{x_0\}$ is closed in \mX, then the inclusion $A \ra X$ is a pointed cofibration iff 
$i_0X \cup IA/I\{x_0\}$ is a retract of $I(X,x_0)$.  
Observe that for this it is not necessary that \mA itself be closed.

Let $(X,A,x_0)$ be a pointed pair $-$then a  
\un{Str{\o}m structure}
\index{Str{\o}m structure (pointed pair)} 
on $(X,A,x_0)$ consists of a continuous function 
$\phi:X \ra [0,1]$ such that $A \subset \phi^{-1}(0)$, a continuous function $\psi:X \ra [0,1]$ such that 
$\{x_0\} = \psi^{-1}(0)$, and a homotopy $\Phi:IX \ra X$ of $\id_X \ \rel A$ such that $\Phi(x,t) \in A$ whenever 
$\min\{t,\psi(x)\} > \phi(x)$.

[Note: \ $\Phi$ is therefore a pointed homotopy.]\\

\index{Theorem: Pointed Cofibration Characerization Theorem}
\index{Pointed Cofibration Characerization Theorem}
\textbf{\small POINTED COFIBRATION CHARACTERIZATION THEOREM}
%\endgroup %%------------------------------------<<
\quadx 
Let $x_0 \in A \subset X$ and suppose that $\{x_0\}$ is a zero set in \mX 
$-$then the inclusion $A \ra X$ is a pointed cofibration iff the 
pointed pair $(X,A,x_0)$ admits a Str{\o}m structure.

[Necessity: Fix 
$\psi \in C(X,[0,1])$ : $\{x_0\} = \psi^{-1}(0)$ and let $X \overset{p}{\lla} IX \overset{q}{\lra} [0,1]$ 
be the projections.  
Put $Y = \{(x,t) \in i_0X \cup IA: t \leq \psi(x)\}$.  
Define a continuous function $f:i_0X \cup IA \ra Y$ by 
$f(x,t) = (x,\min\{t,\psi(x)\}$ and let $F:IX \ra Y$ be some continuous extension of $f$.  
Consider 
$\phi(x) = \sup\limits_{0 \leq t \leq 1}\abs{\min\{t,\psi(x)\} - qF(x,t)}$, $\Phi(x,t) = pF(x,t)$.

Sufficiency: Given a Str{\o}m structure $(\phi,\psi,\Phi)$ on $(X,A,x_0)$, define a retraction 
$r:I(X,x_0) \ra i_0X \cup IA/I\{x_0\}$ by 
\[
r(x,t) = 
\begin{cases}
\ (\Phi(x,t),0) \hspace{2.9cm} (t\psi(x) \leq \phi(x))\\
\ (\Phi(x,t),t - \phi(x)/\psi(x)) \hspace{0.75cm}   (t\psi(x) > \phi(x))
\end{cases}
.]
\]
\\

\textbf{\small LEMMA}  \  
Let $(X,A,x_0)$ be a pointed pair.  Suppose that the inclusions 
$
\begin{cases}
\ \{x_0\} \ra A\\[-.15cm]
\ \{x_0\} \ra X
\end{cases}
$
are closed cofibrations and that the inclusion $A \ra X$ is a pointed cofibration $-$then the pair $(X,x_0)$ has a 
Str{\o}m structure $(f,F)$ for which $F(IA) \subset A$.

[Fix a  Str{\o}m structure $(f_X,F_X)$ on $(X,x_0)$.  Choose a Str{\o}m structure $(\phi,\psi,\Phi)$ on $(X,A,x_0)$ such that $\phi \leq \psi = f_X$.  
Fix a  Str{\o}m structure $(f_A,F_A)$ on $(A,x_0)$.  Extend the pointed homotopy 
$i \circ F_A:IA \ra A \overset{i}{\ra} X$ to a pointed homotopy $\ov{F}:IX \ra X$ with $\ov{F} \circ i_0 = \id_X$.  Put
\[
\ov{f}(x)
\ = \ 
\begin{cases}
\ (1 - \phi(x)/\psi(x)) f_A(\Phi(x,1)) + \phi(x) \hspace{0.75cm} (\phi(x) < \psi(x))\\[-.15cm]
\ \psi(x) \hspace{5.75cm} (\phi(x) = \psi(x))
\end{cases}
.
\]
%%----------------------------------------------------------------------------------------------30
Then $\ov{f} \in C(X,[0,1])$, $\restr{\ov{f}}{A} = f_A$, and $\ov{f}^{-1}(0) = \{x_0\}$.  Consider 
$f(x) = \min\{1,\ov{f}(x) + f_X(\ov{F}(x,1))\}$, 
\[
F(x,t) = 
\begin{cases}
\ \ov{F}(x,t/\ov{f}(x)) \hspace{2.25cm} (t < \ov{f}(x))\\[-.15cm]
\ F_X(\ov{F}(x,1), t - \ov{f}(x)) \hspace{0.75cm}  (t \geq \ov{f}(x))
\end{cases}
.]
\]
\\

\begin{proposition} \ %19
Let $(X,A,x_0)$ be a pointed pair.  Suppose that the inclusions 
$
\begin{cases}
\ \{x_0\} \ra A\\[-.15cm]
\ \{x_0\} \ra X
\end{cases}
$
are closed cofibrations $-$then the inclusion $A \ra X$ is a cofibration iff it is a pointed cofibration.
\end{proposition}

[To establish the nontrivial assertion, take $(f,F)$ as in the lemma and choose a Str{\o}m structure 
$(\ov{\phi},\ov{\psi},\ov{\Phi})$ on $(X,A,x_0)$ with $\ov{\phi} \leq \ov{\psi} = f$.  Define  a Str{\o}m structure $(\phi,\Phi)$ on $(X,A)$ by $\phi(x) = \ov{\phi}(x) - \ov{\psi}(x) + \sup\limits_{0 \leq t \leq 1} \ov{\psi}(\ov{\Phi}(x,t))$, 
\[
\Phi(x,t) = F(\ov{\Phi}(x,t), \min\{t,\ov{\phi}(x)/\ov{\psi}(x)\}) \hspace{0.75cm} (x \neq x_0)
\]
and $\Phi(x_0,t) = x_0$.]\\

So, under conditions commonly occurring in practice, the pointed and unpointed notions of cofibration are equivalent.\\

Let $X \overset{f}{\la} Z \overset{g}{\ra} Y$ be a pointed 2-source $-$then there is an embedding $M_{*,*} \ra M_{f,g}$ 
and the quotient $M_{f,g}/M_{*,*}$ is the pointed double mapping cylinder of $f, g$.  
Here $M_{*,*}$ is the double mapping cylinder of the 2-source $* \la * \ra *$, which, being $* \times [0,1]$, is contractible.  
Thus if \mX, \mY, and \mZ are wellpointed, then 
$M_{f,g}/M_{*,*}$ is wellpointed and the projection $M_{f,g} \ra M_{f,g}/M_{*,*}$ is a homotopy equivalence 
(cf. p. \pageref{3.20}).

[Note: \ The pointed mapping torus of a pair $u,v:X \ra Y$ of pointed continuous functions is 
$T_{u,v}/T_{*,*}$, where $T_{*,*}$ is $* \times \bS^1$, which is not contractible.]\\

\begingroup%%----------------------------------->>
\fontsize{9pt}{11pt}\selectfont
The \cd 
\begin{tikzcd}[ sep=large]
{Iz_0} \ar{d} 
&{z_0 \amalg z_0} \ar{d} \ar{l} \ar{r}
&{x_0 \amalg y_0} \ar{d}\\
{IZ} 
&{Z \amalg Z} \ar{l}[swap]{i_0}\ar{l}{i_1} \ar{r}[swap]{f \amalg g}
&{X \amalg Y}
\end{tikzcd}
leads to an induced map of pushouts $Iz_0 \ra M_{f,g}$ which we claim is a cofibration.  Thus, since 
$
\begin{cases}
\ X\\
\ Y
\end{cases}
$
are wellpointed, the arrow $x_0 \amalg y_0 \ra X \amalg Y$ is a cofibration.  
On the other hand, the pushout of the 2-source 
$Iz_0 \la z_0 \ds\amalg z_0 \ra Z \ds\amalg Z$ can be identified with $i_0Z \cup Iz_0 \cup i_1Z$ 
(even though $z_0$ is not assumed to be closed) and the inclusion 
$i_0Z \cup Iz_0 \cup i_1Z \ra IZ$ is a cofibration (cf. p. \pageref{3.21}).  
The claim is then seen to be a consequence of the proof of Proposition 4 in $\S 12$ 
(which depends only on the fact that cofibrations are pushout stable (cf. Proposition 2)).  Consideration of
%%----------------------------------------------------------------------------------------------31
the pushout square 
\begin{tikzcd}[ sep=large]
{Iz_0} \ar{d}  \ar{r} &{*} \ar{d}\\
{M_{f,g}} \ar{r} &{M_{f,g}/M_{*,*}}
\end{tikzcd}
now implies that $M_{f,g}/M_{*,*}$ is wellpointed.  Finally, one can view $M_{f,g}$ itself as a wellpointed space (take 
$[z_0,1/2]$ as the base point).  The projection $M_{f,g} \ra M_{f,g}/M_{*,*}$ is therefore a homotopy equivalence between wellpointed spaces, hence is actually a pointed homotopy equivalence (cf. p. \pageref{3.22}).\\
\label{3.23}
\label{3.24}
\endgroup %%------------------------------------<<

In particular: There are pointed versions of $\Gamma X$ and $\Sigma X$ of the cone and suspension of a pointed space \mX.  Each is a quotient of its unpointed counterpart (and has the same homotopy type if \mX is wellpointed).  $\Sigma X$ is a 
cogroup in $\bHTOP_*$.  In terms of the smash product, $\Gamma X = X \# [0,1]$ (0 the base point of $[0,1]$) and 
$\Sigma X = X \# \bS^1$ ($(1,0)$ the base point of $\bS^1$).  
Example: $\Gamma(X \vee Y) = \Gamma X \vee \Gamma Y$ and $\Sigma(X \vee Y) = \Sigma X \vee \Sigma Y$.  
The 
\un{mapping space functor}
\index{mapping space functor} 
$\Theta:\bTOP_* \ra \bTOP_*$ is the functor that sends 
$(X,x_0)$ to the subspace of $C([0,1),X)$ consisting of those $\sigma$ such that $\sigma(0) = x_0$ and the 
\un{loop space functor}
\index{loop space functor} 
$\Omega:\bTOP_* \ra \bTOP_*$ is the functor that sends 
$(X,x_0)$ to the subspace of $C([0,1),X)$ consisting of those $\sigma$ such that $\sigma(0) = x_0 = \sigma(1)$, 
the base point in either case being the constant path $[0,1] \ra x_0$.  $\Omega X$ is a group object in $\bHTOP_*$.  
$(\Gamma,\Theta)$ and $(\Sigma,\Omega)$ are adjoint pairs.  Both drop to 
$\bHTOP_*:[\Gamma X,Y] \approx [X,\Theta Y]$ and $[\Sigma X,Y] \approx [X,\Omega Y]$.

[Note: \ If \mX is wellpointed, then so are $\Theta  X$ and $\Omega X$.]\\

\begingroup%%----------------------------------->>
\fontsize{9pt}{11pt}\selectfont
The mapping space $\Theta X$ is contractible and there is a pullback square
\begin{tikzcd}[ sep=large]
{\Omega X} \ar{d}  \ar{r} &{\Theta X} \ar{d}{p_1}\\
{x_0} \ar{r} &{X}
\end{tikzcd}
in \bTOP, hence in $\bTOP_*$.\\
\endgroup %%------------------------------------<<

\label{14.72}
\index{The Moore Loop Space (example)}
\begingroup%%----------------------------------->>
\fontsize{9pt}{11pt}\selectfont
\textbf{\small EXAMPLE  \ (\un{The Moore Loop Space})} \  
Given a pointed space $(X,x_0)$, let $\Omega_M X$ be the set of all pairs $(\sigma,r_\sigma)$: $\sigma \in C([0,r_\sigma],X)$ 
$(0 \leq r_\sigma < \infty)$ and $\sigma(0) = x_0 = \sigma(r_\sigma)$.  
Attach to each $(\sigma,r_\sigma) \in \Omega_M X$ the function $\ov{\sigma}(t) = \sigma(\min\{t,r_\sigma\})$ on $\R_{\geq 0})$ $-$then the assignment 
$(\sigma,r_\sigma) \ra (\ov{\sigma},r_\sigma)$ injects $\Omega_M X$ into $C(\R_{\geq 0},X) \times \R_{\geq 0}$.  
Equip $\Omega_M X$ with the induced topology from the product (compact open topology on $C(\R_{\geq 0},X)$).  
Define an associative multiplication on $\Omega_M X$ by writing 
$
(\tau + \sigma)(t) = 
\begin{cases}
\ \sigma(t) \hspace{1.05cm}   (0 \leq t \leq r_\sigma)\\
\ \tau(t - r_\sigma) \hspace{0.35cm} (r_\sigma \leq t \leq r_{\tau + \sigma})
\end{cases}
\hspace{-.2cm},
$
where $r_{\tau + \sigma} = r_\tau + r_\sigma$, the unit thus being $(0,0)$ $(0 \ra x_0)$.  Since ``+'' is continuous, 
$\Omega_M X$ is a monoid in \bTOP, the \un{Moore loop space} of \mX, and $\Omega_M$ is a functor 
$\bTOP_* \ra \bMON_{\bTOP}$.  The inclusion $\Omega X \ra \Omega_M X$ is an embedding (but it is not a pointed map).
\endgroup %%------------------------------------<< %not sure why had to resort tp this break

\begingroup%%----------------------------------->>
\fontsize{9pt}{11pt}\selectfont
Claim: $\Omega X$ is a deformation retract of $\Omega_M X$.

[Consider the homotopy $H:I\Omega_M X \ra \Omega_M X$ defined as follows.  
The domain of $H((\sigma,r_\sigma),t)$ is the interval $[0,(1 - t)r_\sigma + t]$ and there
%%----------------------------------------------------------------------------------------------32
\[
H((\sigma,r_\sigma),t)(T) = \sigma\left(\frac{Tr_\sigma}{(1 - t)r_\sigma + t}\right)
\]
if $r_\sigma > 0$, otherwise $H((0,0),t)(T) = x_0$.]

One can also introduce $\Theta_M X$, the 
\un{Moore mapping space}
\index{Moore mapping space} 
of \mX.  Like 
$\Theta X$, $\Theta_M X$ is contractible and evaluation at the free end defines a Hurewicz fibration $\Theta_M X \ra X$ whose fiber over the base point is $\Omega_M X$.\\
\endgroup %%------------------------------------<<

Let $f:X \ra Y$ be a pointed continuous function, $C_f$ its pointed mapping cone.\\

\textbf{\small LEMMA}  \  
If $f$ is a pointed cofibration, then the projection $C_f \ra Y/f(X)$ is a pointed homotopy equivalence.\\

\label{4.53}
In general, there is a pointed cofibration $j:Y \ra C_f$ and an arrow $C_f \ra \Sigma X$.  Iterate to get a pointed cofibration $j^\prime:C_f \ra C_j$ $-$then the triangle 
\begin{tikzcd}%[ sep=small]
{C_f} \ar{rd} \ar{r} &{C_j}\ar{d}\\
&{\Sigma X}
\end{tikzcd}
commutes and by the lemma, the vertical arrow is a pointed homotopy equivalence.  
Iterate again to get a pointed cofibration 
$j\pp:C_j \ra C_{j^\prime}$ $-$then the triangle 
\begin{tikzcd}%[ sep=small]
{C_j} \ar{rd} \ar{r} &{C_{j^\prime}}\ar{d}\\
&{\Sigma Y}
\end{tikzcd}
commutes and by the lemma, the vertical arrow is a pointed homotopy equivalence.  
Example: Given pointed spaces 
$
\begin{cases}
\ X\\
\ Y
\end{cases}
, \ 
$
let $X \ov{\#} Y$\label{3.31} be the pointed mapping cone of the inclusion $f:X \vee Y \ra X \times Y$ $-$then in $\bHTOP_*$, 
$C_j \approx \Sigma(X \vee Y)$ and $C_{j^\prime} \approx \Sigma(X \times Y)$.\\

Let $f:X \ra Y$ be a pointed continuous function $-$then the 
\un{pointed mapping cone} \un{sequence}
\index{pointed mapping cone sequence} 
associated with $f$ is given by 
$X \overset{f}{\ra}$ 
$Y \ra $
$C_f \ra $ 
$\Sigma X \ra $ 
$\Sigma Y \ra $ 
$\Sigma C_f \ra $ 
$\Sigma^2 X \ra $ 
$\cdots$.  
Example:  When $f = 0$, this sequence becomes 
$X \overset{0}{\ra}$ 
$Y \ra $
$Y \vee \Sigma X \ra $ 
$\Sigma X \ra $ 
$\Sigma Y \ra $ 
$\Sigma Y \vee \Sigma^2 X \ra $ 
$\Sigma^2 X \ra $ 
$\cdots$. 

[Note: \ If the diagram 
\begin{tikzcd}[ sep=large]
{X} \ar{d} \ar{r}{f} &{Y} \ar{d}\\
{X^\prime}\ar{r}[swap]{f^\prime} &{Y^\prime}
\end{tikzcd}
commutes in $\bHTOP_*$ and if the vertical arrows are pointed homotopy equivalences, then the pointed mapping cone sequences of $f$ and $f^\prime$ are connected by a commutative ladder in $\bHTOP_*$, all of whose vertical arrows are pointed homotopy equivalences.]\\

%%----------------------------------------------------------------------------------------------33
\index{Theorem: Replication Theorem}
\index{Replication Theorem}
\textbf{\small REPLICATION THEOREM} \ \  
Let $f:X \ra Y$ be a pointed continuous function $-$then for any pointed space \mZ, there is an exact sequence
\[
\cdots \ra [\Sigma Y,Z] \ra [\Sigma X,Z] \ra [C_f,Z] \ra [Y,Z] \ra [X,Z]
\]
in $\bSET_*$.

[Note: \ A sequence of pointed sets and pointed functions 
$(X,x_0) \overset{\phi}{\ra} $ 
$(Y,y_0) \overset{\psi}{\ra} $ 
$(Z,z_0)$ is said to be 
\un{exact}
\index{exact (sequence of pointed sets and pointed functions)} 
in $\bSET_*$ if the range of $\phi$ is equal to the kernel of $\psi$.]\\

\begingroup%%----------------------------------->>
\fontsize{9pt}{11pt}\selectfont
\label{5.32}
\label{5.55}
\textbf{\small EXAMPLE} \quadx 
Let $f:X \ra Y$ be a pointed continuous function, \mZ a pointed space.  Given pointed continuous functions 
$\alpha:\Sigma X \ra Z$, $\phi:C_f \ra Z$, write
$
(\alpha \cdot \phi)[x,t] =
\begin{cases}
\ \alpha(x,2t) \hspace{0.9cm} (0 \leq t \leq 1/2)\\
\ \phi(x,2t-1) \hspace{0.35cm} (1/2 \leq t \leq 1)
\end{cases}
(x \in X)
$
$\&$ $(\alpha \cdot \phi)(y) = \phi(y)$ $(y \in Y)$ $-$then this prescription defines a left action of 
$[\Sigma X,Z]$ on $[C_f,Z]$ and the orbits are the fibers of the arrow $[C_f,Z] \ra [Y,Z]$.\\
\endgroup %%------------------------------------<<

\label{3.30}
\label{4.52}
\label{14.43}
\label{14.123a}
\begingroup%%----------------------------------->>
\fontsize{9pt}{11pt}\selectfont
\textbf{\small FACT} \   
Given a pointed continuous function $f:X \ra Y$ and a pointed space \mZ, put $f_Z = f \# \id_Z$ $-$then there is a commutative ladder 
\[
\begin{tikzcd}[ sep=large]
{X \# Z} \ar{d}[swap]{\id} \ar{r}
&{Y \# Z} \ar{d}[swap]{\id} \ar{r}
&{C_{f_Z}} \ar{d} \ar{r}
&{\Sigma(X \# Z)}\ar{d} \ar{r} 
&{\Sigma(Y \# Z)}\ar{d} \ar{r}
&{\cdots}\\
{X \# Z}  \ar{r}
&{Y \# Z} \ar{r}
&{C_f \# Z} \ar{r}
&{\Sigma X \# Z} \ar{r}
&{\Sigma Y \# Z} \ar{r}
&{\cdots}
\end{tikzcd}
\]
in $\bHTOP_*$, all of whose vertical arrows are pointed homotopy equivalences.
\\ \indent
[Show that there are mutually inverse pointed homotopy equivalences
$
\begin{cases}
\ \phi:C_f \# Z \ra C_{f_Z}\\
\ \psi:C_{f_Z} \ra C_f \# Z 
\end{cases}
$
for which the triangles
\[
\begin{tikzcd}%[ sep=small]
&{C_f \# Z}\ar{dd}{\phi}\\
{Y \# Z} \ar{ru} \ar{rd}\\
&{C_{f_Z}}
\end{tikzcd}
\indent\indent
\begin{tikzcd}%[ sep=small]
&{C_f \# Z}\\
{Y \# Z} \ar{ru} \ar{rd}\\
&{C_{f_Z}} \ar{uu}[swap]{\psi}
\end{tikzcd}
\]
commute.]\\
\endgroup %%------------------------------------<<

Given a pointed space $(X,x_0)$, let $\widecheck X$ be the mapping cylinder of the inclusion 
$\{x_0\} \ra X$ and denote by $\widecheck{x}_0$ the image of $x_0$ under the embedding $i:\{x_0\} \ra \widecheck{X}$ 
$-$then 
$(\widecheck{X},\widecheck{x}_0)$ is wellpointed and $\{\widecheck{x}_0\}$ is closed in $\widecheck{X}$ 
(cf. p. \pageref{3.23}).  
The embedding $j:X \ra \widecheck{X}$ is a closed cofibration (cf. p. \pageref{3.24}).  
It is not a pointed map but the retraction 
$r:\widecheck{X} \ra X$ is both a pointed map and a homotopy equivalence.  
We shall term $(X,x_0)$ 
\un{nondegenerate}
\index{nondegenerate (pointed space)} 
if $r:\widecheck{X} \ra X$ is a pointed homotopy equivalence.

[Note: \ Consider $X \vee [0,1]$, where $x_0 = 0$ $-$then $\widecheck{X}$ is homeomorphic to $X \vee [0,1]$ with 
$\widecheck{x}_0 \leftrightarrow 1$.]\\

%%----------------------------------------------------------------------------------------------34
\begingroup%%----------------------------------->>
\fontsize{9pt}{11pt}\selectfont
\textbf{\small FACT} \   
Suppose that 
$
\begin{cases}
\ (X,x_0)\\
\ (Y,y_0)
\end{cases}
$
are nondegenerate.  
Assume: 
$
\begin{cases}
\ X\\
\ Y
\end{cases}
$
are numerably contractible $-$then $X \vee Y$ and $X \# Y$ are numerably contractible.

[To discuss $X \# Y$, take
$
\begin{cases}
\ (X,x_0)\\
\ (Y,y_0)
\end{cases}
$
wellpointed with 
$
\begin{cases}
\ \{x_0\} \subset X\\
\ \{y_0\} \subset Y
\end{cases}
$
closed.  The mapping cone of the inclusion $X \vee Y \ra X \times Y$ is numerably contractible (cf. p. \pageref{3.25}) and has the homotopy type of $X \times Y/X \vee Y = X \# Y$, which is therefore numerably contractible.]\\
\endgroup %%------------------------------------<<

\begingroup%%----------------------------------->>
\fontsize{9pt}{11pt}\selectfont
\textbf{\small FACT} \   
Suppose that 
$
\begin{cases}
\ (X,x_0)\\
\ (Y,y_0)
\end{cases}
$
are nondegenerate.  Let $f \in C(X,x_0;Y,y_0)$ $-$then the pointed mapping cone $C_f$ is numerably contractible provided that \mY is numerably contractible.
\\ \indent
[Consider the \cd
\begin{tikzcd}[ sep=large]
{X \vee [0,1]} \ar{d} \ar{r}{f \vee \id} &{Y \vee [0,1]} \ar{d}\\
{X} \ar{r}[swap]{f} &{Y}
\end{tikzcd}
.  
By hypothesis, the vertical arrows are pointed homotopy equivalences, so $C_{f \vee \id}$ and $C_f$ have the same pointed homotopy type.  Look at the unpointed mapping cone of $f \vee \id$.]\\
\endgroup %%------------------------------------<<

\begingroup%%----------------------------------->>
\fontsize{9pt}{11pt}\selectfont
Application: The pointed suspension of any nondegenerate space is numerably contractible.\\
\endgroup %%------------------------------------<<

A pointed space $(X,x_0)$ is said to satisfy 
\un{Puppe's condition}
\index{Puppe's condition} 
provided that there exists a halo \mU of $\{x_0\}$ in \mX and a homotopy 
$\Phi:IU \ra X$ of the inclusion $U \ra X \ \rel\{x_0\}$ such that 
$\Phi \circ i_1(U) = \{x_0\}$.  
Every wellpointed space satisfies Puppe's condition.\\

\textbf{\small LEMMA}  \  
Let $(X,A,x_0)$ be a pointed pair.  Suppose that there exists a pointed homotopy $H:IX \ra X$ of $\id_X$ such that 
$H \circ i_1(A) = \{x_0\}$ and 
$H \circ i_t(A) \subset A$ $(0 \leq t \leq 1)$ $-$then the projection $X \ra X/A$ is a pointed homotopy equivalence.\\

\begin{proposition} \ 
Let $(X,x_0)$ be a pointed space $-$then $(X,x_0)$ is nondegenerate iff it satisfies Puppe's condition.
\end{proposition}

[Necessity: Let $\rho:X \ra \cX$ be a pointed homotopy inverse for $r$.  
Fix a homotopy $H:IX \ra X$ of $\id_X \  \rel\{x_0\}$ such that 
$H \circ i_1 = r \circ \rho$.  
Put $U = \rho^{-1}(\{x_0\} \times ]0,1])$ 
$-$then \mU is a halo of $\{x_0\}$ in \mX with haloing function $\pi$ the composite 
$X \overset{\rho}{\ra} \cX \ra \cX/X = [0,1]$.  Consider $\Phi = \restr{H}{IU}$.

Sufficiency: One can assume that \mU is closed (cf. p. \pageref{3.26}).  Set
\[
\Phi^\prime(x,t) = 
\begin{cases}
\ \Phi(x,2t) \hspace{0.75cm}  (\in X \subset \cX) \hspace{1.13cm}   (0 \leq t \leq 1/2)\\
\ 2t - 1 \hspace{1.05cm} (\in [0,1] \subset \cX) \hspace{0.75cm} (1/2 \leq t \leq 1)
\end{cases}
\ \  (x \in U).
\]
Define a pointed homotopy $H:I\cX \ra \cX$ by
\[
(H \circ \restr{i_t}{X})(x) =
\begin{cases}
\ x \hspace{2.25cm} (x \notin U)\\
\ \Phi^\prime(x,t\pi(x)) \hspace{0.5cm} (x \in U)
\end{cases}
\]
%%----------------------------------------------------------------------------------------------35
and
\[
(H \circ \restr{i_t}{[0,1]})(T) =
\begin{cases}
\ T \hspace{3.42cm} (0 \leq t \leq 1/2)\\
\ 1 - (1 - T)(2 - 2t) \hspace{0.5cm} (1/2 \leq t \leq 1)
\end{cases}
.
\]
The lemma implies that $r:\cX \ra \cX/[0,1] = X$ is a pointed homotopy equivalence.]\\

\begingroup%%----------------------------------->>
\fontsize{9pt}{11pt}\selectfont
\textbf{\small EXAMPLE} \quadx 
Take $X = [0,1]^\kappa$ $(\kappa > \omega)$ and let $x_0 = 0_\kappa$, the ``origin'' in \mX 
$-$then $(X,x_0)$ is not wellpointed (cf. p. \pageref{3.27}) but is nondegenerate.\\
\endgroup %%------------------------------------<<

\label{5.27e}
\begingroup%%----------------------------------->>
\fontsize{9pt}{11pt}\selectfont
\textbf{\small FACT} \   
A pointed space $(X,x_0)$ is nondegenerate iff it has the same pointed homotopy type as $(\widecheck{X},\widecheck{x}_0)$.\\
\endgroup %%------------------------------------<<

\label{9.96}
\begingroup%%----------------------------------->>
\fontsize{9pt}{11pt}\selectfont
Application: Nondegeneracy is a pointed homotopy type invariant.
\\ \indent
[Note: \ Compare this with the remark on p. \pageref{3.28}.]\\
\endgroup %%------------------------------------<<

\label{4.23}
\begingroup%%----------------------------------->>
\fontsize{9pt}{11pt}\selectfont
\textbf{\small FACT} \   \
Suppose that 
$
\begin{cases}
\ (X,x_0)\\
\ (Y,y_0)
\end{cases}
$
are nondegenerate.  Let $f \in C(X,x_0;Y,y_0)$ $-$then $f$ is a homotopy equivalence in \bTOP iff $f$ is a homotopy equivalence in 
$\bTOP_*$.\\
\endgroup %%------------------------------------<<

\begingroup%%----------------------------------->>
\fontsize{9pt}{11pt}\selectfont
\index{The Moore Loop Space (example)}
\textbf{\small EXAMPLE  \ (\un{The Moore Loop Space})}  \  
Suppose that the pointed space \mX is nondegenerate $-$then $\Omega X$ and $\Omega_M X$ are nondegenerate.  Since the retraction of $\Omega_M X$ onto $\Omega X$ is not only a homotopy equivalence in \bTOP but a pointed map as well, it follows that $\Omega X$ and $\Omega_M X$ have the same pointed homotopy type.\\
\endgroup %%------------------------------------<<

\label{5.0o}
\begin{proposition} \ %21
Let $(X,x_0)$ be a pointed space $-$then $(X,x_0)$ is wellpointed and $\{x_0\}$ is closed in \mX iff $(X,x_0)$ is nondegenerate and $\{x_0\}$ is a zero set in \mX.
\end{proposition}

[This is a consequence of Propositions 10 and 20.]\\

As noted above, nondegeneracy is a pointed homotopy type invariant.  It is also a relatively stable property: \mX nondegenerate $\implies$ $\Gamma X$, $\Sigma X$, $\Theta X$, $\Omega X$ nondegenerate and \mX, \mY nondegenerate 
$\implies$ $X \times Y$, $X \vee Y$, $X \# Y$ nondegenerate.\\

\begingroup%%----------------------------------->>
\fontsize{9pt}{11pt}\selectfont
To illustrate, consider $X \# Y$.  In $\bHTOP_*$, $X \# Y \approx \cX \# \cY$, and since 
$
\begin{cases}
\ \{\widecheck{x}_0\} \ra \cX\\
\ \{\widecheck{y}_0\} \ra \cY
\end{cases}
$
are closed cofibrations, $\cX \# \cY$ is wellpointed (cf. p. \pageref{3.29}), hence a fortiori, nondegenerate.  Thus the same is true of $X \# Y$.\\
\endgroup %%------------------------------------<<

Given pointed spaces $(X_1,x_1), \ldots, (X_n,x_n)$ , write $X_1\Delta \cdots \Delta X_n$ for the subspace
\[
(\{x_1\} \times X_2 \times \cdots \times X_n) \cup \cdots \cup (X_1 \times \cdots \times  X_{n-1} \times \{x_n\})
\]
%%----------------------------------------------------------------------------------------------36
of $X_1 \times \cdots \times X_n$ and let $X_1 \# \cdots \# X_n$ be the quotient 
$X_1 \times \cdots \times X_n/X_1\Delta \cdots \Delta X_n$.\\

\begin{proposition} \ %22
Let \mX, \mY, \mZ be nondegenerate $-$then $(X \# Y) \# Z$ and $X\#(Y\#Z)$ have the same pointed homotopy type.
\end{proposition}

[There is a pointed 2-source
$(X\# Y)\# Z \la X \# Y \# Z \ra X \# (Y \# Z)$ arising from the identity.  Both arrows are continuous bijections and it will be enough to show that they are pointed homotopy equivalences.  For this purpose, consider instead the pointed 2-source 
$(\cX \# \cY) \# \cZ \la \cX \# \cY \# \cZ \ra \cX \# (\cY \# \cZ)$ and, to be specific, work on the left, calling the arrow $\phi$.  
Define pointed continuous functions 
$
\begin{cases}
\ u:\cX \ra \cX\\[-.15cm]
\ v:\cY \ra \cY
\end{cases}
$
by 
$
\begin{cases}
\ (\restr{u}X)(x) = x\\[-.15cm]
\ (\restr{v}Y)(y) = y
\end{cases}
$
$\&$ 
$
\begin{cases}
\ (\restr{u}{[0,1]})(t) = \max\{0,2t - 1\}\\[-.15cm]
\ (\restr{v}{[0,1]})(t) = \max\{0,2t - 1\}
\end{cases}
$ 
$-$then $u \times v \times \id_Z$ induces a pointed function 
$\psi:(\cX \# \cY) \# \cZ  \ra \cX \# \cY \# \cZ$.  
To check that $\psi$ is continuous, introduce closed subspaces 
$
\begin{cases}
\ A\\[-.15cm]
\ B
\end{cases}
$
of $\cX \# \cY$: Points of \mA are represented by pairs $(x,y)$, where $x \geq 1/2$ $(y \in \cY)$ or $y \geq 1/2$ $(x \in \cX)$, and points of \mB are represented by pairs $(x,y)$, where 
$
\begin{cases}
\ x \in X\\[-.15cm]
\ y \in Y
\end{cases}
$ 
or 
$
\begin{cases}
\ x \leq 1/2 \ (y \in Y) \\[-.15cm]
\ y \leq 1/2 \  (x \in X)
\end{cases}
$
 or $x \leq 1/2$ $\&$ $y \leq 1/2$.  Since the projection $(\cX \# \cY) \times \cZ \ra (\cX \# \cY) \# \cZ$ is closed, the images 
 $
\begin{cases}
\ A_Z\\[-.15cm]
\ B_Z
\end{cases}
$
of 
$
\begin{cases}
\ A \times \cZ\\[-.15cm]
\ B \times \cZ
\end{cases}
$
 in $(\cX \# \cY) \# \cZ$ are closed and their union fills out $(\cX \# \cY) \# \cZ$.  The continuity of $\psi$ is a consequence of the continuity of $\restr{\psi}{A_Z}$ and 
$\restr{\psi}{B_Z}$ ($B_Z$ is homeomorphic to $B \times \cZ/B \times\{\widecheck{z_0}\}$ and 
$B \times \cZ$ is closed in both $(\cX \# \cY) \times \cZ$ and $\cX \times \cY \times \cZ$).  
To see that 
$
\begin{cases}
\ \phi\\[-.15cm]
\ \psi
\end{cases}
$
 are mutually invervse pointed homotopy equivalences, define pointed homotopies 
 $
\begin{cases}
\ H:I\cX \ra \cX\\[-.15cm]
\ G:I\cY \ra \cY
\end{cases}
$
by 
$
\begin{cases}
\ (H \circ \restr{i_t}{X})(x) = x\\[-.15cm]
\ (G \circ \restr{i_t}{Y})(y) = y
\end{cases}
$
$\&$ 
$
\begin{cases}
\ (H \circ \restr{i_t}{[0,1]})(T)\\[-.15cm]
\ (G \circ \restr{i_t}{[0,1]})(T)
\end{cases}
$
$= \ds\max\left\{0, \frac{2T - t}{2 - t}\right\}.$
  \mH and \mG combine with $\id_Z$ to define a pointed homotopy on $\cX \times \cY \times \cZ$ which  
(i) induces a pointed homotopy on $\cX \# \cY \# \cZ$ between the identity and $\psi \circ \phi$ and 
(ii) induces a pointed homotopy on $(\cX \# \cY) \# \cZ$ between the identity and $\phi \circ \psi$
.]\\

Application: If \mX and \mY are nondegenerate $-$then in $\bHTOP_*$, 
$\Sigma (X \# Y) \approx$ 
$\Sigma X \# Y \approx$ 
$X \# \Sigma Y$.

[Note: \ Nondegeneracy is not actually necessary for the truth of this conclusion (cf. p. \pageref{3.30}).]\\

\begingroup%%----------------------------------->>
\fontsize{9pt}{11pt}\selectfont
Within the class of nondegenerate spaces, associativity of the smash product is natural, i.e., if $f:
%%----------------------------------------------------------------------------------------------37
X \ra X^\prime$, $g:Y \ra Y^\prime$, $h:Z \ra Z^\prime$ are pointed continuous functions, then the diagram
\[
\begin{tikzcd}[ sep=large]
{(X \#  Y) \# Z} \ar{d}[swap]{(f\# g)\# h} \ar{r}
&{X \# (Y \# Z)} \ar{d}{f \# (g\# h)}\\
{(X^\prime \# Y^\prime)\# Z^\prime} \ar{r}
&{X^\prime \# (Y^\prime \# Z^\prime)}
\end{tikzcd}
\]
commutes in $\bHTOP_*$.
\\ \indent
[Note: \ The horizontal arrows are the pointed homotopy equivalences figuring in the proof of Proposition 22.]\\
\endgroup %%------------------------------------<<

\begin{proposition} \ %23
Suppose that \mX and \mY are nondegenerate $-$then the projection $X \ov{\#} Y \ra X \# Y$ is a pointed homotopy equivalence.
\end{proposition}

[Consider the \cd
\begin{tikzcd}[ sep=large]
{\cX \ov{\#}\cY} \ar{d} \ar{r}&{\cX\#\cY} \ar{d}\\
{X \ov{\#} Y} \ar{r} &{X\# Y}
\end{tikzcd}
.  The upper horizontal arrow and the two vertical arrows are pointed homotopy equivalences, thus so is the lower horizontal arrow.]\\

Given pointed spaces 
$
\begin{cases}
\ X\\[-.15cm]
\ Y
\end{cases}
, \ 
$
the pointed mapping cone sequence associated with the inclusion $f:X \vee Y \ra X \times Y$. reads:
$X \vee Y \overset{f}{\ra} $
$X \times Y \ra $ 
$X \ov{\#} Y \ra $ 
$\Sigma (X \vee Y) \ra $
$\Sigma (X \times Y) \ra $ 
$\cdots$.\\

\textbf{\small LEMMA}  \  
The arrow $F:X \ov{\#} Y \ra \Sigma (X \vee Y)$ is nullhomotopic.

[There is a pointed injection $X \ov{\#} Y \ra \Gamma (X \times Y)$.  It is continuous (but not necessarily an embedding).  
Write $\Sigma (X \vee Y) = \Sigma X \vee \Sigma Y$ to realize 
$
F:
\begin{cases}
\ F[x,y_0,t] = [x,t] \in \Sigma X\\[-.15cm]
\ F[x_0,y,t] = [y,t] \in \Sigma Y
\end{cases}
$
$\&$ $F[x,y,1] = *$, the base point.  Put 
$
\begin{cases}
\ \ov{\Sigma}X = \Sigma X/\{[x,t]:x \in X, t \leq 1/2\}\\[-.15cm]
\ \un{\Sigma}Y = \Sigma Y/\{[y,t]:y \in Y, t \geq 1/2\}
\end{cases}
$
$-$then the arrows 
$
\begin{cases}
\ \Sigma X \ra \ov{\Sigma}X \\[-.15cm]
\ \Sigma Y \ra \un{\Sigma} Y
\end{cases}
$
are pointed homotopy equivalences, hence the same holds for their wedge: 
$\Sigma X \vee \Sigma Y \ra \ov{\Sigma}X \vee \un{\Sigma} Y$.  
The assignment 
$
[x,y,t] \ra 
\begin{cases}
\ [x,t] \ (t \geq 1/2)\\[-.15cm]
\ [y,t] \ (t \leq 1/2)
\end{cases}
$
defines a pointed continuous function $\Gamma (X \times Y) \ra \ov{\Sigma}X \vee \un{\Sigma} Y$.  
The composite 
$X \ov{\#} Y \ra \Gamma(X \times Y) \ra \ov{\Sigma}X \vee \un{\Sigma} Y$ is equal to the composite 
$X \ov{\#} Y \overset{F}{\ra} \Sigma X \vee \Sigma Y \ra \ov{\Sigma}X \vee \un{\Sigma} Y$.  
But the first composite is nullhomotopic.  
Therefore the second composite is nullhomotopic and this implies that $F \simeq 0$.]\\

\index{Puppe Formula} 
\textbf{\small PUPPE FORMULA}  \ \  
Suppose that \mX and \mY are nondegenerate $-$then in $\bHTOP_*$, 
$\Sigma (X \times Y) \approx \Sigma X \vee \Sigma Y \vee \Sigma(X \# Y)$.

%
%%----------------------------------------------------------------------------------------------38
[The third term of the pointed mapping cone sequence of 
$0 \ra X \ov{\#} Y \ra \Sigma(X \vee Y)$ is $\Sigma(X \vee Y) \vee \Sigma(X \ov{\#} Y)$, so from the lemma, 
$C_F \approx \Sigma(X \vee Y) \vee \Sigma(X \ov{\#} Y)$.  
Using now the notation of p. \pageref{3.31}, there is a commutative triangle
\begin{tikzcd}[ sep=large]
{X \ov{\#} Y} \ar{rd}[swap]{F} \ar{r}{j^\prime} &{C_j} \ar{d}\\
&{\Sigma (X \vee Y)}
\end{tikzcd}
in which the vertical arrow is a pointed homotopy equivalence, thus $C_{j^\prime} \approx C_F$ or still, 
$\Sigma (X \times Y) \approx $
$\Sigma (X \vee Y) \vee \Sigma(X \ov{\#} Y) \approx$ 
$\Sigma X \vee \Sigma Y \vee \Sigma(X \# Y)$ (cf. Proposition 23).]\\

Thanks to Proposition 22, this result can be iterated.  Let $X_1, \ldots, X_n$ be nondegenerate $-$then 
$\Sigma (X_1 \times \cdots \times X_n)$ has the same pointed homotopy type as 
$\ds\bigvee\limits_N \Sigma\left(\underset{i \in N}{\#} X_i\right)$, 
where $N$ runs over the nonempty subsets of $\{1, \ldots, n\}$.  
Example: 
$\ds\Sigma(\bS^{k_1} \times \cdots \times \bS^{k_n}) \approx \ds\bigvee\limits_N \bS^N$, $\bS^N$ a sphere of dimension 
$1 + \ds\sum\limits_{i \in N} k_i$.\\
\vspace{0.25cm}

\begingroup%%----------------------------------->>
\fontsize{9pt}{11pt}\selectfont
\index{Whitehead Products (example)}
\textbf{\small EXAMPLE (\un{Whitehead Products})} \quadx 
Let 
$
\begin{cases}
\ X\\
\ Y
\end{cases}
$
be nondegenerate $-$then for any pointed space \mE, there is a short exact sequence of groups
\[
0 \ra [\Sigma(X \# Y),E] \ra [\Sigma(X \times Y),E] \ra [\Sigma(X \vee Y),E] \ra 0.
\]
Here, composition is written additively even though the groups involved may not be abelian.  
This data generates a pairing 
$[\Sigma X,E] \times [\Sigma Y,E] \ra [\Sigma(X \# Y),E]$.  
Take 
$
\begin{cases}
\ \alpha \in [\Sigma X, E]\\
\ \beta \in [\Sigma Y, E]
\end{cases}
$
and use the embeddings 
$
\begin{cases}
\ [\Sigma X, E]\\
\ [\Sigma Y, E]
\end{cases}
\ra [\Sigma(X \times Y),E]
$
to form the commutator $\alpha + \beta - \alpha - \beta$ in $[\Sigma (X \times Y),E]$.  
Because it lies in the kernel of the homomorphism 
$[\Sigma (X \times Y),E] \ra [\Sigma (X \vee Y),E]$, by exactness there exists a unique element 
$[\alpha,\beta] \in [\Sigma (X \# Y),E]$ with image $\alpha + \beta - \alpha - \beta$.  $[\alpha,\beta]$ is called the 
\un{Whitehead product}
\index{Whitehead product} 
of $\alpha$, $\beta$.  $[\alpha,\beta]$ and $[\beta,\alpha]$ are connected by the relation $[\alpha,\beta] + [\beta,\alpha] \circ \Sigma\Tee = 0$, where $\Tee:X \# Y \ra Y \# X$ is the interchange.  
Of course, $[\alpha,0] = [0,\beta] = 0$.  
In general, $[\alpha,\beta]  = 0$ if \mE is an H-space (since then 
$[\Sigma(X \times Y),E]$ is abelian), hence, always $\ds\Sigma[\alpha,\beta] = 0$ (look at the arrow $E \ra \Omega\Sigma E$).  
There are left actions
\endgroup

\begingroup
\fontsize{9pt}{11pt}\selectfont
\[
\begin{cases}
\ [\Sigma X,E] \times [\Sigma (X \# Y),E] \ra [\Sigma (X \# Y),E]\\
\ [\Sigma Y,E] \times [\Sigma (X \# Y),E] \ra [\Sigma (X \# Y),E] 
\end{cases}
:
\begin{cases}
\ (\alpha,\xi) \ra \alpha \cdot \xi = \alpha + \xi - \alpha\\
\ (\beta,\xi) \ra \beta \cdot \xi = \beta + \xi - \beta
\end{cases}
\text{ (abuse of notation).}
\]
One has
$
\begin{cases}
\ [\alpha + \alpha^\prime,\beta] = \alpha\cdot[\alpha^\prime,\beta] + [\alpha,\beta]\\
\ [\alpha,\beta+\beta^\prime] = [\alpha,\beta] + \beta\cdot[\alpha,\beta^\prime]
\end{cases}
\hspace{-.2cm}. \ 
$
These relations simplify if the cogroup objects 
$
\begin{cases}
\ \Sigma X\\
\ \Sigma Y
\end{cases}
$ are commutative (as would be the case, e.g., when 
$
\begin{cases}
\ X = \Sigma X^\prime\\
\ Y = \Sigma Y^\prime
\end{cases}
$
for nondegenerate
$
\begin{cases}
\ X^\prime\\
\ Y^\prime
\end{cases}
\hspace{-.2cm}\big). \ 
$
Indeed, under this assumption, $[\Sigma (X \# Y),E]$ is abelian.  
Therefore the 
$
\begin{cases}
\ \alpha\cdot[\alpha^\prime,\beta] - [\alpha^\prime,\beta]\\
\ \beta\cdot[\alpha,\beta^\prime] - [\alpha,\beta^\prime]
\end{cases}
$
must vanish (``being commutative''), implying that
$
\begin{cases}
\ [\alpha + \alpha^\prime,\beta] = [\alpha,\beta] + [\alpha^\prime,\beta]\\
\ [\alpha,\beta+\beta^\prime] = [\alpha,\beta] + [\alpha,\beta^\prime]
\end{cases}
\hspace{-.2cm}. \ 
$
The Whitehead product also satisfies a form of the
%%----------------------------------------------------------------------------------------------39
Jacobi identity.  Precisely: Suppose given nondegenerate \mX, \mY, \mZ whose associated cogroup objects $\Sigma X$, 
$\Sigma Y$, $\Sigma Z$  are commutative 
$-$then
\[
[[\alpha,\beta],\gamma] + [[\beta,\gamma],\alpha] \circ \Sigma\sigma + [[\gamma,\alpha],\beta] \circ \Sigma\tau = 0
\]
in the group $[\Sigma(X \#Y \# Z),E]$, where 
$
\begin{cases}
\ \sigma:X \#Y \# Z \ra Y \# Z \# X\\
\ \tau:X \#Y \# Z \ra Z \# X \# Y\
\end{cases}
$
(cf. Proposition 22).  
The verification is a matter of manipulating commutator identities.]\\
\endgroup %%------------------------------------<<

\begingroup%%----------------------------------->>
\fontsize{9pt}{11pt}\selectfont
A 
\un{graded Lie algebra}
\index{graded Lie algebra} 
over a commutative ring \mR with unit is a graded \mR-module 
$L = \ds\bigoplus\limits_{n \geq 0} L_n$ together with bilinear pairings 
$[ \ ,\ ]:  L_n \times L_m \ra L_{n+m}$ such that 
$[x,y] = (-1)^{\abs{x}\abs{y}+1}[y,x]$ and 
\[
(-1)^{\abs{x}\abs{z}}[[x,y,],z] + (-1)^{\abs{y}\abs{x}}[[y,z],x] + (-1)^{\abs{z}\abs{y}}[[z,x],y] = 0.
\]
\mL is said to be 
\un{connected}
\index{connected graded Lie algebra} 
if $L_0 = 0$.  
Example: Let $A = \ds\bigoplus\limits_{n \geq 0} A_n$ be a graded $R$-algebra.  
For $x \in A_n$, $y \in A_m$, put 
$[x,y] = xy - (-1)^{\abs{x}\abs{y}}yx$ $-$then with this definition of the bracket, \mA is a graded Lie algebra over \mR.
\\ \indent
[Note: \ As usual, an absolute value sign stands for the degree of a homogeneous element in a graded $R$-module.]\\
\endgroup %%------------------------------------<<

\begingroup%%----------------------------------->>
\fontsize{9pt}{11pt}\selectfont
\textbf{\small EXAMPLE} \quadx 
Let \mX be a path connected topological space.  Given 
$
\begin{cases}
\ \alpha \in \pi_n(X)\\
\ \beta \in \pi_m(X)
\end{cases}
, \ 
$
the Whitehead product $[\alpha,\beta] \in \pi_{n+m-1}(X)$. One has 
$[\alpha,\beta] = (-1)^{nm + n + m}[\beta,\alpha]$.  
Moreover, if $\gamma \in \pi_r(X)$, then 
\[
(-1)^{nr + m}[[\alpha,\beta],\gamma] + (-1)^{mn + r}[[\beta,\gamma],\alpha] + (-1)^{rm + n}[[\gamma,\alpha],\beta] = 0.
\]
Assume now that \mX is simply connected.  
Consider the graded $\Z$-module 
$\pi_*(\Omega X) = \ds\bigoplus\limits_{n \geq 0} \pi_n(\Omega X)$.  Since $\pi_{n+1}(X) = \pi_n(\Omega X)$, 
the Whitehead product determines a bilinear pairing 
$[\ ,\ ]:\pi_n(\Omega X) \times \pi_m(\Omega X) \ra \pi_{n+m}(\Omega X)$ 
with respect to which $\pi_*(\Omega X)$ acquires the structure of a connected graded Lie algebra over $\Z$.\\
\endgroup %%------------------------------------<<

\begingroup%%----------------------------------->>
\fontsize{9pt}{11pt}\selectfont
\textbf{\small FACT} \   
Suppose that \mX is simply connected $-$then the Hurewicz homomorphism $\pi_*(\Omega X) \ra H_*(\Omega X)$ 
is  a morphism of graded Lie algebras, i.e., preserves the brackets.
\\ \indent
[Note: \ Recall that $H_*(\Omega X)$ is a graded $\Z$-algebra (Pontryagin product), hence can be regarded as a graded Lie algebra over $\Z$.]\\
\endgroup %%------------------------------------<<

\label{13.16} %dmc this seems out of place - but G lists it as p. 3-39 - still seems out of place

A pair $(X,A)$ is said to be 
\un{$n$-connected}
\index{n-connected (pair)} 
$(n \geq 1)$ if each path component of \mX meets \mA and $\pi_q(X,A,x_0) = 0$ $(1 \leq q \leq n)$ for all $x_0 \in A$ or,  equivalently, if every map 
$(\bD^q,\bS^{q-1}) \ra (X,A)$ is homotopic rel $\bS^{q-1}$ to a map $\bD^q \ra A$ $(0 \leq q \leq n)$.  
If $A$ is path connected, 
%%----------------------------------------------------------------------------------------------40
then $\forall \ x_0^\prime$, $x_0\pp \in A$, $\pi_n(X,A,x_0^\prime) \approx \pi_n(X,A,x_0\pp)$ $(n \geq 1)$.  
Examples: 
(1) $(\bD^{n+1},\bS^n)$ is $n$-connected; 
(2) $(\bB^{n+1},\bB^{n+1} - \{0\})$ is $n$-connected.

\label{4.41}
\label{4.54}
[Note: \ Take $A = \{x_0\}$ $-$then $\pi_q(X,\{x_0\},x_0) = \pi_q(X,x_0)$, so \mX is 
\un{$n$-connected}
\index{n-connected} 
$(n \geq 1)$ provided that \mX is path connected and $\pi_q(X) = 0$ $(1 \leq q \leq n)$.  
Example: $\bS^{n+1}$ is $n$-connected.]\\

\begingroup%%----------------------------------->>
\fontsize{9pt}{11pt}\selectfont
\label{4.56}
\label{5.30}
\textbf{\small EXAMPLE} \  
If \mX is $n$-connected and \mY is $m$-connected, then $X*Y$ is $((n+1) + (m+1))$-connected.
\\ \indent
[Note: \ If \mX is path connected and \mY is nonempty but arbitrary, then $X*Y$ is 1-connected.]\\
\endgroup %%------------------------------------<<

\begingroup%%----------------------------------->>
\fontsize{9pt}{11pt}\selectfont
\textbf{\small EXAMPLE} \  
Suppose that 
$
\begin{cases}
\ X\\
\ Y
\end{cases}
$
are nondegenerate and \mX is $n$-connected and \mY is $m$-connected $-$then $X \# Y$ is $(n + m + 1)$-connected.\\
\endgroup %%------------------------------------<<

\begingroup%%----------------------------------->>
\fontsize{9pt}{11pt}\selectfont
\textbf{\small FACT} \   
Let $f : \bS^n \ra A$ be a continuous function.  
Put $X = \bD^{n+1} \sqcup_f A$ $-$then $(X,A)$ is $n$-connected.\\
\endgroup %%------------------------------------<<

\begingroup%%----------------------------------->>
\fontsize{9pt}{11pt}\selectfont
\textbf{\small EXAMPLE} \  
The pair $(\bS^n \times \bS^m,\bS^n \vee \bS^m)$ is $n + m - 1$ connected.\\
\endgroup %%------------------------------------<<

\index{Theorem: Homotopy Excision Theorem}
\index{Homotopy Excision Theorem}
\textbf{\small HOMOTOPY EXCISION THEOREM} \ \ 
Suppose that 
$
\begin{cases}
\ X_1\\
\ X_2
\end{cases}
$
are subspaces of \mX with $X = \itr X_1 \cup \itr X_2$.  Assume: 
$
\begin{cases}
\ (X_1,X_1 \cap X_2)\\
\ (X_2,X_2 \cap X_1)
\end{cases}
$ 
\hspace{-.2cm} is \ 
$
\begin{cases}
\ \text{$n$-connected}\\
\ \text{$m$-connected}
\end{cases}
$
$-$then the arrow $\pi_q(X_1,X_1 \cap X_2) \ra \pi_q(X_1 \cup X_2,X_2)$ induced by the inclusion 
$(X_1,X_1 \cap X_2) \ra (X_1 \cup X_2,X_2)$ is bijective for $1 \leq q < n + m$ and surjective for $q = n + m$.

[This is dealt with at the end of the $\S$.]\\

\textbf{\small LEMMA}  \  
Let \mX be a strong deformation retract of \mY and let $A \subset X$ be a strong deformation retract of $B \subset Y$ $-$then 
$\forall \ n \geq 1$, $\pi_n(X,A) \approx \pi_n(Y,B)$.

[Use the exact sequence for a pair and the five lemma.]\\

\begin{proposition} \ %24
Let 
$
\begin{cases}
\ A\\[-.15cm]
\ B
\end{cases}
$
be closed subspaces of \mX with $X = A \cup B$.  
Put $C = A \cap B$.  
Assume: The inclusions
$
\begin{cases}
\ C \ra A\\[-.15cm]
\ C \ra B
\end{cases}
$
are cofibrations and 
$
\begin{cases}
\ (A,C)\\[-.15cm]
\ (B,C)
\end{cases}
$
is 
$
\begin{cases}
\ \text{$n$-connected}\\[-.15cm]
\ \text{$m$-connected}
\end{cases}
$
$-$then the arrow $\pi_q(A,C) \ra \pi_q(X,B)$ is bijective for $1 \leq q < n + m$ and surjective for $q = n + m$.
\end{proposition}

[Set $\ov{X} = i_0 A \cup IC \cup i_1B$, 
$
\begin{cases}
\ \ov{X}_1 = i_0A \cup IC\\[-.15cm]
\ \ov{X}_2 = IC \cup i_1B
\end{cases}
\hspace{-.2cm} :
$
$\ov{X}_1 \cap \ov{X}_2 = IC$ and 
$
\begin{cases}
\ \itr \ov{X}_1 \supset \ov{X} - i_1 B\\[-.15cm]
\ \itr \ov{X}_2 \supset \ov{X} - i_0 A
\end{cases}
$
$\implies$ $\ov{X} = \itr \ov{X}_1 \ \cup \  \itr \ov{X}_2$.  
From the lemma, 
$
\begin{cases}
\ \pi_q(A,C) \approx \pi_q(\ov{X}_1,IC)\\[-.15cm]
\ \pi_q(B,C) \approx \pi_q(\ov{X}_2,IC)
\end{cases}
$
$\implies$ 
$
\begin{cases}
\ (\ov{X}_1,IC)\\[-.15cm]
\ (\ov{X}_2,IC)
\end{cases}
$
is
$
\begin{cases}
\ \text{$n$-connected}\\[-.15cm]
\ \text{$m$-connected}
\end{cases}
\hspace{-.2cm} , \ 
$
thus the homotopy excision theorem is applicable to the triple $(\ov{X},\ov{X}_1,\ov{X}_2)$.  
Because the inclusions 
$
\begin{cases}
\ C \ra A\\[-.15cm]
\ C \ra B
\end{cases}
$
are cofibrations, $i_0A \cup IC$ is a strong deformation retract
%%----------------------------------------------------------------------------------------------41
of $IA$ and $IC \cup i_1B$ is a strong deformation retract of $IB$ (cf. p. \pageref{3.32}.  
Therefore $\ov{X}$ is a strong deformation retract of $IA \cup IB = IX$, so 
$\pi_q(\ov{X},\ov{X}_2) \approx $ 
$\pi_q(IX,IB) \approx $ 
$\pi_q(X,B) $ .]\\

\textbf{\small LEMMA}  \  
Let $f:(X,A) \ra (Y,B)$ be a homotopy equivalence in $\bTOP^2$ $-$then $\forall \ x_0 \in A$ and any $q \geq 1$, the induced map $f_*:\pi_q(X,A,x_0) \ra \pi_q(Y,B,f(x_0))$ is bijective.\\

\begin{proposition} \ %25
Let \mA be a nonempty closed subspace of \mX.  
Assume: The inclusion $A \ra X$ is a  cofibration and \mA is $n$-connected, 
$(X,A)$ is $m$-connected 
$-$then the arrow $\pi_q(X,A) \ra \pi_q(X/A,*)$ is bijective for $1 \leq q \leq n + m$ and surjective for 
$q = n + m + 1$.  
\end{proposition}

[Denote by $C_i$ the unpointed mapping cone of the inclusion $i:A \ra X$.  There are closed cofibrations 
$
\begin{cases}
\ \Gamma A \ra C_i\\[-.15cm]
\ X \ra C_i
\end{cases}
$
and $C_i = \Gamma A \cup X$, with $\Gamma A \cap X = A$.  Since the pair $(\Gamma A,A)$ 
is $(n+1)$-connected, it follows from Proposition 24 that the arrow 
$\pi_q(X,A) \ra \pi_q(C_i,\Gamma A)$ 
is bijective for 
$1 \leq q \leq n + m$ and surjective for $q = n + m + 1$.  But $\Gamma A$ is contractible, hence the projection 
$(C_i,\Gamma A) \ra (C_i/\Gamma A,*)$ is a homotopy equivalence in $\bTOP^2$ (cf. Proposition 14).  Taking into account the lemma, it remains only to observe that $X/A$ can be identified with $C_i/\Gamma A$.]\\

\index{Theorem: Freudenthal Suspension Theorem}
\index{Freudenthal Suspension Theorem}
\textbf{\small FREUDENTHAL SUSPENSION THEOREM} \ \ 
Suppose that \mX is nondegenerate and $n$-connected $-$then the suspension homomorphism 
$\pi_q(X) \ra \pi_{q+1}(\Sigma X)$ is bijective for $0 \leq q \leq 2n$ and surjective for $q = 2n+1$.

[Take \mX wellpointed with a closed base point and, for the moment, work with its unpointed suspension $\Sigma X$.  
Using the notation of p. \pageref{3.33}, write $\Sigma X = \Gamma^- X \cup \Gamma^+ X$ $-$then $\forall \ q$, 
$\pi_q(X) \approx \pi_q(\Gamma^- X \cap \Gamma^+ X)$ 
$\approx \pi_{q+1}(\Gamma^- X,\Gamma^- X \cap \Gamma^+ X)$.  
On the other hand, Proposition 25 implies that the arrow 
$\pi_{q+1}(\Gamma^- X,\Gamma^- X \cap \Gamma^+ X) \ra \pi_{q+1}(\Sigma X)$ is a bijection for $1 \leq q + 1 \leq 2n + 1$ and a surjection for $q + 1 = 2n + 2$.  
Moreover, \mX is wellpointed, therefore its pointed and unpointed suspensions have the same homotopy type.]

[Note: \ This result is true if \mX is merely path connected, i.e., $n = 0$ is admissible (inspect the proof of Proposition 25),]\\

Application:  Suppose that $n \geq 1$ $-$then 
(i)  $\pi_q(\bS^n) = 0$ $(0 \leq q < n)$;
(ii)  $\pi_q(\bS^n) \approx \pi_{q+1}(\bS^{n+1})$ $(0 \leq q \leq 2n - 2)$; 
(iii)  $\pi_n(\bS^n) \approx \Z$.

[As regards the last point, note that in the sequence 
$\pi_1(\bS^1) \ra $ 
$\pi_2(\bS^2) \ra $ 
$\pi_3(\bS^3) \ra $ 
$\cdots$, 
the first homomorphism is an epimorphism, the others are isomorphisms, and 
$\pi_1(\bS^1) \approx \Z$, 
$\pi_2(\bS^2) \approx \Z$ (a piece of the exact sequence associated with the Hopf map $\bS^3 \ra \bS^2$ is 
$\pi_2(\bS^3) \ra $ 
$\pi_2(\bS^2) \ra $
$\pi_1(\bS^1) \ra \pi_1(\bS^3)$).]\\

%%----------------------------------------------------------------------------------------------42
The infinite cyclic group $\pi_n(\bS^n)$ is generated by $[\iota_n]$, $\iota_n$ the identity $\bS^n \ra \bS^n$.  
Form 
the Whitehead product $[\iota_n,\iota_n] \in \pi_{2n-1}(\bS^n)$ $-$then the kernel of the suspension homomorphism 
$\pi_{2n-1}(\bS^n) \ra \pi_{2n}(\bS^{n+1})$ is generated by $[\iota_n,\iota_n]$ 
(Whitehead\footnote[2]{\textit{Elements of Homotopy Theory}, Springer Verlag (1978), 549.}
).\\

\begingroup%%----------------------------------->>
\fontsize{9pt}{11pt}\selectfont
The proof of the homotopy excision theorem is elementary but complicated.  This is the downside.  The upside is that the highpowered approaches are cluttered with unnecessary assumptions, hence do not go as far.\\
\endgroup %%------------------------------------<<

\index{Theorem: Open Homotopy Excision Theorem}
\index{Open Homotopy Excision Theorem}
\textbf{\small OPEN HOMOTOPY EXCISION THEOREM} \  \ 
\begingroup%%----------------------------------->>
\fontsize{9pt}{11pt}\selectfont
Suppose that 
$
\begin{cases}
\ X_1\\
\ X_2
\end{cases}
$
are open subspaces of \mX with $X = X_1 \cup X_2$.  Assume: 
$
\begin{cases}
\ (X_1,X_1 \cap X_2)\\
\ (X_2,X_2 \cap X_1)
\end{cases}
$
is 
$
\begin{cases}
\ \text{$n$-connected}\\
\ \text{$m$-connected}
\end{cases}
$
$-$then the arrow $\pi_q(X_1,X_1 \cap X_2) \ra$ $\pi_q(X_1 \cup X_2, X_2)$ induced by the inclusion 
$(X_1,X_1 \cap X_2) \ra (X_1 \cup X_2,X_2)$ is bijective for $1 \leq q < n + m$ and surjective for $q = n + m$.
\\ \indent
[Note: 
Goodwillie\footnote[3]{\textit{Memoirs Amer. Math. Soc.} \textbf{431} (1990), 1-317.}
has extended the open homotopy excision theorem to ``(N+1)-ads''.]\\
\endgroup %%------------------------------------<<

\begingroup%%----------------------------------->>
\fontsize{9pt}{11pt}\selectfont
Admit the open homotopy excision theorem.\\
\endgroup %%------------------------------------<<

\index{Theorem: CW Homotopy Excision Theorem}
\index{CW Homotopy Excision Theorem}
\textbf{\small CW HOMOTOPY EXCISION THEOREM} \  \ 
\begingroup%%----------------------------------->>
\fontsize{9pt}{11pt}\selectfont
Suppose that 
$
\begin{cases}
\ K_1\\
\ K_2
\end{cases}
$
are subcomplexes of a CW complex \mK with $K = K_1 \cup K_2$.  Assume: 
$
\begin{cases}
\ (K_1,K_1 \cap K_2)\\
\ (K_2,K_2 \cap K_1)
\end{cases}
$
is 
$
\begin{cases}
\ \text{$n$-connected}\\
\ \text{$m$-connected}
\end{cases}
$
$-$then the arrow 
$\pi_q(K_1,K_1 \cap K_2) \ra$ $\pi_q(K_1 \cup K_2, K_2)$ induced by the inclusion 
$(K_1,K_1 \cap K_2) \ra (K_1 \cup K_2,K_2)$ is bijective for $1 \leq q < n + m$ and surjective for $q = n + m$.

[Fix a neighborhood 
$
\begin{cases}
\ U\\
\ V
\end{cases}
$
of $K_1 \cap K_2$ in 
$
\begin{cases}
\ K_1\\
\ K_2
\end{cases}
$
such that $K_1 \cap K_2$ is a strong deformation retract of 
$
\begin{cases}
\ U\\
\ V
\end{cases}
$
and put 
$
\begin{cases}
\ K_1^\prime = K_1 \cup V\\
\ K_2^\prime = K_2 \cup U
\end{cases}
\hspace{-.2cm}. \ 
$
Write
$
\begin{cases}
\ U = O \cap K_1\\
\ V = P \cap K_2
\end{cases}
\hspace{-.2cm}, \ 
$
where 
$
\begin{cases}
\ O\\
\ P
\end{cases}
$
are open in \mK $-$then 
$
\begin{cases}
\ K_1^\prime = P \cup (K - K_2)\\
\ K_2^\prime = O \cup (K - K_1)
\end{cases}
\hspace{-.2cm}, \ 
$
hence
$
\begin{cases}
\ K_1^\prime\\
\ K_2^\prime
\end{cases}
$
are open in \mK and $K = K_1^\prime \cup K_2^\prime$.  Since 
$
\begin{cases}
\ K_1 \ \& \ V\\
\ K_2 \ \& \ U
\end{cases}
$
are closed in 
$
\begin{cases}
\ K_1^\prime\\
\ K_2^\prime
\end{cases}
\hspace{-.2cm}, \ 
$
the homotopy deforming 
$
\begin{cases}
\ V\\
\ U
\end{cases}
$
into $K_1 \cap K_2$ can be extended to all of 
$
\begin{cases}
\ K_1^\prime\\
\ K_2^\prime
\end{cases}
$
in the obvious way, so 
$
\begin{cases}
\ K_1\\
\ K_2
\end{cases}
$
is a strong deformation retract of 
$
\begin{cases}
\ K_1^\prime\\
\ K_2^\prime
\end{cases}
\hspace{-.2cm}. \ 
$
On the other hand, $K_1^\prime \cap K_2^\prime = U \cup V$ and 
$
\begin{cases}
\ U\\
\ V
\end{cases}
$
is closed in $U \cup V$, thus the union of the deforming homotopies is continuous and $K_1 \cap K_2$ is a strong deformation retract of $K_1^\prime \cap K_2^\prime$.  Therefore 
$
\begin{cases}
\ (K_1^\prime,K_1^\prime \cap K_2^\prime)\\
\ (K_2^\prime,K_2^\prime \cap K_1^\prime)
\end{cases}
$
is
$
\begin{cases}
\ \text{$n$-connected}\\
\ \text{$m$-connected}
\end{cases}
$
and the open homotopy excision theorem is applicable to the triple $(K,K_1^\prime,K_2^\prime)$.  
Consider the commutative triangle.  
%%----------------------------------------------------------------------------------------------43
\[
\begin{tikzcd}[ sep=small]
{\pi_q(K_1,K_1 \cap K_2)}\ar{rddd} \ar{rr} &&{{\pi_q(K_1^\prime,K_1^\prime \cap K_2^\prime)}}\ar{lddd}\\
\\
\\
&{\pi_q(K_1 \cup K_2,K_2)} 
\end{tikzcd}
.]
\]
\\
\endgroup %%------------------------------------<<

\begingroup%%----------------------------------->>
\fontsize{9pt}{11pt}\selectfont
The CW homotopy excision theorem implies the homotopy excision theorem.  
For choose a CW resolution 
$L \ra X_1 \cap X_2$.  
There exist: 
(1) A CW complex $K_1 \supset L$ and a CW resolution $f_1:K_1 \ra X_1$ such that the square
\begin{tikzcd}[ sep=large]
{K_1} \ar{r}  &{X_1}\\
{L} \ar{u} \ar{r} &{X_1 \cap X_2} \ar{u}
\end{tikzcd}
commutes; 
(2) A CW complex $K_2 \supset L$ and a CW resolution $f_2:K_2 \ra X_2$ such that the square
\begin{tikzcd}[ sep=large]
{K_2} \ar{r}  &{X_2}\\
{L} \ar{u} \ar{r} &{X_2 \cap X_1} \ar{u}
\end{tikzcd}
commutes.  
Note that 
$
\begin{cases}
\ (K_1,L)\\
\ (K_2,L)
\end{cases}
$
is 
$
\begin{cases}
\ \text{$n$-connected}\\
\ \text{$m$-connected}
\end{cases}
. \ 
$
Define a CW complex \mK by the pushout square 
\begin{tikzcd}[ sep=large]
{L} \ar{d} \ar{r} &{K_2} \ar{d}\\
{K_1} \ar{r} &{K}
\end{tikzcd}
: $K = K_1 \cup K_2$ $\&$ $L = K_1 \cap K_2$ $-$then there is an arrow $f:K \ra X$ determined by 
$
\begin{cases}
\ f_1\\
\ f_2
\end{cases}
\hspace{-.25cm}
,\ 
$
viz. 
$
\begin{cases}
\ \restr{f}{K_1} = f_1\\
\ \restr{f}{K_2} = f_2
\end{cases}
. \ 
$
\\[0.25cm]
\endgroup %%------------------------------------<<

\textbf{\small LEMMA}  \  
\begingroup%%----------------------------------->>
\fontsize{9pt}{11pt}\selectfont
$f$ is a weak homotopy equivalence.

[Set $\ov{K} = i_0 K_1 \cup I L \cup i_1 K_2$: 
$
\begin{cases}
\ U_1 = \ov{K} - i_1K_2\\
\ U_2 = \ov{K} - i_0K_1
\end{cases}
$
$-$then 
$
\begin{cases}
\ U_1\\
\ U_2
\end{cases}
$
are open in $\ov{K}$ and $\ov{K} = U_1 \cup U_2$.  Let $\ov{p}:\ov{K} \ra K$ be the restriction of the projection 
$p:IK \ra K$ and denote by $\ov{f}$ the composite $f \circ \ov{p}$:
$
\begin{cases}
\ \ov{f}(U_1) \subset X_1\\
\ \ov{f}(U_2) \subset X_2
\end{cases}
$
and 
$
\begin{cases}
\ \restr{\ov{f}}{U_1}\\
\ \restr{\ov{f}}{U_2}
\end{cases}
$
$\&$ $\restr{\ov{f}}{U_1 \cap U_2}$ are weak homotopy equivalences.  
But by assumption $X = \itr X_1 \cup \itr X_2$.  
Therefore $\ov{f}$ is a weak homotopy equivalence (cf. p. \pageref{3.34}).  
The inclusions 
$
\begin{cases}
\ K_1 \ra K\\
\ K_2 \ra K
\end{cases}
$
are closed cofibrations (cf. p. \pageref{3.35}), hence $\ov{K}$ is a strong deformation retract of $IK$.  
Consequently, 
$\ov{p}$ is a homotopy equivalence, so $f$ is a weak homotopy equivalence.]\\
\endgroup %%------------------------------------<<

\begingroup%%----------------------------------->>
\fontsize{9pt}{11pt}\selectfont
The CW homotopy excision theorem is applicable to the triple $(K,K_1,K_2)$.  Examination of the commutative square
\[
\begin{tikzcd}[ sep=large]
{\pi_q(K_1,K_1 \cap K_2)} \ar{d}\ar{r}&{\pi_q(K_1 \cup K_2,K_2)} \ar{d}\\
{\pi_q(X_1,X_1 \cap X_2)} \ar{r} &{\pi_q(X_1 \cup X_2,X_2)}
\end{tikzcd}
\]	
thus justifies the claim.  
Accordingly, it is the open homotopy excision theorem which is the heart of the matter.
\\[-.2cm]
\endgroup %%------------------------------------<<

%%----------------------------------------------------------------------------------------------44
\begingroup%%----------------------------------->>
\fontsize{9pt}{11pt}\selectfont
Given a $p$-dimensional cube \mC in $\R^q$ $(q \geq 1, 0 \leq p \leq q)$, 
denote by $\sk_d \ C$ its $d$-dimensional skeleton, i.e., the set of its $d$-dimensional faces.  
Put $\dot{C} = \bigcup \sk_{p-1} C$ $-$then the inclusion 
$\dot{C} \ra C$ is a closed cofibration.  
Analytically, \mC is specified by a point $(c_1 \ldots c_q) \in \R^q$, a positive number 
$\delta$, and a subset \mP of $\{1, \ldots, q\}$ of cardinality $p$ : \mC is the set of $x \in \R^q$ such that 
$c_i \leq x_i \leq c_i + \delta$ $(i \in P)$ $\&$ $x_i = c_i$ $(i \notin P)$.  
Here, if $P = \emptyset$, then 
$C = \{(c_1, \ldots, c_q)\}$.  For $1 \leq d \leq q$, let 
$
\begin{cases}
\ \ds{ K_d(C) = \{x\in C: x_i < c_i + \ds\frac{\delta}{2} \text{ for at least $d$ indices $i \in P$}\}}\\
\\
\  \ds{L_d(C) = \{x\in C: x_i > c_i + \ds\frac{\delta}{2} \text{ for at least $d$ indices $i \in P$}\}}
\end{cases}
.
$
When $d > p$, it is understood that 
$
\begin{cases}
\ K_d(C) = \emptyset\\
\ L_d(C) = \emptyset
\end{cases}
.
$
\\[0.25cm]
\endgroup %%------------------------------------<<

\index{Compression Lemma}
\textbf{\small COMPRESSION LEMMA} \  \ 
\begingroup%%----------------------------------->>
\fontsize{9pt}{11pt}\selectfont
Fix a $p$-dimensional cube \mC in $\R^q$ $(q \geq 1, 1 \leq p \leq q)$, a positive integer $d \leq p$, and a pair $(X,A)$.  
Suppose that $f:C \ra X$ is a continuous function such that $\forall \ D \in \sk_{p-1}\ C$, 
$f^{-1}(A) \cap D \subset K_d(D)$ $(L_d(D))$ $-$then there exists a continuous function $g:C \ra X$ with 
$f \simeq g \ \rel\dot{C}$ and $g^{-1}(A) \subset K_d(C)$ $(L_d(C))$.
\\ \indent
[Take $p = q$, $C = [0,1]^q$, and put 
$x_0 = (1/4,\ldots, 1/4)$.  Given an $x \in [0,1]^q$, let $\ell((x_0,x)$ be the ray that starts at $x_0$ and passes through $x$.  
Denote by $P(x)$ the intersection of $\ell(x_0,x)$ with the frontier of 
$[0,1/2]^q$, $Q(x)$ the intersection of $\ell(x_0,x)$ with the frontier of $[0,1]^q$.  
Let $\phi:[0,1]^q \ra [0,1]^q$ be the continuous function that sends the line segment joining 
$P(x)$ and $Q(x)$ 
to the point $Q(x)$ 
and maps the line segment joining $x_0$ and $P(x)$ linearly onto the line segment joining $x_0$ and $Q(x)$.  
Note that $\phi \simeq \id_{[0,1]^q} \  \rel \fr[0,1]^q$.  
Now set \  $g = f \circ \phi$.  \ 
Assume: $x \in g^{-1}(A)$.  \ 
Case 1:  $x_i < 1/2$ \ $(\forall \ i)$ $\implies$ $x \in K_q([0,1]^q) \subset K_d([0,1]^q)$.  \ 
Case 2: $x_i \geq 1/2$ $(\exists \ i)$ $\implies$ $\phi(x) \in \fr([0,1]^q$ $\implies$ $\phi(x) \in D$ 
$(\exists \ D \in \sk_{p-1}\ [0,1]^q)$ $\implies$ $\phi(x) \in K_d(D)$ $\implies$ 
$1/2 > \phi(x)_i = 1/4 + t(x_i - 1/4)$ for at least $d$ indices 
$i$ $\implies$ $1/2 > \phi(x)_i \geq x_i$ $(t \geq 1)$ for at least $d$ indices $i$ $\implies$ $x \in K_d([0,1]^q)$.]
\\ \indent
[Note: \ The parenthetical assertion is analogous.]\\
\endgroup %%------------------------------------<<

\begingroup%%----------------------------------->>
\fontsize{9pt}{11pt}\selectfont
Notation: Put
$I^q = [0,1]^q$, 
$\dot{I}^q = \fr  [0,1]^q$, 
$I_0^{q-1} = I^{q-1} \times \{0\}$, $(q > 1)$ $\&$ 
$I_0^0 = \{0\}$ $(q = 1)$, 
$J^{q-1} = \dot{I}^{q-1} \times I \cup I^{q-1} \times \{1\}$ $(q > 1)$, $\&$ 
$J^0 = \{1\}$ $(q = 1)$, so 
$\dot{I}^q  = I_0^{q-1} \cup J^{q-1}$ and 
$\dot{I}_0^{q-1} = I_0^{q-1} \cap J^{q-1}$ $-$then for any pointed pair $(X,A,x_0)$, 
$\pi_q(X,A,x_0) = [I^q,\dot{I}^q,J^{q-1};X,A,x_0]$.
\\ \indent
[Note: A continuous function $f:(I^q,\dot{I}^q,J^{q-1}) \ra (X,A,x_0)$ represents 0 in $\pi_q(X,A,x_0)$ iff there 
exists a continuous function $g:I^q \ra A$ such that $f \simeq g \ \rel \dot{I}^q$.]\\
\endgroup %%------------------------------------<<

\begingroup%%----------------------------------->>
\fontsize{9pt}{11pt}\selectfont
There are two steps in the proof of the open homotoy excision theorem: 
(1) Surjectivity in the range $1 \leq q \leq n + m$;
(2) Injectivity in the range $1 \leq q < n + m$.  
The argument in either situation is founded on the same iterative principle.
\\ \indent
Starting with surjectivity, let $\alpha \in \pi_q(X_1 \cup X_2, X_2,x_0)$, $x_0 \in X_1 \cap X_2$ the ambient base point.  
Represent $\alpha$ by an $f:(I^q,\dot{I}^q,J^{q-1}) \ra (X_1 \cup X_2, X_2,x_0)$.  It will be shown that $\exists$ 
$F \in \alpha$: 
$\pr(F^{-1}(X - X_1)) \cap \pr(F^{-1}(X - X_2)) = \emptyset$, 
$\pr:I^q \ra I^{q-1}$ the projection.  Granted this, choose a continuous function 
$\phi:I^{q-1} \ra [0,1]$ which is 1 on 
$\pr(F^{-1}(X - X_1))$ and 0 on $\dot{I}^{q-1} \cup \pr(F^{-1}(X - X_2))$.  Define
%%----------------------------------------------------------------------------------------------45
$\Phi:I^q \ra I^q$ by $\Phi(x_1, \ldots, x_q) = (x_1, \ldots x_{q-1},t + (1-t)x_q)$, where 
$t = \phi(x_1, \ldots, x_{q-1})$, and put $g = F \circ \Phi$ $-$then 
$g:(I^q,\dot{I}^q,J^{q-1}) \ra (X_1, X_1 \cap X_2,x_0)$ is a continuous function whose class 
$\beta \in \pi_q(X_1, X_1 \cap X_2,x_0)$ is sent to $\alpha$ under the inclusion.
\\ \indent
There remains the task of producing \mF.  
Since 
$\{f^{-1}(X_1), f^{-1}(X_2)\}$ is an open covering of $I^q$, one can subdivide $I^q$ into a collection $\sC$ of $q$-dimensional 
cubes \mC such that either $f(C) \subset X_1$ or $f(C) \subset X_2$.  
Enumerate the elements in $\sk_d \hsx C$ $(C \in \sC, d = 0, 1, \ldots, q)$: $\sD = \{D\}$.  
In $\sD$, distinguish two subcollections 
$
\begin{cases}
\ \{D_k:k = 1, \ldots, r\}: f(D_k) \subset X_2\\
\ \{D_l:l = 1, \ldots, s\}: f(D_l) \subset X_1
\end{cases}
$
but 
$
\begin{cases}
\ f(D_k) \not\subset X_1\\
\ f(D_l) \not\subset X_2
\end{cases}
\hspace{-.2cm}, \ 
$
arranging the indexing so that $\dim D_j \leq \dim D_{j+1}$.\\
\indent\indent ($\mu$) \ There exist continuous functions 
$\mu_0 = f$, 
$\mu_k:I^q \ra X$ $(k = 1, \ldots, r)$ such that $\forall \ k$: $\mu_k \simeq \mu_0$ (as a map of triples), 
$\mu_k^{-1}(X_2 - X_1 \cap X_2) \cap D_j$ $\subset$ $K_{n+1}(D_j)$ $(j \leq k)$, and $\forall \ D \in \sD$ :
$
\begin{cases}
\ \mu_0(D) \subset X_1\\
\ \mu_0(D) \subset X_2
\end{cases}
$
$\implies$
$
\begin{cases}
\ \mu_k(D) \subset X_1\\
\ \mu_k(D) \subset X_2
\end{cases}
$
or $\mu_0(D) \subset X_1 \cap X_2$ $\implies$ $\mu_k(D) \subset X_1 \cap X_2$.  
This is seen via induction on $k$, $\mu_0 = f$ being the initial step.  
Assume that $\mu_{k-1}$ has been constructed.
\\ \indent
Claim: $\exists$ a homotopy $h_k:ID_k \ra X_2 \ \rel\dot{D}_k$ such that 
$h_k \circ i_0 = \restr{\mu_{k-1}}{D_k}$ and $(h_k \circ i_1)^{-1}(X_2 - X_1 \cap X_2)$ $\subset$ $K_{n+1}(D_k)$.
\\ \indent
[Case 1: \ $\dim D_k = 0$.  Here, $K_{n+1}(D_k) = \emptyset$ and the point $\mu_{k-1}(D_k) \in X_2$ can be joined by a path in $X_2$ to some point of $X_1 \cap X_2$.  
Case 2: \ $0 < \dim D_k < n + 1$.  Here, $K_{n+1}(D_k) = \emptyset$ and the induction hypothesis forces the containment 
$\mu_{k-1}(\dot{D}_k) \subset X_1 \cap X_2$, hence $\restr{\mu_{k-1}}{D_k}$ represents an element of 
$\pi_{d_k}(X_2,X_1 \cap X_2)  = 0$ $(d_k = \dim D_k)$.
Case 3: \ $\dim D_k \geq n + 1$.  Apply the compression lemma.]
\\ \indent
Extend $h_k$ to a homotopy 
$H_k:I^q \times I \ra X$ of 
$\mu_{k-1} \ \rel \cup \{D:f(D) \subset X_1\} \cup \ds\bigcup\limits_{j=1}^{k-1} D_j$ such that
$\ds\bigcup\limits_{j=k+1}^r H_k(ID_j) \subset X_2$.  Complete the induction by taking $\mu_k = H_k \circ i_1$.
\\
\indent\indent ($\nu$) \  There exist continuous functions $\nu_0 = \mu_r$, $\nu_l: I^q \ra X$ $(l = 1, \ldots s)$ such that 
$\forall \ l$: $\nu_l \simeq \nu_0 \  \rel \cup \{D:f(D) \subset X_2\}$, 
$\nu_l^{-1}(X_1 - X_1 \cap X_2) \cap D_j$ $\subset$ $L_{m+1}(D_j)$ $(j \leq l)$, and $\forall \ D \in \sD$: 
$
\begin{cases}
\ \nu_0(D) \subset X_1\\
\ \nu_0(D) \subset X_2
\end{cases}
$
$\implies$
$
\begin{cases}
\ \nu_l(D) \subset X_1\\
\ \nu_l(D) \subset X_2
\end{cases}
$
or $\nu_0(D) \subset X_1 \cap X_2$ $\implies$ $\nu_l(D) \subset X_1 \cap X_2$.  
As above, this is seen via induction on $l$, $\nu_0 = \mu_r$ being the initial step.  
Observe that $\cup\{D:f(D) \subset X_2\} \supset \dot{I}^q \supset J^{q-1}$.
\\ \indent
Definition: $F = \nu_s$ $(\implies F \in \alpha)$.  If $\pr(F^{-1}(X - X_1)) \cap \pr(F^{-1}(X - X_2))$ were nonempty, then there would exist an $x \in I^{q-1}$ and a cube $D \subset I^{q-1}$: 
$
\begin{cases}
\ x \in K_n(D)\\
\ x \in L_m(D)
\end{cases}
, \ 
$
an impossibility since $q - 1 < n + m$.
\\ \indent
Turning to injectivity, let 
$f,g:(I^q,\dot{I}^q,J^{q-1}) \ra (X_1, X_1 \cap X_2, x_0)$ 
be continuous functions such that 
$u \circ f \simeq u \circ g$ as maps of triples 
$u:(X_1,X_1\cap X_2,x_0) \ra (X_1 \cup X_2, X_2 ,x_0)$ the inclusion.  
Fix a homotopy 
$h:(I^q,\dot{I}^q,J^{q-1}) \times I \ra (X_1 \cup X_2,X_2,x_0)$:
$
\begin{cases}
\ h \circ i_0 = u \circ f\\
\ h \circ i_1 = u \circ g
\end{cases}
. \ 
$
Using the techniques employed in the proof of surjectivity, one can replace $h$ by another homotopy \mH such that 
$\pr \times \id_I(H^{-1}(X - X_1)) \cap \pr \times \id_I(H^{-1}(X - X_2)) = \emptyset$.  
It is this extra dimension that accounts for the restriction $q < n + m$.  Choose a continuous function 
$\phi:I^{q-1} \times I \ra [0,1]$ which is 1 on $\pr \times \id_I(H^{-1}(X - X_1))$ and 0 on 
$(\dot{I}^{q-1} \times I) \cup (I^{q-1} \times \dot{I}) \cup \pr \times \id_I(H^{-1}(X - X_2))$.  
Define $\Phi:I^q \times I \ra I^q \times I$ by 
$\Phi(x_1, \ldots, x_q,x_{q+1}) = $
%%----------------------------------------------------------------------------------------------46
$(x_1, \ldots, x_{q-1},t+(1-t)x_q,x_{q+1})$, where $t = \phi(x_1,\ldots,x_{q-1},x_{q+1})$ $-$then the composite 
$H \circ \Phi$ is a homotopy between $f$ and $g:H \circ \Phi(\dot{I}^q \times I) \subset X_1 \cap X_2$ $\&$ 
$H \circ \Phi(J^{q-1} \times I) = \{x_0\}$.\\
\endgroup %%------------------------------------<<

\begingroup%%----------------------------------->>
\fontsize{9pt}{11pt}\selectfont
Given a pair $(X,A)$, let $\pi_0(X,A)$ be the quotient $\pi_0(X)/\sim$, where $\sim$ means that the path components of \mX which meet \mA are identified.  With this agreement, $\pi_0(X,A)$ is a pointed set.  If 
$f:(X,A) \ra (Y,B)$ is a map of pairs, then 
$f_*:\pi_0(X,A) \ra \pi_0(Y,B)$  is a morphism of pointed sets and the sequence 
$* \ra \pi_0(X,A) \ra \pi_0(Y,B)$ is exact in $\bSET_*$ iff 
$(f_*)^{-1} \im(\pi_0(B) \ra \pi_0(Y)) = \im(\pi_0(A) \ra \pi_0(X))$.\\
\endgroup %%------------------------------------<<

\begingroup%%----------------------------------->>
\fontsize{9pt}{11pt}\selectfont
\textbf{\small LEMMA}  \  
Let $f:(X,A) \ra (Y,B)$ be a continuous function.  Fix $q \geq 0$ $-$then $\forall \ x_0 \in A$, 
$f_*:\pi_q(X,A,x_0) \ra \pi_q(Y,B,f(x_0))$ is injective and 
$f_*:\pi_{q+1}(X,A,x_0) \ra \pi_{q+1}(Y,B,f(x_0))$ is surjective iff in any diagram
\begin{tikzcd}[ sep=large]
{(X,A)} \ar{r}{f}&{(Y,B)}\\
{(J^q,\dot{I}_0^q)} \ar{u}{\phi} \ar{r} &{(I^{q+1},I_0^q)} \ar{u}[swap]{\psi}
\end{tikzcd}
, where $f \circ \phi \simeq \psi$ on $J^q$ by 
$h:(J^q,\dot{I}_0^q) \times I \ra (Y,B)$, there exists a $\Phi:(I^{q+1},I_0^q) \ra (X,A)$ such that 
$\restr{\Phi}{(J^q,\dot{I}_0^q)} = \phi$ and an 
$H:(I^{q+1},I_0^q) \times I \ra (Y,B)$ such that 
$\restr{H}{(J^q,\dot{I}_0^q) \times I} = h$ and $f \circ \Phi \simeq \psi$ on $I^{q+1}$ by \mH.
\\ \indent
[Note: When $q = 0$, replace injectivity by the statement ``$* \ra \pi_0(X,A) \ra \pi_0(Y,B)$'' is exact.  
Observe that $f \circ \phi = \psi$ on $J^q$ is permissible ($h$ = constant homotopy) and implies by specialization the direct assertion.  In addtion, if $\Phi$ $\&$ \mH exist in this case, then $\Phi$ $\&$ \mH exist in general.  Thus the point is to show that the direct assertion entails the existence of $\Phi$ $\&$ \mH under the assumption that $f \circ \phi = \psi$ on $J^q$.]\\
\endgroup %%------------------------------------<<

\label{4.59}
\begingroup%%----------------------------------->>
\fontsize{9pt}{11pt}\selectfont
\textbf{\small FACT} \   
Suppose that 
$
\begin{cases}
\ X_1\\
\ X_2
\end{cases}
$
$\&$
$
\begin{cases}
\ Y_1\\
\ Y_2
\end{cases}
$
are open subspaces of 
$
\begin{cases}
\ X\\
\ Y
\end{cases}
$
with 
$
\begin{cases}
\ X = X_1 \cup X_2\\
\ Y = Y_1 \cup Y_2
\end{cases}
. \ 
$
Let $f:X \ra Y$ be a continuous function such that 
$
\begin{cases}
\ X_1 = f^{-1}(Y_1)\\
\ X_2 = f^{-1}(Y_2)
\end{cases}
. \ 
$
Fix $n \geq 1$.  Assume: The sequence 
$* \ra \pi_0(X_i,X_1\cap X_2) \ra \pi_0(Y_i,Y_1 \cap Y_2)$ is exact $(i = 1,2)$ and that 
$f_*:\pi_q(X_i,X_1\cap X_2) \ra \pi_q(Y_i,Y_1 \cap Y_2)$ is bijective for $1 \leq q < n$ and surjective for $q = n$ $(i = 1, 2)$ 
$-$then the sequence 
$* \ra \pi_0(X,X_i) \ra \pi_0(Y,Y_i)$ is exact $(i = 1,2)$ and 
$f_*:\pi_q(X,X_i) \ra \pi_q(Y,Y_i)$ is bijective for $1 \leq q < n$ and surjective for $q = n$ $(i = 1,2)$.
\\ \indent
[Fix $i_0 \in \{1,2\}$, $0 \leq q < n$, and maps 
$\phi:(J^q,\dot{I}_0^q) \ra (X,X_{i_0})$, 
$\psi:(I^{q+1},I_0^q) \ra (Y,Y_{i_0})$, satisfying $f \circ \phi = \psi$ on $J^q$.  
In view of the lemma, it suffices to exhibit an extension 
$\Phi : (I^{q+1},I_0^q) \ra (X,X_{i_0})$ of $\phi$ and a homotopy 
$H:(I^{q+1},I_0^q) \times I \ra (Y,Y_{i_0})$ such that $\restr{H}{(J^q,\dot{I}_0^q) \times I}$ is the constant homotopy at 
$f \circ \phi$ and $f \circ \Phi \simeq \psi$ on $I^{q+1}$ by \mH.  
Subdivide $I^{q+1}$ into a collection $\sC$ of $(q+1)$-dimensional cubes $C:\forall \ C \in \sC$, $\exists$ $i_C \in \{1,2\}$: 
$\phi(C \cap J^q) \subset X_{i_C}$ and $\psi(C) \subset Y_{i_C}$ (possible, 
$
\begin{cases}
\ \phi^{-1}(X - X_1) \cup \psi^{-1}(Y - Y_1)\\
\ \phi^{-1}(X - X_2) \cup \psi^{-1}(Y - Y_2)
\end{cases}
$
being disjoint and closed). \   
Regard $I^{q+1}$ as $I^q \times I$ $-$then $\sC$ restricts to a subdivision of $I^q$ and induces a partition of \mI into subintervals $I_k = [a_{k-1},a_k]$: $0 = a_0 < a_1 < \cdots < a_r = 1$.  Break the subdivision of $I^q$ into its skeletal constituents \mD.  Construct $\Phi$ on $D \times I_k$ $\&$ \mH on $I(D \times I_k)$ via downward induction on $k$ and for fixed $k$, via upward induction on $\dim D$.  Arrange matters so that: 
(1) $\psi(D \times I_k) \subset Y_i$ $\implies$ 
%%----------------------------------------------------------------------------------------------47
$\Phi(D \times I_k) \subset X_i$ $\&$ $H(I(D \times I_k)) \subset Y_i$; 
(2) $\psi(D \times \{a_{k-1}\}) \subset Y_1 \cap Y_2$ $\implies$ $\Phi(D \times \{a_{k-1}\}) \subset X_1 \cap X_2$ $\&$ 
$H(I(D \times \{a_{k-1}\})) \subset Y_1 \cap Y_2$.  
The first condition plus the second when $k = 1$ yield 
$\Phi(I_0^q) \subset X_{i_0}$ $\&$ $H(I_0^q \times I) \subset Y_{i_0}$.  
At each stage, the induction hypothesis secures $\Phi$ on $\dot{D} \times I_k \cup D \times \{a_k\}$ $\&$ \mH on 
$I(\dot{D} \times I_k \cup D \times \{a_k\})$.  
Case 1:  If either $\psi(D \times \{a_{k-1}\})$ is not contained in $Y_1 \cap Y_2$ or $\psi(D \times I_k)$ is contained in $Y_1 \cap Y_2$, 
use the fact that $\dot{D} \times I_k \cup D \times \{a_k\}$ is a strong deformation retract of $D \times I_k$ to specify 
$\Phi$ on $D \times I_k$ $\&$ \mH on $I(D \times I_k)$.
Case 2: If $\psi(D \times \{a_{k-1}\})$ is contained in $Y_1 \cap Y_2$ and $\psi(D \times I_k)$ is contained in just one of the 
$Y_i$, realize 
$\Phi:(\dot{D} \times I_k \cup D \times \{a_k\},\dot{D} \times \{a_{k-1}\}) \ra (X_i,X_1 \cap X_2)$ $\&$ 
$H:(\dot{D} \times I_k \cup D \times \{a_k\},\dot{D} \times \{a_{k-1}\}) \times I \ra (Y_i,Y_1 \cap Y_2)$.  
Apply the lemma to produce the required extension of $\Phi$ to $D \times I_k$ $\&$ \mH to $I(D \times I_k)$.  
Here, of course, the assumption on $f$ comes in.]\\
\endgroup %%------------------------------------<<

%%%%%%%%%%%%%%%%%%%%%%%%%%%%%%%%%%%%%%
%%%%%%%%%%%%%%%%%%%%%%%%%%%%%%%%%%%%%%
%%%%%%%%%%%%%%%%%%%%%%%%%%%%%%%%%%%%%%

\begin{center}
$\S \ 3$
\\[0.5cm]
$\mathcal{REFERENCES}$\\
\end{center}

\[
\textbf{BOOKS}
\]

\begingroup
\fontsize{9pt}{11pt}\selectfont
\setlength\parindent{0 cm}

[1] \quad Brown, R., \textit{Topology}, Ellis Horwood (1988).
\\[-.2cm]

[2] \quad tom Dieck, T., Kamps, K., and Puppe, D., \textit{Homotopietheorie}, Springer Verlag (1970).
\\[-.2cm]

[3] \quad James, I., \textit{General Topology and Homotopy Theory}, Springer Verlag (1984).
\\[-.2cm]

[4] \quad James, I., \textit{Fibrewise Topology}, Cambridge University Press (1989).
\\[-.2cm]

[5] \quad Piccinini, R., \textit{Lectures on Homotopy Theory}, North Holland (1992).
\\[-.2cm]
 
%[6]  {\cyr Postniknov M., Lektsii po Algebraicheskoy Topologii}[{\cyr Osmpvy Teorii Gomotopiy}] {\cyr Nauka} (1984).\\
[6] 
\hspace{.25cm}
{\fontencoding{OT2}\selectfont
Postnikov M., Lektsii po Algebraicheskoy Topologii [Osmpvy Teorii Gomotopi\u], Nauka} (1984).
\\
\endgroup

\[
\textbf{ARTICLES}
\]

\begingroup
\fontsize{9pt}{11pt}\selectfont
\setlength\parindent{0 cm}

[1] \quad Cockcroft, W. and Jarvis, T., An Introduction to Homotopy Theory and Duality, 
\textit{Bull. Soc. Math.}

\hspace{.8cm}\textit{Belgique} \textbf{16} (1964), 407-427, and \textbf{17} (1965), 3-26.
\\[-.2cm]

[2] \quad Nomura, Y., An Application of the Path Space Technique to the Theory of Triads, 
\textit{Nagoya Math. J.} 

\hspace{.8cm}\textbf{22} (1963), 169-188.
\\[-.2cm]

[3] \quad Puppe, D., Homotopiemengen und ihre Induzierten Abbildungen I, 
\textit{Math. Zeit.} \textbf{69} (1958), 299-344.
\\[-.2cm]

[4] \quad Spanier, E., The Homotopy Excision Theorem, 
\textit{Michigan Math. J.} \textbf{14} (1967), 245-255.
\\[-.2cm]

[5] \quad Str{\o}m, A., Note on Cofibrations, 
\textit{Math. Scand.} \textbf{19} (1966), 11-14.
\\[-.2cm]

[6] \quad Str{\o}m, A., Note on Cofibrations II, 
\textit{Math. Scand.} \textbf{22} (1968), 130-142.

\setlength\parindent{2em}

\endgroup

\chapter{
$\boldsymbol{\S}$\textbf{4}.\quadx  FIBRATIONS}
\setlength\parindent{2em}
\setcounter{proposition}{0}
%%----------------------------------------------------------------------------------------------01
$\text{ }$\\[-1.5cm]

The technology developed below, like that of the preceding $\S$, underlies the foundations of homotopy theory in \bTOP or $\bTOP_*$.

Let $B$ be a toplogical space.  
An object  in $\bTOP/B$ is a topological space $X$ together with a continuous function $p:X \ra B$ called the 
\un{projection}.
\index{projection (fibration)}   
For $O \subset B$, put $X_O = p^{-1}(O)$, which is therefore an object  in 
$\bTOP/O$ (with projection $p_O = \restr{p}{X_O}$).  
The notation $\restr{X}{O}$ is also used.  In particular: $X_b = p^{-1}(b)$ is the 
\un{fiber}
\index{fiber} 
over $b \in B$.  
A morphism in $\bTOP/B$ is a continuous function $f:X \ra Y$ over $B$, i.e., an $f \in C(X,Y)$ such that the triangle
\begin{tikzcd}
& X \arrow{rr}{f} \arrow{dr}[swap]{p} &&Y \arrow{dl}{q}\\
&& B
\end{tikzcd}
commutes.  Notation: $f \in C_B(X,Y)$, $f_O = \restr{f}{X_O} \ (O \subset B)$.  
The base space $B$ is an object in $\bTOP/B$, where $p = \id_B$.  An element $s \in C_B(B,X)$ is called a 
\un{section}
\index{section} 
of $X$, written $s \in \text{sec}_B(X)$.

[Note: \  The product of 
$
\begin{cases}
\ p:X \ra B\\[-.15cm]
\ q:Y \ra B
\end{cases}
$
in $\bTOP/B$ is the fiber product: $X \times_B Y$.  If $B^\prime$ is a topological space and if $\Phi^\prime \in C(B^\prime, B)$, 
then $\Phi^\prime$ determines a functor $\bTOP/B \ra \bTOP/B^\prime$ that sends $X$ to $X^\prime = B^\prime \times_B X$.  
Obviously, $(X \times_B Y)^\prime = X^\prime \times_{B^\prime} Y^\prime$.]\\

\begingroup%%----------------------------------->>
\fontsize{9pt}{11pt}\selectfont
\textbf{\small EXAMPLE}  \ 
Let $X$ be in $\bTOP/B$ $-$then the assignment $O \ra \text{sec}_O(X_O)$, $O$ open in \mB, defines a sheaf of sets on $B$, the 
\un{sheaf of sections}
\index{sheaf of sections} 
$\Gamma_X$ of \mX.

[Note: \  Recall that for any sheaf of sets $\mathcal{F}$ on $B$, there exists an $X$ in 
$\bTOP/B$ with $p:X \ra B$ a local homeomorphism such that $\mathcal{F}$ is isomorphic to $\Gamma_X$.  
In fact, the category of sheaves of sets on $B$ is equivalent to the full subcategory of 
$\bTOP/B$ whose objects are those \mX for which $p:X \ra B$ is a local homeomorphism.]\\
\endgroup%%------------------------------------<<

\begingroup%%----------------------------------->>
\fontsize{9pt}{11pt}\selectfont
\textbf{\small FACT}  \ 
Let $X$ be in $\bTOP/B$ $-$then the projection $p:X \ra B$ is a local homeomorphism iff both it and the diagonal embedding $X \ra X \times_B X$ are open maps.\\
\endgroup%%------------------------------------<<

\begingroup%%----------------------------------->>
\fontsize{9pt}{11pt}\selectfont
\textbf{\small FACT}  \ 
Let $X$ be in $\bTOP/B$.  Assume: $X$ $\&$ $B$ are path connected Hausdorff spaces and the projection $p:X \ra B$ is a local homeomorphism $-$then $p$ is a homeomorphism iff $p$ is proper and $p_*:\pi_1(X) \ra \pi_1(B)$ is surjective.\\
\endgroup%%------------------------------------<<

There is a functor $\bTOP \ra \bTOP/B$ that sends a topological space \mT to $B \times T$ (product topology) with projection $B \times T \ra B$.  
An \mX in $\bTOP/B$ is said to be 
\un{trivial}
\index{trivial (fibration)} 
if there exists a \mT in $\bTOP$ such that \mX is homeomorphic over $B$ to $B \times T$, 
\un{locally trivial}
\index{locally trivial} 
if there exists an open covering 
$\{O\}$ of \mB such that $\forall \ O, X_O$ is trivial over $O$.

%%----------------------------------------------------------------------------------------------02
[Note: \  Spelled out, local triviality means that $\forall \  O$ there exists a topological space $T_O$ and a homeomorphism $X_O \ra O \times T_O$ over $O$.  If $T_O$ can be chosen independent of $O$, so  $\forall \  O$, $T_O = T$, then \mX is said to be 
\un{locally trival with fiber \mT}.
\index{locally trival with fiber \mT}  
When \mB is connected, this can always be arranged.]\\

\begingroup%%----------------------------------->>
\fontsize{9pt}{11pt}\selectfont
\textbf{\small FACT}  \ 
Let $X$ be in $\bTOP/IB$.  Suppose that $\restr{X}{(B \times [0,1/2])}$ and $\restr{X}{(B \times [1/2,1])}$ are trivial 
$-$then \mX is trivial.\\
\endgroup%%------------------------------------<<

\begingroup%%----------------------------------->>
\fontsize{9pt}{11pt}\selectfont
\textbf{\small EXAMPLE}  \ 
Let \mX be in $\bTOP/[0,1]^n$ $(n \geq 1)$.  
Suppose that \mX is locally trivial $-$then \mX is trivial.\\
\endgroup%%------------------------------------<<

A 
\un{fiber homotopy}
\index{fiber homotopy} 
is a homotopy over 
$B: f \underset{B}{\simeq} g$ $(f,g \in C_B(X,Y))$.  
Isomorphisms in the associated homotopy category are the fiber homotopy equivalences and any two 
$
\begin{cases}
\ X\\[-.15cm]
\ Y
\end{cases}
$
in $\bTOP/B$ for which there exists a fiber homotopy equivalence $X \ra Y$ have the same fiber homotopy type.  
The fiber homotopy type of $X \times_B Y$ depends only on the fiber homotopy types of \mX and \mY.  
The objects in $\bTOP/B$ that have the same fiber homotopy type of \mB itself are said to be 
\un{fiberwise contractible}.
\index{fiberwise contractible}   
Example:  The path space $PB$ with projection $p_0$ is in $\bTOP/B$ and is fiberwise contractible 
(consider the fiber homotopy $H: IPB \ra PB$ defined by $H(\sigma, t)(T) = \sigma(tT))$.

[Note: \  A fiber homotopy with domain $IB$ is called a 
\un{vertical homotopy}.]
\index{vertical homotopy}
\\

%%%%
\textbf{\small LEMMA}  \  Let \mX be in $\bTOP/B$.  
Assume: \mX is fiberwise contractible $-$then for any 
$\Phi^\prime \in C(B^\prime,B)$, $X^\prime$ is fiberwise contractible.\\
%%%%%

\begingroup%%----------------------------------->>
\fontsize{9pt}{11pt}\selectfont
Let $f:X \ra Y$ be a continuous function.  
View its mapping cylinder $M_f$ as an object in $\bTOP/Y$ with projection $r:M_f \ra Y$ $-$then $j \in \text{sec}_Y(M_f)$ and $M_f$ is fiberwise contractible.\\
\endgroup%%------------------------------------<<

Let \mX,\mY be in $\bTOP/B$ $-$then a fiber preserving function $f:X \ra Y$ is said to be 
\un{fiberwise constant}
\index{fiberwise constant} 
if $f = t \circx p$ for some section $t:B \ra Y$.  
Elements of $C_B(X,Y)$ that are fiber homotopic to a fiberwise constant function are 
\un{fiberwise inessential}.
\index{fiberwise inessential}\\

\begingroup%%----------------------------------->>
\fontsize{9pt}{11pt}\selectfont
Suppose that \mB is not in \bCG $-$then the identity map $kB \ra B$ is continuous and constant on fibers but not fiberwise constant.\\
\endgroup%%------------------------------------<<

\textbf{\small LEMMA}  \  Let \mX be in $\bTOP/B$ $-$the \mX is fiberwise contractible iff $\id_X$ is fiberwise inessential.\\

\begingroup%%----------------------------------->>
\fontsize{9pt}{11pt}\selectfont
\textbf{\small EXAMPLE}  \ 
Take $X = ([0,1] \times \{0,1\}) \cup (\{0\} \times [0,1])$, $B = [0,1]$, and let $p$ be the vertical projection $-$then \mX is contractible but not fiberwise contractible.\\
\endgroup%%------------------------------------<<

%%----------------------------------------------------------------------------------------------03

\begingroup%%----------------------------------->>
\fontsize{9pt}{11pt}\selectfont
\textbf{\small EXAMPLE}  \ 
Let \mX be a subspace of $B \times \R^n$ and suppose that there exists an $s \in \sec_B(X)$, say $b \ra (b,s(b))$, such that 
$\forall \ b \in B$, $\forall \ x \in X_b$, $\{(b,(1-t)s(b) + tx):0 \leq t \leq 1\} \subset X_b$ $-$then \mX is fiberwise contractible.\\
\endgroup%%------------------------------------<<

\begingroup%%----------------------------------->>
\fontsize{9pt}{11pt}\selectfont
\textbf{\small FACT}  \ 
Let \mX be in $\bTOP/B$; let $f, g \in C_B(X,X)$.  Suppose that $\{O,P\}$ is a numerable covering of \mB for which 
$
\begin{cases}
\ f_O\\
\ g_P
\end{cases}
$
are fiberwise inessential $-$then $g \circx f$ is fiberwise inessential.

[Fix fiber homotopies 
$
\begin{cases}
\ K:IX_O \ra X_O\\
\ L:IX_P \ra X_P
\end{cases}
$
between
$
\begin{cases}
\ f_O \ \& \ k \circx p_O\\
\ g_P \ \& \ l \circx p_P
\end{cases}
, \ 
$
where 
$
\begin{cases}
\ k \in \sec_O(X_O)\\
\ l \in \sec_P(X_P)
\end{cases}
. \ 
$
Through reparametrization, it can be assumed that 
$
\begin{cases}
\ K \circx i_t\\
\ L \circx i_t
\end{cases}
$
are independent of $t$ when $0 \leq t \leq 1/4$, $3/4 \leq t \leq 1$.  Choose 
$
\begin{cases}
\ \mu\\
\ \nu
\end{cases}
\in C(B,[0,1]) \ :
$
$
\begin{cases}
\ \sptx \mu \subset O\\
\ \sptx \nu \subset P
\end{cases}
$
$\&$ $\mu + \nu = 1$.  Let $\Delta$ be the triangle in $\R^2$ with the vertexes $(0,0)$, $(1,0)$, $(0,1)$.  
Note that the transformation $(\xi,\eta) \ra (\xi,(1-\xi)\eta)$ takes $I[0,1] - I\{1\}$ homeomorphically onto $\Delta - \{(1,0)\}$.  
The continuous fiber preserving function $\Phi:I^2 X_{O \cap P} \ra X_{O \cap P}$ defined by $\Phi(x,(\xi,\eta)) = L(K(x,\eta),\xi)$ is independent of $\eta$ when $\xi = 1$, 
thus it induces a continuous fiber preserving function 
$\Phi_\Delta:X_{O \cap P} \times \Delta \ra X_{O \cap P}$.  On $X_{O \cap P} \times \fr\Delta$, one has 
$\Phi_\Delta(x,(t,1-t)) = L(k(p(x)),t)$, 
$\Phi_\Delta(x,(0,t)) = g(K(x,t))$, 
$\Phi_\Delta(x,(t,0)) =L(f(x),t)$.  Write
$
s(b) = 
\begin{cases}
\ L(k(b),\nu(b)) \hspace{0.35cm} \ \  (b \in O \cap P)\\ 
\ g(k(b))  \hspace{1.13cm} \  \ (b \in O - P)\\
\ l(b) \hspace{1.85cm} (b \in P - O)
\end{cases}
$
$-$then $s \in \sec_B(X)$ and $g \circx f$ is fiber homotopic to $s \circx p$ via 
\endgroup%%------------------------------------<<

\begingroup%%----------------------------------->>
\fontsize{9pt}{11pt}\selectfont
\[
H(x,t) = 
\begin{cases}
\ \Phi_\Delta(x,t(\nu(b),\mu(b))) \hspace{0.5cm}  \ \ (b \in O \cap P)\\
\ g(K(x,t)) \hspace{2.01cm} (b \in O - P)\\
\ L(f(x),t) \hspace{1.98cm} \  (b \in P - O)
\end{cases}
(x \in X_b).]
\]
\\[-.5cm]
\endgroup%%------------------------------------<<

\begingroup%%----------------------------------->>
\fontsize{9pt}{11pt}\selectfont
Consequently, if $f_1, \ldots, f_n \in C_B(X,X)$ and if $O_1, \ldots, O_n$ is a numerable covering of \mB such that 
$\forall \ i$, $f_{O_i}$ is fiberwise inessential, then $f_1 \circx \cdots \circx f_n$ is fiberwise inessential.  
Example: $X_{O_i}$ is fiberwise contractible $(i = 1, \ldots, n)$ $\implies$ \mX is fiberwise contractible (cf. p. \pageref{4.1}).\\ 
\endgroup%%------------------------------------<<

Let \mX be in $\bTOP/B$ $-$then \mX is said to have the 
\un{section extension property}
\index{section extension property} 
(SEP) 
\index{SEP} 
provided that for each $A \subset B$, 
every section $s_A$ of $X_A$ which admits an extension $s_O$ to a halo \mO of \mA in \mB can be extended to a section $s$ of \mX: $\restr{s}{A} = s_A$.

[Note: \ If \mX has the SEP, then $\sec_B(X)$ is nonempty (take $A = \emptyset = O$.]\\

\begingroup%%----------------------------------->>
\fontsize{9pt}{11pt}\selectfont

Let $X$ be in $\bTOP/B$ and suppose that \mX has the SEP.  
Let $s$ be a section of $\restr{X}{\phi^{-1}(]0,1])}$, where 
$\phi \in C(B,[0,1])$ $-$then $\forall \ \epsilon$, $0 <  \epsilon < 1$, 
$\restr{s}{\phi^{-1}([\epsilon,1])}$ can be extended to a 
section $s_\epsilon$ of \mX but it is false in general that $s$ can be so extended.\\
\endgroup%%------------------------------------<<

\begingroup%%----------------------------------->>
\fontsize{9pt}{11pt}\selectfont
\textbf{\small EXAMPLE}  \ 
Suppose that \mB is a CW complex of combinatorial dimension \ $\leq n + 1$ and \mT is $n$-connected $-$then $B \times T$ has the SEP.\\
\endgroup%%------------------------------------<<

%%----------------------------------------------------------------------------------------------04
\begin{proposition} \ %01
Let \mX, \mY be in $\bTOP/B$ and suppose that \mY has the SEP.  Assume: $\exists$ 
$
\begin{cases}
\ f \in C_B(X,Y)\\
\ g \in C_B(Y,X)
\end{cases}
: g \circx f \underset{B}{\simeq} \id_X
$
$-$then \mX has the SEP.
\end{proposition}

[Fix a fiber homotopy $H:IX \ra X$ between $\id_X$ and $g \circx f$.  
Given $A \subset B$, let $s_A$ be a section of $X_A$ which admits an extension $s_O$ to a halo \mO of \mA in \mB.  
Choose a closed halo \mP of \mA in \mB: 
$A \subset P \subset O$ and \mO a halo of \mP in \mB (cf. HA$_2$, p. \pageref{4.2}).
Since \mY has the SEP, there exists a section $t$ of \mY: $\restr{t}{P} = f \circx \restr{s_O}{P}$.  
With $\pi$ a haloing function of \mP, define 
$s:B \ra X$ by
$
s(b) = 
\begin{cases}
\ g \circx t(b) \hspace{2.37cm} (b \in \pi^{-1}(0))\\
\ H(s_O(b),1 - \pi(b)) \hspace{0.5cm} (b \in P)
\end{cases}
$
to get a section $s$ of \mX: $\restr{s}{A} = s_A$.]
\\

Application: Fiberwise contractible spaces have the SEP.
\\

\textbf{\small LEMMA}  \ 
Let \mX be in $\bTOP/B$ and suppose that \mX has the SEP.  Let \mO be a cozero set in \mB $-$then $X_O$ has the SEP.

[There is no loss of generality in assuming that $A = f^{-1}(]0,1])$, where $f \in C(O,[0,1])$.  
Accordingly, given a section $s_A$ of $X_A$, it will be enough to construct a section $s$ of $X_O$ which agrees with $s_A$ on $f^{-1}(1)$.  
Fix $\phi \in C(B,[0,1])$ : $O = \phi^{-1}(]0,1])$.  
Claim: There exist sections $s_2, s_3, \ldots$ of \mX such that $s_{n+1}(b) = s_n(b)$ 
$\bigl(\phi(b) > \ds\frac{1}{n}\bigr)$ and 
$s_n(b) = s_A(b)$ $\bigl(f(b) > 1 - \ds\frac{1}{n} \ \& \ \phi(b) > \ds\frac{1}{n+1}\bigr)$.  
Granted the claim, we are done.  
Put 
$
F(b) = 
\begin{cases}
\ f(b) \phi(b) \hspace{0.4cm} (b \in O)\\
\ 0 \hspace{1.6cm} (b \in B - O)
\end{cases}
\hspace{-.2cm}:
$
$F \in C(B,[0,1])$.  
Since \mX has the SEP and $s_A$ is defined on $F^{-1}(]0,1])$, a halo of $F^{-1}([1/6,1])$ in \mB, there exists a section of \mX that agrees with $s_A$ on $f^{-1}(]1/2,1]) \cap \phi^{-1}(]1/3,1])$.  
Call it $s_2$, thus setting the stage for induction.
Choose continuous functions $\mu_n$, $\nu_n:[0,1] \ra [0,1]$ subject to 
$\ds\frac{1}{n+3} < \nu_n(x) < \mu_n(x) \leq \ds\frac{1}{n}$ with 
$\mu_n(x) \leq \ds\frac{1}{n+2}$ $\bigl(x \geq 1 - \ds\frac{1}{n+1}\bigr)$  and 
$\nu_n(x) \geq \ds\frac{1}{n+1}$ $\bigl(x \leq 1 - \ds\frac{1}{n}\bigr)$ $(n = 2, 3, \ldots)$.  
Let 
$A_n = \{b \in O: \phi(b) > \mu_n(f(b))\}$, $O_n = \{b \in O:\phi(b) > \nu_n(f(b))\}$ $-$then $O_n$ is a halo of $A_n$ in \mB, a haloing function being 1 on $\{b \in O: \mu_n(f(b)) \leq \phi(b)\}$, 
\[
\frac{\phi(b) - \nu_n(f(b))}{\mu_n(f(b)) - \nu_n(f(b))}
\quad \text{ on } \quad
\{b \in O: \nu_n(f(b)) \leq \phi(b) \leq \mu_n(f(b))\},
\]
and 0 on $\{b \in O: \phi(b) \leq \nu_n(f(b))\} \cup B - O$.  
To pass from $n$ to $n + 1$, note that the prescription 
$
b \ra 
\begin{cases}
\ s_n(b) \hspace{0.5cm} \bigl(\phi(b) > \ds\frac{1}{n + 1}\bigr)\\
\ s_A(b) \hspace{0.5cm} \bigl(f(b) > 1 - \ds\frac{1}{n}\bigr)
\end{cases}
$
defines a section of $X_{O_n}$.  
Its restriction to $A_n$ can therefore be extended to a section $s_{n+1}$ of \mX with the required properties.]
\\

\index{Section Extension Theorem}
\index{Theorem: Section Extension Theorem}
%\textbf{\small SECTION EXTENSION THEOREM}
\textbf{\small SECTION EXTENSION THEOREM}
Let \mX be in $\bTOP/B$.  Suppose that $\sO = \{O_i:i \in I\}$ is a numerable covering of \mB such that $\forall \ i$, 
$X_{O_i}$ has the SEP $-$then \mX has the SEP.

%%----------------------------------------------------------------------------------------------05
[Given $A \subset B$, let $s_A$ be a section of $X_A$ which admits an extension $s_O$ to a halo \mO of \mA in \mB.  
Fix a haloing function $\pi$ for \mO and let $\{\pi_i:i \in I\}$ be a partition of unity on \mB subordinate to $\sO$.  
Put 
$\Pi_S = \sum\limits_{i \in S} (1 - \pi)\pi_i + \pi$ $(S \subset I)$.  
Consider the set $\sS$ of all pairs $(S,s)$:  
$s$ is a section of $\restr{X}{\Pi_S^{-1}(]0,1])}$ $\&$ $\restr{s}{A} = s_A$: $\sS$ is nonempty (take $S = \emptyset$, 
$s = \restr{s_O}{\pi^{-1}(]0,1])})$.  
Order $\sS$ by stipulating that 
$(S^\prime,s^\prime) \leq (S\pp,s\pp)$ iff $S^\prime \subset S\pp$ and $s^\prime(b) = s\pp(b)$ when 
$\Pi_{S^\prime}(b) = \Pi_{S\pp}(b) > 0$.  
One can check that every chain in $\sS$ has an upper bound, so by Zorn, $\sS$ has a maximal element $(S_0,s_0)$.  
Since $\Pi_I = 1$, to finish it need only be shown that $S_0 = I$.  
Suppose not.  Select an $i_0 \in I - S_0$, set 
$\Pi_0 = \Pi_{S_0}$ $\&$ $\pi_0 = (1 - \pi)\pi_{i_0}$, and define a continuous function 
$\phi_0:\pi_0^{-1}(]0,1]) \ra [0,1]$ by $\phi_0(b) = \min\{1,\Pi_0(b)/\pi_0(b)\}$.  
Owing to the lemma, $\restr{X}{\pi_0^{-1}(]0,1])}$ has the SEP $(\pi_0^{-1}(]0,1]) \subset O_{i_0})$.  
On the other hand, 
$\phi_0^{-1}(]0,1])$ is a halo of 
$\phi_0^{-1}(1)$ in $\pi_0^{-1}(]0,1])$ and 
$\restr{s_0}{\phi_0^{-1}(1)}$ 
admits an extension to 
$\phi_0^{-1}(]0,1])$, viz. 
$\restr{s_0}{\phi_0^{-1}(]0,1])}$.  Therefore 
$\restr{s_0}{\phi_0^{-1}(1)}$  can be 
extended to a section 
$s_{i_0}$ of 
$\restr{X}{\pi_0^{-1}(]0,1])}$.  Let $T = S_0 \cup \{i_0\}$ and write
$
t(b) = 
\begin{cases}
\ s_0(b) \quadx (\pi_0(b) \leq \Pi_0(b))\\
\ s_{i_0}(b) \quadx (\pi_0(b) \geq \Pi_0(b))
\end{cases}
(\Pi_T(b) > 0)
$
$-$then $(T,t) \in \sS$ and $(S_0,s_0) < (T,t)$, contradicting the maximality of $(S_0,s_0)$.]\\

\begingroup%%----------------------------------->>
\fontsize{9pt}{11pt}\selectfont
\label{3.3}
\textbf{\small FACT}  \ 
Let \mA be a subspace of \mX.  
Suppose that there exists a numerable covering $\sU = \{U_i:i \in I\}$ of \mX such that 
$\forall \ i$, the inclusion $A \cap U_i \ra U_i$ is a cofibration $-$then the inclusion $A \ra X$ is a cofibration.
\\ \indent
[Let $\{\kappa_i:i \in I\}$ be a partition of unity on \mX subordinate to $\sU$.  The lemma on p. \pageref{4.3} implies that 
$\forall \ i$, the inclusion 
$A \cap \kappa_i^{-1}(]0,1]) \ra \kappa_i^{-1}(]0,1])$ is a cofibration.  Therefore one can assume that $\sU$ is numerable 
and open.  Fix a topological space \mY and a pair $(F,h)$ of continuous functions 
$
\begin{cases}
\ F:X \ra Y\\
\ h:IA \ra Y
\end{cases}
$
such that $\restr{F}{A} = h \circx i_0$.  
Define a sheaf of sets $\sF$ on \mX by assigning to each open set \mU the set of all continuous functions $H:IU \ra Y$ such that $\restr{F}{U} = H \circx i_0$ and $\restr{H}{I(A \cap U)} = \restr{h}{I(A \cap U)}$.  
Choose a topological space $E$ and a local homeomorphism $p:E \ra X$ for which $\sF(U) = \sec_U(E_U)$ at each \mU.  
Show that $\forall \ i$, 
$E_{U_i}$ has the SEP.  The section extension theorem then says that $\exists$ $H \in \sF(X)$.]\\
\endgroup%%------------------------------------<<

Let \mX be in $\bTOP/B$.  Let \mE be in \bTOP; let $\phi \in C(E,B)$ $-$then a continuous function 
$\Phi:E \ra X$ is a \un{lifting}
\index{lifting} 
of $\phi$ provided that $p \circx \Phi = \phi$.  
Example: Every $s \in \sec_B(X)$ is a lifting of $\id_B$.\\

\begingroup%%----------------------------------->>
\fontsize{9pt}{11pt}\selectfont
\textbf{\small FACT}  \ 
Suppose that \mX is fiberwise contractible.  Let $\phi \in C(E,B)$ $-$then for any halo \mU of any \mA in \mE and all 
$\psi \in C(U,X)$: $p \circx \psi = \restr{\phi}{U}$, there exists a lifting $\Phi$ of $\phi:\restr{\Phi}{A} = \restr{\psi}{A}$.
\\ \indent
[Note: \ The condition is also characteristic.  
First take $E = B$, $A = \emptyset = U$, and $\phi = \id_B$ to see that 
$\exists$ $s \in \sec_B(X)$.  
Next let $E = IX$, $A = i_0 X \cup i_1X$, 
$U = X \times [0,1/2[ \ \cup \hsx X \times \hsx ]1/2,1]$, and 
define $\phi:IX \ra B$ by $\phi(x,t) = p(x)$, $\psi:U \ra X$ by 
$
\psi(x,t) = 
\begin{cases}
\ x \hspace{1.4cm} (t < 1/2)\\
\ s \circx p(x) \hspace{0.5cm} (t > 1/2)
\end{cases}
\hspace{-.2cm}. \ 
$
Since \mU is a halo of \mA in $IX$, every lifting $\Phi$ of $\phi$ with $\restr{\Phi}{A} = \restr{\psi}{A}$ is a fiber homotopy between 
$\id_X$ and $s \circx p$, i.e., \mX is fiberwise contractible.]\\
\endgroup%%------------------------------------<<

%%----------------------------------------------------------------------------------------------06
\indent\indent (HLP) \ Let \mY be a topological space $-$then the projection $p:X \ra B$ is said to have the 
\un{homotopy lifting property with respect to \mY}
\index{homotopy lifting property with respect to \mY}
(HLP w.r.t \mY) 
\index{HLP w.r.t \mY}
if given continuous functions 
$
\begin{cases}
\ F:Y \ra X\\[-.15cm]
\ h:IY \ra B
\end{cases}
$
such that $p \circx F = h \circx i_0$, there is a continuous function $H:IY \ra X$ such that $F = H \circx i_0$ and $p \circx H = h$.\\

\begingroup%%----------------------------------->>
\fontsize{9pt}{11pt}\selectfont
If $p:X \ra B$ has the HLP w.r.t. \mY and if 
$
\begin{cases}
\ f \in C(Y,B)\\
\ g \in C(Y,B)
\end{cases}
$
are homotopic, then $f$ has a lifting $F \in C(Y,X)$ iff $g$ has a lifting $G \in C(Y,X)$.\\
\endgroup%%------------------------------------<<

\begingroup%%----------------------------------->>
\fontsize{9pt}{11pt}\selectfont
\textbf{\small EXAMPLE}  \ 
Take $X = [0,1] \amalg *$, $B = [0,1]$ and define $p:X \ra B$ by $p(t) = t$, $p(*) = 0$.  
Fix a nonempty \mY and let $f$ be 
the constant map $Y \ra 0$ $-$then the constant map $Y \ra *$ is a lifting $F \in C(Y,X)$ of $f$.  
Put $h(y,t)  = t$, so 
$h:IY \ra B$.  
Obviously, $p \circx F = h \circx i_0$ but there does not exist $H \in C(IY,X)$: $F = H \circx i_0$ and 
$p \circx H = h$.\\
\endgroup%%------------------------------------<<

Let \mX be in \ $\bTOP/B$.  \ Given a topological space \mY and continuous functions \ 
$
\begin{cases}
\ F:Y \ra X\\[-.15cm]
\ h:IY \ra B
\end{cases}
$
such that $p \circx F = h \circx i_0$, let \mW be the subspace of $Y \times PX$ consisting of the pairs 
$(y,\sigma)$: $F(y) = \sigma(0)$ $\&$ $h(y,t) = p(\sigma(t))$ $(0 \leq t \leq 1)$.  View \mW as an object in $\bTOP/Y$ with 
projection $(y,\sigma) \ra y$.

\label{4.7}
\textbf{\small LEMMA}  \ 
The \cd 
\begin{tikzcd}[sep=large]
{Y} \ar{d}[swap]{i_0} \ar{r}{F} &{X} \ar{d}{p}\\
{IY} \ar{r}[swap]{h} &{B}
\end{tikzcd}
admits a filler $H:IY \ra X$ iff $\sec_Y(W) \neq \emptyset$.\\

\begin{proposition} \ %02
Suppose that $p:X \ra B$ has the HLP w.r.t \mY $-$then $\forall$ pair $(F,h)$, \mW has the SEP.
\end{proposition}

[Fix $A \subset Y$ and let \mV be a halo of \mA in \mY for which there exists a homotopy $H_V:IV \ra X$ such that
$\restr{F}{V} = H_V \circx i_0$ and $p \circx H_V = \restr{h}{IV}$.  
To construct a homotopy $H:IY \ra X$ such that $F = H \circx i_0$ and 
$p \circx H = h$, with $\restr{H}{IA} = \restr{H_V}{IA}$, take \mV closed 
(cf. HA$_2$, p. \pageref{4.4}) and using a haloing function $\pi$, 
put $\ov{h}(y,t) = h(y,\min\{1,\pi(y) + t\})$, so $\ov{h}:IY \ra B$.  
Define $\ov{H}_V:i_0Y \cup IV \ra X$ by 
$
\begin{cases}
\ \ov{H}_V(y,0) = F(y)\\
\ \ov{H}_V(y,t) = H_V(y,t)
\end{cases}
$
and define $\ov{F}:Y \ra X$ by $\ov{F}(y) = \ov{H}_V(y,\pi(y))$.  
Since $p \circx \ov{F} = \ov{h} \circx i_0$, there is a 
continuous function $\ov{H}:IY \ra X$ such that $\ov{F} = \ov{H} \circx i_0$ and $p \circx \ov{H} = \ov{h}$.  
The rule
\[
H(y,t) = 
\begin{cases}
\ \ov{H}_V(y,t) \hspace{1.75cm} (0 \leq t \leq \pi(y))\\
\ \ov{H}(y,t - \pi(y)) \hspace{0.8cm}  (\pi(y) \leq t \leq 1)
\end{cases}
\]
then specifies a homotopy  $H:IY \ra X$ having the properties in question.]\\

%%----------------------------------------------------------------------------------------------07
Let $\sY$ be a class of topological spaces $-$then $p:X \ra B$ is said to be a 
\un{$\sY$ fibration}
\index{fibration ($\sY$)} 
if $\forall \ Y \in \sY$, $p:X \ra B$ has the HLP w.r.t. \mY.

\indent\indent (H) \ Take for $\sY$ the class of topological spaces $-$then a $\sY$ fibration $p:X \ra B$ is called a 
\un{Hurewicz fibration}.
\index{Hurewicz fibration}

\indent\indent (S) \ Take for $\sY$ the class of CW complexes $-$then a $\sY$ fibration $p:X \ra B$ is called a 
\un{Serre fibration}.
\index{Serre fibration}

Every Hurewicz fibration is a Serre fibration.  
The converse is false (cf. p. \pageref{4.5}).
\label{16.35}

Observation: Let $Y \in \sY$ and suppose that $p:X \ra B$ is a $\sY$ fibration $-$then any inessential $f \in C(Y,B)$ admits a lifting $F \in C(Y,X)$.

[Note: \ It is thus a corollary that if $B \in \sY$ is contractible, then $\sec_B(X)$ is nonempty.]\\

\label{12.11}
\label{18.24}
\begingroup%%----------------------------------->>
\fontsize{9pt}{11pt}\selectfont
Other possibilities suggest themselves.  For example, one could consider $p:X \ra B$, where both \mX and \mB are in \bCG, and work with the class $\sY$ of compactly generated spaces.  This leads to the notion of 
\un{\bCG fibration}.
\index{fibration (\bCG)}  
\label{13.8}
\label{13.39}
Any \bCG fibration is a Serre fibration.  In general, if $p:X \ra B$ is a Hurewicz 
fibration, then $kp:kX \ra kB$ is a \bCG fibration.  
Another variant would be to consider pointed spaces and pointed homotopies.  
Via the artifice of adding a disjoint base point 
(cf. p. \pageref{4.6}), one sees that every pointed Hurewicz fibration is a Hurewicz fibration.  
In the opposite direction, an $f \in C_B(X,Y)$ is said to be a 
\un{fiberwise Hurewicz fibration}
\index{fiberwise Hurewicz fibration}
if it has the fiber homotopy lifting property with repsect to all \mE in $\bTOP/B$.  
Of course, if $f$ is a Hurewicz fibration, then $f$ is a fiberwise Hurewicz fibration.  
On the other hand, for any \mX in $\bTOP/B$, the projection $p:X \ra B$ is always a fiberwise Hurewicz fibration.\\
\endgroup%%------------------------------------<<

\begingroup%%----------------------------------->>
\fontsize{9pt}{11pt}\selectfont
\textbf{\small FACT}  \ 
Suppose that $p:X \ra B$ is a Hurewicz fibration.  Let \mE be a topological space with the homotopy type of a compactly generated space $-$then a $\phi \in C(E,B)$ has a lifting $E \ra X$ iff $k\phi \in C(kE,kB)$ has a lifting $kE \ra kX$.
\\ \indent
[The identity map $kE \ra E$ is a homotopy equivalence.]\\
\endgroup%%------------------------------------<<

\begingroup%%----------------------------------->>
\fontsize{9pt}{11pt}\selectfont
\textbf{\small EXAMPLE}  \ 
For any topological space \mT, the projection $B \times T \ra B$ is a Hurewicz fibration.  Take, e.g., $T = \bD^n$, let 
$X_0 \subset B \times \bS^{n-1}$, and put $X = B \times \bD^n - X_0$ $-$then the restriction to \mX of the projection 
$B \times \bD^n \ra B$ is a Hurewicz fibration.\\
\endgroup%%------------------------------------<<

\begingroup%%----------------------------------->>
\fontsize{9pt}{11pt}\selectfont
\label{5.0s}
\index{Covering Spaces (example)}
\textbf{\small EXAMPLE \  (\un{Covering Spaces})} \ 
A continuous function $p:X \ra B$ is said to be a 
\un{covering} \un{projection}
\index{covering projection} 
if each $b \in B$ 
has a neighborhood \mO such that $X_O$ is trivial with discrete fiber.  Every covering projection is a Hurewicz fibration.
\\ \indent
[Note: \ A sheaf of sets $\sF$ on \mB is 
\un{locally constant}
\index{locally constant} 
provided that each $b \in B$ 
has a basis $\sB$ of neighborhoods such that whenever $U$, $V \in \sB$ with $U \subset V$, the restriction map 
$\sF(V) \ra \sF(U)$ is a bijection.  
If $p:X \ra B$ is a covering projection, then its sheaf of sections $\Gamma_X$ is locally constant.  
Moreover, every locally constant sheaf of sets $\sF$ on \mB can be so realized.]\\
\endgroup%%------------------------------------<<

%%----------------------------------------------------------------------------------------------08
\begingroup%%----------------------------------->>
\fontsize{9pt}{11pt}\selectfont
\textbf{\small EXAMPLE}  \ 
Let \mX be the triangle in $\R^2$ with vertexes $(0,0)$, $(1,0)$, $(0,1)$ $-$then the vertical projection $p:X \ra [0,1]$ is a Hurewicz fibration but \mX is not locally trivial.
\\ \indent
[Note: \ 
Ferry\footnote[2]{\textit{Trans. Amer. Math. Soc.} \textbf{327} (1991), 201-219; 
see also Husch, \textit{Proc. Amer. Math. Soc.} \textbf{61} (1976), 155-156.}
has constructed an example of a Hurewicz fibration $p:X \ra [0,1]$ whose fibers are connected $n$-manifolds but such that \mX is not locally trivial.]\\
\endgroup%%------------------------------------<<

\label{12.18}
\textbf{\small LEMMA}  \ 
Let \mX be in $\bTOP/B$ $-$then $p:X \ra B$ is a Serre fibration iff it has the HLP w.r.t. the $[0,1]^n$ $(n \geq 0)$.\\

\label{4.5}
\label{4.8}
\begingroup%%----------------------------------->>
\fontsize{9pt}{11pt}\selectfont
\textbf{\small EXAMPLE}  \ 
Take $X = \{(x,-x): 0 \leq x \leq 1\} \cup \ds\bigcup\limits_1^\infty ([0,1] \times \{1/n\})$, $B = [0,1]$, 
and let $p$ be the vertical projection $-$then $p$ is a Serre fibration but not a Hurewicz fibration.
\\ \indent
[Note: \  $p^{-1}(0)$ and $p^{-1}(1)$ do not have the same homotopy type.]\\
\endgroup%%------------------------------------<<

\begingroup%%----------------------------------->>
\fontsize{9pt}{11pt}\selectfont
\textbf{\small EXAMPLE}  \ 
Let \mB be a topological space which is not compactly generated $-$then $\Gamma B$ is not compactly generated  and the identity map $k\Gamma B \ra \Gamma B$ is a Serre fibration but not a Hurewicz fibration.
\\ \indent
[For any compact Hausdorff space \mK, the arrow $C(K,k\Gamma B) \ra C(K,\Gamma B)$ is a bijection.]\\
\endgroup%%------------------------------------<<

\label{4.62}
\begingroup%%----------------------------------->>
\fontsize{9pt}{11pt}\selectfont
\textbf{\small EXAMPLE}  \ 
Let $B = [0,1]^\omega$, the Hilbert cube.  Put $X = B \times B - \Delta_B$ and let $p$ be the vertical projection, $q$ the horizontal projection $-$then $p:X \ra B$ is a Serre fibration.  
Moreover, \mB is an AR as are the $X_b$ (each being homeomophic to 
$B \times [0,1[$) but $p:X \ra B$ is not a Hurewicz fibration.
\\ \indent
[If so, then there would exist an $s \in \sec_B(X)$.  
Consider $q \circx s$: It is a continous function $B \ra B$ without a fixed point, contradicting Brouwer.]\\
\endgroup%%------------------------------------<<

\begingroup%%----------------------------------->>
\fontsize{9pt}{11pt}\selectfont
Ungar\footnote[2]{\textit{Pacific J. Math.} \textbf{30} (1969), 549-553.}
has shown that if \mX and \mB are compact ANRs of finite topological dimension, then a Serre fibration $p:X \ra B$ is necessarily a Hurewicz fibration.\\
\endgroup%%------------------------------------<<

The projection $p:X \ra B$ is a Hurewicz fibration iff the commutative diagram 
\begin{tikzcd}[sep=large]
{PX} \ar{d}[swap]{Pp} \ar{r}{p_0} &{X} \ar{d}{p}\\
{PB} \ar{r}[swap]{p_0} &{B}
\end{tikzcd}
is a weak pullback square.  Homeomorphisms are Hurewicz fibrations.  Maps with an empty domain are Hurewicz fibrations.  The composite of two Hurewicz fibrations is a Hurewicz fibration.\\

\begin{proposition} \ %03
Let 
$
\begin{cases}
\ p_1:X_1 \ra B_1\\
\ p_2:X_2 \ra B_2
\end{cases}
$
be Hurewicz fibrations $-$then $p_1 \times p_2:X_1 \times X_2 \ra B_1 \times B_2$ is Hurewicz fibration.\\
\end{proposition}

%%----------------------------------------------------------------------------------------------09
\begin{proposition} \  %04
Let 
\begin{tikzcd}[sep=large]
{X^\prime} \ar{d}[swap]{p^\prime} \ar{r} &{X} \ar{d}{p}\\
{B^\prime} \ar{r} &{B}
\end{tikzcd}
be a pullback square.  
Suppose that $p$ is a Hurewicz fibration $-$then $p^\prime$ is a Hurewicz fibration.\\
\end{proposition}

Application: Let $p:X \ra B$ be a Hurewicz fibration $-$then $\forall \ O \subset B$, $p_O:X_O \ra O$ is a Hurewicz fibration.\\

\begin{proposition} \  %05
Let $p:X \ra B$ be a Hurewicz fibration $-$then for any LCH space \mY, the postcomposition arrow 
$p_*:C(Y,X) \ra C(Y,B)$ is a Hurewicz fibration (compact open topology).
\end{proposition}

[Convert
\[
\begin{tikzcd}%[sep=large]
{E} \ar{d} \ar{r} &{C(Y,X)} \ar{d}\\
{IE} \ar[dashed]{ru}\ar{r} &{C(Y,B)}
\end{tikzcd}
\text{\quadx to \quadx}
\begin{tikzcd}%[sep=large]
{E \times Y} \ar{d} \ar{r} &{X} \ar{d}\\
{I(E \times Y)} \ar[dashed]{ru}\ar{r} &{B}
\end{tikzcd}
.]
\]
\vspace{0.2cm}

Application: Let $p:X \ra B$ be a Hurewicz fibration $-$then $Pp:PX \ra PB$ is a Hurewicz fibration.\\

\begin{proposition} \ %06
Let $i:A \ra X$ be a closed cofibration, where \mX is a LCH space $-$then for any topological space \mY, the precomposition 
arrow $i^*:C(X,Y) \ra C(A,Y)$ is a Hurewicz fibration (compact open topology).\
\end{proposition}

[Convert
\[
\begin{tikzcd}%[sep=large]
{E} \ar{d} \ar{r} &{C(X,Y)} \ar{d}\\
{IE} \ar[dashed]{ru} \ar{r} &{C(A,Y)}
\end{tikzcd}
\text{\quadx to \quadx}
\begin{tikzcd}%[sep=large]
{E \times X} \ar{d} \ar{r} &{Y} \\
{I(E \times X)} \ar[dashed]{ru} &{I(E \times A)} \ar{l} \ar{u}
\end{tikzcd}
.]
\]
\vspace{0.2cm}

Application: Let \mX be a topological space $-$then $p_t:PX \ra X$ $(0 \leq t \leq 1)$ is a Hurewicz fibration.\\

\label{9.33}
\begingroup%%----------------------------------->>
\fontsize{9pt}{11pt}\selectfont
\textbf{\small EXAMPLE}  \ 
Let $i:A \ra X$ be a closed cofibration, where \mX is a LCH space.  Fix $a_0 \in A$ and put $x_0 = i(a_0)$ $-$then for any pointed topological space $(Y,y_0)$, the precomposition arrow $i^*:C(X,x_0;Y,y_0) \ra C(A,a_0;Y,y_0)$ is a Hurewicz fibration (compact open topology).
\\ \indent
[The \cd 
\begin{tikzcd}[sep=large]
{C(X,x_0;Y,y_0)} \ar{d} \ar{r} &{C(X,Y)}  \ar{d}\\
{C(A,a_0;Y,y_0)}  \ar{r} &{C(A,Y)}
\end{tikzcd}
is a pullback square.]\\
\endgroup%%------------------------------------<<

%%----------------------------------------------------------------------------------------------10
\label{12.25}
\begingroup%%----------------------------------->>
\fontsize{9pt}{11pt}\selectfont
\textbf{\small FACT}  \ 
Let \mX be a topological space $-$then 
$
\Pi:
\begin{cases}
\ PX \ra X \times X\\
\ \sigma \ra (\sigma(0),\sigma(1))
\end{cases}
$
is a Hurewicz fibration.  Moreover, \mX is locally path connected iff $\Pi$ is open.
\\ \indent
[Note: \ Fix $x_0 \in X$ $-$then the fiber of $\Pi$ over $(x_0,x_0)$ is $\Omega X$, the loop space of $(X,x_0)$.]\\
\endgroup%%------------------------------------<<

\index{Stacking Lemma}
%\textbf{\small STACKING LEMMA} \quadx
\textbf{\small STACKING LEMMA} \ 
Given a topological space \mY, let $\{P_i:i \in I\}$ be a numerable covering of $IY$ $-$then there exists a numerable covering 
$\{Y_j:j \in J\}$ of \mY and positive real numbers $\epsilon_j$ $(j \in J)$ such that $\forall \ t^\prime, \ t\pp \in [0,1]$ with 
$t^\prime \leq t\pp$ $\&$ $t\pp - t^\prime < \epsilon_j$, $\exists \ i \in I$: $Y_j \times [t^\prime,t\pp] \subset P_i$.\\
\label{4.12}

[Let $\{\rho_i:i \in I\}$ be a partition of unity on $IY$ subordinate to $\{P_i:i \in I\}$.  
Put 
$J = \ds\bigcup\limits_1^\infty I^r$.  Take $j \in J$, say, $j = (i_1, \ldots, i_r) \in I^r$, define $\pi_j \in C(Y,[0,1])$ by 
\[
\pi_j(y) = \prod\limits_{k=1}^r \min\left\{\rho_{i_k}(y,t) \ :\ t \in \left[\frac{k-1}{r+1},\frac{k+1}{r+1}\right] \right\}
\]
and set $Y_j = \pi_j^{-1}(]0,1])$, $\epsilon_j = 1/2r$.  Since 
$Y_j \subset \ds\bigcap\limits_{k=1}^r \left\{y:\{y\} \times 
\left[\ds\frac{k-1}{r+1},\frac{k+1}{r+1}\right] \subset P_{i_k}\right\}$, 
the $\epsilon_j$ will work.  
Moreover, due to the compactness of $[0,1]$, for each $y \in Y$ there is: 
(1) An index $j \in I^r$ such that  
$\{y\} \times \ds\left[\ds\frac{k-1}{r+1},\frac{k+1}{r+1}\right] \subset \rho_{i_k}^{-1}(]0,1])$ 
$(k = 1, \ldots, r)$ and 
(2) A neighborhood $V$ of $y$ such that $IV$ meets but a finite number of the $\rho_i^{-1}(]0,1])$.  
Therefore $\{Y_j:j \in J\} = \ds\bigcup\limits_1^\infty \{Y_j:j \in I^r\}$ is a $\sigma$-neighborhood finite cozero set covering 
of \mY, hence is numerable.]\\

\index{Local-Global Principle}
%\textbf{\small LOCAL-GLOBAL PRINCIPLE} \quadx
\textbf{\small LOCAL-GLOBAL PRINCIPLE} \ 
Let \mX be in \bTOP/\mB.  
Suppose that $\sO = \{O_i:i \in I\}$ is a numerable covering of \mB such that $\forall \ i$, 
$p_{O_i}:X_{O_i} \ra O_i$ is a Hurewicz fibration $-$then $p:X \ra B$ is a Hurewicz fibration.

[Fix a topological space \mY and a pair $(F,h)$ of continuous functions 
$
\begin{cases}
\ F:Y \ra X\\
\ h:IY \ra B
\end{cases}
$
such that $p \circx F = h \circx i_0$.  
To establish the existence of an $H:IY \ra X$ such that 
$F = H \circx i_0$ and 
$p \circx H = h$ 
is equivalent to proving that $\sec_Y(W) \neq \emptyset$ (cf. p. \pageref{4.7}).  
For this, we shall use the section extension theorem and show that $W$ has the SEP, which suffices.  
Set $P_i = h^{-1}(O_i)$: $\{P_i:i \in I\}$ is a numerable covering of $IY$ and the stacking lemma is applicable.  
Given $j$, put $W_j = \restr{W}{Y_j}$, choose $t_k: 0 = t_0 < t_1 < \cdots < t_n = 1$, $t_k - t_{k-1} < \epsilon_j$, and select $i$ 
accordingly: $h(Y_j \times [t_{k-1},t_k]) \subset O_i$. \   
The claim is that $W_j$ has the SEP.  
So let $A \subset Y_j$, let $V$ 
be a halo of \mA in $Y_j$, and let $H_V:IV \ra X$ be a homotopy such that $\restr{F}{V} = H_V \circx i_0$ and 
$p \circx H_V = h_{IV}$.  
With $\pi$ a haloing function of $V$, put $A_k = \pi^{-1}([t_k,1])$ $(k = 1, \ldots,n)$: $A_k$ is 
a halo of $A_{k+1}$ in $Y_j$ and $V$ is a halo 
%%----------------------------------------------------------------------------------------------11
of $A_1$ in $Y_j$.  
Owing to Proposition 2, there exist homotopies $H_k:Y_j \times [t_{k-1},t_k] \ra X$ having the following properties: 
$p \circx H_k = \restr{h}{Y_j \times [t_{k-1},t_k]}$, 
$H_k(y,t_{k-1}) = H_{k-1}(y,t_{k-1})$ $(k > 1)$, 
$H_1(y,0) = F(y)$, 
$\restr{H_k}{A_k \times [t_{k-1},t_k]} = \restr{H_V}{A_k \times [t_{k-1},t_k]}$.  
The $H_k$ thus combine to determine a homotopy 
$H:IY_j \ra X$ such that 
$\restr{F}{Y_j} = H \circx i_0$, 
$p \circx H = \restr{h}{IY_j}$, and 
$\restr{H}{IA} = \restr{H_V}{IA}$.]\\

Application: Suppose that \mB is a paracompact Hausdorff space.  Let \mX be in $\bTOP/B$.  Assume: \mX is locally 
trivial $-$then $p:X \ra B$ is a Hurewicz fibration.\\

\begingroup%%----------------------------------->>
\fontsize{9pt}{11pt}\selectfont
\textbf{\small EXAMPLE}  \ 
Let $B = L^+$, the long ray.  Put $X = \{(x,y) \in L^+ \times L^+: x < y\}$ and let $p$ be the vertical projection $-$then 
\mX is locally trivial but $p:X \ra B$ is not a Hurewicz fibration.\\
\endgroup%%------------------------------------<<

\label{4.66}
\label{13.7}
\textbf{\small FACT}  \ 
\begingroup%%----------------------------------->>
\fontsize{9pt}{11pt}\selectfont
Let \mX be in \bTOP/\mB.  
Suppose that $\sO = \{O_i:i \in I\}$ is an open covering of \mB such that $\forall \ i$, 
$p_{O_i}:X_{O_i} \ra O_i$ is a Hurewicz fibration $-$then the projection $p:X \ra B$ is a $\sY$ fibration, where $\sY$ is 
the class of paracompact Hausdorff spaces.
 \\ \indent
[Given $Y \in \sY$ and continuous functions 
 $
\begin{cases}
\ F:Y \ra X\\
\ h:IY \ra B
\end{cases}
$
such that $p \circx F = h \circx i_0$, consider the pullback square 
\begin{tikzcd}[sep=large]
{IY \times_B X} \ar{d} \ar{r} &{X} \ar{d}{p}\\
{IY} \ar{r}[swap]{h} &{B}
\end{tikzcd}
, observing that $IY \in \sY$.]
 \\ \indent
[Note: \ It follows that $p:X \ra B$ is a Serre fibration.]\\
\endgroup%%------------------------------------<<

Let $f:X \ra Y$ be a continuous function $-$then the 
\un{mapping track}
\index{mapping track} 
$W_f$ of $f$ is defined 
by the pullback square 
\begin{tikzcd}[sep=large]
{W_f} \ar{d} \ar{r} &{P Y} \ar{d}{p_0}\\
{X} \ar{r}[swap]{f} &{Y}
\end{tikzcd}
.  \ 
Special case: $\forall \ y_0 \in Y$, the mapping track of the inclusion $\{y_0\} \ra Y$ is the mapping space $\Theta Y$ of 
$(Y,y_0)$.  
There is a projection $p:W_f \ra X$, a homotopy $G:W_f \ra PY$, and a unique continuous function $s:X \ra W_f$ such that $p \circx s = \id_X$ and $G \circx s = j \circx f$ $(j:Y \ra PY)$.  
One has $s \circx p \underset{X}{\simeq} \id_{W_f}$.  
The composition $p_1 \circx G$ is a projection $q:W_f \ra Y$ and $f = q \circx s$.

[Note: \ The mapping track is a functor $\bTOP(\ra) \ra \bTOP$.]\\

\textbf{\small LEMMA}  \ 
$p$ is a Hurewicz fibration and $W_f$ is fiberwise contractible over \mX.\\

\textbf{\small LEMMA}  \ 
$q$ is a Hurewicz fibration.\\

%%----------------------------------------------------------------------------------------------12
[To construct a filler for 
\begin{tikzcd}[sep=large]
{E} \ar{d}[swap]{i_0} \ar{r}{\Phi} &{W_f} \ar{d}{q}\\
{IE} \ar{r}[swap]{h} &{Y}
\end{tikzcd}
, write 
$
\Phi(e) = (x_e, \tau_e):
\begin{cases}
\ x_e \in X\\
\ \tau_e \in PY
\end{cases}
$
$\&$ $f(x_e) = \tau_e(0)$, and define $H:IE \ra W_f$ by $H(e,t) = (x_e,\ov{h}(e,t))$, where
\[
\ov{h}(e,t)(T) = 
\begin{cases}
\ \tau_e(2T(2 - t)^{-1}) \hspace{0.75cm}  (T \leq 1 - t/2)\\
\ h(e,2T + t - 2) \hspace{0.75cm} (T \geq 1 - t/2)
\end{cases}
.]
\]
\vspace{0.1cm}

\begin{proposition} \  %07
Every morphism in \bTOP can be written as the composite of a homotopy equivalence and a Hurewicz fibration.\\
\end{proposition}

\label{12.5}
\begingroup%%----------------------------------->>
\fontsize{9pt}{11pt}\selectfont
\label{4.27}
\textbf{\small FACT}  \ 
Let $f:X \ra Y$ be a continuous function $-$then $f$ can be factored as 
$
f = 
\begin{cases}
\ \Phi \circx k\\
\ \Psi \circx l
\end{cases}
$
, where 
$
\begin{cases}
\ \Phi\\
\ \Psi
\end{cases}
$
is a Hurewicz fibration, 
$
\begin{cases}
\ k\\
\ l
\end{cases}
$
is a closed cofibration, and 
$
\begin{cases}
\ k\\
\ \Psi
\end{cases}
$
is a homotopy equivalence.
\\ \indent
[Per Proposition 7, write $f = q \circx s$, form $S = Is(X) \cup W_f \times ]0,1] \subset IW_f$, and let 
$\omega:IW_f \ra [0,1]$ be the projection.  The restriction to \mS of the Hurewicz fibration $IW_f \ra W_f$ is a Hurewicz fibration, call it $p$.  
Proof: Given continous functions
$
\begin{cases}
\ F:Y \ra S\\
\ h:IY \ra W_f
\end{cases}
$
such that $p \circx F = h \circx i_0$, consider $H:IY \ra S$, where $H(y,t) = (h(y,t),t + (1 - t)\omega(F(y)))$.  
Next, if $k:X \ra S$ is defined by $k(x) = (s(x),0)$, 
then $k(X)$ is both a strong deformation retract of \mS and a zero set in \mS (being $(\restr{\omega}{S})^{-1}(0))$.  
Therefore $k$ is a closed cofibration (cf. $\S 3$, Proposition 10).  And: $f = q \circx p \circx k$.  
To derive the other factorization, write $f = r \circx i$ (cf. $\S 3$, Proposition 16) and decompose $r$ as above.]\\
\endgroup%%------------------------------------<<

Let \mX be in $\bTOP/B$.  Define $\lambda:PX \ra W_p$ by $\sigma \ra (\sigma(0),p \circx \sigma)$.\\

\begin{proposition} \  %08
The projection $p:X \ra B$ is a Hurewicz fibration iff $\lambda$ has a right inverse $\Lambda$.
\end{proposition}

[Note: \ $\Lambda$ is called a 
\un{lifting function}.]
\index{lifting function}\\

\begingroup%%----------------------------------->>
\fontsize{9pt}{11pt}\selectfont
\label{4.19}
\textbf{\small FACT}  \ 
Let $p:X \ra B$  be a Hurewicz fibration.  Suppose that \mA is a subspace of \mX for which there exists a fiber preserving retraction $r:X \ra A$ $-$then the restriction of $p$ to \mA is a Hurewicz fibration $A \ra B$.\\
\endgroup%%------------------------------------<<

\begingroup%%----------------------------------->>
\fontsize{9pt}{11pt}\selectfont
\textbf{\small EXAMPLE}  \ 
Let \mX be a nonempty compact subspace of $\R^n$.  Realize $\Gamma X$ in $\R^{n+1}$ by writing 
$\Gamma X = \ds\bigcup\limits_x \{(t,tx):0 \leq t \leq 1\}$, so $\Gamma^2 X$ is 
$\ds\bigcup\limits_x \{(s,st,stx): 0 \leq s \leq 1 \ \& \ 0 \leq t \leq 1\}$, a subspace of $\R^{n+2}$.  
Claim: The projection 
$
p: 
\begin{cases}
\ \Gamma^2 X \ra [0,1]\\
\ (s,st,stx) \ra s
\end{cases}
$
is a Hurewicz fibration.  To see this, consider 
$[0,1] \times \Gamma X = \ds\bigcup\limits_x \{(s,t,tx): 0 \leq s \leq 1 \ \& \ 0 \leq t \leq 1\}$ with projection 
$(s,t,tx) \ra s$ and define a fiber preserving retraction 
%%----------------------------------------------------------------------------------------------13
$r:[0,1] \times \Gamma X \ra \Gamma^2 X$ by 
$
r(s,t,tx) = 
\begin{cases}
\ (s,s,sx) \quadx (t \geq s)\\
\ (s,t,tx) \quadx \ (t \leq s)
\end{cases}
\hspace{-.2cm}
. \ 
$
The fibers of $p$ over the points in $]0,1]$ can be identified with $\Gamma X$, while $p^{-1}(0) = *$.
\\ \indent
[Note: \ If \mX is the Cantor set, then $\Gamma X$ is not an ANR.]\\
\endgroup%%------------------------------------<<

Let \mX be in $\bTOP/B$ $-$then there is a morphism 
\begin{tikzcd}%[sep=large]
{X} \ar{rd}[swap]{p} \ar{rr}{\gamma} &&{W_p} \ar{ld}{q}\\
&{B}
\end{tikzcd}
.  
Here, in a change of notation, $\gamma$ sends $x$ to $(x,j(p(x)))$, $j:B \ra PB$ the embedding.\\

\begin{proposition} \  %09
Suppose that $p:X \ra B$ is a Hurewicz fibration $-$then $\gamma:X \ra W_p$ is a fiber homotopy equivalence.
\end{proposition}

[Choose a lifting function $\Lambda:W_p: \ra PX$.  Define a fiber homotopy 
$H:IX \ra X$ by $H(x,t) = \Lambda(\gamma(x))(t)$ and a fiber homotopy $G:IW_p \ra W_p$ by 
$G((x,\tau),t) = (\Lambda(x,\tau)(t),\tau_t)$ $(\tau_t(T) = \tau(t + T - tT))$ $-$then it is clear that the assignment 
$(x,\tau) \ra \Lambda(x,\tau)(1)$ is a fiber homotopy inverse for $\gamma$.]\\

\label{4.61}
\label{13.35}
Application: The fibers of a Hurewicz fibration over a path connected base have the same homotopy type.

[Note: \ This need not be true if ``Hurewicz'' is replaced by ``Serre'' (cf. p. \pageref{4.8}).  
It can also fail if ``path connected'' is weakened to ``connected''.  Indeed, for a connected \mB whose path components are singletons, every 
$p:X \ra B$ is a Hurewicz fibration.]\\

\begingroup%%----------------------------------->>
\fontsize{9pt}{11pt}\selectfont
A Hurewicz fibration $p:X \ra B$ is said to be 
\un{regular}
\index{regular (Hurewicz fibration)} 
if the morphism 
\begin{tikzcd}%[sep=large]
{X} \ar{rd}[swap]{p} \ar{rr}{\gamma} &&{W_p} \ar{ld}{q}\\
&{B}
\end{tikzcd}
has a left inverse $\Gamma$ in $\bTOP/B$.\\
\endgroup%%------------------------------------<<

\begingroup%%----------------------------------->>
\fontsize{9pt}{11pt}\selectfont
\textbf{\small FACT}  \ 
The Hurewicz fibration $p:X \ra B$ is regular iff there exists a lifting function $\Lambda_0:W_p \ra PX$ with the property that 
$\Lambda_0(x,\tau) \in j(X)$ whenever $\tau \in j(B)$.
\\ \indent
[Given a left inverse $\Gamma$ for $\gamma$, consider the lifting function $\Lambda_0:W_p \ra PX$ defined by 
$\Lambda_0(x,\tau)(t) = \Gamma(x,\tau_t)$, where $\tau_t(T) = \tau(tT)$.]\\
\endgroup%%------------------------------------<<

\begingroup%%----------------------------------->>
\fontsize{9pt}{11pt}\selectfont
\textbf{\small FACT}  \ 
The Hurewicz fibration $p:X \ra B$ is regular iff every commutative diagram  
\begin{tikzcd}[sep=large] %\ds
{Y} \ar{d}[swap]{i_0} \ar{r}{F} &{X} \ar{d}{p}\\
{IY} \ar{r}[swap]{h} &{B}
\end{tikzcd}
admits a filler $H:IY \ra X$ such that $H$ is stationary with $h$, 
i.e., $\restr{h}{I\{y_0\}}$ constant $\implies$ $\restr{H}{I\{y_0\}}$ constant.
\\ \indent
[Note: \ The local-global principle is valid in the regular situation (work with a suitable subspace of 
\mW to factor in the stationary condition).]\\
\endgroup%%------------------------------------<<

%%----------------------------------------------------------------------------------------------14
\begingroup%%----------------------------------->>
\fontsize{9pt}{11pt}\selectfont
A sufficient condition for the regularity of the Hurewicz fibration $p:X \ra B$ is that $j(B)$ be a zero set in $PB$.  
Thus let $\phi \in C(PB,[0,1])$: $j(B) = \phi^{-1}(0)$.  
Define $\Phi \in C(PB,PB)$ by 
$
\Phi(\tau)(t) =
\begin{cases}
\ \tau(t/\phi(\tau)) \hspace{0.5cm} (t < \phi(\tau))\\
\ \tau(1) \hspace{1.25cm} (\phi(\tau) \leq t \leq 1)
\end{cases}
\hspace{-.2cm}. \ 
$
Take any lifting function $\Lambda$ and put $\Lambda_0(x,\tau)(t) = \Lambda(x,\Phi(\tau))(\phi(\tau)t)$ to get a lifting function $\Lambda_0:W_p \ra PX$ with the property that $\Lambda_0(x,\tau) \in j(X)$ whenever $\tau \in j(B)$.  
Example: $j(B)$ is a zero set in $PB$ if $\Delta_B$ is a zero set in $B \times B$, e.g., if the inclusion 
$\Delta_B \ra B \times B$ is a closed cofibration, a condition satisfied by a CW complex or a metrizable topological manifold 
(cf. p. \pageref{4.9}).\\
\endgroup%%------------------------------------<<

\begingroup%%----------------------------------->>
\fontsize{9pt}{11pt}\selectfont
\textbf{\small EXAMPLE}  \ 
Let $B = [0,1]/[0,1[$ $-$then the Hurewicz fibration $p_0:PB \ra B$ is not regular.\\
\endgroup%%------------------------------------<<

\label{4.35}
\begingroup%%----------------------------------->>
\fontsize{9pt}{11pt}\selectfont
\textbf{\small FACT}  \ 
Suppose that $p:X \ra B$ is a regular Hurewicz fibration $-$then $\forall \ x_0 \in X$, $p:(X,x_0) \ra (B,b_0)$ is a pointed 
Hurewicz fibration ($b_0 = p(x_0)$).\\
\endgroup%%------------------------------------<<

\label{4.64}
\label{6.18}
\begingroup%%----------------------------------->>
\fontsize{9pt}{11pt}\selectfont
Let \mX be in $\bTOP/B$ $-$then the projection $p:X \ra B$ is said to have the 
\un{slicing structure property}
\index{slicing structure property} 
if there exists an open covering $\sO = \{O_i:i \in I\}$ of \mB 
and continuous functions $s_i:O_i \times X_{O_i} \ra X_{O_i}$ $(i \in I)$ such that $s_i(p(x),x) = x$ and $p \circx s_i(b,x) = b$.  Note that $p$ is necessarily open.  
Example: \mX locally trivial $\implies$ $p:X \ra B$ has the slicing structure property (but not conversely).
\\ \indent
Observation: Suppose that $p:X \ra B$ has the slicing structure property $-$then $\forall \ i$, 
$p_{O_i}:X_{O_i} \ra O_i$ is a regular Hurewicz fibration.
\\ \indent
[Consider the lifting function $\Lambda_i$ defined by $\Lambda_i(x,\tau)(t) = s_i(\tau(t),x)$.]
\\ \indent
So, if $p:X \ra B$ has the slicing structure property, then $p:X \ra B$ must be a Serre fibration and is even a regular 
Hurewicz fibration provided that \mB is a paracompact Hausdorff space.\\
\endgroup%%------------------------------------<<

\begingroup%%----------------------------------->>
\fontsize{9pt}{11pt}\selectfont
\textbf{\small FACT}  \ 
Let \mX be in $\bTOP/B$, where \mB is uniformly locally contractible.  Assume: The projection $p:X \ra B$ is a regular Hurewicz fibration $-$then $p$ has the slicing structure property.\\
\endgroup%%------------------------------------<<

\label{6.24}
\begingroup%%----------------------------------->>
\fontsize{9pt}{11pt}\selectfont
Application: Suppose that \mB is a uniformly locally contractible paracompact Hausdorff space.  Let \mX be in 
$\bTOP/B$ $-$then the projection $p:X \ra B$ is a regular Hurewicz fibration iff $p$ has the slicing structure property.
\\ \indent
[Note: \ It therefore follows that if \mB is a CW complex or a metrizable topological manifold, then the Hurewicz fibrations with base \mB are precisely the $p:X \ra B$ which have the slicing structure property.]\\
\endgroup%%------------------------------------<<

\begingroup%%----------------------------------->>
\fontsize{9pt}{11pt}\selectfont
\textbf{\small FACT}  \ 
Let $p:X \ra B$ be a Serre fibration, where \mX and \mB are CW complexes $-$then $p$ is a \bCG fibration.
\\ \indent
[An open subset of a CW complex is homeomorphic to a retract of a CW complex (cf. p. \pageref{4.10}).]
\\ \indent
[Note: \ If $X \times B$ is compactly generated, the $p$ is a Hurewicz fibration.]\\
\endgroup%%------------------------------------<<

Cofibrations are embeddings (cf. p. \pageref{4.11}).  By analogy, one might expect that surjective Hurewicz fibrations are quotient maps.  However, this is not true in general.  Example:
%%----------------------------------------------------------------------------------------------15
Take $X = \Q$ (discrete topology), $B = \Q$ (usual topology), $p = \id_{\Q}$ $-$then $p:X \ra B$ is a surjective Hurewicz fibration but is not a quotient map.\\

\begin{proposition} \  %10
Let $p:X \ra B$ be a Hurewicz fibration.  
Assume: $p$ is surjective and $B$ is locally path connected $-$then $p$ is a quotient map.
\end{proposition}

[Consider the \cd
\begin{tikzcd}[sep=large]
{PX} \ar{d}[swap]{p_1} \ar{r}{\lambda} &{W_p} \ar{d}{q}\\
{X} \ar{r}[swap]{p} &{B}
\end{tikzcd}
.  Since $\lambda$ and $p_1$ have right inverses, they are quotient, so $p$ is quotient iff $q$ is quotient.  Take a nonempty subset 
$O \subset B$: $W_O$ is open in $W_p$.  Fix $b \in O$, $x \in X_b$, and choose a neighborhood $O_b$ of 
$b: (\{x\} \times PO_b) \cap W_p \subset W_O$.  The path component $O_0$ of $O_b$ containing $b$ is open.  
Given $b_0 \in O_0$, $\exists \ \tau \in PO_b$ connecting $b$ and $b_0$.  
But $(x,\tau) \in W_O$ $\implies$ 
$b_0 = q(x,\tau) \in O$ $\implies$ $O_0 \subset O$.  
Therefore $O$ is open in \mB, hence $q$ is quotient.]\\

Application: Every connected locally path connected nonempty space \mB is the quotient of a contractible space.

[Fix $b_0 \in B$ and consider the mapping space $\Theta B$ of $(B,b_0)$ with projection $\tau \ra \tau(1)$.]\\

Let $p:X \ra B$ be a Hurewicz fibration $-$then for any path component \mA of \mX, $p(A)$ is a path component of \mB and $A \ra p(A)$ is a Hurewicz fibration.  
Therefore $p(X)$ is a union of path components of \mB.  
So, if \mB is path connected and \mX is nonempty, then $p$ is surjective.\\

\begingroup%%----------------------------------->>
\fontsize{9pt}{11pt}\selectfont
\textbf{\small FACT}  \ 
Let $p:X \ra B$ be a Hurewicz fibration.  Assume: $B$ is path connected and $X_b$ is path connected for some $b \in B$ $-$then 
\mX is is path connected.
\\ \indent
[Note: \ The fibers of a Hurewicz fibration $p:X \ra B$ need not be path connected but if \mX is path connected, 
then any two path components of a given fiber have the same homotopy type.]\\
\endgroup%%------------------------------------<<

\label{4.31}
\begingroup%%----------------------------------->>
\fontsize{9pt}{11pt}\selectfont
\textbf{\small FACT}  \ 
Suppose that \mB is path connected $-$then $B$ is locally path connected iff every Hurewicz fibration $p:X \ra B$ is open.\\
\endgroup%%------------------------------------<<

\begin{proposition} \ %11
Let $p:X \ra B$ be a Hurewicz fibration.  Suppose that the inclusion $O \ra B$ is a closed cofibration $-$then the inclusion 
$X_O \ra X$ is a closed cofibration.
\end{proposition}

[Fix a Str{\o}m structure $(\phi,\Phi)$ on $(B,O)$.  Let $H:IX \ra X$ be a filler for the commutative diagram 
\begin{tikzcd}[sep=large]
{X} \ar{d}[swap]{i_0} \ar{r}{\id_X} &{X} \ar{d}{p}\\
{IX} \ar{r}[swap]{h} &{B}
\end{tikzcd}
, where $h = \Phi \circx Ip$.  Define a Str{\o}m structure $(\psi,\Psi)$ on $(X,X_O)$ by $\psi = \phi \circx p$, 
$\Psi(x,t) = H(x, \min\{t,\psi(x)\})$.]\\

%%----------------------------------------------------------------------------------------------16
\label{4.14}
Application: Let $p:X \ra B$ be a Hurewicz fibration.  
Let \mA be a subspace of \mX and suppose that the inclusion 
$A \ra X$ is a closed cofibration.  
View \mA as an object in \bTOP/\mB with projection $p_A = \restr{p}{A}$ $-$then the inclusion 
$W_{p_A} \ra W_p$ is a closed cofibration.\\

\begingroup%%----------------------------------->>
\fontsize{9pt}{11pt}\selectfont
\label{4.36}
\textbf{\small EXAMPLE}  \ \ 
Let $(X,x_0)$ be a pointed space.  \ 
Assume: The inclusion $\{x_0\} \ra X$ is a closed cofibration $-$then Proposition 11 implies 
that the inclusion $j:\Omega X \ra \Theta X$ is a closed cofibration.  
Call $\theta$ the continuous function 
$\Gamma \Omega X \ra \Theta X$ that sends $[\sigma,t]$ to $\sigma_t$, where $\sigma_t(T) = \sigma(tT)$.  
The arrow i: 
$
\begin{cases}
\ \Omega X \ra \Gamma \Omega X\\
\ \sigma \ra [\sigma,1]
\end{cases}
$ 
is a closed cofibration and $\theta \circx i = j$.  Consider the commutative diagram
\begin{tikzcd}[sep=large]
%{\Omega X} \ar[equals]{d} \ar{r}{i} 
{\Omega X} \arrow[d,shift right=0.5,dash] \arrow[d,shift right=-0.5,dash]  \ar{r}{i}
&{\Gamma \Omega X} \ar{d}{\theta}\\
{\Omega X} \ar{r}[swap]{j} &{\Theta X}
\end{tikzcd}
.  \ Because $\Gamma \Omega X$ and $\Theta X$ are contractible, it follows from $\S 3$, Proposition 14 that the arrow 
$(\id_{\Omega X},\theta)$ is a homotopy equivalence in $\bTOP(\ra)$.\\
\endgroup%%------------------------------------<<

\textbf{\small LEMMA}  \ 
Let $\phi \in C(Y,[0,1])$ :such that $A = \phi^{-1}(0)$ is a strong deformation retract of \mY.  Suppose that $p:X \ra B$ is a 
Hurewicz fibration $-$then every commutative diagram 
\begin{tikzcd}[sep=large]
{A} \ar{d}[swap]{i} \ar{r}{g} &{X} \ar{d}{p}\\
{Y} \ar{r}[swap]{f} &{B}
\end{tikzcd}
\hspace{0.5cm}
has a filler $F:Y \ra X$.

[Fix a retraction $r:Y \ra A$ and a homotopy $\Phi:IY \ra Y$ between $i \circx r$ and $\id_Y$ rel \mA.  
Define a homotopy $h:IY \ra Y$ by 
$
h(y,t) = 
\begin{cases}
\ \Phi(y,t/\phi(y)) \hspace{0.5cm} (t < \phi(y))\\
\ \Phi(y,1) \hspace{1.25cm} \ (t \geq \phi(y))
\end{cases}
\hspace{-.25cm}.
$
Since $p$ is a Hurewicz fibration, there exists a homotopy $H:IY \ra X$ such that $g \circx r = H \circx i_0$ and 
$p \circx H = f \circx h$.  
Take for $F:Y \ra X$ the continuous function $y \ra H(y,\phi(y))$.]

[Note: \ The hypotheses on \mA are realized when the inclusion $i:A \ra Y$ is both a homotopy equivalence and a closed cofibration (cf. $\S 3$, Proposition 5).]\\

\label{12.1}
\label{12.6}
\begingroup%%----------------------------------->>
\fontsize{9pt}{11pt}\selectfont
\textbf{\small FACT}  \ 
Let $i:A \ra Y$ be a continuous function with a closed image $-$then $i$ is both a homotopy equivalence and a closed cofibration 
iff every commutative diagram 
\begin{tikzcd}[sep=large]
{A} \ar{d}[swap]{i} \ar{r} &{X} \ar{d}{p}\\
{Y} \ar{r} &{B}
\end{tikzcd}
, where $p$ is a Hurewicz fibration, has a filler $Y \ra X$.
\\ \indent
[First take $X = PB$, $p = p_0$ to see that $i$ is a closed cofibration.  
Next, identify \mA with $i(A)$ and produce a 
retraction $r:Y \ra A$ from a filler for 
\begin{tikzcd}[sep=large]
{A} \ar{d}[swap]{i} \ar{r}{\id_A} &{A} \ar{d}\\
{Y} \ar{r} &{*}
\end{tikzcd}
.  Finally, consider 
\begin{tikzcd}[sep=large]
{A} \ar{d}[swap]{i} \ar{r}{j} &{PY} \ar{d}{\Pi}\\
{Y} \ar{r}[swap]{\rho} &{Y \times Y}
\end{tikzcd}
where $\rho(y) = (y,r(y))$ ($\Pi$ as on p. \pageref{4.12}).]\\
\endgroup%%------------------------------------<<

%%----------------------------------------------------------------------------------------------17
\begingroup%%----------------------------------->>
\fontsize{9pt}{11pt}\selectfont
\textbf{\small FACT}  \ 
Let $p:X \ra B$ be a continuous function $-$then $p$ is a Hurewicz fibration iff every commutative diagram
\begin{tikzcd}[sep=large]
{A} \ar{d}[swap]{i} \ar{r} &{X} \ar{d}{p}\\
{Y} \ar{r} &{B}
\end{tikzcd}
, where $i$ is both a homotopy equivalence and a closed cofibration, has a filler $Y \ra X$.\\
\endgroup%%------------------------------------<<

\label{5.11}
\label{12.2}
\label{12.7}
\begingroup%%----------------------------------->>
\fontsize{9pt}{11pt}\selectfont
\textbf{\small FACT}  \ 
Let 
\begin{tikzcd}[sep=large]
{X_0} \ar{d} &{X_1} \ar{d} \ar{l} &{\cdots} \ar{l}\\
{Y_0} &{Y_1} \ar{l} &{\cdots} \ar{l}
\end{tikzcd}
be a commutative ladder of topological spaces.  Assume: $\forall \ n$, the horizontal arrows 
$
\begin{cases}
\ X_n \la X_{n+1}\\
\ Y_n \la Y_{n+1}
\end{cases}
$
are Hurewicz fibrations and the vertical arrows $\phi_n:X_n \ra Y_n$ are homotopy equivalences $-$then the induced map 
$\phi:\lim X_n \ra \lim Y_n$ is a homotopy equivalence.
\\ \indent
[The mapping cylinder is a functor $\bTOP(\ra) \ra \bTOP$, so there is an arrow 
$\pi_n:M_{\phi_{n+1}} \ra M_{\phi_{n}}$.  Use $\S 3$, Proposition 17 to construct a commutative triangle 
\begin{tikzcd}[sep=large]
{X_{0}} \ar{d}[swap]{i} \ar{r}{\id} &{X_{0}}\\
{M_{\phi_{0}}} \ar{ru}[swap]{r_0}
\end{tikzcd}
.  \ The lemma then provides a filler $r_1:M_{\phi_{1}} \ra X_1$ for 
\begin{tikzcd}[sep=large]
{X_{1}} \ar{d}[swap]{i} \ar{r}{\id} &{X_{1}} \ar{d}\\
{M_{\phi_{1}}} \ar{r}[swap]{r_0 \circx \pi_0} &{X_0}
\end{tikzcd}
, \ hence, by induction, a filler $f_{n+1}:M_{\phi_{n+1}} \ra X_{n+1}$ for 
\begin{tikzcd}[sep=large]
{X_{n+1}} \ar{d}[swap]{i} \ar{r}{\id} &{X_{n+1}} \ar{d}\\
{M_{\phi_{n+1}}} \ar{r}[swap]{r_n \circx \pi_n} &{X_n}
\end{tikzcd}
.  \ 
Give the composite $Y_n \overset{j}{\ra} M_{\phi_{n}} \overset{r_n}{\lra} X_n$ a name, say $\psi_n$, and take limits to get 
a left homotopy inverse $\psi$ for $\phi$.]\\
\endgroup%%------------------------------------<<

\begin{proposition} \  %12
Let \mA be a closed subspace of \mY and assume that the inclusion $A \ra Y$ is a cofibration.  Suppose that $p:X \ra B$ is a 
Hurewicz fibration $-$then every commutative diagram 
\begin{tikzcd}[sep=large]
{i_0Y \cup IA} \ar{d} \ar{r}{F} &{X} \ar{d}{p}\\
{IY} \ar{r}[swap]{h} &{B}
\end{tikzcd}
\ 
has a filler $H:IY \ra X$.
\end{proposition}

[Quote the lemma: $i_0Y \cup IA$ is a strong deformation retract of $IY$ (cf. p. \pageref{4.13}) and $i_0Y \cup IA$ is a zero set in $IY$.]\\

Application: Let  $p:X \ra B$ be a Hurewicz fibration, where \mB is a LCH space.  Suppose that the inclusion $O \ra B$ 
is a closed cofibration $-$then the arrow of restriction $\sec_B(X) \ra \sec_O(X_O)$ is a Hurewicz fibration.\\

%%----------------------------------------------------------------------------------------------18
\begingroup%%----------------------------------->>
\fontsize{9pt}{11pt}\selectfont
\index{Vertical Homotopies (example)}
\textbf{\small EXAMPLE  \  (\un{Vertical Homotopies})} \ 
Let $p:X \ra B$ be a Hurewicz fibration.  Suppose that $s^\prime$, $s\pp \in \sec_B(X)$ are homotopic $-$then 
$s^\prime$, $s\pp$ are vertically homotopic.
\\ \indent
[Take any homotopy \ $H:IB \ra X$ between $s^\prime$ and $s\pp$.  \ 
Define $G:IB \ra X$ by 
$G(b,t) \ = \  $
$
\begin{cases}
\ H(b,2t) \hspace{2.35cm} (0 \leq t \leq 1/2) \\
\ s\pp \circx p \circx H(b,2 - 2t) \hspace{0.35cm} \ \ (1/2 \leq t \leq 1)
\end{cases}
$
\hspace{-.2cm}.  \ 
Since $p \circx G(b,t) = p \circx G(b,1 - t)$, it follows that $p \circx G$ is homotopic rel $B \times \{0,1\}$ to the projection 
$B \times [0,1] \ra B$.]\\
\endgroup%%------------------------------------<<

\textbf{\small LEMMA}  \ 
Let \mA be a closed subspace of \mY and assume that the inclusion $A \ra Y$ is a cofibration.  
Suppose that $p:X \ra B$ is a Hurewicz fibration.  Let $F:i_0 Y \cup IA \ra X$ be a continuous function such that $\forall \ a \in A$: 
$p \circx F(a,t) = p \circx F(a,0)$ $(0 \leq t \leq 1)$ $-$then there exists a continuous function $H:IY \ra X$ which extends $F$ such that $\forall \ y \in Y$: $p \circx H(y,t) = p \circx H(y,0)$ $(0 \leq t \leq 1)$.

[Choose $\phi \in C(Y,[0,1])$: $A = \phi^{-1}(0)$ and fix a retraction $r:IY \ra i_0 Y \cup IA$.  
Put $f = p \circx F \circx r$.  Define $G \in C(IY,PB)$ as follows: 
%G(y,t)(T) = 
%(i) 
%(ii) 
%(iii) 
%(iv) $f(y,0)$ $(t\phi(y) \leq T \leq 1
\[
G(y,t)(T) =
\begin{cases}
\ f(y,(t\phi(y) - T(2 - \phi(y)))/\phi(y)) \hspace{0.3cm}  (0 \leq T \leq t\phi(y)/2  \ \& \ \phi(y) \neq 0)\\
\ f(y,t) \hspace{4.7cm} (0 \leq T \leq t\phi(y)/2 \ \& \ \phi(y) = 0)\\
\ f(y,t\phi(y) - T) \hspace{3.22cm}  (t\phi(y)/2 \leq T \leq t\phi)y))\\
\ f(y,0) \hspace{4.65cm} (t\phi(y) \leq T \leq 1)
\end{cases}
\hspace{-.2cm}.
\]
Take a lifting function $\Lambda:W_p \ra PX$ and set $H(y,t) = \Lambda(F\circx r(y,t)$, $G(y,t))(t\phi(y))$.]\\

\index{Lifting Principle}
%\textbf{\small LIFTING PRINCIPLE} \quadx
\textbf{\small LIFTING PRINCIPLE} \ 
Let $p:X \ra B$ be a Hurewicz fibration.  
Let \mA be a subspace of \mX and suppose that the inclusion $A \ra X$ is a closed cofibration.  
View \mA as an object in \bTOP/\mB with projection $p_A = \restr{p}{A}$ and assume that $p_A:A \ra B$ is a 
Hurewicz fibration.  Let $\Lambda_A:W_{p_A} \ra PA$ be a lifting function $-$then there exists a lifting function 
$\Lambda_X:W_p \ra PX$ such that $\Lambda_\restr{X}{W_{p_A}} = \Lambda_A$.

[The inclusion $W_{p_A} \ra W_p$ is a closed cofibration (cf. p. \pageref{4.14}).  Therefore the inclusion 
$i_0W_p \cup IW_{p_A} \ra IW_p$ is a closed cofibration (cf. p. \pageref{4.15} or $\S 3$, Proposition 7).  
Fix a lifting function 
$\Lambda:W_p \ra PX$. \ 
Define a continuous function 
$F:i_0 IW_p \cup I(i_0W_p \cup IW_{p_A}) \ra X$ by 
\[
F((x,\tau),t,T) = 
\begin{cases}
\ \Lambda(x,\tau)(t) \hspace{3.55cm} \  (T = 0 \ \& \ (x,\tau) \in W_p)\\
\ x \hspace{5.05cm} (t = 0 \ \& \ (x,\tau) \in W_p)\\
\ \Lambda_A(a,\tau)(t) \hspace{3.4cm}  (0 \leq t \leq T \ \& \ (a,\tau) \in W_{p_A})\\
\ \Lambda(\Lambda_A(a,\tau)(T),\tau * T)(t - T) \hspace{0.5cm} (T \leq t \leq 1 \ \& \ (a,\tau) \in W_{p_A})
\end{cases}
\hspace{-.2cm}.
\]
Here, 
$
\tau * T(t)  =  
\begin{cases}
\ \tau(t + T) \hspace{0.5cm} (t \leq 1 - T)\\
\ \tau(1) \hspace{1.2cm}  (t \geq 1 - T)
\end{cases}
\hspace{-.2cm} .
$
Apply the lemma to get a continuous function $H:I^2W_p \ra X$ which extends $F$ such that 
$\forall \ ((x,\tau),t) \in IW_p$: 
$p \circx H((x,\tau), t, T) = p \circx H((x,\tau), t, 0)$.  
Put $\Lambda_X(x,\tau)(t) = H((x,\tau),t,1)$ $-$then 
$\Lambda_X:W_p \ra PX$ is a lifting function that restricts to $\Lambda_A$.]\\

\begin{proposition} \  %13
Let \mX be in \bTOP/\mB.  Suppose that $X = A_1 \cup A_2$, where 
$
\begin{cases}
\ A_1\\
\ A_2
\end{cases}
$
are closed and the inclusions 
$A_0 = A_1 \cap A_2 \ra 
\begin{cases}
\ A_1\\
\ A_2
\end{cases}
$
are cofibrations.  Assume: 
$
\begin{cases}
\ p_1 = p_{A_1}:A_1 \ra B\\
\ p_2 = p_{A_2}:A_2 \ra B
\end{cases}
$ 
$\&$ $p_0 = p_{A_0}:A_0 \ra B$ are Hurewicz fibrations $-$then $p:X \ra B$ is
%%----------------------------------------------------------------------------------------------19
a Hurewicz fibration.
\end{proposition}

[Choose a lifting function $\Lambda_0:W_{p_0} \ra PA_0$.  Use the lifting principle to secure lifting functions 
$
\begin{cases}
\ \Lambda_1:W_{p_1} \ra PA_1\\
\ \Lambda_2:W_{p_2} \ra PA_2
\end{cases}
$
such that 
$
\begin{cases}
\ \restr{\Lambda_1}{W_{p_0}} = \Lambda_0\\
\ \restr{\Lambda_2}{W_{p_0}} = \Lambda_0
\end{cases}
\hspace{-.25cm}.
$
Define a lifting function $\Lambda:W_p \ra PX$ by 
$
\Lambda(x,\tau) = 
\begin{cases}
\ \Lambda_1(x,\tau) \hspace{0.5cm} ((x,\tau) \in W_{p_1})\\
\ \Lambda_2(x,\tau)  \hspace{0.5cm} ((x,\tau) \in W_{p_2})
\end{cases}
$
and cite Proposition 8.]\\

\label{5.46}
\begingroup%%----------------------------------->>
\fontsize{9pt}{11pt}\selectfont
\index{Mayer-Vietoris Condition (example)} 
\textbf{\small FACT \ (\un{Mayer-Vietoris Condition})} \ 
Suppose that $B = B_1 \cup B_2$, where 
$
\begin{cases}
\ B_1\\
\ B_2
\end{cases}
$
are closed and the inclusions $B_0 = B_1 \cap B_2 \ra$ 
$
\begin{cases}
\ B_1\\
\ B_2
\end{cases}
$
are cofibrations.  Let 
$
\begin{cases}
\ X_1 \ra B_1\\
\ X_2 \ra B_2
\end{cases}
$
be Hurewicz fibrations.  Assume: 
$
\begin{cases}
\ \restr{X_1}{B_0}\\
\ \restr{X_2}{B_0}
\end{cases}
$
have the same fiber homotopy type $-$then there exists a Hurewicz fibration $X \ra B$ such that 
$
\begin{cases}
\ X_1 \ \& \ \restr{X}{B_1}\\
\ X_2 \ \& \ \restr{X}{B_2}\
\end{cases}
$
have the same fiber homotopy type.\\
\endgroup%%------------------------------------<<

\begingroup%%----------------------------------->>
\fontsize{9pt}{11pt}\selectfont
\textbf{\small FACT}  \ 
Let 
\begin{tikzcd}[sep=large]
{X_0} \ar{d} \ar{r}{p_0} &{B_0} \ar{d} &{Y_0} \ar{l}[swap]{q_0} \ar{d}\\
{X} \ar{r}[swap]{p} &{B} &{Y} \ar{l}{q}
\end{tikzcd}
be a \cd in which the vertical arrows are inclusions and closed cofibrations.  Assume that the projections 
$
\begin{cases}
\ p_0\\
\ p
\end{cases}
$
are Hurewicz fibrations $-$then the induced map $X_0 \times_{B_0} Y_0 \ra X \times_B Y$ is a closed cofibration.
\\ \indent
[The inclusion $p^{-1}(B_0) \ra X$ is a closed cofibration (cf. Proposition 11).  
Since $X_0$ is contained in 
$p^{-1}(B_0)$ and since the inclusion $X_0 \ra X$ is a closed cofibration, the inclusion $X_0 \ra p^{-1}(B_0)$ is a closed cofibration (cf. $\S 3$, Proposition 9).  
Proposition 13 then implies that the arrow 
$i_0 p^{-1}(B_0) \cup IX_0 \ra B_0$ is a Hurewicz fibration.  
Consequently, (cf. Proposition 12), the commutative diagram 
\endgroup%%------------------------------------<<

\begingroup%%----------------------------------->>
\fontsize{9pt}{11pt}\selectfont
\[
\begin{tikzcd}[sep=large]
{i_0 p^{-1}(B_0) \cup IX_0} \ar{d} \ar{r}{\id} &{i_0 p^{-1}(B_0) \cup IX_0} \ar{d}\\
{Ip^{-1}(B_0)} \ar{r} &{B_0}
\end{tikzcd}
\]
has a filler $r:Ip^{-1}(B_0) \ra i_0p^{-1}(B_0) \cup IX_0$.  
Therefore the inclusion 
$X_0 \times_{B_0} Y_0 \ra p^{-1}(B_0) \times_B Y_0$ is a closed cofibration.  On the other hand, the projection
$X \times_B Y \ra Y$ is a Hurewicz fibration (cf. Proposition 4) and the inclusion $Y_0 \ra Y$ is a closed cofibration, 
so the inclusion $p^{-1}(B_0) \times_B Y_0 \ra X \times_B Y$ is a closed cofibration (cf. Proposition 11).]\\
\endgroup%%------------------------------------<<

\begingroup%%----------------------------------->>
\fontsize{9pt}{11pt}\selectfont
Application: Consider the 2-sink
$X \overset{p}{\ra} B \overset{q}{\la} Y$, where $p:X \ra B$ is a Hurewicz fibration.  Assume: The inclusions 
$\Delta_X \ra X \times X$, 
$\Delta_B \ra B \times B$, 
$\Delta_Y \ra Y \times Y$
are closed cofibrations $-$then the diagonal embedding 
$X \times_B Y \ra (X \times_B Y) \times (X \times_B Y)$ is a closed cofibration.\\
\endgroup%%------------------------------------<<

%%----------------------------------------------------------------------------------------------20
\label{4.20}
\label{4.34}
Let 
$X \overset{p}{\ra} B \overset{q}{\la} Y$ be a 2-sink $-$then the 
\un{fiber join}
\index{fiber join} 
$X *_B Y$ 
\index{$X *_B Y$} 
is the double mapping cylinder of the 2-source 
$X \overset{\xi}{\la} X *_B Y \overset{\eta}{\ra} Y$.  The fiber homotopy type of $X *_B Y$ depends only on the fiber homotopy types of \mX and \mY.  There is a projection 
$X *_B Y \ra B$ and the fiber over $b$ is $X_b * Y_b$.  
Examples: 
(1) The fiber join of $X \overset{p}{\ra} B \la B \times \{0\}$ is $\Gamma_B X$, the 
\un{fiber cone}
\index{fiber cone} 
of \mX;
(2) The fiber join of $X \overset{p}{\ra} B \la B \times \{0,1\}$ is $\Sigma_B X$, the 
\un{fiber suspension}
\index{fiber suspension} 
of $X$;
(3) The fiber join of $B \times T_1 \ra B \la B \times T_2$ is $B \times (T_1 * T_2)$; 
(4) The fiber join of  $\{b_0\} \ra B \overset{p}{\la} X$ is the mapping cone $C_{b_0}$ of the inclusion 
$X_{b_0} \ra X$.\\

\begingroup%%----------------------------------->>
\fontsize{9pt}{11pt}\selectfont
Let \mX be in \bTOP/\mB $-$then $\Gamma_B X$ can be identified with the mapping cylinder $M_p$ and 
$\Sigma_B X$ can be identified with the double mapping cylinder $M_{p,p}$.\\
\endgroup%%------------------------------------<<

\textbf{\small LEMMA}  \ 
Let $f \in C_B(X,Y)$.  Suppose that
$
\begin{cases}
\ p:X \ra B\\[-.15cm]
\ q:Y \ra B
\end{cases}
$ 
are Hurewicz fibrations $-$then the projection $\pi:M_f \ra B$ is a Hurewicz fibration.

[Fix lifting functions 
$
\begin{cases}
\ \Lambda_X:W_p \ra PX\\[-.15cm]
\ \Lambda_Y:W_q \ra PY
\end{cases}
\hspace{-.25cm}. \ 
$
Define a lifting function $\Lambda:W_\pi \ra PM_f$ as follows: Given $((x,t),\tau) \in IX \times_B PB$, put
\[
\Lambda((x,t),\tau)(T) = 
\begin{cases}
\ (\Lambda_X(x,\tau)(T),(t-1/2)(1 + T) + (1 - T)/2) \hspace{0.28cm} \  (1/2 \leq t \leq 1)\\
\ (\Lambda_X(x,\tau)(T),t - T/2) \hspace{3.7cm}\ \ (0 \leq t \leq 1/2 \ \& \ T \leq 2t)\\
\ \Lambda_Y(f(\Lambda_X(x,\tau)(2t)),\tau_{2t})(T - 2t) \hspace{2.15cm} \ (0 \leq t \leq 1/2 \ \& \ T \geq 2t)
\end{cases}
\hspace{-.4cm},
\]
where $\tau_{2t}(T) = \tau (\min\{2t + T,1\})$, and given $(y,\tau) \in Y \times_B PB$, put 
$\Lambda(y,\tau) = \Lambda_Y(y,\tau)$.]\\

\begin{proposition} \  %14
Suppose that 
$
\begin{cases}
\ p:X \ra B\\[-.15cm]
\ q:Y \ra B
\end{cases}
$
are Hurewicz fibrations $-$then the projection $X *_B Y \ra B$ is a Hurewicz fibration.
\end{proposition}

[Consider the pushout square 
\begin{tikzcd}[sep=large]
{X \times_B Y} \ar{d} \ar{r} &{M_\eta} \ar{d}\\
{M_\xi} \ar{r} &{X *_B Y}
\end{tikzcd}
(cf. p. \pageref{4.16}).  Here, the arrows 
$
X \times_B Y \ra 
\begin{cases}
\ M_\eta\\[-.15cm]
\ M_\xi
\end{cases}
\ra X *_B Y
$
are closed cofibrations and the projections
$X \times_B Y \ra B$, \ 
$
\begin{cases}
\ M_\eta\\[-.15cm]
\ M_\xi
\end{cases}
\ra B
$
are Hurewicz fibrations.  That the projection $X *_B Y \ra B$ is a Hurewicz fibration is therefore a consequence of 
Proposition 13.]\\

Application: Let $p:X \ra B$ be a Hurewicz fibration $-$then the projections 
$
\begin{cases}
\ \Gamma_B X \ra B\\[-.15cm]
\ \Sigma_B X \ra B
\end{cases}
$
are Hurewicz fibrations.\\

\begingroup%%----------------------------------->>
\fontsize{9pt}{11pt}\selectfont
Let 
$X \overset{p}{\ra} B \overset{q}{\la} Y$ be a 2-sink, where $p$ is a Hurewicz fibration.  There is a commutative diagram
%%----------------------------------------------------------------------------------------------21
\begin{tikzcd}[sep=large]
\vspace{0.25cm}
{X} \arrow[d,shift right=0.5,dash] \arrow[d,shift right=-0.5,dash] \ar{r}{p} 
&{B} \arrow[d,shift right=0.5,dash] \arrow[d,shift right=-0.5,dash] 
&{Y} \ar{l}[swap]{q} \ar{d}{\gamma}\\
{X} \ar{r}[swap]{p} &{B} &{W_q} \ar{l}
\end{tikzcd}
and $\gamma$ is a homotopy equivalence, thus the induced map 
$X \times_B Y \ra X \times_B W_q$ is a homotopy equivalence (cf. p. \pageref{4.17}).  
Consideration of 
\begin{tikzcd}[sep=large]
{X} \arrow[d,shift right=0.5,dash] \arrow[d,shift right=-0.5,dash] 
&{X \times_B Y} \ar{l} \ar{d} \ar{r} 
&{Y} \arrow[d,shift right=0.5,dash] \arrow[d,shift right=-0.5,dash]\\
{X} &{X \times_B W_q} \ar{l} \ar{r} &{Y} 
\end{tikzcd}
then leads to a homotopy equivalence $X *_B Y \ra X *_B W_q$ (cf. p. \pageref{4.18}).  
Example: $\forall \ b_0 \in B$, $X *_B \Theta B$ and $C_{b_0}$ have the same homotopy type.\\

\label{12.3}
\label{12.8}
Assume in addition that $q$ is a closed cofibration and define \mP by the pushout square   
\begin{tikzcd}%[sep=large]
{X \times_B Y} \ar{d}[swap]{\xi} \ar{r}{\eta} &{Y} \ar{d}\\
{X} \ar{r} &{P}
\end{tikzcd}
$-$then Proposition 11 implies that $\xi$ is a closed cofibration.  Therefore the arrow 
$X *_B Y \ra P$ of $\S 3$, Proposition 18 is a homotopy equivalence.  
Example: $\forall \ b_0 \in B$ such that the inclusion $\{b_0\} \ra B$ is a closed cofibration, 
$\Theta B *_B \Theta B$ and $\Theta B / \Omega B$ have the same homotopy type.\\
\endgroup%%------------------------------------<<

\begin{proposition} \  %15
Suppose that 
$
\begin{cases}
\ p:X \ra B\\[-.15cm]
\ q:Y \ra B
\end{cases}
$
are Hurewicz fibrations.  Let $\phi \in C_B(X,Y)$.  
Assume that $\phi$ is a homotopy equivalence $-$then $\phi$ is a homotopy equivalence in \bTOP/\mB.
\end{proposition}

[This is the analog of $\S 3$, Proposition 13.  It is a special case of Proposition 16 below.]\\

Application: Let $p:X \ra B$ be a homotopy equivalence $-$then $W_p$ is fiberwise contractible.

[Write $p = q \circx \gamma$: $p$ and $\gamma$ are homotopy equivalences, thus so is $q$.]

[Note: \ Similar reasoning leads to another proof of Proposition 9.]\\

\begingroup%%----------------------------------->>
\fontsize{9pt}{11pt}\selectfont
\textbf{\small EXAMPLE}  \ 
Let $p:X \ra B$ be a Hurewicz fibration.  View $PX$ as an object in $\bTOP/W_p$ with projection 
$\lambda:PX \ra W_p$ $-$then $PX$ is fiberwise contractible.\\
\endgroup%%------------------------------------<<

\begingroup%%----------------------------------->>
\fontsize{9pt}{11pt}\selectfont
\textbf{\small FACT}  \ 
Let $p:X \ra B$ be a continuous function $-$then $p$ is both a homotopy equivalence and a Hurewicz fibration iff every commutative diagram 
\begin{tikzcd}[sep=large]
{A} \ar{d}[swap]{i} \ar{r} &{X} \ar{d}{p}\\
{Y} \ar{r} &{B}
\end{tikzcd}
, where $i$ is a closed cofibration, has a filler $Y \ra X$.
\\ \indent
[To discuss the necessity, use Proposition 12, noting that \mX is fiberwise contractible, hence $\exists$ 
$s \in \sec_B(X)$: $s \circx p \underset{B}{\simeq} \id_X$.]\\
\endgroup%%------------------------------------<<

%%----------------------------------------------------------------------------------------------22
\label{12.4} %dmc these first 9 may well need tuning 
\label{12.9}
\label{14.101}
\begingroup%%----------------------------------->>
\fontsize{9pt}{11pt}\selectfont
Application: Let 
\begin{tikzcd}[sep=large]
{X^\prime} \ar{d}[swap]{p^\prime} \ar{r} &{X} \ar{d}{p}\\
{B^\prime} \ar{r} &{B}
\end{tikzcd}
be a pullback square.  
Suppose that $p$ is a Hurewicz fibration and a homotopy equivalence $-$then $p^\prime$ is a Hurewicz fibration and a homotopy equivalence.\\
\endgroup%%------------------------------------<<

\begingroup%%----------------------------------->>
\fontsize{9pt}{11pt}\selectfont
\textbf{\small FACT}  \ 
Let $i:A \ra Y$ be a continuous function $-$then $i$ is a closed cofibration iff every commutative diagram 
\begin{tikzcd}[sep=large]
{A} \ar{d}[swap]{i} \ar{r} &{X} \ar{d}{p}\\
{Y} \ar{r} &{B}
\end{tikzcd}
, where $p$ is both a homotopy equivalence and a Hurewicz fibration, has a filler $Y \ra X$.
\\ \indent
[To establish the sufficiency, first consider \ 
\begin{tikzcd}[sep=large]
{A} \ar{d}[swap]{i} \ar{r} &{PX} \ar{d}{p_0}\\
{Y} \ar{r} &{X}
\end{tikzcd}
\ 
to see that $i$ is a cofibration.  \ 
Taking $i$ to be an inclusion, put $X = IA \cup Y \times \ ]0,1]$ $-$then the restriction to \mX 
of the Hurewicz fibration $IY \ra Y$ is a Hurewicz fibration (cf. p. \pageref{4.19}), call it $p$.  \ 
Since $p$ is also a homotopy equivalence, 
the commutative diagram 
\begin{tikzcd}[sep=large]
{A} \ar{d}[swap]{i} \ar{r} &{X} \ar{d}{p}\\
{Y} \arrow[r,shift right=0.5,dash] \arrow[r,shift right=-0.5,dash] &{Y}
\end{tikzcd}
\ 
has a filler $f:Y \ra X$ ($a \ra (a,0)$ $(a \in A)$), therefore \mA is a zero set in \mY, thus is closed.]\\
\endgroup%%------------------------------------<<

\begingroup%%----------------------------------->>
\fontsize{9pt}{11pt}\selectfont
\textbf{\small FACT}  \ 
Let $X \overset{p}{\ra} B \overset{q}{\la} Y$ be a 2-sink, where $p:X \ra B$ is a Hurewicz fibration.  Denote by $W_*$ the 
mapping track of the projection $X *_B Y \ra B$ $-$then $X *_B W_q$ and $W_*$ have the same fiber homotopy type.\\
\endgroup%%------------------------------------<<

\textbf{\small LEMMA}  \ 
Suppose that $\xi \in C_B(X,E)$ is a fiberwise Hurewicz fibration.  
Let $f \in C(X,X)$: $\xi \circx f = \xi$ $\&$ 
$f \underset{B}{\simeq} \id_X$ $-$then $\exists$ $g \in C(X,X)$ : 
$\xi \circx g = \xi$ $\&$ $f \circx g \underset{E}{\simeq} \id_X$.

[Let $H:IX \ra X$ be a fiber homotopy with 
$H \circx i_0 = f$ and 
$H \circx i_1 = \id_X$;  
let
$G:IX \ra X$ be a fiber homotopy with 
$G \circx i_0 = \id_X$ and $\xi \circx G = \xi \circx H$.  
Define $F:IX \ra X$ by 
$
F(x,t) = 
\begin{cases}
\ f \circx G(x,1-2t) \hspace{0.5cm} (0 \leq t \leq 1/2)\\
\ H(x,2t-1) \hspace{1.15cm} (1/2 \leq t \leq 1)
\end{cases}
$
and put
\[
k((x,t),T) = 
\begin{cases}
\ \xi \circx G(x,1-2t(1-T)) \hspace{1.65cm}  (0 \leq t \leq 1/2)\\
\ \xi \circx H(x,1-2(1-t)(1-T)) \hspace{0.5cm} \ (1/2 \leq t \leq 1)
\end{cases}
\]
to get a fiber homotopy $k:I^2X \ra E$ with $\xi \circx F = k \circx i_0$.  Choose a fiber homotopy 
$K:I^2X \ra X$ such that $F = K \circx i_0$ and $\xi \circx K = k$.  Write $K_{(t,T)}:X \ra X$ for the function 
$x \ra K((x,t),T)$.  Obviously, 
$K_{(0,0)} \simeq K_{(0,1)} \simeq$ $K_{(1,1)} \simeq$ $K_{(1,0)}$ all fiber homotopies being over $E$.  
Set 
$g = G \circx i_1$ $-$then 
$f \circx g =$ $F \circx i_0 =$ $K_{(0,0)} \underset{E}{\simeq}$ $K_{(1,0)} =$ $F \circx i_1 =$ $\id_X$.]

[Note: \ Take $B = *$, $E = B$, $\xi = p$, so $p:X \ra B$ is a Hurewicz fibration $-$then the lemma asserts that 
$\forall \ f \in C_B(X,X)$, with $f \simeq \id_X$, $\exists$ $g \in C_B(X,X)$: $f \circx g \underset{B}{\simeq} \id_X$.]\\

%%----------------------------------------------------------------------------------------------23
\begin{proposition} \  %16
Suppose that 
$
\begin{cases}
\ \xi \in C_B(X,E)\\
\ \eta \in C_B(Y,E)
\end{cases}
$
are fiberwise Hurewicz fibrations.  Let $\phi \in C(X,Y)$ : $\eta \circx \phi = \xi$.  
Assume that $\phi$ is a homotopy equivalence in 
\bTOP/\mB $-$then $\phi$ is a homotopy equivalence in \bTOP/\mE.
\end{proposition}

[Since $\xi$ is a fiberwise Hurewicz fibration, there exists a fiber homotopy inverse $\psi:Y \ra X$ for $\phi$ with 
$\xi \circx \psi = \eta$, thus, from the lemma, 
$\exists$ $\psi^\prime \in C(Y,Y)$ : 
$\eta \circx \psi^\prime = \eta$ 
$\&$ $\phi \circx \psi \circx \psi^\prime \underset{E}{\simeq} \id_Y$.  
This says that $\phi^\prime = \psi \circx \psi^\prime$ is a homotopy right inverse for $\phi$ over $E$.  
Repeat the argument with $\phi$ replaced by $\phi^\prime$ to conclude that $\phi^\prime$ has a right homotopy inverse 
$\phi\pp$ over $E$, hence that $\phi^\prime$ is a homotopy equivalence in \bTOP/\mE or still, that $\phi$ is a homotopy equivalence in \bTOP/\mE.]

[Note: \ To recover Proposition 15, take $B = *$, $E = B$, $\xi = p$, and $\eta = q$.]\\

\begin{proposition} \  %17
Suppose given a \cd
\begin{tikzcd}%[sep=large]
{X} \ar{d}[swap]{\phi} \ar{r}{p} &{B} \ar{d}{\psi}\\
{Y} \ar{r}[swap]{q} &{A}
\end{tikzcd}
in which 
$
\begin{cases}
\ p\\[-.15cm]
\ q
\end{cases}
$
are Hurewicz fibrations and 
$
\begin{cases}
\ \phi\\[-.15cm]
\ \psi
\end{cases}
$
are homotopy equivalences $-$then $(\phi,\psi)$ is a homotopy equivalence in $\bTOP(\ra)$.
\end{proposition}

[This is the analog of $\S 3$, Proposition 14.]\\

Let $X \overset{f}{\ra} Z \overset{g}{\la} Y$ be a 2-sink $-$then the 
\un{double mapping track}
\index{double mapping track} 
$W_{f,g}$ 
\index{$W_{f,g}$} of $f, g$ 
is defined by the pullback square
\begin{tikzcd}%[sep=large]
{W_{f,g}} \ar{d} \ar{r} &{PZ} \ar{d}[swap]{p_0} \ar{d}{p_1}\\
{X \times Y} \ar{r}[swap]{f \times g} &{Z \times Z}
\end{tikzcd}
.  \ 
The homotopy type of $W_{f,g}$ depends only on the homotopy classes of $f$ and $g$ and $W_{f,g}$ is homeomorphic to 
$W_{g,f}$.  
There are Hurewicz fibrations 
$
\begin{cases}
\ p:W_{f,g} \ra X\\
\ q:W_{f,g} \ra Y
\end{cases}
\hspace{-.25cm}. \ 
$
The diagram 
\begin{tikzcd}%[sep=large]
{W_{f,g}} \ar{d}[swap]{p} \ar{r}{q} &{Y} \ar{d}{g}\\
{X} \ar{r}[swap]{f} &{Z}
\end{tikzcd}
is homotopy commutative and if the diagram 
\begin{tikzcd}%[sep=large]
{W} \ar{d}[swap]{\xi} \ar{r}{\eta} &{Y} \ar{d}{g}\\
{X} \ar{r}[swap]{f} &{Z}
\end{tikzcd}
is homotopy commutative, then there exists a $\phi:W \ra W_{f,g}$ such that 
$
\begin{cases}
\ \xi = p \circx \phi\\[-.15cm]
\ \eta = q \circx \phi
\end{cases}
\hspace{-.25cm}.
$
\label{12.23}

[Note: \ The \cd 
\begin{tikzcd}%[sep=large]
{W_{f,g}} \ar{d} \ar{r} &{Y} \ar{d}{g}\\
{W_f} \ar{r}[swap]{q} &{Z}
\end{tikzcd}
is a pullback square $(f = q \circx s)$.]\\

%%----------------------------------------------------------------------------------------------24
\begingroup%%----------------------------------->>
\fontsize{9pt}{11pt}\selectfont
\textbf{\small FACT}  \ 
Let $X \overset{f}{\ra} Z \overset{g}{\la} Y$ be a 2-sink $-$then the assignment 
$(x,y,\tau) \ra \tau(1/2)$ defines a Hurewicz fibration 
$W_{f,g} \ra Z$.
\\ \indent
[Let 
$
\begin{cases}
\ W_f^+ = \{(x,\tau):f(x) = \tau(0), \tau \in C([0,1/2],Z)\}\\
\ W_g^- = \{(y,\tau):g(y) = \tau(1), \tau \in C([1/2,1],Z)\}
\end{cases}
$
\hspace{-.25cm}.  \ The projections 
$
\begin{cases}
\ W_f^+ \ra Z\\
\ (x,\tau) \ra \tau(1/2)
\end{cases}
$
,
$
\begin{cases}
\ W_g^- \ra Z\\
\ (y,\tau) \ra \tau(1/2)
\end{cases}
$
are Hurewicz fibrations and the commutative diagram 
\begin{tikzcd}[sep=large]
{W_{f,g}} \ar{d} \ar{r} &{W_g^-} \ar{d}\\[-.15cm]
{W_f^+} \ar{r} &{Z}
\end{tikzcd}
is a pullback square.]\\
\endgroup%%------------------------------------<<

Every 2-sink $X \overset{f}{\ra} Z \overset{g}{\la} Y$ determines a pullback square 
\begin{tikzcd}[sep=large]
{P} \ar{d}[swap]{\xi} \ar{r}{\eta} &{Y} \ar{d}{g}\\
{X} \ar{r}[swap]{f} &{Z}
\end{tikzcd}
and there is an arrow $\phi:P \ra W_{f,g}$ characterized by the conditions 
$
\begin{cases}
\ \xi = p \circx \phi\\[-.15cm]
\ \eta = q \circx \phi
\end{cases}
$
$\&$ $P \overset{\phi}{\ra} W_{f,g} \ra PZ =$
$
\begin{cases}
\begin{tikzcd}
%{j \circx f \circx \xi} \ar[equals]{d}\\[-.15cm]
{j \circx f \circx \xi} \arrow[d,shift right=0.5,dash] \arrow[d,shift right=-0.5,dash]\\[-.15cm]
{j \circx g \circx \eta}
\end{tikzcd}
\end{cases}
$
\hspace{-.25cm}.\\
\vspace{0.25cm}

\begin{proposition} \  %18
If $f$ is a Hurewicz fibration, then $\phi:P \ra W_{f,g}$ is a homotopy equivalence in \bTOP/\mY.
\end{proposition}

[Use Proposition 9 and the fact that the pullback of a fiber homotopy equivalence is a fiber homotopy equivalence.]\\

Application: Let $p:X \ra B$ is a Hurewicz fibration.  Suppose that 
$
\begin{cases}
\ \Phi_1^\prime\\[-.15cm]
\ \Phi_2^\prime
\end{cases}
\in C(B^\prime,B)
$
are homotopic $-$then 
$
\begin{cases}
\ X_1^\prime\\[-.15cm]
\ X_2^\prime
\end{cases}
$
have the same homotopy type over $B^\prime$.\\

\label{4.21}
\label{5.0au1}
\label{6.28}
\label{13.31}
For example, under the assumption that $p:X \ra B$ is a Hurewicz fibration, if 
$\Phi^\prime:B^\prime \ra B$ is homotopic to the constant map $B^\prime \ra b_0$, then $X^\prime$ 
is fiber homotopy equivalent to $B^\prime \times X_{b_0}$.\\

\label{4.33}
\label{5.0t}
\begingroup%%----------------------------------->>
\fontsize{9pt}{11pt}\selectfont
\textbf{\small FACT}  \ 
Suppose that $p:X \ra B$ is a Hurewicz fibration.  Let $\Phi^\prime:B^\prime \ra B$ be a homotopy equivalence $-$then 
the arrow $X^\prime \ra X$ is a homotopy equivalence.\\
\endgroup%%------------------------------------<<

\begingroup%%----------------------------------->>
\fontsize{9pt}{11pt}\selectfont
Denote by $\abs{\id,\Delta}_{\bTOP}$ the comma category corresponding to the identity functor id on 
$\bTOP \times \bTOP$ and the diagonal functor $\Delta:\bTOP \ra \bTOP \times \bTOP$.  So, an object in 
 $\abs{\id,\Delta}_{\bTOP}$ is a 2-sink $X \overset{f}{\ra} Z \overset{g}{\la} Y$ and a morphism of 2-sinks is a 
 commutative diagram
\begin{tikzcd}[sep=large]
{X} \ar{d} \ar{r}{f} &{Z} \ar{d} &{Y} \ar{l}[swap]{g} \ar{d}\\
{X^\prime} \ar{r}[swap]{f^\prime} &{Z^\prime} &{Y^\prime} \ar{l}{g^\prime}
\end{tikzcd}
.  The double mapping
%%----------------------------------------------------------------------------------------------25
track is a functor $\abs{\id,\Delta}_{\bTOP} \ra \bTOP$.  It has a left adjoint $\bTOP \ra \abs{\id,\Delta}_{\bTOP}$, viz. 
the functor that sends \mX to the 2-sink $X \overset{i_0}{\ra} IX \overset{i_1}{\la} X$.\\
\endgroup%%------------------------------------<<

\label{6.19}
\begingroup%%----------------------------------->>
\fontsize{9pt}{11pt}\selectfont
\textbf{\small FACT}  \ 
Let 
\begin{tikzcd}[sep=large]
{X} \ar{d} \ar{r}{f} &{Z} \ar{d} &{Y} \ar{l}[swap]{g} \ar{d}\\
{X^\prime} \ar{r}[swap]{f^\prime} &{Z^\prime} &{Y^\prime} \ar{l}{g^\prime}
\end{tikzcd}
be a \cd in which the vertical arrows are homotopy equivalences $-$then the arrow 
$W_{f,g} \ra W_{f^\prime,g^\prime}$ is a homotopy equivalence.\\
\endgroup%%------------------------------------<<

\label{4.17}
\begingroup%%----------------------------------->>
\fontsize{9pt}{11pt}\selectfont
Application: 
Suppose that 
$
\begin{cases}
\ p:X \ra B\\
\ p^\prime:X^\prime \ra B^\prime
\end{cases}
$
are Hurewicz fibrations.  Let 
$
\begin{cases}
\ g:Y \ra B\\
\ g^\prime:Y^\prime \ra B^\prime
\end{cases}
$ 
be continuous functions.  Assume that the diagram 
\begin{tikzcd}[sep=large]
{X} \ar{d} \ar{r}{p} &{B} \ar{d} &{Y} \ar{l}[swap]{g} \ar{d}\\
{X^\prime} \ar{r}[swap]{p^\prime} &{B^\prime} &{Y^\prime} \ar{l}{g^\prime}
\end{tikzcd}
commutes and that the vertical arrows are homotopy equivalences $-$then the induced map 
$X \times_B Y \ra X^\prime \times_{B^\prime} Y^\prime$
is a homotopy equivalence.\\
\endgroup%%------------------------------------<<

\label{4.44}
\begingroup%%----------------------------------->>
\fontsize{9pt}{11pt}\selectfont
\textbf{\small EXAMPLE}  \ 
Suppose given a \cd 
\begin{tikzcd}[sep=large]
{X} \ar{d}[swap]{\phi} \ar{r}{p} &{B} \ar{d}{\psi}\\
{Y} \ar{r}[swap]{q} &{A}
\end{tikzcd}
in which 
$
\begin{cases}
\ p\\
\ q
\end{cases}
$
are Hurewicz fibrations and 
$
\begin{cases}
\ \phi\\
\ \psi
\end{cases}
$
are homotopy equivalences $-$then $\forall \ b \in B$, the induced map 
$X_b \ra Y_{\psi(b)}$ is a  homotopy equivalence.
\label{4.22}
\\ \indent
[Note: Let $f:X \ra Y$ be a homotopy equivalence, fix $x_0 \in X$ and put $y_0 = f(x_0)$, form the commutative 
diagram 
\begin{tikzcd}[sep=large]
{\Theta X} \ar{d} \ar{r}{p_1} &{X} \ar{d} &{\{x_0\}} \ar{l} \ar{d}\\
{\Theta Y} \ar{r}[swap]{p_1} &{Y} &{\{y_0\}} \ar{l}
\end{tikzcd}
, and conclude that the arrow $\Omega X \ra \Omega Y$ is a homotopy equivalence.]\\
\endgroup%%------------------------------------<<

\label{5.28}
\begingroup%%----------------------------------->>
\fontsize{9pt}{11pt}\selectfont
Given a 2-sink $X \overset{p}{\ra} B \overset{q}{\la} Y$, let $X \bboxsub_B Y$ be the double mapping cylinder of the 2-source $X \la W_{p,q} \ra Y$.  
It is an object in \bTOP/\mB with projection 
$
\begin{cases}
\ x \ra p(x)\\
\ y \ra q(y)
\end{cases}
, ((x,y,\tau),t) \ra \tau(t).
$
\\
\endgroup%%------------------------------------<<

\label{4.55}
\begingroup%%----------------------------------->>
\fontsize{9pt}{11pt}\selectfont
\textbf{\small FACT}  \ 
There is a  homotopy equivalence $X \bboxsub_B Y \overset{\phi}{\ra} W_p *_B W_q$.
\\ \indent
[Define $\phi$ by 
$
\begin{cases}
\ \phi(x) = \gamma(x)\\
\ \phi(y) = \gamma(y)
\end{cases}
$
$\&$ $\phi((x,y,\tau),t) = ((x,\tau_t),(y,\ov{\tau}_t),t)$, where 
$\tau_t(T) = \tau(tT)$ and 
$\ov{\tau}_t(T) = \tau(tT + 1 - T)$.]
\\ \indent
[Note: \ More is true if $p:X \ra B$ is a Hurewicz fibration: $X \bboxsub_B Y$ and $X *_B Y$ have the same homotopy type.  
Indeed, $W_p *_B W_q$ has the same fiber homotopy type as $X *_B W_q$ which in turn has the same homotopy type as $X *_B Y$ (cf. p. \pageref{4.20} ff.).]\\
\endgroup%%------------------------------------<<

\begingroup%%----------------------------------->>
\fontsize{9pt}{11pt}\selectfont
Application: $\forall \ b_0 \in B$, $\Sigma\Omega B$ and $\Theta B *_B \Theta B$ have the same homotopy type.
\\ \indent
%%----------------------------------------------------------------------------------------------26
[Note: \ The suspension is taken in \bTOP, not $\bTOP_*$.]\\
\endgroup%%------------------------------------<<

Given $f \in C_B(X,Y)$, let \mW be the subspace of $X \times PY$ consisting of the pairs 
$(x,\tau)$: $f(x) = \tau(0)$ and $p(x) = q(\tau(t))$ $(0 \leq t \leq 1)$ $-$then \mW is in \bTOP/\mY with projection 
$(x,\tau) \ra \tau(1)$ and is fiberwise contractible if $f$ is a fiber homotopy equivalence (cf. Proposition 16).

[Note: \ $W$ is an object in \bTOP/\mB with projection $(x,\tau) \ra p(x)$. 
 Viewed as an object in \bTOP/\mY, its projection $(x,\tau) \ra \tau(1)$ is therefore a morphism in \bTOP/\mB and as such, is a fiberwise Hurewicz fibration.]\\

\textbf{\small LEMMA}  \ 
$f$ admits a right fiber homotopy inverse iff $\sec_Y(W) \neq \emptyset$.\\

\begin{proposition} \  %19
Let $f \in C_B(X,Y)$.  Suppose that there exists a numerable covering $\sO = \{O_i: i \in I\}$ of \mB such that 
$\forall \ i$, $f_{O_i}: X_{O_i} \ra Y_{O_i}$ is a fiber homotopy equivalence $-$then $f$ is a fiber homotopy equivalence.
\end{proposition}

[It need only be shown that $\sec_Y(W) \neq \emptyset$.  
For then, by the lemma, $f$ has a right fiber homotopy inverse $g$ and, repeating the argument, $g$ has a right fiber homotopy inverse $h$, which means that $g$ is a fiber homotopy equivalence, thus so is $f$.  
This said, work with $f_{O_i} \in C_{O_i}(X_{O_i},Y_{O_i})$ and, as above, form 
$W_{O_i} \subset X_{O_i} \times PY_{O_i}$.  
Obviously, $\restr{W}{Y_{O_i}} = W_{O_i}$.  
The assumption that $f_{O_i}$ is a fiber homotopy equivalence implies that $W_{O_i}$ is fiberwise contractible, hence has the SEP.  
But $\{Y_{O_i}:i \in I\}$ is a numerable covering of \mY.  
Therefore, on the basis of the section extension theorem, \mW has the SEP.  In particular: $\sec_Y(W) \ne \emptyset$.]\\

\label{4.1}
Application: Let \mX be in \bTOP/\mB.  
Suppose that there exists a numerable covering $\sO = \{O_i:i \in I\}$ of $B$ such that $\forall \ i$, $X_{O_i}$ is fiberwise contractible $-$then \mX is fiberwise contractible.\\

\begin{proposition} \  %20
Let 
$
\begin{cases}
\ p:X \ra B\\[-.15cm]
\ q:Y \ra B\
\end{cases}
$
be Hurewicz fibrations, where \mB is numerably contractible.  
Suppose that $f \in C_B(X,Y)$ has the property that 
$f_b:X_b \ra Y_b$ is a homotopy equivalence at one point $b$ in each path component of \mB $-$then $f:X \ra Y$ is a fiber homotopy equivalence.
\end{proposition}

[Fix a numerable covering $\sO = \{O_i: i \in I\}$ of \mB for which the inclusions $O_i \ra B$ are inessential, say homotopic to $O_i \ra b_i$, where $f_{b_i}:X_{b_i} \ra Y_{b_i}$ is a homotopy equivalence $-$then $\forall \ i$, 
$f_{O_i}:X_{O_i} \ra Y_{O_i}$ is a fiber homotopy equivalence (cf. p. \pageref{4.21}), so Proposition 19 is applicable.]\\

\begingroup%%----------------------------------->>
\fontsize{9pt}{11pt}\selectfont
\textbf{\small EXAMPLE}  \ 
Take 
$B = \{0\} \cup \{1/n:n = 1, 2, \ldots\}$, 
$T = B \cup \{n:n = 1, 2, \ldots\}$, and put $X = B \times T$.  Observe that \mB is not numerably contractible.  
Let $k = 1, 2, \ldots, \infty$, $l = 0, 1, 2, \ldots$, and define
%%----------------------------------------------------------------------------------------------27
$f \in C_B(X,X)$ as follows: 
%(i) 
%(ii) $f(1/k,1/l) = (1/k,1/l)$ $(0 < l < k)$, $(1/k 
$
(i) \ f(1/k,l) =
\begin{cases}
\ (1/k,l) \hspace{0.85cm} (l < k)\\
\ (1/k,1/k)  \hspace{0.45cm} (l = k \neq 1)\\
\ (1/k,l-1) \hspace{0.27cm} (l > k)
\end{cases}
$
\hspace{-.35cm};
$
(ii) \ f(1/k,1/l) =
\begin{cases}
\ (1/k,1/l) \hspace{1.0cm} (0 < l < k)\\
\ (1/k,1/(l+1))  \hspace{0.17cm}  (l \geq k)
\end{cases}
$
$-$then $f$ is bijective and $\forall \ b \in B$, $f_b:X_b \ra X_b$ is a homeomorphism $(X_b = \{b\} \times T)$.  
Nevertheless, $f$ is not a fiber homotopy equivalence.  
For if it were, then $f$ would have to be a homeomorphism, an impossibility ($f^{-1}$ is not continuous at $(0,0)$).\\
\endgroup%%------------------------------------<<

\begingroup%%----------------------------------->>
\fontsize{9pt}{11pt}\selectfont
\label{6.17}
\label{14.46}
\label{14.53}
\label{14.60}
\label{14.73}
\label{14.92}
\index{Delooping Homotopy Equivalences (example)}
\textbf{\small EXAMPLE \ (\un{Delooping Homotopy Equivalences})} \ 
Suppose that 
$
\begin{cases}
\ X\\
\ Y
\end{cases}
$
are path connected and numerably contractible.  Let $f:X \ra Y$ be a continuous function.  Fix $x_0 \in X$ and put 
$y_0 = f(x_0)$ $-$then $f:X \ra Y$ is a homotopy equivalence  iff $\Omega f: \Omega X \ra \Omega Y$ is a homotopy equivalence.  In fact, the necessity is true without any restriction on \mX or \mY (cf. p. \pageref{4.22}).  
Turning to the sufficiency, write $f = q \circx s$, where $q:W_f \ra Y$.  
Since $s$ is a homotopy equivalence, one need only deal with $q$.  Form the pullback square
$
\begin{tikzcd}[sep=large]
{X \times_Y \Theta Y} \ar{d}\ar{r} &{\Theta Y} \ar{d}{p_1}\\
{X} \ar{r}[swap]{f} &{Y}
\end{tikzcd}
.  
$
The map 
$
\begin{cases}
\ \Theta X \ra X \times_Y \Theta Y\\
\ \sigma \ra (\sigma(1),f \circx \sigma)
\end{cases}
$
is a morphism in \bTOP/\mX which, when restricted to the fibers over $x_0$, is $\Omega f$, thus is a fiber homotopy equivalence (cf. Proposition 20).  In particular: $X \times_Y \Theta Y$ is contractible.  Consider now the commutative triangle 
$
\begin{tikzcd}%[sep=large]
{W_f} \ar{rd}[swap]{q} \ar{rr}&&{PY} \ar{ld}{p_1}\\
&{Y}
\end{tikzcd}
.
$
\ 
The fiber of $p_1$ over $y_0$ is contractible; on the other hand, the fiber of $q$ over $y_0$ is homeomorphic to 
$X \times_Y \Theta Y$ (parameter reversal).  
The arrow $W_f \ra PY$ is therefore a homotopy equivalence (cf. Proposition 20).  
But $p_1$ is a homotopy equivalence, hence so is $q$.\\
\endgroup%%------------------------------------<<

\begingroup%%----------------------------------->>
\fontsize{9pt}{11pt}\selectfont
\label{5.0m}
\index{H-Groups (example)}
\textbf{\small EXAMPLE \ (\un{H-Groups})} \ 
In any H-group (= cogroup object in $\bHTOP_*$), the operations of left and right translation are homotopy equivalences 
(so all path components have the same homotopy type).  
Conversely, let $(X,x_0)$ be a nondegenerate homotopy associative H-space with the property that 
the operations of left and right translation are homotopy equivalences.  
Assume: \mX is numerably contractible $-$then \mX admits a homotopy inverse, thus is an H-group.  
To see this, consider the shearing map
$
sh:
\begin{cases}
\ X \times X \ra X \times X\\
\ (x,y) \ra (x,xy)
\end{cases}
\hspace{-.25cm}. \ 
$
Agreeing to view $X \times X$ as an object in \bTOP/\mX via the first projection, Proposition 20 implies that sh is a homotopy equivalence over \mX.  Therefore sh is a homotopy equivalence or still, sh is a pointed homotopy equivalence, 
$(X \times X,(x_0,x_0))$ being nondegenerate (cf. p. \pageref{4.23}).  
Consequently, \mX is an H-group.

[Note: \ If $(X,x_0)$ is a homotopy associative H-space and if $\pi_0(X)$ is a group, then the operations of left and right translation are homotopy equivalences.]

Example: Let \mK be a compact ANR.  
Denote by $HE(K)$ the subspace of $C(K,K)$ (compact open topology) 
consisting of the homotopy equivalences $-$then $HE(K)$ is open in $C(K,K)$, hence is an ANR (cf. $\S 6$, Proposition 6).  
In particular: $(HE(K),\id_K)$ is wellpointed (cf. p. \pageref{4.24}) and numerably contractible
%%----------------------------------------------------------------------------------------------28
(cf. p. \pageref{4.25}).  
Because $HE(K)$ is a topological semigroup with unit under composition and $\pi_0(HE(K))$ is a group, it follows that $HE(K)$ is an H-group.\\
\endgroup%%------------------------------------<<

\begingroup%%----------------------------------->>
\fontsize{9pt}{11pt}\selectfont
\label{5.0am}
\label{5.45}
\index{Small Skeletons (example)}
\textbf{\small EXAMPLE \ (\un{Small Skeletons})} \ 
In algebraic topology, it is often necessary to determine whether a given category has a small skeleton.  
For instance,
if \mB  is a connected, locally path connected, locally simply connected space, then the full subcategory of \bTOP/\mB 
whose objects are the covering projections $X \ra B$ has a small skeleton.  
Here is a less apparent example.  
Fix a nonempty topological space $F$.  
Given a numerably contractible topological space \mB, 
let $\bFIB_{B,F}$ be the category whose objects are the Hurewicz fibrations $X \ra B$ such that $\forall \ b \in B$, $X_b$ has the homotopy type of $F$, 
and whose morphisms $X \ra Y$ are the fiber homotopy classes $[f]:X \ra Y$.  
The functor $\bFIB_{B,F} \ra \bFIB_{B^\prime,F}$ determined by a homotopy equivalence 
$\Phi^\prime:B^\prime \ra B$ induces a bijection 
$\Ob \ov{\bFIB}_{B,F} \ra \Ob\ov{\bFIB}_{B^\prime,F}$, hence $\bFIB_{B,F}$ has a small skeleton iff this is the case of 
$\bFIB_{B^\prime,F}$.
\\ \indent
Claim: Consider a 2-source $B_1 \overset{\phi_1}{\la} B_0 \overset{\phi_2}{\ra} B_2$, where $B_0$, 
$
\begin{cases}
\ B_1\\
\ B_2
\end{cases}
$
are numerably contractible.  Suppose that $\bFIB_{B_0,F}$, 
$
\begin{cases}
\ \bFIB_{B_1,F}\\
\ \bFIB_{B_2,F}
\end{cases}
$
have small skeletons $-$then $\bFIB_{M_{\phi_1,\phi_2},F}$ has a small skeleton.
\\ \indent
[Observing that the double mapping cylinder $M_{\phi_1,\phi_2}$ is numerably contractible, write
$
\begin{cases}
\ \phi_1 = r_1 \circx i_1\\
\ \phi_2 = r_2 \circx i_2
\end{cases}
\hspace{-.25cm},
$
where
$
\begin{cases}
\ r_1\\
\ r_2
\end{cases}
$
are homotopy equivalences and 
$
\begin{cases}
\ i_1\\
\ i_2
\end{cases}
$
are closed cofibrations (cf. $\S 3$, Proposition 16).   There is a commutative diagram
\begin{tikzcd}[sep=large]
{M_{i_1}} \ar{d}[swap]{r_1} 
&{B_0} \ar{l}[swap]{i_1} \arrow[d,shift right=0.5,dash] \arrow[d,shift right=-0.5,dash] \ar{r}{i_2} 
&{M_{i_2}} \ar{d}{r_2}\\
{B_1}   &{B_0}  \ar{l}{\phi_1}  \ar{r}[swap]{\phi_2} &{B_2}
\end{tikzcd}
and the arrow  $M_{i_1,i_2} \ra M_{\phi_1,\phi_2}$ is a homotopy equivalence  (cf. p. \pageref{4.26}).  
Thus one can assume that 
$
\begin{cases}
\ \phi_1\\
\ \phi_2
\end{cases}
$
are closed cofibrations.  But then if \mB is defined by the pushout square 
\begin{tikzcd}[sep=large]
{B_0} \ar{d}[swap]{\phi_1} \ar{r}{\phi_2} &{B_2} \ar{d}\\
{B_1} \ar{r} &{B}
\end{tikzcd}
, the arrow $M_{\phi_1,\phi_2} \ra B$ is a homotopy equivalence (cf. $\S 3$, Proposition 18).  
So, with $B_0 = B_1 \cap B_2$, 
$
\begin{cases}
\ B_1 \subset B\\
\ B_2 \subset B
\end{cases}
$
\hspace{-.25cm}, take an \mX in $\bFIB_{B,F}$ and put $X_0 = \restr{X}{B_0}$, 
$
\begin{cases}
\ X_1 = \restr{X}{B_1}\\
\ X_2 = \restr{X}{B_2}
\end{cases}
$
to get a commutative diagram 
\begin{tikzcd}[sep=large]
{X_1} \ar{d} &{X_0} \ar{l}[swap]{\psi_1} \ar{d} \ar{r}{\psi_2} &{X_2} \ar{d}\\
{B_1}   &{B_0} \ar{l}{\phi_1} \ar{r}[swap]{\phi_2} &{B_2}
\end{tikzcd}
in which
$
\begin{cases}
\ \psi_1\\
\ \psi_2
\end{cases}
$
are closed cofibrations (cf. Proposition 11).  In the skeletons of $\bFIB_{B_0,F}$, 
$
\begin{cases}
\ \bFIB_{B_1,F}\\
\ \bFIB_{B_2,F}
\end{cases}
\hspace{-.25cm} ,
$
choose objects $Y_0$, 
$
\begin{cases}
\ Y_1\\
\ Y_2
\end{cases}
$
and fiber homotopy equivalences $f_0:Y_0 \ra X_0$, 
$
\begin{cases}
\ f_1:Y_1 \ra X_1\\
\ f_2:Y_2 \ra X_2
\end{cases}
\hspace{-.25cm}:
$
$
\begin{cases}
\ p_1 \circx f_1 = q_1\\
\ p_2 \circx f_2 = q_2
\end{cases}
$
(obvious notation).  Let 
$
\begin{cases}
\ g_1:X_1 \ra Y_1\\
\ g_2:X_2 \ra Y_2
\end{cases}
$
be a fiber homotopy inverse for 
$
\begin{cases}
\ f_1\\
\ f_2
\end{cases}
\hspace{-.25cm}.
$
Set 
$
\begin{cases}
\ F_1 = g_1 \circx \psi_1 \circx f_0\\
\ F_2 = g_2 \circx \psi_2 \circx f_0\
\end{cases}
\ :\ 
$
$
\begin{cases}
\ f_1 \circx F_1 \simeq \psi_1 \circx f_0\\
\ f_2 \circx F_2 \simeq \psi_2 \circx f_0
\end{cases}
\hspace{-.25cm}.
$
Write
$
\begin{cases}
\ F_1 = \Psi_1 \circx l_1\\
\ F_2 = \Psi_2 \circx l_2
\end{cases}
\hspace{-.25cm}, \ 
$
%%----------------------------------------------------------------------------------------------29
where 
$
\begin{cases}
\ \Psi_1\\
\ \Psi_2
\end{cases}
$
are Hurewicz fibrations and homotopy equivalences and 
$
\begin{cases}
\ l_1\\
\ l_2
\end{cases}
$
are closed cofibrations (cf. p. \pageref{4.27}), say 
$
\begin{cases}
\ l_1:Y_0 \ra \ov{Y}_1 \ \& \  \Psi_1:\ov{Y}_1 \ra Y_1\\
\ l_2:Y_0 \ra \ov{Y}_2 \ \& \  \Psi_2:\ov{Y}_2 \ra Y_2
\end{cases}
\hspace{-.25cm}. \ 
$
Here
$
\begin{cases}
\ \ov{Y}_1\\
\ \ov{Y}_2
\end{cases}
$ 
is an object in 
$
\begin{cases}
\ \bTOP / B_1\\
\ \bTOP / B_2
\end{cases}
$
with projection 
$
\begin{cases}
\ q_1 \circx \Psi_1\\
\ q_2 \circx \Psi_2
\end{cases}
$
and 
$
\begin{cases}
\ f_1 \circx \Psi_1:\ov{Y}_1 \ra X_1\\
\ f_2 \circx \Psi_2:\ov{Y}_2 \ra X_2
\end{cases}
$
is a fiber homotopy equivalence (cf. Proposition 15).  Change
$
\begin{cases}
\ f_1 \circx \Psi_1\\
\ f_2 \circx \Psi_2
\end{cases}
$
by a homotopy over 
$
\begin{cases}
\ B_1\\
\ B_2
\end{cases}
$
into a map
$
\begin{cases}
\ G_1\\
\ G_2
\end{cases}
$
such that
$
\begin{cases}
\ G_1 \circx l_1 = \psi_1 \circx f_0\\
\ G_2 \circx l_2 = \psi_2 \circx f_0\
\end{cases}
$
\hspace{-.25cm}.  Form the pushout square
\begin{tikzcd}[sep=large]
{Y_0} \ar{d}[swap]{l_1} \ar{r}{l_2} &{\ov{Y}_2} \ar{d}\\
{\ov{Y}_1} \ar{r} &{Y}
\end{tikzcd}
$-$then \mY is in \bTOP/\mB and there is a fiber homotopy equivalence $f:Y \ra X$, i.e., this process picks up all the isomorphism classes in $\bFIB_{B,F}$.]
\\ \indent
Example: Let \mB be a CW complex $-$then \mB is numerably contractible (cf. p. \pageref{4.28}) and $\bFIB_{B,F}$ has a small skeleton.  In fact, $B = \colimx B^{(n)}$, so by induction, $\bFIB_{B^{(n)},F}$ has a small skeleton $\forall \ n$.  On the other hand, \mB and $\tel B$ have the same homotopy type (cf. p. \pageref{4.29}) and $\tel B$ is a double mapping cylinder calculated on the $B^{(n)}$ (cf. p. \pageref{4.30}).\\
\endgroup%%------------------------------------<<

\begingroup%%----------------------------------->>
\fontsize{9pt}{11pt}\selectfont
\textbf{\small FACT}  \ 
Let \mX be in \bTOP/\mB.  
Suppose that $\sU = \{U_i:i \in I\}$ is a numerable covering of \mX such that for every nonempty finite subset $F \subset I$, 
the restriction of $p$ to $\ds\bigcap\limits_{i \in F} U_i$ is a Hurewicz fibration $-$then $p:X \ra B$ is a Hurewicz fibration.

[Equip $I$ with a well ordering $<$ and use the Segal-Stasheff construction to produce a lifting function 
$\Lambda:W_p \ra PX$.  Compare this result with Proposition 13 when $I = \{1,2\}$.]\\
\endgroup%%------------------------------------<<

The property of being a Hurewicz fibration is not a fiber homotopy type invariant, i.e., if \mX and \mY have the same fiber homotopy type and if $p:X \ra B$ is a Hurewicz fibration, then $q:Y \ra B$ need not be a Hurewicz fibration.  
Example: Take 
$X = [0,1] \times [0,1]$, 
$Y = ([0,1] \times \{0\}) \cup (\{0\} \times [0,1])$, 
$B = [0,1]$, and let $p, q$ be the vertical projections $-$then \mX and \mY are fiberwise contractible and $p:X \ra B$ is a Hurewicz fibration but $q:Y \ra B$ is not a Hurewicz fibration.  This difficulty can be circumvented by introducing still another notion of ``fibration''.

Let \mX be in \bTOP/\mB.  
Let \mY be in \bTOP $-$then the projection $p:X \ra B$ is said to have the HLP w.r.t \mY 
\un{up to homotopy}
\index{HLP up to homotopy} 
if given continuous functions
$
\begin{cases}
\ F:Y \ra X\\[-.15cm]
\ h:IY \ra B
\end{cases}
$
such that $p \circx F = h \circx i_0$, there is a continuous function $H:IY \ra X$ such that 
$F \underset{B}{\simeq} H \circx i_0$ and $p \circx H = h$.

[Note: \ To interpret the condition $F \underset{B}{\simeq} H \circx i_0$, view \mY as an object in \bTOP/\mB with projection $p \circx F$.]\\

\textbf{\small LEMMA}  \ 
The projection $p:X \ra B$ has the HLP w.r.t \mY up to homotopy iff given 
%%----------------------------------------------------------------------------------------------30
continuous functions
$
\begin{cases}
\ F:Y \ra X\\[-.15cm]
\ h:IY \ra B
\end{cases}
$
such that $p \circx F = h \circx i_t$ $(0 \leq t \leq 1/2)$, there is a continuous function $H:IY \ra X$ such that 
$F = H \circx i_0$ and $p \circx H = h$.\\

Let \mX in \bTOP/\mB $-$then $p:X \ra B$ is said to be a 
\un{Dold fibration}
\index{Dold fibration} 
if it has the HLP w.r.t. \mY up to homotopy for every \mY in \bTOP.  Obviously, Hurewicz $\implies$ Dold, but 
Dold $\nRightarrow$ Serre and  Serre $\nRightarrow$ Dold.  
The pullback of a Dold fibration is a Dold fibration and the local-global principle remains valid.\\
\label{4.46}

\begin{proposition} \  %21
Let \mX, \mY be in \bTOP/\mB  and suppose that $q:Y \ra B$ is a Dold fibration.  
Assume: $\exists$ 
$
\begin{cases}
\ f \in C_B(X,Y)\\[-.15cm]
\ g \in C_B(Y,X)
\end{cases}
: g \circx f \underset{B}{\simeq} \id_X 
$
$-$then $p:X \ra B$ is a Dold fibration.
\end{proposition}

[Fix a topological space $E$ and continuous functions 
$
\begin{cases}
\ \Phi : E \ra X\\[-.15cm]
\ \Psi : IE \ra B
\end{cases}
$
such that $p \circx \Phi = \psi \circx i_0$.  Since $q \circx f = p$, $\exists$ $G:IE \ra Y$ with 
$f \circx \Phi \underset{B}{\simeq} G \circx i_0$, and $q \circx G = \psi$.  Put 
$\Psi = g \circx G$: $\Phi \underset{B}{\simeq} g \circx f \circx \Phi$ 
$\underset{B}{\simeq} \Psi \circx i_0$ $\&$ 
$p \circx \Psi = $
$p \circx g \circx G =$
$q \circx G = \psi$.]\\

The property of being a Dold fibration is therefore a fiber homotopy type invariant.  
Example: Take $X = ([0,1] \times \{0\}) \cup (\{0\} \times [0,1])$, $B = [0,1]$, and let $p$ be the vertical projection 
$-$then $p:X \ra B$ is a Dold fibration but not a Hurewicz fibration (nor is $p$ an open map (cf. p. \pageref{4.31})).\\

\begingroup%%----------------------------------->>
\fontsize{9pt}{11pt}\selectfont
\textbf{\small EXAMPLE}  \ 
Define $f:[-1,1] \ra [-1,1]$ by $f(x) = 2 \abs{x} - 1.$  Put $X = I[-1,1]/\sim$, where $(x,0) \sim (f(x),1)$, and let 
$p:X \ra \bS^1$ be the projection $-$then $p$ is an open map and a Dold fibration but not a Hurewicz fibration.\\
\endgroup%%------------------------------------<<

\begingroup%%----------------------------------->>
\fontsize{9pt}{11pt}\selectfont
\textbf{\small FACT}  \ 
Suppose that $B$ is numerably contractible, so $B$ admits a numerable covering $\{O\}$ for which each inclusion 
$O \ra B$ is inessential.  
Let \mX be in \bTOP/\mB $-$then the projection $p:X \ra B$ is a Dold fibration iff 
$\forall \ O$ there exists a topological space $T_O$ and a fiber homotopy equivalence $X_O \ra O \times T_O$ over $O$.\\
\endgroup%%------------------------------------<<

The homotopy theory of Hurewicz fibrations carries over to Dold fibrations.  
The proofs are only slightly more complicated.  
Specifically, Propositions 15, 17, 18, and 20 are true if ``Hurewicz'' is replaced by ``Dold''.\\

\begin{proposition} \  %22
Let \mX be in \bTOP/\mB $-$then \mX is fiberwise contractible iff $p:X \ra B$ is a Dold fibration and a homotopy equivalence.
\end{proposition}

[The necessity is a consequence of Proposition 21 and the sufficiency is a consequence of Proposition 15.]\\

%%----------------------------------------------------------------------------------------------31
\begin{proposition} \  %23
Let \mX be in \bTOP/\mB $-$then $p:X \ra B$ is a Dold fibration  iff $\gamma:X \ra W_p$ is a fiber homotopy equivalence.
\end{proposition}

[Bearing in mind that $q:W_q \ra B$ is a Hurewicz fibration, the reasoning is the same as that used in the proof of Proposition 22.]\\

\label{4.47}
Application: The fibers of a Dold fibration over a path connected base have the same homotopy type.

[Note: \ Take $X = ([0,1] \times \{0,1\}) \cup (\{0\} \times [0,1])$, $B = [0,1]$, and let $p$ be the vertical projection 
$-$then $p:X \ra B$ is not a Dold fibration.]\\

\begingroup%%----------------------------------->>
\fontsize{9pt}{11pt}\selectfont
\textbf{\small EXAMPLE}  \ 
Let
$
\begin{cases}
\ p:X \ra B\\
\ q:Y \ra B
\end{cases}
$
be Hurewicz fibrations $-$then the projection $X \bboxsub_B Y \ra B$ is a Dold fibration, hence $X *_B Y$ and 
$X \bboxsub_B Y$ have the same fiber homotopy type.\\
\endgroup%%------------------------------------<<

\begingroup%%----------------------------------->>
\fontsize{9pt}{11pt}\selectfont
\textbf{\small EXAMPLE}  \ 
Let \mX be a topological space.  
Fix a numerable covering $\sU = \{U_i: i \in I\}$ of \mX 
$-$then, in the notation of p. \pageref{4.32}, the projection $p_{\sU}:B\sU \ra X$ is a Dold fibration (for $B\sU$ , as an object in \bTOP/\mX, is fiberwise contractible).\\
\endgroup%%------------------------------------<<

\begingroup%%----------------------------------->>
\fontsize{9pt}{11pt}\selectfont
Notation: Given \ $b_0 \in B$, put \ $B_0 = B - \{b_0\}$ \ and for \mX, \mY in \bTOP/\mB, write $X_0$, $Y_0$ in place of $X_{b_0}$, $Y_{b_0}$.\\
\endgroup%%------------------------------------<<

\begingroup%%----------------------------------->>
\fontsize{9pt}{11pt}\selectfont
\index{Expansion Principle}
\textbf{\small FACT \ (\un{Expansion Principle})} \ 
Let \mX be in \bTOP/\mB.  
Suppose that $p_{B_0}:X_{B_0} \ra B_0$ is a Dold fibration and $b_0$ has a halo $O \subset B$  
contractible to $b_0$, with $O - \{b_0\}$ numerably contractible.  
Assume: $r:X_O \ra X_0$ is a homotopy equivalence 
which $\forall \ b \in O$ induces a homotopy equivalence $r_b:X_b \ra X_0$ $-$then there exists a \mY in \bTOP/\mB and an embedding $X \ra Y$ over $B$ such that $q:Y \ra B$ is a Dold fibration and 
$
\begin{cases}
\ X\\
\ X_0
\end{cases}
$
is a strong deformation retract of 
$
\begin{cases}
\ Y\\
\ Y_0
\end{cases}
$
\hspace{-.25cm}.
\\ \indent
[The \cd 
\begin{tikzcd}[sep=large]
{X_0^\prime} \ar{d} \ar{r} &{X_0} \ar{d}\\
{O} \ar{r} &{b_0}
\end{tikzcd}
$(X_0^\prime = O \times X_0)$ is a pullback square.  
Since $O \ra b_0$ is a homotopy equivalence, $X_0^\prime \ra X_0$ is a homotopy equivalence (cf. p. \pageref{4.33}), thus the arrow $r^\prime:X_O \ra X_0^\prime$ defined by 
$x \ra (p(x),r(x))$ is a homotopy equivalence.  
Let $Y$ be the double mapping cylinder of the 2-source
$X \lla X_O \overset{r^\prime}{\lra} X_0^\prime$ : \mY is in \bTOP/\mB and there is an embedding $X \ra Y$ over $B$.  
It is a closed cofibration.  
$Y_O$ is the mapping cylinder of $r^\prime$, so $X_O$ is a strong deformation retract of 
$Y_O$ (cf. $\S 3$, Proposition 17).  
Therefore \mX is a strong deformation retract of \mY (cf. $\S 3$, Proposition 3).  
Similar remarks apply to $X_0$ and $Y_0$.  Finally, to see that $q$ is a Dold fibration, note that $\{O,B_0\}$ is a numerable covering of $B$.  
Accordingly, taking into account the local-global principle, it is enough to verify that 
$q_O:Y_O \ra O$ and $q_{B_0}:Y_{B_0} \ra B_0$ are Dold fibrations.  
Consider, e.g., the latter.  
The hypotheses on $r$, in conjunction
%%----------------------------------------------------------------------------------------------32
with Proposition 20, imply that the embedding $X_{B_0} \ra Y_{B_0}$ is a fiber homotopy equivalence.  
But $p_{B_0}$ is a Dold fibration, hence the same holds for $q_{B_0}$.]\\
\endgroup%%------------------------------------<<

\label{5.29}
Let $f:X \ra Y$ be a pointed continuous function $-$then the 
\un{mapping fiber}
\index{mapping fiber ($E_f$)} $E_f$ 
\index{$E_f$} of $f$ is defined by the pullback square 
\begin{tikzcd}%[sep=large]
{E_f} \ar{d} \ar{r} &{W_f} \ar{d}{q}\\
{\{x_0\}} \ar{r}[swap]{f} &{Y}
\end{tikzcd}
, i.e., $E_f$ is the double mapping track of the 2-sink $X \overset{f}{\ra} Y \la \{y_0\}$.  
Example: The mapping fiber $E_0$ of $0$ : $X \ra Y$ is $X \times \Omega Y$.\\

\begingroup%%----------------------------------->>
\fontsize{9pt}{11pt}\selectfont
\label{4.40}
\textbf{\small EXAMPLE}  \ 
Let $f:X \ra Y$ be a pointed continuous function.  
Assume: $f$ is a Hurewicz fibration.  
Denote by $C_{y_0}$ the mapping cone of the inclusion $X_{y_0} \ra X$ $-$then the mapping fiber of 
$C_{y_0} \ra Y$ has the same homotopy type as $X_{y_0} * \Omega Y$ (cf. p. \pageref{4.34} ff.).\\
\endgroup%%------------------------------------<<

\begingroup%%----------------------------------->>
\fontsize{9pt}{11pt}\selectfont
\textbf{\small FACT}  \ 
Let $X \overset{p}{\ra} B \overset{q}{\la} Y$ be a 2-sink.  
Denote by $W_{\Box}$ the mapping track of the projection 
$X \bboxsub_B Y \ra B$ $-$then $W_p *_B W_q$ and $W_{\Box}$ has the same fiber homotopy type.\\
\endgroup%%------------------------------------<<

\begingroup%%----------------------------------->>
\fontsize{9pt}{11pt}\selectfont
Application: The mapping fiber of the projection $X \bboxsub_B Y \ra B$ has the same homotopy type as 
$E_p * E_q$.\\ 
\endgroup%%------------------------------------<<

\label{18.23}
Let $f:X \ra Y$ be a pointed continuous function $-$then $W_f$ and $E_f$ are pointed spaces, the base point in either case being $(x_0,j(y_0))$.  The pointed homotopy type of $W_f$ or $E_f$ depends only on the pointed homotopy class of $f$.  The projection $q:W_f \ra Y$ is a pointed Hurewicz fibration and the restriction $\pi$ of the projection 
$p:W_f \ra X$ to $E_f$ is a pointed Hurewicz fibration with $\pi^{-1}(x_0) = \Omega Y$.  By construction, 
$f \circx \pi$ is nullhomotopic and for any $g:Z \ra X$ with $f \circx g$ nullhomotopic, there is a $\phi:Z \ra E_f$ such that 
$g = \pi \circx \phi$.\\

\begingroup%%----------------------------------->>
\fontsize{9pt}{11pt}\selectfont
When is a pointed continuous function which is a Hurewicz fibration actually a pointed Hurewicz fibration?  
Regularity, suitably localized, is what is relevant.  
Thus let $p:X \ra B$ be a Hurewicz fibration taking $x_0$ to $b_0$.  
Assume: $\exists$ a lifting function $\Lambda$ such that $\Lambda(x_0,j(b_0)) = j(x_0)$ $-$then $p$ is a pointed Hurewicz fibration.
\\ \indent
[Note: \ For this, it is sufficient that $\{b_0\}$ be a zero set in $B$, 
any Hurewicz fibration $p:X \ra B$ automatically becoming a pointed Hurewicz fibration $\forall \ x_0 \in X_{b_0}$ (argue as on p. \pageref{4.35}).  
The condition is satisfied if the inclusion $\{b_0\} \ra B$ is a closed cofibration.]\\
\endgroup%%------------------------------------<<

\label{4.39}
\begingroup%%----------------------------------->>
\fontsize{9pt}{11pt}\selectfont
\textbf{\small LEMMA}  \ 
Let \mX, \mY, \mZ be pointed spaces; let
$
\begin{cases}
\ f:X \ra Z\\
\ g:Y \ra Z
\end{cases}
$
be pointed continuous functions $-$then  the projections 
$
\begin{cases}
\ W_{f,g} \ra X\\
\ W_{f,g} \ra Y
\end{cases}
$
$\&$ $W_{f,g} \ra X \times Y$ are pointed Hurewicz fibrations, the base point of $W_{f,g}$ being the triple 
$(x_0,y_0,j(z_0))$.
\\ \indent
%%----------------------------------------------------------------------------------------------33
[To deal with $p:W_{f,g} \ra X$, define a lifting function 
$\Lambda:W_p \ra PW_{f,g}$ by $\Lambda((x,y,\tau),\sigma)(t) =$
$(\sigma(t),y,\tau_t)$, where 
\[
\tau_t(T) = 
\begin{cases}
\ f \circx \sigma (t - 2T) \hspace{0.75cm} (0 \leq T \leq t/2)\\
\ \tau\left(\ds\frac{2T - t}{2 - t}\right) \hspace{0.80cm}   \ (t/2 \leq T \leq 1)
\end{cases}
.
\]
Obviously, 
$\Lambda((x_0,y_0,j(z_0)),j(x_0)) =$
$j(x_0,y_0,j(z_0))$, so $p:W_{f,g} \ra X$ is a pointed Hurewicz fibration.]\\
\endgroup%%------------------------------------<<

\begin{proposition} \  %24
Consider the pullback square
\begin{tikzcd}%[sep=large]
{X^\prime} \ar{d} \ar{r} &{X} \ar{d}{p}\\
{B^\prime} \ar{r}[swap]{\Phi^\prime} &{B}
\end{tikzcd}
, where $p$ is a Hurewicz fibration.  Suppose that 
$
\begin{cases}
\ X\\[-.15cm]
\ B
\end{cases}
\& \ B^\prime
$
are wellpointed, that the inclusions 
$
\begin{cases}
\ \{x_0\} \ra X\\[-.15cm]
\ \{b_0\} \ra B
\end{cases}
$
$\&$ $\{b_0^\prime\} \ra B^\prime$ are closed, and that $p(x_0) = b_0 = \Phi^\prime(b_0^\prime)$.  
Put $x_0^\prime = (b_0^\prime,x_0)$ 
$-$then  the inclusion 
$\{x_0^\prime\} \ra X^\prime$ is a closed cofibration.
\end{proposition}

[The arrow $X_{b_0} \ra X$ is a closed cofibration (cf. Proposition 11).  
Therefore the composite 
$X_{b_0^\prime}^\prime \ra X^\prime \ra X$ is a closed cofibration.  
On the other hand, the composite 
$\{x_0^\prime\} \ra X^\prime_{b_0^\prime} \ra X^\prime \ra X$ is a closed cofibration.  
Therefore the inclusion $\{x_0^\prime\} \ra X^\prime_{b_0^\prime}$ is a closed cofibration (cf. $\S 3$, Proposition 9).  
But the arrow $X^\prime_{b_0^\prime} \ra X^\prime$ is a closed cofibration (cf. Proposition 11), thus the inclusion 
$\{x_0^\prime\} \ra X^\prime$ is a closed cofibration.]\\

Application: Let $f:X \ra Y$ be a pointed continuous function.  Assume: 
$
\begin{cases}
\ X\\[-.15cm]
\ Y
\end{cases}
$
are wellpointed with closed base points $-$then $W_f$ and $E_f$ are wellpointed with closed base points.

[$PY$ is wellpointed with a closed base point (cf. $\S 3$, Proposition 6).]\\

\begingroup%%----------------------------------->>
\fontsize{9pt}{11pt}\selectfont
\label{5.27f}
\textbf{\small FACT}  \ 
Let $f:X \ra Y$ be a pointed continuous function.  Suppose that 
$\phi:X^\prime \ra X$ $(\psi:Y \ra Y^\prime)$ is a pointed homotopy equivalence $-$then the arrow 
$E_{f \circx \phi} \ra E_f$ $(E_f \ra E_{\psi \circx f})$ is a pointed homotopy equivalence.\\
\endgroup%%------------------------------------<<

\label{9.118}
\begingroup%%----------------------------------->>
\fontsize{9pt}{11pt}\selectfont
Application: Let $X$ be wellpointed with $\{x_0\} \subset X$ closed $-$then the mapping fiber of the diagonal embedding $X \ra X \times X$ has the same pointed homotopy type as $\Omega X$.
\\ \indent
[The embedding $j:X \ra PX$ is a pointed homotopy equivalence and 
$
\Pi:
\begin{cases}
\ PX \ra X \times X\\
\ \sigma \ra (\sigma(0),\sigma(1))
\end{cases}
$
is a pointed Hurewicz fibration.]\\
\endgroup%%------------------------------------<<

\begingroup%%----------------------------------->>
\fontsize{9pt}{11pt}\selectfont
\textbf{\small EXAMPLE}  \ 
Let 
$
\begin{cases}
\ X\\
\ Y
\end{cases}
$
be wellpointed with 
$
\begin{cases}
\ \{x_0\} \subset X\\
\ \{y_0\} \subset Y
\end{cases}
$
closed.\\
\indent\indent (1) \ The mapping fiber of the inclusion  $X \vee Y \ra X \times Y$ has the same pointed homotopy type as 
$\Omega X * \Omega Y$.\\
%%----------------------------------------------------------------------------------------------34
\indent\indent (2) \ The mapping fiber of the projection $X \vee Y \ra Y$ has the same pointed homotopy type as 
$X \times \Omega Y / \{x_0\} \times \Omega Y$.

[In both situations, replace $\Theta$ by $\Gamma \Omega$ as on p. \pageref{4.36}.]\\
\endgroup%%------------------------------------<<

\label{14.152} %dmc mnft
\begingroup%%----------------------------------->>
\fontsize{9pt}{11pt}\selectfont
\textbf{\small FACT}  \ 
Let 
$
\begin{cases}
\ f:X \ra Y\\
\ g:Y \ra Z
\end{cases}
$
be pointed continuous functions $-$then there is a homotopy equivalence $E_{g \circx f} \ra W$, where $W$ is the double mapping track of the 2-sink $X \overset{f}{\ra} Y \overset{\pi}{\la} E_g$.
\\ \indent
[Consider the diagram
\begin{tikzcd}[sep=large]
{E_{g \circx f}} \ar{d} \ar{r} &{E_g} \ar{d} \ar{r} &{*} \ar{d}\\
{X} \ar{r} &{Y} \ar{r} &{Z}
\end{tikzcd}
.]\\
\endgroup%%------------------------------------<<

\label{4.51}
Let $f:X \ra Y$ be a pointed continuous function, $E_f$ its mapping fiber.\\

\textbf{\small LEMMA}  \ 
If $f$ is a pointed Hurewicz fibration, then the embedding $X_{y_0} \ra E_f$ is a pointed homotopy equivalence.\\

\label{4.49}
In general, there is a pointed Hurewicz fibration $\pi:E_f \ra X$ and an embedding $\Omega Y \ra E_f$.  
Iterate to get a pointed Hurewicz fibration $\pi^\prime:E_\pi \ra E_f$ $-$then the triangle 
\begin{tikzcd}%[sep=large]
{E_\pi} \ar{r} &{E_f}\\
{\Omega Y} \ar {u} \ar{ru}
\end{tikzcd}
commutes and by the lemma, the vertical arrow is a pointed homotopy equivalence.  Iterate again to get a pointed Hurewicz fibration $\pi\pp:E_{\pi^\prime} \ra E_\pi$ $-$then the triangle 
\begin{tikzcd}%[sep=large]
{E_{\pi^\prime}} \ar{r} &{E_\pi}\\
{\Omega X} \ar {u} \ar{ru}
\end{tikzcd}
commutes and by the lemma, the vertical arrow is a pointed homotopy equivalence.  
Example:  Given pointed spaces 
$
\begin{cases}
\ X\\[-.15cm]
\ Y
\end{cases}
,
$
let $X \flat Y$ be the mapping fiber of the inclusion $f:X \vee Y \ra X \times Y$ $-$then in $\bHTOP_*$, 
$E_\pi \approx \Omega (X \times Y)$ and $E_{\pi^\prime} \approx \Omega(X \vee Y)$.\\

\begingroup%%----------------------------------->>
\fontsize{9pt}{11pt}\selectfont
\textbf{\small LEMMA}  \ 
Let 
$
\begin{cases}
\ X\\
\ Y
\end{cases}
$
be wellpointed with 
$
\begin{cases}
\ \{x_0\} \subset X\\
\ \{y_0\} \subset Y
\end{cases}
$
closed.  Denote by $S$ the subspace of $X *Y$ consisting of the 
$
\begin{cases}
\ [x,y_0,t]\\
\ [x_0,y,t]
\end{cases}
$
$-$then $X * Y/S = \Sigma(X\#Y)$ and the projection $X *Y \ra X*Y /  S$ is a pointed homotopy equivalence.
\\ \indent
[Note: \ The base point of $X * Y$ is $[x_0,y_0,1/2]$ and $\Sigma$ is the pointed suspension.]\\
\endgroup%%------------------------------------<<

\label{4.50}
\begingroup%%----------------------------------->>
\fontsize{9pt}{11pt}\selectfont
Application: Let 
$
\begin{cases}
\ X\\
\ Y
\end{cases}
$
be wellpointed with 
$
\begin{cases}
\ \{x_0\} \subset X\\
\ \{y_0\} \subset Y
\end{cases}
$
closed $-$then $X \flat Y$ has the same pointed homotopy type as $\Sigma(\Omega X \# \Omega Y)$.\\
\endgroup%%------------------------------------<<

\begingroup%%----------------------------------->>
\fontsize{9pt}{11pt}\selectfont
\textbf{\small EXAMPLE}  \  
Suppose that \mX and \mY are nondegenerate $-$then the Puppe formula says that in $\bHTOP_*$, 
$\Sigma(\Omega X \times \Omega Y) \approx$ 
$\Sigma\Omega X \vee \Sigma\Omega Y \vee \Sigma(\Omega X \# \Omega Y)$, and by the above, 
$\Sigma(\Omega X \# \Omega Y) \approx$ 
$X \flat Y$.\\
\endgroup%%------------------------------------<<

%%----------------------------------------------------------------------------------------------35
\begingroup%%----------------------------------->>
\fontsize{9pt}{11pt}\selectfont
\textbf{\small EXAMPLE \ (\un{The Flat Product})} \ 
In contrast to the smash product $\#$ (or its modification $\ov{\#}$), the flat product $\flat$ does not possess the properties that one might expect to hold by analogy.  
Specifically, for nondegenerate spaces, it is generally false that in
$\bHTOP_*$: 
(1) $(X \flat Y)\flat Z \approx X \flat (Y \flat Z)$; 
(2) $(X \times Y)\flat Z \approx (X \flat Z) \times (Y \flat Z)$; 
(3) $\Omega(X \flat Y ) \approx \Omega X \flat Y$.
Counterexamples: 
(1) Take $X = Y = \bP^\infty (\C)$, $Z = \bP^\infty (\bH)$; 
(2) Take $X = Y = Z = \bP^\infty (\C)$;
(3) Take $X = Y = \bP^\infty (\C)$.
Look, e.g., at (1).  
Using the fact that 
$\Omega \bP^\infty (\C) \approx \bS^1$, 
$\Omega \bP^\infty (\bH)  \approx \bS^3$, 
compute: 
$\bP^\infty (\C) \flat \bP^\infty (\C)  \approx$ 
$\Omega \bP^\infty (\C) * \Omega \bP^\infty (\C) \approx$ 
$\bS^1 * \bS^1 \approx$ 
$\bS^3$ $\&$
$\bS^3 \flat \bP^\infty (\bH) \approx$ 
$\Omega \bS^3 * \bS^3 \approx$ 
$\Sigma (\Omega \bS^3 \# \bS^3) \approx$ 
$\Sigma \Omega \bS^3 \# \bS^3 \approx$ 
$\Sigma \Omega \bS^3 \# \Sigma^3\bS^0 \approx$ 
$\Sigma^4 \Omega \bS^3 \# \bS^0 \approx$
$\Sigma^4 \Omega \bS^3$ 
$\implies$ 
$(\bP^\infty (\C) \flat \bP^\infty (\C))\flat \bP^\infty (\bH) \approx$
$\Sigma^4 \Omega \bS^3$ 
Similarly, 
$\bP^\infty (\C) \flat (\bP^\infty (\C)\flat \bP^\infty (\bH)) \approx$
$\Sigma^2 \Omega \bS^5$.  
The singular homology functor $H_8(-;\Z)$ distinguishes these spaces: 
$H_8(\Sigma^4 \Omega \bS^3;\Z) \approx \Z$, 
$H_8(\Sigma^2 \Omega \bS^5;\Z) \approx 0$.\\
\endgroup%%------------------------------------<<

Let $f:X \ra Y$ be a pointed continuous function $-$then the \un{mapping fiber sequence} 
\index{mapping fiber sequence} associated with $f$ is given by \ 
$\cdots \ra \Omega^2 Y \ra \Omega E_f \ra \Omega X \ra \Omega Y \ra$ $E_f \ra X \overset{f}{\ra} Y$.  
Example: When $f = 0$, this sequence becomes 
$\cdots \ra \Omega^2 Y \ra \Omega X \times \Omega^2 Y \ra \Omega X \ra$ $\Omega Y \ra$ 
$X \times \Omega Y \ra$ 
$X \overset{0}{\ra} Y$.

[Note: \ If the diagram 
\begin{tikzcd}%[sep=large]
{X} \ar{d} \ar{r}{f} &{Y} \ar{d}\\
{X^\prime} \ar{r}[swap]{f^\prime} &{Y^\prime}
\end{tikzcd}
commutes in $\bHTOP_*$ and if the vertical arrows are pointed homotopy equivalences, then the mapping fiber sequences of $f$ and $f^\prime$ are connected by a commutative ladder in $\bHTOP_*$, all of whose vertical arrows are pointed homotopy equivalences.]\\

\label{6.18a}
\begingroup%%----------------------------------->>
\fontsize{9pt}{11pt}\selectfont
\textbf{\small FACT}  \ 
Let $f:X \ra Y$ be a pointed Hurewicz fibration.  
Assume: The inclusion $X_{y_0} \ra X$ is nullhomotopic $-$then 
$\Omega Y$ has the same pointed homotopy type as $X_{y_0} \times \Omega X$.
\\ \indent
[For $\pi:E_f \ra X$ is nullhomotopic, thus in $\bHTOP_*$: $E_\pi \approx E_f \times \Omega X$ 
$\implies$ $\Omega Y \approx X_{y_0} \times \Omega X$.]\\
\endgroup%%------------------------------------<<

\index{Replication Theorem}
\index{Theorem: Replication Theorem}
\textbf{\small REPLICATION THEOREM} \quadx
Let $f:X \ra Y$ be a pointed continuous function $-$then for any pointed space \mZ, there is an exact sequence
\[
\cdots \ra [Z,\Omega X] \ra [Z,\Omega Y] \ra [Z,E_f] \ra [Z,X] \ra [Z,Y] 
\]
in $\bSET_*$.\\

If $f:X \ra Y$ is a pointed Dold fibration or if $f:X \ra Y$ is a Dold fibration and \mZ is nondegenerate, then in the replication theorem one can replace $E_f$ by $X_{y_0}$ (cf. p. \pageref{4.37}).  This replacement can also be made if 
$f:X \ra Y$ is a Serre fibration provided that \mZ is a CW complex (cf. infra).  In particular, when $f:X \ra Y$ is either a Dold fibration or a Serre fibration, there is an exact sequence 
\[
\cdots \ra \pi_2(Y) \ra \pi_1(X_{y_0}) \ra \pi_1(X) \ra \pi_1(Y) \ra \pi_0(X_{y_0}) \ra \pi_0(X) \ra \pi_0(Y).
\]

%%----------------------------------------------------------------------------------------------36
\begingroup%%----------------------------------->>
\fontsize{9pt}{11pt}\selectfont
\textbf{\small LEMMA}  \ 
Let $f:X \ra Y$ be a pointed continuous function.  Assume: $f$ is a Serre fibration $-$then for every pointed CW complex $Z$, the arrow $[Z,X_{y_0}] \ra [Z,E_f]$ is a pointed bijection.
\\ \indent
[Proposition 12 is true for Serre fibrations if the ``cofibration data'' is restricted to CW complexes.]\\
\endgroup%%------------------------------------<<

Examples: Suppose that $f:X \ra Y$ is either a Dold fibration or a Serre fibration, where 
$
\begin{cases}
\ X \ne \emptyset\\[-.15cm]
\ Y \ne \emptyset
\end{cases}
$
.  
(1) If $X_{y_0}$ is simply connected, then $\forall \ x_0 \in X_{y_0}$, $\pi_1(X,x_0) \approx \pi_1(Y,y_0)$;
(2) If \mX is simply connected, then $\forall \ y_0 \in f(X)$, there is a bijection $\pi_1(Y,y_0) \ra \pi_0(X_{y_0})$;
(3) \ If \mX is path connected and if \mY is simply connected, then $\forall \ y_0 \in Y$, $\pi_0(X_{y_0}) = *$; 
(4) \ If \mY is path connected and $X_{y_0}$ is path connected, then \mX is path connected.\\

\begingroup%%----------------------------------->>
\fontsize{9pt}{11pt}\selectfont
\textbf{\small LEMMA}  \ 
Let $f:X \ra Y$ be a Hurewicz fibration.  Fix $y_0 \in f(X)$ $\&$ $x_0 \in X_{y_0}$ and let $(Z,z_0)$ be wellpointed with 
$\{z_0\} \subset Z$ closed $-$then there is a left action 
$\pi_1(X,x_0) \times [Z,z_0;X_{y_0},x_0]  \ra [Z,z_0;X_{y_0},x_0]$.
\\ \indent
[Represent $\alpha \in \pi_1(X,x_0)$ by a loop $\sigma \in \Omega(X,x_0)$.  Given 
$\phi:(Z,z_0) \ra (X_{y_0},x_0)$, consider the commutative diagram
\begin{tikzcd}[sep=large]
{i_0Z \cup I\{z_0\}} \ar{d} \ar{r}{F} &{X} \ar{d}{f}\\
{IZ} \ar{r}[swap]{h} &{Y}
\end{tikzcd}
, where 
$
F(z,t) = 
\begin{cases}
\ (i \circx \phi)(z) \hspace{0.5cm} (t = 0)\\
\ \sigma(t) \hspace{0.97cm} \quadx (z = z_0)
\end{cases}
$
($i$ the inclusion $X_{y_0} \ra X$) and $h(z,t) = (f \circx \sigma)(t)$.  Proposition 12 says that this diagram has a filler 
$H:IZ \ra X$.  Put $\psi(z) = H(z,1)$ to get a pointed continuous function $\psi:(Z,z_0) \ra (X_{y_0},x_0)$.  
Definition: $\alpha \cdot [\phi] = [\psi]$.]
\\ \indent
[Note: \ There is a left action 
$\pi_1(X,x_0) \times [Z,z_0;X,x_0] \ra [Z,z_0;X,x_0]$ 
and a left action 
$\pi_1(X_{y_0},x_0) \times [Z,z_0;X_{y_0},x_0] \ra [Z,z_0;X_{y_0},x_0]$ (cf. p. \pageref{4.38}).  The arrow
$[Z,z_0;X_{y_0},x_0] \ra [Z,z_0;X,x_0]$ induced by the inclusion $X_{y_0} \ra X$ is a morphism of $\pi_1(X,x_0)$-sets 
and the operation of $\pi_1(X_{y_0},x_0)$ on $[Z,z_0;X_{y_0},x_0]$ coincides with that defined via the homomorphism 
$\pi_1(X_{y_0},x_0) \ra \pi_1(X,x_0)$.]\\
\endgroup%%------------------------------------<<

\label{5.16}
\begingroup%%----------------------------------->>
\fontsize{9pt}{11pt}\selectfont
\textbf{\small EXAMPLE}  \ 
Let $f:X \ra Y$ be a Hurewicz fibration.  Fix $y_0 \in f(X)$ $\&$ $x_0 \in X_{y_0}$ and $n \geq 1$ $-$then there is a left action 
$\pi_1(X,x_0) \times \pi_n(X,x_0) \ra \pi_n(X,x_0)$, 
and a left action 
$\pi_1(X,x_0) \times \pi_n(Y,y_0) \ra \pi_n(Y,y_0)$, 
and a left action 
$\pi_1(X,x_0) \times \pi_n(X_{y_0},x_0) \ra \pi_n(X_{y_0},x_0)$.  All the homomorphisms in the exact sequence 
\[
\cdots \ra \pi_{n+1}(Y,y_0) \ra \pi_n(X_{y_0},x_0) \ra \pi_n(X,x_0) \ra \pi_n(Y,y_0) \ra \cdots
\]
are $\pi_1(X,x_0)$-homomorphisms.
\\ \indent
[Note: \ Suppose that $X_{y_0}$ is path connected $-$then there is a left action 
$\pi_1(Y,y_0) \times \pi_n^*(X_{y_0},x_0) \ra \pi_n^*(X_{y_0},x_0)$, where 
$\pi_n^*(X_{y_0},x_0)$ is $\pi_n(X_{y_0},x_0)$ modulo the (normal) subgroup generated by the 
$\alpha \cdot \xi - \xi$ ($\alpha \in \pi_1(X_{y_0},x_0)$, $\xi \in \pi_n(X_{y_0},x_0)$).]\\
\endgroup%%------------------------------------<<

\begingroup%%----------------------------------->>
\fontsize{9pt}{11pt}\selectfont
\textbf{\small EXAMPLE}  \ 
Let $f:X \ra Y$ be a Hurewicz fibration.  
Fix $y_0 \in f(X)$ $\&$ $x_0 \in X_{y_0}$ $-$then $\pi_1(Y,y_0)$ operates to the left on $\pi_0(X_{y_0})$ 
and the orbits are the fibers of the arrow $\pi_0(X_{y_0}) \ra \pi_0(X)$.\\
\endgroup%%------------------------------------<<

\begingroup%%----------------------------------->>
\fontsize{9pt}{11pt}\selectfont
\textbf{\small FACT}  \ 
Let $f:X \ra Y$ be a Hurewicz fibration.  Fix $y_0 \in f(X)$ $\&$ $x_0 \in X_{y_0}$ $-$then $\forall \ n \geq 1$, 
$\pi_1(X_{y_0},x_0)$ operates trivially on $\ker(\pi_n(X_{y_0},x_0) \ra \pi_n(X,x_0))$.\\
\endgroup%%------------------------------------<<

%%----------------------------------------------------------------------------------------------37
\label{5.17}
\label{9.19}
\label{18.25}
\index{Mayer-Vietoris Sequence (example)}
\begingroup%%----------------------------------->>
\fontsize{9pt}{11pt}\selectfont
\textbf{\small EXAMPLE \ (\un{Mayer-Vietoris Sequence})} \ 
Let \mX, \mY, \mZ be pointed spaces; let 
$
\begin{cases}
\ f:X \ra Z\\
\ g:Y \ra Z
\end{cases}
$
be pointed continuous functions $-$then the projection $W_{f,g} \ra X \times Y$ is a pointed Hurewicz fibration 
(cf. p. \pageref{4.39}) and there is a long exact sequence 
$\cdots \ra $
$\pi_{n+1}(Z) \ra$ 
$\pi_{n}(W_{f,g}) \ra$ 
$\pi_{n}(X) \times \pi_{n}(Y) \ra$ 
$\pi_{n}(Z) \ra$ 
$\cdots \ra$
$\pi_{2}(Z) \ra$
$\pi_{1}(W_{f,g}) \ra$ 
$\pi_{1}(X) \times \pi_{1}(Y) \ra$ 
$\pi_1(Z) \ra $ 
$\pi_{0}(W_{f,g}) \ra$ 
$\pi_{0}(X \times Y)$.
\\ \indent
[Note: \ It follows that if \mX and \mY are path connected and if every $\gamma \in \pi_1(Z)$ has the form 
$\gamma = f_*(\alpha) \cdot g_*(\beta)$ $(\alpha \in \pi_1(X)$, $\beta \in \pi_1(Y))$, then 
$W_{f,g}$ is path connected.]\\
\endgroup%%------------------------------------<<

If $f:X \ra Y$ is either a Dold fibration or a Serre fibration, then the homotopy groups of \mX and \mY are related to those of the fibers by a long exact sequence.  
As for the homology groups, there is still a connection but it is intricate and is best expressed in terms of a spectral sequence.

[Note: \ In the simplest case, viz. that of a projection $Y \times T \ra Y$, the K\"unneth formula computes the homology of $Y \times T$ in terms of the homology of \mY and the homology of $T$.]\\

\begingroup%%----------------------------------->>
\fontsize{9pt}{11pt}\selectfont
\label{4.48}
\textbf{\small EXAMPLE}  \ 
Let $f:X \ra Y$ be a Hurewicz fibration, where \mX is nonempty and \mY is path connected.  
Fix $y_0 \in Y$ $-$then 
$\forall \ q \geq 1$, the projection $(X,X_{y_0}) \ra (Y,y_0)$ induces a bijection 
$\pi_q(X,X_{y_0}) \ra \pi_q(Y,y_0)$.  
The analog of this in homology is false.  
Consider, e.g., the Hopf map 
$\bS^{2n+1} \ra \bP^n(\C)$ with fiber $\bS^1:H_q(\bS^{2n+1},\bS^1) = 0$ $(2 < q \leq 2n)$ $\&$ 
$H_{2q}(\bP^n(\C)) \approx \Z$ $(1 < q \leq n)$.  
However, a partial result holds in that if $X_{y_0}$ is $n$-connected 
and \mY is $m$-connected, then the arrow 
$H_q(X,X_{y_0}) \ra H_q(Y,y_0)$ induced by the projection 
$(X,X_{y_0}) \ra (Y,y_0)$ is bijective  for $1 \leq q < n + m + 2$ and surjective for $q = n + m + 2$.  
Consequently, under these conditions, there is an exact sequence
\endgroup%%------------------------------------<< % y ask y
\\[-1.35cm]

\begingroup%%----------------------------------->>
\fontsize{9pt}{11pt}\selectfont
\begin{align*}
&H_{n+m+1}(X_{y_0}) \ra H_{n+m+1}(X) \ra H_{n+m+1}(Y) \ra H_{n+m}(X_{y_0}) \ra \cdots \\
&\hspace{1.75cm} \ra H_2(Y) \ra H_1(X_{y_0}) \ra H_1(X) \ra H_1(Y).
\end{align*}
\\[-.85cm]
\indent
[One can assume that the inclusion $\{y_0\} \ra Y$ is a closed cofibration (pass to a CW resolution $K \ra Y$), 
hence that the inclusion $X_{y_o} \ra X$ is a closed cofibration (cf. Proposition 11).  
The mapping cone of the latter is path connected and the mapping fiber of $C_{y_0} \ra Y$ has the same homotopy type as 
$X_{y_0} * \Omega Y$ 
(cf. p. \pageref{4.40}), which is $(n + m + 1)$-connected (cf. p. \pageref{4.41}).  
Thus the arrow 
$C_{y_0} \ra Y$ is an $(n + m + 2)$-equivalence, so the Whitehead theorem implies that the induced map 
$H_q(C_{y_0}) \ra H_q(Y)$ is bijective for $0 \leq q < n + m + 2$ and surjective for $q = n + m + 2$.  
But the projection $C_{y_0} \ra X / X_{y_0}$ is a homotopy equivalence (cf. p. \pageref{4.42}) and 
$H_q(X,X_{y_0}) \approx H_q(X / X_{y_0},*)$ (cf. p. \pageref{4.43}).]\\
\endgroup%%------------------------------------<<

\begingroup%%----------------------------------->>
\fontsize{9pt}{11pt}\selectfont
Application: Suppose that $X$ is $(n+1)$-connected $-$then 
$H_q(X) \approx H_{q-1}(\Omega X)$ $(2 \leq q \leq 2n+2)$.  
\label{5.13}
\\ \indent
[Note: \ It is a corollary that if \mX is nondegenerate and $n$-connected, then the arrow of adjunction 
$e:X \ra \Omega\Sigma X$ induces an isomorphism $H_q(X) \ra H_q(\Omega\Sigma X)$ for $0 \leq q \leq 2n+1$.  
Therefore, by the Whitehead theorem, the suspension homomorphism 
$\pi_q(X) \ra \pi_{q+1}(\Sigma X)$ is bijective for $0 \leq q \leq 2n$ and surjective for $q = 2n + 1$ (Freudenthal).]\\
\endgroup%%------------------------------------<<

%%----------------------------------------------------------------------------------------------38
\label{13.125}
Let \mX be a topological space, $\sin X$ its singular set $-$then $\sin X$ can be regarded as a category: 
\begin{tikzcd}%[sep=large]
{\dpm} \ar{rd} \ar{rr}{\Delta^\alpha} &&{\dpn} \ar{ld}\\
&{X}
\end{tikzcd}
$(\alpha \in \Mor([m],[n]))$.  The objects of $[(\sin X)^\OP,\bAB]$ are called 
\un{coefficient systems}
\index{coefficient systems} 
on \mX.  
Given a coefficient system $\sG$, the singular homology $H_*(X;\sG)$ of \mX with coefficients in $\sG$ is by definition the homology of the chain complex
\[
\bigoplus\limits_{\sigma_0 \in \sin_0 X} \sG\sigma_0 \overset{\partial}{\lla} 
\bigoplus\limits_{\sigma_1 \in \sin_1 X} \sG\sigma_1 \overset{\partial}{\lla} 
\bigoplus\limits_{\sigma_2 \in \sin_2 X} \sG\sigma_2 \overset{\partial}{\lla} 
\cdots, 
\]
where $\partial = \ds\sum\limits_0^n (-1)^i \ds\bigoplus\limits_{\sigma_n \in \sin_n X} \sG d_i$.

[Note: \ To interpret $\sG d_i$, recall that there are arrows
$d_i:\sin_n X \ra \sin_{n-1} X$ corresponding to the face operators 
$\delta_i:[n-1] \ra [n]$ $(0 \leq i \leq n)$.  
So, $\forall \ \sigma \in \sin_n X$, 
$\sG d_i:\sG(\dpn \overset{\sigma}{\ra} X) \ra \sG(\Delta^{n-1} \overset{d_i\sigma}{\lra} X)$.]

Example: Fix an abelian group \mG and define $\sG_G$ by 
$
\begin{cases}
\ \sG_G \sigma = G\\[-.15cm]
\ \sG_G \Delta^\alpha = \id_G
\end{cases}
$
$-$then $H_*(X;\sG_G) = H_*(X;G)$, the singular homology of \mX with coefficients in \mG.

A coefficient system $\sG$ is said to be 
\un{locally constant}
\index{locally constant (coefficient system)}
provided that $\forall \ \alpha$, $\sG \Delta^\alpha$ is invertible.  
$\bLCCS_X$ 
\index{$\bLCCS_X$} 
is the full subcategory of $[(\sin X)^\OP,\bAB]$ whose objects are the locally constant coefficient systems on \mX.

[Note: \ A coefficient system $\sG$ is said to be 
\un{constant}
\index{constant (coefficient system)}
if for some abelian group \mG, $\sG$ is isomorphic to $\sG_G$.]\\

\begingroup%%----------------------------------->>
\fontsize{9pt}{11pt}\selectfont
Suppose that \mX is locally path connected and locally simply connected $-$then the category of locally constant coefficient systems on \mX is equivalent to the category of locally constant sheaves of abelian groups on \mX.\\
\endgroup%%------------------------------------<<

\begin{proposition} \  
$\bLCCS_X$ is equivalent to $[(\Pi X)^\OP,\bAB]$.
\end{proposition}

[We shall define a functor $\sG \ra \sG_\Pi$ from 
$\bLCCS_X$ to $[(\Pi X)^\OP,\bAB]$ and a functor $\sG \ra \sG_{\sin}$ from 
$[(\Pi X)^\OP,\bAB]$ to $\bLCCS_X$ such that 
$
\begin{cases}
\ (\sG_\Pi)_{\sin} \approx \sG\\[-.15cm]
\ (\sG_{\sin})_\Pi \approx \sG
\end{cases}
\hspace{-.25cm}. \ 
$
\\

Definition of $\sG_\Pi$: Given $x \in X$, put $\sG_{\Pi}x = \sG\sigma_x$, where $\sigma_x \in \sin_0 X$ with 
$\sigma_x(\dpz) = x$.  
Given a morphism $[\sigma]:x \ra y$, put 
$\sG_\Pi[\sigma] = (\sG d_1) \circ(\sG d_0)^{-1}$, where $\sigma \in \sin_1 X$ with 
$
\begin{cases}
\ d_1 \sigma = x\\[-.15cm]
\ d_0 \sigma = y
\end{cases}
\hspace{-.25cm}. \ 
$
In other words, $\sG_\Pi[\sigma] $ is the composite 
$\sG y \ra \sG \sigma \ra \sG x$.  
Note that $\sG_\Pi[\sigma] $ is welldefined.  Indeed, if 
$
\begin{cases}
\ \sigma^\prime\\[-.15cm]
\ \sigma\pp
\end{cases}
\in \sin_1 X
$
with 
$
\begin{cases}
\ d_1\sigma^\prime = x = d_1\sigma\pp\\[-.15cm]
\ d_0\sigma^\prime = y = d_0\sigma\pp
\end{cases}
$
and $[\sigma^\prime] = [\sigma\pp]$, then there exists a $\tau \in \sin_2 X$ such that 
$
\begin{cases}
\ d_1\tau = \sigma^\prime\\[-.15cm]
\ d_2 \tau = \sigma\pp
\end{cases}
$
and $s_0d_0\sigma^\prime = d_0\tau = s_0d_0\sigma\pp$.

%%----------------------------------------------------------------------------------------------39
Definition of $\sG_{\sin}$: Given $\sigma \in \sin_n X$, put 
$\sG_{\sin}\sigma = \sG(e_n\sigma(\dpz))$, where $e_n:\sin_n X \ra \sin_0 X$ 
is the arrow associated with the vertex operator $e_n:[0] \ra [n]$ that sends 0 to $n$.
Given a morphism
\begin{tikzcd}%[sep=large]
{\dpm} \ar{rd}[swap]{r} \ar{rr}{\Delta^\alpha} &&{\dpn} \ar{ld}{\sigma}\\
&{X}
\end{tikzcd}
, \ put $\sG_{\sin}\Delta^\alpha = \sG(\sigma \circx \Delta^\theta)$, where $\theta:[1] \ra [n]$ is defined by 
$
\begin{cases}
\ \theta(0) = \alpha(m)\\[-.15cm]
\ \theta(1) = n
\end{cases}
\hspace{-.25cm}:
$
$\sigma \circx \Delta^\theta$ is a path in \mX which begins at $e_m\tau(\dpz)$ and ends at $e_n\sigma(\dpz)$.]\\

\label{5.0ao}
\label{5.34b}
Because of this result, one can always pass back and forth between locally constant coefficient systems on \mX and cofunctors $\Pi X \ra \bAB$.  The advantage of dealing with the latter is that in practice a direct description is sometimes available.  For example, fix $n \geq 2$ and assign to each $x \in X$ the homotopy group $\pi_n(X,x)$ $-$then every morphism $[\sigma]:x \ra y$ determines an isomorphism $\pi_n(X,y) \ra \pi_n(X,x)$ and there is a cofunctor 
$\pi_n X: \Pi X \ra \bAB$.

[Note: \ Suppose that $\sG$ is in $[(\Pi X)^\OP,\bAB]$ $-$then $\forall \ x_0 \in X$, the fundamental group 
$\pi_1(X,x_0)$ operates to the right on 
$\sG x_0:\sG x_0 \times \pi_1(X,x_0) \ra \sG x_0$.  
Conversely, if \mX is path connected and if $G_0$ is an abelian group on which $\pi_1(X,x_0)$ operates to the right, then there exists a $\sG$ in $[(\Pi X)^\OP,\bAB]$, unique up to isomorphism, with $\sG x_0 = G_0$ and inducing the given operation of $\pi_1(X,x_0)$ on $G_0$.]\\

Application: On a simply connected space, every locally constant coefficient system is isomorphic to a constant coefficient system.\\

\begingroup%%----------------------------------->>
\fontsize{9pt}{11pt}\selectfont
\label{4.45}
\textbf{\small EXAMPLE}  \ 
Let $f:X \ra Y$ be a Hurewicz fibration $-$then $\forall \ q \geq 0$, there is a cofunctor 
$\sH_q(f):\Pi Y \ra \bAB$ that assigns to each $y \in Y$ the singular homology group $H_q(X_y)$ of the fiber $X_y$.  
Thus let $[\tau]:y_0 \ra y_1$ be a morphism.  
Case 1: 
$
\begin{cases}
\ y_0 \\
\ y_1
\end{cases}
\hspace{-.25cm}\notin f(X).
$
In this situation, $X_{y_0}$ $\&$ $X_{y_1}$ are empty, hence $H_q(X_{y_0}) = 0 = H_q(X_{y_1})$.  
Definition: $\sH_q(f)[\tau]$ is the zero morphism.
Case 2: 
$
\begin{cases}
\ y_0 \\
\ y_1
\end{cases}
\hspace{-.25cm} \in f(X).
$
Fix a homotopy $\Lambda:IX_{y_0} \ra X$ such that 
$
\begin{cases}
\ f \circx \Lambda(x,t) = \tau(t)\\
\ \Lambda(x,0) = x
\end{cases}
$
$-$then the arrow 
$
\begin{cases}
\ X_{y_0} \ra X_{y_1}\\
\ x \ra \Lambda(x,1)
\end{cases}
$
is a homotopy equivalence.  
Definition: $\sH_q(f)[\tau]$ is the inverse of the induced isomorphism $H_q(X_{y_0}) \ra H_q(X_{y_1})$ (it is independent of the choices).\\
\endgroup%%------------------------------------<<

\textbf{\small LEMMA}  \ 
Suppose that \mX is path connected.  Given a locally constant coefficient system $\sG$, fix $x_0 \in X$, put 
$G_0 = \sG x_0$, and let $H_0$ be the subgroup of $G_0$ generated by the 
$g - g\cdot \alpha$ $(g \in G_0, \ \alpha \in \pi_1(X,x_0))$ $-$then $H_0(X;\sG) \approx G_0 /  H_0$.\\

%%----------------------------------------------------------------------------------------------40
\label{14.102}
Let $f:X \ra Y$ and $f^\prime:X^\prime \ra Y^\prime$ be a pair of continuous functions.  Call
$\mathcal{H}om(f^\prime,f)$ the simplicial set specified by taking for 
$\mathcal{H}om(f^\prime,f)_n$ the set of all
$
\begin{cases}
\ u \in C(\dpn \times X^\prime,X)\\[-.15cm]
\ v \in C(\dpn \times Y^\prime,Y)
\end{cases}
$
such that the diagram
\begin{tikzcd}[sep=large]
{\dpn \times X^\prime} \ar{d}[swap]{\id \times f^\prime} \ar{r}{u} &{X} \ar{d}{f}\\
{\dpn \times Y^\prime} \ar{r}[swap]{v} &{Y}
\end{tikzcd}
commutes and define
$
\begin{cases}
\ d_i\\[-.15cm]
\ s_i
\end{cases}
$
in the obvious way.

Now specialize, putting $Y^\prime = \dpz$, so $f^\prime:X^\prime \ra \dpz$ is the constant map, and write 
$\mathcal{H}om(X^\prime,f)$ in place of $\mathcal{H}om(f^\prime,f)$.  In succession, let 
$X^\prime = \dpz, \dpo, \ldots$ to obtain a sequence of simplicial sets and simplicial maps:
\[
\begin{tikzcd}[ sep=small]
{\mathcal{H}om(\dpz,f)} 
&{\mathcal{H}om(\dpo,f)} \ar[shift left]{l} \ar[shift right]{l} 
&{\mathcal{H}om(\Delta^2,f) \ \cdots}  \ar[l, shift left=0] \ar[l, shift left=2] \ar[l, shift right=2]
%&{\mathcal{H}om(\Delta^2,f) \ \cdots}  \ar[l, shift left=1] \ar[l, shift left=3] \ar[l, shift right=1] \ar[l, shift right=3]
\end{tikzcd}
\]
Here, the arrows come from the face operators 
\begin{tikzcd}[ sep=small]
{[0]} \ar[shift left]{r} \ar[shift right]{r} &{[1],\ [1]}
\ar{r} \ar[shift left=2]{r} \ar[shift right=2]{r} &{[2], \ldots \ .}
\end{tikzcd}
This data generates a double chain complex 
$K_{\bullet \bullet} = \{K_{n,m}:n \geq 0, m \geq 0\}$ of abelian groups if we write 
$K_{n,m} = F_{\ab}(\mathcal{H}om(\dpn,f)_m)$ and define 
$
\begin{cases}
\ \partial_I: K_{n,m} \ra K_{n-1,m}\\
\ \partial_{II}: K_{n,m} \ra K_{n,m-1}
\end{cases}
$
as follows.\\
\indent\indent $(\partial_I)$  The arrows %%%%%%%%%%%%%%%%%%%
\begin{tikzcd}[ sep=small] 
{\mathcal{H}om(\dpn,f)_m}  
    \arrow[r, draw=none, "\raisebox{+1.5ex}{\vdots}" description]
    \arrow[r,shift right=3]
    \arrow[r, shift left=3]
&{\mathcal{H}om(\Delta^{n-1},f)_m}
\end{tikzcd}
lead to arrows\\
$
\begin{tikzcd}[ sep=small] 
{K_{n,m}}  
    \arrow[r, draw=none, "\raisebox{+1.5ex}{\vdots}" description]
    \arrow[r,shift right=3]
    \arrow[r, shift left=3]
&{K_{n-1,m}}
\end{tikzcd}
.\ 
$
Take for $\partial_I$ their alternating sum multiplied by $(-1)^m$.\\
\indent\indent $(\partial_{II})$ The arrows %%%%%%%%%%%%%%%%%%%
\begin{tikzcd}[ sep=small] 
{\mathcal{H}om(\dpn,f)_m}  
    \arrow[r, draw=none, "\raisebox{+1.5ex}{\vdots}" description]
    \arrow[r,shift right=3]
    \arrow[r, shift left=3]
&{\mathcal{H}om(\Delta^{n},f)_{m-1}}
\end{tikzcd}
lead to arrows \\
$
\begin{tikzcd}[ sep=small] 
{K_{n,m}}  
    \arrow[r, draw=none, "\raisebox{+1.5ex}{\vdots}" description]
    \arrow[r,shift right=3]
    \arrow[r, shift left=3]
&{K_{n,m-1}}
\end{tikzcd}
.  \ 
$
Take for $\partial_{II}$ their alternating sum.

One can check that 
$\partial_I \circx \partial_I = 0 = \partial_{II} \circx \partial_{II}$ and 
$\partial_I \circx \partial_{II} + \partial_{II} \circx \partial_I = 0$.  
Form the total chain complex 
$K_{\bullet} = \{K_p\}$: $K_p = \bigoplus\limits_{n+m = p} K_{n,m}$, 
where
$\partial = \partial_I + \partial_{II}$ $-$then there are first quadrant spectral sequences
\[
\begin{cases}
\ _IE_{p,q}^2 = _IH_p(_{II}H_q(K_{\bullet\bullet})) \implies H_{p+q}(K_{\bullet})\\
\ _{II}E_{p,q}^2 = _{II}H_p(_{I}H_q(K_{\bullet\bullet})) \implies H_{p+q}(K_{\bullet})
\end{cases}
.
\]

\textbf{\small LEMMA}  \ 
%\[
$
_IE_{p,q}^2 \approx 
\begin{cases}
\ H_q(X) \hspace{0.575cm} (p = 0)\\[-.15cm]
\ 0 \hspace{1.5cm}  (p > 0)
\end{cases}
$
\hspace{-.25cm}.

[From the definitions, $\sin X = \mathcal{H}om(\dpz,f)$.  
On the other hand, each projection 
$\dpn \ra \dpz$ 
is a homotopy equivalence and induces an arrow 
$\sin X \ra \mathcal{H}om(\dpn,f)$.  
Since there are $n + 1$ commutative diagrams 
\begin{tikzcd}%[sep=large] 
{\sin X} \ar{d} \ar{r}{\id} &{\sin X} \ar{d}\\
{\mathcal{H}om(\dpn,f)} \ar{r} &{\mathcal{H}om(\Delta^{n-1},f)}
\end{tikzcd}
, passing to homology per $\partial_{II}$ gives  
\[
\begin{tikzcd}%[sep=large] 
{\underset{(p = 0)}{H_q(X)}} 
&{\underset{(p = 1)}{H_q(X)}} \ar{l}[swap]{0}
&{\underset{(p = 2)}{H_q(X)}} \ar{l}[swap]{\id}
&{\cdots} \ar{l}[swap]{0}
\end{tikzcd}
.]
\]

%%----------------------------------------------------------------------------------------------41
Thus the first spectral sequence $_IE$ collapses and $H_*(K_{\bullet}) \approx H_*(X)$.  
To explicate the second spectral sequence $_{II}E$, given $\tau \in \sin_n Y$, 
let $X_\tau$ be the fiber over $\tau$ of the induced map $\sin_n X \ra \sin_n Y$, i.e., $X_\tau = \{\sigma:f \circx \sigma = \tau\}$.  
View $X_\tau$ as a subspace of $\sin_n X$ (compact open topology).  
Put $\sH_q(f)\tau = H_q(X_\tau)$ and $\forall \ \alpha$, let $\sH_q(f)\Delta^\alpha$ be the homomorphism on homology defined by the arrow 
$X_\tau \ra X_{\tau \circx \Delta^\alpha}$ $-$then $\sH_q(f)$ is in 
$[(\sin Y)^\OP,\bAB]$ or still, is a coefficient system on \mY.

[Note: \ $\forall \ y \in Y$, $\sH_q(f)_{\tau_y} = H_q(X_y)$, where $\tau_y \in \sin_0 Y$ with $\tau_y(\dpz) = y$.]\\

\textbf{\small LEMMA}  \ 
$_{II}E_{p,q}^2 \approx  H_p(Y;\sH_q(f))$.

[$_IH_q(K_{\bullet\bullet})$ can be identified with the chain complex on which the homology of $\sH_q(f)$ is computed.]\\

\begin{proposition} \  %26
Suppose that $f:X \ra Y$ is a Hurewicz fibration $-$then $\sH_q(f)$ is locally constant.
\end{proposition}

[Fix $\alpha \in \Mor([m],[n])$ $-$then $\alpha$ determines arrows 
$
\begin{cases}
\ C(\dpn,X) \ra C(\dpm,X)\\[-.15cm]
\ C(\dpn,Y) \ra C(\dpm,Y)
\end{cases}
$
and there is a commutative diagram 
\begin{tikzcd}%[sep=large] 
{C(\dpn,X)} \ar{d} \ar{r}{f_*} &{C(\dpn,Y)} \ar{d}\\
{C(\dpm,X)} \ar{r}[swap]{f_*} &{C(\dpm,Y)}
\end{tikzcd}
.  
According to Proposition 5, the horizontal arrows are Hurewicz fibrations.  
But the vertical arrows are homotopy equivalences, thus $\forall \ \tau \in C(\dpn,Y)$ the induced map 
$X_\tau \ra X_{\tau \circx \Delta^\alpha}$ is a homotopy equivalence 
(cf. p. \pageref{4.44}), so 
$\sH_q(f)\Delta^\alpha:H_q(X_\tau) \ra H_q(X_{\tau \circx \Delta^\alpha})$ is an isomorphism.]

[Note: \ Retaining the assumption that $f:X \ra Y$ is a Hurewicz fibration, one may apply the procedure figuring in the proof of Proposition 25 to the locally constant coefficient system $\sH_q(f)$.  The result is the cofunctor 
$\sH_q(f):\Pi Y \ra \bAB$ defined in the example on p. \pageref{4.45}.]\\

Proposition 26 is also true if $f:X \ra Y$ is either a Dold fibration or a Serre fibration.\\

\label{9.97} %dmc mnft
\label{9.98} %dmc mnft
\begingroup%%----------------------------------->>
\fontsize{9pt}{11pt}\selectfont
Consider first the case when $f$ is Dold $-$then Proposition 5 still holds and the validity of the relevant homotopy theory has already been mentioned (cf. p. \pageref{4.46}).  
As for the case when $f$ is Serre, note that the arrow $C(\dpn,X) \ra C(\dpn,Y)$ is again Serre (as can be seen from the proof of Proposition 5).  
Therefore, thanks to the Whitehead theorem, the lemma below suffices to complete the argument.\\
\endgroup%%------------------------------------<<

\begingroup%%----------------------------------->>
\fontsize{9pt}{11pt}\selectfont
\textbf{\small LEMMA}  \ 
Suppose given a \cd 
\begin{tikzcd}[sep=large] 
{X} \ar{d}[swap]{\phi} \ar{r}{p} &{B} \ar{d}{\psi}\\
{Y} \ar{r}[swap]{q} &{A}
\end{tikzcd}
in which 
$
\begin{cases}
\ p\\
\ q
\end{cases}
$
are Serre fibrations and 
$
\begin{cases}
\ \phi\\
\ \psi
\end{cases}
$
are weak homotopy equivalences $-$then $\forall \ b \in B$, the induced map $X_b \ra Y_{\psi(b)}$ is a weak homotopy
%%----------------------------------------------------------------------------------------------42
equivalence.
\\ \indent
[If $X_b$ is empty, then so is $Y_{\psi(b)}$ and the assertion is trivial.  Otherwise, let $a = \psi(b)$ and apply the five lemma to the commutative diagram
\[
\begin{tikzcd}[sep=large]
{\cdots} \ar{r}
&{\pi_{q+1}(B)} \ar{d} \ar{r}
&{\pi_{q}(X_b)} \ar{d} \ar{r}
&{\pi_{q}(X)} \ar{d} \ar{r}
&{\pi_{q}(B)} \ar{d} \ar{r}
&{\cdots}\\
{\cdots} \ar{r}
&{\pi_{q+1}(A)} \ar{r}
&{\pi_{q}(Y_a)} \ar{r}
&{\pi_{q}(Y)} \ar{r}
&{\pi_{q}(A)} \ar{r}
&{\cdots}
\end{tikzcd}
,
\]
with the usual caveat at the $\pi_0$ and $\pi_1$ level.]\\
\endgroup%%------------------------------------<<

The coefficient system $\sH_q(f)$ is defined in terms of the integral singular homology of the fibers.  Embelish the notation and denote it by $\sH_q(f;\Z)$.  One may then replace $\Z$ by any abelian group \mG: $\sH_q(f;G)$, a coefficient system which is locally constant if $f:X \ra Y$ is either a Dold fibration or a Serre fibration.\\

\index{Spectral Sequence:\\ Fibration Spectra Sequence}
\index{Fibration Spectra Sequence}
\textbf{\small FIBRATION SPECTRAL SEQUENCE} \ 
Let $f:X \ra Y$ be either a Dold fibration or a Serre fibration $-$then for any abelian group \mG, there is a first quadrant spectral sequence 
$E = \{E_{p,q}^r,d^r\}$ such that 
$E_{p,q}^2 \approx H_p(Y;\sH_q(f;G))$ $\implies$ $H_{p+q}(X;G)$ and $\forall \ n$, $H_n(X;G)$ admits an increasing filtration
\[
0 = H_{-1,n+1} \subset H_{0,n} \subset \cdots \subset H_{n-1,1} \subset H_{n,0} = H_n(X;G)
\]
by subgroups $H_{i,n-i}$, where 
$E_{p,q}^\infty \approx H_{p,q} / H_{p-1,q+1}$.

[Note: \ The fibration spectral sequence is natural, i.e., if the diagram
\begin{tikzcd}[sep=large]
{X} \ar{d} \ar{r}{f} &{Y} \ar{d}\\
{X^\prime} \ar{r}[swap]{f^\prime} &{Y^\prime}
\end{tikzcd}
commutes, then there is a morphism $\mu:E \ra E^\prime$ of spectral sequences such that 
$\mu_{p,q}^2$ coincides with the homomorphism 
$H_p(Y;\sH_q(f;G)) \ra H_p(Y^\prime;\sH_q(f^\prime;G))$ induced by the arrow 
$\sH_q(f;G) \ra \sH_q(f^\prime;G)$.]\\

\begingroup%%----------------------------------->>
\fontsize{9pt}{11pt}\selectfont
\index{Wang Homology Sequence}
\textbf{\small WANG HOMOLOGY SEQUENCE} \ 
Take $Y = \bS^{n+1}$ $(n \geq 1)$ and let $f:X \ra Y$ be a Hurewicz fibration with path connected fibers $X_y$ $-$then 
there is an exact sequence
\[
\cdots \ra H_q(X) \ra H_{q - n - 1}(X_y) \ra H_{q-1}(X_y) \ra H_{q-1}(X) \ra \cdots .
\]
\endgroup%%------------------------------------<<

\begingroup%%----------------------------------->>
\fontsize{9pt}{11pt}\selectfont
\textbf{\small EXAMPLE}  \ 
Suppose that $n \geq 1$ $-$then $H_{kn}(\Omega \bS^{n+1}) \approx \Z$ $(k = 0, 1, \ldots)$, while 
$H_q(\Omega \bS^{n+1}) = 0$ otherwise.  Moreover, the Pontryagin ring 
$H_*(\Omega \bS^{n+1})$ is isomorphic to $\Z[t]$, where t generates $H_n(\Omega \bS^{n+1})$.\\
\endgroup%%------------------------------------<<

%%----------------------------------------------------------------------------------------------43
\begingroup%%----------------------------------->>
\fontsize{9pt}{11pt}\selectfont
As formulated, the fibration spectral sequence applies to singular homology.  There is also a companion result for singular cohomology (with additional multiplicative structure when the coefficient group \mG is a commutative ring).\\
\endgroup%%------------------------------------<<

\index{Wang Cohomology Sequence}
\begingroup%%----------------------------------->>
\fontsize{9pt}{11pt}\selectfont
\textbf{\small WANG COHOMOLOGY SEQUENCE} \ 
Take $Y = \bS^{n+1}$ $(n \geq 1)$ and let $f:X \ra Y$ be a Hurewicz fibration with path connected fibers $X_y$ $-$then 
there is an exact sequence
\[
\cdots \ra H^q(X) \ra H^{q}(X_y) \overset{\theta}{\ra} H^{q-n}(X_y) \ra H^{q+1}(X) \ra \cdots .
\]
\\ \indent
[Note: \ In the graded ring $H^*(X_y)$, 
$\theta(\alpha\cdot\beta) = \theta(\alpha) \cdot \beta + (-1)^{n\abs{\alpha}} \alpha \cdot \theta(\beta)$.]\\
\endgroup%%------------------------------------<<

\begingroup%%----------------------------------->>
\fontsize{9pt}{11pt}\selectfont
\textbf{\small EXAMPLE}  \ 
Suppose that $n \geq 1$ $-$then 
$\theta:H^{kn}(\Omega \bS^{n+1}) \ra H^{(k-1)n}(\Omega \bS^{n+1})$ $(k \geq 1)$ 
is an isomorphism and $H^0(\Omega \bS^{n+1})$ is the infinite cyclic group generated by 1.  
Put $\alpha_0 = 1$ and define $\alpha_k$ $(k \geq 1)$ inductively through the relation 
$\theta(\alpha_k) = \alpha_{k-1}$.  
Case 1:  $n$ even.  One has $k!\alpha_k = \alpha_1^k$, therefore $H^*(\Omega \bS^{n+1})$ 
is the divided polynomial algebra generated by $\alpha_1, \alpha_2, \ldots$.
Case 2: $n$ odd.  One has 
$\alpha_1^2 = 0$, 
$\alpha_1\alpha_{2k} = \alpha_{2k+1}$, 
$\alpha_1\alpha_{2k+1}  = 0$, and 
$\alpha_2^k = k!\alpha_{2k}$, thus $\alpha_1$ generates an exterior algebra isomorphic to 
$H^*(\bS^n)$ and 
$\alpha_2, \alpha_4, \ldots$ generate a divided polynomial algebra isomorphic to 
$H^{*}(\Omega \bS^{2n+1})$, so
$H^{*}(\Omega \bS^{n+1}) \approx$ 
$H^*(\bS^n) \otimes H^*(\Omega \bS^{2n+1})$.\\
\endgroup%%------------------------------------<<

In what follows, we shall assume that \mX is nonempty and \mY is path connected.

[Note: \ If $f$ is Dold, then the $X_y$ have the same homotopy type (cf. p. \pageref{4.47}), 
while if $f$ is Serre, then $X_y$ have the same weak homotopy type (cf. Proposition 31).]\\
\indent\indent (ED$_\tH$) \ Let $e_H:E_{p,0}^\infty \ra E_{p,0}^2$ be the edge homomorphism on the horizontal axis.  The arrow of augmentation $H_0(X_y,G) \ra G$ is independent of $y$, so there is a homomophism 
$H_p(Y;\sH_0(f;G)) \ra H_p(Y;G)$.  The composite 
$H_p(X;G) \ra H_{p,0} /  H_{p-1,1} \approx$ 
$E_{p,0}^\infty \overset{e_H}{\lra}$ 
$E_{p,0}^2 \approx$
$H_p(Y;\sH_0(f;G)) \ra$ 
$H_p(Y;G)$ is the homomorphism on homology induced by $f:X \ra Y$.\\
\indent\indent (ED$_\tV$) \ Let $e_V:E_{0,q}^2 \ra E_{0,q}^\infty$ be the edge homomorphism on the vertical axis.  
Fix $y \in Y$ $-$then there is an arrow
$H_q(X_y;G) \ra H_0(Y;\sH_q(f;G))$.  The composite 
$H_q(X_y;G) \ra$ 
$H_0(Y;\sH_q(f;G)) \approx$ 
$E_{0,q}^2 \overset{e_V}{\lra}$
$E_{0,q}^\infty \lra$ 
$H_q(X;G)$ is the homomorphism on homology induced by the inclusion $X_y \ra X$.

Keeping to the preceding hypotheses, $f:X \ra Y$ is said to be 
\un{$G$-orientable}
\index{G-orientable ($f:X \ra Y$)} 
provided that the $X_y$ are path connected and $\forall \ q$, $\sH_q(f;G)$ is constant, so $\forall \ y$ the right action 
$H_q(X_y;G) \times \pi_1(Y,y) \ra H_q(X_y;G)$ is trivial.

[Note: \ If $f:X \ra Y$ is $G$-orientable, then by the universal coefficient theorem, 
$E_{p,q}^2 \approx$ 
$H_p(Y;\sH_q(X_y;G)) \approx$ 
$H_p(Y) \otimes H_q(X_y;G) \otimes \Tor(H_{p-1}(Y),H_q(X_y;G))$.]\\

%%----------------------------------------------------------------------------------------------44

\label{5.0g}
\label{5.10}
\begingroup%%----------------------------------->>
\fontsize{9pt}{11pt}\selectfont
\textbf{\small EXAMPLE}  \ 
Let $f:X \ra Y$ be $G$-orientable.  Assume: $H_i(X_{y_0};G) = 0$ $(0 < i \leq n)$ and 
$H_j(Y;\Z) = 0$ $(0 < j \leq m)$ $-$then there is an exact sequence
\\[-.5cm]
\begin{align*}
H_{n+m+1}(X_{y_0};G) \ra H_{n+m+1}(X;G) \ra H_{n+m+1}(Y;G) \ra H_{n+m}(X_{y_0};G) \ra \cdots\\
\indent\indent \ra H_2(Y;G) \ra H_1(X_{y_0};G) \ra H_1(X;G) \ra H_1(Y;G).
\end{align*}
\\[-.5cm]
\indent
[For $2 \leq r < n + m + 2$, combine the exact sequence
\[
0 \ra E_{r,0}^\infty \ra E_{r,0}^r \overset{d^r}{\ra} E_{0,r-1}^r \ra E_{0,r-1}^\infty \ra 0
\]
with the exact sequence
\[
0 \ra E_{0,r}^\infty \ra H_r(X;G) \ra E_{r,0}^\infty \ra 0 
\]
observing that 
$H_r(Y;G) \approx$ 
$E_{r,0}^2 \approx$ 
$E_{r,0}^r$ 
and 
$H_{r-1}(X_{y_0};G) \approx$ 
$E_{0,r-1}^2 \approx$ 
$E_{0,r-1}^r $, the arrow $H_r(Y;G) \ra H_{r-1}(X_{y_0};G)$ being the transgression.]
\\ \indent
[Note: \ The above assumptions are less stringent than those imposed earlier in the case $G = \Z$ (cf. p. \pageref{4.48}).]\\
\endgroup%%------------------------------------<<

\label{5.0aw}
\label{8.4}
\label{8.41} %dmc mnft
\label{9.79}
\begingroup%%----------------------------------->>
\fontsize{9pt}{11pt}\selectfont
\textbf{\small EXAMPLE}  \ 
Let $f:X \ra Y$ be $\Lambda$-orientable, where $\Lambda$ is a principal ideal domain $-$then the arrow 
$H_*(X;\Lambda) \ra H_*(Y;\Lambda)$ is an isomorphism iff $\forall \ q > 0$, the $H_q(X_{y_0};\Lambda) = 0$ 
and the arrow 
$H_*(X_{y_0};\Lambda) \ra$ 
$H_*(X;\Lambda)$ is an isomorphism iff $\forall \ q > 0$, $H_q(Y;\Lambda) = 0$.
\\ \indent
[Note: \ The formulation is necessarily asymmetric (take \mY simply connected and consider $\Theta Y \ra Y$).]\\
\endgroup%%------------------------------------<<

\label{5.0ap}
\label{5.8d}
\begingroup%%----------------------------------->>
\fontsize{9pt}{11pt}\selectfont
\textbf{\small FACT}  \ 
Suppose that $f:X \ra Y$ is $\Z$-orientable $-$then any two of the following conditions imply the third: 
(1) $\forall \ p$, $H_p(Y)$ is finitely generated;
(2) $\forall \ q$, $H_q(X_{y_0})$ is finitely generated;
(3) $\forall \ n$, $H_n(X)$ is finitely generated.\\
\endgroup%%------------------------------------<<

\begingroup%%----------------------------------->>
\fontsize{9pt}{11pt}\selectfont
\textbf{\small FACT}  \ 
Suppose that $f:X \ra Y$ is $\Z$-orientable $-$then any two of the following conditions imply the third: 
(1) $\forall \ p > 0$, $H_p(Y)$ is finite;
(2) $\forall \ q > 0$, $H_q(X_{y_0})$ is finite;
(3) $\forall \ n > 0$, $H_n(X)$ is finite.\\
\endgroup%%------------------------------------<<

Given pointed spaces 
$
\begin{cases}
\ X\\[-.15cm]
\ Y
\end{cases}
, \ 
$
the mapping fiber sequence associated with the inclusion $f:X \vee Y \ra X \times Y$ reads: 
$\cdots \ra \Omega(X \vee Y) \ra \Omega(X \times Y) \ra X \flat Y \ra X \vee Y \ra X \times Y$.

[Note: \ The homology of $\Omega(X \vee Y)$ can be calculated in terms of the homology of $\Omega X$ and 
$\Omega Y$ 
(Aguade-Castellet\footnote[2]{\textit{Collect. Math.} \textbf{29} (1978), 3-6; 
see also Dula-Katz, \textit{Pacific J. Math.} \textbf{86} (1980), 451-461.}).]\\

\textbf{\small LEMMA}  \ 
The arrow $F:\Omega(X \times Y) \ra X \flat Y$ is nullhomotopic.

%%----------------------------------------------------------------------------------------------45
[Put
$
\begin{cases}
\ \ov{\Omega} X = \{\sigma: \sigma([1/2,1]) = x_0\}\\
\ \un{\Omega} Y = \{\tau: \tau([0,1/2]) = y_0\}
\end{cases}
$
$-$then the inclusions 
$
\begin{cases}
\ \ov{\Omega} X \ra \Omega X\\
\ \un{\Omega} Y \ra \Omega Y
\end{cases}
$ 
are pointed homotopy equivalences, hence the same holds for their product:
$\ov{\Omega} X \times \un{\Omega} Y \ra \Omega X \times \Omega Y =$ 
$\Omega (X \times Y)$.  Use two parameter reversals to see that the composite 
$\ov{\Omega} X \times \un{\Omega} Y \ra$ 
$\Theta (X \vee Y) \ra$ 
$X \flat Y$ is equal to the composite 
$\ov{\Omega} X \times \un{\Omega} Y \ra$ 
$\Omega (X \times Y) \overset{F}{\lra}$ $X \flat Y$, from which $F \simeq 0$.]\\

\index{Ganea-Nomura Formula}
%\textbf{\small GANEA-NOMURA FORMULA} \quadx
\textbf{\small GANEA-NOMURA FORMULA} \ 
Suppose that \mX and \mY are nondegenerate $-$then in $\bHTOP_*$, 
$\Omega(X \vee Y) \approx \Omega X \times \Omega Y \times \Omega\Sigma(\Omega X \# \Omega Y)$.

[The mapping fiber of 0: 
$\Omega (X \times Y) \ra X \flat Y$ is 
$\Omega (X \times Y)  \times \Omega(X \flat Y)$ and by the lemma, 
$E_F \approx \Omega (X \times Y)  \times \Omega(X \flat Y)$.  
Employing the notation of p. \pageref{4.49}, there is a commutative triangle 
\begin{tikzcd}%[sep=large]
{E_\pi} \ar{r}{\pi^\prime} &{X \flat Y}\\
{\Omega (X \times Y)} \ar{u} \ar{ur}[swap]{F}
\end{tikzcd}
.  \ The vertical arrow is a pointed homotopy equivalence, thus 
$E_{\pi^\prime} \approx E_F$ or still, 
$\Omega(X \vee Y) \approx$
$\Omega (X \times Y)  \times \Omega(X \flat Y) \approx$ 
$\Omega X \times \Omega Y \times \Omega\Sigma(\Omega X \# \Omega Y)$ (cf. p. \pageref{4.50}).]\\
\vspace{0.15cm}

Given pointed spaces 
$
\begin{cases}
\ X\\[-.15cm]
\ Y
\end{cases}
, \ 
$
the mapping fiber sequence associated with the projection $f:X \vee Y \ra Y$ reads:
$\cdots \ra \Omega(X \vee Y) \ra \Omega Y \ra E_f \ra X \vee Y \ra Y$.\\

\textbf{\small LEMMA}  \ 
The arrow $F:\Omega Y \ra E_f$ is nullhomotopic.

[Define $g:Y \ra X \vee Y$ by $g(y) = (x_0,y)$, so $f \circx g = \id_Y$.  
Let \mZ be any pointed space $-$then in view of the replication theorem, there is an exact sequence 
$[Z,\Omega(X \vee Y)] \ra [Z,\Omega Y] \ra [Z,E_f]$.  
Since $\Omega f$ has a right inverse, the arrow 
$[Z,\Omega(X \vee Y)] \ra [Z,\Omega Y]$ is surjective.  This means that the arrow 
$[Z,\Omega Y] \ra [Z,E_f]$ is the zero map, therefore $F$ is nullhomotopic.]\\

\index{Gray-Nomura Formula}
\textbf{\small GRAY-NOMURA FORMULA} \ 
Suppose that \mX and \mY are nondegenerate $-$then in $\bHTOP_*$, 
$\Omega(X \vee Y) \approx \Omega Y \times \Omega(X \times \Omega Y / \{x_0\} \times \Omega Y)$.

[Argue as in the proof of the Ganea-Nomura formula ($E_f$ is determined on p. \pageref{4.51}).]\\

\begin{proposition} \  
Let \mX, \mY be pointed spaces $-$then $\Sigma X \times Y/\{x_0\} \times Y$ has the same pointed homotopy type as 
$\Sigma X \vee (\Sigma X \# Y)$.
\end{proposition}

[$\Sigma X \times Y/\{x_0\} \times Y \approx$ 
$\Sigma X \# Y_+ \approx$
$\Sigma X \# (\bS^0 \vee Y) \approx$ 
$X \# \Sigma(\bS^0 \vee Y) \approx$ 
$X \# (\bS^1 \vee \Sigma Y) \approx$ 
$(X \# \bS^1) \vee (X \# \Sigma Y) \approx$ 
$\Sigma X \vee (\Sigma X \#Y)$.]

[Note: \ Recall that in $\bHTOP_*$, 
$\Sigma (X \# Y) \approx$ 
$\Sigma X \# Y \approx$ 
$X \# \Sigma Y$ for arbitrary pointed \mX and \mY (cf. p. \pageref{4.52}).]\\

%%----------------------------------------------------------------------------------------------46
So, if $X$ is the pointed suspension of a nondegenerate space, then the Gray-Nomura formula can be simplified: 
$\Omega(X \vee Y) \approx$ 
$\Omega Y \times \Omega(X \vee (X \# \Omega Y))$.  Consequently, for all nondegenerate \mX and \mY, 
\[
\Omega\Sigma (X \vee Y) \approx
\begin{cases}
\ \Omega\Sigma X \times \Omega\Sigma (Y \vee (Y \# \Omega\Sigma X))\\
\ \Omega\Sigma Y \times \Omega\Sigma (X \vee (X \# \Omega\Sigma Y))
\end{cases}
.
\]

\begingroup%%----------------------------------->>
\fontsize{9pt}{11pt}\selectfont
Suppose that 
$
\begin{cases}
\ X\\
\ Y
\end{cases}
,Z
$
are wellpointed with 
$
\begin{cases}
\ \{x_0\} \subset X\\
\ \{y_0\} \subset Y
\end{cases}
, \{z_0\} \subset Z
$
closed.  
Let $f:X \ra Y$ be a pointed continuous function, $C_f$ its pointed mapping cone.  Let $p:Z \ra C_f$ be a pointed continuous function, $Z_0$ its fiber over the base point.  
Assume: $p$ is a Hurewicz fibration $-$then $p$ is a pointed Hurewicz fibration.  Form the pullback square 
\begin{tikzcd}[sep=large]
{P} \ar{d} \ar{r} &{Z} \ar{d}{p}\\
{Y} \ar{r}[swap]{j} &{C_f}
\end{tikzcd}
.  Since $j \circx f \simeq 0$, there is a commutative triangle 
\begin{tikzcd}[sep=large]
&{Z} \ar{d}{p}\\
{X} \ar{ru}{k} \ar{r}[swap]{j \circx f} &{C_f}
\end{tikzcd}
and an induced map $e:X \ra P$.\\
\vspace{0.25 cm}\endgroup%%------------------------------------<<

\begingroup%%----------------------------------->>
\fontsize{9pt}{11pt}\selectfont
\textbf{\small FACT}  \ 
The pointed mapping cone of the arrow $C_e \ra Z$ has the pointed homotopy type of $X * Z_0$.\\
\endgroup%%------------------------------------<<

\begingroup%%----------------------------------->>
\fontsize{9pt}{11pt}\selectfont
\textbf{\small EXAMPLE}  \ 
Let \mX be wellpointed with $\{x_0\} \subset X$ closed.  The pointed mapping cone of $X \ra *$ is $\Sigma X$, the pointed suspension of \mX.  Consider the pullback square
\begin{tikzcd}[sep=large]
{\Omega\Sigma X} \ar{d} \ar{r} &{\Theta\Sigma X} \ar{d}{p_1}\\
{*} \ar{r} &{\Sigma X}
\end{tikzcd}
Here, $e:X \ra \Omega\Sigma X$ is the arrow of adjunction and the pointed mapping cone 
$C_e \ra \Theta\Sigma X$ has the same pointed homotopy type as $C_e \ra *$, thus in $\bHTOP_*$, 
$\Sigma C_e \approx X * \Omega\Sigma X$.\\
\endgroup%%------------------------------------<<

Given a pointed space \mX, the pointed mapping cone sequence associated with the arrow of adjunction 
$e:X \ra \Omega\Sigma X$ reads: 
$X \overset{e}{\ra} \Omega\Sigma X \ra C_e \ra \Sigma X \ra \Sigma\Omega\Sigma X \ra \cdots$.\\

\begin{proposition} \  
Let \mX be nondegenerate $-$then $\Sigma\Omega\Sigma X$ has the same pointed homotopy type as 
$\Sigma X \vee \Sigma(X \# \Omega \Sigma X)$.
\end{proposition}

[Because the evaluation map $r:\Sigma\Omega\Sigma X \ra \Sigma X$ exhibits $\Sigma X$ as a retract of 
$\Sigma\Omega\Sigma X$, the replication theorem of $\S 3$ implies that the arrow $F:C_e \ra \Sigma X$ 
is nullhomotopic, hence $C_F \approx \Sigma X \vee \Sigma C_e$.  Reverting to the notation of p. \pageref{4.53}, there is a commutative triangle 
\begin{tikzcd}%[sep=large]
{C_e} \ar{rd}[swap]{F} \ar{r}{j^\prime} &{C_j} \ar{d}\\
&{\Sigma X}
\end{tikzcd}
in which the vertical arrow is a pointed homotopy equivalence.  Accordingly, 
$C_{j^\prime} \approx C_F$ $\implies$ 
$\Sigma\Omega\Sigma X \approx \Sigma X \vee \Sigma C_e \approx$ 
$\Sigma X \vee \Sigma (X \# \Omega\Sigma X)$, the last step by the preceding example.]\\

%%----------------------------------------------------------------------------------------------47
Assume: \mX and \mY are nondegenerate.  
Put $X^{[0]} = \bS^0$, $X^{[n]} = X \# \cdots \# X$ ($n$ factors).  
Starting from the formula
$\Omega\Sigma (X \vee Y) \approx \Omega\Sigma X \times \Omega\Sigma(Y \vee (Y \# \Omega\Sigma X))$, 
successive application of Proposition 28 gives:
\[
\Omega\Sigma (X \vee Y) \approx \Omega\Sigma X \times \Omega\Sigma
\left(\bigvee\limits_0^N Y \# X^{[n]}  \vee (Y \# X^{[N]} \# \Omega\Sigma X)\right).
\]

\begingroup%%----------------------------------->>
\fontsize{9pt}{11pt}\selectfont
\textbf{\small FACT}  \ 
Let 
$
\begin{cases}
\ X\\
\ Y
\end{cases}
$
be nondegenerate and path connected $-$then $\forall \ q > 0$, 
$\pi_q(\Sigma X \vee \Sigma Y) \approx$ 
$\pi_q(\Sigma X) \oplus \pi_q \bigl(\ds\Sigma\left(\bigvee\limits_0^\infty Y \# X^{[n]}\right)\bigr)$.
\\ \indent
[By the above, $\pi_q(\Sigma X \vee \Sigma Y)$ is isomorphic to 
\[
\pi_q(\Sigma X) \oplus \pi_q (\ds\Sigma\left(\bigvee\limits_0^N Y \# X^{[n]} \vee (Y \# X^{[N]} \# \Omega\Sigma X)\right)).
\]
Since $\Sigma(Y \#X^{[N]} \# \Omega\Sigma X)$ is $(N + 2)$-connected (cf. p. \pageref{4.54}), it follows that 
$\forall \ q \leq N + 2$: 
$\pi_q(\Sigma X \vee \Sigma Y) \approx$ 
$\pi_q(\Sigma X) \oplus \pi_q(\Sigma\left(\bigvee\limits_0^N Y \# X^{[n]}\right))$.  
But 
$\Sigma \left(\bigvee\limits_{n > N} Y \# X^{[n]}\right)$ 
is also $(N + 2)$-connected.  Therefore, $\forall \ q > 0$ : 
$\pi_q(\Sigma X \vee \Sigma Y) \approx$ 
$\pi_q(\Sigma X) \oplus \pi_q \bigl(\Sigma\left(\bigvee\limits_0^\infty Y \# X^{[n]}\right)\bigr)$. ]\\
\endgroup%%------------------------------------<<

A continuous function $f:X \ra Y$ is said to be an 
\un{$n$-equivalence}
\index{n-equivalence} 
$(n \geq 1)$ provided that $f$ induces a one-to-one correspondence between the path components of 
$
\begin{cases}
\ X\\
\ Y
\end{cases}
$
and $\forall \ x_0 \in X$, $f_*:\pi_q(X,x_0) \ra \pi_q(Y,f(x_0))$ is bijective for $1 \leq q < n$ and surjective for $q = n$.  
Example: A pair $(X,A)$ is $n$-connected iff the inclusion $A \ra X$ is an $n$-equivalence.

[Note: \ $f$ is an $n$-equivalence iff the pair $(M_f,i(X))$ is $n$-connected.]\\

\begingroup%%----------------------------------->>
\fontsize{9pt}{11pt}\selectfont
\textbf{\small FACT}  \ 
Let $X \overset{p}{\lra} B \overset{q}{\lla} Y$ be a 2-sink.  Suppose that 
$
\begin{cases}
\ p \text{ is an $n$-equivalence}\\
\ q \text{ is an $m$-equivalence}
\end{cases}
$
$-$then the projection $X \bboxsub_B Y \ra B$ is an $(n + m + 1)$-equivalence.
\\ \indent
[There is an arrow $X \bboxsub_B Y \overset{\phi}{\lra} W_p *_B W_q$ that commutes with the projections and is a homotopy equivalence (cf. p. \pageref{4.55}), thus one can assume that 
$
\begin{cases}
\ p\\
\ q
\end{cases}
$
are Hurewicz fibrations and work instead with $X *_B Y$ (the connectivity of the join is given on p. \pageref{4.56}).]\\
\endgroup%%------------------------------------<<

A continuous function $f:X \ra Y$ is said to be a 
\un{weak homotopy equivalence}
\index{weak homotopy equivalence} 
if $f$ is an $n$-equivalence $\forall \ n \geq 1$.  
Example: Consider the coreflector $k:\bTOP \ra \bCG$ $-$then for every topological space $X$, the identity map 
$kX \ra X$ is a weak homotopy equivalence.  

[Note: \ When $X$ and $Y$ are path connected, $f$ is a \whe provided that at some $x_0 \in X$, 
$f_*:\pi_q(X,x_0) \ra \pi_q(Y,f(x_0))$ is bijective $\forall \ q \geq 1$.]\\

%%----------------------------------------------------------------------------------------------48
\label{18.22}
\label{9.104}
\label{14.26}
\label{14.78}
\label{14.111}
\label{14.136}
\label{14.150}
Example: Let 
\begin{tikzcd}%[sep=large]
{X} \ar{d} \ar{r}{f} &{Z} \ar{d} &{Y} \ar{l}[swap]{g} \ar{d}\\
{X^\prime}  \ar{r}[swap]{f^\prime} &{Z^\prime} &{Y^\prime} \ar{l}{g^\prime} 
\end{tikzcd}
be a \cd in which the vertical arrows are homotopy equivalences $-$then the arrow $W_{f,g} \ra W_{f^\prime,g^\prime}$ is a weak homotopy equivalence.

[Compare Mayer-Vietoris sequences (use an ad hoc argument to establish that 
$\pi_0(W_{f,g}) \approx \pi_0(W_{f^\prime,g^\prime})$).]

Example: Let 
$
\begin{cases}
\ X^0 \subset X^1 \subset \cdots\\[-.15cm]
\ Y^0 \subset Y^1 \subset \cdots
\end{cases}
$
be an expanding sequence of topological spaces.  Assume: $\forall \ n$, the inclusions 
$
\begin{cases}
\ X^n \ra X^{n+1}\\[-.15cm]
\ Y^n \ra Y^{n+1}
\end{cases}
$
are closed cofibrations.  Suppose given a sequence of continuous functions $\phi^n:X^n \ra Y^n$ such that $\forall \ n$, the diagram 
\begin{tikzcd}%[sep=large]
{X^n} \ar{d}[swap]{\phi^n} \ar{r} &{X^{n+1}} \ar{d}{\phi^{n+1}}\\
{Y^n}  \ar{r} &{Y^{n+1}} 
\end{tikzcd}
commutes $-$then $\phi^\infty:X^\infty \ra Y^\infty$ is a weak homotopy equivalence if this is the case of the $\phi^n$.

\label{14.23b}
[Consider the commutative diagram  
\begin{tikzcd}%[sep=large]
{\tel X^\infty} \ar{d}[swap]{\tel \phi} \ar{r} &{X^\infty} \ar{d}{\phi^\infty}\\
{\tel Y^\infty}  \ar{r} &{Y^\infty} 
\end{tikzcd}
(cf. p. \pageref{4.57}).  
Since the horizontal arrows are homotopy equivalences, it suffices to prove that $\tel \phi$ is a weak homotopy equivalence.  
To see this, recall that there are projections 
$
\begin{cases}
\ \tel X^\infty \ra [0,\infty[\\[-.15cm]
\ \tel Y^\infty \ra [0,\infty[
\end{cases}
\hspace{-.25cm}, \ 
$
thus a compact subset of 
$
\begin{cases}
\ \tel X^\infty\\[-.15cm]
\ \tel Y^\infty
\end{cases}
$
must lie in 
$
\begin{cases}
\ \telsub_n X^\infty\\[-.15cm]
\ \telsub_n Y^\infty
\end{cases}
$
$(\exists \ n \gg 0)$.  
But $\forall \ n$, the arrow $\telsub_n X^\infty \ra \telsub _n Y^\infty$ is a weak homotopy equivalence.]

[Note: \ Here is a variant.  
Let 
$
\begin{cases}
\ X^0 \subset X^1 \subset \cdots\\[-.15cm]
\ Y^0 \subset Y^1 \subset \cdots
\end{cases}
$
be expanding sequences of topological spaces.  
Assume: $\forall \ n$, 
$
\begin{cases}
\ X^n\\[-.15cm]
\ Y^n
\end{cases}
$
is $\tT_1$.   
Suppose given a sequence of continuous functions $\phi^n:X^n \ra Y^n$ such that $\forall \ n$, the diagram
\begin{tikzcd}%[sep=large]
{X^n} \ar{d}[swap]{\phi^n} \ar{r} &{X^{n+1}} \ar{d}{\phi^{n+1}}\\
{Y^n}  \ar{r} &{Y^{n+1}} 
\end{tikzcd}
commutes $-$then $\phi^\infty:X^\infty \ra Y^\infty$ is a weak homotopy equivalence if this is the case of the $\phi^n$.]\\ 
\vspace{0.15cm}

\begingroup%%----------------------------------->>
\fontsize{9pt}{11pt}\selectfont
\textbf{\small EXAMPLE}  \ 
Given pointed spaces $X$ and $Y$, let $X \bowtie Y$ be the double mapping track of the 2-sink  
$X \ra X \vee Y \la Y$.  The projection $X \bowtie Y \ra X \times Y$ is a pointed Hurewicz fibration.  Its fiber over $(x_0,y_0)$ is 
$\Omega(X \vee Y)$ and the composite 
$\Omega(X \flat Y) \ra \Omega(X \vee Y) \ra X \bowtie Y$ defines a weak homotopy equivalence 
$\Omega(X \flat Y) \ra X \bowtie Y$.\\
\endgroup%%------------------------------------<<

\begingroup%%----------------------------------->>
\fontsize{9pt}{11pt}\selectfont
Assume $X$ and $Y$ are nondegenerate $-$then the argument used to establish that
%%----------------------------------------------------------------------------------------------49
\[
\Omega\Sigma (X \vee Y) \approx \Omega\Sigma X \times 
\Omega\Sigma \left( \bigvee\limits_0^N Y \#X^{[n]} \vee \left(Y \#X^{[N]} \# \Omega\Sigma X \right) \right)
\]
does not explicity produce a pointed homotopy equivalence between either side but such precision is possible.  
Let 
$
\begin{cases}
\ \iota_{\Sigma X}\\
\ \iota_{\Sigma Y}
\end{cases}
$
be the inclusions 
$
\begin{cases}
\ \Sigma X \ra \Sigma X \vee \Sigma Y\\
\ \Sigma Y \ra \Sigma X \vee \Sigma Y
\end{cases}
\hspace{-.25cm}. \ 
$
With $w_0 = \iota_{\Sigma Y}$, inductively define 
$w_1 = [w_0,\iota_{\Sigma Y}], \ldots , w_n = [w_{n-1},\iota_{\Sigma Y}]$, the bracket being the Whitehead product, 
so 
$w_1: \Sigma (Y \#X) \ra \Sigma X \vee \Sigma Y,$ $\ldots$,  
$w_n: \Sigma (Y \# X^{[n]}) \ra \Sigma X \vee \Sigma Y$.  Write
$\Omega(\iota_{\Sigma Y}) + \Omega\left(\bigvee\limits_0^N w_n \vee [w_N,\iota_{\Sigma X} \circx r] \right)$ 
for the composite 
\[
\Omega\Sigma X \times 
\Omega\Sigma \left( \bigvee\limits_0^N Y \#X^{[n]} \vee \left(Y \#X^{[N]} \# \Omega\Sigma X \right) \right)
\lra 
\Omega\Sigma (X \vee Y) \times \Omega\Sigma (X \vee Y)  \overset{+}{\lra} \Omega\Sigma (X \vee Y).
\]
Then Spencer\footnote[2]{\textit{J. London Math. Soc.} \textbf{4} (1971), 291-303.}
has shown that 
$\Omega(\iota_{\Sigma X}) + \Omega\left( \ds\bigvee\limits_0^N w_n \vee [w_N,\iota_{\Sigma X} \circx r]\right)$
is a pointed homotopy equivalence.\\
\endgroup%%------------------------------------<<

\begingroup%%----------------------------------->>
\fontsize{9pt}{11pt}\selectfont
\textbf{\small EXAMPLE}  \ 
Let 
$
\begin{cases}
\ X\\
\ Y
\end{cases}
$
be nondegenerate and path connected $-$then the map
\[
\Omega(\iota_{\Sigma X}) + \Omega\left( \bigvee\limits_0^\infty w_n\right): 
\Omega\Sigma X \times \Omega\Sigma \left(\bigvee\limits_0^\infty Y \#X^{[n]}\right) \lra \Omega\Sigma (X \vee Y)
\]
is a weak homotopy equivalence.\\
\endgroup%%------------------------------------<<

\begingroup%%----------------------------------->>
\fontsize{9pt}{11pt}\selectfont
Let $L$ be the free Lie algebra over $\Z$ on two generators $t_1$, $t_2$.  
The basic commutators of weight one are 
$t_1$ and $t_2$. Put $e(t_1) = 0$, $e(t_2) = 0$.  
Proceeding inductively, suppose that the basic commutators of weight less than $n$ have been defined and ordered as 
$t_1, \ldots t_p$ and that a function $e$ from $\{1, \ldots, p\}$ to the nonnegative 
integers has been defined: $\forall \ i$, $e(i) < i$.  
Take for the basic commutators of weight $n$ the $[t_i,t_j]$, where weight 
$t_i + $ weight $t_j = n$ and $e(i) \leq j < i$.  
Order these commutators in any way and label them 
$t_{p+1} \ldots t_{p+q}$.  
Complete the construction by setting $e([t_i,t_j]) = j$.  Let $B$ be the set of basic commutators thus obtained $-$then $B$ is an additive basis for $L$, the 
\un{Hall basis}.
\index{Hall basis}\\
\endgroup%%------------------------------------<<

\begingroup%%----------------------------------->>
\fontsize{9pt}{11pt}\selectfont
\index{Hilton-Milnor Formula}
\textbf{\small EXAMPLE \ (\un{Hilton-Milnor Formula})} \ 
Let 
$
\begin{cases}
\ X\\
\ Y
\end{cases}
$
be nondegenerate and path connected.  
Put \ 
$
\begin{cases}
\ Z(t_1) = X\\
\ Z(t_2) = Y
\end{cases}
$
and let \ 
$
\begin{cases}
\ \zeta_1:\Sigma Z(t_1) \ra \Sigma X \vee \Sigma Y\\
\ \zeta_2:\Sigma Z(t_2) \ra \Sigma X \vee \Sigma Y
\end{cases}
$
be the inclusions.  For $t \in B$ of weight $n > 1$, write uniquely $t = [t_i,t_j]$, where 
weight $t_i$ $+$ weight $t_j = n$. \  
Via recursion on the weight, put $Z(t) = Z(t_i) \#Z(t_j)$ and let 
$\zeta_t:\Sigma Z(t) \ra \Sigma X$ $\vee$ $\Sigma Y$ be the Whitehead product $[\zeta_i,\zeta_j]$, where
$
\begin{cases}
 \zeta_i:\Sigma Z(t_i) \ra \Sigma X \vee \Sigma Y\\
 \zeta_j:\Sigma Z(t_j) \ra \Sigma X \vee \Sigma Y
\end{cases}
\hspace{-.25cm}. \ 
$
The $\zeta_t$ combine to define a continuous function 
$\zeta = \ds\sum\limits_{t \in B} \Omega \zeta_t$ from 
$(w) \ds\prod\limits_{t \in B} \Omega \Sigma Z(t)$ 
(cf. p. \pageref{4.58}) to $\Omega \Sigma (X \vee Y)$.
Claim: $\zeta$ is a weak homotopy equivalence.  To see this, attach to each $N = 1, 2, \ldots,$ a ``remainder'' 
$R_N = \ds\bigvee\limits_{\substack{i \geq N\\ e(i) < N}} Z(t_i)$.  
Applying the preceding example to 
$\Omega\Sigma \bigl(Z(t_N) \vee \ds\bigvee\limits_{\substack{i > N\\ e(i) < N}} Z(t_i)\bigr)$, 
it follows that the map
\[
\ds\sum\limits_{i = 1}^N \Omega\zeta_i + \Omega \bigl( \bigvee\limits_{\substack{i > N\\ e(i) \leq N}} \zeta_i \bigr):
\prod\limits_{i = 1}^N \Omega\Sigma Z(t_i) \times \Omega\Sigma (R_{N+1}) \ra \Omega\Sigma (X \vee Y)
\]
%%----------------------------------------------------------------------------------------------50
is a weak homotopy equivalence.  
To finish, let $N \ra \infty$ (justified, since the connectivity of $R_{N+1}$ tends to $\infty$ with \mN).
\\ \indent
[Note: \ The isomorphism 
$\zeta_*:\underset{t \in B}{\oplus} \pi_*(\Omega \Sigma Z(t)) \ra \pi_*(\Omega \Sigma (X \vee Y))$ depends on the 
choice of the Hall basis \mB.  
Consult 
Goerss\footnote[2]{\textit{Quart. J. Math.} \textbf{44} (1993), 43-85.}
for an intrinsic description.]\\
\endgroup%%------------------------------------<<

\label{9.99}
A nonempty path connected topological space $X$ is said to be 
\un{homotopically trivial}
\index{homotopically trivial} 
if $X$ is $n$-connected for all $n$, i.e., provided that $\forall \ q > 0$, $\pi_q(X) = 0$.  
Example: A contractible space is homotopically trivial.

Example: Let $X \overset{f}{\ra} Z \overset{g}{\la} Y$ be a 2-sink.  
Assume: $X$ $\&$ $Z$ are homotopically trivial $-$then the arrow $W_{f,g} \ra Y$ is a weak homotopy equivalence.\\

\index{The Wedge of the Broom}
\begingroup%%----------------------------------->>
\fontsize{9pt}{11pt}\selectfont
\textbf{\small EXAMPLE}  \ 
A homotopy equivalence is a weak homotopy equivalence but the converse is false.
\\
\indent\indent (1) \ \ (\un{The Wedge of the Broom}) \ \ 
Consider the subspace $X$ of $\R^2$ consisting of the line segments joining 
$(0,1)$ to $(0,0)$ $\&$ \ $(1/n,0)$\  $(n = 1, 2, \ldots)$ \ $-$then $X$ is contractible, thus it and its base point $(0,0)$ have the same homotopy type.  
But in the pointed homotopy category, $(X,(0,0))$ and $(\{(0,0)\},(0,0))$ are not equivalent.  \ 
\hspace{0.5cm}
\parbox{4cm}{
\begin{tikzpicture}[scale=2.5]%[scale=0.5,shift={(-5,-3)}]
\draw[violet]
(0,1) -- (1,0);
\draw[teal]
(0,1) -- (0.5,0);
\draw[violet]
(0,1) -- (0.333333,0);
\draw[teal]
(0,1) -- (0.25,0);
\draw[violet]
(0,1) -- (0.2,0);
\draw[teal]
(0,1) -- (0.166667,0);
\draw[violet]
(0,1) -- (0.142857,0);
\draw[teal]
(0,1) -- (0.125,0);
\draw[violet]
(0,1) -- (0.111111,0);
\draw[teal]
(0,1) -- (0.1,0);
\draw[violet]
(0,1) -- (0.0909091,0);
\draw[teal]
(0,1) -- (0.0833333,0);
\draw[violet]
(0,1) -- (0.0769231,0);
\draw[teal]
(0,1) -- (0.0714286,0);
\draw[violet]
(0,1) -- (0.0666667,0);
\draw[teal]
(0,1) -- (0.0625,0);
\draw[violet]
(0,1) -- (0.0588235,0);
\draw[black]
(0,1) -- (0,0);
\draw[white] (0, -.1) -- (0,-.11);
\end{tikzpicture}
}
\setlength\parindent{0em}
%\hspace{0.5cm}
\hspace{-1.5cm} Consider $X$ $\vee$ $X$, the subspace of $\R^2$ consisting of the line segments joining \ 
$
\begin{cases}
\ (0,1) \text{ to } (0,0) \ \& \ (1/n,0)\\
\ (0,-1) \text{ to } (0,0) \ \& \ (-1/n,0)
\end{cases}
(n = 1, 2, \ldots), 
$
\hspace{0.5cm}
\parbox{4cm}{
\begin{tikzpicture}[scale=2.5]%[scale=0.5,shift={(-5,-3)}]
\draw[violet]
(0,1) -- (1,0);
\draw[teal]
(0,1) -- (0.5,0);
\draw[violet]
(0,1) -- (0.333333,0);
\draw[teal]
(0,1) -- (0.25,0);
\draw[violet]
(0,1) -- (0.2,0);
\draw[teal]
(0,1) -- (0.166667,0);
\draw[violet]
(0,1) -- (0.142857,0);
\draw[teal]
(0,1) -- (0.125,0);
\draw[violet]
(0,1) -- (0.111111,0);
\draw[teal]
(0,1) -- (0.1,0);
\draw[violet]
(0,1) -- (0.0909091,0);
\draw[teal]
(0,1) -- (0.0833333,0);
\draw[violet]
(0,1) -- (0.0769231,0);
\draw[teal]
(0,1) -- (0.0714286,0);
\draw[violet]
(0,1) -- (0.0666667,0);
\draw[teal]
(0,1) -- (0.0625,0);
\draw[violet]
(0,1) -- (0.0588235,0);
\draw[black]
(0,1) -- (0,0);
\draw[violet]
(0,-1) -- (-1,0);
\draw[teal]
(0,-1) -- (-0.5,0);
\draw[violet]
(0,-1) -- (-0.333333,0);
\draw[teal]
(0,-1) -- (-0.25,0);
\draw[violet]
(0,-1) -- (-0.2,0);
\draw[teal]
(0,-1) -- (-0.166667,0);
\draw[violet]
(0,-1) -- (-0.142857,0);
\draw[teal]
(0,-1) -- (-0.125,0);
\draw[violet]
(0,-1) -- (-0.111111,0);
\draw[teal]
(0,-1) -- (-0.1,0);
\draw[violet]
(0,-1) -- (-0.0909091,0);
\draw[teal]
(0,-1) -- (-0.0833333,0);
\draw[violet]
(0,-1) -- (-0.0769231,0);
\draw[teal]
(0,-1) -- (-0.0714286,0);
\draw[violet]
(0,-1) -- (-0.0666667,0);
\draw[teal]
(0,-1) -- (-0.0625,0);
\draw[violet]
(0,-1) -- (-0.0588235,0);
\draw[black]
(0,-1) -- (0,0);
\draw[white] (0, -.1) -- (0,-.11);
\end{tikzpicture}
}
\setlength\parindent{0em}
\hspace{0.5cm}

$-$then $X \vee X$ is path connected and homotopically trivial.  
However, $X \vee X$ is not contractible, so the map that sends $X \vee X$ to $(0,0)$ is a weak homotopy equivalence but not a homotopy equivalence.
\\[-.2cm]
\setlength\parindent{2em}
%$\vspace{0.5cm}$

\index{The Warsaw Circle}
\indent\indent (2) \ (\un{The Warsaw Circle}) 
Consider the subspace $X$ of $\R^2$ consisting of the union of \\
\vspace{0.1cm}
$
%\left\{(x,y):
\{(x,y):
\begin{cases}
\ x = 0, -2 \leq y \leq 1\\
\ 0 \leq x \leq 1, y = -2\\
\ x = 1, -2 \leq y \leq 0
\end{cases}
%\right\}
\
$
and $\{(x,y): 0 < x \leq 1, y = \sin(2\pi/x)\}$ 
\qquad
\parbox{4cm}{
\begin{tikzpicture}[x=3cm,xscale=2]%[scale=5]%[scale=0.5,shift={(-5,-3)}]
\draw[blue,domain=0.01:0.63661977,samples=5000] plot (\x-0.01,{sin((1/\x)r)});
\draw[black]
(0,1) -- (0,-2);
\draw[black]
(0,-2) -- (0.63661977,-2);
\draw[black]
(0.62661977,-2) -- (0.62661977,1);
\end{tikzpicture}
}
$-$then $X$ is path connected and homotopically trivial.  However, $X$ is not contractible, so the map that sends $X$ to $(0,0)$ is a weak homotopy equivalence but not a homotopy equivalence.\\ 
\endgroup%%------------------------------------<<

\begingroup%%----------------------------------->>
\fontsize{9pt}{11pt}\selectfont
\textbf{\small FACT}  \ 
Let $p:X \ra B$ be a Hurewicz fibration, where $X$ and $B$ are path connected and $X$ is nonempty.  
Suppose that $[p]$ is both a monomorphism and an epimorphism in \bHTOP $-$then $p$ is a weak homotopy equivalence.\\
\endgroup%%------------------------------------<<

A continuous function 
$f:(X,A) \ra (Y,B)$ 
is said to be a 
\un{relative $n$-equivalence}
\index{relative $n$-equivalence} 
$(n \geq 1)$ provided that the sequence 
$* \ra \pi_0(X,A) \ra \pi_0(Y,B)$ 
is exact and $\forall \ x_0 \in A$, 
$f_*:\pi_q(X,A,x_0) \ra \pi_q(Y,B,f(x_0))$ 
is bijective for $1 \leq q < n$ and surjective for $q = n$.\\

\begin{proposition} \quad %29
Suppose that 
\hspace{0.25cm}
$
\begin{cases}
\ X_1\\[-.15cm]
\ X_2
\end{cases}
$
$\&$ 
\hspace{0.25cm}
$
\ 
\begin{cases}
\ Y_1\\[-.15cm]
\ Y_2
\end{cases}
$
are open subspaces of 
\hspace{0.25cm}
$
\begin{cases}
\ X\\[-.15cm]
\ Y
\end{cases}
$ 
with 
\hspace{0.5cm} 
%%----------------------------------------------------------------------------------------------51
$
\begin{cases}
\ X = X_1 \cup X_2\\[-.15cm]
\ Y = Y_1 \cup Y_2
\end{cases}
\hspace{-.25cm}. \ 
\hspace{0.5cm}
$
Let $f:X \ra Y$ be a continuous function such that 
\ 
$
\begin{cases}
\ X_1 = f^{-1}(Y_1)\\[-.15cm]
\ X_2 = f^{-1}(Y_2)
\end{cases}
\hspace{-.25cm}. \ 
$
Fix $n \geq 1$.  
Assume: $f:(X_i,X_1 \cap X_2) \ra (Y_i,Y_1 \cap Y_2)$ is a relative $n$-equivalence $(i = 1,2)$ $-$then 
$f:(X,X_i) \ra (Y,Y_i)$ is a  relative $n$-equivalence $(i = 1,2)$.
\end{proposition}

[This is the content of the result on p. \pageref{4.59}.]\\

A continuous function $f:(X,A) \ra (Y,B)$ is said to be a 
\un{relative weak homotopy} \un{equivalence}
\index{relative weak homotopy equivalence} 
if $f$ is a relative $n$-equivalence $\forall \ n \geq 1$.  
Example: Let $p:X \ra B$ be a Serre fibration, where $B$ is path connected and $X$ is nonempty $-$then $\forall \ b \in B$, 
the arrow $(X,X_b) \ra (B,b)$ is a relative weak homotopy equivalence.\\

\label{5.0k}
\textbf{\small LEMMA}  \ 
Let $f:(X,A) \ra (Y,B)$ be a continuous function.  
Assume: $f:A \ra B$ and $f:X \ra Y$ are weak homotopy equivalences $-$then $f:(X,A) \ra (Y,B)$ is a relative weak homotopy equivalence.\\

\begin{proposition} \quad %30
Suppose that 
\ 
$
\begin{cases}
\ X_1\\[-.15cm]
\ X_2
\end{cases}
\& \ 
$
\ 
$
\begin{cases}
\ Y_1\\[-.15cm]
\ Y_2
\end{cases}
$
are open subspaces of 
\ 
$
\begin{cases}
\ X\\[-.15cm]
\ Y
\end{cases}
$ 
with 
$
\begin{cases}
\ X = X_1 \cup X_2\\[-.15cm]
\ Y = Y_1 \cup Y_2
\end{cases}
\hspace{-.2cm}. \ 
$
Let $f:X \ra Y$ be a continuous function such that  
$
\begin{cases}
\ X_1 = f^{-1}(Y_1)\\[-.15cm]
\ X_2 = f^{-1}(Y_2)
\end{cases}
\hspace{-.2cm}. \ 
$
Assume: 
$
\begin{cases}
\ f:X_1 \ra Y_1\\[-.15cm]
\ f:X_2 \ra Y_2
\end{cases}
$ 
$\&$ $f:X_1 \cap X_2 \ra Y_1 \cap Y_2$ are weak homotopy equivalences $-$then $f:X \ra Y$ is a weak homotopy equivalence.
\end{proposition}

[The lemma implies that $f(X_i,X_1 \cap X_2) \ra (Y_i,Y_1 \cap \ Y_2)$ is a relative weak homotopy equivalence $(i = 1, 2)$.  Therefore, on the basis of Proposition 29, $f:(X,X_i) \ra (Y,Y_i)$ is a relative weak homotopy equivalence $(i = 1, 2)$.  Since a given $x \in X$ belongs to at least one of the $X_i$, this suffices (modulo low dimensional details).]\\

\label{12.22}
\label{14.23a}
\label{14.77}
\label{14.135}
\label{14.151}
Application: Let 
\begin{tikzcd}%[sep=large]
{X} \ar{d}  &{Z} \ar{l}[swap]{f} \ar{d} \ar{r}{g} &{Y} \ar{d}\\
{X^\prime} &{Z^\prime} \ar{l}{f^\prime} \ar{r}[swap]{g^\prime} &{Y^\prime}
\end{tikzcd}
be a \cd in which the vertical arrows are weak homotopy equivalences $-$then the arrow 
$M_{f,g} \ra M_{f^\prime,g^\prime}$ is a weak homotopy equivalence.

[Note: If in addition 
$
\begin{cases}
\ f\\[-.15cm]
\ f^\prime
\end{cases}
$
are closed cofibrations, then the arrow 
$X \sqcup_g Y \ra X^\prime \sqcup_{g^\prime} Y^\prime$ is a weak homotopy equivalence 
(cf. $\S 3$, Proposition 18).]\\

\begingroup%%----------------------------------->>
\fontsize{9pt}{11pt}\selectfont
\textbf{\small FACT}  \ 
Let 
$
\begin{cases}
\ X\\
\ Y
\end{cases}
$
be topological spaces and let $f:X \ra Y$ be a continuous function.  
Assume: $\sV = \{V\}$ is an open covering of $Y$ which is closed under finite intersections such that 
$\forall \ V \in \sV$, $f:f^{-1}(V) \ra V$ is a weak homotopy equivalence $-$then
$f:X \ra Y$  is a weak homotopy equivalence.
\\ \indent
%%----------------------------------------------------------------------------------------------52
\label{4.63}
[Use Zorn on the collection of subspaces $B$ of $Y$ that have the following properties: $B$ is a union of elements of $\sV$, $f:f^{-1}(B) \ra B$ is a weak homotopy equivalence, and $\forall \ V \in \sV$, 
$f:f^{-1}(B \cap V) \ra B \cap V$ is a weak homotopy equivalence.  Order this collection by inclusion and fix a maximal element $B_0$.  
Claim: $B_0 = Y$.  If not, choose $V \in \sV$: $V \not\subset B_0$ and consider $B_0 \cup V$.]\\
\endgroup%%------------------------------------<<

\begingroup%%----------------------------------->>
\fontsize{9pt}{11pt}\selectfont
\textbf{\small SUBLEMMA} \quadx
Let $f \in C(X,Y)$ and suppose given continuous functions 
$
\begin{cases}
\ \phi:\bS^{n-1} \ra X\\
\ \psi:\bD^n \ra Y
\end{cases}
$ 
with $f \circx \phi = \restr{\psi}{\bS^{n-1}}$ $-$then there exists a neighborhood \mU of $\bS^{n-1}$ in $\bD^n$ and continuous functions 
$
\begin{cases}
\ \ov{\phi}:U \ra X\\
\ \ov{\psi}:\bD^n \ra Y
\end{cases}
$
such that $\restr{\ov{\phi}}{\bS^{n-1}} = \phi$ and $f \circx \ov{\phi} = \restr{\ov{\psi}}{U}$, where 
$\psi \simeq \ov{\psi} \  \rel \  \bS^{n-1}$.

[Let $U = \{x: 1/2 < \norm{x} \leq 1\}$ and put $\ov{\phi}(x) = \phi(x/\norm{x})$ $(x \in U)$.  
Write
$
v(x) = 
\begin{cases}
\ x   \hspace{1.15cm} (\norm{x} \leq 1)\\
\ x/\norm{x} \hspace{0.35cm}  (\norm{x} \geq 1) 
\end{cases}
\hspace{-.25cm}. \ 
$
Define $H:I\bD^n \ra Y$ by $H(x,t) = \psi(v((1 + t)x))$ and take $\ov{\psi} = H \circx i_1$.]\\
\endgroup%%------------------------------------<<

\label{3.34}
\begingroup%%----------------------------------->>
\fontsize{9pt}{11pt}\selectfont
\textbf{\small LEMMA}  \ 
Suppose that 
$
\begin{cases}
\ X_1\\
\ X_2
\end{cases}
\& \ 
$
\ 
$
\begin{cases}
\ Y_1\\
\ Y_2
\end{cases}
$
are subspaces of 
$
\ 
\begin{cases}
\ X\\
\ Y
\end{cases}
$
with \\
$
\begin{cases}
\ X = \itr X_1 \cup \itr X_2\\
\ Y = \itr Y_1 \cup \itr Y_2
\end{cases}
. \ 
$
Let $f:X \ra Y$ be a continuous function such that 
$
\begin{cases}
\ f(X_1) \subset Y_1\\
\ f(X_2) \subset Y_2
\end{cases}
\hspace{-.25cm}.
$
Assume: 
$
\begin{cases}
\ f:X_1 \ra Y_1\\
\ f:X_2 \ra Y_2
\end{cases}
$
$\&$ $f:X_1 \cap X_2 \ra Y_1 \cap Y_2$ are weak homotopy equivalences 
$-$then $f:X \ra Y$ is a weak homotopy equivalence.
\\ \indent
[In the notation employed at the end of $\S 3$, given continuous functions 
$
\begin{cases}
\ \phi:\dot{I}^q \ra X\\
\ \psi:{I}^q \ra Y
\end{cases}
$
such that $f \circx \phi = \restr{\psi}{\dot{I}^q}$, it is enough to find a continuous function $\Phi:I^q \ra X$ such that 
$\restr{\Phi}{\dot{I}^q} = \phi$ and $f \circx \Phi \simeq \psi \rel \dot{I}^q$.  
This can be done by a subdivision argument.  The trick is to consider 
$
\begin{cases}
\ \phi^{-1}(X - \itr X_1) \cup \ov{\psi^{-1}(Y - Y_1)}\\
\ \phi^{-1}(X - \itr X_2) \cup \ov{\psi^{-1}(Y - Y_2)}
\end{cases}
\hspace{-.25cm}. \ 
$
These sets are closed.  However, they need not be disjoint and the point of the sublemma is to provide an escape for this difficulty.]\\
\endgroup%%------------------------------------<<

\begingroup%%----------------------------------->>
\fontsize{9pt}{11pt}\selectfont
\textbf{\small EXAMPLE}  \ 
In the usual topology, take $Y = \R$, $Y_1 = \Q$, $Y_2 = \PP$; in the discrete topology, take 
$X = \R$, $X_1 = \Q$, $X_2 = \PP$ $-$then the identity map $X \ra Y$ is not a weak homotopy equivalence, yet the restrictions 
$
\begin{cases}
\ X_1 \ra Y_1\\
\ X_2 \ra Y_2
\end{cases}
, \ 
$
$X_1 \cap X_2 \ra Y_1 \cap Y_2$ are weak homotopy equivalences.\\
\endgroup%%------------------------------------<<

\begingroup%%----------------------------------->>
\fontsize{9pt}{11pt}\selectfont
\textbf{\small FACT}  \ 
Let
$
\begin{cases}
\ X\\
\ Y
\end{cases}
$
be topological spaces and let $f:X \ra Y$ be a continuous function.  Suppose that 
$
\begin{cases}
\ \sU = \{U_i: i \in I\}\\
\ \sV = \{V_i: i \in I\}
\end{cases}
$
are open coverings of 
$
\begin{cases}
\ X\\
\ Y
\end{cases}
$
such that $\forall \ i$: $f(U_i) \subset V_i$.  
Assume: For every nonemtpy finite subset $F \subset I$, the induced map 
$\ds\bigcap\limits_{i \in F} U_i \ra \ds\bigcap\limits_{i \in F} V_i$ is a weak homotopy equivalence $-$then $f$ is a weak homotopy equivalence.\\
\endgroup%%------------------------------------<<

Topological spaces 
$
\begin{cases}
\ X\\[-.15cm]
\ Y
\end{cases}
$
are said to have the same 
\un{weak homotopy type}
\index{weak homotopy type} 
if there exists a topological space $Z$ and weak homotopy equivalences 
$
\begin{cases}
\ f:Z \ra X\\[-.15cm]
\ g:Z \ra Y
\end{cases}
\hspace{-.25cm}. \ 
$
The relation of having the same weak homotopy type is an equivalence relation.

%%----------------------------------------------------------------------------------------------53
[Note: \ One can always replace $Z$ by a CW resolution $K \ra Z$ , hence 
$
\begin{cases}
\ X\\[-.15cm]
\ Y
\end{cases}
$
have the same weak homotopy type iff they admit CW resolutions 
$
\begin{cases}
\ K \ra X\\[-.15cm]
\ K \ra Y
\end{cases}
$
by the same CW complex $K$.]\\

\begingroup%%----------------------------------->>
\fontsize{9pt}{11pt}\selectfont
Transitivity is the only issue.  For this, let $X_1$, $X_2$, $X_3$, be topological spaces, let $K, L$ be CW complexes, and consider the diagram
\begin{tikzcd}%[sep=large]
&{K} \ar{ld}[swap]{f_1}\ar{rd}{f_2} &&{L}\ar{ld}[swap]{g_2} \ar{rd}{g_3}\\
{X_1}  &&{X_2} &&{X_3}
\end{tikzcd}
, where 
$
\begin{cases}
\ f_1\\
\ f_2
\end{cases}
,
$
$
\begin{cases}
\ g_2\\
\ g_3
\end{cases}
$
are weak homotopy equivalences.  Since $(K,f_2)$ and $(L,g_2)$ are both CW resolutions of $X_2$, there is a homotopy equivalence $\phi:K \ra L$ such that $f_2 \simeq g_2 \circx \phi$ (cf. p. \pageref{4.60}).  Thus $g_3 \circx \phi:K \ra X_3$ is a weak homotopy equivalence, so $X_1$ and $X_3$ have the same weak homotopy type.\\
\endgroup%%------------------------------------<<

\begingroup%%----------------------------------->>
\fontsize{9pt}{11pt}\selectfont
\textbf{\small EXAMPLE}  \ 
Two aspherical spaces having the same isomorphic fundamental groups have the same weak homotopy type.

[Note: \ A path connected topological space $X$ is said to be 
\un{aspherical}
\index{aspherical (path connected topological space)}
provided that $\forall \ q > 1$, $\pi_q(X) = 0$.  
Example: If $X$ is path connected and metrizable with $\dim X = 1$, then $X$ is aspherical.]\\
\endgroup%%------------------------------------<<

Let $X$ be in \bTOP/\mB.  Assume that the projection $p:X \ra B$ is surjective $-$then $p$ is said to be a 
\un{quasifibration}
\index{quasifibration} 
if $\forall \ b \in B$, the arrow $(X,X_b) \ra (B,b)$ is a relative weak homotopy equivalence.  
If $p:X \ra B$ is a quasifibration, then $\forall \ b_0 \in B$, $\forall \ x_0 \in X_{b_0}$, there is an exact sequence
\[
\cdots \ra \pi_2(B) \ra \pi_1(X_{b_0}) \ra \pi_1(X) \ra \pi_1(B) \ra \pi_0(X_{b_0}) \ra \pi_0(X) \ra \pi_0(B).
\]
\vspace{0.15cm}

\textbf{\small LEMMA}  \ 
Let $p: X \ra B$ be a Serre fibration.  Suppose that $B$ is path connected and $X$ is nonempty $-$then $p$ is a quasifibration.\\

\begingroup%%----------------------------------->>
\fontsize{9pt}{11pt}\selectfont
\textbf{\small EXAMPLE}  \ 
Take $X = ([-1,0] \times \{1\}) \cup (\{0\} \times [0,1]) \cup ([0,1] \times \{0\})$, $B = [-1,1]$, and let $p$ be the vertical projection $-$then $p$ is a quasifibration ($X$ and $B$ are contractible, as are all the fibers) but $p$ is neither a Serre fibration nor a Dold fibration.
\\ \indent
[Note: \ The pullback of a Serre fibration is a Serre fibration, i.e., Proposition 4 is valid with ``Hurewicz'' replaced by ``Serre''.  
This fails for quasifibrations.  Let $B^\prime = [0,1]$ and define $\Phi^\prime: B^\prime \ra B$ by 
$
\Phi(t) = 
\begin{cases}
\ t \sin(1/t) \hspace{0.35cm} \ \  (t > 0)\\
\ \ 0 \hspace{1.57cm} (t = 0)
\end{cases}
$
$-$then the projection $p^\prime :X^\prime  \ra B^\prime$ is not a quasifibration (consider $\pi_0$).]\\
\endgroup%%------------------------------------<<

\begin{proposition} \  %31
Let $p:X \ra B$ be a quasifibration, where $B$ is path connected $-$then the fibers of $p$ have the same weak homotopy type.
\end{proposition}

%%----------------------------------------------------------------------------------------------54
[Using the mapping track $W_p$, factor $p$ as $q \circx \gamma$ and note that $\forall \ b \in B$, $\gamma$ induces a weak homotopy equivalence $X_b \ra q^{-1}(b)$.  But $q:W_p \ra B$ is a Hurewicz fibration and since $B$ is path connected, the fibers of $q$ have the same weak homotopy type (cf. p. \pageref{4.61}).]\\

\begingroup%%----------------------------------->>
\fontsize{9pt}{11pt}\selectfont
\textbf{\small EXAMPLE}  \ 
Let $B = [0,1]^n$ $(n \geq 1)$.  Put $X = B \times B - \Delta_B$ and let $p$ be the vertical projection $-$then $p$ is not a quasifibration (cf. p. \pageref{4.62}).\\
\endgroup%%------------------------------------<<

\textbf{\small LEMMA}  \ 
Let $p:X \ra B$ be a continuous function.  Suppose that $O \subset B$ and $p_O:X_O \ra O$ is a quasifibration $-$then the arrow 
$(X,X_O) \ra (B,O)$ is a relative weak homotopy equivalence iff $\forall \ b \in O$, the arrow 
$(X,X_b) \ra (B,b)$ is a relative weak homotopy equivalence.\\

\begin{proposition} \  %32
Let $X$ be in \bTOP/\mB.  Suppose that 
$
\begin{cases}
\ O_1\\[-.15cm]
\ O_2
\end{cases}
$
are open subspaces of $B$ with $B = O_1 \cup O_2$.  Assume: 
$
\begin{cases}
\ p_{O_1}:X_{O_1} \ra O_1\\[-.15cm]
\ p_{O_2}:X_{O_2} \ra O_2\
\end{cases}
\& \ p_{O_1 \ \cap \ O_2}:X_{O_1 \ \cap \ O_2} \ra O_1\  \cap \ O_2
$
are quasifibrations $-$then $p:X \ra B$ is a quasifibration.

[From the lemma, the arrows 
$(X_{O_i}, X_{O_1 \cap O_2}) \ra (O_i, O_1 \ \cap\  O_2)$ are relative weak homotopy equivalences $(i = 1, 2)$. 
 Therefore the arrow $(X,X_{O_i}) \ra (B,O_i)$ is a relative weak homotopy equivalence $(i = 1,2)$ (cf. Proposition 29).  
Since $p$ is clearly surjective, another appeal to the lemma completes the proof.]\\
\end{proposition}

Application: Let $X$ be in \bTOP/\mB.  
Suppose that $\sO = \{O_i:i \in I\}$ is an open covering of $B$ which is closed under finite intersections.  
Assume: $\forall \ i$, $p_{O_i}:X_{O_i} \ra O_i$ is a quasifibration $-$then $p:X \ra B$ is a quasifibration.

[The argument is the same as that indicated on p. \pageref{4.63} for weak homotopy equivalences.]

[Note: \ This is the local-global principle for quasifibrations.  Here, numerability is irrelevant.]\\

\begingroup%%----------------------------------->>
\fontsize{9pt}{11pt}\selectfont
\textbf{\small EXAMPLE}  \ 
Let $X = \R^2$ be equipped with the following topology.  Basic neighborhoods of $(x,y)$, where 
$
\begin{cases}
\ x \leq 0\\
\ x \geq 1
\end{cases}
\& -\infty < y < \infty
$
or
$
\begin{cases}
\ 0 < x < 1 \ \& \ y > 0\\
\ 0 < x < 1 \ \&\  y < 0
\end{cases}
\hspace{-.25cm},
$
are the usual neighborhoods but the basic neighborhoods of $(x,0)$, where $0 < x < 1$, are the open semicircles centered at 
$(x,0)$ of radius $< \min\{x,1-x\}$ that lie in the closed upper half plane.  
\[
\begin{tikzpicture}[scale=2.5]%[x=5cm]%[scale=5]%[scale=0.5,shift={(-5,-3)}]
\draw[black]
(-1,0) -- (2,0);

\draw[black]
(0,-.5) -- (0,1);

\draw[black]
(1,-.25) -- (1,.75);

\draw[red]
(0,0) -- (1,0);

\begin{scope}
\clip (0,0) rectangle (1, 1); 
\draw[violet]
(0.5,0) circle(0.5);
\draw[violet]
(0.5,0) circle(0.375);
\draw[violet]
(0.5,0) circle(0.25);
\draw[violet]
(0.5,0) circle(0.125);

\draw[green]
(0.25,0) circle(0.25);
\draw[green]
(0.25,0) circle(0.1875);
\draw[green]
(0.25,0) circle(0.125);
\draw[green]
(0.25,0) circle(0.0625);

\draw[green]
(0.75,0) circle(0.25);
\draw[green]
(0.75,0) circle(0.1875);
\draw[green]
(0.75,0) circle(0.125);
\draw[green]
(0.75,0) circle(0.0625);

\end{scope}

%\begin{axis}
%xmin = -1,
%xmax = 1,
%ymin = -1,
%ymax = 1,
%axis lines = center,
%xlabel = $x$, 
%ylabel = $y$, 
%\end{axis}
\end{tikzpicture}
\]
Take $B = \R^2$ (usual topology) $-$then the identity map $p:X \ra B$ is not a quasifibration (since 
$\pi_1(B) = 0$, $\pi_1(X) \neq 0)$ and the fibers are points).  
Put
$
\begin{cases}
\ O_1 = \{(x,y):x > 0\}\\
\ O_2 = \{(x,y):x < 1\}
\end{cases}
:
$
$
\begin{cases}
\ O_1\\
\ O_2
\end{cases}
$
are open subspaces of $B$ with $B = O_1 \cup O_2$.  Moreover 
$
\begin{cases}
\ X_{O_1}\\
\ X_{O_2}
\end{cases}
$
are contractible, thus 
$
\begin{cases}
\ p_{O_1}:X_{O_1} \ra O_1\\
\ p_{O_2}:X_{O_2} \ra O_2
\end{cases}
$
are quasifibrations.  However, $p_{O_1 \ \cap \ O_2}: X_{O_1 \ \cap \ O_2} \ra O_1 \cap O_2$ is not a quasifibration.\\
\endgroup%%------------------------------------<<

%%----------------------------------------------------------------------------------------------55
\begingroup%%----------------------------------->>
\fontsize{9pt}{11pt}\selectfont
\textbf{\small FACT}  \ 
Let $p:X \ra B$ be a surjective continuous function, where $B = \colimx B^n$ is $\tT_1$.  Assume: $\forall \ n$, 
$p^{-1}(B^n) \ra B^n$ is a quasifibration $-$then $p$ is a quasifibration.\\
\endgroup%%------------------------------------<<

Let $A$ be a subspace of $X$, $i:A \ra X$ the inclusion.

\indent\indent (WDR) $A$ is said to be a 
\un{weak deformation retract}
\index{weak deformation retract}
of $X$ if there is a homotopy $H:IX \ra X$ such that 
$H \circx i_0 = \id_X$, 
$H \circx i_t(A) \subset A$ $(0 \leq t \leq 1)$, and 
$H \circx i_1(X) \subset A$.  

[Note: \ Define $r:X \ra A$ by $i \circx r = H \circx i_1$ $-$then $i \circx r \simeq \id_X$ and $r \circx i \simeq \id_A$.]

A strong deformation retract is a weak deformation retract.  The comb is a weak deformation retract of $[0,1]^2$ 
(consider the homotopy $H((x,y),t) = (x,(1 - t)y))$ but the comb is not a retract of $[0,1]^2$.

[Note: \ A pointed space $(X,x_0)$ is contractible to $x_0$ in $\bTOP_*$ iff $\{x_0\}$ is a weak (or strong) deformation retract of $X$.  
The broom with base point $(0,0)$ is an example of a pointed space which is contractible in \bTOP but not in 
$\bTOP_*$.  Therefore a deformation retract need not be a weak deformation retract.]\\

\begingroup%%----------------------------------->>
\fontsize{9pt}{11pt}\selectfont
On a subspace $A$ of $X$ such that the inclusion $A \ra X$ is a cofibration, ``strong'' = ``weak''.\\
\endgroup%%------------------------------------<<

\begin{proposition} \  %33
Let $p:X \ra B$ be a surjective continuous function.  Suppose that $O$ is a subspace of $B$ for which 
$p_O:X_O \ra O$ is a quasifibration and 
$
\begin{cases}
\ O\\[-.15cm]
\ X_O
\end{cases}
$
is a weak deformation retract of 
$
\begin{cases}
\ B\\[-.15cm]
\ X
\end{cases}
 \hspace{-.2cm}
,
$
say 
$
\begin{cases}
\ \rho:B \ra O\\[-.15cm]
\ \tau:X \ra X_O
\end{cases}
\hspace{-.25cm}. \ 
$
Assume: $p \circx r = \rho \circx p$ and $\forall \ b \in B$, $\restr{r}{X_b}$ is a weak homotopy equivalence 
$X_b \ra X_{\rho(b)}$ $-$then $p:X \ra B$ is a quasifibration.
\end{proposition}

[Given $b \in B$, $r:(X,X_b) \ra (X_O,X_{\rho(b)})$, as a map of pairs, is a relative weak homotopy equivalence and, by the assumption, the diagram
\begin{tikzcd}%[sep=large]
{(X,X_b)} \ar{d} \ar{r} &{(X_O,X_{\rho(b)})}\ar{d}\\
{(B,b)} \ar{r} &{(O,\rho(b))}
\end{tikzcd}
commutes.]\\

Application: Let
\begin{tikzcd}%[sep=large]
{X} \ar{d}  &{Z} \ar{l}[swap]{f} \ar{d} \ar{r}{g} &{Y} \ar{d}\\
{X^\prime} &{Z^\prime} \ar{l}{f^\prime} \ar{r}[swap]{g^\prime} &{Y^\prime}
\end{tikzcd}
be a \cd in which the vertical arrows are quasifibrations.  
Assume: $\forall \ z^\prime \in Z^\prime$, 
$
\begin{cases}
\ \restr{f}{Z_{z^\prime}}\\[-.15cm]
\ \restr{g}{Z_{z^\prime}}
\end{cases}
$
is a weak homotopy equivalence
$
\begin{cases}
\ Z_{z^\prime} \ra X_{f^\prime(z^\prime)}\\[-.15cm]
\ Z_{z^\prime} \ra Y_{g^\prime(z^\prime)}
\end{cases}
$
$-$then the arrow $M_{f,g} \ra M_{f\prime,g^\prime}$ is a quasifibration.\\
\vspace{0.25 cm}

%%----------------------------------------------------------------------------------------------56
\begin{proposition} \  %34
Let
\begin{tikzcd}%[sep=large]
{X} \ar{d}  &{Z} \ar{l}[swap]{f} \ar{d} \ar{r}{g} &{Y} \ar{d}\\
{X^\prime} &{Z^\prime} \ar{l}{f^\prime} \ar{r}[swap]{g^\prime} &{Y^\prime}
\end{tikzcd}
be a \cd in which the left vertical arrow is a surjective Hurewicz fibration and the right vertical arrow is a quasifibration.  
Assume: 
\begin{tikzcd}%[sep=large]
{X} \ar{d}  &{Z} \ar{l}[swap]{f} \ar{d}\\
{X^\prime} &{Z^\prime} \ar{l}{f^\prime}
\end{tikzcd}
is a pullback square, 
$
\begin{cases}
\ f\\[-.15cm]
\ f^\prime
\end{cases}
$ 
are closed cofibrations, and $\forall \ z^\prime \in Z^\prime$, $\restr{g}{Z_{z^\prime}}$ is a weak homotopy equivalence 
$Z_{z^\prime} \ra Y_{g^\prime(z^\prime)}$ $-$then the induced map 
$X \sqcup_g Y \ra X^\prime \sqcup_{g^\prime} Y^\prime$ is a quasifibration.
\end{proposition}

[Consider the \cd
\begin{tikzcd}%[sep=large]
{M_{f,g}} \ar{d}[swap]{\phi} \ar{r}{\mu} &{M_{f^\prime,g^\prime}} \ar{d}{\phi^\prime}\\
{X \sqcup_g Y} \ar{r}[swap]{\nu} &{X^\prime \sqcup_{g^\prime } Y^\prime }
\end{tikzcd}
.  \ Since 
$
\begin{cases}
\ f\\[-.15cm]
\ f^\prime
\end{cases}
$ 
are cofibrations, 
$
\begin{cases}
\ \phi\\
\ \phi^\prime
\end{cases}
$ 
are homotopy equivalences (cf. $\S 3$, Proposition 18) and, by the above, $\mu$ is a quasifibration.  
Thus it need only be shown that $\forall \ m^\prime \in M_{f^\prime,g^\prime}$, the arrow 
$\mu^{-1}(m^\prime) \ra \nu^{-1}(\phi^\prime(m^\prime))$ is a weak homotopy equivalence, which can be done by examining cases.]
\\

The conclusion of Proposition 34 cannot be strengthened to ``Hurewicz fibration''.  
To see this, take 
$X = [-1,0] \times [0,1]$, 
$Y = [0,2] \times [0,2]$, 
$Z = \{0\} \times [0,1]$, 
$X^\prime = [-1,0]$, 
$Y^\prime  = [0,2]$, 
$Z^\prime  = \{0\}$, 
let 
$
\begin{cases}
\ f:Z \ra X\\[-.15cm]
\ g:Z \ra Y
\end{cases}
\hspace{-.25cm}, \ 
$ 
$
\begin{cases}
\ f^\prime:Z^\prime \ra X^\prime\\[-.15cm]
\ g^\prime:Z^\prime \ra Y^\prime
\end{cases}
$ 
be the inclusions, and let $X \ra X^\prime$, $Z \ra Z^\prime$, $Y \ra Y^\prime$ be the vertical projections $-$then
$X \sqcup_g Y = X \cup Y$, 
$X^\prime \sqcup_{g^\prime} Y^\prime = X^\prime \cup Y^\prime$, 
and the induced map 
$X \cup Y \ra X^\prime \cup Y^\prime$ is the vertical projection.  But it is not a Hurewicz fibration since it fails to have the slicing structure property (cf. p. \pageref{4.64}).\\

\begingroup%%----------------------------------->>
\fontsize{9pt}{11pt}\selectfont
\label{14.75}
\index{Cone Construction}
\textbf{\small EXAMPLE \ (\un{Cone Construction})} \ 
Fix nonempty topological spaces $X$, $Y$ and let $\phi:X \times Y \ra Y$ be a continuous function.  
Define $E$ by the pushout square
\begin{tikzcd}[sep=large]
{X \times Y}  \ar{d} \ar{r}{\phi} &{Y} \ar{d}\\
{\Gamma X \times Y} \ar{r} &{E}
\end{tikzcd}
.  \ 
Assume: $\forall \ x \in X$, $\phi_x:\{x\} \times Y \ra Y$ is a weak homotopy equivalence.  Consider the commutative diagram 
\begin{tikzcd}[sep=large]
{\Gamma X \times Y} \ar{d}  &{X \times Y} \ar{l} \ar{d} \ar{r}{\phi} &{Y} \ar{d}\\
{\Gamma X} &{X} \ar{l} \ar{r} &{*}
\end{tikzcd}
.  \ 
Since the arrows 
$X \ra \Gamma X$, $X \times Y \ra \Gamma X \times Y$ are closed cofibrations, all the hypotheses of Proposition 34 are met.  
Therefore the induced map 
$E \ra \Sigma X$ is a quasifibration.
\\ \indent
[Note: \ The same construction can be made in the pointed category provided that $(X,x_0)$ is wellpointed
%%----------------------------------------------------------------------------------------------57
with $\{x_0\} \subset X$ closed.]\\
\endgroup%%------------------------------------<<

\begingroup%%----------------------------------->>
\fontsize{9pt}{11pt}\selectfont
\index{Dold-Lashof Construction}
\textbf{\small EXAMPLE \  (\un{Dold-Lashof Construction})} \ 
Let $G$ be a topological semigroup with unit in which the operations of left and right translation are homotopy equivalences.  
Let $p:X \ra B$ be a quasifibration.  Assume: There is a right action 
$
\begin{cases}
\ X \times G \ra X\\
\ (x,g) \ra x \cdot g
\end{cases}
$
such that $p(x\cdot g) = p(x)$ and the arrow 
$
\begin{cases}
\ G \ra X_{p(x)}\\
\ g \ra x \cdot g
\end{cases}
$
is a weak homotopy equivalence.  
Define $\ov{X}$ by the pushout square 
\begin{tikzcd}[sep=large]
{X \times G} \ar{d} \ar{r} &{X}\ar{d}\\
{\Gamma X \times G} \ar{r} &{\ov{X}}
\end{tikzcd}
and put $\ov{B} = C_p$.  Since the diagram 
\begin{tikzcd}[sep=large]
{\Gamma X \times G} \ar{d}  &{X \times G} \ar{l} \ar{d} \ar{r} &{X} \ar{d}\\
{\Gamma X} &{X} \ar{l} \ar{r} &{B}
\end{tikzcd}
commutes, Proposition 34 implies that $\ov{p}:\ov{X} \ra \ov{B}$ is a quasifibration.
Represent a generic point of  $\ov{X}(\ov{B})$  by the symbol $[x,t,g]$ $([x,t])$ 
(with the obvious understanding at $t = 0$ or $t = 1$), 
so $\ov{p}[x,t,g] = [x,t]$.  The assignment
$
\begin{cases}
\ \ov{X} \times G \ra \ov{X}\\
\ ([x,t,g],h) = [x,t,gh]
\end{cases}
$
is unambiguous and satisfies the algebraic conditions for a right action of $G$ on $\ov{X}$ but it is not necessarily continuous.
The resolution is to place a smaller topology on $\ov{X}$. 
Let 
$t:\ov{X} \ra [0,1]$ be the function $[x,t,g] \ra t$; 
let 
$x:t^{-1}(]0,1[) \ra X$ be the function $[x,t,g] \ra x$; 
let 
$g:t^{-1}([0,1[) \ra G$ be the function $[x,t,g] \ra g$; 
let 
$x\cdot g:t^{-1}(]0,1]) \ra X$ be the function $[x,t,g] \ra x \cdot g$.  
Definition: The 
\un{coordinate topology}
\index{coordinate topology}
on $\ov{X}$ is the initial topology determined by $t, x, g, x \cdot g$.  The injection 
$
\begin{cases}
\ X \ra \ov{X}\\
\ x \ra [x,1,e]
\end{cases}
$
is an embedding, as is the injection 
$
\begin{cases}
\ G \ra \ov{X}\\
\ g \ra [x,t,g]
\end{cases}
(t \neq 0,1). 
$
Moreover, $G$ acts continuously and $\forall \ \ov{x} \in \ov{X}$, the arrow 
$
\begin{cases}
\ G \ra \ov{X}_{\ov{p}(\ov{x})}\\
\ g \ra \ov{x} \cdot g
\end{cases}
$
is a weak homotopy equivalence.  
Now equip $\ov{B}$ with its coordinate topology (cf. p. \pageref{4.65}) $-$then 
$\ov{p}: \ov{X}\ra \ov{B}$ is continuous and remains a quasifibration (apply Propositions 32 and 33 to 
$
\begin{cases}
\ O_1 = \{[x,t]:0 < t \leq 1\}\\
\ O_2 = \{[x,t]:0 \leq t < 1\}
\end{cases}
).
$
In other words, $(\ov{X},\ov{B})$ satisfies the same conditions as $(X,B)$ and there is a commutative diagram 
$
\begin{tikzcd}[sep=large]
{X} \ar{d} \ar{r} &{\ov{X}} \ar{d}\\
{B} \ar{r} &{\ov{B}}
\end{tikzcd}
,
$
where $X \ra \ov{X}$ is inessential (consider 
$
H:
\begin{cases}
\ IX \ra \ov{X}\\
\ (x,t) \ra [x,t,e]
\end{cases}
).
$
\\ \indent
Example: Let $G$ be a topological group $-$then $\ov{G}$ (coordinate topology) is homeomorphic to 
$G *_c G$ (coarse join).\\
\endgroup%%------------------------------------<<

Let $G$ be a topological group, $X$ a topological space.  Suppose that $X$ is a right $G$-space:
$
\begin{cases}
\ X \times G \ra X\\
\ (x,g)  \ra x\cdot g
\end{cases}
$
$-$then the projection 
$X \ra X/G$ is an open map and $X/G$ is Hausdorff iff $X \times_{X/G} X$ is closed in $X \times X$.  
The continuous function 
$\theta:X \times G \ra X \times_{X/G} X$ defined by $(x,g) \ra (x,x\cdot g)$ is surjective.  
It is injective iff the action is free,  
\index{free action}
i.e., iff $\forall \ x \in X$, the stabilizer $G_x = \{g: x\cdot g = x\}$ of $x$ in $G$ is trivial.  
A free right $G$-space $X$ is said
%%----------------------------------------------------------------------------------------------58
to be 
\un{principal}
\index{principal (free right $G$-space)}
provided  that $\theta$ is a homeomorphism or still, that the division function
$
d:
\begin{cases}
\ X \times_{X/G} X \ra G\\
\ (x,x\cdot g) \ra g
\end{cases}
$
is continuous.

Let $X$ be in $\bTOP/B$ $-$then $X$ is said to be a 
\un{principal $G$-space over $B$}
\index{principal $G$-space over $B$}
if $X$ is a principal $G$-space, $B$ is a trivial \mG-space, the projection $p:X \ra B$ is open, surjective, and equivariant, and $G$ operates transitively on the fibers.  
There is a commutative triangle
\begin{tikzcd}%[sep=large]
&{X} \ar{ld}\ar{d}\\
{X/G} \ar{r}  &{B}
\end{tikzcd}
and the arrow $X/G \ra B$ is a homeomorphism.  
$\bPRIN_{B,G}$ 
\index{$\bPRIN_{B,G}$} 
is the category whose objects are the principal $G$-spaces over $B$ and whose morphisms are the equivariant continuous functions over $B$.  
If $\Phi^\prime \in C(B^\prime,B)$, then for every $X$ in $\bPRIN_{B,G}$ there is a pullback square 
\begin{tikzcd}%[sep=large]
{X^\prime} \ar{d} \ar{r}{f^\prime} &{X} \ar{d}\\
{B^\prime} \ar{r}[swap]{\Phi^\prime}  &{B}
\end{tikzcd}
with $X^\prime = B^\prime \times_B X$ in $\bPRIN_{B^\prime,G}$ and $f^\prime$ equivariant.\\

\textbf{\small LEMMA}  \ 
Every morphism in $\bPRIN_{B,G}$ is an isomorphism.

[Note: \ The objects in $\bPRIN_{B,G}$  which are isomorphic to $B \times G$ (product topology) are said to be 
\un{trivial}
\index{trivial (principal $G$-spaces over $B$)}.  
It follows from the lemma that the trivial objects are precisely those that admit a section.]\\

Application: Let 
$
\begin{cases}
\ X^\prime\\[-.15cm]
\ X
\end{cases}
$
be in 
$
\begin{cases}
\ \bPRIN_{B^\prime,G}\\[-.15cm]
\ \bPRIN_{B,G}
\end{cases}
\hspace{-.25cm};
$
let $f^\prime \in C(X^\prime,X)$, $\Phi^\prime \in C(B^\prime,B)$.  
Assume: $f^\prime$ is equivariant and 
$p \circx f^\prime = \Phi^\prime \circx p^\prime$ $-$then the commutative diagram 
\begin{tikzcd}%[sep=large]
{X^\prime} \ar{d} \ar{r}{f^\prime} &{X} \ar{d}\\
{B^\prime} \ar{r}[swap]{\Phi^\prime}  &{B}
\end{tikzcd}
is a pullback square.

[Compare this diagram with the pullback square defining the fiber product.]\\

Let $X$ be in $\bTOP/B$ $-$then $X$ is said to be a 
\un{$G$-bundle over $B$}
\index{G-bundle over $B$}
if $X$ is a free right $G$-space, $B$ is a trivial $G$-space, the projection $p:X \ra B$ is open, surjective, and equivariant, 
and there exists an open covering $\sO = \{O_i: i \in I\}$ of $B$ such that $\forall \ i$, $X_{O_i}$ is equivariantly homeomorphic to 
$O_i \times G$ over $O_i$.  
Since the division function is necessarily continuous and $G$ operates transitively on the fibers, $X$ is a principal $G$-space over $B$.  
If $\sO$ can be chosen numerable, then $X$ is said to be a 
\un{numerable $G$-bundle over $B$}
\index{numerable $G$-bundle over $B$}
(a condition that is automatic when $B$ is a paracompact Hausdorff space, e.g., a CW complex).  
$\bBUN_{B,G}$ 
\index{$\bBUN_{B,G}$ } 
is the full subcategory of $\bPRIN_{B,G}$ whose objects are the numerable $G$-bundles over $B$.  Each $X$ in $\bBUN_{B,G}$ is numerably locally trivial with fiber $G$ and the local-global principal implies 
%%----------------------------------------------------------------------------------------------59
that the projection $X \ra B$ is a Hurewicz fibration.  There is a functor 
$I:\bBUN_{B,G} \ra \bBUN_{IB,G}$ that sends $p:X \ra B$ to $Ip:IX \ra IB$, where 
$(x,t)\cdot g = (x \cdot g,t)$.\\

\begingroup%%----------------------------------->>
\fontsize{9pt}{11pt}\selectfont
\textbf{\small EXAMPLE}  \ 
A $G$-bundle over $B$ need not be numerable.  For instance, take $G = \R$ $-$then every object in $\bBUN_{B,\R}$ 
admits a section ($\R$ being contractible), hence is trivial.  
Let now $X$ be the subset of $\R^3$ defined by the equation $x_1x_3 + x_2^2 = 1$ and let $\R$ act on $X$ via 
$(x_1,x_2,x_3)\cdot t = (x_1, x_2+tx_1,x_3-2tx_2-t^2x_1)$.  
$X$ is an $\R$-bundle over $X/\R$ but it is not numerable.  For if it were, then there would exist a section 
$X/\R \ra X$, an impossibility since $X/\R$ is not Hausdorff.\\
\endgroup%%------------------------------------<<

\begingroup%%----------------------------------->>
\fontsize{9pt}{11pt}\selectfont
\textbf{\small FACT}  \ 
Suppose that $X$ is a $G$-bundle over $B$ $-$then the projection $p:X \ra B$ is a Serre fibration (cf. p. \pageref{4.66}) which is $\Z$-orientable if $B$ and $G$ are path connected.\\
\endgroup%%------------------------------------<<

Let 
$
\begin{cases}
\ X^\prime\\[-.15cm]
\ X
\end{cases}
$
be in 
$
\begin{cases}
\ \bBUN_{B^\prime,G}\\[-.15cm]
\ \bBUN_{B,G}
\end{cases}
\hspace{-.25cm}. \ 
$
Write $X^\prime \times_G X$ for the orbit space $(X^\prime \times X)/G$ $-$then there is a commutative diagram 
\begin{tikzcd}%[sep=large]
{X^\prime \times X} \ar{d} \ar{r}&{X^\prime} \ar{d}\\
{X^\prime \times_G X}\ar{r} &{B^\prime}
\end{tikzcd}
which is a pullback square.  
As an object in $\bTOP/B^\prime$, $X^\prime \times_G X$ is numerably locally trivial with fiber $X$ so, e.g., has the SEP if $X$ is contractible.  
The $s^\prime \in \sec_{B^\prime}(X^\prime \times_G X)$ correspond bijectively to the equivariant 
$f^\prime \in C(X^\prime,X)$.  
As an object in $\bTOP/B^\prime \times B$, $X^\prime \times_G X$ is numerably locally trivial with fiber \mG.  
Given $\Phi^\prime \in C(B^\prime,B)$, there exists an equivariant $f^\prime \in C(X^\prime,X)$ rendering the diagram 
\begin{tikzcd}%[sep=large]
{X^\prime} \ar{d} \ar{r}{f^\prime} &{X} \ar{d}\\
{B^\prime} \ar{r}[swap]{\Phi^\prime}  &{B}
\end{tikzcd}
commutative iff the arrow 
$
\begin{cases}
\ B^\prime \ra B^\prime \times B\\
\ b^\prime \ra (b^\prime,\Phi^\prime(b^\prime))
\end{cases}
$
admits a lifting 
\begin{tikzcd}%[sep=large]
&{X^\prime \times_G X} \ar{d}\\
{B^\prime} \ar[dashed]{ru} \ar{r}  &{B^\prime \times B}
\end{tikzcd}
.\\

\index{Theorem: Covering Homotopy Theorem}
\index{Covering Homotopy Theorem}
\textbf{\small COVERING HOMOTOPY THEOREM} \quadx
Let 
$
\begin{cases}
\ X^\prime\\[-.15cm]
\ X
\end{cases}
$
be in 
$
\begin{cases}
\ \bBUN_{B^\prime,G}\\[-.15cm]
\ \bBUN_{B,G}
\end{cases}
\hspace{-.25cm}. \ 
$
Suppose that $f^\prime:X^\prime \ra X$ is an equivariant continuous function and 
$h:IB^\prime \ra B$ is a homotopy with $p \circx f^\prime = h \circx i_0 \circx p^\prime$ $-$then there exists an equivariant homotopy 
$H:IX^\prime \ra X$ such that $H \circx i_0 = f^\prime$ and for which the diagram 
\begin{tikzcd}%[sep=large]
{IX^\prime} \ar{d} \ar{r}{H} &{X} \ar{d}\\
{IB^\prime} \ar{r}[swap]{h}  &{B}
\end{tikzcd}
commutes.\\

[Take $\Phi^\prime = h \circx i_0$ to get a lifting 
\begin{tikzcd}%[sep=large]
&{X^\prime \times_G X} \ar{d}\\
{B^\prime} \ar[dashed]{ru} \ar{r}  &{B^\prime \times B}
\end{tikzcd}
and a commutative diagram
%%----------------------------------------------------------------------------------------------60
\begin{tikzcd}%[sep=large]
{B^\prime} \ar{d}[swap]{i_0} \ar{r} &{IX^\prime \times_G X} \ar{d}\\
{IB^\prime} \ar{r}  &{IB^\prime \times B}
\end{tikzcd}
.  The projection $IX^\prime \times_G X \ra IB^\prime \times B$ is a Hurewicz fibration, thus the diagram has a filler 
$IB^\prime \ra IX^\prime \times_G X$ and this guarantees the existence of \mH.]\\

Application: Let $X$ be in $\bBUN_{B,G}$.  
Suppose that 
$
\begin{cases}
\ \Phi_1^\prime\\[-.15cm]
\ \Phi_2^\prime
\end{cases}
\in C(B^\prime,B)
$
are homotopic $-$then
$
\begin{cases}
\ X_1^\prime\\[-.15cm]
\ X_2^\prime
\end{cases}
$
are isomorphic in $\bBUN_{B^\prime,G}$.\\
\vspace{0.25cm}

\begingroup%%----------------------------------->>
\fontsize{9pt}{11pt}\selectfont
\textbf{\small FACT}  \ 
The functor $I:\bBUN_{B,G} \ra \bBUN_{IB,G}$ has a representative image.\\
\endgroup%%------------------------------------<<

\label{5.48}

The relation ``isomorphic to'' is an equivalence relation on $\Ob\bBUN_{B,G}$.  
Call $k_G B$ the ``class'' of equivalence classes arising therefrom $-$then $k_G B$ is a ``set'' (see below).  
Since for any $\Phi^\prime \in C(B^\prime,B)$ and each $X$ in $\bBUN_{B,G}$, the isomorphism class $[X^\prime]$ 
of $X^\prime$ in $\bBUN_{B^\prime,G}$ depends only on the homotopy class $[\Phi^\prime]$ of $\Phi^\prime$, 
$k_G$ is a cofunctor $\bHTOP \ra \bSET$.  
A topological space $B_G$ is said to be a 
\un{classifying space}
\index{classifying space (for $G$ of a topological space)} 
for $G$ 
if $B_G$ represents $k_G$, i.e., if there exists a natural isomorphism $\Xi:[-,B_G] \ra k_G$, an 
$X_G \in \Xi_{B_G}(\id_{B_G})$ being a 
\un{universal}
\index{universal numerable $G$-bundle}
numerable $G$-bundle over $B_G$.  From the definitions, 
$\forall \ \Phi \in C(B,B_G)$, $\Xi_B[\Phi] = [X]$, where $X$ is defined by the pullback square
\begin{tikzcd}%[sep=large]
{X} \ar{d} \ar{r} &{X_G} \ar{d}\\
{B} \ar{r}[swap]{\Phi}  &{B_G}
\end{tikzcd}
and $\Phi$ is the 
\un{classifying map}.
\index{classifying map (for $G$ of a topological space)}
\\

\indent\indent (UN) \ Assume that
$
\begin{cases}
\ \Xi^\prime \ra [-,B_G^\prime] \ra k_G\\[-.15cm]
\ \Xi\pp \ra [-,B_G\pp] \ra k_G
\end{cases}
$
are natural isomorphisms $-$then there exist mutually inverse homotopy equivalences 
$
\begin{cases}
\ \Phi^\prime: B_G^\prime \ra B_G\pp\\[-.15cm]
\ \Phi\pp: B_G\pp \ra B_G^\prime
\end{cases}
$
such that \\
$
\begin{cases}
\ k_G[\Phi^\prime]([X_G\pp]) \ = [X_G^\prime]\\[-.15cm]
\ k_G[\Phi\pp]([X_G^\prime]]) = [X_G\pp]
\end{cases}
\hspace{-.25cm}.
$
\\
\vspace{0.25cm}

\begingroup%%----------------------------------->>
\fontsize{9pt}{11pt}\selectfont
Recall that the members of a class are sets, therefore $k_G B$ is not a class but rather a conglomerate.  
Still, $\bBUN_{B,G}$ has a small skeleton $\ov{\bBUN}_{B,G}$.  
Indeed, any $X$ in $\bBUN_{B,G}$ is isomorphic to $B \times G$.  
Here, the topology on $B \times G$ depends on $X$ and is in general not the product topology but the action is the same
$((b,g)\cdot h = (b,gh))$.  
Thus one can modify the definition of $k_G$ and instead take for $k_G B$ the set 
$\Ob\ov{\bBUN}_{B,G}$.\\
\endgroup%%------------------------------------<<

\begin{proposition} \  %35
Suppose that there exists a $B_G$ in \bTOP and an $X_G$ in $\bBUN_{B_G,G}$ such that $X_G$ is contractible $-$then 
$k_G$ is representable.
\end{proposition}

[Define a natural transformation $\Xi:[-,B_G] \ra k_G$ by assigning to a given homotopy class $[\Phi]$ 
($\Phi \in C(B,B_G)$) the isomorphism class $[X]$ of the numerable $G$-bundle $X$ over \mB 
%%----------------------------------------------------------------------------------------------61
defined by the pullback square
\begin{tikzcd}%[sep=small]
{X} \ar{d} \ar{r} &{X_G} \ar{d}\\
{B} \ar{r}[swap]{\Phi}  &{B_G}
\end{tikzcd}
\hspace{-.2cm}
.  The claim is that $\forall \ B$, $\Xi_B:[B,B_G] \ra k_GB$ is bijective.

Surjectivity: Take any $X$ in $\bBUN_{B,G}$ and form $X \times_G X_G$.  
Since $X_G$ is contractible, $X \times_G X_G$ has the SEP, thus $\sec_B(X \times_G X_G)$ is nonempty, so there exists an equivariant $f \in C(X,X_G)$.  
Determine $\Phi \in C(B,B_G)$ from the commutative diagram
\begin{tikzcd}%[sep=large]
{X} \ar{d} \ar{r}{f} &{X_G} \ar{d}\\
{B} \ar{r}[swap]{\Phi}  &{B_G}
\end{tikzcd}
$-$then $\Xi_B[\Phi] = [X]$.

Injectivity: Let $\Phi^\prime$, $\Phi\pp \in C(B,B_G)$ and assume that 
$\Xi_B[\Phi^\prime] = \Xi_B[\Phi\pp]$, say $[X^\prime] = [X\pp]$, where 
\begin{tikzcd}[sep=small]
{X^\prime} \ar{rdd} \ar{rr}{\phi} &&{X\pp} \ar{ldd}\\
\\
&{B}
\end{tikzcd}
, with $\phi$ equivariant.  There are pullback squares
$
\begin{tikzcd}%[sep=large]
{X^\prime} \ar{d} \ar{r}{f^\prime} &{X_G} \ar{d}\\
{B} \ar{r}[swap]{\Phi^\prime}  &{B_G}
\end{tikzcd}
,
$
$
\begin{tikzcd}%[sep=large]
{X\pp} \ar{d} \ar{r}{f\pp} &{X_G} \ar{d}\\
{B} \ar{r}[swap]{\Phi\pp}  &{B_G}
\end{tikzcd}
.
$
\ 
Put $B_0 = B \times ([0,1/2[ \ \cup \ ]1/2,1])$ and define $H_0:\restr{IX^\prime}{B_0} \ra X_G$ by 
$
H_0(x^\prime,t) = 
\begin{cases}
\ f^\prime (x^\prime) \hspace{0.99cm} (t < 1/2)\\ 
\ f\pp \circx \phi(x^\prime)  \hspace{0.25cm}  (t > 1/2)
\end{cases}
\hspace{-.2cm}:\ 
$
$H_0$ is equivariant, hence corresponds to a section $s_0$ of $(IX^\prime \times_G \restr{X_G)}{B_0}$.
Since $B_0$ is a halo of $i_0B \cup i_1B$ in $IB$ and since $IX^\prime \times_G X_G$ has the SEP, $\exists$ 
$s \in \sec_{IB}(IX^\prime \times_G X_G)$ : 
$\restr{s}{B \times (\{0\} \cup \{1\})} = \restr{s_0}{B \times (\{0\} \cup \{1\})}$.  
Translated, this means that there exists an equivariant homotopy $H:IX^\prime \ra X_G$.  
Determine $h:IB \ra B_G$ from the commutative diagram 
\begin{tikzcd}%[sep=large]
{IX^\prime} \ar{d} \ar{r}{H} &{X_G} \ar{d}\\
{IB} \ar{r}[swap]{h}  &{B_G}
\end{tikzcd}
$-$then
$
\begin{cases}
\ h \circx i_0 = \Phi^\prime\\ 
\ h \circx i_1 = \Phi\pp
\end{cases}
$
$\implies [\Phi^\prime] = [\Phi\pp]$.]\\

\begingroup%%----------------------------------->>
\fontsize{9pt}{11pt}\selectfont
The converse to Proposition 35 is also true:  In order that $k_G$ be representable, it is necessary that $X_G$ be contractible.  Thus let $X_G^\infty$ be the numerable $G$-bundle over $B_G^\infty$ produced by the Milnor construction $-$then 
$X_G^\infty$ is contractible, so $\Xi^\infty$ is a natural isomorphism.  As the same holds for $\Xi$ by assumption, there are pullback squares 
\begin{tikzcd}[sep=large]
{X_G} \ar{d} \ar{r}{f} &{X_G^\infty} \ar{d}\\
{B_G} \ar{r}[swap]{\Phi}  &{B_G^\infty}
\end{tikzcd}
,
\begin{tikzcd}[sep=large]
{X_G^\infty} \ar{d} \ar{r}{f^\infty} &{X_G} \ar{d}\\
{B_G^\infty} \ar{r}[swap]{\Phi^\infty}  &{B_G}
\end{tikzcd}
and $\Phi^\infty \circx \Phi \simeq \id_{B_G}$.  Owing to the covering homotopy theorem, $f^\infty \circx f$ is equivariantly homotopic to an isomorphism 
\begin{tikzcd}[sep=small]
{X_G} \ar{rdd} \ar{rr}{\phi} &&{X_G} \ar{ldd}\\
\\
&{B_G}
\end{tikzcd}
.  But $\phi$ is
%%----------------------------------------------------------------------------------------------62
necessarily inessential, $X_G^\infty$ being contractible.\\
\endgroup%%------------------------------------<<

\begingroup%%----------------------------------->>
\fontsize{9pt}{11pt}\selectfont
\textbf{\small EXAMPLE}  \ 
Let $E$ be an infinite dimensional Hilbert space $-$then its general linear group $\bGL(E)$ is contractible 
(cf. p. \pageref{4.66a}).
Any compact Lie group $G$ can be embedded as a closed subgroup of $\bGL(E)$.  So, if 
$X_G = \bGL(E)$, $B_G = \bGL(E)/G$, then $B_G$ is a classifying space for \mG and $X_G$ is universal.
\\ \indent
[$B_G$ is a paracompact Hausdorff space.  Local triviality of $X_G$ is a consequence of a generality due to Gleason, 
viz: 
Suppose that $G$ is a compact Lie group and $X$ is a Hausdorff principal $G$-space which is completely regular $-$then $X$, as an object in \bTOP/\mB ($B = X/G$), is a $G$-bundle.]\\
\endgroup%%------------------------------------<<

\begingroup%%----------------------------------->>
\fontsize{9pt}{11pt}\selectfont
\textbf{\small EXAMPLE}  \ 
Let $G$ be a noncompact connected semisimple Lie group with a finite center, $K \subset G$ a maximal compact subgroup.  
The coset space $K\backslash G$ is contractible, being diffeomorphic to some $\R^n$.  Let $\Gamma$ be a discrete subgroup of $G$.  
Assume: $\Gamma$ is cocompact and torsion free $-$then $\Gamma$ operates on 
$K\backslash G$ by right translation and $K\backslash G$  is a numerable $\Gamma$-bundle over $K\backslash G/\Gamma$.  
So if 
$X_\Gamma = K\backslash G$, $B_\Gamma = K\backslash G/\Gamma$, then $B_\Gamma$ is a classifying space for 
$\Gamma$ and $X_\Gamma$ is universal.
\\ \indent
[Note: \ $B_\Gamma$ is a compact riemannian manifold.  Its universal covering space is $X_\Gamma$, thus 
$B_\Gamma$ is aspherical and of homotopy type $(\Gamma,1)$.]\\
\endgroup%%------------------------------------<<

\index{Milnor Construction}
\textbf{\small MILNOR CONSTRUCTION} \ 
Let $G$ be a topological group.  Consider the subset of 
$([0,1] \times G)^\omega$ made up of the strings $\{(t_i,g_i)\}$ for which 
$\ds\sum\limits_i t_i = 1$ $\&$ $\#\{i:t_i = \neq 0\} < \omega$.  
Write
$\{(t_i^\prime,g_i^\prime)\} \sim \{(t_i\pp,g_i\pp)\}$ iff $\forall \ i$, $t_i^\prime = t_i\pp$ and at those $i$ such that
$t_i^\prime = t_i\pp$ is positive, $g_i^\prime = g_i\pp$.  
Call $X_G^\infty$ the resulting set of equivalence classes.  
Define coordinate functions $t_i$ and $g_i$ by 
$
t_i= 
\begin{cases}
\ X_G^\infty \ra [0,1]\\[-.15cm]
\ x \ra t_i(x)
\end{cases}
$
and 
$
g_i= 
\begin{cases}
\ t_i^{-1}(]0,1]) \ra G\\[-.15cm]
\ x \ra g_i(x)
\end{cases}
, \ 
$
where $x = [(t_i(x),g_i(x))]$.  
The 
\un{Milnor topology}
\index{Milnor topology} 
on $X_G^\infty$ 
is the initial topology determined by the $t_i$ and $g_i$.  Thus topologized, $X_G^\infty$ is a right $G$-space:
$
\begin{cases}
\ X_G^\infty \times G \ra X_G^\infty\\[-.15cm]
\ (x,g) \ra x \cdot g
\end{cases}
. \ 
$
Here, $t_i(x \cdot g) = t_i(x)$ and $g_i(x \cdot g) = g_i(x)g$.  Let $B_G^\infty$ be the orbit space $X_G^\infty/G$.

[Note: \ Put 
$X_G^0 = G$, $X_G^n = G *_c \cdots *_c  G$, the $(n+1)$-fold coarse join of $G$ with itself.  
One can identify $X_G^n$ with $\{x: \forall \ i \geq n+1, t_i(x) = 0\}$.  
Each $X_G^n$ is a zero set in $X_G^\infty$ and there is an equivariant embedding $X_G^n \ra X_G^{n+1}$.  
So, $X_G^0 \subset X_G^1 \subset \cdots $ is an expanding sequence of topological spaces and the colimit in \bTOP associated with this data is $X_G^\infty$ equipped with the final topology determined by the inclusions $X_G^n \ra X_G^\infty$.  
The colimit topology is finer than the Milnor topology and in general, there is no guarantee that the 
$G$-action $(x,g) \ra x\cdot g$ remains continuous.]\\
\indent\indent (M) \ $X_G^\infty$  is a numerable $G$-bundle over $B_G^\infty$.

[It is clear that $X_G^\infty$  is a principal $G$-space.  Write $O_i$ for the image of $t_i^{-1}(]0,1])$ 
under the projection $X_G^\infty \ra B_G^\infty$ $-$then $\{O_i\}$ is a countable cozero set covering of $B_G^\infty$, hence is numerable (cf. p. \pageref{4.67}).  On the other hand, $\forall \ i$, $\sec_{O_i}(\restr{X_G^\infty}{O_i})$ is nonempty.  To see this, 
%%----------------------------------------------------------------------------------------------63
define a continuous fiber preserving function 
$f_i:\restr{X_G^\infty}{O_i} \ra \restr{X_G^\infty}{O_i}$ by $f_i(x) = x \cdot g_i(x)^{-1}$: $\forall \ g \in G$, $f_i(x \cdot g) = f_i(x)$.  
Consequently, $f_i$ drops to a section $s_i:O_i \ra \restr{X_G^\infty}{O_i}$, therefore $\restr{X_G^\infty}{O_i}$  is trivial.]\\
\indent\indent (D) \ $X_G^\infty$ is contractible.

[Let $\Delta_G^\infty$ be the subset of $X_G^\infty$ consisting of those $x$ such that 
$g_i(x) = e$ if $t_i(x) > 0$ $-$then $\Delta_G^\infty$ is contractible, so one need only construct a homotopy 
$H:IX_G^\infty \ra X_G^\infty$ such that 
$H \circx i_0 = \id_{X_G^\infty}$ and 
$H \circx i_1(X_G^\infty) \subset \Delta_G^\infty$.  
Put $U_k = \tau_k^{-1}(]0,1])$ and $A_k = \tau_k^{-1}(1)$, where 
$\tau_k = \ds\sum\limits_{i \leq k} t_i$.  
Define $H_k^\prime:IU_k \ra U_k$ by 
\[
t_i(H_k^\prime(x,t)) = 
\begin{cases}
\ \ds\frac{t + (1 - t)\tau_k(x)}{\tau_k(x)} \hspace{0.9cm} (i \leq k)\\
\ (1 - t)t_i(x) \hspace{1.7cm}    (i > k)
\end{cases}
\]
and $g_i(H_k^\prime(x,t)) = g_i(x)$ when $t_i(H_k^\prime(x,t)) > 0$.  
Note that $H_k^\prime(x,0) = x$, 
$H_k^\prime(x,1) \in A_k$, and 
$x \in \Delta_G^\infty$ $\implies$ $H_k^\prime(x,t) \in \Delta_G^\infty$ $(0 \leq t \leq 1)$.  
Define $H_k\pp:IA_k \ra A_{k+1}$ by 
\[
t_i(H_k\pp(x,t)) = 
\begin{cases}
\ (1 - t)t_i(x)  \hspace{0.5cm} \ \ (i \leq k)\\
\ t \hspace{2.5cm}  (i = k + 1)\\
\ 0 \hspace{2.45cm} (i > k+1)
\end{cases}
\]
and 
$
g_i(H_k\pp(x,t)) = 
\begin{cases}
\ g_i(x)  \hspace{0.5cm}  (i \leq k)\\[-.15cm]
\ e \hspace{0.85cm} \ \ (i = k+1)
\end{cases}
$
when $t_i(H_k\pp(x,t) > 0$.  Note that 
$H_k\pp(x,0) = x$, 
$H_k\pp(x,1) \in \Delta_G^\infty$, and 
$x \in \Delta_G^\infty$ $\implies$ $H_k\pp(x,t) \in \Delta_G^\infty$ $(0 \leq t \leq 1)$.  
Combine
$
\begin{cases}
\ H_k^\prime\\[-.15cm]
\ H_k\pp
\end{cases}
$
and obtain a homotopy $H_k:IU_k \ra U_{k+1}$ such that 
$H_k(x,0) = x$, 
$H_k(x,1) \in \Delta_G^\infty$, and  
$x \in \Delta_G^\infty$ $\implies$ $H_k(x,t) \in \Delta_G^\infty$ $(0 \leq t \leq 1)$.  
Proceeding recursively, write $G_1 = H_1$ and 
\[
G_{k+1}(x,t) = 
\begin{cases}
\ G_{k}(x,t) \hspace{4.22cm} \  \  (2/3 \leq \tau_k(x) \leq 1)\\
\ H_{k+1}(G_{k}(x,t),2t(2-3\tau_k(x))) \hspace{0.75cm} (1/2 \leq \tau_k(x) \leq 2/3)\\
\ H_{k+1}(G_{k}(x,2t(3\tau_k(x) - 1)),t) \hspace{0.75cm} (1/3 \leq \tau_k(x) \leq 1/2)\\
\ H_{k+1}(x,t) \hspace{3.80cm}  \ \  (0 \leq \tau_k(x) \leq 1/3)
\end{cases}
\]
to get a sequence of homotopies $G_k:IU_k \ra U_{k+1}$ such that 
$\restr{G_{k+1}}{I\tau_k^{-1}(]2/3,1])} = \restr{G_k}{I\tau_k^{-1}(]2/3,1])}$ and 
$G_k(x,0) = x$, $G_k(x,1) \in \Delta_G^\infty$ .  Take for $H$ the homotopy 
$IX_G^\infty \ra X_G^\infty$ that agrees on $I\tau_k^{-1}(]2/3,1])$ with $G_k$.]

[Note: \ The argument shows that $\Delta_G^\infty$ is a weak deformation retract of $X_G^\infty$.]\\

\begingroup%%----------------------------------->>
\fontsize{9pt}{11pt}\selectfont
\index{Borel Construction}
\textbf{\small FACT \ (\un{Borel Construction})} \ 
Let $X$ be in $\bBUN_{B,G}$.  There is a pullback square
\begin{tikzcd}[sep=large]
{X} \ar{d} \ar{r}{f} &{X_G^\infty} \ar{d}\\
{B} \ar{r}[swap]{\Phi}  &{B_G^\infty}
\end{tikzcd}
and since $f$ is equivariant, the continuous function 
$
\begin{cases}
\ X \ra X \times X_G^\infty\\
\ x \ra (x,f(x))
\end{cases}
$
induces a map $B \ra X \times_G X_G^\infty$, which is a homotopy equivalence 
(cf. p. \pageref{4.67a}).\\
\endgroup%%------------------------------------<<

%%----------------------------------------------------------------------------------------------64
\begingroup%%----------------------------------->>
\fontsize{9pt}{11pt}\selectfont
\textbf{\small FACT}  \ 
Let $\alpha:G \ra K$ be a continuous homomorphism $-$then $\alpha$ determines a continuous function 
$f_\alpha:X_G^\infty \ra X_K^\infty$ such that $f_\alpha(x\cdot g) = f_\alpha(x) \cdot \alpha(g)$.  
There is a commutative diagram 
\begin{tikzcd}[sep=large]
{X_G^\infty} \ar{d} \ar{r}{f_\alpha} &{X_K^\infty} \ar{d}\\
{B_G^\infty} \ar{r}[swap]{\Phi_\alpha}  &{B_K^\infty}
\end{tikzcd}
and $\Phi_\alpha$ is a homotopy equivalence iff $\alpha$ is a homotopy equivalence.\\
\endgroup%%------------------------------------<<

\index{Theorem: Classification Theorem}
\index{Classification Theorem}

\textbf{\small CLASSIFICATION THEOREM} \ \ 
For any topological group $G$, the functor $k_G$ is representable.

[This follows from Proposition 35 and the Milnor construction.]\\

The isomorphism classes of numerable $G$-bundles over $B$ are therefore in a one-to-one correspondence with the elements of $[B,B_G^\infty]$.  By comparison, recall that on general grounds the isomorphism classes of $G$-bundles over $B$ are in a one-to-one correspondence with the elements of the cohomology set $H^1(B;\bG)$ (\bG the sheaf of $G$-valued continuous functions on $B$).\\

\begingroup%%----------------------------------->>
\fontsize{9pt}{11pt}\selectfont
\textbf{\small LEMMA}  \ 
Suppose that $G$ is metrizable $-$then the Milnor topology on $X_G^\infty$ is metrizable.
\\ \indent
[Fix a metric $d_G$ on $G$: $d_G \leq 1$.  Define a metric $d$ on $X_G^\infty$ by
\[
d(x,y) = \sum\limits_i \min \{t_i(x),t_i(y)\} d_G(g_i(x),g_i(y)) + \bigl( 1 - \sum\limits_i \min \{t_i(x),t_i(y)\} \bigr).
\]
To check the triangle inequality, consider 
$\ds\frac{1}{2}\abs{t_i(x) - t_i(y)} +  \min \{t_i(x),t_i(y)\} d_G(g_i(x),g_i(y))$ and distinguish two cases: 
$t_i(z) \geq \min \{t_i(x),t_i(y)\}$ $\&$ 
$t_i(z) < \min \{t_i(x),t_i(y)\}$.  
In the metric topology, the coordinate functions are continuous, thus the metric topology is finer than the Milnor topology.  
To go the other way, let $\{x_n\}$ be a net in $X_G^\infty$ such that $x_n \ra x$ in the Milnor topology.  
Claim: $x_n \ra x$ in the metric topology.  Fix $\epsilon > 0$.  Since
$\ds\sum\limits_i t_i(x) = 1$, $\exists \ N$: $\ds\sum\limits_1^N t_i(x) > 1 - \ds\frac{\epsilon}{4}$.  
Choose $n_0$: $\forall \ n \geq n_0$ $\&$ $1 \leq i \leq N$, 
$\abs{t_i(x_n) - t_i(x)} < \ds\frac{\epsilon}{4N}$ and 
$t_i(x) > 0$ $\implies$ $t_i(x_n) > 0$ with 
$d_G(g_i(x_n),g_i(x)) < \ds\frac{\epsilon}{4N}$, from which
%\[
%d(x_n,x) \leq \ds\sum\limits_1^N \min\{t_i(x_n),t_i(x)\} d_G(g_i(x_n),g_i(x)) + 
%\bigl(1 - \ds\sum\limits_1^N \min\{t_i(x_n),t_i(x)\} \bigr) 
%\leq \ds\frac{\epsilon}{4} + 1 - \left(1 - \ds\frac{\epsilon}{2}\right) < \epsilon.]
%\]
\begin{align*}
d(x_n,x) \ 
&\leq 
\sum\limits_1^N \min\{t_i(x_n),t_i(x)\} d_G(g_i(x_n),g_i(x)) + \big(1 - \ds\sum\limits_1^N \min\{t_i(x_n),t_i(x)\} \big) 
\\[11pt]
&\leq 
\frac{\epsilon}{4} + 1 - \bigg(1 - \frac{\epsilon}{2}\bigg)
\\[11pt]
&<
\epsilon.]
\end{align*}
\label{11.5}
\\[-.75cm] \indent
[Note: \ $B_G^\infty$ is also metrizable.  For this, it need only be shown that $B_G^\infty$ is locally metrizable and paracompact (cf. p. \pageref{4.68}).  Local metrizability follows from the fact that 
$\restr{X_G^\infty}{O_i}$ is homeomorphic to 
$O_i \times G$.  Since a metrizable space is paracompact and since $\{O_i\}$ is numerable, $B_G^\infty$ admits a neighborhood finite closed covering by paracompact subspaces, hence is a paracompact Hausdorff space (cf. p. \pageref{4.69}).]\\
\endgroup%%------------------------------------<<

\begingroup%%----------------------------------->>
\fontsize{9pt}{11pt}\selectfont
\textbf{\small EXAMPLE}  \ 
$X_G^\infty$ in the colimit topology is contractible.  This is because $\forall \ n$, the inclusion 
$X_G^n \ra X_G^{n+1}$ is a cofibration (cf. p. \pageref{4.70}) and inessential, thus the result on p. \pageref{4.71} can be applied.  
Consequently, if the underlying topology on $G$ is locally compact and Hausdorff (e.g., if $G$ is Lie), then 
$\colim (X_G^n \times G) = (\colim X_G^n) \times G$, so  $X_G^\infty$ in the colimit topology is a right $G$-space.  
As such, it is a numerable $G$-bundle
%%----------------------------------------------------------------------------------------------65
over $B_G^\infty$, which is therefore a classifying space for $G$ (cf. Proposition 35).  
While the topology on $B_G^\infty$ arising in this fashion is finer than that produced by the Milnor construction, it has the advantage of being ``computable''.  
For example, let $G$ be $\bS^0$, $\bS^1$, or $\bS^3$, the multiplicative group of elements of norm one in 
$\R$, $\C$, or $\HH$ $-$then $X_G^n = \bS^n$, $\bS^{2n+1}$, or $\bS^{4n+3}$, hence $X_G^\infty = \bS^\infty$ 
and factoring in the action, 
$B_G^\infty = \bP^\infty(\R)$, $\bP^\infty(\C)$, or $\bP^\infty(\HH)$.  
As a colimit of the $\bS^n$, $\bS^\infty$ is not first countable.  
However, the three topologies on its underlying set coming from the Milnor construction are metrizable, in particular first countable.
\\ \indent
\label{9.24}
[Note: \ Here is another model for $X_G$ and $B_G$ when $G = \bS^0$, $\bS^1$, or $\bS^3$.  
Take an infinite dimensional Banach space \mE over $\R$, $\C$, or $\HH$ and let $S$ be its unit sphere $-$then $S$ is an AR 
(cf. p. \pageref{4.72}), hence contractible (cf. p. \pageref{4.73}), so $X_G = S$ is universal and $B_G = S/G$ is classifying.]\\
\endgroup%%------------------------------------<<

\begingroup%%----------------------------------->>
\fontsize{9pt}{11pt}\selectfont
Let $G$ be a compact Lie group $-$then 
Notbohm\footnote[2]{\textit{J. London Math. Soc.} \textbf{52} (1995), 185-198.}
has shown that the homotopy type of $B_G^\infty$ determines the Lie group isomorphism class of $G$.\\
\endgroup%%------------------------------------<<

\label{5.14a}
Consider $G$ as a pointed space with base point $e$.  Let
$x_G^\infty = [(1,e),(0,e),\ldots]$ be the base point in $X_G^\infty$, 
$b_G^\infty = x_G^\infty \cdot G$ the base point in $B_G^\infty$ $-$then $\forall \ q \geq 0$, 
$\pi_q(G) \approx \pi_{q+1}(B_G^\infty)$.  
Choose a homotopy 
$H:IX_G^\infty \ra X_G^\infty$ such that 
$
\begin{cases}
\ H(x,0) = x_G^\infty\\
\ H(x,1) = x
\end{cases}
. \ 
$
Taking adjoints and projecting leads to a map $X_G^\infty \ra \Theta B_G^\infty$.  
The triangle
\begin{tikzcd}[sep=small]
{X_G^\infty} \ar{rdd} \ar{rr} &&{\Theta B_G^\infty} \ar{ldd}{p_1}\\
\\
&{B_G^\infty}
\end{tikzcd}
commutes, thus there is an arrow $G \ra \Omega B_G^\infty$.\\

\begin{proposition} \ 
The arrow  $G \ra \Omega B_G^\infty$ is a homotopy equivalence.
\end{proposition}

[The map $X_G^\infty \ra \Theta B_G^\infty$ is a homotopy equivalence (by contractibility).  But the projections 
$X_G^\infty \ra B_G^\infty$, 
$\Theta B_G^\infty \overset{p_1}{\ra} B_G^\infty$ are Hurewicz fibrations.  Therefore the map 
$X_G^\infty \ra \Theta B_G^\infty$ is a fiber homotopy equivalence (cf. Proposition 15).]\\

\begingroup%%----------------------------------->>
\fontsize{9pt}{11pt}\selectfont
\textbf{\small EXAMPLE}  \ 
Take $B = \bS^n$ $(n \geq 1)$ $-$then 
$k_G \bS^n \approx$ 
$[\bS^n,B_G^\infty] \approx$ 
$\pi_1(B_G^\infty,b_G^\infty)\backslash[\bS^n,s_n;B_G^\infty,b_G^\infty] \approx$ 
$\pi_1(B_G^\infty,b_G^\infty)\backslash \pi_n(B_G^\infty,b_G^\infty) \approx$ 
$\pi_0(G,e) \backslash \pi_{n-1}(G,e)$, 
i.e., in brief: 
$k_G \bS^n \approx \pi_0(G) \backslash \pi_{n-1}(G)$.\\
\endgroup%%------------------------------------<<

\begingroup%%----------------------------------->>
\fontsize{9pt}{11pt}\selectfont
\textbf{\small LEMMA}  \ 
Suppose that $G$ is an ANR $-$then $X_G^\infty$ and $B_G^\infty$ are ANRs (cf. p. \pageref{4.74}) and the arrow 
$G \ra \Omega B_G^\infty$ is a pointed homotopy equivalence.
\\ \indent
[Being ANRs, $(G,e)$ $\&$ 
$
\begin{cases}
\ (X_G^\infty,x_G^\infty)\\
\ (B_G^\infty,b_G^\infty)
\end{cases}
$
are wellpointed (cf. p. \pageref{4.75}).  
Therefore $X_G^\infty$ is contractible to $x_G^\infty$ in $\bTOP_*$ and the arrow 
$G \ra \Omega B_G^\infty$ is a pointed map.  
But $(\Omega B_G^\infty,j(b_G^\infty))$ is wellpointed (cf. p. 
%%----------------------------------------------------------------------------------------------66
\pageref{4.76}) (actually $\Omega B_G^\infty$ is an ANR (cf. $\S 6$, Proposition 7)), so the arrow 
$G \ra \Omega B_G^\infty$ is a pointed homotopy equivalence (cf. p. \pageref{4.77}).]\\
\endgroup%%------------------------------------<<

\begingroup%%----------------------------------->>
\fontsize{9pt}{11pt}\selectfont
\textbf{\small EXAMPLE}  \ 
Let $G$ be a Lie group $-$then $G$ is an ANR (cf. p. \pageref{4.78}).  Consider $k_G\Sigma B$, where $(B,b_0)$ is nondegenerate and $\Sigma B$ is the pointed suspension.  Thus 
$k_G\Sigma B \approx$ 
$[\Sigma B,B_G^\infty] \approx$ 
$\pi_1(B_G^\infty,b_G^\infty)\backslash [B,b_0;\Omega B_G^\infty,j(b_G^\infty)] \approx$ 
$\pi_0(G,e) \backslash [B,b_0;G,e]$, which, when $G$ is path connected, simplifies to 
$[B,b_0;G,e]$ or still, 
$[B,G]$ (the action of $\pi_1(G,e)$ on $[B,b_0;G,e]$ is trivial).
\\ \indent
[Note: \ Suppose that $G$ is an arbitrary path connected topological group $-$then again 
$k_G\Sigma B \approx$ 
$[B,b_0;\Omega B_G^\infty,j(b_G^\infty)]$.  
However, $\Omega B_G^\infty$ is a path connected H-group, hence
$[B,b_0;\Omega B_G^\infty,j(b_G^\infty)] \approx$ 
$[B,\Omega B_G^\infty]$ and, by Proposition 36, 
$[B,\Omega B_G^\infty] \approx$ 
$[B,G]$.]\\
\endgroup%%------------------------------------<<
%%%%%%%%%%%%%%%%%%%%%%%%%%%%%%%%%%%%%%
%%%%%%%%%%%%%%%%%%%%%%%%%%%%%%%%%%%%%%
%%%%%%%%%%%%%%%%%%%%%%%%%%%%%%%%%%%%%%

\begin{center}
$\S \ 4$
\\[0.5cm]
$\mathcal{REFERENCES}$\\
\end{center}

\[
\textbf{BOOKS}
\]

\begingroup
\fontsize{9pt}{11pt}\selectfont
\setlength\parindent{0 cm}

[1] \quad tom Dieck, T., Kamps, K., and Puppe, D., \textit{Homotopietheorie}, Springer Verlag (1970).
\\[-.2cm]

[2] \quad Husemoller, D., \textit{Fiber Bundles}, Springer Verlag (1994).
\\[-.2cm]

[3] \quad James, I., \textit{General Topology and Homotopy Theory}, Springer Verlag (1984).
\\[-.2cm]

[4] \quad James, I., \textit{Fibrewise Topology}, Cambridge University Press (1989).
\\[-.2cm]

[5] \quad Piccinini, R., \textit{Lectures on Homotopy Theory}, North Holland (1992).
\\[-.2cm]

%[6]  {\cyr Postniknov M., Lektsii po Algebraicheskoy Topologii}[{\cyr Osmpvy Teorii Gomotopiy}] {\cyr} (1984).\\
[6] 
\hspace{.25cm}
{\fontencoding{OT2}\selectfont
Postnikov M. Lektsii po Algebraicheskoy Topologii [Osmpvy Teorii Gomotopi\u], Nauka} (1984).
\\
\endgroup

\[
\textbf{ARTICLES}
\]

\begingroup
\fontsize{9pt}{11pt}\selectfont
\setlength\parindent{0 cm}

%[1] \quad Booth, P., Classifying Spaces for a General Theory of Fibrations.\\
[1] \quad Booth, P., Fibrations and Classifying Spaces: An Axiomatic Approach I, II, 
\textit{Cahiers de Topologie et} 

\hspace{.8cm}\textit{G\'eom\'etrie Diff\'erentielle Cat\'egoriques}, (1998), 83-116 and 181-203.
\\[-.2cm]

[2] \quad Dold, A., Partitions of Unity in the Theory of Fibrations, 
\textit{Ann. of Math.} \textbf{78} (1963), 223-255.
\\[-.2cm]

[3] \quad Dold, A. and Thom, R., Quasifaserungen und Unendliche Symmetrische Produkte, 
\textit{Ann. of Math.} \textbf{67} 

\hspace{.8cm}(1958), 239-281.
\\[-.2cm]

[4] \quad Dold, A. and Lashof, R., Principal Quasifibrations and Fibre Homotopy Equivalences, 
\textit{Illinois J. Math.} 

\hspace{.8cm}\textbf{3} (1959), 285-305.
\\[-.2cm]

[5] \quad Earle, C. and Eells, J., A Fibre Bundle Description of Teichm\"uller Theory, 
\textit{J. Differential Geom.} \textbf{3} 

\hspace{.8cm}(1969), 19-43.
\\[-.2cm]

[6] \quad Eells, J., Fiber Bundles, In: \textit{Global Analysis and its Applications}, vol. I, P. Harpe, et al. (ed.),

\hspace{.8cm}International Atomic Energy Agency (1974), 53-82.
\\[-.2cm]

[7] \quad Fuchs, M., A Modified Dold-Lashof Construction that does Classify H-Principal Fibrations, 
\textit{Math.} 

\hspace{.8cm}\textit{Ann.} \textbf{192} (1971), 328-340.
\\[-.2cm]

[8] \quad Holm, P., The Microbundle Representation Theorem, 
\textit{Acta Math.} \textbf{117} (1967), 191-213.
\\[-.2cm]

[9] \quad Holm, P., Microbundles and S-Duality, 
\textit{Acta Math.} \textbf{118} (1967), 271-296.
\\[-.2cm]

[10] \quad Hurewicz, W., On the Concept of Fiber Space, 
\textit{Proc. Nat. Acad. Sci.} \textbf{41} (1955), 956-961.
\\[-.2cm]

[11] \quad Lashof, R. and Rothenberg, M., G-Smoothing Theory, 
\textit{Proc. Symp. Pure Math.} \textbf{32} (1978), 211-266.
\\[-.2cm]

[12] \quad Martino, J., Classifying Spaces and their Maps, 
\textit{Contemp. Math.} \textbf{188} (1995), 153-190.
\\[-.2cm]

[13] \quad May, J., Classifying Spaces and Fibrations, 
\textit{Memoirs Amer. Math. Soc.} \textbf{155} (1975). 1-98.
\\[-.2cm]

[14] \quad Milnor, J., Construction of Universal Bundles I and II, 
\textit{Ann. of Math.} \textbf{63} (1956), 272-284 and 430-436.
\\[-.2cm]

[15] \quad Milnor, J., Microbundles, 
\textit{Topology} \textbf{3} (1964), 53-80.
\\[-.2cm]

[16] \quad Morgan, C. and Piccinini, R., Fibrations, 
\textit{Expo. Math.} \textbf{4} (1986), 217-242.
\\[-.2cm]

[17] \quad Nomura, Y., On Mapping Sequences, 
\textit{Nagoya Math. J.} \textbf{17} (1960), 111-145.
\\[-.2cm]

[18] \quad Pacati, C., Pave\u sic, P., and Piccinini, R., The Dold-Lashof Construction Revisited, 
\textit{Univ. Ljubljana} 

\hspace{.95cm}\textbf{32} (1994), 1-32.
\\[-.2cm]

[19] \quad Schwarz, G., Lifting Smooth Homotopies of Orbit Spaces, 
\textit{Publ. Math. I.H.E.S.} \textbf{51} (1980), 37-135.
\\[-.2cm]

[20] \quad Serre, J-P., Homologie Singuli\`ere des Espaces Fibr\'es, 
\textit{Ann. of Math.} \textbf{54} (1951), 425-505.
\\[-.2cm]

[21] \quad Stasheff, J., H-Spaces and Classifying Spaces: Foundations and Recent Developments, 
\textit{Proc. Symp.} 

\hspace{.8cm}\textit{Pure Math.} \textbf{22} (1971), 247-272.
\\[-.2cm]

[22] \quad Steenrod, N., Milgram's Classifying Space of a Topological Group, 
\textit{Topology} \textbf{7} (1968), 349-368.

\setlength\parindent{2em}

\endgroup

\chapter{
$\boldsymbol{\S}$\textbf{5}.\quadx  VERTEX SCHEMES AND CW COMPLEXES}
\setlength\parindent{2em}
\setcounter{proposition}{0}
%%----------------------------------------------------------------------------------------------01
$\text{ }$\\[-1.9cm]

Vertex schemes and CW complexes pervade algebraic topology.  
What follows is an account of their basic properties.  
All the relevant facts will be stated with precision but I shall only provide proofs for those that are not readily available in the standard treatments.

A \un{vertex scheme} $K$ is a pair $(V,\Sigma)$ consisting of a set $V = \{v\}$ and a subset $\Sigma = \{\sigma\} \subset 2^V$ subject to: (1) $\forall \ \sigma: \sigma \neq \emptyset$ $\&$ $\#(\sigma) < \omega$; 
(2) $\forall \ \sigma: \emptyset \neq \tau \subset \sigma \implies \tau \in \Sigma$; (3) $\forall \ v: \{v\} \in \Sigma$.  
The elements $v$ of $V$ are called the \un{vertexes} of $K$ and the elements $\sigma$ of $\Sigma$ are called the \un{simplexes} of $K$, the nonempty $\tau \subset \sigma$ being termed the \un{faces} of $\sigma$.  
A \un{vertex map} $f: K_1 = (V_1,\Sigma_1) \ra K_2 = (V_2,\Sigma_2)$ is a function $f:V_1 \ra V_2$ such that $\forall \ \sigma_1 \in \Sigma_1$, $f(\sigma_1) \in \Sigma_2$.  \textbf{VSCH} is the category whose objects are the vertex schemes and whose morphisms are the vertex maps.\\
\index{vertex scheme}
\index{vertexes}
\index{simplexes}
\index{faces}
\index{vertex map}
\index{\textbf{VSCH}}

\begingroup%%----------------------------------->>
\fontsize{9pt}{11pt}\selectfont
\textbf{\small EXAMPLE} \ 
Let $X$ be a set; let $\sS = \{S\}$ be a collection of subsets of $X$ $-$then the 
\un{nerve}
\index{nerve (set)} 
of $\sS$, written $N(\sS)$, is the vertex scheme whose vertexes are the nonempty elements of $\sS$ and whose simplexes are the nonempty finite subsets of $\sS$ with nonempty intersection.\\
\endgroup%%------------------------------------<<

Let $K = (V,\Sigma)$ be a vertex scheme.  If $\#(\Sigma) < \omega$ $(\leq \omega)$, then $K$ is said to be 
\un{finite}
\index{vertex scheme (finite)} 
(\un{countable}).
\index{vertex scheme (countable)}  
If $\forall \ v$, $\#\{\sigma: v \in \sigma\} < \omega$, then $K$ is said to be 
\un{locally finite}.
\index{vertex scheme (locally finite)}  
A 
\un{subscheme}
\index{subscheme} 
of $K$ is a vertex scheme $K^\prime = (V^\prime,\Sigma^\prime)$ such that 
$
\begin{cases}
\ V^\prime \subset V\\[-.1cm]
\ \Sigma^\prime \subset \Sigma
\end{cases}
\hspace{-.26cm}. \ 
$
An 
\un{$n$-simplex}
\index{n-simplex} 
is a simplex of cardinality $n+1$ $(n \geq 0)$.  The 
\un{$n$-skeleton}
\index{n-skeleton} 
of $K$ is the subscheme $K^{(n)} = (V^{(n)},\Sigma^{(n)})$ of $K$ defined by putting $V^{(n)} = V$ and letting $\Sigma^{(n)} \subset \Sigma$ be the set of $m$-simplexes of $K$ with $m \leq n$.  
The 
\un{combinatorial dimension}
\index{combinatorial dimension (vertex scheme)} 
of $K$, written $\dim K$, is $-1$ if $K$ is empty, otherwise is $n$ if $K$ contains an $n$-simplex but no $(n+1)$-simplex and is $\infty$ if $K$ contains $n$-simplexes for all $n \geq 0$.  
If $K$ is finite, then $\dim K$ is finite.  The converse is trivially false.\\

\begingroup%%----------------------------------->>
\fontsize{9pt}{11pt}\selectfont
\textbf{\small EXAMPLE} \ 
In the plane, take $V = \{(0,0)\} \cup \{(1,1/n): n \geq 1\}$.  Let $K = (V,\Sigma)$ 
be any vertex scheme having for its 1-simplexes the sets $\sigma_n = \{(0,0),(1,1/n)\}$ $(n \geq 1)$ $-$then $K$ is not locally finite.\\
\endgroup%%------------------------------------<<

Given a vertex scheme $K = (V,\Sigma)$, let $\abs{K}$ be the set of all functions 
$\phi:V \ra [0,1]$ such that $\phi^{-1}(]0,1]) \in \Sigma$ \  $\&$ \ $\sum\limits_v \phi(v) = 1$.  
Assign to each $\sigma$ the sets
%$
%\begin{equation*}
%\begin{cases}
%\ \langle\sigma\rangle \ = \{\phi \in \abs{K}: \phi^{-1}(]0,1]) = \sigma\}\\
%\ \abs{\sigma} \ \  = \{\phi \in \abs{K}: \phi^{-1}(]0,1]) \subset \sigma\}
%\end{cases}
%\end{equation*}
%.\ 
%$
%\[
%\begin{cases}
%\ \langle\sigma\rangle \ = \{\phi \in \abs{K}: \phi^{-1}(]0,1])  = \sigma\}\\
%\ \abs{\sigma} \ \  = \{\phi \in \abs{K}: \phi^{-1}(]0,1])  \subset \sigma\}
%\end{cases}
%\hspace{-.26cm}.
%\]
$
\begin{cases}
\ \langle\sigma\rangle \ = \{\phi \in \abs{K}: \phi^{-1}(]0,1]) \\[-.1cm]
\ \abs{\sigma} \ \  = \{\phi \in \abs{K}: \phi^{-1}(]0,1]) 
\end{cases}
$
$
\begin{aligned}
 = \sigma\}\\
 \subset \sigma\}
\end{aligned}
. \ 
$
So, $\forall \ \sigma: \langle\sigma\rangle \  \subset \abs{\sigma}$ and 
$\abs{K} = \ds\bigcup\limits_\sigma \langle\sigma\rangle$, a disjoint union.  
Traditionally, there are two ways to topologize $\abs{K}$.

\indent \indent (WT) \   If $\sigma$ is an $n$-simplex, then $\abs{\sigma}$ can be viewed as a compact Hausdorff space: $\abs{\sigma} \leftrightarrow \Delta^n$.  
This said, the 
\un{Whitehead topology}
\index{Whitehead topology} 
on $\abs{K}$ is the final topology determined 
%%----------------------------------------------------------------------------------------------02
by the inclusions $\abs{\sigma} \ra \abs{K}$.  $\abs{K}$ is a perfectly normal paracompact Hausdorff space.  
Moreover, $\abs{K}$ is 
$
\begin{cases}
\ \text{compact} \\[-.1cm]
\ \text{locally compact}
\end{cases}
$
\hspace{-.26cm} iff $K$ is 
$
\begin{cases}
\ \text{finite} \\[-.1cm]
\ \text{locally finite}
\end{cases}
$
\hspace{-.26cm}.\\

\indent \indent (BT) \  There is a map
$
\begin{cases}
\ V \ra [0,1]^{\abs{K}} \\[-.1cm]
\ v \mapsto b_v: b_v(\phi) = \phi(v)
\end{cases}
\hspace{-.26cm}. \ 
$
The $b_v$ are called the 
\un{barycentric} \un{coordinates}, 
\index{barycentric coordinates}
the initial topology on $\abs{K}$ determined by them being the 
\un{barycentric topology}, 
\index{barycentric topology} 
a topology that is actually metrizable: $d(\phi,\psi) = \sum\limits_v \abs{b_v(\phi) - b_v(\psi)}$.

To keep things straight, denote by $\abs{K}_b$ the set $\abs{K}$ equipped with the barycentric topology $-$then the identity map $i: \abs{K} \ra \abs{K}_b$ is continuous, thus the Whitehead topology is finer than the barycentric topology.  The two agree iff $K$ is locally finite.

[Note: \  A vertex map $f:K_1 = (V_1,\Sigma_1) \ra K_2 = (V_2,\Sigma_2)$ induces a map $\abs{f}:$
$
\begin{cases}
\ \abs{K_1} \ra \abs{K_2} \\
\ \phi_1 \mapsto \phi_2
\end{cases}
\hspace{-.26cm}, \ 
$
where $\phi_2(v_2) = \ds\sum\limits_{f(v_1) = v_2} \phi_1(v_1)$.  
Topologically, $\abs{f}$ is continuous in either the Whitehead topology or the barycentric topology.  
Consequently, there are two functors from \textbf{VSCH} to \bTOP, connected by the obvious natural transformation.]\\

\begingroup%%----------------------------------->>
\fontsize{9pt}{11pt}\selectfont
\textbf{\small EXAMPLE} \ 
Let $E$ be a vector space over $\R$.  Let $V$ be a basis for \mE; let $\Sigma$ be the set of nonempty finite subsets of \mV.  
Call $K(E)$ the associated vertex scheme.  
Equip $E$ with the finite topology $-$then $\abs{K(E)}$ can be identified with the convex hull of $V$ in \mE.  
But $\abs{K(E)}$ and $\abs{K(E)}_b$ are homeomorphic iff $E$ is finite dimensional.
\vspi
[Note: \  Let $K = (V,\Sigma)$ be a vertex scheme.  
Take for $E$ the free $\R$-module on \mV, equipped with the finite topology $-$then $\abs{K}$ can be embedded in $\abs{K(E)}$.]\\
\endgroup%%------------------------------------<<

\begin{proposition} \ %01
The identity map $i: \abs{K} \ra \abs{K}_b$ is a homotopy equivalence.
\end{proposition}

[The collection $\{b_v^{-1}(]0,1])\}$ is an open covering of $\abs{K}_b$, hence has a precise neighborhood finite open refinement $\{U_v\}$.  Choose a partition of unity $\{\kappa_v\}$ on $\abs{K}_b$ subordinate to $\{U_v\}$.  Let $j:\abs{K}_b \ra \abs{K}$ be the map that sends $\psi$ to the function
$
\begin{cases}
\ V \ra [0,1]\\[-.1cm]
\ v \mapsto \kappa_v(\psi)
\end{cases}
$
\hspace{-.26cm}.  
Consider the homotopies
$
\begin{cases}
\ H:I\abs{K} \ra \abs{K}\\[-.1cm]
\ G:I\abs{K}_b \ra \abs{K}_b
\end{cases}
$
defined by
$
\begin{cases}
\ H(\phi,t) = t\phi + (1 - t)j \circx i(\phi)\\[-.1cm]
\ G(\psi,t) = t\psi + (1 - t)i \circx j(\psi)
\end{cases}
\hspace{-.26cm}.]
$
\\
\vspace{0.25cm}

\begingroup%%----------------------------------->>
\fontsize{9pt}{11pt}\selectfont
Let $X$ be a topological space $-$then two continuous functions
$
\begin{cases}
\ f:X \ra \abs{K}\\
\ g:X \ra \abs{K}
\end{cases}
$
are said to be 
\un{contiguous}
\index{contiguous (functions)} 
if $\forall \ x \in X$ $\exists \ \sigma \in \Sigma:$ $\{f(x),g(x)\} \subset \abs{\sigma}$.\\
\endgroup%%------------------------------------<<

\begingroup%%----------------------------------->>
\fontsize{9pt}{11pt}\selectfont
\textbf{\small FACT} \ 
Suppose that 
$
\begin{cases}
\ f:X \ra \abs{K}\\
\ g:X \ra \abs{K}
\end{cases}
$
are contiguous $-$then $f \simeq g$.
\vspi
[Define a homotopy $H:IX \ra \abs{K}_b$ between $i \circx f$ and $i \circx g$ by writing $b_v(H(x,t)) = (1 - t)b_v(f(x)) + tb_v(g(x))$ and apply Proposition 1.]\\
\endgroup%%------------------------------------<<

%%----------------------------------------------------------------------------------------------03
\label{6.14}
\label{19.29}
\begingroup%%----------------------------------->>
\fontsize{9pt}{11pt}\selectfont
\textbf{\small EXAMPLE} \ 
Let $X$ be a topological space; let $\sU = \{U\}$ be a numerable open covering of $X$ $-$then a 
\un{$\sU$-map}
\index{U-map (numerable open covering)} %\index{$\sU$-map (numerable open covering)}
is a continuous function $f:X \ra \abs{N(\sU)}$ such that $\forall \ U \in \sU:$ $(b_U \circx f)^{-1}(]0,1]) \subset U$.  
Every partition of unity on \mX subordinate to $\sU$ defines a $\sU$-map and any two $\sU$-maps are contiguous, hence homotopic.\\
\endgroup%%------------------------------------<<

\begingroup%%----------------------------------->>
\fontsize{9pt}{11pt}\selectfont
\textbf{\small FACT} \ 
Let $X$ be a topological space.  Suppose that
$
\begin{cases}
\ f:X \ra \abs{K}\\
\ g:X \ra \abs{K}
\end{cases}
$
are two continuous functions such that $\forall \ x \in X$ $\exists \ v \in V:$ $\{f(x),g(x)\} \subset b_v^{-1}(]0,1])$ $-$then $f \simeq g$.\\
\endgroup%%------------------------------------<<

\index{Theorem Adjunction}
\textbf{\small ADJUNCTION THEOREM} \ 
Let $K$ and $L^\prime$ be vertex schemes.  Let $K^\prime$ be a subscheme of $K$ and let $f:K^\prime \ra L^\prime$ be a vertex map $-$then there exists a vertex scheme $L$ containing $L^\prime$ as a subscheme and a homeomorphism $\abs{K} \sqcup_{\abs{f}} \abs{L^\prime} \ra \abs{L}$ whose restriction to $\abs{L^\prime} $ is the identity map.\\

A topological space $X$ is said to be a 
\un{polyhedron}
\index{polyhedron} 
if there exists a vertex scheme $K$ and a homeomorphism $f: \abs{K} \ra X$ ($\abs{K}$ in the Whitehead topology).  
The ordered pair $(K,f)$ is called a 
\un{triangulation}
\index{triangulation} 
of $X$.  
Put $f_v = b_v \circx f^{-1}$ $-$then the collection $\sT_K = \{f_v^{-1}(]0,1])\}$ is a numerable open covering of $X$ and Whitehead's\footnote[2]{\textit{Proc. London Math. Soc.} \textbf{45} (1939), 243-327.} 
``Theorem 35'' says:  
For any open covering $\sU$ of $X$, there exists a triangulation $(K,f)$ of $X$ such that $\sT_K$ refines $\sU$ .

\label{14.106}
Every polyhedron is a perfectly normal paracompact Hausdorff space.  
A polyhedron is metrizable iff it is locally compact.  
Every open subset of a polyhedron is a polyhedron.\\

\begingroup%%----------------------------------->>
\fontsize{9pt}{11pt}\selectfont
Let $X$ be a topological space $-$then a closure preserving closed covering $\sA = \{A_j: j \in J\}$ of $X$ is said to be 
\un{absolute}
\index{absolute (closure preserving closed covering)} 
if for every subset $I \subset J$, the subspace 
$X_I = \ds\bigcup\limits_i A_i$ has the final topology with respect to the inclusions $A_i \ra X_I$.  
Example:  Every neighborhood finite closed covering of $X$ is absolute.
\vspi
[Note: \  Let $K$ be a vertex scheme $-$then 
$\{\abs{\sigma}\}$ is an absolute closure preserving closed covering of $\abs{K}$ but, in general, 
is only a closure preserving closed covering of $\abs{K}_b$.]\\
\endgroup%%------------------------------------<<

\begingroup%%----------------------------------->>
\fontsize{9pt}{11pt}\selectfont
\textbf{\small EXAMPLE} \ 
Take $X = [0,1]$, put $X_1 = [0,1]$, $X_n = \{0\} \cup [1/n,1]$ $(n > 1)$ $-$then 
$\{X_n\}$ is a closure preserving closed covering of $X$ but $\{X_n\}$ is not absolute since 
$X = \ds\bigcup\limits_{n>1} X_n$ does not have the final topology with respect to the inclusions $X_n \ra X$ $(n > 1)$.\\
\endgroup%%------------------------------------<<

\begingroup%%----------------------------------->>
\fontsize{9pt}{11pt}\selectfont
\textbf{\small LEMMA} \ 
Let  $\sA = \{A_j: j \in J\}$ be an absolute closure preserving closed covering of $X$ $-$then for any compact Hausdorff space $K$, $\sA  \times K = \{A_j \times X; j \in J\}$ is an absolute closure preserving closed covering of $X \times K$.\\
\endgroup%%------------------------------------<<

%%----------------------------------------------------------------------------------------------04
\label{1.17} 
\label{4.69}
\begingroup%%----------------------------------->>
\fontsize{9pt}{11pt}\selectfont
\textbf{\small FACT} \ 
If $X$ is a topological space and if $\sA = \{A_j: j \in J\}$ be an absolute closure preserving closed covering of $X$ 
such that each $A_j$ is a normal (normal and countably paracompact, perfectly normal, collectionwise normal, paracompact) Hausdorff space, then $X$ is a normal (normal and countably paracompact, perfectly normal, collectionwise normal, paracompact) Hausdorff space.
\vspi
[In every case, $X$ is $\tT_1$.  And: $\tT_1 + $ normal $\implies$ Hausdorff.
\\
\indent\indent (Normal) \ 
Let $A$ be a closed subset of $X$, take an $f \in C(A,[0,1])$ and let $\sF$ be the set of continuous functions $F$ that are extensions of $f$ and have domains of the form $A \cup X_I$, where $X_I = \ds\bigcup\limits_i A_i$ $(I \subset J)$.  
Order $\sF$ by writing $F^\prime \leq F^{\prime\prime}$ iff $F^{\prime\prime}$ is an extension of $F^\prime$.  
Every chain in $\sF$ has an upper bound, so by Zorn, $\sF$ has a maximal element $F_0$.  
But the domain of $F_0$ is necessarily all of $X$ and $\restr{F_0}{A} = f$.
\\
\indent\indent (Normal and Countably Paracompact) \ 
First recall that a normal Hausdorff space is countably paracompact iff its product with $[0,1]$ is normal.  Since $\sA \times [0,1] = \{A_j \times [0,1]: j \in J\}$ is an absolute closure preserving closed covering of $X \times [0,1]$, it follows that $X \times [0,1]$ is normal, thus $X$ is countably paracompact.
\\
\label{6.32}
\indent\indent (Perfectly Normal) \ 
Fix a closed subset $A$ of $X$.  To prove that $A$ is a zero set in $X$, equip $J$ with a well ordering $<$.  
Given $j \in J$, put $X(j) = \ds\bigcup\limits_{i \leq j} A_i$.  Inductively construct continuous functions $f_j:X(j) \ra [0,1]$ such that 
$\restr{f_{j^{\prime\prime}}}{X(j^\prime)} =$ 
$f_{j^{\prime}}$ if $j^\prime < j^{\prime\prime}$ and $Z(f_j) = A \cap X(j)$.
\\
\indent\indent (Collectionwise Normal) \ 
Let $A$ be a closed subset of $X$, $E$ any Banach space $-$then it suffices to show that every $f \in C(A,E)$ admits an extension $F \in C(X,E)$ (cf. p. \pageref{5.0}).  This can be done by imitating the argument used to establish normality.
\\
\indent\indent (Paracompact) \  
Tamano's theorem says that a normal Hausdorff space $X$ is paracompact iff $X \times \beta X$ is normal, which enables one to proceed as in the proof of countable paracompactness.]\\
\endgroup%%------------------------------------<<

\begingroup%%----------------------------------->>
\fontsize{9pt}{11pt}\selectfont
\textbf{\small EXAMPLE} \ 
The ordinal space $[0,\Omega[$ is not paracompact but $\{[0,\alpha]: \alpha < \Omega\}$ is a covering of $[0,\Omega[$ by compact Hausdorff spaces and $[0,\Omega[$ has the final topology with respect to the inclusions $[0,\alpha] \ra [0,\Omega[$.\\
\endgroup%%------------------------------------<<

\begingroup%%----------------------------------->>
\fontsize{9pt}{11pt}\selectfont
\textbf{\small FACT} \ 
Let $X$ be a topological space; let $\sA = \{A_j: j \in J\}$ be an absolute closure preserving closed covering of $X$. 
Suppose that each $A_j$ can be embedded as a closed subspace of a polyhedron $-$then $X$ 
can be embedded as a closed subspace of a polyhedron.
\vspi
[For every $j$ there is a vertex scheme $K_j$, a vector space $E_j$ over $\R$, and a closed embedding $f_j:A_j \ra \abs{K_j}$ $(\subset E_j)$.  
Write $E$ for the direct sum of the $E_j$ and give $E$ the finite topology.  
Let $E_I$ stand for the direct sum of the $E_i$ $(i \in I)$ and put $K_I = K(E_I)$ $-$then $\abs{K_I} \subset \abs{K(E)}$.  Here, as above, $I$ is a subset of $J$.  
Consider the set $\sP$ of all pairs $(I,f_I)$ , where $f_I:X_I \ra \abs{K_I}$ is a closed embedding.  
Order $\sP$ by stipulating that 
$(I^\prime, f_{I^\prime}) \leq (I^{\prime\prime}, f_{I^{\prime\prime}})$ iff $I^\prime \subset I^{\prime\prime}$ and 
(1) 
$\restr{f_{I^{\prime\prime}}}{X_{I^\prime}} = f_{I^\prime}$ $\&$ 
(2) 
$f_{I^{\prime\prime}}(X_{I^{\prime\prime}} - X_{I^\prime}) \cap \abs{K_{I^\prime}} =$ $\emptyset$.  
Every chain in $\sP$ has an upper bound, so by Zorn, $\sP$ has a maximal element 
$(I_0,f_{I_0})$.  Verify that $X_{I_0} = X$.]\\
\endgroup%%------------------------------------<<

\begingroup%%----------------------------------->>
\fontsize{9pt}{11pt}\selectfont
Application:  Let $X$ be a paracompact Hausdorff space.  
Suppose that $X$ admits a covering $\sU$ by open
%%----------------------------------------------------------------------------------------------05
sets $U$, each of which is homeomorphic to a closed subspace of a polyhedron $-$then $X$ is homeomorphic to a closed subspace of a polyhedron.\\
\endgroup%%------------------------------------<<

\label{5.5}
\begingroup%%----------------------------------->>
\fontsize{9pt}{11pt}\selectfont
The embedding theorem of dimension theory implies that every second countable compact Hausdorff space of finite topological dimension can be embedded in some euclidean space (cf. p. \pageref{5.1}).  
It therefore follows that if a topological space $X$ has an absolute closure preserving closed covering made up of metrizable compacta of finite topological dimension, then $X$ can be embedded as a closed subspace of a polyhedron.  
This setup is realized, e.g., by the CW complexes (cf. p. \pageref{5.2}).\\
\endgroup%%------------------------------------<<

The product $X \times Y$ of polyhedrons $X$ and $Y$ need not be a polyhedron 
(cf. p. 
\pageref{5.3}), 
although this will be the case if one of the factors is locally compact.\\

\begingroup%%----------------------------------->>
\fontsize{9pt}{11pt}\selectfont
\label{5.0l}
\textbf{\small FACT} \ 
Let $X$ and $Y$ be polyhedrons $-$then $X \times Y$ has the homotopy type of a polyhedron.
\vspi
[Consider a product $\abs{K} \times \abs{L}$.  Since 
$
\begin{cases}
\ \abs{K}  \ \& \  \abs{K}_b\\
\ \abs{L}  \  \&  \ \abs{L}_b
\end{cases}
$
have the same homotopy type, it need only be shown that  $\abs{K}_b \times \abs{L}_b$ has the homotopy type of a polyhedron.  Let
$
\begin{cases}
\ \sU\\
\ \sV
\end{cases}
$
be the cozero set covering of 
$
\begin{cases}
\ \abs{K}_b\\
\ \abs{L}_b
\end{cases}
$
associated with the barycentric coordinates $-$then
$
\begin{cases}
\ K\\
\ L
\end{cases}
$
can be identified with the corresponding nerve
$
\begin{cases}
\ N(\sU)\\
\ N(\sV)
\end{cases}
.\ 
$
Put $\sU\times \sV = \{U \times V: U \in \sU, V \in \sV\}$.  
Claim:  There is a homotopy equivalence $\abs{N(\sU \times \sV)}_b \ra \abs{N(\sU)}_b \times  \abs{N(\sV)}_b$.  
Indeed, the projections
$
\begin{cases}
\ \sU \times \sV \ra \sU \quadx (U \times V \ra U)\\
\ \sU \times \sV \ra \sV \quadx (U \times V \ra V)
\end{cases}
$
define vertex maps
$
\begin{cases}
\ p_\sU: N(\sU \times \sV) \ra N(\sU)\\
\ p_\sV: N(\sU \times \sV) \ra N(\sV)
\end{cases}
, \ 
$
from which $p:\abs{N(\sU \times \sV)}_b \ra \abs{N(\sU)}_b \times  \abs{N(\sV)}_b$, where $p = \abs{p_\sU} \times \abs{p_\sV}$.  A homotopy inverse $q:\abs{N(\sU)}_b \times  \abs{N(\sV)}_b \ra \abs{N(\sU \times \sV)}_b$ to $p$ is given in terms of barycentric coordinates by $b_{U \times V}(q(\phi,\psi)) = b_U(\phi)b_V(\psi)$.]\\
\endgroup%%------------------------------------<<

Let $X$ be a topological space; let $A$ be a closed subspace of $X$ $-$then $X$ is said to be obtained from $A$ by
\un{attaching $n$-cells} 
\index{attaching $n$-cells} 
if there exists an indexed collection of continuous functions $f_i:\bS^{n-1} \ra A$ 
such that $X$ is homeomorphic to the adjunction space $\bigl(\coprod\limits_i\bD^{n}\bigr) \sqcup_f A$ $(f = \coprod\limits_i f_i)$.  
When this is so, $X - A$ is homeomorphic to $\coprod\limits_i \left(\bD^{n} - \bS^{n-1}\right) = \coprod\limits_i \bB^{n} $, 
a decomposition that displays its path components as a collection of $n$-cells.\\

\begingroup%%----------------------------------->>
\fontsize{9pt}{11pt}\selectfont
\textbf{\small EXAMPLE} \ 
Put $s_n = (1, 0, \ldots , 0) \in \R^{n+1}$ $(n \geq 1)$.  
Let $I$ be a set indexing a collection of copies of the pointed space $(\bS^n,s_n)$ 
$-$then the wedge $\ds\bigvee\limits_I \bS^n$ is a pointed space with basepoint $*$.  
Since the quotient $\bD^{n}/\bS^{n-1}$ can be identified with $\bS^n$, $\ds\bigvee\limits_I \bS^n$ is obtained from $*$ by attaching $n$-cells.\\
\endgroup%%------------------------------------<<

%%----------------------------------------------------------------------------------------------06
Let $X$ be a topological space $-$then a 
\un{CW structure}
\index{CW structure} 
on $X$ is a sequence $X^{(0)}, X^{(1)},\ldots$
of closed subspaces $X^{(n)}:$
$
\begin{cases}
\ X = \bigcup\limits_0^\infty X^{(n)}\\[-.1cm]
\ X^{(n)} \subset X^{(n+1)}
\end{cases}
$
and subject to:\\
\indent \indent (CW$_1$) \quadx $X^{(0)}$ is discrete.\\
\indent \indent (CW$_2$) \quadx $X^{(n)}$ is obtained from $X^{(n-1)}$ by attaching $n$-cells $(n > 0)$.\\
\indent \indent (CW$_3$) \quadx $X$ has the final topology determined by the inclusions $X^{(n)} \ra X$.

A 
\un{CW complex}
\index{CW complex} 
is a topological space $X$ equipped with a CW structure.  
Just as a polyhedron may have more than one triangulation, a CW complex may have more than one CW structure.  
Every CW complex is a perfectly normal paracompact Hausdorff space.

[Note: \   Let $K$ be a vertex scheme.  
Consider $\abs{K}$ (Whitehead topology) $-$then $\abs{K^{(0)}}$ is discrete and 
$\abs{K^{(n)}}$ is obtained from $\abs{K^{(n-1)}}$ by attaching $n$-cells 
$(n > 0): \abs{\sigma} - \langle\sigma\rangle \ra \abs{K^{(n-1)}}$, $\sigma$ an $n$-simplex.  
Since $\abs{K}$ has the final topology determined by the inclusions 
$\abs{K^{(n)}} \ra \abs{K}$, it follows that the sequence $\{\abs{K^{(n)}}\}$ is a CW structure on $\abs{K}$.]

\bCW 
\index{\bCW} 
is the full subcategory of \bTOP whose objects are the CW complexes and 
\bHCW
\index{\bHCW}  is the associated homotopy category.\\

\begingroup%%----------------------------------->>
\fontsize{9pt}{11pt}\selectfont
\textbf{\small EXAMPLE} \ 
Equip $\R^\infty$ with the finite topology. Let $\bS^\infty = \ds\bigcup\limits_0^\infty \bS^n$ and give it the induced topology or, 
what amounts to the same, the final topology determined by the inclusions $\bS^n \ra \bS^\infty$.  
The sequence $\{\bS^n\}$ is a CW structure on $\bS^\infty$.  
Indeed, $\bS^n$ is obtained from $\bS^{n-1}$ by attaching two $n$-cells $(n > 0)$ (seal the upper and lower hemispheres at the equator).  
On the other hand, $\R^n$ is not obtained from $\R^{n-1}$ by attaching $n$-cells.  
Therefore the sequence $\{\R^n\}$ is not a CW structure on $\R^\infty$.  But $\R^\infty$ is obviously a polyhedron.  
A less apparent aspect is this.  Put $s_\infty = (1, 0, \ldots)$ $-$then it can be shown that $\bS^\infty$ and $\bS^\infty - \{s_\infty\}$ are homeomorphic.  
Since stereographic projection from $s_\infty$ defines a homeomorphism $\bS^\infty - \{s_\infty\} \ra \R^\infty$, 
the conclusion is that $\bS^\infty$ and $\R^\infty$ are actually homeomorphic.
\vspi
[Note: \  The sequence $\{\bD^n\}$ is not a CW structure for $\bD^\infty = \ds\bigcup\limits_0^\infty \bD^n$.  
However, $\bD^n \cup \bS^n$ can be obtained from $\bD^{n-1} \cup \hsx \bS^{n-1}$ by attaching $n$-cells $(n > 0)$, so the sequence $\{\bD^n \cup \hsx \bS^n\}$ is a CW structure for $\bD^\infty$.]\\
\endgroup%%------------------------------------<<

Let $X$ be a CW complex with CW structure $\{X^{(n)}\}:X^{(n)}$ is the 
\un{$n$-skeleton}
\index{CW complex $n$-skeleton} 
of $X$. The inclusion $X^{(n)} \ra X$ is a closed cofibration (cf. p. \pageref{5.0a}) and $\forall \ n \geq 1$, the pair $(X,X^{(n)})$ is $n$-connected.  
Put $\sE_0 = X^{(0)}$ and denote by $\sE_n$ the set of path components of 
$X^{(n)} - X^{(n-1)}$ $(n > 0)$.  
Let $\sE = \bigcup\limits_0^\infty \sE_n$ $-$then an element  $e$ of $\sE$ is said to be a 
\un{cell}
\index{CW complex cell}in $X$, $e$ being termed an 
\un{$n$-cell}
\index{CW complex $n$-cell}if $e \in \sE_n$.  
Set theoretically, $X$ is the disjoint union of its cells.  
On the basis of the definitions, for every $e \in \sE_n$, there exists a continuous function 
$\Phi_e:\bD^n \ra$ 
$e \cup X^{(n-1)}$, the 
\un{characteristic map}
\index{CW complex characteristic map} 
of $e$, such that $\restr{\Phi_e}{\bB^n}$ is an embedding and 
(i) $\Phi_e(\bB^n) = e$; 
(ii) $\Phi_e(\bS^{n-1}) \subset X^{(n-1)}$; 
(iii) $\Phi_e(\bD^n) = \bar{e}$.  $X$ has the final topology determined by the $\Phi_e$.
%%----------------------------------------------------------------------------------------------07
A subspace $A \subset X$ is called a 
\un{subcomplex}
\index{CW complex subcomplex} 
if there exists a subset $\sE_A \subset \sE: A = \bigcup \sE_A$ $\&$ $\forall \ e \in \sE_A \cap \sE_n$, $\Phi_e(\bD^n) \subset A$.  
A subcomplex $A$ of $X$ is itself a CW complex with CW structure $\{A^{(n)} = A \cap X^{(n)}\}$.  
The inclusion $A \ra X$ is a closed cofibration and for every open $U \supset A$ 
there exists an open $V \supset A$ with $V \subset U$ such that \mA is a strong deformation retract of \mV. 
If $\sE^\prime \subset \sE$, then $\bigcup \sE^\prime$ is a subcomplex iff  $\bigcup \sE^\prime$ is closed.  
Arbitrary unions and intersections of subcomplexes are subcomplexes.  
In general, the $\bar{e}$ are not subcomplexes, although this will be the case if all characteristic maps are embeddings.  
The 
\un{combinatorial dimension}
\index{CW complex combinatorial dimension}
of $X$, written $\dim X$, is $-1$ is $X$ is empty, otherwise is the smallest value of $n$ such that $X = X^{(n)}$ 
(or $\infty$ if there is no such $n$).  
It is a fact that $\dim X$ is equal to the topological dimension of $X$ (cf. p. \pageref{5.0b}), therefore is independent of the CW structure.\\

\begingroup%%----------------------------------->>
\fontsize{9pt}{11pt}\selectfont
Let $X$ be a CW complex $-$then the collection $\overline{\sE} = \{\bar{e}: e \in \sE\}$ is a closed covering of $X$ and $X$ has the final topology determined by the inclusions $\bar{e} \ra X$ but 
$\overline{\sE}$ need not be closure preserving.\\
\endgroup%%------------------------------------<<

\label{13.2}
\index{simplicial sets (example)}
\begingroup%%----------------------------------->>
\fontsize{9pt}{11pt}\selectfont
\textbf{\small EXAMPLE \ \un{(Simplicial Sets)}} \ 
Let $X$ be a simplicial set $-$then its geometric realization $\abs{X}$ is a CW complex with CW structure $\{|X^{(n)}|\}$.  In fact $|X^{(0)}|$ is discrete and, using the notation of p. \pageref{5.0c}, the commutative diagram
\begin{tikzcd}[sep=large]
{X_n^\# \cdot \dot{\Delta}[n]} \arrow{r} \arrow{d} &X^{(n-1)} \arrow{d}\\
{X_n^\# \cdot {\Delta}[n]} \arrow{r} &X^{(n)}
\end{tikzcd}
is a pushout square in \bSISET.  Since the geometric realization functor $\abs{?}$ is a left adjoint, it preserves colimits.  
Therefore the commutative diagram
\begin{tikzcd}[sep=large]
{X_n^\# \cdot \dot{\Delta}^n} \arrow{r} \arrow{d} &{|X^{(n-1)}|} \arrow{d}\\
{X_n^\# \cdot {\Delta}^n} \arrow{r} &{|X^{(n)}|}
\end{tikzcd}
is a pushout square in \bTOP, which means that $|X^{(n)}|$ is obtained from $|X^{(n-1)}|$ by attaching 
$n$-cells $(n > 0)$.  Moreover, 
$X = \colimx X^{(n)}$ $\implies$ 
$\abs{X} = \colimx |X^{(n)}|$ , 
so $\abs{X}$ has the final topology determined by the inclusions $|X^{(n)}| \ra \abs{X}$.  
Denoting now by $G$ the identity component of the homeomorphism group of $[0,1]$, there is a left action 
$G \times \abs{X} \ra \abs{X}$ and the orbits of \mG are the cells of $\abs{X}$.
\vspi
[Note: \ If \mY is a simplicial subset of $X$, then $\abs{Y}$ is a subcomplex of $\abs{X}$, thus the inclusion 
$\abs{Y} \ra$ 
$\abs{X}$ is a closed cofibration.]\\
\endgroup%%------------------------------------<<

\begingroup%%----------------------------------->>
\fontsize{9pt}{11pt}\selectfont
It is true but not obvious that if $X$ is a simplicial set, then $\abs{X}$ is actually a polyhedron (cf. p. \pageref{5.0d}).\\
\endgroup%%------------------------------------<<

A \un{CW pair} 
\index{CW pair} 
is a pair $(X,A)$, where $X$ is a CW complex and $A \subset X$ is a subcomplex.  
$\bCW^2$ is the full subcategory of $\bTOP^2$ whose objects are the CW pairs and $\bHCW^2$ is the associated homotopy category.\\

%%----------------------------------------------------------------------------------------------08
A \un{pointed CW complex} 
\index{pointed CW complex} 
is a pair $(X,x_0)$, where $X$ is a CW complex and $x_0 \in X^{(0)}$.  $\bCW_*$ is the full subcategory of $\bTOP_*$ whose objects are the pointed CW complexes and 
$\bHCW_*$ is the associated homotopy category.

[Note: \  If $(X,x_0)$ is a pointed CW complex, then $\forall \ q \geq 1$, 
$\pi_q(X,x_0) \approx$ 
$\colim \ \pi_q(X^{(n)},x_0)$.]\\

\begingroup%%----------------------------------->>
\fontsize{9pt}{11pt}\selectfont
Let $X$ be a CW complex $-$then $\forall x_0 \in X$, the inclusion $\{x_0\} \ra X$ is a cofibration 
(cf. p. \pageref{5.0e}), thus $(X,x_0)$ is wellpointed.  Of course, a given $x_0$ need not be in $X^{(0)}$ but there is always some CW structure on \mX having $x_0$ as a 0-cell.\\
\endgroup%%------------------------------------<<

Let $X$ be a topological space, $A \subset X$ a closed subspace$-$then a 
\un{relative CW} \un{structure}
\index{relative CW structure} \ 
on \ \ $(X,A)$\ \   is \  a  \ sequence \  $(X,A)^{(0)}, (X,A)^{(1)},\ldots$ \ 
of closed subspaces $(X,A)^{(n)}: \ $ 
$
\begin{cases}
\ X = \bigcup\limits_0^\infty (X,A)^{(n)}\\[-.1cm]
\ (X,A)^{(n)} \subset (X,A)^{(n+1)}
\end{cases}
$
and subject to:\\
\indent \indent (RCW$_1$) \ \  $(X,A)^{(0)}$ is obtained from \mA by attaching 0-cells.\\
\indent \indent (RCW$_2$) \ \  $(X,A)^{(n)}$ is obtained from $(X,A)^{(n-1)}$ by attaching $n$-cells $(n > 0)$.\\
\indent \indent (RCW$_3$) \ \ \  $X$ has the final topology determined by the inclusions $(X,A)^{(n)} \ra X$.

[Note: \ $(X,A)^{(0)}$ is the coproduct of $A$ and a discrete space, so when $A = \emptyset$ the definition reduces to that of a CW structure.]

A 
\un{relative CW complex}
\index{relative CW complex} 
is a topological space $X$ and a closed subspace $A$ equipped with a relative CW structure.

[Note: \   If $(X,A)$ is a relative CW complex, then the inclusion $A \ra X$ is a closed cofibration and $X/A$ is a CW complex.  
On the other hand, if $X$ is a CW complex and if $A \subset X$ is a subcomplex, then $(X,A)$ is a relative CW complex.]

Example:  Suppose that $(X,A)$ is a relative CW complex $-$then $(IX,IA)$ is a relative CW complex, where $(IX,IA) ^{(n)} = i_0(X,A)^{(n)} \cup (I(X,A)^{(n-1)} \cup IA) \cup i_1(X,A)^{(n)}$.

Let $(X,A)$ be a relative CW complex with relative CW structure $\{(X,A)^{(n)}\}$ : $(X,A)^{(n)}$ is the 
\un{$n-$skeleton}
\index{n-skeleton(relative CW complex)} 
of $X$ relative to $A$.  
The inclusion $(X,A)^{(n)} \ra X$ is a closed cofibration (cf. p. \pageref{5.0f}) and $\forall \ n \geq 1$, the pair $(X,(X,A)^{(n)})$ is $n$-connected.  The 
\un{relative combinatorial} \un{dimension}
\index{relative combinatorial dimension (relative CW complex)} 
of $(X,A)$, written $\dim (X,A)$, is $-1$ if $X$ is empty, otherwise is the smallest value of $n$ such that 
$X = (X,A)^{(n)}$ (or $\infty$ if there is no such $n$).  
Obviously, $\dim (X,A) = \dim (X/A)$ provided that $X$ is nonempty.\\

\textbf{\small LEMMA} \ 
Let $(X,A)$ be a relative CW complex $-$then for every compact subset $K \subset X$ there exists an index $n$ such that $K \subset (X,A)^{(n)}$.

[Consider the image of $K$ under the projection $X \ra X/A$, bearing in mind that $X/A$ is a CW complex.]\\

%%----------------------------------------------------------------------------------------------09
Application:  Let $(X,A,x_0)$ be a pointed pair.  Assume: $(X,A)$ is a relative CW complex 
$-$then $\forall \ q \geq 1$, $\pi_q(X,x_0) \approx \colimx \pi_q((X,A)^{(n)},x_0)$.\\

\index{Theorem, Hopf Extension}
\textbf{\small HOPF EXTENSION THEOREM} \quad 
Let $(X,A)$ be a relative CW complex with $\dim (X,A)$ $\leq$ $n +1$ $(n \geq 1)$.  Suppose that $f \in C(A,\bS^n)$ $-$then $\exists \ F \in C(X,\bS^n): \restr{F}{A} = f$ iff $f^*(H^n(\bS^n)) \subset i^*(H^n(X))$, $i:A \ra X$ the inclusion.\\

\index{Theorem, Hopf Classification}
\textbf{\small HOPF CLASSIFICATION THEOREM} \quad 
Let $(X,A)$ be a relative CW complex with $\dim (X,A) \leq n$ $(n \geq 1)$.  Fix a generator  $\iota \in H^n(\bS^n,s_n;\Z)$ $-$then the assignment $[f] \ra f^* \iota$ defines a bijection $[X,A;\bS^n,s_n] \ra H^n(X,A;\Z)$,\\

\label{5.0ax}
\begingroup%%----------------------------------->>
\fontsize{9pt}{11pt}\selectfont
\textbf{\small EXAMPLE} \ 
The unit tangent bundle of $\bS^{2n}$ can be identified with the Stiefel manifold $\bV_{2n+1,2}$.  
It is $(2n - 2)$-connected with euclidean dimension $4n-1$.  
One has 
$H_q(\bV_{2n+1,2}) \approx \bZ$ $(q = 0, 4n - 1)$, 
$H_{2n-1}(\bV_{2n+1,2}) \approx \Z/2\Z$, and 
$H_q(\bV_{2n+1,2}) = 0$ otherwise.  
By the Hopf classification theorem, 
$[\bV_{2n+1,2},\bS^{4n-1}] \approx$ 
$H^{4n-1}(\bV_{2n+1,2})$, so there is a map 
$f:\bV_{2n+1,2} \ra \bS^{4n-1}$ such that $f^*$ induces an isomorphism 
$H^{4n-1}(\bS^{4n-1}) \ra$ 
$H^{4n-1}(\bV_{2n+1,2})$.  
Consequently, under $f_*$, 
$H_*(\bV_{2n+1,2});\Q) \approx$ 
$H_*(\bS^{4n-1};\Q)$, 
thus the mapping fiber $E_f$ of $f$ is rationally acyclic, i.e., 
$\widetilde{H}_*(E_f;\Q) = 0$ (cf. \pageref{5.0g}).\\
\endgroup%%------------------------------------<<

Let
$
\begin{cases}
\ X\\
\ Y
\end{cases}
$
be CW complexes with CW structures
$
\begin{cases}
\ \{X^{(n)}\}\\
\ \{Y^{(n)}\}
\end{cases}
$
$-$then a 
\un{skeletal map}
\index{skeletal map} 
is a continuous function $f:X \ra Y$ such that $\forall \ n:$ $f(X^{(n)}) \subset Y^{(n)}$.

[Note: \  A CW complex is filtered by its skeletons, 
so the term ``skeletal map'' is just the name used for ``filtered map'' in the CW context.]\\

\begingroup%%----------------------------------->>
\fontsize{9pt}{11pt}\selectfont
\index{simplicial sets (example)}
\textbf{\small EXAMPLE \ \un{Simplicial Sets}} \ 
If $f:X \ra Y$ is a simplicial map, then $\abs{f}:\abs{X} \ra \abs{Y}$ is a skeletal map and transforms cells of $\abs{X}$ onto cells of $\abs{Y}$.\\
\endgroup%%------------------------------------<<

\index{Theorem, Skeletal Approximation}
\textbf{\small SKELETAL APPROXIMATION THEOREM} \ 
Let $X$ and $Y$ be CW complexes.  Suppose that $A$ is a subcomplex of $X$ $-$then for any continuous function $f:X \ra Y$  such that $\restr{f}{A}$ is skeletal there exists a skeletal map $g:X \ra Y$  such that $\restr{f}{A} = \restr{g}{A}$ and $f \simeq g$ rel $A$.

[Note: \  In particular, every continuous function $f:X \ra Y$ is homotopic to a skeletal map 
$g:X \ra Y$.]\\

Let
$
\begin{cases}
\ (X,A)\\[-.1cm]
\ (Y,B)
\end{cases}
$
be relative CW complexes with relative CW structures
$
\begin{cases}
\ \{(X,A)^{(n)}\}\\[-.1cm]
\ \{(Y,B)^{(n)}\}
\end{cases}
$
$-$then a 
\un{relative skeletal map}
\index{relative skeletal map}
is a continuous function $f:(X,A) \ra (Y,B)$ such that $\forall \ n:$ $ f((X,A)^{(n)}) \subset (Y,B)^{(n)}$.\\

%%----------------------------------------------------------------------------------------------10
\index{Theorem, Relative Skeletal Approximation}
\textbf{\small RELATIVE SKELETAL APPROXIMATION THEOREM} \ 
Let $(X,A)$ and $(Y,B)$ be relative CW complexes $-$then every continuous function $f:(X,A) \ra (Y,B)$ is homotopic rel $A$ to a relative skeletal map $g:(X,A) \ra (Y,B)$. \\

Here is a summary of the main topological properties of CW complexes.\\
$
\begin{array}{lll}
\indent\indent \text{(TCW$_1$)} \quadx &\text{Every CW complex is compactly generated.}\\
\indent\indent \text{(TCW$_2$}) \quadx &\text{Every CW complex is stratifiable, hence is hereditarily paracompact.}\\
\indent\indent \text{(TCW$_3$)} \quadx &\text{Every CW complex is uniformly locally contractible, therefore locally } \\
\text{contractible.}\\
\label{14.48}
\label{14.61}
\indent\indent \text{(TCW$_4$)} \quadx &\text{Every CW complex is numerably contractible.}\\
\indent\indent \text{(TCW$_5$)} \quadx &\text{Every CW complex is locally path connected.}\\
\indent\indent \text{(TCW$_6$)} \quadx &\text{Every CW complex is the coproduct of its path components and these }\\
\text{are subcomplexes.}\\
\indent\indent \text{(TCW$_7$)} \quadx &\text{Every connected CW complex is path connected.}\\
\indent\indent \text{(TCW$_8$)} \quadx &\text{Every connected CW complex has a universal covering space.}\\
\end{array}
$
%\begin{array}{lll}
%\indent &(TCW$_1$) \quadx &Every CW complex is compactly generated.\\
%\indent &(TCW$_2$) \quadx &Every CW complex is stratifiable, hence is hereditarily paracompact.\\
%\indent &(TCW$_3$) \quadx &Every CW complex is uniformly locally contractible, therefore locally \\
%&contractible.\\
%\indent &(TCW$_4$) \quadx &Every CW complex is numerably contractibe.\\
%\indent &(TCW$_5$) \quadx &Every CW complex is locally path connected.\\
%\indent &(TCW$_6$) \quadx &Every CW complex is the coproduct of its path components and these \\
%&are subcomplexes.\\
%\indent &(TCW$_7$) \quadx &Every CW complex is path connected.\\
%\indent &(TCW$_8$) \quadx &Every CW complex has a universal covering space.\\
%\end{array}

[Note: \  If $X$ is a connected CW complex with CW structure $\{X^{(n)}\}$ and if $p:\widetilde{X} \ra X$ is a covering projection, then the sequence $\{\widetilde{X}^{(n)} = p^{-1}(X^{(n)})\}$ is a CW structure on $\widetilde{X}$ with respect to which $p$ is skeletal.]

If $(X,A)$ is a relative CW complex, then certain topological properties of $A$ are automatically transmitted to $X$.  
For example, if $A$ is in \bCG, $\dcg$, or \textbf{CGH}, then the same holds for $X$.  
Analogous remarks apply to a Hausdorff \mA which is normal, perfectly normal, paracompact, etc.
\\[-.2cm]

\indent\indent (F) \  A CW complex $X$ is said to be 
\un{finite}
\index{finite (CW complex)} 
if $\#(\sE) < \omega$.  Every finite CW complex is compact and conversely.  
A compact subset of a CW complex is contained in a finite subcomplex.\\
\indent\indent (C) \  A CW complex $X$ is said to be 
\un{countable}
\index{countable (CW complex)} 
if $\#(\sE) \leq \omega$.  A CW complex is countable iff it does not contain an uncountable discrete set.  Every countable CW complex is Lindel\"of and conversely.

[Note: \  The homotopy groups of a countable connected CW complex are countable.]\\
\indent\indent (LF) \  A CW complex $X$ is said to be 
\un{locally finite}
\index{locally finite (CW complex)} 
if each $x \in X$ has a neighborhood $U$ such that $U$ is contained in a finite subcomplex of $X$.  Every locally finite CW complex is locally compact and conversely.  Every locally finite CW complex is metrizable and conversely.  
A locally finite connected CW complex is countable.\\

\begingroup%%----------------------------------->>
\fontsize{9pt}{11pt}\selectfont
What spaces carry a CW structure?  There is no known characterization but the foregoing conditions impose a priori limitations.  For example, a nonmetrizable LCH space cannot be equipped with a CW
%%----------------------------------------------------------------------------------------------11
structure.  
On the other hand, the Cantor set and the Hilbert cube are metrizable compact Hausdorff spaces but neither supports a CW structure.
\vspi
[Note: \  Every compact differentiable manifold can be triangulated 
but examples are known of compact topological manifolds that cannot be triangulated, i.e., that are not polyhedrons 
(David-Januszkiewicz\footnote[2]{\textit{J. Differential Geom.} \textbf{34} (1991), 347-388.}).]\\
\endgroup%%------------------------------------<<

\label{19.24}
\label{2.4}
\index{Sorgenfrey Line}
\begingroup%%----------------------------------->>
\fontsize{9pt}{11pt}\selectfont
\textbf{\small EXAMPLE \  (\un{The Sorgenfrey Line})} \ 
Topologize $X = \R$ by choosing for the basic neighborhoods of a given $x$ all sets of the form $[x,y[$ $(x < y)$.  
In this topology, the line is a perfectly normal paracompact Hausdorff space but it is not locally compact.  
While not second countable, $X$ is first countable (and separable), therefore is compactly generated.  
However, $X$ is not locally connected, thus carries no CW structure.
\vspi
[Note: \  The square of the Sorgenfrey line is not normal (apply Jones' lemma).]\\
\endgroup%%------------------------------------<<

\index{Niemytzki Plane}
\begingroup%%----------------------------------->>
\fontsize{9pt}{11pt}\selectfont
\textbf{\small EXAMPLE \  (\un{The Niemytzki Plane})} \ 
Let \mX be the closed upper half plane in $\R^2$.  Topologize \mX as follows: The basic neighborhoods of $(x,y)$ $(y > 0)$ are as usual but the basic neighborhoods of $(x,0)$ are the 
$\{(x,0)\} \cup B$, where $B$ is an open disk in the upper half plane with horizontal tangent at $(x,0)$.  \mX is a compactly generated CRH space.  In addition, $X$ is Moore, hence is perfect.  And \mX is connected, locally path connected and even contractible (consider the homotopy 
$H((x,y),t) =$ 
$
\begin{cases}
\ (x,y) + t(0,1) \indent \indent (0 \leq t \leq 1/2)\\
\ t(0,1) + 2(1 - t)(x,y)  \ \  \ \  (1/2 \leq t \leq 1)
\end{cases}
\big). \ 
$ 
 However, $X$ is not normal, thus carries no CW structure.
 \vspi
 [Note: \ $X$ is neither  countably paracompact  nor metacompact but is countably metacompact.]\\
\endgroup%%------------------------------------<<

\begingroup%%----------------------------------->>
\fontsize{9pt}{11pt}\selectfont
\textbf{\small EXAMPLE} \ 
An open subset of a polyhedron is a polyhedron but an open subset of a CW complex need not be a CW complex.  
To see this, fix an enumeration $\{q_n\}$ of $\Q  \ \cap \ ]0,1[$\ .  
Consider the CW complex $X$ defined as follows: $X^{(0)} = \{0,1\}$, $X^{(1)} = [0,1]$
$
\begin{cases}
\ 0 \ra 0\\
\ 1 \ra 1
\end{cases}
$ 
and at each point $q_n$ attach a 2-cell by taking for $f_n:\bS^1 \ra X^{(1)}$ the constant map $f_n = q_n$.  
Choose a point $x_n \in e_n$ $(\in \sE_2)$ and put $A = \{x_n\}$ $-$then \mA is closed and $U = X - A$ carries no CW structure.
\vspi
[Otherwise: (a) $[0,1] \subset U^{(1)}$; (b) $\forall \ n$, $U^{(1)} \cap e_n \neq \emptyset$; (c) $\forall \ n$ $q_n \in U^{(0)}$.]\\
\endgroup%%------------------------------------<<

\begin{proposition} \ %02
Every CW complex has the homotopy type of a polyhedron.
\end{proposition}

[Let $X$ be a CW complex with CW structure $\{X^{(n)}\}: X = \colimx  X^{(n)}$.  
Taking into account $\S 3$, Proposition 15, it will be enough to construct a sequence of vertex schemes $K_{(n)}$ such that $\forall \ n$, $K_{(n-1)}$ is a subscheme of $K_{(n)}$ and a sequence of homotopy equivalences $\phi_n:X^{(n)} \ra \abs{K_{(n)}}$ such that $\forall \ n$, $\restr{\phi_n}{X^{(n-1)}} = \phi_{n-1}$.  
Proceeding by induction, make the obvious choices when $n = 0$ and then assume that $K_{(0)}, \ldots, K_{(n-1)}$ and $\phi_0, \ldots, \phi_{n-1}$ have been defined.  At level $n$ there is an index set $I_n$ and a pushout
%%----------------------------------------------------------------------------------------------12
square
\begin{tikzcd}%[sep=large]
{I_n \cdot \dot{\Delta}^n} \arrow{r}{f} \arrow{d} &X^{(n-1)} \arrow{d}\\
{I_n \cdot \Delta}^n\arrow{r} &X^{(n)}
\end{tikzcd}
$(f = \coprod\limits_i f_i)$.  
Given $i \in I_n$, use the simplicial approximation theorem to produce a vertex scheme $K_i$ and 
a vertex map 
$g_i:K_i \ra K_{(n-1)}$ with 
$\abs{K_i} = \dot{\Delta}^n$ and 
$\abs{g_i} \simeq \phi_{n-1} \circx f_i$.  
Combine the $K_i$ and put
$\abs{g} = \coprod\limits_i \abs{g_i}$.  
The adjunction theorem implies that there exists a vertex scheme $K_{(n)}$ containing $K_{(n-1)}$ as a subscheme and a homeomorphism $I_n \cdot \Delta^n \sqcup_{\abs{g}} \abs{K_{(n-1)}} \ra \abs{K_{(n)}}$ 
whose restriction to 
$\abs{K_{(n-1)}}$ is the identity map.  The triangle
\begin{tikzcd}%[sep=large]
{I_n \cdot \dot{\Delta}^n} \arrow{r}{f} \arrow{rd}[swap]{\abs{g}} &{X^{(n-1)}} \arrow{d}{\phi_{n-1}}\\
&{\abs{K_{(n-1)}}}
\end{tikzcd}
is homotopy commutative: $\abs{g} \simeq \phi_{n-1} \circx f$.  
Since $\phi_{n-1}$ is a homotopy equivalence, one can find a homotopy equivalence 
$\phi_n:I_n \cdot \Delta^n \sqcup_f X^{(n-1)} \ra I_n \cdot \Delta^n \sqcup_{\abs{g}} \abs{K_{(n-1)}}$ such that $\restr{\phi_n}{X^{(n-1)}} = \phi_{n-1}$ (cf. p. \pageref{5.4}), which completes the induction.]

[Note: \  Similar methods lead to the expected analogs in $\bCW^2$ or $\bCW_*$.  
Consider e.g., a CW pair $(X,A)$ with relative CW structure $\{(X,A)^{(n)}\}:$ $(X,A)^{(n)} = X^{(n)} \cup A$.  
Choose a vertex scheme $L$ and a homotopy equivalence 
$\phi:A \ra \abs{L}$ $-$then there is a vertex scheme $K_{(0)}$ containing \mL as a subscheme and a 
homotopy equivalence of pairs $((X,A)^{(0)},A) \ra (\abs{K_{(0)}},\abs{L})$ so, 
arguing as above, there is a vertex scheme $K$ containing $L$ as a subscheme and a homotopy equivalence $\Phi:X \ra \abs{K}$ such that $\restr{\Phi}{A} = \phi$.  
Conclusion:  In $\bHTOP^2$, $(X,A) \approx (\abs{K},\abs{L})$ (cf. $\S 3$, Proposition 14).]\\

\label{6.15}
\begin{proposition} \ %03
Let $X$ be a CW complex.  
Assume: (i) $X$ is finite (countable) or (ii) $\dim X \leq n$ 
$-$then there exists a vertex scheme $K$ such that $X$ has the homotopy type of 
$\abs{K}$, where (i) $K$ is finite (countable) or (ii) $\dim K \leq n$.
\end{proposition}

[This is implicit in the proof of the preceding proposition.]\\

\label{5.2}
\begingroup%%----------------------------------->>
\fontsize{9pt}{11pt}\selectfont
Let $X$ be a CW complex; let $\sA$ be the collection of finite subcomplexes of $X$ 
$-$then $\sA$ is an absolute closure preserving closed covering of $X$.  
Since every finite subcomplex of $X$ is a second countable compact Hausdorff space of finite topological dimension, 
it follows that $X$ can be embedded as a closed subspace of a polyhedron (cf. p.\pageref{5.5}).\\
\endgroup%%------------------------------------<<

\label{4.10}
\label{5.7}
\label{6.43}
\begingroup%%----------------------------------->>
\fontsize{9pt}{11pt}\selectfont
\textbf{\small FACT} \ 
Every CW complex is the retract of a polyhedron, hence every open subset of a CW complex is the retract of a polyhedron.\\
\endgroup%%------------------------------------<<

\begingroup%%----------------------------------->>
\fontsize{9pt}{11pt}\selectfont
\textbf{\small EXAMPLE} \ 
Every polyhedron is a CW complex but there exist CW complexes that cannot be triangulated.  
Thus let $f(t) = t \sin (\pi/2t)$ $(0 < t \leq 1)$ and set $f(0) = 0$.  
Denote by $m$ the absolute minimum 
%%----------------------------------------------------------------------------------------------13
of $f$ on $[0,1]$ (so $-1 < m < 0$).  
Take for $X$ the image of the square $[0,1] \times [0,1]$ under the map $(u,v) \mapsto (u, uv, f(v))$.  
The following subspaces constitute a CW structure on \mX:
\begin{align*}
X^{(0)} &= \{(0,0,0),(1,0,0),(0,0,1),(1,1,1),(0,0,m)\},\\
X^{(1)} &= \{(u,0,0): 0 \leq u \leq 1\} \cup \{(u,u,1): 0 \leq u \leq 1\} \cup
\begin{cases}
\hsx \{(0,0,v): m \leq v \leq 0\}\\
\hsx \{(0,0,v): 0 \leq v \leq 1\}
\end{cases}
\hspace{-.4cm}\cup \{(1,v,f(v)): 0 \leq v \leq 1\},
\end{align*}
and $X^{(2)} = X$.  
Using the fact that $f$ has a sequence $\{M_n\}$ of relative maxima: $M_1 > M_2 > \cdots $ $(1 > M_1)$, 
look at the $(0,0,M_n)$ and deduce that \mX is not a polyhedron.\\
\endgroup%%------------------------------------<<

\label{13.12}
\begingroup%%----------------------------------->>
\fontsize{9pt}{11pt}\selectfont
\textbf{\small FACT} \ 
Let $X$ be a CW complex.  Suppose that all the characteristic maps are embeddings $-$then $X$ is a polyhedron.\\
\endgroup%%------------------------------------<<

There are two other issues.\\
\index{CW complexes - products}
\indent \indent (Products)  Let
$
\begin{cases}
\ X \\[-.1cm]
\ Y
\end{cases}
$
be CW complexes with CW structures
$
\begin{cases}
\ \{X^{(n)}\}\\[-.1cm]
\ \{Y^{(n)}\}
\end{cases}
$
\hspace{-.26cm}.  
Put $(X \times_k Y)^{(n)} = \bigcup\limits_{p+q=n} X^{(p)} \times_k Y^{(q)}$.  
Consider $X \times_k Y$ $-$then the sequence $\{(X \times_k Y)^{(n)}\}$ satisfies 
$\text{CW}_1$, $\text{CW}_2$, and $\text{CW}_3$ above, meaning it is a CW structure on $X \times_k Y$.  
When can ``$\times_k$'' be replaces by ``$\times$''?  
Useful sufficient conditions to ensure this are that one of the factors be locally finite or that both of the factors be countable 
(necessary conditions have been discussed by 
Tanaka\footnote[2]{\textit{Proc. Amer. Math. Soc.} \textbf{86} (1982), 503-507.}).\\

\index{Dowker's Product}
\begingroup%%----------------------------------->>
\fontsize{9pt}{11pt}\selectfont
\textbf{\small EXAMPLE \  (\un{Dowker's Product})} \ \ 
Suppose that $X$ and $Y$ are CW complexes $-$then the product $X \times Y$ need not be compactly generated, hence, 
when this happens, $X \times Y$ is not a CW complex.  
Here is an illustration.  
Definition of \mX:  Put $X^{(0)} = \N^\N \cup \{0\}$ (discrete topology), let $f_s:\{0,1\} \ra X^{(0)}$ be the map
$
\begin{cases}
\ 0 \ra 0\\
\ 1 \ra s
\end{cases}
$
$(s \in \N^\N)$, write $X^{(1)}$ for the space thereby obtained from $X^{(0)}$ by attaching 1-cells,
 and take $X = X^{(0)} \cup X^{(1)}$.  
 Definitition of \mY:  Put $Y^{(0)} = \N \cup \{0\}$ (discrete topology), 
 let $f_n:\{0,1\} \ra Y^{(0)}$ be the map
$
\begin{cases}
\ 0 \ra 0\\
\ 1 \ra n
\end{cases} 
$
$(n \in \N)$, write $ Y^{(1)}$ for the space thereby obtained from $Y^{(0)}$ by attaching 1-cells, 
and take $Y = Y^{(0)}  \cup Y^{(1)} $.  Let $\Phi_s$ ($\Phi_n$) be the characteristic map of the 1-cell corresponding to the $s \in \N^\N$ ($n \in \N$).  
Consider the following subset of $X \times Y$: $K = \{(\Phi_s(1/s_n),\Phi_n(1/s_n)): (s,n) \in \N^\N \times \N\}$.  
Evidently $K$ is a closed subset of $X \times_k Y$.  But $K$ is not a closed subset of $X \times Y$.  
For if it were, $X \times Y - K$ would be open and since the point $(0,0) \in X \times Y - K$, 
there would be a basic neighborhood $U \times V: (0,0) \in U \times V \subset X \times Y - K$.  
Given $s \in \N^\N$, $\exists$ a real number $a_s: 0 < a_s \leq 1$ such that $U \supset \{\Phi_s(p): p < a_s\}$ and given $n \in \N$, $\exists$ a real number $b_n: 0 < b_n \leq 1$ such that $V \supset \{\Phi_n(q): q < b_n\}$.  
Define $\ov{s} \in \N^\N$ by $\ov{s}_n = 1 + [\max \{n,1/b_n\}]$ (so $\ov{s}_n > n$ $\&$ $\ov{s}_n > 1/b_n$); 
define $\ov{n} \in \N$ by $\ov{n} = 1 + [1/a_{\ov{s}}]$ (so
%%----------------------------------------------------------------------------------------------14
$\ov{n} > 1/a_{\ov{s}}$) $-$then the pair $(\Phi_{\ov{s}}(1/\ov{s}_{\ov{n}}), \Phi_{\ov{n}}(1/\ov{s}_{\ov{n}})$ is in both $U \times V$ and $K$.  Contradiction.  
Incidentally, one can show that the projections
$
\begin{cases}
\ X \times_k Y \ra X\\
\ X \times_k Y \ra Y
\end{cases} 
$
are not Hurewicz fibrations (although, of course, they are \bCG fibrations).
\vspi
[Note: \  This construction has an obvious interpretation in terms of cones.  
Observe too that $X$ and $Y$ are polyhedrons.  
\label{5.3}
\label{5.53}
Corollary:  The square of a polyhedron need not be a polyhedron.]\\
\endgroup%%------------------------------------<<

\begingroup%%----------------------------------->>
\fontsize{9pt}{11pt}\selectfont
\textbf{\small FACT} \ 
Every countable CW complex has the homotopy type of a locally finite countable CW complex.
\\ \indent
[Let $X$ be a countable CW complex.  
Fix an enumeration $\{e_k\}$ of its cells.  
Given $e_k$, denote by $X(e_k)$ the intersection of all subcomplexes of $X$ containing $e_k$ $-$then $X(e_k)$ is a finite subcomplex of $X$.  
Put $X^n = \ds\bigcup\limits_0^n X(e_k): X^0 \subset X^1 \subset \cdots $ 
is an expanding sequence of topological spaces with $X^\infty = X$.  
The telescope $\tel X^\infty$ of $X^\infty$ has the same homotopy type as $X^\infty = X$ 
(cf. p. \pageref{5.0h}) and is a CW complex.  
In fact, $\tel X^\infty$ is the subcomplex of $X \times_k[0,\infty[ \ =X \times [0,\infty[$  made up of the cells $e \times \{n\}$, 
$e \times \hsx ]n,n+1[$, where $e$ is a cell of $X^m$ $(m \leq n)$, a description which makes it clear that $\tel X^\infty$ is locally finite.]
\\ \indent
[Note: \  Suppose that $X$ is a locally finite countable CW complex 
$-$then there exists a sequence of finite subcomplexes $X_n$ such that 
$\forall \ n$, $X_n \subset \itry X_{n+1}$, with $X = \ds\bigcup\limits_n X_n$.]\\
\endgroup%%------------------------------------<<
%% -------------------<

\index{Adjunctions (CW complexes)}
\indent \indent (Adjunctions) \quadx Let
$
\begin{cases}
\ X \\[-.1cm]
\ Y
\end{cases}
$
be CW complexes with CW structures
$
\begin{cases}
\ \{X^{(n)}\}\\[-.1cm]
\ \{Y^{(n)}\}
\end{cases}
\hspace{-.25cm}. \ 
$
Suppose that \mA is a subcomplex of $X$. Let $f:A \ra Y$ be a skeletal map $-$then the adjunction space $X \sqcup_f Y$ is a CW complex, the CW structure being 
$\{X^{(n)} \sqcup_{f^{(n)}} Y^{(n)}\}$ $(f^{(n)} = \restr{f}{A^{(n)}})$.  
Examples: 
(1)  If $X$ is a CW complex and if $A \subset X$ is a subcomplex, then the quotient $X/A$ is a CW complex;
(2)  If $X$ is a CW complex , then its cone $\Gamma X$ and its suspension $\Sigma X$ are CW complexes;
(3)  If $X$ and $Y$ are CW complexes and if $f:X \ra Y$ is a skeletal map, then the mapping cylinder $M_f$ of $f$ is a CW complex, containing both $X$ and $Y$ as embedded subcomplexes; 
(4)  If $X$ and $Y$ are CW complexes and if $f:X \ra Y$ is a skeletal map, then the mapping cone $C_f$ of $f$ is a CW complex containing $Y$ as an embedded subcomplex.

[Note: \  There are also pointed analogs of these results.  For example, if 
$
\begin{cases}
\ (X,x_0) \\
\ (Y,y_0)
\end{cases}
$
are pointed CW complexes, then the smash product $X \#_k Y$ is a pointed CW complex.]\\

\label{6.13}
\begingroup%%----------------------------------->>
\fontsize{9pt}{11pt}\selectfont
Let $X$ and $Y$ be CW complexes.  Let $A$ be a subcomplex of $X$ and let $f:A \ra Y$ be a continuous function $-$then $X \sqcup_f Y$ has the homotopy type of a CW complex:  Proof:  By the skeletal approximation theorem, there exists a skeletal map $g:A \ra Y$ such that $f \simeq g$, so  $X \sqcup_f Y$ has the same homotopy type as $X \sqcup_g Y$ (cf. p. \pageref{5.0i}).\\
\endgroup%%------------------------------------<<

\begingroup%%----------------------------------->>
\fontsize{9pt}{11pt}\selectfont
\textbf{\small FACT} \ 
A CW complex is path connected iff its $1$-skeleton is path connected.\\
\endgroup%%------------------------------------<<

%%----------------------------------------------------------------------------------------------15
\label{5.0ai}
\index{trees (example)}
\begingroup%%----------------------------------->>
\fontsize{9pt}{11pt}\selectfont
\textbf{\small EXAMPLE  \ (\un{Trees})} \ 
Let $X$ be a nonempty connected CW complex $-$then a 
\un{tree}
\index{tree} 
in $X$ is a nonempty simply connected subcomplex $T$ of $X$ with $\dim T \leq 1$.  
Every tree in $X$ is contractible and contained in a maximal tree.  A tree is maximal iff it contains $X^{(0)}$.  
If $T$ is a maximal tree in $X$, then $X/T$ is a connected CW complex with exactly one 0-cell and the projection $X \ra X/T$ is a homotopy equivalence (cf. p. \pageref{5.6}).\\
\endgroup%%------------------------------------<<

\label{5.0y}
\label{5.49}
\label{13.21}
\index{WHE Criterion}
%\indent\textbf{\small WHE CRITERION} \quadx
\textbf{\small WHE CRITERION} \ 
Let
$
\begin{cases}
\ X \\[-.1cm]
\ Y
\end{cases}
$
be topological spaces, $f:X \ra Y$ a continuous function $-$then $f$ is a weak homotopy equivalence if for any finite CW pair $(K,L)$ and any diagram
\begin{tikzcd}%[sep=large]
L \arrow{r}{\phi} \arrow{d} &X \arrow{d}{f}\\
K \arrow{r}[swap]{\psi} &Y
\end{tikzcd}
, where $f \circx \phi = \restr{\psi}{L}$, there exists a $\Phi:K \ra X$ such that $\restr{\Phi}{L} = \phi$ and 
$f \circx \Phi \simeq \psi$ rel $L$.

[Indeed, diagrams of the form
\begin{tikzcd}%[sep=large]
s_n \arrow{r} \arrow{d} &X \arrow{d}{f}\\
\bS^n \arrow{r} &Y
\end{tikzcd}
,
\begin{tikzcd}%[sep=large]
\bS^n \arrow{r} \arrow{d} &X \arrow{d}{f}\\
\bD^{n+1} \arrow{r} &Y
\end{tikzcd}
evidently suffice.]\\

\textbf{\small LEMMA} \ 
Suppose that $f:X \ra Y$ is an $n$-equivalence $-$then in any diagram
\begin{tikzcd}%[sep=large]
\bS^{n-1} \arrow{r}{\phi} \arrow{d} &X \arrow{d}{f}\\
\bD^{n} \arrow{r}[swap]{\psi} &Y
\end{tikzcd}
,
 where $f \circx \phi \simeq \psi$ on $\bS^{n-1}$ by $h: I\bS^{n-1} \ra Y$, there exists a $\Phi:\bD^{n} \ra X$ such that 
$\restr{\Phi}{\bS^{n-1}} = \phi$ and $H:I \bD^{n} \ra Y$ such that $\restr{H}{I \bS^{n-1}} = h$ and $f \circx \Phi \simeq \psi$ on $\bD^{n}$ by $H$.\\
\vspace{0.25cm}

\index{Homotopy Extension Lifting Property}
\textbf{\small HOMOTOPY EXTENSION LIFTING PROPERTY} \ 
Suppose that $f:X \ra Y$ is a weak homotopy equivalence.  Let $(K,L)$ be a relative CW complex $-$then in any diagram
\begin{tikzcd}%[sep=large]
L \arrow{r}{\phi} \arrow{d} &X \arrow{d}{f}\\
K \arrow{r}[swap]{\psi} &Y
\end{tikzcd}
, 
where $f \circx \phi \simeq \psi$ on $L$ by $h:IL \ra Y$, there exists a $\Phi:K \ra X$ such that $\restr{\Phi}{L} = \phi$ and an $H:IK \ra Y$ such that $\restr{H}{IL} = h$ and $f \circx \Phi \simeq \psi$ on $K$ by $H$.\\

\label{4.47}
Application: Let $f:X \ra Y$ be a weak homotopy equivalence $-$then for any CW complex $K$, the arrow $f_*:[K,X] \ra [K,Y]$ is bijective.

[To see that $f_*$ is surjective (injective), apply the homotopy extension lifting property to 
$(K,\emptyset)$ $((IK,i_0K \cup i_1K))$.]

[Note: \  The condition is also characteristic.  
Thus first take $K = *$ and reduce to
%%----------------------------------------------------------------------------------------------16
when
$
\begin{cases}
\ X \\[-.15cm]
\ Y
\end{cases}
$
are path connected.  
Next take $K = \ds\bigvee\limits_I \bS^1$ ($I$ a suitable index set) to get that $\forall \ x \in X$, 
$f_*:\pi_1(X,x) \ra \pi_1(Y,f(x))$ is surjective.  Finish by taking $K = \bS^n$ (cf. p. \pageref{5.0j}).]\\

\label{5.0ama}
\label{9.108}
\begingroup%%----------------------------------->>
\fontsize{9pt}{11pt}\selectfont
\textbf{\small EXAMPLE} \ 
Let 
$
\begin{cases}
\ (X,x_0) \\
\ (Y,y_0)
\end{cases}
$
be pointed connected CW complexes.  
Suppose that $f \in C(X,x_0,Y,y_0)$ has the property that $\forall \ n > 1$, $f_*:\pi_n(X,x_0) \ra \pi_n(Y,y_0)$ is bijective 
$-$then for any pointed simply connected CW complex $(K,k_0)$, the arrow $f_*:[K,k_0;X,x_0] \ra [K,k_0;Y,y_0]$ is bijective.\\
\endgroup%%------------------------------------<<
%% -------------------<

\label{9.109}
\label{12.20}
\begingroup%%----------------------------------->>
\fontsize{9pt}{11pt}\selectfont
\textbf{\small FACT} \ 
Let $p:X \ra B$ be a continuous function $-$then $p$ is both a weak homotopy equivalence and a Serre fibration iff for any relative CW complex $(K,L)$ and any diagram
\begin{tikzcd}[sep=large]
L \arrow{r}{\phi} \arrow{d} &X \arrow{d}{p}\\
K \arrow{r}[swap]{\psi} &B
\end{tikzcd}
, where $p \circx \phi = \restr{\psi}{L}$, there exists a $\Phi:K \ra X$ such that $\restr{\Phi}{L} = \phi$ and $p \circx \Phi = \psi$.
\vspi
\label{12.19}
[Note: \  The characterization can be simplified:  
A continuous function $p:X \ra B$ is both a weak homotopy equivalence and a Serre fibration iff every commutative diagram
\begin{tikzcd}[sep=large]
\bS^{n-1} \arrow{r} \arrow{d} &X \arrow{d}\\
\bD^n \arrow{r}&B
\end{tikzcd}
$(n \geq 0)$ admits a filler $\bD^n \ra X$.]\\
\endgroup%%------------------------------------<<
%% -------------------<

\begingroup%%----------------------------------->>
\fontsize{9pt}{11pt}\selectfont
Application:  Let 
\begin{tikzcd}[sep=large]
X^\prime \arrow{r} \arrow{d}[swap]{p^\prime} &X \arrow{d}{p}\\
B^\prime \arrow{r} &B
\end{tikzcd}
be a pullback square.  Suppose that $p$ is a Serre fibration and a weak homotopy equivalence $-$then $p^\prime$ is a Serre fibration and a weak homotopy equivalence.\\
\endgroup%%------------------------------------<<
%% -------------------<

\label{9.110}
A continuous function $f:(X,A) \ra (Y,B)$ is said to be a 
\un{weak homotopy equivalence} 
\un{of pairs}
\index{weak homotopy equivalence of pairs}
provided that $f:X \ra Y$ and $f:A \ra B$ are weak homotopy equivalences.

[Note: \  A weak homotopy equivalence of pairs is a relative weak homotopy equivalence (cf. p. \pageref{5.0k}) but not conversely.]\\

Application:  Let $f:(X,A) \ra (Y,B)$ be a weak homotopy equivalence of pairs $-$then for any CW pair $(K,L)$, the arrow $f_*:[K,L;X,A] \ra [K,L;Y,B]$ is bijective.

[Note: \  The condition is also characteristic. 
For $[K,\emptyset;X,A] \approx [K,\emptyset;Y,B]$ $\implies$ $[K,X] \approx [K,Y]$ and $[IK, i_0K;X,A] \approx [IK,i_0K;Y,B]$ 
$\implies$ $[K,A] \approx [K,B]$.]\\

\index{Theorem, Realization Theorem}
\textbf{\small REALIZATION THEOREM} \ 
Suppose that $X$ and $Y$ are CW complexes.  Let $f:X \ra Y$ be a weak homotopy equivalence $-$then $f$ is a homotopy equivalence.

[Note: \  It is a corollary that the result remains true when $X$ and $Y$ have the homotopy type of CW complexes.]\\

%%----------------------------------------------------------------------------------------------17
\label{5.0q}
Application:  A connected CW complex is contractible iff it is homotopically trivial.\\

\label{5.0r}
\begingroup%%----------------------------------->>
\fontsize{9pt}{11pt}\selectfont
\textbf{\small EXAMPLE} \ 
Let $X$ and $Y$ be CW complexes $-$then the identity map $X \times_k Y \ra X \times Y$ is a homotopy equivalence.
\vspi
[A priori, the identity map $X \times_k Y \ra X \times Y$ is a weak homotopy equivalence.  
However, $X$ and $Y$ each have the homotopy type of a polyhedron (cf. Proposition 2), thus the same holds for their product $X \times Y$ 
(cf. p. \pageref{5.0l}).]\\
\endgroup%%------------------------------------<<

\index{H-Groups (example)}
\begingroup%%----------------------------------->>
\fontsize{9pt}{11pt}\selectfont
\textbf{\small EXAMPLE \  (\un{H-Groups}}) \ 
Let $(X,x_0)$ be a nondegenerate homotopy associative H-space.  
Assume:  $X$ is path connected $-$then the shearing map sh: $
\begin{cases}
\ X \times X \ra X \times X \\
\ (x,y) \mapsto (x,xy)
\end{cases}
$
is a weak homotopy equivalence, thus $X$ is an H-group if $X$ carries a CW structure (cf. p. \pageref{5.0m}).\\
\endgroup%%------------------------------------<<

\label{5.44a}
The pointed version of the realization theorem says that if 
$
\begin{cases}
\ X \\[-.1cm]
\ Y
\end{cases}
$
are CW complexes and if $f:X \ra Y$ is a weak homotopy equivalence, then $f$ is a pointed homotopy equivalence for any choice of 
$
\begin{cases}
\ x_0 \in X \\[-.1cm]
\ y_0 \in Y
\end{cases}
$
with $f(x_0) = y_0$.  Proof: By the realization theorem, $f$ is a homotopy equivalence, so $f$ is actually a pointed homotopy equivalence, 
$
\begin{cases}
\ (X,x_0) \\[-.1cm]
\ (Y,y_0)
\end{cases}
$
being wellpointed (cf. p. \pageref{5.0n}).\\

\index{Theorem, Relative Realization Theorem}
\textbf{\small RELATIVE REALIZATION THEOREM} \ 
Suppose that $(X,A)$ and $(Y,B)$ are CW pairs.  Let $f:(X,A) \ra (Y,B)$ be a weak homotopy equivalence of pairs 
$-$then $f$ is a homotopy equivalence of pairs.

[Note: \  This result need not be true if one merely assumes that $f$ is a relative weak homotopy equivalence.  Example: Take $X$ path connected, fix a point $a_0 \in A$, and consider the projection 
$(X \times A, a_0 \times A) \ra (X, a_0)$.  It is a relative weak homotopy equivalence but the induced map on relative singular homology is not necessarily an isomorphism.]\\

\begingroup%%----------------------------------->>
\fontsize{9pt}{11pt}\selectfont
The relative realization theorem is a consequence of the following assertion.  
Suppose that $(X,A)$ and $(Y,B)$ are relative CW complexes.  Let $f:(X,A) \ra (Y,B)$ be a weak homotopy equivalence of pairs with 
$\restr{f}{A}:A \ra B$ a homotopy equivalence $-$then $f$ is a homotopy equivalence of pairs.\\
\endgroup%%------------------------------------<<

\begingroup%%----------------------------------->>
\fontsize{9pt}{11pt}\selectfont
\textbf{\small EXAMPLE} \ 
Let $(K,L)$ be a relative CW complex.  
Assume:  The inclusion $L \ra K$ is a weak homotopy equivalence $-$then the inclusion $L \ra K$ is a homotopy equivalence.  
Proof:  Consider the arrow $(L,L) \ra (K,L)$.\\
\endgroup%%------------------------------------<<

\begin{proposition} \ %04
Let $(Y,B)$ and $(Y^\prime,B^\prime)$ be pairs and let 
$h:(Y,B)  \ra (Y^\prime,B^\prime)$ be a continuous function; 
let $(X,A)$ and $(X^\prime,A^\prime)$ be CW pairs and let 
$f:(X,A) \ra (Y,B)$ 
%%----------------------------------------------------------------------------------------------18
$\&$ $f^\prime: (X^\prime,A^\prime) \ra (Y^\prime,B^\prime)$ be continuous functions.  
Assume: $f^\prime$ is a weak homotopy equivalence of pairs $-$then there exists a continuous function $g:(X,A) \ra (X^\prime,A^\prime)$, unique up to homotopy of pairs, such that the diagram
\begin{tikzcd}%[sep=large]
(X,A)  \arrow{r}{g} \arrow{d}[swap]{f} &(X^\prime,A^\prime) \arrow{d}{f^\prime}\\
(Y,B) \arrow{r}[swap]{h} &(Y^\prime,B^\prime)
\end{tikzcd}
commutes up to homotopy of pairs.
\end{proposition}

[The arrow $f_*^\prime:[X,A;X^\prime,A^\prime] \ra [X,A;Y^\prime,B^\prime]$ is bijective.]\\

\label{4.60}
Given a topological space $X$, a 
\un{CW resolution}
\index{CW resolution} 
for \mX is an ordered pair $(K,f)$, where \mK is a CW complex and $f:K \ra X$ is a weak homotopy equivalence.  The homotopy type of a CW resolution is unique.  
Proof: Let $f:K \ra X$ $\&$ $f^\prime:K^\prime \ra X$ be CW resolutions 
of $X$ $-$then by Proposition 4, there exists a continuous function $g:K \ra K^\prime$ such that the diagram 
\begin{tikzcd}%[sep=large]
{K} \ar{d}[swap]{f} \ar{r}{g} &{K^\prime} \ar{d}{f^\prime}\\
{X} \arrow[r,shift right=0.45,dash] \arrow[r,shift right=-0.45,dash]  &{X}
\end{tikzcd}
is homotopy commutative: $f \simeq f^\prime \circx g$.  
Therefore $g$ is a weak homotopy equivalence, hence a homomotopy equivalence (via the realization theorem).\\

\index{Theorem, Resolution Theorem}
\textbf{\small RESOLUTION THEOREM} \ 
Every topological space $X$ admits a CW resolution $f:K \ra X$.

[Note: \  If $X$ is path connected ($n$-connected), then one can choose $K$ path connected with $K^{(0)}$ ($K^{(n)}$) a singleton.]\\

Application: Suppose that $X$ is homotopically trivial $-$then for any CW complex $K$, the elements of $C(K,X)$ are inessential.\\

Given a pair $(X,A)$, a 
\un{relative CW resolution}
\index{relative CW resolution} 
for $(X,A)$ is an ordered pair $((K,L),f)$, where $(K,L)$ is a CW pair and $f:(K,L) \ra (X,A)$ is a weak homotopy equivalence of pairs.  A relative CW resolution is unique up to homotopy of pairs (cf. Proposition 4).\\

\index{Theorem, Relative Resolution Theorem}
\textbf{\small RELATIVE RESOLUTION THEOREM} \ 
Every pair $(X,A)$ admits a relative CW resolution $f:(K,L) \ra (X,A)$.

[Fix CW resolutions
$
\begin{cases}
\ \phi:L \ra A \\
\ \psi:K \ra X
\end{cases}
$
and let $i:A \ra X$ be the inclusion.  Using Proposition 4, choose a $g:L \ra K$ such that $\psi \circx g \simeq i \circx \phi$.  
Owing to the skeletal approximation theorem, one can assume that $g$ is skeletal, thus its mapping cylinder $M_g$ is a CW complex containing $L$ and $K$ as embedded subcomplexes.  
If $r:M_g \ra K$ is the usual retraction, then $r$ is a homotopy equivalence and $\psi \circx \restr{r}{L} \simeq i \circx \phi$.  
Since the inclusion $L \ra M_g$ is a 
%%----------------------------------------------------------------------------------------------19
cofibration, $\psi \circx r$ is homotopic to a map $f:M_g \ra X$ such that $\restr{f}{L} = i \circx \phi$.  
Change the notation to conclude the proof.]

[Note: \  If $(X,A)$ is $n$-connected, then one can choose $K$ with $K^{(n)} \subset L$.]\\

\begingroup%%----------------------------------->>
\fontsize{9pt}{11pt}\selectfont
It follows from the proof of the relative resolution theorem that given $(X,A)$ and a CW resolution $g:L \ra A$, there exists a relative CW resolution $f:(K,L) \ra (X,A)$ extending $g$.\\
\endgroup%%------------------------------------<<

Let $X$ and $Y$ be topological spaces $-$then $X$ is said to be 
\un{dominated in homotopy}
\index{dominated in homotopy} 
by $Y$ if there exist continuous functions
$
\begin{cases}
\ f:X \ra Y \\
\ g:Y \ra X
\end{cases}
$
such that $g \circx f \simeq \text{id}_X$.   
Example:  A topological space is contractible iff it is dominated in homotopy by a one point space.

[Note: \  Let $f:X \ra Y$ be a continuous function, $M_f$ its mapping cylinder $-$then $f$ admits a left homotopy inverse $g:Y \ra X$ iff $i(X)$ is a retract of $M_f$.  By comparison, $f$ is a homotopy equivalence iff $i(X)$ is a strong deformation retract of $M_f$ (cf. $\S 3$, Proposition 17).]\\

\begingroup%%----------------------------------->>
\fontsize{9pt}{11pt}\selectfont
\textbf{\small EXAMPLE} \ 
Let $X$ be a topological space which is dominated in homotopy by a compact connected $n$-manifold $Y$.  
Assume:  $H^n(X;\Z_2) \neq 0$ $-$then 
Kwasik\footnote[2]{\textit{Canad. Math. Bull.} \textbf{27} (1984), 448-451.} 
has shown that $X$ and $Y$ have the same homotopy type.\\
\endgroup%%------------------------------------<<

\label{9.106}
\begingroup%%----------------------------------->>
\fontsize{9pt}{11pt}\selectfont
\textbf{\small FACT} \ 
If $X$ is dominated in homotopy by a CW complex, then the path components of $X$ are open.\\
\endgroup%%------------------------------------<<

\index{Theorem, Domination Theorem}
\textbf{\small DOMINATION THEOREM} \ 
Let $X$ be a topological space $-$then $X$ has the homotopy type of a CW complex iff $X$ is dominated in homotopy by a CW complex.

[Suppose that $X$ is dominated in homotopy by a CW complex $Y$:
$
\begin{cases}
\ f:X \ra Y \\[-.1cm]
\ g:Y \ra X
\end{cases}
$
$\& \ g \circx f \simeq \text{id}_X$.  Fix a CW resolution $h:K \ra X$.  Using Proposition 4, choose continuous functions
$
\begin{cases}
\ f^\prime:K \ra Y \\[-.1cm]
\ g^\prime:Y \ra K
\end{cases}
$
such that the diagram
\begin{tikzcd}%[sep=large]
K \arrow{d}[swap]{h} \arrow{r}{f^\prime} &Y \arrow[equal]{d} \arrow{r}{g^\prime} &K \arrow{d}{h}\\
X \arrow{r}[swap]{f} &Y \arrow{r}[swap]{g} &X
\end{tikzcd}
is homotopy commutative.  Claim:  $h$ is a homotopy equivalence with homotopy inverse $g^\prime \circx f$.  In fact:  
$(g \circx f) \circx h$ $\simeq$ $g \circx f^\prime$ $\simeq$ $h \circx (g^\prime \circx f^\prime)$ 
$\&$ $(g \circx f) \circx h$ $\simeq$ $h \circx \text{id}_K$ 
$\implies$ 
$g^\prime \circx f^\prime$ $\simeq$ $\text{id}_K$ (cf. Proposition 4), so $(g^\prime \circx f) \circx h$ $\simeq$  $g^\prime \circx f^\prime$ $\simeq$ $\text{id}_K$ 
$\&$ $h \circx (g^\prime \circx f)$ $\simeq$ $g \circx f$ $\simeq$ $\text{id}_X$.]\\

Application:  Every retract of a CW complex has the homotopy type of a CW complex.

%%----------------------------------------------------------------------------------------------20
\label{5.0u}
[Note: \  Consequently, every open subset of a CW complex has the homotopy type of a CW complex (cf. p. )\pageref{5.7}.]\\

\index{Theorem, Countable Domination Theorem}
\textbf{\small COUNTABLE DOMINATION THEOREM} \ 
Let $X$ be a topological space $-$then $X$ has the homotopy type of a countable CW complex iff $X$ is dominated in homotopy by a countable CW complex.

[Suppose that $X$ is dominated in homotopy by a countable CW complex $Y$:
$
\begin{cases}
\ f:X \ra \\[-.1cm]
\ g:Y \ra
\end{cases}
$
$
\begin{aligned}
 X\\
 Y
\end{aligned}
$
$\& \ g \circx f \simeq \id_X$.  Using the notation of the preceding proof, consider the image $g^\prime(Y)$ of $Y$ in $K$.  
Claim:  $g^\prime(Y)$ is contained in a countable subcomplex $L_0$ of $K$.  
Indeed, for any cell $e$ of $Y$, $g^\prime(\bar{e})$ is  compact, thus is contained in a finite subcomplex of $K$ and a countable union of finite subcomplexes is a countable subcomplex.  
Fix a homotopy $H:IK \ra K$ between $g^\prime \circx f \circx h$ and $\id_K$.  
Since $IL_0$ is a countable CW complex, there exists a countable subcomplex $L_1 \subset K$: $H(IL_0) \subset L_1$.  
Iteration then gives a sequence $\{L_n\}$ of countable subcomplexes $L_n$ of $K$: $\forall \ n$, $H(IL_n) \subset L_{n+1}$.  
The union $L = \bigcup\limits_n L_n$ is a countable CW complex whose homotopy type is that of $X$.]\\

Application:  Every Lindel\"of space having the homotopy type of a CW complex has the homotopy type of a countable CW complex.

[The subcomplex generated by a Lindel\"of subspace of a CW complex is necessarily countable.]\\

\begingroup%%----------------------------------->>
\fontsize{9pt}{11pt}\selectfont
Is it true that if $X$ is dominated in homotopy by a finite CW complex, then $X$ has the homotopy type of a finite CW complex?  The answer is ``no'' in general but ``yes'' under certain assumptions.
\vspi
Notation:  Given a group $G$, let $\Z[G]$ be its integral group ring and write $\widetilde{K}_0(G)$ for the reduced Grothendieck group attached to the category of finitely generated projective $\Z[G]$-modules.
\vspi
The following results are due to Wall\footnote[2]{\textit{Ann. of Math.} \textbf{81} (1965), 56-69.}.\\
\endgroup%%------------------------------------<<
%% -------------------<

\begingroup%%----------------------------------->>
\fontsize{9pt}{11pt}\selectfont
\index{Theorem, Obstruction Theorem}
\textbf{\small OBSTRUCTION THEOREM} \ 
Suppose that $X$ is path connected and dominated in homotopy by a finite CW complex $-$then there exists an element $\widetilde{w}(X) \in \widetilde{K}_0(\pi_1(X))$ such that $\widetilde{w}(X) = 0$ iff $X$ has the homotopy type of a finite CW complex.\\
\endgroup%%------------------------------------<<

\begingroup%%----------------------------------->>
\fontsize{9pt}{11pt}\selectfont

One calls $\widetilde{w}(X)$ 
\un{Wall's obstruction to finiteness}.  
\index{Wall's obstruction to finiteness}
Example:  If $X$ is simply connected and dominated in homotopy by a finite CW complex, then $X$ has the homotopy type of a finite CW complex.\\
\endgroup%%------------------------------------<<

%%----------------------------------------------------------------------------------------------21
\begingroup%%----------------------------------->>
\fontsize{9pt}{11pt}\selectfont
\index{Fulfillment Lemma}
\textbf{\small FULFILLMENT LEMMA} \ 
Let $G$ be a finitely presented group $-$then given any $\alpha \in \widetilde{K}_0(G)$, there exists a connected CW complex $X_\alpha$ which is dominated in homotopy by a finite CW complex such that $\pi_1(X_\alpha) = G$ and $\widetilde{w}(X) = \alpha$.\\
\endgroup%%------------------------------------<<

\begingroup%%----------------------------------->>
\fontsize{9pt}{11pt}\selectfont
Let $A$ be a Dedekind domain, e.g., the ring of algebraic integers in an algebraic number field $-$then the reduced Grothendieck group of $A$ is isomorphic to the ideal class group of $A$.  
This fact, in conjunction with the fulfillment lemma, can be used to generate examples.  
Thus fix a prime $p$, put $\omega_p = \exp (2\pi\sqrt{-1}/p)$, and consider $\Z[\omega_p]$, the ring of algebraic integers in $\Q(\omega_p)$.  It is known that $\widetilde{K}_0(\Z/p\Z)$ is isomorphic to the reduced Grothendieck group of $\Z[\omega_p]$.  But the ideal class group of $\Z[\omega_p]$ is nontrivial for $p > 19$ (Montgomery).  
Moral:  There exist connected CW complexes which are dominated in homotopy by a finite CW complex, yet do not have the homotopy type of a finite CW complex.\\
\endgroup%%------------------------------------<<

\begingroup%%----------------------------------->>
\fontsize{9pt}{11pt}\selectfont
\textbf{\small EXAMPLE} \ 
Every path connected compact Hausdorff space $X$ which is dominated in homotopy by a CW complex is automatically dominated in homotopy by a finite CW complex.  
Is  $\widetilde{w}(X) = 0$?  Every connected compact ANR (in particular, every connected compact topological manifold) has the homotopy type of a CW complex 
(cf. p. \pageref{5.8}), thus is dominated in homotopy by a finite CW complex and one can prove that its Wall obstruction to finiteness must vanish, so such an $X$ does have the homotopy type of a finite CW complex.  
Still, some restriction on $X$ is necessary.  
This is because Ferry\footnote[2]{\textit{Topology} \textbf{19} (1980), 101-110; 
see also \textit{SLN} \textbf{870} (1981), 1-5 and 73-81.} 
has shown that any Hausdorff space which is dominated in homotopy by a second countable compact Hausdorff space must itself have the homotopy type of a second countable compact Hausdorff space and since there exist connected CW complexes with a nonzero Wall obstruction to finiteness, it follows that there exist path connected metrizable compacta which are dominated in homotopy by a finite CW complex, yet do not have the homotopy type of a finite CW complex.\\
\endgroup%%------------------------------------<<

\begingroup%%----------------------------------->>
\fontsize{9pt}{11pt}\selectfont
\textbf{\small EXAMPLE} \   
Suppose that $X$ is path connected and dominated in homotopy by a finite CW complex $-$then 
Gersten\footnote[3]{\textit{Amer. J. Math.} \textbf{88} (1966), 337-346; 
see also Kwasik, \textit{Comment. Math. Helv.} \textbf{58} (1983), 503-508.}
has shown that for any connected CW complex $K$ of zero Euler characteristic, the product $X \times K$ has the homotopy type of a finite CW complex, i.e., multiplication by $K$ kills Wall's obstruction to finiteness.  
For example, one can take $K = \bS^{2n+1}$.  
In particular: $X \times \bS^1$ is homotopy equivalent to a finite CW complex $Y$, say $f:X \times \bS^1 \ra Y$.  
Since $X$ is homotopy equivalent to $X \times \R$ and $X \times \R$ is the covering space of $X \times \bS^1$ determined by 
$\pi_1(X) \subset \pi_1(X \times \bS^1)$, it follows that $X$ is homotopy equivalent to the covering space 
$\widetilde{Y}$ of $Y$ determined by the subgroup $f_*(\pi_1(X))$ of $\pi_1(Y)$.  
Conclusion:  $X$ has the homotopy type of a finite dimensional CW complex.\\
\endgroup%%------------------------------------<<

A (pointed) topological space is said to be a
\un{(pointed) CW space}
\index{CW space} 
\index{pointed CW space} 
if it has the (pointed) homotopy type of a (pointed) CW complex.  
\bCWSP ($\bCWSP_*$) 
\index{\bCWSP}
\index{$\bCWSP_*$} is the full subcategory 
%%----------------------------------------------------------------------------------------------22
of \bTOP ($\bTOP_*$) whose objects are the CW spaces (pointed CW spaces) and 
\bHCWSP ($\bHCWSP_*$) 
\index{\bHCWSP}
\index{$\bHCWSP_*$} 
is the associated homotopy category.  
Example:  Suppose that $(X,A)$ is a relative CW complex, where $A$ is a CW space $-$the $X$ is a CW space.

\label{9.94}
[Note: \  If $(X,x_0)$ is a pointed CW space, then $(X,x_0)$ is nondegenerate 
(cf. p. \pageref{5.0o}).]

Every CW space is numerably contractible 
(cf. p. \pageref{5.0p}).  Every connected CW space is path connected.  Every totally disconnected CW space is discrete.  Every homotopically trivial CW space is contractible 
(cf. p. \pageref{5.0q}).

[Note: \  A CW space need not be locally path connected.]

The product $X \times Y$ of CW spaces
$
\begin{cases}
\ X \\[-.1cm]
\ Y
\end{cases}
$
is a CW space.  Proof:  There exist CW complexes
$
\begin{cases}
\ K\\[-.1cm]
\ L
\end{cases}
$
such that in \bHTOP, 
$
\begin{cases}
\ X \approx K \\[-.1cm]
\ Y \approx L
\end{cases}
$
\hspace{-.26cm}$\implies$ $X \times Y \approx K \times L \approx K \times_k L$ 
(cf. p. \pageref{5.0r}) and $K \times_k L$ is a CW complex.

\label{9.95}
A CW space need not be compactly generated.  
Example:  Suppose that $X$ is not in \bCG $-$then $\Gamma X$ is not in \bCG but $\Gamma X$ is a CW space.  
However, for any CW space $X$, the identity map $kX \ra X$ is a homotopy equivalence.\\

\begin{proposition} \ %05
Let $X$ be a connected CW space $-$then $X$ has a simply connected covering space $\widetilde{X}$ which is universal.  Moreover, every simply connected covering space of $X$ is homeomorphic over $X$ to $\widetilde{X}$.
\end{proposition}

[Fix a CW complex $K$ and a homotopy equivalence $\phi:X \ra K$.  Let $\widetilde{K}$ be a universal covering space of $K$ and define $\widetilde{X}$ by the pullback square
\begin{tikzcd}%[sep=large]
\widetilde{X}  \arrow{d} \arrow{r}{\widetilde{\phi}} &\widetilde{K} \arrow{d}\\
X  \arrow{r}[swap]{\phi} &K
\end{tikzcd}
.  
Since the covering projection $\widetilde{K} \ra K$ is a Hurewicz fibration 
(cf. p. \pageref{5.0s}), $\widetilde{\phi}$ is a homotopy equivalence 
(cf. p. \pageref{5.0t}), so $\widetilde{X}$ is a simply connected covering space of $X$. To see that $\widetilde{X}$ is universal, 
let $\widetilde{X}^\prime$ be some other connected covering space of $X$ $-$then the claim is that there is an arrow 
%\begin{tikzcd}%[sep=large]
%\widetilde{X}   \arrow{r}{f} &\widetilde{X}^\prime
%\end{tikzcd}
$\widetilde{X}   \overset{f}{\ra} \widetilde{X}^\prime$
and a commutative triangle
\begin{tikzcd}[sep=small]
\widetilde{X} \arrow{rdd}  \arrow{rr}{f} &&\widetilde{X}^\prime \arrow{ldd}\\
\\
&X
\end{tikzcd}
.  For this, form the pullback square
\begin{tikzcd}%[sep=large]
\widetilde{K}^\prime  \arrow{d} \arrow{r}{\widetilde{\psi}} &\widetilde{X}^\prime \arrow{d}\\
K  \arrow{r}[swap]{\psi} &X
\end{tikzcd}
, 
$\psi$ a homotopy inverse for $\phi$.  Due to the universality of $\widetilde{K}$, there is an arrow
%\begin{tikzcd}%[sep=large]
%\widetilde{K}   \arrow{r}{g} &\widetilde{K}^\prime
%\end{tikzcd}
$\widetilde{K}   \overset{g}{\ra} \widetilde{K}^\prime$
and a commutative triangle
\begin{tikzcd}[sep=small]
\widetilde{K} \arrow{rdd}  \arrow{rr}{g} &&\widetilde{K}^\prime \arrow{ldd}\\
\\
&K
\end{tikzcd}
.  
Consider the diagram
%%----------------------------------------------------------------------------------------------23
\[
\begin{tikzcd}%[sep=large]
\widetilde{X}  \arrow{d}[swap]{p} \arrow{r}{\widetilde{\phi}} 
&\widetilde{K}  \arrow{d}  \arrow{r}{g}   
&\widetilde{K}^\prime  \arrow{d} \arrow{r}{\widetilde{\psi}} 
&\widetilde{X}^\prime \arrow{d}{p^\prime}\\
X  \arrow{r}[swap]{\phi}    
&K  \arrow[r,shift right=0.45,dash] \arrow[r,shift right=-0.45,dash]                 
&K  \arrow{r}[swap]{\psi} &X
\end{tikzcd}
.
\]
From the definitions, $p^\prime \circx \widetilde{\psi} \circx g \circx \widetilde{\phi} = \psi \circx \phi \circx p \simeq p$, 
thus $\exists \ f \in C_X(\widetilde{X}, \widetilde{X}^\prime): f \simeq \widetilde{\psi} \circx g \circx \widetilde{\phi}$.  
Finally, if $\widetilde{X}^\prime$ is simply connected, then $\widetilde{K}^\prime$ is simply connected and one can assume that $g$ is a homeomorphism.   
Therefore $f$ is a fiber homotopy equivalence (cd. $\S 4$, Proposition 15).  
Because the fibers are discrete, it follows that $f$ is also an open bijection, hence is a homeomorphism.]\\

\begingroup%%----------------------------------->>
\fontsize{9pt}{11pt}\selectfont
\textbf{\small EXAMPLE} \ 
The Cantor set is not a CW space.  The topologist's sine curve $C = A \cup B$, where \\
\vspace{0.05cm}
$
\begin{cases}
\ A = \{(0,y): -1 \leq y \leq 1\}\\
\ B = \{(x,\sin (2\pi /x)): 0 <x \leq 1\}
\end{cases}
, \ 
$
is not a CW space.  
The wedge of the broom is not a CW space but the broom, being contractible, is a CW space, although it carries no CW structure.   
The product $\ds\prod\limits_1^\infty \bS^n$ is not a CW space.\\
\endgroup%%------------------------------------<<

\label{11.16}
\label{19.39}
\begingroup%%----------------------------------->>
\fontsize{9pt}{11pt}\selectfont
\textbf{\small FACT} \ 
Suppose that $X$ is a connected CW space.  
Assume:  $\pi_1(X)$ is finite and $\forall \ q > 1$, $\pi_q(X)$ is finitely generated $-$then there exists a homotopy equivalence $f:K \ra X$, where $K$ is a CW complex such that $\forall \ n$, $K^{(n)}$ is finite.\\
\endgroup%%------------------------------------<<

\begingroup%%----------------------------------->>
\fontsize{9pt}{11pt}\selectfont
Dydak\footnote[2]{\textit{Proc. Amer. Math. Soc.} \textbf{116} (1992), 1171-1173; 
see also Dyer-Roitberg, \textit{Topology Appl.} \textbf{46} (1992), 119-124.} 
has shown that the full subcategory of $\bHCWSP_*$ whose objects are the pointed connected CW spaces is balanced.\\
\endgroup%%------------------------------------<<

Every open subset of a CW complex is a CW space 
(cf. p. \pageref{5.0u}).  Every open subset of a metrizable topological manifold is a CW space 
(cf. p. \pageref{5.0v}).\\

\begin{proposition} \ %06
Let $U$ be an open subset of a normed linear space $E$ $-$then $U$ is a CW space.
\end{proposition}

[Fix a countable neighborhood basis at zero in $E$ consisting of convex balanced sets $U_n$ such that $U_{n+1} \subset U_n$.  
Assuming that $U$ is nonempty, for each $x \in U$, there exists an index $n(x): x + 2U_{n(x)} \subset U$.  
Since $U$ is paracompact, the open covering $\{x + U_{n(x)} : x \in U\}$ 
has a neighborhood finite open refinement $\sO = \{O\}$.  
So, $\forall \ O \in \sO$ $\exists \ x_O \in U: O \subset x_O + U_{n(O)}$ $(n(O) = n(x_O))$.  
Let $\{\kappa_O: O \in \sO\}$ be a partition of unity on $U$ subordinate to $\sO$.
%%----------------------------------------------------------------------------------------------24
Consider $N(\sO)$, the nerve of $\sO$.  If $\{O_1, \ldots, O_k\}$ is a simplex of $N(\sO)$ and if $n(O_1) \leq \cdots \leq n(O_k)$, then the convex hull of $\{x_{O_1}, \ldots, x_{O_k}\}$ is contained in 
$x_{O_1} + 2U_{n(O_1)}$ $\subset$ $U$.  Define continuous functions
$
\begin{cases}
\ f:U \ra \abs{(N(\sO)}\\[-.1cm]
\ g:\abs{(N(\sO)} \ra U
\end{cases}
$
by 
$
\begin{cases}
\ f(x) = \sum\limits_O \kappa_O(x) \chi_O\\[-.1cm]
\ g(\phi) = \sum\limits_O \phi(O) x_O
\end{cases}
$
and put $H(x,t) = tx + (1-t) \sum\limits_O \kappa_O(x) x_O$ to get a homotopy $H:IU \ra U$ between $g \circx f$ and id$_U$.  This shows that $U$ is dominated in homotopy by $\abs{(N(\sO)}$, hence, by the domination theorem, has the homotopy type of a CW complex.]

[Note: \  If $E$ is second countable, then $U$ has the homotopy type of a countable CW complex.  Reason:  Every open covering of a second countable metrizable space has a countable star finite refinement 
(cf. p. \pageref{5.0w}).]\\

\begingroup%%----------------------------------->>
\fontsize{9pt}{11pt}\selectfont
\textbf{\small FACT} \ 
Let $E$ be a normed linear space.  Suppose that $E_0$ is a dense linear subspace of $E$.  Equip $E_0$ with the finite topology $-$then for every open subset $U$ of $E$, the inclusion $U \cap E_0 \ra U$ is a weak homotopy equivalence.\\
\endgroup%%------------------------------------<<

\begingroup%%----------------------------------->>
\fontsize{9pt}{11pt}\selectfont
\textbf{\small FACT} \ 
Let $E$ be a normed linear space.  Suppose that $E^0 \subset E^1 \subset \cdots$ is an increasing sequence of finite dimensional linear subspaces of $E$ whose union is dense in $E$.  Given an open subset $U$ of $E$, put 
$U^n = U \cap E^n$ $-$then $U^0 \subset U^1 \subset \cdots$ is an expanding sequence of topological spaces and the inclusion $U^\infty \ra U$ is a homotopy equivalence.\\
\endgroup%%------------------------------------<<

\begin{proposition} \ %07
Let $A \ra X$ be a closed cofibration and let $f: A \ra Y$ be a continuous function.  
Assume: $A$, $X$, and $Y$ are CW spaces $-$then $X \sqcup_f Y$ is a CW space.
\end{proposition}

[There is a CW pair $(K,L)$ and a commutative diagram
\begin{tikzcd}%[sep=large]
K \arrow{d} &L \arrow{l} \arrow{r}{g} \arrow{d} &Y \arrow[equal]{d}\\
X                &A \arrow{l} \arrow{r}[swap]{f}                 &Y 
\end{tikzcd}
, where the vertical arrows are homotopy equivalences and $g$ is the composite.  Accordingly, $K \sqcup_g Y \approx X \sqcup_f Y$ in \bHTOP (cf. p. \pageref{5.0x} ff.) and $K \sqcup_g Y$ is a CW space (cf. p. \pageref{5.0y}).]\\

Application:  Let 
%\begin{tikzcd}
%X &Z \arrow{l}[swap]{f} \arrow{r}{g} &Y 
%\end{tikzcd}
$X \overset{f}{\lla} Z \overset{g}{\lra} Y$ 
be a 2-source.  Assume: $X$, $Y$, and $Z$ are CW spaces $-$then $M_{f,g}$ is a CW space.

[Note: \  One can establish an analogous result for the double mapping track of a 2-sink in \bCWSP (cf. $\S 6$, Proposition 8).  For example, given a nonempty CW space $X$, $\forall \ x_0 \in X$, $\Omega(X,x_0)$ is a CW space (consider the 2-sink $* \ra X \leftarrow *$).]\\

\begingroup%%----------------------------------->>
\fontsize{9pt}{11pt}\selectfont
\textbf{\small EXAMPLE} \ 
Suppose that $X$ and $Y$ are CW spaces $-$then their join $X * Y$ is a CW space.
\vspi
%%----------------------------------------------------------------------------------------------25
[Note: \  The double mapping cylinder of $X \leftarrow X \times Y \ra Y$ defines the join.  If $X$ and $Y$ are CW  complexes, then $X * Y$ is a CW complex provided that $X \times Y = X \times_k Y$.  Otherwise, consider $X *_k Y$, the double mapping cylinder of $X \leftarrow X \times_k Y \ra Y$.]\\
\endgroup%%------------------------------------<<

\label{14.90}
\textbf{\small LEMMA} \ 
Let $X^0 \subset X^1 \subset \cdots$ be an expanding sequence of topological spaces.  Assume:  $\forall \ n$, $X^n$ is a CW complex containing $X^{n-1}$ as a subcomplex $-$then $X^\infty$ is a CW complex containing $X^n$ as a subcomplex.\\

\begingroup%%----------------------------------->>
\fontsize{9pt}{11pt}\selectfont
\label{5.44b}
\index{mapping telescope (example)}
\textbf{\small EXAMPLE \ (\un{The Mapping Telescope})} \ 
Let 
$
\begin{cases}
\ (\bX,\bff) \\
\ (\bY,\bg) 
\end{cases}
$
be objects in $\bFIL(\bTOP)$.  
Suppose that $\phi:(\bX,\bff) \ra (\bY,\bg)$ is a homotopy morphism, i.e., $\forall \ n$, the diagram
\begin{tikzcd}[sep=large]
X_n \arrow{d}[swap]{\phi_n} \arrow{r}{f_n} &X_{n+1} \arrow{d}{\phi_{n+1}}\\
Y_n \arrow{r}[swap]{g_n} &Y_{n+1}
\end{tikzcd}
is homotopy commutative. $-$then there is an arrow $\tel \phi:\tel(\bX,\bff) \ra \tel(\bY,\bg)$ such that $\forall \ n$, the diagram
\begin{tikzcd}[sep=large]
X_n \arrow{d} &\telsub_n(\bX,\bff) \arrow{l} \arrow{d} \arrow{r} &\tel(\bX,\bff) \arrow{d}\\
Y_n \arrow{r} &\telsub_n(\bY,\bg) \arrow{l} \arrow{r} &\tel(\bY,\bg)
\end{tikzcd}
is homotopy commutative and $\tel \phi$ is a homotopy equivalence if each $\phi_n$ is a homotopy equivalence.  
Thanks to the skeletal approximation theorem and the lemma, it then follows that for any object (\bX,\bff) in 
$\bFIL(\bCW)$, there exists another object (\bX,\bg) in 
$\bFIL(\bCW)$ such that tel(\bX,\bff) and $\tel (\bX,\bg)$ have the same homotopy type and $\tel(\bX,\bg)$ is a CW complex.
\vspi
[The mapping telescope is a double mapping cylinder 
(cf. p. \pageref{5.0z}).  Use the fact that a homotopy morphism of 2-sources, i.e., a homotopy commutative diagram
\begin{tikzcd}[sep=large]
X \arrow{d}                       &Z             \arrow{l}[swap]{f} \arrow{d} \arrow{r}{g}              &Y \arrow{d}\\
X^\prime      &Z^\prime  \arrow{l}{f^\prime}     \arrow{r}[swap]{g^\prime} &Y^\prime
\end{tikzcd}
, gives rise to an arrow $M_{f,g} \ra M_{f^\prime,g^\prime}$ which is a homotopy equivalence if this is the case of the vertical arrows (cf. p. \pageref{5.0aa}).]\\
\endgroup%%------------------------------------<<

\begin{proposition} \ %08
Let $X^0 \subset X^1 \subset \cdots$ be an expanding sequence of topological spaces.  
Assume:  $\forall \ n$, $X^n$ is a CW space and the inclusion $X^n \ra X^{n+1}$ is a cofibration $-$then $X^\infty$ is a CW space.
\end{proposition}

[There is a commutative ladder
\begin{tikzcd}%[sep=large]
K^0  \arrow{d} \arrow{r} &K^1 \arrow{d} \arrow{r} &\cdots\\
X^0  \arrow{r} &X^1 \arrow{r} &\cdots
\end{tikzcd}
, where the vertical arrows $K^n \ra X^n$ are homotopy equivalences and 
$K^0 \subset K^1 \subset \cdots$ is an expanding sequence of CW complexes such that $\forall n$, $(K^n,K^{n-1})$ is a CW pair.  
The induced map $K^\infty \ra X^\infty$ is a homotopy equivalence (cf. $\S 3$, Proposition 15) and, 
by the lemma, $K^\infty$ is a CW complex.]\\

%%----------------------------------------------------------------------------------------------26
Application:  Let (\bX,\bff)  be an object $\bFIL(\bTOP)$.   
Assume:  $\forall \ n$, $X_n$ is a CW space $-$then tel(\bX,\bff)  is a CW space.\\

\begingroup%%----------------------------------->>
\fontsize{9pt}{11pt}\selectfont
\textbf{\small FACT} \ 
Let $X$ be a topological space.  
Suppose that $\sU = \{U_i:i \in I\}$ is a numerable covering of $X$ with the property that for every nonempty finite subset $F \subset I$, $\ds\bigcap\limits_{i \in F} U_i$ is a CW space $-$then $X$ is a CW space.
\vspi
[In the notation of the Segal-Stasheff construction, show that $B\mathcal{U}$ is a CW space.]\\

Application:  Let $X$ be a topological space.  
Suppose that $\sU = \{U_i:i \in I\}$ is a numerable covering of 
$X$ with the property that for every nonempty finite subset $F \subset I$, 
$\ds\bigcap\limits_{i \in F} U_i$ is either empty or contractible $-$then $X$ is a CW space.
\vspi
[Note: \  One can be more precise:  $X$ and $\abs{N(\sU)}$ have the same homotopy type.  
Example:  Every paracompact open subset of a locally convex topological vector space is a CW space (cf. Proposition 6).]\\

\textbf{\small EXAMPLE} \ 
Let $X$ be the Cantor set.  In $\Sigma X$, let $U_1$ be the image of $X \times [0,2/3[$ and let $U_2$ be the image of $X \times ]1/3,1]$ $-$then $\{U_1,U_2\}$ is a numerable covering of $\Sigma X$.  Both $U_1$ and $U_2$ are contractible, hence are CW spaces.  But $\Sigma X$ is not a CW space.  In this connection, observe that $U_1 \cap U_2$ has the same homotopy type as $X$, thus is not a CW space.\\
\endgroup%%------------------------------------<<

A sequence of groups $\pi_n$ $(n \geq 1)$ is said to be a 
\un{homotopy system}
\index{homotopy system} 
if $\forall \ n > 1: \pi_n$ is abelian and there is a left action $\pi_1 \times \pi_n \ra \pi_n$.\\

\index{Theorem, Homotopy System}
\textbf{\small HOMOTOPY SYSTEM THEOREM} \ 
Let $\{\pi_n: n \geq 1\}$ be a homotopy system $-$then there exists a pointed connected CW complex $(X,x_0)$ 
and $\forall \ n \geq 1$, an isomorphism $\pi_n(X,x_0) \ra \pi_n$ such that the action of 
$\pi_1(X,x_0)$ on $\pi_n(X,x_0)$ corresponds to the action of $\pi_1$ on $\pi_n$.

[Note: \   One can take $X$ locally finite if all the $\pi_n$ are countable.]\\

Let $\pi$ be a group and let $n$ be an integer $\geq 1$, 
where $\pi$ is abelian if $n > 1$ $-$then a pointed path connected space $(X,x_0)$ is said to have 
\un{homotopy type $(\pi,n)$} 
\index{homotopy type $(\pi,n)$} 
if $\pi_n(X,x_0)$ is isomorphic to $\pi$ and $\pi_q(X,x_0) = 0$ $(q \neq n)$.  
An \un{Eilenberg-MacLane space}
\index{Eilenberg-MacLane space}of type $(\pi,n)$ 
is a pointed connected CW space $(X,x_0)$ of homotopy type $(\pi,n)$.  
Notation: $(X,x_0) =$ $(K(\pi,n),k_{\pi,n})$.  Two spaces of homotopy type $(\pi,n)$ have the same weak homotopy type and two Eilenberg-MacLane spaces of type $(\pi,n)$ have the same pointed homotopy type.  
Every Eilenberg-MacLane space is nondegenerate, therefore the same is true of its loop space which, moreover, is a pointed CW space (cf. p. \pageref{5.0ab}).  Example:  $\Omega K(\pi,n + 1) = K(\pi,n)$, $\pi$ abelian.\\

%%----------------------------------------------------------------------------------------------27

\begingroup%%----------------------------------->>
\fontsize{9pt}{11pt}\selectfont
\textbf{\small EXAMPLE} \ 
A model for $K(G,1)$, $G$ a discrete topological group, is $B_G^\infty$ (cf. p. \pageref{5.0ac}).\\
\endgroup%%------------------------------------<<

Upon specializing the homotopy system theorem, it follows that for every $\pi$, $(K(\pi,n)$,$k_{\pi,n})$ exists as a pointed CW complex.  
If in addition $\pi$ is abelian, then $(K(\pi,n),k_{\pi,n})$ carries the structure of a homotopy commutative H-group, unique up to homotopy, and the assignment $(X,A) \ra [X,A;K(\pi,n),k_{\pi,n}]$ defines a cofunctor 
$\bTOP^2 \ra \bAB$.\\

\begingroup%%----------------------------------->>
\fontsize{9pt}{11pt}\selectfont
\textbf{\small EXAMPLE} \ 
A model for $K(\Z^n,1)$ is $\bT^n$.
\vspi
[Note: \  Suppose that $X$ is a homotopy commutative H-space with the pointed homotopy type of a finite connected CW complex $-$then 
Hubbuck\footnote[2]{\textit{Topology} \textbf{8} (1969), 119-126.\vspace{0.11cm}}
 has shown that in $\bHTOP_*$, $X \approx \bT^n$ for some $n \geq 0$.]\\
\endgroup%%------------------------------------<<

\begingroup%%----------------------------------->>
\fontsize{9pt}{11pt}\selectfont
\textbf{\small EXAMPLE} \ 
A model for $K(\Z/n\Z,1)$ is the orbit space $\bS^\infty/\Gamma$, where $\Gamma$ is the subgroup of $\bS^1$ generated by a primitive $n^\text{th}$ root of unity.
\vspi
[Note: \  Recall that $\bS^\infty$ is contractible (cf. p. \pageref{5.0ad}).]\\
\endgroup%%------------------------------------<<

\label{5.9}
\begingroup%%----------------------------------->>
\fontsize{9pt}{11pt}\selectfont
\textbf{\small EXAMPLE} \ 
A model for $K(\Q,1)$ is the pointed mapping telescope of the sequence $\bS^1 \ra \bS^1 \ra \cdots$, the $k^{\thx}$ map having degree $k$.
\vspi
[Note: \  
Shelah\footnote[3]{\textit{Proc. Amer. Math. Soc.} \textbf{103} (1988), 627-632.\vspace{0.11cm}}
 has shown that if $X$ is a compact metrizable space which is path connected and locally path connected, then $\pi_1(X)$ cannot be isomorphic to $\Q$.]\\
\endgroup%%------------------------------------<<

\begingroup%%----------------------------------->>
\fontsize{9pt}{11pt}\selectfont
The homotopy type of $\ds\prod\limits_{q=1}^N K(\Z,2q)$ or $\ds\prod\limits_{q=1}^N K(\Z /n\Z,2q)$ admits an interpretation in terms of the theory of algebraic cycles 
(Lawson\footnote[6]{\textit{Ann. of Math.} \textbf{129} (1989), 253-291.}).\\

\endgroup%%------------------------------------<<

\label{5.0an}
$(\pi,1)$\quadx
Suppose that $(X,x_0)$ has homotopy type $(\pi,1)$ $-$then for any pointed connected CW complex $(K,k_0)$, the assignment $[f] \ra f_*$ defines a bijection 
$[K,k_0;X,x_0] \ra$ 
$\Hom(\pi_1(K,k_0),\pi_1(X,x_0))$.  \ \ 
Since \ $(K,k_0)$ \ is wellpointed, the orbit space \ 
$\pi_1(X,x_0)\backslash$ $[K,k_0;X,x_0] $ can be identified with $[K,X]$ 
(cf. p. \pageref{5.0ae}), thus there is a bijection 
$[K,X] \ra \pi_1(X,x_0)\backslash$$\Hom(\pi_1(K,k_0),\pi_1(X,x_0))$, 
the set of conjugacy classes of homomorphisms $\pi_1(K,k_0) \ra \pi_1(X,x_0)$.  
If $\pi$ is abelian, then 
$\text{Hom}(\pi_1(K,k_0),\pi_1(X,x_0)) \approx$ 
$\text{Hom}(H_1(K,k_0),$ $\pi_1(X,x_0)) \approx$ 
$H^1(K,k_0;\pi_1(X,x_0))$ and the forgetful function $[K,k_0;X,x_0]  \ra [K,X]$ is bijective.

Example:  Fix a pointed connected CW complex $(K,k_0)$ $-$then the functor $\bGR\ra \bSET$ that sends $\pi$ to $[K,k_0;K(\pi,1),k_{\pi,1}]$ is represented by $\pi_1(K,k_0)$.\\

%%----------------------------------------------------------------------------------------------28
\begingroup%%----------------------------------->>
\fontsize{9pt}{11pt}\selectfont
\textbf{\small EXAMPLE} \ 
Take $X = K(\pi,1)$, $x_0 = k_{\pi,n}$ and realize $(X,x_0)$ as a pointed CW complex.  
Assume:  $X$ is locally finite and finite dimensional.  
Write $HE(X,x_0)$ $(HE(X))$ for the space of homotopy equivalences of $(X,x_0)$ $(X)$ equipped with the compact open topology $-$then $\pi_0(HE(X,x_0))$ $(\pi_0(HE(x)))$ is the isomorphism group of $(X,x_0)$ $(X)$ viewed as an object in 
$\bHTOP_*$ (\bHTOP).  By the above, $\pi_0(HE(X,x_0)) \approx \Aut\pi$ $(\pi_0(HE(X)) \approx \Out\pi)$.  The evaluation
$
\begin{cases}
\ HE(X) \ra X \\
\ f \mapsto f(x_0) 
\end{cases}
$
is a Hurewicz fibration (cf. $\S 4$ Proposition 6) and its fiber over $x_0$ is $HE(X,x_0)$.  With id$_X$ as the base point, one has $\pi_q(HE(X,x_0),\id_X) = 0$ $(q > 0)$, $\pi_q(HE(X),\id_X) = 0$ $(q > 1)$, and 
$\pi_1(HE(X),\id_X) \approx \Cen\pi$, the center of $\pi$.  The homotopy sequence of the evaluation thus reduces to 
$1 \ra$ 
$\pi_1(HE(X),\id_X)  \ra$ 
$\pi_1(X,x_0) \ra$ 
$\pi_0(HE(X,x_0),\id_X) \ra$ 
$\pi_0(HE(X),\id_X) \ra 1$, 
%\[
%1 \ra \pi_1(HE(X),\id_X)  \ra \pi_1(X,x_0) \ra \pi_0(HE(X,x_0),\id_X) \ra \pi_0(HE(X),\id_X) \ra 1,
%\]
i.e., to
%\[
%1 \ra \Cen \ \pi \ra \pi \ra \Aut \ \pi \ra \Out \ \pi \ra 1.
%\]
%
$1 \ra$ 
$\Cen \pi \ra$ 
$\pi \ra$ 
$\Aut \pi \ra$ 
$\Out \pi \ra 1$.\\
\endgroup%%------------------------------------<<

\begingroup%%----------------------------------->>
\fontsize{9pt}{11pt}\selectfont
\textbf{\small EXAMPLE} \ 
Let $p:X \ra B$ be a Hurewicz fibration, where $B = K(G,1)$.  Suppose that $\forall \ b \in B$, $X_b$ is a $K(\pi,1)$ ($\pi$ abelian) $-$then the only nontrivial part of the homotopy sequence for $p$ is the short exact sequence 
$1 \ra \pi \ra \pi_1(X) \ra G \ra 1$.  Therefore $\pi_1(X)$ is an extension of $\pi$ by $G$ 
and $X$ is a $K(\pi_1(X),1)$ (cf. $\S 6$, Proposition 11).  
Algebraically, there is a left action $G \times \pi \ra \pi$ and geometrically, there is a left action $G \times \pi \ra \pi$.  These two actions are identical.\\
\endgroup%%------------------------------------<< 

\label{5.43}
\begingroup%%----------------------------------->>
\fontsize{9pt}{11pt}\selectfont
\textbf{\small EXAMPLE} \ 
Consider a 2-source $\pi^\prime \leftarrow G \ra \pi^{\prime\prime}$ in \bGR, where the arrows are monomorphisms.  
Define $\pi$ by the pushout square
\begin{tikzcd}[sep=large]
G \arrow{d} \arrow{r} &\pi^{\prime\prime}  \arrow{d}\\
\pi^\prime \arrow{r}   &\pi
\end{tikzcd}
, i.e., $\pi = \pi^\prime *_G \pi^{\prime\prime}$ 
$-$then there exists a pointed CW complex $X = K(\pi,1)$ and pointed subcomplexes
$
\begin{cases}
\ X^\prime = K(\pi^\prime,1) \\
\ X^{\prime\prime} =  K(\pi^{\prime\prime},1)
\end{cases}
, \ 
$
$Y = K(G,1)$ such that $X = X^\prime \cup X^{\prime\prime}$ and $Y = X^\prime \cap X^{\prime\prime}$ .\\
\endgroup%%------------------------------------<< 

\label{14.56}
\begingroup%%----------------------------------->>
\fontsize{9pt}{11pt}\selectfont
\textbf{\small EXAMPLE} \ 
Let $X$ and $Y$ be connected CW complexes.  
Suppose that $f:X \ra Y$ is a continuous function such that for every finite connected CW complex $K$, 
the induced map $[K,X] \ra [K,Y]$ is bijective 
$-$then $f$ is a homotopy equivalence iff $\forall \ x \in X$, $f_*:\pi_1(X,x) \ra \pi_1(Y,f(x))$ is surjective 
(cf. p. \pageref{5.0af}) but this condition is not automatic.  
To construct an example, 
let $S_\infty$ be the subgroup of the symmetric group of $\N$ consisting of those permutations that have finite support.  
Each injection $\iota:\N \ra \N$ determines a homomorphism $\iota_\infty:S_\infty \ra S_\infty$ viz.
$
\begin{cases}
\ \restr{\iota_\infty(\sigma)}{(\N - \iota(\N))} = \id\\
\ \restr{\iota_\infty(\sigma)}{\iota(\N)} = \iota \circx \sigma \circx \iota^{-1}
\end{cases}
, \ 
$
and on any finite product, 
$\ds\prod \iota_\infty: S_\infty\backslash \ds\prod S_\infty \ra S_\infty\backslash \prod S_\infty$ 
is bijective. 
Here the action of $S_\infty$ on $\ds\prod S_\infty$ is by conjugation.  
Choose $\phi:K(S_\infty,1) \ra K(S_\infty,1)$ such that $\phi_* = \iota_\infty$ on $S_\infty$ $-$then for every finite connected CW complex $K$, the induced map 
$[K,K(S_\infty,1)] \ra [K,K(S_\infty,1)]$ is bijective (consider first a finite wedge of circles).  
However, $\phi$ is not a homotopy equivalence unless $\iota$ is surjective.
\vspi
[Note: \  There are various conditions on $\pi_1(X)$ (or $\pi_1(Y)$) which guarantee that $f_*$ is surjective (under the given assumptions).  For example, any of the following will do:  (1) $\pi_1(X)$ (or $\pi_1(Y)$) nilpotent; (2) $\pi_1(X)$ (or $\pi_1(Y)$) finitely generated; (3) $\pi_1(X)$ (or $\pi_1(Y)$) free.]\\
\endgroup%%------------------------------------<< 

%%----------------------------------------------------------------------------------------------29

\begingroup%%----------------------------------->>
\fontsize{9pt}{11pt}\selectfont
\textbf{\small EXAMPLE} \ 
Let $\pi$ be a group $-$then $K(\pi,1)$ can be realized by a path connected metrizable topological manifold 
(cf. p. \pageref{5.0ag}) iff $\pi$ is countable and has finite cohomological dimension 
(Johnson\footnote[2]{\textit{Proc. Camb. Phil. Soc.} \textbf{70} (1971), 387-393.}).
\vspi
[Note: \  Under these circumstances, the cohomological dimension of $\pi$ cannot exceed the euclidean dimension of $K(\pi,1)$, there being equality iff $K(\pi,1)$ is compact.]\\
\endgroup%%------------------------------------<< 

\begingroup%%----------------------------------->>
\fontsize{9pt}{11pt}\selectfont
\textbf{\small EXAMPLE} \ 
The homotopy type of an aspherical compact topological manifold is completely determined by its fundamental group.  
Question:  If $X$ and $Y$ are aspherical compact topological manifolds and if $\pi_1(X) \approx \pi_1(Y)$, 
is it then true that $X$ and $Y$ are homeomorphic?  Borel has conjectured that the answer is ``yes''.  To get an idea of the difficulty of this problem, a positive resolution easily leads to a proof of the Poincar\'e conjecture (modulo a result of Milnor).  Additional information and references can be found in 
Farrell-Jones\footnote[3]{\textit{CBMS Regional Conference} \textbf{75} (1990), 1-54; 
see also Conner-Raymond, \textit{Bull. Amer. Math. Soc.} \textbf{83} (1977), 36-85.}.\\
\endgroup%%------------------------------------<< 

\label{5.31}
$(\pi,n)$ \quadx
Suppose that $(X,x_0)$ has the homotopy type $(\pi,n)$, where $\pi$ is abelian.  
Let $\iota \in H^n(X,x_0;\pi_n(X,x_0))$ be the fundamental class 
$-$then for any pointed connected CW complex $(K,k_0)$, the assignment $[f] \ra f^* \iota$ defines a bijection 
$[K,k_0;X,x_0] \ra$ 
$H^n(K,k_0;$ $\pi_n(X,x_0))$.

Assuming that 
$\pi^\prime$ and  $\pi^{\prime\prime}$ are abelian, 
$[K(\pi^\prime,n),k_{\pi^\prime,n};K(\pi^{\prime\prime},n),k_{\pi^{\prime\prime},n}] \approx$ 
$[K(\pi^\prime,n),$ $K(\pi^{\prime\prime},n)] \approx$ 
$\text{Hom}(\pi^\prime,\pi^{\prime\prime})$.  
Example: Suppose that $0 \ra \pi^\prime \ra \pi \ra \pi^{\prime\prime} \ra 0$ is a short exact sequence of abelian groups $-$then
(1) The mapping fiber of the arrow $K(\pi,n) \ra K(\pi^{\prime\prime},n)$ is a $K(\pi^\prime,n)$; 
(2)  The mapping fiber of the arrow $K(\pi^\prime,n+1) \ra K(\pi,n+1)$ is a $K(\pi^{\prime\prime},n)$; 
(3) The mapping fiber of the arrow $K(\pi^{\prime\prime},n) \ra K(\pi^\prime,n+1)$ is a $K(\pi,n)$.

[Note: \  $\bCWSP_*$ is closed under the formation of mapping fibers (cf. $\S 6$, Proposition 8).]\\

\begingroup%%----------------------------------->>
\fontsize{9pt}{11pt}\selectfont
\textbf{\small EXAMPLE} \ 
A model for $K(\Z,2)$ is $\bP^\infty(\C)$.  Fix $n > 1$ and choose a map $\bP^\infty(\C) \ra K(\Z,2n)$ representing a generator of $H^{2n}(\bP^\infty(\C);\Z) \approx \Z$.  Put $Y = \bP^\infty(\C)$ and define $X$ by the pullback square
\begin{tikzcd}[sep=large]
X \arrow{d} \arrow{r} &{\Theta K(\Z,2n)} \arrow{d}\\
Y \arrow{r}   &K(\Z,2n)
\end{tikzcd}
.  The fiber $X_{y_0}$ is a $K(\Z,2n - 1)$.  Since $2n - 1 \geq 3$, there is an isomorphism $\pi_{2n-1}(X_{y_0}) \approx \pi_{2n-1}(X)$ but the corresponding arrow in homology $H_{2n-1}(X_{y_0}) \ra H_{2n-1}(X)$ is not even one-to-one.\\
\endgroup%%------------------------------------<< 

Let $(X,A)$ be a relative CW complex $-$then for any abelian group $\pi$, there is a bijection $[X,A;K(\pi,n),k_{\pi,n}] \ra H^n(X,A;\pi)$ which, in fact, is an isomorphism of abelian groups,
%%----------------------------------------------------------------------------------------------30
natural in $(X,A)$.  This applies in particular when $A = \emptyset$, thus there is an isomorphism 
$[X;K(\pi,n)] \ra H^n(X;\pi)$ of abelian groups, natural in $X$.  So, on \bHCW the cofunctor $H^n(-;\pi)$ is representable by $K(\pi,n)$.  But on \bHTOP itself, this is no longer true in that the relation $[X,K(\pi,n)] \approx H^n(X;\pi)$ can fail if $X$ is not a CW complex.\\

\begingroup%%----------------------------------->>
\fontsize{9pt}{11pt}\selectfont
\textbf{\small EXAMPLE} \ 
Let $X$ be the Warsaw circle and take $\pi = \Z$ $-$then $H^1(X,\Z) = 0$, while $[X,K(\Z,1)] \approx \Z$ or still, 
$[X,K(\Z,1)] \approx$ $\check{H}^1(X;\Z)$.\\
\endgroup%%------------------------------------<< 

\label{20.1}
In general, for an arbitrary abelian group $\pi$ and an arbitrary pair $(X,A$), there is a natural isomorphism $[X,A;K(\pi,n),k_{\pi,n}] \ra$ $\check{H}(X,A;\pi)$ 
(cf. p. \pageref{5.0ah}).  
Moral: It is $\check{C}$ech cohomology rather than singular cohomology that is the representable theory.\\

\label{6.30}
\begingroup%%----------------------------------->>
\fontsize{9pt}{11pt}\selectfont
Suppose that $(X,x_0)$ is a pointed connected CW complex.  
Equip $C(X,K(\pi,n))$ with the compact open topology $-$then $[X,K(\pi,n)] = \pi_0(C(X,K(\pi,n)))$, $X$ being a compactly generated Hausdorff space.  
Because the forgetful function $[X,x_0;K(\pi,n),k_{\pi,n}] \ra [X,K(\pi,n)]$ is surjective, every path component of $C(X,K(\pi,n))$ contains a pointed map $f_0: f_0(x_0) = k_{\pi,n}$.\\
\endgroup%%------------------------------------<< 

\label{9.29}
\begingroup%%----------------------------------->>
\fontsize{9pt}{11pt}\selectfont
\textbf{\small EXAMPLE} \ 
Let $(X,x_0)$ be a pointed connected CW complex.  
Assume: $X$ is locally finite $-$then for any abelian group $\pi$, 
$\pi_q(C(X,K(\pi,n)),f_0) \approx$
$
\begin{cases}
\ H^{n-q}(X;\pi) \ \ \  (1 \leq q \leq n)\\
\ 0 \indent\indent\quadx  (q > n)
\end{cases}
.
$
\vspi
[Since $K(\pi,n)$ is an H-group, all the path components of $C(X,K(\pi,n))$ have the same homotopy type.  
Let $f_0$ be the constant map $X \ra k_{\pi,n}$, $C_0(X,K(\pi,n))$ its path component.  
To compute 
$\pi_q(C_0(X,K(\pi,n)),f_0)$, consider the  Hurewicz fibration 
$C_0(X,K(\pi,n)) \ra K(\pi,n)$ which sends $f$ to $f(x_0)$ (cf. $\S 4$, Proposition 6), bearing in mind that 
$\pi_1(C_0(X,K(\pi,n)),f_0)$ is abelian.]
\vspi
[Note: \ Suppose in addition that $X$ is finite $-$then $C(X,K(\pi,n))$ (compact open topology) is a CW space 
(cf. p. \pageref{5.0aha}) \ 
and there is a decomposition 
$H^n(C(X,K(\pi,n)) \times X;\pi) \approx$ 
$\ds \bigoplus\limits_{q = 0}^n H^q (C(X,K(\pi,n));$ $H^{n-q}(X;\pi)).$  
Let 
$\ev:C(X,K(\pi,n)) \times X \ra K(\pi,n)$ be the evaluation.  
Take the fundamental class 
$\iota \in H^n(K(\pi,n);\pi)$ and write 
$\ev^* \iota = \ds \bigoplus\limits_{q = 0}^n \mu_q$, 
where
$\mu_q \in H^q(C(X,K(\pi,n)); H^{n-q}(X;\pi))$.  
Let 
$[f_q] \in [C(X,K(\pi,n)),K(H^{n-q}(X;\pi),q)]$ 
correspond to $\mu_q$ (conventionally, 
$K(H^n(X;\pi),0)$ is $H^n(X;\pi)$ (discrete topology)).  
The $f_q$ determine an arrow 
$C(X,K(\pi,n)) \ra \ds \prod\limits_{q=0}^n K(H^{n-q}(X;\pi),q)$.  
It is a weak homotopy equivalence, hence, by the realization theorem, a homotopy equivalence.]\\
\endgroup%%------------------------------------<< 

\begingroup%%----------------------------------->>
\fontsize{9pt}{11pt}\selectfont
\textbf{\small EXAMPLE} \ 
Let $(X,x_0)$ be a pointed connected CW complex.  
Assume: $X$ is locally finite and finite dimensional $-$then for any group $\pi$, $\pi_q(C(X,K(\pi,1)),f_0) \approx$
$
\begin{cases}
\ \text{Cen}(\pi,f_0) \ \ \ (q = 1)\\
\ 0 \indent\indent \  (q > 1)
\end{cases}
. \ 
$
Here, $\Cen(\pi,f_0)$ is the centralizer of $(f_0)_*(\pi_1(X,x_0))$ in $\pi_1(K(\pi,1),k_{\pi,1}) \approx \pi$.  Special case:  Suppose that $(X,x_0)$  is aspherical, let $\pi = \pi_1(X,x_0)$, take $f_0 = \id_X$, and conclude that the path component of the identity in 
%%----------------------------------------------------------------------------------------------31
$C(X,X)$ has homotopy type $(\Cen\pi,1)$, $\Cen\pi$ the center of $\pi$.  Example:  $\Cen\pi$ is trivial if $X$ is a compact connected riemannian manifold whose sectional curvatures are $< 0$.
\vspi
[Reduce to when $X^{(0)} = \{x_0\}$ 
(cf. p. \pageref{5.0ai}), 
observe that 
$\pi_q(C(X,K(\pi,1)),f_0) \approx$ 
$\pi_q(C(X^{(1)},K(\pi,1)),$ $\restr{f_0}{X^{(1)}})$, 
and use the fact that $X^{(1)}$ is a wedge of circles.]
\vspi
[Note: \  It can happen that $\pi$ is finitely generated but  $\Cen(\pi,f_0)$ is infinitely generated even if $X = \bS^1$ (Hansen\footnote[2]{\textit{Compositio Math.} \textbf{28} (1974), 33-36.}).]\\
\endgroup%%------------------------------------<<

\label{9.86}
A 
\un{compactly generated group}
\index{compactly generated group} 
is a group $G$ equipped with a compactly generated topology in which inversion $G \ra G$ is continuous and multiplication 
$G \times_k G \ra G$ is continuous.  
Since multiplication is not required to be continuous on $G \times G$ (product topology), 
a compactly generated group is not necessarily a topological group, 
although this will be the case if $G$ is a LCH space or if $G$ is first countable.  
Example:  Let $G$ be a simplicial group $-$then its geometric realization $\abs{G}$ is a compactly generated group 
(cf. p. \pageref{5.0aj}).

[Note: \  If $G$ is a topological group, then $kG$ is a compactly generated group but $kG$ need not be a topological group 
(cf. p. \pageref{5.0ak}).  
A compactly generated group is $\tT_0$ iff it is \dsp.  
Therefore any \dsp compactly generated group which is not Hausdorff cannot be a topological group.]

Suppose that $\pi$ is abelian 
$-$then it is always possible to realize $K(\pi,n)$ as a pointed CW complex carrying the structure 
of an abelian compactly generated group on which $\Aut\pi$ operates to the right by 
base point preserving skeletal homeomorphisms such that $\forall \ \phi \in \Aut\pi$, there is a commutative square
\begin{tikzcd}%[sep=large]
\pi_n(K(\pi,n)) \arrow{d}[swap]{\phi_*} \indent\approx &\pi  \arrow{d}{\phi}\\
\pi_n(K(\pi,n)) \indent\approx   &\pi
\end{tikzcd}
(Adem-Milgram\footnote[3]{\textit{Cohomology of Finite Groups}, Springer Verlag (1994), 51.}) $(0 = k_{\pi,n})$.   
With this understanding, let $G$ be a group, assume that $\pi$ is a right $G$-module, 
and denote by $\chi:G \ra \Aut\pi$ the associated homomorphism.  
Calling $\widetilde{K}(G,1)$ the universal covering space of $K(G,1)$, 
form the product $\widetilde{K}(G,1) \times K(\pi,n)$ and write 
$K(\pi,n;\chi)$ for the orbit space $(\widetilde{K}(G,1) \times K(\pi,n))/G$.  
As an object in $\bTOP/(K(G,1)$, $K(\pi,n;\chi)$  is locally trivial with fiber $K(\pi,n)$, thus the projection 
$p_\chi:K(\pi,n;\chi) \ra$ 
$K(G,1)$ is a Hurewicz fibration (local-global principle) and $K(\pi,n;\chi)$ is a CW space (cf. $\S 6$, Proposition 11).  
The inclusion $\widetilde{K}(G,1) \times \{0\} \ra \widetilde{K}(G,1) \times K(\pi,n)$ 
defines a section $s_\chi:K(G,1) \ra K(\pi,n;\chi)$, so $K(\pi,n;\chi)$ is an object in $\bTOP(K(G,1))$ 
(cf. p. \pageref{5.0al}).  
Example:  Take  $G = \Aut\pi:\ $ 
$
\begin{cases}
\ \pi \times \Aut\pi \ra \pi \\
\ (\alpha, \phi) \mapsto \phi^{-1}(\alpha)
\end{cases}
$
$-$then the associated homomorphism $\Aut\pi \ra \Aut\pi$ is $\id_{\Aut\pi} \equiv \chi_{_\pi}$.

%%----------------------------------------------------------------------------------------------32
[Note: \  Given $G$, consider the trivial action $\pi \times G \ra \pi$, where $\chi: \ $
$
\begin{cases}
\ G \ra \Aut\pi\\
\ g \mapsto \id_\pi
\end{cases}
\hspace{-.25cm}.
$ 
In this case, $K(\pi,n;\chi)$  reduces to the product $K(G,1) \times K(\pi,n)$.]

Example:  Take $\pi = \Z$, $G = \Z/2\Z$ and let $\chi: G \ra \Aut\pi$  be the nontrivial homomorphism $-$then $K(\Z,2;\chi)$ ``is'' $B_{\bO(2)}$.\\

\label{11.7}
\begingroup%%----------------------------------->>
\fontsize{9pt}{11pt}\selectfont
\label{5.47}
\textbf{\small EXAMPLE} \ 
The homotopy sequence for $p_\chi$ breaks up into a collection of split short exact sequences $0 \ra \pi_q(K(\pi,n)) \ra \pi_q(K(\pi,n;\chi)) \ra \pi_q(K(G,1))  \ra 0$.  
Case 1:  $n \geq 2$.  Here 
$\pi_q(K(\pi,n;\chi)) \approx $
$
\begin{cases}
\ \pi \quadx (q = n) \\
\ G  \quadx (q = 1)
\end{cases}
$
and $\pi_q(K(\pi,n;\chi)) = 0$ otherwise.  The algebraic right action 
$\pi \times G \ra \pi$ corresponds to an algebraic left action 
$G \times \pi \ra \pi$ and this is the same as the geometric left action 
$G \times \pi \ra \pi$
Case 2:  $n = 1$.  In this situation, $\pi_1(K(\pi,n;\chi))$ is a split extension of $\pi$ by $G$ and the higher homotopy groups are trivial.  If $\Theta_{s,p}K(\pi,n;\chi)$ is the subspace of $PK(\pi,n;\chi)$ made up of those $\sigma$ such that $\sigma (0) \in s_\chi(K(G,1))$ and $p_\chi(\sigma(t)) = p_\chi(\sigma(0))$ $(0 \leq t \leq 1)$, then the projection $\Theta_{s,p}K(\pi,n;\chi) \ra K(\pi,n;\chi)$ sending $\sigma$ to $\sigma(1)$ is a Hurewicz fibration whose fiber over the base point is $\Omega K(\pi,n)$.  Specialize and take 
$G = \Aut\pi$ (so $\chi = \chi_\pi$).  
Let $B$ be a connected CW complex.  
The ``class'' of fiber homotopy classes of Hurewicz fibrations $X \ra B$ with fiber $K(\pi,n)$ is a ``set'' 
(cf. p. \pageref{5.0am} ff.).  
As such, it is in a one-to-one correspondence with the set of homotopy classes 
$[B,K(\pi,n+1;\chi_\pi)]:[X] \leftrightarrow [\Phi]$, 
$\Phi:B \ra K(\pi,n+1;\chi_\pi)$ the 
\un{classifying map}, 
\index{classifying map} 
where $X$ is define by the pullback square
\begin{tikzcd}[sep=large]
X\arrow{d} \arrow{r} &{\Theta_{s,p}K(\pi,n+1;\chi_\pi)}  \arrow{d}\\
B \arrow{r}   &{K(\pi,n+1;\chi_\pi)}
\end{tikzcd}
.  
For example, if $X$ is a connected CW space with two nonzero homotopy groups 
$\pi_1(X) = G$ and $\pi_n(X) = \pi$ $(n > 1)$, then the geometry furnishes a right action 
$\pi \times G \ra \pi$ and an associated homomorphism 
$\chi:G \ra \Aut\pi$.  
To construct $X$ up to homotopy, fix a map 
$f:X \ra K(G,1)$ 
which induces the identity on $G$, pass to the mapping track $W_f$, 
and consider the Hurewicz fibration $W_f \ra K(G,1)$.  
There is an arrow 
$\Phi: K(G,1) \ra K(\pi,n+1;\chi_\pi)$ such that 
$\chi = \Phi_*:G \ra \Aut\pi$ and $[W_f] \leftrightarrow [\Phi]$.
\vspi
[Note: \  Suppose that $B$ is a pointed simply connected CW complex $-$then the set of fiber homotopy classes of Hurewicz fibrations $X \ra B$ with fiber $K(\pi,n)$ is in a one-to-one correspondence with 
$\Aut\pi\backslash H^{n+1}(B;\pi)$.  
Proof:  The set of homotopy classes $[B,K(\pi,n+1;\chi_\pi)]$ can be identified with the set of pointed homotopy classes $[B,K(\pi,n+1;\chi_\pi)] \mod \pi_1(K(\pi,n+1,\chi_\pi))$, i.e., with the set of pointed homotopy classes 
$[B,K(\pi,n+1;\chi_\pi)] \mod \Aut\pi$, i.e., with the set of homotopy classes $[B,K(\pi,n+1)] \mod \Aut\pi$ 
(cf. p. \pageref{5.0ama}), i.e., with 
 $\Aut\pi \backslash H^{n+1}(B;\pi)$.  Translated, this means that in the simply connected case, one can use
\begin{tikzcd}%[sep=large]
\Theta K(\pi,n+1) \arrow{d} \\
K(\pi,n+1)
\end{tikzcd}
to carry out the classification but then it is also necessary to build in the action of $\Aut\pi$.]\\
\endgroup%%------------------------------------<<

\begingroup%%----------------------------------->>
\fontsize{9pt}{11pt}\selectfont
\textbf{\small EXAMPLE} \ 
Let $G$ be a group; let 
$
\begin{cases}
\ \chi^\prime: G \ra \Aut\pi^\prime \\
\ \chi^{\prime\prime}: G \ra \Aut\pi^{\prime\prime}
\end{cases}
$
be homomorphisms, where
$
\begin{cases}
\ \pi^\prime\\
\ \pi^{\prime\prime}
\end{cases}
$
are abelian $-$then $[K(\pi^\prime,n+1;\chi^\prime),K(\pi^{\prime\prime},n+1;\chi^{\prime\prime})]_G$ $\approx$ $\text{Hom}_G(\pi^\prime,\pi^{\prime\prime})$, $[\ ,\ ]_G$ standing for homotopy in
%%----------------------------------------------------------------------------------------------33
$\bTOP(K(G,1))$.\\
\endgroup%%------------------------------------<<

Notation:  
Given $X$ in \bTOP/\mB and $\phi \in C(E,B)$, let $\text{lif}_\phi(E,X)$ be the set of liftings $\Phi:E \ra X$ of $\phi$.  
Relative to a choice of base points $b_0 \in B$ , $x_0 \in X_{b_0}$, and $e_0 \in E$, where $\phi(e_0) = b_0$, let $\text{lif}_\phi(E,e_0;X,x_0)$ be the subset of $\text{lif}_\phi(E,X)$ consisting of those $\Phi$ such that $\Phi(e_0) = x_0$.  
Write $[E,X]_\phi$ for the set of fiber homotopy classes in $\text{lif}_\phi(E,X)$ and $[E,e_0;X,x_0]_\phi$ for the set of pointed fiber homotopy classes in $\text{lif}_\phi(E,e_0;X,x_0)$.\\

\textbf{\small LEMMA} \ 
If $(B,b_0)$, $(E,e_0)$ are wellpointed with $\{b_0\} \subset B$, $\{e_0\} \subset E$ closed, then the fundamental group $\pi_1(X_{b_0},x_0)$ operates to the left on $[E,e_0;X,x_0]_\phi$ and the forgetful function $[E,e_0;X,x_0]_\phi \ra [E,X]_\phi$ passes to the quotient to define an injection $\pi_1(X_{b_0},x_0)\backslash [E,e_0;X,x_0]_\phi \ra [E,X]_\phi$ which, when $X_{b_0}$ is path connected, is a bijection.\\

\begingroup%%----------------------------------->>
\fontsize{9pt}{11pt}\selectfont
Let $G$ and $\pi$ be groups.  Given $\chi \in \text{Hom}(G,\Aut\pi)$, denote by $\text{Hom}_\chi(G,\pi)$ the set of 
\un{crossed} 
\un{homomorphisms}
\index{crossed homomorphisms}
per $\chi$, so $f:G \ra \pi$ is in $\text{Hom}_\chi(G,\pi)$ iff $f(g^\prime,g^{\prime\prime}) = f(g^\prime)(\chi(g^\prime)f(g^{\prime\prime}))$.  
There is a left action $\pi \times \text{Hom}_\chi(G,\pi) \ra \text{Hom}_\chi(G,\pi)$, viz. $(\alpha \cdot f)(g) = \alpha f(g)(\chi(g)\alpha^{-1})$.
\vspi
[Note: \ The elements of $\Hom_\chi(G,\pi)$ correspond bijectively to the sections 
$s:G \ra \pi \rtimes_\chi G$, where $\pi \rtimes_\chi G$ is the semidirect product (cf. p. \pageref{5.0amb}).]\\
\endgroup%%------------------------------------<<

\begingroup%%----------------------------------->>
\fontsize{9pt}{11pt}\selectfont
\textbf{\small EXAMPLE} \ 
Suppose that $B$ is a connected CW complex.  
Fix a group $\pi$ and a Hurewicz fibration $p:X \ra B$ with fiber $K(\pi,1)$.  
Assume: $\text{sec}_B(X) \neq \emptyset$, say $s \in \text{sec}_B(X)$.  
Choose $b_0 \in B$ and put $x_0 = s(b_0)$.  Let $(E,e_0)$ be a pointed connected CW complex, 
$\phi:E \ra B$ a pointed continuous function.  
There is a split short exact sequence $1 \ra \pi_1(X_{b_0},x_0) \ra \pi_1(X,x_0) \ra\pi_1(B,b_0) \ra 1$, 
from which a left action of $G = \pi_1(E,e_0)$ on $\pi = \pi_1(X_{b_0},x_0)$ or still, 
a homomorphism $\chi:G \ra \Aut\pi$, $\chi(g)$ thus being conjugation by $(s \circx \phi)_*(g)$.  
Attach to $\Phi \in \text{lif}_\phi(E,e_0;X,x_0)$ an element $f_\Phi \in \text{Hom}_\chi(G,\pi)$ 
via the prescription $f_\Phi(g) = \Phi_*(g)(s \circx \phi)_*(g)^{-1}$ $-$then the assignment $\Phi \ra f_\Phi$ induces a bijection 
$[E,e_0;X,x_0]_\phi \ra$ 
$\text{Hom}_\chi(G,\pi)$, so $[E,X]_\phi \approx$ 
$\pi\backslash [E,e_0;X,x_0]_\phi \approx$ 
$\pi\backslash \text{Hom}_\chi(G,\pi)$.
\vspi
[Note: \  The considerations on 
p. \pageref{5.0an} 
are recovered by taking $B = *$ and $X = K(\pi,1)$.]\\
\endgroup%%------------------------------------<<

\index{locally constant coefficients}
\indent\indent (Locally Constant Coefficients)  
Let $(X,x_0)$ be a pointed connected CW complex.  Assume given a homomorphism 
$\chi_\sG:\pi_1(X,x_0) \ra G$ and a homomorphism 
$\chi:G \ra \Aut\pi$, where $\pi$ is abelian.  Let 
$\sG:\Pi X \ra \bAB$ be the cofunctor determined by the composite 
$\chi \circx \chi_\sG$ 
(cf. p. \pageref{5.0ao}).  
Choose a pointed continuous function $f_\sG:X \ra K(G,1)$ corresponding to $\chi_\sG$ and put $k_{\pi,n;\chi} = s_\chi(k_{G,1})$ $-$then $[X,x_0;K(\pi,n;\chi),k_{\pi,n;\chi}]_{f_\sG} \approx H^n(X,x_0;\sG)$.  
So, if $n = 1$, 
$H^1(X,x_0;\sG) \approx \text{Hom}_{\chi \circx \chi_\sG}(\pi_1(X,x_0),\pi)$ (see the preceding example) 
$\implies$ $H^1(X;\sG)$ $\approx$ 
$\pi\backslash H^1(X,x_0;\sG)$ $\approx$ 
$\pi\backslash \text{Hom}_{\chi \circx \chi_\sG}(\pi_1(X,x_0),\pi)$ $\approx$ 
$[X;K(\pi,1;\chi)]_{f_\sG}$ 
but if $n > 1$, 
$H^n(X,x_0;\sG)$ $\approx $ 
$H^n(X;\sG)$ $\approx $ 
$[X,K(\pi,n;\chi)]_{f_\sG}$.

%%----------------------------------------------------------------------------------------------34
\label{9.2} %dmc mnft
\label{9.45} %dmc mnft
\label{9.51} %dmc mnft
[Note: \  The cohomology of any cofunctor $\sG:\Pi X \ra \bAB$ fits into this scheme.  
Simply take $\pi = \sG_{x_0}$, $G = \Aut\pi$, $\chi = \chi_{_\pi}$, 
and let $\chi_{_\sG}:\pi_1(X,x_0) \ra \Aut\pi$ be the homomorphism derived from the right action 
$\pi \times \pi_1(X,x_0) \ra \pi$ (of course, $H^0(X,\sG)$ is $\text{fix}_{\chi_\sG}(\pi)$, 
the subgroup of $\pi$ whose elements are fixed by $\chi_\sG$).  
When $\chi_\sG$ is trivial, one can choose $f_\sG$ as the map to the base point of 
$K(\Aut\pi,1)$ and recover the fact that $[X,K(\pi,n)] \approx H^n(X;\pi)$.]\\

\textbf{\small LEMMA} \ 
Fix a set of representatives $f_i$ for 
$[X,x_0;K(G,1),k_{G,1}]$ $-$then 
$[X,x_0;K(\pi,n;\chi),$ $k_{\pi,n;\chi}]$ 
is in a one-to-one correspondence with the union 
$\bigcup\limits_i [X,x_0;K(\pi,n;\chi),k_{\pi,n;\chi}]_{f_i}$ (which is necessarily disjoint).\\

\label{5.27a}
Application:  There is a one-to-one correspondence between the set of pointed homotopy classes of pointed continuous functions 
$f:X \ra K(\pi,n;\chi)$ such that $\pi_1(f) = \chi_\sG$ and the elements of $H^n(X;\sG)$ $(n > 1)$.\\

\begingroup%%----------------------------------->>
\fontsize{9pt}{11pt}\selectfont
\textbf{\small FACT} \ 
Let
$
\begin{cases}
\ (X,x_0) \\
\ (Y,y_0)
\end{cases}
$
be pointed connected CW complexes; let $f \in C(X,x_0;Y,y_0)$.  Assume given a homomorphism $\chi_\sG:\pi_1(Y,y_0) \ra G$ and a homomorphism $\chi:G \ra \Aut\pi$.  Put $\chi_{f^*\sG} = \chi_\sG \circx \pi_1(f)$ and suppose that $f^*:[Y,y_0;K(\pi,n;\chi),k_{\pi,n;\chi}] \ra [X,x_0;K(\pi,n;\chi),k_{\pi,n;\chi}]$ is bijective $-$then $H^n(Y;\sG) \approx H^n(X;f^*\sG)$.\\
\endgroup%%------------------------------------<< 

The singular homology and cohomology groups of an Eilenberg-MacLane space of type $(\pi,n)$ with coefficients in $G$ depend only on $(\pi,n)$ and $G$.  
Notation:  $H_q(\pi,n;G)$, $H^q(\pi,n;G)$
(or $H_q(\pi,n)$, $H^q(\pi,n)$ if $G = \Z$).  
Example:  $H_n(\pi,n) \approx \pi/[\pi,\pi]$.

[Note: \  There are isomorphisms $H_*\pi \approx H_*(\pi,1)$ ($H^*\pi \approx H^*(\pi,1)$), where $H_*\pi$ ($H^*\pi$) is the homology (cohomology) of $\pi$.  
In general, if $G$ is a right $\pi$-module and if $\sG$ is the locally constant coefficient system on $K(\pi,1)$ associated with $G$, then $H_*(\pi,G)$ ($H^*(\pi,G)$) is isomorphic to 
$H_*(K(\pi,1);\sG)$ ($H^*(K(\pi,1);\sG)$).]\\

\begingroup%%----------------------------------->>
\fontsize{9pt}{11pt}\selectfont
\textbf{\small EXAMPLE} \ 
If $\pi$ is abelian, then $\forall \ n \geq 2$, $H_{n+1}(\pi,n) = 0$ but this can fail if $n = 1$ since, e.g., 
$H_2(\Z/2\Z \oplus \Z/2\Z, 1) \approx H_1(\Z/2\Z, 1) \otimes H_1(\Z/2\Z, 1) \approx \Z/2\Z$.  
When does $H_2(\pi,1)$ vanish?  
To formulate the answer, let $0 \ra \pi_\text{tor} \ra \pi \ra \Pi \ra 0$ be the short exact sequence in which 
$\pi_\text{tor}$ is the torsion subgroup of $\pi$ and denote by 
$\pi_\text{tor}(p)$ the $p$-primary component of $\pi_\text{tor}$ $-$then 
Varadarajan\footnote[2]{\textit{Ann. of Math.} \textbf{84} (1966), 368-371.} 
has shown that $H_2(\pi,1) = 0$ iff rank 
$\Pi \leq 1$ plus $\forall \ p: (p_1)$ 
$\pi_\text{tor}(p) \otimes \Pi = 0$ $\&$ $(p_2)$ 
$\pi_\text{tor}(p)$ is the direct sum of a divisible group and a cyclic group.  
Example:  Assume that $\pi$ is finite $-$then $H_2(\pi,1) = 0$ iff $\pi$ is cyclic.  Other examples include $\pi = \Z$, $\pi = \Q$, and $\pi = \Z/p^\infty \Z$ (the $p$-primary component of $\Q/\Z$).\\
\endgroup%%------------------------------------<< 

%%----------------------------------------------------------------------------------------------35
\label{5.0av}
\label{5.27b}
\label{8.37}
\label{9.65}
\begingroup%%----------------------------------->>
\fontsize{9pt}{11pt}\selectfont
\textbf{\small EXAMPLE} \ 
Let$(X,x_0)$ be a pointed path connected space.  
Denote by $\text{hur}_n(X)$ the image in $H_n(X)$ of $\pi_n(X)$ under the Hurewicz homomorphism.\\
\indent\indent $(\pi,1)$ \quadx 
Set $\pi = \pi_1(X)$ and assume that $\pi_q(X) = 0$ for $1 < q < n$ $-$then $H_q(X) \approx H_q(\pi,1)$ $(q < n)$ and $H_n(X) /\text{hur}_n(X) \approx H_n(\pi,1)$.
\vspi
[Note: \  In particular, there is an exact sequence $\pi_2(X) \ra H_2(X) \ra H_2(\pi,1) \ra 0$.]\\
\indent\indent $(\pi,n)$ \quadx 
Set $\pi = \pi_n(X)$ $(n > 1)$ and assume that $\pi_q(X) = 0$ for $1 \leq q < n$ $\&$ $\pi_q(X) = 0$ for $n < q < N$ $-$then $H_q(X) \approx H_q(\pi,n)$ $(q < N)$ and $H_N(X)/\text{hur}_N(X) \approx H_N(\pi,n)$.

[Note: \ Take $N = n + 1$ to see that under the stated conditions the Hurewicz homomorphism 
$\pi_{n+1}(X) \ra H_{n+1}(X)$ is surjective.]\\
\endgroup%%------------------------------------<< 

\label{5.8b}
\label{5.8c}
\label{5.14}
\begingroup%%----------------------------------->>
\fontsize{9pt}{11pt}\selectfont
\textbf{\small EXAMPLE} \ 
Let $\pi$ be a finitely generated (finite) abelian group $-$then $\forall \ q \geq 1$, $H_q(\pi,n)$ is finitely generated (finite).  
The $H_q(\pi,1)$ are handled by computation.  
Simply note that $H_q(\Z,1) = $
$
\begin{cases}
\ \Z \quadx (q = 1)\\
\ 0 \quadx (q > 1)
\end{cases}
$
$\&$ $H_q(\Z/k\Z,1) = $
$
\begin{cases}
\ \Z/k\Z \quadx \ \ (q \text{ odd})\\
\ 0 \quadx\indent  (q \text{ even})
\end{cases}
$
and use the K\"unneth formula.  
To pass inductively from $n$ to $n+1$, apply the generalities on 
p.  \pageref{5.0ap} 
to the $\Z$-orientable Hurewicz fibration $\Theta K(\pi,n+1) \ra K(\pi,n+1)$.  
One can, of course, say much more.  Indeed, 
Cartan\footnote[3]{\textit{Collected Works}, vol. III, Springer Verlag (1979), 1300-1394; 
see also Moore, \textit{Ast\'erisque \textbf{32-33} (1976), 173-212.}} 
has explicitly calculated $H_q(\pi,n;G)$, $H^q(\pi,n;G)$ for any finitely generated abelian $G$.  
However, there are occasions when a qualitative description suffices.  
To illustrate, recall that $H^*(\Z,n;\Q)$ is an exterior algebra on one generator of degree 
$n$ if $n$ is odd and a polynomial algebra on one generator of degree $n$ if $n$ is even.  
Therefore, if $n$ is odd, then 
$H_q(\Z,n;\Q) = \Q$ for $q = 0$ $\&$ $q = n$ with 
$H_q(\Z,n;\Q) = 0$ otherwise and if $n$ is even, then 
$H_q(\Z,n;\Q) = \Q$ for $q = kn$ $(k = 0, 1, \ldots)$ with 
$H_q(\Z,n;\Q) = 0$ otherwise.  
So, by the above, if $n$ is odd, then 
$H_q(\Z,n) $ is finite for $q \neq 0$ $\&$ $q \neq n$ and if $n$ is even, then 
$H_q(\Z,n) $ is finite unless $q = kn$ $(k = 0, 1, \ldots)$, 
$H_{kn}(\Z,n)$ being the direct sum of a finite group and an infinite cyclic group.\\
\endgroup%%------------------------------------<< 

\begingroup%%----------------------------------->>
\fontsize{9pt}{11pt}\selectfont
\textbf{\small EXAMPLE} \ 
If $\pi^\prime$ and $\pi^{\prime\prime}$ are finitely generated abelian groups and $\F$ is a field, then the algebra 
$H^*(\pi^\prime \otimes \pi^{\prime\prime}),n;\F)$ is isomorphic to the tensor product over $\F$ of the algebras 
$H^*(\pi^\prime,n;\F)$ and $H^*(\pi^{\prime\prime},n;\F)$.  
Specialize and take $\F = \F_2$ $-$then for $\pi$ a finitely generated abelian group, the determination of $H^*(\pi,n;\F_2)$ reduces to the determination of 
$H^*(\pi,n;\F_2)$ when $\pi = \Z/2^k\Z$, $\pi = \Z/p^\ell\Z$ \ (p = odd prime), or $\pi = \Z$.  
The second possibility is easily dispensed with: 
$H^q(\Z/p^\ell\Z,n; \F_2) = 0$ $\forall \ q > 0$, so 
$H^*(\Z/p^\ell\Z,n; \F_2) =$ $\F_2$.  
The outcome in the other cases involves the Steenrod squares $Sq^i$ and their iterates $Sq^I$.  To review the definitions, a sequence $I = (i_1, \ldots, i_r)$ of positive integers is termed admissible provided that 
$i_1 \geq 2i_2, \ldots i_{r-1} \geq 2i_r$, its excess $e(I)$ being the difference 
$(i_1 - 2i_2) + \cdots + (i_{r-1} - 2i_r) + i_r$.  $Sq^I$ is the composite 
$Sq^{i_1} \circx \cdots \circx Sq^{i_r}$ ($Sq^I$ = id if $e(I) = 0$).
\\
%%----------------------------------------------------------------------------------------------36
\indent\indent $(\pi = \Z/2^k\Z)$ \quadx  Let $u_n$ be the unique nonzero element of $H^n(\Z/2^k\Z,n;\F_2)$.
\vspi
$(k=1)$ \quadx $H^*(\Z/2\Z,1;\F_2) = \F_2[u_1]$, the polynomial algebra with generator $u_1$.  
For $n > 1$, $H^*(\Z/2\Z,n;$ $\F_2) =$ 
$\F_2[(Sq^I u_n)]$, the polynomial algebra with generators the $Sq^I u_n$, 
where $I$ runs through all admissible sequences of excess $e(I) < n$.
\vspi
$(k>1)$ \quadx%%%
$H^*(\Z/2^k\Z,1;\F_2) = \ds\bigwedge (u_1) \otimes \F_2[v_2]$, the tensor product of the exterior algebra with generator $u_1$ and the polynomial algebra with generator $v_2$.  
Here, $v_2$ is the image of the fundamental class under the Bockstein operator 
$H^1(\Z/2^k\Z,1;\F_2) \ra H^2(\Z/2^k\Z,1;\F_2)$ corresponding to the exact sequence 
$0 \ra \Z/2\Z \ra \Z/2^{k+1}\Z \ra \Z/2^k\Z \ra 0$.  
Using this, extend the definition and let $v_n$ be the image of the fundamental class under the Bockstein operator 
$H^n(\Z/2^k\Z,n;\F_2) \ra H^{n+1}(\Z/2^k\Z,n;\F_2)$.  
Write $\overline{Sq}^I u_n = Sq^I u_n$ if $i_r > 1$ and 
$\overline{Sq}^I u_n = Sq^{i_1} \circx \cdots \circx Sq^{i_{r-1}} v_n$ if $i_r = 1$ $-$then for $n > 1$, 
$H^*(\Z/2^k\Z,n;\F_2) =\F_2[\overline{(Sq}^I u_n)]$, 
the polynomial algebra with generators the 
$\overline{Sq}^I u_n$, 
where $I$ runs through all admissible sequences of excess $e(I) < n$.
\\
\indent\indent $(\pi = \Z)$ \quadx%%%
Let $u_n$ be the unique nonzero element of $H^n(\Z,n;\F_2)$ $-$then 
$H^*(\Z,1;\F_2) = \ds\bigwedge (u_1)$, 
the exterior algebra with generator $u_1$, and for $n > 1$, $H^*(\Z,n;\F_2) = \F_2[(Sq^I u_n)]$, the polynomial algebra with generators $Sq^I u_n$,  where $I$ runs through all admissible sequences of excess $e(I) < n$ and $i_r > 1$.
\vspi
Let $\pi$ be a finitely generated abelian group $-$then, as vector spaces over $\F_2$, the $H^q(\pi,n;\F_2)$ are finite dimensional, so it makes sense to consider the associated Poincar\'e series:
\endgroup%%------------------------------------<<

\begingroup%%----------------------------------->>
\fontsize{9pt}{11pt}\selectfont
\[  
P(\pi,n;t) \ = \ \ds\sum\limits_{q=0}^\infty \dim(H^q(\pi,n;\F_2)) \cdot t^q.  
\]
\endgroup%%------------------------------------<<

\begingroup%%----------------------------------->>
\fontsize{9pt}{11pt}\selectfont
Obviously, 
$P(\pi^\prime \oplus \pi^{\prime\prime},n;t) =$ 
$P(\pi^\prime,n;t) \cdot P(\pi^{\prime\prime},n;t)$.  
Examples:  (1) $P(\Z/2\Z,1;t) = \ds\sum\limits_0^\infty t^q$; (2) $P(\Z,1;t) = 1 + t$.
\\
\indent\indent (PS$_1$) \quadx %%%%%%%%%%%%%%%%%%%%%%%% PS1
$P(\pi,n;t)$ converges in the interval $0\leq t < 1$.
\vspi
[It suffices to treat the cases $\pi = \Z/2^k\Z$, $\pi = \Z/p^\ell\Z$ ($p$ = odd prime), $\pi = \Z$.  
The second case is trivial: $P(\Z/p^\ell\Z,n;t) = 1$.\\
\indent\indent $(\pi = \Z/2^k\Z)$ \quadx
In view of what has been said above 
$H^*(\Z/2^k\Z,n;\F_2)$ and 
$H^*(\Z/2\Z,n;\F_2)$ are isomorphic as vector spaces over $\F_2$, 
thus one need only examine the situation when $k = 1$ and $n > 1$.  
Given an admissible $I$, let 
$\abs{I} = i_1 + \cdots + i_r$ 
$(\implies e(I) = 2i_1 - \abs{I})$ $-$then 
$P(\Z/2\Z,n;t) = \ds\prod\limits_{e(I) < n} \frac{1}{1 - t^{n+\abs{I}}}$.  
Since the number of admissible $I$ with 
$e(I) < n$ such that $n + \abs{I} = N$ is equal to the number of decompositions of $N$ of the form 
$N = 1 + 2^{h_1} + \cdots +  2^{h_{n-1}}$, where $0 \leq h_1 \leq \cdots \leq h_{n-1}$, it follows that 
\endgroup%%------------------------------------<<

\begingroup%%----------------------------------->>
\fontsize{9pt}{11pt}\selectfont
\[
P(\Z/2\Z,n;t) \ = \prod\limits_{0 \leq h_1 \leq \cdots \leq h_{n-1}} \ \frac{1}{1 - t^{1 + 2^{h_1} + \cdots +  2^{h_{n-1}}}}.
\]
\endgroup%%------------------------------------<<

\begingroup%%----------------------------------->>
\fontsize{9pt}{11pt}\selectfont
The associated series 
$\ds\sum\limits_{0 \leq h_1 \leq \cdots \leq h_{n-1}} \ t^{1 + 2^{h_1} + \cdots +  2^{h_{n-1}}}$ is convergent if $0\leq t < 1$.
\\ 
\indent\indent $(\pi = \Z)$ \quadx %%%
Assuming that $n > 1$, the extra condition $i_r > 1$ is incorporated by the requirement $h_{n-1} = h_{n-2}$.  Consequently, $P(\Z,n;t) = P(\Z/2\Z,n-1;t)/P(\Z,n-1;t)$ or still,
\endgroup%%------------------------------------<<

\begingroup%%----------------------------------->>
\fontsize{9pt}{11pt}\selectfont
\[
P(\Z,n;t) \ = \ \frac{P(\Z/2\Z,n-1;t) \cdot P(\Z/2\Z,n-3;t) \cdots}{P(\Z/2\Z,n-2;t) \cdot P(\Z/2\Z,n-4;t) \cdots}
\]
\endgroup%%------------------------------------<<

\begingroup%%----------------------------------->>
\fontsize{9pt}{11pt}\selectfont
%%----------------------------------------------------------------------------------------------37
via iteration of the data.]
\vspi
Put $\Phi(\pi,n;x) = \log_2 P(\pi,n;1 - 2^{-x})$ $(0 \leq x < \infty)$.\\
\indent\indent (PS$_2$) \quadx %%%%%%%%%%%%%%%%%%%%%%%% PS2
Suppose that $\pi$ is the direct sum of $\mu$ cyclic groups of order a power of 2, 
a finite group of odd order, and $\nu$ cyclic groups of infinite order $-$then: 
(i) $\mu \geq 1$ $\implies$ $\ds\Phi(\pi,n;x) \sim \frac{\mu x^n}{n!}$; 
(ii) $\mu = 0$ $\&$ $\nu \geq 1$ $\implies$ $\ds\Phi(\pi,n;x) \sim \frac{\nu x^{n-1}}{(n-1)!}$; 
(iii) $\mu = 0$ $\&$ $\nu = 0$ $\implies$ $\Phi(\pi,n;x) = 0$.\\
\label{5.0ar}
\label{5.0as}
\vspi
[The essential point is the asymptotic relation 
$\ds\Phi(\Z/2\Z,n;x) \sim \frac{x^n}{n!}$, everything else being a corollary.  
Observe first that 
$P(\Z/2\Z,1;t) = \ds\frac{1}{1 - t}$ $\implies$ $\Phi(\Z/2\Z,1;x) = x$.  
Proceeding by induction on $n$, introduce the abbreviations $P_n(t) = P(\Z/2\Z,n;t)$, $\Phi_n(x) = \Phi(\Z/2\Z,n;x)$, and the auxiliary functions 
$Q_n(t) =$ 
$\ds\prod\limits_{0 \leq h_1 \leq \cdots \leq h_{n-1}} \ \frac{1}{1 - t^{2^{h_1} + \cdots + 2^{h_{n-1}}}}$, 
$\Psi_n(x) \ = \ $ 
$\log_2 Q_n(1 - 2^{-x})$ $-$then 
$Q_n(t)/P_{n-1}(t) \leq P_n(t) \leq Q_n(t)$ $(0 \leq t < 1)$ $\implies$ 
$\Psi_n(x) - \Phi_{n-1}(x) \leq \Phi_n(x) \leq \Psi_n(x)$ $(0 \leq x < \infty)$.  
Because 
$\Phi_{n-1}(x) \sim$ 
$\ds\frac{x^{n-1}}{(n-1)!}$ (induction hypothesis), one need only show that 
$\Psi_n(x) \sim$ 
$\ds\frac{x^n}{n!}$.  But from the definitions 
$Q_n(t)/P_{n-1}(t) = Q_n(t^2)$, hence 
$\Psi_n(x) =$ 
$\Phi_{n-1}(x) + \Psi_n(x - 1 - \log_2(1 - 2^{-x-1}))$.  So, 
$\forall \ \epsilon > 0$, $\exists \ x_\epsilon > 0$: $\forall \ x > x_\epsilon$,
\endgroup%%------------------------------------<<

\begingroup%%----------------------------------->>
\fontsize{9pt}{11pt}\selectfont
\[
\Psi_n(x-1) +  \frac{(1 - \epsilon)}{(n-1)!}x^{n-1} \leq \Psi_n(x) \leq \Psi_n(x-1 + \epsilon) + \frac{(1 + \epsilon)}{(n-1)!}x^{n-1}.
\]
\endgroup%%------------------------------------<<

\begingroup%%----------------------------------->>
\fontsize{9pt}{11pt}\selectfont
Claim:  Given \mA and $n \geq 1$, 
there exists a polynomial $F_n(x)$ of degree $n$ with leading term $\ds\frac{Ax^n}{n!}$ such that $F_n(x) =$
$\ds F_n(x-1) + \frac{Ax^{n-1}}{(n-1)!}$.
\vspi
[Use induction on $n$:  Put $F_1(x) = Ax$ and consider 
$F_n(x) =$ 
$\ds\frac{Ax^n}{n!} + \sum\limits_{k=2}^n \ \frac{(-1)^k}{k!} F_{n-k+1}(x)$.]
\vspi
Claim:  Let $f \in C([0,\infty[)$.  Assume: 
$f(x) \leq$ 
$\ds f(x-1) + \frac{Ax^{n-1}}{(n-1)!}$ $\left(f(x) \geq f(x-1) + \ds\frac{Ax^{n-1}}{(n-1)!}\right)$ 
$-$then there exists a constant $C^\prime$ $(C^{\prime\prime})$ such that 
$f(x) \leq$ 
$\ds F_n(x) + C^\prime$  ($f(x) \geq F_n(x) + C^{\prime\prime}$).
\vspi
[Let $C^\prime = \max \{f(x) - F_n(x): 0 \leq x \leq 1\}: f(x) \leq F_n(x) + C^\prime$ $(0 \leq x \leq 1)$ and by induction on $N: N \leq x \leq N + 1$ $\implies$ 
$f(x) \leq$ 
$\ds f(x-1) + \frac{Ax^{n-1}}{(n-1)!} \leq$ 
$\ds F_n(x - 1) + C^\prime + \frac{Ax^{n-1}}{(n-1)!} =$ 
$F_n(x) +  C^\prime$.]
\vspi
These generalities allow one to say that $\forall \ \epsilon > 0$, there exist polynomials $R_\epsilon^\prime$ and $R_\epsilon^{\prime\prime}$ of degree $< n: \forall \ x \gg 0$, 
\endgroup%%------------------------------------<<

\begingroup%%----------------------------------->>
\fontsize{9pt}{11pt}\selectfont
\[
(1 - \epsilon) \frac{x^n}{n!} + R_\epsilon^\prime (x) \leq \Psi_n(x) \leq 
\left(\frac{1 + \epsilon}{1 - \epsilon}\right)\frac{x^n}{n!} + R_\epsilon^{\prime\prime}(x).
\]
\endgroup%%------------------------------------<<

\begingroup%%----------------------------------->>
\fontsize{9pt}{11pt}\selectfont
Since $\epsilon$ is arbitrary, this means that 
$\Psi_n(x) \sim$ 
$\ds \frac{x^n}{n!}$.]\\
\endgroup%%------------------------------------<<

\label{5.34c}
\textbf{\small LEMMA} \ 
Suppose that $A$ is path connected $-$then $\forall \ n \geq 1$ there exists a path connected space $X \supset A$ which is obtained from $A$ by attaching $(n+1)$-cells such that $\pi_n(X) = 0$ and, under the inclusion $A \ra X$, $\pi_q(A) \approx \pi_q(X)$ $(q < n)$.

[Let $\{\alpha\}$ be a set of generators for $\pi_n(A)$.  Represent $\alpha$ by 
$f_\alpha:\bS^n \ra A$ and put 
$X = \bigl(\coprod\limits_\alpha \bD^{n+1}\bigr) \sqcup_f A$ \ 
$(f = \coprod\limits_\alpha f_\alpha)$.]\\

%%----------------------------------------------------------------------------------------------38

Let $X$ be a pointed path connected space.  Fix $n \geq 0$ $-$then an 
\un{$n^\text{th}$ Postnikov} 
\un{approximate}
\index{n$^\text{th}$ Postnikov approximate} 
to $X$ is a pointed path connected space $X[n] \supset X$, 
where $(X[n],X)$ is a relative CW complex whose cells in $X[n] - X$ have dimension 
$> n+1$, such that $\pi_q(X[n]) = 0$ $(q > n)$ and, under the inclusion 
$X \ra X[n]$, 
$\pi_q(X) \approx$ 
$\pi_q(X[n])$ $(q \leq n)$.

[Note: \  $X[0]$ is homotopically trivial  and $X[1]$ has the homotopy type $(\pi_1(X),1)$.]\\

\begin{proposition} \ %09
Every pointed path connected space $X$ admits an $n^{th}$ Postnikov approximate $X[n]$.
\end{proposition}

[Using the lemma, construct a sequence 
$X = X_0 \subset X_1 \subset \cdots$ 
of pointed connected spaces $X_k$ such that $\forall \ k > 0$, $X_k$ is obtained from $X_{k-1}$ by attaching 
$(n + k + 1)$-cells, 
$\pi_{n+k}(X_k) = 0$, and, under the inclusion 
$X_{k-1} \ra X_k$, 
$\pi_q(X_{k-1}) \approx$ 
$\pi_q(X_k)$ $(q < n + k)$.  
Consider $X[n] = \colim \ X_k$.]

[Note: \ If \mX is a pointed connected CW space, then the $X[n]$ are pointed connected CW spaces.]\\

\begingroup%%----------------------------------->>
\fontsize{9pt}{11pt}\selectfont
\label{5.44}
\textbf{\small EXAMPLE} \ 
Let $\pi$ be a group and let $n$ be an integer $\geq 1$, where $\pi$ is abelian if $n > 1$ 
$-$then a pointed connected CW space $X$ is said to be a 
\un{Moore space}
\index{Moore space} 
of type $(\pi,n)$ provided that $\pi_n(X)$ is isomorphic to $\pi$ and
$
\begin{cases}
\ \pi_q(X) = 0 \hspace{0.60cm} (q < n) \\
\ H_q(X) = 0 \hspace{0.5cm} (q > n)
\end{cases}
\hspace{-.25cm}. \  
$
Notation: $X = M(\pi,n)$.  
If $n = 1$, then $M(\pi,n)$ exists iff $H_2(\pi,1) = 0$ but if $n > 1$, then $M(\pi,n)$ always exists.  
If $n = 1$ and $H_2(\pi,1) = 0$, then the pointed homotopy type of $M(\pi,1)$ is not necessarily unique (e.g., when $\pi = \Z)$ but if $n > 1$, then the pointed homotopy type of $M(\pi,n)$ is unique.  
In any event, $M(\pi,n)[n] = K(\pi,n)$.\\
\endgroup%%------------------------------------<< 

\label{9.119}
\label{11.33}
\begingroup%%----------------------------------->>
\fontsize{9pt}{11pt}\selectfont
\textbf{\small FACT} \ 
Suppose that $X$ is a pointed path connected space.  
Fix $n \geq 1$ $-$then there exists a pointed $n$- connected space $\widetilde{X}_n$ in $\bTOP/X$ such that the projection $\widetilde{X}_n \ra X$ is a pointed Hurewicz fibration and induces an isomorphism $\pi_q(\widetilde{X}_n) \ra \pi_q(X)$ $\forall \ q > n$.
\vspi
[Consider the mapping fiber of the inclusion $X \ra X[n]$.]\\
\endgroup%%------------------------------------<< 

\begingroup%%----------------------------------->>
\fontsize{9pt}{11pt}\selectfont
\textbf{\small EXAMPLE} \ 
Take $X = \bS^3$ $-$then the fibers of the projection $\widetilde{X}_3 \ra X$ have homotopy type $(\Z,2)$ and $\forall \ q \geq 1$, $H_q(\widetilde{X}_3) = $
$
\begin{cases}
\ 0 \hspace{1.65cm} (q \text{ odd})\\
\ \Z/(q/2)\Z \hspace{0.50cm}  (q \text{ even})
\end{cases}
$
\hspace{-.25cm}.
\vspi
[Use the Wang cohomology sequence and the fact that $H^*(\Z,2)$ is the polynomial algebra over $\Z$ generated by an element of degree 2.]
\vspi
[Note: \  Given a prime $p$, let $\sC$ be the class of finite abelian groups with order prime to $p$ 
$-$then from the above, $H_n(\widetilde{X}_3) \in \sC$ $(0 < n < 2p)$, 
so by the mod $\sC$ Hurewicz theorem, $\pi_n(\widetilde{X}_3)  \in \sC$ 
$(0 < n < 2p)$ and the Hurewicz homomorphism 
$\pi_{2p}(\widetilde{X}_3) \ra H_{2p}(\widetilde{X}_3)$ is $\sC-$bijective.  
Therefore the $p$-primary component of $\pi_n(\bS^3)$ is 0 if $n < 2p$ and is $\Z/p\Z$ if $n = 2p$.]\\
\endgroup%%------------------------------------<< 

%%----------------------------------------------------------------------------------------------39
\begingroup%%----------------------------------->>
\fontsize{9pt}{11pt}\selectfont
Put $W_1 = \widetilde{X}_1$.  Let $W_2$ be the mapping fiber of the inclusion $\widetilde{X}_1 \ra \widetilde{X}_1[2]$ $-$then the mapping fiber of the projection $W_2 \ra W_1$ has the homotopy type $(\pi_2(X),1)$.  
Iterate:  The result is a sequence of pointed Hurewicz fibrations $W_n \ra W_{n-1}$, where the mapping fiber has homotopy type $(\pi_n(X),n-1)$ and $W_n$ is $n$-connected with $\pi_q(W_n) \approx \pi_q(X)$ $(\forall \ q > n)$.  
The diagram $X = W_0 \leftarrow W_1 \leftarrow \cdots$ is called ``the'' 
\un{Whitehead tower} 
\index{Whitehead tower} 
of $X$.
\vspi
[Note: \  If $X$ is a pointed connected CW space, then the $W_n$ are pointed connected CW spaces and the mapping fiber of the projection $W_n \ra W_{n-1}$ is a $K(\pi_n(X),n-1)$.]\\
\endgroup%%------------------------------------<< 

\label{5.0au}
\begingroup%%----------------------------------->>
\fontsize{9pt}{11pt}\selectfont
\textbf{\small EXAMPLE} \ 
Let $X$ be a pointed simply connected CW complex which is finite and noncontractible.  
Assume:  $\exists \ i > 0$ such that $H_i(X;\F_2) \neq 0$ $-$then $\pi_q(X)$ contains a subgroup isomorphic to $\Z$ or $\Z/2\Z$ for infinitely many $q$.
\vspi
[Because the $H_q(X)$ are finitely generated $\forall \ q$, the same is true for the $\pi_q(X)$ 
(cf. p. \pageref{5.0aq}).  
The set of positive integers $n$ such that $\pi_n(X) \otimes \Z/2\Z \neq 0$ is nonempty.  
To get a contradiction, suppose that there is a largest integer $N$.  
Working with the Whitehead tower of $X$, let 
$P_n(t) =$ 
$\ds\sum\limits_{q=0}^\infty \ \dim(H^q(W_n;\F_2)) \cdot t^q$, the mod 2 Poincar\'e series of 
$H^*(W_n;\F_2)$ (meaningful, the $H^q(W_n;\F_2)$ being finite dimensional over $\F_2$).  
In particular:  $P_N(t) = 1$, $P_{N-1}(t) = P(\pi_N(X),N;t)$, $P_1(t) = P_X(t)$, the Poincar\'e series of $H^*(X;\F_2)$.  
On general grounds, there is a majorization $P_n(t) \prec P_{n-1}(t)$ $\cdot P(\pi_N(X),n-1;t)$, 
where the symbol $\prec$ means that each coefficient of the formal power series on the left is $\leq$ the corresponding coefficient of the formal power series on the right.  
So, starting with $n = N-1$ and multiplying out, one finds that 
$P(\pi_N(X),N;t) \prec$ 
$\ds P_X(t) \cdot \prod\limits_{1 < i < N} P(\pi_i(X),i-1;t)$.  Since $P_X(t)$ is a polynomial, hence is bounded on $[0,1]$, $\exists \ C > 0$: 
$\ds P(\pi_N(X),N;t) \leq$ 
$\ds C  \cdot \prod\limits_{1 < i < N} P(\pi_i(X),i-1;t)$, or still, in the notation of 
p. \pageref{5.0ar}, 
$\Phi(\pi_N(X),N;x) \leq$ 
$\ds \log_2 C + \sum\limits_{1 < i < N} \Phi(\pi_i(X),i-1;x)$ $(0 \leq x < \infty)$.  Comparing the asymptotics of either side leads to an immediate contradiction 
(cf. p. \pageref{5.0as}).]
\vspi
[Note: \  This analysis is due to 
Serre\footnote[2]{\textit{Comment. Math. Helv.} \textbf{27} (1953), 198-232.}.  
It has been extended to all odd primes by 
Umeda\footnote[3]{\textit{Proc. Japan Acad.} \textbf{35} (1959), 563-566; 
see also McGibbon-Neisendorfer, \textit{Comment. Math. Helv.} \textbf{59} (1984), 253-257.}.  
Accordingly, if $X$ is a pointed simply connected CW complex which is finite and noncontractible, then $\pi_q(X)$ is nonzero for infinitely many $q$.  Proof:  If $\forall \ p \in \Pi$ $\&$ $\forall \ i > 0$, $H_i(X;\F_p) = 0$, then the arrow $X \ra *$ is a homology equivalence 
(cf. p. \pageref{5.0at}), 
thus by the Whitehead theorem, $X$ is contractible.]\\
\endgroup%%------------------------------------<< 

\textbf{\small LEMMA} \ 
Let $(X,A,x_0)$ be a pointed pair.  
Assume: $(X,A)$ is a relative CW complex whose cells in $X - A$ have dimension $> n + 1$.  
Suppose that $(Y,y_0)$ is a pointed space such that $\pi_q(Y,y_0) = 0$ $\forall \ q > n$ $-$then every pointed continuous function $f:A \ra Y$ has a pointed continuous extension $F:X \ra Y$.\\

%%----------------------------------------------------------------------------------------------40
It follows from the lemma that if $X$ and $Y$ are pointed path connected spaces and if $f:X \ra Y$ is a pointed continuous function, 
then for $m \leq n$ there exists a pointed continuous function $f_{n,m}:X[n] \ra Y[m]$ rendering the diagram
\begin{tikzcd}%[sep=large]
X  \arrow{r}{f} \arrow{d} &Y \arrow{d}\\
X[n] \arrow{r}[swap]{f_{n,m}} &Y[m]
\end{tikzcd}
commutative, any two such being homotopic $\rel X$.  
Proof:  Let $F:X \ra Y[m]$ be the composite
%\begin{tikzcd}
%X \arrow{r}{f} &Y \arrow{r} &Y[m]
%\end{tikzcd}
$X \overset{f}{\ra} Y \ra Y[m]$ 
.  
To establish the existence of $f_{n,m}$, consider any filler for 
\begin{tikzcd}%[sep=large]
X[n]  \arrow[dashed]{dr}\\
X \arrow{u} \arrow{r}[swap]{F} &Y[m]
\end{tikzcd}
and to establish the uniqueness of $f_{n,m}$ $\rel X$, take two extenstion $f_{n,m}^\prime$ $\&$ $f_{n,m}^{\prime\prime}$, define $\Phi:i_0X[n] \cup IX \cup i_1X[n] \ra Y[m]$ by
$
\begin{cases}
\ \Phi(x,0) = f_{n,m}^\prime(x)\\
\ \Phi(x,1) = f_{n,m}^{\prime\prime}(x)
\end{cases}
, \ 
$
$\Phi(x,t) = F(x)$, and consider any filler for
\begin{tikzcd}%[sep=large]
IX[n]  \arrow[dashed]{dr}\\
{i_0X[n] \cup IX \cup i_1X[n]} \arrow{u} \arrow{r}[swap]{\Phi} &Y[m]
\end{tikzcd}
.\\
\vspace{0.25cm}

Application:  Let $X^\prime[n]$ and $X^{\prime\prime}[n]$ be the $n^\text{th}$ Postnikov approximates to $X$ $-$then in $\bHTOP^2$, $(X^\prime[n],X) \approx (X^{\prime\prime}[n],X)$.\\

\begingroup%%----------------------------------->>
\fontsize{9pt}{11pt}\selectfont
\textbf{\small EXAMPLE} \ \ 
Let $X$ and $Y$ be pointed connected CW spaces $-$then it can happen that $X[n]$ and $Y[n]$ have the same pointed homotopy type for all $n$, yet $X$ and $Y$ are not homotopy equivalent.  
To construct an example, let $K$ be a pointed simply connected CW complex.  
Assume:  $K$ is finite and noncontractible.  Put 
$X =$ 
$\ds (w)\prod\limits_0^\infty K[n]$, $Y = X \times K$ 
$-$then $\forall \ n$, $X[n] \approx Y[n]$ in $\bHTOP_*$.  
However, it is not true that 
$X \approx Y$ in \bHTOP.  For if so, $K$ would be dominated in homotopy by $X$ or still, by 
$K[0] \times \cdots \times K[n]$ $(\exists \ n)$, thus $\forall \ q$, $\pi_q(K)$ would be a direct summand of 
$\pi_q(K[0] \times \cdots \times K[n])$.  
But this is impossible:  The $\pi_q(K)$ are nonzero for infinitely many $q$ (cf. \pageref{5.0au}).
\vspi
[Note: \ This subject has its theoretical aspects as well.  McGibbon-M$\o$ller
\footnote[2]{\textit{Topology} \textbf{31} (1992), 177-201; 
see also Dror-Dwyer-Kan, \textit{Proc. Amer. Math. Soc.} \textbf{74} (1979), 183-186.}).]\\
\endgroup%%------------------------------------<< 

Let $X$ be a pointed path connected space.  Given a sequence $X[0], X[1], \ldots$  \ of Postnikov approximates to $X$, 
$\forall \ n \geq 1$ there is a pointed continuous function $f_n:X[n] \ra X[n-1]$ such that the triangle
\begin{tikzcd}[ sep=small]
&X \arrow{ldd} \arrow{rdd} \\
\\
X[n] \arrow{rr}[swap]{f_n} &&X[n-1]
\end{tikzcd}
commutes.  Put $P_0X = X[0]$, let
%%----------------------------------------------------------------------------------------------41
$s_0$ be the identity map, and denote by $P_1X$ the mapping track of $f_1:$
$
\begin{tikzcd}%[sep=large]
X[1]  \arrow{r}{f_1} \arrow{d}[swap]{s_1} &X[0]\arrow{d}{s_0}\\
P_1X \arrow{r}[swap]{p_1} &P_0X
\end{tikzcd}
.
$
\  
Recall that $s_1$ is a pointed homotopy equivalence, while $p_1$ is the usual pointed Hurewicz fibration associated with this setup.  Repeat the procedure, taking for $P_2X$ the mapping  track of $s_1 \circx f_2:$
\begin{tikzcd}%[sep=large]
X[2]  \arrow{r}{f_2} \arrow{d}[swap]{s_s} &X[1]\arrow{d}{s_1}\\
P_2X \arrow{r}[swap]{p_2} &P_1X
\end{tikzcd}
.  
\label{9.11}
The upshot is that the $f_n$ can be converted to pointed Hurewicz fibrations $p_n$, where at each stage there is a commutative triangle
$
\begin{tikzcd}[ sep=small]
&X \arrow{ldd} \arrow{rdd} \\
\\
P_nX \arrow{rr}[swap]{p_n} &&P_{n-1}X\\
\end{tikzcd}
.
$
\ 
The diagram $P_0X \leftarrow P_1X \leftarrow \cdots$ 
of pointed Hurewicz fibrations is called ``the'' 
\vspace{0.1cm}
\un{Postnikov tower}
\index{Postnikov tower} 
of $X$.  Obviously, $\pi_q(P_nX) = 0$ $(q > n)$, $\pi_q(X) \approx \pi_q(P_nX)$ $(q \leq n)$, 
and $\pi_q(P_nX) \approx \pi_q(P_{n-1}X)$ $(q \neq n)$.  
Therefore the mapping fiber of $p_n$ has homotopy type $(\pi_n(X),n)$.

[Note: \  If $X$ is a pointed connected CW space, then the $P_nX$ are pointed connected CW spaces, so the mapping fiber of $p_n$ is a $K(\pi_n(X),n)$.]\\

\label{9.63}
\begingroup%%----------------------------------->>
\fontsize{9pt}{11pt}\selectfont
\textbf{\small EXAMPLE} \ 
Let $X$ be a pointed path connected space.  Fix $n > 1$ $-$then $\pi_n(X)$ defines a locally constant coefficient system on $P_{n-1}X$ and there is an exact sequence
\begin{align*}
H_{n+2}(P_nX) \ra &H_{n+2}(P_{n-1}X) \ra H_{1}(P_{n-1}X;\pi_n(X)) \ra H_{n+1}(P_{n}X) \ra H_{n+1}(P_{n-1}X)\\
&\ra H_{0}(P_{n-1}X;\pi_n(X)) \ra H_{n}(P_{n}X) \ra H_{n}(P_{n-1}X) \ra 0.]
\end{align*}
\indent
[Work with the fibration spectral sequence of $p_n:P_nX \ra P_{n-1}X$, noting that $E_{p,q}^r = 0$ if $0 < q < n$ or $q = n+1$.]\\
\endgroup%%------------------------------------<< 

A nonempty path connected topological space $X$ is said to be 
\un{abelian}
\index{abelian (space)} 
if 
$\pi_1(X)$ is abelian and if $\forall \ n > 1$, $\pi_1(X)$ operates trivially on $\pi_n(X)$.  
Every simply connected space is abelian as is every path connected H-space or every path connected 
compactly generated semigroup with unit (obvious definition).

[Note: \  If $X$ is abelian, then $\forall \ x_0 \in X$, the forgetful function  
$[\bS^n;s_n;X, x_0] \ra$ 
$[\bS^n,X]$ is bijective (cf. p. \pageref{5.1}).\\

\begingroup%%----------------------------------->> 
\fontsize{9pt}{11pt}\selectfont
\textbf{\small EXAMPLE} \ 
$\bP^n(\R)$ is abelian iff $n$ is odd.\\
\endgroup%%------------------------------------<< 

%%----------------------------------------------------------------------------------------------42
\label{5.20}
Let $X$ be a pointed connected CW space.  
Assume: $X$ is abelian.  
There is a commutative triangle
\begin{tikzcd}[ sep=small]
&X \arrow{ldd} \arrow{rdd}\\ 
\\
X[n+1]   \arrow{rr}[swap]{f_{n+1}} &&X[n]
\end{tikzcd}
and an embedding $IX \ra M_{f_{n+1}}$.  
Define $\widehat{X}[n]$ by the pushout square
\begin{tikzcd}%[sep=large]
IX \arrow{d} \arrow{r}{p} &X  \arrow{d}\\
M_{f_{n+1}} \arrow{r}   &\widehat{X}[n]
\end{tikzcd}
$-$then  $\widehat{X}[n]$ contains $X[n]$ as a strong deformation retract, hence $\pi_q(\widehat{X}[n]) \approx \pi_q(X[n])$ $(q \geq 1)$.  Using the exact sequence
\[
\cdots 
\hspace{-.05cm}
\ra 
\hspace{-.05cm}
\pi_{q+1}(X[n+1]) 
\hspace{-.05cm}
\ra 
\hspace{-.05cm}
\pi_{q+1}(\widehat{X}[n]) 
\hspace{-.05cm}
\ra 
\hspace{-.05cm}
\pi_{q+1}(\widehat{X}[n],X[n+1]) 
\ra\pi_{q}(X[n+1]) 
\hspace{-.05cm}
\ra 
\hspace{-.05cm}
\pi_{q}(\widehat{X}[n]) 
\hspace{-.05cm}
\ra 
\hspace{-.05cm}
\cdots,
\]
one finds that $\pi_q(\widehat{X}[n],X[n+1]) = 0$ 
$(q \neq n+2)$ and 
$\pi_{n+2}(\widehat{X}[n],X[n+1])$ $\approx$ 
$\pi_{n+1}(X[n+1])$ $\approx$ $\pi_{n+1}(X)$.  
Thus the relative Hurewicz homomorphism 
$\text{hur}:\pi_{n+2}(\widehat{X}[n],X[n+1]) \ra$ 
$H_{n+2}(\widehat{X}[n],X[n+1])$ 
is bijective, so the composite
\begin{tikzcd}%[sep=large]
{\kappa_{n+2}:H_{n+2}(\widehat{X}[n],X[n+1])} \arrow{r}{\text{hur}^{-1}}  
&{}%&{\pi_{n+2}(\widehat{X}[n],X[n+1])}
\end{tikzcd}
${\pi_{n+2}(\widehat{X}[n],X[n+1])}$
%  \arrow{r} 
%&{\pi_{n+1}(X)}
%$\end{tikzcd}
$\ra {\pi_{n+1}(X)}$ 
is an isomorphism.  Since $H_{n+1}(\widehat{X}[n],X[n+1]) = 0$, the universal coefficient theorem implies that $H^{n+2}(\widehat{X}[n],X[n+1]; \pi_{n+1}(X))$ can be identified with $
\Hom(H_{n+2}(\widehat{X}[n],X[n+1]); \pi_{n+1}(X))$, therefore 
$\kappa_{n+2}$ 
corresponds to a cohomology class in 
$H^{n+2}(\widehat{X}[n],X[n+1]; \pi_{n+1}(X))$ 
whose image 
$\textbf{k}^{n+2}$ $(= \textbf{k}^{n+2}(X))$ in $H^{n+2}(X[n];\pi_{n+1}(X))$ is the 
\label{5.21}
\un{Postnikov invariant}
\index{Postnikov invariant} 
of $X$ in dimension $n+2$.  Put $K_{n+2} = K(\pi_{n+1}(X),n+2)$, let $k_{n+2}:X[n] \ra K_{n+2}$ be the arrow associated with $\textbf{k}^{n+2}$, and define $W[n+1]$ by the pullback square
\begin{tikzcd}%[sep=large]
W[n+1] \arrow{d} \arrow{r} &{\Theta K_{n+2}}  \arrow{d}\\
X[n] \arrow{r}[swap]{k_{n+2}}   &{K_{n+2}}
\end{tikzcd}
$-$then $W[n+1]$ is a CW space (cf. $\S 6$, Proposition 9) and there is a lifting
\begin{tikzcd}%[sep=large]
&W[n+1]  \arrow{d}\\ 
X[n+1]   \arrow{ur}{\Lambda_{n+1}} \arrow{r}[swap]{f_{n+1}} &X[n]
\end{tikzcd}
of $f_{n+1}$ which is a weak homotopy equivalence or still, a homotopy equivalence (realization theorem).  
The restriction of $\Lambda_{n+1}$ to $X$ is an embedding and 
$\Lambda_{n+1}:(X[n+1],X) \ra (W[n+1],X)$ is a homotopy equivalence of pairs.

[Note: \  $\Lambda_{n+1}$ is constructed by considering a specific factorization of $k_{n+2}$ as a composite 
$X[n] \ra$ 
$\widehat{X}[n]/X[n+1] \ra$ 
$K_{n+2}$ $(k_{n+2}$ 
is determined only up to homotopy.)]\\

\index{Theorem, Invariance}
\textbf{\small INVARIANCE THEOREM} \ 
Let
$
\begin{cases}
\ X\\
\ Y
\end{cases}
$
be pointed CW spaces.  Assume:
$
\begin{cases}
\ X\\
\ Y
\end{cases}
$
are abelian.  Suppose that $\phi:X \ra Y$ is a pointed continuous function.   Fix a pointed $\phi_n:X[n] \ra Y[n]$ such that the diagram
\begin{tikzcd}%[sep=large]
X \arrow{d} \arrow{r}{\phi} &Y \arrow{d} \\
X[n] \arrow{r}[swap]{\phi_n} &Y[n]
\end{tikzcd}
commutes $-$then 
$\forall \ n$, $\phi_n^*\textbf{k}^{n+2}(Y) = \phi_{\text{co}}\textbf{k}^{n+2}(X)$ in $H^{n+2}(X[n];\pi_{n+1}(Y))$.

[Note: \   Here $\phi_\text{co}$ is the coefficient group homomorphism 
$H^{n+2}(X[n];\pi_{n+1}(X)) \ra$ 
$H^{n+2}(X[n];\pi_{n+1}(Y))$.]\\

%%----------------------------------------------------------------------------------------------43
\index{Theorem, Nullity}
\textbf{\small NULLITY THEOREM} \ 
Let $X$ be a pointed CW space.  
Assume: $X$ is abelian $-$then 
$\textbf{k}^{n+1} = 0$ iff the Hurewicz homomorphism $\pi_n(X) \ra H_n(X)$ is split injective.\\

\label{14.95a}
\begingroup%%----------------------------------->>
\fontsize{9pt}{11pt}\selectfont
\textbf{\small EXAMPLE} \ 
Suppose that $\textbf{k}^{n+1} = 0$ $-$then $W[n]$ is fiber homotopy equivalent to 
$X[n-1] \times$ $K(\pi_n(X),n)$ 
(cf. p. \pageref{5.0au1}), 
hence $X[n] \approx X[n-1] \times K(\pi_n(X),n)$.  
Therefore $X$ has the same pointed homotopy type as the weak product  \ 
$\ds (w)\prod\limits_1^\infty K(\pi_n(X),n)$ 
\ 
provided that the Hurewicz homomorphism 
$\pi_n(X) \ra H_n(X)$ is split injective for all $n$.  
This condition can be realized.  
In fact, 
Puppe\footnote[2]{\textit{Math. Zeit.} \textbf{68} (1958), 367-421.}
 has shown that if $G$ is a path connected abelian compactly generated semigroup with unit, then $\forall \ n$, the Hurewicz homomorphism $\pi_n(G) \ra H_n(G)$ is split injective, thus  
 $G \approx$ 
 $\ds (w)\prod\limits_1^\infty K(\pi_n(G),n)$ when $G$ is in addition a CW space.
 \vspi
[Note: \  Analogous remarks apply if $G$ is a path connected abelian topological semigroup with unit.  
Reason:  The identity map $kG \ra G$ is a weak homotopy equivalence.]\\ 
\endgroup%%------------------------------------<<

\begingroup%%----------------------------------->>
\fontsize{9pt}{11pt}\selectfont
\label{11.15}
\index{Theorem, Abelian Obstruction}
\textbf{\small ABELIAN OBSTRUCTION THEOREM} \ 
Let $(X,A)$ be a relative CW complex; let $Y$ be a pointed abelian CW space.  Suppose that $\forall \ n > 0$, 
$H^{n+1}(X,A;\pi_n(Y)) =$ $0$ 
$-$then  every $f \in C(A,Y)$ admits an extension $F \in C(X,Y)$, 
any two such being homotopic rel $A$ provided that $\forall \ n > 0$, $H^{n}(X,A;\pi_n(Y)) = 0$.\\
\endgroup%%------------------------------------<<

\begingroup%%----------------------------------->>
\fontsize{9pt}{11pt}\selectfont
\textbf{\small EXAMPLE} \ 
Let $(X,x_0)$ be a pointed CW complex; let $(Y,y_0)$ be a pointed simply connected CW complex.  Assume:  $\forall \ n > 0$, $H^{n}(X;\pi_n(Y)) = 0$ $-$then $[X,x_0;Y,y_0] = *$.
\vspi
[In fact, $H^{n}(X,x_0;\pi_n(Y,y_0)) \approx H^{n}(X;\pi_n(Y)) = 0 \implies [X,x_0;Y,y_0]  = *$ $(\implies [X,Y] = *$ (cf. p. \pageref{5.8a})).]\\
\endgroup%%------------------------------------<<

\begin{proposition} \ %10
Let $X$ be a pointed abelian CW space.  Assume:  The $H_q(X)$ are finitely generated $\forall \ q$ $-$then $\forall \ n$, the $H_q(X[n])$ are finitely generated $\forall \ q$.
\end{proposition}

[The assertion is trivial if $n = 0$.  
Next, $X[1]$ is a $K(\pi_1(X),1)$, hence 
$\pi_1(X) \approx H_1(X)$, which is finitely generated.  
For $q > 1$, 
$H_q(X[1]) \approx$ 
$H_q(\pi_1(X),1)$ and these too are 
%%----------------------------------------------------------------------------------------------44
finitely generated (cf. p. \pageref{5.8b}).  Proceeding by induction, suppose that the $H_q(X[n])$ are finitely generated 
$\forall \ q$ $-$then the $H_q(X[n],X)$ are finitely generated $\forall \ q$.  In particular, $H_{n+2}(X[n],X)$ is finitely generated.  Since 
$\pi_{n+1}(X[n]) = \pi_{n+2}(X[n]) = 0$, the arrow 
$\pi_{n+2}(X[n],X) \ra \pi_{n+1}(X)$ is an isomorphism.  But \mX is abelian, so from the relative Hurewicz theorem, 
$\pi_{n+2}(X[n],X) \approx$ 
$H_{n+2}(X[n],X)$.  
Therefore $\pi_{n+1}(X)$ is finitely generated.  
Consider now the mapping track $W_{n+2}$ of $k_{n+2}:X[n] \ra K_{n+2}$.  
The fiber of the $\Z$-orientable Hurewicz fibration 
$W_{n+2} \ra K_{n+2}$ over the base point is homeomorphic to $W[n+1]$ (parameter reversal).  The 
$H_q(K_{n+2}) = H_q(\pi_{n+1}(X),n+2)$ are finitely generated $\forall \ q$ (cf. p. \pageref{5.8c}), as are the 
$H_q(W_{n+2})$ (induction hypothesis), thus the $H_q(W[n+1])$ are finitely generated $\forall \ q$ 
(cf. p. \pageref{5.8d}).  
Because $X[n+1]$ and $W[n+1]$ have the same homotopy type, this completes the passage from $n$ to $n + 1$.]\\

\label{5.12}
\label{5.0aq}
Application:  Let $X$ be a pointed abelian CW space.  Assume:  The $H_q(X)$ are finitely generated $\forall \ q$ $-$then the $\pi_q(X)$ are finitely generated $\forall \ q$.

[Note: \  This result need not be true for a nonabelian $X$.  
Example:  Take $X = \bS^1 \vee \bS^2$ 
$-$then the $H_q(X)$  are finitely generated $\forall \ q$ and $\pi_1(X) \approx \Z$.  
On the other hand, $\pi_2(X) \approx H_2(\widetilde{X})$, 
$\widetilde{X}$ the universal covering space of $X$, 
i.e., the real line with a copy of $\bS^2$ attached at each integral point.  
Therefore $\pi_2(X)$ is free abelian on countably many generators.]\\

\begin{proposition} \ %11
Let $X$ be a pointed abelian CW space.  Assume:  The $H_q(X)$ are finite $\forall \ q > 0$ $-$then $\forall \ n$, the $H_q(X[n])$ are finite $\forall \ q > 0$.\\
\end{proposition}

\label{5.26}
Application:  Let $X$ be a pointed abelian CW space.  Assume:  The $H_q(X)$ are finite $\forall \ q > 0$ $-$then the $\pi_q(X)$ are finite $\forall \ q > 0$.\\

\begingroup%%----------------------------------->>
\fontsize{9pt}{11pt}\selectfont
\label{9.17}
\label{17.14}
\index{Homotopy Groups of Spheres}
\textbf{\small EXAMPLE \ (\un{Homotopy Groups of Spheres})} \ 
The $\pi_q(\bS^{2n+1})$ of the odd dimensional spheres are finite for $q > 2n + 1$ and the $\pi_q(\bS^{2n})$ of the even dimensional sphere are finite for $q > 2n$ except that $\pi_{4n - 1}(\bS^{2n})$ is the direct sum of $\Z$ and a finite group.  Here are the details.\\
\indent \indent $(2n + 1)$ \quadx
Fix a map 
$f:\bS^{2n + 1} \ra K(\Z,2n + 1)$ classifying a generator of 
$H^{2n + 1}(\bS^{2n + 1})$ $-$then $f_*$ induces an isomorphism 
$H_*(\bS^{2n + 1};\Q) \ra H_*(K(\Z,2n + 1); \Q)$
 (cf. p. \pageref{5.0av}), 
 so $\forall \ q > 0$, $H_q(E_f; \Q) = 0$ 
 (cf. p. \pageref{5.0aw}).  
 Accordingly, $\forall \ q > 0$, $H_q(E_f)$ is finite (being finitely generated).  
 Therefore all the homotopy groups of $E_f$ are finite.  
 But $\pi_q(E_f) \approx \pi_q(\bS^{2n + 1})$ if $q > 2n + 1$.\\
\indent \indent $(2n)$ \quadx
The even dimensional case requires a double application of the odd dimensional case.  
First, consider the Stiefel manifold $\bV_{2n +1, 2}$ and the map 
$f:\bV_{2n +1, 2} \ra \bS^{4n - 1}$ defined on p. \pageref{5.0ax}.  
As noted there, $\forall \ q > 0$, $H_q(E_f;\Q) = 0$, hence the $\pi_q(E_f)$ are finite and this means that the 
$\pi_q(\bV_{2n +1, 2})$ are finite save for 
$\pi_{4n -1}(\bV_{2n +1, 2})$ which is the direct sum of $\Z$ and a finite group. 
 Second, examine the homotopy sequence of the Hurewicz fibration 
$\bV_{2n +1, 2} \ra \bS^{2n}$, noting that its fiber is $\bS^{2n - 1}$.\\
\endgroup%%------------------------------------<<

%%----------------------------------------------------------------------------------------------45
\label{8.39} %dmc mnft
Given a category \bC, the 
\un{tower category}
\index{tower category} 
\textbf{TOW(C)} of \bC is the functor category $[[\N]^{\OP},\bC]$.  
Example:  The Postnikov tower of a pointed path connected space is an object in \textbf{TOW(TOP$_*$)}.

Take \bC = \bAB $-$then an object (\bG,\bff) in \textbf{TOW(AB)} is a sequence 
$\{G_n,f_n:G_{n+1} \ra G_n\}$, where $G_n$ is an abelian group and 
$f_n:G_{n+1} \ra$ $G_n$ is a homomorphism, a morphism 
$\phi:(\bG^\prime,\bff^\prime) \ra (\bG^{\prime\prime},\bff^{\prime\prime)})$ in 
$\textbf{TOW(AB)}$ being a sequence 
$\{\phi_n\}$, where $\phi_n:G_n^\prime \ra G_n^{\prime\prime}$ is a homomorphism and 
$\phi_n \circx f_n^\prime = f_n^{\prime\prime} \circx \phi_{n+1}$.  
\textbf{TOW(AB)} is an abelian category.   
As such, it has enough injectives.

\label{11.30}
[Note: \  Equip $[\N]$ with the topology determined by $\leq$, i.e., regard $[\N]$ as an A space 
$-$then \textbf{TOW(AB)}  is equivalent to the category of sheaves of abelian groups on $[\N]$.]

The functor $\lim: \textbf{TOW(AB)}  \ra \bAB$ 
that sends \bG  to $\lim \bG$ is left exact (being a right adjoint) but it need not be exact.  
The right dervied functors $\lim^i$ of $\lim$ live only in dimensions 0 and 1,
 i.e., the $\lim^i$ $(i > 1)$ necessarily vanish.  
To compute 
$\lim^i \bG$, form $G = \prod\limits_n G_n$ and define 
$d:G \ra G$ by $d(x_0,x_1, \ldots) = (x_0 - f_0(x_1), x_1 - f_1(x_2), \ldots)$ $-$then 
$\ker d = \lim \bG$ and $\text{coker } d = \lim^1 \bG$.  
Example:  Suppose that $\forall \ n$, $G_n$ is finite, $-$then $\lim^1 \bG = 0$.

[Note: \  Translated to sheaves, $\lim^i$ corresponds to the $i^{th}$ 
right derived functor of the global section functor.]\\

\begingroup%%----------------------------------->>
\fontsize{9pt}{11pt}\selectfont
The fact that the $\lim^i$ $(i > 1)$ vanish is peculiar to the case at hand.  
Indeed, if $(I,\leq)$ is a directed set and if \bI is the associated filtered category, then for a suitable choice of $I$, one can exhibit a \bG in \textbf{[I$^{\OP}$,AB]} such that $\lim^i \bG \neq 0$ $\forall \ i > 0$ 
(Jensen\footnote[2]{\textit{SLN} \textbf{254} (1972), 51-52.}).\\
\endgroup%%------------------------------------<< 

\label{8.40}
\begingroup%%----------------------------------->>
\fontsize{9pt}{11pt}\selectfont
\textbf{\small EXAMPLE} \ 
Let $\mu \neq \nu$ be relatively prime natural numbers $> 1$.  
Define $\bG(\mu$) in \textbf{TOW(AB)}  by $G(\mu)_n = \Z \ \forall \ n$ $\&$
$
\begin{cases}
\ G(\mu)_{n+1} \ra G(\mu)_n\\
\ 1 \ra \mu
\end{cases}
$
and $\phi \in \Mor(\bG(\mu),\bG(\mu))$ by $\phi_n(1) = \nu$ 
$-$then the cokernel of $\phi$ is isomorphic to the constant tower on $[\N]$ with value $\Z/\nu\Z$.  
Applying $\lim$ to the exact sequence 
%\begin{tikzcd}
%0 \arrow{r} &{\bG(\mu)} \arrow{r}{\phi} &{\bG(\mu)} \arrow{r} &{\text{coker}\phi} \arrow{r} &0
%\end{tikzcd}
$ 0 \ra $
$\bG(\mu) \overset{\phi}{\ra}$
$\bG(\mu) \ra$ 
$\coker \ \phi \ra 0$
and noting that $\lim \bG(\mu) = 0$, one obtains a sequence 
$0 \ra 0 \ra 0 \ra \Z/\nu\Z \ra 0$ which is not exact.  
On the other hand, the sequence 
%\begin{tikzcd}
%0 \arrow{r} &{\Z/\nu\Z} \arrow{r}{\phi} &{\lim^1 \bG(\mu)} \arrow{r}{\lim^1} &{\bG(\mu)} \arrow{r} &0
%\end{tikzcd}
$ 0 \ra $
$\Z/\nu\Z \ra$ 
\begin{tikzcd}%[sep=large]
{\lim^1 \bG(\mu)} \arrow{r}{\lim^1 \phi} &{\bG(\mu)}
\end{tikzcd}
$\ra 0$
is exact, so $\lim^1 \bG(\mu)$ contains a copy of $\Z/\nu\Z$ $\forall \ \nu:$ $(\mu,\nu) = 1$.\\
\endgroup%%------------------------------------<< 

To extend the applicability of the preceding considerations, replace \bAB by \bGR.  
Again, there is a functor $\lim:\bTOW(\bGR) \ra \bGR$ that sends \bG to $\lim \bG$.  
As for $\lim^1$\bG, it is the quotient $\prod\limits_n G_n/\sim$, where 
$
\begin{cases}
\ x^\prime  = \{x_n^\prime\}\\[-.1cm]
\ x^{\prime\prime} = \{x_n^{\prime\prime}\}
\end{cases}
$
are equivalent iff $\exists \ x = \{x_n\}$ such that 
%%----------------------------------------------------------------------------------------------46
$\forall \ n:$ $x_n^{\prime\prime} = x_n x_n^\prime f_n (x_{n+1}^{-1})$. 
 While not necessarily a group, $\lim^1$\bG is a pointed set with base point the equivalence class of $\{e_n\}$ and it is clear that 
$\lim^1:\bTOW(\bGR) \ra \bSET_*$ is a functor.

[Note: \  \ Put $X \ =  \ \prod\limits_n G_n$  \ $-$then the assignment 
$((g_0,g_1, \ldots),(x_0,x_1,\ldots))$ \  $\lra$ \ 
$(g_0 x_0 f_0(g_1^{-1}), g_1 x_1 f_1 (g_2^{-1}), \ldots)$ 
defines a left action of the group 
$\prod\limits_n G_n$  on the pointed set $X$.  
The stabilizer of the base point is $\lim \bG$ and the orbit space
 $\prod\limits_n G_n\backslash X$ is $\lim^1 \bG$.  
 For the definition and properties of $\lim^1$ ``in general'', consult 
 Bousfield-Kan\footnote[2]{\textit{SLN} \textbf{304} (1972), 305-308.}.]\\

\textbf{\small LEMMA} \ 
Let $* \ra \bG^{\prime} \ra \bG \ra \bG^{\prime\prime} \ra  *$ 
be an exact sequence in $\bTOW(\bGR)$ $-$then there is a natural exact sequence of groups and pointed sets
\[
* \ra \lim \bG^{\prime} \ra \lim \bG \ra \lim \bG^{\prime\prime} \ra  
{\lim}^1 \bG^\prime \ra {\lim}^1 \bG \ra {\lim}^1 \bG\pp \ra *.
\]

[Note: \  Specifically, the assumption is that $\forall \ n$, the sequence 
$* \ra G_n^{\prime} \ra$ 
$G_n \ra$ 
$G_n^{\prime\prime} \ra  *$ is exact in $\bGR$.]\\

\begingroup%%----------------------------------->>
\fontsize{9pt}{11pt}\selectfont
\textbf{\small EXAMPLE} \ 
Suppose that $\{G_n\}$ is a tower of fintely generated abelian groups $-$then $\lim^1 G_n$ is isomorphic to a group of the form 
Ext$(G,\Z)$, where $G$ is countable and torsion free.  
To see this, write $G_n^\prime$ for the torsion subgroup of $G_n$ and call 
$G_n^{\prime\prime}$ the quotient $G_n/G_n^\prime$.  
Since each $G_n^\prime$ is finite, 
$\lim^1 G_n^\prime = *$ $\implies$ 
$\lim^1 G_n \approx \lim^1 G_n^{\prime\prime}$.  
Assume, therefore, that the $G_n$ are torsion free.  
Let  \ 
$K_n = \ds\bigoplus\limits_{i \leq n} G_i = G_n \oplus K_{n-1}$ and define \ 
$K_n \ra K_{n-1}$ by $G_n \ra G_{n-1} \ra K_{n-1}$ 
\ on the first factor and by the identity on the second factor.  \ 
So, $\forall \ n$, $K_n \ra K_{n-1}$ is surjective, thus the sequence \ 
$0 \ra$ $\lim G_n \ra \lim K_n\ra \lim K_n/G_n \ra \lim^1 G_n \ra 0$ is exact.  \ 
Because $G_n$, $K_n$ and $K_n/G_n$ are free abelian, the sequence \ 
$0 \ra$ $\text{Hom}(K_n/G_n,\Z) \ra$ 
$\text{Hom}(K_n,\Z) \ra$ 
$\text{Hom}(G_n,\Z)$ $\ra 0$ is exact $\implies$ the sequence \ \ 
$0 \ra \text{colim Hom}(K_n/G_n,\Z) \ra$ 
$\text{colim Hom}(K_n,\Z) \ra$ 
$\text{colim Hom}(G_n,\Z)$ $\ra 0$ is exact 
$\implies$ the sequence \ 
$0 \ra$ $\text{Hom(colim Hom}(G_n,\Z),\Z) \ra$ 
$\text{Hom(colim Hom}(K_n,\Z),\Z) \ra$ \ 
$\text{Hom(colim Hom}(K_n/G_n,\Z),\Z) \ra$ \  
$\text{Ext(colim Hom}(G_n,\Z),\Z)   \ra$ 
$\text{Ext(colim Hom}(K_n,\Z),\Z) $ is exact
$\implies$ the sequence
$0 \ra \lim G_n \ra \lim K_n \ra$ 
$\lim K_n/G_n \ra$ $\text{Ext(colim Hom}(G_n,\Z),\Z)   \ra 0$
is exact (for colim Hom$(K_n,\Z) \approx$ $\ds\bigoplus\limits_n \text{Hom}(G_n,\Z),$ which is free).  
Consequently, 
$\lim^1 G_n \approx \text{Ext(colim Hom}(G_n,\Z),\Z)$, where $\colimx \Hom (G_n,\Z)$ is countable and torsion free.
\vspi
[Note: \  It follows that $\lim^1 G_n$ is divisible, hence if 
$\lim^1 G_n \neq *$, then on general grounds, there exist cardinals 
$\alpha$ and $\gamma (p)$ 
$(p \in \Pi)$: 
$\lim^1 G_n \approx$ 
$\ds\alpha \cdot \Q \hsx \oplus \hsx \bigoplus\limits_p \hsx \gamma(p) \cdot (\Z/p^\infty\Z)$.  
But here one can say more, viz $\alpha = 2^\omega$ and $\forall \ p$, $\gamma(p)$ is finite or $2^\omega$.]
\vspi
%%----------------------------------------------------------------------------------------------47
Huber-Warfield\footnote[3]{\textit{Arch. Math.} \textbf{33} (1979), 430-436.}
have shown that an abelian \mG is isomorphic to a $\lim^1 \bG$ for some \bG in $\bTOW(\bAB)$ iff $\Ext(\Q,G) = 0$.\\
\endgroup%%------------------------------------<< 

When is $\lim^1 \bG = *$?  An obvious sufficient condition is that the $f_n:G_{n+1} \ra G_n$ be surjective for every n.  
More generally, \bG is said to be 
\un{Mittag-Leffler}
\index{Mittag-Leffler} if 
$\forall \ n$ $\exists \ n^\prime \geq n$: 
$\forall \ n^{\prime\prime} \geq n^\prime$, 
$\text{im } (G_{n^{\prime}} \ra G_n) =$ 
$\text{im } (G_{n\pp} \ra G_n)$.\\

\index{Mittag-Leffler criterion}
\textbf{\small MITTAG-LEFFLER CRITERION} \ 
Suppose that \bG is Mittag-Leffler $-$then $\lim^1 \bG = *$.

[Note: \  There is a partial converse, viz.  if $\lim^1 \bG = *$ and if the $G_n$ are countable, then \bG is Mittag-Leffler (Dydak-Segal\footnote[2]{\textit{SLN} \textbf{688} (1978), 78-80.}).]\\

\begingroup%%----------------------------------->>
\fontsize{9pt}{11pt}\selectfont
\textbf{\small EXAMPLE} \ 
Fix a sequence $\mu_0 <\mu_1 \cdots$ of natural numbers $(\mu_0 > 1)$.  
Put $G_n = \prod\limits_{k \geq n} \Z/\mu_k\Z$ and let $G_{n+1} \ra G_n$ be the inclusion $-$then \bG is not Mittag-Leffler, yet $\lim^1 \bG = *$.\\
\endgroup%%------------------------------------<< 

\begingroup%%----------------------------------->>
\fontsize{9pt}{11pt}\selectfont
\textbf{\small FACT} \ 
Assume: $\lim^1 \bG \neq *$ and the $G_n$ are countable $-$then $\lim^1 \bG$ is uncountable.\\
\endgroup%%------------------------------------<< 

\label{8.1}
\begingroup%%----------------------------------->>
\fontsize{9pt}{11pt}\selectfont
\textbf{\small EXAMPLE} \ 
Let $X$ be a CW complex.  Suppose that 
$X_0 \subset X_1 \subset \cdots$ 
is an expanding sequence of subcomplexes of $X$ such that 
$X = \ds\bigcup\limits_n X_n$.  Fix a cofunctor 
$\sG: \Pi X \ra$  
$\bAB$ and put $\sG_n =$ 
$\restr{\sG}{X_n}$ $-$then $\forall \ q \geq 1$, there is an exact sequence 
$0 \ra \lim^1 H^{q-1}(X_n;\sG_n) \ra$ 
$H^q(X;\sG) \ra \lim H^q(X_n;\sG_n) \ra 0$ 
of abelian groups 
(Whitehead\footnote[3]{\textit{Elements of Homotopy Theory}, Springer Verlag (1978), 273-274.}).  
To illustrate, take $X = K(\Q,1)$ (realized as on p. \pageref{5.9}) and let 
$\sG: \Pi X \ra \bAB$ be the cofunctor corresponding to the usual action of 
$\Q$ on $\Q[\Q]$ 
(cf. p. \pageref{5.10}).  This data generates a short exact sequence
$0 \ra$ 
$\lim^1 H^1(\Z;\Q[\Q]) \ra H^2(\Q;\Q[\Q]) \ra$ 
$\lim H^2(\Z;\Q[\Q]) \ra 0$.  
The tower 
$H^1(\Z;\Q[\Q]) \leftarrow$ 
$H^1(\Z;\Q[\Q]) \leftarrow$ $\cdots$ is not Mittag-Leffler but 
$H^1(\Z;\Q[\Q])$ is countable, therefore 
$\lim^1 H^1(\Z;\Q[\Q])$ is uncountable.  
In particular: $H^2(\Q;\Q[\Q]) \neq 0$.\\
\endgroup%%------------------------------------<< 

\begingroup%%----------------------------------->>
\fontsize{9pt}{11pt}\selectfont
\textbf{\small FACT} \ 
Let $\{G_n\}$ be a tower of nilpotent groups.  Assume: $\forall \ n$, 
$\#(G_n) \leq \omega$ $-$then 
$\lim^1 G_n = *$ iff 
$\lim^1 G_n/[G_n,G_n] = *$.
\vspi
[For as noted above, in the presence of countability, 
$\lim^1 G_n/[G_n,G_n] = *$ $\implies$ 
$\{G_n/[G_n,G_n]\}$ is Mittag-Leffler.]\\
\endgroup%%------------------------------------<< 

\begin{proposition} \ %12
Let
$
\begin{cases}
\ \{X_n\} \\[-.1cm]
\ \{Y_n\}
\end{cases}
$
be two sequences of pointed spaces.  Suppose given pointed continuous functions 
$
\begin{cases}
\ \phi_n:X_n \ra X_{n+1}\\[-.1cm]
\ \psi_n:Y_{n+1} \ra Y_{n}
\end{cases}
$
\hspace{-.26cm}.  
Assume: The $\phi_n$ are closed cofibrations and the $\psi_n$ are pointed Hurewicz fibrations $-$then there is an exact sequence
%%----------------------------------------------------------------------------------------------48
%\begin{tikzcd}
%* \arrow[r]  &\lim^1[X_n,\Omega Y_n]  \arrow{r}{\iota} &{[\colimx X_n,\lim Y_n]} \arrow{r} &{\lim[X_n,Y_n]} %\arrow{r} &*
%\end{tikzcd}\\
\[
* \ra {\lim}^1[X_n,\Omega Y_n] \overset{\iota}{\ra} [\colim \ X_n,\lim Y_n] \ra \lim[X_n,Y_n] \ra *
\]
in $\bSET_*$ and $\iota$ is an injection.
\end{proposition}

[Write $X_\infty = \colimx X_n$ $\&$ $Y_\infty = \lim Y_n$.  Embedded in the data are arrows
$
\begin{cases}
\ \Phi_n:X_n \ra\\[-.1cm]
\ \Psi_n:Y_\infty \ra 
\end{cases}
$
$
\begin{aligned}
 X_\infty\\
 Y_{n}
\end{aligned}
$
with 
$
\begin{cases}
\ \Phi_{n+1} \circx \phi_n = \Phi_n\\[-.1cm]
\ \psi_n \circx \Psi_{n+1} = \Psi_n
\end{cases}
$
 and $\forall \ n$, an arrow 
 $[X_{n+1},Y_{n+1}] \ra $ 
 $[X_n,Y_n]$, viz. 
 $\abs{f} \ra$ 
 $[\psi_n \circx f \circx \phi_n]$.\\
 
Define $\xi_n:[X_\infty,Y_\infty] \ra [X_n,Y_n]$ by $\xi_n([f]) = [\Psi_n \circx f \circx \Phi_n]$.  
 Because the collection $\{\xi_n:[X_\infty,Y_\infty] \ra [X_n,Y_n]\}$ is a natural source, there exists a unique pointed map 
 $\xi_\infty:[X_\infty,Y_\infty] \ra \lim[X_n,Y_n]$ such that $\forall \ n$, the triangle
\begin{tikzcd}%[sep=large]
{[X_\infty,Y_\infty]} \arrow{r}{\xi_\infty} \arrow{rd}[swap]{\xi_n} &{\lim [X_n,Y_n]} \arrow{d}\\
&{[X_n,Y_n]}
\end{tikzcd}
 commutes.  
To prove that $\xi_\infty$ is surjective, take 
$\{[f_n]\} \in \lim [X_n,Y_n]$ 
$-$then $\forall \ n$, $\psi_n \circx f_{n+1} \circx \phi_n \simeq f_n$.  
Set $\bar{f}_0 = f_0$ and, proceeding inductively, assume that 
$\bar{f}_1 \in [f_1],\ldots, \bar{f}_n \in [f_n]$ have been found with 
$\psi_{k-1} \circx \bar{f}_k \circx \phi_{k-1} = \bar{f}_{k-1}$ $(1 \leq k \leq n)$.  
Choose a pointed homotopy $h_n:IX_n \ra Y_n$: 
$
\begin{cases}
\ h_n \circx i_0 = \psi_n \circx f_{n+1} \circx \phi_n\\[-.1cm]
\ h_n \circx i_1 = \bar{f}_n
\end{cases}
. \ 
$
Since $\psi_n$ is a pointed Hurewicz fibration, the commutative diagram
\begin{tikzcd}%[sep=large]
X_n \arrow{rr}{f_{n+1}\circ\phi_n} \arrow{d}[swap]{i_0} &&Y_{n+1} \arrow{d}{\psi_n}\\
{IX_n} \arrow{rr}[swap]{h_n} &&Y_n
\end{tikzcd}
admits a pointed filler $H_n:IX_n \ra Y_{n+1}$.  
Fix a retraction $r_n:IX_{n+1} \ra i_0X_{n+1} \cup I\phi_n(X_n)$ 
(cf. $\S 3$, Proposition 1) and specify a pointed continuous function 
$F_{n+1}:i_0X_{n+1} \cup I\phi_n(X_n) \ra Y_{n+1}$ by the prescription
 $
\begin{cases}
\ F_{n+1}(x_{n+1},0) = f_{n+1}(x_{n+1})\\[-.1cm]
\ F_{n+1}(\phi_n(x_n),t) = H_n(x_n,t)
\end{cases}
. \ 
$
Put $\bar{h}_n = \psi_n \circx F_{n+1} \circx r_n$ to get a commutative diagram
\begin{tikzcd}%[sep=large]
{i_0X_{n+1} \cup I\phi_n(X_n)} \arrow{r}{F_{n+1}} \arrow{d} &{Y_{n+1}} \arrow{d}{\psi_n}\\
{IX_{n+1}} \arrow{r}[swap]{\bar{h}_n} &Y_n
\end{tikzcd}
. \  Bearing in mind that $\phi_n$ is a closed cofibration, this diagram has a pointed filler $\overline{H}_{n+1}:IX_{n+1} \ra Y_{n+1}$ (cf. $\S 4$, Proposition 12).  Finally, to push the induction forward, let $\bar{f}_{n+1} = \overline{H}_{n+1} \circx i_1$.  Conclusion:  There exists a pointed continuous function $\bar{f}_\infty: X_\infty \ra Y_\infty$ such that $\xi_\infty([\bar{f}_\infty]) = \{[f_n]\}$, i.e., $\xi_\infty$ is surjective.
 
As for the kernel of $\xi_\infty$, it consists of those $[f]:$ $\forall \ n$, $\Psi_n \circx f \circx \Phi_n$ is nullhomotopic.  
Thus there are pointed homotopies 
$\Xi_n:IX_n \ra$ $Y_n$ such that 
$\Xi_n \circx i_0 = 0_n$ $\&$
 $\Xi_n \circx i_1 =$ $\Psi_n \circx f \circx \Phi_n$ with
$\psi_n \circx \Xi_{n+1} \circx I\phi_n \circx i_0 = 0_n$ $\&$ $\psi_n \circx \Xi_{n+1} \circx I\phi_n \circx i_1 = \Psi_n \circx f \circx \Phi_n$, where $0_n$ is the zero morphism $X_n \ra Y_n$.  To define $\eta_\infty: \ker \xi_\infty \ra \lim^1[X_n,\Omega Y_n]$, let $\sigma_{n,f}:X_n \ra \Omega Y_n$ be the pointed continuous function given by
\[
\sigma_{n,f}(x_n,t) = 
\begin{cases}
\ \Xi_n(x_n,2t) \hspace{3.35cm}   (0 \leq t \leq 1/2)\\[-.1cm]
\ \psi_n \circx \Xi_{n+1}(\phi_n(x_n),2 - 2t) \hspace{0.7cm} (1/2 \leq t \leq 1).
\end{cases}
\hspace{-.26cm}.
\]
%%----------------------------------------------------------------------------------------------49
The $\sigma_{n,f}$ determine a string in \ $\prod\limits_n[X_n,\Omega Y_n]$ or still, an element of \ $\lim^1[X_n,\Omega Y_n]$, call it $[\sigma_f]$.  Definition:  $\eta_\infty([f]) = [\sigma_f]$.  One can check that $\eta_\infty$ does not depend on the choice of $\Xi_n$ and is independent of the choice of $f \in [f]$.  Claim:  $\eta_\infty$  is bijective.  To verify, e.g., injectivity, suppose that $\eta_\infty ([f^\prime]) = \eta_\infty ([f^{\prime\prime}])$ $-$then there exists a string $\{[\sigma_n]\} \in \prod\limits_n [X_n,\Omega Y_n]$: $\forall \ n$,\\
\[
\begin{cases}
\ \sigma_n(x_n,3t) \hspace{3.5cm}  (0 \leq t \leq 1/3)\\[-.1cm]
\ \sigma_{n,f^\prime}(x_n,3t - 1) \hspace{2.5cm} (1/3 \leq t \leq 2/3)\\[-.1cm]
\ \psi_n \circx \sigma_{n+1}(\phi(x_n),3- 3t) \hspace{1.05cm} (2/3 \leq t \leq 1)
\end{cases}
\]
%^
represents $\sigma_{n,f\pp}$.  In addition, the formulas\\
%^
\[
\begin{cases}
\ \Xi_n^\prime (x_n,1 - 3t) \hspace{3.0cm} (0 \leq t \leq 1/3)\\[-.1cm]
\ \sigma_n(x_n,2 - 3t) \hspace{3.05cm} (1/3 \leq t \leq 2/3)\\[-.1cm]
\ \Xi_n^{\prime\prime}(x_n,3t -2) \hspace{3.02cm} (2/3 \leq t \leq 1)
\end{cases}
\]
%^
define a pointed homotopy $H_n:IX_n \ra Y_n$ having the property that 
$H_n \circx i_0 =$ $\Psi_n \circx f^\prime \circx \Phi_n$ $\&$ 
$H_n \circx i_1 =$ $\Psi_n \circx f^{\prime\prime} \circx \Phi_n$.  
Arguing as before, construct pointed homotopies 
$\overline{H}_n:IX_n \ra$ $Y_n$ such that 
$\overline{H}_n \circx i_0 =$ $\Psi_n \circx f^\prime \circx \Phi_n$ $\&$ 
$\overline{H}_n \circx i_1 =$ $\Psi_n \circx f^{\prime\prime} \circx \Phi_n$ with 
$\psi_n \circx \overline{H}_{n+1} \circx I\phi_n =$ $\overline{H}_n$.
The $\overline{H}_n$ combine and induce a pointed homotopy 
$\overline{H}_\infty: IX_\infty \ra$ $Y_\infty$ between 
$f^\prime$ and $f^{\prime\prime}$, i.e., $\eta_\infty$ is injective.]\\

\label{11.12}
Application:
Let $\{X_n\}$ be a sequence of pointed spaces.  Suppose given pointed continuous functions $\phi_n:X_n \ra X_{n+1}$ such that $\forall \ n$, $\phi_n$ is a closed cofibration $-$then for any pointed space $Y$, there is an exact sequence
% \begin{tikzcd}
%&* \arrow{r} &{\lim^1 [\Sigma X_n,Y]} \arrow{r}{\iota} &{[\text{colim}X_n,Y]} \arrow{r} &{\lim [X_n,Y]} \arrow{r} &*\\
%\end{tikzcd}
\[
* \ra {\lim}^1 [\Sigma X_n,Y] \overset{\iota}{\ra} [\colim \ X_n,Y] \ra \lim [X_n,Y] \ra *
\]
in \bSET$_*$ and $\iota$ is an injection.\\

\label{9.36}
\begingroup%%----------------------------------->>
\fontsize{9pt}{11pt}\selectfont
\textbf{\small EXAMPLE} \ 
Fix an abelian group $\pi$.  Let $(X,x_0)$ be a pointed CW complex.  
Suppose that $x_0 \in X_0 \subset X_1 \subset \cdots$ is an expanding sequence of subcomplexes of $X$ such that 
$X = \ds\bigcup\limits_n X_n$ $-$then $\forall \ q \geq 1$,  there is an exact sequence
$0 \ra$ 
$\lim^1 \widetilde{H}^{q-1}(X_n;\pi) \ra$ 
$\widetilde{H}^{q}(X_;\pi) \ra$ 
$\lim \widetilde{H}^{q}(X_n;\pi)\ra 0$ of abelian groups.  
Example:  $\forall \ q \geq 1$, ${H}^{q}(\Z/p^\infty \Z,n) \approx$ $\lim {H}^{q}(\Z/p^k \Z,n)$.
\vspi
[In the above, substitute $Y = K(\pi,q)$.]\\
\endgroup%%------------------------------------<<

\label{5.54}
\label{11.3}
\begingroup%%----------------------------------->>
\fontsize{9pt}{11pt}\selectfont
\textbf{\small LEMMA} \ 
Let $X$ be a pointed finite CW complex.  Let $K$ be a pointed connected CW complex.  Assume:  The homotopy groups of $K$ are finite $-$then the pointed set $[X,K]$ is finite.
\vspi
[This result is contained in obstruction theory but one can also give a direct inductive proof.]\\
\endgroup%%------------------------------------<<

%%----------------------------------------------------------------------------------------------50
\label{9.12} %dmc mnft
\begingroup%%----------------------------------->>
\fontsize{9pt}{11pt}\selectfont
\textbf{\small EXAMPLE} \ 
Let $(X,x_0)$ be a pointed CW complex.   
Suppose that 
$x_0 \in X_0 \subset X_1 \subset \cdots$ 
is an expanding sequence of finite subcomplexes of $X$ such that 
$X = \ds\bigcup\limits_n X_n$.  
Let $K$ be a pointed connected CW complex.   
Assume:  The homotopy groups of $K$ are finite $-$then the natural map 
$\pi_X:[X,K] \ra$ 
$\lim[X_n,K]$ is bijective.  
In fact, surjectivity is automatic, so injectivity is what's at issue.  
For this, consider the natural map 
$\pi_{IX}:[IX,K] \ra$ 
$\lim[i_0X \cup IX_n \cup i_1 X,K]$ 
and the obvious arrows 
$i_0,i_1:\lim[i_0 X \cup IX_n \cup i_1 X,K] \ra$ 
$[X,K]$.  
Since $i_0 \circx \pi_{IX} = i_1 \circx \pi_{IX}$  and since $\pi_{IX}$ is surjective, $i_0 = i_1$.  
That $\pi_{IX}$ is injective is thus a consequence of the following claim.
\vspi
Claim:  If $\pi_X([f_0]) = \pi_X([f_1])$, then there exists $[F] \in \lim[i_0X \cup IX_n \cup i_1 X,K]$: 
$
\begin{cases}
\ \abs{f_0} = i_0([F])\\
\ \abs{f_1} = i_1([F])\\
\end{cases}
.
$
\vspi
[Let 
$i_0^n,i_1^n:[i_0X \cup IX_n \cup i_1 X,K] \ra [X,K]$ be the obvious arrows.  
For each $n$, there is at least one 
$[F_n] \in [i_0X \cup IX_n \cup i_1 X,K]$:
$
\begin{cases}
\ [f_0] = i_0^n([F_n])\\
\ [f_1] = i_1^n([F_n])\\
\end{cases}
. \ 
$
Denote by $I_n$ the subset of 
$[i_0X \cup IX_n \cup i_1 X,K]$ consisting of all such $[F_n]$ $-$then, from the lemma, $I_n$ is finite, hence $\lim I_n \neq \emptyset$.]
\vspi
[Note: \  The $\Sigma X_n$ are finite CW complexes, therefore the $[\Sigma X_n,K]$ are finite groups, so $\lim^1 [\Sigma X_n,K] = *$.  
But this only means that the kernel of $\pi_X$ is $[0]$.]\\
\endgroup%%------------------------------------<<

\label{5.55a}

Application:  Let $\{Y_n\}$ be a sequence of pointed spaces.  Suppose given pointed continuous functions 
$\psi_n:Y_{n+1} \ra Y_n$ such that 
$\forall \ n$, $\psi_n$ is a pointed Hurewicz fibration $-$then for any pointed space $X$, there is an exact sequence\\
%
%\begin{tikzcd}%[sep=large]
%\indent * \arrow[r]  &\lim^1[X_n,\Omega Y_n]  \arrow{r}{\iota} &{[X_n,\lim Y_n]} \arrow{r} &{\lim[X,Y_n]} \arrow{r} &*
%\end{tikzcd}\\
%
\[
* \ra {\lim}^1[X_n,\Omega Y_n]  \overset{\iota}{\ra} [X,\lim Y_n] \ra \lim[X,Y_n] \ra *
\]
in $\bSET_*$ and $\iota$ is an injection.

[Note: \   The exact sequence
%\begin{tikzcd}%[sep=large]
%* \arrow[r]  &\lim^1\pi_{q+1}(Y_n) \arrow{r}{\iota} &\pi_q(\lim Y_n) \arrow{r} &\lim \pi_q(Y_n) \arrow{r} &*
%\end{tikzcd}
$* \ra \lim^1\pi_{q+1}(Y_n) \overset{\iota}{\ra} \pi_q(\lim Y_n) \ra \lim \pi_q(Y_n) \ra *$
of pointed sets is a special case (take $X = \bS^q$).]\\

\begingroup%%----------------------------------->>
\fontsize{9pt}{11pt}\selectfont
\textbf{\small EXAMPLE} \ 
For each $n$, put $Y_n = \bS^1$ and let $\psi_n:Y_{n+1} \ra Y_n$ be the squaring map
$
\begin{cases}
\ \bS^1 \ra \bS^1\\
\ s \ \ \mapsto s^2\\
\end{cases}
$
$-$then $\lim \pi_1(Y_n) = 0$ but $\lim^1(Y_n) \approx \widehat{\Z}_2/\Z$, the 2-adic integers mod $\Z$.\\
\endgroup%%------------------------------------<<

\begingroup%%----------------------------------->>
\fontsize{9pt}{11pt}\selectfont
\textbf{\small EXAMPLE} \ 
Let $\bpi = \{\pi_n\}$ be a tower of abelian groups.  
Assume: $\bpi$ is Mittag-Leffler $-$then $\forall \ q \geq 1$, 
$K(\lim \bpi,q) = \lim K(\pi_n,q)$, 
so for any pointed CW complex $(X,x_0)$, there is an exact sequence
$0 \ra$ 
$\lim^1 \widetilde{H}^{q-1} (X;\pi_n) \ra$ 
$\widetilde{H}^q(X;\lim \bpi) \ra$ 
$\lim \widetilde{H}^q(X;\pi_n) \ra 0$
of abelian groups.\\
\endgroup%%------------------------------------<<

\indent Given a pointed path connected space $X$, let $P_\infty X = \lim P_n X$ $-$then $\forall \ q \geq 0$, $\pi_q(P_\infty X) \approx \lim \pi_q(P_n X) \approx \pi_q(P_q X)$.  
Proof: The relevant $\lim^1$ term vanishes.\\

\begin{proposition} \ %13
The canonical arrow $X \ra P_\infty X$ is a weak homotopy equivalence.
\end{proposition}

%%----------------------------------------------------------------------------------------------51
[For each $n$, there is an inclusion $X \ra X[n]$, a projection $P_\infty X \ra P_n X$, and a pointed homotopy equivalence $X[n] \ra P_n X$.  Consider the associated commutative diagram
\begin{tikzcd}%[sep=large]
X \arrow{r} \arrow{d} &{P_\infty X} \arrow{d}\\
X[n] \arrow{r} &P_n(X)
\end{tikzcd}
, recalling that $\pi_n(X) \approx \pi_n(X[n])$.]\\
\vspace{0.25cm}

\label{9.64}
\begingroup%%----------------------------------->>
\fontsize{9pt}{11pt}\selectfont
\textbf{\small FACT} \ 
Let $\{X_n, f_n: X_{n+1} \ra X_n\}$ be a tower in \bTOP.  
Assume:  The $X_n$ are CW spaces and the $f_n$ are Hurewicz fibrations $-$then $\lim X_n$ is a CW space iff all but finitely many of the $f_n$ are homotopy equivalences.
\vspi
[Necessity:  
If infinitely many of the $f_n$ are not homotopy equivalences, then $\lim X_n$ is not numerably contractible.
\vspi
Sufficiency:  
If all of the $f_n$ are homotopy equivalences, then $X_0$ and $\lim X_n$ have the same homotopy type (cf. p. \pageref{5.11}).]\\
\endgroup%%------------------------------------<< 

\begingroup%%----------------------------------->>
\fontsize{9pt}{11pt}\selectfont
Application:  Suppose that $X$ is a pointed connected CW space $-$then the canonical arrow $X \ra P_\infty X$ is a homotopy equivalence iff $X$ has finitely many nontrivial homotopy groups.\\
\endgroup%%------------------------------------<< 

\index{Theorem Whitehead}
\textbf{\small WHITEHEAD THEOREM} \ 
Suppose that $X$ and $Y$ are path connected topological spaces.\\
\indent \indent (1) \quadx Let $f:X \ra Y$ be an $n$-equivalence $-$then $f_*:H_q(X) \ra H_q(Y)$ is bijective for $1 \leq q < n$ and surjective for $q = n$.\\
%^
\indent \indent (2) \quadx Suppose in addition that $X$ and $Y$ are simply connected.   Let $f:X \ra Y$ be a continuous function such that  $f_*:H_q(X) \ra H_q(Y)$ is bijective for $1 \leq q < n$ and surjective for $q = n$ 
$-$then $f$ is an $n$-equivalence.

[The condition on $f_*$ amounts to requiring that $H_q(M_f,i(X)) = 0$ for $q \leq n$, thus the result follows from the relative Hurewicz theorem.]\\
\vspace{0.25cm}

\begingroup%%----------------------------------->>
\fontsize{9pt}{11pt}\selectfont
\label{5.23}
\textbf{\small EXAMPLE} \ 
Let $X$ be a pointed connected CW space $-$then the inclusion $X \ra X[n]$ is an $(n+1)$-equivalence, hence there are bijections $H_q(X) \approx H_q(X[n])$ $(q \leq n)$ and a surjection $H_{n+1}(X) \ra H_{n+1}(X[n])$.  
So, if $X$ is abelian and if the $\pi_q(X)$ are finitely generated $\forall \ q$, 
then the $H_q(X)$ are finitely generated $\forall \ q$ (cf. p. \pageref{5.12}).\\
\endgroup%%------------------------------------<< 

\index{Theorem, Suspension Theorem}
\begingroup%%----------------------------------->>
\fontsize{9pt}{11pt}\selectfont
\textbf{\small EXAMPLE \ (\un{Suspension Theorem})} \ 
Suppose that $X$ is nondegenerate and $n$-connected.  Let $K$ be a pointed CW complex 
$-$then the suspension map $[K,X] \ra [\Sigma K, \Sigma X]$ is bijective if $\dim K \leq 2n$ and surjective if $\dim K \leq 2n + 1$.  In fact, the arrow of adjunction 
$e:X \ra \Omega\Sigma X$ induces an isomorphism 
$H_q(X) \ra H_q(\Omega\Sigma X)$ for 
$0 \leq q \leq 2n +1$ (cf. p. \pageref{5.13}), therefore by the Whitehead theorem $e$ is a $(2n+1)$-equivalence.  So, if $\dim K$ is finite and if $n \geq 2 + \dim K$, then 
$[\Sigma^n K,\Sigma^n X] \approx$ 
$[\Sigma^{n+1}K,\Sigma^{n+1}X]$.\\
\endgroup%%------------------------------------<< 

%%----------------------------------------------------------------------------------------------52
A continuous function $f:X \ra Y$ is said to be a 
\un{homology equivalence}
\index{homology equivalence} 
if \ $\forall \ n \geq 0$, \qquad
$f_*:H_n(X) \ra H_n(Y)$ 
is an isomorphism.  \ 
Example:  Consider the coreflector 
$k: \bTOP \ra \bCG$ $-$then for every topological space $X$, the identity map 
$kX \ra X$ is a homology equivalence.\\

\begingroup%%----------------------------------->>
\fontsize{9pt}{11pt}\selectfont
\label{5.33}
\textbf{\small EXAMPLE} \ 
A homology equivalence $f:X \ra Y$ need not be a weak homotopy equivalence.  One can take, e.g., \mX to be Poincar\'e's homology 3-sphere $\bS^3/\bSL(2,5)$ and $Y = \bS^3$.  There is a homology equivalence $f:X \ra Y$ obtained by collapsing the 2-skeleton of \mX to a point which, though, is not a weak homotopy equivalence, the fundamental group of \mX being 
$\bSL(2,5)$.  Eight different descriptions of \mX have been examined by 
Kirby-Scharlemann\footnote[2]{In: \textit{Geometric Topology}, J. Cantrell (ed.), Academic Press (1979), 113-146.}.\\
\endgroup%%------------------------------------<< 

\index{Theorem Whitehead (bis)}
\textbf{\small WHITEHEAD THEOREM (bis)} \ 
Suppose that $X$ and $Y$ are path connected topological spaces.\\
\indent \indent (1) \quadx Let $f:X \ra Y$ be a weak homotopy equivalence $-$then $f$ is a homology equivalence.

[Note: \  It is a corollary that in general a weak homotopy equivalence is a homology equivalence.]

\indent \indent (2) \quadx Suppose in addition that $X$ and $Y$ are simply connected.   Let $f:X \ra Y$ be a homology equivalence $-$then $f$ is a weak homotopy equivalence.\\

Consequently, if $X$ and $Y$ are simply connected topological spaces that are dominated in homotopy by CW complexes, then a continuous function $f:X \ra Y$ is a homotopy equivalence iff it is a homology equivalence.\\

\begingroup%%----------------------------------->>
\fontsize{9pt}{11pt}\selectfont
The following familiar remarks serve to place this result in perspective.\\
\indent \indent (1) \ There exist path connected topological spaces $X$ and $Y$ such that $\forall \ n$: $\pi_n(X)$ is isomorphic to $\pi_n(Y)$ but $\exists \ n$: $H_n(X)$ is not isomorphic to $H_n(Y)$.
\\
\indent \indent (2) \ There exist simply connected topological spaces $X$ and $Y$ such that $\forall \ n$: $H_n(X)$ is isomorphic to $H_n(Y)$ but $\exists \ n$: $\pi_n(X)$ is not isomorphic to $\pi_n(Y).$
\\
\indent \indent (3) \ There exist path connected topological spaces $X$ and $Y$ admitting a homology equivalence $f:X \ra Y$ 
with the property that $f_*:\pi_1(X) \ra \pi_1(Y)$ is an isomorphism, yet $f$ is not a weak homotopy equivalence.
\vspi
[Note: \  Recall too that there exist topological spaces $X$ and $Y$ such that 
$\forall \ n$: $H_n(X)$ is isomorphic to $H_n(Y)$ and 
$\forall \ n$: $\pi_n(X,x_0)$ is isomorphic to $\pi_n(Y,y_0)$ 
$(\forall \ x_0 \in X$, $\forall \ y_0 \in Y)$, yet $X$ and $Y$ do not have the same homotopy type.  
Example:
$
\begin{cases}
\ X = \{0\} \cup \{1/n: n \geq 1\} \\
\ Y= \{0\} \cup \{n: n \geq 1\} 
\end{cases}
.]
$
\\
\endgroup%%------------------------------------<< 

%%----------------------------------------------------------------------------------------------53
\begingroup%%----------------------------------->>
\fontsize{9pt}{11pt}\selectfont
\textbf{\small EXAMPLE} \ 
There exists a sequence $X_1, X_2, \ldots$ of simply connected CW complexes $X_n$ having isomorphic integral singular cohomology rings such that $\forall \ n^\prime \neq n^{\prime\prime}$, 
the homotopy types of $X_{n^\prime}$ $\&$ $X_{n^{\prime\prime}}$ are distinct 
(Body-Douglas\footnote[2]{\textit{Topology} \textbf{13} (1974), 209-214.}).\\
\endgroup%%------------------------------------<< 

\begingroup%%----------------------------------->>
\fontsize{9pt}{11pt}\selectfont
\textbf{\small EXAMPLE} \ 
Let $X$ be a pointed connected CW space  $-$then $\Sigma X$ is contractible iff 
$H_1(\pi,1) =$ 
$0 =$ 
$H_2(\pi,1)$ $(\pi = \pi_1(X))$ and $H_q(X) = 0$ $(q \geq 2)$.\\
\endgroup%%------------------------------------<< 

\index{stable splitting}
\begingroup%%----------------------------------->>
\fontsize{9pt}{11pt}\selectfont
\textbf{\small EXAMPLE \ (\un{Stable Splitting})} \ 
Let \mG be a finite abelian group $-$then there exist positive integers $T$ and $t$ such that 
$\Sigma^T K(G,1)$ has the pointed homotopy type of a wedge 
$X_1 \vee \cdots \vee X_t$, where the $X_i$ are pointed simply connected CW spaces.  
For let 
$G =$ 
$G(p_1) \oplus \cdots \oplus G(p_n)$ be the primary decomposition of $G$.  
Since the arrow 
$K(G(p_1),1) \vee \cdots \vee K(G(p_n),1) \ra$ 
$K(G(p_1),1) \times \cdots \times K(G(p_n),1) =$ $K(G,1)$ 
is a homology equivalence, its suspension is a pointed homotopy equivalence, thus one can assume that $G$ is 
$p$-primary, say 
$G =$ 
$\Z/p^{e_1}\Z \oplus \cdots \oplus \Z/p^{e_r}\Z$, so 
$K(G,1) =$ 
$\ds \prod\limits_1^r K(\Z/p^{e_i}\Z,1)$.  
Accordingly, thanks to the Puppe formula and the fact that 
$\Sigma(X \# Y) \approx$ 
$\Sigma X \# Y \approx$ 
$X \# \Sigma Y$, it suffices to consider $K(\Z/p^e\Z,1)$.
\vspi
Claim: There exist pointed simply connected CW spaces 
$X_1, \ldots, X_{p-1}$ and a pointed homotopy equivalence 
$\Sigma K(\Z/p^e\Z,1) \ra$ $X_1 \vee \cdots \vee X_{p-1}$.
\vspi
[A generator of the multiplicative group of units in $\Z/p\Z$ defines a pointed homotopy equivalence 
$K(\Z/p^e\Z,1) \ra K(\Z/p^e\Z,1))$.]\\
\endgroup%%------------------------------------<< 

The rather restrictive assumption that
$
\begin{cases}
\ \pi_1(X) = 0 \\[-.1cm]
\ \pi_1(Y) = 0
\end{cases}
$
is not necessary in order to guarantee that a homology equivalence $f:X \ra Y$ is a weak homotopy equivalence.  For example,
$
\begin{cases}
\ X \\[-.1cm]
\ Y
\end{cases}
$
abelian will do and in fact one can get away with considerably less.

Notation:  Given a group $G$, let $\Z[G]$ be its integral group ring and $I[G] \subset \Z[G]$ the augmentation ideal.  
Given a $G$-module $M$, let $M_G$ be its group of coinvariants, i.e., the quotient $M/I[G] \cdot M$ or still, $H_0(G;M)$.

[Note: \  In this context, ``$G$-module'' means left $G$-module.  If $K$ is a normal subgroup of $G$, then the action of $G$ on $M$ induces an action $G/K$ on $M_K$ and $M_G \approx (M_K)_{G/K}$.]\\

\index{Fundamental Exact Sequence}
\textbf{\small FUNDAMENTAL EXACT SEQUENCE} \ 
Fix a $G$-module $M$.  Let $K$ be a normal subgroup of $G$ $-$then there is an exact sequence
\[
H_2(G;M) \ra H_2(G/K;M_K) \ra H_1(K;M)_{G/K} \ra H_1(G;M) \ra H_1(G/K;M_K) \ra 0.
\]

[The LHS spectral sequence reads: 
$E_{p,q}^2 \approx$ $H_p(G/K;H_q(K;M))$ $\implies$ 
$H_{p+q}(G;M)$.  Explicate the associated five term exact sequence
%\begin{tikzcd}[ sep=small]
%H_2(G;M) \ar{r} &E_{2,0}^2  \ar{r}{d^2} &E_{0,1}^2  \ar{r} &H_1(G;M)  \ar{r} &E_{1,0}^2 \ar{r} &0
%\end{tikzcd}\\
$H_2(G;M) \ra$ 
$E_{2,0}^2  \overset{d^2}{\lra}$  
$E_{0,1}^2  \ra$ 
$H_1(G;M)  \ra$ 
$E_{1,0}^2 \ra 0$
.]\\

%%----------------------------------------------------------------------------------------------54
\label{11.32}
Application:  Let $K$ be a normal subgroup of $G$ $-$then there is an exact sequence 
$H_2(G) \ra  H_2(G/K) \ra  K/[G,K] \ra  H_1(G) \ra H_1(G/K) \ra  0$.

[Specialize the fundamental exact sequence and take $M = \Z$ (trivial $G$-action).  
Observe that the arrows
$
\begin{cases}
\ H_1(G) \ra H_1(G/K) \\
\ H_2(G) \ra H_2(G/K)
\end{cases}
$
are induced by the projection $G \ra G/K$.]\\

\begingroup%%----------------------------------->>
\fontsize{9pt}{11pt}\selectfont
Using a superscript to denote the ``invariants'' functor, the fundamental exact sequence in cohomology is
$0 \ra H^1(G/K;M^K) \ra$ 
$H^1(G;M) \ra$ 
$H^1(K;M)^{G/K} \ra$ 
$H^2(G/K;M^K)  \ra$ 
$H^2(G;M)$.\\
\endgroup%%------------------------------------<< 

Notation:  Given a group $G$, let \ $\Gamma^0(G) \supset \Gamma^1(G) \supset \ldots$ \ be its descending central series, so $\Gamma^{i+1}(G) = [G,\Gamma^i(G)].$  In particular:  $\Gamma^0(G) = G$, $\Gamma^1(G) = [G,G]$ and $G$ is \un{nilpotent}
\index{nilpotent} 
if there exists a $d$ : $\Gamma^d(G) = \{1\}$, the smallest such $d$ being its 
\un{degree of nilpotency}: 
\index{degree of nilpotency}
$\nil G$. \\

\label{8.18}
\begingroup%%----------------------------------->>
\fontsize{9pt}{11pt}\selectfont
\label{5.25}
\textbf{\small FACT} \ 
Let $G$ be a nilpotent group $-$then $G$ is finitely generated iff $G/[G,G]$ is finitely generated.\\
\endgroup%%------------------------------------<< 

\label{8.21}
\begingroup%%----------------------------------->>
\fontsize{9pt}{11pt}\selectfont
\textbf{\small EXAMPLE} \ 
Let $G$ be a nilpotent group $-$then $G$ is finitely generated iff $\forall \ q \geq 1$, $H_q(G)$ is finitely generated.  
For suppose that $G$ is finitely generated.  
Case 1:  $\nil G \leq 1$.  In this situation, $G$ is abelian and the assertion is true 
(cf. p. \pageref{5.14}).  
Case 2:  $\nil G > 1$.  Argue by induction, using the LHS spectral sequence 
$E_{p,q}^2 \approx$ 
$H_p(G/\Gamma^i(G);H_q(\Gamma^i(G)/\Gamma^{i+1}(G)))$ $\implies$ 
$H_{p+q}(G/\Gamma^{i+1}(G))$.  
To discuss the converse, note that $H_1(G) \approx G/[G,G]$ and quote the preceding result.\\
\endgroup%%------------------------------------<< 

\begingroup%%----------------------------------->>
\fontsize{9pt}{11pt}\selectfont
It is false in general that a subgroup of a finitely generated group is finitely generated.  Example:  Let $G$ be the free group on two symbols and consider $[G,G]$.\\
\endgroup%%------------------------------------<< 

\label{8.24} %dmc mnft
\begingroup%%----------------------------------->>
\fontsize{9pt}{11pt}\selectfont
\textbf{\small FACT} \ 
Suppose that $G$ is a finitely generated nilpotent group $-$then every subgroup of $G$ is finitely generated.\\
\endgroup%%------------------------------------<< 

\begingroup%%----------------------------------->>
\fontsize{9pt}{11pt}\selectfont
\textbf{\small FACT} \ 
Suppose that $G$ is a finitely generated nilpotent group $-$then $G$ is finitely presented.
[The class of finitely presented groups is closed with respect to the formation of extensions.]\\
\endgroup%%------------------------------------<< 

\begingroup%%----------------------------------->>
\fontsize{9pt}{11pt}\selectfont
Notation:  Given a group $G$, $G_{\tor}$ is its subset of elements of finite order.

\label{8.12}
[Note: \  $G_{\tor}$ need not be a subgroup of $G$ (consider $G = \Z/2\Z * \Z/2\Z)$ but will be if $G$ is nilpotent 
(since $\nil G \leq d$ and $y^m = e$ $\implies$ $(xy)^{m^d} = x^{m^d}$).]\\
\endgroup%%------------------------------------<< 

\begingroup%%----------------------------------->>
\fontsize{9pt}{11pt}\selectfont
\textbf{\small FACT} \ 
Suppose that $G$ is a finitely generated nilpotent group.  
Assume: $G$ is torsion $-$then $G$ is finite.\\
\endgroup%%------------------------------------<< 

\begingroup%%----------------------------------->>
\fontsize{9pt}{11pt}\selectfont
Application:  If $G$ is a finitely generated nilpotent group, the $G_{\tor}$ is a finite nilpotent normal subgroup.\\
\endgroup%%------------------------------------<< 

%%----------------------------------------------------------------------------------------------55
\begin{proposition} \ %14
Let $f:G \ra K$ be a homomorphism of groups.  Assume: 
(i) $f_*:H_1(G) \ra H_1(K)$ is bijective and 
(ii) $f_*:H_2(G) \ra H_2(K)$ is surjective. $-$then 
$\forall \ i \geq 0$, the induced map $G/\Gamma^i(G) \ra K/\Gamma^i(K)$ is an isomorphism.
\end{proposition}

[The assertion is trivial if $i = 0$ and holds by assumption if $i = 1$.  Fix $i > 1$ and proceed by induction.  There is a commutative diagram\\
\[
\begin{tikzcd}[sep=small]
H_2(G) \ar{dd} \ar{r} 
&H_2(G/\Gamma^i(G)) \ar{dd} \ar{r} 
&\Gamma^i(G)/\Gamma^{i+1}(G)  \ar{dd} \ar{r} 
&H_1(G) \ar{dd} \ar{r} 
&H_1(G/\Gamma^i(G)) \ar{r}\ar{dd} &0 \ar{dd}\\
\\
H_2(K) \ar{r}           
& H_2(K/\Gamma^i(K)) \ar{r}            
&\Gamma^i(K)/\Gamma^{i+1}(K)   \ar{r}          
&H_1(K) \ar{r}            
& H_1(K/\Gamma^i(K)) \ar{r} 
&0\\
\end{tikzcd}
\]
with exact rows, hence, by the five lemma, $\Gamma^i(G)/\Gamma^{i+1}(G) \approx \Gamma^i(K)/\Gamma^{i+1}(K)$.  But then from
\[
\begin{tikzcd}[sep=large]
1  \ar{r} 
&\Gamma^i(G)/\Gamma^{i+1}(G)  \ar{d} \ar{r} 
&G/\Gamma^{i+1}(G) \ar{d} \ar{r} 
&G/\Gamma^i(G) \ar{r}\ar{d} 
&1 \\
1\ar{r}            
&\Gamma^i(K)/\Gamma^{i+1}(K)   \ar{r}           
&K/\Gamma^{i+1}(K) \ar{r}          
 &K/\Gamma^i(K) \ar{r} &1
\end{tikzcd},
\]
one concludes that $G/\Gamma^{i+1}(G) \approx K/\Gamma^{i+1}(K)$.]\\

\label{5.15}
\label{8.22}
\label{8.23}  %dmc mnft
\label{8.33}  %dmc mnft
Application:  Let $f:G \ra K$ be a homomorphism of nilpotent groups.  
Assume: 
(i) $f_*:H_1(G) \ra H_1(K)$ is bijective and 
(ii) $f_*:H_2(G) \ra H_2(K)$ is surjective. $-$then $f$ is an isomorphism.\\

\label{5.27c}
\label{8.27}
\label{8.28} %dmc mnft
Let $G$ and $\pi$ be groups.  
Suppose that $G$ operates on $\pi$, i.e., suppose given a homomorphism $\chi:G \ra \Aut\pi$.  
Put 
$\Gamma_\chi^0(\pi) = \pi$ and, via recursion, write 
$\Gamma_\chi^{i+1}(\pi)$ for the subgroup of $\pi$ generated by the 
$\alpha(\chi(g)\alpha_i)\alpha^{-1}\alpha_i^{-1}$ 
$(\alpha \in \pi, \alpha_i \in \Gamma_\chi^{i}(\pi))$, where $g \in G$ $-$then 
$\Gamma_\chi^{i}(\pi)$ is a $G$-stable normal subgroup of $\pi$ containing 
$\Gamma_\chi^{i+1}(\pi)$.  
The quotient 
$\Gamma_\chi^{i}(\pi) / \Gamma_\chi^{i+1}(\pi)$ is abelian and the induced action of $G$ is trivial.  
One says that $G$ operates nilpotently 
\index{operate nilpotently} 
on $\pi$ or that $\pi$ is 
\un{$\chi$-nilpotent} 
\index{chi $\chi$-nilpotent} 
if there exists a 
$d: \Gamma_\chi^{d}(\pi) = \{1\}$, the smallest such $d$ being its 
\un{degree of nilpotency}: 
\index{degree of nilpotency} 
nil$_\chi\pi$.  
Example:  Take $G = \pi$ and let 
$\chi:\pi \ra \Aut \pi$ 
be the representation of $\pi$ by inner automorphisms $-$then $\pi$ is $\chi$-nilpotent iff $\pi$ is nilpotent.
 
[Note: \  From the definitions, for any 
$\chi$, $\Gamma^i(\pi) \subset \Gamma_\chi^i(\pi)$, thus if $\pi$ is $\chi$-nilpotent, then $\pi$ must be nilpotent.]

Let $\Pi$ and $\pi$ be groups, where $\Pi \subset \Aut \pi$.  
Suppose that 
$\pi = \pi_0 \supset$ $\pi_1 \supset \cdots$ $\supset \pi_d =$ $\{1\}$ 
is a finite filtration on $\pi$ by $\Pi -$stable normal subgroups such that $\Pi$ operates trivially on the 
$\pi_i/\pi_{i+1}$ $-$then there is a lemma in group theory that says that $\Pi$ must be nilpotent 
(Suzuki\footnote[2]{\textit{Group Theory}, vol II, Springer Verlag (1986), 19-20.}).  
So, given $\chi:G \ra \Aut \pi$, im $\chi$ is nilpotent provided that $\pi$ is $\chi$-nilpotent.\\

%%----------------------------------------------------------------------------------------------56
\label{8.29}
\begingroup%%----------------------------------->>
\fontsize{9pt}{11pt}\selectfont
\label{5.0amb}
\label{5.14b}
\textbf{\small FACT} \ 
Given a homomorphism 
$\chi:G \ra \Aut \pi$, consider the semidirect product 
$\pi \rtimes_\chi G$, 
i.e., the set of all ordered pairs $(\alpha, g) \in \pi \times G$ with the law of composition 
$(\alpha^\prime, g^\prime)(\alpha^{\prime\prime}g^{\prime\prime}) =$ 
$(\alpha^\prime(\chi(g^\prime)\alpha^{\prime\prime}),g^\prime g^{\prime\prime})$ $-$then 
$\pi \rtimes_\chi G$ is nilpotent iff $\pi$ is $\chi$-nilpotent and $G$ is nilpotent.\\
\endgroup%%------------------------------------<< 

\begingroup%%----------------------------------->>
\fontsize{9pt}{11pt}\selectfont
\textbf{\small EXAMPLE} \ 
Every finite $p$-group is nilpotent.  
Since the semidirect product of two finite $p$-groups is a finite $p$-group, it follows that if $G$ and $\pi$ are finite $p$-groups and if $G$ operates on $\pi$, then $G$ actually operates nilpotently on $\pi$.\\
\endgroup%%------------------------------------<< 

\begingroup%%----------------------------------->>
\fontsize{9pt}{11pt}\selectfont
\textbf{\small FACT} \ 
Suppose that $G$ operates on $\pi$ $-$then $G$ operates nilpotently on $\pi$ iff $\pi$ is nilpotent and $G$ operates nilpotently on $\pi/[\pi,\pi]$.\\
\endgroup%%------------------------------------<< 

\label{8.16}
\begingroup%%----------------------------------->>
\fontsize{9pt}{11pt}\selectfont
\textbf{\small EXAMPLE} \ 
Let $1 \ra G^\prime \ra G \ra G^{\prime \prime} \ra 1$ be a short exact sequence of groups.  
Obviously: $G$ nilpotent  $\implies$
$
\begin{cases}
\ G^\prime \\
\ G^{\prime \prime}
\end{cases}
$
nilpotent.  The converse is false (consider $A_3 \subset S_3$).  However, there is a characterization:  \mG is nilpotent iff
$
\begin{cases}
\ G^\prime \\
\ G^{\prime \prime}
\end{cases}
$
are nilpotent and the action of $G^{\prime \prime}$ on $G^\prime /  [G^\prime,G^\prime]$ is nilpotent.\\
\endgroup%%------------------------------------<< 

%%=================
Example:  Suppose that $\pi = M$ is a $G$-module.  Since $M$ is abelian, it is nilpotent but it needn't be $\chi$-nilpotent.  In fact, $\Gamma_\chi^i(M) = (I[G])^i\cdot M$, therefore $M$ is $\chi$-nilpotent iff $(I[G])^d\cdot M = 0$ for some $d$.  
When this is so, $M$ is referred to as a 
\un{nilpotent $G$-module}.
\index{nilpotent $G$-module}
\\
%%=================

\begingroup%%----------------------------------->>
\fontsize{9pt}{11pt}\selectfont
\textbf{\small EXAMPLE} \ 
Let $\pi$ be a nilpotent $G$-module.  Fix $n \geq 1$ $-$then $\forall \ q \geq 0$, $H_q(\pi,n)$ is a nilpotent $G$-module.
\vspi
[\mG operates nilpotently on the 
$\Gamma_\chi^i(\pi)$ and 
$\forall \ i$, there is a short exact sequence 
$0 \ra \Gamma_\chi^{i+1}(\pi) \ra$ 
$\Gamma_\chi^{i}(\pi) \ra$ 
$\Gamma_\chi^{i}(\pi)/\Gamma_\chi^{i+1}(\pi) \ra 0$ 
of $G$-modules, the action of $G$ on 
$\Gamma_\chi^{i}(\pi) / \Gamma_\chi^{i+1}(\pi) $ being trivial.  
The mapping fiber of the arrow 
$K(\Gamma_\chi^{i}(\pi),n) \ra$ 
$K(\Gamma_\chi^i(\pi)) / \Gamma_\chi^{i+1}(\pi),n)$ is a 
$K(\Gamma_\chi^{i+1}(\pi),n)$.  
Consider the associated fibration spectral sequence, noting that by induction, $G$ operates nilpotently on $E_{p,q}^2$.]\\
\endgroup%%------------------------------------<< 

\label{5.24}
\label{5.27}
\label{11.28}
\begingroup%%----------------------------------->>
\fontsize{9pt}{11pt}\selectfont
\textbf{\small FACT} \ 
Suppose that \mG is a finitely generated nilpotent group.  Let $M$ be a nilpotent $G$-module.\\
\indent \indent (1) \quadx If $M$ is finitely generated, then $\forall \ q \geq 0$, $H_q(G;M)$ is finitely generated.\\
\indent \indent (2) \quadx If $M$ is not finitely generated, then $H_0(G;M)$ is not finitely generated.\\
\endgroup%%------------------------------------<< 

%%------------
\label{5.22}
A nonempty path connected topological space $X$ is said to be 
\un{nilpotent}
\index{nilpotent (space)} 
if $\pi_1(X)$ is nilpotent and if $\forall \ n > 1$, $\pi_1(X)$ operates nilpotently on $\pi_n(X)$.  
Examples:  
(1)  Every abelian topological space is nilpotent; 
(2) Every path connected topological space whose homotopy groups are finite $p$-groups is nilpotent (cf. supra);  
(3)  Take for $X$ the Klein bottle $-$then $\pi_1(X)$ is not a nilpotent group;  
(4)  Take for $X$ the real projective plane $-$then 
$\pi_1(X) \approx$ 
$\Z/2\Z$, $\pi_2(X) \approx$ 
$\Z$ and the action of $\pi_1(X)$ on $\pi_2(X)$ is the inversion $n \ra -n$,
%%----------------------------------------------------------------------------------------------57
thus $\pi_1(X)$ does not operate nilpotently on $\pi_2(X)$;  
(5)  Take for $X$ the torus $\bS^1 \times \bS^1$ $-$then $X$ is nilpotent but its $1-$skeleton 
$X^{(1)} = \bS^1 \vee \bS^1$ is not nilpotent.\\

\begingroup%%----------------------------------->>
\fontsize{9pt}{11pt}\selectfont
\textbf{\small EXAMPLE} \ 
Let \mG be a topological group with base point $e$ and denote by $G_0$ the path component of $e$ $-$then 
$\pi_0(G) = G/G_0$ can be identified with $\pi_1(B_G^\infty)$ and $\pi_n(G) = \pi_n(G_0)$ can be identified with 
$\pi_{n+1}(B_G^\infty)$ 
(cf. p. \pageref{5.14a}).  
These identifications are compatible in that the homomorphism 
$\chi_n:\pi_0(G) \ra \Aut\pi_n(G_0)$ arising from the operation of $G$ on itself by inner automorphisms corresponds to the action of $\pi_1(B_G^\infty)$ on $\pi_{n+1}(B_G^\infty)$.  
Accordingly, $B_G^\infty$ is a nilpotent topological space iff $\pi_0(G)$ is a nilpotent group and $\forall \ n \geq 1$, 
$\pi_n(G_0)$ is $\chi_n$-nilpotent or still $\forall \ n \geq 1$, 
the semidirect product $\pi_n(G_0) \rtimes_{\chi_n} \pi_0(G)$ is nilpotent 
(cf. p. \pageref{5.14b}).  
The forgetful function 
$[\bS^n,s_n;G_0,e] \ra [\bS^n,G_0]$ is bijective, hence 
$[\bS^n,G_0] \approx \pi_n(G_0)$.  
In addition, 
$[\bS^n,G_0]$ is isomorphic to $\pi_n(G_0) \rtimes_{\chi_n} \pi_0(G)$.  
To see this, let 
$f:\bS^n \ra G$ be a continuous function.  Choose 
$g_f \in G: f(\bS^n) \subset G_0 g_f$, put $f_0 = f \cdot g_f^{-1}$ and consider the assignment 
$[f] \ra ([f_0],g_f G_0)$.  It therefore follows that $B_G^\infty$ is a nilpotent topological space iff 
$\forall \ n \geq 1$, $[\bS^n,G]$ is a nilpotent group.  
Example:  $B_{\bO(2n+1)}^\infty$ is nilpotent but $B_{\bO(2n)}^\infty$ is not nilpotent.
\vspi
[Note: \   Here is another illustration.  The higher homotopy groups of a connected nilpotent Lie group are trivial.  So, if $G$ is an arbitrary nilpotent Lie group, then $B_G^\infty$ is a nilpotent topological space.]\\ 
\endgroup%%------------------------------------<< 

\begingroup%%----------------------------------->>
\fontsize{9pt}{11pt}\selectfont
\textbf{\small FACT} \ 
Let $G$ be a topological group.  
Assume: $\forall \ n \geq 1$, $[\bS^n,G]$ is a nilpotent group $-$then for any finite CW complex $K$, $[K,G]$ is a nilpotent group.
\vspi
[Take $K$ connected and argue by induction on the number of cells.]\\
\endgroup%%------------------------------------<< 

\begingroup%%----------------------------------->>
\fontsize{9pt}{11pt}\selectfont
\textbf{\small EXAMPLE} \ 
Let $X$ be a nilpotent CW space $-$then Mislin\footnote[2]{\textit{Ann. of Math.} \textbf{103} (1976), 547-556.}
 has shown that $X$ is dominated in homotopy by a finite CW complex iff 
 the $H_q(X)$ are finitely generated $\forall \ q$ and there exists $q_0: \forall \ q > q_0$, $H_q(X) = 0$.  
 Moreover, under these conditions, Wall's obstruction to finiteness is zero provided that $\pi_1(X)$ is infinte but this can fail if $\pi_1(X)$ is finite (Mislin\footnote[3]{\textit{Topology} \textbf{14} (1975), 311-317.}).\\
\endgroup%%------------------------------------<< 

\index{Theorem: Dror's Whitehead}
\textbf{\small DROR'S WHITEHEAD THEOREM} \ 
Suppose that $X$ and $Y$ are nilpotent topological spaces.  
Let $f:X \ra Y$ be a homology equivalence $-$then $f$ is a weak homotopy equivalence.

[To prove that $f$ is a weak homotopy equivalence amounts to proving that for every $n$, 
the pair$(M_f,i(X))$ is $n$-connected, where, a priori, $H_*(M_f,i(X)) = 0$.  
Consider the commutative diagram
\begin{tikzcd}[]
X \ar{d} \ar{r}{f} &{Y} \ar{d}\\
X[1] \ar{r}[swap]{f_{1,1}} &{Y[1]}
\end{tikzcd}
.  
Since the vertical arrows are 2-equivalences, $f_{1,1}$
%%----------------------------------------------------------------------------------------------58
 induces a bijection $H_1(X[1]) \ra H_1(Y[1])$ and a surjection $H_2(X[1]) \ra H_2(Y[1])$.
 But 
$
\begin{cases}
\ X[1] \\[-.1cm]
\ Y[1]
\end{cases}
$
has the homotopy type 
$
\begin{cases}
\ (\pi_1(X),1) \\[-.1cm]
\ (\pi_1(Y),1)
\end{cases}
$
and 
$
\begin{cases}
\ \pi_1(X) \\[-.1cm]
\ \pi_1(Y)
\end{cases}
$
are nilpotent groups, thus $f_*:\pi_1(X) \ra \pi_1(Y)$ is an isomorphism 
(cf. p. \pageref{5.15}) 
and so $(M_f,i(X))$ is 1-connected.  
Noting that here $\pi_2(M_f,i(X))$ is abelian, fix $n > 1$ and assume inductively that 
$\pi_q(M_f,i(X)) = 0 $ for $q < n$ $-$then, from the relative Hurewicz theorem, 
$\pi_n(M_f,i(X))_{\pi_1(X)} = 0$, 
i.e., 
$\pi_n(M_f,i(X)) = I[\pi_1(X)] \cdot \pi_n(M_f,i(X))$.  
On the other hand, there is an exact sequence 
$\pi_n(M_f) \ra  \pi_n(M_f,i(X)) \ra \pi_{n-1}(i(X))$ of $\pi_1(X)$-modules. 
 Because the flanking terms are, by hypothesis, nilpotent $\pi_1(X)$-modules, the same must be true of $\pi_n(M_f,i(X))$.  
 Conclusion: $\pi_n(M_f,i(X)) = 0$.]\\

\begin{proposition} \ %15
Let $f:X \ra Y$ be a Hurewicz fibration, where $X$ and $Y$ are path connected.  Assume: $X$ is nilpotent $-$then 
$\forall \ y_0 \in Y$, the path components of $X_{y_0}$ are nilpotent.
\end{proposition}

[Fix $x_0 \in X_{y_0}$ and take $X_{y_0}$ path connected.  The homomorphisms in the homotopy sequence
\[
\cdots \ra \pi_{n+1}(Y,y_0)  \ra \pi_n(X_{y_0},x_0) \ra \pi_n(X,x_0) \ra \pi_{n}(Y,y_0) \ra \cdots
\]
of $f$ are $\pi_1(X,x_0)$-homomorphisms (cf. p. \pageref{5.16}).  
Of course $\pi_1(X,x_0)$ operates on $\pi_n(Y,y_0)$ through $f_*$ and if 
$i:X_{y_0} \ra X$ is the inclusion, then 
$\alpha \cdot \xi = (i_* \alpha) \cdot \xi$ 
$(\alpha \in \pi_1(X_{y_0},x_0)$, $\xi \in \pi_n(X_{y_0},x_0))$.  
Since the base points will play no further role, drop them from the notation.\\
%^
\indent \indent $(n = 1)$ \quadx 
To see that $\pi_1(X_{y_0})$ is nilpotent, consider the short exact sequence associated with the exact sequence
\begin{tikzcd}[ sep=small]
{\pi_2(Y)} \ar{r}{\partial} &{\pi_1(X_{y_0})} \ar{r}{i_*} &{\pi_1(X),}
\end{tikzcd}
noting that $\im \partial$ is contained in the center of $\pi_1(X_{y_0})$.\\
\indent \indent $(n > 1)$ \quadx 
There is an exact sequence
\begin{tikzcd}[ sep=small]
{\pi_{n+1}(Y)} \ar{r}{\partial} &{\pi_n(X_{y_0})} \ar{r}{i_*} &{\pi_n(X)}
\end{tikzcd}
and by assumption, $\exists \ d$ : $(I[\pi_1(X)])^d \cdot \pi_n(X) = 0$.
Claim:  $(I[\pi_1(X)])^{d+1} \cdot \pi_n(X_{y_0}) = 0$.
For let 
$\alpha \in (I[\pi_1(X_{y_0})])^d$, $\xi \in \pi_n(X_{y_0})$: 
$i_*(\alpha \cdot \xi) =$ $i_* \alpha \cdot i_* \xi =$ $0$ $\implies$ 
$\alpha \cdot \xi =$ $\partial \eta$ $(\eta \in \pi_{n+1}(Y))$. 
And: 
$\forall \ \beta \in \pi_1(X_{y_0})$, 
$(i_* \beta - 1) \cdot \eta =$ $(f_* i_* \beta - 1) \cdot \eta = 0$, so 
$0 = \partial((i_* \beta - 1) \cdot \eta )$ 
$= (i_* \beta - 1) \cdot \partial \eta$ $=$ 
$((\beta - 1) \alpha) \cdot \xi$.  Hence the claim.]\\

\label{9.25}
Application:  Let $X$ and $Y$ be pointed path connected spaces.  
Assume: $X$ is nilpotent $-$then for every pointed continuous function $f:X \ra Y$, 
the path componenets of the mapping fiber $E_f$ of $f$ are nilpotent.\\

\label{9.26}
\begingroup%%----------------------------------->>
\fontsize{9pt}{11pt}\selectfont
\textbf{\small EXAMPLE} \ 
Let $(K,k_0)$ be a pointed connected CW complex.  
Assume: $K$ is finite $-$then for any pointed path connected space $(X,x_0)$ the path components of $C(K,k_0;X,x_0)$ are nilpotent.  In 
%%----------------------------------------------------------------------------------------------59
particular, the fundamental group of the path component of the constant map $K \ra x_0$ is nilpotent, 
thus $[K,k_0;\Omega X,j(x_0)]$ is a nilpotent group.  
Observe that the base points play a role here: $[\bS^1,\Omega \textbf{P}^2(\R)]$ is a group but it is not nilpotent.\\
\endgroup%%------------------------------------<< 

\label{9.18}
\begingroup%%----------------------------------->>
\fontsize{9pt}{11pt}\selectfont
\textbf{\small FACT} \ 
Let $f:X \ra B$ be a Hurewicz fibration.  Given $\Phi^\prime \in C(B^\prime,B)$ define $X^\prime$ by the pullback square
\begin{tikzcd}%[ sep=small]
{X^\prime} \ar{d} \ar{r} &X \ar{d}{f}\\
{B^\prime} \ar{r}[swap]{\Phi^\prime} &B
\end{tikzcd}
.\ 
Assume:
$
\begin{cases}
\ X \\
\ B
\end{cases}
$
$\&$ $B^\prime$ are nilpotent $-$then the path components of $X^\prime$ are nilpotent.
\vspi
[Work with the Mayer-Vietoris sequence (cf. p. \pageref{5.17}).]\\
\endgroup%%------------------------------------<< 

\label{9.30}
\begingroup%%----------------------------------->>
\fontsize{9pt}{11pt}\selectfont
\textbf{\small EXAMPLE} \ 
The preceding result implies that nilpotency behaves well with respect to pullbacks but the situation for pushouts is not as satisfactory since nilpotency is not ordinarily inherited (consider $\bS^1 \vee \bS^2$).  
For example, suppose that $f:X \ra Y$ is a continuous function, where $X$ and $Y$ are nonempty path connected CW spaces.  Assume: $Y$ is nilpotent $-$then Rao\footnote[2]{\textit{Proc. Amer. Math Soc.} \textbf{87} (1983), 335-341.} 
has shown that the mapping cone $C_f$ of $f$ is nilpotent iff one of the following conditions is satisfied: 
(i) \ $f_*:\pi_1(X) \ra \pi_1(Y)$ is surjective;
(ii) \ $\forall \ q > 0$, $H_q(X) = 0$;
(iii) \  $\exists$ a prime $p$ such that $\pi_1(C_f)$ is a finite $p$-group and $\forall \ q > 0$, $H_q(X)$ is a $p$-group of finite exponent.
Example:  If $f:X \ra Y$ is a closed cofibration, then under (i), (ii), or (iii), $Y/f(X)$ is nilpotent 
(cf. p. \pageref{5.18}).  
Moreover, under (ii), the projection $Y \ra Y/f(X)$ is a homology equivalence 
(cf. p. \pageref{5.19}), hence by Dror's Whitehead theorem is a homotopy equivalence.\\
\endgroup%%------------------------------------<< 

Let
$
\begin{cases}
\ X \\[-.15cm]
\ Y
\end{cases}
$
be pointed connected CW spaces.  Suppose that $f:X \ra Y$ is a pointed continuous function. $-$then $f$ is said to admit a \un{principal refinement of order $n$}
\index{principal refinement of order $n$} 
if $f$ can be written as a composite
\begin{tikzcd}[ sep=small]
X \ar{r}{\Lambda} 
&{W_N} \ar{r}{q_N} 
&{W_{N-1}} \ar{r} 
&{\cdots} \ar{r} 
&{W_1} \ar{r}{q_1} 
&{W_0} %\quadx = &{Y}
\end{tikzcd}
$= Y,$ \ 
where $\Lambda$ is a pointed homotopy equivalence and each $q_i:W_i \ra W_{i-1}$ is a pointed Hurewicz fibration for which there is an abelian group $\pi_i$ and a pointed continuous function
$\Phi_{i-1}:W_{i-1} \ra K(\pi_i,n+1)$ such that the diagram
\begin{tikzcd}%[ sep=small]
{W_i} \ar{d}[swap]{q_i} \ar{r} &{\Theta K(\pi_i,n+1)} \ar{d}\\
{W_{i-1}} \ar{r}[swap]{\Phi_{i-1}}  &{K(\pi_i,n+1)}
\end{tikzcd}
is a pullback square.

[Note: \  $W_i$ is a pointed connected CW space homeomorphic to $E_{\Phi_{i-1}}$ (parameter reversal).]

Example:  If $X$ is a pointed abelian CW space, then $\forall \ n$, the arrow $f_n:X[n] \ra X[n-1]$ admits a principal refinement of order $n$:
\begin{tikzcd}%[ sep=small]
&{W[n]}\ar{d}\\
{X[n]} \ar{ru} \ar{r}  &{X[n-1]}
\end{tikzcd}
(cf. p. \pageref{5.20}), with $N = 1$.\\
\vspace{0.5cm}

%%----------------------------------------------------------------------------------------------60

\begingroup%%----------------------------------->>
\fontsize{9pt}{11pt}\selectfont
\index{central extensions}
\textbf{\small EXAMPLE \ (\un{Central Extensions})} \ 
Let $\pi$ and $G$ be groups, where $\pi$ is abelian $-$then the isomorphism classes of central extensions
$1 \ra \pi \ra \Pi \ra G \ra 1$ 
of $\pi$ by $G$ are in a one-to-one correspondence with the elements of 
$H^2(G,1;\pi)$ or still, with the elements of $[K(G,1),K(\pi,2)]$.  
Therefore $G$ is nilpotent iff the constant map $K(G,1) \ra *$ admits a principal refinement of order 1.
\vspi
[Any nilpotent $G$ generates a finite sequence of central extensions
$1 \ra$ 
$\Gamma^i(G)/\Gamma^{i+1}(G) \ra$ 
$G/\Gamma^{i+1}(G)  \ra$ 
$G/\Gamma^i(G) \ra 1$
.]\\
\endgroup%%------------------------------------<< 

Let $X$ be a pointed connected CW space $-$then, in view of the preceding example, the arrow $f_1:X[1] \ra X[0]$ admits a principal refinement of order 1 iff $\pi_1(X)$ is nilpotent.\\

\begin{proposition} \ %16
Let $X$ be a pointed connected CW space.  Fix $n > 1$ $-$then the arrow $f_n:X[n] \ra X[n-1]$ admits a principal refinement of order $n$ iff $\pi_1(X)$ operates nilpotently on $\pi_n(X)$.
\end{proposition}

[Necessity: \  Suppose that $f_n$ factors as a composite
\begin{tikzcd}[ sep=small]
{X[n]} \ar{r}{\Lambda} 
&{W_N} \ar{r}{q_N} 
&{W_{N-1}} \ar{r} &{}
\end{tikzcd} 
$\cdots$
\begin{tikzcd}[ sep=small]
\ar{r} 
&{W_1} 
\ar{r}{q_1}
&{W_0} %\ar{r}{=} &{X[n-1]}
\end{tikzcd}
$ = X[n-1]$
, where $\Lambda$ and the $q_i$ are as in the definition.  
Obviously $\pi_1(X) \approx \pi_1(W_i)$ for all $i$.  
Since $\pi_n(W_0) = \pi_n(X[n-1]) = 0$, $\pi_1(x)$ operates nilpotently on $\pi_n(W_0)$.  
Claim:  $\pi_1(x)$ operates nilpotently on $\pi_n(W_1)$.  
Thus let $\overline{W_0}$ be the mapping track of $\Phi_0$ and define $\overline{W_1}$ by the pullback square
\begin{tikzcd}%[sep=large]
{\overline{W_1}} \ar{d} \ar{r}  &{\Theta K(\pi_1,n+1)} \ar{d}\\
{\overline{W_0}} \ar{r} &{K(\pi_1,n+1)}
\end{tikzcd}
 $-$then there is a pointed homotopy equivalence $W_1 \ra \overline{W_1}$ and, from the proof of the ``$n > 1$'' part of Proposition 15,  
 $\pi_1(x)$ operates nilpotently on $\pi_n(\overline{W_1})$.  
Iterate to conclude that $\pi_1(X)$ operates nilpotently on $\pi_n(W_N) \approx \pi_n(X)$.

Sufficiency: \ 
One can copy the argument employed in the abelian case to construct the Postnikov invariant 
(cf. p. \pageref{5.21}).  
At the first stage, the only difference is that after replacing $n$ by $n-1$, the coefficient group for cohomology is not $\pi_n(X)$ but $\pi_n(X)_{\pi_1(X)} = H_0(\pi_1(X);\pi_n(X))$.  
Because the initial lifting 
\begin{tikzcd}%[ sep=small]
&{W_1}\ar{d}{q_1}\\
{X[n]} \ar{ru}{\Lambda_1} \ar{r}[swap]{f_n}  &{X[n-1]}
\end{tikzcd}
of $f_n$ is a pointed homotopy equivalence iff $I[\pi_1(X)] \cdot \pi_n(X) = 0$, 
it is in general necessary to repeat the procedure, which will then terminate after finitely many steps.]\\

\label{11.9}
Application:  Let $X$ be a pointed connected CW space $-$then $X$ is nilpotent iff $\forall \ n$, the arrow $f_n:X[n] \ra X[n-1]$ admits a principal refinement of order $n$.

[Note: \   If $X$ is nilpotent and if $\chi_n:\pi_1(X) \ra \Aut\pi_n(X)$ is the homomorphism corresponding to the action of 
$\pi_1(X)$ on $\pi_n(X)$, then a choice for the abelian groups figuring in the principal refinement of the arrow 
$X[n] \ra X[n-1]$ are the $\Gamma_{\chi_n}^i(\pi_n(X))/\Gamma_{\chi_n}^{i+1}(\pi_n(X))$.]\\
%%----------------------------------------------------------------------------------------------61

\begingroup%%----------------------------------->>
\fontsize{9pt}{11pt}\selectfont
\textbf{\small EXAMPLE} \ 
Let $K$ be a finite CW complex $-$then for any pointed nilpotent CW space $X$, 
the path components of $C(K,X)$ are nilpotent.
\vspi
[Bearing in mind $\S 4$, Proposition 5, use Proposition 15 and induction to show that $\forall \ n$, 
the path components of $C(K,X[n])$ are nilpotent.]\\
\endgroup%%------------------------------------<< 

\label{9.28}
\begingroup%%----------------------------------->>
\fontsize{9pt}{11pt}\selectfont
\textbf{\small EXAMPLE} \ 
Let $(K,k_0)$ be a pointed CW complex.  
Assume: $K$ is finite $-$then for any pointed nilpotent CW space $(X,x_0)$, the path components of $C(K,k_0;X,x_0)$ are nilpotent.  
Indeed, $C(K,k_0;X,x_0) = C(K_0,k_0;X,x_0) \times C(K_1,X) \times \cdots \times C(K_n,X)$, where $K_0$, $K_1$, $\ldots$, $K_n$ are the path components of $K$ and $k_0 \in K_0$.\\
\endgroup%%------------------------------------<< 

\index{Theorem Nilpotent Obstruction}
\begingroup%%----------------------------------->>
\fontsize{9pt}{11pt}\selectfont
\textbf{\small NILPOTENT OBSTRUCTION THEOREM} \ 
Let $(X,A)$ be a relative CW complex; let $Y$ be a pointed nilpotent CW space.  
Suppose that $\forall \ n > 0$ $\&$ $\forall \ i \geq 0$, 
$H^{n+1}(X,A;\Gamma_{\chi_n}^i(\pi_n(Y))/\Gamma_{\chi_n}^{i+1}(\pi_n(Y))) = 0$ 
$-$then every $f \in C(A,Y)$ admits an extension $F \in C(X,Y)$, any two such being homotopic rel $A$ provided that $\forall \ n > 0$ $\&$ $\forall \ i \geq 0$, 
$H^{n}(X,A;\Gamma_{\chi_n}^i(\pi_n(Y))/\Gamma_{\chi_n}^{i+1}(\pi_n(Y))) = 0$.\\
\endgroup%%------------------------------------<<

\begin{proposition} \ %17
Let $X$ be a pointed connected CW space, $\widetilde{X}$ its universal covering space.  Assume $\pi_1(X)$ is nilpotent $-$then $X$ is nilpotent iff $\forall \ n \geq 1$, $\pi_1(X)$ operates nilpotently on $H_n(\widetilde{X})$.
\end{proposition}

[$\widetilde{X}$ exists and is a pointed connected CW space (cf. Proposition 5).

Necessity:  Consider the Postnikov tower of $\widetilde{X}$, so $\widetilde{p}_n:P_n\widetilde{X} \ra P_{n-1}\widetilde{X}$.  
Suppose inductively that $\pi_1(X)$ operates nilpotently on the homology of $P_{n-1}\widetilde{X}$.
Since $X$ is nilpotent, the $H_q(\pi_n(X),n)$  are nilpotent $\pi_1(X)$-modules (cf. p. \pageref{5.22}), i.e., $\pi_1(X)$ operates nilpotently on the homology of the mapping fiber of $\widetilde{p}_n$.  
Therefore, by the universal coefficient theorem, the $E_{p,q}^2 \approx H_p(P_{n-1}\widetilde{X};H_q(\pi_n(X),n))$ in the fibration spectral sequence of $\widetilde{p}_n$ are nilpotent $\pi_1(X)$-modules, thus the same is true of the $H_i(P_n\widetilde{X})$.  
But the arrow $\widetilde{X} \ra P_n\widetilde{X}$ induces an isomorphism of 
$\pi_1(X)$-modules $H_i(\widetilde{X}) \ra H_i(P_n\widetilde{X})$ for $i \leq n$.

Sufficiency:  Introduce the Whitehead tower of $\widetilde{X}$ and argue as above.]\\

\begin{proposition} \ %18
Let $X$ be a pointed connected CW space.  
Assume:  $X$ is nilpotent $-$then the $\pi_q(X)$ are finitely generated $\forall \ q$ iff the $H_q(X)$ are finitely generated $\forall \ q$. 
\end{proposition}

[Suppose that the $\pi_q(X)$ are finitely generated $\forall \ q$ $-$then, $\widetilde{X}$ being simply connected, hence abelian, the $H_q(\widetilde{X})$ are finitely generated $\forall \ q$ 
(cf. p. \pageref{5.23}).  
On the other hand, according to Proposition 17, $\pi_1(X)$ operates nilpotently on the $H_q(\widetilde{X})$.  
Consequently, the $H_p(\pi_1(X);H_q(\widetilde{X}))$ are finitely generated 
(cf. p. \pageref{5.24}).  However, these terms are precisely the $E_{p,q}^2$ in the spectral sequence of the covering projection $\widetilde{X} \ra X$ (see below), so $\forall \ i$, $H_i(X)$ is finitely generated.

%%----------------------------------------------------------------------------------------------62
Suppose that the $H_q(X)$ are finitely generated $\forall \ q$ 
$-$then, since $\pi_1(X)/[\pi_1(X),\pi_1(X)]$ $\approx$ $H_1(X)$, 
the nilpotent group $\pi_1(X)$ is finitely generated 
(cf. p. \pageref{5.25}).  
As for the $\pi_q(X)$ $(q > 1)$, their finite generation will follow if it can be shown that the $H_q(\widetilde{X})$ are finitely generated 
(cf. p. \pageref{5.26}).  
Proceeding by contradiction, fix an $i_0$ such that $H_{i_0}(\widetilde{X})$ is not finitely generated and take $i_0$ minimal.  
The 
$E_{p,q}^2 \approx$ 
$H_p(\pi_1(X);H_q(\widetilde{X}))$ are finitely generated if $q < i_0$ but 
$E_{0,i_0}^2 \approx$ 
$H_0(\pi_1(X);H_{i_0}(\widetilde{X}))$ is not finitely generated 
(cf. p. \pageref{5.27}), thus $E_{0,i_0}^\infty$ is not finitely generated.  
Therefore $H_{i_0}(X)$ contains a subgroup which is not finitely generated.]

[Note: \  A finitely generated nilpotent group is finitely presented and its integral group ring is (left and right) noetherian.  This said, it then follows that under the equivalent conditions of the proposition, $X$ necessarily has the pointed homotopy type of a pointed CW complex with a finite $n$-skeleton $\forall \ n$ 
(Wall\footnote[2]{\textit{Ann. of Math.} \textbf{81} (1965), 56-69.}).]\\

\label{7.2}
The spectral sequence 
$E_{p,q}^2 \approx H_p(\pi_1(X);H_q(\widetilde{X}))$ $\Rightarrow$ $H_{p+q}(X)$ 
of the covering projection $\widetilde{X} \ra X$ is an instance of a fibration spectral sequence.  In fact, consider the inclusion 
$i:X \ra X[1] =$ $K(\pi_1(X),1)$ and pass to its mapping track 
$W_i \ra$ $K(\pi_1(X),1)$ $-$then $E_i$ has the same pointed homotopy type as $\widetilde{X}$.  
Moreover, 
$H_p(\pi_1(X);H_q(\widetilde{X})) \approx$ 
$H_p(K(\pi_1(X),1);\sH_q(\widetilde{X}))$, 
where 
$\sH_q(\widetilde{X})$ is the locally constant coefficient system on 
$K($ $\pi_1(X),1)$ 
determined by $H_q(\widetilde{X})$ (cf. p. \pageref{5.27a}).\\

\label{5.34a}
\label{5.40}
\label{9.15} %dmc mnft
\begingroup%%----------------------------------->>
\fontsize{9pt}{11pt}\selectfont
\textbf{\small FACT} \ 
Suppose that 
$
\begin{cases}
\ X \\
\ Y
\end{cases} 
$
are pointed connected CW spaces.  Let $f:X \ra Y$ be a pointed Hurewicz fibration with $\pi_0(X_{y_0}) = *$ $-$then $\pi_1(X)$ operates nilpotently on the $\pi_q(X_{y_0})$ $\forall \ q$ iff $X_{y_0}$ is nilpotent and $\pi_1(Y)$  operates nilpotently on the $H_q(X_{y_0})$ $\forall \ q$.
\\
\endgroup%%------------------------------------<< 

\begingroup%%----------------------------------->>
\fontsize{9pt}{11pt}\selectfont
\textbf{\small EXAMPLE} \ 
Suppose that 
$
\begin{cases}
\ X \\
\ Y
\end{cases}
$
are pointed connected CW spaces.  
Let $f:X \ra Y$ be a pointed Hurewicz fibration with $\pi_0(X_{y_0}) = *$ $-$then any two of the following conditions imply the third and the third implies that $X_{y_0}$ is nilpotent:
(i) \ $X$ is nilpotent;
(ii) \ $Y$ is nilpotent;
(iii) \ $\pi_1(X)$ operates nilpotently on the $\pi_q(X_{y_0})$ $\forall \ q$.  Assume now that $\pi_1(Y)$  operates nilpotently on the $H_q(X_{y_0})$ $\forall \ q$.  
Claim: $X$ is nilpotent iff both $Y$ and $X_{y_0}$ are nilpotent.  
For $X$ nilpotent $\implies$  $X_{y_0}$ is nilpotent (cf. Proposition 15) $\implies$ $\pi_1(X)$ operates nilpotently on the 
$\pi_q(X_{y_0})$ $\forall \ q$ $\implies$ $Y$ nilpotent, and conversely.\\
\endgroup%%------------------------------------<< 

\index{Theorem Hilton-Roitberg Comparison}
\begingroup%%----------------------------------->>
\fontsize{9pt}{11pt}\selectfont
\textbf{\small HILTON-ROITBERG\footnote[3]{\textit{Quart. J. Math.} \textbf{27} (1976), 433-444; 
see also Sch\"on, \textit{Quart. J. Math.} \textbf{32} (1981), 235-237.}
 COMPARISON THEOREM} \ 
Suppose that 
$
\begin{cases}
\ X \\[-.1cm]
\ Y
\end{cases}
$
\hspace{-.25cm}
$\&$
$
\begin{cases}
\ X^\prime \\[-.1cm]
\ Y^\prime
\end{cases}
$
are pointed connected CW spaces.  Let $f:X \ra Y$ and $f^\prime:X^\prime \ra Y^\prime$ be pointed Hurewicz fibrations such that $E_f$ 
%%----------------------------------------------------------------------------------------------63
and $E_f^\prime$ are path connected and 
$
\begin{cases}
\ \pi_1(Y) \\
\ \pi_1(Y^\prime)
\end{cases}
$
operates nilpotently on the 
$
\begin{cases}
\ H_q(E_f)\\
\ H_q(E_{f^\prime})
\end{cases}
$
\hspace{-.25cm}
$\forall \ q$.  
Suppose there is a commutative diagram
\begin{tikzcd}%[sep=large]
X \ar{d} \ar{r}{f}  &Y\ar{d}\\
X^\prime \ar{r}[swap]{f^\prime} &Y^\prime
\end{tikzcd}
, where $\pi_1(Y) \approx \pi_1(Y^\prime)$ or $\pi_1(Y)$ $\&$ $\pi_1(Y^\prime)$ are nilpotent $-$then, assuming that all isomorphisms are induced, any two of the following conditions imply the third:
(1)\ $\forall \ p$, $H_p(Y) \approx H_n(Y^\prime)$;
(2) \ $\forall \ q$, $H_q(E_f) \approx H_q(E_{f^\prime})$;
(3) \ $\forall \ n$, $H_n(X) \approx H_n(X^\prime)$.
\\
\endgroup%%------------------------------------<< 

A nonempty path connected topological space $X$ is said to be 
\un{acyclic}
\index{acyclic (space)} 
provided that $\forall \ q > 0$, $H_q(X) = 0$.  
So: $X$ acyclic $\implies$ 
$\pi = [\pi,\pi]$ and $H_1(\pi,1) = 0 = H_2(\pi,1)$ 
(cf. p. \pageref{5.27b}), 
where $\pi = \pi_1(X)$.  
Example:  Every nilpotent acyclic space is homotopically trivial (quote Dror's Whitehead theorem).\\

\index{example (acyclic groups)}
\index{acyclic (group)}
\begingroup%%----------------------------------->>
\fontsize{9pt}{11pt}\selectfont
\textbf{\small EXAMPLE \ \un{(Acyclic Groups)}} \   
A group \mG is said to be \un{acyclic} if $\forall \ n > 0$, $H_n(G) = 0$ or, equivalently, if $K(G,1)$ is an acyclic space.  Nontrivial finite groups are never acyclic (Swan\footnote[2]{\textit{Proc. Amer. Math. Soc.} \textbf{11} (1960), 885-887.}).  
However, there are plenty of concretely defined infinite acyclic groups.  A list of examples has been compiled by 
Harpe-McDuff\footnote[3]{\textit{Comment. Math. Helv.} \textbf{58} (1983), 48-71; 
see also Berrick, In: \textit{Group Theory}, K. Cheng and Y. Leong (ed.), Walter deGruyter (1989), 253-266.}.  
They include: 
(1) The symmetric group on an infinite set; 
(2) The group of invertible linear transformations of an infinite dimensional vector space; 
(3) The group of invertible bounded linear transformations of an infinite dimensional Hilbert space; 
(4) The automorphism group of the measure algebra of the unit interval; 
(5) The group of compactly supported homeomorphisms of $\R^n$.\\
\endgroup%%------------------------------------<< 

\begingroup%%----------------------------------->>
\fontsize{9pt}{11pt}\selectfont
\textbf{\small FACT} \ 
Let $G$ be a group which is the colimit of subgroups $G_n$ $(n \in \N)$ with the property that 
$\forall \ n$, there exists a nontrivial $g_n \in G_{n+1}$ and a homomorphism 
$\phi_n:G_n \ra \text{Cen}_{G_{n+1}}(G_n)$ 
such that $\forall g \in G_n$, $g = [g_n,\phi_n(g)]$ $-$then $G$ is acyclic.
\vspi
[It suffices to work with coefficients in an arbitrary field \bk.  
Since $H_*(G;\bk) \approx \text{colim}H_*(G_n;\bk)$, 
one need only show that $\forall \ n \geq 1$ $\&$ $N \geq 1$, the morphism 
$H_q(G_n;\bk) \ra H_q(G_{n+N};\bk)$ induced by the inclusion 
$G_n \ra G_{n+N}$ is trivial when $1 \leq q < 2^N$.  
For this, fix $n$ and use induction on $N$.  
Recall that conjugation induces the identity on homology and apply the K\"unneth formula.]
\vspi
[Note: \  It is clear that $\phi_n$ is injective ($\implies g_n \in G_{n+1} - G_n$).  
Observe too that it is not necessary to assume that $\phi_n(G_n)$ is contained in the centralizer of $G_n$ in $G_{n+1}$ as this is implied by the other condition.  
Proof:  $\forall \ g,h \in G_n$: 
$[g_n,\phi_n(gh)] =$ 
$[g_n,\phi_n(g)] \cdot [\phi_n(g),[g_n,\phi_n(h)]] \cdot [g_n,\phi_n(h)]$ $\implies$ 
$gh = g[\phi_n(g),h]h$ $\implies$ 
$e = [\phi_n(g),h]$.]\\
\endgroup%%------------------------------------<< 

%%----------------------------------------------------------------------------------------------64
\label{5.42}
\begingroup%%----------------------------------->>
\fontsize{9pt}{11pt}\selectfont
\textbf{\small EXAMPLE} \ 
Let $H_c(\Q)$ be the set of bijections of $\Q$ that are the identity outside some finite interval.  
Given a group $G$, let $F_c(\Q,G)$ be the set of functions 
$\Q \ra G$ that send all elements outside some finite interval to the identity.  
Both $H_c(\Q)$ and $F_c(\Q,G)$ are groups and there is a homomorphism 
$\chi:H_c(\Q) \ra$ 
$\Aut F_c(\Q,G)$, 
viz. 
$\chi(\beta)\alpha(q) =$ $\alpha(\beta^{-1}(q))$.  The 
\un{cone}
\index{cone (of a group)} 
of $G$ is the associated semidirect product:
$\Gamma G = F_c(\Q,G) \rtimes_\chi H_c(\Q)$.  
The assignment
$
\begin{cases}
\ G \ra \Gamma G \\
\ g \ra \alpha_g
\end{cases}
$
:  $\alpha_g(q) = $
$
\begin{cases}
\ g \quadx (q = 0) \\
\ e \quadx (q \neq 0)
\end{cases}
$
is a monomorphism of groups and $\Gamma G$ is acyclic.
\vspi
[Let 
$\Gamma G_n = \{(\alpha,\beta): \spt \alpha \cup \spt \beta \subset [-n,n]\}$ 
and construct a homomorphism 
$\phi_n:\Gamma G_n \ra$ 
$\Cen_{\Gamma G_{n+1}}(\Gamma G_n)$ in terms of a bijection 
$\beta_n \in H_c(\Q)$: $\spt \beta_n \subset [-n-1,n+1]$ \ 
$\&$ \  $\forall \ k$ : $\beta_n^k[-n,n] \cap [-n,n] = \emptyset$.]\\
\endgroup%%------------------------------------<< 

\begingroup%%----------------------------------->>
\fontsize{9pt}{11pt}\selectfont
\textbf{\small FACT} \ 
Every group can be embedded in an acyclic simple group.
\vspi
[By the above, every group can be embedded in an acyclic group.  On the other hand, every group can be embedded in a simple group (Robinson\footnote[2]{\textit{Finiteness Conditions and Generalized Solvable Groups, vol. I}, Springer Verlag (1972), 144.}).  
So given $G$, there is a sequence 
$G \subset G_1 \subset G_2 \subset \cdots$, 
where $G_n$ is acyclic if $n$ is odd and simple if $n$ is even.  Consider $\ds\bigcup\limits_n G_n$.]\\
\endgroup%%------------------------------------<< 

\label{5.35}
Recall that a group $G$ is said to be 
\un{perfect}
\index{perfect (group)} 
if $G = [G,G]$.  
Examples: 
(1)  Every acyclic group is perfect; 
(2) Every nonabelian simple group is perfect.

[Note: \  The fundamental group of an acyclic space is perfect.]

The homomorphic image of a perfect group is perfect.  Therefore, if $G$ is perfect and $\pi$ is nilpotent, then $G$ operates nilpotently on $\pi$ iff $G$ operates trivially on $\pi$ 
(cf. p. \pageref{5.27c}).  
Proof:  A perfect nilpotent group is trivial.

Every group $G$ has a unique maximal perfect subgroup $G_\text{per}$, the 
\un{perfect radical}
\index{perfect radical (group)} 
of $G$.  The automorphisms of $G$ stabilize $G_\text{per}$, thus $G_\text{per}$ is normal.\\
\indent \indent ($\text{P}_1$) \quadx Let $f:G \ra K$ be a homomorphism of groups $-$then $f(G_\text{per}) \subset K_\text{per}$.\\
\indent \indent ($\text{P}_2$) \quadx Let $f:G \ra K$ be a homomorphism of groups, where $K_\text{per} = \{1\}$ $-$then $G_\text{per}\subset \ker f$.\\

\label{8.8}
\begingroup%%----------------------------------->>
\fontsize{9pt}{11pt}\selectfont
\textbf{\small FACT} \ 
A locally free group is acyclic iff it is perfect.
\vspi
[Note: \  A group is said to be 
\un{locally free}
\index{locally free} 
if its finitely generated subgroups are free.]\\
\endgroup%%------------------------------------<< 

\textbf{\small LEMMA} \ 
Let $f:G \ra K$ be an epimorphism of groups.  Put $N = \ker f$ $-$then $f(G_\text{per}) = K_\text{per}$ 
provided that $\exists \ n: N^{(n)} \subset G_\text{per}$.

[Note: \  $N^{(n)}$ is the $n^\text{th}$ derived group of $N: N^{(0)} = N$, $ N^{(i+1)} = [N^{(i)}, N^{(i)}]$.  
Obviously, 
$N^{(0)} \subset G_\text{per}$ if $N$ is perfect and $N^{(1)} \subset G_\text{per}$ if $N$ is central.]\\

%%----------------------------------------------------------------------------------------------65
\label{5.34f}
\label{5.38}
Application: Let $N$ be a perfect normal subgroup of $G$ $-$then the perfect radical of 
$G/N$ is the quotient $G_\text{per}/N$, hence the perfect radical of $G/N$ is trivial iff 
$N = G_\text{per}$.\\

\begingroup%%----------------------------------->>
\fontsize{9pt}{11pt}\selectfont
\textbf{\small EXAMPLE} \ 
Let $A$ be a ring with unit.  Agreeing to employ the usual notation of algebraic K-theory, denote by 
$\bGL(A)$ the infinite general linear group of $A$ and write $\bE(A)$ for the subgroup of $\bGL(A)$ consisting of the elementary matrices $-$then, according to the Whitehead lemma, 
$\bE(A) =$ 
$[\bE(A),\bE(A)] =$ 
$[\bGL(A),\bGL(A)]$, thus $\bE(A)$ is the perfect radical of $\bGL(A)$.  
Let now $\bST(A)$ be the Steinberg group of $A$: $\bST(A)$ is perfect and there is an epimorphism 
$\bST(A) \ra$ 
$\bE(A)$ of groups whose kernel is the center of $\bST(A)$.
\vspi
[Note: \ On occasion, it is necessary to consider rings which may not have a unit (pseudorings).  
Given a pseudoring $A$, let $\ov{A}$ be the set of all functions 
$X:\N \times \N \ra$ $A$ such that $\#\{(i,j):X_{ij} \neq 0\} < \omega$ $-$then 
$\ov{A}$ is again a pseudoring (matrix operations).  
The law of composition 
$X * Y = $ $X + Y + X \times Y$ equips $\ov{A}$ with the structure of a semigroup with unit.  
Definition: $\ov{\bGL}(A)$ is the group of units of $(\ov{A},*)$.  
Therefore, using obvious notation, 
$\ov{\bE}(A) =$ 
$[\ov{\bE}(A),\ov{\bE}(A)] =$ 
$[\ov{\bGL}(A),\ov{\bGL}(A)]$.  
Every bijection 
$\phi:\N \ra$ $\N \times \N$ defines an isomorphism of pseudorings: 
$\ov{\ov{A}} \approx$ $\ov{A}$, hence 
$\ov{\bGL}(\ov{A}) \approx$ 
$\ov{\bGL}(A)$.  
In the event that $A$ has a unit, the assignment
$
\begin{cases}
\ \ov{\bGL}(A) \ra \bGL(A) \\
\ X \ra X + I
\end{cases}
$
is an isomorphism of groups ($\implies$ 
$\ov{\bGL}(\ov{A}) \approx$ $\bGL(A)$).]\\
\endgroup%%------------------------------------<< 

\begingroup%%----------------------------------->>
\fontsize{9pt}{11pt}\selectfont
\index{universal central extensions}
\label{5.36}
\label{5.37}
\textbf{\small EXAMPLE \ \un{(Universal Central Extensions)}} \ 
Let $G$ be a group $-$then a central extension
$1 \ra N \ra U \ra G \ra 1$ is said to be \un{universal} if for any other central extension
$1 \ra \pi \ra \Pi \ra G \ra 1$ there is a unique homomophism
\begin{tikzcd}[sep=small]
U \ar{rr} \ar{rdd} &&\Pi \ar{ldd}\\
\\
&G
\end{tikzcd}
over \mG.  
A central extension 
$1 \ra N \ra U \ra G \ra 1$ is universal iff $H_1(U) = 0 = H_2(U)$.  
On the other hand, a universal central extension 
$1 \ra N \ra U \ra G \ra 1$ exists iff $G$ is pefect.  
To identify $N$ in terms of $G$, use a portion of the fundamental exact sequence:
$H_2(U) \ra$ 
$H_2(G) \ra$ 
$N/[U,N] \ra$ 
$H_1(U)$ or still,
$0 \ra$ \
$H_2(G) \ra$ 
$N/[U,N] \ra 0$ $\implies$ 
$H_2(G) \approx N$.  
Example;  Take 
$G = \bE(A)$ $-$then 
$H_1(\textbf{ST}(A)) =$ $0 = H_2(\textbf{ST}(A))$ and there is a universal extension
$1 \ra$ 
$H_2(\bE(A)) \ra$ 
$\textbf{ST}(A)) \ra$ 
$\bE(A) \ra 1$.\\
\endgroup%%------------------------------------<<

\begingroup%%----------------------------------->>
\fontsize{9pt}{11pt}\selectfont
\textbf{\small EXAMPLE} \ 
Let \textbf{ACYGR} be the full subcategory of \bGR whose objects are the acyclic groups $-$then 
Berrick\footnote[2]{\textit{J. Pure Appl. Algebra} \textbf{44} (1987), 35-43.}
 has defined a functor $\alpha:\bAB \ra \textbf{ACYGR}$ such that $\forall \ G$, the center of $\alpha G$ is naturally isomorphic to $G$.  
The quotient $\beta G = \alpha G / \text{Cen }G$ is a perfect group and the central extension 
$1 \ra$ 
$G \ra \alpha G \ra$ 
$\beta G \ra 1$ is universal, so $G \approx$ $H_2(\beta G)$. 
\vspi
[Note: \  By contrast, the cone construction defines a functor $\Gamma:\bGR \ra$ $\textbf{ACYGR}$.]\\
\endgroup%%------------------------------------<<

\begingroup%%----------------------------------->>
\fontsize{9pt}{11pt}\selectfont
\textbf{\small FACT} \ 
Let
$
\begin{cases}
\ G_1 \\
\ G_2
\end{cases}
$
be groups $-$then the perfect radical of $G_1 \times G_2$ is $(G_1)_{\text{per}} \times (G_2)_{\text{per}}$.\\
\endgroup%%------------------------------------<<

%%----------------------------------------------------------------------------------------------66

\begingroup%%----------------------------------->>
\fontsize{9pt}{11pt}\selectfont
\textbf{\small FACT} 
Let
$
\begin{cases}
\ G_1 \\
\ G_2
\end{cases}
$
be groups with trivial perfect radicals $-$then the perfect radical of their free product $G_1 * G_2$ is trivial.
\vspi
[A theorem of Kurosch says that any subgroup $G$ of $G_1 * G_2$ has the form 
%$F* (\underset{i}{\ast} G i)$, 
$F*(
\underset{i}{
\begingroup%%X
\fontsize{18pt}{18pt}\selectfont
\text{$*$} 
\endgroup%%X 
}
G i)
$
where $F$ is a free group and $\forall \ i$, $G_i$ is isomorphic to a subgroup of either 
$G_1$ or $G_2$.  Put 
$X =$ 
$K(F,1) \vee \ds\bigvee\limits_i K(G_i,1)$: $\pi_1(X) \approx$ $G$.  
If $G$ is perfect, then 
$0 =$ 
$H_1(X) \approx$ 
$H_1(F) \oplus \ds\bigoplus\limits_i H_1(G_i)$, 
and it follows that $F$ and the $G_i$ are perfect, hence trivial.]\\
\endgroup%%------------------------------------<<

Let
$
\begin{cases}
\ X \\[-.1cm]
\ Y
\end{cases}
$
be pointed connected CW spaces.  
Suppose that $f:X \ra Y$ is a pointed continuous function $-$then $f$ is said to be 
\un{acyclic} 
\index{acylic} 
if its mapping fiber $E_f$ is acyclic.  
For this, it is therefore necessary that $\pi_0(E_f) = *$.

[Note: \  Using the mapping cylinder $M_f$, write $f = r \circx i$ 
(cf. p. \pageref{5.27d})) 
$-$then $(M_f,i(x_0))$ is nondegenerate, thus $r:M_f \ra Y$ is a pointed homotopy equivalence
(cf. p. \pageref{5.27e})) 
which implies that the arrow 
$E_i \ra$ 
$E_{r \circx i} =$ 
$E_f$ is a pointed homotopy equivalence 
(cf. p. \pageref{5.27f})).   
Conclusion: $f:X \ra Y$ is acyclic iff $i:X \ra M_f$ is acyclic.]

Observation:  Suppose that $f:X \ra Y$ is acyclic $-$then 
$f_*:\pi_1(X) \ra \pi_1(Y)$ is surjective and its kernel is a perfect normal subgroup of $\pi_1(X)$.

[Inspect the exact sequence $\pi_2(Y) \ra \pi_1(E_f) \ra \pi_1(X) \ra \pi_1(Y) \ra\pi_0(E_f)$.]\\

\begin{proposition} \ %19
Let
$
\begin{cases}
\ X \\[-.1cm]
\ Y
\end{cases}
$
be pointed connected CW spaces, $f:X \ra Y$ a pointed continuous function $-$then $f$ is a pointed homotopy equivalence iff $f$ is acyclic and $f_*: \pi_1(X) \ra \pi_1(Y)$ is an isomorphism.
\end{proposition}

[The necessity is clear.  As for the sufficiency, the arrow $\pi_2(Y) \ra \pi_1(E_f)$ is surjective, hence $\pi_1(E_f)$ is both abelian and perfect.  But this means that $\pi_1(E_f)$ must be trivial, so, being a pointed connected CW space, $E_f$ is contractible.]\\

\label{9.68}
\begingroup%%----------------------------------->>
\fontsize{9pt}{11pt}\selectfont
Let $P$ be a set of primes.  
Fix an abelian group $G$ $-$then $G$ is said to be 
\un{$P$-primary}
\index{P-primary} 
if \ 
$\forall \ g \in G$, \ 
$\exists \ F \subset P$ $(\#(F) < \omega)$ $\&$ $n \in \N$: 
$\ds\bigl(\prod\limits_{p \in F} p\bigr)^n g = 0$ \ 
$\ds\bigl(\prod\limits_\emptyset = 1\bigr)$ 
and $G$ is said to be 
\un{uniquely $P$-divisible}
\index{uniquely $P$-divisible} 
if $\forall \ g \in G$, $\forall \ p \in P$, $\exists ! \ h \in G: ph = g$.
\vspi
[Note: \  If $P$ is empty, then the only 
$P$-primary abelian group is the trivial group and every abelian group is uniquely $P$-divisible.]\\
\endgroup%%------------------------------------<< 

\begingroup%%----------------------------------->>
\fontsize{9pt}{11pt}\selectfont
\textbf{\small LEMMA} \ 
Let $\sC$ be a class of abelian groups containing 0.  Assume:  $\sC$ is closed under the formation of direct sums and five term exact sequences, i.e., for any exact sequence $G_1 \ra G_2 \ra G_3 \ra G_4 \ra G_5$ of abelian groups
$
\begin{cases}
\ G_1, G_2 \\
\ G_4, G_5
\end{cases}
$
$\in \sC$ $\implies G_3 \in \sC$ $-$then there exists a set of primes $P$ such that $\sC$  is either the class of $P$-primary abelian groups or the class of uniquely $P$-divisible abelian groups.
\vspi
%%----------------------------------------------------------------------------------------------67
[The hypotheses imply that $\sC$  is colimit closed.  Given a set $P$ of primes, 
it follows that if $\Z/p\Z \in \sC$  $\forall \ p \in P$, 
then every $P$-primary abelian group is in 
$\sC$  or if $\Q \in \sC$ and 
$\Z/p\Z \in \sC$ $\forall \ p \notin P$, 
then every uniquely $P$-divisible abelian group is in $\sC$.  On the other hand, if some $G \in \sC$ is not uniquely $P$-divisible, then $\Z/p\Z \in \sC$ (consider
\begin{tikzcd}[ sep=small]
G \ar{r}{p} &{G)}
\end{tikzcd}
and if some $G \in \sC$ is not torsion, then $\Q \in \sC$ 
(consider $\Q \otimes G =$ 
$\colimx (\cdots \ra G \overset{n}{\ra} G \ra \cdots)$).
%\begin{tikzcd}[ sep=small]
%{(\cdots} \ar{r} &G \ar{r}{n} &G \ar{r} &{\cdots )).}
%\end{tikzcd}
To summarize: 
(1) If $\Q \notin \sC$ and  $\Z/p\Z \in \sC$ exactly for  $p \in P$, then $\sC$ consists of the $P$-primary abelian groups;
(2)  If $\Q \in \sC$ and $\Z/p\Z \in \sC$ exactly for  $p \notin P$, then $\sC$ consists of the uniquely $P$-divisible abelian groups.]\\
\endgroup%%------------------------------------<< 

\begingroup%%----------------------------------->>
\fontsize{9pt}{11pt}\selectfont
Application:
Fix abelian groups
$
\begin{cases}
\ A \\
\ B
\end{cases}
$
$-$then $A \otimes B = 0 = \Tor(A,B)$ iff there exists a set $P$ of primes such that one of the groups is $P$-primary and the other is uniquely $P$-divisible.
\vspi
[Supposing that $A \otimes B = 0 = \Tor(A,B)$, the class of abelian groups $G$ for which 
$G \otimes B =$ 
$0 =$ 
$\Tor(G,B)$ satisfies the assumptions of the lemma.]\\
\endgroup%%------------------------------------<< 

\begingroup%%----------------------------------->>
\fontsize{9pt}{11pt}\selectfont
\textbf{\small EXAMPLE} \ 
Given a 2-sink
$X \overset{p}{\lra} B \overset{q}{\lla} Y$, 
where
$
\begin{cases}
\ X \\
\ Y
\end{cases}
$
$\&$ $B$ are pointed connected CW spaces, form $X \square_BY$ 
(cf. p. \pageref{5.28}).  
Let $r:X \square_B Y \ra B$ be the projection $-$then the following conditions are equivalent:
(i) \ $r$ is a pointed homotopy equivalence;
(ii) \ $E_r$ is acyclic;
(iii) \ $\exists \ P$ such that one of 
$
\begin{cases}
\ \widetilde{H}_*(E_p) = \bigoplus\limits_i \widetilde{H}_i(E_p)\\
\ \widetilde{H}_*(E_q) = \bigoplus\limits_j \widetilde{H}_j(E_q)
\end{cases}
$
is $P$-primary and the other is uniquely $P$-divisible.  
To see this, recall that $E_r \approx E_p * E_q$ 
(cf. p. \pageref{5.29}) 
and, on general grounds, 
$\widetilde{H}_{k+1}(E_p * E_q) \approx$ 
$\ds\bigoplus\limits_{i+j=k}\widetilde{H}_i(E_p) \otimes \widetilde{H}_j(E_q) \ \oplus$ 
$\ds\bigoplus\limits_{i+j=k-1} \Tor(\widetilde{H}_i(E_p),\widetilde{H}_j(E_q))$.
In particular:  $E_r$ acyclic $\implies$ 
$0 = \widetilde{H}_1(E_r) =$ 
$\widetilde{H}_0(E_p) \otimes \widetilde{H}_0(E_q)$, 
so at least one of $E_p$ and $E_q$ is path connected, thus $E_p * E_q$ is simply connected 
(cf. p. \pageref{5.30}) 
or still, $E_r$ is contractible and $r$ is a pointed homotopy equivalence.  
Therefore (i) and (ii) are equivalent.  To check (ii) $\Leftrightarrow$ (iii), use the algebra developed above.\\
\endgroup%%------------------------------------<< 

\begingroup%%----------------------------------->>
\fontsize{9pt}{11pt}\selectfont
\label{9.27}
\textbf{\small EXAMPLE} \ 
Let
$
\begin{cases}
\ X \\
\ Y
\end{cases}
$
be pointed connected CW spaces, $f:X \ra Y$ a pointed continuous function.  
Denote by $C_\pi$ the mapping cone of the Hurewicz fibration $\pi:E_f \ra X$ $-$then, specializing the preceding example, the projection $C_\pi \ra Y$ is a pointed homotopy equivalence iff $\exists$ $P$ such that one of 
$
\begin{cases}
\ \widetilde{H}_*(E_f) = \bigoplus\limits_i \widetilde{H}_i(E_f)\\
\ \widetilde{H}_*(\Omega Y) = \bigoplus\limits_j \widetilde{H}_j(\Omega Y)
\end{cases}
$
is $P$-primary and the other is uniquely $P$-divisible.  
To illustrate the situation when $P$ is the set of all primes, consider the short exact sequence
$0 \ra \Z \ra \Q \ra \Q/\Z \ra 0$ $-$then the mapping fiber of the arrow
$K(\Z,n+1) \ra K(\Q,n+1)$ is a $K(\Q/\Z,n)$ 
(cf. p. \pageref{5.31}).  
Furthermore, $\Omega K(\Q,n+1) = K(\Q,n)$ and 
$\widetilde{H}_*(\Q,n)$ is a uniquely divisible abelian group (being a vector space over $\Q$), while 
$\widetilde{H}_*(\Q/\Z,n)$ is a torsion abelian group 
(cf. p. \pageref{5.31a}).  
When $P = \emptyset$, there are two possibilities:
(1) \ $\widetilde{H}_*(E_f) = 0$;
(2) \ $\widetilde{H}_*(\Omega Y) = 0$.  
In the first case, $f$ is acyclic and in the second case, $Y$ is contractible and 
$\pi:E_f \ra X$ is a pointed homotopy equivalence.  
Consequently, if $\pi_1(Y) \neq 0$, then $f$ is acyclic iff the projection 
$C_\pi \ra Y$ is a pointed homotopy equivalence.
\vspi
%%----------------------------------------------------------------------------------------------68
[Note: \  A priori, $C_\pi$ is calculated in \bTOP but is viewed as an object in $\bTOP_*$.  As such, it has the same pointed homotopy type as the pointed mapping cone of $\pi$.]\\
\endgroup%%------------------------------------<< 

\begingroup%%----------------------------------->>
\fontsize{9pt}{11pt}\selectfont
\textbf{\small FACT} \ 
Suppose that $f:X \ra Y$ is acyclic.  Let $Z$ be any pointed space $-$then the arrow 
$[Y,Z] \ra$ $[X,Z]$ is injective.
\vspi
[The orbits of the action of $[\Sigma E_f,Z]$ on $[C_\pi,Z]$ are the fibers of the arrow $[C_\pi,Z] \ra [X,Z]$ 
(cf. p. \pageref{5.32}).  
But $\Sigma E_f$ is contractible in $\bTOP_*$, hence 
$[\Sigma E_f,Z]$ is the trivial group and, as noted above, one can replace $C_\pi$ by \mY.]\\
\endgroup%%------------------------------------<< 

%% ---------------------------<
\begin{proposition} \ %20
Let
$
\begin{cases}
\ X \\[-.1cm]
\ Y
\end{cases}
$
be pointed connected CW spaces.  
Suppose that $f:X \ra Y$ is a pointed continuous function with 
$\pi_0(E_f) = *$ $-$then $f$ is acyclic 
iff $f$ is a homology equivalence and $\pi_1(Y)$ operates nilpotently on the $H_q(E_f)$ $\forall \ q$.
\end{proposition}

[Consider the commutative diagram
\begin{tikzcd}%[sep=large]
W_f \ar{r} \ar{d} \ar{r} &Y \arrow[d,shift right=0.5,dash] \arrow[d,shift right=-0.5,dash] \\
Y  \arrow[r,shift right=0.5,dash] \arrow[r,shift right=-0.5,dash]  &Y
\end{tikzcd}
and apply the Hilton-Roitberg comparison theorem.]\\
%\ar[equal, shift left]{r} \ar[equal, shift right]{r} \\
%%-------------------------------

\begingroup%%----------------------------------->>
\fontsize{9pt}{11pt}\selectfont
\textbf{\small EXAMPLE} \ 
Take $X = \bS^3/\bSL(2,5)$, $Y = \bS^3$ $-$then the arrow $X \ra Y$ is an acyclic map 
(cf. p. \pageref{5.33}).\\
\endgroup%%------------------------------------<< 

\begingroup%%----------------------------------->>
\fontsize{9pt}{11pt}\selectfont
\textbf{\small FACT} \ 
Let
$
\begin{cases}
\ X \\
\ Y
\end{cases}
$
be pointed connected CW spaces, $f:X \ra Y$ a pointed continuous function.  
Denote by $C_f$ its mapping cone $-$then $f$ acyclic $\implies C_f$ contractible and $C_f$ contractible $\implies$ 
$f$ acyclic provided that $\pi_1(Y) = 0$.

[If $C_f$ is contractible and $Y$ is simply connected, then $f$ is a homology equivalence 
(cf. p. \pageref{5.34}) 
and $\pi_1(Y)$ operates trivially on the $H_q(E_f)$ $\forall \ q$, so Proposition 20 can be cited.]\\
\endgroup%%------------------------------------<< 

\begingroup%%----------------------------------->>
\fontsize{9pt}{11pt}\selectfont
\textbf{\small FACT} \ 
Let
$
\begin{cases}
\ X \\
\ Y
\end{cases}
$
be pointed connected CW spaces, $f:X \ra Y$ a pointed continuous function.  
Assume:  $X$ is acyclic and $f_*:\pi_1(X) \ra \pi_1(Y)$ is trivial $-$then $f$ is nullhomotopic.
\vspi
[Take $X$ to be a pointed connected CW complex, consider a lifting 
$\widetilde{f}:X \ra \widetilde{Y}$ of $f$, and show that 
$\widetilde{Y} \ra C_{\widetilde{f}}$ is an acyclic map.]
\vspi
\label{8.11}
[Note: \  It is a corollary that if $X$ is acyclic and $\text{Hom}(\pi_1(X,x_0),\pi_1(Y,y_0)) = *$, then $C(X,x_0;Y,y_0)$ is homotopically trivial.]\\
\endgroup%%------------------------------------<< 

\label{5.34e}
\begingroup%%----------------------------------->>
\fontsize{9pt}{11pt}\selectfont
Application: 
Let
$
X \ \& \
\begin{cases}
\ Y \\
\ Y^\prime
\end{cases}
$
be pointed connected CW spaces.  Suppose 
$f:X \ra Y$ \ $\&$ \ 
$f^\prime:X \ra Y^\prime$ 
are pointed continuous functions with $f$ acyclic $-$then there exists a pointed continuous function 
$g:Y \ra Y^\prime$ such that
$g \circx f \simeq f^\prime$ iff 
$\ker \pi_1(f) \subset \ker \pi_1(f^\prime)$.
\vspi
[Note: \  Up to pointed homotopy, $g$ is unique.]\\
\endgroup%%------------------------------------<< 

%%----------------------------------------------------------------------------------------------69
\label{9.49} %dmc mnft
\begin{proposition} \ %21
Let
$
\begin{cases}
\ X \\
\ Y
\end{cases}
$
be pointed connected CW spaces.  
Suppose that $f:X \ra Y$ is a pointed continuous function with $\pi_0(E_f) = *$ 
$-$then $f$ is a pointed homotopy equivalence 
iff $f$ is a homology equivalence and $\pi_1(X)$ operates nilpotently on the $\pi_q(E_f)$ $\forall \ q$.
\end{proposition}

[The stated condition on $\pi_1(X)$ implies that $\pi_1(Y)$ operates nilpotently on the 
$H_q(E_f)$ $\forall \ q$ 
(cf. p. \pageref{5.34a}), 
thus, by Proposition 20, $E_f$ is acyclic.  
But $E_f$ is also nilpotent.  
Therefore $E_f$ is contractible and $f:X \ra Y$ is a pointed homotopy equivalence.]\\

\label{11.18}
It will be convenient to insert here a technical addendum to the fibration spectral sequence.

Notation: A continuous function 
$f:X \ra Y$ 
induces a functor 
$f^*:\bLCCS_Y \ra$ $\bLCCS_X$ or still, a functor 
$f^*:[(\Pi Y)^\OP,\bAB] \ra$ $[(\Pi X)^\OP,\bAB]$ (cf. $\S 4$, Proposition 25).  
If \mX is a subspace of \mY and $f$ is the inclusion, one writes 
$\restr{\sG}X$ instead of $f^*\sG$.

Let
$f:X \ra Y$ be a Hurewicz fibration, where 
$
\begin{cases}
\ X \\[-.1cm]
\ Y
\end{cases}
$
and the $X_y$ are path connected.  
Fix a cofunctor 
$\sG:\Pi Y$ $\ra$ $\bAB$ $-$then $\forall \ y \in Y$, the projection 
$X_y \ra Y$ is inessential, hence 
$f^* \restr{\sG}{X_y}$ is constant.  So, $\forall \ q \geq 0$, there is a cofunctor 
$\sH_q(f;\sG):\Pi Y$ $\ra$ $\bAB$ that assigns to each $y \in Y$ the singular homology group 
$H_q(X_y;f^*\restr{\sG}{X_y})$ and the fibration spectral sequence assumes the form 
$E_{p,q}^2 \approx$ 
$H_p(Y;\sH_q(f;\sG)) \Rightarrow$ 
$H_{p+q}(X;f^*\sG)$.

[Note: \ A morphism 
$[\tau]:y_0 \ra y_1$ determines a homotopy equivalence 
$X_{y_0} \ra$ $X_{y_1}$ 
(cf. p. \pageref{5.34b}) 
and an isomorphism 
$\sG[\tau]:\sG y_1 \ra$ $\sG y_0$, thus 
$\sH_q(f;\sG)[\tau]$ is the composite 
$H_q(X_{y_1};\sG y_1) \ra$ 
$H_q(X_{y_0};\sG y_1) \ra$
$H_q(X_{y_0};\sG y_0)$.]\\

\begin{proposition} \ %22
Let
$
\begin{cases}
\ X \\[-.1cm]
\ Y
\end{cases}
$
be pointed connected CW spaces, $f:X \ra Y$ a pointed continuous function $-$then $f$ is acyclic iff for every locally constant coefficient system $\sG$ on $Y$, the induced map $f_*:H_*(X;f^*\sG) \ra H_*(Y;\sG)$ is an isomorphism.
\end{proposition}

[Upon passing to the mapping track, one can assume that $f$ is a pointed Hurewicz fibration.

Necessity: 
$\forall \ y \in Y$, $X_y$ is acyclic, thus from the universal coefficient theorem, 
$\forall \ q > 0$, 
$H_q(X_y;f^* \restr{\sG}{X_y}) = 0$.  Accordingly, the edge homomorphism 
$e_H:E_{p,0}^\infty \ra$ $E_{p,0}^2$ is an isomorphism, so 
$\forall \ p \geq 0$, 
$H_p(X;f^*\sG) \approx$ $H_p(Y;\sG)$.

Sufficiency: The integral group ring $\Z[\pi_1(Y)]$ is a right $\pi_1(Y)$-module.  
Viewed as a locally constant coefficient system on $Y$, its homology is that of $\widetilde{Y}$.  
Form the pullback square 
\begin{tikzcd}%[sep=large]
{X \times_Y \widetilde{Y}} \ar{d} \ar{r}{f^\prime} &{\widetilde{Y}} \ar{d}\\
{X} \ar{r}[swap]{f} &{Y}
\end{tikzcd}
$-$then 
$H_*(X \times_Y \widetilde{Y}) \approx$ 
$H_*(X;f^*(\Z[\pi_1(Y)]))$ and 
$f_*^\prime:H_*(X \times_Y 
%%----------------------------------------------------------------------------------------------70
\widetilde{Y}) \ra$ $H_*(\widetilde{Y})$ is the composite 
$H_*(X \times_Y \widetilde{Y}) \ra$ 
$H_*(X;f^*(\Z[\pi_1(Y)])) \overset{f_*}{\lra}$ 
$H_*(Y;\Z[\pi_1(Y)]) \ra$ 
$H_*(\widetilde{Y})$.  
By hypothesis, $f_*$ is an isomorphism, hence $f_*^\prime$ is too.  
Since $\widetilde{Y}$ is simply connected, $E_{f^\prime}$ is path connected.  Consider the commutative diagram
\begin{tikzcd}%[sep=large]
{X \times_Y \widetilde{Y}} \ar{d}[swap]{f^\prime} \ar{r}{f^\prime} &{\widetilde{Y}} \ar{d}{\id.}\\
{\widetilde{Y}} \ar{r}[swap]{\id} &{\widetilde{Y}}
\end{tikzcd}
Owing to the Hilton-Roitberg comparison theorem, the projection 
$E_{f^\prime} \ra *$ is a homology equivalence.  Therefore $E_f$ is acyclic.]\\

\label{5.34d}
Application: Let $X$, $Y$, $Z$ be pointed connected CW spaces.  Suppose that
$
\begin{cases}
\ f:X \ra Y \\[-.1cm]
\ g:Y \ra Z
\end{cases}
$
are pointed continuous functions.  Assume: $f$ is acyclic $-$then $g$ is acyclic iff $g \circx f$ is acyclic.\\

\begingroup%%----------------------------------->>
\fontsize{9pt}{11pt}\selectfont
\textbf{\small FACT} \ 
Let $X \overset{f}{\lla} Z \overset {g}{\lra} Y$ be a pointed 2-source, where 
$
\begin{cases}
\ X\\
\ Y
\end{cases}
\& \ Z
$
are pointed connected CW spaces.  Consider the pushout square 
\begin{tikzcd}[sep=large]
{Z} \ar{d}[swap]{f} \ar{r}{g} &{Y} \ar{d}{\eta}\\
{X} \ar{r}[swap]{\xi} &{P}
\end{tikzcd}
. \ 
Assume: $f$ is a cofibration $-$then $f$ (or $g$) acyclic $\implies$ $\eta$ (or $\xi$) acyclic.\\
\endgroup%%------------------------------------<< 

\index{Plus Construction}
\textbf{\small PLUS CONSTRUCTION} \ 
Fix a pointed connected CW space $X$.  
Let $N$ be a perfect normal subgroup of $\pi_1(X)$ $-$then there exists a pointed connected CW space $X_N^+$ 
and an acyclic map $f_N^+:X \ra X_N^+$ such that $\ker \pi_1(f_N^+) = N$ $(\implies \pi_1(X_N^+) \approx \pi_1(X)/N)$.  
Moreover, the pointed homotopy type of $X_N^+$ is unique, i.e., if $g_N^+:X \ra Y_N^+$ is acyclic and if $\ker \pi_1(g_N^+) = N$, 
then there is a pointed homotopy equivalence $\phi:X_N^+ \ra Y_N^+$ such that $\phi \circx f_N^+ \simeq g_N^+$.

[Existence:  We shall first deal with the case when $N = \pi_1(X)$.  
Thus let $\{\alpha\}$ be a set of generators for $\pi_1(X)$.  
Represent 
$\alpha$ by $f_\alpha:\bS^1 \ra X$ and put 
$X_1 = \bigl(\coprod\limits_\alpha \bD^2\bigr) \sqcup_f X$ 
$\bigl(f = \coprod f_\alpha\bigr)$ to obtain a relative CW complex $(X_1,X)$ with $\pi_1(X) = 0$ 
(cf. p. \pageref{5.34c}).  
Consider the exact sequence 
$H_2(X_1) \ra H_2(X_1,X) \ra H_1(X):$ 
(a) $\pi_2(X_1) \approx H_2(X_1)$; 
(b) $H_2(X_1,X)$ is free abelian on generators $\omega_\alpha$, say; 
(c) $H_1(X) = 0$.  
Given $\alpha$, choose a continuous function 
$g_\alpha:\bS^2 \ra X_1$ 
such that the homotopy class 
$[g_\alpha]$ maps to $\omega_\alpha$ 
under the composite 
$\pi_2(X_1) \ra H_2(X_1) \ra H_2(X_1,X)$.  Put 
$X_N^+ = \bigl(\coprod\limits_\alpha \bD^3\bigr) \sqcup_g X_1$ 
$\bigl(g = \coprod\limits_\alpha g_\alpha \bigr)$ $-$then the pair $(X_N^+,X_1)$ 
is a relative CW complex with $\pi_1(X_N^+) = 0$.  
The inclusion $X \ra X_N^+$ is a closed cofibration.  
In addition, it is a homology equivalence (for 
$H_*(X_n^+,X)$ $=$ $0)$, hence is an acyclic map (cf. Proposition 20).  
Turning to the general case, let 
$\widetilde{X}_N$ be the covering space of $X$ corresponding to $N$ 
(so $\pi_1(\widetilde{X}_N) \approx N$).  
Apply the foregoing procedure to $\widetilde{X}_N$ to get an acyclic closed cofibration 
$\widetilde{f}_N^+:\widetilde{X}_N \ra \widetilde{X}_N^+$, where 
$\widetilde{X}_N^+$ is simply connected.  Define
%%----------------------------------------------------------------------------------------------71
$\widetilde{X}_N^+$ by the pushout square
\begin{tikzcd}%[sep=large]
\widetilde{X}_N \arrow{r}{\widetilde{f}_N^+} \arrow{d} &\widetilde{X}_N^+ \arrow{d}\\
X \arrow{r}[swap]{f_N^+} &X_N^+
\end{tikzcd}
.  
Thanks to Proposition 7, $X_N^+$ is a pointed connected CW space.  
And: $f_N^+$ is an acyclic closed cofibration 
(cf. p. \pageref{5.34d}).  
Finally, the Van Kampen theorem implies that $\pi_1(X_N^+) \approx \pi_1(X) / N$.

Uniqueness:  Since \mN =
$
\begin{cases}
\ \ker \pi_1(f_N^+) \\[-.1cm]
\ \ker \pi_1(g_N^+) 
\end{cases}
$
\hspace{-.25cm}, 
there exists a pointed continuous function $\phi:X_N^+ \ra Y_N^+$ such that $\phi \circx f_N^+ \simeq g_N^+$ 
(cf. p. \pageref{5.34e}).  
But
$
\begin{cases}
\ f_N^+ \\[-.1cm]
\ g_N^+
\end{cases}
$
acyclic $\implies \phi$ acyclic and $\phi_*:\pi_1(X_N^+) \ra \pi_1(Y_N^+)$ is necessarily an isomorphism.  Therefore $\phi$ is a pointed homotopy equivalence (cf. Proposition 19).]

[Note: \  $X_N^+$ is called the 
\un{plus construction}
%\index{plus construction} 
with respect to $N$.  
Like an Eilenberg-MacLane space, $X_N^+$  is really a pointed homotopy type, thus, while a given representative may have a certain property, it need not be true that all representatives do.  
As for $\phi$, if $f_N^+$ is an acyclic closed cofibration and if $g_N^+$ is another such, then matters can be arranged so that there is commutativity on the nose: $\phi \circx f_N^+ = g_N^+$.  
This in turn means that $\phi$ is a homotopy equivalence in $X\backslash\bTOP$ (cf. $\S3$, Proposition 13).]\\

\begingroup%%----------------------------------->> 
\fontsize{9pt}{11pt}\selectfont
One can interpret $X_N^+$ as a  representing object of the functor on the homotopy category of pointed connected CW spaces which assigns to each $Y$ the set of all $[f] \in [X,Y]: \ker \pi_1(f) \supset N$.\\
\endgroup%%------------------------------------<< 

Different notation is used when $N = \pi_1(X)_\text{per}$, the perfect radical of 
$\pi_1(X): X_N^+$ is replaced by $X^+$ and $f_N^+:X \ra X_N^+$ is replaced by $i^+:X \ra X^+$.  
Example:  $X$ acyclic $\implies X^+$ contractible.

[Note: \  The perfect radical of  $\pi_1(X^+)$ is trivial 
(cf. p. \pageref{5.34f}).]

Examples: Let
$
\begin{cases}
\ X \\[-.1cm]
\ Y
\end{cases}
$
be pointed connected CW spaces $-$then 
(1) $X^+ \times Y^+$ is a model for $(X \times Y)^+$; 
(2) $X^+ \vee Y^+$ is a model for $(X \vee Y)^+$; 
(3) $X^+ \# Y^+$ is a model for $(X \# Y)^+$.\\

\begingroup%%----------------------------------->> 
\fontsize{9pt}{11pt}\selectfont
\index{homology spheres}
\textbf{\small EXAMPLE \  (\un{Homology Spheres})} \ 
Fix $n > 1$.  Suppose that $X$ is a pointed connected CW space such that $\widetilde{H}_q(X) =
\begin{cases}
\ \Z \hspace{0.5cm} (q = n) \\
\  0 \hspace{0.5cm} (q \neq n)
\end{cases}
$
$-$then $\pi_1(X)$ is perfect and $X^+$ has the same pointed homotopy type as $\bS^n$.\\
\endgroup%%------------------------------------<< 

\begingroup%%----------------------------------->> 
\fontsize{9pt}{11pt}\selectfont
\textbf{\small FACT} \ 
Let $X$ be a pointed connected CW space $-$then for any pointed acyclic CW space $Z$, the arrow 
$[Z,E_{i^+}] \ra [Z,X]$ is bijective.
\vspi
[Note: \  The central extension 
$1 \ra \text{im }\pi_2(X^+) \ra \pi_1(E_{i^+}) \ra\pi_1(X)_\text{per} \ra 1$ is universal.]\\
\endgroup%%------------------------------------<< 

%%----------------------------------------------------------------------------------------------72
Convention: Henceforth it will be assumed that $i^+:X \ra X^+$ is an acyclic closed cofibration.\\
\label{5.36a}

\textbf{\small LEMMA} \ 
Let
$
\begin{cases}
\ X \\[-.1cm]
\ Y
\end{cases}
$
be pointed connected CW spaces.  Suppose that $f:X \ra Y$ is a pointed continuous function $-$then there is a pointed continuous function $f^+:X^+ \ra Y^+$ rendering the diagram 
\begin{tikzcd}%[sep=large]
X \arrow{r}{f} \arrow{d} &Y \arrow{d}\\
X^+ \arrow{r}[swap]{f^+} &Y^+
\end{tikzcd}
commutative, $f^+$ being unique up to pointed homotopy.\\

\label{5.41}
Application: Let
$
\begin{cases}
\ X \\[-.1cm]
\ Y
\end{cases}
$
be pointed connected CW spaces.  Assume: $X$ and $Y$ have the same pointed homotopy type $-$they $X^+$ and $Y^+$ have the same pointed homotopy type.\\

\begin{proposition} \ %23
Let $X$ be a pointed connected CW space.  
Denote by 
$\widetilde{X}_N$ the covering space of $X$ corresponding to $N$, where $N$ is a normal subgroup of $\pi_1(X)$ containing 
$\pi_1(X)_\text{per}$ $-$then $\widetilde{X}_N^+$ 
has the same pointed homotopy type as the covering space of $X^+$ corresponding to the normal subgroup 
$N/\pi_1(X)_\text{per}$  of 
$\pi_1(X^+) \approx$ 
$\pi_1(X)/\pi_1(X)_\text{per}$.
\end{proposition}

[The pointed homotopy type of $\widetilde{X}_N$ can be calculated as the mapping fiber of the composite 
$X \ra X[1] = K(\pi_1(X),1) \ra K(\pi_1(X)/N,1)$.  
This arrow factors through $X^+$ and $\pi_1(X)/N \approx (\pi_1(X)/\pi_1(X)_\text{per})/(N/\pi_1(X)_\text{per})$.]\\

\label{9.100}
\begingroup%%----------------------------------->> 
\fontsize{9pt}{11pt}\selectfont
Notation:  Given a group $G$, put $BG = K(G,1)$.\\
\endgroup%%------------------------------------<< 

\begingroup%%----------------------------------->> 
\fontsize{9pt}{11pt}\selectfont
\textbf{\small EXAMPLE} \ 
$BG_\text{per}$ is the covering space of $BG$ corresponding to $G_\text{per}$.  There is an arrow $BG_\text{per}^+ \ra BG^+$ and $BG_\text{per}^+$ ``is'' the universal covering space of $BG^+$.\\
\endgroup%%------------------------------------<< 

\begingroup%%----------------------------------->> 
\fontsize{9pt}{11pt}\selectfont
\textbf{\small EXAMPLE} \ 
Let $A$ be a ring with unit $-$then the fundamental group of the mapping fiber of $B\bGL(A) \ra B\bGL(A)^+$ is isomorphic to $\textbf{ST}(A)$.\\
\endgroup%%------------------------------------<< 

\begin{proposition} \ %24
Let
$
\begin{cases}
\ X \\
\ Y
\end{cases}
$
be pointed connected CW spaces.  Suppose that $f:X \ra Y$ is a pointed continuous function with $\pi_0(E_f) = *$ $-$then $\pi_0(E_{f^+}) = *$ and the perfect radical of $\pi_1(E_{f^+})$ is trivial.
\end{proposition}

[Note: \  It follows that there is a commutative triangle
\begin{tikzcd}%[sep=large]
{E_f} \arrow{r} \ar{d} &{E_f^+} \ar{ld}\\
{E_{f^+}}
\end{tikzcd}
.]\\[0.5cm]

%%----------------------------------------------------------------------------------------------73
\begingroup%%----------------------------------->> 
\fontsize{9pt}{11pt}\selectfont
\textbf{\small FACT} \ 
Let
$
\begin{cases}
\ X \\[-.1cm]
\ Y
\end{cases}
$
be pointed connected CW spaces.  
Suppose that $f:X \ra Y$ is a pointed continuous function with $\pi_0(E_f) = *$ 
$-$then the arrow $E_f^+ \ra E_{f^+}$ is a pointed homotopy equivalence if $\pi_1(Y)_\text{per}$ 
is trivial or if $E_f^+$ is nilpotent and $\pi_1(Y)_\text{per}$ operates nilpotently on the $H_q(E_f)$ $\forall \ q$.
\vspi
[Note: \  $\pi_1(Y)_\text{per}$ operates nilpotently on the $H_q(E_f)$ $\forall \ q$ iff $\pi_1(Y)_\text{per}$ operates trivially on the $H_q(E_f)$ $\forall \ q$ 
(cf. p. \pageref{5.35}).]
\\[.1cm]
\endgroup%%------------------------------------<< 

\label{5.39}
\label{5.37a}
\index{central extensions}
\begingroup%%----------------------------------->> 
\fontsize{9pt}{11pt}\selectfont
\textbf{\small EXAMPLE \ \un{(Central Extensions)}} \ 
Let $\pi$ and $G$ be groups, where $\pi$ is abelian.  Consider a central extension 
$1 \ra \pi \ra \Pi \ra G \ra 1$ 
$-$then $B\pi$ can be identified with the mapping fiber of the arrow $B\Pi^+ \ra BG^+$.
\vspi
[Since $\pi$ is abelian, $B\pi = B\pi^+$ and $G$ $(= \pi_1(BG))$ operates trivially on $\pi$, 
hence operates trivially on the $H_q(B\pi)$ $\forall \ q$.]
\\[.1cm]
\endgroup%%------------------------------------<< 

\begingroup%%----------------------------------->> 
\fontsize{9pt}{11pt}\selectfont
\textbf{\small EXAMPLE} \ 
Let $G$ be an abelian group $-$then there is a universal central extension 
$1 \ra$ 
$G \ra$ 
$\alpha G \ra$ \
$\beta G \ra 1$ 
(cf. p. \pageref{5.36}).  
Specializing the preceding example, the mapping fiber of the arrow 
$K(\alpha G,1)^+ \ra K(\beta G,1)^+$ is a $K(G,1)$ and $K(\beta G,1)^+$ is a $K(G,2)$.
\vspi
[Recall that $\alpha G$ is acyclic, thus $K(\alpha G,1)^+$ is contractible.]\\
\endgroup%%------------------------------------<< 

\begin{proposition} \ %25
Let
$
\begin{cases}
\ X \\[-.1cm]
\ Y
\end{cases}
$
be pointed connected CW spaces.  Suppose that $f:X \ra Y$ is a pointed continuous function for which the normal closure of $f_*(\pi_1(X)_\text{per})$ is $\pi_1(Y)_\text{per}$ $-$then the adjunction space $X^+ \sqcup_f Y$ represents $Y^+$.
\end{proposition}

[Since $i^+:X \ra X^+$ is an acyclic closed cofibration, the same is true of the inclusion 
$Y \ra X^+ \sqcup_f Y$ 
(cf. p. \pageref{5.36a}).  
On the other hand, by Van Kampen, the fundamental group of 
$X^+ \sqcup_f Y$ is isomorphic to $\pi_1(Y)$ modulo the normal closure of 
$f_*(\pi_1(X)_\text{per})$, i.e., to 
$\pi_1(Y)/\pi_1(Y)_\text{per}$.]\\

\begingroup%%----------------------------------->> 
\fontsize{9pt}{11pt}\selectfont
\label{18.31}
\index{Algebraic K-Theory (example)}
\textbf{\small EXAMPLE  \ \un{(Algebraic K-Theory)}} \ 
Let $A$ be a ring with unit $-$then by definition, 
$K_0(A)$ is the Grothendieck group attached to the category of finitely generated projective $A$-modules and for 
$n \geq 1$, $K_n(A)$ is taken to be the homotopy group $\pi_n(B\bGL(A)^+)$.  
While it is immediate that $K_0$ is a functor from 
\textbf{RG} to \bAB, 
the plus construction requires some choices, so to guarantee that $K_n$ is a functor one has to fix the data.  
Thus first construct $B\bGL(\Z)^+)$.  This done, define 
$B\bGL(A)^+)$ by the pushout square
\begin{tikzcd}[sep=large]
{B\bGL(\Z)} \ar{r} \ar{d} &{B\bGL(A)} \ar{d}\\
{B\bGL(\Z)^+} \ar{r} &{B\bGL(A)^+}
\end{tikzcd}
.  Here, Proposition 25 comes in (the normal closure of $\text{im}(\bE(\Z) \ra \bE(A))$ is $\bE(A)$).  
Observe that the $K_n$ preserve products: 
$K_n(A^\prime \times A^{\prime\prime}) \approx K_n(A^\prime) \times K_n(A^{\prime\prime})$.\\
\indent\indent $(n = 1)$ \quadx 
$K_1(A) = \pi_1(B\bGL(A)^+) \approx$ 
$\pi_1(B\bGL(A))/\pi_1(B\bGL(A))_\text{per} \approx$ 
$\bGL(A)/[\bGL(A),\bGL(A)]$ $=$ $H_1(\bGL(A))$.
\\
\indent\indent $(n = 2)$ \quadx 
$K_2(A) = \pi_2(B\bGL(A)^+) \approx \pi_2(B\bE(A)^+) \approx H_2(B\bE(A)^+) \approx 
H_2(B\bE(A)) \approx H_2(\bE(A))$.
\vspi
%%----------------------------------------------------------------------------------------------74
[Note: \  The central extension 
$1 \ra K_2(A) \ra \textbf{ST}(A) \ra \bE(A) \ra 1$
is universal (cf. p. \pageref{5.37}) 
and $BK_2(A)$ can be identified with the mapping fiber of the arrow 
$B\textbf{ST}(A)^+ \ra B\bE(A)^+$.]
\\
\indent\indent $(n = 3)$ \quadx 
$K_3(A) = \pi_3(B\bGL(A)^+) \approx$ 
$\pi_3(B\bE(A)^+) \approx$ 
$\pi_3(B\textbf{ST}(A)^+) \approx$ 
$H_3(B\textbf{ST}(A)^+) \approx$ 
$H_3(B\textbf{ST}(A)) =$ $H_3(\textbf{ST}(A))$.
\vspi
There is no known homological interpretation of $K_4$ and beyond.\\
\endgroup%%------------------------------------<< 

\index{relative Algebraic K-Theory (example)}
\begingroup%%----------------------------------->> 
\fontsize{9pt}{11pt}\selectfont
\textbf{\small EXAMPLE \ \un{(Relative Algebraic K$-$Theory)}} \ 
Let \mA be a ring with unit, $I \subset A$ a two sided ideal.  
Write 
$\widehat{\bGL}(A/I)$ for the image of $\bGL(A)$ in $\bGL(A/I)$ 
$-$then 
$\widehat{\bGL}(A/I) \supset {\bE}(A/I)$, thus 
$\widehat{\bGL}(A/I)$ is normal and 
$G_{A,I} = {\bGL}(A/I) / \widehat{\bGL}(A/I)$ is abelian.  
Since $B\widehat{\bGL}(A/I)^+$ can be identified with the mapping fiber of the arrow 
$B{\bGL}(A/I)^+ \ra BG_{A,I}^+ (= BG_{A,I})$ 
(cf. p. \pageref{5.37a}), 
it follows that 
$\pi_n(B\widehat{\bGL}(A/I)^+) \approx$ $\pi_n(B{\bGL}(A/I)^+)$ $(n > 1)$ but 
$\pi_1(B\widehat{\bGL}(A/I)^+) \approx$ $\text{im}(K_1(A) \ra K_1(A/I))$ 
and there is a short exact sequence 
$0 \ra$  
$\pi_1(B\widehat{\bGL}(A/I)^+) \ra$ 
$K_1(A/I) \ra$  
$G_{A,I} \ra 0$.  
If $\sK(A,I)$ is the mapping fiber of the arrow 
$B{\bGL}(A)^+ \ra$ 
$B\widehat{\bGL}(A/I)^+$, then 
$\sK(A,I)$ is path connected, so letting 
$K_n(A,I) = \pi_n(\sK(A,I))$ $(n \geq 1)$, 
one obtains a functorial long exact sequence \ 
$\cdots \ra K_{n+1}(A/I) \ra$ 
$K_n(A,I)$ $\ra$ $K_n(A)$ 
$\ra$ $K_n(A/I) \ra$ 
$\cdots \ra K_1(A,I) \ra$ 
$K_1(A) \ra K_1(A/I)$.\\
\endgroup%%------------------------------------<< 

\begin{proposition} \ %26
Let $X$ be pointed connected CW space.  
Put 
$\pi = \pi_1(X)$ and denote by $\widetilde{X}_\text{per}$ the mapping fiber of the composite 
$X \ra K(\pi,1) \ra K(\pi/\pi_\text{per},1)$.  
Assume: $\pi/\pi_\text{per}$ is nilpotent and $\pi/\pi_\text{per}$ operates nilpotently on the 
$H_q(\widetilde{X}_\text{per})$ $\forall \ q$ $-$then $X^+$ is nilpotent.
\end{proposition}

[Since $(\pi/\pi_\text{per})$ is trivial 
(cf. p. \pageref{5.38}), 
$\widetilde{X}_\text{per}^+$ can be identified with the mapping fiber of the composite 
$X^+ \ra K(\pi,1)^+ \ra K(\pi/\pi_\text{per},1)^+$ 
(cf. p. \pageref{5.39}).  
By construction, $\widetilde{X}_\text{per}^+$ is simply connected (cf. Proposition 23), hence nilpotent.  
But $K(\pi/\pi_\text{per},1)^+ = K(\pi/\pi_\text{per},1)$ is also nilpotent.  
Therefore, bearing in mind that the inclusion 
$\widetilde{X}_\text{per} \ra \widetilde{X}_\text{per}^+$ is a homology equivalence, it follows that 
$X^+$  is nilpotent 
(cf. p. \pageref{5.40}).]\\

\begingroup%%----------------------------------->> 
\fontsize{9pt}{11pt}\selectfont
\textbf{\small FACT} \ 
Let $G$ be a group.  Fix $\phi \in \Aut G$.  
Assume:  Given $g_1, \ldots, g_n \in G$, $\exists \ g \in G$: $\phi(g_i) = gg_ig^{-1}$ 
$(1 \leq i \leq n)$ $-$then $\phi_*:H_*(G) \ra H_*(G)$ is the identity.\\
\endgroup%%------------------------------------<<

\begingroup%%----------------------------------->> 
\fontsize{9pt}{11pt}\selectfont
Application:  Let $G$ be a group. Let \mK be a normal subgroup of \mG which is the colimit of subgroups 
$K_n$ $(n \in \N)$ such that $\forall \ n$, $G = K\cdot \Cen_G(K_n)$ $-$then \mG operates trivially on $H_*(K)$.\\
\endgroup%%------------------------------------<<

\label{8.9}
\label{9.89}
\begingroup%%----------------------------------->> 
\fontsize{9pt}{11pt}\selectfont
\textbf{\small EXAMPLE} \ 
Let $A$ be a ring with unit $-$then 
$B\bGL(A)^+$ is nilpotent.  To see this, consider the short exact sequence
$1 \ra$ 
$\bE(A) \ra$ 
$\bGL(A) \ra$ 
$\bGL(A)/\bE(A) \ra 1$.  
Here, 
$\bE(A) = \bGL(A)_\text{per}$ and $B\bE(A)$ is the mapping fiber of the arrow 
$B\bGL(A) \ra$ $K(\bGL(A)/\bE(A),1)$.  
The quotient 
$\bGL(A)/\bE(A)$ is abelian, hence nilpotent.  
On the other hand, if $\bE(n,A)$ is the subgroup of $\bGL(n,A)$ consisting of the elementary matrices, then 
$\bE(A) = \colim \ \bE(n,A)$ and $\forall \ n$, 
$\bGL(A) = \bE(A) \cdot \Cen_{\bGL(A)}(\bE(n,A))$, so 
$\bGL(A)$ operates trivially on $H_*(\bE(A))$.  That 
$B\bGL(A)^+$ is nilpotent is therefore a consequence of Proposition 26.
\vspi
%%----------------------------------------------------------------------------------------------75
[Note: \ More is true.  Thus define a homomorphism \ 
$\oplus:\bGL(A) \times \bGL(A) \ra \bGL(A)$ by $(X,Y) \ra X \oplus Y$, where 
$
(X \oplus Y)_{ij} =
\begin{cases}
\ x_{kl} \quadx (i = 2k - 1, j = 2l - 1)\\
\ y_{kl} \quadx (i = 2k, j = 2l)
\end{cases}
\& \ 0
$
otherwise $-$then Loday\footnote[2]{\textit{Ann. Sci. \'Ecole Norm. Sup.} \textbf{9} (1976), 309-377.}
has shown that the composite 
$B\bGL(A)^+ \times B\bGL(A)^+ \ra B(\bGL(A) \times \bGL(A))^+$
$\ra B\bGL(A)^+$ serves to equip $B\bGL(A)^+$ with the structure of a homotopy commutative H-group.  
In particular: $B\bGL(A)^+$ is abelian.]\\
\endgroup%%------------------------------------<<

\label{14.185}
\begingroup%%----------------------------------->> 
\fontsize{9pt}{11pt}\selectfont
\textbf{\small EXAMPLE} \ 
Let $A$ be a ring with unit.  
Write $\bU\bT(A)$ for the ring of upper triangular 2-by-2 matrices with entries in $A$ $-$then 
the projection 
$p:\bU\bT(A) \ra A \times A$ 
$\bigl(
p
\begin{pmatrix}
a_1 &a\\
0 &a_2\\
\end{pmatrix}
= (a_1,a_2)\bigr)
$
induces an epimorphism 
$p:\bGL(\bU\bT(A)) \ra \bGL(A \times A)$.  
Its kernel is not perfect, therefore 
$Bp:B\bGL(\bU\bT(A)) \ra$ 
$B\bGL(A \times A)$ is not acyclic.  
Nevertheless, $Bp$ is a homology equivalence.  
Consider now the commutative diagram
\begin{tikzcd}[sep=large]
{B\bGL(\bU\bT(A))} \ar{d}[swap]{Bp} \ar{r} &{B\bGL(\bU\bT(A))^+} \ar{d}{{Bp}^+}\\
{B\bGL(A \times A)} \ar{r}  &{B\bGL(A \times A)^+}
\end{tikzcd}
.  
Since the horizontal arrows are homology equivalences, ${Bp}^+$ is a pointed homotopy equivalence, so 
$\forall \ n \geq 1$, 
$K_n(\bU\bT(A)) \approx$ $K_n(A) \times K_n(A)$.
\vspi
[Note: \ ${Bp}^+$ is acyclic (cf. Proposition 19), thus the composite 
$B\bGL(\bU\bT(A)) \overset{Bp}{\lra}$ 
$B\bGL(A \times A) \ra$ 
$B\bGL(A \times A)^+$ is acyclic even though $Bp$ is not.]\\
\endgroup%%------------------------------------<<

\begingroup%%----------------------------------->> 
\fontsize{9pt}{11pt}\selectfont
\textbf{\small FACT} \ 
Let \mG be a group.  Assume:
\\
\indent\indent $(\oplus)$ \quadx There is a homomorphism 
$\oplus:G \times G \ra G$ such that for any finite set 
$\{g_1, \ldots, g_n\} \subset G$, $\exists$ 
$
\begin{cases}
\ u \\
\ v
\end{cases}
\in G:
$
$
\begin{cases}
\ u(g_i \oplus e)u^{-1} = g_i \\
\ v(e \oplus g_i) v^{-1} = g_i
\end{cases}
(i = 1, \ldots, n).
$
\\
\indent\indent $(\frak{p})$ There is a homomorphism 
$\frak{p}:G \ra G$ such that for any finite set 
$\{g_1, \ldots, g_n\} \subset G$, $\exists$ $\rho \in G$: 
$\rho(g_i \oplus \frak{p} g_i)\rho^{-1} = g_i$ $(i = 1, \ldots, n)$.
\vspi
Then $G$ is acyclic.
\vspi
[Fix a field of coefficients \bk.  
Let $\Delta:G \ra G \times G$ be the diagonal map $-$then 
$\frak{p}$ and 
$\oplus \circx (\id \times \frak{p}) \circx \Delta$ 
operate in the same way on homology.  
Since $H_1(G;\bk) = 0$, one can take $n > 1$ and assume inductively that 
$H_q(G;\bk) = 0$ $(0 < q < n)$.  
Let $x \in H_n(G;\bk)$: 
$\frak{p}_*(x) =$ 
$(\oplus \circx (\id \times \frak{p}) \circx \Delta)_*(x) =$ 
$\oplus_*(x \otimes 1 + 1 \otimes \frak{p}_*(x)) =$ 
$x + \frak{p}_*(x)$ $\implies$ $x = 0$.]\\
\endgroup%%------------------------------------<<

\index{delooping Algebraic K-Theory}
\begingroup%%----------------------------------->> 
\fontsize{9pt}{11pt}\selectfont
\textbf{\small EXAMPLE \ \un{(Delooping Algebraic K-Theory)}} \ 
Let \mA be a ring with unit.  Denote by $\Gamma A$ the set of all functions 
$X:\N \times \N \ra A$ such that 
$\forall \ i$, $\#\{j:X_{ij} \neq 0\} < \omega$ and 
$\forall \ j$, $\#\{i:X_{ij} \neq 0\} < \omega$ 
$-$then 
$\Gamma A$ is a ring with unit containing $\ov{A}$ as a two sided ideal.  
$\Gamma A$ is called the 
\un{cone}
\index{cone (of a ring)} 
of $A$ and the quotient 
$\Sigma A =$ $\Gamma A /\ov{A}$ is called the 
\un{suspension}
\index{suspension (of a ring)} 
of $A$.  
Define a homomorphism 
$\oplus:\Gamma A \times \Gamma A \ra$ $\Gamma A$ by 
$(X,Y) \ra$ $X \oplus Y$, where 
$X \oplus Y)_{ij} =$
$
\begin{cases}
\ x_{kl} \ \ (i = 2k - 1, j = 2l - 1) \\
\ y_{kl} \ \ (i = 2k, j = 2l) 
\end{cases}
\& \ 0
$
otherwise and define a homomorphism 
$\frak{p}:\Gamma A \ra$ $\Gamma A$ by
$\frak{p}(X)_{ij} = X_{mn}$ 
if
$
\begin{cases}
\ i = 2^k(2m - 1) \\
\ j = 2^k(2n - 1) 
\end{cases}
$
for some $k, m, n \ \& \ 0$ otherwise.  
Evidently, 
$X \oplus \frak{p} X = \frak{p}X$ for all 
%%----------------------------------------------------------------------------------------------76
$X \in \Gamma A$ and 
$
\begin{cases}
\ \otimes \\
\ \frak{p} 
\end{cases}
$
induce homomorphisms 
$\otimes : \bGL(\Gamma A) \times \bGL(\Gamma A) \ra \bGL(\Gamma A)$ $\&$ 
$\frak{p}:\bGL(\Gamma A) \ra \bGL(\Gamma A)$ satifying the preceeding assumptions.  
Therefore $\bGL(\Gamma A)$ is acyclic, so 
$\bGL(\Gamma A) = \bE(\Gamma A)$.  
Taking into account the exact sequences 
$1 \ra \ov{\bGL}(\ov{A}) \ra \ov{\bGL}(\Gamma A) \ra \ov{\bGL}(\Sigma A)$,
$\ov{\bE}(\Gamma A) \ra \ov{\bE}(\Sigma A) \ra 1$, it follows that there is an exact sequence 
$1 \ra \bGL(A) \ra \bGL(\Gamma A) \ra \bE(\Sigma A) \ra 1$.  
The mapping fiber of the arrow 
$B\bGL(\Gamma A)^+ \ra B\bE(\Sigma A)^+$ is $B\bGL(A)^+$.  Since $B\bGL(\Gamma A)^+$ is contractible, this means that in 
$\bHTOP_*$, 
$B\bGL(A)^+ \approx \Omega B\bGL(\Sigma A)^+$.  
Consequently, $\forall \ n \geq 1$, 
$K_n(A) = \pi_n(B\bGL(A)^+) \approx$ 
$\pi_n(\Omega B\bE(\Sigma A)^+) \approx$ 
$\pi_{n+1}(B\bE(\Sigma A)^+) \approx$ 
$\pi_{n+1}(B\bGL(\Sigma A)^+) =$ 
$K_{n+1}(\Sigma A)$.  
It is also true that $K_0(A) \approx K_1(\Sigma A)$ 
(Farrell-Wagoner\footnote[2]{\textit{Comment. Math. Helv.} \textbf{47} (1972), 474-501.\vspace{0.11cm}}).  
Let $\Omega_0 B\bGL(\Sigma A)^+$ be the path component of $\Omega B\bGL(\Sigma A)^+$ 
containing the constant loop $-$then in $\bHTOP_*$, 
$\Omega B\bE(\Sigma A)^+ \approx \Omega_0 B\bGL(\Sigma A)^+$ (cf. p. \pageref{5.41}).  
But $\pi_1(B\bGL(\Sigma A)^+) = K_1(\Sigma A)$, hence 
$K_0(A) \times B\bGL(\Sigma A)^+ \approx \Omega B\bGL(\Sigma A)^+$.
\vspi
[Note: \  Additional information can be found in 
Wagoner\footnote[3]{\textit{Topology} \textbf{11} (1972), 349-370.}
There it is shown that by fixing the data, the pointed homotopy equivalence 
$K_0(A) \times B\bGL(A)^+ \approx \Omega B\bGL(\Sigma A)^+$ can be made natural, i.e., if 
$f:A^\prime \ra A\pp$ is a morphism of rings, then the diagram
\begin{tikzcd}%[sep=large]
{K_0(A^\prime) \times B\bGL(A^\prime)^+} \ar{d} \indent \approx &{\Omega B\bGL(\Sigma A^\prime)^+} \ar{d}\\
{K_0(A\pp) \times B\bGL(A\pp)^+} \indent \approx &{\Omega B\bGL(\Sigma A\pp)^+}
\end{tikzcd}
is pointed homotopy commutative.]\\
\endgroup%%------------------------------------<<
\vspace{0.2cm}

\begingroup%%----------------------------------->> 
\fontsize{9pt}{11pt}\selectfont
\textbf{\small EXAMPLE} \ 
Let \mA be a ring with unit $-$then 
$\Sigma \bU\bT(A) \approx \bU\bT(\Sigma A)$ $\implies$ 
$K_0(\bU\bT(A)) \approx $
$K_1(\Sigma \bU\bT(A))$ $\approx$ 
$K_1(\bU\bT(\Sigma A))$ $\approx$ 
$K_1(\Sigma A) \times K_1(\Sigma A)$ $\approx$ 
$K_0(A) \times K_0(A)$.\\
\endgroup%%------------------------------------<<
\vspace{0.2cm}

\index{Theorem Kan-Thurston}
\textbf{\small KAN$-$THURSTON THEOREM} \ 
Let $X$ be a pointed connected CW space $-$then there exists a group $G_X$ and an acyclic map 
$\kappa_X:K(G_X,1) \ra X$.

[Because of Proposition 2, one can take for $X$ a pointed connected CW complex with all characteristic maps embeddings.  
Moreover, it will be enough to deal with finite $X$, the transition to infinite $X$ being straightforward 
(given the naturality built into the argument).  
Since $\dim X \leq 1$ $\implies$ $X$ is aspherical, we shall assume that $\dim X > 1$ and proceed by induction on $\#(\sE)$, 
supposing that the construction has been carried out in such a way that if $X_0$ is a connected subcomplex of $X$, 
then $K(G_{X_0},1) = \kappa_X^{-1}(X_0)$ and $G_{X_0} \ra G_X$ is injective.  
To execute the inductive step, consider the pushout square
\begin{tikzcd}%[sep=large]
{\bS^{n-1}} \ar{d} \ar{r} &X \ar{d}\\
{\bD^{n}} \ar{r} &Y
\end{tikzcd}
$(n \geq 2)$, where the horizontal arrows are embeddings and 
$
\begin{cases}
\ X_0 = \text{im}(\bS^{n-1} \ra X) \\[-.1cm]
\ Y_0 = \text{im}(\bD^{n} \ra Y) 
\end{cases}
$
are connected
%%----------------------------------------------------------------------------------------------77
subcomplexes of 
$
\begin{cases}
\ X\\[-.1cm]
\ Y
\end{cases}
$
\hspace{-.26cm}, so
\begin{tikzcd}%[sep=large]
X_0 \ar{d} \ar{r} &X \ar{d}\\
Y_0 \ar{r} &Y
\end{tikzcd}
is a pushout square.  Recalling that there is a monomorphism $G_{X_0} \ra \Gamma G_{X_0}$ of groups 
(cf. p. \pageref{5.42}), 
define $G_Y$ by the pushout square
\begin{tikzcd}%[sep=large]
G_{X_0} \ar{d} \ar{r} &G_X \ar{d}\\
{\Gamma G_{X_0}} \ar{r} &G_Y
\end{tikzcd}
and realize $K(G_Y,1)$ by the pushout square
\begin{tikzcd}%[sep=large]
{K(G_{X_0},1)} \ar{d} \ar{r} &K(G_X,1) \ar{d}\\
{K(\Gamma G_{X_0},1)} \ar{r} &K(G_Y,1)
\end{tikzcd}
(cf. p. \pageref{5.43}).  
Extend $\kappa_X:K(G_X,1) \ra X$ to $\kappa_Y:K(G_Y,1) \ra Y$  in the obvious way (thus $\kappa_YK(\Gamma G_{X_0},1) \subset Y_0$ and the diagram
\begin{tikzcd}%[sep=large]
{K(G_{X},1)} \ar{d} \ar{r}{\kappa_X} &X \ar{d}\\
{K(G_{Y},1)}  \ar{r}[swap]{\kappa_Y} &Y
\end{tikzcd}
commutes).  
The induction hypothesis implies that $\kappa_X$ and $\kappa_{X_0}$ are acyclic.  
In addition, $K(\Gamma G_{X_0},1)$ is an acyclic space and $Y_0$ is contractible, hence 
$\restr{\kappa_Y}{K(\Gamma G_{X_0},1)}$ is acyclic (cf. Proposition 20).  
Therefore, by comparing Mayer-Vietoris sequences and applying the five lemma, it follows that $\kappa_Y$ is acyclic (cf Proposition 22).  
Finally, the condition on connected subcomplexes passes on to $Y$.]

[Note: \  Put $N = \ker \pi_1(\kappa_X)$ $-$then $X$ is a model for $K(G_X,1)_N^+$.]\\

Application:  Every nonempty path connected topological space has the homology of a $K(G,1)$.\\

\label{9.6}
\begingroup%%----------------------------------->> 
\fontsize{9pt}{11pt}\selectfont
\textbf{\small EXAMPLE} \ 
Suppose given two sequences 
$\pi_n$ $(n \geq 2)$ $\&$ $G_q$ $(q \geq 1)$ 
of abelian groups $-$then there exists a pointed connected CW space $Z$ such that 
$\forall \ n \geq 2$: $\pi_n(Z) \approx \pi_n$ $\&$ $\forall \ q \geq 1$: 
$H_q(Z) \approx G_q$.  Thus choose $X$: 
$\pi_{n+1}(X) \approx \pi_n$ $(n \geq 2)$ (homotopy system theorem) and put 
$Y = \ds\bigvee\limits_1^\infty M(G_q,q)$ 
(cf. p. \pageref{5.44}): 
$H_q(Y) \approx G_q$ $(q \geq 1)$.  Using Kan-Thurston, form 
$
\begin{cases}
\ \kappa_X:K(G_X,1) \ra X \\
\ \kappa_Y:K(G_Y,1) \ra Y
\end{cases}
$
and consider 
$Z =$ 
$E_{\kappa_X} \times K(G_Y,1)$, 
the mapping fiber of the arrow 
$K(G_X \times G_Y,1) = K(G_X,1) \times K(G_Y,1) \ra X$.  
Example:  If $G_q$ $(q \geq 1)$ is any sequence of abelian groups, then there exists a group $G$ such that 
$\forall \ q \geq 1$: $H_q(G) \approx G_q$.
\vspi
[Note: \  $Z$ also has the property that $\pi_1(Z)$ operates trivially on $\pi_n(Z)$ $\forall \ n \geq 2$.]\\
\endgroup%%------------------------------------<< 

The homotopy categories of algebraic topology are not complete (or cocomplete), 
a circumstance that precludes application of the representable functor theorem and the general adjoint functor theorem (or their duals).  
However, there is still a certain amount of structure.  
For instance, consider \bHTOP.  
It has products and the double mapping track furnishes weak pullbacks.  
Therefore \bHTOP is weakly complete, i.e., every diagram $\Delta:\bI \ra \bHTOP$ has a weak limit 
(meaning: ``existence without uniqueness'').  \bHTOP is also weakly cocomplete.  
In fact, \bHTOP has coproducts, while weak pushouts are furnished
%%----------------------------------------------------------------------------------------------78
by the double mapping cylinder.  
Example:  Let (\bX,\bff) be an object in $\bFIL(\bHTOP)$ $-$then tel(\bX,\bff) is a weak colimit of (\bX,\bff).

[Note:   The discussion of $\bHTOP_*$ is analogous.  
Example: Let $f:X \ra Y$ be a pointed continuous function, $C_f$ its pointed mapping cone $-$then 
$C_f$ is a weak cokernel of [f].]\\

\begingroup%%----------------------------------->> 
\fontsize{9pt}{11pt}\selectfont

\textbf{\small EXAMPLE} \ 
For each $n$, put 
$Y_n = \bS^3$ and let $[\psi_n]:Y_{n+1} \ra Y_n$ 
be the homotopy class of maps of degree 2 $-$then $Y = \lim Y_n$ does not exist in \bHTOP.  
To see this, assume the contrary, thus 
$\forall \ X$, $[X,Y] \approx \lim[X,Y_n]$, so, in particular, $Y$ must be 3-connected.  
Form the adjunction space 
$\bD^3 \sqcup_f \bS^2$, where $f:\bS^2 \ra \bS^2$ 
is skeletal of degree 3.  
Since $\dim (\bD^3 \sqcup_f \bS^2) \leq 3$, of necessity 
$[\bD^3 \sqcup_f \bS^2,Y] = *$.  But according to the Hopf classification theorem, 
$[\bD^3 \sqcup_f \bS^2,\bS^3] \approx$ 
$H^3(\bD^3 \sqcup_f \bS^2;\Z)$, which is $\Z/3\Z$, and in the limit, 
$[\bD^3 \sqcup_f \bS^2,Y] \approx \Z/3\Z$.\\
\endgroup%%------------------------------------<< 

\begingroup%%----------------------------------->> 
\fontsize{9pt}{11pt}\selectfont
\textbf{\small EXAMPLE} \ 
Working in $\bHTOP_*$, let $f:X \ra Y$ be a pointed Hurewicz fibration, where \mX and \mY are path connected.  
Suppose that $K = \ker [f]$ exists, say $[\kappa]:K \ra X$.  
If $\pi$ is the projection $E_f \ra X$, then $f \circx \pi \simeq 0$, so there exists a pointed continous function 
$\phi:E_f \ra K$ such that $\kappa \circx \phi \simeq \pi$ and by construction, 
$f \circx \kappa \simeq 0$, so there exists a pointed continuous function 
$\psi:K \ra E_f$ such that $\kappa \simeq \pi \circx \psi$.  
Thus 
$\kappa \circx \phi \circx \psi \simeq \kappa$ $\implies$ 
$\phi \circx \psi \simeq \id_K$, $[\kappa]$ being a monomorphism in $\bHTOP_*$.  
Take now 
$X = \bSO(3)$, 
$Y = \bSO(3)/\bSO(2)$, and let $f:X \ra Y$ be the canonical map $-$then 
$\pi_1(E_f) \approx$ 
$\Z$, 
$\pi_1(K) \approx$ 
$\Z/2\Z$ and $\Z/2\Z$ is not a direct summand of $\Z$.
\vspi
[Note: \ Similar examples show that cokernels do not exist in $\bHTOP_*$.]\\
\endgroup%%------------------------------------<< 

Let \bC be a category with products and weak pullbacks $-$then every diagram in \bC has a weak limit.  Any functor $F:\bC \ra \bSET$ that preserves products and weak pullbacks necessarily preserves weak limits.\\

\begin{proposition} \ %27
Let \bC be a category with products and weak pullbacks.  Assume: Ob\bC contains a set $\sU = \{U\}$ with the following properties.\\
%^
\indent \indent $(\sU_1)$ \quadx A morphism $f:X \ra Y$ is an isomorphism provided that 
$\forall \ U \in \sU$, the arrow $\Mor(Y,U) \ra \Mor(X,U)$ is bijective.\\
%^
\indent \indent $(\sU_2)$ \quadx Each object $(\bX,\bff)$ in \textbf{TOW(C)} has a weak limit $X_\infty$ such that 
$\forall \ U \in \sU$, the arrow 
$\colimx \Mor(X_n,U) \ra \Mor(X_\infty,U)$ is bijective.

Then a functor $F:\bC \ra \bSET$ is representable iff it preserves products and weak pullbacks.
\end{proposition}

[The condition is certainly necessary.  As for the sufficiency, introduce the comma category $\abs{*,F}$.  
Recall that an object of $\abs{*,F}$ is a pair $(x,X)$ $(x \in FX, X \in \Ob\bC)$, while a morphism 
$(x,X) \ra (y,Y)$ is an arrow $f:X \ra Y$ such that $(Ff)x = y$.  The assumptions imply that 
$\abs{*,F}$ has products and weak pullbacks, hence is weakly complete, and $F$ is
%%----------------------------------------------------------------------------------------------79
representable iff  $\abs{*,F}$  has an initial object.  Let $\sU_F$ be the subset of Ob$\abs{*,F}$ 
consisting of the pairs $(u,U)$ $(u \in FU, U \in \sU)$.

Claim:  $\forall \ (x,X) \in \text{Ob}\abs{*,F}$ 
$\exists \ (\bar{x},\overline{X}) \in \text{Ob}\abs{*,F}$ 
and a morphism 
$(\bar{x},\overline{X})  \ra (x,X)$ 
such that 
$\forall \ (u,U) \in \sU_F$ there is a unique morphism $(\bar{x},\overline{X})  \ra (u,U)$.

[Define an object (\bX,\bff) in \textbf{TOW}$(\abs{*,F})$ by setting $(x_0,X_0) = (x, X) \times \prod (u,U)$ and inductively choose $(x_{n+1},X_{n+1}) \ra (x_n,X_n)$ to equalize all pairs of morphisms
%\begin{tikzcd}[ sep=small]
%{(x_n, X_n)} \ar[r, shift left] \ar[r,shift right] &{(u,U)}
%\end{tikzcd}
$(x_n, X_n)$ $\rightrightarrows$ $(u,U)$ 
$((u,U) \in \sU_F)$.  
Any weak limit of (\bX,\bff) created via $\sU_2$ is a candidate for $(\bar{x},\overline{X})$.]

The existence of an initial object in $(\abs{*,F})$ is then a consequence of observing that for all $(x,X)$ $\&$ $(y,Y)$:
(i)  Every morphism $(\bar{x},\overline{X}) \ra (\bar{y},\overline{Y})$ is an isomorphism (apply the claim and $\sU_1$);
(ii)  There is at least one morphism $(\bar{x},\overline{X}) \ra (\bar{y},\overline{Y})$ (the composite
$\overline{(\bar{x},\overline{X}) \times (u,Y)} \ra (\bar{x},\overline{X}) \times (y,Y) \ra (\bar{x},\overline{X})$ is an isomorphism;
(iii)  There is at most one morphism $(\bar{x},\overline{X}) \ra (y,Y)$ (form the equalizer $(z,Z)$ of
\begin{tikzcd}[ sep=small]
{(\bar{x},\overline{X})} \ar[r, shift left] \ar[r,shift right] &{(y,Y)}
\end{tikzcd}
and consider the composite$(\bar{z},\overline{Z}) \ra (z,Z) \ra (\bar{x},\overline{X})$).]

\label{16.16}
[Note: \  Proposition 27 can also be formulated in terms of a category \bC that has coproducts and weak pushouts together with a set $\sU = \{U\}$ of objects satisfying the following conditions.\\
%^
\indent \indent $(\sU_1)$ \quadx A morphism $f:X \ra Y$ is an isomorphism provided that $\forall \ U \in \sU$, the arrow $\Mor(U,X) \ra \Mor(U,Y)$ is bijective.\\
%^
\indent \indent $(\sU_2)$ \quadx Each object $(\bX,\bff)$ in \textbf{FIL(C)} has a weak colimit $X_\infty$ such that $\forall \ U \in \sU$, the arrow $\colimx \Mor(U,X_n) \ra \Mor(U,X_\infty)$ is bijective.

Under these hypotheses, the conclusion is that a cofunctor $F:\bC \ra \bSET$ is representable iff it converts coproducts into products and weak pushouts into weak pullbacks.]\\

%% ------------------->

\label{15.38}
\begingroup%%----------------------------------->> 
\fontsize{9pt}{11pt}\selectfont
\textbf{\small EXAMPLE} \ 
Let \bC be a category with coproducts and weak pushouts whose representable cofunctors are precisely those that convert coproducts into products and weak pushouts into weak pullbacks.  
Suppose that $\bT = (T,m,\epsilon)$ is an idempotent triple in \bC and let $S \subset \Mor\bC$ be the class consisting of those $f$ such that $Tf$ is an isomorphism $-$then
(1) \ \mS admits a calculus of left fractions;
(2) \ \mS is saturated;
(3) \ \mS satisfies the solution set condition;
(4) \ \mS is coproduct closed, i.e., 
$s_i:X_i \ra Y_i$ in \mS $\forall \ i \in I$ $\implies$ 
$\ds\coprod\limits_i s_i: \ds\coprod\limits_i X_i \ra \ds\coprod\limits_i Y_i$ in \mS.  
Conversely, any class $S \subset \Mor\bC$ with properties $(1) - (4)$ is generated by an idempotent triple, thus $S^\perp$ is the object class of a reflective subcategory of \bC.
\vspi
[The functor $L_S :\bC \ra S^{-1}\bC$ preserves coproducts and weak pushouts.  
So, for fixed 
$Y \in \text{Ob}S^{-1}\bC$, $\Mor(L_S-,Y)$ is a cofunctor 
$\bC \ra \bSET$ 
which converts coproducts into products and weak pushouts into weak pullbacks., hence is representable: $\Mor(L_SX,Y) \approx \Mor(X,Y_S)$.  
Use the assignment $Y \ra Y_S$ to define a functor 
$S^{-1}\bC \ra \bC$ and take for $T$ the composite 
$\bC \ra S^{-1}\bC \ra \bC$.  
Let 
$\epsilon_X \in \Mor(X,TX)$ correspond to $\id_{L_SX}$ under the bijection 
$\Mor(L_SX,L_SX) \approx \Mor(X,TX)$ $-$then 
$\epsilon: \id_\bC \ra T$ is a natural transformation, 
$\epsilon T = T \epsilon$ is a natural isomorphism, and $Tf$ is an isomorphism iff $f \in S$.]\\
\endgroup%%------------------------------------<< 

%%----------------------------------------------------------------------------------------------80
Notation: $\bCONCW_*$ is the full subcategory of $\bCW_*$  whose objects are the pointed connected CW complexes and $\bHCONCW_*$ is the associated homotopy category.\\

\textbf{\small LEMMA} \ 
$\bHCONCW_*$ has coproducts and weak pushouts.

[If 
$X \overset{f}{\lla} Z \overset{g}{\lra} Y$
is a 2-source in $\bCONCW_*$, then using the skeletal approximation theorem, one can always arrange that $M_{f,g}$ remains in $\bCONCW_*$.]\\

\index{Brown Representability Theorem}
\textbf{\small BROWN REPRESENTABILITY THEOREM} \ 
A cofunctor $F:\bHCONCW_* \ra \bSET$ is representable iff it converts coproducts into products and weak pushouts into weak pullbacks.

[Take for $\sU$ the set $\{(\bS^n,s_n):n \in \N\}$ $-$then $\sU_1$ holds since in $\bCONCW_*$ 
a pointed continuous function $f:X \ra Y$ is a pointed homotopy equivalence iff it is a weak homotopy equivalence 
(cf. p. \pageref{5.44a}) 
and $\sU_2$ holds since one can take for a weak colimit of an object 
(\bX,\bff) in $\bFIL(\bHCONCW_*)$ the pointed mapping telescope constructed using pointed skeletal maps 
(cf. p. \pageref{5.44b}).]

[Note: \  Since $F$ converts coproducts into products, $F$ takes an initial object to a terminal object: $F* = *$ and $X \ra *$ $\implies$ $* = F* \ra FX$, thus $FX$ has a natural base point.]\\

Spelled out, here are the conditions on $F$ figuring in the Brown representability theorem.

\indent\indent (Wedge Condition) \quadx 
For any collection $\{X_i: i \in I\}$ in $\bCONCW_*$, 
$F(\bigvee\limits_i X_i) \approx$ 
$\prod\limits_i FX_i$.

\indent\indent (Mayer-Vietoris Condition) \quadx
For any weak pushout square
\begin{tikzcd}%[sep=large]
Z \ar{d}[swap]{f} \ar{r}{g} &Y \ar{d}{\eta}\\
X  \ar{r}[swap]{\xi} &P
\end{tikzcd}
in $\bHCONCW_*$,
\begin{tikzcd}%[sep=large]
FP \ar{d}[swap]{F\xi} \ar{r}{F\eta} &FY \ar{d}{Fg}\\
FX  \ar{r}[swap]{Ff} &FZ
\end{tikzcd}
is a weak pullback square in \bSET, so $\forall$
$
\begin{cases}
\ x \in FX\\[-.1cm]
\ y \in FY
\end{cases}
:
$
$(Ff)x = (Fg)y$, $\exists \ p \in FP:$
$
\begin{cases}
\ (F\xi)p = x\\[-.1cm]
\ (F\eta)p = y
\end{cases}
\hspace{-.25cm}.
$

[Note: \  It is not necessary to make the verification for an arbitrary weak pushout square.  In fact, it is sufficient to consider pointed double mapping cylinders calculated relative to skeletal maps, thus it is actually enough to consider diagrams of the form
$
\begin{tikzcd}%[sep=large]
C \ar{d} \ar{r} &B \ar{d}\\
A  \ar{r} &X
\end{tikzcd}
, 
$
where $X$ is a pointed connected CW complex and 
$
\begin{cases}
\ A\\[-.1cm]
\ B
\end{cases}
$
$\& \ C$ are pointed connected subcomplexes such that $X = A \cup B$, $C = A \cap B$.]

%%----------------------------------------------------------------------------------------------81
Examples: \  
(1) \
Fix a pointed path connected space $(X,x_0)$ $-$then $[-;X,x_0]$ is a cofunctor on $\bHCONCW_*$ satisfying the wedge and Mayer-Vietoris conditions, hence there exists a pointed connected CW complex $(K,k_0)$ and a natural isomorphism $\Xi: [-;K,k_0] \ra [-;X,x_0]$, each $f \in \Xi_{K,k_o}([\id_K])$ being a weak homotopy equivalence $K \ra X$, thus the Brown representability theorem implies the resolution theorem;
(2) \
Fix $n \in \N$ and an abelian group $\pi$ $-$then the cofunctor $H^n(-;\pi)$ (singular cohomology) satisfies the wedge and Mayer-Vietoris conditions, hence there exists a pointed connected CW complex $(K(\pi,n),k_{\pi,n})$ and a natural isomorphism $\Xi: [-;K(\pi,n),k_{\pi,n}] \ra$ 
$H^n(-;\pi)$, thus the Brown representability theorem implies the existence of Eilenberg-MacLane spaces of type $(\pi,n)$ ($\pi$ abelian);
(3) \
Fix a group $\pi$ $-$then the cofunctor that assigns to a pointed connected CW complex $(K,k_0)$ the set of homomorphisms $\pi_1(K,k_0) \ra \pi$ satisfies the wedge and Mayer-Vietoris conditions, hence there exists a pointed connected CW complex   $(K(\pi,1),k_{\pi,1})$ and a natural isomorphism 
$\Xi: [-;K(\pi,1),k_{\pi,1}] \ra$ 
$\text{Hom}(\pi_1-;\pi)$, thus the Brown representability theorem implies the existence of Eilenberg-MacLane spaces of type $(\pi,1)$ ($\pi$ arbitrary).

[Note: \   Both $\bHCW_*$ and \bHCW have coproducts and weak pushouts but Brown representability can fail.  
Indeed, Matveev
\footnote[2]{\textit{Math. Notes} \textbf{39} (1986), 471-474.\vspace{0.11cm}} 
has given an example of a nonrepresentable cofunctor $F:\bHCW_* \ra \bSET$ which converts coproducts into products and weak pushouts into weak pullbacks and 
Heller\footnote[3]{\textit{J. London Math. Soc.} \textbf{23} (1981), 551-562.} 
has given an example of a nonrepresentable cofunctor $F:\bHCW \ra \bSET$ 
which converts coproducts into products and weak pushouts into weak pullbacks.]\\

\label{17.11}
\begingroup%%----------------------------------->> 
\fontsize{9pt}{11pt}\selectfont
\textbf{\small EXAMPLE} \ 
Let $U:\bGR \ra \bSET$ be the forgetful functor.
\\
\indent\indent ($\bHCW_*$) \quadx Suppose that $F: \bHCW_* \ra \bGR$ is a cofunctor such that $U \circx F$ converts coproducts into products and weak pushouts into weak pullbacks $-$then $U \circx F$ is representable.
\vspi
[Represent the composite 
$\bHCONCW_* \ra \bHCW_* \ra \bGR  \ra \bSET$
by $K$.  Put $G = F\bS^0$ and equip it with the discrete topology.
\vspi
Claim:  For any $X$ in $\bCONCW_*$, $U \circx F(X_+) \approx [X_+,K \times G]$.
\vspi
[There is a split short exact sequence 
$1 \ra FX \ra FX_+ \ra F\bS^0 \ra 1$, hence $U \circx F(X_+) \approx U \circx F(X) \times
%%----------------------------------------------------------------------------------------------82
G$ $\approx$ $[X,K] \times G$ or, reinstating the base points: 
$U \circx F(X_+) \approx [X,x_0;K,k_0] \times G$.  And: $[X,x_0;K,k_0] \approx [X,K]$ $\implies$ $[X,x_0;K,k_0] \times G$ 
$\approx$ $[X,K] \times G$ $\approx$ $[X,K] \times [X,G]$ $\approx$ $[X,K \times G]$ $\approx$  $[X_+,K \times G]$.]
\vspi
Given $(X,x_0)$ in $\bCW_*$, let $X_{i_0}$, $X_i$ $(i \in I)$ be its set of path components, where 
$x_0 \in X_{i_0}$ $-$then 
$X = X_{i_0} \vee \ds\bigvee\limits_i X_{i+}$, so 
$U \circx F(X) \approx$ 
$U \circx F(X_{i_0}) \times \ds\prod\limits_i U \circx F(X_{i_+}) \approx$ 
$[X_{i_0},K] \times \ds\prod\limits_i[X_{i+},K \times G] \approx$ 
$[X_{i_0},K \times G] \times \ds\prod\limits_i[X_{i+},K \times G] \approx$ 
$[X,K \times G]$.]
\\
\indent\indent (\bHCW) \quadx Suppose that 
$F:\bHCW \ra \bGR$ is a cofunctor such that $U \circx F$ converts coproducts into products and weak pushouts into weak pullbacks $-$then $U \circx F$ is representable.
\vspi
[Let $F_*$ be the composite 
$\bHCONCW_* \ra$ 
$\bHCONCW \ra$ 
$\bHCW \ra$ 
$\bGR \ra$ 
$\bSET$.
\vspi
Claim: If $F* = *$, then $F_*$ is representable.
\vspi
[The assumption on $F$ implies that $FA = *$ for any discrete topological space $A$.  
To check that $F_*$ satisfies the wedge condition, put 
$X = \ds \coprod\limits_i X_i$ and let $A \subset X$ be the set made up of the base points 
$x_i \in X_i$ $-$then 
$F(X/A) \approx FX$.  
But 
$X/A = $
$\ds \bigvee\limits_i X_i$ $\implies$ 
$F_*\bigl(\bigvee\limits_i X_i\bigr) \approx$ 
$U \circx F(X) \approx$ 
$\ds\prod\limits_i F_*X_i$.  
As $F_*$ necessarily satisfies the Mayer-Vietoris condition, $F_*$ is representable: 
$[-,K_*] \approx$ $F_*$.]
\vspi
Claim: If $F* = *$, then $U \circx F$ is representable.
\vspi
[If $X$ is in \bCW and if 
$X = \ds \coprod\limits_i X_i$ 
is its decomposition into path components, then 
$U \circx F(X) \approx$ 
$\ds\prod\limits_i U \circx F(X_i) \approx$ 
$\ds\prod\limits_i F_* X_i \approx$ 
$\ds\prod\limits_i [X_i,K_*] \approx$ 
$\ds[\coprod\limits_i X_i,K_*] \approx$ 
$[X,K_*]$.]
\vspi
Given $X$ in \bCW, view $\pi_0(X)$ as a discrete topological space $-$then $U \circx F \circx \pi_0$ is represented by $F*$ (discrete topology).  
On the other hand, $F$ is the semidirect product of $F \circx \pi_0$ and the kernel $F_0$ of $F \ra F \circx \pi_0$ induced by the embedding $\pi_0(X) \ra X$.  Moreover, $U \circx F \approx U \circx F_0 \times U \circx F \circx \pi_0$ and $F_0* = *$ $\implies$ $U \circx F_0$ is representable.]\\
\endgroup%%------------------------------------<< 

\begingroup%%----------------------------------->> 
\fontsize{9pt}{11pt}\selectfont
Given a small, full subcategory $\bC_0$ of $\bHCW_*$, denote by $\overline{\bC_0}$ the full subcategory of $\bHCW_*$ whose objects are those $Y$ such that 
$g:Y \ra$ 
$Z$ is an isomorphism (= pointed homotopy equivalence) if 
$g_*:[X_0,Y] \ra$ 
$[X_0,Z]$ is bijective for all $X_0 \in \text{Ob}\bC_0$.\\
\endgroup%%------------------------------------<< 

\begingroup%%----------------------------------->> 
\fontsize{9pt}{11pt}\selectfont
\textbf{\small FACT} \ 
Suppose that $F: \bHCW_* \ra \bSET$ is a cofunctor which converts coproducts into products and weak pushouts into weak pullbacks 
$-$then there exists an object $X_F$ in $\bHCW_*$ and a natural transformation $\Xi:[-, X_F] \ra F$ 
such that $\forall \ X_0 \in \text{Ob}\bC_0$, $\Xi_{X_0}:[X_0;X_F] \ra FX_0$ is bijective.\\
\endgroup%%------------------------------------<< 

\begingroup%%----------------------------------->> 
\fontsize{9pt}{11pt}\selectfont
\textbf{\small FACT} \ 
Suppose that $F: \bHCW_* \ra \bSET$ is a cofunctor which converts coproducts into products and weak pushouts into weak pullbacks 
$-$then \mF is representable if for some $\bC_0$, $X_F \in \text{Ob}\overline{\bC}_0$.
\vspi
[With $\Xi$ as above, put $x_F = \Xi_{X_F}([\id_{X_F}])$, so that $\forall \ X \in \text{Ob}\bHCW_*$, $\Xi_X([f]) = F[f]_{X_F}$ $([f] \in [X,X_F])$.
\vspi
Surjectivity:  
Given $X \in \text{Ob}\bHCW_*$, call $\bC_0^\prime$ the full subcategory of $\bHCW_*$ obtained by adding $X$ and $X_F$ to $\bC_0$.  
Determine $X_F^\prime$ and $\Xi^\prime:[-,X_F^\prime] \ra F$ accordingly.  
In particular, $\Xi_{X_F}^\prime:[X_F,X_F^\prime] \ra FX_F$ is surjective, thus $\exists \ [f] \in [X_F,X_F^\prime]$: $x_F = F[f]_{x_F^\prime}$.  
From the definitions, $\forall \ X_0 \in \text{Ob}{\bC}_0$, $f_*:[X_0,X_F] \ra
%%----------------------------------------------------------------------------------------------83
[X_0,X_F^\prime]$ is bijective.  
Therefore $f$ is an isomorphism.  
Let $x \in FX$ and choose 
$[g] \in [X,X_F^\prime]$: $\Xi_X^\prime([g]) = x$ $-$then 
$\Xi_X([f^{-1}] \circx [g]) =$ 
$F([f^{-1}] \circx [g])x_F =$ 
$F[g](F[f^{-1}] x_F) =$ 
$F[g]x_F^\prime = x$.
\vspi
Injectivity: Given $X \in \Ob\bHCW_*$, let 
$u,v: X \ra X_F$ be a pair of morphisms: 
$\Xi_X([u]) =$ 
$\Xi_X([v])$, i.e.,  
$F[u]x_F =$ 
$F[v]x_F$.  
Fix a weak coequalizer 
$f:X_F \ra Z$ of $u, v$ and choose $z \in FZ$: 
$F[f]z =$ $x_F$.  
Since
$\Xi_Z:[Z,X_F] \ra FZ$ 
is surjective, 
$\exists \ g:Z \ra X_F$ such that $\Xi_Z([g]) = z$, hence $x_F = F[g \circx f]x_F$.  
From the definitions, $\forall \ X_0 \in \text{Ob}{\bC}_0$, $(g \circx f)_*:[X_0,X_F] \ra [X_0,X_F]$ is bijective.
Therefore $g  \circx f$ is an isomorphism.  
Finally, $f \circx u \simeq f \circx v$ $\implies$ $g \circx f \circx u \simeq g \circx f \circx v$ $\implies$ $u \simeq v$.]\\
\endgroup%%------------------------------------<< 

\begingroup%%----------------------------------->> 
\fontsize{9pt}{11pt}\selectfont
Application:   Let $\bC_0$ be the full subcategory of $\bHCW_*$ consisting of the $(\bS^n,s_n)$ $(n \geq 0)$, so $\overline{\bC_0} = \bHCONCW_*$ $-$then a cofunctor $F: \bHCW_* \ra \bSET$ which converts coproducts into products and weak pushouts into weak pullbacks is representable provided that $\#(F\bS^0) = 1$.
\vspi
[In fact, $\pi_0(X_F) = [\bS^0,X_F] = F\bS^0$, thus $X_F$ is connected.]\\
\endgroup%%------------------------------------<< 

\begingroup%%----------------------------------->> 
\fontsize{9pt}{11pt}\selectfont
\textbf{\small EXAMPLE} \ 
Fix a nonempty topological space $F$.  
Given a CW complex $B$, let $k_FB$ be the set 
$\text{Ob}\overline{\textbf{FIB}}_{B,F}$ where $\overline{\textbf{FIB}}_{B,F}$ is the skeleton of $\textbf{FIB}_{B,F}$ 
(cf. p. \pageref{5.45}) 
$-$then $k_F$ is a cofunctor 
$\bHCW \ra \bSET$ 
which converts coproducts into products and weak pushouts into weak pullbacks 
(cf. p. \pageref{5.46}).  
However, $k_F$ is not automatically representable since Brown representability can fail in  \bHCW.  
To get around this difficulty, one employs a subterfuge.  
Thus given a pointed CW complex $(B,b_0)$, let 
$\textbf{FIB}_{B,F_{;*}}$ be the category whose objects are the pairs $(p,i)$, where 
$p:X \ra B$ is a Hurewicz fibration such that 
$\forall \ b \in B$, $X_b$ has the homotopy type of $F$ and 
$i:F \ra p^{-1}(b_0)$ is a homotopy equivalence, and whose morphisms 
$(p,i) \ra (q,j)$ are the fiber homotopy classes 
$[f]:X \ra Y$ and the homotopy classes 
$[\phi]:F \ra F$ such that
$f_{b_0} \circx i \simeq j \circx \phi$.  
As in the unpointed case, 
$\textbf{FIB}_{B,F_{;*}}$ has a small skeleton and there is a cofunctor 
$k_{F_{;*}}:\bHCW_* \ra \bSET$ which converts coproducts into products and weak pushouts into weak pullbacks.  
Since $\#(k_{F_{;*}}\bS^0) = 1$, it follows from the above that $k_{F_{;*}}$ is representable: 
$[-; B_F,b_F] \approx$ $k_{F_{;*}}$, $(B_F,b_F)$ a pointed connected CW complex.  
If now $B$ is a CW complex, then the functor 
$\textbf{FIB}_{B,F} \ra$ 
$\bFIB_{B_+,F_{;*}}$  that assigns to $p:X \ra B$ the pair 
$(p \coprod c, \id_F)$ $(c:F \ra *)$ induces a bijection
$\Ob\overline{\bFIB}_{B,F} \ra$ $\Ob\overline{\bFIB}_{B_+,F_{;*}}$, so 
$k_FB \approx$ 
$k_{F_{;*}}B_+ \approx$ 
$[B_+,*;B_F,b_F] \approx$ 
$[B,B_F]$, i.e., $B_F$ represents $k_F$.
Example:  Take $F = K(\pi,n)$ ($\pi$ abelian) $-$then 
$B_F$ has the same pointed homotopy type as $K(\pi,n+1;\chi_\pi)$ 
(cf. p. \pageref{5.47}) 
$(K(\pi,n+1;\chi_\pi)$ is not necessarily a CW complex).\\
\endgroup%%------------------------------------<< 

\begingroup%%----------------------------------->> 
\fontsize{9pt}{11pt}\selectfont
Example \ 
Consider the Hurewicz fibration 
$p_1:\Theta \bS^n \ra \bS^n$ $(n \geq 2)$.  Let 
$i:\Omega \bS^n \ra \Omega\bS^n$ be the identity and 
$\iota:\Omega \bS^n \ra \Omega\bS^n$ 
the inversion $-$then the pairs $(p_1,i)$ and $(p_1,\iota)$ are not isomorphic in 
$\textbf{FIB}_{\bS^n,\Omega \bS^n;*}$.\\
\endgroup%%------------------------------------<< 

\begingroup%%----------------------------------->> 
\fontsize{9pt}{11pt}\selectfont
Let $G$ be a topological group $-$then in the notation of 
p. \pageref{5.48}, 
the restriction $\restr{k_G}{\bHCW}$ is a cofunctor $\bHCW \ra \bSET$ which converts coproducts into products and weak pushouts into weak pullbacks.  
To ensure that it is representable, one can introduce the pointed analog of $\bBUN_{B,G}$, say $\bBUN_{B,G;*}$ and proceed as above.  
The upshot is that the classifying space $B_G$ is now a CW complex but this need not be true of the universal space $X_G$.  
To clarify the situation, consider the pullback square
\begin{tikzcd}[sep=large]
X_G \ar{d} \ar{r} &X_G^\infty \ar{d}\\
B_G  \ar{r} &B_G^\infty
\end{tikzcd}
.  Since
%%----------------------------------------------------------------------------------------------84
for any CW complex \mB, $[B,B_G] \approx k_GB \approx [B,B_G^\infty]$, the arrow 
$B_G \ra B_G^\infty$ is a weak homotopy equivalence 
(cf. p. \pageref{5.49} ff.).  
Therefore the arrow $X_G \ra X_G^\infty$ is a weak homotopy equivalence, so $X_G$ is homotopically trivial ($X_G^\infty$ being contractible).\\
\endgroup%%------------------------------------<< 

\begingroup%%----------------------------------->> 
\fontsize{9pt}{11pt}\selectfont
\textbf{\small LEMMA} \ 
$X_G$ is contractible iff $G$ is a CW space.
\vspi
[Necessity:  For then $X_G$ is a CW space and because the fibers of the Hurewicz fibration 
$X_G \ra B_G$ are homeomorphic to $G$, it follows that $G$ is a CW space 
(cf. p. \pageref{5.50}).
\vspi
[Sufficiency:  Due to $\S 6$, Proposition 11, $X_G$ is a CW space.  
But a homotopically trivial CW space is contractible.]
\\
\endgroup%%------------------------------------<< 

\begingroup%%----------------------------------->> 
\fontsize{9pt}{11pt}\selectfont
Moral:  When $G$ is a CW space, $k_G$ can be represented by a CW complex (cf. $\S 4$, Proposition 35).
\vspi
[Note: \  Under these conditions, $B_G$ and $B_G^\infty$ have the same homotopy type (representing objects are isomorphic), thus $B_G^\infty$ is a CW space (see p. \pageref{5.51} for another argument).]\\
\endgroup%%------------------------------------<< 

Notation:  $\bFCONCW_*$
\index{$\bFCONCW_*$}
is the full subcategory of $\bCONCW_*$ whose objects are the pointed finite connected CW complexes and 
$\bHFCONCW_*$ 
\index{$\bHFCONCW_*$}
 is the associated homotopy category.
 
[Note: \  Any skeleton $\overline{\bHFCONCW}_*$ of $\bHFCONCW_*$  is countable 
(cf. p. \pageref{5.52}).]\\

A cofunctor $F:\bHFCONCW_* \ra \bSET$ is said to be 
\un{representable in the large}
\index{representable in the large} 
if there exists a pointed connected CW complex $X$ and a natural isomorphism $[-,X] \ra F$.

[Note: \  In this context, $[-,X] $ stands for the restriction to $\bHFCONCW_*$ of the representable cofunctor determined by $X$.  
Observe that in general it is meaningless to consider $FX$.]

Example:  The restriction to $\bHFCONCW_*$ of any cofunctor 
$\bHCONCW_* \ra \bSET$ satisfying the wedge and Mayer-Vietoris conditions is representable in the large.

Let $F:\bHCONCW_* \ra \bSET$ be a cofunctor.\\
\indent\indent (Finite Mayer-Vietoris Condition)  \ \ 
For any weak pushout square
\begin{tikzcd}%[sep=large]
Z \ar{d}[swap]{f} \ar{r}{g} &Y \ar{d}{\eta}\\
X  \ar{r}[swap]{\xi} &P
\end{tikzcd}
in $\bHCONCW_*$, where $Z$ is finite,
\begin{tikzcd}%[sep=large]
FP \ar{d}[swap]{F\xi} \ar{r}{F\eta} &FY \ar{d}{Fg}\\
FX  \ar{r}[swap]{Ff} &FZ
\end{tikzcd}
is a weak pullback square in \bSET, so $\forall$
$
\begin{cases}
\ x \in FX\\[-.1cm]
\ y \in FY
\end{cases}
$
\hspace{-.26cm}: $(Ff)x = (Fg)y$, $\exists \ p \in FP:$
$
\begin{cases}
\ (F\xi)p = x\\[-.1cm]
\ (F\eta)p = y
\end{cases}
\hspace{-.26cm}.
$
\\
\indent\indent (Limit Condition) \ \ 
For any pointed connected CW complex $X$ and for any collection $\{X_i: i \in I\}$ of pointed connected subcomplexes of $X$ such that $X = \text{colim}X_i$,
%%----------------------------------------------------------------------------------------------85
where $I$ is directed and the $X_i$ are ordered by inclusion, the arrow $FX \ra \lim FX_i$ is bijective.\\

\textbf{\small SUBLEMMA} \ 
Let $F:\bHCONCW_* \ra \bSET$  be a cofunctor satisfying the wedge and finite Mayer-Vietoris conditions.  
Fix an $X$ in $\bCONCW_*$ and choose $x \in FX$.  Suppose that 
$X  \overset{f}{\lla} K \overset{g}{\lra} X$
is a pointed 2-source, where $K$ is in $\bFCONCW_*$ and 
$
\begin{cases}
\ f\\[-.1cm]
\ g
\end{cases}
$
are skeletal with $(Ff)x = (Fg)x$ $-$then there is a $Y$ in $\bCONCW_*$ containing $X$ as an embedded pointed subcomplex, say $i:X \ra Y$, such that $i \circx f \simeq i \circx g$ and a $y \in FY$ such that $(Fi)y = x$.

[Consider the weak pushout square \ 
\begin{tikzcd}%[sep=large]
{K \vee K} \ar{d}[swap]{\nabla_K} \ar{r}{f \vee g} &X \ar{d}\\
K  \ar{r} &Y
\end{tikzcd}
, \ where $Y$ is the pointed double mapping cylinder of the folding map $\nabla_K$ and the wedge $f \vee g$.  
By construction, $Y$ is a pointed weak coequalizer of 
$
\begin{cases}
\ f\\[-.1cm]
\ g
\end{cases}
$
and the existence of $y \in FY$ follows from the assumptions.]\\

\textbf{\small LEMMA} \ 
Let $F:\bHFCONCW_* \ra \bSET$ be a cofunctor  satisfying the wedge, finite Mayer-Vietoris, and limit conditions.  
Fix an $X$ in $\bCONCW_*$ and choose $x \in FX$ $-$then there is a $Y$ in $\bCONCW_*$ containing $X$ as an embedded pointed subcomplex, say 
$i:X \ra Y$, such that $i \circx f \simeq i \circx g$ for any pointed 2-source
$X \overset{f}{\lla} K \overset{g}{\lra} X$, 
where $K$ is in $\bFCONCW_*$ and 
$
\begin{cases}
\ f\\[-.1cm]
\ g
\end{cases}
$
are skeletal with $(Ff)x = (Fg)x$ and a $y \in FY$ such that $(Fi)y = x$.

[Since it is enough to let $K$ run over the objects in $\overline{\bHFCONCW_*}$, one need only deal with a set
$\{X \overset{f_s}{\lla} K_s \overset{g_s}{\lra} X: s \in S\}$
of pointed 2-sources.  Given any $T \subset S$, proceed as in the proof of the sublemma and form the weak pushout square
$
\begin{tikzcd}%[sep=large]
{\bigvee\limits_t(K_t \vee K_t)} \ar{d} \ar{r} &X \ar{d}{i_T}\\
{\bigvee\limits_t K_t} \ar{r} &Y_T
\end{tikzcd}
, 
$
so for $T^\prime \subset T^{\prime\prime}$ there is a commutative triangle
$
\begin{tikzcd}[sep=small]
&X \ar{ldd} \ar{rdd}\\
\\
Y_{T^\prime} \ar{rr}[swap]{j} &&Y_{T^{\prime\prime}}
\end{tikzcd}
.  
$  
\ 
Consider the set $\sT$ of pairs $(T,y_T)$ $(y_T \in FY_T)$: $(Fi_T)y_T = x$.  Order $\sT$ by writing $(T^\prime,y_{T^\prime}) \leq (T^{\prime\prime},y_{T^{\prime\prime}})$ iff $T^\prime \subset T^{\prime\prime}$ and $(Fj)y_{T^{\prime\prime}} = y_{T^\prime}$ $-$then the limit condition implies that every chain in $\sT$ has an upper bound, thus $\sT$ has a maximal element $(T_0,y_{T_0})$ (Zorn).  Thanks to the sublemma, $T_0 = S$, therefore one can take $Y = Y_S$, $y = y_S$.]\\

\begin{proposition} \ %28
Let $F:\bHCONCW_* \ra \bSET$ be a cofunctor satisfying the
%%----------------------------------------------------------------------------------------------86
wedge, finite Mayer-Vietoris, and limit conditions $-$then the restriction of \mF to $\bHFCONCW_*$ is representable in the large.
\end{proposition}

[Put $X^0 = \bigvee\limits_{K,k} K$, 
where $K$ runs over the objects in $\overline{\bHFCONCW_*}$ and for each $K$, $k$ runs over $FK$.  
Using the wedge condition, choose $x^0 \in FX^0$ such that the associated natural transformation 
$\Xi^0:[-,X^0] \ra F$ has the property that $\Xi_K^0:[K,X^0] \ra FK$ is surjective for all $K$.  
Per the lemma, construct $X^0 \subset X^1$ $\&$ $x^1 \in FX^1$ 
and continue by induction to obtain an expanding sequence 
$X^0 \subset X^1 \subset \cdots $ of topological spaces and elements 
$x^0 \in FX^0, x^1 \in FX^1, \ldots$ such that $\forall \ n$, 
$X^n$ is a pointed connected CW complex containing $X^{n-1}$ as a pointed subcomplex and 
$x^n \ra x^{n-1}$ under $X^{n-1} \ra X^n$.  
Put $X = X^\infty$ $-$then $X$ is a pointed connected CW complex containing 
$X^n$ as a pointed subcomplex (cf. p. \pageref{5.53}).  
Let $x \in FX$ be the element corresponding to $\{x^n\}$ via the limit condition and let 
$\Xi:[-,X] \ra F$ be the associated natural transformation.  
That $\Xi_K$ is surjective for all $K$ is automatic.  
But $\Xi_K$ is also injective for all $K$: 
$\Xi_K([f]) =$ $\Xi_K([g])$,
 i.e., $(Ff)x = (Fg)x$ ($f$, $g$ skeletal) $\implies$ 
 $(Ff)x^n =$ $(Fg)x^n$ $(\exists \ n)$ $\implies$ 
 $i \circx f \simeq i \circx g$ $(i:X^n \ra X^{n+1})$.]
 \\

Given a cofunctor  $F:\bHFCONCW_* \ra \bSET$, for $X$ in $\bCONCW_*$, let $\overline{F}X =$ 
$\lim FX_k$, 
where $X_k$ runs over the pointed finite connected subcomplexes of $X$ ordered by inclusion $-$then 
$\overline{F}$ is the object function of a cofunctor $\bHCONCW_* \ra \bSET$ whose restriction to
 $\bHFCONCW_*$ is (naturally isomorphic) to $F$.  
On the basis of the definitions, $\overline{F}$ satisfies the limit condition.  
Moreover, $\overline{F}$ satisfies the wedge condition provided that $F$ converts finite coproducts into finite products so, 
in order to conclude that $F$ is representable in the large, 
it need only be shown that $\overline{F}$ satisfies the finite Mayer-Vietoris condition (cf. Proposition 28).  
Assume, therefore, that $F$ converts weak pushouts into weak pullbacks.  
Consider the diagram
\begin{tikzcd}%[sep=large]
C \ar{d} \ar{r} &B \ar{d}\\
A \ar{r} &X
\end{tikzcd}
\hspace{-.26cm}, where $X$ is a pointed connected CW complex and
$
\begin{cases}
\ A\\[-.1cm]
\ B
\end{cases}
$
$\&\ C$ are pointed connected subcomplexes such that $X = A \cup B$, $C = A \cap B$ with $C$ finite.  
To prove that 
\begin{tikzcd}%[sep=large]
{\overline{F}X} \ar{d} \ar{r} &{\overline{F}B} \ar{d}\\
{\overline{F}A}  \ar{r} &{\overline{F}C}
\end{tikzcd}
is a weak pullback square, let
$
\begin{cases}
\ K_i\\[-.1cm]
\ L_j
\end{cases}
$
run over the pointed finite connected subcomplexes of
$
\begin{cases}
\ A\\[-.1cm]
\ B
\end{cases}
$
which contain $C$ and using obvious notation, let
$
\begin{cases}
\ \bar{a} \in \overline{F}A\\[-.1cm]
\ \bar{b} \in \overline{F}B
\end{cases}
$
\hspace{-.26cm}: $\restr{\bar{a}}{C} = \restr{\bar{b}}{C}$ $-$then the question is whether there exists $\bar{x} \in \overline{F}$X:
$
\begin{cases}
\ \restr{\bar{x}}{A} = \bar{a}\\[-.1cm]
\ \restr{\bar{x}}{B} = \bar{b}
\end{cases}
\hspace{-.26cm}. \ 
$
For this, note first that 
$
\begin{cases}
\ \overline{F}A = \lim FK_i\\[-.1cm]
\ \overline{F}B = \lim FL_j
\end{cases}
$
and $\overline{F}X = \lim FX_{ij}$ $(X_{ij} = K_i \cup L_j)$.  Represent
%%----------------------------------------------------------------------------------------------87
$
\begin{cases}
\ \bar{a}\\[-.1cm]
\ \bar{b}
\end{cases}
$ by 
$
\begin{cases}
\ \{a_i\} \quadx (a_i \in FK_i)\\[-.1cm]
\ \{b_j\} \quadx (b_j \in FL_j)
\end{cases}
$
and let $S_{ij}$ be the set of $x_{ij} \in FX_{ij}:$
$
\begin{cases}
\ \restr{x_{ij}}{K_i} = a_i\\[-.1cm]
\ \restr{x_{ij}}{L_j} = b_j
\end{cases}
$
\hspace{-.26cm}.  
Since $S_{ij}$ is nonempty and $\lim S_{ij}$ is a subset of $\lim FX_{ij}$, 
it suffices to prove that $\lim S_{ij}$ is nonempty as any $\bar{x} \in \lim S_{ij}$ will work.  
However, this is a subtle point that has been resolved only by placing restrictions on the range of $F$.\\

\begingroup%%----------------------------------->> 
\fontsize{9pt}{11pt}\selectfont
\label{5.55b}
\textbf{\small EXAMPLE} \ 
Let $U:\bCPTHAUS \ra \bSET$ be the forgetful functor.  
Suppose that 
$F:\bHFCONCW_*$ $\ra$ $\bCPTHAUS$ is a cofunctor such that $U \circx F$ 
converts finite coproducts into finite products and weak pushouts into weak pullbacks $-$then 
$U \circx F$ is representable in the large.  
In fact, if $T_{ij}$ is the subspace of $FX_{ij}$ such that 
$UT_{ij} = S_{ij}$, then $T_{ij}$ 
is closed and $\lim T_{ij}$ is calculated over a cofiltered category, hence $\lim T_{ij}$ is a nonempty compact Haudorff space.  
But $U$ preserves limits, therefore $\lim S_{ij} = U(\lim T_{ij})$ is also nonempty.
\vspi
[Note: \  More is true: $\overline{U \circx F}$ satisfies the Mayer-Vietoris condition, hence is representable.  
Example: If $Y$ is a pointed connected CW complex whose homotopy groups are finite, then for every pointed finite connected CW complex $X$, $[X,Y]$ is finite 
(cf. p. \pageref{5.54}), thus is a compact Hausdorff space (discrete topology) and so $\overline{[-,Y]}$ is representable.]\\
\endgroup%%------------------------------------<< 

\index{Theorem Replication Theorem}
\textbf{\small REPLICATION THEOREM} \ 
Let $f:K \ra L$ be a pointed skeletal map, where
$
\begin{cases}
\ K\\[-.1cm]
\ L
\end{cases}
$
are in $\bFCONCW_*$ $-$then for any cofunctor $F:\bHFCONCW_* \ra \bSET$ which converts finite coproducts into finite products and weak pushouts into weak pullbacks, there is an exact sequence
\[
\cdots \ra F \Sigma L \ra F \Sigma K \ra FC_f \ra FL \ra FK
\]
in \bSET$_*$.

[Note: \  $F$ takes (abelian) cogroup objects to (abelian) group objects, so all the arrows to the left of $F \Sigma K$ are homomorphisms of groups.  In addition, $F \Sigma K$ operates to the left on $FC_f$ and the orbits are the fibers of the arrow $FC_f \ra FL$ 
(cf. p. \pageref{5.55}).]
\\[0.1cm]

Application:  There is an exact sequence
\[
F \Sigma K_i \times F \Sigma L_j  \ra F \Sigma C \ra FX_{ij} \ra FK_i \times FL_j
\]
in \bSET$_*$.

[The pointed mapping cone of the arrow $K_i \vee L_j \ra X_{ij}$ has the same pointed homotopy type as $\Sigma C$.]
\\[0.1cm]

Let $(I,\leq)$ be a nonempty directed set, \bI the associated filtered category.  
Suppose that 
$\Delta:\bI^\text{op} \ra \bSET$ is a diagram, 
where $\forall \ i \in \Ob \bI$, $\Delta_i \neq \emptyset$ and $\forall \ \delta \in \Mor \bI$, $\Delta \delta$ is surjective.
%%----------------------------------------------------------------------------------------------88
In $I$, write $i \sim j$ iff there exists a bijective map $f:\Delta_i \ra \Delta_j$ and a $k$ with 
$
\begin{cases}
\ i\\[-.1cm]
\ j
\end{cases}
$
\hspace{-.26cm} $\leq k$ such that the triangle
\begin{tikzcd}[sep=small]
&{\Delta_k} \ar{ldd} \ar{rdd}\\
\\
{\Delta_i}  \ar{rr}[swap]{f} &&{\Delta_j} 
\end{tikzcd}
commutes.\\
\vspace{0.5cm}

\textbf{\small LEMMA} \ 
If $\#(I/\sim) \leq \omega$, then $\lim \Delta$ is nonempty.
\\

\index{Adams Representability Theorem}
\textbf{\small ADAMS REPRESENTABILITY THEOREM} \ 
Let $U:\bGR \ra \bSET$ be the forgetful functor.  
Suppose that $F:\bHFCONCW_* \ra \bGR$ is a cofunctor such that $U \circx F$ converts finite coproducts into finite products and weak pushouts into weak pullbacks $-$then $U \circx F$ is representable in the large.

[The arrow $S_{i^\prime,j^\prime} \ra S_{ij}$ is surjective if 
$
\begin{cases}
\ K_i \subset K_{i^\prime}\\[-.1cm]
\ L_j \subset L_{j^\prime}
\end{cases}
$
\hspace{-.26cm}.  
This is because $F \Sigma C$ acts transitively to the left on
$
\begin{cases}
\ S_{i^\prime,j^\prime}\\[-.1cm]
\ S_{ij}
\end{cases}
$
and $S_{i^\prime,j^\prime} \ra S_{ij}$ is equivariant.  Claim: $\#(\{ij\}/\sim) \leq \omega$.  For one can check that $ij \sim i^\prime j^\prime$ iff $F \Sigma K_i \times F \Sigma L_j \ra F \Sigma C$ $\&$ $F \Sigma K_{i^\prime} \times F \Sigma L_{j^\prime} \ra F \Sigma C$ have the same image, of which there are at most a countable number of possibilities.  The lemma thus implies that $\lim S_{i j}$ is nonempty.]\\

Working in $\bCONCW_*$, two pointed continuous functions $f,g:X \ra Y$ are said to be \un{prehomotopic} if for any pointed finite connected CW complex $K$ and any pointed continuous function $\phi:K \ra X$, $f \circx \phi \simeq g \circx \phi$.  Homotopic maps are prehomotopic but the converse is false since, e.g., there are phantom maps that are not nullhomotopic (see below).

Notation: $\bPREHCONCW_*$ is the quotient category of $\bCONCW_*$ 
defined by the congruence of prehomotopy, $[X,Y]_\text{pre}$ being the set of morphisms from $X$ to $Y$.

If $F:\bHFCONCW_* \ra \bSET$ is a cofunctor, then $\overline{F}$ can be viewed as a cofunctor $\bPREHCONCW_* \ra \bSET$.  
Given $X$ in $\bCONCW_*$, there is a bijection 
Nat$([-,X]_\text{pre},\overline{F})$ $\ra$ $\overline{F}X$ (Yoneda).  
On the other hand, there is a bijection
Nat$([-,X],F) \ra \overline{F}X$, viz. 
$\Xi \ra \{\Xi_{X_k}([i_k])\}$, 
$i_k:X_k \ra X$ the inclusion.  
Example:  Take $F = [-,X]$, so $\overline{[X,X]} = \lim [X_k,X]$, and put $\iota_X = \{[i_k]\}$ $-$then $id_{[-,X]} \leftrightarrow \iota_X$.\\

\begin{proposition} \ %29
Let $Y$ be in $\bCONCW_*$.  
Assume: $\overline{[-,Y]}$ satisfies the finite Mayer-Vietoris condition $-$then for all $X$ in $\bCONCW_*$, the natural map 
$[X,Y]_\text{pre} \ra$ 
$\lim[X_k,Y]$ is bijective.
\end{proposition}

[Injectivity is immediate.  Turning to surjectivity, note that by definition 
$\lim [X_k,Y] = \overline{[X,Y]}$.  
Fix $x_0 \in \overline{[X,Y]}$ and let $y_0 = \iota_Y$ $(\in \overline{[Y,Y]})$.  Put $Z_0 = X \vee Y$ and write $z_0 = (x_0,y_0) \in\overline{[Z_0,Y]} \approx \overline{[X,Y]} \times \overline{[Y,Y]}$.  Imitating the argument used in the proof of Proposition
%%----------------------------------------------------------------------------------------------89
28, construct a $Z$ in $\bCONCW_*$ containing $Z_0$ as an embedded pointed subcomplex and an element $z \in \overline{[Z,Y]}$ which restricts to $z_0$ such that the associated natural transformation 
$[K,Z] \ra [K,Y]$ is a bijection for all $K$.  
Specialize and take $K = \bS^n$ $(n \in \N)$ to see that the inclusion $j:Y \ra Z$ is a pointed homotopy equivalence (realization theorem) and then compose the inclusion $i:X \ra Z$ with a homotopy inverse for $j$ to get a pointed continuous function $f_0:X \ra Y$ whose prehomotopy class is sent to $x_0$.]\\

\label{11.4} %dmc mnft
\label{11.8} %dmc mnft
\begingroup%%----------------------------------->> 
\fontsize{9pt}{11pt}\selectfont

\textbf{\small FACT} \ 
If $Y$ is a pointed connected CW complex whose homotopy groups are countable, then $\overline{[-,Y]}$ satisfies the finite Mayer-Vietoris condition.
\vspi
[Note: \  Under this assumption on $Y$, it follows that for all $X$ in $\bCONCW_*$, 
the natural map $[X,Y] \ra \lim [X_k,Y]$ is surjective (and even bijective provided that the homotopy groups of $Y$ are finite 
(cf. p. \pageref{5.55a} $\&$ p. \pageref{5.55b})).]\\
\endgroup%%------------------------------------<< 

\begin{proposition} \ %30
Suppose that $F:\bHFCONCW_* \ra \bSET$  
is a cofunctor which converts finite coproducts into finite products and weak pushouts into weak pullbacks.  
Assume: \ $\overline{F}$ \ satisfies \ the \ finite \ Mayer-Vietoris \ condition \ $-$then the cofunctor \ 
$\overline{F}$ : $\bPREHCONCW_* \ra \bSET$  is representable.
\end{proposition}

[By Proposition 28, there is an $X$ in $\bCONCW_*$ and a natural isomorphism  
$\Xi: [-,X] \ra F$.  Repeating the reasoning used in the proof of Proposition 29, one finds that the extension 
$\overline{\Xi}: [-,X]_\text{pre} \ra \overline{F}$ is a natural isomorphism as well.]\\

\begin{proposition} \ %31
Suppose that $F, F^\prime:\bHFCONCW_* \ra \bSET$ are cofunctors which convert finite coproducts into finite products and weak pushouts into weak pullbacks.  
Assume $\overline{F}$ and $\overline{F}^\prime$ satisfy the finite Mayer-Vietoris condition.  
Fix natural isomorphisms $\Xi: [-,X] \ra F$, $\Xi^\prime: [-,X^\prime] \ra F^\prime$, where $X$, $X^\prime$ are pointed connected CW complexes.  
Let $T:F \ra F^\prime$ be a natural transformation $-$then there is a pointed continuous function $f:X \ra X^\prime$, unique up to prehomotopy, such that the diagram
%\begin{tikzcd}[ sep=small]
%[K,X] \arrow{r}{f_*} \arrow{d}{\Xi_k} &[K,X^\prime] \arrow{d}{\Xi_K^\prime}\\
%FK \arrow{r}{T_K} &F^\prime K\\
%\end{tikzcd}
\begin{tikzcd}%[sep=large]
{[K,X]} \arrow{r}{f_*} \arrow{d}[swap]{\Xi_K} &{[K,X^\prime]} \arrow{d}{\Xi_K^\prime}\\
FK \arrow{r}[swap]{T_K} &F^\prime K
\end{tikzcd}
commutes for  all K.
\end{proposition}

[Note: \  If $F = F^\prime$ and $T$ is the identity, 
then $f:X \ra X^\prime$ is a pointed homotopy equivalence.]\\

\begin{proposition} \ %32
Any representing object in the Adams representability theorem is a group object in $\bPREHCONCW_*$ and all such have the same pointed homotopy type.\\
\end{proposition}

%%----------------------------------------------------------------------------------------------90
\begingroup%%----------------------------------->> 
\fontsize{9pt}{11pt}\selectfont
\textbf{\small FACT} \ 
Let $F:\bHFCONCW_* \ra \bSET$ be a cofunctor which converts finite coproducts into finite products and weak pushouts into weak pullbacks.  Assume: $\forall \ K$, $\#(FK) \leq \omega$ $-$then \mF is representable in the large.
\vspi
[Note: \  It is unknown whether the cardinality assumption can be dropped.]\\
\endgroup%%------------------------------------<<

Given pointed connected CW complexes
$
\begin{cases}
\ X \\[-.11cm]
\ Y
\end{cases}
\hspace{-.25cm},
$
a pointed continuous function $f:X \ra Y$ is said to be a 
\un{phantom map}
\index{phantom map} 
if it is prehomotopic to 0.  Let $\Ph(X,Y)$ be the set of pointed homotopy classes of phantom maps from $X$ to $Y$ $-$then there is an exact sequence
\[
* \ra \Ph(X,Y) \ra [X,Y] \ra \lim[X_k,Y]
\]
in $\bSET_*$.  
Of course $[0] \in \Ph(X,Y)$ but $\#(\Ph(X,Y)) > 1$ is perfectly possible.  
Example:  Take $X = K(\Q,3)$, $Y = K(\Z,4)$ $(\implies$ 
$[X,Y] \approx$ 
$H^4(\Q,3) \approx$ 
$\text{Ext}(\Q,\Z) \approx \R)$, realize $X$ as the pointed mapping telescope of the sequence 
$\bS^3 \ra \bS^3 \ra \cdots ,$ the $k^{\thx}$ map having degree $k$, and note that up to homotopy, every 
$\phi:K \ra X$ factors through $\bS^3$ $(\implies$ $\Ph(X,Y) =$ $[X,Y])$.\\

\begingroup%%----------------------------------->> 
\fontsize{9pt}{11pt}\selectfont
Is the arrow $[X,Y] \ra \lim[X_k,Y]$ always surjective?  While the answer is ``yes'' under various assumptions on $X$ or $Y$, what happens in general has yet to be decided.
\vspi
[Note: \   
By contrast, there is a bijection $\Ph(X,Y) \ra \lim^1[\Sigma X_k,Y]$ of pointed sets 
(Gray-McGibbon\footnote[2]{\textit{Topology} \textbf{32} (1993), 371-394.\vspace{0.11cm}}).]\\
\endgroup%%------------------------------------<<

\begingroup%%----------------------------------->> 
\fontsize{9pt}{11pt}\selectfont
\textbf{\small EXAMPLE} \ 
Meier\footnote[3]{\textit{Quart. J. Math.} \textbf{29} (1978), 469-481.\vspace{0.11cm}}  
has shown that $\Ph(K(\Z,n),\bS^{n+1}) \approx \text{Ext}(\Q,\Z)$ for all positive even $n$.  
Special case:  $\Ph(\bP^\infty(\C),\bS^3) \approx \text{Ext}(\Q,\Z)$.
\vspi
[Note: \  Suppose that \mG is an abelian group which is countable and torsion free $-$ then 
$\exists \ X$ $\&$ $Y: \Ph(X,Y) \approx \text{Ext}(G,\Z)$ 
(Roitberg\footnote[6]{\textit{Topology Appl.} \textbf{59} (1994), 261-271.}).]\\ 
\endgroup%%------------------------------------<<

\label{9.22}
\label{15.24}
\begingroup%%----------------------------------->> 
\fontsize{9pt}{11pt}\selectfont
\index{universal phantom maps}
\textbf{\small EXAMPLE \ (\un{Universal Phantom Maps})} \ 
Let $X$ be a pointed connected CW complex.  
Assume: $X$ has a finite number of cells in each dimension 
$-$then it is clear that $f:X \ra Y$ is a phantom map iff $\forall \ n > 0$, $\restr{f}{X^{(n)}}$ is nullhomotopic.  
Denote by $\tel^+ X$ the pointed telescope of $X$ which starts at $X^{(1)}$ rather than $X^{(0)}$.  
Recall that the projection $p:\text{tel}^+X \ra X$ is a pointed homotopy equivalence 
(cf. p. \pageref{5.55c}).  
Now collapse each integral joint of $\text{tel}^+X$ to a point, i.e., mod out by 
$\ds\bigvee\limits_{n > 0}  X^{(n)}$.  The resulting quotient can be identified with 
$\ds\bigvee\limits_{n > 0} \Sigma X^{(n)}$ and the arrow 
$\Theta:\text{tel}^+X \ra$ 
$\ds\bigvee\limits_{n > 0} \Sigma X^{(n)}$ is a phantom map.  It is universal
%%----------------------------------------------------------------------------------------------91
in the sense that if 
$f:X \ra Y$ is a phantom map and if 
$\ov{f} = f \circx p$, then there is a pointed continuous function 
$F:\ds\bigvee\limits_{n > 0} \Sigma X^{(n)} \ra Y$ such that $\ov{f} \simeq F \circx \Theta$.  
This is because the inclusion $i: \ds\bigvee\limits_{n > 0}  X^{(n)} \ra \telsub^+X$ is a closed cofibration, hence 
$C_i \approx$ 
$\ds\bigvee\limits_{n > 0} \Sigma X^{(n)}$ 
(cf. p. \pageref{5.55d}).  
Corollary:  All phantom maps out of $X$ are nullhomotopic iff $\Theta$ is nullhomotopic.
\vspi
[Note: \  Here is an application.  Suppose that $
\begin{cases}
\ X \\
\ Y
\end{cases}
$
are pointed connected CW complexes with a finite number of cells in each dimension.  Claim:  If 
$f:X \ra Y$ and 
$g:Y \ra Z$ are phantom maps, then 
$g \circx f: X \ra Z$ is nullhomotopic.  To see this, observe that the composite 
$\ds\bigvee\limits_{n > 0} \Sigma X^{(n)} \overset{F}{\lra}$ 
$Y \overset{p^{-1}}{\lra}$ 
$\telsub^+Y \overset{\Theta}{\lra}$ 
$\ds\bigvee\limits_{n > 0} \Sigma Y^{(n)}$ 
%\begin{tikzcd}
%\bigvee\limits_{n > 0} \Sigma X^{(n)}  \arrow{r}{F} &Y \arrow{r}{p^{-1}} &\text{tel}^+Y \arrow{r}{\Theta} &\bigvee\limits_{n > 0} \Sigma Y^{(n)}
%\end{tikzcd}
is a phantom map.  
Accordingly, its restriction to each 
$\Sigma X^{(n)}$ is nullhomotopic, so actually 
$\Theta \circx p^{-1} \circx F \simeq 0$.  Therefore 
$g \circx f \simeq (\ov{g} \circx p^{-1}) \circx (\ov{f} \circx p^{-1}) \simeq$ 
$(G \circx \Theta \circx p^{-1}) \circx (F \circx \Theta \circx p^{-1}) \simeq$ 
$G \circx (\Theta \circx p^{-1} \circx F) \circx \Theta \circx p^{-1} \simeq 0$.]\\
\endgroup%%------------------------------------<<

%%%%%%%%%%%%%%%%%%%%%%%%%%%%%%%%%%%%%%
%%%%%%%%%%%%%%%%%%%%%%%%%%%%%%%%%%%%%%

\begin{center}
$\S \ 5$
\\[0.5cm]
$\mathcal{REFERENCES}$\\
\end{center}

\[
\textbf{BOOKS}
\]

\begingroup
\fontsize{9pt}{11pt}\selectfont
\setlength\parindent{0 cm}

[1] \quad Baues, H., \textit{Obstruction Theory}, Springer Verlag (1977).
\\[-.2cm]

[2] \quad Berrick, A., \textit{An Approach to Algebraic K-Theory}, Pitman (1982).
\\[-.2cm]

[3] \quad Brown, K., \textit{Cohomology of Groups}, Springer Verlag (1982).
\\[-.2cm]

[4] \quad Cohen, M., \textit{A Course in Simple Homotopy Theory}, Springer Verlag (1973).
\\[-.2cm]

[5] \quad Cooke, G. and Finney, R., \textit{Homology of Cell Complexes}, Princeton University Press (1967).
\\[-.2cm]

[6] \quad tom Dieck, T., \textit{Topologie}, Walter de Gruyter (1991).
\\[-.2cm]
 
[7] \quad Dold, A., \textit{Halbezakte Homotopiefunktoren}, Springer Verlag (1966).
\\[-.2cm]

[8] \quad Fritsch, R. and Piccinini R., \textit{Cellular Structures in Topology}, Cambridge University Press (1990).
\\[-.2cm] 

[9] \quad Gray, B., \textit{Homotopy Theory}, Academic Press (1975).
\\[-.2cm] 

[10] \quad Hilton, P., \textit{Nilpotente Gruppen und Nilpotente R\"aume}, Springer Verlag (1984).
\\[-.2cm]

[11] \quad Hilton, P., Mislin, G., and Roitberg, J., \textit{Localization of Nilpotent Groups and Spaces}, North Holland 

\hspace{0.95cm}(1975).
\\[-.2cm] 

[12] \quad Hughes, B. and Ranicki, A., \textit{Ends of Complexes}, Cambridge University Press (1996).
\\[-.2cm] 

[13] \quad Lundell, A. and Weingram, S., \textit{The Topology of CW Complexes}, Van Nostrand (1969).
\\[-.2cm] 

%Postnikov M. Lectures in Algebraic Topology,  Homotopy Theory of Cell Spaces, Nauka 1985\\[-.2cm]
%[14]{\cyr Postniknov M., Lektsii po Algebraicheskoy Topologii}[{\cyr Teoriya Gomotopiy Kletochnykh }{\cyr Prostranstv]} (1985).\\[-.2cm]
[14] \quad {\fontencoding{OT2}\selectfont
Postnikov M., Lektsii po Algebraicheskoy Topologii, [Teoriya Gomotopi\u i Kletochnykh Pros-}

\hspace{0.875cm} {\fontencoding{OT2}\selectfont transtv], Nauka (1985)}.
\\[-.2cm]

[15] \quad Schubert, H., \textit{Topology}, Allyn and Bacon (1968).
\\[-.2cm]

[16] \quad Varadarajan, K., \textit{The Finiteness Obstruction of C.T.C. Wall}, John Wiley (1989).
\\[-.2cm] 

[17] \quad Whitehead, G., \textit{Elements of Homotopy Theory}, Springer Verlag (1978).
\\
\endgroup

\[
\textbf{ARTICLES}
\]

\begingroup
\fontsize{9pt}{11pt}\selectfont
\setlength\parindent{0 cm}

[1] \quad Adams, J., A Variant of E.H. Brown's Representability Theorem, 
\textit{Topology} \textbf{10} (1971), 185-198.
\\[-.2cm]

[2] \quad Baumslag, G., Dyer, E., and Heller, A., The Topology of Discrete Groups, 
\textit{J. Pure Appl. Algebra} \textbf{16} 

\hspace{0.8cm}(1980), 1-47.
\\[-.2cm]

[3] \quad Brown, E., Abstract Homotopy Theory, 
\textit{Trans. Amer. Math. Soc.} \textbf{119} (1965), 79-85.
\\[-.2cm]

[4] \quad Cauty, R., Sur les Sous-Espaces des Complexes Simpliciaux, 
\textit{Bull. Soc. Math. France} \textbf{100} (1972), 

\hspace{0.8cm}129-155.
\\[-.2cm]

[5] \quad Dror, E., A Generalization of the Whitehead Theorem, 
\textit{SLN} \textbf{249} (1971), 13-22.
\\[-.2cm]

[6] \quad Hausmann, J. and Husemoller, D., Acylic Maps, 
\textit{Enseign. Math.} \textbf{25} (1979), 53-75.
\\[-.2cm]

[7] \quad Kan, D. and Thurston, W., Every Connected Space has the Homology of a $K(\pi,1)$, 
\textit{Topology} \textbf{15} 

\hspace{0.8cm}(1976), 253-258.
\\[-.2cm]

[8] \quad Mather, M., Pullbacks in Homotopy Theory, 
\textit{Canad. J. Math.} \textbf{28} (1976), 225-263.
\\[-.2cm]

[9] \quad McGibbon, C., Phantom Maps, In: \textit{Handbook of Algebraic Topology}, I. James (ed.), North Holland 

\hspace{0.8cm}(1995), 1209-1257.
\\[-.2cm]

[10] \quad  Mislin, G., Wall's Finiteness Obstruction, In: \textit{Handbook of Algebraic Topology}, I. James (ed.), North 

\hspace{0.8cm} Holland (1995), 1259-1291.
\\[-.2cm]

[11] \quad Palais, R., Homotopy Theory of Infinite Dimensional Manifolds, 
\textit{Topology} \textbf{5} (1966), 1-16.
\\[-.2cm]

[12] \quad Robinson, C., Moore-Postnikov Systems for Nonsimple Fibrations, 
\textit{Illinois J. Math.} \textbf{16} (1972), 234-

\hspace{0.95cm}242.
\\[-.2cm]

[13] \quad Roitberg, J., Computing Homotopy Classes of Phantom Maps, 
\textit{CRM Proc.} \textbf{6} (1994), 141-168.
\\[-.2cm]

[14] \quad Sklyarenko, E., Some Applications of the Functor $\underset{\lla}{lim}^1$, 
\textit{Math. Sbornik} \textbf{51} (1985), 367-387.
\\[-.2cm]

[15] \quad Steenrod, N., Cohomology Operations and Obstructions to Extending Continuous Functions, 
\textit{Adv.} 

\hspace{0.95cm}\textit{Math.} \textbf{8} (1972), 371-416.
\\[-.2cm]

[16] \quad Whitehead, J., Combinatorial Homotopy I and II, 
\textit{Bull. Amer. Math. Soc.} \textbf{55} (1949), 213-245 and 

\hspace{0.95cm}453-496.
\\[-.2cm]

[17] \quad Whitehead, J., A Certain Exact Sequence, 
\textit{Ann. of Math.} \textbf{52} (1950), 51-110.

\setlength\parindent{2em}

\endgroup

\chapter{
$\boldsymbol{\S}$\textbf{6}.\quadx  ABSOLUTE NEIGHBORHOOD RETRACTS}
\setlength\parindent{2em}
\setcounter{proposition}{0}
\setcounter{chapter}{6}

%%----------------------------------------------------------------------------------------------01
$\text{ }$\\[-1.25cm]

From the point of view of homotopy theory, the central result of this $\S$ is the CW-ANR theorem which says that a topological space has the homotopy type of a CW complex iff it has the homotopy type of an ANR.  
But absolute neighborhood retracts also have a life of their own.  For example, their theory is an essential component of infinite dimensional topology.

Consider a pair $(X,A)$, i.e., a topological space $X$ and a subspace $A \subset X$.  
Let $Y$ be a topological space.  
Suppose given a continuous function $f:A \ra Y$ $-$then the extension question is:  
Does there exist a continuous function $F:X \ra Y$ such that $\restr{F}{A} = f$?  
While this is a complex multifaceted issue, there is an evident connection with the theory of retracts.  
For if we take $Y = A$, then the existence of a continuous extension $r:X \ra A$ of the identity map id$_A$ amounts to saying that $A$ is a retract of $X$.  
Every retract of a Hausdorff space $X$ is necessarily closed in \mX.  
On the other hand, if $A$ is closed in $X$, then with no assumptions on $X$, a continuous function $f:A \ra Y$ has a continuous extension $F:X \ra Y$ iff $Y$ is a retract of the adjunction space $X \sqcup_f Y$.  
The opposite end of the spectrum is when $A$ is dense in $X$.  In this case, one can be quite specific and we shall start with it.

Let $(X,A)$ be a pair with $A$ dense in $X$.  
Write $\tau_X$ and $\tau_A$ for the corresponding topologies.
Define a map $\text{ex}:\tau_A \ra \tau_X$ by $\text{ex}(O) = X - \overline{A - O}$, the bar denoting closure in 
$X$ $-$then $\text{ex}(O) \cap A = O$ and $\text{ex}(O) = \bigcup \{U: U \in \tau_X \ \& \  U \cap A = O\}$.  Obviously, 
$
\begin{cases}
\ \text{ex}(\emptyset) = \emptyset\\
\ \text{ex}(A)  = X
\end{cases}
$
and $\forall \ O, P \in \tau_A$: $\text{ex}(O \cap P) = \text{ex}(O) \cap \text{ex}(P)$.  
Put $\text{ex}(\sO) = \{\text{ex}(O): O \in \sO\}$ $(\sO  \subset \tau_A)$.\\

\begin{proposition} \ %01
Let $A$ be a dense subspace of a topological space $X$; 
let $Y$ be a regular Hausdorff space $-$then a given $f \in C(A,Y)$ admits a continuous extension 
$F \in C(X,Y)$ iff $X = \bigcup \text{ex}(f^{-1}(\sV))$ for every open covering $\sV$ of $Y$.
\end{proposition}

[The condition is clearly necessary.  As for the sufficiency, suppose that $X \neq \emptyset$ and $\#(Y) > 1$.  
Call
$
\begin{cases}
\ \tau_X\\
\ \tau_Y
\end{cases}
$
the topologies on $X$ and $Y$.\\
\indent\indent $(F^*)$  Define a map $F^*: \tau_Y \ra \tau_X$ by
\[
F^*(V) = \bigcup \left\{\text{ex}(f^{-1}(V^\prime)): V^\prime \in \tau_Y \ \& \ \overline{V^\prime} \subset V\right\}.
\]
Note that 
$
\begin{cases}
\ F^*(\emptyset) = \emptyset\\
\ F^*(Y) = X
\end{cases}
\text{and } \forall \ V_1, V_2 \in \tau_Y: F^*(V_1 \cap V_2) \hsx = \hsx F^*(V_1) \cap F^*(V_2).
$
Let $\{V_j\} \subset \tau_Y$ $-$then $F^*(\bigcup\limits_j V_j) \supset \bigcup\limits_j F^*(V_j)$ and in fact equality prevails.  
To see this, write 
%%----------------------------------------------------------------------------------------------02
$\bigcup\limits_j V_j = \cup \sV$, where $\sV$ is the set of all $V \in \tau_Y$: $\overline{V} \subset V_j$ $(\exists \ j)$.
Take a 
$V^\prime \in \tau_Y$: $\overline{V^\prime} \subset \bigcup\limits_j V_j$.  
Since $Y = (Y - \overline{V^\prime}) \bigcup (\cup \sV)$, 
$X = \text{ex}(f^{-1}(Y - \overline{V^\prime} )) \bigcup (\cup \text{ex}(f^{-1}(\sV)))$.  
But $\emptyset = \text{ex}(f^{-1}(V^\prime)) \cap \text{ex}(f^{-1}(Y - \overline{V^\prime} ))$ $\implies$ 
$\text{ex}(f^{-1}(V^\prime))$ $\subset$ $\text{ex}(f^{-1}(\sV))$ $\subset$ 
$\bigcup\limits_j F^*(V_j)$, from which it follows that 
$F^*(\bigcup\limits_j V_j) \subset \bigcup\limits_j F^*(V_j)$.

\indent\indent $(F_*)$ Define a map $F_*: \tau_X \ra \tau_Y$ by
\[
F_*(U) = \bigcup \{V: V \in \tau_Y \ \& \ F^*(V) \subset U\}.
\]
Note that $\forall \ U \in \tau_X$ and $\forall \ V \in \tau_Y$: $V \subset F_*(U)$ $\iff$ $F^*(V) \subset U$.  Indeed, $F^*$ respects arbitrary unions.  
We claim now that $\forall \ x \in X$ $\exists$ $y \in Y$ : $F_*(X - \overline{\{x\}}) = Y - \{y\}$.
Let $F_*(X - \overline{\{x\}}) = Y - B_x$.  
Case 1: $B_x = \emptyset$.  Here, $X = F^*(Y) \subset X - \overline{\{x\}}$, an impossibility.  
Case 2: $\#(B_x) > 1$.  Choose $y_1, y_2 \in B_x$ : $y_1 \neq y_2$.  Choose $V_1, V_2 \in \tau_Y$: $V_1 \cap V_2 = \emptyset$ $\&$ 
$
\begin{cases}
\ y_1 \in V_1\\
\ y_2 \in V_2
\end{cases}
$
$-$then $V_1 \cap V_2 \subset F_*(X - \overline{\{x\}})$ $\implies$ $F^*(V_1 \cap V_2) \subset X - \overline{\{x\}}$, i.e., 
$F^*(V_1) \cap F^*(V_2)  \subset X - \overline{\{x\}}$, thus either $F^*(V_1)$ or $F^*(V_2)$ is contained in $X - \overline{\{x\}}$ and so either $V_1$ or $V_2$ is contained in $F_*(X - \overline{\{x\}}) = Y - B_x$, a contradiction.

\indent\indent $(F)$ Define a map $F:X \ra Y$ by stipulating that $F(x) = y$ iff $F_*(X - \overline{\{x\}}) = Y - \{y\}$.  The definitions imply that
$
\begin{cases}
\ F^{-1}(V) = F^*(V)\\
\ F^{-1}(V) \cap A = f^{-1}(V)
\end{cases}
(V \in \tau_Y), \ 
$
therefore $F \in C(X,Y)$ and $\restr{F}{A} = f$.]\\
%\vspace{0.25cm}

Retain the assumption that $A$ is dense in $X$ and $Y$ is regular Hausdorff.  
Assign to each $x \in X$ the collection $\sU(x)$ of all its neighborhoods $-$then a continuous function $f:A \ra Y$ has a continuous extension  $F:X \ra Y$ iff $\forall \ x$ the filter base $f(\sU(x) \cap A)$ converges.
The nontrivial part of this assertion is a simple consequence of the preceding result.  
For suppose that for some open covering $\sV$ of $Y :$ $X \neq \bigcup \text{ex}(f^{-1}(\sV))$.  
Choose $x \in X :$ $x \notin \bigcup \text{ex}(f^{-1}(\sV))$, so $\forall \ U \in \sU(x)$  and $\forall \ V \in \sV :$ $U \cap A \not\subset f^{-1}(V)$ or still, $f(U \cap A) \not\subset V$.  
But $f(\sU(x) \cap A)$ converges to $y \in Y$.  Accordingly, there is
(i) \ $V_0 \in \sV :$ $y \in V_0$ and
(ii) \ $U_0 \in \sU(x) :$ $f(U_0 \cap A) \subset V_0$.  Contradiction.

Here are two applications.

\indent\indent (C) \ 
Suppose that $Y$ is compact Hausdorff $-$then a continuous function $f:A \ra Y$ has a continuous extension  $F:X \ra Y$ iff for every finite open covering $\sV$ of $Y$ there exists a finite open covering $\sU$ of $X$ such that $\sU \cap A$ is a refinement of $f^{-1}(\sV)$.

In this statement, one can replace ``compact'' by ``Lindel\"of'' if ``finite'' is replaced by ``countable''.  
More is true:  It suffices to assume that $Y$ is merely $\R$-compact (recall that every Lindel\"of regular Hausdorff space is 
$\R$-compact).

%%----------------------------------------------------------------------------------------------03
\indent\indent (R-C) \ 
Suppose that $Y$ is $\R$-compact $-$then a continuous function $f:A \ra Y$ has a continuous extension  $F:X \ra Y$ iff for every countable open covering $\sV$ of $Y$ there exists a countable open covering $\sU$ of $X$ such that $\sU \cap A$ is a refinement of $f^{-1}(\sV)$.

[There is a closed embedding $Y \ra \prod \R$.  
Postcompose $f$ with a generic projection $\prod\R \ra \R$ and extend it to $X$.  Form the associated diagonal map $F:X \ra \prod \R$ $-$then $F$ is continuous and $\restr{F}{A} = f$ (viewed as a map $A \ra \prod \R$).  
Conclude by remarking that $F(X) = F(\overline{A})$ $\subset$ $\overline{F(A)} \subset$ $\overline{Y} = Y$.]

[Note: \  The $\R$-compactness of $Y$ is essential.  Consider $X = [0,\Omega]$, $A = Y = [0,\Omega[$, and let $f = \text{id}_A$ ($Y$ is not $\R$-compact, being countably compact but not compact).]\\

\begingroup%%------------------------------------>>
\fontsize{9pt}{11pt}\selectfont
\textbf{\small EXAMPLE} \ 
The proposition can fail if the assumption ``$Y$ regular Hausdorff'' is weakened to ``$Y$ Hausdorff''.  Let $X$ be the set of nonnegative real numbers.  Put $D = \{1/n: n = 1, 2, \ldots\}$ $-$then the collection of all sets of the form $U \cup (V -D)$, where $U$ and $V$ are open in the usual topology on $X$, is also a topology, call the resulting space $Y$.  Observe that $Y$ is Hausdorff but not regular.  Let $A = X - D$ and define $f \in C(A,Y)$ by $f(x) = x$.  
It is clear that there is no $F \in C(X,Y) :$ $\restr{F}{A} = f$, 
yet for every open covering $\sV$ of $Y$, $X = \bigcup \text{ex}(f^{-1}(\sV))$.\\
\endgroup%%------------------------------------<<

\begingroup%%------------------------------------>>
\fontsize{9pt}{11pt}\selectfont
\textbf{\small FACT} \ 
Let $A$ be a dense subspace of a topological space \mX; let $Y$ be a regular Hausdorff space $-$then a given $f \in C(A,Y)$ admits a continuous extension 
$F \in C(X,Y)$ iff $\forall \ x \in X - A$ $\exists \ f_x \in C(A \cup \{x\},Y)$: $\restr{f_x}{A} = f$.\\
\endgroup%%------------------------------------<<

Let $X$ and $Y$ be topological spaces.

\index{EP}
\index{extension property with respect to Y}
\indent\indent (EP) \ A subspace $A \subset X$ is said to have the \un{extension property with respect} 
\un{ to $Y$} (EP. w.r.t. \mY) if $\forall \ f \in C(A,Y)$ $\exists \ F \in C(X,Y)$: $\restr{F}{A} = f$.

\index{NEP}
\index{neighborhood extension property with respect to Y}
\indent\indent (NEP) \  \ A subspace $A \subset X$ is said to have the 
\un{neighborhood extension property} 
\un{with respect to $Y$} (NEP. w.r.t. \mY) if $\forall \ f \in C(A,Y)$ $\exists$ 
$
\begin{cases}
\ U \supset A\\
\ F \in C(U,Y)
\end{cases}
(U \text{ open):} \restr{F}{A} = f.
$

\index{retract}
\index{neighborhood retract}
[Note: \  In this terminology, $A$ is a retract (neighborhood retract) of $X$ iff $A$ has the EP (NEP) w.r.t $Y$ for every $Y$.]

Two related special cases of importance are when $Y = \R$ or $Y = [0,1]$.  If $A$ has the EP w.r.t $\R$, then $A$ has the EP w.r.t $[0,1]$.  
Reason: If $f \in C(A,[0,1])$ and if $F \in C(X,\R)$ is a continuous extension of $f$, then $\min\{1, \max\{0,F\}\}$ is a continuous extension of $f$ with range a subset of $[0,1]$.  
The converse is trivially false.  
Example: Let $X$ be a CRH space $-$then $X$, as a subspace of $\beta X$, has the EP w.r.t. $[0,1]$ but $X$ has the EP w.r.t $\R$ iff $X$ is pseudocompact (of course in general $X$, as a subspace of $vX$, has the EP w.r.t. $\R$).  
Bear in mind that a CRH space is compact iff it is both $\R$-compact and pseudocompact.

%%----------------------------------------------------------------------------------------------04
[Note: \  Suppose that $X$ is Hausdorff $-$then $X$ is normal iff every closed subspace has the EP w.r.t. $\R$ (or, equivalently, $[0,1]$.]\\

\begingroup%%------------------------------------>>
\fontsize{9pt}{11pt}\selectfont
Suppose that $X$ is a CRH space.  Let $A$ be a subspace of $X$.
\\
\indent\indent ($\beta$) \ If $A$ has the EP w.r.t. $[0,1]$, then the closure of $A$ in $\beta X$ is $\beta A$ and conversely.\\
\indent\indent ($\nu$) \ If $A$ has the EP w.r.t. $\R$, then the closure of $A$ in $\nu X$ is $\nu A$ and conversely provided that $X$ is in addition normal.
\\ \indent
[Note: \  The Niemytzki plane is a nonnormal hereditarily $\R$-compact space, so the unconditional converse is false.]\\
\endgroup%%------------------------------------<<

Two subsets $A$ and $B$ of a topological space $X$ are said to be 
\un{completely separated}
\index{completely separated} 
in $X$ if $\exists \ \phi \in C(X,[0,1])$: 
$
\begin{cases}
\ \restr{\phi}{A} = 0\\
\ \restr{\phi}{B} = 1
\end{cases}
\hspace{-.25cm}.\ 
$
For this to be the case, it is necessary and sufficient that $A$ and $B$ are contained in disjoint zero sets.  Example: Suppose that $X$ is a CRH space $-$then any two disjoint closed subsets of $X$, one of which is compact, are completely separated in $X$ (no compactness assumption being necessary if $X$ is in addition normal).

[Note: \  It is enough to find a function $f \in C(X)$: 
$
\begin{cases}
\ \restr{f}{A} \leq 0\\
\ \restr{f}{B} \geq 1
\end{cases}
\hspace{-.25cm}. \ 
$
Reason: Take $\phi = \min\{1, \max\{0,f\}\}$.  
Moreover, 0 and 1 can be replaced by any real numbers $r$ and $s$ with $r < s$.]\\

\begin{proposition} \ %02
Let $A \subset X$ $-$then $A$ has the EP w.r.t. $[0,1]$ iff any two completely separated subsets of $A$ are completely separated in $X$.
\end{proposition}

[Assume that $A$ has the stated property.  
Fix an $f \in C(A,[0,1])$.  
To construct an extension $F \in C(X,[0,1])$ of $f$, we shall first define by recursion two sequences $\{f_n\}$ and $\{g_n\}$ subject to: 
$f_n \in BC(A)$ $\&$ $\norm{f_n} \leq 3r_n$ and 
$g_n \in BC(X)$ $\&$ $\norm{g_n} \leq r_n$, where $r_n = (1/2)(2/3)^n$ 
(so $\ds\sum\limits_1^\infty r_n = 1$).  Set $f_1 = f$.  Given $f_n$, let
$
\begin{cases}
\ S_n^- = \{x \in A: f_n(x) \leq -r_n\}\\
\ S_n^+ = \{x \in A: f_n(x) \geq r_n\}
\end{cases}
\hspace{-.25cm}. \ 
$
Since
$
\begin{cases}
\ S_n^-\\
\ S_n^+
\end{cases}
$
are completely separated in $A$, they are, by hypothesis, completely separated in $X$.  
Choose $g_n \in BC(X)$:
$
\begin{cases}
\ \restr{g_n}{S_n^-} = -r_n\\
\ \restr{g_n}{S_n^+} = r_n
\end{cases}
\& \ \norm{g_n} \leq r_n. 
$
Push the recursion forward by setting $f_{n+1} = f_n - \restr{g}{A}$.  
The series $\ds\sum\limits_1^\infty g_n$ is uniformly convergent on $X$, thus its sum $G$ is a continuous function on 
$X$ : $\restr{G}{A} = f$.  Take $F = \max\{0,G\}$.]\\

\label{19.31}
Application:  Suppose that $X$ is a CRH space $-$then any compact subset of $X$ has the EP w.r.t $[0,1]$ 
(cf. p. \pageref{6.1}).\\

%%----------------------------------------------------------------------------------------------05
\begingroup%%------------------------------------>>
\fontsize{9pt}{11pt}\selectfont
\textbf{\small FACT} \ 
Let $A \subset X$; let $f \in BC(A)$ $-$then $\exists$ $F \in BC(X)$ : $\restr{F}{A} = f$ iff $\forall \ a, b \in \R$: $a < b$, the sets
$
\begin{cases}
\ f^{-1}(]-\infty,a])\\
\ f^{-1}([b,+\infty[)
\end{cases}
$
are completely separated in \mX.\\
\endgroup%%------------------------------------<<

\begin{proposition} \ %03
Let $A \subset X$ $-$then $A$ has the EP w.r.t. $\R$ iff $A$ has the EP w.r.t $[0,1]$ and is completely separated from any zero set in $X$ disjoint from it.
\end{proposition}

[Necessity: Let $Z$ be a zero set in $X$ disjoint from $A$ : $Z = Z(g)$, where $g \in C(X,[0,1])$.  
Put $f = \restr{(1/g)}{A}$.  Choose $h \in C(X)$ : $\restr{h}{A} = f$. Consider $gh$.

Sufficiency:  Fix an $f \in C(A)$.  Because $\text{arctan} \circx f \in C(A,[-\pi/2,\pi/2])$, it has an extension 
$G \in C(X,[-\pi/2,\pi/2])$.  Let $B = G^{-1}(\pm \pi/2)$ $-$then $B$ is a zero set in $X$ disjoint from $A$, 
so there exists $\phi \in C(X,[0,1])$ :
$
\begin{cases}
\ \restr{\phi}{A} = 1\\
\ \restr{\phi}{B} = 0
\end{cases}
\hspace{-.25cm}. \ 
$
Put $F = \tan (\phi G)$ : $F \in C(X)$ $\&$ $\restr{F}{A} = f$.]\\

Consequently, every zero set in $X$ that has the EP w.r.t $[0,1]$ actually has the EP w.r.t $\R$.  On the other hand, a zero set in $X$ need not have the EP w.r.t $[0,1]$.  Examples:
(1) Take for $X$ the Isbell$-$Mr\'owka space $\Psi(\N)$ $-$then $A = \sS$ is a zero set in $X$ but $\sS$ does not have the EP w.r.t $[0,1]$;
(2) Take for $X$ the Niemytzki plane $-$then $A = \{(x,y): y = 0\}$ is a zero set in $X$ but $A$ does not have the EP w.r.t $[0,1]$.\\

\begingroup%%------------------------------------>>
\fontsize{9pt}{11pt}\selectfont
\index{Katet\"ov Space}
\textbf{\small EXAMPLE \ (\un{Katet\"ov Space})} \ 
As a subspace of $\R$, $\N$ has the EP w.r.t $[0,1]$, so the closure of $\N$ in $\beta \R$ is $\beta \N$.  
Let $X = \beta \R - (\beta \N - \N)$ $-$then $\beta X = \beta \R$ and $X$ is a LCH space which is actually pseudocompact 
(an unbounded continuous function on $X$ would be unbounded on a closed subset of $\R$ disjoint from $\N$).  
However, $X$ is not countably compact, thus is not normal (cf. $\S 1$, Proposition 5).  
As a subspace of $X$, $\N$ has the EP w.r.t $[0,1]$, but does not have the EP w.r.t $\R$.
\vspi
[Note: \ $\N$ is a closed $G_\delta$ but is not a zero set in $X$.]\\
\endgroup%%------------------------------------<<

\begingroup%%------------------------------------>>
\fontsize{9pt}{11pt}\selectfont
\index{Z-embedded}%\index{$\sZ$-embedded}
A subspace $A \subset X$ is said to be \un{$\sZ$-embedded} in $X$ if every zero set in $A$ is the intersection of $A$ with a zero set in $X$. 
Example:  Any cozero set in $X$ is $\sZ$-embedded in $X$.  If $A$ has the EP w.r.t $[0,1]$, 
then $A$ is $\sZ$-embedded in $X$ (but not conversely), so, e.g., any retract of $X$ is $\sZ$-embedded in $X$.
Examples: Suppose $X$ is Hausdorff $-$then
(1) \ Every subspace of a perfectly normal $X$ is $\sZ$-embedded in $X$;
(2) \ Every $F_\sigma$-subspace of a normal $X$ is $\sZ$-embedded in $X$;
(3) \ Every Lindel\"of subspace of a completely regular $X$ is $\sZ$-embedded in $X$.\\
\endgroup%%------------------------------------<<

\begingroup%%------------------------------------>>
\fontsize{9pt}{11pt}\selectfont
\textbf{\small FACT} \ 
Let $A \subset X$ $-$then $A$ has the EP w.r.t $\R$ iff $A$ is $\sZ$-embedded in $X$ and is completely separated from any zero set in $X$ disjoint from it.
\vspi
[Note: \  It is a corollary that if $A$ is a zero set in $X$, then $A$ has the EP w.r.t. $\R$ iff $A$ is $\sZ$-embedded in $X$.  
Both the Isbell-Mr\'owka space and the Niemytzki plane contain zero sets that are not $\sZ$-embedded.]
\\
\endgroup%%------------------------------------<<

%%----------------------------------------------------------------------------------------------06
\begingroup%%------------------------------------>>
\fontsize{9pt}{11pt}\selectfont
Application:  Suppose that $X$ is a Hausdorff space $-$then $X$ is normal iff every closed subset of $X$ is $\sZ$-embedded in $X$.\\
\endgroup%%------------------------------------<<

\begin{proposition} \ %04
Let $A \subset X$ $-$then $A$ has the EP w.r.t. $[0,1]$ ($\R$) iff for every finite (countable) numerable open covering $\sO$ of $A$ there exists a finite (countable) numerable open covering $\sU$ of $X$ such that $\sU \cap A$ is a refinement of $\sO$.
\end{proposition}

[The proof of necessity is similar to but simpler than the proof of the sufficiency so we shall deal just with it, 
assuming only that there exists a numerable open covering $\sU$ of $X$ such that $\sU \ \cap \ A$ 
is a refinement of $\sO$, thereby omitting the cardinality assumption on $\sU$.

\indent\indent ([0,1]) Let
$
\begin{cases}
\ S^{\prime}\\
\ S^{\prime\prime}
\end{cases}
$
be two completely separated subspaces of $A$; let
$
\begin{cases}
\ Z^{\prime}\\
\ Z^{\prime\prime}
\end{cases}
$
%^
be two disjoint zero sets in $A$ :
$
\begin{cases}
\ S^{\prime} \subset Z^{\prime}\\
\ S^{\prime\prime} \subset Z^{\prime\prime}
\end{cases}
\hspace{-.25cm}. \ 
$
Let $\sO = \{A - Z^\prime,A - Z^{\prime\prime}\}$.  
Take $\sU$ per $\sO$ and choose a neighborhood finite cozero set covering $\sV$ of $X$ such that $\sV$ is a star refinement of $\sU$ (cf. $\S 1$, Proposition 13).  
Put
$
\begin{cases}
\ W^{\prime} = X - \bigcup \{V \in \sV: V \cap Z^{\prime} = \emptyset\}\\
\ W^{\prime\prime} = X - \bigcup \{V \in \sV: V \cap Z^{\prime\prime} = \emptyset\}
\end{cases}
$
$-$then
$
\begin{cases}
\ W^{\prime}\\
\ W^{\prime\prime}
\end{cases}
$
are disoint zero sets in $X$ :
$
\begin{cases}
\ Z^{\prime} \subset W^{\prime}\\
\ Z^{\prime\prime} \subset W^{\prime\prime}
\end{cases}
\hspace{-.25cm}. \ 
$
Therefore $S^{\prime}$ and $S^{\prime\prime}$ are completely separated in $X$, thus, by Proposition 2, $A$ has the EP w.r.t. [0,1].

\indent\indent ($\R$)  Let $Z$ be a zero set in $X$ : $A \cap Z = \emptyset$, say $Z = Z(f)$, where $f \in C(X,[0,1])$.  The collection 
$\sO = \{f^{-1}(]1/n,1]) \cap A\}$ is a countable cozero set covering of $A$, hence is numerable 
(cf. p. \pageref{6.2}).  
Take $\sU$ per $\sO$ and choose a neighborhood finite cozero set covering $\sV = \{V_j : j \in J\}$ of $X$ and a zero set covering  $\sZ = \{Z_j: j \in J\}$ of $X$ such that $\sV$ is a refinement of $\sU$ with $Z_j \subset V_j$ ($\forall \  j$) 
(cf. p. \pageref{6.3}).  
Given $j$, $\exists \ n_j$: $Z_j \cap A \subset f^{-1}(]1/n_j,1]) \cap A$.  
Put 
$W = \bigcup\limits_j Z_j \cap f^{-1}(]1/n_j,1])$ 
$-$then $W$ is a zero set in $X$ containing $A$ and disjoint from $Z$, so $A$ and $Z$ are completely separated in $X$.  
Since the first part of the proof implies that $A$ necessarily has the EP w.r.t. [0,1], it follows from Proposition 3 that $A$ has the EP w.r.t $\R$.]\\

\begingroup%%------------------------------------>>
\fontsize{9pt}{11pt}\selectfont
\textbf{\small FACT} \ 
Let $A \subset X$$-$then $A$ is $\sZ$-embedded  in $X$ iff for every finite numerable open covering $\sO$ of $A$ there exists  a cozero set $U$ containing $A$ and a finite numerable open covering $\sU$ of $U$ such that $\sU \cap A$ is a refinement of $\sO$.\\
\endgroup%%------------------------------------<<

\textbf{\small LEMMA} \ 
Let $(X,d)$ be a metric space; let $A$ be a nonempty closed proper subspace of $X$ $-$then there exists a subset $\{a_i: i \in I\}$ of $A$ and a neighborhood finite open covering $\{U_i: i \in I\}$ of $X - A$ such that $\forall \ i$: $x \in U_i \implies d(x,a_i) \leq 2d(x,A)$.

[Assign to each $x \in X - A$ the open ball $B_x$ of radius $d(x,A)/4$.  The collection $\{B_x: x \in X - A\}$ is an open covering of $X - A$, thus by paracompactness has a neighborhood finite open refinement $\{U_i: i \in I\}$.  Each $U_i$ determines a point $x_i \in X - A$: $U_i \subset B_{x_i}$,
%%----------------------------------------------------------------------------------------------07
from which a point $a_i \in A$: $d(x_i, a_i) \leq (5/4)d(x_i,A)$.  Obviously $\forall \ x \in U_i$: $d(x,a_i) \leq (3/2)d(x_i,A)$ and $d(x_i,A) \leq (4/3)d(x,A)$.]\\

\index{Dugundji Extension Theorem}
\index{Theorem (Dugundji Extension Theorem)}
\textbf{\small DUGUNDJI EXTENSION THEOREM} \ 
Let $(X,d)$ be a metric space; let $A$ be a closed subspace of $X$.  
Let $E$ be a locally convex topological vector space.  Equip
$
\begin{cases}
\ C(A,E) \\
\ C(X,E)
\end{cases}
$
with the compact open topology $-$then there exists a linear embedding 
$\text{ext}: C(A,E) \ra C(X,E)$ such that $\forall \ f \in C(A,E)$, $\restr{\text{ext}(f)}{A} = f$ 
and the range of $\text{ext}(f)$ is contained in the convex hull of the range of $f$.

[Assume that $A$ is nonempty, proper and, using the notation of the lemma, choose a partition of unity 
$\{\kappa_i: i \in I\}$ on $X - A$ subordinate to $\{U_i: i \in I\}$.  
Given $f \in C(A,E)$, let
\[
\text{ext}(f)(x) = \ 
\begin{cases}
\ f(x) \hspace{2.43cm} (x \in A) \\
\ \ds\sum\limits_i \kappa_i(x)f(a_i) \hspace{1cm}  (x \in X - A)
\end{cases}
.
\]
Then $\restr{\text{ext}(f)}{A} = f$ and it is clear that $\text{ext}(f)(X)$ is contained in the convex hull of $f(A)$.  
The continuity of $\text{ext}(f)$ is built in at the points of $X - A$.
As for the points of $A$, fix $a_0 \in A$ and let $N$ be a balanced convex neighborhood of zero in $E$.  \ 
Choose a $\delta > 0$ : $d(a,a_0) \leq \delta$ $\implies$ $f(a) - f(a_0) \in N$ $(a \in A)$.  \ 
Suppose that 
$
\begin{cases}
\ x \in X- A \\
\ d(x,a_0) < \delta/3
\end{cases}
\hspace{-.25cm}. \ 
$
If $\kappa_i(x) > 0$, then, from the lemma, $d(x,a_i) \leq 2 d(x,A)$, hence $d(a_i,a_0) \leq 3 d(x,a_0) < \delta$.  
Consequently,
\[
\text{ext}(f)(x) - \text{ext}(f)(a_0) 
= 
\sum\limits_i \kappa_i(x)(f(a_i) - f(a_0)) 
\in 
\sum\limits_i \kappa_i(x)N 
\subset 
N.
\]
Therefore $\text{ext}(f) \in C(X,E)$.  
By construction, ext is linear and one-to-one, so the only remaining issue is its continuity.  
Take a nonempty compact subset $K$ of $X$ and let $O(K,N) = \{F \in C(X,E): F(K) \subset N\}$.  
Put $K_A = K \cap A \cup \{a_i \in A: K \cap U_i \neq \emptyset\}$.
Let $O(K_A,N) = \{f \in C(A,E): f(K_A) \subset N\}$.  
Plainly, $f \in O(K_A,N) \implies$ $\text{ext}(f) \in O(K,N)$.
Claim: $K_A$ is compact.  
To see this, let $\{x_n\}$ be a sequence in $K_A$.  
Since $K \cap A$ is compact, we can suppose that $\{x_n\}$ has no subsequence in $K \cap A$, 
%dmc previous sentence???
thus without loss of generality,  
$x_n = a_{i_n}$ for some $i_n$ : $K \cap U_{i_n} \neq \emptyset$.  
Pick $y_n \in K \cap U_{i_n}$ and assume that $y_n \ra y \in K$.
Case 1: $y \in K \cap A$.  
Here, $d(x_n,y) = d(a_{i_n},y) \leq 3d(y_n,y) \ra 0$.
Case 2: \ $y \in K \cap (X - A)$.  
There is a neighborhood of $y$ that meets finitely many of the $U_i$ and once $y_n$ is in this neighborhood, the index $i_n$ is constrained to a certain finite subset of $I$, which means that $\{x_n\}$ has a constant subsequence.]

[Note: \  Suppose that $E$ is a normed linear space $-$then the image of $\restr{\text{ext}}{BC(A,E)}$ is contained in $BC(X,E)$ and, per the uniform topology, $\text{ext}: BC(A,E) \ra BC(X,E)$ is a linear isometric embedding: $\forall \ f \in BC(A,E)$, $\norm{f} = \norm{\ext(f)}$.]\\

%%----------------------------------------------------------------------------------------------08
In passing, observe that if the $a_i$ are chosen from some given dense subset 
$A_0 \subset A$, then the range of $\ext(f)$ is contained in the union of $f(A)$ and the convex hull of $f(A_0)$.\\

\begingroup%%------------------------------------>>
\fontsize{9pt}{11pt}\selectfont
The Dugundji extension theorem has many applications.  To mention one, it is a key ingredient in the proof of a theorem of Milyutin to the effect that if $X$ and $Y$ are uncountable metrizable compact Hausdorff spaces, then $C(X)$ and $C(Y)$ are linearly homeomorphic 
(Pelczynski\footnote[2]{\textit{Dissertationes Math.} \textbf{58} (1968), 1-92; 
see also Semadeni, \textit{Banach Spaces of Continuous Functions}, PWN (1971), \textbf{379}.\vspace{0.11cm}}).  
Extensions to the case of noncompact $X$ and $Y$ have been given by 
Etcheberry\footnote[3]{\textit{Studia Math.} \textbf{53} (1975), 103-127; 
see also Hess, \textit{SLN} \textbf{991} (1983), 103-110.\vspace{0.11cm}}.
\\ \indent
[Note: \  The Banach$-$Stone theorem states that if $X$ and $Y$ are compact Hausdorff spaces, then $X$ and $Y$ are homeomorphic provided that the Banach spaces $C(X)$ and $C(Y)$ are isometrically isomorphic 
(Behrends\footnote[6]{\textit{SLN} \textbf{736} (1979), 138-140.\vspace{0.11cm}}).]\\
\endgroup%%------------------------------------<<

Is Dugundji's extension theorem true for an arbitrary topological vector space $E$?  In other words, can the ``locally convex'' supposition on $E$ be dropped?  The answer is ``no'', even if $E$ is a linear metric space 
(cf. p. \pageref{6.4}).

[Note: \  A topological vector space $E$ is said to be a 
\un{linear metric space}
\index{linear metric space} 
if it is metrizable.  Every linear metric space $E$ admits a translation invariant metric (Kakutani) but $E$ need not be normable.]\\

\begingroup%%------------------------------------>>
\fontsize{9pt}{11pt}\selectfont
Let $X$ be a CRH space; let $A$ be a nonempty closed subspace of $X$.  
Let $E$ be a locally convex topological vector space (normed linear space) 
$-$then a linear operator $T : C(A,E) \ra C(X,E)$ ($T: BC(A,E) \ra BC(X,E)$) continuous for the compact open topology (uniform topology) is said to be a 
\un{linear extension operator}
\index{linear extension operator} 
if for all $f$ in $C(A,E)$ ($BC(A,E)$) : $\restr{Tf}{A} = f$.  
Write LEO$(X,A;E)$ (LEO$_b$$(X,A;E)$) for the set of linear extension operators associated with 
$C(A,E)$ $(BC(A,E))$.  
Assuming that $X$ is metrizable, the Dugundji extension theorem asserts:  
$\forall$ $A$, $C(A,E)$ ($BC(A,E)$) possesses a linear extension operator 
(and even more in that the ``same'' operator works for both).  
Question: What conditions on $X$ or $A$ serve to ensure that LEO$(X,A;E)$ (LEO$_b$$(X,A;E)$) is not empty?\\
\endgroup%%------------------------------------<<

\label{19.24a}
\index{Michael Line}
\begingroup%%------------------------------------>>
\fontsize{9pt}{11pt}\selectfont
\textbf{\small EXAMPLE \  (\un{The Michael Line})} \ 
Take the set $\R$ and topologize it by isolating the points of $\PP$, leaving the points of $\Q$ with their usual neighborhoods.  
The resulting space $X$ is Hausdorff and hereditarily paracompact but not locally compact.  
And $A = \Q$ is a closed subspace of $X$ which, however, is not a $G_\delta$ in $X$.  
Let $E = C(\PP)$, $\PP$ in its usual topology, 
$-$then $E$ is a locally convex topological vector space (compact open topology). 
Claim:  LEO$(X,A;E)$ is empty.  For this, it suffices to exhibit an $f \in C(A,E)$
%%----------------------------------------------------------------------------------------------09
that cannot be extended to an $F \in C(X,E)$.  
If $\PP$ has its usual topology, then the continuous function
$
\begin{cases}
\ A \times \PP \ra \R\\
\ (x,y) \mapsto 1/(y - x)
\end{cases}
\vspace{0.2cm}
$
has no continuous extension $X \times \PP \ra \R$ (thus $X \times \PP$ is not normal).  
Defining $f \in C(A,E)$ by $f(x)(y) = 1/(y - x)$, it follows that $f$ has no extension $F \in C(X,E)$.\\
\endgroup%%------------------------------------<<

\begingroup%%------------------------------------>>
\fontsize{9pt}{11pt}\selectfont
A Hausdorff space is said to be 
\un{submetrizable}
\index{submetrizable}
if its topology contains a metrizable topology.  
Examples:
(1) \ The Michael line is submetrizable and normal but not perfect; 
(2) \ The Niemytzki plane is submetrizable and perfect but not normal.\\
\endgroup%%------------------------------------<<

\begingroup%%------------------------------------>>
\fontsize{9pt}{11pt}\selectfont
\textbf{\small FACT} \ 
Let $X$ be a submetrizable CRH space.  
Suppose that $A$ is a nonempty closed subspace of $X$ with a compact frontier 
$-$then $\forall \ E$, $\text{LEO}(X,A;E)$ ($\text{LEO}_b(X,A;E)$) is not empty.
\vspi
[Note: \  In view of the preceding example, the hypothesis on $A$ is not superfluous.]\\
\endgroup%%------------------------------------<<

\begingroup%%------------------------------------>>
\fontsize{9pt}{11pt}\selectfont
When $E = \R$, denote by LEO$(X,A)$ (LEO$_b(X,A)$) the set of linear extension operators for $C(A)$ ($BC(A)$).\\
\endgroup%%------------------------------------<<

\begingroup%%------------------------------------>>
\fontsize{9pt}{11pt}\selectfont
\textbf{\small EXAMPLE} \ 
$\text{LEO}_b(X,A)$ can be empty, even if $X$ is a compact Hausdorff space.  
For a case in point, take $X = \beta \N$ $\&$ $A = \beta \N - \N$.
Claim: $\text{LEO}(X,A)$ ($= \text{LEO}_b(X,A)$) is empty.
Suppose not and let $T:C(A) \ra C(X)$ be a linear extension operator.  
Fix an uncountable collection $\sU = \{U_i: i \in I\}$ of nonempty pairwise disjoint open subsets of $A$.  
Pick an $a_i \in U_i$ and choose $f_i \in C(A,[0,1])$ : 
$
\begin{cases}
\ f_i(a_i) = 1 \\
\ \restr{f_i}{(A - U_i)} = 0
\end{cases}
.
$
Let $O_i = \{x \in X: Tf_i(x) > 1/2\}$.  
Since $X$ is separable, there exists an uncountable subset $I_0$ of $I$ and a point 
$x_0 \in X$ : $x_0 \in \ds\bigcap\limits_{i \in I_0} O_i$.  
Let $n$ be some integer $> \norm{T}$.  
Select distinct indices $i_k$ $(k = 1, \ldots, 2n)$ in $I_0$.  
Put $f = \ds\sum\limits_1^{2n} f_{i_k}$, so $\norm{f} = 1$.  A contradiction then results by writing
\[
n 
\hsx = \hsx 
n\norm{f} 
\hsx \geq \hsx
\norm{Tf} 
\geq 
Tf(x_0) 
\hsx =  \hsx
\sum\limits_1^{2n}T f_{i_k}(x_0) 
\hsx > \hsx
2n \cdot \frac{1}{2} 
\hsx = \hsx
n.
\]
\indent
[Note: \  Let $X$ be a compact Hausdorff space; 
let $A$ be a nonempty closed subspace of $X$.  
Set 
$\rho(X,A) = \inf\{\norm{T}: T \in \text{LEO}(X,A)\}$ (where $\rho(X,A) = \infty$ if $\text{LEO}(X,A)$ is empty).  
Of course, $\rho(X,A) \geq 1$ and 
Benyamini\footnote[2]{\textit{Israel J. Math.} \textbf{16} (1973), 258-262.\vspace{0.11cm}} 
has shown that $\forall \ r$ : $1 \leq r < \infty$, there exists a pair $(X,A)$ : $\rho(X,A) = r$.]\\
\endgroup%%------------------------------------<<

\begingroup%%------------------------------------>>
\fontsize{9pt}{11pt}\selectfont
The space $X$ figuring in the preceding example is not perfect 
(no point of $\beta \N - \N$ is a $G_\delta$ in $\beta \N$).
Can one get a positive result if perfection is assumed?  The answer is ``no''.  Indeed, van 
Douwen\footnote[2]{\textit{General Topology Appl.} \textbf{5} (1975), 297-319.} 
has constructed an example of a CRH space $X$ that is simultaneously perfect and paracompact, 
yet contains a nonempty closed subspace $A$ for which $\text{LEO}_b(X,A) = \emptyset$.
\endgroup%%------------------------------------<<

%%----------------------------------------------------------------------------------------------10
\begingroup%%------------------------------------>>
\fontsize{9pt}{11pt}\selectfont
The assumption that $\text{LEO}_b(X,A)$ is not empty $\forall \ A$ has implications for the topology of $X$.  
To quantify the situation, given $r: 1 \leq r < \infty$, let $b_r$ be the condition: 
$\forall \ A$, $\{T \in \text{LEO}_b(X,A): \norm{T} \leq r\} \neq \emptyset$.  
Claim: If $b_r$ is in force, then for any discrete collection $\sA = \{A_i: i \in I\}$ of nonempty closed subsets of $X$ there is a collection $\sU = \{U_i: i \in I\}$ of open subsets of $X$ such that
(1) $A_i \subset U_i$ $\&$ $i \neq j$ $\implies$ $U_i \cap A_j = \emptyset$ and 
(2) $\text{ord}(\sU) \leq [r]$.  Thus put $A = \cup \sA$, let $\chi_i:A \ra [0,1]$ be the characteristic function of $A_i$, choose $T \in \text{LEO}_b(X,A)$ : $\norm{T} \leq r$, and consider $\sU = \{U_i: i \in I\}$, where 
$U_i = \{x \in X: T_{\chi_i}(x) > r/[r] + 1\}$.  
Example: Suppose that $X$ satisfies $b_r$ for some $r < 2$ $-$then $X$ is collectionwise normal.
\endgroup%%------------------------------------<<

\begingroup%%------------------------------------>>
\fontsize{9pt}{11pt}\selectfont
[Note: \  Let $X$ be the Michael line $-$then one can show that $X$ satisfies $\text{b}_1$,\  yet \ 
$\text{LEO}(X,A) = \emptyset$\  if $A = \Q$.]\\
\endgroup%%------------------------------------<<

\begingroup%%------------------------------------>>
\fontsize{9pt}{11pt}\selectfont
\textbf{\small FACT} \ 
Let $X$ be a Moore space.  
Assume: $X$ satisfies $\text{b}_r$ for some $r$ $-$then $X$ is normal and metacompact.\\
\endgroup%%------------------------------------<<

\label{6.6}
Let $X$ be a nonempty topological space $-$then an 
\un{equiconnecting structure}
\index{equiconnecting structure} 
on $X$ is a continuous function $\lambda:IX^2 \ra X$ such that $\forall \ x,y \in X$ and $\forall \ t \in [0,1]$:
$
\begin{cases}
\ \lambda(x,y,0) = x\\
\ \lambda(x,y,1) = y
\end{cases}
$
$\&$ $\lambda(x,x,t) = x$.  A subset $A \subset X$ for which $\lambda(IA^2)$ is called 
\un{$\lambda$-convex}.
\index{lambda-convex, $\lambda$-convex}  
In order that $X$ have an equiconnecting structure, 
it is necessary that $X$ be both contractible and locally contractible but these conditions are not sufficient as can be seen by considering Borsuk's cone 
(cf. p. \pageref{6.5}).  
Example: Suppose that $X$ is a contractible topological group.  
Let $H:IX \ra X$ be a homotopy contracting $X$ to its unit element $e$ 
$-$then the prescription $\lambda(x,y,t) = H(e,t)^{-1}H(xy^{-1},t)y$ defines an equiconnecting structure on $X$.  
In particular, if $X$ is a topological vector space, then $H(x,t) = (1 - t)x$ will do.

\label{6.11}
[Note: \  Let $E$ be an infinite dimensional Banach space.  
Consider $\bGL(E)$, the group of invertible bounded linear transformations $T:E \ra E$.  
Equip $\bGL(E)$ with the topology induced by the operator norm $-$then $\bGL(E)$ is a topological group and, 
being an open subset of a Banach space, has the homotopy type of a CW complex (cf. $\S 5$, Proposition 6).  
If $E$ is actually a Hilbert space, then $\bGL(E)$ is contractible 
\label{4.66a}
(Kuiper\footnote[3]{\textit{Topology} \textbf{3} (1965), 19-30.\vspace{0.11cm}}) 
but this need not be true in general (even if $E$ is reflexive), although it is the case for certain specific spaces, 
e.g., $C([0,1])$ or $L^p([0,1])$ $(1 \leq p \leq \infty)$.  
See 
Mityagin\footnote[6]{\textit{Russian Math. Surveys} \textbf{25} (1970), 59-103.\vspace{0.11cm}}
for proofs and other remarks.]\\

\begingroup%%------------------------------------>>
\fontsize{9pt}{11pt}\selectfont
\textbf{\small FACT} \ 
A nonempty topological space $X$ has an equiconnecting structure iff the diagonal $\Delta_X$ is a strong deformation retract of $X \times X$.
\vspi
[Necessity: Given $\lambda$, consider the homotopy $H:IX^2 \ra X^2$ defined by $H((x,y),t) = (\lambda(x,y,t),y)$.
\vspi
%%----------------------------------------------------------------------------------------------11
[Sufficiency: Given $H$, consider the equiconnecting structure  $\lambda:IX^2 \ra X$  defined by
\[
\lambda(x,y,t) = \ 
\begin{cases}
\ p_1(H((x,y),2t)) \hspace{0.95cm} \ \  (0 \leq t \leq 1/2)\\
\ p_2(H((x,y),2 - 2t)) \hspace{0.5cm} \  (1/2 \leq t \leq 1)
\end{cases}
,
\]
where $p_1$ and $p_2$ are the projections onto the first and second factors.]\\
\endgroup%%------------------------------------<<

\begingroup%%------------------------------------>>
\fontsize{9pt}{11pt}\selectfont
\textbf{\small FACT} \ 
Suppose that $X$ is a nonempty topological space for which  the inclusion $\Delta_X \ra X \times X$ is a cofibration $-$then $X$ has an equiconnecting structure iff $X$ is contractible.
\vspi
[Choose a homotopy $H:IX \ra X$ contracting $X$ to $x_0$:
$
\begin{cases}
\ H(x,0) = x\\
\ H(x,1) = x_0
\end{cases}
$
and then define $\Lambda:IX^2 \ra X^2$ by $\Lambda((x,y),t) = (H(x,t),H(y,t))$ to see that $\Delta_X$ is a weak deformation retract of $X \times X$.]\\
\endgroup%%------------------------------------<<

\label{6.7}
\begingroup%%------------------------------------>>
\fontsize{9pt}{11pt}\selectfont
A nonempty topological space $X$ is said to be 
\un{locally convex}
\index{locally convex} 
if it admits an equiconnecting structure $\lambda$ such that every $x \in X$ has a neighborhood basis comprised of $\lambda$-convex sets.  The convex subsets of a locally convex topological vector space are therefore locally convex, where $\lambda(x,y,t) = (1 - t)x + ty$.    
On the other hand, the long ray $L^+$ is not locally convex.\\
\endgroup%%------------------------------------<<

\label{6.41}
\begingroup%%------------------------------------>>
\fontsize{9pt}{11pt}\selectfont
\textbf{\small EXAMPLE} \ 
Let $K = (V,\Sigma)$ be a vertex scheme.  
Suppose that $K$ is 
\un{full}
\index{full (vertex scheme)}, 
i.e., if $F \subset V$ is finite and nonempty, 
then $F \in \Sigma$.  Claim: $\abs{K}$ is locally convex.  
Thus fix a point $* \in V$.  Let $\phi \in \abs{K}$ $-$then 
$\phi = \sum\limits_{v \neq *} b_v(\phi)\chi_v + \bigl(1 - \sum\limits_{v \neq *} b_v(\phi)\bigr)\chi_*$.
Here, $\chi_v$, ($\chi_*$) is the characteristic function of $\{v\}$ ($\{*\}$).  
Define $\beta:\abs{K} \times \abs{K} \ra \abs{K}$ 
by 
$\beta(\phi,\psi) = 
\sum\limits_{v \neq *} \beta(\phi,\psi)_v \chi_v + \bigl(1  -  \sum\limits_{v \neq *} \beta(\phi,\psi)_v \bigr)\chi_*$, where $\beta(\phi,\psi)_v = \min\{b_v(\phi),b_v(\psi)\}$.  
The assignment
\[
\lambda(\phi,\psi,t) = \ 
\begin{cases}
\ (1 - 2t)\phi + 2t\beta(\phi,\psi) \hspace{1.6cm}  (0 \leq t \leq 1/2)\\
\ (2 - 2t)\beta(\phi,\psi)  + (2t - 1)\psi \hspace{0.75cm}  (1/2 \leq t \leq 1)
\end{cases}
\]
is an equiconnecting structure on $\abs{K}$ relative to which $\abs{K}$ is locally convex.\\
\endgroup%%------------------------------------<<

\begingroup%%------------------------------------>>
\fontsize{9pt}{11pt}\selectfont
\textbf{\small FACT} \ 
Let $A \subset X$, where $X$ is metrizable and $A$ is closed $-$then $A$ 
has the EP w.r.t any locally convex topological space.\\
\endgroup%%------------------------------------<<

\index{Placement Lemma}
\textbf{\small PLACEMENT LEMMA} \ 
Every metric space $(X,d)$ can be isometrically embedded as a closed subspace of a normed linear space $E$, 
where $\text{wt} E = \omega \text{wt} X$.

[Denote by $\Sigma$ the collection of all nonempty finite subsets of $X$.  
Give $\Sigma$ the discrete topology.  
Fix a point $x_0 \in X$.  Attach to each $x \in X$ a function
$
f_x :
\begin{cases}
\ \Sigma \ra \R\\
\ \sigma \mapsto d(x,\sigma) - d(x_0,\sigma)
\end{cases}
$
$-$then $f_x \in BC(\Sigma)$ and the assignment
$
\iota :  
\begin{cases}
\ X \ra BC(\Sigma)\\
\ x \mapsto f_x
\end{cases}
$
is an isometric embedding.  
Note that $f_{x_0} \equiv 0$.  Let $E$ be the linear span of $\iota(X)$ in $BC(\Sigma)$.  
To see that $\iota(X)$ is closed in $E$, take a 
$\phi \in E - \iota(X)$, say 
$\phi = \ds\sum\limits_0^n r_i f_{x_i}$ ($r_i$ real), put $\sigma = \{x_0 \ldots, x_n\}$ and choose $\delta$ 
%%----------------------------------------------------------------------------------------------12
positive and less than $(1/2) \min\norm{\phi - f_{x_i}}$.  
Claim: No element of $\iota(X)$ can be within $\delta$ of $\phi$.  
Suppose not, so $\exists$ $x \in X$: $\norm{\phi - f_x} < \delta$.  Since $\iota$ is an isometry,
\[
d(x,x_i) = \norm{f_{x_i} - f_x} \geq \norm{\phi - f_{x_i}} - \norm{\phi - f_x} > 2\delta - \delta = \delta,
\]
from which $\norm{\phi - f_x}  \geq \abs{\phi(\sigma) - f_x(\sigma)} = d(x,\sigma) \geq \delta$, a contradiction.  There remains the assertion on the weights.  For this, let $D$ be a dense subset of $\iota(X)$ of cardinality $\leq \kappa$: $f_{x_0} \in D$ $-$then the linear span of $D$ is dense in $E$ and contains a dense subset of cardinality $\leq \omega\kappa$.]

\label{6.38}
[Note: \  One can obviously arrange that $E$ is complete provided this is the case of $(X,d)$.]\\

\begingroup%%------------------------------------>>
\fontsize{9pt}{11pt}\selectfont
\textbf{\small FACT} \ 
Every CRH space $X$ can be embedded as a closed subspace of a locally convex topological vector space $E$.\\
\endgroup%%------------------------------------<<

Let $Y$ be a nonempty metrizable space.\\
\indent\indent (AR) \ $Y$ is said to be an 
\un{absolute retract}
\index{absolute retract} 
(AR) if under any closed embedding $Y \ra Z$ into a metrizable space $Z$, the image of $Y$ is a retract of $Z$.\\
\indent\indent (ANR) \ $Y$ is said to be an 
\un{absolute neighborhood retract}
\index{absolute neighborhood retract} 
(ANR) if under any closed embedding $Y \ra Z$ into a metrizable space $Z$, the image of $Y$ is a neighborhood retract of $Z$.

[Note: \  There is no map from a nonempty set to the empty set, thus $\emptyset$ cannot be an AR, but there is a map from the empty set to the empty set, so we shall extend the terminology to agree that $\emptyset$ is an ANR.]\\

\begin{proposition} \ %05
Let $Y$ be a nonempty metrizable space $-$then $Y$ is an AR (ANR) iff for every pair $(X,A)$, where $X$ is metrizable and $A \subset X$ is closed, $A$ has the EP (NEP) w.r.t. $Y$.
\end{proposition}

[The indirect assertion is obvious.  
Turning to the direct assertion, in view of the placement lemma, $Y$ can be realized as a closed subspace of a normed linear space $E$.  
Assuming that $Y$ is an AR, fix a retraction $r:E \ra Y$.   If now $f:A \ra Y$ is a continuous function, then by the Dugundji extension theorem, $\exists$ $F \in C(X,E)$: $\restr{F}{A} = f$.  Consider $r \circx F$.]\\

\label{6.4}

\begingroup%%------------------------------------>>
\fontsize{9pt}{11pt}\selectfont
\textbf{\small EXAMPLE} \ 
Cauty\footnote[2]{\textit{Fund. Math.} \textbf{146} (1994), 85-99.} 
has given an example of a linear metric space $E$ which is not an absolute retract. 
So, for this $E$, the Dugundji extension theorem must fail.
\vspi
%%----------------------------------------------------------------------------------------------13
[Note: \  Therefore a metrizable space that has an equiconnecting structure need not be an AR.]\\
\endgroup%%------------------------------------<<

A countable product of nonempty metrizable spaces is an AR iff all the factors are ARs.  
Example: $[0,1]^n$, $\R^n$, $[0,1]^\omega$, and $\R^\omega$ are absolute retracts.  
A countable product of nonempty metrizable spaces is an ANR iff all the factors are ANRs and all but finitely many of the factors are ARs.
Example: $\bS^n$ and $\textbf{T}^n$ are absolute neighborhood retracts but
$
\begin{cases}
\ \bS^n \times \bS^n \times \cdots \\
\ \textbf{T}^n \times \textbf{T}^n \times \cdots
\end{cases}
$
($\omega$ factors) are not absolute neighborhood retracts.

Every retract (neighborhood retract) of an AR (ANR) is an AR (ANR).  
An open subspace of an ANR is an ANR.\\

\begingroup%%------------------------------------>>
\fontsize{9pt}{11pt}\selectfont
\label{4.72}
\textbf{\small EXAMPLE} \ 
Let $E$ be a normed linear space $-$then every nonempty convex subset of $E$ is an AR and every open subset of $E$ is an ANR.  
Assume in addition that $E$ is infinite dimensional.  
Let $S$ be the unit sphere in $E$ $-$then $S$ is an AR.  
To establish this, it need only be shown that $S$ is a retract of $D$, the closed unit ball in $E$.  
Fix a proper dense linear subspace $E_0 \subset E$ (the kernel of a discontinuous linear functional will do).  
In the notation of the Dugundji extension theorem, work with the pair $(D,S)$, picking the points defining ext in $S \cap E_0$, 
and let $f = \text{id}_S$ 
$-$then there exists a continuous function 
$\text{ext}(f):D \ra E$ such that 
$\restr{\text{ext}(f)}{S} = \text{id}_S$, with $\text{ext}(f)(D)$ contained in 
$S \cup (D \cap E_0)$, a proper subset of $D$.  
Choose a point $p$ in the interior of $D: p \notin \text{ext}(f)(D)$, 
let $r:D - \{p\} \ra S$ be the corresponding radial retraction and consider $r \circx \text{ext}(f)$.  
Corollary: Not every continuous function $D \ra D$ has a fixed point.
\\ \indent
[Note: \ 
There is another way to argue.  
Klee\footnote[2]{\textit{Proc. Amer. Math. Soc.} \textbf{7} (1956), 673-674.} 
has shown that if $E$ is an infinite dimensional normed linear space and if $K \subset E$ is compact, 
then $E$ and $E - K$ are homeomorphic.  
In particular, $E - \{0\}$ is homeomorphic to $E$, thus is an AR, and so $S$, being a retract of $E - \{0\}$ is an AR.  
Matters are trivial if $E$ is an infinite dimensional Banach space, since then $E$ is actually homeomorphic to $S$.]\\
\endgroup%%------------------------------------<<

\begingroup%%------------------------------------>>
\fontsize{9pt}{11pt}\selectfont
\textbf{\small EXAMPLE} \ 
Let $Y$ be any set lying between $]0,1[^n$ and $[0,1]^n$ $-$then $Y$ is an AR.  
Thus let $f$ be a closed embedding $Y \ra Z$ of $Y$ into a metrizable space $Z$.  
Call $j$ the inclusion $Y \ra [0,1]^n$, so $j \circx f^{-1} \in C(f(Y),[0,1]^n)$.  
Choose a $g \in C(Z,[0,1]^n)$: $\restr{g}{f(Y)} = j \circx f^{-1}$.  
Fix a compatible metric $d$ on $Z$ and define a continuous function 
$h:Z \ra [0,1]^n \times [0,1]$ by sending $z$ to $(g(z),\min\{1,d(z,f(Y))\})$.  
The range of $h$ is therefore a subset of $i_0Y \cup [0,1]^n \times ]0,1]$.  Let
$r:i_0Y \cup [0,1]^n \times ]0,1] \ra i_0Y$ be the retraction determined by projecting from the point $(1/2, \ldots, 1/2,-1) \in \R^{n+1}$ 
and let $p:i_0Y \ra Y$ be the canonical map.  
That$ f(Y)$ is a retract of $Z$ is then seen by considering the composite $f \circx p \circx r \circx h$.\\
\endgroup%%------------------------------------<<

\begingroup%%------------------------------------>>
\fontsize{9pt}{11pt}\selectfont
\textbf{\small FACT} \ 
Let $Y$ be an AR; let $B$ be a nonempty closed subspace of $Y$ $-$then $B$ is an AR iff $B$ is a strong deformation retract of $Y$.
\vspi
%%----------------------------------------------------------------------------------------------14
[To see that the condition is necessary, fix a retraction $r:Y \ra B$ and define a continuous function 
$h:i_0Y \cup IB \cup i_1Y \ra Y$ by
$
h(y,t) = 
\begin{cases}
\ y \quadx \ \ (y \in Y, t = 0) \\
\ y \quadx \ \ (y \in B, 0 \leq t \leq 1) \\
\ f(y) \ (y \in Y, t = 1) \\
\end{cases}
\hspace{-.2cm}. \ 
$
Since $i_0Y \cup IB \cup i_1Y$ is a closed subspace of $IY$ and since $B$ is an AR, 
it follows from Proposition 5 that $h$ has a continuous extension $H:IY \ra Y$.]\\
\endgroup%%------------------------------------<<

\label{14118.}
\begingroup%%------------------------------------>>
\fontsize{9pt}{11pt}\selectfont
Let $Y$ be an ANR $-$then $Y$ is homeomorphic to its diagonal $\Delta_Y$ which is therefore a strong deformation retract of 
$Y \times Y$ and this means that $Y$ has an equiconnecting structure 
(cf. p. \pageref{6.6}).

[Note: \  A metrizable locally convex topological space is an AR 
(cf. p. \pageref{6.7} 
and Proposition 5) but not every AR is locally convex.]\\
\endgroup%%------------------------------------<<

\label{6.20} %dmc mnft
\label{6.21} %dmc mnft
\begingroup%%------------------------------------>>
\fontsize{9pt}{11pt}\selectfont
\textbf{\small FACT} \ 
Let $Y$ be an ANR; let $B$ be a closed subspace of $Y$ $-$then $B$ is an ANR iff the inclusion $B \ra Y$ is a cofibration.
\vspi
[If $B$ is an ANR, then so is $i_0Y \cup IB$ 
(cf. p. \pageref{6.7a} 
(NES$_4$)), 
thus there exists a neighborhood \mO of 
$i_0Y \cup IB$  in $IY$ and a retraction $r:O \ra i_0Y \cup IB$.  
Choose a neighborhood $V$ of $B$ in $Y$: $IV \subset O$ and fix  
$
\phi \in C(Y,[0,1]):
\begin{cases}
\ \restr{\phi}{B} = 1\\
\ \restr{\phi}{Y-V} = 0
\end{cases}
\hspace{-.2cm}. \ 
$
Consider the map
$ 
\begin{cases}
\ IY \ra i_0Y \cup IB \\
\ (y,t) \mapsto r(y,\phi(y)t)
\end{cases}
$
\hspace{-.2cm}.]\\
\endgroup%%------------------------------------<<

\label{4.24}
\label{4.75}
\begingroup%%------------------------------------>>
\fontsize{9pt}{11pt}\selectfont
Let $Y$ be an ANR $-$then $Y$ is homeomorphic to its diagonal $\Delta_Y$, hence the inclusion $\Delta_Y \ra Y \times Y$ is a cofibration.  
Consequently, $Y$ is uniformly locally contractible 
(cf. p. \pageref{6.8}) 
and $\forall \ y_0 \in Y$, $(Y,y_0)$ is wellpointed 
(cf. p. \pageref{6.9}).
\vspi
[Note: \  It is unknown whether every metrizable uniformly locally contractible space is an ANR.  Any counterexample would necessarily have infinite topological dimension (cf. infra).]\\
\endgroup%%------------------------------------<<

\label{4.73}
Thanks to the placement lemma and the fact that a retract of a contractible (locally contractible) space is contractible (locally contractible), every AR (ANR) is contractible (locally contractible).  Both the broom and the cone over the Cantor set are contractible but, failing to be locally contractible, neither is an ANR.\\

\textbf{\small LEMMA} \ 
Suppose that $Y$ is a contractible ANR $-$then $Y$ is an AR.\\

A locally path connected topological space $X$ is said to be 
\un{locally $n$-connected}
\index{locally $n$-connected} 
$(n \geq 1)$ provide that for any $x \in X$ and any neighborhood $U$ of $x$ there exists a neighborhood 
$V \subset U$ of $x$ such that the arrow $\pi_q(V,x) \ra \pi_q(U,x)$ induced by the inclusion 
$V \ra U$ is the trivial map $(1 \leq q \leq n)$.  
If $X$ is locally $n$-connected for all $n$, then $X$ is called 
\un{locally homotopically trivial}.
\index{locally homotopically trivial}  
Example: A locally contractible space is locally homotopic\-ally trivial.\\

\begingroup%%------------------------------------>>
\fontsize{9pt}{11pt}\selectfont
\textbf{\small EXAMPLE} \ 
Working in $\ell^2$, let $p_k = (r_k(2k+1),0,\ldots)$, where $r_k  = 1/2k(k+1)$ $(k = 1, 2, \ldots)$, and put $p_0 = \lim\limits_k p_k$ $(= (0, 0, \ldots))$.  
Denote by $X_k(n)$ the set consisting of those points $x = \{x_i\}$: $x_i = 0$ 
%%----------------------------------------------------------------------------------------------15
$(i > n+1)$ and whose distance from $p_k$ is $r_k$.  
The union $\{p_0\} \cup \ds\bigcup\limits_{k= 1}^\infty X_k(n+1)$ 
is locally $n$-connected but not locally $(n+1)$-connected, while the union  
$\{p_0\} \ds\cup \bigcup\limits_{k= 1}^\infty X_k(k)$ is locally homotopically trivial but not locally contractible.\\
\endgroup%%------------------------------------<<

\label{19.49}
\label{19.55}

\begingroup%%------------------------------------>>
\fontsize{9pt}{11pt}\selectfont
Let $Y$ be a nonempty metrizable space.
\\
\indent\indent (LC$^n$) \ $Y$ is locally $n$-connected iff for every pair $(X,A)$, where $X$ is metrizable and $A \subset X$ is closed with $\dim(X-A) \leq n+1$, $A$ has the NEP w.r.t $Y$.
\\
\indent\indent (C$^n +$ LC$^n$) \ $Y$ is $n$-connected and locally $n$-connected iff for every pair $(X,A)$, where $X$ is metrizable and $A \subset X$ is closed with $\dim(X-A) \leq n+1$, $A$ has the EP w.r.t $Y$.
\vspi
Let $Y$ be a nonempty metrizable space of topological dimension $\leq n$.
\\
\indent\indent (LC$^n +$ $\dim \leq n$) \ $Y$ is locally $n$-connected iff $Y$ is locally contractible iff $Y$ is an ANR.
\\
\indent\indent (C$^n +$ LC$^n+$ $\dim \leq n$) \ $Y$ is $n$-connected and locally $n$-connected iff $Y$ is contractible 
and locally contractible iff $Y$ is an AR.
\vspi
The proofs of these results can be found in 
Dugundji\footnote[2]{\textit{Compositio Math.} \textbf{13} (1958), 229-246; 
see also Kodama, \textit{Proc. Japan Acad. Sci.} \textbf{33} (1957), 79-83.}
.
\vspi
\label{6.23}
[Note: \  It follows that a metrizable space of finite topological dimension is uniformly locally contractible iff it is an ANR and has an equiconnecting structure iff it is an AR.]\\
\endgroup%%------------------------------------<<

\label{6.5}
\begingroup%%------------------------------------>>
\fontsize{9pt}{11pt}\selectfont
\index{Borsuk's Cone}
\textbf{\small EXAMPLE \ (\un{Borsuk's Cone})} \ 
There exists a contractible, locally contractible compact metrizable space that is not an ANR.  
Choose a sequence: $0 = t_0 < t_1 < \cdots < 1$, $\lim t_n = 1$.  
Inside the product $\ds\prod\limits_0^\infty [0,1]$, for $n = 1, 2, \ldots$, 
form $Y_n = [t_{n-1},t_n] \times [0,1]^n \times 0 \times \cdots$, put 
$Y_\infty = 1 \times \ds\prod\limits_1^\infty [0,1]$, and let 
$Y = \bigl(\ds\bigcup\limits_1^\infty \text{fr}Y_n\bigr) \cup Y_\infty$ 
$-$then $Y$ is a compact connected metrizable space which we claim is locally contractible yet has nontrivial singular homology in every dimension, thus is not an ANR 
(cf. p. \pageref{6.10}).  
Local contractibility at the points of $Y - Y_\infty$ being obvious, let $y_\infty = (1, y_1, \ldots) \in Y_\infty$ and fix a neighborhood $U$ of $y_\infty$.  
There is no loss of generality in assuming that $U$ is the intersection of $Y$ with a set
$[a_0,1] \times [a_1,b_1] \times \cdots \times [a_k,b_k] \times [0,1] \times\cdots$.
Consider a neighborhood $V$ of $y_\infty$ that is the intersection of $Y$ with a set
$[a_0,1] \times [a_1,b_1] \times \cdots \times [a_k,b_k] \times [a_{k+1},b_{k+1}] \times [0,1] \times\cdots$, 
where $b_{k+1} - a_{k+1} < 1$.
There are two cases 
$ 
\begin{cases}
\ 1 \notin [a_{k+1},b_{k+1}] \\
\ 0 \notin [a_{k+1},b_{k+1}]
\end{cases}
\hspace{-.25cm}. \ 
$
As both are handled in a similar manner, suppose, e.g., that $1 \notin [a_{k+1},b_{k+1}]$ and define a homotopy $H:IV \ra U$ between the inclusion $V \ra U$ and the constant map $V \ra y_\infty$ by letting $H(v,t)$ be consecutively
\endgroup%%

\begingroup%%why ask why
\fontsize{9pt}{11pt}\selectfont
\[
\begin{cases}
\ (v_0, v_1, \ldots, v_k, (1-3t)v_{k+1},v_{k+2},\ldots) \\
\ (3t-1 + (2-3t)v_0, v_1, \ldots, v_k, 0, v_{k+2},\ldots) \\
(1, y_1 - 3(1-t)(y_1 - v_1), y_2 - 3(1-t)(y_2-v_2),\ldots)
\end{cases}
.
\]
\endgroup%%

\begingroup%%why ask why
\fontsize{9pt}{11pt}\selectfont
%%----------------------------------------------------------------------------------------------16
Here, $v = (v_0, v_1, \ldots) \in V$ and
$ 
\begin{cases}
\ 0 \leq t \leq 1/3 \\
\ 1/3 \leq t \leq 2/3 \\
\ 2/3 \leq t \leq 1
\end{cases}
\hspace{-.25cm}. \ 
$
That $Y$ is not an ANR is seen by remarking that fr$Y_n$ is a retract of $Y$, hence $H_n(\text{fr}Y_n) \approx \Z$ is isomorphic to a direct summand of $H_n(Y)$.  The cone $\Gamma Y$ of $Y$ is a contractible, locally contractible compact metrizable space.  
And $Y$, as a closed subspace of $\Gamma Y$, is a neighborhood retract of $\Gamma Y$.  Therefore $\Gamma Y$ is not an ANR.  
Finally, $Y$ is not uniformly locally contractible, so $\Gamma Y$ does not have an equiconnecting structure.
\vspi
[Note: \  Other, more subtle examples of this sort are known 
(Daverman-Walsh\footnote[2]{\textit{Michigan Math. J.} \textbf{30} (1983), 17-30.})
.]\\
\endgroup%%------------------------------------<<

\begingroup%%------------------------------------>>
\fontsize{9pt}{11pt}\selectfont
\textbf{\small FACT} \ 
Let $Y \subset \R^n$ $-$then $Y$ is a neighborhood retract of $\R^n$ iff $Y$ is locally compact and locally contractible.\\
\endgroup%%------------------------------------<<

\label{6.40}
\begingroup%%------------------------------------>>
\fontsize{9pt}{11pt}\selectfont
Haver\footnote[3]{\textit{Proc. Amer. Math. Soc.} \textbf{40} (1973), 280-284.} 
has shown that if a locally contractible metrizable space $Y$ can be written as a countable union of compacta of finite topological dimension, 
then $Y$ is an ANR.  
Example: Every metrizable CW complex $X$ is an ANR.  
Indeed, for this one can assume that $X$ is connected (cf. Proposition 12).  
But then $X$, being locally finite, is necessarily countable, hence can be written as a countable union of finite subcomplexes.\\
\endgroup%%------------------------------------<<

Certain function spaces or automorphism groups that arise ``in nature'' turn out to be ARs or, equivalently, contractible ANRs.  
Example:  Let $E$ be an infinite dimensional Hilbert space $-$then \textbf{GL}($E$) is contractible 
(cf. p. \pageref{6.11}).  
However, \textbf{GL}($E$) is an open subset of a Banach space, thus is an ANR.  Conclusion: \textbf{GL}($E$) is an AR.\\

\begingroup%%------------------------------------>>
\fontsize{9pt}{11pt}\selectfont
\index{measurable functions (example)}
\textbf{\small EXAMPLE \ (\un{Measurable Functions})} \ 
Let $Y$ be a nonempty metrizable space.  
Denote by $M_Y$ the set of equivalence classes of Borel measurable functions 
$f:[0,1] \ra Y$ equipped with the topology of convergence in measure 
$-$then $M_Y$ is metrizable, a compatible metric being given by the assignment
$(f,g) \lra \ds\int\limits_0^1 d(f(x),g(x)) dx$, where $d$ is a compatible metric on $Y$ bounded by 1.  
Nhu\footnote[6]{\textit{Fund. Math.} \textbf{124} (1984), 243-254.}
has shown that $M_Y$ is an ANR.  
Claim: $M_Y$ is contractible.  
To see this, fix a point $y_0 \in Y$ and consider the homotopy $H(f,t)(x) = $
$
\begin{cases}
\ f(x) \quadx (x > t)\\
\ y_0 \quadx \ \ \  (x \leq t)
\end{cases}
. \ 
$
Therefore $M_Y$ is an AR.
\\ \indent
[Note: \quadx Take $Y = \R$ $-$then $M_\R$ is a linear metric space.  But its dual $M_\R^*$ is trivial, hence $M_\R$ is not locally convex.]\\
\endgroup%%------------------------------------<<

\index{measurable transformations (example)}
\begingroup%%------------------------------------>>
\fontsize{9pt}{11pt}\selectfont
\textbf{\small EXAMPLE \  (\un{Measurable Transformations})} \ 
Let $\Gamma$ be the set of equivalence classes of measure preserving Borel measurable bijections $\gamma:[0,1] \ra [0,1]$, i.e., 
let $\Gamma$ be the automorphism group of the 
%%----------------------------------------------------------------------------------------------17
measure algebra \textbf{A} of the unit interval.  
Equip $\Gamma$ with the topology of pointwise convergence on \bA 
$-$then a subbasis for the neighborhoods at each fixed $\gamma_0 \in \Gamma$ is the collection of all sets of the form 
$\{\gamma: \abs{\gamma A \Delta \gamma_0 A} < \epsilon\}$ $(A \in \textbf{A} \ \& \ \epsilon > 0)$, $\Delta$ being symmetric difference.  
With respect to this topology, $\Gamma$ is a first countable Hausdorff topological group, so $\Gamma$ is metrizable.  
Nhu\footnote[5]{\textit{Proc. Amer. Math. Soc.} \textbf{110} (1990), 515-522.\vspace{0.11 cm}} 
has shown that $\Gamma$ is an ANR.  
Claim: $\Gamma$ is contractible.  
To see this, let $B$ be the complement of $A$ in $[0,1]$ and assign to each pair $(A,\gamma)$ its return partition, viz. the sequence $\{\Omega_n\}$, where $\Omega_0 = B$, $\Omega_1 = A \cap \gamma^{-1}A$, and for $n \geq 2$,
$\Omega_n = A \cap \gamma^{-1}B \cap \cdots \cap \gamma^{-(n-1)}B \cap \gamma^{-n}A$.  
Define $\gamma_A \in \Gamma$ by $\gamma_A(x) = \gamma^n(x)$ $(x \in \Omega_n)$, check that the map
$
\begin{cases}
\ \textbf{A} \times \Gamma \ra \Gamma\\
\ (A,\gamma) \mapsto \gamma_A
\end{cases}
$
is continuous, and consider the homotopy $H(t,\gamma) = \gamma_{[t,1]}$.  Therefore $\Gamma$ is an AR.
\\ \indent
[Note: \  Confining the discussion to the unit interval is not unduly restrictive since the Halmos-von Neumann theorem says that every separable, non atomic, normalized measure algebra is isomorphic to \textbf{A}.]\\
\endgroup%%------------------------------------<<

Let $X$ be a second countable topological manifold of euclidean dimension $n$.  Denote by $H(X)$ the set of all homeomorphims $X \ra X$ endowed with the compact open topology $-$then $H(X)$ is a topological group 
(cf. p. \pageref{6.12}).  
Moreover, $H(X)$ is metrizable and one can ask: Is $H(X)$ an ANR?  
If $X$ is not compact, then the answer is ``no'' since there are examples where $H(X)$ is not even locally contractible 
(Edwards-Kirby\footnote[2]{\textit{Ann. of Math.} \textbf{93} (1971), 63-88.\vspace{0.11 cm}}).  
If $X$ is compact, then $H(X)$ is locally contractible 
(\u Cernavski\u i\footnote[3]{\textit{Math. Sbornik} \textbf{8} (1969), 287-333; 
see also Rushing, \textit{Topological Embeddings}, Academic Press (1973), 270-293.\vspace{0.11 cm}}) 
and there is some evidence to support a conjecture that $H(X)$ might be an ANR.

[Note: \  If $X$ is not compact but is homeomorphic to the interior of a compact topological manifold with boundary, then $H(X)$ is locally contractible (\u Cernavski\u i (ibid.)).  
\label{2.5}
Example: $H(\R^n)$ is locally contractible.]\\

\begingroup%%------------------------------------>>
\fontsize{9pt}{11pt}\selectfont
\textbf{\small EXAMPLE} \ 
Take $X = [0,1]$ $-$then $H([0,1])$ is homeomorphic to $\R^\omega \times \{0,1\}$ (thus is an ANR).  In other words, the claim is that the identity component $H_e([0,1])$ of $H([0,1])$ is homeomorphic to $\R^\omega$.  
Form the product 
$\ds\prod\limits_{n = 0}^\infty \hspace{0.1cm} \ds\prod\limits_{i = 1}^{2^n} \hspace{0.1cm} ]0,1[_{n,i}$ and define a homeomorphism between it and $H_e([0,1])$ by assigning to a typical string $(x_{n,i})$ an order preserving homeomorphism $\phi:[0,1] \ra [0,1]$ via the following procedure.  
Suppose that $n$ is given and that there have been defined two sets of points
\[
\begin{cases}
\ A_n = \{0 = a(n,0) < a(n,1) < \cdots < a(n,2^n) = 1\}\\
\ B_n = \{0 = b(n,0) < b(n,1) < \cdots < b(n,2^n) = 1\}
\end{cases}
,
\]
with $\phi(a(n,i)) = b(n,i)$.  To extend the definition of $\phi$ to an order preserving bijection $A_{n+1} \ra B_{n+1}$, where
$
\begin{cases}
\ A_{n+1} \supset A_n\\
\ B_{n+1} \supset B_n
\end{cases}
$
and both have cardinality $2^{n+1} + 1$, distinguish two cases.  
Case 1: $n$ is odd.  Let $\alpha_i$
%%----------------------------------------------------------------------------------------------18
be the midpoint of $[a(n,i-1),a(n,i)]$ and set $\beta_i = \phi(\alpha_i) = x_{n,i}(b(n,i) - b(n,i-1)) + b(n,i-1)$.  
Case 2: $n$ is even.  Let $\beta_i$ be the midpoint of 
$[b(n,i-1),b_n,i)]$ and set $\alpha_i = \phi^{-1}(\beta_i) = x_{n,i}(a(n,i) - a(n,i-1)) + a(n,i-1)$.  Define 
$
\begin{cases}
\ A_{n+1} = A_n \cup \{\alpha_i: i = 1, \ldots, 2^n\}\\
\ B_{n+1}  =  B_n \cup \{\beta_i: i = 1, \ldots, 2^n\}
\end{cases}
$
, so that in the obvious notation
\[
\begin{cases}
\ A_{n+1} = \{0 = a(n+1,0) < a(n+1,1) < \cdots < a(n+1,2^{n+1}) = 1\}\\
\ B_{n+1} = \{0 = b(n+1,0) < b(n+1,1) < \cdots < b(n+1,2^{n+1}) = 1\}
\end{cases}
\hspace{-.25cm}, 
\]
with $\phi(a(n+1,i)) = b(n+1,i)$.  If now
$
\begin{cases}
\ A = \ds\bigcup\limits_1^\infty A_n\\
\ B = \ds\bigcup\limits_1^\infty B_n
\end{cases}
\hspace{-.25cm}, 
$
then
$
\begin{cases}
\ A\\
\ B
\end{cases}
$
are dense in $[0,1]$ and $\phi:A \ra B$ is an order preserving bijection, hence admits an extension to an order preserving homeomorphism $\phi:[0,1] \ra [0,1]$.
\vspi
[Note: \  $H([0,1])$ and $H(]0,1[)$ are homeomorphic.  In fact, the arrow of restriction $H([0,1]) \ra H(]0,1[)$ is continuous and has for its inverse the arrow of extension $H(]0,1[) \ra H([0,1])$, which is also continuous.  Corollary: $H(]0,1[)$  is an ANR.  Corollary $H(\R)$ is an ANR.]\\
\endgroup%%------------------------------------<<

\begingroup%%why ask why
\fontsize{9pt}{11pt}\selectfont
\textbf{\small EXAMPLE} \ 
Take $X = \bS^1$ $-$then $H(\bS^1)$ is homeomorphic to $\R^\omega \times \bS^1 \times \{0,1\}$ (thus is an ANR).  To see this, it suffices to observe that $H(\bS^1)$ is homeomorphic to $G \times \bS^1$, where $G$ is the subgroup of $H(\bS^1)$ consisting of those $\phi$ which fix $(1,0)$.\\
\endgroup%%------------------------------------<<

\begingroup%%why ask why
\fontsize{9pt}{11pt}\selectfont
Therefore, if $X$ is a compact 1-manifold, then $H(X)$ is an ANR.  Ths same conclusion obtains if $X$ is a compact 2-manifold 
(Luke-Mason\footnote[2]{\textit{Trans. Amer. Math. Soc.} \textbf{164} (1972), 275-285.}) 
but if $n > 2$, then it is unknown whether $H(X)$ is an ANR.\\
\endgroup%%------------------------------------<<

\begingroup%%why ask why
\fontsize{9pt}{11pt}\selectfont
\textbf{\small EXAMPLE} \ 
Take $X = [0,1]^\omega$, the Hilbert cube $-$then $H(X)$ (compact open topology) is metrizable and 
Ferry\footnote[3]{\textit{Ann. of Math.} \textbf{106} (1977), 101-119.} 
has shown that $H(X)$ is an ANR.\\
\endgroup%%------------------------------------<<

\textbf{\small LEMMA} \ 
Let $K = (V,\Sigma)$ be a vertex scheme $-$then $\abs{K}_b$ is an ANR.

[There are three steps to the proof.

\indent\indent (I) Fix a point $* \notin V$ and put $V_* = V \cup \{*\}$.  Let $\Sigma_*$ be the set of all nonempty finite subsets of $V_*$.  
Call $K_*$ the associated vertex scheme.  
Claim: $\abs{K_*}_b$ is an AR.  
Indeed, the inclusion $\abs{K_*}_b \ra \ell^1(V_*)$ is an isometric embedding with a convex range.

\indent\indent (II) Let $\Gamma_*$ be the subspace of $\abs{K_*}_b$ consisting of $\chi_*$, the characteristic function of 
$\{*\}$, and those $\phi \neq \chi_*$ : $\phi^{-1}(]0,1]) \cap V \in \Sigma$.  
Claim: $\Gamma_*$ is an AR.  To establish this, it suffices to exhibit a retraction 
$r:\abs{K_*}_b \ra \Gamma_*$.  Take a $\phi \in \abs{K_*}_b$.  
Case 1: $\phi = \chi_*$.  There is no choice here: $r(\chi_*) = \chi_*$.  
Case 2: $\phi \neq \chi_*$.  
Suppose that $\phi^{-1}(]0,1]) - \{*\} = \{v_0, \ldots, v_n\}$.  
%%----------------------------------------------------------------------------------------------19
Order the vertexes $v_i$ so that $\phi(v_0) \geq \cdots \geq \phi(v_n)$.  
Denote by $k$ the maximal index: $\{v_0, \ldots, v_k\} \in \Sigma$ and define $r(\phi)$ by the following formulas:\\
\[
\begin{cases}
\ r(\phi)(*) = 1 - \sum\limits_{v \in V} r(\phi)(v)\\
\ r(\phi)(v) = 0 \ (v \ \in V - \{v_0, \ldots, v_k\})
\end{cases}
\]
and
\[
\begin{cases}
\ k = n: r(\phi)(v_i) = \phi(v_i) \hspace{2.23cm} (0 \leq i \leq k)\\
\ k < n: r(\phi)(v_i) = \phi(v_i) - \phi(v_{k+1}) \hspace{0.5cm} (0 \leq i \leq k)
\end{cases}
.
\]
One can check that $r$ is welldefined and continuous.

\indent\indent (III) Since $\Gamma_* - \{\chi_*\}$ is open in $\Gamma_*$, it is an ANR.  
Claim: $\abs{K}_b$ is a retract of $\Gamma_* - \{\chi_*\}$, hence is an ANR.  
To see this, consider the map 
$\phi \ra \ds\frac{\phi - \phi(*)\chi_*}{1 - \phi(*)}$.]\\

A topological space is said to be a 
(\un{finite}, 
\un{countable}) 
\un{CW space}
\index{CW space}
\index{finite CW space}
\index{countable CW space}
if it has the homotopy type of a (finite, countable) CW complex.  The following theorems characterize these classes in terms of ANRs.\\

\label{5.8}
\index{Theorem (CW-ANR Theorem)}
\textbf{\small CW-ANR THEOREM} \ 
Let $X$ be a topological space $-$then $X$ has the homotopy type of a CW complex iff $X$ has the homotopy type of an ANR.

[If $X$ has the homotopy type of a CW complex, then there exists a vertex scheme $K$ such that $X$ has the homotopy type of $\abs{K}$ (cf. $\S 5$, Proposition 2) or still, the homotopy type of $\abs{K}_b$ (cf. $\S 5$, Proposition 1) and, by the lemma, $\abs{K}_b$ is an ANR.  
Conversely, if $X$ has the homotopy type of an ANR $Y$, use the placement lemma to realize $Y$ as a closed subspace of a normed linear space $E$.  
Fix an open $U \subset E$ : $U \supset Y$ and a retraction $r:U \ra Y$.  
Since $U$ has the homotopy type of a CW complex (cf. $\S 5$, Proposition 6), the domination theorem implies that the same is true of $Y$.]\\

\index{Theorem (Countable CW-ANR Theorem)}
\textbf{\small COUNTABLE CW-ANR THEOREM} \ 
Let $X$ be a topological space $-$then $X$ has the homotopy type of a countable CW complex iff $X$ has the homotopy type of a second countable ANR.

[If $X$ has the homotopy type of a countable CW complex, then there exists a countable locally finite vertex scheme $K$ such that $X$ has the homotopy type of $\abs{K}$ (cf. $\S 5$, Proposition 3 and 
p. \pageref{6.13}).  
Therefore, $\abs{K} = \abs{K}_b$ is Lindel\"of, hence second countable, and, by the lemma, $\abs{K}_b$ is an ANR.  Conversely, if $X$ has the homotopy type of a second countable ANR $Y$, then the ``$E$'' figuring in the preceding argument is second countable, therefore the ``$U$'' has the homotopy type of a countable CW complex (cf. $\S 5$, Proposition 6) and the countable domination theorem can be applied.]\\

%%----------------------------------------------------------------------------------------------20
\index{Theorem (Finite CW-ANR Theorem)}
\textbf{\small FINITE CW-ANR THEOREM} \ 
Let $X$ be a topological space $-$then $X$ has the homotopy type of a finite CW complex iff $X$ has the homotopy type of a compact ANR.

[One direction is easy: If $X$ has the homotopy type of a finite CW complex, then there exists a finite vertex scheme $K$ such that $X$ has the homotopy type of $\abs{K} = \abs{K}_b$ (cf. $\S 5$, Proposition 3), which, by the lemma, is an ANR.  The converse, however, is difficult: Its proof depends on an application of a number of theorems from infinite dimensional topology 
(West\footnote[2]{\textit{Ann. of Math.} \textbf{106} (1977), 1-18.\vspace{0.11cm}}).]\\

\label{6.10}
Application: The singular homology groups of a compact ANR are finitely generated and vanish beyond a certain point and the fundamental group of a compact connected ANR is finitely presented.\\

\begingroup%%------------------------------------>>
\fontsize{9pt}{11pt}\selectfont
According to the CW-ANR theorem, if $Y$ is an ANR, then it and each of its open subsets has the homotopy type of a CW complex.  
On the other hand, it can be shown that every metrizable space with this property is an ANR 
(Cauty\footnote[2]{\textit{Fund. Math.} \textbf{144} (1994), 11-22.\vspace{0.11cm}})
.\\
\endgroup%%------------------------------------<<

\begingroup%%------------------------------------>>
\fontsize{9pt}{11pt}\selectfont
\textbf{\small FACT} \ 
Let $Y$ be a nonempty metrizable space $-$then $Y$ is an AR iff $Y$ is a homotopically trivial ANR.
\vspi
[A connected CW complex is homotopically trivial iff it is contractible.  Quote the CW-ANR theorem.]\\
\endgroup%%------------------------------------<<

\begingroup%%------------------------------------>>
\fontsize{9pt}{11pt}\selectfont
Let $X$ and $Y$ be topological spaces.  Let $\sO = \{O\}$ be an open covering of $Y$ $-$then two continuous functions
$
\begin{cases}
\ f:X \ra Y\\
\ g:X \ra Y
\end{cases}
$
are said to be 
\un{$\sO$-contiguous}
\index{o-contiguous} %\index{$\sO$-contiguous}
if $\forall \ x \in X$ $\exists$ $O \in \sO$: $\{f(x),g(x)\} \subset O$.\\
\endgroup%%------------------------------------<<

\begingroup%%------------------------------------>>
\fontsize{9pt}{11pt}\selectfont
\textbf{\small LEMMA} \ 
Suppose that $Y$ is an ANR $-$then there exists an open covering $\sO = \{O\}$ of $Y$ 
such that for any topological space 
$
X : 
\begin{cases}
\ f \in C(X,Y) \\
\ g \in C(X,Y)
\end{cases}
$
$\sO-$contiguous $\implies f \simeq g$.
\vspi
[Choose a normed linear space $E$ containing $Y$ as a closed subspace. Fix a neighborhood $U$ of $Y$ in $E$ and a retraction $r:U \ra Y$.  Let $\sC = \{C\}$ be a covering of $U$ by convex open sets.  Put $\sO = \sC  \cap Y$.  Take two $\sO$-contiguous functions $f$ and $g$.  Define $h:IX \ra E$ by $h(x,t) = (1 - t)f(x) + t(g(x)$ $-$then $h(IX) \subset U$, so $H = r \circx h$ is a homotopy $IX \ra Y$ between $f$ and $g$.]\\ 
\endgroup%%------------------------------------<<

\begingroup%%------------------------------------>>
\fontsize{9pt}{11pt}\selectfont
Let $X$ be a topological space, $\sU = \{U\}$ an open covering of $X$.  Let $K = (V,\Sigma)$ be a vertex scheme $-$then a function $f:\abs{K^{(0)}} \ra X$ is said to be 
\un{confined}
\index{confined} 
by $\sU$ if $\forall \ \sigma \in \Sigma$ $\exists$ $U \in \sU$: $f(\abs{\sigma} \cap \abs{K^{(0)}}) \subset U$.\\
\endgroup%%------------------------------------<<

%%----------------------------------------------------------------------------------------------21
\begingroup%%------------------------------------>>
\fontsize{9pt}{11pt}\selectfont
\textbf{\small LEMMA} \ 
Suppose that $Y$ is an ANR. Let $\sO = \{O\}$ be an open covering of $Y$ 
$-$then there exists an open refinement $\sP = \{P\}$ of $\sO$ such that for every vertex scheme $K = (V,\Sigma)$ 
and every function $f:\absx{K^{(0)}} \ra Y$ confined by $\sP$ there exists a continuous function 
$F:\absx{K} \ra Y$ such that 
$\restr{F}{\absx{K^{(0)}}} = f$ and 
$\forall \ \sigma \in \Sigma$, $\forall \ P \in \sP$: 
$f(\abs{\sigma} \cap \absx{K^{(0)}}) \subset P$ $\implies$ 
$\exists$ $O \in \sO$: $F(\absx{\sigma}) \cup P \subset O$.
\vspi
[Choose a normed linear space $E$ containing $Y$ as a closed subspace.  
Fix a neighborhood $U$ of $Y$ in $E$ and a retraction $r:U \ra Y$.  
Let $\sC = \{C\}$ be a refinement of $r^{-1}(\sO)$ consisting of convex open sets.  
Put $\sP = \sC \cap Y$ $-$then $\sP$ is an open refinement of $\sO$ which we claim has the properties in question.  
Thus let $K = (V,\Sigma)$ be a vertex scheme.  
Take a function $f:\absx{K^{(0)}} \ra Y$ confined by $\sP$.  
Given $\sigma \in \Sigma$, write $C_\sigma$ for the convex hull of $f(\abs{\sigma} \cap \absx{K^{(0)}})$, 
itself a subset of some element $C \in \sC$.
Construct by induction continuous functions 
$\Phi_n:\absx{K^{(n)}} \ra U$ subject to 
$\Phi_0 = f$, $\restr{\Phi_{n+1}}{\absx{K^{(n)}}} = \Phi_n$, and 
$\forall \ \sigma \in \Sigma$, 
$\Phi_n(\abs{\sigma} \cap \absx{K^{(n)}}) \subset C_\sigma$.
Here the point is that if $\Phi_n$ has been constructed and if $\sigma$ is an (n+1)-simplex, then 
$\abs{\sigma} - \langle \sigma\rangle \subset \absx{K^{(n)}}$, therefore the restriction of 
$\Phi_n$ to  $\abs{\sigma} - <\sigma>$ can be continuously extended to $\abs{\sigma}$, $C_\sigma$ being an AR.  
This done, define $\Phi:\absx{K} \ra U$ by $\restr{\Phi}{\absx{K^{(n)}}} = \Phi_n$.  
Since each $\Phi_n$ is continuous, so is $\Phi$.  Consider $F = r \circx \Phi$.]\\
\endgroup%%------------------------------------<<

\begingroup%%------------------------------------>>
\fontsize{9pt}{11pt}\selectfont
These lemmas can be used to prove that if $Y$ is an ANR of topological dimension $\leq n$, then $Y$ is dominated in homotopy by $\abs{K}$, where $K$ is a vertex scheme: $\dim K \leq n$, a result not directly implied by the CW-ANR theorem.  In succession, let $\sO$ be an open covering of $Y$ per the first lemma, let $\sP$ be an open refinement of $\sO$ per the second lemma, and let $\sQ$ be a neighborhood finite star refinement of $\sP$ (cf. $\S 1$, Proposition 13) $-$then $\sQ$ has a precise open refinement $\sV$ of order $\leq n+1$ (cf. $\S 19$, Proposition 6).  
Obviously $\dim N(\sV) \leq n$, $N(\sV)$ the nerve of $\sV$.  Fix a point $y_V$ in each $V \in N(\sV)^{(0)}$.  Define $f:\abs{N(\sV)^{(0)}} \ra Y$ by $f(\chi_V) = y_V$.
Claim: $f$ is confined by $\sP$.  For suppose that $\sigma = \{V_1, \ldots, V_k\}$ is a simplex of $N(\sV)$.  
Since $V_1 \cap \cdots \cap V_k \neq \emptyset$ and since $\sV$ is a star refinement of $\sP$, 
there exists $P \in \sP$: $V_1 \cup \cdots \cup V_k \subset P$ 
$\implies$ $f(\abs{\sigma} \cap \abs{N(\sV)^{(0)}}) \subset P$.
Now take $F:\abs{N(\sV)} \ra Y$ as above and choose a $\sV$-map $G:Y \ra \abs{N(\sV)}$ 
(cf. p. \pageref{6.14}).  
One can check that $F \circx G$ and $\text{id}_Y$ are $\sO$-contiguous, hence homotopic.
\\ \indent
[Note: \  By analogous arguments, if $Y$ is a compact (connected) ANR of topological dimension $\leq n$, then $Y$ is dominated in homotopy by $\abs{K}$, where $K$ is a vertex scheme: $\dim K \leq n$ and $\abs{K}$ is compact (connected).]\\
\endgroup%%------------------------------------<<

\label{19.54}
\begingroup%%------------------------------------>>
\fontsize{9pt}{11pt}\selectfont
Application: Let $Y$ be an ANR of topological dimension $\leq n$ $-$then the singular homology groups of $Y$ vanish in all dimensions $> n$.\\
\endgroup%%------------------------------------<<

\label{19.53}
\begingroup%%------------------------------------>>
\fontsize{9pt}{11pt}\selectfont
\textbf{\small EXAMPLE} \ 
Suppose that $Y$ is a compact connected ANR: $\dim Y = 1$, $\&$ $\pi_1(Y) \neq 1$ $-$then $\pi_1(Y)$ is finitely generated and free.  Consequently, $Y$ has the homotopy type of a finite wedge of 1-spheres.\\
\endgroup%%------------------------------------<<

There are two variants of the CW-ANR theorem.

\indent\indent (Paired Version) \ A 
\un{CW pair}
\index{CW pair} 
is a pair $(X,A)$, where $X$ is a CW complex and $A \subset X$ is a subcomplex; an 
\un{ANR pair}
\index{ANR pair} 
is a pair $(Y,B)$, where $Y$ is an ANR and $B \subset Y$ is 
%%----------------------------------------------------------------------------------------------22
closed and an ANR.  Working then in the category of pairs of topological spaces, 
the result is that an arbitrary object in this category has the homotopy type of a CW pair iff it has the homotopy type of an ANR pair.

\label{6.16}
\indent\indent (Pointed Version) \ A 
\un{pointed CW complex}
\index{pointed CW complex} 
is a pair $(X,x_0)$, where $X$ is a CW complex and $x_0 \in X^{(0)}$; a 
\un{pointed ANR}
\index{pointed ANR} 
is a pair $(Y,y_0)$, where $Y$ is an ANR and $y_0 \in Y$.  
Working then in the category of pointed topological spaces, the result is that an arbitrary object in this category has the homotopy type of a pointed CW complex iff it has the homotopy type of a pointed ANR.

[Note: \  There is also a CW-ANR theorem for the category of pointed pairs of topological spaces.]\\

\begingroup%%------------------------------------>>
\fontsize{9pt}{11pt}\selectfont

In \textbf{HTOP}$^2$, the relevant reduction is that if $(X,A)$ is a CW pair, then there exists a vertex scheme $K$ and a subscheme $L$ such that 
$(X,A) \approx (\abs{K},\abs{L})$, while in \textbf{HTOP}$_*$, 
the relevant reduction is that if $(X,x_0)$ is a pointed CW complex, then there exists a vertex scheme $K$ and a vertex 
$v_0 \in V$ such that $(X,x_0) \approx (\abs{K},\abs{v_0})$ 
(cf. p. \pageref{6.15}).\\
\endgroup%%------------------------------------<<

Convention:  The function spaces encountered below carry the compact open topology.\\

\textbf{\small LEMMA} \ 
Let $X$, $Y$, and $Z$ be topological spaces.

\indent\indent (i) \ Let $f \in C(X,Y)$ $-$then the homotopy class of the precomposition arrow 
$f^*:C(Y,Z) \ra C(X,Z)$ depends only on the homotopy class of $f$.

\indent\indent (ii) \ Let $g \in C(Y,Z)$ $-$then the homotopy class of the postcomposition arrow 
$g_*:C(X,Y) \ra C(X,Z)$ depends only on the homotopy class of $g$.\\

\label{9.35}
\label{9.92}
\label{9.103}
Application: The homotopy type of $C(X,Y)$ depends only on the homotopy types of $X$ and $Y$.

[Note: \  By the same token, in \textbf{TOP}$^2$ the homotopy type of $(C(X,A;Y,B)$, 
$C(X,B))$ depends only on the homotopy types of $(X,A)$ and $(Y,B)$, 
whereas in \textbf{TOP}$_*$ the homotopy type of $C(X,x_0;Y,y_0)$ depends only on the homotopy types of $(X,x_0)$ and $(Y,y_0)$.]\\

\begin{proposition} \ %06
Let $K$ be a nonempty compact metrizable space; let $Y$ be a metrizable space $-$then $C(K,Y)$ is an ANR iff $Y$ is an ANR.
\end{proposition}

[Necessity: Assuming that $Y$ is nonempty, embed $Y$ in $C(K,Y)$ via the assignment $y \mapsto j(y)$, 
where $j(y)$ is the constant map $K \ra y$.  
Fix a point $k_0 \in K$ and denote by $e_0:C(K,Y) \ra Y$ the evaluation $\phi \ra \phi(k_0)$.  
Because $j \circx e_0$ is a retraction of $C(K,Y)$ onto $j(Y)$, it follows that if $C(K,Y)$ is an ANR, then so is $Y$.

%%----------------------------------------------------------------------------------------------23
[Sufficiency: Let $(X,A)$ be a pair, where $X$ is metrizable and $A \subset X$ is closed.  
Let $f:A \ra C(K,Y)$ be a continuous function.  Define a continuous function 
$\phi:A \times K \ra Y$ by setting 
$\phi(a,k) = f(a)(k)$.  Since $Y$ is an ANR, there is a neighborhood $O$ of $A \times K$ in $X \times K$ 
and a continuous function $\Phi:O \ra Y$ with $\restr{\Phi}{A \times K} = \phi$.  
Fix a neighborhood $U$ of $A$ in $X$ : $U \times K \subset O$.  
Define a continuous function $F:U \ra C(K,Y)$ by setting $F(u)(k) = \Phi(u,k)$.  
Obviously, $\restr{F}{A} = f$, thus $C(K,Y)$ is an ANR (cf. Proposition 5).]\\

\label{6.27}
Keeping to the above notation, the compactness of $K$ implies that $\pi_0(C(K,Y)) = [K,Y]$.  
Assume in addition that $Y$ is separable $-$then $C(K,Y)$ is separable.  But $C(K,Y)$ is also an ANR, hence its path components are open.  
Conclusion: $\#[K,Y] \leq \omega$.\\

\label{5.0aha}
Here is another corollary.  Suppose that $X$ is a finite CW space $-$then, on the basis of the CW-ANR theorem, for any CW space $Y$, $C(X,Y)$ has the homotopy type of an ANR, hence is again a CW space.

[Note: \  Some assumption on $X$ is necessary.  Example: Give $\{0,1\}$ the discrete topology and consider $\{0,1\}^\omega$.]\\

\begingroup%%------------------------------------>>
\fontsize{9pt}{11pt}\selectfont
\textbf{\small EXAMPLE} \ 
Let $X$ be a topological space $-$then the 
\un{free loop space}
\index{free loop space} 
$\Lambda X$ of $X$  is defined by the pullback square
\begin{tikzcd}[ sep=large]
{\Lambda X} \arrow{d} \arrow{r} &{PX} \arrow{d}{\Pi}\\
{X} \arrow{r} &{X \times X}
\end{tikzcd}
, where $\Pi$ is the Hurewicz fibration $\sigma \mapsto (\sigma(0),\sigma(1))$ and $X \ra X \times X$ is the diagonal embedding.  
The arrow $\Lambda X \ra X$ is a Hurewicz fibration and its fiber over $x_0$ is $\Omega(X,x_0)$, 
so if $X$ is path connected, then the homotopy type of $\Omega(X,x_0)$ is independent of the choice of $x_0$.  
Since $\Lambda X$ can be identified with $C(\bS^1,X)$ (compact open topology), the free loop space of $X$ is a CW space when $X$ is a CW space.
\vspi
[Note: Given a topological group \mG, define $W_G^\infty$ by the pullback square
\begin{tikzcd}[ sep=large]
{W_G^\infty} \arrow{d} \arrow{r} &{PX_G^\infty} \arrow{d}{\Pi}\\
{X_G^\infty \times G} \arrow{r}[swap]{\Phi} &{X_G^\infty \times X_G^\infty}
\end{tikzcd}, 
where $\Phi(x,g) = (x,x\cdot g)$ $-$then $W_G^\infty /G$ can be identified with 
$\Lambda B_G^\infty$ and there is a weak homotopy equivalence 
$\Lambda BG^\infty \ra (X_G^\infty \times G)/G$ (the action of $G$ on itself being by conjugation).]\\
\endgroup%%------------------------------------<<

\begingroup%%------------------------------------>>
\fontsize{9pt}{11pt}\selectfont
\textbf{\small EXAMPLE} \ 
Suppose that $X$ and $Y$ are path connected CW spaces for which there exists an $n$ such that 
(i) $X$ has the homotopy type of a locally finite CW complex with a finite $n$-skeleton and 
(ii) $\pi_q(Y) = 0$ $(\forall \ q > n)$ $-$then $C(X,Y)$ is a CW space.
\vspi
[Take $X$ to be a locally finite CW complex with a finite $n$-skeleton $X^{(n)}$.  
One can assume that $n$ is $> 0$ because when $n = 0$, $Y$ is contractible and the result is trivial.  
Consider the inclusion $i:X^{(n)} \ra X$ $-$then the precomposition arrow 
$i^*:C(X,Y) \ra C(X^{(n)},Y)$ is a Hurewicz fibration (cf. $\S 4$, Proposition 6)
%%----------------------------------------------------------------------------------------------24
and, in view of the assumption on $Y$, its fibers are either empty or contractible.  
But $C(X^{(n)},Y)$ is a CW space, thus so is $C(X,Y)$ (cf. Proposition 11).]\\
\endgroup%%------------------------------------<<

\begin{proposition} \ %07
Let $K$ be a nonempty compact metrizable space, $L \subset K$ a nonempty closed subspace; 
let $Y$ be a metrizable space, $Z \subset Y$ a closed subspace.  
Suppose that $Y$ is an ANR $-$then $C(K,L;Y,Z)$ is an ANR iff $Z$ is an ANR.
\end{proposition}

[Assuming that $Z$ is nonempty, one may proceed as in the proof of Proposition 6 and show that $Z$ is homeomorphic to a retract of $C(K,L;Y,Z)$, from which the necessity.  
Consider now a pair $(X,A)$, where $X$ is metrizable and $A \subset X$ is closed.  Let $f:A \ra C(K,L;Y,Z)$ be a continuous function.  
Define a continuous function $\phi:A \times L \ra Z$ by setting $\phi(a,\ell) = f(a)(\ell)$.  
Since $Z$ is an ANR, there is a neighborhood $O$ of $A \times L$ in $X \times L$ and a continuous function 
$\Phi:O \ra Z$ with $\restr{\Phi}{A \times L} = \phi$.  
Fix a neighborhood $U$ of $A$ in $X$: 
$\overline{U} \times L \subset O$.  
Define a continuous function $\psi:A \times K \cup \hsy \overline{U} \times L \ra Y$ by setting
$
\begin{cases}
\ \psi(a,k) = f(a)(k)\\
\ \psi(u,\ell) = \Phi(u,\ell)
\end{cases}
\hspace{-.25cm}
. \ 
$
Since $Y$ is an ANR, there is a neighborhood $P$ of $A \times K \cup \hsy \overline{U} \times L$ in $X \times K$ 
and a continuous function $\Psi:P \ra Y$ with 
$\restr{\Psi}{A \times K \cup \hsy \overline{U} \times L} = \psi$.  
Fix a neighborhood $V$ of $A$ in $X$: 
$ V \times K \subset P$ $\&$ $V \subset U$.  
Define a continuous function $F:V \ra C(K,L;Y,Z)$ by setting $F(v)(k) = \Psi(v,k)$.  
Obviously, $\restr{F}{A} = f$, thus $C(K,L;Y,Z)$ is an ANR (cf. Proposition 5).]\\

Take, e.g., $(K,L) = (\bS^n,s_n)$ $(s_n = (1, 0, \ldots, 0) \in \R^{n+1}, n \geq 1)$ and let $y_0 \in Y$ $-$then $\pi_n(Y,y_0) = \pi_0(C (\bS^n,s_n;Y,y_0))$.  
Accordingly, if $Y$ is separable, then $\pi_n(Y,y_0)$ is countable.
Example:  The homotopy groups of a countable connected CW complex are countable.\\

\label{5.0ab}
\index{Theorem (Loop Space Theorem)}
\textbf{\small LOOP SPACE THEOREM} \ 
Let $(X,x_0)$ be a pointed CW space $-$then the loop space $\Omega (X,x_0)$ is a pointed CW space.

[Fix a pointed ANR $(Y,y_0) $ with pointed homotopy type of $(X,x_0)$ 
(cf. p. \pageref{6.16}) 
$-$then $\Omega (Y,y_0) = C (\bS^1,s_1;Y,y_0)$ is a pointed ANR (cf. Proposition 7), so 
$\Omega(X,x_0) = C (\bS^1,s_1;X,x_0))$ is a pointed CW space.]\\

\begingroup%%------------------------------------>>
\fontsize{9pt}{11pt}\selectfont
\textbf{\small EXAMPLE} \ 
Suppose that $(X,x_0)$ is path connected and numerably contractible.  
Assume: $\Omega X$ is a CW space $-$then $X$ is a CW space.  Thus let $f:K \ra X$ be a pointed CW resolution.  
Owing to the loop space theorem, $\Omega K$ is a CW space.  
But the arrow $\Omega f: \Omega $K$ \ra \Omega X$ is a weak homotopy equivalence and since $\Omega X$ is a CW space, it follows from the realization theorem that $\Omega f$ is a homotopy equivalence.  
Therefore $f$ is a homotopy equivalence 
(cf. p. \pageref{6.17}).
\vspi
[Note: Let $X$ be the Warsaw circle $-$then $X$ is not a CW space.  On the other hand, there exists a continuous bijection $\phi:[0,1[ \ra X$ which is a regular Hurewicz fibration.  As this implies that $\phi$ is a pointed
%%----------------------------------------------------------------------------------------------25
Hurewicz fibration 
(cf. p. \pageref{6.18}), 
$\Omega X$ has the same pointed homotopy type as $\Omega [0,1[$ 
(cf. p. \pageref{6.18a}), 
hence is a CW space, so $X$ is not numerably contractible.]\\
\endgroup%%------------------------------------<<

\label{5.0ac}
\label{5.51}
\index{classifying spaces (example)}
\begingroup%%------------------------------------>>
\fontsize{9pt}{11pt}\selectfont
\textbf{\small EXAMPLE \ (\un{Classifying Spaces})} \ 
Let $G$ be a topological group $-$then $B_G^\infty$ is path connected and numerably contractible (inspect the Milnor construction).  
Moreover, according to $\S 4$, Proposition 36, $G$ and $\Omega B_G^\infty$ have the same homotopy type.  
Taking into account the preceding example, it follows that if $G$ is a CW space, then the same is true of $B_G^\infty$.  
Corollary: Any classifying space for $G$ is a CW space provided that $G$ itself is a CW space.\\
\endgroup%%------------------------------------<<

\textbf{\small LEMMA} \ 
Let $X \overset{f}{\ra} Z \overset{g}{\leftarrow} Y$ be a 2-sink.  Assume: $X,Y$, and $Z$ are ANRs $-$then $W_{f,g}$ is an ANR.\\

\begin{proposition} \ %08
Let $X \overset{f}{\ra} Z \overset{g}{\leftarrow} Y$ be a 2-sink.  Assume: $X,Y$, and $Z$ are CW spaces $-$then $W_{f,g}$ is a CW space.
\end{proposition}

[Fix ANRs 
$
\begin{cases}
\ X^\prime\\
\ Y^\prime
\end{cases}
\hspace{-.25cm}, \ 
$
homotopy equivalences
$
\begin{cases}
\ \phi:X^\prime \ra X\\
\ \psi:Y^\prime \ra Y
\end{cases}
\hspace{-.25cm}, \ 
$
and put
$
\begin{cases}
\ f^\prime = f \circx \phi\\
\ g^\prime = f \circx \psi
\end{cases}
$ $-$then there is a commutative diagram
\begin{tikzcd}%[ sep=large]
{X^\prime} \ar{d}[swap]{\phi} \arrow{r}{f^\prime} &{Z} \arrow[d,shift right=0.5,dash] \arrow[d,shift right=-0.5,dash] 
&{Y^\prime} \ar{d}{\psi} \ar{l}[swap]{g^\prime}\\
{X} \ar{r}[swap]{f} &{Z} &{Y} \ar{l}{g}
\end{tikzcd}
, thus the arrow $W_{f^\prime,g^\prime} \ra W_{f,g}$ is a homotopy equivalence 
(cf. p. \pageref{6.19}).  
Choose a homotopy equivalence $\zeta:Z \ra Z^\prime$, where $Z^\prime$ is an ANR.  There is an arrow 
 $W_{f^\prime,g^\prime} \ra W_{\zeta \circx f^\prime,\zeta \circx g^\prime}$ and it too is a homotopy equivalence.  But from the lemma, $W_{\zeta \circx f^\prime,\zeta \circx g^\prime}$ is an ANR.]\\

For a case in point, let $X$ and $Y$ be CW spaces $-$then $\forall \ f \in C(X,Y)$, $W_f$ is a CW space, and $\forall \ f \in C(X,x_0;Y,y_0)$, $E_f$ is a CW space.\\

\begingroup%%------------------------------------>>
\fontsize{9pt}{11pt}\selectfont
\textbf{\small FACT} \ 
Let $p:X \ra B$ be a regular Hurewicz fibration.  Assume: $\exists$ $b_0 \in B$ such that $\Omega(B,b_0)$ and $X_{b_0}$ are CW spaces $-$then $\forall \ x_0 \in X_{b_0}$, $\Omega(X,x_0)$ is a CW space.
\vspi
[By regularity, there is a lifting function 
$\Lambda_0:W_p \ra PX$ with the property that 
$\Lambda_0(x,\tau) \in j(X)$ whenever $\tau \in j(B)$.  
Define 
$f:\Omega(B,b_0) \ra X_{b_0}$ by $f(\tau) = \Lambda_0(x_0,\tau)(1)$, so $f(j(b_0)) = x_0$.  
The mapping fiber $E_f$ of $f$ has the same homotopy type as $\Omega(X,x_0)$.]\\
\endgroup%%------------------------------------<<

\begin{proposition} \ %09
Suppose that $p:X \ra B$ is a Hurewicz fibration and let $\Phi^\prime \in C(B^\prime,B)$.  Assume: $X$, $B$, and $B^\prime$ are CW spaces $-$then $X^\prime = B^\prime \times_B X$ is a CW space.
\end{proposition}

[In view of the preceding proposition, this follows from $\S 4$, Proposition 18.]\\

\label{5.50}
Application:  Let $p:X \ra B$ be a Hurewicz fibration, where $X$ and $B$ are CW spaces 
$-$then $\forall \ b \in B$, $X_b$ is a CW space.

%%----------------------------------------------------------------------------------------------26
[Note: \  Let $X$ be a CW space.  Relative to a base point, work first with 
\begin{tikzcd}[ sep=small]
{P X} \arrow{r}{p_0} &X
\end{tikzcd}
to see that $\Theta X$ is a CW space and then consider
\begin{tikzcd}[ sep=small]
{\Theta X} \arrow{r}{p_1} &X
\end{tikzcd}
to see that $\Omega X$ is a CW space, thereby obtaining an unpointed variant of the loop space theorem.]\\

\begin{proposition} \ %10
Suppose that $p:X \ra B$ is a Hurewicz fibration and let $O \subset B$.  
Assume: $X$ is an ANR, $B$ is metrizable, and the inclusion $O \ra B$ is a closed cofibration $-$then $X_O$ is an ANR.
\end{proposition}

[The inclusion $X_O \ra X$ is a closed cofibration (cf. $\S 4$, Proposition 11), a condition which is characteristic 
(cf. p. \pageref{6.20}).]\\

Application:  Let $p:X \ra B$ be a Hurewicz fibration, where $X$ and $B$ are ANRs 
$-$then $\forall \ b \in B$, $X_b$ is an ANR.

[Given $b \in B$, the inclusion $\{b\} \ra B$ is a closed cofibration 
(cf. p. \pageref{6.21}).]\\

\begingroup%%------------------------------------>>
\fontsize{9pt}{11pt}\selectfont
\textbf{\small EXAMPLE} \ 
Let $(Y,B,b_0)$ be a pointed pair.  Assume: $Y$ and $B$ are ANRs, with $B \subset Y$ closed.  Let $\Theta (Y,B)$ be the subspace of $\Theta Y$ consisting of those $\tau$ such that $\tau(1) \in B$ $-$then $\Theta (Y,B)$ is an ANR.  In fact, $\Theta Y$ is an ANR and there is a pullback square
\begin{tikzcd}[ sep=large]
{\Theta (Y,B)} \arrow{d} \arrow{r} &{\Theta Y} \arrow{d}{p_1}\\
{B} \arrow{r} &{Y}
\end{tikzcd}
.\\
\endgroup%%------------------------------------<<

\begingroup%%------------------------------------>>
\fontsize{9pt}{11pt}\selectfont
\textbf{\small EXAMPLE} \ 
Take $Y = \bS^n \times \bS^n \times \cdots $ ($\omega$ factors), $y_0 = (s_n, s_n, \ldots)$ $-$then $Y$ is not an ANR.   
Nevertheless, for every pair $(X,A)$, where $X$ is metrizable and $A \subset X$ is closed, $A$ has the HEP w.r.t $Y$
(cf. p. \pageref{6.22}).  
Therefore $\Theta Y$ is an AR.  Still, $\Omega Y$ is not an ANR.  
Indeed, none of the fibers of the Hurewicz fibration $p_1:\Theta Y \ra Y$ is an ANR.\\
\endgroup%%------------------------------------<<

\begin{proposition} \ %11
Suppose that $p:X \ra B$ is a Hurewicz fibration.  Assume: $B$ is a CW space and $\forall \ b\in B$, $X_b$ is a CW space $-$then $X$ is a CW space.
\end{proposition}

[Fix a CW resolution $f:K \ra X$.  Consider the Hurewicz fibration $q:W_f \ra X$ $(f = q \circx s)$.  
Since $s:K \ra W_f$ is a homotopy equivalence, $W_f$ is a CW space.  
Moreover, $q$ is a weak homotopy equivalence and the composite $p \circx q:W_f \ra B$ is a Hurewicz fibration.
The fibers $(p \circx q)^{-1}(b) = q^{-1}(X_b)$ are therefore CW spaces.  
Comparison of the homotopy sequences of $p \circx q$ and $p$ shows that the arrow 
$q_b:q^{-1}(X_b) \ra X_b$ is a weak homotopy equivalence, hence a homotopy equivalence.  
Because $B$ is numerably contractible (being a CW space), one can then apply 
$\S 4$, Proposition 20 to conclude that $q:W_f \ra X$ is a homotopy equivalence.]

[Note: \  If $p:X \ra B$ is a Hurewicz fibration and if $X$ and the $X_b$ are CW spaces, then it need not be true that $B$ is a CW space (consider the Warsaw circle).]\\
%%----------------------------------------------------------------------------------------------27

\begingroup%%------------------------------------>>
\fontsize{9pt}{11pt}\selectfont
Let $p:X \ra B$ be a Hurewicz fibration, where $X$ is metrizable and $B$ and the $X_b$ are ANRs.  Question: Is $X$ an ANR?  While the answer in unknown in general, the following lemma implies that the answer is ``yes'' provided that the topological dimension of $X$ is finite 
(cf. p. \pageref{6.23}).  
Infinite dimensional results can be found in 
Ferry\footnote[2]{\textit{Pacific J. Math.} \textbf{75} (1978), 373-382.}.\\
\endgroup%%------------------------------------<<

\begingroup%%------------------------------------>>
\fontsize{9pt}{11pt}\selectfont
\textbf{\small LEMMA} \ 
Suppose that $p:X \ra B$ is a Hurewicz fibration.  
Assume: $B$ is an ANR and $\forall \ b \in B$, $X_b$ is locally contractible $-$then $X$ is locally contractible.
\vspi
[Fix $x_0 \in X$, put $b_0 = p(x_0)$, and let $U$ be any neighborhood of $x_0$.  
Since $p$ has the slicing structure property 
(cf. p. \pageref{6.24}), 
it is an open map.  
Accordingly, one can assume at the outset that there is a contininuous function $\Phi:p(U) \ra PB$ such that
$
\begin{cases}
\ \Phi(b)(0) = b\\
\ \Phi(b)(1) = b_0
\end{cases}
\& \ \Phi(b_0)(t) = b_0 \ (0 \leq t \leq 1).  
$
Using the local contractibility of $X_{b_0}$, choose a neighborhood 
$O_0 \subset U \cap X_{b_0}$ of $x_0$ in $X_{b_0}$ and a homotopy 
$\phi:IO_0 \ra U \cap X_{b_0}$ satisfying 
$
\begin{cases}
\ \phi(x,0) = x\\
\ \phi(x,1)= u_0
\end{cases}
(u_0 \in U \cap X_{b_0}).
$
Fix a neighborhood $U_0$ of $x_0$: $U_0 \subset U$ and $O_0 = U_0 \cap X_{b_0}$.  Let $\Lambda_0:W_p \ra PX$ be a lifting function with the property that $\Lambda_0(x,\tau) \in j(X)$ whenever $\tau \in j(B)$.  
Define $F \in C(U,PX)$ by $F(x) = \Lambda_0(x,\Phi(p(x)))$.  
Because $F(x_0) =j(x_0) \in$ $\{\sigma \in PX:\sigma([0,1]) \subset U_0\}$, there is a neighborhood $V \subset U_0$ of $x_0$ such that $\forall \ x \in V$, $F(x)(t) \in U_0$ $(0 \leq t \leq 1)$.  If now $H:IV \ra U$ is the homotopy $H(x,t) = $ 
$
\begin{cases}
\ F(x)(2t)  \hspace{1.65cm}   (0 \leq t \leq 1/2)\\
\ \phi(F(x)(1),2t - 1) \hspace{0.35cm} (1/2 \leq t \leq 1)
\end{cases}
\hspace{-.25cm}, \ 
$
then
$
\begin{cases}
\ H(x,0) = x\\
\ H(x,1) = u_0
\end{cases}
\hspace{-.25cm}, \ 
$
i.e., the inclusion $V \ra U$ is inessential.]\\
\endgroup%%------------------------------------<<

Let $Y$ be a metrizable space.  Suppose that $Y$ admits a covering $\sV$ by pairwise disjoint open sets $V$, each of which is an ANR $-$then $Y$ is an ANR.  
To see this, assume that $Y$ is realized as a closed subspace of a metrizable space $Z$.  
Fix a compatible metric $d$ on $Z$.  
Given a nonempty $V \in \sV$, put $O_V = \{z: d(z,V) < d(z,Y - V)\}$ $-$then $O_V$ is open in $Z$ and $O_V \cap Y = V$.  
Moreover, the $O_V$ are pairwise disjoint.  
By hypothesis, there exists an open subset $U_V$ of $O_V$ containing $V$ and a retraction $r_V:U_V: \ra V$.  
Form $U = \bigcup\limits_V U_V$, a neighborhood of $Y$ in $Z$, and define a retraction 
$r:U \ra Y$ by $\restr{r}{U_V} = r_V$.

What is less apparent is that the same assertion is still true if the $V$ are not pairwise disjoint.\\

\textbf{\small LEMMA} \ 
Let $Y$ be a metrizable space.  
Suppose that $Y = Y_1 \cup Y_2$, where $Y_1$ and $Y_2$ are open and ANRs $-$then $Y$ is an ANR.

[This is proved in a more general context on 
p. \pageref{6.25} 
(cf. NES$_3$).]\\

\begin{proposition} \ %12
Let $Y$ be a metrizable space.  Suppose that $Y$ admits a covering $\sV$ by open sets \mV, each of which is an ANR $-$then $Y$ is an ANR.
\end{proposition}

%%----------------------------------------------------------------------------------------------28
[Use the domino principle 
(cf. p. \pageref{6.26}).]\\

\label{3.9}
\label{4.78}
\label{5.0v}
Application: Every metrizable topological manifold is an ANR, hence by the CW-ANR theorem has the homotopy type of a CW complex.\\

\label{5.0ag}
In particular, every compact topological manifold is an ANR, hence by the finite CW-ANR theorem has the homotopy type of a finite CW complex.  If $X$ and $Y$ are finite CW complexes, then $\#[X,Y] \leq \omega$ 
(cf. p. \pageref{6.27}).  
Specializing to the attaching process (and recalling that the inclusion $\bS^{n-1} \ra \textbf{D}^{n}$ is a closed cofibration), it follows that the set of homotopy types of compact topological manifolds is countable.

[Note: \  One can even prove that the set of homeomorphism types of compact topological manifolds is countable 
(Cheeger-Kister\footnote[2]{\textit{Topology} \textbf{9} (1970), 149-151.\vspace{0.11cm}}
).]\\
\label{5.52}

\begingroup%%------------------------------------>>
\fontsize{9pt}{11pt}\selectfont
The use of the term ``set'' in the above is justified by remarking that the full subcategory of \textbf{TOP} whose objects are the compact topological manifolds has a small skeleton.\\
\endgroup%%------------------------------------<<

\begingroup%%------------------------------------>>
\fontsize{9pt}{11pt}\selectfont
\textbf{\small EXAMPLE} \ 
Let $p:X \ra B$ be a covering projection.  Suppose that $X$ is metrizable and $B$ is an ANR $-$then $X$ is an ANR.
\vspi
[Note: \  The assumption that $X$ is metrizable is superfluous.]\\
\endgroup%%------------------------------------<<

\begingroup%%------------------------------------>>
\fontsize{9pt}{11pt}\selectfont
\textbf{\small EXAMPLE} \ 
Let $p:X \ra B$ be a Hurewicz fibration.  
Assume: $X$ is an ANR and $B$ is a path connected, numerably contractible, paracompact Hausdorff space $-$then $B$ is an ANR.  
For let $O$ be an open subset of $B$ with the property that the inclusion 
$O \ra B$ is inessential, say homotopic to $O \ra b.$  
Since $X_O$ is fiber homotopy equivalent to $O \times X_b$ 
(cf. p. \pageref{6.28}), 
$\text{sec}_O(X_O)$ is nonempty (cf. $\S 4$, Proposition 1), so $O$ is homeomorphic to a retract of $X_O$, an ANR.  
Therefore $B$ is locally an ANR, hence an ANR (recall that locally metrizable + paracompact $\implies$ metrizable; 
cf. p. \pageref{6.29}).\\
\endgroup%%------------------------------------<<

\begingroup%%------------------------------------>>
\fontsize{9pt}{11pt}\selectfont
\textbf{\small EXAMPLE} \ 
Let $X$ be an aspherical compact topological manifold.  
Assume: $\chi(X) \neq 0$ $-$then the path component of the identity in $C(X,X)$ is contractible.
\vspi
[Since $C(X,X)$ is an ANR (cf. Proposition 6), the  path component of the identity in 
$C(X,X)$ is a $K(\Cen\pi,1)$ 
(cf. p. \pageref{6.30} ff.), 
where $\pi = \pi_1(X)$.  
On the other hand, the assumption $\chi(X) \neq 0$ implies that $\Cen \pi$ is trivial.]\\
\endgroup%%------------------------------------<<

\begingroup%%------------------------------------>>
\fontsize{9pt}{11pt}\selectfont
Let $X$ and $Y$ be metrizable spaces. \  Let $A$ be a closed subspace of $X$ and let \ $f:A \ra Y$ be a continuous function $-$then 
Borges\footnote[3]{\textit{Proc. Amer. Math. Soc.} \textbf{24} (1970), 446-451.\vspace{0.11cm}} 
has shown that $X \sqcup_f Y$ is metrizable iff every point of $X \sqcup_f Y$ belongs
%%----------------------------------------------------------------------------------------------29
to a compact subset of countable character, i.e., having countable neighborhood basis in $X$.  
In particular, this condition is satisfied if $X \sqcup_f Y$ is first countable or if $A$ is compact.
\vspi
[Note: \  In any event, $X \sqcup_f Y$ is a perfectly normal paracompact Hausdorff space (AD$_5$ 
(cf. p. \pageref{6.31})).]\\
\endgroup%%------------------------------------<<

\label{19.57}
\begingroup%%------------------------------------>>
\fontsize{9pt}{11pt}\selectfont
\textbf{\small LEMMA} \ 
Let $B$ be a closed subspace of a metrizable space $Y$ such that the inclusion $B \ra Y$ is a cofibration.  
Suppose that $B$ and $Y - B$ are ANRs $-$then $Y$ is an ANR.
\vspi
[Fix a Str{\o}m structure $(\psi,\Psi)$ on $(Y,B)$ and put $V = \psi^{-1}([0,1[)$.  Show that $V$ is an ANR.]\\
\endgroup%%------------------------------------<<

\begingroup%%------------------------------------>>
\fontsize{9pt}{11pt}\selectfont
\textbf{\small FACT} \ 
Let $X$ and $Y$ be ANRs.  Let $A$ be a closed subspace of $X$ and let $f:A \ra Y$ be a continuous function.  Suppose that $A$ is an ANR $-$then $X \sqcup_f Y$ is an ANR provided that it is metrizable.\\
\endgroup%%------------------------------------<<

\begingroup%%------------------------------------>>
\fontsize{9pt}{11pt}\selectfont
\textbf{\small LEMMA} \ 
Let $B$ be a closed subspace of a metrizable space $Y$ such that the inclusion $B \ra Y$ is a cofibration.  
Suppose that $B$ is an AR and $Y - B$ is an ANR$-$then $Y$ is an AR if $B$ is a strong deformation retract of $Y$.
\vspi
[It follows from the previous lemma that $Y$ is an ANR.  But $Y$ and $B$ have the same homotopy type and $B$ is contractible.]\\
\endgroup%%------------------------------------<<

\begingroup%%------------------------------------>>
\fontsize{9pt}{11pt}\selectfont
\textbf{\small FACT} \ 
Let $X$ and $Y$ be ARs.  Let $A$ be a closed subspace of $X$ and let $f:A \ra Y$ be a continuous function.  Suppose that $A$ is an AR $-$then $X \sqcup_f Y$ is an AR provided that it is metrizable.\\
\endgroup%%------------------------------------<<

\begingroup%%------------------------------------>>
\fontsize{9pt}{11pt}\selectfont
\textbf{\small EXAMPLE} \ 
Take $X = [0,1]^2$, $A = [1/4,3/4] \times \{1/2\}$, $Y = [0,1]^3$ and let $f:A \ra Y$ be a continuous surjective map $-$then $X \sqcup_f Y$ is a compact AR of topological dimension 3, yet it is not homeomorphic to any CW complex.\\
\endgroup%%------------------------------------<<

\begingroup%%------------------------------------>>
\fontsize{9pt}{11pt}\selectfont
Let $(X,A)$ be a CW pair.  Is it true that $A$ has the EP w.r.t. any locally convex topological vector space?  A priori, this is not clear since CW complexes are not metrizable in general.  
There is, however, a class of topologically significant spaces, encompassing both the class of metrizable spaces and the class of CW complexes for which a satisfactory extension theory exists.
\\ \indent
Let $X$ be a Hausdorff space; let $\tau$ be the topology on $X$ $-$then $X$ is said to be 
\un{stratifiable}
\index{stratifiable} 
if there exists a function $\text{ST}_X:\N \times \tau \ra \tau$, termed a 
\un{stratification}
\index{stratification}, 
such that 
(a) \ $\forall \ U \in \tau$, $\overline{\text{ST}_X(n,U)} \subset U$;
(b) \ $\forall \ U \in \tau$, $\bigcup\limits_n \text{ST}_X(n,U) = U$;
(c) \ $\forall \ U,V \in \tau$ : $U \subset V \implies \text{ST}_X(n,U) \subset \text{ST}_X(n,V)$.
A stratifiable space is perfectly normal and every subspace of a stratifiable space is stratifiable.  
A finite or countable product of stratifiable spaces is stratifiable.  
A stratifiable space need not be compactly generated and a compactly generated space need not be stratifiable, even if it is regular and countable 
(Foged\footnote[2]{\textit{Proc. Amer. Math. Soc.} \textbf{81} (1981), 337-338; 
see also \textit{Proc. Amer. Math. Soc.} \textbf{92} (1984), 470-472.\vspace{0.11cm}}).
\label{3.8}
%%----------------------------------------------------------------------------------------------30
Example: Every metrizable space is stratifiable.
Example: The Sorgenfrey line, the Niemytzki plane, and the Michael line are not stratifiable.
\\ \indent
[Note: \  
Junnila\footnote[3]{\textit{Colloq. Math. Soc. J\'anos Bolyai} \textbf{23} (1980), 689-703; 
see also Harris, \textit{Pacific J. Math.} \textbf{91} (1980), 95-104.\vspace{0.11cm}}  
has shown that every topological space is the open image of a stratifiable space.]\\
\endgroup%%------------------------------------<<

\begingroup%%------------------------------------>>
\fontsize{9pt}{11pt}\selectfont
\textbf{\small FACT} \ 
Let $X$ be a topological space; let $\sA = \{A_j: j \in J\}$ be an absolute closure preserving closed covering of $X$.  Suppose that each $A_j$ is stratifiable $-$then $X$ is stratifiable.
\vspi
[$X$ is necessarily a perfectly normal Hausdorff space 
(cf. p. \pageref{6.32}).  
As for stratifiability, consider the set $\sP$ of all pairs $(I,\text{ST}_I)$, 
where $I \subset J$ and $\text{ST}_I$ is a stratification of $X_I = \bigcup\limits_i A_i$.
Order $\sP$ by stipulating that $(I^\prime,\text{ST}_{I^{\prime}}) \leq (I^{\prime\prime},\text{ST}_{I^{\prime\prime}})$ iff $I^\prime \subset I^{\prime\prime}$ and for each open subset $U$ of $X_{I\pp}:$
\endgroup%%------------------------------------<<

\begingroup%%------------------------------------>> %dmc so adhoc
\fontsize{9pt}{11pt}\selectfont
\vspace{-.35cm}
\[
\text{ST}_{I^{\prime\prime}}(n,U) \cap X_{I^{\prime}} = \text{ST}_{I^{\prime}}(n,U \cap X_{I^{\prime}} ) \quadx\& \quadx \overline{\text{ST}_{I^{\prime\prime}}(n,U)} \cap X_{I^{\prime}} = \overline{\text{ST}_{I^{\prime}}(n,U \cap X_{I^{\prime}})}.
\]
\vspace{-.65cm}

\noindent Every chain in $\sP$ has an upper bound, so by Zorn, $\sP$ has a maximal element $(I_0,\text{ST}_{I_0})$.  Verify that $X_{I_0} = X$.]\\
\endgroup%%------------------------------------<<

\begingroup%%------------------------------------>>
\fontsize{9pt}{11pt}\selectfont
\label{3.6}
Application: Every CW complex is stratifiable.
\vspi
[The collection of finite subcomplexes of a CW complex $X$ is an absolute closure preserving closed covering of $X$.]
\vspi
Application: Let $E$ be a vector space over $\R$.  Equip $E$ with the finite topology $-$then $E$ is stratifiable.
\vspi
[Fix a basis $\{e_i: i \in I\}$ for $E$.  Assign to each finite subset of $I$ the span of the corresponding $e_i$.  The resulting collection of linear subspaces is an absolute closure preserving closed covering of $E$.]\\[.1cm]
\endgroup%%------------------------------------<<

\begingroup%%------------------------------------>>
\fontsize{9pt}{11pt}\selectfont
\textbf{\small FACT} \ 
Suppose that $X$ and $Y$ are stratifiable $-$then the coarse joine $X *_c Y$ is stratifiable.\\[.1cm]
\endgroup%%------------------------------------<<

\begingroup%%------------------------------------>>
\fontsize{9pt}{11pt}\selectfont
Application:  Let $G$ be a stratifiable topological group $-$then $\forall \ n$, $X_G^n$ is stratifiable.\\[.1cm]
\endgroup%%------------------------------------<<

\begingroup%%------------------------------------>>
\fontsize{9pt}{11pt}\selectfont
\textbf{\small LEMMA} \ 
Let $X = \ds\bigcup\limits_0^\infty X_n$ be a topological space, where 
$X_n \subset X_{n+1}$ and $X_n$ is stratifiable and a zero set in $X$, say $X_n = \phi_n^{-1}(0)$ $(\phi_n \in C(X,[0,1])$.  
Suppose that there is a retraction 
$r_n: \phi_n^{-1}([0,1[) \ra X_n$ such that 
$\forall \ x \in X_n - X_{n-1}$ $(X_{-1} = \emptyset)$, the sets
$r_n^{-1}(U) \cap \phi_n^{-1}([0,t[)$ form a neighborhood basis of $x$ in $X$ 
($U$ a neighborhood of $x$ in $X_n$ and $0 < t \leq 1$) $-$then $X$ stratifiable.
\vspi
[The assumptions imply that $X$ is Hausdorff.  
To construct $\text{ST}_X$, fix a stratification  $\text{ST}_{X_n}$ of $X_n$: $\text{ST}_{X_n}(k,U) \subset \text{ST}_{X_n}(k+1,U)$.  
Given an open subset $U$ of $X$, denote by $U(n,k)$ the interior of
\[
\{x \in X_n: r_n^{-1}(x) \cap \phi_n^{-1}([0,1/(k+1)[) \subset U\}
\]
in $X_n$ and for $N = 1, 2, \ldots$, put
\[
\text{ST}_{X_n}(N,U) \ = \ 
\bigcup\limits_{n,k \leq N} r_n^{-1}(\text{ST}_{X_n}(N,U(n,k))) \cap \phi_n^{-1}([0,1/(k+2)[).]
\] 
\endgroup%%------------------------------------<<

%%----------------------------------------------------------------------------------------------31
\label{6.45}
\index{classifying spaces (example)}
\begingroup%%------------------------------------>>
\fontsize{9pt}{11pt}\selectfont
\textbf{\small EXAMPLE \ (\un{Classifying Spaces})} \ 
Let $G$ be a stratifiable topological group $-$then $X_G^\infty$ and $B_G^\infty$ are stratifiable.
\vspi
[Since the $X_G^n$ are stratifiable, the lemma can be used to establish the stratifiability of $X_G^\infty$.  
As for $B_G^\infty$, in the notation of the Milnor construction, $\restr{X_G^\infty}{O_i}$ is homeomorphic to $O_i \times G$, 
thus $O_i$ is stratifiable and so $B_G^\infty$ admits a neighborhood finite closed covering by stratifiable subspaces, hence is stratifiable.]\\
\endgroup%%------------------------------------<<

\begingroup%%------------------------------------>>
\fontsize{9pt}{11pt}\selectfont
\textbf{\small FACT} \ 
Let $X$ and $Y$ be stratifiable.  
Let $A$ be a closed subspace of $X$ and let $f:A \ra Y$ be a continuous function $-$then $X \sqcup_f Y$ is stratifiable.\\
\endgroup%%------------------------------------<<

\begingroup%%------------------------------------>>
\fontsize{9pt}{11pt}\selectfont
Application:  Suppose that $(X,A)$ is a relative CW complex.  
Assume: $A$ is stratifiable $-$then $X$ is stratifiable.\\
\endgroup%%------------------------------------<<

\begingroup%%------------------------------------>>
\fontsize{9pt}{11pt}\selectfont
Let $X$ be a topological space; 
Let $\sS$ and $\sT$ be collections of subsets of $X$ $-$then $\sS$ is said to be 
\un{cushioned}
\index{cushioned} 
in $\sT$ if there exists a function $\Gamma : \sS \ra \sT$ such that 
$\forall \  \sS_0 \subset  \sS$: $\overline{\bigcup \{S: S \in  \sS_0\}}$ 
$\subset$ 
${\bigcup \{\Gamma (S): S \in  \sS_0}\}$.  
For example, if  $\sS$ is closure preserving, then $\sS$ is cushioned in $\overline{\sS}$.
A collection $\sS$ which is the union of a countable number of subcollections $\sS_n$, each of which is cushioned in $\sT$, is said to be 
\un{$\sigma$-cushioned}
\index{sigma-cushioned, $\sigma$-cushioned} 
in $\sT$.
\\ \indent
Michael\footnote[2]{\textit{Proc. Amer. Math. Soc.} \textbf{10} (1959), 309-314.}  
has shown that a CRH space $X$ is paracompact iff every open covering of $X$ has a $\sigma$-cushioned open refinement 
(cf. p. \pageref{6.33}).   
This result can be used to prove that stratifiable spaces are paracompact.  
For suppose that $\sU = \{U\}$ is an open covering of $X$.  
Put $\sU_n = \{\text{ST}_X(n,U): U \in \sU\}$.  
Let $\sU_0 \subset \sU$ $-$then $\forall \ U \in \sU_0$, 
$\text{ST}_X(n,U) \subset$ $\text{ST}_X(n,\cup \ \sU_0) \subset$ 
$\overline{\text{ST}_X(n,\cup \ \sU_0)} \subset$ $\cup \ \sU_0$, from which 
$\overline{\cup \{\text{ST}_X(n,U) : U \in \sU_0\}} \subset$ 
$\cup \ \sU_0$, thus $\sU_n$ is cushioned in 
$\sU$ and so $\sU$ has a $\sigma$-cushioned open refinement.  
Therefore $X$ is paracompact.  
Example:  A nonmetrizable Moore space is not stratifiable (Bing 
(cf. p. \pageref{6.34})).
\\ \indent
[Note: \  Another way to argue is to show that every stratifiable space is collectionwise normal and subparacompact (cf. $\S 1$, Proposition 10 and the ensuing remark).]\\
\endgroup%%------------------------------------<<

\begingroup%%------------------------------------>>
\fontsize{9pt}{11pt}\selectfont
Let $X$ be a CRH space $-$then $X$ is said to satisfy 
\un{Arhangel'ski\u i's condition}
\index{Arhangel'ski\u i's condition} 
if there exists a sequence $\{\sU_n\}$ of collections of open subsets of $\beta X$ such that each $\sU_n$ covers $X$ and $\forall \ x \in X$: $\ds\bigcap\limits_n \text{st}(x,\sU_n) \subset X$.  
Example: Every topologically complete CRH space $X$ satisfies Arhangel'ski\u i's condition.  
In fact $X$ is a $G_\delta$ in $\beta X$, thus 
$X = \ds\bigcap\limits_1^\infty U_n$ ($U_n$ open in $\beta X$) and so we can take $\sU_n = \{U_n\}$.  
Example: Every Moore space satisfies Arhangel'ski\u i's condition.\\
\endgroup%%------------------------------------<<

\begingroup%%------------------------------------>>
\fontsize{9pt}{11pt}\selectfont
\textbf{\small FACT} \ 
Let $X$ be a CRH space.  Suppose that $X$ satisfies Arhangel'ski\u i's condition $-$then $X$ is compactly generated.\\
\endgroup%%------------------------------------<<

%%----------------------------------------------------------------------------------------------32
\begingroup%%------------------------------------>>
\fontsize{9pt}{11pt}\selectfont
Let $X$ be a CRH space $-$then 
Kullman\footnote[2]{\textit{Proc. Amer. Math. Soc.} \textbf{27}  (1971), 154-160.} 
has shown that $X$ is Moore iff $X$ is submetacompact, has a $G_\delta$ diagonal, and satisfies Arhangel'ski\u i's condition.  
Since a stratifiable space is paracompact and has a perfect square, 
it follows that every stratifiable space satisfying Arhangel'ski\u i's condition is metrizable (Bing 
(cf. p. \pageref{6.35})).  
Consequently, a nonmetrizable stratifiable space cannot be embedded in a topologically complete stratifiable space.
Example: Every stratifiable LCH space is metrizable.\\
\endgroup%%------------------------------------<<

\begingroup%%------------------------------------>>
\fontsize{9pt}{11pt}\selectfont
A Hausdorff space $X$ is said to satisfy 
\un{Ceder's condition}
\index{Ceder's condition} 
if $X$ has a $\sigma $-closure preserving basis. 
Example:  Suppose that $X$ is metrizable $-$then $X$ satisfies Ceder's condition.
Reason: The Nagata-Smirnov metrization theorem says that a regular Hausdorff space $X$ is metrizable iff 
$X$ has a $\sigma$-neighborhood finite basis.  
On the other hand, every CW complex satisfies Ceder's condition (cf. infra) and a CW complex is not in general metrizable.\\
\endgroup%%------------------------------------<<

\begingroup%%------------------------------------>>
\fontsize{9pt}{11pt}\selectfont
\textbf{\small FACT} \ 
Let $X$ be a Hausdorff space.  
Suppose that $X$ is the closed image of a metrizable space $-$then $X$ satisfies Ceder's condition.\\
\endgroup%%------------------------------------<<

\begingroup%%------------------------------------>>
\fontsize{9pt}{11pt}\selectfont
Any $X$ that satisfies Ceder's condition is stratifiable.  
Proof: Let $\sO = \bigcup\limits_n \sO_n$ be a $\sigma$-closure preserving basis for $X$, 
attach to each closed set $A \subset X$:  
$O(n,A) = X - \bigcup \{\overline{O}: O \in \sO_n \  \& \ A \cap \overline{O} = \emptyset\}$ 
and then define 
$\text{ST}_X :\N \times \tau \ra \tau$ by setting $\text{ST}_X(n,U) = X - \overline{O(n,X - U)}$.
\vspi
[Note: \  It is unknown whether the converse holds.]\\
\endgroup%%------------------------------------<<

\begingroup%%------------------------------------>>
\fontsize{9pt}{11pt}\selectfont
\index{\tM complexes}
\textbf{\small EXAMPLE \ (\un{M complexes})} \ 
A topological space is said to be an 
\un{$\tM_0$ space}
\index{space, $M_0$ space}  
if it is metrizable and, recursively, a topological space is said to be an 
\un{$\tM_{n+1}$ space}
\index{space, $\tM_{n+1}$ space} 
if it is homeomorphic to an adjunction $X \sqcup_f Y$, where $X$ is an $\tM_0$ space and $Y$ is an $\tM_n$ space.  
An \un{$\tM_\infty$ space} is a topological space that is an $\tM_n$ space for some $n$.
\\ \indent
A topological space $X$ is said to be an 
\un{M complex}
\index{M complex} 
if there exists a sequence of closed $\tM_\infty$ subspaces $A_j:$
$
\begin{cases}
\ X = \ds\bigcup\limits_j A_j \\
\ A_j \subset A_{j+1} \\
\end{cases}
$
and the topology on $X$ is the final topology determined by the inclusions $A_j \ra X$.
Example: Every CW complex is an M complex.  Since an M complex is the quotient of a metrizable space, an M complex is necessarily compactly generated.  Therefore a subspace of an M complex is an M complex iff it is compactly generated.  
Every M complex satisfies Ceder's condition, hence is stratifiable.
\vspi
[Note: \  Not every CW complex is the closed image of a metrizable space.]\\
\endgroup%%------------------------------------<<

\index{Dugundji Extension Theorem! (Stratifiable Space)}
\index{Theorem (Dugundji Extension Theorem! (Stratifiable Space))}
\begingroup%%------------------------------------>>
\fontsize{9pt}{11pt}\selectfont
\textbf{\small DUGUNDJI EXTENSION THEOREM} \ \ 
Let $X$ be a stratifiable space; let $A$ be a closed subspace of $X$.  
Let $E$ be a locally convex topological vector space.  
Equip
$
\begin{cases}
\ C(A,E) \\
\ C(X,E)
\end{cases}
$
with the compact 
%%----------------------------------------------------------------------------------------------33
open topology $-$then there exists a linear embedding $\ext : C(A,E) \ra C(X,E)$ such that $\forall \ f \in C(A,E)$, 
$\restr{\ext(f)}{A} = f$ and the range of ext(f) is contained in the convex hull of the range of $f$.
\\ \indent
[Normalize $\ST_X:$
$
\begin{cases}
\ \ST_X(n,X) = X \\
\ \ST_X(1,X - \{x\}) = \emptyset
\end{cases}
$
$\& \ \ST_X(n,U) \subset \ST_X(n+1,U)$.  
Given $x \in U$, let $n(x,U)$ be the smallest integer $n$ : $x \in \ST_X(n,U)$.  
Put
$U(x) = \ST_X(n(x,U),U) - \overline{\ST_X(n(x,U),X - \{x\})}$, a neighborhood of $x$.  
Plainly, $U(x) \cap V(y) \neq \emptyset$ $\&$ $n(x,U) \leq n(y,V)$ $\implies$ $y \in U$.  
On the other hand, 
$
\begin{cases}
\ n(x,X) = 1 \\
\ X(x) = X
\end{cases}
$
$\implies \{U: y \in U(x)\} \neq \emptyset$.  
Assuming that $A$ is nonempty and proper, attach to each 
$x \in X-A$: $n(x) = \max\{n(a,O) (O \in \tau): a \in A \ \& \ x \in O(a)\}$ $-$then $n(x) < n(x,X -A)$.  
Since every subspace of $X$ is stratifiable, $X - A$ is, in particular, paracompact.  
Thus the open covering $\{(X - A)(x): x \in X - A\}$ has a neighborhood finite open refinement $\{U_i: i \in I\}$ .  
Each $U_i$ determines a point $x_i \in X - A$: $U_i \subset (X - A)(x_i)$, from which a point $a_i \in A$ and a neighborhood $O_i$ of $a_i$: $x_i \in O_i(a_i)$ $\&$ $n(x_i) = n(a_i,O_i)$.  
Choose a partition of unity $\{\kappa_i: i \in I\}$ on $X - A$ subordinate to $\{U_i: i \in I\}$.  
Given $f \in C(A,E)$, let 
\[
\text{ext}(f)(x) =
\begin{cases}
\ f(x) \hspace{2.35cm}  (x \in A) \\
\ \ds\sum\limits_i \kappa_i(x) f(a_i) \hspace{1cm}   (x \in X - A)
\end{cases}
.
\]
Referring back to the proof of the Dugundji extension theorem in the metrizable case and eschewing the obvious, 
it is apparent that there are two nontrivial claims.
\\ \indent
Claim 1: \ ext($f$) is continuous at the points of \mA.
\\ \indent
[Let $a \in A$; let $N$ be a convex neighborhood of $f(a)$ in $E$.  
By continuity of $f$, there exists a neighborhood $O$ of $a$ in $X$:$f(A \cap O) \subset N$.  
Assertion: $\text{ext}(f)(O(a)(a)) \subset N$.  
Case 1: \ $x \in A \cap O(a)(a)$.  
Here, $x \in A \cap O$ and $\text{ext}(f)(x) = f(x) \in N$.
Case 2: \ $x \in (X - A) \cap O(a)(a)$.  
Take any index $i: \kappa_i(x) \neq 0$ $(\implies x \in U_i)$ 
$-$then $\emptyset \neq U_i \cap O(a)(a) \subset (X - A)(x_i) \cap O(a)$ $\implies$ $x_i \in O(a)$ 
$\implies$ 
$n(a,O) \leq n(x_i) =$ $n(a_i,O_i)$ $\implies$ $a_i \in O$ $\implies$ $f(a_i) \in N$ $\implies$ $\text{ext}(f)(x) \in N$.]
\\ \indent
Claim 2: \ $\text{ext} \in \text{LEO}(X,A;E)$.
\\ \indent
[Define a function $\phi:X \ra 2^A$ by the rule
$
\begin{cases}
\ \phi(a) = \{a\} \hspace{2.2cm} (a \in A)\\
\ \phi(x) = \{a_i: i \in I_x\} \hspace{1cm}  (x \in X - A)
\end{cases}
\hspace{-.25cm},
$
$I_x$ the set $\{i \in I: x \in \text{spt } \kappa_i\}$.  
Given a nonempty compact subset $K$ of $X$, put $K_A = \ds\bigcup\limits_{x \in K} \phi(x)$.
Assertion: $K_A$ is compact.  
Since the $\phi(x)$ are finite, hence compact, it will be enough to show that for every $x \in X$ 
and for every open subset $V$ of $A$ containing $\phi(x)$ there exists an open subset $U$ of $X$ containing $x$ such that 
$\cup\  \phi(U) \subset V$.  
Case 1: \ $x \in X - A$.  
Here one need only remark that there exists a neighborhood $U$ of $x$ in 
$X - A$: $y \in U$ $\implies$ $\phi(y) \subset \phi(x)$.  
Case 2: \ $a \in A$.  
Let $O$ be an open subset of $X$: 
$\phi(a) = \{a\} \subset O$.  If $x \in A \cap O(a)(a)$, then 
$\phi(x) = \{x\} \subset O$, while if $x \in (X - A) \cap O(a)(a)$, then arguing as in the first claim, 
$\forall \ i \in I_x$, $a_i \in O$.
Conclusion: $\cup \ \phi(O(a)(a)) \subset A \cap O$.]]
\\ \indent
[Note: \  Suppose that $E$ is a normed linear space 
$-$then the image of $\restr{\text{ext}}{BC(A,E)}$ is contained in $BC(X,E)$ and, per the uniform topology, 
$\text{ext}: BC(A,E) \ra BC(X,E)$ is a linear isometric embedding: $\forall \ f \in BC(A,E)$, $\norm{f} = \norm{\text{ext}(f)}$.]\\
\endgroup%%------------------------------------<<

%%----------------------------------------------------------------------------------------------34
\label{6.42}
\begingroup%%------------------------------------>>
\fontsize{9pt}{11pt}\selectfont
\textbf{\small FACT} \ 
Let $A \subset X$, where $X$ is stratifiable and $A$ is closed $-$then $A$ has the EP w.r.t. any locally convex topological space.\\
\endgroup%%------------------------------------<<

\begingroup%%------------------------------------>>
\fontsize{9pt}{11pt}\selectfont
Is it true that if $K$ is a compact Hausdorff space and $X$ is stratifiable, then $C(K,X)$ is stratifiable?  The answer is ``no'' even if $K = [0,1]$.\\
\endgroup%%------------------------------------<<

\begingroup%%------------------------------------>>
\fontsize{9pt}{11pt}\selectfont
\textbf{\small EXAMPLE} \ 
Let $X$ be the closed upper half plane in $\R^2$.  
Topologize $X$ as follows: 
The basic neighborhoods of $(x,y)$ $(y > 0)$ are as usual but the basic neighborhoods of $(x,0)$ are the ``butterflies'' 
$N_\epsilon(x)$ $(\epsilon > 0)$, where $N_\epsilon(x)$ is the point $(x,0)$ together with all points in the open upper half plane having distance $< \epsilon$ from $(x,0)$ and lying beneath the union of the two rays emanating from $(x,0)$ with slopes $\pm\epsilon$.  
Thus topologized, $X$ is stratifiable (and satisfies Ceder's condition).  
Moreover, $X$ is first countable and separable.  
But $X$ is not second countable, so $X$ is not metrizable.  
Therefore $X$ carries no CW structure (since for a CW complex, metrizability is equivalent to first countability).  
Claim: $C([0,1]),X)$ is not stratifiable.  
To see this, assign to each $r \in \R$ an element $f_r \in C([0,1]),X)$ by putting $f_r(1/2) = (r,0)$ and then laying down $[0,1]$ symmetrically around the circle of radius 1 centered at $(r,1)$.  The set $\{f_r\}$ is a closed discrete subspace of $C([0,1]),X)$ of cardinality $2^\omega$.  
Construct a closed separable subspace of $C([0,1]),X)$ containing $\{f_r\}$ and finish by quoting Jones' lemma.
\\ \indent
[Note: \  $X$ is compactly generated (being first countable).  
However, $C([0,1]),X)$ is not compactly generated.]\\
\endgroup%%------------------------------------<<

\begingroup%%------------------------------------>>
\fontsize{9pt}{11pt}\selectfont
Cauty\footnote[2]{\textit{Arch. Math. (Basel)} \textbf{27} (1976), 306-311; 
see also Guo, \textit{Tsukuba J. Math.} \textbf{18} (1994), 505-517.}
has shown that if $X$ is a CW complex, then for any compact Hausdorff space $K$, $C(K,X)$ is stratifiable, hence is perfectly normal and paracompact.\\
\endgroup%%------------------------------------<<

Let $\kappa$ be an infinte cardinal.  A Hausdorff space $X$ is said to be 
\un{$\kappa -$collectionwise normal}
\index{kappa-collectionwise normal! $\kappa-$collectionwise normal} 
if for every discrete collection $\{A_i: i \in I\}$ of closed subsets of $X$ with $\#(I) \leq \kappa$ there exists a pairwise disjoint collection $\{U_i: i \in I\}$ of open subsets of $X$ such that $\forall \ i \in I$: $A_i \subset U_i$.  
So: $X$ is collectionwise normal iff $X$ is $\kappa$-collectionwise normal for every $\kappa$.

[Note: \  Recall that every paracompact Hausdorff space is collectionwise normal (cf. $\S 1$, Proposition 9).]\\

\begingroup%%------------------------------------>>
\fontsize{9pt}{11pt}\selectfont
\textbf{\small EXAMPLE} \ 
If $X$ is normal, then $X$ is $\omega$-collectionwise normal  
(cf. p. \pageref{6.36}) 
and conversely.\\
\endgroup%%------------------------------------<<

Let $\kappa$ be an infinite cardinal; let $I$ be a set: $\#(I) = \kappa$.  
Assuming that $0 \notin I$, let $V = \{0\} \cup I$ and put $\Sigma = \{\{0\},\{i\}(i \in I)\} \cup \{\{0,i\}(i \in I)\}$ 
$-$then $K = (V, \Sigma)$ is a vertex scheme.  
Equipping $I$ with the discrete topology, one may view $\abs{K}$ as the cone
%%----------------------------------------------------------------------------------------------35
$\Gamma I$.  
Therefore $\abs{K}$ is contractible, hence so is $\abs{K}_b$ (cf. $\S$ 5, Proposition 1), 
the latter being by definition the \un{star space}
\index{star space} 
\textbf{S($\kappa$)} 
\index{\textbf{S($\kappa$)}} 
corresponding to $\kappa$.
It is clear that \textbf{S($\kappa$)} is completely metrizable of weight $\kappa$.  
The elements of \textbf{S($\kappa$)} are equivalence classes $[i,t]$ of pairs $(i,t)$, where 
$(i^\prime,t^\prime) \sim (i\pp,t\pp)$ iff $t^\prime = 0 = t\pp$ or $i^\prime = i\pp$ $\&$ $t^\prime = t\pp$.
There is a continuous map 
$
\pi_\kappa: 
\begin{cases}
\ \textbf{S($\kappa$)} \ra [0,1]\\
\ [i,t] \mapsto t
\end{cases}
$
and $\forall \ i \in I$ there is an embedding $e_i:$
$
\begin{cases}
\ [0,1] \ra \textbf{S($\kappa$)}\\
\ t \mapsto [i,t] 
\end{cases}
.
$
The point $e_i(0)$ is independent of $i$ and will be denoted by $0_\kappa$.\\

\begin{proposition} \ %13
Let $X$ be a Hausdorff space $-$then $X$ is $\kappa$-collectionwise normal iff every closed subspace $A$ of $X$ has the EP w.r.t. \bf{S($\kappa$)}
\end{proposition}

[Necessity: \ Fix an \ $f \in C(A,\textbf{S($\kappa$)})$ and let \ $\Phi:X \ra [0,1]$ \ be a continuous extension of \ 
$\pi_\kappa \circx f$.  \ 
Put \ $A_i = f^{-1}(\{[i,t]: 0 < t \leq 1\})$ : $\{A_i: i \  \in \  I\}$ \ 
is a discrete collection of closed subsets of $\Phi^{-1}(]0,1])$.  \ 
Since $\Phi^{-1}(]0,1])$ is an $F_\sigma$, it too is $\kappa$-collectionwise normal, 
thus there exists a pairwise disjoint collection \ $\{U_i$ : $i \ \in\  I\}$ \  of open subsets of $X$ such that $\forall \ i \in I$ : $A_i \subset U_i$.   
Define a continuous function  $g : A \cup (X \ -  \ \bigcup\limits_i U_i) \lra [0,1]$ \ by the conditions
$
\begin{cases}
\ \restr{g}A = \pi_\kappa \circx f\\
\ \restr{g}{X - \bigcup\limits_i U_i} = 0
\end{cases}
$
and extend it to a continuous function $G:X \ra [0,1]$.  Set $F(x) = $
$
\begin{cases}
\ e_i \circx G(x) \hspace{0.5cm} (x \in U_i)\\
\ 0_\kappa \hspace{1.67cm} (x \in X - \bigcup\limits_i U_i)
\end{cases}
$
$-$then $F \in C(X,\textbf{S($\kappa$)})$ and $\restr{F}{A} = f$.

Sufficiency: Let $\{A_i: i \in I\}$ be a discrete collection of closed subsets of $X$ with $\#(I) = \kappa$.  
Put $A = \bigcup\limits_i A_i$ $-$then $A$ is a closed subspace of $X$.  
Define $f \in C(A,\textbf{S($\kappa$)})$ piecewise: $\restr{f}{A_i} = [i,1]$.  
Extend $f$ to $F \in C(X,\textbf{S($\kappa$)})$ and consider the collection $\{U_i:i \in I\}$, where 
$U_i = F^{-1}(\{[i,t]: 1/2 < t \leq 1\})$.]\\

Application:  The star space $\bS(\kappa$) is an AR.\\

\label{6.44}

\begingroup%%------------------------------------>>
\fontsize{9pt}{11pt}\selectfont
\textbf{\small EXAMPLE} \ 
Let $\kappa$ be an infinite cardinal $-$then there exists a $\kappa$-collectionwise normal space $X$ which is not $\kappa^+$-collectionwise normal, $\kappa^+$ the cardinal successor to $\kappa$.  For this, fix a set $I^+$ of cardinality $\kappa^+$  and equip $I^+$ with the discrete topology.  
There is an embedding
$I^+ \ra \ds\prod \textbf{S($\kappa$)}$, 
the terms of the product being indexed by elements of $C(I^+,\textbf{S($\kappa$)})$.  
Let $X$ be the result of retopologizing $\ds\prod \textbf{S($\kappa$)}$ 
by isolating the points of $\ds\prod \textbf{S($\kappa$)} - I^+$.
\vspi
Claim: $X$ is $\kappa$-collectionwise normal.
\vspi
[Let $\{A_i: i \in I\}$ be a discrete collection of closed subsets of $X$ with $\#(I) = \kappa$.  
Since $X - I^+$ is discrete, there is no loss of generality in assuming that the $A_i$ are contained in $I^+$.  
Define a continuous function 
$f:\ds\bigcup\limits_i A_i \ra \textbf{S($\kappa$)}$ by $\restr{f}{A_i} = [i,1]$ 
and then, using Proposition 13, extend $f$ to an element 
$F \in C(I^+,\textbf{S($\kappa$)})$, determining a projection 
$p_F:\ds\prod \textbf{S($\kappa$)} \ra \textbf{S($\kappa$)}$ such that $\restr{p_F}{I^+} = F$.  
Consider the collection $\{U_i:i \in I\}$, where $U_i = p_F^{-1}(\{[i,t]: 1/2 < t \leq 1\})$.]
\vspi
%%----------------------------------------------------------------------------------------------36
Claim: $X$ is not $\kappa^+$-collectionwise normal.
\vspi
[If $X$ were $\kappa^+$-collectionwise normal, then it would be possible to separate the points of $I^+$ by a collection of nonempty pairwise disjoint open subsets of $X$ of cardinality $\kappa^+$.  
Taking into account how $X$ is manufactured from $\ds\prod \textbf{S($\kappa$)}$, one arrives at a contradiction to an obvious corollary of the Hewitt-Pondiczery theorem.]
\vspi
[Note: Give $I^+ \times \{0\} \cup \ds\bigcup\limits_1^\infty (X - I^+) \times \{1/n\}$ the topology induced by the product $X \times [0,1]$ $-$then this space is perfectly normal and $\kappa$-collectionwise normal but is not 
$\kappa^+$-collectionwise normal.  
And: It is not a LCH space 
(cf. p. \pageref{6.37}).]\\
\endgroup%%------------------------------------<<

\index{Kowalsky's Lemma}
\textbf{\small KOWALSKY'S LEMMA} \ 
Let $\kappa$ be an infinite cardinal.  Let $Y$ be an AR of weight $\kappa$ $-$then every metrizable space $X$ of weight $\leq \kappa$ can be embedded in $Y^\omega$.

[Let $\sU = \bigcup\limits_n \sU_n$ be a $\sigma$-discrete basis for $X$ : $\sU_n = \{U_n(i) : i \in I_n\}$, 
where 
$I = \coprod\limits_n I_n$ and $\#(I) \leq \kappa$.  
Write $\cup \ \sU_n = \bigcup\limits_m A_{mn}$, $A_{mn}$ closed in $X$.  
Fix distinct points $a, b$ which do not belong to $I$.  Since wt $Y$ $= \kappa$, there exists in $Y$ a collection of nonempty pairwise disjoint open sets $V_j$ $(j \in I \cup \{a,b\})$.  
Choose a point $y_j \in V_j$.  
Given $n$, define a continuous function 
$f_n: \cup \ \ov{\sU}_n \ra Y$ by $\restr{f_n}{\overline{U_n(i)}} = y_i$ $(i \in I_n)$ and extend $f_n$ to a continuous function
$F_n:X \ra Y$.  
Given $mn$, define a continuous function $f_{mn}:A_{mn} \cup (X - \cup \ \sU_n) \ra Y$ by 
$
\begin{cases}
\ \restr{f_{mn}}{A_{mn}} = y_a\\
\ \restr{f_{mn}}{X - \cup \ \sU_n} = y_b
\end{cases}
$
and extend $f_{mn}$ to a continuous function $F_{mn}:X \ra Y$.  Let $\Phi_{mn}:X \ra Y^2$ be the diagonal of $F_n$ and $F_{mn}$.  
Let $\Phi$ be the diagonal of the $\Phi_{mn}$, so $\Phi:X \ra (Y^2)^{\omega^2} \equiv Y^\omega$  
$-$then $\Phi$ is an embedding.]

[Note: \  Suppose that $Y$ is not compact $-$then every completely metrizable space $X$ of weight $\leq \kappa$ can be embedded in $Y^\omega$ as a closed subspace.  
For $X$, as a subspace of $Y^\omega$, is a $G_\delta$ (being completely metrizable), thus on elementary grounds is homeomorphic to a closed subspace of $Y^\omega \times \R^\omega$: 
Take a compatible metric $d$ on $Y^\omega$, represent the complement $Y^\omega - X$ as a countable union 
$\bigcup\limits_j B_j$ of closed subsets $B_j$, let $d_j:Y^\omega \ra \R$ be the function $y \ra d(y,B_j)$, and consider the graph of the diagonal of the $d_j$.  
Claim: There is a closed embedding $\R \ra Y^\omega$.  To see this, fix a closed discrete subset $\{y_n: n \in \Z\}$ in $Y$.  
Let
$
\begin{cases}
\ S = \bigcup\limits_{-\infty}^\infty [2n,2n+1]\\
\ T = \bigcup\limits_{-\infty}^\infty [2n+1,2n+2]
\end{cases}
$
and define continuous functions
$
\begin{cases}
\ f:S \ra Y\\
\ g:T \ra Y
\end{cases}
$
by 
$
\begin{cases}
\ \restr{f}{[2n,2n+1]} = y_n\\
\ \restr{g}{[2n+1,2n+2]} = y_n\\
\end{cases}
\hspace{-.25cm}. \ 
$
Extend 
$
\begin{cases}
\ f\\
\ g\\
\end{cases}
$
to a continuous function 
$
\begin{cases}
\ F:\R \ra Y\\
\ G:\R \ra Y
\end{cases}
$
and let $H:\R \ra Y^2$ be the diagonal of $F$ and $G$.  If $\Phi:\R \ra Y^\omega$ is any embedding, then the diagonal of 
$\Phi$ and $H$ is a closed embedding $\R \ra Y^\omega \times Y^2 \equiv Y^\omega$.]\\

%%----------------------------------------------------------------------------------------------37
\label{19.22}
Application: Every metrizable space $X$ of weight $\leq \kappa$ can be embedded in $\bS(\kappa)^\omega$.\\

\label{19.52}
Let $\kappa$ be an infinite cardinal.  Let $X$ be a topological space 
$-$then a subspace $A \subset X$ is said to have the 
\un{extension property with respect to $\sB(\kappa)$} 
\index{extension property with respect to $\sB(\kappa)$} 
\index{EP w.r.t. $\sB(\kappa)$}
(EP w.r.t. $\sB(\kappa)$) if it has the EP w.r.t. every Banach space of weight $\leq \kappa$.  
Since every completely metrizable AR can be realized as a closed subspace of a Banach space 
(cf. p. \pageref{6.38}), 
it is clear that $A$ has the EP w.r.t. $\sB(\kappa)$ iff it has the EP w.r.t every completely metrizable AR of weight $\leq \kappa$.\\

\begin{proposition} \ %14
Fix a pair $(X,A)$.  
Suppose that for some noncompact AR $Y$ of weight $\kappa$, $A$ has the EP w.r.t $Y$ $-$then $A$ has the EP w.r.t $\sB(\kappa)$ .
\end{proposition}

[Let $E$ be a Banach space of weight $\leq \kappa$.  
Owing to Kowalsky's lemma, $E$ can be realized as a closed subspace of $Y^\omega$.  Let $f \in C(A,E)$.  
By hypothesis, $f$ has a continuous extension $F \in C(X,Y^\omega)$.  Consider $r \circx F$, where $r:Y^\omega \ra E$ is a retraction.]\\

One conclusion that can be drawn from this is that $A$ has the EP w.r.t. $\R$ iff $A$ has the EP w.r.t $\sB(\omega)$.  So: If $X$ is a Hausdorff space, then $X$ is normal iff every closed subspace $A$ of $X$ has the EP w.r.t every separable Banach space.

\label{5.0}
\label{20.2}
Another conclusion is that $A$ has the EP w.r.t $\mathbf{S}(\kappa)$ iff $A$ has the EP w.r.t $\sB(\kappa)$.  
Consequently, if $X$ is a Hausdorff space, then $X$ is $\kappa$-collectionwise normal iff every closed subspace $A$ of $X$ has the EP w.r.t $\sB(\kappa)$. (cf. Proposition 13).  
Corollary: A Hausdorff space $X$ is collectionwise normal iff every closed subspace $A$ of $X$ has the EP w.r.t every Banach space.\\

\begingroup%%------------------------------------>>
\fontsize{9pt}{11pt}\selectfont
\textbf{\small FACT} \ 
Let $A \subset X$ $-$then $A$ has the EP w.r.t $\R$ iff $IA \subset IX$ has the EP w.r.t $[0,1]$.\\
\endgroup%%------------------------------------<<

\label{19.46}
Let $X$ be a topological space.  Let $\{\sU_n\}$ be a sequence of open coverings of $X$ $-$then $\{\sU_n\}$ is said to be a 
\un{star sequence}
\index{star sequence} 
if $\forall \ n$, $\sU_{n+1}$ is a star refinement of $\sU_n$.  
By means of a standard construction from metrization theory, 
one can associate with a given star sequence $\{\sU_n\}$ a continuous pseudometric $\delta$ on $X$ such that $\delta(x,y) = 0$ iff 
$y \in \bigcap\limits_1^\infty \text{st}(x,\sU_n)$, a subset $U \subset X$ being open in the topology generated by $\delta$ iff $\forall \ x \in U$ $\exists \ n$: $\text{st}(x,\sU_n) \subset U$.  
Let $X_\delta$ be the metric space obtained from $X$ by identifying points at zero distance from one another and write $p:X \ra X_\delta$ for the projection.\\

\begin{proposition} \ %15
Let $A \subset X$ $-$then $A$ has the EP w.r.t $\sB(\kappa)$ iff for every numerable open covering $\sO$ of $A$ of cardinality $\leq \kappa$ there exists a numerable open covering $\sU$ of $X$ of cardinality $\leq \kappa$ such that $\sU \cap A$ is a refinement of $\sO$.
\end{proposition}

%%----------------------------------------------------------------------------------------------38
[Necessity: Let $\sO = \{O_i: i \in I\}$ be a numerable open covering of $A$ with $\#(I) \leq \kappa$.
Choose a partition of unity $\{\kappa_i: i \in I\}$ on $A$ subordinate to $\sO$.  
Form the Banach space
$\ell^1(I): r = (r_i) \in \ell^1(I)$ iff $\norm{r} = \sum\limits_i \abs{r_i} < \infty$.  
The assignment
$
\begin{cases}
\ A \ra \ell^1(I)\\
\ a \mapsto (\kappa_i(a))
\end{cases}
$
defines a continuous function $f$ whose range is contained in 
$S^+ = \{r: \norm{r} = 1\} \cap \{r: \forall \ i, r_i \geq 0\}$, a closed convex subset of $\ell^1(I)$.  
Therefore $f$ has a continuous extension $F:X \ra S^+$.  Let $p_i$ be the projection
$
\begin{cases}
\ \ell^1(I) \ra \R\\
\ r \ra r_i
\end{cases}
\hspace{-.25cm}; \ 
$
let $\sigma_i = p_i \circx F$ $-$then $\restr{\sigma_i}{A} = \kappa_i$ and 
$\ds\sum\limits_i \sigma_i(x) = 1$ $(\forall \ x \in X)$.  
Put $U_i = \sigma_i^{-1}(]0,1])$ and apply NU 
(cf. p. \pageref{6.39}) 
to see that the collection $\sU = \{U_i: i \in I\}$ is a numerable open covering of $X$ of cardinality $\leq \kappa$.  
And by construction, $\sU \cap A$ is a refinement of $\sO$.

Sufficiency: Let $E$ be a Banach space of weight $\leq \kappa$.  
Fix a dense subset $E_0$ in $E$ of cardinality $\leq \kappa$ and let $\sE_n$ be the open covering of $E$ consisting of the open balls of radius $1/3^n$ centered at the points of $E_0$.  
Suppose that $f:A \ra E$ is continuous 
$-$then $\forall \ n$, $f^{-1}(\sE_n)$ is a numerable open covering of $A$ of cardinality $\leq \kappa$, 
so there exists a star sequence $\{\sU_n\}$ of open coverings of $X$ of cardinality $\leq \kappa$ such that $\forall \ n$, $\sU_n \cap A$ is a refinement of $f^{-1}(\sE_n)$.
Viewed as a map from $A$ endowed with the topology induced by the pseudometric $\delta$ associated with $\{\sU_n\}$, $f$ is continuous, 
thus passes to the quotient to give a continuous $f_\delta:A_\delta \ra E$, where $A_\delta = p(A)$.  
Because $f_\delta$ is actually uniformly continuous, 
there exists a continuous extension 
$\overline{f}_\delta:\overline{A}_\delta \ra E$ of $f_\delta$ 
to the closure 
$\overline{A}_\delta$ of $A_\delta$ in $X_\delta$.  
Choose $F_\delta \in C(X_\delta,E)$ : $\restr{F_\delta}{\overline{A}_\delta} = \overline{f}_\delta$ and consider $F = F_\delta \circx p$.]\\

Examples: Let $X$ be a CRH space $-$then $\forall \ \kappa$ 
(1) Every compact subspace of $X$ has the EP w.r.t $\sB(\kappa)$; 
(2) Every pseudocompact subspace of $X$ which has the EP w.r.t $[0,1]$ has the EP w.r.t $\sB(\kappa)$; 
(3) Every Lindel\"of subspace of $X$ which has the EP w.r.t $\R$ has the EP w.r.t $\sB(\kappa)$.\\

Suppose that $X$ is collectionwise normal.  
Let $A$ be a closed subspace of $X$; let $\sO = \{O_i: i \in I\}$ be a neighborhood finite open covering of $A$  
$-$then Proposition 15 implies that there exists a neighborhood finite open covering 
$\sU = \{U_i: i \in I\}$ of $X$ such that $\forall \ i \in I$, $U_i \cap A \subset O_i$.  
Question: Is it possible to arrange matters so that $\forall \ i \in I$, $U_i  \cap  A = O_i$?  
The answer is ``no'' since Rudin's Dowker space fails to admit this improvement 
(Przymusi\' nski-Wage\footnote[2]{\textit{Fund. Math.} \textbf{109} (1980), 175-187.\vspace{0.11cm}}) 
but ``yes'' if $X$ is in addition countably paracompact. 
(Kat\u etov\footnote[3]{\textit{Colloq. Math.} \textbf{6} (1958), 145-151.\vspace{0.11cm}}).\\

%%----------------------------------------------------------------------------------------------39

\begingroup%%------------------------------------>>
\fontsize{9pt}{11pt}\selectfont
Let $(X,\delta)$ be a pseudometric space; 
let $A$ be a closed subspace of $X$ $-$then $A$ has the EP w.r.t every AR $Y$.  
Proof: Let $X_\delta$ be the metric space obtained from $X$ by identifying points at zero distance from one another, 
write $p$ for the projection $X \ra X_\delta$, and put $A_\delta = p(A)$, a closed subspace of $X_\delta$.  
Each $f \in C(A,Y)$ passes to the quotient to give an $f_\delta \in C(A_\delta,Y)$ for which there exists an extension $F_\delta \in C(X_\delta,Y)$.  
Consider $F = F_\delta \circx p$.\\
\endgroup%%------------------------------------<<

\begingroup%%------------------------------------>>
\fontsize{9pt}{11pt}\selectfont
The weight of a pseudometric is the weight of its associated topology.\\
\endgroup%%------------------------------------<<

\begingroup%%------------------------------------>>
\fontsize{9pt}{11pt}\selectfont
\textbf{\small LEMMA} \ 
Let $A \subset X$ $-$then $A$ has the EP w.r.t $\sB(\kappa)$ iff every continuous pseudometric on $A$ of weight $\leq \kappa$ can be extended to a continuous pseudometric on $X$.
\vspi
[Necessity: Let $\delta$ be a continuous pseudometric on $A$ of weight $\leq \kappa$.  Let $A_\delta$ 
be the metric space obtained from $A$ by identifying points at zero distance from one another.  
Embed $A_\delta$ isometrically into a Banach space $E$ of weight $\leq \kappa$ 
$-$then the projection $A \ra A_\delta \subset E$ has a continuous extension $\Phi:X \ra E$ and the assignment $\Delta:$
$
\begin{cases}
\ X \times X \ra \R\\
\ (x^\prime,x^{\prime\prime}) \ra \norm{\Phi(x^{\prime}) - \Phi(x^{\prime\prime})}
\end{cases}
$
is a continuous extension of $\delta$.
\vspi
Sufficiency: Let $E$ be a Banach space of weight $\leq \kappa$; let $f \in C(A,E)$.  
Define a pseudometric $\delta$ on $A$ by 
$\delta(a^\prime,a^{\prime\prime}) = \norm{f(a^\prime) - f(a^{\prime\prime})}$ $-$then $\delta$ is continuous of weight $\leq \kappa$, hence admits a continuous extension $\Delta$.  Call $X(\Delta)$ the set $X$ equipped with the topology determined by $\Delta$.  
Let $A(\Delta)$ be the closure of $A$ in $X(\Delta)$.  
Extend $f$ continuously to a function $f(\Delta):A(\Delta) \ra E$ and note that $A(\Delta) \subset X(\Delta)$ has the EP w.r.t. $E$.]\\
\endgroup%%------------------------------------<<

\begingroup%%------------------------------------>>
\fontsize{9pt}{11pt}\selectfont
\textbf{\small FACT} \ 
Let $A$ be a zero set in $X$.  Suppose that $A$ has the EP w.r.t $\sB(\kappa)$ $-$then $A$ has the EP w.r.t every AR $Y$ of weight $\leq \kappa$.
\vspi
[Choose a $\phi \in C(X,[0,1])$: $A = \phi^{-1}(0)$.  
Fix a compatible metric $d$ on $Y$.  
Given $f \in C(A,Y)$, define a pseudometric $\delta$ on $A$ by $\delta(a^\prime,a^{\prime\prime}) = d(f(a^\prime),f(a^{\prime\prime}))$.  Let $\Delta$ be a continuous extension of $\delta$ to $X$ and consider the sum of $\Delta(x^\prime,x^{\prime\prime}) $ and $\abs{\phi(x^{\prime}) - \phi(x^{\prime\prime})}$.]\\
\endgroup%%------------------------------------<<

\begingroup%%------------------------------------>>
\fontsize{9pt}{11pt}\selectfont
Let $X$ be a CRH space.  Suppose that $X$ is perfectly normal and collectionwise normal $-$then it follows that every closed subspace $A$ of $X$ has the EP w.r.t every AR.\\
\endgroup%%------------------------------------<<

\begingroup%%------------------------------------>>
\fontsize{9pt}{11pt}\selectfont
\textbf{\small FACT} \ 
Let $X$ be a submetrizable CRH space.  
Suppose that $A \subset X$ has the EP with respect to every normed linear space $-$then $A$ is a zero set in $X$.
\vspi
[Note: \  Take for $X$ the Michael line and let $A = \Q$ $-$then $X$ is a paracompact Hausdorff space, so $A$ has the EP w.r.t every Banach space.  
On the other hand, $X$ is submetrizable but $A$ is not a $G_\delta$.  
Therefore $A$ does not have the EP w.r.t. every normed linear space.]\\
\endgroup%%------------------------------------<<

\textbf{\small LEMMA} \ 
Fix a pair $(X,A)$.  Suppose that $A$ has the EP w.r.t $\sB(\kappa)$ $-$then every continuous function 
$\phi:i_0X \cup IA \ra \bS(\kappa)$ has a continuous extension $\Phi:IX \ra \bS(\kappa)$.

%%----------------------------------------------------------------------------------------------40
[The restriction \  $\psi$ of $\phi$ to $IA$ determines a continuous function $A \ra C([0,1],\bS(\kappa))$. \ 
But $ C([0,1],\bS(\kappa))$ is a completely metrizable AR (cf. the proof of Proposition 6), 
the weight of which is $\leq \kappa$, so our assumption on $A$ guarantees that this function has a continuous extension $X \ra  C([0,1],\bS(\kappa))$, 
leading thereby to a continuous function $\Psi:IX \ra \bS(\kappa)$ whose restriction to $IA$ is $\psi$.  
Choose an $f \in C(X,[0,1])$ : $f^{-1}(0) = \{x: \phi(x,0) = \Psi(x,0)\}$.  
Let $F$ be the function
$
\begin{cases}
\ X \ra \bS(\kappa)\\
\ x \ra \Psi(x,f(x))
\end{cases}
\hspace{-.25cm}. \ 
$
Because $\bS(\kappa)$ is contractible, there is a homotopy $H:IX \ra \bS(\kappa)$ such that
$
\begin{cases}
\ H(x,0) = \phi(x,0)\\
\ H(x,1) = F(x)
\end{cases}
\hspace{-.25cm}. \ 
$
Consider the function $\Phi:IX \ra \bS(\kappa)$ defined by $\Phi(x,t) = $
$
\begin{cases}
\ \Psi(x,t) \hspace{1.45cm}(t \geq f(x))\\
\ H(x,t/f(x)) \hspace{0.5cm}  (t < f(x))
\end{cases}
\hspace{-.25cm}.]
$
\\
\vspace{0.25cm}

\begin{proposition} \ 
Let $A \subset X$ $-$then $A$ has the EP w.r.t $\sB(\kappa)$ iff $i_0X \cup IA$, as a subspace of $IA$, has the EP w.r.t every completely metrizable ANR $Y$ of weight $\leq \kappa$.
\end{proposition}

[Necessity: Let $f:i_0X \cup IA \ra Y$ be continuous.  Using Kowalsky's lemma, realize $Y$ as a closed subspace of $\bS(\kappa)^\omega$ and let $r:O \ra Y$ be a retraction ($O$ open in $\bS(\kappa)^\omega$).  
Given a projection $p:\bS(\kappa)^\omega \ra \bS(\kappa)$, let $\phi_p = p \circx f$ $-$then by what has been said above, $\phi_p$ has a continuous extension $\Phi_p:IX \ra \bS(\kappa)$.  
Therefore $f$ has a continuous extension $\Phi :IX \ra \bS(\kappa)^\omega$.  Set $P = \Phi^{-1}(O)$.  
Since $P$ is a cozero set in $IX$ containing $IA$ and since the projection $IX \ra X$ takes zero sets to zero sets, there is a cozero set $U$ in $X$ such that $A \subset U$ and $IU \subset P$.  
On the other hand, $A$ has the EP w.r.t. $\R$, so it follows from Proposition 3 that $\exists \ \phi \in C(X,[0,1])$ :
$
\begin{cases}
\ \restr{\phi}{A} = 1\\
\ \restr{\phi}{X - U} = 0
\end{cases}
\hspace{-.25cm}. \ 
$
Define $F \in C(IX,Y)$ by $F(x,t) = r(\Phi(x,\phi(x)t))$ : $F$ is a continuous extension of $f$.

Sufficiency: Let $\sO = \{O_i: i \in I\}$ be a neighborhood finite cozero covering of $A$ with $\#(I) \leq \kappa$.  Put
\[
\sP = \{O_i \times ]1/3,1]: i \in I\} \cup \{i_0X \cup A \times [0,2/3[\}.
\]
Then $\sP$ is a neighborhod finite cozero set covering of $i_0X \cup IA$ of cardinality $\leq \kappa$, thus Proposition 15 implies that there exists a numerable open covering $\sV$ of $IX$ of cardinality $\leq \kappa$ such that $\sV \cap (i_0X \cup IA)$ is a refinement of $\sP$.
Let $\sU = \sV \cap (i_1X)$: $\sU$ is a numerable open covering of $i_1X$ such that 
$\sU \cap (i_1A)$ is a refinement of $\sP \cap (i_1A) = i_1 \sO$.  Finish by quoting Proposition 15.]\\

\label{19.45}
\begingroup%%------------------------------------>>
\fontsize{9pt}{11pt}\selectfont
\textbf{\small EXAMPLE} \ 
Suppose that the inclusion $A \ra X$ is a cofibration $-$then $i_0 X \cup IA$ is a retract of $IX$ (cf. $\S 3$, Proposition 1), so Proposition 16 implies that $A$ has the EP w.r.t. every Banach space.
\vspi
[Note: \  This applies in particular to a relative CW complex $(X,A)$.]\\
\endgroup%%------------------------------------<<

\label{6.22}
\label{19.10}
\label{19.40}
\label{19.51}
Let $X$ and $Y$ be topological spaces.

%%----------------------------------------------------------------------------------------------41
\indent\indent (HEP) \ A subspace $A \subset X$ is said to have the 
\un{homotopy extension property with} \un{respect to $Y$} 
\index{homotopy extension property with respect to $Y$} 
\index{HEP}
(HEP w.r.t $Y$) if given continuous functions
$
\begin{cases}
\ F:X \ra Y\\
\ h:IA \ra Y
\end{cases}
$
such that $\restr{F}{A} = h \circx i_0$, there is a continuous function $H:IX \ra Y$ such that $F = H \circx i_0$ and 
$\restr{H}{IA} = h$.

[Note: \  In this terminology, the inclusion $A \ra X$ is a cofibration iff $A$ has the HEP w.r.t $Y$ for every $Y$.]

Suppose that $A$ has the HEP w.r.t $Y$.  Let 
$
\begin{cases}
\ f \in C(A,Y)\\
\ g \in C(A,Y)
\end{cases}
$
be homotopic.  
Assume: $f$ has a continuous extension $F \in C(X,Y)$ $-$then $g$ has a continuous extension $G \in C(X,Y)$ and $F \simeq G$.  
Therefore, under these circumstances, the extension question for continuous functions $A \ra Y$ is a problem in the homotopy category.

If $A \subset X$ is closed and if $i_0X \cup IA$, as a subspace of $IX$, has the EP w.r.t $Y$, then it is clear that $A$ has the HEP w.r.t. $Y$.  
Conditions ensuring that this is so are provided by Proposition 16.  
Here are two illustrations.

\indent\indent (1) Every closed subspace $A$ of a normal Hausdorff space $X$ has the HEP w.r.t. every second countable completely metrizable ANR $Y$.

\indent\indent (2) Every closed subspace $A$ of a collectionwise normal Hausdorff space $X$ has the HEP w.r.t. every completely metrizable ANR $Y$.

[Note: \  Historically, these results were obtained by imposing in addition a countable paracompactness assumption on $X$.  
Reason: If $X$ is a normal Hausdorff space, then the product $IX$ is normal iff $X$ is countably paracompact.]

If $A \subset X$ and if $A$ has the EP w.r.t $\sB(\kappa)$, then $A$ has the HEP w.r.t every completely metrizable ANR $Y$ of weight $ \leq \kappa$. 
Proof: Take a pair of continuous functions
$
\begin{cases}
\ F:X \ra Y\\
\ h:IA \ra Y
\end{cases}
$
such that $\restr{F}{A} = h \circx i_0$ and define $\phi:i_0X \cup IA \ra Y$ by
$
\begin{cases}
\ \phi(x,0) = F(x)\\
\ \phi(a,t) = h(a,t)
\end{cases}
\hspace{-.25cm}.
$
In view of Proposition 16, the only issue is the continuity of $\phi$.  
To see this, embed $Y$ in a Banach space $E$ of weight $\leq \kappa$.  
Since $IA$, as a subspace of $IX$, has the EP w.r.t. $\sB(\kappa)$, $h$ has a continuous extension 
$\overline{h}: I\overline{A} \ra E$.  Define $\overline{\phi}:i_0 X \cup I\overline{A} \ra E$ by
$
\begin{cases}
\ \overline{\phi}(x,0) = F(x)\\
\ \overline{\phi}(\overline{a},t) = \overline{h}(\overline{a},t)
\end{cases}
$
$-$then $\overline{\phi}$ is a welldefined continuous function which agrees with $\phi$ on $i_0 X \cup IA$.\\

\begingroup%%------------------------------------>>
\fontsize{9pt}{11pt}\selectfont
\textbf{\small EXAMPLE} \ 
The product $Y = \bS^n \times \bS^n \times \cdots$ ($\omega$ factors) is not an ANR.  
But if $X$ is normal and $A \subset X$ is closed, then $A$ has the HEP w.r.t. $Y$.\\
\endgroup%%------------------------------------<<

\begingroup%%------------------------------------>>
\fontsize{9pt}{11pt}\selectfont
\textbf{\small FACT} \ 
Suppose that $X$ is Hausdorff.  Let $A$ be a zero set in $X$.
\\
\indent\indent (1) If $X$ is normal, then $A$ has the HEP w.r.t. every second countable ANR $Y$.
\\
\indent\indent (2) If $X$ is collectionwise normal, then $A$ has the HEP w.r.t. every ANR $Y$.\\
\endgroup%%------------------------------------<<

%%----------------------------------------------------------------------------------------------42
\begingroup%%------------------------------------>>
\fontsize{9pt}{11pt}\selectfont
\textbf{\small FACT} \ 
Let $Y$ be a nonempty metrizable space.  
Suppose that $Y$ is locally contractible $-$then $Y$ is an ANR iff for every pair $(X,A)$, 
where $X$ is metrizable and $A \subset X$ is closed, $A$ has the HEP w.r.t. $Y$.\\
\endgroup%%------------------------------------<<

Let $\sX$ be a homeomorphism invariant class of normal Hausdorff spaces that is closed hereditary, 
i.e., if $X \in\sX$  and if $A \subset X$ is closed, then $A \in \sX$ .\\

\begingroup%%------------------------------------>>
\fontsize{9pt}{11pt}\selectfont
Let $\sX$ be the class consisting of the Hausdorff spaces satisfying Ceder's condition $-$then it is unknown whether $\sX$  is closed hereditary.\\
\endgroup%%------------------------------------<<

A nonempty topological space $Y$ is said to be an 
\un{extension space}
\index{extension space} 
for $\sX$ if every closed subspace of every element of $\sX$ has the EP w.r.t $Y$.  
Denote by ES($\sX$) 
\index{ES($\sX$)} 
the class of extension spaces for $\sX$.  
Obviously, if $\sX^\prime \subset \sX\pp$, then 
$\text{ES}(\sX\pp) \subset \text{ES}(\sX^\prime)$, so $\forall \ \sX$: ES(normal)  $\subset \text{ES}(\sX)$.

\indent\indent (ES$_1$) The class $\text{ES}(\sX)$ is closed under the formation of products.

\indent\indent (ES$_2$) Any retract of an extension space for $\sX$ is in $\text{ES}(\sX)$.

\indent\indent (ES$_3$) Suppose that $Y = Y_1 \cup \hsx Y_2$, where $Y_1$, and $Y_2$ are open and 
$
\begin{cases}
\ Y_1\\
\ Y_2
\end{cases}
\hspace{-.25cm}
$
$\in \text{ES}(\sX)$ $\&$ $Y_1 \cap Y_2 \in \text{ES}(\sX)$ 
$-$then $Y \in \text{ES}(\sX)$.

\indent\indent (ES$_4$) Assume: The elements of $\sX$ are hereditarily normal.  Suppose that $Y = Y_1 \cup Y_2$, where $Y_1$, and $Y_2$ are closed and 
$
\begin{cases}
\ Y_1\\
\ Y_2
\end{cases}
\hspace{-.25cm}
\in \text{ES}(\sX) \ \& \ Y_1 \cap Y_2 \in \text{ES}(\sX) 
$
$-$then $Y \in \text{ES}(\sX)$.

\indent\indent (ES$_5$) Suppose that $Y = Y_1 \cup Y_2$, where $Y_1$ and $Y_2$ are closed $-$then $Y \in \text{ES}(\sX)$ $\&$ $Y_1 \cap Y_2  \in \text{ES}(\sX)$ $\implies$ 
$
\begin{cases}
\ Y_1\\
\ Y_2
\end{cases}
\hspace{-.25cm}
\in \text{ES}(\sX)$.\\

\begingroup%%------------------------------------>>
\fontsize{9pt}{11pt}\selectfont
\textbf{\small EXAMPLE} \ 
A nonempty topological space $Y$ is an extension space for the class of metrizable spaces iff it is an extension space for the class of M complexes.\\
\endgroup%%------------------------------------<<

A nonempty topological space $Y$ is said to be a 
\un{neighborhood extension space}
\index{neighborhood extension space} 
for $\sX$ if every closed subspace of every element of $\sX$ has the NEP w.r.t. $Y$.  
Denote by $\text{NES}(\sX)$ 
\index{NES$(\sX)$} 
the class of neighborhood extension spaces for $\sX$.
Obviously, if $\sX^\prime \subset \sX\pp$, then 
$\text{NES}(\sX\pp) \subset \text{NES}(\sX^\prime)$, 
so $\forall \ \sX$ : NES(normal)  $\subset \text{NES}(\sX)$. \ 
Of course, $\text{ES}(\sX) \subset \text{NES}(\sX)$.  
In the other direction, every contractible element of $\text{NES}(\sX)$ is in $\text{ES}(\sX)$.

[Note: \  It is convenient to agree that $\emptyset \in \text{NES}(\sX)$.  So, if $Y \in \text{NES}(\sX)$ and if $V \subset Y$ is open, then $V \in \text{NES}(\sX)$.]

\indent\indent (NES$_1$) The class $\text{NES}(\sX)$ is closed under the formation of finite products.

\indent\indent (NES$_2$) Any neighborhood retract of a neighborhood extension space for $\sX$ is in $\text{NES}(\sX)$.
 
%%----------------------------------------------------------------------------------------------43
\indent\indent (NES$_3$) \label{6.25} Suppose that $Y = Y_1 \cup Y_2$, where $Y_1$, and $Y_2$ are open and 
$
\begin{cases}
\ Y_1\\
\ Y_2
\end{cases}
$
$\in \text{NES}(\sX)$
$-$then $Y \in \text{NES}(\sX)$.

\indent\indent (NES$_4$) 
\label{6.7a}
Assume: The elements of $\sX$ are hereditarily normal.  Suppose that $Y = Y_1 \cup Y_2$, where $Y_1$ and $Y_2$ are closed and 
$
\begin{cases}
\ Y_1\\
\ Y_2
\end{cases}
\hspace{-.25cm}
\in \text{NES}(\sX) \ \& \  Y_1 \cap Y_2 \in \text{NES}(\sX) 
$
$-$then $Y \in \text{NES}(\sX)$.

\indent\indent (NES$_5$)  Suppose that $Y = Y_1 \cup Y_2$, where $Y_1$, and $Y_2$ are closed $-$then $Y \in \text{NES}(\sX)$ $\&$ $Y_1 \cap Y_2  \in \text{NES}(\sX)$ $\implies$ 
$
\begin{cases}
\ Y_1\\
\ Y_2
\end{cases}
\hspace{-.25cm}
\in \text{NES}(\sX)$.

[Note: \  There is a slight difference between the formulation of ES$_3$ and  NES$_3$.  Reason: An empty intersection is permitted in NES$_3$ but not in ES$_3$ (consider $X = [0,1]$, $A = Y = \{0,1\}$).]\\

\begingroup%%------------------------------------>>
\fontsize{9pt}{11pt}\selectfont
\label{20.7}
\textbf{\small EXAMPLE \  (\un{CW Complexes})} \ 
Metrizable CW complexes are ANRs 
(cf. p. \pageref{6.40}).
\\
\indent\indent (1) Every finite CW complex is in NES(normal).
\\
\indent\indent (2) Every CW complex is in NES(compact) (but it is not true that every CW complex is in NES(paracompact)).
\\
\label{3.7}
\indent\indent (3) Every CW complex is in NES(stratifiable).
\vspi
[First, if $K$ is a full vertex scheme, then $\abs{K}$ is a locally convex topological vector space 
(cf. p. \pageref{6.41}), 
so $\abs{K} \in $ ES(stratifiable) 
(cf. p. \pageref{6.42}).  
Second, if $K$ is a vertex scheme and if $L$ is a subscheme, then $\abs{L}$ is a neighborhood retract of $\abs{K}$.  
Third, if $X$ is a CW complex, then $X$ is a retract of a polyhedron 
(cf. p, \pageref{6.43}).]\\
\endgroup%%------------------------------------<<

\label{20.5}
\begingroup%%------------------------------------>>
\fontsize{9pt}{11pt}\selectfont
\textbf{\small FACT} \ 
Every CW complex has the homotopy type of an ANR which is in NES(paracompact).\\
\endgroup%%------------------------------------<<

\label{19.14}
\begingroup%%------------------------------------>>
\fontsize{9pt}{11pt}\selectfont
\textbf{\small EXAMPLE} \ \ 
Suppose that $X = Y \cup Z$ is metrizable.  \ Let $K$ and $L$ be finite CW complexes.  \ 
Assume: Every closed subspace of 
$
\begin{cases}
\ Y\\
\ Z
\end{cases}
$
has the EP w.r.t
$
\begin{cases}
\ K\\
\ L
\end{cases}
$ 
$-$then every closed subspace of $X$ has the EP w.r.t $K*L$.\\
\endgroup%%------------------------------------<<

The ``ES'' arguments are similar to but simpler than the ``NES'' arguments.  
Of the latter, the most difficult is the one for NES$_3$, which runs as follows.  
Take an $X$ in $\sX$ and let $A \subset X$ be closed 
$-$then the claim is that $\forall \ f \in C(A,Y)$ there exists an open $U \supset A$ and an $F \in C(U,Y)$: $\restr{F}{A} = f$.  
Since $X$ is covered by open sets
$
\begin{cases}
\ f^{-1}(Y_1) \cup (X - A)\\
\ f^{-1}(Y_2) \cup (X - A)
\end{cases}
$
and since $X$ is normal, there exist closed sets
$
\begin{cases}
\ X_1 \subset X\\
\ X_2 \subset X
\end{cases}
$
which cover $X$ with 
$
\begin{cases}
\ X_1 \subset f^{-1}(Y_1) \cup (X - A)\\
\ X_2 \subset f^{-1}(Y_2) \cup (X - A)
\end{cases}
\hspace{-.25cm}
. \ 
$
Put
$
\begin{cases}
\ A_1 = X_1 \cap A\\
\ A_2 = X_2 \cap A
\end{cases}
\hspace{-.25cm}
. \ 
$
There are now two cases, depending on whether $Y_1 \cap Y_2$ is empty or not.  The second possibility is more involved than
%%----------------------------------------------------------------------------------------------44
the first so we shall look only at it.  
Because  $Y_1 \cap Y_2 \in \text{NES}(\sX)$, the restriction $\restr{f}{A_1 \cap A_2}$ has an extension 
$f_{12} \in C(O,Y_1 \cap Y_2)$, where $O$ is some open subset of $X_1 \cap X_2$ containing $A_1 \cap A_2$. 
Choose an open subset $P$ of $X_1 \cap X_2$ : $A_1 \cap A_2 \subset P \subset$ $\overline{P} \subset O$.  
Observing that $A \cap \overline{P} = A_1 \cap A_2$, define $g \in C(A \cup \overline{P},Y)$ by g(x) = 
$
\begin{cases}
\ f(x) \quadx \  (x \in A)\\
\ f_{12}(x) \ \ (x \in \overline{P}
\end{cases}
\hspace{-.25cm}
.
$
Because
$
\begin{cases}
\ Y_1 \in \text{NES}(\sX) \\
\ Y_2 \in \text{NES}(\sX)
\end{cases}
\hspace{-.25cm}
,
$
the restriction 
$
\begin{cases}
\ \restr{g}{A_1 \cup \overline{P}}\\
\ \restr{g}{A_2 \cup \overline{P}}
\end{cases}
$
has an extension 
$
\begin{cases}
\ G_1 \in C(O_1,Y_1)\\
\ G_2 \in C(O_2,Y_2)
\end{cases}
\hspace{-.25cm}
, \ 
$
where 
$
\begin{cases}
\ O_1\\
\ O_2
\end{cases}
$
is some open subset of 
$
\begin{cases}
\ X_1\\
\ X_2
\end{cases}
$
containing
$
\begin{cases}
\ A_1 \cup \overline{P}\\
\ A_2 \cup \overline{P}
\end{cases}
\hspace{-.25cm}
. \ 
$
Choose an open subset
$
\begin{cases}
\ P_1 \text{ of } X_1\\
\ P_2 \text{ of } X_2
\end{cases}
\hspace{-.25cm}
$
:
$
\begin{cases}
\ A_1 \cup \overline{P} \subset P_1 \subset \overline{P}_1 \subset O_1\\
\ A_2 \cup \overline{P} \subset P_2 \subset \overline{P}_2 \subset O_2
\end{cases}
$
and an open subset $V \subset X$: $A \subset V$ $\&$ $(X_1 \cap X_2 - P) \cap \overline{V} = \emptyset$.
Let
$
\begin{cases}
\ B_1 = (\overline{P_1 - X_2} \cap \overline{V}) \cup \overline{P}\\
\ B_2 = (\overline{P_2 - X_1} \cap \overline{V}) \cup \overline{P}
\end{cases}
\hspace{-.25cm}
. \ 
$
It is clear that
$
\begin{cases}
\ B_1 \subset O_1\\
\ B_2 \subset O_2
\end{cases}
\hspace{-.25cm}
, \ 
$
with $B_1 \cap B_2 = \overline{P}$, so the prescription 
$
G(x) = 
\begin{cases}
\ G_1(x) \ (x \in B_1)\\
\ G_2(x) \ (x \in B_2)
\end{cases}
$
is a continuous extension of $f$ to $B_1 \cup B_2 \supset A$.  The set $(P_1 - X_2) \cup (P_2 - X_1) \cup P$ is open in $X$.  
Denote by $U$ its intersection with $V$ and let $F = \restr{G}{U}$.

[Note: \ \  To reduce  \ NES$_4$ to NES$_3$,  \ put instead
$
\begin{cases}
\ A_1 = f^{-1}(Y_1)\\
\ A_2 = f^{-1}(Y_2)
\end{cases}
\hspace{-.25cm}
. \ 
$
\qquad
Since
$
\begin{cases}
\ \overline{A_1 - A_2} \cap (A_2 - A_1) = \emptyset\\
\ (A_1 - A_2) \cap \overline{A_2 - A_1}= \emptyset
\end{cases}
$
and since $X$ is hereditarily normal, there exists an open set 
$U_0 \subset X$ : $A_1 - A_2 \subset U_0 \subset \overline{U}_0 \subset X - (A_2 - A_1)$.  
Setting 
$
\begin{cases}
\ X_1 = \overline{U}_0 \cup (A_1 \cap A_2)\\
\ X_2 = (X - U_0) \cup (A_1 \cap A_2)
\end{cases}
\hspace{-.25cm}
,
$
the argument then proceeds as before.]\\

\begingroup%%------------------------------------>>
\fontsize{9pt}{11pt}\selectfont
Why work with classes of normal Hausdorff spaces?  
Answer: If the class $\sX$ contains a space that is not normal, then every nonempty Hausdorff space $Y \in \text{NES}(\sX)$ is necessarily a singleton.\\
\endgroup%%------------------------------------<<

\label{3.10}
\begingroup%%------------------------------------>>
\fontsize{9pt}{11pt}\selectfont
\textbf{\small FACT} \ 
Suppose that $Y$ is an AR (ANR).
\\
\indent\indent (1) Let $\sX$ be the class of perfectly normal paracompact Hausdorff spaces $-$then $Y \in \text{ES}(\sX)$ $(\text{NES}(\sX))$.
\\
\indent\indent (2) Let $\sX$ be the class of perfectly normal Hausdorff spaces $-$then $Y \in \text{ES}(\sX)$ $(\text{NES}(\sX))$ iff $Y$ is second countable.
\\ \indent
[For the necessity, remark that every collection of nonempty pairwise disjoint open subsets of $Y$ is countable.  
Reason: The construction on 
p. \pageref{6.44} ff. 
furnishes a perfectly normal Hausdorff space $X$ containing an uncountable closed discrete subspace $A$, the points of which cannot be separated by a collection of nonempty pairwise disjoint open subsets of $X$.]
\\
\indent\indent (3) Let $\sX$ be the class of paracompact Hausdorff spaces $-$then $Y \in \text{ES}(\sX)$ $(\text{NES}(\sX))$ iff $Y$ is completely metrizable.
\\ \indent
%%----------------------------------------------------------------------------------------------45
[To establish the necessity, assume, e.g., that $Y$ is an AR.  Let $X$ be the result of retopologizing $\beta Y$ by isolating the points of $\beta Y - Y$.  Every open covering of $X$ has a $\sigma$-discrete open refinement, hence $X$ is a paracompact Hausdorff space.  
Since $Y$ sits inside $X$ as a closed subspace, there is a retraction $r:X \ra Y$.
On the other hand, $Y$ is metrizable, thus is Moore, so $Y$ satisfies Arhangel'ski\u i's condition.  
Fix a sequence $\{\sV_n\}$ of collections of open subsets of $\beta Y$ such that each $\sV_n$ covers $Y$ and $\forall \ y \in Y$: $\ds\bigcap\limits_n \text{st}(y,\sV_n) \subset Y$.
Assign to a given $ V \in \sV_n$ the open subset $P_V \subset V$ determined by intersecting $V$ with the interior in $\beta Y$ of $r^{-1}(V \cap Y)$.  
Put $P_n = \ds\bigcup \{P_V: V \in \sV_n\}$ : $P_n \supset Y$ 
$\&$ $Y = \ds\bigcap\limits_n P_n$, therefore $Y$ is topologically complete or still, is completely metrizable.]
\\
\indent\indent (4) Let $\sX$ be the class of normal Hausdorff spaces $-$then $Y \in \text{ES}(\sX)$ $(\text{NES}(\sX))$ iff $Y$ is second countable and completely metrizable.\\
\endgroup%%------------------------------------<<

\begingroup%%------------------------------------>>
\fontsize{9pt}{11pt}\selectfont
\textbf{\small FACT} \ 
 Let $\sX$ be the class consisting of the Hausdorff spaces that can be realized as a closed subspace of a product of a compact Hausdorff space and a metrizable space (the elements of $\sX$ are precisely those paracompact Hausdorff spaces satisfying Arhangel'ski\u i's condition) $-$then every AR (ANR) is in $\text{ES}(\sX)$ $(\text{NES}(\sX))$.
\\ \indent
[Suppose that $X \in \sX$ is closed in $K \times Z$, where $K$ is compact Hausdorff and $Z$ is metrizable.  The projection $K \times Z \ra Z$ is closed and has compact fibers, thus the same is true of its restriction $p$ to $X$.  
Fix a closed subspace $A \subset X$.  Take an AR $Y$ of weight $\leq \kappa$ and let $f \in C(A,Y)$.  
Embed $Y$ in $\bS(\kappa)^\omega$ and apply Proposition 13 to produce a continuous extension 
$\phi:X \ra \bS(\kappa)^\omega$ of $f$.
Write $\Phi$ for the diagonal of $\phi$ and $p$ $-$then $\Phi(A)$ is closed in $\bS(\kappa)^\omega \times p(X)$.  
Therefore the restriction to $\Phi(A)$ of the projection 
$\psi:\bS(\kappa)^\omega \times p(X) \ra \bS(\kappa)^\omega$ has a continuous extension
$\Psi:\bS(\kappa)^\omega \times p(X) \ra Y$.  
Put $F = \Psi \circx \Phi$ : $F \in C(X,Y)$ $\&$ $\restr{F}{A} = f$.]\\
\endgroup%%------------------------------------<<

\begingroup%%------------------------------------>>
\fontsize{9pt}{11pt}\selectfont
Application: If $K$ is a compact Hausdorff space and if $Y$ is an ANR, then $C(K,Y)$ is an ANR (so for any CW complex $X$, $C(K,X)$ is a CW space).
\vspi
[Inspect the proof of Proposition 6, keeping in mind the preceding result.]\\
\endgroup%%------------------------------------<<

\begingroup%%------------------------------------>>
\fontsize{9pt}{11pt}\selectfont
Suppose that $G$ is a stratifiable topological group $-$then $X_G^\infty$ and $B_G^\infty$ are stratifiable 
(cf. p. \pageref{6.45}) and 
Cauty\footnote[2]{\textit{Arch. Math (Basel)} \textbf{28} (1977), 623-631.} 
has shown that if $G$ is also in NES(stratifiable), then the same holds true for $X_G^\infty$ and $B_G^\infty$.
\label{4.74}
\label{11.6}
Example: If $G$ is an ANR, then $X_G^\infty$ and $B_G^\infty$ are ANRs 
(cf. p.4-65 \pageref{4.46}).\\
\endgroup%%------------------------------------<<

\textbf{\small LEMMA} \ 
Let $Y$ be a topological space.  
Suppose that $Y$ admits a covering $\sV$ by pairwise disjoint open sets $V$, each of which is in NES(collectionwise normal) 
$-$then $Y$ is in NES(collectionwise normal).

%%----------------------------------------------------------------------------------------------46
[Let $X$ be collectionwise normal, $A \subset X$ closed, and let $f \in C(A,Y)$.  
Put $A_V = f^{-1}(V)$, $f_V = \restr{f}{A_V}$ $-$then there exists a neighborhood $O_V$ of $A_V$ in $X$ and an $F_V \in C(O_V,V)$: $\restr{F_V}{A_V} = f_V$.  
Since $\{A_V\}$ is a discrete collection of closed subsets of $X$, there exists a pairwise disjoint collection $\{U_V\}$ of open subsets of $X$ such that $\forall \ V$: $A_V \subset U_V$.  
Set $U = \bigcup\limits_V (O_V \hspace{0.05cm} \cap \hspace{0.05cm} U_V)$ 
and define 
$F:U \ra Y$ by $\restr{F}{O_V \hspace{0.05cm} \cap \hspace{0.05cm} U_V}$ $=$ 
$\restr{F_V}{O_V \hspace{0.05cm} \cap \hspace{0.05cm} U_V}$
 to get a continuous extension of $f$ to $U$.]\\[-.1cm]

Let $Y$ be a topological space.  
Suppose that $Y$ admits a numerable covering $\sV$ by open sets $V$, 
each of which is in NES(collectionwise normal) 
$-$then, from the proof of Proposition 12, it follows that $Y$ is in NES(collectionwise normal).\\[-.1cm]

\begingroup%%------------------------------------>>
\fontsize{9pt}{11pt}\selectfont
\textbf{\small FACT} \ 
Let $Y$ be a topological space.  
Suppose that $Y$ admits a covering $\sV$ by open sets $V$, each of which is in NES(paracompact) $-$then $Y$ is in NES(paracompact).\\[-.1cm]
\endgroup%%------------------------------------<<

\begingroup%%------------------------------------>>
\fontsize{9pt}{11pt}\selectfont
Application: Every topological manifold is in NES(paracompact).
\vspi
[Note: \  This applies in particular to the Pr\"ufer manifold, which is not metrizable and contains a closed submanifold that is not a neighborhood retract.]
\\[-.1cm]
\endgroup%%------------------------------------<<

Assume: $I\sX \subset \sX$.  
Let $Y \in \text{NES}(\sX)$ $-$then for every pair $(X,A)$, where $X \in \sX$ and $A \subset X$ is closed, $A$ has the HEP w.r.t. $Y$.  
Proof: $i_0 X \cup IA$, as a closed subspace of $IX$, has the EP w.r.t $Y$.\\[-.1cm]

\begingroup%%------------------------------------>>
\fontsize{9pt}{11pt}\selectfont
\label{20.6}
\textbf{\small EXAMPLE \  (\un{CW complexes})} \ 
If $X$ is stratifiable and $A \subset X$ is closed, then $A$ has the HEP w.r.t. any CW complex.\\[-.1cm]
\endgroup%%------------------------------------<<

\begin{proposition} \ %17
Assume: $I\sX \subset \sX$.  Let $Y \in \text{NES}(\sX)$ and suppose that $Y$ is homotopy equivalent to a 
$Z \in \text{ES}(\sX)$ $-$then $Y \in \text{ES}(\sX)$. 
\end{proposition}

[Choose continuous functions $\phi:Y \ra Z$, $\psi:Z \ra Y$ such that 
$\psi \circx \phi \simeq \text{id}_Y$, $\phi \circx \psi \simeq \text{id}_Z$.  
Take an $X$ in $\sX$ and let $A \subset X$ be closed.  
Given $f \in C(A,Y)$, $\exists \ F \in C(X,Z)$: $F \circx i = \phi \circx f$, where $i:A \ra X$ is the inclusion.  
But $A$ has the HEP w.r.t. $Y$ and $\psi \circx F \circx i \simeq f$, so $f$ admits a continuous extension to $X$.]\\[-.1cm]

\begingroup%%------------------------------------>>
\fontsize{9pt}{11pt}\selectfont
\label{20.8}
\textbf{\small FACT} \ 
Suppose that $X$ is an ANR.  Let $Y$ be a topological space such that every closed subset $A \subset X$ has the EP w.r.t $Y$.  Fix a weak homotopy equivalence $K \ra Y$, where $K$ is a CW complex $-$then every closed subset $A \subset X$ has the EP w.r.t $K$.
\vspi
[Owing to the CW-ANR theorem, the induced map $[X,K] \ra [X,Y]$ is bijective 
(cf. p. \pageref{4.47}).  
On the other hand, every closed subset $A \subset X$ has the HEP w.r.t. $K$ (metrizable $\implies$ stratifiable).]\\
\endgroup%%------------------------------------<<

%%%%%%%%%%%%%%%%%%%%%%%%%%%%%%%%%%%%%%
%%%%%%%%%%%%%%%%%%%%%%%%%%%%%%%%%%%%%%
%%%%%%%%%%%%%%%%%%%%%%%%%%%%%%%%%%%%%%

\begin{center}
$\S \ 6$
\\[0.5cm]
$\mathcal{REFERENCES}$\\
\end{center}

\[
\textbf{BOOKS}
\]

\begingroup
\fontsize{9pt}{11pt}\selectfont
\setlength\parindent{0 cm}

[1] \quad Al\`o, R. and Shapiro, H., \textit{Normal Topological Spaces}, Cambridge University Press (1974).
\\[-.2cm]

[2] \quad Borsuk, K., \textit{Theory of Retracts}, PWN (1967).
\\[-.2cm]

[3] \quad Hu, S.-T., \textit{Theory of Retracts}, Wayne State University Press (1965).
\\[-.2cm]
\endgroup

\[
\textbf{ARTICLES}
\]

\begingroup
\fontsize{9pt}{11pt}\selectfont
\setlength\parindent{0 cm}

[1] \quad Bogatyi, S. and Fedorchuk, V., Theory of Retracts and Infinite Dimensional Manifolds, 
\textit{J. Soviet Math.} 

\hspace{0.8cm}\textbf{44} (1989), 372-423.
\\[-.2cm]

[2] \quad Bryant, J., Ferry, S., Mio, W., and Weinberger, S., Topology of Homology Manifolds, 
\textit{Ann. of Math.} 

\hspace{0.8cm}\textbf{143} (1996), 435-467.
\\[-.2cm]

[3] \quad Chigogidze, A., Noncompact Absolute Extensors in Dimension $n$, $n$-Soft Mappings, and Their Applic-

\hspace{0.8cm}ations, 
\textit{Math. Izvestiya} \textbf{28} (1987), 151-174.
\\[-.2cm]

[4] \quad van Douwen, E., Lutzer, D., and Przymusi\'nski, T., Some Extensions of the Tietze-Urysohn Theorem, 
%\textit{Amer.}

\hspace{0.8cm}\textit{Amer. Math. Monthly} \textbf{84} (1977), 435-441.
\\[-.2cm]

[5] \quad Dranishnikov, A., Absolute Extensors in Dimension $n$ and Dimension Raising $n$-Soft Mappings, 
\textit{Rus-}

\hspace{0.8cm}\textit{sian Math. Surveys} \textbf{39} (1984), 63-111.
\\[-.2cm]

[6] \quad Dydak, J., Extension Theory: The Interface between Set-Theoretic and Algebraic Topology, 
\textit{Topology}

\hspace{0.95cm}\textit{Appl.} \textbf{74} (1996), 225-258.
\\[-.2cm]

[7] \quad Hanner, O., Retraction and Extension of Mappings of Metric and Nonmetric Spaces, 
\textit{Ark. Mat.} \textbf{2} 

\hspace{0.8cm}(1952), 315-360.
\\[-.2cm]

[8] \quad Gruenhage, G., Generalized Metric Spaces, In: \textit{Handbook of Set Theoretic Topology}, K. Kunen and J. 

\hspace{0.8cm}Vaughan (ed.), North Holland (1984), 423-501.
\\[-.2cm]

[9] \quad Hoshina, T., Extensions of Mappings, In: \textit{Topics in General Topology}, K. Morita and J. Nagata (ed.), 

\hspace{0.8cm}North Holland (1989), 41-80.
\\[-.2cm]

[10] \quad Milnor, J., On Spaces Having the Homotopy Type of a CW complex, 
\textit{Trans. Amer. Math. Soc.} \textbf{90} 

\hspace{0.95cm}(1959), 272-280.
\\[-.2cm]

[11] \quad Morita, K., Extensions of Mappings, In: op. cit. [9], 1-39.
\\[-.2cm]

[12] \quad Price, T., Towards Classifying All Manifolds, 
\textit{Math. Chronicle} \textbf{7} (1978), 1-47.
\\[-.2cm]

[13] \quad Repov\u s, D., Detection of Higher Dimensional Topological Manifolds among Topological Spaces, In:

\hspace{0.95cm}\textit{Seminari di Geometria}, M. Ferri (ed.), Bologna (1992), 113-143.
\\[-.2cm]

[14] \quad West, J., Open Problems in Infinite Dimensional Topology, In: \textit{Open Problems in Topology}, J. van 

\hspace{0.95cm}Mill and G. Reed (ed.), North Holland (1990), 523-597.
\\[-.2cm]

[15] \quad Zarichnyi, M., and Fedorchuk, V., Covariant Functors in Categories of Topological Spaces, 
\textit{J. Soviet}

\hspace{0.95cm}\textit{Math.} \textbf{53} (1991), 147-176.

\setlength\parindent{2em}

\endgroup

\chapter{
$\boldsymbol{\S}$\textbf{7}.\quadx  $\mathcal{C}-$THEORY}
\setlength\parindent{2em}
\setcounter{proposition}{0}

%%----------------------------------------------------------------------------------------------01
$\text{ }$\\[-1.25cm]

A classical technique in algebraic topology is to work modulo a Serre class of abelian groups.  
I shall review these matters here, supplying proofs of the less familiar facts.

Let $\sC \subset \Ob\bAB$ be a nonempty class of abelian groups $-$then $\sC$ is said to be a 
\un{Serre class}
\index{Serre class} 
provided that for any short exact sequence  
$0 \ra$ 
$G^\prime \ra$ 
$G \ra$ 
$G\pp \ra 0$ in \bAB, $G \in \sC$ iff 
$
\begin{cases}
\ G^\prime \\
\ G\pp
\end{cases}
\hspace{-.25cm}
\in \sC
$
or, equivalently, for any exact sequence 
$G^\prime \ra$ 
$G \ra$ 
$G\pp$ in \bAB
$
\begin{cases}
\ G^\prime \\
\ G\pp
\end{cases}
\hspace{-.25cm}
\in \sC
$
$\implies$ $G \in \sC$.

[Note: \ To show that a nonempty class $\sC \subset \Ob\bAB$ is a Serre class, it is usually simplest to check that 
$\sC$ is closed under subgroups, homomorphic images, and extensions.]

Example: For any Serre class $\sC$, the subclass $\sC_\tor$ of torsion groups in $\sC$ is a Serre class.

[Note: \ A Serre class $\sC$ is said to be 
\un{torsion}
\index{torsion (Serre class)} 
if $\sC = \sC_\tor$.]\\

\label{7.1}
\begingroup%%----------------------------------->>
\fontsize{9pt}{11pt}\selectfont
\index{p-Primary Abelian Groups}
\textbf{\small EXAMPLE \ (\un{$p$-Primary Abelian Groups})} \  
An abelian $p$-group \mG is said to be 
\un{$p$-primary}.
\index{p-primary ($p$-primary abelian group)}  
The 
\un{rank}
\index{rank(abelian group)} 
$r(G)$ of a $p$-primary \mG is the cardinality of a maximal independent system in 
\mG.  If $G[p] = \{g:pg = 0\}$, then $G[p]$ is a vector space over $\F_p$ and $\dim G[p] = r(G)$.  
The 
\un{final rank}
\index{final rank($p$-primary abelian group)} 
$r_f(G)$ of a $p$-primary \mG is the infimum of the $r(p^n G)$ $(n \in \N)$.  
Every $p$-primary \mG can be written as $G = G^\prime \oplus G\pp$, where 
$G^\prime$ is bounded and $r(G\pp) = r_f(G\pp)$ 
(Fuchs\footnote[2]{\textit{Infinite Abelian Groups}, vol. I, Academic Press (1970), 152.}).  
Fix now a symbol $\infty$, considered to be larger than all the cardinals.  
Given a Serre class $\sC$ of $p$-primary abelian groups, let 
$\Phi(\sC)$ be the smallest cardinal number $> r(G)$ $\forall \ G \in \sC$ if such a number exists, otherwise put 
$\Phi(\sC) = \infty$, and let 
$\Psi(\sC)$ be the smallest cardinal number $> r_f(G)$ $\forall \ G \in \sC$ if such a number exists, otherwise put 
$\Psi(\sC) = \infty$.  
Obviously, $\Phi(\sC) \geq \Psi(\sC)$, 
$
\begin{cases}
\ \Phi(\sC) = 1 \text{\quad or \quad} \Phi(\sC) \geq \omega\\
\ \Psi(\sC)= 1 \text{\quad or \quad} \Psi(\sC) \geq \omega
\end{cases}
\hspace{-.25cm}
. \ 
$
And: $\sC$ is precisely the class of $p$-primary \mG for which $r(G) < \Phi(\sC)$ $\&$ $r_f(G) < \Psi(\sC)$.  
On the other hand, suppose that 
$
\begin{cases}
\ \Phi\\
\ \Psi
\end{cases}
$
are cardinal numbers or $\infty$ with $\Phi \geq \Psi$, 
$
\begin{cases}
\ \Phi = 1 \text{\quad or \quad} \Phi \geq \omega\\
\ \Psi = 1 \text{\quad or \quad} \Psi \geq \omega
\end{cases}
\hspace{-.25cm}
. \ 
$
Let $\sC$ be the class of $p$-primary \mG for which $r(G) < \Phi$ $\&$ $r_f(G) < \Psi$ $-$then $\sC$ is a Serre class such that 
$\Phi(\sC) = \Phi$ $\&$ $\Psi(\sC) = \Psi$.  
Thus the conclusion is that there is a one-to-one correspondence between the conglomerate of Serre classes of $p$-primary abelian groups and the conglomerate of ordered pairs $(\Phi,\Psi)$, where 
$
\begin{cases}
\ \Phi\\
\ \Psi
\end{cases}
$
are cardinal numbers or $\infty$: $\Phi \geq \Psi$, 
$
\begin{cases}
\ \Phi = 1 \text{\quad or \quad} \Phi \geq \omega\\
\ \Psi = 1 \text{\quad or \quad} \Psi \geq \omega
\end{cases}
$
\hspace{-.25cm}
.
\\ \indent
[Note: \ If $\sC$ is a Serre class and if $\sC(p)$ is the subclass of $\sC$ consisting of the $p$-primary $G$ in $\sC$, then 
$\sC(p)$ is a Serre class.]\\
\endgroup %%------------------------------------<<

%%----------------------------------------------------------------------------------------------02
Notation: Given a Serre class $\sC$, $\tf(\sC)$ is the subclass of $\sC$ made up of the torsion free groups in $\sC$.\\

\begin{proposition} \ %01
Let $\sC$ be a Serre class.  Assume: $\tf(\sC)$ contains a group of infinite rank $-$then either $\tf(\sC)$ is the class of all torsion free abelian groups or $\tf(\sC)$ is the class of all torsion free abelian groups of cardinality $< \kappa$, where 
$\kappa > \omega$.
\end{proposition}

[Any torsion free abelian group \mG of infinite rank contains a free abelian group of rank $= \#(G)$.]\\

\begingroup%%----------------------------------->>
\fontsize{9pt}{11pt}\selectfont
\textbf{\small EXAMPLE} \quad 
Fix a cardinal number $\kappa > \omega$.  Let $\sT_\kappa$ be the class of torsion abelian groups of cardinality $< \kappa$; let $\sF_\kappa$ be the class of torsion free abelian groups of cardinality $< \kappa$.  Take any Serre class $\sT$ of torsion abelian groups: $\sT \supset \sT_\kappa$ $-$then the class $\sC$ consisting of all abelian groups \mG which are extensions of a group in $\sT$ by a group in $\sF_\kappa$ is a Serre class such that $\sC_\tor = \sT$ and $\tf(\sC) = \sF_\kappa$.\\
\endgroup %%------------------------------------<<

A 
\un{characteristic}
\index{characteristic} 
is a sequence $\chi = \{\chi_p:p \in \bPi\}$, where each $\chi_p$ is a nonnegative integer or $\infty$.  Given characteristics
$
\begin{cases}
\ \chi^\prime\\
\ \chi\pp
\end{cases}
\hspace{-.25cm}
$
, write 
$\chi^\prime \sim \chi\pp$ iff $\#\{p:\chi_p^\prime \neq \chi_p\pp\} < \omega$ and 
$\chi_p^\prime = \infty$ $\Longleftrightarrow$ $\chi_p\pp = \infty$ $-$then $\sim$ is an equivalence relation on the set of characteristics, an equivalence class \bt being called a 
\un{type}.
\index{type (equivalence class on characteristics)}  
The sum $\bt^\prime + \bt\pp$ of types 
$
\begin{cases}
\ \bt^\prime\\
\ \bt\pp
\end{cases}
$
is the type containing the characteristic $\{\chi_p^\prime + \chi_p\pp: p \in \bPi\}$ and $\bt^\prime \leq \bt\pp$ provided that 
$\chi_p^\prime \leq \chi_p\pp$ for almost all $p$, 
$\bt\pp - \bt^\prime$ being the largest type \bt such that 
$\bt + \bt^\prime \leq \bt\pp$.

\indent\indent (Rational Groups) \ 
A nonzero abelian group \mG is said to be 
\un{rational}
\index{rational (abelian group)} 
if it is isomorphic to a subgroup of $\Q$ or still, is torion free of rank 1.  
Such groups can be classified.  For assume that \mG is rational, say $G \subset \Q$.  Take 
$g \in G$: $g \neq 0$.  
Given $p \in \bPi$, consider the set $S_p(g)$ of nonnegative integers $n$ such that the equation 
$p^n x = g$ has a solution in \mG.  
Put 
$\chi_p(g) = \sup S_p(g)$, the 
\un{$p$-height}
\index{p-height} of $g$ $-$then 
$\chi(g) = \{\chi_p(g): p \in \bPi\}$ is a characteristic.  
Moreover, distinct nonzero elements of \mG determine equivalent characteristics.  
Definition: The 
\un{type}
\index{type (of a rational group)} 
$\bt(G)$ of \mG is the type of the characteristic of any nonzero element of \mG.  
Every type \bt can be realized by a rational group, i.e., $\bt = \bt(G)$ $(\exists \ G)$ and rational 
$
\begin{cases}
\ G^\prime\\
\ G\pp
\end{cases}
$
are isomorphic iff $\bt(G^\prime) = \bt(G\pp)$ (in general, $G^\prime$ is isomorphic to a subgroup of $G\pp$ iff 
$\bt(G^\prime) \leq \bt(G\pp)$).

Example: Suppose that $\Z \subset G \subset \Q$ $-$then 
$G/\Z \approx \bigoplus\limits_p \Z/p^{\chi_p}\Z$, $\{\chi_p:p \in \Pi\}$ the characteristic of 1, and $\Hom(G,G)$ is isomorphic to the subring of $\Q$ generated by 1 and the $p^{-1}$: $pG = G$.\\

\begingroup%%----------------------------------->>
\fontsize{9pt}{11pt}\selectfont
\textbf{\small FACT} \quad 
If \mG and \mK are rational, then $G \otimes K$ is rational and $\bt(G \otimes K) = \bt(G) + \bt(K)$.\\
\endgroup %%------------------------------------<<

%%----------------------------------------------------------------------------------------------03
\begingroup%%----------------------------------->>
\fontsize{9pt}{11pt}\selectfont
\textbf{\small FACT} \quad 
If \mG and \mK are rational, then $\Hom(G,K) = 0$ if $\bt(G) \not\leq \bt(K)$, but is rational if $\bt(G) \leq \bt(K)$ with 
$\bt(\Hom(G,K)) = \bt(K) - \bt(G)$.\\
\endgroup %%------------------------------------<<

Notation: \bT is a nonempty set of types such that 
(i) $\bt_0 \in \bT$ $\&$ $\bt \leq \bt_0$ $\implies$ $\bt \in \bT$ and 
(ii) $\bt^\prime, \bt\pp \ \in \bT$ $\implies$ $\bt^\prime + \bt\pp \in \bT$, $\bT(\bAB)$ 
being the class of abelian groups \mG which admit a monomorphism 
$G \ra \bigoplus\limits_1^n G_i$, where the $G_i$ are rational ($n$ depending on \mG) and the 
$\bt(G_i) \in \bT$.\\

\begingroup%%----------------------------------->>
\fontsize{9pt}{11pt}\selectfont
\textbf{\small FACT} \quad 
A torsion free abelian group \mG of finite rank is in $\bT(\bAB)$ iff for each nonzero homomorphism 
$\phi:G \ra \Q$, $\bt(\phi(G)) \in \bT$.\\
\endgroup %%------------------------------------<<

\begin{proposition} \ %02
Let $\sC$ be a Serre class.  Assume: $\tf(\sC)$ contains only groups of finite rank and at least one group of positive rank $-$then 
$\tf(\sC) = \bT(\bAB)$ for some \bT.
\end{proposition}

[Let \bT be the set of types \bt such that a rational group of type \bt is in $\tf(\sC)$.  If 
$G_1, \ldots, G_n$ are rational and if 
$\bt(G_1), \ldots, \bt(G_n)$ belong to \bT, then 
$\bigoplus\limits_1^n G_i \in \tf(\sC)$ and every subgroup of 
$\bigoplus\limits_1^n G_i$ is in $\tf(\sC)$.  
On the other hand, for any $G \neq 0$ in $\tf(\sC)$, there are rational
$G_1, \ldots, G_n$ and a monomorphism $G \ra \bigoplus\limits_1^n G_i$.  
Upon restricting to homomorphic images, one can arrange that the $G_i \in \tf(\sC)$, so the $\bt(G_i) \in \bT$.  Since $\sC$ is closed under subgroups, \bT satisfies condition (i) above.  
As for condition (ii), let 
$
\begin{cases}
\ \bt^\prime\\
\ \bt\pp
\end{cases}
\hspace{-.25cm}
\in \bT.
$
\  Choose 
$
\begin{cases}
\ G^\prime\\
\ G\pp
\end{cases}
\hspace{-.25cm}
:  
$
$\Z \subset$
$
\begin{cases}
\ G^\prime\\
\ G\pp
\end{cases}
\hspace{-.25cm}
$
$\subset \Q$ $\&$
$
\begin{cases}
\ \bt^\prime = \bt(G^\prime)\\
\ \bt\pp = \bt(G\pp)
\end{cases}
$
is represented by the characteristic
$
\begin{cases}
\ \chi^\prime\\
\ \chi\pp
\end{cases}
$
corresponding to 1.  Suppose first that $\forall \ p$, 
$
\begin{cases}
\ \chi_p^\prime\\
\ \chi_p\pp
\end{cases}
$
is finite.  Let $\Z \subset G \subset \Q$: $\chi(1) = \chi^\prime + \chi\pp$.  
Fix an isomorphim 
$\phi:G^\prime/\Z \ra G/G\pp$ and let \mK be the subgroup of $G^\prime \oplus G$ composed of the 
$(g^\prime,g):\phi(g^\prime + \Z) = g + G\pp$ $-$then there is a short exact sequence 
$0 \ra$ 
$G\pp \ra$ 
$K \ra$ 
$G^\prime \ra 0$, hence $K \in \sC$.  
But there is also an epimorphism $K \ra G$, thus $G \in \sC$ and $\bt^\prime + \bt\pp \in \bT$.  
Passing to the general case, write
$
\begin{cases}
\ \chi^\prime = \chi_f^\prime + \chi_{0,\infty}^\prime\\
\ \chi\pp = \chi_f\pp + \chi_{0,\infty}\pp
\end{cases}
\hspace{-.25cm}
,
$
where 
$
\begin{cases}
\ \chi_f^\prime \\
\ \chi_f\pp 
\end{cases}
$
take finite values and 
$
\begin{cases}
\ \chi_{0,\infty}^\prime\\
\ \chi_{0,\infty}\pp
\end{cases}
$
have values 0 or $\infty$.  
Let $\Z \subset G_f \subset \Q$ : $\chi_f(1) = \chi_f^\prime + \chi_f\pp$; let 
$\Z \subset$
$
\begin{cases}
\ G_{0,\infty}^\prime\\
\ G_{0,\infty}\pp
\end{cases}
\hspace{-.25cm}
$
$\subset \Q:$ 
$
\begin{cases}
\ \chi_{0,\infty}^\prime (1) = \chi_{0,\infty}^\prime\\
\ \chi_{0,\infty}\pp (1) = \chi_{0,\infty}\pp
\end{cases}
\hspace{-.25cm}
. \ 
$
From the foregoing, $G_f \in \sC$; in addition, 
$
\begin{cases}
\ G_{0,\infty}^\prime\\
\ G_{0,\infty}\pp
\end{cases}
$
is isomorphic to a subgroup of 
$
\begin{cases}
\ G^\prime\\
\ G\pp
\end{cases}
\hspace{-.25cm}
\in \sC.
$
Therefore
$G_f \oplus G_{0,\infty}^\prime \oplus G_{0,\infty}^{\prime\prime} \in \sC$ and 
$G_f + G_{0,\infty}^\prime + G_{0,\infty}^{\prime\prime} \subset \Q$ 
has type $\bt^\prime + \bt\pp$.]\\

\begingroup%%----------------------------------->>
\fontsize{9pt}{11pt}\selectfont
\textbf{\small EXAMPLE} \quad 
Given \bT, let $\sT$ be a Serre class of torsion abelian groups with the property that the
%%----------------------------------------------------------------------------------------------04
type determined by a characteristic $\chi$ belongs to \bT iff 
$\ds\bigoplus\limits_p \Z/p^{\chi_P} \Z \in \sT$ $-$then the class $\sC$ consisting of all abelian groups \mG 
which are extensions of a group in $\sT$ by a group in $\bT(\bAB)$ is a Serre class such that 
$\sC_\tor = \sT$ and $\tf(\sC) = \bT(\bAB)$.\\
\endgroup %%------------------------------------<<

Every torsion abelian group \mG contains a 
\un{basic subgroup}
\index{basic subgroup (of a torsion abelian group)} 
\mB, i.e., \mB is a direct sum of cyclic groups, \mB is pure in \mG, and $G/B$ is divisible.  If
$
\begin{cases}
\ G^\prime\\
\ G\pp
\end{cases}
$
are torsion and if 
$
\begin{cases}
\ B^\prime \subset G^\prime\\
\ B\pp \subset G\pp
\end{cases}
$
are basic, then 
$G^\prime \otimes G\pp \approx$ $B^\prime \otimes B\pp$.  
Corollary: The tensor product of two torsion abelian groups is a direct sum of cyclic groups.\\

\textbf{\small LEMMA} \quad 
Let $0 \ra G^\prime \ra G \ra G\pp \ra 0$ be a short exact sequence of abelian groups.  Suppose that the image of $G^\prime$ in \mG is pure $-$then for any \mK, the sequence 
$0 \ra$ 
$G^\prime \otimes K \ra$
$G \otimes K \ra$
$G\pp \otimes K \ra 0$ is exact and the image of $G^\prime \otimes K$ in $G \otimes K$ is pure.

[Note: \ Under the same assumptions, the sequence 
$0 \ra$ 
$\Tor(G^\prime,K) \ra$ 
$\Tor(G,K) \ra$ 
$\Tor(G\pp,K) \ra 0$ 
is exact and the image of 
$\Tor(G^\prime,K)$ in $\Tor(G,K)$ 
is pure.]\\

A Serre class $\sC$ is said to be a 
\un{ring}
\index{ring (Serre class)} 
if $G, K \in \sC$ $\implies$ $G \otimes K \in \sC$, 
$\Tor(G,K) \in \sC$.

[Note: \ $\sC$ is a ring provided that $\forall \ G \in \sC$: $G \otimes G \in \sC$, $\Tor(G,G) \in \sC$.  This is because 
$G, K \in \sC$ $\implies$ 
$G \otimes K \subset (G \oplus K) \otimes (G \oplus K)$,
$\Tor(G,K) \subset \Tor(G \oplus K,G \oplus K)$.]\\

\begingroup%%----------------------------------->>
\fontsize{9pt}{11pt}\selectfont
\textbf{\small EXAMPLE} \quad 
Let $\sC$ be a ring.  Fix a group \mG $-$then $G/[G,G] \in \sC$ iff $\forall \ i$, $\Gamma^i(G)/\Gamma^{i+1}(G) \in \sC$.  
\vspi
[The iterated commutator map 
$\otimes^{i+1} (G/[G,G]) \ra \Gamma^i(G)/\Gamma^{i+1}(G)$ is surjective.]\\
\endgroup %%------------------------------------<<

\label{7.4}
\begingroup%%----------------------------------->>
\fontsize{9pt}{11pt}\selectfont
\textbf{\small EXAMPLE} \quad 
Let $\sC$ be a ring.  Fix a group \mG such that $\forall \ n > 0$, $H_n(G) \in \sC$.  Let $M \in \sC$ be a nilpotent 
$G$-module $-$then $\forall \ n \geq 0$, $H_n(G;M) \in \sC$.
\vspi
[Since the $(I[G])^i \cdot M/(I[G])^{i+1} \cdot M \in \sC$, it suffices to look at the case when the action of \mG on \mM is trivial.]\\
\endgroup %%------------------------------------<<

\begingroup%%----------------------------------->>
\fontsize{9pt}{11pt}\selectfont
\textbf{\small FACT} \quad 
Let $\sC$ be a Serre class.  Suppose that $G \in \sC$ $-$then for any finitely generated \mK, $G \otimes K$ and $\Tor(G,K)$ belong to $\sC$.\\
\endgroup %%------------------------------------<<

\begin{proposition} \ %03
Let $\sC$ be a Serre class $-$then $\sC$ is a ring iff $\sC_\tor$ is a ring.
\end{proposition}

[Setting aside the trivial case when $\sC$ is the class of all abelian groups, let us assume that $\sC_\tor \neq \sC$ is a ring.  Fix $G \in \sC - \sC_\tor$: 
$\Tor(G,G) \approx$ 
$\Tor(G_\tor,G_\tor) \in \sC_\tor$, $G_\tor$ the torsion subgroup of \mG.  To deal with $G \otimes G$, put
$\tf(G) = G/G_\tor$ and consider the exact sequences
\[
\begin{cases}
\ 0 \ra G_\tor \otimes G \ra G \otimes G \ra \tf(G) \otimes G \ra 0\\
\ 0 \ra G_\tor \otimes G_\tor \ra G_\tor \otimes G \ra G_\tor \otimes \tf(G) \ra 0\\
\ 0 \ra \tf(G) \otimes G_\tor \ra \tf(G) \otimes G \ra \tf(G) \otimes \tf(G) \ra 0
\end{cases}
.
\]
%%----------------------------------------------------------------------------------------------05
Because $G_\tor \otimes G_\tor \in \sC_\tor$, it will be enough to prove that 
$G_\tor \otimes \tf(G)$ and
$\tf(G) \otimes \tf(G)$ are in $\sC$.

\indent\indent (I) \ 
Suppose that $\tf(G)$ contains a group of infinite rank.  
Choose $\kappa > \omega$ as in Proposition 1 (so $\sC$ contains all abelian groups of cardinality $< \kappa$): $\#(\tf(G)) < \kappa$ $\implies$ $\#(\tf(G) \otimes \tf(G)) < \kappa$ 
$\implies$ 
$\tf(G) \otimes \tf(G) \in \sC$.  
There is a free group \mF in $\sC$ and an epimorpism $F \ra \tf(G) \ra 0$, where $\rank F < \kappa$.  Let \mB be a basic subgroup of $G_\tor$ and form the exact sequence 
$0 \ra$ 
$B \otimes F \ra$
$G_\tor \otimes F \ra$
$G_\tor/B \otimes F \ra 0$.  Using the fact that \mB is a direct sum of cyclic groups, 
$B \otimes F \approx$ 
$B \otimes B_\kappa$: $\#(B_\kappa) < \kappa$ $\implies$ $B \otimes F \in \sC$.  Analogously, by an application of the structure theorem for divisible abelian groups, $G_\tor/B \otimes F \in \sC$.  
Conclusion: $G_\tor \otimes F \in \sC$ $\implies$ $G_\tor \otimes \tf(G) \in \sC$.

\indent\indent (II) \ 
Suppose that $\tf(\sC) = \bT(\bAB)$ (cf. Proposition 2).  Let \mF be the free abelian group generated by a maximal independent system in $\tf(G)$ $-$then there is an exact sequence 
$0 \ra$ 
$F \ra$ 
$\tf(G) \ra$ 
$\tf(G)/F \ra 0$ and $\tf(G)/F \in \sC_\tor$.  
Tensor this sequence with $G_\tor$ to get another exact sequence 
$F \otimes G_\tor \ra$ 
$\tf(G) \otimes G_\tor \ra$ 
$\tf(G)/F \otimes G_\tor$.  Of course, $\tf(G)/F \otimes G_\tor \in \sC_\tor$; moreover 
$F \otimes G_\tor \in \sC$ which implies that 
$\tf(G) \otimes G_\tor$ itself is in $\sC$.  Finally, the sequence 
$0 \ra$ 
$F \otimes \tf(G) \ra$ 
$\tf(G) \otimes \tf(G) \ra$ 
$\tf(G)/F \otimes \tf(G) \ra 0$ is exact.  Obviously, $F \otimes \tf(G) \in \sC$ and, repeating the preceding argument, 
$\tf(G)/F \otimes \tf(G) \in \sC$, hence $\tf(G) \otimes \tf(G) \in \sC$.]\\

In what follows, $\alpha$ and $\gamma$ are functions having cardinal numbers as values, the domain of $\alpha$ being 
$\bPi \times \N$ and the domain of $\gamma$ being $\bPi$.

Examples: \ 
(1) Let \mG be a torsion abelian group.  Assume: \mG is a direct sum of cyclic groups $-$then 
$G \approx$ 
$\bigoplus\limits_p \bigoplus\limits_n \alpha(p,n) \cdot (\Z/p^n \Z)$;
(2) Let \mG be a torsion abelian group.  Assume: \mG is divisible $-$then 
$G \approx$ 
$\bigoplus\limits_p  \gamma(p)\cdot (\Z/p^\infty \Z)$;
(3) Let \mG be a torsion abelian group.  Assume: \mG is $p$-primary and satisfies the descending chain condition on subgroups $-$then 
$G \approx$ 
$\bigoplus\limits_n \alpha(p,n) \cdot (\Z/p^n \Z) \oplus \gamma(p) \cdot (\Z/p^\infty \Z)$, where 
$\sum\limits_n \alpha(p,n) < \omega$ and $\gamma(p)$ is finite.

[Note: \ For use below, recall that $\Z/p^\infty \Z$ is a homomorphic image of 
$\bigoplus\limits_n \Z/p^n \Z$ 
(in fact, every countable $p$-primary \mG is a homomorphic image  of 
$\bigoplus\limits_n \Z/p^n \Z$).]

Notation: Given a torsion Serre class $\sC$, 
$\balpha(\sC) = \{\alpha: \bigoplus\limits_p \bigoplus\limits_n \alpha(p,n) \cdot (\Z/p^n \Z) \in \sC\}$ and 
$\bgamma(\sC) = \{\gamma: \bigoplus\limits_p \gamma(p) \cdot (\Z/p^\infty \Z) \in \sC\}$.

Observations: 
(i) $\gamma_0 \in \bgamma(\sC)$ $\&$ $\gamma \leq \gamma_0$ $\implies$ $\gamma \in \bgamma(\sC)$ and 
(ii) $\gamma^\prime, \gamma\pp \in \bgamma(\sC)$ $\implies$  $\gamma^\prime + \gamma\pp$ $\in \bgamma(\sC)$.

Suppose that $\sC$ is a torsion Serre class.  Let $G \in \sC$ $-$then 
$G \approx$ 
$\bigoplus\limits_p G(p)$, $G(p)$ the $p$-primary component of \mG.  Denote by $\sC_0$ the subclass of $\sC$ comprised of those \mG such that 
%%----------------------------------------------------------------------------------------------06
each $G(p)$ is bounded, so $\forall \ p$, $\exists$ $M(p)$: $p^{M(p)} G(p) = 0$, and put $\balpha_0(\sC) = \balpha(\sC_0)$ 
(meaningful, $\sC_0$ being Serre).\\

\index{Cardinal Lemma}
\textbf{\small CARDINAL LEMMA} \quad 
Let $\sC$ be a torsion Serre class $-$then $\forall \ \alpha \in \balpha(\sC)$, $\exists$ $\alpha_0 \in \balpha_0(\sC)$ $\&$ 
$\gamma \in \bgamma(\sC)$ such that 
$\alpha(p,n) \leq \alpha_0(p,n) + \gamma(p)$, where $\gamma(p) \geq \omega$ or $\gamma(p) = 0$.

[Set 
$\sigma(p,n) = \ds\sum\limits_{m=n}^\infty \alpha(p,m)$ and choose $M(p)$ such that 
$\sigma(p,n) =$ $\sigma(p,n+1) = \cdots$ $(n \geq M(p))$.  
Define $\alpha_0$ by 
$
\alpha_0(p,n) = 
\begin{cases}
\ \alpha(p,n) \hspace{0.5cm} (n < M(p))\\
\ 0 \hspace{1.35cm} \  (n \geq M(p))
\end{cases}
$
and define $\gamma$ by 
$\gamma(p) = \sigma(p,M(p))$: $\alpha \leq \alpha_0 + \gamma$ and 
$\alpha_0 \in \balpha_0(\sC)$, thus the issue is whether $\gamma \in \bgamma(\sC)$.  To see this, it need only be shown that
$\forall  \ p$, $\gamma(p)\cdot (\Z/p^\infty\Z)$ is a homomorphic image of 
$\bigoplus\limits_n \alpha(p,n)\cdot (\Z/p^n\Z)$.  
Case 1: $\gamma(p) = \omega$.  Here $\#\{n:\alpha(p,n) \neq 0\} = \omega$ and there are epimorphisms 
$\bigoplus\limits_n \alpha(p,n)\cdot (\Z/p^n\Z) \ra$
$\bigoplus\limits_n (\Z/p^n\Z) \ra$ 
$\gamma(p) \cdot (\Z/p^\infty\Z)$.  
Case 2: $\gamma(p) > \omega$.  
Put $N_\infty = \{n:\alpha(p,n) > \omega\}$: $\#(N_\infty) = \omega$ and there are epimorphisms 
$\bigoplus\limits_n \alpha(p,n)\cdot (\Z/p^n\Z) \ra$
$\bigoplus\limits_{n \in N_\infty} n\alpha(p,n)\cdot (\Z/p^n\Z) \ra$
$\bigoplus\limits_{n \in N_\infty} \alpha(p,n)\cdot (\Z/p\Z \oplus \cdots \oplus \Z/p^n\Z) \ra$
$\gamma(p) \cdot (\bigoplus\limits_n \Z/p^n\Z) \ra$ 
$\gamma(p) \cdot (\Z/p^\infty\Z)$.]\\

Given a torsion Serre class $\sC$, let $\sC^*$ be the subclass of those \mG such that each $G(p)$ satisfies the descending chain condition on subgroups.  
Note that $\sC^*$ is Serre.\\

\begin{proposition} \ %04
Let $\sC$ be a torsion Serre class $-$then $\sC$ is a ring iff $\sC^*$ is a ring.
\end{proposition}

[Straightforward computations establish the necessity.  As for the sufficiency, fix $G \in \sC$ and let \mB be a basic subgroup of \mG.  Applying the cardinal lemma, one finds that $B \otimes B \in \sC$.  But 
$G \otimes G \approx B \otimes B$, thus $G \otimes G \in \sC$.  
The verification that $\Tor(G,G)  \in \sC$ hinges on a preliminary remark.

Claim: Suppose that $\sC^*$ is a ring $-$then $\forall \ \gamma \in \bgamma(\sC)$, $\gamma^2 \in \bgamma(\sC)$.

[Write 
$\gamma = \gamma^\prime + \gamma\pp$, where $\forall \ p$, $\gamma^\prime(p)$ is finite and 
$\gamma\pp(p) \geq \omega$ or $\gamma\pp(p) = 0$, so 
$\gamma^2 = (\gamma^\prime)^2 + \gamma\pp$.  Since $\sC^*$ is a ring, 
$(\gamma^\prime)^2 \in \bgamma(\sC)$, hence $\gamma^2 \in \bgamma(\sC)$.]

Consider the exact sequences
\[
\begin{cases}
\ 0 \ra \Tor(B,G) \ra \Tor(G,G)\ra \Tor(G/B,G) \ra 0\\
\ 0 \ra \Tor(B,B) \ra \Tor(G,B)\ra \Tor(G/B,B) \ra 0\\
\ 0 \ra \Tor(B,G/B) \ra \Tor(G,G/B) \ra \Tor(G/B,G/B) \ra 0\\
\end{cases}
.
\]
Owing to the claim, $\Tor(G/B,G/B) \in \sC$.  Proof: 
$G/B \approx \bigoplus\limits_p \gamma(p) \cdot (\Z/p^\infty\Z)$ $\implies$ 
$\Tor(G/B,G/B) \approx$ $\bigoplus\limits_p \gamma^2(p) \cdot (\Z/p^\infty\Z)$.  
In addition, 
$\Tor(B,B) \approx B \otimes B$ $\in \sC$.  Therefore everything comes down to showing that 
$\Tor(B,G/B) \in \sC$ or still, that 
$\bigoplus\limits_p \gamma(p) \cdot B(p) \in \sC$.  Using 
%%----------------------------------------------------------------------------------------------07
the cardinal lemma, represent \mB by $B_0 \oplus B_\infty$ with 
$B_0(p) = \bigoplus\limits_n \alpha_0(p,n) \cdot (\Z/p^n\Z)$ and 
$B_\infty(p) = \bigoplus\limits_n \alpha_\infty(p,n) \cdot (\Z/p^n\Z)$, subject to 
$(\alpha_0)$ $\forall \ p$, $\exists$ $M(p)$: $n \geq M(p)$ $\implies$ $\alpha_0(p,n) = 0$ and 
$(\alpha_\infty)$ $\exists$ $\gamma_\infty \in \bgamma(\sC)$: $\forall \ p$, $\forall \ n$, 
$\alpha_\infty(p,n) \leq \gamma_\infty(p)$, where $\gamma_\infty(p) \geq \omega$ or $\gamma_\infty(p) = 0$.  
From the definitions, 
$\bigoplus\limits_p \gamma(p) \cdot B_0(p) \approx$ 
$B_0 \otimes \bigl(\bigoplus\limits_p \gamma(p) \cdot (\Z/p^{M(p)}\Z)\bigr) \in \sC$.  
Turning to $B_\infty$, for each $p$, there is a monomorphism 
$\gamma(p) \cdot B_\infty(p) \ra$ $(\gamma(p) + \gamma_\infty(p)) \cdot (\Z/p^\infty\Z)$.  
Because $\gamma + \gamma_\infty \in \bgamma(\sC)$, it follows that 
$\bigoplus\limits_p \gamma(p) \cdot B_\infty(p) \in \sC$.]\\

Application: Let $\sC$ be a Serre class.  
Assume $\tf(\sC)$ contains a free group of infinite rank $-$then $\sC$ is a ring.\\

\begingroup%%----------------------------------->>
\fontsize{9pt}{11pt}\selectfont
\textbf{\small EXAMPLE} \quad 
Not every Serre class is a ring.  For instance, let $\sC$ be the class of all torsion abelian groups \mG such that $\forall \ p$, 
$G(p)$ is finite, so 
$G(p) = \ds\bigoplus\limits_p \alpha(p,n) \cdot (\Z/p^n/\Z)$, where 
$r(G(p)) = \ds\sum\limits_n \alpha(p,n) < \omega$ 
(cf. p. \pageref{7.1}).  
Enumerate $\bPi$: $p_1 < p_2 < \cdots$ $-$then the subclass of $\sC$ consisting of those \mG for which the sequence 
$\{r(G(p_k))/k\}$ is bounded is a Serre class but it is not a ring (consider 
$G = \ds\bigoplus\limits_k k \cdot (\Z/p_k\Z)$).
\\ \indent
[Note: \ $\sC$ is a Serre class and it is a ring.]\\
\endgroup %%------------------------------------<<

A Serre class $\sC$  is said to be 
\un{acyclic}
\index{acyclic (Serre class)} 
if $\forall \ G \in \sC$, $H_n(G) \in \sC$ $(n > 0)$.\\

\begingroup%%----------------------------------->>
\fontsize{9pt}{11pt}\selectfont
\textbf{\small FACT} \quad 
Let $\sC$ be a Serre class.  
Suppose that $G \in \sC$ is finitely generated $-$then $H_n(G) \in \sC$ $(n > 0)$.\\
\endgroup %%------------------------------------<<

If \mG is a torsion abelian group and if 
$G \approx \bigoplus\limits_p G(p)$ is its primary decomposition, then $\forall \ n > 0$, the $H_n(G)$ are torsion and 
$\forall \ p$, 
$H_n(G)(p) \approx$ 
$H_n(G(p))$ ($\implies$ 
$H_n(G) \approx$ 
$\bigoplus\limits_p H_n(G(p))$).

[Note: \ $\forall \ n > 0$, $G(p)$ bounded $\implies$ $H_n(G(p))$ bounded (in fact, 
$p^{M(p)}G(p) = 0$ $\implies$ 
$p^{M(p)}H_n(G(p)) = 0$.]

Example: 
$\Q/\Z \hsx \approx \hsx \bigoplus\limits_p \Z/p^\infty\Z$ 
$\implies$ 
$H_n(\Q/\Z) \hsx \approx \hsx$ 
$\bigoplus\limits_p H_n(\Z/p^\infty\Z)$, where for $n > 0$, 
$H_n(\Z/p^\infty\Z) = \colimx H_n(\Z/p^k\Z) =$ 
$
\begin{cases}
\ \Z/p^\infty\Z \ \hspace{0.5cm} (\text{$n$ odd})\\
\ 0 \hspace{1.65cm}  (\text{$n$ even})
\end{cases}
\hspace{-.25cm}
.
$
\\

\begingroup%%----------------------------------->>
\fontsize{9pt}{11pt}\selectfont
\textbf{\small FACT} \quad 
Fix a prime $p$.  For $k = 1, 2, \ldots$, let $G_k$ be a direct sum of $k$ copies of $\Z/p\Z$ $-$then by the K\"unneth formula, 
$\forall \ n > 0$, $H_n(G_k) = G_{d(n,k)}$, where $d(1,k) = k$ and 
$d(n,k+1) = \ds\sum\limits_{i = 1}^n d(i,k) + (1 - (-1)^n)/2$ 
(hence $d(n,k) \leq k^n)$.\\
\endgroup %%------------------------------------<<

\begingroup%%----------------------------------->>
\fontsize{9pt}{11pt}\selectfont
\textbf{\small FACT} \quad 
Fix a prime $p$.  For $k = 1, 2, \ldots$, let $G_k$ be a direct sum of $k$ copies of $\Z/p^\infty\Z$ 
$-$then by the K\"unneth formula, $\forall \ n > 0$, $H_n(G_k) \ = \ G_{d(n,k)}$, where $d(n,k) = 0$ ($n$ even) and \ 
$\ds d(n,k)  \ = \  \ds\binom{k + \frac{n-1}{2}}{\frac{n+1}{2}}$ \ \ ($n$ odd) (hence $d(n,k) \leq k^n)$.\\
\endgroup %%------------------------------------<<

%%----------------------------------------------------------------------------------------------08
\textbf{\small LEMMA} \quad 
Suppose that $\sC$ is a Serre class.  Let $0 \ra K \ra G \ra G/K \ra 0$ be a short exact sequence in $\sC$ $-$then for $n > 0$, 
$H_n(G) \in \sC$ provided that the $H_p(G/K;H_q(K)) \in \sC$ $(p + q > 0)$.

[Apply the LHS spectral sequence.]

[Note: \ By the universal coefficient theorem, 
$H_p(G/K;H_q(K)) \approx$ 
$H_p(G/K) \otimes H_q(K) \oplus \Tor(H_{p-1}(G/K),H_q(K))$.]\\

\index{Theorem: Theorem of Balcerzyk}
\index{Theorem of Balcerzyk}
\textbf{\small THEOREM OF BALCERZYK} \quad 
Let $\sC$ be a Serre class $-$then $\sC$ is acyclic iff $\sC$ is a ring.

[Suppose that $\sC$ is acyclic.  Since $G$ torsion $\implies$ $H_n(G)$ torsion $(n > 0)$, $\sC_\tor$ is acyclic, thus one can assume that $\sC$ is torsion (cf. Proposition 3) and then, taking into account Proposition 4, work with $\sC^*$ (which is acyclic).  So, let $G \in \sC^*$:
$G \approx$ $\bigoplus\limits_p G(p)$ $\implies$ 
$G \otimes G \approx$ 
$\bigoplus\limits_p G(p) \otimes G(p)$ and 
$\#(G(p) \otimes G(p)) < \omega$ $\implies$ 
$G(p) \otimes G(p) \approx$ 
$H_2(G(p)) \oplus H_2(G(p)) \oplus G(p)$ $\implies$ 
$G \otimes G \approx$ 
$\bigoplus\limits_p (H_2(G(p)) \oplus H_2(G(p)) \oplus G(p)) \approx$ 
$H_2(G) \oplus H_2(G) \oplus G \in \sC^*$.  
To check that $\Tor(G,G) \in \sC^*$, it is obviously enough to look at the case when 
$G \approx$ 
$\bigoplus\limits_p \gamma(p) \cdot (\Z/p^\infty\Z)$, where $\forall \ p$, $\gamma(p) < \omega$.  \ Thus: 
$H_3(G) \approx$ \ 
$\bigoplus\limits_p H_3(\gamma(p) \cdot (\Z/p^\infty\Z))$ 
$\bigoplus\limits_p \ds\binom{\gamma(p)  + 1}{2} \cdot (\Z/p^\infty\Z)$ (cf. supra) and 
$2 \ds\binom{\gamma(p)  + 1}{2} \geq$ 
$\gamma(p)^2$
$\implies$  
$\gamma^2 \in \bgamma(\sC)$ $\implies$ $\Tor(G,G) \in \sC^*$.

Suppose that $\sC$ is a ring.  Let $G \in \sC$ $-$then there is a short exact sequence 
$0 \ra$ 
$G_\tor \ra$ 
$G \ra$ 
$G/G_\tor \ra 0$.  
Accordingly, in view of the lemma, to prove that $H_n(G) \in \sC$ $(n > 0)$, it suffices to prove that 
$H_p(G/G_\tor;H_q(G_\tor)) \in \sC$ $(p + q > 0)$.  But 
$H_p(G/G_\tor;H_q(G_\tor)) \approx$ 
$H_p(G/G_\tor) \otimes H_q(G_\tor) \oplus \Tor(H_{p-1}(G/G_\tor),H_q(G_\tor))$ 
and the verification that 
$H_n(G) \in \sC$ $(n > 0)$ reduces to when 
(i) \mG is torsion free or 
(ii) \mG is torsion.

\indent\indent (Torsion Free) \ If $\tf(\sC)$ is the class of all torsion free abelian groups of cardinality $< \kappa$ 
$(\kappa > \omega)$ (cf. Proposition 1), then $G \in \tf(\sC)$ $\implies$ 
$\#(H_n(G)) < \kappa$ $\implies$  $H_n(G) \in \sC$ $(n > 0)$.  
The other possibility is that $\tf(\sC) = \bT(\bAB)$ for some \bT (cf. Proposition 2).  
Under these circumstances, a given 
$G \in \tf(\sC)$ contains a free subgroup $F \approx r \cdot \Z$ of finite rank such that the sequence 
$0 \ra F \ra G \ra G/F \ra 0$ is exact.  Here, 
$G/F \approx$ 
$\bigoplus\limits_1^r T_i$ is torsion and the $p$-primary components of each $T_i$ are isomorphic to 
$\Z/p^{n_i}\Z$ or $\Z/p^\infty\Z$.  
Therefore 
$H_n(T_i) \approx$
$
\begin{cases}
\ T_i \quad \text{($n$ odd)}\\
\ 0 \  \quad \text{($n$ even)}
\end{cases}
\hspace{-.25cm}
\in \sC \  (n > 0)
$
$\implies$ $H_n(T) \in \sC$ $(n > 0)$ (K\"unneth).  
On the other hand, 
$
H_n(F) \approx
\begin{cases}
\ \ds\binom{r}{n} \cdot \Z \hspace{0.5cm} (n \leq r) \\[8pt]
\ \ \ 0 \hspace{1.15cm} \  \ (n > r)
\end{cases}
\hspace{-.25cm} \in \sC (n > 0).
$
The lemma now implies that $H_n(G) \in \sC$ $(n > 0)$.

%%----------------------------------------------------------------------------------------------09
\indent\indent (Torsion) \ Let $G \in \sC_\tor$.  Choose a basic subgroup \mB of \mG: 
$0 \ra$ 
$B \ra$ 
$G \ra$ 
$G/B \ra 0$ $-$then thanks to the lemma, one need only consider $H_n(B)$ and $H_n(G/B)$ $(n > 0)$.  
Using the cardinal lemma, represent \mB by $B_0 \oplus B_\omega \oplus B_\infty$ with 
$B_0(p) =$ $\bigoplus\limits_n \alpha_0(p,n) \cdot (\Z/p^n\Z)$, 
$B_\omega(p) =$ $\bigoplus\limits_n \alpha_\omega(p,n) \cdot (\Z/p^n\Z)$, 
and
$B_\infty(p) =$ $\bigoplus\limits_n \alpha_\infty(p,n) \cdot (\Z/p^n\Z)$, subject to 
$(\alpha_0)$ $\forall \ p$, $\sum\limits_n \alpha_0(p,n) < \omega$, 
$(\alpha_\omega)$ $\forall \ p$, 
$\exists \ M(p)$: $n \geq M(p)$ $\implies$ $\alpha_\omega(p,n) = 0$ $\&$ $\forall \ n$:
$\alpha_\omega(p,n) \geq \omega$ or $\alpha_\omega(p,n) = 0$, and 
$(\alpha_\infty)$ $\exists$ $\gamma_\infty \in \bgamma(\sC)$: $\forall \ p$, $\forall \ n$,  
$\alpha_\infty(p,n) \leq \gamma_\infty(p)$, where $\gamma_\infty(p) \geq \omega$ or $\gamma_\infty(p) = 0$.  
That $H_n(B) \in \sC$ $(n > 0)$ results from the following observations (modulo K\"unneth): 
(O$_0$) \ $\forall \ p$, $\#(B_0(p)) < \omega$, hence there is a monomorphism 
$H_n(B_0(p)) \ra \otimes^n B_0(p)$; 
(O$_\omega$) \ $\forall \ p$, $\forall \ \alpha \geq \omega$, 
$H_n(\alpha \cdot (\Z/p^k\Z)) \approx$ 
$\alpha \cdot (\Z/p^k\Z)$; 
(O$_\infty$) \ $\forall \ p$, $\#(B_\infty(p)) \leq \gamma_\infty(p)$, hence there is a monomorphism 
$H_n(B_\infty(p)) \ra \gamma_\infty(p) \cdot (\Z/p^\infty\Z)$.  
Finally, write 
$G/B \approx$ 
$\bigoplus\limits_p \gamma(p) \cdot (\Z/p^\infty\Z)$ and fix $n > 0$.  
Case 1: $n$ even $\implies$ $H_n(G/B) = 0$.  
Case 2: $n$ odd.  If $\gamma(p) \geq \omega$, then 
$H_n(\gamma(p) \cdot (\Z/p^\infty\Z)) \approx$ 
$\gamma(p) \cdot (\Z/p^\infty\Z)$, while if $\gamma(p) < \omega$, then 
$H_n(\gamma(p) \cdot (\Z/p^\infty\Z)) \approx$ 
$\ds\binom{\gamma(p) + \frac{n-1}{2}}{\frac{n+1}{2}} \cdot (\Z/p^\infty\Z)$ (cf supra).  However, 
$\ds\binom{\gamma(p) + \frac{n-1}{2}}{\frac{n+1}{2}} \leq$ 
$(\gamma(p))^n$ and $\gamma^n \in \bgamma(\sC)$.]\\

\label{7.3}
\begingroup%%----------------------------------->>
\fontsize{9pt}{11pt}\selectfont
\textbf{\small EXAMPLE} \quad 
Let $\sC$ be a ring.  Fix a nilpotent group \mG such that $G/[G,G] \in \sC$ $-$then $\forall \ n > 0$, $H_n(G) \in \sC$.\\
\endgroup %%------------------------------------<<

\label{5.31a}
\begingroup%%----------------------------------->>
\fontsize{9pt}{11pt}\selectfont
\textbf{\small FACT} \quad 
Let $\sC$ be a ring.  Suppose that \mX is simply connected $-$then $H_q(X) \in \sC$ $\forall \ q > 0$ iff 
$H_q(\Omega X) \in \sC$ $\forall \ q > 0$.\\
\endgroup %%------------------------------------<<

\begingroup%%----------------------------------->>
\fontsize{9pt}{11pt}\selectfont
Application: Let $\sC$ be a ring.  Fix $\pi \in \sC$ $-$then the $H_q(\pi,n) \in \sC$ $(q > 0)$.\\
\endgroup %%------------------------------------<<

\label{8.2}
If $\sC$ is a Serre class, then a homomorphism $f:G \ra K$ of abelian groups is said to be 
\un{$\sC$-injective}
\index{C-injective} %\index$\sC$-injective} 
(\un{$\sC$-surjective})
\index{C-surjective} %\index$\sC$-surjective} 
if the kernel (cokernel) of $f$ is in $\sC$, $f$ being 
\un{$\sC$-bijective}
\index{C-bijective} %\index$\sC$-bijective} 
provided that it is both $\sC$-injective and $\sC$-surjective.\\

\index{Theorem: Mod $\sC$ Hurewicz Theorem}
\index{Mod $\sC$ Hurewicz Theorem}
\textbf{\small MOD $\sC$ HUREWICZ THEOREM} \quad 
Let $\sC$ be a Serre class.  Assume: $\sC$ is a ring.  Suppose that \mX is abelian $-$then if $n \geq 2$, the condition 
$\pi_q(X) \in \sC$ $(1 \leq q < n)$ is equivalent to the condition $H_q(X) \in \sC$ $(1 \leq q < n)$ and either implies that the Hurewicz homomorphism $\pi_n(X) \ra H_n(X)$ is $\sC$-bijective.\\

\begingroup%%----------------------------------->>
\fontsize{9pt}{11pt}\selectfont
\textbf{\small EXAMPLE} \quad 
Let $\sC$ be a ring.  Suppose that \mX is a pointed connected CW space which is nilpotent.  Agreeing to write 
$\pi_1(X) \in \sC$ if $\pi_1(X) /[\pi_1(X) ,\pi_1(X) ] \in \sC$, fix $n \geq 2$ $-$then the following conditions are 
%%----------------------------------------------------------------------------------------------10
equivalent:
(i)  $\pi_q(X) \in \sC$ $(1 \leq q < n)$; 
(ii) $H_q(X) \in \sC$  $(1 \leq q < n)$; 
(iii) $\pi_1(X) \in \sC$ $\&$ $H_q(\widetilde{X}) \in \sC$ $(1 \leq q < n)$.  
Furthermore, under (i), (ii), or (iii), the Hurewicz homomorphism $\pi_n(X) \ra H_n(X)$ induces a $\sC$-bijection 
$\pi_n(X)_{\pi_1(X)} \ra H_n(X)$.
\\ \indent
 [To illustrate the line of argument, assume (iii) and consider the spectral sequence 
$E_{p,q}^2 \approx H_p(\pi_1(X);$ $H_q(\widetilde{X}))$ $\Rightarrow$ $H_{p+q}(X)$ of the covering projection 
$\widetilde{X} \ra X$ 
(cf. p. \pageref{7.2}).  
Since $\pi_1(X) \in \sC$ is nilpotent, 
$E_{p,0}^2 \in \sC$ $(p > 0)$ 
(cf. p. \pageref{7.3}).  In addition, the $H_q(\widetilde{X})$ $(q > 0)$ are nilpotent 
$\pi_1(X)$-modules (cf. $\S 5$, Proposition 17), thus 
$E_{p,q}^2 \in \sC$ $(p \geq 0, 1 \leq q < n)$, 
(cf. p. \pageref{7.4}) $\implies$ $H_q(X) \in \sC$ $(1 \leq q < n)$ and there is a $\sC$-bijection 
$H_n(\widetilde{X})_{\pi_1(X)} \ra H_n(X)$.  
Owing to the mod $\sC$ Hurewicz theorem, $\pi_q(X) \approx \pi_q(\widetilde{X}) \in \sC$ $(2 \leq q < n)$ and the 
Hurewicz homomorphism $\pi_n(\widetilde{X}) \ra H_n(\widetilde{X})$ is $\sC$-bijective.  But then the arrow 
$\pi_n(\widetilde{X})_{\pi_1(X)} \ra H_n(\widetilde{X})_{\pi_1(X)}$ is also $\sC$-bijective, $\pi_n(\widetilde{X})$ and 
$H_n(\widetilde{X})$ being nilpotent $\pi_1(X)$-modules.]\\
\endgroup %%------------------------------------<<

A Serre class $\sC$ is said to be an 
\un{ideal}
\index{ideal (Serre class)}
if $G \in \sC$ $\implies$ $G \otimes K \in \sC$, $\Tor(G,K) \in \sC$ for all \mK in \bAB.\\

\textbf{\small LEMMA} \quad 
Let $\sC$ be a Serre class $-$then $\sC$ is an ideal iff $\forall \ G \in \sC$, $\bigoplus\limits_i G_i \in \sC$, where 
$\bigoplus\limits_i$ is taken over any index set and $\forall \ i$, $G_i \approx G$.\\

Example: Let $\sC$ be an ideal.  Suppose that $\sG \in [(\sin X)^\OP,\bAB]$ is a coefficient system on \mX such that 
$\forall \ \sigma$, $\sG \sigma \in \sC$ $-$then $\forall \ n \geq 0$, $H_n(X;\sG) \in \sC$.\\

\begingroup%%----------------------------------->>
\fontsize{9pt}{11pt}\selectfont
\textbf{\small EXAMPLE} \quad 
The conglomerate of torsion Serre classes which are ideals is codable by a set.  For in the set of subsets of 
$F(\N,\Z_{\geq 0} \cup \{\infty\})$, write $S \sim T$ iff each sequence in $S$ is $ \leq$ a finite sum of sequences in \mT and each sequence in \mT is $\leq$ a finite sum of sequences in \mS.  
Let \mE be the resulting set of equivalence classes.  
Claim: The conglomerate of torsion ideals is in a one-to-one corrspondence with \mE.  
Thus given a torsion ideal $\sC$, assign to 
$G \in \sC$  the sequence $\{x_n(G)\} \in F(\N,\Z_{\geq 0} \cup \{\infty\})$ by letting $x_n(G)$ be the least uppper bound of the exponents of the elements in $G(p_n)$, where $\forall \ n$, $p_n < p_{n+1}$.  
Put $S_{\sC} = \{\{x_n(G)\}:G \in \sC\}$ and call $[S_{\sC}] \in E$ the equivalence class corresponding to $S_{\sC}$.  
To go the other way, take an \mS and let $\sC_S$ be the class of torsion abelian groups \mG with the property that there exists a finite number of sequences in \mS such that the $n^\text{th}$ term of their sum is an upper bound on the exponents of the elements in $G(p_n)$ $-$then $\sC_S$ is an ideal and $S \sim T$ $\implies$ $\sC_S = \sC_T$, so $\sC_{[S]}$ makes sense.  
One has 
$\sC \ra [S_{\sC}] \ra \sC_{[S_{\sC}]} = \sC$ and 
$[S] \ra \sC_{[S]} \ra [S_{\sC_{[S]}}] = [S]$.
\\ \indent
[Note: \ It is sufficient to consider torsion ideals since any ideal containing a nonzero torsion free group is necessarily the class of all abelian groups.]\\
\endgroup %%------------------------------------<<

\index{Theorem: Mod $\sC$ Whitehead Theorem}
\index{Mod $\sC$ Whitehead Theorem}
\textbf{\small MOD $\sC$ WHITEHEAD THEOREM} \quad 
Let $\sC$ be a Serre class.  Assume: $\sC$ is an ideal.  Suppose that \mX and \mY are abelian and $f:X \ra Y$ is a continuous function 
$-$then if $n \geq 2$, the condition $f_*:\pi_q(X) \ra \pi_q(Y)$ is $\sC$-bijective for $1 \leq q < n$ and $\sC$-surjective for
%%----------------------------------------------------------------------------------------------11
$q = n$ is equivalent to the condition $f_*:H_q(X) \ra H_q(Y)$ is $\sC$-bijective for $1 \leq q < n$ and $\sC$-surjective for 
$q = n$.\\

\begingroup%%----------------------------------->>
%\setstretch{1.}
\fontsize{9pt}{11pt}\selectfont
\textbf{\small EXAMPLE} \quad 
Let 
$
\begin{cases}
\ X\\
\ Y
\end{cases}
$
be abelian.  
Assume: $\forall \ q$, 
$
\begin{cases}
\ H_q(X)\\
\ H_q(Y)
\end{cases}
$
is finitely generated $(\implies$ $\forall \ q$, 
$
\begin{cases}
\ \pi_q(X)\\
\ \pi_q(Y)
\end{cases}
$
is finitely generated).
\\ 
\indent\indent $(\text{char \bk}= 0)$ \ 
Let $f:X \ra Y$ be a continuous function.  
Fix a field \bk of characteristic 0 and denote by $\sF$ the class of finite abelian groups, 
$\sT$ the class of torsion abelian groups $-$then if $n \geq 2$, the following conditions are equivalent: 
(1) $f_*: H_q(X) \ra H_q(Y)$ is $\sF$-bijective for $1 \leq q < n$ and $\sF$-surjective for $q = n$;
(2) $f_*: H_q(X) \ra H_q(Y)$ is $\sT$-bijective for $1 \leq q < n$ and $\sT$-surjective for $q = n$;
(3) $f_*: H_q(X;\bk) \ra H_q(Y;\bk)$ is bijective for $1 \leq q < n$ and surjective for $q = n$;
(4) $f^*: H^q(Y;\bk) \ra H^q(X;\bk)$ is bijective for $1 \leq q < n$ and injective for $q = n$.
\\ 
\indent\indent $(\text{char \bk} = p)$ \ 
Let $f:X \ra Y$ be a continuous function.  
Fix a field \bk of characteristic $p$ and denote by $\sF_p$ the class of finite abelian groups with order prime to $p$, 
$\sT_p$ the class of torsion abelian groups with trivial $p$-primary 
component $-$then if $n \geq 2$, the following conditions are equivalent: 
(1) $f_*: H_q(X) \ra$ $H_q(Y)$ is $\sF_p$-bijective for $1 \leq q < n$ and $\sF_p$-surjective for $q = n$;
(2) $f_*: H_q(X) \ra H_q(Y)$ is $\sT_p$-bijective for $1 \leq q < n$ and $\sT_p$-surjective for $q = n$;
(3) $f_*: H_q(X;\bk) \ra H_q(Y;\bk)$ is bijective for $1 \leq q < n$ and surjective for $q = n$;
(4) $f^*:H^q(Y;\bk) \ra H^q(X;\bk)$ is bijective for $1 \leq q < n$ and injective for $q = n$.
\\ \indent
Example: If $\forall \ n$, $f_*$ induces an isomorphism $H_n(X;\F_p) \ra H_n(Y;\F_p)$, then $\forall \ n$, $f_*$ induces an isomorphism $\pi_q(X)(p) \ra \pi_n(Y)(p)$ of $p$-primary components.\\
\endgroup %%------------------------------------<<

\begingroup%%----------------------------------->>
\fontsize{9pt}{11pt}\selectfont
\textbf{\small FACT} \quad 
Let \mX be a CW complex.  Assume: \mX  is finite and $n$-connected $-$then the Hurewicz homomorphism 
$\pi_q(X) \ra H_q(X)$ is $\sC$-bijective for $q < 2n + 1$, where $\sC$ is the class of finite abelian groups.\\
\endgroup %%------------------------------------<<

%%%%%%%%%%%%%%%%%%%%%%%%%%%%%%%%%%%%%%
%%%%%%%%%%%%%%%%%%%%%%%%%%%%%%%%%%%%%%
%%%%%%%%%%%%%%%%%%%%%%%%%%%%%%%%%%%%%%

\begin{center}
$\S \ 7$
\\[0.5cm]
$\mathcal{REFERENCES}$\\
\end{center}

\[
\textbf{BOOKS}
\]

\begingroup
\fontsize{9pt}{11pt}\selectfont
\setlength\parindent{0 cm}

[1] \quad Hu, S., \textit{Homotopy Theory}, Academic Press (1959).
\\[-.2cm]

[2] \quad Mosher, R., and Tangora, M., \textit{Cohomology Operations and Applications in Homotopy Theory}, 

\hspace{0.8cm}Harper $\&$
Row (1968).
\\[-.2cm]
\endgroup

\[
\textbf{ARTICLES}
\]

\begingroup
\fontsize{9pt}{11pt}\selectfont
\setlength\parindent{0 cm}

[1] \quad Balcerzyk, S., On Classes of Abelian Groups, 
\textit{Fund. Math.} \textbf{51} (1962), 149-178; see also \textit{Fund.} 

\hspace{0.8cm}\textit{Math.} \textbf{56} (1964), 199-202.
\\[-.2cm]

[2] \quad Goncalves, D., Generalized Classes of Groups, $\C$-Nilpotent Spaces, and "The Hurewicz Theo-

\hspace{0.8cm}rem", 
\textit{Math. Scand.} \textbf{53} (1983), 39-61.
\\[-.2cm]

[3] \quad Serre, J-P., Groupes d'Homotopie et Classes de Groupes Ab\'eliens, 
\textit{Ann. of Math.} \textbf{58} (1953), 

\hspace{0.8cm}258-294.
\\[-.2cm]

[4] \quad Walker, C., and Walker, E., Quotient Categories and Rings of Quotients, 
\textit{Rocky Mountain J.}

\hspace{0.8cm}\textit{Math.} \textbf{2} (1972), 513-555.

\setlength\parindent{2em}

\endgroup

\chapter{
$\boldsymbol{\S}$\textbf{8}.\quadx  LOCALIZATION OF GROUPS}
\setlength\parindent{2em}
\setcounter{proposition}{0}
\setcounter{chapter}{8}

%%----------------------------------------------------------------------------------------------01
$\text{ }$\\[-1.25cm]

The algebra of this section is the point of departure for the developments in the next $\S$.  
While the primary focus is on the ``abelian-nilpotent'' theory, part of the material is presented in a more general setting.  
I have also included some topological applications that will be of use in the sequel.

The Serre classes in \bAB that are closed under the formation of coproducts (and hence colimits) 
are in a one-to-one correspondence with the Giraud subcategories of \bAB.  
Under this correspondece, the class of all abelian groups corresponds to the class of trivial groups.  
The remaining classes are necessarily torsion ideals and their determination is embedded in abelian localization theory.

[Note: \ Not every torsion ideal is closed under the formation of coproducts (consider, e.g., the class of bounded abelian groups).]

Notation: \mP is a set of primes, $\ov{P}$ its complement in the set of all primes.
 
Given \mP, put $S_P = \{1\} \cup \{n > 1:p \in P \implies p \dnd n \}$ 
$-$then $\Z_P = S_P^{-1}\Z$ is the localization of $\Z$ at $P$ and the inclusion $\Z \ra \Z_P$ is an epimorphism in \bRG.  
$\Z_P$ is a principal ideal domain.  
Moreover, $\Z_P$ is a subring of $\Q$ and every subring of $\Q$ is a $\Z_P$ for a suitable $P$.  
The characteristic of 1 in $\Z_P$ is
$
\begin{cases}
\ 0 \quad (p \in P)\\
\ \infty \  (p \in \ov{P})
\end{cases}
$
$\implies$ $\Z_P/\Z \approx \bigoplus\limits_{p \in \ov{P}} \Z/p^\infty \Z$.  
Examples: 
(1) Take $P = \emptyset$: $\Z_P = \Q$; 
(2) Take $P = \bPi$: $\Z_P = \Z$; 
(3) Take $P = \bPi - \{p\}$ : $\Z_P = \Z\left[\ds\frac{1}{p}\right]$;
(4) Take $P = \bPi - \{2,5\}$ : $\Z_P =$ all rationals whose decimal expansion is finite.

[Note: \ Write $\Z_p$ in place of $\Z_{\{p\}}$: $\Z_p$ is a local ring and its residue field is isomorphic to $\F_p$.]\\

\label{8.5}
\begingroup%%----------------------------------->>
\fontsize{9pt}{11pt}\selectfont
\textbf{\small EXAMPLE} \  
Suppose that $P \neq \emptyset$ $-$then as vector spaces over $\Q$, $\Ext(\Q,\Z_P) \approx \R$.\\
\endgroup %%------------------------------------<<

Equip $S_P$ with the structure of a directed set by stipulating that $n^\prime \leq n\pp$ iff $n^\prime| n\pp$.  
View $(S_P,\leq)$ as a filtered category $\bS_P$ and let $\Delta_P:\bS_P \ra \bAB$ be the diagram that sends an object $n$ to $\Z$ and a morphism $n^\prime \ra n\pp$ to the multiplication $\ds\frac{n\pp}{n^\prime}:\Z \ra \Z$ $-$then the homomorphism 
$\colimx \Delta_P \ra \Z_P$ is an isomorphism.  
Example: $\Z_P \otimes \Z/p^n\Z = $ 
$
\begin{cases}
\ 0 \hspace{1.5cm} (p \in \ov{P})\\
\ \Z/p^n\Z \ \hspace{0.5cm}  (p \in P)
\end{cases}
\hspace{-.25cm}.
$
\\

\label{8.32}
\label{9.3}
\begingroup%%----------------------------------->>
\fontsize{9pt}{11pt}\selectfont
\textbf{\small EXAMPLE} \  
Fix $P \neq \bPi$ $-$then there is a short exact sequence 
$0 \ra$ 
$\lim^1 H^1(\Z;\Q([\Z_P]) \ra$ 
$H^2(\Z_P;\Q([\Z_P]) \ra$ 
$\lim H^2(\Z;\Q([\Z_P]) \ra 0$.  Here, $H^2(\Z_P;\Q([\Z_P]) \neq 0$ (in fact, is uncountable 
(cf. p. \pageref{8.1})).\\ 
\endgroup %%------------------------------------<<

%%----------------------------------------------------------------------------------------------02
\textbf{\small LEMMA} \  
Let $P^\prime$ and $P\pp$ be two sets of primes $-$then 
(i) $\Z_{P^\prime} + \Z_{P\pp} = \Z_{P^\prime \cap P\pp}$ and 
(ii) $\Z_{P^\prime} \cap \Z_{P\pp} = \Z_{P^\prime \cup P\pp}$, the sum and intersection being as subgroups of $\Q$.  Furthermore, the biadditive function 
$
\begin{cases}
\ \Z_{P^\prime} \times \Z_{P\pp} \ra \Z_{P^\prime \cap P\pp}\\
\ (z^\prime,z\pp) \ra z^\prime z\pp
\end{cases}
$
defines an isomorphism of rings: 
$\Z_{P^\prime} \otimes \Z_{P\pp} \approx$ 
$\Z_{P^\prime \cap P\pp}$ $(\implies \Z_P \otimes \Z_P \approx \Z_P)$.\\

\begingroup%%----------------------------------->>
\fontsize{9pt}{11pt}\selectfont
\textbf{\small FACT} \  
There is a \cd \ 
\begin{tikzcd}[sep=large]
{\Z_{P^\prime \cup P\pp}} \ar{d}[swap]{i^\prime} \ar{r}{i\pp} &{\Z_{P\pp}} \ar{d}{j\pp}\\
{\Z_{P^\prime}} \ar{r}[swap]{j^\prime} &{\Z_{P^\prime \cap P\pp}}
\end{tikzcd}
and a short exact sequence \ \ 
$0 \ra$ 
$\Z_{P^\prime \cup P\pp} \overset{\mu}{\ra}$ 
$\Z_{P^\prime} \oplus \Z_{P\pp} \overset{\nu}{\ra}$ 
$\Z_{P^\prime \cap P\pp} \ra 0$ 
$(\mu(z) = (i^\prime(z),i\pp(z))$ $\&$ $\nu(z^\prime,z\pp) = j^\prime(z^\prime) - j\pp(z\pp))$, thus the square is simultaneously a pullback and a pushout in \bAB.\\
\endgroup %%------------------------------------<<

An abelian group \mG is said to be 
\un{$S_P$-torsion}
\index{S$_P$-torsion} 
if $\forall \ g \in G$, $\exists$ 
$n \in S_P$ : $ng = 0$.  
Denote by $\sC_P$ the class of $S_P$-torsion abelian groups $-$then $\sC_P$ is a Serre class which is closed under the formation of coproducts and every torsion Serre class with this property is a $\sC_P$ for some \mP.  
Examples: 
(1) Take $P = \emptyset$ : $\sC_P$ is the class of torsion abelian groups; 
(2) Take $P = \bPi$ : $\sC_P$ is the class of trivial groups;  
(3) Take $P = \{p\}$ : $\sC_P$ is the class of torsion abelian groups with trivial $p$-primary component; 
(4) Take $P = \bPi - \{p\}$ : $\sC_P$ is the class of abelian $p$-groups.

[Note: \ \mG is $S_P$-torsion iff \mG is $\ov{P}$-primary or still, iff $\Z_P \otimes G = 0$.]

Let $f:G \ra K$ be a homomorphism of abelian groups $-$then $f$ is said to be 
\un{$P$-injective}
\index{P-injective} 
(\un{$P$-surjective}
\index{P-surjective})
if the kernel (cokernel) of $f$ is $S_P$-torsion, $f$ being
\un{$P$-bijective}
\index{P-bijective} 
provide that it is both $P$-injective and $P$-surjective.

[Note: \ This is the terminlogy on 
p. \pageref{8.2}, specialized to the case $\sC = \sC_P$.]\\

\index{Five Lemma}
\textbf{\small FIVE LEMMA} \quad 
Let 
\[
\begin{tikzcd}%[ sep=small]
{G_1} \ar{d}{f_1} \ar{r} 
&{G_2} \ar{d}{f_2} \ar{r}
&{G_3} \ar{d}{f_3} \ar{r}
&{G_4} \ar{d}{f_4} \ar{r}
&{G_5} \ar{d}{f_5}\\
{K_1} \ar{r} 
&{K_2} \ar{r}
&{K_3} \ar{r}
&{K_4} \ar{r}
&{K_5}
\end{tikzcd}
\]
be a commutative diagram of abelian groups with exact rows.\\
\indent\indent (1) \quad If $f_2$ and $f_4$ are $P$-surjective and $f_5$ is $P$-injective, then $f_3$ is $P$-surjective.\\
\indent\indent (2) \quad If $f_2$ and $f_4$ are $P$-injective and $f_1$ is $P$-surjective, then $f_3$ is $P$-injective.\\

\begingroup%%----------------------------------->>
\fontsize{9pt}{11pt}\selectfont
The definition of ``$S_P$-torsion'' carries over without changes to \bGR, as does the definition of ``$P$-injective'' but it is best to modify the definition of ``$P$-surjective''.  
Thus let $f:G \ra K$ be a homomorphism of groups $-$then $f$ is said to be 
\un{$P$-surjective}
\index{P-surjective (homomorphism of groups)} 
if $\forall \ k \in K$, $\exists$ $n \in S_P$: $k^n \in \im f$ 
(when \mG and \mK are nilpotent, this is equivalent to requiring that $\coker f$ be $S_P$-torsion).  Assigning to the term 
``$P$-bijective''
%%----------------------------------------------------------------------------------------------03
the obvious interpretation, the five lemma retains its validity under the following additional assumptions: 
$(1)_+$ $\im (K_2 \ra K_3) \subset \Cen \ K_3$ or 
$(2)_+$ $\im(G_1 \ra G_2) \subset \Cen \ G_2$ (no extra conditions are needed in the nilpotent case).\\
\endgroup %%------------------------------------<<

\label{9.20} %dmc mnft
Given an abelian group \mG, 
\un{the localization of \mG at \mP}
\index{localization of \mG at \mP} 
is the tensor product 
$G_P = \Z_P \otimes G$.  The functor 
$Z_P \otimes -: \bAB \ra \Z_P\text{-}\bMOD$ preserves colimits and is exact.  
Examples: 
(1)  Suppose that \mG is finitely generated, say 
$G \approx \bigoplus\limits_1^r \Z \oplus \bigoplus\limits_p \bigoplus\limits_n \alpha(p,n) \cdot (\Z/p^n\Z)$ $-$then 
$G_P \approx \bigoplus\limits_1^r \Z_P \oplus \bigoplus\limits_{p \in P} \bigoplus\limits_n \alpha(p,n) \cdot (\Z/p^n\Z)$; 
(2)  Suppose that \mG is torsion, say 
$G = \bigoplus\limits_p G(p)$ $-$then 
$G_P = \bigoplus\limits_{p \in P} G(p)$.

[Note: \ $G_{\Q} = \Q \otimes G$ is the 
\un{rationalization}
\index{rationalization} 
of \mG.  
Example: $\Q \otimes \Z^\omega \neq \Q^\omega$.  $G_p = \Z_p \otimes G$ is the 
\un{$p$-localization}
\index{p-localization} of \mG.  
Example: $(\Q/\Z)_p = \Z/p^\infty\Z$.]\\

\begingroup%%----------------------------------->>
\fontsize{9pt}{11pt}\selectfont
\textbf{\small FACT} \  
Let \mG be an abelian group $-$then the \cd 
\begin{tikzcd}%[ sep=small]
{G} \ar{d} \ar{r} &{G_{\ov{P}}} \ar{d}\\
{G_P} \ar{r} &{G_{\Q}}
\end{tikzcd}
is simultaneously a pullback square and a pushout square in \bAB and the arrow 
$
\begin{cases}
\ G_P \ra G_{\Q}\\
\ G_{\ov{P}} \ra G_{\Q}
\end{cases}
$
is a 
$
\begin{cases}
\ \ov{P} \text{-bijection}\\
\ P \text{-bijection}
\end{cases}
$
.\\
\endgroup %%------------------------------------<<
\vspace{0.25cm}

\begingroup%%----------------------------------->>
\fontsize{9pt}{11pt}\selectfont
\textbf{\small FACT} \  
Let \mG be an abelian group $-$then \mG is finitely generated iff 
$
\begin{cases}
\ G_P\\
\ G_{\ov{P}}
\end{cases}
$
are finitely generated 
$
\begin{cases}
\ \Z_P\\
\ \Z_{\ov{P}}
\end{cases}
\text{-modules.}
$
\vspi
[Note: \ $\bigl(\ds\bigoplus\limits_p \Z/p\Z\bigr)_q$ is a finitely generated $\Z_q$-module for every prime $q$ but 
$\ds\bigoplus\limits_p \Z/p\Z$ is not a finitely generated abelian group.]\\
\endgroup %%------------------------------------<<

\begingroup%%----------------------------------->>
\label{17.51}
\fontsize{9pt}{11pt}\selectfont
\textbf{\small FACT} \  
Let \mG be an abelian group $-$then $G = 0$ iff $\forall \ p$, $G_p = 0$.\\
\endgroup %%------------------------------------<<

\begingroup%%----------------------------------->>
\fontsize{9pt}{11pt}\selectfont
\textbf{\small FACT} \  
Let 
$
\begin{cases}
\ G\\
\ K
\end{cases}
$
be finitely generated abelian groups.  Assume: $\forall \ p$, $G_p \approx K_p$ $-$then $G \approx K$.
\vspi
[Note: \ To see the failure of this conclusion when one of \mG and \mK is not finitely generated, take $G = \Z$ and let \mK be the additive subgroup of $\Q$ consisting of those rationals of the form $m/n$, 
where $n$ is square free $-$then $\forall \ p$, 
$G_p \approx K_p$, yet $G \not\approx K$.  
Replacing ``square free'' by ``$k^\text{th}$ -power free'', it follows that there exist infinitely many mutually nonisomorphic abelian groups whose $p$-localization is isomorphic to $\Z_p$ at every prime $p$.]\\
\endgroup %%------------------------------------<<

\begingroup%%----------------------------------->>
\fontsize{9pt}{11pt}\selectfont
\textbf{\small FACT} \  
Let $f:G \ra K$ be a homomorphism of abelian groups $-$then $f$ is injective (surjective) iff $\forall \ p$, $f_p:G_p \ra K_p$ is injective (surjective).
\vspi
[Localization is an exact functor, hence preserves kernels and cokernels.]\\ 
\endgroup %%------------------------------------<<

%%----------------------------------------------------------------------------------------------04
\begingroup%%----------------------------------->>
\fontsize{9pt}{11pt}\selectfont
\textbf{\small FACT} \  
Let $f, g:G \ra K$ be homomorphisms of abelian groups.  Assume: $\forall \ p$, $f_p = g_p$, $-$then $f = g$.

[The vertical arrows in the \cd
\begin{tikzcd}[ sep=large]
{G} \ar{d} \ar{r}{f}\ar{r}[swap]{g} &{K} \ar{d}\\
{\prod G_p}  \ar{r}{\prod f_p}\ar{r}[swap]{\prod g_p} &{\prod K_p}
\end{tikzcd}
are one-to-one.]\\
\endgroup %%------------------------------------<<

\textbf{\small LEMMA} \  
Let $G_\tor$ be the torsion subgroup of \mG $-$then $(G_\tor)_P$ is the torsion subgroup of $G_P$.\\

\begingroup%%----------------------------------->>
\fontsize{9pt}{11pt}\selectfont
\textbf{\small EXAMPLE} \  
Take $G = \ds\prod\limits_p \Z/p\Z$ $-$then 
$G_\tor \approx \ds\bigoplus\limits_p \Z/p\Z$ $\implies$ $(G_p)_\tor \approx \Z/p\Z$, so $\forall \ p$, $(G_p)_\tor$ is a direct summand of $G_p$, yet $G_\tor$ is not a direct summand of \mG.\\
\endgroup %%------------------------------------<<

Let \mG be an abelian group $-$then one may attach to \mG a sink 
$\{r_p:G_p \ra G_\Q\}$ and a source 
$\{l_p:G \ra G_P\}$, where 
$
\forall \ 
\begin{cases}
\ p\\[-.1cm]
\ q
\end{cases}
, \ r_p \circx l_p = r_q \circx l_q.
$
\\

\label{8.17}
\index{Fracture Lemma}
\textbf{\small FRACTURE LEMMA} \  
Suppose that \mG is a finitely generated abelian group $-$then the source $\{l_p:G \ra G_p\}$ is the multiple pullback of the sink $\{r_p:G_p \ra G_\Q\}$.

[It suffices to look at two cases: 
(i) $G = \Z/p^n\Z$ and 
(ii) $G = \Z$.]\\

\begingroup%%----------------------------------->>
\fontsize{9pt}{11pt}\selectfont
\textbf{\small EXAMPLE} \  
Take $G = \ds\bigoplus\limits_p \Z/p\Z$ $-$then $G_p = \Z/p\Z$ and $G_{\Q} = 0$, the final object in \bAB.  
Accordingly, the multiple pullback of the sink $\{\Z/p\Z \ra 0\}$ is the source $\bigg\{\ds\prod\limits_p \Z/p\Z \ra \Z/p\Z \bigg\}$.\\
\vspace{0.25cm}
\endgroup %%------------------------------------<<

An abelian group \mG is said to be 
\un{$P$-local}
\index{P-local (abelian group)} 
if the map 
$
\begin{cases}
\ G \ra G\\[-.1cm]
\ g \ra ng
\end{cases}
$
is bijective $\forall \ n \in S_P$.  $\bAB_P$ is the full subcategory of \bAB whose objects are the $P$-local abelian groups.  
$\bAB_P$ is a Giraud subcategory of \bAB with exact reflector 
$
L_P: 
\begin{cases}
\ \bAB \ra \bAB_P\\
\ G \ra G_P
\end{cases}
$
and arrow of localization $l_P:G \ra G_P$.  
Therefore \mG is $P$-local iff $l_P$ is an isomorphism.  
In general, the kernel and cokernel of $l_P:G \ra G_P$ are $S_P$-torsion, i.e., $l_P$ is $P$-bijective.

[Note: \ The objects of $\bAB_P$ are the uniquely $\ov{P}$-divisible abelian groups.  
Changing the notation momentarily, let $S_P \subset \Mor \bAB$ be the class consisting of those $s$ such that 
$\ker s \in \sC_P$ and 
$\coker s \in \sC_P$ $-$then the localization $S_P^{-1}\bAB$ is equivalent to $\bAB_P$ and the endomorphism ring of $\Z$, considered as an object in $S_P^{-1}\bAB$, is isomorphic to $\Z_P$.  
Moreover, a homomorphism $f:G \ra K$ of abelian groups is $P$-bijective iff $f_P:G_P \ra K_P$ is bijective.]\\

\label{9.4}
\index{Recognition Principal}
\textbf{\small RECOGNITION PRINCIPAL} \  
Let \mG be an abelian group $-$then \mG is $P$-local iff it carries the structure of a $\Z_P$-module or satisfies one of the following equivalent conditions.\\
%%----------------------------------------------------------------------------------------------05
\indent\indent (REC$_1$) \quad $\Z_P/\Z \otimes G = 0$ $\&$ $\Tor(\Z_P/\Z,G) = 0$.\\
\indent\indent (REC$_2$) \quad $\forall \ n \in S_P$, $\Z/n\Z \otimes G = 0$ $\&$ $\Tor(\Z/n\Z,G) = 0$.\\
\indent\indent (REC$_3$) \quad $\Hom(\Z_P/\Z,G) = 0$ $\&$ $\Ext(\Z_P/\Z,G) = 0$.\\
\indent\indent (REC$_4$) \quad $\forall \ n \in S_P$, $\Hom(\Z/n\Z,G) = 0$ $\&$ $\Ext(\Z/n\Z,G) = 0$.

[Note: \ In REC$_2$ or REC$_4$, one can just as well work with $p \in \ov{P}$.]\\

\begingroup%%----------------------------------->>
\fontsize{9pt}{11pt}\selectfont
\textbf{\small FACT} \  
Let \mG be an abelian group.  Suppose that \mG is isomorphic to a subgroup of a 
$P$-local abelian group and a quotient group of a 
$P$-local abelian group $-$then \mG is $P$-local.\\
\endgroup %%------------------------------------<<

\begingroup%%----------------------------------->>
\fontsize{9pt}{11pt}\selectfont
\textbf{\small FACT} \  
Let $0 \ra G^\prime \ra G \ra G\pp \ra 0$ be a short exact sequence of abelian groups.  
Assume: Two of the groups are $P$-local $-$then so is the third.
\vspi
[Note: \ $\bAB_P$ is closed with respect to the formation of five term exact sequences but this need not be true of three term exact sequences unless \mP is the set of all primes, this being the only case when $\bAB_P$ is a Serre class.]\\
\endgroup %%------------------------------------<<

\begingroup%%----------------------------------->>
\fontsize{9pt}{11pt}\selectfont
\textbf{\small EXAMPLE} \  
The homology groups attached to a chain complex of $P$-local abelian groups are $P$-local.\\
\endgroup %%------------------------------------<<

\begingroup%%----------------------------------->>
\fontsize{9pt}{11pt}\selectfont
\textbf{\small EXAMPLE} \  
Let $f:X \ra Y$ be a Dold fibration or a Serre fibration.  Assume: 
$
\begin{cases}
\ X\\
\ Y
\end{cases}
$
and the $X_y$ are path connected and 
$
\begin{cases}
\ \pi_1(X)\\
\ \pi_1(Y)
\end{cases}
$
and the $\pi_1(X_y)$ are abelian.  Fix $y_0 \in Y$ $\&$ $x_0 \in X_{y_0}$ $-$then there is an exact sequence 
$\cdots \ra$ 
$\pi_{n+1}(Y,y_0) \ra$ 
$\pi_{n}(X_{y_0},x_0) \ra$ 
$\pi_n(X,x_0) \ra$ 
$\pi_n(Y,y_0) \ra$ 
$\cdots$ and if any two of 
$\{\pi_n(X_{y_0},x_0)\}$, $\{\pi_n(X,x_0)\}$, $\{\pi_n(Y,y_0)\}$ are $P$-local, so is the third.\\
\endgroup %%------------------------------------<<

\textbf{\small LEMMA} \  
$L_P:\bAB \ra \bAB_P$ preserves finite limits.\\

\begingroup%%----------------------------------->>
\fontsize{9pt}{11pt}\selectfont
\textbf{\small EXAMPLE} \  
$L_P$ need not preserve arbitrary limits.  For instance, take $P = \bPi - \{2\}$ and define \bG in $\bTOW(\bAB)$ by 
$G_n = \Z$ $\forall \ n$ and 
$
\begin{cases}
\ G_{n+1} \ra G_n\\
\ 1 \ra 2
\end{cases}
$
$-$then $\lim \bG = 0$ but 
$\lim \bG_P = \ds\Z\left[\frac{1}{2}\right]$.\\
\endgroup %%------------------------------------<<

\label{8.20}
Let $f:G \ra K$ be a homomorphism of abelian groups $-$then $f$ is 
\un{$P$-localizing}
\index{P-localizing (morphism of abelian groups)} 
if $\exists$ an isomorphism $\phi:G_P \ra K$ such that 
$f = \phi \circx l_P$ 
(cf. p. \pageref{8.3}).\\

\textbf{\small LEMMA} \  
Let $f:G \ra K$ be a homomorphism of abelian groups $-$then $f$ is $P$-localizing iff 
$f$ is $P$-bijective and \mK is $P$-local.\\

Example: Let 
$
\begin{cases}
\ X\\
\ Y
\end{cases}
$
be path connected topological spaces, $f:X \ra Y$ a continuous function $-$then by the universal coefficient theorem, 
$f_*:H_n(X) \ra H_n(Y)$ is $P$-localizing 
%%----------------------------------------------------------------------------------------------06
$\forall \ n \geq 1$ iff 
$f_*:H_n(X; \Z_P) \ra H_n(Y; \Z_P)$ is an isomorphism $\forall \ n \geq 1$ and $H_n(Y)$ is $P$-local $\forall \ n \geq 1$.

\label{9.5}
Example: Let \mX be a path conncected topological space $-$then by the universal coefficient theorem,  
$H_n(X)$ is $P$-local $\forall \ n \geq 1$ iff $\forall \ p \in \ov{P}$, $H_n(X;\Z/p\Z) = 0$ $\forall \ n \geq 1$.\\

\label{9.21}
\begingroup%%----------------------------------->>
\fontsize{9pt}{11pt}\selectfont
\textbf{\small FACT} \  
Let 
\[
\begin{tikzcd}[ sep=large]
{G_1} \ar{d}{f_1} \ar{r} 
&{G_2} \ar{d}{f_2} \ar{r}
&{G_3} \ar{d}{f_3} \ar{r}
&{G_4} \ar{d}{f_4} \ar{r}
&{G_5} \ar{d}{f_5}\\
{K_1} \ar{r} 
&{K_2} \ar{r}
&{K_3} \ar{r}
&{K_4} \ar{r}
&{K_5}
\end{tikzcd}
\]
be a commutative diagram of abelian groups with exact rows.  Assume: $f_1, f_2, f_4, f_5$ are $P$-localizing $-$then $f_3$ is 
$P$-localizing.\\
\endgroup %%------------------------------------<<

\begingroup%%----------------------------------->>
\fontsize{9pt}{11pt}\selectfont
\textbf{\small EXAMPLE} \  
Let 
$
\begin{cases}
\ G\\
\ K
\end{cases}
$
be abelian groups $-$then 
$(G \otimes K)_P \approx$ 
$G_P \otimes K \approx$ 
$G \otimes K_P \approx$ 
$G_P \otimes K_P$ and 
$\Tor(G,K)_P \approx$ 
$\Tor(G_P,K) \approx$ 
$\Tor(G,K_P) \approx$ 
$\Tor(G_P,K_P)$ .\\
\endgroup %%------------------------------------<<

\begingroup%%----------------------------------->>
\fontsize{9pt}{11pt}\selectfont
\textbf{\small EXAMPLE} \  
Let 
$
\begin{cases}
\ G\\
\ K
\end{cases}
$
be abelian groups.
\\
\indent\indent (R) \ Assume: \mG is finitely generated $-$then 
$\Hom(G,K)_P \approx \Hom(G,K_P)$ and 
$\Ext(G,K)_P \approx \Ext(G,K_P)$.
\\
\indent\indent (L) \ Assume: \mK is $P$-local $-$then 
$\Hom(G_P,K) \approx \Hom(G,K)$ and 
$\Ext(G_P,K) \approx \Ext(G,K)$.
\vspi
[An injective $\Z_P$-module is also injective as an abelian group.]\\
\endgroup %%------------------------------------<<

\begingroup%%----------------------------------->>
\fontsize{9pt}{11pt}\selectfont
\textbf{\small FACT} \  
Let \mG be an abelian group $-$then $\forall \ n \geq 1$, $H_n(l_P):H_n(G) \ra H_n(G_P)$ is $P$-localizing.  In particular: 
\mG $P$-local $\implies$ $H_n(G)$ $P$-local $(\forall \ n \geq 1)$ and convsersely.

[There are three steps: 
(1) $G = \Z/p^n\Z$ or $G = \Z$ (direct verification); 
(2) \mG finitely generated (K\"unneth); 
(3) \mG arbitrary (take colimits).]

[Note: \ It is a corollary that for any abelian group \mG, 
$H_n(G;\Z_P) \approx H_n(G_P;\Z_P)$ $(n \geq 1)$.  This is also true if \mG is nilpotent (cf. Proposition 8) but is false in general.  
Example: Take $G = S_3$, $P = \{3\}$ $-$then 
$H_3(G;\Z_P) \neq 0$ $\&$ $H_3(G_P;\Z_P) = 0$.]\\
\endgroup %%------------------------------------<<

\begin{proposition} \ %01
Let $f:X \ra Y$ be either a Dold fibration or a Serre fibration such that $\forall \ p \in \ov{P}$, $f$ is $\Z/p\Z$-orientable $-$then any two of the following conditions imply the third: 
(1) $\forall \ k \geq 1$, $H_k(Y)$ is $P$-local; 
(2) $\forall \ l \geq 1$, $H_l(X_{y_0})$ is $P$-local; 
(3) $\forall \ n \geq 1$, $H_n(X)$ is $P$-local.
\end{proposition}

[In the notation of 
p. \pageref{8.4}, take $\Lambda = \Z/p\Z$.  By what has been said there, 
$\widetilde{H}_*(-,\Lambda) = 0$ for any two of \mY, $X_{y_0}$, and \mX entails $\widetilde{H}_*(-,\Lambda) = 0$ for the third.]\\

Application: Let $\pi$ be a $P$-local abelian group $-$then $\forall \ q \geq 1$, $H_q(\pi,n)$ is $P$-local.

%%----------------------------------------------------------------------------------------------07
[As recorded above, this is true when $n = 1$.  To treat the general case, proceed by induction, bearing in mind that the mapping fiber of the projection 
$\Theta K(\pi,n+1) \ra K(\pi,n+1)$ is a $K(\pi,n)$.]

[Note: \ If $\pi$ is any abelian group, then the arrow of localization $l_P:\pi \ra \pi_P$ induces a map 
$l_P:K(\pi,n) \ra K(\pi_P,n)$ and $\forall \ q \geq 1$, $H_q(l_P):H_q(\pi,n) \ra H_q(\pi_P,n)$ is $P$-localizing.  
In fact, $H_q(l_P)$ is $P$-bijective (mod $\sC_P$ Whitehead theorem) and $H_q(\pi_P,n)$ is $P$-local.]\\

\begingroup%%----------------------------------->>
\fontsize{9pt}{11pt}\selectfont
\textbf{\small FACT} \  
Let \mX be a pointed connected CW space.  Assume: \mX is simply connected $-$then $\forall \ n \geq 1$, 
$\pi_n(X)$ is $P$-local iff $\forall \ n \geq 1$, $H_n(X)$ is $P$-local.
\vspi
[Pass from homotopy to homology via the Postnikov tower of \mX and pass from homology to homotopy via the Whitehead tower of \mX.]\\
\endgroup %%------------------------------------<<

\label{9.8}
\begingroup%%----------------------------------->>
\fontsize{9pt}{11pt}\selectfont
\textbf{\small FACT} \  
Let 
$
\begin{cases}
\ X\\
\ Y
\end{cases}
$
be pointed connected CW spaces, $f:X \ra Y$ a pointed continuous function.  Assume: \mX $\&$ \mY are simply connected $-$then $\forall \ n \geq 1$, 
$f_*:\pi_n(X) \ra \pi_n(Y)$ is $P$-localizing iff $\forall \ n \geq 1$, 
$f_*:H_n(X) \ra H_n(Y)$ is $P$-localizing.
\vspi
[Taking into account the preceding fact, this follows from the mod $\sC_P$ Whitehead theorem.]\\
\endgroup %%------------------------------------<<

If \mG and \mK are $P$-local abelian groups, then $\Hom(G,K)$, $\Ext(G,K)$, $G \otimes K$, $\Tor(G,K)$ are $P$-local and 
$\Z_P$-isomorphic to their $\Z_P$ counterparts, hence can be identified with them.\\

\textbf{\small LEMMA} \  
Suppose that $P \neq \emptyset$ and let \mG be $P$-local.  
Assume: $\Hom(G,\Z_P) = 0$ $\&$ $\Ext(G,\Z_P) = 0$ $-$then 
$G = 0$.

[To begin with, 
$\Hom(\Tor(\Q,G),\Z_P) \oplus \Ext(\Q \otimes G,\Z_P) \approx$ 
$\Hom(\Q,\Ext(G,\Z_P)) \oplus \Ext(\Q, \Hom(G,\Z_P))$ $\implies$ 
$\Ext(\Q \otimes G,\Z_P) = 0$.  On the other hand, the condition 
$\Ext(G,\Z_P) = 0$ implies that \mG is torsion free, so if $G \neq 0$, 
then $\Q \otimes G$ is a nontrivial vector space over $\Q$ : 
$\Q \otimes G \approx I \cdot \Q$ $(\#(I) \geq 1)$ $\implies$ 
$\Ext(\Q \otimes G,\Z_P) \approx$ 
$\Ext(\Q,\Z_P)^I \approx \R^I$ 
(cf. p. \pageref{8.5}).  Contradiction.]\\

\begin{proposition} \ %02
Let 
$
\begin{cases}
\ X\\[-.1cm]
\ Y
\end{cases}
$
be path connected topological spaces, $f:X \ra Y$ a continuous function $-$then 
$f_*:H_*(X;\Z_P) \ra H_*(Y;\Z_P)$ is an isomorphism iff 
$f^*:H^*(Y;\Z_P) \ra H^*(X;\Z_P)$ is an isomorphism.
\end{proposition}

[There is an exact sequence
\[ 
\cdots \ra 
\widetilde{H}_n(X;\Z_P) \ra 
\widetilde{H}_n(Y;\Z_P) \ra
\widetilde{H}_n(C_f;\Z_P) \ra
\widetilde{H}_{n-1}(X;\Z_P) \ra 
\widetilde{H}_{n-1}(Y;\Z_P) \ra
\cdots 
\]
in homology and there is an exact sequence
\[ 
\cdots \ra 
\widetilde{H}^{n-1}(Y;\Z_P) \ra 
\widetilde{H}^{n-1}(X;\Z_P) \ra
\widetilde{H}^n(C_f;\Z_P) \ra
\widetilde{H}^{n}(Y;\Z_P) \ra 
\widetilde{H}^{n}(X;\Z_P) \ra
\cdots 
\]
%%----------------------------------------------------------------------------------------------08
in cohomology.  Accordingly, it need only be shown that 
$\widetilde{H}_*(C_f;\Z_P) = 0$ iff 
$\widetilde{H}^*(C_f;\Z_P) = 0$.  
Case 1: $P = \emptyset$.  Here, 
$\widetilde{H}^n(C_f;\Q) \approx \Hom(\widetilde{H}_n(C_f;\Q),\Q)$ and the assertion is obvious.  
Case 2: $P \neq \emptyset$.  Since 
$\widetilde{H}^n(C_f;\Z_P) \approx$ 
$\Hom(\widetilde{H}_n(C_f;\Z_P),\Z_P) \oplus \Ext(\widetilde{H}_{n-1}(C_f;\Z_P),\Z_P)$,  
it is clear that 
$\widetilde{H}_*(C_f;\Z_P) = 0$ $\implies$ $\widetilde{H}^*(C_f;\Z_P) = 0$, while if 
$\widetilde{H}^*(C_f;\Z_P) = 0$, then $\forall \ n$, 
$\Hom(\widetilde{H}_n(C_f;\Z_P),\Z_P) = 0$ \ $\&$ \  $\Ext(\widetilde{H}_n(C_f;\Z_P),\Z_P) = 0$, 
thus from the lemma, 
$\widetilde{H}_n(C_f;\Z_P) = 0$.]\\

\begin{proposition} \ %03
Let 
$
\begin{cases}
\ X\\[-.1cm]
\ Y
\end{cases}
$
be path connected topological spaces, $f:X \ra Y$ a continuous function $-$then 
$f_*:H_*(X;\Z_P) \ra H_*(Y;\Z_P)$ is an isomorphism iff 
$f_*:H_*(X;\Q) \ra H_*(Y;\Q)$ is an isomorphism and $\forall \ p \in P$, 
$f_*:H_*(X;\Z/p\Z) \ra H_*(Y;\Z/p\Z)$ is an isomorphism.
\end{proposition}

[Introducing again the mapping cone, it suffices to prove that 
$\widetilde{H}_*(C_f;\Z_P) = 0$ iff 
$\widetilde{H}_*(C_f;\Q) = 0$ and $\forall \ p \in P$, 
$\widetilde{H}_*(C_f;\Z/p\Z) = 0$.  
If first 
$\widetilde{H}_*(C_f;\Z_P) = 0$, then 
$\widetilde{H}_*(C_f;\Q) \approx$ 
$\Q \otimes \widetilde{H}_*(C_f) \approx$ 
$\Q \otimes (\Z_P \otimes \widetilde{H}_*(C_f)) \approx$ 
$\Q \otimes \widetilde{H}_*(C_f;\Z_P)  = 0$ and because $p \in P$ $\implies$ 
$\Z_P \otimes \Z/p\Z = \Z/p\Z$, $\forall \ n$, 
$\widetilde{H}_n(C_f;\Z/p\Z ) \approx$ 
$\widetilde{H}_n(C_f) \otimes \Z/p\Z \oplus \Tor(\widetilde{H}_{n-1}(C_f),\Z/p\Z) \approx$ 
$\widetilde{H}_n(C_f) \otimes( \Z_P \otimes \Z/p\Z) \oplus \Tor(\widetilde{H}_{n-1}(C_f),\Z_P \otimes\Z/p\Z  \approx$ 
$(\widetilde{H}_n(C_f) \otimes \Z_P) \otimes \Z/p\Z \oplus \Tor(\Z_P \otimes \widetilde{H}_{n-1}(C_f),\Z/p\Z) \approx$ 
$\widetilde{H}_n(C_f;\Z_P) \otimes \Z/p\Z \oplus \Tor(\widetilde{H}_{n-1}(C_f;\Z_P),\Z/p\Z) = 0$.  
As for the implication in the opposite direction, 
$\widetilde{H}_*(C_f;\Z_P) = 0$ iff $\widetilde{H}_*(C_f)$ is $S_P$-torsion, so 
$\widetilde{H}_*(C_f;\Q) = 0$ $\implies$ $\widetilde{H}_*(C_f)$ is torsion and $\forall \ n$, 
$\widetilde{H}_n(C_f;\Z_P) = 0$ $\implies$ 
$\Tor(\widetilde{H}_n(C_f),\Z/p\Z) = 0$ $\implies$  
$\widetilde{H}_n(C_f)(p) = 0$ $(p \in P)$, 
i.e., $\widetilde{H}_*(C_f)$ is $S_P$-torsion.]\\

\label{5.0at}
Application: 
Let 
$
\begin{cases}
\ X\\
\ Y
\end{cases}
$
be path connected topological spaces, $f:X \ra Y$ a continuous function $-$then 
$f_*:H_n(X) \ra H_n(Y)$ is $P$-localizing $\forall \ n \geq 1$ iff 
$f_*:H_n(X;\Q) \ra H_n(Y;\Q)$ is an isomorphism $\forall \ n \geq 1$ and $\forall \ p \in P$, 
$f_*:H_n(X;\Z/p\Z) \ra H_n(Y;\Z/p\Z)$ is an isomorphism $\forall \ n \geq 1$ and $\forall \ p \in \ov{P}$, 
$H_n(Y;\Z/p\Z) = 0$ $\forall \ n \geq 1$.

[Note: \ When $P = \bPi$, ``$P$-localizing'' = ``homology equivalence'' and the last condition is vacuous.]\\

\begingroup%%----------------------------------->>
\fontsize{9pt}{11pt}\selectfont
\textbf{\small FACT} \  Let 
$
\begin{cases}
\ X\\
\ Y
\end{cases}
$
be path connected topological spaces, $f:X \ra Y$ a continuous function.  Assume: $\forall \ n$, 
$
\begin{cases}
\ H_n(X)\\
\ H_n(Y)
\end{cases}
$
is finitely generated $-$then for $P \neq \emptyset$, $f_*:H_*(X;\Z_P) \ra H_*(Y;\Z_P)$ is an isomorphism iff 
$\forall \ p \in P$, 
$f_*:H_*(X;\Z/p\Z) \ra H_*(Y;\Z/p\Z)$ is an isomorphism.\\
\endgroup %%------------------------------------<<

The theory set forth below has been developed by a number of mathematicians and can be approached in a variety of ways.  What follows is an account of the bare essentials.

%%----------------------------------------------------------------------------------------------09
A group \mG is said to be 
\un{$P$-local}
\index{P-local (group)} 
if the map 
$
\begin{cases}
\ G \ra G\\
\ g \ra g^n
\end{cases}
$
is bijective $\forall \ n \in S_P$.  $\bGR_P$ is the full subcategory of \bGR whose objects are the $P$-local groups.  
On general grounds 
(cf. p. \pageref{8.6}),  $\bGR_P$ is a reflective subcategory of \bGR with reflector 
$
L_P: 
\begin{cases}
\ \bGR \ra \bGR_P\\
\ G \ra G_P
\end{cases}
$
and arrow of localization $l_P:G \ra G_P$.

[Note: \ If \mG is abelian, then the restriction of $L_P$ to \bAB ``is'' the $L_P$ introduced earlier.]

Example: Fix $P \neq \bPi$ $-$then no nontrivial free group is $P$-local.\\

\label{9.38}
\begingroup%%----------------------------------->>
\fontsize{9pt}{11pt}\selectfont
\textbf{\small EXAMPLE} \  
Let \mX be a pointed connected CW space $-$then $\pi_1(X)$ and the $\pi_q(X) \rtimes \pi_1(X)$ $(q \geq 2)$ are $P$-local iff  
$\forall \ n \in S_P$, the arrrow 
$
\begin{cases}
\ \Omega X \ra \Omega X\\
\ \sigma \ra \sigma^n
\end{cases}
$
is a pointed homotopy equivalence.
\vspi
[For $[\bS^{q-1},\Omega X]$ (no base points) is isomorphic to $\pi_q(X) \rtimes \pi_1(X)$ $(q \geq 2)$.]\\
\endgroup %%------------------------------------<<

The kernel of $l_P:G \ra G_P$ contains the set of $S_P$-torsion elements of \mG but is ordinarily much larger.  
Definition: An element $g \in G$ is said to be of 
\un{type $S_P$}
\index{type $S_P$} 
if $\exists$ $a, b \in G$ and $n \in S_P$ : 
$g = ab^{-1}$ $\&$ $a^n = b^n$.  
The subset of \mG consisting of the elements of type $S_P$ is closed under inversion and conjugation and is annihilated by $l_P$.  
Proceeding recursively, construct a sequence 
$\{1\} = \Lambda_0 \subset \Lambda_1 \subset \cdots$ of normal subgroups of \mG by letting $\Lambda_{k+1}/\Lambda_k$ be the subgroup of $G/\Lambda_k$ generated by the elements of type $S_P$.  Put 
$\Lambda_P(G) = \bigcup\limits_k \Lambda_k$: $\Lambda_P(G)$ is a normal subgroup of \mG and it is clear that if 
$f:G \ra K$ is a homomorphism of groups, then $f(\Lambda_P(G)) \subset \Lambda_P(K)$.  On the other hand, \mG $P$-local 
$\implies$ $\Lambda_P(G) = \{1\}$, so $\ker l_P \supset \Lambda_P(G)$.  
The containment can be proper since there are examples where $\Lambda_P(G)$ is trivial but $\ker l_P$ is not trivial 
(Berrick-Casacuberta\footnote[2]{\textit{SLN} \textbf{1509} (1992), 20-29.}).  
However, for certain \mG, $\ker l_P$ is always trivial, e.g. when \mG is locally free 
(cf. p. \pageref{8.7})

Observation: $\Lambda_P(G) = \{1\}$ iff $\forall \ n \in S_P$, the map
$
\begin{cases}
\ G \ra G\\
\ g \ra g^n
\end{cases}
$
is injective.\\

\begingroup%%----------------------------------->>
\fontsize{9pt}{11pt}\selectfont
\index{generically trivial groups (example)}
\textbf{\small EXAMPLE \  (\un{Generically Trivial Groups})} \  
A group \mG is said to be 
\un{generically trivial}
\index{generically trivial (group)} provided that $\forall \ p$, $G_p = 1$.  
Example: The infinite alternating group is generically trivial.  
The homomorphic image of a generically trivial group is generically trivial, 
so generically trivial groups are perfect (but not conversely as there exist perfect groups which are locally free 
(cf. p. \pageref{8.8})).  
Since a perfect nilpotent group is trivial, the only generically trivial nilpotent group is the trivial group and since a finite $p$-group is nilpotent, 
a perfect finite group is generically trivial.  
Example: Let \mA be a ring with unit $-$then $\bST(A)$ is generically trivial 
(Berrick-Miller\footnote[3]{\textit{Math. Proc. Cambridge Philos. Soc.} \textbf{111} (1992), 219-229.}), hence $\bE(A)$ is too 
($\implies$ $\bGL(\Gamma A) = \bE(\Gamma A)$ is acyclic and generically trivial 
(cf. p. \pageref{8.9} ff.)).
\endgroup %%------------------------------------<<

%%----------------------------------------------------------------------------------------------10
\begingroup%%----------------------------------->>
\fontsize{9pt}{11pt}\selectfont
[Note: \ In the same paper it is shown that if $\{G_n:n \geq 2\}$ is a sequence of abelian groups, then there exists a generically trivial group \mG such that $H_n(G) \approx G_n$ $(n \geq 2)$.]\\
\endgroup %%------------------------------------<<

\begingroup%%----------------------------------->>
\fontsize{9pt}{11pt}\selectfont
\index{separable groups (example)}
\textbf{\small EXAMPLE \ (\un{Separable Groups})} \  
A group \mG is said to be 
\un{separable}
\index{separable (group)} 
provided that the arrow $G \ra \ds\prod\limits_p G_p$ is injective.  
The class of separable groups is closed under the formation of products and subgroups, thus is the object class of  an epireflective subcategory \bGR 
(cf. p. \pageref{8.10}).  Every nilpotent group is separable as is every locally free group.\\
\endgroup %%------------------------------------<<

\begingroup%%----------------------------------->>
\fontsize{9pt}{11pt}\selectfont
\textbf{\small FACT} \  
A group \mG is generically trivial iff every homomorphism $f:G \ra K$, where \mK is separable, is trivial.\\
\endgroup %%------------------------------------<<

\begingroup%%----------------------------------->>
\fontsize{9pt}{11pt}\selectfont
\textbf{\small EXAMPLE} \  
Let \mX be a pointed connected CW space.  
Assume: \mX is acyclic and $\pi_1(X)$ is generically trivial 
$-$then for every pointed connected CW space \mY such that $\pi_1(Y)$ is separable, $C(X,x_0;Y,y_0)$ is homotopically trivial 
(cf. p. \pageref{8.11}).\\
\endgroup %%------------------------------------<<

\label{8.14}
\label{9.120} %dmc mnft
\textbf{\small LEMMA} \  
Suppose that \mG is torsion $-$then \mG is $P$-local iff \mG is $S_{\ov{P}}$-torsion.

[Necessity: Given $g \in G$, $\exists \ t$ : $g^t = e$.  
Write $t = n\bar{n}$ ($n \in S_P$, $\bar{n} \in S_{\ov{P}}$) : 
$(g^{\ov{n}})^n = e$ $\implies$ $g^{\ov{n}} = e$.  Therefore \mG is $S_{\ov{P}}$-torsion.

Sufficiency: Fix $n \in S_P$.  For each $\bar{n} \in S_{\ov{P}}$, choose $k, l$: $kn + l\bar{n} = 1$, hence 
(i) Given $g \in G$, $\exists$ $\bar{n} \in S_{\ov{P}}$: 
$g^{\ov{n}} = e$ 
$\implies$ 
$g = g^{k n + l \bar{n}} = (g^k)^n$
(ii) Given $g_1, g_2 \in G$, $\exists$ $\bar{n} \in S_{\ov{P}}$: 
$g_1^{\ov{n}} = e = g_2^{\ov{n}}$, so $g_1^n = g_2^n$ 
$\implies$
$g_1 = (g_1^n)^k (g_1^{\bar{n}})^l =$ 
$(g_2^n)^k(g_2^{\bar{n}})^l = g_2$.]\\

\textbf{\small LEMMA} \  
Suppose that \mG is torsion $-$then $l_P:G \ra G_P$ is surjective and $\ker l_P$ is generated by the 
$S_P$-torsion elements of \mG.

[Let $\Lambda$ be the subgroup of \mG generated by the $S_P$-torsion elements of \mG.  Since \mG is torsion, 
$G/\Lambda$ is $S_{\ov{P}}$-torsion, thus $P$-local.  In addition, for every homomorphism $f:G \ra K$, where \mK is $P$-local, 
$f(\Lambda) = \{1\}$.]\\

\begingroup%%----------------------------------->>
\fontsize{9pt}{11pt}\selectfont
\textbf{\small FACT} \  
Let $1 \ra G^\prime \ra G \ra G\pp \ra 1$ be a short exact sequence of groups.  
Assume: $G^\prime$ is $P$-local and $G\pp$ is 
$S_{\ov{P}}$-torsion $-$then \mG is $P$-local.\\
\endgroup %%------------------------------------<<

\begingroup%%----------------------------------->>
\fontsize{9pt}{11pt}\selectfont
\textbf{\small EXAMPLE} \  
Let \mX be a pointed connected CW space.  
Assume: $\pi_1(X)$ is $S_{\ov{P}}$-torsion and $\forall \ q \geq 2$, $\pi_q(X)$ is $P$-local 
$-$then $\forall \ n \ \in S_P$, the arrow 
$
\begin{cases}
\ \Omega X \ra \Omega X\\
\ \sigma \ra \sigma^n
\end{cases}
$
is a pointed homotopy equivalence.\\
\endgroup %%------------------------------------<<

\label{8.30}
\begingroup%%----------------------------------->>
\fontsize{9pt}{11pt}\selectfont
\textbf{\small FACT} \  
Let $1 \ra G^\prime \ra G \ra G\pp \ra 1$ be a short exact sequence of groups.  
Assume: $G\pp$ is $S_{\ov{P}}$-torsion $-$then the sequence 
$1 \ra G_P^\prime \ra G_P \ra G_P\pp \ra 1$ is exact.\\
\endgroup %%------------------------------------<<

%%----------------------------------------------------------------------------------------------11
\begingroup%%----------------------------------->>
\fontsize{9pt}{11pt}\selectfont
\textbf{\small EXAMPLE} \  
Let $\pi$ be the fundamental group of the Klein bottle $-$then there is a short exact sequence 
$1 \ra$ 
$\Z \oplus \Z \ra$ 
$\pi \ra$ 
$\Z/2\Z \ra 1$ so if $2 \in P$, there is a short exact sequence 
$1 \ra$ 
$\Z_P \oplus \Z_P \ra$ 
$\pi_P \ra$ 
$\Z/2\Z \ra 1$ and $l_P:\pi \ra \pi_P$ is injective (but this is false if $2 \notin P$).\\
\endgroup %%------------------------------------<<

\label{9.76}
\label{10.9}
\begingroup%%----------------------------------->>
\fontsize{9pt}{11pt}\selectfont
\index{finite groups (example)}
\textbf{\small EXAMPLE  \ (\un{Finite Groups})} \ 
Let \mG be a finite group $-$then $l_P:G \ra G_P$ is surjective and $\ker l_P$ is the subgroup of \mG generated by the Sylow $p$-subgroups 
$(p \in \ov{P})$, so e.g. if \mG is a $p$-group, 
$
G_P = 
\begin{cases}
\ G \quad (p \in P)\\
\ 1 \quad \  (p \in \ov{P}
\end{cases}
. \ 
$
Therefore \mG is $P$-local iff $\#(G) \in S_{\ov{P}}$.\\
\endgroup %%------------------------------------<<

\begingroup%%----------------------------------->>
\fontsize{9pt}{11pt}\selectfont
\textbf{\small FACT} \  
Let \mG be a finite group $-$then \mG is $P$-local iff $\forall \ n \geq 1$, $H_n(G)$ is $P$-local.

[Given a nontrivial subgroup $K \subset G$, the homomorphism $H_n(K) \ra H_n(G)$ is nonzero for infinitely many n (Swan\footnote[2]{\textit{Proc. Amer. Math. Soc.} \textbf{11} (1960), 885-887.}).  
Since 
$H_n(G) \approx \ds\bigoplus\limits_{p|\#(G)} H_n(G)(p)$, it follows that $\forall \ p|\#(G)$, $H_n(G)(p) \neq 0$ for infinitely many $n$.]\\
\endgroup %%------------------------------------<<

\label{9.31} %dmc mnft

\begingroup%%----------------------------------->>
\fontsize{9pt}{11pt}\selectfont
\textbf{\small FACT} \  
Let \mG be a finite group $-$then  
$H_1(l_P):H_1(G;\Z_P) \ra H_1(G_P;\Z_P)$ is bijective and 
$H_2(l_P):H_2(G;\Z_P) \ra H_2(G_P;\Z_P)$ is surjective.

[The short exact sequence 
$1 \ra \ker l_P \ra G \ra G_P \ra 1$ leads to an exact sequence 
$H_2(G;\Z_P) \ra$ 
$H_2(G_P;\Z_P) \ra$ 
$\Z_P \otimes \ker l_P/[G,\ker l_P] \ra$ 
$H_1(G;\Z_P) \ra$ 
$H_1(G_P;\Z_P) \ra 0$ in which the middle term is zero.]\\
\endgroup %%------------------------------------<<

\begingroup%%----------------------------------->>
\fontsize{9pt}{11pt}\selectfont
\textbf{\small FACT} \  
Let \mG be a finite group $-$then $\forall \ n \geq 1$, $H_n(l_P):H_n(G) \ra H_n(G_P)$ is $P$-localizing iff $\ker l_P$ is $S_P$-torsion.\\
\endgroup %%------------------------------------<<

\begingroup%%----------------------------------->>
\fontsize{9pt}{11pt}\selectfont
Application: Let \mG be finite group.  
Suppose that $\forall \ p$ $\&$ $\forall \ n \geq 1$,  $H_n(l_p):H_n(G) \ra H_n(G_p)$ 
is $P$-localizing $-$then \mG is nilpotent.
\vspi
[The Sylow subgroups of \mG are normal.]\\
\endgroup %%------------------------------------<<

A subgroup \mK of a group \mG is said to be 
\un{$P$-isolated}
\index{P-isolated} 
if $\forall \ g \in G$, 
$\forall \ n \in S_P$ : $g^n \in K$ $\implies$ $g \in K$.  
The intersection of a collection of $P$-isolated subgroups of \mG is $P$-isolated.  
Therefore every nonempty subset $X \subset G$ is contained in a unique minimal $P$-isolated subgroup of \mG, the
 \un{$P$-isolator}
\index{P-isolator} 
of \mX, writtten $I_P(G,X)$.  To describe $I_P(G,X)$, let $X_1 = X$, 
$I_1 = \langle X_1 \rangle$, and define $X_{i+1}, I_{i+1}$ inductively by setting 
$X_{i+1} =$ $\{g:g^n \in I_i \ (\exists \  n \in S_P)\}$, $I_{i+1} = \langle X_{i+1} \rangle$ $-$then 
$I_P(G,X) = \bigcup\limits_i I_i$.  
Corollary: \mX conjugation invariant $\implies$ $I_P(G,X)$ normal.

[Note: \ A $P$-isolated subgroup of a $P$-local group is $P$-local.]

%%----------------------------------------------------------------------------------------------12
\label{8.13}
\label{8.15}
Example: For any \mG, $G_P = I_P(G_P,l_P(G))$.

[Note: \ More generally, if $f:G \ra K$ is a homomorphism of groups, then 
$f_P(G_P) =$ 
$I_P(K_P,l_P(f(G)))$.  
Corollary: $f$ surjective $\implies$ $f_P$ surjective.]\\

\begingroup%%----------------------------------->>
\fontsize{9pt}{11pt}\selectfont
\textbf{\small EXAMPLE} \  
Fix a prime $p$ $-$then $\Z/p^\infty \Z$ is isomorphic to $I_{\ov{p}}(\Q,\Z)/\Z$.\\
\endgroup %%------------------------------------<<

\label{9.7}
\begingroup%%----------------------------------->>
\fontsize{9pt}{11pt}\selectfont
\textbf{\small EXAMPLE} \  
Let \mF be a free group on $n > 1$ generators $-$then $F/[F,F] \approx n\cdot \Z$.  
By contrast, 
Baumslag\footnote[2]{\textit{Acta Math.} \textbf{104} (1960), 217-303 (cf. 253-254 $\&$ 291-293).}
has shown that 
$F_P/I_P(F_P,[F_P,F_P]) \approx n \cdot \Z_P$, while 
$F_P/[F_P,F_P] \approx n \cdot \Z_P \oplus \ds\bigoplus\limits_{p \in \ov{P}} \omega \cdot (\Z/p^\infty\Z)$.]
\vspi
[Note: \ Since 
$\ds\bigoplus\limits_{p \in \ov{P}} \omega \cdot (\Z/p^\infty\Z)$ is $S_P$-torsion, $H_1(F_P)$ is not $P$-local if 
$P \neq \bPi$.  
This example also shows that in \bGR, the operations 
$G \ra$ 
$G/[G,G] \ra$ 
$(G/[G,G])_P$, 
$G \ra$ 
$G_P \ra$ 
$G_P/[G_P,G_P]$ need not coincide.]\\
\endgroup %%------------------------------------<<

\begingroup%%----------------------------------->>
\fontsize{9pt}{11pt}\selectfont
\textbf{\small FACT} \  
If \mG is a nilpotent group and if \mK is a subgroup of \mG, then $\{g:g^n \in K \  (\exists \ n \in S_P)\}$ 
is a subgroup of \mG, hence equals 
$I_P(G,K)$.
\vspi
[Assuming that \mG is generated by the $g$ such that for some $n \in S_P$, $g^n \in K$, one can argue inductively on 
$d = \nil G > 1$ and suppose that $\forall \ g \in G$, $\exists \ n \in S_P$ $\&$ $h \in \Gamma^{d-1}(G)$, $k \in K$ : $g^n = h k$.  
On the other hand, 
$\otimes^d ([G,G] \cdot K/[G,G]) \ra$ 
$\otimes^d (G/[G,G]) \ra$ 
$\Gamma^{d-1}(G)$, so $\exists \ m \in S_P$: $h^m \in K$.  
But $h$ is central, thus $g^{nm} = h^m k^m \in K$.]
\vspi
[Note: \ In particular, the set of $S_P$-torsion elements in a nilpotent group is a subgroup 
(cf. p. \pageref{8.12}).]\\
\endgroup %%------------------------------------<<

\index{Commutator Formula}
\textbf{\small COMMUTATOR FORMULA} \  
Suppose that $\Lambda_P(G) = \{1\}$.  Let $
\begin{cases}
\ K\\
\ L
\end{cases}
$
be subgroups of \mG $-$then $[K,L] = \{1\}$ $\implies$ $[I_P(G,K),I_P(G,L)] = \{1\}$.

[Given 
$
\begin{cases}
\ x \in I_P(G,K)\\
\ y \in I_P(G,L)
\end{cases}
\!\!\!\!\!, \ 
$
the claim is that $xyx^{-1} = y$.  
This is trivial if 
$
\begin{cases}
\ x \in I_1(G,K)\\
\ y \in I_1(G,L)
\end{cases}
\!\!\!\!\!\! ,\ 
$
so argue by induction on $i$, assuming that 
$
\begin{cases}
\ x \in I_{i+1}(G,K)\\
\ y \in I_{i+1}(G,L)
\end{cases}
$
with 
$
\begin{cases}
\ x^n \in I_i(G,K)\\
\ y^n \in I_i(G,L)
\end{cases}
\!\!\!\!\! (\exists \ n \in S_P)
$
$-$then $(y^{-n}xy^n)^n =$ 
$y^{-n}x^ny^n = x^n$ $\implies$ 
$y^{-n}xy^n = x$ $\implies$
$xy^nx^{-1} = y^n$  $\implies$ 
$(xyx^{-1})^n = y^n$  $\implies$ 
$xyx^{-1} = y$.]\\

\quad Application: Suppose that $\Lambda_P(G) = \{1\}$.  Let $g_1, g_2$, be elements of \mG such that 
$[g_1^{n_1},g_2^{n_2}] = 1$, where $n_1, n_2 \in S_P$ $-$then $[g_1,g_2] = 1$.\\

\textbf{\small LEMMA} \  
Suppose that $\Lambda_P(G) = \{1\}$.  Let \mK be a $P$-isolated central subgroup of \mG $-$then $\Lambda_P(G/K) = \{1\}$.

%%----------------------------------------------------------------------------------------------13
[Consider an element of type $S_P$ in $G/K$, say 
$gK = (aK)(b^{-1}K)$ $\&$ $a^nK = b^nK$ $(\exists \ n \in S_P)$.  So: 
$a^n = b^n k$ $(\exists \ k \in K)$ $\implies$ $[a^n,b^n] = 1$ $\implies$ $[a,b] = 1$ $\implies$ $(ab^{-1})^n \in K$ 
$\implies$ $ab^{-1} \in I_P(G,K) = K$ $\implies$ $aK = bK$.]\\

\index{Transmission of Nilpotency}
\textbf{\small TRANSMISSION OF NILPOTENCY} \quad Suppose that $\Lambda_P(G) = \{1\}$.  Let \mK be a nilpotent subgroup of \mG $-$then $I_P(G,K)$ is nilpotent with $\nil I_P(G,K) = \nil K$.

[The assertion is obvious if \mK consists of the identity alone.  \ 
Assume next that \mK is abelian and nontrivial: 
$[K,K] = \{1\}$ $\implies$ $[I_P(G,K),I_P(G,K)] = \{1\}$ $\implies$ $I_P(G,K)$ is abelian and nontrivial.  
Induction hypothesis: The assertion is true whenever \mL is a nilpotent subgoup of \mH provided that $\Lambda_P(H) = \{1\}$ and 
$\nil L \leq d - 1$, where $d = \nil K > 1$.  
Let \mZ be the center of \mK $-$then $[K,Z] = \{1\}$ $\implies$ 
$[I_P(G,K),I_P(G,Z)] = \{1\}$, thus $I_P(G,Z)$ is a $P$-isolated central subgroup of $I_P(G,K)$, so by the lemma, 
$\Lambda_P(I_P(G,K)/I_P(G,Z)) = \{1\}$.  
Now put $X = K \cdot I_P(G,Z)$: $I_1 = X_1$ is a group $(X_1 = X)$ and 
$I_1/I_P(G,Z) \approx$ 
$K/K \cap I_P(G,Z) \approx$ 
$K/Z$.  
Since $\nil K/Z = \nil K - 1$, it follows that $I_1$ is nilpotent with $\nil I_1 = d$.  
Write, as above, 
$I_P(G,X) = \bigcup\limits_i I_i$.   
Assume that $I_i$ is nilpotent with  $\nil I_i = d$ \ $\forall \ i \leq i_0$.  
Fix a well ordering of the elements $X_{i_0 + 1}$: $\{x_\beta: 0 \leq \beta < \alpha\}$.  
Let $W_\gamma$ be the subgroup of \mG generated by $I_{i_0}$ and $\{x_\beta: 0 \leq \beta < \gamma\}$ $-$then 
$I_{i_0 + 1} = \bigcup\limits_\gamma W_\gamma$ and the claim is that $\forall \ \gamma$, $W_\gamma$ is nilpotent with 
$\nil W_\gamma = d$, hence that $I_{i_0 + 1}$ is nilpotent with $\nil I_{i_0 + 1} = d$.  
Consider $W_1$: $I_{i_0}/I_P(G,Z)$ is a nilpotent subgroup of $W_1/I_P(G,Z)$ with $\nil I_{i_0}/I_P(G,Z) = d - 1$.  
Therefore the induction hypothesis implies that 
$W_1/I_P(G,Z) = I_P(W_1/I_P(G,Z),I_{i_0}/I_P(G,Z))$ is nilpotent with $\nil W_1/I_P(G,Z) = d - 1$.  
This means that $W_1$ is nilpotent with $\nil W_1 = d$, which sets the stage for an evident transfinite recursion.  
Conclusion: $\forall \ i$, $I_i$ is nilpotent with $\nil I_i = d$, i.e., $I_P(G,X)$ is nilpotent with 
$\nil I_P(G,X) = d$ or still, $I_P(G,K)$ is nilpotent with 
$\nil I_P(G,K) = d$.]\\

\begin{proposition} \ %04
Let \mG be a nilpotent group $-$then $G_P$ is nilpotent and $\nil G_P \leq \nil G$.
\end{proposition}

[In fact, $G_P = I_P(G_P,l_P(G)))$ and transmission of nilpotency ensures that $G_P$ is nilpotent with 
$\nil G_P = \nil l_P(G) \leq \nil G$.]\\

Notation: \bNIL is the category of nilpotent groups and $\bNIL^d$ is the category of nilpotent groups with degree of nilpotency $\leq d$.

Thanks to Proposition 4, $L_P$ respects \bNIL : \mG nilpotent $\implies$ $G_P$ nilpotent, thus $\bNIL_P$, the full subcategory of \bNIL whose objects are the $P$-local nilpotent groups, is a reflective subcategory of \bNIL.  More is true: 
$L_P$ respects $\bNIL^d$ and there is a commutative diagram
%%----------------------------------------------------------------------------------------------14
\begin{tikzcd}%[ sep=small]
{\bNIL^{d+1}} \ar{r}{L_P} &{\bNIL_P^{d+1}}\\
{\bNIL^d} \ar{u} \ar{r}[swap]{L_P} &{\bNIL_P^d} \ar{u}
\end{tikzcd}
(obvious notation).\\[.25cm]

\begingroup%%----------------------------------->>
\fontsize{9pt}{11pt}\selectfont
\textbf{\small FACT} \  
Let \mG be a group.  Assume: \mG is locally nilpotent $-$then $G_P$ is locally nilpotent.
\vspi
[Note: \ A group is said to be 
\un{locally nilpotent}
\index{locally nilpotent (group)} 
if its finitely generated subgroups are nilpotent.]\\
\endgroup %%------------------------------------<<

\begingroup%%----------------------------------->>\
\fontsize{9pt}{11pt}\selectfont
\textbf{\small FACT} \  
Let \mG be a group.  Assume: \mG is virtually nilpotent $-$then $G_P$ is virtually nilpotent.
\vspi
[Note: \ A group is said to be 
\un{virtually nilpotent}
\index{virtually nilpotent (group)} 
if it contains a nilpotent subgroup of finite index.]\\
\endgroup %%------------------------------------<<

\label{10.8}
Given a set of primes \mP, a group \mG is said to be 
\un{residually finite \mP}
\index{residually finite \mP (group)} 
if $\forall \ g \neq e$ in \mG there is a finite $S_{\ov{P}}$-torsion group $X_g$ and an epimorphism $\phi_g:G \ra X_g$ such that $\phi_g(g) \neq e$.

[Note: \ When $P = \bPi$, the term is
\index{residually finite (group)}  
\un{residually finite}. 
Example: $\Q$ is not residually finite but $\Z_P$ $(P \neq \emptyset)$ is residually finite $p \ \forall \ p \in P$.]

Examples: 
(1) \  (Iwasawa) Every free group is residually finite $p$ for all primes $p$; 
(2) \  (Hirsch) Every polycyclic group is residually finite ($\implies$ every finitely generated nilpotent group is residually finite); 
(3) \  (Gruenberg) Every finitely generated torsion free nilpotent group is residually finite $p$ for all primes $p$; 
(4) \  (Hall) Every finitely generated abelian-by-nilpotent group is residually finite.

[Note: Proofs of these results can be found in 
Robinson\footnote[2]{\textit{Finiteness Conditions and Generalized Soluble Groups}, vol II, Springer Verlag (1972); 
see also Magnus, \textit{Bull. Amer. Math. Soc.} \textbf{75} (1969), 305-316.}.]\\

\textbf{\small LEMMA} \  Let \mG be a finitely generated nilpotent group.  
Assume: All the torsion in \mG is 
$S_{\ov{P}}$-torsion, where $P \neq \emptyset$ $-$then \mG is residually finite \mP.

[Fix $g \neq e$ in \mG.  
Case 1: $g \notin G_\tor$.  
According to Gruenberg, $\forall \ p$, $G/G_\tor$ is residually finite $p$, so a fortiori $G/G_\tor$ is residually finite \mP.  
Case 2: $g \in G_\tor$.  
According to Hirsch, there is a finite nilpotent group $X_g$ and an epimorphism $\phi_g:G \ra X_g$ such that $\phi_g(g) \neq e$.  Write $X_g = \prod\limits_p X_g(p)$, $X_g(p)$ the Sylow $p$-subgroup of $X_g$.  Let $\pi_P$ be the projection $X_g \ra \prod\limits_{p \in P} X_g(p)$ and consider the composite $\pi_P \circx \phi_g$.]\\

\begin{proposition} \ %05
Let \mG be a nilpotent group $-$then $l_P:G \ra G_P$ is $P$-bijective.
\end{proposition}

[Since $G_P$ is nilpotent, $\{g_P:g_P^n \in l_P(G) (\ \exists \ n \in S_P)\}$ equals $I_P(G_P,l_P(G))$ 
(cf. p. \pageref{8.13}) or still, $G_P$, thus $l_P$ is $P$-surjective.  
To verify that $l_P$ is $P$-injective, suppose first that $P$ is nonempty.  
Because the kernel of $l_P$ contains the $S_P$-torsion, one can assume that
%%----------------------------------------------------------------------------------------------15
all the torsion in \mG is $S_{\ov{P}}$-torsion.  The claim in this situation is that $l_P$ is injective.  
If to begin with \mG is finitely generated, then on the basis of the lemma, there is an embedding $G \ra \prod\limits_{g \neq e} X_g$, where each $X_g$ is a finite $S_{\ov{P}}$-torsion group, hence $P$-local 
(cf. p. \pageref{8.14}).  Therefore $\prod\limits_{g \neq e} X_g$ is 
$P$-local, so $l_P$ is necessarily injective.  
To see that $l_P$ is injective in general, express \mG as the colimit of its finitely generated subgroups $G_i$ and compute the kernel of $G \ra G_P$ as the colimit of the kernels of the $G_i \ra G_{i,P}$.  
There remains the possibility that $P$ is empty.  
To finesse this, choose $P: P \neq \emptyset$ $\&$ $\ov{P} \neq \emptyset$ and note that the arrow 
$(G_P)_{\ov{P}} \ra G_{\emptyset}$ $(=G_\Q)$ is an isomorphism which implies that $G_{\tor} = \ker l_{\Q}$.]\\

Application: Every torsion free nilpotent group embeds in its rationalization.\\

\textbf{\small LEMMA} \  
Let $f:G \ra K$ be a homomorphism of nilpotent groups.  Assume: $f$ is injective (surjective) $-$then $f_P:G_P \ra K_P$ is injective (surjective).

[It will be enough to establish injectivity 
(see p. \pageref{8.15} for surjectivity).  Suppose that 
$f_P(g_P) = e$ $(g_P \in G_P)$.  
Since $l_P$ is $P$-surjective, $\exists \ g \in G$ $\&$ $n \in S_P$: $l_P(g) = g_P^n$ 
$\implies$ $l_P(f(g)) = e$ $\implies$ $\exists$ $m \in S_P$ : $f(g)^m = e$ $\implies$ $g^m = e$ $\implies$ 
$g \in \ker l_P$ $\implies$ $g_P^n = e$ $\implies$ $g_P = e$, $G_P$ being $P$-local.]\\

\begin{proposition} \ %06
$L_P:\bNIL \ra \bNIL_P$ is exact, i.e., if 
$1 \ra G^\prime \ra G \ra G\pp \ra 1$ is a short exact sequence in \bNIL, then 
$1 \ra G_P^\prime \ra G_P \ra G_P\pp \ra 1$ is a short exact sequence in $\bNIL_P$.
\end{proposition}

[It is straightforward to check that $\im(G_P^\prime \ra G_P) = \ker(G_P \ra G_P\pp)$.]\\

\textbf{\small LEMMA} \  
Let \mG be a nilpotent group.  Suppose that \mK is a central subgroup of \mG $-$then $K_P$ is a central subgroup of $G_P$.

[In fact, $[G,K] = \{1\}$ $\implies$ $[l_P(G),l_P(K)] = \{1\}$, so by the commutator formula, 
$[I_P(G_P,l_P(G)),I_P(G_P,l_P(K))] = \{1 \}$ $\implies$ $[G_P,K_P] = \{1\}$.]\\

\begin{proposition} \ %07
$L_P:\bNIL \ra \bNIL_P$ preserves central extensions.\\
\end{proposition}

\textbf{\small LEMMA} \  
Let $G_1 \ra G_2 \ra G_3 \ra G_4 \ra G_5$ be an exact sequence of nilpotent groups.  Assume: 
$
\begin{cases}
\ G_1, G_2\\
\ G_4, G_5
\end{cases}
$
are $P$-local $-$then $G_3$ is $P$-local.\\[.25cm]

Application: Let $1 \ra G^\prime \ra G \ra G\pp \ra 1$ be a short exact sequence of nilpotent groups.  Assume: Two of the groups are $P$-local $-$then so is the third.\\

%%----------------------------------------------------------------------------------------------16
\label{9.43} %dmc mnft
\begingroup%%----------------------------------->>
\fontsize{9pt}{11pt}\selectfont
\textbf{\small EXAMPLE} \  
Let \mX be a pointed connected CW space.  Assume: \mX is nilpotent and $\forall \ q \geq 1$, $\pi_q(X)$ is $P$-local $-$then 
$\forall \ n \in S_P$, the arrrow 
$
\begin{cases}
\ \Omega X \ra \Omega X\\
\ \sigma \ra \sigma^n
\end{cases}
$
is a pointed homotopy equivalence.
\vspi
[There is a split short exact sequence 
$1 \ra$ 
$\pi_q(X) \ra$ 
$\pi_q(X) \rtimes \pi_1(X) \ra$ 
$\pi_1(X) \ra 1$, where $\pi_q(X) \rtimes \pi_1(X)$ $(q \geq 2)$ is nilpotent 
(cf. p. \pageref{8.16}), hence $P$-local.]
\\[-.2cm]
\endgroup %%------------------------------------<<

\begingroup%%----------------------------------->>
\fontsize{9pt}{11pt}\selectfont
If $f, g:G \ra K$ are homomorphisms of nilpotent groups such that $\forall \ p$, $f_p = g_p$, then $f = g$.  
In other words, morphisms in \bNIL (as in \bAB) are determined by their localizations.  
For finitely generated objects the situation is different.  
Definition: Given a finitely generated nilpotent group \mG, the 
\un{genus}
\index{genus (of a finitely generated nilpotent group)}
$\gen G$  of $G$ is the conglomerate of the isomorphism classes of finitely generated nilpotent groups \mK such that $\forall \ p$, 
$G_p \approx K_p$.  
By contrast to what obtains in \bAB, it can happen that $\#(\gen G) > 1$ although always $\#(\gen G) < \omega$ 
(Pickel\footnote[2]{\textit{Trans. Amer. Math. Soc.} \textbf{160} (1971), 327-341;
 see also Mislin, \textit{SLN} \textbf{418} (1974), 103-120 
 and Warfield, \textit{J. Pure Appl. Algebra} \textbf{6} (1975), 125-132.}).
\vspi
[Note: \ If \mG is a finitely generated abelian group and if \mK is a a finitely generated nilpotent group such that 
$\forall \ p$, $G_p \approx K_p$, then $G \approx K$ ($\implies$ $\gen \ G = \{[G]\}$).\\
\endgroup %%------------------------------------<<

\begingroup%%----------------------------------->>
\fontsize{9pt}{11pt}\selectfont
\textbf{\small FACT} \  
Let \mG be a nilpotent group $-$then two elements of \mG are conjugate iff their images in every $G_p$ are conjugate.\\ 
\endgroup %%------------------------------------<<

Let \mG be a nilpotent group $-$then one may attach to \mG a sink 
$\{r_p:G_p \ra G_\Q\}$ and a source 
$\{l_p:G \ra G_p\}$, where 
$
\forall
\begin{cases}
\\[-.95cm]
\ p\\[-.1cm]
\ q\\[-.1cm]
\end{cases}
, \ r_p \circx l_p = r_q \circx l_q.
$
\\

\textbf{\small LEMMA} \  
Let $1 \ra G^\prime \ra G \ra G\pp \ra 1$ be a short exact sequence of nilpotent groups.  
Assume: The source 
$
\begin{cases}
\ \{l_p:G^\prime \ra G_p^\prime\}\\
\ \{l_p:G\pp \ra G_p\pp \}
\end{cases}
$
is the multiple pullback of the sink 
$
\begin{cases}
\ \{r_p:G_p^\prime \ra G_\Q^\prime \}\\
\ \{r_p:G_{p}\pp \ra G_\Q\pp \}
\end{cases}
$
$-$then the source $\{l_p:G \ra G_p\}$ is the multiple pullback of the sink 
\\[.2cm]$\{r_p:G_p \ra G_\Q\}$.

[The verification is a diagram chase, using the exactness of 
$1 \ra G_p^\prime \ra G_p \ra G_p\pp \ra 1$.  
Precisely: Given elements $g_p \in G_p$ \ $\&$ \  $g_\Q \in G_\Q$ : $\forall \ p$, 
$r_p(g_p) = g_\Q$, $\exists$! $g \in G$ : $\forall \ p$, $l_p(g) = g_p$.]\\

\label{9.37}
\index{Fracture Lemma}
\textbf{\small FRACTURE LEMMA} \  
Suppose that \mG is a finitely generated nilpotent group $-$then the source $\{l_p:G \ra G_p\}$ is the multiple pullback of the sink $\{r_p:G_p \ra G_\Q\}$.

[Proceed by induction on $\nil G$.  The assertion is true if $\nil G \leq 1$ 
(cf. p. \pageref{8.17}).  
Assume therefore that $\nil G > 1$ and consider the short exact sequence
$1 \ra$ 
$\Gamma^1(G) \ra$ 
$G \ra$ 
%%----------------------------------------------------------------------------------------------17
$G/\Gamma^1(G) \ra 1$ 
of nilpotent groups.  Since \mG is finitely generated, $\Gamma^1(G)$ is finitely generated 
(cf. p. \pageref{8.18}), as is 
$G/\Gamma^1(G)$.  
Furthermore, 
$\nil \Gamma^1(G) < \nil G$ and 
$\nil G/\Gamma^1(G) = 1$, thus the lemma is applicable.]\\

Let $f:G \ra K$ be a homomorphism of nilpotent groups $-$then $f$ is said to be 
\un{$P$-localizing}
\index{P-localizing} 
if $\exists$ an isomorphism $\phi:G_P \ra K$ such that $f  = \phi \circx l_P$ 
(cf. p. \pageref{8.19}).\\

\textbf{\small LEMMA} \  
Let $f:G \ra K$ be a homomorphism of nilpotent groups 
$-$then $f$ is $P$-localizing iff $f$ is $P$-bijective and \mK is $P$-local.

[Note: \ A homomorphism $f:G \ra K$ of nilpotent groups is $P$-bijective iff $f_P:G_P \ra K_P$ is bijective 
(cf. Proposition 5).]\\

\begingroup%%----------------------------------->>
\fontsize{9pt}{11pt}\selectfont
\textbf{\small FACT} \  
Let 
\[
\begin{tikzcd}[sep=large]
{G_1} \ar{d}{f_1} \ar{r} 
&{G_2} \ar{d}{f_2} \ar{r}
&{G_3} \ar{d}{f_3} \ar{r}
&{G_4} \ar{d}{f_4} \ar{r}
&{G_5} \ar{d}{f_5}\\
{K_1} \ar{r} 
&{K_2} \ar{r}
&{K_3} \ar{r}
&{K_4} \ar{r}
&{K_5}
\end{tikzcd}
\]
be a commutative diagram of nilpotent groups with exact rows.  Assume: $f_1, f_2, f_4, f_5$ are $P$-localizing 
$-$then $f_3$ is $P$-localizing.\\
\endgroup %%------------------------------------<<

\begin{proposition} \ %08
Let \mG be a nilpotent group $-$then $\forall \ n \geq 1$, $H_n(l_P):H_n(G) \ra H_n(G_P)$ is $P$-localizing.
\end{proposition}

[This is true if $\nil G \leq 1$, so argue by induction on $\nil G > 1$.  
There is a commutative diagram
\[
\begin{tikzcd}%[sep=large]
{1}\ar{r} 
&{\Cen G} \ar{d} \ar{r}
&{G} \ar{d} \ar{r}
&{G/\Cen G} \ar{d} \ar{r}
&{1}\\
{1} \ar{r} 
&{(\Cen G)_P} \ar{r}
&{G_P} \ar{r}
&{(G/\Cen G)_P} \ar{r}
&{1}
\end{tikzcd}
\]
of central extesnsions (cf. Proposition 7) and a morphism 
$\{E_{p,q}^2 \approx$ 
$H_p(G/\Cen G; H_q(\Cen G))\}$ $\ra$ 
$\{\ov{E}_{p,q}^2 \approx H_p((G/\Cen G)_P;H_q((\Cen G)_P))\}$ of LHS spectral sequences.  Since 
$\nil \Cen G \leq 1$ and 
$\nil G/\Cen G \leq \nil G - 1$, it follows from the induction hypothesis and the universal coefficient theorem that the arrow 
$E_{p,q}^2 \ra \ov{E}_{p,q}^2$ is $P$-localizing $(p + q > 0)$.  
However, the homology groups attached to a chain complex of $P$-local abelian groups are $P$-local 
(cf. p. \pageref{8.20}), 
thus this conclusion persists through the spectral sequence and in the end, it is seen that the arrow 
$E_{p,q}^\infty \ra \ov{E}_{p,q}^\infty$ is $P$-localizing $(p + q > 0)$.  Fix now an $n \geq 1$.  Consider the commutative diagram 
\[
\begin{tikzcd}%[sep=large]
{0}  \ar{r} 
&{H_{p-1,q+1}} \ar{d} \ar{r}
&{H_{p,q}} \ar{d} \ar{r}
&{E_{p,q}^\infty} \ar{d} \ar{r}
&{0} \\
{0} \ar{r} 
&{\ov{H}_{p-1,q+1}} \ar{r}
&{\ov{H}_{p,q}} \ar{r}
&{\ov{E}_{p,q}^\infty} \ar{r}
&{0}
\end{tikzcd}
,
\]
%%----------------------------------------------------------------------------------------------18
where $p + q = n$ $-$then the obvious recursion argument allows one to say that the arrow 
$H_{p,q} \ra \ov{H}_{p,q}$ is $P$-localizing, therefore $H_n(l_P): H_n(G) \ra H_n(G_P)$ is $P$-localizing.]\\

Application: Let \mG be a nilpotent group $-$then $\forall \ n \geq 1$, $H_n(G)_P \approx H_n(G_P)$.\\

\begingroup%%----------------------------------->>
\fontsize{9pt}{11pt}\selectfont
\textbf{\small FACT} \  
Suppose that \mG and \mK are finitely generated nilpotent groups.  Assume: $\gen G = \gen K$ $-$then $\forall \ n \geq 1$, 
$H_n(G) \approx H_n(K)$.

[The point here is that $H_n(G)$ and $H_n(K)$ are finitely generated 
(cf. p. \pageref{8.21}).]\\

\endgroup %%------------------------------------<<

\begin{proposition} \  %9
Let \mG be a nilpotent group .  Assume: $\forall \ n \geq 1$, $H_n(G)$ is $P$-local $-$then \mG is $P$-local.
\end{proposition}

[According to Proposition 8, 
$H_n(l_P):H_n(G) \ra H_n(G_P)$ is $P$-localizing or still, is an isomorphism, $H_n(G)$ being $P$-local.  But this means that 
$l_P:G \ra G_P$ is an isomorphism 
(cf. p. \pageref{8.22}).]\\

\begin{proposition} \ %10
Let $f:G \ra K$ be a homomorphism of nilpotent groups $-$then $f$ is $P$-localizing iff $\forall \ n \geq 1$, 
$H_n(f):H_n(G) \ra H_n(K)$ is $P$-localizing.
\end{proposition}

[Necessity: By definition, $\exists$ an isomorphism $\phi:G_P \ra K$ such that $f = \phi \circx l_P$, so 
$H_n(f) = H_n(\phi) \circx H_n(l_P)$, where $H_n(\phi)$ is an isomorphism and $H_n(l_P)$ is $P$-localizing 
(cf. Proposition 8).

Sufficiency: Since $\forall \ n \geq 1$, $H_n(K)$ is $P$-local, Proposition 9 implies that \mK is $P$-local, hence by universality, $\exists$ a homomorphism $\phi:G_P \ra K$ such that $f = \phi \circx l_P$.  
Claim: $\phi$ is an isomorphism.  In fact, $H_n(f) = H_n(\phi) \circx H_n(l_P)$, where $H_n(f)$ and $H_n(l_P)$ are 
$P$-localizing, thus $\forall \ n \geq 1$, $H_n(\phi)$ is an isomorphism, from which the claim 
(cf. p. \pageref{8.23}).]

[Note: \ Similar considerations show that if $f:G \ra K$ is a homomorphism of nilpotent groups, 
then $f$ is $P$-bijective iff $\forall \ n \geq 1$, $H_n(f):H_n(G;\Z_P) \ra H_n(K;\Z_P)$ is bijective.]\\

\begin{proposition} \ %11
Let $f:G \ra K$ be a homomorphism of nilpotent groups.  Assume: $f$ is $P$-localizing $-$then $\forall \ i \geq 0$, 
$\Gamma^i(f):\Gamma^i(G) \ra \Gamma^i(K)$ is $P$-localizing.
\end{proposition}

[On the basis of the \cd
\[
\begin{tikzcd}%[ sep=large]
{1} \ar{r} 
&{\Gamma^i(G)} \ar{d} \ar{r}
&{G} \ar{d} \ar{r}
&{G/\Gamma^i(G)} \ar{d} \ar{r}
&{1}\\
{1} \ar{r} 
&{\Gamma^i(K)} \ar{r}
&{K} \ar{r}
&{K/\Gamma^i(K)} \ar{r}
&{1}
\end{tikzcd}
,
\]
it need only be shown that $\forall \ i$, the induced map $f_i$ is $P$-localizing.  
This can be done by induction on $i$.  Indeed, the assertion is trivial if $i = 0$ and a consequence of Proposition 10 if $i = 1$, so to pass from $i$ to $i+1$, it suffices to remark that the arrow 
$\Gamma^i(G)/\Gamma^{i+1}(G) \ra$ 
$\Gamma^i(K)/\Gamma^{i+1}(K)$ is $P$-localizing (inspect the proof of Proposition 14 in $\S 5$).]\\

%%----------------------------------------------------------------------------------------------19
Application: Let \mG be a nilpotent group $-$then $\forall \ i \geq 0$, $\Gamma^i(G)_P \approx \Gamma^i(G_P)$.\\

\textbf{\small LEMMA} \  
Let 
$
\begin{cases}
\\[-.95cm]
\ \phi:G \ra K\\[-.1cm]
\ \psi:H \ra K\\[-.1cm]
\end{cases}
\vspace{0.2cm}
$
be homomorphisms of nilpotent groups $-$then $f:G \times_K H \ra$ $G_P \times_{K_P} H_P$ is $P$-localizing.

[For $f$ is clearly $P$-injective, being the restriction to $G \times_K H$ of the $P$-bijection 
$l_P \times l_P:G \times H \ra G_P \times H_P$.  
To show that $f$ is $P$-surjective, take $(g_P,h_P) \in G_P \times_{K_P} H_P$, so 
$\phi_P(g_P) = \psi_P(h_P)$.  
Choose
$
\begin{cases}\\[-.95cm]
\ g \in G\\[-.1cm]
\ h \in H\\[-.1cm]
\end{cases}
\vspace{0.2cm}
$
$\&$
$
\begin{cases}\\[-.95cm]
\ m\\[-.1cm]
\ n\\[-.1cm]
\end{cases}
\in S_P:
$
$
\begin{cases}\\[-.95cm]
\ g_P^m = l_P(g)\\[-.1cm]
\ h_P^m = l_P(h)\\[-.1cm]
\end{cases}
$
$\implies$ 
$l_P \circx \phi(g^n) =$ 
$\phi_P \circx l_P(g^n) =$ 
$\phi_P(g_P^{mn}) =$ 
$\psi_P(h_P^{mn}) =$ 
$\psi_P \circx l_P(h^m) =$
$l_P \circx \psi(h^m)$ 
$\implies$ 
$\phi(g^n) = \psi(h^m)k$ $(k \in \ker l_P)$.  
Choose $t \in S_P$ : $k^t = e$.  
Fix $d$ : $\nil K \leq d$ $-$then 
$\phi(g^n)^{t^d} =$ 
$(\psi(h^m)k)^{t^d} =$ 
$(\psi(h^m)^{t^d}$   
(cf. p. \pageref{8.24}) 
$\implies$  
$(g^{nt^d},h^{mt^d}) \in G \times_K H$ 
$\implies$ 
$(g_P,h_P)^{mnt^d} =$
$f(g^{nt^d},h^{mt^d})$ 
$\implies$ 
$(g_P,h_P)^{mnt^d} \in \im f$, i.e., $f$ is $P$-surjective.  
Since $G_P \times_{K_P} H_P$ is necessarily $P$-local, it follows that $f$ is $P$-localizing.]\\

\textbf{\small LEMMA} \  
Let 
$
\begin{cases}
\\[-.95cm]
\ \phi:G \ra K\\[-.1cm]
\ \psi:G \ra K\\[-.1cm]
\end{cases}
\vspace{0.2cm}
$
be homomorphisms of nilpotent groups $-$then $f:\eq(\phi,\psi) \ra \eq(\phi_P,\psi_P)$ is $P$-localizing.

[Imitate the argument used in the preceding proof.]\\

\begin{proposition} \ %12
$L_P:\bNIL \ra \bNIL_P$ preserves finite limits.
\end{proposition}

[Combine the foregoing lemmas.]\\

\begingroup%%----------------------------------->>
\fontsize{9pt}{11pt}\selectfont
\textbf{\small EXAMPLE} \  
Let \mG be a nilpotent group; let 
$
\begin{cases}
\ G^\prime\\
\ G\pp
\end{cases}
$
be subgroups of \mG $-$then $(G^\prime \cap G\pp)_P \approx G_P^\prime \cap G_{P}\pp$.\\
\endgroup %%------------------------------------<<

\label{8.26}
\begingroup%%----------------------------------->>
\fontsize{9pt}{11pt}\selectfont
\textbf{\small FACT} \  
Let \mG be a nilpotent group, $\{g_\mu\}$ a subset of \mG.  
Fix $n \in \N$.  Assume: 
(1) The set $\{g_\mu[G,G]\}$ generates $G/[G,G]$; 
(2) Each $g_\mu$ is a product of $n^\text{th}$ powers $-$then the map 
$
\begin{cases}
\ G \ra G\\
\ g \ra g^n
\end{cases}
$
is surjective.
\vspi
[$\Gamma^i(G)/\Gamma^{i+1}(G)$ has $n^\text{th}$ roots (consider the arrow 
$\otimes^{i+1} (G/[G,G]) \ra$ 
$\Gamma^i(G)/\Gamma^{i+1}(G)$, thus 
$G/\Gamma^{i+1}(G)$ has $n^\text{th}$ roots (consider the central extension 
$1 \ra$ 
$\Gamma^i(G)/\Gamma^{i+1}(G) \ra$ 
$G/\Gamma^{i+1}(G) \ra$ 
$G/\Gamma^{i}(G) \ra 1$).]\\
\endgroup %%------------------------------------<

\begingroup%%----------------------------------->>
\fontsize{9pt}{11pt}\selectfont
\textbf{\small EXAMPLE} \  
Let \mG be a nilpotent group; let \mK be a subgroup of \mG.  Write 
$\norx_{G} K$ for the normal closure of \mK in \mG, 
$\norx_{G_P} K_P$ for the normal closure of $K_P$ in $G_P$ $-$then 
$(\norx_{G} K)_P \approx \norx_{G_P} K_P$.\\
\endgroup %%------------------------------------<<

\begingroup%%----------------------------------->>
\fontsize{9pt}{11pt}\selectfont
\textbf{\small EXAMPLE} \  
Let \mG be a nilpotent group; let 
$
\begin{cases}
\ G^\prime\\
\ G\pp
\end{cases}
$
be subgroups of \mG.  Write 
$\langle G^\prime,G\pp \rangle$ for the subgroup of \mG generated by $G^\prime \cup G\pp$, 
$\langle G_P^\prime,G_P\pp \rangle$ for the subgroup of $G_P$ generated by $G_P^\prime \cup G_P\pp$ $-$then 
$\langle G^\prime,G\pp \rangle_P \approx \langle G_P^\prime,G_P\pp \rangle$.\\
\endgroup %%------------------------------------<<

%%----------------------------------------------------------------------------------------------20
Notation: Given groups 
$
\begin{cases}
\\[-.95cm]
\ G^\prime\\[-.1cm]
\ G\pp\\[-.1cm]
\end{cases}
\hspace{-.2cm}, \ 
$
the kernel $\car(G^\prime,G\pp)$ of the epimorphism 
$G^\prime * G\pp \ra G^\prime \times G\pp$ is the 
\un{cartesian subgroup}
\index{cartesian subgroup} 
of $G^\prime * G\pp$.  
It is freely generated by $\{[g^\prime,g\pp]$:$g^\prime \neq e \ \& \ g\pp \neq e\}$.  If 
$
\begin{cases}
\\[-.95cm]
\ G^\prime\\[-.1cm]
\ G\pp\\[-.1cm]
\end{cases}
$
are subgroups of \mG, then $\nabla_G(\car(G^\prime,G\pp)) = [G^\prime,G\pp]$, where 
$\nabla_G:G * G \ra G$ is the folding map.

Suppose that 
$
\begin{cases}
\\[-.95cm]
\ G^\prime\\[-.1cm]
\ G\pp\\[-.1cm]
\end{cases}
$
are in $\bNIL^d$.  Put $G^\prime *_d G\pp = G^\prime * G\pp/\Gamma^d(G^\prime * G\pp)$ $-$then 
$G^\prime *_d G\pp$ is the coproduct in $\bNIL^d$.  Call $\car_d(G^\prime,G\pp)$ the kernel of the epimorphism 
$G^\prime *_d G\pp \ra G^\prime \times G\pp$, so 
$\car_d(G^\prime,G\pp) \approx$ 
$\car(G^\prime,G\pp)/\Gamma^d(G^\prime *_d G\pp)$ .\\

\begingroup%%----------------------------------->>
\fontsize{9pt}{11pt}\selectfont
\textbf{\small FACT} \  
$\bNIL^d$ is a reflective subcategory of \bGR, hence is complete and cocomplete.
\vspi
[Note: \ \bNIL is finitely complete but not finitely cocomplete.]\\
\endgroup %%------------------------------------<<

\begingroup%%----------------------------------->>
\fontsize{9pt}{11pt}\selectfont
\textbf{\small FACT} \  
Let \mG be a nilpotent group $-$then the \cd 
\begin{tikzcd}[ sep=large]
{G} \ar{d} \ar{r} &{G_{\ov{P}}} \ar{d}\\
{G_P} \ar{r} &{G_{\Q}}
\end{tikzcd}
is simultaneously a pullback square and a pushout square in \bNIL and the arrow 
$
\begin{cases}
\ G_P \ra G_{\Q}\\
\ G_{\ov{P}} \ra G_{\Q}
\end{cases}
$
is a 
$
\begin{cases}
\ \ov{P} \text{-bijection}\\
\ P \text{-bijection}
\end{cases}
. \ 
$
\vspace{0.35cm}
\endgroup %%------------------------------------<<

\begin{proposition} \ %13
Let \mG be a nilpotent group; let
$
\begin{cases}
\ G^\prime\\
\ G\pp
\end{cases}
$
be subgroups of \mG $-$then $l_P:G \ra G_P$ restricts to an arrow 
$f:[G^\prime,G\pp] \ra [G_P^\prime,G_P\pp]$ which is $P$-localizing.
\end{proposition} 

[Trivially, $f$ is $P$-injective.  To check that $f$ is $P$-surjective, look first at the commutative diagram
\[
\begin{tikzcd}%[sep=large]
{1} \ar{r} 
&{\car_d(G^\prime,G\pp)} \ar{d} \ar{r}
&{G^\prime *_d G\pp} \ar{d} \ar{r}
&{G^\prime \times G\pp} \ar{d} \ar{r}
&{1}\\
{1} \ar{r} 
&{\car_d(G_P^\prime,G_P\pp)} \ar{r}
&{G_P^\prime *_d G_P\pp} \ar{r}
&{G_P^\prime \times G_P\pp} \ar{r}
&{1}
\end{tikzcd}
,
\]
it being assumed that $\nil G \leq d$.  Since $L_P$ preserves colimits, 
$(G^\prime *_d G\pp)_P \approx$ 
$(G_P^\prime *_d G_P\pp)_P$ 
(cf. p. \pageref{8.25}).  Therefore the arrow 
$G^\prime *_d G\pp \ra G_P^\prime *_d G_P\pp$ is $P$-bijective, thus the same is true of the arrow 
$\car_d(G^\prime,G\pp) \ra \car_d(G_P^\prime,G_P\pp)$. \ 
Consequently, upon forming the commutative square
\begin{tikzcd}%[sep=large]
{\car_d(G^\prime,G\pp)} \ar{d} \ar{r} &{[G^\prime,G\pp]} \ar{d}{f}\\
{\car_d(G_P^\prime,G_P\pp)} \ar{r} &{[G_P^\prime,G_P\pp]}
\end{tikzcd}
in which the horizontal arrows are the epimorphisms induced by the folding maps, it is seen that $f$ is $P$-surjective.  
Turning to the verification that $[G_P^\prime,G_P\pp]$ is $P$-local, there is no $S_P$-torsion and $\forall \ n \in S_P$,  
$G_P^\prime *_d G_P\pp$ has $n^\text{th}$ roots (consider generators 
(cf. p. \pageref{8.26})), so $\forall \ n \in S_P$, 
$\car_d(G_P^\prime,G_P\pp)$ has $n^\text{th}$ roots and this suffices.]\\

%%----------------------------------------------------------------------------------------------21
Application:  Let \mG be a nilpotent group; let
$
\begin{cases}
\ G^\prime\\
\ G\pp
\end{cases}
$
be subgroups of \mG $-$then $[G^\prime,G\pp]_P \approx [G_P^\prime,G_P\pp]$.\\

\label{9.10} %dmc mnft
Let \mG and $\pi$ be groups.  Suppose that \mG operates on $\pi$, i.e., suppose given a homomorphism 
$\chi:G \ra \Aut \pi$ $-$then $\chi$ determines a homomorphism $\chi_P:G \ra \Aut \pi_P$, thus \mG operates on $\pi_P$.\\

\begingroup%%----------------------------------->>
\fontsize{9pt}{11pt}\selectfont
\textbf{\small FACT} \  
If \mG operates on $\pi$ and if $\pi$ is nilpotent, then 
$\Gamma_\chi^i(\pi)_P \approx \Gamma_{\chi_P}^i(\pi_P)$ (here the notation is that of 
p. \pageref{8.27}).  In particular: 
$\pi$ $\chi$-nilpotent $\implies$ $\pi_P$ $\chi_P$-nilpotent.
\vspi
[Use induction and Proposition 13, so that 
$[\pi,\Gamma_\chi^i(\pi)]_P \approx$ 
$[\pi_P,\Gamma_\chi^i(\pi)_P] \approx$ 
$[\pi_P,\Gamma_{\chi_P}^i(\pi_P)]$.]\\
\endgroup %%------------------------------------<<

Given groups \mG and $\pi$, let $\Hom_\nil(G,\Aut \pi)$ be the subset of $\Hom(G,\Aut \pi)$ consisting of those 
$\chi$ such that $\pi$ is $\chi$-nilpotent.

[Note: \ In order that $\Hom_\nil(G,\Aut \pi)$ be nonempty, it is necessary that $\pi$ be nilpotent 
(cf. p. \pageref{8.28}).]

 Suppose that \mG and $\pi$ are nilpotent.

\indent\indent ($\nil_1$) \quad The arrow 
$\Hom(G,\Aut \pi) \ra \Hom(G,\Aut \pi_P)$ restricts to an arrow \\
$\Hom_\nil(G,\Aut \pi)$ $\ra$ $\Hom_\nil(G,\Aut \pi_P)$.

[For, as noted above, $\pi$ $\chi$-nilpotent $\implies$ $\pi_P$ $\chi_P$-nilpotent.]

\indent\indent ($\nil_2$) \quad There is an arrow 
$\Hom_\nil(G,\Aut \pi) \ra \Hom_\nil(G_P,\Aut \pi_P)$ that sends $\chi$ to $\ov{\chi}_P$, where 
$\ov{\chi}_P \circx l_P = \chi_P$.

[This semidirect product $\Pi = \pi \rtimes_{\chi} G$ is nilpotent 
(cf. p. \pageref{8.29}).  
Localize the split short exact sequence $1 \ra \pi \ra \Pi \ra G \ra 1$ and consider the associated action of $G_P$ on 
$\pi_P:\Pi_P = \pi_P \rtimes_{\ov{\chi}_P} G_P$.]

\indent\indent ($\nil_3$) \quad The arrow 
$\Hom(G_P,\Aut \pi_P) \ra \Hom(G,\Aut \pi_P)$ restricts to an arrow 
$\Hom_\nil(G_P,\Aut \pi_P) \ra \Hom_\nil(G,\Aut \pi_P)$ which is bijective.  
If $\hbar$ is its inverse, then $\forall \ \chi$, 
$\hbar(\chi_P) = \ov{\chi}_P$.

[Implicit in the construction of $\hbar$ is the relation 
$\Gamma_{\chi_P}^i(\pi_P) \approx \Gamma_{\ov{X}_P}^i (\pi_P)$.]\\

\begingroup%%----------------------------------->>
\fontsize{9pt}{11pt}\selectfont
\textbf{\small FACT} \  
Suppose that \mG operates nilpotently on $\pi$ and $\pi$ is abelian $-$then for any half exact functor $F:\bAB \ra \bAB$, \mG operates nilpotently on $F \pi$.\\
\endgroup %%------------------------------------<<

\begingroup%%----------------------------------->>
\fontsize{9pt}{11pt}\selectfont
\textbf{\small EXAMPLE} \  
Fix a path connected topological space \mX and let $\pi$ be a nilpotent $G$-module $-$then $\forall \ n \geq 0$, $H_n(X;\pi)$ is a 
nilpotent $G$-module.\\
\endgroup %%------------------------------------<<

\begin{proposition} \ %14
Let \mG be a nilpotent group, \mM a nilpotent $G$-module $-$then $\forall \ n \geq 0$, the arrow 
$H_n(G;M) \ra H_n(G_P;M_P)$ is $P$-localizing.
\end{proposition}

%%----------------------------------------------------------------------------------------------22
[From the definitions, \ 
$H_0(G;M) \approx M/\Gamma_\chi^1(M)$ and 
$H_0(G_P;M_P) \approx M_P/\Gamma_{\ov{\chi}_P}^1(M_P)$.  \ 
Accordingly, \ since $L_P$ is exact, \ 
$(M/\Gamma_\chi^1(M))_P$ $\approx$ 
$M_P/\Gamma_\chi^1(M)_P$ $\approx$ 
$M_P/\Gamma_{{\chi}_P}^1(M_P)$ $\approx$ 
$M_P/$ $\Gamma_{\ov{\chi}_P}^1(M_P)$, 
thereby dispensing with the case $n = 0$.  
Assume henceforth that $n \geq 1$.  
Matters are plain when $\nil_\chi M = 0$.  
If $\nil_\chi M = 1$, i.e., if \mG operates trivially on \mM, then $G_P$ operates trivially on $M_P$ and one can apply the universal coefficient theorem, 
in conjuction with Proposition 10, to derive the desired conclusion.  
Arguing inductively, suppose that $\nil_\chi M \leq d$ $(d > 1)$ and that the assertion holds for operations having degree of nilpotency $\leq d - 1$.  
Consider the short exact sequence 
$0 \ra$ 
$\Gamma_\chi^1 (M) \ra$ 
$M \ra$ 
$M/\Gamma_\chi^1 (M) \ra 0$.  
The degree of nilpotency of the induced action of \mG on $\Gamma_\chi^1(M)$ is $\leq d - 1$, 
while that of \mG on $M/\Gamma_\chi^1 (M)$ is $\leq 1$.  
Comparison of the long exact sequence 
$\cdots \ra$ 
$H_{n+1}(G;M/\Gamma_\chi^1 (M)) \ra$ 
$H_n(G;\Gamma_\chi^1(M)) \ra$ 
$H_n(G;M) \ra$ 
$H_{n}(G;M/\Gamma_\chi^1 (M)) \ra$ 
$H_{n-1}(G;\Gamma_\chi^1(M)) \ra \cdots$ 
with its local companion terminates the proof.]\\

Application: Let \mG be a nilpotent group, \mM a nilpotent $G$-module $-$then $\forall \ n \geq 0$, 
$H_n(G;M)_P \approx H_n(G_P;M_P)$.\\

\label{9.39}
Given a group \mG, $G$-\bACT is the category whose objects are the groups on which \mG operates to the left and whose morphisms are the equivariant homomorphisms.  
An object $\pi$ in $G$-\bACT is really a pair$(\chi,\pi)$, where 
$\chi:G \ra \Aut \pi$.  
One says that $\pi$ is 
\un{$P$-local}
\index{P-local (object in $G$-\bACT)} 
or that \mG operates $P$-locally on $\pi$ if $\forall \ n \in S_P$ $\&$ $\forall \ g \in G$, the map $\pi \ra \pi$ that sends $\alpha$ to 
$\alpha(\chi(g)\alpha) \cdots (\chi(g^{n-1})\alpha)$ is bijective, so $\pi$ is necessarily a $P$-local group.  
Denote by 
$G\text{-}\bACT_P$ the full subcategory of $G$-\bACT whose objects are the $P$-local $\pi$ $-$then $G\text{-}\bACT_P$ is a reflective subcategory of $G$-\bACT with reflector $L_{G,P}$.  This can be seen by applying the relative subcategory theorem.  Thus let 
$F_G$ be the free $G$-group on one generator $*$, i.e., the free group on the symbols $g \cdot *$ $(g \in G)$ with the obvious left action.  Write $S_{G,P}$ for the set of $G$-maps
$
\begin{cases}
\ F_G \ra F_G\\
\ * \ra \rho_g^n(*)
\end{cases}
(n \in S_P),
$
where $\rho_g^n(*) = *(g\cdot*) \cdots (g^{n-1}\cdot *)$.  Working through the definitions, one finds that 
$\Ob G\text{-}\bACT_P = S_{G,P}^\perp$.

\label{8.31} %dmc mnft
Example: $\pi \rtimes_\chi G$ is a $P$-local group iff \mG operates $P$-locally on $\pi$ and \mG is a $P$-local group.

[Note: \ It is a corollary that if \mG is $S_{\ov{P}}$-torsion, then every $P$-local group in $G$-\bACT is actually in 
$G\text{-}\bACT_P$.  Proof: Consider the short exact sequence 
$1 \ra \pi \ra \pi \rtimes_\chi G \ra G \ra 1$ and quote the generality on 
p. \pageref{8.30}.]\\

\begingroup%%----------------------------------->>
\fontsize{9pt}{11pt}\selectfont
\textbf{\small FACT} \  
Let $f:G \ra K$ be a homomorphism of groups $-$then the functor $f^*:K\text{-}\bACT \ra G\text{-}\bACT$ has a left adjoint 
$f_*:G\text{-}\bACT \ra K\text{-}\bACT$.
\vspi
%%----------------------------------------------------------------------------------------------23
[Let 
$
\begin{cases}
\ \ov{\pi}_G\\
\ \ov{\pi}_K
\end{cases}
$
be the normal closure of $\pi$ in 
$
\begin{cases}
\ \pi * G\\
\ \pi * K
\end{cases}
$
.  There are pushout squares
$
\begin{tikzcd}[ sep=large]
{\pi} \ar{d} \ar{r} &{\pi * G} \ar{d}\\
{*} \ar{r} &{G}
\end{tikzcd}
,
$
$
\begin{tikzcd}[ sep=large]
{\pi} \ar{d} \ar{r} &{\pi * K} \ar{d}\\
{*} \ar{r} &{K}
\end{tikzcd}
,
$
short exact sequences 
$1 \ra \ov{\pi}_G \ra \pi * G \ra G \ra 1$, 
$1 \ra \ov{\pi}_K \ra \pi * K \ra K \ra 1$, 
and a commutative diagram 
$
\begin{tikzcd}[ sep=large]
{\pi * G} \ar{d}[swap]{\id * f} \ar{r} &{G} \ar{d}{f}\\
{\pi * K} \ar{r} &{K}
\end{tikzcd}
.
$
\ 
Let $\pi_{\chi,G}$ be the normal closure in $\pi * G$ of the words $g \alpha  g^{-1} (\chi(g) \alpha)^{-1}$, 
$f(\pi_{\chi,G})$ 
the normal closure in $\pi * K$ of the words 
$\id * f(g \alpha g^{-1} (\chi(g)\alpha)^{-1})$ 
$-$then $\pi_{\chi,G}$ is a normal subgroup of $\ov{\pi}_G$, the quotient 
$\ov{\pi}_G/\pi_{\chi,G}$ is equivariantly isomorphic to $\pi$, and $f(\pi_{\chi,G}) \subset \ov{\pi}_K$.  
Definition: $f_*(\pi) = \ov{\pi}_K/f_(\pi_{\chi,G})$ the action of \mK being conjugation.  
Note that the arrow 
$\pi \ra f^*f_*(\pi)$ is equivariant.]\\
\endgroup %%------------------------------------<<

\begingroup%%----------------------------------->>
\fontsize{9pt}{11pt}\selectfont
\textbf{\small EXAMPLE} \  
For any homomorphism $f:G \ra K$ of groups, the composite $L_{K,P} \circx f_*$ is a functor 
$G\text{-}\bACT \ra K\text{-}\bACT  \ra K\text{-}\bACT_P$.  
Specialize and take $K = G_P$, $f = l_P$.  
Given $\pi$ in $G\text{-}\bACT$, form 
$\pi \rtimes_\chi G$ $-$then its localization $(\pi \rtimes_\chi G)_P$ is isomorphic to a semidirect product 
$? \rtimes G_P$ and ? can be identified with $L_{G_P,P} \circx l_{P,*}(\pi)$.\\
\endgroup %%------------------------------------<<

\label{9.44}
\label{9.47}
\label{9.48}
Given a group \mG, a 
\un{$P$-local $G$-module}
\index{P-local $G$-module} 
is a $G$-module on which $G$ operates $P$-locally.  Every $P$-local $G$-module is a $P$-local abelian group.

[Note: \ If $(\Z[G])_{S_P}$ is the localization of $\Z[G]$ at the multiplicative closure of the 
$1 + g + \cdots + g^{n-1}$ $(n \in S_P)$, then the $P$-local $G$-modules are the $(\Z[G])_{S_P}$-modules.  
When \mG is trivial, $(\Z[G])_{S_P}$ reduces to $\Z_P$.]\\

\begin{proposition} \ %15
Suppose that \mG is $S_P$-torsion $-$then every $P$-local $G$-module is trivial.
\end{proposition}

[In $\Z[G]$, consider the identity $g^n - 1 = (g - 1) (1 + g + \cdots + g^{n-1})$.]\\

\begingroup%%----------------------------------->>
\fontsize{9pt}{11pt}\selectfont
\textbf{\small FACT} \  
Let $0 \ra M^\prime \ra M \ra M\pp \ra 0$ be a short exact sequence of $G$-modules.  
Assume: Two of the modules are $P$-local 
$-$then so is the third.\\
\endgroup %%------------------------------------<<

\begingroup%%----------------------------------->>
\fontsize{9pt}{11pt}\selectfont
\textbf{\small EXAMPLE} \  
Suppose that \mM is a $P$-local $G$-module and \mN is a nilpotent $G$-module $-$then $M \otimes N$, 
$\Tor(M,N)$, $\Hom(N,M)$, $\Ext(N,M)$ are $P$-local $G$-modules.\\
\endgroup %%------------------------------------<<

\begingroup%%----------------------------------->>
\fontsize{9pt}{11pt}\selectfont
\textbf{\small FACT} \  
Let $G \ra \pi$ be a homomorphism of groups $-$then every $P$-local $\pi$-module is a $P$-local $G$-module.\\
\endgroup %%------------------------------------<<

%%----------------------------------------------------------------------------------------------24
\begingroup%%----------------------------------->>
\fontsize{9pt}{11pt}\selectfont
\textbf{\small EXAMPLE} \  
Suppose that $1 \ra G^\prime \ra G \ra G\pp \ra 1$ is a central extension of groups.  Let \mM be a $P$-local $G$-module $-$then 
$\forall \ n \geq 0$, the action of $G\pp$ on $H_n(G^\prime;M)$ and $H^n(G^\prime;M)$ is $P$-local.\\
\endgroup %%------------------------------------<<

Given a group \mG, a 
\un{$P[G]$-module}
\index{P[G]-module} 
is a $P$-local $G_P$-module.  
Every $P[G]$-module is a $P$-local $G$-module via $l_P:G \ra G_P$.

[Note: \ A $G_P$-module \mM is a $P[G]$-module iff the corresponding semidirect product $M \rtimes G_P$ is a 
$P$-local group 
(cf. p. \pageref{8.31}).]

Example: Suppose that \mG is nilpotent.  Let \mM be a nilpotent $G_P$-module which is $P$-local as an abelian group 
$-$then \mM is a $P[G]$-module.\\

\label{9.46}
\begingroup%%----------------------------------->>
\fontsize{9pt}{11pt}\selectfont
\textbf{\small FACT} \  
Let \mM be a $P[G]$-module $-$then $H^0(G_P;M) = H^0(G;M)$, i.e., the $G_P$-invariants in \mM are equal to the $G$-invariants in \mM.
\vspi
[Let $m \in M^G$.  Define homomorphisms $\phi,\psi:G_P \ra M \rtimes G_P$ by the rules 
$\phi(g) = (g \cdot m - m, g)$, 
$\psi(g) = (0,g)$: $\phi \circx l_P = \psi \circx l_P$ $\implies$ $\phi = \psi$, $M \rtimes G_P$ being $P$-local, 
i.e., $m \in M^{G_P}$.]\\
\endgroup %%------------------------------------<<

\begin{proposition} \ %16
Let \mG be a nilpotent group, \mM a $P[G]$-module $-$then $\forall \ n \geq 0$, 
$H_n(G;M) \approx H_n(G_P;M)$.
\end{proposition}

[It suffices to treat the case of an abelian \mG.  There are short exact sequences 
$0 \ra$ 
$\ker l_P \ra$ 
$G \ra $ 
$\im l_P \ra 0,$ \ 
$0 \ra$ 
$\im l_P \ra$ 
$G_P \ra $ 
$\coker l_P \ra 0$ \ 
and associated LHS spectral sequences.  Since $\ker l_P$  is $S_P$-torsion, $H_q(\ker l_P) \in \sC_P$ $(q > 0)$.  
But the action of $\ker l_P$ on \mM is by definition trivial, and as an abelian group, \mM is $P$-local, thus the universal coefficient theorem implies that $H_q(\ker l_P;M) = 0$ $(q > 0)$.  So, $\forall \ n \geq 0$, 
$H_n(G;M) \approx H_n(\im l_P;M)$.  
On the other hand, from the above, the action of $\coker l_P$ on the 
$H_q(\im l_P;M)$ is $P$-local, hence  trivial (cf. Proposition 15).  
Appealing once again to the universal coefficient theorem, it follows that 
$H_p(\coker l_P;H_q(\im l_P;M)) = 0$ $(p > 0)$.  So, $\forall \ n \geq 0$, $H_n(\im l_P;M) \approx H_n(G_P;M)$.]\\

\begingroup%%----------------------------------->>
\fontsize{9pt}{11pt}\selectfont
\textbf{\small FACT} \  
Let \mG be a nilpotent group, \mM a $P[G]$-module $-$then $\forall \ n \geq 0$, $H^n(G_P;M) \approx H^n(G;M)$.\\
\endgroup %%------------------------------------<<

\begingroup%%----------------------------------->>
\fontsize{9pt}{11pt}\selectfont
\textbf{\small EXAMPLE} \  
The preceding result can fail if \mM is not a $P[G]$-module.  Thus fix $P \neq \bPi$ and take $G = \Z$ : $H^2(\Z;\Q[\Z_P]) = 0$ (since $\Z$ has cohomological dimension one) but $H^2(\Z_P;\Q[\Z_P]) \neq 0$ 
(cf. p. \pageref{8.32}).\\
\endgroup %%------------------------------------<<

\label{9.57}
\begingroup%%----------------------------------->>
\fontsize{9pt}{11pt}\selectfont
\textbf{\small FACT} \  
Let \mG be a finite group $-$then $\ker l_P$ is $S_P$-torsion iff $\forall \ n \ \geq 0$, 
$H_n(G;M) \approx$ 
$H_n(G_P;M)$, where \mM is any $P[G]$-module.\\
\endgroup %%------------------------------------<<

There is another reflective subcategory of \bGR that one can attach to a given $P \subset \bPi$ whose definition is homological in character.  The associated reflector agrees with $L_P$ on \bNIL but differs from $L_P$ on \bGR.\\

%%----------------------------------------------------------------------------------------------25
\index{Colimit Lemma}
\textbf{\small COLIMIT LEMMA} \  
Let \bC be a cocomplete category with the property that there exists a set $S_0 \subset \Ob\bC$ 
such that each object in \bC is a filtered colimt of objects in $S_0$.  
Let $F:\bC \ra \bSET_*$ be a functor which preserves filtered colimts $-$then there exists a set 
$K_0 \subset \ker F$ such that each $X \in \ker F$ is a filtered colimit of objects in $K_0$.

[Note: \ As the notation suggests, $\ker F = \{X:FX = *\}$.]\\

\label{9.58}
Let \mA be an abelian group $-$then a homomorphism $f:G \ra K$ of groups is said to be an 
\un{$HA$-homomorphism}
\index{HA-homomorphism (of abelian groups)} 
if 
$f_*: H_1(G;A) \ra$ $H_1(K;A)$ is bijective and 
$f_*: H_2(G;A) \ra$ $H_2(K;A)$ is surjective.  
Example: An $H\Z$-homomorphism of nilpotent groups is an isomorphism 
(cf. p. \pageref{8.33}).\\
\indent\indent ($HA$-localization) \quad Let $S_{HA} \subset \Mor\bGR$ be the class of $HA$-homomorphisms $-$then 
$S_{HA}^\perp$ is the object class of a reflective subcategory $\bGR_{HA}$ of \bGR.  
The reflector 
$
L_{HA}: 
\begin{cases}
\ \bGR \ra \bGR_{HA}\\
\ G \ra G_{HA}
\end{cases}
$
is called 
\un{$HA$-localization}
\index{HA-localization} 
and the objects in $\bGR_{HA}$ are called the 
\un{$HA$-local}
\index{HA-local(groups)} 
groups.

[In order to apply the reflective subcategory theorem, it suffices to exhibit a set $S_0 \subset S_{HA}$: 
$S_0^\perp = S_{HA}^\perp$.  
For this purpose, put $\bC = \bGR(\ra)$ $(\approx [\btwo,\bGR])$ and let 
$F:\bC \ra \bSET_*$ be the functor that sends 
$f:G \ra K$ to $\ker_1 \oplus \hspace{0.07cm} \cokersub_1 \oplus \cokersub_2$, where 
$\ker_1$ is the kernel of 
$f_*:H_1(G;A) \ra H_1(K;A)$ and $\cokersub_i$ is the cokernel of 
$f_*:H_i(G;A) \ra H_i(K;A)$ $(i = 1, 2)$.  
Owing to the colimit lemma, there exists a set $S_0 \subset S_{HA}$ such that each element of 
$S_{HA}$ is a filtered colimit of elements in $S_0$, so $S_0^\perp = S_{HA}^\perp$.]

[Note: \ In general, the containment $S_{HA} \subset S_{HA}^{\perp\perp}$ is strict (see below).]

When $A = \Z_P$, the ``$\Z$'' is dropped from the notation, thus one writes $S_{HP}$ for the class of 
$HP$-homomophisms and 
$
L_{HP}:
\begin{cases}
\ \bGR \ra \bGR_{HP}\\
\ G \ra G_{HP}
\end{cases}
$
for the associated reflector, the objects in $\bGR_{HP}$ then being referred to as the 
\un{$HP$-local}
\index{HP-local(groups)} 
groups.  
Example: Every abelian $P$-local group is $HP$-local.

[Note: \ In the two extreme cases, viz. $P = \emptyset$ or $P = \bPi$, $HP$ is replaced by $H\Q$ or $H\Z$.]\\

\begin{proposition} \ %17
Every $HP$-local group is $P$-local.
\end{proposition}

[The homomorphisms 
$
\begin{cases}
\ \Z \ra \Z \\
\ 1 \ra n
\end{cases}
(n \in S_P)
$
are $HP$-homomorphisms, thus $S_P \subset S_{HP}$ $\implies$ 
$S_{HP}^\perp = \Ob\bGR_{HP} \subset \Ob\bGR_P = S_P^\perp$.]\\

Consequently, there is a natural transformation $L_P \ra L_{HP}$.

[Note: \ For any \mG, the arrow of localization $l_P:G \ra G_P$ is an $HP$-homomophism 
(cf. p. \pageref{8.34}).  As regards $l_{HP}:G \ra G_{HP}$, it too is an $HP$-homomorphism 
(cf. p. \pageref{8.35}ff), although a priori it can only be said that $l_{HP} \in S_{HP}^{\perp\perp}$.]\\

%%----------------------------------------------------------------------------------------------26
\begin{proposition} \ %18
Let $f:G \ra K$ be an $HP$-homomorphism $-$then $\forall \ i \geq 0$, the induced map 
$(G/\Gamma^i(G))_P \ra (K/\Gamma^i(K))_P$ is an isomorphism.
\end{proposition}

[Taking into account Propositions 6 and 8, one has only to repeat the proof of Proposition 14 in $\S 5$.]\\

\textbf{\small LEMMA} \  
Let $1 \ra G^\prime \ra G \ra G\pp \ra 1$ be a central extension of groups.  Assume: $G^\prime$ is $P$-local $-$then in any commutative diagram 
\begin{tikzcd}%[ sep=large]
{K} \ar{d}[swap]{f} \ar{r} &{G} \ar{d}\\
{L} \ar{r} &{G\pp}
\end{tikzcd}
of groups, where $f:K \ra L$ is an $HP$-homomophism, there is a unique lifting 
\begin{tikzcd}%[ sep=large]
{K} \ar{d}[swap]{f} \ar{r} &{G} \ar{d}\\
{L} \ar[dashed]{ru} \ar{r} &{G\pp}
\end{tikzcd}
rendering the triangles commutative.

[Suppose that
$
\begin{cases}
\ \phi \\
\ \psi
\end{cases}
$
are liftings and $\lambda:L \ra G^\prime$ is a homomorphism such that $\phi(l) = \psi(l)\lambda(l)$ $(l \in L)$.  
Since $\lambda \circx f$ is trivial and 
$\Z_P \otimes (K/[K,K]) \approx$ 
$\Z_P \otimes (L/[L,L])$, it follows that $\lambda$ is trivial, hence $\phi = \psi$, which settle uniqueness.  
Existence can be established by passing to Eilenberg-MacLane spaces and using obstruction theory 
(cf. p. \pageref{8.36}).]\\

\begin{proposition} \ %19
Let $1 \ra G^\prime \ra G \ra G\pp \ra 1$ be a central extension of groups.  Assume: $G^\prime$ is $P$-local and $G\pp$ is 
$HP$-local $-$then \mG is $HP$-local.
\end{proposition}

[The claim is that $f \perp G$ for every $HP$-homomorphism $f:K \ra L$.  This, however, is obviously implied by the lemma.]

[Note: \ Changing the assumption to $G\pp$ is $P$-local changes the conclusion to \mG is $P$-local (but, of course, the proof is different).]\\

\label{10.2}
Application: If \mG is nilpotent, then $G_P \approx G_{HP}$ and 
$\restr{L_P}{\bNIL} \approx \restr{L_{HP}}{\bNIL}$.

[Note: \ It is not necessary to use Proposition 19 to make this deduction.  
Thus let \mG be a nilpotent $P$-local group with $\nil G \leq d$ $-$then for any $HP$-homomorhpism $K \ra L$, $\Hom(L,G) \approx \Hom(K,G)$.  
Proof: $\bNIL^d$ is a reflective subcategory of \bGR, hence 
$\Hom(L,G) \approx \Hom(L/\Gamma^d(L),G)$, 
$\Hom(K,G) \approx \Hom(K/\Gamma^d(K),G)$,
and $\bNIL_P^d$ is a reflective subcategory of $\bNIL^d$, hence
$\Hom(L/\Gamma^d(L),G) \approx \Hom((L/\Gamma^d(L))_P,G)$, 
$\Hom(K/$ $\Gamma^d(K)),G) \approx \Hom((K/\Gamma^d(K))_P,G)$, 
And: 
$(K/\Gamma^d(K))_P \approx (L/\Gamma^d(L))_P$ (cf. Proposition 18).]\\

\begingroup%%----------------------------------->>
\fontsize{9pt}{11pt}\selectfont
\textbf{\small FACT} \  
Suppose that \mG is a group such that for some $i$, $\Gamma^i(G)/\Gamma^{i+1}(G)$ is $S_P$-torsion $-$then 
$G_{HP} \approx$ 
$(G/\Gamma^i(G))_P$.
\vspi
[The short exact sequence 
$1 \ra \Gamma^i(G) \ra G \ra G/\Gamma^i(G) \ra 1$ leads to an exact sequence 
$H_2(G;\Z_P) \ra$ 
$H_2(G/\Gamma^i(G);\Z_P) \ra$
$\Z_P \otimes (\Gamma^i(G)/\Gamma^{i+1}(G)) \ra$ 
$H_1(G;\Z_P) \ra$ 
$H_1(G/\Gamma^i(G);\Z_P) \ra 0$.  
Therefore the arrow 
%%----------------------------------------------------------------------------------------------27
$G \ra$ 
$G/\Gamma^i(G)$ is an $HP$-homomorphism $\implies$ 
$G_{HP} \approx$ 
$(G/\Gamma^i(G))_{HP}$ or still, 
$G_{HP} \approx$ 
$(G/\Gamma^i(G))_{P}$, 
$G/\Gamma^i(G)$ being nilpotent.]\\
\endgroup %%------------------------------------<<

\label{9.62} %dmc mnft

\begingroup%%----------------------------------->>
\fontsize{9pt}{11pt}\selectfont
\textbf{\small EXAMPLE} \  
The $HP$-localization of every finite group is nilpotent.\\
\endgroup %%------------------------------------<<

\begingroup%%----------------------------------->>
\fontsize{9pt}{11pt}\selectfont
\textbf{\small EXAMPLE} \  
The $HP$-localization of every perfect group is trivial.  So, if \mG is perfect and if 
$H_2(G;\Z_P) \neq 0$, then the arrow $* \ra G$ is in $S_{HP}^{\perp\perp}$ but not in $S_{HP}$.\\
\endgroup %%------------------------------------<<

\begingroup%%----------------------------------->>
\fontsize{9pt}{11pt}\selectfont
\textbf{\small FACT} \  
The class of $HP$-homomorphisms admits a calculus of left fractions.\\
\endgroup %%------------------------------------<<

\index{Theorem: Kan Factorization Theorem}
\index{Kan Factorization Theorem}
\begingroup%%----------------------------------->>
\fontsize{9pt}{11pt}\selectfont
\textbf{\small KAN\footnote[2]{In: \textit{Algebra, Topology, and Category Theory}, A. Heller and M. Tierney (ed.), Academic Press (1976) 95-99.}
 FACTORIZATION THEOREM} \  
Let 
$
\begin{cases}
\ X\\
\ Y
\end{cases}
$
be pointed connected CW spaces, $f:X \ra Y$ a pointed continuous function.  
Assume: 
$f_*:H_q(X;\Z_P) \ra H_q(Y;\Z_P)$ is bijective for $1 \leq q < n$ and surjective for $q = n$ $-$then there exists a pointed connected CW space $X_f$ and pointed continuous functions 
$\phi_f:X \ra X_f$, $\psi_f:X_f \ra Y$ with $f = \psi_f \circx \phi_f$ such that 
$H_*(\phi_f):H_*(X;\Z_P) \ra H_*(X_f;\Z_P)$ is an isomorphism and $\psi_f:X_f \ra Y$ is an $n$-equivalence.
\\ \indent
[The case when $n = 1$ is handled by appropriately attaching 1-cells and 2-cells.  
In general, one iterates the following statement (which can be established by appropriately attaching $(n+1)$-cells and $(n+2)$-cells).
\\
\indent\indent (ST$_n$) \ Let $
\begin{cases}
\ X\\
\ Y
\end{cases}
$
be pointed connected CW spaces, $f:X \ra Y$ a pointed continuous function.  
Assume: $f$ is an $n$-equivalence and
$f_*:H_q(X;\Z_P) \ra H_q(Y;\Z_P)$ is bijective for $1 \leq q \leq n$ and surjective for $q = n + 1$ 
$-$then there exists a pointed connected CW space $X_f$ and pointed continuous functions 
$\phi_f:X \ra X_f$, $\psi_f:X_f \ra Y$ with $f = \psi_f \circx \phi_f$ such that 
$H_*(\phi_f):H_*(X;\Z_P) \ra H_*(X_f;\Z_P)$ is an isomorphism and $\psi_f:X_f \ra Y$ is an $(n + 1)$-equivalence.]\\
\endgroup %%------------------------------------<<

\label{9.67}
\label{10.4}
\label{10.5}
\begingroup%%----------------------------------->>
\fontsize{9pt}{11pt}\selectfont
Application: Let $f:G \ra K$ be a homomorphism of groups.  Assume: 
$f_*:H_1(G;\Z_P) \ra H_1(K;\Z_P)$ is surjective $-$then there exists a factorization 
$G \overset{\phi_f}{\lra} G_f \overset{\psi_f}{\lra} K$ of $f$ with $\phi_f$ an $HP$-homomorphism and $\psi_f$ surjective.
\vspi
[Recall that for any pointed connected space \mX, there is a surjection 
$H_2(X;\Z_P) \ra H_2(\pi_1(X);\Z_P)$ 
(cf. p. \pageref{8.37}).]\\
\endgroup %%------------------------------------<<

\begingroup%%----------------------------------->>
\fontsize{9pt}{11pt}\selectfont
\textbf{\small EXAMPLE} \  
Let $f:G \ra K$ be a homomorphism of $HP$-local groups $-$then $f$ is surjective iff 
$f_*:H_1(G;\Z_P) \ra H_1(K;\Z_P)$ is surjective.
\vspi
[To check sufficiency, note that the commutative diagram \ 
\begin{tikzcd}[sep=large]
{G} \ar{d}[swap]{\phi_f} \arrow[r,shift right=0.5,dash] \arrow[r,shift right=-0.5,dash] &{G} \ar{d}{f}\\
{G_f} \ar{r}[swap]{\psi_f} &{K}
\end{tikzcd}
\ 
has a filler $G_f \ra G$ rendering the triangles commutative.]\\
\endgroup %%------------------------------------<<

%%----------------------------------------------------------------------------------------------28
\begingroup%%----------------------------------->>
\fontsize{9pt}{11pt}\selectfont
\textbf{\small FACT} \  
Let $f:G \ra K$ be a homomorphism of $HP$-local groups $-$then $\im f$ is $HP$-local.\\
\endgroup %%------------------------------------<<

Let \mA be a ring with unit.  Fix a right $A$-module \mR $-$then a homomorphism $f:M \ra N$ of left $A$-modules is said to be an 
\un{$HR$-homomorphism}
\index{HR-homomorphism} provided that 
$R \otimes_A M \ra R \otimes_A N$ is an isomorphism and 
$\Tor_1^A(R,M) \ra \Tor_1^A(R,N)$ is an epimorphism.

\indent\indent ($HR$-Localization) \ Let $S_{HR} \subset \Mor \AbMOD$ be the class of $HR$-homomorphisms $-$then 
$S_{HR}^\perp$ is the object class of a reflective subcategory $\AbMOD_{HR}$ of \AbMOD.  
The reflector 
$
L_{HR}:
\begin{cases}
\ \AbMOD \ra \AbMOD_{HR}\\
\ M \ra M_{HR}
\end{cases}
$
is called 
\un{$HR$-localization}
\index{HR-localization} 
and the objects in $\AbMOD_{HR}$ are called 
\un{$HR$-local}
\index{HR-local ((left) $A$-modules)}
(left) $A$-modules.

[Each object in \AbMOD is $\kappa$-definite for some $\kappa$.  
Accordingly, due to the reflective subcategory theorem, 
one has only to find a set $S_0 \subset S_{HR}$ : $S_0^\perp = S_{HR}^\perp$, 
which can be done by using the colimit lemma.]\\

\begin{proposition} \ %20
$L_{HR}:\AbMOD \ra \AbMOD_{HR}$ is an additive functor.\\
\end{proposition}

Let \mG be a group, $A = \Z[G]$ and write \GbMOD in place of $\Z[G]\text{-}\bMOD$.  
Take $R = \Z$ (trivial $G$-action) $-$then 
a homomorphism $f:M \ra N$ of $G$-modules is an $H\Z$ homomorphism iff 
$f_*:H_0(G;M) \ra H_0(G;N)$ is bijective and 
$f_*:H_1(G;M) \ra H_1(G;N)$ is surjective.  The reflector 
$
L_{H\Z}:
\begin{cases}
\ \GbMOD \ra \GbMOD_{H\Z}\\
\ M \ra M_{H\Z}
\end{cases}
$
is called 
\un{$H\Z$-localization}
\index{H$\Z$-localization} 
and the objects in $\GbMOD_{H\Z}$ are called 
\un{$H\Z$-local}
\index{H$\Z$-local ((left) $G$-modules)}
(left) $G$-modules.  
Example: Every trivial $G$-module is $H\Z$-local.

\label{9.59}
[Note: \ The arrow of localization $l_{H\Z}:M \ra M_{H\Z}$ is an $H\Z$-homomorphism 
(cf. p. \pageref{8.38}), i.e., 
$l_{H\Z} \in S_{H\Z} \subset S_{H\Z}^{\perp\perp}$.]\\

\begin{proposition} \ %21
The $H\Z$-localization of any $M$ in $G$-\bMOD which is $P$-local as an abelian group is again $P$-local: 
$M = \Z_P \otimes M$ $\implies$ $M_{H\Z} = \Z_P \otimes M_{H\Z}$.
\end{proposition}

[This is because $L_{H\Z}$ is an additive functor (cf. Proposition 20).]\\

\textbf{\small SUBLEMMA} \quad 
Suppose that 
\begin{tikzcd}%[ sep=large]
{M} \ar{d}[swap]{f} \ar{r} &{N} \ar{d}{g}\\
{P} \ar{r} &{Q}
\end{tikzcd}
is a pushout square in $G$-\bMOD.  
Assume: $f$ is an $H\Z$-homomorphism $-$then $g$ is an $H\Z$-homomorphism.

[There is a \cd 
\begin{tikzcd}%[ sep=large]
{M} \ar{rd}[swap]{f} \ar{r}{\pi} &{\ov{M}}\ar{d} \ar{r} &{N} \ar{d}{g}\\
&{P} \ar{r} &{Q}
\end{tikzcd}
, where $\pi$ is surjective and the square is simultaneously a pullback and a pushout in $G$-\bMOD.  
Observing that the arrow $\ov{M} \ra P$ is an $H\Z$-homomorphism, consider the long exact sequence 
$H_1(G;\ov{M}) \ra$ 
%%----------------------------------------------------------------------------------------------29
$H_1(G;N) \oplus H_1(G;P) \ra$ 
$H_1(G;Q) \ra $ 
$H_0(G;\ov{M}) \ra$ 
$H_0(G;N) \oplus H_0(G;P) \ra$ 
$H_0(G;Q) \ra$ $0$.]\\

\label{9.66} %dmc this seems out of place

\textbf{\small LEMMA} \  
Let 
$0 \ra M^\prime \ra M \ra M\pp \ra 0$ be a short exact sequence of $G$-modules.  Assume: $M^\prime$ is $H\Z$-local $-$then 
in any commutative diagram 
\begin{tikzcd}%[sep=large]
{P} \ar{d}[swap]{f} \ar{r} &{M} \ar{d}\\
{Q} \ar{r} &{M\pp}
\end{tikzcd}
of $G$-modules, where $f:P \ra Q$ is an $H\Z$-homomorphism, there is a unique lifting 
\begin{tikzcd}%[sep=large]
{P} \ar{d}[swap]{f} \ar{r} &{M} \ar{d}\\
{Q} \ar[dashed]{ru}\ar{r} &{M\pp}
\end{tikzcd}
rendering the triangles commutative.

[Uniqueness is elementary, so we shall deal only with existence.  
Define \mN by the pushout square 
\quad
\begin{tikzcd}%[sep=large]
{P} \ar{d}[swap]{f} \ar{r} &{M} \ar{d}\\
{Q} \ar{r} &{N}
\end{tikzcd}
\quad
and \ display \ the \  data \  in a \ commutative diagram
$
\begin{tikzcd}%[sep=large]
{P} \ar{d}[swap]{f} \ar{r} &{M} \ar{d} \arrow[r,shift right=0.5,dash] \arrow[r,shift right=-0.5,dash] &{M} \ar{d}\\
{Q} \ar{r} &{N} \ar{r}[swap]{\pi} &{M\pp}
\end{tikzcd}
.
$
\ 
Put $N^\prime = \ker \pi$, define $\ov{N}$ by the pushout square
$ 
\begin{tikzcd}%[sep=large]
{N^\prime} \ar{d} \ar{r} &{N} \ar{d}\\
{N_{H\Z}^\prime} \ar{r} &{\ov{N}}
\end{tikzcd}
, 
$ 
and pass to 
\[
\begin{tikzcd}%[sep=large]
{0}  \ar{r} 
&{M^\prime} \ar{d} \ar{r}
&{M} \ar{d} \ar{r}
&{M\pp} \arrow[d,shift right=0.5,dash] \arrow[d,shift right=-0.5,dash]  \ar{r}
&{0}\\
{0} \ar{r} 
&{N^\prime} \ar{d} \ar{r}
&{N} \ar{d} \ar{r}
&{M\pp}\arrow[d,shift right=0.5,dash] \arrow[d,shift right=-0.5,dash]   \ar{r}
&0{}\\
{0} \ar{r} 
&{N_{H\Z}^\prime} \ar{r}
&{\ov{N}} \ar{r}
&{M\pp} \ar{r}
&{0}
\end{tikzcd}
.
\]
According to the sublemma, the arrows 
$M \ra N$, $N \ra \ov{N}$ are $H\Z$-homomorphisms, thus the composite 
$M^\prime \ra N^\prime \ra N_{H\Z}^\prime$ is an $H\Z$-homomorphism, hence is an isomorphism (since $M^\prime$ and 
$N_{H\Z}^\prime$ are $H\Z$-local).  Therefore the composite
$M \ra N \ra \ov{N}$ is an isomorphism.  Precompose its inverse with the arrow $N \ra \ov{N}$ to get a lifting 
$
\begin{tikzcd}%[sep=large]
{M} \ar{d} \arrow[r,shift right=0.5,dash] \arrow[r,shift right=-0.5,dash] &{M} \ar{d}\\
{N} \ar[dashed]{ru}\ar{r}[swap]{\pi} &{M\pp}
\end{tikzcd}
, 
$ 
which may then be precomposed with the arrow $Q \ra N$ to get a lifting
$
\begin{tikzcd}%[sep=large]
{P} \ar{d}[swap]{f} \ar{r} &{M} \ar{d}\\
{Q} \ar[dashed]{ru} \ar{r} &{M\pp}
\end{tikzcd}
,
$ 
as desired.]\\

\begin{proposition} \ 
Let 
$0 \ra M^\prime \ra M \ra M\pp \ra 0$ be a short exact sequence of $G$-modules.  Assume: $M^\prime$ and $M\pp$ are 
$H\Z$-local $-$then \mM is $H\Z$-local.\\
\end{proposition}

Application: Every nilpotent $G$-module is $H\Z$-local.

%%----------------------------------------------------------------------------------------------30
[Note: \ More generally, if \mM is a $G$-module such that for some $i$, 
$(I[G])^i\cdot M = (I[G])^{i+1} \cdot M$, then 
$M_{H\Z} \approx M/(I[G])^i \cdot M$.  
Proof: It follows from the exact sequence 
$H_1(G;M) \ra$ 
$H_1(G;M/(I[G])^i \cdot M) \ra$ 
$(I[G])^i \cdot M / (I[G])^{i+1} \cdot M \ra$ 
$H_0(G;M) \ra$ 
$H_0(G;M/(I[G])^i \cdot M) \ra 0$
that the arrow 
$M \ra M/(I[G])^i \cdot M$ is an $H\Z$-homomorphism.  
On the other hand, 
$M/(I[G])^i \cdot M$ is a nilpotent $G$-module.  
As for the realizability of the condition, recall that 
$G/[G,G] \approx I[G]/I[G]^2$, hence \mG perfect $\implies$ $I[G] = I[G]^2$ and 
$G/[G,G]$ divisible + torsion $\implies$ $I[G]^2 = I[G]^3 = \cdots$.]\\

\begingroup%%----------------------------------->>
\fontsize{9pt}{11pt}\selectfont
\textbf{\small FACT} \  
The class of $H\Z$-homomorphisms admits a calculus of left fractions.\\
\endgroup %%------------------------------------<<

\label{9.71} %dmc mnft

\textbf{\small LEMMA} \  
Let $f:M \ra N$ be a homomorphism of $G$-modules.  
Assume: $f_*:H_0(G;M) \ra H_0(G;N)$ is surjective $-$then there exists a factorization 
$M \overset{\phi_f}{\lra} M_f \overset{\psi_f}{\lra} N$ of $f$ with $\phi_f$ an $H\Z$-homomorphism and $\psi_f$ surjective.

[Choose a free $G$-module \mP and a surjection $\mu:M \oplus P \ra N$ such that $\restr{\mu}{M} = f$.  
Since the composite 
$H_0(G;\ker \mu) \ra H_0(G;M \oplus P) \ra H_0(G;P)$ is surjective and $H_0(G;P)$ is free abelian, one can find a free $G$-module 
\mQ and a homomorphism $\nu:Q \ra \ker \mu$ such that $H_0(G;Q) \approx H_0(G;P)$ through 
$Q \overset{\nu}{\ra} \ker \mu \ra M \oplus P \ra P$.  
Factor $f$ as 
$M \overset{\phi_f}{\lra} (M \oplus P)/\nu(Q) \overset{\psi_f}{\lra} N$, where $\phi_f$ is induced by the inclusion 
$M \ra M \oplus P$ and $\psi_f$ is induced by $\mu$.]\\

\begin{proposition} \ %23
Let $f:M \ra N$ be a homomorphism of $H\Z$-local $G$-modules $-$then $f$ is surjective iff $f_*:H_0(G;M) \ra H_0(G;N)$ is surjective.
\end{proposition}

[To check sufficiency, note that the commutative diagram 
\begin{tikzcd}%[sep=large]
{M} \ar{d}[swap]{\phi_f} \arrow[r,shift right=0.5,dash] \arrow[r,shift right=-0.5,dash] &{M} \ar{d}{f}\\
{M_f} \ar{r}[swap]{\psi_f} &{N}
\end{tikzcd}
has a filler $M_f \ra M$ rendering the triangles commutative.]\\

\begin{proposition} \ %24
Let $f:M \ra N$ be a homomorphism of $H\Z$-local $G$-modules $-$then $\im f$ is $H\Z$-local.
\end{proposition}

[Let $\ov{N} \supset f(M)$ be the largest $G$-submodule of \mN for which the induced map 
$H_0(G;$ $f(M)) \ra H_0(G;\ov{N})$ is surjective.  There is a commutative triangle 
\begin{tikzcd}%[sep=large]
{M} \ar{d}[swap]{\ov{f}} \ar{rd}{f} \\
{\ov{N}} \ar{r}[swap]{j} &{N}
\end{tikzcd}
and a factorization 
$M \overset{\phi_{\ov{f}}}{\lra} M_{\ov{f}} \overset{\psi_{\ov{f}}}{\lra} \ov{N}$ of 
$\ov{f}$ with $\phi_{\ov{f}}$ an $H\Z$-homomorphism and $\psi_{\ov{f}}$ surjective.  
Consider any lifting 
$M_{\ov{f}} \ra M$ of $j \circx \psi_{\ov{f}}$ to see that 
$\ov{N} = f(M)$.  But $\ov{N}$ is $H\Z$-local.]\\

\begin{proposition} \ %25
Let $f:M \ra N$ be a homomorphism of $H\Z$-local $G$-modules $-$then $\coker f$ is $H\Z$-local.
\end{proposition}

%%----------------------------------------------------------------------------------------------31
[Since $\im f$ is $H\Z$-local (cf. Proposition 24), one can assume that $f$ is injective, the claim thus being that 
$N/M$ is $H\Z$-local.  There is a commutative diagram
\[
\begin{tikzcd}%[sep=large]
{0}  \ar{r} 
&{M} \ar{d} \ar{r}
&{N} \arrow[d,shift right=0.5,dash] \arrow[d,shift right=-0.5,dash] \ar{r}
&{N/M} \ar{d} \ar{r}
&{0}\\
{0} \ar{r} 
&{K} \ar{r}
&{N} \ar{r}
&{(N/M)_{H\Z}} \ar{r}
&{0}
\end{tikzcd}
\]
of short exact sequences, where the kernel \mK is $H\Z$-local.  The arrow $M \ra K$ is obviously injective.  
That it is also surjective can be seen by comparing the exact sequence 
$H_1(G;N) \ra H_1(G;N/M) \ra H_0(G;M) \ra H_0(G;N) \ra H_0(G;N/M)$ from the first row with its analog from the second row and applying the five lemma: 
$H_0(G;M) \ra H_0(G;K)$ surjective $\implies$ $M \ra K$ surjective (cf. Proposition 23).  Conclusion: 
$N/M \approx (N/M)_{H\Z}$.]

[Note: \ A priori, cokernels in $G$-$\bMOD_{H\Z}$ are calculated first in $G$-\bMOD and then reflected back into 
$G$-$\bMOD_{H\Z}$.  The point of the proposition is that the second step is not needed.]\\

Application: $G$-$\bMOD_{H\Z}$ is an abelian category and the reflector 
$L_{H\Z}:G$-$\bMOD \ra$ $G$-$\bMOD_{H\Z}$ is right exact.\\

\label{9.69}
\label{9.70}
\begingroup%%----------------------------------->>
\fontsize{9pt}{11pt}\selectfont
\textbf{\small EXAMPLE} \  
Let \mM be an $H\Z$-local $G$-module $-$then $\forall \ n$, $\Z/n\Z \otimes M$ is $H\Z$-local.\\
\endgroup %%------------------------------------<<

\begingroup%%----------------------------------->>
\fontsize{9pt}{11pt}\selectfont
\textbf{\small EXAMPLE} \  
Let \bM be a tower in $G\text{-}\bMOD_{H\Z}$ $-$then $\lim \bM$ and $\lim^1 \bM$ are $H\Z$-local 
(cf. p. \pageref{8.39}).\\
\endgroup %%------------------------------------<<

\begingroup%%----------------------------------->>
\fontsize{9pt}{11pt}\selectfont
\textbf{\small FACT} \  
Let \bM be a tower in $G\text{-}\bMOD_{H\Z}$ Assume: \mG is finitely generated $-$then $\lim^1 \bM = 0$ iff 
$\lim^1 H_0(G;\bM) = 0$.
\vspi
[Here, $H_0(G;\bM)$ stands for the tower determined by the arrows 
$H_0(G;M_{n+1}) \ra H_0(G;M_n)$.  Use Proposition 23 and the fact that \mG finitely generated $\implies$ 
$H_0\bigl(G;\ds\prod\limits_n M_n \bigr) \approx \ds\prod\limits_n H_0(G;M_n)$ 
(Brown\footnote[2]{\textit{Comment. Math. Helv.} \textbf{50} (1975), 129-135; see also 
Strebel, \textit{Math. Zeit.} \textbf{151} (1976), 263-275.}).]\\
\endgroup %%------------------------------------<<

\begin{proposition} \ %26
Let $G \ra \pi$ be a homomorphism of groups $-$then every $H\Z$-local $\pi$-module is an $H\Z$-local $G$-module.
\end{proposition}

[The forgetful functor $\pi$-$\bMOD \ra G$-$\bMOD$ has a left adjoint $G$-$\bMOD \ra \pi$-$\bMOD$ that sends \mM to $\Z[\pi] \otimes_{\Z[G]} M$.  
Thanks to the change of rings spectral sequence, the homomorphism 
$H_i(G;M) \ra H_i(\pi;\Z[\pi] \otimes_{\Z[G]} M)$ is bijective for $i = 0$ and surjective for $i = 1$.  
Therefore an 
$H\Z$-homomorphism of $G$-modules goes over to an $H\Z$-homomorphism
%%----------------------------------------------------------------------------------------------32
of $\pi$-modules.  Suppose now that \mP is an $H\Z$-local $\pi$-module.  Let $M \ra N$ be an $H\Z$-homomorphism of 
$G$-modules $-$then the bijectivity of the arrow 
$\Hom(N,P) \ra \Hom(M,P)$ follows from the bijectivity of the arrow 
$\Hom(\Z[\pi] \otimes_{\Z[G]} N, P) \ra \Hom(\Z[\pi] \otimes_{\Z[G]} M,P)$.]\\

\label{9.75}
\begingroup%%----------------------------------->>
\fontsize{9pt}{11pt}\selectfont
\textbf{\small EXAMPLE} \  
Let \mM be an $H\Z$-local $G_{HP}$-module $-$then \mM is an $H\Z$-local $G$-module.\\
\endgroup %%------------------------------------<<

Although one can consider $HA$-localization for an arbitrary abelian group \mA, apart from $A = \Z_P$ the other case of topological significance is when $A = \F_p$.  The general aspects of the $H\F_p$-theory are similar to those of the $HP$-theory.  
\label{10.1}
For instance, the analog of Proposition 19 says that if 
$1 \ra G^\prime \ra G \ra G\pp \ra 1$ is a central extension of groups with $G^\prime$ an $\F_p$-module and $G\pp$ 
$H\F_p$-local, then \mG is $H\F_p$-local.

[Note: \ An abelian group is a $\Z_P$-module iff it is $P$-local iff it is $HP$-local.  To perfect the analogy, one can relax the assumption on $G^\prime$ and suppose only that $G^\prime$ is $H\F_p$-local (cf. Proposition 33).]\\

\begin{proposition} \ %27
Every $H\F_p$-local group is $p$-local.\\
\end{proposition}

\begingroup%%----------------------------------->>
\fontsize{9pt}{11pt}\selectfont
\textbf{\small EXAMPLE} \  
Let \mG be a finite group $-$then $G_{H\F_p} \approx G_p$.\\
\endgroup %%------------------------------------<<

\label{9.72}
\begingroup%%----------------------------------->>
\fontsize{9pt}{11pt}\selectfont
The Kan factorization theorem remains valid if $\Z_P$ is replaced by $\F_p$.  Therefore a homomorphism 
$f:G \ra K$ of $H\F_p$-local groups is surjective iff $f_*:H_1(G;\F_p) \ra H_1(K;\F_p)$ is surjective.\\
\endgroup %%------------------------------------<<

The class of $H\F_p$-local abelian groups turns out to be the same as the class of $p$-cotorsion abelian groups (cf. Proposition 30).  It will therefore be convenient to review the theory of the latter starting with the global situation.

An abelian group $G$ is said to be 
\un{cotorsion}
\index{cotorsion (abelian group)} 
if 
$\Hom(\Q,G) = 0$ $\&$ 
$\Ext(\Q,G) = 0$.  Taking into account the exact sequence 
$\Hom(\Q,G) \ra$ 
$\Hom(\Z,G) \ra$ 
$\Ext(\Q/\Z,G) \ra$ 
$\Ext(\Q,G)$ 
and making the identification $G \approx \Hom(\Z,G)$, it follows that \mG is cotorsion iff the arrow 
$G \ra \Ext(\Q/\Z,G)$ is an isomorphism.

[Note: \ One motivation for the terminology is that if 
$0 \ra A \ra B \ra C \ra 0$ is a short exact sequence of abelian groups, then the sequence 
$0 \ra$ 
$\Hom(K,A) \ra$ 
$\Hom(K,B) \ra$ 
$\Hom(K,C) \ra 0$ is exact for all torsion groups \mK iff the sequence 
$0 \ra$ 
$\Hom(C,L) \ra$
$\Hom(B,L) \ra$
$\Hom(A,L) \ra 0$ is exact for all cotorsion groups $L$.]\\

\begingroup%%----------------------------------->>
\fontsize{9pt}{11pt}\selectfont
Let $0 \ra A \ra B \ra C \ra 0$ be a short exact sequence of abelian groups $-$then
$0 \ra$ 
$\Hom(K,A) \ra$ 
$\Hom(K,B) \ra$ 
$\Hom(K,C) \ra 0$ is exact $\forall$ torsion \mK iff 
$0 \ra$ 
$\Hom(K,A) \ra$ 
$\Hom(K,B) \ra$ 
$\Hom(K,C) \ra 0$ is exact $\forall$ finite cyclic \mK iff 
$0 \ra$
$A \ra$ 
$B \ra$ 
$C \ra 0$ is pure short exact iff 
$0 \ra$ 
$\Hom(C,L) \ra$
$\Hom(B,L) \ra$ 
%%----------------------------------------------------------------------------------------------33
$\Hom(A,L) \ra 0$ is exact $\forall$ finite cyclic \mL iff 
$0 \ra$ 
$\Hom(C,L) \ra$ 
$\Hom(B,L) \ra$ 
$\Hom(A,L) \ra 0$ is exact $\forall$ cotorsion \mL.\\ 
\endgroup %%------------------------------------<<

\textbf{\small LEMMA} \  
For any abelian group \mG, $\Ext(\Q/\Z,G)$ is cotorsion.

[Given $A, B, C$ in \bAB, there are isomorphisms
\begin{align*}
\Ext(A,\Ext(B,C)) &\approx \Ext(\Tor(A,B),C),\\
\Ext(A,\Hom(B,C)) \oplus \Hom(A,\Ext(B,C))  &\approx \Ext(A \otimes B,C) \oplus \Hom(\Tor(A,B),C).]
\end{align*}
\vspace{0.05cm}

\textbf{\small LEMMA} \  
For any abelian group \mG, $\Ext(\Q/\Z,\Ext(\Q/\Z,G)) \approx \Ext(\Q/\Z,G)$.\\

Consequently, the full subcategory of \bAB whose objects are the cotorsion groups is a reflective subcategory of \bAB, the arrow of reflection being $G \ra \Ext(\Q/\Z,G)$.

[Note: \ By comparison, the full subcategory of \bAB whose objects are the torsion groups is a coreflective subcategory of \bAB, the arrow of coreflection being $\Tor(\Q/\Z,G) \ra G$.]\\

\begingroup%%----------------------------------->>
\fontsize{9pt}{11pt}\selectfont
\textbf{\small EXAMPLE} \  
$\Z/n\Z$ is cotorsion but $\Z$ is not cotorsion.\\
\endgroup %%------------------------------------<<

A cotorsion group \mG is said to be 
\un{adjusted}
\index{adjusted (cotorsion group)} 
if \mG has no torsion free direct summand or, equivalently, if $G/G_\tor$ is divisible.\\

\index{Cotorsion Structure Lemma}
\textbf{\small COTORSION STRUCTURE LEMMA} \  
Suppose that \mG is cotorsion $-$then there is a split short exact sequence 
$0 \ra K \ra G \ra L \ra 0$, where 
$K \approx \Ext(\Q/\Z,G_\tor)$ is adjusted cotorsion and 
$L \approx \Ext(\Q/\Z,G/G_\tor)$ is torsion free cotorsion.

[Note: \ In the opposite direction, recall that every abelian group is split by its maximal divisible subgroup and the associated quotient is reduced.]\\

\index{Theorem: Harrison's First Theorem}
\index{Harrison's First Theorem}
\textbf{\small HARRISON'S\footnote[2]{\textit{Ann. of Math.} \textbf{69} (1959), 366-391.}
 FIRST THEOREM} \  
Let \bC be the full subcategory of \bAB whose objects are the torsion free cotorsion groups; let \bD be the full subcategory of \bAB whose objects are the divisible torsion groups.  Define 
$\Phi:\bC \ra \bD$ by 
$\Phi G = \Q/\Z \otimes G$; define
$\Psi:\bD \ra \bC$ by 
$\Psi G = \Hom(\Q/\Z,G)$ $-$then the pair $(\Phi,\Psi)$ is an adjoint equivalence of categories.\\

\index{Theorem: Harrison's Second Theorem}
\index{Harrison's Second Theorem}
\textbf{\small HARRISON'S$^
%\footnote[2]{} %% symbol no ref since same as above
%dmc missing ref
\dagger$
 SECOND THEOREM} \  
 Let \bC be the full subcategory of \bAB whose objects are the adjusted cotorsion groups; let \bD be the full subcategory of \bAB
 whose
%%----------------------------------------------------------------------------------------------34
objects are the reduced torsion groups.  Define 
$\Phi:\bC \ra \bD$ by 
$\Phi G = \Tor(\Q/\Z,G)$; define
$\Psi:\bD \ra \bC$ by 
$\Psi G = \Ext(\Q/\Z,G)$ $-$then the pair $(\Phi,\Psi)$ is an adjoint equivalence of categories.\\

An abelian group \mG is said to be 
\un{$p$-cotorsion}
\index{$p$-cotorsion (abelian group)} 
if 
$\Hom\big(\Z\left[\ds\frac{1}{p}\right],G\big) = 0$ $\&$
$\Ext\bigl(\Z\left[\ds\frac{1}{p}\right],G\bigr)$ $= 0$.  
Taking into account the exact sequence 
$\Hom\bigl(\Z\left[\ds\frac{1}{p}\right],G\bigr)$ $\ra$ 
$\Hom(\Z,G) \ra$ 
$\Ext(\Z/p^\infty\Z,$ $G) \ra$ 
$\Ext\bigl(\Z\left[\ds\frac{1}{p}\right],G\bigr)$ 
and making the identification 
$G \approx \Hom(\Z,G)$, it follows that \mG is $p$-cotorsion iff the arrow 
$G \ra \Ext(\Z/p^\infty\Z,G)$ is an isomorphism.  
Example: $\forall \ n$, $\Z/p^n\Z$ is $p$-cotorsion.

[Note: \ The full subcategory of \bAB whose objects are the $p$-cotorsion groups is a reflective subcategory of \bAB with arrow of reflection $G \ra \Ext(\Z/p^\infty\Z,G)$ and there are evident variants of Harrison's first and second theorems.]\\

\begingroup%%----------------------------------->>
\fontsize{9pt}{11pt}\selectfont
\textbf{\small EXAMPLE} \  
If $G = \widehat{\Z}_p$, the $p$-adic integers, then 
$\widehat{\Z}_p \approx \Ext(\Z/p^\infty \Z,\widehat{\Z}_p)$, hence $\widehat{\Z}_p$ is $p$-cotorsion.
\vspi
[Note: \ A subgroup of $\widehat{\Z}_p$ is $p$-cotorsion iff it is an ideal.]\\
\endgroup %%------------------------------------<<

\begingroup%%----------------------------------->>
\fontsize{9pt}{11pt}\selectfont
\textbf{\small EXAMPLE} \  
The following abelian groups are not $p$-cotorsion: 
$\Z/p^\infty \Z,$ 
$\ds\bigoplus\limits_n \Z/p^n \Z$, 
$\widehat{\Z}_p \otimes \widehat{\Z}_p$.\\
\endgroup %%------------------------------------<<

\label{9.73}
\label{9.74}
\begingroup%%----------------------------------->>
\fontsize{9pt}{11pt}\selectfont
\textbf{\small EXAMPLE} \  
For any abelian group \mG, $\Hom(\Z/p^\infty\Z,G)$ is $p$-cotorsion.  In fact, 
$\Hom\bigl(\Z\left[\ds\frac{1}{p}\right],$ $\Hom(\Z/p^\infty\Z,G)\bigr) \approx$
$\Hom\bigl(\Z\left[\ds\frac{1}{p}\right] \otimes \Z/p^\infty\Z,G\bigr) \approx$ 
$\Hom(0,G) = 0$ and 
$\Ext\bigl(\Z\left[\ds\frac{1}{p}\right],\Hom(\Z/p^\infty\Z,G)\bigr) \approx$
$\Ext\bigl(\Tor\bigl(\Z\left[\ds\frac{1}{p}\right],\Z/p^\infty\Z\bigr),G\bigr) \approx$
$\Ext(0,G) = 0$.\\
\endgroup %%------------------------------------<<

\begingroup%%----------------------------------->>
\fontsize{9pt}{11pt}\selectfont
\textbf{\small FACT} \  
Let \mG be a group and let \mM be a $G$-module.  
Assume: \mM is $H\Z$-local $-$then $\Ext(\Z/p^\infty\Z,M)$ is $H\Z$-local.
\label{10.7}
\vspi
[The arrow 
$\Ext(\Z/p^\infty\Z,M) \ra \lim\Ext(\Z/p^n\Z,M)$ is surjective and its kernel can be identified with 
$\lim^1\Ext(\Z/p^n\Z,M)$ 
(Weibel\footnote[2]{\textit{An Intoduction to Homological Algebra}, Cambridge University Press (1994), 85; 
see also Jensen, \textit{SLN} \textbf{254} (1972), 35-37.}), i.e., there is a short exact sequence 
$0 \ra$ 
$\lim^1\Hom(\Z/p^n\Z,M) \ra$ 
$\Ext(\Z/p^\infty\Z,M) \ra$ 
$\lim\Ext(\Z/p^n\Z,M)$.  Since 
$\Ext(\Z/p^n\Z,M) \approx M/p^n M$ and $M/p^n M$ is $H\Z$-local (cf. Proposition 25), 
$\lim\Ext(\Z/p^n\Z,M)$ must be $H\Z$-local too ($G\text{-}\bMOD_{H\Z}$ is limit closed).  
Similar remarks imply that 
$\lim^1\Hom(\Z/p^n\Z,M)$  is $H\Z$-local (it is a cokernel 
(cf. p. \pageref{8.40})). Now quote Proposition 22.]\\
\endgroup %%------------------------------------<<

\begingroup%%----------------------------------->>
\fontsize{9pt}{11pt}\selectfont
\textbf{\small FACT} \  
For any abelian group \mG, the arrow of reflection $G \ra \Ext(\Z/p^\infty\Z,G)$ induces an isomorphism 
$\F_p \otimes G \ra \F_p \otimes \Ext(\Z/p^\infty\Z,G)$ and an epimorphism 
$\Tor(\F_p,G) \ra \Tor(\F_p,\Ext(\Z/p^\infty\Z,G))$.

%%----------------------------------------------------------------------------------------------35
[To check the first assertion,  observe that 
$\F_p \otimes G$
\hsx $\approx$ \hsx
$\Ext(\F_p,G)$ 
\hsx $\approx$ \hsx
$\Ext(\Tor(\F_p,\Z/p^\infty\Z),G)$ 
\hsx $\approx$ \hsx
$\Ext(\F_p,\Ext(\Z/p^\infty\Z,G))$ 
$\approx$ 
$\F_p \otimes \Ext(\Z/p^\infty\Z,G)$.]\\
\endgroup %%------------------------------------<<

Notation: Given an abelian group \mG, $\xdiv G$ is the maximal divisible subgroup of \mG and $\divx_p G$ is the maximal $p$-divisible subgroup of \mG.

[Note: \ The kernel of the arrow of reflection $G \ra \Ext(\Z/p^\infty\Z,G)$ is $\divx_p G$.]\\

\begin{proposition} \ %28
Suppose that \mG is cotorsion $-$then $G \approx \prod\limits_p G_p$, 
where $G_p = \ds\bigcap\limits_{q \neq p} \divx_q G$ is the maximal $p$-cotorsion subgroup of \mG.
\end{proposition}

[The point here is that $\Ext(\Q/\Z,G) \approx \prod\limits_p \Ext(\Z/p^\infty\Z,G)$.]

[Note: \ This result is the analog for a cotorsion group of the primary decomposition of a torsion group.]\\

\label{9.82}
\textbf{\small LEMMA} \  
If \mA and \mG are abelian groups with \mG $p$-cotorsion, then 
(i) $A \otimes \F_p = 0$ $\implies$ $\Hom(A,G) = 0$ and 
(ii) $\Tor(A,\F_p) = 0$ $\implies$ $\Ext(A,G) = 0$.

[To check the second assertion, observe that 
$\Ext(A,G) \approx$ 
$\Ext(A,\Ext(\Z/p^\infty\Z,G)) \approx$ 
$\Ext(\Tor(A,\Z/p^\infty\Z),G) \approx$ 
$\Ext(0,G) = 0$.]\\

\begin{proposition} \ %29
Let 
$
\begin{cases}
\ X\\
\ Y
\end{cases}
$
be path connected topological spaces, $f:X \ra Y$ a continuous function $-$then 
$f_*:H_*(X;\F_p) \ra H_*(Y;\F_p)$ is an isomorphism iff $f^*:H^*(Y;G) \ra H^*(X;G)$ is an isomorphism for all $p$-cotorsion abelian groups \mG.
\end{proposition}

[By passing to the mapping cylinder, one can assume that $f$ is an inclusion.  
If $\forall \ n \geq 1$, 
$H_n(Y,X;\F_p) = 0$, then $\forall \ n \geq 1$, 
$H_n(Y,X) \otimes \F_p = 0$ and 
$\Tor(H_n(Y,X),\F_p) = 0$.  
So, from the lemma, for any $p$-cotorsion \mG, 
$\Hom(H_n(Y,X),G) = 0$ and 
$\Ext(H_n(Y,X),G) = 0$ $\forall \ n \geq 1$, thus $H^n(Y,X;G) = 0$ $\forall \ n \geq 1$.  
To reverse the argument, specialize and take $G = \F_p$.]\\

In the context of $HR$-localization, take $A = \Z$ and $R = \F_p$ $-$then the object class of the corresponding reflective subcategory of $\Z$-$\bMOD \approx \bAB$ is the class of $p$-cotorsion groups.\\

\begin{proposition} \ %30
Let \mG be an abelian group $-$then \mG is $H\F_p$-local iff \mG is $p$-cotorsion.
\end{proposition}

[Let $S_1 \subset \Mor \bAB$ be the class of homomorphisms $f:A \ra B$ such that $A \otimes \F_p \ra B \otimes \F_p$ 
is an isomorphism and $\Tor(A,\F_p) \ra \Tor(B,\F_p)$ is an epimorphism (thus $S_1^\perp$ is the class of $p$-cotorsion groups) and let $S_2 \subset \Mor \bAB$ be the class of homomorphisms $f:A \ra B$ such that  
$f_*:H_1(A;\F_p) \ra H_1(B;\F_p)$ is bijective and  
$f_*:H_2(A;\F_p) \ra H_2(B;\F_p)$ is
%%----------------------------------------------------------------------------------------------36
surjective (thus $S_2^\perp$ is the class of abelian $H\F_p$-local groups) (cf. infra).  
Claim: $S_1 = S_2$.  
For, in either case, $A/pA \approx B/pB$.  This said, consider the commutative diagram
\[
\begin{tikzcd}%[sep=large]
{0}  \ar{r} 
&{H_2(A) \otimes \F_p} \ar{d} \ar{r}
&{H_2(A;\F_p)} \ar{d} \ar{r}
&{\Tor(A,\F_p)} \ar{d}\ar{r}
&{0}\\
{0} \ar{r} 
&{H_2(B) \otimes \F_p} \ar{r}
&{H_2(B;\F_p)} \ar{r}
&{\Tor(B,\F_p)} \ar{r}
&{0}
\end{tikzcd}
\]
of short exact sequences.  Since
$
\begin{cases}
\ H_2(A) \otimes \F_p \approx \wedge^2 (A/pA)\\
\ H_2(B) \otimes \F_p \approx \wedge^2 (B/pB)
\end{cases}
$
(Brown\footnote[2]{\textit{Cohomology of Groups}, Springer Verlag (1982), 126.}), 
the five lemma implies that if $\Tor(A,\F_p) \ra \Tor(B,\F_p)$ is an epimorphism, then 
$f_*:H_2(A;\F_p) \ra H_2(B;\F_p)$ is surjective.  The converse is trivial.]\\

\begingroup%%----------------------------------->>
\fontsize{9pt}{11pt}\selectfont
The reflective subcategory theorem is applicable to \bAB, so one can define the notion of ``abelian $H\F_p$-local group'' internally.  That this is the same as ``abelian $+ H\F_p$-local'' is a consequence of the following lemma.\\
\endgroup %%------------------------------------<<

\begingroup%%----------------------------------->>
\fontsize{9pt}{11pt}\selectfont
\textbf{\small LEMMA} \  
An $H\F_p$-homomorphism $G \ra K$ of groups induces an $H\F_p$-homomorphism 
$G/[G,G] \ra K/[K,K]$ of abelian groups.\\
\endgroup %%------------------------------------<<

Given a group \mG, let $\rho_p:G^\omega \ra G^\omega$ be the function defined by 
$\rho_p(g_0,g_1, \ldots) =$ $(g_0g_1^{-p},g_1g_2^{-p},\ldots$).\\

\begin{proposition} \ %31
Suppose that \mG is abelian $-$then $\rho_p$ is a homomorphism and 
$\ker \rho_p \approx \lim \bG_p \approx \Hom\bigl(\Z\left[\ds\frac{1}{p}\right],G)\bigr)$, 
$\coker \rho_p \approx \lim^1 \bG_p \approx \Ext\bigl(\Z\left[\ds\frac{1}{p}\right],G\bigr)$, 
where $\bG_p$ is the tower 
$\cdots \la G \overset{p}{\la} G \la \cdots$.
\end{proposition}

[Representing $\Z\left[\ds\frac{1}{p}\right]$ as a colimit $\cdots \ra \Z \overset{p}{\ra} \Z  \ra \cdots$ gives 
$\lim\bG_p \approx \Hom\bigl(\Z\left[\ds\frac{1}{p}\right],G\bigr)$ and, from the short exact sequence 
$0 \ra$ 
$\lim^1 \Hom(\Z,G) \ra$
$\Ext\bigl(\Z\left[\ds\frac{1}{p}\right],G\bigr) \ra \lim \Ext(\Z,G) \ra 0$ 
(Weibel\footnote[3]{\textit{An Introduction to Homological Algebra}, Cambridge University Press (1994), 85; 
see also Jensen, \textit{SLN} \textbf{254} (1972), 35-37.}), 
one has $\lim^1 \bG_p \approx \Ext\bigl(\Z\left[\ds\frac{1}{p}\right],G\bigr)$.]\\

Application:  An abelian group \mG is $p$-cotorsion ($= H\F_p$-local) iff $\lim \bG_p = 0$ $\&$ $\lim^1 \bG_p = 0$, i.e., 
iff $\rho_p$ is bijective.\\

%%----------------------------------------------------------------------------------------------37
Let \mG be a group $-$then \mG is said to be 
\un{$p$-cotorsion}
\index{p-cotorsion (group)} 
provided that 
$\rho_p$ is bijective.  
Claim: The full subcategory of \bGR whose objects are the $p$-cotorsion groups is a reflective subcategory of \bGR.  To see this, let 
$F_\omega$ be the free group on generators $x_0,x_1, \ldots$, define a homomorphism 
$f:F_\omega \ra F_\omega$ by $f(x_i) = x_ix_{i+1}^{-p}$ and consider $f^\perp$ (reflective subcategory theorem).\\

\label{8.42}
\begingroup%%----------------------------------->>
\fontsize{9pt}{11pt}\selectfont
\textbf{\small FACT} \  
Suppose that \mG is $p$-cotorsion $-$then $\Cen G$ is $p$-cotorsion.\\
\endgroup %%------------------------------------<<

\begin{proposition} \ 
Every $H\F_p$-local group is $p$-cotorsion.
\end{proposition}

[It is enough to prove that $f:F_\omega \ra F_\omega$ is an $H\F_p$-homomorphism.  But 
$f_*:H_1(F_\omega;\F_p) \ra H_1(F_\omega;\F_p)$ is the identity 
$\omega \cdot \F_p \ra \omega \cdot \F_p$ and 
$H_2(F_\omega;\F_p) \approx$ 
$H_2(F_\omega) \otimes \F_p \oplus \Tor(H_1(F_\omega),\F_p)$ vanishes.]\\

\label{9.78} %dmc mnft
The abelian $p$-cotorsion theory has been extended to \bNIL by 
Huber-Warfield\footnote[2]{\textit{J. Algebra} \textbf{74} (1982), 402-442.}.  
Thus the full subcategory of \bNIL whose objects are the $p$-cotorsion groups is a reflective subcategory of \bNIL.  It is traditional to denote the arrow of reflection by 
$G \ra$ 
$\Ext(\Z/p^\infty \Z,G)$ even though the ``Ext'' has no a priori connection with extensions of \mG by $\Z/p^\infty \Z$.  
One reason for this is that each short exact sequence 
$1 \ra$
$G^\prime \ra$ 
$G \ra$ 
$G\pp \ra 1$ of nilpotent groups gives rise to an exact sequence 
$0 \ra$ 
$\Hom(\Z/p^\infty \Z,G^\prime) \ra$
$\Hom(\Z/p^\infty \Z,G) \ra$ 
$\Hom(\Z/p^\infty \Z,G\pp) \ra$ 
$\Ext(\Z/p^\infty \Z,G^\prime) \ra$ 
$\Ext(\Z/p^\infty \Z,G) \ra$ 
$\Ext(\Z/p^\infty \Z,G\pp) \ra 0$.

[Note: \ It is reasonable to conjecture that the $p$-cotorsion reflector in \bGR extends the $p$-cotorsion reflector in \bNIL but I know of no proof.]

The $p$-cotorsion reflector in \bNIL respects $\bNIL^d$ : $\nil\ \Ext(\Z/p^\infty \Z,G) \leq \nil G$, 
hence its restriction to 
\bAB ``is'' the $p$-cotorsion reflector in \bAB.

Notation: Given a nilpotent group \mG, $\xdiv G$ is the maximal divisible subgroup of \mG and div$_p$ $G$ is the maximal 
$p$-divisible subgroup of \mG.

[Note: \ The kernel of the arrow of reflection $G \ra \Ext(\Z/p^\infty \Z,G)$ is div$_p$ $G$.]\\

\textbf{\small LEMMA} \  
For any nilpotent group \mG, $\Hom(\Z/p^\infty\Z,G)$ is a torsion free $p$-cotorsion abelian group.

[Let $G_\tor(p)$ be the maximal $p$-torsion subgroup of \mG $-$then $\xdiv G_\tor (p)$ is abelian and the range of every homomorphism $f:\Z/ p^\infty\Z \ra G$ is contained in $\xdiv G_\tor (p)$.]

[Note: \ Therefore \mG $p$-cotorsion $\implies$ $\Hom(\Z/ p^\infty\Z,G) = 0$.]\\

\label{9.77}
\begingroup%%----------------------------------->>
\fontsize{9pt}{11pt}\selectfont
\textbf{\small FACT} \  
Let \mG be a nilpotent group $-$then the arrow $g \ra g^p$ is bijective iff 
$\Hom(\Z/p^\infty\Z,G) = 0$ $\&$ 
$\Ext(\Z/p^\infty\Z,G) = 0$ or still, iff $\forall \ n > 0$, 
$H_n(G;\F_p) = 0$.\\
\endgroup %%------------------------------------<<

%%----------------------------------------------------------------------------------------------38
\begingroup%%----------------------------------->>
\fontsize{9pt}{11pt}\selectfont
\textbf{\small EXAMPLE} \  
There is a short exact sequence 
$0 \ra$ 
$\Z \ra$ 
$\widehat{\Z}_p \ra$ 
$\widehat{\Z}_p/\Z \ra 0$ and $\widehat{\Z}_p/\Z$ is uniquely $p$-divisible, hence 
$H_*(\Z,\F_p) \approx H_*(\widehat{\Z}_p;\F_p)$ 
(cf. p. \pageref{8.41}).\\
\endgroup %%------------------------------------<<

\begingroup%%----------------------------------->>
\fontsize{9pt}{11pt}\selectfont
\textbf{\small FACT} \  
Let 
$1 \ra$
$G^\prime \ra$ 
$G \ra$ 
$G\pp \ra 1$ be a short exact sequence of nilpotent groups.  Assume: Two of the groups are $p$-cotorsion $-$then so is the third.\\
\endgroup %%------------------------------------<<

\begingroup%%----------------------------------->>
\fontsize{9pt}{11pt}\selectfont
\textbf{\small EXAMPLE} \  
Suppose that \mG is nilpotent and $p$-cotorsion $-$then $G/\Cen G$ is $p$-cotorsion.
\vspi
[$\Cen G$ is necessarily $p$-cotorsion 
(cf. p. \pageref{8.42}).]\\
\endgroup %%------------------------------------<<

\label{8.36}
\textbf{\small LEMMA} \  
Let 
$1 \ra$
$G^\prime \ra$ 
$G \ra$ 
$G\pp \ra 1$ be a central extension of groups.  Assume: $G^\prime$ is $H\F_p$-local $-$then for any commutative diagram 
\begin{tikzcd}%[sep=large]
{K} \ar{d}[swap]{f} \ar{r} &{G} \ar{d}\\
{L} \ar{r} &{G\pp}
\end{tikzcd}
of groups, where $f:K \ra L$ is an $H\F_p$-homomorphism, there is a unique lifting 
\begin{tikzcd}%[sep=large]
{K} \ar{d}[swap]{f} \ar{r} &{G} \ar{d}\\
{L} \ar[dashed]{ru}\ar{r} &{G\pp}
\end{tikzcd}
rendering the triangles commutative.

[Put 
$
\begin{cases}
\ X = K(K,1)\\
\ Y = K(L,1)
\end{cases}
$
and consider the diagram 
$
\begin{tikzcd}%[sep=large]
{X} \ar{d}[swap]{f} \ar{r} &{K(G,1)} \ar{d}\\
{Y} \ar[dashed]{ru} \ar{r} &{K(G\pp,1)}
\end{tikzcd}
.  
$ 
Supposing, as we may, that $f$ is an inclusion, the obstruction to lifting lies in $H^2(Y,X;G^\prime)$.  
Claim: $H^2(Y,X;G^\prime)  = 0$.  
To verify this, look at the short exact sequence 
$0 \ra$ 
$\Ext(H_1(Y,X),G^\prime) \ra$ 
$H^2(Y,X;G^\prime) \ra$ 
$\Hom(H_2(Y,X),G^\prime) \ra 0$.  Since 
$f_*:H_1(X;\F_p) \ra H_1(Y;\F_p)$ is bijective and 
$f_*:H_2(X;\F_p) \ra H_2(Y;\F_p)$ is surjective, 
$H_2(Y,X) \otimes \F_p = 0$ and 
$\Tor(H_1(Y,X),\F_p) = 0$.  
But $G^\prime$ is $H\F_p$-local or still, $p$-cotorsion (cf. Proposition 30), thus 
$\Hom(H_2(Y,X),G^\prime) =0$ and 
$\Ext(H_1(Y,X),G^\prime) = 0$ \ (see \hsx the \hsx lemma preceding Proposition 29).  
Therefore 
$H^2(Y,X;G^\prime) = 0$ and the lifting exists.  
As for its uniqueness, of necessity 
$H_1(Y,X;\F_p) = 0$, i.e., 
$H_1(Y,X) \otimes \F_P = 0$, thus 
$H^1(Y,X;G^\prime) \approx$ 
$\Hom(H_1(Y,X),G^\prime) = 0$.]\\

\begin{proposition} \ %33
Let $1 \ra G^\prime \ra G \ra G\pp \ra 1$ be a central extension of groups.  
Assume: $G^\prime$ is $H\F_p$-local and 
$G\pp$ is $H\F_p$-local $-$then \mG is $H\F_p$-local.
\end{proposition}

[The proof is the same as that of Proposition 19.]\\

Application: If \mG is nilpotent and $p$-cotorsion, then \mG is $H\F_p$-local.

[In fact, $\Cen G$ and $G /\Cen G$ are $p$-cotorsion, so one can proceed by induction.]\\

\begin{proposition} \ %34
Let \mG be a $p$-cotorsion nilpotent group $-$then there exists a central series 
$G = C^0(G) \supset C^1(G) \supset \cdots$ having the same length as the descending central series of \mG such that 
$\forall \ i$, $C^i(G)/C^{i+1}(G)$ is a $p$-cotorsion abelian group.
\end{proposition}

%%----------------------------------------------------------------------------------------------39
[Define $C^i(G)$ to be the kernel of the composite 
 $G \ra G/\Gamma^i(G) \ra$ 
 $\Ext(\Z/p^\infty\Z,$ $G/\Gamma^i(G))$.]
 
[Note: \ Here is a variant.  Let \mG be a group.  
 Let \mM be a nilpotent $G$-module, $\chi:G \ra \Aut M$ the associated homomorphism.  
 Assume: \mM is $p$-cotorsion $-$then there exists a finite filtration 
 $M = C_\chi^0(M) \supset C_\chi^1(M) \supset \cdots \supset C_\chi^d(M) = \{0\}$ of \mM by $G$-submodules 
 $C_\chi^i(M)$ such that $\forall \ i$, \mG operates trivially on $C_\chi^i(M)/C_\chi^{i+1}(M)$ and 
 $C_\chi^i(M)/C_\chi^{i+1}(M)$ is $p$-cotorsion.]\\

%%%%%%%%%%%%%%%%%%%%%%%%%%%%%%%%%%%%%%
%%%%%%%%%%%%%%%%%%%%%%%%%%%%%%%%%%%%%%
%%%%%%%%%%%%%%%%%%%%%%%%%%%%%%%%%%%%%%

\begin{center}
$\S \ 8$
\\[0.5cm]
$\mathcal{REFERENCES}$\\
\end{center}

\[
\textbf{BOOKS}
\]

\begingroup
\fontsize{9pt}{11pt}\selectfont
\setlength\parindent{0 cm}

[1] \quad Griffith, P., \textit{Infinite Abelian Group Theory}, University of Chicago Press (1970).
\\[-.2cm]

[2] \quad Hilton, P., \textit{Nilpotente Gruppen und Nilpotente R\"aume}, Springer Verlag (1984).
\\[-.2cm]

[3] \quad Hilton, P., Mislin, G., and Roitberg, J., \textit{Localization of Nilpotent Groups and Spaces}, North Holland 

\hspace{0.8cm}(1975).
\\[-.2cm]

[4] \quad Warfield, R., \textit{Nilpotent Groups}, Springer Verlag (1976).
\\[-.2cm]
\endgroup

\[
\textbf{ARTICLES}
\]

\begingroup
\fontsize{9pt}{11pt}\selectfont
\setlength\parindent{0 cm}

[1] \quad Baumslag, G., Some Aspects of Groups with Unique Roots, 
\textit{Acta Math.} \textbf{104} (1960), 217-303.
\\[-.2cm]

[2] \quad Baumslag, G., Lecture Notes on Nilpotent Groups, 
\textit{CBMS Regional Conference} \textbf{2} (1971), 1-73.
\\[-.2cm]

[3] \quad Bousfield, A., Homological Localization Towers for Groups and $\Pi$-Modules, 
\textit{Memoirs Amer. Math.}

\hspace{0.8cm}\textit{Soc.} \textbf{186} (1977), 1-68.
\\[-.2cm]

[4] \quad Bousfield, A., Constructions of Factorization Systems in Categories, 
\textit{J. Pure Appl. Algebra} \textbf{9} (1977), 

\hspace{0.8cm}207-220.
\\[-.2cm]

[5] \quad Hilton, P., Localization and Cohomology of Nilpotent Groups, 
\textit{Math. Zeit.} \textbf{132} (1973), 263-286.
\\[-.2cm]

[6] \quad Hilton, P., Remarks on the Localization of Nilpotent Groups, 
\textit{Comm. Pure Appl. Math.} \textbf{24} (1973), 

\hspace{0.8cm}703-713.
\\[-.2cm]

[7] \quad Ribenboim, P., Equations in Groups, with Special Emphasis on Localization and Torsion I and II,

\hspace{0.8cm}\textit{Atti Accad. Naz. Lincei Mem. Cl. Sci. Fis. Mat. Natur. Sez.} \textbf{19} (1987), 23-60 
and \textit{Portugal. Math.} 

\hspace{0.8cm}\textbf{44} (1987), 417-445.
\\[-.2cm]

[8] \quad Stallings, J., Homology and Central Series of Groups, 
\textit{J. Algebra} \textbf{2} (1965), 170-181.

\setlength\parindent{2em}

\endgroup

\chapter{
$\boldsymbol{\S}$\textbf{9}.\quadx  HOMOTOPICAL LOCALIZATION}
\setlength\parindent{2em}
\setcounter{proposition}{0}

%%----------------------------------------------------------------------------------------------01
$\text{ }$\\[-1.25cm]

Localization at a set of primes is a powerful tool in commutative algebra and group theory, thus it should come as no surprise that the transcription of this process to algebraic topology is of fundamental importance.  
More generally, one can interpret ``localization'' as the search for and construction of reflective subcategories in a homotopy category.\\

\label{9.90}
\label{11.24}
\begingroup%%----------------------------------->>
\fontsize{9pt}{11pt}\selectfont
\textbf{\small EXAMPLE} \  
\bHCW is not a reflective subcategory of \bHTOP.  
Reason: \bHCW is not isomorphism closed.  
\bHCWSP is not a reflective subcategory of \bHTOP.  
Reason: \bHCWSP is not limit closed (e.g., the product $\ds\prod\limits_1^\infty \bS^n$ is not a CW space).  
On the other hand, \bHCWSP is a coreflective subcategory of \bHTOP, 
the coreflector being the functor that assigns to each topological space \mX the geometric realization of its singular set (the arrow of adjunction $\as X \ra X$ is a weak homotopy equivalence (Giever-Milnor theorem)).  
In particular: \bHCWSP has products, viz. the product of $\{X_i\}$ in \bHCWSP is 
%$\as {\ds\prod\limits_i X_i}$, 
$\bigl | \sin {\ds\prod\limits_i X_i} \bigr |$, 
where $\ds\prod\limits_i X_i$ is the product in \bHTOP (or still, the product in \bTOP).
\\ \indent
[Note: \ Analogous remarks apply in the pointed setting.  So, e.g., the $n^{th}$ homotopy group of 
$\ds\prod\limits_i X_i$ 
(taken in $\bHCWSP_*$) is isomorphic to $\ds\prod\limits_i \pi_n(X_i)$.]\\
\endgroup %%------------------------------------<<

Notation: $\bCONCWSP_*$ 
\index{$\bCONCWSP_*$} 
is the full subcategory of $\bCWSP_*$ whose objects are the pointed connected CW spaces and $\bHCONCWSP_*$
\index{$\bHCONCWSP_*$} 
is the associated homotopy category.\\

\label{9.101}
\label{9.113}
\begingroup%%----------------------------------->>
\fontsize{9pt}{11pt}\selectfont
\textbf{\small EXAMPLE} \  
Write $\bHCONCWSP_*[n]$ for the full subcategory of $\bHCONCWSP_*$ whose objects have trivial homotopy groups in dimension $> n$ $(n \geq 0)$ $-$then $\bHCONCWSP_*[n]$ is a reflective subcategory of $\bHCONCWSP_*$, the reflector being the functor that assigns to each \mX its $n^\text{th}$ Postnikov approximate $X[n]$.  
Example: The fundamental group functor $X \ra \pi_1(X)$ sets up an equivalence between $\bHCONCWSP_*[1]$ and \bGR.
\\ \indent
[Note: The data generates an orthogonal pair $(S,D)$.   
Here $[f]:X \ra Y$ is in \mS iff $f_* : \pi_q(X) \ra \pi_q(Y)$ is bijective for $q \leq n$.]\\
\endgroup %%------------------------------------<<

\begingroup%%----------------------------------->>
\fontsize{9pt}{11pt}\selectfont
\textbf{\small EXAMPLE} \  
Write $\bHSCONCWSP_*$ for the full subcategory of $\bHCONCWSP_*$ whose objects are simply connected $-$then 
$\bHSCONCWSP_*$ is not a reflective subcategory of $\bHCONCWSP_*$.  
For suppose it were and, to get a contradiction, take 
$X = \bP^2(\R)$.  Consider, in the notation of 
p. \pageref{9.1}, $\epsilon_X:X \ra TX$.  By definition, 
$\epsilon_X \perp K(\Z,2)$ $\implies$ 
$H^2(TX) \approx$ 
$H^2(X) \approx \Z/2\Z$.  But $H_1(TX) = 0$ $\implies$ $H^2(TX) \approx \Hom(H_2(TX),\Z)$, which is torsion free.
\\ \indent
[Note: \ Let $f:\bS^1 \ra *$ $-$then $f^\perp$ is the object class of $\bHSCONCWSP_*$.]\\
\endgroup %%------------------------------------<<

Given a set of primes $P$, a pointed connected CW space \mX is said to be 
\un{$P$-local in} \un{homotopy}
\index{P-local in homotopy} 
if $\forall \ n \geq 1$, $\pi_n(X)$ is $P$-local.\\

%%----------------------------------------------------------------------------------------------02

\label{9.40}
\label{9.87}
\label{18.37}
\label{9.114}
\begingroup%%----------------------------------->>
\fontsize{9pt}{11pt}\selectfont
\textbf{\small EXAMPLE} \  
Fix $P \neq \bPi$ $-$then the full subcategory of $\bHCONCWSP_*$ whose objects are $P$-local in homotopy is not the object class of a reflective subcategory of $\bHCONCWSP_*$.  To see this, suppose the opposite and consider $\bS^1$.  Calling its localization $\bS_P^1$, for any $P$-local group \mG, the universal arrow $l_P:\bS^1 \ra \bS_P^1$ necessarily induces a bijection 
$[\bS_P^1,K(G,1)] \approx [\bS^1,K(G,1)]$ 
$\implies$ $\Hom(\pi_1(\bS_P^1),G) \approx \Hom(\pi_1(\bS^1),G)$.  Since $\pi_1(\bS_P^1)$ is by definition $P$-local, it follows that $\pi_1(\bS_P^1) \approx \Z_P$.  Form now 
$K(\Q[\Z_P],2;\chi)$, where $\chi:\Z_P \ra \Aut \Q[\Z_P]$ is the homomorphism corresponding to the action of $\Z_P$ on 
$\Q[\Z_P]$.  
Since $K(\Q[\Z_P],2;\chi)$ is $P$-local, the bijection 
$[\bS_P^1,K(\Q[\Z_P],2;\chi)] \approx [\bS^1,K(\Q[\Z_P],2;\chi)]$ 
restricts to an isomorphism 
$H^2(\bS_P^1;\Q[\Z_P])$ $\approx$ $
H^2(\bS^1;\Q[\Z_P])$ 
(cf. p. \pageref{9.2}) (locally constant coefficients), thus
$H^2(\bS_P^1;\Q[\Z_P]) = 0$.  
But 
$H^2(\pi_1(\bS_P^1);\Q[\Z_P])$ embeds in \hsx
$H^2(\bS_P^1;\Q[\Z_P])$ 
(consider the spectral sequence 
\hsx $E_2^{p,q} \hsx \approx \hsx H^p(\pi_1(\bS_P^1);H^q(\widetilde{\bS}_P^1;\Q[\Z_P]))$ %dmc the tilde over the 1 or not?
\hsx $\Rightarrow$ \hsx $H^{p+q}(\bS_P^1;\Q[\Z_P]))$, 
which contradicts the fact that $H^2(\Z_P;\Q[\Z_P]) \neq 0$ 
(cf. p. \pageref{9.3}).
\\ \indent
[Note: \ Let $\rho_n^q:\bS^q \ra \bS^q$ $(q \geq 1)$ be a map of degree $n$ $(n \in S_P)$.  Working in 
$\bHCONCWSP_*$, put $S_0 = \{[\rho_n^q]\}$ $-$then $S_0^\perp$ is the class of objects in $\bHCONCWSP_*$ which are 
$P$-local in homotopy.]\\
\endgroup %%------------------------------------<<

Given integers $k,n > 1$, let $k:\bS^{n-1} \ra \bS^{n-1}$ be a map of degree $k$ $-$then the adjunction space 
$\bP^n(k) = \bD^n \sqcup_k \bS^{n-1}$ is a Moore space of type $(\Z/k\Z,n-1)$ and 
$\Sigma \bP^n (k) = \bP^{n+1}(k)$.

Given a pointed connected CW space \mX, the 
\un{$n^\text{th}$ mod $k$ homotopy group}
\index{n$^\text{th}$ mod $k$ homotopy group} 
of \mX is 
$[\bP^n(k),X]$, the set of pointed homotopy classes of pointed continuous functions $\bP^n(k) \ra X$.  
Notation: $\pi_n(X;\Z/k\Z)$.  Here, the language is slightly deceptive.  
While it is true that $\pi_n(X;\Z/k\Z)$ is a group if 
$n > 2$ (which is abelian if $n > 3$), $\pi_2(X;\Z/k\Z)$ is merely a pointed set (but there is a left action 
$\pi_2(X) \times \pi_2(X;\Z/k\Z) \ra \pi_2(X;\Z/k\Z)$.  In the event that $\pi_1(X)$ is abelian, put 
$\pi_1(X;\Z/k\Z)  = \pi_1(X) \otimes \Z/k\Z$.

[Note: \ When \mX is an H-space, $\pi_2(X;\Z/k\Z)$ is a group (and $\pi_n(X;\Z/k\Z)$ is abelian if $n > 2$).]\\

A pointed continuous function $f:X \ra Y$ between pointed connected CW spaces induces a map 
$f_*: \pi_n(X;\Z/k\Z) \ra \pi_n(Y; \Z/k\Z)$.  
It is a homomorphism if $n > 2$ and respects the action of $\pi_2$ if $n = 2$.\\

\index{Theorem: Universal Coefficient Theorem}
\index{Universal Coefficient Theorem}
\textbf{\small UNIVERSAL COEFFICIENT THEOREM} \  
For each $n > 1$, there is a functorial exact sequence 
$0 \ra \pi_n(X) \otimes \Z/k\Z \ra \pi_n(X;\Z/k\Z) \ra \Tor(\pi_{n-1}(X),\Z/k\Z) \ra 0$.

[The arrows
$\bS^{n-1} \overset{k}{\lra} \bS^{n-1}$, 
$\bS^{n-1} \lra \bP^n(k) \lra \bS^n$, 
$\bS^n \overset{k}{\lra} \bS^n$ generate a functorial exact sequence 
$\pi_n(X) \overset{k}{\lra} \pi_n(X) \lra \pi_n(X;\Z/k\Z) \lra \pi_{n-1}(X)$ 
$\overset{k}{\lra} \pi_{n-1}(X)$.]

[Note: \ If $n = 2$, interpret exactness in $\bSET_*$ and if $\pi_1(X)$ is not abelian, interpret 
$\Tor(\pi_1(X),\Z/k\Z)$ as the kernel of $\pi_1(X) \overset{k}{\lra} \pi_1(X)$.]\\

Example: Let \mX be a pointed connected CW space $-$then \mX is $P$-local in homotopy
%%----------------------------------------------------------------------------------------------03
iff $\pi_1(X)$ is $P$-local and $\forall \ p \in \ov{P}$, $\pi_n(X;\Z/p\Z) = 0$ $\forall \ n > 1$.

[Apply REC$_2$ of the recognition principle 
(cf. p. \pageref{9.4} ff.).]\\

\begingroup%%----------------------------------->>
\fontsize{9pt}{11pt}\selectfont
Neisendorfer\footnote[2]{\textit{Memoirs Amer. Math. Soc.} \textbf{232} (1980), 1-67.}
has established a mod $k$ analog of the Hurewicz theorem.\\
\endgroup %%------------------------------------<<

\begingroup%%----------------------------------->>
\fontsize{9pt}{11pt}\selectfont
\index{Theorem: Mod $k$ Hurewicz Theorem}
\index{Mod $k$ Hurewicz Theorem}
\textbf{\small MOD $k$ HUREWICZ THEOREM} \quad 
Suppose that \mX is a pointed abelian CW space $-$then if $n \geq 2$, the condition 
$\pi_q(X;\Z/k\Z) = 0$ $(1 \leq q < n)$ is equivalent to the condition $H_q(X;\Z/k\Z) = 0$ $(1 \leq q < n)$ and either implies 
that the Hurewicz map $\pi_n(X;\Z/k\Z) \ra H_n(X;\Z/k\Z)$ is bijective.
\\ \indent
[Note: \ The arrow $\bP^n(k) \ra \bS^n$ induces as isomorphism 
$H_n(\bP^n(k);\Z/k\Z) \ra H_n(\bS^n;\Z/k\Z)$, so there is a generator of $H_n(\bP^n(k);\Z/k\Z)$ that is sent to the canonical generator of $H_n(\bS^n;\Z/k\Z)$, from which the Hurewicz map 
$\pi_n(X;\Z/k\Z) \ra H_n(X;\Z/k\Z)$ (it is a homomorphism if $n > 2$).]\\
\endgroup %%------------------------------------<<

\begingroup%%----------------------------------->>
\fontsize{9pt}{11pt}\selectfont
The mod $k$ analog of the Whitehead theorem is also true (consult 
Suslin\footnote[3]{\textit{J. Pure Appl. Algebra} \textbf{34} (1984), 301-318.}
for a variant with applications to algebraic K-theory).\\
\endgroup %%------------------------------------<<

Given a set of primes \mP, a pointed connected CW space \mX is said to be 
\un{$P$-local in} \un{homology}
\index{P-local in homology} 
if $\forall \ n \geq 1$, $H_n(X)$ is $P$-local.

[Note: \ \mX is $P$-local in homology iff $\forall \ p \in \ov{P}$, $H_n(X;\Z/p\Z) = 0$ $\forall \ n \geq 1$ 
(cf. p. \pageref{9.5}).]\\

\begingroup%%----------------------------------->>
\fontsize{9pt}{11pt}\selectfont
\label{9.41}
\textbf{\small EXAMPLE} \  
Fix $P \neq \bPi$ $-$then there exists a pointed connected CW space \mX such that $\forall \ n \geq 2$, 
$\pi_n(X) \approx \Z$ and $\forall \ n \geq 1$, $H_n(X) \approx \Z_P$ 
(cf. p. \pageref{9.6}), so $P$-local in homology need not imply $P$-local in homotopy.
\vspi
[Note: \ In the other direction, $P$-local in homotopy need not imply $P$-local in homology.  Reason: There exists a $P$-local group \mG such that $G/[G,G] (\approx H_1(G))$ has an $S_P$-torsion direct summand 
(cf. p. \pageref{9.7}), 
e.g., $G = (\Z*\Z)_P$.]\\
\endgroup %%------------------------------------<<

\begin{proposition} \ %01
Let 
$
\begin{cases}
\ X\\
\ Y
\end{cases}
$
be pointed nilpotent CW spaces, $f:X \ra Y$ a pointed continuous function.  Assume: $\forall \ n \geq 1$, 
$f_*:\pi_n(X) \ra \pi_n(Y)$ is $P$-localizing $-$then $\forall \ n \geq 1$, $f_*: H_n(X) \ra H_n(Y)$ is $P$-localizing.
\end{proposition}

[There is a \cd
\begin{tikzcd}%[sep=large]
{\widetilde{X}} \ar{d}[swap]{\widetilde{f}} \ar{r} &{X} \ar{d}{f}\\
{\widetilde{Y}} \ar{r} &{Y}
\end{tikzcd}
and a morphism 
$\{E_{p,q}^2 \approx$ 
$H_p(\pi_1(X);$ $H_q(\widetilde{X}))\} \ra$ 
$\{\ov{E}_{p,q}^2 \approx H_p(\pi_1(Y);H_q(\widetilde{Y}))\}$ 
of fibration spectral sequences.  Since
$
\begin{cases}
\ \widetilde{X}\\
\ \widetilde{Y}
\end{cases}
$
are
%%----------------------------------------------------------------------------------------------04
simply connected, $\forall \ q \geq 1$, $\widetilde{f}_*:H_q(\widetilde{X}) \ra H_q(\widetilde{Y})$ is $P$-localizing 
(cf. p. \pageref{9.8}).  In addition, $\forall \ q \geq 1$, 
$
\begin{cases}
\ \pi_1(X)\\
\ \pi_1(Y)
\end{cases}
$
operates nilpotently on 
$
\begin{cases}
\ H_q(\widetilde{X})\\
\ H_q(\widetilde{Y})
\end{cases}
$
(cf. $\S 5$, Proposition 17), thus $\forall \ q \geq 1$, the arrow 
$E_{p,q}^2 \ra \ov{E}_{p,q}^2$ is $P$-localizing (cf. $\S 8$, Proposition 14).  
Recalling that $\forall \ p \geq 1$, the arrow 
$H_p(\pi_1(X)) \ra H_p(\pi_1(Y))$ is $P$-localizing (cf. $\S 8$, Proposition 10), one can pass through the spectral sequence to see that $\forall \ q \geq 1$, $f_*: H_q(X) \ra H_q(Y)$ is $P$-localizing.]\\

\label{9.14}
Application: Let \mX be a pointed nilpotent CW space.  Assume: \mX is $P$-local in homotopy $-$then \mX is $P$-local in homology.

[Note: \ The converse is also true 
(cf. p. \pageref{9.9}).]\\

\begin{proposition} \ %02
Let 
$
\begin{cases}
\ X\\
\ Y
\end{cases}
$
be pointed nilpotent CW spaces, $f:X \ra Y$ a pointed continuous function.  
Assume: $\forall \ n \geq 1$, $f_*: H_n(X) \ra H_n(Y)$ is $P$-localizing 
$-$then for any pointed nilpotent CW space \mZ which is $P$-local in homotopy, the precomposition arrow 
$f^*:[Y,Z] \ra [X,Z]$ is bijective.
\end{proposition}

[There is no loss in generality in supposing that 
$
\begin{cases}
\ X\\
\ Y
\end{cases}
$
are pointed nilpotent CW complexes with \mX a pointed subcomplex of \mY (take $f$ skeletal and replace \mY by the pointed mapping cylinder of $f$).  
Because the inclusion $X \ra Y$ is a cofibration, this reduction converts the problem into one that can be treated by obstruction theory.  
Thus given a pointed continuous function $\phi:X \ra Z$, the obstructions to extending $\phi$ to a pointed continuous function $\Phi:Y \ra Z$ and the obstructions to any two such being homotopic rel \mX (hence pointed homotopic) lie in the 
$H^p(Y,X;\Gamma_{\chi_q}^i(\pi_q(Z))/\Gamma_{\chi_q}^{i+1}(\pi_q(Z)))$ for certain $p$ and $q$ (nilpotent obstruction theorem).  
The claim is that these groups are trivial.  
But, by hypothesis, $\forall \ n \geq 1$, 
$f_*:H_n(X;\Z_P) \ra H_n(Y;\Z_P)$ is an isomorphism, hence $\forall \ n \geq 1$, $H_n(Y,X;\Z_P) = 0$.  Since 
$\Z_P$ is a principal ideal domain and since the 
$\Gamma_{\chi_q}^i(\pi_q(Z))/\Gamma_{\chi_q}^{i+1}(\pi_q(Z))$ are $\Z_P$-modules 
(cf. p. \pageref{9.10}), the universal coefficient theorem implies that the obstructions to existence and uniqueness do indeed vanish.]

[Note: \ Otherwise said, under the stated conditions, $[f] \perp Z$ for any pointed nilpotent CW space \mZ which is $P$-local in homotopy.]\\

\label{9.42}
Notation: $\bNILCWSP_*$
\index{$\bNILCWSP_*$} 
is the full subcategory of $\bCWSP_*$ whose objects are the pointed nilpotent CW spaces and $\bHNILCWSP_*$
\index{$\bHNILCWSP_*$} 
is the associated homotopy category, while 
$\bNILCWSP_{*,P}$
\index{$\bNILCWSP_{*,P}$} 
is the full subcategory of $\bNILCWSP_*$ whose objects are the pointed nilpotent CW spaces which are $P$-local in homotopy and $\bHNILCWSP_{*,P}$
\index{$\bHNILCWSP_{*,P}$} 
is the associated homotopy category.\\

%%----------------------------------------------------------------------------------------------05
\index{Theorem: Nilpotent $P$-Localization Theorem}
\index{Nilpotent $P$-Localization Theorem}
\textbf{\small NILPOTENT $P$-LOCALIZATION THEOREM} \ \  
$\bHNILCWSP_{*,P}$ is a reflective subcategory of $\bHNILCWSP_{*}$.

[On general grounds, it is a question of assigning to each \mX in $\bHNILCWSP_*$ an object in $X_P$ in 
$\bHNILCWSP_{*,P}$ and a pointed homotopy class $[l_P]:X \ra X_P$ with the property that for any pointed homotopy class 
$[f]:X \ra Y$, where \mY is in $\bHNILCWSP_{*,P}$, there exists a unique pointed homotopy class 
$[\phi]: X_P \ra Y$ such that $[f] = [\phi] \circx [l_P]$.  
In view of Propositions 1 and 2, it will be enough to construct a pair $(X_P,l_P)$: $\forall \ q \geq 1$, 
$\pi_q(l_P):\pi_q(X) \ra \pi_q(X_P)$ is $P$-localizing.  
For this, we shall work first with the $n^\text{th}$ Postnikov approximate $X[n]$ of \mX and produce $(X[n]_P,l_P)$ inductively.  
Matters being plain if $n = 0$ ($X[0]$ is contractible), take $n > 0$.  
Consider a principal refinement of order $n$ of the arrow 
$X[n] \ra X[n-1]$, i.e., a factorization 
$X[n] \overset{\Lambda}{\ra}$ 
$W_N \overset{q_N}{\ra}$ 
$W_{N-1}  \ra $
$\cdots \ra$ 
$W_1 \overset{q_1}{\ra}$ 
$W_0 = X[n-1]$, where $\Lambda$ is a pointed homotopy equivalence and each $q_i:W_i \ra W_{i-1}$ is a pointed Hurewicz fibration for which there is an abelian group $\pi_i$ and a pointed continuous function 
$\Phi_{i-1}:W_{i-1} \ra K(\pi_i,n+1)$ such that the diagram
\begin{tikzcd}%[sep=large]
{W_i} \ar{d}[swap]{q_i} \ar{r} &{\Theta K(\pi_i,n+1)} \ar{d}\\
{W_{i-1}} \ar{r}[swap]{\Phi_{i-1}} &{K(\pi_i,n+1)}
\end{tikzcd}
is a pullback square.  
To exhibit pairs $(W_{i,P},l_P)$ 
(and hence produce ($X[n]_P,l_P)$), one can proceed via recursion on $i > 0$, the existence of $(W_{0,P},l_P)$ 
being secured by the induction hypothesis.  Choose a filler 
$\Phi_{i-1,P}:W_{i-1,P} \ra K(\pi_{i,P},n+1)$ for 
\begin{tikzcd}%[sep=large]
{W_{i-1}} \ar{d} \ar{r} &{K(\pi_i,n+1)} \ar{d}\\
{W_{i-1,P}} \ar[dashed]{r} &{K(\pi_{i,P},n+1)}
\end{tikzcd}
and define $W_{i,P}$ by the pullback square
\begin{tikzcd}%[sep=large]
{W_{i,P}} \ar{d} \ar{r} &{\Theta K(\pi_{i,P},n+1)} \ar{d}\\
{W_{i-1,P}}  \ar{r}[swap]{\Phi_{i-1,P}} &{K(\pi_{i,P},n+1)}
\end{tikzcd}
.  
Since the composite
$W_i \ra W_{i-1} \ra W_{i-1,P} \ra K(\pi_{i,P},n+1)$ is nullhomotopic, there is a filler $l_P:W_i \ra W_{i,P}$ for
\begin{tikzcd}%[sep=large]
{W_i} \ar[dashed]{d} \ar{r} &{W_{i-1}} \ar{d}\\
{W_{i,P}} \ar{r} &{W_{i-1,P}}
\end{tikzcd}
.  
From the definitions, 
$
\begin{cases}
\ W_i\\
\ W_{i,P}
\end{cases}
$
is a pointed connected CW space homeomorphic to 
$
\begin{cases}
\ E_{\Phi_{i-1}}\\
\ E_{\Phi_{i-1},P}
\end{cases}
$
(parameter reversal).  Moreover, 
$
\begin{cases}
\ W_i\\
\ W_{i,P}
\end{cases}
$
is nilpotent (cf. $\S 5$, Proposition 15) and by comparing the homotopy sequences of 
$
\begin{cases}
\ \Phi\\
\ \Phi_{i,P}
\end{cases}
$
one finds that $\forall \ q \geq 1$, $\pi_q(l_P):\pi_q(W_i) \ra \pi_q(W_{i,P})$ is $P$-localizing.  Recall now that $\forall \ n$, there is a pointed homotopy equivalence $X[n] \ra P_nX$ and a pointed Hurewicz fibration $P_nX \ra P_{n-1}X$ 
(cf. p. \pageref{9.11}).  
Passing to mapping tracks and changing $l_P$ within its pointed homotopy class, one can always arrange that $\forall \ n$, the arrow $(P_nX)_P \ra (P_{n-1}X)_P$ is a pointed Hurewicz
%%----------------------------------------------------------------------------------------------06
fibration and the diagram
\begin{tikzcd}%[sep=large]
{P_n X} \ar{d} \ar{r} &{P_{n-1} X} \ar{d}\\
{(P_n X)_P} \ar{r} &{(P_{n-1} X)_P}
\end{tikzcd}
commutes.  
So, 
$\lim l_P:\lim P_n X \ra \lim (P_n X)_P$ exists and $\forall \ q \geq 1$, 
$\pi_q(\lim l_P):\pi_q(\lim P_n X) \ra \pi_q(\lim(P_n X)_P)$
is $P$-localizing 
(cf. p. \pageref{9.12}).  
Fix a CW resolution $X_P \ra \lim (P_n X)_P$ and let $l_P:X \ra X_P$ be a filler for
\begin{tikzcd}%[sep=large]
{X} \ar[dashed]{d} \ar{r} &{\lim P_nX} \ar{d}\\
{X_P} \ar{r} &{\lim (P_nX)_P}
\end{tikzcd}
(cf. $\S 5$, Proposition 4).  
Because the arrow $X \ra \lim P_n X$ is a weak homotopy equivalence (cf. $\S 5$, Proposition 13), it follows that $\forall \ q \geq 1$, $\pi_q(l_P):\pi_q(X) \ra \pi_q(X_P)$ is $P$-localizing.]\\

The reflector $L_P$ figuring in the nilpotent $P$-localization theorem sends \mX to $X_P$ (special cases: $X_{\Q}$, $X_p$ 
($p \in \bPi$)) with arrow of localization $[l_P]:X \ra X_P$.  
Brackets are often omitted, e.g., given $f:X \ra Y$, there is a diagram
\begin{tikzcd}%[sep=large]
{X} \ar{d} \ar{r}{f} &{Y} \ar{d}\\
{X_P} \ar{r}[swap]{f_P} &{Y_P}
\end{tikzcd}
, commutative up to pointed homotopy.

[Note: $L_P$ respects the ``abelian subcategory'' and the ``simply connected subcategory''.]

Let $[f]:X \ra Y$ be a morphism in $\bHNILCWSP_*$ $-$then $[f]$ (or $f$) is said to be 
\un{$P$-localizing}
\index{P-localizing (f a morphism in $\bHNILCWSP_*$)} 
if $\exists$ an isomorphism 
$[\phi]:X_P \ra Y$ such that $[f] = [\phi] \circx [l_P]$ 
(cf. p. \pageref{9.13}).\\

\begin{proposition} \ %03
Let 
$
\begin{cases}
\ X\\
\ Y
\end{cases}
$
be pointed nilpotent CW spaces, $f:X \ra Y$ a pointed continuous function $-$then $f$ is $P$-localizing iff $\forall \ n \geq 1$, 
$f_*:\pi_n(X) \ra \pi_n(Y)$ is $P$-localizing.
\end{proposition}

[This is implicit in the proof of the nilpotent $P$-localization theorem.]\\

\label{9.56}
Example: For any nilpotent group \mG, $K(G,1)_P \approx K(G_P,1)$.\\

\begin{proposition} \ %04
Let 
$
\begin{cases}
\ X\\
\ Y
\end{cases}
$
be pointed nilpotent CW spaces, $f:X \ra Y$ a pointed continuous function $-$then $f$ is $P$-localizing iff $\forall \ n \geq 1$, 
$f_*:H_n(X) \ra H_n(Y)$ is $P$-localizing.
\end{proposition}

[The point behind the sufficiency is that $\forall \ n \geq 1$, $H_n(Y)$ is $P$-local, therefore Dror's Whitehead theorem implies that 
$l_P:Y \ra Y_P$ is a pointed homotopy equivalence, thus \mY is $P$-local in homotopy.]\\

\label{9.9}
\label{9.32}
Application: Let \mX be a pointed nilpotent CW space.  Assume: \mX is $P$-local in homology $-$then \mX is $P$-local in homotopy.

%%----------------------------------------------------------------------------------------------07
[Note: \ The converse is also true 
(cf. p. \pageref{9.14}).]\\

\begingroup%%----------------------------------->>
\fontsize{9pt}{11pt}\selectfont
\textbf{\small FACT} \  
Let $P^\prime$ and $P\pp$ be two sets of primes $-$then for any pointed nilpotent CW space \mX, 
$(X_{P^\prime})_{P\pp} \approx (X_{P\pp})_{P^\prime}$.
\vspi
[The left hand side computes $X_{P^\prime \cap P\pp}$ and the right hand side computes 
$X_{P\pp \cap P^\prime}$.]\\
\endgroup %%------------------------------------<<

The nilpotent $P$-localization theorem has been relativized by 
Llerena\footnote[2]{\textit{Math. Zeit.} \textbf{188} (1985), 397-410.}.  
In fact, suppose that \ 
$
\begin{cases}
\ X\\
\ Y
\end{cases}
\& \ Z 
$
are pointed connected CW spaces.  
Let $f:X \ra Y$ be a pointed Hurewicz fibration with $E_f$ nilpotent $-$then there exists a pointed connected CW space $X(P)$, a pointed Hurewicz fibration $f(P):X(P) \ra Y$ with $E_{f(P)}$ nilpotent and $P$-local in homotopy, and a pointed continuous function $l(P):X \ra X(P)$ over \mY such that the induced map $E_f \ra E_{f(P)}$ is $P$-localizing: 
$(E_f)_P \approx E_{f(P)}$.  
In addition, for any pointed Hurewicz fibration $g:Z \ra Y$ with $E_g$ nilpotent and $P$-local in homotopy, $[f(P),g] \approx [f,g]$ in the sense of pointed fiber homotopy, i.e., given a commutative triangle 
\begin{tikzcd}[sep=small]
{X} \ar{rdd}[swap]{f} \ar{rr}{\phi} &&{Z} \ar{ldd}{g}\\
\\
&{Y}
\end{tikzcd},
there is a commutative triangle 
\begin{tikzcd}[sep=small]
{X(P)} \ar{drd}[swap]{f(P)} \ar{rr}{\phi(P)} &&{Z} \ar{ldd}{g}\\
\\
&{Y}
\end{tikzcd}
: $[\phi] = [\phi(P)] \circx [l(P)]$, $\phi(P)$ being unique up to pointed fiber homotopy.\\

\begingroup%%----------------------------------->>
\fontsize{9pt}{11pt}\selectfont
\textbf{\small EXAMPLE} \  
Let \mX be a pointed connected CW space $-$then the diagram\\
\begin{tikzcd}[sep=large]
{\widetilde{X}} \ar{d}[swap]{l_P} \ar{r} 
&{X} \ar{d}\ar{r} 
&{K(\pi_1(X),1)} \arrow[d,shift right=0.5,dash] \arrow[d,shift right=-0.5,dash]\\
{\widetilde{X}_p} \ar{r} 
&{X(P)} \ar{r}
&{K(\pi_1(X),1)}
\end{tikzcd}
commutes in $\bHCONCWSP_*$ 
(cf. p. \pageref{9.15}).  
Here, $\pi_1(X) \approx \pi_1(X(P))$ and $\forall \ n \geq 2$, the arrow $\pi_n(X) \ra \pi_n(X(P))$ is $P$-localizing.\\
\endgroup %%------------------------------------<<

Nilpotent $P$-localization is compatible with homotopy and homology in that $\forall \ n \geq 1$, 
$\pi_n(X)_P \approx \pi_n(X_P)$ and 
$H_n(X)_P \approx H_n(X_P)$ but this is false for cohomology.  
Example: Take $X = \bS^n$: $\bS_P^n = M(\Z_P,n)$ $\implies$ 
$H^{n+1}(\bS_P^n) \approx \Ext(\Z_P,\Z) \neq 0$ $(P \neq \bPi)$.

[Note: \ By contrast, taking coefficients in $\Z_P$, $\forall \ n \geq 1$, $H^n(X_P;\Z_P) \approx H^n(X;\Z_P)$ 
(cf. $\S 8$, Proposition 2).]

Let $[f]:X \ra Y$ be a morphism in $\bHNILCWSP_*$ $-$then $[f]$ (or $f$) is said to be a 
\un{$P$-equivalence}
\index{P-equivalence} 
if $f_P:X_P \ra Y_P$ is a pointed homotopy equivalence.  
With regard to the underlying orthogonal pair $(S,D)$, $[f]$ is a $P$-equivalence iff $[f] \in S$, so $[f]$ is $P$-localizing iff $[f] \in S$ $\&$ 
$Y \in D$ 
(cf. p. \pageref{9.16}).

%%----------------------------------------------------------------------------------------------08
[Note: \ When $P = \emptyset$, the term is 
\un{rational equivalence}
\index{rational equivalence}.  
Examples: 
(1) There is a rational equivalence 
$\bS^3 \ra K(\Z,3)$ 
but there is no rational equivalence
$K(\Z,3) \ra \bS^3$; 
(2) There are rational equivalences 
$\bS^3 \vee \bS^5 \ra$ 
$\bS^3 \vee K(\Z,5),\ \bS^3 \vee K(\Z,5) \ra$ 
$K(\Z,3) \vee K(\Z,5),\ \bS^3 \vee \bS^5 \ra$ 
$K(\Z,3) \vee \bS^5$, \ 
$K(\Z,3) \vee \bS^5 \ \ra$ 
$K(\Z,3) \vee K(\Z,5)$ 
but there are no rational equivalences 
$\bS^3 \vee K(\Z,5) \ra$ 
$K(\Z,3) \vee \bS^5$, 
$K(\Z,3) \vee \bS^5 \ra$ 
$\bS^3 \vee K(\Z,5)$.]\\

\begin{proposition} \ %05
Let 
$
\begin{cases}
\ X\\
\ Y
\end{cases}
$
be pointed nilpotent CW spaces, $f:X \ra Y$ a pointed continuous function $-$then $f$ is a $P$-equivalence iff 
$f_*:H_*(X;\Z_P) \ra H_*(Y;\Z_P)$ is an isomorphism.
\end{proposition}

[Note: \ This holds iff $f_*:H_*(X;\Q) \ra H_*(Y;\Q)$ is an isomorphism and $\forall \ p \in P$, 
$f_*:H_*(X;\Z/p\Z) \ra H_*(Y;\Z/p\Z)$ is an isomorphism (cf. $\S 8$, Proposition 3).]\\

Example: Fix a positive integer $d$.  Let $P_d$ be the set of primes that do not divide $d$ $-$then 
$\bS^n \overset{d}{\lra} \bS^n$ is a $P_d$-equivalence.\\

\begingroup%%----------------------------------->>
\fontsize{9pt}{11pt}\selectfont
\index{local spheres (example)}
\textbf{\small EXAMPLE \ (\un{Local Spheres})} \  
Given $P$, let $p_1 < p_2 < \cdots$ be an enumeration of the elements of $\ov{P}$ and put 
$d_k = p_1^k \cdots p_k^k$ $(k = 1, 2\ldots)$ $-$then a model for $\bS_P^n$ is the pointed mapping telescope of the sequence 
$\bS^n \ra \bS^n \ra \cdots$, the $k^\text{th}$ map having degree $d_k$.  Since $\Q$ is $P$-local, 
$H^*(\bS_P^n;\Q) \approx H^*(\bS^n;\Q)$.  
Accordingly, $\bS_P^n$ cannot be an H-space if $n$ is even (Hopf).  
As for what happens when $n$ is odd, 
Adams\footnote[2]{\textit{Quart. J. Math.} \textbf{12} (1961), 52-60.}
has shown that if $2 \notin P$, then $\bS_P^n$ is an H-space while if $2 \in P$, then $\bS_P^n$ is an H-space iff 
$n = 1, 3$, or 7.\\
\endgroup %%------------------------------------<<

\begingroup%%----------------------------------->>
\fontsize{9pt}{11pt}\selectfont
\index{rational spheres (example)}
\textbf{\small EXAMPLE \ (\un{Rational Spheres})} \  
If $n$ is odd, then $\bS_{\Q}^n = K(\Q,n)$ but if $n$ is even, then $\bS_{\Q}^n =$ $E_f$, where 
$f:K(\Q,n) \ra K(\Q,2n)$ corresponds to 
$t^2 \in H^{2n}(\Q,n;\Q)$ ($H^*(\Q,n;\Q) = \Q[t]$, $\abs{t} = n$).  
Consequently, if $n$ is odd, then 
$
\Q \otimes \pi_q(\bS^n) = 
\begin{cases}
\ \Q \quad (q = n)\\
\ 0 \quad (q \neq n)
\end{cases}
$
but if $n$ is even, then $\Q \otimes \pi_q(\bS^n) =$ 
$
\begin{cases}
\ \Q \quad (q = n, 2n - 1)\\
\ 0 \quad (q \neq n, 2n - 1)
\end{cases}
$
(cf. p. \pageref{9.17}).\\
\endgroup %%------------------------------------<<

\begin{proposition} \ %06
Let 
$
\begin{cases}
\ X\\
\ Y
\end{cases}
$
be pointed nilpotent CW spaces, $f:X \ra Y$ a pointed continuous function.  
Suppose that $f$ is a $P$-equivalence $-$then for any $P^\prime \subset P$, $f$ is a $P^\prime$-equivalence.\\
\end{proposition}

\begin{proposition} \ %07
Let 
$
\begin{cases}
\ X\\
\ Y
\end{cases}
$
be pointed nilpotent CW spaces, $f:X \ra Y$ a pointed continuous function.  Suppose that $f$ is a $P^\prime$-equivalence and a 
$P\pp$-equivalence $-$then $f$ is a $(P^\prime \cup P\pp)$-equivalence.\\
\end{proposition}

%%----------------------------------------------------------------------------------------------09
\begingroup%%----------------------------------->>
\fontsize{9pt}{11pt}\selectfont
\textbf{\small FACT} \  
Let \mX be a pointed nilpotent CW space.  Fix \mP $-$then for any $P^\prime \subset P$, the canonical arrow 
$X_P \ra X_{P^\prime}$ is a $(P^\prime \cup \ov{P})$-equivalence.\\
\endgroup %%------------------------------------<<

\begingroup%%----------------------------------->>
\fontsize{9pt}{11pt}\selectfont
\textbf{\small EXAMPLE} \  
Let 
$
\begin{cases}
\ X\\
\ Y
\end{cases}
$
be pointed nilpotent CW spaces.  Assume: $\exists$ a pointed homotopy equivalence $\phi:X_{\Q} \ra Y_{\Q}$ $-$then there is a pointed nilpotent CW space \mZ such that 
$
\begin{cases}
\ Z_P \approx X_P\\
\ Z_{\ov{P}} \approx Y_{\ov{P}}
\end{cases}
$
.
\vspi
[Choose $r_P:X_P \ra X_{\Q}$ 
$\hsx \& \hsx$ 
$r_{\ov{P}}:Y_{\ov{P}} \ra Y_{\Q}$ 
: 
$l_{\Q}\simeq r_P \circx l_P$ 
$(l_{\Q}:X \ra X_{\Q}$) 
$\hsx \& \hsx$ 
$l_{\Q} \simeq r_{\ov{P}} \circx l_{\ov{P}}$ $(l_{\Q}:Y \ra Y_{\Q})$.  
The double mapping track \mZ of the pointed 2-sink 
$X_P \overset{\phi\circx r_P}{\lra} Y_{\Q} \overset{r_{\ov{P}}}{\lla} Y_{\ov{P}}$
is a pointed CW space (cf. $\S 6$, Proposition 8).  
To check that \mZ is path connected (hence nilpotent) 
(cf. p. \pageref{9.18})), fix $\gamma \in \pi_1(Y_{\Q})$.  Since $\phi \circx r_P$ is a $\ov{P}$-equivalence and 
$r_{\ov{P}}$ is a $P$-equivalence, 
$\exists$ $m \in S_{\ov{P}}$: $\gamma^m = (\phi \circx r_P)_*(\alpha)$ 
$(\alpha \in \pi_1(X_P))$ $\&$ $\exists$ $n \in S_P$: $\gamma^n = (r_{\ov{P}})_*(\beta)$ $(\beta \in \pi_1(Y_{\ov{P}}))$.  
But $m$ and $n$ are relatively prime, so $\exists \ k$ and $l: km + ln = 1$ $\implies$ 
$\gamma = (\phi \circx r_P)_*(\alpha^k)\cdot(r_{\ov{P}})_*(\beta^l)$, which means that \mZ is path connected 
(cf. p. \pageref{9.19}).  And:
$
\begin{cases}
\ Z \ra X_P\\
\ Z \ra Y_{\ov{P}}
\end{cases}
$
is a 
$
\begin{cases}
\ \text{$P$-equivalence}\\
\ \text{$\ov{P}$-equivalence}
\end{cases}
.]
$
\\
\endgroup %%------------------------------------<<

\begin{proposition} \ %08
Let 
$
\begin{cases}
\ X\\
\ Y
\end{cases}
$
be pointed nilpotent CW spaces, $f:X \ra Y$ a pointed continuous function $-$then $f$ is a pointed homotopy equivalence provided that 
$\forall \ p$, $f_p:X_p \ra Y_p$ is a pointed homotopy equivalence.
\end{proposition}

[In fact, $\forall \ p$, $H_*(f)_p:H_*(X)_p \ra H_*(Y)_p$ is an isomorphism.  
Therefore $f$ is a homology equivalence 
(cf. p. \pageref{9.20}) and Dror's Whitehead theorem is applicable.]\\

In the simply connected situation, there is another approach to $P$-localization which depends on Proposition 2 but not on Proposition 1.  Thus let \mX be a pointed simply connected CW space $-$then it will be enough to construct a pair 
$(X_P,l_P)$: $\forall \ q \geq 1$, $H_q(l_P):H_q(X) \ra H_q(X_P)$ is $P$-localizing and for this one can assume that \mX is a pointed simply connected CW complex.

Observation: A model for $X_P$, where $X = \bigvee\limits_I \bS^n$ $(n > 1)$, is a Moore space of type 
$(I \cdot \Z_P,n)$: $X_P = \bigvee\limits_I M(\Z_P,n)$.

\indent\indent $(\dim X < \infty)$ \ If $\dim X = 2$, then \mX has the pointed homotopy type of a wedge 
$\bigvee\limits_I \bS^2$, hence $(X_P,l_P)$ exists in this case.  
Proceeding by induction on the dimension, suppose that 
$(X_P,l_P)$ has been constructed for all \mX with $\dim X \leq n$ $(n \geq 2)$ and consider an \mX with $\dim X = n + 1$.  
Up to pointed homotopy type, \mX is the pointed mapping cone $C_f$ of a pointed continuous function 
$f:\bigvee\limits_I \bS^n \ra X^{(n)}$ $(\#(I) = \#(\sE_{n+1}))$ and the pointed cofibration $j:X^{(n)} \ra C_f$ is a cofibration 
(cf. $\S 3$, Proposition 19).  Choose a filler
%%----------------------------------------------------------------------------------------------10
$f_P: \bigvee\limits_I \bS_P^n \ra X_P^{(n)}$ for 
\begin{tikzcd}%[sep=large]
{\bigvee\limits_I \bS^n} \ar{d} \ar{r} &{X^{(n)}} \ar{d}\\
{\bigvee\limits_I \bS_P^n} \ar[dashed]{r} &{X_P^{(n)}}
\end{tikzcd}
.  \ 
Since the composite 
$\bigvee\limits_I \bS^n \ra$
$X^{(n)} \ra$ 
$X_P^{(n)} \ra C_{f_P}$ 
is nullhomotopic, there is a filler $l_P:C_f \ra C_{f_P}$ for 
$
\begin{tikzcd}%[sep=large]
{X^{(n)}} \ar{d} \ar{r} &{C_f} \ar[dashed]{d}\\
{X_P^{(n)}} \ar{r} &{C_{f_P}}
\end{tikzcd}
.  
$ \ 
Assembling the data leads to a commutative diagram 
\[
\begin{tikzcd}%[sep=large]
{\widetilde{H}_{q}(\bigvee\limits_I \bS^n)} \ar{d} \ar{r} 
&{\widetilde{H}_{q}(X^{(n)})} \ar{d} \ar{r} 
&{\widetilde{H}_{q}(C_f)} \ar{d} \ar{r} 
&{\widetilde{H}_{q-1}(\bigvee\limits_I \bS^n)} \ar{d} \ar{r} 
&{\widetilde{H}_{q-1}(X^{(n)})}\ar{d}\\
{\widetilde{H}_{q}(\bigvee\limits_I \bS_P^n)} \ar{r} 
&{\widetilde{H}_{q}(X_P^{(n)})} \ar{r} 
&{\widetilde{H}_{q}(C_{f_P})}\ar{r}
&{\widetilde{H}_{q-1}(\bigvee\limits_I \bS_P^n)}\ar{r} 
&{\widetilde{H}_{q-1}(X_P^{(n)})}
\end{tikzcd}
\]
of abelian groups with exact rows, where both vertical arrows on either side of the arrow 
$\widetilde{H}_q(C_f) \ra \widetilde{H}_q(C_{f_P})$ are $P$-localizing.  
But this means that 
$\widetilde{H}_q(C_f) \ra \widetilde{H}_q(C_{f_P})$ is $P$-localizing as well 
(cf. p. \pageref{9.21}).\\
\indent\indent $(\dim X = \infty)$ \ One can arrange matters in such a way that $\forall \ n$, the diagram 
\begin{tikzcd}[sep=large]
{X^{(n)}} \ar{d} \ar{r} &{X^{(n+1)}} \ar{d}\\
{X_P^{(n)}} \ar{r} &{X_P^{(n+1)}}
\end{tikzcd}
is commutative and the arrow $X_P^{(n)} \ra X_P^{(n+1)}$ is a cofibration.  Put 
$X_P = \colim X_P^{(n)}$ (cf. $\S 5$, Proposition 8) and define $l_P:X \ra X_P$ in the obvious fashion.\\
\vspace{0.25cm}

\begingroup%%----------------------------------->>
\fontsize{9pt}{11pt}\selectfont
\textbf{\small FACT} \  
Let \mX $\&$ 
$
\begin{cases}
\ Y\\
\ Z
\end{cases}
$
be pointed simply connected CW spaces with finitely generated homotopy groups.  Suppose that $g:Y \ra Z$ is a rational equivalence $-$then $g$ induces a bijection $[X_{\Q},Y] \ra [X_{\Q},Z]$.

[Assuming that \mX is a pointed connected CW complex, construct $X_{\Q}$ as above, and show by induction that 
$\forall \ n$, $[X_{\Q}^{(n)},Y] \approx [X_{\Q}^{(n)},Z]$.]\\
\endgroup %%------------------------------------<<

\begingroup%%----------------------------------->>
\fontsize{9pt}{11pt}\selectfont
\index{phantom maps (example)}
\textbf{\small EXAMPLE \ (\un{Phantom Maps})} \  
The notion of phantom map, as defined on 
p. \pageref{9.22} for pointed connected CW complexes, extends to pointed connected CW spaces
$
\begin{cases}
\ X\\
\ Y
\end{cases}
\hspace{-.25cm}
: \Ph(X,Y).
$
This said, let
$
\begin{cases}
\ X\\
\ Y
\end{cases}
$
be pointed simply conncected CW spaces with finitely generated homotopy groups $-$then 
$\Ph(X,Y) = l_{\Q}^*[X_{\Q},Y]$ $\subset [X,Y]$ 
(cf. p. \pageref{9.23}).  
For instance, take $X = \Omega \bS^3$, $Y = \bS^3$.  
To compute $[\Omega \bS^3,\bS^3]$, note first that 
$\Sigma\Omega \bS^3 \approx$ 
$\Sigma\Omega\Sigma \bS^2 \approx$ 
$\Sigma \bigl(\ds\bigvee\limits_{n \geq 1} \bS^{2n}\bigr) \approx$ 
$\bigvee\limits_{n \geq 1} \bS^{2n+1}$ (cf. $\S 4$, Proposition 28 and subsequent discussion) and 
$\bS^3 \approx \Omega B_{\bS^3}^\infty$ 
(cf. p. \pageref{9.24}), hence 
$[\Omega \bS^3,\bS^3] \approx$ 
$[\Omega \bS^3,\Omega B_{\bS^3}^\infty] \approx$ 
$[\Sigma\Omega \bS^3,B_{\bS^3}^\infty] \approx$
$\bigl[\ds\bigvee\limits_{n \geq 1} \bS^{2n+1},B_{\bS^3}^\infty\bigr] \approx$ 
%%----------------------------------------------------------------------------------------------11
$\ds\prod\limits_{n \geq 1}[\bS^{2n + 1},B_{\bS^3}^\infty] \approx$ 
$\ds\prod\limits_{n \geq 1}[\bS^{2n},\bS^3]$.  
By the same token,
$[(\Omega \bS^3)_{\Q},\bS^3] \approx$ 
$[\Omega (\bS^3_{\Q}),\bS^3] \approx$ 
$\ds\prod\limits_{n \geq 1}[\bS_{\Q}^{2n},\bS^3]$ 
or still, 
$\approx \ds\prod\limits_{n \geq 1}[\bS_{\Q}^{2n},K(\Z,3)]$, the arrow $\bS^3 \ra K(\Z,3)$ being a rational equivalence.  
Conclusion: $\Ph(\Omega \bS^3,\bS^3) = 0$.\\
\endgroup %%------------------------------------<<

\begingroup%%----------------------------------->>
\fontsize{9pt}{11pt}\selectfont
\textbf{\small LEMMA}  \  
Let 
$
\begin{cases}
\ X\\
\ Y
\end{cases}
$
be pointed connected CW spaces, $f:X \ra Y$ a pointed Hurewicz fibration with $\pi_0(X_{y_0}) = *$ $-$then there is an exact sequence 
$\cdots \ra$ 
$\pi_{n+1}(Y;\Z/k\Z) \ra$ 
$\pi_n(X_{y_0};\Z/k\Z) \ra$
$\pi_n(X;\Z/k\Z) \ra$
$\pi_n(Y;\Z/k\Z) \ra$
$\cdots \ra \pi_2(Y;\Z/k\Z)$.\\
\endgroup %%------------------------------------<<

\begingroup%%----------------------------------->>
\fontsize{9pt}{11pt}\selectfont
\textbf{\small EXAMPLE} \  
Let 
$
\begin{cases}
\ X\\
\ Y
\end{cases}
$
be pointed simply connected CW spaces with finitly generated homotopy groups, $f:X \ra Y$ a pointed continuous function 
$-$then $f$ is a $p$-equivalence iff $\forall \ n \geq 2$, 
$f_*:\pi_n(X;\Z/p\Z) \ra \pi_n(Y;\Z/p\Z)$ is bijective.\\
\endgroup %%------------------------------------<<

\indent\indent (Products) \ 
Let 
$
\begin{cases}
\ X\\
\ Y
\end{cases}
$
be pointed nilpotent CW spaces $-$then $(X \times Y)_P \approx X_P \times Y_P$.\\

\begingroup%%----------------------------------->>
\fontsize{9pt}{11pt}\selectfont
\index{H-Spaces (example)}
\textbf{\small EXAMPLE \ (\un{H-Spaces})} \  
Suppose that \mX is a a path connected H-space $-$then $X_P$ is a path connected H-space and the arrow of localization 
$l_P:X \ra X_P$ is an H-map.\\
\endgroup %%------------------------------------<<

\indent\indent (Mapping Fibers) \ 
Let 
$
\begin{cases}
\ X\\
\ Y
\end{cases}
$
be pointed nilpotent CW spaces, $f:X \ra Y$ a pointed continuous function.  
Assume: $E_f$ is nilpotent $-$then 
$(E_f)_P \approx E_{f_P}$.

[Since $\pi_0(E_f) = *$, the arrow $\pi_1(X) \ra \pi_1(Y)$ is surjective, thus the same is true of the arrow 
$\pi_1(X)_P \ra \pi_1(Y)_P$ or still, of the arrow $\pi_1(X_P) \ra \pi_1(Y_P)$.  
Therefore $\pi_0(E_{f_P}) = *$ and 
$E_{f_P}$ is nilpotent 
(cf. p. \pageref{9.25}).  Compare the long exact sequences in homotopy.]\\

Application: Let $(K,k_0)$ be a pointed finite connected CW complex.  Suppose that $f:X \ra Y$ is $P$-localizing $-$then 
for any pointed continuous function $\phi:K \ra X$, the arrow $C(K,k_0;X,x_0:\phi) \ra C(K,k_0; Y,y_0:f \circx \phi)$ is 
$P$-localizing.

[Note: \ $C(\cdots:\phi)$, $C(\cdots:f \circx \phi)$ stand for the path component to which $\phi$, $f \circx \phi$ belong 
(cf. p. \pageref{9.26} ff.).]

Example: Given a pointed nilpotent CW space \mX, 
$(\Omega_0 X)_P \approx \Omega_0(X_P)$, where $\Omega_0$? is the path component of $\Omega$? containing the constant loop.\\

\begingroup%%----------------------------------->>
\fontsize{9pt}{11pt}\selectfont
\textbf{\small EXAMPLE} \  
Let \mX be a pointed nilpotent CW space.  Denote by $C_{\pi_P}$ the mapping cone of the pointed Hurewicz fibration 
$\pi_P:E_{l_P} \ra X$ $-$then the projection $C_{\pi_P} \ra X_P$ is a pointed homotopy equivalence iff \mX is $P$-local or 
$X_P$ is simply connected 
(cf. p. \pageref{9.27}).\\
\endgroup %%------------------------------------<<

%%----------------------------------------------------------------------------------------------12
\begingroup%%----------------------------------->>
\fontsize{9pt}{11pt}\selectfont
\textbf{\small FACT} \  
Let \mK be a finite CW complex; let \mX be a pointed nilpotent CW space.  
Fix a continuous function $\phi:K \ra X$.  
Denote by 
$C(K,X:\phi)$, $C(K,X_P : l_P \circx \phi)$ the path component of $C(K,X)$, $C(K,X_P)$ containing $\phi$, $\l_P \circx \phi$ 
$-$then $C(K,X:\phi)$ is nilpotent 
(cf. p. \pageref{9.28}) and 
$C(K,X:\phi)_P \approx$ 
$C(K,X_P:l_P \circx \phi)$.
\vspi
[Reduce to when \mK is connected and work with the Postnikov tower of \mX.]\\
\endgroup %%------------------------------------<<

\begingroup%%----------------------------------->>
\fontsize{9pt}{11pt}\selectfont
\textbf{\small EXAMPLE} \  
Let $X = \bS^{2m} \times \bS^{2n+1}$ $(m,n > 0)$ $-$then 
$C(X,X:\id_X)_{\Q} \approx$ 
$\ds\prod\limits_{i=1}^{4m-1} K(\Q^{d_i},i) \times \ds\prod\limits_{j=1}^{2n+1} K(H^{2n+1-j}(X;\Q),j)$, 
where 
$d_i = \dim_{\Q}H^{4m-1-i}(X;\Q) - \dim_{\Q}H^{2m-1-i}(X;\Q)$ 
(cf. p. \pageref{9.29}).\\
\endgroup %%------------------------------------<<

\indent\indent (Mapping Cones) \  Let 
$
\begin{cases}
\ X\\
\ Y
\end{cases}
$
be pointed nilpotent CW spaces, $f:X \ra Y$ a pointed continuous function.  Assume: $C_f$ is nilpotent $-$then 
$(C_f)_P \approx C_{f_P}$.

[$C_{f_P}$ is path connected and by Van Kampen, $\pi_1(C_{f_P}) \approx (\pi_1(C_f))_P$.  
But why is $C_{f_P}$ nilpotent?  For this, it is necessary to use the result of Rao mentioned on 
p. \pageref{9.30} (and transferred to the pointed setting).  
Take, e.g., the third possibility: $\exists$ a prime $p$ such that $\pi_1(C_f)$ is a finite $p$-group and $\forall \ q > 0$, 
$H_q(X)$ is a $p$-group of finite exponent.  
Case 1: $p \notin P$.  Here, 
$(\pi_1(C_f))_P = 1$ 
(cf. p. \pageref{9.31}) and $C_{f_P}$ is simply connected.  
Case 2: $p \in P$.  \mX is then $P$-local in homology, hence is $P$-local in homotopy 
(cf. p. \pageref{9.32}), i.e., 
$X \approx X_P$, and $\pi_1(C_f) \approx \pi_1(C_{f_P})$.  
Therefore $C_{f_P}$ is nilpotent.  Comparing the long exact sequences in homology finishes the proof.]

\label{11.14}
Example: Given a pointed nilpotent CW space \mX, $(\Sigma X)_P \approx \Sigma X_P$.\\

\begingroup%%----------------------------------->>
\fontsize{9pt}{11pt}\selectfont
\textbf{\small EXAMPLE} \ 
Let 
$
\begin{cases}
\ X\\
\ Y
\end{cases}
$
be pointed simply connected CW spaces $-$then 
$(X \# Y)_P \approx X_P \# Y_P$.
\vspi
[Observing that $(X \vee Y)_P \approx X_P \vee Y_P$, identify $X \# Y$ with the pointed mapping cone 
$X \ov{\#} Y$ of the inclusion $X \vee Y \ra X \times Y$ (cf. $\S 3$, Proposition 23).]\\
\endgroup %%------------------------------------<<

Every nilpotent group \mG is separable, i.e., the arrow $G \ra \prod\limits_p G_p$ is injective.  
The following result is the homotopy theoretic analog.\\

\begin{proposition} \ %09
Let \mX be a pointed nilpotent CW space $-$then for any pointed finite connected CW complex \mK, the arrow 
$[K,X] \ra \prod\limits_p [K,X_p]$ is injective.
\end{proposition}

[The assertion is certainly true if \mK is a finite wedge of circles.  Arguing inductively, consider the pushout square
\begin{tikzcd}%[sep=large]
{\bS^{n-1}} \ar{d} \ar{r}{f} &{L} \ar{d}\\
{\bD^{n}} \ar{r} &{K}
\end{tikzcd}
$(n \geq 2)$ and suppose that the arrow
%%----------------------------------------------------------------------------------------------13
$[L,X] \ra \prod\limits_p [L,X_p]$ is injective.  
Taking $f$ skeletal, there is a factorization
\begin{tikzcd}%[sep=large]
{\bS^{n-1}} \ar{d}[swap]{i} \ar{r}{f} &{L}\\
{M_f} \ar{ru}[swap]{r}
\end{tikzcd}
, 
where $L \approx M_f$ and $K \approx C_f \approx C_i$, so one can assume that $f$ is a closed cofibration.  
Restoring the base points, the corresponding arrow of restriction $f^*:C(L,l_0;X,x_0) \ra C(\bS^{n-1},s_{n-1};X,x_0)$ is then a Hurewicz fibration 
(cf. p. \pageref{9.33}) and the fiber of $f^*$ over 0 is homeomorphic to 
$C(L/\bS^{n-1},*_{\bS^{n-1}};X,x_0)$, $*_{\bS^{n-1}}$ the image of $\bS^{n-1}$ in $L/\bS^{n-1}$.  
But the projection 
$C_f \ra L/\bS^{n-1}$ is a pointed homotopy equivalence 
(cf. p. \pageref{9.34}), thus 
$C(K,k_0;X,x_0) \approx$ 
$C(L/\bS^{n-1},*_{\bS^{n-1}};X,x_0)$ 
(cf. p. \pageref{9.35}).  This said, given $\phi \in C(K,k_0;X,x_0)$, put 
$\psi = \restr{\phi}{L}$, let 
$
\begin{cases}
\ (C,\phi) = C(K,k_0;X,x_0:\phi)\\
\ (C,\psi) = C(L,l_0;X,x_0:\psi)
\end{cases}
$
and call 
$
\begin{cases}
\ [K,X]_\phi\\
\ [L,X]_\psi
\end{cases}
$
the pointed set
$
\begin{cases}
\ [K,X]\\
\ [L,X]
\end{cases}
$
with 
$
\begin{cases}
\ \phi\\
\ \psi
\end{cases}
$
as the base point.  Noting that 
$\pi_1(C(\bS^{n-1},s_{n-1};X,x_0),0) \approx$ 
$\pi_n(X)$, a portion of the homotopy sequence of our fibration reads: 
$\pi_1(C,\psi) \ra$ 
$\pi_n(X) \ra$ 
$[K,X]_\phi \ra$ 
$[L,X]_\psi$.  Here, $\pi_n(X)$ operates on $[K,X]_\phi$ and the orbit of $\phi$ consists of those maps which are pointed homotopic to $\psi$ when restricted to \mL, the stabilizer of $\phi$ being precisely $\im\pi_1(C,\psi)$.  Collect the data and display it in a commutative diagram
\[
\begin{tikzcd}%[sep=large]
{\pi_1(C,\psi)} \ar{d} \ar{r} 
&{\pi_n(X)} \ar{d} \ar{r} 
&{[K,X]_\phi} \ar{d} \ar{r} 
&{[L,X]_\psi}\ar{d}\\
{\prod\limits_p\pi_1(C_p,l_p \circx \psi)} \ar{r} 
&{\prod\limits_p\pi_n(X_p)} \ar{r}
&{\prod\limits_p [K,X_p]_{l_p \circx \phi}} \ar{r}
&{\prod\limits_p [L,X_p]_{l_p \circx \psi}} 
\end{tikzcd}
.
\]
The components of the first and second vertical arrows are $p$-localizing and by hypothesis, the fourth vertical arrow is injective.  
As for the third vertical arrow, its injectivity amounts to showing that if 
$\phi^\prime:K \ra X$ and if $\forall \ p$, $l_p \circx \phi^\prime \simeq l_p \circx \phi$, then $\phi^\prime \simeq \phi$.  
To begin, $\forall \ p$, 
$l_p \circx \psi^\prime \simeq l_p \circx \psi$ $\implies$ $\psi^\prime \simeq \psi$, hence $\phi^\prime$ lies on the 
$\pi_n(X)$-orbit of $\phi$, i.e., $\exists$! 
$\alpha \in \pi_n(X)/\im\pi_1(C,\psi)$: $[\phi^\prime] = \alpha \cdot [\phi]$.  
Claim: $\alpha$ is trivial.  In fact, $\forall \ p$, $l_p(\alpha)$ is trivial in 
$\pi_n(X_p)/\im\pi_1(C_p,l_p \circx \psi)$ and the arrow 
$\pi_n(X)/\im\pi_1(C,\psi) \ra $ 
$\prod\limits_p(\pi_n(X_p)/\im\pi_1(C_p,l_p \circx \psi))$ is one-to-one.]\\

Application: Let \mK be a pointed finite nilpotent CW complex; let \mX be a pointed nilpotent CW complex.  Suppose that 
$f,g:K \ra X$ are pointed continuous functions.  Assume: $\forall \ p$, $f_p \simeq g_p$ $-$then $f \simeq g$.\\

\begingroup%%----------------------------------->>
\fontsize{9pt}{11pt}\selectfont
\textbf{\small EXAMPLE} \  
Suppose that $P \neq \emptyset$ $\&$ $\ov{P} \neq \emptyset$.  Define \mK by the pushout square
\begin{tikzcd}[sep=large]
{\bS^n} \ar{d} \ar{r}{f} &{\bS_P^n \vee \bS_{\ov{P}}^n} \ar{d}\\
{\bD^{n+1}} \ar{r} &{K}
\end{tikzcd}
$(n \geq 2)$, where 
$f = (1,1) \in \pi_n(\bS_P^n \vee \bS_{\ov{P}}^n) \approx$ 
$\Z_P \oplus \Z_{\ov{P}}$.  
Let $\phi:K \ra \bS^{n+1}$ be the collapsing map $-$then $\forall \ p$, $l_p \circx \phi \simeq 0$ but $[\phi] \neq [0]$.  
Therefore, even when \mX is a sphere, Proposition 9 can fail if \mK is not finite (but Proposition 9 does imply that 
$\phi \in \Ph(K,\bS^{n+1})$).\\
\endgroup %%------------------------------------<<

%%----------------------------------------------------------------------------------------------14
\begingroup%%----------------------------------->>
\fontsize{9pt}{11pt}\selectfont
\textbf{\small FACT} \  
Let \mX be a pointed nilpotent CW space $-$then for any pointed finite connected CW complex \mK, the commutative diagram 
\begin{tikzcd}[sep=large]
{[K,X]} \ar{d} \ar{r} &{[K,X_{\ov{P}}]} \ar{d}\\
{[K,X_P]} \ar{r} &{[K,X_{\Q}]}
\end{tikzcd}
is a pullback square in $\bSET_*$.
\\ \indent
[Note: \ \mX ``is'' the double mapping track of the pointed 2-sink $X_P \ra X_{\Q} \la X_{\ov{P}}$.]\\
\endgroup %%------------------------------------<<

\begingroup%%----------------------------------->>
\fontsize{9pt}{11pt}\selectfont
\textbf{\small EXAMPLE} \  
The assumption on \mK plays a role in the preceding result.  
Thus suppose that 
$P \neq \emptyset$ $\&$ $\ov{P} \neq \emptyset$ $-$then the commutative diagram 
\begin{tikzcd}[sep=large]
{[\bP^\infty(\C),\bS^3]} \ar{d} \ar{r} &{[\bP^\infty(\C),\bS_{\ov{P}}^3]} \ar{d}\\
{[\bP^\infty(\C),\bS_P^3]} \ar{r} &{[\bP^\infty(\C),\bS_{\Q}^3]}
\end{tikzcd}
is not a pullback square in $\bSET_*$.
\\ \indent
[Show that the arrow 
$\lim^1[\Sigma \bP^n(\C),\bS^3] \ra$ 
$\lim^1[\Sigma \bP^n(\C),\bS_P^3] \oplus \lim^1[\Sigma \bP^n(\C),\bS_{\ov{P}}^3]$ 
is not one-to-one 
(cf. p. \pageref{9.36}).]\\
\endgroup %%------------------------------------<<

\begingroup%%----------------------------------->>
\fontsize{9pt}{11pt}\selectfont
\textbf{\small FACT} \  
Let \mX be a pointed nilpotent CW space $-$then for any finite CW complex \mK, the arrow 
$[K,X] \ra \ds\prod\limits_p [K,X_p]$ is injective.
\\ \indent
[Note: \ In this context, the brackets refer to homotopy classes of maps, not to pointed homotopy classes of pointed maps.]\\
\endgroup %%------------------------------------<<

Let \mX be a pointed nilpotent CW space $-$then one may attach to \mX a sink $\{r_p:X_p \ra X_{\Q}\}$ and a source 
$\{l_p:X \ra X_p\}$, where $\forall$ 
$
\begin{cases}
\ p\\
\ q
\end{cases}
\hspace{-.25cm}
, \  r_p \circx l_p \simeq r_q \circx l_q.
$
\\

\begin{proposition} \ %10
Let \mX be a pointed nilpotent CW space with finitely generated homotopy groups.  
Suppose given a pointed finite connected CW complex \mK and pointed continuous functions $\phi(p):K \ra X_p$ such that $\forall$ 
$
\begin{cases}
\ p\\
\ q
\end{cases}
\hspace{-.25cm}
,
$
$r_p \circx \phi(p) \simeq r_q \circx \phi(q)$ 
$-$then there is a pointed continuous function $\phi:K \ra X$ such that $\forall \ p$, $l_p \circx \phi \simeq \phi(p)$.
\end{proposition}

[The fracture lemma on 
p. \pageref{9.37} implies that the result holds if \mK is a finite wedge of circles.  
Proceeding via induction, consider the pushout square
\begin{tikzcd}%[sep=large]
{\bS^{n-1}} \ar{d} \ar{r}{f} &{L} \ar{d}\\
{\bD^{n}} \ar{r} &{K}
\end{tikzcd}
$(n \geq 2)$ and assume that there is a pointed continuous function $\psi:L \ra X$ such that $\forall \ p$, 
$l_p \circx \psi \simeq \psi(p)$, where $\psi(p) = \restr{\phi(p)}{L}$.  
Since $\forall \ p$, $\psi(p) \circx f \simeq 0$, from Proposition 9, 
$\psi \circx f \simeq 0$, so $\exists$ a pointed continuous function $\phi^\prime:K \ra X$ which restricts to $\psi$.  
Taking $f$ to be a closed cofibration and following the proof of Proposition 9, form the commutative diagram
\[
\begin{tikzcd}%[sep=large]
{\pi_1(C,\psi)} \ar{d} \ar{r} 
&{\pi_n(X)} \ar{d} \ar{r} 
&{[K,X]_{\phi^\prime}} \ar{d} \ar{r} 
&{[L,X]_\psi}\ar{d}\\
{\pi_1(C_p,l_p \circx \psi)} \ar{r} 
&{\pi_n(X_p)} \ar{r}
&{[K,X_p]_{l_p \circx \phi^\prime}} \ar{r}
&{[L,X_p]_{l_p \circx \psi}} 
\end{tikzcd}
.
\]
%%----------------------------------------------------------------------------------------------15
Because $\restr{\phi(p)}{L} \simeq \l_p \circx \restr{\phi^\prime}{L}$, $\phi(p)$ 
must be on the $\pi_n(X_p)$-orbit of $l_p \circx \phi^\prime$, 
i.e., $\exists$! $\alpha(p) \in \pi_n(X_p)/\im\pi_1(C_p,l_p \circx \psi)$: $[\phi(p)] = \alpha(p) \cdot [l_p \circx \phi^\prime]$.  
However, the $\alpha(p)$ all rationalize to the same element of 
$(\pi_n(X)/\im\pi_1(C,\psi))_{\Q}$, thus $\exists$! $\alpha \in \pi_n(X)/\im\pi_1(C,\psi)$: $\forall \ p$, 
$l_p(\alpha) = \alpha(p)$.  Put $\phi = \alpha \cdot \phi^\prime$: 
$l_p \circx \phi \simeq$ 
$l_p \circx (\alpha \cdot \phi^\prime) \simeq$ 
$l_p(\alpha) \cdot (l_p \circx \phi^\prime) \simeq$ 
$\alpha(p) \cdot (l_p \circx \phi^\prime) \simeq$ 
$\phi(p)$.]\\

\begingroup%%----------------------------------->>
\fontsize{9pt}{11pt}\selectfont
\textbf{\small FACT} \  
Let \mX be a pointed nilpotent CW space with finitely generated homotopy groups.  Suppose given a finite CW complex \mK and pointed continuous functions $\phi(p):K \ra X_p$ such that $\forall$ 
$
\begin{cases}
\ p\\
\ q
\end{cases}
, \ 
$
$r_p \circx \phi(p) \simeq r_q \circx \phi(q)$ 
$-$then there is a continuous function $\phi:K \ra X$ such that $\forall \ p$, $l_p \circx \phi \simeq \phi(p)$.\\ 
\endgroup %%------------------------------------<<

\index{Hasse Principle}
\textbf{\small HASSE PRINCIPLE} \quad 
Let \mX be a pointed nilpotent CW space with finitely generated homotopy groups $-$then for any pointed finite conncected CW complex \mK, the source
$\{[K,X] \ra [K,X_p]\}$ is the multiple pullback of the sink 
$\{[K,X_p] \ra [K,X_{\Q}]\}$.

[This is a consequence of Propositions 9 and 10.]\\

\begingroup%%----------------------------------->>
\fontsize{9pt}{11pt}\selectfont
Given a pointed nilpotent CW space \mX with finitely generated homotopy groups, the 
\un{genus}
\index{genus (of a pointed nilpotent CW space with finitely generated homotopy groups)} 
$\gen X$ of \mX is the conglomerate of pointed homotopy types $[Y]$, 
where \mY is a pointed nilpotent CW space with finitely generated homotopy groups such that 
$\forall \ p$, $X_p \approx Y_p$.  
The members of $\gen X$ have isomorphic higher homotopy groups (but their fundamental groups are not necessarily isomorphic) 
and isomorphic integral singular homology groups (but their integral singular cohomology rings are not necessarily isomorphic) .  
\\ \indent
Examples: 
(1) $\gen \bS^n = \{[\bS^n]\}$; 
(2) $\gen K(\pi,n) = \{[K(\pi,n)]\}$, $\pi$ a finitely generated abelian group; 
(3) $\gen M(\pi,n) = \{[M(\pi,n)]\}$, $\pi$ a finitely generated abelian group $(n \geq 2)$.\\
\endgroup %%------------------------------------<<

\begingroup%%----------------------------------->>
\fontsize{9pt}{11pt}\selectfont
\textbf{\small EXAMPLE} \  
Fix a generator $\alpha \in \pi_6(\bS^3) \approx \Z/12\Z$.  Put 
$X = \bD^7 \sqcup_\alpha \bS^3$, 
$Y = \bD^7 \sqcup_{5\alpha} \bS^3$ 
$-$then $\forall \ p$, $X_p \approx Y_p$ but \mX and \mY do not have the same pointed homotopy type.\\
\endgroup %%------------------------------------<<

\begingroup%%----------------------------------->>
\fontsize{9pt}{11pt}\selectfont
\textbf{\small EXAMPLE} \  
It has been shown by Wilkerson\footnote[2]{\textit{Topology} \textbf{15} (1976), 111-130.} 
that if \mX is a pointed finite simply connected CW complex, then $\#(\gen X) < \omega$ but this can fail when \mX is not finite.  
For instance, take $X = \bP^\infty(\HH)$ $-$then $\gen X$ is in a one-to-one correspondence with the set of all functions 
$\bPi \ra \{\pm 1\}$ 
(Rector\footnote[3]{\textit{SLN} \textbf{249} (1971), 99-105.}), hence has cardinality $2^\omega$.
\\ \indent
[Note: \ It is unknown whether $\#(\gen X) < \omega$ for an arbitrary pointed finite nilpotent CW complex \mX.]\\
\endgroup %%------------------------------------<<

\begingroup%%----------------------------------->>
\fontsize{9pt}{11pt}\selectfont
\textbf{\small EXAMPLE} \  
Let 
$
\begin{cases}
\ X\\
\ Y
\end{cases}
$
be pointed nilpotent CW spaces $-$then \mX and \mY are said to be 
\un{clones}
\index{clones (CW spaces)} 
if 
(i) $\forall \ n$, $X[n] \approx Y[n]$ and 
(ii) $\forall \ p$, $X_p \approx Y_p$. 
While neither (i) nor (ii) alone suffices to imply that $X \approx Y$, one
%%----------------------------------------------------------------------------------------------16
can ask whether this is the case of their conjunction.  
In other words, if \mX and \mY are clones, does it follow that \mX and \mY have the same pointed homotopy type?  The answer is ``no''.  
Take $X = \bS^3 \times K(\Z,3)$ $-$then, up to pointed homotopy type, the number of distinct clones of \mX is uncountable 
(McGibbon\footnote[6]{\textit{Comment. Math. Helv.} \textbf{68} (1993), 263-277.}).\\
\endgroup %%------------------------------------<<

Given a set of primes \mP, a pointed CW space \mX is said to be 
\un{$P$-local}
\index{P-local (pointed CW space )} 
if $\forall \ n \in S_p$, the arrow
$
\begin{cases}
\ \Omega X \ra \Omega X\\
\ \sigma \ra \sigma^n
\end{cases}
$
is a pointed homotopy equivalence.

[Note: \ \mX is $P$-local iff $\pi_1(X)$ and the $\pi_q(X) \rtimes \pi_1(X)$ $(q \geq 2)$ are $P$-local groups 
(cf. p. \pageref{9.38}) or still, iff $\pi_1(X)$ is a $P$-local group and the $\pi_q(X)$ $(q \geq 2)$ are $P$-local 
$\pi_1(X)$-modules 
(cf. p. \pageref{9.39}).  
Therefore a $P$-local space is $P$-local in homotopy (but not conversely 
(cf. p. \pageref{9.40}).]

Example: For any $P$-local group \mG, $K(G,1)$ is a $P$-local space.

[Note: \ Accordingly, a $P$-local space is not necessarily $P$-local in homology 
(cf. p. \pageref{9.41}).]

Notation: $\bCONCWSP_{*,P}$
\index{$\bCONCWSP_{*,P}$} 
is the full subcategory of $\bCONCWSP_{*}$ whose objects are the pointed connected CW spaces which are $P$-local and 
$\bHCONCWSP_{*,P}$
\index{$\bHCONCWSP_{*,P}$} 
is the associated homotopy category.

[Note: \ This notation is a consistent extension of that introduced on 
p. \pageref{9.42} for the nilpotent category, i.e., a pointed nilpotent CW space which is $P$-local in homotopy is $P$-local 
(cf. p. \pageref{9.43}).]

Observation: Set 
$\bS_T^q = \bS^1$ $(q = 1)$, 
$\bS_T^q = (\bS^{q-1} \amalg *)\#\bS^1$ $(q \geq 2)$ 
and let 
$\rho_n^q = \rho_n$ $(q = 1)$, 
$\rho_n^q = \id\#\rho_n$ $(q \geq 2)$, 
where $\rho_n:\bS^1 \ra \bS^1$ is a map of degree $n$ $(n \in S_P)$.  
Working in $\bHCONCWSP_*$, put 
$S_0 = \{[\rho_n^q]\}$ $-$then $S_0^\perp$ is the object class of $\bHCONCWSP_{*,P}$.

[In fact, 
$[\bS_T^1,X] \approx \pi_1(X)$, 
$[\bS_T^q,X] \approx \pi_q(X) \rtimes \pi_1(X)$ $(q \geq 2)$ and 
$(\rho_n^q)^*:[\bS_T^q,X] \ra [\bS_T^q,X]$ is the $n^\text{th}$ power map $\forall \ q \geq 1$.]

Let $[f]:X \ra Y$ be a morphism in $\bHCONCWSP_{*}$ $-$then $\abs{f}$ (or $f$) is said to be a 
\un{$P$-equivalence}
\index{P-equivalence} 
if $\abs{f}$ is orthogonal to every $P$-local pointed connected CW space.

[Note: \ This terminology does not conflict with that used earlier in the nilpotent category (cf. Proposition 12).]

Convention: Given a pointed connected CW space \mX, a 
\un{$P[X]$-module} 
\index{P[X]-module} 
is a $P[\pi_1(X)]$-module.

[Note: \ If 
$
\begin{cases}
\ X\\
\ Y
\end{cases}
$
are pointed connected CW spaces, and if $f:X \ra Y$ is a pointed continuous function, then every $P[Y]$-module can be construed as a $P[X]$-module 
(cf. p. \pageref{9.44}).]\\
%%----------------------------------------------------------------------------------------------17
\begin{proposition} \ %11
Let 
$
\begin{cases}
\ X\\
\ Y
\end{cases}
$
be pointed connected CW spaces, $f:X \ra Y$ a pointed continuous function $-$then $f$ is a $P$-equivalence iff 
$\pi_1(f)_P:\pi_1(X)_P \ra \pi_1(Y)_P$ is bijective and for every locally constant coefficient system $\sG$ on \mY arising from a 
$P[Y]$-module, $H^n(Y;\sG) \approx H^n(X;f^*\sG)$ $\forall \ n$.
\end{proposition}

[Necessity: \ Given a $P$-local group \mG, 
$[f] \perp K(G,1)$ $\implies$
$[Y,K(G,1)] \approx [X,K(G,1)]$ $\implies$
$\Hom(\pi_1(Y),G) \approx \Hom(\pi_1(X),G)$ $\implies$
$\pi_1(f)  \perp G$ $\implies$ 
$\pi_1(f)_P:\pi_1(X)_P \approx \pi_1(Y)_P$.  
To check the cohomological assertion, fix a right $P[Y]$-module $\pi$ and let 
$\chi:\pi_1(Y)_P \ra \Aut\pi$ be the associated homomorphism.  
Denote by $\sG:\Pi Y \ra \bAB$ the cofunctor corresponding to the composite $\chi \circx l_P$, where $l_P:\pi_1(Y) \ra \pi_1(Y)_P$.  
Since for positive $n$, $K(\pi,n;\chi)$ is $P$-local, 
$[f] \perp K(\pi,n;\chi)$ $\implies$
$[Y,K(\pi,n;\chi)] \approx [X,K(\pi,n;\chi)]$ $\implies$
$H^n(Y;\sG) \approx H^n(X;f^*\sG)$  
(cf. p. \pageref{9.45}), $n > 0$.  
There remains the claim that 
$H^0(Y;\sG) \approx H^0(X;f^*\sG)$, i.e., that the $\pi_1(Y)$-invariants in $\pi$ equal the $\pi_1(X)$-invariants in $\pi$.  
To see this, consider the commutative diagram
\begin{tikzcd}%[sep=large]
{\pi_1(X)} \ar{d} \ar{r} &{\pi_1(Y)} \ar{d}\\
{\pi_1(X)_P} \ar{r} &{\pi_1(Y)_P}
\end{tikzcd}
.  
From what has been said above, the arrow $\pi_1(X)_P \ra \pi_1(Y)_P$ is an isomorphism.  
The claim thus follows from the fact that the $\pi_1(Y)_P$-invariants in $\pi$ are equal to the the $\pi_1(Y)$-invariants in $\pi$ 
(cf. p. \pageref{9.46}).

Sufficiency: \ In order to apply the machinery of full blown obstruction theory (locally constant coefficients 
(Olum\footnote[2]{\textit{Ann. of Math.} \textbf{52} (1950), 1-50; see also Baues, \textit{Obstruction Theory}, Springer Verlag (1977).})), take
$
\begin{cases}
\ X\\
\ Y
\end{cases}
$
to be pointed connected CW complexes with \mX a pointed subcomplex of \mY, so $f$ is the inclusion $X \ra Y$,
Fix a pointed continuous function $\phi:X \ra Z$, where \mZ is $P$-local $-$then $\pi_1(f) \perp \pi_1(Z)$ $\implies$ $\exists$! 
$\theta \in \Hom(\pi_1(Y),\pi_1(Z))$: $\pi_1(\phi) = \theta \circx \pi_1(f)$.  
By restriction of scalars, i.e., using the filler for
\begin{tikzcd}%[sep=large]
{\pi_1(Y)} \ar{d} \ar{r}{\theta} &{\pi_1(Z)}\\
{\pi_1(Y)_P} \ar[dashed]{ru}
\end{tikzcd}
, the $\pi_n(Z)$ $(n \geq 2)$ become $P[Y]$-modules and there is a long exact sequence 
$H^1(Y;\pi_n(Z)) \ra$ 
$H^1(X;\pi_n(Z)) \ra$
$H^2(Y,X;\pi_n(Z)) \ra$
$H^2(Y;\pi_n(Z)) \ra$
$H^2(X;\pi_n(Z)) \ra$
$\cdots$.

\indent\indent (Existence) \ One can find a pointed continuous function $\psi:(Y,X)^{(2)} \ra Z$ such that 
$\restr{\psi}{X} = \phi$ and 
$\pi_1(\psi) = \theta$ ($(Y,X)^{(2)} = Y^{(2)} \cup X$ and 
$\pi_1((Y,X)^{(2)}) \approx \pi_1(Y)$).  
On the other hand, the higher order obstructions to the existence of a pointed continuous function 
$\Phi:Y \ra Z$ such that $\restr{\Phi}{X} = \phi$ ($\implies$ $\pi_1(\Phi) = \theta$) lie in the 
$H^{n+1}(Y,X;\pi_n(Z))$ $(n \geq 2)$.  
As these groups necessarily vanish, the precomposition arrow 
$f^*:[Y,Z] \ra [X,Z]$ is surjective.

\indent\indent (Uniqueness) \ Suppose that $\Phi^\prime, \Phi\pp:Y \ra Z$ are pointed continuous functions
%%----------------------------------------------------------------------------------------------18
with 
$
\begin{cases}
\ \restr{\Phi^\prime}{X} = \phi\\
\ \restr{\Phi\pp}{X} = \phi
\end{cases}
$
$-$then the claim is that $\Phi^\prime$ and $\Phi\pp$ are pointed homotopic.  
Indeed, 
$\pi_1(\Phi^\prime) = \theta = \pi_1(\Phi\pp)$ $\implies$ 
$\restr{\Phi^\prime}{(Y,X)^{(1)}} \simeq {\Phi\pp}{(Y,X)^{(1)}}$ $\rel X$ and since the 
$H^n(Y,X;\pi_n(Z))$ $(n \geq 2)$ are trivial, $\Phi^\prime$ and $\Phi\pp$ are homotopic $\rel X$.]\\

\textbf{\small LEMMA}  \  
Let 
$
\begin{cases}
\ X\\
\ Y
\end{cases}
$
be pointed connected CW spaces, $f:X \ra Y$ a pointed continuous function.  
Fix a group \mG and a ring \mA with unit.  
Suppose given a homomorphism $\pi_1(Y) \ra G$ and a homomorphism $\Z[G] \ra A$.  
Let $\sA$ be the locally constant coefficient system on \mY corresponding to \mA.  
Assume: $\forall \ n \geq 0$, $H_n(X;f^*\sA) \approx H_n(Y;\sA)$ $-$then for every locally constant coefficient system $\sM$ on \mY corresponding to a 
$
\begin{cases}
\ \text{right}\\
\ \text{left}
\end{cases}
$
$A$-module \mM, 
$
\begin{cases}
\ H_n(X;f^*\sM) \approx H_n(Y;\sM)\\
\ H^n(Y;\sM) \approx H^n(X;f^*\sM)
\end{cases}
\forall \ n \geq 0.
$

[It suffices to work with pointed connected CW complexes \mX and \mY, where \mX is a pointed subcomplex of \mY 
($f$ becoming the inclusion).  
Put $\pi = \pi_1(Y)$ and let $C_*(\widetilde{Y},\widetilde{X})$ be the associated relative skeletal chain complex 
(Whitehead\footnote[2]{\textit{Elements of Homotopy Theory}, Springer Verlag (1978), 287-288.}), 
so each $C_n(\widetilde{Y},\widetilde{X})$ is a free left $\Z[\pi]$-module and $\forall \ n \geq 0$, 
$H_n(Y,X;\sA) = H_n(A \otimes_{\Z[\pi]} C_*(\widetilde{Y},\widetilde{X}))$.  
Here, however, $\forall \ n \geq 0$, $H_n(Y,X;\sA) = 0$, and this means that 
$A \otimes_{\Z[\pi]} C_*(\widetilde{Y},\widetilde{X})$ is a free resolution of 0 as an $A$-module.
Therefore, for any right $A$-module \mM, 
$H_n(Y,X;\sM) \approx$ 
$H_n(M \otimes_{\Z[\pi]} C_*(\widetilde{Y},\widetilde{X})) \approx$
$H_n(M \otimes_{A} A \otimes_{\Z[\pi]} C_*(\widetilde{Y},\widetilde{X})) \approx$ 
$\Tor_n^A(M,0) = 0$ $\forall \ n \geq 0$ and for any left $A$-module \mM, 
$H^n(Y,X;\sM) \approx$ 
$H^n(\Hom_{\Z[\pi]} (C_*(\widetilde{Y},\widetilde{X}),M)) \approx$
$H^n(\Hom_{\Z[\pi]} (C_*(\widetilde{Y},\widetilde{X}),\Hom_A(A,M))) \approx$
$H^n(\Hom_A(A \otimes_{\Z[\pi]} C_*(\widetilde{Y},\widetilde{X}),M)) \approx$ 
$\Ext_A^n(0,M) = 0$ $\forall \ n \geq 0$.]

[Note: \ Recall that when dealing with modules over a group ring, there is no essential distinction between ``left'' and 
``right''.  In particular: The $C_n(\widetilde{Y},\widetilde{X})$ are both left and right free $\Z[\pi]$-modules.]\\

It is a corollary that $f$ is a $P$-equivalence provided that 
$\pi_1(f)_P:\pi_1(X)_P \ra \pi_1(Y)_P$ is bijective and for every locally constant coefficient system $\sG$ on \mY arising from a 
$P[Y]$-module, 
$H_n(X;f^*\sG) \approx H_n(Y;\sG)$ $\forall \ n \geq 0$.  In fact, to pass from homology to cohomology, one may apply the lemma, taking $G = \pi_1(Y)_P$ and $A = (\Z[G])_{S_P}$ 
(cf. p. \pageref{9.47}).\\

\begingroup%%----------------------------------->>
\fontsize{9pt}{11pt}\selectfont
\textbf{\small EXAMPLE} \  
Let \mX be a pointed connected CW space.  Suppose that \mN is a perfect normal subgroup of $\pi_1(X)$ which is contained in the kernel of the arrow of localization $\pi_1(X) \ra \pi_1(X)_P$ $-$then $f_N^+:X \ra X_N^+$ is a $P$-equivalence.
\\ \indent
%%----------------------------------------------------------------------------------------------19
[The assumption on \mN guarantees that $\pi_1(X)_P \approx \pi_1(X_N^+)_P$.  But $f_N^+$ is acyclic, so for every locally constant coefficient system $\sG$ on $X_N^+$, 
$H_*(X;(f_N^+)^*\sG) \approx H_*(X_N^+;\sG)$ (cf. $\S 5$, Proposition 22) and the lemma can be quoted.]
\\ \indent
[Note: \ It is not really necessary to use the lemma.  This is because acyclic maps can equally well be characterized in terms of cohomology with locally constant coefficients.]\\
\endgroup %%------------------------------------<<

\begin{proposition} \ %12
Let 
$
\begin{cases}
\ X\\
\ Y
\end{cases}
$
be pointed nilpotent CW spaces, $f:X \ra Y$ a pointed continuous function.  Assume 
$f_*:H_*(X;\Z_P) \ra H_*(Y;\Z_P)$ is an isomorphism $-$then for every locally constant coefficient system $\sG$ on \mY arising from a $P[Y]$-module, $H_n(X;f^*\sG) \approx H_n(Y;\sG)$ $\forall \ n \geq 0$.
\end{proposition}

[According to Proposition 5, $f_P: X_P \ra Y_P$ is a pointed homotopy equivalence, so there is no loss of generality in supposing that $Y = X_P$, $f = l_P$.  Consider the diagram
\begin{tikzcd}%[sep=large]
{\widetilde{X}} \ar{d}[swap]{\widetilde{f}} \ar{r}{p} 
&{X} \ar{r} \ar{d}
&{K(\pi_1(X),1)} \ar{d}\\
{\widetilde{Y}} \ar{r}[swap]{q}
&{Y} \ar{r}
&{K(\pi_1(Y),1)}
\end{tikzcd}
.  It commutes up to pointed homotopy and because 
$
\begin{cases}
\ \widetilde{X}\\
\ \widetilde{Y}
\end{cases}
$
are simply connected, 
$H_n(\widetilde{X};p^*f^*\sG) \approx$ 
$H_n(\widetilde{X};\widetilde{f}^*q^*\sG) \approx$
$H_n(\widetilde{X};G)$, 
$H_n(\widetilde{Y};q^*\sG) \approx$
$H_n(\widetilde{Y};G)$, \mG the underlying $P$-local $\pi_1(Y)$-module.  
Bearing in mind that \mG is, in particular, a $\Z_P$-module, the fact that
$H_*(\widetilde{X};\Z_P) \approx H_*(\widetilde{Y};\Z_P)$, 
in conjunction with the universal coefficient theorem then gives
$H_*(\widetilde{X};G) \approx H_*(\widetilde{Y};G)$.  
Pass now to the morphism 
$\{E_{p,q}^2 \approx H_p(\pi_1(X);H_q(\widetilde{X};G))\} \ra$ 
$\{\ov{E}_{p,q}^2 \approx H_p(\pi_1(Y);H_q(\widetilde{Y};G))\}$ 
of fibration spectral sequences.  Since the action of 
$\pi_1(Y)$ on the 
$H_q(\widetilde{Y})$ is nilpotent (cf. $\S 5$, Proposition 17), each 
$H_q(\widetilde{Y};G)$ is a $P$-local $\pi_1(Y)$-module 
(cf. p. \pageref{9.48}), i.e., is a $P[X]$-module 
$(\pi_1(X)_P = \pi_1(Y))$.  
Therefore, $\forall \ p$ $\&$ $\forall \ q$, 
$H_p(\pi_1(X);H_q(\widetilde{X};G)) \approx$ 
$H_p(\pi_1(Y);H_q(\widetilde{Y};G))$ (cf. $\S 8$, Proposition 16), which serves to complete the proof  
(cf. p. \pageref{9.49}).]

[Note: \ In the nilpotent category, the term ``$P$-equivalence'' has two possible interpretations.  
The point of the proposition is that they coincide (cf. $\S 8$, Proposition 2).]\\

If \mS is the class of $P$-equivalences and if \mD is the class of $P$-local spaces, then $(S,D)$ is an orthogonal pair.  
Proof: $S = D^\perp$ (by definition) and 
$S_0^\perp = D$ $\implies$ 
$S_0^{\perp\perp} = S$ $\implies$ 
$D = S_0^{\perp\perp\perp}  = S^\perp$.  
Consequently, \mS has the closure properties (1)-(3) formulated on 
p. \pageref{9.50}.  
It will also be necessary to know the interplay between $P$-equivalences, wedges, and certain weak colimits.

\indent\indent (Wedges) \quad
Let 
$
\begin{cases}
\ X_i\\
\ Y_i
\end{cases}
(i \in I)
$
be pointed connected CW spaces.  Suppose that $\forall \ i$, $f_i:X_i \ra Y_i$ is a $P$-equivalence $-$then
$\bigvee\limits_i f_i:\bigvee\limits_i X_i \ra \bigvee\limits_i Y_i$ is a $P$-equivalence.

%%----------------------------------------------------------------------------------------------20
[By assumption, $\forall \ i$, 
$(\pi_1(X_i) \ra \pi_1(Y_i)) \perp \Ob\bGR_P$, hence 
$(\underset{i}{*} \pi_1(X_i) \ra \underset{i}{*} \pi_1(Y_i))$ $\perp \Ob\bGR_P$, i.e., 
$(\pi_1(\bigvee\limits_i X_i) \ra \pi_1(\bigvee\limits_i Y_i))$ $\perp \Ob\bGR_P$.  Let $\sG$ be a locally constant coefficient system on $\bigvee\limits_i Y_i$ arising from a $P[\bigvee\limits_i Y_i]$-module.  Employing the notation used in the proof of Proposition 11, 
$[\bigvee\limits_i Y_i, K(\pi,n;\chi)] \approx$ 
$\prod\limits_i [Y_i, K(\pi,n;\chi)] \approx$
$\prod\limits_i [X_i, K(\pi,n;\chi)] \approx$
$[\bigvee\limits_i X_i, K(\pi,n;\chi)]$ 
$\implies$
$H^n(\bigvee\limits_i Y_i;\sG) \approx$ 
$H^n(\bigvee\limits_i X_i;(\bigvee\limits_i f_i)^*\sG)$ 
(cf. p. \pageref{9.51}), $n > 0$. 
Finally, the 
$\pi_1(\bigvee\limits_i Y_i)$-invariants in $\pi$ equal $\bigcap\limits_i \pi^{\pi_1(Y_i)}$ and the 
$\pi_1(\bigvee\limits_i X_i)$-invariants in $\pi$ equal $\bigcap\limits_i \pi^{\pi_1(X_i)}$.  
And: $\forall \ i$, $\pi^{\pi_1(Y_i)} = \pi^{\pi_1(X_i)}$.]

\indent\indent (Double Mapping Cylinders) \quad 
Let $X \overset{f}{\la} Z \overset{g}{\ra} Y$ be a pointed 2-source, where 
$
\begin{cases}
\ X\\
\ Y
\end{cases}
\& \ Z
$
are pointed connected CW spaces and $f$ is a $P$-equivalence.  Form the pointed double mapping cylinder $M_{f,g}$ of $f, g$ 
$-$then the arrow $Y \ra M_{f,g}$ is a $P$-equivalence.

[Assuming that 
$
\begin{cases}
\ X\\
\ Y
\end{cases}
\& \ Z
$
are pointed connected CW complexes and 
$
\begin{cases}
\ f\\
\ g
\end{cases}
$
are skeletal, pass from 
\begin{tikzcd}%[sep=large]
{Z} \ar{d}[swap]{f} \ar{r}{g} &{Y} \ar{d}\\
{X} \ar{r} &{M_{f,g}}
\end{tikzcd}
to 
\begin{tikzcd}%[sep=large]
{Z} \ar{d} \ar{r} &{M_g} \ar{d}{j}\\
{M_f} \ar{r}[swap]{i} &{M_{f,g}}
\end{tikzcd}
(cf. p. \pageref{9.52}), noting that the arrow $Z \ra M_f$ is a $P$-equivalence.  
Thanks to Van Kampen, the commutative diagram 
\begin{tikzcd}%[sep=large]
{\pi_1(Z)} \ar{d} \ar{r} &{\pi_1(M_g)} \ar{d}\\
{\pi_1(M_f)} \ar{r} &{\pi_1(M_{f,g})}
\end{tikzcd}
is a pushout square \bGR, so 
$(\pi_1(Z) \ra \pi_1(M_f))$ $\perp \Ob \bGR_P$ $\implies$ 
$(\pi_1(M_g) \ra \pi_1(M_{f,g}))$ $\perp \Ob \bGR_P$.  
Let $\sG$ be a locally constant coefficient system on $M_{f,g}$ arising from a $P[M_{f,g}]$-module.  
On general grounds (excision), 
$H^n(M_{f,g},M_g;\sG) \approx$ 
$H^n(M_f,\Z;\restr{\sG}{M_f})$ $\forall \ n \geq 0$, 
thus 
$H^n(M_{f,g},M_g;\sG) = 0$ $\forall \ n \geq 0$ $\implies$
$H^n(M_{f,g};\sG) \approx$ 
$H^n(M_{g};\restr{\sG }{M_g})$ $\forall \ n \geq 0$.  
That the arrow $M_g \ra M_{f,g}$ is a $P$-equivalence is therefore implied by Proposition 11.]

\indent\indent (Mapping Telescopes) \quad
Let $\{X_k,f_k\}$ be a sequence, where $X_k$ is a pointed connected CW space and $f_k:X_k \ra X_{k+1}$ is a $P$-equivalence.  
Form the pointed mapping telescope $\telx(\bX,\bff)$ of $(\bX,\bff)$ 
$-$then the arrow $X_0 \ra \telx(\bX,\bff)$ is a 
$P$-equivalence.

[Assuming that the $X_k$ are pointed connected CW complexes and the $f_k$ are skeletal, there is a commutative diagram
\begin{tikzcd}%[sep=large]
{\telx_k(\bX,\bff)} \ar{d} \ar{r} &{\telx_{k+1}(\bX,\bff)}\ar{d}\\
{X_k} \ar{r} &{X_{k+1}}
\end{tikzcd}
in which the vertical arrows are pointed homotopy equivalences 
(cf. p. \pageref{9.53}).  By hypothesis, $\forall \ k$, 
$(\pi_1(\telx_0(\bX,\bff)) \ra$
$\pi_1(\telx_k(\bX,\bff))) \perp \Ob\bGR_P$, so 
$(\pi_1(\telx_0(\bX,\bff)) \ra$
$\colimx \pi_1(\telx_k(\bX,\bff))) \perp \Ob\bGR_P$ $\implies$ 
$(\pi_1(\telx_0(\bX,\bff)) \ra$ 
$\pi_1(\telx(\bX,\bff)))$ 
$\perp \Ob\bGR_P$.  Let $\sG$ be a locally constant coefficient system 
%%----------------------------------------------------------------------------------------------21
on $\telx(\bX,\bff)$ arising from a $P[\telx(\bX,\bff)]$-module and put 
$\sG_k = \restr{\sG}{\telx_k(\bX,\bff)}$ 
$-$then $\forall \ n \geq 0$, 
$H^n(\telx_k(\bX,\bff);\sG_k) \approx$ 
$H^n(\telx_0(\bX,\bff);\sG_0)$ $\implies$ 
$\lim H^n(\telx_k(\bX,\bff);\sG_k) \approx$ 
$H^n(\telx_0(\bX,\bff);\sG_0)$.  
Since 
$\pi_1(\telx(\bX,\bff)) \approx$ 
$\colimx (\pi_1(\telx_k(\bX,\bff))$, 
$H^0(\telx(\bX,\bff);\sG) \approx$ 
$\lim H^0(\telx_k(\bX,\bff))$.  
Moreover, $\forall \ n \geq 1$, there is an exact sequence 
$0 \ra $ 
$\lim^1 H^{n-1}(\telx_k(\bX,\bff);\sG_k) \ra$
$H^n(\telx(\bX,\bff);\sG) \ra$ 
$\lim H^n(\telx_k(\bX,\bff);\sG_k) \ra 0$ 
of abelian groups 
(Whitehead\footnote[2]{\textit{Elements of Homotopy Theory}, Springer Verlag (1978), 273-274.}). 
 But here the $\lim^1$ terms vanish, so $\forall \ n \geq 1$, 
$H^n(\telx(\bX,\bff);\sG) \approx$ 
$\lim H^n(\telx_k(\bX,\bff); \sG_k)$.]\\

%\begingroup%%----------------------------------->>
%\fontsize{9pt}{11pt}\selectfont
\index{Theorem: Homotopical $P$-Localization Theorem}
\index{Homotopical $P$-Localization Theorem}
\textbf{\small HOMOTOPICAL $P$-LOCALIZATION THEOREM} \quad 
$\bHCONCWSP_{*,P}$ is a reflective subcategory of $\bHCONCWSP_*$.

[The theorem will follow provided that one can show that it is possible to assign to each pointed connected CW space \mX a $P$-local pointed connected CW space $X_P$ and a $P$-equivalence $l_P:X \ra X_P$.  Let $M_n^q$ be the pointed double mapping cylinder of the pointed 2-source 
$\bS_T^q \overset{\rho_n^q}{\lla} \bS_T^q \overset{\rho_n^q}{\lra} \bS_T^q$ $-$then the diagram
\begin{tikzcd}%[sep=large]
{\bS_T^q} \ar{d}[swap]{\rho_n^q} \ar{r}{\rho_n^q} &{\bS_T^q} \ar{d}{j_n^q}\\
{\bS_T^q} \ar{r}[swap]{i_n^q} &{M_n^q}
\end{tikzcd}
is pointed homotopy commutative and 
$
\begin{cases}
\ i_n^q\\
\ j_n^q
\end{cases}
$
are $P$-equivalences.  Choose pointed continuous functions 
$\phi_n^q:M_n^q \ra \bS_T^q$ such that 
$
\begin{cases}
\ \phi_n^q \circx i_n^q\\
\ \phi_n^q \circx j_n^q
\end{cases}
= \id.
$
\ 
We shall now construct a sequence $\{X_k,f_k\}$ such that $X_0 = X$ and $f_k:X_k \ra X_{k+1}$ is a $P$-equivalence.  
Thus, arguing by recursion, assume that $X_k$ has been constructed.  
Consider the set of morphisms $[f] \in [\bS_T^q,X_k]$ which 
cannot be factored through $\rho_n^q$ (failure of surjectivity of $(\rho_n^q)^*)$ and the set of morphisms 
$[g] \in [M_n^q,X_k]$ which cannot be factored through $\phi_n^q$ (failure of injectivity of  $(\rho_n^q)^*)$).  
If $\forall \ q$ $\&$ $\forall \ n$, these two sets are empty, then $X_k$ is $P$-local, so one can let $X_P = X_k$ and take for 
$l_P:X \ra X_P$ the composite 
$X_0 \ra$ 
$X_1 \ra$ 
$\cdots \ra X_k$.  Otherwise, form the pointed 2-source
\[
X_k \overset{\vee}{\lla} \bigvee\limits_{q,n} \bigl(\bigl(\bigvee\limits_f \bS_T^q\bigr) \vee 
\bigl(\bigvee\limits_g M_n^q\bigr)\bigr) \overset{h}{\lra} 
\bigvee\limits_{q,n} \bigl(\bigl(\bigvee\limits_f \bS_T^q\bigr) \vee 
\bigl(\bigvee\limits_g  \bS_T^q\bigr)\bigr)
\]
and call $X_{k+1}$ the pointed double mapping cylinder of $\vee$, $h$.  
Since $h$ is a $P$-equivalence (being a wedge of 
$P$-equivalences), the same is true of the arrow $X_k \ra X_{k+1}$, thereby completing the transition from $k$ to $k + 1$.  
Definition: $X_P = \telx(\bX,\bff)$.  
Accordingly, $l_P:X \ra X_P$ is a $P$-equivalence.  To prove that $X_P$ is $P$-local, it suffices to show that $X_P$ is orthogonal to the $\rho_n^q$.  Due to the compactness of 
$
\begin{cases}
\ \bS_T^q\\
\ M_n^q
\end{cases}
$
\hspace{-.3cm}, 
matters may be arranged
%%----------------------------------------------------------------------------------------------22
in such a way that any continuous function $
\begin{cases}
\ \bS_T^q \ra X_P\\
\ M_n^q \ra X_P
\end{cases}
$
factors through some $X_k$ 
(cf. p. \pageref{9.54}), hence the very construction of $X_P$ guarantees that every triangle
\begin{tikzcd}%[sep=large]
{\bS_T^q} \ar{d}[swap]{s} \ar{r}{\rho_n^q} &{\bS_T^q}\\
{X_P}
\end{tikzcd}
has a unique filler $\bS_T^q \overset{t}{\ra} X_P$.]\\
%\endgroup %%------------------------------------<<

\label{8.34}
The reflector $L_P$ produced by the homotopical $P$-localization theorem, when restricted to $\bHNILCWSP_*$, ``is'' the $L_P$ produced by the nilpotent $P$-localization theorem.  
Therefore, the idempotent triple corresponding to 
$P$-localization in $\bHCONCWSP_*$ is an extension of the idempotent triple corresponding to $P$-localization in 
$\bHNILCWSP_*$ 
(cf. p. \pageref{9.55a}) 
(however it is not the only such extension 
(cf. p. \pageref{9.55})).

Remarks: \ 
(1) $\forall \ X$, $\pi_1(X)_P \approx \pi_1(X_P)$; 
(2) $\forall \ X$ $\&$ $\forall \ n \geq 1$, the arrow $H_n(X) \ra H_n(X_P)$ is $P$-bijective but $H_n(X_P)$ need not be 
$P$-local (unless \mX is nilpotent); 
(3) $\forall \ X$ $\&$ $\forall \ n > 1$, $\pi_n(X_P)$ is $P$-local but the arrow $\pi_n(X) \ra \pi_n(X_P)$ need not be $P$-bijective (unless \mX is nilpotent).\\

\begingroup%%----------------------------------->>
\fontsize{9pt}{11pt}\selectfont
\textbf{\small EXAMPLE} \  
Let \mG be a group $-$then $K(G,1)_P \approx K(G_P,1)$ if \mG is nilpotent 
(cf. p. \pageref{9.56}) but this is false in  
general ($K(G,1)_P$ will ordinarily have nontrivial higher homotopy groups).  To illustrate, suppose that \mG is finite.  
Claim: $K(G,1)_P \approx K(G_P,1)$ iff $\ker l_P$ is $S_P$-torsion, $l_P:G \ra G_P$ the arrow of localization.  In fact, 
$K(G_P,1)$ is $P$-local, so the question is whether the arrow $K(G,1) \ra K(G_P,1)$ is a $P$-equivalence, which is the case iff 
$\ker l_P$ is $S_P$-torsion 
(cf. p. \pageref{9.57}).
\\ \indent
[Note: $K(S_3,1)_3$ is simply connected but $\pi_3(K(S_3,1)_3) \approx \Z/3\Z$.]\\
\endgroup %%------------------------------------<<

\begingroup%%----------------------------------->>
\fontsize{9pt}{11pt}\selectfont
\textbf{\small FACT} \  
For any \mG, the arrow of localization $l_P:G \ra G_P$ is an $HP$-homomorphism.
\\ \indent
[The triangle
\begin{tikzcd}[sep=large]
{K(G,1)} \ar{d} \ar{r} &{K(G,1)_P} \ar{ld}\\
{K(G_P,1)}
\end{tikzcd}
commutes in $\bHCONCWSP_*$.  In addition,\\
$H_*(K(G,1);\Z_P)$ $\approx$ 
$H_*(K(G,1)_P;\Z_P)$ and 
$\pi_1(K(G,1)_P) \approx G_P$.]\\
\endgroup %%------------------------------------<<

\begingroup%%----------------------------------->>
\fontsize{9pt}{11pt}\selectfont
The methods used in the proof of the homotopical $P$-localization theorem are of a general character and can easily be abstracted.  
What follows isolates the essentials.
\\ \indent
Fix a category \bC with coproducts.  
Let $S \subset \Mor \bC$ be a class of morphisms containing the isomorphisms of \bC which is closed under composition and cancellable.  
Problem: Find additional conditions on $S$ that will ensure that $S^\perp$ is the object class of a reflective subcategory of \bC.  
For this, assume that \mS is closed under coproducts and that for every 2-source $B \overset{f}{\la} A \ra A^\prime$, where $f \in S$, there is a weak pushout square
\begin{tikzcd}[sep=large]
{A} \ar{d}[swap]{f} \ar{r} &{A^\prime} \ar{d}{f^\prime}\\
{B} \ar{r} &{B^\prime}
\end{tikzcd}
, 
where $f^\prime \in S$.  Suppose further that there is a set $S_0 \subset S$ : $S_0^\perp = S^\perp$ 
%%----------------------------------------------------------------------------------------------23
and a regular cardinal $\kappa$ such that $\forall$ limit ordinal $\lambda \leq \kappa$, every diagram 
$\Delta:[0,\lambda[ \ra \bC$ in which the $\Delta_0 \ra \Delta_\alpha$ $(\alpha < \lambda)$ are in \mS admits a weak colimit 
$\Delta_\lambda$ such that $\Delta_0 \ra \Delta_\lambda$ is in \mS and when $\lambda = \kappa$, for each 
$f:A \ra B$ in $S_0$ 
(i)  $\forall$ $\phi \in \Mor(A,\Delta_\kappa)$, $\exists$ $\alpha < \kappa$ $\&$ $\phi_\alpha \in \Mor(A,\Delta_\alpha)$ : 
\begin{tikzcd}[sep=large]
{A} \ar{d}[swap]{\phi_\alpha} \ar{r}{\phi} &{\Delta_\kappa}\\
{\Delta_\alpha} \ar{ru}
\end{tikzcd}
commutes and 
(ii)  $\forall$ $\psi^\prime, \psi\pp \in \Mor(B,\Delta_\kappa)$ : $\psi^\prime \circx f = \psi\pp \circx f$, $\exists$ 
$\alpha < \kappa$ $\&$ $\psi_\alpha^\prime, \psi_\alpha\pp \in \Mor(B,\Delta_\alpha)$ :
\begin{tikzcd}[sep=large]
{B} \ar{d}[swap]{\psi_\alpha^\prime} \ar{r}{\psi^\prime} &{\Delta_\kappa}\\
{\Delta_\alpha} \ar{ru}
\end{tikzcd}
,
\begin{tikzcd}[sep=large]
{B} \ar{d}[swap]{\psi_\alpha\pp} \ar{r}{\psi\pp} &{\Delta_\kappa}\\
{\Delta_\alpha} \ar{ru}
\end{tikzcd}
commute $\&$ $\psi_\alpha^\prime \circx f = \psi_\alpha\pp \circx f$.
\\ \indent
Conclusion: $S = S^{\perp\perp}$ and $S^\perp$ is the object class of a reflective subcategory of \bC.
\\ \indent
[The verification proceeds by transfinite recursion, the only new wrinkle being that a limit ordinal $\lambda \leq \kappa$, 
$X_\lambda$ is taken to be the weak colimit of the $\{X_\alpha: \alpha < \lambda\}$ (as predicated per the hypotheses).  
Therefore, in the usual notation, $TX \equiv X_\kappa$.  
It is automatic that the arrow $\epsilon_X:X \ra TX$ is in \mS.  
Since $TX \in S_0^\perp = S^\perp$, what remains to be shown is that $S = S^{\perp\perp}$.  
Thus let $f:A \ra B$ be orthogonal to $S^\perp$.  Since $\epsilon_A:A \ra TA$ is in \mS and $TB \in S^\perp$, there is a unique filler 
$Tf:TA \ra TB$ for the diagram 
\begin{tikzcd}[sep=large]
{A} \ar{d}[swap]{\epsilon_A} \ar{r}{f} &{B} \ar{d}{\epsilon_B}\\
{TA} \ar[dashed]{r} &{TB}
\end{tikzcd}
.  
On the other hand, $\epsilon_B \circx f$ is orthogonal to $TA$, so one can find an arrow $TB \ra TA$ inverting $Tf$.  
It follows that $Tf$ is an isomorphism, hence $Tf \in S$ $\implies$ $\epsilon_B \circx f \in S$ $\implies$ $f \in S$, \mS being cancellable.]
\\ \indent
\label{10.6}
[Note: \ If \bC is cocomplete, then the statement simplifies.  
Example: The reflective subcategory theorem is a special case of these considerations 
(Ad\'amek-Rosicky\footnote[2]{\textit{Locally Presentable and Accessible Categories}, Cambridge University Press (1994), 30-35; 
see also Borceux, \textit{Handbook of Categorical Algebra 1}, Cambridge University Press (1994), 193-209.}).  
Applied to \bGR, one sees, e.g., that the $P$-localization of a countable group is countable.]\\
\endgroup %%------------------------------------<<

\label{8.35} %dmc mnft
\label{8.38} %dmc mnft indeed both of these seem not clear - near but not clearly here
\begingroup%%----------------------------------->>
\fontsize{9pt}{11pt}\selectfont
There are situations where the preceding remarks are not applicable since the assumption of cancellability on \mS may not be satisfied.  
The point is that cancellable means right cancellable and left cancellable, i.e.,  
$g \circx f \in S$ $\&$ $f \in S$ $\implies$ $g \in S$ and 
$g \circx f \in S$ $\&$ $g \in S$ $\implies$ $f \in S$.  
Let us drop the supposition that \mS is left cancellable (but retain everything else, including right cancellable) $-$then the argument above still implies that it is possible to assign to each object $X \in \Ob\bC$ another object $TX \in \Ob\bC$ and a morphism $\epsilon_X:X \ra TX$ in \mS.  Again, $TX \in S_0^\perp = S^\perp$, thus $S^\perp$ is the object class of a reflective subcategory of \bC but now the containent $S \subset S^{\perp\perp}$ can be strict (left cancellable is used to get 
$S = S^{\perp\perp}$).\\
\endgroup %%------------------------------------<<

\begingroup%%----------------------------------->>
\fontsize{9pt}{11pt}\selectfont
\textbf{\small EXAMPLE} \  
Let \bC be a cocomplete category, each object of which is $\kappa$-definite for some $\kappa$.  Let 
%%----------------------------------------------------------------------------------------------24
$S \subset \Mor\bC$ be a class of morphisms containing the isomorphisms of \bC which is closed under composition and right cancellable.  Assume that if 
\begin{tikzcd}[sep=large]
{A} \ar{d}[swap]{f} \ar{r} &{A^\prime} \ar{d}{f^\prime}\\
{B} \ar{r} &{B^\prime}
\end{tikzcd}
is a pushout square, then $f \in S$ $\implies$ $f^\prime \in S$ and if $\Xi \in \Nat(\Delta,\Delta^\prime)$, where 
$\Delta,\Delta^\prime:\bI \ra \bC$, then $\Xi_i \in S$ $(\forall \ i)$ $\implies$ $\colimx\Xi \in S$.  
Finally, suppose that there is a set $S_0 \subset S$: $S_0^\perp = S^\perp$.  
Accordingly, $S^\perp$ is the object class of a reflective subcategory of \bC and 
$\forall \ X$, the arrow $\epsilon_X:X \ra TX$ is in \mS.  
Examples: 
(1) Take $\bC = \bGR$ $-$then the class of $HP$-homomorphisms satisfies these conditions, 
hence $\forall$ \mG, the arrow of localization 
$l_{HP}:G \ra G_{HP}$ is in $S_{HP}$ 
(cf. p. \pageref{9.58}); 
(2) Take $\bC = G\text{-}\bMOD$ $-$then the class of $H\Z$-homomorphisms satisfies these conditions, hence $\forall \ M$, the arrow of localization $l_{H\Z}:M \ra M_{H\Z}$ is in $S_{H\Z}$ 
(cf. p. \pageref{9.59}).\\
\endgroup %%------------------------------------<<

\begingroup%%----------------------------------->>
\fontsize{9pt}{11pt}\selectfont
The role played in the theory by ``closure'' properties can be pinned down.
\vspi
Given a category \bC, let $S \subset \Mor\bC$ be a class of morphisms containing the isomorphisms of \bC and closed under composition with them.  Definition: \mS is said to be a 
\un{localization class}
\index{localization class} 
provided that it is possible to assign to each object $X \in \Ob\bC$ another object 
$TX \in \Ob\bC$ and a morphism $\epsilon_X:X \ra TX$ in \mS  with the following universal property: For every $f:A \ra B$ in 
\mS and for every $g:A \ra X$ there is a unique $t:B \ra TX$ such that $\epsilon_X \circx g = t \circx f$.  So, for any arrow 
$X \ra Y$, there is a commutative diagram
\begin{tikzcd}[sep=large]
{X} \ar{d} \ar{r}{\epsilon_X} &{TX} \ar{d}\\
{Y} \ar{r}[swap]{\epsilon_Y} &{TY}
\end{tikzcd}
, thus \mT defines a functor $\bC \ra \bC$ and $\epsilon:\id_{\bC} \ra T$ is a natural transformation.  Here, 
$\epsilon T = T \epsilon$ is not necessarily a natural isomorphism (it is if \mS is closed under composition).\\[.1cm]
\endgroup %%------------------------------------<<

\index{Theorem: Theorem of Korostenski-Tholen}
\index{Theorem of Korostenski-Tholen}
\begingroup%%----------------------------------->>
\fontsize{9pt}{11pt}\selectfont
\textbf{\small THEOREM OF KOROSTENSKI-THOLEN\footnote[2]{\textit{Comm. Algebra} \textbf{14} (1986), 741-766.}} \quad 
Let \mS be a localization class in a category \bC $-$then $S = S^{\perp\perp}$ iff \mS is closed under composition and left cancellable.  In addition, the assignment $S \ra S^\perp$ sets up a one-to-one correspondence between those localization classes \mS such that $S = S^{\perp\perp}$ and the conglomerate of reflective subcategories of \bC.\\
\endgroup %%------------------------------------<<

Let $[f]:X \ra Y$ be a morphism in $\bHCONCWSP_*$ $-$then $[f]$ (or f) is said to be an 
\un{$HP$-equivalence}
\index{HP-equivalence} 
if $\forall \ n \geq 0$, 
$f_*:H_n(X;\Z_P) \ra H_n(Y;\Z_P)$ is an isomorphism.

[Note: \ In the two extreme cases, viz. $P = \emptyset$ or $P = \bPi$, $HP$ is replaced by $H\Q$ or $H\Z$.]

\indent\indent (Wedges) Let 
$
\begin{cases}
\ X_i\\
\ Y_i
\end{cases}
(i \in I)
$
be pointed connected CW spaces.  Suppose that $\forall \ i$, $f_i:X_i \ra Y_i$ is an $HP$-equivalence $-$then
$\bigvee\limits_i f_i:\bigvee\limits_i X_i \ra \bigvee\limits_i Y_i$ is an $HP$-equivalence.

%%----------------------------------------------------------------------------------------------25
[This is because $\forall \ n \geq 1$, 
$H_n(\bigvee\limits_i X_i;\Z_P) \approx$ $\bigoplus\limits_i H_n(X_i;\Z_P)$ and 
$H_n(\bigvee\limits_i Y_i;\Z_P) \approx$ \\$\bigoplus\limits_i H_n(Y_i;\Z_P)$.]\\
\indent\indent (Pushouts) Suppose that 
$
\begin{cases}
\ X\\
\ Y
\end{cases}
$
are pointed connected CW spaces, $A \subset X$ a pointed connected CW subspace, and $f:A \ra Y$ a pointed continuous function.  
Assume: The inclusion $A \ra Y$ is a closed cofibration and an $HP$-equivalence $-$then the arrow $Y \ra X \sqcup_f Y$ is an $HP$-equivalence.

[The adjunction space $X \sqcup_f Y$ is a pointed connected CW space (cf. $\S 5$, Proposition 7) and it has the same pointed homotopy type as the pointed double mapping cylinder of the pointed 2-source 
$X \la A \ra Y$ (cf. $\S 3$, Proposition 18).]\\

\begin{proposition} \ %13
Every $P$-equivalence $f:X \ra Y$ is an $HP$-equivalence.
\end{proposition}

[Specializing Proposition 11, one can say that $\forall \ n \geq 0$, 
$f^*:H^n(Y;\Z_P) \ra H^n(X;\Z_P)$ is an isomorphism, hence $\forall \ n \geq 0$, 
$f_*:H_n(X;\Z_P) \ra H_n(Y;\Z_P)$ is an isomorphism (cf. $\S 8$, Proposition 2).]

[Note: \ An $HP$-equivalence need not be a $P$-equivalence.  For instance, take $P = \bPi$ $-$then 
$HP$-equivalence = homology equivalence and $P$-equivalence = pointed homotopy equivalence.]\\

Given a set of primes $P$, a pointed connected CW space \mX is said to be 
\un{$HP$-local}
\index{HP-local} 
provided that $[f] \perp X$ for every $HP$-equivalence $f$.\\

\textbf{\small SUBLEMMA} \quad 
Let \mK be a pointed connected CW complex, $L \subset K$ $(L \neq K)$ a pointed connected subcomplex such that 
$H_*(K,L;\Z_P) = 0$ $-$then there exists a pointed countable connected subcomplex $A \subset K$ such that 
$A \not\subset L$ and $H_*(A,A \cap L;\Z_P) = 0$.

[We shall construct an expanding sequence of pointed countable connected subcomplexes $A_1, A_2, \ldots$ of \mK such that 
$\forall \ n$, $A_n \not\subset L$ and the arrow 
$H_*(A_n,A_n \cap L;\Z_P) \ra$ 
$H_*(A_{n+1},A_{n+1} \cap L;\Z_P)$ is the zero map.  Thus fix $A_1$: $A_1 \not\subset L$.  Given $A_n$, for each element 
$x \in H_*(A_n,A_n \cap L;\Z_P)$ choose a pointed finite connected subcomplex $K_x \subset K$ such that $x$ goes to zero in 
$H_*(A_n \cup K_x, (A_n \cup K_x) \cap L; \Z_P)$.  
Let $A_{n+1}$ be the union of the $A_n$ and the $K_x$ and put 
$A = \bigcup\limits_n A_n$.]

[Note: \ $A \cap L$ is necessarily connected.]\\

\label{13.60}
\textbf{\small LEMMA}  \  
Let \mZ be a pointed connected CW space.  Suppose that for any CW pair $(K,L)$, where \mK is a pointed countable connected CW complex and $L \subset K$ $(L \neq K)$ is a pointed connected subcomplex such that 
$H_*(K,L;\Z_P) = 0$, the arrow $[K,Z] \ra [L,Z]$ is surjective $-$then \mZ is $HP$-local.

%%----------------------------------------------------------------------------------------------26
[The claim is that for every $HP$-equivalence $f:X \ra Y$, the precomposition arrow 
$f^*:[Y,Z] \ra [X,Z]$ is bijective.  Since it is clear that the class of $HP$-equivalences admits a calculus of left fractions 
(cf. p. \pageref{9.60}), it need only be shown that $f^*:[Y,Z] \ra [X,Z]$  is surjective.  For this purpose, one can make the usual adjustments and take 
$
\begin{cases}
\ X\\
\ Y
\end{cases}
$
to be pointed connected CW complexes and $f:X \ra Y$ the inclusion with $X \neq Y$.  
Build now a transfinite sequence of pointed connected subcomplexes $X_\alpha$ of \mY via the following procedure.  Set $X_0 = X$.  
Owing to the sublemma, there exists a pointed countable connected subcomplex $A_0 \subset Y$ such that $A_0 \not\subset X_0$ and 
$H_*(A_0,A_0 \cap X_0;\Z_P) = 0$.   Set $X_1 = A_0 \cup X_0$.  
Case 1: $X_1 = Y$.  In this situation, the arrow $[Y,Z] \ra [X,Z]$ is surjective.  For let $\phi:X \ra Z$ be a pointed continuous function.  Since the inclusion $A_0 \cap X_0 \ra A_0$ is a cofibration, our assumptions imply that the restriction of $\phi$ to 
$A_0 \cap X_0$ extends to a pointed continuous function $A_0 \ra Z$, thus $\phi$ extends to a pointed continuous function 
$\Phi:Y \ra Z$.
Case: 2 $X_1 \neq Y$.  
Utilizing excision, $H_*(X_1,X_0;\Z_P) = 0$, so from the exact sequence of the triple 
$(Y,X_1,X_0)$, $H_*(Y,X_1;\Z_P) = 0$.  
Therefore the sublemma is applicable to the pair $(Y,X_1)$, hence there exists a pointed countable connected subcomplex $A_1 \subset Y$ such that $A_1 \not\subset X_1$ and 
$H_*(A_1,A_1 \cap X_1;\Z_P) = 0$.  
Set $X_2 = A_1 \cup X_1$.  
Continue on out to a sufficiently large regular cardinal 
$\kappa$ (if necessary) , taking $X_\lambda = \bigcup\limits_{\alpha < \lambda} X_\alpha$ at a limit ordinal 
$\lambda \leq \kappa$ (observe that $H_*(Y,X_\lambda;\Z_P) = 0$), where $X_\kappa = Y$.]\\

Notation: $\bCONCWSP_{*,HP}$
\index{$\bCONCWSP_{*,HP}$} 
is the full subcategory of $\bCONCWSP_{*}$ whose objects are the pointed connected CW spaces which are $HP$-local and $\bHCONCWSP_{*,HP}$
\index{$\bHCONCWSP_{*,HP}$} 
is the associated homotopy category.\\

\index{Theorem: Homotopical $HP$-Localization Theorem}
\index{Homotopical $HP$-Localization Theorem}
\textbf{\small HOMOTOPICAL $HP$-LOCALIZATION THEOREM} \quad 
$\bHCONCWSP_{*,HP}$ is a reflective subcategory of $\bHCONCWSP_*$.

[The theorem will follow provided that one can show that it is possible to assign to each pointed connected CW space \mX an $HP$-local pointed connected CW space $X_{HP}$ and an $HP$-equivalence $l_{HP}: X \ra X_{HP}$.  
The full subcategory $\bHCW_*$ whose objects are the pointed countable connected CW complexes has a small skeleton.  
One can therefore choose a set of CW pairs $(K_i,L_i)$, where $K_i$ is a pointed countable connected CW complex and $L_i \subset K_i$, $(L_i \neq K_i)$ is a pointed connected subcomplex such that $H_*(K_i,L_i;\Z_P) = 0$, which contains up to isomorphism all such CW pairs with these properties.  
Assuming that \mX is a pointed connected CW complex, construct an expanding transfinite sequence 
$X = X_0 \subset X_1 \subset \cdots \subset X_\alpha \subset X_{\alpha + 1} \subset \cdots \subset X_\Omega$ of pointed connected CW complexes by setting $X_\lambda = \bigcup\limits_{\alpha < \lambda} X_\alpha$ at a limit ordinal 
$\lambda \leq \Omega$ and defining $X_{\alpha + 1}$ by the pushout square
%%----------------------------------------------------------------------------------------------27
\begin{tikzcd}%[sep=large]
{\bigvee\limits_i\bigvee\limits_f L_i} \ar{d} \ar{r} 
&{X_\alpha} \ar{d}\\
{\bigvee\limits_i\bigvee\limits_f K_i}  \ar{r} 
&{X_{\alpha + 1}}
\end{tikzcd}
.  
Here, $f$ runs over a set of skeletal representatives in $[L_i,X_\alpha]$ and the arrow $X_\alpha \ra X_{\alpha + 1}$ is an 
$HP$-equivalence.  Put $X_{HP} = X_\Omega$ $-$then $\forall \ i$, $[K_i,X_{HP}] \ra [L_i,X_{HP}]$ is surjective, thus by the lemma, $X_{HP}$ is $HP$-local.  That the inclusion $X \ra X_{HP}$ is an $HP$-equivalence is automatic.]\\

\label{9.55}
The reflector $L_{HP}$ produced by the homotopical $HP$-localization theorem, when restricted to $\bHNILCWSP_*$, 
``is'' the $L_P$ produced by the nilpotent $P$-localization theorem.  
Proof: If \mX is nilpotent and $P$-local, then \mX is $HP$-local, as can be seen by appealing to the preceding lemma and using the nilpotent obstruction theorem (cf. Proposition 2).  
Therefore the idempotent triple corresponding to $HP$-localization in 
$\bHCONCWSP_*$ is an extension of the idempotent triple corresponding to localization in $\bHNILCWSP_*$ 
(cf. p. \pageref{9.61}).  On the other hand, Proposition 13 implies that every $HP$-local space is $P$-local, so there is a natural transformation 
$L_P \ra L_{HP}$.\\

\begin{proposition} \ %14
Let $[f]:X \ra Y$ be a morphism in $\bHCONCWSP_*$.  
Assume: $[f]$ is orthogonal to every $HP$-local pointed connected CW space $-$then $[f]$ is an $HP$-equivalence.
\end{proposition}

[By hypothesis, for every $HP$-local \mZ, $[Y,Z] \approx [X,Z]$.  Specialize and substitute in $Z = K(\Z_P,n)$ (which is 
$HP$-local) to get $H^n(Y;\Z_P) \approx H^n(X;\Z_P)$ $\forall \ n \geq 1$ or still, 
$H_n(X;\Z_P) \approx H_n(Y;\Z_P)$ $\forall \ n \geq 1$ (cf. $\S 8$, Proposition 2).]

[Note: \ Thus, in the homotopy category, the class of $HP$-equivalences is ``saturated'' but the group theoretic analog of this is false 
(cf. p. \pageref{9.62}).]\\

In the $P$-local situation, one starts with an intrinsic definition of the $P$-local objects and defines the $P$-equivalences via orthogonality, while in the $HP$-local situation, one starts with an intrinsic definition of the $HP$-equivalences and defines the 
$HP$-local objects via orthogonality.  The $P$-equivalences are characterized in Proposition 11, so to complete the picture, it is necessary to characterize the $HP$-local objects.

A pointed connected CW space \mX is said to satisfy 
\un{Bousfield's condition}
\index{Bousfield's condition} if 
$\forall \ n \geq 1$, $\pi_n(X)$ is an $HP$-local group and $\forall \ n \geq 2$ $\pi_n(X)$ is an $H\Z$-local $\pi_1(X)$-module.

[Note: \ Recall that an abelian group is $P$-local iff it is $HP$-local.]\\

\index{Lemma B}
\textbf{\small LEMMA B} \quad  Let \mX be a pointed connected CW space.  Fix $n > 1$ and suppose that 
$\phi:\pi_n(X) \ra M$ is a homomorphism of $\pi_1(X)$-modules $-$then 
$\phi_P:\pi_n(X)_P \ra M_P$ is an $H\Z$-homomorphism iff there exists a pointed connected CW space \mY and a pointed continuous
%%----------------------------------------------------------------------------------------------28
function $f:X \ra Y$ such that 
$H_*(f):H_*(X;\Z_P) \approx H_*(Y;\Z_P)$, 
$\pi_q(f):\pi_q(X) \approx \pi_q(Y)$ $(q < n)$, and 
$\pi_n(f) \approx \phi$ in $\pi_n(X)\backslash \pi_1(X)$-\bMOD.

[To establish the sufficiency, compare the exact sequence \ 
$H_{n+2}(P_{n-1}X;\Z_P)$  \ $\ra$
$H_1(P_{n-1}X;\pi_n(X)_P)$ $\ra$ 
$H_{n+1}(P_nX;\Z_P) \ra $
$H_{n+1}(P_{n-1}X;\Z_P) \ra$
$H_0(P_{n-1}X;\pi_n(X)_P) \ra$ 
$H_n(P_nX;\Z_P) \ra$
$H_n(P_{n-1}X;\Z_P) \ra 0$ on 
p. \pageref{9.63} with its analog for \mY, noting that
$H_1(P_{n-1}X;\pi_n(X)_P) \approx$ 
$H_1(\pi_1(X);\pi_n(X)_P)$, 
$H_0(P_{n-1}X;\pi_n(X)_P) \approx$ 
$H_0(\pi_1(X);\pi_n(X)_P)$.  
Indeed, there are bijections  \ 
$H_q(P_nX;\Z_P) \approx$ 
$H_q(P_nY;\Z_P)$ $(q \leq n)$  \ and a surjection  \ 
$H_{n+1}(P_nX;\Z_P) \ra H_{n+1}(P_nY;\Z_P)$ 
(cf. p. \pageref{9.64}).

To establish the necessity, attach certain $n$-cells and $(n+1)$-cells to \mX so as to produce a relative CW complex 
$(\ov{X},X)$ and an isomorphism $\pi_n(\ov{X}) \ra M$ such that $X[n-1] \approx \ov{X}[n-1]$ and the triangle
\begin{tikzcd}[sep=small]
&{\pi_n(X)} \ar{ldd} \ar{rdd}{\phi}\\
\\
{\pi_n(\ov{X})} \ar{rr} &&{M}
\end{tikzcd}
commutes.  The composite $X \ra X[n] \ra \ov{X}[n]$ induces an arrow 
$H_q(X;\Z_P) \ra H_q(\ov{X}[n];\Z_P)$ which is bijective for $q \leq n$ and surjective for $q = n+1$.  Apply the Kan factorization theorem.]\\

\begin{proposition} \ %15
Let 
$
\begin{cases}
\ X\\
\ Y
\end{cases}
$
be pointed connected CW spaces, $f:X \ra Y$ a pointed continuous function.  Assume:
$
\begin{cases}
\ X\\
\ Y
\end{cases}
$
satisfy Bousfield's condition and $f$ is an $HP$-equivalence $-$then $f$ is a pointed homotopy equivalence.
\end{proposition}

[Obviously, $\Z_P \otimes \pi_1(X)/[\pi_1(X),\pi_1(X)] \approx \Z_P \otimes \pi_1(Y)/[\pi_1(Y),\pi_1(Y)]$.  
Furthermore, the horizontal arrows in the commutative diagram
\begin{tikzcd}%[sep=large]
{H_2(X;\Z_P)} \ar{d} \ar{r} &{H_2(\pi_1(X);\Z_P)} \ar{d}\\
{H_2(Y;\Z_P)} \ar{r} &{H_2(\pi_1(Y);\Z_P)}
\end{tikzcd}
are surjective 
(cf. p. \pageref{9.65}) and 
$H_2(X;\Z_P) \approx H_2(Y;\Z_P)$.  Therefore 
$f_*:\pi_1(X) \ra \pi_1(Y)$ is an $HP$-homomorphism.  But this means that $f_*$ is an isomorphism, 
$
\begin{cases}
\ \pi_1(X)\\
\ \pi_1(Y)
\end{cases}
$
being $HP$-local.  Next, consider the commutative diagram
$
\begin{tikzcd}%[sep=large]
{\pi_2(X)} \ar{d} \ar{r}{f_*} &{\pi_2(Y)}\ar{d}\\
{\pi_2(X)_P} \ar{r}[swap]{(f_*)_P} &{\pi_2(Y)_P}
\end{tikzcd}
.  
$
The vertical arrows are isomorphisms and $(f_*)_P$ is an $H\Z$-homomorphism (cf. Lemma B).  Consequently, 
$f_*:\pi_2(X) \ra \pi_2(Y)$ is an $H\Z$-homomorphism between $H\Z$-local $\pi_1(X)$-modules, hence is an isomorphism.  That $f$ is a weak homotopy equivalence then follows by iteration.]\\

\textbf{\small LEMMA}  \  
For any pointed connected CW space \mX, there exists a pointed connected CW space $X_B$ which satisfies Bousfield's condition and an $HP$-equivalence $l_B:X \ra X_B$, where $\pi_1(X)_{HP} \approx \pi_1(X_B)$.

%%----------------------------------------------------------------------------------------------29
[Fix a pointed continuous function $\phi:X \ra K(\pi_1(X)_{HP},1)$ such that $\phi_* = l_{HP}$, where 
$l_{HP}:\pi_1(X) \ra \pi_1(X)_{HP}$ is the arrow of localization.  Since $l_{HP}$ is an $HP$-homomorphism, 
the Kan factorization theorem implies that there exists a pointed connected CW space $X_1$ and pointed continuous functions 
$f_1:X \ra X_1$, $\psi_1:X_1 \ra K(\pi_1(X)_{HP},1)$ with $\phi = \psi_1 \circx f_1$ such that $f_1$ is an $HP$-equivalence and 
$\pi_1(\psi_1):\pi_1(X_1) \ra \pi_1(X)_{HP}$ is an isomorphism.  
Continuing, construct a pointed connected CW space $X_2$, a pointed continuous function $f_2:X_1 \ra X_2$, 
and an isomorphism $\pi_2(X_2) \ra (\pi_2(X_1)_P)_{H\Z}$ such that $f_2$ is an 
$HP$-equivalence, $\pi_1(f_2):\pi_1(X_1) \ra \pi_1(X_2)$ is an isomorphism, and the composite 
$\pi_2(X_1) \ra$ 
$\pi_2(X_2) \ra$
$(\pi_2(X_1)_P)_{H\Z}$ 
equals the composite 
$\pi_2(X_1) \ra$ 
$\pi_2(X_1)_P \ra$
$(\pi_2(X_1)_P)_{H\Z}$ 
(cf. Lemma B and $\S 8$, Proposition 21).  
This gives $X \ra X_1 \ra X_2$.  Proceed from here inductively and let $X_B$ be the pointed mapping telescope of the sequence thereby obtained.]

[Note: \ It is apparent from the construction of $X_B$ that if $\pi_q(X)$ is an $HP$-local group for $1 \leq q \leq n$ and if 
$\pi_q(X)$ is an $H\Z$-local $\pi_1(X)$-module for $2 \leq q \leq n$, then $\forall \ q \leq n$, 
$\pi_q(X) \approx \pi_q(X_B)$.]\\

\begin{proposition} \ %16
Let \mX be a pointed connected CW space $-$then \mX is $HP$-local iff \mX satisfies Bousfield's condition.
\end{proposition}

[Suppose that \mX satisfies Bousfield's condition.  
Bearing in mind that the class of $HP$-equivalences admits a calculus of left fractions, 
to  prove that \mX is $HP$-local, it suffices to show that every $HP$-equivalence $f:X \ra Y$ has a left inverse 
$g:Y \ra X$ in $\bHCONCWSP_*$, i.e., $g \circx f \simeq \id_X$.  
For this purpose, apply the lemma to get 
$l_B:Y \ra Y_B$ $-$then the composite $l_B \circx f:X \ra Y_B$ is a pointed homotopy equivalence (cf. Proposition 15), so 
$\exists$ $h:Y_B \ra X$ such that $h \circx l_B \circx f \simeq \id_X$ and we can take $g = h \circx l_B$.  
Conversely, suppose that \mX is $HP$-local.  By what has just been said, $X_B$ is $HP$-local, thus $l_B:X \ra X_B$ is a pointed homotopy equivalence.]\\

Application: $\forall \ X$, $\pi_1(X)_{HP} \approx \pi_1(X_{HP})$.\\

\begingroup%%----------------------------------->>
\fontsize{9pt}{11pt}\selectfont
\textbf{\small EXAMPLE} \  
Take $X = \bS^1 \vee \bS^1$: $\pi_1(X)_P$ is countable but $\pi_1(X)_{HP}$ is uncountable if $2 \in P$.\\
\endgroup %%------------------------------------<<

\begingroup%%----------------------------------->>
\fontsize{9pt}{11pt}\selectfont
\textbf{\small EXAMPLE} \  
When \mP is the set of all primes, every space is $P$-local.  However, not every space is $H\Z$-local and in fact the effect of 
$H\Z$-localization on the higher homotopy groups can be drastic even if the fundamental group is nilpotent.  Thus let \mX be a pointed connected CW space and for $q > 1$, put $\widehat{\pi}_q(X) = \lim \pi_q(X)/(I[\pi_1(X)])^i \cdot \pi_q(X)$.  
Note that $\widehat{\pi}_q(X)$ is an $H\Z$-local $\pi_1(X)$-module, being the limit of nilpotent $\pi_1(X)$-modules 
(cf. p. \pageref{9.66}).  Assume now that $\pi_1(X)$ is a finitely generated nilpotent group.  Suppose further that 
(i) $\pi_q(X)$ is a nilpotent $\pi_1(X)$-module $(1 < q \leq n)$ and
(ii) $\pi_q(X)$ is a finitely
%%----------------------------------------------------------------------------------------------30
generated $\pi_1(X)$-module $(n < q \leq 2n)$ $(n \geq 1)$ $-$then 
Dror-Dwyer\footnote[2]{\textit{Illinois J. Math.} \textbf{21} (1977), 675-684.}
have shown that 
(i) $\pi_q(X_{H\Z}) \approx \pi_q(X)$ $(1 < q \leq n)$ and 
(ii) $\pi_q(X_{H\Z}) \approx \widehat{\pi}_q(X)$ $(n < q \leq 2n)$ $(n \geq 1)$.  
In this situation, the first conclusion is actually automatic, so the impact lies in the second.  
Example: Take $X = \bP^2(\R)$ and $n = 1$ to see that $\pi_2(X_{H\Z}) \approx \widehat{\Z}_2$, the 2-adic integers.\\
\endgroup %%------------------------------------<<

\index{Theorem: \textit{HP} Whitehead Theorem}
\index{HP Whitehead Theorem}%\index{\textit{HP} Whitehead Theorem}
\textbf{\small \textit{HP} WHITEHEAD THEOREM} \quad 
Suppose that \mX and \mY are $HP$-local and let $f:X \ra Y$ be a pointed continuous function.  Assume: 
$f_*:H_q(X;\Z_P) \ra H_q(Y;\Z_P)$ is bijective for $1 \leq q < n$ and surjective for $q = n$ $-$then $f$ is an $n$-equivalence.

[If $n = 1$, the claim is that 
$f_*:H_1(\pi_1(X);\Z_P) \ra H_1(\pi_1(Y);\Z_P)$ surjective $\implies$ $f_*:\pi_1(X) \ra \pi_1(Y)$ surjective, which is true 
(cf. p. \pageref{9.67}).  
If $n > 1$, use the Kan factorization theorem to write $f = \psi_f \circx \phi_f$, where $\phi_f:X \ra X_f$ is an $HP$-equivalence and $\psi_f:X_f \ra Y$ is an $n$-equivalence.  Since \mX is $HP$-local, $X \approx (X_f)_{HP}$ and since \mY is $HP$-local, 
$\pi_q(X_f) \approx \pi_q(Y)$ $(1 \leq q < n)$ $\implies$ 
$\pi_q(X_f) \approx \pi_q((X_f)_{HP})$ $(1 \leq q < n)$.  Therefore the arrow $\pi_q(X) \ra \pi_q(Y)$ is bijective for 
$1 \leq q < n$ and surjective for $q = n$, i.e., $f$ is an $n$-equivalence.]

[Note: \ Taking $\Z_P = \Z$ and 
$
\begin{cases}
\ X\\
\ Y
\end{cases}
$
nilpotent leads to a refinement of Dror's Whitehead theorem (which, of course, can also be derived directly).]\\

\begingroup%%----------------------------------->>
\fontsize{9pt}{11pt}\selectfont
\textbf{\small EXAMPLE} \  
Let \mX be a pointed connected CW space.  Assume: $\widetilde{H}_*(X;\Z_P) = 0$, i.e., \mX is $\Z_P$-acyclic $-$then 
$X_{HP}$ is contractible.\\
\endgroup %%------------------------------------<<

Given an abelian group \mG, one can introduce the notion of ``$HG$-equivalence'' and play the tape again.  So, employing obvious notation, the upshot is that $\bHCONCWSP_{*,HG}$ is a reflective subcategory of $\bHCONCWSP_{*}$, with reflector $L_{HG}$ which sends \mX to $X_{HG}$.

[Note: \ The CW pairs $(K,L)$ that intervene when testing for ``$HG$-local'' have the property that the cardinality of the set of cells in \mK is $\leq \#(G)$ if $\#(G)$ is infinite and $\leq \omega$ if $\#(G)$ is finite.]

While the number of distinct homological localizations appears to be large, the reality is that all the possibilities can be described in a simple way.  
Definition: $L_{HG^\prime}$ and $L_{HG\pp}$ have the \un{same acyclic spaces} if 
$\widetilde{H}_*(X;G^\prime) = 0$ $\Leftrightarrow$ 
$\widetilde{H}_*(X;G\pp) = 0$ or still, if the $HG^\prime$-equivalences are the same as $HG\pp$-equivalences, hence that 
$L_{HG^\prime}$ and $L_{HG\pp}$ are naturally isomorphic.

\label{17.31}
Given an abelian group \mG, call $S(G)$ the class of abelian groups \mA such that $A \otimes G = 0 = \Tor(A,G)$.\\

%%----------------------------------------------------------------------------------------------31
\begin{proposition} \ %17
Let $\Acy_G$ be the class of $G$-acyclic spaces $-$then $\sS(G) = \{\widetilde{H}_n(X): n \geq 0 \ \& \ X \in \Acy_G\}$.
\end{proposition}

[This follows from the universal coefficient theorem and the existence of Moore spaces.]\\

Application: $\sS(G^\prime) = \sS(G\pp)$ iff $\Acy_{G^\prime} = \Acy_{G\pp}$\\

\label{9.122}
Given an abelian group \mG, let $P_G$ be the set of primes $p$ such that \mG is not uniquely divisible by $p$ and put  
$
S_G \hsx = \hsx 
\begin{cases}
\ \bigoplus\limits_{p \in P_G} \Z/p\Z \text{\quad if \ } \Q \otimes G = 0\\
\ \Z_{P_G} \indent\quad \ \text{ if \ }  \Q \otimes G \neq 0
\end{cases}
$
$-$then $S(G) = S(S_G)$ 
(cf. p. \pageref{9.68} ff.).  
Corollary: $L_{HG} \approx L_{HS_G}$.  
Therefore, besides the $L_{HP}$, the only other homological localizations that need be considered are those corresponding to 
$\bigoplus\limits_{p \in P} \Z/p\Z$ for some \mP.\\

\begingroup%%----------------------------------->>
\fontsize{9pt}{11pt}\selectfont
\textbf{\small FACT} \quad
Let \mX be a pointed connected CW space $-$then \hsx 
$\widetilde{H}_*(X;\ds\bigoplus\limits_{p \in P} \Z/p\Z) \hsx = \hsx 0$ \hsx iff \hsx
$\widetilde{H}_*(X;\ds\prod\limits_{p \in P} \Z/p\Z)$ $= 0$.\\
\endgroup %%------------------------------------<<

The ``$\Z/p\Z$-theory'' (= ``$\F_p$-theory''), in its general aspects, runs parallel to the ``$\Z_P$-theory'' but there are some differences in detail.

A pointed connected CW space \mX is said to satisfy 
\un{Bousfield's condition mod $p$}
\index{Bousfield's condition mod $p$} 
if $\forall \ n \geq 1$, $\pi_n(X)$ is an $H\F_p$-local group and $\forall \ n \geq 2$, $\pi_n(X)$ is an $H\Z$-local $\pi_1(X)$-module.

[Note: \ Recall that an abelian group is $H\F_p$-local iff it is $p$-cotorsion.]\\

\index{Lemma B mod $p$}
\textbf{\small LEMMA B mod $p$} \quad 
Let \mX be a pointed connected CW space .  
Fix $n > 1$ and suppose that $\phi:\pi_n(X) \ra M$ is a homomorphism of 
$\pi_1(X)$-modules.  Consider the following conditions.\\
\indent\indent $(C_1)$ \quad $\id \otimes \phi:\F_p \otimes \pi_n(X) \ra \F_p \otimes M$ is an $H\Z$-homomorphism.\\
\indent\indent $(C_2)$ \quad $\phi_*:H_0(\pi_1(X);\Tor(\F_p,\pi_n(X))) \ra H_0(\pi_1(X);\Tor(\F_p,M))$ is surjective.\\
\indent\indent $(C_3)$ \quad $\id \otimes \phi:\F_p \otimes \pi_n(X) \ra \F_p \otimes M$ is an isomorphism.\\ 
Then $C_1 + C_2 \implies$\\
\indent\indent ($E$) \  There exists a pointed connected CW space \mY and a pointed continuous function $f:X \ra Y$ such that 
$H_*(f):H_*(X;\F_p) \approx H_*(Y;\F_p)$, 
$\pi_q(f): \pi_q(X) \approx \pi_q(Y)$ $q < n)$, and 
$\pi_n(f) \approx \phi$ in $\pi_n(X)\backslash \pi_1(X)$-\bMOD.

Conversely, $E \implies  C_1$ and $E + C_3 \implies C_2$.\\

\begin{proposition} \ %18
Let 
$
\begin{cases}
\ X\\
\ Y
\end{cases}
$
be pointed connected CW spaces, $f:X \ra Y$ a pointed continuous function.  Assume:
$
\begin{cases}
\ X\\
\ Y
\end{cases}
$
satisfy Bousfield's condition mod $p$ and $f$ is an $H\F_p$-equivalence $-$then $f$ is a pointed homotopy equivalence.
\end{proposition}

%%----------------------------------------------------------------------------------------------32
[Arguing as in the proof of Proposition 15, one finds that $f_*:\pi_1(X) \ra \pi_1(Y)$ is an isomorphism.  
To discuss 
$f_*:\pi_2(X) \ra \pi_2(Y)$, define $M, N$ in $\pi_1(X)$-\bMOD by the exact sequence 
$0 \ra M \ra \pi_2(X) \ra \pi_2(Y) \ra N \ra 0$.  
The claim is that $M = 0 = N$, hence that 
$f_*:\pi_2(X) \ra \pi_2(Y)$ is an isomorphism.  
For this, it need only be shown that $\F_p \otimes M = 0 = \F_p \otimes N$ 
(both \mM and \mN are $H\F_p$-local).  Since $f$ is an $H\F_p$-equivalence, 
$\id \otimes f_*: \F_p \otimes \pi_2(X) \ra \F_p \otimes \pi_2(Y)$ is an $H\Z$-homomorphism $(E \implies C_1)$.  
But 
$
\begin{cases}
\ \F_p \otimes \pi_2(X)\\
\ \F_p \otimes \pi_2(Y)
\end{cases}
$
are $H\Z$-local 
(cf. p. \pageref{9.69}), so 
$\F_p \otimes \pi_2(X) \approx \F_p \otimes \pi_2(Y)$, from which $\F_p \otimes N = 0$.  Using now the exact sequence 
$\Tor(\F_p,\pi_2(X)) \ra$
$\Tor(\F_p,\pi_2(Y)) \ra$ 
$\F_p \otimes M \ra 0$, $E + C_3 \implies C_2$ gives $H_0(\pi_1(X);\F_p \otimes M) = 0$.  However, \mM is $H\Z$-local 
(being a kernel), thus $\F_p \otimes M$ is $H\Z$-local 
(cf. p. \pageref{9.70}).  And: 
$\F_p \otimes M = I[\pi_1(X)]\cdot (\F_p \otimes M)$ $\implies$ 
$(\F_p \otimes M)_{H\Z} = 0$  
(cf. p. \pageref{9.71}) $\implies$ $\F_p \otimes M = 0$.  
That $f$ is a weak homotopy equivalence then follows by iteration.]\\

\textbf{\small LEMMA}  \  
For any pointed connected CW space \mX, there exists a pointed connected CW space $X_B$ which satisfies Bousfield's condition mod $p$ and an $H\F_p$-equivalence $l_B:X \ra X_B$, where $\pi_1(X)_{H\F_p} \approx \pi_1(X_B)$.\\

[Construct $f_1:X \ra X_1$ as before (the Kan factorization theorem holds mod $p$ 
(cf. p. \pageref{9.72})).  
Continuing, construct a pointed connected CW space $X_1^\prime$, a pointed continuous function 
$f_1^\prime:X_1 \ra X_1^\prime$, and an isomorphism 
$\pi_2(X_1^\prime) \ra \pi_2(X_1)_{H\Z}$ such that $f_1^\prime$ is an $H\Z$-equivalence, 
$\pi_1(f_1^\prime):\pi_1(X_1) \ra \pi_1(X_1^\prime)$ is an isomorphism, and the composite 
$\pi_2(X_1) \ra$ 
$\pi_2(X_1^\prime) \ra$ 
$\pi_2(X_1)_{H\Z}$ is the arrow 
$\pi_2(X_1) \ra \pi_2(X_1)_{H\Z}$ 
(cf. Lemma B ($P = \bPi$)).  This gives 
$X \ra X_1 \ra X_1^\prime$.  
Next, construct a pointed connected CW space $X_2$, a pointed continuous function 
$f_1\pp:X_1^\prime \ra X_2$, and an isomorphism 
$\pi_2(X_2) \ra \Ext(\Z/p^\infty\Z,\pi_2(X_1^\prime))$ such that $f_1\pp$ is an $H\F_p$-equivalence, 
$\pi_1(f_1\pp):\pi_1(X_1^\prime) \ra \pi_1(X_2)$ is an isomorphism, and the composite 
$\pi_2(X_1^\prime) \ra$ 
$\pi_2(X_2) \ra$ 
$\Ext(\Z/p^\infty\Z,\pi_2(X_1^\prime))$ is the arrow 
$\pi_2(X_1^\prime) \ra$ 
$\Ext(\Z/p^\infty\Z,\pi_2(X_1^\prime))$ 
(cf. Lemma B mod $p$ and 
p. \pageref{9.73}).  
To justify that application of $C_1 + C_2 \implies E$, note that the arrow 
$\F_p \otimes \pi_2(X_1^\prime) \ra$ 
$\F_p \otimes \Ext(\Z/p^\infty\Z,\pi_2(X_1^\prime))$ is bijective and the arrow 
$\Tor(\F_p,\pi_2(X_1^\prime)) \ra$ 
$\Tor(\F_p, \Ext(\Z/p^\infty\Z,\pi_2(X_1^\prime))$ is surjective 
(cf. p. \pageref{9.74}).  This gives 
$X \ra X_1 \ra X_1^\prime \ra X_2$.  
Proceed from here inductively and let $X_B$ be the pointed mapping telescope of the sequence thereby obtained.]

[Note: \ It is apparent from the construction of $X_B$ that if $\pi_q(X)$ is an $H\F_p$-local group for $1 \leq q \leq n$ and if $\pi_q(X)$ is an $H\Z$-local $\pi_1(X)$-module for $2 \leq q \leq n$, then $\forall \ q \leq n$, 
$\pi_q(X) \approx \pi_q(X_B)$.]\\

\begin{proposition} \ %19
Let \mX be a pointed connected CW space $-$then \mX is $H\F_p$-local iff \mX satisfies Bousfield's condition mod $p$.
\end{proposition}

[The proof is the same as that of Proposition 16.]\\

%%----------------------------------------------------------------------------------------------33
Application: $\forall \ X$, $\pi_1(X) _{H\F_p} \approx \pi_1(X_{H\F_p})$.\\

\begingroup%%----------------------------------->>
\fontsize{9pt}{11pt}\selectfont
\textbf{\small EXAMPLE} \  
Let \mX be a pointed connected CW space.  Assume: The homotopy groups of \mX are finite $-$then $\forall \ n \geq 1$ 
$\pi_n(X_{H\F_p})$ is a finite $p$-group, thus $X_{H\F_p}$ is nilpotent.
\\ \indent
[For here $\pi_1(X)_{H\F_p} \approx \pi_1(X)_p$ 
(cf. p. \pageref{9.75}), which is a finite $p$-group 
(cf. p. \pageref{9.76}).]\\
\endgroup %%------------------------------------<<

\begingroup%%----------------------------------->>
\fontsize{9pt}{11pt}\selectfont
\textbf{\small EXAMPLE} \  
Every $H\F_p$-local space is $p$-local (cf. Proposition 13 and $\S 8$, Proposition 3), so there is a natural transformation 
$L_p \ra L_{H\F_p}$.  
If \mG is a finite group, then $K(G,1)_p \approx K(G,1)_{H\F_p}$ but if \mG is infinite, this is false 
(consider $G = \Z$).\\
\endgroup %%------------------------------------<<

\begingroup%%----------------------------------->>
\fontsize{9pt}{11pt}\selectfont
\textbf{\small EXAMPLE} \  
Suppose that \mX is a pointed nilpotent CW space $-$then $X_{H\F_p}$ is nilpotent and $\forall \ n \geq 1$, there is a split short exact sequence 
$0 \ra$ 
$\Ext(\Z/p^\infty\Z,\pi_n(X)) \ra$ 
$\pi_n(X_{H\F_p}) \ra$ 
$\Hom(\Z/p^\infty\Z,\pi_{n-1}(X)) \ra 0$ (see below).  
Therefore, even in the nilpotent case, it need not be true that 
$\pi_n(X)_{H\F_p}$ ``is'' $\pi_n(X_{H\F_p})$ when $n > 1$.\\
\endgroup %%------------------------------------<<

\label{11.1}
\index{Theorem: \textit{H}$\F_p$ Whitehead Theorem}
\index{H$\F_p$ Whitehead Theorem}%\index{\textit{H}$\F_p$ Whitehead Theorem}
\textbf{\small \textit{H}$\F_p$ WHITEHEAD THEOREM} \quad 
Suppose that \mX and \mY are $H\F_p$-local and let $f:X \ra Y$ be a pointed continuous function.  
Assume: 
$f_*:H_q(X;\F_p) \ra H_q(Y;\F_p)$ is bijective for $1 \leq q < n$ and surjective for $q = n$ 
$-$then $f$ is an $n$-equivalence.

[The proof is the same as that of the $HP$ Whitehead theorem.]\\

\label{9.121}
\begingroup%%----------------------------------->>
\fontsize{9pt}{11pt}\selectfont
\textbf{\small EXAMPLE} \  
Let \mX be a pointed connected CW space.  Assume: $\widetilde{H}_*(X;\F_p) = 0$, i.e., \mX is $\F_p$-acyclic $-$then 
$X_{H\F_p}$ is contractible.
\\ \indent
[Note: \ A pointed nilpotent CW space \mX is $\F_p$-acyclic iff $\forall \ n \geq 1$, 
$\Hom(\Z/p^\infty\Z,\pi_{n}(X)) = 0$ $\&$ 
$\Ext(\Z/p^\infty\Z,\pi_n(X)) = 0$ 
(cf. p. \pageref{9.77}).]\\
\endgroup %%------------------------------------<<

\begin{proposition} \ %20
Let \mZ be a pointed nilpotent CW space $-$then \mZ is $H\F_p$-local iff $\forall \ n \geq 1$, $\pi_n(Z)$ is $p$-cotorsion.
\end{proposition}

[Necessity: Since \mZ satisfies Bousfield's condition mod $p$ (cf. Proposition 19), the $\pi_n(Z)$ are $H\F_p$-local, 
hence are $p$-cotorsion (cf. $\S 8$, Proposition 32).

Sufficiency: The claim is that for every $H\F_p$-equivalence $f:X \ra Y$, the precomposition arrow 
$f^*:[Y,Z] \ra [X,Z]$ is bijective.  For this, one can assume that 
$
\begin{cases}
\ X\\
\ Y
\end{cases}
$
are pointed connected CW complexes with \mX a pointed subcomplex of \mY and argue as in the proof of Proposition 2.  
However, it is no longer possible to work with the 
$\Gamma_{\chi_q}^i(\pi_q(Z))/\Gamma_{\chi_q}^{i+1}(\pi_q(Z))$ (since they need not be $p$-cotorsion).  
Instead, one uses the 
$C_{\chi_q}^i(\pi_q(Z))/C_{\chi_q}^{i+1}(\pi_q(Z))$ (which are $p$-cotorsion) (cf. $\S 8$, Proposition 34).  
Thus now, 
$\forall \ n \geq 1$, $H_n(Y,X;\F_p) = 0$ $\implies$ 
$H^p(Y,X;C_{\chi_q}^i(\pi_q(Z))/C_{\chi_q}^{i+1}(\pi_q(Z))) = 0$ (cf. $\S 8$ Proposition 29) and the obvious modification 
of the nilpotent obstruction theorem is applicable.]\\

%%----------------------------------------------------------------------------------------------34

\begingroup%%----------------------------------->>
\fontsize{9pt}{11pt}\selectfont
\textbf{\small EXAMPLE} \  
Fix a prime $p$ $-$then every $H\F_p$-local pointed nilpotent CW space is $\F_q$-acyclic for all primes $q \neq p$.
\\ \indent
[A $p$-cotorsion nilpotent group is uniquely $q$-divisible for all primes $q \neq p$ 
(cf. p. \pageref{9.78}).]\\
\endgroup %%------------------------------------<<

\textbf{\small LEMMA}  \  
Let \mF be a free abelian group $-$then the arrow 
$K(F,n) \ra K(\widehat{F}_p,n)$ is an $H\F_p$-equivalence.

[Since $\widehat{F}_p/F$ is uniquely $p$-divisible, $K(\widehat{F}_p/F,n)$ is $\F_p$-acyclic.  
On the other hand, 
$K(F,n)$ is the mapping fiber of the arrow 
$K(\widehat{F}_p,n) \ra K(\widehat{F}_p/F,n)$, so 
$H_*(F,n;\F_p) \approx H_*(\widehat{F}_p,n;\F_p)$ 
(cf. p. \pageref{9.79}).]

[Note: \ $\widehat{F}_p$ is the $p$-adic completion of \mF.  Since \mF is torsion free, 
$\Ext(\Z/p^\infty \Z,F) \approx \widehat{F}_p$ 
(cf. p. \pageref{9.80}).]\\

Let \mG be an abelian group.  Fix a presentation 
$0 \ra R \ra F \ra G \ra 0$ of \mG, i.e., a short exact sequence with \mR and \mF free abelian $-$then there is an exact sequence 
$0 \ra$ 
$\Hom(\Z/p^\infty\Z,G) \ra$ 
$\Ext(\Z/p^\infty\Z,R) \ra$
$\Ext(\Z/p^\infty\Z,F) \ra$
$\Ext(\Z/p^\infty\Z,G) \ra 0$ or still, an exact sequence
$0 \ra$ 
$\Hom(\Z/p^\infty\Z,G) \ra$ 
$\widehat{R}_p \ra$ 
$\widehat{F}_p \ra$ 
$\Ext(\Z/p^\infty\Z,G) \ra 0$.  Consider the following diagram
\[
\begin{tikzcd}[sep=large]
{K(F,n)} \ar{d} \ar{r} 
&{K(G,n)} \ar[dashed]{d} \ar{r}
&{K(R,n+1)} \ar{d} \ar{r}
&{K(F,n+1)} \ar{d}\\
{K(\widehat{F}_p,n)} \ar{r} 
&{K(G,n)_{H\F_p}} \ar{r} 
&{K(\widehat{R}_p,n+1)} \ar{r} 
&{K(\widehat{F}_p,n+1)}
\end{tikzcd}
,
\]
where by definition $K(G,n)_{H\F_p}$ is the mapping fiber of the arrow 
$K(\widehat{R}_p,n+1) \ra$ 
$K(\widehat{F}_p,n+1)$.  To justify the notation, first note that $K(G,n)_{H\F_p}$ has two nontrivial homotopy groups, namely 
$\pi_n(K(G,n)_{H\F_p}) \approx \Ext(\Z/p^\infty\Z,G)$ and 
$\pi_{n+1}(K(G,n)_{H\F_p}) \approx \Hom(\Z/p^\infty\Z,G)$.  
Since both of these groups are $p$-cotorsion, Proposition 20 implies that $K(G,n)_{H\F_p}$ is $H\F_p$-local.  
Taking into account the lemma, standard spectral sequence generalities allow one to infer that the filler 
\begin{tikzcd}[sep=small]
{K(G,n)} \ar[dashed]{r} &{K(G,n)_{H\F_p}}
\end{tikzcd}
is an $H\F_p$-equivalence.  
Therefore $K(G,n)_{H\F_p}$ is the $H\F_p$-localization of $K(G,n)$.  
Example: $K(\Q,n)_{H\F_p} \approx *$.\\

\begingroup%%----------------------------------->>
\fontsize{9pt}{11pt}\selectfont
\textbf{\small EXAMPLE} \  
Suppose that $G = \Z/p^\infty\Z$ $-$then 
$K(\Z/p^\infty\Z,n)_{H\F_p} \approx K(n+1,\widehat{Z}_p)$  
(cf. p. \pageref{9.81}).\\
\endgroup %%------------------------------------<<

Let \mX be a pointed nilpotent CW space.  
Thanks to the preceding considerations, one can copy the proof of the nilpotent 
$P$-localization theorem to see that $X_{H\F_p}$ is nilpotent.  
In so doing, one finds that there is a short exact sequence 
$0 \ra$ 
$\Ext(\Z/p^\infty\Z,\pi_n(X)) \ra$ 
%%----------------------------------------------------------------------------------------------35
$\pi_n(X_{H\F_p}) \ra \Hom(\Z/p^\infty\Z,\pi_{n-1}(X)) \ra 0$ which necessarily splits 
(Ext (torsion free, $p$-cotorsion) = 0 
(cf. p. \pageref{9.82})).  Moreover, the triangle
\begin{tikzcd}[sep=large]
{\pi_n(X)} \ar{rd} \ar{r} &{\Ext(\Z/p^\infty\Z,\pi_n(X)) } \ar{d}\\
&{\pi_n(X_{H\F_p})}
\end{tikzcd}
commutes.  
When the homotopy groups of \mX are finitely generated, it is common to write $\widehat{X}_p$ in place of 
$X_{H\F_p}$ and refer to $\widehat{X}_p$ as the 
\un{$p$-adic completion}
\index{p-adic completion} 
of \mX, the rationale being that in this case, $\forall \ n$, 
$\pi_n(\widehat{X}_p) \approx$ 
$\pi_n(X)\overset{\raisebox{-0.05cm}{$\widehat{ \ }$}}{_p}$ %repair orig to notation of para 10.  ie here orig has a wedge - not good
(cf. p. \pageref{9.83}).

Observation: Let \mX be a pointed nilpotent CW space $-$then $\forall \ p \in P$, 
$(X_P)_{H\F_p} \approx X_{H\F_p}$ and $\forall \ p \notin P$, $(X_P)_{H\F_p} \approx *$.\\

\begingroup%%----------------------------------->>
\fontsize{9pt}{11pt}\selectfont
\textbf{\small EXAMPLE} \  
Given $n \geq 1$, 
%$\widehat{\bS} {}_p^n$
$[\widehat{\bS} {}_p^n,\widehat{\bS} {}_p^n] \approx$ 
$[\bS^n,\widehat{\bS} {}_p^n] \approx \pi_n(\widehat{\bS} {}_p^n) \approx 
\widehat{\Z}_p$, the $p$-adic integers.  
This correspondence is an isomorphism of rings, thus a pointed homotopy equivalence 
$\widehat{\bS} {}_p^n \ra \widehat{\bS} {}_p^n$ determines a $p$-adic unit (i.e., in the notation of 
p. \pageref{9.84}, an element of $\widehat{\bU}_p$) and vice versa.
\\ \indent
[Note: \ $\bS_P^n = M(\Z_p,n)$ but $\widehat{\bS} {}_p^n \neq M(\widehat{\Z}_p,n)$.]\\
\endgroup %%------------------------------------<<

\begingroup%%----------------------------------->>
\fontsize{9pt}{11pt}\selectfont
\textbf{\small LEMMA}  \  
Let \mG be a finite group whose order is prime to $p$.  Suppose that \mX is a path connected free right $G$-space $-$then 
$H^*(X/G;\F_p) \approx H^*(X;\F_p)^G$.\\
\endgroup %%------------------------------------<<

\begingroup%%----------------------------------->>
\fontsize{9pt}{11pt}\selectfont
\index{Sullivan's Loop Space (example)}
\textbf{\small EXAMPLE\  (\un{Sullivan's Loop Space})} \  
Assume that $p$ is odd and that $n$ divides $p-1$ $-$then 
$\widehat{\bS} {}_p^{2n-1}$ has the pointed homotopy type of a loop space.  
This is seen as follows.  
Since 
$\widehat{\bU}_p \approx$ 
$\Z/(p-1) \Z \oplus \widehat{Z}_p$ 
(cf. p. \pageref{9.85}), $\Z/n\Z$ $(\subset \Z/(p-1)\Z)$ operates on $\widehat{\Z}_p$ 
(but the action is not nilpotent).  
Realize $K(\widehat{\Z}_p,2)$ per 
p. \pageref{9.86} and form 
$K(\widehat{\Z}_p,2;\chi) =$ 
$(\widetilde{K}(\Z/n\Z,1) \times K(\widehat{\Z}_p,2))/(\Z/n\Z)$, where 
$\chi:\Z/n\Z \ra \Aut \widehat{\Z}_p$ (thus 
$\pi_1(K(\widehat{\Z}_p,2;\chi)) \approx \Z/n\Z$ and 
$\pi_2(K(\widehat{\Z}_p,2;\chi)) \approx \widehat{\Z}_p$).  
Since 
$H^*(\widehat{\Z}_p,2;\F_p) \approx \F_p[t]$ $(\abs{t} = 2)$, the lemma implies that 
$H^*(\widehat{\Z}_p,2;\chi;\F_p) \approx \F_p[t]$ $(\abs{t} = 2n)$.  
Fix a pointed continuous function 
$f:\bP^2(n) \ra K(\widehat{\Z}_p,2;\chi)$ which induces an isomorphism of fundamental groups 
$(\bP^2(n) = M(\Z/n\Z,1)$ 
(cf. p. \pageref{9.87}) $-$then $C_f$ is simply connected (Van Kampen) and the arrow 
$K(\widehat{\Z}_p,2;\chi) \ra C_f$ is an $H\F_p$-equivalence, hence 
$K(\widehat{\Z}_p,2;\chi)_{H\F_p} \approx (C_f)_{H\F_p} \equiv B$.
\\ \indent
Claim: \mB is $(2n-1)$-connected.
\\ \indent
[$H^q(B;\F_p) = 0$ $(1 \leq q < 2n)$ $\implies$ 
$H_q(B) \otimes \F_p = 0$ $(1 \leq q < 2n)$ $\&$ $\pi_1(B) = *$ $\implies$ 
$\pi_2(B) \approx H_2(B)$ (Hurewicz $\implies$ $\pi_2(B) = 0$ ($\pi_2(B)$ is $p$-cotorsion and $p$-divisible), so by iteration, \
$\pi_q(B) = 0$ $(1 \leq q < 2n)$.]
\\ \indent
The cohomology algebra $H^*(\Omega B;\F_p)$ is an exterior algebra on one generator of degree $2n - 1$ and there is an 
$H\F_p$-equivalence $\bS^{2n-1} \ra \Omega B$.  Accordingly, 
$\widehat{\bS} {}_p^{2n-1} \approx \Omega B$, $\Omega B$ being $H\F_p$-local 
(cf. p. \pageref{9.88}).\\
\endgroup %%------------------------------------<<

\label{11.2}
\begingroup%%----------------------------------->>
\fontsize{9pt}{11pt}\selectfont
\textbf{\small EXAMPLE} \  
Let \mA be a ring with unit $-$then $B\bGL(A)^+$ is nilpotent (in fact, abelian 
(cf. p. \pageref{9.89} ff.).  Supposing that the $K_n(A)$ are finitely generated, $\forall \ n \geq 1$,
$\pi_n(B\bGL(A)_{H\F_p}^+) \approx$ 
$\Ext(\Z/p^\infty\Z,K_n(A)) \approx$ 
$\widehat{Z}_p \otimes K_n(A)$.
\\ \indent
%%----------------------------------------------------------------------------------------------36
[Note: \ This assumption is in force whenever \mA is a finite field 
(Quillen\footnote[2]{\textit{Ann. of Math.} \textbf{96} (1972), 552-586.}) or the ring of integers in an algebraic number field 
(Quillen\footnote[3]{\textit{SLN} \textbf{341} (1973), 179-198.}).]\\
\endgroup %%------------------------------------<<

\begingroup%%----------------------------------->>
\fontsize{9pt}{11pt}\selectfont
\textbf{\small FACT} \  \ 
Suppose that \mX is a pointed simply connected CW space which is $H\F_p$-local $-$then $H^n(X;\widehat{\Z}_p)$ is a finite 
$p$-group $\forall \ n \geq 1$ iff $\pi_n(X)$ is a finite $p$-group $\forall \ n \geq 1$.

[Since $X_{\Q}$ is $\F_p$-acyclic, the projection $E_{l_{\Q}} \ra X$ is an $H\F_p$-equivalence, so 
$(E_{l_{\Q}})_{H\F_p} \approx X$.  
In addition, the homotopy groups of \mX are $p$-cotorsion, thus are uniquely $q$-divisible for all primes $q \neq p$.  Therefore the $\pi_n(E_{l_{\Q}})$ are $p$-primary.  
The mod $\sC$ Hurewicz theorem then implies that 
$\forall \ n \geq 1$, $H_n(E_{l_{\Q}})$ is $p$-primary ($E_{l_{\Q}}$ is abelian).  
Finally, if the homotopy groups of either 
$E_{l_{\Q}}$ or \mX are finite $p$-groups, then $E_{l_{\Q}} \approx X$.]\\
\endgroup %%------------------------------------<<

\begin{proposition} \ %21
Let $[f]:X \ra Y$ be a morphism in $\bHCONCWSP_*$.  Assume: $[f]$ is orthogonal to every $H\F_p$-local pointed connected CW space $-$then $[f]$ is an $H\F_p$-equivalence.
\end{proposition}

[This is the $H\F_p$ version of Proposition 14 and is proved in the same way (cf. $\S 8$, Proposition 29).]\\

Given a set of primes \mP, put $\F_P = \bigoplus\limits_{p \in P} \F_p$.\\

\begin{proposition} \ %22
Let \mX be a pointed nilpotent CW space $-$then $\forall \ P$, $X_{H\F_p}$ is nilpotent and 
$X_{H\F_p} \approx \prod\limits_{p \in P} X_{H\F_p}$.
\end{proposition}

[Extending the algebra of $p$-cotorsion abelian or nilpotent groups to a \mP-cotorsion theory is a formality.  
The other point is that the product may be infinite, hence has to be interpreted as on 
p. \pageref{9.90}.]\\

\begingroup%%----------------------------------->>
\fontsize{9pt}{11pt}\selectfont
\index{arithmetic square}
\textbf{\small EXAMPLE \  (\un{Arithmetic Square})} \  
Suppose that \mX is a pointed nilpotent CW space $-$then for any \mP, the diagram 
\begin{tikzcd}[sep=large]
{X_P} \ar{d} \ar{r} &{X_{H\F_p}} \ar{d}\\
{(X_P)_{\Q}} \ar{r} &{(X_{H\F_p})_{\Q}}
\end{tikzcd}
is pointed homotopy commutative and $X_P$ ``is'' the double mapping track of the pointed 2-sink 
$(X_P)_{\Q} \ra (X_{H\F_p})_{\Q} \la X_{H\F_p}$
(Dror-Dwyer-Kan\footnote[6]{\textit{Illinois J. Math.} \textbf{21} (1977), 242-254.}).
\\ \indent
[Note: \ When $P = \bPi$, the result asserts that \mX ``is'' the double mapping track of the pointed 2-sink 
$X_{\Q} \ra \bigl(\ds\prod\limits_p X_{H\F_p}\bigr)_{\Q} \la \ds\prod\limits_p X_{H\F_p}$.  
Replacing the $X_{H\F_p}$ by the $X_p$, it can also be shown that \mX ``is'' the 
%%----------------------------------------------------------------------------------------------37
double mapping track of the pointed 2-sink 
$X_{\Q} \ra \bigl(\ds\prod\limits_p X_p\bigr)_{\Q} \la \ds\prod\limits_p X_p$.  
(Hilton-Mislin\footnote[5]{\textit{Comment. Math. Helv.} \textbf{50} (1975), 477-491.}).]\\
\endgroup %%------------------------------------<<

\begin{proposition} \ %23
Let \mG be an abelian group.  Suppose that 
$
\begin{cases}
\ f\\
\ g
\end{cases}
$
are $HG$-equivalences $-$then so is $f \times g$.\\
\end{proposition}

Application: Let
$
\begin{cases}
\ X\\
\ Y
\end{cases}
$
be pointed connected CW spaces $-$then $(X \times Y)_{HG} \approx X_{HG} \times Y_{HG}$.

[Note: The product $X_{HG} \times Y_{HG}$ is, a priori, $HG$-local.]\\

\begin{proposition} \ %24
Let \mG be an abelian group.  Suppose that  $X \overset{f}{\ra} Z \overset{g}{\la} Y$ is a pointed 2-sink, where 
$
\begin{cases}
\ X\\
\ Y
\end{cases}
\& \ Z
$
are $HG$-local pointed connected CW spaces $-$then the path component $W_0$ of $W_{f,g}$ which contains the base point 
$(x_0,y_0,j(z_0))$ is $HG$-local.
\end{proposition}

[It suffices to prove that if \mK is a pointed connected CW complex and $L \subset K$ $(L \neq K)$ is a pointed connected subcomplex such that $H_*(K,L;G) = 0$, then any pointed continuous function $\phi:L \ra W_0$ admits a pointed continuous extension $\Phi:K \ra W_0$.  Thus write $\phi = (x_\phi,y_\phi,\tau_\phi)$ and view $\tau_\phi$ as a pointed homotopy 
$I(L,l_0) \ra Z$ between $f \circx x_\phi$ and $g \circx y_\phi$ (note that $\phi(l_0) = (x_0,y_0,j(z_0))$).  Fix pointed continuous functions 
$
\begin{cases}
\ x_\Phi:K \ra X\\
\ y_\Phi:K \ra Y
\end{cases}
$
extending 
$
\begin{cases}
\ x_\phi\\
\ y_\phi
\end{cases}
$
and define $H:i_0K \cup I(L,l_0) \cup i_1K \ra Z$ accordingly 
$
\bigl(
\begin{cases}
\ X\\
\ Y
\end{cases}
\text{are } HG\text{-local}\bigr).
$
Since the inclusion 
$i_0K \cup I(L,l_0) \cup i_1K \ra I(K,k_0)$ is an $HG$-equivalence and \mZ is $HG$-local, \mH can be extended to 
$\tau_\phi:I(K,k_0) \ra Z$.  Therefore one can take $\Phi = (x_\Phi, y_\Phi, \tau_\Phi)$.]\\

\label{9.88}
Application: For any $HG$-local pointed connected CW space \mX, the path component $\Omega_0 X$ of $\Omega X$ which contains the constant loop is $HG$-local.\\

Notation: \ Given compactly generated Hausdorff spaces
$
\begin{cases}
\ X\\
\ Y
\end{cases}
\hspace{-.25cm}
, \ 
$
put $\map(X,Y) = kC(X,Y)$, where $C(X,Y)$ carries the compact open topology 
(cf. p. \pageref{9.91}).

[Note: \ If 
$
\begin{cases}
\ (X,x_0)\\
\ (Y,y_0)
\end{cases}
$
are pointed compactly generated Hausdorff spaces, then $\map_*(X,Y)$ is the closed subpace of $\map(X,Y)$ consisting of the base point preserving continuous functions.]\\

%%----------------------------------------------------------------------------------------------38
\label{11.13}
Let 
$
\begin{cases}
\ (X,x_0)\\
\ (Y,y_0)
\end{cases}
$
be pointed connected CW spaces.  Consider $C(X,x_0;Y,y_0)$ (compact open topology) $-$then the pointed homotopy type of 
$C(X,x_0;Y,y_0)$ depends only on the pointed homotopy types of $(X,x_0)$ and $(Y,y_0)$ 
(cf. p. \pageref{9.92}).  
Therefore, when dealing with questions involving the pointed homotopy type of $C(X,x_0;Y,y_0)$, one can always assume that 
$(X,x_0)$ and $(Y,y_0)$ are pointed connected CW complexes, hence are wellpointed compactly generated Hausdorff spaces.  
Of course, the homotopy type of $\map_*(X,Y)$ is not necessarily that of $C(X,x_0;Y,y_0)$ but the arrow 
$\map_*(X,Y) \ra C(X,x_0;Y,y_0)$ is at least a weak homotopy equivalence 
(cf. p. \pageref{9.93})

[Note: \ The evaluation map$ f \ra f(x_0)$ defines a \bCG fibration $\map(X,Y) \ra Y$ whose fiber over $y_0$ is 
$\map_*(X,Y)$.]

Observation: \ If $\pi_q(\map_*(X,Y))$ is computed on the path component containing the constant map, then 
$\pi_q(\map_*(X,Y)) \approx [\Sigma^qX,Y]$.

Examples: 
(1) $\forall$ $HP$-local \mX, 
$\pi_q(\map_*(\bS_{HP}^n,X)) \approx$ 
$\pi_{n+q}(X)$ $(\Sigma^q \bS_{HP}^n \approx \bS_{HP}^{n+q})$; 
(2) $\forall$ $H\F_p$-local \mX, 
$\pi_q(\map_*(\bS_{H\F_p}^n,X)) \approx$ 
$\pi_{n+q}(X)$ $((\Sigma^q \bS_{H\F_p}^n)_{H\F_p} \approx \bS_{H\F_p}^{n+q})$.\\

\begingroup%%----------------------------------->>
\fontsize{9pt}{11pt}\selectfont
Let $(X,x_0)$ be a pointed connected CW space $-$then  $(X,x_0)$ is nondegenerate 
(cf. p. \pageref{9.94}), thus satisfies Puppe's condition (cf. $\S 3$, Proposition 20).  On the other hand, the identity map $kX \ra X$ is a homotopy equivalence 
(cf. p. \pageref{9.95}).  Moreover $(kX,x_0)$ satisifes Puppe's condition.  Therefore $(kX,x_0)$ is nondegenerate 
(cf. $\S 3$, Proposition 20) and the identity map $kX \ra X$ is a pointed homotopy equivalence 
(cf. p. \pageref{9.96}).\\
\endgroup %%------------------------------------<<

\begin{proposition} \ %25
Fix an abelian group \mG.  Let 
$
\begin{cases}
\ (X,x_0)\\
\ (Y,y_0)
\end{cases}
$
be pointed connected CW spaces, $f:X \ra Y$ a pointed continuous function.  
Assume: $f$ is an $HG$-equivalence $-$then for any $HG$-local pointed connected CW space $(Z,z_0)$ the precomposition arrow 
$f^* : C(Y,y_0;Z,z_0) \ra C(X,x_0;Z,z_0)$ is a weak homotopy equivalence.
\end{proposition}

[Make the transition spelled out above and consider instead 
$f^*:\map_*(Y,Z) \ra$ $\map_*(X,Z)$, there being no loss of generality in supposing that $f$ is an inclusion.  Since 
$
\begin{cases}
\ \map(Y,Z) \ra Z\\
\ \map(X,Z) \ra Z
\end{cases}
$
are \bCG fibrations, thus are Serre, and since the diagram\\
\begin{tikzcd}%[sep=large]
{\map_*(Y,Z)} \ar{d} \ar{r} 
&{\map(Y,Z)} \ar{d} \ar{r}
&{Z}\ar[equals]{d}\\
{\map_*(X,Z)} \ar{r}
&{\map(X,Z)} \ar{r}
&{Z}
\end{tikzcd}
commutes, it need only be shown that 
%$f^*:\map_*(Y,Z) \ra \map_*(X,Z)$ %dmc Garth has infra - but is the correct?
$f^*:\map(Y,Z) \ra \map(X,Z)$
 is a weak homotopy equivalence 
(cf. p. \pageref{9.97} ff,.).  
Claim: $\forall$ finite connected CW pair $(K,L)$, the diagram 
\begin{tikzcd}%[sep=large]
{L} \ar{d} \ar{r}{\phi} &{\map(Y,Z)} \ar{d}{f^*}\\
{K} \ar{r}[swap]{\psi} &{\map(X,Z)}
\end{tikzcd}
admits a 
%%----------------------------------------------------------------------------------------------39
filler $\Phi:K \ra \map(Y,Z)$ such that $\restr{\Phi}{L} = \phi$ and $f^* \circx \Phi = \psi$.  
For this, convert to 
\ 
$
\begin{tikzcd}%[sep=large]
{K \times X \cup L \times Y} \ar{rd} \ar{r}{i} &{K \times Y} \ar[dashed]{d}\\
&{Z}
\end{tikzcd}
.  \ 
$
Because $i$ is a cofibration (cf. $\S 3$, Proposition 7) and an $HG$ equivalence (Mayer-Vietoris), there exists an arrow 
$K \times Y \ra Z$ rendering the triangle strictly commutative.  Now quote the WHE criterion.]

[Note: \ The fact that \mZ is $HG$-local gives $[Y,Z] \approx [X,Z]$, i.e., 
$\pi_0(\map_*(Y,Z)) \approx \pi_0(\map_*(X,Z))$, so $f^*$ automatically induces a bijection of path components.]\\

Application: Fix an abelian group \mG.  Let 
$
\begin{cases}
\ (X,x_0)\\
\ (Y,y_0)
\end{cases}
$
be pointed connected CW spaces.  Assume: \mX is $HG$-acyclic and \mY is $HG$-local $-$then $C(X,x_0;Y,y_0)$ is homotopically trivial.

[The constant map $X \ra x_0$ is an $HG$-equivalence.]\\

\label{14.50}
\label{14.85}
\label{16.36}
\textbf{\small LEMMA}  \  
Let 
$
\begin{cases}
\ X\\
\ Y
\end{cases}
$
be topological spaces, $f:X \ra Y$ a continuous function.  
Assume: $f$ is a weak homotopy equivalence $-$then for any CW complex \mZ, the postcomposition arrow 
$f_*:C(Z,X) \ra C(Z,Y)$ is a weak homotopy equivalence.

[Given a finite CW pair $(K,L)$ convert
\[
\begin{tikzcd}%[sep=large]
{L} \ar{d} \ar{r} &{C(Z,X)} \ar{d}{f_*}\\
{K} \ar[dashed]{ru} \ar{r} &{C(Z,Y)}
\end{tikzcd}
\text{\indent to \indent}
\begin{tikzcd}%[sep=large]
{L \times Z} \ar{d} \ar{r} &{X} \ar{d}{f}\\
{K \times Z} \ar[dashed]{ru} \ar{r} &{Y}
\end{tikzcd}
\quad .
\]
This is permissible:
$
\begin{cases}
\ L \times Z\\
\ K \times Z
\end{cases}
$
are CW complexes, hence are compactly generated Hausdorff spaces.  Accordingly, the arrows
$
\begin{cases}
\ C(L \times Z,X) \ra C(L,C(Z,X))\\
\ C(K \times Z,Y) \ra C(K,C(Z,Y))
\end{cases}
$
are homeomorphisms (compact open topology) 
(Engleking\footnote[2]{\textit{General Topology}, Heldermann Verlag (1989), 160.}.)]

\label{13.75}
[Note: \ Let 
$
\begin{cases}
\ X\\
\ Y
\end{cases}
$
be compactly generated Hausdorff spaces, $f:X \ra Y$ a continuous function.  Assume: $f$ is a weak homotopy equivalence $-$then for any CW complex \mZ, the postcomposition arrow 
$f_*:\map(Z,X) \ra \map(Z,Y)$ is a weak homotopy equivalence (same argument).  
When 
$
\begin{cases}
\ X\\
\ Y
\end{cases}
\& \ f
$
are pointed, consideration of 
%\begin{tikzcd}%[sep=large]
%{\map_*(Z,X)} \ar{d} \ar{r} 
%&{\map(Z,X)} \ar{d} \ar{r} 
%&{X} \ar{d}\\
%{\map_*(Z,Y)} \ar{r} 
%&{\map(Z,Y)} \ar{r}
%&{Y}
%\end{tikzcd}
\begin{tikzcd}%[sep=large]
{\map_*(Z,X)} \ar{d} \ar{r} 
&{\map(Z,X)} \ar{d}\\
{\map_*(Z,Y)} \ar{r} 
&{\map(Z,Y)} 
\end{tikzcd}
\begin{tikzcd}%[sep=large]
{}\ar{r} 
&{X} \ar{d}\\
{}\ar{r} 
&{Y}
\end{tikzcd}
implies that 
$f_*:\map_*(Z,X) \ra \map_*(Z,Y)$ is a weak homotopy equivalence 
(cf. p. \pageref{9.98} ff.).]\\

%%----------------------------------------------------------------------------------------------40
\label{9.116}
\begingroup%%----------------------------------->>
\fontsize{9pt}{11pt}\selectfont
\textbf{\small EXAMPLE} \  
Fix a prime $p$.  Let \mK be a pointed connected CW complex; let \mX be a pointed nilpotent CW complex.  Assume: \mK is
$\Z\left[\ds\frac{1}{p}\right]$-acyclic, i.e., $\widetilde{H}_*\left(K;\Z\left[\ds\frac{1}{p}\right]\right) = 0$ $-$then the arrow of localization 
$l_{H\F_p}:X \ra X_{H\F_p}$ induces a weak homotopy equivalence 
$\map_*(K,X) \ra \map_*(K,X_{H\F_p})$.
\endgroup %%------------------------------------<<

\begingroup%%----------------------------------->>
\fontsize{9pt}{11pt}\selectfont
[Every pointed nilpotent CW complex \mZ which is either rational or $H\F_p$-local $(q \neq p)$ is necessarily 
$H\Z\left[\ds\frac{1}{p}\right]$-local.  Therefore $\map_*(K,Z)$ is homotopically trivial.  This said, work in the compactly generated category and consider the arithmetic square 
\begin{tikzcd}[sep=large]
{X} \ar{d} \ar{r} &{L} \ar{d}\\
{X_\Q} \ar{r} &{L_\Q}
\end{tikzcd}
, where $L = X_{H\F_{\bPi}}$ $(P = \bPi)$.  Since \mX can be identified with the double mapping track of the pointed 2-sink 
$X_\Q \ra L_\Q \la L$, $\map_*(K,X)$ is the double mapping track of the pointed 2-sink 
$\map_*(K,X_\Q) \ra  \map_*(K,L_\Q) \la \map_*(K,L)$.  Because $\map_*(K,X_\Q)$ and $\map_*(K,L_\Q)$ are both homotopically trivial, the arrow $\map_*(K,X) \ra \map_*(K,L)$ is a weak homotopy equivalence 
(cf. p. \pageref{9.99}).  
However, by definition, there is a weak homotopy equivalence 
$L \ra X_{H\F_p} \times_k \ds\prod\limits_{q \neq p} X_{H\F_q}$, so from the above, the arrow 
$\map_*(K,L) \ra$ 
$\map_*(K,X_{H\F_p}) \times_k \ds\prod\limits_{q \neq p} \map_*(K,X_{H\F_q})$ is a weak homotopy equivalence.  
But 
$\ds\prod\limits_{q \neq p} \map_*(K,X_{H\F_q})$ is homotopically trivial, thus the projection 
$\map_*(K,L) \ra \map_*(K,X_{H\F_q})$ is a weak homotopy equivalence.]\\
\endgroup %%------------------------------------<<

\begingroup%%----------------------------------->>
\fontsize{9pt}{11pt}\selectfont
\textbf{\small EXAMPLE} \  
Let \mG be a finite $p$-group $-$then $BG (= K(G,1)$ 
(cf. p. \pageref{9.100})) is $\Z\left[\ds\frac{1}{p}\right]$-acyclic 
(Brown\footnote[2]{\textit{Cohomology of Groups}, Springer Verlag (1982), 84.}).  
So, for any pointed nilpotent CW space \mX, $[BG,X] \approx [BG,X_{H\F_p}]$.
\\ \indent
[Note: \ If \mX is a simply connected CW space and if the homotopy groups of \mX are finite $p$-groups, then \mX is 
$\Z\left[\ds\frac{1}{p}\right]$-acyclic.  Proof: $\forall \ n > 0$, $H_n(X)$ is a finite $p$-group (mod $\sC$ Hurewicz), hence 
$\forall n > 0$, $H_n\left(X;\Z\left[\ds\frac{1}{p}\right]\right) = \Z\left[\ds\frac{1}{p}\right] \otimes H_n(X) = 0$.]\\
\endgroup %%------------------------------------<<

\label{9.115 9}
\begingroup%%----------------------------------->>
\fontsize{9pt}{11pt}\selectfont
\textbf{\small EXAMPLE} \  
Fix a prime $p$.  Let \mX be a pointed nilpotent CW complex $-$then the arrow of localization $l_p:X \ra X_p$ induces a weak homotopy equivalence 
$\map_*(B\Z/p\Z,X) \ra \map_*(B\Z/p\Z,X_p)$.
\\ \indent
[The point is that $X_{H\F_p}$ can be identified with $(X_p)_{H\F_p}$.]\\
\endgroup %%------------------------------------<<

If 
\setstretch{1.}
$
\begin{cases}
\ A\\
\ B
\end{cases}
$
\setstretch{1.25}
are pointed connected CW complexes and if $\rho:A \ra B$ is a pointed continuous function, 
then $\rho^\perp$ need not be the object class of a reflective subcategory of $\bHCONCWSP_*$ 
(cf. p. \pageref{9.101}).  Of course, $Z \in \rho^\perp$ iff 
$\rho^*:\pi_0(C(B,b_0;Z,z_0)) \ra \pi_0(C(A,a_0;Z,z_0))$ 
is bijective and it is a fundamental point of principle that the class of \mZ for which 
$\rho^*:C(B,b_0;Z,z_0) \ra C(A,a_0;Z,z_0)$ is a weak homotopy equivalence 
is the object class of a reflective subcategory of 
$\bHCONCWSP_*$ 
(cf. p. \pageref{9.102}).  This means that the ``orthogonal subcategory problem'' 
in $\bHCONCWSP_*$ has a positive solution if the notion of 
%%----------------------------------------------------------------------------------------------41
``orthogonality'' is strengthened so as to include not just $\pi_0$ but all the $\pi_n$ $(n > 0)$ as well (cf. Proposition 25 (and its proof)).

\label{9.111}
The formalities are best handled by working in $\bCGH_*$.  In fact, it is actually more convenient to work in \bCGH.  Thus let 
$
\setstretch{1.}
\begin{cases}
\ A\\
\ B
\end{cases}
\setstretch{1.25}
$
be CW complexes, $\rho:A \ra B$ a continuous function $-$then an object \mZ in \bCGH is said to be 
\un{$\rho$-local}
\index{rho-local, $\rho$-local (object in \bCGH)} 
if 
$\rho^*:\map(B,Z) \ra \map(A,Z)$ is a weak homotopy equivalence.

[Note: \ Since the diagram
\begin{tikzcd}%[sep=large]
{\map(B,Z)} \ar{d} \ar{r} &{\map(A,Z)} \ar{d}\\
{C(B,Z)} \ar{r} &{C(A,Z)}
\end{tikzcd}
commutes and the vertical arrows are weak homotopy equivalences, $Z$ is $\rho$-local iff 
$\rho^*:C(B,Z) \ra C(A,Z)$ is a weak homotopy equivalence.]

Notation: $\rho$-\bloc is the full subcategory of \bCGH whose objects are $\rho$-local.

[Note: If 
$
\setstretch{1.}
\begin{cases}
\ \rho_1\\
\ \rho_2
\end{cases}
\setstretch{1.25}
$
are homotopic, then the same holds for 
$
\setstretch{1.}
\begin{cases}
\ \rho_1^*\\
\ \rho_2^*
\end{cases}
\setstretch{1.25}
$
(cf. p. \pageref{9.103}).  Therefore \mZ is in $\rho_1$-\bloc iff \mZ is in $\rho_2$-\bloc.]

$\rho$-\bloc is closed under the formation of products in \bCGH and is invariant under homotopy equivalence.\\

\textbf{\small LEMMA}  \  
Let
$
\setstretch{1.}
\begin{cases}
\ A\\
\ B
\end{cases}
\setstretch{1.25}
$
be pointed CW complexes, $\rho:A \ra B$ a pointed continuous function.  
Suppose that \mZ is a pointed compactly generated Hausdorff space 
$-$then $\rho^*:\map_*(B,Z) \ra \map_*(A,Z)$ is a weak homotopy equivalence if \mZ is 
$\rho$-local and conversely if $\pi_0(Z) = *$.\\

\begingroup%%----------------------------------->>
\fontsize{9pt}{11pt}\selectfont
\textbf{\small EXAMPLE} \  
Take 
$
\begin{cases}
\ A = W\\
\ B = *
\end{cases}
, \ 
$ 
where \mW is path connected, and let $\rho:W \ra *$ $-$then the $\rho$-local objects are said to be 
\un{$W$-null}.
\index{W-null}  
So, \mZ is $W$-null iff the arrow $Z \ra \map(W,Z)$ is a weak homotopy equivalence.  On the other hand, relative to some choice of a base point in \mW, a pointed path connected \mZ is $W$-null iff the arrow 
$* \ra \map_*(W,Z)$ is a weak homotopy equivalence or still, iff $\map_*(W,Z)$ is homotopically trivial, i.e., iff $\forall \ q \geq 0$, 
$[\Sigma^q W,Z] = 0$.  
Example: When $W = \bS^{n+1}$ $(n \geq 0)$, a pointed path connected \mZ is $W$-null iff $\pi_q(Z) = 0$ $(q > n)$.\\
\endgroup %%------------------------------------<<

\begingroup%%----------------------------------->>
\fontsize{9pt}{11pt}\selectfont
\textbf{\small FACT} \  
Let $f:X \ra Y$ be a \bCG fibration, where \mY is path connected.  Fix $y_0 \in Y$ and assume that $X_{y_0}$ $\&$ \mY are $W$-null $-$then \mX is $W$-null.
\\ \indent
[Observing that the arrow 
$\map(W,X) \ra \map(W,Y)$ is a \bCG fibration, consider the commutative diagram
\begin{tikzcd}[sep=large]
{X_{y_0}} \ar{d} \ar{r} 
&{X} \ar{d} \ar{r} 
&{Y} \ar{d}\\
{\map(W,X_{y_0})} \ar{r} 
&{\map(W,X)} \ar{r}
&{\map(W,Y)}
\end{tikzcd}
.]
\\ \indent
[Note: \ By the same token, \mX $\&$ \mY $W$-null $\implies$ $X_{y_0}$ $W$-null .]\\
\endgroup %%------------------------------------<<

%%----------------------------------------------------------------------------------------------42
\begin{proposition} \ %26
Let
$
\begin{cases}
\ A\\
\ B
\end{cases}
$
be CW complexes, $\rho:A \ra B$ a continuous function.  Suppose that \mZ is $\rho$-local $-$then $\forall \ Y$ in \bCW, 
$\map(Y,Z)$ is $\rho$-local.
\end{proposition}

[The arrow 
$\map(B,\map(Y,Z)) \ra \map(A,\map(Y,Z))$ is a weak homotopy equivalence iff the arrow 
$\map(B \times_k Y,Z) \ra \map(A \times_k Y,Z)$ is a weak homotopy equivalence, i.e., iff the arrow 
$\map(Y,\map(B,Z)) \ra \map(Y,\map(A,Z)$ is a weak homotopy equivalence.]\\

\textbf{\small LEMMA}  \  
Given \mX in \bCGH,
$
\begin{cases}
\ Y\\
\ Z
\end{cases}
$
in $\bCGH_*$, $\map(X,\map_*(Y,Z))$ is homeomorphic to $\map_*(Y,\map(X,Z))$.

[
$\map(X,\map_*(Y,Z)) \  \approx \  $
$\map_*(X_+,\map_*(Y,Z)) \ \approx \ $ 
$\map_*(X_+ \#_k Y,Z) \  \approx \ $
$\map_*(Y,$ 
$\map_*(X_+,Z))$ $\approx$ $\map_*(Y,\map(X,Z))$.]\\

\begin{proposition} \ %27
Let
$
\begin{cases}
\ A\\
\ B
\end{cases}
$
be pointed CW complexes, $\rho:A \ra B$ a pointed continuous function.  Suppose that \mZ is pointed and $\rho$-local $-$then $\forall \ Y$ in $\bCW_*$, 
$\map_*(Y,Z)$ is $\rho$-local.
\end{proposition}

[The arrow 
$\map(B,\map_*(Y,Z)) \ra \map(A,\map_*(Y,Z))$ is a weak homotopy equivalence iff the arrow 
$\map_*(Y,\map(B,Z)) \ra \map_*(Y,\map(A,Z))$ is a weak homotopy equivalence.]\\

\begingroup%%----------------------------------->>
\fontsize{9pt}{11pt}\selectfont
Given a pointed compactly generated Hausdorff space \mX, put 
$\Sigma_k X = X \#_k \bS^1$, 
$\Omega_k X = \map_*(\bS^1,X)$ $-$then the assignments $X \ra \Sigma_k X$, $X \ra \Omega_k X$ define functors 
$\bCGH_* \ra \bCGH_*$ and $(\Sigma_k,\Omega_k)$ is an adjoint pair.\\
\endgroup %%------------------------------------<<

\begingroup%%----------------------------------->>
\fontsize{9pt}{11pt}\selectfont
\textbf{\small EXAMPLE} \  
Let
$
\begin{cases}
\ A\\
\ B
\end{cases}
$
be pointed CW complexes, $\rho:A \ra B$ a pointed continuous function.  Suppose that \mZ is pointed and $\rho$-local $-$then 
$\Omega_k Z$ is $\rho$-local.  Therefore the arrow 
$\map_*(B,\Omega_k Z) \ra$ $\map_*(A,\Omega_k Z)$ is a weak homotopy equivalence, i.e., the arrow 
$\map_*(\Sigma_k B,Z) \ra$ $\map_*(\Sigma_k A,Z)$ is a weak homotopy equivalence, so \mZ is $\Sigma_k \rho$-local provided that \mZ is path connected.\\
\endgroup %%------------------------------------<<

\begin{proposition} \ %28
Let
$
\begin{cases}
\ A\\
\ B
\end{cases}
$
be CW complexes, $\rho:A \ra B$ a continuous function.  Suppose that $X \ra Z \la Y$ is a 2-sink of compactly generated Hausdorff spaces.  Assume: 
$
\begin{cases}
\ X\\
\ Y
\end{cases}
\& \ Z
$
are $\rho$-local $-$then the compactly generated double mapping track \mW is $\rho$-local.
\end{proposition}

[The vertical arrows in the commutative diagram 
\begin{tikzcd}%[sep=large]
{\map(B,X)} \ar{d} \ar{r} 
&{\map(B,Z)} \ar{d} & \ar{l} \\
{\map(A,X)}  \ar{r} 
&{\map(A,Z)} & \ar{l} 
\end{tikzcd}
\begin{tikzcd}%[sep=large]
&{\map(B,Y)} \ar{d}\\
&{\map(A,Y)}
\end{tikzcd}
are weak homotopy equivalences, thus the arrow  
$\map(B,W) \ra \map(A,W)$ is
%%----------------------------------------------------------------------------------------------43
a weak homotopy equivalence 
(cf. p. \pageref{9.104}).]\\

\begin{proposition} \ %29
Let
$
\begin{cases}
\ A\\
\ B
\end{cases}
$
be CW complexes, $\rho:A \ra B$ a continuous function.  Suppose that \mW is a retract of \mZ, where \mZ is $\rho$-local $-$then \mW is $\rho$-local.
\end{proposition}

[There is a \cd 
\begin{tikzcd}%[sep=large]
{\map(B,W)} \ar{d} \ar{r} 
&{\map(B,Z)} \ar{d} \ar{r} 
&{\map(B,W)} \ar{d}\\
{\map(A,W)}  \ar{r} 
&{\map(A,Z)} \ar{r} 
&{\map(A,W)}
\end{tikzcd}
in which the composite of the horizontal arrows across the top and the bottom is the respective identity map, i.e., the arrow 
$\map(B,W) \ra \map(A,W)$ is is a retract of the arrow 
$\map(B,Z) \ra \map(A,Z)$ 
(cf. p. \pageref{9.105}).  But the retract of a weak homotopy equivalence is a weak homotopy equivalence.]\\

\begingroup%%----------------------------------->>
\fontsize{9pt}{11pt}\selectfont
\textbf{\small EXAMPLE} \  
If \mZ is $\rho$-local and a CW space, then any nonempty union of its path components is again $\rho$-local.

[\mZ is the coproduct of its path components 
(cf. p. \pageref{9.106}).]\\
\endgroup %%------------------------------------<<

\label{13.82b}
\indent\indent ($(A,B)$ Construction) \quad 
Let
$
\begin{cases}
\ A\\
\ B
\end{cases}
$
be CW complexes, $\rho:A \ra B$ a continuous function.  
Because the objects in $\rho\text{-}\bloc$ depend only on $[\rho]$, there is no loss of generality in taking $\rho$ skeletal.  
The mapping cylinder $M_\rho$ of $\rho$ is then a CW complex and it is clear that the 
$\rho$-local spaces are the same as the $i$-local spaces, $i:A \ra M_\rho$ the embedding.  
One can therefore assume that \mA is a subcomplex of \mB and $\rho:A \ra B$ the inclusion (which is a closed cofibration).  
Let $(K,L)$ be 
$(\bD^n,\bS^{n-1})$ $(n \geq 0)$.  
Given an \mX in \bCGH, put $X_0 = X$ and with $f$ running over 
$\map(K \times A \cup L \times B,X_0)$, define $X_1$ by the pushout square
\[
\begin{tikzcd}%[sep=large]
{\coprod\limits_{(K,L)}\coprod\limits_f K \times A \cup L \times B} \ar{d} \ar{r} 
&{X_0}\ar{d}\\
{\coprod\limits_{(K,L)}\coprod\limits_f K \times B} \ar{r} 
&{X_1}
\end{tikzcd}
.
\]
Since $K \times A \cup L \times B \ra K \times B$ 
is a closed cofibration (cf. $\S 3$, Proposition 7), $X_0 \ra X_1$ is a closed cofibration and $X_1$ is in \bCGH 
(cf. p. \pageref{9.107}).  Proceeding, construct an expanding transfinite sequence 
$X = X_0 \subset$ 
$X_1 \subset$ 
$\cdots \subset$
$X_\alpha \subset$ 
$X_{\alpha + 1} \subset$ 
$\cdots \subset X_\kappa$ of compactly generated Hausdorff spaces by setting 
$X_\lambda = \ds \bigcup\limits_{\alpha < \lambda} X_\alpha$ at a limit ordinal $\lambda \leq \kappa$ and defining 
$X_{\alpha + 1}$ by the pushout square
\[
\begin{tikzcd}%[sep=large]
{\coprod\limits_{(K,L)}\coprod\limits_f K \times A \cup L \times B} \ar{d} \ar{r} 
&{X_\alpha}\ar{d}\\
{\coprod\limits_{(K,L)}\coprod\limits_f K \times B} \ar{r} 
&{X_{\alpha + 1}}
\end{tikzcd}
,
\]
%%----------------------------------------------------------------------------------------------44
where $f$ runs over 
$\map(K \times A \cup L \times B,X_\alpha)$.  
Here, it is understood that each $X_\lambda$ has the final topology per the $X_\alpha \ra X_\lambda$ $(\alpha < \lambda)$.  
Transfinite induction then implies that all the $X_\alpha$ $(\alpha \leq \kappa)$ are in \bCGH and every embedding 
$X_\alpha \ra X_\beta$ $(\alpha < \beta \leq \kappa)$ is a closed cofibration.  
As for $\kappa$, choose it to be a regular cardinal 
$> \sup\limits_{(K,L)} \#(K \times A \cup L \times B)$ (thus $\kappa$ is independent of \mX).  
Now fix a pair $(K,L)$.  
Claim: The arrow of restriction 
$\map(K \times B,X_\kappa) \ra \map(K \times A \cup L \times B,X_\kappa)$ is surjective.  
To see this, let 
$f:K \times A \cup L \times B \ra X_\kappa$.  
Given $x \in K \times A \cup L \times B$, $\exists$ 
$\alpha_x < \kappa$: $f(x) \in X_{\alpha_x}$ $\implies$ $\alpha = \sup\limits_x \alpha_x < \kappa$, so $f$ factors through 
$X_\alpha$, hence the claim.  
Consequently, 
$\rho^*:\map(B,X_\kappa) \ra \map(A,X_\kappa)$ is a weak homotopy equivalence 
(cf. p. \pageref{9.108}) (the arrow 
$\map(B,X_\kappa) \ra \map(A,X_\kappa)$ is a \bCG fibration (cf. $\S 4$, Proposition 6)), i.e., $X_\kappa$ is $\rho$-local.

Definition: Given an \mX in \bCGH, put $L_\rho X = X_\kappa$ $-$then this assignment defines a functor 
$L_\rho:\bCGH \ra \bCGH$ and there is a natural transformation $\id \ra L_\rho$.

[Note: \ The very construction of $L_\rho$ guarantees that the embedding $l_\rho:X \ra L_\rho X$ is a closed cofibration.]

Remarks: 
(1) $
\begin{cases}
\ A\\
\ B
\end{cases}
\& \ X
$
path connected $\implies$ $L_\rho X$ path connected; 
(2) \mX in \bCWSP $\implies$ $L_\rho X$ in \bCWSP.\\

\begin{proposition} \ %30
Let
$
\begin{cases}
\ A\\
\ B
\end{cases}
$
be CW complexes, $\rho:A \ra B$ a continuous function.  Suppose that \mZ is $\rho$-local $-$then $\forall \ X$, the arrow 
$\map(L_\rho X,Z) \ra \map(X,Z)$ is a weak homotopy equivalence.
\end{proposition}

[By definition, $L_\rho X = \underset{\alpha < \kappa}{\colimx} X_\alpha$, hence 
$\map(L_\rho X, Z) \approx \lim\limits_{\alpha < \kappa}\map(X_\alpha,Z)$ 
(homeomorphism of compactly generated Hausdorff spaces) (limit in \bCGH).  
On the other hand, the arrows in the ``long'' tower 
$\map(X_0,Z) \la $ 
$\map(X_1,Z) \la $ 
$\cdots \la $
$\map(X_\alpha,Z) \la $ 
$\map(X_{\alpha + 1},Z) \la \cdots$ \ 
are \bCG fibrations and at a limit ordinal \ $\lambda$, \ 
$\map(X_\lambda,Z) \approx \lim\limits_{\alpha < \lambda}\map(X_\alpha,Z)$, 
so it will be enough to prove that $\forall \ \alpha$, 
$\map(X_{\alpha + 1},Z) \ra \map(X_\alpha,Z)$ is a weak homotopy equivalence.  
But the commutative diagram
\[
\begin{tikzcd}%[sep=large]
{\map(X_{\alpha + 1},Z)} \ar{d} \ar{r} 
&{\map(\coprod\limits_{(K,L)}\coprod\limits_f (K \times B,Z) \ar{d}{p}}\\
{\map(X_\alpha,Z)} \ar{r} 
&{\map(\coprod\limits_{(K,L)}\coprod\limits_f (K \times A \cup L \times B,Z)}
\end{tikzcd}
\]
is a pullback square in \bCGH and $p$ is a \bCG fibration, thus one has only to show that $p$ is a weak homotopy equivalence 
(cf. p. \pageref{9.109}).  To this end, fix a pair $(K,L)$ and consider
%%----------------------------------------------------------------------------------------------45
the triangle 
\[
\begin{tikzcd}%[sep=large]
{\map(K \times B,Z)} \ar{rd} \ar{r} &{\map(K \times A \cup L \times B,Z)} \ar{d}\\
&{\map(K \times A,Z)}
\end{tikzcd}
.
\]
According to Proposition 26, the oblique arrow is a weak homotopy equivalence.  In addition, the commutative diagram
\[
\begin{tikzcd}%[sep=large]
{\map(K \times A \cup L \times B,Z)} \ar{d} \ar{r} 
&{\map(L \times B,Z)}\ar{d}\\
{\map(K \times A,Z)} \ar{r} 
&{\map(L \times A,Z)}
\end{tikzcd}
\]
is a pullback square in \bCGH and another appeal to Proposition 26 says that the \bCG fibration 
$\map(L \times B,Z) \ra \map(L \times A,Z)$ is a weak homotopy equivalence.  Therefore the arrow 
$\map(K \times A \cup L \times B,Z) \ra \map(K \times A,Z)$ is a weak homotopy equivalence 
(cf. p. \pageref{9.110}).  
Finally, then, the arrow 
$\map(K \times B,Z) \ra \map(K \times A \cup L \times B,Z)$ is a weak homotopy equivalence and our assertion follows.]\\

Application: Suppose that \mZ is $\rho$-local $-$then every diagram 
\begin{tikzcd}%[sep=large]
{X} \ar{d}[swap]{l_\rho} \ar{r}{\phi} &{Z}\\
{L_\rho X}
\end{tikzcd}
has a filler $\Phi:L_\rho X \ra Z$ in the homotopy category: $\phi \simeq \Phi \circx l_\rho$.  
And: $\Phi$ is unique up to homotopy.\\

Because $L_\rho$ is a functor $\bCGH \ra \bCGH$, given $f,g \in \map(X,Y)$, there are commutative diagrams 
\begin{tikzcd}%[sep=large]
{X} \ar{d} \ar{r}{f} &{Y}\ar{d}\\
{L_\rho X} \ar{r}[swap]{L_\rho f} &{L_\rho Y}
\end{tikzcd}
,
\begin{tikzcd}%[sep=large]
{X} \ar{d} \ar{r}{g} &{Y}\ar{d}\\
{L_\rho X} \ar{r}[swap]{L_\rho g} &{L_\rho Y}
\end{tikzcd}
.  
If further $f \simeq g$, then $L_\rho f \simeq L_\rho g$, i.e., $L_\rho$ respects the homotopy congruence.\\

\index{Theorem: Homotopical $\rho$-Localization Theorem}
\index{Homotopical $\rho$-Localization Theorem}
\textbf{\small HOMOTOPICAL $\brho$-LOCALIZATION THEOREM} \quad
Let
$
\begin{cases}
\ A\\
\ B
\end{cases}
$
be CW complexes, $\rho:A \ra B$ a continuous function.  
Let \bC be either the homotopy category of compactly generated Hausdorff spaces or the homotopy category of compactly generated CW Hausdorff spaces $-$then the full subcategory of \bC whose objects are $\rho$-local is reflective.

[Note: \ Analogous conclusions can be drawn in the path connected situation provided that 
$
\begin{cases}
\ A\\
\ B
\end{cases}
$
are themselves path connected.]\\

%%----------------------------------------------------------------------------------------------46
Let $f \in \map(X,Y)$ $-$then $f$ is said to be a 
\un{$\rho$-equivalence}
\index{rho-equivalence, $\rho$-equivalence (f $\in \map(X,Y)$} 
if $L_\rho f : L_\rho X \ra L_\rho Y$ is a homotopy equivalence.  
On general grounds, $f$ is a $\rho$-equivalence iff 
$\forall$ $\rho$-local \mZ, 
$f^*:[Y,Z] \ra [X,Z]$ is bijective.  More is true: $f$ is a $\rho$-equivalence iff 
$\forall$ $\rho$-local \mZ, 
$f^*:\map(Y,Z) \ra \map(X,Z) $ is a weak homotopy equivalence.  
Proof: Consider the commutative diagram 
\begin{tikzcd}%[sep=large]
{\map(L_\rho Y,Z)} \ar{d} \ar{r} &{\map(L_\rho X,Z)}\ar{d}\\
{\map(Y,Z)} \ar{r} &{\map(X,Z)}
\end{tikzcd}
.\\
\vspace{0.25cm}

\begingroup%%----------------------------------->>
\fontsize{9pt}{11pt}\selectfont
In the special case when $\rho:W \ra *$, and where \mW is path connected, homotopical $\rho$-localization is referred to as 
\un{\mW-nullification}
\index{\mW-nullification} 
and one writes $l_W:X \ra L_W X$ in place of 
$l_\rho:X \ra L_\rho X$, the $\rho$-equivalences being termed 
\un{\mW-equivalences}.
\index{\mW-equivalences}\\
\endgroup %%------------------------------------<<

\begin{proposition} \ %31
Let 
$
\begin{cases}
\ X\\
\ Y
\end{cases}
$
be compactly generated CW Hausdorff spaces $-$then 
$L_\rho(X \times_k Y) \approx L_\rho X \times_k L_\rho Y$.
\end{proposition}

[The product $L_\rho X \times_k L_\rho Y$ is necessarily $\rho$-local, thus it suffices to prove that the arrow 
$X \times_k Y \ra L_\rho X \times_k L_\rho Y$ is a $\rho$-equivalence.  
To see this, let \mZ be $\rho$-local.  Thanks to Proposition 26, $\map(L_\rho Y,Z)$ and $\map(X,Z)$ are $\rho$-local.  Consider the composite 
$\map(L_\rho X \times_k L_\rho Y,Z) \ra$ 
$\map(L_\rho X,\map(L_\rho Y,Z)) \ra$ 
$\map(X,\map(L_\rho Y,Z)) \ra$ 
$\map(X \times_k ,L_\rho Y,Z) \ra$ 
$\map(L_\rho Y, \map(X,Z)) \ra$ 
$\map(Y,\map(X,Z)) \ra$ 
$\map(X \times_k Y,Z)$.]

[Note: \ $L_\rho$ need not preserve arbitrary products.]\\

As it stands, base points play no role in the homotopical $\rho$-localization theorem but they can be incorporated.

Let 
$
\begin{cases}
\ A\\
\ B
\end{cases}
$
be pointed CW complexes, $\rho:A \ra B$ a pointed continuous function.  
Since $l_\rho:X \ra L_\rho X$ is a closed cofibration, $X$ wellpointed $\implies$ $L_\rho X$ wellpointed.  
Accordingly, for any $\rho$-local, wellpointed \mZ, the arrow 
$\map_*(L_\rho X,Z) \ra \map_*(X,Z)$ is a weak homotopy equivalence.  
Therefore if \bC is either the homotopy category of wellpointed compactly generated Hausdorff spaces or the homotopy category of wellpointed compactly generated CW Hausdorff spaces, then the full subcategory of \bC whose objects are $\rho$-local is reflective.

[Note: \ While the data is pointed, $\rho$-local is defined in terms of map, not $\map_*$ (but one can use $\map_*$ for path connected objects 
(cf. p. \pageref{9.111})).]

\label{9.102}
Let 
$
\begin{cases}
\ A\\
\ B
\end{cases}
$
be pointed connected CW complexes, $\rho:A \ra B$ a pointed continuous function $-$then an object \mZ in $\bCONCWSP_*$ is said to be 
\un{$\rho$-local}
\index{rho-local, $\rho$-local} 
if $\rho^*:C(B,b_0;Z,z_0) \ra C(A,a_0;Z,z_0)$ is a weak homotopy equivalence.\\

\index{Theorem: Localization Theorem of Dror Farjoun}
\index{Localization Theorem of Dror Farjoun}
\textbf{\small LOCALIZATION THEOREM OF DROR FARJOUN} \quad 
The $\rho$-local \mZ constitute the object class of a reflective subcategory of $\bHCONCWSP_*$.

%%----------------------------------------------------------------------------------------------47
[It is a question of assigning to each \mX a $\rho$-local object $L_\rho X$ and an arrow $l_\rho:X \ra L_\rho X$ such that $\forall$ $\rho$-local \mZ, $l_\rho^*$ induces a bijection $[L_\rho X,Z] \ra [X,Z]$.  
Fix a pointed CW complex 
$(\ov{X},\ov{x}_0)$ and a pointed homotopy equivalence 
$(X,x_0) \ra (\ov{X},\ov{x}_0)$.  
Definition: 
$L_\rho X = L_\rho \ov{X}$, $l_\rho:X \ra L_\rho X$ being the composite $X \ra \ov{X} \ra L_\rho \ov{X}$.

Claim: $L_\rho X$ is $\rho$-local.

[Setting $Y = L_\rho X$, by construction, the arrow $\map(B,Y) \ra \map(A,Y)$ is a weak homotopy equivalence.  
Therefore the arrow $\map_*(B,Y) \ra \map_*(A,Y)$ is a weak homotopy equivalence, so inspection of 
\begin{tikzcd}%[sep=large]
{\map_*(B,Y)} \ar{d} \ar{r} &{\map_*(A,Y)} \ar{d}\\
{C(B,b_0;Y,y_0)} \ar{r} &{C(A,a_0;Y,y_0)}
\end{tikzcd}
shows that $L_\rho X$ is $\rho$-local.]

Given a $\rho$-local \mZ, choose a pointed CW complex $(\ov{Z},\ov{z_0})$ and a pointed homotopy equivalence 
$(Z,z_0) \ra (\ov{Z},\ov{z_0})$.  
Consideration of 
\begin{tikzcd}%[sep=large]
{C(B,b_0;\ov{Z},\ov{z}_0)} \ar{d} \ar{r} &{C(A,a_0;\ov{Z},\ov{z}_0)} \ar{d}\\
{C(B,b_0;Z,z_0)} \ar{r} &{C(A,a_0;Z,z_0)}
\end{tikzcd}
and 
\begin{tikzcd}%[sep=large]
{\map_*(B,\ov{Z})} \ar{d} \ar{r} &{\map_*(A,\ov{Z})} \ar{d}\\
{C(B,b_0;\ov{Z},\ov{z_0})} \ar{r} &{C(A,a_0;Z,z_0)}
\end{tikzcd}
allows one to infer that the arrow 
$\map_*(B,\ov{Z}) \ra \map_*(A,\ov{Z})$ is a weak homotopy equivalence.  
In turn, this means that the arrow 
$\map(B,\ov{Z})$ $\ra$ $\map(A,\ov{Z})$ is a weak homotopy equivalence $(\pi_0(\ov{Z}) = *)$.  
Take now any $\phi:X \ra Z$ and chase the diagram 
\begin{tikzcd}[sep=small]
{X} \arrow[d,shift right=-1] \ar{r}{\phi}
&{Z}\arrow[r,shift right=-1]
&{\ov{Z}} \arrow[l,shift right=-1]\\
{\ov{X}}\ar{d}\arrow[u,shift right=-1]\\
{L_\rho\ov{X}}
\end{tikzcd}
to see that up to pointed homotopy, there exists a unique $\Phi:L_\rho X \ra Z$ such that $\phi \simeq \Phi \circx l_\rho$.]

[Note: \ If 
$
\begin{cases}
\ A\\
\ B
\end{cases}
$
are $n$-connected, then $\pi_q(l_\rho):\pi_q(X) \ra \pi_q(L_\rho X)$ is an isomorphism for $q \leq n$ 
(cf. p. \pageref{9.112}).]\\

\begingroup%%----------------------------------->>
\fontsize{9pt}{11pt}\selectfont
\textbf{\small EXAMPLE} \  
Consider $L_{\bS^{n+1}}$, the nullification functor corresponding to $\bS^{n+1} \ra *$ $(n \geq 0)$ $-$then, in this situation, one recovers the fact that $\bHCONCWSP_*[n]$ is a reflective subcategory of $\bHCONCWSP_*$ 
(cf. p. \pageref{9.113}), where $\forall \ X$, $L_{\bS^{n+1}} X \approx X[n]$.\\
\endgroup %%------------------------------------<<

\begingroup%%----------------------------------->>
\fontsize{9pt}{11pt}\selectfont
\textbf{\small EXAMPLE} \  
Fix a set of primes \mP.  Given a pointed connected CW space \mX, its loop space $\Omega X$ 
is a pointed CW space (loop space theorem), thus the arrow 
$
\begin{cases}
\ \Omega X \ra \Omega X\\
\ \sigma \ra \sigma^n
\end{cases}
(n \in S_P)
$
is a pointed homotopy equivalence iff it is a weak homotopy equivalence.  
To interpret this, put 
$\rho = \ds\bigvee\limits_n \rho_n$, where $\rho_n:\bS^1 \ra \bS^1$ is a map of degree $n$ $(n \in S_P)$ $-$then the $\rho$-local objects in $\bCONCWSP_*$ are precisely the objects of
%%----------------------------------------------------------------------------------------------48
$\bCONCWSP_{*,P}$ and the homotopical $P$-localization theorem is seen to be a special case of the localization theorem of 
Dror Farjoun.
\\ \indent
[Note: \ The full subcategory of $\bHCONCWSP_*$ whose objects are $P$-local in homotopy is not the object class of a reflective subcategory of $\bHCONCWSP_*$ 
(cf. p. \pageref{9.114}).  However, the full subcategory of $\bHCONCWSP_*$ whose objects are $P$-local in ``higher homotopy'' is the object class of a reflective subcategory of $\bHCONCWSP_*$.  
Proof: Consider the pointed suspension of $\rho$.  
Therefore $L_{\Sigma\rho}$ induces an isomorphism of fundamental groups and $P$-localizes the higher homotopy groups.]\\
\endgroup %%------------------------------------<<

\begingroup%%----------------------------------->>
\fontsize{9pt}{11pt}\selectfont
\textbf{\small EXAMPLE} \  
Fix an abelian group \mG.  Choose a set of CW pairs $(K_i,L_i)$, where $K_i$ is a pointed connected CW complex and 
$L_i \subset K_i$ $(L_i \neq K_i)$ is a pointed connected subcomplex such that 
$H_*(K_i,L_i;G) = 0$ subject to the restriction that the cardinality of the set of cells in $K_i$ is $\leq \#(G)$ if $\#(G)$ is infinite and 
$\leq \omega$ if $\#(G)$ if finite, which contains up to isomorphism all such CW pairs with these properties.  Let 
$\rho:\bigvee\limits_i L_i \ra \ds\bigvee\limits_i K_i$ $-$then a pointed connected CW space is $HG$-local iff it is $\rho$-local, proving once again that $\bHCONCWSP_{*,HG}$ is a  a reflective subcategory of $\bHCONCWSP_*$.
\\ \indent
[Note: \ Take $G = \Z$ and let \mW be the pointed mapping cone of $\rho$ $-$then the nullification functor $L_W$ assigns to each \mX its plus construction $X^+$.]\\
\endgroup %%------------------------------------<<

\begingroup%%----------------------------------->>
\fontsize{9pt}{11pt}\selectfont
\textbf{\small EXAMPLE} \  
Fix a prime $p$.  Let $W = M(\Z/p\Z,1)$ be the ``standard'' Moore space of type $(\Z/p\Z,1)$ $-$then a simply connected \mZ is $W$-null iff $\forall \ n \geq 2$, $\pi_n(Z)$ is $\ov{\rho}$-local.\\
\endgroup %%------------------------------------<<

\begingroup%%----------------------------------->>
\fontsize{9pt}{11pt}\selectfont
\textbf{\small EXAMPLE} \  
Fix a prime $p$.  Let $W = M\bigl(\Z\left[\ds\frac{1}{p}\right],1\bigr)$ be the ``standard'' Moore space of type 
$\bigl(\Z\left[\ds\frac{1}{p}\right],1\bigr)$ 
$-$then a simply connected \mZ is $W$-null iff $\forall \ n \geq 2$, $\pi_n(Z)$ is $p$-cotorsion.\\
\endgroup %%------------------------------------<<

\begingroup%%----------------------------------->>
\fontsize{9pt}{11pt}\selectfont
\textbf{\small EXAMPLE} \  
Fix a prime $p$.  Put $W = B\Z/p\Z$ $-$then a nilpotent \mZ is $W$-null iff $Z_p$ is $W$-null iff $Z_{H\F_p}$ is $W$-null 
(cf. p. \pageref{9.115 9}).  In general, a $W$-null \mZ is $W_k$-null, where $W_k = B\Z/p^k\Z$ $(1 \leq k < \infty)$
(consider the short exact sequence 
$0 \ra$ 
$Z/p\Z \ra$ 
$Z/p^{k+1}\Z \ra$
$Z/p^{k}\Z \ra 0$, show that the pointed mapping cone of 
$BZ/p\Z \ra$ 
$BZ/p^{k+1}\Z$
is $B\Z/p^k\Z$, and use induction (replication theorem)), hence \mZ is $W_\infty$-null, where $W_\infty = B\Z/p^\infty\Z$.
\\ \indent
[Note: The arrow $W \ra *$ is a $\ov{\rho}$-equivalence, so every $\ov{\rho}$-local space is $W$-null.  
Example: $K(\Z\left[\ds\frac{1}{p}\right],1)$ is $W$-null.]\\
\endgroup %%------------------------------------<<

\begingroup%%----------------------------------->>
\fontsize{9pt}{11pt}\selectfont
\textbf{\small LEMMA}  \  
Let 
$
\begin{cases}
\ A\\
\ B
\end{cases}
$
be pointed connected CW complexes, $\rho:A \ra B$ a pointed continuous function.  Assume: 
$\pi_1(\rho):\pi_1(A) \ra \pi_1(B)$ is surjective $-$then for any $\rho$-local \mZ, its universal covering space $\widetilde{Z}$ is 
$\rho$-local.
\\ \indent
[Note: \ Therefore $\pi_1(Z) = *$ $\implies$ $\pi_1(L_\rho Z) = *$.]\\
\endgroup %%------------------------------------<<

%%----------------------------------------------------------------------------------------------49
\begingroup%%----------------------------------->>
\fontsize{9pt}{11pt}\selectfont
\textbf{\small EXAMPLE} \  
Fix a prime $p$.  Put $W = B\Z/p\Z$ $-$then \mZ $W$-null $\implies$ $\widetilde{Z}$ $W$-null.  
Suppose now that \mX is a simply connected CW space.  
Assume: The homotopy groups of $X_{H\F_p}$ are finite $p$-groups.  
Claim: $L_W X_{H\F_p}$ is contractible if $(L_W X)_{H\F_p}$ is contractible.  
For let \mZ be $W$-null.  
Since $X_{H\F_p}$ is simply connected  and $\widetilde{Z}$ is $W$-null, one need only show that 
$[X_{H\F_p},\widetilde{Z}] = *$.  But 
$X_{H\F_p}$ is $\Z\left[\ds\frac{1}{p}\right]$-acyclic 
(cf. p. \pageref{9.116}) and $\widetilde{Z}_{H\F_p}$ is 
$W$-null (cf. supra), hence 
$[X_{H\F_p},\widetilde{Z}] \approx$ 
$[X_{H\F_p},\widetilde{Z}_{H\F_p}] \approx$ 
$[X,\widetilde{Z}_{H\F_p}] \approx$
$[L_W X,\widetilde{Z}_{H\F_p}] \approx$
$[(L_W X)_{H\F_p},\widetilde{Z}_{H\F_p}] \approx$
$[*,\widetilde{Z}_{H\F_p}] = *$.\\
\endgroup %%------------------------------------<<

\begingroup%%----------------------------------->>
\fontsize{9pt}{11pt}\selectfont
\textbf{\small LEMMA}  \  
Let $\pi$ be a group $-$then for any pointed connected CW space \mX, the path components of 
$C(X,x_0;K(\pi,1),k_{\pi,1})$ are homotopically trivial.\\
\endgroup %%------------------------------------<<

\begingroup%%----------------------------------->>
\fontsize{9pt}{11pt}\selectfont
\textbf{\small EXAMPLE} \  
Let 
$
\begin{cases}
\ A\\
\ B
\end{cases}
$
be pointed connected CW complexes, $\rho:A \ra B$ a pointed continuous function $-$then the precomposition arrow 
$\Hom(\pi_1(B),\pi) \ra \Hom(\pi_1(A),\pi)$ determined by $\pi_1(\rho)$ is bijective iff $K(\pi,1)$ is $\rho$-local.\\
\endgroup %%------------------------------------<<

\begingroup%%----------------------------------->>
\fontsize{9pt}{11pt}\selectfont
\textbf{\small EXAMPLE} \  
Fix a prime $p$.  Put $W = B\Z/p\Z$ $-$then $K(\pi,1)$ is $W$-null iff $\pi$ has no $p$-torsion.  
Example: \mZ is $W$-null provided that $\pi_1(Z)$ has no $p$-torsion and $\widetilde{Z}$ is $W$-null.\\
\endgroup %%------------------------------------<<

\begingroup%%----------------------------------->>
\fontsize{9pt}{11pt}\selectfont
\textbf{\small FACT} \  
Fix a pointed connected CW complex \mW $-$then \mW is acyclic iff $\forall \ Z$, $l_W:X \ra L_W X$ is a homology equivalence.

[Note: \ Assuming that \mW is acyclic, \mX is $W$-null iff $[W,X] = 0$.]\\
\endgroup %%------------------------------------<<

\textbf{\small LEMMA}  \  
Given $\rho_1$ $\&$ $\rho_2$, suppose that $\rho_2$ is a $\rho_1$-equivalence $-$then there exists a natural transformation 
$L_{\rho_2} \ra L_{\rho_1}$ in  $\bHCONCWSP_*$ and the class of $\rho_2$-equivalences is contained in the class of $\rho_1$-equivalences.\\

Let 
$
\begin{cases}
\ A\\
\ B
\end{cases}
$
be pointed connected CW complexes, $\rho:A \ra B$ a pointed continuous function.\\

\label{9.112}
Application: If 
$
\begin{cases}
\ A\\
\ B
\end{cases}
$
are $n$-connected, then $\pi_q(l_\rho):\pi_q(X) \ra \pi_q(L_\rho X)$ is an isomorphism for $q \leq n$.

[The class of $\rho_{n+1}$-equivalences, where $\rho_{n+1}:\bS^{n+1} \ra *$, is the class of maps $X \ra Y$ inducing isomorphisms in homotopy up to degree $n$.  But $\rho$ is a $\rho_{n+1}$-equivalence and $X \ra L_\rho X$ is a 
$\rho$-equivalence.]\\

\begingroup%%----------------------------------->>
\fontsize{9pt}{11pt}\selectfont
\textbf{\small FACT} \  
If \mW is $n$-connected, then $\pi_{n+1}(l_W):\pi_{n+1}(X) \ra \pi_{n+1}(L_W X)$ is surjective.\\
\endgroup %%------------------------------------<<

%%----------------------------------------------------------------------------------------------50
Localization theory has been developed in extenso by 
Bousfield\footnote[2]{\textit{J. Amer. Math. Soc.} \textbf{7} (1994), 831-873.}
and Dror Farjoun\footnote[3]{\textit{Cellular Spaces, Null Spaces and Homotopy Localization}, Springer Verlag (1996).}
.  
While I shall not pursue these developments in detail, let us at least set up some of the machinery without proof and see how it is used to make computations.\\

\begingroup%%----------------------------------->>
\fontsize{9pt}{11pt}\selectfont
 The simplest situation is that of $W$-nullification, where \mW is a pointed connected CW complex.\\

\index{Fibration Rule}
\textbf{\small FIBRATION RULE} \quad 
Let 
$
\begin{cases}
\ X\\
\ Y
\end{cases}
$
be pointed connected CW spaces, $f:X \ra Y$ a pointed continuous function with $\pi_0(E_f) = *$.  Suppose that 
$L_WE_f$ is contractible $-$then $f$ is a $W$-equivalence, i.e., the arrow $L_W f:L_W X \ra L_W Y$ is a pointed homotopy equivalence.\\
\endgroup %%------------------------------------<<

\begingroup%%----------------------------------->>
\fontsize{9pt}{11pt}\selectfont
\textbf{\small EXAMPLE} \  
Fix a prime $p$.  Put $W = B\Z/p\Z$ $-$then the arrow $W \ra *$ is a $W$-equivalence, thus 
$L_WK(\Z/p\Z,1)$ is contractible.  So, $\forall \ k$, 
$L_WK(\Z/p^k\Z,1)$ is contractible and this implies that
$L_WK(\Z/p^\infty\Z,1)$ is contractible.  
Examples: 
(1) From the short exact sequence 
$0 \ra$ 
$\Z \ra$ 
$\Z\left[\ds\frac{1}{p}\right] \ra$ 
$\Z/p^\infty\Z \ra 0$, $\forall \ n \geq 2$, 
$L_WK(\Z,n) \approx K\bigl(\Z\left[\ds\frac{1}{p}\right],n\bigr)$;
(2) From the short exact sequence 
$0 \ra$ 
$\widehat{\Z}_p \ra$ 
$\widehat{\Q}_p \ra$ 
$\Z/p^\infty\Z \ra 0$ 
(cf. p. \pageref{9.117}), $\forall \ n \geq 2$, 
$L_WK(\widehat{\Z}_p,n) \approx K(\widehat{\Q}_p,n)$.
\\ \indent
[Note: \ $L_WK(\pi,1)$ is contractible if $\pi$ is a finite $p$-group and, when $\pi$ is in addition abelian, 
$L_WK(\pi,n)$ is contractible as can be checked by considering 
$K(\pi,n-1) \ra \Theta K(\pi,n) \ra K(\pi,n)$.]\\
\endgroup %%------------------------------------<<

\begingroup%%----------------------------------->>
\fontsize{9pt}{11pt}\selectfont
\index{Zabrodsky Lemma}
\textbf{\small ZABRODSKY LEMMA} \quad 
Let 
$
\begin{cases}
\ X\\
\ Y
\end{cases}
\& \ Z
$
be wellpointed compactly generated connected CW Hausdorff spaces, $f:X \ra Y$ a pointed continuous function with $\pi_0(E_f) = *$.  
Assume: $\map_*(E_f, Z)$ and $\map_*(X,Z)$ are homotopically trivial $-$then $\map_*(Y,Z)$ is homotopically trivial.

[Note: \ In this setting, $E_f$ is the compactly generated mapping track.  
Its base point is $(x_0,j(y_0))$ and the inclusion 
$\{(x_0,j(y_0))\} \ra E_f$ is a closed cofibration 
(cf. p. \pageref{9.118}).]\\
\endgroup %%------------------------------------<<

\begingroup%%----------------------------------->>
\fontsize{9pt}{11pt}\selectfont
\textbf{\small EXAMPLE} \  
Miller\footnote[6]{\textit{Ann. of Math.} \textbf{120} (1984), 39-87.}
has shown that if \mG is a locally finite group, then every pointed finite dimensional connected CW complex \mZ is \mW-null, where $W = BG$.  
Using the Zabrodsky lemma, it follows by induction that for any locally finite abelian group $\pi$, all such \mZ are $K(\pi,n)$-null.
\\ \indent
[Note: \ A group is said to be 
\un{locally finite}
\index{locally finite (group)} 
if its finitely generated subgroups are finite.  
Example: Let \mX be a pointed simply connected CW space with finitely generated homotopy groups $-$then the homotopy groups of $E_{l_{\Q}}(l_{\Q}:X \ra X_{\Q})$ are locally finite.]\\
\endgroup %%------------------------------------<<

\begingroup%%----------------------------------->>
\fontsize{9pt}{11pt}\selectfont
\textbf{\small EXAMPLE} \  
Suppose that \mG is a locally finite group with the property that 
$\#\{n:H_n(G) \neq 0\} < \omega$ $-$then \mG is acyclic.
\\ \indent
%%----------------------------------------------------------------------------------------------51
[$\Sigma BG$ has the pointed homotopy type of a pointed finite dimensional connected CW complex, so by Miller, 
$[\Sigma BG,\Sigma BG] = *$.  Therefore $\Sigma BG$ is contractible, thus \mG is acyclic.]\\
\endgroup %%------------------------------------<<

\begingroup%%----------------------------------->>
\fontsize{9pt}{11pt}\selectfont
\textbf{\small EXAMPLE} \  
Miller (ibid.) has shown that if \mZ is a pointed nilpotent CW space such that $H_n(Z;\F_p) = 0$ for 
$n \gg 0$, then \mZ is $W$-null, where $W = B\Z/p\Z$.  
Example: $\forall \ n > 0$, $\bS^n$ and $\widehat{\bS}{}_p^n$ are $W$-null.\\
\endgroup %%------------------------------------<<

\begingroup%%----------------------------------->>
\fontsize{9pt}{11pt}\selectfont
\index{Preservation Rule}
\textbf{\small PRESERVATION RULE} \quad 
Let 
$
\begin{cases}
\ X\\
\ Y
\end{cases}
$
be pointed connected CW spaces, $f:X \ra Y$ a pointed continuous function with $\pi_0(E_f) = *$.  
Suppose that 
\mY is $W$-null $-$then the arrow $L_WE_f\ra E_{L_W}f$ is a pointed homotopy equivalence.
\\ \indent
[Note: \ The assumption that \mY is $W$-null can be weakened to $L_{\Sigma W} Y \approx L_W Y$.]\\[-.1cm]
\endgroup %%------------------------------------<<

\begingroup%%----------------------------------->>
\fontsize{9pt}{11pt}\selectfont
\textbf{\small EXAMPLE} \  
Let \mX be a pointed connected CW complex.  
Assume: \mX is finite and $\pi_2(X)$ is torsion $-$then $\forall \ n \geq 2$, 
$(L_W\widetilde{X}_n)_{H\F_p} = X_{H\F_p}$ ($\widetilde{X}_n$ as on 
\pageref{9.119}), where $W = B\Z/p\Z$.
\\ \indent
[Let \mE be the mapping fiber of the pointed Hurewicz fibration $\widetilde{X}_n \ra X$.  \ 
According to Miller's theorem, \mX is 
$W$-null, so 
$L_WE$ can be identified with the mapping fiber of the arrow $L_W \widetilde{X}_n \ra X$, hence 
$(L_WE)_{H\F_p}$ can be identified with the mapping fiber of the arrow $(L_W \widetilde{X}_n)_{H\F_p} \ra X_{H\F_p}$.  
Let $\ov{E}$ be the mapping fiber of the arrow of localization $l_{\ov{p}}:E \ra E_{\ov{p}}$.  Since $\pi_2(X) \approx \pi_1(E)$ and 
$\pi_2(X)$ is torsion, $\pi_1(E)$ maps onto $\pi_1(E)_{\ov{p}}$ 
(cf. p. \pageref{9.120}).  Therefore $\ov{E}$ is path connected.  
On the other hand, the nonzero homotopy groups of $\ov{E}$ are finite in number and each of them is a locally finite 
$p$-group.  From this it follows that $L_W \ov{E}$ is contractible, thus 
$L_W E \approx$ 
$L_W E_{\ov{p}} \approx$ 
$E_{\ov{p}}$.  
But the homotopy groups of $E_{\ov{p}}$ are uniquely $p$-divisible which means that $E_{\ov{p}}$ is $\F_p$-acyclic or still, that $(E_{\ov{p}})_{H\F_p}$ is contractible 
(cf. p. \pageref{9.121}).  Consequently, $(L_W E)_{H\F_p}$ is contractible and 
$(L_W \widetilde{X}_n)_{H\F_p} \approx X_{H\F_p}$.]
\\ \indent
[Note: \ Here is a numerical illustration.  
Take $X = \bS^3$ $-$then the fibers of the projection $\widetilde{X}_3 \ra X$ have the homotopy type $(\Z,2)$ and 
$L_W K(\Z,2) \approx K\bigl(\Z\left[\ds\frac{1}{p}\right],2\bigr)$, the mapping fiber of the arrow $L_W \widetilde{X}_3 \ra X$.  
The potentially nonzero homotopy groups of 
$K\bigl(\Z\left[\ds\frac{1}{p}\right],2\bigr)_{H\F_p}$ are 
$\Ext\bigl(\Z/p^\infty\Z,\Z\left[\ds\frac{1}{p}\right]\bigr)$ and 
$\Hom\bigl(\Z/p^\infty\Z,\Z\left[\ds\frac{1}{p}\right]\bigr)$, which in fact vanish, 
$\Z\left[\ds\frac{1}{p}\right]$ being uniquely $p$-divisible.  
Therefore 
$(L_W K(\Z,2))_{H\F_p}$ is contractible.  
Observe too that the mapping fiber of the arrow 
$(\widetilde{X}_3)_{H\F_p} \ra X_{H\F_p}$ is a $K(\widehat{\Z}_p,2)$.  
Because $X_{H\F_p}$ is $W$-null, 
$L_W K(\widehat{\Z}_p,2) \approx K(\widehat{\Q}_p,2)$ can be identified with the mapping fiber of the arrow 
$L_W((\widetilde{X}_3)_{H\F_p}) \ra X_{H\F_p}$.]\\[-.1cm]
\endgroup %%------------------------------------<<

\begingroup%%----------------------------------->>
\fontsize{9pt}{11pt}\selectfont
Given abelian groups \mG and \mA, call \mA 
\un{$G$-null}
\index{G-null} 
if $\Hom(G,A) = 0$.  Every abelian group \mA has a maximal $G$-null quotient $A//G$.\\[-.1cm]
\endgroup %%------------------------------------<<

\begingroup%%----------------------------------->>
\fontsize{9pt}{11pt}\selectfont
\textbf{\small EXAMPLE} \  
Fix an abelian group \mG.  Put $W = M(G,n)$ $(n \geq 2)$ and let $P_G$ be the set of primes $p$ such that \mG is uniquely divisible by $p$ ($P_G$ has the opposite meaning on 
p. \pageref{9.122}).  Let \mX be a pointed connected CW space $-$then 
$\pi_q(L_W X) \approx$ 
$\pi_q(X)$ $(q < n)$ and 
$\pi_n(L_W X) \approx \pi_n(X)//G$.  Moreover, for 
%%----------------------------------------------------------------------------------------------52
$q > n$, $\pi_q(L_W X) \approx \Z_{P_G} \otimes \pi_q(X)$ if $\Q \otimes G = 0$, while if $\Q \otimes G \neq 0$, 
there is a split short exact sequence 
$0 \ra$ 
$\ds\prod\limits_{p \in P_G} \Ext(\Z/p^\infty \Z,\pi_q(X)) \ra$ 
$\pi_q(L_W X) \ra$ 
$\ds\prod\limits_{p \in P_G} \Hom(\Z/p^\infty \Z,\pi_{q-1}(X))$ $\ra$ $0$.
\\ \indent
[\mZ is $W$-null iff $\Hom(G,\pi_q(Z)) = 0 = \Ext(G,\pi_q(Z))$ $\forall \ q > n$ and 
$\Hom(G,\pi_n(Z)) = 0$.  This said, reduce to when \mX is $(n-1)$-connected and show first that 
$\pi_n(L_W X) \approx \pi_n(X)//G$.  
Next, set 
$
S_G = 
\begin{cases}
\ \Z_{P_G} \indent\quad \  \text{ if } \Q \otimes G = 0\\
\ \bigoplus\limits_{p \in P_G} \Z/p\Z \ \ \   \text{ if } \Q \otimes G \neq 0
\end{cases}
\hspace{-.25cm}
. \ 
$
Since $\widetilde{H}_*(W;S_G) = 0$, each $HS_G$-local space is $W$-null, thus there is a natural transformation 
$L_W \ra L_{HS_G}$.  Deduce from this that 
$\pi_q(L_W X) \approx \pi_q(X_{HS_G})$ for $q > n$.]
\endgroup %%------------------------------------<<

%%%%%%%%%%%%%%%%%%%%%%%%%%%%%%%%%%%%%%
%%%%%%%%%%%%%%%%%%%%%%%%%%%%%%%%%%%%%%
%%%%%%%%%%%%%%%%%%%%%%%%%%%%%%%%%%%%%%

\begin{center}
$\S \ 9$
\\[0.5cm]
$\mathcal{REFERENCES}$\\
\end{center}

\[
\textbf{BOOKS}
\]

\begingroup
\fontsize{9pt}{11pt}\selectfont
\setlength\parindent{0 cm}

[1] \quad Adams, J., \textit{Localization and Completion}, University of Chicago (1975).
\\[-.2cm]

[2] \quad Arkowitz, M., \textit{Localization and H Spaces}, Aarhus Universitet (1976).
\\[-.2cm]

[3] \quad Bousfield, A., and Kan, D., \textit{Homotopy Limits, Completions and Localizations}, Springer Verlag (1972).
\\[-.2cm]

[4] \quad Dror Farjoun, E., \textit{Cellular Spaces, Null Spaces and Homotopy Localization}, Springer Verlag (1996).
\\[-.2cm]

[5] \quad Hilton, P., \textit{Nilpotente Gruppen und Nilpotente R\"aume}, Springer Verlag (1984).
\\[-.2cm]

[6] \quad Hilton, P., Mislin, G., and Roitberg, J., \textit{Localization of Nilpotent Groups and Spaces}, North Holland 

\hspace{0.8cm}(1975).
\\[-.2cm]

%[7] \quad Hirschhorn, P., Localization, Cellularization, and Homotopy Colimits,\\ Orig version
[7] \quad Hirschhorn, P., \textit{Model Categories and Their Localizations}, Amer. Math. Soc. (2003).
\\[-.2cm]

[8] \quad Schiffman, S., \textit{Ext $p$-Completion in the Homotopy Category}, Ph.D. Thesis, Dartmouth College, Hanover 

\hspace{0.8cm}(1974).
\\[-.2cm]

[9] \quad Sullivan, D., \textit{Geometric Topology}, MIT (1970).
\\[-.2cm]
\endgroup

\[
\textbf{ARTICLES}
\]

\begingroup
\fontsize{9pt}{11pt}\selectfont
\setlength\parindent{0 cm}

[1] \quad Bousfield, A., The Localization of Spaces with respect to Homology, 
\textit{Topology} \textbf{14} (1975), 133-150.
\\[-.2cm]

[2] \quad Bousfield, A., Homological Localization Towers for Groups and $\Pi$-Modules, 
\textit{Memoirs Amer. Math.}

\hspace{0.8cm}\textit{Soc.} \textbf{186} (1977), 1-68.
\\[-.2cm]

[3] \quad Bousfield, A., Localization and Periodicity in Unstable Homotopy Theory, 
\textit{J. Amer. Math. Soc.} \textbf{7} 

\hspace{0.8cm}(1994), 831-873.
\\[-.2cm]

[4] \quad Bousfield, A., Unstable Localization and Periodicity, In: 
\textit{New Trends in Localization and Periodicity}, 

\hspace{0.8cm}C. Broto et al. (ed.), Birkh\"auser (1996), 33-50.
\\[-.2cm]

[5] \quad Bousfield, A., Homotopical Localizations of Spaces, 
\textit{Amer. J. Math.} \textbf{119} (1997), 1321-1354.
\\[-.2cm]

[6] \quad Casacuberta, C., Recent Advances in Unstable Localization, 
\textit{CRM Proc.} \textbf{6} (1994), 1-22.
\\[-.2cm]

[7] \quad Casacuberta, C., and Peschke, G., Localizing with respect to Selfmaps of the Circle, 
\textit{Trans. Amer.} 

\hspace{0.8cm}\textit{Math. Soc.} \textbf{339} (1993), 117-140.
\\[-.2cm]

[8] \quad Dror Farjoun, E., Homotopy Localization and $v_1$-Periodic Spaces, 
\textit{SLN} \textbf{1509} (1992), 104-113.
\\[-.2cm]

[9] \quad Hilton, P., Localization in Topology, 
\textit{Amer. Math. Monthly} \textbf{82} (1975), 113-131.
\\[-.2cm]

[10] \quad May, J., Fibrewise Localization and Completion, 
\textit{Trans. Amer. Math. Soc.} \textbf{258} (1980), 127-146.
\\[-.2cm]

[11] \quad McGibbon, C., The Mislin Genus of a Space, 
\textit{CRM Proc.} \textbf{6} (1994), 75-102.
\\[-.2cm]

[12] \quad Postnikov, M., Localization of Topological Spaces, 
\textit{Russian Math. Surveys} \textbf{32} (1977), 121-184.
\\[-.2cm]

[13] \quad Sullivan, D., Genetics of Homotopy Theory and the Adams Conjecture, 
\textit{Ann. of Math.} \textbf{100} (1974), 

\hspace{0.95cm}1-79.

\setlength\parindent{2em}

\endgroup

\chapter{
$\boldsymbol{\S}$\textbf{10}.\quadx  COMPLETION OF GROUPS}
\setlength\parindent{2em}
\setcounter{proposition}{0}

%%----------------------------------------------------------------------------------------------01
$\text{ }$\\[-1.25cm]

There are many ways to ``complete'' a group.  
While various procedures are related by a web of interconnections, the theory is less systematic than that of $\S 8$, 
one reason for this being that completion functors are generally not idempotent.  
Still, the material is more or less standard, so I shall omit the details and settle for a survey of what is relevant.

Let \mG be a topological group.  
Assume: The left and right uniform structures on \mG coincide $-$then the 
\un{completion}
\index{completion (of a group)} 
$\widehat{G}$ of \mG
is the uniform completion of $G/\ov{\{e\}}$.  
Therefore $\widehat{G}$ is a uniformly complete Hausdorff topological group which is universal with respect to continuous homomorphisms $G \ra K$, where \mK is a uniformly complete Hausdorff topological group: 
\begin{tikzcd}%[sep=small]
{G} \ar{d} \ar{r} &{K}\\
{\widehat{G}} \ar[dashed]{ru}
\end{tikzcd}
.

[Note: \ The assumption is automatic if \mG is abelian.  
In this case, $\widehat{G}$ is also abelian.  
Example: Each prime $p$ determines a metrizable topology on $\Q$ and a corresponding completion $\widehat{\Q}_p$, the field of 
\un{$p$-adic numbers}.
\index{$p$-adic numbers}  
It is homeomorphic to 
$\coprod\limits_1^\infty C,$ \mC the Cantor set.]\\

\begingroup%%----------------------------------->>
\fontsize{9pt}{11pt}\selectfont
\textbf{\small EXAMPLE} \quad 
Let \mG be a group and let $\{G_i\}$ be a collection of normal subgroups of \mG directed by inclusion 
(i.e., $i \leq j$ $\Leftrightarrow G_j \subset G_i$).  
Equip \mG with the structure of a topological group by stipulating that the 
$G_i$ are to be a fundamental system of neighborhoods of $e$, thus the underlying topology is Hausdorff iff 
$\ds\bigcap\limits_i G_i = \{e\}$.  Because the $G_i$ are normal, the left and right uniform structures on \mG coincide.  
On the other hand, the $G/G_i$ are discrete, 
therefore $\lim G/G_i$ is a uniformly complete Hausdorff topological group and the canonical arrow 
$\widehat{G} \ra \lim G/G_i$ is an isomorphism of topological groups.\\
\endgroup %%------------------------------------<<

Let \mG be a group $-$then by a 
\un{filtration}
\index{filtration (on a group)} 
on \mG we understand a sequence 
$\{G_n\}$ of normal subgroups of \mG such that $\forall \ n$, $G_n \supset G_{n+1}$.  
The filtration is said to be 
\un{exhaustive}
\index{exhaustive (filtration)} 
provided that $\bigcup\limits_n G_n = G$.  If \mK is a subgroup of \mG, 
$\{K \cap G_n\}$ is a filtration on \mK (the 
\un{induced}
\index{induced (filtration)} 
filtration) and if \mK is a normal subgroup of \mG, 
$\{K \cdot G_n/K\}$ is a filtration on $G/K$ (the 
\un{quotient}
\index{quotient (filtration)} 
filtration).

[Note: \ The $n$ run over $\Z$ but in practice it often happens that $G_0 = G$.]

Let \mG be a group with a filtration, i.e., a filtered group.  
Endow \mG with the structure of a topological group in which the $G_n$ 
become a fundamental system of neighborhoods of $e$ $-$then the canonical arrow 
$\widehat{G} \ra \lim G/G_n$ is an isomorphism of topological groups (cf. supra).  
More is true: $\widehat{G}_n$ can be identified with the closure of the image of $G_n$ in $\widehat{G}$ and the 
$\widehat{G}_n$ form a fundamental system of neighborhoods of $e$ in $\widehat{G}$, hence are normal 
%%----------------------------------------------------------------------------------------------02
open subgroups of $\widehat{G}$.  
The topology on $\widehat{G}$ is defined by the filtration $\{\widehat{G}_n\}$.  
In addition: 
$G/G_n \approx$ $\widehat{G}/\widehat{G}_n$ $\implies$ 
$\lim G /G_n \approx$ 
$\lim \widehat{G}/\widehat{G}_n$ $\implies$
$\widehat{G} \approx$ 
$\widehat{\widehat{G}\hspace{0.05cm}}$.

[Note: \ If \mK is a subgroup of \mG, the induced topology on \mK is the topology defined by the induced filtration and if \mK is a normal subgroup of \mG, the quotient topology on $G/K$ is the topology defined by the quotient filtration.]\\

\begingroup%%----------------------------------->>
\fontsize{9pt}{11pt}\selectfont
\textbf{\small EXAMPLE} \quad 
Let \mG be a filtered abelian group $-$then $\forall \ n$, there is a short exact sequence 
$0 \ra$ 
$G_n \ra$ 
$G \ra$ 
$G/G_n \ra 0$.  Since $\lim^1 G = 0$, it follows that there is an exact sequence 
$0 \ra$ 
$\lim G_n \ra$ 
$G \ra$ 
$\lim G/G_n \ra$ 
$\lim^1 G_n \ra 0$, hence $\lim^1 G_n \approx \widehat{G}/G$ provided that $\ds\bigcap\limits_n G_n = 0$.\\
\endgroup %%------------------------------------<<

\label{9.83}
\label{17.33} %dmc mnft
\indent\indent ($p$-Adic Completions) \quad Fix a prime $p$.  
Given a group \mG, let $G^{p^n}$ $(n \geq 0)$ be the subgroup of \mG generated by the $g^{p^n}$ $(g \in G)$ 
(take $G^{p^n} = G$ for $n < 0$) and set 
$G^{p^\omega}$ $\hsx = \hsx$ $\bigcap \limits_1^\infty G^{p^n}$ $-$then the $G^{p^n}$ filter \mG, thus one can form 
$\widehat{G}_p = \lim G/G^{p^n}$, the 
\un{p-adic completion}
\index{p-adic completion (of a group)} 
of \mG.  
The assignment 
$G \ra \widehat{G}_p$ defines a functor $\bGR \ra \bGR$ and this data generates a triple in \bGR.  
In general, 
$\widehat{G}_p \not\approx$ 
$(\widehat{G}_p)\overset{\raisebox{-0.05cm}{$\widehat{ \ }$}}{_p}$
but if \mG is nilpotent, then $\widehat{G}_p$ is nilpotent with 
$\nil \widehat{G}_p = \nil G/G^{p^\omega}$ and 
$\widehat{G}_p \approx$ 
$(\widehat{G}_p)\overset{\raisebox{-0.05cm}{$\widehat{ \ }$}}{_p}$
(the kernel of the projection $\widehat{G}_p \ra G/G^{p^n}$ is $(\widehat{G}_p)^{p^n}$) 
(Warfield\footnote[2]{\textit{SLN} \textbf{513} (1976), 59-60.}).  
Accordingly, $p$-adic completion restricts to a functor $\bNIL \ra \bNIL$ and 
$\bNIL_p^{\widehat{}}$, 
the full subcategory of \bNIL whose objects are Hausdorff and complete in the $p$-adic topology, is a reflective subcategory of \bNIL.  
Every object in 
$\bNIL_p^{\widehat{}}$ 
is $p$-cotorsion.

[Note: \ On a subgroup of \mG, the induced $p$-adic topology need not agree with the intrinsic $p$-adic topology.  
Moreover, the image of \mG in $\widehat{G}_p$ need not be normal and 
$(\widehat{G}_p)$\raisebox{0.1cm}{$\widehat{ \ }$} \ 
is conceptually distinct from 
$(\widehat{G}_p)\overset{\raisebox{-0.05cm}{$\widehat{ \ }$}}{_p}$
.]\\

Example: \ Take $G = \Z$ $-$then $\widehat{G}_p = \lim \Z/p^n\Z$ is $\widehat{Z}_p$, the (ring of) 
\un{$p$-adic integers}.
\index{p-adic integers}

[Note: \ $\widehat{\Z}_p$ is homeomorphic to the Cantor set, hence is uncountable.  
A 
\un{$p$-adic module}
\index{p-adic module} 
is a $\widehat{\Z}_p$-module.  
Example: Let \mG be an abelian group $-$then \mG is a $p$-adic module if \mG is $p$-primary or $p$-cotorsion.]\\

\label{9.80}
\label{11.27}
\label{11.31}
\begingroup%%----------------------------------->>
\fontsize{9pt}{11pt}\selectfont
\index{nilpotent groups (example)}
\textbf{\small EXAMPLE \  (\un{Nilpotent Groups})} \  
Suppose that \mG is nilpotent $-$then there is a short exact sequence 
$1 \ra$ 
$\Ext(\Z/p^\infty \Z ,G)^{p^\omega} \ra$ 
$\Ext(\Z/p^\infty \Z ,G) \ra$ 
$\widehat{G}_p \ra 1$, hence 
$\Ext(\Z/p^\infty \Z ,G)_p^{\widehat{}}$ 
$\approx \widehat{G}_p$.  Here, 
$\Ext(\Z/p^\infty \Z ,G)^{p^\omega}$ is a $p$-cotorsion abelian group.  
It is trivial if $G_\tor(p)$ has finite exponent, in particular, if \mG is finitely generated or torsion free.  
When \mG is abelian, $\Ext(\Z/p^\infty \Z ,G)^{p^\omega}$ can be alternatively described as 
$\Pur \Ext(\Z/p^\infty \Z ,G)$ (the subgroup of $\Ext(\Z/p^\infty \Z ,G)$ which classifies the pure extensions of \mG by 
$\Z/p^\infty \Z$) or as $\lim^1 \Hom(\Z/p^n \Z,G)$.
\\ \indent
%%----------------------------------------------------------------------------------------------03
[Note: \ Proofs of the above assertions can be found in 
Huber-Warfield\footnote[2]{\textit{J. Algebra} \textbf{74} (1982), 402-442.}.  They also show that if
$1 \ra$
$G^\prime \ra$
$G \ra$
$G\pp \ra 1$ is a short exact sequence of nilpotent groups and if $G_{\tor}\pp(p)$ has finite exponent, then the sequence 
$1 \ra$
$\widehat{G}_p^\prime \ra$
$\widehat{G}_p \ra$
$\widehat{G}_p\pp \ra 1$ is short exact.]\\
\endgroup %%------------------------------------<<

\label{9.81}
\label{9.117}
\begingroup%%----------------------------------->>
\fontsize{9pt}{11pt}\selectfont
\index{p-Adic Integers (example)}
\textbf{\small EXAMPLE \  (\un{$p$-Adic Integers})} \  
$\widehat{\Z}_p$ is a principal ideal domain.  It is the closure of $\Z$ in $\widehat{\Q}_p$ and 
$\Q \otimes \widehat{\Z}_p \approx$ 
$\widehat{\Q}_p$.  $\widehat{\Z}_p$ is a local ring with unique maximal ideal $p\widehat{\Z}_p$ and 
$\widehat{\Z}_p/p\widehat{\Z}_p \approx \F_p$.  
Examples: 
(1) $\Hom(\widehat{\Z}_p,\widehat{\Z}_p) \approx \widehat{\Z}_p$; \ 
(2) $\Hom(\Z/p^\infty \Z,\Z/p^\infty \Z) \approx \widehat{\Z}_p$; \ 
(3) $\widehat{\Q}_p/\widehat{\Z}_p \approx \Z/p^\infty \Z$; \ 
(4) $\widehat{\Z}_p \otimes \widehat{\Z}_p \approx \widehat{\Z}_p \oplus 2^\omega \cdot \Q$; \ 
(5) 
$\widehat{\Z}_p^\omega \approx (2^\omega \cdot$ 
$\widehat{\Z}_p)\overset{\raisebox{-0.05cm}{$\widehat{ \ }$}}{_p}$; \ 
(6) $\Z^\omega /\omega \cdot \Z \approx 2^\omega \cdot \Q \oplus \ds\prod\limits_p \widehat{\Z}_p^\omega$; \ 
(7) $\Ext(\Z/p^\infty \Z,\omega \cdot \Z) \approx$ 
$(\omega \cdot 
\Z)\overset{\raisebox{-0.05cm}{$\widehat{ \ }$}}{_p}$; \  
(8) $\Ext(\widehat{\Z}_p,\Z) \approx$ $\Z/p^\infty \Z \oplus \Q^{2^\omega}$.\\
\endgroup %%------------------------------------<<

\begingroup%%----------------------------------->>
\fontsize{9pt}{11pt}\selectfont
\textbf{\small EXAMPLE} \quad 
The \cd 
\begin{tikzcd}[sep=large]
{\Z_p} \ar{d} \ar{r} &{\widehat{\Z}_p} \ar{d}\\
{\Q} \ar{r} &{\widehat{\Q}_p}
\end{tikzcd}
is simultaneously a pullback and a pushout in \bAB.\\
\endgroup %%------------------------------------<<

\begingroup%%----------------------------------->>
\fontsize{9pt}{11pt}\selectfont
\textbf{\small FACT} \quad 
The $p$-adic completion functor on \bAB is not right exact.  Its $0^\text{th}$ left derived functor is $\Ext(\Z/p^\infty\Z,-)$ and its $1^\text{st}$ left derived functor is $\Hom(\Z/p^\infty\Z,-)$.\\
\endgroup %%------------------------------------<<

\indent\indent ($\F_p$-Completions) \quad Fix a prime $p$.  Given a group \mG, let 
$G = \Gamma_p^0(G) \supset \Gamma_p^1(G)$ $\supset$ $\cdots$ be its 
\un{descending $p$-central series},
\index{descending $p$-central series} 
so $\Gamma_p^{i+1}(G)$ is the subgroup of \mG generated by $[G,\Gamma_p^i(G)]$ and the $g^p$ $(g \in \Gamma_p^i(G))$.  
Note that $\Gamma_p^{i}(G)/\Gamma_p^{i+1}(G)$ is central in $G/\Gamma_p^{i+1}(G)$ and 
$\Gamma_p^{i}(G)/\Gamma_p^{i+1}(G)$ is an $\F_p$-module.  Moreover, 
$H_1(G;\F_p) \approx$ 
$\F_p \otimes (G/[G,G]) \hsx \approx \hsx$ 
$G/\Gamma_p^1(G)$.  
Definition: $\F_p G = \lim G/\Gamma_p^{i}(G)$ is the 
\un{$\F_p$-completion}
\index{F$_p$-completion} 
of \mG.  
The assignment $G \ra \F_p G$ defines a functor $\bGR \ra \bGR$ and this data generates a triple in \bGR.  In general, 
$\F_p G \not\approx \F_p \F_p G$ but 
Bousfield\footnote[3]{\textit{Memoirs Amer. Math. Soc.} \textbf{186} (1977), 1-68.} 
has shown that if $H_1(G;\F_p)$ is a finitely generated $\F_p$-module, then 
$\F_p G \approx \F_p \F_p G$.  Therefore $\F_p$-completion is idempotent on the class of finitely generated groups or the class of perfect groups.\\

\begingroup%%----------------------------------->>
\fontsize{9pt}{11pt}\selectfont
\textbf{\small LEMMA} \quad 
A group \mG has a finite central series whose factors are elementary abelian $p$-groups iff $\exists \ i$: 
$\Gamma_p^i(G) = \{1\}$ or still, iff \mG is nilpotent and $\exists \ n$: $G^{p^n} = \{1\}$.\\
\endgroup %%------------------------------------<<

\begingroup%%----------------------------------->>
\fontsize{9pt}{11pt}\selectfont
\index{nilpotent groups (example)}
\textbf{\small EXAMPLE (\un{Nilpotent Groups})} \quad 
For any group \mG, $G^{p^i} \subset \Gamma_p^i(G)$ $\forall \ i$, thus there is an arrow 
$\widehat{G}_p \ra \F_p G$.  If in addition \mG is nilpotent, then $\forall \ n$, $G/G^{p^n}$ is nilpotent and 
$(G/G^{p^n})^{p^n} = \{1\}$, hence by the lemma $\exists \ i$: $\Gamma_p^i(G/G^{p^n}) = \{1\}$ $\implies$ 
$\Gamma_p^i(G) \subset G^{p^n}$ $\implies$ $\widehat{G}_p \approx \F_p G$.  
Corollary: \mG nilpotent $\implies$ $\F_p G \approx \F_p \F_p G$.
\vspi
%%----------------------------------------------------------------------------------------------04
Recall that if $1 \ra G^\prime \ra G \ra G\pp \ra 1$ is a central extension of groups with $G^\prime$ an 
$\F_p$-module and $G\pp$ $H\F_p$-local, then \mG is $H\F_p$-local 
(cf. p. \pageref{10.1}).  
Consequently, given any \mG, it follows by induction that $\forall \ i$, $G/\Gamma_p^i(G)$ is $H\F_p$-local which means that $\F_p G$ is $H\F_p$-local as well (for, being reflective in \bGR, $\bGR_{H\F_p}$ is limit closed).  
Accordingly, there is a commutative triangle
\begin{tikzcd}[sep=large]
{G} \ar{d} \ar{r} &{\F_p G}\\
{G_{H\F_p}} \ar{ru}
\end{tikzcd}
and the arrow $G_{H\F_p} \ra \F_p G$ is an isomorphism iff $G \ra \F_p G$ is an $H\F_p$-homomorphism.  
Example: Suppose that \mG is a nilpotent group for which $G_\tor(p)$ has finite exponent $-$then 
$G_{H\F_p} \approx \F_p G$.  Proof: 
$G_{H\F_p} \approx$ 
$\Ext(\Z/p^\infty \Z,G) \approx$ 
$\widehat{G}_p \approx$ 
$\F_p G$.\\
\endgroup %%------------------------------------<<

\begingroup%%----------------------------------->>
\fontsize{9pt}{11pt}\selectfont
\textbf{\small EXAMPLE} \quad 
Take $G = \ds\bigoplus\limits_1^\infty \Z/p^n\Z$ $-$then the arrow $G_{H\F_p} \ra \F_p G$ is not an isomorphism.

[Show that the induced map $H_2(G;\F_p) \ra H_2(\F_p G; \F_p)$ is not surjective, hence that $G \ra \F_p G$ is not an 
$H\F_p$-homomorphism.]\\
\endgroup %%------------------------------------<<

\begingroup%%----------------------------------->>
\fontsize{9pt}{11pt}\selectfont
\textbf{\small FACT} \quad 
Let $f:G \ra K$ be an $H\F_p$-homomorphism $-$then $\forall \ i \geq 0$, the induced map 
$G/\Gamma_p^i(G) \ra K/\Gamma_p^i(K)$ is an isomorphism.
\vspi
[Note: \ Compare this result with Proposition 18 in $\S 8$.]\\
\endgroup %%------------------------------------<<

\begingroup%%----------------------------------->>
\fontsize{9pt}{11pt}\selectfont
Fix a set of primes \mP.  Given a group \mG, its 
\un{$P$-completion}
\index{P-completion (of a group)} 
$PG$ is 
$\lim (G/\Gamma^i(G))_P$.  The assignment $G \ra PG$ defines a functor $\bGR \ra \bGR$ and this data generates a triple in 
\bGR.  In general, $PG \not\approx PPG$ but  
Bousfield\footnote[2]{\textit{Memoirs Amer. Math. Soc.} \textbf{186} (1977), 1-68.} 
has shown that if $H_1(G;\Z_p)$ is a finitely generated $\Z_p$-module, then 
$P G \approx P P G$.  Therefore $P$-completion is idempotent on the class of finitely generated groups or the class of perfect groups.
\\ \indent
[Note: \ It is clear that $PG \approx PPG$ if \mG is nilpotent.]
\\ \indent
$P$-completion is related to $HP$-localization in the same way that $\F_p$-completion is related to $H\F_p$-localization.  In fact, since $G/\Gamma^i(G)$ is nilpotent, 
$(G/\Gamma^i(G))_P \approx$ 
$(G/\Gamma^i(G))_{HP}$ 
(cf. p. \pageref{10.2}) $\implies$ $PG$ is 
$HP$-local.  
Thus there is a commutative triangle 
\begin{tikzcd}[sep=large]
{G} \ar{d} \ar{r} &{P G}\\
{G_{HP}} \ar{ru}
\end{tikzcd}
and the arrow $G_{HP} \ra PG$ is an isomorphism iff $G \ra PG$ is an $HP$-homomorphism.\\
\endgroup %%------------------------------------<<

\begingroup%%----------------------------------->>
\fontsize{9pt}{11pt}\selectfont
\textbf{\small EXAMPLE} \quad 
Let $\pi$ be the fundamental group of the Klein bottle $-$then the arrow $\pi_{HP} \ra P_\pi$ is not an isomorphism if 
$2 \in P$.
\vspi
[By definition, $\pi \ra \pi_{HP}$ is an $HP$-homomorphism, so $H_2(\pi_{HP}; \Q) = 0$.  
On the other hand, there is a short exact sequence 
$1 \ra \Z_P \oplus \widehat{\Z}_2 \ra P_\pi \ra \Z/2\Z \ra 1$ and, from the LHS spectral sequence, 
$H_2(P_\pi;\Q) \approx$ 
$H_2(\widehat{\Z}_2;\Q) \approx$ 
$\ds\bigwedge$$_{_\Q}^{^2} (\widehat{\Z}_2 \otimes \Q)$, which is uncountable.]\\
\endgroup %%------------------------------------<<

%%----------------------------------------------------------------------------------------------05
\begingroup%%----------------------------------->>
\fontsize{9pt}{11pt}\selectfont
Notation: Given a category \bC, $\bTRI_{\bC}$ is the metacategory whose objects are the triples in \bC and 
$\bIDTRI_{\bC}$ is the full submetacategory of $\bTRI_{\bC}$ whose objects are the idempotent triples in \bC.
\vspi
[Note: \ Recall that a morphism of triples is a morphism in the metacategory $\bMON_{[\bC,\bC]}$ 
(cf. p. \pageref{10.3}).]\\
\endgroup %%------------------------------------<<

\index{Theorem: Theorem of Fakir}
\index{Theorem of Fakir}
\begingroup%%----------------------------------->>
\fontsize{9pt}{11pt}\selectfont
\textbf{\small THEOREM OF FAKIR\footnote[2]{\textit{C. R. Acad. Sci. Paris} \textbf{270} (1970), 99-101.}} \quad 
Let \bC be a category.  Assume: \bC is complete and wellpowered $-$then $\bIDTRI_{\bC}$ is a monocoreflective submetacategory of $\bTRI_{\bC}$.
\vspi
[Note: \ The coreflector sends $\bT = (T,m,\epsilon)$ to its 
\un{idempotent modification}
\index{idempotent modification} 
$\bT^\infty = (T^\infty ,m^\infty ,\epsilon^\infty)$.  
In addition: 
(1) $\forall \ \bT$, \bT and $\bT^\infty$ have the same equivalences, i.e., a morphism is rendered invertible by \mT iff it is 
rendered invertible by $T^\infty$; 
(2) $\forall \ \bT$, $\epsilon^\infty T:T \ra T^\infty \circ T$ is a natural isomorphism.]\\
\endgroup %%------------------------------------<<

\begingroup%%----------------------------------->>
\fontsize{9pt}{11pt}\selectfont
Let us take $\bC = \bGR$ and apply this result to the triple determined by $P$-completion.  
Thus, in obvious notation, 
$P^\infty G$ is the idempotent modification of $PG$, so $P^\infty G$ embeds in $PG$ while 
$PG \approx P P^\infty G$ (by 
(1)) $\&$ $PG \approx P^\infty P G$ (by (2)).  
Of course, those \mG for which the arrow $G \ra P^\infty G$ is an isomorphism constitute the object class of a reflective subcategory of \bGR.   
Moreover, $P^\infty G$ is $HP$-local, hence there is a commutative diagram 
$
\begin{tikzcd}[sep=large]
&{G}  \ar{ld} \ar{d} \ar{rd}\\
{G_P} \ar{r} &{G_{HP}} \ar{r} &{P^\infty G}
\end{tikzcd}
.  
$
When restricted to \bNIL, $L_P$, $L_{HP}$, and $P^\infty$ are naturally isomorphic but on \bGR, these functors are distinct (see below).\\
\endgroup %%------------------------------------<<

\begingroup%%----------------------------------->>
\fontsize{9pt}{11pt}\selectfont
\textbf{\small FACT} \quad 
The arrow $PG \ra PPG$ is surjective iff the induced map $H_1(G;\Z_P) \ra H_1(PG:\Z_P)$ is surjective.
\vspi
Claim: $\forall \ G$, $PG$ embeds in $PPG$.
\vspi
[For $P^\infty G$ embeds in $PG$ $\implies$ $P^\infty G$ embeds in $PPG$, i.e., $PG$ embeds in $PPG$.]
\vspi
Therefore $PG \approx PPG$ iff the induced map $H_1(G;\Z_P) \ra H_1(PG:\Z_P)$ is surjective.  This can be rephrased: 
$PG \approx PPG$ iff the arrow $G_{HP} \ra PG$ is surjective.  
Proof: Since $G_{HP}$ and $PG$ are $HP$-local, the arrow $G_{HP} \ra PG$ is surjective iff the induced map 
$H_1(G_{HP};\Z_P) \ra H_1(PG:\Z_P)$ is surjective 
(cf. p. \pageref{10.4}).\\
\endgroup %%------------------------------------<<

\begingroup%%----------------------------------->>
\fontsize{9pt}{11pt}\selectfont
\textbf{\small EXAMPLE} \quad 
Let $\pi$ be the fundamental group of the Klein bottle $-$then $\pi$ is finitely generated, 
hence $P \pi \approx P P \pi$ and the arrow $\pi_{HP} \ra P \pi$ is surjective but, as seen above, it is not an isomorphism if $2 \in P$.\\
\endgroup %%------------------------------------<<

\begingroup%%----------------------------------->>
\fontsize{9pt}{11pt}\selectfont
\textbf{\small FACT} \quad 
Let $f:G \ra K$ be a homomorphism of groups $-$then the following conditions are equivalent: 
(1) $P^\infty f: P^\infty G \ra P^\infty K$ is an isomorphism; 
(2) $Pf:PG \ra PK$ is an isomorphism; 
(3) $f \perp P X$ for every group \mX; 
(4) $f_*:(G/\Gamma^i(G))_P \ra (K/\Gamma^i(K))_P$ is an isomorphism $\forall \ i$.\\
\endgroup %%------------------------------------<<

%%----------------------------------------------------------------------------------------------06
\begingroup%%----------------------------------->>
\fontsize{9pt}{11pt}\selectfont
Application: $\forall \ G$, $H_1(G;\Z_P) \approx H_1(P^\infty G;\Z_P)$.\\
\endgroup %%------------------------------------<<

\begingroup%%----------------------------------->>
\fontsize{9pt}{11pt}\selectfont
Thus, as a consequence, $\forall \ G$, the induced map 
$H_1(G_{HP};\Z_P) \ra H_1(P^\infty G;\Z_P)$ is an isomorphism which means that the arrow 
$G_{HP} \ra P^\infty G$ is surjective 
(cf. p. \pageref{10.5}).  
Corollary: The range of the arrow $G_{HP} \ra PG$ is $P^\infty G$.
\\ \indent
[Note: \ Accordingly, 
$P^\infty G \approx PG \Leftrightarrow PG \approx PPG \Leftrightarrow H_1(G;\Z_P) \approx H_1(PG;\Z_P)$.]\\
\endgroup %%------------------------------------<<

\begingroup%%----------------------------------->>
\fontsize{9pt}{11pt}\selectfont
\textbf{\small EXAMPLE} \quad 
Let $\pi$ be the fundamental group of the Klein bottle $-$then for any \mP, $\pi_P$ is countable 
(cf. p. \pageref{10.6}).  
If now $2 \in P$, then $P^\infty \pi \approx P \pi$ is uncountable, so $\pi_P \not\approx \pi_{HP}$.  On the other hand, 
$\pi_{HP} \not\approx P^\infty \pi$.\\
\endgroup %%------------------------------------<<

\label{8.7}
\begingroup%%----------------------------------->>
\fontsize{9pt}{11pt}\selectfont
\textbf{\small FACT} \quad 
Suppose that \mG is a free group $-$then the arrow of localization $l_P:G \ra G_P$ is one-to-one.
\\ \indent
[Since \mG is free, the quotients $G/\Gamma^i(G)$ are torsion free nilpotent groups and the intersection 
$\ds\bigcap\limits_i \Gamma^i(G)$ is trivial.]\\
\endgroup %%------------------------------------<<

\indent\indent ($I$-Adic Completions) \quad Let \mA be a ring with unit, $I \subset A$ a two sided ideal.  Put 
$A_n = I^n$ $(n \geq 0)$, $A_n = A$ $(n < 0)$ $-$then $\{A_n\}$ is an exhaustive filtration on \mA, the associated topology being the 
\un{$I$-adic topology}.  
\index{I-adic topology}  
\mA is a topological ring in the $I$-adic topology.  Moreover, 
$\widehat{A}$ is a topological ring but in general, $(\widehat{I})^n \neq \widehat{I}^n$ and the $\widehat{I}$-adic topology on $\widehat{A}$ need not agree with the filtration topology.

[Note: \ Given a left $A$-module $M$, put $M_n = I^n \cdot M$ $(n \geq 0)$, $M_n = M$ $(n < 0)$ $-$then $\{M_n\}$ is an exhaustive filtration on \mM, the associated topology being the \un{$I$-adic topology}.  \mM is a topological left $A$-module in the $I$-adic topology.  Moreover, $\widehat{M}$ is a topological left $\widehat{A}$-module and 
$\widehat{M}_n = \widehat{I}^n \cdot \widehat{M} = \widehat{I}^n \cdot \im M$ $\forall \ n$ provided that \mM is finitely generated (in which case $\widehat{M}$ is finitely generated).  
Example: Take \mA commutative and \mI finitely generated: 
$\widehat{I}^n = I^n \cdot \widehat{A}$ $\implies$ 
$\widehat{I} = I \cdot \widehat{A}$ $\implies$ 
$(\widehat{I})^n = I^n \cdot \widehat{A} = \widehat{I}^n$, so, in this situation, the $\widehat{I}$-adic topology on 
$\widehat{A}$ agrees with the filtration topology.]\\

\begingroup%%----------------------------------->>
\fontsize{9pt}{11pt}\selectfont
Let \mA be a left Noetherian ring with unit, $I \subset A$ a two sided ideal $-$then \mI is said to have the 
\un{left Artin-Rees  property}
\index{left Artin-Rees  property} 
if for every finitely generated left $A$-module \mM and every left submodule $N \subset M$, the $I$-adic topology on $N$ is the restriction of the $I$-adic topology to \mM.  
Example: $I$ has the left Artin-Rees property if $\forall \ M$, $N$, $\exists \ i$: $I^i \cdot M \cap N \subset I \cdot N$.
\\ \indent
[Note: \ The theory has been surveyed by 
Smith\footnote[2]{\textit{SLN} \textbf{924} (1982), 197-240.}.]\\
\endgroup %%------------------------------------<<
%%----------------------------------------------------------------------------------------------07

\begingroup%%----------------------------------->>
\fontsize{9pt}{11pt}\selectfont
\textbf{\small EXAMPLE} \quad 
Fix a group \mG.  Definition: \mG is said to have the 
\un{Artin-Rees  property}
\index{Artin-Rees  property} 
if 
$\Z[G]$ is noetherian and $I[G]$ has the Artin-Rees  property.  Here, it is not necessary to distinguish between 
``left'' and ``right''.  
Example: Every finitely generated nilpotent group \mG has the Artin-Rees  property.\\
\endgroup %%------------------------------------<<

\begingroup%%----------------------------------->>
\fontsize{9pt}{11pt}\selectfont
Let \mA be a ring with unit, $I \subset A$ a two sided ideal $-$then there is a homomorphism of rings 
$A \ra \widehat{A}$, hence $\widehat{A}$ can be viewed as an $A$-bimodule.  Given a left $A$-module \mM, its 
\un{formal completion}
\index{formal completion (of a left module)} 
is the left $\widehat{A}$-module obtained from \mM by extension of scalars, i.e., the tensor product $\widehat{A} \otimes_A M$.
\\ \indent
[Note: \ A homomorphism \ $f:M \lra \ N$ of left \ $A$-modules leads to a commutative diagram
\begin{tikzcd}%[sep=small]
{\widehat{A} \otimes_A M} \ar{d} \ar{r} &{\widehat{A} \otimes_A N} \ar{d}\\
{\widehat{M}} \ar{r}[swap]{\widehat{f}} &{\widehat{N}}
\end{tikzcd}
of left $\widehat{A}$-modules.]
\\ \indent
Assume again that \mA is left noetherian and \mI has the left Artin-Rees property $-$then, like in the commutative case, the functor $M \ra \widehat{M}$ is exact on the category of finitely generated left $A$-modules and for all such \mM, the arrow 
$\widehat{A} \otimes_A M \ra \widehat{M}$ is bijective.  Moreover, $\widehat{A}$, as a right $A$-module, is flat.\\
\endgroup %%------------------------------------<<

\begingroup%%----------------------------------->>
\fontsize{9pt}{11pt}\selectfont
\textbf{\small FACT} \quad 
Suppose that \mA is left and right noetherian and \mI has the left and right Artin-Rees property.  
Let \mM be a left $A$-module $-$then 
$\Tor_*^A(A/I,M) \approx \Tor_*^A(A/I,\widehat{A} \otimes_A M)$.\\
\endgroup %%------------------------------------<<

\begingroup%%----------------------------------->>
\fontsize{9pt}{11pt}\selectfont
\textbf{\small EXAMPLE} \quad 
Fix a group \mG with the Artin-Rees property.  Let \mM be a finitely generated $G$-module $-$then 
$H_*(G;M) \approx H_*(G;\widehat{M})$.  Consequently, a homomorphism $f:M \ra N$ of finitely generated $G$-modules is an 
$H\Z$-homomorphism iff $\widehat{f}:\widehat{M} \ra \widehat{N}$ is an isomorphism.\\
\endgroup %%------------------------------------<<

\begingroup%%----------------------------------->>
\fontsize{9pt}{11pt}\selectfont
\textbf{\small FACT} \quad 
Suppose that \mG is a finitely generated nilpotent group.  Let \mM be a finitely generated $G$-module $-$then 
$\widehat{M}$ is $H\Z$-local and the arrow of completion $M \ra \widehat{M}$ is an $H\Z$-homomorphism, thus 
$M_{H\Z} \approx \widehat{M}$.\\
\endgroup %%------------------------------------<<

\begingroup%%----------------------------------->>
\fontsize{9pt}{11pt}\selectfont
\textbf{\small EXAMPLE} \quad 
Take $G = \Z/2\Z$ and for any abelian group \mM, let \mG operate on \mM by ``negation''.  In this situation, 
$M_{H\Z} \approx \Ext(\Z/2^\infty \Z,M)$ and there is a short exact sequence 
$0 \ra$ 
$\lim^1 \Hom(\Z/2^n\Z,M) \ra$ 
$\Ext(\Z/2^\infty\Z,M) \ra$ 
$\widehat{M} \ra 0$ 
(cf. p. \pageref{10.7}).  And: The epimorphism 
$\Ext(\Z/2^\infty\Z,M) \ra \widehat{M}$ has a nonzero kernel if 
$M = \ds\bigoplus\limits_1^\infty \Z/2^n\Z$.\\
\endgroup %%------------------------------------<<

\label{10.10}
\label{10.11}
A Hausdorff topological group \mG is said to be 
\un{profinite}
\index{profinite (Hausdorff topological group)} 
if it is compact and totally disconnected or, equivalently, that 
$G \approx \lim G_i$, where $i$ runs over a directed set and $\forall \ i$, $G_i$ is a finite group (discrete topology).

[Note: \ If \mG is profinite, then $G \approx \lim G/U$, $U$ open and normal.]\\

\begingroup%%----------------------------------->>
\fontsize{9pt}{11pt}\selectfont
\textbf{\small EXAMPLE} \quad 
Let \mG be a Hausdorff topological group.  Assume: \mG is compact and torsion $-$then \mG is profinite.\\
\endgroup %%------------------------------------<<
%%----------------------------------------------------------------------------------------------08

\begingroup%%----------------------------------->>
\fontsize{9pt}{11pt}\selectfont
\textbf{\small EXAMPLE} \quad 
Let \mG be an abelian group $-$then \mG is algebraically isomorphic to a profinite abelian group iff \mG is algebraically isomorphic to a product 
$\ds\prod\limits_p \bigl[\widehat{\Z}_p^{\kappa_p} \times \ds\prod\limits_{i \in I_p} \Z/p^{n_i}\Z\bigr]$.  
Here, $\kappa_p$ is a cardinal number (possibly zero), $I_p$ is an index set (possibly empty), and $n_i$ is a positive integer.\\
\endgroup %%------------------------------------<<

\begingroup%%----------------------------------->>
\fontsize{9pt}{11pt}\selectfont
\textbf{\small EXAMPLE} \quad 
Let $k$ be a field, \mK a Galois extension of $k$.  Put $G = \Gal(K/k)$ $-$then \mG is a profinite group.  In fact, 
$G \approx \lim G_i$, where $G_i = \Gal(K_i/k)$, $K_i$ a finite Galois extension of $k$.
\\ \indent
[Note: The quotient $G/\ov{[G,G]}$ can be identified with $\Gal(k^\text{ab}/k)$, $k^\text{ab}$ the maximal abelian extension of $k$ in \mK.]\\
\endgroup %%------------------------------------<<

Given a group \mG, the 
\hsx
\un{profinite completion} 
\index{profinite completion (of a group)} 
\hsx 
$\pro G$ of \mG is 
$\lim G/U$, the limit being taken over the normal subgroups of finite index in \mG.  The assignment 
$G \ra \pro G$ defines a functor $\bGR \ra \bGR$ and this data generates a triple in \bGR which, however, is not idempotent.

Example: Take $G = \Z$ $-$then $\pro \Z = \lim \Z/n\Z$ is $\widehat{\Z}$, the (ring of) 
\un{$\Pi$-adic integers}
\index{Pi-adic integers, $\Pi$-adic integers}.\\

\begingroup%%----------------------------------->>
\fontsize{9pt}{11pt}\selectfont
\textbf{\small EXAMPLE} \quad 
Every residually finite group embeds in its profinite completion.  This said, 
Evans\footnote[2]{\textit{J. Pure Appl. Algebra} \textbf{65} (1990), 101-104.}
has shown that for each prime $p$, there exists a countable, torsion free, residually finite group \mG such that $\pro G$ contains an element of order $p$.\\
\endgroup %%------------------------------------<<

\begingroup%%----------------------------------->>
\fontsize{9pt}{11pt}\selectfont
\textbf{\small EXAMPLE} \quad 
Let $k = \F_p$ $-$then $\Gal(\bar{k}/k) \approx \widehat{\Z}$.  
Moreover, the infinite cyclic group generated by the Frobenius is dense in $\Gal(\bar{k}/k)$.\\
\endgroup %%------------------------------------<<

\begingroup%%----------------------------------->>
\fontsize{9pt}{11pt}\selectfont
\textbf{\small EXAMPLE} \quad 
It follows from the positive solution to the congruence subgroup problem for $\bSL(n,\Z)$ $(n > 2)$ that 
$\pro \bSL(n,\Z) \approx \ds\prod\limits_p \bSL(n,\widehat{\Z}_p)$.\\
\endgroup %%------------------------------------<<

\begingroup%%----------------------------------->>
\fontsize{9pt}{11pt}\selectfont
\textbf{\small EXAMPLE} \quad 
Define a homomorphism 
$\chi:\widehat{\Z} \ra \Aut \widehat{\Z}$ by 
$\chi(\widehat{n}) = id_{\widehat{\Z}}$ if 
$\widehat{n} \in 2\widehat{\Z}$ and 
$\chi(\widehat{n}) = -id_{\widehat{\Z}}$ if 
$\widehat{n} \notin 2\widehat{\Z}$.  $-$then the semidirect product 
$\widehat{\Z} \rtimes_\chi \widehat{\Z}$ is isomorphic to $\pro \pi$, $\pi$ the fundamental group of the Klein bottle.\\
\endgroup %%------------------------------------<<

\label{11.11}
\begingroup%%----------------------------------->>
\fontsize{9pt}{11pt}\selectfont
\textbf{\small EXAMPLE} \quad 
Let \mG be a finitely generated nilpotent group $-$then $\pro G$ is nilpotent and $\nil G = \nil \pro G$.  
Proof: \mG is residually finite 
(cf. p. \pageref{10.8}), hence embeds in $\pro G$.
\\ \indent
\label{11.17}
[Note: Blackburn's\footnote[3]{\textit{Proc. Amer. Math. Soc.} \textbf{16} (1965), 143-148.}
theorem says that two elements of \mG are conjugate iff their images in every finite quotient of \mG are conjugate, i.e., two elements of \mG are conjugate iff they are conjugate in $\pro G$.]\\
\endgroup %%------------------------------------<<

%%----------------------------------------------------------------------------------------------09
\label{11.10}
\label{11.21}
\begingroup%%----------------------------------->>
\fontsize{9pt}{11pt}\selectfont
\textbf{\small EXAMPLE} \quad 
If 
$1 \ra G^\prime \ra G \ra G\pp \ra 1$ is short exact, then 
$1 \ra \pro G^\prime \ra \pro G \ra \pro G\pp \ra 1$ need not be short exact even when the data is abelian 
(e.g., pro turns 
$0 \ra$ 
$\Z \ra$ 
$\Q \ra$ 
$\Q/\Z \ra 0$ into 
$0 \ra$
$\widehat{\Z} \ra$
$0 \ra$
$0 \ra 0$).  
However, there are positive results.  For instance 
Schneebeli\footnote[2]{\textit{Arch. Math.} \textbf{31} (1978), 244-253.}
has shown that pro preserves short exact sequences in the class of polycyclic groups, thus in the class of finitely generated nilpotent groups.\\
\endgroup %%------------------------------------<<

\begingroup%%----------------------------------->>
\fontsize{9pt}{11pt}\selectfont
\textbf{\small FACT} \quad 
Suppose that \mG is a finitely generated nilpotent group $-$then $\forall \ i \geq 0$, 
$\pro \Gamma^i(G) \approx \Gamma^i(\pro G)$.\\
\endgroup %%------------------------------------<<

\begingroup%%----------------------------------->>
\fontsize{9pt}{11pt}\selectfont
\textbf{\small FACT} \quad 
Suppose that \mG is a finitely generated nilpotent group $-$then every normal subgroup of $\pro G$ of finite index is open.
\\ \indent
[Note: \ This can fail if \mG is not finitely generated (consider a discontinuous homomorphism 
$(\Z/p\Z)^\omega \ra \Z/p\Z)$.]\\
\endgroup %%------------------------------------<<

\label{11.19}
\begingroup%%----------------------------------->>
\fontsize{9pt}{11pt}\selectfont
A group \mG is said to have 
\un{property S}
\index{property S} 
if for any $\pro G$-module \mM which is finite as an abelian group, $H^n(\pro G;M) \approx H^n(G;M)$ $\forall \ n$.  
Example: Every cyclic group has property S.\\
\endgroup %%------------------------------------<<

\begingroup%%----------------------------------->>
\fontsize{9pt}{11pt}\selectfont
\textbf{\small FACT} \quad 
Suppose that \mG is a finitely generated nilpotent group $-$then \mG has property S.
\endgroup %%------------------------------------<<

\begingroup%%----------------------------------->>
\fontsize{9pt}{11pt}\selectfont
[Consider first the case of a central extension 
$1 \ra K \ra G \ra G/K \ra 1$, where \mK is cyclic and assume that the assertion holds for $G/K$.  
Claim: The assertion hold for \mG.  Indeed, since \mG is a finitely generated nilpotent group, the sequence 
$1 \ra \pro K \ra \pro G \ra \pro G/K \ra 1$ is exact (cf. supra), so there is a morphism of LHS spectral sequences 
\endgroup

\begingroup
\fontsize{9pt}{11pt}\selectfont
\[
\begin{tikzcd}%[sep=small]
{H^p(\pro G/K;H^q(\pro K;M))} \ar{d} &{\implies}  &{H^{p+q}(\pro G;M)} \ar{d} \\
{H^p(G/K;H^q(K;M))} &{\implies} &{H^{p+q}(G;M)}
\end{tikzcd}
\]
which is an isomorphism on the $E_2^{p,q}$.  In general, one can find a central series 
$G = G^0 \supset \cdots \supset G^n = \{1\}$, where $\forall \ i$, $G_i$ is normal in \mG and $G^i/G^{i+1}$ is cyclic.  
Proceed from here inductively to see that the $G/G_i$ have property S.]\\
\endgroup %%------------------------------------<<

Although profinite completion is not an idempotent functor on \bGR, it is idempotent on \bTOPGR, the category of topological groups.  
Thus let \mG be a topological group $-$then its 
\un{continuous profinite completion}
\index{continuous profinite completion} 
$\proc G$ is $\lim G/U$, the limit being taken over the open normal subgroups of finite index in \mG.  
With this understanding, 
$\proc G \approx \proc\proc G$.

\label{11.20}
[Note: \ Given a group \mG, $\pro G \approx \pro \pro G$ iff every normal subgroup of $\pro G$ of finite index is open.  
Corollary: $\pro G \approx \pro\pro G$ iff every homomorphism $G \ra F$, 
%
%%----------------------------------------------------------------------------------------------10
where \mF is finite, can be extended uniquely to a homomorphism $\pro G \ra F$ (in general, 
$\Hom_c(\pro G,F) \approx \Hom(G,F)$, the subscript standing for ``continuous''.  
Example: pro is idempotent on the class of finitely generated nilpotent groups.]\\

\begingroup%%----------------------------------->>
\fontsize{9pt}{11pt}\selectfont
\textbf{\small FACT} \quad 
Let $f:G \ra K$ be a homomorphism of groups $-$then $\pro f: \pro G \ra \pro K$ is an isomorphism of topological groups iff 
$\forall$ finite group \mF, $\Hom(K,F) \approx \Hom(G,F)$.
\vspi
[Note: \ pro is not a conservative functor 
(Platonov-Tavgen\footnote[2]{\textit{K-Theory} \textbf{4} (1990), 89-101.}).]\\
\endgroup %%------------------------------------<<

Let \mG be a profinite group $-$then \mG is said to be 
\un{$p$-profinite}
\index{p-profinite (group)} 
if \mG is $p$-local.  In this connection, recall that a finite group is a $p$-group iff it is $p$-local 
(cf. p. \pageref{10.9}).  
Upon representing \mG as $\lim G_i$ 
(cf. p. \pageref{10.10}), it follows that \mG is $p$-profinite iff $\forall \ i$, $G_i$ is 
$p$-local.

[Note: \ Let \mG be a finite group $-$then \mG is $p$-local iff $\forall \ q \neq p$, the arrow $g \ra g^q$ is surjective.]\\

\label{9.84}
\label{9.85}
\begingroup%%----------------------------------->>
\fontsize{9pt}{11pt}\selectfont
\index{p-Adic Units (example)}
\textbf{\small EXAMPLE \ (\un{$p$-Adic Units})} \  
Put $\widehat{\bU}_p = \lim (\Z/p^n\Z)^\times$ $-$then $\widehat{\bU}_p$ is $p$-profinite.  
It is the group of units in 
$\widehat{\Z}_p$.  
Using the ''exp-log'' correspondence, one shows that 
$\widehat{\bU}_p \approx \Z/(p-1) \Z \oplus \widehat{\Z}_p$ if $p$ is odd, while 
$\widehat{\bU}_2 \approx \Z/2\Z \oplus \widehat{\Z}_2$.\\
\endgroup %%------------------------------------<<

\begingroup%%----------------------------------->>
\fontsize{9pt}{11pt}\selectfont
\textbf{\small EXAMPLE} \quad 
Let $\Q^{\cy}$ be the field generated over $\Q$ by the roots of unity in $\ov{\Q}$.  
For each prime $p$, choose $\omega_n$ subject to $\omega_n^{p^n} = 1$ $\&$ $\omega_{n+1}^{p^n} = \omega_n$ $(n \geq 1)$.  Let $K_p$ be the field generated over $\Q$ by the roots of unity in $\ov{\Q}$ whose order is a power of $p$ $-$then 
$K_p = \ds\bigcup\limits_n \Q(\omega_n)$ $\implies$
$\Gal(K_p/\Q) \approx$ $\lim \Gal(\Q(\omega_n)/\Q)$.  But
$\Gal(\Q(\omega_n)/\Q) \approx$ $(\Z/p^n\Z)^\times$ $\implies$
$\Gal(K_p/\Q) \approx$ $\widehat{\bU}_p$ $\implies$ 
$\Gal(\Q^{\cy}/\Q) \approx$ 
$\ds\prod\limits_p \widehat{\bU}_p \approx$ $\widehat{\Z}^\times$.
\\ \indent
[Note: \ It follows from global class field theory that $\Q^{\cy}$ is the maximal abelian extension $\Q^{\ab}$ of $\Q$ in $\ov{\Q}$.]\\
\endgroup %%------------------------------------<<

\begingroup%%----------------------------------->>
\fontsize{9pt}{11pt}\selectfont
\textbf{\small EXAMPLE} \quad 
Suppose that \mG is $p$-profinite.  Assume: \mG is torsion $-$then 
Zelmanov\footnote[3]{\textit{Israel J. Math.} \textbf{77} (1992), 83-95.}
has shown that \mG is locally finite.\\
\endgroup %%------------------------------------<<

\begingroup%%----------------------------------->>
\fontsize{9pt}{11pt}\selectfont
Platonov had conjectured that every Hausdorff topological group which is compact and torsion is locally finite (such a group is necessarily profinite 
(cf. p. \pageref{10.11})).  
Wilson\footnote[6]{\textit{Monatsh. Math.} \textbf{96} (1983), 57-66.}
reduced this to the $p$-profinite case which was then disposed of by Zelmanov.\\
\endgroup %%------------------------------------<<

Given a group \mG, the 
\un{$p$-profinite completion}
\index{p-profinite completion (of a group)} 
$\prop G$ of \mG is $\lim G/U$, the limit being taken over the normal subgroups of finite index in \mG subject to 
$[G:U] \in \{p^n\}$.  The
%%----------------------------------------------------------------------------------------------11
assignment $G \ra \prop G$ defines a functor $\bGR \ra \bGR$ and this data generates a triple in \bGR which, however, is not idempotent.

[Note: \ Since $\prop G$ is $p$-local, there is a commutative triangle
\begin{tikzcd}%[sep=small]
{G} \ar{d} \ar{r} &{\prop G}\\
{G_p} \ar{ru}
\end{tikzcd}
and a natural transformation $L_p \ra \prop$.]

Example: Take $G = \Z$ $-$then $\prop \Z = \lim \Z/p^n \Z$ is $\widehat{\Z}_p$, the (ring of) 
\un{$p$-adic integers}.
\index{p-adic integers}\\

\begingroup%%----------------------------------->>
\fontsize{9pt}{11pt}\selectfont
\textbf{\small EXAMPLE} \quad 
Define a homomorphism 
$\chi:\widehat{\Z}_2 \ra$ $\Aut \widehat{\Z}_2$ by 
$\chi(\widehat{n}) =$ $\id_{\widehat{\Z}_2}$ if 
$\widehat{n} \hspace{0.03cm}  \in 2\widehat{\Z}_2$ and 
$\chi(\widehat{n})$ $=$ $-\id_{\widehat{\Z}_2}$ if 
$\widehat{n} \notin 2\widehat{\Z}_2$  $-$then the semidirect product 
$\widehat{\Z}_2 \rtimes_\chi \widehat{\Z}_2$ is isomorphic to $\prox_2 \pi$, $\pi$ the fundamental group of the Klein bottle.
\\ \indent
[Note: \ For $p$ odd, $\prox_p \pi \approx \widehat{\Z}_p$.  Therefore a nonabelian group can have an abelian $p$-profinite completion.]\\
\endgroup %%------------------------------------<<

\begingroup%%----------------------------------->>
\fontsize{9pt}{11pt}\selectfont
\textbf{\small LEMMA} \quad 
Suppose tht $G/\Gamma_p^1(G)$ is finite $-$then $\forall \ i > 1$, $G/\Gamma_p^i(G)$ is a finite $p$-group.\\
\endgroup %%------------------------------------<<

\begingroup%%----------------------------------->>
\fontsize{9pt}{11pt}\selectfont
Application: $\dim H_1(G;\F_p) < \omega$ $\implies$ $\pro_p G \approx \F_p G$.\\[-.1cm]
\endgroup %%------------------------------------<<

\begingroup%%----------------------------------->>
\fontsize{9pt}{11pt}\selectfont
\textbf{\small EXAMPLE} \quad 
Let \mF be a free group on $n > 1$ generators $-$then $\prop F \approx \F_p F$ and 
Bousfield\footnote[2]{\textit{Trans. Amer. Math. Soc.} \textbf{331} (1992), 335-359.}
has shown that $H_1(\prop F;\F_p) \approx n \cdot \F_p$ but for some $q > 1$, $H_q(\prop F;\F_p)$ is uncountable.
\\ \indent
[Note: \ If $F^k$ is the subgroup of \mF generated by the $k^\text{th}$ powers, 
then it follows from the negative solution to the Burnside problem that $F/F^k$ is infinite provided that $k \gg 0$ 
(Ivanov\footnote[3]{\textit{Bull. Amer. Math. Soc.} \textbf{27} (1992), 257-260.}).  
This circumstance makes it difficult to compare $\widehat{F}_p$ and $\prop F$.]\\[-.1cm]
\endgroup %%------------------------------------<<

\label{11.22}
\begingroup%%----------------------------------->>
\fontsize{9pt}{11pt}\selectfont
\textbf{\small EXAMPLE} \quad 
For any \mG there is an arrow $\widehat{G}_p \ra \prop G$.  It is an isomorphism if \mG is finitely generated and nilpotent but not in general (consider $\omega \cdot (\Z/p\Z)$).\\[-.1cm]
\endgroup %%------------------------------------<<

\label{11.23}
\begingroup%%----------------------------------->>
\fontsize{9pt}{11pt}\selectfont
\textbf{\small FACT} \quad 
Suppose that \mG is a finitely generated nilpotent group $-$then the arrow $\pro G \ra \ds\prod\limits_p \prop G$ is an isomorphism.
\\ \indent
[Note: \ This can fail if \mG is not nilpotent (Consider $S_3$).]\\[-.1cm]
\endgroup %%------------------------------------<<

\begingroup%%----------------------------------->>
\fontsize{9pt}{11pt}\selectfont
\textbf{\small FACT} \quad 
Suppose that \mG is a finitely generated nilpotent group.  Let \mK be a subgroup of \mG $-$then the $p$-profinite topology on \mK is the restriction of the $p$-profinite topology on \mG.
\endgroup %%------------------------------------<<
%%%%%%%%%%%%%%%%%%%%%%%%%%%%%%%%%%%%%%
%%%%%%%%%%%%%%%%%%%%%%%%%%%%%%%%%%%%%%
%%%%%%%%%%%%%%%%%%%%%%%%%%%%%%%%%%%%%%

\begin{center}
$\S \ 10$
\\[0.5cm]
$\mathcal{REFERENCES}$\\
\end{center}

\[
\textbf{BOOKS}
\]

\begingroup
\fontsize{9pt}{11pt}\selectfont
\setlength\parindent{0 cm}

[1] \quad Bourbaki, N., \textit{Commutative Algebra}, Addison-Wesley (1972).
\\[-.2cm]

[2] \quad Serre, J-P., \textit{Cohomologie Galoisienne}, Springer Verlag (1994).
\\[-.2cm]

[3] \quad Shatz, S., \textit{Profinite Groups, Arithmetic, and Geometry}, Princeton University Press (1972).
\\[-.2cm]
\endgroup

\[
\textbf{ARTICLES}
\]

\begingroup
\fontsize{9pt}{11pt}\selectfont
\setlength\parindent{0 cm}

[1] \quad Bousfield, A., Homological Localization Towers for Groups and $\Pi$-Modules, 
\textit{Memoirs Amer. Math.} 

\hspace{0.8cm}\textit{Soc.} \textbf{186} (1977), 1-68.
\\[-.2cm]

[2] \quad Brown, K., and Dror, E., The Artin-Rees Property and Homology, 
\textit{Israel J. Math.} \textbf{22} (1975), 93-109.
\\[-.2cm]

[3] \quad Dwyer, W., Homological Localization of $\Pi$-Modules, 
\textit{J. Pure Appl. Algebra} \textbf{10} (1977), 135-151.

\setlength\parindent{2em}

\endgroup

\chapter{
$\boldsymbol{\S}$\textbf{11}.\quadx  HOMOTOPICAL COMPLETION}
\setlength\parindent{2em}
\setcounter{proposition}{0}

%%----------------------------------------------------------------------------------------------01

$\text{ }$\\[-1.25cm]

In homotopy theory, completion appeared on the scene before localization and, to a certain extent, has been superceded by it.  Because of this, a semiproofless account will suffice.

One approach to completing a space at a prime $p$ is due to 
Bousfield-Kan\footnote[2]{\textit{SLN} \textbf{304} (1972); 
see also Iwase, \textit{Trans. Amer. Math. Soc.} \textbf{320} (1990), 77-90.\vspace{0.15 cm}}.  
It is the analog of the $\F_p$-completion process for groups.  Thus there is a functor $X \ra \F_p X$ on 
$\bHCONCWSP_*$ called 
\un{$\F_p$-completion}
\index{F$_p$-completion ($\bHCONCWSP_*$)} 
which is part of a triple.  
It is not idempotent but $\F_p X$ is $H\F_p$-local so there is a triangle 
\begin{tikzcd}%[ sep=small]
{X} \ar{d} \ar{r} &{\F_p X}\\
{X_{H\F_p}} \ar{ru}
\end{tikzcd}
, commutative up to pointed homotopy.
Definition: \mX is said to be 
\un{$\F_p$-good}
\index{F$_p$-good} 
provided that the arrow 
$X_{H\F_p} \ra \F_p X$ is a pointed homotopy equivalence; otherwise, \mX is said to be 
\un{$\F_p$-bad}.
\index{F$_p$-bad}  
For \mX to be $\F_p$-good, it is necessary and sufficient that the arrow \
$\F_p X \ra \F_p \F_p X$ be a pointed homotopy equivalence.  Therefore $\F_p$-completion is idempotent on the class of 
$\F_p$-good spaces.

[Note: \ \mX is $\F_p$-good iff the arrow $X \ra \F_p X$ is an $H\F_p$-equivalence.]

Examples: 
(1) Let \mX be a pointed connected CW space $-$then \mX is $\F_p$-good if 
(i) \mX is nilpotent or 
(ii) $\pi_1(X)$ is finite or 
(iii) $H_1(X;\F_p)$ is trivial; 
(2) Let \mF be a free group $-$then 
$\F_p K(F,1) \approx K(\F_p F,1)$ but $K(F,1)$ is $\F_p$-bad if \mF is free on two generators, i.e., $\bS^1 \vee \bS^1$ is 
$\F_p$-bad 
(Bousfield\footnote[3]{\textit{Trans. Amer. Math. Soc.} \textbf{331} (1992), 335-359.\vspace{0.15 cm}}).\\

\begingroup%%----------------------------------->>
\fontsize{9pt}{11pt}\selectfont
As a  heuristic guide, $H\F_p$-localization can be thought of as the ``idempotent modification'' of $\F_p$-completion.  
Reason: $f:X \ra Y$ is an $H\F_p$-equivalence iff $\F_p f:\F_p X \ra \F_p Y$ is a pointed homotopy equivalence, thus 
$H\F_p$-localization and $\F_p$-completion have the same equivalences (cf. $\S 9$, Proposition 21).

[Note: \ In a sense that can be made precise, the $\F_p$-completion of a space is but an initial step along the transfinite road to its $H\F_p$-localization 
(Dror-Dwyer\footnote[6]{\textit{Comment. Math. Helv.} \textbf{52} (1977), 185-210; 
see also \textit{Israel J. Math.} \textbf{29} (1978), 141-154.}).]\\

\endgroup %%------------------------------------<<

\begingroup%%----------------------------------->>
\fontsize{9pt}{11pt}\selectfont
\index{Theorem: Fiber Theorem}
\index{Fiber Theorem}
\textbf{\footnotesize FIBER THEOREM} \quad 
Let 
$
\begin{cases}
\ X\\
\ Y
\end{cases}
$
be pointed connected CW spaces, $f:X \ra Y$ a pointed continuous function with $\pi_0(E_f) = *$.  Assume: The action of 
$\pi_1(Y)$ on the $H_n(E_f;\F_p)$ is nilpotent $\forall \ n$ $-$then $\F_p E_f$ can be identified with the mapping fiber of the arrow $\F_p X \ra \F_p Y$.
%%----------------------------------------------------------------------------------------------02

[Note: \ The action of $\pi_1(\F_p Y)$ on the $H_n(\F_p E_f;\F_p)$ is nilpotent $\forall \ n$ if $E_f$ is $\F_p$-good, thus if $E_f$ and \mY are both $\F_p$-good, then so is \mX.]\\ 

\endgroup %%------------------------------------<<

\begingroup%%----------------------------------->>
\fontsize{9pt}{11pt}\selectfont
\textbf{\small EXAMPLE} \quad 
Suppose that \mX is a pointed connected CW space with the property that $\pi_1(X)$ operates nilpotently on the 
$H_n(\widetilde{X};\F_p)$ $\forall \ n$ $-$then \mX is $\F_p$-good if in addition $\pi_1(X)$ is nilpotent.\\

\endgroup %%------------------------------------<<

\begingroup%%----------------------------------->>
\fontsize{9pt}{11pt}\selectfont
\index{Theorem: $\F_p$ Whitehead Theorem}
\index{F$_p$ Whitehead Theorem}
\textbf{\footnotesize $\F_p$ WHITEHEAD THEOREM} \quad 
Let 
$
\begin{cases}
\ X\\
\ Y
\end{cases}
$
be pointed connected CW spaces, $f:X \ra Y$ a pointed continuous function.  Assume: 
$f_*:H_q(X;\F_p) \ra H_q(Y;\F_p)$ is bijective for $1 \leq q < n$ and surjective for $q = n$ 
$-$then $\F_p f$ is an $n$-equivalence.

[Note: \ To explain the difference in formulation between the $\F_p$ Whitehead theorem and the $H\F_p$ Whitehead theorem 
(cf. p. \pageref{11.1}), one has only to recall that the arrows 
$
\begin{cases}
\ X \ra X_{H\F_p}\\
\ Y \ra Y_{H\F_p}
\end{cases}
$
are $H\F_p$-equivalences.]\\

Application: \mX $n$-connected $\implies$ $\F_p X$ $n$-connected.\\

\endgroup %%------------------------------------<<

\label{11.25}
\label{11.29}
\begingroup%%----------------------------------->>
\fontsize{9pt}{11pt}\selectfont
\textbf{\small EXAMPLE} \quad 
Define functors $L_n^p:\bGR \ra \bGR$ by writing 
$L_n^p G = \pi_{n+1}(\F_p K(G,1))$ $(n \geq 0)$.  So, e.g., for any pointed connected CW space \mX, 
$\pi_1(\F_pX) \approx L_0^p \pi_1(X)$ ($\F_p$ Whitehead theorem).  Since 
$\F_p K(G,1)$ is $H\F_p$-local, $L_n^p G$ is abelian $p$-cotorsion $(n \geq 1)$.  
Examples: 
(1) If \mG is free, then $L_0^p G \approx \F_p G$ and $L_n^p G = 0$ $(n \geq 1)$; 
(2) If \mG is nilpotent, then 
$L_0^p G \approx \Ext(\Z/p^\infty\Z,G)$,
$L_1^p G  \approx \Hom(\Z/p^\infty\Z,G)$, 
and 
$L_n^p G = 0$ $(n \geq 2)$; 
(3) If \mG is finite, then $L_n^p G$ is a finite $p$-group which is trivial when $p$ and $\#(G)$ are relatively prime.

[Note: \ $\forall \ G$, there is a surjection $L_0^p G \ra \F_p G$ 
(Bousfield\footnote[2]{\textit{Memoirs Amer. Math. Soc.} \textbf{186} (1977), 1-68 (cf. 66).\vspace{0.15 cm}})
which is a bijection whenever 
$H_1(G;\F_p)$ and 
$H_2(G;\F_p)$ are finite dimensional, e.g., if \mG is finitely presented 
(Brown\footnote[3]{\textit{Cohomology of Groups}, Springer Verlag (1982), 197-198.\vspace{0.15 cm}}).]\\

\endgroup %%------------------------------------<<

\begingroup%%----------------------------------->>
\fontsize{9pt}{11pt}\selectfont
\textbf{\small EXAMPLE} \quad 
Let \mA be a ring with unit $-$then the arrow 
$B\bGL(A) \ra B\bGL(A)^+$ is a homology equivalence, hence it is an $H\F_p$-equivalence.  Therefore 
$L_n^p\bGL(A) \approx \pi_{n+1}(\F_pB\bGL(A)^+)$, so if the $K_n(A)$ are finitely generated, 
$L_n^p\bGL(A) \approx$ $\widehat{\Z}_p \otimes K_{n+1}(A)$ 
(cf. p. \pageref{11.2}).\\

\endgroup %%------------------------------------<<

\begingroup%%----------------------------------->>
\fontsize{9pt}{11pt}\selectfont
Here is a final point.  
Fix a set of primes \mP $-$then Bousfield-Kan (ibid.) have shown that the $P$-completion process for groups can be imitated in the homotopy category, i.e., there is a functor $X \ra PX$ on $\bHCONCWSP_*$ called 
\un{$P$-completion}
\index{P-completion ($\bHCONCWSP_*$)} 
which is part of a triple.  
Its formal properties are identical to those of the $\F_p$-completion and its ``idempotent modification'' is $HP$-localization.  
Example: $\bP^2(\R)$ is $P$-bad if $2 \in P$ but $\bP^2(\R)$ is $\F_p$-good $\forall \ p$ 
(since $\pi_1(\bP^2(\R)) \approx \Z/2\Z$ is finite).\\ 

\endgroup %%------------------------------------<<

Another approach to completing a space at a prime $p$ is due to 
Sullivan\footnote[4]{\textit{Ann. of Math.} \textbf{100} (1974), 1-79.}.  
In this
%%----------------------------------------------------------------------------------------------03
context, there is also an analog of the profinite completion process for groups and we shall consider it first.

Notation: $\bF_*$ is the full subcategory of $\bCONCW_*$ whose objects are the pointed connected CW complexes with finite homotopy groups and $\bH\bF_*$ is the associated homotopy cateogory.

[Note: \ Any skeleton $\ov{\bH\bF}_*$ of $\bH\bF_*$ is small.]\\

\textbf{\small LEMMA} \quad For every pointed connected CW complex \mX, the category 
$X\backslash \ov{\bH\bF}_*$ is cofiltered.

[This is because  $\bH\bF_*$ has finite products and weak pullbacks.]

[Note: \ The objects of $X\backslash \ov{\bH\bF}_*$ are the pointed homotopy classes of maps $X \ra K$ and the morphisms 
$(X \ra K) \ra (X \ra L)$ are the pointed homotopy commutative triangles
\begin{tikzcd}[ sep=small]
&{X} \ar{ldd}\ar{rdd}\\
\\
{K} \ar{rr} &&{L}
\end{tikzcd}
.]\\
\vspace{0.5cm}

In what follows, $\lim\limits_X$ stands for a limit calculated over $X\backslash \ov{\bH\bF}_*$.\\

\begin{proposition} \ %1
For every pointed connected CW complex \mX, the cofunctor 
$F_X:\bHCONCW_* \ra \bSET$ defined by $F_X Y = \lim\limits_X [Y,K]$ is representable.
\end{proposition}

[It \  is \ a \ question \ of \ applying \ the \ Brown \ representability \ theorem.  \  
That $F_X$ satisfies the wedge condition is automatic.  \ 
Turning to the Mayer-Vietoris condition, if $Y_k$ is a pointed finite connected subcomplex of \mY, then $[Y_k,K]$ is finite 
(cf. p. \pageref{11.3}).  \ 
Give it the discrete topology and form $\lim [Y_k,K]$, a nonempty compact Hausdorff space.  Since 
$[Y,K] \approx \lim [Y_k,K]$ 
(cf. p. \pageref{11.4}), it follows that there is a factorization \qquad
$
\begin{tikzcd}%[ sep=small]
{\bHCONCW_*}\ar{rd}[swap]{F_X} \ar[dashed]{r} &{\bCPTHAUS} \ar{d}{U}\\
&{\bSET}
\end{tikzcd}
, \ 
$
where \mU is the forgetful functor.  The verification that $F_X$ satisfies the Mayer-Vietoris condition is now straightforward.]\\

\index{profinite completion (space)}
The \un{profinite completion} of \mX, denoted $\pro X$, is an object that represents $F_X$.  
There is a natural transformation $[-,X] \lra [-,\pro X]$ and an arrow $\prox_X:X \lra \pro X$ (Yoneda).

[Note: \ Profinite completion generates a triple in $\bHCONCW_*$ or ($\bHCONCWSP_*$) which, however, is not idempotent.]\\

\begingroup%%----------------------------------->>
\fontsize{9pt}{11pt}\selectfont
\textbf{\small EXAMPLE} \quad 
Let \mG be a topological group.  Assume: \mG is Lie and $\#(\pi_0(G)) < \omega$ $-$then $B_G^\infty$ is metrizable 
(cf. p. \pageref{11.5}) ($B_G^\infty$ is even an ANR 
(cf. p. \pageref{11.6})), in particular, $B_G^\infty$ is a compactly generated 
%%----------------------------------------------------------------------------------------------04
Hausdorff space.  And: For every pointed finite dimensional connected CW complex \mX, 
$\map_*(B_G^\infty,\pro X)$ is homotopically trivial 
(Friedlander-Mislin\footnote[2]{\textit{Invent. Math.} \textbf{83} (1986), 425-436.}).

[Note: \ Taking $G = \bS^1$, the Zabrodsky lemma and induction imply that $\forall \ n \geq 2$, 
$\map_*(K(\Z,n),\pro X)$ is homotopically trivial.]\\

\endgroup %%------------------------------------<<

\begingroup%%----------------------------------->>
\fontsize{9pt}{11pt}\selectfont
\textbf{\small FACT} \quad 
Let \mX  be a pointed connected CW complex $-$then for any CW complex \mY, the arrow 
$[Y,\pro X] \ra \lim\limits_X[Y,K]$ is bijective.

[Note: \ In this context, the bracket refers to homotopy classes of maps, not to pointed homotopy classes of pointed maps.]\\

\endgroup %%------------------------------------<<

The homotopy groups of $\pro X$ are profinite: 
Proof: $\pi_n(\pro X) \approx$
$[\bS^n,\pro X] \approx$ 
$\lim\limits_X[\bS^n,K]$ and the $[\bS^n,K]$ are finite.

[Note: \  It follows that $\forall \ n$, there is a commutative triangle\\
%
%
%$
\[
\begingroup%%----------------------------------->>
%\fontsize{9pt}{11pt}\selectfont
%dmc - for some reason this is problematic
\begin{tikzcd}[sep=tiny]
&{\pi_n(X)} \ar{ldddddd}\ar{rdddddd}\\
\\
\\
\\
\\
\\
{\pro \pi_n(X)} \ar{rr} &&{\pi_n(\pro X)}
\end{tikzcd}
\endgroup
.]
\]
%$
%\\
\vspace{0.25cm}

\begin{proposition} \ %2
Let \mX be a pointed connected CW complex $-$then $\pi_1(\pro X) \approx$ $\pro \pi_1(X)$.
\end{proposition}

[The full subcategory of $X\backslash \ov{\bH\bF}_*$ consisting of those objects $X \ra K$ such that the induced map 
$\pi_1(X) \ra \pi_1(K)$ is surjective is an initial subcategory.  To see this, let 
$\widetilde{K} \ra K$ be the covering of \mK corresponding to $\im \pi_1(X)$, and consider 
\begin{tikzcd}%[ sep=small]
&{\widetilde{K}} \ar{d}\\
{X} \ar{ru} \ar{r} &{K}
\end{tikzcd}
.  On the other hand, for any normal subgroup \mG of $\pi_1(X)$ of finite index, there is an arrow 
$X \ra K(\pi_1(X)/G,1)$.]\\

\begingroup%%----------------------------------->>
\fontsize{9pt}{11pt}\selectfont
\textbf{\small EXAMPLE} \quad 
The arrow $\pro \pi_n(X) \ra \pi_n(\pro X)$ is not necessarily bijective when $n > 1$.  Thus take 
$X = \bS^1 \vee \Sigma \bP^2 (\R)$ $-$then 
$\pi_1(X) \approx \Z$, 
$\pi_2(X) \approx \omega \cdot (\Z/2\Z)$ and 
$\pi_1(\pro X) \approx \widehat{\Z}$, 
$\pi_2(\pro X) \approx (\Z/2\Z)^\omega$  but 
$\pro \pi_2(X) \approx \Hom((\Z/2\Z)^\omega,\Z/2\Z)$.\\

\endgroup %%------------------------------------<<

\textbf{\small LEMMA} \quad 
Suppose that \mG is a finitely generated abelian group $-$then $\pro K(G,n) \approx$ $K(\pro G,n)$.\\

\begingroup%%----------------------------------->>
\fontsize{9pt}{11pt}\selectfont
\textbf{\small EXAMPLE} \quad 
$\pro K(\Z,n) \approx K(\widehat{\Z},n)$ but $\pro K(\Q/\Z,n) \approx K(\widehat{\Z},n+1)$.\\

\endgroup %%------------------------------------<<

%%----------------------------------------------------------------------------------------------05

\label{11.25}
\begingroup%%----------------------------------->>
\fontsize{9pt}{11pt}\selectfont

\textbf{\small EXAMPLE} \quad 
Consider $K(\Z,2;\chi)$, where $\chi:\Z/2\Z \ra \Aut \Z$ is the nontrivial homomorphism (so 
$K(\Z,2;\chi) \approx$ 
$B_{\bO(2)}$ 
(cf. p. \pageref{11.7}) $-$then $\chi$ extends to a homomorphism 
$\widehat{\chi}:\Z/2\Z \ra \Aut \widehat{\Z}$ and $\pro K(\Z,2;\chi) \approx K(\widehat{\Z},2;\widehat{\chi})$.\\ 

\endgroup %%------------------------------------<<

\begingroup%%----------------------------------->>
\fontsize{9pt}{11pt}\selectfont

\textbf{\small FACT} \quad 
Let \mX be a pointed connected CW complex $-$then $\forall \ q$, 
$H^q(X;\widehat{\Z}) \approx \lim_n H^q(X;\Z/n\Z)$.

[$H^q(X;\widehat{\Z}) \approx [X,K(\widehat{\Z},q)] \approx [X,\pro K(\Z,q)] \approx
\lim_n [X,K(\Z/n\Z,q)] \approx \lim_n H^q(X;\Z/n\Z)$.]\\

\endgroup %%------------------------------------<<

\begingroup%%----------------------------------->>
\fontsize{9pt}{11pt}\selectfont

\textbf{\small FACT} \quad 
Let \mX be a pointed connected CW complex $-$then $\forall \ q$, 
$H^q(X;\widehat{\Z}) \approx \lim H^q(X_k;\widehat{\Z})$, where $X_k$ runs over the pointed finite connected subcomplexes of \mX.

[$H^q(X;\Z/n\Z) \approx [X,K(\Z/n\Z,q)] \approx 
\lim[X_k,K(\Z/n\Z,q)] \approx \lim H^q(X_k;\Z/n\Z)$ 
(cf. p. \pageref{11.8}) $\implies$ 
$H^q(X;\widehat{\Z}) \approx$
$\lim_n H^q(X;\Z/n\Z) \approx$ 
$\lim_n \lim H^q(X_k;\Z/n\Z) \approx$
$\lim \lim_n H^q(X_k;\Z/n\Z) \approx$
$\lim  H^q(X_k;\widehat{\Z})$.]\\
 
\endgroup %%------------------------------------<<

In general, it is difficult to relate the higher homotopy groups of $\pro X$ to those of \mX itself except under the most favorable conditions.\\

\begin{proposition} \ %3
Let \mX be a pointed nilpotent CW space with finitely generated homotopy groups $-$then $\forall \ n$, 
$\pi_n(\pro X) \approx \pro \pi_n(X)$.
\end{proposition}

[Note: \ Recall that a particular choice for the abelian groups figuring in a principal refinement of order $n$ of 
$X[n] \ra X[n-1]$ are the 
$\Gamma_{\chi_n}^i(\pi_n(X))/\Gamma_{\chi_n}^{i+1} (\pi_n(X))$ 
(cf. p. \pageref{11.9}).  
Since the $\pi_n$ are finitely generated, there is a unique continuous nilpotent action of $\pro \pi_1(X)$ on 
$\pro \pi_n(X)$ compatible with the action of $\pi_1(X)$ on $\pi_n(X)$.  This said, 
Hilton-Roitberg\footnote[2]{ \textit{J. Algebra} \textbf{60}  (1979), 289-306.}
have shown that, in obvious notation 
(i) $\nil_{\chi_n} \pi_n(X) = \nil_{\pro\chi_n} \pro\pi_n(X)$ and 
(ii) $\pro(\Gamma_{\chi_n}^i(\pi_n(X))/\Gamma_{\chi_n}^{i+1} (\pi_n(X))) \approx 
\Gamma_{\pro\chi_n}^i(\pro\pi_n(X))/\Gamma_{\pro\chi_n}^{i+1} (\pro\pi_n(X))$.  
Since profinite completion preserves short exact sequences of finitely generated nilpotent groups 
(cf. p. \pageref{11.10}), 
the conclusion is that the arrow 
$(\pro X)[n] \ra (\pro X)[n-1]$ admits a ``canonical'' principal refinement of order $n$, viz. apply pro to the ``canonical'' principal refinement of order $n$ of $X[n] \ra X[n-1]$.  
Corollary: Under the stated assumptions on \mX, $\pro X$ is nilpotent (but the unconditional assertion 
``\mX nilpotent $\implies$ $\pro X$ nilpotent'' is seemingly in limbo.]

Example: $\bS^n = M(\Z,n)$ but $\pro \bS^n \neq M(\pro \Z,n)$.\\

\begingroup%%----------------------------------->>
\fontsize{9pt}{11pt}\selectfont
\textbf{\small FACT} \quad 
Let \mX be a pointed nilpotent CW space with finitely generated homotopy groups $-$then 
for every pointed finite connected CW complex \mK, the arrow 
$[K,X] \ra [K,\pro X]$ is injective.

[Note: \ As a reality check, take $K = \bS^1$ and $X = K(G,1)$, where \mG is a finitely generated nilpotent group, and observe that the injectivity of the arrow 
$[\bS^1,K(G,1)] \ra [\bS^1,K(\pro G,1)]$ is equivalent to the assertion that \mG embeds in $\pro G$ 
(cf. p. \pageref{11.11}).]\\

\endgroup %%------------------------------------<<

%%----------------------------------------------------------------------------------------------06
\begingroup%%----------------------------------->>
\fontsize{9pt}{11pt}\selectfont
Application: Let \mY be a pointed nilpotent CW space with finitely generated homotopy groups $-$then for every pointed connected CW space \mX, $\Ph(X,Y)$ is the kernel of the arrow $[X,Y] \ra [X,\pro Y]$.\\

\endgroup %%------------------------------------<<

\begingroup%%----------------------------------->>
\fontsize{9pt}{11pt}\selectfont
\textbf{\small LEMMA} \quad 
Let $\{G_n,f_n:G_{n+1} \ra G_n\}$ be a tower in \bGR.  Assume: $\forall \ n$, $G_n$ is a compact Hausdorff topological group and $f_n$ is a continuous homomorphism $-$then $\lim^1 G_n = *$.

[Note: The result is false if the ``Hausdorff'' hypothesis is dropped.]\\

\endgroup %%------------------------------------<<

\begingroup%%----------------------------------->>
\fontsize{9pt}{11pt}\selectfont
\textbf{\small EXAMPLE} \quad 
Let \mX be a pointed connected CW complex with a finite number of cells in each dimension; let \mY be a pointed nilpotent CW space with finitely generated homotopy groups $-$then $\forall \ n$, 
$[\Sigma X^{(n)}, \pro Y]$ is a compact Hausdorff topological group and the arrow 
$[\Sigma X^{(n+1)}, \pro Y] \ra$ $[\Sigma X^{(n)}, \pro Y]$ is a continuous homomorphism.  So, by the lemma, 
$\lim^1[\Sigma X^{(n)}, \pro Y] = *$, i.e., $\Ph(X,\pro Y) = *$, 
(cf. p. \pageref{11.12}).

Claim: A pointed continuous function $f:X \ra Y$ is a phantom map iff $\text{pro$_{Y}$} \circx f \simeq 0$.

[Necessity: \ $f \in \Ph(X,Y)$ $\implies$ 
$\text{pro$_{Y}$} \circx f \in \Ph(X, \pro Y)$ $\implies$ 
$\text{pro$_{Y}$} \circx f \simeq 0$.

Sufficiency: Let $\phi:K \ra X$ be a pointed continuous function, where \mK is a pointed finite connected CW complex 
$-$then 
$\text{pro$_{Y}$} \circx f \circx \phi \simeq 0$ $\implies$ 
$f \circx \phi \simeq 0$, the arrow 
$[K,Y] \ra [K,\pro Y]$ being one-to-one.]\\

\endgroup %%------------------------------------<<

\begingroup%%----------------------------------->>
\fontsize{9pt}{11pt}\selectfont
\textbf{\small LEMMA} \quad 
Let 
$
\begin{cases}
\ X\\
\ Y
\end{cases}
$
be pointed simply connected CW spaces with finitely generated homotopy groups $-$then the function space of pointed continuous functions $X_{\Q} \ra \pro Y$ is homotopically trivial (compact open topology).

[Adopt the conventions on 
p. \pageref{11.13} and work with 
$\map_*(X_{\Q}, \pro Y)$.  Since $\Sigma^n X_{\Q} \approx (\Sigma^n X)_{\Q}$ 
(cf. p. \pageref{11.14}), 
$\widetilde{H}_*(\Sigma^n X_{\Q} ;\F_p) = 0$ $\forall \ p$, thus 
$\widetilde{H}^*(\Sigma^n X_{\Q} ;\pi_q(\pro Y)) \ = \ 0$ $\forall \ q$ \ (the $\pi_q(\pro Y)$ are cotorsion).  
Accordingly, by obstruction theory 
(cf. p. \pageref{11.15}), $\forall \ n \geq 0$, 
$[\Sigma^n X_{\Q},\pro Y] = *$.]\\

\endgroup %%------------------------------------<<

\label{9.23}
\begingroup%%----------------------------------->>
\fontsize{9pt}{11pt}\selectfont
\textbf{\small EXAMPLE} \quad 
Let 
$
\begin{cases}
\ X\\
\ Y
\end{cases}
$
be pointed simply connected CW spaces with finitely generated homotopy groups $-$then 
$\Ph(X,Y) = l_{\Q}^*[X_{\Q},Y] \subset [X,Y]$.

[There is no loss of generality is supposing that \mX is a pointed simply connected CW complex with a finite number of cells in each dimension 
(cf. p. \pageref{11.16}).

\indent\indent $(\Ph(X,Y) \ \subset \  l_{\Q}^*[X_{\Q},Y])$ \quad Fix an $f \in \Ph(X,Y)$.  From the above, 
$\text{pro$_Y$} \circx f \simeq 0$, so $\exists$ a $g:X \ra E$ such that $f = \pi \circx g$, \mE the mapping fiber of 
$\text{pro$_Y$}$ and $\pi:E \ra Y$ the projection.  
Since \mE is rational (each of its homotopy groups is a direct sum of copies of $\widehat{\Z}/\Z$), $\exists$ an $h:X_{\Q} \ra E$ such that $g \simeq h \circx l_{\Q}$, thus $f \simeq f_{\Q} \circx l_{\Q}$, where $f_{\Q} = \pi \circx h: X_{\Q} \ra Y$.

\indent\indent $(l_{\Q}^*[X_{\Q},Y])  \ \subset \ \Ph(X,Y)$ \quad Assume that $f \simeq f_{\Q} \circx l_{\Q}$, where 
$f_{\Q}:X_{\Q} \ra Y$.  Thanks to the lemma, the composite $\text{pro$_Y$} \circx f_{\Q}$ is nullhomotopic, hence 
$\text{pro$_Y$} \circx f$ is too.]\\

\endgroup %%------------------------------------<<

\begingroup%%----------------------------------->>
\fontsize{9pt}{11pt}\selectfont
\textbf{\small FACT} \quad 
Let \mX be a pointed nilpotent CW space with finitely generated homotopy groups $-$then for every finite CW complex \mK, 
the arrow 
$[K,X] \ra [K,\pro X]$ is injective.
%%----------------------------------------------------------------------------------------------07

[Note: \ In this context, the brackets refer to homotopy classes of maps, not to pointed homotopy classes of pointed maps.]\\

\endgroup %%------------------------------------<<

\begingroup%%----------------------------------->>
\fontsize{9pt}{11pt}\selectfont
\textbf{\small EXAMPLE} \quad 
The preceding result has content even when \mK is connected.  Thus, restoring the base points, it follows that the arrow 
$\pi_1(X)\backslash[K,k_0;X,x_0] \ra$ 
$\pi_1(\pro X)\backslash[K,k_0;\pro X,\pro x_0]$ is one-to-one.  Specializing this to 
$K = \bS^1$, $X = K(G,1)$, where \mG is a finitely generated nilpotent group, one recovers Blackburn's theorem 
(cf. p. \pageref{11.17}).\\

\endgroup %%------------------------------------<<

\begin{proposition} \ %4
Let \mX be a pointed nilpotent CW space with finitely generated homotopy groups $-$then for every locally constant coefficient system $\sG$ on $\pro X$ arising from a finite $\pro \pi_1(X)$-module, 
$H^*(\pro X;\sG) \approx H^*(X; \text{pro}_X^*\sG)$.
\end{proposition}

[The main idea here is to proceed inductively, playing off 
$K(\pi_n(X),n) \ra$ 
$P_n X \ra$ 
$P_{n-1} X$ 
against 
$\pro K(\pi_n(X),n) \ra$ 
$\pro P_n X \ra$ 
$\pro P_{n-1} X$ 
(use the cohomological version of the fibration spectral sequence formulated on 
p. \pageref{11.18}).  To get the induction off the ground, one has to deal with $K(\pi_1(X),1)$, the point being that $\pi_1(X)$ has property S 
(cf. p. \pageref{11.19}).]\\

\textbf{\small LEMMA} \quad 
Let 
$
\begin{cases}
\ X\\
\ Y
\end{cases}
\& \ Z
$
be pointed connected CW spaces, $f:X \ra Y$ a pointed continuous function $-$then the precomposition arrow 
$f^*:[Y,Z] \ra [X,Z]$ is bijective whenever \mZ has finite homotopy groups iff
\[
\text{and}
\begin{array}{l}
\indent \text{(A$_1$) \quad $\Hom(\pi_1(Y),F) \approx \Hom(\pi_1(X),F)$ for any finite group \mF}\\
\indent \text{(A$_2$) \quad $H^n(Y;\sG) \approx H^n(X;f^*\sG)$ $\forall \ n$ for any locally constant coefficient system}
\end{array}
\]
$\sG$ on \mY arising from a finite $\pi_1(Y)$-module.

[Tailor the proof of Proposition 11 in $\S 9$ to the setup at hand.]\\

\begin{proposition} \ %5
Let \mX be a pointed nilpotent CW space with finitely generated homotopy groups $-$then every pointed continuous function 
$\phi:X \ra K$, where \mK is a pointed connected CW complex with finite homotopy groups, admits a continuous extension 
$\pro \phi:\pro X \ra K$ which is unique up to pointed homotopy.
\end{proposition}

[Each homomorphism $\pi_1(X) \ra F$, where \mF is finite, can be extended uniquely to a homomorphism 
$\pro \pi_1(X) \ra F$ 
(cf. p. \pageref{11.20} ff.), therefore A$_1$ holds.  That A$_2$ holds is the content of Proposition 4.]\\

Application: pro is idempotent on the class of pointed nilpotent CW spaces with finitely generated homotopy groups.\\

Fix a prime $p$ $-$then upon replacing ``finite group'' by ``finite $p$-group'' in the foregoing, one arrives at the 
\un{$p$-profinite completion}
\index{p-profinite completion (spaces)} 
$\prox_p X$ of \mX.  Modulo minor changes, the theory
%%----------------------------------------------------------------------------------------------08
carries over in the expected way.  Consider, e.g., Proposition 4.  
There it is necessary to look only at those $\sG$ whose underlying $\prox_p  \pi_1(X)$-module 
\mG is a finite abelian $p$-group such that the associated homomorphism 
$\prox_p \pi_1(X) \ra \Aut G$ factors through a $p$-subgroup of $\Aut G$.  
Another point to bear in mind is that $p$-adic completion preserves short exact sequences of finitely generated nilpotent groups 
(cf. p. \pageref{11.21}) and $p$-adic completion = $p$-profinite completion in the class of finitely generated nilpotent groups 
(cf. p. \pageref{11.22}).\\

\begingroup%%----------------------------------->>
\fontsize{9pt}{11pt}\selectfont
\textbf{\small EXAMPLE} \quad 
Let \mX be a pointed simply connected CW complex with a finite number of cells in each dimension.  Denote by 
$\text{pro$_{p,T }$} X$ the pointed mapping telescope of the sequence 
$\{\prox_p X^{(n)} \ra$ $\prox_p X^{(n+1)}\}$ $-$then $\forall \ n$, 
$\pi_n(\text{pro$_{p,T }$} X) \approx \widehat{\Z}_p \otimes \pi_n(X)$ $\implies$ 
$\text{pro$_{p,T }$} X \approx \prox_p X$.\\

\endgroup %%------------------------------------<<

It is clear that $\forall \ p$, there is an arrow 
$\pro X \ra \prox_p X$ from which the arrow $\pro X \ra \prod\limits_p \prox_p X$ (product in $\bHTOP_*$).  
And: 
$[\bS^n,\pro X] \ra [\bS^n,\prod\limits_p \prox_p X]$ $\implies$ $\pi_n(\pro X) \ra \prod\limits_p \pi_n(\prox_p X)$.\\

\begin{proposition} \ %6
Let \mX be a pointed nilpotent CW space with finitely generated homotopy groups $-$then the arrow 
$\pro X \ra \prod\limits_p \prox_p X$ is a weak homotopy equivalence. 
\end{proposition}

[In this situation, $\forall \ n$, 
$\pi_n(\pro X) \approx \pro \pi_n(X)$ $\&$ 
$\pi_n(\prox_p X) \approx \prox_p \pi_n(X)$.  
Moreover, for any finitely generated nilpotent group \mG, the arrow 
$\pro G \ra \prod\limits_p \prox_p G$ is an isomorphism  
(cf. p. \pageref{11.23}).

[Note: \ If the product is taken over $\bHCWSP_*$ 
(cf. p. \pageref{11.24}), then the arrow 
$\pro X \ra \prod\limits_p \prox_p X$ is a pointed homotopy equivalence.]\\

\begingroup%%----------------------------------->>
\fontsize{9pt}{11pt}\selectfont
\textbf{\small EXAMPLE} \quad 
Let $X = B_{\bO(2)}$ 
(cf. p. \pageref{11.25}) $-$then, in obvious notation, 
$\text{pro}_2 X \approx K(\widehat{\Z}_2,2;\widehat{\chi}_2)$ but at an odd prime $p$, $\prox_p X$ is simply connected and in fact $\Omega \prox_p X \approx 
%\widehat{\bS}_p^3 %dmc note the difference - and comp with G -note neither mathc exactly but he tends the hat only over the 
\widehat{\bS_p^3}$ %S - which actually seems to look better - whereas the second is really the intent
.  
Thus here, it is false that the arrow 
$\pro X \ra \prod\limits_p \prox_p X$ is a weak homotopy equivalence.\\

\endgroup %%------------------------------------<<

Let \mX be a pointed nilpotent CW space $-$then $\prox_p X$ and $X_{H\F_p}$ ($= \F_p X$) are, in general, not the ``same''.  
Reason: $\prop$ fails to be idempotent.  However, when the homotopy groups of \mX are finitely generated, 
$\pi_n(\prox_p X) \approx$ 
$\prox_p \pi_n(X) \approx$ 
$\Ext(\Z/p^\infty\Z,\pi_n(X)) \approx$ 
$\pi_n(X_{H\F_p})$.  Therefore $\prox_p X$ is $H\F_p$-local (cf. $\S 9$, Proposition 20) ($\prox_p X$ is nilpotent) and in this case, 
$\prox_p X \approx$ $X_{H\F_p}$ $(= \widehat{X}_p)$.

[Note: \ It is a fact that for nilpotent \mX, $\prox_p X \approx X_{H\F_p}$ under the sole hypothesis that $\forall \ n$, 
$H^n(X;\F_p)$ is finite dimensional 
(cf. p. \pageref{11.26}).  In this connection, recall that if the homotopy groups of a nilpotent \mX are finitely generated, then the $H_n(X)$ are finitely generated (cf. $\S 5$, Proposition 18), hence $\forall \ n$, 
$H^n(X;\F_p)$ is finite dimensional.]\\

%%----------------------------------------------------------------------------------------------09
\begin{proposition} \ %7
Let \mX be a path connected topological space $-$then the following conditions are equivalent:

\indent\indent \textnormal{(CO$_1$)} \quad $\forall \ n$, $H^n(X;\F_p)$ is finite dimensional; 

\indent\indent \textnormal{(HO$_1$)} \quad $\forall \ n$, $H_n(X;\F_p)$ is finite dimensional; 

\indent\indent \textnormal{(CO$_2$)} \quad $\forall \ n$, $H^n(X;\widehat{\Z}_p)$  is finitely generated over $\widehat{\Z}_p$; 

\indent\indent \textnormal{(HO$_2$)} \quad $\forall \ n$, $H_n(X;\widehat{\Z}_p)$ is finitely generated over $\widehat{\Z}_p$; 

\indent\indent \textnormal{(CO$_3$)} \quad $\forall \ n$, $H^n(X;\Z_p)$ is finitely generated over $\Z_p$;  

\indent\indent \textnormal{(HO$_3$)} \quad $\forall \ n$, $H_n(X;\Z_p)$ is finitely generated over $\Z_p$; 

\indent\indent \textnormal{(CO$_4$)} \quad $\forall \ n$, $H^n(X;\Q)$ is finite dimensional and $H^n(X;\Z)_\tor(p)$ is finite; 

\indent\indent \textnormal{(HO$_4$)} \quad $\forall \ n$, $H_n(X;\Q)$ is finite dimensional and $H_n(X;\Z)_\tor(p)$ is finite.\\
\end{proposition}

\begingroup%%----------------------------------->>
\fontsize{9pt}{11pt}\selectfont
\textbf{\small EXAMPLE} \quad 
Suppose that \mX is a pointed simply connected CW space which is $H\F_p$-local $-$then $H^n(X;\F_p)$ is finite dimensional 
$\forall \ n$ iff $\pi_n(X)$ is a finitely generated $\widehat{\Z}_p$-module $\forall \ n$.

[Note: \ $\pi_n(X)$ is $p$-cotorsion, hence is a $p$-adic module 
(cf. p. \pageref{11.27}).]\\

\endgroup %%------------------------------------<<

A group \mG is said to be 
\un{$\F_p$-finite}
\index{F$_p$-finite (group)} 
provided that 
$H^1(G;\F_p)$ and $H^2(G;\F_p)$ are finite dimensional.  
Example: Every finitely generated nilpotent group is $\F_p$-finite 
(cf. p. \pageref{11.28}).

[Note: \ Let \mG be an abelian group $-$then \mG is $\F_p$-finite iff $G \otimes \F_p$ and $\Tor(G,\F_p)$ are finite or still, \mG is $\F_p$-finite iff 
$H^n(G,n;\F_p)$ and $H^{n+1}(G,n;\F_p)$ are finite dimensional.]\\

\label{11.34}
\begingroup%%----------------------------------->>
\fontsize{9pt}{11pt}\selectfont
\textbf{\small EXAMPLE} \quad 
Suppose that \mG is $\F_p$-finite $-$then $H_1(G;\F_p)$ and $H_2(G;\F_p)$ are finite dimensional.  
Therefore, $L_0^pG \approx \F_p G$ 
(cf. p. \pageref{11.29}).  In particular, for any nilpotent $\F_p$-finite group \mG, 
$\Ext(\Z/p^\infty\Z,G) \approx$ 
$\F_p G \approx$ 
$\widehat{G}_p \approx$ 
$\prox_p G$.

[Note: \ In the abelian case, one may proceed directly.  Thus observe first that if \mG is abelian and $\F_p$-finite, then 
$\forall \ n$, $\Tor(G,\Z/p^n\Z)$ is finite (argue by induction, using the coefficient sequence associated with the short exact sequence 
$0 \ra$
$\Z/p \Z \ra$ 
$\Z/p^{n+1}\Z \ra$ 
$\Z/p^n\Z \ra 0$).  Accordingly, $\forall \ n$, 
$\Hom(\Z/p^n\Z,G)$ is finite $\implies$ $\lim^1 \Hom(\Z/p^n\Z,G) = 0$ 
(cf. p. \pageref{11.30}) $\implies$ 
$\Ext(\Z/p^\infty\Z,G) \approx \widehat{G}_p$
(cf. p. \pageref{11.31}).]\\

\endgroup %%------------------------------------<<

\begingroup%%----------------------------------->>
\fontsize{9pt}{11pt}\selectfont
\textbf{\small EXAMPLE} \quad 
Any abelian group in any of the following four classes is $\F_p$-finite:
(C$_1$) \  The finite abelian $p$-groups; 
(C$_2$) \  The free abelian groups of finite rank; 
(C$_3$) \  The uniquely $p$-divisible abelian groups; 
(C$_4$) \  The $p$-primary divisible abelian groups satisfying the descending chain condition on subgroups.  
Moreover, every $\F_p$-finite abelian group \mG admits a composition series 
$G = G^0 \supset G^1 \supset \cdots \supset$ $G^n = \{0\}$ such that $\forall \ i$, $G^i/G^{i+1}$ is in one of these four classes.

[Given an $\F_p$-finite abelian group \mG, $\exists$ a short exact sequence 
$0 \ra$ 
$G^\prime \ra$ 
$G \ra$ 
$G\pp \ra 0$, where 
$
\begin{cases}
\ G^\prime\\
\ G\pp
\end{cases}
$
are $\F_p$-finite with $G^\prime$ finitely generated and $G\pp$ $p$-divisible.  
Proof: One may take $G\pp = G/G^\prime$, where $G^\prime$ is a finitely generated subgroup of \mG mapping onto 
$G/pG$.]\\

\endgroup %%------------------------------------<<

%%----------------------------------------------------------------------------------------------10

\begingroup%%----------------------------------->>
\fontsize{9pt}{11pt}\selectfont
\textbf{\small FACT} \quad 
Let \mG be an abelian group.  Assume: \   \mG is $\F_p$-finite $-$then $\forall \ n$, $H^n(G;\F_p)$ is finite dimensional.\\

\endgroup %%------------------------------------<<

\begin{proposition}
Let \mG be an $\F_p$-finite nilpotent group $-$then $\forall \ n$, $H^n(G;\F_p)$ is finite dimensional.
\end{proposition}

[This is true if \mG is abelian (cf. supra).  Since in general, the iterated commutator map 
$\otimes^{i+1} (G/[G,G]) \ra \Gamma^i(G)/\Gamma^{i+1}(G)$ is surjective, 
$H^1(\Gamma^i(G)/\Gamma^{i+1}(G);\F_p)$ is finite dimensional $\forall \ i$.  In particular: 
$H^1(\Gamma^{d-1}(G);\F_p)$ is finite dimensional ($d = \nil G > 1$).  Put 
$K = \Gamma^{d-1}(G)$ and consider the central extension 
$1 \ra$ 
$K \ra$ 
$G \ra$ 
$G/K \ra 1$.  The associated LHS spectral sequence is 
$H^p(G/K;H^q(K;\F_p)) \Rightarrow H^{p+q}(G;\F_p)$, so it need only be shown that the $E_2^{p,q}$ are finite dimensional.  Specialized to the present situation, the fundamental exact sequence in cohomomolgy reads 
$0 \ra$ 
$H^1(G/K;\F_p) \ra$ 
$H^1(G;\F_p) \ra$ 
$H^1(K;\F_p) \ra$
$H^2(G/K;\F_p) \ra$ 
$H^2(G;\F_p)$  
(cf. p. \pageref{11.32}).  Therefore 
$H^1(G/K;\F_p)$ and $H^2(G/K;\F_p)$ are finite dimensional, hence by induction, $\forall \ n$, $H^n(G/K;\F_p)$ is finite dimensional.  
Claim:  $H^2(K;\F_p)$ is finite dimensional.  To see this, suppose the contrary.  
Because 
$\dim E_2^{2,1} < \omega$, $E_3^{0,2}$ (the kernel of the differential $E_2^{0,2} \ra E_2^{2,1}$) would be infinite dimensional.  
But $\dim E_2^{3,0} < \omega$ $\implies$ $\dim E_3^{3,0} < \omega$, which means that $E_4^{0,2}$ (the kernel of the differential $E_3^{0,2} \ra E_3^{3,0}$) would be infinite dimensional.  This, however, is untenable: 
$E_4^{0,2} \ra E_\infty^{0,2}$ and $H^2(G;\F_p)$ is finite dimensional.  Thus the conclusion is that \mK is $\F_p$-finite and, being abelian, $H^n(K;\F_p)$ is finite dimensional $\forall \ n$.  It now follows that $\forall \ p$ $\&$ $\forall \ q$,  
$E_2^{p,q}$ is finite dimensional.]\\

Application: Let \mG be an $\F_p$-finite nilpotent group $-$then $\forall \ i$, $\Gamma^i(G)/\Gamma^{i+1}(G)$ is an 
$\F_p$-finite abelian group.\\

\begingroup%%----------------------------------->>
\fontsize{9pt}{11pt}\selectfont
\textbf{\small FACT} \quad 
Let \mG be a group, \mM a nilpotent $G$-module.  Assume: $H^1(G;\F_p)$ is finite dimensional and \mM is $\F_p$-finite $-$then $\forall \ i$,  $\Gamma\chi^i(M)/\Gamma\chi^{i+1}(M)$  is $\F_p$-finite.\\

\endgroup %%------------------------------------<<

\textbf{\small LEMMA} \quad 
Let \mG be a group, \mM a nilpotent $G$-module which is a vector space over $\F_p$.  Assume; $H^1(G;\F_p)$ is finite dimensional and $H^0(G;M)$ is finite dimensional $-$then \mM is finite dimensional.

[The assertion is clear if \mG operates trivially on \mM.  
Agreeing to argue inductively on $d = \nil_\chi M > 1$, put 
$N = \Gamma\chi^{d-1}(M)$ and consider the exact sequence 
$0 \ra$ 
$H^0(G;N) \ra$ 
$H^0(G;M) \ra$ 
$H^0(G;M/N) \ra$ 
$H^1(G;N) \ra \cdots$.  \ 
Since \mG operates trivially on \mN, $H^0(G;N) = N$, thus \mN is finite dimensional.  Consequently, 
$H^1(G;N)$ is finite dimensional, so $H^0(G;M/N)$ is finite dimensional.  Owing to the induction hypothesis, $M/N$ is finite dimensional, hence the same holds for \mM itself.]\\

%%----------------------------------------------------------------------------------------------11
\begin{proposition}
Let \mX be a pointed nilpotent CW space $-$then $\forall \ n$, $H^n(X;\F_p)$ is finite dimensional iff $\forall \ n$, 
$\pi_n(X)$ is $\F_p$-finite.
\end{proposition}

[We shall prove that the condition on the homotopy groups is necessary, the verifiction that it is also sufficient being similar.  For this, consider the 5-term exact sequence 
$0 \ra$ 
$E_2^{1,0} \ra$ 
$H^1(X;\F_p) \ra$ 
$E_2^{0,1} \ra$ 
$E_2^{2,0} \ra$ 
$H^2(X;\F_p)$ associated with the fibration spectral sequence 
$H^p(\pi_1(X);H^q(\widetilde{X};\F_p)) \Rightarrow H^{p+q}(X;\F_p)$ to see that 
$H^1(\pi_1(X);\F_p)$ and $H^2(\pi_1(X);\F_p)$ are finite dimensional, i.e., that $\pi_1(X)$ is $\F_p$-finite.
Since $\pi_1(X)$ operates nilpotently on the $H_n(\widetilde{X})$ (cf. $\S 5$, Proposition 17), 
$H_n(\widetilde{X};\F_p)$ is a nilpotent $\pi_1(X)$-module, as is its dual $H^n(\widetilde{X};\F_p)$.  
Taking into account Proposition 8, one finds from the lemma that $H^2(\widetilde{X};\F_p)$ is finite dimensional and then by iteration that $H^n(\widetilde{X};\F_p)$ is finite dimensional $\forall \ n$.  This sets the stage for the discussion of $\pi_2(X)$.  
Thus, in the notation of 
p. \pageref{11.33}, consider 
$\widetilde{X}_2 \ra$
$\widetilde{X}_1 \ra$
$K(\pi_2(X),2)$ $(\widetilde{X}_1 \approx \widetilde{X})$.  
Once again, there is a fibration spectral sequence 
$H^p(K(\pi_2(X),2);H^q(\widetilde{X}_2;\F_p)) \Rightarrow$ 
$H^{p+q}(\widetilde{X}_1;\F_p)$ and a low degree exact sequence 
$H^2(\pi_2(X),2;\F_p) \ra$
$H^2(\widetilde{X}_1;\F_p) \ra$
$H^2(\widetilde{X}_2;\F_p) \ra$
$H^3(\pi_2(X),2;\F_p) \ra$
$H^3(\widetilde{X}_1;\F_p)$.  
Because $H^2(\widetilde{X}_1;\F_p)$ and $H^3(\widetilde{X}_1;\F_p)$ are finite dimensional and 
$H^2(\widetilde{X}_2;\F_p) = 0$, it follows that 
$H^2(\pi_2(X),2;\F_p)$ and $H^3(\pi_2(X),2;\F_p)$ are finite dimensional.  Therefore $\pi_2(X)$ is $\F_p$-finite and the process can be continued.]\\

\begingroup%%----------------------------------->>
\fontsize{9pt}{11pt}\selectfont
\textbf{\small FACT} \quad 
Let \mG be an $\F_p$-finite nilpotent group $-$then $\prox_p G$ operates nilpotently on the $L_n^p G$.\\

\endgroup %%------------------------------------<<

\begingroup%%----------------------------------->>
\fontsize{9pt}{11pt}\selectfont
\textbf{\small FACT} \quad 
Let \mG be an $\F_p$-finite nilpotent group, \mM an $\F_p$-finite nilpotent $G$-module $-$then $\prox_p G$ operates nilpotently on the $L_n^p M$.\\

\endgroup %%------------------------------------<<

\label{11.26}
\begingroup%%----------------------------------->>
\fontsize{9pt}{11pt}\selectfont
\index{coincidence criterion}
\textbf{\small COINCIDENCE CRITERION} \quad 
Let \mX be a pointed nilpotent CW space such that $\forall \ n$, $H^n(X;\F_p)$ is finite dimensional $-$then $\forall \ n$, there is a split short exact sequence 
$0 \ra$
$\prox_p \pi_n(X) \ra$ 
$\pi_n(\prox_p X) \ra$ 
$\Hom(\Z/p^\infty\Z,\pi_{n-1}(X)) \ra 0$, hence $\prox_p X \approx X_{H\F_p}$.

[Note: \ Recall that here, 
$\Ext(\Z/p^\infty\Z,\pi_n(X)) \approx \F_p \pi_n(X) \approx \pi_n(X)_p^{\widehat{}} \approx \prox_p \pi_n(X)$ 
(cf. p. \pageref{11.34}).]\\

\endgroup %%------------------------------------<<

\begingroup%%----------------------------------->>
\fontsize{9pt}{11pt}\selectfont
\textbf{\small EXAMPLE} \quad 
Let \mX be a pointed nilpotent CW space such that $\forall \ n$, $H^n(X;\F_p)$ is finite dimensional.  
Let $\sA_p$ be the 
mod $p$ Steenrod algebra $-$then $H^*(X;\F_p)$ is an unstable $\sA_p$-module and 
Lannes-Schwartz\footnote[2]{\textit{Topology} \textbf{28} (1989), 153-169.}
have shown that \mX is $W$-null, where $W = B\Z/p\Z$, iff every cyclic submodule of $H^*(X;\F_p)$ is finite.\\

\endgroup %%------------------------------------<<
%%%%%%%%%%%%%%%%%%%%%%%%%%%%%%%%%%%%%%
%%%%%%%%%%%%%%%%%%%%%%%%%%%%%%%%%%%%%%
%%%%%%%%%%%%%%%%%%%%%%%%%%%%%%%%%%%%%%

\begin{center}
$\S \ 11$
\\[0.5cm]
$\mathcal{REFERENCES}$\\[-.2cm]
\end{center}

\[
\textbf{BOOKS}
\]

\begingroup
\fontsize{9pt}{11pt}\selectfont
\setlength\parindent{0 cm}

[1] \quad Adams, J., \textit{Localization and Completion}, University of Chicago (1975).
\\[-.2cm]

[2] \quad Artin, M. and Mazur, B., \textit{Etale Homotopy}, Springer Verlag (1969).
\\[-.2cm]

[3] \quad Bousfield, A. and Kan, D., \textit{Homotopy Limits, Completions and Localizations}, Springer Verlag (1987).
\\[-.2cm]

[4] \quad Schwartz, L., \textit{Unstable Modules over the Steenrod Algebra and Sullivan's Fixed Point Set Conjecture}, 

\hspace{0.8cm}University of Chicago (1994).
\\[-.2cm]

[5] \quad Sullivan, D., \textit{Geometric Topology}, MIT (1970).
\\[-.2cm]
\endgroup

\[
\textbf{ARTICLES}
\]

\begingroup
\fontsize{9pt}{11pt}\selectfont
\setlength\parindent{0 cm}

[1] \quad Greenlees, J. and May, J., Derived Functors of $I$-adic Completion and Local Homology, 
\textit{J. Algebra} 

\hspace{0.8cm}\textbf{149} (1992), 438-453.
\\[-.2cm]

[2] \quad Greenlees, J. and May, J., Completions in Algebra and Topology, In: 
\textit{Handbook of Algebraic Topology}, 

\hspace{0.8cm}I. James (ed.), North Holland (1995), 255-276.
\\[-.2cm]

[3] \quad Lannes, J., Sur les Espaces Fonctionnels dont la Source est le Classifiant d'un $p$-Groupe Ab\'elien 

\hspace{0.8cm}\'El\'ementaire, \textit{Publ. Math. I.H.E.S.} \textbf{75} (1992), 135-244.
\\[-.2cm]

[4] \quad Morel, F., Quelques Remarques sur la Cohomologie Modulo $p$ Continue des Pro-$p$-Espaces et les 

\hspace{0.8cm}R\'esultats de J. Lannes Concernant les Espaces Fonctionnels hom$(BV,X)$, 
\textit{Ann. Sci. \'Ecole Norm.}

\hspace{0.8cm}\textit{Sup.} \textbf{26} (1993), 309-360.
\\[-.2cm]

[5] \quad Quillen, D., An Application of Simplicial Profinite Groups, 
\textit{Comment. Math. Helv.} \textbf{44} (1969), 45-60.
\\[-.2cm]

[6] \quad Rector, D., Homotopy Theory of Rigid Profinite Spaces, 
\textit{Pacific J. Math.} \textbf{85} (1979), 413-445.
\\[-.2cm]

[7] \quad Sullivan, D., Genetics of Homotopy Theory and the Adams Conjecture, 
\textit{Ann. of Math.} \textbf{100} (1974), 

\hspace{0.8cm}1-79.

\endgroup

\chapter{
$\boldsymbol{\S}$\textbf{12}.\quadx  MODEL CATEGORIES}
\setcounter{proposition}{0}
\setlength\parindent{2em}

%%----------------------------------------------------------------------------------------------01
$\text{ }$\\[-1.25cm]

Of the various proposals that have been advanced for the development of abstract homotopy theory, perhaps the most widely used and successful axiomization is Quillen's.  
The resulting unification is striking and the underlying techniques are applicable not only in topology but also in algebra.

Let $i:A \ra Y$, $p:X \ra B$ be morphisms in a category \bC $-$then $i$ is said to have the 
\un{left lifting property with respect to $p$}
\index{left lifting property with respect to $p$} 
(LLP w.r.t. $p$)
\index{LLP} 
and $p$ is said to have the 
\un{right lifting property with respect to $i$}
\index{right lifting property with respect to $i$} 
(RLP w.r.t. $i$)
\index{RLP} 
if for all $u: A \ra X$,  $v: Y \ra B$ such that $p \circ u = v \circ i$, there is a $w: Y \ra X$ such that $w \circ i = u$, $p \circ w = v$.\\

\begingroup%%----------------------------------->>
\fontsize{9pt}{11pt}\selectfont
For instance, take \bC = \bTOP $-$then $i:A \ra Y$ is a cofibration iff $\forall \ X$, $i$ has the LLP w.r.t. $p_0:PX \ra X$, i.e., 
\begin{tikzcd}[sep=large]
A  \ar{d} \ar{r} &PX \ar{d}{p_0}\\
Y \ar[dashed]{ur}\ar{r} &X
\end{tikzcd}
and $p:X \ra B$ is a Hurewicz fibration iff $\forall \ Y$, $p$ has the RLP w.r.t. $i_0:Y \ra IY$, i.e., 
\begin{tikzcd}[sep=large]
Y  \ar{d}[swap]{i_0} \ar{r} &X \ar{d}{p}\\
IY \ar[dashed]{ur} \ar{r} &B
\end{tikzcd}
.\\
\endgroup%%------------------------------------<<

\label{9.105}
Consider a category \bC equipped with three composition closed classes of morphisms termed 
\un{weak equivalences} 
\index{weak equivalences (model category)} 
(denoted $\overset{\sim}{\ra}$), 
\un{cofibrations} 
\index{cofibrations (model category)} 
(denoted $\rightarrowtail$), 
and 
\un{fibrations} 
\index{fibrations (model category)} 
(denoted $\twoheadrightarrow$), 
each containing the isomorphisms of \bC.  
Agreeing to call a morphism which is both a weak equivalence and a cofibration (fibration) an 
\un{acyclic cofibration}
\index{acyclic cofibrations (model category)}
(\un{fibration})
\index{acyclic fibrations (model category)}
\bC is said to be a 
\un{model category} 
\index{model category}
provided that the following axioms are satisfied.

\qquad (MC$-$1) \quadx 
\bC is finitely complete and finitely cocomplete.

\qquad (MC$-$2) \quadx 
Given composable morphisms $f$, $g$, if any two of $f$, $g$, $g \circ f$ are weak equivalences, so is the third.

\qquad(MC$-$3) \quadx 
Every retract of a weak equivalence, cofibration, or fibration is again a weak equivalence, cofibration, or fibration.
\label{15.36}

[Note: \ To say that $f:X \ra Y$ is a 
\un{retract}
\index{retract (model category)} 
of $g:W \ra Z$ means that there exist morphisms 
$i:X \ra W$, 
$r:W \ra X$, 
$j:Y \ra Z$, 
$s:Z \ra Y$ with 
$g \circ i = j \circ f$, 
$f \circ r = s \circ g$, 
$r \circ i = \id_X$, 
$s \circ j = \id_Y$.  
A retract of an isomorphism is an isomorphism.]

\qquad (MC$-$4) \quadx 
Every cofibration has the LLP w.r.t. every acyclic fibration and every fibration has the RLP w.r.t. every acyclic cofibration.

\qquad (MC$-$5) \quadx 
Every morphism can be written as the composite of a cofibration and an acyclic fibration and the composite of an acyclic cofibration and a fibration.

[Note: \  In proofs, the axioms for a model category are often used without citation.]

%%----------------------------------------------------------------------------------------------02

Remark: A weak equivalence which is a cofibration and a fibration is an isomorphism.

A model category \bC has an initial object (denoted $\emptyset$) and a final object (denoted $*$).  
An object $X$ in \bC is said to be 
\un{cofibrant}
\index{cofibrant} 
if $\emptyset \ra X$ is a cofibration and 
\un{fibrant}
\index{fibrant} 
if $X \ra *$ is a fibration.\\

\label{12.41}

\begingroup%%----------------------------------->>
\fontsize{9pt}{11pt}\selectfont
\textbf{\small FACT} \ 
Suppose that \bC is a model category.  Let $X \in \Ob\bC$ $-$then $X$ is cofibrant iff every acyclic fibration $Y \ra X$ has a right inverse and $X$ is fibrant iff every acyclic cofibration $X \ra Y$ has a left inverse.\\
\endgroup%%------------------------------------<<

Example:  Take \bC  = \bTOP  $-$then \bTOP  is a model category if weak equivalence = homotopy equivalence, cofibration = closed cofibration, fibration = Hurewicz fibration.  All objects are cofibrant and fibrant.

[MC$-$1 is clear, as is MC$-$2.  That MC$-$4 obtains is implied by what can be found on 
p. \pageref{12.1} $\&$ 
p. \pageref{12.2}, 
p. \pageref{12.3} $\&$ 
p. \pageref{12.4} 
and that MC$-$5 obtains is implied by what can be found on 
p. \pageref{12.5}.
There remains the verification of MC$-$3.  That MC$-$3 obtains for closed cofibrations or Hurewicz fibrations is implied by what can be found on 
p. \pageref{12.6} $\&$ 
p. \pageref{12.7}, 
p. \pageref{12.8} $\&$ 
p. \pageref{12.9}.  
Finally, suppose that $f$ is the retract of a homotopy equivalence $-$then $\abs{f}$ is the retract of an isomorphism in \bHTOP, so $\abs{f}$ is an isomorphism in \bHTOP, i.e., $f$ is a homotopy equivalence.]

[Note: \  We shall refer to this structure of a model category on \bTOP as the 
\un{standard} \un{structure}.]
\index{standard structure (model category on \bTOP)}

\label{14.64}
Addendum: \bCG has a standard model category structure, viz. weak equivalence = homotopy equivalence, cofibration = closed cofibration, fibration = \bCG fibration.

[The verification of MC-4 for \bCG is essentially the same as it is for \bTOP.  To check MC-5, note that $k$ preserves homotopy equivalences, sends closed cofibrations to closed cofibrations 
(cf. p. \pageref{12.10}), 
and takes Hurewicz fibrations to \bCG fibrations 
(cf. p. \pageref{12.11}).  
Therefore, if
$
\begin{cases}
\ X\\
\ Y
\end{cases}
$
are in \bCG and if $f:X \ra Y$ is a continuous function, one can first factor $f$ in \bTOP and then apply $k$ to get the desired factorization of $f$ in \bCG.]\\

\begingroup%%----------------------------------->>
\fontsize{9pt}{11pt}\selectfont
\textbf{\small  EXAMPLE} \ 
Let \bA be an abelian category.  
Write \textbf{CXA}
\index{\textbf{CXA}} 
for the abelian category of chain complexes over \bA.  
Given a morphism $f:X \ra Y$ in \textbf{CXA}, call $f$ a weak equivalence if $f$ is a chain homotopy equivalence, a cofibration if $\forall \ n$, $f_n:X_n \ra Y_n$ has a left inverse, and a fibration if $\forall \ n$, $f_n:X_n \ra Y_n$ has a right inverse $-$then \textbf{CXA} is a model category.  Every object is cofibrant and fibrant.\\
\endgroup%%------------------------------------<<

\begingroup%%----------------------------------->>
\fontsize{9pt}{11pt}\selectfont
\textbf{\small  EXAMPLE} \ 
Let \bA be an abelian category with enough projectives.  Write 
\textbf{CXA}$_{\geq 0}$
\index{\textbf{CXA}$_{\geq 0}$} 
for the full subcategory of \textbf{CXA} whose objects $X$ have the property that $X_n = 0$ if $n < 0$.  
Given a morphism $f:X \ra Y$ in \textbf{CXA}$_{\geq 0}$, call $f$ a weak equivalence if $f$ is a homology equivalence, a cofibration if $\forall \ n$, $f_n:X_n \ra Y_n$ is a monomorphism with a projective cokernel, and a fibration if $\forall \ n > 0$, $f_n:X_n \ra Y_n$ is
%%----------------------------------------------------------------------------------------------03
an epimorphism $-$then \textbf{CXA}$_{\geq 0}$ is a model category.  
Every object is fibrant and the cofibrant objects are those $X$ such that $\forall \ n$, $X_n$ is projective.\\
\endgroup%%------------------------------------<<

\begingroup%%----------------------------------->>
\fontsize{9pt}{11pt}\selectfont
There are lots of other ``algebraic'' examples of model categories, many of which figure prominently in rational homotopy theory (specifics can be found in the references at the end of the $\S$).\\
\endgroup%%------------------------------------<<

\label{13.91}
Given a model category \bC, $\bC^\OP$ acquires the structure of a model category by stipulating 
that $f^\OP$ is a weak equivalence in $\bC^\OP$ iff $f$ is a weak equivalence in \bC, 
that $f^\OP$ is a cofibration in $\bC^\OP$,  iff $f$ is a fibration in \bC,  and 
that $f^\OP$ is a fibration in $\bC^\OP$ ,  iff $f$ is a cofibration in \bC.

\label{14.123}
\label{16.37}
Given a model category \bC and objects $A$, $B$ in \bC, the categories $A\backslash \bC$, \bC/\mB are again model categories, a morphism in either case being declared a weak equivalence, cofibration, or fibration if it is such when viewed in \bC alone.
\label{12.30}

Example: Take \bC  = \bTOP  (standard structure) $-$then an object $(X,x_0)$ in $\bTOP_*$ is cofibrant iff $* \ra (X,x_0)$ is a closed cofibration (in \bTOP), i.e., iff $(X,x_0)$ is wellpointed with $\{x_0\} \subset X$ closed.\\

\label{12.32}
\begin{proposition} \ 
Let \bC be a model category.\\
\indent\indent (1) \ \  The cofibrations in \bC are the morphisms that have the LLP w.r.t. acyclic fibrations.\\
\indent\indent (2) \ \  The acyclic cofibrations in \bC are the morphisms that have the LLP w.r.t. fibrations.\\
\indent\indent (3) \ \  The fibrations in \bC are the morphisms that have the RLP w.r.t. acyclic cofibrations.\\ 
\indent\indent (4) \ \  The acyclic fibrations in \bC are the morphisms that have the RLP w.r.t. cofibrations.
\end{proposition}

[Statements (3) and (4) follow from statements (1) and (2) by duality.  
The proofs of (1) and (2) being analogous, consider (1).  
Thus suppose that $i:A \ra Y$ has the LLP w.r.t. acyclic fibrations.  Using MC-5, write $i = p \circ j$, where $j:A \ra X$ is a cofibration and $p:X \ra Y$ is an acyclic fibration.  
By hypothesis, $\exists$ a $w$ such that $w \circ i = j$, $p \circ w = \id_Y$, and this implies that $i$ is a retract of $j$, so $i$ is a cofibration.]\\

Example: Take \bC  = \bCG  (standard structure) $-$then an arrow $A \ra Y$ that has the LLP w.r.t. acyclic \bCG fibrations must be a closed cofibration.\\

\label{12.36}
\label{12.39}
\label{13.103}
\label{13.104}
\label{16.29}
\label{16.30}
\label{16.38}
\label{16.42}
\label{16.48}
\begingroup%%----------------------------------->>
\fontsize{9pt}{11pt}\selectfont
\textbf{\small  EXAMPLE} \ 
Let \bC and \bD be  model categories.  Suppose that
$
\begin{cases}
\ F: \bC \ra \bD\\
\ G: \bD \ra \bC
\end{cases}
$
are functors and $(F,G)$ is an adjoint pair $-$then $F$ preserves cofibrations and acyclic cofibrations iff $G$ preserves fibrations and acyclic fibrations.
\vspi
%%----------------------------------------------------------------------------------------------04
[Note: \  Either condition is equivalent to requiring that $F$ preserves cofibrations and $G$ preserves fibrations.]\\
\endgroup%%------------------------------------<<

\label{12.35} %dmc mnft

In a model category \bC, the classes of cofibrations and fibrations possess a number of ``closure'' properties (all verifications are simple consequences of Proposition 1).

\index{coproducts}
\indent\indent (Coproducts) 
If $\forall \ i$, $f_i:X_i \ra Y_i$ is a cofibration (acyclic cofibration), then 
$\ds\coprod\limits_i f_i: \coprod\limits_i X_i \ra \coprod\limits_i Y_i$ is a cofibration (acyclic cofibration).

\index{products}
\indent\indent (Products)  
If $\forall \ i$, $f_i:X_i \ra Y_i$ is a fibration (acyclic fibration), then 
$\ds\prod\limits_i f_i: \prod\limits_i X_i \ra \prod\limits_i Y_i$ is a cofibration (acyclic fibration).

\index{pushouts}
\indent\indent (Pushouts)  
Given a 2-source
\begin{tikzcd}[sep=small]
X  &Z \ar{l}[swap]{f}  \ar{r}{g} &Y,
\end{tikzcd}
define \mP by the pushout square
\begin{tikzcd}[sep=large]
Z  \ar{d}[swap]{f} \ar{r}{g} &Y \ar{d}{\eta}\\
X \ar{r}[swap]{\xi} &P
\end{tikzcd}
.  
Assume: $f$ is a cofibration (acyclic cofibration) $-$then $\eta$ is a cofibration (acyclic cofibration).

\index{pullbacks}
\indent\indent (Pullbacks)  
Given a 2-sink
\begin{tikzcd}[sep=small]
X  \ar{r}{f} &Z   &Y \ar{l}[swap]{g},
\end{tikzcd}
define \mP by the pullback square
\begin{tikzcd}[sep=large]
P  \ar{d}[swap]{\xi} \ar{r}{\eta} &Y \ar{d}{g}\\
X \ar{r}[swap]{f} &Z
\end{tikzcd}
.  
Assume: $g$ is a fibration (acyclic fibration) $-$then $\xi$ is a fibration (acyclic fibration).

\index{sequential colimits}
\indent\indent (Sequential Colimits)  
If $\forall \ n$, $f_n:X_n \ra X_{n+1}$  is a cofibration (acyclic cofibration), then $\forall \ n$, $i_n:X_n \ra \colimx X_n$ is a cofibration (acyclic cofibration).

\index{sequential limits}
\indent\indent (Sequential Limits)  
If $\forall \ n$, $f_n:X_{n+1} \ra X_n$  is a fibration (acyclic fibration), then $\forall \ n$, $p_n:\lim X_n \ra X_n$ is a fibration (acyclic fibration).

[Note: \  It is assumed that the relevant coproducts, products, sequential colimits, and sequential limits exist.]\\

\begingroup%%----------------------------------->>
\fontsize{9pt}{11pt}\selectfont
\index{pushouts (example)}
\textbf{\small EXAMPLE (\un{Pushouts})} \quadx
Fix a model category \bC.  Let \bI be the category
${1 \ \bullet} \overset{a}{\la} \underset{3}{\bullet} \overset{b}{\ra} {\bullet \ 2}$ 
(cf. p. \pageref{12.12}) 
$-$then the functor category $[\bI,\bC]$ is again a model category.  
Thus an object of $[\bI,\bC]$ is a 2-source
${X} \overset{f}{\la} {Z} \overset{g}{\ra} {Y}$
 and a morphism $\Xi$ of 2-sources is a commutative diagram
 \begin{tikzcd}[sep=large]
{X}  \ar{d}   &{Z} \ar{d} \ar{l}[swap]{f} \ar{r}{g}  &{Y} \ar{d} \\
{X^\prime}   &{Z^\prime}  \ar{l}{f^\prime} \ar{r}[swap]{g^\prime} &{Y^\prime} 
\end{tikzcd}
 .  
 Stipulate that $\Xi$ is a weak equivalence or a fibration if this is the case of each of its vertical constituents.  
 Define now $P_L$, $P_R$ by the pushout squares
 \begin{tikzcd}[sep=large]
{X}  \ar{d}   &{Z} \ar{d} \ar{l}[swap]{f}  \ar{d} \\
{P_L}   &{Z^\prime}  \ar{l}
\end{tikzcd}
,
\begin{tikzcd}[sep=large]
{Z} \ar{d}  \ar{r}{g}  &{Y} \ar{d} \\
{Z^\prime} \ar{r} &{P_R} 
\end{tikzcd}
, let $\rho_L:P_L \ra X^\prime$, $\rho_R : P_R \ra Y^\prime$ be the 
%%----------------------------------------------------------------------------------------------05
induced morphisms, and call $\Xi$ a cofibration provided that $Z \ra Z^\prime$, $\rho_L$, and $\rho_R$ are cofibrations.  
With these choices,  $[\bI,\bC]$ is a model category.  
The fibrant objects  
${X} \overset{f}{\lla} {Z} \overset{g}{\lra} {Y}$
in $[\bI,\bC]$ are those for which $X$, $Y$, and $Z$ are fibrant.  The cofibrant objects
${X} \overset{f}{\lla} {Z} \overset{g}{\lra} {Y}$
in $[\bI,\bC]$ are those for which $Z$ is cofibrant and 
$
\begin{cases}
\ f:Z \ra X\\
\ g:Z \ra Y
\end{cases}
$
are cofibrations.
\vspi
[Note: \  The story for pullbacks is analogous.]\\
\endgroup%%------------------------------------<<

\label{12.38}
\label{13.96}
\begingroup%%----------------------------------->>
\fontsize{9pt}{11pt}\selectfont
\textbf{\small  EXAMPLE} \ 
Fix  model category \bC $-$then $\bFIL(\bC)$ is again a model category.  
Thus let $\phi:(\bX,\bff) \ra (\bY,\bg)$ be a morphism in $\bFIL(\bC)$.  
Stipulate that $\phi$ is a weak equivalence or a fibration if this is the case of each $\phi_n$.  
Define now $P_{n+1}$ by the pushout square 
\begin{tikzcd}[sep=large]
{X_n}  \ar{d}[swap]{\phi_n} \ar{r}{f_n} &{X_{n+1}} \ar{d}\\
{Y_n} \ar{r}[swap]{g_n} &{P_{n+1}}
\end{tikzcd}
, 
let $\rho_{n+1}:P_{n+1} \ra Y_{n+1}$ be the induced morphism, and call $\phi$ a cofibration provided that $\phi_0$ and all the $\rho_{n+1}$ are cofibrations (each $\phi_n$ $(n > 0)$ is then a cofibration as well).  
With these choices, $\bFIL(\bC)$ is a model category.  
The fibrant objects $(\bX,\bff)$ in $\bFIL(\bC)$ are those for which $X_n$ is fibrant $\forall$ $n$.  
The cofibrant objects $(\bX,\bff)$ in $\bFIL(\bC)$ are those for which $X_0$ is cofibrant and 
$\forall \ n$, $f_n:X_n \ra X_{n+1}$ is a cofibration.
\vspi
[Note: \  The story for $\bTOW(\bC)$ is analogous.]\\
\endgroup%%------------------------------------<<

\begingroup%%----------------------------------->>
\fontsize{9pt}{11pt}\selectfont
\textbf{\small FACT} \ 
Let \bC be a model category.  Suppose that 
\begin{tikzcd}[sep=large]
{A}  \ar{d}[swap]{i} \ar{r}{u} &{X} \ar{d}{p}\\
{Y} \ar{r}[swap]{v} &{B}
\end{tikzcd}
is a comutative diagram in \bC, where $i$ is a cofibration, $p$ is a weak equivalence, 
and $X$ is fibrant $-$then $\exists$ a $w:Y \ra X$ such that $w \circ i = u$.
\vspi
[Note: \  There is a similar assertion for fibrations and cofibrant objects.]\\
\endgroup%%------------------------------------<<

Given a model category \bC, objects $X^{\prime}$ and $X^{\prime\prime}$ are said to be 
\un{weakly equivalent} 
\index{weakly equivalent} 
if there exists a path beginning at $X^{\prime}$ and ending at 
$X\pp$: 
$X^{\prime} = X_0 \ra X_1 \leftarrow \cdots \ra X_{2n-1}$ $\la X_{2n} = X\pp$,
where all the arrows are weak equivalences.  
Example: Take \bC = \bTOP (standard structure) 
$-$then $X^{\prime}$ and $X^{\prime\prime}$ are weakly equivalent iff they have the same homotopy type.\\

\begingroup%%----------------------------------->>
\fontsize{9pt}{11pt}\selectfont
\textbf{\small  EXAMPLE} \ 
The arrow category \bC($\rightarrow$) of a model category \bC is again a model category 
(cf. p. \pageref{12.13}).  
Therefore it makes sense to consider weakly equivalent morphisms.  Example: Every morphism in \bC is weakly equivalent to a fibration with a fibrant domain and codomain.\\
\endgroup%%------------------------------------<<

\textbf{\small COMPOSITION LEMMA} \ 
Consider the commutative diagram
\hspace{-.3cm}
\begin{tikzcd}[sep=large]
{\bullet}  \ar{d} \ar{r}  &{\bullet} \ar{d} \ar{r}  &{\bullet} \ar{d} \\
{\bullet}  \ar{r} &{\bullet}  \ar{r} &{\bullet} 
\end{tikzcd}
in a category \bC.  
Suppose that both the squares are pushouts$-$then the rectangle is a pushout.  
Conversely, if  the rectangle and the first square are pushouts, then the second square is a pushout.\\

%%----------------------------------------------------------------------------------------------06
Application:  Consider the commutative cube
\begin{tikzcd}[sep=small]
&{\bullet} \ar{dd} \ar{rr} &&{\bullet} \ar{dd} \\
  {\bullet} \ar{dd} \ar{rr}  \ar{ru} &&{\bullet} \ar{dd}\ar{ru} \\
&{\bullet}\ar{rr}   &&{\bullet} \\
  {\bullet} \ar{rr}  \ar{ru} &&{\bullet}  \ar{ru}
\end{tikzcd}
in a category \bC.  Suppose that the top and the left and right hand sides are pushouts $-$then the bottom is a pushout.\\

\begin{proposition} \ 
Let \bC be a model category.  Given a $2-$source
$X \overset{f}{\lla} Z \overset{g}{\lra} Y$, 
define $P$ by the pushout square
\begin{tikzcd}[sep=large]
Z  \ar{d}[swap]{f} \ar{r}{g} &Y \ar{d}{\eta}\\
X \ar{r}[swap]{\xi} &P
\end{tikzcd}
Assume: $f$ is a cofibration and $g$ is a weak equivalence $-$then $\xi$ is a weak equivalence provided that $Z$ $\&$ Y are cofibrant.
\end{proposition}

[Introduce the cylinder object $IZ$ for $Z$ 
(cf. p. \pageref{12.14}) 
and define $M_g$ by the pushout square
\begin{tikzcd}[sep=large]
{Z  \coprod Z} \ar{d}[swap]{\iota} \ar{r}{g \coprod \id_Z} &{Y  \coprod Z} \ar{d}\\
IZ \ar{r} &M_g
\end{tikzcd}
(cf. p. \pageref{12.15}).  
Noting that
\begin{tikzcd}[sep=small]
&{Z  \coprod Z} \ar{dl} \ar{r} &{Y  \coprod Z} \ar{dd}[swap]{(\id_Y,g)}\\
IZ \ar{rd}\\
&Z \ar{r}[swap]{g} &Y
\end{tikzcd}
commutes, choose $r:M_g \ra Y$ accordingly, so $g = r \circ i$ and $r \circ j = \id_Y$, where $i:Z \ra M_g$ is the composite $Z \ra Y \coprod Z \ra M_g$ and $j:Y \ra M_g$ is the composite $Y \ra Y\coprod Z \ra M_g$.  Since $\iota$ is a cofibration and
\begin{tikzcd}[sep=large]
{\emptyset}  \ar{d} \ar{r} &Z \ar{d}\\
Y \ar{r} &{Y  \coprod Z}
\end{tikzcd}
is a pushout square, $i$ and $j$ are cofibrations.  
Moreover, $j$ is acyclic.  
This is because $i_0:Z \ra IZ$ is an acyclic cofibration and $j$ is obtained from $i_0$ via
\begin{tikzcd}%[sep=large]
Z  \ar{d}[swap]{\text{in}_0} \ar{r}{g} &Y \ar{d}\\
{Z  \coprod Z}  \ar{d}[swap]{\iota} \ar{r} &{Y  \coprod Z} \ar{d}\\
IZ \ar{r} &M_g
\end{tikzcd}
$(i_0 = \iota \circ \text{in}_0)$.  
Therefore $r$ is a weak equivalence.  
But, by assumption, $g$ is a weak equivalence.  
Therefore $i$ is a weak equivalence.  
Define $\ov{I}$ by the pushout square
\begin{tikzcd}[sep=large]
Z  \ar{d}[swap]{f} \ar{r}{i} &M_g  \ar{d}{\ov{f}}\\
X \ar{r}[swap]{\ov{i}} &{\ov{I}}
\end{tikzcd}
Since $f$ is a cofibration and $i$ is an acyclic cofibration, $\ov{f}$ is a cofibration and $\bar{i}$ is an acyclic cofibration.  The commutative diagram
\begin{tikzcd}[sep=large]
M_g  \ar{d}[swap]{\ov{f}} \ar{r}{r} &Y \ar{d}{\eta}\\
{\ov{I}} \ar{r}[swap]{r_i} &P
\end{tikzcd}
is a pushout square and $\xi = r_i \circ \bar{i}$.  Define $\ov{J}$ by the pushout square
\begin{tikzcd}[sep=large]
Y  \ar{d}[swap]{\eta} \ar{r}{j} &M_g \ar{d}{\ov{\eta}}\\
P \ar{r}[swap]{\ov{j}} &{\ov{J}}
\end{tikzcd}
.  Since $\eta$ is a cofibration and $j$ is an acyclic cofibration, $\bar{\eta}$ is a cofibration and $\bar{j}$ is an 
%%----------------------------------------------------------------------------------------------07
acyclic cofibration.  
The commutative diagram
\begin{tikzcd}[sep=large]
M_g  \ar{d}[swap]{\ov{\bar{\eta}}} \ar{r}{r} &Y \ar{d}{\eta}\\
{\ov{J}} \ar{r}[swap]{r_j} &P
\end{tikzcd}
is a pushout square and $\id_P = r_j \circ \bar{j}$.  
Therefore $r_j$ is a weak equivalence.  Define $Z_0$, $Z_1$ by the pushout squares
\begin{tikzcd}[sep=large]
Z  \ar{d}[swap]{f} \ar{r}{i_0} &IZ \ar{d}{f_0}\\
X \ar{r}{\sim} &{Z_0}
\end{tikzcd}
,
\begin{tikzcd}[sep=large]
Z  \ar{d}[swap]{f} \ar{r}{i_1} &IZ \ar{d}{f_1}\\
X \ar{r}{\sim} &{Z_1}
\end{tikzcd}
.  The composites
\begin{tikzcd}[sep=small]
Z  \ar{r}{i_0} &IZ \ar{r}{\sim} &Z,
\end{tikzcd}
\begin{tikzcd}[sep=small]
Z  \ar{r}{i_1} &IZ \ar{r}{\sim} &Z
\end{tikzcd}
being $\id_Z$, there are weak equivalences $\zeta_0:Z_0 \ra X$, $\zeta_1:Z_1 \ra X$ and factorizations
\begin{tikzcd}[sep=small]
X \ar{r}{\sim} &Z_0 \ar{r}{\zeta_0} &X,
\end{tikzcd}
\begin{tikzcd}[sep=small]
X \ar{r}{\sim} &Z_1 \ar{r}{\zeta_1} &X 
\end{tikzcd}
of $\id_X$.  
Define $W$ by the pushout square
\begin{tikzcd}[sep=large]
IZ  \ar{d}[swap]{f_0} \ar{r}{f_1} &Z_1 \ar{d}\\
Z_0 \ar{r} &{W}
\end{tikzcd}
and determine $\zeta:W \ra X$ so that $\zeta_0$ is the composite
\begin{tikzcd}[sep=small]
Z_0 \ar{r} &W \ar{r}{\zeta} &X
\end{tikzcd}
and $\zeta_1$  is the composite
%$
%\begin{tikzcd}[sep=small]
%Z_1\ar{r} &W \ar{r}{\zeta} &X
%\end{tikzcd}
%.
%$ 
$Z_1 \lra W \overset{\zeta}{\lra} X$.
Decompose $\zeta$ per 
$W \rightarrowtail \ov{W} \overset{\sim}{\twoheadrightarrow} X$ 
%\begin{tikzcd}[sep=small]
%W \ar{r} &\ov{W} \ar{r}{\sim} &X
%\end{tikzcd}
$-$then the composites $Z_0 \ra \ov{W}$, $Z_1 \ra \ov{W}$ are acyclic cofibrations.  
To go from $Z$ to $\ov{I}$ through 
\begin{tikzcd}[sep=small]
Z \ar{r}{i_1} &IZ \ar{r}{\sim} &Z \ar{r}{i} &{M_g} \ar{r}{\bar{f}} &\ov{I}
\end{tikzcd}
is the same as going from $Z$ to $\ov{I}$ through 
$
\begin{tikzcd}[sep=small]
Z \ar{r}{f} &X \ar{r}{\bar{i}} &\ov{I}
\end{tikzcd}
.
$  
Consequently, there is an arrow $\bar{i}_1:Z_1 \ra \ov{I}$ such that the composite
\begin{tikzcd}[sep=small]
X \ar{r}{\sim} &Z_1 \ar{r}{\bar{i}_1} &\ov{I}
\end{tikzcd}
is $\bar{i}$ and the commutative diagram
\begin{tikzcd}%[sep=large]
IZ  \ar{d}[swap]{f_1} \ar{r} &M_g \ar{d}\\
Z_1 \ar{r}[swap]{\bar{i}_1} &{\ov{I}}
\end{tikzcd}
is a pushout square.  
But $\bar{i}$ is a weak equivalence.  
Therefore $\bar{i}_1$ is a weak equivalence.  Define $\ov{K}$ by the pushout square
$
\begin{tikzcd}%[sep=large]
Z_1  \ar{d} \ar{r}{\bar{i}_1} &\ov{I} \ar{d}\\
\ov{W} \ar{r} &{\ov{K}}
\end{tikzcd}
.
$
Since $\bar{i}_1$ is a weak equivalence, the same holds for $\ov{W} \ra \ov{K}$.  
To go from $Z$ to $\ov{J}$ through
\begin{tikzcd}[sep=small]
Z \ar{r}{i_0} &IZ \ar{r}{\sim} &Z \ar{r}{g} &Y \ar{r}{j} &M_g \ar{r}{\bar{\eta}} &\ov{J}
\end{tikzcd}
 is the same as going from $Z$ to $\ov{J}$ through
\begin{tikzcd}[sep=small]
Z \ar{r}{f} &X \ar{r}{\xi} &P  \ar{r}{\bar{j}} &\ov{J.}
\end{tikzcd}
%.  
Consequently, there is an arrow 
$\bar{j}_0:Z_0 \ra \ov{J}$ such that the composite
$X \overset{\sim}{\ra} Z_0 \overset{\ov{j}_0}{\lra} \ov{J}$ is $\ov{j} \circ \xi$ 
%\begin{tikzcd}[sep=small]
%X \ar{r}{\sim} &Z_0 \ar{r} &\ov{W}  \ar{r}{\sim} &\ov{K}
%\end{tikzcd}
and the commutative diagram
\begin{tikzcd}%[sep=large]
IZ  \ar{d}[swap]{f_0} \ar{r} &M_g \ar{d}\\
Z_0 \ar{r}[swap]{\bar{j}_0} &{\ov{J}}
\end{tikzcd}
is a pushout square.  
To go from $IZ$ to $\ov{K}$ by $IZ \ra M_g \ra \ov{I} \ra \ov{K}$ is the same as going from $IZ$ to $\ov{K}$ by
\begin{tikzcd}[sep=small]
IZ \ar{r} &Z_0 \ar{r} &\ov{W}  \ar{r}{\sim} &\ov{K},
\end{tikzcd}
thus there is an arrow $\ov{J} \ra \ov{K}$ and a commutative diagram
\begin{tikzcd}%[sep=large]
Z_0  \ar{d}[swap]{f_0} \ar{r}{\bar{j}_0} &\ov{J} \ar{d}\\
\ov{W} \ar{r} &\ov{K}
\end{tikzcd}
which is a pushout square:
%\begin{tikzcd}%[sep=small]
%& & &{M_g} \ar{d} \ar{r} &{\ov{I}} \ar{d} \\
%& & &{\ov{J}}  \ar{r} &{\ov{K}}\\
%{IZ}\ar{r}  \ar{rrruu}  \ar{d} &{Z_1}  \ar{rrruu} \ar{d} \\
%{Z_0} \ar{r}   \ar{rrruu} &{\ov{W}}  \ar{rrruu}
%\end{tikzcd}
\begin{tikzcd}%[sep=small]
&{M_g} \ar{d} \ar{r} &{\ov{I}} \ar{d} \\
&{\ov{J}}  \ar{r} &{\ov{K}}\\
{IZ}\ar{r}  \ar{ruu}  \ar{d} &{Z_1}  \ar{ruu} \ar{d} \\
{Z_0} \ar{r}   \ar{ruu} &{\ov{W}}  \ar{ruu}
\end{tikzcd}
.
\\
%%----------------------------------------------------------------------------------------------08
It follows that $\bar{j}_0$ is a weak equivalence and this implies that $\bar{j} \circ \xi$ is a weak equivalence.  
Finally, $\xi = \id_P \circ \xi = r_j \circ \bar{j} \circ \xi$ is a weak equivalence.  ]

[Note: \  There is a parallel statement for fibrations and pullbacks.]\\

\begingroup%%----------------------------------->>
\fontsize{9pt}{11pt}\selectfont
\textbf{\small  EXAMPLE} \ 
Working in \bC = \bTOP (standard structure), suppose that $A \ra X$ is a closed cofibration.  Let $f:A \ra Y$ be a homotopy equivalence $-$then the arrow $X \ra X \sqcup_f Y$ is a homotopy equivalence 
(cf. p. \pageref{12.16}).\\
\endgroup%%------------------------------------<<

\begin{proposition} \ 
Let \bC be a model category.  Suppose given a commutative diagram
\begin{tikzcd}[sep=large]
{X}  \ar{d}   &Z \ar{d}  \ar{l}[swap]{f}             \ar{r}{g}             &{Y} \ar{d}\\
{X^\prime}  &{Z^\prime} \ar{l}{f^\prime} \ar{r}[swap]{g^\prime} &{Y^\prime}
\end{tikzcd}
, where 
$
\begin{cases}
\ f\\
\ f^\prime
\end{cases}
$
are cofibrations and the vertical arrows are weak equivalences $-$then the induced morphism $P \ra P^\prime$ of pushouts is a weak equivalence provided that
$
\begin{cases}
\ Z \ \& \ Y\\
\ Z^\prime \ \& \ Y^\prime
\end{cases}
$
are cofibrant.
\end{proposition}

[We shall first treat the special case when $g$ is a cofibration.  In this situation, the arrow
$Y \ra Z^\prime \underset{Z}{\sqcup} Y$ is a weak equivalence (cf. Proposition 2) and $Z^\prime \underset{Z}{\sqcup} Y$ is cofibrant.  Form the pushout square
\begin{tikzcd}[sep=large]
{Y}  \ar{d} \ar{r}              &{Z^\prime \underset{Z}{\sqcup} Y}  \ar{d}\\
{X \underset{Z}{\sqcup} Y} \ar{r}  &{X \underset{Z}{\sqcup} (Z^\prime \underset{Z}{\sqcup} Y)}
\end{tikzcd}
and apply Proposition 2 once again to see that the arrow
$X \underset{Z}{\sqcup} Y \ra {X \underset{Z}{\sqcup} (Z^\prime \underset{Z}{\sqcup} Y)}$ is a weak equivalence.
Next write
${X \underset{Z}{\sqcup} (Z^\prime \underset{Z}{\sqcup} Y)}$ $\approx$ $(X \underset{Z}{\sqcup} Z^\prime )
\underset{Z^\prime}{\sqcup}  (Z^\prime \underset{Z}{\sqcup} Y)$
and note that the arrow 
 $X \underset{Z}{\sqcup} Z^\prime  \ra X^\prime$ is a weak equivalence (cf. Proposition 2).  
Consider now the commutative diagram
\[
\begin{tikzcd}[sep=large]
{Z^\prime}  \ar{d} \ar{r}        
&{X \underset{Z}{\sqcup} Z^\prime} \ar{d} \ar{r}                                             
&{X^\prime} \ar{d}\\
{Z^\prime \underset{Z}{\sqcup} Y} \ar{r}  
&{(X \underset{Z}{\sqcup} Z^\prime) \underset{Z^\prime}{\sqcup} (Z^\prime \underset{Z}{\sqcup} Y)} \ar{r} 
&{X^\prime \underset{{Z^\prime}}{\sqcup} (Z^\prime \underset{Z}{\sqcup} Y)}
\end{tikzcd}
\]
in which both the squares and the rectangle are pushouts.  
Since 
$Z^\prime \rightarrowtail Z^\prime \underset{Z}{\sqcup} Y$ 
$\implies$
$X \underset{Z}{\sqcup} Z^\prime \rightarrowtail 
(X \underset{Z}{\sqcup} Z^\prime) \underset{Z^\prime}{\sqcup} (Z^\prime \underset{Z}{\sqcup} Y)$
 \ and \ 
$X \underset{Z}{\sqcup} Z^\prime $ 
is cofibrant, still another application of Proposition 2 implies that the arrow \ 
$(X \underset{Z}{\sqcup} Z^\prime) \underset{Z^\prime}{\sqcup} (Z^\prime \underset{Z}{\sqcup} Y)
\ra X^\prime \underset{Z^\prime}{\sqcup} (Z^\prime \underset{Z}{\sqcup} Y)$ \ 
is a weak equivalence.  
Repeating the reasoning with
\[
\begin{tikzcd}[sep=large]
{Z^\prime}  \ar{d} \ar{r}        
&{Z^\prime \underset{Z}{\sqcup} Y} \ar{d} \ar{r}                                             
&{Y^\prime} \ar{d}\\
{X^\prime}  \ar{r}  
&{X^\prime \underset{Z\prime}{\sqcup} (Z^\prime \underset{Z}{\sqcup} Y)} \ar{r} 
&{X^\prime \underset{{Z^\prime}}{\sqcup} Y^\prime}
\end{tikzcd}
\]
%%----------------------------------------------------------------------------------------------09
leads to the conclusion that the arrow
$X^\prime \underset{Z^\prime}{\sqcup} (Z^\prime \underset{Z}{\sqcup} Y) \rightarrow
X^\prime \underset{{Z^\prime}}{\sqcup} Y^\prime$
is a weak equivalence.  
We have therefore built a weak equivalence from \mP to $P^\prime$.  
To proceed in general, factor $g$ as 
$Z \rightarrowtail \ov{Y} \overset{\sim}{\twoheadrightarrow} Y$.  
%\begin{tikzcd}[sep=small]
%Z  \ar{r}{\sim} &\ov{Y} \ar{r}{\sim} \ar[r, two heads] &{Y.}
%\end{tikzcd}
Define $\ov{X}^\prime$, $\ov{Y}^\prime$ by the pushout squares
\begin{tikzcd}[sep=large]
Z  \ar{d} \ar{r} &X \ar{d}\\
Z^\prime \ar{r} &{\ov{X}^\prime }
\end{tikzcd}
, 
\begin{tikzcd}[sep=large]
Z  \ar{d} \ar{r} &\ov{Y} \ar{d}\\
Z^\prime \ar{r} &\ov{Y}^\prime
\end{tikzcd}
$-$then there are weak equivalences 
$\ov{X^\prime} \ra X^\prime$, 
$\ov{Y^\prime} \ra Y^\prime$.
The 2-sources 
$X \leftarrow Z \ra \ov{Y}$,  
$\ov{X}^\prime \leftarrow Z^\prime \ra \ov{Y}^\prime$
generate pushouts $\ov{P}$, $\ov{P^\prime}$.
Since the arrows on the ``right'' are cofibrations, the induced morphisms 
$\ov{P} \ra P$, 
$\ov{P} \ra \ov{P^\prime}$, 
$\ov{P^\prime} \ra P^\prime$, 
are weak equivalences.
The assertion thus follows from the fact that the diagram
\begin{tikzcd}[sep=large]
\ov{P}  \ar{d} \ar{r} &P \ar{d}\\
\ov{P}^\prime \ar{r} &P^\prime
\end{tikzcd}
commutes.]

[Note: \  There is a parallel statement for fibrations and pullbacks.]\\

\begingroup%%----------------------------------->>
\fontsize{9pt}{11pt}\selectfont
\textbf{\small  EXAMPLE} \ 
Working in \bC = \bTOP (standard structure), suppose that 
$
\begin{cases}
\ A \ra X\\
\ A^\prime \ra X^\prime
\end{cases}
$
are closed cofibrations.  Let 
$
\begin{cases}
\ f:A \ra Y\\
\ f^\prime: A^\prime \ra Y^\prime
\end{cases}
$
be continuous functions.    Assume that the diagram 
\begin{tikzcd}%[sep=large]
{X}  \ar{d}  &{A} \ar{d} \ar{l} \ar{r}{f} &{Y} \ar{d}\\
{X^\prime}  &{A^\prime} \ar{l} \ar{r}[swap]{f^\prime} &{Y^\prime}
\end{tikzcd}
commutes and that the vertical arrows are homotopy equivalences $-$then the induced map
$X \sqcup_f Y \ra X^\prime \sqcup_{f^\prime} Y^\prime$
is a homotopy equivalence. 
(cf. p. \pageref{12.17} ff.).\\ 
\endgroup%%------------------------------------<<

\begin{proposition} \ %04
Let \bC be a model category.  Suppose given a commutative diagram
\begin{tikzcd}%[sep=large]
{X}  \ar{d}  &{Z} \ar{d} \ar{l}[swap]{f} \ar{r}{g} &{Y} \ar{d}\\
{X^\prime}  &{Z^\prime} \ar{l}{f^\prime} \ar{r}[swap]{g^\prime} &{Y^\prime}
\end{tikzcd}
, where $Y \ra Y^\prime$ and  
${X \underset{{Z}}{\sqcup} Z^\prime} \ra X^\prime$ are cofibrations (acyclic cofibrations) 
$-$then the induced morphism  $P \ra P^\prime$ of pushouts is a  cofibration (acyclic cofibration).
\end{proposition}

[Each morphism in the string
$P = X \underset{Z}{\sqcup} Y \ra {X \underset{{Z}}{\sqcup} Y^\prime}
\approx ({X \underset{{Z}}{\sqcup} Z^\prime}) \underset{Z^\prime}{\sqcup} Y^\prime
\ra  X^\prime \underset{Z^\prime}{\sqcup} Y^\prime
= P^\prime$
is a  cofibration (acyclic cofibration).]

[Note: \ There is a parallel statement for fibrations and pullbacks.]\\

\begingroup%%----------------------------------->>
\fontsize{9pt}{11pt}\selectfont
In the topological setting, Proposition 4 is related to but does not directly imply the lemma on 
p. \pageref{12.17a} ff.\\
\endgroup%%------------------------------------<<

\indent\indent (Small Object Argument) \  
Suppose that \bC is a cocomplete category.  
Let $S_0 = \{L_i \overset{\phi_i}{\ra} K_i \ (i \in I)\}$ be a set of morphisms in \bC.  
Given a morphism $f:X \ra Y$, consider the set of pairs of morphisms $(g,h)$ such that the diagram
\begin{tikzcd}[sep=large]
{L_i}  \ar{d}[swap]{\phi_i} \ar{r}{g}  &{X} \ar{d}{f}\\
{K_i}  \ar{r}[swap]{h} &{Y}
\end{tikzcd}
commutes.  Put
%%----------------------------------------------------------------------------------------------10
$X_0 = X$ and define $X_1$ by the pushout square
\begin{tikzcd}%[sep=large]
{\coprod\limits_i \coprod\limits_{(g,h)} L_i}  \ar{d} \ar{r}  &{X_0} \ar{d}\\
{\coprod\limits_i \coprod\limits_{(g,h)} K_i}  \ar{r} &{X_1}
\end{tikzcd}
.  
Observing that the data furnishes a commutative triangle
\begin{tikzcd}[sep=small]
{X_0}  \ar{ddr} \ar{rr} &&{X_1} \ar{ddl}\\
\\
&{Y}
\end{tikzcd}
, one may proceed and construct a sequence $X = X_0 \ra X_1 \ra \cdots \ra X_\omega$ of objects in \bC, taking $X_\omega = \colimx X_n$.  There is a commutative triangle
\begin{tikzcd}[sep=small]
{X}  \ar{ddr}[swap]{f} \ar{rr}{i_\omega}  &&{X_\omega} \ar{ddl}{f_\omega}\\
\\
&{Y}
\end{tikzcd}
and if $\forall \ i$, $L_i$ is $\omega$-definite, then the conclusion is that $f_\omega:X_\omega \ra Y$ has the RLP w.r.t. each $\phi_i$.

[Note: \  All that's really required of the $L_i$ is that the arrow $\colimx \Mor(L_i,X_n) \ra \Mor(L_i,X_\omega)$ be surjective $\forall \ i$.]

\label{16.11}
Example: Take $\bC  = \bTOP$ $-$then  \bTOP is a model category if weak equivalence = weak homotopy equivalence, fibration = Serre fibration, cofibration = all continuous functions which have the LLP w.r.t Serre fibrations that are weak homotopy equivalences.  Every object is fibrant and every CW complex is cofibrant.  Every object is weakly equivalent to a CW complex.

[Axioms MC-1, MC-2 and MC-3 are immediate.

Claim: Every continuous function $f:X \ra Y$ can be written as a composite $f_\omega \circ i_\omega$, where $i_\omega:X \ra X_\omega$ is a weak homotopy equivalence and has the LLP w.r.t Serre fibrations and $f_\omega:X_\omega \ra Y$ is a Serre fibration.

[Serre fibrations can be characterized by the property that they have the RLP w.r.t the embeddings
$i_0:[0,1]^n \ra I[0,1]^n$ $(n \geq 0)$ 
(cf. p. \pageref{12.18}).  
Accordingly, in the small object argument, take 
$S_0 = \{[0,1]^n \overset{i_0}{\lra} I[0,1]^n$ $(n \geq 0)\}$ $-$then $\forall \ k$, the arrow $X_k \ra X_{k+1}$ is a homotopy equivalence and has the LLP w.r.t. Serre fibrations.  
Consider the factorization of $f$ arising from the small object argument
\begin{tikzcd}[sep=small]
{X}  \ar{ddr}[swap]{f} \ar{rr}{i_\omega}  &&{X_\omega} \ar{ddl}{f_\omega}\\
\\
&{Y}
\end{tikzcd}
.  
It is clear that $i_\omega$ has the LLP w.r.t. Serre fibrations.  
On the other hand, since the points of $X_\omega - i_\omega(X)$ are closed,  every compact subset of $X_\omega$ lies in some $X_k$, thus the arrow
$\colimx C([0,1]^n,X_k) \ra C([0,1]^n,X_\omega)$ is surjective $\forall \ n$.  Therefore $f_\omega$ has the RLP w.r.t each $i_0:[0,1]^n \ra I[0,1]^n$, hence is a Serre fibration.  And: $i_\omega$ is a homotopy equivalence (cf. $\S 3$, Proposition 15), hence is a weak homotopy equivalence.]

%%----------------------------------------------------------------------------------------------11
Claim:  Every continuous function $f:X \ra Y$ can be written as a composite $f_\omega \circ i_\omega$, where $i_\omega:X \ra X_\omega$ has the LLP w.r.t. Serre fibrations that are weak homotopy equivalences and $f_\omega$ is both a weak homotopy equivalence and a Serre fibration.

[Serre fibrations that are weak homotopy equivalences can be characterized by the property that they have the RLP w.r.t. the inclusions
$\textbf{S}^{n-1} \ra \bD^{n}$ $(n \geq 0)$ 
(cf. p. \pageref{12.19}).  
Accordingly, in the small object argument, take
$S_0 = \{\textbf{S}^{n-1} \ra \bD^{n} \  (n \geq 0)\}$ and reason as above.]

Combining the claims gives MC-5.  Turning to the nontrivial half of MC-4, viz. that ``every fibration has the RLP w.r.t. every acyclic cofibration'', suppose that $f:X \ra Y$ is an ``acyclic cofibration''.  Decompose $f$ per the first claim: $f = f_\omega \circ i_\omega$.  Since $f$ and $i_\omega$ are weak homotopy equivalences, the same is true of $f_\omega$, so $\exists$ a $g:Y \ra X_\omega$ such that $g \circ f = i_\omega$, $f_\omega \circ g = \id_Y$.  This means that $f$ is a retract of $i_\omega$.  
But the class of maps which have the LLP w.r.t Serre fibrations is closed under the formation of retracts.]

[Note: \  We shall refer to this structure of a model category on \bTOP as the 
\un{singular} \un{structure}.]
\index{singular structure (model category)}

Remark: If $(K,L)$ is a relative CW complex, then the inclusion $L \ra K$ has the LLP w.r.t Serre fibrations that are weak homotopy equivalences 
(cf. p. \pageref{12.20}), 
hence is a cofibration in the singular structure.

[Note: \  Every cofibration in the singular structure is a cofibration in the standard structure, thus is a closed cofibration.  In fact there is a characterization: A continuous function is a cofibration in the singular structure iff it is a retract of a ``countable composition'' 
$X_0 \ra X_1 \ra \cdots \ra X_\omega$, where $\forall \ k$ the arrow $X_k \ra X_{k+1}$ is defined by the pushout square
\begin{tikzcd}%[sep=large]
{\coprod\limits_{n \geq 0} \coprod \textbf{S}^{n-1}}  \ar{d} \ar{r}  &{X_k} \ar{d}\\
{\coprod\limits_{n \geq 0} \coprod \bD^{n}}  \ar{r} &{X_{k+1}}
\end{tikzcd}
.]

\label{13.69}
Addendum: \bCG, $\dcg$, and \bCGH have a singular model category structure, viz. 
weak equivalence = weak homotopy equivalence, 
fibration = Serre fibration, 
cofibration = all continuous functions which have the LLP w.r.t Serre fibrations that are weak homotopy equivalences.

[In fact, if $f:X \ra Y$ is a continuous function, where
$
\begin{cases}
\ X\\
\ Y
\end{cases}
$
are in \bCG, $\dcg$, or \bCGH, then the $X_\omega$ that figures in either of the small object arguments used above is again in \bCG, $\dcg$, and \bCGH.]\\

\begingroup%%----------------------------------->>
\fontsize{9pt}{11pt}\selectfont
\textbf{\small  EXAMPLE} \ 
Take \bC = \bTOP (singular structure) $-$then any cofibrant $X$ is a CW space.  Thus fix a CW resolution $f:K \ra X$.  Factor $f$ as 
$K \underset{i}{\rightarrowtail} L \underset{p}{\overset{\sim}{\twoheadrightarrow}} X,$
%\begin{tikzcd}[sep=large]
%{K}  \ar{r}{f} \ar{r}  &{L} \ar[two heads]{r}{\sim} \ar[two heads]{r}[swap]{p}  &{X}
%\end{tikzcd}
where \mL is a cofibrant CW space (that this
%%----------------------------------------------------------------------------------------------12
is possible is implicit in the relevant small object argument).  
Since $X$ is cofibrant, $\exists$ an $s:X \ra L$ such that $p \circ s = \id_X$.  Fix a $j:L \ra K$ for which 
$
\begin{cases}
\ i \circ j \simeq \id_L\\
\ j \circ i \simeq \id_K
\end{cases}
$
($i$ is a weak homotopy equivalence, hence a homotopy equivalence (realization theorem)).  
So: $f \circ (j \circ s) =$ $
(p \circ i) \circ (j \circ s) \simeq$ $p \circ s = \id_X$.  
Therefore $X$ is dominated in homotopy by $K$, thus by the domination theorem is a CW space.

[Note: \  $L$ is a compactly generated Hausdorff space and $s:X \ra L$ is a closed embedding.  
Conclusion: Every cofibrant $X$ is in \bCGH.  
Example: $[0,1]/[0,1[$ is compactly generated (and contractible) but not Hausdorff, hence not cofibrant.]\\
\endgroup%%------------------------------------<<

A model category \bC is said to be 
\un{proper}
\index{proper model category} 
provided that the following axiom is satisfied.

\indent\indent (PMC) \quadx Given a 2-source $X \overset{f}{\lla} Z \overset{g}{\lra} Y$, define $P$ by the pushout square
\begin{tikzcd}%[sep=large]
{Z}  \ar{d}[swap]{f} \ar{r}{g} &{Y} \ar{d}{\eta}\\
{X} \ar{r}[swap]{\xi} &{P}
\end{tikzcd}
.  
Assume: $f$ is a cofibration and $g$ is a weak equivalence $-$then $\xi$ is a weak equivalence.  
Given a 2-sink $X \overset{f}{\lra} Z \overset{g}{\lla} Y$,  define $P$ by the pullback square
\begin{tikzcd}%[sep=large]
{P}  \ar{d}[swap]{\xi} \ar{r}{\eta} &{Y} \ar{d}{g.}\\
{X} \ar{r}[swap]{f} &{Z}
\end{tikzcd}  
Assume: $g$ is a fibration and $f$ is a weak equivalence $-$then $\eta$ is a weak equivalence.

Remark: In a proper model category, Proposition 2 becomes an axiom (no cofibrancy conditions), which suffices to ensure the validity of Proposition 3 (no cofibrancy conditions).\\

\begin{proposition} \ 
Let \bC be a model category.  Assume: All the objects of \bC are cofibrant and fibrant $-$then \bC is proper.
\end{proposition}

[This follows from Proposition 2.]

[Note: \  Not every model category is proper 
(cf. p. \pageref{12.21}).]\\

Example: \bTOP (or \bCG), in its standard structure, is a proper model category.\\

\label{13.54}
\begingroup%%----------------------------------->>
\fontsize{9pt}{11pt}\selectfont
\textbf{\small  EXAMPLE} \ 
\bTOP (or \bCG, $\dcg$, \bCGH), in its singular structure is a proper model category.  
In fact, since every object is fibrant, half of Proposition 5 is immediately applicable.  
However, not every object is cofibrant so for this part an ad hoc argument is necessary.  
Thus consider the commutative diagram
\begin{tikzcd}[sep=large]
{X}  \arrow[d,shift right=0.5,dash] \arrow[d,shift right=-0.5,dash] 
&{Z}  \arrow[d,shift right=0.5,dash] \arrow[d,shift right=-0.5,dash]  \ar{l}[swap]{f} \ar{r}{\id_Z} &{Z}  \ar{d}{g}\\
{X} &{Z} \ar{l}{f}  \ar{r}[swap]{g} &{Z}
\end{tikzcd}
, where $f$ is a cofibration in the singular structure and $g$ is a weak homotopy equivalence 
$-$then $f$ is a closed cofibration, therefore, $\xi:X \ra P$ is a weak homotopy equivalence 
(cf. p. \pageref{12.22}).
\vspi
%%----------------------------------------------------------------------------------------------13
[Note: \  Let $X$ be a topological space which is not compactly generated 
$-$then $\Gamma X$ is not compactly generated and the identity map $k\Gamma X \ra \Gamma X$ is an acyclic Serre fibration, 
so $\Gamma X$ is not cofibrant (but $\Gamma X$ is a CW space).]\\
\endgroup%%------------------------------------<<

\label{13.132}
Let \bC be a proper model category $-$then a commutative diagram
\begin{tikzcd}%[sep=large]
{W}  \ar{d} \ar{r}  &{Y} \ar{d}{g}\\
{X} \ar{r}[swap]{f} &{Z}
\end{tikzcd}
in \bC is said to be a 
\un{homotopy pullback} 
\index{homotopy pullback} 
if for some factorization 
\begin{tikzcd}[sep=small]
{Y}   \ar{r}{\sim}  &{\ov{Y}}  \ar[r, two heads] &{Z}
\end{tikzcd}
of $g$, the induced morphism $W \ra X \times_Z \ov{Y}$ is a weak equivalence.  
This definition is essentially independent of the choice of the factorization of $g$ since any two such factorizations
$
\begin{cases}
{Y}   \overset{\sim}{\ra} {\ov{Y}^\prime}  \twoheadrightarrow {Z}\\
{Y}   \overset{\sim}{\ra} {\ov{Y}^{\prime\prime}}  \twoheadrightarrow {Z}
\end{cases}
$
lead to a commutative diagram
\begin{tikzcd}%[sep=large]
&{X \times_Z \ov{Y}^\prime} \\
{W} \ar{r} \ar{ru} \ar{rd} &{\bullet} \ar{u}[swap]{\sim} \ar{d}{\sim}\\
&{X \times_Z \ov{Y}^{\prime\prime}}
\end{tikzcd}
and it does not matter whether one factors $g$ or $f$ (see below).  
Example: A pullback square
\begin{tikzcd}%[sep=large]
{P}  \ar{d}[swap]{\xi} \ar{r}{\eta}  &{Y} \ar{d}{g}\\
{X} \ar{r}[swap]{f} &{Z}
\end{tikzcd}
is a homotopy pullback provided that $g$ is a fibration.

[Note: \  The dual notion is 
\un{homotopy pushout} 
\index{homotopy pushout}.]\\
\vspace{0.25cm}

\begingroup%%----------------------------------->>
\fontsize{9pt}{11pt}\selectfont
Take two factorizations
$
\begin{cases}
{Y}   \overset{\sim}{\ra} {\ov{Y}^\prime}  \twoheadrightarrow {Z}\\
{Y}   \overset{\sim}{\ra} {\ov{Y}^{\prime\prime}}  \twoheadrightarrow {Z}
\end{cases}
$
of $g$, form the pullback \ $\ov{Y}^\prime \times_Z \ov{Y}^{\prime\prime},$ and note that the projections
$\ov{Y}^\prime \times_Z \ov{Y}^{\prime\prime} \ra \ov{Y}^\prime$,  
$\ov{Y}^\prime \times_Z \ov{Y}^{\prime\prime} \ra \ov{Y}^{\prime\prime}$
are fibrations.  Factor the arrow 
$Y \ra \ov{Y}^\prime \times_Z \ov{Y}^{\prime\prime}$
as
$Y \overset{\sim}{\ra} \ov{W} \twoheadrightarrow \ov{Y}^\prime \times_Z \ov{Y}^{\prime\prime}$. 
Since the diagram
\begin{tikzcd}[sep=large]
{Y}  \ar{d} \ar{r}  \ar{rd} &{\ov{Y}^{\prime\prime}}\\
{\ov{Y}^\prime} &{\ov{W}}  \ar{l} \ar{u}
\end{tikzcd}
commutes, the arrows 
$\ov{W} \ra  \ov{Y}^\prime$, 
$\ov{W} \ra  \ov{Y}^{\prime\prime}$
are weak equivalences.  
Consider the commutative diagrams
\begin{tikzcd}[sep=large]
{X}  \arrow[d,shift right=0.5,dash] \arrow[d,shift right=-0.5,dash] \ar{r}  \ar{r} 
&{Z} \arrow[d,shift right=0.5,dash] \arrow[d,shift right=-0.5,dash] 
&\ov{W} \ar{l} \ar{d}\\
{X} \ar{r} &{Z} &{\ov{Y}^\prime} \ar{l}
\end{tikzcd}
,
$
\begin{tikzcd}[sep=large]
{X}  \arrow[d,shift right=0.5,dash] \arrow[d,shift right=-0.5,dash] \ar{r}  \ar{r} 
&{Z} \arrow[d,shift right=0.5,dash] \arrow[d,shift right=-0.5,dash] 
&\ov{W} \ar{l} \ar{d}\\
{X} \ar{r} &{Z}  &{\ov{Y}^{\prime\prime}} \ar{l}
\end{tikzcd}
.
$
Because the arrows 
$\ov{W} \ra  Z$,
$\ov{Y}^\prime \ra  Z$, 
$\ov{Y}^{\prime\prime} \ra  Z$ 
are fibrations, Proposition 3 implies that the induced morphisms
$X \times_Z \ov{W} \ra X \times_Z \ov{Y}^\prime$, 
$X \times_Z \ov{W} \ra X \times_Z \ov{Y}^{\prime\prime}$
are weak equivalences.  Therefore one may put 
$\bullet = X \times_Z \ov{W}$ 
in the above.
\vspi
[Note: \  Take a factorization 
${Y} \overset{\sim}{\ra} {\ov{Y}}  \twoheadrightarrow {Z}$ 
of $g$ and a factorization 
${X} \overset{\sim}{\ra} {\ov{X}}  \twoheadrightarrow {Z}$ 
of $f$.  
Claim: The induced morphism
$W \ra X \times_Z \ov{Y}$
is a weak equivalence iff the induced morphism 
$W \ra \ov{X} \times_Z {Y}$
is a weak equivalence.  Proof: The diagram
\begin{tikzcd}[sep=large]
{W}  \ar{d} \ar{r}  \ar{r} &{X \times_Z \ov{Y}} \ar{d}\\
{\ov{X} \times_Z {Y}} \ar{r} &{\ov{X} \times_Z \ov{Y}}  
\end{tikzcd}
commutes and the arrows
$X \times_Z \ov{Y} \ra \ov{X} \times_Z \ov{Y}$, 
$\ov{X} \times_Z {Y} \ra \ov{X} \times_Z \ov{Y}$ 
are weak equivalences (cf. Proposition 3).]\\
\endgroup%%------------------------------------<<

%%----------------------------------------------------------------------------------------------14

Example:  In a proper model category \bC, a commutative diagram
\begin{tikzcd}[sep=large]
{W}  \ar{d} \ar{r}  &{Y} \ar{d}{g}\\
{X} \ar{r}[swap]{f} &{Z}
\end{tikzcd}
, where $f$ is a weak equivalence, is a homotopy pullback iff the arrow $W \ra Y$ is a weak equivalence.\\

\textbf{\small COMPOSITION LEMMA} \ 
Consider the commutative diagram
\hspace{-.3cm}
\begin{tikzcd}[sep=large]
{\bullet}  \ar{d} \ar{r}  &{\bullet} \ar{d} \ar{r}  &{\bullet} \ar{d} \\
{\bullet}  \ar{r} &{\bullet}  \ar{r} &{\bullet} 
\end{tikzcd}
in a proper model category \bC.  Suppose that both the squares are homotopy pullbacks $-$then the rectangle is a homotopy pullback.  Conversely, if  the rectangle and the second square are homotopy pullbacks, then the first square is a homotopy pullback.\\

\label{13.128}
\label{13.129}

\begingroup%%----------------------------------->>
\fontsize{9pt}{11pt}\selectfont
\textbf{\small  EXAMPLE} \ 
Take \bC = \bTOP (standard structure) $-$then the commutative diagram
\begin{tikzcd}[sep=large]
{W}  \ar{d} \ar{r}  &{Y} \ar{d}{g}\\
{X} \ar{r}[swap]{f} &{Z}
\end{tikzcd}
 is a homotopy pullback iff the arrow $W \ra W_{f,g}$ is a homotopy equivalence.  Proof: The commutative diagram
\begin{tikzcd}[sep=large]
{W_{f,g}}  \ar{d} \ar{r}  &{Y} \ar{d}{g}\\
{W_f} \ar{r}[swap]{q} &{Z}
\end{tikzcd}
is a pullback square $(f = q \circ s)$ 
(cf. p. \pageref{12.23}).  
One may therefore take this condition as the definition of homotopy pullback in \bTOP.  \ Example:  A pullback square
\begin{tikzcd}[sep=large]
{P}  \ar{d}[swap]{\xi} \ar{r}{\eta}  &{Y} \ar{d}{g}\\
{X} \ar{r}[swap]{f} &{Z}
\end{tikzcd}
is a homotopy pullback provided that $g$ is a Dold fibration (cf. $\S 4$, Propostion 18 (with ``Hurewicz'' replaced by ``Dold'')).
\vspi
[Note: \  Let $W$ be a topological space; let
$
\begin{cases}
\ X\\
\ Y
\end{cases}
$
be pointed topological spaces, $f:X \ra Y$ a pointed continuous function $-$then 
$W  \lra X \overset{f}{\lra} {Y}$
is said to be a 
\un{fibration up to homotopy} 
\index{fibration up to homotopy} 
(or a 
\un{homotopy fiber sequence})
\index{homotopy fiber sequence} if the diagram
\begin{tikzcd}[sep=large]
{W}  \ar{d} \ar{r}  &{\{y_0\}} \ar{d}\\
{X} \ar{r}[swap]{f} &{Y}
\end{tikzcd}
commutes and the induced map $W \ra E_f$ is a homotopy equivalence.  
Because $E_f$ is the double mapping track of the 2-sink 
$X  \overset{f}{\lra} Y \lla \{y_0\}$
, a sequence
$W  \lra X \overset{f}{\lra} {Y}$
is a fibration up to homotopy if the composite $W \ra Y$ is the constant map $W \ra y_0$
%%----------------------------------------------------------------------------------------------15
and the commutative diagram
\begin{tikzcd}[sep=large]
{W}  \ar{d} \ar{r}  &{\{y_0\}} \ar{d}\\
{X} \ar{r}[swap]{f} &{Y}
\end{tikzcd}
is a homotopy pullback.]\\
\endgroup%%------------------------------------<<

\label{14.22}
\label{14.23}
\label{14.74}
\begingroup%%----------------------------------->>
\fontsize{9pt}{11pt}\selectfont
\textbf{\small FACT} \ 
Let
\begin{tikzcd}[sep=large]
{X}  \ar{d}  &{Z} \ar{d} \ar{l}[swap]{f}  \ar{r}{g} &{Y}  \ar{d}\\
{X^\prime}  &{Z^\prime} \ar{l}{f^\prime} \ar{r}[swap]{g^\prime} &{Y^\prime}
\end{tikzcd}
be a commutative diagram of topological spaces in which the squares are homotopy pullbacks $-$then in the commutative diagram
$
\begin{tikzcd}[sep=large]
{X}  \ar{d} \ar{r}  &{M_{f,g}} \ar{d} &{Y} \ar{l} \ar{d}\\
{X^\prime} \ar{r} &{M_{f^\prime,g^\prime}} &{Y^\prime} \ar{l} 
\end{tikzcd}
,
$
the squares are homotopy pullbacks.\\
\endgroup%%------------------------------------<<

\begingroup%%----------------------------------->>
\fontsize{9pt}{11pt}\selectfont
Application: Suppose that 
$
\begin{cases}
\ A \ra X\\
\ A^\prime \ra X^\prime
\end{cases}
$
are closed cofibrations.  Let
$
\begin{cases}
\ f: A \ra Y\\
\ f^\prime:A^\prime \ra Y^\prime
\end{cases}
$
be continuous functions.  Assume that the diagram
\begin{tikzcd}[sep=large]
{X}  \ar{d}  &{A} \ar{l}  \ar{d} \ar{r}{f} &{Y}  \ar{d}\\
{X^\prime} &{A^\prime} \ar{l}  \ar{r}[swap]{f^\prime} &{Y^\prime} 
\end{tikzcd}
commutes and that the squares are homotopy pullbacks $-$then in the commutative diagram
\begin{tikzcd}[sep=large]
{X}  \ar{d} \ar{r}  &{X \sqcup_f Y} \ar{d} &{Y} \ar{l} \ar{d}\\
{X^\prime} \ar{r} &{X^\prime \sqcup_{f^\prime} Y^\prime}  &{Y^\prime} \ar{l}
\end{tikzcd}
, the squares are homotopy pullbacks.\\
\endgroup%%------------------------------------<<

\begingroup%%----------------------------------->>
\fontsize{9pt}{11pt}\selectfont
\textbf{\small FACT} \ 
Let
$
\begin{cases}
\ (\textbf{X},\bff)\\
\ (\bY,\bg)
\end{cases}
$
be objects in $\bFIL(\bTOP)$, $\phi:(\bX,\bff)\ \ra (\bY,\bg)$ a morphism.  Assume: $\forall \ n$, 
\begin{tikzcd}[sep=large]
{X_n}  \ar{d}[swap]{\phi_n} \ar{r}{f_n}  &{X_{n+1}} \ar{d}{\phi_{n+1}} \\
{Y_n} \ar{r}[swap]{g_n} &{Y_{n+1}}
\end{tikzcd}
is a homotopy pullback, $-$then $\forall \ n$, 
\begin{tikzcd}[sep=large]
{X_n}  \ar{d} \ar{r}  &{\tel(\bX,\bff)} \ar{d} \\
{Y_n} \ar{r} &{\text{tel}(\bY,\bg)}
\end{tikzcd}
is a homotopy pullback.\\
\endgroup%%------------------------------------<<

\begingroup%%----------------------------------->>
\fontsize{9pt}{11pt}\selectfont
Application:  Let 
\begin{tikzcd}[sep=large]
{X^0}  \ar{d}\ar{r}  &{X^1} \ar{d} \ar{r} &{\cdots} \\
{Y^0} \ar{r} &{Y^1} \ar{r} &{\cdots}
\end{tikzcd}
be a commutative ladder connecting two expanding sequences of topological spaces.  
Assume: $\forall \ n$, the inclusions
$
\begin{cases}
\ X^n \ra X^{n+1}\\
\ Y^n \ra Y^{n+1}
\end{cases}
$
are cofibrations and 
\begin{tikzcd}[sep=large]
{X^n}  \ar{d} \ar{r} &{X^{n+1}} \ar{d} \\
{Y^n} \ar{r}&{Y^{n+1}}
\end{tikzcd}
is a homotopy pullback $-$then $\forall \ n$, 
\begin{tikzcd}[sep=large]
{X^n}  \ar{d} \ar{r} &{X^\infty} \ar{d} \\
{Y^n} \ar{r}&{Y^\infty}
\end{tikzcd}
is a homotopy pullback.\\
\endgroup%%------------------------------------<<
\vspace{0.5cm}

%%----------------------------------------------------------------------------------------------16
\label{16.44} %dmc mnft
Let \bC be a model category $-$then a morphism $g:Y \ra Z$ in \bC is said to be a 
\un{homotopy fibration} 
\index{homotopy fibration} 
if in any commutative diagram
$
\begin{tikzcd}[sep=large]
{X^\prime \times_Z Y}  \ar{d} \ar{r}{\Phi} &{X \times_Z Y}  \ar{d} \ar{r} &Y \ar{d}{g}\\
{X^\prime} \ar{r}[swap]{\phi} &{X} \ar{r}[swap]{f} &{Z}
\end{tikzcd}
, 
$
$\Phi$ is a weak equivalence whenever $\phi$ is a weak equivalence.  Example: Every fibration in a proper model category is a homotopy fibration.\\

\label{13.131}
\label{18.14}
\textbf{\small LEMMA} \ 
Let \bC be a proper model category.  Suppose that $g:Y \ra Z$ is a homotopy fibration $-$then the pullback square
\begin{tikzcd}%[sep=large]
{X \times_Z Y}  \ar{d} \ar{r} &Y \ar{d}{g}\\
{X} \ar{r}[swap]{f } &{Z}
\end{tikzcd}
is a homotopy pullback.

[Fix factorizations
\begin{tikzcd}[sep=small]
Y  \ar{r}{\sim} &\ov{Y} \ar[r,two heads] &{Z,\ }
%\end{tikzcd}
%, 
%\begin{tikzcd}[sep=small]
X  \ar{r}{\sim} &\ov{X} \ar[r,two heads] &Z
\end{tikzcd}
of $g$, $f$ and form the commutative diagram
\[
\begin{tikzcd}%[sep=large]
{X \times_Z Y}  \ar{d} \ar{r} &{\ov{X} \times_Z Y}  \ar{d} \ar{r} &Y \ar{d}{\sim}\\
{X \times_Z \ov{Y}}  \ar{d} \ar{r} &{\ov{X} \times_Z \ov{Y}}  \ar{d} \ar{r} &\ov{Y} \ar[d,two heads]\\
{X} \ar{r}[swap]{\sim} &{\ov{X}} \ar[r,two heads] &{Z} 
\end{tikzcd}
.
\]
Isolate the upper left hand corner:
\begin{tikzcd}%[sep=large]
{X \times_Z Y}  \ar{d}[swap]{\Phi} \ar{r} &{\ov{X} \times_Z Y}  \ar{d}\\
{X \times_Z \ov{Y}} \ar{r} &{\ov{X} \times_Z \ov{Y}}
\end{tikzcd}
.  From the assumptions, the three unlabeled arrows are weak equivalences.  Therefore $\Phi$ is a weak equivalence.]\\

\begingroup%%----------------------------------->>
\fontsize{9pt}{11pt}\selectfont
\textbf{\small FACT} \ 
The class of homotopy fibrations is closed under composition and the formation of retracts and is pullback stable.\\
\endgroup%%------------------------------------<<

In a model category \bC, one can introduce two notions of ``homotopy'', which are defined respectively via ``cylinder objects'' and ``path objects''.  These considerations then lead to the construction of the homotopy category \bHC of \bC.

\label{12.14}
\indent\indent (CO) \ 
A 
\un{cylinder object} 
\index{cylinder object} 
for $X$ is an object $IX$ in \bC together with a diagram
$X \coprod X \overset{\iota}{\rightarrowtail} IX \overset{\sim}{\ra} X$
that factors the folding map $X \coprod X \ra X$.  Write
$
\begin{cases}
\ i_0:X \ra IX\\
\ i_1:X \ra IX\\
\end{cases}
$
for the arrows
$
\begin{cases}
\ \iota \circ \text{in}_0\\
\ \iota \circ \text{in}_1
\end{cases}
. \ 
$
Since id$_X$ factors as
$
\begin{cases}
\ X \ra IX \overset{\sim}{\ra} X\\
\ X \ra IX \overset{\sim}{\ra} X
\end{cases}
, \ 
$
$
\begin{cases}
\ i_0\\
\ i_1\\
\end{cases}
$
are weak equivalences.  If $X$ is in addition cofibrant, then
$
\begin{cases}
\ i_0\\
\ i_1
\end{cases}
$
are cofibrations.  Proof: $X \coprod X$ is defined by the 
%%----------------------------------------------------------------------------------------------17
pushout square
\begin{tikzcd}%[sep=large]
{\emptyset}  \ar{d} \ar{r} &{X}  \ar{d}{\text{in}_1}\\
{X} \ar{r}[swap]{\text{in}_0} &{X \coprod X}
\end{tikzcd}
, so 
$
\begin{cases}
\ \text{in}_0\\
\ \text{in}_1
\end{cases}
$
are cofibrations and the class of cofibrations is composition closed.

\indent\indent (PO) \ 
A 
\un{path object} 
\index{path object} 
for $X$ is an object $PX$ in \bC together with a diagram
$X \overset{\sim}{\ra} PX \overset{\Pi}{\twoheadrightarrow} X \times X$
that factors the diagonal map $X \ra X \times X$.  Write
$
\begin{cases}
\ p_0:PX \ra X\\
\ p_1:PX \ra X\\
\end{cases}
$
for the arrows
$
\begin{cases}
\ \pr_0 \circ \Pi\\
\ \pr_1 \circ \Pi
\end{cases}
. \ 
$
Since id$_X$ factors as
$
\begin{cases}
\ X \overset{\sim}{\ra} PX \ra X\\
\ X \overset{\sim}{\ra} PX \ra X
\end{cases}
, \ 
$
$
\begin{cases}
\ p_0\\
\ p_1
\end{cases}
$
are weak equivalences.  If $X$ is in addition fibrant, then
$
\begin{cases}
\ p_0\\
\ p_1
\end{cases}
$
are fibrations.  Proof: $X \times X$ is defined by the pullback square
\begin{tikzcd}%[sep=large]
{X \times X}  \ar{d}[swap]{\pr_0} \ar{r}{\pr_1} &{X}  \ar{d}\\
{X} \ar{r} &{*}
\end{tikzcd}
, so 
$
\begin{cases}
\ \pr_0\\
\ \pr_1
\end{cases}
$
are fibrations and the class of fibrations is composition closed.

[Note: \  Cylinder objects and path objects exist (cf. MC-5).]\\

\begingroup%%----------------------------------->>
\fontsize{9pt}{11pt}\selectfont
\textbf{\small  EXAMPLE} \ 
Take \bC  = \bTOP  (standard structure) $-$then a choice for $IX$ is $X \times [0,1]$ 
(cf. p. \pageref{12.24}) 
and a choice for $PX$ is $C([0,1],X)$ 
(cf. \pageref{12.25}).\\
\endgroup%%------------------------------------<<

\begingroup%%----------------------------------->>
\fontsize{9pt}{11pt}\selectfont
\textbf{\small  EXAMPLE} \ 
Take \bC  = \bTOP  (singular structure) $-$then a choice for $IX$ is $X \times [0,1]$ if $X$ is a CW complex (but not in general).  However, for any $X$, a choice for $PX$ is $C([0,1],X)$.

[Note: \  Let $X$ be the Warsaw circle $-$then the inclusion 
$i_0X \cup i_1 X \ra X \times [0,1]$ is not a cofibration in the singular structure.  
Thus consider
\begin{tikzcd}[sep=large]
{i_0X \cup i_1 X}  \ar{d} \ar{r}{f} &{X}  \ar{d}\\
{X \times [0,1]} \ar{r} &{*}
\end{tikzcd}
, where 
$
\begin{cases}
\ f(x,0) = x\\
\ f(x,1) = x_0
\end{cases}
.
$
Since $X \ra *$ is a Serre fibration and a weak homotopy equivalence, the existence of a filler for this diagram would mean that $X$ is contractible which it isn't.]\\
\endgroup%%------------------------------------<<

\label{13.33}
\begingroup%%----------------------------------->>
\fontsize{9pt}{11pt}\selectfont
\textbf{\small LEMMA} \ 
Let $(K,L)$ be a relative CW complex, where $K$ is a LCH space.  
Suppose that $X \ra B$ is a Serre fibration $-$then the arrow 
$C(K,X) \ra C(L,X) \times_{C(L,B)} C(K,B)$ is a Serre fibration which is a weak homotopy equivalence if this is the case of 
$L \ra K$ or $X \ra B$.
\vspi
[Note: \  Dropping the assumption that $(K,L)$ is a relative CW complex and supposing only that $L \ra K$ is a closed cofibration (with $K$ a LCH space), the result continues to hold if ``Serre'' is replaced by ``Hurewicz'' and weak homotopy equivalence by homotopy equivalence.]\\
\endgroup%%------------------------------------<<

\begingroup%%----------------------------------->>
\fontsize{9pt}{11pt}\selectfont
Application: Let $(K,L)$ be a relative CW complex, where $K$ is a LCH space.   
Suppose that $A \ra Y$ is a cofibration in the singular structure $-$then the arrow
$L \times Y \cup K \times A \ra K \times Y$ is a cofibration in the singular structure which is a weak homotopy equivalence if this is the case of $L \ra K$ or $A \ra Y$.\\
\endgroup%%------------------------------------<<

%%----------------------------------------------------------------------------------------------18

\label{12.27}
\begingroup%%----------------------------------->>
\fontsize{9pt}{11pt}\selectfont
\textbf{\small  EXAMPLE} \ 
Take $L = \{0,1\}$, $K = [0,1]$ $-$then for any cofibration $A \ra Y$ in the singular structure, the inclusion 
$i_0 Y \cup A \times [0,1] \cup i_1 Y \ra Y \times [0,1]$ is a cofibration in the singular structure 
(cf. p. \pageref{12.26}).  
In particular, $\forall$ cofibrant $X$, a choice for $IX$ is $X \times [0,1]$.\\
\endgroup%%------------------------------------<<

\indent\indent (LH) \ Morphisms $f, g:X \ra Y$ in \bC are said to be 
\un{left homotopic} 
\index{left homotopic} 
if $\exists$ a cylinder object $IX$ for $X$ and a morphism 
$H:IX \ra Y$ such that 
$H \circ i_0 = f$, 
$H \circ i_1 = g$.
One calls $H$ a 
\un{left homotopy} 
\index{left homotopy} 
between $f$ and $g$.  
Notation: $f \underset{l}{\simeq}g$.  
\index{$\underset{l}{\simeq}$}  
If $Y$ is fibrant and if $f \underset{l}{\simeq}g$, then $\exists$ a cylinder object $I^\prime X$ for $X$ with 
$X \coprod X \overset{\iota^\prime}{\rightarrowtail} I^\prime X \overset{\sim}{\twoheadrightarrow} X$
and a left homotopy $H^\prime:I^\prime X \ra Y$ between $f$ and $g$.
Proof: Factor $IX \overset{\sim}{\ra} X$ as 
$IX \overset{\sim}{\rightarrowtail} I^\prime X \overset{\sim}{\twoheadrightarrow} X$ 
and consider a filler $H^\prime:I^\prime X \ra Y$ for the commutative diagram
$
\begin{tikzcd}%[sep=large]
{IX}  \ar{d} \ar{r}{H} &{Y}  \ar{d}\\
{I^\prime X} \ar{r} &{*}
\end{tikzcd}
.
$

[Note: \  Suppose that $f \underset{l}{\simeq}g$ $-$then $f$ is a weak equivalence iff $g$ is a weak equivalence.]

\indent\indent (RH) \ Morphisms $f, g:X \ra Y$ in \bC are said to be 
\un{right homotopic} 
\index{right homotopic} 
if $\exists$ a path object $PY$ for $Y$ and a morphism $G:X \ra PY$ such that 
$p_0 \circ G = f$, 
$p_1 \circ G = g$.
One calls \mG a 
\un{right homotopy}
\index{right homotopy} 
between $f$ and $g$.  
Notation: $f \underset{r}{\simeq}g$.  
\index{$\underset{l}{\simeq}$}  
If $X$ is cofibrant and if $f \underset{r}{\simeq}g$, then $\exists$ a path object $P^\prime Y$ for $Y$ with 
$Y \overset{\sim}{\rightarrowtail} P^\prime Y \overset{\Pi^\prime}{\twoheadrightarrow} Y \times Y$
and a right homotopy 
$G^\prime:X \ra P^\prime Y$ between $f$ and g.
Proof: Factor $Y \overset{\sim}{\ra} PY$ as 
$Y \overset{\sim}{\rightarrowtail} P^\prime Y \overset{\sim}{\twoheadrightarrow} PY$ 
and consider a filler $G^\prime:X \ra P^\prime Y$ for the commutative diagram
$
\begin{tikzcd}%[sep=large]
{\emptyset}  \ar{d} \ar{r} &{P^\prime Y}  \ar{d}\\
{X} \ar{r}[swap]{G} &{PY}
\end{tikzcd}
.
$

[Note: \  Suppose that $f \underset{r}{\simeq}g$ $-$then $f$ is a weak equivalence iff $g$ is a weak equivalence.]\\

Notation: Given $X$, $Y$ $\in $ Ob\bC, let
$
\begin{cases}
\ [X,Y]_l\\
\ [X,Y]_r
\end{cases}
$
be the set of equivalence classes in $\Mor(X,Y)$ under the equivalence relation generated by 
$
\begin{cases}
\ \text{left}\\
\ \text{right}
\end{cases}
$
homotopy.

[Note: \  The relations of
$
\begin{cases}
\ \text{left}\\
\ \text{right}
\end{cases}
$
homotopy are reflexive and symmetric but not necessarily transitive.  
Elements of 
$
\begin{cases}
\ [X,Y]_l\\
\ [X,Y]_r
\end{cases}
$
are denoted by
$
\begin{cases}
\ [f]_l\\
\ [f]_r
\end{cases}
$
and referred to as
$
\begin{cases}
\ \text{left}\\
\ \text{right}
\end{cases}
$
homotopy classes of morphisms.]
\\
\vspace{0.5cm}

\begingroup%%----------------------------------->>
\fontsize{9pt}{11pt}\selectfont
Left homotopy is reflexive.  Proof: Given $f:X \ra Y$, take for $H$ the composition
$IX \overset{\sim}{\ra} X \overset{f}{\ra} Y$.
\vspi
Left homotopy is symmetric.  Proof: Given $f,g:X \ra Y$, and $H:IX \ra Y$ such that 
$H \circ i_0 = f$, $H \circ i_1 = g$, let 
$\mathsf{T}:X \coprod X \ra X \coprod X$ be the interchange, note that
$X \coprod X \overset{\iota \circ \mathsf{T}}{\ra} IX \overset{\sim}{\ra} X$
factors the folding map $X \coprod X  \ra X$, and 
$H \circ (\iota \circ \mathsf{T}) \circ \text{in}_0 = g$, 
$H \circ (\iota \circ \mathsf{T}) \circ \text{in}_1 = f$.\\
\endgroup%%------------------------------------<<

%%----------------------------------------------------------------------------------------------19

\begin{proposition} \ 
Left homotopy is an equivalence relation on $\Mor(X,Y)$ if $X$ is cofibrant and right homotopy is an equivalence relation on $\Mor(X,Y)$ if $Y$ is fibrant.\
\end{proposition}

[To check transitivity in the case of left homotopy, suppose that $f \underset{l}{\simeq} g$ $\&$ $g \underset{l}{\simeq} h$, say
$
\begin{cases}
\ H \circ i_0 = f\\
\ H \circ i_1 = g
\end{cases}
$
$\&$
$
\begin{cases}
\ H \circ i_0^\prime = g\\
\ H \circ i_1^\prime = h
\end{cases}
. \ 
$
Define $I^{\prime\prime} X$ by the pushout square
\begin{tikzcd}%[sep=large]
{X}  \ar{d}[swap]{i_1} \ar{r}{i_0^\prime} &{I^\prime X}  \ar{d}{j_1}\\
{IX} \ar{r}[swap]{j_0^\prime} &{I^{\prime\prime}X}
\end{tikzcd}
$-$then $I^{\prime\prime} X$ is a cylinder object for $X$ 
(specify $\iota^{\prime\prime}:X \coprod X \ra I^{\prime\prime} X$ by
$
\begin{cases}
\ \iota^{\prime\prime} \circ \text{in}_0 = j_0^\prime \circ i_0\\
\ \iota^{\prime\prime} \circ \text{in}_1 = j_1 \circ i_1^\prime
\end{cases}
). 
$
Moreover, $H \circ i_1 = H^\prime \circ i_0^\prime$ $\implies$ $\exists$ 
$H^{\prime\prime}:I^{\prime\prime} X \ra Y$:
$
\begin{cases}
\ H^{\prime\prime} \circ \iota_0^{\prime\prime} = f\\
\ H^{\prime\prime} \circ \iota_1^{\prime\prime} = h
\end{cases}
.]
$

[Note: \  Here is the verification that $\iota^{\prime\prime}$ is a cofibration.  Form the commutative diagram
\begin{tikzcd}%[sep=large]
{X}  \ar{d}[swap]{i_0}  &{\emptyset} \ar{l} \ar{d} \ar{r} &{X} \ar{d}{i_1^\prime}\\
{IX}  &{X} \ar{l}{i_1} \ar{r}[swap]{i_0^\prime}  &{I^{\prime}X}
\end{tikzcd}
and apply Proposition 4.]\\

\begin{proposition} \ %07
If $X$ is cofibrant and $p:Y \ra Z$ is an acyclic fibration, then the postcomposition arrow $p_*:[X,Y]_l \ra [X,Z]_l$ is bijective, while if $Z$ is fibrant and $i:X \ra Y$ is an acyclic cofibration, then the precomposition arrow $i^*:[Y,Z]_r \ra [X,Z]_r$ is bijective.
\end{proposition}

[In either case, the arrows are welldefined.  That $p_*$ is surjective follows from the fact that, generically,
\begin{tikzcd}%[sep=large]
{\emptyset}  \ar{d} \ar{r} &{Y}  \ar{d}{p}\\
{X} \ar{r} &{Z}
\end{tikzcd}
has a filler $X \ra Y$.  Assume now that $p \circ f \underset{l}{\simeq} p \circ g$, where $f$, $g \in \Mor(X,Y)$.  
Choose $H:IX \ra Z$ with 
$
\begin{cases}
\ H \circ i_0 = p \circ f\\
\ H \circ i_1 = p \circ g
\end{cases}
$
$-$then any filler $IX \ra Y$ in 
\begin{tikzcd}%[sep=large]
{X \coprod X}  \ar{d}[swap]{\iota} \ar{r}{f \coprod g} &{Y}  \ar{d}{p}\\
{IX} \ar{r}[swap]{H} &{Z}
\end{tikzcd}
is a left homotopy between $f$ and $g$.  Therefore $p_*$ is injective.]\\
\vspace{0.25cm}

\begingroup%%----------------------------------->>
\fontsize{9pt}{11pt}\selectfont
\textbf{\small FACT} \ 
Suppose that 
$
\begin{cases}
\ Y\\
\ Z
\end{cases}
$
are fibrant and $p:Y \ra Z$ is a weak equivalence $-$then for any $X$, the postcomposition arrow $p_*:[X,Y]_r \ra [X,Z]_r$ is injective.\\
\endgroup%%------------------------------------<<

\begingroup%%----------------------------------->>
\fontsize{9pt}{11pt}\selectfont
\textbf{\small FACT} \ 
Suppose that 
$
\begin{cases}
\ X\\
\ Y
\end{cases}
$
are cofibrant and $i:X \ra Y$ is a weak equivalence $-$then for any \mZ, the precomposition arrow $i^*:[Y,Z]_l \ra [X,Z]_l$ is injective.\\
\endgroup%%------------------------------------<<

\textbf{\small LEMMA (LH)} \quadx
Let $f$, $g \in \Mor(X,Y)$ be left homotopic.  
Assume: $Y$ is fibrant $-$then $\forall \ \phi: X^\prime \ra X$, $f \circ \phi \underset{l}{\simeq} g \circ \phi$.

%%----------------------------------------------------------------------------------------------20
[Since $Y$ is fibrant, one can arrange that the left homotopy $H:IX \ra Y$ between $f$ and $g$ is computed per 
$X \coprod X \overset{\iota}{\rightarrowtail} IX \overset{\sim}{\twoheadrightarrow} X$ (cf. LH).  This said, form the commutative diagram
\begin{tikzcd}%[sep=large]
{X^\prime \coprod X^\prime}  \ar{d}[swap]{\iota^\prime} \ar{r}{\phi \coprod \phi} 
&{X \coprod X}  \ar{r}{\iota} &{IX} \ar{d}\\
{IX^\prime} \ar{r} &{X^\prime} \ar{r}[swap]{\phi} &{X}
\end{tikzcd}
, choose a filler $\Phi:IX^\prime \ra IX$, and note that $H \circ \Phi$ is a left homotopy between 
$f \circ \phi$ and $g \circ \phi$.]\\

\begin{proposition} \ %08
\textbf{(LH)} \ 
Suppose that $Y$ is fibrant $-$then the composition in $\Mor\bC$ induces a map
$[X^\prime,X]_l \times [X,Y]_l \ra [X^\prime,Y]_l$.
\end{proposition}

[The contention is that $[f]_l = [g]_l$ ($f$, $g \in \Mor(X,Y)$) $\&$ $[\phi]_l = [\psi]_l$ 
($\phi, \psi \in \Mor(X^\prime,X)$) $\implies$ $[f \circ \phi]_l =$ $[g \circ \psi]_l$.  
From the definitions, 
$\exists$ $f_1, \ldots, f_n \in \Mor(X,Y)$ : $f_1 = f$, $f_n = g$ 
with $f_i \underset{l}{\simeq} f_{i+1}$, hence by the lemma, 
$f_i \circ \phi \underset{l}{\simeq} f_{i+1} \circ \phi$ $\forall \ i$ $\implies$ $[f \circ \phi]_l = [g \circ \phi]_l$.  
But trivially, $[g \circ \phi]_l = [g \circ \psi]_l$.]\\

\textbf{\small LEMMA (RH)} \quadx
Let $f$, $g \in \Mor(X,Y)$ be right homotopic.  
Assume: $X$ is cofibrant $-$then $\forall \ \psi:Y \ra Y^\prime$, $\psi \circ f \underset{r}{\simeq} \psi \circ g$.\\

\setcounter{proposition}{7}
\begin{proposition} \ 
\textbf{(RH)} \ 
Suppose that $X$ is cofibrant $-$then  composition in $\Mor\bC$ induces a map
$[X,Y]_r \times [Y,Y^\prime]_r \ra [X,Y^\prime]_r$.\\
\end{proposition}

\begingroup%%----------------------------------->>
\fontsize{9pt}{11pt}\selectfont
\textbf{\small FACT} \ 
Let $f$, $g \in \Mor(X,Y)$ be left homotopic.  
Suppose that $\phi:X^\prime \ra X$ is an acyclic fibration $-$then $f \circ \phi \underset{l}{\simeq} g \circ \phi$.\\
\endgroup%%------------------------------------<<

\begingroup%%----------------------------------->>
\fontsize{9pt}{11pt}\selectfont
\textbf{\small FACT} \ 
Let $f$, $g \in \Mor(X,Y)$ be right homotopic.  
Suppose that $\psi:Y \ra Y^\prime $ is an acyclic cofibration $-$then $\psi\circ f \underset{r}{\simeq} \psi \circ g$.\\
\endgroup%%------------------------------------<<

\begin{proposition} \ 
Let $f$, $g \in \Mor(X,Y)$ $-$then
(i) $X$ cofibrant $\&$ $f \underset{l}{\simeq} g$ $\implies$ $f \underset{r}{\simeq} g$ and 
(ii) $Y$ fibrant $\&$ $f \underset{r}{\simeq} g$ $\implies$ $f \underset{l}{\simeq} g$.
\end{proposition}

[We shall prove (i), the proof of (ii) being analogous.  
Choose a left homotopy $H:IX \ra Y$ between $f$ and $g$ and let
$p:IX \ra X$ be the abient weak equivalence.  
Fix a path object $PY$ for $Y$ and let $j:Y \ra PY$ be the ambient weak equivalence.  
Since $X$ is cofibrant, $i_0$ is an acyclic cofibration, thus the commutative diagram
\begin{tikzcd}%[sep=large]
{X}  \ar{d}[swap]{i_0} \ar{r}{j \circ f} &{PY}  \ar{d}{\Pi}\\
{IX} \ar{r}[swap]{(f \circ p,H)} &{Y \times Y} 
\end{tikzcd}
has a filler $\rho:IX \ra PY$ and the composite $G = \rho \circ i_1$ is a right homotopy between $f$ and g.]\\
%

%%----------------------------------------------------------------------------------------------21
\label{12.29}
Notation: Given a cofibrant $X$ and a fibrant $Y$, write $\simeq$ for 
$\underset{l}{\simeq} \ = \  \underset{r}{\simeq}$, call this equivalence relation 
\un{homotopy},
\index{homotopy (model category)} 
and let $[X,Y]$ be the set of homotopy classes of morphisms in $\Mor(X,Y)$, a typical element being [$f$].

[Note: \  If $f \simeq g$, then $f$ is a weak equivalence iff $g$ is a weak equivalence.]

Observation: Suppose that $X$ is cofibrant and $Y$ is fibrant.  Let $f$, $g \in \Mor(X,Y)$ $-$then the following conditions are equivalent: 
(1) $f$ and $g$ are left homotopic; 
(2) $f$ and $g$ are right homotopic with respect to a fixed choice of path object; 
(3) $f$ and $g$ are right homotopic ;
(4) $f$ and $g$ are left homotopic with respect to a fixed choice of cylinder object.\\

\begingroup%%----------------------------------->>
\fontsize{9pt}{11pt}\selectfont
\textbf{\small FACT} \ 
Let 
\begin{tikzcd}[sep=large]
{X}  \ar{d}[swap]{\phi} \ar{r}{f} &{Y}  \ar{d}{\psi}\\
{W} \ar{r}[swap]{g} &{Z}
\end{tikzcd}
be a diagram in \bC, where $X$ is cofibrant and $Z$ is fibrant.  Assume: $\psi \circ f \simeq g \circ \phi$ $-$then if $W$ is fibrant and $g$ is a fibration, $\exists \ \widetilde{\phi}:X \ra W$ such that $\phi \simeq \widetilde{\phi}$ $\&$ $g \circ \widetilde{\phi} =$ $\psi \circ f$ and if $Y$ is cofibrant and $f$ is a cofibration, $\exists \ \widetilde{\psi}:Y \ra Z$ such that 
$\psi \simeq \widetilde{\psi}$ $\&$ $\widetilde{\psi} \circ f =$ $g \circ \phi$.\\
\endgroup%%------------------------------------<<

\begin{proposition} \ %10
Suppose that 
$
\begin{cases}
\ X\\
\ Y
\end{cases}
$
are cofibrant and fibrant.  Let $f \in \Mor(X,Y)$ $-$then $f$ is a weak equivalence iff $f$ has a homotopy inverse, i.e., iff there exists a $g \in \Mor(Y,X)$ such that $g \circ f \simeq \id_X$ $\&$ $f \circ g \simeq \id_Y$.
\end{proposition}

[Necessity:  Write $f = p \circ i$, where $i:X \ra Z$ is an acyclic cofibration and $p:Z \ra Y$ is a fibration.  Note that $Z$ is both cofibrant and fibrant and $p$ is a weak equivalence.  
Fix a filler $r:Z \ra X$ for
\begin{tikzcd}%[sep=large]
{X}  \ar{d}[swap]{i} \ar[equal]{r} &{X}  \ar{d}\\
{Z} \ar{r}[swap]{H} &{*}
\end{tikzcd}
.  
Since $i^*([i \circ r]) =$ $[i \circ r \circ i] =$ $[i] = [\id_Z \circ i] =$ $i^*([\id_Z])$, it follows that  $i \circ r \simeq \id_Z$ (cf. Proposition 7).  
Therefore $r$ is a homotopy inverse for $i$.  Similarly, $p$ admits a homotopy inverse $s$.  
Put $g = r \circ s$ $-$then $g:Y \ra X$ is a homotopy inverse for $f$.

Sufficiency: Decompose $f$ as above: $f = p \circ i$.  
Because $i$ is a weak equivalence, one has only to prove that $p$ is a weak equivalence.  
Let $g:Y \ra X$ be a homotopy inverse for $f$.  
Fix a left homotopy $H:IY \ra Y$ between $f \circ g$ and id$_Y$ and choose a filler $H^\prime:IY \ra Z$ in 
\begin{tikzcd}%[sep=large]
{Y}  \ar{d}[swap]{i_0} \ar{r}{i \circ g} &{Z}  \ar{d}{p}\\
{IY} \ar{r}[swap]{H} &{Y}
\end{tikzcd}
.  
Set $s = H^\prime \circ i_1$ ($\implies p \circ s = \id_Y$).  If $r:Z \ra X$ is a homotopy inverse for $i$, then 
$p \simeq f \circ r$ $\implies$ $s \circ p \simeq i \circ g \circ p \simeq$ 
$i \circ g \circ f \circ r \simeq$ $i \circ r \simeq$ id$_Z$, so $s \circ p$ is a weak equivalence.  
But $p$ is a retract of $s \circ p$, hence it too is a weak equivalence.]\\

\begingroup%%----------------------------------->>
\fontsize{9pt}{11pt}\selectfont
\textbf{\small  EXAMPLE} \ 
Take \bC = \bTOP (singular structure) and let $X$, $Y$ be cofibrant, e.g., CW complexes $-$then Proposition 10 says that a weak homotopy equivalence $f:X \ra Y$ is a homotopy equivalence, which, 
%%----------------------------------------------------------------------------------------------22
when specialized to $X$,$Y$ CW complexes, is the realization theorem.
\vspi
[Note: \  Bear in mind that a cylinder object for a cofibrant $X$, $Y$ is $IX$, $IY$ 
(cf. p. \pageref{12.27}).]\\
\endgroup%%------------------------------------<<

%dmc come back and maybe add infra to a symbol index
Notation:   $\bC_{\bc}$ is the full subcategory of \bC whose objects are cofibrant, $\bC_\bff$
is the full subcategory of \bC whose objects are fibrant, and $\bC_{\bcf}$ is the full subcategory of \bC whose objects are cofibrant and fibrant.
$\bH_{\br}\bC_{\bc}$ is the category with Ob$\bH_{\br}\bC_{\bc} =$ Ob$\bC_{\bc}$ and Mor$\bH_{\br}\bC_{\bc} = $ right homotopy classes of morphisms (cf. Proposition 8 (RH)), 
$\bH_{\bl}\bC_\bff$ is the category with Ob$\bH_{\bl}\bC_\bff =$ Ob$\bC_\bff$ and Mor$\bH_{\bl}\bC_\bff = $ left homotopy classes of morphisms (cf. Proposition 8 (LH)), 
and 
$\bH\bC_{\bcf}$ is the category with 
Ob$\bH\bC_{\bcf} =$ Ob$\bC_{\bcf}$ and 
Mor$\bH\bC_{\bcf} = $ homotopy classes of morphisms (cf. Proposition 9).

[Note: \  Write $\bH\bC_{\bc}$  ($\bH\bC_\bff$) for $\bH\bC_{\bcf}$ if all objects are fibrant (cofibrant).]
\label{12.31}
\label{14.33}

Given $X \in $ Ob\bC, use MC-5 to factor 
$\emptyset \ra X$ as $\emptyset \rightarrowtail \sL X \overset{\sim}{\twoheadrightarrow} X$ 
and 
$X \ra *$ as $X \overset{\sim}{\rightarrowtail} \sR X \twoheadrightarrow *$,
thus $\pi_X:\sL X \ra X$ is an acyclic fibration and $\iota_X: X \ra \sR X$ is an acyclic cofibration.

[Note: \  $\sL X$ is cofibrant and $\sR X$ is fibrant.  If $X$ is cofibrant, take $\sL X = X$ $\&$ $\pi_X = \id_X$ and if $X$ is fibrant, take $\sR X = X$ $\&$ $\iota_X = \id_X$.]\\

\index{Lemma \textbf{$\sL$}}
\textbf{\small  LEMMA $\sL$} \quadx
Fix
$
\begin{cases}
\ X\\
\ Y
\end{cases}
$
$\in$ $\Ob\bC$ and let $f \in  \Mor(X,Y)$ $-$then there exists $\sL f \in \Mor(\sL X,\sL Y)$ such that the diagram
\begin{tikzcd}%[sep=large]
{\sL X}  \ar{d}[swap]{\pi_X} \ar{r}{\sL f} &{\sL Y}  \ar{d}{\pi_Y}\\
{X} \ar{r}[swap]{f} &{Y}
\end{tikzcd}
commutes.  $\sL f$ is uniquely determined up to left homotopy and is a weak equivalence iff $f$ is.  Moreover, for fibrant $Y$, 
$\sL f$ is uniquely determined up to left homotopy by $[f]_l$.

[To establish the existence of $\sL f$, consider any filler $\sL X \ra \sL Y$ for
$
\begin{tikzcd}%[sep=large]
{\emptyset}  \ar{d} \ar{r}&{\sL Y}  \ar{d}{\pi_Y}\\
{\sL X} \ar{r}[swap]{f \circ \pi_X} &{Y}
\end{tikzcd}
.  
$
Since $\sL X$ is cofibrant and $\pi_Y$ is an acyclic fibration, the postcomposition arrow 
$[\sL X,\sL Y]_l \ra [\sL X,Y]_l$ determined by $\pi_Y$ is bijective (cf. Proposition 7).  
This implies that $\sL f$ is unique up to left homotopy.  
The weak equivalence assertion is clear.  
Finally, if $Y$ is fibrant, then composition in $\Mor\bC$ induces a map 
$[\sL X,X]_l \times [X,Y]_l \ra [\sL X,Y]_l$ (cf. Proposition 8 (LH)).  Therefore 
$[f]_l = [g]_l$ $\implies$ $[f\circ \pi_X]_l = [g \circ \pi_X]_l$ $\implies$ $[\pi_Y \circ \sL f]_l = [\pi_Y \circ \sL g]_l$ $\implies$ $\sL f \underset{l}{\simeq} \sL g$ (cf. Proposition 7).]\\

Application: 
$\sL\id_X \underset{l}{\simeq} \id_{\sL X}$ $\implies$  
$\sL\id_X \underset{r}{\simeq} \id_{\sL X}$ 
and 
$\sL(g \circ f) \underset{l}{\simeq} \sL g \circ \sL f$ $\implies$ 
$\sL(g \circ f) \underset{r}{\simeq} \sL g \circ \sL f$ 
(cf. Proposition 9), thus there is a functor 
$\sL:\bC \ra \bH_\br\bC_\bc$ that takes $X$ to $\sL X$ and $f:X \ra Y$ to 
$[\sL f]_r \in [\sL X,\sL Y]_r$.\\

%%----------------------------------------------------------------------------------------------23
\label{13.79}
\index{Lemma \textbf{$\sR$}}
\textbf{\small  LEMMA $\sR$} \quadx
Fix
$
\begin{cases}
\ X\\
\ Y
\end{cases}
$
$\in \Ob\bC$ and let $f \in \Mor(X,Y)$ $-$then there exists $\sR f \in \Mor(\sR X,\sR Y)$ such that the diagram
\begin{tikzcd}%[sep=large]
{X}  \ar{d}[swap]{\iota_X} \ar{r}{f} &{Y}  \ar{d}{\iota_Y}\\
{\sR X} \ar{r}[swap]{\sR f} &{\sR Y}
\end{tikzcd}
commutes.  $\sR f$ is uniquely determined up to right homotopy and is a weak equivalence iff $f$ is.  
Moreover, for cofibrant $X$, 
$\sR f$ is uniquely determined up to right homotopy by $[f]_r$.\\

Application: 
$\sR\id_X \underset{r}{\simeq} \id_{\sR X}$ $\implies$  
$\sR\id_X \underset{l}{\simeq} \id_{\sR X}$
and 
$\sR(g \circ f) \underset{r}{\simeq} \sR g \circ \sR f$ $\implies$ 
$\sR(g \circ f) \underset{l}{\simeq} \sR g \circ \sR f$ 
(cf. Proposition 9), thus there is a functor 
$\sR:\bC \ra \bH_\bl\bC_\bff$ that takes $X$ to $\sR X$ and $f:X \ra Y$ to 
$[\sR f]_l \in [\sR X,\sR Y]_l$.\\

\begingroup%%----------------------------------->>
\fontsize{9pt}{11pt}\selectfont
\index{Reedy's Lifting Lemma}
\textbf{\small REEDY'S LIFTING LEMMA} \ 
Suppose that 
$
\begin{cases}
\ X\\
\ Y
\end{cases}
$
are cofibrant.  Let $\phi \in \Mor(X,Y)$ $-$then $\phi$ is a weak equivalence iff given any commutative diagram
\begin{tikzcd}[sep=large]
{X}  \ar{d}[swap]{\phi} \ar{r}{u} &{U}  \ar{d}{\Phi}\\
{Y} \ar{r}[swap]{v} &{V}
\end{tikzcd}
, where $\Phi$ is a fibration, $\exists$ $w:Y \ra U$ $\&$ $H:IX \ra U$ such that $\Phi \circx w = v$,
$
\begin{cases}
\ H \circx i_0 = u\\
\ H \circx i_1 = w \circx \phi
\end{cases}
, \ 
$
and $\Phi \circx H = v \circx \phi \circx p$, $p:IX \ra X$ the projection.
\vspi
[Necessity: Write $\phi = \eta \circx \xi$, where $\xi:X \ra Z$ is an acyclic cofibration and $\eta:Z \ra Y$ is an acyclic fibration.  Define $IZ$ by the pushout square
\begin{tikzcd}[sep=large]
{X \coprod X}  \ar{d} \ar{r}{\iota} &{IX}  \ar{d}\\
{Z \coprod Z} \ar{r} &{IZ}
\end{tikzcd}
to get a cylinder object for $Z$ compatible with that for $X$ in the sense that there is a commutative diagram
\begin{tikzcd}[sep=large]
{IX}  \ar{d}[swap]{I\xi} \ar{r}{p} &{X}  \ar{d}{\xi}\\
{IZ} \ar{r}[swap]{p} &{Z}
\end{tikzcd}
.  \ 
Since $Y$ is cofibrant, one can find an $s:Y \ra Z$ such that $\eta \circx s = \id_Y$.  
Therefore $\eta \circx \id_Z = \eta \circx (s \circx \eta)$ $\implies$ $\exists$ $h:IZ \ra Z$ such that
$
\begin{cases}
\ h \circx i_0 = \id_Z\\
\ h \circx i_1 = s \circx \eta
\end{cases}
$
and $\eta \circx h = \eta \circx p$:
\begin{tikzcd}[sep=large]
{IZ}  \ar{d}[swap]{p} \ar{r}{h} &{Z}  \ar{d}{\eta}\\
{Z} \ar{r}[swap]{\eta} &{Y}
\end{tikzcd}
(cf. Proposition 7 and its proof).  Choose now a filler $\sigma:Z \ra U$ for 
\begin{tikzcd}[sep=large]
{X}  \ar{d}[swap]{\xi} \ar{r}{u} &{U}  \ar{d}{\Phi}\\
{Z} \ar{r}[swap]{v\circx \eta} &{V}
\end{tikzcd}
.  \ 
Definition: $w = \sigma \circx s$ $\&$ $H = \sigma \circx h \circx I\xi$.  
So, e.g., 
$\Phi \circx H = \Phi \circx \sigma \circx h \circx I\xi = $ $v \circx \eta \circx h \circx I\xi = $ 
$v \circx \eta \circx p \circx I\xi = $ 
$v \circx \eta \circx \xi \circx p = $ $ v \circx \phi \circx p$.
\vspi
Sufficiency: If $\phi:X \ra Y$ has the stated property, then for every fibrant \mZ, $\phi^*:[Y,Z]_l \ra [X,Z]_l$
%%----------------------------------------------------------------------------------------------24
is surjective and $\phi^*:[Y,Z]_r \ra [X,Z]_r$ is injective, hence $\phi^*:[Y,Z] \ra [X,Z]$ is bijective.  Because the horizontal arrows in the commutative diagram
\begin{tikzcd}[sep=large]
{[\sR Y,\sL Z]}  \ar{d} \ar{r} &{[Y,Z]}  \ar{d}\\
{[\sR X,\sL Z]} \ar{r} &{[X,Z]}
\end{tikzcd}
are bijective, $(\sR \phi)^*:[\sR Y,\sL Z] \ra [\sR X,\sL Z]$ is also bijective for every fibrant \mZ.  
Take
$Z = \sR\sL X$: $\sL Z = \sL\sR\sL X$ $= \sR\sL X$ $= \sR X$ $\implies$ $\exists$ $\psi:\sR Y \ra \sR X$ such that 
$(\sR\phi)^*([\psi]) = [\id_{\sR X}]$, i.e., $\psi \circx \sR\phi \simeq \id_{\sR X}$.
Working next with $Z = \sR\sL Y$, it follows that 
$\psi^*:[\sR X,\sR Y] \ra [\sR Y,\sR Y]$ is the inverse to the bijection
$(\sR\phi)^*:[\sR Y,\sR Y] \ra [\sR X,\sR Y]$, thus
$(\sR\phi)^*([\id_{\sR Y}]) = [\sR\phi]$ $\implies$
$\psi^*([\sR\phi]) = [\id_{\sR Y}]$ $\implies$ 
$\sR\phi \circx \psi \simeq \id_{\sR Y}$.
In other words, $\sR\phi$ has a homotopy inverse and this means that $\sR\phi$ is a weak equivalence (cf. Proposition 10) or still, $\phi$ is a weak equivalence.]\\
\endgroup%%------------------------------------<<

\begingroup%%----------------------------------->>
\fontsize{9pt}{11pt}\selectfont
The proof of Proposition 2 can be shortened by using Reedy's lifting lemma.  Thus consider the pushout square
$
\begin{tikzcd}[sep=large]
{Z}  \ar{d}[swap]{g} \ar{r}{f} &{X}  \ar{d}{\xi}\\
{Y} \ar{r}[swap]{\eta} &{P}
\end{tikzcd}
, 
$
where $f$ is a cofibration, $g$ is a weak equivalence, and 
$
\begin{cases}
\ Z\\
\ Y
\end{cases}
$
are cofibrant $-$then the claim is that $\xi$ is a weak equivalence.  First define $M_f$ by the pushout square
\begin{tikzcd}[sep=large]
{Z}  \ar{d}[swap]{i_0} \ar{r}{f} &{X}  \ar{d}\\
{IZ} \ar{r} &{M_f}
\end{tikzcd}
(cf. p. \pageref{12.28}) 
and construct a cylinder object $IX$ for $X$ with the property that the arrow $M_f \ra IX$ is an acyclic cofibration.  This done, fix a commutative diagram
\begin{tikzcd}[sep=large]
{X}  \ar{d}[swap]{\xi} \ar{r}{u} &{U}  \ar{d}{\Phi}\\
{P} \ar{r}[swap]{v} &{V}
\end{tikzcd}
(note that \mP is cofibrant).  Since $g$ is a weak equivalence, $\exists$ $\ov{w}:Y \ra U$ $\&$ $\ov{H}:IZ \ra U$ such that $\Phi \circx \ov{w} = v \circx \eta$, 
$
\begin{cases}
\ \ov{H} \circx i_0 = u \circx f\\
\ \ov{H} \circx i_1 = \ov{w} \circx g
\end{cases}
, \ 
$
and $\Phi \circx \ov{H} = v \circx \eta \circx g \circx p$, $p:IZ \ra Z$ the projection.  Choose a filler $H:IX \ra U$ for
\begin{tikzcd}[sep=large]
{M_f}  \ar{d} \ar{r}{(\ov{H},u)} &{U}  \ar{d}{\Phi}\\
{IX} \ar{r}[swap]{v \circx \xi \circx p} &{V}
\end{tikzcd}
($p:IX \ra X$) and then determine $w:P\ra U$ from the commutativity of
$
\begin{tikzcd}[sep=large]
{Z}  \ar{d}[swap]{g} \ar{r}{f} &{X}  \ar{d}{H \circx i_1}\\
{Y} \ar{r}[swap]{\ov{w}} &{U}
\end{tikzcd}
.
$
\\
\endgroup%%------------------------------------<<

\begin{proposition} \ %11
The restriction of the functor $\sL:\bC \ra \bH_\br\bC_\bc$
to $\bC_\bff$ induces a functor 
$H_\sL:\bH_\bl\bC_\bff \ra \bH\bC_{\bcf}$, 
while the restriction of the functor $\sR:\bC \ra \bH_\bl\bC_\bff$ 
to $\bC_\bc$ induces a functor 
$H_{\sR}:\bH_\br\bC_\bc \ra \bH\bC_{\bcf}$.\\
\end{proposition}

\label{13.57}
\label{16.14}
\label{16.28}
Definition: Let \bC be a model category $-$then the 
\un{homotopy category}
\index{homotopy category (model category)} 
$\bH\bC$
\index{$\bH\bC$} 
of \bC is the category whose underlying object class is the same as that of \bC, the morphism set $[X,Y]$
%%----------------------------------------------------------------------------------------------25
of X,Y being $[\sR\sL X,\sR\sL Y]$.

[Note: \  $[\sR\sL X,\sR\sL Y]$ is the morphism set of 
$H_{\sR} \circx \sL(X)$, $H_{\sR} \circx \sL(Y)$ in the category $ \bH\bC_{\bcf}$.  Of course, the situation is symmetrical in that one could just as well work with $H_\sL \circx \sR$.]

Denote by $Q$ the functor $\bC \ra \bH\bC$ which is the identity on objects and sends $f:X \ra Y$ to $H_{\sR} \circx \sL(f) = [\sR\sL f]$. \label{13.76}\\

\label{16.43} %dmc mnft
\begingroup%%----------------------------------->>
\fontsize{9pt}{11pt}\selectfont
\textbf{\small FACT} \ 
Let $f$, $g \in \Mor(X,Y)$ $-$then $\sR\sL f \simeq \sR\sL g$ iff 
$\iota_Y \circx f \circx \pi_X \simeq \iota_Y \circx g \circx \pi_X$.\\
\endgroup%%------------------------------------<<

\begin{proposition} \ %12
Let $f \in \Mor(X,Y)$ $-$then $Qf$ is an isomorphism iff $f$ is a weak equivalence.
\end{proposition}

[This follows from Proposition 10 and the fact that $f$ is a weak equivalence iff $\sR\sL f$ is a weak equivalence.]
\\

Application: Weakly equivalent objects in \bC are isomorphic in $\bH\bC$.\\

\begin{proposition} \ %13
The inclusion $\bH\bC_{\bcf} \ra \bH\bC$ is an equivalence of categories.
\end{proposition}

[The inclusion is obviously full and faithful.  
On the other hand, a given $X \in \Ob\bC$ is weakly equivalent to 
$\sR\sL X: X \overset{\pi_X}{\lla} \sL X \overset{\iota_{\sL X}}{\lra} \sR\sL X$
, thus the inclusion has a representative image.]\\

\textbf{\small LEMMA} \ 
Let \bC be a model category.  Suppose that $F:\bC \ra \bD$ is a functor which sends weak equivalences to isomorphisms $-$then
$
\begin{cases}
\ f \underset{l}{\simeq} g\\
\ f \underset{r}{\simeq} g
\end{cases}
$
$\implies$ $Ff = Fg$.

[Consider the case of left homotopy:
$
\begin{cases}
\ H \circx i_0 = f\\
\ H \circx i_1 = g
\end{cases}
$
and let $p:IX \overset{\sim}{\ra} X$ be the projection:
$
\begin{cases}
\ p \circx i_0\\
\ p \circx i_1
\end{cases}
$
= $\id_X$ $\implies$ $Fp \circx F_{i_0} = Fp \circx F_{i_1}$ 
$\implies$ 
$F_{i_0} = F_{i_1}$ $\implies$ $F f = FH \circx F_{i_0}$ $=$ $FH \circx F_{i_1}$ $=$ $F g$.]\\

\label{16.15} %dmc mnft
Given a cofibrant $X$ and a fibrant $Y$, the symbol $[X,Y]$ has two possible interpretations.  If $\Mor(X,Y)/\simeq$ is the quotient of $\Mor(X,Y)$ modulo homotopy (the meaning of $[X,Y]$ on 
p. \pageref{12.29}), 
then the lemma implies that $Q$ induces a map $\Mor(X,Y)/\simeq \ \ra \  [X,Y]$, which is in fact bijective.\\

\begingroup%%----------------------------------->>
\fontsize{9pt}{11pt}\selectfont
\textbf{\small FACT} \ 
Let $p:Y \ra Z$ be a weak equivalence, where 
$
\begin{cases}
\ Y\\
\ Z
\end{cases}
$
are fibrant $-$then for any cofibrant $X$ and any $f:X \ra Z$, $\exists$ a $g:X \ra Y$ such that $p \circx g \simeq f$, $g$ being unique up to homotopy.\\
\endgroup%%------------------------------------<<

%%----------------------------------------------------------------------------------------------26
\index{Theorem Q}
\textbf{\small THEOREM Q} \ 
Let \mS be the class of weak equivalences $-$then $S^{-1}\bC = \bHC$, 
i.e., the pair ($\bH\bC$,Q) is a localization of \bC at \mS.

[Proposition 12 implies that \mQ sends weak equivalences to isomorphisms.  
Suppose now that \bD is a metacategory and $F:\bC \ra \bD$ is a functor such that $\forall \ s \in S$, $Fs$ is an isomorphism.  
Claim: There exists a unique functor $F^\prime:\bHC \ra \bD$ such that $F = F^\prime \circx Q$.  
Thus take $F^\prime = F$ on objects and given $[f] \in [X,Y]$, represent $[f]$ by 
$\phi \in \Mor(\sR\sL X,\sR\sL Y)$ and let $F^\prime[f]$ be the filler $FX \ra FY$ in the diagram
\[
\begin{tikzcd}[sep=large]
{F\sR\sL X}  \ar{d}[swap]{F \phi} &{F\sL X}  \ar{l}[swap]{F {\iota_{\sL X}}} \ar{r}{F {\pi_X}} &{FX} \ar[dashed]{d}\\
{F\sR\sL Y}   &{F\sL Y}  \ar{l}{F {\iota_{\sL Y}}} \ar{r}[swap]{F {\pi_Y}} &{FY}
\end{tikzcd}
.]
\]
\\

Example: Let \bC be a finitely complete and finitely cocomplete category 
$-$then \bC is a model category if 
weak equivalence = isomorphism, 
cofibration = any morphism, 
fibration = any morphism and \bHC = \bC.
\label{12.13}

Example: Consider the arrow category $\bC(\rightarrow)$ of a model cateogory \bC $-$then $\bC(\rightarrow)$ can be equipped with two distinct model category structures.  Thus let 
$(\phi,\psi): (X,f,Y) \ra (X^\prime,f^\prime,Y^\prime)$ be a morphism in $\bC(\ra)$, so
\begin{tikzcd}%[sep=large]
{X}  \ar{d}[swap]{\phi} \ar{r}{f} &{Y}  \ar{d}{\psi}\\
{X^\prime} \ar{r}[swap]{f^\prime} &{Y^\prime}
\end{tikzcd}
commutes.  
In the first structure, call $(\phi,\psi)$ a weak equivalence if $\phi$ $\&$ $\psi$ are weak equivalences, 
a cofibration if $\phi$ and $X^\prime \underset{X}{\sqcup} \ra Y^\prime$ are cofibrations, 
a fibration if $\phi$ $\&$ $\psi$ are fibrations and, in the second structure, 
call $(\phi,\psi)$ a weak equivalence if $\phi$ $\&$ $\psi$ are weak equivalences, a cofibration if $\phi$ $\&$ $\psi$ are cofibrations, a fibration if $\psi$ and $X \ra X^\prime \times_{Y^\prime} Y$ are fibrations.  
The weak equivalences in either structure are the same, thus both lead to the same homotopy category 
$\bH\bC(\rightarrow)$.\\

\begingroup%%----------------------------------->>
\fontsize{9pt}{11pt}\selectfont
\textbf{\small  EXAMPLE} \ 
Take \bC = \bTOP (standard structure) $-$then \bHTOP ``is''  \bHTOP but the pointed situation is different.  
Thus let $\bTOP_{\textbf{*c}}$ be the full subcategory of $\bTOP_{*}$ whose objects are the $(X,x_0)$ such that $* \ra (X,x_0)$ is a closed cofibration, i.e., whose objects are cofibrant relative to the model category structure on 
$\bTOP_*$ inherited from \bTOP 
(cf. p. \pageref{12.30}).  
The corresponding homotopy category of $\bTOP_*$ is equivalent to $\bHTOP_{\textbf{*c}}$ (cf. Proposition 13).  
Here, the ``\bH'' has its usual interpretation since for $X$ in $\bTOP_{\textbf{*c}}$, the inclusion 
$X \vee X \ra I(X,x_0)$ is a closed cofibration, 
so a homotopy between objects in $\bTOP_{\textbf{*c}}$ preserves the base points.  
However, $\bHTOP_{\textbf{*c}}$ is not equivalent to $\bHTOP_*$ 
if this symbol is assigned its customary meaning.  
Reason: The isomorphism closure in $\bHTOP_*$ of objects in $\bTOP_{*c}$ is the class of nondegenerate spaces, therefore the inclusion
$\bHTOP_{\textbf{*c}} \ra \bHTOP_*$ does not have a representative 
%%----------------------------------------------------------------------------------------------27
image.  Of course the explanation is that the machine is rendering invertible not just pointed homotopy equivalences between pointed spaces but also homotopy equivalences between pointed spaces.
\vspi
[Note: \  $\bTOP_{\textbf{*c}}$ itself satisfies all the axioms for a model category except the first.]\\
\endgroup%%------------------------------------<<

\begingroup%%----------------------------------->>
\fontsize{9pt}{11pt}\selectfont
\textbf{\small  EXAMPLE} \ 
Take \bC = \bTOP (singular structure) $-$then \bHC is equivalent to \bHCW.\\
\endgroup%%------------------------------------<<

Let \bC be a model category.  Given a category \bD and a functor $F:\bC \ra \bD$, a 
\un{left derived functor} 
\index{left derived functor} 
for $F$ is a pair $(LF,l)$ consisting of a functor 
$LF :\bHC \ra \bD$ and a natural transformation $l:L F \circx Q \ra F$, $(LF,l)$ being final among all pairs having this property, i.e., for any pair $(F^\prime,\Xi^\prime)$ where $F^\prime \in \Ob[\bHC,\bD]$, $\&$ $\Xi^\prime \in \text{Nat}(F^\prime \circx Q,F)$, there exists a unique natural transformation $\Xi:F^\prime \ra LF$ such that $\Xi^\prime = l \circx \Xi Q$.  
Left derived functors, if they exist, are unique up to natural isomorphism.

[Note: \ A 
\un{right derived functor} 
\index{right derived functor} 
for $F$ is a pair $(RF,r)$ consisting of a functor 
$RF :\bHC \ra \bD$ and a natural transformation 
$r:F \ra R F \circx Q$, $(RF,r)$ being initial among all pairs having this property, i.e., for any pair 
$(F^\prime,\Xi^\prime)$, where 
$F^\prime \in \Ob[\bHC,\bD]$ $\&$ $\Xi^\prime \in \text{Nat}(F,F^\prime \circx Q)$, 
there exists a unique natural transformation $\Xi:RF \ra F^\prime$ such that $\Xi^\prime = \Xi Q \circx r$.]

Example: Suppose that $F:\bC \ra \bD$ sends weak equivalences to isomorphisms $-$then by Theorem Q, there exists a unique functor $F^\prime:\bHC \ra \bD$ with $F = F^\prime \circx Q$, so one can take $LF = F^\prime$ and $l = \id_F$.\\

\begingroup%%----------------------------------->>
\fontsize{9pt}{11pt}\selectfont
\textbf{\small FACT} \ 
Let 
$
\begin{cases}
\ F\\
\ G
\end{cases}
$
be functors $\bHC \ra \bD$.  Suppose that $\Xi:F \circx Q \ra G \circx Q$ is a natural transformation $-$then $\Xi$ induces a natural transformation $F \ra G$.\\
\endgroup%%------------------------------------<<

\textbf{\small LEMMA} \ 
Let \bC be a model category.  
Suppose that $F:\bC_\bc \ra \bD$ is a functor which sends acyclic cofibrations to isomorphisms $-$then 
$f \underset{r}{\simeq} g$ $\implies$ $Ff = Fg$.

[Fix a path object $PY$ for $Y$ with 
$Y \overset{\sim}{\rightarrowtail} PY \overset{\Pi}{\twoheadrightarrow} Y \times Y$ 
and a right homotopy $G:X \ra PY$ between $f$ and $g$ (cf. $\bR\bH$ (\mX is cofibrant)).  Calling $j$ the acyclic cofibration $Y \ra PY$, 
$Fj$ is an isomorphism.  Therefore
$
\begin{cases}
\ p_0 \circx j = \id_Y\\
\ p_1 \circx j = \id_Y
\end{cases}
$
$\implies$ $F_{p_0} \circx Fj = F_{p_1} \circx Fj$ $\implies$ $F_{p_0} = F_{p_1}$ $\implies$ 
$Ff = F_{p_0} \circx FG = F_{p_1} \circx FG = Fg$.]\\

\begin{proposition} \ \  %14
Let \bC be a model category.  \ Given a category \bD and a functor $F:\bC \ra \bD$, suppose that $F$ sends weak equivalences between cofibrant objects to isomorphisms $-$then a left derived functor $(LF,l)$ of $F$ exists and $\forall$ cofibrant $X$, 
$l_X:LFX \ra FX$ is an isomorphism.
\end{proposition}

[The lemma implies that $F$ induces a function $\ov{F}:\bH_{\br}\bC_{\bc} \ra \bD$.  
In addition, there is a functor
$\sL:\bC \ra \bH_{\br}\bC_{\bc}$ that takes $X$ to $\sL X$ and $f:X \ra Y$ to 
$[\sL f]_r \in [\sL X,\sL Y]_r$ 
%%----------------------------------------------------------------------------------------------28 
(cf. p. \pageref{12.31}).  
Since the composite $\ov{F} \circx \sL$  sends weak equivalences to isomorphisms, it follows from Theorem Q that there exists a unique functor $LF:\bHC \ra \bD$ such that $LF \circx Q = \ov{F} \circx \sL$.  
Define a natural transformation $l:LF \circx Q \ra F$ by assigning to each $X \in \Ob\bC$ 
the element $l_X = F\pi_X \in \Mor(F\sL X,FX)$ $-$then $X$ cofibrant 
$
\implies
\begin{cases}
\ \sL X = X\\
\ \pi_X = \id_X
\end{cases}
\implies l_X = F\id_X = \id_{FX}.
$
It remains to prove that the pair $(LF,l)$ is final.  So fix a pair $(F^\prime,\Xi^\prime)$ as above.  
Define a natural transformation
$\Xi:F^\prime \ra LF$ by assigning to each $X \in \Ob\bHC$ the element 
$\Xi_X \in \Mor(F^\prime X,LFX)$ determined from 
\begin{tikzcd}[sep=large]
{F^\prime X}   \ar{rr}{F^\prime (Q \pi_X)^{-1}} 
&&{F^\prime \sL X }  \ar{rr}{\Xi_{\sL X}^\prime}
&&{F \sL X = LFX.}
\end{tikzcd}
Bearing in mind that $\forall \ X$, $QX = X$ and $\sL X$ is cofibrant, the commutativity of
\[
\begin{tikzcd}[sep=large]
{F^\prime \sL X}  \ar{d}[swap]{F^\prime Q \pi_X} \ar{r}{\Xi_{\sL X}} 
&{L F \sL X}  \ar[equals]{d}  \ar{r}{l_{\sL X}} 
&{F \sL X}  \ar{d}{F\pi_X}\\
{F^\prime X} \ar{r}[swap]{\Xi_X} &{L F X} \ar{r}[swap]{l_X} &{FX}
\end{tikzcd}
\]
ensures the uniqueness of $\Xi$.]\\

Given model categories 
$
\begin{cases}
\ \bC\\
\ \bD
\end{cases}
$
and a functor $F:\bC \ra \bD$, a 
\un{total left derived functor}
\index{total left derived functor} 
for $F$ is a functor
$\bL F:\bHC \ra \bHD$ which is a left derived functor for the composite $Q \circx F: \bC \ra \bHD$.  
Total left derived functors, if they exist, are unique up to natural isomorphism.

[Note: A 
\un{total right derived functor}
\index{total right derived functor} 
for $F$ is a functor 
$\bR F:\bHC \ra \bHD$ which is a right derived functor for the composite $Q \circx F: \bC \ra \bHD$.]

Remark: The substitute for the failure of
\begin{tikzcd}%[sep=large]
{\bC}  \ar{d} \ar{r}{F} &{\bD}  \ar{d}\\
{\bHC} \ar{r}[swap]{\bL F} &{\bHD}
\end{tikzcd}
to commute is the natural transformation $l:\bL F \circx Q \ra Q \circx F$.

Example: Suppose that $F:\bC \ra \bD$ sends weak equivalences between cofibrant objects to weak equivalences $-$then by Proposition 14, $\bL F$ exists and $\forall$ cofibrant 
$X$, $l_X:\bL FX \ra FX$ is an isomorphism.\\

\label{13.94}
\label{16.49}
\label{16.50}

\textbf{\small LEMMA} \ 
Let $F:\bC \ra \bD$ be a functor between model categories.  
Suppose that $F$ sends acyclic cofibrations between cofibrant objects to weak equivalences $-$then $F$ preserves weak equivalences between cofibrant objects.

[Let $f:X \ra Y$ be a weak equivalence, where $X \ \& \ Y$ are cofibrant.  Factor 
$f \coprod \id_Y: X \coprod Y \ra Y$ as $p \circx i$, where $i: X \coprod Y \ra Z$ is a cofibration and $p:Z \ra Y$ is an acyclic fibration.  Since $X \ \& \ Y$ are cofibrant, the composites
$
\begin{cases}
\ i \circx \text{in}_0:X \ra Z\\
\ i \circx \text{in}_1:Y \ra Z
\end{cases}
$
are cofibrations.  In addition
$
\begin{cases}
\ p \circx i \circx \text{in}_0\\
\ p \circx i \circx \text{in}_1
\end{cases}
$
are weak equivalences, hence
$
\begin{cases}
\ i \circx \text{in}_0\\
\ i \circx \text{in}_1
\end{cases}
$
are weak equivalences.
%%----------------------------------------------------------------------------------------------29
Therefore
$
\begin{cases}
\ F( i \circx \text{in}_0)\\
\ F(i \circx \text{in}_1)
\end{cases}
$
are weak equivalences.  But $Fp \circx F(i \circx \text{in}_1) = \id_{FY}$, thus $Fp$ is a weak equivalence and so $Ff = Fp \circx F( i \circx \text{in}_0)$ is a  weak equivalence.]\\

\label{13.59a}
\label{13.83}
\label{16.39}
\index{Theorem TDF Theorem}
\textbf{\small TDF THEOREM} \ 
Let \bC and \bD be model categories.  Suppose that 
$
\begin{cases}
\ F:\bC\ra \bD\\
\ G:\bD\ra \bC
\end{cases}
$
are functors and $(F,G)$ is an adjoint pair.  Assume: $F$ preserves cofibrations and $G$ preserves fibrations $-$then
$
\begin{cases}
\ \bL F:\bHC\ra \bHD\\
\ \bR G:\bHD\ra \bHC
\end{cases}
$
exist and $(\bL F,\bR G)$ is an adjoint pair.\\

[The existence of $\bL F$ follows from the fact that $F$ preserves acyclic cofibrations 
(cf. p. \pageref{12.32} ff.), 
thus by the lemma, $F$ preserves weak equivalences between cofibrant objects, and Proposition 14 is applicable (the argument for $\bR G$ is dual).  
Because $F$ is a left adjoint and $G$ is a right adjoint, $F$ preserves initial objects and $G$ preserves final objects.  
Therefore $F$ sends cofibrant objects to cofibrant objects and $G$ sends fibrant objects to fibrant objects.  
Consider now the bijection of adjunction $\Xi_{X,Y}:\Mor(FX,Y) \ra \Mor(X,GY)$ 
(cf. p. \pageref{12.33}).  If
$
\begin{cases}
\ X \in \Ob\bC_{\bc}\\
\ Y \in \Ob\bD_\bff
\end{cases}
, \ 
$
then $\Xi_{X,Y}$ respects the relation of homotopy and induces a bijection $[FX,Y] \ra [X,GY]$.  
Using the definitions, for arbitrary
$
\begin{cases}
\ X \in \Ob\bC\\
\ Y \in \Ob\bD
\end{cases}
$
\
this leads to functorial bijections 
$[\bL FX,Y] \approx$ $[F\sL X,\sR Y] \approx$ $[\sL X,G \sR Y] \approx$ $[X,\bR GY]$.]
\label{16.18}
\label{16.31}

[Note: \  Suppose that $\forall$ 
$
\begin{cases}
\ X \in \Ob\bC_{\bc}\\
\ Y \in \Ob\bD_\bff
\end{cases}
, \ 
$
$\Xi_{X,Y}$ maps the weak equivalences in $\Mor(FX,Y)$ onto weak equivalences in $\Mor(X,GY)$ $-$then the pair 
$(\bL F,\textbf{R}G)$ is an adjoint equivalence of categories.]\\

Implicit in the proof of the TDF theorem is the fact that $\forall$ $X$, $\bL FX$ is isomorphic 
(in \bHD) to $FX^\prime$, where $X^\prime$ is any cofibrant object which is weakly equivalent to $X$.\\

\begingroup%%----------------------------------->>
\fontsize{9pt}{11pt}\selectfont
\index{Pushouts (Example)}
\textbf{\small EXAMPLE (\un{Pushouts})} \quadx
Fix  a model category \bC.  Let \bI be the category 
$1 \ \bullet \overset{a}{\leftarrow} \underset{3}{\bullet} \overset{b}{\ra} \bullet \ 2$ 
(cf. p. \pageref{12.34}) 
$-$then the functor category [\bI,\bC] is again a model category 
(cf. p. \pageref{12.35} ff.).  
Given a 2-source
$X \overset{f}{\leftarrow} Z \overset{g}{\ra} {Y},$
define $P$ by the pushout square
\begin{tikzcd}[sep=large]
{Z}  \ar{d}[swap]{f} \ar{r}{g} &{Y}  \ar{d}{\eta}\\
{X} \ar{r}[swap]{\xi} &{P}
\end{tikzcd}
and put colim$(X \overset{f}{\leftarrow} Z \overset{g}{\ra} {Y}) = P$ to get a functor
colim: $[\bI,\bC] \ra \bC$ which is left adjoint to the constant diagram functor 
$K:\bC \ra [\bI,\bC]$.  
Since $K$ preserves fibrations and acyclic fibrations, the hypotheses of the TDF theorem are satisified 
(cf. p. \pageref{12.36} ff.).  
Therefore $\bL \colim$ and $\bR K$ exist and $(\bL \colim, \bR K )$ is an adjoint pair.  
Moreover, according to the theory, 
$\bL\colimx(X \overset{f}{\leftarrow} Z \overset{g}{\ra} {Y})$ is isomorphic 
(in \bHC) to $\colimx(X \overset{f}{\leftarrow} Z \overset{g}{\ra} {Y})$ 
whenever 
$(X \overset{f}{\leftarrow} Z \overset{g}{\ra} {Y})$ 
is cofibrant, i.e., whenever $Z$ is cofibrant and 
$
\begin{cases}
\ f:Z \ra X\\
\ g:Z \ra Y
\end{cases}
$
are cofibrations.  
For instance, by way of illustration, let us take \bC = \bTOP (standard structure).  
Claim: 
$\bL \colim(X \overset{f}{\leftarrow} Z \overset{g}{\ra} {Y})$
and 
$M_{f,g}$
have the same homotopy 
%%----------------------------------------------------------------------------------------------30
type.  To see this, consider the 2-source $M_f \leftarrow Z \ra M_g$.  It is cofibrant and the vertical arrows in the commutative diagram
\begin{tikzcd}[sep=large]
{M_f}  \ar{d}   &{Z} \arrow[d,shift right=0.5,dash] \arrow[d,shift right=-0.5,dash] \ar{l} \ar{r} &{M_g}  \ar{d}\\
{X}  &{Z}  \ar{l} \ar{r} &{Y} 
\end{tikzcd}
are homotopy equivalences 
(but $M_f \leftarrow Z \ra M_g$ is not
$\sL(X \overset{f}{\lla} Z \overset{g}{\lra} Y)$), so 
$\bL\colim(X \overset{f}{\lla} Z \overset{g}{\lra} Y)$ $\approx$ 
$\colim(M_f \la Z \ra M_g$) $\approx$ $M_{f,g}$ 
(cf. p. \pageref{12.37}).
\vspi
[Note: \  The story for pullbacks is analogous (work with $\bR\lim$).]\\
\endgroup%%------------------------------------<<

\label{13.97}
\begingroup%%----------------------------------->>
\fontsize{9pt}{11pt}\selectfont
\textbf{\small  EXAMPLE} \ 
Fix  a model category \bC $-$then $\bFIL(\bC)$ is again a model category 
(cf. p. \pageref{12.38}).  
Assuming that \bC admits sequential colimits, there is a functor
$\colim:\bFIL(\bC) \ra \bC$ which is left adjoint to the constant diagram functor 
$K:\bC \ra \bFIL(\bC)$.  
Since $K$ preserves fibrations and acyclic fibrations, the hypotheses of the TDF theorem are satisified 
(cd. p. \pageref{12.39} ff.).  
Therefore $\bL \colim$ and $\bR K$ exist and $(\bL \colim,\bR K)$ is an adjoint pair.  
Moreover, according to the theory, 
$\bL\colimx(\bX,\bff)$ is isomorphic (in \bHC) to $\colimx(\bX,\bff)$ whenever $(\bX,\bff)$ is cofibrant, i.e., whenever $X_0$ is cofibrant and $\forall \ n$, $f_n:X_n \ra X_{n+1}$ is a cofibration.  
If $\bC  = \bTOP$  (standard structure), $\bL \colimx(\bX,\bff)$ and $\tel(\bX,\bff)$ have the same homotopy type 
(cf. p. \pageref{12.40}).  
In general, $\colim:\bFIL(\bC) \ra \bC$ preserves weak equivalences between cofibrant objects, a fact which specialized to the topological setting recovers Proposition 15 in $\S 3$ provided that the cofibrations are closed.
\vspi
[Note: \  The story for $\bTOW(\bC)$ is analogous (work with $\bR \lim$).]\\
\endgroup%%------------------------------------<<

The axioms defining a model category interlock cofibrations and fibrations in such a way that certain canonical examples are excluded.  This difficulty can be circumvented by simply weakening the assumptions and concentrating on either the cofibrations or the fibrations.

Consider a category \bC equipped with two composition closed classes of morphisms termed 
\un{weak equivalences} 
\index{weak equivalences (cofibration category)} 
(denoted $\overset{\sim}{\ra}$) and 
\un{cofibrations}
\index{cofibration (cofibration category)}   
(denoted $\rightarrowtail$), each containing the isomorphisms of \bC.  
Agreeing to call a morphism which is both a weak equivalence and a cofibration an 
\un{acyclic cofibration}, 
\index{acyclic cofibration (cofibration category)}  \bC is said to be a 
\un{cofibration category} 
\index{cofibration category}  
provided that the following axioms are satisfied.

\qquad (CC-1) \quadx \bC has an initial object $\emptyset$.

\qquad (CC-2) \quadx Given composable morphisms $f$, $g$, if any two of $f$, $g$, $g \circx f$ are weak equivalences, so is the third.

\qquad (CC-3) \quadx Every 2-source 
$X \overset{f}{\leftarrow} Z \overset{g}{\ra} {Y}$, 
where $f$ is a cofibration (acyclic cofibration), admits a pushout 
$X \overset{\xi}{\ra} Z \overset{\eta}{\leftarrow} {Y},$ 
where $\eta$ is a cofibration (acyclic cofibration).

\qquad (CC-4) \quadx Every morphism can be written as the composite of a cofibration and a weak equivalence.

[Note: \  The axioms defining a 
\un{fibration category} 
\index{fibration category} 
are dual.]\\

%%----------------------------------------------------------------------------------------------31
Let \bC be a cofibration category $-$then an $X \in $ Ob\bC is said to be 
\un{cofibrant} 
\index{cofibrant} if $\emptyset \ra X$ is a cofibration and 
\un{fibrant} 
\index{fibrant} 
if every acyclic cofibration $X \ra Y$ has a left inverse 
(cf. p. \pageref{12.41}).

\index{Fibrant Embedding Axiom} \index{FEA}
\indent\indent (Fibrant Embedding Axiom) (FEA) \quadx Given an object $X$ in \bC, there is an acyclic cofibration 
$\iota_X:X \ra \sR X$, where $\sR X$ is fibrant.

[Note: \  The FEA is trivially met if all objects are fibrant.]

Example: The cofibrant objects in a model category are the object class of a cofibration category satisyfing the FEA.\\

\begingroup%%----------------------------------->>
\fontsize{9pt}{11pt}\selectfont
\textbf{\small  EXAMPLE} \ 
Take \bC = \bTOP $-$then \bTOP is a cofibration category if weak equivalence = homotopy equivalence, cofibration = cofibration.  All objects are cofibrant and fibrant.\\
\endgroup%%------------------------------------<<

\begingroup%%----------------------------------->>
\fontsize{9pt}{11pt}\selectfont
\textbf{\small  EXAMPLE} \ 
Take \bC = $\bTOP_*$ $-$then $\bTOP_*$ is a cofibration category if weak equivalence = pointed homotopy equivalence, cofibration = pointed cofibration.  All objects are cofibrant and fibrant.
\vspi
[Note: \  This is the ``internal'' structure of a cofibration category on $\bTOP_*$.  
An ``external'' structure is obtained by letting the weak equivalences be the pointed maps which are homotopy equivalences in \bTOP and the cofibrations be the pointed maps which are cofibrations in \bTOP.  
Here, all objects are fibrant and the cofibrant objects are the wellpointed spaces.  
Another ``external' structure arises by requiring that the cofibrations be closed, 
which reduces the number of cofibrant objects.]\\
\endgroup%%------------------------------------<<

\begingroup%%----------------------------------->>
\fontsize{9pt}{11pt}\selectfont
\textbf{\small  EXAMPLE} \ 
Take for \bC the category whose objects are pairs $(X,N_X)$, where $X$ is a pointed connected CW space and $N_X$ is a perfect normal subgroup of $\pi_1(X)$, and whose morphisms $f:(X,N_X) \ra (Y,N_Y)$ are pointed continuous functions $f:X \ra Y$ such that $f_*(N_X) \subset N_Y$.  
Stipulate that $f$ is a weak equivalence if $f_*:\pi_1(X)/N_X \approx \pi_1(Y)/N_Y$ and 
$f_*:H_*(X;f^*\sG) \approx H_*(Y;\sG)$ for every locally constant coefficient system $\sG$ on $Y$ arising from a $\pi_1(Y)/N_Y$-module.  If by cofibration one understands a pointed continuous function which is a closed cofibration in \bTOP, then \bC is a cofibration category satisfying the FEA.
\vspi
[CC-1, CC-2, and CC-4 are clear.  As for CC-3, given a 2-source $X \overset{f}{\leftarrow} Z \overset{g}{\ra} Y$, where $f$ is a cofibration, define \mP by the pushout square
\begin{tikzcd}[sep=large]
{Z}  \ar{d}[swap]{f} \ar{r}{g} &{Y}  \ar{d}{\eta}\\
{X} \ar{r}[swap]{\xi} &{P}
\end{tikzcd}
and let $N_P$ be the normal subgroup of $\pi_1(P) = \pi_1(X) *_{\pi_1(Z)} \pi_1(Y)$ generated by $N_X$ $\&$ $N_Y$.  To check the FEA assertion, fix a pair $(X,N_X)$.  
Thanks to the plus construction, there is a pair $(X_{N_X}^+,0)$ and a cofibration
$(X,N_X) \ra (X_{N_X}^+,0)$ which is a weak equivalence (cf. $\S 5$, Proposition 22).  
Claim: $(X_{N_X}^+,0)$ is fibrant.  
For suppose given $(X_{N_X}^+,0) \overset{\sim}{\rightarrowtail} (Y,N_Y)$.  
Denote by $f$ the composite
$(X_{N_X}^+,0) \overset{\sim}{\rightarrowtail}$ $(Y,N_Y) \overset{\sim}{\rightarrowtail} (Y_{N_Y}^+,0)$, 
so 
$f_*:\pi_1(X_{N_X}^+) \approx$ $\pi_1(Y)/N_Y \approx$ $\pi_1(Y_{N_Y}^+)$.
Since $f$ is acyclic (as a map) and a cofibration, one may now invoke $\S 5$, Proposition 19 and $\S 3$, Proposition 5.]\\
\endgroup%%------------------------------------<<

%%----------------------------------------------------------------------------------------------32
\begingroup%%----------------------------------->>
\fontsize{9pt}{11pt}\selectfont
\textbf{\small  EXAMPLE} \ 
Take for \bC the category whose objects are the pointed connected CW spaces.  
Fix an abelian group \mG $-$then \bC = $\bCONCWSP_*$ is a cofibration category if weak equivalence = $HG$-equivalence, cofibration = closed cofibration in \bTOP and this structure satisfies the FEA.
\vspi
[Note: \  The fibrant objects are the $HG$-local spaces.]\\
\endgroup%%------------------------------------<<

\label{15.7}
The formal ``one sided'' results in model category theory carry over to cofibration categories, e.g., Propositions 2, 3, and 4.  Assuming in addition that \bC satisfies the FEA, one can also show that the inclusion 
$\bHC_{\bcf} \ra \bHC$ is an equivalence of categories (cf. Proposition 13) and $S^{-1}$\bC = \bHC, where $S$ is the class of weak equivalences (cf. Theorem Q).\\

\begingroup%%----------------------------------->>
\fontsize{9pt}{11pt}\selectfont
\textbf{\small  EXAMPLE} \ 
Take \bC = $\bTOP_*$ $-$then $\bHTOP_*$ ``is'' $\bHTOP_*$ if $\bTOP_*$ carries its ``internal'' structure of a cofibration category.\\
\endgroup%%------------------------------------<<

\begingroup%%----------------------------------->>
\fontsize{9pt}{11pt}\selectfont
\textbf{\small  EXAMPLE} \ 
The homotopy category of the cofibration category evolving from the plus construction is equivalent to $\bHCONCWSP_*$.\\
\endgroup%%------------------------------------<<

Let \bC be a category.  Suppose given a composition closed class $S \subset \Mor \bC$ containing the isomorphisms of \bC such that for composable morphisms $f$, $g$, if any two of $f$, $g$, $g \circx f$ are in $S$, so is the third.  
Problem: Does $S^{-1}$\bC exist as a category?  
The assumption that $S$ admits a calculus of left or right fractions does not suffice to resolve the issue.  
However, one strategy that will work is to somehow place on \bC the structure of a model category (or a cofibration category) in which $S$ appears as the class of weak equivalences.  
For then $S^{-1}$\bC ``is'' \bHC and \bHC is a category.\\

\label{13.105}
\label{18.1}
\begingroup%%----------------------------------->>
\fontsize{9pt}{11pt}\selectfont
\textbf{\small  EXAMPLE} \ 
Let \bC be a model category.  Assume: \bC is complete and cocomplete.  Suppose that \bI is a small category and let 
$S \subset  \Mor[\bI,\bC]$ be the class of levelwise weak equivalences $-$then it has been shown by 
%Dwyer-Kan\footnote[2]{Model Categories and General Abstract Homotopy Theory,???}
Dwyer-Hirschhorn-Kan-Smith\footnote[2]{\textit{Homotopy limit functors on model categories, and homotopical categories, Mathematical Surveys and Monographs}, Amer. Math Soc. \textbf{113} (2004).}
that 
$S^{-1}[\bI,\bC]$ exists as a category even though $[\bI,\bC]$ need not carry the structure of a model category having 
$S$ for its class of weak equivalences.
\vspi
[Note: \  Given a functor $[\bI,\bC] \ra \bC$ or $\bC \ra [\bI,\bC]$, one can define in the obvious way its total left (right) derived functor.  
In particular:  
$\colim: [\bI,\bC] \ra \bC$ ($\lim : [\bI,\bC] \ra \bC$) is a left (right) adjoint for the constant diagram functor 
$K:\bC \ra [\bI,\bC]$.  
Moreover, $\bL\colim$ and $\bR K$ ($\bL K$ and $\bR\lim$) exist and 
$(\bL\colim,\bR K)$ 
($(\bL K,\bR\lim)$) is an adjoint pair (Dwyer-Hirschhorn-Kan-Smith (ibid.)).]\\
\endgroup%%------------------------------------<<

%%%%%%%%%%%%%%%%%%%%%%%%%%%%%%%%%%%%%%
%%%%%%%%%%%%%%%%%%%%%%%%%%%%%%%%%%%%%%
%%%%%%%%%%%%%%%%%%%%%%%%%%%%%%%%%%%%%%

\begin{center}
$\S \ 12$
\\[0.5cm]
$\mathcal{REFERENCES}$\\[-.2cm]
\end{center}

\[
\textbf{BOOKS}
\]

\begingroup
\fontsize{9pt}{11pt}\selectfont
\setlength\parindent{0 cm}

[1] \quad Baues, H., \textit{Algebraic Homotopy}, Cambridge University Press (1989).
\\[-.2cm]

[2] \quad Gunnarsson, T., \textit{Abstract Homotopy Theory and Related Topics}, Ph.D. Thesis, Chalmers Tekniska 

\hspace{0.8cm}H\"ogskola, 
G\"oteborg (1978).
\\[-.2cm]

[3] \quad Quillen, D., \textit{Homotopical Algebra}, Springer Verlag (1967).
\\[-.2cm]

[4] \quad Tanr\'e, D., \textit{Homotopie Rationnelle: Mod\`eles de Chen, Quillen, Sullivan}, Springer Verlag (1983).
\\[-.2cm]

[5] \quad Warner,G., \textit{Categorical Homotopy Theory}, Preprint. 

\hspace{0.8cm}\url{https://sites.math.washington.edu//~warner/CHT_Warner.pdf}

\endgroup

\[
\textbf{ARTICLES}
\]

\begingroup
\fontsize{9pt}{11pt}\selectfont
\setlength\parindent{0 cm}

[1] \quad Anderson, D., Fibrations and Geometric Realizations, 
\textit{Bull. Amer. Math. Soc.} \textbf{84} (1978), 765-788.
\\[-.2cm]

[2] \quad Avramov, L., Local Algebra and Rational Homotopy, 
\textit{Ast\'erisque} \textbf{113-114} (1984), 15-43.
\\[-.2cm]

[3] \quad Avramov, L. and Halperin, S., Through the Looking Glass: A Dictionary Between Rational Homotopy 

\hspace{0.8cm}Theory and Local Algebra, 
\textit{SLN} \textbf{1183} (1986), 1-27.
\\[-.2cm]

[4] \quad Bousfield, A. and Gugenheim, V., On PL DeRham Theory and Rational Homotopy Type, 
\textit{Memoirs}

\hspace{0.8cm}\textit{Amer. Math. Soc.} \textbf{179} (1976), 1-94.
\\[-.2cm]

%[5] \quad Casacuberta, C., Paricio, L., and Rodriguez, J., Models for Torsion Homotopy Types,\\[-.2cm]
[5] \quad Casacuberta, C., Paricio, L., and Rodriguez, J., Models for Torsion Homotopy Types, 
\textit{Isr J. Math.} 

\hspace{0.8cm}\textbf{107} (1998), 301-318.
\\[-.2cm]
%www.ub.edu.topologia.casacuberta.articles.crhp.pdf %(verify)\\[-.2cm]

[6] \quad Doeraene, J-P., L.S.-Category in a Model Category, 
\textit{J. Pure Appl. Algebra} \textbf{84} (1993), 215-261.
\\[-.2cm]

[7] \quad Dold, A., Zur Homotopietheorie der Kettenkomplexe, 
\textit{Math. Ann.} \textbf{140} (1960), 278-298.
\\[-.2cm]

[8] \quad Dupont, N., Problems and Conjectures in Rational Homotopy Theory, 
\textit{Expo. Math.} \textbf{12} (1994), 323-

\hspace{0.8cm}352.
\\[-.2cm]

[9] \quad Dwyer, W. and Spalinski, J., Homotopy Theories and Model Categories, In: 
\textit{Handbook of Algebraic}

\hspace{0.8cm}\textit{Topology}, I. James (ed.), North Holland (1995), 73-126.
\\[-.2cm]

[10] \quad Halperin, S., Lectures on Minimal Models, 
\textit{M\'emoires  Soc. Math. France} \textbf{111} (1983), 1-261.
\\[-.2cm]

[11] \quad Jardine, J., Homotopy and Homotopical Algebra, In: 
\textit{Handbook of Algebra}, M. Hazewinkel (ed.), 

\hspace{0.95cm}North Holland (1996), 639-669.
\\[-.2cm]

[12] \quad Lehmann, D., Th\'eorie Homotopique des Formes Diff\'erentielles, 
\textit{Ast\'erisque} \textbf{45} (1990/1977), 1-145.
\\[-.2cm]

[13] \quad Quillen, D., Rational Homotopy Theory, 
\textit{Ann. of Math.} \textbf{90} (1969), 205-295.
\\[-.2cm]

[14] \quad Reedy, C., Homotopy Theory of Model Categories, \textit{Preprint}.
\\[-.2cm]

[15] \quad Roig, A., Mod\`eles Minimaux et Foncteurs D\'eriv\'es, 
\textit{J. Pure Appl. Algebra} \textbf{91} (1994), 231-254.
\\[-.2cm]

[16] \quad Sullivan, D., Infinitesimal Computations in Topology, 
\textit{Publ. Math. I.H.E.S.} \textbf{47} (1977), 269-331.

\setlength\parindent{2em}

\endgroup

\chapter{
$\boldsymbol{\S}$\textbf{13}.\quadx  SIMPLICIAL SETS}
\setlength\parindent{2em}
\setcounter{proposition}{0}
%%----------------------------------------------------------------------------------------------01
$\text{ }$\\[-1.25cm]

It is possible to develop much of algebraic topology entirely within the context of simplicial sets.  
However, I shall not go down that road.  
Instead, the focus will be on the simplicial aspects of model categories which, for instance, is the homotopical basis for the algebraic K-theory of rings or spaces.

\bSISET (= \textbf{$\widehat{\bDelta}$}) is complete and cocomplete, wellpowered and cowellpowered, and cartesian closed 
(cf. \pageref{13.1}).

\label{13.106}
[Note: \  \bSISET admits an involution $X \ra X^\OP$, where $d_i^\OP = d_{n-i}$, $s_i^\OP = s_{n-i}$.  Example: $\forall$ small category \bC, 
$\ner\bC^\OP = (\ner\bC)^\OP$.]

Notation: $\emptyset$ stands for an initial object in \bSISET (e.g., $\ddz$) and $*$ stands for a final object in \bSISET 
(e.g., $\dz$).\\

\begingroup%%----------------------------------->>
\fontsize{9pt}{11pt}\selectfont
The four exponential objects associated with $\emptyset$ and $*$ are 
$\emptyset^\emptyset = *$, 
$*^\emptyset = *$, 
$\emptyset^* = \emptyset$, 
$*^* = *$.\\
\endgroup%%------------------------------------<< 

Let $X$ be a simplicial set $-$then $\abs{X}$ is a CW complex 
(cf. p. \pageref{13.2}), 
thus is a compactly generated Hausdorff space.  Therefore ``geometric realization'' can be viewed as a functor $\bSISET \ra \bCGH$.\\

\begingroup%%----------------------------------->>
\fontsize{9pt}{11pt}\selectfont
$\abs{?}:\bSISET \ra \bTOP$ preserves colimits (being a left adjoint) and it is immaterial whether the colimit is taken in \bTOP or \bCGH.  Reason: A colimit in \bCGH is calculated by taking the maximal Hausdorff quotient of the colimit calculated in \bTOP.\\
\endgroup%%------------------------------------<<

\begingroup%%----------------------------------->>
\fontsize{9pt}{11pt}\selectfont
\textbf{\small EXAMPLE} \ 
The pushout square 
\begin{tikzcd}[sep=large]
{\ddn} \ar{d} \ar{r} &{\Delta[0]} \ar{d}\\
{\Delta[n]} \ar{r} &{\bS[n]}
\end{tikzcd}
defines the 
\un{simplicial $n$-sphere}
\index{simplicial $n$-sphere} $\bS[n]$.  
Its geometric realization is homeomorphic to $\bS^n$.\\
\endgroup%%------------------------------------<<

\textbf{\small LEMMA} \ 
$\abs{?}:\bSISET \ra \bCGH$ preserves equalizers.

[Let $X$ and $Y$ be simplicial sets; let $u, \ v:X \ra Y$ be a pair of simplicial maps $-$then 
$Z = \eq(u,v)$ is a simplicial subset of $X$ and $\abs{Z}$ is a subcomplex of $\abs{X}$ which is contained in 
$\eq(\abs{u},\abs{v})$.  
Take now a point $[x,t] \in \eq(\abs{u},\abs{v})$, say $x \in X_n^{\#}$ $\&$ 
$t \in \mdpn$ 
(cf. p. \pageref{13.3}).  
Write
$
\begin{cases}
\ u(x) = (Y\alpha)y_u\\
\ v(x) = (Y\beta)y_v
\end{cases}
,
$
where $y_u$, $y_v \in Y$ are nondegenerate and $\alpha, \ \beta \in \Mor\bDelta$ are epimorphisms.  
By assumption, $\abs{u}([x,t]) = \abs{v}([x,t]) $; moreover, 
$
\begin{cases}
\ \abs{u}([x,t]) = [u(x),t] = [(Y\alpha)y_u,t] = [y_u,\Delta^\alpha(t)]\\
\ \abs{v}([x,t]) = [v(x),t] = [(Y\beta)y_v,t] = [y_v,\Delta^\beta(t)]
\end{cases}
,
$
so $y_u = y_v$ and $\Delta^\alpha(t) = \Delta^\beta(t)$ (because the issue is one of
%%----------------------------------------------------------------------------------------------02
\label{13.17}
\label{13.20}
epimorphisms, interior points go to interior points).  But $\Delta^\alpha(t) = \Delta^\beta(t)$ $\implies$ $\alpha = \beta$, hence $u(x) = v(x)$ or still, $x \in Z$ $\implies$ $[x,t] \in \abs{Z}$.]\\

\label{13.67}
\textbf{\small LEMMA} \ 
$\abs{?}:\bSISET \ra \bCGH$ preserves finite products.

[Let $X$ and $Y$ be simplicial sets.  Write
$
\begin{cases}
\ X = \colim_i \Delta[m_i]\\
\ Y = \colim_j \Delta[n_j]
\end{cases}
$
(cf. p. \pageref{13.4}).  
Since \bSISET is cartesian closed, products commute with colimits.  
Therefore 
$\abs{X \times Y}$ $\approx$ \\
$\abs{\colim_{i,j} \Delta[m_i] \times \Delta[n_j]}$ 
from which 
$\abs{X \times Y} \approx$ 
$\colim_{i,j}\abs{\Delta[m_i] \times \Delta[n_j]} \approx$ 
$\colim_{i,j}(\abs{\Delta[m_i]} \times_k \abs{\Delta[n_j]})$, 
the arrow 
$\abs{\Delta[m_i] \times \Delta[n_j]} \ra$  
$\abs{\Delta[m_i]} \times \abs{\Delta[n_j]} \equiv$ 
$\abs{\Delta[m_i]} \times_k \abs{\Delta[n_j]}$
being a homeomorphism 
(cf. p. \pageref{13.5}).  
But \bCGH is also cartesian closed 
(cf. p. \pageref{13.6}), 
thus once again products commute with colimits.  This gives 
$\abs{X \times Y} \approx$ 
$\colim_i\abs{\Delta[m_i]} \times_k \colim_j\abs{\Delta[n_j]} \approx$ 
$\abs{X} \times_k \abs{Y}$, 
i.e., the arrow
$\abs{X \times Y} \ra$ $\abs{X} \times_k \abs{Y}$ is a homeomorphism.]

[Note: \  While the arrow 
$\abs{X \times Y} \ra$ $\abs{X} \times \abs{Y}$ is a set theoretic bijection, it need not be a homeomorphism when 
$\abs{X} \times \abs{Y}$ has the product topology.]\\

\begin{proposition} \ 
$\abs{?}:\bSISET \ra \bCGH$ preserves finite limits.
\end{proposition}

[This is implied by the lemmas.]

[Note: \  $\abs{?}:\bSISET \ra \bCGH$ does not preserve arbitrary limits.  
Example: The arrow $\abs{\dw^\omega} \ra \abs{\dw}^\omega$ is not a homeomorphism.]\\

\label{13.17}
\label{13.20}
Example: The composite $\abs{?} \circx \sin$ preserves homotopies 
($f \simeq g $ $\implies$ $\abs{\sin f} \simeq \abs{\sin g}$).

[For any topological space $X$, 
$\abs{\sin X} \times \Delta^1 \approx$ 
$\abs{\sin X} \times \abs{\dw}$ $\approx$ 
$\abs{\sin X \times \dw}$ $\lra$ \\
$\abs{\sin X \times \sin \abs{\dw}}$ $\approx$ 
$\abs{\sin (X \times \Delta^1)}$, 
$\ra$ being the geometric realization of $\id_{\sin X}$ times the arrow of adjunction $\dw \ra \sin\abs{\dw}$.  
So, if $H:X \times \Delta^1 \ra Y$ is a homotopy, then 
$\abs{\sin X} \times \Delta^1 \lra $
\begin{tikzcd}%[sep=small]
%{\abs{\sin X} \times \Delta^1} \ar{r} 
{\abs{\sin (X \times \Delta^1)}} \ar{r}{\abs{\sin H}}
&{\abs{\sin Y}}
\end{tikzcd}
%$\abs{\sin X} \times \Delta^1 \ra \abs{\sin (X \times \Delta^1)} \overset{\abs{\sin Y}}{\longrightarrow} \abs{\sin Y}$ 
is a homotopy.]\\

\begingroup%%----------------------------------->>
\fontsize{9pt}{11pt}\selectfont
\label{5.0aj}
\textbf{\small EXAMPLE} \ 
Let $G$ be a simplicial group $-$then $\abs{G}$ is a compactly generated group.
\vspi
[Note: \  $\abs{G}$ is a topological group if $\abs{G}$ is countable, i.e., if 
$\forall \ n$, $\#(G_n^\#) \leq \omega$.]\\
\endgroup%%------------------------------------<<

\begingroup%%----------------------------------->>
\fontsize{9pt}{11pt}\selectfont
\textbf{\small FACT} \ 
Let $X$ and $Y$ be simplicial sets, $\Pi X$ and $\Pi Y$ their fundamental groupoids $-$then 
$\Pi(X \times Y) \approx$ $\Pi X \times \Pi Y$.
\vspi
[Note: \  The functor $\Pi:\bSISET \ra \bGRD$ does not preserve equalizers.  
Example: Define $X$ by the pushout square
\begin{tikzcd}[sep=large]
{\dot\Delta[2]} \ar{d} \ar{r} &{\Delta[2]} \ar{d}{v}\\
{\Delta[2]} \ar{r}[swap]{u} &{X}
\end{tikzcd}
: $\Pi\dot\Delta[2] = \Pi \eq(u,v) \neq \eq(\Pi u,\Pi v) = \Pi\Delta[2]$.]\\
\endgroup%%------------------------------------<<

Let $\langle 2n\rangle$ be the category whose objects are the integers in the interval $[0,2n]$ and whose morphisms, apart from identities, are depicted by 
$\underset{0}{\bullet} \ra \underset{1}{\bullet} \la \ldots \ra \underset{2n-1}{\bullet}\la \underset{2n}{\bullet}$.  
Put  
%%----------------------------------------------------------------------------------------------03
$I_{2n} = \ner\langle 2n\rangle$: $\abs{I_{2n}}$ is homeomorphic to $[0,2n]$.  
Given a simplicial set $X$, a 
\un{path} 
\index{path (in a simplicial set)} 
in $X$ is a simplicial map 
$\sigma:I_{2n}  \ra X$.  
One says that $\sigma$ \un{begins} at $\sigma(0)$ and \un{ends} at $\sigma(2n)$.  
Write $\pi_0(X)$ for the quotient of $X_0$ with respect to the equivalence relation obtained by declaring that 
$x^\prime \sim x\pp$ iff there exists a path in $X$ which begins at $x^\prime$ and ends at $x\pp$ $-$then the assignment 
$X \ra \pi_0(X)$ defines a functor $\pi_0:\bSISET \ra \bSET$ which preserves finite products and is a left adjoint for the functor 
$\si:\bSET \ra \bSISET$ that sends $X$ to siX, the 
\un{constant simplicial set}
\index{constant simplicial set} 
on $X$, i.e., 
$\si X([n]) = X$ $\&$ 
$
\begin{cases}
\ d_i = \id_X\\
\ s_i = \id_X
\end{cases}
(\forall \ n).
$

[Note: \  The geometric realization of siX is $X$ equipped with the discrete topology.]\\

\begingroup%%----------------------------------->>
\fontsize{9pt}{11pt}\selectfont
Let $X$ be a simplicial set, $\Pi X$ its fundamental groupoid $-$then there is a canonical surjection  \ 
$\ds\bigcup\limits_0^\infty \Nat (I_{2n},X)$ $\ra$ $\Mor \Pi X$ 
compatible with the composition of morphisms.  
Thus fix $n$ and call 
$\ini_i: \dw \ra I_{2n}$ the injection corresponding to $i$.  
Attach to 
$\sigma:I_{2n} \ra X$ an element $x_i \in X_1$ by setting 
$x_i = \sigma \circx \ini_i(\id_{[1]})$: $\sigma \ra \pi_\sigma \in \Mor \Pi X$, where 
$\pi_\sigma = x_{2n}^{-1} \circx x_{2n-2} \circx \cdots \circx x_2^{-1} \circx x_1$.  
Corollary: $\pi_0(X) \leftrightarrow \pi_0(\Pi X)$.
\vspi
[Note: \  $\Pi X$ and $\Pi \abs{X}$ are equivalent but, in general, not isomorphic.]\\
\endgroup%%------------------------------------<<

\label{14.58}
\label{14.59}
\begingroup%%----------------------------------->>
\fontsize{9pt}{11pt}\selectfont
\textbf{\small FACT} \ 
Let $X$ be a simplicial set; let 
$
\begin{cases}
\ d_1: X_1 \ra X_0\\
\ d_0: X_1 \ra X_0
\end{cases}
$
$-$then $\pi_0(X) \approx \coeq(d_1,d_0)$.\\
\endgroup%%------------------------------------<<

Given a simplicial set $X$, the decomposition of $X_0$ into equivalence classes determines a partition of $X$ into simplicial subsets $X_i$.  The $X_i$ are called the 
\un{components} 
\index{components (simplicial set)} 
of $X$ and $X$ is 
\un{connected} 
\index{connected (simplicial set)} 
if it has exactly one component.

[Note: \  $X = \coprod\limits_i X_i$ $\implies$ $\abs{X} = \coprod\limits_i \abs{X_i}$, 
$\abs{X_i}$ running through the components of $\abs{X}$, so $\pi_0(X) \leftrightarrow \pi_0(\abs{X})$.]\\

\begingroup%%----------------------------------->>
\fontsize{9pt}{11pt}\selectfont
\textbf{\small EXAMPLE} \ 
A small category \bC is connected iff its nerve $\ner\bC$ is connected or, equivalently, iff its classifying space $B\bC$ is connected (= path connected).\\
\endgroup%%------------------------------------<<

\label{13.102}
Let $B$ be a simplicial set. An object in $\bSISET/B$ is a simplicial set $X$ together with a simplicial map 
$p:X \ra B$ called the 
\un{projection}.
\index{projection (\bSISET/B)}    
Given $b \in B_n$, define $X_b$ by the pullback square
\begin{tikzcd}[sep=large]
{X_b} \ar{d} \ar{r} &{X} \ar{d}{p}\\
{\dn} \ar{r}[swap]{\Delta_b} &{B}
\end{tikzcd}
$-$then $X_b$ is the 
\un{fiber} 
\index{fiber (simplicial set)} of $p$ over $b$ if $b \in B_0$.\\

There is a functor $\bSISET \ra \bSISET/B$ that sends a simplicial set $T$ to $B \times T$ with projection 
$B \times T \ra B$.  An $X$ in $\bSISET/B$ is said to be 
\un{trivial} 
\index{trivial ($\bSISET/B$)} 
if there exists a $T$ in \bSISET such that $X$ is isomorphic over $B$ to $B \times T$, 
\un{locally trivial} 
\index{locally trivial ($\bSISET/B$)} 
if $\forall \ n$ $\&$ $\forall \ b \in B_n$, $X_b$ is trival over $\dn$, say $X_b \approx \dn \times T_b$.

%%----------------------------------------------------------------------------------------------04
[Note: \  If for some \mT, $T_b \approx T$ $\forall \ n$ $\&$ $\forall \ b \in B_n$, then $X$ is said to be 
\un{locally trivial with} 
\un{fiber \mT}.] 
\index{locally trivial with fiber $T$ ($\bSISET/B$)}

Notation: Given $b \in B_n$, let $b_0, b_1, \ldots, b_n$ be its vertex set, i.e., 
$b_i = (B\epsilon_i)b$, $\epsilon_i:[0] \ra [n]$ the $i^{th}$ vertex operator $(i = 0, 1, \ldots, n)$.\\

\label{13.49}

\textbf{\small SUBLEMMA} \ 
Let $X$ be in $\bSISET/B$.  
Assume $X$ is locally trivial $-$then $\forall \ b \in B_n$, $T_b$ is isomorphic to $X_{b_i}$ $(i = 0, 1, \ldots, n)$.

[Take $i = 0$ and consider the \cd
\begin{tikzcd}[sep=large]
{X_{b_0}} \ar{d} \ar{r} &{X_b} \ar{d} \ar{r} &{X} \ar{d}{p}\\
{\dz} \ar{r} &{\dn} \ar{r}[swap]{\Delta_b} &{B}
\end{tikzcd}
.  Here
\begin{tikzcd}[sep=large]
{X_{b_0}} \ar{d} \ar{r} &{X_b} \ar{d}\\
{\dz} \ar{r} &{\dn}
\end{tikzcd}
is a pullback square.  
But $X_b$, viewed as an object in $\bSISET/\dn$, is isomorphic to $\dn \times T_b$, so $X_{b_0}$ is isomorphic to $T_b$.]\\

\textbf{\small LEMMA} \ 
Let $X$ be in $\bSISET/B$.  
Assume $X$ is locally trivial and $B$ is connected $-$then $X$ is locally trivial with fiber \mT.

[The sublemma implies that $\forall$
$
\begin{cases}
\ b^\prime \in B_{n^\prime}\\
\ b\pp \in B_{n\pp}
\end{cases}
,
$
$
\begin{cases}
\ T_{b^\prime} \approx X_{b_0^\prime}\\
\ T_{b\pp} \approx X_{b_0\pp}
\end{cases}
$
and $\forall \ b \in B_1$, $X_{b_0} \approx X_{b_1}$.]\\

The terms ``trivial'', ``locally trivial'', and ``locally trivial with fiber \mT'' as used in \bTOP are also used in \bCGH, the only difference being that the products are taken in \bCGH.\\

\begin{proposition} \ %02
Let $X$ be a locally trivial object in $\bSISET/B$ $-$then $\abs{X}$ is a locally trivial object in \bCGH/$\abs{B}$.
\end{proposition}

[There is no loss in generality in assuming that $\abs{B}$ is connected, hence that $B$ is connected.  
So, thanks to the lemma, $X$ is locally trivial with fiber $T$ and the contention is that 
$\abs{X}$ is locally trivial with fiber $\abs{T}$.  
Fix a point $[b,t] \in \abs{B}$ with $b \in B_n^\#$, $t \in \overset{\circx \hspace{0.2cm}}{\Delta^n}$ 
$-$then the asssociated $n$-cell $e_b$ is an open subset of $\abs{B^{(n)}} = \abs{B}^{(n)}$.  
Employing a standard collaring procedure, one can find an expanding sequence 
$e_b = O_n \subset O_{n+1} \subset \cdots $ of subsets of $\abs{B}$ such that 
$O_\infty = \colimx O_m$ is open in $\abs{B}$ and contains $e_b$ as a strong deformation retract.  
In this connection, recall that 
$O_{m-1} = \abs{B^{(m-1)}} \hspace{0.05cm} \cap \hspace{0.05cm}  O_m$, 
$O_m$ is open in $\abs{B^{(m)}}$, and there is a pushout square
\begin{tikzcd}%[sep=large]
{\coprod\limits_{x \in B_m^\#} \dot O_x} \ar{d} \ar{r} &{O_{m-1}} \ar{d}\\
{\coprod\limits_{x \in B_m^\#}  O_x} \ar{r} &{O_m}
\end{tikzcd}
, where $\forall \ x$, 
$
\begin{cases}
\ \dot O_x \subset \dot\Delta^m\\
\ O_x \subset \Delta^m
\end{cases}
$
and $\dot O_x \ra O_x$ is a closed cofibration, thus $O_{m-1} \ra O_m$ is a closed cofibration.  It will, of course, be enough to prove that $\abs{p}^{-1}(O_\infty) \approx O_\infty \times_k \abs{T}$.  
One can go further.  Indeed $O_\infty \times_k \abs{T} = \colim(O_m \times_k \abs{T})$ 
%%----------------------------------------------------------------------------------------------05
and $\abs{p}^{-1}(O_\infty) = \colim \abs{p}^{-1}(O_m)$, which reduces the problem to constructing a compatible sequence of homeomorphisms 
\begin{tikzcd}[sep=small]
{\abs{p}^{-1}(O_m)} \ar{rddd} \ar{rr} &&{O_m \times_k \abs{T}} \ar{lddd}\\
\\
\\
&{O_m}
\end{tikzcd}
.

\indent\indent $(m = n)$ \ 
Applying $\abs{?}$ to the pullback square
\begin{tikzcd}%[sep=large]
{X_b} \ar{d} \ar{r} &{X} \ar{d}{p}\\
{\dn} \ar{r}[swap]{\Delta_b} &{B}
\end{tikzcd}
in \bSISET gives a pullback square
\begin{tikzcd}%[sep=large]
{\abs{X_b}} \ar{d} \ar{r} &{\abs{X}} \ar{d}{\abs{p}}\\
{\dn} \ar{r}[swap]{\abs{\Delta_b}} &{\abs{B}}
\end{tikzcd}
in \bCGH (cf. Proposition 1).  
On the other hand, $X_b \approx \dn \times T$ and 
$\abs{\Delta_b}:\mdpn \ra e_b$ is a homeomorphism.

\indent\indent $(m > n)$ \ 
Suppose that the homeomorphism 
\begin{tikzcd}[sep=small]
{\abs{p}^{-1}(O_{m-1})} \ar{rdd} \ar{rr} &&{O_{m-1} \times_k \abs{T}} \ar{ldd}\\
\\
&{O_{m-1}}
\end{tikzcd}
has been constructed.  There is a pushout square
\[
\begin{tikzcd}[sep=large]
{\coprod\limits_{x \in B_m^\#} \dot O_x \times_{\abs{B}}\abs{X}} \ar{d} \ar{r} 
&{\abs{p}^{-1}(O_{m-1})} \ar{d}\\
{\coprod\limits_{x \in B_m^\#}  O_x  \times_{\abs{B}}\abs{X}} \ar{r} 
&{\abs{p}^{-1}(O_m)}
\end{tikzcd}
, 
\]homeomorphisms
\[
\begin{tikzcd}[sep=small]
{\dot O_x \times_{\abs{B}}\abs{X}} \ar{rddd} \ar{rr} &&{\dot O_x \times_k\abs{T}} \ar{lddd}\\
\\
\\
&{\dot O_x}
\end{tikzcd}
\ , \ 
\begin{tikzcd}[sep=small]
{O_x \times_{\abs{B}}\abs{X}} \ar{rddd} \ar{rr} &&{ O_x \times_k\abs{T}} \ar{lddd}\\
\\
\\
&{O_x}
\end{tikzcd}
\]
and a \cd
\[
\begin{tikzcd}[sep=large]
{\coprod\limits_{x \in B_m^\#}  O_x \times_{\abs{B}}\abs{X}} \ar{d} 
&{\coprod\limits_{x \in B_m^\#} \dot O_x \times_{\abs{B}}\abs{X}} \ar{d} \ar{l}\ar{r}
&{\abs{p}^{-1}(O_{m-1})} \ar{d}\\
{\coprod\limits_{x \in B_m^\#} O_x \times_k\abs{T}}
&{\coprod\limits_{x \in B_m^\#} \dot O_x \times_k\abs{T}} \ar{l} \ar{r}
&{O_{m-1} \times_k \abs{T}}
\end{tikzcd}
\]
compatible with the projections.  Accordingly, the induced map 
$\abs{p}^{-1}(O_m) \ra O_m \times_k \abs{T}$ is a homeomorphism over $O_m$.]\\

\label{13.37}
Application:  Let $X$ be in $\bSISET/B$. 
Assume: $X$ is locally trivial $-$then
$\abs{p}:\abs{X} \ra \abs{B}$ is a \bCG fibration 
(cf. p. \pageref{13.7}), 
thus is Serre 
(cf. p. \pageref{13.8}).\\

%%----------------------------------------------------------------------------------------------06

\begingroup%%----------------------------------->>
\fontsize{9pt}{11pt}\selectfont
The following lemma has been implicitly used in the proof of Proposition 2.\\
\endgroup%%------------------------------------<< 

\begingroup%%----------------------------------->>
\fontsize{9pt}{11pt}\selectfont
\textbf{\small LEMMA} \ 
Fix $B$ in \bCGH, $X$ in $\bCGH/B$,and let $\Delta:\bI \ra \bCGH/B$ be a diagram.  Assume: The colimit of $\Delta$ calculated in 
\bTOP is Hausdorff $-$then the arrow 
$\colim(\Delta_i \times_B X) \ra (\colim \Delta_i) \times_B X$ is a homeomorphism of compactly generated Hausdorff spaces.\\
\endgroup%%------------------------------------<< 

\begingroup%%----------------------------------->>
\fontsize{9pt}{11pt}\selectfont
Let \mX be in $\bSISET/B$ $-$then $p:X \ra B$ is said to be a 
\un{covering projection} 
\index{covering projection}
if \mX is locally trivial and $\forall \ b \in B_0$, $X_b$ is discrete, i.e., $X_b = X_b^{(0)}$.\\
\endgroup%%------------------------------------<< 

\begingroup%%----------------------------------->>
\fontsize{9pt}{11pt}\selectfont
\textbf{\small FACT} \ 
A simplicial map $p:X \ra B$ is a covering projection iff every commutative diagram 
\begin{tikzcd}[sep=large]
{\dz} \ar{d} \ar{r} &{X} \ar{d}{p}\\
{\dn} \ar{r} &{B}
\end{tikzcd}
has a unique filler.\\
\endgroup%%------------------------------------<< 

\begingroup%%----------------------------------->>
\fontsize{9pt}{11pt}\selectfont
\textbf{\small EXAMPLE} \ 
A covering projection in \bSISET is sent by $\abs{?}$ to a covering projection in \bTOP and a covering projection in \bTOP is sent by sin to a covering projection in \bSISET.\\
\endgroup%%------------------------------------<< 

\begingroup%%----------------------------------->>
\fontsize{9pt}{11pt}\selectfont
\textbf{\small EXAMPLE} \ 
Let \bC be a small category $-$then the category of covering spaces of $B\bC$ is equivalent to the functor category 
$[\pi_1(\bC),\bSET]$, $\pi_1(\bC)$ the fundamental groupoid of \bC 
(cf. p. \pageref{13.9}).\\
\endgroup%%------------------------------------<< 

\begin{proposition} \ %03
Let $\Phi$, $\Psi: \textbf{$\bDelta$} \ra \bSISET$ be functors; let $\Xi \in \text{Nat}(\Phi, \Psi)$.  
Assume: $\forall \ n$, $\abs{\Xi_{[n]}}: \abs{\Phi{[n]}} \ra \abs{\Psi{[n]}}$ is a homotopy equivalence $-$then $\forall$ simplicial set $X$, the geometric realization of the arrow $\Gamma_\Phi X \ra \Gamma_\Psi X$ is a homotopy equivalence provided that $\Gamma_\Phi$, $\Gamma_\Psi$ preserve injections.
\end{proposition}

[$\Gamma_\Phi$, $\Gamma_\Psi$ are the realization functors corresponding to $\Phi$, $\Psi$, so 
$\Gamma_\Phi \circx \Delta = \Phi$, 
$\Gamma_\Psi \circx \Delta = \Psi$ 
(cf. p. \pageref{13.10}), 
thus the assertion is true if $X = \dn$, thus too if $X = \coprod \dn$.  
In general there are pushout squares\\
\[
\begin{tikzcd}%[sep=large]
{X_n^\# \cdot \Gamma_\Phi \ddn} \ar{d}\ar{r} &{\Gamma_\Phi X^{(n-1)}} \ar{d}\\
{X_n^\# \cdot \Gamma_\Phi \dn} \ar{r} &{\Gamma_\Phi X^{(n)}}
\end{tikzcd}
\ , \ 
\begin{tikzcd}%[sep=large]
{X_n^\# \cdot \Gamma_\Psi \ddn} \ar{d}\ar{r} &{\Gamma_\Psi X^{(n-1)}} \ar{d}\\
{X_n^\# \cdot \Gamma_\Psi \dn} \ar{r} &{\Gamma_\Psi X^{(n)}}
\end{tikzcd}
,
\]
where, by hypothesis, the vertical arrows on the left are injective simplicial maps.  Consider now the \cd\\
\[
\begin{tikzcd}%[sep=large]
{X_n^\# \cdot \abs{\Gamma_\Phi \dn}} \ar{d}
&{X_n^\# \cdot \abs{\Gamma_\Phi \ddn}} \ar{l}\ar{d} \ar{r}
&{\abs{\Gamma_\Phi X^{(n-1)}}} \ar{d}\\
{X_n^\# \cdot \abs{\Gamma_\Psi \dn}} 
&{X_n^\# \cdot \abs{\Gamma_\Psi \ddn}} \ar{l} \ar{r}
&{\abs{\Gamma_\Psi X^{(n-1)}}}
\end{tikzcd}
.
\]
%%----------------------------------------------------------------------------------------------07
Since the geometric realization of an injective simplicial map is a closed cofibration and since inductively the arrows 
$\abs{\Gamma_\Phi \ddn} \ra \abs{\Gamma_\Psi \ddn}$, 
$\abs{\Gamma_\Phi X^{(n-1)}} \ra \abs{\Gamma_\Psi X^{(n-1)}}$ 
are homotopy equivalences, the induced map 
$\abs{\Gamma_\Phi X^{(n)}} \ra \abs{\Gamma_\Psi X^{(n)}}$ 
of pushouts is a homotopy equivalence 
(cf. p. \pageref{13.11} ff.).  
Finally, 
$
\begin{cases}
\ \Gamma_\Phi X = \colim \Gamma_\Phi X^{(n)}\\
\ \Gamma_\Psi X = \colim \Gamma_\Psi X^{(n)}
\end{cases}
$
$\implies$ 
$
\begin{cases}
\ \abs{\Gamma_\Phi X} = \colim \abs{\Gamma_\Phi X^{(n)}}\\
\ \abs{\Gamma_\Psi X} = \colim \abs{\Gamma_\Psi X^{(n)}}
\end{cases}
,
$
which leads to the desired conclusion (cf. $\S 3$, Proposition 15).]\\

\begingroup%%----------------------------------->>
\fontsize{9pt}{11pt}\selectfont
\textbf{\small EXAMPLE} \ 
Let $\Phi:\bDelta \ra \bSISET$ be a functor such that $\forall \ n$, $\abs{\Phi[n]}$ is contractible.  
Assume given a natural transformation $\Phi \ra Y_{\bDelta}$ $-$then $\forall$ simplicial set $X$, 
$\abs{\Gamma_\Phi X} \ra \abs{X}$ is a homotopy equivalence whenever $\Gamma_\Phi$ preserves injections.\\
\endgroup%%------------------------------------<<

Let $M_{\bDelta}$ be the set of monomorphisms in $\Mor \bDelta$; 
let $E_{\bDelta}$ be the set of epimorphisms in $\Mor \bDelta$ 
$-$then every $\alpha \in \Mor \bDelta$ can be written uniquely in the form 
$\alpha = \alpha^\sharp \circx \alpha^\flat$, where $\alpha^\sharp \in M_{\bDelta}$ and $\alpha^\flat \in E_{\bDelta}$.

[Note: \  Every $\alpha \in E_{\bDelta}$ has a ``maximal'' right inverse $\alpha^+ \in M_{\bDelta}$, viz. 
$\alpha^+(i) = \max \alpha^{-1}(i)$.]

Notation: $\bDelta_M$ is the category with $\Ob\bDelta_M = \Ob\bDelta$ and $\Mor\bDelta_M = M_{\bDelta}$, 
$\iota_M:\bDelta_M \ra \bDelta$ being the inclusion and $\Delta_M:\bDelta_M \ra \widehat{\bDelta}_M$ being the Yoneda embedding.

Write \bSSISET 
\index{\bSSISET} 
for the functor category $[\bDelta_M^\OP,\bSET]$ 
$-$then an object in \bSSISET is called a 
\un{semisimplicial set} 
\index{semisimplicial set} 
and a morphism in \bSSISET is called a 
\un{semisimplicial map}.
\index{semisimplicial map}  
There is a commutative triangle
\begin{tikzcd}[sep=large]
{\bDelta_M} \ar{d}[swap]{\Delta_M} \ar{r}{\bDelta \circx \iota_M} &{\widehat{\bDelta}} \\
{\widehat{\bDelta}_M} \ar{ru}[swap]{\Gamma_{\bDelta \circx \iota_M}}
\end{tikzcd}
, where $\Gamma_{\bDelta \circx \iota_M}$ is the realization functor corresponding to $\bDelta \circx \iota_M$.  
It assigns to a semisimplicial set $X$ a simplicial set $PX$, the 
\un{prolongment} 
\index{prolongment (of a semisimplicial set)} 
of $X$.  
Explicitly, the elements of $(PX)_n$ are all pairs $(x,\rho)$ with $x \in X_p$ and $\rho:[n] \ra [p]$ an epimorphism, thus 
$(PX\alpha)(x,\rho) = ((X(\rho \circx \alpha)^\sharp)x, (\rho \circx \alpha)^\flat)$ if the codomain of $\alpha$ is $[n]$.  
And: $P$ assigns to a semisimplicial map $f:X \ra Y$ the simplicial map 
$
Pf: 
\begin{cases}
\ PX \ra PY\\
\ (x,\rho) \mapsto (f(x),\rho)
\end{cases}
. \ 
$
The prolongment functor is a left adjoint for the forgetful functor 
$U:\widehat{\bDelta} \ra \widehat{\bDelta}_M$ (the singular functor in this setup.)

[Note: \   The Kan extension theorem implies that $U$ is also a left adjoint.  In particular: $U$ preserves colimits.]

Definition: $\abs{?}_M = \abs{?} \circx P$.  So, $(\abs{?}_M, U \circx \sin)$ is an adjoint pair and $\abs{?}_M$ is the realization functor determined by the composite $\Delta^? \circx \iota_M$.

[Note: \ $\abs{?}_M:\bSSISET \ra \bCGH$ does not preserve finite products.]\\

\label{14.17}

\begin{proposition} \ %04
For any simplicial set $X$, the arrow $\abs{UX}_M \ra \abs{X}$ is a homotopy
%%----------------------------------------------------------------------------------------------08
equivalence.
\end{proposition}

[In the notation of Proposition 3, take $\Phi = P \circx U \circx \Delta$, $\Psi = \Delta$, and let
$\Xi \in \Nat(\Phi,\Psi)$ be the natural transformation arising from the arrow of adjunction $P \circx U \ra \id$ 
via precomposition.  
Because $\Gamma_\Phi$, $\Gamma_\Psi$ preserve injections, it need only be shown that $\forall \ n$, 
the arrow $\abs{PU\dn} \ra \abs{\dn}$ is a homotopy equivalence or still, that $\forall \ n$, $\abs{PU\dn}$ is contractible.  
Suppose first that $n = 0$.  
In this case $\abs{PU\dz} = \coprod\limits_n \Delta^n/\sim$, the equivalence relation being generated by writing
$(t_0,\ldots,t_{i-1},0,t_{i+1}, \ldots, t_n) \sim (t_0,\ldots,t_{i-1},t_{i+1}, \ldots, t_n)$.  
Therefore $\abs{PU\dz}$ is the infinite dimensional ``dunce hat'' \mD.  As such, it is contractible.  
For positive $n$, let $D * \cdots * D$ be the quotient of $D \times \cdots \times D \times \Delta^n$ 
with respect to the relations 
$(d_0^\prime, \ldots, d_n^\prime,(t_0, \ldots, t_n)) \sim (d_0\pp, \ldots , d_n\pp,(t_0, \ldots, t_n))$ 
iff $d_i^\prime = d_i\pp$ when $t_i \neq 0$ $-$then up to homeomorphism, $\abs{PU\dn}$ is $D * \cdots * D$, a contractible space.]\\

Given $n$, let $\ov{\Delta}[n]$ be the simplicial set defined by the following conditions.

\indent\indent (Ob) \  $\ov{\Delta}[n]$ assigns to an object $[p]$ the set $\ov{\Delta}[n]_p$ of all finite sequences 
$\mu = (\mu_0, \ldots, \mu_p)$ of monomorphisms in $\bDelta$ having codomain $[n]$ such that $\forall \ i, \ j$ 
$(0 \leq i \leq j \leq p)$ there is a monomorphism $\mu_{ij}$ with $\mu_i = \mu_j \circx \mu_{ij}$.

\indent\indent (Mor) \  $\ov{\Delta}[n]$ assigns to a morphism $\alpha:[q] \ra [p]$ the map
$\ov{\Delta}[n]_p \ra \ov{\Delta}[n]_q$ taking $\mu$ to $\mu \circx \alpha$, i.e., 
$(\mu_0, \ldots, \mu_p) \ra (\mu_{\alpha(0)}, \ldots, \mu_{\alpha(q)})$.

Call $\ov{\Delta}$ the functor $\bDelta \ra \widehat{\bDelta}$ that sends $[n]$ to $\ov{\Delta}[n]$ and 
$\alpha:[m] \ra [n]$ to $\ov{\Delta}[\alpha]:\ov{\Delta}[m] \ra \ov{\Delta}[n]$, where
$\ov{\Delta}[\alpha]\nu = ((\alpha \circx \nu_0)^\sharp, \ldots, (\alpha \circx \nu_p)^\sharp)$.  
The associated realization functor $\Gamma_{\ov{\Delta}}$ is a functor $\bSISET \ra \bSISET$ such that
$\Gamma_{\ov{\Delta}} \circx \Delta = \ov{\Delta}$.  
It assigns to a simplicial set $X$ a simplicial set 
$\Sd X = \ds\int^{[n]} X_n \cdot \ov{\Delta}[n]$, 
the 
\un{subdivision}
\index{subdivision (simplicial set)} 
of $X$, and to a simplicial map $f:X \ra Y$ a simplicial map 
$\Sd f: \Sd X \ra \Sd Y$, the 
\un{subdivision} 
\index{subdivision (simplicial map)} 
of $f$.  
In particular, $\Sd \dn = \ov{\Delta}[n]$ and $\Sd \Delta[\alpha] = \ov{\Delta}[\alpha]$.  
On the other hand, the realization functor $\Gamma_\Delta$ associated with the Yoneda embedding $\Delta$ is naturally isomorphic to the identity functor id on 
\bSISET: $X = \ds\int^{[n]} X_n \cdot \dn$. 
If 
$\td_n: \ov{\Delta}[n] \ra \dn$ is the simplicial map that sends 
$\mu = (\mu_0, \ldots, \mu_p) \in \ov{\Delta}[n]_p$ to 
$\td_n\mu \in \dn_p$ : 
$\td_n(\mu(i) = \mu_i(m_i)$ 
$(\mu_i:[m_i] \ra [n])$, then the 
$\td_n$ determine a natural transformation 
$\td:\ov{\Delta} \ra \Delta$ which, by functoriality, leads to a natural transformation 
$\td:\Gamma_{\ov{\Delta}} \ra \Gamma_\Delta$.  
Thus, $\forall \ X, \ Y$ and $\forall \ f:X \ra Y$ there is a commutative diagram
\begin{tikzcd}%[sep=large]
{\Sd X} \ar{d}[swap]{\Sd f} \ar{r}{\td_X} &{X} \ar{d}{f} \\
{\Sd Y} \ar{r}[swap]{\td_Y} &{Y}
\end{tikzcd}
.  
It will be shown below that $\abs{\td_X}:\abs{\Sd X} \ra \abs{\Sd Y}$ is a homotopy equivalence (cf. Proposition 5).

Given $n$, write 
$\ovdn$
%$\ov{\Delta}^n$ 
%$\overset{   \underset{  \ov{\hspace{0.3cm}}} {\hspace{0.3cm}}   \hspace{0.2cm}     }{\Delta^n}$ 
%$\overset{   \raisebox{0.05cm}{{\un{\hspace{0.3cm}}}}   \hspace{0.2cm}     }{\Delta^n}$ 
for 
$\abs{\ov{\Delta}[n]}$ and 
$\ovdalpha$ for $\abs{\ov{\Delta}[\alpha]}$.  
The elements of $\ovdn$  are equivalence classes 
$[\mu,t]$.  
Any two representatives of $[\mu,t]$ are related by a finite chain of ``elementary
%%----------------------------------------------------------------------------------------------09
equivalences'' involving omission of $\mu_i$ and $t_i$ if $t_i = 0$ and replacement of $t_i$ and $t_{i+1}$ by 
$t_i + t_{i+1}$ if $\mu_{i+1} = \mu_i$. 
Every $[\mu,t]$ has a 
\un{canonical representative},
\index{canonical representative ($\ov{\Delta}^n$)} 
meaning that $[\mu,t]$ can be represented by a pair 
$(\mu,t)$: $\mu = (\mu_0, \ldots, \mu_n) \in \ov{\Delta}[n]_n$ with $\mu_i:[i] \ra [n]$ $(0 \leq i \leq n)$ and 
$t = (t_0, \dots, t_n) \in \Delta^n$.  So, $\mu_n = \id_{[n]}$ and there exists a permutation $\pi$ of $\{0, 1, \ldots, n\}$ such that $\forall \ i$, $\mu_i([i]) = \{\pi(0),\pi(1),\ldots,\pi(i)\}$.

Notation: Given $\alpha \in M_{\bDelta}$, say $\alpha:[m] \ra [n]$, put 
$b(\alpha) = \ds\frac{1}{m+1} \sum\limits_0^m e_{\alpha(i)}\in \R^{n+1}$.\\

\textbf{\small LEMMA} \ 
For each $n \geq 0$, the assignment 
$[u,t] \ra \ds\sum\limits_0^p t_ib(\mu_i)$ is a (welldefined) homeomorphism $h_n:\ovdn \ra \Delta^n$.

[Note: \  Geometrically, $\ovdn$ is ``barycentric subdivision'' of $\Delta^n$.]\\

\begingroup%%----------------------------------->>
\fontsize{9pt}{11pt}\selectfont
The homeomorphisms $h_n$ do not determine a natural transformation  
$\abs{?}\circx \ov{\Delta} \ra \abs{?}\circx \Delta$.  In fact, it is impossible for these functors to be naturally isomorphic.  
To see this, suppose to the contrary that there exists a natural isomorphism 
$\Xi:\abs{?}\circx \ov{\Delta} \ra \abs{?}\circx \Delta$.  There would then be homeomorphisms 
$
\begin{cases}
\ \Xi_m:\ovdm \ra \Delta^m\\
\ \Xi_n:\ovdn \ra \Delta^n
\end{cases}
$
such that for any $\alpha:[m] \ra [n]$ the diagram 
\begin{tikzcd}[sep=large]
{\ovdm} \ar{d}[swap]{\ovdalpha} \ar{r}{\Xi_m} &{\Delta^m} \ar{d}{\Delta^\alpha}\\
{\ovdn} \ar{r}[swap]{\Xi_n} &{\Delta^n}
\end{tikzcd}
commutes.  Take $m = 2$, $n = 1$ and trace the effect on the pair $(\id_{[2]},1)$ when $\alpha$ is in succession 
$\sigma_0:[2] \ra [1]$, $\sigma_1:[2] \ra [1]$.
\vspi
[Note: \ If $\alpha:[m] \ra [n]$ is a monomorphism, then the diagram
\begin{tikzcd}[sep=large]
{\ovdm} \ar{d}[swap]{\ovdalpha} \ar{r}{h_m} &{\Delta^m} \ar{d}{\Delta^\alpha}\\
{\ovdn} \ar{r}[swap]{h_n} &{\Delta^n}
\end{tikzcd}
commutes.]\\
\endgroup%%------------------------------------<<

\index{Theorem: Subdivision Theorem (Simplicial Sets)}
\textbf{\small SUBDIVISION THEOREM} \ 
Let $X$ be a simplicial set $-$then there is a homeomorphism $h_X:\abs{\Sd X} \ra \abs{X}$.

[Before proceeding to the details, I shall first outline the argument.  
In order to define a continuous function 
$h_X:\abs{\Sd X} \ra \abs{X}$, it is enough to define a continuous function 
$\coprod\limits_n X_n \times \ovdn \ra \coprod\limits_n X_n \times \Delta^n$
that respects the relations defining 
$\abs{\Sd X} = \ds \int^{[n]} X_n \cdot \ovdn$ 
and 
$\aX = \ds \int^{[n]} X_n \cdot \Delta^n$.  
This amounts to exhibiting a collection of continuous functions 
$h_x: \ovdn \ra \dpn$ $(x \in X_n, n \geq 0)$ such that for all $\alpha:[m] \ra [n]$, the diagram
\begin{tikzcd}%[sep=large]
{\ovdm} \ar{d}[swap]{\ovdalpha} \ar{r}{h_y} &{\Delta^m} \ar{d}{\Delta^\alpha}\\
{\ovdn} \ar{r}[swap]{h_x} &{\Delta^n}
\end{tikzcd}
commutes.  Here, $y = (X\alpha)x$.  
To ensure that $h_X$ is a homeomorphism, one need only arrange that if $x \in X_n^\#$ $(n \geq 0)$, then $h_x$ restricts to a homeomorphism $h_n^{-1}(\mdpn) \ra \mdpn$.

%%----------------------------------------------------------------------------------------------10
Let $x \in X_n$.   
Consider a pair $(\mu,t)$, with $\mu = (\mu_0, \ldots,  \mu_p)  \in \ \ov{\Delta}[n]_p$ and 
$t = (t_0, \ldots, t_p) \in \Delta^p$.  
Write $(X \mu_i)x = (X\alpha_i)x_i$, where $\alpha_i$ is an epimorphism and $x_i$ is nondegenerate.  Put 
$\gamma_{ij} = (\alpha_j \circx \mu_{ij})^\flat$, 
$b_{ij} = b(\mu_j \circx \gamma_{ij}^+)$ $(0 \leq i \leq j \leq p)$.  
Definition: 
\[
h_x([\mu,t]) = t_pb_{pp} + \sum\limits_{0 \leq i < p} t_i(1 - t_p - \cdots - t_{i+1})b_{ii} + \sum\limits_{0 \leq i < j \leq p} t_i t_j b_{ij}.
\]
This expression is a convex combination of points in $\dpn$, hence is in $\dpn$.  
Moreover, its value depends only on the class $[\mu,t]$ and not on a specific representative $(\mu,t)$.  
Therefore $h_x:\ovdn \ra \dpn$ makes sense.  
Because there exist finitely many nondegenerate $\mu$ such that 
$\bigcup\limits_\mu \abs{\Delta_\mu}(\dpn) = \ovdn$, $h_x$ is continuous.  
Turning to compatibility, fix $\alpha: [m] \ra [n]$ $-$then the claim is that 
$\Delta^\alpha \circx h_y = h_x \circx \ovdalpha$.  
Given 
$\nu = (\nu_0, \ldots, \nu_p) \in \ov{\Delta}[m]_p$, let $\mu = \ov{\Delta}[\alpha]\nu \in \ov{\Delta}[n]_p$ and construct 
$\beta_i$, $y_i$, $\delta_{ij}$ per $\nu$ and $y$ exactly like 
$\alpha_i$, $x_i$, $\gamma_{ij}$ are constructed per $\mu$ and $x$.  
From the definitions, $\alpha \circx \nu_i \circx \delta_{ij}^+ = \mu_i \circx \gamma_{ij}^+$ and this implies that 
$\Delta^\alpha$ matches barycenters, which suffices.

Let $x \in X_n^\#$.  
Pick a canonical representative $(\mu,t)$ for $[\mu,t]$ $-$then $\forall \ i$, 
$\gamma_{in} =$ $\gamma_{in}^+ =$ $\id_{[i]}$ and $[\mu,t] \in h_n^{-1}(\mdpn)$ iff $t_n > 0$.  
Since each of the coordinates of $h_x([\mu,t]) \in \dpn$ is bounded from below by $t_n/(n+1)$, it follows that 
$h_x(h_n^{-1}(\mdpn)) \subset \mdpn$.  
To address the issue of injectivity, suppose that $[\mu^\prime,t^\prime]$, $[\mu\pp,t\pp] \in h_n^{-1}(\mdpn)$ and 
$h_x([\mu^\prime,t^\prime]) =$ $(t_0, \ldots, t_n) =$ $h_x([\mu\pp,t\pp])$.  
In terms of canonical representatives, one has to prove that 
$\forall \ i$, $t_i^\prime = t_i\pp$ and 
$\mu_i^\prime = \mu_i\pp$ if $t_i^\prime$ $\&$ $t_i\pp$ are $> 0$.  
This will be done by decreasing induction on $i$.  
Let
$
\begin{cases}
\ \pi^\prime\\
\ \pi\pp
\end{cases}
$
be the permutations attached to 
$
\begin{cases}
\ \mu^\prime\\
\ \mu\pp
\end{cases}
. \ 
$
Looking at $t_{\pi^\prime(n)} = t_n^\prime/(n+1)$ and $t_{\pi\pp(n)} = t_n\pp/(n+1)$ yields $t_n^\prime = t_n\pp$, starting the induction.  
Assume that $k < n$ and that the assertion is true $\forall \ i > k$.  
Define $T^\prime = (T_0^\prime, \ldots, T_n^\prime)$ by
\[
\sum\limits_{0 \leq i \leq k} t_i^\prime (1 - t_n^\prime - \cdots - t_{i+1}^\prime) b_{ii}^\prime + \sum\limits_{\substack{0 \leq i \leq k, \\ i < j \leq n}} t_i^\prime t_j^\prime b_{ij}^\prime.
\]
Define $T\pp = (T_0\pp, \ldots, T_n\pp)$ analogously $-$then, from the induction hypothesis, $T^\prime = T\pp$.  
Case 1: $\mu_k^\prime \neq \mu_k\pp$.  Choose $l \in [n]$ : $l \in \mu_k^\prime([k])$ $\&$ $l \notin \mu_k\pp([k])$ 
$\implies$ $t_k^\prime t_n^\prime /(k+1) \leq$ $T_l^\prime  =$ $T_l\pp = 0$ $\implies$ $t_k^\prime  = 0$.  
Similarly, $t_k\pp = 0$.  
Case 2: $\mu_k^\prime = \mu_k\pp$.  Take $T^\prime$ and split off
\[
(1 - t_n^\prime - \cdots - t_{k+1}^\prime) b_{kk}^\prime + \sum\limits_{k < j \leq n} t_j^\prime b_{kj}^\prime
\]
to get 
$S^\prime = (S_0^\prime, \ldots, S_n^\prime)$.  
Do the same with $T\pp$ to get 
$S\pp = (S_0\pp, \ldots, S_n\pp)$ 
$-$then from the induction hypothesis, $S^\prime = S\pp$.  
Set $l = \pi^\prime(k)$ and compute: 
$t_k^\prime S_l^\prime = $ $T_l^\prime = $ $T_l\pp \geq$
%%----------------------------------------------------------------------------------------------11
$t_k\pp S_l\pp= $ $t_k\pp S_l^\prime$.  
But $S_l^\prime \geq t_n^\prime/(k+1) > 0$ $\implies$ $t_k^\prime \geq t_k\pp$.  
Similarly, $t_k\pp \geq t_k^\prime$.  Thus the induction is complete.  
Owing to the theorem of invariance of domain, $h_x(h_n^{-1}(\mdpn))$ is open in $\mdpn$ and the restriction 
$h_n^{-1}(\mdpn) \ra h_x(h_n^{-1}(\mdpn))$ is a homeomorphism.  However, 
$h_x(\ovdn - h_n^{-1}(\mdpn)) \subset \ddpn$, so
$h_x(h_n^{-1}(\mdpn)) = \mdpn \cap h_x(\ovdn)$ is closed in $\mdpn$.  
Being nonempty, 
$h_x(h_n^{-1}(\mdpn))$ must be equal to $\mdpn$.]\\

\begingroup%%----------------------------------->>
\fontsize{9pt}{11pt}\selectfont
\index{Barratt's Lemma}
\textbf{\small BARRATT'S LEMMA} \ 
Let $\Delta$ be a simplex, $\Delta_1$ a proper face of $\Delta$, $\Delta_0$ a proper face of $\Delta_1$.  Let $r:\Delta_1 \ra \Delta_0$ be an affine retraction, i.e., a retraction induced by the composition of a linear map and a translation mapping vertexes onto vertexes.  Define $X$ by the pushout square
\begin{tikzcd}[sep=large]
\Delta_1 \ar{d} \ar{r}{r} &{\Delta_0} \ar{d}\\
\Delta \ar{r} &X
\end{tikzcd}
$-$then there exists a homeomorphism $\phi:X \ra \Delta$ such that the triangle
\begin{tikzcd}[sep=large]
{\Delta_0} \ar{d} \ar{r} &{\Delta} \\
X \ar{ru}[swap]{\phi}
\end{tikzcd}
commutes.
\vspi
[Supposing that $n + 1 = \dim \Delta$, normalize the situation as follows.  Take for $\Delta$ the one point compactification of 
$\{(x_0, \ldots, x_n): x_n \geq 0\}$, let $\Delta_1$ be the convex hull of 
$\{0,e_0, \ldots, e_m\}$, 
let $\Delta_0$ be the convex hull of 
$\{0,e_0, \ldots, e_k\}$, 
and let $P$ be the orthogonal projection onto the span of 
$\{e_0, \ldots, e_k, e_{m+1}, \ldots, e_n\}$, 
so 
$\restr{P}{\Delta_1} = r$ and $X = \Delta/\sim$, where $x \sim y$ iff $x = y \notin \Delta_1$ or $r(x) = r(y)$ $(x,y \in \Delta_1)$.  Let $d(x)$ be the distance of $x$ from $\Delta_1$, $f(x) = \min\{1,d(x)\}$, and put 
$\phi(x) = f(x)x + (1 - f(x))P(x)$ (thus $\phi(\infty) = \infty$ and $\restr{\phi }{\Delta_1} = r$).
\vspi
Claim: $\phi: \Delta \ra \Delta$ is surjective and $\restr{\phi}{\Delta - \Delta_1}$ is injective.
\vspi
[Given $x = (x_0, \ldots, x_n)$, set 
$x(t) = (x_0, \ldots, x_k,tx_{k+1}, \ldots, tx_m, x_{m+1}, \ldots, x_n)$.  
Obviously, $x_{k+1} =$ $\cdots =$ $x_m = 0$ $\implies$ $\phi(x) = x$.  
On the other hand, if some $x_i \neq 0$ $(k < i \leq m)$, then $t \ra \infty$ $\implies$ $x(t) \ra \infty$ $\implies$ 
$f(x(t)) = 1$ $(t \gg 0)$.  
However, 
$\phi(x(t)) = (x_0, \ldots, x_k,tf(x(t))x_{k+1}, \ldots, tf(x(t))x_m, x_{m+1}, \ldots, x_n)$ 
and the intermediate value theorem guarantees that $\exists$ $t: tf(x(t)) = 1$.  
Assume now that $x, \ y \in \Delta_1$ with $\phi(x) = \phi(y)$: $x_i = y_i$ $(i \leq k \ \& \ i > m)$, 
$f(x)x_i = f(y)y_i$ $(k < i \leq m)$ $\implies$ $y = x\ds\left(\frac{f(x)}{f(y)}\right)$.  
But $t \ra \phi(x(t))$ is one-to-one $(\implies x = y)$.  
To see this, it need only be shown that $t \ra d(x(t))$ is nondecreasing.  Proceeding by contradiction, suppose that $d(x(t^\prime)) < d(x(t))$ $(\exists \ t^\prime > t)$ 
and choose $u: d(x(t^\prime)) < u < d(x(t))$ $\implies$ $u > d(x(0))$, i.e., $x(0)$, $x(t^\prime) \in d^{-1}([0,u])$, 
$x(t) \notin d^{-1}([0,u])$, an impossibility, $d^{-1}([0,u])$ being convex.]
\vspi
Therefore $\phi$ determines a continuous bijection $X \ra \Delta$ between compact Hausdorff spaces with the stated property.]\\
\endgroup%%------------------------------------<<

\label{5.0d}
\begingroup%%----------------------------------->>
\fontsize{9pt}{11pt}\selectfont
\textbf{\small FACT} \ 
Let $X$ be a simplicial set $-$then $\abs{\Sd X}$ is a polyhedron, hence $\aX$ can be triangulated.
\vspi
[Using Barratt's lemma, apply the criterion on 
p. \pageref{13.12} 
to $\abs{\Sd X}$, observing that $\forall$ nondegenerate $x$ in $(\Sd X)_n$ there is a pushout square
\begin{tikzcd}[sep=large]
{\Delta[n-1]} \ar{d}[swap]{\Delta[\delta_n]} \ar{r} &{\langle d_nx\rangle} \ar{d} \\
{\Delta[n]} \ar{r} &{\langle x \rangle}
\end{tikzcd}
, where (?) equals ``generated simplicial subset''.]\\
\endgroup%%------------------------------------<<

%%----------------------------------------------------------------------------------------------12

\begin{proposition} \ %05
 Let $X$ be a simplicial set $-$then $\abs{\td_X}:\abs{\Sd X} \ra \abs{X}$ is a homotopy equivalence.
\end{proposition}

[One can define $\abs{\td_X}$ by a collection of continuous functions $\td_x:\ovdn \ra \dpn$ satisfying the same compatibility conditions as the $h_x:\ovdn \ra \dpn$ that figure in the proof of the subdivision theorem.  
Introduce 
$H_x:\ovdn \times [0,1] \ra \dpn$ by writing $H_x(u,t) = (1 - t)h_x(u) + t\td_x(u)$ $-$then, in total, the $H_x$ define a homotopy $\abs{\Sd X} \times [0,1] \ra \aX$  between $h_X$ and $\abs{\td_X}$.]

[Note: \  $h_X$ is not natural but is homotopic to $\abs{\td_X}$ which is natural.  The fact that $\abs{\td_X}$ is a homotopy equivalence can also be seen directly.  
Proof: $\forall \ n$, $\abs{\ov{\Delta}[n]} = \ovdn$ is contractible and 
$\Gamma_{\ov{\Delta}} = \Sd$ preserves injections, thus the example following Proposition 3 is applicable.]\\

\begingroup%%----------------------------------->>
\fontsize{9pt}{11pt}\selectfont
\textbf{\small EXAMPLE} \ 
Let $X$ be a simplicial set $-$then $\aX$ is homeomorphic to $B(c\Sd^2 X)$ 
(Fritsch-Latch\footnote[2]{\textit{Math. Zeit.} \textbf{177} (1981), 147-179.}).  
Therefore the geometric realization  of a simplicial set is homeomorphic to the classifying space of a small category.
\vspi
[Note: \  The homeomorphism is not natural.]\\
\endgroup%%------------------------------------<<

\begingroup%%----------------------------------->>
\fontsize{9pt}{11pt}\selectfont
Sd is the realization functor $\Gamma_{\ov{\Delta}}$.  
The associated singular functor $S_{\ov{\Delta}}$ is denoted by $\Ex$ 
\index{Ex (extension functor)} and referred to as 
\un{extension} 
\index{extension (functor)}.  
Since $(\Sd,\Ex)$ is an adjoint pair, there is a bijective map $\Xi_{X,Y}:\Nat(\Sd X,Y) \ra \Nat(X,\Ex Y)$ which is functorial in $X$ and $Y$ 
(cf. p. \pageref{13.13}).  
Put $e_X = \Xi_{X,X}(\td_X)$ $-$then $e_X:X \ra \Ex X$ is the simplicial map given by 
$e_X(x) = \Delta_x \circx \td_n$ $(x \in X_n)$, hence $e_X$ is injective.\\
\endgroup%%------------------------------------<<

\label{13.45}
\label{13.73}
\label{18.30}
\begingroup%%----------------------------------->>
\fontsize{9pt}{11pt}\selectfont
\textbf{\small LEMMA} \ 
For every simplicial set $X$, $\abs{e_X}:\aX \ra \abs{\Ex X}$ is a homotopy equivalence 
(cf. p. \pageref{13.14}).
\vspi
[Note: \  Since $e_X$ is injective, $\aX$ can be considered as a strong deformation retract of $\abs{\Ex X}$ (cf. $\S 3$, Proposition 5).]\\
\endgroup%%------------------------------------<<

\label{13.27}
\begingroup%%----------------------------------->>
\fontsize{9pt}{11pt}\selectfont
Denote by $\Exx^\infty$ 
\index{$\Exx^\infty$} 
the colimit of 
$\id \ra \Ex \ra \Exx^2 \ra \cdots$ $-$then $\Exx^\infty$ is a functor $\bSISET \ra \bSISET$ 
and for any simplicial set $X$, there is an arrow 
$e_X^\infty:X \ra \Exx^\infty X$.  
Claim: $\abs{e_X^\infty}:\aX \ra \abs{\Exx^\infty X}$ is a homotopy equivalence.  
In fact, $\abs{\Exx^n X}$ embeds in $\abs{\Exx^{n+1} X}$ as a strong deformation retract and 
$\abs{\Exx^\infty X} = \colim \abs{\Exx^n X}$.  
Therefore $\aX$ is a strong deformation retract of $\abs{\Exx^\infty X}$ 
(cf. p. \pageref{13.15}).\\
\endgroup%%------------------------------------<<

The subdivision functor can also be introduced in the semisimplicial setting.  
It is compatible with prolongment in that there is a commutative diagram
\begin{tikzcd}%[sep=large]
{\bSSISET} \ar{d}[swap]{P} \ar{r}{\Sd} &{\bSSISET} \ar{d}{P} \\
{\bSISET} \ar{r}[swap]{\Sd} &{\bSISET}
\end{tikzcd}
%%----------------------------------------------------------------------------------------------13
and, in contradistinction to what happens in the simplicial setting, the homeomorphism 
$h_{PX}:\abs{\Sd X}_M \ra \aX_M$ is natural, as is the homotopy between $h_{PX}$ and $\abs{\td_{PX}}$.

Put $S = U \circx \sin$ $-$then $S:\bTOP \ra \bSSISET$ and $(\abs{?}_M,S)$ is an adjoint pair.  
Given a topological space $X$, postcompose $h_{PSX}:\abs{\Sd SX}_M \ra \abs{S X}_M$ with the arrow 
$\abs{S X}_M \ra X$ to get a continuous function $\abs{\Sd SX}_M \ra X$ which by adjointness corresponds to a 
semisimplicial map $g_{SX}:\Sd SX \ra SX$.  
Definition: 
$b_X = \abs{Pg_{SX}} \circx h_{PSX}^{-1} \in C(\abs{SX}_M,\abs{SX}_M)$.  
Using Proposition 4, one can check that $b_X$ is naturally homotopic to $\id_{\abs{SX}_M}$.  In effect, the triangle
\begin{tikzcd}[sep=small]
{\abs{SX}_M} \ar{rddd} \ar{rr}{b_X} &&{\abs{SX}_M} \ar{lddd} \\
\\
\\
&{\abs{\sin X}}
\end{tikzcd}
commutes up to homotopy.\\
\vspace{0.5cm}

\index{Theorem Simplicial Excision}
\index{Simplicial Excision Theorem}
\textbf{\small SIMPLICIAL EXCISION THEOREM} \ 
Let $X$ be a topological space.  Suppose that 
$
\begin{cases}
\ X_1\\
\ X_2
\end{cases}
$
are subspaces of $X$ with $X = \text{int}X_1 \cup \text{int}X_2$ $-$then the geometric realization of 
$\text{sin}X_1 \cup \text{sin}X_2$ is a strong deformation retract of $\abs{\sin X}$.

[The inclusion $\abs{\sin X_1 \cup \sin X_2} \ra \abs{\sin X}$ is a closed cofibration, thus it will be enough to prove that it is a homotopy equivalence (cf. $\S 3$, Proposition 5).  According to Proposition 4, the vertical arrows in the commutative diagram 
\begin{tikzcd}%[sep=large]
{\abs{S X_1 \cup S X_2}_M} \ar{d} \ar{r} &{\abs{S X}_M} \ar{d} \\
{\abs{\sin X_1 \cup \sin X_2}} \ar{r} &{\abs{\sin X}}
\end{tikzcd}
are homotopy equivalences, which reduces the problem to showing that the inclusion 
$\abs{S X_1 \cup S X_2}_M \ra \abs{S X}_M$ is a homotopy equivalence or still, a weak homotopy equivalence.  
To this end, fix $n \geq 0$ and let $f:\bD^n \ra \abs{S X}_M$ be a continuous function such that 
$f(\bS^{n-1}) \subset \abs{S X_1 \cup S X_2}_M$.  
Since the image of $f$ is contained in the union of a finite number of cells of $\abs{S X}_M$, 
$\exists$ $k \gg 0$: $b_X^k \circx f$ factors through $\abs{S X_1 \cup S X_2}_M$ 
(the ``excisive'' consequence of the assumption that $X = \itr X_1 \cup \itr X_2$).  
On the other hand, by naturality, $b_X(\abs{S X_1 \cup S X_2}_M) \subset \abs{S X_1 \cup S X_2}_M$ and the same is true of the homotopy between $b_X$ and $\id_{\abs{S X}_M}$, hence too for the $k^{th}$ iterate of $b_X^k$.  
Therefore $f$ is homotopic $\rel \bS^{n-1}$ to a continuous function $g:\bD^n \ra \abs{S X}_M$ with 
$g(\bD^n) \subset \abs{S X_1 \cup S X_2}_M$.  
These considerations suffice to imply that the inclusion 
$\abs{S X_1 \cup S X_2}_M \ra \abs{S X}_M$ is a weak homotopy equivalence 
(cf. p. \pageref{13.16}).]\\

Let $\sC$ be a class of topological spaces $-$then $\sC$ is said to be 
\un{homotopy cocomplete} 
\index{homotopy cocomplete} 
provided that the following conditions are satisfied.\\
\indent\indent (HOCO$_1$) \  If $X \in \sC$ and if $Y$ has the same homotopy type as $X$, then $Y \in \sC$.\\
\indent\indent (HOCO$_2$) \  $\sC$ is closed under the formation of coproducts.\\
\indent\indent (HOCO$_3$) \  If $X \overset{f}{\la} Z \overset{g}{\ra} Y$ is a 2-source with 
$
\begin{cases}
\ X\\
\ Y
\end{cases}
\& \ Z \in \sC,
$
then $M_{f,g} \in \sC$.

Examples: 
(1) The class of CW spaces is homotopy cocomplete; 
(2) The class of numerably contractible spaces is homotopy cocomplete.\\

%%----------------------------------------------------------------------------------------------14

\begin{proposition} \ %06
The class of topological spaces for which the arrow of adjunction $\abs{\text{sin} X} \ra X$ is a homotopy equivalence is homotopy cocomplete.
\end{proposition}

[If $f:X \ra Y$ is a homotopy equivalence, then 
$\abs{\sin f}:\abs{\sin X} \ra \abs{\sin Y}$ is a homotopy equivalence 
(cf. p. \pageref{13.17}).  
Since the diagram
\begin{tikzcd}%[sep=large]
{\abs{\sin X}} \ar{d} \ar{r}{\abs{\sin f}} &{\abs{\sin Y}} \ar{d} \\
X \ar{r}[swap]{f} &{Y}
\end{tikzcd}
commutes, HOCO$_1$ obtains.  That HOCO$_2$ holds is clear, so it remains to deal with HOCO$_3$.  
Viewing $M_{f,g}$ as a quotient of $X \amalg IZ \amalg Y$, let $\ov{X}$  be the image of 
$X \amalg Z \times [0,2/3]$, 
let $\ov{Y}$ be the image of $Z \times [1/3,1] \amalg Y$ and put 
$\ov{Z} = \ov{X} \cap \ov{Y}$ $-$then 
$M_{f,g} = \itr\ov{X} \cup \itr\ov{Y}$ and there are homotopy equivalences $\ov{X} \ra X$, $\ov{Y} \ra Y$, 
$\ov{Z} \ra Z$.  
Because \mX, \mY, \mZ are in our class, the same is true of $\ov{X}, \ov{Y},\ov{Z}$.  
To establish that the arrow $\abs{\sin M_{f,g}} \ra M_{f,g}$ is a homotopy equivalence, consider the commutative diagram 
$
 \begin{tikzcd}%[sep=large]
{\abs{\sin \ov{X}}} \ar{d}   & {\abs{\sin \ov{Z}}}  \ar{l} \ar{d} \ar{r}  &{\abs{\sin \ov{Y}}} \ar{d} \\
{\ov{X}}  &{\ov{Z}} \ar{l} \ar{r} &{\ov{Y}}
\end{tikzcd}
. 
$
The horizontal arrows are closed cofibrations, hence the induced map of pushouts is a homotopy equivalence 
(cf. p. \pageref{13.18} ff.).  
The pushout arising from the 2-source on the bottom is $M_{f,g}$, while the pushout arising from the 2-source on the top is $\abs{\sin \ov{X} \cup \sin \ov{Y}}$ which, by the simplicial excision theorem, is a strong deformation retract of 
$\abs{\sin M_{f,g}}$.  Inspection of the triangle 
\begin{tikzcd}%[sep=large]
{\abs{\sin \ov{X} \cup \sin \ov{Y}}} \ar{dr} \ar{r} &{\abs{\sin M_{f,g}}} \ar{d} \\
&{M_{f,g}}
\end{tikzcd}
finishes the argument.]

[Note: \  $\forall \ X$, $\abs{\sin X}$ is a CW complex, thus $X$ is a CW space if the arrow of adjunction 
$\abs{\sin X} \ra X$ is a homotopy equivalence.\\

Any homotopy cocomplete class of topological spaces that contains a one point space necessarily contains the class of CW spaces.  But $\#(X) = 1$ $\implies$ $\#(\abs{\sin X}) = 1$, therefore the class of CW spaces is precisely the class of topological spaces for which the arrow of adjunction $\abs{\sin X} \ra X$ is a homotopy equivalence.\\

\index{Theorem Giever-Milnor}
\index{Giever-Milnor Theorem}
\textbf{\small GIEVER-MILNOR THEOREM} \ 
Let $X$ be a topological space $-$then the arrow of adjunction $\abs{\text{sin} X} \ra X$ is a weak homotopy equivalence.

[The adjoint pair $(\abs{?},\sin)$ determines a cotriple in \bTOP 
(cf. p. \pageref{13.19}), 
which induces a cotriple in 
\bHTOP ($\abs{?} \circx \sin$ preserves homotopies 
(cf. p. \pageref{13.20})).  
On general grounds, $\forall \ Y$, the 
postcomposition arrow 
$[\abs{\sin Y},\abs{\sin X}] \ra [\abs{\sin Y},X]$ is surjective.  However here it is also injective.  Reason: $\forall \ Z$, 
the arrow of adjucntion $\abs{\sin \abs{\sin Z}} \ra \abs{\sin Z}$ is a homotopy equivalence, i.e., is an isomophism in 
\bHTOP.  It therefore follows that for every CW complex $K$, the postcomposition arrow 
$[K, \abs{\sin X}] \ra [K,X]$ is bijective and
%%----------------------------------------------------------------------------------------------15
this means that the arrow of adjunction $\abs{\sin X} \ra X$ is a weak homotopy equivalence 
(cf. p. \pageref{13.21} ff.).]\\

\label{13.26}
\label{13.46}
\label{13.58}
\label{13.108}
\label{13.127}

Application: Let $X$ be a simplicial set $-$then the geometric realization of the arrow of adjunction 
$X \ra \sin \aX$ is a homotopy equivalence.

[The triangle
\begin{tikzcd}%[sep=large]
{\aX} \ar{dr}[swap]{\id_{\aX}} \ar{r} &{\abs{\sin \aX}} \ar{d} \\
&{\aX}
\end{tikzcd}
commutes.]\\
\vspace{0.25cm}

\begingroup%%----------------------------------->>
\fontsize{9pt}{11pt}\selectfont
\textbf{\small EXAMPLE} \ 
Consider the adjoint situation $(F,G,\mu,\nu)$, where $F = \abs{?}$, $G = \sin$ $-$then in the notation of 
p. \pageref{13.22}, 
$
\begin{cases}
\ S^{-1}\bSISET\\
\ T^{-1}\bTOP
\end{cases}
$
are equivalent to \bHCW.\\
\endgroup%%------------------------------------<<

\label{14.44}

Given simplicial sets $X$ and $Y$, write $\map(X,Y)$ 
\index{$\map(X,Y)$} in place of $Y^X$ 
(cf. p. \pageref{13.23}).  
The elements of $\map(X,Y)_0 \approx \Nat(X,Y)$ are the simplicial maps $X \ra Y$, two such being termed 
\un{homotopic}
\index{homotopic (simplicial maps)}
 if they belong to the same component of $\map(X,Y)$.  
In other words, simplicial maps $f,g \in \Nat(X,Y)$ are homotopic $(f \simeq g)$ provided $\exists \ n \geq 0$ and a simplicial map 
$H:X \times I_{2n} \ra Y$ such that if
$
\begin{cases}
\ H \circx i_0 \ : X \approx X \times \dz  \ \overset{\id_X \times e_0}{\longrightarrow} X \times I_{2n} \overset{H}{\lra} Y\\
\ H \circx i_{2n}: X \approx X \times \dz \underset{\id_X \times e_{2n}}{\longrightarrow} X \times I_{2n} \underset{H}{\lra} Y
\end{cases}
, \ 
$
then 
$
\begin{cases}
\ H \circx i_0 = f\\
\ H \circx i_{2n} = g
\end{cases}
, \ 
$
where 
$
\begin{cases}
\ e_0:\dz \ra I_{2n}\\
\ e_{2n}:\dz \ra I_{2n}
\end{cases}
$
are the vertex inclusions per 
$
\begin{cases}
\ 0\\
\ 2n
\end{cases}
. \ 
$

[Note: \  Paths $I_{2n} \ra \map(X,Y)$ correspond to homotopies $H:X \times I_{2n} \ra Y$.]

Given simplicial sets $X$ and $Y$, simplicial maps \ $f$, $g \in \Nat(X,Y)$ \ are said to be 
%\un{simplicially} \un{homotopic}
\un{simplicially homotopic} 
\index{simplicially homotopic} 
$(f \underset{s}{\simeq} g)$ 
\index{$\underset{s}{\simeq}$} 
provided $\exists$ a simplicial map $H:X \times \dw \ra Y$ such that if 
$
\begin{cases}
\ H \circx i_0: X \approx X \times \dz \ \overset{\id_X \times e_0}{\lra} \ X \times \dw \overset{H}{\lra} Y\\
\ H \circx i_{1}: X \approx X \times \dz \ \underset{\id_X \times e_{1}}{\lra} \ X \times \dw \underset{H}{\lra} Y
\end{cases}
, \ 
$
then \quad
$
\begin{cases}
\ H \circx i_0 = f\\
\ H \circx i_1 = g
\end{cases}
, \ 
$
where 
$
\begin{cases}
\ e_0:\dz \ra \dw\\
\ e_1:\dz \ra \dw
\end{cases}
$
are the vertex inclusions per 
$
\begin{cases}
\ 0\\
\ 1
\end{cases}
$
.  The relation $\underset{s}{\simeq}$ is reflexive but it needn't be symmetric or transitive.

[Note: \   Elements of $\map(X,Y)_1$ correspond to simplicial homotopies $H:X \times \dw \ra Y$.]

\label{14.36}
\label{18.11}
Example: Suppose that 
$
\begin{cases}
\ \bC\\
\ \bD
\end{cases}
$
are small categories.  Let $F, \ G:\bC \ra \bC$ be functors, $\Xi: F \ra G$ a natural transformations $-$then $\Xi$ defines a functor 
$\Xi_H:\bC \times [1] \ra \bD$, 
hence 
$\ner\Xi_H:\ner(\bC \times [1] ) \ra \ner\bD$, i.e., 
$\ner\Xi_H:\ner\bC \times \dw \ra \ner\bD$ is a simplicial homotopy between $\ner F$ and $\ner G$.  So, e.g., 
$
\begin{cases}
\ B\bC\\
\ B\bD
\end{cases}
$
have the same homotopy type if there is a functor $\bC \ra \bD$ which admits a left or right adjoint.  In particular: The classifying space of a small category having either an initial object or a final object is contractible.  
\label{13.112}
\label{13.121}
\label{14.97}
Example: $B\bDelta$ is contractible.\\

%%----------------------------------------------------------------------------------------------16

\label{13.36}
\begingroup%%----------------------------------->>
\fontsize{9pt}{11pt}\selectfont
\textbf{\small EXAMPLE} \  
Take $X = Y = \dn$ $(n > 0)$.  Let $C_0:\dn \ra \dn$ be the projection of $\dn$ onto the $0^{\text{th}}$ vertex, i.e., send 
$(\alpha_0, \ldots, \alpha_p) \in \dn_p$ to $(0, \ldots, 0) \in \dn_p$.  
Claim: $C_0 \underset{s}{\simeq} \id_{\dn}$.  
To see this, consider the simplicial map $H:\dn \times \dw \ra \dn$ 
defined by 
$H((\alpha_0, \ldots, \alpha_p),(0, \ldots, 0,1,\ldots,1)) = $
$(0, \ldots, 0,\alpha_{i+1}, \ldots, \alpha_p)$ \ so that \ 
$H((\alpha_0, \ldots, \alpha_p), (0, \ldots, 0))$ \ $=$ \ 
$(0, \ldots, 0), H((\alpha_0, \ldots, \alpha_p), (1, \ldots, 1)) = $\ 
$(\alpha_0, \ldots, \alpha_p)$ $-$then \mH is a simplicial homotopy between $C_0$ and $\id_{\dn}$.  
On the other hand, there is no simplicial homotopy \mH between $\id_{\dn}$ and $C_0$.  For suppose that
$H((1,1),(0,1)) = (\mu,\nu) \in \dn_1$.  Apply $d_1$ $\&$ $d_0$ to get $\mu = 1$ $\&$ $\nu = 0$, an impossibility.
\vspi
[Note: Let $C_k:\dn \ra \dn$ be the projection of $\dn$ onto the $k^{\text{th}}$ vertex, i.e., send 
$(\alpha_0, \ldots, \alpha_p) \in \dn_p$ to $(k, \ldots, k) \in \dn_p$ $(0 \leq k \leq n)$ $-$then 
$\id_{\dn}  \underset{s}{\simeq} C_n$ but 
$\id_{\dn}  \underset{s}{\not\cong} C_k$ $(0 \leq k < n)$.  
Still, $\forall \ k$, $\exists$ a homotopy $H_k:\dn \times I_2 \ra \dn$ such that $H_k \circx e_0 = \id_{\dn}$ and 
$H_k \circx e_2 = C_k$.]\\
\endgroup%%------------------------------------<<

\label{13.42}
\begingroup%%----------------------------------->>
\fontsize{9pt}{11pt}\selectfont
\textbf{\small FACT} \ 
Suppose that $f, \ g:X \ra Y$ are simplicially homotopic $-$then $\Ex f, \ \Ex g: \Ex X \ra \Ex Y$ are simplicially homotopic.
\vspi
[Ex is a right adjoint, hence preserves products.]\\
\endgroup%%------------------------------------<<

The equivalence relation generated by $\underset{s}{\simeq}$ is $\simeq$.  
Given simplicial sets $X$ and $Y$, put $[X,Y]_0 = \Nat(X,Y)/\simeq$, so $[X,Y]_0  = \pi_0(\map(X,Y))$ $-$then
\bHZEROSISET 
\index{\bHZEROSISET} 
is the category whose objects are the simplicial sets and whose morphisms are the homotopy classes of simplicial maps.

[Note: \  The symbol \bHSISET is reserved for a different role 
(cf. p. \pageref{13.24}).]\\

\begingroup%%----------------------------------->>
\fontsize{9pt}{11pt}\selectfont
To check that the relation of homotopy is compatible with composition, let $X$, $Y$, and $Z$ be simplicial sets.  
Define a simplicial map $C_{X,Y,Z}:\map(X,Y) \times \map(Y,Z) \ra \map(X,Z)$ by assigning to a pair $(f,g)$ in 
$\map(X,Y)_n \times \map(X,Z)_n$ the composite 
\begin{tikzcd}[sep=large]
{X \times \dn}  \ar{r}{\id \times \di} &{X \times (\dn \times \dn)}
\end{tikzcd}
$\overset{A}{\lra}$
${(X \times \dn) \times \dn)}$
%
%$\overset{f \times \id}{\lra} Y \times \dn$ 
$\overset{f \times \id}{\xrightarrow{\hspace*{1cm}}} Y \times \dn$ 
%
%
%\begin{tikzcd}[sep=large]
%\ar{r}{f \times \id} &{Y \times \dn }
%\end{tikzcd}
%
$\overset{g}{\ra} Z$
in $\map(X,Z)_n$.
At level 0, $C_{X,Y,Z}$ is composition of simplicial maps.  Since 
$\pi_0(\map(X,Y) \times \map(Y,Z)) \approx$ 
$\pi_0(\map(X,Y)) \times \pi_0(\map(Y,Z))$, 
$C_{X,Y,Z}$ induces an arrow 
$[X,Y]_0 \times [Y,Z]_0 \ra [X,Z]_0$ with the requisite properties.
\vspi
[Note: \  \bHZEROSISET has finite products.  In addition, 
$\map(X \times Y,Z) \approx$ $\map(X,\map(Y,Z))$ $\implies$ 
$[X \times Y,Z]_0 \approx$ $[X,\map(Y,Z)]_0$ so \bHZEROSISET is cartesian closed.]\\
\endgroup%%------------------------------------<<

\begingroup%%----------------------------------->>
\fontsize{9pt}{11pt}\selectfont
\textbf{\small EXAMPLE} \ 
Geometric realization preserves homotopies but $\abs{?}:\bHZEROSISET \ra \bHTOP$ is not conservative.
\vspi
[Take $X = \dz$, $Y = \ner(\infty)$, where $(\infty)$ is the zig-zag on the set of nonnegative integers: 
$0<1>2<3>4 \ldots$, and consider the inclusion $X \ra Y$ corresponding to $0 \ra 0$.]\\
\endgroup%%------------------------------------<<

Notation: Given a simplicial set $X$, write $IX$ in place of $X \times \Delta[1]$.\\

%%----------------------------------------------------------------------------------------------17
\begingroup%%----------------------------------->>
\fontsize{9pt}{11pt}\selectfont
The obvious composite $X \coprod X \ra IX \ra X$ factors the folding map $X \coprod X \ra X$ and \bSISET carries the structure of a model category in which $IX$ is a cylinder object 
(cf. p. \pageref{13.25}).\\
\endgroup%%------------------------------------<<

\label{13.74}
\label{13.116} %dmc mnft

A simplicial map $f:X \ra Y$ is said to be a 
\un{weak homotopy equivalence}
\index{weak homotopy equivalence (simplicial map)} 
if its geometric realization 
$\abs{f}:\abs{X} \ra \abs{Y}$ is a weak homotopy equivalence (= homotopy equivalence).  
Example: $\forall \ X$, the projection $IX \ra X$ is a weak homotopy equivalence.

[Note: \  A homotopy equivalence in \bSISET is a weak homotopy equivalence (but not conversely).]\\

\begingroup%%----------------------------------->>
\fontsize{9pt}{11pt}\selectfont
\label{13.50}
\label{13.59}
\textbf{\small EXAMPLE} \ 
Suppose that 
$
\begin{cases}
\ X\\
\ Y
\end{cases}
$
are topological spaces and $f:X \ra Y$ is a continuous function $-$then there is a \cd
\begin{tikzcd}[sep=large]
{\abs{\sin X}} \ar{d} \ar{r}{\abs{\sin f}} &{\abs{\sin Y}} \ar{d}\\
{X} \ar{r}[swap]{f} &{Y}
\end{tikzcd}
, thus $f$ is a weak homotopy equivalence iff $\sin f$ is a weak homotopy equivalence (Giever-Milnor theorem).]\\
\endgroup%%------------------------------------<<
\label{14.93}

\index{simplicial groups (example)}
\begingroup%%----------------------------------->>
\fontsize{9pt}{11pt}\selectfont
\textbf{\small EXAMPLE \ (\un{Simplicial Groups})} \ 
Given a simplicial group $G$, put $N_nG = \bigcap\limits_{i > 0} \ker d_i$ $(n > 0)$ $(N_0G =G_0)$ and let 
$\partial_n:N_nG \ra N_{n-1}G$ be the restriction $\restr{d_0}{N_nG}$ $(n > 0)$ $(\partial_0:N_0G \ra 0)$ $-$then 
$\im \partial_{n+1}$ is a normal subgroup of $\ker \partial_n$.  
Definition: The 
\un{homotopy groups} 
\index{homotopy groups (Simplicial Groups)} of $G$ are the quotients 
$\pi_n(G) = \ker \partial_n / \im \partial_{n+1}$.   
Justification: $\forall \ n \geq 0$, $\pi_n(G) \approx \pi_n(\abs{G}),e)$.  
Since a homomorphism $f:G \ra K$ of simplicial groups induces a morphism $Nf:NG \ra NK$ of chain complexes, 
thus a homomorphism $\pi_*(f):\pi_*(G) \ra \pi_*(K)$ in homotopy, 
it follows that $f$ is a weak homotopy equivalence iff $\pi_*(f)$ is bijective.
\vspi
[Note: \  A short exact sequence $1 \ra G^\prime \ra G \ra G\pp \ra 1$ of simplicial groups gives rise to a short exact sequence $1 \ra NG^\prime \ra NG \ra NG\pp \ra 1$ of chain complexes and a long exact sequence 
$\cdots \ra$ $\pi_{n+1}(G\pp) \ra$ $\pi_n(G^\prime) \ra$ 
$\pi_n(G) \ra$ $\pi_n(G\pp) \ra$ $\pi_{n-1}(G^\prime) \ra \cdots$ 
of homotopy groups.]\\
\endgroup%%------------------------------------<<

\label{13.100}
\label{13.126}
\begingroup%%----------------------------------->>
\fontsize{9pt}{11pt}\selectfont
\index{simplex categories (example)}
\textbf{\small EXAMPLE\  (\un{Simplex Categories})} \ 
Let $X$ be a simplicial set $-$then $X$ is a cofunctor $\bDelta \ra \bSET$, thus one can form the Grothendieck construction 
$\gro_{\bDelta}X$ on $X$.  
So: The objects of $\gro_{\bDelta}X$ are the $([n],x)$ $(x \in X_n)$ and the morphisms 
$([n],x) \ra ([m],y)$ are the $\alpha:[n] \ra [m]$ such that $(X\alpha)y = x$.  \ 
One calls $\gro_{\bDelta}X$ the 
\un{simplex category} 
\index{simplex category} 
of $X$.  It is isomorphic to the comma category 
$\abs{Y_{\bDelta},K_X}$: 
$
\begin{tikzcd}[sep=small]
{\Delta[n]} \ar{rdd} \ar{rr} &&{\Delta[m]} \ar{ldd}\\
%{\Delta[n]} \ar{rdd} \ar{rr} &&{\Delta[m]} \ar{ldd}{\qquad \cdot}\\
\\
&{X}
\end{tikzcd}
.
$
There is a natural weak homotopy equivalence $\ner(\gro_{\bDelta}X) \ra X$, viz. the rule 
$\nersub_p(\gro_{\bDelta}X) \ra X_p$ that sends 
$([n_0],x_0) \overset{\alpha_0}{\longrightarrow} \cdots \overset{\alpha_{p-1}}{\longrightarrow} ([n_p],x_p)$ 
to $(X\alpha)x_p$, where $\alpha:[p] \ra [n_p]$ is defined by 
$\alpha(i) = \alpha_{p-1} \circx \cdots \circx \alpha_i(n_i)$ $(0 \leq i \leq p)$ $(\alpha(p) = n_p)$.
\vspi
[First check the assertion when $X = \dn$.]\\
\endgroup%%------------------------------------<<

%%----------------------------------------------------------------------------------------------18
A simplicial map $f:X \ra Y$ is said to be a 
\un{cofibration} 
\index{cofibration (simplicial map)} 
if its geometric realization 
$\abs{f}:\abs{X} \ra \abs{Y}$ is a cofibration.  
Example: $\forall \ X$, the arrows 
$
\begin{cases}
\ i_0:X \ra IX\\
\ i_1:X \ra IX
\end{cases}
$
are cofibrations and weak homotopy equivalences.\\

\textbf{\small LEMMA} \ 
The cofibrations in \bSISET are the injective simplicial maps.\\

Example:  Let $X$ be a simplicial set $-$then the arrow of adjunction $X \ra \sin \abs{X}$ is a cofibration and a weak homotopy equivalence 
(cf. p. \pageref{13.26}).\\

\begingroup%%----------------------------------->>
\fontsize{9pt}{11pt}\selectfont
\textbf{\small EXAMPLE} \ 
Let $X$ be a simplicial set $-$then $e_X:X \ra \Ex X$ is a cofibration, as is $e_X^\infty:X \ra \Exx^\infty X$ and both are weak homotopy equivalences 
(cf. p. \pageref{13.27}).\\
\endgroup%%------------------------------------<<

\begin{proposition} \ 
Let $p:X \ra B$ be a simplicial map $-$then $p$ has the RLP w.r.t the incusions $\dot\Delta[n] \ra \Delta[n]$ $(n \geq 0)$ iff $p$ has the RLP w.r.t all cofibrations.
\end{proposition}

[Let $i:A \ra Y$ be an injective simplicial map.  To construct a filler for 
$
\begin{tikzcd}%[sep=large]
{A} \ar{d}[swap]{i} \ar{r}{u} &{X} \ar{d}{p}\\
{Y} \ar{r}[swap]{v} &{B}
\end{tikzcd}
,
$
take $i$ to be an inclusion and call $(Y,A)_n^\#$ the subset of $Y_n^\#$ consisting of those elements which do not belong to $A$ $-$then $\forall \ n$, there is a pushout square 
$
\begin{tikzcd}%[sep=large]
{(Y,A)_n^\# \cdot \ddn} \ar{d} \ar{r} &{Y^{(n-1)} \cup A} \ar{d}\\
{(Y,A)_n^\# \cdot \dn}  \ar{r} &{Y^{(n)} \cup A}
\end{tikzcd}
,
$
so one can construct the arrow $Y \ra X$ by induction.]\\

Given $n \geq 1$, the 
\un{$k^{th}$-horn} 
\index{k$^{th}$-horn}  $\Lambda[k,n]$ 
\index{$\Lambda[k,n]$} 
of $\dn$ $(0 \leq k \leq n)$ is the simplicial subset of $\dn$ defined by the condition that $\Lambda[k,n]_m$ is the set of 
$\alpha:[m] \ra [n]$ whose image does not contain the set $[n] - \{k\}$.  
So: $\abs{\Lambda[k,n]} = \Lambda^{k,n}$ is the subset of $\abs{\dn} = \Delta^n$ consisting of those 
$(t_0, \ldots, t_n)$: $t_i = 0$ $(\exists \ i \neq k)$, thus $\Lambda^{k,n}$ is a strong deformation retract of $\Delta^n$.

Example: 
Let 
$
\begin{cases}
\ X\\
\ Y
\end{cases}
$
be topological spaces, $f:X \ra Y$ a continuous function $-$then $f$ is a Serre fibration iff $f$ has the RLP w.r.t. the inclusions 
$\Lambda^{k,n} \ra \Delta^n$ $(0 \leq k \leq n, n \geq 1)$.\\

\label{13.134}
\begingroup%%----------------------------------->>
\fontsize{9pt}{11pt}\selectfont
The representation of $\dot\Delta[n]$ as a coequalizer can be modified to exhibit $\Lambda[k,n]$ as a coequalizer (in the notation of 
p. \pageref{13.28}, 
replace 
$\coprod\limits_{0 \leq i \leq n} \Delta[n-1]_i$ by 
$\coprod\limits_{\substack{0 \leq i \leq n, \\ i \neq k}}  \Delta[n-1]_i$).  
A corollary is that for every simplicial set $X$, $\Nat(\Lambda[k,n],X)$ is in a one-to-one correspondence with the set of finite sequences 
$(x_0, \ldots, \hat{x}_k, \ldots, x_n)$ of elements of $X_{n-1}$ such that 
$d_ix_j = d_{j-1}x_i$ $(i < j \ \& \ i,\ j \neq k)$.\\
\endgroup%%------------------------------------<<

A retract invariant, composition closed class of injective simplicial maps is said to be 
\un{replete}
\index{replete (class of injective simplicial maps)} 
if it contains the isomorphisms and is stable under formation of coproducts, 
%%----------------------------------------------------------------------------------------------19
pushouts, and sequential colimits.  
The 
\un{repletion} 
\index{repletion (of a set of injective simplicial maps)} 
of a set $S_0$ of injective simplicial maps is $\cap M$, $M$ replete with $S_0 \subset M$.

Specialize to $S_0 = \{\Lambda[k,n] \ra \Delta[n] \ (0 \leq k \leq n, n \geq 1)\}$ $-$then the repletion of $S_0$ is the class of 
\un{anodyne extensions} 
\index{anodyne extensions}.
Examples:
(1) The injections $\Delta[\delta_i]:\Delta[n-1] \ra \Delta[n]$ are anodyne extensions;
(2) The inclusions $\Delta[m] \times \Lambda[k,n] \cup \dot\Delta[m] \times \Delta[n] \ra \Delta[m] \times \Delta[n]$ 
are anodyne extensions.\\

\begin{proposition} \ 08
Let $f:X \ra Y$ be an anodyne extension $-$then $\abs{f}(\abs{X})$ is a strong deformation retract of $\abs{Y}$.
\end{proposition}

[The class of injective simplicial maps with this property is replete (cf. $\S 3$, Proposition 3 and 
p. \pageref{13.29}) 
and contains $S_0$.]\\

Application: Every anodyne extension is a weak homotopy equivalence.\\

\label{13.34} %dmc ?

\begin{proposition} \ %09
Let 
$
\begin{cases}
\ A\\
\ B
\end{cases}
$
be a simplicial subset of 
$
\begin{cases}
\ X\\
\ Y
\end{cases}
. \ 
$
Suppose that the inclusion $B \ra Y$ is an anodyne extension $-$then the inclusion 
$X \times B \cup A \times Y \ra X \times Y$ 
is an anodyne extension.
\end{proposition}

[The class of injective simplicial maps $B^\prime \ra Y^\prime$ for which the arrow 
$X \times B^\prime \underset{A \times B^\prime}{\sqcup} A \times Y^\prime$ $\ra X \times Y^\prime$ 
is an anodyne extension is replete.  On the other hand, an induction shows that the inclusions 
$X \times \Lambda[k,n] \cup A \times \dn \ra X \times \dn$ are anodyne.]\\

\begingroup%%----------------------------------->>
\fontsize{9pt}{11pt}\selectfont
\textbf{\small EXAMPLE} \ 
The inclusion $\Sd\Lambda[k,n] \ra \Sd\dn$ is an anodyne extension.
\vspi
[Note: \  In general, Sd preserves anodyne extensions 
(cf. p. \pageref{13.30}).]\\
\endgroup%%------------------------------------<<

\begingroup%%----------------------------------->>
\fontsize{9pt}{11pt}\selectfont
\textbf{\small FACT} \ 
The class of homotopy classes of anodyne extensions admits a calculus of left fractions.
\vspi
[The point is to show that if $f,g:X \ra Y$ are simplicial maps and if $s:X^\prime \ra X$ is an anodyne extension with 
$f \circx s \simeq g \circx s$, then $\exists$ an anodyne extension $t:Y \ra Y^\prime$ with $t \circx f \simeq t \circx g$.]\\
\endgroup%%------------------------------------<<

Let $p:X \ra B$ be a simplicial map $-$then $p$ is said to be a 
\un{Kan fibration} 
\index{Kan fibration (simplicial map)} 
if it has the RLP w.r.t. the inclusions $\Lambda[k,n] \ra \dn$ $(0 \leq k \leq n, n \geq 1)$.

[Note: \  Let $p:X \ra B$ be a Kan fibration $-$then for any component $A$ of $X$, $p(A)$ is a component of $B$ and $A \ra p(A)$ is a Kan fibration.  Therefore $p(X)$ is a union of components of $B$.  So, if $B$ is connected and $X$ is nonempty, then $p$ is surjective.]

Example: Let 
$
\begin{cases}
\ X\\
\ Y
\end{cases}
$
be topological spaces, $f:X \ra Y$ a continuous function $-$then $f$ is a Serre fibration iff 
$\sin f: \sin X \ra \sin Y$ is a Kan fibration.\\

%%----------------------------------------------------------------------------------------------20
\begingroup%%----------------------------------->>
\fontsize{9pt}{11pt}\selectfont
In ``parameters'', the condition that $p$ be a Kan fibration is equivalent to requiring that if 
$(x_0, \ldots,$ $\hat{x}_k, \ldots, x_n)$ is a finite sequence of elements of $X_{n-1}$ such that 
$d_ix_j = d_{j-1}x_i$ $(i < j \ \& \ i,j \neq k)$ and $p(x_i) = d_ib$ $(b \in B_n)$, then $\exists$ $x \in X_n$: 
$d_ix = x_i  \ (i \neq k)$ with $p(x) = b$.\\
\endgroup%%------------------------------------<<

\begin{proposition} \ %10
Let $p:X \ra B$ be a simplicial map $-$then $p$ is a Kan fibration iff it has the RLP w.r.t every anodyne extension.
\end{proposition}

[The class of injective simplicial maps that have the LLP w.r.t $p$ is replete.]\\

Application: Let $A$ be a simplicial subset of $Y$.  Suppose that $p:X \ra B$ is a Kan fibration $-$then every \cd
\begin{tikzcd}%[sep=large]
{i_0Y \cup IA} \ar{d} \ar{r}{F} &{X} \ar{d}{p}\\
{IY} \ar{r}[swap]{h} &{B}
\end{tikzcd}
has a filler $H:IY \ra X$ (cf. $\S 4$, Proposition 12).

[The vertex inclusion $e_0:\dz \ra \dw$ is anodyne.]\\

\label{13.41}
\label{13.48}
\begingroup%%----------------------------------->>
\fontsize{9pt}{11pt}\selectfont
\textbf{\small FACT} \ 
Let $p:X \ra B$ be a  Kan fibration $-$then $\Ex p:\Ex X \ra \Ex B$ is a Kan fibration.\\
\endgroup%%------------------------------------<<

A simplicial set $X$ is said to be 
\un{fibrant} 
\index{fibrant (simplicial set)} 
if the arrow $X \ra *$  is a Kan fibration.  The fibrant objects are therefore those $X$ such that every simplicial map $f:\Lambda[k,n] \ra X$ 
can be extended to a simplicial map $F:\dn \ra X$ $(0 \leq k \leq n, n \geq 1)$.

[Note: \  The components of a fibrant $X$ are fibrant.]

Example: Let $X$ be a topological space $-$then $\sin X$ is fibrant.\\

\textbf{\small LEMMA} \ 
Suppose that $X$ is fibrant.  
Assume: $\exists \ n_0 \geq 1$ such that $\#(X_{n_0}^\#) \geq 1$ $-$then $\forall \ n \geq n_0$, 
$\#(X_{n}^\#) \geq 1$.

[Fix $x \in X_{n_0}^\#$ and choose $y \in X_{n_0+1}$ such that 
$d_0y = x$, $d_1y = s_0d_0x$.  
Claim: $y \in X_{n_0+1}^\#$.  
Suppose not, so $y = s_iz$ $(\exists \ i)$.  
Case 1: $i \geq 1$: $x = d_0y =$ $d_0s_iz =$ $s_{i-1}d_0z$, an impossibility.  
Case 2: $i = 0$: $x = d_0y =$ $d_0s_0z = z$ $\implies$ $x = z$ $\implies$ $y = s_0x$ 
$\implies$ $d_1y = d_1s_0x$ $\implies$ $x = d_1s_0x = s_0d_0x$, an impossibility.]\\

Application: $\Delta[n]$ $(n \geq 1)$ is not fibrant.\\

Remark: Let $Y$ be a simplicial set $-$then the arrow $Y \ra *$ is a homotopy fibration.\\
Proof: Take any \cd 
\begin{tikzcd}%[sep=large]
{X^\prime \times Y} \ar{d} \ar{r}{\Phi} &{X \times Y} \ar{d} \ar{r} &{Y} \ar{d}\\
{X^\prime} \ar{r}[swap]{\phi} &{X} \ar{r} &{*} 
\end{tikzcd}
, where $\phi$ is a 
%%----------------------------------------------------------------------------------------------21
weak homotopy equivalence, and apply $\abs{?}$ to get a \cd\\
\begin{tikzcd}%[sep=large]
{\abs{X^\prime} \times_k \abs{Y}} \ar{d} \ar{r}{\abs{\Phi}} 
&{\abs{X} \times_k \abs{Y}} \ar{d} \ar{r} &{\abs{Y}} \ar{d}\\
{\abs{X^\prime}} \ar{r}[swap]{\abs{\phi}} &{\abs{X}} \ar{r} &{*} 
\end{tikzcd}
in \bCGH (cf. Proposition 1).  Since the projection 
$\abs{X} \times_k \abs{Y} \ra \abs{X}$ is a \bCG fibration and $\abs{\phi}$ 
is a homotopy equivalence, $\abs{\Phi}$ is a homotopy equivalence 
(cf. p. \pageref{13.31}), 
i.e., $\Phi$ is a homotopy equivalence.

[Note: \  See 
p. \pageref{13.32} 
for the \mc structure on \bSISET.]\\

\label{13.55} %\dmc ?
\begingroup%%----------------------------------->>
\fontsize{9pt}{11pt}\selectfont
\textbf{\small EXAMPLE} \ 
The underlying simplicial set of a simplicial group $G$ is fibrant.
\vspi
[Let $(x_0, \ldots, \hat{x}_k, \ldots, x_n)$ be a finite sequence of elements of $G_{n-1}$ such that 
$d_ix_j = d_{j-1}x_i$ $(i < j \ \& \ i, j \neq k)$.  
Claim: $\exists$ elements $g_{-1},g_0, \ldots \in G_n$ such that $d_ig_r = x_i$ $(i \leq r, i \neq k)$.  
Thus put 
$g_{-1} = e \in G_n$ and assume that $g_{r-1} \in G_n$ has been constructed.  
Case 1: $r = k$.  Take $g_r = g_{r-1}$.
Case 2:  $r \neq k$.  Take $g_r = g_{r-1}(s_rh_r)^{-1}$, where $h_r = x_r^{-1}(d_rg_{r-1})$.]
\vspi
[Note: \  A homomorphism $f:G \ra K$ of simplicial groups is a Kan fibration iff $N_nf:N_nG \ra N_nK$ is surjective $\forall \ n > 0$.  
Therefore a surjective homomorphism of simplicial groups is a Kan fibration.]\\
\endgroup%%------------------------------------<<

\begingroup%%----------------------------------->>
\fontsize{9pt}{11pt}\selectfont
\textbf{\small EXAMPLE} \ 
Let \bC be a small category $-$then $\ner\bC$ is fibrant iff \bC is a groupoid.
\vspi
[Note: \  It is a corollary that $\dn \ (n \geq 1)$ is not fibrant.]\\
\endgroup%%------------------------------------<<

\begingroup%%----------------------------------->>
\fontsize{9pt}{11pt}\selectfont
\textbf{\small LEMMA} \ 
Put $\text{d}_{k,n} = \text{d}_{\Lambda[k,n]}$ $(0 \leq k \leq n, n \geq 1)$ $-$then there is a simplicial map 
$\tD_{k,n}:\Sd^2\dn \ra \Sd\Lambda[k,n]$ 
such that 
$\restr{\tD_{k,n}}{\Sd^2\Lambda[k,n]} = \Sd \hspace{0.03cm} \text{d}_{k,n}$.\\
\endgroup%%------------------------------------<<

\label{13.51}
\label{13.73a}
\begingroup%%----------------------------------->>
\fontsize{9pt}{11pt}\selectfont
\textbf{\small FACT} \ 
For any simplicial set $X$, $\Exx^\infty X$ is fibrant.
\vspi
[Suppose given a simplicial map $f:\Lambda[k,n] \ra \Exx^\infty X$.  
Choose an $r$ such that $f$ factors through 
$\Exx^r X$ and let $g$ be the composite 
$\Lambda[k,n] \ra \Exx^r X \ra \Ex\Exx^r X$ $-$then, under 
$\Nat(\Lambda[k,n],\Ex\Exx^r X) \approx$ $\Nat(\Sd \hspace{0.03cm} \Lambda[k,n],\Exx^r X)$, $g$ corresponds to 
$h:\Sd\Lambda[k,n] \ra \Exx^r X$ and an extension $F:\dn \ra \Exx^\infty X$ of $f$ can be constructed by working with the 
``double adjoint'' of $h \circx \tD_{k,n}$ (it being a simplicial map from $\dn$ to $\Exx^2\Exx^r X$).]\\
\endgroup%%------------------------------------<<

The class of Kan fibrations is pullback stable.  In particular: The fibers of a Kan fibration are fibrant objects.\\

\begin{proposition} \ %11
Let $p:X \ra B$ be a Kan fibration $-$then $B$ fibrant $\implies$ $X$ fibrant and $X$ fibrant + $p$ surjective $\implies$ $B$ fibrant.\\
\end{proposition}

%%----------------------------------------------------------------------------------------------22
\label{13.44}
\label{13.44a}
\label{13.44b}
\label{13.110}
\begin{proposition} \ %12
Suppose that $L \ra K$ is an inclusion of simplicial sets and $X \ra B$ is a Kan fibration $-$then the arrow 
$\map(K,X) \ra \map(L,X) \times_{\map(L,B)} \map(K,B)$ 
is a Kan fibration.
\end{proposition}

[Pass from 
\[
\begin{tikzcd}%[sep=large]
{\Lambda[k,n]} \ar{d} \ar{r} &{\map(K,X)} \ar{d}\\
{\dn} \ar{r} &{\map(L,X) \times_{\map(L,B)} \map(K,B)}
\end{tikzcd}
,
\]
to 
\[
\begin{tikzcd}%[sep=large]
{\Lambda[k,n] \times K \cup \dn \times L} \ar{d}[swap]{i} \ar{r} &{X} \ar{d}\\
{\dn \times K} \ar{r} &{B}
\end{tikzcd}
\]
and note that $i$ is anodyne (cf. Proposition 9).]

[Note: \  Compare this result with its topological analog on 
p. \pageref{13.33}.]\\

\label{13.72}
Application: Let $p:X \ra B$ be a Kan fibration $-$then for any simplicial set $Y$, the postcomposition arrow 
$p_*:\map(Y,X) \ra \map(Y,B)$ is a Kan fibration (cf. $\S 4$, Proposition 5).

[Note: \  Take $B = *$ to see that $X$ fibrant $\implies$ $\map(Y,X)$ fibrant $\forall \ Y$.]\\

Application: Let $i:A \ra X$ be a cofibration $-$then for any fibrant $Y$, the precomposition arrow 
$i^*:\map(X,Y) \ra \map(A,Y)$ is a Kan fibration  (cf. $\S 4$, Proposition 6).\\

\begingroup%%----------------------------------->>
\fontsize{9pt}{11pt}\selectfont
\textbf{\small FACT} \ 
Let $L \ra K$ be an anodyne extension $-$then $\forall$ fibrant $Z$, the arrow 
$[K,Z]_0 \ra [L,Z]_0$ is bijective.
\vspi
[Since $Z$ is fibrant, the arrow $[K,Z]_0 \ra [L,Z]_0$ is surjective, hence bijective 
(cf. p. \pageref{13.34}).]\\
\endgroup%%------------------------------------<< 

\begingroup%%----------------------------------->>
\fontsize{9pt}{11pt}\selectfont
Application: Let $L \ra K$ be an anodyne extension $-$then $\forall$ fibrant $Z$, the arrow 
$\map(K,Z) \ra \map(L,Z)$ is a homotopy equivalence.
\vspi
[For any simplicial set $X$, the inclusion $X \times L \ra X \times K$ is anodyne (cf. Proposition 9).  But 
$[X,\map(K,Z)]_0 \ra [X,\map(LZ)]_0$ is bijective iff $[X \times K,Z]_0 \ra [X \times L,Z]_0$ is bijective.]\\
\endgroup%%------------------------------------<< 

\begin{proposition} \ 
Let $p:X \ra B$ be a Kan fibration.  Suppose that $b^\prime, b\pp \in B_0$ are in the same component of $B$ $-$then the fibers 
$X_{b^\prime}$, $X_{b\pp}$ have the same homotopy type.
\end{proposition}

[Note: \  Compare this result with its topological analog on 
p. \pageref{13.35}.]\\

\textbf{\small LEMMA} \ 
For any fibrant $X$, simplicial homotopy of simplicial maps $\Delta[0] \ra X$ is an equivalence relation.

%
%%----------------------------------------------------------------------------------------------23
[The relation is reflexive: $\forall \ x \in X_0$, $d_1s_0x = x = d_0s_0x$.

The relation is transitive: For suppose that $x \underset{s}{\simeq} y$ $\&$ $y \underset{s}{\simeq} z$ 
$(x, y, z \in X_0)$, say
$
\begin{cases}
\ d_1u = x\\
\ d_0u = y
\end{cases}
$
$(u \in X_1)$ $\&$ 
$
\begin{cases}
\ d_1v = y\\
\ d_0v = z
\end{cases}
$
$(v \in X_1)$.  
The pair $(v,u)$ determines a simplicial map $\Lambda[1,2] \ra X$.  
Extend it to a simplicial map $F:\Delta[2] \ra X$ and put
$w = d_1 F \in X_1$ : 
$d_1w = d_1d_1 F = $ $d_1 d_2 F = x$ 
$\&$ 
$d_0w = d_0d_1 F =$ $d_0d_0F = z$ $(F \leftrightarrow F(\id_{[2]}))$.

The relation is symmetric.  For suppose that $x \underset{s}{\simeq} y$ $(x, y \in X_0)$, say
 $
\begin{cases}
\ d_1u = x\\
\ d_0u = y
\end{cases}
$
$(u \in X_1)$. \ 
The pair $(s_0x, u)$ determines a simplicial map $\Lambda[0,2] \ra X$.  
Extend it to a simplicial map 
$G:\Delta[2] \ra X$ and put 
$v = d_0G$ : $d_1v = d_1d_0G = d_0d_2G = y$ 
$\&$ 
$d_0v = d_0d_0G = d_0d_1G = x$ $(G \leftrightarrow G(\id_{[2]}))$.]

[Note: \  It is a corollary that $\dn$ $(n \geq 1)$ is not fibrant.]\\

Application:  For any fibrant $X$ and any $Y$, simplicial homotopy of simplicial maps $Y \ra X$ is an equivalence relation, so homotopy = simplicial homotopy in this situation.

[In fact, $X$ fibrant $\implies$ $\map(Y,X)$ fibrant $\forall \ Y$ (cf. supra).]\\

Denote by $\iota_n$ the inclusion $\ddn \ra \dn$.  Given a Kan fibration $p:X \ra B$, put 
$\map(\iota_n,p) = \map(\ddn,X) \times_{\map(\ddn,B)} \map(\dn,B)$ 
and let $\iota_n/p$ be the arrow map 
$(\dn,X) \ra \map(\iota_n,p)$ $-$then $\iota_n/p$ is a Kan fibration (cf. Proposition 12).  
Definition: Elements $x^\prime, x\pp \in X_n$ are said to be 
\un{$p$-connected}
\index{p-connected (elements of a simplicial set)} 
$(x^\prime \ \underset{p}{\simeq} x\pp)$ 
if 
$\Delta_{x^\prime}$, $\Delta_{x\pp} \in \map(\dn,X)_0$ 
belong to the same component of the same fiber of $\iota_n/p$.  
Since an element of $\map(\iota_n,p)_0$ is a pair $(f,F)$, where 
$f:\ddn \ra X$, $F:\dn \ra B$ and $p \circx f = F \circx \iota_n$, an element $\Delta_x \in \map(\dn,X)_0$ lies on the fiber 
$\map(\dn,X)_{(f,F)}$ of $\iota_n/p$ over $(f,F)$ if 
$
\begin{cases}
\ p \circx \Delta_x = F\\
\ \Delta_x \circx \iota_n = f
\end{cases}
. \ 
$  
Accordingly, elements $x^\prime, x\pp \in X_n$ with 
$
\begin{cases}
\ p \circx \Delta_{x^\prime}\\
\ p \circx \Delta_{x\pp}
\end{cases}
= F
$
$\& \ $ 
$
\begin{cases}
\ \Delta_{x^\prime} \circx \iota_n\\
\ \Delta_{x\pp} \circx \iota_n
\end{cases}
= f
$
are $p$-connected if 
$\exists$ $H:I\dn \ra X$: $H \circx i_0 = \Delta_{x^\prime}$, $H \circx i_1 = \Delta_{x\pp}$, $p \circx H = F \circx \pr$,  
$\restr{H}{I\ddn} = f \circx \pr$ 
or still, 
$\exists$ $H^\prime, H\pp:I\dn \ra X$: 
$
\begin{cases}
\ H^\prime \circx i_0 = \Delta_{x^\prime}\\
\ H\pp \circx i_0 = \Delta_{x\pp}
\end{cases}
\&
$
$H^\prime \circx i_1 = H\pp \circx i_1$ (or 
$
\begin{cases}
\ H^\prime \circx i_1 = \Delta_{x^\prime}\\
\ H\pp \circx i_1 \Delta_{x\pp}
\end{cases}
\&
$
$H^\prime \circx i_0 = H\pp \circx i_0$), 
$p \circx H^\prime = p \circx H\pp$, 
$\restr{H^\prime}{I\ddn} = \restr{H\pp}{I\ddn}$.

[Note: \   The relation $\underset{p}{\simeq}$ is an equivalence relation on $X_n$.]\\

\textbf{\small LEMMA} \ 
Let $X$ be a simplicial set.  Suppose that $x^\prime$, $x\pp \in X_n$ are degenerate $-$then 
$d_ix^\prime = d_ix\pp$ $(0\leq i \leq n)$ $\implies$ $x^\prime = x\pp$.

[Write 
$x^\prime = s_ky^\prime$, $x\pp = s_ly\pp$.  
Case 1: $k = l$.  
Here, 
$y^\prime =$ 
$d_kx^\prime =$ 
$d_k x\pp =$ 
$y\pp$ $\implies$ $x^\prime = x\pp$.\\
{Case 2:} $k \neq l$, say $k < l$.  
(1) $y^\prime =$ $d_kx^\prime =$ $d_kx\pp =$ $d_ks_ly\pp =$ $s_{l-1}d_ky\pp$;
(2) $x^\prime =$ $s_ky^\prime =$ $s_ks_{l-1}d_ky\pp =$ $s_ls_kd_ky\pp$; 
(3) $y\pp =$ $d_lx\pp =$ $d_lx^\prime =$ $d_ls_ls_kd_ky\pp =$ $s_kd_ky\pp$; 
(4) $x^\prime = s_ly\pp = x\pp$.]\\

%%----------------------------------------------------------------------------------------------24
\label{13.38}
Application: Given a Kan fibration $p:X \ra B$, degenerate elements $x^\prime, x\pp \in X_n$ are $p$-connected iff they are equal.\\

A Kan fibration $p:X \ra B$ is said to be 
\un{minimal} 
\index{minimal (Kan fibraton)} 
if $\forall \ n$, 
$\forall$ $x^\prime, x\pp \in X_n$: $x^\prime \underset{p}{\simeq} x\pp$ $\implies$ $x^\prime = x\pp$.

[Note: \  A fibrant $X$ is minimal when $X \ra *$ is minimal.]\\

\begingroup%%----------------------------------->>
\fontsize{9pt}{11pt}\selectfont
\textbf{\small FACT} \ 
Suppose that $X$ is fibrant $-$then $X$ is minimal iff $\forall \ n$, $\forall \ x^\prime$, $x\pp \in X_n$: $d_i x^\prime = d_i x\pp$ 
$(\forall i \neq j)$ $\implies$ $d_j x^\prime = d_j x\pp$ $(0 \leq i, j \leq n)$.\\
\endgroup%%------------------------------------<<

\begingroup%%----------------------------------->>
\fontsize{9pt}{11pt}\selectfont
\textbf{\small EXAMPLE} \  Let $G$ be a simplicial group $-$then $G$ is minimal iff the chain complex $(NG, \partial)$ is minimal, i.e., iff $\forall \ n$, $\partial_n:N_nG \ra N_{n-1}G$ is the zero homomorphism.\\
\endgroup%%------------------------------------<<

The class of minimal Kan fibrations is pullback stable.  In particular: The fibers of a minimal Kan fibration are minimal fibrant objects.\\

\begin{proposition} \ 
A minimal Kan fibration $p:X \ra B$ is locally trivial.
\end{proposition}

[The claim is that $\forall \ n$ $\&$ $\forall \ b \in B_n$, $X_b$ is trivial over $\dn$.  
Therefore it will be enough to prove that every minimal Kan fibration $p:X \ra \dn$ is trivial.  
To this end, let 
$C_0:\dn \ra \dn$ be the projection onto the $0^{\text{th}}$ vertex and choose a simplicial homotopy 
$H:I\dn \ra \dn$ between 
$C_0$ and $\id_{\dn}$ 
(cf. p. \pageref{13.36}).  
Call \mA the fiber of $p$ over the $0^{\text{th}}$ vertex $-$then there is a retraction 
$r:X \ra A$ and a simplicial homotopy $\ov{H}:IX \ra X$ between $X \overset{r}{\ra} A \ra X$ and $\id_X$ with 
$p \circx \ov{H} = H \circx (p \times \id_{\dw})$.  
Define a simplicial map 
$f:X \ra \dn \times A$ over $\dn$ by $f(x) = (p(x),r(x))$.
To establish that $f$ is an isomorphism, we shall proceed by induction on $k$, 
taking $X_{-1} = \emptyset$ and assuming that 
$\restr{f}{X_l}$ is bijective $(l < k, k \geq 0)$.

Injectivity: \ Suppose that $f(x^\prime) = f(x\pp)$, where $x^\prime$, $x\pp \in X_k$.  
Put $H^\prime(\alpha,t) = \ov{H}((X\alpha)x^\prime,t)$, $H\pp(\alpha,t) = \ov{H}((X\alpha)x\pp,t)$ 
to get simplicial homotopies $H^\prime$, $H\pp: I\dk \ra X$ such that
$
\begin{cases}
\ H^\prime \circx i_1= \Delta_{x^\prime}\\
\ H\pp \circx i_1 = \Delta_{x\pp}
\end{cases}
$
$\&$ $H^\prime \circx i_0 = H\pp \circx i_0$, $p \circx H^\prime = p \circx H\pp$, 
$\restr{H^\prime}{I\ddk} = \restr{H\pp}{I\ddk}$, thus $x^\prime \underset{p}{\simeq} x\pp$, 
so minimality forces $x^\prime = x\pp$.

Surjectivity: \ Let $ (\alpha_0,a_0) \in (\dn \times A)_k$.  
The induction hypothesis, coupled with the injectivity of $f$, ensures the existence of a simplicial map 
$g:\ddk \ra X$ such that $\forall \ \alpha \in \ddk$, $f \circx g(\alpha) = (\alpha_0 \circx \alpha, (X \alpha) a_0)$.  
In addition, one can find a simplicial homotopy 
$G:I\dk \ra X$ satisfying $G \circx i_0 = \Delta_{a_0}$, $\restr{G}{I\ddk} = \ov{H} \circx (g \times \id_{\Delta[1]})$, 
$p \circx G(\alpha,t) = H(\alpha_0 \circx \alpha,t)$.  
Write $\bar{a}_0 = r(x_k)$ $(x_k = G(\id_{[k]},1))$ and set $\ov{G} = \ov{H}\circx (G \circx  i_1 \times \id_{\Delta[1]})$ 
$-$then
$
\begin{cases}
\ G \circx i_0 = \Delta_{a_0}\\
\ \ov{G} \circx i_0 = \Delta_{\bar{a}_0}
\end{cases}
$
%%----------------------------------------------------------------------------------------------25
$\&$ $G \circx i_1 = \ov{G} \circx i_1$, 
$p \circx  G = p \circx \ov{G}$, 
$\restr{G}{I\ddk} = \restr{\overline{G}}{I \ddk}$.
Therefore $a_0 \underset{p}{\simeq} \bar{a}_0$ $\implies$ $a_0 = \bar{a}_0$ $\implies$ $f(x_k) = (\alpha_0,a_0)$.]\\

\label{13.40}
Application: The geometric realization of a minimal Kan fibration is a Serre fibration 
(cf. p. \pageref{13.37}).\\

Let $p:X \ra B$ be a Kan fibration; let $A$ be a simplicial subset of $X$, $i:A \ra X$ the inclusion.

\indent\indent (DR) \  \mA is said to be a 
\un{deformation retract} 
\index{deformation retract (simplicial set)} 
of $X$ over $B$ if there is a simplicial map 
$r:X \ra A$ over $B$ and a simplicial homotopy $H:IX \ra X$ over $B$ such that 
$r \circx i = \id_A$ 
and 
$H \circx i_0 = i \circx r$, 
$H \circx i_1 = \id_X$.

\indent\indent (SDR) \  \mA is said to be a 
\un{strong deformation retract} 
\index{strong deformation retract (simplicial set)} 
of $X$ over $B$ if there is a simplicial map 
$r:X \ra A$ over $B$ and a simplicial homotopy $H:IX \ra X$ over $B$ such that 
$r \circx i = \id_A$ 
and 
$H \circx i_0 = i \circx r$, 
$H(a,t) = a$ $(a \in A)$, 
$H \circx i_1 = \id_X$.

[Note: \  Taking $B = *$ leads to the corresponding absolute notions for fibrant objects.]

If $p:X \ra B$ is Kan and $A \subset X$ is a retract of $X$ over \mB, then the restriction $p_A = \restr{p}{A}$ is Kan.\\

\begingroup%%----------------------------------->>
\fontsize{9pt}{11pt}\selectfont
\textbf{\small FACT} \ 
Let $p:X \ra B$ be a Kan fibration.  Suppose that $A \subset X$ is a deformation retract of $X$ over $B$ $-$then $p$ has the RLP w.r.t. every cofibration that has the LLP w.r.t $p_A$.\\
\endgroup%%------------------------------------<<

\begin{proposition} \ %15
Let $p:X \ra B$ be a Kan fibration $-$then there is a simplicial subset $A \subset X$ which is a strong deformation retract of $X$ over $B$ such that $p_A$ is a minimal Kan fibration.
\end{proposition}

[Let $E$ be a set of representatives for the equivalence classes per $\underset{p}{\simeq}$ containing the degenerate elements of $X$ 
(cf. p. \pageref{13.38}).  
Choose a simplicial subset $A \subset X$ maximal with respect to $A \subset E$: $p_A$ will be minimal if it is Kan.  
Consider the set $\sY$ of all pairs $(Y,G)$, where $A \subset Y \subset X$ and $G:IY \ra X$ is a simplicial homotopy over $B$ such that $G(i_0(Y)) \subset A$, $G(a,t) = a$ $(a \in A)$, $G \circx i_1 = Y \ra X$.
Example: $(A,IA \overset{\pr}{\lra} A \ra X) \in \sY$.  
Order $\sY$ by stipulating that
$(Y^\prime,G^\prime) \leq (Y\pp,G\pp)$ iff $Y^\prime \subset Y\pp$ $\&$ $\restr{G\pp}{IY^\prime} = G^\prime$.  
Every chain in $\sY$ has an upper bound, so by Zorn, $\sY$ has a maximal element $(Y_0,G_0)$.  
Claim: $Y_0 = X$.  
Supposing this is false, take $x \in X_n$: $x \notin Y_0$, with $n$ minimal.  
Note that $x$ is nondegenerate.  
Call $Y_x$ the smallest simplicial subset of $X$: $Y_0 \subset Y_x$ $\&$ $x \in Y_x$.  
Since $\restr{\Delta_x}{\ddn}$ factors through $Y_0$, there is a pushout square 
$
\begin{tikzcd}%[sep=large]
{\ddn} \ar{d} \ar{r} &{Y_0} \ar{d}\\
{\dn} \ar{r}[swap]{\Delta_x} &{Y_x}
\end{tikzcd}
.
$
Fix a simplicial homotopy $H_x:I\dn \ra X$ over 
%%----------------------------------------------------------------------------------------------26
\mB such that $H_x \circx i_1 = \Delta_x$ and 
$\restr{H_x}{I\ddn} = G_0 \circx (\restr{\Delta_x}{\ddn} \times \id_{\dw})$.  
Put $x\pp = H_x(\id_{[n]},0)$ and define $x^\prime \in E$ via $x^\prime \underset{p}{\simeq} x\pp$: $d_ix\pp \in A$ 
$(0 \leq i \leq n)$ $\implies$ $x^\prime \in A$.  
Fix a simplicial homotopy 
$H:I\dn \ra X \rel \ddn$ over $B$ such that $H \circx i_0 = \Delta_{x^\prime}$, $H \circx i_1 = \Delta_{x\pp}$.  
Determine a simplicial map $K:I^2\dn \ra X$ satisfying $p \circx K(\alpha,t, T) = p((X\alpha)x)$,
$
\begin{cases}
\ K(\alpha,t, 1) = H_x(\alpha, t)\\
\ K(\alpha, 0,T) = H(\alpha, T)
\end{cases}
, \ 
$
$K(\alpha,t,T) = G_0((X\alpha)x,t)$ $(\alpha \in \ddn)$, and $K(\alpha,1,T) = (X\alpha)x$.  
Extend $G_0$ to a simplicial homotopy $G_x:IY_x \ra X$ $(G_x(x,t) = K(\id_{[n]}, t, 0))$ : $(Y_x,G_x) \in \sY$.  
Contradiction.]\\

\label{13.43} %dmc mnft
\begingroup%%----------------------------------->>
\fontsize{9pt}{11pt}\selectfont
\textbf{\small LEMMA} \ 
Let $f,\  g:X \ra Y$ be simplicial maps, where $f \underset{s}{\simeq} g$ $\&$ 
$
\begin{cases}
\ X\\
\ Y
\end{cases}
$
are fibrant.  Assume: $f$ is an isomorphism and $Y$ is minimal $-$then $g$ is an isomorphism.\\
\endgroup%%------------------------------------<<

\begingroup%%----------------------------------->>
\fontsize{9pt}{11pt}\selectfont
Application: A simplicial homotopy equivalence between minimal fibrant objects is an isomorphism.\\
\endgroup%%------------------------------------<<

\begingroup%%----------------------------------->>
\fontsize{9pt}{11pt}\selectfont
Consequently, if $X$ is fibrant and if 
$
\begin{cases}
\ A^\prime\\
\ A\pp
\end{cases}
$
are deformation retracts of $X$ that are minimal, then 
$
\begin{cases}
\ A^\prime\\
\ A\pp
\end{cases}
$
are isomorphic.\\[.25cm]
\endgroup%%------------------------------------<<

\label{13.52}
A simplicial map $p:X \ra B$ which has the RLP w.r.t. the inclusions $\ddn \ra \dn$  $(n \geq 0)$ is a Kan fibration (cf. Proposition 7).  Moreover, $p$ is a simplicial homotopy equivalence.  
Proof: $p$ admits a section $s:B \ra X$ and 
\begin{tikzcd}%[sep=large]
{X \times \ddw} \ar{d} \ar{r}{u} &{X} \ar{d}{p}\\
{X \times \dw} \ar{r} &{B}
\end{tikzcd}
admits a filler $H:X \times \dw \ra X$.  Here, $u(x,0) = s(p(x))$, $u(x,1) = x$.\\

\begin{proposition} \ %16
Let $p:X \ra B$ be a Kan fibration $-$then $p$ can be written as the composite of a simplicial map which has the RLP w.r.t. the inclusions  $\dot\Delta[n] \ra \Delta[n]$ $(n \geq 0)$ and a minimal Kan fibration.
\end{proposition}

[Using the notation of Proposition 15, write $p = p_A \circx r$, $r:X \ra A$ the retraction.  Suppose given a \cd
\begin{tikzcd}%[sep=large]
{\ddn} \ar{d} \ar{r}{u} &{X} \ar{d}{r}\\
{\dn} \ar{r}[swap]{v} &{A}
\end{tikzcd}
.  
Since $A$ is a strong deformation retract of $X$ over $B$, there is a simplicial homotopy $H:IX \ra X$ over $B$ such that 
$H \circx i_0 = i \circx r$, $H(a,t) = a$ $(a \in A)$, $H \circx i_1 = \id_X$.  
Choose a simplicial homotopy 
$G:I\dn \ra X$ subject to 
$G(\alpha,0) = v(\alpha)$,  
$\restr{G}{I\ddn} = H \circx (u \times \id_{\dw})$,  
$p \circx G(\alpha,t) = p(v(\alpha))$.
Let 
$\ov{G}(\alpha,t) = H((X\alpha)\ov{x},t)$, where $\ov{x} = G(\id_{[[n]},1)$.  
Put 
$
\begin{cases}
\ a^\prime = v(\id_{[n]})\\
\ a\pp = r(\bar{x})
\end{cases}
, \ 
$
$
\begin{cases}
\ G^\prime = r \circx G\\
\ G\pp = r \circx \overline{G}
\end{cases}
$
%%----------------------------------------------------------------------------------------------27
$-$then \ 
$
\begin{cases}
\ G^\prime \circx i_0 = \Delta_{a^\prime}\\
\ G\pp \circx i_0 = \Delta_{a\pp}
\end{cases}
$
$\&$ $G^\prime \circx i_1 = G\pp \circx i_1$, $p_A \circx G^\prime = p_A \circx G\pp$, 
$\restr{G^\prime}{I\ddn} = \restr{G\pp}{I\ddn}$.  
So: 
$a^\prime \underset{p_A}{\simeq} a\pp$ $\implies$ $a^\prime = a\pp$ (by minimality), hence 
$\Delta_{\bar{x}}:\dn \ra X$ is our filler.]\\

\textbf{\small LEMMA} \  
Suppose that $p:X \ra B$ has the RLP w.r.t the inclusions $\dot\Delta[n] \ra \Delta[n]$ $(n \geq 0)$ $-$then 
$\abs{p}: \abs{X} \ra \abs{B}$ is a \bCG fibration, thus is Serre 
(cf. p. \pageref{13.39}).

[Consider a filler $X \times B \ra X$ for 
\begin{tikzcd}%[sep=large]
{X} \ar{d} \ar{r}{\id_X} &{X} \ar{d}{p}\\
{X \times B} \ar{r}&{B}
\end{tikzcd}
, bearing in mind that $\abs{X \times B} \approx \abs{X} \times_k \abs{B}$.]\\

\begin{proposition} \ 
The geometric realization of a Kan fibration is a Serre fibration.
\end{proposition}

[This follows from Proposition 16, the lemma, and the fact that the geometric realization of a minimal Kan fibration is a Serre fibration 
(cf. p. \pageref{13.40}).]

[Note: \  The argument proves more: The geometric realization of a Kan fibration is a \bCG fibration.]\\

\label{13.53}
\begingroup%%----------------------------------->>
\fontsize{9pt}{11pt}\selectfont
For instance, suppose that $p:X \ra B$ is Kan and a weak homotopy equivalence.  Let $B^\prime \ra B$ be a simplicial map and define 
$X^\prime$ by the pullback square 
\begin{tikzcd}[sep=large]
{X^\prime} \ar{d}[swap]{p^\prime} \ar{r} &{X} \ar{d}{p}\\
{B^\prime} \ar{r}&{B}
\end{tikzcd}
$-$then $p^\prime$ is Kan and a weak homotopy equivalence.\\
\endgroup%%------------------------------------<< 

Suppose that $X$ is fibrant $-$then $X$ is said to be 
\un{simplicially contractible} 
\index{simplicially contractible} 
if the projection $X \ra *$ is a simplicial homotopy equivalence.\\

\begingroup%%----------------------------------->>
\fontsize{9pt}{11pt}\selectfont
\textbf{\small EXAMPLE} \ 
Let $X$ be fibrant $-$then Ex $X$ is fibrant 
(cf. p. \pageref{13.41}) 
and is simplicially contractible if this is so of X.
\vspi
[Recall that Ex preserves simplicial homotopy equivalences 
(cf. p. \pageref{13.42}).]\\
\endgroup%%------------------------------------<< 

\begin{proposition} \ %18
A fibrant $X$ is simplicially contractible iff every simplicial map $f:\dot\Delta[n] \ra X$ can be extended to a simplicial map 
$F:\Delta[n] \ra X$ $(n \geq 0)$.
\end{proposition}

[The stated extension property implies that $X$ is fibrant and simplicially contractible 
(cf. p. \pageref{13.43}).  
To deal with the converse, fix a section $s:\dz \ra X$ for $p:X \ra \dz$ and a simplicial homotopy $H:IX \ra X$ between 
$s \circx p$ and $\id_X$.  Given $f:\ddn \ra X$, choose $G:I\dn \ra X$ such that 
$G \circx i_0 = s \circx (\dn \ra \dz)$, $\restr{G}{I\ddn} = H \circx (f \times \id_{\dw})$ and put 
$F = G \circx i_1$ $-$then $\restr{F}{\ddn} = f$.]\\

%%----------------------------------------------------------------------------------------------28
A 
\un{simplicial pair} 
\index{simplicial pair} 
is a pair $(X,A)$, where $X$ is a simplicial set and $A \subset X$ is a simplicial subset.  
Example: Fix $x_0 \in X_0$ and, in an abuse of notation, let $x_0$ be the simplicial subset of $X$ generated by $x_0$ so that
$(x_0)_n = \{s_{n-1} \cdots s_0x_0\}$ $(n \geq 1)$ $-$then $(X,x_0)$ is a simplicial pair.

A 
\un{pointed simplicial set}
\index{pointed simplicial set} 
is a simplicial pair $(X,x_0)$.  
A 
\un{pointed simplicial map}
\index{pointed simplicial map} 
is a base point preserving simplicial map $f:X \ra Y$, i.e., a simplicial map $f:X \ra Y$ for which the triangle 
\begin{tikzcd}[sep=small]
&{\dz} \ar{ldd}[swap]{\Delta_{x_0}} \ar{rdd}{\Delta_{y_0}}\\
\\
{X} \ar{rr}[swap]{f} &&{Y}
\end{tikzcd}
commutes or, in brief, $f(x_0) = y_0$.

$\bSISET_*$  is the category whose objects are the pointed simplicial sets and whose morphisms are the pointed simplicial maps.  
Thus $\bSISET_* = [\bDelta^\OP,\bSET_*]$ and the forgetful functor $\bSISET_* \ra \bSISET$ has a left adjoint that sends a simplicial set $X$ to the pointed simplicial set $X_+ = X \coprod *$.

[Note: \  The vertex inclusion $e_0: \Delta[0] \ra \Delta[1]$ defines the base point of $\Delta[1]$, hence of $\dot\Delta[1]$.]\\

\begingroup%%----------------------------------->>
\fontsize{9pt}{11pt}\selectfont
$\Delta[0]$ is a zero object in $\bSISET_*$ and $\bSISET_*$ has the obvious products and coproducts.  In addition, the pushout square
\begin{tikzcd}[sep=large]
{X \vee Y} \ar{d}\ar{r} &{\dz} \ar{d}\\
{X \times Y} \ar{r} &{X\#Y}
\end{tikzcd}
defines the 
\un{smash product}
\index{smash product (in $\bSISET_*$)} 
$X \# Y$.  
Therefore $\bSISET_*$ is a closed category if
$X \otimes Y = X \# Y$ and $e = \ddw$.  
Here, the internal hom functor sends $(X,Y)$ to $\map_*(X,Y)$, the simplicial subset of $\map(X,Y)$ 
whose elements in degree $n$ are the $f:X \times \dn \ra Y$ with $f(x_0 \times \dn) = y_0$, 
i.e., the pointed simplicial maps $X\#\dn_+ \ra Y$, the zero morphism $0_{XY}$ being the base point.\\
\endgroup%%------------------------------------<< 

\begingroup%%----------------------------------->>
\fontsize{9pt}{11pt}\selectfont
\textbf{\small FACT} \ 
Let $i:A \ra X$ be a pointed cofibration $-$then for any pointed fibrant $Y$, the precomposition arrow 
$i^*:\map_*(X,Y) \ra \map_*(A,Y)$ is a Kan fibration.
\vspi
[Consider the pullback square
\begin{tikzcd}[sep=large]
{\map_*(X,Y)} \ar{d} \ar{r} &{\map(X,Y)} \ar{d}\\
{\map_*(A,Y)} \ar{r}&{\map(A,Y)}
\end{tikzcd}
, \ recalling that the arrow $\map(X,Y) \ra \map(A,Y)$ is a Kan fibration 
(cf. p. \pageref{13.44}).]\\
\endgroup%%------------------------------------<< 

\begingroup%%----------------------------------->>
\fontsize{9pt}{11pt}\selectfont
Application: Fix a pointed fibrant $Y$ $-$then $\forall$ pointed $X$, $\map_*(X,Y)$ is fibrant.\\
\endgroup%%------------------------------------<< 

Suppose that $X$ is fibrant.  
Fix $x_0 \in X_0$ $-$then the 
\un{mapping space} 
\index{mapping space (simplicial sets)} 
$\Theta X$ of the pointed simplicial set $(X,x_0)$ is defined by the pullback square
$
\begin{tikzcd}%[sep=large]
{\Theta X} \ar{d} \ar{r} &{\map(\dw,X)} \ar{d}{e_0^*}\\
{\dz} \ar{r}[swap]{\Delta_{x_0}} &{\map(\dz,X)} \approx {X}
\end{tikzcd}
.$
%%----------------------------------------------------------------------------------------------29
Since $X$ is fibrant, $e_0^*$ is a Kan fibration 
(cf. p. \pageref{13.44a}), 
hence $\Theta X$ is fibrant.  
Furthermore, the composite 
$\Theta X \ra \map(\dw,X) \overset{e_1^*}{\lra} \map(\dz,X) \approx X$ is a Kan fibration, call it $p_1$.  
Proof: Consider the pullback square
\begin{tikzcd}%[sep=large]
{\Theta X} \ar{d}[swap]{p_1} \ar{r} &{\map(\dw,X)} \ar{d}{i^*}\\
{X}  \ar[bend right]{rr}[swap]{(\Delta_{x_0},\id_X)} &{\map(\ddw,X)} \ \ \approx &{X \times X}\\
\end{tikzcd}
, noting that $i^*$ is a Kan fibration 
(cf. p. \pageref{13.44b}).\\
\vspace{0.25cm}

\begingroup%%----------------------------------->> 
\fontsize{9pt}{11pt}\selectfont
$\Theta X$ can be identified with $\map_*(\dw,X)$, thus is a pointed simplicial set.  
The fiber of $p_1:\Theta X \ra X$ over the base point is the 
\un{loop space} 
\index{loop space (simplicial set)} $\Omega X$, i.e., 
$\map_*(\bS[1],X)$, $\bS[1] = \bDelta[1]/\dot\bDelta[1]$, the simplicial circle.  
\label{14.170}
Example: $\forall$ pointed topological space $X$, there are natural isomorphisms 
$\Theta (\sin X) \approx \sin\Theta X$, 
$\Omega(\sin X) \approx \sin\Omega X$.\\
\endgroup%%------------------------------------<< 

\textbf{\small LEMMA} \ 
$e_0^*:\map(\dw,X) \ra \map(\dz,X)$ has the RLP w.r.t. the inclusions $\ddn \ra \dn$ $(n \geq 0)$.

[Convert 
\begin{tikzcd}%[sep=large]
{\ddn} \ar{d} \ar{r} &{\map(\dw,X)} \ar{d}{e_0^*}\\
{\dn} \ar{r} &{\map(\dz,X)}
\end{tikzcd}
\ \ to \ \ 
\begin{tikzcd}%[sep=large]
{\dn \times \dz \cup \ddn \times \dw} \ar{d} \ar{r} &{X} \ar{d}\\
{\dn \times \dw} \ar{r} &{*}
\end{tikzcd}
\vspace{0.15cm}
, bearing in mind that $e_0:\dz \ra \dw$ is anodyne.]\\

\begin{proposition} \ %19
Suppose that $X$ is fibrant $-$then $\Theta X$ is simplicially contractible.
\end{proposition}

[In view of the lemma, this is a consequence of Proposition 18.]\\

\begingroup%%----------------------------------->>
\fontsize{9pt}{11pt}\selectfont
\label{13.14}
\textbf{\small LEMMA} \ 
For every simplicial set $X$, $\abs{e_X}:\abs{X} \ra \abs{\Ex X}$ is a homotopy equivalence 
(cf. p. \pageref{13.45}).\vspi

[Show that $\abs{e_X}$ is bijective on $\pi_0$ and $\pi_1$ and, using an acyclic models argument, that $\abs{e_X}$ is a homology equivalence.  
To handle the higher homotopy groups, define $\Theta X$ by the pullback square
\begin{tikzcd}[sep=large]
{\Theta X} \ar{d} \ar{r} &{\Theta \sin\abs{X}} \ar{d}{p_1}\\
{X} \ar{r} &{\sin\abs{X}}
\end{tikzcd}
.  
Since $X \ra \sin\abs{X}$ is a weak homotopy equivalence
 (cf. p. \pageref{13.46}), 
 the same is true of $\Theta X \ra \Theta \sin\abs{X}$ 
(cf. p. \pageref{13.47}). 
But $\Theta \sin\abs{X}$ is simplicially contractible (cf. Proposition 19), thus $\Theta X \ra *$ is a weak homotopy equivalence and so 
$\Ex \Theta X \ra *$ is a weak homotopy equivalence.  
In addition: $\Theta X \ra X$ Kan $\implies$ $\Ex \Theta X \ra \Ex X$ Kan 
(cf. p. \pageref{13.48}).  
Compare the homotopy sequences of the associated Serre fibrations and use induction.]\\
\endgroup%%------------------------------------<<

\index{Simplicial Extension Theorem}
\textbf{\small SIMPLICIAL EXTENSION THEOREM} \ 
Let $(K,L)$ be a simplicial pair, $p:X \ra B$ a Kan fibration.  Suppose given a \cd
\begin{tikzcd}%[sep=large]
{L} \ar{d} \ar{r}{g} &{X} \ar{d}{p}\\
{K} \ar{r}[swap]{f} &{B}
\end{tikzcd}
$-$the $\forall \ \phi \in
%%----------------------------------------------------------------------------------------------30
C(\abs{K},\abs{X})$ such that $\restr{\phi}{\abs{L}} = \abs{g}$ and $\abs{p} \circx \phi = \abs{f}$ there is a simplicial map 
$F:K \ra X$ with $\restr{F}{L} = g$, $p \circx F = f$, and 
$\abs{F} \underset{\abs{B}}{\simeq} \phi \ \rel \abs{L}$.

[It will be enough to consider the case when 
$
\begin{cases}
\ K = \dn\\
\ L = \ddn
\end{cases}
(n \geq 0).  
$
Identify
$
\begin{cases}
\ X\\
\ B
\end{cases}
$
with its image in 
$
\begin{cases}
\ \sin \abs{X}\\
\ \sin \abs{B}
\end{cases}
$
under the arrow of adjunction 
$
\begin{cases}
\ X \ra \sin \abs{X}\\
\ B \ra \sin \abs{B}
\end{cases}
$
, so that $\phi \in C(\dpn,\abs{X}) = \sin_n \abs{X}$, $d_i\phi \in X$ $(0 \leq i \leq n)$, 
$b_\phi = \abs{p} \circx \phi \in B_n$.  
The assertion can thus be recast: $\exists$ $x \in X_n$ such that $x \underset{\sin\abs{p}}{\simeq} \phi$.  
This being clear if $n = 0$, take $n > 0$, write $b_\phi= (B\beta)b$, where $\beta$ is an epimorphism and $b$ is nondegenerate, and argue inductively on $n$ and on the finite set of epimorphisms having domain $[n]$ 
(viz., $\beta^\prime \leq \beta\pp$ iff $\forall \ i$, $\beta^\prime(i) \leq \beta\pp(i)$).

[Note: \  $p$ Kan $\implies$ $\abs{p}$ Serre (cf. Proposition 17) $\implies$ $\sin\abs{p}$ Kan.]

\indent\indent (I) \  $\beta:[n] \ra [0]$. Here, $b \in B_0$ and $d_i\phi \in X_b$ $(0 \leq i \leq n)$.  
View $X_b$ (which is fibrant) as a pointed simplicial set with base point $\phi_0$ 
(the $0^{\text{th}}$ element in the vertex set of $\phi$ 
(cf. p. \pageref{13.49})).  
Put $Y = X_b$, $W = \Theta Y$, $q = p_1$, and choose a finite sequence 
$(w_0, \ldots, w_{n-1}, \widehat{w}_n)$ of elements of $W_{n-1}$ such that $d_iw_j = d_{j-1}w_i$ $(i < j \ \& \ i, j \neq n)$ 
with $q(w_i) = d_i\phi$ $(0 \leq i \leq n-1)$ ($q$ maps $W$ surjectively onto the component of $Y$ containing the base point).  
Encode the data in the \cd 
\begin{tikzcd}%[sep=large]
{\Lambda[n,n]} \ar{d} \ar{r} &{\sin\abs{W}} \ar{d}\\
{\dn} \ar{r}[swap]{\Delta_\phi} &{\sin\abs{Y}}
\end{tikzcd}
to produce a $\psi \in \sin_n\abs{W}$ : $\sin\abs{q}(\psi) = \phi$.  
The induction hypothesis furnishes a $w_n \in W_{n-1} \ :$ $w_n \underset{\sin\abs{q}}{\simeq} d_n\psi$.  
On the other hand, $W$ is simplicially contractible (cf. Proposition 19), so one can find a 
$w \in W_n:$ $d_iw = w_i$ $(0 \leq i \leq n)$ (cf. Proposition 18).  
Claim:  $x \underset{\sin\abs{p}}{\simeq} \phi$, where $x = q(w)$.  
To see this, fix a simplicial homotopy $H:I\Delta[n-1] \ra \sin\abs{W}\rel\dot\Delta[n-1]$ over $\sin\abs{Y}$ 
such that 
$H \circx i_0 = \Delta_{w_n}$, 
$H \circx i_1 = \Delta_{d_n\psi}$.  \ 
Define a simplicial map 
$\ov{H}:\dn \times \dot{\Delta}[1]$ $\cup$ $\ddn \times \dw \ra \sin\abs{W}$ by the recipe 
$\ov{H} \circx i_0 = \Delta_{{w}}$,
$\ov{H} \circx i_1 = \Delta_\psi$, 
$\ov{H}(d_i\id_{[n]},t) = w_i$ $(0 \leq i \leq n-1)$, 
$\ov{H}(d_{n} \id_{[n]},t) = H(\id_{[n-1]},t)$.  
Using the fact that $(\abs{?},\sin)$ is an adjoint pair, $\ov{H}$ determines a continuous function 
$\ov{G}:i_0\Delta^n \cup i_1 \Delta^n \cup I \dot\Delta^n \ra \abs{W}$ 
which can then be extended to a continuous function 
$\widetilde{G}:I\Delta^n \ra \abs{W}$ ($\abs{W}$ is contractible).  
Pass back to get a simplicial homotopy $\widetilde{H}:I\dn \ra \sin\abs{W}$ extending $\ov{H}$.  
Consider the composite $\sin \abs{q} \circx \widetilde{H}$ followed by the inclusion $\sin\abs{Y} \ra \sin\abs{X}$.

\indent\indent (II) \  $\beta:[n] \ra [m]$ $(m > 0)$. Let $k = \min\limits_{0 \leq i \leq n} i: \beta(i) \neq \beta(i+1)$.  
Choose $\ov{x} \in X_n$: $d_i\ov{x} = d_i \phi$ $(0 \leq i \leq n-1)$ with $p(\ov{x}) = b_\phi$ and choose 
$\psi \in \sin_{n+1}\abs{X}$: $d_k\psi = \ov{x}$, $d_{k+1}\psi = \phi$, $d_i\psi = d_is_k\phi$ 
$(0 \leq i \leq n, i \neq k, k+1)$ with $\abs{p} \circx \psi = s_k b_\phi$ $-$then $\exists$ $\ov{y} \in X_n$ : 
$\ov{y} \underset{\sin\abs{p}}{\simeq} d_{n+1}\psi$ (induction).  
Choose $\ov{w} \in X_{n+1}$ : 
$d_k\ov{w} = \ov{x}$, 
$d_{n+1}\ov{w} = \ov{y}$, 
$d_i\ov{w} = d_is_k\phi$  $(0 \leq i \leq n, i \neq k , k + 1)$ with
$p(\ov{w}) = s_kb_\phi$.  
Fix a simplicial homotopy $H:I\dn \ra \sin\abs{X} \rel \ddn$ over $\sin\abs{B}$ such
%%----------------------------------------------------------------------------------------------31
that $H \circx i_0 = \Delta_{\bar{y}}$, $H \circx i_1 = \Delta_{d_{n+1}\psi}$ 
and incorporate the choices into a simplicial homotopy 
$\ov{H}:I \Delta[n+1] \ra \sin\abs{X}$ satisfying 
$\ov{H} \circx i_0 = \Delta_{\ov{w}}$,
$\ov{H} \circx i_1 = \Delta_\psi$, 
$\ov{H}(d_i\id_{[n+1]},t) = d_i\bar{w}$ $(0 \leq i \leq n, i \neq k + 1)$, 
$\ov{H}(d_{n+1} \id_{[n+1]},t) = H(\id_{[n]},t),\abs{p} \circx \ov{H}(\id_{[n+1]},t) = s_kb_\phi$.  
Put $x = d_{k+1}\ov{w}$ and examine 
$\ov{H} \circx (\Delta[\delta_{k+1}] \times \id_{\dw}):I\dn \ra \sin\abs{X}$ 
to conclude that $x \underset{\sin\abs{p}}{\simeq} \phi$.]\\

\begingroup%%----------------------------------->> 
\fontsize{9pt}{11pt}\selectfont
Specialized to $B = *$, one can say that if $(K,L)$ is a simplicial pair and $X$ is fibrant, then given a simplicial map 
$g:L \ra X$ and a continuous extension $\phi:\abs{K} \ra \abs{X}$ of $\abs{g}$, there exists a simplicial extension 
$F:K \ra X$ of $g$ such that $\abs{F} \simeq \phi \ \rel\abs{L}$.  
Conversely, every simplicial set $X$ with this property is fibrant.  
Proof: The geometric realization of a simplicial map $\Lambda[k,n] \ra X$ can be extended to a continuous function 
$\Delta^n \ra \abs{X}$.\\
\endgroup%%------------------------------------<<

Example: Suppose that $X$ is fibrant $-$then $X$ is a strong deformation retract of $\sin\abs{X}$.
[Apply the simplicial extension theorem to the \cd
$
\begin{tikzcd}%[sep=large]
{X} \ar{d} \ar{r}{\id_X} &{X} \ar{d}\\
{\sin\abs{X}} \ar{r} &{*}
\end{tikzcd}
,
$
taking for $\phi \in C(\abs{\sin\abs{X}},\abs{X})$ the arrow of adjunction $\abs{\sin\abs{X}} \ra \abs{X}$.]\\

\begingroup%%----------------------------------->>
\fontsize{9pt}{11pt}\selectfont
\textbf{\small EXAMPLE} \ 
Let 
$
\begin{cases}
\ X\\
\ Y
\end{cases}
$
be simplicial sets.  Assume: $Y$ is fibrant $-$then there is a weak homotopy equivalence 
$\abs{\map(X,Y)} \ra \map(\abs{X},\abs{Y})$.
\vspi
[Since $Y$ is fibrant, the arrow of adjunction $Y \ra \sin\abs{Y}$ is a simplicial homotopy equivalence, thus the arrow 
$\map(X,Y) \ra$ $\map(X,\sin\abs{Y})$ is a simplicial homotopy equivalence.  But 
$\map(X,\sin\abs{Y}) \approx$ $\sin\map(\abs{X},\abs{Y})$ 
and the arrow of adjunction 
$\abs{\sin\map(\abs{X},\abs{Y})} \ra \map(\abs{X},\abs{Y})$ is a weak homotopy equivalence (Giever-Milnor theorem).]\\
\endgroup%%------------------------------------<<

\begin{proposition} \ %20
Let 
$
\begin{cases}
\ X\\
\ Y
\end{cases}
$
 be fibrant $-$then a simplicial map $f:X \ra Y$ is a simplicial homotopy equivalence iff its geometric realization 
$\abs{f}:\abs{X} \ra \abs{Y}$ is a homotopy equivalence.
\end{proposition}

[In general, geometric realization takes simplicial homotopy equivalences to homotopy equivalences.  The fibrancy of $X$ $\&$ $Y$ is used to go the other way.  Thus fix a homotopy inverse 
$g:\abs{Y} \ra \abs{X}$ for $\abs{f}$ and let $r:\sin\abs{X} \ra X$ be a simplicial homotopy inverse for 
$X \ra \sin\abs{X}$ (cf. supra) $-$then the composite 
%$Y \ra \sin\abs{Y} \overset{\sin g}{\longrightarrow} \sin\abs{X} \overset{r}{\ra} X$ 
$Y \ra \sin\abs{Y}$ 
$\overset{\sin g}{\xrightarrow{\hspace*{1cm}}}$ 
$\sin\abs{X} \overset{r}{\ra} X$
%$X$ 
is a simplicial homotopy inverse for $f$.]

[Note: \  It is a corollary that a fibrant $X$ is simplicially contractible iff $\abs{X}$ is contractible.]\\

Application: Suppose that
$
\begin{cases}
\ X\\
\ Y
\end{cases}
$
are topological spaces and $f:X \ra Y$ is a continuous
%%----------------------------------------------------------------------------------------------32
function $-$then $f$ is a weak homotopy equivalence iff $\sin f: \sin X \ra \sin Y$ is a simplicial homotopy equivalence.

[If $f$ is a weak homotopy equivalence, then $\abs{\sin f}$ is a weak homotopy equivalence 
(cf. p. \pageref{13.50}) 
or still, a homotopy equivalence.  But this means that $\sin f$ is a simplicial homotopy equivalence, 
$
\begin{cases}
\ \sin X\\
\ \sin Y
\end{cases}
$
being fibrant.]\\
\vspace{0.25cm}

\begingroup%%----------------------------------->> 
\fontsize{9pt}{11pt}\selectfont
A simplicial set $X$ is said to be 
\un{finite} 
\index{finite (simplicial set)} if $\abs{X}$ is finite.
\vspi
[Note: \  A finite simplicial set is a simplicial object in the category of finite sets (but not conversely).]\\
\endgroup%%------------------------------------<<

\index{Simplicial Approximation Theorem}
\begingroup%%----------------------------------->> 
\fontsize{9pt}{11pt}\selectfont
\textbf{\small SIMPLICIAL APPROXIMATON THEOREM} \ 
Let 
$
\begin{cases}
\ X\\
\ Y
\end{cases}
$
be simplicial sets with $X$ finite.  
Fix $\phi \in C(\abs{X},\abs{Y})$ $-$then $\exists$ $n > 0$ and a simplicial map
$f:\Sd^nX \ra Y$ such that $\abs{f} \simeq \phi \circx \abs{\td_X^n}$.
\vspi
[Since $\Exx^\infty Y$ is fibrant 
(cf. p. \pageref{13.51}), 
it follows from the simplicial extension theorem that there exists a simplicial map $F:X \ra \Exx^\infty Y$ such that $\abs{F} \simeq \abs{e_Y^\infty} \circx \phi$.  
But $X$ is finite, so $F$ factors through $\Exx^nY$ for some $n$.]
\vspi
[Note: \ The natural transformations $\td^n:\Sd^n \ra \id$ are defined inductively by 
$\td_X^0 = \id_X$, $\td_X^{n+1} = \td_X^n \circx \td_{\Sd^nX}$.]\\
\endgroup%%------------------------------------<< 

\begin{proposition} \ %21
Let $p:X \ra B$ be a simplicial map $-$then $p$ is a Kan fibration and a weak homotopy equivalence iff $p$ has the RLP w.r.t the inclusions $\ddn \ra \dn$ $(n \geq 0)$.
\end{proposition}

[That the condition is sufficient has been noted on 
p. \pageref{13.52}.  
As for the necessity, one can assume that $p$ is minimal (cf. Proposition 16).  To construct a filler $\dn \ra X$ for 
\begin{tikzcd}%[sep=large]
{\ddn} \ar{d} \ar{r} &{X} \ar{d}{p}\\
{\dn} \ar{r}[swap]{\Delta_b} &{B}
\end{tikzcd}
$(b \in B_n)$, it suffices to construct a filler $\dn \ra X_b$ for 
$
\begin{tikzcd}%[sep=large]
{\ddn} \ar{d} \ar{r} &{X_b} \ar{d}\\
{\dn} \arrow[r,shift right=0.5,dash] \arrow[r,shift right=-0.5,dash]  &{\dn}
\end{tikzcd}
.  
$
But the projection $X_b \ra \dn$ is a weak homotopy equivalence 
(cf. p. \pageref{13.53}) 
and $X_b$ is trivial over $\dn$ (cf. Proposition 14), say $X_b \approx \dn \times T_b$, where $T_b$ is fibrant.  
Therefore $\abs{T_b}$ is contractible, hence $T_b$ is simplicially contractible.  Now quote Proposition 18.]\\

Recall that \bCGH in its singular structure is a proper model category 
(cf. p. \pageref{13.54}).\\

\label{13.32}
\index{Fundamental Theorem of Simplicial Homotopy Theory}
\textbf{\small FUNDAMENTAL THEOREM OF SIMPLICIAL HOMOTOPY THEORY} \ \ \ 
\bSISET is a proper model category if weak equivalence = weak homotopy equivalence, cofibration = injective simplicial map, fibration = Kan fibration.  Every object is cofibrant and the fibrant objects are the fibrant simplicial sets.

%%----------------------------------------------------------------------------------------------33
[Axioms MC-1, MC-2, and MC-3 are immediate.

Claim: Every simplicial map $f:X \ra Y$ can be written as a composite $f_w \circx i_w$, where $i_w:X \ra X_w$ is an anodyne extension and $f_w:X_w \ra Y$ is a Kan fibration.

[In the small object argument, take $S_0 = \{\Lambda[k,n] \ra \Delta[n]$ $(0 \leq k \leq n, n \geq 1)\}$.]

Claim: Every simplicial map $f:X \ra Y$ can be written as a composite $f_w \circx i_w$, where $i_w:X \ra X_w$ is a cofibration and $f_w:X_w \ra Y$ is both a weak homotopy equivalence and a Kan fibration.

[In the small object argument, take $S_0 = \{\dot{\Delta}[n] \ra \Delta[n] \ (n \geq 0)\}$.]

Combining the claims gives MC-5.  Turning to MC-4, consider a commutative diagram
\begin{tikzcd}%[sep=large]
{A}  \ar{d}[swap]{i}  \ar{r}{u}  &{X} \ar{d}{p}\\
{Y}   \ar{r}[swap]{v} &{B}
\end{tikzcd}
, where $i$ is a cofibration and $p$ is a Kan fibration.  
If $p$ is a weak homotopy equivalence, then the existence of a filler 
$w:Y \ra X$ is implied by Proposition 7 and Proposition 21.  
On the other hand, if $i$ is a weak homotopy equivalence, 
then by the first claim $i = q \circx j$, where $j:A \ra Z$ is anodyne and 
$q:Z \ra Y$ is a Kan fibration which is necessarily a weak homotopy equivalence, 
so $\exists$ $f:Y \ra Z$ such that $f \circx i = j$, $q \circx f = \id_Y$.  
Consequently, $i$ is a retract of $j$, thus is itself anodyne.

There remains the verification of PMC.  Since all objects are cofibrant, half of this is automatic (cf. $\S 12$, Proposition 5).  Employing the usual notation, consider a pullback square
\begin{tikzcd}%[sep=large]
{P}  \ar{d}[swap]{\xi}  \ar{r}{\eta}  &{Y} \ar{d}{g}\\
{X}   \ar{r}[swap]{f} &{Z}
\end{tikzcd}
in \bSISET.  
Assume: $g$ is a Kan fibration and $f$ is a weak homotopy equivalence $-$then $\eta$ is a weak homotopy equivalence .  
Proof: 
\begin{tikzcd}%[sep=large]
{\abs{P}}  \ar{d}[swap]{\abs{\xi}}  \ar{r}{\abs{\eta}}  &{\abs{Y}} \ar{d}{\abs{g}}\\
{\abs{X}}   \ar{r}[swap]{\abs{f}} &{\abs{Z}}
\end{tikzcd}
is a pull back square in \bCGH (cf. Proposition 1), $\abs{g}$ is a Serre fibration (cf. Proposition 17), and $\abs{f}$ is a weak homotopy equivalence.  Therefore $\abs{\eta}$ is a weak homotopy equivalence.]

[Note: \  It is a corollary that $\bSISET_*$ (= $\Delta[0]\backslash$\bSISET) is a proper model category.]\\

\label{13.47} %dmc ?

\index{simplicial groups (example)}
\begingroup%%----------------------------------->>
\fontsize{9pt}{11pt}\selectfont
\textbf{\small EXAMPLE \  (\un{Simplicial Groups})} \ 
The free group functor $F_\text{gr}:\textbf{SET}\ra \bGR$ extends to a functor 
$F_\text{gr}:\bSISET \ra \textbf{SIGR}$ which is a left adjoint to the forgetful functor
$U: \textbf{SIGR} \ra \bSISET$.  Call a homomorphism $f:G \ra K$ of simplicial groups a weak equivalence if $U f$ is a weak homotopy equivalence, a fibration if $U f$ is a Kan fibration, and a cofibration if $f$ has the LLP w.r.t. acyclic fibrations $-$then with these choices, \textbf{SIGR} is a model category.  Here the point is that $f:G \ra K$ is a fibration (acyclic
%%----------------------------------------------------------------------------------------------34
 fibration) iff it has the RLP w.r.t the arrows $F_\text{gr} \Lambda[k,n] \ra F_\text{gr} \Delta[n]$ $(0 \leq k \leq n, n \geq 1)$ $(F_\text{gr} \dot{\Delta}[n] \ra F_\text{gr} \Delta[n]$ $(n \geq 0))$.  
Since $F_\text{gr}$ preserves cofibrations and $U$ preserves fibrations, the TDF theorem implies that 
$\textbf{L} F_\text{gr}:\bHSISET \ra \textbf{HSIGR}$ and 
$\textbf{R} U:\textbf{HSIGR} \ra \bHSISET$ 
exist and constitute an adjoint pair.
\vspi
[Note: Every object in \textbf{SIGR} is fibrant 
(cf. p. \pageref{13.55}) 
but not every object in \textbf{SIGR} is cofibrant.  Definition: A simplicial group $G$ is said to be 
\un{free}
\index{free (simplicial group)} 
if $\forall \ n$, $G_n$ is a free group with a specified basis $B_n$ such that $s_iB_n \subset B_{n+1}$ $(0 \leq i \leq n)$.  Every free simplicial group is cofibrant and every cofibrant simplicial group is the retract of a free simplicial group.]\\
\endgroup%%------------------------------------<<

\index{groupoids (example)}
\begingroup%%----------------------------------->>
\fontsize{9pt}{11pt}\selectfont
\textbf{\small EXAMPLE \ (\un{Groupoids})} \ 
\textbf{GRD} acquires the structure of a model category when one stipulates that a functor $F$ is a weak equivalence if 
$F$ is an equivalence of categories, a cofibration if $F$ is injective on objects, and a fibration if $\ner F$ is a Kan fibration.  
All objects are cofibrant and fibrant.\\
\endgroup%%------------------------------------<<

\label{13.71}
\index{G-Sets (example)}
\begingroup%%----------------------------------->>
\fontsize{9pt}{11pt}\selectfont
\textbf{\small EXAMPLE \  (\un{$G-$Sets})} \ 
Fix a group $G$.  
Denote by \bG the groupoid having a single object $*$ with Mor$(*,*) = G$ $-$then the category 
$\bSET_G$ of right $G$-sets is the functor category 
$[\bG^\OP,\bSET]$ and the category of simplicial right $G$-sets 
$\bSISET_G$ is the functor category 
$[\bDelta^\OP,[\bG^\OP,\textbf{SET}]]$ $\approx$ $[(\bDelta \times \bG)^\OP,\bSET]$.
Claim: $\bSISET_G$ is a model category.  Thus let
$U:\bSISET_G \ra \bSISET$ be the forgetful functor and declare that a morphism 
$f:X \ra Y$ of simplicial right $G$-sets is 
a weak equivalence if $U f$ is a weak homotopy equivalence, 
a fibration if $U f$ is a Kan fibration, 
and a cofibration if $f$ has the LLP w.r.t. acyclic fibrations.
\\
\indent\indent (CO) \  An object $X$ in $\bSISET_G$ is cofibrant iff $\forall \ n$, $X_n $ is a free $G$-set.
\vspi
Fix a cofibrant $XG$ in $\bSISET_G$ such that $XG \ra *$ is an acyclic fibration.  
Put $BG = XG/G$ $-$then $XG$ is simplicially contractible and locally trivial with fiber $G$ (i.e., si$G$), 
the projection $XG \ra BG$ is a Kan fibration, $BG$ is fibrant, and $\abs{BG}$ is a $K(G,1)$.  
Explicit models for $(XG,BG)$ can be found, e.g., in the notation of 
p. \pageref{13.56}, 
$XG = \text{bar}(*; \textbf{G};G)$ ($\approx$ $\ner \tran G$), $BG  = \text{bar}(*; \textbf{G};*)$ ($\approx$ $\ner \bG$).
\vspi
[Note: \ $U$ has a left adjoint $F_G$ which sends $X$ to $X \times$ $\si G$.  And, thanks to the TDF theorem, 
$(\textbf{L}F_G,\textbf{R}U)$ is an adjoint pair.]\\
\endgroup%%------------------------------------<<

Remark: The class of anodyne extensions is precisely the class of acyclic cofibrations.\\

\label{13.30}

\begingroup%%----------------------------------->>
\fontsize{9pt}{11pt}\selectfont
Claim: Sd preserves anodyne extensions.  
For suppose that $f:X \ra Y$ is anodyne and form the commutative diagram \ 
$
\begin{tikzcd}[sep=large]
{\text{Sd }X}  \ar{d}[swap]{\td_X}  \ar{r}{\text{Sd }f}  &{\text{Sd }Y} \ar{d}{\td_Y}\\
{X}   \ar{r}[swap]{f} &{Y}
\end{tikzcd}
.  \ 
$
Since Sd preserves injections, Sd $f$ is a cofibration.  But $\td_X$ $\&$ $\td_Y$ are weak homotopy equivalences 
(cf. Proposition 5), thus Sd $f$ is an acyclic cofibration, i.e., is anodyne.\\
\endgroup%%------------------------------------<<

\begin{proposition} \ %22
Suppose that $L \ra K$ is an inclusion of simplicial sets and $X \ra B$ is a Kan fibration $-$then the arrow
$\text{map}(K,X) \ra \text{map}(L,X) \times_{\text{map}(L,B)} \text{map}(K,B)$ is a
%%----------------------------------------------------------------------------------------------35
Kan fibration (cf. Proposition 12) which is a weak homotopy equivalence if this is the case of $L \ra K$ or $X \ra B$.
\end{proposition}

[Owing to Proposition 21, the problem is to produce a filler $\Delta[n] \times K \ra X$ for 
\begin{tikzcd}%[sep=large]
{\dot\Delta [n] \times K \cup \Delta [n] \times L}  \ar{d}[swap]{i}  \ar{r}  &{X} \ar{d}\\
{\Delta [n] \times K}   \ar{r} &{B}
\end{tikzcd}
.  
If $L \ra K$ is an acyclic cofibration, then, as pointed out above, $L \ra K$ is anodyne.  
Therefore $i$ is anodyne (cf. Proposition 9) and the filler exists.  
If $X \ra B$ is an acyclic Kan fibration, then the existence of the filler is guaranteed by MC-4.]\\

\label{13.24}
\label{13.25}
\bHSISET 
\index{\bHSISET} 
is the homotopy category of \bSISET 
(cf. \pageref{13.57} ff.).  
In this situation, $IX = X \times \Delta[1]$ serves as a cylinder object while $PX = \text{map}(\Delta[1],X)$ is a path object when $X$ is fibrant but not in general 
(Berger\footnote[2]{\textit{Bull. Soc. Math. France} \textbf{123} (1995), 1-32.}).  
Since all objects are cofibrant, $\mathcal{L}X = X$ $\forall \ X$ and there are canonical choices for $\mathcal{R}X$, e.g., $\sin \abs{X}$ or $\text{Ex}^\infty X$.  If $X$ is cofibrant and $Y$ is fibrant, then left homotopy = right homotopy or still, simplicial homotopy: $[X,Y] \approx [X,Y]_0$.  \bHSISET has finite products.  
And: \bHSISET is cartesian closed.  
Proof: 
$[X \times Y, Z]$ $\approx$ 
$[X \times Y, \sin\abs{Z}]$ $\approx$ 
$[X \times Y, \sin\abs{Z}]_0$ $\approx$ 
$[X,\text{map}(Y,\sin \abs{Z})]_0$ $\approx$ 
$[X,\text{map}(Y,\sin \abs{Z})]$.

[Note: \  Recall too that the inclusion $\bHSISET_{\bff} \ra \bHSISET$ is an equivalence of categories 
(cf. $\S 12$, Proposition 13).]

Example: $X$ and $X^\OP$ are naturally isomorphic in \bHSISET.\\

\begingroup%%----------------------------------->>
\fontsize{9pt}{11pt}\selectfont
\textbf{\small FACT} \ 
Let $S \subset \Mor \bH_0\bSISET$ be the class of homotopy classes of anodyne extensions  
$-$then $S^{-1}\bH_0\bSISET$ is equivalent to \bHSISET.\\
\endgroup%%------------------------------------<<

\index{Theorem: Comparison Theorem (simplicial sets)}
\textbf{\small COMPARISON THEOREM} \ 
The adjoint pair $(\abs{?},\sin)$ induces an adjoint equivalence of categories between \bHSISET and \bHTOP (singular structure).

[In the TDF theorem, take $F = \abs{?}$, $G = \sin$ $-$then $F$ preserves cofibrations and $G$ preserves fibrations, thus
$
\begin{cases}
\ \textbf{L}F\\
\ \textbf{R}G
\end{cases}
$
exist and $(\bL F,\bR G)$ is an adjoint pair.  
Consider now the bijection of adjunction 
$\Xi_{X,Y}:C(\abs{X},Y) \ra \text{Nat}(X,\sin Y)$ so $\Xi_{X,Y}f$ is the composition 
$X \ra \sin \abs{X}$ 
%\overset{\sin f}{\longrightarrow} 
$\overset{\sin f}{\xrightarrow{\hspace*{1cm}}}$ 
$\sin Y$.  
Since the arrow $X \ra \sin \abs{X}$ is a weak homotopy equivalence 
(cf. p. \pageref{13.58}), 
$\Xi_{X,Y}f$ is a weak homotopy equivalence 
iff $\sin f$ is a weak homotopy equivalence, 
i.e., 
iff $f$ is a weak homotopy equivalence 
(cf. p. \pageref{13.59}).  
Therefore the pair $(\bL F,\bR G)$ is an adjoint equivalence of categories 
(cf. p. \pageref{13.59a}).]\\

Application: \bHSISET is equivalent to \bHCW.

%%----------------------------------------------------------------------------------------------36
[Note: \   Analogously, $\bHSISET_*$ is equivalent to $\bHCW_*$.]\\

Are there model other categories \bC whose associated homotopy category \textbf{HC} is equivalent to \bHCW?  
The answer is ``yes''.\\

\label{17.18}

\begingroup%%----------------------------------->>
\fontsize{9pt}{11pt}\selectfont
\textbf{\small EXAMPLE} \ 
Take \bC = \textbf{CAT} and call a morphism $f$ a weak equivalence if Ex$^2 \circx $ ner $f$ is a weak homotopy equivalence, a fibration if Ex$^2 \circx $ ner $f$ is a Kan fibration, and a cofibration if $f$ has the LLP w.r.t. all fibrations that are weak equivalences $-$then 
Thomason\footnote[2]{\textit{Cahiers Topologie G\'eom. Diff\'erentielle} \textbf{21} (1980), 305-324.} 
has shown that \textbf{CAT} is a proper model category.  Put 
$
\begin{cases}
\ F = c \circx \text{Sd}^2\\
\ G = \text{Ex}^2 \circx ner
\end{cases}
$
: \ $(F,G)$ is an adjoint pair with the property that $F$ preserves cofibrations and $G$ preserves fibrations, thus 
$
\begin{cases}
\ \textbf{L}F\\
\ \textbf{R}G
\end{cases}
$
exist and $(\bL F,\bR G)$ is an adjoint pair (TDF theorem).  
Moreover, the arrow
$X \ra \text{Ex}^2 \circx \ner \circx c \circx \text{Sd}^2 X$ is a weak homotopy equivalence of simplicial sets, so the pair $(\bL F,\bR G)$ is an adjoint equivalence of cateogories.  
It therefore follows that \bHSISET, \bHCAT, and \bHCW are equivalent.
\vspi
[Note: 
Latch\footnote[3]{\textit{J. Pure Appl. Algebra} \textbf{9} (1977), 221-237.} 
proved that $\ner: \bCAT \ra \bSISET$ induces an equivalence $\bHCAT \ra \bHSISET$ 
(but the adjoint pair $(c,\ner)$ does not induce an adjoint equivalence).]\\
\endgroup%%------------------------------------<<

\begingroup%%----------------------------------->>
\fontsize{9pt}{11pt}\selectfont
\textbf{\small EXAMPLE} \ 
The category of simplicial groupoids is a model category and its homotopy category is equivalent to \bHSISET, hence to \bHCW
(Dwyer-Kan\footnote[6]{\textit{Nederl. Akad. Wetensch. Indag. Math.} \textbf{46} (1984), 379-385; 
see also Heller, \textit{Illinois J. Math.} \textbf{24} (1980), 576-605.}).
\vspi
\label{13.62}
[Note: \  A 
\un{simplicial groupoid}
\index{simplicial groupoid} 
\bG is a category object $(M,O)$ in \bSISET, where \mO is a constant simplicial set, equipped with a simplicial map $\chi:M \ra M$ such that
$s = t \circx \chi$, $t = s \circx \chi$, $c \circx (\chi \times \id_M) = e \circx s$, $c \circx(\id_M \times \chi) = e \circx t$.
So, $\forall \ n$, $\bG_n$ is a groupoid and $\Ob\bG_n = \Ob\bG_0$.  
Introducing the obvious notion of morphism, the simplicial groupoids are seen to constitute a category which is complete and cocomplete.  
Its model category structure is derived from (1)-(3) below.
\\
\indent\indent (1) \  A morphism $F:\bG \ra \bK$ of simplicial groupoids is a weak equivalence if $F$ restricts to a bijection on components and $\forall \ X \in O$, the induced morphism $\bG(X) \ra \bK(FX)$ of simplicial groups is a weak equivalence.
\\
\indent\indent (2) \  A morphism $F:\bG \ra \bK$ of simplicial groupoids is a fibration if $F_0:\bG_0 \ra \bK_0$ is a fibration of groupoids and $\forall \ X \in O$, the induced morphism $\bG(X) \ra \bK(FX)$ of simplicial groups is a fibration.
\\
\indent\indent (3) \ A morphism $F:\bG \ra \bK$ of simplicial groupoids is a cofibration if it has the LLP w.r.t acyclic fibrations.]\\
\endgroup%%------------------------------------<<

%%----------------------------------------------------------------------------------------------37
\label{13.87}
Fix a small category \bI $-$then the functor category [\bI,\bSISET] admits two proper model category structures.  However, the weak equivalences in either structure are the same, so both give rise to the same homotopy category \textbf{H[I,SISET]}.

\label{13.70}
\label{13.92}
\label{13.107a}
\indent\indent (L) \  Given functors $F,G:\bI \ra \bSISET$, call $\Xi \in \text{Nat}(F,G)$ a weak equivalence if $\forall \ i$, $\Xi_i:F_i \ra G_i$ is a weak homotopy equivalence, a fibration if $\forall \ i$, $\Xi_i:F i \ra G i$ is a Kan fibration, a cofibration if $\Xi$ has the LLP w.r.t acyclic fibrations.

\indent\indent (R) \  Given functors $F,G:\bI \ra \bSISET$, call $\Xi \in \Nat(F,G)$ a weak equivalence if $\forall \ i$, $\Xi_i:F_i \ra G_i$ is a weak homotopy equivalence, a cofibration if $\forall \ i$, $\Xi_i:F i \ra G i$ is an injective simplicial map, a fibration if $\Xi$ has the RLP w.r.t acyclic cofibrations.

In practice, both structures are used but for theoretical work, structure L is generally the preferred choice.

[Note: \  When \bI is discrete, structure L = structure R (all data is levelwise).]\\

\begingroup%%----------------------------------->>
\fontsize{9pt}{11pt}\selectfont
Since the arguments are dual, it will be enough to outline the proof in the case of structure L.
\vspi
Notation: Let $f:X \ra Y$ be a simplicial map $-$then $f$ admits a functorial factorization
$X \overset{i_f}{\lra} \mathcal{L}_f \overset{\pi_f}{\lra} Y$,
where $i_f$ is a cofibration and $\pi_f$ is an acyclic Kan fibration, and a functorial factorization
$X \overset{\iota_f}{\lra} \mathcal{R}_f \overset{p_f}{\lra} Y$,
where $\iota_f$ is an acyclic cofibration and $p_f$ is a Kan fibration.
\vspi
Observation: These factorizations extend levelwise to factorizations of $\Xi:F \ra G$, viz.
$F \overset{i_\Xi}{\lra} \mathcal{L}_\Xi \overset{\pi_\Xi}{\lra} G$
and 
$F \overset{\iota_\Xi}{\lra} \mathcal{R}_\Xi \overset{p_\Xi}{\lra} G$.

Write $\bI_\text{dis}$ for the discrete category underlying \bI $-$then the forgetful functor
%$U:\textbf{[I,SISET]} \ra \textbf{[I$_\text{dis}$,SISET]}$ 
$U:[\bI,\bSISET] \ra [\bI_\text{dis},\bSISET]$
has a left adjoint that sends $X$ to $\fr X$, where
$\fr X j = \ds\coprod\limits_{i \in \Ob\bI} \text{Mor}(i,j) \cdot X i$.\\
\endgroup%%------------------------------------<<

\begingroup%%----------------------------------->>
\fontsize{9pt}{11pt}\selectfont
\textbf{\small LEMMA} \ 
Fix an $F$ in [\bI,\bSISET].  \ 
Suppose that $\Phi:UF \ra X$ is a cofibration in 
$[\bI_\text{dis},\bSISET]$ and 
\begin{tikzcd}%[sep=large]
{\fr UF}  \ar{d}[swap]{\nu_F}  \ar{r}{\fr\Phi}  &{\fr X} \ar{d}{u}\\
{F}   \ar{r} &{G}
\end{tikzcd}
is a pushout square in [\bI,\bSISET] $-$then the composite 
$Uu \circx \mu_X:X \overset{\mu_X}{\lra} U\text{fr}X \overset{Uu}{\lra} UG$ is a cofibration in 
$[\bI_\text{dis},\bSISET]$.
\vspi
[The commutative diagram %note - still need to shift right the top line of the next diagram
\[
\begin{tikzcd}[sep=huge]
&&{X j} \ar{d} \arrow[r,shift right=0.5,dash] \arrow[r,shift right=-0.5,dash]  &{X j} \ar{d}{(\mu_X)_j}\\
&{\bigl(\ds\coprod\limits_{\substack{i \overset{\delta}{\ra} j\\ \delta \neq \id_j}} F i \bigr) \amalg F j}                          \ar{d} \ar{r}
&{\bigl(\ds\coprod\limits_{\substack{i \overset{\delta}{\ra} j\\ \delta \neq \id_j}} F i \bigr) \amalg {X j}}                        
\ar{d} \ar{r}
&{\bigl(\ds\coprod\limits_{\substack{i \overset{\delta}{\ra} j\\ \delta \neq \id_j}} X i \bigr) \amalg X j}
\ar{d}{u_j}\\
&{F j} \ar{r}[swap]{\Phi_j} &{X j} \ar{r}[swap]{\mu_j \circx (\mu_X)_j} &{G j}
\end{tikzcd}
\]
%%----------------------------------------------------------------------------------------------38
tells the tale.  Indeed, the middle row is a factorization of $(\fr\Phi)_j$ (suppression of ``U''), the bottom square on the right is a pushout, and a coproduct of cofibrations is a cofibration.]
\vspi
[Note: \  As usual,
$
\begin{cases}
\ \mu\\
\ \nu
\end{cases}
$
are the ambient arrows of adjunction.]\\
\endgroup%%------------------------------------<<

\label{13.107b}
\begingroup%%----------------------------------->>
\fontsize{9pt}{11pt}\selectfont
Consider any $\Xi:F \ra G$.  Claim: $\Xi$ can be written as the composite of a cofibration and an acyclic fibration.  Thus define $F_1$ by the pushout square
\begin{tikzcd}[sep=large]
{\text{fr}UF}  \ar{d}[swap]{\nu F}  \ar{r}{\text{fr}Ui_\Xi}  &{\fr U\mathcal{L}_\Xi} \ar{d}\\
{F}   \ar{r} &{F_1}
\end{tikzcd}
$-$then there is a commutative diagram
\[
\begin{tikzcd}[sep=huge]
&{\fr U F}  \ar{d}[swap]{\nu_F}  \ar{r}{\fr U i_\Xi}  
&{\fr U \mathcal{L}_\Xi} \ar{d} \ar{r}{\fr U \pi_\Xi} 
&{\fr U G} \ar{d}{\nu_G}\\
&{F} \ar{rd}   \ar{r} &{F_1}  \ar{d} \ar{r} &{G}\\
& &{\mathcal{L}_\Xi} \ar{ru}
\end{tikzcd}
\]
in which 
$\fr U\mathcal{L}_\Xi \ra F_1 \ra \mathcal{L}_\Xi$ is $\nu_{\mathcal{L}_\Xi}$.
Putting $F_0 = F$ (and $\Xi_0 = \Xi$), iterate the construction to obtain a sequence 
$F = F_0 \ra F_1 \ra \cdots \ra F_\omega$ of objects in [\bI,\bSISET], taking 
$F_\omega = \colimx F_n$.  This leads to a commutative triangle
\begin{tikzcd}[sep=small]
{F}  \ar{ddr}[swap]{\Xi}  \ar{rr}{\i_\omega}  &&{F_\omega} \ar{ddl}{\Xi_\omega}\\
\\
&{G}
\end{tikzcd}
.  
Here, $i_\omega$ is a cofibration (since the $F_n \ra F_{n+1}$ are).  
Moreover, $i_\omega$ is a weak equivalence whenver $\Xi$ is a weak equivalence and in that situation, 
$i_\omega$ has the LLP w.r.t. all fibrations.  
To see that $\Xi_\omega$ is an acyclic fibration, look at the interpolation\\
\vspace{0.2cm}
\[
\begin{tikzcd}[sep=large]
{UF_0} \ar{d} \ar{r} &{U\mathcal{L}_{\Xi_0}} \ar{d} \ar{r} &{UF_1} \ar{r} \ar{d} \ar{r}&{U\mathcal{L}_{\Xi_1}}\ar{d} \ar{r}  &{\cdots} \\
{UG}  \arrow[r,shift right=0.5,dash] \arrow[r,shift right=-0.5,dash]  
&{UG} \arrow[r,shift right=0.5,dash] \arrow[r,shift right=-0.5,dash]  
&{UG} \arrow[r,shift right=0.5,dash] \arrow[r,shift right=-0.5,dash]  
&{UG} \arrow[r,shift right=0.5,dash] \arrow[r,shift right=-0.5,dash] 
&{\cdots}
\end{tikzcd}
\]
in $[\bI_\dis,\bSISET]$.  
Thanks to the lemma, the horizontal arrows in the top row are cofibrations.  
On the other hand, the arrows $U\mathcal{L}_{\Xi_n} \ra UG$ are acyclic fibrations.  
But then $U\Xi_\omega$ is an acyclic fibration per $[\bI_\dis,\bSISET]$, i.e., $\Xi_\omega$ is an acyclic fibration per [\bI,\bSISET].  
Hence the claim.
\vspi
To finish the verification of MC-5, one has to establish that $\Xi$ can be written as the composite of an acyclic cofibration and a fibration.  
This, however, is immediate: 
Apply the claim to $\iota_\Xi$.  MC-4 is equally clear.  
For if $\Xi$ is a cofibration, then $\Xi$ is a retract of $i_\omega$, so if $\Xi$ is an acyclic cofibration, then $\Xi$ has the LLP w.r.t all fibrations.  PMC is obvious.\\
\endgroup%%------------------------------------<<
\label{13.120}

\begingroup%%----------------------------------->>
\fontsize{9pt}{11pt}\selectfont
\textbf{\small EXAMPLE} \ 
Definition: A functor $F:\bI \ra \bSISET$ is said to be 
\un{free} 
\index{free (functor)} 
if $\exists$ functors 
$B_n:\bI_\dis \ra \bSET$ $(n \geq 0)$ such that 
$\forall \ j \in \Ob\bI$ : $B_nj \subset (F j)_n$ \ $\&$ \ $s_iB_nj \subset B_{n+1}j$ $(0 \leq i \leq n)$, 
with \ $\fr B_n \approx F_n$ \  ($F_nj = (F j)_n$).  
Every free functor is cofibrant in structure L and every cofibrant functor in structure L is the retract of a free functor.  
Example: $\ner(\bI/-)$ is a free functor, hence is cofibrant in structure L.\\
\endgroup%%------------------------------------<<

%%----------------------------------------------------------------------------------------------39
Fix an abelian group \mG.  Let $f:X \ra Y$ be a simplicial map $-$then $f$ is said to be an 
\un{$HG$-equivalence} 
\index{equivalence! $HG$-equivalence}
if $\forall \ n \geq 0$, $\abs{f}_*:H_n(\abs{X};G) \ra H_n(\abs{Y};G)$ is an isomorphism.  
Agreeing that an 
\un{$HG$-cofibration} 
\index{cofibration! $HG$-cofibration}
is an injective simplicial map, an 
\un{$HG$-fibration}
\index{fibration! $HG$-fibration}
is a simplicial map which has the RLP w.r.t all $HG$-cofibrations that are $HG$-equivalences.  
Every $HG$-fibration is a Kan fibration.  
Proof: $\Lambda^{k,n}$ is a strong deformation retract of $\Delta^n$.\\

\begin{proposition} \ %23
Let $p:X \ra B$ be a simplicial map $-$then $p$ is an $HG$-fibration and an $HG$-equivalence iff $p$ is a Kan fibration and a weak homotopy equivalence.
\end{proposition}

[Necessity: Write $p = q \circx j$, where $j:X \ra Y$ is a cofibration and $q:Y \ra B$ is an acyclic Kan fibration.  
Since $p$ is an $HG$-equivalence, the same is true of $j$, thus the commutative diagram 
\begin{tikzcd}%[sep=large]
{X}  \ar{d}[swap]{j} \arrow[r,shift right=0.5,dash] \arrow[r,shift right=-0.5,dash]  &{X} \ar{d}{p}\\
{Y} \ar{r}[swap]{q}  &{B}  
\end{tikzcd}
admits a filler $g:Y \ra X$.  
Therefore $p$ is a retract of $q$, hence is an acyclic Kan fibration.

Sufficiency: Apply Proposition 7 and Proposition 21.]\\

Notation: Given a simplicial set $X$, write $\#(X)$ for $\#(\sE)$, the cardinality of the set of cells in $\abs{X}$.

[Note: \  $\forall$ set $X$, $\#(\si X) = \#(X)$, the cardinality of X.]\\

\begin{proposition} \ %24
Let $p:X \ra B$ be a simplicial map which has the RLP w.r.t. every inclusion 
$A \ra Y$, where $H_*(\abs{Y},\abs{A};G) = 0$ and $\#(Y)$ is $\leq \#(G)$ if $\#(G)$ is infinite and $\leq \omega$ if $\#(G)$ is finite $-$then $p$ is an $HG$-fibration.
\end{proposition}

[It suffices to prove that $p$ has the RLP w.r.t every inclusion $L \ra K$ $(L \neq K)$ with 
$H_*(\abs{K},\abs{L};G) = 0$.  This can be established by using Zorn's lemma.  Indeed, $\exists$ a simplicial subset 
$A \subset K$ $(A \not\subset L)$ such that $H_*(\abs{A},\abs{A \cap L};G) = 0$ subject to the restriction that 
$\#(A)$ is $\leq \#(G)$ if $\#(G)$ is infinite and $\leq \omega$ if $\#(G)$ is finite 
(cf. p. \pageref{13.60}).]\\

\index{Prefactorization Lemma}
\textbf{\small PREFACTORIZATION LEMMA} \ 
Suppose that $\kappa$ is an infinite cardinal.  Let $f:X \ra Y$ be a simplicial map $-$then $f$ can be written as a composite $f = p_f \circx i_f$, where $i_f:X \ra X_f$ is an injection with $H_*(\abs{X_f},\abs{X};G) = 0$, such that every commutative diagram
\begin{tikzcd}%[sep=large]
{L} \ar{d} \ar{r} &{X} \ar{r}{i_f} &{X_f} \ar{d}{p_f}\\
{K} \ar{rr} &&{Y} 
\end{tikzcd}
has a filler $K \ra X_f$, $(K,L)$ being any simplicial pair with $\#(K) \leq \kappa$ and $H_*(\abs{K},\abs{L};G) = 0$.

[Choose a set of simplicial pairs $(K_i,L_i)$ with $\#(K_i) \leq \kappa$ and 
$H_*(\abs{K_i},\abs{L_i};G) = 0$ which contains up to isomorphism all such simplicial pairs.  Consider the set of pairs 
%%----------------------------------------------------------------------------------------------40
of morphisms $(g,h)$ such that the diagram
\begin{tikzcd}%[sep=large]
{L_i}  \ar{d} \ar{r}{g} &{X} \ar{d}{f}\\
{K_i} \ar{r}[swap]{h}  &{Y}  
\end{tikzcd}
commutes, define $X_f$ by the pushout square 
$
\begin{tikzcd}%[sep=large]
{\coprod\limits_i\coprod\limits_{(g,h)}L_i}  \ar{d} \ar{r} &{X} \ar{d}{i_f}\\
{\coprod\limits_i\coprod\limits_{(g,h)}K_i} \ar{r}  &{X_f}  
\end{tikzcd}
, 
$
and let $p_f:X_f \ra Y$ be the induced simplicial map.]\\
\vspace{0.25cm}

\index{Homological Model Category Theorem}
\textbf{\small HOMOLOGICAL MODEL CATEGORY THEOREM} \ 
Fix an abelian group $G$ $-$then \bSISET is a model category if 
weak equivalence = $HG$-equivalence, 
cofibration = $HG$-cofibration, 
fibration = $HG$-fibration.

[On the basis of Proposition 23, one has only to show that every simplicial map $f:X \ra Y$ can be written as a composite 
$p \circx i$, where $i$ is an acyclic $HG$-cofibration and $p$ is an $HG$-fibration.  This can be done by a transfinite lifting argument, using the prefactorization lemma with $\kappa$ a regular cardinal $> \#(G)$ (cf. Proposition 24).]

[Note: \  The fibrant objects in this structure are the 
\un{$HG$-local}
\index{local!$HG$-local (objects)} 
objects, i.e., those $X$ such that $X \ra *$ is an $HG$-fibration.]\\

\begin{proposition} \ %25
Suppose that $L \ra K$ is an inclusion of simplicial sets and $X \ra B$ is an $HG$-fibration $-$then the arrow 
$\text{map}(K,X) \ra \text{map}(L,X)  \times_{\text{map}(L,B) } \text{map}(K,B)$
is an $HG$-fibration which is an $HG$-equivalence if this is the case of $L \ra K$ or $X \ra B$.\\
\end{proposition}

\label{12.21}
\begingroup%%----------------------------------->>
\fontsize{9pt}{11pt}\selectfont
\textbf{\small EXAMPLE} \ 
The \mc structure on \bSISET provided by the homological model category theorem is generally not proper.  
Thus factor 
$X \ra *$ as $X \ra X_{HG} \ra *$, where $X \ra X_{HG}$ is an acyclic $HG$-cofibration and 
$X_{HG} \ra *$ is an $HG$-fibration.  
Assuming that $X$ is fibrant and connected, define $E_{HG}$ by the pullback square
\begin{tikzcd}[sep=large]
{E_{HG}}  \ar{d} \ar{r} &{\Theta X_{HG}} \ar{d}\\
{X} \ar{r}  &{X_{HG}}  
\end{tikzcd}
$-$then the arrow $E_{HG} \ra \Theta X_{HG}$ is not necessarily an  $HG$-equivalence.\\
\endgroup%%------------------------------------<<

\begingroup%%----------------------------------->>
\fontsize{9pt}{11pt}\selectfont
\textbf{\small FACT} \ 
Suppose given simplicial maps $f:X \ra Y$, $g:Y \ra Z$, where $f$ is a Kan fibration and $g$, $g \circx f$ are $HG$-fibrations $-$then $f$ is an $HG$-fibration.\\
\endgroup%%------------------------------------<<

\begingroup%%----------------------------------->>
\fontsize{9pt}{11pt}\selectfont
Application: If $f:X \ra Y$ is a Kan fibration and if 
$
\begin{cases}
\ X\\
\ Y
\end{cases}
$
are $HG$-local, then $f$ is an $HG$-fibration.\\
\endgroup%%------------------------------------<<

\begingroup%%----------------------------------->>
\fontsize{9pt}{11pt}\selectfont
\textbf{\small EXAMPLE} \ 
The $HG$-local objects in \bSISET are closed under the formation of products and $\map(X,Y)$ is $HG$-local $\forall \ X$ provided that $Y$ is $HG$-local.  Given a 2 sink 
$X \overset{f}{\ra} Z \overset{g}{\la} Y$ 
of $HG$-local objects
%%----------------------------------------------------------------------------------------------41
with $f$ a Kan fibration, the pullback $X \times_Z Y$ is  $HG$-local.  Finally, for any tower
$X_0 \la X_1 \la \cdots $ of Kan fibrations and  $HG$-local $X_n$, the limit $\lim X_n$ is $HG$-local.\\
\endgroup%%------------------------------------<<

A 
\un{simplicial category}
\index{simplicial category} 
is a \bSISET-category.  
So, to specify a simplicial category one must specify a class of objects \mO and a function that assigns to each ordered pair $X,Y \in O$ a simplicial set $\HOM(X,Y)$ plus simplicial maps $C_{X,Y,Z}:\text{HOM}(X,Y) \times \text{HOM}(Y,Z) \ra \text{HOM}(X,Z)$,
$I_X: \Delta[0] \ra \text{HOM}(X,X)$ satisfying \bSISET-cat$_1$ and \bSISET-cat$_2$ 
(cf. p. \pageref{13.61}).
Here is an equivalent description.  
Fix a class \mO.  
Consider the metacategory $\mathcal{CAT}_O$ whose objects are the categories with object class \mO, the morphisms being the functors which are the identity on objects $-$then a simplicial category with object class \mO is a simplicial object in $\mathcal{CAT}_O$.\\

\begingroup%%----------------------------------->>
\fontsize{9pt}{11pt}\selectfont
A category object $(M,O)$ in \bSISET, where \mO is a constant simplicial set, is a simplicial category.  
In particular: A simplicial groupoid is a simplicial category 
(cf. p. \pageref{13.62}).\\
\endgroup%%------------------------------------<<

\begingroup%%----------------------------------->>
\fontsize{9pt}{11pt}\selectfont
\textbf{\small EXAMPLE} \ 
There is a functor $\bDelta^\OP \ra \bSISET$ which sends [n] to $\Delta[1]^n$ and 
$
\begin{cases}
\ \delta_i \text{ to } d_i\\
\ \sigma_i \text{ to } s_i
\end{cases}
$
, where
\[
d_i(\alpha_1, \ldots, \alpha_n) = 
\begin{cases}
\ (\alpha_2, \ldots, \alpha_n) \qquad\qquad\qquad\qquad\  \ \ \ (i = 0)\\
\ (\alpha_1, \ldots, \max(\alpha_{i+1},\alpha_i), \ldots, \alpha_n) \ \ \  (0 < i < n),\\
\ (\alpha_1, \ldots, \alpha_{n-1}) \qquad\qquad\qquad\qquad\  (i = n)
\end{cases}
\]
$s_i(\alpha_1, \ldots, \alpha_n) = (\alpha_1, \ldots, \alpha_i, 0, \alpha_{i+1}, \ldots, \alpha_n)$.  \ 
Now fix a small category \bC.  \ Given $X, \ Y \in \Ob\bC$, let $C = C(X,Y)$ be the cosimplicial set defined by taking for 
$C(X,Y)^n$ the set of all functors 
$F:[n+1] \ra \bC$ with $F_0 = X$, $F_{n+1} = Y$ and letting 
$C\delta_i:C^n \ra C^{n+1}$, $C\sigma_i:C^n \ra C^{n-1}$ be the assignments 
$(f_0, \ldots, f_n) \ra (f_0, \ldots, f_{i-1}, \id, f_i \ldots, f_n)$, 
$(f_0, \ldots, f_n) \ra (f_0, \ldots, f_{i+1} \circx f_i, \ldots, f_n)$.  
Definition: 
$\HOM(X,Y) = \ds\int^{[n]} \Delta[1]^n \times C(X,Y)^n$.  \ 
Since 
$\HOM(X,Y)_m = \ds\int^{[n]} \Delta[1]_m^n \times C(X,Y)^n$, one can introduce a ``composition'' rule and a ``unit'' rule satisfying the axioms.  
The upshot, therefore, is a simplicial category \bFRC with $O = \Ob\bC$.
\vspi
[Note: \  The abstract interpretation of \bFRC is this.  
Observe first that the forgetful functor from \bCAT to the category of small graphs with distinguished loops at the vertexes has a left adjoint.  
Consider the associated cotriple in \bCAT $-$then the standard resolution of \bC is \bFRC and the underlying category \bUFRC is the free category on $\Ob\bC$ having one generator for each nonidentity morphism in \bC.]\\
\endgroup%%------------------------------------<<

\label{13.64}
Let \bC be a category.  \ 
Suppose that $X,\ Y$ are simplicial objects in \bC and let $K$ be a simplicial set $-$then a 
\un{formality}
\index{formality (simplicial objects)} \ 
$f:X \bbox K \ra Y$ is a collection of morphisms $f_n(k): X_n \ra Y_n$ in \bC, one for each $n \geq 0$ and $k \in K_n$, such that $Y\alpha \circx f_n(k) = f_m((K\alpha)k)\circx X\alpha$ $(\alpha:[m] \ra [n])$.  
Notation: For$(X\bbox K, Y)$.  
Example: For$(X \bbox \Delta[0],Y)$ can be identified with $\Nat(X,Y)$.

%%----------------------------------------------------------------------------------------------42
[Note: \  As it stands, $X \bbox K$ is just a symbol, not an object in \bSIC (but see below).]\\

\begin{proposition} \ %26 dmcxx
Let \bC be a category $-$then the class of simplicial object in \bC is the object class of a simplicial category \textbf{SIMC}.
\end{proposition}

[Define: $\HOM(X,Y)$ by letting $\HOM(X,Y)_n$ be For$(X \bbox \Delta[n],Y)$.]

[Note: \  \bSIC is isomorphic to the underlying category of \bSIMC.]\\

A 
\un{simplicial functor}
\index{simplicial functor} 
is a \bSISET-functor.  Example:  If 
$
\begin{cases}
\ \bC\\
\ \bD
\end{cases}
$
are categories and $F:\bC \ra \bD$ is a functor, then $F$ extends to a simplicial functor $SF: \bSIMC \ra \bSIMD$.\\

\begingroup%%----------------------------------->>
\fontsize{9pt}{11pt}\selectfont
\textbf{\small EXAMPLE} \ 
\bCAT is cartesian closed, hence can be viewed as a \bCAT-category.  Since 
$\ner:\bCAT \ra \bSISET$ is a morphism of symmetric monoidal categories, 
$\ner_*\bCAT$ 
is a simplicial category whose object class is the class of small categories, $\HOM(\bC,\bD)$ being $\ner[\bC,\bD]$ 
(cf. p. \pageref{13.63}).  
One may therefore interpret ner as a simplicial functor 
$\ner_*\bCAT \ra \bSISET$ (for $\ner[\bC,\bD] \approx \map(\ner \bC, \ner \bD)$).\\
\endgroup%%------------------------------------<<

Given a category \bC, a \un{simplicial action} on \bC is a functor 
$\Box:\bC \times \bSISET \ra \bC$, 
together with natural isomorphisms $R$ and $A$, where 
$R_X:X \bbox \Delta[0] \ra X$, 
$A_{X,K,L}:X \bbox (K \times L) \ra (X \bbox K) \bbox L$, 
subject to the following assumptions.

\indent\indent (SA$_1$) \ The diagram\\
\[
\begin{tikzcd}[sep=small] %dmcxx
{X \bbox(K \times (L \times M))}  \ar{dd}[swap]{\id \bbox A} \ar{r}{A}
&{(X \bbox K) \bbox (L \times M)} \ar{r}{A}
&{((X \bbox K) \bbox L) \bbox M}\\
\\
{X \bbox ((K \times L) \times M)}  \ar{rr}[swap]{A}
&&{(X \bbox (K \times L)) \bbox M} \ar{uu}[swap]{A \bbox \id}\\
\end{tikzcd}
\]
commutes.

\indent\indent (SA$_2$) \ The diagram\\

\begin{tikzcd}%[sep=large]
&&&{X \bbox(\Delta[0] \times K)}  \ar{d}[swap]{\id \bbox L} \ar{r}{A}
&{(X \bbox \Delta[0]) \bbox K} \ar{d}{R \bbox \id}\\
&&&{X \bbox K}  \arrow[r,shift right=0.5,dash] \arrow[r,shift right=-0.5,dash] 
&{X \bbox K}
\end{tikzcd}\\
commutes.

[Note: \  Every category admits a simplicial action, viz. the trivial simplicial action.]\\

\begingroup%%----------------------------------->>
\fontsize{9pt}{11pt}\selectfont
It is automatic that the diagram\\
\begin{tikzcd}[sep=large]
&&&{X \bbox(K \times \Delta[0])}  \ar{d}[swap]{\id \bbox R} \ar{r}{A}
&{(X \bbox K) \bbox \Delta[0]} \ar{d}{R}\\
&&&{X \bbox K}  \arrow[r,shift right=0.5,dash] \arrow[r,shift right=-0.5,dash] 
&{X \bbox K}\\
\end{tikzcd}\\
%%----------------------------------------------------------------------------------------------43
commutes.\\
\endgroup%%------------------------------------<<

\begingroup%%----------------------------------->>
\fontsize{9pt}{11pt}\selectfont
\textbf{\small EXAMPLE} \ 
If $\Box$ is a simplicial action on \bC, then for every small category \bI, the composition 
$[\bI,\bC] \times \bSISET \ra$ 
$[\bI,\bC] \times [\bI,\bSISET] \approx$ 
%$[\bI,\bC \times \bSISET] \overset{[\bI,\Box]}{\longrightarrow}[\bI,\bC]$ 
\begin{tikzcd}[sep=large]
{[\bI,\bC \times \bSISET] } \ar{r}{[\bI,\Box]} &{[\bI,\bC]}
\end{tikzcd}
is a simplicial action on $[\bI,\bC]$.\\
\endgroup%%------------------------------------<<

\begin{proposition} \ %27
Let \bC be a category.  Assume: \bC admits a simplicial action 
$\Box$ $-$then there is a simplicial category  
$\Box\bC$ such that \bC is isomorphic to the underlying category $\bU\Box\bC$.
\end{proposition}

[Put \ $O = \Ob\bC$ \ and assign to each ordered pair \ $X, \ Y \in O$ \ the simplicial set 
$\HOM(X,Y)$ defined by $\HOM(X,Y)_n = \Mor(X \bbox \dn, Y)$ $(n \geq 0)$.

\indent\indent (Composition) \  Given $X, \ Y, \ Z$, let 
$C_{X,Y,Z}:\ \HOM(X,Y) \times \HOM(Y,Z) \ \ra\ $ $\ \HOM(X,Z)$ be the simplicial map that sends 
$
\begin{cases}
\ f:X \bbox \dn \ra Y\\
\ g:Y \bbox \dn \ra Z
\end{cases}
$
to the composite 
\begin{tikzcd}%[sep=large]
{X \bbox \dn} \ar{r}{\id \Box \di} 
&{X \Box (\dn \times \dn)} \ar{r}{A}
&{(X \Box \dn) \Box \dn} \ar{r}{f \Box \id}
&{Y \Box \dn} \ar{r}{g}
&{Z.}
\end{tikzcd}
\\

\indent\indent (Unit) \ Given $X$, let 
$I_X:\Delta[0] \ra \HOM(X,X)$ be the simplicial map that sends $[n] \ra [0]$ to 
$X \bbox \dn \ra X \bbox \Delta[0] \overset{R}{\ra} X$.

Call $\Box \bC$ the simplicial category arising from this data.  That \bC is isomorphic to the underlying category 
$\bU \bbox \bC$ can be seen by considering the functor which is the identity on objects and sends a morphism 
$f:X \ra Y$ in \bC to $X \bbox \Delta[0] \overset{R}{\ra} X \overset{f}{\ra} Y$, an element of 
$\Mor(X \bbox \Delta[0],Y) = \HOM(X,Y)_0 \approx \Nat(\Delta[0],\HOM(X,Y))$.]

[Note: \  $\HOM:\bC^\OP \times \bC \ra \bSISET$ is a functor and the simplicial set $\HOM(X,Y)$ is called the \un{simplicial mapping space}
\index{simplicial mapping space} 
between $X$ and \mY.  
Example: Take for $\Box$ the trivial simplical action $-$then in this case, $\HOM(X,Y) = \si\Mor(X,Y)$.]\\
%\vspace{0.75cm}

\label{13.66}
Examples:
(1) \bSISET admits a simplicial action: $K \bbox L = K \times L$ (so $\HOM(K,L) = \map(K,L)$); 
(2) \bCGH admits a simplicial action: $X \bbox K = X \times_k \abs{K}$ (so $\HOM(X,Y)_n = $ all continuous functions $X \times_k \Delta^n \ra Y$); 
(3) $\bSISET_*$ admits a simplicial action: $K \bbox L = K \# L_+$ (so $\HOM(K,L) = \map_*(K,L)$); 
(4) $\bCGH_*$ admits a simplicial action: $X \bbox K = X \#_k \abs{K}_+$ (so $\HOM(X,Y)_n = $ all pointed continuous functions $X \#_k \Delta_+^n \ra Y$).

[Note: \  If $X, \ Y$ are in \bCGH, then 
$\HOM(X,Y) \approx \sin(\map(X,Y))$ 
and if $X, \ Y$ are in $\bCGH_*$, then 
$\HOM(X,Y) \approx \sin(\map_*(X,Y))$.  
In either situation, $\HOM(X,Y)$ is fibrant.]\\

\begingroup%%----------------------------------->>
\fontsize{9pt}{11pt}\selectfont
Neither \bTOP nor $\bTOP_*$ fits into the preceding framework (products or smash products are preserved in general only if the compactly generated category is used).  
This difficulty can be circumvented by restricting the definition of simplicial action to the full subcategory of \bSISET whose objects are the finite simplicial sets.  
It is therefore still the case that \bTOP ($\bTOP_*$) is isomorphic to the underlying category of 
%%----------------------------------------------------------------------------------------------44
a simplicial category with $\HOM(X,Y)_n = $ all continuous functions $X \times \Delta^n \ra Y$ (all pointed continuous functions 
$X \# \Delta_+^n \ra Y$).\\
\endgroup%%------------------------------------<<

\label{13.68}
\label{13.82a} %dmc mnft
\label{13.99} %dmc mnft
Example: Let \bC be a category.  Assume: \bC has coproducts $-$then \bSIC admits a simplicial action $\Box$ such that $\Box\bSIC$ is isomorphic to \bSIMC (cf. Proposition 26).

[Define $X \bbox K$ by $(X \bbox K)_n = K_n \cdot X_n$ (thus for $\alpha:[m] \lra [n]$, 
$K_n \cdot X_n \overset{X_\alpha}{\lra}$ 
$K_n \cdot X_m \overset{K_\alpha}{\lra} K_m \cdot X_m$).  
The symbol $X \bbox K$ also has another connotation 
(cf. p. \pageref{13.64}).  
To reconcile the ambiguity, note that there is a formality 
$\ini: X \bbox K \ra X \bbox K $, where 
$\ini_n(k):X_n \ra (X \bbox K)_n$ 
is the injection from $X_n$ to $K_n \cdot X_n$ corresponding to $k \in K_n$ 
(cf. p. \pageref{13.65}).  
Moreover, 
$\ini^*:\Nat(X \bbox K,Y) \ra \For(X \bbox K,Y)$ is bijective and functorial.  
Therefore $\Box\bSIC$ and $\bSIMC$ are isomorphic.]

[Note: \  $\Box$ is the 
\un{canonical} 
\index{simplicial action (canonical)} 
simplicial action on \bSIC.]\\

\begingroup%%----------------------------------->>
\fontsize{9pt}{11pt}\selectfont
\textbf{\small EXAMPLE} \ 
Let \bI be a small category $-$then there is an induced simplicial action on $[\bI,\bSISET]$ 
($(F \bbox K)_i = F i \times K$ 
(cf. p. \pageref{13.66})).  
And:
$\HOM(F,G) \approx \ds\int_i\map(F i,G i)$.  In fact, 
$\HOM(F,G)_n \approx$ 
$\Nat(F \bbox \dn,G) \approx$
$\ds\int_i \Nat(F i \times \dn,G i) \approx$
$\ds\int_i \Nat(\dn,\map(F i,G i)) \approx$
$\Nat(\dn,\ds\int_i\map(F i,G i)) \approx$
$\bigl(\ds\int_i \map(F i,G i)\bigr)_n$.\\
\endgroup%%------------------------------------<<

A simplicial action $\Box$ on a category \bC is said to be 
\un{cartesian}
\index{cartesian (simplicial action)} 
if $\forall \ X \in \Ob\bC$, the functor
$X \bbox -: \bSISET \ra \bC$ has a right adjoint.

Example: Let \bC be a category.  Assume: \bC has coproducts $-$then the canonical simplicial action $\Box$ on \bSIC is cartesian.

[Let \ $\HOM(X,Y)$ \ be the simplicial set figuring in the definition of \  \bSIMC, so 
$\HOM(X,Y)_n = \For(X \bbox \dn,Y)$ (cf. Proposition 26).  Define 
$\ev \in \For(X \bbox \HOM(X,Y),Y)$ by $\ev_n(f) = f_n(\id_{[n]}):X_n \ra Y_n$.  
Viewing ev as ``evaluation'', there is an induced functorial bijection 
$\Nat(K,\HOM(X,Y)) \ra \For(X \bbox K,Y)$.  However, 
$\For(X \bbox K,Y) \approx \Nat(X \bbox K,Y)$ (cf. supra), hence $\Box$ is cartesian.]\\

\begin{proposition} \ %28
Suppose that the simplicial action $\Box$ on \bC is cartesian $-$then $\forall \ X \in \Ob\bC$, 
$\HOM(X, -):\bC \ra \bSISET$ is a right adjoint for $X\Box -$.
\end{proposition}

[Given a simplicial set $K$, write \ $K = \colim_i \Delta[n_i]$ : \  
$\Mor(X \bbox K,Y)$ $\approx$ 
$\lim_i \Mor(X \bbox$ $\Delta[n_i],Y)$ $\approx$ 
$\lim_i \HOM(X,Y)_{n_i}$ $\approx$ 
$\lim_i \Nat(\Delta[n_i],\HOM(X,Y)) \approx$ 
$\Nat(K,\HOM(X,Y))$.]\\

A simplicial action $\Box$ on a category \bC is said to be 
\un{closed} 
\index{simplicial action (closed)} 
provided that it is cartesian and each of the functors $-\Box K:\bC \ra \bC$ has a right adjoint $X \ra \texttt{HOM}(K,X)$, so 
$\Mor(X \bbox K,Y) \approx \Mor(X,\texttt{HOM}(K,Y))$.

%%----------------------------------------------------------------------------------------------45
\label{13.85} %dmc mnft
[Note: \  The above defined simplicial actions on \bSISET, \bCGH, $\bSISET_*$, and $\bCGH_*$ are closed.]

If \bC admits a closed simplicial action, then $\bC^\OP$ admits a closed simplicial action.

Example: \bGRD admits a closed simplicial action:
$\bG \bbox K = $
$\bG \times \bPi K (\texttt{HOM}(K,\bG) =$ 
$[\bPi K,\bG])$.

[Note: \  Recall that $\bPi:\bSISET \ra \bGRD$ preserves finite products 
(cf. p. \pageref{13.67}).]\\

\begingroup%%----------------------------------->>
\fontsize{9pt}{11pt}\selectfont
\textbf{\small EXAMPLE} \ 
If $\bbox$ is a closed simplicial action on \bC, then for every small category \bI, the composition 
$[\bI,\bC] \times \bSISET \ra [\bI,\bC] \times [\bI,\bSISET] \approx$ 
$[\bI,\bC \times \bSISET]$ 
$\overset{[\bI,\bbox]}{\xrightarrow{\hspace*{1cm}}}$ 
$[\bI,\bC]$ 
is a closed simplicial action on 
$[\bI,\bC]$.\\
\endgroup%%------------------------------------<<

\begin{proposition} \ %29
Suppose that the simplicial action $\Box$ on \bC is closed $-$then 
$\text{HOM}(X \Box K,Y) \approx \text{map}(K,\HOM(X,Y)) \approx \text{HOM}(X, \texttt{HOM}(K,Y))$.\\
%%%% !!!!!!!!!!!!!!!!!!!!!!!!!!!!!!!!!!!!!!!!!!!!!!!!!!!!!!!!!!!!!!!!!!! note in the text the last HOM is written differently  ie small - find out why - it happens elsewhere as well - I have a notion - note the arguments - K sisets
\end{proposition}

\label{13.71a}

\begingroup%%----------------------------------->>
\fontsize{9pt}{11pt}\selectfont
\textbf{\small FACT} \ 
Let
$
\begin{cases}
\ \bC\\
\ \bD
\end{cases}
$
be categories equipped with closed simplicial actions.  Suppose that
$
\begin{cases}
\ F:\bC \ra \bD\\
\ G:\bD \ra \bC
\end{cases}
$
are functors and $(F,G)$ is an adoint pair.  
Assume: $\forall \ X$, $\forall \ K$: $F(X\bbox K) \approx FX \bbox K$ $-$then 
$\HOM(FX,FY)$ $\approx$ $\HOM(X,GY)$ 
and
$G\texttt{HOM}(K,Y) \approx \texttt{HOM}(K,GY)$.
\endgroup%%------------------------------------<<

Notation: Given a category \bC and a simplicial object $X$ in \bC, write $h_X$ for the cofunctor 
$\bC \ra \bSISET$ defined by $(h_X A)_n = \Mor(A,X_n)$.

[Note: \  For all $X, \ Y$ in \bSIC, $\Nat(X,Y) \approx \Nat(h_X,h_Y)$ (simplicial Yoneda).]\\

\begin{proposition} \ %30
Let \bC be a category.  
Assume: \bC has coproducts and is complete $-$then the canonical simplicial action 
$\Box$ on \textbf{SIC} is closed ($\Box$ is necessarily cartesian 
(cf. p. \pageref{13.68})).
\end{proposition}

[Given a simplicial set $K$, write \ $K \times \dn = \colim_i\Delta[n_i]$ : 
$\Nat(K \times \dn,h_Y A)$ $\approx$ 
$\lim_i\Nat(\Delta[n_i],h_Y A)$ $\approx$ 
$\lim_i \Mor(A,Y_{n_i})$ $\approx$ 
$\Mor(A,\lim_i Y_{n_i})$ $\approx$ 
$\Mor(A,\texttt{HOM}(K,Y)_n)$, \ 
where by definition \ 
$\texttt{HOM}(K,Y)_n = \lim_i Y_{n_i}$.  \ 
In other words, $\texttt{HOM}(K,Y)_n$ represents 
$A \ra \Nat(K \times \dn,h_Y A)$.  
Varying $n$ yields a simplicial object $\texttt{HOM}(K,Y)$ in \bC with 
$h_{\texttt{HOM}(K,Y)} \approx \map(K,h_Y)$.  
Agreeing to let $h_X \bbox K$ be the cofunctor $\bC \ra \bSISET$ that sends $A$ to 
$h_X A \times K$, we have
$\Nat(X \bbox K,Y)$ $\approx$ 
$\Nat(h_{X \bbox K},h_Y)$ $\approx$ 
$\Nat(h_X \bbox K,h_Y)$ $\approx$ 
$\Nat(h_X,\map(K,h_Y))$ $\approx$ 
$\Nat(h_X,h_{\texttt{HOM}(K,Y)}) \approx$ 
$\Nat(X,\texttt{HOM}(K,Y))$, 
which proves that $\Box$ is closed.]\\

Example: The canonical simplicial action $\Box$ on \bSIGR or \bSIAB is closed.\\

\begingroup%%----------------------------------->>
\fontsize{9pt}{11pt}\selectfont
\index{\mG-Sets (example)}
\textbf{\small EXAMPLE \  (\un{$G-$Sets})} \ 
Fix a group $G$ $-$then $\bSISET_G$ admits a canonical simplicial action $\Box$, viz. 
$X \bbox K = X \times K$, with trivial operations on $K$.  
In addition, $\bbox$ is closed, $\texttt{HOM}(K,X)$ being $\map(K,X)$ 
(operations in the target).  
Obviously, $F_G(X \bbox K) \approx F_G(X) \bbox K$.\\
\endgroup%%------------------------------------<<

%%----------------------------------------------------------------------------------------------46
A 
\un{simplicial model category}
\index{simplicial model category} 
is a model category \bC equipped with a closed simplicial action $\Box$ satisfying

\indent\indent (SMC) \  Suppose that $A \ra Y$ is a cofibration and $X \ra B$ is a fibration $-$then the arrow
$\text{HOM}(Y,X) \ra \text{HOM}(A,X)  \times_{\text{HOM}(A,B) } \text{HOM}(Y,B)$ 
is a Kan fibration which is a weak homotopy equivalence if $A \ra Y$ or $X \ra B$ is acyclic.

Observation: It is clear that SMC $\implies$ MC-4.  Indeed, the commutative diagram 
\begin{tikzcd}%[sep=large]
{A}  \ar{d} \ar{r} &{X} \ar{d}\\
{Y} \ar{r}  &{B}  
\end{tikzcd}
is a vertex of $\HOM(A,X) \times_{\HOM(A,B)}\HOM(Y,B)$, a filler $Y \ra X$ is a preimage in $\HOM(Y,X)_0$, and acyclic Kan fibrations are surjective.

Example: \bSISET, \bCGH, $\bSISET_*$, $\bCGH_*$ are simplicial model categories.

[Note: \  \bCGH and $\bCGH_*$ are taken in their singular structures. 
(cf. p. \pageref{13.69}).]\\

\label{13.82}
\label{13.88}

\begingroup%%----------------------------------->>
\fontsize{9pt}{11pt}\selectfont
\textbf{\small EXAMPLE} \ 
Fix a small category \bI $-$then the functor category [\bI,\bSISET] is a simplicial model category (use structure L 
(cf. p. \pageref{13.70})).\\
\endgroup%%------------------------------------<<

\label{13.98}
\begingroup%%----------------------------------->>
\fontsize{9pt}{11pt}\selectfont
\textbf{\small EXAMPLE} \ 
Fix an abelian group $G$ and take \bSISET in the \mc structure furnished by the homological \mc theorem.  Since every  $HG$-fibration is a Kan fibration, it follows from Propositions 23 and 25 that \bSISET is a simplicial model category.\\
\endgroup%%------------------------------------<<

\index{\mG-Sets (example simplicial model category)}
\begingroup%%----------------------------------->>
\fontsize{9pt}{11pt}\selectfont
\textbf{\small EXAMPLE \  (\un{$G-$Sets})} \ 
Fix a group $G$ $-$then $\bSISET_G$ is a simplicial model category 
(cf. p. \pageref{13.71}).\\
\endgroup%%------------------------------------<<

In a simplicial model category \bC: \ 
(1) $X \bbox \Delta[0] \approx X$;
(2) $\texttt{HOM}(\Delta[0],X) \approx X$; 
(3) $\emptyset \bbox K \approx \emptyset$;
(4) $\texttt{HOM}(K,*) \approx *$; 
(5) $\HOM(\emptyset,X) \approx \Delta[0]$;
(6) $\HOM(X,*) \approx \Delta[0]$; 
(7) $X \bbox \emptyset \approx \emptyset$; 
(8) $\texttt{HOM}(\emptyset,X) \approx *$.\\

\begin{proposition} \ %31
Suppose that $\Box$ is a closed simplicial action on a model category \bC $-$then \bC is a simplicial model category  iff whenever $A \ra Y$ is a cofibration in \bC and $L \ra K$ is an inclusion of simplicial sets, the arrow 
$A \Box K \underset{A \Box L}{\sqcup} Y \Box L \ra Y \Box K$ 
is a cofibration which is acyclic if $A \ra Y$ or $L \ra K$ is acyclic.\\
\end{proposition}
%\vspace{0.75cm}

\label{13.77}
\label{13.80}
\label{13.95}
\label{13.95a}

Application: Let \bC be a simplicial model category.

\indent\indent (i) Suppose that $A \ra Y$ is a cofibration in \bC $-$then for every simplicial set $K$, the arrow %dmcxx
$A \bbox K \ra Y \bbox K$ is a cofibration which is acyclic if $A \ra Y$ is acyclic.

\indent\indent (ii) Suppose that $Y$ is cofibrant and $L \ra K$ is an inclusion of simplicial sets $-$then the arrow 
$Y \bbox L \ra Y \bbox K$ is a cofibration which is acyclic if $L \ra K$ is acyclic.

[Note: \  In particular, $Y$ cofibrant $\implies$ $Y \bbox K$ cofibrant.]\\ 
%%----------------------------------------------------------------------------------------------47

\label{13.86}

\begingroup%%----------------------------------->>
\fontsize{9pt}{11pt}\selectfont
\textbf{\small FACT} \ 
Suppose that $\Box$ is a closed simplicial action on a \mc \bC $-$then \bC is a simplicial \mc iff whenever $A \ra Y$ is a cofibration in \bC, the arrows 
$A \bbox \Delta[n] \underset{A \bbox \dot\Delta[n]}{\sqcup} Y \bbox \dot\Delta[n] \ra  Y \bbox \Delta[n]$ $(n \geq 0)$ 
are cofibrations which are acyclic if $A \ra Y$ is acyclic and the arrows 
$A \bbox \Delta[1] \underset{A \bbox \Lambda[i,1]}{\sqcup} Y \bbox \Lambda[i,1]\ra  Y \bbox \Delta[1]$ $(i = 0, 1)$ 
are acyclic cofibrations.\\
\endgroup%%------------------------------------<<

\begin{proposition} \ 
Suppose that $\Box$ is a closed simplicial action on a \mc \bC $-$then \bC is a simplicial \mc iff whenever $L \ra K$ is an inclusion of simplicial sets and $X \ra B$ is a fibration in \bC, the arrow 
$\texttt{HOM}(K,X) \ra \texttt{HOM}(L,X) \times_{\texttt{HOM}(L,B)} \texttt{HOM}(K,B)$ 
is a fibration which is acyclic if $L \ra K$ or $X \ra B$ is acyclic.\\
\end{proposition}

\label{13.78}
\label{13.81}
Application: Let \bC be a simplicial model category.\\
\indent\indent (i) \  Suppose that $L \ra K$ is an inclusion of simplicial sets and $X$ is fibrant $-$then the arrow 
$\texttt{HOM}(K,X) \ra \texttt{HOM}(L,X)$ is a fibration which is acyclic if $L \ra K$ is acyclic.\\
\indent\indent (ii) \  Suppose that $X \ra B$ is a fibration in \bC $-$then for every simplicial set $K$, the arrow 
$\texttt{HOM}(K,X) \ra \texttt{HOM}(K,B)$ is a fibration which is acyclic if $X \ra B$ is acyclic.

[Note: \  In particular, $X$ fibrant $\implies$ $\texttt{HOM}(K,X)$ fibrant.]\\

\begingroup%%----------------------------------->>
\fontsize{9pt}{11pt}\selectfont
\textbf{\small FACT} \ 
Suppose that $\Box$ is a closed simplicial action on a \mc \bC $-$then \bC is a simplicial \mc iff 
whenever \ $X \ra B$ \ is a fibration in \bC, \ the arrows \  
$\texttt{HOM}(\Delta[n],X) \ra$ 
$\texttt{HOM}(\dot\Delta[n],X)$ $\times_{\texttt{HOM}(\dot\Delta[n],B)} \texttt{HOM}(\Delta[n],B)$ $(n \geq 0)$ 
are fibrations which are acyclic if $X \ra B$ is acyclic and the arrows \ 
$\texttt{HOM}(\Delta[1],X) \ra \texttt{HOM}(\Lambda[i,1],X) \times_{\texttt{HOM}(\Lambda[i,1],B)}\texttt{HOM}(\Delta[1],B)$ $(i = 0, 1)$ 
are acyclic fibrations.\\
\endgroup%%------------------------------------<<

\label{13.119}

Example: Let \bC be a category.  Assume: \bC is complete and cocomplete and there is an adjoint pair $(F,G)$ where
$
\begin{cases}
\ F:\bSISET \ra \bSIC\\
\ G:\bSIC \ra \bSISET
\end{cases}
, \ 
$
subject to the requirement that $G$ preserves filtered colimits.  
Call a morphism $f:X \ra Y$ a \we if $G f$ is a weak homotopy equivalence, 
a fibration if $G f$ is a Kan fibration, a cofibration if $f$ has the LLP w.r.t. acyclic fibrations 
$-$then \bSIC is a \mc provided that every cofibration with the LLP w.r.t. fibrations is a \we (cf. infra).  
Claim: \bSIC is a simplicial model category ($\Box = $ canonical simplicial action (cf. Proposition 30)).  
To see this, note first that 
$F(X \times K) \approx F X \bbox K$, hence 
$G \texttt{HOM}(K,Y) \approx \map(K,GY)$ 
(cf. p.  \pageref{13.71a}).  
Let now $L \ra K$ be an inclusion of simplicial sets and $X \ra B$ be a fibration in \bSIC.  Apply $G$ to the arrow 
$\texttt{HOM}(K,X) \ra \texttt{HOM}(L,X) \times_{\texttt{HOM}(L,B)} \texttt{HOM}(K,B)$ 
to get
$G\texttt{HOM}(K,X) \ra G\texttt{HOM}(L,X) \times_{G\texttt{HOM}(L,B)} G\texttt{HOM}(K,B)$
or still, 
$\map(K,GX) \ra \map(L,GX) \times_{\map(L,GB)} \map(K,GB)$.
Taking into account Proposition 22 and the definitions, the claim thus follows from Proposition 32.

%%----------------------------------------------------------------------------------------------48
[Note: \  Typically, such a setup is realized in ``algebraic'' situations (consider, e.g., \bC = \textbf{GR}).  Consult 
Crans\footnote[2]{\textit{J. Pure. Appl. Algebra} \textbf{101} (1995), 35-57.} 
for a variation on the overall procedure with applications to simplicial sheaves.]\\

\begingroup%%----------------------------------->>
\fontsize{9pt}{11pt}\selectfont
The model category structure on \bSIC is produced by a small object argument.  
Thus one works with the 
$F \ddn \ra F \dn$ $(n \geq 0)$ 
to show that every $f$ can be written as the composite of a cofibration and an acyclic fibration and one works with the 
$F \Lambda[k,n] \ra F\dn$ $(0 \leq k \leq n, n \geq 1)$ to show that every $f$ can be written as the composite of a cofibration that has the LLP w.r.t fibrations and a fibration.  This leads to MC-5 under the assumption that every cofibration with the LLP w.r.t fibrations is a weak equivalence, which is also needed to establish the nontrivial half of MC-4.  In practice, this condition can be forced.\\
\endgroup%%------------------------------------<<

\begingroup%%----------------------------------->>
\fontsize{9pt}{11pt}\selectfont
\textbf{\small SUBLEMMA} \ 
Let 
$
\begin{cases}
\ X\\
\ Y
\end{cases}
$
be topological spaces, $f:X \ra Y$ a continuous function; let $\phi:X^\prime \ra X$, $\psi:Y \ra Y^\prime$ be continuous functions.  
Assume: $f \circx \phi$, $\psi \circx f$ are weak homotopy equivalences $-$then $f$ is a weak homotopy equivalence.\\
\endgroup%%------------------------------------<<

\begingroup%%----------------------------------->>
\fontsize{9pt}{11pt}\selectfont
\textbf{\small LEMMA} \ 
Suppose that there is a functor $T: \bSIC \ra \bSIC$ and a natural transformation 
$\epsilon:\id_{\bSIC} \ra T$ such that $\forall \ X$, $\epsilon_X:X \ra TX$ is a \we and $TX \ra *$ is a fibration 
$-$then every cofibration with the LLP w.r.t. fibrations is a \we.

[Let $i:A \  \ra \  Y$ be a cofibration with the stated properties.  
Fix a filler $w:Y \ra TA$ for
$
\begin{tikzcd}[sep=large]
{A}  \ar{d}[swap]{i} \ar{r}{\epsilon_A} &{TA} \ar{d}\\
{Y} \ar{r}  &{*}  
\end{tikzcd}
. \ 
$  
Consider the commutative diagram
$
\begin{tikzcd}[sep=large]
{A}  \ar{d}[swap]{i} \ar{r}{f} &{\texttt{HOM}(\Delta[1],TY)} \ar{d}{\Pi}\\
{Y} \ar{r}[swap]{g}  &{\texttt{HOM}(\dot\Delta[1],TY)}  
\end{tikzcd}
,  \ 
$
where $f$ is the arrow 
$A \overset{i}{\lra} Y  \overset{\epsilon_Y}{\lra}TY \approx \texttt{HOM}(\Delta[0],TY) \ra \texttt{HOM}(\Delta[1],TY)$ 
and $g$ is the arrow
$
\begin{cases}
\ Y  \overset{\epsilon_Y}{\ra}TY\\
\ Y  \overset{w}{\ra}TA \overset{T i}{\ra} TY
\end{cases}
$
$(\texttt{HOM}(\dot\Delta[1],TY) \approx TY \times TY)$.  
Since $G T Y$ is fibrant and 
$
\begin{cases}
\ G\texttt{HOM}(\Delta[1],TY) \approx \map(\Delta[1],GTY)\\
\ G\texttt{HOM}(\dot\Delta[1],TY) \approx \map(\dot\Delta[1],GTY)\\
\end{cases}
$
, $\Pi$ is a fibration 
(cf. p. \pageref{13.72}), 
thus our diagram admits a filler $Y \ra \texttt{HOM}(\Delta[1],TY)$.
This in turn implies that $Ti \circx w$ is a weak equivalence, i.e., $\abs{G Ti} \circx \abs{G w}$ is a weak homotopy equivalence.  
Assemble the data: 
$\abs{G A} \overset{\abs{G i}}{\lra} \abs{G Y} \overset{\abs{G w}}{\lra}  \abs{G T A} \overset{\abs{G T i}}\lra  \abs{G T Y}$.
Because $\abs{G w} \circx \abs{G i} = \abs{G \epsilon_A}$ is a weak homotopy equivalence, one can apply the sublemma and conclude that $\abs{G w}$ is a weak homotopy equivalence.  Therefore $\abs{G i}$ is a weak homotopy equivalence which means by definition that $i$ is a weak equivalence.]\\
\endgroup%%------------------------------------<<

\begingroup%%----------------------------------->>
\fontsize{9pt}{11pt}\selectfont
\textbf{\small EXAMPLE} \ 
The hypotheses of the lemma are trivially met if $\forall \ X$, $X \ra *$ is a fibration.  So, for instance, \bSIC is a simplicial model category when 
\bC = \bGR, \bAB, or \bAMOD, $G$ being the forgetful functor.\\
\endgroup%%------------------------------------<<

%%----------------------------------------------------------------------------------------------49

\begingroup%%----------------------------------->>
\fontsize{9pt}{11pt}\selectfont
Retaining the supposition that \bC is complete and cocomplete, let us assume in addition that \bC has a set of separators and is cowellpowered.  Given a simplicial object $X$ in \bC, the cofunctor $\bC \ra \bSET$ defined by 
$A \ra (\Ex \HOM(A,X))_n$ $(n \geq 0)$ is representable (view $A$ as a constant simplicial object).  Indeed, 
$\HOM(-,X)$ converts colimits into limits and Ex preserves limits.  The assertion is then a consequence of the special adjoint functor theorem.  Accordingly, $\exists$ an object $(\Ex X)_n$ in \bC and a natural isomorphism 
$\Mor(A,(\Ex X)_n) \approx (\Ex \HOM(A,X))_n$.  
Thus there is a functor 
$\Ex:\bSIC \ra \bSIC$, where $\forall \ X$, 
$\Ex X([n]) = (\Ex X)_n$ $(n \geq 0)$, with 
$\HOM(A,\Ex X) \approx \Ex\HOM(A,X)$ (since 
$\HOM(A,\Ex X)_n \approx$
$\Nat(A \bbox \Delta[n],\Ex X) \approx$
$\Mor(A,(\Ex X)_n) \approx$
$(\Ex \HOM(A,X))_n) $.  
Iterate to arrive at
$\Exx^\infty:\bSIC \ra \bSIC$
and 
$\epsilon^\infty:\id_{\bSIC} \ra \Exx^\infty$.\\
\endgroup%%------------------------------------<<

\begingroup%%----------------------------------->>
\fontsize{9pt}{11pt}\selectfont
\index{Small Object Contstruction}
\textbf{\small SMALL OBJECT CONSTRUCTION} \ 
Fix a $P \in \Ob\bC$ such that $\Mor(P,-):\bC \ra \bSET$ preserves filtered colimits.  
Viewing $P$ as a constant simplicial object, define 
$G:\bSIC \ra \bSISET$ 
by 
$GX \ $=$ \  \HOM(P,X)$ 
$-$then $G$ has a left adjoint $F$, viz. $FX$ $=$ $P \ \Box \ K$, and $G$ preserves filtered colimits \  (for 
$(G \hspace{0.03cm}\colimx X_i)_n \approx$ 
$\HOM(P,\colimx X_i)_n \approx$ 
$\Nat(P \ \Box \ \Delta[n], \colimx X_i)$ $\approx$ 
$\Mor(P,(\colimx X_i)_n)$ $\approx$ \ 
$\Mor(P,\colimx(X_i)_n)$ \ $\approx$ \ 
%%%
$\colimx \Mor(P,(X_i)_n)$ \ $\approx$ \ 
$\colimx \Nat(P \ \Box \ \Delta[n],X_i)$ $\approx$ \ 
$\colimx \HOM(P,X_i)_n$ $\approx$ \ 
$(\colimx G X_i)_n$).  
In the lemma, take \ 
$T = \Exx^\infty$, $\epsilon = \epsilon^\infty$.  \ \
Because \ 
$\HOM(P,\Exx^\infty X)$ $\approx$ \ 
$\HOM(P,\colimx\Exx^n X)$ $\approx$ \ 
$\colimx \HOM(P,\Exx^n X)$ $\approx$ 
$\Exx^\infty \HOM(P,X)$, 
it follows that $\forall \ X$, 
$\epsilon_X^\infty:X \ra \Exx^\infty X$ is a weak equivalence 
(cf. p. \pageref{13.73}) 
and $\Exx^\infty X \ra *$ is a fibration 
(cf. p. \pageref{13.73a}).  
Therefore \bSIC admits the structure of a simplicial model category in which a morphism $f:X \ra Y$ is a weak equivalence or a fibration if this is the case of the simplicial map $f_*:\HOM(P,X) \ra \HOM(P,Y)$.\\
\endgroup%%------------------------------------<<

\begingroup%%----------------------------------->>
\fontsize{9pt}{11pt}\selectfont
\textbf{\small EXAMPLE} \ 
In the small object construction, take \bC = \bSISET $-$then every finite simplicial set $P$ determines a simplicial model category structure on $[\bDelta^\OP,\bSISET]$.\\
\endgroup%%------------------------------------<<

\begin{proposition} \ %33
Let $X$, $Y$, and $Z$ be objects in a simplicial model category \bC.
\indent\indent (i) \  If $f:X \ra Y$ is an acyclic cofibration and $Z$ is fibrant, then 
$f^*:\text{HOM}(Y,Z) \ra \text{HOM}(X,Z)$ is a weak homotopy equivalence.\\
\indent\indent (ii) \ If $g:Y \ra Z$ is an acyclic fibration and $X$ is cofibrant, then 
$g_*:\text{HOM}(X,Y) \ra \text{HOM}(X,Z)$ is a weak homotopy equivalence.\\
\end{proposition}

\begin{proposition} \ %34
Let $X$, $Y$, and $Z$ be objects in a simplicial model category \bC.
\indent\indent (i) \ If $f:X \ra Y$ is a weak equivalence between cofibrant objects and $Z$ is fibrant, then 
$f^*:\text{HOM}(Y,Z) \ra \text{HOM}(X,Z)$ is a weak homotopy equivalence.\\
\indent\indent (ii) \ If $g:Y \ra Z$ is a weak equivalence between fibrant objects and $X$ is cofibrant, then 
$g_*:\text{HOM}(X,Y) \ra \text{HOM}(X,Z)$ is a weak homotopy equivalence.
\end{proposition}

[Use Proposition 33 and the lemma prefacing the proof of the TDF theorem.]\\

\begingroup%%----------------------------------->>
\fontsize{9pt}{11pt}\selectfont
\textbf{\small EXAMPLE} \ 
Take $\bC = \bCGH$ (singular structure) $-$then all objects are fibrant, so if $g:Y \ra Z$ is a weak homotopy equivalence and $X$ is cofibrant, 
$g_*:\HOM(X,Y) \ra \HOM(X,Z)$ is a weak homotopy 
%%----------------------------------------------------------------------------------------------50
equivalence.  But $\HOM(X,Y) \approx \sin(\map(X,Y))$, $\HOM(X,Z) \approx \sin(\map(X,Z))$, thus 
$g_*:\map(X,Y) \ra \map(X,Z)$ is a weak homotopy equivalence 
(cf. p. \pageref{13.74}).
\vspi
[Note: Contrast this approach with that used on 
p. \pageref{13.75}.]\\
\endgroup%%------------------------------------<<

Let $i:A \ra Y$, $p:X \ra B$ be morphisms in a simplicial model category \bC.  
Assume: $i$ is a cofibration and $p$ is a fibration $-$then $i$ is said to have the 
\un{homotopy left lifting } \un{property with respect to $p$}
\index{homotopy left lifting property with respect to $p$} 
(HLLP w.r.t. $p$) 
\index{HLLP w.r.t. $p$} and $p$ is said to have the 
\un{homotopy right} \un{lifting property with respect to $i$}
\index{homotopy right lifting property with respect to i} (HRLP w.r.t $i$) 
\index{HRLP w.r.t $i$} 
if the arrow 
$\HOM(Y,X) \ra$ $\HOM(A,X)$ $\times_{\HOM(A,B)}$ $\HOM(Y,B)$ is a weak homotopy equivalence.\\

\begingroup%%----------------------------------->>
\fontsize{9pt}{11pt}\selectfont
\textbf{\small FACT} \ 
Given a cofibration $i:A \ra Y$ and a fibration $p:X \ra B$ in a simplicial model category \bC, 
each of the following conditions is equivalent to $i$ having the HLLP w.r.t $p$ and $p$ having that HRLP w.r.t. $i$.
\\
\indent\indent (1) \ If $L \ra K$ is an inclusion of simplicial sets, then $p$ has the RLP w.r.t. the arrow 
$A \bbox K \underset{A \Box L}{\sqcup} Y \bbox L \ra Y \bbox K$.
\\
\indent\indent (2) \ The fibration $p$ has the RLP w.r.t. the arrows 
$A \bbox \dn \underset{A \Box \ddn}{\sqcup} Y \bbox \ddn \ra Y \bbox \dn$ $(n \geq 0)$.
\\
\indent\indent (3) \ If $L \ra K$ is an inclusion of simplicial sets, then $i$ has the LLP w.r.t. the arrow 
$\texttt{HOM}(K,X) \ra \texttt{HOM}(L,X) \times_{\texttt{HOM}(L,B)} \texttt{HOM}(K,B)$.
\\
\indent\indent (4) \ The cofibration $i$ has the LLP w.r.t. the arrows 
$\texttt{HOM}(\dn,X) \ra \texttt{HOM}(\ddn,X) \times_{\texttt{HOM}(\ddn,B)}\texttt{HOM}(\dn,B)$ $(n \geq 0)$.\\
\endgroup%%------------------------------------<<

Let \bC be a simplicial model category.  
Agreeing to identify $\Mor(X,Y)$ and $\HOM(X,Y)_0$,  one may transfer from \bSISET to \bC the notions of 
\un{homotopic} 
\index{homotopic (in a simplicial model category)} 
$(f \simeq g)$ and \un{simplicially} \un{homotopic} 
\index{simplicially homotopic (in a simplicial model category)} 
$(f \underset{s}{\simeq} g)$ leading thereby to \bHZEROC 
\index{\bHZEROC} 
\ 
(thus $[X,Y]_0 = \Mor(X,Y)/\simeq $ \ $(\equiv \ \pi_0(\HOM(X,$ $Y))$).

[Note: \  $\Mor(X\Box I_{2n},Y)$ \  $\approx$  \ 
$\Nat(I_{2n},\HOM(X,Y)) \approx$ \  
$\Mor(X,\texttt{HOM}(I_{2n},Y))$ \ 
and \ 
$\Mor(X\Box\dw,Y) \approx \Nat(\dw,\HOM(X,Y)) \approx \Mor(X,\texttt{HOM}(\dw,Y))$.]

Example: Suppose that $i:A \ra Y$ is a cofibration and $p:X \ra B$ is a fibration.  
Assume: $i$ has the HLLP w.r.t. $p$ $-$then every commutative diagram 
\begin{tikzcd}%[sep=large]
{A}  \ar{d} \ar{r} &{X} \ar{d}\\
{Y} \ar{r}  &{B}
\end{tikzcd}
has a filler and any two such are homotopic.\\

\begin{proposition} \ %35
Let \bC be a simplicial model category.  Suppose that $f \simeq g$ $-$then $f$, $g$ are left homotopic and right homotopic.
\end{proposition}

[Note: \  Therefore $Qf = Qg$ 
(cf. p. \pageref{13.76}).  
Corollary: A homotopic equivalence in \bC is a weak equivalence (but not conversely).]\\
%%----------------------------------------------------------------------------------------------51

\begin{proposition} \ %36
Let \bC be a simplicial model category.  
Assume: $X$ is cofibrant and $Y$ is fibrant $-$then the relations of homotopy, simplicial homotopy, left homotopy, and right homotopy on 
$\text{Mor}(X,Y)$ coincide and are equivalence relations.  
Therefore ``homotopy is homotopy'' and $[X,Y]_0 \leftrightarrow [X,Y]$.
\end{proposition}

[Note: \  $\text{HOM}(X,Y)$ is necessarily fibrant (cf. SMC).]\\

\begingroup%%----------------------------------->>
\fontsize{9pt}{11pt}\selectfont
\textbf{\small EXAMPLE} \ 
Under the assumption that $X$ is cofibrant and $Y$ is fibrant, 
$[X \ \Box \ K,Y]$ $\approx$ $[K,$ $\HOM(X,Y)]$ $\approx$ $[X,\texttt{HOM}(K,Y)]$. %%check on homs sizes ??????????????
\vspi
[Note: \  Bear in mind that $X \ \Box K$ is cofibrant 
(cf. p. \pageref{13.77}) 
and $\texttt{HOM}(K,Y)$ is fibrant 
(cf. p. \pageref{13.78}).]\\
\endgroup%%------------------------------------<<

\begin{proposition} \ %37
Let $X$, $Y$, and $Z$ be objects in a simplicial model category \bC.\\
\indent\indent (i)  \ Let $f \in \Mor(X,Y)$ $-$then the homotopy class of the precomposition arrow 
$f^*:\HOM(Y,Z) \ra \HOM(X,Z)$ depends only on the homotopy class of $f$.

[Note: \   Thus \ $f^*$  \ is \  a \ homotopy \ equivalence \ of \ simplicial \ sets \  
if $f$\  is a \ homotopy equivalence.]\\
\indent\indent (ii) \ Let $g \in \Mor(Y,Z)$ $-$then the homotopy class of the postcomposition arrow 
$g_*:\HOM(X,Y) \ra \HOM(X,Z)$ depends only on the homotopy class of $g$.
\end{proposition}

[Note: \  Thus $g_*$ is a homotopy equivalence of simplicial sets if $g$ is a homotopy  equivalence.]\\

\begin{proposition} \ %38
Suppose that \bC is a simplicial model category.  Let $f \in \Mor(X,Y)$.  Assume: The precomposition arrows 
$
\begin{cases}
\ \HOM(Y,X) \ra \HOM(X,X)\\
\ \HOM(Y,Y) \ra \HOM(X,Y)
\end{cases}
$
are weak homotopy equivalences $-$then $f$ is a homotopy equivalence.
\end{proposition}

[Note: \   \ The result can also be  formulated in terms of  postcomposition arrows \\
$
\begin{cases}
\ \HOM(X,X) \ra \HOM(X,Y)\\
\ \HOM(Y,X) \ra \HOM(Y,Y)
\end{cases}
.]
$
\\
\vspace{0.25cm}

\begin{proposition} \ %39
Let \bC be a simplicial model category $-$then a morphism $f:X \ra Y$ is a \we if $\forall$ fibrant $Z$, the precomposition arrow 
$f^*:\HOM(Y,Z) \ra \HOM(X,Z)$ is a weak homotopy equivalence.
\end{proposition}

[Using the notation of Lemma $\sR$ 
(cf. p. \pageref{13.79}), 
consider the commutative diagram
\begin{tikzcd}%[sep=large]
{X}  \ar{d}[swap]{\iota_X} \ar{r}{f} &{Y} \ar{d}{\iota_Y}\\
{\sR X} \ar{r}[swap]{\sR f}  &{\sR Y}  
\end{tikzcd}
and apply $\HOM(-,Z)$ to get
\begin{tikzcd}%[sep=large]
{\HOM(X,Z)}  &{\HOM(Y,Z)} \ar{l}\\
{\HOM(\sR X,Z)} \ar{u}  &{\HOM(\sR Y,Z)}  \ar{l} \ar{u}
\end{tikzcd}
(Z fibrant).  Since 
$
\begin{cases}
\ \iota_X\\
\ \iota_Y
\end{cases}
$
are acyclic cofibrations, the vertical arrows are weak homotopy equivalences (cf. Proposition 33).  
Taking into account the hypothesis, it follows that
%%----------------------------------------------------------------------------------------------52
$(\sR f)^*:\HOM(\sR Y,Z) \ra \HOM(\sR X,Z)$
is a weak homotopy equivalence.  But
$
\begin{cases}
\ \sR X\\
\ \sR Y
\end{cases}
$
are fibrant, so one can let $Z = \sR X, \ \sR Y$ and conclude that $\sR f$ is a homotopy equivalence (cf. Proposition 38), hence a weak equivalence (cf. Proposition 35).  Therefore $f$ is a weak equivalence (cf. Lemma $\sR$).]

[Note: \  The result can also be formulated in terms of the postcomposition arrows 
$f_*:\HOM(Z,X) \ra \HOM(Z,Y)$ (\mZ cofibrant).]\\

\label{13.95b}
Application: Let \bC be a simplicial model category.  Suppose that $f:X \ra Y$ is a weak equivalence between cofibrant objects $-$then $\forall$ $K$, 
$f \ \Box \id_K: X \ \Box \ K \ra Y \ \Box \ K$ is a weak equivalence between cofibrant objects 
(cf. p. \pageref{13.80}).

[Take any fibrant $Z$ and consider the arrow
$\HOM(Y \ \Box \ K, Z) \ra \HOM(X \ \Box \ K, Z)$
or still, the arrow
$\HOM(Y,\texttt{HOM}(K,Z)) \ra \HOM(X,\texttt{HOM}(K,Z))$.  
Because $\texttt{HOM}(K,Z)$ is fibrant 
(cf. p. \pageref{13.81}), 
the latter is a weak homotopy equivalence (cf. Proposition 34), so by the above, the arrow 
$X \ \Box \ K \ra Y \ \Box \ K$ is a \we.]\\

\label{13.123}

\begingroup%%----------------------------------->>
\fontsize{9pt}{11pt}\selectfont
\textbf{\small EXAMPLE} \ 
Fix a small category \bI and view the functor category 
$[\bI^\OP,\bSISET]$ as a simplicial model category 
(cf. p. \pageref{13.82}).  
Suppose that $L \ra K$ is a weak equivalence, where 
$L,\ K:\bI^\OP \ra \bSISET$ are cofibrant $-$then $\forall \ F:\bI \ra \bSISET$, the induced map
$\ds\int^i F i \times L i \ra \ds\int^i F i \times K i$
of simplicial sets is a weak homotopy equivalence.  To see this, use Proposition 39.  Thus take any fibrant $Z$ and consider the arrow
$\map(\ds\int^i F i \times K i,Z) \ra \map(\ds\int^i F i \times L i, Z)$ 
i.e., the arrow 
$\ds\int_i \map(Fi \times Ki,Z) \ra \ds\int_i \map(Fi \times Li,Z)$, 
i.e., the arrow 
$\ds\int_i \map(Ki, \map(Fi,Z)) \ra \ds\int_i\map(Li, \map(Fi,Z))$
i.e., the arrow 
$\HOM(K,\map(F,Z)) \ra \HOM(L,\map(F,Z))$ 
(cf. p. \pageref{13.82a}), 
which is a weak homotopy equivalence (cf. Proposition 34).
\vspi
[Note: \  Here, $\map(F,Z)$ is the functor 
$\bI^\OP \ra \bSISET$ defined by $i \ra \map(Fi,Z)$ thus $\map(F,Z)$ is a fibrant object in $[\bI^\OP,\bSISET]$.]\\
\endgroup%%------------------------------------<<

Let $\rho:A \ra B$ be an inclusion of simplicial sets $-$then a fibrant object $Z$ in \bSISET is said to be 
\un{$\rho$-local} 
\index{local! $\rho$-local (fibrant object)} 
if $\rho^*:\map(B,Z) \ra \map(A,Z)$ is a weak homotopy equivalence.

[Note: \  Since $Z$ is fibrant, $\rho^*$ is actually a simplicial homotopy equivalence (cf. Proposition 20).]

Imitating the $(A,B)$ construction in $\S 9$ 
(cf. p. \pageref{13.82b} ff.), 
one can show that there is a functor 
$L_\rho:\bSISET \ra \bSISET$ and a natural transformation $\id \ra L_\rho$, where $\forall \ X$, $L_\rho X$ is $\rho$-local and 
$l_\rho:X \ra L_\rho X$ is a cofibration such that for all $\rho$-local $Z$, the arrow  
$\map(L_\rho X,Z) \ra \map(X,Z)$ is a weak homotopy equivalence.  
Consequently, the full subcategory of \bHZEROSISET whose objects are $\rho$-local is reflective.

%%----------------------------------------------------------------------------------------------53
[Note: \  Observe that it is necessary to work not only with 
$A \times \Delta[n] \underset{A \times \dot\Delta[n]}{\sqcup} B \times \dot\Delta[n] \ra B \times \Delta[n]$ $(n \geq 0)$ 
but also with the 
$\Lambda[k,n] \ra \Delta[n]$ $(0 \leq k \leq n, n \geq 1)$ (this to insure that $L_\rho X$ is fibrant).]\\

\textbf{\small LEMMA} \ 
Let $f:X \ra Y$ be a cofibration in \bSISET.  
Assume: $\forall$ $\rho$-local $Z$, $f^*:\map(Y,Z) \ra \map(X,Z)$ is a weak homotopy equivalence $-$then 
$L_\rho f: L_\rho X \ra L_\rho Y$ is a homotopy equivalence.

[Pass from 
\begin{tikzcd}%[sep=large]
{X}  \ar{d} \ar{r}{f} &{Y} \ar{d}\\
{L_\rho X} \ar{r}[swap]{L_\rho f}  &{L_\rho Y}  
\end{tikzcd}
to 
\begin{tikzcd}%[sep=large]
{\map(X,Z)}   &{\map(Y,Z)} \ar{l}\\
{\map(L_\rho X,Z)} \ar{u}  &{\map(L_\rho Y,Z)} \ar{l} \ar{u}
\end{tikzcd}
(Z $\rho$-local), take $Z = L_\rho X, \ L_\rho Y$, and quote Proposition 38.]\\

Application: Suppose that $f:X \ra Y$ is an acyclic cofibration $-$then 
$L_\rho f: L_\rho X \ra L_\rho Y$ is a homotopy equivalence.

[Note: \   Therefore $ L_\rho: \bSISET \ra \bSISET$ preserves weak homotopy equivalences 
(cf. p. \pageref{13.83}) 
(all objects are cofibrant), hence
$\bL L_\rho: \bHSISET \ra \bHSISET$ exsists (cf. $\S 12$, Proposition 14).]\\

\begingroup%%----------------------------------->>
\fontsize{9pt}{11pt}\selectfont
\textbf{\small EXAMPLE} \ 
Fix an inclusion $\rho:A \ra B$ of simplicial sets.  
Let $f:X \ra Y$ be a simplicial map $-$then $f$ is said to be a 
\un{$\rho$-equivalence} 
\index{equivalence! $\rho$-equivalence (simplicial set)} 
if $L_\rho f: L_\rho X \ra L_\rho Y$ is a homotopy equivalence (or just a weak homotopy equivalence (cf. Proposition 20)).  
Agreeing  that a 
\un{$\rho$-cofibration}
\index{cofibration! $\rho$-cofibration (simplicial set)} 
is an injective simplicial map, a 
\un{$\rho$-fibration}
\index{fibration! $\rho$-fibration (simplicial set)} 
is a simplicial map which has the RLP w.r.t all $\rho$-cofibrations that are $\rho$-equivalences.  
Every $\rho$-fibration is a Kan fibration (cf. supra).  
This said, \bSISET acquires the structure of a simplicial model category by letting 
weak equivalence = $\rho$-equivalence, 
cofibration = $\rho$-cofibration, 
fibration =$\rho$-fibration.
\vspi
[Note: \  The fibrant objects in this structure are the $\rho$-local objects.]\\
\endgroup%%------------------------------------<<

Let \bC be a complete and cocomplete category $-$then in the notation of 
p. \pageref{13.84}, 
the truncation 
$\tr^{(n)}:\bSIC \ra \bSIC_n$ 
has a left adjoint 
$\sk^{(n)}:\bSIC_n \ra \bSIC$, 
where $\forall$ $X$ in $\bSIC_n$, 
%$(\sk^{(n)} X)_m = \colim\limits_{\substack{[m] \ra [k] \\k \leq n}} X_k$, 
$(\sk^{(n)} X)_m = \underset{\substack{[m] \ra [k] \\k \leq n}}{\colim} X_k$,
and a right adjoint 
$\cosk^{(n)}:\bSIC_n \ra \bSIC$, 
where $\forall$ $X$ in $\bSIC_n$,
$(\cosk^{(n)} X)_m = \lim\limits_{\substack{[k] \ra [m]\\k \leq n}} X_k$.

[Note: \  The colimit and limit are taken over a comma category.]\\

\index{Extension Principle  (Objects)}
\textbf{\small EXTENSION PRINCIPLE \  (OBJECTS)} \ 
Let $X$ be an object in $\bSIC_{n}$ $-$then a factorization 
$(\sk^{(n)} X)_{n+1} \ra X_{n+1} \ra (\cosk^{(n)} X)_{n+1}$
of the arrow 
$(\sk^{(n)} X)_{n+1} \ra (\cosk^{(n)} X)_{n+1}$
determines an extension of $X$ to an object in $\bSIC_{n+1}$.\\

%%----------------------------------------------------------------------------------------------54
\index{Extension Principle (Morphisms)}
\textbf{\small EXTENSION PRINCIPLE \  (MORPHISMS)} \ \ 
Let 
$
\begin{cases}
\ X\\
\ Y
\end{cases}
$
be objects in $\bSIC_{n+1}$; \ let \ \ 
$f:\restr{X}{{\bDelta_n^\OP}} \ra \restr{Y}{{\bDelta_n^\OP}}$ be a morphism $-$then the arrow 
$X_{n+1} \ra Y_{n+1}$ determines an extension $F:X \ra Y$ of $f$ provided that 
\begin{tikzcd}[sep=large]
{(\sk^{(n)} X)_{n+1}} \ar{d} \ar{r}  &{X_{n+1}}  \ar{d} \ar{r} &{(\cosk^{(n)} X)_{n+1}} \ar{d} \\
{(\sk^{(n)} Y)_{n+1}} \ar{r}  &{Y_{n+1}} \ar{r} &{(\cosk^{(n)} Y)_{n+1}} 
\end{tikzcd}
commutes in \bC.\\

Let $X$ be a simplicial object in \bC.  Recall that 
$sk^{(n)} X = \sk^{(n)}$ ($\tr^{(n)}X$) 
and 
$cosk^{(n)} X = \cosk^{(n)}$ ($\tr^{(n)} X$) 
(cf. p. \pageref{13.84a}).

\label{13.90}
\indent\indent (L) \ The 
\un{latching object} 
\index{latching object} 
of $X$ at $[n]$ is $L_n X = (sk^{(n-1)} X)_n$ and the 
\un{latching} \un{morphism} 
\index{latching morphism} 
is the arrow $L_n X \ra X_n$.

\indent\indent (M) \ The 
\un{matching object} 
\index{matching object} 
of $X$ at $[n]$ is $M_n X = (cosk^{(n-1)} X)_n$ and the 
\un{matching} \un{morphism} 
\index{matching morphism}
is the arrow $X_n \ra M_n X$.

[Note: \  The 
\un{connecting morphism} 
\index{connecting morphism (simplicial object)} 
of $X$ at $[n]$ is the composite 
$L_n X \ra X_n \ra M_n X$.]

In particular: $L_0 X$ is an initial object in \bC and $M_0 X$ is a final object in \bC.\\

\begin{proposition} \ %40
Let \bC be a complete and cocomplete \mc.  Suppose that $f:X \ra Y$ is a morphism in \bSIC such that $\forall \ n$, the arrow 
$X_n \underset{L_n X}{\sqcup} L_n Y \ra Y_n$ is a cofibration (acyclic cofibration) in \bC $-$then 
$\forall \ n$, $L_n f:L_n X \ra L_n Y$ is a cofibration (acyclic cofibration) in \bC.
\end{proposition}

[One checks by induction that $L_n f$ has the LLP w.r.t. acyclic fibrations (fibrations) in \bC.]

[Note: \  There is a parallel statement for fibrations (acyclic fibrations) involving the arrows $X_n \ra M_n X \times _{M_n Y} Y_n$.]\\

\begin{proposition} \ %41
Let \bC be a complete and cocomplete \mc.  
Suppose that $f:X \ra Y$ is a morphism in \bSIC such that $\forall \ n$, the arrow 
$X_n \underset{L_n X}{\sqcup} L_n Y \ra Y_n$ 
($X_n \ra M_n X \times _{M_n Y} Y_n$) is a cofibration (fibration) in \bC 
$-$then $\forall \ n$, $f_n:X_n \ra Y_n$ is a cofibration (fibration) in \bC.
\end{proposition}

[Consider the pushout square
\begin{tikzcd}%[sep=large]
{L_n X}  \ar{d} \ar{r} &{L_n Y} \ar{d}\\
{X_n} \ar{r}  &{X_n \underset{L_n X}{\sqcup} L_n Y} 
\end{tikzcd}
Owing to Proposition 40, the arrow $L_n X \ra L_n Y$ is a cofibration.  Therefore the arrow 
$X_n \ra X_n \underset{L_n X}{\sqcup} L_n Y$ is a cofibration.  But $f_n$ is the composite 
$X_n \ra X_n \underset{L_n X}{\sqcup} L_n Y \ra Y_n$.]\\

\begin{proposition} \ %42
Let \bC be a complete and cocomplete \mc.  Suppose that $f:X \ra Y$ is a morphism in \bSIC such that $\forall \ n$, 
$f_n:X_n \ra Y_n$ is a weak equivalence 
%%----------------------------------------------------------------------------------------------55
in \bC and the arrow 
$X_n \underset{L_n X}{\sqcup} L_n Y \ra Y_n$ is a cofibration in \bC $-$then $\forall \ n$, the arrow 
$X_n \underset{L_n X}{\sqcup} L_n Y \ra Y_n$ is an acyclic cofibration in \bC.
\end{proposition}

[One checks by induction that $L_n f$ has the LLP w.r.t. fibrations in \bC.]

[Note: \  There is a parallel statement for fibrations involving arrows $X_n \ra M_n X \times _{M_n Y} Y_n$.]\\

\label{13.133}
Let \bC be a complete and cocomplete model category.  
Given a morphism $f:X \ra Y$ in \bSIC, call $f$ a \we if $\forall \ n$, $f_n:X_n \ra Y_n$ is a weak equivalence in \bC, 
a cofibration if $\forall \ n$, the arrow $X_n \underset{L_nX}{\sqcup} \L_nY \ra Y_n$ is a cofibration in \bC, 
a fibration if $\forall \ n$, the arrow $X_n \ra M_nX \times_{M_nY} Y_n$ is a fibration in \bC.  
This structure is the 
\un{Reedy structure}
\index{Reedy structure} 
on \bSIC.\\
\label{14.1}

\index{Theorem: Reedy Model Category Theorem}
\textbf{\small REEDY MODEL CATEGORY THEOREM} \ 
Let \bC be a complete and cocomplete (proper) model category $-$then \bSIC in the Reedy structure is a (proper) model category.

[The crux of the matter is the verification of MC-4 and MC-5. 
However, due to the extension principle, 
the requisite liftings and factorizations can be constructed via induction using Propositions 40, 41, and 42.]

[Note: \  Suppose further that \bC is a simplicial model category $-$then \bSIC is a simplicial model category.  
In fact, \bSIC admits a closed simplicial action derived from that on \bC 
(cf. p. \pageref{13.85}), 
so it suffices to verify that SMC holds.  
For this, it is convenient to employ Proposition 31.  
Thus let $X \ra Y$ be a cofibration in \bSIC and $L \ra K$ an inclusion of simplicial sets.  Claim:  The arrow
$X \Box K \underset{X \Box L}{\sqcup} Y \Box L \ra Y \Box K$ 
is a cofibration which is acyclic if $X \ra Y$ or $L \ra K$ is acyclic.  
Fix $n$ and consider the arrow
$(X \Box K \underset{X \Box L}{\sqcup} Y \Box L)_n 
\sqcup_{L_n(X \Box K \underset{X \Box L}{\sqcup} Y \Box L)} L_n(Y \Box K) \ra (Y \Box K)_n$ 
or, equivalently, the arrow
$(X_n \underset{L_nX}{\sqcup} L_nY) \Box K \sqcup_{(X_n \underset{L_nX}{\sqcup} L_nY)\Box L} Y_n 
\Box L \ra Y_n \Box K$.  
On the other hand, the canonical simplicial action $\Box$ on \bSIC need not be compatible with the Reedy structure on \bSIC.  Thus let $X \ra Y$ be a cofibration in \bSIC and consider the arrows 
$X \bbox \Delta[1] \underset{X \Box \Lambda[i,1]}{\sqcup} Y \bbox \Lambda[i,1] \ra Y \bbox \Delta[1]$ $(i = 0,1)$ 
(cf. p. \pageref{13.86}).  
While cofibrations, they need not be weak equivalences.]\\

\begingroup%%----------------------------------->>
\fontsize{9pt}{11pt}\selectfont
\textbf{\small EXAMPLE} \ 
Take \bC = $\bTOP_*$ (singular structure) $-$then according to 
Dwyer-Kan-Stover\footnote[2]{\textit{J. Pure Appl. Algebra} \textbf{90} (1993), 137-152; 
see also \textit{J. Pure Appl. Algebra} \textbf{103} (1995), 167-188.} 
there is a \mc structure on $\bSITOP_*$ having for its weak equivalences those $f:X \ra Y$ such that 
$\forall \ n \geq 1$, $f_*:\pi_n(X) \ra \pi_n(Y)$ is a weak equivalence of simplicial groups.  
Obviously, every weak equivalence in the Reedy structure is a weak equivalence in this structure (but not conversely).\\
\endgroup%%------------------------------------<<

%%----------------------------------------------------------------------------------------------56
The functor category $[\bDelta^\text{op},\bSISET]$ carries two other proper model category structures 
(cf. p. \pageref{13.87}).  
Every cofibration in the Reedy structure is a cofibration in structure R and every fibration in the Reedy structure is a fibration in structure L (cf. Proposition 41).  
Therefore every fibration in structure R is a fibration in the Reedy structure and every cofibration in structure L is a cofibration in the Reedy structure.

[Note: \  In reality, the cofibrations in the Reedy structure are precisely the levelwise injective simplicial maps, thus the Reedy structure is structure R.]\\

\label{14.69a}
\label{14.122a}

\begingroup%%----------------------------------->>
\fontsize{9pt}{11pt}\selectfont
$\bGamma$ is the category whose objects are the finite sets $\bn$ $\equiv \{0, 1, \ldots, n\}$ $(n \geq 0)$ with base point 0 and whose morphisms are the base point preserving maps.
\vspi
[Note: \  Suppose that $\gamma: \bm \ra \bn$ is a morphism in $\bGamma$ $-$then the partition 
$\ds\coprod\limits_{0 \leq j \leq n} \gamma^{-1}(j) = \bm$ of \bm determines a permutation 
$\theta: \bm \ra \bm$ such that $\gamma \circx \theta$ is order preserving.  
Therefore $\gamma$ has a unique factorization of the form $\alpha \circx \sigma$, where $\alpha:\bm \ra \bn$ 
is order preserving and $\sigma:\bm \ra \bm$ is a base point preserving permutation which is order preserving 
in the fibers of $\gamma$.]
\vspi
\label{14.177}
Notation: Write $\bGamma\bSISET_*$ for the full subcategory of $[\bGamma,\bSISET_*]$ whose objects are the 
$X: \bGamma \ra \bSISET_*$ such that $X_0 = *$, ($X_n = X(\bn)$).\\
\endgroup%%------------------------------------<<

\label{14.166}
\begingroup%%----------------------------------->>
\fontsize{9pt}{11pt}\selectfont
\textbf{\small EXAMPLE} \ 
Let $G$ be an abelian semigroup with unit.  
Using additive notation, view $G^n$ as the set of base point preserving functions 
$\bn \ra G$ $-$then the rule $X_n = \text{si}G^n$ defines an object in $\bGamma\bSISET_*$.  
Here the arrow $G^m \ra G^n$ attached to $\gamma: \bm \ra \bn$ sends $(g_1, \ldots g_m)$ to 
$(\ov{g}_1, \ldots, \ov{g}_n)$, where 
$\ov{g}_j = \sum\limits_{\gamma(i) = j} g_i$ if $\gamma^{-1}(j) \neq \emptyset$, 
$\ov{g}_j = 0$ if $\gamma^{-1}(j) = \emptyset$.\\
\endgroup%%------------------------------------<<

\begingroup%%----------------------------------->>
\fontsize{9pt}{11pt}\selectfont
Let \ $S_n(\bSISET_*$) \ be the category whose objects are the pointed simplicial left \ $S_n$-sets 
$-$then $S_n(\bSISET_*$) is a simplicial model category 
(cf. p. \pageref{13.88}).
\vspi
[Note: \  The group of base point preserving permutations $\bn \ra \bn$ is $S_n$ for any $X$ in 
$\bGamma\bSISET_*$, $X_n$ is a pointed simplicial left $S_n$-set.]
\vspi
Let $\bGamma_n$ be the full subcategory of $\bGamma$ whose objects are the $\bm$ $(m \leq n)$.  \ 
Assigning to the symbol \ $\bGamma_n\bSISET_*$ \ the obvious interpretation, one can follow the usual procedure and introduce tr$^{(n)}:\bGamma\bSISET_* \ra \bGamma_n\bSISET_*$ and its left (right) adjoint sk$^{(n)}$ (cosk$^{(n)}$) 
(cf. p. \pageref{13.89}).  
Put
$sk^{(n)} = \text{sk}^{(n)}\circx \text{tr}^{(n)}$
(the \un{n-skeleton}),
\index{n-skeleton simplicial model category} 
$cosk^{(n)} = \text{cosk}^{(n)}\circx \text{tr}^{(n)}$
(the \un{n-coskeleton}).
\index{n-coskeleton simplicial model category}\\
\endgroup%%------------------------------------<<

\index{Extension Principle (Objects)}
\begingroup%%----------------------------------->>
\fontsize{9pt}{11pt}\selectfont
\textbf{\small EXTENSION PRINCIPLE (OBJECTS)} \ 
Let $X$ be an object in $\bGamma_n\bSISET_*$ then a factorization 
$(\text{sk}^{(n)}X)_{n+1} \ra X_{n+1} \ra (\text{cosk}^{(n)}X)_{n+1}$ 
of the arrow
$(\text{sk}^{(n)}X)_{n+1} \ra (\text{cosk}^{(n)}X)_{n+1}$
in $S_{n+1}(\bSISET_*)$ determines an extension of $X$ to an object in $\bGamma_{n+1}\bSISET_*$.\\
\endgroup%%------------------------------------<<

\index{Extension Principle (Morphisms)}
\begingroup%%----------------------------------->>
\fontsize{9pt}{11pt}\selectfont
\textbf{\small EXTENSION PRINCIPLE (MORPHISMS)} \ 
Let
$
\begin{cases}
\ X\\
\ Y
\end{cases}
$
be objects in $\bGamma_{n+1}\bSISET_*$; let $f:\restr{X}{{\bGamma_n}} \ra \restr{Y}{{\bGamma_n}}$ be a morphism $-$then an $S_{n+1}$-equivariant arrow $X_{n+1} \ra Y_{n+1}$ determins an 
%%----------------------------------------------------------------------------------------------57
extension $F:X \ra Y$ of $f$ provided that
\begin{tikzcd}[sep=large]
{(\text{sk}^{(n)}X)_{n+1}} \ar{d} \ar{r} &{X_{n+1}} \ar{d} \ar{r} &{(\text{cosk}^{(n)}X)_{n+1}} \ar{d}\\
{(\text{sk}^{(n)}Y)_{n+1}} \ar{r}           &{Y_{n+1}} \ar{r}           &{(\text{cosk}^{(n)}Y)_{n+1}}
\end{tikzcd}
commutes in $S_{n+1}(\bSISET_*)$.\\
\endgroup%%------------------------------------<<

\begingroup%%----------------------------------->>
\fontsize{9pt}{11pt}\selectfont
Given an $X$ in $\bGamma\bSISET_*$ , write 
$L_nX = (sk^{(n-1)}X)_n$, $M_nX = (cosk^{(n-1)}X)_n$ for the latching, matching objects of $X$ at $\bn$ 
(cf. p. \pageref{13.90}).
\vspi
Given a morphism $f:X \ra Y$, call $f$ a \we if $\forall \ n \geq 1$, $f_n:X_n \ra Y_n$ is a \we in $S_n(\bSISET_*)$, a cofibration if $\forall \ n \geq 1$, the arrow $X_n \underset{L_nX}{\cup} L_nY \ra Y_n$ is a cofibration in $S_n(\bSISET_*)$, a fibration if $\forall \ n \geq 1$, the arrow $X_n \ra M_nX \times_{M_nY} Y_n$ is a fibration  in $S_n(\bSISET_*)$.  This structure is the 
\un{Reedy structure}
\index{Reedy structure ($\bGamma\bSISET_*$)} 
on  $\bGamma\bSISET_*$.\\
\endgroup%%------------------------------------<<

\index{Bousfield-Friedlander Model Category Theorem}
\begingroup%%----------------------------------->>
\fontsize{9pt}{11pt}\selectfont
\textbf{\small BOUSFIELD-FRIEDLANDER MODEL CATEGORY THEOREM} \quad
$\bGamma\bSISET_*$ in the Reedy structure is a proper simplicial model category.\\
\endgroup%%------------------------------------<<

Observation: The opposite of a \mc is a \mc 
(cf. p. \pageref{13.91}).  
So, if \bC is a complete and cocomplete \mc, then by the above $[\bDelta^\OP,\bC^\OP]$ is a \mc.  
Therefore $[\bDelta^\OP,\bC^\OP]^\OP$ is a \mc, i.e. \bCOSIC is a \mc (Reedy structure).\\

\begingroup%%----------------------------------->>
\fontsize{9pt}{11pt}\selectfont
\textbf{\small EXAMPLE} \ 
Take \bC = \bSISET $-$then the class of \wes in $[\bDelta,\bSISET]$ (Reedy Structure) is the same as the class of \wes in $[\bDelta,\bSISET]$ (structure \bL 
(cf. p. \pageref{13.92})) 
but the class of cofibrations is larger.  
Example: $Y_{\bDelta} \equiv \bDelta$ 
(cf. p. \pageref{13.93}) 
is a cosimplicial object in $\widehat{\bDelta}$ which is cofibrant in the Reedy structure but not in structure \bL.\\
\endgroup%%------------------------------------<<

\begin{proposition} \ %43
Let \bC be a complete and cocomplete model category. Equip \textbf{SIC} with its Reedy structure $-$then the functor $L_n:\textbf{SIC} \ra \bC$ preserves weak equivalences between cofibrant objects.
\end{proposition}

[Inspect the proof of Proposition 42 and quote the lemma on 
p. \pageref{13.94}.]\\

Let \bC be a simplicial model category.  Assume: \bC is complete and cocomplete.  
Given an $X$ in \textbf{SIC}, put \ 
${\displaystyle\abs{X} =  \int^{[n]} X_n \Box \Delta [n]}$ 
$-$then $\abs{X}$ is the 
\un{realization} 
\index{realization (in \textbf{SIC}} 
of $X$ and the assignment 
$X \ra \abs{X}$ is a functor $\textbf{SIC} \ra \bC$.  \  
$\abs{?}$ is a left adjoint for $\sin:\bC \ra \textbf{SIC}$, where 
$\sin_nY =\texttt{HOM}(\Delta[n],Y)$.  
 In fact, 
${\ds\Mor(\abs{X},Y) \approx}$ 
${\ds\Mor(\int^{[n]} X_n \Box \Delta [n],Y) \approx}$ 
${\ds \int_{[n]} \text{Mor}(X_n \Box \Delta[n],Y) \approx}$
${\ds \int_{[n]} \text{Mor}(X_n,\texttt{HOM}(\Delta[n],Y)) \approx}$
${\ds \int_{[n]} \text{Mor}(X_n, \sin_n Y) \approx}$ 
${\Nat(X,\sin Y)}$.\\
\vspace{0.25cm}

%%----------------------------------------------------------------------------------------------58

\begingroup%%----------------------------------->>
\fontsize{9pt}{11pt}\selectfont
\textbf{\small EXAMPLE} \ 
Take \bC = \bSISET and let $X$ be a simplicial object in \bC.  \ 
One can fix [m] and form $\abs{X_m^h}$, the geometric realization of $[n] \ra X([n],[m])$ and one can fix [n] and form $\abs{X_n^v}$, the geometric realization of $[m] \ra  X([n],[m])$.  The assignments
$
\begin{cases}
\ [m] \ra \abs{X_m^h}\\
\ [n] \ra \abs{X_n^v}
\end{cases}
$
define simplicial objects
$
\begin{cases}
\ X^h\\
\ X^v
\end{cases}
$
in \bCGH  and their realizations 
$
\begin{cases}
\  \abs{X^h}\\
\ \abs{X^v}
\end{cases}
$
are homeomorphic to the geometric realization of $\abs{X}$.\\
\endgroup%%------------------------------------<<
\vspace{0.25cm}

\textbf{\small LEMMA} \ 
Let $X$ be a simplicial object in \bC $-$then $\abs{X} \approx \colim\abs{X}_n$, where
%$\abs{X}_n = \displaystyle \int^{[k]} X_k \Box \Delta[k]^{(n)}$\\
%$\abs{X}_n = {\displaystyle \int}^{[k]} X_k \Box \Delta[k]^{(n)}$\\
$\abs{X}_n = {\ds\int}^{[k]} X_k \ \Box \  \Delta[k]^{(n)}$.  Moreover, $\forall \ n > 0$ there is a pushout square
\[
\begin{tikzcd}%[sep=huge]
{L_nX \ \Box \ \Delta[n] \underset{L_nX \ \Box \ \dot\Delta[n]}{\sqcup} X_n  \ \Box \ \dot\Delta[n]}
 \ar{d} \ar{r} &{\abs{X}_{n-1}} \ar{d}\\
{X_n \ \Box \ \Delta[n]} \ar{r} &{\abs{X}_n}
\end{tikzcd}
.
\]

[The functors $X_n \ \Box \  -$ are left adjoints, hence preserve colimits, so
$\abs{X} = {\ds\int}^{[n]} X_n \  \Box \  \Delta[n]$  
$\approx$ $\ds\int^{[n]} X_n \  \Box \  \colim_k \Delta[n]^{(k)}$ 
$\approx$ $\ds\int^{[n]} \colim_k X_n \bbox \Delta[n]^{(k)}$ 
$\approx \colim_n \ds\int^{[k]} X_k  \Box \  \Delta[k]^{(n)}$  
$= \colim_n \abs{X}_n$.  
And: Relative to the inclusion $\bDelta_n \ra \bDelta$, the left Kan extension of $[m] \ra \Delta[m]$ $(m \leq n)$ is 
$[k] \ra \Delta[k]^{(n)}$, thus $\abs{X}_n$ can be identified with 
$\ds\int^{[m]} X_m \  \Box \  \Delta[m]$ $(m \leq n)$.]\\

If $X$ is a cofibrant object in \bSIC (Reedy structure), then the latching morphism $L_nX \ra X_n$ is a cofibration in \bC.  
Therefore the arrow 
$L_n X \ \Box \ \Delta[n] \underset{L_n X \ \Box \ \dot\Delta[n]}{\sqcup} X_n \ \Box \ \dot\Delta[n] \ra X_n \ \Box \ \Delta[n]$
is a cofibration in \bC (cf. Proposition 31).  
Consequently, the arrow $\abs{X}_{n-1} \ra \abs{X}_n$ is a cofibration in \bC.

[Note: \  It follows from Proposition 40 that $L_nX$ is a cofibrant object in \bC, hence $X_n$ is a cofibrant object in \bC.  This means that $L_n X \ \Box \ \dot\Delta[n]$, $L_n X \ \Box \ \Delta[n]$,  and $X_n \ \Box \ \dot\Delta[n]$  are cofibrant objects in \bC, so
$L_n X \ \Box \ \Delta[n] \underset{L_n X \ \Box \ \dot\Delta[n]}{\sqcup} X \ \Box \ \dot\Delta[n]$ is a cofibrant object in \bC 
(cf. p. \pageref{13.95}).]\\

\textbf{\small LEMMA} \ 
Let \bC be a simplicial model category. Assume: \bC is complete and cocomplete.  Suppose that
$
\begin{cases}
\ X\\
\ Y
\end{cases}
$
are cofibrant objects in \bSIC (Reedy structure) and $f:X \ra Y$ is a \we $-$then the arrow 
\[
L_n X \ \Box \ \Delta[n] \underset{L_n X \ \Box \ \dot\Delta[n]}{\sqcup} X_n \ \Box \ \dot\Delta[n] \ra 
L_n Y \ \Box \ \Delta[n] \underset{L_n Y \ \Box \ \dot\Delta[n]}{\sqcup} Y_n \ \Box \ \dot\Delta[n]
\]
is a \we in \bC.

%%----------------------------------------------------------------------------------------------59
[Consider the commutative diagram
\[
\begin{tikzcd}%[sep=large]
{L_n X \ \Box \ \Delta[n]} \ar{d}  &{L_n X \ \Box \ \dot\Delta[n]} \ar{l} \ar{r} \ar{d} &{X_n \ \Box \ \dot\Delta[n]} \ar{d}\\
{L_n Y \ \Box \ \Delta[n]}            &{L_n Y \ \Box \ \dot\Delta[n]} \ar{l} \ar{r}           &{Y_n \ \Box \ \dot\Delta[n]}
\end{tikzcd}
.  
\]
The horizontal arrows are cofibrations 
(cf. p. \pageref{13.95a}) 
and the vertical arrows are weak equivalences (cf. Proposition 43 and 
p. \pageref{13.95b}).  
Therefore Proposition 3 in $\S 12$ is applicable.]\\

\label{13.118}

\begin{proposition} \ %44
Let \bC be a simplicial model category.  Assume: \bC is complete and cocomplete.  Suppose that
$
\begin{cases}
\ X\\
\ Y
\end{cases}
$
are cofibrant objects in \bSIC (Reedy structure) and $f:X \ra Y$ is a \we $-$then $\abs{f}:\abs{X} \ra \abs{Y}$ is a \we.
\end{proposition}

[Since 
$
\begin{cases}
\ \abs{X}_0 = X_0\\
\ \abs{Y}_0 = Y_0
\end{cases}
$
and $\forall \ n$, 
$
\begin{cases}
\ \abs{X}_n \ra \abs{X}_{n+1}\\
\ \abs{Y}_n \ra \abs{Y}_{n+1}
\end{cases}
$
is a cofibration in \bC, one may view 
$
\begin{cases}
\ \{\abs{X}_n: n \geq 0\}\\
\ \{\abs{Y}_n: n \geq 0\}
\end{cases}
$
as cofibrant objects in $\bFIL(\bC)$ 
(cf. p. \pageref{13.96}).  
So, to prove that $\abs{f}:\abs{X} \ra \abs{Y}$ is a \we , it need only be shown that $\forall \ n$, $\abs{f}_n:\abs{X}_n \ra \abs{Y}_n$ is a \we 
(cf. p. \pageref{13.97}).  
For this, work with 
\[
\begin{tikzcd}%[sep=large]
&&{X_n \ \Box \ \Delta[n]} \ar{d}  
&{L_n X \ \Box \ \Delta[n] \underset{L_n X \ \Box \ \dot\Delta[n]}{\sqcup} X_n \ \Box \ \dot\Delta[n]}  \ar{l} \ar{r} \ar{d} 
&{\abs{X}_{n-1}} \ar{d}\\
&&{Y_n \ \Box \ \Delta[n]}            
&{L_n Y \ \Box \ \Delta[n] \underset{L_n Y \ \Box \ \dot\Delta[n]}{\sqcup} Y_n \ \Box \ \dot\Delta[n]}  \ar{l} \ar{r}            
&{\abs{Y}_{n-1}}
\end{tikzcd}\\
\]
and use induction (cf. $\S 12$, Proposition 3).]\\

\begingroup%%----------------------------------->>
\fontsize{9pt}{11pt}\selectfont
\textbf{\small EXAMPLE} \ 
Take \bC = \bSISET and suppose that $f:X \ra Y$ is a weak equivalence, i.e., 
$\forall \ n$, $f_n:X_n \ra Y_n$ is a weak equivalence $-$then 
$\abs{f}:\abs{X} \ra \abs{Y}$ is a weak homotopy equivalence.
\vspi
[All simplicial objects in $\widehat{\bDelta}$ are cofibrant in the Reedy structure.]
\vspi
[Note: \  Fix an abelian group $G$ and consider \bSISET in the homological model category structure determined by $G$ $-$then \bSISET is a simplicial model category 
(cf. p. \pageref{13.98}), 
hence $\abs{f}:\abs{X} \ra \abs{Y}$ is an  $HG$-equivalence if $\forall \ n$, 
$f_n:X_n \ra Y_n$ is an  $HG$-equivalence.]\\
\endgroup%%------------------------------------<<

\label{13.117} %dmc mnft
\begingroup%%----------------------------------->>
\fontsize{9pt}{11pt}\selectfont
\textbf{\small EXAMPLE} \ 
Suppose that \bC is a \smc which is complete and cocomplete.  Let $X$ be a cofibrant object in \bSIC (Reedy structure).  Assume: $\forall \ \alpha$, $X\alpha$ is a \we $-$then the arrow $\abs{X}_0 \ra \abs{X}$ is a \we.\\
\endgroup%%------------------------------------<<

Let \bC be a simplicial model category.  Assume: \bC is complete and cocomplete.  
Given an $X$ in \bCOSIC, put  
$\tot X = \ds\int_{[n]} \shom(\Delta[n],X_n)$ $-$then $\tot X$ is the 
\un{totalization}
\index{totalization}
%%----------------------------------------------------------------------------------------------60
of $X$ and the assignment $X \ra \text{tot} X$ is a functor $\bCOSIC \ra \bC$.  tot is a right adjoint for 
cosin$:\bC \ra \bCOSIC$, where cosin$_nY = Y_n \ \Box \ \Delta[n]$.  In fact, 
$\Mor(Y,\tot X) \approx$ 
$\Mor(Y,\ds\int_{[n]} \shom(\Delta[n],X_n))$ $\approx$
$\ds\int_{[n]}\Mor(Y,\shom(\Delta[n],X_n))$ $\approx$ \ 
$\ds\int_{[n]}\Mor(Y \ \Box \ \Delta[n],X_n)$ $\approx$ \ 
$\ds\int_{[n]}\Mor(\text{cosin}_n Y,X_n)$ $\approx$ \ 
$\Nat(\text{cosin } Y,X)$.

Example: Take \bC = \bSISET $-$then $\tot X = \HOM(Y_{\bDelta},X)$ 
(cf. p. \pageref{13.99}).

Example: Let $X$ be a simplicial set.  
Given a cosimplicial object $Y$ in $\widehat{\bDelta}$, the functor 
$\bDelta \ra \bSISET$ that sends [n] to $\map(X,Y_n)$ defines another  cosimplicial object in $\widehat{\bDelta}$, call it 
$\map(X,Y)$.  And: 
$\text{tot } \map(X,Y) \approx$
$\ds\int_{[n]}\map(\Delta[n],\map(X,Y_n))$ $\approx$
$\ds\int_{[n]}\map(X,\map(\Delta[n],$ $Y_n))$ $\approx$
$\map(X, \ds\int_{[n]}\map(\Delta[n],Y_n))$ $\approx$
$\map(X,\tot Y).$\\
\vspace{0.25cm}

\begingroup%%----------------------------------->>
\fontsize{9pt}{11pt}\selectfont
\textbf{\small EXAMPLE} \ 
Given a simplicial set $K$ and a compactly  generated Hausdorff space $X$, let $X^K$ be the cosimplicial object in \bCGH with 
$(X^K)_n = X^{K_n}$ $-$then $\map(\abs{K},X) \approx \tot X^K$.\\
\endgroup%%------------------------------------<<

\begingroup%%----------------------------------->>
\fontsize{9pt}{11pt}\selectfont
\textbf{\small EXAMPLE} \ 
Fix a prime $p$ $-$then there is a forgetful functor from the category of simplicial vector spaces over $\F_p$ to \bSISET.  
It has a left adjoint, thus this data determines a triple in \bSISET.  
Write res$_pX$ for the standard resolution of \mX: res$_pX$ is therefore a cosimplicial object in $\widehat{\bDelta}$ 
and $\tot \res_pX$ is the 
\un{$\F_p$-completion} 
\index{completion! $\F_p$-completion} 
$\F_pX$ of $X$ 
(Bousfield-Kan\footnote[2]{\textit{SLN} \textbf{304} (1972).}).\\
\endgroup%%------------------------------------<<

\begin{proposition} \ %45
Let \bC be a simplicial model category.  Assume: \bC is complete and cocomplete.  Suppose that
$
\begin{cases}
\ X\\
\ Y
\end{cases}
$
are fibrant objects in \textbf{COSIC} (Reedy Structure) and $f:X \ra Y$ is a weak equivalence $-$then 
$\text{tot} f: \tot X \ra \tot Y$ is a weak equivalence.
\end{proposition}

[The proof is dual to that of Proposition 44.  Of course, $\tot X$ $\approx$ $\lim \text{tot}_nX$ (obvious notation).]\\

The simplex category $\gro_{\bDelta} K$ of a simplicial set $K$ can be viewed as a comma category:
\begin{tikzcd}[sep=small]
{\Delta [n]} \ar{ddr} \ar{rr} &&{\Delta [m]} \ar{ddl}\\
\\
&{K}
\end{tikzcd}
(cf. p. \pageref{13.100}).  
Call this interpretation of $\bDelta K$, $\bDelta^\OP K$ being its opposite.  There is a forgetful functor
$\Delta K: \textbf{$\bDelta K$} \ra \bSISET$ and $K \approx \colimx \Delta K$ 
(cf. p. \pageref{13.101}).\\

\begingroup%%----------------------------------->>
\fontsize{9pt}{11pt}\selectfont
\textbf{\small FACT} \ 
The fundamental groupoid of $\bDelta K$ is equivalent to the fundamental groupoid of \mK.\\
\endgroup%%------------------------------------<<

%%----------------------------------------------------------------------------------------------61
Given a category \bC, write $K$-\bSIC 
\index{$K$-\bSIC} 
for the functor category $[\bDelta^\OP K,\bC]$ and $K$-\bCOSIC  
\index{$K$-\bCOSIC} 
for the functor category $[\bDelta K,\bC]$ $-$then by definition, a 
\un{$K$-simplicial object}
\index{K-simplicial object} 
in \bC is an object in $K$-\bSIC and a 
\un{$K$-cosimplicial object}
\index{K-cosimplicial object} 
in \bC is an object in $K$-\bCOSIC.

[Note: \  Take $K = \Delta[0]$ to recover \bSIC and \bCOSIC.]

The preceding results can now be generalized.  Thus if \bC is a complete and cocomplete \mc, one can again introduce latching objects and matching objects and use them to equip $K$-\bSIC (dually, $K$-\bCOSIC) with the structure of a \mc (Reedy structure).  Assuming in addition that \bC is a simplicial \mc, there is a realization functor 
$\abs{?}_K:K\text{-}\bSIC \ra \bC$ that sends $X$ to 
$\abs{X}_K = \ds\int^{\bDelta K} X \ \Box \ \Delta K$, 
where 
$ X \ \Box \ \Delta K: \bDelta^\OP K \times \bDelta K \ra \bC$ 
is the composite
%$\bDelta^\OP K \times \bDelta K
\begin{tikzcd}[sep=large]
{\bDelta^\OP K \times \bDelta K}  \ar{r}{X \times \Delta K} &{\bC \times \bSISET}
\end{tikzcd}
$\overset{\Box}{\lra} \bC$.
So in the notation of the Kan extension theorem, 
$\abs{?}_K = \abs{?} \circx \lan$, 
i.e., the diagram
\begin{tikzcd}[sep=large]
{K\text{-}\bSIC} \ar{dr}[swap]{\abs{?}_K} \ar{r}{\lan} &{\bSIC} \ar{d}{\abs{?}}\\
&{\bC}
\end{tikzcd}
commutes.  
Here, lan is computed from the arrow 
$\bDelta^\OP K \ra \bDelta^\OP$
induced by the projection $K \ra \Delta[0]$.  $\abs{?}_K$ is a left adjoint for $\sin_K:\bC \ra K\text{-}\bSIC$.  On the other hand, there is a totalization functor 
$\tot_K:K\text{-}\bCOSIC \ra \bC$ that sends $X$ to 
$\tot_K X = \ds\int_{\bDelta K} \texttt{HOM}(\Delta K,X)$, 
where 
$\texttt{HOM}(\Delta K,X): \bDelta^\OP K \times \bDelta K \ra \bC$ 
is the composite
$\bDelta^\OP K \times \bDelta K$ 
$\overset{\bDelta^\OP K \times X}{\xrightarrow{\hspace*{1.5cm}}}$ 
$\bSISET^\OP \times \bC$ 
$\overset{\texttt{HOM}}{\xrightarrow{\hspace*{0.75cm}}}$ 
$\bC$.
So, in the notation of the Kan extension theorem, $\tot_K = \tot \circx \ran$, i.e., the diagram 
\begin{tikzcd}%[sep=large]
{K\text{-}\bCOSIC} \ar{dr}[swap]{\tot_K} \ar{r}{\ran} &{\bCOSIC} \ar{d}{\tot}\\
&{\bC}
\end{tikzcd}
commutes.  
Here, ran is computed from the arrow 
$\bDelta K \ra \bDelta$ induced by the projection 
$K \ra \Delta[0]$.  
$\tot_K$ is a right adjoint for 
$\cosin_K: \bC \ra K\text{-}\bCOSIC$.\\
\vspace{0.25cm}

\begingroup%%----------------------------------->>
\fontsize{9pt}{11pt}\selectfont
To check the claimed factorization of \ $\abs{?}_K$, \ represent \ $\abs{X}_K$ \ as the coequalizer of the diagram 
$\ds\coprod\limits_{k \ra l} X_l \xbox \Delta Kk \rightrightarrows \ds\coprod\limits_k X_k \xbox \Delta Kk$.  \ 
Noting \ that \ \ 
$(\lan X)_n = \ds\coprod\limits_{k \in K_n} X_k$, \ \ 
we have \ 
$\ds\coprod\limits_{k \ra l} X_l \ \xbox \ \Delta K k \ \approx$  
$\ds\coprod\limits_{n,m \geq 0}$ $\ds\coprod\limits_{[n] \ra [m]} \ds\coprod\limits_{l \in K_m}  X_l \ \xbox \ \Delta[n] \approx$ 
$\ds\coprod\limits_{n,m \geq 0} \ds\coprod\limits_{[n] \ra [m]} (\lan X)_m \ \xbox \ \Delta[n]$ 
and 
$\ds\coprod\limits_k X_k \ \xbox \ \Delta K k \approx $ 
$\ds\coprod\limits_{n \geq 0} \ds\coprod\limits_{k \in K_n} X_k \ \xbox \ \Delta[n]$ $\approx$ 
$\ds\coprod\limits_{n \geq 0}  (\lan X)_n \ \xbox \ \Delta[n]$, 
\ \ 
i.e., $\abs{X}_K$ \ is \ naturally \ isomorphic \  to \ the coequalizer \ of \ the \ diagram 
$\ds\coprod\limits_{n,m \geq 0} \coprod\limits_{[n] \ra [m]}$ $(\lan X)_m \ \xbox \ \Delta[n] \rightrightarrows 
\ds\coprod\limits_{n \geq 0} (\lan X)_n \ \xbox \ \Delta[n]$,
i.e., to $\abs{\lan X}$.
\vspi
Example: Take \bC = \bSISET $-$then $\abs{*}_K = K$.\\
\endgroup%%------------------------------------<<

\label{13.113}
\begingroup%%----------------------------------->>
\fontsize{9pt}{11pt}\selectfont
\textbf{\small EXAMPLE} \ 
Let $B$ be a simplicial set.  Fix an $X$ in \bSISET/\mB $-$then $\forall \ n$ $\&$ $\forall \ b \in B_n$, there is a
%%----------------------------------------------------------------------------------------------62
pullback square
\begin{tikzcd}[sep=large]
{X_b} \ar{d} \ar{r} &{X} \ar{d}{p}\\
{\Delta[n]} \ar{r}[swap]{\Delta_b} &{B}
\end{tikzcd}
(cf. p. \pageref{13.102}).  
This data thus determines a \mB-cosimplicial object $X_B$ in \bSISET.  
One has $X \approx \colim X_B$ and $X_B$ cofibrant in the Reedy structure.\\
\endgroup%%------------------------------------<<

%%%%%%%%%%%%%%%%%%%\begin{proposition} \ 
\begingroup%%----------------------------------->>
\fontsize{11pt}{11pt}\selectfont
\textbf{\small PROPOSITION 44 \  (K)} \ 
Let \bC be a simplicial model category.  Assume: \bC is complete and cocomplete.  Suppose that 
$
\begin{cases}
\ X\\
\ Y
\end{cases}
$
are cofibrant objects in $K$-\bSIC (Reedy structure) and $f:X \ra Y$ is a \we $-$then
$\abs{f}_K: \abs{X}_K \ra \abs{Y}_K$ is a weak equivalence.\\
\endgroup%%------------------------------------<<

\begingroup%%----------------------------------->>
\fontsize{11pt}{11pt}\selectfont
\textbf{\small PROPOSITION 45 \  (K)} \ 
Let \bC be a simplicial model category.  Assume: \bC is complete and cocomplete.  Suppose that 
$
\begin{cases}
\ X\\
\ Y
\end{cases}
$
are fibrant objects in $K$-\bCOSIC (Reedy structure) and $f:X \ra Y$ is a \we $-$then
$\tot_K f: \tot_K X \ra \tot_K Y$ is a weak equivalence.\\
\endgroup%%------------------------------------<<

\begingroup%%----------------------------------->>
\fontsize{9pt}{11pt}\selectfont
\textbf{\small FACT} \ 
$\sin_K$ preserves fibrations and acyclic fibrations.
\vspi
[Note: \  Therefore $\abs{?}_K$ preserves cofibrations and acyclic cofibrations.  
(cf. p. \pageref{13.103} ff.).]\\
\endgroup%%------------------------------------<<

\label{108a}
\begingroup%%----------------------------------->>
\fontsize{9pt}{11pt}\selectfont
\textbf{\small FACT} \ 
$\cosin_K$ preserves cofibrations and acyclic cofibrations.
\vspi
[Note: \  Therefore $\tot_K$ preserves fibrations and acyclic fibrations.  
(cf. p. \pageref{13.104} ff.).]\\
\endgroup%%------------------------------------<<

Notation: Let \bI be a small category.  Put $\bDelta \bI = \bDelta\ner\bI$ and call it the 
\un{simplex category}
\index{simplex category} 
of \bI $-$then $\bDelta \bI$  
is isomorphic to the comma category $\abs{\iota,K_\bI}$:
\begin{tikzcd}[sep=small]
{[n]} \ar{rdd}   \ar{rr}  &&{[m]} \ar{ldd}\\
\\
&{\bI}
\end{tikzcd}
$(\iota: \bDelta \ra \textbf{CAT})$.
There is a projection $\pi_\bI: \bDelta\bI \ra \bI$ that sends an object $[n] \overset{f}{\ra} \bI$ to 
$fn \in \Ob\bI$.  Example: $\bDelta\textbf{1} = \bDelta$.

[Note: \  $\bDelta^\OP\bI$ is the opposite of $\bDelta\bI$.  
Example: $\bDelta^\OP\textbf{1} = \bDelta^\OP$.  
Replacing \bI by $\bI^\OP$, there is a projection 
$\pi_\bI^\OP: \bDelta^\OP\bI^\OP \ra \bI$ that sends an object  $[n] \overset{f}{\ra} \bI^\OP$ 
to $fn \in \Ob \bI$.]\\

\label{13.109}

\begingroup%%----------------------------------->>
\fontsize{9pt}{11pt}\selectfont
\textbf{\small EXAMPLE} \ 
Let \bC be a complete and cocomplete \mc.  Suppose that $F:\bI \ra \bC$ is a functor such that $\forall \ i$, $F i$ is cofibrant (fibrant) $-$then 
$F \circx \pi_{\bI}^\OP$ ($F \circx \pi_{\bI}$) is a cofibrant (fibrant) object in 
$[\bDelta^\OP \bI^\OP,\bC]$ ($[\bDelta \bI,\bC]$) (Reedy structure).\\
\endgroup%%------------------------------------<<

Let \bI be a small category and \bC a simplicial model category.  
Assume: \bC is complete and cocomplete $-$then the functor
$\colim:[\bI,\bC] \ra \bC$ ($\lim:[\bI,\bC] \ra \bC$) need not preserve levelwise weak equivalences between levelwise cofibrant (fibrant) objects.  
To rememdy 
%%----------------------------------------------------------------------------------------------63
this defect, one introduces the notion of 
\un{homotopy colimit}
\index{homotopy colimit} 
(
\un{limit})
\index{homotopy limit}.  
Thus define a functor 
$\hocolim_{\bI}:[\bI,\bC] \ra \bC$
by 
$\hocolim_{\bI} F$ (or $\hocolimx F$) $= \ds\int^{\bI^\OP} F \xbox \ner(-\bs\bI)^\OP$
and define a functor 
$\holim_{\bI}:[\bI,\bC] \ra \bC$
by 
$\holim_{\bI} F$ (or $\holimx F$) $= \ds\int_{\bI} \texttt{HOM}(\ner(\bI/-),F)$.

[Note: \   One has \ 
$\HOM(\hocolim_{\bI} F,Y)$ \ $\approx $ \ 
$\HOM ( \ds\int^i F i \xbox \ner(i\bs\bI)^\OP,Y)$ \ $\approx $  
$\ds\int_i \HOM(F i$ $\xbox \ner(i\bs\bI)^\OP,Y)$  $\approx $ 
$\ds\int_i \map(\ner(i\bs\bI)^\OP,\HOM(F i,Y))$  $\approx $ 
$\ds\int_i \map(\ner(\bI^\OP/i),$ $\HOM(F i,Y))$ $\approx $
$\holim_{\bI^\OP} \HOM(F,Y)$
where 
$\HOM(F,Y):\bI^\OP \ra \bSISET$ sends $i$ to $\HOM(Fi,Y)$.]

Remark: The functor $\hocolim$ has a right adjoint, viz. 
$\texttt{HOM}(\ner(-\bs\bI)^\OP,-)$, 
and the functor $\holim$ has a left adjoint, viz. 
$- \xbox \ner(\bI/-)$.

Remark:  There are natural transformations $\hocolim \ra \colim$, and $\lim \ra \holim$.

[Note: \  It can be shown that $\bL\hocolim$, and $\bR\holim$ exist and that there are natural isomorphisms 
$\bL\hocolim \ra \bL\colim$,
$\bR\lim \ra \bR\holim$
(Dwyer-Kan\footnote[2]{\textit{Model Categories and General Abstract Homotopy Theory}, Preprint.}) 
%% maybe should be this
%(Dwyer-Hirschhorn-Kan-Smith\footnote[2]{\textit{Homotopy limit functors on model categories, and homotopical categories, %Mathematical Surveys and Monographs}, Amer. Math Soc. \textbf{113} (2004).})
(cf. p. \pageref{13.105}).]

Example: Take \bC = \bSISET, \bCGH, $\bSISET_*$, $\bCGH_*$ $-$then 
$F i \xbox \ner(i\bs\bI)^\OP = F i \times \ner(i\bs\bI)^\OP$, 
$F i \times_k  B(i\bs\bI)^\OP$, 
$F i \# \ner(i\bs\bI)_+^\OP$, 
$F i \#_k \ B(i\bs\bI)_+^\OP$, 
and 
$\texttt{HOM}(\ner(\bI/i),F i) = \map(\ner(\bI/i),F i)$, 
$\map(B(\bI/i),F i)$, 
$\map_*(\ner(\bI/i)_+,F i)$, 
$\map_*(B(\bI/i)_+,F i)$.

[Note: \  Consider $\ds\int^i F i \xbox \ner(i\bs\bI)$ and $\ds\int^i F i \xbox \ner(i\bs\bI)^\OP$.  
When \bC = \bSISET or $\bSISET_*$, they are simplicial opposites of one another 
(cf. p. \pageref{13.106}), 
hence are naturally weakly equivalent, and when \bC = \bCGH or $\bCGH_*$, they are related by a natural homeomorphism (since $\forall \ i$, 
$B(i\bs\bI) \approx B(i\bs\bI)^\OP$ 
(cf. p. \pageref{13.107})).]\\

\begingroup%%----------------------------------->>
\fontsize{9pt}{11pt}\selectfont
Place on $[\bI^\OP,\bSISET]$ and [\bI,\bSISET] structure L 
(cf. p. \pageref{13.107a}) 
$-$then $i \ra \ner(i\bs\bI)^\OP$ is a cofibrant object in $[\bI^\OP,\bSISET]$ and $i \ra \ner(\bI/i)$ is a cofibrant object in [\bI,\bSISET] 
(cf. p. \pageref{13.107b}).  
Observe too that $\forall \ i \in \Ob \bI$, the classifying spaces $B(i\bs\bI)^\OP$ and $B(\bI/i)$ are contractible 
(cf. p. \pageref{13.108}).\\
\endgroup%%------------------------------------<<

\begingroup%%----------------------------------->>
\fontsize{9pt}{11pt}\selectfont
\textbf{\small EXAMPLE} \ 
Let $F$ be the functor $\bI \ra \bSISET$ that sends $i \in \Ob \bI$ to $F i = \Delta[0]$ $-$then 
$\hocolim F \approx \ner \bI^\OP$, 
i.e., 
$\ds\int^i \Delta[0] \times \ner(i\bs\bI)^\OP \approx \ner \bI^\OP$ 
or still, 
$\ds\int^i \Delta[0] \times \ner(\bI^\OP/i) \approx \ner \bI^\OP$.  
Similarly, 
$\ds\int^i \Delta[0] \times \ner(i\bs\bI) \approx \ner \bI$
and 
$\ds\int^i \Delta[0] \times \ner(i\bs\bI^\OP) \approx \ner \bI^\OP$.  
In addition, 
$\ds\int^i \Delta[0] \times \ner(i\bs\bI^\OP)^\OP \approx \ner \bI$ 
or still, 
$\ds\int^i \Delta[0] \times \ner(\bI/i) \approx \ner \bI$.\\
\endgroup%%------------------------------------<<

\begingroup%%----------------------------------->>
\fontsize{9pt}{11pt}\selectfont
\textbf{\small EXAMPLE} \ 
Let $U:\bCGH_* \ra \bCGH$ be the forgetful functor and consider a functor $F:\bI \ra \bCGH_*$.  
Question: What is the relation between $\hocolimx F$ $\&$ $\hocolimx U \circx F$ and $\holimx F$ $\&$ $\holimx U \circx F$?  
The 
%%----------------------------------------------------------------------------------------------64
answer for homotopy limits is that there is essentially no difference (since $\map_*(X_+,Y) \approx \map(X,UY)$).  
Turning to homotopy colimits, assume that $\forall \ i$, $F i$ is cofibrant $-$then there is a cofibration 
$B\bI^\OP \ra \hocolimx U \circx F$ 
and a homeomorphism 
$\hocolimx U \circx F/B\bI^\OP \ra \hocolimx F$.
\vspi
[Note: \  If $B\bI^\OP$ is contractible, the projection 
$\hocolimx U \circx F \ra \hocolimx F$ is a weak homotopy equivalence.  Proof: Consider the pushout square 
\begin{tikzcd}[sep=large]
{B\bI^\OP} \ar{d}   \ar{r}  &{*} \ar{d}\\
{\hocolimx U \circx F} \ar{r}&{\hocolimx F}
\end{tikzcd}
,
bearing in mind that \bCGH (singular structure) is a proper model category.]\\
\endgroup%%------------------------------------<<

\label{13.130}
\label{18.16}
\begingroup%%----------------------------------->>
\fontsize{9pt}{11pt}\selectfont
\textbf{\small LEMMA} \ 
Let $X \ra B$ be a simplicial map.  \ Suppose that for every commutative diagram \qquad
$
\begin{tikzcd}[sep=large]
{X_{b^\prime}} \ar{d}   \ar{r}  &{X_b} \ar{d} \ar{r} &{X} \ar{d}{p}\\
{\Delta[n^\prime]} \ar{r}&{\Delta[n]} \ar{r}[swap]{\Delta_b} &{B}
\end{tikzcd}
,
$
the arrow $X_{b^\prime} \ra X_b$ is a weak homotopy equivalence $-$then $p$ is a homotopy fibration.\\
\endgroup%%------------------------------------<<

\label{13.110a}
\label{13.114}
\label{13.111}  %dmc mnft
\label{13.124}
\begingroup%%----------------------------------->>
\fontsize{9pt}{11pt}\selectfont
\textbf{\small FACT} \ 
Let $F:\bI \ra \bSISET$ be a functor $-$then the arrow $\hocolimx F \ra \ner \bI^\OP$ is a homotopy fibration iff 
$\forall \ \delta \in \Mor \bI$, $F \delta$ is a weak homotopy equivalence.\\
\endgroup%%------------------------------------<<

\begin{proposition} \ %46
Fix $F \in \Ob[\bI,\bC]$ $-$then
\[
\hocolimx F \approx \int^{\bDelta \bI^\OP} F \circx \pi_{\bI}^\OP \  \xbox \ \ner \bI^\OP \hspace{0.5cm} (= \abs{F \circx \pi_{\bI}^\OP }_{\ner \bI^\OP})
\]
and 
\[
\holimx F \approx \int_{\bDelta \bI} \texttt{HOM}(\Delta \ner \bI,F \circx \pi_{\bI}) \hspace{1.5cm} (= \tot_{\ner\bI} F \circx \pi_{\bI}). \hspace{0.25cm}
\]
\\
\end{proposition}

Application:  Let $F, \ G:\bI \ra \bC$ be functors and let $\Xi:F \ra G$ be a natural transformation.  
Assume: $\forall \ i$, $\Xi_i:F i \ra G i$ is a \we $-$then 
$\hocolimx \Xi: \hocolimx F \ra \hocolimx G$ is a \we provided that $\forall \ i$, 
$
\begin{cases}
\ F i\\
\ G i
\end{cases}
$
is cofibrant and 
$\holimx \Xi:\holimx F \ra \holimx G$ is a \we  provided that $\forall \ i$, 
$
\begin{cases}
\ F i\\
\ G i
\end{cases}
$
is fibrant.

[In view of the above result and the example on 
p. \pageref{13.109}, 
this follows from Propositions 44(K) and 45(K).]\\

\label{14.31}
\begingroup%%----------------------------------->>
\fontsize{9pt}{11pt}\selectfont
\textbf{\small EXAMPLE} \ 
Let $F:\bI \ra \bSISET$ be a functor $-$then there is a natural homeomorphism $\abs{\hocolimx F} \ra \hocolimx \abs{F}$ of compactly generated Hausdorff spaces.
\vspi
[Geometric realization is a left adjoint, hence preserves colimits.]\\
\endgroup%%------------------------------------<<

%%----------------------------------------------------------------------------------------------65
\begingroup%%----------------------------------->>
\fontsize{9pt}{11pt}\selectfont
\textbf{\small EXAMPLE} \ 
Let $F: \bI \ra \bCGH$ be a functor such that $\forall \ i$, $Fi$ is cofibrant $-$then there is a natural weak homotopy equivalence 
$\hocolimx \sin F \ra \sin \hocolimx F$.
\vspi
[Consider the natural transformation $\abs{\sin F} \ra F$.  Thanks to the Giever-Milner theorem, $\forall \ i$,  
$\abs{\sin Fi} \ra F i$ is a weak homotopy equivalence, thus the arrow 
$\hocolimx\abs{\sin F} \ra \hocolimx F$ is a weak homotopy equivalence (cf. supra).  But from the preceding example, 
$\abs{\hocolimx \sin F} \approx \hocolimx\abs{\sin F}$, so taking adjoints leads to the conclusion.]\\
\endgroup%%------------------------------------<<

\begingroup%%----------------------------------->>
\fontsize{9pt}{11pt}\selectfont
\textbf{\small EXAMPLE} \ 
Let $F: \bI \ra \bCGH$ be a functor $-$then there is a natural isomorphism 
$\sin \holim F \ra \holim \sin F$ of simplicial sets.\\
\endgroup%%------------------------------------<<

\begingroup%%----------------------------------->>
\fontsize{9pt}{11pt}\selectfont
\textbf{\small EXAMPLE} \ 
Let $F: \bI \ra \bSISET$ be a functor such that $\forall \ i$, $F i$ is fibrant $-$then there is a natural weak homotopy equivalence 
$\abs{\holimx F} \ra \holimx \abs{F}$.\\
\endgroup%%------------------------------------<<

Another corollary to Proposition 46 is the fact that 
$\hocolimx F \approx \abs{\lan F \circx \pi_\bI^\OP}$ 
and 
$\holimx F \approx \tot \ran F \circx \pi_\bI$.\\

\index{Simplicial Replacement Lemma}
\textbf{\small SIMPLICIAL REPLACEMENT LEMMA} \ 
Fix $F \in \Ob[\bI,\bC]$.  Define $\coprod F$ in \bSIC by 
$\left(\coprod F\right)_n = \ds\coprod\limits_{\overset{f}{[n] \ra \bI^\OP}} F f n$ 
$-$then 
$\coprod F \approx \lan F \circx \pi_\bI^\OP$.\\
\vspace{0.25cm}

\index{Cosimplicial Replacement Lemma}
\textbf{\small COSIMPLICIAL REPLACEMENT LEMMA} \ 
Fix $F \in \Ob[\bI,\bC]$.  Define $\prod F$ in \bCOSIC by 
$\left(\prod F\right)_n = \ds\prod\limits_{\overset{f}{[n] \ra \bI}} F f n$ 
$-$then 
$\prod F \approx \ran F \circx \pi_\bI$.\\
\vspace{0.25cm}

\begingroup%%----------------------------------->>
\fontsize{9pt}{11pt}\selectfont
\textbf{\small FACT} \ 
Let $F, G: \bI \ra \bSISET$ be functors and let $\Xi:F \ra G$ be a natural transformation.  Assume: 
$\forall \ i$, $\Xi_i:F i \ra G i$ is a Kan fibration $-$then 
$\holimx \Xi: \holimx F \ra \holimx G$ is a Kan fibration.
\vspi
[The arrow $\ds\prod \Xi: \ds\prod F \ra \ds\prod G$ is a fibration in $[\bDelta,\bSISET]$ (Reedy structure).  But 
$\tot:[\bDelta,\bSISET]$ $\ra$ $\bSISET$ preserves fibrations 
(cf. p. \pageref{108a}).]\\
\endgroup%%------------------------------------<<

\begingroup%%----------------------------------->>
\fontsize{9pt}{11pt}\selectfont
Application: Let $F:\bI \ra \bSISET$ be a functor.  Assume: $\forall \ i$, $F i$ is fibrant $-$then 
$\holimx F$ is fibrant.\\
\endgroup%%------------------------------------<<

\begingroup%%----------------------------------->>
\fontsize{9pt}{11pt}\selectfont
\textbf{\small EXAMPLE} \ 
Let $\rho:A \ra B$ be an inclusion of simplicial sets.  Suppose that
$F:\bI \ra \bSISET$ is a functor such that $\forall \ i$, $F i$ is $\rho$-local $-$then $\holimx F$ is $\rho$-local.

[Each $F i$ is fibrant, so $\holimx F$ is fibrant.  Denote by 
$
\begin{cases}
\ \map(A,F) \\
\ \map(B,F)
\end{cases}
$
the functor $\bI \ra \bSISET$ that sends $i$ to 
$
\begin{cases}
\ \map(A,F i) \\
\ \map(B,F i)
\end{cases}
, \ 
$
which are fibrant 
(cf. p. \pageref{13.110}).  
Since
$
\begin{cases}
\ \map(A,\holimx F) \approx \holimx \map(A,F)\\
\ \map(B,\holimx F) \approx \holimx \map(B,F)\\
\end{cases}
$
and each $F i$ is $\rho$-local, the arrow
$\map(B,\holimx F) \ra \map(A,\holimx F)$ is a weak homotopy equivalence  
(cf. p. \pageref{13.110a}).]\\
\endgroup%%------------------------------------<<

%%----------------------------------------------------------------------------------------------66

\begingroup%%----------------------------------->>
\fontsize{9pt}{11pt}\selectfont
\textbf{\small EXAMPLE} \ 
Let $\rho: \A \ra B$ be an inclusion of simplicial sets.   
Suppose that $F,G: \bI \ra \bSISET$ are functors and $\Xi:F \ra G$ is a natural transformation.  
Assume: $\forall \ i$, $\Xi_i:F i \ra G i$ is a $\rho$-equivalence $-$then 
$\hocolimx \Xi: \hocolimx F \ra \hocolim G$ is a $\rho$-equivalence.
\vspi
[It is a question of proving that the arrow 
$\map(\hocolimx(G,Z) \ra \map(\hocolimx F,Z)$ is a weak homotopy equivalence $\forall$ $\rho$-local $Z$ or still, that the arrow $\holimx \map(G,Z) \ra \holimx \map(F,Z)$ is a weak homotopy equivalence, which is true 
(cf. p. \pageref{13.111}).]\\
\endgroup%%------------------------------------<<

\begin{proposition} \ %47
For any cofibrant object $F$ in $[\bI,\bSISET]$ (structure L), the arrow 
$\hocolimx F \ra \colimx F$ is a weak homotopy equivalence.
\end{proposition}

 [It suffices to show that $\forall$ fibrant $Z$, the arrow $\map(\colimx F,Z) \ra \map(\hocolimx F,Z)$ is a weak homotopy equivalence (cf. Proposition 39).  
Since 
$\hocolimx F \ra \colimx F$ 
is induced by the projection 
$\ner(-\bs\bI)^\OP \ra *$, 
one need only consider the arrow 
$\HOM(F,\map(*,Z)) \ra \HOM(F,\map(\ner(-\bs\bI)^\OP,Z))$.  
But $F$ is a cofibrant object in $[\bI,\bSISET]$ and 
$\map(*,Z) \ra \map(\ner(-\bs\bI)^\OP,Z)$
is a weak equivalence between fibrant objects in $[\bI,\bSISET]$, thus the assertion is a consequence of Proposition 34.]\\

\begingroup%%----------------------------------->>
\fontsize{9pt}{11pt}\selectfont
\textbf{\small FACT} \ 
Suppose that \bI is filtered $-$then $\forall$ $F$ in $[\bI,\bSISET]$, the arrow $\hocolimx F \ra \colimx F$ is a weak homotopy equivalence.\\
\endgroup%%------------------------------------<<

\label{14.55}
\begingroup%%----------------------------------->>
\fontsize{9pt}{11pt}\selectfont
\textbf{\small EXAMPLE} \ 
$\forall$ $F$ in $\bFIL(\bSISET)$, the arrow $\hocolimx F \ra \colimx F$ is a weak homotopy equivalence.  
Therefore $\abs{\colim F}$ is contractible iff $\forall \ n$, $\abs{F n}$ is contractible.
\vspi
[The arrow $\hocolimx F \ra \ner [\N]^\OP$ is a weak homotopy equivalence.  
And: $[\N]^\OP$ has a final object, hence $B[\N]^\OP$ is contractible 
(cf. p. \pageref{13.112}).]\\
\endgroup%%------------------------------------<<

\textbf{\small LEMMA} \ 
If $X$ is a cofibrant $K$-simplicial ($K$-cosimplicial) object in \bSISET, then $\forall$ fibrant $Y$ in \bSISET, $\map(X,Y)$ is a fibrant $K$-cosimplicial ($K$-simplicial) object in \bSISET.\\

\begin{proposition} \ %48
For any cofibrant $K$-simplicial ($K$-cosimplicial) object $X$ in \bSISET, the arrow $\hocolimx X \ra \colimx X$ is a weak homotopy equivalence.\\
\end{proposition}

\begingroup%%----------------------------------->>
\fontsize{9pt}{11pt}\selectfont
\textbf{\small EXAMPLE} \ 
Let $B$ be a simplicial set.  
Fix an $X$ in \bSISET/\mB and determine the cofibrant \mB-cosimplicial object $X_B$ in \bSISET as on 
p. \pageref{13.113} ff. 
$-$then the arrow $\hocolimx X_B \ra \colim X_B$ $(\approx X)$ is a weak homotopy equivalence.
\vspi
[Note: \  Suppose given 
\begin{tikzcd}[sep=small]
{X} \ar{rdd}[swap]{p}   \ar{rr}{f}  &&{Y} \ar{ldd}{q}\\
\\
&{B}
\end{tikzcd}
such that $\forall \ n$ $\&$ $\forall \ b \in B_n$, $X_b \ra Y_b$ is a weak homotopy equivalence $-$then $\hocolimx X_B \ra \hocolimx Y_B$ is a weak homotopy equivalence 
(cf. p. \pageref{13.114}).  
Since there
%%----------------------------------------------------------------------------------------------67
is a \cd 
\begin{tikzcd}[sep=large]
{\hocolim X_B} \ar{d}   \ar{r}  &{X} \ar{d}{f}\\
{\hocolim Y_B} \ar{r} &{Y}
\end{tikzcd}
, it follows that $f$ is a weak homotopy equivalence.\\
Example: $p$ is a weak homotopy equivalence if the $\abs{X_b}$ are contractible.]\\
\endgroup%%------------------------------------<<

Given a category \bC, write \bBISIC for the functor category 
$[(\bDelta \times \bDelta)^\OP,\bC]$ (i.e., $[\bDelta^\OP,\bSIC]$) 
$-$then by definition, a 
\un{bisimplicial object} 
\index{bisimplicial object} 
in \bC is an object in \bBISIC (i.e., a simplicial object in \bSIC).  
Example: Assuming that \bC has finite products, if 
$
\begin{cases}
\ X \\
\ Y
\end{cases}
$
are simplicial objects in \bC, the assignment $([n],[m]) \ra X_n \times Y_m$ 
defines a bisimplicial object $X \un{\times} Y$ in \bC.

Specialize to \bC = \bSET $-$then an object in \bBISISET ($= \reallywidehat{\bDelta \times \bDelta}$) is called a 
\un{bisimplicial} 
\un{set} 
\index{bisimplicial set} 
and a morphism in \ \bBISISET is called a \ 
\un{bisimplicial map}.  
\index{bisimplicial map}  \ \ 
Given a bisimplicial set $X$, put $X_{n,m} = X([n],[m])$ $(= X_n[m]))$ $-$then there are horizontal operators 
$
\begin{cases}
\ d_i^h:X_{n,m} \ra X_{n-1,m}\\
\ s_i^h:X_{n,m} \ra X_{n+1,m}
\end{cases}
(0 \leq i \leq n)
$
and vertical operators
$
\begin{cases}
\ d_j^v:X_{n,m} \ra X_{n,m-1}\\
\ s_j^v:X_{n,m} \ra X_{n,m+1}
\end{cases}
(0 \leq j \leq m).
$
\ 
The horizontal operators commute with the vertical operators, the simplicial identities are satisfied horizontally and vertically, 
and thanks to the Yoneda lemma, 
$\Nat(\Delta[n,m],X)$ $\approx X_{n,m}$, 
where $\Delta[n,m] = \Delta[n] \un{\times} \Delta[m]$.

[Note: \  Every simplicial set $X$ can be regarded as a bisimplicial set by trivializing its structure in either the horizontal or vertical direction, i.e., $X_{n,m} = X_m$ or $X_{n,m} = X_n$.]\\

\begingroup%%----------------------------------->>
\fontsize{9pt}{11pt}\selectfont
\textbf{\small EXAMPLE} \ 
Any functor $T:\bDelta \ra \bCAT$ gives rise to a functor $X_T:\bCAT \ra \bBISISET$ by writing 
$X_T \bI([n],[m]) = \nersub_n(T[m],\bI)$ 
$(\approx \Nat([n],[T[m],\bI) \approx$ 
$\Nat(T[m],[[n],\bI]) \approx$ $(S_T[[n],\bI])_m$, \ 
$S_T$ the singular functor 
(cf. p. \pageref{13.115})).\\
\endgroup%%------------------------------------<<

\label{18.7}
\begingroup%%----------------------------------->>
\fontsize{9pt}{11pt}\selectfont
\textbf{\small EXAMPLE} \ 
Let \bC be a double category, i.e., a category object in \bCAT $-$then $\ner \bC$ is a simplicial object in \bCAT, hence 
$\ner(\ner \bC)$ is a bisimplicial set.\\
\endgroup%%------------------------------------<<

Viewing [$n$] as a small category, one may form its simplex category $\bDelta[n]$ 
($=\bDelta \ner[n] = \bDelta\Delta[n] = \bDelta/[n]$).  The assignments
$[n] \ra \ner\bDelta[n]$, $[n] \ra \Delta[n]$ 
define cosimplicial objects $\textbf{Y}_{\bDelta}$, $Y_{\bDelta}$ in \bSISET 
which are cofibrant in the Reedy structure and there is a weak equivalence $\textbf{Y}_{\bDelta} \ra Y_{\bDelta}$ 
(cf. p. \pageref{13.116}).

Let $X$ be a bisimplicial set $-$then 
$\hocolimx X = $
$\ds\int^{[n]} X_n \times \ner([n]\bs\bDelta^\OP)^\OP = $
$\ds\int^{[n]} X_n \times \ner(\bDelta/[n]) = $
$\ds\int^{[n]} X_n \times \ner \bDelta[n] \ra $
$\ds\int^{[n]} X_n \times \Delta[n] = $
$\abs{X}$.\\
\vspace{0.25cm}

\begin{proposition} \ %49
The arrow $\hocolimx X \ra \abs{X}$ is a weak homotopy equivalence.
\end{proposition}

%%----------------------------------------------------------------------------------------------68
[Bearing in mind Proposition 39, take a fibrant $Z$ and consider the arrow 
$\map(\abs{X},Z)$ $\ra$ $\map(\hocolimx X,Z)$
or still, the arrow 
$\HOM(X,\map(Y_{\bDelta},Z)) \ra \HOM(X,\map(\textbf{Y}_{\bDelta},Z))$.  
In the Reedy structure, $X$ is necessarily cofibrant while 
$\map(Y_{\bDelta},Z) \ra \map(\textbf{Y}_{\bDelta},Z)$ 
is a weak equivalence between fibrant objects (see the lemma prefacing Proposition 48).  
One may therefore quote Proposition 34.]\\

Using the notation of the Kan extension theorem, take 
$\bC = \bDelta^\OP$,
$\bD = \bDelta^\OP \times \bDelta^\OP$,
$\bS = \bSET$, 
and let $K$ be the diagonal 
$ \bDelta^\OP \ra \bDelta^\OP \times \bDelta^\OP$
$-$then the functor 
$[K,\bS] \equiv \di:\bBISISET \ra \bSISET$ 
has both a right and left adjoint.  
One calls di 
\index{diagonal functor} the 
\un{diagonal}: 
$(\di X)_n = X([n],[n])$, the operators being
$
\begin{cases}
\ d_i = d_i^h d_i^v = d_i^v d_i^h\\
\ s_i = s_i^h s_i^v = s_i^v s_i^h
\end{cases}
.
$
Example: $\di(X \un{\times} Y) = X \times Y$ ($\implies$ $\di\Delta[n,m] = \Delta[n] \times \Delta[m]$).\\

\begin{proposition} \  %50
Up to natural isomorphism, $\di$ and $\abs{?}$ are the same.
\end{proposition}

[It suffices to prove that di is a left adjoint for 
$\sin:\Nat(\di X,Y) \approx \Nat(X,\sin Y)$.  
But
$X \approx \colim_{i,j} \Delta[n_i,m_j]$ 
and one has 
$\Nat(\Delta[n,m],\sin Y) \approx $
$\map(\Delta[n],Y)_m \approx $
$\Nat(\Delta[n] \times \Delta[m],Y) \approx $
$\Nat(\di \Delta[n,m],Y)$.]\\

Application: $\forall$ bisimplicial set $X$, there is a weak homotopy equivalence $\hocolimx X \ra \di X$.\\

\begingroup%%----------------------------------->>
\fontsize{9pt}{11pt}\selectfont
\textbf{\small EXAMPLE} \ 
Let $F:\bI \ra \bSISET$ $-$then in the notation of the simplicial replacement lemma, $F$ determines a bisimplicial set $\coprod F$ by the rule 
$\left(\coprod F\right)_n = \ds\coprod\limits_{[n]\overset{f}{ \ra}\bI^\OP} F f n$. \ 
And:
$\hocolimx F \approx $
$\abs{\coprod F}  \approx $
$\di \coprod F$.\\
\endgroup%%------------------------------------<<

\label{14.32}
\begingroup%%----------------------------------->>
\fontsize{9pt}{11pt}\selectfont
\textbf{\small EXAMPLE} \ 
Place on \bCGH its singular structure and equip $[\bDelta^\OP,\bCGH]$ with the corresponding Reedy structure.  
Take an $X$ in \bSICGH 
which is both Reedy fibrant and Reedy cofibrant and let $UX$ be the simplicial set obtained from $X$ by forgetting the topologies $-$then the arrow $\abs{UX} \ra \abs{X}$ is a weak homotopy equivalence.  
To see this, let $\sin X$ be the bisimplicial set defined by 
$(\sin X)_n = \sin X_n$ and write $\sin^\Tee X$ for the ``transpose'' of $\sin X$, i.e., 
$(\sin^\Tee X)_{n,m} = (\sin X)_{m,n}$ $(\implies (\sin^\Tee X)_{0,*} \approx UX)$.  
Since $\sin X$ is Reedy fibrant $\forall \ \alpha$, $\sin^\Tee X(\alpha)$ is a weak homotopy equivalence.  
Therefore the arrow
$\abs{\sin^\Tee X}_0 \ra \abs{\sin^\Tee X}$ is a weak homotopy equivalence 
(cf. p. \pageref{13.117}).  
Write $\abs{\sin X}$ for the simplicial object in \bCGH with 
$\abs{\sin X}_n = \abs{\sin X_n}$.  
Because $\abs{\sin X}$ is Reedy cofibrant, in view of the Giever-Milnor theorem, the arrow 
$\norm{\sin X} \ra \abs{X}$ is a weak homotopy equivalence (cf. Proposition 44).  
So, putting everything together gives
$\abs{UX} \approx $
$||\sin^\Tee X|_0| \overset{\sim}{\ra} $
$\norm{\sin^\Tee X} \approx $
$\abs{\di \sin^\Tee X} = $
$\abs{\di \sin X} \approx $
$\norm{\sin X} \overset{\sim}{\ra} $
$\abs{X}$.\\
\endgroup%%------------------------------------<<

%%----------------------------------------------------------------------------------------------69
\begin{proposition} \ %51
Suppose that $f:X \ra Y$ is a bisimplicial map.  Assume: $\forall \ n$, $f_n:X_n \ra Y_n$ is a weak homotopy equivalence $-$then 
$\di f:\di X \ra \di Y$ is a weak homotopy equivalence.
\end{proposition}

[Since all simplicial objects in $\widehat{\bDelta}$ are cofibrant in the Reedy structure, this is a consequence of Propositions 44 and 50.]

[Note: \  In both the statement and the conclusion, one can replace ``weak homotopy equivalence'' by ``HG-equivalence'' 
(cf. p. \pageref{13.118}).]\\

\begingroup%%----------------------------------->>
\fontsize{9pt}{11pt}\selectfont
Let $X \overset{f}{\ra} Z \overset{g}{\la} Y$ be a 2-sink in \bSISET $-$then a \cd 
\begin{tikzcd}[sep=large]
{W} \ar{d}   \ar{r}  &{Y} \ar{d}{g}\\
{X} \ar{r}[swap]{f} &{Z}
\end{tikzcd}
is said to be a 
\un{pullback up to homotopy} 
\index{pullback up to homotopy} 
if the arrow $W \ra X \times_Z Y$ is a weak homotopy equivalence.  
Example: 
\begin{tikzcd}[sep=large]
{\dot\Delta[1]} \ar{d}   \ar{r}  &{\Lambda[1,2]} \ar{d}\\
{\Delta[1]} \ar{r} &{\Delta[2]}
\end{tikzcd}
is not a homotopy pullback but is a pullback up to homotopy.\\
\endgroup%%------------------------------------<<

\begingroup%%----------------------------------->>
\fontsize{9pt}{11pt}\selectfont
\textbf{\small FACT} \ 
Let $f:X \ra Y$ be a bisimplicial map.  Assume $\forall \ m,\ n$ $\&$ $\forall \ \alpha:[m] \ra [n]$, the \cd
\begin{tikzcd}[sep=large]
{X_n} \ar{d}[swap]{f_n}   \ar{r}{X_\alpha}  &{X_m} \ar{d}{f_m}\\
{Y_n} \ar{r}[swap]{Y_\alpha} &{Y_m}
\end{tikzcd}
is a pullback up to homotopy $-$then $\forall \ n$, 
\begin{tikzcd}[sep=large]
{X_n \times \Delta[n]} \ar{d}   \ar{r}  &{\di X} \ar{d}\\
{Y_n \times \Delta[n]} \ar{r} &{\di Y}
\end{tikzcd}
is a pullback up to homotopy.\\
\endgroup%%------------------------------------<<

\begin{proposition} \  %52
\bBISISET carries a proper model category structure in which a bisimplicial map $f:X \ra Y$ is a \we if $\di f$ is a weak homotopy equivalence, a fibration if $\di f$ is a Kan fibration, and a cofibration if $f$ has the LLP w.r.t acyclic fibrations.
\end{proposition}

[This is an instance of the generalities on 
p. \pageref{13.119}, 
the essential point being that di (which plays the role of ``\mG'') has both a right and left adjoint.  In particular: di preserves filtered colimits.  The stage is thus set for a small category argument.  
Let \mD be the left adjoint of di normalized by the condition $D\Delta[n] = \Delta[n,n]$.  Put
$\dot\Delta[n,n] = D\dot\Delta[n]$, 
$\Lambda[k,n,n] = D\Lambda[k,n]$ 
$-$then the arrow $\dot\Delta[n,n] \ra \Delta[n,n]$ is a cofibration and the arrow $\Lambda[k,n,n] \ra \Delta[n,n]$ is an acyclic cofibration ($\abs{\di \Lambda[k,n,n]}$ is contractible).  
The requisite factorizations can therefore be established in the usual way.  
Let us note only that every $f$ admits a decomposition of the form $f = p \circx i$, where $p$ is a fibration and $i$ is an acyclic cofibration that has the LLP w.r.t. fibrations (specifically, $i$ is a sequential colimit of pushouts of coproducts of inclusions 
$\Lambda[k,n,n] \ra \Delta[n,n]$).  
As for  properness, the part 
%%----------------------------------------------------------------------------------------------70
of PMC concerning pullbacks is obvious while the part concering pushouts follows from the observation that a cofibration is necessarily an injective bisimplicial map.]\\

\begingroup%%----------------------------------->>
\fontsize{9pt}{11pt}\selectfont
\textbf{\small FACT} \  Take \bBISISET in the model category structure supplied by Proposition 52 $-$then the adjoint pair $(D,\di)$ induces an adjoint equivalence of categories between \bHSISET and \bHBISISET.\\
\endgroup%%------------------------------------<<

\label{14.28a}
\label{14.66}
For certain purposes, it is technically more convenient to use a modification of the homotopy colimit in order to minimize the proliferation of opposites.  
Definition: Given $F \in \Ob [\bI,\bC]$, put $\ohc_\bI F$ 
\index{$\ohc$} 
(or $\ohc F)$ $= \ds\int^{\bI^\OP} F \ \xbox \ \ner(-\bs \bI)$. 
The formal properties of $\ohc$ are the same as those of hocolim, the primary difference being that 
$\ohc F \approx \abs{\coprod F}$, where now $(\coprod F)_n = \ds\coprod\limits_{[n] \overset{f}{\ra}\bI } F f 0$.\\

\label{14.28}

\begingroup%%----------------------------------->>
\fontsize{9pt}{11pt}\selectfont
\textbf{\small EXAMPLE} \ 
Let $F:\bI \ra \bCAT$ be a functor $-$then the 
\un{Grothendieck construction} 
\index{Grothendieck construction} 
on $F$ is the category $\gro_\bI F$ 
\index{$\gro_\bI F$} 
whose objects are the pairs $(i,X)$, where $i \in \Ob\bI$ and $X \in \Ob Fi$, and whose morphisms are the arrows $(\delta,f):(i,X) \ra (j,Y)$, where $\delta \in \Mor(i,j)$ and $f \in \Mor((F\delta)X,Y)$ 
(composition is given by 
$(\delta^\prime,f^\prime) \circx (\delta,f) =$ $(\delta^\prime \circx \delta,f^\prime \circx (F\delta^\prime)f)$).  
Put $NF = \ner \circx F$, so $NF:\bI \ra \bSISET$.  One can thus form $\ohc NF$ and 
Thomason\footnote[2]{\textit{Math. Proc. Cambridge Philos. Soc.} \textbf{85} (1979), 91-109; 
see also Heggie, \textit{Cahiers Topologie G\'eom. Diff\'erentielle Cat\'egoriques} 34 (1993), 13-36.}
has shown that there is a natural weak homotopy equivalence $\eta:\ohc NF \ra \ner\gro_\bI F$.  
The situation for homotopy limits is simpler.  
Indeed, 
$\holim NF \approx$ 
$\ds\int_i \map(\ner(\bI/i),(\ner \circx F)i) \approx$ 
$\ds\int_i \ner[\bI/i,Fi] \approx$ 
$\ner\bigl(\ds\int_i [\bI/i,Fi]\bigr)$.
\vspi
[Note: \  Here is the definition of $\eta$.  
Representing $\ohc NF$ as $\di \ds\coprod NF$, fix $n$ and consider a typical string
$(i_0 \overset{\delta_0}{\ra} i_1 \ra \cdots \ra i_{n-1} \overset{\delta_{n-1}}{\ra} i_n, X_0, \ra X_1, \ra \cdots \ra X_{n-1} \ra X_n)$, where the $X_k \in \Ob F_{i_0}$ $(0 \leq k \leq n)$ $-$then $\eta_n$ takes it to the element of 
$\nersub_n\gro_\bI F$ given by 
$(i_0,X_0) \ra (i_1,(F\delta_0)X_1) \ra \cdots \ra (i_n(F\delta_{n-1} \circx \cdots \circx F\delta_0)X_n)$.]\\
\endgroup%%------------------------------------<<

Let \bI and \bJ be small categories, $\nabla :\bJ \ra \bI$ a functor.

Notation: Given $i \in \Ob\bI$, write $i\backslash\nabla$ for the comma category $\abs{K_i,\nabla}$.

[Note: \  Dually, $\nabla/ i$ stands for the comma category $\abs{\nabla,K_i}$.]

Observation: The commutative diagram
\begin{tikzcd}%[sep=large]
{i\backslash\nabla} \ar{d}   \ar{r}  &{\bJ} \ar{d}\\
{i\backslash\bI} \ar{r} &{\bI}
\end{tikzcd}
is a pullback square in \bCAT.

[Note: \  The 
\un{fiber} 
\index{fiber ($\nabla$ over $i$)} 
of $\nabla$ over $i$ is defined by the pullback square
$
\begin{tikzcd}%[sep=large]
{\nabla^{-1}(i)} \ar{d}   \ar{r}  &{\bJ} \ar{d}{\nabla}\\
{\textbf{1}} \ar{r}[swap]{K_i} &{\bI}
\end{tikzcd}
.
$
So:
%%----------------------------------------------------------------------------------------------71
$\nabla^{-1}(i)$ is the subcategory of \bJ having objects $j$ such that $\nabla j = i$, morphisms $\delta$ such that $\nabla \delta = \id_i$, and there is a commutative diagram
\begin{tikzcd}%[sep=large]
{\nabla^{-1}(i)} \ar{d} \ar{dr}   \ar{r}  &{i\backslash\nabla} \ar{d}\\
{\nabla/i} \ar{r} &{\bJ}
\end{tikzcd}
.]\\
\vspace{0.25cm}

\begingroup%%----------------------------------->>
\fontsize{9pt}{11pt}\selectfont
\textbf{\small EXAMPLE} \ 
The arrow 
$\colimx \ner(-\bs\nabla) \ra \ner \bJ$ 
is an isomorphism.  \ 
Viewed as an object in $[\bI^\OP,\bSISET]$ (structure L), \ 
$\ner(-\bs\nabla)$ is free, hence cofibrant 
(cf. p. \pageref{13.120}).  
Therefore the arrow 
$\ohc \ner (-\bs\nabla) \ra \colimx \ner (-\bs\nabla)$ $(\approx \ner \bJ)$ 
is a weak homotopy equivalence (cf. Proposition 47).  
\vspi
[Note: \  Take \bI = \bJ and $\nabla = \id_\bI$ $-$then the arrow 
$\ohc \ner (-\bs\bI) = $
$\ds\int^i \ner(i\bs\bI) \times \ner(i\bs\bI^\OP) \ra$ 
$\ds\int^i \ner(i\bs\bI) \times \Delta[0] \approx$ 
$\ner \bI$ 
is a weak homotopy equivalence, as is the arrow
$\ohc \ner (-\bs\bI) = $
$\ds\int^i \ner(i\bs\bI) \times \ner(i\bs\bI^\OP) \ra$ 
$\ds\int^i \Delta[0] \times \ner(i\bs\bI^\OP) \approx$ 
$\ner \bI^\OP$.]\\
\endgroup%%------------------------------------<<

\textbf{\small LEMMA} \ 
Let \bI and \bJ be small categories, $\nabla:\bJ \ra \bI$ a functor 
$-$then $\forall \ F$ in $[\bI^\OP,\bSISET]$, 
$\ds\int_i \map(\ner(i\bs\nabla),Fi) \approx$ 
$\ds\int_j \map(\ner(j\bs\bJ),(F \circx \nabla^\OP)j)$, 
i.e., 
$\HOM(\ner(-\bs\nabla),$ $F) \approx$ 
$\HOM(\ner(-\bs\bJ),F \circx \nabla^\OP)$.

[The left Kan extension of $\ner(-\bs\bJ)$ along $\nabla^\OP$ is $\ner(-\bs\nabla)$.]\\

\begin{proposition} \ %53
Let \bI and \bJ be small categories, $\nabla:\bJ \ra \bI$ a functor $-$then $\forall \ F$ in $[\bI,\bSISET]$, the arrow
$\ds\int^j (F \circx \nabla)j \times \ner(j\bs\bJ) \ra \ds\int^i Fi \times \ner(i\bs\nabla)$
is a weak homotopy equivalence.
\end{proposition}

[This is yet another application of Proposition 39.  Thus fix a fibrant $Z$ and pass to 
$\map(\ds\int^i Fi \times \ner(i\bs\nabla),Z) \ra$ 
$\map(\ds\int^i (F \circx \nabla)j \times \ner(j\bs\bJ),Z)$, 
i.e., to
$\ds\int_i\map(\ner(i\bs\nabla),\map$ $(Fi,Z)) \ra$
$\ds\int_j \map(\ner(j\bs\bJ),\map((F \circx \nabla)j,Z))$, 
i.e., to 
$\HOM(\ner(-\bs\nabla),\map(F,Z)) \ra$ 
$\HOM(\ner(-\bs\bJ),\map(F,Z) \circx \nabla^\OP)$, 
which by the lemma is an isomorphism, hence a fortiori, a weak homotopy equivalence.]\\

A small category is 
\un{contractible} 
\index{contractible (small category)} 
if its classifying space is contractible.  Example: Every filtered category is contractible.\\

\begingroup%%----------------------------------->>
\fontsize{9pt}{11pt}\selectfont
\textbf{\small EXAMPLE} \ 
Let \bC be a small category $-$then the 
\un{cone} 
\index{cone (small category)} $\Gamma \bC$ 
\index{$\Gamma \bC$}  of \bC is the small category with $\Ob \Gamma \bC = \Ob \bC \ds\coprod \{\emptyset\}$, 
where $\emptyset$ is an adjoined initial object.  
Example: $\Gamma \textbf{0} = \textbf{1}$.  So $\Gamma \bC$ is contractible 
(cf. p. \pageref{13.121}) and 
$B \Gamma \bC \approx \Gamma B \bC$.
\vspi
%%----------------------------------------------------------------------------------------------72
[Note: \  Given small categories 
$
\begin{cases}
\ \bC \\
\ \bD
\end{cases}
, \ 
$
their 
\un{join} 
\index{join (small categories)} 
$\bC * \bD$ %\index{$\bC * \bD$} 
is the full subcategory of $\Gamma \bC \times \Gamma \bD$ with 
$\Ob\bC * \bD = \Ob\bC \times \Ob\bD \ds\coprod \Ob\bC \times \{\emptyset\} \ds\coprod \{\emptyset\} \times \Ob\bD$.
Under the join, \bCAT is a symmetric monoidal category (\textbf{0} is the unit).  One has 
$B(\bC * \bD) \approx B\bC *_k B\bD$].\\
\endgroup%%------------------------------------<<

Given small categories 
$
\begin{cases}
\ \bI\\
\ \bJ
\end{cases}
, \ 
$
a functor $\nabla:\bJ \ra \bI$ is said to be 
\un{strictly final} 
\index{strictly final (functor)} 
provided that for every $i \in \Ob\bI$, the comma category $\abs{K_i,\nabla}$ is contractible.  
A strictly final functor is final.  
In particular: $\nabla:\bJ \ra \bI$ strictly final $\implies$ 
$\colimx \Delta \circx \nabla \approx$ $\colimx\Delta$, where $\Delta:\bI \ra \bSISET$ 
(cf. p. \pageref{13.122}).

[Note: \  A subcategory of a small category is 
\un{strictly final}
\index{strictly final (subcategory)} 
if the inclusion is a strictly final functor.]\\

\begin{proposition} \ %54
Let \bI and \bJ be small categories, $\nabla:\bJ \ra \bI$ a strictly final functor 
$-$then $\forall \ F$ in $[\bI,\bSISET]$, the arrow
%$\overline{\text{hocolim}} F \circx \nabla \rightarrow \overline{\text{hocolim}} F$ is a weak homotopy equivalence.\\
$\ohc F \circx \nabla \ra \ohc F$ is a weak homotopy equivalence.
\end{proposition}

[According to Proposition 53, the arrow 
$\ohc F\circx\nabla = $
$\ds\int^j (F \circx\nabla)_j \times \ner(j\bs\bJ) \ra$ 
$\ds\int^i Fi \times \ner(i\bs\nabla)$ 
is a weak homotopy equivalence.  Claim: The arrow
$\ds\int^i Fi \times \ner(i\bs\nabla) \ra$ 
$\ds\int^i Fi \times \ner(i\bs\bI) =$ 
$\ohc F$ 
is a weak homotopy equivalence.  
Indeed: $\ner(-\bs\nabla)$, $\ner(-\bs\bI)$ are cofibrant objects in 
$[\bI^\OP,\bSISET]$ and since $\nabla$ is strictly final, the arrow 
$\ner(-\bs\nabla) \ra \ner(-\bs\bI)$ is a \we.  
Therefore one may appeal to the example on 
p. \pageref{13.123}.]\\

\begingroup%%----------------------------------->>
\fontsize{9pt}{11pt}\selectfont
\textbf{\small FACT} \ 
Let \bI and \bJ be small categories, $\nabla:\bJ \ra \bI$ a functor.  
Assume: $\ner \nabla:\ner\bJ \ra \ner\bI$ is a  weak homotopy equivalence.  Suppose that 
$F:\bI \ra \bSISET$ sends the morphisms in \bI to weak homotopy equivalences 
$-$then the arrow $\ohc F \circx \nabla \ra \ohc F$ is a weak homotopy equivalence.
\vspi
[The commutative diagram
\begin{tikzcd}[sep=large]
{\ds\coprod F \circx \nabla} \ar{d} \ar{d}   \ar{r}  &{\ds\coprod F} \ar{d}\\
{\ds\coprod *} \ar{r} &{\ds\coprod *}
\end{tikzcd}
of bisimplicial sets is a pullback square, therefore the commutative diagram
\begin{tikzcd}[sep=large]
{\di\ds\coprod F \circx \nabla} \ar{d} \ar{d}   \ar{r}  &{\di\ds\coprod F} \ar{d}\\
{\di\ds\coprod *} \ar{r} &{\di\ds\coprod *}
\end{tikzcd}
of simplicial sets is a pullback square (di is a right adjoint).  Accordingly, in \bSISET, the \cd
\begin{tikzcd}[sep=large]
{\ohc F \circx \nabla} \ar{d} \ar{d}   \ar{r}  &{\ohc F} \ar{d}\\
{\ner \bJ} \ar{r} &{\ner \bI}
\end{tikzcd}
is a 
%%----------------------------------------------------------------------------------------------73
pullback square.  The result thus follows from the fact that the arrow $\ohc F \ra \ner \bI$ is a homotopy fibration 
(cf. p. \pageref{13.124}).]\\
\endgroup%%------------------------------------<<

\begingroup%%----------------------------------->>
\fontsize{9pt}{11pt}\selectfont
\textbf{\small EXAMPLE} \ 
If \bI is contractible and if $F:\bI \ra \bSISET$ sends the morphisms in \bI to weak homotopy equivalences, then $\forall \ i \in \Ob\bI$, the arrow $Fi \ra \ohc F$ is a weak homotopy equivalence.\\
\endgroup%%------------------------------------<<

\index{homotopy pushdowns}
\begingroup%%----------------------------------->>
\fontsize{9pt}{11pt}\selectfont
\textbf{\small FACT \  (\un{Homotopy Pushdowns})} \ 
Let \bI and \bJ be small categories, $\nabla:\bJ \ra \bI$ a functor.  
Given a functor $G:\bJ \ra \bSISET$, define an object 
$\ohc_\nabla G$ in $[\bI,\bSISET]$ by $(\ohc_\nabla G)i = \ohc_{\nabla/i} G \circx U_i$, where 
$U_i:\nabla/i \ra \bJ$ is the forgetful functor $-$then the arrow 
$\ohc_{\bI} \ohc_\nabla G \ra \ohc_{\bJ} G$ 
is a weak homotopy equivalence.\\
\endgroup%%------------------------------------<<

\index{Theorem Quillen's Theorem A}
\textbf{\small QUILLEN'S THEOREM A} \ \ 
Suppose that \bI  and \bJ are small categories and $\nabla:\bJ \ra \bI$ is a strictly final functor $-$then 
$\ner\nabla:\ner\bJ \ra \ner\bI$ 
is a weak homotopy equivalence, hence 
$B\nabla:B\bJ \ra B\bI$ 
is a homotopy equivalence.

[In Proposition 54, let $F$ be the functor $\bI \ra \bSISET$ that sends 
$i \in \Ob\bI$ to $Fi = \Delta[0]$.]

[Note: \  The same conclusion obtains if $\nabla$ is ``strictly initial''.]\\

\begingroup%%----------------------------------->>
\fontsize{9pt}{11pt}\selectfont
\textbf{\small EXAMPLE} \ 
Let $X$ be a topological space, sin $X$ its singular set $-$then sin $X$ can be regarded as a category:
\begin{tikzcd}[sep=small]
{\Delta^m} \ar{rdd}   \ar{rr}{\Delta^\alpha}  &&{\Delta^n} \ar{ldd}\\
\\
&{X}
\end{tikzcd}
($\alpha \in \Mor([m],[n])$) 
(cf. p. \pageref{13.125}).  
This category is isomorphic to $\Delta/X \equiv \text{gro}_\Delta \sin X$ and there is a natural weak homotopy equivalence $\ner \Delta/X \ra \sin X$ 
(cf. p. \pageref{13.126}), 
which thus gives a natural weak homotopy equivalence $B\Delta/X \ra X$ (Giever-Milnor theorem). 
Let \bC be any small full subcategory of \bTOP/\mX containing $\Delta/X$ as a subcategory.  Assume: $\forall \ Y \ra X$ in \bC, Y is homotopically trivial $-$then the arrow $B\iota:B\Delta/X \ra B\bC$ induced by the inclusion $\iota:\Delta/X \ra \bC$ is a homotopy equivalence.  
To see this, one can suppose that $X$ is nonempty and appeal to Quillen's theorem A.  
Claim: $\iota$ is a strictly initial functor i.e., $\forall \ Y \ra X$ in \bC, the comma category $\iota/Y \ra X$ is contractible.  
Indeed, $\iota/Y \ra X$ is simply $\Delta/Y$ and the arrow $B\Delta/Y \ra *$ is a weak homotopy equivalence, hence a homotopy equivalence.\\
\endgroup%%------------------------------------<<

Let \bC be a category $-$then the 
\un{twisted arrow category} 
\index{twisted arrow category} 
\bC$(\leadsto)$ of \bC is the category whose objects are the arrows $f:X \ra Y$ of \bC and whose morphisms $f \ra f^\prime$ are the pairs
$
(\phi,\psi): 
\begin{cases}
\ \phi \in \text{Mor}(X^\prime,X)\\
\ \psi \in \text{Mor}(Y,Y^\prime)
\end{cases}
$
for which the square
\begin{tikzcd}%[sep=large]
{X}    \ar{r}{f}  &{Y} \ar{d}{\psi}\\
{X^\prime}  \ar{u}{\phi} \ar{r}[swap]{f^\prime} &{Y^\prime}
\end{tikzcd}
commutes.  Denote by 
$
\begin{cases}
\ s\\
\ t
\end{cases}
$
the canonical projections
$
\begin{cases}
\ \bC(\leadsto) \ra \bC^\OP \\
\ \bC(\leadsto) \ra \bC
\end{cases}
. 
$
\\
\vspace{0.5cm}

%%----------------------------------------------------------------------------------------------74

\begingroup%%----------------------------------->>
\fontsize{9pt}{11pt}\selectfont
\textbf{\small EXAMPLE} \ 
Suppose that \bC is a small category \ $-$then
$\ner s: \ner\bC(\leadsto) \ra \ner \bC^\OP$, 
$\ner t: \ner\bC(\leadsto) \ra \ner \bC$
are weak homotopy equivalences.
\vspi
[To discuss $\ner s$, observe that $\forall$ $X$, the functor $X\bs\bC \ra s/X$ that sends $X \overset{f}{\ra} Y$ to  
$(X \overset{f}{\ra} Y,\id_X)$ 
(so $s(X \overset{f}{\ra} Y) \overset{\id_X}{\lra} X)$ has a left adjoint.  
Since $X\bs\bC$ is contractible, $s/X$ must be too 
(cf. p. \pageref{13.127}), 
i.e., $s$ is strictly initial, thus by Quillen's theorem A, $\ner s$ is a weak homotopy equivalence.]
\vspi
[Note: \  It is a corollary that ner \bC and $\ner \bC^\OP$ are naturally weakly equivalent.]\\
\endgroup%%------------------------------------<<

Let \bI and \bJ be small categories, $\nabla:\bJ \ra \bI$ a functor, then by $\nabla(\leadsto)$ 
\index{$\nabla(\leadsto)$ } 
we shall understand the category whose 
objects are the triples $(i, \delta, j)$, where $\delta:i \ra \nabla_j$, and whose 
morphisms $(i, \delta, j) \ra (i^\prime, \delta^\prime, j^\prime)$ are the pairs $(\phi,\psi)$ : 
$
\begin{cases}
\ \phi \in \Mor(i^\prime,i)\\
\ \psi \in \Mor(j,j^\prime)
\end{cases}
$
for which the square
\begin{tikzcd}%[sep=large]
{i}   \ar{r}{\delta}  &{\nabla j} \ar{d}{\nabla\psi}\\
{i^\prime}  \ar{u}{\phi} \ar{r}[swap]{\delta^\prime} &{\nabla j^\prime} 
\end{tikzcd}
commutes.  Example:  $\id_{\bI}(\leadsto) = \bI(\leadsto)$.\\
\vspace{0.25cm}

\index{Theorem Quillen's Theorem B}
\textbf{\small QUILLEN'S THEOREM B} \ \ 
Suppose that \bI and \bJ are small categories and $\nabla:\bJ \ra \bI$ is a functor with the property that for every morphism 
$i^\prime \ra i\pp$ in \bI, the arrow $\ner(i\pp\bs\nabla) \ra \ner(i^\prime\bs\nabla)$ is a weak homotopy equivalence 
$-$then $\forall \ i \in \Ob \bI$, the pullback square \qquad
\begin{tikzcd}%[sep=large]
{\ner(i\bs\nabla)} \ar{d}  \ar{r}  &{\ner \bJ} \ar{d}\\
{\ner(i\bs\bI)}  \ar{r} &{\ner \bI} 
\end{tikzcd}
is a homotopy pullback.

[Each of the squares in the commutative diagram
\begin{tikzcd}%[sep=large]
{i\bs\nabla} \ar{d}  \ar{r}  &{\nabla(\leadsto)} \ar{d}  \ar{r} &{\bJ} \ar{d}\\
{i\bs\bI}      \ar{d}  \ar{r}  &{\bI(\leadsto)}      \ar{d}  \ar{r} &{\bI}\\
{\textbf{1}} \ar{r} &{\bI^\OP}
\end{tikzcd}
are pullback squares in \bCAT, hence each of the squares in the commutative diagram
\[
\begin{tikzcd}%[sep=large]
{\ner(i\bs\nabla)} \ar{d}  \ar{r}  &{\ner\nabla(\leadsto)} \ar{d}  \ar{r} &{\ner\bJ} \ar{d}\\
{\ner(i\bs\bI)}      \ar{d}  \ar{r}  &{\ner\bI(\leadsto)}      \ar{d}  \ar{r} &{\ner\bI}\\
{\Delta[0]} \ar{r} &{\ner\bI^\OP}
\end{tikzcd}
\]
are pullback squares  in \  \bSISET \  (ner is a right adjoint).  \ 
And, from the definitions, 
$\ohc \ner(-\bs\nabla)$ \ $\approx$ \  $\ner \nabla(\leadsto)$, \ 
$\ohc \ner(-\bs\bI)$ \  $\approx$  \ $\ner \bI(\leadsto)$. \ \ 
Since the arrows 
$\ohc \ner(-\bs\bI) \ra \ner \bI^\OP$, \ 
$\ner(i\bs\bI) \ra \Delta[0]$ 
are weak homotopy equivalences, the commutative diagram
\begin{tikzcd}%[sep=large]
{\ner(i\bs\bI)} \ar{d}  \ar{r} &{\ohc \ner(-\bs \bI)} \ar{d}\\
{\Delta[0]}                 \ar{r} &{\ner \bI^\OP} 
\end{tikzcd}
is a homotopy pullback 
%%----------------------------------------------------------------------------------------------75
(cf. p. \pageref{13.128}); 
since the arrows
$\ohc \ner(-\bs\nabla) \ra \ner \bJ$,
$\ohc \ner(-\bs\bI) \ra \ner \bI$ 
are weak homotopy equivalences, the commutative diagram \quad
\begin{tikzcd}%[sep=large]
{\ohc \ner(-\bs\nabla)} \ar{d}  \ar{r} &{\ner \bJ} \ar{d}\\
{\ohc \ner(-\bs\bI)}                \ar{r} &{\ner \bI} 
\end{tikzcd}
is a \ homotopy \ pullback 
(cf. p. \pageref{13.129}).  \ 
Owing  to  our \ assumption on  $\nabla$, the arrow \quad
$\ohc \ner(-\bs\nabla)$ $\ra$ $\ner \bI^\OP$ 
is a homotopy fibration 
(cf. p. \pageref{13.130}).  \quad
Accordingly, the pullback square
\begin{tikzcd}%[sep=large]
{\ner(i\bs\nabla)} \ar{d} \ar{r} &{\ohc \ner(-\bs\nabla)} \ar{d}\\
{\Delta[0]}                \ar{r} &{\ner \bI^\OP} 
\end{tikzcd}
is a homotopy pullback 
(cf. p. \pageref{13.131}).  
The \ composition \ \  lemma \  \ therefore \ \ implies\  \ that \ \  the\  \ commutative \ diagram \quad
\begin{tikzcd}%[sep=large]
{\ner(i\bs\nabla)} \ar{d} \ar{r} &{\ohc \ner(-\bs\nabla)} \ar{d}\\
{\ner(i\bs\bI)}                \ar{r} &{\ohc \ner(-\bs\bI)}
\end{tikzcd}
is a homotopy pullback. \ \  Finally, then, by another application of the composition lemma, \ one concludes that the commutative diagram \quad
\begin{tikzcd}%[sep=large]
{\ner(i\bs\nabla)} \ar{d} \ar{r} &{\ner \bJ} \ar{d}\\
{\ner(i\bs\bI)}                \ar{r} &{\ner \bI} 
\end{tikzcd}
is a homotopy pullback.]

[Note: \  One can also formulate the result in terms of $\nabla/i$.]\\

\label{14.24}
\label{14.57}
\label{18.15}
\textbf{\small LEMMA} \ 
If 
\begin{tikzcd}%[sep=large]
{W} \ar{d} \ar{r} &{Y} \ar{d}{g}\\
{X}  \ar{r}[swap]{f} &{Z} 
\end{tikzcd}
is a homotopy pullback in \bSISET, then 
\begin{tikzcd}%[sep=large]
{\abs{W}} \ar{d} \ar{r} &{\abs{Y}} \ar{d}{\abs{g}}\\
{\abs{X}}  \ar{r}[swap]{\abs{f}} &{\abs{Z}} 
\end{tikzcd}
is a homotopy pullback in \bCGH (singular structure) and the arrow 
$\abs{W} \ra W_{\abs{f},\abs{g}}$ is a homotopy equivalence 
(compactly generated double mapping track).

[In the notation 
p. \pageref{13.132}, write
$Y \overset{\sim}{\ra} \ov{Y} \twoheadrightarrow Z$ 
$-$then $\ov{Y} \ra Z$ Kan $\implies$ $\abs{\ov{Y}} \ra \abs{Z}$ Serre and 
$W \ra X \times_{\ov{Y}} Z$ goes to 
$\abs{W} \ra \abs{X \times_{\ov{Y}} Z} =$ $\abs{X} \times_{\abs{\ov{Y}}} \abs{Z}$ (cf. Proposition 1), so
\begin{tikzcd}%[sep=large]
{\abs{W}} \ar{d} \ar{r} &{\abs{Y}} \ar{d}{\abs{g}}\\
{\abs{X}}  \ar{r}[swap]{\abs{f}} &{\abs{Z}} 
\end{tikzcd}
is a homotopy pullback in \bCGH.  \ 
The  double mapping track of the 2-sink
$\abs{X} \overset{\abs{f}}{\ra} \abs{Z} \overset{\abs{g}}{\la} \abs{Y}$ 
calculated in \bTOP is a CW space (cf. $\S 6$, Proposition 8).  
Its image under $k$ is $W_{\abs{f},\abs{g}}$, thus $W_{\abs{f},\abs{g}}$ is a CW space.  
Therefore the arrow $\abs{W} \ra W_{\abs{f},\abs{g}}$, which is a priori a weak homotopy equivalence, is actually a homotopy equivalence.]
\\

Consequently, under the conditions of Quillen's theorem B, $\forall \ i \in \Ob\bI$: $\nabla^{-1}(i) \neq 0$, there is a homotopy equivalence $B(i\bs\nabla) \ra E_{B\nabla}$ (compactly generated mapping fiber), so $\forall \ j \in \nabla^{-1}(i)$, there is an exact sequence 
\[
\cdots \ra \pi_{q+1}(B\bI,i) \ra \pi_q(B(i\bs\nabla),(j,\id_i)) \ra \pi_q(B\bJ,j) \ra \pi_q(B\bI,i) \ra \cdots.
\]

%%----------------------------------------------------------------------------------------------76
Remark: It is thus a corollary that theorem B $\implies$ theorem A.\\

\begingroup%%----------------------------------->>
\fontsize{9pt}{11pt}\selectfont
Waldhausen\footnote[2]{\textit{CMS Conf. Proc.} \textbf{2} (1982), 141-184.} 
has extended Quillen's theorems A and B from \textbf{CAT} to $[\Delta^\OP,\textbf{CAT}]$.\\
\endgroup%%------------------------------------<<

Fix an abelian group $G$ $-$then a commutative diagram 
\begin{tikzcd}%[sep=large]
{W} \ar{d} \ar{r} &{Y} \ar{d}{g}\\
{X}  \ar{r}[swap]{f} &{Z} 
\end{tikzcd} \ 
of simplicial sets is said to be an 
\un{$HG$-pullback} 
\index{HG-pullback (simplicial sets)} if for some factorization 
$Y \overset{\sim}{\ra} \ov{Y} \twoheadrightarrow Z$ of $g$, the induced simplicial map 
$W \ra X \times_Z \ov{Y}$ is an  $HG$-equivalence.  
Here, the factorization of $g$ is in the usual model category structure on \bSISET and not in that of the homological model category theorem, hence the choice of the factorization of $g$ is immaterial and one can work with either $g$ or $f$.  
Example: A homotopy pullback is an  $HG$-pullback.

[Note: When $G = \Z$, the term is 
\un{homology pullback}.]
\index{homology pullback (simplicial sets)}\\

Example: A commutative diagram
\begin{tikzcd}%[sep=large]
{W} \ar{d} \ar{r} &{Y} \ar{d}{g}\\
{X}  \ar{r}[swap]{f} &{Z} 
\end{tikzcd}
of simplicial sets, where $f$ is a weak homotopy equivalence, is an  $HG$-pullback iff the arrow $W \ra Y$ is an  $HG$-equivalence.\\

\index{Composition Lemma (Simplicial Spaces)}
\textbf{\small COMPOSITION LEMMA} \ 
Consider the commutative diagram
\begin{tikzcd}%[sep=large]
{\bullet}  \ar{d} \ar{r}  &{\bullet} \ar{d} \ar{r}  &{\bullet} \ar{d} \\
{\bullet}  \ar{r} &{\bullet}  \ar{r} &{\bullet} 
\end{tikzcd}
in \bSISET.  Assume: The square on the right is a homotopy pullback $-$then the rectangle is an  $HG$-pullback iff the square on the left is an  $HG$-pullback.\\

Rappel: \bSISET is a topos, so $\forall \ B$, $\bSISET/B$ is a topos 
(MacLane-Moerdijk\footnote[2]{\textit{Sheaves in Geometry and Logic}, Springer Verlag (1992), 190.}), 
thus is cartesian closed.

[Note: \  Similar remarks apply to \bBISISET.]\\

\begin{proposition} \ %55
Let $F:\bI \ra \bSISET$ be a functor.  
Assume: $\forall \ \delta \in \Mor \bI$, $F \delta$ is an  $HG$-equivalence $-$then $\forall \ i \in \Ob \bI$, the pullback square
\begin{tikzcd}[sep=large]
{Fi} \ar{d} \ar{r} &{\ohc F} \ar{d}\\
{\Delta[0]}  \ar{r}[swap]{\Delta_i} &{\ner \bI} 
\end{tikzcd}
is an  $HG$-pullback.
\end{proposition}

%%----------------------------------------------------------------------------------------------77
[Factor 
$\Delta[0] \underset{\Delta_i}{\lra} \ner \bI$ 
as
$\Delta[0] \underset{\Delta_x}{\rat} X \thra \ner \bI$, 
where $\Delta_x$ is a weak homotopy equivalence, the claim being that the arrow
$Fi \ra X \times_{\ner \bI} \ohc F$ 
is an  $HG$-equivalence.  In view of the small object argument, one can suppose that $\Delta_x$ is a sequential colimit of pushouts of coproducts of inclusions $\Lambda[k,n] \ra \Delta[n]$.  Because of this and the fact that the functor 
$- \times_{\ner \bI} \ohc F$
preserves colimits, it is obviously enough to prove that every diagram of the form
\begin{tikzcd}%[sep=large]
&&{\ohc F} \ar{d}\\
{\Lambda[k,n]}  \ar{r} &{\Delta[n]} \ar{r}[swap]{\Delta_f} &{\ner \bI} 
\end{tikzcd}
leads to an  $HG$-equivalence 
$\Lambda[k,n] \times_{\ner \bI}  \ohc F \ra$ 
$\Delta[n] \times_{\ner \bI} \ohc F$.  
To begin with, 
$\Delta[n] \times_{\ner \bI} \ohc F \approx$ $\ohc F \circx f$ $(f:[n] \ra \bI)$.
Furthermore, the initial object $0 \in [n]$ defines a natural transformation $F \circx f(0) \ra F \circx f$, so there is a \cd
\[
\begin{tikzcd}[sep=large]
&&{\ds\coprod\limits_{\substack{\alpha_0 \ra \cdots \ra \alpha_m\\ \in \Lambda[k,n]}}F \circx f(0)} \ar{d} \ar{r} 
&{\ds\coprod\limits_{\substack{\alpha_0 \ra \cdots \ra \alpha_m\\ \in \Delta[n]}}F \circx f(0)}  \ar{d}\\
&&{\ds\coprod\limits_{\substack{\alpha_0 \ra \cdots \ra \alpha_m\\ \in \Lambda[k,n]}}F \circx f(\alpha_0)} \ar{r}
&{\ds\coprod\limits_{\substack{\alpha_0 \ra \cdots \ra \alpha_m\\ \in \Delta[n]}}F \circx f(\alpha_0)}
\end{tikzcd}
\]
of bisimplicial sets. \  The hypothesis on $F$, in conjunction with the appended note to Proposition 51, implies that the diagonal of either vertical arrow is an $HG$-equivalence.  
But the diagonal of the top horizontal arrow is the weak homotopy equivalence 
$\Lambda[k,n] \times F \circx f(0) \ra \Delta[n] \times F \circx f(0)$, therefore the diagonal of the bottom horizontal arrow is an  $HG$-equivalence, i.e., 
$\Lambda[k,n] \times_{\ner \bI} \ohc F \ra \Delta[n] \times_{\ner \bI} \ohc F$  is an $HG$-equivalence.]\\

\begin{proposition} \ %56
Suppose that \bI and \bJ are small categories and $\nabla:\bJ \ra \bI$ is a functor with the property that for every morphism 
$i^\prime \ra i\pp$ in \bI, the arrow $\ner(i\pp\bs\nabla) \ra \ner(i^\prime\bs\nabla)$ is an  $HG$-equivalence $-$then $\forall \ i \in \Ob\bI$, the pullback square
\begin{tikzcd}%[sep=large]
{\ner(i\bs\nabla)} \ar{d} \ar{r} &{\ner \bJ} \ar{d}\\
{\ner(i\bs\bI)}                \ar{r} &{\ner \bI} 
\end{tikzcd}
is an $HG$-pullback.
\end{proposition}

[One has only to trace the proof of Quillen's theorem B, using Proposition 55 to establish that the pullback square
\begin{tikzcd}%[sep=large]
{\ner(i\bs\nabla)} \ar{d} \ar{r} &{\ohc\ner (-\bs\nabla)} \ar{d}\\
{\Delta[0]}                \ar{r} &{\ner \bI^\OP} 
\end{tikzcd}
is an $HG$-pullback.]

[Note: \  It follows that $\forall \ i \in \Ob\bI$: 
$\nabla^{-1}(i) \neq 0$, the arrow $B(i\backslash\nabla) \ra E_{B\nabla}$ is an $HG$-equivalence 
(compactly generated mapping fiber).]\\

%%----------------------------------------------------------------------------------------------78

\begingroup%%----------------------------------->>
\fontsize{9pt}{11pt}\selectfont
Proposition 56 is the homological analog of Quillen's theorem B.  
The same style of argument can also be used for it 
(in Proposition 55, replace ``$HG$-equivalence'' by ``weak homotopy equivalence'' and ``$HG$-pullback'' by ``homotopy pullback'').\\
\endgroup%%------------------------------------<<

Let $(M,O)$ be a category object in \bSISET.  
Suppose that $Y$ is a left \bM-object and $\tran Y$ is the associated translation category $-$then the projection $T:Y \ra O$ gives rise to an internal functor $\tran Y \ra \bM$ from which a morphism  
$\ner\tran Y \ra \ner \bM$ of simplicial objects in $\widehat{\bDelta}$ or still, a bisimplicial map.  
Each $x \in O_0$ determines a pullback square
\begin{tikzcd}%[sep=large]
{Y_x} \ar{d} \ar{r} &{Y} \ar{d}{T}\\
{\Delta[0]}  \ar{r}[swap]{\Delta_x} &{O} 
\end{tikzcd}
in \bSISET and through $e:O \ra M$, arrows $\Delta[0] \underset{\Delta_x}{\ra} \nersub_n \bM$, 
thus there is a pullaback square 
\begin{tikzcd}%[sep=large]
{Y_x} \ar{d} \ar{r} &{\ner \tran Y} \ar{d}\\
{\Delta[0]}  \ar{r}[swap]{\Delta_x} &{\ner \bM} 
\end{tikzcd}
in \bBISISET (abuse of notation).

[Note: \  $\forall \ f \in M_0$, 
$
\begin{cases}
\ sf\\
\ tf
\end{cases}
\in O_0
$
and $\lambda:M \times_O Y \ra Y$ defines an arrow $Y_{sf} \ra Y_{tf}$.]\\
\vspace{0.25cm}

\begin{proposition} \ 
If $\forall \ f \in M_0$, the arrow $Y_{sf} \ra Y_{tf}$ is an $HG$-equivalence, then the pullback square
\begin{tikzcd}%[sep=large]
{Y_x} \ar{d} \ar{r} &{\abs{\ner \  \tran Y}} \ar{d}\\
{\Delta[0]}  \ar{r}[swap]{\abs{\Delta_x}} &{\abs{\ner \  \textbf{M}}} 
\end{tikzcd}
(cf. Proposition 50) is an $HG$-pullback  provided that \mO is a constant simplicial set.
\end{proposition}

[Use the model category structure on \bBISISET furnished by Proposition 52 to factor 
$\Delta[0] \underset{\Delta_x}{\ra} \ner \bM$ as $p \circx i$, where $p$ is a fibration and $i$ is an acyclic cofibration representable as a sequential colimit of pushouts of coproducts of inclusions $\Lambda[k,n,n] \ra \Delta[n,n]$.  
Reasoning as in the proof of Proposition 55, it suffices to show that for any diagram of the form
\begin{tikzcd}%[sep=large]
&&{\ner  \tran Y} \ar{d}\\
{\Lambda[k,n,n]}  \ar{r} &{\Delta[n,n]} \ar{r}[swap]{\Delta_f} &{\ner  \bM} 
\end{tikzcd}
, $\abs{\Lambda[k,n,n]} \times_{\abs{\ner \bM}} \abs{\ner \tran Y} \ra \abs{\Delta[n,n]} 
\times_{\abs{\ner  \bM}} \abs{\ner \tran Y}$
is an $HG$-equivalence.  The arrow
$\Delta_f:\Delta[n,n] \ra \ner \bM$ corresponds to 
$x_0 \overset{f_0}{\ra} x_1 \ra \cdots \ra x_{n-1} \overset{f_{n-1}}{\ra} x_n$,
where the $x_i \in O_n$ $(= O)$ and the $f_i \in M_n$.  This said, consider the \cd
\[
\begin{tikzcd}%[sep=large]
{\Lambda[k,n,n] \times Y_{x_0}} \ar{d} \ar{r} 
&{\Delta[n,n] \times Y_{x_0}} \ar{d}\\
{\Lambda[k,n,n] \times_{\ner \bM} \ner \tran Y}  \ar{r}
&{\Delta[n,n] \times_{\ner \bM} \ner \tran Y}  
\end{tikzcd}
\]
%%----------------------------------------------------------------------------------------------79
which results from piecing together the definitions.  The diagonal of the top horizontal arrow is an $HG$-equivalence 
($\abs{\di\Lambda[k,n,n]}$ is contractible), as is the diagonal of the two vertical arrows.]

[Note: \  Changing the assumption to ``weak homotopy equivalence'' changes the conclusion to ``homotopy pullback''.]\\

\label{13.136}
\label{14.54}
\label{14.55a}
\label{14.63}

\begingroup%%----------------------------------->>
\fontsize{9pt}{11pt}\selectfont
\textbf{\small EXAMPLE} \ 
Let $(M,O)$ be a category object in \bSISET with $O \approx \Delta[0]$.  
So: \mM is a simplicial monoid or, equivalently, \mM is a simplicial object in $\bMON_\bSET$.  Let $Y$ be a left \bM-object.  
Assume: $\forall \ m \in M_0$, 
$m_*:H_*(\abs{Y};G) \ra H_*(\abs{Y};G)$ is an isomorphism $-$then the pullback square
\begin{tikzcd}[sep=large]
{Y} \ar{d} \ar{r} &{\abs{\text{bar}(*;\bM;Y)}} \ar{d}\\
{\Delta[0]}  \ar{r} &{\abs{\text{bar}(*;\bM;*)}} 
\end{tikzcd}
is an $HG$-pullback.\\
\endgroup%%------------------------------------<<

Let
\begin{tikzcd}%[sep=large]
{W} \ar{d} \ar{r} &{Y} \ar{d}\\
{X}  \ar{r} &{Z} 
\end{tikzcd}
be a \cd of bisimplicial sets.  
Problem: Find conditions which ensure that 
\begin{tikzcd}%[sep=large]
{\di W} \ar{d} \ar{r} &{\di Y} \ar{d}\\
{\di X}  \ar{r} &{\di Z} 
\end{tikzcd}
is a homotopy pullback.  
To this end, assume that $\forall \ n$,
\begin{tikzcd}%[sep=large]
{W_n} \ar{d} \ar{r} &{Y_n} \ar{d}\\
{X_n}  \ar{r} &{Z_n} 
\end{tikzcd}
is a homotopy pullback.  
Using the Reedy structure on $[\bDelta^\OP,\bSISET]$, construct a \cd
\begin{tikzcd}%[sep=large]
{W} \ar{d} \ar{r} &{Y} \ar{d} \ar{r} &{\ov{Y}} \ar{d}\\
{X}  \ar{r} &{Z}  \ar{r} &{\ov{Z}}
\end{tikzcd}
,where 
$
\begin{cases}
\ Y \ra \ov{Y}\\
\ Z \ra \ov{Z}
\end{cases}
$
are levelwise weak homotopy equivalences, 
$
\begin{cases}
\ \ov{Y}\\
\ \ov{Z}
\end{cases}
$
are Reedy fibrant, and $\ov{Y} \ra \ov{Z}$ is a Reedy fibration $-$then $\forall \ n$,
\begin{tikzcd}%[sep=large]
{W_n}  \ar{d} \ar{r} &{\ov{Y}_n} \ar{d}\\
{X_n}  \ar{r} &{\ov{Z}_n}
\end{tikzcd}
is a homotopy pullback.  
Form the \cd
\begin{tikzcd}%[sep=large]
{\di W} \ar{d} \ar{r} &{\di Y} \ar{d} \ar{r} &{\di \ov{Y}} \ar{d}\\
{\di X}  \ar{r} &{\di Z}  \ar{r} &{\di \ov{Z}}
\end{tikzcd}
.  
The square 
\begin{tikzcd}%[sep=large]
{\di Y} \ar{d} \ar{r} &{\di \ov{Y}} \ar{d}\\
{\di Z}  \ar{r} &{\di \ov{Z}}
\end{tikzcd}
is a homotopy pullback (cf. Proposition 51), so by the composition lemma, \qquad
\begin{tikzcd}%[sep=large]
{\di W} \ar{d} \ar{r} &{\di Y} \ar{d}\\
{\di X}  \ar{r} &{\di Z}
\end{tikzcd}
will be a homotopy pullback if this is the case of 
$
\begin{tikzcd}%[sep=large]
{\di W} \ar{d} \ar{r} &{\di \ov{Y}} \ar{d}\\
{\di X}  \ar{r} &{\di \ov{Z}}
\end{tikzcd}
.  
$
Since $\forall \ n$, $\ov{Y}_n \ra \ov{Z}_n$ is a Kan fibration (cf. Proposition
%%----------------------------------------------------------------------------------------------80
41), the induced map $W \ra X \times_{\ov{Z}} \ov{Y}$ of bisimplicial sets is a levelwise weak homotopy equivalence, thus
$\di W \ra \di X \times_{\di \ov{Z}} \di \ov{Y}$ is a weak homotopy equivalence (cf Proposition 51).  
Therefore the central issue is whether
$\di \ov{Y} \ra \di \ov{Z}$ is a Kan fibration.  However it is definitely not automatic that di takes Reedy fibrations to Kan fibrations, meaning that conditions have to be imposed.\\

\begingroup%%----------------------------------->>
\fontsize{9pt}{11pt}\selectfont
\textbf{\small EXAMPLE} \ 
Let 
$
\begin{cases}
\ X\\
\ Y
\end{cases}
$
be simplicial sets, $f:X \ra Y$ a simplicial map.  Extend 
$
\begin{cases}
\ X\\
\ Y
\end{cases}
$
to bisimplicial sets by rendering them trivial in the vertical direction $-$then the associated bisimplicial map is a fibration in the Reedy  structure and its diagonal is $f$ but, of course, $f$ need not be Kan.\\
\endgroup%%------------------------------------<<

\begin{proposition} \ %58
Let 
$
\begin{cases}
\ X\\
\ Y
\end{cases}
$
be bisimplicial sets, $f:X \ra Y$ a Reedy fibration.  Assume: $\forall \ m$, the arrow $X_{*,m} \ra Y_{*,m}$ is a Kan fibration $-$then $\di f: \di X \ra \di Y$ is a Kan fibration.
\end{proposition}

[Convert\  the \  lifting \ problem \qquad
\begin{tikzcd}%[sep=large]
\Lambda[k,n] \ar{d} \ar{r} &{\di X} \ar{d}\\
\Delta[n] \ar[dashed]{ru}  \ar{r} &{\di Y}
\end{tikzcd}
\qquad to \ the \ lifting \ problem \qquad
\begin{tikzcd}%[sep=large]
\Lambda[k,n,n] \ar{d} \ar{r} &{X} \ar{d}\\
\Delta[n,n] \ar[dashed]{ru}  \ar{r} &{Y}
\end{tikzcd}
(notation as in the proof of Proposition 52) and factor the inclusion 
$\Lambda[k,n,n] \ra \Delta[n,n]$ 
as
$\Lambda[k,n,n] \ra \Lambda[k,n] \un{\times} \Delta[n] \ra \Delta[n,n]$.  
Since $f$ is Reedy, it has the RLP w.r.t the first inclusion and since $f$ is horizontally Kan, it has the RLP w.r.t the second inclusion.]\\

Let $K$ be a simplicial set.  Given a bisimplicial set $X$, the 
\un{matching space}
\index{matching space} 
of $X$ at $K$ is the simplicial set
$M_KX$ defined by the end  
$\ds\int_{[n]} X_n^{K_n}$. 
So:
$M_KX([m])  \approx$ 
$\Nat(\Delta[m],\ds\int_{[n]} X_n^{K_n})$ $\approx$  
$\ds\int_{[n]} \Nat(\Delta[m], X_n^{K_n})$ $\approx$ 
$\ds\int_{[n]} \Nat(\Delta[m], X_n)^{K_n}$ $\approx$ 
$\ds\int_{[n]}  X_{n,m}^{K_n}$ \ $\approx$ 
$\ds\int_{[n]}  \Mor(K_n,X_{n,m})$ $\approx$ 
$\Nat(K,X_{*,m})$.  
Obviously, $M_KX$ is functorial, covariant in $X$ and contravariant in \mK.

[Note: \  The functor $X \ra M_KX$ is a right adjoint for the functor $L \ra K \un{\times} L$.]

Examples:\  
(1) $M_{\Delta[n]}X([m]) \approx \Nat(\Delta[n],X_{*,m}) \approx X_{n,m}$ $\implies$ 
$M_{\Delta[n]}X \approx X_{*,m}$ $(\equiv X_n)$; 
(2) $M_{\dot\Delta[n]}X([m]) \approx$ $\Nat(\dot\Delta[n],X_{*,m}) \approx$ $\Nat(sk^{(n-1)}\Delta[n],X_{*,m}) \approx$ 
$(cosk^{(n-1)}X)_n$ $\implies$ $M_{\dot\Delta[n]}X \approx M_nX$.

[Note: \  The inclusion $\dot\Delta[n] \ra \Delta[n]$ leads to an arrow
$M_{\Delta[n]}X \ra M_{\dot\Delta[n]}X$ or still, to an arrow $X_n \ra M_nX$, which is precisely the matching morphism.]

One can use an analogous definition for the matching space of $X$ at $K$ if $X$ is a simplicial set rather than a bisimplicial set: 
$M_KX \approx \ds\int_{[n]} X_n^{K_n}$ $(\approx \Nat(K,X))$.

[Note: \  Suppose that $X$ is a bisimplicial set $-$then $M_K X_{*,m} \approx (M_KX)_m$.]

Put $M_{k,n}X = M_{\Lambda[k,n]}X$ $(0 \leq k \leq n, n \geq 1)$.  
Because $\Lambda[k,n] \subset \dot\Delta[n]$, there are arrows $X_n \ra M_nX \ra M_{k,n}X$ natural in X.\\

%%----------------------------------------------------------------------------------------------81
\label{13.135} %dmc mnft
\textbf{\small LEMMA} \ 
A simplicial map $K \ra L$ is a Kan fibration iff the arrows
$K_n \ra M_{k,n}K \times_{M_{k,n}L} L_n$ are surjective $(0 \leq k \leq n, \ n \geq 1)$.

[Note: \   A simplicial map $K \ra L$ is a Kan fibration and a weak homotopy equivalence iff the arrows
$K_n \ra M_nX \times_{M_nL} L_n$ are surjective $(n \geq 0)$.]\\

\begin{proposition} \ %59
Let 
$
\begin{cases}
\ X\\
\ Y
\end{cases}
$
be bisimplicial sets, $f:X \ra Y$ a Reedy fibration.  Suppose that the arrows 
$\pi_0(X_{n,*}) \ra \pi_0(M_{k,n}X \times_{M_{k,n}Y} Y_{n,*})$
arising from the squares
\begin{tikzcd}%[sep=large]
{X_{n,*}} \ar{d} \ar{r} &{Y_{n,*}} \ar{d}\\
{M_{k,n}X}   \ar{r} &{M_{k,n}Y}
\end{tikzcd}
are surjective $(0 \leq k \leq n, n \geq 1)$ $-$then $\di f$ is a Kan fibration.
\end{proposition}

[Since \bSISET satisfies SMC, so does \bBISISET (Reedy structure) 
(cf. p. \pageref{13.133}).  
Applying this to the cofibration
$\Lambda[k,n] \un{\times} \Delta[0] \ra \Delta[n] \un{\times} \Delta[0]$,
it follows that the arrow
$\HOM(\Delta[n] \un{\times} \Delta[0],x) \ra
\HOM(\Lambda[k,n] \un{\times}\Delta[0],X) \times_{\HOM(\Lambda[k,n] \un{\times} \Delta[0],Y)} \HOM(\Delta[n] \un{\times} \Delta[0],Y)$
is a Kan fibration.  Therefore the arrow
$X_{n,*} \ra M_{k,n}X \times_{M_{k,n}Y} Y_{n,*}$ 
is a Kan fibration.  It is surjective by the assumption on $\pi_0$.  
The lemma thus implies that $f$ is horizontally Kan, from which the assertion (cf. Proposition 58).]\\

Convention: The homotopy groups of a pointed simplicial set are those of its geometric realization.\\

\begingroup%%----------------------------------->>
\fontsize{9pt}{11pt}\selectfont
Homotopy groups commute with finite products.  Homotopy groups also commute with infinite products if the data is fibrant but not in general (consider $\pi_1(\bS[1]^\omega)$).\\
\endgroup%%------------------------------------<<

Let $X$ be a bisimplicial set $-$then for every $n, q \geq 1$ and $x \in X_{n,0}$, there are homomorphisms 
$(d_i^h)_*: \pi_q(X_{n,*},x) \ra \pi_q(X_{n-1,*},d_i^hx)$ $(0 \leq i \leq n)$.\\ 
\indent\indent ($\pi_q$) \ \   $X$ satisfies the 
\un{$\pi_q$-Kan condition} 
\index{Kan condition! $\pi_q$-Kan condition} at $x \in X_{n,0}$ if for every finite sequence 
$(\alpha_0, \ldots, \widehat{\alpha_k}, \ldots, \alpha_n)$, where $\alpha_i \in \pi_q(X_{n-1,*},d_i^hx)$ and 
$(d_i^h)_*\alpha_j = (d_{j-1}^h)_*\alpha_i$ $(i < j \ \& \ i, j \neq k)$, $\exists \ \alpha \in \pi_q(X_{n,*},x): (d_i^h)_*\alpha = \alpha_i$ $(i \neq k)$.

[Note: \  If $x^\prime, x\pp \in X_{n,0}$ are in the same component of $X_n$, then $X$ satisfies the $\pi_q$-Kan condition at $x^\prime$ iff $X$ satisfies the $\pi_q$-Kan condition at $x\pp$.]

\label{13.137}
Definition: A bisimplicial set $X$ satisfies the 
\un{$\pi_*$-Kan condition} 
\index{Kan condition! $\pi_*$-Kan condition} if $\forall \ n,q \geq 1$, $X$ satisfies the 
$\pi_q$-Kan condition at each $x \in X_{n,0}$.

Example: Bisimplicial groups satisfy the $\pi_*$-Kan condition.\\

\begingroup%%----------------------------------->>
\fontsize{9pt}{11pt}\selectfont
\textbf{\small EXAMPLE} \ 
Let $X$ be a bisimplicial set such that $\forall \ n$, $X_n$ is connected $-$then $X$ satisfies the $\pi_*$-Kan condition.
\vspi
[Consider the $\pi_q$-Kan condition at $x = s_{n-1}^h \cdots s_0^h x_0 \ (x_0 \in X_{0,0})$.]\\
\endgroup%%------------------------------------<<

%%----------------------------------------------------------------------------------------------82
\textbf{\small LEMMA} \ 
Let 
$
\begin{cases}
\ X\\
\ Y
\end{cases}
$
be bisimplicial sets, $f:X \ra Y$ a bisimplicial map.  Assume: $f$ is a levelwise weak homotopy equivalence $-$then $X$ satisfies the $\pi_*$-Kan condition iff $Y$ satisfies the $\pi_*$-Kan condition.\\

One can describe $\dot\Delta[n]$ as the simplicial subset of $\Delta[n]$ generated by the $d_i\id_{[n]}$ $(0 \leq i \leq n)$ and one can describe the $\Lambda[k,n]$ as the simplicial subset of $\Delta[n]$ generated by the $d_i\id_{[n]}$ $(0 \leq i \leq n, i \neq k)$.  In general, if $t_0, \ldots, t_r$ are integers such that $0 \leq t_0 < \cdots < t_r \leq n$, let 
$\Delta_n^{(t_0, \ldots, t_r)}$ be  the simplicial subset of $\Delta[n]$ generated by the 
$d_{t_0}\id_{[n]}, \ldots, d_{t_r}\id_{[n]}$ $-$then there is a pushout square
\[
\begin{tikzcd}%[sep=large]
{\Delta_{n-1}^{(t_0, \ldots, t_{r-1)}}} \ar{d} \ar{rr}{\Delta[\delta_{t_{r-1}}^n]} 
&&{\Delta_n^{(t_0, \ldots, t_{r-1})}} \ar{d}\\
{\Delta[n-1]}   \ar{rr}[swap]{\Delta[\delta_{t_r}^n]} 
&&{\Delta_n^{(t_0, \ldots, t_r)}}
\end{tikzcd}
.
\]

[Note: \  $\Delta_n^{(t_0, \ldots, t_r)}$ is a simplicial subset of $\Lambda[k,n]$ provided that $k \neq t_i$ $(i = 0, \ldots, r)$.]

Given a bisimplicial set $X$, write $M_n^{(t_0, \ldots, t_r)}X$ for the matching space of $X$ at $\Delta_n^{(t_0, \ldots, t_r)}$.  There are arrows
$X_n \ra M_nX \ra M_n^{(t_0, \ldots, t_r)}X$ natural in $X$.  
Example: $M_n^{(0, \ldots, , \widehat{k}, \ldots, n)}X = M_{k,n}X$.

[Note: \  $M_n^{(t_0, \ldots, t_r)}X([m])$ consists of the set of finite sequences $(x_{t_0}, \ldots, x_{t_r})$ of elements of $X_{n-1,m}$ such that $d_i^hx_j = d_{j-1}^hx_i$ for all $i < j$ in $\{t_0, \ldots, t_r\}$ 
(cf. p. \pageref{13.134}).  
Moreover, the arrow $X_n \ra M_n^{(t_0, \ldots, t_r)}X$ sends $x \in X_{n,m}$ to $(d_{t_0}^hx, \ldots, d_{t_r}^hx)$ and it is Kan if $X$ is Reedy fibrant.]\\

\textbf{\small LEMMA} \ 
Let $X$ be a bisimplicial set.  Assume: $X$ is Reedy fibrant and satisfies the $\pi_*$-Kan condition.  Suppose that
$x = (x_{t_0}, \ldots, x_{t_r}) \in M_n^{(t_0, \ldots, t_r)}X([0])$ $-$then $\forall \ q \geq 1$, the map
$\pi_q(M_n^{(t_0, \ldots, t_r)}X,x) \ra \pi_q(X_{n-1,*},x_{t_0}) \times \cdots \times \pi_q(X_{n-1,*},x_{t_r})$ 
is injective and its range is the set of finite sequences $(\alpha_{t_0}, \ldots, \alpha_{t_r})$ in the product such that
$(d_i^h)_*\alpha_j = (d_{j-1}^h)_*\alpha_i$ for all $i < j$ in $\{t_0, \ldots, t_r\}$.

[Work inductively with the pullback squares

\[
\begin{tikzcd}%[sep=large]
{M_n^{(t_0, \ldots, t_r)}X} \ar{d} \ar{rr} &&{X_{n-1,*}} \ar{d}\\
{M_n^{(t_0, \ldots, t_{r-1})}X}   \ar{rr} &&{M_{n-1}^{(t_0, \ldots, t_{r-1})}X}
\end{tikzcd}
.]
\]

[Note: \  The result also holds for $q = 0$.]\\

Given a bisimplicial set $X$, define a simplicial set $\pi_0(X)$ by $\pi_0(X)_n = \pi_0(X_n)$ $(= \pi_0(X_{n,*}))$.  
Example: Suppose that $X$ is Reedy fibrant and satisfies the $\pi_*$-Kan condition $-$then 
$\pi_0(M_{k,n}X) \approx M_{k,n}\pi_0(X)$.\\
%%----------------------------------------------------------------------------------------------83

\begingroup%%----------------------------------->>
\fontsize{9pt}{11pt}\selectfont
\textbf{\small EXAMPLE} \ 
Let $X$ be a bisimplicial set such that $\forall \ n$, the path components of $\abs{X_n}$ are abelian.  Write 
$[\bS^q,X]$ for the simplicial set with $[\bS^q,X]_n = [\bS^q,\abs{X_n}]$ $-$then $X$  satisfies the $\pi_*$-Kan condition if the simplicial map $[\bS^q,X] \ra \pi_0(X)$ is a Kan fibration $\forall \ q \geq 1$.\\
\endgroup%%------------------------------------<<

\label{13.138}
\label{13.139}
\begingroup%%----------------------------------->>
\fontsize{9pt}{11pt}\selectfont
\textbf{\small FACT} \ 
Let $X$ be a bisimplicial set such that $\forall \ n$ $X_n$ is connected $-$then $\di X$ is connected.
\vspi
[There is a coequalizer diagram 
%$\pi_0(X_1)  \underset{\underset{d_0}{\lra}}{\overset{d_1}{\lra}}  \pi_0(X_0) \ra \pi_0(diX)$.]\\
\begin{tikzcd}[sep=small]
\pi_0(X_1) 
%\arrow[r,shift right=1]{d_0} \arrow[r,shift right=-1] 
\ar[shift right=-1]{r}{d_1} 
\ar[shift right=1]{r}[swap]{d_0}
&{\pi_0(X_0)} \ar{r}
&{\pi_0(\di X)}
\end{tikzcd}
.]\\
\endgroup%%------------------------------------<<

\begin{proposition} \ %60
Let 
$
\begin{cases}
\ X\\
\ Y
\end{cases}
$
be bisimplicial sets, $f:X \ra Y$ a Reedy fibration with $f_*:\pi_0(X) \ra \pi_0(Y)$ a Kan fibration.  Assume: 
$
\begin{cases}
\ X\\
\ Y
\end{cases}
$
are Reedy fibrant and satisfy the $\pi_*$-Kan condition $-$then $\di f$ is a Kan fibration.
\end{proposition}

[According to Proposition 59, \ it suffices to show that the arrows 
$\pi_0(X_{n,*}) \ \ra$ 
$\pi_0(M_{k,n}X$ $\times_{M_{k,n}Y} Y_{n,*})$ \ are \ surjective \ $(0 \leq k \leq n, n \geq 1)$. 
Consider the square 
\begin{tikzcd}%[sep=large]
{X_{n,*}} \ar{d} \ar{r} &{Y_{n,*}} \ar{d}\\
{M_{k,n}X}  \ar{r} &{M_{k,n}Y} 
\end{tikzcd}
$-$then 
$\pi_0(M_{k,n}X \times_{M_{k,n}Y} Y_{n,*}) \approx \pi_0(M_{k,n}X) \times_{\pi_0(M_{k,n}Y)} \pi_0(Y_{n,*})$.  
In fact, $Y_{n,*} \ra M_{k,n}Y$ is a Kan fibration and the lemma implies that $\forall \ y \in Y_{n.0}$, $Y_{n,*} \ra M_{k,n}Y$ induces a surjection of fundamental groups (cf. infra).  But 
$\pi_0(M_{k,n}X \times_{\pi_0(M_{k,n}Y)} \pi_0(Y_{n,*})$ 
$\approx$ 
$M_{k,n}\pi_0(X) \times_{M_{k,n}\pi_0(Y)} \pi_0(Y_{n,*})$ 
and 
$\pi_0(X_{n,*}) \ra  M_{k,n}\pi_0(X) \times_{M_{k,n}\pi_0(Y)} \pi_0(Y_{n,*})$ 
is surjective, $\pi_0(X) \ra \pi_0(Y)$ being Kan by assumption 
(cf. p. \pageref{13.135}).]\\

\begingroup%%----------------------------------->>
\fontsize{9pt}{11pt}\selectfont
\textbf{\small LEMMA} \ 
Let 
\begin{tikzcd}[sep=large]
{X^\prime} \ar{d} \ar{r} &{X} \ar{d}{p}\\
{B^\prime}  \ar{r} &{B} 
\end{tikzcd}
be a pullback square of topological spaces, where $p:X \ra B$ is a Serre fibration.  Assume: $\forall \ x \in X$, the homomorphism 
$\pi_1(X,x) \ra \pi_1(B,p(x))$ is surjective $-$then the arrow 
$\pi_0(X^\prime) \ra \pi_0(B^\prime) \times_{\pi_0(B)} \pi_0(X)$ is bijective.
\vspi
[Injectivity is a consequence of the $\pi_1$-hypothesis.]\\
\endgroup%%------------------------------------<<

\index{Theorem Bousfield-Friedlander}
\textbf{\small THEOREM OF BOUSFIELD-FRIEDLANDER} \ 
Let
\begin{tikzcd}%[sep=large]
{W} \ar{d} \ar{r} &{Y} \ar{d}\\
{X}  \ar{r} &{Z} 
\end{tikzcd}
be a commutative diagram of bisimplicial sets such that $\forall \ n$, 
\begin{tikzcd}%[sep=large]
{W_n} \ar{d} \ar{r} &{Y_n} \ar{d}\\
{X_n}  \ar{r} &{Z_n} 
\end{tikzcd}
is a homotopy pullback.  Assume: $\pi_0(Y) \ra \pi_0(Z)$ is a Kan fibration and $Y$, $Z$ satisfy the $\pi_*$-Kan condition $-$then
\begin{tikzcd}%[sep=large]
{\di W} \ar{d} \ar{r} &{\di Y} \ar{d}\\
{\di X}  \ar{r} &{\di Z} 
\end{tikzcd}
is a homtopoy pullback.

[Proceed as on 
p. \pageref{13.136} ff.: 
$\di \ov{Y} \ra \di \ov{Z}$ is a Kan fibration (cf. Proposition 60).]

%%----------------------------------------------------------------------------------------------84
[Note: \  When $Y_n$, $Z_n$ are connected $\forall \ n$, $\pi_0(Y) \ra \pi_0(Z)$ is trivially Kan and $Y$, $Z$  necessarily satisfy the $\pi_*$-Kan condition 
(cf. p. \pageref{13.137}).]\\

\begingroup%%----------------------------------->>
\fontsize{9pt}{11pt}\selectfont
Let $K$ be a simplicial set.  
Given a bisimplicial set $X$, define a bisimplicial set $\map(K,X)$ by 
$\map(K,X)_n = \map(K,X_n)$.\\
\endgroup%%------------------------------------<<

\begingroup%%----------------------------------->>
\fontsize{9pt}{11pt}\selectfont
\textbf{\small LEMMA} \ 
There is a canonical arrow $\abs{\map(K,X)} \ra \map(K, \abs{X})$.
\vspi
[The evaluation $K \times \map(K,X_n) \ra X_n$ defines a bisimplicial map $K \times \map(K,X) \ra X$ or still, a simplicial map $\abs{K \times \map(K,X)} \ra \abs{X}$.  However multiplication by $K$ in \bBISISET is a left adjoint, hence 
$\abs{K \times \map(K,X)} \approx K \times \abs{\map(K,X)}$.]\\
\endgroup%%------------------------------------<<

\begingroup%%----------------------------------->>
\fontsize{9pt}{11pt}\selectfont
A bisimplicial set $X$ is said to be 
\un{pointed} 
\index{pointed (bisimplicial set)} 
if an $x \in X_{0,0}$ has been fixed and each $X_n$ is equipped with the base point $s_{n-1}^h \cdots s_0^h x_0$.\\
\endgroup%%------------------------------------<<

\begingroup%%----------------------------------->>
\fontsize{9pt}{11pt}\selectfont
\textbf{\small EXAMPLE} \  Let $X$ be a Reedy fibrant pointed bisimplicial set such that $\forall \ n$, 
$X_n$ is connected $-$then $X$ is $\pi_*$-Kan, thus $\abs{X} \approx \di X$ is fibrant (cf. Proposition 60).  
Denote by $\Theta X$ $(\Omega X)$ the bisimplicial set which takes [n] to $\Theta X_n$ $(\Omega X_n)$ (it follows from Proposition 41 that $\forall \ n$, $X_n$ is fibrant).  
Specializing the lemma to $K = \Delta[1]$ provides us with the canonical arrows 
$\abs{\Theta X} \ra \Theta \abs{X}$ $(\abs{\Omega X} \ra \Omega \abs{X})$ ($\abs{?}$ preserves pullbacks) and a \cd
\begin{tikzcd}[sep=large]
{\abs{\Omega X}} \ar{d} \ar{r} &{\abs{\Theta X}} \ar{d} \ar{r} 
&{\abs{X}} \arrow[d,shift right=0.5,dash] \arrow[d,shift right=-0.5,dash]\\
{\Omega \abs{X}}  \ar{r} &{\Theta \abs{X}} \ar{r} &{\abs{X}}
\end{tikzcd}
.  On the other hand, the theorem of Bousfield-Friedlander says that
\begin{tikzcd}[sep=large]
{\abs{\Omega X}} \ar{d} \ar{r} &{\abs{\Theta X}} \ar{d}\\
{\Delta[0]}  \ar{r}  &{\abs{X}}
\end{tikzcd}
is a homotopy pullback.  Because the geometric realization of $\abs{\Theta X}$ is contractible (cf. Proposition 51), the conclusion is that the canonical arrow $\abs{\Omega X} \ra \Omega \abs{X}$ is a weak homotopy equivalence.\\
\endgroup%%------------------------------------<<

\begingroup%%----------------------------------->>
\fontsize{9pt}{11pt}\selectfont
\textbf{\small EXAMPLE} \  Let $X$ be a pointed bisimplicial set such that $\forall \ n$, $\abs{X_n}$ is simply connected $-$then $\abs{\di X}$ is simply connected.
\vspi
[For this, one can suppose that $X$ is Reedy fibrant.  On general grounds, $\abs{\di X}$ is path connected 
(cf. p. \pageref{13.138}) 
and by the preceding example, $\pi_0(\di \Omega X) \approx \pi_0(\Omega \di X)$.  But $\forall \ n$, $\Omega X_n$ is connected, thus $\di \Omega X$ is connected 
(cf. p. \pageref{13.139}) 
and so $\abs{\di X}$ is simply connected.
\vspi
[Note: \  It is clear that the argument can be iterated: $\abs{X_n}$ $k$-conneted $\forall \ n$ $\implies$ 
$\abs{\di X}$ $k$-connected.]\\
\endgroup%%------------------------------------<<

%%%%%%%%%%%%%%%%%%%%%%%%%%%%%%%%%%%%%%
%%%%%%%%%%%%%%%%%%%%%%%%%%%%%%%%%%%%%%
%%%%%%%%%%%%%%%%%%%%%%%%%%%%%%%%%%%%%%

\begin{center}
$\S \ 13$
\\[0.5cm]
$\mathcal{REFERENCES}$\\
\end{center}

\[
\textbf{BOOKS}
\]

\begingroup
\fontsize{9pt}{11pt}\selectfont
\setlength\parindent{0 cm}

[1] \quad Dwyer, Hirschhorn, Kan, Smith, 
\textit{Homotopy Limit Functors on Model Categories and Homotopical}

\hspace{0.8cm}\textit{Categories}, Mathematical Surveys and Monographs, Amer. Math Soc. 113 (2004).
\\[-.2cm]

[2] \quad Fritsch, R. and Piccinini, R., 
\textit{Cellular Structures in Topology}, Cambridge University Press (1990).
\\[-.2cm]

[3] \quad Gabriel, P. and Zisman, M., 
\textit{Calculus of Fractions and Homotopy Theory}, Springer Verlag (1967).
\\[-.2cm]

[4] \quad Goerss, P. and Jardine, J., \textit{Simplicial Homotopy Theory}, Progress in Mathematics, volume 174, 

\hspace{0.8cm}Birkh\"auser, (1999).
\\[-.2cm]

[5] \quad Heggie, M., 
\textit{Tensor Products in Homotopy Theory}, Ph.D. Thesis, McGill University, Montreal (1986).
\\[-.2cm]

%[5] \quad Hirschhorn, P., Localization, Cellularization, and Homotopy Colimits,\\ Orig version
[6] \quad Hirschhorn, P., 
\textit{Model Categories and Their Localizations}, Amer. Math. Soc. (2003).
\\[-.2cm]

[7] \quad Lamotke, K., 
\textit{Semisimpliziale Algebraische Topologie}, Springer Verlag (1968).
\\[-.2cm]

[8] \quad May, J., 
\textit{Simplicial Objects in Algebraic Topology}, Van Nostrand (1967).
\\[-.2cm]

[9] \quad Quillen, D., 
\textit{Homotopical Algebra}, Springer Verlag (1967).
\\[-.2cm]

[10] \quad Reedy, C., 
\textit{Homology of Algebraic Theories}, Ph.D. Thesis, University of California, San Diego (1974).
\\[-.2cm]

[11] \quad Thomason, R., 
\textit{Homotopy Colimits in \bCAT with Applications to Algebraic K-Theory and Loop Space}

\hspace{0.95cm}\textit{Theory}, Ph.D. Thesis, Princeton Universtiy, Princeton (1977).
\\
\endgroup

\[
\textbf{ARTICLES}
\]

\begingroup
\fontsize{9pt}{11pt}\selectfont
\setlength\parindent{0 cm}

[1] \quad Anderson, D., Fibrations and Geometric Realizations, 
\textit{Bull. Amer. Math. Soc.} \textbf{84} (1978), 765-788.
\\[-.2cm]

[2] \quad Bousfield, A. and Friedlander, E., Homotopy Theory of $\bGamma$-Spaces, Spectra, and Bisimplicial Sets, 
\textit{SLN} 

\hspace{0.8cm}\textbf{658} (1978), 80-130.
\\[-.2cm]

[3] \quad Curtis, E., Simplicial Homotopy Theory, 
\textit{Adv. Math.} \textbf{6} (1971), 107-209.
\\[-.2cm]

[4] \quad Duskin, J., Simplicial Methods and the Interpretation of "Triple" Cohomology, 
\textit{Memoirs Amer. Math.}

\hspace{0.8cm}\textit{Soc.} \textbf{163} (1975), 1-135.
\\[-.2cm]

[5] \quad Dwyer, W., Hopkins, M., and Kan, D., The Homotopy Theory of Cyclic Sets, 
\textit{Trans. Amer. Math.}

\hspace{0.8cm}\textit{Soc.} \textbf{291} (1985), 281-289.
\\[-.2cm]

[6] \quad Dwyer, W. and Kan, D., Equivalences Between Homotopy Theories of Diagrams, In: 
\textit{Algebraic Topol-}

\hspace{0.8cm}\textit{ogy and Algebraic K-Theory}, W. Browder (ed.), Princeton University Press (1987), 180-205.
\\[-.2cm]

[7] \quad Fritsch, R., Zur Unterteilung Semisimplizialer Mengen I and II, 
\textit{Math. Zeit.} \textbf{108} (1969), 329-367 and 

\hspace{0.8cm}\textbf{109} (1969), 131-152.
\\[-.2cm]

[8] \quad Fritsch, R. and Puppe, D., Die Hom\"oomorphie der Geometrischen Realisierungen einer Semisim

\hspace{0.8cm}plizialen Menge und ihrer Normalunterteilung, 
\textit{Arch. Math.} \textbf{18} (1967), 508-512.
\\[-.2cm]

%this is now a book and should be moved accordingly - this messes up infra
% moved to books [9] \quad Goerss, P. and Jardine, J., Simplicial Homotopy Theory.\\
%Birkh\"auser (1999).\\

[9] \quad Goerss, P. and Jardine, J., Localization Theories for Simplicial Presheaves, \textit{Can. J. Math.} \textbf{50} (1998), 

\hspace{0.8cm}1048-1089.
\\[-.2cm]

[10] \quad Gray, J., Closed Categories, Lax Limits and Homotopy Limits, 
\textit{J. Pure Appl. Algebra} \textbf{19} (1980), 

\hspace{0.95cm}127-158.
\\[-.2cm]

[11] \quad Heller, A., Homotopy Theories, 
\textit{Memoirs Amer. Math. Soc.} \textbf{383} (1988), 1-78.
\\[-.2cm]

[12] \quad Jardine, J., The Homotopical Foundations of Algebraic K-Theory, 
\textit{Contemp. Math.} \textbf{83} (1989), 57-82.
\\[-.2cm]
 
[13] \quad Joyal, A. and Tierney, M., Classifying Spaces for Sheaves of Simplicial Groupoids, 
\textit{J. Pure Appl.}

\hspace{0.95cm}\textit{Algebra} \textbf{89} (1993), 135-161.
\\[-.2cm]

[14] \quad Joyal, A. and Tierney, M., On the Homotopy Theory of Sheaves of Simplicial Groupoids, 
\textit{Math. Proc.}

\hspace{0.95cm}\textit{Cambridge Philos. Soc.} \textbf{120} (1996), 263-290.
\\[-.2cm]

[15] \quad Kan, D., On c.s.s. Complexes, 
\textit{Amer. J. Math.} \textbf{79} (1957), 449-476.
\\[-.2cm]

[16] \quad Milnor, J., The Geometric Realization of a Semisimplicial Complex, 
\textit{Ann. of Math.} \textbf{65} (1957), 

\hspace{0.95cm}357-362.
\\[-.2cm]

[17] \quad Moerdijk, I., Bisimplicial Sets and the Group Completion Theorem, In: 
\textit{Algebraic K-Theory: Con-}

\hspace{0.95cm}\textit{nections with Geometry and Topology}, J. Jardine and V Snaith (ed.), Kluwer Academic Publishers 

\hspace{0.95cm}(1989), 225-240.
\\[-.2cm]

[18] \quad Quillen, D., Higher Algebraic K-Theory: I, 
\textit{SLN} \textbf{341} (1973), 85-147.
\\[-.2cm]

[19] \quad Reedy, C., Homotopy Theory of Model Categories, \textit{Preprint}.
\\[-.2cm]

[20] \quad Schwede, S., Stable Homotopy of Algebraic Theories, \textit{Topology} \textbf{40} (2001), 1-41.
\\[-.2cm]

[21] \quad Zisman, M., Quelques Propri\'etes des Fibr\'es au Sens de Kan, 
\textit{Ann. Inst. Fourier (Grenoble)} \textbf{10} 

\hspace{0.95cm}(1960), 345-457.

\setlength\parindent{2em}

\endgroup

\chapter{
$\boldsymbol{\S}$\textbf{14}.\quadx  SIMPLICIAL SPACES}
\setlength\parindent{2em}
\setcounter{proposition}{0}
\setcounter{chapter}{14}

%%----------------------------------------------------------------------------------------------01
$\text{ }$\\[-1.25cm]

After working through the foundations of the theory, various applications will be given, e.g., the James construction and infinite symmetric products.  I have also included some material on operads and delooping procedures.

A 
\un{simplicial space} 
\index{simplicial space} 
is a simplicial object in \bTOP and a 
\un{simplicial map} 
\index{simplicial map} 
is a morphism of simplicial spaces.  
\bTOP, in its standard structure, is a model category, thus \bSITOP is a model category (Reedy structure) 
(cf. p. \pageref{14.1}).  
This fact notwithstanding, it will be simplest to proceed from first principles.

There is a forgetful functor $\bSITOP \ra \bSISET$ and it has a left and right adjoint 
(cf. p \pageref{14.2}).

[Note: \ The purely set theoretic properties of simplicial spaces are the same as  those of simplicial sets.]

Given an $X$ in \bSITOP, put 
$\abs{X} = \ds\int^{[n]} X_n \times \dn$ $-$then $\abs{X}$ is the 
\un{geometric realization} 
\index{geometric realization (in \bSITOP)} 
of $X$ and the assignment $X \ra \abs{X}$ is a functor 
$\bSITOP \ra \bTOP$.  $\abs{?}$ has a right adjoint $\bTOP \ra \bSITOP$ (compact open topology on the singular set).\\

\label{14.38}
\begingroup%%----------------------------------->>
\fontsize{9pt}{11pt}\selectfont
\index{Star Construction (example)}
\textbf{\small EXAMPLE \  (\un{Star Construction})} \ 
Let $X$ be a nonempty topological space.  Define a simplicial space $\Lambda X$ by the prescription 
$(\Lambda X)_n = X \times \cdots \times X$ (n+1 factors) with 
$d_i(x_0, \ldots x_n) = (x_0, \ldots, \hat{x_i}, \ldots ,x_n)$, 
$s_i(x_0, \ldots x_n) = (x_0, \ldots, x_i,x_i, \ldots ,x_n)$.  
Represent $\dpn$ as the set of points 
$(t_1, \ldots, t_n)$ in $\R^n$ such that $0 \leq t_1 \leq \cdots \leq t_n \leq 1$ (which entails a change in the formulas defining the simplicial operators).  
Form $X^*$ as on 
p. \pageref{14.3} 
and let 
$\lambda_n:X^{n+1} \times \dpn \ra X^*$ be the continuous function that sends 
$((x_0, \ldots, x_n),(t_1, \ldots, t_n))$ 
to the right continuous step function $[0,1[ \ra X$ which is equal to $x_i$ on $[t_i,t_{i+1}[$ $(t_0 = 0, t_{n+1} = 1)$ 
$-$then
the $\lambda_n$ combine to give a continuous bijection 
$\lambda: \abs{\Lambda X} \ra X^*$.  
Since $X$ $\tT_2$ $\implies$ $X^*$ $\tT_2$, $\abs{\Lambda X}$ is Hausdorff whenever $X$ is and in this situation, the composite
$X \ra X \times \Delta^0 \ra  \abs{\Lambda X}$ is a closed embedding.
\vspi
[Note: \  Like $X^*$, $\abs{\Lambda X}$ is contractible 
(cf. p. \pageref{14.3a}).]\\
\endgroup %%------------------------------------<<

\index{Hausdorff (simplicial space)} \index{compactly generated (simplicial space)}
A simplicial space $X$ is said to be 
\un{Hausdorff}, 
\un{compactly generated} 
$\ldots$ if $\forall \ n$, $X_n$ is Hausdorff, compactly generated $\ldots$, i.e., if $X$ is a simplicial object in \bHAUS, \bCG, $\ldots$.  
On general grounds, the geometric realization of a compactly generated simplicial space is automatically compactly generated but there is no a priori guarantee that the geometric realization of a Hausdorff simplicial space is Hausdorff.

Observation: If $X$ is a simplicial space and if $\alpha:[m] \ra [n]$ is an epimorphism, then 
$X\alpha:X_n \ra X_m$ is an embedding and $(X\alpha)X_n$ is a retract of $X_m$.

%%----------------------------------------------------------------------------------------------02
Let $X$ be a simplicial space $-$then $X$ is said to satisfy the 
\un{embedding condition} 
\index{embedding condition (simplicial space)} 
if $\forall \ n$ $\&$ $\forall \ i$, $s_i:X_{n-1} \ra X_n$ is a closed embedding.  
Examples: 
(1) A Hausdorff simplicial space  satisfies the embedding condition;  
(2) A $\Delta$-separated compactly generated simplicial space satisfies the embedding condition.
\\

\textbf{\small LEMMA} \ 
Suppose given a diagram
\begin{tikzcd}%[sep=small]
{X^\prime} \ar{d}[swap]{p^\prime}  &{X} \ar{d}{p}\\
{B^\prime} \ar{r}[swap]{i}  &{B}
\end{tikzcd}
of topological spaces and continuous functions, where $p$ is quotient and $i$ is one-to-one.  
Assume: $\exists$ a neighborhood finite collection $\{A_j\}$ of closed subsets of $X$ and continuous functions 
$f_j:A_j \ra X^\prime$ such that 
$p^{-1}(i(B^\prime)) \ = \  \bigcup\limits_j A_j$ 
with \ 
$\restr{p}{A_j} = i \circx p^\prime \circx f_j$ $\forall \ j$ $-$then $p^\prime$ is quotient and $i$ is a closed embedding.\\

If \mX is a simplicial space, then $\aX$ can be identified with the quotient 
$\coprod\limits_n X_n \times \dpn / \sim$, the equivalence relation being generated by writing 
$((X\alpha)x,t) \sim (x,\Delta^\alpha t)$.  Let $p:\coprod\limits_n X_n \times \dpn \ra \aX$ be the projection and put 
$\aX_n = p\bigl(\coprod\limits_{m \leq n} X_m \times \Delta^m\bigr)$.\\

\begin{proposition} \ %01
Let $X$ be a simplicial space.  
Assume: $X$ satisifies the embedding condition $-$then $\forall \ n$, $\aX_n$ is a closed subspace of $\aX$ and 
$\aX = \colimx \aX_n$.
\end{proposition}

[Fix $n^\prime$ and consider
$
\begin{tikzcd}%[sep=large]
{\coprod\limits_{m \leq n^\prime} X_m \times \Delta^m} \ar{d}[swap]{p^\prime}  
&{\coprod\limits_n X_n \times \dpn} \ar{d}{p}\\
{\aX_{n^\prime}} \ar{r}[swap]{i}  &{\aX}
\end{tikzcd}
. \ 
$
For each $m \leq n^\prime$ and $n$, there are but finitely many diagrams of the form 
$[m] \overset{\beta}{\lla} [k] \overset{\alpha}{\lra}[n]$, where $\alpha$ is a monomorphism and $\beta$ is an epimorphism.  
Put 
$A_{\alpha, \beta} = (X\alpha)^{-1}(X \beta) X_m \times \Delta^\alpha \Delta^k \subset$ $X_n \times \Delta^n$, 
define 
$f_{\alpha,\beta}:A_{\alpha, \beta} \ra X_m \times \Delta^m$ 
by 
$f_{\alpha,\beta}(x,t) = (y,\Delta^\beta u)$ $(t = \Delta^\alpha u$ $(\exists! \ u \in \Delta^k)$, 
$(X\alpha)x = (X\beta)y$ $(\exists! \ y \in X_m))$, and apply the lemma.]\\

\label{14.40}
\begingroup%%----------------------------------->>
\fontsize{9pt}{11pt}\selectfont
\textbf{\small FACT} \ 
Suppose that $X$ is a simplicial space satisfying the embedding condition.  Define a simplicial set $\pi_0(X)$ by 
$\pi_0(X)_n = \pi_0(X_n)$ $-$then $\pi_0(\aX) \approx \pi_0\abs{\pi_0(X)}$.
\vspi
[Every point in $\aX$ can be joined by a path in $\aX$ to a point in $X_0 = \aX_0$.  On the other hand, given $x \in X_1$, $\sigma(t) = [x,(1-t,t)]$ $(0 \leq t \leq 1)$ is a path in $\aX$ which begins at $d_1x$ and ends at $d_0x$.]
\vspi
[Note: \  Therefore $\aX$ is path connected if $X_0$ is path connected.]\\
\endgroup %%------------------------------------<<

Notation: Given an $X$ in \bSITOP, write $sX_{n-1}$ for the union 
$s_0X_{n-1} \cup \cdots \cup s_{n-1}X_{n-1}$.\\

%%----------------------------------------------------------------------------------------------03
\begin{proposition} \ %02
Let $X$ be a simplicial space.  
Assume: $X$ satisfies the embedding condition $-$then $\forall \ n$, there is a pushout square
$
\begin{tikzcd}%[sep=small]
{X_n \times \ddpn \cup sX_{n-1} \times \dpn} \ar{d}  \ar{r} &{\aX_{n-1}} \ar{d}\\
{X_n \times \dpn} \ar{r}  &{\aX_n}
\end{tikzcd}
.
$

[The arrow $X_n \times \dpn \ra \aX_n$ is quotient.  
To see this, form
$
\begin{tikzcd}%[sep=small]
{X_n \times \dpn} \ar{d}  &{\coprod\limits_{m \leq n} X_m \times \Delta^m} \ar{d}{\quad \cdot}\\
{\aX_n} \arrow[r,shift right=0.5,dash] \arrow[r,shift right=-0.5,dash]  &{\aX_n}
\end{tikzcd}
%.  
$
Taking into account the lemma, let 
$f_n:X_n \times \dpn \ra X_n \times \dpn$ be the identity.  
To define 
$f_m:X_m \times \Delta^m \ra X_n \times \dpn$ if $m < n$, fix 
a monomorphism $\alpha:[m] \ra [n]$, 
an epimorphism $\beta:[n] \ra [m]$ such that $\beta \circx \alpha = \id_{[m]}$, and put 
$f_m(x,t) = ((X\beta)x, \Delta^\alpha t)$.]\\
\end{proposition}

Application: Suppose that $X$ is a $\Delta$-separated compactly generated simplicial space $-$then $\aX$ is a $\Delta$-separated compactly generated space.

[$\aX_n$ is a $\Delta$-separated compactly generated space (AD$_6$ 
(cf. p. \pageref{14.4})), 
thus the assertion follows from the fact that $\aX = \colimx \aX_n$ 
(cf. p. \pageref{14.5}).]\\

Let $X$ be a simplicial space $-$then $X$ is said to satisfy the 
\un{cofibration condition} 
\index{cofibration condition (simplicial space)}
if $\forall \ n$ $\&$ $\forall \ i$, $s_i:X_{n-1} \ra X_n$ is a closed cofibration.  
Since the commutative diagram 
\begin{tikzcd}%[sep=small]
{X_{n-1}} \ar{d}[swap]{s_j}  \ar{r}{s_i} &{X_n} \ar{d}{s_{j+1}}\\
{X_n} \ar{r}[swap]{s_i}  &{X_{n+1}}
\end{tikzcd}
is a pullback square $(0 \leq i \leq j \leq n-1)$, one can use Proposition 8 in $\S 3$ to see that the cofibration condition implies that the $sX_{n-1} \ra X_n$ are closed cofibrations.

Example:  Given a topological space $X$, denote by $\si X$ the 
\un{constant simplicial set} 
\index{constant simplicial set} 
on $X$, i.e., 
$\si X([n]) = X \ \&$
$
\begin{cases}
\ d_i = \id_X\\
\ s_i = \id_X
\end{cases}
(\forall \ n)
$
$-$then $\si X$ satisfies the cofibration condition and $\abs{\si X} \approx X$.\\

\begingroup%%----------------------------------->>
\fontsize{9pt}{11pt}\selectfont
Since $L_nX$ can be identified with $sX_{n-1}$, every $X$ which satisfies the cofibration condition is necessarily cofibrant (Reedy structure).\\
\vspi
\textbf{\small FACT} \ 
Suppose that $X$ is a simplicial space satisfying the embedding condition $-$then $X$ satisfies the cofibration condition iff $X$ is Reedy cofibrant.\\
\endgroup %%------------------------------------<<

\begin{proposition} \ %03
Let $X$ be a simplicial space.  Assume: $X$ satisifies the cofibration condition $-$then $\forall \ n$, the arrow 
$\aX_{n-1} \ra \aX_n$ is a closed cofibration.
\end{proposition}

%%----------------------------------------------------------------------------------------------04
[The arrow 
$X_n \times \ddpn \cup sX_{n-1} \times \dpn \ra X_n \times \dpn$ 
is a closed cofibration (cf. $\S 3$, Proposition 7).  Now quote Proposition 2 (cf. $\S 3$, Proposition 2).]\\

Application: Let $X$ be a compactly generated simplicial space satisfying the cofibration condition.  \ 
Assume: $\forall \ n$, $X_n$ is Hausdorff $-$then $\abs{X}$ is a compactly generated Hausdorff space.

[This follows from the lemma on 
p. \pageref{14.6} 
and condition B on 
p. \pageref{14.7}.]\\

\label{14.19}
Application: Let $X$ be a compactly generated simplicial space satisfying the cofibration condition.  
Assume: $\forall \ n$, $X_n$ is Hausdorff $-$then $\aX$ is a compactly generated Hausdorff space.

[This follows from the lemma on 
p. \pageref{14.8} 
and condition B on 
p. \pageref{14.9}.]\\

Application: Let $X$ be a simplicial space satisfying the cofibration condition.  
Assume: $\forall \ n$, $X_n$ is numerably contractible $-$then $\aX$ is numerably contractible.

[It suffices to show that the $\aX_n$ are numerably contractible 
(cf. p. \pageref{14.10}).  
But inductively, the double mapping cylinder of the 2-source 
$X_n \times \dpn \la$ $X_n \times \ddpn \cup sX_{n-1} \times \dpn \ra$ $\aX_{n-1}$ is numerably contractible and numerable contractibility is a homotopy type invariant 
(cf. p. \pageref{14.11}).]\\

\begingroup%%----------------------------------->>
\fontsize{9pt}{11pt}\selectfont
\textbf{\small EXAMPLE} \ 
Let $X$ be a Hausdorff simplicial space.  Assume: $\forall \ n$, the inclusions $\Delta_{X_n} \ra X_n \times X_n$ is a cofibration $-$then $X$ satisfies the cofibration condition.
\vspi
[$\forall \ i$, $s_iX_{n-1}$ is a retract of $X_n$, hence the inclusion $s_iX_{n-1} \ra X_n$ is a closed cofibration 
(cf. p. \pageref{14.12}).]\\
\endgroup %%------------------------------------<<

\label{14.99}
\label{14.115}
\label{14.131}
\label{14.147}
\begingroup%%----------------------------------->>
\fontsize{9pt}{11pt}\selectfont
\textbf{\small LEMMA} \ 
Let
\begin{tikzcd}[sep=large]
{X^0} \ar{d}  \ar{r} &{X^1} \ar{d}  \ar{r} &{\cdots}\\
{Y^0} \ar{r}  &{Y^1}  \ar{r} &{\cdots}
\end{tikzcd}
be a commutative ladder connecting two expanding sequences of topological spaces.  
Assume: $\forall \ n$, the inclusions 
$
\begin{cases}
\ X^n \ra X^{n+1}\\
\ Y^n \ra Y^{n+1}
\end{cases}
$
are closed cofibrations, 
\begin{tikzcd}[sep=large]
{X^n} \ar{d} \ar{r} &{X^{n+1}} \ar{d}\\
{Y^n} \ar{r}  &{Y^{n+1}}
\end{tikzcd}
is a pullback square, and the vertical arrows $\phi^n:X^n \ra Y^n$ are closed cofibrations $-$then the induced map 
$\phi^\infty:X^\infty \ra Y^\infty$ is a closed cofibration.
\vspi
[Take any arrow $Z \ra B$ which is both a homotopy equivalence and a Hurewicz fibration and construct a filler 
$Y^\infty \ra Z$ for 
\begin{tikzcd}[sep=large]
{X^\infty} \ar{d} \ar{r} &{Z} \ar{d}\\
{Y^\infty} \ar{r}  &{B}
\end{tikzcd}
via induction, noting that
$Y^n \underset{X^n}{\sqcup} X^{n+1} \ra Y^{n+1}$ is a closed cofibration (cf. $\S 3$, Proposition 8).]\\
\endgroup %%------------------------------------<<

\label{14.114}
\begingroup%%----------------------------------->>
\fontsize{9pt}{11pt}\selectfont
Application: Let $X^0 \subset X^1 \subset \cdots $ be an expanding sequence of topological spaces.  
Assume: $\forall \ n$, $X^n$ is in $\bDelta\text{-}\bCG$, $X^n \ra X^{n+1}$ is a cofibration, and 
$\Delta_{X^n} \ra X^n \times_k X^n$ is a cofibration $-$then 
$\Delta_{X^\infty} \ra X^\infty \times_k X^\infty$ is a cofibration.\\
\endgroup %%------------------------------------<<

%%----------------------------------------------------------------------------------------------05
\begingroup%%----------------------------------->>
\fontsize{9pt}{11pt}\selectfont
\textbf{\small EXAMPLE} \ 
Let $X$ be a $\Delta$-separated compactly generated simplicial space.  
Assume: $\forall \ n$,  
$\Delta_{X^n} \ra X_n \times_k X_n$ is a cofibration $-$then $X$ satisfies the cofibration condition 
(cf. p. \pageref{14.13}) 
and 
$\Delta_{\aX_n} \ra \aX_n \times_k \aX_n$ is a cofibration 
(cf. p. \pageref{14.14}).  
Therefore 
$\Delta_{\aX} \ra \aX \times_k \aX$ is a cofibration.\\
\endgroup %%------------------------------------<<

\begingroup%%----------------------------------->>
\fontsize{9pt}{11pt}\selectfont
\textbf{\small FACT} \ 
Let 
$
\begin{cases}
\ X\\
\ Y
\end{cases}
$
be $\Delta$-separated compactly generated simplicial spaces satisfying the cofibration condition.  Suppose that
$f:X \ra Y$ is a simplicial map such that $\forall \ n$, $f_n:X_n \ra Y_n$ is a cofibration $-$then 
$\abs{f}:\aX \ra \abs{Y}$ is a cofibration.
\vspi
[Use the lemma on 
p. \pageref{14.15} ff. 
to conclude that $\forall \ n$, $\abs{f}_n:\aX_n \ra \abs{Y}_n$ is a cofibration.  
And: 
\begin{tikzcd}[sep=large]
{\aX_{n-1}} \ar{d} \ar{r} &{\aX_n} \ar{d}\\
{\abs{Y}_{n-1}} \ar{r}  &{\abs{Y}_n}
\end{tikzcd}
is a pullback square.]\\
\vspace{0.25cm}
\endgroup %%------------------------------------<<

\label{14.144} %dmc mnft

\begin{proposition} \ %04
Suppose that
$
\begin{cases}
\ X\\
\ Y
\end{cases}
$
are simplicial spaces satisfying the cofibration condition and let $f:X \ra Y$ be a simplicial map.  
Assume: $\forall \ n$, 
$f_n:X_n \ra Y_n$ is a homotopy equivalence $-$then $\abs{f}:\aX \ra \abs{Y}$ is a homotopy equivalence.
\end{proposition}

[Since 
$
\begin{cases}
\ \aX = \colimx \aX_n\\
\ \abs{Y} = \colim\abs{Y}_n
\end{cases}
$
and the 
$
\begin{cases}
\ \aX_{n-1} \ra \aX_n\\
\ \abs{Y}_{n-1} \ra \abs{Y}_n
\end{cases}
$
are closed cofibrations, it need only be shown that the $\aX_n \ra \abs{Y}_n$ are homotopy equivalences (cf. $\S 3$, Proposition 15).  This is done by induction, the point being that $sX_{n-1} \ra sY_{n-1}$ is a homotopy equivalence.]\\

\begingroup%%----------------------------------->>
\fontsize{9pt}{11pt}\selectfont
\textbf{\small EXAMPLE} \ 
Let $X$ be a simplicial space such that $\forall \ n$, $X_n$ has the homotopy type of a compactly generated space $-$then the arrow 
$\abs{kX} \ra \aX$ is a homotopy equivalence if $X$ satisfies the cofibration condition.
\vspi
[$\forall \ n$, $kX_n \ra X_n$ is a homotopy equivalence and $kX$ satisfies the cofibration condition 
(cf. p. \pageref{14.16}).]\\
\endgroup %%------------------------------------<<

Given an $X$ in \bSITOP, the 
\un{homotopic realization}
\index{homotopic realization (\bSITOP)} 
of $X$ is the quotient 
$\HR X = \coprod\limits_n X_n \times \dpn / \sim$, where $\sim$ is restricted to the monomorphisms in $\bDelta$, i.e., 
$((X\alpha)x,t) \sim (x,\Delta^\alpha t)$ $(\alpha \in M_{\bDelta})$.  
Write (HR$X)_n$ for the image of 
$\coprod\limits_{m \leq n} X_m \times \Delta^m$ under the projection 
$\coprod\limits_n X_n \times \dpn / \sim \ \ra \  \HR X$.

Example: Viewing a simplicial set $X$ as a ``discrete'' simplicial space, $\HR X = \abs{UX}_M$ 
(cf. p. \pageref{14.17}).

Example: $\abs{*} = *$ but  $\HR* = $ ``a large contractible space''.
\\

\begin{proposition} \ %05
Let $X$ be a simplicial space $-$then $\forall \ n$, $(\HR X)_n$ is a closed subspace of 
$\HR X$ and $\HR X= \colim(\HR X)_n$.\\
\end{proposition}

%%----------------------------------------------------------------------------------------------06
\begin{proposition} \ %06
Let $X$ be a simplicial space $-$then $\forall \ n$, there is a pushout square
\begin{tikzcd}%[sep=small]
{X_n \times \ddpn} \ar{d} \ar{r} &{(\HR X)_{n-1}} \ar{d}\\
{X_n \times \dpn} \ar{r}  &{(\HR X)_n}
\end{tikzcd}
and the arrow $(\HR X)_{n-1} \ra (\HR X)_n$ is a closed cofibration.\\
\end{proposition}
\vspace{0.25cm}

\label{14.20}
\begingroup%%----------------------------------->>
\fontsize{9pt}{11pt}\selectfont
\textbf{\small FACT} \ 
Let $X$ be a simplicial space.  Assume: $X_0$ is numerably contractible $-$then $\HR X$ is numerably contractible.
\vspi
[It suffices to show that the $(\HR X)_n$ are numerably contractible 
(cf. p. \pageref{14.18}).  
This is done by induction on $n$, starting from $(\HR X)_0 = X_0$.  
Suppose, therefore, that $n$ is positive and 
$(\HR X)_{n-1}$ is numerably contractible.  
Choose distinct points $u, v \ \in \mathring{\Delta}^n$.  
Because the arrow 
$X \times \dpn \ra (\HR X)_n$ is surjective, $(\HR X)_n = U \cup V$, where 
$U = \im(X_n \times \dpn - \{u\})$, 
$V = \im(X_n \times \dpn - \{v\})$.  
But $\{U,V\}$ is a numerable covering of $(\HR X)_n$ and the retractions
$\dpn - \{u\} \ra \ddpn$, 
$\dpn - \{v\} \ra \ddpn$, induce homotopy equivalences 
$U \ra (\HR X)_{n-1}$, 
$V \ra (\HR X)_{n-1}$.]\\
\endgroup %%------------------------------------<<

It follows from Propositions 5 and 6 that the homotopic realization of a Hausdorff simplicial space is a Hausdorff space and the homotopic realization of a ($\Delta$-separated, Hausdorff) compactly generated simplicial space is a ($\Delta$-separated, Hausdorff) compactly generated space.

[Note: \  Another corollary is that if $\forall \ n$, $X_n$ is a CW space, then $\HR X$ is a CW space (cf. $\S 5$ Propositions 7 and 8).]\\

\begingroup%%----------------------------------->>
\fontsize{9pt}{11pt}\selectfont
Notation UW is the semisimplicial set defined by 
$\UW_n = \{(i_0, \ldots, i_n): i_j \in \Z_{\geq 0} \ \& \  i_0 < \cdots < i_n\}$, where 
$d_j: \UW_n \ra \UW_{n-1}$ sends $(i_0, \ldots, i_n)$ to $(i_0, \ldots, \widehat{i}_j, \ldots, i_n)$.
\vspi
Let $X$ be a simplicial space $-$then the 
\un{unwinding}
\index{unwinding (simplicial space)} 
$\UW X$ is the ``homotopic realization'' of the cofunctor $\bDelta_M \ra \bTOP$ which takes $[n]$ to 
$X_n \times \UW_n$ $(= \ds\coprod\limits_{i_0 < \cdots < i_n} X_n)$.  
Example: $\UW *$ is the ``infinite dimensional simplex'' (Whitehead topology).\\
\endgroup %%------------------------------------<<

\begingroup%%----------------------------------->>
\fontsize{9pt}{11pt}\selectfont
\textbf{\small EXAMPLE} \ 
Let $G$ be a topological group, \bG the topological groupoid having a single object $*$ with $\Mor(*,*) = G$ $-$then 
$\ner \bG$ is a simplicial space and there is a canonical continuous bijection $\UW\ner \bG \ra B_G^\infty$.
\vspi
[Note: \  This arrow is not a homeomorphism (consider $G = *$) but it is a homotopy equivalence.]\\
\endgroup %%------------------------------------<<

\begingroup%%----------------------------------->>
\fontsize{9pt}{11pt}\selectfont
\textbf{\small FACT} \ 
For every simplicial space $X$, the projection $\UW X \ra \HR X$ is a homotopy equivalence.\\
\endgroup %%------------------------------------<<

\begin{proposition} \ %07
Let $X$ be a simplicial space.  Assume: $X$ satisfies the cofibration condition $-$then the arrow 
$\HR X \ra \aX$ is a homotopy equivalence.
\end{proposition}

[The argument is similar to that used in the proof of Proposition 4 in $\S 13$.]\\

%%----------------------------------------------------------------------------------------------07
Application: Let $X$ be a simplicial space.  Assume: $\forall \ n$, $X_n$ is a CW space $-$then $\aX$ is a CW space whenever $X$ satisfies the cofibration condition.\\

\label{14.41}
\label{14.52}
\begingroup%%----------------------------------->>
\fontsize{9pt}{11pt}\selectfont
\textbf{\small EXAMPLE} \ 
Let $X$ be a simplicial space satisfying the cofibration condition.  
Assume: $X_0$ is numerably contractible $-$then $\aX$ is numerably contractible 
(cf. p. \pageref{14.19}).
\vspi
[$\HR X$ is numerably contractible  
(cf. p. \pageref{14.20}) 
and numerable contractibility is a homotopy type invariant 
(cf. p. \pageref{14.21}).]\\
\endgroup %%------------------------------------<<

\label{14.45}
\begingroup%%----------------------------------->>
\fontsize{9pt}{11pt}\selectfont
\textbf{\small FACT} \ 
Equip \bTOP with its standard structure.  Let $f:X \ra Y$ be a simplicial map.  Assume: 
$\forall \ m,n \ \& \ \forall \  \alpha:[m] \ra [n]$, the commutative diagram 
\begin{tikzcd}[sep=large]
{X_n} \ar{d}[swap]{f_n} \ar{r}{X_\alpha} &{X_m} \ar{d}{f_m} \\
{Y_n} \ar{r}[swap]{Y\alpha}  &{Y_m} 
\end{tikzcd}
is a homotopy pullback $-$then $\forall \ n$, 
\begin{tikzcd}[sep=large]
{X_n \times \dpn} \ar{d} \ar{r} &{\HR X} \ar{d} \\
{Y_n \times \dpn} \ar{r}  &{\HR Y} 
\end{tikzcd}
is a homotopy pullback.
\vspi
[One first shows by induction that $\ \forall \ n$, \quad
\begin{tikzcd}[sep=large]
{X_n \times \dpn} \ar{d} \ar{r} &{(\HR X)_n} \ar{d}\\
{Y_n \times \dpn} \ar{r}  &{(\HR Y)_n}
\end{tikzcd}
\quad
is a homotopy pullback.  \ 
To carry out the passage from $n-1$ to $n$, observe that the squares in the commutative diagram\\
\begin{tikzcd}[sep=large]
{X_n \times \dpn} \ar{d} &{X_n \times \ddpn} \ar{l}  \ar{d}  \ar{r} &{(\HR X)_{n-1}}  \ar{d}\\
{Y_n \times \dpn}           &{Y_n \times \ddpn} \ar{l}             \ar{r} &{(\HR Y)_{n-1}} 
\end{tikzcd}
are homotopy pullbacks, thus the squares in the commutative diagram
\begin{tikzcd}[sep=large]
{X_n \times \dpn} \ar{d} \ar{r} &{(\HR X)_n} \ar{d}  &{(\HR X)_{n-1}}  \ar{l} \ar{d}\\
{Y_n \times \dpn} \ar{r}  &{(\HR Y)_n}  &{(\HR Y)_{n-1}} \ar{l}
\end{tikzcd}
are homotopy pullbacks 
(cf. p. \pageref{14.22}).  
So $\forall \ n$, 
\begin{tikzcd}[sep=large]
{(\HR X)_n} \ar{d} \ar{r} &{\HR X} \ar{d}\\
{(\HR Y)_n} \ar{r} &{\HR Y}
\end{tikzcd}
is a homotopy pullback 
(cf. p. \pageref{14.23}).  
Accordingly, both squares in the commutative diagram 
\begin{tikzcd}[sep=large]
{X_n \times \dpn} \ar{d} \ar{r} &{(\HR X)_n} \ar{d} \ar{r} &{\HR X} \ar{d}\\
{Y_n \times \dpn} \ar{r}  &{(\HR Y)_n} \ar{r} &{\HR Y}
\end{tikzcd}
are homotopy pullbacks, hence by the composition lemma, 
\begin{tikzcd}[sep=large]
{X_n \times \dpn} \ar{d} \ar{r} &{\HR X} \ar{d} \\
{Y_n \times \dpn} \ar{r}  &{\HR Y} 
\end{tikzcd}
is a homotopy pullback.]\\
\endgroup %%------------------------------------<<

%%----------------------------------------------------------------------------------------------08
\begingroup%%----------------------------------->>
\fontsize{9pt}{11pt}\selectfont
It follows from Proposition 7 that this result remains valid if 
$
\begin{cases}
\ \HR X\\
\ \HR Y
\end{cases}
$
are replaced by 
$
\begin{cases}
\ \aX\\
\ \abs{Y}
\end{cases}
$
provided that 
$
\begin{cases}
\ X\\
\ Y
\end{cases}
$
satisfy the cofibration condition.  Proof: Consider the commutative diagram \\
\begin{tikzcd}%[sep=small]
{X_n \times \dpn} \ar{d} \ar{r} &{\HR X} \ar{d} \ar{r} &{\aX} \ar{d}\\
{Y_n \times \dpn} \ar{r}  &{\HR Y} \ar{r} &{\abs{Y}}
\end{tikzcd}
.\\
\endgroup %%------------------------------------<<

\label{18.29} %dmc mnft
\begin{proposition} \ 
Suppose that
$
\begin{cases}
\ X\\
\ Y
\end{cases}
$
are simplicial spaces and let $f:X \ra Y$ be a simplicial map.  Assume: $\forall \ n$, $f_n:X_n \ra Y_n$ is a homotopy equivalence $-$then $\HR f: \HR X \ra \HR Y$ is a homotopy equivalence.\\
\end{proposition}

\begin{proposition} \ %9
Suppose that
$
\begin{cases}
\ X\\
\ Y
\end{cases}
$
are simplicial spaces and let $f:X \ra Y$ be a simplicial map.  
Assume: $\forall \ n$, $f_n:X_n \ra Y_n$ is a weak homotopy equivalence $-$then $\HR f: \HR X \ra \HR Y$ is a weak homotopy equivalence.
\end{proposition}

[If the vertical arrows in the commutative diagram
\begin{tikzcd}%[sep=small]
{X_n \times \dpn} \ar{d} &{X_n \times \ddpn} \ar{l} \ar{r} \ar{d} &{}\\
{Y_n \times \dpn} &{Y_n \times \ddpn} \ar{l} \ar{r}  &{}
\end{tikzcd}
\begin{tikzcd}%[sep=small]
{(\HR X)_{n-1}} \ar{d}\\
{(\HR Y)_{n-1}}
\end{tikzcd}
are weak homotopy equivalences, then the induced map 
$(\HR X)_n \ra (\HR Y)_n$ is a weak homotopy equivalence 
(cf. p. \pageref{14.23a}).  
Pass now to colimits via the result on 
p. \pageref{14.23b}.]\\

\label{14.25}
\label{14.29} 
\label{14.42}
\label{14.129}
\label{14.142}
\label{14.153}
\label{14.178}
Application: Let 
$
\begin{cases}
\ X\\
\ Y
\end{cases}
$
be simplicial spaces satisfying the cofibration condition.  Suppose that $f:X \ra Y$ is a simplicial map such that $\forall \ n$, $f_n:X_n \ra Y_n$ is a weak homotopy equivalence $-$then 
$\abs{f}:\aX \ra \abs{Y}$ is a weak homotopy equivalence.\\

\label{14.30}
\label{14.51}
Example: Let $X$ be a simplicial space satisfying the cofibration condition.  Consider the commutative triangle
\begin{tikzcd}%[sep=small]
{k\aX} \ar{r} \ar{r} &{\aX}\\
&{\abs{kX}} \ar{lu} \ar{u}
\end{tikzcd}
.  By the above, $\abs{kX} \ra \aX$ is a weak homotopy equivalence.  
Since the same is true of 
$k\aX \ra \aX$, it follows that the arrow $\abs{kX} \ra k\aX$ is a weak homotopy equivalence.\\

\label{14.29}
\label{14.49}
\label{14.80}
\begingroup%%----------------------------------->>
\fontsize{9pt}{11pt}\selectfont
\textbf{\small EXAMPLE} \ 
Given an $X$ in \bSITOP, denote by $\abs{\sin X}$ the simplicial space which takes $[n]$ to $\abs{\sin X_n}$.  
Thanks to the Giever-Milnor theorem, the arrow of adjunction $\abs{\sin X_n} \ra X_n$ is a weak homotopy   
%%----------------------------------------------------------------------------------------------09
equivalence.  
On the other hand, $\abs{\sin X}$ satisfies the cofibration condition.  
Consequently, the arrow 
$\norm{\sin X} \ra \aX$ is a weak homotopy equivalence if $X$ satisfies the cofibration condition.
\vspi
[Note: \  $\sin X$ is a bisimplicial set and $\abs{\di\sin X} \approx \norm{\sin X}$.]\\
\endgroup %%------------------------------------<<

\label{18.21}
\begingroup%%----------------------------------->>
\fontsize{9pt}{11pt}\selectfont
\index{Homotopy Pullbacks}
\textbf{\small EXAMPLE \ (\un{Homotopy Pullbacks})} \ 
Equip \bCG with its singular structure and suppose given a commutative diagram
\begin{tikzcd}[sep=large]
{W} \ar{d}\ar{r} \ar{r} &{Y} \ar{d}{g}\\
{X}   \ar{r}[swap]{f}  &{Z}
\end{tikzcd}
of compactly generated simplicial spaces such that
\begin{tikzcd}[sep=large]
{W_n} \ar{d} \ar{r} \ar{r} &{Y_n} \ar{d}\\
{X_n}   \ar{r}  &{Z_n}
\end{tikzcd}
is a homotopy pullback in \bCG $\forall \ n$, where $Y_n$, $Z_n$ are path connected.  The associated commutative diagram 
\begin{tikzcd}[sep=large]
{\sin W} \ar{d}\ar{r} \ar{r} &{\sin Y} \ar{d}\\
{\sin X}   \ar{r}  &{\sin Z}
\end{tikzcd}
of bisimplicial sets then has the property that $\forall \ n$, 
\begin{tikzcd}[sep=large]
{\sin W_n} \ar{d}\ar{r} \ar{r} &{\sin Y_n} \ar{d}\\
{\sin X_n}   \ar{r}  &{\sin Z_n}
\end{tikzcd}
is a homotopy pullback in \bSISET with $\sin Y_n$ and $\sin Z_n$ connected.  Accordingly, 
\begin{tikzcd}[sep=large]
{\di\sin W} \ar{d}\ar{r} \ar{r} &{\di\sin Y} \ar{d}\\
{\di\sin X}   \ar{r}  &{\di\sin Z}
\end{tikzcd}
is a homotopy pullback in \bSISET (theorem of Bousfield-Friedlander), so 
\begin{tikzcd}[sep=large]
{\abs{\di\sin W}} \ar{d}\ar{r} \ar{r} &{\abs{\di\sin Y}} \ar{d}\\
{\abs{\di\sin X}}   \ar{r}  &{\abs{\di\sin Z}}
\end{tikzcd}
is a homotopy pullback in \bCG 
(cf. p. \pageref{14.24}).  
Therefore 
\begin{tikzcd}[sep=large]
{\abs{W}}\ar{d} \ar{r} &{\abs{Y}} \ar{d}{\abs{g}}\\
{\aX} \ar{r}[swap]{\abs{f}} &{\abs{Z}}
\end{tikzcd}
is a homotopy pullback in \bCG if $W$, $X$, $Y$, $Z$ satisfy the cofibration condition.
\vspi
[Note: \ Equip \bTOP with its singular structure and suppose given a commutative diagram 
\begin{tikzcd}[sep=large]
{W} \ar{d}\ar{r} \ar{r} &{Y} \ar{d}{g}\\
{X}   \ar{r}[swap]{f}  &{Z}
\end{tikzcd}
of simplicial spaces such that
\begin{tikzcd}[sep=large]
{W_n} \ar{d} \ar{r} \ar{r} &{Y_n} \ar{d}\\
{X_n}   \ar{r}  &{Z_n}
\end{tikzcd}
is a homotopy pullback in \bTOP $\forall \ n$, where $Y_n$, $Z_n$ are path connected $-$then
\begin{tikzcd}[sep=large]
{\abs{W}}\ar{d} \ar{r} &{\abs{Y}} \ar{d}{\abs{g}}\\
{\aX} \ar{r}[swap]{\abs{f}} &{\abs{Z}}
\end{tikzcd}
is a homotopy pullback in \bTOP if $W$, $X$, $Y$, $Z$ satisfy the cofibration condition.  To see this, observe that
\begin{tikzcd}[sep=large]
{\abs{kW}} \ar{d} \ar{r} &{\abs{kY}} \ar{d}{\abs{kg}}\\
{\abs{kX}}   \ar{r}[swap]{\abs{kf}}  &{\abs{kZ}}
\end{tikzcd}
is a homotopy pullback in \bCG, thus the arrow
%%----------------------------------------------------------------------------------------------10
$\abs{kW} \ra W_{\abs{kf},\abs{kg}}$ is a weak homotopy equivalence.  
In the commutative diagram
$
\begin{tikzcd}[sep=large]
{\abs{W}} \ar{r} &{W_{\abs{f},\abs{g}}}\\
{k\abs{W}} \ar{u}  \ar{r}  &{W_{k\abs{f},k\abs{g}}} \ar{u} 
\end{tikzcd}
,
$ 
the vertical arrow on the left is a weak homotopy equivalence as is the vertical arrow on the right.  
Therefore 
$\abs{W} \ra W_{\abs{f},\abs{g}}$
is a weak homotopy equivalence iff 
$k\abs{W} \ra W_{k\abs{f},k\abs{g}}$ 
is a weak homotopy equivalence.  Working in the compactly generated category, form 
\begin{tikzcd}[sep=large]
{\abs{kX}} \ar{d} \ar{r}{\abs{kf}} &{\abs{kZ}} \ar{d} &{\abs{kY}} \ar{d} \ar{l}[swap]{\abs{kg}}\\
{k\aX} \ar{r}[swap]{k\abs{f}} &{k\abs{Z}} &{k\abs{Y}} \ar{l}{k\abs{g}}
\end{tikzcd}
.  
The vertical arrows are weak homotopy equivalences 
(cf. p. \pageref{14.25}), 
so 
$W_{\abs{kf},\abs{kg}} \ra W_{k\abs{f},k\abs{g}}$
is a weak homotopy equivalence 
(cf. p. \pageref{14.26}).  
Examination of 
\begin{tikzcd}[sep=large]
{W_{\abs{kf},\abs{kg}}} \ar{r} &{W_{k\abs{f},k\abs{g}}}\\
{\abs{kW}} \ar{u}  \ar{r}  &{k\abs{W}} \ar{u} 
\end{tikzcd}
then implies that
$k\abs{W} \ra W_{k\abs{f},k\abs{g}}$ 
is a weak homotopy equivalence.]\\
\endgroup %%------------------------------------<<

\label{14.81} %dmc mnft
\begin{proposition} \ 
If 
$
\begin{cases}
\ X\\
\ Y
\end{cases}
$
are Hausdorff simplicial spaces and if $f:X \ra Y$ is a simplicial map such that $\forall \ n$, 
$f_n:X_n \ra Y_n$ is a homology equivalence, then 
$\HR f:\HR X \ra \HR Y$ is a homology equivalence, thus so is 
$\abs{f}: \aX \ra \abs{Y}$ 
subject to the cofibration condition on 
$
\begin{cases}
\ X\\
\ Y
\end{cases}
.
$
\end{proposition}

[By Mayer-Vietoris and the five lemma, the arrow 
$(\HR X)_n \ra (\HR Y)_n$ is a homology equivalence $\forall \ n$.]

[Note: \ The Hausdorff assumption can be replaced by \dsep and compactly generated.]\\

Notation: Given an $X$ in \bSITOP, put $IX = X \times \si[0,1]$, so $\forall \ n$, 
$(IX)_n = IX_n$.  \\

\textbf{\small LEMMA} \ 
For every simplicial space $X$, $\abs{IX} \approx I\aX$.

[The functor $-\times[0,1]:\bTOP \ra \bTOP$ has a right adjoint, thus preserves colimits, in particular, coends.]\\

Application: Let $X$, $Y$, be simplicial spaces, $H:IX \ra Y$ a simplicial map $-$then 
$\abs{H \circx i_0} \simeq \abs{H \circx i_1}$. \\

Example: Suppose that $X$ is a simplicial space.  Define simplicial spaces 
$\Gamma X$, $\Sigma X$ by 
$(\Gamma X)_n = \Gamma X_n$, $(\Sigma X)_n = \Sigma X_n$
$-$then
$\abs{\Gamma X} \approx \Gamma \aX$,
$\abs{\Sigma X} \approx \Sigma \aX$.

[The diagrams
\begin{tikzcd}%[sep=large]
{X_n} \ar{d} \ar{r} &{*} \ar{d} \\
{{IX_n}} \ar{r}  &{{\Gamma X_n}} 
\end{tikzcd}
,
\begin{tikzcd}%[sep=large]
{X_n \amalg X_n} \ar{d} \ar{r} &{* \amalg *} \ar{d} \\
{{IX_n}} \ar{r}  &{{\Sigma X_n}} 
\end{tikzcd}
determine pushout squares
%%----------------------------------------------------------------------------------------------11
in $[\bDelta^\OP,\bTOP]$, thus the diagrams
\begin{tikzcd}%[sep=large]
{\aX} \ar{d} \ar{r} &{*} \ar{d} \\
{\abs{IX}} \ar{r}  &{\abs{\Gamma X}} 
\end{tikzcd}
,
\begin{tikzcd}%[sep=large]
{\aX \amalg \aX} \ar{d} \ar{r} &{* \amalg *} \ar{d} \\
{\abs{IX}} \ar{r}  &{\abs{\Sigma X}} 
\end{tikzcd}
are pushout squares in \bTOP and, from the lemma, $\abs{IX} \approx I\aX$.]

[Note: \  When dealing with a pointed simplicial space $X$, one can work with either its unpointed geometric realization 
$\ds\int^{[n]} X_n \times \dpn$ 
or its pointed geometric realization 
$\ds\int^{[n]} X_n \# \dpn_+$ .  
However, both give the ``same'' result (consider right adjoints).  
Therefore if one defines pointed simplicial spaces  
$\Gamma X$, $\Sigma X$ by 
$(\Gamma X)_n = \Gamma X_n$, $(\Sigma X)_n = \Sigma X_n$ (pointed cone, pointed suspension), then it is still the case that
$\abs{\Gamma X} \approx \Gamma \aX$,
$\abs{\Sigma X} \approx \Sigma \aX$ (unpointed geometric realization ).]\\

\label{14.130}
\label{14.158}
\label{14.159}
\label{14.161}
\label{14.184}
\begingroup%%----------------------------------->>
\fontsize{9pt}{11pt}\selectfont
\textbf{\small EXAMPLE} \ 
Let $X$ be a pointed simplicial space satisfying the cofibration condition (give $\aX$ the base point $x_0 \in X_0 = \aX_0$).  
Assume: $\forall \ n$, $X_n$ is path connected.  Denote by $\Theta X$ ($\Omega X$) the simplicial space which takes 
$[n]$ to $\Theta X_n$ ($\Omega X_n$) $-$then $\Theta X$ ($\Omega X$) satisfies the cofibration condition (inspect the proof of Proposition 6 in $\S 3$) hence, 
\begin{tikzcd}[sep=large]
{\abs{\Omega X}} \ar{d} \ar{r} &{\abs{\Theta X}} \ar{d} \\
{\{x_0\}} \ar{r}  &{\aX} 
\end{tikzcd}
is a homotopy pullback in \bTOP (singular structure).  
Because there is a commutative diagram
\begin{tikzcd}[sep=large]
{\abs{\Omega X}} \ar{d} \ar{r} &{\abs{\Theta X}} \ar{d} \ar{r} 
&{\aX} \arrow[d,shift right=0.5,dash] \arrow[d,shift right=-0.5,dash]  \\
{\Omega \aX} \ar{r}  &{\Theta\aX}  \ar{r} &{\aX}
\end{tikzcd}
and $\abs{\Theta X}$ is contractible, it follows that the arrow $\abs{\Omega X} \ra \Omega\aX$ is a weak homotopy equivalence.\\
\endgroup %%------------------------------------<<

\label{14.180}
\begingroup%%----------------------------------->>
\fontsize{9pt}{11pt}\selectfont
\textbf{\small FACT} \ 
Let $X$ be a pointed simplicial space satisfying the cofibration condition 
(give $\aX$ the base point $x_0 \in X_0 = \aX_0$) 
$-$then $X_0$ $n$-connected, $X_1$ $(n-1)$-connected, 
$\ldots, X_{n-1}$ 1-connected $\implies$ $\aX$ $n$-connected.
\vspi
[If $n = 1$, one can suppose that $\forall \ m > 1$, $X_m$ is path connected, 
thus $\abs{\Omega X}$ is path connected and 
$* = \pi_0(\abs{\Omega X}) \approx$ 
$\pi_0(\Omega \aX) \approx$ 
$\pi_1(\aX)$.  
If $n > 1$, show that $H_q(\aX) = 0$ $(q \leq n)$ and quote Hurewicz.]\\
\endgroup %%------------------------------------<<

Recall that if $X$ is a locally compact space and $g:Y \ra Z$ is quotient, then $\id_X \times g: X \times Y \ra X \times Z$
is quotient (cf. $\S 2$, Proposition 1 (\mX is cartesian)).  
Here is a variant in which $X$ is allowed to be arbitrary.\\

\index{Whitehead Lemma}
\textbf{\small WHITEHEAD LEMMA} \quad
Let $g:Y \ra Z$ be quotient.  
Assume: $\forall \ z \in Z$ and $\forall$ neighborhood $V$ of $z$, there exists an open subset $U \subset Y$ with $\ov{U}$ compact and contained in $g^{-1}(V)$ such that $g(U)$ is a neighborhood of $z$ $-$then for any $X$, 
$\id_X \times g: X \times Y \ra X \times Z$ is quotient.

%%----------------------------------------------------------------------------------------------12
[Writing $p = \id_X \times g$, the claim is that a subset $O \subset X \times Z$ having the property that 
$p^{-1}(O)$ is open in $X \times Y$ is itself open in $X \times Z$.  
Fix $(x_0,z_0) \in O$ and choose an open 
$Y_0 \subset Y$: $\{x_0\} \times Y_0 = (\{x_0\} \times Y) \cap p^{-1}(O)$.  
If $V_0 = g(Y_0)$, then 
$Y_0 = g^{-1}(V_0)$, so $V_0$ is open in $Z$.  
Per $z_0$ $\&$ $V_0$, take $U_0$ as in the assumption and let
$X_0 = \{x: \{x\} \times \ov{U_0} \subset p^{-1}(O)\}$.  
Since $X_0$ is open in $X$ and 
$(x_0,z_0) \in X_0 \times g(U_0) \subset O$, it follows that $O$ is open in $X \times Z$.]

[Note: \  The argument goes through for any arrow $X \ra W$ which is quotient.]\\

\label{14.141} %dmc mnft

\label{14.37}
Application: For every topological space $X$, $\abs{\si X \times \Delta[1]} \approx X \times [0,1]$.
\\

\textbf{\small LEMMA} \ 
For every smplicial space $X$, $\abs{X \times \dw} \approx \aX \times [0,1]$.

[
$\abs{X \times \dw} \approx$ 
$\ds\int^{[n]} \ds\int^{[m]} X_n \times \dw_m \times \dpn \times \Delta^m \approx$ 
$\ds\int^{[n]} \left(\ds\int^{[m]} X_n \times \dw_m \times \Delta^m \right) \times \dpn \approx$ 
$\ds\int^{[n]} X_n \times [0,1] \times \dpn \approx$ 
$\left( \ds\int^{[n]} X_n \times \dpn \right) \times [0,1]  \approx$ 
$\aX \times [0,1]$.]\\
\vspace{0.25cm}

\begingroup%%----------------------------------->>
\fontsize{9pt}{11pt}\selectfont
\textbf{\small FACT} \ 
Let $X$, $Y$, be simplicial spaces and let $f,g:X \ra Y$ be simplicial maps.  Suppose that $\forall \ n$, there are continuos functions $h_i:X_n \ra Y_{n+1}$ $(0 \leq i \leq n)$ such that 
$d_0 \circx h_0 = f_n$, 
$d_{n+1} \circx h_n = g_n$ 
and
\[
d_{i} \circx h_j = 
\begin{cases}
\ h_{j-1} \circx d_{i} \quad \ (i < j)\\
\ d_{i} \circx h_{i-1} \quad \ (i = j > 0)\\
\ h_j \circx d_{i-1} \quad\  (i > j+1)
\end{cases}
, \qquad s_{i} \circx h_j = 
\begin{cases}
h_{j+1} \circx s_{i} \quad  \ (i \leq j)\\
h_j \circx s_{i-1} \quad \  (i > j)
\end{cases}
.
\]
Then $\abs{f} = \abs{g}$ in the homotopy category.\\
\endgroup %%------------------------------------<<

\label{14.128} %dmc  ? syncs with G but ?
\begingroup%%----------------------------------->>
\fontsize{9pt}{11pt}\selectfont
\textbf{\small EXAMPLE} \ 
Given a triple $\bT = (T,m,\varepsilon)$ in \bTOP, $\forall$ \bT-algebra $X$, $\abs{\barr(T;\bT;X)}$ and $X$ have the same homotopy type 
(cf. p. \pageref{14.27} ff.).]\\
\endgroup %%------------------------------------<<

\label{14.39}
\label{14.47}
\label{18.17}
\label{18.18}
\begingroup%%----------------------------------->>
\fontsize{9pt}{11pt}\selectfont
\textbf{\small EXAMPLE} \ 
Let $X$ be a simplicial space $-$then the 
\un{translate} 
\index{translate (simplicial space)} 
$TX$ of $X$ is the simplicial space with 
$T_nX = X_{n+1}$, where if $\alpha:[m] \ra [n]$, 
$TX(\alpha):T_nX \ra T_mX$ is 
$X(T\alpha):X_{n+1} \ra X_{m+1}$, 
$T\alpha:[m+1] \ra [n+1]$ being the rule that sends 0 to 0 and $i$ to $\alpha(i - 1) + 1$ $(i > 0)$.  
There are simplicial maps 
$\si X_0 \ra TX$, $TX \ra \si X_0$, viz. 
$s_0^{n+1}:X_0 \ra X_{n+1}$,
$d_1^{n+1}:X_{n+1} \ra X_0$, 
and the composition 
$\si X_0 \ra TX \ra \si X_0$ is the identity.  
On the other hand, if 
$h_i:T_nX \ra T_{n+1}X$ is defined by  
$h_i = s_0^{i+1} \circx d_1^i$ $(0 \leq i \leq n)$, then 
$d_1 \circx h_0 = \id$, 
$d_{n+2} \circx h_n = s_0^{n+1} \circx d_1^{n+1}$ and 
\[
d_{i+1} \circx h_j =
\begin{cases}
\ h_{j-1} \circx d_{i+1} \hspace{0.5cm} (i < j)\\
\ d_{i+1} \circx h_{i-1} \hspace{0.5cm} (i = j > 0)\\
\ h_j \circx d_i \hspace{1.15cm} (i > j+1)
\end{cases}
\hspace{-.25cm}, 
\qquad s_{i+1} \circx h_j = 
\begin{cases}
h_{j+1} \circx s_{i+1} \hspace{0.5cm} (i \leq j)\\
h_j \circx s_i \hspace{1.05cm} \ (i > j)
\end{cases}
\hspace{-.25cm}.
\]
Therefore $\abs{TX}$ and $X_0$ have the same homotopy type. 
 In particular: $X_0$ contractible $\implies$ $\abs{TX}$ contractible.\\
\endgroup %%------------------------------------<<

While the general theory of simplicial spaces does not require a compactly generated hypothesis, one can say more with it than without it.  
A key point here is that \bCG
%%----------------------------------------------------------------------------------------------13
admits a closed simplicial action, viz.
$X \bbox K = X \times_k \abs{K}$, relative to which \bCG satisfies SMC in either its standard or singular model category structure.  Note, 
however, that the formal definition of, e.g., 
$\ohc_{\bI}:[\bI,\bCG] \ra \bCG$ depends only on 
$\Box \hthree(\ohc_{\bI}- = \ds\int^i -i \times_k B(i\backslash \bI)$ 
(cf. p. \pageref{14.28})) 
and not on the underlying simplicial model category structure.\\

\begingroup%%----------------------------------->>
\fontsize{9pt}{11pt}\selectfont
\textbf{\small LEMMA} \ 
Let 
$F,G:\bI \ra \bCG$ 
be functors and let
$\Xi:F \ra G$ 
be a natural transformation.  
Assume: $\forall \ i$, 
$\Xi_i:F_i \ra Gi$ is a weak homotopy equivalence $-$then
$\ohc \Xi: \ohc F \ra  \ohc G$ is a weak homotopy equivalence.
\vspi
[One has
$
\begin{cases}
\ \ohc \ F \approx \abs{\ds\coprod F}\\
\ \ohc \ G \approx \abs{\ds\coprod G}
\end{cases}
$
(cf. p. \pageref{14.28a}) and 
$
\begin{cases}
\ \ds\coprod F\\
\ \ds\coprod G
\end{cases}
$
satisfy the cofibration condition.  In addition, $\forall \ n$, 
$\left(\ds\coprod \Xi\right)_n: \left(\ds\coprod F\right)_n \ra \left(\ds\coprod G\right)_n$ 
is a weak homotopy equivalence.  Therefore
$\abs{\ds\coprod \Xi} \abs{\ds\coprod F}  \ra \abs{\ds\coprod G}$ is a weak homotopy equivalence 
(cf. p. \pageref{14.29}).]
\vspi
[Note: \  Changing the assumption to ``homotopy equivalence'' changes the conclusion to ``homotopy equivalence'' (cf. Proposition 4).]\\
\endgroup %%------------------------------------<<

\label{14.34}
\begingroup%%----------------------------------->>
\fontsize{9pt}{11pt}\selectfont
\textbf{\small EXAMPLE} \ 
For any compactly generated space $X$, $\ohc X$ and $\HR X$ have the same weak homotopy type.  
To see this, consider $\abs{\sin X}$ 
(cf. p. \pageref{14.30} ff.) 
$-$then the arrow 
$\ohc \abs{\sin X} \ra \ohc X$ is a weak homotopy equivalence (by the lemma) and the arrow 
$\HR \abs{\sin X} \ra \HR X$ is a weak homotopy equivalence (cf. Proposition 9).  But 
$\abs{\ohc \sin X}$ is homeomorphic to $\ohc \abs{\sin X}$ 
(cf. p. \pageref{14.31}) 
and the homotopy type of 
$\abs{\ohc \sin X}$ is the same as that of 
$\abs{\di \sin X} \approx \norm{\sin X}$ 
(cf. p. \pageref{14.32}), 
the homotopy type of the latter being that of 
$\HR \abs{\sin X}$ (cf. Proposition 7).
\vspi
[Note: \  More is true: $\ohc \ X$ and $\HR X$ have the same homotopy type.  
Thus take \bCG in its standard structure and equip \bSICG with the corresponding Reedy structure $-$then $\forall$ Reedy cofibrant $X$, the arrow 
$\ohc \ X \ra \aX$ is a homotopy equivalence (cf. $\S 13$, Proposition 49) and $\aX$ has the same homotopy type as 
$\HR X$ (cf. Proposition 7).  
To handle an arbitrary $X$, pass to $\sL X$ 
(cf. p. \pageref{14.33}).  
Because the arrow 
$\sL X \ra X$ is a levelwise a homotopy equivalence, 
$\ohc \ \sL X$ and $\ohc \ X$ have the same homotopy type (cf. supra).  However $\sL X$ is Reedy cofibrant, so 
$\ohc \ \sL X$ has the same homotopy type as $\HR \sL X$, i.e., as $\HR X$ (cf. Proposition 8).]\\
\endgroup %%------------------------------------<<

Let 
$
\begin{cases}
\ \bC\\
\ \bD
\end{cases} 
$
and \bI be small categories.\\

\indent\indent $(\otimes_{\bI})$ \quad This is the functor 
$[\bC \times \bI^{\OP},\bCG] \times [\bI \times \bD,\bCG] \ra [\bC \times \bD,\bCG]$ 
given by $(F \otimes_{\bI} G)_{X,Y} = \ds\int^i F(X,i) \times_k G(i,Y)$.\\
\indent\indent $(\Hom_{\bI})$ \quad This is the functor 
$[\bC \times \bI,\bCG]^{\OP} \times [\bI \times \bD,\bCG] \ra [\bC^{\OP} \times \bD,\bCG]$ 
given by $\Hom_{\bI}(F,G)_{X,Y} = \ds\int_i G_{i,Y}^{F_{X,i}}$.

%%----------------------------------------------------------------------------------------------14
[Note: \  In either situation one can, of course, take
$
\begin{cases}
\ \bC\\
\ \bD
\end{cases}
= \bone. 
$
Special cases: 
$* \otimes - \ \approx \  \colim_{\bI}-, \  \Hom_{\bI}(*,-)  \ \approx  \ \lim_{\bI} -$.]

Examples: \quad 
(1) \ $(F \otimes_{\bI} G) \otimes_{\bJ} H \approx F \otimes_{\bI} (G \otimes_{\bJ} H)$;
(2) \ $\Hom_{\bJ} (F \otimes_{\bI} G, H) \approx \Hom_{\bI^\OP}(F,\Hom_{\bJ}(G,H))$.

\label{18.6}
Example: Suppose that $X$ is a compactly generated simplicial space $-$then 
$X \otimes_{\bDelta} \Delta^? = \aX$.

[Note: \  $\Delta^?:\bDelta \ra \bCG$ sends $[n]$ to $\dpn$.]

Example: Suppose that $X$ is a compactly generated simplicial space $-$then 
$X_M \otimes_{\bDelta_M} \Delta_M^? = \HR X$.

[Note: \  $X_M$ is the restriction of $X$ to $\bDelta_M$ and 
$\Delta_M^?:\bDelta_M \ra \bCG$ sends $[n]$ to $\dpn$.]\\

\begingroup%%----------------------------------->>
\fontsize{9pt}{11pt}\selectfont
Given $Y$, $Z$ in $[\bI,\bCG]$, put 
$Z^Y \approx \Hom_{\bI}(\Mor \times Y,Z)$, where 
$\Mor \times Y:\bI^{\OP} \times \bI \ra \bCG$ sends $(j,i)$ to $\Mor(j,i) \times Y_i$.  
So, e.g., $\Hom_{\bI}(\Mor,Z)_j = \ds\int_i Z_i^{\Mor(j,i)} = Z_j$ (integral Yoneda).\\
\vspi
\textbf{\small FACT} \ 
The functor category $[\bI,\bCG]$ is cartesian closed.
\vspi
[Let $X$, $Y$, $Z$ be in $[\bI,\bCG]$  $-$then 
$\Nat(X \times Y, Z) \approx$ 
$\Nat(X \otimes_{\bI^{\OP}}(\Mor \times Y),Z)  \approx$ 
$\Nat(X,\Hom_{\bI}(\Mor \times Y,Z)) \approx$ 
$\Nat(X,Z^Y)$.]\\
\endgroup %%------------------------------------<<

\textbf{\small LEMMA} \ 
Let \bI and \bJ be small categories, $\nabla:\bJ \ra \bI$ a functor $-$then 
$F \circx \nabla^{\OP} \otimes_{\bJ} G \approx F \otimes_{\bI} \lan G$.\\

Notation: Given a small category \bI and functors
$
\begin{cases}
\ F\\
\ G
\end{cases}
:\bI \ra \bCG.
$
write $F \un{\times}_k G$ for the functor $\bI \times \bI \ra \bCG$ that sends $(i,j)$ to $Fi \times_k Gj$.\\

\textbf{\small LEMMA} \ 
Relative to the diagonal 
$\bDelta \ra \bDelta \times \bDelta$, 
$\lan \Delta^? \approx \Delta^?  \un{\times}_k \Delta^?$.\\

\begin{proposition} \ %11
If $X$ and $Y$ are compactly generated simplicial spaces, then 
$\abs{X \times_k Y} \approx \aX \times_k \aY$.
\end{proposition}

[One has 
$\abs{X \times_k Y} \approx$ 
$(X \times_k Y) \otimes_{\bDelta} \Delta^? \approx$ 
$(X \un{\times}_k Y) \otimes_{\bDelta \times \bDelta} \Delta^? \un{\times}_k \Delta^? \approx$ 
$(X \otimes_{\bDelta} \Delta^?) \times_k (Y \otimes_{\bDelta} \Delta^?) \approx$ 
$\aX \times_k \aY$.]

[Note: \  Therefore $\abs{?}$ preserves finite products as long as one works in $[\bDelta^{\OP},\bCG]$.]\\

\begingroup%%----------------------------------->>
\fontsize{9pt}{11pt}\selectfont
It is not true that HR preserves finite products.  However 
$\ohc (X \times_k Y)$ 
and 
$\ohc X \times_k \ohc Y$ 
are homeomorphic, thus 
$\HR(X \times_k Y)$ and 
$\HR X \times_k \HR Y$ 
have the same homotopy type 
(cf. p. \pageref{14.34}).\\
\endgroup %%------------------------------------<<

\begingroup%%----------------------------------->>
\fontsize{9pt}{11pt}\selectfont
\textbf{\small FACT} \ 
Let $X$ be a simplicial object in $\bCG/B$; let $Y$ be an object in $\bCG/B$.  Assume: $B$ is \dsep $-$then
$\abs{X \times_{\si B} \si Y} \approx \aX \times_B Y$.
\vspi
%%----------------------------------------------------------------------------------------------15
[Since $B$ is $\Delta$-separated, the functor $- \times_B Y$ has a right adjoint 
(cf. p. \pageref{14.35}).]\\
\endgroup %%------------------------------------<<

\begingroup%%----------------------------------->>
\fontsize{9pt}{11pt}\selectfont
\textbf{\small FACT} \ 
$\aq:[\bDelta^{\OP},\bDelta\text{-}\bCG] \ra \bDelta\text{-}\bCG$ preserves finite limits.
\vspi
[It suffices to deal with equalizers.  For this, let 
$u,v:X \ra Y$ be a pair of simplicial maps $-$then 
$\abs{\eq(u,v)}$ is closed in $\aX$, which is enough.]\\
\endgroup %%------------------------------------<<

Let \bC be a small category $-$then \bC is said to be \un{compactly generated} 
\index{compactly generated (small category)}
if $O = \Ob \bC$ and $M = \Mor \bC$ are compactly generated topological spaces and the four structure functions 
$s:M \ra O$, 
$t:M \ra O$, 
$e:O \ra M$, 
$c:M \times_O M \ra M$ 
are continuous.  
One appends the term \dsep or Hausdorff 
\index{compactly generated  (small category)}
\index{compactly generated Hausdorff (small category)}
when $O$ and $M$ are, in addtion, \dsep or Hausdorff.  
Example: Every compactly generated semigroup with unit (= monoid in \bCG) determines a compactly generated category.

[Note: \  Any small category can be regarded as a compactly generated category by equipping its objects and morphisms with the discrete topology.]

If \bC, \bD are compactly generated categories, then a functor $F:\bC \ra \bD$ is said to be 
\un{continuous} 
\index{continuous (functor of compactly generated categories)}
provided that the functions 
$
\begin{cases}
\ \Ob \bC \ra \Ob \bD\\
\ X \ra FX
\end{cases}
$
, \ 
$
\begin{cases}
\ \Mor \bC \ra \Mor \bD\\
\ f \ra F f
\end{cases}
$
are continuous.

If \bC, \bD are compactly generated categories and if $F, G:\bC \ra \bD$ are continuous functors, 
then a natural transformation 
$\Xi:F \ra G$ is said to be 
\un{continuous} 
\index{continuous (natural transformation of compactly generated categories)}
provided that the function
$
\begin{cases}
\ \Ob \bC \ra \Mor \bD\\
\ X \ra \Xi_X
\end{cases}
$
is continuous.\\

\begingroup%%----------------------------------->>
\fontsize{9pt}{11pt}\selectfont
In other words, per \bCG, 
compactly generated category = internal category, 
continuous functor = internal functor, 
continuous natural transformation = internal natural transformation.
\vspi
[Note: \  If $(M,O)$ is a category object in \bSISET, then $(\abs{M},\abs{O})$ is a category object in \bCG.  
Conversely, if $(M,O)$ is a category object in \bCG, then $(\sin M,\sin O)$  is a category object in \bSISET.]\\
\endgroup %%------------------------------------<<

Let \bC be a compactly generated category $-$then $\ner \bC$ is a compactly generated simplicial space: 
$\ner_0 \bC = O$, $\ner_1 \bC = M, \ldots, \ner_n \bC = M \times_O \cdots \times_O M$ ($n$ factors) (fiber product in \bCG), 
an $n$-tuple $(f_{n-1}, \ldots, f_0)$ corresponding to 
$X_0 \overset{f_0}{\lra} X_1 \ra \cdots \lra X_{n-1} \overset{f_{n-1}}{\lra} X_n$.  
Thus one can form either the geometric realization or the homotopic realization of $\ner \bC$.  
These two spaces are necessarily compactly generated and they have the same homotopy type if $\ner \bC$ satisfies the cofibration condition (cf. Proposition 7).

[Note: \ 
Meyer\footnote[2]{\textit{Israel J. Math.} \textbf{48} (1984), 331-339.} 
has established versions of Quillen's theorems A and B for compactly generated categories.]\\

%%----------------------------------------------------------------------------------------------16
\begingroup%%----------------------------------->>
\fontsize{9pt}{11pt}\selectfont
\textbf{\small EXAMPLE} \ 
Let \bC be a compactly generated category, where \mO has the discrete topology $-$then \bC is a \bCG-category and 
$\forall \ X,\ Y, \ \Mor(X,Y)$ is a clopen subset of $M$, so $\ner \bC$ satisfies the cofibration condition provided that
$\forall \ X$, the inclusion $\{\id_X\} \ra \Mor(X,X)$ is a closed cofibration.\\
\endgroup %%------------------------------------<<

\begingroup%%----------------------------------->>
\fontsize{9pt}{11pt}\selectfont
\textbf{\small EXAMPLE} \ 
Let \bC be a compactly generated category.  View $M$ as an object in $\bCG/O \times_k O$ via 
$
\begin{cases}
\ s: M \ra O\\
\ t:M \ra O
\end{cases}
. \ 
$
Assume: The \bCG embedding $e:O \ra M$ is a closed cofibration over $O \times_k O$ $-$then $\ner \bC$ satisfies the cofibration condition.\\
\endgroup %%------------------------------------<<

Example: Given an internal category \bM in \bCG, and a right \bM-object $X$ and a left \bM-object $Y$, consider 
$\barr(X;\bM;Y)$, the bar construction on $(X,Y)$.  So: 
$\barr(X;\bM;Y) \approx \ner \bM_{X,Y}$, 
where $\bM_{X,Y} = \tran(X,Y)$, is the translation category of $(X,Y)$.

[Note: \ Suppose that \bI is a small category.  Let 
$F:\bI^\OP \ra \bCG$, $G:\bI \ra \bCG$ be functors 
$-$then \mF determines a right \bI-object $X_F$, \mG  determines a left \bI-object $Y_G$, and there is a canonical arrow 
$\abs{\barr(X_F;\bI;Y_G)} \ra F \otimes_{\bI} G$.]\\

\label{14.104}
\begingroup%%----------------------------------->>
\fontsize{9pt}{11pt}\selectfont
\label{14.175}
To simplify the notation, write $\barr(F;\bI,G)$ in place of $\barr(X_F;\bI;Y_G)$.
\vspi
Examples: 
(1) The assignment $j \ra \abs{\barr(\Mor(-,j);\bI;G)}$ defines a functor $PG:\bI \ra \bCG$ and the arrow of evalutation 
$(PG)j \ra Gj$ is a homotopy equivalence; 
(2) The assignment $i \ra \abs{\barr(F;\bI;\Mor(i,-))}$ defines a functor $PF:\bI^\OP \ra \bCG$ and the arrow of evalutation 
$(PF)i \ra Fi$ is a homotopy equivalence.
\vspi
Observation: $\abs{\barr(F;\bI,G)} \approx PF \otimes_{\bI} G \approx F \otimes_{\bI}PG$.\\
\endgroup %%------------------------------------<<

\begingroup%%----------------------------------->>
\fontsize{9pt}{11pt}\selectfont
\textbf{\small EXAMPLE} \ 
$\ohc G \approx$ 
$B(-\backslash \bI) \otimes_{\bI} G \approx$ 
$P* \otimes_{\bI} G \approx$  
$* \otimes_{\bI} PG \approx$  
$\colimx PG$.\\
\endgroup %%------------------------------------<<

\begingroup%%----------------------------------->>
\fontsize{9pt}{11pt}\selectfont
Working with the unit interval, one can define a notion of homotopy $(\simeq)$ in the functor category $[\bI,\bCG]$ that formally extends the special case $\bI = \bone$.  This leads to a quotient category 
$[\bI,\bCG]/ \simeq$.  
Agreeing to call a morphism in $[\bI,\bCG]$ a 
\un{homotopy equivalence} 
\index{homotopy equivalence ($[\bI,\bCG]$)}
if its image in $[\bI,\bCG]/\simeq$ is an isomorphism, it is seen by the usual argument that $[\bI,\bCG]/\simeq$ is the localization of $[\bI,\bCG]$ at the class of homotopy equivalences.
\vspi
[Note: \  The functor $P:[\bI,\bCG] \ra [\bI,\bCG]$ respects homotopy congruence.]\\
\endgroup %%------------------------------------<<

\begingroup%%----------------------------------->>
\fontsize{9pt}{11pt}\selectfont
\textbf{\small LEMMA} \ 
Let $G^\prime, \ G\pp: \bI \ra \bCG$ be functors and let $\Xi:G^\prime \ra G\pp$ be a natural transformation.  
Assume: $\forall \ j$, $\Xi_j:G^\prime j \ra G\pp j$ is a homotopy equivalence $-$then 
$P\Xi:PG^\prime \ra PG\pp$ is a homotopy equivalence.\\
\endgroup %%------------------------------------<<

\begingroup%%----------------------------------->>
\fontsize{9pt}{11pt}\selectfont
Application: $\forall \ G$, the arrow of evaluation $PPG \ra PG$ is a homotopy equivalence.\\
\endgroup %%------------------------------------<<

\begingroup%%----------------------------------->>
\fontsize{9pt}{11pt}\selectfont
Application: Assume: $\forall \ j, \ G^\prime j, \ G\pp j$ are contractible $-$then there is a homotopy equivalence 
$PG^\prime \ra PG\pp$.
\vspi
%%----------------------------------------------------------------------------------------------17
[The arrows 
$PG^\prime \ra P*$, $PG\pp \ra P*$ are homotopy equivalences.]
\vspi
[Note: \  There is only one homotopy class of arrows \ $PG^\prime \ra PG\pp$.  \ 
Thus suppose that $\Phi$, \ $\Psi:PG^\prime \ra PG\pp$ are not homotopic and form the commutative diagrams
$
\begin{tikzcd}[sep=large]
{PPG^\prime} \ar{d} \ar{r}{P\Phi} &{PPG\pp} \ar{d} \ar{r}{PT} &{P*} \ar{d} \\
{PG^\prime} \ar{r}[swap]{\Phi}  &{PG\pp} \ar{r}[swap]{T} &{*}
\end{tikzcd}
,
$
$
\begin{tikzcd}[sep=large]
{PPG^\prime} \ar{d} \ar{r}{P\Psi} &{PPG\pp} \ar{d} \ar{r}{PT} &{P*} \ar{d} \\
{PG^\prime} \ar{r}[swap]{\Psi}  &{PG\pp} \ar{r}[swap]{T} &{*}
\end{tikzcd}
. \ 
$
Since the vertical arrows in the squares on the left are homotopy equivalences, $P \Phi$, $P \Psi$ are not homotopic.  
On the other hand, 
$T \circx \Phi = T \circx \Psi$ $\implies$ 
$PT \circx P \Phi = PT \circx P\Psi$ $\implies$ 
$P \Phi \simeq P\Psi$ (\mT is a levelwise homotopy equivalence, hence $PT$ is a homotopy equivalence).  Contradiction.]\\
\endgroup %%------------------------------------<<

\label{14.121}

\begin{proposition} \ %12
Suppose that
$
\begin{cases}
\ \bC\\
\ \bD
\end{cases}
$
are compactly generated categories.  Let $F,\ G:\bC \ra \bD$ be continuous functors, $\Xi:F \ra G$ a continuous natural transformation $-$then 
$\abs{\ner F},\ \abs{\ner G}:\abs{\ner \bC} \ra \abs{\ner \bD}$ are homotopic via $\abs{\ner \Xi_H}$ 
(cf. p. \pageref{14.36}).
\end{proposition}

[Note: \  A 
\un{topological category} 
\index{topological category}
is a category object in \bTOP.  And: The analog of Proposition 12 is true in this setting as well 
(since 
$\abs{? \times \dw} \approx \aq \times [0,1]$ 
(cf. p. \pageref{14.37})).]\\

\begingroup%%----------------------------------->>
\fontsize{9pt}{11pt}\selectfont
\textbf{\small EXAMPLE} \ 
Let $X$ be a nonempty compactly generated space.  View $\grd X$ as a compactly generated category $-$ then 
$\abs{\ner \grd X}$ is contractible.
\vspi
\label{14.3a}
[Note: \   For any nonempty topological space $X$, $\grd X$ is a topological category and 
$\abs{\ner \grd X}$ ($= \abs{\Lambda X}$ 
(cf. p. \pageref{14.38})) 
is contractible.]\\
\endgroup %%------------------------------------<<

Given a monoid $G$ in \bCG with the property that the inclusion $\{e\} \ra G$ is a closed cofibration, write \bG for the associated compactly generated category and put 
$XG = \abs{\barr(*;\bG;G)}$ ($(XG)_n = \abs{\barr(*;\bG;G)}_n$), 
$BG = \abs{\barr(*;\bG;*)}$ ($(BG)_n = \abs{\barr(*;\bG;*)}_n$) $-$then 
there are projections 
$XG \ra BG$ ($(XG)_n \ra (BG)_n$) and closed cofibrations 
$G \ra XG$, $\{e\} \ra BG$.

[Note: \  The assumption on $G$ implies that 
$\barr(*;\bG;G)$, $\barr(*;\bG;*)$ satisfy the cofibration condition.]\\

\begingroup%%----------------------------------->>
\fontsize{9pt}{11pt}\selectfont
\textbf{\small EXAMPLE} \ 
$\barr(*;\bG;G)$ is isomorphic to $T\barr(*;\bG;*)$, the translate of $\barr(*;\bG;*)$ 
(cf. p. \pageref{14.39}).

[Use the transposition 
$\barr(*;\bG;G) \overset{\Tee}{\lra} T \barr(*;\bG;*)$ 
defined by 
$\barr_n$
$(*;\bG;G) \overset{\Tee_n}{\lra} T_n \barr(*;\bG;*)$ 
where
$\Tee_n(g_0, \ldots, g_{n-1},g_n) = (g_n, g_0, \ldots, g_{n-1})$.]\\
\endgroup %%------------------------------------<<

\textbf{\small LEMMA} \ 
$XG$ is contractible.\\

%%----------------------------------------------------------------------------------------------18
[Consider the compactly generated category $\tran G$.  It has an initial object, viz. $e$ (the unique morphism from $e$ to $g$ is $(g,e)$).  But the assignment 
$
\begin{cases}
\ G \ra G \times_k G\\
\ g \ra (g,e)
\end{cases}
$
is continuous.  Therefore $\abs{\barr(*;\bG;G)}$ is contractible (cf. Proposition 12).]

[Note: \  $XG$ is a right $G$-space.]\\

\textbf{\small LEMMA} \ 
$BG$ is path connected 
(cf. p. \pageref{14.40}) 
and numerably contractible 
(cf. p. \pageref{14.41}).

[Note: \ $BG$ is called the 
\un{classifying space} 
\index{classifying space} 
of $G$ but I shall pass in silence on just what $BG$ classifies (for an abstract approach to this question, see 
Moerdijk\footnote[2]{\textit{SLN} \textbf{1616} (1995).}).]\\

Remark: \ $XG$ and $BG$ are abelian monoids in \bCG provided that $G$ is abelian.\\

\label{14.84}
\begingroup%%----------------------------------->>
\fontsize{9pt}{11pt}\selectfont
The formation of $\abs{\barr(X;\bG;Y)}$ is functorial in the sense that if 
$\phi:G \ra G^\prime$ is a continuous homomorphism and 
$
\begin{cases}
\ X \ra X^\prime\\
\ Y \ra Y^\prime
\end{cases}
$
are $\phi$-equivariant, then there is an arrow 
$\abs{\barr(X;\bG;Y)} \ra \abs{\barr(X^\prime;\bG^\prime;Y^\prime)}$.  
In particular: $\phi$ induces arrows 
$XG \ra XG^\prime$, 
$BG \ra BG^\prime$.
\vspi
The formation of $\abs{\barr(X;\bG;Y)}$ is product preserving in the sense that the projections define a natural homeomorphism 
$\abs{\barr(X \times_k X^\prime; \bG \times_k \bG^\prime;Y \times_k Y^\prime)} \ra$ 
$\abs{\barr(X;\bG;Y)} \times_k \abs{\barr(X^\prime;\bG^\prime;Y^\prime)}$.
\vspi
[Note:  \ In the compactly generated category, 
$B(G \times_k G^\prime) \approx BG \times_k BG^\prime$ 
but in the topological category all one can say is that the arrow 
$B(G \times G^\prime) \ra BG \times BG^\prime$ 
is a homotopy equivalence 
(Vogt\footnote[3]{\textit{Math. Zeit.} \textbf{153} (1977), 59-82.}).]\\
\endgroup %%------------------------------------<<

\begingroup%%----------------------------------->>
\fontsize{9pt}{11pt}\selectfont
\textbf{\small EXAMPLE} \ 
Let $G$ be a compactly generated group with $\{e\} \ra G$ a closed cofibration $-$then $XG$ is a compactly generated group containing $G$ as a closed subgroup, the action 
$XG \times_k G \ra XG$ agrees with the product in $XG$, $BG$ is the homogeneous space $XG/G$, and $XG$ is a numerable $G$-bundle over $BG$ (in the compactly generated category).\\
\endgroup %%------------------------------------<<

A 
\un{cofibered monoid} 
\index{cofibered monoid} 
is a monoid $G$ in \bCG, for which the inclusion 
$\{e\} \ra G$ is a closed cofibration.\\

\textbf{\small LEMMA} \ 
Let $G,$ $K$ be cofibered monoids in \bCG, $f:G \ra K$ a continuous homomorphism.  Assume: $f$ is a weak homotopy equivalence $-$then $Bf:BG \ra BK$ is a weak homotopy equivalence.

[Apply the criterion on 
p. \pageref{14.42} 
to $\barr \ f: \barr(*;\bG;*) \ra \barr(*;\bK;*)$.]\\

%%----------------------------------------------------------------------------------------------19
\begingroup%%----------------------------------->>
\fontsize{9pt}{11pt}\selectfont
Let $G$ be a monoid in \bCG.  
If the inclusion $\{e\} \ra G$ is not a closed cofibration, consider 
%$\widecheck{G}$ apparently - not check but vee
$\overset{\vee}{G}$ 
(cf. p \pageref{14.43}) 
$-$then by construction, the inclusion 
$\{\overset{\vee}{e}\} \ra \overset{\vee}{G}$ is a closed cofibration.  
Moreover, $\overset{\vee}{G}$ is a monoid in \bCG: Take for the product in $[0,1]$ the usual product and extend the product in $G$ by writing $gt = g = tg$ $(g \in G, 0 \leq t \leq 1)$.  
The retraction $r:\overset{\vee}{G} \ra G$ is a morphism of monoids and a homotopy equivalence.\\
\endgroup %%------------------------------------<<

\label{14.62}
\begingroup%%----------------------------------->>
\fontsize{9pt}{11pt}\selectfont
\textbf{\small EXAMPLE \  (\un{Wreath Products})} \ 
Let $G$ be a cofibered monoid in \bCG $-$then the 
\un{wreath product} 
\index{wreath product (cofibered monoid in \bCG)} 
$S_n \ds \int G$ is the cofibered monoid in \bCG with 
$S_n \ds \int G = S_n \times G^n$ as a set, multiplication being 
$(\sigma,(g_1, \ldots, g_n)) \cdot (\tau(h_1, \ldots, h_n)) =$ 
$(\sigma \tau, (g_{\tau(1)}h_1, \ldots, g_{\tau(n)}h_n)$ (so $(\id,(e, \ldots,e))$ is the unit).  
Generalizing the fact that 
$BS_n \approx XS_n / S_n$, one has $B\bigl( S_n \ds\int G \bigr) \approx XS_n \times_{S_n} (BG)^n$.
\vspi
[Note: \  Embedding $S_n$ in $S_{n+1}$ as the subgroup fixing the last letter and embedding $G^n$ in $G^{n+1}$ as $G^n \times \{e\}$ serves to fix an embedding of $S_n \ds \int G$ in $S_{n+1} \ds \int G$ and 
$S_\infty \ds \int G$ is by definition $\ds\bigcup\limits_n S_n \ds\int G$ (colimit topology).  
Another point is that if $X$ is a compactly generated space on which $G$ operates to the right, then $X^n$ is a compactly generated space on which $S_n \ds \int G$ operates to the right: 
$(x_1, \ldots, x_n) \cdot (\sigma,(g_1, \ldots, g_n)) =$ 
$(x_{\sigma(1)} \cdot g_1, \ldots, x_{\sigma(n)} \cdot g_n)$.]\\
\endgroup %%------------------------------------<<

\label{14.164}
A 
\un{discrete monoid} 
\index{discrete monoid} 
is a monoid $G$ in \bSET equipped with the discrete topology.  If $G$ is a discrete monoid, then $G$ is a cofibered monoid and $BG = B\bG$.  
Example: Suppose that $G$ is a discrete group $-$then $BG$ is a $K(G,1)$.\\

\begingroup%%----------------------------------->>
\fontsize{9pt}{11pt}\selectfont
\textbf{\small EXAMPLE} \ 
Let $G$ be a discrete monoid; let $\phi,\ \psi: G \ra G$ be homomorphisms $-$then $\phi,\ \psi$ correspond to functors 
$\Phi, \ \Psi:\bG \ra \bG$ and there exists a natural transformation $\Xi:\Phi \ra \Psi$ iff $\phi,\ \psi$ are 
\un{semiconjugate}\index{semiconjugate (functors discrete monoid)} 
in the sense that \ $\xi \phi \ = \  \psi \xi$ \  for some $\xi \in G$.  \ 
Semiconjugate homomophisms lead to homotopic maps at the classifying space level 
(cf. p. \pageref{14.44}).  \ 
To illustrate, suppose that $X$ is an infinite set and let $M_X$ be the monoid of one-to-one functions $X \ra X$.  
Fix $\iota \in M_X$: $\#(\iota(X)) = \#(X - \iota(X))$.  \ 
Define a homomorphism \ 
$\phi:M_X \ra M_X$ by $\phi(f)(x) = \iota(f(\iota^{-1}(x)))$ if $x \in \iota(X)$
$\phi(f)(x) = x$ if $x \notin \iota(X)$.  
Obviously, $\iota \hspace{0.03cm}\id_{M_X} = \phi \iota$.  Fix an injection 
$i:X \ra X - \iota(X)$ and let 
$
C_{\id_X} = 
\begin{cases}
\ M_X \ra M_X\\
\ f \ra \id_X
\end{cases}
$
, so $i C_{\id_X} = \phi i$.  
Conclusion: $BM_X$ is contractible.\\
\endgroup %%------------------------------------<<

\begingroup%%----------------------------------->>
\fontsize{9pt}{11pt}\selectfont
\textbf{\small EXAMPLE} \ 
Every nonempty path connected topological space has the weak homotopy type of the classifying space of a discrete monoid (McDuff\footnote[2]{\textit{Topology} \textbf{18} (1979), 313-320.}).  
Consequently, if $G$ is a discrete monoid, then the $\pi_q(BG)$ can be anything at all.
\vspi
[Note: \  Compare this result with the Kan-Thurston theorem.]\\
\endgroup %%------------------------------------<<

\label{14.103}
\begin{proposition} \ 
Let $G$ be a cofibered monoid in \bCG.  Assume: $G$ admits a homotopy inverse $-$then the sequence $G \ra XG \ra BG$ is a fibration up to homotopy (per \bCG (standard structure)).
\end{proposition}

%%----------------------------------------------------------------------------------------------20
\label{14.173}
[The fact that $G$ has a homotopy inverse implies that $\forall \ m, \ n$ $\&$ $\forall \ \alpha$: 
$[m] \ra [n]$, the commutative diagram
\begin{tikzcd}%[sep=small]
{G^{n+1}} \ar{d} \ar{r} &{G^{m+1}} \ar{d}\\
{G^n} \ar{r} &{G^m}
\end{tikzcd}
is a homotopy pullback, which suffices 
(cf. p. \pageref{14.45} ff.).]

[Note: \  If the inclusion $\{e\} \ra G$ is a closed cofibration, $\pi_0(G)$ is a group, and $G$ is numerably contractible, then $G$ admits a homotopy inverse 
(cf. p. \pageref{14.46}).]\\

Notation:  Given a pointed compactly generated space $X$, put 
$\Theta_k X = X^{[0,1]}$, 
$\Omega_k X = X^{\bS^1}$ 
(pointed exponential objects in $\bCG_*$) (dispense with the ``sub $k$'' if there is no question as to the context).

\label{14.79}
Returning to $G$, there is a morphism of H-spaces $G \ra \Omega B G$ which sends $g$ to the loop 
$\sigma_g:[0,1] \ra BG$ defined by $\sigma_g(t) = [g,(1-t,t)]$ $(0 \leq t \leq 1)$.

[Note: \  The base point of $BG$ is $[e,1]$.]\\

\begin{proposition} \ %14
Let $G$ be a cofibered monoid in \bCG.  
Assume: $G$ admits a homotopy inverse $-$then the arrow $G \ra \Omega BG$ is a pointed homotopy equivalence.
\end{proposition}

[There is an arrow $XG \lra \Theta BG$ and a commutative diagram 
$
\begin{tikzcd}%[sep=small]
{G} \ar{d} \ar{r} &{XG} \ar{d}\\ %\ar{r} &{}\\
{\Omega BG} \ar{r} &{\Theta BG} %\ar{r} &{}
\end{tikzcd}
$
$
\begin{tikzcd}%[sep=small]
\ar{r} &{BG} \arrow[d,shift right=0.5,dash] \arrow[d,shift right=-0.5,dash] \\
\ar{r} &{BG}
\end{tikzcd}
.  \ 
$  
Since $XG$ is contractible, the arrow from the compactly generated mapping fiber of $XG \lra BG$ to the compactly generated mapping fiber of $\Theta BG \lra BG$,  i.e., to \ $\Omega BG$, is a homotopy equivalence.  \ 
Therefore by Proposition 13, the arrow 
$G \lra \Omega BG$ is a homotopy equivalence or still, a pointed homotopy equivalence, both spaces being wellpointed.]
\\

Example: Let $G$ be an abelian group $-$then $BG$ is an abelian compactly generated group, so 
$B^{(2)}G \equiv BBG$ is a $K(G,2)$ and by iteration, $B^{(n)}G$ is a $K(G,n)$.\\

\begingroup%%----------------------------------->>
\fontsize{9pt}{11pt}\selectfont
Let $X$ be a pointed compactly generated simplicial space.  Given $n \geq 1$, there are maps 
$\pi_i: [1] \ra [n]$ $(i = 1, \ldots, n)$ where $\pi_i(0) = i - 1$, $\pi_i(1) = i$.  
Definition: $X$ is said to be 
\un{monoidal} 
\index{monoidal (pointed compactly generated simplicial space)} 
if $X_0 = *$ and $\forall \ n \geq 1$, the arrow 
$X_n \ra X_1 \times_k \cdots \times_k X_1$ determined by the $\pi_i$ is a pointed homotopy equivalence.  
Example: Let $G$ be a monoid in \bCG $-$then $\ner \  \bG$ is monoidal.\\
\endgroup %%------------------------------------<<

\begingroup%%----------------------------------->>
\fontsize{9pt}{11pt}\selectfont
\textbf{\small EXAMPLE} \ 
There is a functor 
$\spx:\bCG_* \ra [\bDelta^\OP,\bCG_*]$ that assigns to each pointed compactly generated space $(X,x_0)$ a monoidal compactly generated simplicial space $\spx X$, where, suitably topologized, 
%%----------------------------------------------------------------------------------------------21
$\spx_nX$ is the set of continuous functions $\dpn \ra X$ which carry the vertexes $v_i$ of $\dpn$ to the base point $x_0$ of $X$.  
In particular: $\spx_1X = \Omega X$.
\vspi
[Consider $[0,n]$ as the segmented interval consisting of the edges of $\dpn$ connecting the vertexes 
$v_0, \ldots v_n$ $-$then $[0,n]$ is a  strong deformation retract of $\dpn$ and a continuous function 
$f:[0,n] \ra X$ such that $f(v_i) = x_0$ can be identified with a sequence of $n$ loops in \mX.]
\vspi
[Note: \  $\spx X$ generally does not satisfy the cofibration condition.]\\
\endgroup %%------------------------------------<<

\begingroup%%----------------------------------->>
\fontsize{9pt}{11pt}\selectfont
If $X$ is monoidal, then $X_1$ is a homotopy associative H-space: 
$X_1 \times_k X_1 \ra X_2 \overset{d_1}{\lra} X_1$ 
(relative to some choice of homotopy inverse for $X_2 \ra X_1 \times_k X_1$), thus $\pi_0(X_1)$ is a monoid.  
Moreover, one has an arrow 
$\Sigma X_1 \ra \aX$ ($\Sigma =$ pointed suspension), hence, by adjunction, an arrow 
$X_1 \ra \Omega \aX$ (which is a morphism of H-spaces).\\
\endgroup %%------------------------------------<<

\begingroup%%----------------------------------->>
\fontsize{9pt}{11pt}\selectfont
\textbf{\small FACT} \ 
Let $X$ be a monoidal compactly generated simplicial space.  Assume: $X$ satisfies the cofibration condition and $X_1$ admits a homotopy inverse $-$then the arrow $X_1 \ra \Omega \aX$ is a pointed homotopy equivalence.
\vspi
[The role of $XG$ in the above is played here by the contractible space $\abs{TX}$, where $TX$ is the translate of $X$
(cf. p. \pageref{14.47}), 
and the sequence 
$X_1 \ra \abs{TX} \ra \aX$ is a fibration up to homotopy (per \bCG (standard structure)).]
\vspi
[Note: \  The $d_0:X_{n+1} \ra X_n$ define a simplicial map $TX \ra X$.]\\
\endgroup %%------------------------------------<<

\label{14.179}
\begingroup%%----------------------------------->>
\fontsize{9pt}{11pt}\selectfont
Remark: If \bC is a pointed category with finite products and if $X$ is a monoidal simplicial object in \bC (obvious definition), then $X_1$ is a monoid object in \bC.\\
\endgroup %%------------------------------------<<

\index{Dold-Lashof Theorem}
\index{Theorem Dold-Lashof}
\textbf{\small DOLD-LASHOF THEOREM} \quad
Let $G$ be a cofibered monoid in \bCG $-$then the arrow $G \ra \Omega BG$ is a weak homotopy equivalence iff $\pi_0(G)$ is a group.

[The necessity is clear.  
To establish the sufficiency, note that $\abs{\sin G}$ is a cofibered monoid in \bCG.  Form now the commutative diagram
$
\begin{tikzcd}%[sep=small]
{\abs{\sin G}} \ar{d} \ar{r} &{\Omega B\abs{\sin G}} \ar{d}\\
{G} \ar{r} &{\Omega BG}
\end{tikzcd}
.
$  Thanks to the Giever-Milnor theorem, the arrow of adjunction $\abs{\sin G} \ra G$ is a weak homotopy equivalence.  
Because $\pi_0(\abs{\sin G})$ is a group and $\abs{\sin G}$ is a CW complex, hence numerably contractible 
(cf. p. \pageref{14.48}
 (TCW$_4$)), the arrow $\abs{\sin G} \ra \Omega B \abs{\sin G}$ is, in particular, a  weak homotopy equivalence (cf. supra).  Finally, $B\abs{\sin G} \ra BG$ is a weak homotopy equivalence 
 (cf. p. \pageref{14.49}), 
 thus 
$\Omega B\abs{\sin G} \ra \Omega BG$ is a weak homotopy equivalence 
(cf. p. \pageref{14.50}).  
Therefore the arrow 
$G \ra \Omega BG$ is a weak homotopy equivalence.]
\\

\label{14.88}
Example:  Let $G,\ K$ be path connected cofibered monoids in \bCG, $f:G \ra K$ a continuous homomorphism.  
Assume: $Bf:BG \ra BK$ is a weak homotopy equivalence $-$then $f$ is a weak homotopy equivalence.

%%----------------------------------------------------------------------------------------------22
[Consider the \cd
\begin{tikzcd}%[sep=small]
{G} \ar{d} \ar{r} &{\Omega BG} \ar{d} \\
{K} \ar{r} &{\Omega BK}
\end{tikzcd}
and apply Dold-Lashof.]\\

\begingroup%%----------------------------------->>
\fontsize{9pt}{11pt}\selectfont
Modulo obvious changes in the definitions, Propositions 13 and 14 are valid for cofibered monoids in \bTOP.  
The same holds for the Dold-Lashof theorem.  
Indeed, if $G$ is a cofibered monoid in \bTOP, then $kG$ is a cofibered monoid in \bCG and the arrow $kG \ra G$ is a weak homotopy equivalence.  Suppose in addition that $\pi_0(G)$ is a group $-$then $\pi_0(kG)$ is a group, so the arrow $kG \ra \Omega_k BkG$ is a weak homotopy equivalence.  
On the other hand, $BkG \ra BG$ is a weak homotopy equivalence 
(cf. p. \pageref{14.51}), 
thus 
$\Omega BkG \ra \Omega BG$ is a weak homotopy equivalence, as is $\Omega_k BkG \ra \Omega BkG$.  Since the diagram
\begin{tikzcd}[sep=large]
{kG} \ar{d} \ar{r} &{\Omega_k BkG}\ar{d} \\
{G}  \ar{r} &{\Omega BG}
\end{tikzcd}
commutes, it follows that the arrow $G \ra \Omega BG$ is a weak homotopy equivalence.\\
\endgroup %%------------------------------------<<

\label{14.83} %dmd mnft
\label{14.86}
\begingroup%%----------------------------------->>
\fontsize{9pt}{11pt}\selectfont
\textbf{\small EXAMPLE \ (\un{The Moore Loop Space})} \ 
Let $(X,x_0)$ be a pointed topological space $-$then $\Omega_MX$ is a monoid in \bTOP.  As such, it admits a homotopy inverse and there is a canonical arrow 
$B\Omega_MX \ra X$ such that the composite 
$\Omega X \ra \Omega_MX \ra \Omega B \Omega_MX \ra \Omega X$ is the identity.  Assume now that $X$ is path connected, numerably contractible, and the inclusion $\{x_0\} \ra X$ is a closed cofibration (so $\Omega_MX$ is cofibered (cf. $\S 3$, Proposition 21)).  Owing to Proposition 14, the arrow 
$\Omega_MX \ra \Omega B \Omega_MX$ is a homotopy equivalence.  But the retraction 
$\Omega_MX \ra \Omega X$ is a homotopy equivalence.  Therefore the arrow 
$\Omega B \Omega_MX \ra \Omega X$ is a homotopy equivalence.
Since $B \Omega_MX$ is numerably contractible 
(cf. p. \pageref{14.52}), 
the delooping criterion on 
p. \pageref{14.53} 
then says that the arrow 
$B \Omega_MX \ra X$ is a homotopy equivalence.
\vspi
[Note: \  The same reasoning shows that
$B \Omega_MX \ra X$ is a weak homotopy equivalence provided that $X$ is path connected and the inclusion $\{x_0\} \ra X$ is a closed cofibration.]\\
\endgroup %%------------------------------------<<

\textbf{\small LEMMA} \ 
Let \mM be a simplicial monoid, $Y$ a left \bM-object $-$then $\aY$ is a left $\abs{\bM}$-object and the geometric realization of 
$\abs{\barr(*;\bM;Y)}$ can be identified with 
$\abs{\barr(*;\abs{\bM};\abs{Y})}$.

[One has 
$\barr_n(*;\bM;Y) = M \times \cdots \times M \times Y$.  The geometric realization of 
$[m] \ra \barr_n(*;\bM;Y)_m = M_m \times \cdots \times M_m \times Y_m$ is 
$\abs{M}^n \times_k \abs{Y} = \barr_n(*;\abs{\bM};\abs{Y})$, which, when realized with respect to $[n]$, gives 
$\abs{\barr(*;\abs{\bM};\abs{Y})}$.]

[Note: \  As a special case, 
$\norm{\barr(*;\bM;*)}$ ($= \norm{\ner \bM}$) $\approx B \abs{M}$.  
Alternatively, 
$|[m] \ra |[n] \ra$  $\barr_n(*;\bM;*)_m||$ $\approx$ 
$\abs{[m] \ra \abs{\ner \bM_m}}$ $\approx$ 
$\abs{[m] \ra B\bM_m} \approx$ 
$B\abs{M}$.]\\

\label{14.67}
\label{18.32}
\begingroup%%----------------------------------->>
\fontsize{9pt}{11pt}\selectfont
\textbf{\small EXAMPLE \ (\un{Algebraic K-Theory})} \ 
Let $A$ be a ring with unit.  Put 
$M(A) = \ds\coprod\limits_{n \geq 0} \ner \bGL(n,A)$ ($\ner \bGL(0,A) = \dz$) $-$then 
$M(A)_k = \ds\coprod\limits_{n \geq 0} \bGL(n,A)^k$, 
thus $M(A)$ acquires the structure of a simplicial monoid from matrix addition, i.e., if
$
\begin{cases}
\ (g_1, \ldots, g_k) \in \bGL(n,A)^k\\
\ (h_1, \ldots, h_k) \in \bGL(m,A)^k
\end{cases}
, \ 
$
$(g_1, \ldots, g_k) \cdot (h_1, \ldots, h_k) =$ 
%%----------------------------------------------------------------------------------------------23
$(g_1 \oplus h_1, \ldots, g_k\oplus h_k)$, where 
$g_i \oplus h_i =
\begin{pmatrix}
g_i &0\\
0 &h_i
\end{pmatrix}
\in$ $\bGL(n+m,A)$ $(i = 1, \ldots, k)$.  
Right multiplication by the vertex $1 \in \ner_0 \bGL(1,A)$ determines a simplicial map 
$- \oplus 1: M(A) \ra M(A)$ whose restriction to $\ner \bGL(n,A)$ is the arrow 
$\ner \bGL(n,A) \ra \ner \bGL(n+1,A)$ induced by the canonical inclusion 
$\bGL(n,A) \ra \bGL(n+1,A)$.  The colimit of the diagram
\begin{tikzcd}%[sep=small]
{M(A)}  \ar{r}{-\oplus 1} &{M(A)} \ar{r}{-\oplus 1}  &{\cdots}
\end{tikzcd}
is isomorphic to the simplicial set 
$Y(A) = \ds\coprod\limits_\Z \ner \bGL(A)$.  
It is a left $\bM(A)$-object and the pullback square 
\begin{tikzcd}%[sep=small]
{Y(A)} \ar{d} \ar{r} &{\abs{\barr(*;\bM(A);Y(A))}}\ar{d} \\
{\dz}  \ar{r} &{\abs{\barr(*;\bM(A);*)}}
\end{tikzcd}
is a homology pullback 
(cf. p. \pageref{14.54}).  
In fact, left multiplication by a vertex $n \in M(A)$ shifts the vertexes of $Y(A)$ 
(the term indexed by $z \in \Z$ is sent to the term indexed by $n+z$) and the corresponding map of simplicial sets 
$\ner \bGL(A) \ra \ner \bGL(A)$ is induced by the homomorphism
$
\begin{cases}
\ \bGL(A) \ra \bGL(A)\\
\ g \ra I_n \oplus g
\end{cases}
$
$(I_n$ = rank $n$ identity matrix), so 
$n_*:H_*(\abs{Y(A)}) \ra H_*(\abs{Y(A)})$ is an isomophism.  
But 
$\barr(*;\bM(\bA);Y(A)) \approx$ 
$\colim_{[\N]} \barr(*;\bM(\bA);M(A))$ $\implies$ 
$\abs{\barr(*;\bM(\bA);Y(A))} \approx$ 
$\colim_{[\N]} \abs{\barr(*;\bM(\bA);M(A))}$ and, by the lemma, the geometric realization of 
$\abs{\barr(*;\bM(\bA);M(A))}$ is 
$\abs{\barr(*;\abs{\bM(\bA)};\abs{M(A)})} = X \abs{M(A)}$,  
a contractible space.  
Therefore the geometric realization of 
$\abs{\barr(*;\bM(\bA);Y(A))}$ is contractible 
(cf. p. \pageref{14.55}).  
Consequently, 
$\abs{Y(A)} = \ds\coprod\limits_{\Z} B \bGL(A)$ has the homology of $\Omega B \abs{M(A)}$
($ \abs{M(A)} = \ds\coprod\limits_{n \geq 0} B \bGL(n,A)$) and a model for $B \bGL(A)^+$ is the path component of 
$\Omega B \abs{M(A)}$ containing the constant loop.
\vspi
\label{14.68a}
[Note: \  An analogous discussion can be given for the simplicial monoid 
$M_\infty = \ds\coprod\limits_{n \geq 0} \ner S_n$ that one obtains from the symmetric groups $S_n$.  
Spelled out, if $S_\infty$ is as on 
p. \pageref{14.56}, 
$\ds\coprod\limits_{\Z} B S_\infty$ has the homology of 
$\Omega B \abs{M_\infty}$ ($\abs{M_\infty} = \ds\coprod\limits_{n \geq 0} B S_n$) and a model for 
$B S_\infty^+$ is the path component of $\Omega B \abs{M_\infty}$ containing the constant loop.]\\
\endgroup %%------------------------------------<<

A left \bG-object $Y$ is a compactly generated space on which $G$ operates to the left and there is a commutative diagram
$
\begin{tikzcd}%[sep=small]
{Y} \ar{d} \ar{r} &{\abs{\barr(*;\bG;Y)}}\ar{d} \\
{*}  \ar{r} &{\abs{\barr(*;\bG;*)} = BG}
\end{tikzcd}
.
$
\\

\begin{proposition} \ %15
Let $G$ be a cofibered monoid in \bCG.  Let $Y$ be a left \bG-object such that $\forall \ g \in G$, the arrow 
$y \ra g\cdot y$ is a weak homotopy equivalence $-$then the sequence
$Y \ra \abs{\barr(*;\bG;Y)} \ra BG$ is a fibration up to homotopy (per \bCG (singular structure)).
\end{proposition}

[Pass to the simplicial monoid $\sin G$, noting that $\sin Y$ is a left $\textbf{sin} \hsx \bG$-object.  Since every 
$g \in \sin_0 G$ induces a weak homotopy equivalence $\sin Y \ra \sin Y$, the pullback square 
\begin{tikzcd}%[sep=small]
{\sin Y} \ar{d} \ar{r} &{\abs{\barr(*;\textbf{sin} \hsx \bG;\sin Y)}} \ar{d} \\
{\dz}  \ar{r} &{\abs{\barr(*;\textbf{sin} \hsx \bG;*)}}
\end{tikzcd}
is a homotopy pullback 
(cf. p. \pageref{14.55a}).  
Therefore, taking into account the lemma, the sequence 
$\abs{\sin Y} \ra \abs{\barr(*;\abs{\textbf{sin} \bG};\abs{\sin Y})} \ra B \abs{\sin G}$
%%----------------------------------------------------------------------------------------------24
is a fibration up to homotopy (per \bCG (singular structure)) 
(cf. p. \pageref{14.57}).  
The obvious comparison then implies 
that the same is true for the sequence 
$Y \ra \abs{\barr(*;\bG;Y)} \ra BG$.]

[Note: \  Similar methods lead to a homological version of this proposition.]\\

\begingroup%%----------------------------------->>
\fontsize{9pt}{11pt}\selectfont
\textbf{\small EXAMPLE} \ 
Given a cofibered monoid $G$ in \bCG, let $UG$ be the associated discrete monoid $-$then the mapping fiber of the arrow 
$BUG \ra BG$ at the base point has the weak homotopy type of $\abs{\barr(*;\textbf{UG};G)}$ whenever $\pi_0(G)$ is a group.\\
\endgroup %%------------------------------------<<

The forgetful functor from the category of groups to the category of monoids has a left adjoint that sends a monoid $G$ to its 
\un{group completion}
\index{group completion (of a monoid)} 
$\ov{G}$.  
Example: Let $G$ be any monoid with a zero element ($0g = g0 = 0$ $\forall \ g \in G$), e.g., 
$G = \Z_2^\times$ $-$then $\ov{G} = *$, the trivial group.

[Note: \  $G$ abelian $\implies$ $\ov{G}$ abelian.]\\

\textbf{\small LEMMA} \ 
The functor $G \ra \overline{G}$ preserves finite products.\\

\begingroup%%----------------------------------->>
\fontsize{9pt}{11pt}\selectfont
\textbf{\small EXAMPLE} \ 
Suppose that $G$ is a discrete abelian monoid.  In this situarion, a model for $\ov{G}$ is the quotient of $G \times G$ by the equivalence relation 
$(g^\prime,h^\prime) \sim (g\pp,h\pp)$ iff $\exists \ k^\prime, \ k\pp \in G$ such that 
$(g^\prime k^\prime, h^\prime k^\prime) = (g\pp k\pp,h\pp k\pp)$, the morphism $G \ra \ov{G}$ being induced by 
$g \ra (g,e)$.  Let $G$ operate on $G \times G$ via the diagonal and form 
$\abs{\barr(*;\bG;G \times G)}$ $-$then $\pi_0(\abs{\barr(*;\bG;G \times G)})$ is the coequalizer of 
$
\begin{cases}
\ d_1: G \times (G \times G) \ra G \times G\\
\ d_0: G \times (G \times G) \ra G \times G
\end{cases}
$
(cf. p. \pageref{14.58}), 
which, from the definitions, is precisely $\ov{G}$.
\vspi
[Note: \  There is an arrow 
$\abs{\barr(*;*;G)} \ra \abs{\barr(*;\bG;G \times G)}$ corresponding to 
$(*,*,g) \ra (*,e,(g,e))$ and 
$G \approx \pi_0(G)$ $\approx \pi_0(\abs{\barr(*;*;G)})$.]\\
\endgroup %%------------------------------------<<

\begingroup%%----------------------------------->>
\fontsize{9pt}{11pt}\selectfont
\textbf{\small FACT} \ 
Let \mM be a simplicial monoid, $\ov{M}$ its simplicial group completion $-$then the arrow 
$\pi_0(M) \ra \pi_0(\ov{M})$ is a morphism of monoids and 
$\ov{\pi_0(M)} \approx \pi_0(\ov{M})$.
\vspi
[Representing $\pi_0(M)$ as $\coeq(d_1,d_0)$ 
(cf. p. \pageref{14.59}), 
one has 
$\ov{\pi_0(M)} = \ov{\coeq(d_1,d_0)}$ $\approx$ $\coeq(\ov{d}_1,\ov{d}_0) = $ $\pi_0(\ov{M})$.]\\
\endgroup %%------------------------------------<<

\textbf{\small LEMMA} \ 
Let $X$ be a pointed simplicial set.  
Assume: $X_0 = *$ $-$then $cX$ is a monoid and $\pi_1(X) \approx \ov{cX}$.\\

\label{14.187}
Application: Let $M$ be simplicial monoid $-$then 
$c\abs{\ner \bM} \approx \pi_0(M)$, hence 
$\pi_1(\abs{\ner \bM}) \approx \ov{\pi_0(M)}$ or still, 
$\pi_1(B \abs{M}) \approx \ov{\pi_0(M)}$.\\

\begin{proposition} \ %16
Let $G$ be a cofibered monoid in \bCG $-$then $\pi_1(BG) \approx \ov{\pi_0(G)}$.
\end{proposition}

[In the above, take $M = \sin G$ to get $\pi_1(B \abs{\sin G}) \approx \ov{\pi_0(\sin G)}$.]

%%----------------------------------------------------------------------------------------------25
[Note: \  If $G$ is a discrete monoid, then 
$\pi_1(BG) \approx \ov{G} \approx \pi_1(B\ov{G})$.]\\

\label{14.94} %dmc mnft
\begingroup%%----------------------------------->>
\fontsize{9pt}{11pt}\selectfont
Let \ $M$ be a simplicial monoid, \ $\ov{M}$ its simplicial group completion \ $-$then \ 
$\ov{\pi_0(M)} \approx$ 
$\pi_0(\ov{M})$, so 
$\pi_1(B\abs{M}) \approx$
$\pi_1(B\abs{\ov{M}})$.  
When $\pi_0(M)$ is a group, $\abs{M}$ and $\abs{\ov{M}}$ admit a homotopy inverse 
(cf. p. \pageref{14.60}) 
(CW complexes are numerably contractible 
(cf. p. \pageref{14.61} 
(TCW$_4$))), thus the rows in the commutative diagram
\begin{tikzcd}[sep=large]
{\abs{M}} \ar{d} \ar{r} &{X\abs{M}} \ar{d} \ar{r} &{B\abs{M}} \ar{d}\\
{\abs{\ov{M}}} \ar{r} &{X\abs{\ov{M}}} \ar{r} &{B\abs{\ov{M}}}
\end{tikzcd}
are fibrations up to homotopy per \bCG (standard structure) (cf. Proposition 13).  
Therefore the arrow 
$\abs{M} \ra \abs{\ov{M}}$ is a pointed homotopy equivalence iff the arrow 
$B\abs{M} \ra B\abs{\ov{M}}$ is a pointed homotopy equivalence, i.e., iff the arrow 
$B\abs{M} \ra B\abs{\ov{M}}$ is acyclic (cf. $\S 5$ Proposition 19).  Of course, the arrow 
$\abs{M} \ra \abs{\ov{M}}$ cannot be a pointed homotopy equivalence if $\pi_0(M)$ is not a group.  
Since the fundamental groups of 
$B\abs{M}$ and $B\abs{\ov{M}}$ are isomorphic, the general question is whether the arrow 
$B\abs{M} \ra B\abs{\ov{M}}$ is acyclic and for this one has the criterion provided by Proposition 22 in $\S 5$.\\
\endgroup %%------------------------------------<<

\begingroup%%----------------------------------->>
\fontsize{9pt}{11pt}\selectfont
\textbf{\small EXAMPLE} \ 
Suppose that $G$ is a discrete monoid $-$then the arrow 
$BG \ra B\ov{G}$ is a pointed homotopy equivalence iff \ 
$\Tor_*^{\Z[G]}(\Z,\Z[\ov{G}]) \approx$ 
$\Tor_*^{\Z[\ov{G}]}(\Z,\Z[\ov{G}])$\  i.e.,\  iff \ 
$\Tor_q^{\Z[G]}(\Z,\Z[\ov{G}]) = 0$ $\forall \ q \geq 1$ and 
$\Z \otimes_{\Z[G]} Z[\ov{G}] \approx \Z$.  For instance, this will be true if $G$ is abelian.  It also holds when $G$ is free 
(Cartan-Eilenberg\footnote[2]{\textit{Homological Algebra}, Princeton University Press (1956), 192.}).
\vspi
[Note: \  
$\Tor_0^{\Z[G]}(\Z,\Z[\ov{G}]) \approx$ 
$\Z \otimes_{\Z[G]} \Z[\ov{G}] \approx$ 
$(\Z[G] / I[G]) \otimes_{\Z[G]} \Z[\ov{G}] \approx$ 
$ \Z[\ov{G}] / I[G] \cdot \Z[\ov{G}] \approx$ $\Z$, 
$I[G] \cdot  \Z[\ov{G}]$ being $I[\ov{G}]$.]\\
\endgroup %%------------------------------------<<

\begingroup%%----------------------------------->>
\fontsize{9pt}{11pt}\selectfont
\textbf{\small FACT} \ 
Let $M$ be a simplicial monoid, $\ov{M}$ its simplicial group completion.  Suppose that $\forall \ n$, the arrow 
$BM_n \ra B \ov{M}_n$ is a pointed homotopy equivalence $-$then the arrow 
$B\abs{M} \ra B\abs{\ov{M}}$ is a pointed homotopy equivalence.
\vspi
[Given a $\pi_0(\ov{M})$-module $\ov{A}$, compare the spectral sequence 
$E_{n,m}^1 \approx \Tor_m^{\Z[M_n]}(\Z,\ov{A}) \implies $ $H_{n+m}(B\abs{M},\ov{A})$ 
with the spectral sequence 
$E_{n,m}^1 \approx \Tor_m^{\Z[\ov{M}_n]}(\Z,\ov{A}) \implies $ $H_{n+m}(B\abs{\ov{M}},\ov{A})$.]\\
\endgroup %%------------------------------------<<

\begingroup%%----------------------------------->>
\fontsize{9pt}{11pt}\selectfont
Application: If $\forall \ n$, $M_n$ is abelian or free, then the arrow 
$B\abs{M} \ra B\abs{\ov{M}}$ is a pointed homotopy equivalence.\\
\endgroup %%------------------------------------<<

According to the Dold-Lashof theorem, for a cofibered monoid $G$ in \bCG, the arrow $G \ra \Omega B G$ is a weak homotopy equivalence iff $\pi_0(G)$ is a group.  
What happens in general?  
To give an answer, one replaces ``homotopy'' by ``homology'', the point being that the arrow $G \ra \Omega B G$ is a morphism of H-spaces, thus the arrow $H_*(G) \ra H_*(\Omega B G)$  is a morphism of 
%%----------------------------------------------------------------------------------------------26
Pontryagin rings.  
Viewing $\pi_0(G)$ as a multiplicative subset of $H_*(G)$, the image of $\pi_0(G)$ in 
$H_*(\Omega B G)$ consists of units (since $\pi_0(\Omega B G)$ 
is a group) and under certain conditions, $H_*(\Omega B G)$ represents the localization of 
$H_*(G)$ at $\pi_0(G)$.\\

\index{Theorem: Group Completion Theorem}
\index{Group Completion Theorem}
\textbf{\small GROUP COMPLETION THEOREM} \quad
Let $G$ be a cofibered monoid in \bCG.  Assume: $\pi_0(G)$ is in the center of $H_*(G) $ $-$then 
$H_*(G) [\pi_0(G)^{-1}] \approx H_*(\Omega B G)$.

[Note: \  The diagram
\begin{tikzcd}%[sep=large]
{\Z[\pi_0(G)]} \ar{d} \ar{r} &{\Z[\ov{\pi_0(G)}]} \ar{d}\\
{H_*(G)} \ar{r}  &{H_*(\Omega B G)}
\end{tikzcd}
is therefore a pushout square in the category of graded associative $\Z$-algebras.]\\

\begingroup%%----------------------------------->>
\fontsize{9pt}{11pt}\selectfont

\textbf{\small EXAMPLE} \ \ 
The group completion theorem is false for an arbitrary cofibered monoid in \bCG.  
Thus choose a discrete monoid $G$ whose classifying space $BG$ has the weak homotopy type of $\bS^n$ $(n > 1)$ 
(cf. p. \pageref{14.62}) 
$-$then if the group completion theorem held for \mG, one would have 
$H_*(\Omega \bS^n) \approx H_0(\Omega \bS^n)$, 
an absurdity.\\
\endgroup %%------------------------------------<<

To eleminate topological technicalities, we shall work with $\abs{\sin G}$ and argue simplicially.\\

\textbf{\small LEMMA} \ 
Let $A$ be a ring with unit.  
Suppose that $S$ is a countable multiplicative subset of $A$ which is contained in the center of $A$ 
$-$then $A[S^{-1}]$ is isomorphic as a (left or right) \mA-module to the colimit of 
%$A \overset{\rho_{s_1}}{\lra} A \overset{\rho_{s_2}}{\lra} \cdots$, where $\rho_{s_i}$ 
$A$
$\overset{\rho_{s_1}}{\xrightarrow{\hspace*{0.75cm}}}$ 
$A$
$\overset{\rho_{s_2}}{\xrightarrow{\hspace*{0.75cm}}}$ 
$\cdots$, 
where $\rho_{s_i}$
is right multiplication by $s_i$ 
and $\{s_i\}$ is an enumeration of the elements of $S$, each element being repeated infinitely often.\\

\begin{proposition} \ %17
Let \mM be a simplicial monoid $-$then $H_*(\abs{M})[\pi_0(\abs{M})^{-1}] \approx H_*(\Omega B \abs{M})$ 
provided that $\pi_0(\abs{M})$ is contained in the center of $H_*(\abs{M})$.
\end{proposition}

[As functors of $M$, both sides of the purported relation commute with filtered colimits.  
Because $M$ can be written as a filtered colimit of countable simplicial submonoids $M_k$ such that $\pi_0(\abs{M_k})$ is contained in the center of 
$H_*(\abs{M_k})$, one can assume that $M$ is countable.  
Pick a vertex in each component of $M$ and, with an eye to the lemma, arrange them in a sequence $\{m_i\}$ subject to the proviso that every choice appears an infinity of times.  Consider
%$M \overset{\rho_{m_1}}{\lra} M \overset{\rho_{m_2}}{\lra} \cdots$, 
$M$
$\overset{\rho_{m_1}}{\xrightarrow{\hspace*{0.75cm}}}$ 
$M$
$\overset{\rho_{m_2}}{\xrightarrow{\hspace*{0.75cm}}}$ 
$\cdots$, 
where $\rho_{m_i}:M \ra M$ is 
right multiplication by $m_i$.  This sequence defines an object in $\bFIL(\bSISET)$.  Form its colimit to get  a left \bM-object $Y$ such that the geometric realization of $\abs{\barr(*;\bM;Y)}$ is contractible (compare the discussion in the example preceeding Proposition 15).  By construction, 
$H_*(\abs{Y}) \approx H_*(\abs{M})[\pi_0(\abs{M})^{-1}]$, hence $\forall \ m \in M_0$, 
$m_*:H_*(\abs{Y}) \ra H_*(\abs{Y})$ is an isomorphism.  This means that the
%%----------------------------------------------------------------------------------------------27
pullback square
\begin{tikzcd}%[sep=large]
{Y} \ar{d} \ar{r} &{\abs{\barr(*;\bM;Y)}} \ar{d}\\
{\dz} \ar{r}  &{\abs{\barr(*;\bM;*)}}
\end{tikzcd}
is a homology pullback 
(cf. p. \pageref{14.63}), 
so the arrow from $\abs{Y}$ to the mapping fiber $E$ of 
$\abs{\barr(*;\abs{\bM};\abs{Y})} \ra B\abs{M}$ over the base point is a homology equivalence.  
Working with the standard model category structure on \bCG 
(cf. p. \pageref{14.64}), 
factor the projection 
$X \abs{M} \ra B\abs{M}$ into an acyclic closed cofibration $X\abs{M} \ra X$ followed by a \bCG fibration 
$X \ra B\abs{M}$ to get the commutative diagram
\begin{tikzcd}%[sep=small]
{\abs{M}} \ar{d} \ar{r} &{X \abs{M}} \ar{d} \ar{r} 
&{B\abs{M}} \arrow[d,shift right=0.5,dash] \arrow[d,shift right=-0.5,dash]\\
{F} \ar{r}  &{X} \ar{r}  &{B\abs{M}} 
\end{tikzcd}
, $F$ denoting the fiber.  Choose a filler $X \ra \Theta B \abs{M}$ for
\begin{tikzcd}%[sep=small]
{X \abs{M}} \ar{d} \ar{r} &{\Theta B\abs{M}} \ar{d}\\
{X} \ar{r}  &{B\abs{M}} 
\end{tikzcd}
$-$then
\begin{tikzcd}%[sep=small]
{\abs{M}} \ar{d} \ar{r} &{X \abs{M}} \ar{d} \ar{r} 
&{B\abs{M}} \arrow[d,shift right=0.5,dash] \arrow[d,shift right=-0.5,dash]\\
{F}           \ar{d} \ar{r} &{X}              \ar{d} \ar{r} 
&{B\abs{M}} \arrow[d,shift right=0.5,dash] \arrow[d,shift right=-0.5,dash]\\
{\Omega B\abs{M}}    \ar{r} &{\Theta B\abs{M}}   \ar{r} &{B\abs{M}}
\end{tikzcd}
commutes, the composite $\abs{M} \ra F \ra \Omega B \abs{M}$ being our morphism of H-spaces.  There is also a commutative diagram
\begin{tikzcd}%[sep=small]
{\abs{M}} \ar{d} \ar{r} &{X \abs{M}} \ar{d} \ar{r} 
&{B\abs{M}} \arrow[d,shift right=0.5,dash] \arrow[d,shift right=-0.5,dash]\\
{\abs{Y}} \ar{r}  &{Y\abs{M}} \ar{r}  &{B\abs{M}} 
\end{tikzcd}
, where 
$Y\abs{M} = \abs{\barr(*;\abs{\bM};\abs{Y})}$.  
Putting everything together leads finally to the commutative diagram
$
\begin{tikzcd}%[sep=small]
&{\abs{M}} \ar{ld} \ar{d}\ar{r} &{\abs{Y}} \ar{d}\\
{\Omega B\abs{M}}   &{F} \ar{l}\ar{r}  &{E} 
\end{tikzcd}
.
$  Since the arrows $\Omega B \abs{M} \la F \ra E$ are homotopy equivalences, the result then falls out by applying $H_*$.]\\

Upon forming the commutative diagram 
\begin{tikzcd}%[sep=small]
{\abs{\sin G}} \ar{d} \ar{r} &{\Omega B \abs{\sin G}} \ar{d}\\
{G} \ar{r}  &{\Omega B G}
\end{tikzcd}
, the group completion theorem is seen to follow from Proposition 17.

[Note: \  The centrality hypothesis on $\pi_0(G)$ is automatic if $G$ is homotopy commutative.]\\

\begingroup%%----------------------------------->>
\fontsize{9pt}{11pt}\selectfont
The group completion theorem remains in force when $\Z$ is replaced by any commutative ring \bk with unit as long as 
$\pi_0(G)$ is in the center of $H_*(G;\bk)$.\\
\endgroup %%------------------------------------<<

\index{Strict Monoidal Categories (example)}
\begingroup%%----------------------------------->>
\fontsize{9pt}{11pt}\selectfont
\textbf{\small EXAMPLE \  (\un{Strict Monoidal Categories})} \ 
\bCAT is a monidal category $(\otimes = \times, e = 1)$ and a monoid therein is a strict monoidal category (strict in the sense that multiplication is literally associative (not just up to natural isomorphism) and the unit is a two sided identity).  A strict monoidal category is therefore a category object in \bCAT with object element \bone.  When considered as a discrete cateogory,
%%----------------------------------------------------------------------------------------------28
every monoid in \bSET becomes a strict monoidal category.  
Fix now a strict monoidal category \bM.  Viewing \bM as an internal category in \bCAT, 
one can form $\barr(\bone;\bM,\bone)$ 
(cf. p. \pageref{14.65}), 
which is a simplicial object in \bCAT.  On the other hand, viewing \bM as a small category (= internal category in \bSET), one can form $\ner \bM$ (a simplicial monoid) and $B\bM$ (a cofibered monoid in \bCG).  
Bearing in mind that $\barr(\bone;\bM,\bone):\bDelta^\OP \ra \bCAT$, put 
$\bGM = \gro_{\bDelta^\OP}\barr(\bone;\bM,\bone)$ 
$-$then there is a weak homotopy equivalence 
$\ohc N \barr(\bone;\bM,\bone) \ra \ner \bGM$ 
(cf. p. \pageref{14.66}).  
But there is also a weak homotopy equivalence 
$\ohc N \barr(\bone;\bM,\bone) \ra \abs{N \barr(\bone;\bM,\bone)}$ (cf. $\S 13$, Proposition 49).  Since 
$N \barr(\bone;\bM,\bone) \approx \barr(*;\ner \bM;*)$ 
and
$\norm{\barr(*;\ner \bM;*)} \approx B\abs{\ner \bM}$, 
it follows that $B\abs{\ner \bM}$ and $B\bGM$ have the same homotopy type.  
Therefore
$H_*(B\bM)[\pi_0(BM)^{-1}] \approx$ $H_*(\Omega B\bGM)$ if \bM is in addition symmetric 
(for this condition implies that $B\bM$ is homotopy commutative).
\vspi
\label{14.96}
\label{14.181}
\label{16.57}
[Note: \  A symmetric strict monoidal category is said to be 
\un{permutative}.
\index{permutative (strict monoidal category)}  
Every small symmetric monoidal category is equivalent to a permutative category 
(Isbell\footnote[2]{\textit{J. Algebra} \textbf{13} (1969), 299-307.}).  
Examples: 
(1) $\bGamma$ is a permutative category under wedge sum.  
Thus $\bm \vee \bn = \bm + \bn$ in blocks 
(the empty wedge sum is \bzero) and for 
$
\begin{cases}
\ \gamma: \bm \ra \bn\\
\ \gamma^\prime:\bm^\prime \ra \bn^\prime
\end{cases}
\hspace{-.25cm},
$
$(\gamma \vee \gamma^\prime(k) = 0$ if $\gamma(k) = 0$ or $\gamma^\prime(k) = 0$, otherwise
$
(\gamma \vee \gamma^\prime)(k) \ = \ 
\begin{cases}
\ \gamma(k) \hspace{1.5cm} (1 \leq k \leq m)\\
\ \gamma^\prime(k - m) + n \  (m < k \leq m + m^\prime)
\end{cases}
\hspace{-.25cm};
$
(2) $\bGamma$ is a permutative category under smash product.  Thus $\bm \# \bn = \bm\bn$ via lexicographic ordering of pairs (the empty smash product is \bone) and for 
$
\begin{cases}
\ \gamma: \bm \ra \bn\\
\ \gamma^\prime:\bm^\prime \ra \bn^\prime
\end{cases}
,
$
$(\gamma \# \gamma^\prime)((i-1)m^\prime + i^\prime) = 0$ if $\gamma(i) = 0$ or $\gamma^\prime(i^\prime) = 0$, 
otherwise 
$(\gamma \# \gamma^\prime)((i-1)m^\prime + i^\prime) = $
$(\gamma(i) - 1)n^\prime + \gamma^\prime(i^\prime)$ $(1 \leq i \leq m, 1 \leq i^\prime \leq m^\prime)$.]\\
\endgroup %%------------------------------------<<

\index{Algebraic K-Theory (example)}
\begingroup%%----------------------------------->>
\fontsize{9pt}{11pt}\selectfont
\textbf{\small EXAMPLE \ (\un{Algebraic K-Theory})} \ 
Let $A$ be a ring with unit.  Denote by $\bM(A)$ the category whose objects are the $A^n$ $(n \geq 0)$, there being no morphism from $A^n$ to $A^m$ unless $n = m$, in which case $\Mor(A^n,A^n) = \bGL(n,A)$ $-$then $\bM(A)$ is a permutative category and 
$\ner \bM(A) =$ $M(A) =$ $\ds\coprod\limits_{n \geq 0} \ner \bGL(n,A)$ 
(cf. p. \pageref{14.67} ff.).  
Here, 
$\Z_{\geq 0} \approx$ $\pi_0(B\bM(A))$, 
$\Z \approx$ $\ov{\pi_0(B\bM(A))} \approx$ $\pi_0(\Omega B\abs{M(A)})$, 
and 
$H_*(B\bM(A))[\pi_0(B\bM(A))^{-1}] \approx$ $H_*(\Omega B\abs{M(A)})$.

[Note: \  Write $\bM_\infty$ for the category whose objects are the finite sets 
$\bn \equiv \{0,1,\ldots,n\}$ $(n \geq 0)$ with base point 0, there being no morphism from \bn to \bm unless 
$n = m$, in which case $\Mor(\bn,\bn) = S_n$ (thus $\bM_\infty = \iso\bGamma$ 
(cf. p. \pageref{14.68})).  
Again, $\bM_\infty$ is permutative and the discussion above can be paralleled 
(cf. p. \pageref{14.68a}).]\\
\endgroup %%------------------------------------<<

\label{14.87}
The compactly generated analog of the ``free topological group'' on $X$ $((X,x_0))$ is meaningful on purely formal grounds 
(cf. p. \pageref{14.69}) 
but this situation is simpler since one has a direct description of the topology on 
$F_{\gr}X$ $(F_{\gr}(X,x_0))$, the 
\un{free compactly generated group} 
\index{free compactly generated group}
on $X$ $((X,x_0))$.  
To be specific, consider an $(X,x_0)$ in $\bCG_*$.  Let $(X^{-1},x_0^{-1})$ be a copy of $(X,x_0)$.  
Put 
$\ov{X} = X \vee X^{-1}$, 
$\ov{X}{}^n = \ov{X} \times_k \cdots \times_k \ov{X}$ ($n$ factors) $-$then with 
%%----------------------------------------------------------------------------------------------29
$F_{\gr}(X,x_0)$ the free group on $X - \{x_0\}$, there is a surjection 
$p:\coprod\limits_n \ov{X}{}^n \ra F_{\gr}(X,x_0)$ sending 
$\ov{X}{}^n$ to $F_{\gr}^n(X,x_0)$, the subset of $F_{\gr}(X,x_0)$ 
consisting of those words of length at most $n$, and $F_{\gr}(X,x_0)$ is equipped with the quotient topology derived from $p$.
When $X$ is $\Delta$-separated, $F_{\gr}(X,x_0)$ is $\Delta$-separated, the arrow of adjunction 
$X \ra F_{\gr}(X,x_0)$ is a closed embedding, $F_{\gr}^n(X,x_0)$ is closed, 
$p_n:\ov{X}{}^n \ra F_{\gr}^n(X,x_0)$ is quotient $(p_n = \restr{p}{\ov{X}{}^n})$, 
$F_{\gr}(X,x_0) = \colimx F_{\gr}^n(X,x_0)$, and the commutative diagram 
\begin{tikzcd}%[sep=small]
{\ov{X}{}_*^{n-1}} \ar{d} \ar{r} &{F_{\gr}^{n-1}(X,x_0)} \ar{d}\\
{\ov{X}{}^n} \ar{r} &{F_{\gr}^{n}(X,x_0)}
\end{tikzcd}
is a pushout square $(\ov{X}{}_*^{n-1} = p_n^{-1}(F_{\gr}^{n-1}(X,x_0)))$.

[Note: \  A reference for this material is 
La Martin\footnote[2]{\textit{Dissertationes Math.} \textbf{146} (1977), 1-36; 
see also Ordman, \textit{General Topology Appl.} \textbf{5} (1975), 205-219.}.  
Incidentally, it is false in general that $k$ applied to the free topological group on $(X,x_0)$ is the free compactly generated group on $(X,x_0)$ but if $X$ is the colimit of an expanding sequence of compact Hausdorff spaces, then the free compactly generated group on $(X,x_0)$ is a topological group, hence is the free topological group on $(X,x_0)$.]\\

\begingroup%%----------------------------------->>
\fontsize{9pt}{11pt}\selectfont
\textbf{\small EXAMPLE} \ 
The structure of $F_{\gr}(X,x_0)$ definitely depends on whether one is working in the topological cateogry or the compactly generated category.  This can be seen by taking $X = \Q$.  For the free topological group on $(\Q,0)$ is not compactly generated and its topology is not the quotient topology associated with the projection 
$\ds\coprod\limits_n \ov{\Q}^n \ra F_{\gr}(\Q,0)$.  
Moreover, $F_{\gr}(\Q,0)$ is not the colimit of the $F_{\gr}^n(\Q,0)$.  
Still, $\forall \ n$, $F_{\gr}^n(\Q,0)$ is closed in $F_{\gr}(\Q,0)$ and every compact subset of 
$F_{\gr}(\Q,0)$ is contained in some $F_{\gr}^n(\Q,0)$.  
Nevertheless, 
$p_n:\ov{\Q}^n \ra F_{\gr}^n(\Q,0)$ is not quotient if $n \gg 0$.
\vspi
[Note: \  Details can be found in 
Fay-Ordman-Thomas\footnote[3]{\textit{General Topology Appl.} \textbf{10} (1979), 33-47}.]\\
\endgroup %%------------------------------------<<

The intent of the preceeding remarks is motivational, our main concern being with the free compactly generated monoids, not free compactly generated groups.  
Thus fix $(X,x_0)$ in $\bCG_*$, call $JX$ the free monoid on $X - \{x_0\}$ and give $JX$ the quotient topology coming from 
$\coprod\limits_n X^n \overset{p}{\ra} JX$.  
Letting $\pi$ be the multiplication in $JX$, consider the commutative diagram
$
\begin{tikzcd}%[sep=large]
{\coprod\limits_n X^n \times_k \coprod\limits_n X^n} \ar{d} \ar{r}{p \times_k p} &{JX \times_k JX} \ar{d}{\pi}\\
{\coprod\limits_n X^n} \ar{r}[swap]{p} &{JX}
\end{tikzcd}
.  \ 
$  
Since $\pi \circx (p \times_k p)$ is continuous and $p \times_k p$ is quotient, $\pi$ is continuous.  Therefore $JX$ is a monoid in $\bCG$.  Suppose now that $G$ is a monoid in \bCG and $f:X \ra G$ is a pointed continuous function.  On algebraic grounds, there exists a unique
%%----------------------------------------------------------------------------------------------30
morphism of monoids $J_f:JX \ra G$ rendering the triangle 
\begin{tikzcd}[sep=small]
{X} \ar{rdd}[swap]{f}  \ar{rr} &&{JX} \ar{ldd}{J_f}\\
\\
&{G}
\end{tikzcd}
commutative.
Claim:  $J_f$ is continuous.  Indeed, there is a continuous function 
$p_f:\coprod\limits_n X^n \ra G$ with $J_f \circx p = p_f$.  
But $p$ is quotient, so $J_f$ is continuous.  Therefore $JX$ is the free compactly generated monoid on $(X,x_0)$.]

[Note: \  $JX$ is the 
\un{James construction} 
\index{James construction} on $(X,x_0)$.]\\
\label{14.82}

\begingroup%%----------------------------------->>
\fontsize{9pt}{11pt}\selectfont
$JX$ can be represented a coend, viz. 
$JX \approx \ds\int^n X^n \times_k J\bn$, $J\bn$ the James construction on the pointed finite set 
$\bn = \{0,1, \ldots, n\}$ 
(cf. p. \pageref{14.69a}).\\
\endgroup %%------------------------------------<<

\begingroup%%----------------------------------->>
\fontsize{9pt}{11pt}\selectfont
\textbf{\small LEMMA} \ 
Let $X$ be a pointed compactly generated simplicial space.  Define a simplicial space $JX$ by
$(JX)_n = JX_n$ $-$then $\abs{JX} \approx J\aX$.
\vspi
[In fact, 
$\abs{JX} = \ds\int^n JX_n \times \dpn \approx$ 
$\ds\int^n \left(\ds\int^{\bm} (X_n)^m \times_k J\bm \right) \times \dpn \approx$ 
$\ds\int^{\bm} \left( \ds\int^n  (X_n)^m \times_k \dpn \right) \times_k J\bm \approx$ 
$\ds\int^\bm \abs{X}^m \times_k J\bm \approx$ $J\aX$.]\\
\endgroup %%------------------------------------<<

Put $J^nX = p(X^n)$ and consider $p^{-1}(J^nX) \cap X^m$.  
Obviously, $m < n \implies p^{-1}(J^nX) \cap X^m = X^m$.  
On the other hand, 
$n < m \implies p^{-1}(J^nX) \cap X^m = \bigcup\limits_S X_S^m$, where for 
$S \subset \{1, \ldots, m\}:\#(S) = m - n$, 
$X_S^m = \{(x_1, \ldots, x_m):x_i = x_0 \ (i \in S)\}$.  
Consequently, $J^nX$ is closed in $JX$ if $\{x_0\}$ is closed in \mX.\\

\textbf{\small LEMMA} \ 
Assume: $\{x_0\}$ is closed in $X$.  Let $A$ be a subset of $J^nX$ such that $p^{-1}(A) \cap X^n$ is closed in $X^n$ $-$then $A$ is closed in $JX$.

[Case 1: $m < n$.  Denoting by $i_{m,n}$ the insertion $X^m \ra X^n$ that sends $(x_1, \ldots, x_m)$ to 
$(x_1, \ldots, x_m,x_0, \ldots, x_0)$, one has 
$p^{-1}(A) \ \cap \ X^m = i_{m,n}^{-1}(p^{-1}(A) \cap X^n)$.  
Case 2: $n < m$.   
Write 
$p^{-1}(A) \ \cap \  X^m = \bigcup\limits_S (p_S^{-1}(p^{-1}(A) \cap X^n))$, 
$p_S:X_S^m \ra X^n$ 
the striking map (i.e., $p_S(x_1, \ldots, x_m)$ retains only those $x_i$, where $i \notin S$).]\\

Accordingly, when $\{x_0\} \subset X$ is closed, the arrow $X^n \ra J^nX$ is quotient and the commutative diagram 
\begin{tikzcd}%[sep=small]
{X_*^n} \ar{d} \ar{r} &{J^{n-1}X} \ar{d}\\
{{X}^n} \ar{r} &{J^nX}
\end{tikzcd}
is a pushout square
$(X_*^n = \bigcup\limits_S X_S^n$ $(\#(S) = 1)$ $\implies$ $X^n / X_*^n \approx X \#_k \cdots \#_k X$ ($n$ factors)).  
It therefore follows that if $X$ is $\Delta$-separated, then each $J^nX$ is \dsep ($\text{AD}_6$ 
(cf. p. \pageref{14.69b})), 
hence 
$JX = \colimx J^nX$ is \dsep 
(cf. p.  \pageref{14.69c}).

[Note: \  The arrow of adjunction $X \ra JX$ is a closed embedding.  Reason: The continuous bijection 
$X \ra J^1X$ is quotient.]\\

%%----------------------------------------------------------------------------------------------31
\begin{proposition} \ 
Let $(X,x_0)$ be a wellpointed compactly generated space with $\{x_0\} \subset X$ closed $-$then $(JX,x_0)$ is a wellpointed compactly generated space with $\{x_0\} \subset JX$ closed, thus is a cofibered monoid in \bCG.
\end{proposition}

[In fact, by the above, $\forall \ n$, $J^{n-1}X \ra J^nX$ is a closed cofibration.]\\

\begingroup%%----------------------------------->>
\fontsize{9pt}{11pt}\selectfont
\textbf{\small LEMMA} \ 
If $(X,x_0)$ is a wellpointed compactly generated Hausdorff space, $-$then $(JX,x_0)$ is a wellpointed compactly generated Hausdorff space.
\vspi
[$\forall \ n$, $J^nX$ is Hausdorff 
(cf. p. \pageref{14.70}) 
and condition B on 
p. \pageref{14.71} 
can be applied.]\\
\endgroup %%------------------------------------<<

\begingroup%%----------------------------------->>
\fontsize{9pt}{11pt}\selectfont
\textbf{\small FACT} \ 
Suppose that $(X,x_0)$ is a pointed CW complex $-$then $(JX,x_0)$ is a pointed CW complex.\\
\endgroup %%------------------------------------<<

If $X$ is a wellpointed compactly generated space with $\{x_0\} \subset X$ closed, then the pointed cone 
$\Gamma X$ and the pointed suspension $\Sigma X$ are wellpointed compactly generated spaces with closed basepoints.
\\

Define $E$ by the pushout square
\begin{tikzcd}%[sep=small]
{X \times_k JX} \ar{d} \ar{r} &{JX} \ar{d}\\
{\Gamma X \times_k JX} \ar{r} &{E}
\end{tikzcd}
, where $X \times_k JX \ra JX$ is multiplication.\\

\textbf{\small LEMMA} \ 
$E$ is contractible.

[Letting $E_n$ be the image of $\Gamma X \times_k J^nX$ in $E$, there is a pushout square
\[
\begin{tikzcd}%[sep=small]
{\Gamma X \times_k J^{n-1}X \cup \{x_0\} \times_k J^nX} \ar{d} \ar{r} &{E_{n-1}} \ar{d}\\
{\Gamma X \times_k J^n X } \ar{r} &{E_n}
\end{tikzcd},
\]
so the arrow $E_{n-1} \ra E_n$ is a closed cofibration.  But 
$E_n /E_{n-1} \approx \Gamma X \#_k (J^nX / J^{n-1}X)$, hence $E_n /E_{n-1}$ is contractible.  
Since $E_0 \approx \Gamma X$, it follows by induction that $E_n$ is contractible 
(cf. p. \pageref{14.70a}).  
Therefore $E = \colimx E_n$ is contractible 
(cf. p. \pageref{14.70b}).]
\\

Notation: Given a pointed compactly generated space $X$, let  $\Theta_{kM}X \ (\Omega_{kM}X)$ be the compactly generated Moore mapping (loop) space of $X$ (dispense with the ``sub $k$'' if there is no question as to context).
\\

\begingroup%%----------------------------------->>
\fontsize{9pt}{11pt}\selectfont
There are two ways to place a compactly generated topology on $\Theta_MX \ (\Omega_MX)$.\\
\indent\indent (1) \quad View $\Theta_MX \ (\Omega_MX)$ as a subsest of $C(\R_{\geq 0},X) \times \R_{\geq 0}$ 
(cf. p. \pageref{14.72} ff.) 
and take the ``$k$-ification'' of the induced topology.\\
\indent\indent (2) \quad Form $kC(\R_{\geq 0},X) \times_k \R_{\geq 0} = kC(\R_{\geq 0},X) \times \R_{\geq 0}$, 
equip $\Theta_MX \ (\Omega_MX)$ with the induced topology, and pass to its ``$k$-ification''.
\vspi
%%----------------------------------------------------------------------------------------------32
Both procedures yield the same compactly generated topology on \ $\Theta_MX \ (\Omega_MX)$, from which 
$\Theta_{kM}X$ \ $(\Omega_{kM}X)$.\\
\endgroup %%------------------------------------<<

\begingroup%%----------------------------------->>
\fontsize{9pt}{11pt}\selectfont
\textbf{\small EXAMPLE} \ 
Let $X$ be a pointed compactly generated space.  Write $\mo X$ for the nerve of the category associated with the compactly generated monoid $\Omega_MX$ $-$then there is a canonical arrow $\mo X \ra \spx X$ which is a levelwise homotopy equivalence.\\
\endgroup %%------------------------------------<<

Let $X$ be a wellpointed compactly generated space with $\{x_0\} \subset X$ closed.  
Choose a continuous function 
$\phi:X \ra [0,1]$ such that $\phi^{-1}(0) = \{x_0\}$ (cf. $\S 3$, Proposition 21) $-$then the 
\un{meridian map} 
\index{meridian map}
$m:X \ra \Omega_M \Sigma X$ is the pointed continuous function specified by the rule 
$m(x)(t) = [x,t/\phi(x)]$ $(0 \leq t \leq \phi(x))$, where $[x_0,0/0]$ is the base point of $\Sigma X$.  Since 
$\Omega_M \Sigma X$ is a monoid in \bCG, $m$ extends to $JX:$
\begin{tikzcd}[sep=tiny]
{X} \ar{rdddd}[swap]{m} \ar{rr} &&{JX}\ar{ldddd}{J_m}\\
\\
\\
\\
&{\Omega_m \Sigma X}
\end{tikzcd} 
, $J_m$ being the 
\un{arrow of James}. 
\index{arrow of James}

[Note: \  The composite
$X \overset{m}{\ra} \Omega_M \Sigma X \ra \Omega \Sigma X$ is $x \ra [x,-]$.]\\

\begingroup%%----------------------------------->>
\fontsize{9pt}{11pt}\selectfont
Ostensibly, the meridian map depends on $\phi$, call it $m_\phi$.  Suppose, however, that $m_\psi$ is the meridian map corresponding to another continuous function $\psi:X \ra [0,1]$ such that $\psi^{-1}(0) = \{x_0\}$ $-$then 
$m_\phi \simeq m_\psi$.
\vspi
[Let $H:IX \ra \Omega_M \Sigma X$ be the homotopy given by 
$H(x,t): [0,(1-t)\phi(x) + t\psi(x)] \ra \Sigma X$, where 
$H(x,t)(T) = [x,T/((1-t)\phi(x) + t\psi(x))]$.  
Write 
$G:X \ra (\Omega_M \Sigma X)^{[0,1]}$ for its adjoint, view $(\Omega_M \Sigma X)^{[0,1]}$ as a monoid in \bCG, 
determine $\ov{G}$ via the commutative triangle
\begin{tikzcd}[sep=tiny]
{X} \ar{rdddd}[swap]{G} \ar{rr} &&{JX} \ar{ldddd}{\ov{G}}\\
\\
\\
\\
&{(\Omega_M \Sigma X)^{[0,1]}}
\end{tikzcd} 
, and consider its adjoint $\ov{H}:IJX \ra \Omega_M \Sigma X$.]\\
\endgroup %%------------------------------------<<

Let \quad $\Gamma_m:\Gamma X \ra \Theta_M\Sigma X$ \quad be the continuous function defined by the prescription 
\quad
$\Gamma_m([x,t])(T) \ \ \ = \ \  \  [x,T/\phi(x)]$ \ \ $(0 \ \leq \ T \ \leq \ t \phi(x))$   \quad \ \ 
$-$then\ \  there\ \  is \ \ an\ \  arrow 
\\
$\Gamma X \times_k JX$
$\overset{\Gamma_m \times_k J_m}{\xrightarrow{\hspace*{1.5cm}}}$ 
$\Theta_M \Sigma X\  \times_k \ \Omega_M \Sigma X$
$\overset{+}{\xrightarrow{\hspace*{0.75cm}}}$ 
$\Theta_M \Sigma X$
%
%$\overset{+}{\lra} \Theta_M \Sigma X$ 
and a commutative diagram \quad
$
\begin{tikzcd}[sep=small]
{X \times_k JX} \ar{dd} \ar{r} &{JX} \ar{dd} \ar{rd}\\
&&{\Omega_M \Sigma X} \ar{ld}\\
{\Gamma X \times_k JX}  \ar{r} &{\Theta_M \Sigma X  }
\end{tikzcd}
. \ 
$
This leads in turn to an arrow $E \ra \Theta_M\Sigma X$ and a commutative triangle 
\begin{tikzcd}[sep=small]
{E} \ar{rdd} \ar{rr} &&{\Theta_M\Sigma X } \ar{ldd}\\
\\
&{\Sigma X}
\end{tikzcd} 
($\Theta_m\Sigma X \ra \Sigma X$ is the \bCG fibration that evaluates a Moore path at its free end).\\

%%----------------------------------------------------------------------------------------------33
\begin{proposition} \ %19
Let $(X,x_0)$ be a wellpointed compactly generated space with $\{x_0\} \subset X$ closed.  
Assume: $X$ is path connected and numerably contractible $-$then the arrow of James 
$JX \ra \Omega_M \Sigma X$ is a pointed homotopy equivalence.
\end{proposition}

[In the commutative diagram 
$
\begin{tikzcd}%[sep=small]
{\Gamma X \times_k JX} \ar{d}
&{X \times_k JX} \ar{l} \ar{d} \ar{r}
&{JX} \ar{d}\\
{\Gamma X} 
&{X} \ar{l}  \ar{r}
&{*}
\end{tikzcd}
, \ 
$ 
the arrows 
$X \times_k JX \ra \Gamma X \times_k JX$,
$X \ra \Gamma X$ 
are closed cofibrations and 
\begin{tikzcd}%[sep=small]
{\Gamma X \times_k JX} \ar{d}
&{X \times_k JX} \ar{l} \ar{d} \\
{\Gamma X} 
&{X} \ar{l}
\end{tikzcd}
is a homotopy pullback, as is 
\begin{tikzcd}%[sep=small]
{X \times_k JX} \ar{d} \ar{r}
&{JX} \ar{d}\\
{X} \ar{r}
&{*}
\end{tikzcd} 
(the shearing map 
$
\text{sh}:
\begin{cases}
\ X \times_k JX \ra\\ %X \times_k JX\\
\ (x,y) \mapsto %(x,xy)
\end{cases}
$
$
\begin{array}{l}
X \times_k JX\\
(x,xy)
\end{array}
$
is a homotopy equivalence 
(cf. p. \pageref{14.73})).  
Consequently, the sequence 
$JX \ra E \ra \Sigma X$ 
is a fibration up to homotopy (per \bCG (standard structure) 
(cf. p. \pageref{14.74})).    
Since \mE is contractible, it remains only to consider the commutative diagram
%$
%\begin{tikzcd}%[sep=small]
%{JX} \ar{d} \ar{r}
%&{E} \ar{d} \ar{r}
%&{\Sigma X} \arrow[d,shift right=0.5,dash] \arrow[d,shift right=-0.5,dash]\\
%{\Omega_M \Sigma X} \ar{r}
%&{\Theta_M \Sigma X}  \ar{r}
%&{\Sigma X}
%\end{tikzcd}
%.]
%$
$
\begin{tikzcd}%[sep=small]
{JX} \ar{d} \ar{r}
&{E} \ar{d} \\
{\Omega_M \Sigma X} \ar{r}
&{\Theta_M \Sigma X} 
\end{tikzcd}
$
$
\begin{tikzcd}%[sep=small]
{} \ar{r}
&{\Sigma X} \arrow[d,shift right=0.5,dash] \arrow[d,shift right=-0.5,dash]\\
{}  \ar{r}
&{\Sigma X}
\end{tikzcd}
.]
$
\\
\vspace{0.25cm}

Application: Under the hypotheses of Proposition 19, the composite
$JX \overset{J_m}{\lra} \Omega_M \Sigma X \ra \Omega \Sigma X$ 
is a pointed homotopy equivalence.\\

\begin{proposition} \ %20
Let $(X,x_0)$ be a wellpointed compactly generated space with $\{x_0\} \subset X$ closed.  
Assume: $X$ is path connected $-$then the arrow of James 
$JX \ra \Omega_M \Sigma X$ is a weak homotopy equivalence.
\end{proposition}

[Thanks to the cone construction 
(cf. p. \pageref{14.75} ff.), 
the arrow $E \ra \Sigma X$ is a quasi-fibration.  
Work with 
\begin{tikzcd}%[sep=small]
{JX} \ar{d} \ar{r}
&{E} \ar{d} \ar{r}
&{\Sigma X} \arrow[d,shift right=0.5,dash] \arrow[d,shift right=-0.5,dash]\\
{\Omega_M \Sigma X} \ar{r}
&{\Theta_M \Sigma X}  \ar{r}
&{\Sigma X}
\end{tikzcd}
and compare the long exact sequences of homotopy groups.]

[Note: \  In the case at hand, $\Sigma X$ is simply connected.]\\

Application: Under the hypotheses of Proposition 20, the composite
$JX \overset{J_m}{\lra} \Omega_M \Sigma X \ra \Omega \Sigma X$ 
is a weak homotopy equivalence.\\

\begingroup%%----------------------------------->>
\fontsize{9pt}{11pt}\selectfont
\textbf{\small EXAMPLE} \ 
Let $X$ be the broom pointed at $(0,0)$ $-$then $X$ is path connected.  But $JX$ and $\Omega\Sigma X$ do not have the same weak homotopy type ($\Sigma X$ is not simply connected).\\
\endgroup %%------------------------------------<<

\begin{proposition} \ %21
Let 
$
\begin{cases}
\ (X,x_0)\\
\ (Y,y_0)
\end{cases}
$
be wellpointed compactly generated spaces with 
$
\begin{cases}
\ {x_0} \subset X\\
\ {y_0} \subset Y
\end{cases}
$
closed and let $f:X \ra Y$ be a pointed continuous function.  Assume: $f$ is a 
%%----------------------------------------------------------------------------------------------34
homotopy equivalence (weak homotopy equivalence) $-$then 
$Jf:JX \ra JY$ is a homotopy equivalence (weak homotopy equivalence).
\end{proposition}

[Arguing by induction from 
$
\begin{tikzcd}%[sep=small]
{X^n} \ar{d} 
&{X_*^n} \ar{l} \ar{d} \ar{r}
&{J^{n-1}X} \ar{d}\\
{Y^n} 
&{Y_*^n} \ar{l} \ar{r}
&{J^{n-1}Y} 
\end{tikzcd}
, \ 
$ 
one finds that $\forall \ n$, 
$J^nX \ra J^nY$ is a homotopy equivalence  
(cf. p. \pageref{14.76} ff.) 
(weak homotopy equivalence) 
(cf. p. \pageref{14.77})), 
hence 
$JX \ra JY$ is a homotopy equivalence (cf. $\S 3$, Proposition 15) (weak homotopy equivalence 
(cf. p. \pageref{14.78})).]\\

Convention: Given a cofibered monoid $G$ in \bCG, 
$\Sigma G \ra BG$ is the adjoint of $G \ra \Omega BG$ 
(cf. p. \pageref{14.79}).\\

\label{14.152} %dmc ??
\textbf{\small LEMMA} \ 
Let $(X,x_0)$ be a wellpointed compactly generated space with $(x_0) \subset X$ closed.  
Assume: $X$ is discrete $-$then the composite 
$\Sigma X \ra \Sigma JX \ra BJX$ is a weak homotopy equivalence.

[Since $X = \bigvee\limits_{X - \{x_0\}} \bS^0$, $JX = \coprod\limits_{X - \{x_0\}} J\bS^0$ \ 
($\coprod$ the coproduct in the category of monoids), 
where $J\bS^0 = \Z_{\geq 0}$, thus it suffices to consider \ %${\Sigma \bS^0}$
\begin{tikzcd}[sep=small]
&{\Sigma \Z_{\geq 0}} \ar{r}  &{B\Z_{\geq 0}} \ar{dd}\\
{\Sigma \bS^0} \ar{ru} \ar{rd}\\
%{} \ar{ru} \ar{rd}\\
&{\Sigma \Z} \ar{r}  &{B\Z}
\end{tikzcd}
.]\\

\begin{proposition} \ %22
Let $(X,x_0)$ be a wellpointed compactly generated space with $\{x_0\} \subset X$ closed $-$then the composite 
$\Sigma X \ra \Sigma JX \ra BJX$ is a weak homotopy equivalence.
\end{proposition}

[The lemma implies that \ $\forall \ n$, the composite \ \ 
$\Sigma \sin_n X \lra \Sigma J \sin_n X \lra BJ\sin_n X$ 
is a weak homotopy equivalence ($\sin_n X$ being supplied with the discrete topology), thus the composite \  
$\abs{n \ra \Sigma \sin_n X} \ra \abs{n \ra \Sigma J\sin_n X} \ra \abs{n \ra BJ \sin_n X}$ \ 
is a weak homotopy equivalence 
(cf. p. \pageref{14.80}).  
But \ 
$\abs{n \ra \Sigma \sin_n X} \approx \Sigma \abs{\sin X}$ 
(cf. p. \pageref{14.81} ff.), 
$\abs{n \ra \Sigma J\sin_n X} \ \approx$ \ 
$\Sigma \abs{n \ra J \sin_n X} \ \approx$ \ 
$\Sigma J \abs{\sin X}$ \ 
(cf. p. \pageref{14.82}), 
$\abs{n \ra BJ \sin_n X} \  \approx$ \ 
$B\abs{n \ra J \sin_n X}$ 
(cf. p. \pageref{14.83}) 
$\approx BJ \abs{\sin X}$ \ 
and there is a commutative diagram 
$
\begin{tikzcd}%[sep=small]
{\Sigma \abs{\sin X}} \ar{d} \ar{r} 
&{\Sigma J\abs{\sin X}} \ar{d} \ar{r}
&{BJ\abs{\sin X}} \ar{d}\\
{\Sigma X} \ar{r} 
&{\Sigma JX} \ar{r} 
&{BJX}
\end{tikzcd}
.
$  
The arrow $\Sigma \abs{\sin X} \ra \Sigma X$ is a weak homotopy equivalence (cf. infra).  
According to Proposition 21, the same holds for the arrow 
$J \abs{\sin X} \ra JX$ or still, for the arrows 
$\Sigma J \abs{\sin X} \ra \Sigma JX$, 
$BJ \abs{\sin X}  \ra BJX$ 
(cf. p. \pageref{14.84}).  
Combining these facts yields the assertion.]\\

\begingroup%%----------------------------------->>
\fontsize{9pt}{11pt}\selectfont
\textbf{\small LEMMA} \ 
Let 
$
\begin{cases}
\ (X,x_0)\\
\ (Y,y_0)
\end{cases}
, (Z,z_0) 
$\
be wellpointed compactly generated spaces with 
$
\begin{cases}
\ \{x_0\} \subset X\\
\ \{y_0\} \subset Y
\end{cases}
%%----------------------------------------------------------------------------------------------35
, % for spaceing (z_0) \subset Z
$
$ \{z_0\} \subset Z$ 
closed and let $f:X \ra Y$ be a pointed continuous function.  Assume: $f$ is a weak homotopy equivalence $-$then 
$f\#_k\id_Z:X\#_kZ \ra Y\#_kZ$ is a weak homotopy equivalence.\\
\endgroup %%------------------------------------<<

\begingroup%%----------------------------------->>
\fontsize{9pt}{11pt}\selectfont
Application: Let 
$
\begin{cases}
\ (X,x_0)\\
\ (Y,y_0)
\end{cases}
$
be wellpointed compactly generated spaces with 
$
\begin{cases}
\ \{x_0\} \subset X\\
\ \{y_0\} \subset Y
\end{cases}
$
closed and let $f:X \ra Y$ be a pointed continuous function.  Assume: $f$ is a weak homotopy equivalence $-$then 
$\Sigma f: \Sigma X \ra \Sigma Y$ is a weak homotopy equivalence.
\vspi
\label{16.63}
[Note: \  Recall too that $\Omega f:\Omega X \ra \Omega Y$ is a weak homotopy equivalence 
(cf. p. \pageref{14.85}).]\\
\endgroup %%------------------------------------<<

\textbf{\small LEMMA} \ 
Let $(X,x_0)$ be a wellpointed compactly generated space with $\{x_0\} \subset X$ closed $-$then there is a canonical arrow 
$B \Omega_M X \ra X$ and a commutative diagram 
$
\begin{tikzcd}%[sep=small]
{\Sigma \Omega_M X} \ar{d} \ar{r} &{B \Omega_M X} \ar{ld}\\
{X}
\end{tikzcd}
.
$

[Note: \  $B \Omega_M X \ra X$ is a weak homotopy equivalence provided that $X$ is path connected 
(cf. p. \pageref{14.86}).]\\

\begin{proposition} \ %23
Let $(X,x_0)$ be a wellpointed compactly generated space with $\{x_0\} \subset X$ closed $-$then the arrow of James 
$J_m:JX \ra \Omega_M \Sigma X$ induces a weak homotopy equivalence 
$BJ_m:BJX \ra B\Omega_M \Sigma X$.
\end{proposition}

[The composite 
$\Sigma X \overset{\Sigma m}{\lra} \Sigma\Omega_M \Sigma X \ra \Sigma X$ is $\id_{\Sigma X}$.  
Proof: $[x,t] \ra [m(x),t] \ra m(x)$ 
$(t\phi(x)) = [x,t\phi(x)/\phi(x)] = [x,t]$.  With this in mind the commutative diagram
\[
\begin{tikzcd}%[sep=large]
&{\Sigma JX}  \ar{d} \ar{r} &{BJX} \ar{d}{BJ_m}\\
{\Sigma X} \ar{ru} \ar{r} \ar{rd} &{\Sigma\Omega_M \Sigma X} \ar{d} \ar{r} &{B\Omega_M \Sigma X} \ar{ld}\\
&{\Sigma X}
\end{tikzcd}
\]
shows that 
%$\Sigma X \ra \Sigma JX \ra BJX \overset{BJ_m}{\lra} B\Omega_M\Sigma X \ra \Sigma X$ is also $\id_{\Sigma X}$.  
$\Sigma X$ 
$\ra$ 
$\Sigma JX$ 
$\ra$ 
$BJX$ 
$\overset{BJ_m}{\xrightarrow{\hspace*{1cm}}}$ 
$B\Omega_M\Sigma X$ 
$\ra$ 
$\Sigma X$
is also $\id_{\Sigma X}$.  
On account of Proposition 22, the composite
$\Sigma X \ra \Sigma JX \ra BJX$ is a weak homotopy equivalence.  
However $\Sigma X$ is path connected, hence 
$B\Omega_M\Sigma X \ra \Sigma X$  is a weak homotopy equivalence.  Therefore 
$BJ_m:BJX \ra  B\Omega_M\Sigma X$ is a weak homotopy equivalence.]

[Note: \ One can view Proposition 23 as the $\pi_0(X) \neq *$ analog of Proposition 20.]\\

\begingroup%%----------------------------------->>
\fontsize{9pt}{11pt}\selectfont
\textbf{\small FACT} \ 
Let $(X,x_0)$ be a wellpointed compactly generated space with $\{x_0\} \subset X$ closed.  
Assume: $X$ is \dsep and $\Delta_X \ra X \times_k X$ a cofibration.  Put $GX = F_{\gr}(X,x_0)$ 
(cf. p. \pageref{14.87}) 
$-$then the arrow $BJX \ra BGX$ is a weak homotopy equivalence.
\vspi
%%----------------------------------------------------------------------------------------------36
[Note: \  It follows that the arrow $JX \ra GX$ is a weak homotopy equivalence whenever $X$ is path connected 
(cf. p. \pageref{14.88}).]\\
\endgroup %%------------------------------------<<

It is also of interest to consider the free abelian compactly generated monoid on $(X,x_0)$, denoted by $SP^\infty X$ 
and referred to as the 
\un{infinite symmetric product} 
\index{infinite symmetric product} 
on $(X,x_0)$.  \ 
Like $JX$, $SP^\infty X$ carries the quotient topology coming from 
$\coprod\limits_n X^n \ra SP^\infty X$.  
Put $SP^n X = p(X^n)$ $-$then if $\{x_0\}$ is closed in $X$, $SP^n X$ is closed in 
$SP^\infty X$ and the arrow $X^n \ra SP^n X$ is quotient, hence $SP^\infty X = \colimx SP^n X$ and 
$X^n /  S_n \approx SP^n X$.  \ 
Example: $SP^\infty \bS^0 \approx \Z_{\geq 0}$.\\

\begingroup%%----------------------------------->>
\fontsize{9pt}{11pt}\selectfont
Under certain conditions, it is possible to indentify $X^n / S_n$.  For instance, $\bS^2 / S_n$ is homeomorphic to 
$\bP^n(\C)$, therefore $SP^\infty \bS^2$ is homeomorphic to $\bP^\infty(\C)$ a $K(\Z,2)$ 
(cf. p. \pageref{14.89}).
\vspi
[Note: \  A survey of this aspect of the theory has been given by 
Wagner\footnote[2]{\textit{Dissertationes Math.} \textbf{182} (1980), 1-52.}.]\\
\endgroup %%------------------------------------<<

\begingroup%%----------------------------------->>
\fontsize{9pt}{11pt}\selectfont
\textbf{\small Example} \ 
Let $X$ be a compact metric space with $\dim X < \infty$.  Assume: $X$ is an ANR $-$then $X^n / S_n$ is an ANR 
(Floyd\footnote[3]{\textit{Duke Math. J.} \textbf{22} (1955), 33-38.}).\\
\endgroup %%------------------------------------<<

\begin{proposition} \ %24
Let $(X,x_0)$ be a wellpointed compactly generated space with $\{x_0\} \subset X$ closed $-$then $(SP^\infty X,x_0)$ is a wellpointed compactly generated space with 
$\{x_0\} \subset SP^\infty X$ closed, thus is an abelian cofibered monoid in \bCG.\\
\end{proposition}

\begingroup%%----------------------------------->>
\fontsize{9pt}{11pt}\selectfont
\textbf{\small LEMMA} \ 
If $(X,x_0)$ is a wellpointed compactly generated Hausdorff space, then $(SP^\infty X,x_0)$ is a wellpointed compactly generated Hausdorff space.\\
\endgroup %%------------------------------------<<

\label{14.95}
\begingroup%%----------------------------------->>
\fontsize{9pt}{11pt}\selectfont
\textbf{\small FACT} \  
Suppose that $(X,x_0)$ is a pointed CW complex $-$then $(SP^\infty X,x_0)$ is a pointed CW complex.
\vspi
[It is enough to place a CW structure on each $SP^n X$ 
in such a way that $SP^{n-1} X$ is a subcomplex of $SP^n X$
(cf. p. \pageref{14.90}).  
For this, it is necessary to alter the CW structure on $X^n$ in order to reflect the action of $S_n$.]\\
\endgroup %%------------------------------------<<

\begin{proposition} \ %25
Let $(X,x_0)$ be a wellpointed compactly generated space with $\{x_0\} \subset X$ closed $-$then  there is an isomorphism 
$BSP^\infty X \approx SP^\infty \Sigma X$ of abelian monoids in \bCG.
\end{proposition}

%%----------------------------------------------------------------------------------------------37
[Analogously, $X S P^\infty X \approx S P^\infty \Gamma X$ and the diagram
\[
\begin{tikzcd}%[sep=small]
{S P^\infty X}   \arrow[d,shift right=0.5,dash] \arrow[d,shift right=-0.5,dash] \ar{r} 
&{X S P^\infty X} \ar{d} \ar{r} 
&{B S P^\infty X}\ar{d} \\
{S P^\infty X} \ar{r} &{S P^\infty \Gamma X} \ar{r}  &{S P^\infty \Sigma X}
\end{tikzcd}
\]
commutes.]\\

\begin{proposition} \ %26
Let $(X,x_0)$ be a wellpointed compactly generated space with $\{x_0\} \subset X$ closed.  
Assume: $X$ is path connected and numerably contractible $-$then the arrow 
$SP^\infty X \ra \Omega B S P^\infty X$ is a pointed homotopy equivalence.
\end{proposition}

[$\forall \ n$, $SP^nX$ is numerably contractible, so $SP^\infty X = \colimx SP^nX$ is numerably contractible 
(cf. p. \pageref{14.91}).  
Since the inclusion $\{x_0\} \ra SP^\infty X$ is a closed cofibration and $SP^\infty X$ is path connected, 
it follows that $SP^\infty X$ admits a homotopy inverse 
(cf. p. \pageref{14.92}).  
Therefore the arrow 
$SP^\infty X \ra \Omega B SP^\infty X$ is a pointed homotopy equivalence (cf. Proposition 14).]\\

Application: Under the hypotheses of Proposition 26, the composite \ 
$SP^\infty X$ $\ra$ $\Omega B SP^\infty X$ $\ra$ $\Omega SP^\infty \Sigma X$ is a pointed homotopy equivalence.\\

\index{Theorem: Dold-Thom Theorem}
\index{Dold-Thom Theorem}
\textbf{\small DOLD-THOM THEOREM} \quad
Suppose that $(X,x_0)$ is a pointed connected CW complex $-$then $\forall \ n > 0$, 
$\pi_n(SP^\infty X) \approx H_n(X)$.

[There are pointed homotopy equivalences 
$\abs{SP^\infty \sin X} \ra SP^\infty \abs{\sin X}$, \ 
$SP^\infty \abs{\sin X} \ra SP^\infty X$.  
One has $\widetilde{H}_*(\abs{\sin X}) \approx \widetilde{H}_*(X)$ 
and, in the notation of 
p. \pageref{14.93}, 
$\pi_*(F_{\ab}(\sin X,x_0)) \approx \widetilde{H}_*(\abs{\sin X})$ 
(Weibel\footnote[2]{\textit{An Introduction to Homological Algebra}, Cambridge University Press (1994), 266-267.}).  
But
$\ov{SP^\infty \sin X} = F_{\ab}(\sin X,x_0)$, thus the arrow 
$\abs{SP^\infty \sin X} \ra \abs{F_{\ab}(\sin X,x_0)}$ \ \ 
is \ a \ pointed \ homotopy \ equivalence \ \ 
(cf, p. \pageref{14.94}).  
Accordingly, \quad
$\pi_*(\abs{SP^\infty \sin X})$ $\approx$ 
$\pi_*(\abs{F_{\ab}(\sin X,x_0)})$ $\approx$ 
$\pi_*(F_{\ab}(\sin X,x_0))$, from which the assertion.]\\

\begingroup%%----------------------------------->>
\fontsize{9pt}{11pt}\selectfont
\textbf{\small EXAMPLE} \ 
Dold-Thom can fail if $X$ is not a CW complex.  
Example: Take for $X$ the Hawaiian earring pointed at $(0,0)$, form its cone 
$\Gamma X$ and consider 
$\Gamma X \vee \Gamma X$ $-$then 
$H_1(\Gamma X \vee \Gamma X) \neq 0$, so either 
$\pi_1(SP^\infty \Gamma X) \neq H_1(\Gamma X)$ or 
$\pi_1(SP^\infty(\Gamma X \vee \Gamma X)) \neq H_1(\Gamma X \vee \Gamma X)$.\\
\endgroup %%------------------------------------<<

Remark: If $(X,x_0)$ is a pointed connected CW complex, then $(SP^\infty X, x_0)$ is a pointed connected CW complex 
(cf. p. \pageref{14.95}) and 
$SP^\infty X \approx (w) \prod\limits_1^\infty K(\pi_n(SP^\infty X),n))$ 
(cf. p. \pageref{14.95a}) 
or still, by the Dold-Thom theorem, 
$SP^\infty X \approx (w) \prod\limits_1^\infty K(H_n(X),n)$.\\

%%----------------------------------------------------------------------------------------------38
\label{14.89}
\begingroup%%----------------------------------->>
\fontsize{9pt}{11pt}\selectfont
\textbf{\small EXAMPLE} \ 
Let $\pi$ be an abelian group and let $X = M(\pi,n)$ (realized as a pointed connected CW complex) $-$then 
$SP^\infty M(\pi,n)$ is a $K(\pi,n)$.  
In particular, $SP^\infty \bS^n$ is a $K(\Z,n)$.\\
\endgroup %%------------------------------------<<

$\bGamma_{\ini}$ is the category whose objects are the finite sets $\bn \equiv \{0, 1, \ldots, n\}$ $(n \geq 0)$ 
with base point 0 and whose morphisms are the base point preserving injective maps.

Example: Let $(X,x_0)$ be a wellpointed compactly generated space with $\{x_0\} \subset X$ closed.  
Viewing $X^n$ as the space of base point preserving continuous functions $\bn \ra X$, define a functor 
$\pow X: \bGamma_{\ini} \ra \bCG_*$ by writing $\pow_n X = X^n$, stipulating that the arrow $X^m \ra X^n$ 
attached to $\gamma:\bm \ra \bn$ sends $(x_1, \ldots, x_m)$ to $(\bar{x}_1, \ldots, \bar{x}_n)$, where
$\bar{x}_j = x_{\gamma^{-1}(j)}$ if $\gamma^{-1}(j) \neq \emptyset$, 
$\bar{x}_j = x_0$ if $\gamma^{-1}(j) = \emptyset$.

[Note: \  $\colimx \pow X$ can be identified with $SP^\infty X$.]\\

\label{14.107}
\label{14.119}
\begingroup%%----------------------------------->>
\fontsize{9pt}{11pt}\selectfont
\textbf{\small EXAMPLE} \ 
For $n > 0$, $\colimx \pow \bn \approx SP^\infty \bn \approx \Z_{\geq 0} \times \cdots \times \Z_{\geq 0}$ ($n$ factors).  
On the other hand, $\ohc \ \pow \bn$ has the homotopy type of 
$\B \bM_\infty \times_k \cdots \times_k B\bM_\infty$ ($n$ factors), $\bM_\infty$ the permutative category on 
p. \pageref{14.96} 
(so $B\bM_\infty = \ds\coprod\limits_{n \geq 0} BS_n$).
\vspi
[Note: \  $\colimx \pow \bzero \approx SP^\infty \bzero \approx \{0\}$ while
$\ohc \ \pow \bzero \approx  B \bGamma_{\ini}$, is a contractible space 
(cf. p. \pageref{14.97}).]\\
\endgroup %%------------------------------------<<

Definition: A 
\un{creation operator} 
\index{creation operator} is a functor 
$\sC:\bGamma_{\ini}^\OP \ra \bCG$ such that $\sC_0 = *$.

[Note: \  $\forall \ n$, $\sC_n$ is a right $S_n$-space.]\\

\label{14.105}
\begingroup%%----------------------------------->>
\fontsize{9pt}{11pt}\selectfont
\textbf{\small EXAMPLE} \ 
Every nonempty compactly generated Hausdorff space $Y$ gives rise to a creation operator $\sC Y$ whose $n^{\thx}$ space is 
$Y^n$ $(Y^0 = *)$, the arrow $Y^n \ra Y^m$ being determined by $\gamma:\bm \ra \bn$ being the map 
$(y_1, \ldots, y_n) \ra (y_{\gamma(1)}, \ldots, y_{\gamma(m)})$.\\
\endgroup %%------------------------------------<<

\label{16.58}
If $\sC$ is a creation operator and if $(X,x_0)$ is a wellpointed compactly generated space  with $\{x_0\} \subset X$ closed, then the 
\un{realization}
\index{realization (creation operator)} 
$\sC[X]$ of $\sC$ at $X$ is 
$\ds\int^{\bn} \sC_n \times_k X^n$ $(= \sC \otimes_{\bGamma_\ini} \pow X)$.  
Example: Suppose that $\sC_n = *$ $\forall \ n$ $-$then 
$\sC[X] \approx *\otimes_{\bGamma_\ini} \pow X$ $\approx$ 
$\colimx \pow X$ $\approx$ $SP^\infty X$.\\

\begingroup%%----------------------------------->>
\fontsize{9pt}{11pt}\selectfont
\textbf{\small EXAMPLE} \ 
Let $\sC_n = S_n$ $\forall \ n$.  Given a morphism $\gamma:\bm \ra \bn$ in $\bGamma_\ini$, specify
$\sC_\gamma:S_n \ra S_m$ as follows: $\forall \ \sigma \in S_n$, there exists a unique order preserving injection 
$\gamma^\prime:\bm \ra \bn$ such that $\gamma^\prime(\bm) = (\sigma \circx \gamma)(\bm)$ and 
$(\sC \gamma)\sigma \in S_m$ is the permutation for which 
$\gamma^\prime \circx (\sC \gamma) \sigma = \sigma \circx \gamma$.  
This data thus defines a creation operator and 
$\forall \ X$, $\sC[X] \approx JX$.\\
\endgroup %%------------------------------------<<

\begin{proposition} \ %27
Suppose that  $(X,x_0)$ is a wellpointed compactly generated space with $\{x_0\} \subset X$ closed and let $\sC$ be a creation operator.  
Denote by $\sC_n[X]$ the image of 
$\ds\coprod\limits_{m \leq n} \sC_m \times_k X^m$ in $\sC[X]$ $-$then $\sC_n[X]$ is a closed subspace of $\sC[X]$ and 
$\sC[X] = \colimx \sC_n[X]$.\\
%
%%----------------------------------------------------------------------------------------------39
In addition, the commutative diagram 
\begin{tikzcd}%[sep=small]
{\sC_n \times_{S_n} X_*^n} \ar{d} \ar{r} &{\sC_{n-1}[X]} \ar{d}\\
{\sC_n \times_{S_n} X^n} \ar{r} &{\sC_{n}[X]}
\end{tikzcd}
is a pushout square and the arrow $\sC_{n-1}[X] \ra \sC_n[X]$ is a closed cofibration.
\end{proposition}

[Note: \  The base point of $\sC[X]$ is $[*,x_0]$ and the inclusion $\{[*,x_0]\} \ra \sC[X]$ is a closed cofibration.]\\

\label{14.169}
Remark: $X$ \dsep $+ \ \sC_n$ \dsep $\forall \ n$ $\implies$ $\sC[X]$ \dsep.\\

\label{14.145} %dmc mnft

\begingroup%%----------------------------------->>
\fontsize{9pt}{11pt}\selectfont
The validation of the above remark depends on Proposition 27 and the following lemma.\\
\endgroup %%------------------------------------<<

\begingroup%%----------------------------------->>
\fontsize{9pt}{11pt}\selectfont
\textbf{\small LEMMA} \ 
Let \mG be a compact Hausdorff topological group.  Suppose that \mX is a $\Delta$-separated right $G$-space $-$then 
$X/G$ is $\Delta$-separated.
\vspi
[It is a matter of proving that 
$\{(x,x\cdot g): x \in X, g \in G\}$ is closed in $X \times_k X$ 
(cf. p. \pageref{14.97a}).  
However, \mG acts to the right on $X \times_k X$, viz. $(x,y)\cdot g = (x,y\cdot g)$, and $\Delta_X$ is closed in 
$X \times_k X$, hence $\Delta_X \cdot G$ is closed in $X \times_k X$, \mG being compact Hausdorff.]\\
\endgroup %%------------------------------------<<

\label{14.127}
\begingroup%%----------------------------------->>
\fontsize{9pt}{11pt}\selectfont
\textbf{\small FACT} \ 
Let 
$
\begin{cases}
\ (X,x_0)\\
\ (Y,y_0)
\end{cases}
$
be wellpointed compactly generated spaces with 
$
\begin{cases}
\ \{x_0\} \subset X\\
\ \{y_0\} \subset Y\\
\end{cases}
$
closed; let $f:X \ra Y$ be a pointed continuous function.  
Assume: $f$ is a closed cofibration $-$then $\forall$ creation operator 
$\sC$, the induced map $\sC[X] \ra \sC[Y]$ is a closed cofibration.
\vspi
[Use the lemma on 
p. \pageref{14.98} ff. 
and the lemma on 
p. \pageref{14.99}.]
\vspi
[Note: \  The conclusion of the lemma on 
p. \pageref{14.100} ff. 
is ``closed cofibration'' rather than just ``cofibration'' 
provided that this is so of the vertical arrow on the right in the hypothesis.  To see this, observe that the argument there can be repeated, testing against any arrow $Z \ra B$ which is both a homotopy equivalence and a Hurewicz fibration 
cf. p. \pageref{14.101}).]\\
\endgroup %%------------------------------------<<

\begin{proposition} \ %28
Let $\phi:\sC \ra \sD$ be a morphism of creation operators.  Assume: $\forall \ n$, $\phi_n:\sC_n \ra \sD_n$ is an $S_n$ equivarariant homotopy equivalence $-$then $\phi$ induces a homotopy equivalence $\sC[X] \ra \sD[X]$.\\
\end{proposition}

By definition, $\ohc \ \pow X \approx B(-\backslash \bGamma_\ini) \otimes_{\bGamma_\ini} \pow X$.  \ 
Problem: Exhibit models for $\ohc \ \pow X$ in the homotopy category.

[Note: \  Strictly speaking, $B(-\backslash\bGamma_\ini)$ is not a creation operator 
(since 
$B(\bzero\backslash\bGamma_\ini) \neq *$).]

A compactly generated paracompact Hausdorff space $X$ is said to be \un{$S_n$-universal} if it is a contractible free right 
$S_n$-space.  The covering projection $X \ra X/S_n$ is then a closed map, 
hence $X/S_n$ is a compactly generated paracompact Hausdorff space.  
Therefore $X/S_n$  is a classifying space for $S_n$ (in the sense of 
p. \pageref{14.102}).  
Examples:
(1) $X_{S_n}^\infty$ is $S_n$-universal; 
(2) $B(\bn\backslash \bGamma_\ini)$ is $S_n$-universal; 
(3) $X S_n$ is $S_n$-universal.\\ 

%%----------------------------------------------------------------------------------------------40
A creation operator $\sC$ is said to be 
\un{universal} 
\index{universal (creation operator)} 
 if $\forall \ n$, $\sC_n$ is $S_n$-universal.

Example: Let $\sC$ be a univeral creation operator $-$then for any cofibered monoid $G$ in \bCG, 
$\sC_n \times_{S_n} (BG)^n$ has the same homotopy type as $B\bigl(S_n \ds\int G \bigr)$ 
(cf. p. \pageref{14.103}).\\

\begin{proposition} \ 
Suppose $\sC$ is a univeral creation operator $-$then there exists an arrow 
$B(-\backslash \bGamma_\ini) \ra \sC$ such that $\forall \ n$, $B(\bn\backslash\bGamma_\ini) \ra \sC_n$ is an $S_n$-equivariant homotopy equivalence.
\end{proposition}

[In the notation of 
p. \pageref{14.104} ff., 
compose the homotopy equivalence $B(-\backslash\bGamma_\ini) \ra P\sC$ and the arrow of evaluation 
$P\sC \ra \sC$.]\\

Application: Let $(X,x_0)$ be a wellpointed compactly generated space with $\{x_0\} \subset X$ closed $-$then $\forall$ 
universal creation operator $\sC$, $\sC[X]$ and $\ohc \ \pow X$ have the same homotopy type.\\

\begingroup%%----------------------------------->>
\fontsize{9pt}{11pt}\selectfont
\textbf{\small FACT} \ 
Let $\phi:\sC \ra \sD$ be a morphism of creation operators.  Assume: $\sC$ and $\sD$ are universal $-$then $\phi$ induces 
a homotopy equivalence $\sC[X] \ra \sD[X]$.\\
\endgroup %%------------------------------------<<

Given a nonempty compactly generated Hausdorff space $Y$, let $F(Y,n)$ be the subspace of $Y^n$ consisting of those 
$n$-tuples $(y_1, \ldots, y_n)$ such that $i \neq j$ $\implies$ $y_i \neq y_j$ $-$then $F(Y,n)$ is open in $Y^n$, hence is a compactly generated Hausdorff space, and $S_n$ operates freely on the right by permuting coordinates.

[Note: \  $F(Y,n)$ is the 
\un{configuration space}
\index{configuration space (compactly generated Hausdorff space)} 
of $n$-tuples of distinct points in $Y$.  
Consult Cohen\footnote[2]{\textit{J. Pure Appl. Algebra} \textbf{100} (1995), 19-42.}
for additional information and references.]

Notation: $\con\  Y$ 
\index{con $Y$} 
is the creation operator that sends $n$ to $F(Y,n)$, the arrow 
$F(Y,n) \ra F(Y,m)$ determined by $\gamma:\bm \ra \bn$ being the map 
$(y_1, \ldots, y_n) \ra (y_{\gamma(1)}, \ldots y_{\gamma(m)})$. 

[Note: \  Therefore $\con\  Y$ is a subfunctor of $\sC Y$ 
(cf. p. \pageref{14.105}).]

Observation: The points of $\con\  Y[X]$ are equivalence classes of pairs $(S,f)$, where $S \subset Y$ is a finite subset of $Y$, $f:S \ra X$ is a function, and $(S,f) \sim (S - \{y\},\restr{f}{S - \{y\}})$ iff $f(y) = x_0$.

[Note: \  All pairs $(S,f)$, where $f(S) = \{x_0\}$, are identified with $(\emptyset,\emptyset)$.]

Examples: 
(1) $\con \ \R^0[X] \approx X$; 
(2) $\con \ Y[\bS^0] \approx \{S \subset Y:\#(S) < \omega\}$.\\

\textbf{\small LEMMA} \ 
$\F(\R^\infty,n)$ is $S_n$-universal.

[$(\R^\infty)^n$ is a polyhedron.  
But $F(\R^\infty,n)$ is an open subset of $(\R^\infty)^n$, thus it too is a polyhedron
(cf. p. \pageref{14.106}).  
Therefore $F(\R^\infty,n)$ is a compactly generated paracompact
%%----------------------------------------------------------------------------------------------41
Hausdorff space.  Contractibility is clear if $n = 0$ or 1, so take $n \geq 2$ and represent $F(\R^\infty,n)$  as 
$\colimx F(\R^q,n)$.  Since for $q \gg 0$, $F(\R^q,n)$ is the complement in $\R^{qn}$ of certain hyperplanes of codimension 
$q$, $F(\R^q,n)$ is $(q - 2)$ connected, and this implies that $F(\R^\infty,n)$ is contractible.]\\

\begin{proposition} \ %30
$\con\  \R^\infty$ is a universal creation operator.\\
\end{proposition}

Application: Let $(X,x_0)$ be a wellpointed compactly generated space with $\{x_0\} \subset X$ closed $-$then 
$\ohc \ \pow X$ and $\con \ \R^\infty[X]$ have the same homotopy type.\\

\begingroup%%----------------------------------->>
\fontsize{9pt}{11pt}\selectfont
\textbf{\small EXAMPLE} \ 
$\con\  \R^\infty[\bS^0] \approx$ 
$\ds\coprod\limits_{n \geq 0} F(\R^\infty,n) /  S_n \approx$ 
$\ds\coprod\limits_{n \geq 0} B S_n$, 
which agrees with the fact that the homotopy type of $\ohc \ \pow \bS^0$ is $B\bM_\infty$ 
(cf. p. \pageref{14.107}).\\
\endgroup %%------------------------------------<<

A $q$-dimensional rectangle $[0,1]^q$ is a product of the form 
$R = [a_1,b_1] \times \cdots \times [a_q,b_q]$, where $0 \leq a_i < b_i \leq 1$.  
Call $R(q)$ the set of such and topologize it as a subspace of $[0,1]^{2q}$.  
Note that there is a closed embedding $R(q) \ra R(q+1)$ defined by multiplication on the right by $[0,1]$ and put 
$R(\infty) = \colimx R(q)$.  
Let $\xBV(R(q),n)$ be the subspace of $F(R(q),n)$ consisting of those $n$-tuples $(R_1, \ldots, R_n)$ with the property that the interior of $R_i$ does not meet the interior of $R_j$ if $i \neq j$ $-$then there is a closed embedding 
$\xBV(R(q),n) \ra \xBV(R(q+1),n)$ and 
$\xBV(R(\infty),n) = \colimx \xBV(R(q),n)$ is a free right $S_n$-space.

Notation $\xBV^\infty$ is the creation operator that sends $n$ to $\xBV(R(\infty),n)$.\\

\textbf{\small LEMMA} \ 
$\xBV(R(\infty),n)$ is $S_n$-universal.

[It follows from condition C on 
p. \pageref{14.108} 
that $\xBV(R(\infty),n)$ is a compactly generated paracompact Hausdorff space.  
Since the closed embedding 
$\xBV(R(q),n) \ra \xBV(R(q+1),n)$ 
is a cofibration, one need only establish that it is also inessential in order to conclude that 
$\xBV(R(\infty),n)$
is contractible 
(cf. p. \pageref{14.109}).  
To define 
$H:I\xBV(R(q),n) \ra \xBV(R(q+1),n)$, 
represent an $n$-tuple $(R_1, \ldots, R_n)$ by a $2n$-tuple 
$(A_1,B_1, \ldots, A_n,B_n)$ of points in $[0,1]^q$.  
Here $R_k \leftrightarrow (A_k,B_k)$ and 
$A_k = (a_{k_1},\ldots,a_{k_q})$,
$B_k = (b_{k_1},\ldots,b_{k_q})$ $(1 \leq k \leq n)$.  
Now write
$H((A_1,B_1, \ldots, A_n,B_n),t) = (A_1(t),B_1(t), \ldots,A_n(t),B_n(t))$, where
\[
A_k(t) = 
\begin{cases}
\ (a_{k_1},\ldots,a_{k_q},2t(k-1)/n) \hspace{2.75cm} (0 \leq t \leq 1/2)\\
\ ((2-2t)a_{k_1}, \ldots,(2-2t)a_{k_q},(k-1)/n) \hspace{0.5cm}  (1/2 \leq t \leq 1)
\end{cases}
\]
and 
\[
B_k(t) = 
\begin{cases}
\ (b_{k_1},\ldots,b_{k_q},1 - 2t(1 - k/n)) \hspace{4.05cm}  (0 \leq t \leq 1/2)\\
\ (2t - 1 + (2-2t)b_{k_1}, \ldots,2t - 1 + (2-2t)b_{k_q},k/n) \hspace{0.5cm} (1/2 \leq t \leq 1).]
\end{cases}
\]

%%----------------------------------------------------------------------------------------------42
[Note: \  At the opposite extreme, each path component of $\xBV(R(1),n)$ is contractible and 
$\pi_0(\xBV(R(1),n)) \approx S_n$.]\\

\begin{proposition} \ %31
$\xBV^\infty$ is a universal creation operator.\\
\end{proposition}

\label{14.116}
Application: Let $(X,x_0)$ be a wellpointed compactly generated space with $\{x_0\} \subset X$ closed $-$then 
$\ohc \ \pow X$ and $\xBV^\infty[X]$ have the same homotopy type.\\

Let $\xBV^q$ be the creation operator that sends $n$ to $\xBV(R(q),n)$ $-$then 
$\xBV^\infty = \colimx  \xBV^q$ $\implies$ $\xBV^\infty[X] = \colimx  \xBV^q[X]$.\\

\begingroup%%----------------------------------->>
\fontsize{9pt}{11pt}\selectfont
\label{14.110}
\textbf{\small FACT} \ 
The arrow $\xBV^q[X] \ra \xBV^{q+1}[X]$ is a closed cofibration.\\
\endgroup %%------------------------------------<<

\begin{proposition} \ %32
The map $\xBV(R(q),n) \ra F(\R^q,n)$ which takes $(R_1, \ldots, R_n)$ to its center is an $S_n$-equivariant homotopy equivalence, hence induces a homotopy equivalence $\xBV^q[X] \ra \con \ \R^q[X]$.\\
\end{proposition}

The elements of $R(q)$ are in a one-to-one correspondence with the functions $[0,1]^q \ra [0,1]^q$ 
of the form 
$R = r_1 \times \cdots \times r_q$, 
where 
$r_i(t) = (b_i - a_i)t + a_i$ $(0 \leq a_i < b_i \leq 1)$.  
Thus $R(q)$ can be viewed as a subspace of $C([0,1]^q,[0,1]^q)$ (compact open topology), there being no ambiguity in so doing since the two interpretations are homeomorphic.

Representing $\bS^q$ as $[0,1]^q/\fr[0,1]^q$, adjust the definitions of $\Sigma^q X$ and $\Omega^q \Sigma^q X$ 
correspondingly $-$then the 
\un{arrow of May}
\index{arrow of May} 
is the continuous function 
$m_q:\xBV^q[X] \ra \Omega^q \Sigma^q X$
specified by the rule
\[
m_q[(R_1,\ldots,R_n),x_1,\ldots, x_n](s) = 
\begin{cases}
\ [x_i,t] \hspace{0.5cm} \text{if } R_i(t) = s\\
\ * \hspace{1.2cm}  \text{if } s \notin \bigcup\limits_i \im R_i
\end{cases}
\hspace{-.25cm}.
\]
\\

\index{Theorem: May's Approximation Theorem}
\index{May's Approximation Theorem}
\textbf{\small MAY'S APPROXIMATION THEOREM} \quad
Let $(X,x_0)$ be a wellpointed compactly generated space with $\{x_0\} \subset X$ closed.  Assume: $X$ is path connected $-$then 
$m_q:\xBV^q[X] \ra \Omega^q \Sigma^q X$ is a weak homotopy equivalence.

[Note: \  If $X$ has the pointed homotopy type of  a pointed connected CW complex, then $\xBV^q[X]$ is a pointed CW space, as is $\Omega^q \Sigma^q X$  (loop space theorem), thus under these circumstances the arrow of May is a pointed homotopy equivalence.]\\

The proof of this result is fairly lengthy and will be omitted.  In principle, the argument is an elaboration of that used in Proposition 20 and can be summarized in a sentence: There
%%----------------------------------------------------------------------------------------------43
is a commutative diagram
\[
\begin{tikzcd}%[sep=small]
{\xBV^q[X]} \ar{d}[swap]{m_q} \ar{r} 
&{E^q(\Gamma X,X)} \ar{d} \ar{r}  
&{\xBV^{q-1}[\Sigma X]} \ar{d}{m_{q-1}} \\
{\Omega \Omega^{q-1} \Sigma^q X} \ar{r} 
&{\Theta \Omega^{q-1} \Sigma^q X} \ar{r}
&{\Omega^{q-1} \Sigma^q X}
\end{tikzcd}
,
\]
where $E^q(\Gamma X,X) \ra \xBV^{q-1}[\Sigma X]$ is a quasifibration with fiber $\xBV^q[X]$ and 
$E^q(\Gamma X, X)$ is contractible, thus one may proceed by induction.  Details are in 
May\footnote[2]{\textit{SLN} \textbf{271} (1972), 50-68.}.

[Note: \  When $q = 1$, $\xBV^0[\Sigma X] = \Sigma X$ and $m_0$ is the identity map.]\\

Notation: Given a pointed \dsep compactly generated space $X$, let $\Oinf\Sinf X = \colimx \Omega^q \Sigma^q X$.

[Note: \  The reason for imposing the $\Delta$-separation condition is that it ensures the validity of the 
\un{repetition principle}: 
\index{repetition principle ($\Oinf \Sinf $)} 
$\Omega \Oinf \Sinf \Sigma X \approx \Oinf \Sinf X$.  
Proof:
$(\Oinf \Sinf \Sigma X)^{\bS^1} \approx$ 
$(\colimx \Omega^q \Sigma^q \Sigma X)^{\bS^1} \approx$ 
$\colimx  (\Omega^q \Sigma^q \Sigma X)^{\bS^1} \approx$ 
$\colimx \Omega^{q+1} \Sigma^{q+1} X \approx$  
$\Oinf \Sinf X$.]\\

\begingroup%%----------------------------------->>
\fontsize{9pt}{11pt}\selectfont
The arrow 
$\Omega^q \Sigma^q X \ra \Omega^{q+1} \Sigma^{q+1} X$ 
is the result of applying $\Omega^q$ to the arrow of adjunction $\Sigma^q X \ra \Omega \Sigma \Sigma^q X$.  
It is a closed embedding but it need not be a closed cofibration even if $X$ is wellpointed (in which case, of course, 
$\Omega^q \Sigma^q X$ is wellpointed $\forall \ q$).\\
\endgroup %%------------------------------------<<

\begingroup%%----------------------------------->>
\fontsize{9pt}{11pt}\selectfont
\textbf{\small EXAMPLE} \ 
Suppose that $X$ and $Y$ are pointed finite CW complexes $-$then 
$\Oinf \Sinf X$ and $\Oinf \Sinf Y$ are homotopy equivalent iff $\Sigma^q X$ and $\Sigma^q Y$ 
are homotopy equivalent for some $q \gg 0$ 
(Bruner-Cohen-McGibbon\footnote[2]{\textit{Quart. J. Math.} \textbf{46} (1995), 11-20.}).\\
\endgroup %%------------------------------------<<

Notation: Given a wellpointed \dsep compactly generated space $X$, put 
$m_\infty = \colimx m_q:\xBV^\infty [X] \ra \Oinf \Sinf X$.

[Note: \  $\xBV^\infty[X]$ is wellpointed (since $\forall \ q$, the arrow 
$\xBV^q[X] \ra \xBV^{q+1}[X]$ is a closed cofibration 
(cf. p. \pageref{14.110})) 
but it is problematic whether this is true of 
$\Oinf \Sinf X$ without additional assumptions on $X$.]\\

\begin{proposition} \ %33
Let $(X,x_0)$ be a wellpointed compactly generated space with $\{x_0\} \subset X$ closed.  
Assume: $X$ is \dsep and path connected $-$then 
$m_\infty:\xBV^\infty[X] \ra \Oinf\Sinf X$ is a weak homotopy equivalence.
\end{proposition}

%%----------------------------------------------------------------------------------------------44
[In the commutative ladder
\begin{tikzcd}%[sep=small]
{\xBV^1[X]} \ar{d} \ar{r} &{\xBV^2[X]}  \ar{d} \ar{r} &{\cdots}\\
{\Omega \Sigma X} \ar{r} &{\Omega^2 \Sigma^2 X} \ar{r} &{\cdots}
\end{tikzcd}
, the vertical arrows are weak homotopy equivalences and the spaces are $\tT_1$, so the generality on 
p. \pageref{14.111} 
can be quoted.]\\

A compactly generated space $X$ is said to be 
\un{$\Delta$-cofibered} 
\index{cofibered! $\Delta$-cofibered}
if the inclusion $\Delta_X \ra X \times_k X$ is a closed cofibration.

[Note: \  It is automatic that $\forall \ x_0 \in X$, $\{x_0\} \ra X$ is a closed cofibration 
(cf. p. \pageref{14.112}).]\\

\begingroup%%----------------------------------->>
\fontsize{9pt}{11pt}\selectfont
\textbf{\small FACT} \ 
Let $K$ be a pointed compact Hausdorff space.  Suppose that $X$ is pointed and $\Delta$-cofibered $-$then the pointed exponential object $X^K $ is $\Delta$-cofibered.\\
\endgroup %%------------------------------------<<

Example: Let $(X,x_0)$ be a pointed compactly generated space.  Assume: $X$ is $\Delta$-cofibered $-$then 
$\Sigma X$ is $\Delta$-cofibered 
(cf. p. \pageref{14.113}), 
as is $\Omega X$.\\

\textbf{\small LEMMA} \ 
Let $(X,x_0)$ be a pointed compactly generated space.  Assume: $X$ is $\Delta$-cofibered $-$then the arrow of adjunction 
$X \ra \Omega \Sigma X$ is a closed cofibration.\\

\label{14.163}
Application: Let $(X,x_0)$ be a pointed compactly generated space.  Assume: $X$ is $\Delta$-cofibered $-$then 
$\forall \ q$, the arrow 
$\Omega^q \Sigma^q X \ra \Omega^{q+1} \Sigma^{q+1}  X$ is a closed cofibration.

\label{14.139}
\label{14.149}
\label{16.65}
[Note: \  It is a corollary that $\Oinf \Sinf X$ is $\Delta$-cofibered 
(cf. p. \pageref{14.114}).]\\

\begingroup%%----------------------------------->>
\fontsize{9pt}{11pt}\selectfont
\textbf{\small LEMMA} \ 
Let $(X,x_0)$ be a pointed compactly generated space.  Assume: $X$ is $\Delta$-cofibered  
$-$then for every pointed $\Delta$-cofibered compact Hausdorff space $K \neq *$, the arrow 
$X \ra (X \#_k K)^K$ 
adjoint to the identity 
$X \#_k K \ra X \#_k K$ 
is a closed cofibration.
\vspi
[Note: \  Specialize and take $K = \bS^1$ to see that the arrow of adjunction $X \ra \Omega \Sigma X$ is a closed cofibration.\\
\endgroup %%------------------------------------<<

\begingroup%%----------------------------------->>
\fontsize{9pt}{11pt}\selectfont
\textbf{\small FACT} \ 
Let 
$
\begin{cases}
\ (X,x_0)\\
\ (Y,y_0)
\end{cases}
$
be pointed compactly generated spaces.  Assume: $X$ is $\Delta$-cofibered  and $Y$ is \dsep 
$-$then for every pointed $\Delta$-cofibered compact Hausdorff space $K \neq *$, 
the arrow $X \ra Y^K$ adjoint to a closed cofibration $X \#_k K \ra Y$ is a closed cofibration.
\vspi
[Factor the arrow $X \ra Y^K$ as the composite 
$X \ra (X \#_k K)^K \ra Y^K$.]\\
\endgroup %%------------------------------------<<

\begingroup%%----------------------------------->>
\fontsize{9pt}{11pt}\selectfont
\textbf{\small FACT} \ 
Suppose that $A \ra X$ is a closed cofibration, where $X$ is $\Delta$-cofibered $-$then $A$ is $\Delta$-cofibered (cf. $\S 3$, Proposition 11) and the arrow
$\Oinf\Sinf A \ra \Oinf\Sinf X$ is a closed cofibration.
\vspi
%%----------------------------------------------------------------------------------------------45
[All the arrows in the pullback square
\begin{tikzcd}[sep=large]
{\Omega^q\Sigma^q A} \ar{d} \ar{r} &{\Omega^{q+1}\Sigma^{q+1} A} \ar{d}\\
{\Omega^q\Sigma^q X} \ar{r}  &{\Omega^{q+1}\Sigma^{q+1} X}
\end{tikzcd}
are closed cofibrations, so one can appeal to the lemma on 
p. \pageref{14.115}.]\\
\endgroup %%------------------------------------<<

\begin{proposition} \ %34
Let $(X,x_0)$ be a pointed compactly generated space.  Assume: $X$ is \dcf and has the pointed homotopy type of a pointed connected CW complex $-$then
$m_\infty: \xBV^\infty[X] \ra \Oinf\Sinf X$ is a pointed homotopy equivalence.
\end{proposition}

[In the commutative ladder 
$
\begin{tikzcd}%[sep=small]
{\xBV^1[X]} \ar{d} \ar{r} &{\xBV^2[X]} \ar{d} \ar{r}  &{\cdots}\\
{\Omega\Sigma X}  \ar{r} &{\Omega^2\Sigma^2 X}  \ar{r}  &{\cdots}
\end{tikzcd}
, 
$
the horizontal arrows are closed cofibrations and the vertical arrows are pointed homotopy equivalenes.  Now cite Proposition 15 in $\S 3$.]\\

\index{Theorem: Homotopy Colimit Theorem}
\index{Homotopy Colimit Theorem}
\textbf{\small HOMOTOPY COLIMIT THEOREM} \quad
Let $(X,x_0)$ be a pointed connected CW complex or a pointed connected ANR $-$then 
$\ohc \ \pow X$ and $\Oinf\Sinf X$ have the same homotopy type.

[One has only to recall that $\ohc \ \pow X$ and $\xBV^\infty[X]$ have the same homotopy type 
(cf. p. \pageref{14.116}).]

[Note: \  For the validity of the condition on the diagonal, 
cf. p. \pageref{14.117} $\&$ 
p. \pageref{14118.}.]\\

\begingroup%%----------------------------------->>
\fontsize{9pt}{11pt}\selectfont
\textbf{\small EXAMPLE} \ 
Connectedness is essential here.  For example, the homotopy type of $\ohc \  \pow \ \bS^0$ is represented by $B\bM_\infty$ 
(cf. p. \pageref{14.119}) 
but the homotopy type of $\Oinf\Sinf \bS^0$ is represented by $\Omega B \abs{M_\infty}$ 
(cf. p. \pageref{14.120}) 
$(\abs{M_\infty} = B\bM_\infty = \ds\coprod\limits_{n \geq 0} B S^n)$.\\
\label{14.165}
\endgroup %%------------------------------------<<

Given a cofunctor $\sC:\iso\bGamma \ra \bCG$, let 
$\widehat{\sC}(\bm,\bn) = \coprod\limits_{\gamma:\bm\ra\bn}\prod\limits_{1 \leq j \leq n} \sC(\#(\gamma^{-1}(j)))$ 
(here $\gamma$ ranges over the morphisms 
$\bm \ra \bn$ in $\bGamma$ and $\widehat{\sC}(\bm,\bzero)$ 
is a point indexed by the unique arrow $\bm \ra \bzero$) $-$then with the obvious choice for the unit, 
$[(\iso\bGamma)^\OP,\bCG]$ aquires the structure of a monoidal category by writing 
$\sC \circx \sD(\bm) = \coprod\limits_{\bn \geq \bzero} \sC(\bn) \times_{S_n} \widehat{\sD}(\bm,\bn)$.\\

\textbf{\small LEMMA} \ 
The functor $-\circx \sD$ has a right adjoint $\shom(\sD,-)$, where 
$\shom(\sD,\sE)(\bn) = \prod\limits_{\bm \geq \bzero} \hom(\widehat{\sD}(\bm,\bn),\sE(\bm))^{S_m}$ 
$(\hom = kC_k$, the internal hom functor in \bCG 
(cf. p. \pageref{14.120a})), 
so 
$\Nat(\sC \circx \sD,\sE) \approx \Nat(\sC,\shom(\sD,\sE))$.\\

An 
\un{operad} 
\index{operad} 
$\sO$ in \bCG is a monoid in the monoidal category $[(\iso\bGamma)^\OP,\bCG]$.  
Examples: 
(1) Let $\sO_n = *  \ \forall \ n$;
(2) Let $\sO_n = S_n  \ \forall \ n$.\\

%%----------------------------------------------------------------------------------------------46

\begingroup%%----------------------------------->>
\fontsize{9pt}{11pt}\selectfont
The definition of an operad makes sense if \bCG is replaced by any symmetric monoidal category \bC which is complete and cocomplete.
\vspi
[Note: \  Agreeing to write \bOPERc for \bMONx, one can show that \bOPERc is complete and cocomplete and that the forgetful functor $\bOPERc \ra \bisox$ has a left adjoint, the free operad functor 
(Getzler-Jones\footnote[2]{\textit{Operads, Homotopy Algebra, and Iterated Integrals for Double Loop Spaces}, Preprint.}).]\\
\endgroup %%------------------------------------<<

Equivalently, an operad $\sO$ in \bCG consists of compactly generated spaces $\sO_n$ equipped with a right action of $S_n$, a point $1 \in \sO_1$ (the \un{unit}) and for each sequence $j_1, \ldots, j_n$ of nonnegative integers, a continuous function 
$\Lambda:\sO_n \times_k (\sO_{j_1} \times_k  \cdots \times_k  \sO_{j_n}) \ra \sO_{j_1 + \cdots + j_n}$ 
satisifying the following conditions.  

\indent\indent (OPER$_1$) \quad Given $\sigma \in S_n$, $\sigma_k \in S_{j_k}$ $(k = 1, \ldots, n)$, and $f \in \sO_n$, 
$g_k \in \sO_{j_k}$, one has 
$\Lambda(f \cdot \sigma;g_1, \ldots, g_n) = \Lambda(f;g_{\sigma^{-1}(1)}, \ldots, g_{\sigma^{-1}(n)}) 
\cdot \sigma(j_1, \ldots, j_n)$ 
($\sigma(j_1, \ldots, j_n)$ the permutation of $S_{j_1 + \cdots + j_n}$ 
that permutes the $n$ blocks of $j_k$ successive integers per $\sigma$, the order within each block staying fixed) and 
$\Lambda(f;g_1\cdot\sigma_1, \ldots, g_n\cdot\sigma_n) = \Lambda(f;g_1, \ldots, g_n) \cdot 
(\sigma_1 \amalg \cdots \amalg \sigma_n)$ ($\sigma_1 \amalg \cdots \amalg \sigma_n$  
the permutation of 
$S_{j_1 + \cdots + j_n}$
 that leaves the $n$ blocks invariant and which restricts to $\sigma_k$ on the $k^{th}$ block).\\
\indent\indent (OPER$_2$) \quad Given $f  \in \sO_n$, $g_k \in \sO_{j_k}$ 
$(k = 1, \ldots, n)$, 
$h_{kl} \in \sO_{i_{kl}}$ 
$(l = 1, \ldots, j_k)$, one has 
$\Lambda(f;\Lambda(g_k;h_{kl})) = \Lambda(\Lambda(f;g_k);h_{kl})$.\\
\indent\indent (OPER$_3$) \quad  Given $f  \in \sO_n$, one has $\Lambda(f;1, \ldots, 1) = f$ and given $g \in \sO_j$, one has 
$\Lambda(1;g) = g$.

Example: \ $\xBV^q$ is an operad in \bCG.  Thus with $\sO_n = \xBV(R(q),n)$, write $(R_1, \ldots, R_n)$.  
$\sigma = (R_{\sigma(1)}, \ldots, R_{\sigma(n)})$ $(\sigma \in S_n)$, take for $1 \in \xBV(R(q),1)$ the identity function, and let 
$\Lambda:\xBV(R(q),n) \times_k (\xBV(R(q),j_1) \times_k \cdots \times_k \xBV(R(q),j_k)) \ra$ $\xBV(R(q),j_1 + \cdots + j_n)$ 
be defined on elements via composition 
$j_1 \cdot [0,1]^q \coprod \cdots \coprod j_n \cdot [0,1]^q \ra n \cdot [0,1]^q \ra$ $[0,1]$.

[Note: \  con $\R^q$ is not an operad in \bCG.]\\

\label{14.133}
\begingroup%%----------------------------------->>
\fontsize{9pt}{11pt}\selectfont
\textbf{\small EXAMPLE} \ 
Let $\sO$ be an operad in \bCG such that $\forall \ n$, $\sO_n \neq \emptyset$.  
Definition: $\grd\sO$ is the operad in \bCG with $\grd_n\sO = \abs{\ner \ \grd\sO_n}$ 
(cf. p. \pageref{14.121}).  
To specify the right action of $S_n$, note that there is a simplicial map $\si S_n \ra \ner\grd S_n$, hence 
$\abs{\ner \ \grd\sO_n} \times S_n \ra  \abs{\ner\grd\sO_n} \times  \abs{\ner \ \grd S_n} \approx$ 
$\abs{\ner \ \grd(\sO_n \times S_n)} \ra$ $ \abs{\ner\grd\sO_n}$.  
Next, 
$\sO_1 =  \abs{\ner \ \grd\sO_1}_0$, so the choice for 1 is clear.
Finally, $\Lambda$ is defined by 
$\abs{\ner \ \grd\sO_n} \times_k \left(\abs{\ner \ \grd\sO_{j_1}} \times_k \cdots \times_k \abs{\ner \ \grd\sO_{j_n}}\right)$ 
$\approx$ $\abs{\ner \ \grd\sO_n \times_k \left(\sO_{j_1} \times_k \cdots \times_k \sO_{j_n}\right))}$ $\ra$ 
$\abs{\ner \ \grd(\sO_{j_1 + \cdots + j_n})}$.  
Example: let $\sO_n = S_n$ $\forall \ n$ $-$then $\grd_n\sO \approx \abs{\ner \ \tran S_n}$ 
(cf. p. \pageref{14.122} ff.), 
i.e., 
$\grd_n\sO \approx X S_n$.\\
\endgroup %%------------------------------------<<

%%----------------------------------------------------------------------------------------------47
In terms of the $\Lambda$, a morphism $\sO \ra \sP$ of operads in \bCG is a sequence of $S_n$-equivariant continuous functions $\sO_n \ra \sP_n$ such that the diagrams
%\begin{tikzcd}%[sep=small]
%{\sO_n \times_k (\sO_{j_1} \times_k \cdots \times_k \sO_{j_n}} \ar{d} \ar{r} 
%&{\sO_{j_1 + \cdots + j_n}}\ar{d}\\
%{\sP_n \times_k (\sP_{j_1} \times_k \cdots \times_k \sP_{j_n}} \ar{r} 
%&{\sP_{j_1 + \cdots + j_n}}
%\end{tikzcd}
\begin{tikzcd}%[sep=small]
{\sO_n \times_k (\sO_{j_1} \times_k \cdots \times_k \sO_{j_n})} \ar{d} \ar{r} 
&{}\\
{\sP_n \times_k (\sP_{j_1} \times_k \cdots \times_k \sP_{j_n})} \ar{r} 
&{}
\end{tikzcd}
\begin{tikzcd}%[sep=small]
{\sO_{j_1 + \cdots + j_n}}\ar{d}\\
{\sP_{j_1 + \cdots + j_n}}
\end{tikzcd}
commute and $\sO_1 \ra \sP_1$ sends 1 to 1.

Example: $\forall \ q$, the arrow $\xBV^q \ra \xBV^{q+1}$ is a morphism of operads in \bCG.\\

\label{14.137}
\begingroup%%----------------------------------->>
\fontsize{9pt}{11pt}\selectfont
\textbf{\small EXAMPLE} \ 
If $\sO$ is an operad in \bCG, then $\sin \sO$ is an operad in \bSISET.  Its geometric realization $\abs{\sin \sO}$ is an operad in \bCG and the arrow $\abs{\sin \sO} \ra \sO$ is a morphism of operads in \bCG.\\
\endgroup %%------------------------------------<<

An operad $\sO$ in \bCG is said to be 
\un{reduced} 
\index{reduced operad} if $\sO_0 = *$.\\

\label{14.172}

\begin{proposition} \ %35
Let $\sO$ be a reduced operad in \bCG $-$then $\sO$ extends to a creation operator in $\bGamma_\ini^\OP \ra \bCG$.
\end{proposition}

[It suffices to define $\sO$ on the order preserving injections 
(cf. p. \pageref{14.122a}) %13-56
or still, for each $n$, on the 
$n+1$ elementary order preserving injections 
$\sigma_i:\bn \ra \bn + \bone$, where
$
\begin{cases}
\ j \ra j \qquad (j \leq i)\\
\ j \ra j + 1 \ (j > i)
\end{cases}
(0 \leq i \leq n),
$
the requisite arrows $\sO_{n+1} \ra \sO_n$ thus being the assignments 
$f \ra \Lambda(f; 1^i, *, 1^{n-i}))$.]\\

Notation: $\bCG_{*\bc}$ is the full subcategory of $\bCG_*$ whose objects are the $(X,x_0)$ such that 
$* \ra (X,x_0)$ is a closed cofibration.

[Note: \  The standard model category structure on $\bCG_*$ is that inherited from the standard model category structure on \bCG 
(cf. p. \pageref{14.123}) %12-3
and the cofibrant objects therein are the objects of $\bCG_{*\bc}$.]

Observation: For any creation operator $\sC$, $\sC[?]$ is a functor 
$\bCG_{*\bc} \ra \bCG_{*\bc}$ (cf. Proposition 27).\\

\begin{proposition} \ %36
Let $\sO$ be a reduced operad in \bCG $-$then $\sO$ determines a triple 
$\bT_{\sO} = (T_{\sO}, m, \epsilon)$ in $\bCG_{*\bc}$.
\end{proposition}

[Take 
$T_{\sO} = \sO[?]$ and for each \mX, define 
$m_X:\sO^2[X] \ra \sO[X]$, 
$\epsilon_X:X \ra \sO[X]$ 
by the formulas 
$m_X[f,[g_1,x_1], \ldots, [g_n,x_n]] = [\Lambda(f; g_1, \ldots, g_n),x_1, \ldots, x_n]$ 
$(f \in \sO_{n}, g_k \in \sO_{j_k}$ $\&$ 
$x_k \in X^{j_k}$ $(1 \leq k \leq n)$), 
$\epsilon_X(x) = [1,x]$ $(x \in X)$.]

[Note: \  A morphism $\sO \ra \sP$ of reduced operads in \bCG leads to a morphism 
$\bT_{\sO} \ra \bT_{\sP}$ of triples in $\bCG_{*\bc}$.]\\

Examples: 
(1) With $\sO_n = *$ $\forall \ n$, $T_{\sO}X = SP^\infty X$; 
(2) With $\sO_n = S_n$ $\forall \ n$, $T_{\sO}X = JX$.\\

%%----------------------------------------------------------------------------------------------48
\label{14.125}
\begingroup%%----------------------------------->>
\fontsize{9pt}{11pt}\selectfont
\textbf{\small FACT} \ 
Let $X$ be a pointed compactly generated simplicial space satisfying the cofibration condition such that $\forall \ n$, $X_n$ 
is in $\bCG_{*\bc}$.  Given a reduced operad $\sO$ in \bCG, define a pointed compactly generated simplicial space 
$\sO[X]$ by  $\sO[X]_n = \sO[X_n]$ $-$then $\abs{\sO[X]} \approx \sO[\aX]$.
\vspi
[Work with the arrow 
$[[f,x_1, \ldots,x_k],t] \ra [f,[x_1,t],\ldots,[x_k,t]]$, where $f \in \sO_k$, $x_j \in X_n$ $(1 \leq j \leq k)$, $t \in \dpn$.]
\vspi
[Note: \  The diagrams
\begin{tikzcd}%[sep=small]
{\abs{\sO^2[X]}} \ar{d}[swap]{\abs{m_X}} \ar{r} &{\sO^2[\aX]} \ar{d}{m_{\aX}} \\
{\abs{\sO[X]}}  \ar{r} &{\sO[\aX]}  
\end{tikzcd}
,
\begin{tikzcd}%[sep=small]
&{\abs{\sO[X]}} \ar{dd}\\
{\aX} \ar{ru}{\abs{\epsilon_X}} \ar{rd}[swap]{\epsilon_{\aX}} \\
&{\sO[\aX]}
\end{tikzcd}
commute.  Consequently, if $X$ is a simplicial $T_{\sO}$-algebra, then $\aX$ is a $T_{\sO}$-algebra (by the composite 
$\sO[\aX] \ra \abs{\sO[x]} \ra \aX$).]\\
\label{14.148}
\label{14.154}
\label{14.155}
\endgroup %%------------------------------------<<

Let $\sO$ be a reduced operad in \bCG $-$then an 
\un{$\sO$-space} 
\index{space! $\sO$-space} 
is an object $(X,x_0)$ in $\bCG_{*\bc}$ 
and continuous functions $\theta_n:\sO_n \times_k X^n \ra X$ $(n \geq 0)$ subject to the following assumptions.\\
\indent\indent ($\sO$-SP$_1$) \quad  Given $\sigma \in S_n$, $f \in \sO_n$, and $x_k \in X$ $(k = 1, \ldots, n)$, 
one has $\theta_n(f\cdot\sigma,x_1, \ldots,x_n) = \theta_n(f,x_{\sigma^{-1}(1)}, \ldots, x_{\sigma^{-1}(n)})$. \\
\indent\indent ($\sO$-SP$_2$) \quad  Given $f \in \sO_n$, $g_k \in \sO_{j_k}$ $(k = 1, \ldots, n)$, 
$x_{kl} \in X$ $(l = 1, \ldots, j_k)$, one has
$
\theta_{j_1 + \cdots + j_n}(\Lambda(f;g_1, \ldots, g_n), x_{11}, \ldots, x_{1j_1}, \ldots, x_{n1}, \ldots, x_{nj_n}) = $
$\theta_n(f,\theta_{j_1}(g_1;x_{11},\ldots,x_{1j_1}),\ldots,$\\ $\theta_{j_n}(g_n;x_{n1}, \ldots, x_{nj_n}))$
.\\
\indent\indent ($\sO$-SP$_3$) \quad  $\theta_0(*) = x_0$ and $\theta_1(1,x) = x \ \forall \ x \in X$. 

\label{14.156}
[Note: \  In practice, one sometimes encounters objects in $\bCG_*$ satisfying all the assumptions that define an $\sO$-space but, strictly speaking, are not $\sO$-spaces because they may not be in $\bCG_{*\bc}$.  
Up to homotopy equivalence, this is not a problem.  
Thus let $X$ be an $\sO$-space in $\bCG_*$ and consider $\overset{\vee}{X}$ 
(cf. p. \pageref{14.123a}).
Define $\overset{\vee}{\theta}_n:\sO_n \times_k \overset{\vee}{X}^n \ra \overset{\vee}{X}$ by 
$\overset{\vee}{\theta}_n(f,\overset{\vee}{x}_1,\ldots,\overset{\vee}{x}_n) = $
$
\begin{cases}
\ {\theta}_n(f,r(\overset{\vee}{x}_1),\ldots,r(\overset{\vee}{x}_n)) \hspace{0.5cm}  \text{if } \overset{\vee}{x}_i \notin [0,1] - \{0\} \ (\exists \ i)\\
\ \overset{\vee}{x}_1 \cdots \overset{\vee}{x}_n  \hspace{2.82cm}   \text{if } \overset{\vee}{x}_i \in [0,1] (\forall \ i)
\end{cases}
$
$-$then $\overset{\vee}{X}$ is an $\sO$-space in $\bCG_{*\bc}$ and the retraction $r:\overset{\vee}{X} \ra X$ is a morphism of $\sO$-spaces.]

Examples: 
(1) If $\sO_n = *$ $\forall \ n$, then the $\sO$-spaces are the abelian cofibered monoids in \bCG; 
(2) $\sO_n = S_n$ $\forall \ n$, then the $\sO$-spaces are the cofibered monoids in \bCG.

Example: $\forall \ X$ in $\bCG_{*\bc}$, $\Omega^qX$ is a $\xBV^q$-space.

[Define $\theta_n:\xBV(R(q),n) \times_k (\Omega^qX)^n \ra \Omega^qX$ by sending 
$((R_1, \ldots, R_n),f_1,\ldots, f_n)$ to that element of $\Omega^qX$ which at $s$ is $f_i(t)$ 
if $R_i(t) = s$ lies in the interior of $R_i$ and is $x_0$ otherwise.]\\

\begingroup%%----------------------------------->>
\fontsize{9pt}{11pt}\selectfont
\textbf{\small EXAMPLE} \ 
Let $\sS$ be the operad in \bCAT with $S_n = \tran S_n$ $\forall \ n$ $-$then in suggestive terminology, a permutative category \bC is an $\sS$-category 
\index{category! $\sS$-category}, 
thus its classifying space $B\bC$ is a B$\sS$-space.
\vspi
[Note: \  $BS_n = B \tran S_n = \abs{\ner \ \tran S_n} = \abs{\barr(*,\bS_n; S_n)} = XS_n$.]\\
\endgroup %%------------------------------------<<

$\sO$-\bSP is the category whose objects are the $\sO$-spaces and whose morphisms $X \ra Y$
%%----------------------------------------------------------------------------------------------49
are the pointed continuous functions $X \ra Y$ such that the diagrams 
\begin{tikzcd}%[sep=small]
{\sO_n \times_k X^n} \ar{d}\ar{r} &{\sO_n \times_k Y^n} \ar{d} \\
{X}  \ar{r} &{Y}  
\end{tikzcd}
commute.

Example: $\sO$-\bSP $= \bCG_{*\bc}$, if $\sO_0 = *$, $\sO_1 = \{1\}$, $\sO_n = \emptyset$ $(n > 1)$.\\

\begingroup%%----------------------------------->>
\fontsize{9pt}{11pt}\selectfont
\textbf{\small EXAMPLE} \ 
If $X$ is an $\sO$-space, then so are $\Omega X$ and $\Theta X$.  Moreover, the inclusion $\Omega X \ra \Theta X$ is a morphism of $\sO$-spaces, as is the \bCG fibration $\Theta X \ra X$.\\
\endgroup %%------------------------------------<<

\label{14.167}
\begin{proposition} \ %37
Let $\sO$ be a reduced operad in \bCG $-$then the categories $\sO$-$\bSP$ and 
$\bT_{\sO}$-$\bALG$ are canonically isomorphic.
\end{proposition}

[There is a one-to-one correspondence between the $\sO$-space structures on $X$ and the $\bT_{\sO}$-algebra structures on $X$, encapsulated in the commutativity of the diagrams
$
\begin{tikzcd}%[sep=small]
{\sO_n \times_k X^n} \ar{rrd}[swap]{\theta_n} \ar{r} &{\sO_n[X]} \ar{r} &{\sO[X]}\ar{d}{\theta}\\
&&{X}  
\end{tikzcd}
$
for all $n$, i.e., the $\theta_n$ combine to define an arrow $\theta:\sO[X] \ra X$ satisfying TA$_1$ and TA$_2$ 
(cf. p. \pageref{14.124} ff.) 
and vice versa).]

[Note:  \ The 
\un{endomorphism operad}
\index{endomorphism operad}
\index{operad: endomorphism operad} 
$\End X$ of $X$ is defined by $(\End X)_n = X^{X^n}$ (pointed exponential object in $\bCG_*$), supplied with the evident operations.  Taking adjoints, the $\bT_{\sO}$-algebra structures on $X$ correspond bijectively to morphisms of operads 
$\sO \ra \End X$ in $\bCG$.]
\\

Example: $\forall \ X$, $\sO[X]$ is a $\bT_{\sO}$-algebra, hence is an $\sO$-space.
\\

\label{14.134} %dmc mnft
\label{14.160}

\begingroup%%----------------------------------->>
\fontsize{9pt}{11pt}\selectfont
\textbf{\small EXAMPLE} \ 
The functors $\Sigma^q:\bCG_* \ra \bCG_*$, $\Omega^q:\bCG_* \ra \bCG_*$ both respect $\bCG_{*\bc}$ and 
$(\Sigma^q,\Omega^q)$ is an adjoint pair, thus $\forall \ X$, there is an arrow of adjunction 
$X \ra \Omega^q\Sigma^q X$.  As noted above, $\Omega^q\Sigma^q X$ is a $\xBV^q$-space or still, a $\bT_{\xBV^q}$-algebra.  The composite 
$\xBV^q[X] \ra \xBV^q[\Omega^q\Sigma^q X] \ra \Omega^q\Sigma^q X$ is $m_q$, the arrow of May.  
It is a morphism of $\bT_{\xBV^q}$-algebras.  On the other hand, $\forall \ X$, there is an arrow of adjunction 
$\Sigma^q\Omega^q X \ra X$, from which $\Omega^q\Sigma^q\Omega^q X \ra \Omega^qX$.  
Viewing the $\xBV^q$-space $\Omega^q X$ as a $\bT_{\xBV^q}$-algebra, its structural morphism 
$\xBV^q[\Omega^q X] \ra \Omega^q X$ is the composite 
$\xBV^q[\Omega^q X] \overset{m_q}{\lra} \Omega^q\Sigma^q\Omega^q X \ra \Omega^qX$.\\
\endgroup %%------------------------------------<<

\begingroup%%----------------------------------->>
\fontsize{9pt}{11pt}\selectfont
\textbf{\small FACT} \ 
Let $X$ be a pointed compactly generated simplicial space satisfying the cofibration condition such that $\forall \ n$, $X_n$ is in 
$\bCG_{*\bc}$ $-$then the arrow $\abs{\Omega^qX} \ra \Omega^q \aX$ is a morphism of $\bT_{\xBV^q}$-algebras.
\vspi
[The structural morphism 
$\xBV^q[\abs{\Omega^q X}] \ra \abs{\Omega^q X}$ is the composite 
$\xBV^q[\abs{\Omega^q X}]  \ra \abs{\xBV^q[\Omega^qX]} \ra
%%----------------------------------------------------------------------------------------------50
\abs{\Omega^qX}$ 
(cf. p. \pageref{14.125}), 
thus one has to check that the diagram
\[
\begin{tikzcd}[sep=large]
{\xBV^q[\abs{\Omega^q X}]}  \ar{d} \ar{rr} &&{\xBV^q[\Omega^q \aX]} \ar{d}\\
{\abs{\xBV^q[\Omega^q X]}} \ar{r} &{\abs{\Omega^q X}} \ar{r} &{\Omega^q \aX}
\end{tikzcd}
\]
commutes.]\\
\endgroup %%------------------------------------<<

\begingroup%%----------------------------------->>
\fontsize{9pt}{11pt}\selectfont
\textbf{\small FACT} \ 
Let $X$ be a pointed compactly generated simplicial space satisfying the cofibration condition such that 
$\forall \ n$, $X_n$ is in 
$\bCG_{*\bc}$ $-$then the diagram
\begin{tikzcd}[sep=large]
{\abs{\xBV^q[X]}}  \ar{d}[swap]{\abs{m_q}} \ar{r} &{\xBV^q[\aX]} \ar{d}{m_q}\\
{\abs{\Omega^q \Sigma^q X}} \ar{r} &{\Omega^q \Sigma^q \aX}
\end{tikzcd}
commutes.\\
\vspace{0.25cm}
\endgroup %%------------------------------------<<

Let $\sO$ be a reduced operad in \bCG, $F:\bCG_{*\bc} \ra \bCG_{*\bc}$ a right $\bT_{\sO}$-functor $-$then for any 
$\bT_{\sO}$-algebra $X$, $\barr(F;\bT_{\sO};X)$ is a simplicial object in $\bCG_{*\bc}$ 
(cf. p. \pageref{14.126}) 
and one writes 
$B(F;\sO;X)$ for its geometric realization (or just $B(\sO;\sO;X)$ if $F = \bT_{\sO} = \sO[?]$).\\

\begin{proposition} \  %38
Let $\sO$ be a reduced operad in \bCG such that $\{1\} \ra \sO_1$ is a closed cofibration.  
Suppose that 
$F:\bCG_{*\bc} \ra \bCG_{*\bc}$ is a right $\bT_{\sO}$-functor which preserves closed cofibrations $-$then $\forall$ 
$\sO$-space $X$, $\barr(F;\bT_{\sO};X)$ satisfies the cofibration condition, hence $B(F;\sO;X)$ is in $\bCG_{*\bc}$.
\end{proposition}

[On general grounds, $\sO[?]$ preserves closed cofibrations 
(cf. p. \pageref{14.127}).  
Moreover, the assumption on the unit of $\sO$ implies that $\epsilon_X:X \ra \sO[X]$ is a closed cofibration $\forall \ X$, so the conclusion follows from the definition of the $s_i$ and the fact that $F$ preserves closed cofibrations.]\\

\begingroup%%----------------------------------->>
\fontsize{9pt}{11pt}\selectfont
\textbf{\small EXAMPLE} \ 
$\Sigma$ is a right $\bT_{\xBV^1}$-functor and preserves closed cofibrations.  If $G$ is a cofibered monoid in \bCG, then $G$ 
acquires the structure of a $\bT_{\xBV^1}$-algebra via the composite 
$\xBV^1[G] \ra JG \ra G$.  
Thus it is meaningful to form 
$\barr(\Sigma;\bT_{\xBV^1};G)$.  
Since 
$\{1\} \ra \xBV(R(1),1)$ is a closed cofibration, $\barr(\Sigma;\bT_{\xBV^1};G)$ 
satisfies the cofibration condition (cf. Proposition 38) and its geometric realization 
$B(\Sigma;\xBV^1;G)$ 
is the 
\un{classifying space} 
\index{classifying space (in the sense of May)}
of \mG in the sense of May.  It is true but not obvious that $B(\Sigma;\xBV^1;G)$ and $BG$ have the same weak homotopy type 
(Thomason\footnote[2]{\textit{Duke Math. J.} \textbf{46} (1979), 217-252; 
see also Fiedorowicz, \textit{Amer. J. Math.} \textbf{106} (1984), 301-350.}).\\
\endgroup %%------------------------------------<<

\begingroup%%----------------------------------->>
\fontsize{9pt}{11pt}\selectfont
\textbf{\small EXAMPLE} \ 
Suppose that $X$ is a path connected $\xBV^q$-space $-$then $X$ has the weak homotopy type of a $q$-fold loop space.  
In fact, $\Sigma^q$ is a right $\bT_{\xBV^q}$-functor, as is $\Omega^q\Sigma^q$, so one can form 
$B(\Sigma^q;\xBV^q;X)$ and $B(\Omega^q\Sigma^q;\xBV^q;X)$, where now $X$ is viewed as a $\bT_{\xBV^q}$-algebra.  
Consider the following diagram
%%----------------------------------------------------------------------------------------------51
in the category of $\bT_{\xBV^q}$-algebras:  
$X \la B(\xBV^q;\xBV^q;X) \ra B(\Omega^q\Sigma^q;\xBV^q;X) \ra \Omega^qB(\Sigma^q;\xBV^q;X)$.
Owing to the generalities on 
p.  \pageref{14.127a} ff., 
the arrow 
$X \la B(\xBV^q;\xBV^q;X)$ is a homotopy equivalence
 (cf. p. \pageref{14.128}).  
 Next, according to May's approximation theorem, 
$\forall \ n,$ $m_q:\xBV^q[(\xBV^q)^n[X]] \ra \Omega^q\Sigma^q(\xBV^q)^n[X]$ is a weak homotopy equivalence.  
Therefore, on account of Proposition 38, the arrow 
$B(\xBV^q;\xBV^q;X) \ra B(\Omega^q\Sigma^q;\xBV^q;X)$ is a weak homotopy equivalence
 (cf. p. \pageref{14.129}).  
As for the arrow 
$B(\Omega^q\Sigma^q;\xBV^q;X) \ra \Omega^qB(\Sigma^q;\xBV^q;X)$
it too is a weak homotopy equivalence.  Indeed, all data is path connected and 
$\barr(\Omega^q\Sigma^q;\bT_{\xBV^q};X) = \Omega^q\barr(\Sigma^q;\bT_{\xBV^q};X)$, 
thus
%note typon in online notes next line ????????????????????????
$\abs{\Omega^q\barr(\Sigma^q;\bT_{\xBV^q};X) }$ $\ra$ 
$\Omega^q\abs{\barr(\Sigma^q;\bT_{\xBV^q};X)}$ 
is a weak homotopy equivalence 
(cf. p. \pageref{14.130}).
\vspi
[Note: \  The composite
$X \ra B(\xBV^q;\xBV^q;X) \ra B(\Omega^q\Sigma^q;\xBV^q;X) \ra \Omega^qB(\Sigma^q;\xBV^q;X)$
is the adjoint of 
$\Sigma^q X \ra B(\Sigma^q;\xBV^q;X)$ 
but it is not a morphism of $\bT_{\xBV^q}$-algebras and one cannot expect to always find a morphism 
$X \ra \Omega^qY$ of $\bT_{\xBV^q}$-algebras which is a weak homotopy equivalence.  
Take, e.g., $q = 1$ and let $X$ be a path connected cofibered monoid in \bCG (thought of as a $\bT_{\xBV^1}$-algebra).  
Claim: The only morphism $X \ra \Omega Y$ of $\bT_{\xBV^1}$-algebras is the constant map $X \ra j(y_0)$.  
Proof: Inspect the commutative diagram
\begin{tikzcd}%[sep=small]
{\xBV(R(1),1) \times_k X} \ar{d} \ar{r} &{\xBV(R(1),1) \times_k \Omega Y} \ar{d}\\
{X} \ar{r} &{\Omega Y}
\end{tikzcd}
.]\\
\endgroup %%------------------------------------<<
\vspace{0.25cm}

\begingroup%%----------------------------------->>
\fontsize{9pt}{11pt}\selectfont
\textbf{\small EXAMPLE} \ 
Let $\sO$ be a reduced operad in \bCG such that $\{1\} \ra \sO_1$ is a closed cofibration.  Assume: $\forall \ n$, 
$\sO_n \ra *$ is an $S_n$-equivariant homotopy equivalence $-$then every $\sO$-space $X$ has the homotopy type of an abelian cofibered monoid in \bCG.  Indeed, $X$ and $B(\sO;\sO;X)$ have the same homotopy type.  
Moreover, $\forall \ n$, the arrow 
$\sO[\sO^n[X]] \ra SP^\infty\sO^n[X]$ is a homotopy equivalence (cf. Proposition 28), so the arrow 
$B(\sO;\sO;X) \ra B(SP^\infty;\sO;X)$
is a homotopy equivalence (cf. Proposition 4 and Proposition 38).  But $B(SP^\infty;\sO;X)$ is an abelian cofibered monoid in \bCG.\\
\endgroup %%------------------------------------<<

Let $\sO$ be a reduced operad in \bCG  $-$then $\sO$ is said to be an \un{$\tE_\infty$ operad} 
\index{E$_\infty$ operad} if $\forall \ n$, $\sO_n$ is a contractible compactly generated Hausdorff space, the action of $S_n$ is free, and the inclusion 
$\{1\} \ra \sO_1$ is a closed cofibration.

Example: $\xBV^\infty = \colimx \xBV^q$ is an $\tE_\infty$ operad, the 
\un{Boardman-Vogt operad}.
\index{Boardman-Vogt operad} \index{operad: Boardman-Vogt}

[In view of Proposition 31, the only thing that has to be checked is the cofibration condition on the unit.  
However, by definition, $\xBV(R(\infty),1) = \colimx \xBV(R(q),1)$ and $\xBV(R(q),1) \ra \xBV(R(q+1),1)$ is a closed cofibration.  
In addition, the diagonal embedding $\xBV(R(q),1) \ra \xBV(R(q),1) \times_k \xBV(R(q),1)$ is a closed cofibration ($\xBV(R(q),1)$ is a polyhedron), thus the diagonal embedding 
$\xBV(R(\infty),1) \ra \xBV(R(\infty),1) \times_k \xBV(R(\infty),1)$ is a closed cofibration 
(cf. p. \pageref{14.131}).  
Therefore the inclusion 
$\{1\} \ra \xBV(R(\infty),1)$ is a closed cofibration 
(cf. p. \pageref{14.132}).] \\

\begingroup%%----------------------------------->>
\fontsize{9pt}{11pt}\selectfont
\textbf{\small EXAMPLE} \ 
Let $\sO_n = S_n$ $\forall \ n$ $-$then $\grd\sO$ is an $\tE_\infty$ operad 
(cf. p. \pageref{14.133}), 
the 
\un{permutation operad} 
\index{permutation operad} 
\index{operad: permutation operad} 
PER.
\vspi
%%----------------------------------------------------------------------------------------------52
[Note: \  In the notation of 
p. \pageref{14.134}, 
$\tPER \approx B\sS$.]\\
\endgroup %%------------------------------------<<

\label{16.32}
Given two real inner product spaces
$
\begin{cases}
\ U\\
\ V
\end{cases}
$
with 
$
\begin{cases}
\ \dim U \leq \omega\\
\ \dim V \leq \omega
\end{cases}
$
, each equipped with the finite topology, let $\sI(U,V)$ be the set of linear isometries $U \ra V$.  
Endow $\sI(U,V)$ with the structure of a compactly generated Hausdorff space by relativising the compact open topology on $C(U,V)$ and taking its 
``$k$-ification''.\\

\label{16.34}
\textbf{\small LEMMA} \ 
Fix a real inner product space $V$ with $\dim V = \omega$ $-$then $\forall$ real inner product space $U$ with 
$\dim U \leq \omega$, $\sI(U,V)$ is contractible.

[Let $\{u_i\}$, $\{v_j\}$ be orthonormal bases for $U, \ V$ and let 
$
\begin{cases}
\ i_1, i_2: U \ra U \oplus U\\
\ j_1, j_2:V \ra V \oplus V
\end{cases}
$
be the inclusions onto the first, second summands.  
Choose a homotopy $F$ through isometries between $i_1$ and $i_2$ and choose a homotopy $\Phi$ through isometries $\id_V$ and $\phi:V \ra V$, where $\phi(v_j) = v_{2j}$.  Let
$h:V \ra V \oplus V$ be the isometry
$
\begin{cases}
\ h(v_{2j}) = (v_j,0)\\
\ h(v_{2j-1}) = (0,v_j)
\end{cases}
\hspace{-.5cm}, \ 
$
fix $f_0 \in \sI(U,V)$, and define $H:I\sI(U,V) \ra \sI(U,V)$  by  
$
H(f,t) = 
\begin{cases}
\ \Phi(2t) \circx f \hspace{3.35cm} (0 \leq t \leq 1/2)\\
\ h^{-1}\circx (f \oplus f_0) \circx F(2t - 1) \hspace{0.5cm} (1/2 \leq t \leq 1)
\end{cases}
$
$-$then 
$H(f,0) = f$, $H(f,1/2) = \phi \circx f =$ $h^{-1} \circx h \circx \phi \circx f =$ $h^{-1} \circx j_1 \circx f =$ 
$h^{-1} \circx (f \oplus f_0) \circx i_1$, and 
$H(f,1) = h^{-1} \circx (f \oplus f_0) \circx i_2 =$ $h^{-1} \circx (f_0 \oplus f_0) \circx i_2$, 
which is independent of $f$.]\\

\label{16.45}
\begingroup%%----------------------------------->>
\fontsize{9pt}{11pt}\selectfont
\textbf{\small FACT} \ 
Suppose that $\dim U < \omega$ and $\dim V = \omega$ $-$then $\sI(U,V)$ is a CW complex, hence the diagonal embedding 
$\sI(U,V) \ra \sI(U,V) \times_k \sI(U,V)$ is a closed cofibration (and, by the lemma, a homotopy equivalence).\\
\endgroup %%------------------------------------<<

\textbf{\small LEMMA} \ 
Fix a real inner product space $V$ with $\dim V = \omega$ $-$then the diagonal embedding 
$\sI(V,V) \ra \sI(V,V) \times_k \sI(V,V)$ is a closed cofibration.

[Write $V = \colimx V_n$, where $\forall \ n$, $\dim V_n = n$ and $V_n \subset V_{n+1} \subset V$.  
Consider the commutative diagram 

\[
\begin{tikzcd}[sep=large]
{\sI(V_n,V)} \ar{d}&{\sI(V_{n+1},V) } \ar{l}  \ar{d} \\
{\sI(V_n,V) \times_k \sI(V_n,V)}  &{\sI(V_{n+1},V) \times_k \sI(V_{n+1},V)} \ar{l} 
\end{tikzcd}
.
\]

Here, the horizontal arrows are \bCG fibrations and the vertial arrows are closed cofibrations and homotopy equivalences.  
Since $\sI(V,V) = \lim \sI(V_n,V)$, the assertion is a consequence of the generality infra.]\\

Application: The inclusion $\{\id_V\} \ra \sI(V,V)$ is a closed cofibration.\\

%%----------------------------------------------------------------------------------------------53
\begingroup%%----------------------------------->>
\fontsize{9pt}{11pt}\selectfont
\textbf{\small LEMMA} \ 
Let
\begin{tikzcd}%[sep=small]
{X_0} \ar{d} &{X_1} \ar{l} \ar{d} &{\cdots} \ar{l} \\
{Y_0} &{Y_1} \ar{l} &{\cdots} \ar{l}
\end{tikzcd}
be a commutative ladder of compactly generated spaces.  
Assume: $\forall \ n$, the horizontal arrows are \bCG fibrations and the vertical arrows are closed cofibrations and homotopy equivalences 
$-$then the induced map 
$\lim X_n \ra \lim Y_n$ is a closed cofibration and a homotopy equivalence.\\
\endgroup %%------------------------------------<<

Example: Let $V$ be a real inner product space with $\dim V = \omega$ and write $V^n$ for the orthogonal direct sum of the 
$n$ copies of $V$ $-$then the assignment $n \ra \sL_n = \sI(V^n,V)$ defines an $\tE_\infty$ operad $\sL$, the 
\un{linear isometries operad}. 
\index{linear isometries operad}\index{operad: linear isometries operad}

[The left action of $S_n$ on $V^n$  by permutations induces a free right action of $S_n$ on $\sL_n$, the unit 
$1 \in \sL_1$ is the identity map $V \ra V$, and 
$\Lambda:\sL_n \times_k (\sL_{j_1} \times_k \cdots \times_k \sL_{j_n}) \ra \sL_{j_1 + \cdots + j_n}$ 
sends 
$(f;g_1, \ldots, g_n)$ to $f \circx (g_1 \oplus \cdots \oplus g_n)$.]\\

\begingroup%%----------------------------------->>
\fontsize{9pt}{11pt}\selectfont
\textbf{\small EXAMPLE} \ 
Take $V = \R^\infty$ $-$then $\Oinf\Sinf \bS^0$ is an $\sL$-space.  
Indeed,  
$\Oinf\Sinf \bS^0 \approx$ $\colimx \Omega^n \bS^n =$ $\colim(\bS^n)^{\bS^n}$ 
and $\forall \ m,n$, there is a smash product pairing 
$(\bS^m)^{\bS^m} \times_k (\bS^n)^{\bS^n} \ra (\bS^m \#_k \bS^n)^{\bS^m \times_k \bS^n}$, 
where $\bS^m \times_k \bS^n = \bS^{m+n}$ 
(cf. p. \pageref{14.133a}).]
\vspi
[Note:  \ 
Boardman-Vogt\footnote[2]{\textit{SLN} \textbf{347} (1973), 207-217; 
see also May, \textit{SLN} \textbf{577} (1977), 9-24.} 
have given a systematic procedure for generating various classes of examples of 
$\sL$-spaces.]\\
\endgroup %%-----------------------------------<<

\textbf{\small LEMMA} \ 
Let $G$ be a finite group and let $X$ be a right $G$-space.  
Assume: Each $x \in X$ has a neighborhood $U$ with the property that $U \cdot g \cap U = \emptyset$ $\forall \ g \neq e$ $-$then the projection $X \ra X/G$ is a covering projection.\\

Application: Let $G$ be a finite group and let $X$ be a right $G$-space.  
Assume: The action of $G$ is free and $X$ is Hausdorff $-$then the projection $X \ra X/G$  is a covering projection.

[Note: \  Subject to these conditions on $X$, given any other right $G$-space $Y$, the product $X \times Y$ satisfies the hypotheses of the lemma, as does $X \times_k X$, hence the projection $X \times Y \ra (X \times Y)/G$  is a covering projection, as is $X \times_k Y \ra (X \times_k Y)/G$.]\\

\begin{proposition} \ %39
Let $\sO \ra \sP$ be a morphism of $\tE_\infty$ operads $-$then $\forall \ X$, the induced map 
$\sO[X] \ra \sP[X]$ is a weak homotopy equivalence.
\end{proposition}

[Consider the \cd
\[
\begin{tikzcd}%[sep=small]
{\sO_n \times_{S_n} X^n} \ar{d} &{\sO_n \times_{S_n} X_*^n} \ar{l} \ar{d} \ar{r} &{\sO_{n-1}[X]} \ar{d} \\
{\sP_n \times_{S_n} X^n}  &{\sP_n \times_{S_n} X_*^n} \ar{l} \ar{r}  &{\sP_{n-1}[X]}
\end{tikzcd}
.
\]
%%----------------------------------------------------------------------------------------------54
Arguing inductively, the arrow $\sO_{n-1}[X] \ra \sP_{n-1}[X]$ is a weak homotopy equivalence.  
But the same is also true of the other two vertical arrows (compare the long exact sequences in the homotopy of the relevant covering projections).  Therefore, since the horizontal arrows on the left are closed cofibrations, it follows that $\sO_n[X] \ra \sP_{n}[X]$ is a weak homotopy equivalence 
(cf. p. \pageref{14.135}), 
thus $\sO[X] \ra \sP[X]$ is a weak homotopy equivalence 
(cf. p. \pageref{14.136}).]\\

Example: Let 
$
\begin{cases}
\ \sO^\prime\\
\ \sO\pp
\end{cases}
$
be $\tE_\infty$ operads $-$then their product $\sO^\prime \times \sO\pp$ is an $\tE_\infty$ operad and $\forall \ X$, the arrows 
$(\sO^\prime \times \sO\pp)[X] \ra $
$
\begin{cases}
\ \sO^\prime[X]\\
\ \sO\pp[X]
\end{cases}
$
induced by the projections 
$\sO^\prime \times \sO\pp \ra $
$
\begin{cases}
\ \sO^\prime\\
\ \sO\pp
\end{cases}
$
are weak homotopy equivalences.

\label{14.151}
Example: Let $\sO$ be an $\tE_\infty$ operad $-$then $\abs{\sin \sO}$ is an $\tE_\infty$ operad 
(cf. p. \pageref{14.137}) 
and 
$\forall \ X$, the arrow $\abs{\sin \sO}[X] \ra \sO[X]$ is a weak homotopy equivalence.

[Note: \  Viewed as a creation operator, $\sO$ need not be universal (but $\abs{\sin \sO}$ is).]\\

\begingroup%%----------------------------------->>
\fontsize{9pt}{11pt}\selectfont
\label{16.64}
\textbf{\small FACT} \ 
Let 
$
\begin{cases}
\ \sC\\
\ \sD
\end{cases}
$
be creation operators, where $\forall \ n$, 
$
\begin{cases}
\ \sC_n\\
\ \sD_n
\end{cases}
$
is a compactly generated Hausdorff space and the action of $S_n$ is free.  Suppose given an arrow $\phi:\sC \ra \sD$ such that
$\forall \ n$, $\phi_n:\sC_n \ra \sD_n$ is a weak homotopy equivalence $-$then $\forall \ X$, $\phi$ induces a weak homotopy equivalence $\sC[X] \ra \sD[X]$.
\vspi
[Note:  \ By the same token, if $f:X \ra Y$ is a weak homotopy equivalence, then $\sC f:\sC[X] \ra \sC[Y]$ is a weak homotopy equivalence provided that $\forall \ n$, $\sC_n$ is a compactly generated Hausdorff space and the action of $S_n$ is free.]\\
\endgroup %%------------------------------------<<

\begin{proposition} \ %40
let $\sO$ be an $\tE_\infty$ operad $-$then every $\sO$-space $X$ is a homotopy associative, homotopy commutative H-space.
\end{proposition}

[To define the product, fix $f_2 \in \sO_2$ and consider $\theta_2(f_2,-|):X^2 \ra X$ (up to homotopy, 
the product is independent of the choice of $f_2 \in \sO_2$ ).]

[Note: \  If $X \ra Y$ is a morphism of $\sO$-spaces, then $X \ra Y$ is a morphism of H-spaces.]\\

\label{14.143}
\begingroup%%----------------------------------->>
\fontsize{9pt}{11pt}\selectfont
\textbf{\small EXAMPLE} \ 
Let $\sO = \tPER \approx B\bS$ and take $f_2 = e \in S_2 \subset XS_2$ $-$then with this choice for the product, every $\sO$-space is a homotopy commutative cofibered monoid in \bCG.\\
\endgroup %%------------------------------------<<

Working in the compactly generated category, let $X$ be a homotopy associative, homotopy commutative H-space $-$then a 
\un{group completion} 
\index{group completion (of a homotopy associative, homotopy commutative H-space)} 
of $X$ is a morphism $X \ra Y$ of H-spaces, where $Y$ is homotopy associative and 
$\pi_0(Y)$ is a group, such that 
$\ov{\pi_0(X)} \approx \pi_0(Y)$ and $H_*(X;\bk)[\pi_0(X)^{-1}] \approx H_*(Y;\bk)$ for every commutative ring \bk with unit.

Example: Let $G$ be a cofibered monoid in \bCG.  
Assume: $G$ is homotopy commutative $-$then according to Proposition 16 and the group completion theorem, 
the arrow $G \ra \Omega B G$ is a group completion.
\\

%%----------------------------------------------------------------------------------------------55
\begingroup%%----------------------------------->>
\fontsize{9pt}{11pt}\selectfont
\textbf{\small EXAMPLE}
Take $X = \Q$ (discrete topology), $Y = \Q$ (usual topology) $-$then the identity map $X \ra Y$ is a group completion but it is not a homotopy equivalence.
\vspi
[Note: \  Suppose that $X \ra Y$ is a group completion, where 
$
\begin{cases}
\ X\\
\ Y
\end{cases}
$
are pointed compactly generated CW spaces $-$then $X \ra Y$ is a weak homotopy equivalence if $\pi_0(X)$ is a group.  
Proof: One has $\pi_0(X) \approx$ $\ov{\pi_0(X)} \approx$ $\pi_0(Y)$ and there are homotopy equivalences 
$
\begin{cases}
\ X \ra X_0 \times \pi_0(X)\\
\ Y \ra Y_0 \times \pi_0(Y)
\end{cases}
$
, where 
$
\begin{cases}
\ X_0\\
\ Y_0
\end{cases}
$ 
is the path component of the identity element, thus the assertion follows from Dror's Whitehead theorem.]\\
\endgroup %%------------------------------------<<

\begingroup%%----------------------------------->>
\fontsize{9pt}{11pt}\selectfont
\textbf{\small EXAMPLE} 
Given a permutative category \bC, let $C^+$ be the simplicial object in \bCAT defined by 
$C_n^+ = \ds\prod\limits_1^{n+2} \bC$, where
\[
d_i(X_0,X_0^\prime,X_1, \ldots,X_{n}) = 
\begin{cases}
\ (X_0\otimes X_1,X_0^\prime \otimes X_1,X_2, \ldots,X_n) \hspace{1.31cm} (i = 0)\\
\ (X_0,X_0^\prime,X_1,\ldots, X_i \otimes X_{i+1}, \ldots,X_n) \hspace{0.75cm}  (0 < i < n)\\
\ (X_0,X_0^\prime,X_1, \ldots,X_{n-1}) \hspace{2.65cm} (i = n)
\end{cases}
,
\]
$s_i(X_0,X_0^\prime,X_1,\ldots,X_n) =$ $(X_0,X_0^\prime,X_1,\ldots, X_i,e,X_{i+1},\ldots,X_n)$ $-$then
there is a functor $\bC \ra \gro_{\bDelta^\OP}C^+$ and 
\label{14.163b}
Thomason\footnote[2]{\textit{Math. Proc. Cambridge Philos. Soc.} \textbf{85} (1979), 91-109.} 
has shown that the arrow 
$B\bC \ra B(\gro_{\bDelta^\OP}C^+)$ is a group completion.\\
\endgroup %%------------------------------------<<

\label{14.176}
\begingroup%%----------------------------------->>
\fontsize{9pt}{11pt}\selectfont
\textbf{\small EXAMPLE}  \ 
Let $X$ be a monoidal compactly generated simplicial space.  Assume: $X$ satisfies the cofibration condition and $X_1$ is homotopy commutative $-$then the arrow $X_1 \ra \Omega \aX$ is a group completion 
(Quillen\footnote[3]{\textit{Memoirs Amer. Math. Soc.} \textbf{529} (1994), 89-105.}).\\
\endgroup %%------------------------------------<<

\label{14.188}
\textbf{\small LEMMA}  \ 
Let $X$ be a homotopy associative, homotopy commutative H-space.  Suppose that $X \ra Y$ is a morphism of H-spaces, where $Y$ is homotopy associative and $\pi_0(Y)$ is a group, such that $\ov{\pi_0(X)} \approx \pi_0(Y)$ and 
$H_*(X;\bk)[\pi_0(X)^{-1}] \approx H_*(Y;\bk)$ for all prime fields \bk $-$then the arrow $X \ra Y$ is a group completion.\\

\begingroup%%----------------------------------->>
\fontsize{9pt}{11pt}\selectfont
\textbf{\small SUBLEMMA}  \ 
Let 
$
\begin{cases}
\ K\\
\ L
\end{cases}
$ 
be pointed CW complexes, $f:K \ra L$ a pointed continuous function.  
Assume: $f$ is a pointed homology equivalence $-$then 
$\Sigma f: \Sigma K \ra \Sigma L$ is a pointed homotopy equivalence.
\vspi
[Given $(X,x_0)$ in $\bCW_*$, let $X_{i_0}$, $X_i$ $(i \in I)$ be its set of path components, where $x_0 \in X_{i_0}$.  
Choose a vertex $x_i$ in each $X_i$ $-$then up to pointed homotopy, 
$\Sigma X = \ds\bigvee\limits_i \Sigma X_i \vee \Sigma \pi_0(X)$.]\\
\endgroup %%------------------------------------<<

\begingroup%%----------------------------------->>
\fontsize{9pt}{11pt}\selectfont
\textbf{\small LEMMA}  \ 
Let 
$
\begin{cases}
\ X\\
\ Y
\end{cases}
, \ Z
$
be \dsp pointed CW spaces in $\bCG_{*\bc}$, $f:X \ra Y$ a pointed homology equivalence.  
Suppose that $Z$ is a homotopy associative H-space such that $\pi_0(Z)$ is a group $-$then the precomposition arrow 
$f^*:[Y,Z] \ra [X,Z]$ is bijective.
\vspi
%%----------------------------------------------------------------------------------------------56
[Take $Z$ path connected and fix a retraction $JZ \ra Z$.  Since $[\Sigma Y,\Sigma Z] \approx [\Sigma X,\Sigma Z]$, 
the arrow 
$[Y,\Omega \Sigma Z] \ra [X,\Omega \Sigma Z]$ is bijective, so the assertion is true for $JZ$ (cf. Proposition 19).  
Now use the commutative diagram
\begin{tikzcd}%[sep=small]
{[Y,JZ]} \ar{r} \ar{d} &{[X,JZ]} \ar{d} \\
{[Y,Z]}  \ar{r} &{[X,Z]}
\end{tikzcd}
to see that the assertion is true for \mZ.]
\vspi
[Note: \  To define a retraction $JZ \ra Z$, make a choice for associating iterated products.  
Continuity is ensured if the homotopy  unit is a strict unit, which can always be arranged (since $Z \vee Z \ra Z \times_k Z$ is a closed cofibration 
(cf. p. \pageref{14.138})).]\\
\endgroup %%------------------------------------<<

\begingroup%%----------------------------------->>
\fontsize{9pt}{11pt}\selectfont
\textbf{\small FACT} \ 
Let 
$
X, \ 
\begin{cases}
\ Y_1\\
\ Y_2
\end{cases}
$
be \dsp pointed CW spaces in $\bCG_{*\bc}$.  Assume: $\pi_0(X) = \Z_{\geq 0}$ and 
$
\begin{cases}
\ X \ra Y_1\\
\ X \ra Y_2
\end{cases}
$
are group completions $-$then $\exists$ a pointed homotopy equivalence $Y_1 \ra Y_2$.
\\
\endgroup %%------------------------------------<<
\vspace{0.25cm}

\index{Theorem: May's Group Completion Theorem}
\index{May's Group Completion Theorem}
\textbf{\small MAY'S GROUP COMPLETION THEOREM} \ 
Let $(X,x_0)$ be a wellpointed compactly generated space with $\{x_0\} \subset X$ closed.  
Assume: $X$ is $\Delta$-separated $-$then 
$m_\infty:\bBV^\infty[X] \ra \Omega^\infty\Sigma^\infty X$
is a group completion.

[Note: \  When specialized to a path connected $X$, one recovers Proposition 33.]
\\

Homological calculations of this sort have their origins in the work of 
Dyer-Lashof\footnote[2]{\textit{Amer. J. Math.} \textbf{84} (1962), 35-88.}.  
Details are in 
May\footnote[3]{\textit{SLN} \textbf{533} (1976), 39-59.}.
\label{14.163a}

Example: $X$ $\Delta$-cofibered $\implies$ $\Omega^\infty\Sigma^\infty X$ $\Delta$-cofibered 
(cf. p. \pageref{14.139}).  
And: $\Omega^\infty\Sigma^\infty X$ is a $\bBV^\infty$-space.  The composite 
$\bBV^\infty[X] \ra \bBV^\infty[\Omega^\infty\Sigma^\infty X] \ra \Omega^\infty\Sigma^\infty X$
is $m_\infty$, the arrow of May.  It is a morphism of $\bT_{\bBV^\infty}$-algebras.
\\

\begin{proposition} \ %41
Let $\sO$ be an $\tE_\infty$ operad $-$then there is a functor 
$G:\sO\text{-}\bSP \ra \bCG_{*\bc}$ and a natural transformation $\id \ra G$ such that for every $\sO$-space $X$, the arrow 
$X \ra GX$ is a group completion.
\end{proposition}

[The product $\sO \times \tPER$ is an $\tE_\infty$ operad and $X$ is an $\sO \times \tPER$-space (through the projection 
$\sO \times \tPER \ra \sO$).  Consider the arrows 
$X \la B(\sO \times \tPER;\sO \times \tPER;X) \ra B(\tPER;\sO \times \tPER;X)$ 
in the category of $\bT_{\sO \times \tPER}\text{-algebras}$.  
The generalities on 
p. \pageref{14.140} ff. 
imply that the arrow 
$X \la B(\sO \times \tPER;\sO \times \tPER;X)$ 
is a homotopy equivalence 
(cf. p. \pageref{14.141}) 
and Propositions 38 and 39 imply that the arrow
$B(\sO \times \tPER;\sO \times \tPER;X) \ra B(\tPER;\sO \times \tPER;X)$ 
is a weak homotopy equivalence 
(cf. p. \pageref{14.142}).  
Since 
$B(\tPER;\sO \times \tPER;X)$ is a $\tPER\text{-space}$, it is a homotopy commutative cofibered monoid in \bCG (cf. 
%%----------------------------------------------------------------------------------------------57
p. \pageref{14.143}).  Put 
$GX = \Omega BB(\tPER; \sO \times \tPER;X)$ and let $X \ra GX$ be the composite 
$X \ra B(\sO \times \tPER; \sO \times \tPER;X) \ra B(\tPER;\sO \times \tPER;X) \ra GX$.]\\

\begingroup%%----------------------------------->>
\fontsize{9pt}{11pt}\selectfont
\textbf{\small FACT} \ 
Let $\sO$ be an $\tE_\infty$ operad.  Suppose that $A \ra X$ is a closed cofibration, where $A, \ X$ are \dsp $\sO$-spaces $-$then $GA \ra GX$ is a closed cofibration.
\vspi
[The arrow 
$B(\tPER;\sO \times \tPER;A) \ra B(\tPER;\sO \times \tPER;X)$ is a closed cofibration 
(cf. p. \pageref{14.144} 
$\&$ 
p. \pageref{14.145}).]\\
\endgroup %%------------------------------------<<

\begin{proposition} \ %42
Let $\sO$ be an $\tE_\infty$ operad such that $\forall \ n$, $\sO_n$ is an $S_n$-CW complex $-$then $\forall$ 
$\Delta$-cofibered $X$, $\sO[X]$ is $\Delta$-cofibered.
\end{proposition}

[By induction, $\forall \ n$, $\sO_n[X]$ is $\Delta$-cofibered 
(cf. p. \pageref{14.146}).  
Therefore 
$\sO[X] = \colimx \sO_n[X]$ is $\Delta$-cofibered 
(cf. p. \pageref{14.147}).]

[Note: \  If $\sO$ is an $\tE_\infty$ operad, then $\abs{\sin \sO}$ is an $\tE_\infty$ operad such that $\forall \ n$,  
$\abs{\sin \sO_n}$ is an $S_n$-CW complex.]
\\

Given an $\tE_\infty$ operad $\sO$, put $\sO^\infty = \sO \times \xBV^\infty$ $-$then every $\sO$-space $X$ is an 
$\sO^\infty$-space.  
On the other hand, $\abs{\sin X}$ is a $\abs{\sin \sO}$-space, hence is a $\abs{\sin \sO^\infty}$-space.  
The arrows 
$\abs{\sin \sO^\infty}[\abs{\sin X}] \ra \abs{\sin \xBV^\infty}[\abs{\sin X}]$, \ 
$\abs{\sin \xBV^\infty}[\abs{\sin X}] \ra \xBV^\infty[\abs{\sin X}]$ 
are weak homotopy equivalences (cf. Proposition 39), thus the composite 
$\abs{\sin \sO^\infty}[\abs{\sin X}] \ra \Omega^\infty \Sigma^\infty \abs{\sin X}$ is a group completion.

[Note: \  The diagram
\[
\begin{tikzcd}%[sep=small]
{\abs{\sin X}} \ar{d}\ar{r}  &{B(\abs{\sin \sO^\infty};\abs{\sin \sO^\infty};\abs{\sin X})}  \ar{d} \\
{X} \ar{r} &{B(\sO^\infty;\sO^\infty;X)}
\end{tikzcd}
\]
commutes.  Here, the horizontal arrows are homotopy equivalences and the vertical arrows are weak homotopy equivalences.]\\

\begin{proposition} \ %43
Let $\sO$ be an $\tE_\infty$ operad.  Suppose that $X$ is an $\sO$-space $-$then the arrow 
$B(\abs{\sin \sO^\infty};\abs{\sin \sO^\infty};\abs{\sin X}) \ra 
B(\Omega^\infty \Sigma^\infty;\abs{\sin \sO^\infty};\abs{\sin X})$ 
is a morphism of $\abs{\sin \sO^\infty}$-spaces 
(cf. p. \pageref{14.148}) 
and a group completion.
\end{proposition}

[Consider the commutative diagram
\[
\begin{tikzcd}%[sep=small]
{\abs{\sin X}} \ar{d}  
&{B(\abs{\sin \sO^\infty};\abs{\sin \sO^\infty};\abs{\sin X})}  \ar{l} \ar{d} \ar{r}
&{B(\Omega^\infty \Sigma^\infty;\abs{\sin \sO^\infty};\abs{\sin X})} \ar{d}\\
{G\abs{\sin X}}
&{B(G\abs{\sin \sO^\infty};\abs{\sin \sO^\infty};\abs{\sin X})}  \ar{l} \ar{r}
&{B(G\Omega^\infty \Sigma^\infty;\abs{\sin \sO^\infty};\abs{\sin X})}
\end{tikzcd}
.
\]
%%----------------------------------------------------------------------------------------------58
The arrow \ 
${\abs{\sin X}} \lla B(\abs{\sin \sO^\infty};\abs{\sin \sO^\infty};\abs{\sin X})$ 
is a homotopy equivalence, as is the arrow \ 
${G\abs{\sin X}} \la B(G\abs{\sin \sO^\infty};\abs{\sin \sO^\infty};\abs{\sin X})$.  
But 
$\abs{\sin X} \ra G\abs{\sin X}$ 
is a group completion, so 
$B(\abs{\sin \sO^\infty};\abs{\sin \sO^\infty};\abs{\sin X}) \ra$ 
$B(G\abs{\sin \sO^\infty};\abs{\sin \sO^\infty};\abs{\sin X})$
is a group completion.  
Since $\Omega^\infty \Sigma^\infty$ preserves closed cofibrations between $\Delta$-cofibered objects 
(cf. p. \pageref{14.149}), 
Proposition 42 implies that 
$\barr(\Omega^\infty \Sigma^\infty;\bT_{\abs{\sin \sO^\infty}};\abs{\sin X})$ 
satisfies the cofibration condition (see the proof of Proposition 38).  
Analogous remarks apply to 
$\barr(G\Omega^\infty \Sigma^\infty;\bT_{\abs{\sin \sO^\infty}};\abs{\sin X})$ 
and 
$\barr(G\abs{\sin \sO^\infty};\bT_{\abs{\sin \sO^\infty}};\abs{\sin X})$.  
Therefore the arrows 
$B(\Omega^\infty \Sigma^\infty;\abs{\sin \sO^\infty};\abs{\sin X}) \ra$ 
$B(G\Omega^\infty \Sigma^\infty;\abs{\sin \sO^\infty};\abs{\sin X})$, \ 
$B(G\abs{\sin \sO^\infty};\abs{\sin \sO^\infty};$\\$\abs{\sin X}) \ra$  
$B(G\Omega^\infty \Sigma^\infty;\abs{\sin \sO^\infty};\abs{\sin X})$
induce isomorphisms in homology $\forall \ k$ (cf. Proposition 10) and the assertion follows.]\\

Maintaining the preceding assumptions, put $\sO^q = \sO \times \xBV^q$.\\

\textbf{\small LEMMA} \ 
Let $\sO$ be an $\tE_\infty$ operad.  Suppose that $X$ is a \dsp $\sO$-space $-$then the arrow 
$B(\Omega^\infty \Sigma^\infty;\abs{\sin \sO^\infty};\abs{\sin X}) \ra B(\Omega^\infty \Sigma^\infty;\sO^\infty;X)$ 
is a weak homotopy equivalence.

[Since \ 
$B(\Omega^\infty \Sigma^\infty;\abs{\sin \sO^\infty};\abs{\sin X})$ \  $\approx$ \ 
$\colimx B(\Omega^q \Sigma^q;\abs{\sin \sO^q};\abs{\sin X})$, \ 
$B(\Omega^\infty \Sigma^\infty;$ $\sO^\infty;X)$ \ $\approx$ \ 
$\colimx B(\Omega^q \Sigma^q;\sO^q;X)$, \ 
where \ 
$B(\Omega^q \Sigma^q;\abs{\sin \sO^q};\abs{\sin X})$ 
$\ra$ 
$B(\Omega^{q+1} \Sigma^{q+1};$ $\abs{\sin \sO^{q+1}};\abs{\sin X})$, 
$B(\Omega^q \Sigma^q;\sO^q;X) \ra B(\Omega^{q+1} \Sigma^{q+1};\sO^{q+1};X)$
are closed embeddings, it will be enough to show that $\forall \ q$, the arrow 
$B(\Omega^q \Sigma^q;\abs{\sin \sO^q};\abs{\sin X}) \ra B(\Omega^q \Sigma^q;\sO^q;X)$
is a weak homotopy equivalence 
(cf. p. \pageref{14.150}).  \ 
However, bearing in mind Proposition 38, $\forall \ n$, 
$\abs{\sin \sO^q}^n[\abs{\sin X}] \ra (\sO^q)^n[X]$ is a weak homotopy equivalence 
(cf. p. \pageref{14.151}), 
hence
$\forall \ n$, $\Omega^q \Sigma^q \abs{\sin \sO^q}^n[\abs{\sin X}] \ra \Omega^q \Sigma^q (\sO^q)^n[X]$ 
is a weak homotopy equivalence 
(cf. p. \pageref{14.152} ff.), 
so the generality on 
p. \pageref{14.153} 
is applicable.]

[Note: \  While $B(\Omega^\infty \Sigma^\infty;\abs{\sin \sO^\infty};\abs{\sin X})$ is in $\bCG_{*\bc}$, 
this is not a prior the case of 
$B(\Omega^\infty \Sigma^\infty;\sO^\infty;X)$ 
(both space are, of course, \dsp).  Still, 
$B(\Omega^\infty \Sigma^\infty;\sO^\infty;X)$ is an $\sO^\infty$-space in $\bCG_{*}$ (see remarks on 
p. \pageref{14.154}).]\\

\begin{proposition} \ %44
Let $\sO$ be an $\tE_\infty$ operad.  Suppose that $X$ is a \dsp $\sO$-space $-$then the arrow 
$B(\sO^\infty;\sO^\infty;X) \ra B(\Omega^\infty \Sigma^\infty;\sO^\infty;X)$ 
is a morphism of $\sO^\infty$-spaces 
(cf. p. \pageref{14.155}) 
and a group completion.
\end{proposition}

[In the \cd
\[
\begin{tikzcd}%[sep=small]
{B(\abs{\sin \sO^\infty};\abs{\sin \sO^\infty};\abs{\sin X})}  \ \ar{d} \ar{r}
&{B(\Omega^\infty \Sigma^\infty;\abs{\sin \sO^\infty};\abs{\sin X})}\ar{d}\\
{B(\sO^\infty;\sO^\infty;X)} \ar{r}
&{B(\Omega^\infty \Sigma^\infty;\sO^\infty;X)} 
\end{tikzcd}
,
\]
%
%%----------------------------------------------------------------------------------------------59
the vertical arrows are weak homotopy equivalences and, by Proposition 43, the top horizontal arrow is a group completion.]

[Note: \  When $X$ is path connected, the arrow 
$B(\sO^\infty;\sO^\infty;X) \ra B(\Omega^\infty \Sigma^\infty;\sO^\infty;X)$ 
is a weak homotopy equivalence (cf. Proposition 33).]\\

A 
\un{spectrum}
\index{spectrum} 
\bX is a sequence of pointed \dsp compactly generated spaces $X_q$ and pointed homeomorphisms 
$X_q \overset{\sigma_q}{\lra} \Omega X_{q+1}$.  
\bSPEC 
\index{\bSPEC} 
is the category whose objects are the spectra and whose morphisms 
$\bff:\bX \ra \bY$ 
are sequences of pointed continuous functions $f_q:X_q \ra Y_q$ such that the diagram
\begin{tikzcd}%[sep=small]
{X_q} \ar{d} \ar{r}{f_q} &{Y_q} \ar{d}\\
{\Omega X_{q+1}} \ar{r}[swap]{\Omega f_{q+1}} &{\Omega Y_{q+1}}
\end{tikzcd}
commutes $\forall \ q$.

[Note: \  The indexing begins at 0.]

\label{16.1}
There is a functor $\bU^\infty:\bSPEC \ra \bDelta$-$\bCG_*$ that sends $\bX = \{X_q\}$ to $X_0$.  It has a left adjoint
$\bQ^\infty:\bDelta\text{-}\bCG_* \ra \bSPEC$ defined by $(\bQ^\infty X)_q = \Omega^\infty \Sigma^\infty \Sigma^q X$.

[Note: \  The repetition principle implies that 
$\Omega \Omega^\infty \Sigma^\infty \Sigma^{q+1} X \approx$ 
$\Omega \Omega^\infty \Sigma^\infty \Sigma \Sigma^{q} X \approx$ 
$\Omega^\infty \Sigma^\infty \Sigma^{q} X$.]

An 
\un{infinite loop space}
\index{infinite loop space} 
is a pointed \dsp compactly generated space in the image of 
$\bU^\infty$.  
Example: $\forall \ X$, $\Omega^\infty \Sigma^\infty X$ is an infinite loop space.  
Every infinite loop space is a $\xBV^\infty$-space (in the extended sense of the word 
(cf. p. \pageref{14.156})).\\

\begingroup%%----------------------------------->>
\fontsize{9pt}{11pt}\selectfont
\textbf{\small EXAMPLE} \ 
If $\bX = \{X_q\}$ is a spectrum such that $X_0$ is wellpointed, then $\forall \ q$, there is an arrow 
$\Omega^q\Sigma^q \Omega^q X_q \ra \Omega^q X_q$, 
from which an arrow 
$\Omega^\infty \Sigma^\infty X_0 \ra X_0$.  Viewing the $\xBV^\infty$-space $X_0$ as a $\bT_{\xBV^\infty}$-algebra, its structural morphism $\xBV^\infty[X_0] \ra X_0$ is the composite 
$\xBV^\infty[X_0] \overset{m_\infty}{\lra}$ $\Omega^\infty \Sigma^\infty X_0 \ra$ $X_0$.\\
\endgroup %%------------------------------------<<

\label{14.169a}
A spectrum \bX is said to be 
\un{connective}
\index{connective (spectrum)} 
if $X_1$ is path connected and $X_q$ is $(q-1)$-connected $(q > 1)$.

Example:  Given an $\tE_\infty$ operad $\sO$ and a \dsp $\sO$-space $X$, the assignment 
$q \ra B_qX = \colimx \Omega^nB(\Sigma^{n+q};\sO^{n+q};X)$ 
specifies a connective spectrum $\bB X$.

[To check that $B_qX$ is $\Delta$-separated, it need only be shown that the arrow 
$\Omega^nB(\Sigma^{n+q};$\\$\sO^{n+q};X) \ra \Omega^{n+1}B(\Sigma^{n+1+q};\sO^{n+1+q};X)$
is a closed embedding 
(cf. p. \pageref{14.157}).  
To see this, note that
$\Sigma B(\Sigma^{n+q};\sO^{n+q};X) \approx$ $B(\Sigma^{n+1+q};\sO^{n+q};X)$ 
(cf. p. \pageref{14.158}) 
and 
$B(\Sigma^{n+1+q};\sO^{n+q};X)$ $\ra$ $B(\Sigma^{n+1+q};\sO^{n+1+q};X)$ is a closed embedding (in fact, a closed cofibration).  Therefore
$B(\Sigma^{n+q};\sO^{n+q};X) \ra$ $\Omega B(\Sigma^{n+1+q};\sO^{n+1+q};X)$ is a closed embedding.  
And: $\Omega^n$ preserves closed embeddings.]

[Note: \  That $\bB X$ is connective is implied by the generalities on 
p. \pageref{14.159}.]\\
Remark: The arrow 
$\colimx B(\Omega^q\Sigma^{q};\sO^{q};X) \ra \colimx \Omega^q B(\Sigma^{q};\sO^{q};X)$ 
is a morphism of \dsp $\sO^\infty$-spaces 
(cf. p. \pageref{14.160} ff.) 
and a weak homotopy equivalence.

%%----------------------------------------------------------------------------------------------60
[In fact, 
$\barr(\Omega^q\Sigma^q;\bT_{\sO^q};X) \ = \  \Omega^q \barr(\Sigma^q;\bT_{\sO^q};X)$, \ 
so \ 
$\abs{\Omega^q \barr(\Sigma^q;\bT_{\sO^q};X)} \ra \Omega^q \abs{\barr(\Sigma^q;\bT_{\sO^q};X)}$
is a weak homotopy equivalence 
(cf. p. \pageref{14.161}).]\\

\begin{proposition} \ %45
Let $\sO$ be an $\tE_\infty$ operad.  Suppose that $X$ is a \dsp $\sO$-space $-$then the composite
$X \ra B(\sO^\infty;\sO^\infty;X) \ra B(\Omega^\infty\Sigma^\infty;\sO^\infty;X) \ra B_0X$ 
is a group completion.
\end{proposition}

[Taking into account Proposition 44, this follows from what has been said above.]

[Note: \  It is not claimed that $B_0X$ is wellpointed.]\\

Therefore every \dsp $\sO$-space $X$ group completes to an infinite loop space.

[Note: \  Consequently, if $X$ is path connected, then $X$ has the weak homotopy type of an infinite loop space.]

Remark: Proposition 45 is true for any \dsp $\sO^\infty$-space (same argument).

[Note: \  Observe that every $\xBV^\infty$-space is an $\sO^\infty$-space.]\\

\begingroup%%----------------------------------->>
\fontsize{9pt}{11pt}\selectfont
\textbf{\small EXAMPLE} \ 
Specializing to $\sO = \tPER$, one sees that the classifying space $B\bC$ of a permutative category \bC group completes to an infinite loop space.\\
\endgroup %%------------------------------------<<

\begin{proposition} \ %46
Let $\sO$ be an $\tE_\infty$ operad.  
Suppose that $\bX = \{X_q\}$ is a spectrum such that $X_0$ is wellpointed $-$then there is a morphism 
$\bb:\bB X_0 \ra \bX$ in \bSPEC such that the diagram
\begin{tikzcd}[sep=large]
{B(\sO^\infty;\sO^\infty;X_0)} \ar{d} \ar{r}&{{B(\Omega^\infty\Sigma^\infty;\sO^\infty;X_0)} } \ar{d}\\
{X_0} &{B_0X_0} \ar{l}{b_0}
\end{tikzcd}
commutes.
\end{proposition}

[Proceeding formally, use the arrow 
$B(\Sigma^{n+q};\sO^{n+q};\Omega^{n+q};X_{n+q}) \ra X_{n+q}$
to define 
$b_q:B_qX_0 \ra X_q$.]

[Note: \  It is a corollary that the composite 
$X_0 \ra B_0X_0 \overset{b_0}{\ra} X_0$ 
is the identity.  Another corollary is that $b_0$ is a weak homotopy equivalence provided that $X_0$ is path connected.]\\

\begin{proposition} \ %47
Let $\sO$ be an $\tE_\infty$ operad $-$then $\forall \ \Delta$-cofibered $X$ in $\bCG_*$, there is a morphism 
$\bff:\bB\sO^\infty[X] \ra \bQ^\infty X$ of spectra such that $\forall \ q$, 
$f_q:B_q\sO^\infty[X] \ra \Omega^\infty\Sigma^\infty\Sigma^q X$ is a pointed homotopy equivalence.
\end{proposition}

[The arrow 
$B(\Sigma^{n+q};\sO^{n+q};\Omega^{n+q}[X]) \ra \Sigma^{n+q}X$ 
is a pointed homotopy equivalence 
(cf. p. \pageref{14.162} ff.).  
Apply $\Omega^n$ and let $n \ra \infty$.  
In this connection, the assumption that $X$ is $\Delta$-cofibered guarantees that 
$\Omega^n\Sigma^{n+q}X \ra \Omega^{n+1}\Sigma^{n+1+q}X$ 
is a closed cofibration 
(cf. p. \pageref{14.163}), 
so Proposition 15 in $\S 3$ is applicable.]

[Note: \  Working through the definitions, one finds that \bff is equal to the composite 
$\bB\sO^\infty[X] \ra \bB \xBV^\infty[X] \ra \bB \Omega^\infty\Sigma^\infty X \overset{\bb}{\ra} \bQ^\infty X$.]\\

%%----------------------------------------------------------------------------------------------61
\label{14.120}
\label{17.7}
\begingroup%%----------------------------------->>
\fontsize{9pt}{11pt}\selectfont
\textbf{\small EXAMPLE} \ 
Take $\sO = \tPER \approx BS$ and let $X = \bS^0$ $-$then 
$\sO[\bS^0] \approx$ $\abs{M_\infty} =$ $\ds\coprod\limits_{n \geq 0} BS_n$ 
and the projection 
$\sO^\infty[\bS^0] \ra \sO[\bS^0]$
is a weak homotopy equivalence.  On the other hand, the composite
$\sO^\infty[\bS^0] \ra B_0 \sO^\infty[\bS^0] \ra \Omega^\infty\Sigma^\infty \bS^0$ 
is a group completion (cf. Propositions 45 and 47), as is the arrow
$\abs{M_\infty} \ra \Omega B \abs{M_\infty}$.
Therefore 
$\Omega^\infty\Sigma^\infty \bS^0$ and $\Omega B \abs{M_\infty}$ 
have the same pointed homotopy type 
(cf. p. \pageref{14.163a}).  
The homotopy groups 
$\pi_*^s$ of $\Omega^\infty\Sigma^\infty \bS^0$ are the stable homotopy groups of spheres.  
Since 
$\Omega B \abs{M_\infty} \approx \Z \times BS_\infty^+$, 
it follows that 
$\pi_*^s \approx \pi_*(BS_\infty^+)$.  
Example: $\pi_1^s \approx \pi_1(BS_\infty^+) =$ $S_\infty/A_\infty \approx$ $\Z/2\Z$.  
There is also a connection with algebraic K-theory.  
Thus $S_\infty \subset \bG\bL(\Z)$, $A_\infty \subset \bE(\Z)$, so there is an arrow 
$B S_\infty^+ \ra B \bG\bL(\Z)^+$.  
The associated homomorphism 
$\pi_n^s \ra K_n(\Z)$ $(= \pi_n(B \bG\bL(\Z)^+))$ can be bijective (e.g., if $n = 1$) but in general is neither injective nor 
surjective 
(see Mitchell\footnote[2]{In: \textit{Algebraic Topology and its Applications}, G. Carlsson et al. (ed.), Springer Verlag (1994), 163-240 (cf. 182-183).} 
for a discussion and more information).
\vspi
[Note: Let $\bC = \bM_\infty = \iso  \bGamma$ $-$then another model for 
$\Omega^\infty \Sigma^\infty \bS^0$ is $B(\gro_{\bDelta\OP} C^+)$ 
(cf. p. \pageref{14.163b}.]\\
\endgroup %%------------------------------------<<

\begingroup%%----------------------------------->>
\fontsize{9pt}{11pt}\selectfont
\textbf{\small EXAMPLE} \ 
Given a discrete group $G$, form 
$S_\infty \ds\int G$ 
(cf. p. \pageref{14.164}) 
$-$then a model for the plus construction on 
$BS_\infty \ds\int G$ is the path component of $\Omega^\infty\Sigma^\infty BG_+$ 
containing the constant loop.  E.g.: When
$G = *$, \ $\Omega^\infty\Sigma^\infty BG_+$ is $\Omega^\infty\Sigma^\infty \bS^0$ 
and when
$G = \Z/2\Z$, \ $\Omega^\infty\Sigma^\infty BG_+$ is $\Omega^\infty\Sigma^\infty \bP^\infty(\R)_+$.\\
\endgroup %%------------------------------------<<

\index{$\bPi$ }
\index{category! $\bPi$}
$\bPi$ is the category whose objects are the finite sets 
$\bn \equiv \{0,1,\ldots,n\}$ $(n \geq 0)$ with base point 0 and whose morphisms are the base point preserving maps $\gamma:\bm \ra \bn$ such that 
$\#(\gamma^{-1}(j)) \leq 1$ $(1 \leq j \leq n)$.  
So: $\bGamma_\ini$ is a subcategory of $\bPi$ and $\bPi$ is a subcategory of $\bGamma$.

[Note: \  Let $(X,x_0)$ be a wellpointed compactly generated space with $\{x_0\} \subset X$ closed $-$then the formulas that define pow $X$ as a functor $\bGamma_\ini \ra \bCG_*$ serve to define pow $X$ as a functor 
$\bPi \ra \bCG_*$.]

A 
\un{category of operators} 
\index{category of operators} 
is a compactly generated category \bC such that $\Ob\bC \ra \Mor\bC$ is a closed cofibration, where $\Ob\bC = \Ob\bGamma$ (discrete topology), subject to the requirement that \bC contains $\bPi$ and admits an augmentation $\epsilon:\bC \ra \bGamma$ which restricts to the inclusion $\bPi \ra \bGamma$.  One writes $\bC(\bm,\bn)$ for the set of morphisms $\bm \ra \bn$.
Example: $\bGamma$ is a category of operators, as is $\bPi$.

Every category of operators is a \bCG-category.

[Note: \ A morphism of categories of operators is a continuous functor $F:\bC \ra \bD$ such that $F \bn \ra \bn$ for all \bn and
\begin{tikzcd}[sep=small]
&{\bGamma}\\
{\bC} \ar{ru}\ar{rr}{F} &&{\bD} \ar{lu}\\
&{\bPi} \ar{lu}\ar{ru}
\end{tikzcd}
commutes.]\\
\vspace{0.25cm}

%%----------------------------------------------------------------------------------------------62
\begingroup%%----------------------------------->>
\fontsize{9pt}{11pt}\selectfont
\textbf{\small FACT} \ 
Let \bC be a category of operators.  Suppose that $X$ is a right \bC-object and $Y$ is a left \bC-object $-$then 
$\barr(X;\bC;Y)$ satisfies the cofibration condition.\\
\endgroup %%------------------------------------<<

A 
\un{cofibered operad}
\index{cofibered operad } in \bCG is a reduced operad $\sO$ in \bCG for which the inclusion $\{1\} \ra \sO_1$ is a 
closed cofibration.  
Example: Every $\tE_\infty$ operad is a cofibered operad in \bCG.\\
\label{14.174}
Notation: Given morphisms $\gamma:\bm \ra \bn$, $\delta:\bn \ra \bp$ in $\bGamma$, let $\sigma_k(\delta,\gamma)$ 
be the permutation on $\#((\delta \circx \gamma)^{-1}(k))$ letters which converts the natural ordering of 
$(\delta \circx \gamma)^{-1}(k)$ to the ordering associated with 
$\bigcup\limits_{\delta(j) = k} \gamma^{-1}(j)$ (all elements of $\gamma^{-1}(j)$ precede all elements of 
$\gamma^{-1}(j^\prime)$ if $j < j^\prime)$ and each $\gamma^{-1}(j)$ has its natural ordering).\\

\begin{proposition} \ %48
Let $\sO$ be a cofibered operad in \bCG $-$then $\sO$ determines a category of operators $\widehat{\sO}$.
\end{proposition}

[Put \ 
$\widehat{\sO}(\bm,\bn) = \coprod\limits_{\gamma:\bm \ra \bn} \prod\limits_{1 \leq j \leq n} \sO(\#(\gamma^{-1}(j)))$ 
(cf. p. \pageref{14.165}).  \ 
Here composition 
$\widehat{\sO}(\bm,\bn) \times \widehat{\sO}(\bn,\bp) \ra \widehat{\sO}(\bm,\bp)$ \ 
is the rule 
$(\delta;g_1, \ldots, g_p) \circx (\gamma,f_1, \ldots, f_n) = (\delta \circx \gamma, h_1, \ldots, h_p)$, 
$h_k$ being $\Lambda(g_k;f_j(\delta(j) = k)) \cdot \sigma_k(\delta,\gamma)$ 
and 
$(\id_\bn;1, \ldots, 1)$ is the identity element in $\widehat{\sO}(\bn,\bn)$.  
The augmentation $\epsilon:\widehat{\sO} \ra \bGamma$ is obvious, viz. $\epsilon(\gamma;f_1, \ldots, f_n) = \gamma$.  
To define the inclusion $\bPi \ra \widehat{\sO}$, send $\gamma:\bm \ra \bn$ to $(\gamma;f_1, \ldots, f_n)$ where
$
\begin{cases}
\ f_j = 1 \quad (j \in \im \gamma)\\
\ f_j = * \quad (j \notin \im \gamma)
\end{cases}
.]
$
\\

Examples: 
(1) Let $\sO_n = *$ $\forall \ n$ $-$then $\Oh = \bGamma$;
(2) Let $\sO_0 = *$, $\sO_1 = \{1\}$, $\sO_n = \emptyset$ $(n > 1)$ $-$then $\Oh = \bPi$.

A
\un{$\bPi$-space}
\index{space! $\bPi$-space} 
is a functor $X:\bPi \ra \bCG_*$ and a 
\un{$\bPi$-map} 
\index{map! $\bPi$-map} 
is a 
natural transformation $f:X \ra Y$.

Given $n \geq 1$, there are projections $\pi_i:\bn \ra \bone$ $(i = 1, \ldots, n)$, where
$
\pi_i(j) = 
\begin{cases}
\ 1 \quad (i = j)\\
\ 0 \quad (i \neq j)
\end{cases}
.\ 
$
A $\bPi$-space $X$ is said to be 
\un{special} 
\index{special ($\bPi$-space)} 
if $X_0 = *$ and $\forall \ n \geq 1$, the arrow 
$X_n \ra X_1 \times_k \cdots \times_k X_1$ determined by the $\pi_i$ is a weak homotopy equivalence.

Given an injection $\gamma:\bm \ra \bn$, let $S_\gamma$ be the subgroup of $S_n$ consisting of those $\sigma$ such that $\sigma(\im \gamma) = \im \gamma$.  
A $\bPi$-space $X$ is said to be 
\un{proper} 
\index{proper ($\bPi$-space)} if $X_0 = *$ and $\forall$ $\gamma: \bm \ra \bn$ 
in $\bGamma_{\ini}$, $X_\gamma:X_m \ra X_n$ is a closed $S_\gamma$-cofibration (cf. infra).  
In particular: $* \ra X_n$ is a closed $S_n$-cofibration, so $\forall \ n$, $X_n$ is in $\bCG_{*\bc}$.

[Note: \  Associated with each $\sigma \in S_\gamma$ is a permutation $\widetilde{\sigma}\in S_m$ such that 
$\sigma \circx \gamma = \gamma \circx \widetilde{\sigma}$ 
and the assignment $\sigma \ra \widetilde{\sigma}$ is a homomorphism $S_\gamma \ra S_m$.  
Thus $X_m$ and $X_n$ are left $S_\gamma$-spaces and $X\gamma:X_m \ra X_n$ is equivariant.]

Example: $\forall \ X$ in $\bCG_{*\bc}$, pow $X$ is a proper special $\bPi$-space.\\

%%----------------------------------------------------------------------------------------------63
\begingroup%%----------------------------------->>
\fontsize{9pt}{11pt}\selectfont
Let $G$ be a finite group.  
Let $A$ and $X$ be left $G$-spaces $-$then an equivariant continuous function 
$i:A \ra X$ is said to be a 
\un{$G$-cofibration}
\index{G-cofibration} 
if it has the following property: 
Given any left $G$-space $Y$ and any pair $(F,h)$ of equivariant continuous functions 
$
\begin{cases}
\ F:X \ra Y\\
\ h:IA \ra Y
\end{cases}
$
such that $F \circx i = h \circx i_0$, there is an equivariant continuous function 
$H:IX \ra Y$ such that $F = H \circx i_9$ and $H \circx Ii = h$.
\vspi
[Note: \  Every $G$-cofibration is an embedding and the induced map $G\bs A \ra G\bs X$ is a cofibration.]
\vspi
The theory set forth in $\S 3$ has an equivariant analog 
(Boardman -Vogt\footnote[2]{\textit{SLN} \textbf{347} (1973), 231-239.}).  
For example, Proposition 1 in $\S 3$ becomes: Let $A$ be an invariant subspace of $X$ $-$then the inclusion 
$A \ra X$ is a $G$-cofibration iff $i_0X \cup IA$ is an equivariant retract of $IX$.  
The notion of an equivariant Str{\o}m structure on $(X,A)$ is clear and there is a $G$-cofibration characterization theorem.
\vspi
[Note:  \ A $G$-cofibration is thus a cofibration.]\\ 
\endgroup %%------------------------------------<<

\begingroup%%----------------------------------->>
\fontsize{9pt}{11pt}\selectfont
\textbf{\small EXAMPLE} \ 
Suppose that $(X,x_0)$ is in $\bCG_{*\bc}$ $-$then the inclusion $X_*^n \ra X^n$ is a closed $S_n$-cofibration.\\
\endgroup %%------------------------------------<<

\begingroup%%----------------------------------->>
\fontsize{9pt}{11pt}\selectfont
\textbf{\small LEMMA} \ 
Let $A$ be an invariant subspace of the left \mG-space $X$.  
Suppose that $A = A_1 \cup \ \cdots \cup A_n$, where each $A_i$ 
is closed in $X$, and suppose that $G$ operates on $\{1, \ldots, n\}$ in such a way that $g \cdot A_i = A_{g\cdot i}$.  
Put $A_S = \ds\bigcap\limits_{i \in S} A_i$ $(S \subset \{1, \ldots, n\})$ $-$then $A \ra X$ is a closed \mG-cofibration if 
$\forall \ S \neq \emptyset$, $A_S \ra X$ is a closed $G_S$-cofibration, $G_S \subset G$ the stabilizer of $S$.
\vspi
[Note: \  Take for $G$ the trivial group to recover Proposition 8 in $\S 3$ (with 2 replaced by $n$).]
\\
\endgroup %%------------------------------------<<

\label{14.168}
\begingroup%%----------------------------------->>
\fontsize{9pt}{11pt}\selectfont
\textbf{\small EXAMPLE} \ 
Let $X$ be a proper special $\bPi$-space.  Put $sX_{n-1} = s_0X_{n-1} \cup \cdots \cup  s_{n-1}X_{n-1}$, 
where $s_i = X\sigma_i$ and 
$
\sigma_i(j) = 
\begin{cases}
\ j \hspace{1.05cm} (j \leq i)\\
\ j + 1 \hspace{0.5cm} (j > i)
\end{cases}
(0 \leq i < n)
$
$-$then the inclusion $sX_{n-1}  \ra X_n$ is a closed $S_n$-cofibration.\\
\endgroup %%------------------------------------<<

Notation: $\pspsp$ is the category of proper special $\bPi$ spaces.
\index{category! $\pspsp$  (category of proper special $\bPi$ spaces)}\\

\begin{proposition} \ %49
Let $L$ be the functor from $\pspsp$ to $\bCG_{*\bc}$ that sends $X$ to $X_1$ and let $R$ be the functor from 
$\bCG_{*\bc}$ to $\pspsp$ that sends $X$ to pow $X$ $-$then $(L,R)$ is an adjoint pair.
\end{proposition}

[Note: \ The arrow of adjunction $LRX \ra X$ is the identity  and the arrow of adjunction $X \ra RLX$ has for its components the map induced by the $\pi_i$.]\\

Let \bC be a category of operators $-$then a 
\un{\bC-space} 
\index{space! \bC-space}  
is a continuous functor 
$X:\bC \ra \bCG_*$ and a 
\un{\bC-map}  
\index{map! \bC-map} 
is a natural transformation $f:X \ra Y$.\\

%%----------------------------------------------------------------------------------------------64
\begingroup%%----------------------------------->>
\fontsize{9pt}{11pt}\selectfont
Continuity in this context means that $\forall \ m,n$ the arrow $\bC(\bm,\bn) \times_k X_m \ra X_n$ is continuous.  
To clarify the matter, let $E = X_n^{X_m}$ (exponential object in \bCG), $E_* = X_n^{X_m}$ (pointed exponential object in $\bCG_*$) $-$then there is a commutative triangle 
$
\begin{tikzcd}%[sep=small]
{\bC(\bm,\bn) }  \ar{rd}\ar{r} &{E_*} \ar{d}\\
&{E} 
\end{tikzcd}
,
$
where $E_* \ra E$ is a \bCG-embedding.  
Thus the arrow $\bC(\bm,\bn) \ra E_*$ is continuous iff the arrow 
$\bC(\bm,\bn) \ra E$ is continuous or still, iff the arrow 
$\bC(\bm,\bn) \times_k X_m \ra X_n$ is continuous.\\
\endgroup %%------------------------------------<<

A \bC-space is said to be 
\un{special}
\index{special (\bC-space)} 
or 
\un{proper}
\index{proper (\bC-space)} 
if its restriction to $\bPi$ is special or proper.

Example: A 
\un{$\bGamma$-space} 
\index{space! $\bGamma$-space}
is an $\widehat{\sO}$-space, where $\sO_n = *$ $\forall \ n$.  
Every abelian monoid $G$ in \bCG gives rise to a special $\bGamma\text{-space}$ 
(cf. p. \pageref{14.166}), 
the
\un{$\bGamma\text{-nerve}$}
\index{nerve! $\bGamma\text{-nerve}$} of 
$G:\bGamma\text{-}\ner G$ 
\index{$\bGamma\text{-}\ner$}
(which is proper if $G$ is cofibered).\\

\textbf{\small LEMMA} \ 
Let $\sO$ be a cofibered operad in \bCG $-$then an \Ohs with underlying  space $\pow X$ determines and is determined by an $\sO$-space structure on \mX.

[To specify an $\sO$-space structure on $X$ is to specify a morphism $\sO \ra \End X$ of operads in \bCG 
(cf. p. \pageref{14.167}), 
from which an \Ohs $\Oh \ra \bCG_*$ with underlying space $\pow X$.  
Conversely, let $\gamma_n:\bn \ra \bone$ be the arrow $j \ra 1$ $(1 \leq j \leq n)$ and view $\sO_n$ as the component of 
$\gamma_n$ in $\Oh(\bn,\bone)$.  
Per an \Ohs with underlying space $\pow X$, restriction of 
$\Oh(\bn,\bone) \ra X^{X^n}$ to $\sO_n$ defines a morphism $\sO \ra \End X$ of operads in \bCG.]\\

Let $\sO$ be a cofibered operad in \bCG $-$then by restriction, $\Oh(-,\bn)$ defines a functor 
$\bPi^\OP \ra \bCG$ $\forall \ n \geq 0$.  Given a $\bPi$-space $X$, put $\Oh_n[X] = \Oh(-,\bn) \otimes_{\bPi}X$ 
(so $\Oh_0[X] = X_0$) and call $\Oh[X]$ the $\bPi$-space which takes $\bn$ to $\Oh_n[X]$.  
Composition in $\Oh$ leads to maps 
$\Oh(\bm,\bn) \times \Oh_m[X] \ra \Oh_n[X]$ or still, to an arrow 
$m_X:\Oh^2[X] \ra \Oh[X]$, while the identities in $\Oh$ induce an arrow $\epsilon_X:X \ra \Oh[X]$.   
Both arrows are natural in $X$ and with $T_{\Oh} = \Oh[?]$, it is seen that 
$\bT_{\Oh} = (T_{\Oh},m,\epsilon)$ is a triple in $\bPi,\bCG_*]$.

Notation: Let $\sE(\bm,\bn)$ be the set of base point preserving maps $\epsilon:\bm \ra \bn$ such that 
$\epsilon^{-1}(0) = \{0\}$ and $i \leq i^\prime$ $\implies$ $\epsilon(i) \leq \epsilon(i^\prime)$.  
Put $S_\epsilon = S_{\epsilon_1} \times \cdots \times S_{\epsilon_n} \subset S_m$, where 
$\epsilon_j = \#(\epsilon^{-1}(j))$.

[Note: \  Let $\sigma \in S_m$ $-$then $\epsilon \circx \sigma \in \sE(\bm,\bn)$ iff $\sigma \in S_\epsilon$.]\\

\begin{proposition} \ %50
Suppose that $X$ is a proper $\bPi$-space. Denote by $\Oh_{m,n}[X]$ the image of 
$\coprod\limits_{m^\prime \leq m} \Oh(\bm^\prime,\bn) \times_k X_{m^\prime}$ in $\Oh_{n}[X]$ $-$then 
$\Oh_{m,n}[X]$ is a closed subspace of $\Oh_{n}[X]$ 
%%----------------------------------------------------------------------------------------------65
and $\Oh_{n}[X]  = \colimx \Oh_{m,n}[X]$.  
In addition, the \cd
\[
\begin{tikzcd}%[sep=small]
{\coprod\limits_{\epsilon \in \sE(\bm,\bn)} \prod\limits_{1 \leq j \leq n} \sO_{\epsilon_j} \times_{S_\epsilon} sX_{m-1}}  
\ar{d} \ar{r}
&\Oh_{m-1,n}[X] \ar{d}\\
{\coprod\limits_{\epsilon \in \sE(\bm,\bn)} \prod\limits_{1 \leq j \leq n} \sO_{\epsilon_j} \times_{S_\epsilon} X_{m}}  
\ar{r}
&\Oh_{m,n}[X] 
\end{tikzcd}
\]
is a pushout square and the arrow $\Oh_{m-1,n}[X] \ra \Oh_{m,n}[X]$ is a closed cofibration.
\end{proposition}

[Note: \  For the definition of ``$s$'', see 
p. \pageref{14.168}.]\\

Remark: $X_n$ \dsp $\forall \ n$ $+$ $\sO_n$ \dsp $\forall \ n$ $\implies$ $\Oh_n[X]$ \dsp $\forall \ n$ 
(cf. p. \pageref{14.169}).\\

\begingroup%%----------------------------------->>
\fontsize{9pt}{11pt}\selectfont
\textbf{\small FACT} \ 
If $X$ is a proper $\bPi$-space, then $\Oh[X]$ is a proper $\bPi$-space and $\epsilon_X:X \ra \Oh[X]$ is a levelwise closed cofibration.\\
\endgroup %%------------------------------------<<

\begin{proposition} \ %51
Fix an $X$ in $\bCG_{*\bc}$ $-$then $L\Oh[RX]$ ($=\Oh_1[\pow X]$ (cf. Proposition 49)) $\approx \Oh[X]$ and 
$\Oh[RX] \approx R\Oh[X]$.\\
\end{proposition}

\textbf{\small LEMMA} \ 
Let $\sO$ be a cofibered operad in \bCG.  \ 
Assume:  $\forall \ n$, $\sO_n$ is a compactly generated Hausdorff space and the action of $S_n$ is free.  
Suppose given a $\bPi$-map $f:X \ra Y$ such that $\forall \ n$, $f_n:X_n \ra Y_n$ is a weak homotopy equivalence 
$-$then $\forall \ n$, $\Oh_nf:\Oh_n[X]$ $\ra$ $\Oh_n[Y]$ is a weak homotopy equivalence provided that $X$ and \mY are proper.

[This is a variant on the argument used in the proof of Proposition 39.]\\

\begin{proposition} \ %52
Let $\sO$ be a cofibered operad in \bCG.  Assume: $\forall \ n$, $\sO_n$ is a compactly generated Hausdorff space and the action of $S_n$ is free.  Suppose that $X$ is a \ps $\bPi$-space $-$then $\Oh[X]$ is a \ps $\bPi$-space.
\end{proposition}

[$\Oh[X]$ is necessarily proper (cf. supra).  To check that $\Oh[X]$ is special, consider the commutative diagram 
$
\begin{tikzcd}%[sep=small]
{\Oh_n[X]} \ar{d} \ar{r} &{\Oh_n[RLX]} \ar{d}\\
{(\Oh_1[X])^n} \ar{r} &{(\Oh_1[RLX])^n}
\end{tikzcd}
, \ 
$
bearing in mind the lemma and the fact that the arrow of adjunction $X \ra RLX$ is a levelwise  weak homotopy equivalence (Proposition 51 supplies an identification $\Oh_n[RLX] \approx (\Oh_1[RLX])^n)$.]\\

Application: Let $\sO$ be an $\tE_\infty$ operad 
$-$then the triple $\bT_{\Oh} = (T_{\Oh},m,\epsilon)$ in $[\bPi,\bCG_*]$ restricts to a triple in $\pspsp$ 
and its associated category of algebras is canonically isomorphic to the category $\psosp$ of proper special $\Oh$-spaces (cf. Proposition 37).\\

%%----------------------------------------------------------------------------------------------66
Suppose that $X$ is a simplicial $\bPi$-space $-$then the 
\un{realization}
\index{realization ($\bPi$-space)}
$\abs{X}$ of $X$ is the $\bPi$-space defined by $\abs{X}(\bn) = \abs{[m] \ra X_m(\bn)}$.

Example: If $\sO$ is an $\tE_\infty$ operad and if $X$ is a \ps \Ohs, then the realization 
$B(\Oh;\Oh;X)$ of $\barr(\bT_{\Oh};\bT_{\Oh};X)$ is a proper special \Ohs.\\

\textbf{\small LEMMA} \ 
Suppose that $F:\bCG_{*\bc} \ra \bV$ is a right $\TO$-functor 
$-$then $F \circx L:\pspsp$ $\ra$ $\bV$ is a right $\TOh$ functor.

[The relevant natural transformation 
$F \circx L \circx \TOh \ra F \circx L$ is the composite 
$FL\Oh[X] \ra$ 
$FL\Oh[RLX] =$ 
$FLR\sO[LX] =$
%$F\sO[LX] \overset{\rho_{LX}}{\longrightarrow} FLX$.]\\
$F\sO[LX]$
$\overset{\rho_{LX}}{\xrightarrow{\hspace*{1cm}}}$ 
$FLX$.]\\

Let $\sO$ be an $\tE_\infty$ operad, $F:\bCG_{*\bc} \ra \bCG_{*\bc} $ a right $\TO$-functor $-$then for any $\TOh$-algebra $X$, $\barr(F \circx L;\TOh;X)$ is a simplicial object in $\bCG_{*\bc}$ and one writes $B(F \circx L;\Oh;X)$ for its geometric realization.

[Note: \  It is clear that there is a version of Proposition 38 applicable to this situation.]\\

\begin{proposition} \ %53
let $\sO$ be an $\tE_\infty$ operad $-$then there is a functor $U$ from $\pOhp$  to $\pOhp$  and a functor $V$ from $\pOhp$  to 
$\sO$-\bSP plus $\Oh$-maps $X \la UX \ra RVX$ natural in $X$ such that $X \la UX$ is a levelwise homotopy equivalence and 
$UX \ra RVX$ is a levelwise weak homotopy equivalence.
\end{proposition}

[Put 
$UX = B(\widehat{\sO};\widehat{\sO};X)$ and 
$VX = B(T_{\sO} \circx L;\Oh;X)$.  
So, in obvious notation 
$RVX = B(R \circx T_{\sO} \circx L; \widehat{\sO};X)$  and the arrow $UX \ra RVX$ is defined in terms of the arrows
$\widehat{\sO}_n[X] \ra (\sO[X_1])^n$, hence is a levelwise weak homotopy equivalence (see the proof of Proposition 52).]

[Note: \  Suppose that $X$ is an $\sO$-space $-$then 
$B(\widehat{\sO};\widehat{\sO};RX) \approx$ 
$RB(\widehat{\sO};\widehat{\sO};X)$ ($\implies$ 
$LB(\widehat{\sO};\widehat{\sO};RX)$ $\approx$ 
$B(\sO;\sO;X)$) and 
$VRX \approx$ $B(\sO;\sO;X)$ (cf. Proposition 51).]\\

Remarks: 
(1) \mX \dsep $\implies$ $UX$, $VX$ \dsep; 
(2) $X \ra UX$ is not an $\widehat{\sO}$-map (but it is a $\Pi$-map).\\

\label{14.191}
\begingroup%%----------------------------------->>
\fontsize{9pt}{11pt}\selectfont
\textbf{\small FACT} \ 
Let $\sO$ be an $\tE_\infty$ operad, $\epsilon:\widehat{\sO} \ra \bGamma$ the augmentation $-$then there are functors 
$\epsilon^*: \psg \ra \pOhp$,
$\epsilon_*: \pOhp \ra \psg$ respecting the 
$\Delta$-separation condition and an $\widehat{\sO}$-map 
$UX \ra \epsilon^*\epsilon_* X$ natural in $X$ which is a levelwise weak homotopy equivalence.\\
\endgroup %%------------------------------------<<

Let $\sO$ be an $\tE_\infty$ operad $-$then there is a functor \bB from the category of \dso to the category of connective spectra 
(cf. p. \pageref{14.169a}) 
and this functor can be extended to the category of \dsep \ps $\Oh$-spaces by writing 
$B_qX = \colimx  \Omega^nB(\Sigma^{n+q}L;\Oh^{n+q};X)$.  To see that this prescription really is an extension, consider
%%----------------------------------------------------------------------------------------------67
any $\Delta$-separated $\sO$-space  
$X:B(\Sigma^{n+q}L;\Oh^{n+q};RX) \approx$ $B(\Sigma^{n+q};\sO^{n+q};X)$
(cf. Proposition 51) $\implies$ $\bB RX \approx \bB X$.\\

\begin{proposition} \ %54
Let $\sO$ be an $\tE_\infty$ operad.  
Suppose that $X$ is a \dsep \ps $\Oh$-space $-$then the composite 
$B(\bT_{\sO^\infty}\circx L;\Oh^\infty;X)$ $\ra$ $B(\Oinf\Sinf L;\Oh^\infty;X) \ra$ $B_0X$ 
is a group completion.
\end{proposition}

[Rework the discussion leading up to Proposition 45.]\\

Let $\sO$ be a cofibered operad in \bCG $-$then an 
\un{infinite loop space machine}
\index{infinite loop space machine} 
on $\Oh$ consists of a functor \bB From the category of \dsep \ps $\Oh$-spaces to the category of connective spectra, 
a functor $K$ from the category of \dsep \ps $\Oh$-spaces to the category of homotopy associative, homotopy commutative 
%\bH 
H-spaces, a natural transformation $L \ra K$ such that $\forall \ X$, the arrow $LX  \ra KX$ is a weak homotopy equivalence, and a natural transformation $K \ra B_0$ such that $\forall \ X$, $KX \ra B_0X$ is a group completion.\\

\begin{proposition} \ %55
Let $\sO$ be an $\tE_\infty$ operad $-$then there exists an infinite loop space machine on $\Oh$, the \un{May machine}. 
\index{May machine}
\end{proposition}

[Take \bB as above and put  \ $KX = B(\bT_{\sO^\infty}\circx  L; \Oh^\infty;X)$.  \  The composite 
$X \ra B(\Oh^\infty;\Oh^\infty;X)$ $\ra$ $RB(\bT_{\sO^\infty}\circx L;\Oh^\infty;X)$
is a levelwise weak homotopy equivalence, hence $LX \ra KX$ is a weak homotopy equivalence.  On the other hand, thanks to Proposition 54, the composite
$KX \ra B(\Oinf\Sinf L;\Oh^\infty;X) \ra B_0X$ 
is a group completion.]\\

Let $\sO$ be an $\tE_\infty$ operad $-$then, using the augmentation $\epsilon:\Oh \ra \bGamma$, a \dsep \ps 
$\bGamma$-space can be regarded as a \dsep \ps $\Oh$-space.  
Therefore an infinite loop space machine on $\Oh$ defines an infinite loop space machine on $\bGamma$.  
However, there is another ostensibly very different method for generating connective spectra from  
\dsep \ps $\bGamma$-spaces which is completely internal and makes no reference to operads.  
The question then arises: Are the spectra thereby produced in some sense the ``same''?  
As we shall see, the answer is ``yes'' 
(cf. Proposition 62), a corollary being that infinite loop space machines associated with distinct $\tE_\infty$ operads $\sO$ and 
$\sP$ attach the ``same'' spectra to a \dsep \ps $\bGamma$-space.\\

\textbf{\small LEMMA} \ 
$\bDelta^\OP$ is isomorphic to the category whose objects are the $\bn_+$ $(j < *, 0 \leq j \leq n)$ 
and whose morphisms are the order preserving maps $\alpha:\bm_+ \ra \bn_+$ such that $\alpha(0) = 0$ and $\alpha(*) = *$.\\

The composite $[n] \ra \bn_+ \ra \bn_+/0 \sim * \equiv \bn$ defines a functor $\bS[1]:\bDelta^\OP \ra \bGamma$.

%%----------------------------------------------------------------------------------------------68
[Note: \  To justify the notation, observe that the pointed simplicial set 
$\bDelta^\OP \ra \bGamma \subset \bSET_*$ thus displayed is in fact a model for the simplicial circle 
(cf. p. \pageref{14.170}).]\\

\begingroup%%----------------------------------->>
\fontsize{9pt}{11pt}\selectfont
\textbf{\small EXAMPLE} \ \ 
Suppose that $\alpha:[n] \ra [m]$ is a morphism in $\bDelta$.  \ 
Put $\gamma = \bS[1]\alpha$ \ (so $\gamma:\bm \ra \bn$ is a morphism in $\bGamma$) $-$then $\gamma$ is given by 
$\gamma^{-1}(j) = \{i: \alpha(j - 1) < i \leq \alpha(j)\}$ $(1 \leq j \leq n)$,  
$\gamma^{-1}(0) = \bm - \ds\bigcup\limits_{j = 1}^n \gamma^{-1}(j)$.  
Examples: 
(1) The $\sigma_i:[n+1] \ra [n]$ of 
p. \pageref{14.171} 
are sent by $\bS[1]$ to the $\sigma_i:\bn \ra \bn + \bone$ of 
p. \pageref{14.172} 
$(n \geq 0, 0 \leq i \leq n)$; 
(2) The $\pi_i:[1] \ra [n]$ of 
p. \pageref{14.173} 
are sent by $\bS[1]$ to the $\pi_i:\bn \ra \bone$ of 
p. \pageref{14.174} 
$(n \geq 1, 1 \leq i \leq n)$.\\
\endgroup %%------------------------------------<<

Notation: Call 
$\ov{\pow}X$ 
\index{$\ov{\pow}X$} 
the functor $\bGamma^\OP \ra \bCG_*$ corresponding to a cofibrant $X$ in $\bCG_*$ 
(standard model category structure).\\

\begingroup%%----------------------------------->>
\fontsize{9pt}{11pt}\selectfont
\textbf{\small EXAMPLE} \ 
Let $Y:\bGamma \ra \bCG$ be a functor $-$then $\forall \ X$, one can form 
$\barr(\ov{\pow}X;\bGamma;Y)$ and denoting by $B(X;\bGamma;Y)$ its geometric realization, there is a canonical arrow 
$B(X;\bGamma;Y) \ra \ov{\pow}X \otimes_{\bGamma} Y$ 
(cf. p. \pageref{14.175}).  
Example: $\forall \ n$, 
$(PY) \bn  \approx B(\bn; \bGamma ; Y)$, 
$Y(\bn) \approx \ov{\pow} \otimes_{\bGamma} Y$  and the arrow of evaluation $(PY) \bn  \ra Y(\bn)$ is a homotopy equivalence.\\
\endgroup %%------------------------------------<<

\label{18.8}
\begingroup%%----------------------------------->>
\fontsize{9pt}{11pt}\selectfont
\textbf{\small EXAMPLE} \ 
Let $\zeta:\bGamma_\ini^\OP \ra \bGamma$ be the functor which is the identity on objects and sends 
$\gamma:\bm \ra \bn$ to $\zeta\gamma:\bn \ra \bm$, where 
$\zeta\gamma(j) = \gamma^{-1}(j)$ if $\gamma^{-1}(j) \neq \emptyset$, 
$\zeta\gamma(j) = 0$ if $\gamma^{-1}(j) = \emptyset$  $-$then for any $X$ in $\bCG_{*c}$, 
$\ov{\pow} \ X \circx \zeta^\OP = \pow \ X$.  
The assignment $\bn \ra \ohc \pow  \bn$ defines a functor $\gamma^\infty:\bGamma \ra \bCG$.  
And: $\ohc  \pow  X \approx \ov{\pow} X \otimes_{\bGamma} \gamma^\infty$.
\vspi
[The left Kan extension of 
$B(-\backslash \bGamma_\ini)$ along $\zeta$ is $\gamma^\infty$, hence
$\ohc  \pow  X \approx$ 
$B(-\backslash \bGamma_\ini) \otimes_{\bGamma_\ini}  \pow  X  \approx$\\
$\pow \ X  \otimes_{\bGamma_\ini^\OP} B(-\backslash \bGamma_\ini) \approx$ 
$\ov{\pow}  X \circx \zeta^\OP \otimes_{\bGamma_\ini^\OP} B(-\backslash \bGamma_\ini) \approx$ 
$\ov{\pow}  X \otimes_{\bGamma} \gamma^\infty$.]
\vspi
[Note: \  Let $X$ be a pointed connected CW complex or a pointed connected ANR $-$then the homotopy colimit theorem says that $\ohc \ \pow \ X$ and $\Oinf \Sinf X$ have the same homotopy type, thus by the above, 
$\ov{\pow} \ X \otimes_{\bGamma} \gamma^\infty$ and $\Oinf \Sinf X$ have the same homotopy type.]\\
\endgroup %%------------------------------------<<

\textbf{\small LEMMA} \ 
Relative to $\bS[1]^\OP:\bDelta \ra \bGamma^\OP$, $\lan\Delta^? \approx \ov{\pow} \bS^1$.
\\

Let $X: \bGamma \ra \bCG$ be a functor $-$then the 
\un{realization} 
\index{realization (simplicial spaces)} 
$\aX_{\bGamma}$ of $X$ is by definition $\abs{X \circx \bS[1]}$, the geometric realization of $X \circx \bS[1]$.  
And:  
$\abs{X \circx \bS[1]} = X \circx \bS[1] \otimes_{\bDelta} \Delta^? \approx$ 
$X \otimes_{\bGamma^\OP} \lan\ \Delta^? \approx$  
$X \otimes_{\bGamma^\OP} \ov{\pow} \bS^1 \approx$ 
$\ov{\pow} \bS^1 \otimes_{\bGamma} X$.

Example: Let $G$ be an abelian cofibered monoid in \bCG $-$then 
$(\bGamma$-$\ner\bG) \circx \bS[1] = \ner \bG$ $\implies$ 
$\abs{\bGamma\text{-}\ner\bG}_{\bGamma} = BG$.\\

\begingroup%%----------------------------------->>
\fontsize{9pt}{11pt}\selectfont
Given an abelian cofibered monoid $G$ in \bCG, let 
$SP^\infty(?;G)$ 
\index{$SP^\infty(?;G)$ } 
be the functor 
$\bCG_{*c} \ra \bCG_{*c}$ that sends $X$ to 
$\ov{\pow} X \otimes_{\bGamma} \bGamma$-$\ner \bG$ $-$then $SP^\infty(X;G)$ is an abelian cofibered monoid in 
\bCG, the \un{infinite symmetric} 
%%----------------------------------------------------------------------------------------------69
\un{product} 
\index{infinite symmetric product (with coefficients in abelian cofibered monoid $G$)}
on $(X,x_0)$ with coefficients in $G$.  
Example: Take $G = \Z_{\geq 0}$ to see that 
$SP^\infty X \approx \ds\int^\bn X^n \times_k SP^\infty \bn \approx$ 
$SP^\infty(X;\Z_{\geq 0})$ (the choice $G = \Z$ leads to the free abelian compactly generated group on $(X,x_0)$).\\
\endgroup %%------------------------------------<<

\begingroup%%----------------------------------->>
\fontsize{9pt}{11pt}\selectfont
\textbf{\small LEMMA} \ 
$\forall \ X$, $Y$, $SP^\infty(X\#_kY;G) \approx$ $SP^\infty(X;SP^\infty(Y;G))$ 
(isomorphism of abelian monoids in \bCG).\\
\endgroup %%------------------------------------<<

\begingroup%%----------------------------------->>
\fontsize{9pt}{11pt}\selectfont
\textbf{\small EXAMPLE} \ 
Let $G$ be an abelian cofibered monoid in \bCG $-$then 
$SP^\infty(\bS^0;G) \approx G$, 
$SP^\infty(\bS^1;G) \approx BG$, and in general 
$SP^\infty(\bS^{n+1};G) \approx B^{(n+1)}G$, where $B^{(n+1)}G = B(B^{(n)}G)$.]
\vspi
[Representing \ $\bS^{n+1}$ as the smash product \ $\bS^n \#_k \bS^1$, the lemma implies that 
$SP^\infty(\bS^{n+1};G)$ $\approx$ \  $SP^\infty(\bS^{n};BG)$.]\\
\endgroup %%------------------------------------<<

Let $X$ be a \ps $\bGamma$-space $-$then $X \circx \bS[1]$ satisifies the cofibration condition.  Moreover, if 
$X \circx \bS[1]$ is monoidal, then $X_1$ is a homotopy associative, homotopy commutative H-space and the arrow 
$X_1 \ra \Omega \aX_{\bGamma}$ is a group completion 
(cf. p. \pageref{14.176}).

[Note: \  $\sin X$ is an object in $\bGamma\bSISET_*$ 
(cf. p. \pageref{14.177}) 
and $\abs{\sin X}$ is a \ps $\bGamma$-space.  The simplicial space $\abs{\sin X} \circx \bS[1]$ is monoidal and there is a commutative diagram 
\begin{tikzcd}%[sep=small]
{\abs{\sin X_1} } \ar{d} \ar{r} &{\Omega \norm{\sin X}_{\bGamma}} \ar{d}\\
{X_1} \ar{r} &{\Omega \aX_{\bGamma}}
\end{tikzcd}
.  Since the vertical arrows are weak homotopy equivalences (Giever-Milnor 
(cf. p. \pageref{14.178} ff.)) 
and since the arrow 
$\abs{\sin X_1} \ra \Omega \norm{\sin X}_{\bGamma}$ is a group completion, it follows that the arrow 
$X_1 \ra \Omega \aX_{\bGamma}$ is a 
\un{weak group completion} 
\index{weak group completion}
($X_1$ is not necessarily an H-space) (but $\forall \ \bk$, $\pi_0(X_1)$ is a central submonoid of $H_*(X_1;\bk)$ and 
$H_*(X_1;\bk)[\pi_0(X_1)^{-1}] \approx H_*(\Omega  \aX_{\bGamma};\bk)$.]\\
\label{14.183}
\label{14.186}

\begingroup%%----------------------------------->>
\fontsize{9pt}{11pt}\selectfont
Remark:  If \bC is a pointed category with finite products and if $X$ is a special $\bGamma$-object in \bC (obvious definition), then $X_1$ is an abelian monoid object in \bC 
(cf. p. \pageref{14.179}).\\
\endgroup %%------------------------------------<<

\begingroup%%----------------------------------->>
\fontsize{9pt}{11pt}\selectfont
\textbf{\small FACT} \ 
Let $X$ be a \ps $\bGamma$-space.  Assume: $\forall \ n \geq 1$, the arrow 
$X_n \ra X_1 \times_k \cdots \times_k X_1$ determined by the $\pi_i$ is an $S_n$-equivariant homotopy equivalence $-$then there exists an abelian cofibered monoid $G$ in \bCG and a levelwise homotopy equivalence $X \ra \bGamma$-$\ner \bG$.\\
\endgroup %%------------------------------------<<

\textbf{\small LEMMA} \ 
Let $X$ be a \ps $\bGamma$-space $-$then $X_1$ path connected $\implies$ $\aX_{\bGamma}$ 
simply connected and $X_1$ $n$-connected  $\implies$ $\aX_{\bGamma}$ $(n+1)$-connected 
(cf. p. \pageref{14.180}).\\

Let 
$\bGamma \overset{\nu_n}{\lra} \bGamma \times \bGamma$ be the functor defined by 
$\bp \ra (\bn,\bp)$ on objects and $\gamma \ra (\id_{\bn},\gamma)$ on morphisms.  
Given a \ps $\bGamma$-space $X$, 
call $\ov{X}_n$ the composite 
$\bGamma \overset{\nu_n}{\lra}$ 
$\bGamma \times \bGamma \overset{\#}{\lra}$
$\bGamma   \overset{X}{\lra}$
$\bCG_*$,
$\#$ being the smash product 
(cf. p. \pageref{14.181}).  
So: 
$\ov{X}_n(\bp) = X_{np}$ and $\ov{X}_n$ is a \ps $\bGamma$-space.

%%----------------------------------------------------------------------------------------------70
[Note: \  Suppose that $\gamma: \bm \ra \bn$ is a morphism in $\bGamma$.  
Set $\gamma_p = \gamma \# \id_\bp:\bm\bp \ra \bn\bp$ $-$then the $\gamma_p$ induce a $\bGamma$-map 
$\ov{X}_m \ra \ov{X}_n$, thus $\ov{X}$ is a functor from $\bGamma$ to $\psg$.]

The 
\un{classifying space} 
\index{classifying space (of a \ps $\bGamma$-space)} 
of a \ps $\bGamma$-space $X$ is the \ps $\bGamma$-space $BX$ which takes \bn to $B_nX = \abs{\ov{X}_n}_{\bGamma}$.  
In particular: $B_1X = \aX_{\bGamma}$ is path connected, hence $B_1X \ra \Omega \abs{BX}_{\bGamma}$ is a weak homotopy equivalence.\\

\begingroup%%----------------------------------->>
\fontsize{9pt}{11pt}\selectfont
\textbf{\small FACT} \ 
Let $G$ be an abelian cofibered monoid in \bCG 
$-$then the classifying space of the $\bGamma$-nerve of $\bG$ is the $\bGamma$-nerve of $\bB \bG$.\\
\endgroup %%------------------------------------<<

Notation:  Given a \ps $\bGamma$-space $X$, write 
$B^{(0)}X = X$,
$B^{(q+1)}X = B(B^{(q)}X)$, and put
$S_0X = \Omega \aX_{\bGamma}$, 
$S_{q+1}X =  \abs{B^{(q)}X}_{\bGamma}$ $(q \geq 0)$.\\

\begingroup%%----------------------------------->>
\fontsize{9pt}{11pt}\selectfont
\textbf{\small EXAMPLE} \ 
Let $X$ be a \ps $\bGamma$-space $-$then $\forall \ q > 0$, $S_qX \approx \ov{\pow} \bS^q \otimes_{\bGamma} X$.\\
\endgroup %%------------------------------------<<

A 
\un{prespectrum}
\index{prespectrum} 
\bX is a sequence of pointed \dsep compactly generated spaces $X_q$ and pointed continuous functions 
$X_q \overset{\sigma_q}{\lra} \Omega X_{q+1}$.  
\bPRESPEC 
\index{\bPRESPEC} 
is the category whose objects are the prespectra and whose morphisms 
$\bff:\bX \ra \bY$ are the sequences of pointed continuous functions $f_q:X_q \ra Y_q$ such that the diagram
\begin{tikzcd}%[sep=small]
{X_q} \ar{d} \ar{r}{f_q} &{Y_q} \ar{d}\\
{\Omega X_{q+1}} \ar{r}[swap]{\Omega f_{q+1}} &{\Omega Y_{q+1}}
\end{tikzcd}
commutes $\forall \ q$.  Every spectrum is a prespectrum.

[Note: \  The indexing begins at 0.]\\

\begingroup%%----------------------------------->>
\fontsize{9pt}{11pt}\selectfont
\textbf{\small EXAMPLE} \ 
Let $\sO$ be an $\tE_\infty$ operad.  Suppose that $X$ is a \dsep \ps $\widehat{\sO}$-space $-$then the assignment 
$q \ra B(\Sigma^qL;\widehat{\sO}^q;X)$ is a prespectrum.\\
\endgroup %%------------------------------------<<

Remark: \bPRESPEC is complete and cocomplete (limits and colimts are calculated levelwise).\\

\begin{proposition} \ %56
Equip $\bDelta$-\bCG$_*$ with its singular structure $-$then \bPRESPEC is a model category if weak equivalences and fibrations are levelwise, a cofibration $\bff:\bX \ra \bY$ being a levelwise cofibration with the additional property that 
$\forall \ q$, the arrow $P_{q+1} \ra Y_{q+1}$ is a cofibration, where $P_{q+1}$ is defined by the pushout square
\begin{tikzcd}%[sep=small]
{\Sigma X_q} \ar{d} \ar{r}&{\Sigma Y_q} \ar{d}\\
{X_{q+1}} \ar{r}&{P_{q+1}}
\end{tikzcd}
.
\end{proposition}

[Note: \  In the presence of the condition on the $P_{q+1} \ra Y_{q+1}$,  
to describe the cofibrations in \bPRESPEC, it suffices to require that $f_0:X_0 \ra Y_0$ be a cofibration.]\\

%%----------------------------------------------------------------------------------------------71
\begingroup%%----------------------------------->>
\fontsize{9pt}{11pt}\selectfont
If \bC is a category and if $F, G: \bC \ra \bPRESPEC$ are functors, then a natural transformation 
$\Xi:F \ra G$ is a function that assigns to each $X \in \Ob\bC$ an element $\Xi_X \in \Mor(FX,GX)$ natural in $X$.  
Using the notation 
$
\begin{cases}
\ FX = \{F_{X,q}\}\\
\ GX = \{G_{X,q}\}
\end{cases}
, \ 
$
$\Xi_X = \{\Xi_{X,q}\}$, the fact that $\Xi_X \in \Mor(FX,GX)$  is expressed by the commutativity of 
\begin{tikzcd}%[sep=small]
{F_{X,q}} \ar{d}[swap]{\sigma_{F,q}} \ar{r}{\Xi_{X,q}} &{G_{X,q}} \ar{d}{\sigma_{G,q}}\\
{\Omega F_{X,q+1}} \ar{r}[swap]{\Omega\Xi_{X,q+1}} &{\Omega G_{X,q+1}}
\end{tikzcd}
$\forall \ q$.  A 
\un{pseudo natural transformation}
\index{pseudo natural transformation} 
$\Xi:F \ra G$ is a function that assigns to each $X \in \Ob\bC$  a sequence of pointed continuous functions 
$\Xi_{X,q}:F_{X,q} \ra G_{X,q}$ natural in $X$ 
and a sequence of pointed homotopies $H_{X,q}$ between 
$\Omega\Xi_{X,q+1} \circx  \sigma_{F,q}$ 
and 
$\sigma_{G,q} \circx \Xi_{X,q}$ 
natural in \mX 
(thus natural $\implies$ pseudo natural (constant homotopies)).  
A pseudo natural homotopy between pseudo natural transformations 
$\Xi_0$, $\Xi_1: F \ra G$ is a pseudo natural transformation $\Upsilon:F\#I_+ \ra G$ such that 
$
\begin{cases}
\ \Upsilon \circx i_0 = \Xi_0\\
\ \Upsilon \circx i_1 = \Xi_1
\end{cases}
, \ 
$
where $(F\#I_+)(X) = FX \# I_+$ $(\{F_{X,q}\#I_+\})$ 
(cf. p. \pageref{14.182}).
\vspi
[Note: \  A natural (pseudo natural) transformation $\Xi$ is called a \ 
\un{natural (pseudo natural) weak} \  \un{equivalence}
\index{natural weak equivalence}
\index{pseudo natural weak equivalence}
if the $\Xi_{X,q}$ are weak homotopy equivalences.]\\
\endgroup %%------------------------------------<<

\label{16.10}
\label{17.9}
\label{17.12}
\index{Cylinder Construction (\bPRESPEC)}
\begingroup%%----------------------------------->>
\fontsize{9pt}{11pt}\selectfont
\textbf{\small EXAMPLE \  (\un{Cylinder Construction})} \ 
There is a functor $M:\bPRESPEC \ra \bPRESPEC$ with the property that $\forall \ \bX$, the arrows 
$(M\bX)_q \ra \Omega(M\bX)_{q+1}$ are closed embeddings.  And: \\
\indent\indent $(M_1)$ \quad $\exists$ a natural transformation $r:M \ra \id$ such that 
$\forall \ \bX$, $r_{\bX,q}:(M\bX)_q \ra X_q$ is a pointed homotopy equivalence.\\
\indent\indent $(M_2)$ \quad $\exists$ a pseudo natural transformation $j:\id \ra M$ such that 
$\forall \ \bX$, $j_{\bX,q}:X_q \ra (M\bX)_q$ is a pointed homotopy equivalence.\\
\indent\indent $(M_3)$ \quad The composite $r \circx j$ is $\id_M$ and the composite $j \circx r$ is pseudo naturally homotopic to $\id_M$.
\vspi
[Construct $M$ by repeated use of pointed mapping cylinders (this forces the definitions of $r$ and $j$).]
\label{17.4}
\label{17.5}
\vspi
[Note: \ $\forall \ \bX$, the rule $q \ra \colimx \Omega^n (M\bX)_{n+q}$ defines a spectrum, call it $eM\bX$.]\\
\endgroup %%------------------------------------<<

\index{Conversion Principle (\bPRESPEC)}
\begingroup%%----------------------------------->>
\fontsize{9pt}{11pt}\selectfont
\textbf{\small FACT \  (\un{Conversion Principle})} \ 
Let \bC be a category and let $F$, $G:\bC \ra \bPRESPEC$ be functors.  
Suppose given a pseudo natural transformation 
$\Xi:F \ra G$ $-$then there exists a natural transformation 
$M \Xi:M \circx F \ra M \circx G$ such that the diagram
\begin{tikzcd}[sep=large]
{MFX} \ar{d}[swap]{r} \ar{r}{M\Xi} &{MGX} \ar{d}{r}\\
{FX} \ar{r}[swap]{\Xi} &{GX}
\end{tikzcd}
is pseudo naturally homotopy commutative.\\
\endgroup %%------------------------------------<<

A prespectrum \bX is said to be 
\un{connective}
\index{connective (prespectrum)} 
if $X_1$ is path connected and $X_q$ is $(q-1)$-connected $(q > 1)$.

Example: Given a \dsep \ps $\bGamma$-space $X$, 
the assignment $q \ra S_qX$ specifies a connective prespectrum $\bS X$.

%%----------------------------------------------------------------------------------------------72
\label{14.192}
[The arrow 
$S_0X \ra \Omega S_1 X$ is the identity map 
$\Omega \aX_{\bGamma} \ra \Omega \aX_{\bGamma}$.  For $q > 0$, the arrow 
$S_q X \ra \Omega S_{q+1} X$ is the weak group completion 
$B_1(B^{(q-1)}X)$ $(=\abs{B^{(q-1)}X}_{\bGamma})$ $\ra \Omega \abs{B^{(q)} X}_{\bGamma}$ of 
p. \pageref{14.183}.]

[Note: \  That $\bS X$is conntective is implied by the generalities on 
p. \pageref{14.184}.]

\label{14.190}
\label{18.35}
A prespectrum \bX is said to be an 
\un{$\Omega$-prespectrum}
\index{prespectrum! $\Omega$-prespectrum} 
if $\forall \ q$, the arrow 
$X_q \overset{\sigma_q}{\lra} \Omega X_{q+1}$ is a weak homotopy equivalence.

Example: Given a \dsep \ps $\bGamma$-space $X$, the assignment $q \ra S_qX$ specifies an $\Omega$-prespectrum 
$\bS X$.\\

\begingroup%%----------------------------------->>
\fontsize{9pt}{11pt}\selectfont
\label{17.60}
\index{Algebraic K-Theory (Example) }
\textbf{\small EXAMPLE \ (\un{Algebraic K-Theory})} \ 
Let $A$ be a ring with unit $-$then the prescription 
$q \ra K_0(\Sigma^q A) \times B\bGL(\Sigma^q A)^+$ attaches to $A$ an $\Omega$-prespectrum  $\bW A$.  
Proof: 
$\Omega(K_0(\Sigma^{q+1} A) \times B\bGL(\Sigma^{q+1} A)^+ \approx$ 
$\Omega B\bGL(\Sigma^{q+1} A)^+$ (trivially) $\approx$ 
$K_0(\Sigma^q A) \times B\bGL(\Sigma^q A) ^+$ 
(cf. p. \pageref{14.185} ff.).
\vspi
[Note: \  As it stands, a morphism $A^\prime \ra A\pp$ of rings does not induce a morphism 
$\bW A^\prime \ra \bW A\pp$ of $\Omega$-prespectra (the relevant diagrams are only pointed homotopy commutative).]\\
\endgroup %%------------------------------------<<

\begin{proposition} \ %57
Let 
$
\begin{cases}
\ \bX\\
\ \bY
\end{cases}
$
be connective $\Omega$-prespectra $-$then a morphism 
$\bff:\bX \ra \bY$ is a weak equivalence provided that $f_0:X_0 \ra Y_0$ is a weak homotopy equivalence.\\
\end{proposition}

\textbf{\small LEMMA} \ 
Let 
$
\begin{cases}
\ X\\
\ Y
\end{cases}
$
be homotopy associative H-spaces such that 
$
\begin{cases}
\ \pi_0(X)\\
\ \pi_0(Y)
\end{cases}
$
is a group under the induced product; let $f:X \ra Y$ be a pointed continuous function such that 
$\pi_0(f): \pi_0(X) \ra \pi_0(Y)$ is bijective $-$then $f$ is a weak homotopy equivalence if $f$ is a homology equivalence.

[Consider the commutative diagram
$
\begin{tikzcd}%[sep=small]
{\abs{\sin X}} \ar{d} \ar{r}{\abs{\sin f}} &{\abs{\sin Y}} \ar{d}\\
{X} \ar{r}[swap]{f} &{Y}
\end{tikzcd}
.  \ 
$
Since the hypotheses on 
$
\begin{cases}
\ X\\
\ Y
\end{cases}
$
and $f$ are also satisfied by 
$
\begin{cases}
\ \abs{\sin X}\\
\ \abs{\sin Y}
\end{cases}
$
and $\abs{\sin f}$ and since there are homotopy equivalences 
$
\begin{cases}
\ \abs{\sin X} \ra \abs{\sin X}_0 \times \pi_0(\abs{\sin X})\\
\ \abs{\sin Y} \ra \abs{\sin Y}_0 \times \pi_0(\abs{\sin Y})
\end{cases}
, \ 
$
where 
$
\begin{cases}
\ \abs{\sin X}_0\\
\ \abs{\sin Y}_0\
\end{cases}
$
is the path component of the identity element, Dror's Whitehead theorem implies that $\abs{\sin f}$ is a homotopy equivalence, hence $f$ is a weak homotopy equivalence (Giever-Milnor).]\\

Example: Suppose that $X \ra Y$ is a group completion $-$then $X \ra Y$ is a weak homotopy equivalence if 
$\pi_0(X)$ is a group.

[Note: \  Let $X$ be a \ps $\bGamma$-space such that $\pi_0(X_1)$ is a group.  Because 
$\pi_0(\abs{\sin X_1})$ is likewise a group, the group completion 
$\abs{\sin X_1} \ra \Omega \norm{\sin X}_{\bGamma}$ is a weak homotopy equivalence, thus the same is true of the weak group completion $X_1 \ra \Omega \abs{X}_{\bGamma}$ 
(cf. p. \pageref{14.186}).]\\

%%----------------------------------------------------------------------------------------------73
\begingroup%%----------------------------------->>
\fontsize{9pt}{11pt}\selectfont
\textbf{\small EXAMPLE} \ 
Let $\sO$ be an $\tE_\infty$ operad.  Suppose that $X$ is a \dsep $\sO$-space.  Assume: $\pi_0(X)$ is a group $-$then $X$ has the weak homotopy type of an infinite loop space.
\vspi
[The group completion $X \ra B_0X$ is a weak homotopy equivalence.]\\
\endgroup %%------------------------------------<<

\begin{proposition} \ %58
Let $
\begin{cases}
\ \bX\\
\ \bY
\end{cases}
$
be connective $\Omega$-prespectra $-$then a morphism $\bff:\bX \ra \bY$ is a \we whenever 
$f_0:X_0 \ra Y_0$ induces a bijection $\pi_0(X_0) \ra \pi_0(Y_0)$ and is a homology equivalence.
\end{proposition}

[There is a commutative diagram 
\begin{tikzcd}%[sep=small]
{X_0} \ar{d} \ar{r}{f_0} &{Y_0} \ar{d}\\
{\Omega X_1} \ar{r}[swap]{\Omega f_1}  &{\Omega Y_1}
\end{tikzcd}
and, in view of the lemma, $\Omega f_1$ is a weak homotopy equivalence.  
So, $f_0$ is a weak homotopy equivalence and one can quote Proposition 57.]\\

\begin{proposition} \ %59
Suppose given an infinite loop space machine on $\bGamma$.  
Let 
$
\begin{cases}
\ X\\
\ Y
\end{cases}
$
be \dsp \ps $\bGamma$-spaces, $f:X \ra Y$ a $\bGamma$-map.  
Assume: $f_1:X_1 \ra Y_1$ is a weak homotopy equivalence or $Kf:KX \ra KY$ is a group completion 
$-$then $\bB f:\bB X \ra \bB Y$ is a weak equivalence.
\end{proposition}

[Work with 
\begin{tikzcd}%[sep=small]
{X_1} \ar{d} \ar{r} &{KX} \ar{d} \ar{r}  &{B_0X} \ar{d}\\
{Y_1} \ar{r} &{KY} \ar{r}  &{B_0Y} 
\end{tikzcd}
and apply Proposition 58.]

[Note: \  There is an evident analog of this result for \bS.]\\

\begin{proposition} \ %60
Suppose given an infinite loop space machine on $\bGamma$.  
Let $
\begin{cases}
\ X\\
\ Y
\end{cases}
$
be \dsp \ps $\bGamma$-spaces $-$then the arrow 
$\bB (X \times Y) \ra \bB X \times \bB Y$ is a weak equivalence.
\end{proposition}

[To begin with, the arrow $K(X \times Y) \ra KX \times_k KY$ is a weak homotopy equivalence (examine 
\begin{tikzcd}%[sep=small]
{K(X \times Y)}  \ar{r} &{KX \times_k KY} \\
{L(X \times Y)} \ar{u} \arrow[r,shift right=0.5,dash] \arrow[r,shift right=-0.5,dash] 
&{LX \times_k LY} \ar{u}
\end{tikzcd}
).  This said, form the commutative diagram 
\begin{tikzcd}%[sep=small]
{K(X \times Y)} \ar{d} \ar{r} &{KX \times_k KY} \ar{d}\\
{B_0(X \times Y)} \ar{r} &{B_0X \times_k B_0Y}
\end{tikzcd}
.  
By definition, $K(X \times Y) \ra B_0(X \times Y)$ is a group completion.  
The same is true of $KX \times_k KY \ra B_0X \times_k B_0 Y$.  
Proof: 
$\ov{\pi_0(KX \times_k KY)} \approx$ 
$\ov{\pi_0(KX) \times \pi_0(KY)} \approx$ 
$\ov{\pi_0(KX)} \times \ov{\pi_0(KY)} \approx$ 
$\pi_0(B_0X) \times \pi_0(B_0Y) \approx$
$\pi_0(B_0X \times_k B_0Y)$ 
(cf. p. \pageref{14.187}) 
and, using the K\"unneth formula, 
$H_*(KX \times_k KY;\bk)[\pi_0(KX \times_k KY)^{-1}] \approx$ 
%%----------------------------------------------------------------------------------------------74
$H_*(B_0X \times_k B_0Y;\bk)$ for all prime fields \bk 
(cf. p. \pageref{14.188}).  
It now follows that 
$\pi_0(B_0(X \times Y)) \approx$ 
$\pi_0(B_0X \times_k B_0Y)$ and 
$H_*(B_0(X \times Y)) \approx$ 
$H_*(B_0X \times_k B_0Y)$, from which the assertion (cf. Proposition 58).]\\

Let $X$ be a \dsep \ps $\bGamma$-space $-$then an infinite loop space machine on $\bGamma$ defines a sequence of functors 
$B_q\ov{X}:\bGamma \ra \bDelta$-$\bCG_*$, viz. $\bn \ra B_q\ov{X}_n$.  
It is not claimed that $B_q\ov{X}_n$ is special.  However $B_q\ov{X}_0$ is homotopically trivial and $\forall \ n \geq 1$, the arrow 
$B_q\ov{X}_n \ra B_q\ov{X}_1 \times_k \cdots \times_k B_q\ov{X}_1$ determined by the $\pi_i$ is a weak homotopy equivalence (cf. Propositions 59 and 60).\\

\label{14.193}
\begingroup%%----------------------------------->>
\fontsize{9pt}{11pt}\selectfont
A $\bGamma$-space \mX is said to be 
\un{semispecial}
\index{semispecial ($\bGamma$-space)} 
or 
\un{semiproper} 
\index{semiproper ($\bGamma$-space)} 
if the requirement $X_0 = *$ is relaxed to $X_0$ homotopically trivial, the other conditions on 
$\restr{X}{\bPi}$ staying the same. 
Example: $\forall \ q \geq 0$, $B_q \overline{X}$ is semispecial.\\
\endgroup %%------------------------------------<<

\begingroup%%----------------------------------->>
\fontsize{9pt}{11pt}\selectfont
\textbf{\small LEMMA} \ 
Suppose that $X$ is a \dsep semispecial $\bGamma$-space $-$then there exists a \dsep semiproper semispecial $\bGamma$-space $WX$ and a $\bGamma$-map $\pi:WX \ra X$ such that $\forall \ n$, $\pi_n:W_nX \ra X_n$ is a weak homotopy equivalence.

[Equip $[0,1]$ with the structure of an abelian cofibered monoid in \bCG by writing 
$st = \min\{s,t\}$.  Put $I = \bGamma$-$\ner[0,1]$, so for $\gamma:\bm \ra \bn$, $I\gamma:I_m \ra I_n$ is the function 
$(s_1, \ldots, s_m) \ra (t_1, \ldots, t_n)$, where $t_j = \min\limits_{\gamma(i) = j} \{s_i\}$ 
(a minimum over the empty set is 1).  
Set $W_0X = X_0$ and define a subfunctor $WX$ of $I \times X$ and a $\bGamma$-map $\pi:WX \ra X$ as follows.  
Given an order preserving injection $\gamma:\bm \ra \bn$, let $[0,1]_\gamma^n$ be the subspace of $[0,1]^n$ consisting of those $(t_1, \ldots, t_n)$, such that $t_j = 0$ if $j \in \im \gamma$, $t_j > 0$ if $j \notin \im \gamma$.  
Now form 
$W_n X = \ds\bigcup\limits_\gamma [0,1]_\gamma^n \times (X\gamma) X_m \subset [0,1]^n \times X_n$ :
$X_n$ embeds in $W_nX$ (consider $\gamma = \id_{\bn}$) and the homotopy 
$H((t_1, \ldots, t_n,x),T) = (t_1T, \ldots, t_nT,x)$ $(0 \leq T \leq 1)$ exhibits 
$X_n$ as a strong deformation retract of $W_nX$ (hence $\pi_n(t_1, \ldots, t_n,x) = (0, \ldots, 0,x)$).  
Therefore the \dsep $\bGamma$-space $WX$ is semispecial.  
To establish that $WX$ is semiproper, one has to show that for each injection $\gamma:\bm \ra \bn$, 
$(WX)\gamma:W_m X \ra W_n X$ is a closed $S_\gamma$-cofibration.  
This can be done by observing that 
$\im(WX)\gamma$ admits the description 
$\{(t_1, \ldots, t_n,x): t_j = 1 \ \forall \ j \notin \im \gamma \ \& \ x \in (X\gamma)X_m\}$.]
\vspi
[Note: \  \mW is functorial and $\pi$ is natural: For any $\bGamma$-map $f:X \ra Y$ between \dsep semispecial $\bGamma$-spaces, the diagram
\begin{tikzcd}[sep=large]
{WX} \ar{d} \ar{r}{Wf} &{WY} \ar{d}\\
{X} \ar{r}[swap]{f} &{Y}
\end{tikzcd}
commutes.]\\
\endgroup %%------------------------------------<<

\begingroup%%----------------------------------->>
\fontsize{9pt}{11pt}\selectfont
Observation: The arrow $W_0X \ra W_nX$ corresponding to $\bzero \ra \bn$ is a closed cofibration.  
Put $\ov{W}_nX = W_nX / W_0X$ $-$then $\ov{W}X$ is a \ps $\bGamma$-space, the projection $WX \ra \ov{W}X$ is a levelwise weak homotopy equivalence, and the diagram 
$X \la WX \ra \ov{W}X$ is natural in $X$.\\
\endgroup %%------------------------------------<<

%%----------------------------------------------------------------------------------------------75
\label{16.19}
Notation: If \bX is a prespectrum, then $\Omega \bX$ is the prespectrum specified by 
 $(\Omega \bX)_q =$ $\Omega X_q$, where $\Omega X_q \ra$ $\Omega\Omega X_{q+1}$ is the composite 
\begin{tikzcd}%[sep=small] 
{\Omega X_q} \ar{r}{\Omega \sigma_q} &{\Omega\Omega X_q} \ar{r}{\Tee} &{\Omega\Omega X_q},
\end{tikzcd}
$\Tee$ being the twist $(\Tee f)(s)(t) = f(t)(s)$.\\

\label{14.189}
\begingroup%%----------------------------------->>
\fontsize{9pt}{11pt}\selectfont
\label{16.3}%dmc mnft
\textbf{\small EXAMPLE} \ 
Let $X$ be a \dsp \ps $\bGamma$-space.  Assume: $X_1$ is path connected $-$then $\forall \ q$, 
$\abs{B^{(q)}X}_{\bGamma}$ is $(q+1)$-connected, hence $\Omega S X$ is a connective $\Omega$-prespectrum.\\
\endgroup %%------------------------------------<<

\textbf{\small LEMMA} \ 
For any \ps $\bGamma$-space $X$, $\Omega X$ is a \ps $\bGamma$-space and there is a canonical arrow 
$\abs{\Omega X}_{\bGamma} \overset{\gamma}{\longrightarrow} \Omega\aX_{\bGamma}$.

[Note: \  Here, of course, $\Omega X$ takes \bn to $\Omega X_n$.]\\

\begin{proposition} \ %61
Let $X$ be a \dsp \ps $\bGamma$-space $-$then there is a morphism 
$\bss:\bS\Omega X \ \ra \Omega \bS X$ in \bPRESPEC such that the triangle
%\begin{tikzcd}%[sep=small]
%&{\Omega X_1} \ar{ld} \ar{rd}\\
%{S \Omega X} \ar{rr}[swap]{s_0}  &&{\Omega S_0 X}
%\end{tikzcd}
\hspace{-.5cm}
\begin{tikzcd}[sep=small]
&{\Omega X_1} \ar{lddd} \ar{rddd}\\
\\
\\
{S_0 \Omega X} \ar{rr}[swap]{s_0}  &&{\Omega S_0 X}
\end{tikzcd}
commutes.
\end{proposition}

[Explicated, the oblique arrow on the left is 
\begin{tikzcd}[sep=small]
&{\Omega X_1} \ar{ld}\\
{\Omega \abs{\Omega X}_{\bGamma}}
\end{tikzcd}
and the composite 
$\Omega X_1 \ra$ 
$\Omega \abs{\Omega X}_{\bGamma} \overset{\Omega_\gamma}{\longrightarrow}$ 
$\Omega\Omega \aX_{\bGamma} \overset{\Tee}{\longrightarrow}$ 
$\Omega\Omega \aX_{\bGamma}$
is $\Omega$ of $X_1 \ra  \Omega\aX_{\bGamma}$,  the oblique arrow 
\begin{tikzcd}[sep=small]
{\Omega X_1} \ar{rd}\\
&{\Omega \Omega \aX_{\bGamma}}
\end{tikzcd}
on the right.    Definition: $s_0 = \Tee \circx \Omega \gamma$.  To force compatibility, take
$s_1 = \gamma: S_1 \Omega X \ra \Omega S_1 X$, 
thereby ensuring that the diagram 
\begin{tikzcd}%[sep=small]
{S_0 \Omega X} \arrow[d,shift right=0.5,dash] \arrow[d,shift right=-0.5,dash] \ar{r}{s_0} 
&{\Omega S_0 X} \ar{d}{\Tee} \\
{\Omega S_1 \Omega X} \ar{r}[swap]{\Omega s_1} &{\Omega \Omega S_1 X}
\end{tikzcd}
commutes.  The arrows
$B_n \Omega X = \abs{\Omega \ov{X}_n}_{\bGamma}$ 
$\overset{\gamma}{\longrightarrow}$ 
$\Omega \abs{\ov{X}_n}_{\bGamma} =$ 
$\Omega B_n X$
yield a $\bGamma$-map $b:B\Omega X \ra \Omega BX$.  Setting $b^{(0)} = \id_{\Omega X}$, let $b^{(q)}$ $(q > 0)$ 
be the composite 
%\begin{tikzcd}%[sep=small]
%{B^{(q)}\Omega X} \ar{r}{B{b^{(q-1)}}} &{B\Omega B^{(q-1)}X} \ar{r}{b} &{\Omega B^{(q)} X}.
%\end{tikzcd}
${B^{(q)}\Omega X}$
$\overset{B{b^{(q-1)}}}{\xrightarrow{\hspace*{1.25cm}}}$ 
$B\Omega B^{(q-1)}X$ 
$\overset{b}{\xrightarrow{\hspace*{.75cm}}}$ 
$\Omega B^{(q)} X$.
Definition: $s_q = \gamma \circx \abs{b^{(q-1)}}_{\bGamma}$ $(q > 1)$.  This makes sense:
$S_q \Omega X = $
\begin{tikzcd}%[sep=small]
{\abs{B^{(q-1)}\Omega X}_{\bGamma}} \ar{r}{\abs{b^{(q-1)}}_{\bGamma}} &{\abs{\Omega B^{(q-1)}X_{\bGamma}}}
\end{tikzcd}
$\overset{\gamma}{\longrightarrow} \Omega \abs{B^{(q-1)}X}_{\bGamma} =$ 
$\Omega S_q X$ 
and the diagram 
\begin{tikzcd}%[sep=small]
{S_q \Omega X} \ar{d}\ar{r}{s_q} &{\Omega S_q X} \ar{d}{\Tee} \\
{\Omega S_{q+1} \Omega X} \ar{r}[swap]{\Omega s_{q+1}} &{\Omega S_{q+1} \Omega X}
\end{tikzcd}
commutes.]

[Note: \  If $X_1$ is path connected, then $\Omega S X$ is a connective $\Omega$-prespectrum 
(cf. p. \pageref{14.189}) 
and $s_0$ is a weak homotopy equivalence 
(cf. p. \pageref{14.190}), 
thus $\bss$ is a weak equivalence (cf. Proposition 57).  
It is also clear that $\bss$ is natural.]\\

%%----------------------------------------------------------------------------------------------76
\begingroup%%----------------------------------->>
\fontsize{9pt}{11pt}\selectfont
\label{14.194}
\textbf{\small LEMMA} \ 
Suppose that $X$ is a \dsep semispecial $\bGamma$-space $-$then there exists a $\bGamma$-map 
$\omega:W \Omega X \ra \Omega W X$ such that the triangle 
\begin{tikzcd}[sep=small]
{W \Omega X} \ar{dddr}[swap]{\pi} \ar{rr}{\omega} &&{\Omega W X} \ar{dddl}{\Omega \pi} \\
\\
\\
&{\Omega X}
\end{tikzcd}
is homotopy commutative, thus $\forall \ n$, $\omega_n:W_n\Omega X \ra \Omega W_n X$ is a weak homotopy equivalence.

[Represent a typical element in 
$W_n \Omega X$ by $(t_1, \ldots, t_n, \sigma)$ $(\sigma \in \Omega_n X = \Omega X_n)$ and let
\[
\omega_n(t_1, \ldots, t_n,\sigma)(t) = 
\begin{cases}
\ (u_1(t), \ldots, u_n(t),\sigma(0)) \hspace{1cm} (0 \leq t \leq 1/3)\\
\ (t_1, \ldots, t_n,\sigma(3t - 1)) \hspace{1.25cm} (1/3 \leq t \leq 2/3),\\
\ (v_1(t), \ldots, v_n(t),\sigma(1)) \hspace{1.11cm} (2/3 \leq t \leq 1)
\end{cases}
\]
where $u_j(t) = 1 - 3t + 3tt_j$, $v_j(t) = 3t - 2 + (3 - 3t)t_j$ $(1 \leq j \leq n)$.  The prescription 
\[
H_\omega((t_1, \ldots, t_n,\sigma),T)(t) = 
\begin{cases}
\ \sigma(0) \hspace{2.22cm} (0 \leq t \leq (1/3)T)\\
\ \sigma\left(\ds\frac{3t - T}{3 - 2T} \right)  \hspace{1cm} ((1/3)T \leq t \leq 1 - (1/3)T)\\
\ \sigma(1) \hspace{2.2cm} (1 - (1/3)T \leq t \leq 1)
\end{cases}
\]
is a homotopy between $\pi$ and $\Omega \pi \circx \omega$.]
\vspi
[Note: \  $\omega$ and $H_\omega$ are natural.]
\\
\endgroup %%------------------------------------<<

\begingroup%%----------------------------------->>
\fontsize{9pt}{11pt}\selectfont
Observation: The diagram
\begin{tikzcd}%[sep=small]
{W \Omega X} \ar{d} \ar{r}{\omega} &{\Omega W X} \ar{d} \\
{\ov{W} \Omega X}    \ar{r}[swap]{\ov{\omega}} &{\Omega \ov{W} X}
\end{tikzcd}
commutes and $\ov{\omega}$ is a levelwise weak homotopy equivalence.\\
\endgroup %%------------------------------------<<

A 
\un{biprespectrum} 
\index{biprespectrum} 
\bX is a sequence of prespectra $\bX_q$ and morphisms 
$\bX_q \overset{\sigma_q}{\longrightarrow} \Omega \bX_{q+1}$ $(q \geq 0)$.  
Spelled out, a biprespectrum is a doubly indexed sequence of pointed \dsp compactly generated speces $X_{q,p}$ and 
pointed continuous functions 
$\sigma_{q,p}:X_{q,p} \ra \Omega X_{q+1,p}$, 
$\sigma_{q,p}:X_{q,p} \ra \Omega X_{q,p+1}$
such that the diagram 
\[
\begin{tikzcd}[sep=huge]
{X_{q,p}} \ar{d}[swap]{\sigma_{q,p}} \ar{rr}{\sigma_{q,p}} &&{\Omega X_{q,p+1}} \ar{d}{\Omega \sigma_{q,p+1}} \\
{\Omega X_{q+1,p}} \ar{r}[swap]{\Omega \sigma_{q+1,p}}
&{\Omega\Omega X_{q+1,p+1}} \ar{r}[swap]{\Tee}
&{\Omega\Omega X_{q+1,p+1}} 
\end{tikzcd}
\]
commutes $\forall \ q, p$.  \bBIPRESPEC 
\index{\bBIPRESPEC} 
is the category whose objects are the biprespectra and whose morphisms 
$\bff:\bX \ra \bY$ are the doubly indexed sequences of pointed continuous functions 
$f_{q,p}:X_{q,p} \ra Y_{q,p}$ such that $f_{q,*}$ $\&$ $f_{*,p}$ are morphisms of prespectra $\forall \ q, p$.\\

\index{Theorem: The Up and Across Theorem}
\textbf{\small THE UP AND ACROSS THEOREM} \quad
Let \bX be a biprespectrum.  
Assume: $\forall \ q$, $\sigma_q$ is a weak equivalence and $\bX_q$ is an $\Omega$-prespectrum 
$-$then the $\Omega$-prespectra 
$
\begin{cases}
\ X_{0,*}\\
\ X_{*,0}
\end{cases}
$ 
are naturally weakly equivalent.

%%----------------------------------------------------------------------------------------------77
[Let \bC be the full subcategory of \bBIPRESPEC whose objects \bX have the property that $\forall \ q$, $\sigma_q$ 
is a weak equivalence and $\bX_q$ is an $\Omega$-prespectrum.  Denote by 
$
\begin{cases}
\ E^\prime\\
\ E\pp
\end{cases}
$ 
the functor $\bC \ra \bPRESPEC$ that sends \bX to 
$
\begin{cases}
\ X_{0,*}\\
\ X_{*,0}
\end{cases}
$ 
$-$then the claim is that
$
\begin{cases}
\ E^\prime \bX\\
\ E\pp \bX
\end{cases}
$ 
are naturally weakly equivalent.  For this, it suffices to construct functors 
$D^\prime$, $D\pp:\bC \ra \bPRESPEC$ and a pseudo natural weak equivalence 
$\Xi_{\bX}:D^\prime \bX \ra D\pp \bX$ 
together with natural weak equivalences
$
\begin{cases}
\ e_{\bX}^\prime:E^\prime \bX \ra D^\prime \bX\\
\ e_{\bX}\pp:E\pp \bX \ra D\pp \bX
\end{cases}
$
.  
Reason:  Consider the diagram
\[
\begin{tikzcd}[sep=large]
&{M D^\prime \bX} \ar{d}[swap]{r} \ar{r}{M\Xi}  &{M D\pp \bX}  \ar{d}{r}\\
{E^\prime \bX} \ar{r} &{D^\prime \bX} \ar{r}[swap]{\Xi}  &{D\pp \bX} &{E\pp \bX} \ar{l}
\end{tikzcd}
\]
furnished by the conversion principle.  
Definition: $D^\prime \bX = \Omega^q X_{q,q} = D\pp \bX$, the arrows of structure
$\sigma_q^\prime:D_q^\prime \bX \ra \Omega D_{q+1}^\prime \bX$, 
$\sigma_q\pp:D_q\pp \bX \ra \Omega D_{q+1}\pp \bX$
being the composites
$\Omega^{q}X_{q,q}$
$\overset{\Omega^{q}\sigma_{q,q}}{\xrightarrow{\hspace*{1cm}}}$ 
$\Omega^{q+1}X_{q,q+1}$
$\overset{\Tee_q}{\xrightarrow{\hspace*{1cm}}}$ 
$\Omega^{q+1}X_{q,q+1}$
$\overset{\Omega^{q+1}\sigma_{q,q+1}}{\xrightarrow{\hspace*{1.5cm}}}$ 
${\Omega^{q+2}X_{q+1,q+1}}$, \ 
$\Omega^{q}X_{q,q}$
$\overset{\Omega^{q}\sigma_{q,q}}{\xrightarrow{\hspace*{1cm}}}$ 
$\Omega^{q+1}X_{q+1,q}$
$\overset{\Tee_q}{\xrightarrow{\hspace*{1cm}}}$ 
$\Omega^{q+1}X_{q+1,q}$
$\overset{\Omega^{q+1}\sigma_{q+1,q}}{\xrightarrow{\hspace*{1.5cm}}}$ 
${\Omega^{q+2}X_{q+1,q+1}}$,
and 
$
\begin{cases}
\ D_q^\prime \bff\\
\ D_q\pp \bff
\end{cases}
= \ \Omega^q f_{q,q}
$ 
, 
where $\bff:\bX \ra \bY$ $(f_{q,q}:X_{q,q} \ra Y_{q,q})$.  
Here $\Tee_q$ is given by twisting the last coordinate past the first q coordinates:
$(\Tee_q(f)(s)(t) = f(t)(s)$ $(s \in \bS^q, \ t \in \bS^1)$.
If $\Xi_{\bX,q}:\Omega^q X_{q,q} \ra \Omega^qX_{q,q}$ 
is the identity for even $q$ and the negative of the identity for odd $q$ 
(i.e., reverse the first coordinate), then there are pointed homotopies $H_{\bX,q}$ between 
$\Omega \Xi_{\bX,q+1} \circx \sigma_q^\prime$ and $\sigma_q\pp \circx \Xi_{\bX,q}$.  
Since the data is natural in \bX, 
$\Xi_{\bX}:D^\prime \bX \ra D\pp \bX$ 
is a pseudo natural weak equivalence.  
Introduce weak homotopy equivalences 
$e_{q,p}^\prime:X_{q,p} \ra \Omega^{p - q}X_{p,p}$ taking $e_{q,q}^\prime = \id$ and inductively letting 
$e_{q,p}^\prime$ $(q < p)$ be the composite 
$X_{q,p}$
$\overset{\sigma_{q,p}}{\xrightarrow{\hspace*{1cm}}}$ 
$\Omega X_{q+1,p}$
$\overset{\Omega e_{q+1,p}^\prime}{\xrightarrow{\hspace*{1cm}}}$ 
$\Omega^{p-q}X_{p,p}$
.  
Call $\omega_{q,p}$ the composite
%%%%
$\Omega^{p-q}X_{p,p}$
$\overset{\Omega^{p-q} \sigma_{p,p}}{\xrightarrow{\hspace*{1.5cm}}}$ 
$\Omega^{p+1-q} X_{p,p+1}$
$\overset{\Tee_{p-q}}{\xrightarrow{\hspace*{1cm}}}$ 
$\Omega^{p+1-q}X_{p,p+1}$
$\overset{\Omega^{p+1-q} \sigma_{p,p+1}}{\xrightarrow{\hspace*{1.75cm}}}$ 
$\Omega^{p+2-q}X_{p+1,p+1}$
$-$then for each $q$, the $e_{q,p}^\prime$ $(q \leq p)$ specify a morphism 
$\{X_{q,p} \overset{\sigma_{q,p}}{\lra} \Omega X_{q,p+1}\} \ra$ 
$\{\Omega^{p-q} X_{p,p} \overset{\omega_{q,p}}{\lra} \Omega^{p+2-q} X_{p+1,p+1}\}$
of prespectra (use induction on $p - q$) (note the shift in the indexing).  
Put $e_{\bX}^\prime = e_{0,*}^\prime$ and define $e_{\bX}\pp$ analogously.]\\

\index{Theorem: Comparison Theorem (infinite loop space)}
\textbf{\small COMPARISON THEOREM} \quad
Suppose given an infinite loop space machine on $\bGamma$ $-$then $\forall$ \dsep \ps $\bGamma$-space $X$, $\bB X$ is naturally weakly equivalent to $\bS X$.

[Note: \ \bS is a functor from the category of \dsep \ps $\bGamma$-spaces to the full subcategory of \bPRESPEC whose objects are the connective $\Omega$-prespectra while \bB is a functor from the category of \dsep \ps $\bGamma$-spaces to the full subcategory of \bPRESPEC whose objects are the connective spectra.  
It is therefore of interest to
%%----------------------------------------------------------------------------------------------78
observe that the proof goes through unchanged if the definition of infinite loop space machine is weakened: It suffices that \bB takes values in the category of connective $\Omega$-prespectra.]\\

Application: Let $\sO$ be an $\tE_\infty$ operad.  Suppose given an infinite loop space machine on $\Oh$ (e.g., the May machine) $-$then $\forall$ \dsep \ps $\bGamma$-space $X$, $\bB(e^*X)$ $(=\bB(X \circx e))$ is naturally equivalent to $\bS X$.\\

\begingroup%%----------------------------------->>
\fontsize{9pt}{11pt}\selectfont
\textbf{\small FACT} \ 
Let $\sO$ be an $\tE_\infty$ operad.  Suppose given an infinite loop space machine on $\Oh$ $-$then $\forall$ \dsep \ps 
$\Oh$-space $X$, $\bB X$, and $\bS(\epsilon_*X)$ are naturally weakly equivalent.
\vspi
[Recalling that 
$\epsilon_*:\pOhp \ra \psg$ respects the $\Delta$-separation condition 
(cf. p. \pageref{14.191}), 
$\bB X$ is naturally weakly equivalent to 
$\bB UX$ or still, is naturally weakly equivalent to $\bB(\epsilon^*\epsilon_*X)$ which is naturally weakly equivalent to  
$\bS(\epsilon_*X)$.]\\
\endgroup %%------------------------------------<<

Heuristics: The proof of the comparison theorem is complicated by a technicality:  
The $B_q\ov{X}$are not necessarily 
\dsep \ps $\bGamma$-spaces (but are \dsep semispecial $\bGamma$-spaces).  
However, let us proceed as if they were 
$-$then one can form the connective $\Omega$-prespectra $\bS B_q \ov{X}$ and there are morphisms 
\begin{tikzcd}%[sep=small]
{\sigma_q:\bS B_q \ov{X}} \ar{r}{\bS \sigma_{q}} 
&{\bS \Omega B_{q+1} \ov{X}} \ar{r}{\bss} &{}
%&{\Omega \bS B_{q+1} \ov{X}}
\end{tikzcd}
%\begin{tikzcd}%[sep=small]
$\Omega \bS B_{q+1} \ov{X}$
%\end{tikzcd}
(cf. Proposition 61).  Since $\forall \ q$, $\sigma_{q}$ is a weak equivalence, 
it follows from the up and across theorem that the connective 
$\Omega$-prespectra 
$\bS B_0 \ov{X}$ $(= \{S_q B_0\ov{X}\})$,
$S_0 \bB \ov{X}$ $(= \{S_0 B_q\ov{X}\})$
are naturally weakly equivalent.  The idea now is to show that $\bS X$ is naturally weakly equivalent to 
$\bS B_0 \ov{X}$ and $\bB X$ is naturally weakly equivalent to $S_0 \bB \ov{X}$.\\
\indent\indent $(\bS B_0 \ov{X})$ \quad $\forall \ n$, there are arrows
$L\ov{X}_n \ra K\ov{X}_n$, 
$K\ov{X}_n \ra B_0\ov{X}_n$,
i.e., there are $\bGamma$-maps 
$L\ov{X} \ra K\ov{X}$,
$K\ov{X} \ra B_0\ov{X}$.  
Because $L\ov{X}_1 \ra K\ov{X}_1$ is a weak homotopy equivalence and 
$K\ov{X}_1 \ra B_0\ov{X}_1$ is a group completion, the arrow 
$\bS L\ov{X} \ra \bS K\ov{X}$ is a weak equivalence, as is the arrow 
$\bS K\ov{X} \ra \bS B_0\ov{X}$ (cf. Proposition 59).  
But $L\ov{X} = X$.\\
\indent\indent $(S_0 \bB \ov{X})$ \quad The weak group completions 
$B_qX = B_q\ov{X}_1 \ra \Omega \abs{B_q \ov{X}}_{\bGamma} = S_0B_q\ov{X}$ 
define a morphism 
$\bB X \ra S_0\bB X$ of connective $\Omega$-prespectra (cf. Proposition 61) which we claim is a weak equivalence.  
In fact, $\pi_0(B_0X)$ is a group, thus 
$B_0X \ra S_0 B_0 \ov{X}$ 
is a weak homotopy equivalence 
(cf. p. \pageref{14.192}), 
so Proposition 57 is applicable.\\

\begingroup%%----------------------------------->>
\fontsize{9pt}{11pt}\selectfont
To establish the comparison theorem in full generality, one first has to extend the basic definitions from the context of proper special $\bGamma$-spaces to that of semiproper semispecial $\bGamma$-spaces.  Thus let $X$ be a semiproper semispecial $\bGamma$-space $-$then there is a closed cofibration $X_0 \ra \abs{X}_{\bGamma}$ and it is best to work with the quotient $\aX_{\ov{\bGamma}} \equiv \aX_{\bGamma} /  X_0$.
Again one has a canonical arrow 
$\Sigma X_1 \ra \aX_{\ov{\bGamma}}$ 
whose adjoint 
$X_1 \ra \Omega \aX_{\ov{\bGamma}}$
is a weak group completion.  
It still makes sense to form $\ov{X}$ and the classifying space $BX$  of $X$ takes \bn to 
$B_n X = \abs{\ov{X}_n}_{\ov{\bGamma}}$.  
The definition of $B^{(q)}X$ is as before but 
$S_0X = \Omega \aX_{\ov{\bGamma}}$, 
$S_{q+1} X = \abs{B^{(q)}X}_{\ov{\bGamma}}$ $(q \geq 0)$.
\vspi
%%----------------------------------------------------------------------------------------------79
Turning to the proof of the comparison theorem, let $X$ be a \dsep \ps $\bGamma$-space $-$then $\forall \ q$, 
$WB_q\ov{X}$ is a \dsep semiproper semispecial $\bGamma$-space 
(cf. p. \pageref{14.193}), 
$\bS WB_q\ov{X}$ 
is a connective $\Omega$-prespectrum, and there are morphisms
$\bsigma_q:\bS WB_q\ov{X}$
$\overset{\bS W\sigma_q}{\xrightarrow{\hspace*{1cm}}}$ 
$\bS W \Omega B_{q+1} \ov{X}$
$\overset{\bS \omega}{\xrightarrow{\hspace*{1cm}}}$ 
$\bS \Omega W B_{q+1} \ov{X}$
$\overset{\bss}{\xrightarrow{\hspace*{1cm}}}$ 
$\Omega \bS W B_{q+1} \ov{X}$
(cf. Proposition 61) ($\omega$ as in the lemma on 
p. \pageref{14.194}).  
Since $\forall \ q$, $\bsigma_q$ is a weak equivalence, 
it follows from the up and across theorem that the connective $\Omega$-prespectra 
$\bS WB_0\ov{X}$ $(=\{S_qWB_0\ov{X}\})$, 
$S_0W\bB\ov{X}$ $(= \{S_0WB_q\ov{X}\})$ 
are naturally weakly equivalent.  
The idea now is to show that $\bS X$ is naturally weakly equivalent to 
$\bS WB_0\ov{X}$ and $\bB X$ is naturally weakly equivalent to $S_0W\bB\ov{X}$.
\\
\indent\indent $(\bS WB_0\ov{X})$ \quad There is a natural weak equivalence $\bS WX \ra \bS X$.  
On the other hand, there are natural weak equivalences
$\bS WL\ov{X} \ra \bS WK\ov{X}$, 
$\bS WK\ov{X} \ra \bS WB_0\ov{X}$
and 
$L\ov{X} = X$.\\
\indent\indent $(S_0 W\bB\ov{X})$ \quad Let $\bW\bB X$ be the connective $\Omega$-prespectrum specified by 
$q \ra W_1B_q\ov{X}$ and 
%\begin{tikzcd}%[sep=small]
%{W_1B_q\ov{X}} \ar{r}{W_1\sigma_q} 
%&{W_1 \Omega B_{q+1} \ov{X}} \ar{r}{\omega_1} 
%&{\Omega W_1 B_{q+1} \ov{X}}
%\end{tikzcd}
$W_1B_q\ov{X}$
$\overset{W_1\sigma_q}{\xrightarrow{\hspace*{1cm}}}$ 
$W_1 \Omega B_{q+1} \ov{X}$
$\overset{\omega_1}{\xrightarrow{\hspace*{1cm}}}$ 
$\Omega W_1 B_{q+1} \ov{X}$
$-$then there is a natural weak equivalence 
$\bW\bB X \ra S_0 W\bB\ov{X}$.  But there is also a pseudo natural weak equivalence $\bW\bB X \ra \bB X$, hence 
$\bB X$ is naturally weakly equivalent to $\bW\bB X$  (conversion principle).\\
\endgroup %%------------------------------------<<

\textbf{\small LEMMA} \ 
Let $X$ be a \dsep \ps $\bGamma$-space $-$then $\Sigma X_1$ is homeomorphic to $(\aX_{\bGamma})_1$, thus the arrow 
$X_1 \ra \Omega \aX_{\bGamma}$ is a closed embedding.\\

Application: Let $X$ be a \dsep \ps $\bGamma$-space $-$then $\forall \ q$, the arrow 
$S_q X \ra \Omega S_{q+1} X$ is a closed embedding.\\

Consequently, if $X$ is a \dsep \ps $\bGamma$-space $-$then the rule $q \ra \colimx \Omega^n S_{n+q}X$ defines a spectrum, call it $\be\bS X$.\\

\begin{proposition} \ %62
Suppose given an infinite loop space machine on $\bGamma$ 
$-$then $\forall$ $\Delta$-separated proper special $\bGamma$-space $X$, $\bB X$ is naturally weakly equivalent to $\beS X$.
\end{proposition}

[There is an obvious natural weak equivalence 
$\bS X \ra \be\bS X$, so the assertion follows from the comparison theorem.]\\

Remark: It is a fact that \bSPEC carries a model category structure in which the weak equivalences are the levelwise weak homotopy equivalences (cf $\S 15$, Proposition 8).  One can therefore interpret Proposition 62 as saying that 
$\bB X$ and $\be\bS X$ are isomorphic in \bHSPEC (a.k.a ``the'' 
\un{stable homotopy category} 
\index{stable homotopy category}).\\

%%%%%%%%%%%%%%%%%%%%%%%%%%%%%%%%%%%%%%
%%%%%%%%%%%%%%%%%%%%%%%%%%%%%%%%%%%%%%
%%%%%%%%%%%%%%%%%%%%%%%%%%%%%%%%%%%%%%

\begin{center}
$\S \ 14$
\\[0.5cm]
$\mathcal{REFERENCES}$\\
\end{center}

\[
\textbf{BOOKS}
\]

\begingroup
\fontsize{9pt}{11pt}\selectfont
\setlength\parindent{0 cm}

[1] \quad Adams, J., 
\textit{Infinite Loop Spaces}, Princeton University Press (1978).
\\[-.2cm]

[2] \quad Boardman, J. and Vogt, R., 
\textit{Homotopy Invariant Algebraic Structures on Topological Spaces}, 

\hspace{0.8cm}Springer Verlag (1973).
\\[-.2cm]

[3] \quad Caruso, J., 
\textit{Configuration Spaces and Mapping Spaces}, Ph.D. Thesis, University of Chicago, Chicago 

\hspace{0.8cm}(1979).
\\[-.2cm]

[4] \quad May, J., 
\textit{The Geometry of Iterated Loop Spaces}, Springer Verlag (1972).
\\[-.2cm]

[5] \quad Milgram, R., 
\textit{Unstable Homotopy from the Stable Point of View}, Springer Verlag (1974).
\\[-.2cm]

[6] \quad Rezk, C., 
\textit{Operads and Spaces of Algebra Structures}, Ph.D. Thesis, MIT, Cambridge (1996).
\\[-.2cm]
\endgroup

\[
\textbf{ARTICLES}
\]

\begingroup
\fontsize{9pt}{11pt}\selectfont
\setlength\parindent{0 cm}

[1] \quad Almgren, F., The Homotopy Groups of the Integral Cycle Groups, 
\textit{Topology} \textbf{1} (1962), 257-299.
\\[-.2cm]

[2] \quad Baues, H., Geometry of Loop Spaces and the Cobar Construction, 
\textit{Memoirs Amer. Math. Soc.} \textbf{230} 

\hspace{0.8cm}(1980), 1-171.
\\[-.2cm]

[3] \quad Berger, C., Op\'erades Cellulaires et Espaces de Lacets It\'er\'es, 
\textit{Ann. Inst. Fourier (Grenoble)} \textbf{46} (1996), 

\hspace{0.8cm}1125-1157.
\\[-.2cm]

[4] \quad B\"odigheimer, C., Stable Splittings of Mapping Spaces, 
\textit{SLN} \textbf{1286} (1987), 174-187.
\\[-.2cm]

[5] \quad Brown, E. and Szczarba, R., Continuous Cohomology and Real Homotopy Type, 
\textit{Trans. Amer. Math.}

\hspace{0.8cm}\textit{Soc.} \textbf{311} (1989), 57-106.
\\[-.2cm]

[6] \quad Carlsson, G. and Cohen, R., The Cyclic Groups and the Free Loop Space, 
\textit{Comment. Math. Helv.} \textbf{62} 

\hspace{0.8cm}(1987), 423-449.
\\[-.2cm]

[7] \quad Carlsson, G. and Milgram, R., Stable Homotopy and Iterated Loop Spaces, In: 
\textit{Handbook of Algebraic} 

\hspace{0.8cm}\textit{Topology}, I. James (ed.), North Holland (1995), 505-583.
\\[-.2cm]

[8] \quad Caruso, J., Cohen, F., May, J., and Taylor, L., James Maps, Segal Maps, and the Kahn-Priddy Theo-

\hspace{0.8cm}rem, 
\textit{Trans. Amer. Math. Soc.} \textbf{281} (1984), 243-283.
\\[-.2cm]

[9] \quad Cohen, F., May, J., and Taylor, L., James Maps, and the $E_n$ Ring Spaces, 
\textit{Trans. Amer. Math. Soc.} 

\hspace{0.8cm}\textbf{281} (1984), 285-295.
\\[-.2cm]

[10] \quad Dold, A. and Thom, R., Quasifaserungen und Unendliche Symmetrische Produkte, 
\textit{Ann. of Math.} 

\hspace{0.95cm}\textbf{67} (1958), 239-281.
\\[-.2cm]

[11] \quad Dunn, G., Uniqueness of $n$-Fold Delooping Machines, 
\textit{J. Pure Appl. Algebra} \textbf{113} (1996), 159-193.
\\[-.2cm]

[12] \quad Floyd, E. and Floyd, W., Actions of the Classical Small Categories of Topology, 
\textit{Preprint}.
\\[-.2cm]

[13] \quad Gajer, P., The Intersection Dold-Thom Theorem, 
\textit{Topology} \textbf{35} (1996), 939-967.
\\[-.2cm]

[14] \quad Hollender, J. and Vogt, R., Modules of Topological Spaces, Applications to Homotopy Limits and 

\hspace{0.95cm}$E_\infty$ Structures, 
\textit{Arch. Math.} \textbf{59} (1992), 115-129.
\\[-.2cm]

[15] \quad James, I., Reduced Product Spaces, 
\textit{Ann. of Math.} \textbf{62} (1955), 170-197.
\\[-.2cm]

[16] \quad Kriz, I. and May, J., Operads, Algebras, Modules, and Motives, 
\textit{Ast\'erisque} \textbf{233} (1995), 1-145.
\\[-.2cm]

[17] \quad May, J., Classifying Spaces and Fibrations, 
\textit{Memoirs Amer. Math. Soc.} \textbf{155} (1975), 1-98.
\\[-.2cm]

[18] \quad May, J., Infinite Loop Space Theory, 
\textit{Bull. Amer. Math. Soc.} \textbf{83} (1977), 456-494.
\\[-.2cm]

[19] \quad May, J., Infinite Loop Space Theory Revisited , 
\textit{SLN} \textbf{741} (1979), 625-642.
\\[-.2cm]

[20] \quad May, J., Applications and Generalizations of the Approximation Theorem, 
\textit{SLN} \textbf{763} (1979), 38-69.
\\[-.2cm]

[21] \quad May, J., Multiplicative Infinite Loop Space Theory, 
\textit{J. Pure Appl. Algebra} \textbf{26} (1982), 1-69.
\\[-.2cm]

[22] \quad May, J. and Thomason, R., The Uniqueness of Infinite Loop Space Machines, 
\textit{Topology} \textbf{17} (1978), 

\hspace{0.95cm}205-224.
\\[-.2cm]

[23] \quad McCord, M., Classifying Spaces and Infinite Symmetric Products, 
\textit{Trans. Amer. Math. Soc.} \textbf{146} 

\hspace{0.95cm}(1969), 273-298.
\\[-.2cm]

[24] \quad Meyer, J-P., Bar and Cobar Constructions I and II, 
\textit{J. Pure Appl. Algebra} \textbf{33} (1984), 163-207 and 

\hspace{0.95cm}\textbf{43} (1986), 179-210.
\\[-.2cm]

[25] \quad Mostow, M., Continuous Cohomology of Spaces with Two Topologies, 
\textit{Memoirs Amer. Math. Soc.} 

\hspace{0.95cm}\textbf{175} (1976), 1-142.
\\[-.2cm]

[26] \quad Segal, G., Categories and Cohomology Theories, 
\textit{Topology} \textbf{13} (1974), 293-312.
\\[-.2cm]

[27] \quad Spanier, E., Infinite Symmetric Products, Function Spaces and Duality, 
\textit{Ann. of Math.} \textbf{69} (1959), 

\hspace{0.95cm}142-198.
\\[-.2cm]

[28] \quad Thomason, R., Uniqueness of Delooping Machines, 
\textit{Duke Math. J.} \textbf{46} (1979), 217-252.
\\[-.2cm]

\setlength\parindent{2em}

\endgroup

\chapter{
$\boldsymbol{\S}$\textbf{15}.\quadx  TRIANGULATED CATEGORIES}
\setlength\parindent{2em}
\setcounter{proposition}{0}
\setcounter{chapter}{15}

%%----------------------------------------------------------------------------------------------01
$\text{ }$\\[-1.25cm]

Because the theory of triangulated categories lies outside the usual categorical experience, an exposition of the basics seems to be in order.  
Topologically, the rationale is that the stable homotopy category is triangulated.

Let \bC be an additive category $-$then an additive functor $\Sigma:\bC \ra \bC$ is said to be a 
\un{suspension functor} 
\index{suspension functor} 
if it is an equivalence of categories.

[Note: \ Thus there is also a functor $\Omega:\bC \ra \bC$ which is simultaneously a right and left adjoint for $\Sigma$ 
and the four arrows of adjunction 
$\Sigma \circx \Omega \overset{\nu}{\ra} \id_{\bC}$, 
$\id_{\bC}  \overset{\mu}{\ra} \Omega \circx \Sigma$, 
$\Omega \circx \Sigma \overset{\mu^{-1}}{\ra} \id_{\bC}$, 
$\id_{\bC} \overset{\nu^{-1}}{\ra} \Sigma \circx \Omega$ 
are natural isomorphisms.]

Let \bC be an additive category, $\Sigma$ a suspension functor $-$then a 
\un{triangle}
\index{triangle (in an additive category)} 
in \bC consists of objects $X,Y,Z$ and morphisms 
$u,v,w$, where 
$X \overset{u}{\ra} Y$, 
$Y \overset{v}{\ra} Z$, 
$Z \overset{w}{\ra} \Sigma X$, a 
\un{morphism of} \un{triangles} 
\index{morphism of triangles (in an additive category)}
being a triple $(f, g, h)$ such that the diagram
\begin{tikzcd}%[sep=large]
{X} \ar{d}{f} \ar{r}{u}
&{Y} \ar{d}{g} \ar{r}{v}
&{Z} \ar{d}{h} \ar{r}{w}
&{\Sigma X} \ar{d}{\Sigma f}\\
{X^\prime} \ar{r}[swap]{u^\prime}
&{Y^\prime} \ar{r}[swap]{v^\prime}
&{Z^\prime}  \ar{r}[swap]{w^\prime}
&{\Sigma X^\prime}
\end{tikzcd}
commutes.

Let \bC be an additive category $-$then a 
\un{triangulation} 
\index{triangulation (additive category)} 
of \bC 
is a pair $(\Sigma,\Delta)$, where $\Sigma$ is a suspension functor and $\Delta$ is a class of triangles 
(the 
\un{exact triangles} 
\index{exact triangles}), 
subject to the following assumptions.

\indent\indent (TR$_1$) \quad Every triangle isomorphic to an exact triangle is exact.

\indent\indent (TR$_2$) \quad For any $X \in \Ob\bC$, the triangle $X \overset{\id_X}{\lra} X \ra 0 \ra \Sigma X$ is exact.

\indent\indent (TR$_3$) \quad Every morphism $X \overset{u}{\ra} Y$ can be completed to an exact triangle 
$X \overset{u}{\ra} $
$Y \overset{v}{\ra} $
$Z \overset{w}{\ra} $
$\Sigma X$.\\
\indent\indent (TR$_4$) \quad The triangle 
$X \overset{u}{\ra} $
$Y \overset{v}{\ra} $
$Z \overset{w}{\ra} $
$\Sigma X$ 
is exact iff the triangle 
$Y \overset{v}{\ra} $
$Z \overset{w}{\ra} $
$\Sigma X \overset{-\Sigma u}{\lra} $
$\Sigma Y$
is exact.\\
\indent\indent (TR$_5$) \quad If 
$X \overset{u}{\ra} $
$Y \overset{v}{\ra} $
$Z \overset{w}{\ra} $
$\Sigma X$
, 
$X^\prime \overset{u^\prime}{\ra} $
$Y^\prime \overset{v^\prime}{\ra} $
$Z^\prime \overset{w^\prime}{\ra} $
$\Sigma X^\prime$
are exact triangles and if in the diagram 
\begin{tikzcd}%[sep=large]
{X} \ar{d}{f} \ar{r}{u}
&{Y} \ar{d}{g} \ar{r}{v}
&{Z}  \ar{r}{w}
&{\Sigma X} \ar{d}{\Sigma f}\\
{X^\prime} \ar{r}[swap]{u^\prime}
&{Y^\prime} \ar{r}[swap]{v^\prime}
&{Z^\prime}  \ar{r}[swap]{w^\prime}
&{\Sigma X^\prime}
\end{tikzcd}
, $g \circx u = u^\prime \circx f$, then there is a morphism $h:Z \ra Z^\prime$ such that $(f,g,h)$ is a morphism of triangles.\\

\begingroup%%----------------------------------->>
\fontsize{9pt}{11pt}\selectfont
\label{15.11}
\textbf{\small EXAMPLE}  \ 
Suppose that 
$X \overset{u}{\ra} $
$Y \overset{v}{\ra} $
$Z \overset{w}{\ra} $
$\Sigma X$ 
is exact.  Let 
$f:X \ra X^\prime$, 
$g:Y \ra Y^\prime$, 
$h:Z \ra Z^\prime$, 
be isomorphisms.  Put
$u^\prime = g \circx u \circx f^{-1}$, 
$v^\prime = h \circx v \circx g^{-1}$, 
$w^\prime = \Sigma f \circx w \circx h^{-1}$ 
$-$then 
$X^\prime \overset{u^\prime}{\ra} $
$Y^\prime \overset{v^\prime}{\ra} $
$Z^\prime \overset{w^\prime}{\ra} $
$\Sigma X^\prime$ 
is exact (cf. TR$_1$).  
Examples: 
(1) 
$X \overset{-u}{\ra} $
$Y \overset{-v}{\ra} $  
$Z \overset{w}{\ra} $ 
$\Sigma X$
is exact; 
(2) 
$Y \overset{-v}{\ra}$ 
$Z \overset{-w}{\ra}$ 
$\Sigma  X \overset{-\Sigma u}{\lra}$
$\Sigma Y$ 
is exact (cf. TR$_4$).\\
\endgroup %%------------------------------------<<

%%----------------------------------------------------------------------------------------------02
\label{15.3}
\label{15.33b} %dmc mnft
\label{16.23} %dmc mnft
\begingroup%%----------------------------------->>
\fontsize{9pt}{11pt}\selectfont
\textbf{\small EXAMPLE}  \ 
$\forall \ X \in \Ob\bC$, the triangle 
$0 \ra X \overset{\id_X}{\lra}  X \ra 0$ $(= \Sigma 0)$ is in $\Delta$ (cf. TR$_2$ $\&$ TR$_4$).\\
\endgroup %%------------------------------------<<

\label{15.4} %dmc may need fine tuning
\label{15.6} %dmc may need fine tuning
\begingroup%%----------------------------------->>
\fontsize{9pt}{11pt}\selectfont
\textbf{\small EXAMPLE}  \ 
Suppose that 
$X \overset{u}{\ra} Y \overset{v}{\ra} Z \overset{w}{\ra} \Sigma X$ is exact $-$then there is a commutative diagram
\begin{tikzcd}[sep=large]
{X} \arrow[d,shift right=0.5,dash] \arrow[d,shift right=-0.5,dash]  \ar{r}{u}  
&{Y} \arrow[d,shift right=0.5,dash] \arrow[d,shift right=-0.5,dash]  \ar{r}{\nu_Z^{-1} \circx v} 
&{\Sigma\Omega Z} \ar{d}{\nu_Z} \ar{r}{w \circx \nu_Z}
&{\Sigma X} \arrow[d,shift right=0.5,dash] \arrow[d,shift right=-0.5,dash]\\
{X}  \ar{r}[swap]{u}  
&{Y} \ar{r}[swap]{v}
&{Z}  \ar{r}[swap]{w}
&{\Sigma X}
\end{tikzcd}
, thus the triangle 
\begin{tikzcd}[sep=large]
{X}  \ar{r}{u}  
&{Y}  \ar{r}{\nu_Z^{-1} \circx v} 
&{\Sigma\Omega Z}  \ar{r}{w \circx v_Z}
&{\Sigma X}
\end{tikzcd}
is exact (cf. TR$_1$) and so, by TR$_4$, the triangle 
\begin{tikzcd}[sep=large]
{\Omega Z}  \ar{rr}{-(\mu_X^{-1} \circx \Omega w)} 
&&{X}  \ar{r}{u}
&{Y} \ar{r}{\nu_Z^{-1} \circx v}
&{\Sigma\Omega Z}
\end{tikzcd}
is exact.
\vspi
[Note: \ Under the bijection of adjunction $\Mor(Z,\Sigma X) \approx \Mor(\Omega Z,X)$, $w$ corresponds to 
$\mu_X^{-1} \circx \Omega w$ and $\Sigma(\mu_X^{-1} \circx \Omega w)$ equals $w \circx \nu_Z$.]\\
\endgroup %%------------------------------------<<

\label{15.13}
\begingroup%%----------------------------------->>
\fontsize{9pt}{11pt}\selectfont
\textbf{\small EXAMPLE}  \ 
Suppose that 
$X \overset{u}{\ra} $
$Y \overset{v}{\ra} $
$Z \overset{w}{\ra} $
$\Sigma X$
is exact $-$then there is a commutative diagram 
\[
\begin{tikzcd}[sep=large]
{\Omega Z} \arrow[d,shift right=0.5,dash] \arrow[d,shift right=-0.5,dash]  \ar{rr}{-(\mu_X^{-1} \circx \Omega w)} 
&&{X} \arrow[d,shift right=0.5,dash] \arrow[d,shift right=-0.5,dash]  \ar{rr}{\nu_Y^{-1} \circx u}
&&{\Sigma\Omega Y} \ar{d}{\nu_Y} \ar{rr}{\Sigma\Omega v}
&&{\Sigma\Omega Z} \arrow[d,shift right=0.5,dash] \arrow[d,shift right=-0.5,dash] \\
{\Omega Z} \ar{rr}[swap]{-(\mu_X^{-1} \circx \Omega w)} 
&&{X} \ar{rr}[swap]{u}
&&{Y} \ar{rr}[swap]{\nu_Z^{-1} \circx v}
&&{\Sigma\Omega Z}
\end{tikzcd}
,\]
thus the triangle 
\begin{tikzcd}[sep=large]
{\Omega Z} \ar{rr}{-(\mu_X^{-1} \circx \Omega w)} 
&&{X}  \ar{rr}{\nu_Y^{-1} \circx u}
&&{\Sigma\Omega Y} \ar{rr}{\Sigma\Omega v}
&&{\Sigma\Omega Z}
\end{tikzcd}
is exact (cf. TR$_1$) and so, by TR$_4$, the triangle 
\begin{tikzcd}[sep=large]
{\Omega Y} \ar{r}{-\Omega v}
&{\Omega Z} \ar{rr}{-(\mu_X^{-1} \circx \Omega w)} 
&&{X} \ar{rr}{\nu_Y^{-1} \circx u}
&&{\Sigma\Omega Y} 
\end{tikzcd}
is exact or still, the triangle 
\begin{tikzcd}[sep=large]
{\Omega Y} \ar{r}{\Omega v}
&{\Omega Z} \ar{rr}{\mu_X^{-1} \circx \Omega w} 
&&{X} \ar{rr}{\nu_Y^{-1} \circx u}
&&{\Sigma\Omega Y} 
\end{tikzcd}
is exact.\\
\endgroup %%------------------------------------<<

A 
\un{triangulated category}
\index{triangulated category} 
is an additive category \bC equipped with a triangulation 
$(\Sigma,\Delta)$.

[Note: \ The opposite of a triangulated category is triangulated.  
In detail: The suspension functor is $\Omega^{\OP}$ and the elements of $\Delta^{\OP}$ are those triangles 
$X \overset{u^{\OP}}{\ra} $ 
$Y \overset{v^{\OP}}{\ra} $ 
$Z \overset{w^{\OP}}{\ra} $
$\Omega^{\OP} X$ 
in $\bC^\OP$ such that 
$\Omega X \overset{-w}{\lra} $ 
$Z \  \overset{v}{\ra}
\begin{tikzcd}[sep=small]
Y \ar{rr}{{\nu_X^{-1} \circx u}} &&{\Sigma\Omega X}
\end{tikzcd}
$
is exact.]

Example: Let \bC be a triangulated category.  Call a triangle 
$X \overset{u}{\ra} Y \overset{v}{\ra} Z \overset{w}{\ra} \Sigma X$ 
\un{antiexact}
\index{antiexact (triangle)}
 if the triangle 
$X \overset{u}{\ra} Y \overset{v}{\ra} Z \overset{-w}{\lra} \Sigma X$ 
is exact $-$then \bC endowed with the class of antiexact triangles is triangulated.\\

\begingroup%%----------------------------------->>
\fontsize{9pt}{11pt}\selectfont
\textbf{\small EXAMPLE}  \ 
Let \bA be an abelian category.  Write $\bC\bX\bA$ for the abelian category of cochain complexes over \bA.  Let 
$\Sigma:\bC\bX\bA \ra \bC\bX\bA$ be the additive functor that sends \mX to $X[1]$, where 
$
\begin{cases}
\ X[1]^n = X^{n+1}\\
\ d_{X[1]}^n = -d_X^{n+1}
\end{cases}
$
$-$then $\Sigma$ is an automorphism of $\bC\bX\bA$, hence is a suspension functor.  
The quotient category $\bK(\bA)$ of 
$\bC\bX\bA$ per cochain homotopy is an additive category and the projection 
$\bC\bX\bA \ra \bK(\bA)$ is an additive functor.  
Moreover, $\Sigma$ induces a suspension functor $\bK(\bA) \ra \bK(\bA)$.  
Definition: A triangle 
$X^\prime \overset{u^\prime}{\ra} $
$Y^\prime \overset{v^\prime}{\ra} $
$Z^\prime \overset{w^\prime}{\ra} $
$\Sigma X^\prime$
in $\bK(\bA)$ is exact if it is isomorphic to a triangle 
$X \overset{f}{\ra} $
$Y \overset{j}{\ra} $
$C_f \overset{\pi}{\ra} $
$\Sigma X$
for some $f$.  Here
%%----------------------------------------------------------------------------------------------03
$C_f$ is the mapping cone of $f:C_f^n = X^{n+1} \oplus Y^n$, 
$
d_{C_f}^n = 
\begin{pmatrix}
d_{\Sigma X}^n &0\\
f^{n+1} &d_Y^n\\
\end{pmatrix}
$
$
\bigl(j^n = 
\begin{pmatrix}
0\\ 
\id_{Y^n}\\
\end{pmatrix}
,\ 
\pi^n = (\id_{X^{n+1}},0)\bigr)$.  
With these choices, one can check by direct computation that $\bK(\bA)$ is triangulated (a detailed explanation can be found in
Kashiwara-Schapira\footnote[2]{\textit{Sheaves on Manifolds}, Springer Verlag (1990), 35-38; 
see also Weibel, \textit{An Introduction to Homological Algebra}, Cambridge University Press (1994), 376.}).\\
\endgroup %%------------------------------------<<

\label{17.80}
\begin{proposition} \ %01
Let \bC be a triangulated category.  Suppose that 
\[
\begin{tikzcd}%[sep=large]
{X} \ar{r}{u} 
&{Y} \ar{d}{g}  \ar{r}{v}  
&{Z} \ar{d}{h} \ar{r}{w}
&{\Sigma X}\\
{X^\prime} \ar{r}[swap]{u^\prime} 
&{Y^\prime}  \ar{r}[swap]{v^\prime}  
&{Z^\prime}  \ar{r}[swap]{w^\prime}
&{\Sigma X^\prime}
\end{tikzcd}
\]
is a diagram with rows in $\Delta$.  Assume: $h \circx v = v^\prime \circx g$ $-$then there is a morphism 
$f:X \ra X^\prime$ such that $(f,g,h)$ is a morphism of triangles.
\end{proposition}

[Bearing in mind TR$_4$, pass to 
\begin{tikzcd}%[sep=large]
{Y} \ar{d}{g}  \ar{r}{v}  
&{Z} \ar{d}{h} \ar{r}{w}
&{\Sigma X} \ar{r}{-\Sigma u}
&{\Sigma Y} \ar{d}{\Sigma g}\\
{Y^\prime}  \ar{r}[swap]{v^\prime}  
&{Z^\prime}  \ar{r}[swap]{w^\prime}
&{\Sigma X^\prime} \ar{r}[swap]{-\Sigma u^\prime}
&{\Sigma Y^\prime}
\end{tikzcd}
\ and apply TR$_5$.]\\

\begin{proposition} \ %2
Let \bC be a triangulated category.  Suppose that 
\[
\begin{tikzcd}%[sep=large]
{X} \ar{d}{f} \ar{r}{u} 
&{Y}   \ar{r}{v}  
&{Z} \ar{d}{h} \ar{r}{w}
&{\Sigma X}\ar{d}{\Sigma f} \\
{X^\prime} \ar{r}[swap]{u^\prime} 
&{Y^\prime}  \ar{r}[swap]{v^\prime}  
&{Z^\prime}  \ar{r}[swap]{w^\prime}
&{\Sigma X^\prime}
\end{tikzcd}
\]
is a diagram with rows in $\Delta$.  Assume: $\Sigma f \circx w = w^\prime \circx h$ $-$then there is a morphism 
$g:Y \ra Y^\prime$ such that $(f,g,h)$ is a morphism of triangles.\\
\end{proposition}

\begin{proposition} \ %3
Let \bC be a triangulated category $-$then for any exact triangle 
$X \overset{u}{\ra}$ 
$Y \overset{v}{\ra}$ 
$Z \overset{w}{\ra}$ 
$\Sigma X$, 
$v \circx u = 0$ and $w \circx v = 0$.
\end{proposition}

[It suffices to prove that $v \circx u = 0$.  But the diagram
\[ 
\begin{tikzcd}%[sep=large]
{X} \arrow[r,shift right=0.5,dash] \arrow[d,shift right=-0.5,dash] \arrow[d,shift right=0.5,dash] \arrow[r,shift right=-0.5,dash]  
&{X} \ar{d}{u}  \ar{r}
&{0} \ar[dashed]{d} \ar{r}{w}
&{\Sigma X}\arrow[d,shift right=-0.5,dash] \arrow[d,shift right=0.5,dash]\\
{X} \ar{r}[swap]{u} 
&{Y}  \ar{r}[swap]{v}  
&{Z}  \ar{r}[swap]{w}
&{\Sigma X}
\end{tikzcd}
\]
must commute (cf. TR$_5$), thus $v \circx u = 0$.]\\

Application: Every morphism $X \overset{u}{\ra} Y$ admits a weak cokernel.

[Thanks to TR$_3$, $\exists$ an exact triangle 
$X \overset{u}{\ra} $
$Y \overset{v}{\ra} $
$Z \overset{w}{\ra} $
$\Sigma X$
and $v \circx u = 0$.  On the other 
%%----------------------------------------------------------------------------------------------04
hand, if $g \circx u = 0$ $(g:Y \ra W)$, then the diagram 
\begin{tikzcd}%[sep=large]
{X} \ar{d} \ar{r}{u}
&{Y} \ar{d}{g} \ar{r}{v}
&{Z} \ar{r}{w}
&{\Sigma X}\ar{d}\\
{0}\ar{r}
&{W}\ar{r}[swap]{\id_W}
&{W}\ar{r}
&{0}
\end{tikzcd}
has a filler $h:Z \ra W$ such that $h \circx v = g$ (cf. TR$_5$)
.]\\

\label{15.26}
\label{15.28} %dmc mnft orig 15-4
\label{15.43} %dmc mnft orig 15-4
\label{15.43} %dmc mnft orig 15-4
\label{15.44} %dmc mnft orig 15-4
\begingroup%%----------------------------------->>
\fontsize{9pt}{11pt}\selectfont
Suppose that a triangulated category \bC has coproducts $-$then \bC has weak pushouts, hence weak co\-limits.  
One can be specific.  
Thus let $\Delta:\bI \ra \bC$ be a diagram.  
Given $\delta \in \Mor \bI$, say $i \overset{\delta}{\ra} j$, put 
$s\delta = i$, 
$t\delta = j$.  
Define an arrow 
$\ds\coprod\limits_{\Mor \bI} \Delta_{s\delta} \ra \ds\coprod\limits_{\Ob\bI} \Delta_i$ by taking the coproduct of the arrows
\begin{tikzcd}[sep=large]
{\Delta_{s\delta}}  \ar{rr}{
\begin{pmatrix}
\id\\
-\Delta\delta\\
\end{pmatrix}
}
&&{\Delta_{s\delta} \amalg \Delta_{s\delta}}
\end{tikzcd}
$-$then a candidate for a weak colimit of $\Delta$ is any completion \mL of 
$\ds\coprod\limits_{\Mor \bI} \Delta_{s\delta} \ra \ds\coprod\limits_{\Ob\bI} \Delta_i$ to an exact triangle (cf. TR$_3$).\\
\endgroup %%------------------------------------<<

Let \bC be a triangulated category, \bD an abelian category $-$then an additive functor (cofunctor) 
$F:\bC \ra \bD$ is said to be 
\un{exact}
\index{exact (additive functor (cofunctor))} 
if for every exact triangle 
$X \overset{u}{\ra} $
$Y \overset{v}{\ra} $
$Z \overset{w}{\ra} $
$\Sigma X,$ 
the sequence
$FX \ra FY \ra FZ$ 
($FZ \ra FY \ra FX$) 
is exact.

[Note: \ An exact functor (cofunctor) generates a long exact sequence involving $\Sigma$ and $\Omega$.]\\

\begin{proposition} \ %4
Let \bC be a triangulated category $-$then $\forall$ $W \in \Ob\bC$, $\Mor(W,-)$ is an exact functor and $\Mor(-,W)$ is an exact cofunctor.
\end{proposition}

[Take any exact triangle 
$X \overset{u}{\ra} $
$Y \overset{v}{\ra} $
$Z \overset{w}{\ra} $
$\Sigma X$ 
and consider 
$\Mor(W,X) \overset{u_*}{\ra} $
$\Mor(W,Y) \overset{v_*}{\ra} $
$\Mor(W,Z)$.
In view of Proposition 3, $\im u_* \subset \ker v_*$.  To go the other way, assume that $v \circx \psi = 0$ 
$(\psi \in \Mor(W,Y))$ $-$then $\exists$ $\phi \in \Mor(W,X)$: $\psi = u \circx \phi$.  
Proof: Examine 
\begin{tikzcd}%[sep=large]
{W} \ar[dashed]{d}{\phi} \ar{r}{\id_W}
&{W} \ar{d}{\psi}  \ar{r}
&{0} \ar{d} \ar{r}
&{\Sigma W}\ar[dashed]{d}{\Sigma \phi}\\
{X} \ar{r}[swap]{u} 
&{Y}  \ar{r}[swap]{v}  
&{Z}  \ar{r}[swap]{w}
&{\Sigma X}
\end{tikzcd}
(cf. Proposition 1).]\\
\vspace{0.5cm}

\label{15.1}
\label{15.35}
\label{15.33c}
Application: If 
\begin{tikzcd}%[sep=large]
{X} \ar{d}{f} \ar{r}{u}
&{Y} \ar{d}{g} \ar{r}{v}
&{Z} \ar{d}{h} \ar{r}{w}
&{\Sigma X} \ar{d}{\Sigma f}\\
{X^\prime} \ar{r}[swap]{u^\prime}
&{Y^\prime} \ar{r}[swap]{v^\prime}
&{Z^\prime}  \ar{r}[swap]{w^\prime}
&{\Sigma X^\prime}
\end{tikzcd}
is a \cd with rows in $\Delta$ and if any two of $f, g, h$ are isomorphisms, then so is the third.

[For instance, suppose that $f$ and $g$ are isomorphisms $-$then the five lemma implies that 
$h_*:\Mor(Z^\prime,Z) \ra \Mor(Z^\prime,Z^\prime)$, 
$h^*:\Mor(Z^\prime,Z) \ra \Mor(Z,Z)$ are isomorphisms so 
$\exists$ $\phi, \ \psi \in \Mor(Z^\prime,Z)$: 
$h \circx \phi = \id_{Z^\prime}$, $\psi \circx h = \id_Z$, i.e., $h$ is an isomorphism.]\\

\begingroup%%----------------------------------->>
\fontsize{9pt}{11pt}\selectfont
\label{15.17}
\textbf{\small EXAMPLE}  \ 
Let \bC be a triangulated category with finite coproducts $-$then $\forall$ $X, \ Y \in \Ob\bC$, the triangle 
$X \ra X \amalg Y \ra Y \overset{0}{\ra} \Sigma X$ is in $\Delta$.
\vspi
%%----------------------------------------------------------------------------------------------05
[According to TR$_3$, the morphism $X \ra X \amalg Y$ can be completed to an exact triangle 
$X \ra X \amalg Y$ 
$\ra Z \ra \Sigma X$.  
Compare it with the exact triangle 
$0 \ra Y$
$\overset{\id_Y}{\lra} Y$
$\ra 0$ 
to get a filler $h:Z \ra Y$ (cf. TR$_5$).  Consideration of 
\[
\begin{tikzcd}[sep=large]
{0} \ar{r} 
&{\Mor(W,\Sigma X)} \ar{r}
&{\Mor(W,\Sigma X \amalg \Sigma Y)} \ar{d} \ar{r}
&{\Mor(W,\Sigma Z)} \ar{d}{(\Sigma h)_*} \ar{r}
&{0}\\
&&{\Mor(W,\Sigma Y)} \arrow[r,shift right=0.5,dash] \arrow[r,shift right=-0.5,dash] 
&{\Mor(W,\Sigma Y)}
\end{tikzcd}
\]
allows one to say that $(\Sigma h)_*$ is an isomorphism $\forall \ W$, hence $\Sigma h$ is an isomorphism or still, $h$ is an isomorphism.]\\
\endgroup %%------------------------------------<<

\label{15.10} %dmc mnft
\label{15.12} %dmc mnft
\label{15.21} %dmc mnft
\label{17.76} %dmc mnft
\begingroup%%----------------------------------->>
\fontsize{9pt}{11pt}\selectfont
\textbf{\small EXAMPLE}  \ 
Let \bC be a triangulated category with finite coproducts $-$then any exact triangle of the form 
$X \overset{u}{\ra}$ 
$Y \overset{v}{\ra}$ 
$Z \overset{0}{\ra}$ 
$\Sigma X$ 
is isomorphic to 
$X \ra$ 
$X \amalg Z \ra$ 
$Z \overset{0}{\ra}$ 
$\Sigma X.$  
Indeed, the triangle 
$Y \overset{v}{\ra}$ 
$Z \overset{0}{\ra}$ 
$\Sigma X \overset{-\Sigma u}{\lra} \Sigma Y$
is exact (cf. TR$_4$) and there is a morphism $Y \ra X \amalg Z$ rendering the diagram
\begin{tikzcd}[sep=large]
{Y} \ar[dashed]{d} \ar{r}{v}
&{Z} \arrow[d,shift right=0.5,dash] \arrow[d,shift right=-0.5,dash]  \ar{r}{0}
&{\Sigma X} \arrow[d,shift right=0.5,dash] \arrow[d,shift right=-0.5,dash]  \ar{r}{-\Sigma u}
&{\Sigma Y} \ar[dashed]{d}\\
{X \amalg Z}  \ar{r}
&{Z} \ar{r}[swap]{0}
&{\Sigma X} \ar{r}
&{\Sigma X \amalg \Sigma Z}
\end{tikzcd}
commutative (cf. Proposition 1).
\vspi
[Note: \ Analogously, an exact triangle of the form 
$X \overset{0}{\ra}$ 
$Y \overset{v}{\ra}$ 
$Z \overset{w}{\ra}$ 
$\Sigma X$ 
is isomorphic to 
$X \overset{0}{\ra}$ 
$Y \ra$ 
$Y \amalg \Sigma X \ra$ 
$\Sigma X$ 
.]\\
\endgroup %%------------------------------------<<

\label{15.8}
\label{15.16}
\begingroup%%----------------------------------->>
\fontsize{9pt}{11pt}\selectfont
\textbf{\small EXAMPLE}  \ 
Let \bC be a triangulated category with finite coproducts.  
Suppose given a morphism $i:X \ra Y$ that admits a left inverse 
$r:Y \ra X$ $-$then there exists an isomorphism $Y \ra X \amalg Z$ and a commutative triangle 
\begin{tikzcd}[sep=large]
{X} \ar{rd} \ar{r}{i} &{Y}\ar{d}\\
&{X \amalg Z} 
\end{tikzcd}
.
\vspi
[Complete $X \overset{i}{\ra} Y$ to an exact triangle 
$X \overset{i}{\ra}$ 
$Y \overset{v}{\ra}$ 
$Z \overset{w}{\ra}$ 
$\Sigma X$ (cf. TR$_3$) and choose a filler 
\begin{tikzcd}[sep=large]
{X} \arrow[d,shift right=0.5,dash] \arrow[d,shift right=-0.5,dash]  \ar{r}{i}
&{Y} \ar{d}{r} \ar{r}{v}
&{Z} \ar[dashed]{d} \ar{r}{w}
&{\Sigma X}\arrow[d,shift right=0.5,dash] \arrow[d,shift right=-0.5,dash] \\
{X} \ar{r}[swap]{\id_X}
&{X} \ar{r}
&{0} \ar{r}
&{\Sigma X}
\end{tikzcd}
(cf. TR$_5$) to see that $w = 0$.]\\
\endgroup %%------------------------------------<<

\label{15.9}
\begingroup%%----------------------------------->>
\fontsize{9pt}{11pt}\selectfont
\textbf{\small EXAMPLE}  \ 
Let \bC be a triangulated category with finite coproducts $-$then the triangles
$X \overset{u}{\ra}$ 
$Y \overset{v}{\ra}$ 
$Z \overset{w}{\ra}$ 
$\Sigma X$, 
$X^\prime \overset{u^\prime}{\ra}$ 
$Y^\prime \overset{v^\prime}{\ra}$ 
$Z^\prime \overset{w^\prime}{\ra}$ 
$\Sigma X^\prime$ 
are exact iff the triangle\\
\vspace{0.15cm}
\begin{tikzcd}[sep=large]
{X \amalg X^\prime} \ar{rr}{
\begin{pmatrix}
u \ 0\\
0 \ u^\prime\\
\end{pmatrix}
}
&&{Y \amalg Y^\prime}
\ar{rr}{
\begin{pmatrix}
v \ 0\\
0 \ v^\prime\\
\end{pmatrix}
}
&&{Z \amalg Z^\prime}
\ar{rr}{
\begin{pmatrix}
w \ 0\\
0 \ w^\prime\\
\end{pmatrix}
}
&&{\Sigma X \amalg \Sigma X^\prime}
\end{tikzcd}
is exact.\\
\endgroup %%------------------------------------<<

\begingroup%%----------------------------------->>
\fontsize{9pt}{11pt}\selectfont
\textbf{\small EXAMPLE}  \ 
Let \bC be a triangulated category with finite coproducts.  Suppose that 
$X \overset{u}{\ra}$ 
$Y \overset{v}{\ra}$ 
$Z \overset{w}{\ra}$ 
$\Sigma X$, 
%%----------------------------------------------------------------------------------------------06
is exact $-$then for any $Y^\prime \in \Ob\bC$ and any $g \in \Mor(Y,Y^\prime)$, the triangle 
\begin{tikzcd}[sep=small]
{Y^\prime \amalg Y} \ar{rr}{
\begin{pmatrix}
\id_{Y^\prime} \ g\\
0 \ \  -v\\
\end{pmatrix}
}
&&{Y^\prime \amalg Z}% \ar{r}{(0 - w)}
\end{tikzcd}
\begin{tikzcd}[sep=small]
{} \ar{rr}{(0 - w)}
&&{\Sigma X} \ar{rr}{
\begin{pmatrix}
\Sigma (g\circx u)\\
-\Sigma u\\
\end{pmatrix}
}
&&{\Sigma Y^\prime \amalg \Sigma Y}
\end{tikzcd}
is exact.\\
\endgroup %%------------------------------------<<

\label{15.2}
\label{15.51} %dmc mnft
\label{15.52} %dmc mnft
\label{15.33a}
\begingroup%%----------------------------------->>
\fontsize{9pt}{11pt}\selectfont
\textbf{\small FACT} \ 
Let \bC be a triangulated category $-$then a morphism $X \overset{u}{\ra} Y$ is an isomorphism iff the triangle 
$X \overset{u}{\ra}$ 
$Y \ra 0$
$\ra \Sigma X$ 
is exact.\\
\endgroup %%------------------------------------<<

\begingroup%%----------------------------------->>
\fontsize{9pt}{11pt}\selectfont
\textbf{\small FACT} \ 
Let \bC be a triangulated category.  Suppose that 
$X \overset{u}{\ra}$ 
$Y \overset{v}{\ra}$ 
$Z \overset{w_i}{\lra}$ 
$\Sigma X$  
$(i = 1,2)$ are exact triangles $-$then $w_1 = w_2$ if $\Mor(\Sigma X, Z) = 0$.\\
\endgroup %%------------------------------------<<

\label{15.15}
\begin{proposition} \ %5
Let \bC be a triangulated category.  Fix a morphism $X \overset{u}{\ra} Y$ in \bC and suppose that 
$X \overset{u}{\ra}$ 
$Y \overset{v}{\ra}$ 
$Z \overset{w}{\ra}$ 
$\Sigma X$, 
$X \overset{u}{\ra}$ 
$Y \overset{v^\prime}{\ra}$ 
$Z^\prime \overset{w^\prime}{\ra}$ 
$\Sigma X$ 
are exact triangles (cf. TR$_3$) $-$then $Z \approx Z^\prime$.
\end{proposition}

[
Any filler for 
\begin{tikzcd}%[sep=large]
{X} \arrow[d,shift right=0.5,dash] \arrow[d,shift right=-0.5,dash]   \ar{r}{u}  
&{Y}\arrow[d,shift right=0.5,dash] \arrow[d,shift right=-0.5,dash]   \ar{r}{v} 
&{Z} \ar[dashed]{d} \ar{r}{w}
&{\Sigma X} \arrow[d,shift right=0.5,dash] \arrow[d,shift right=-0.5,dash] \\
{X}  \ar{r}[swap]{u}  
&{Y} \ar{r}[swap]{v^\prime}
&{Z^\prime}  \ar{r}[swap]{w^\prime}
&{\Sigma X}
\end{tikzcd}
is an isomorphism (cf. p. \pageref{15.1}).]\\

Let \bC be a triangulated category $-$then a full, isomorphism closed subcategory \bD of \bC containing 0 and stable under $\Sigma$ and $\Omega$ is said to be a 
\un{triangulated subcategory}
\index{triangulated subcategory} 
of \bC if $\forall$ $X \overset{u}{\ra} Y$ in $\Mor \bD$, 
there exists an exact triangle 
$X \overset{u}{\ra}$ 
$Y \overset{v}{\ra}$ 
$Z \overset{w}{\ra}$ 
$\Sigma X$
with \mZ in $\Ob\bD.$

[Note: \ \bD is, in its own right, a triangulated category (the suspension functor is the restriction of $\Sigma$ to \bD and the exact triangles 
$X \overset{u}{\ra}$ 
$Y \overset{v}{\ra}$ 
$Z \overset{w}{\ra}$ 
$\Sigma X$
are the elements of $\Delta$ such that $X, \ Y, \ Z \in \Ob\bD$).]\\

\begingroup%%----------------------------------->>
\fontsize{9pt}{11pt}\selectfont
\textbf{\small EXAMPLE}  \ 
Let \bA be an abelian category. 
 Write 
 $\bC\bX\bA^+$ for the full subcategory of $\bC\bX\bA$ consisiting of those \mX which are bounded below 
 $(X^n = 0 (n \ll 0))$, 
write 
$\bC\bX\bA^-$ for the full subcategory of $\bC\bX\bA$ consisiting of those \mX which are bounded above 
$(X^n = 0 (n \gg 0))$, 
and put 
$\bC\bX\bA^\text{b} = \bC\bX\bA^+ \cap \bC\bX\bA^-$ $-$then, in obvious notation, 
$\bK^+(\bA),$ $\bK^-(\bA)$, and $\bK^\text{b}(\bA)$ are triangulated subcategories of $\bK(\bA)$.\\
\endgroup %%------------------------------------<<

\begin{proposition} \ %06
Let \bC be a triangulated category.  Suppose that \mO is the object class of a triangulated subcategory of \bC $-$then for any exact triangle
$X \overset{u}{\ra}$ 
$Y \overset{v}{\ra}$ 
$Z \overset{w}{\ra}$ 
$\Sigma X$, 
if two of \mX, \mY, \mZ are in \mO so is the third.
\end{proposition}

[Assuming that $X,Y \in O$, choose $Z^\prime \in O$: 
$X \overset{u}{\ra}$ 
$Y \overset{v^\prime}{\ra}$ 
$Z^\prime \overset{w^\prime}{\ra}$ 
$\Sigma X$ 
is exact.  
On the basis of Proposition 5, $Z \approx Z^\prime$, hence $Z \in O$ (\mO is isomorphism closed).  
Next assume that 
$Y,Z \in O$ and fix an exact triangle 
$Y \overset{v}{\ra}$ 
$Z \ra W$ 
$\ra \Sigma Y$
with $W \in O$.  By TR$_4$,
%%----------------------------------------------------------------------------------------------07
$Y \overset{v}{\ra}$ 
$Z \overset{w}{\ra}$ 
\begin{tikzcd}%[sep=large]
{\Sigma X} \ar{r}{-\Sigma u} &{\Sigma Y}
\end{tikzcd}
is exact.  Therefore $W \approx \Sigma X$ (cf. Proposition 5) $\implies$ $\Omega W \approx \Omega\Sigma X$.  
But $\Omega W \in O$ $\implies$ $\Omega\Sigma X \in O$
$\implies$ $ X \in O$.  
The argument that $X, Z \in O$ $\implies$ $Y \in O$ is similar.]\\

\begin{proposition} \ %7
Let \bC be a triangulated category.  Suppose given a nonempty class $O \subset \Ob\bC$ $-$then \mO is the object class of a triangulated subcategory of \bC provided that for any exact triangle 
$X \overset{u}{\ra}$ 
$Y \overset{v}{\ra}$ 
$Z \overset{w}{\ra}$ 
$\Sigma X$, 
if two of \mX, \mY, \mZ are in \mO so is the third.
\end{proposition}

[(1) \quad $0 \in O$.  
Proof: $\forall \ X \in O$, $X \overset{\id_X}{\lra} X \ra 0 \ra \Sigma X$ is exact 
(cf. TR$_2$). 
(2) \quad \mO is isomorphism closed.  Proof: If $X \in O$ and if $X \overset{u}{\ra} X^\prime$ is an isomorphism, then the triangle 
$X \overset{u}{\ra}$ 
$X^\prime \ra$ 
$0 \ra \Sigma X$ is exact (cf. p. \pageref{15.2}).  
(3) \ \   $\Sigma O \subset O$. \ Proof: \ For any $X \in O$, 
$X \ra 0 \ra$
\begin{tikzcd}[sep=large]
{\Sigma X} \ar{r}{-\id_{\Sigma X}} &{\Sigma X}
\end{tikzcd}
is exact (cf. TR$_4$), thus $\Sigma X \in O$.  
(4) \  $\Omega O \subset O$. Proof: For any $X \in O$, 
$0 \ra X$ 
$\overset{\id_X}{\lra} X$ 
$\ra 0$ is exact (cf. p. \pageref{15.3}), hence 
$\Omega X \ra 0 \ra$ 
\begin{tikzcd}%[sep=large]
{X} \ar{r}{\nu_X^{-1}} &{\Sigma\Omega X}
\end{tikzcd}
is exact (cf. p. \pageref{15.4}), thus $\Omega X \subset O$.  
The final requirement that \mO must satisfy is clear.]\\

\begingroup%%----------------------------------->>
\fontsize{9pt}{11pt}\selectfont
\label{15.39}
\textbf{\small EXAMPLE}  \ 
Let \bC be a triangulated category, \bD an abelian category.  
Suppose that $F:\bC \ra \bD$ is an exact functor.  
Let $S_F$ be the class of morphisms $X \overset{u}{\ra} Y$ such that $\forall \ n \geq 0$, 
$
\begin{cases}
\ F\Sigma^n u\\
\ F\Omega^n u
\end{cases}
$
is an isomorphism and let $O_F$ be the class of objects \mZ for which there exists an exact triangle 
$X \overset{u}{\ra}$ 
$Y \overset{v}{\ra}$ 
$Z \overset{w}{\ra}$ 
$\Sigma X$
with $u \in S_F$ $-$then $O_F$ is the object class of a triangulated subcategory of \bC.
\vspi
[Note: \ $O_F$ is the class of objects \mZ such that $\forall \ n \geq 0$, 
$
\begin{cases}
\ F\Sigma^n Z = 0\\
\ F\Omega^n Z = 0
\end{cases}
$
.]\\
\endgroup %%------------------------------------<<

\begingroup%%----------------------------------->>
\fontsize{9pt}{11pt}\selectfont
\textbf{\small EXAMPLE}  \ 
Let \bA be an abelian category with a separator.  Suppose that $\sA$ is a Serre class in \bA $-$then $S_{\sA}^{-1}\bA$ exists 
(cf. p. \pageref{15.5}) and the composite 
$\bK(\bA) \overset{H^0}{\lra} \bA \ra S_{\sA}^{-1}\bA$ is exact, hence determines a triangulated subcategory 
$\bK_{\sA}(\bA)$ of $\bK(\bA)$ whose objects \mX are characterized by the condition that $H^n(X) \in \sA$ $\forall \ n$.\\
\endgroup %%------------------------------------<<

Let \bC, \bD be triangulated categories $-$then an additive functor $F:\bC \ra \bD$ is said to be a 
\un{triangulated functor}
\index{triangulated functor} 
if there exists a natural isomorphism 
$\Phi:F \circx \Sigma \ra \Sigma \circx F$ such that 
$X \overset{u}{\ra}$ 
$Y \overset{v}{\ra}$ 
$Z \overset{w}{\ra}$ 
$\Sigma X$
exact $\implies$ 
$FX \overset{Fu}{\ra}$ 
$FY \overset{Fv}{\ra}$ 
%$FZ \overset{\Phi_x \circx Fw}{\lra}$ 
\begin{tikzcd}[sep=large]
{FZ} \ar{r}{\Phi_x \circx Fw} &{\Sigma FX}
\end{tikzcd}
%$\Sigma FX$
exact.
Example: The inclusion functor determined by a triangulated subcategory of a triangulated category is a triangulated functor.\\

\begingroup%%----------------------------------->>
\fontsize{9pt}{11pt}\selectfont
\label{15.42}
\textbf{\small FACT} \ 
Let \bC, \bD be triangulated categories, $F:\bC \ra \bD$ a triangulated functor.  Assume: $G:\bD \ra \bC$ is a left adjoint for \mF $-$then \mG is triangulated.
\vspi
[Note: \ The same conclusion obtains if \mG is a right adjoint for \mF.  Proof: $G^\OP$ is a left adjoint for $F^\OP$, hence $G^\OP$ is triangulated, which implies that \mG is triangulated.]\\
\endgroup %%------------------------------------<<

%%----------------------------------------------------------------------------------------------08
Let \bC, \bD be triangulated categories $-$then a triangulated functor $F:\bC \ra \bD$ is said to be a 
\un{triangulated equivalence}
\index{triangulated equivalence} 
if there exists a triangulated functor $G:\bD \ra \bC$ and 
natural isomorphisms 
$
\begin{cases}
\ \mu:\id_{\bC} \ra G \circx F\\
\ \nu:F \circx G \ra \id_{\bD}
\end{cases}
$
such that the diagrams 
\begin{tikzcd}%[sep=large]
{GF\Sigma X} \ar{rr}{(\Psi F) \circx (G\Phi)} 
&&{\Sigma GFX}\\
{\Sigma X} \ar{u}{\mu_{\Sigma X}} \arrow[rr,shift right=0.5,dash] \arrow[rr,shift right=-0.5,dash] 
&&{\Sigma X} \ar{u}[swap]{\Sigma_{\mu X},}
\end{tikzcd}
\begin{tikzcd}%[sep=large]
{FG\Sigma Y} \ar{d}[swap]{\nu_{\Sigma Y}} \ar{rr}{(\Phi G) \circx (F\Psi)} 
&&{\Sigma FGY} \ar{d}{\Sigma_{\nu Y}}\\
{\Sigma Y}  \arrow[rr,shift right=0.5,dash] \arrow[rr,shift right=-0.5,dash]
&&{\Sigma Y} 
\end{tikzcd}
commute.

[Note: \ $\Phi$ and $\Psi$ are the natural isomorphisms implicit in the definition of \mF and \mG.]\\

\begingroup%%----------------------------------->>
\fontsize{9pt}{11pt}\selectfont
\textbf{\small FACT} \ 
Let \bC, \bD be triangulated categories, $F:\bC \ra \bD$ an additive functor.  
Suppose that there exists a natural transformation 
$\Phi:F \circx \Sigma \ra \Sigma \circx F$ such that 
$X \overset{u}{\ra}$ 
$Y \overset{v}{\ra}$ 
$Z \overset{w}{\ra}$ 
$\Sigma X$
exact $\implies$ 
$FX \overset{Fu}{\lra}$ 
$FY \overset{Fv}{\lra}$ 
\begin{tikzcd}[sep=large]
{FZ} \ar{r}{\Phi_X \circx Fw} &{\Sigma FX}
\end{tikzcd}
exact $-$then $\Phi$ is a natural isomorphism.
\vspi
[For any $X \in \Ob\bC$, the triangle 
$X \ra 0 \ra$
\begin{tikzcd}[sep=large]
{\Sigma X} \ar{r}{-\id_{\Sigma X}} &{\Sigma X}
\end{tikzcd}
is exact.]\\
\endgroup %%------------------------------------<<

\begingroup%%----------------------------------->>
\fontsize{9pt}{11pt}\selectfont
\textbf{\small FACT} \ 
Let \bC, \bD be triangulated categories, $F:\bC \ra \bD$ a triangulated functor.  
Assume: \mF is an equivalence $-$then \mF is a triangulated equivalence.
\vspi
[Given \mG and natural isomorphisms 
$
\begin{cases}
\ \mu:\id_{\bC} \ra G \circx F\\
\ \nu:F \circx G \ra \id_{\bD}
\end{cases}
,
$
consider the inverse of $(G \Sigma \nu) \circx (G \Phi G) \circx (\mu \Sigma G)$.]\\
\endgroup %%------------------------------------<<

Let \bC be a triangulated category $-$then \bC is said to be 
\un{strict}
\index{strict (triangulated category)} 
if its suspension functor $\Sigma$ is an isomorphism (and not just an equivalence).

[Note: \ When \bC is strict, the role of $\Omega$ is played by $\Sigma^{-1}$.]

Example: For any abelian category \bA, $\bK(\bA)$ is a strict triangulated category.\\

\begingroup%%----------------------------------->>
\fontsize{9pt}{11pt}\selectfont
\textbf{\small EXAMPLE}  \ 
Let \bC be a strict triangulated category.  Suppose that 
$X \overset{u}{\ra}$ 
$Y \overset{v}{\ra}$ 
$Z \overset{w}{\ra}$ 
$\Sigma X$
is exact $-$then 
\begin{tikzcd}[sep=large]
{\Sigma^{-1} Z} \ar{r}{-\Sigma^{-1}w} &{X}
\end{tikzcd}
$\overset{u}{\ra} Y$
$\overset{v}{\ra} Z$
is exact (cf. p. \pageref{15.6}).\\
\endgroup %%------------------------------------<<

Given a triangulated category \bC, let $\Z\bC$ be the additive category whose objects are the ordered pairs 
$(n,X)$ $(n \in \Z, \ X \in \Ob\bC)$, the morphisms from $(n,X)$ to $(m,Y)$ being 
$\underset{q \geq n,n}{\colimx} \Mor(\Sigma^{q-n}X,\Sigma^{q-m}Y)$.  
Composition in $\Z\bC$ comes from composition in 
\bC: $\Sigma^{q-n}X \ra \Sigma^{q-m}Y \ra \Sigma^{q-k}Z$.  
To equip $\Z\bC$ with the structure of a strict triangulated category, take for the suspension functor the isomorphism 
$(n,X) \ra (n-1,X)$ and take for the exact triangles the 
$(n,X) \ra (m,Y)$ 
$\ra (k,Z)$
$\ra (n-1,X)$ 
associated with the
$\Sigma^{q-n}X \overset{u}{\ra} \Sigma^{q-m}Y$ 
$\overset{v}{\ra} \Sigma^{q-k}Z$ 
$\overset{w}{\ra} \Sigma\Sigma^{q-n}X$ such that $(u,v,(-1)^qw)$ is exact.\\

\begin{proposition} \ %08
The functor $F:\bC \ra \Z\bC$ that sends \mX to $(0,X)$ is a triangulated equivalence of categories.
\end{proposition}

%%----------------------------------------------------------------------------------------------09
[Note: \ The natural isomorphism $\Phi:F \circx \Sigma \ra \Sigma \circx F$ is defined by letting 
$\Phi_X:(0,\Sigma X) \ra (-1,X)$ be the canonical image of $\id_{\Sigma X}$ in 
$\Mor((0,\Sigma X),(-1,X))$.]\\

\indent\indent (Octahedral Axiom) \ Let \bC be a triangulated category.  Suppose given exact triangles 
$X \overset{u}{\ra} Y \ra Z^\prime \ra \Sigma X$, 
$Y \overset{v}{\ra} Z \ra X^\prime \ra \Sigma Y$, 
$X \overset{v \circx u}{\ra} Z \ra Y^\prime \ra \Sigma X$
$-$then there exists an exact triangle 
$Z^\prime \ra Y^\prime \ra X^\prime \ra \Sigma Z^\prime$
such that the diagram 
\[
\begin{tikzcd}%[sep=large]
{X} \arrow[d,shift right=0.5,dash] \arrow[d,shift right=-0.5,dash] \ar{r}{u} 
&{Y} \ar{d}{v}  \ar{r} 
&{Z^\prime} \ar{d} \ar{r} 
&{\Sigma X} \arrow[d,shift right=0.5,dash] \arrow[d,shift right=-0.5,dash] \\
{X} \ar{d}[swap]{u} \ar{r}[swap]{v \circx u} 
&{Z} \arrow[d,shift right=0.5,dash] \arrow[d,shift right=-0.5,dash]  \ar{r} 
&{Y^\prime} \ar{d} \ar{r} 
&{\Sigma X} \ar{d} \\
{Y} \ar{d} \ar{r}[swap]{v} 
&{Z} \ar{d}  \ar{r} 
&{X^\prime} \arrow[d,shift right=0.5,dash] \arrow[d,shift right=-0.5,dash] \ar{r} 
&{\Sigma Y} \ar{d} \\
{Z^\prime} \ar{r} 
&{Y^\prime} \ar{r} 
&{X^\prime} \ar{r} 
&{\Sigma Z^\prime}
\end{tikzcd}
\]
commutes.

[Note: \ An explanation for the term ``octahedral'' is the diagram
\[
\begin{tikzcd}%[sep=large]
&{Y^\prime} \arrow[ldd,  red] \ar[dashed]{rd}  \\
{Z^\prime} \arrow[d, red] \ar[dashed]{ru} &&{X^\prime} \arrow[ll, red] \arrow[ldd, red]\\
{X} \ar{dr}[swap]{u}\ar{rr}{u \circx v} &&{Z}  \ar{luu} \ar{u}\\
&{Y} \ar{luu} \ar{ru}[swap]{v} 
\end{tikzcd}
\]
%note - marking command not working to put \circ on the arrows - used color instead
%Here $U \arrowcircle[line width=0.5pt] V$ stands for the arrow $U \ra \Sigma V$.]\\
Here 
\begin{tikzcd} [sep=small]
U \ar[r,red] &V 
\end{tikzcd}
stands for the arrow $U \ra \Sigma V$.]

Example: Let \bA be an abelian category $-$then the triangulated category $\bK(\bA)$ satisfies the octahedral axiom.\\

\begingroup%%----------------------------------->>
\fontsize{9pt}{11pt}\selectfont
The stable homotopy category is a triangulated category satisfying the octahedral axiom.\\
\endgroup %%------------------------------------<<

\begingroup%%----------------------------------->>
\fontsize{9pt}{11pt}\selectfont
\textbf{\small EXAMPLE}  \ 
Let \bC be a triangulated category satisfying the octahedral axiom.  Suppose that \mO is the object class of a triangulated subcategory of \bC and write $S_O$ for the class of morphisms $X \overset{u}{\ra} Y$ which can be completed to an exact triangle 
$X \overset{u}{\ra}$ 
$Y \overset{v}{\ra}$ 
$Z \overset{w}{\ra}$ 
$\Sigma X$
with \mZ in \mO $-$then $S_O$ admits a calculus of left and right fractions.

[$S_O$ contains the identities of \bC ($\forall \ X \in \Ob\bC$, 
$X \overset{\id_X}{\lra} X \ra 0 \ra \Sigma X$ is exact and $0 \in O$).  
To check that $S_O$ is closed under composition, let 
$X \overset{u}{\ra}$ 
$Y \ra$ 
$Z^\prime \ra$ 
$\Sigma X$
and 
$Y \overset{v}{\ra}$ 
$Z \ra X^\prime$ 
$\ra \Sigma Y$
be exact triangles
%%----------------------------------------------------------------------------------------------10
with $Z^\prime$, $X^\prime \in O$.  Choose a completion 
$X \overset{v \circx u}{\lra} Z$ to an exact triangle 
$X \overset{v \circx u}{\lra}$ 
$Z \ra Y^\prime$ 
$\ra \Sigma X$ (cf. TR$_3$) $-$then by the octahedral axiom, there exists an exact triangle 
$Z^\prime \ra Y^\prime$
$\ra X^\prime$
$\ra \Sigma Z^\prime$.  
Since $Z^\prime$, $X^\prime \in O$, it follows from Proposition 6 that $Y^\prime \in O$.  
The remaining verifications do not involve the octahedral axiom.]
\vspi
[Note: \ $S_O$ contains the isomorphisms of \bC.]\\
\endgroup %%------------------------------------<<

\begingroup%%----------------------------------->>
\fontsize{9pt}{11pt}\selectfont
\textbf{\small EXAMPLE}  \ 
Let \bC be a triangulated category satisfying the octahedral axiom.  Given classes $O_1, O_2, \subset \Ob\bC$, denote by 
$O_1*O_2$ the class consisting of those \mX which occur in an exact triangle 
$X_1 \ra X \ra X_2 \ra \Sigma X_1$, ($X_1 \in O_1$, $X_2 \in O_2$) $-$then the octahedral axiom implies that the operation 
$*$ is associative.
\vspi
\label{15.20}
[Note: \ Given a class $O \subset \Ob\bC$, an 
\un{extension}
\index{extension (of a class of objects in a triangulated category)}  
of objects of \mO is an element of $\Extx O = \ds\bigcup\limits_{l \geq 0} O * \cdots * O$ ($l$ factors), the elements 
$O * \cdots * O$ being the extensions of objects of \mO of 
\un{length}
\index{length (of an extension of a class of objects in a triangulated category)} 
$l$.]\\
\endgroup %%------------------------------------<<

\begingroup%%----------------------------------->>
\fontsize{9pt}{11pt}\selectfont
\textbf{\small FACT} \ 
Let \bC be a triangulated category with finite coproducts satisfying the octahedral axiom $-$then, in the notation of TR$_5$, 
$\exists$ an $h:Z \ra Z^\prime$ such that $(f, g, h)$ is a morphism of triangles and the triangle

\begin{tikzcd}[sep=large]
{X^\prime \amalg Y} \ar{rr}{
\begin{pmatrix}
u^\prime \ g\\
0 \ -v\\
\end{pmatrix}
}
&&{Y^\prime \amalg Z}\ar{rr}{
\begin{pmatrix}
v^\prime \ h\\
0 \ -w\\
\end{pmatrix}
}
&&{Z^\prime \amalg \Sigma X}\ar{rr}{
\begin{pmatrix}
w^\prime \ \Sigma f\\
0 \ -\Sigma u\\
\end{pmatrix}
}
&&{\Sigma X^\prime \amalg \Sigma Y}
\end{tikzcd}
is exact.\\
\endgroup %%------------------------------------<<

\begin{proposition} \ %09
Let \bC be a triangulated category satisfying the octahedral axiom $-$then every commutative square 
\begin{tikzcd}%[sep=large]
{X} \ar{d}  \ar{r}  &{Y} \ar{d}\\
{X^\prime} \ar{r} &{Y^\prime}
\end{tikzcd}
can be completed to a diagram 
\[
\begin{tikzcd}%[sep=large]
{X} \ar{d}  \ar{r}  
&{Y} \ar{d}  \ar{r}
&{Z} \ar{d}  \ar{r}
&{\Sigma X} \ar{d}\\
{X^\prime} \ar{d}  \ar{r}  
&{Y^\prime} \ar{d}  \ar{r}
&{Z^\prime} \ar{d}  \ar{r}
&{\Sigma X^\prime} \ar{d}\\
{X\pp} \ar{d}  \ar{r}  
&{Y\pp} \ar{d} \ar{r}
&{Z\pp} \ar{d}{\quad \  \ \un{\ \ }} \ar{r}
&{\Sigma X\pp} \ar{d}\\
{\Sigma X} \ar{r}  
&{\Sigma Y} \ar{r}
&{\Sigma Z} \ar{r} 
&{\Sigma^2 X}
\end{tikzcd}
\]
in which the first three rows and the first three columns are exact and all the squares commute except for the one marked with a minus sign which anticommutes.\\
\end{proposition}

\begingroup%%----------------------------------->>
\fontsize{9pt}{11pt}\selectfont
\textbf{\small EXAMPLE}  \ 
Let \bC be a triangulated category satisfying the octahedral axiom.  Suppose that \mO is the object class of a triangulated subcategory of \bC.  Let 
$X \overset{u}{\ra}$ 
$Y \overset{v}{\ra}$ 
$Z \overset{w}{\ra}$ 
$\Sigma X$
, 
$X^\prime \overset{u^\prime}{\ra}$ 
$Y^\prime \overset{v^\prime}{\ra}$ 
$Z^\prime \overset{w^\prime}{\ra}$ 
$\Sigma X^\prime$
be exact triangles.  Assume: There is a diagram 
\begin{tikzcd}[sep=large]
{X} \ar{d}{f} \ar{r}{u}
&{Y} \ar{d}{g} \ar{r}{v}
&{Z} \ar{r}{w}
&{\Sigma X}\ar{d}{\Sigma f}\\
{X^\prime} \ar{r}[swap]{u^\prime}
&{Y^\prime} \ar{r}[swap]{v^\prime}
&{Z^\prime} \ar{r}[swap]{w^\prime}
&{\Sigma X^\prime}
\end{tikzcd}
, where $f, g \in S_O$ and $g \circx u = u^\prime \circx f$ $-$then $\exists$ an $h:Z \ra Z^\prime$ in $S_O$ such that 
$(f, g, h)$ is a morphism of triangles.
\vspi
%%----------------------------------------------------------------------------------------------11
[Note: \ The metacategory $S_O^{-1}\bC$ is triangulated and satisfies the octahedral axiom.  
For instance, consider $\bK(\bA)$, where \bA is an abelian category.  
Let $O = \{X: H^n(X) = 0 \ \forall \ n\}$ $-$then $S_O$ is the class of 
\un{quasiisomorphisms}
\index{quasiisomorphisms (abelian category)} 
of \bA 
(i.e., the $f$ such that $H^n(f)$ is an isomorphism 
$\forall$ $n$ or, equivalently, the $f$ such that $H^n(C_f) = 0 \ \forall \ n$) and the 
\un{derived category}
\index{derived category} 
$\bD(\bA)$ of \bA is the localization $S_O^{-1}\bK(\bA)$.  
But there is a problem with the terminology.  
Reason: A priori, 
$\bD(\bA)$ is only a metacategory.  
However, the assumption that \bA is Grothendieck and has a separator suffices to ensure that 
$\bD(\bA)$ is a category 
(Weibel\footnote[2]{\textit{An Introduction to Homologial Algebra}, Cambridge University Press (1994), 386-387.}
One can also form $\bD^+(\bA)$, $\bD^-(\bA)$ and $\bD^\text{b}(\bA)$.  
Here 
$\bD^+(\bA)$ will be a category if \bA has enough injectives and 
$\bD^-(\bA)$ will be a category if \bA has enough projectives.]\\
\endgroup %%------------------------------------<<

\begingroup%%----------------------------------->>
\fontsize{9pt}{11pt}\selectfont
The derived category $\bD(\bA)$  of 
Freyd's\footnote[3]{\textit{Abelian Categories}, Harper $\&$ Row (1964), 131-132.}
``large'' abelian category \bA is not isomorphic to a category, hence exists only as a metacategory.  
Therefore one cannot find a model category structure on \bA whose weak equivalences are the quasiisomorphisms (cf. p. \pageref{15.7}).\\
\endgroup %%------------------------------------<<

Let \bC be a triangulated category $-$then a subcategory \bD of \bC is said to be 
\un{thick}
\index{thick (subcategory)} 
provided that it is triangulated and for any pair of morphisms 
$i:X \ra Y$, $r:Y \ra X$ with $r \circx i = \id_X$, $Y \in \Ob\bD$ $\implies$ $X \in \Ob\bD$.\\

\begin{proposition} \ %10
Let \bC be a triangulated category with finite coproducts $-$then a triangulated subcategory \bD of \bC is thick iff every object of \bC which is a direct summand of an object of \bD is itself an object of \bD, i.e., $Y \in \Ob\bD$ $\&$ 
$Y \approx  X \coprod Z$ $\implies$ $X \in \Ob\bD$.
\end{proposition}

[Necessity: Since \bD is isomorphism closed, 
$X \coprod Z \in \Ob\bD$, so one only has to consider
\begin{tikzcd}%[sep=large]
{X}  \ar{r}{\ini_X}  &{X \coprod Z} \ar{r}{\pr_X} &{X.}
\end{tikzcd}

Sufficiency: There exists an isomorphism $Y \ra X \amalg Z$ and a commutative diagram 
\begin{tikzcd}%[sep=large]
{X} \ar{dr}  \ar{r}{i}  &{Y} \ar{d}\\
&{X \amalg Z}
\end{tikzcd}
(cf. p. \pageref{15.8}), hence $X \in \Ob\bD$.]
\vspace{0.5cm}

\begin{proposition} \ %11
Let \bC be a triangulated category with finite coproducts satisfying the octahedral axiom $-$then a triangulated subcategory \bD of \bC is thick iff every morphism $X \overset{u}{\ra} Y$ in \bC admitting a factorization 
\begin{tikzcd}[sep=small]
{X} \ar{ddr}[swap]{\phi}  \ar{rr}{u}  &&{Y}\\
\\
&{W} \ar{ruu}[swap]{\psi}
\end{tikzcd}
through an object \mW of \bD and contained in an exact triangle 
$X \overset{u}{\ra} Y$ 
$\overset{v}{\ra} Z$ 
$\overset{w}{\ra} \Sigma X$, where
$Z \in \Ob\bD$, is a morphism in \bD, i.e., $X, \ Y \in \bD$.
\end{proposition}

%%----------------------------------------------------------------------------------------------12
[Necessity: Complete $X \overset{\phi}{\ra} W$ to an exact triangle 
$X \overset{\phi}{\ra} W$ 
$\overset{\omega}{\ra} W^\prime$ 
$\overset{\omega^\prime}{\ra} \Sigma X$ (cf. TR$_3$) $-$then the triangle 
\begin{tikzcd}%[sep=large]
{Y \amalg W} \ar{rr}{
\begin{pmatrix}
\id_Y \ \psi\\
0 \ \  -\omega\\
\end{pmatrix}
}
&&{Y \amalg W^\prime} \ar{r}{(0 - \omega^\prime)}
&{\Sigma X} \ar{r}{
\begin{pmatrix}
\Sigma u\\
-\Sigma \phi\\
\end{pmatrix}
}
&{\Sigma Y \amalg \Sigma W}
\end{tikzcd}
is exact (cf. p. \pageref{15.9} ff.), thus the triangle 
\begin{tikzcd}%[sep=large]
{X}  \ar{r}{
\begin{pmatrix}
-u\\
\phi\\
\end{pmatrix}
}
&{Y \amalg W} \ar{rr}{
\begin{pmatrix}
\id_Y \ \psi\\
0 \ \  -\omega\\
\end{pmatrix}
}
&&{Y \amalg W^\prime} \ar{r}{(0 - \omega^\prime)}
&{\Sigma X}
\end{tikzcd} 
is exact (cf. TR$_4$).  
On the other hand, the triangle 
$Y \amalg W \ra Y$ 
$\overset{0}{\ra} \Sigma W$ 
$\ra \Sigma Y \amalg \Sigma W$ 
is exact 
(cf. p. \pageref{15.10}), %<<-------------
as is the triangle 
$X \overset{-u}{\ra} Y$ 
$\overset{-v}{\ra} Z$ 
$\overset{w}{\ra} \Sigma X$ 
(cf. p. \pageref{15.11}).  %<<-------------
So, in the notation of the octahedral axiom, taking 
$Z^\prime = Y \amalg W^\prime$, 
$X^\prime = \Sigma W$, and 
$Y^\prime = Z$, 
one concludes that there is an exact triangle
$Y \amalg W^\prime \ra Z$ 
$\ra \Sigma W$ 
$\ra \Sigma Y \amalg \Sigma W^\prime$.  
But $Z, \Sigma W \in \Ob\bD$ $\implies$ $Y \amalg W^\prime \in \Ob\bD$ (cf. Proposition 6) 
$\implies$ $Y \in \Ob\bD$ (cf. Proposition 10) 
$\implies$ $X \in \Ob\bD$ (cf. Proposition 6).

Sufficiency: Suppose that $Y \in \Ob\bD$ $\&$ $Y \approx X \amalg Z$ $-$then the triangle
$X \ra X \amalg Z$ 
$\ra Z$ 
$\overset{0}{\ra} \Sigma X$ is exact 
(cf. p. \pageref{15.12}), %<<-----------------------
thus the triangle 
$\Omega Z \overset{0}{\ra}  X$ 
$\ra X \amalg Z$ 
$\ra \Sigma \Omega Z $ is exact 
(cf. p. \pageref{15.13}).  %<<-----------------------
But $0 \in \Ob\bD$ and there is a factorization 
\begin{tikzcd}[sep=small]
{\Omega Z}\ar{ddr} \ar{rr}{0} &&{X}\\
\\
&{0} \ar{ruu}
\end{tikzcd} 
.  
Our assumption implies that $X \in \Ob\bD$, so \bD is thick (cf. Proposition 10).]\\

\begingroup%%----------------------------------->>
\fontsize{9pt}{11pt}\selectfont
\textbf{\small FACT} \ 
Let \bC be a triangulated category with finite coproducts satisfying the octahedral axiom.  
Suppose that \mO is the object class of a thick subcategory of \bC 
$-$then $u \in S_O$ iff $\exists$ $f, g \in \Mor \bC$: $u \circx f \in S_O$, $g \circx u \in S_O$.
\vspi
[Complete $X \overset{u}{\ra} Y$ to an exact triangle 
$X \overset{u}{\ra} Y$ 
$\overset{v}{\ra} Z$ 
$\overset{w}{\ra} \Sigma X$ (cf. TR$_3$), the claim being that $Z \in O$.  
By hypothesis, there are exact triangles 
$X^\prime \overset{u \circx f}{\lra} Y$ 
$\ra Z_f$ 
$\ra \Sigma X^\prime$, 
$X \overset{g \circx u}{\lra} Y^\prime$ 
$\ra Z_g$ 
$\ra \Sigma X$
where $Z_f,\ Z_g \in O$.  
Since $v \circx (u \circx f) = (v \circx u) \circx f = 0$ (cf. Proposition 3) and 
$\Mor(Z_f,Z) \ra $ 
$\Mor(Y,Z) \ra $ 
$\Mor(X^\prime,Z)$ 
is exact (cf. Proposition 4), $\exists$ a factorization 
\begin{tikzcd}[sep=small]
{Y} \ar{ddr} \ar{rr}{v} &&{Z}\\
\\
&{Z_f} \ar{ruu}
\end{tikzcd} 
.  
Complete $Y \overset{g}{\ra} Y^\prime$ to an exact triangle 
$Y \overset{g}{\ra} Y^\prime$ 
$\ra W$ 
$\ra \Sigma Y$ (cf. TR$_3$) and use the octahedral axiom on
$X \overset{u}{\ra} Y$ 
$\ra Z$ 
$\ra \Sigma X$, 
$Y \overset{g}{\ra} Y^\prime$ 
$\ra W$ 
$\ra \Sigma Y$, 
$X \overset{g\circx u}{\lra} Y^\prime$ 
$\ra Z_g$ 
$\ra \Sigma X$, 
to get an exact triangle 
$Z \ra \Sigma Z_g$ 
$\ra W$ 
$\ra \Sigma Z$, 
or still, an exact triangle
$W \ra \Sigma Z$ 
$\ra \Sigma Z_g$ 
$\ra \Sigma W$.  
From the above, the arrow $W \ra \Sigma Z$ factors through $\Sigma Z_f \in O$.  But also $\Sigma Z_g \in O$, thus, as \mO 
is thick, $\Sigma Z \in O$ (cf. Proposition 11), i.e., $Z \in O$.]
\vspi
[Note: \ The condition implies that $S_O$ is saturated: $S_O = \ov{S}_O$ 
(cf. p. \pageref{15.14}), hence 
$X \in O$ iff $L_{S_O}X$ is a zero object.]\\
\endgroup %%------------------------------------<<

\label{15.37}
\begingroup%%----------------------------------->>
\fontsize{9pt}{11pt}\selectfont
Given a triangulated category \bC, call a class $S \subset \Mor \bC$ 
\un{multiplicative}
\index{multiplicative (class of morphisms)} if 
(1) \ \mS admits a calculus of left and right fractions and contains the isomorphisms of \bC; 
(2) \  $u \in S$ $\implies$ $\Sigma u$ $\&$ $\Omega u \in S$; 
(3) \ $f, g \in S$ $\implies$ $\exists \ h \in S$ (data as in TR$_5$); 
(4) \ $u \in S$ iff $\exists$ $f, g \in \Mor \bC$: $u \circx f \in S$, $g \circx u \in S$.
\vspi
%%----------------------------------------------------------------------------------------------13
Example: Let \bC be a triangulated category with finite coproducts satisfying the octahedral axiom $-$then $S_O$ is multiplicative provided that \mO is the object class of a thick subcategory of \bC.  
In fact, the assignment $O \ra S_O$ 
establishes a one-to-one correspondence between the object class of thick subcategories of \bC and the multiplicative classes of morphisms of \bC.
\vspi
[Note: \ To place this conclusion in perspective, recall that in an abelian category there is a one-to-one correspondence between the Serre classes and the saturated morphism classes which admit a calculus of left and right fractions 
(Schubert\footnote[2]{\textit{Categories}, Springer Verlag (1972), 276.}).]\\
\endgroup %%------------------------------------<<

\begin{proposition} \ %12
Let \bC be a triangulated category.  Assume: \bC has coproducts $-$then for any collection 
$\{X_i \ra Y_i \ra Z_i \ra \Sigma X_i\}$ of exact triangles, the triangle 
$\coprod\limits_i X_i \ra \coprod\limits_i Y_i \ra \coprod\limits_i Z_i \ra \coprod\limits_i \Sigma X_i$
is exact.
\end{proposition}

[Note: \ The suspension functor preserves coproducts, so 
$\Sigma \coprod\limits_i X_i \approx \coprod\limits_i \Sigma X_i$.]\\

Let \bC be a triangulated category with coproducts $-$then an $X \in \Ob\bC$ is said to be 
\un{compact}
\index{compact (object in a triangulated category)} 
if $\forall$ collection $\{X_i\}$ of objects in \bC, the arrow 
$\bigoplus\limits_i \Mor(X,X_i) \ra \Mor(X,\coprod\limits_i X_i)$ is an isomorphism.

[Note: \  \mX compact $\implies$ $\Sigma X$ $\&$ $\Omega X$ compact.]\\

\begingroup%%----------------------------------->>
\fontsize{9pt}{11pt}\selectfont
\label{15.32}
\textbf{\small EXAMPLE}  \ 
Let \mA be a commutative ring with unit $-$then the compact objects in $\bD(\bAMOD)$ are those objects which are isomorphic to bounded complexes of finitely generated projective $A$-modules  
(B\"okstedt-Neeman\footnote[3]{\textit{Compositio Math.} \textbf{86} (1993), 209-234.}).\\
\endgroup %%------------------------------------<<

\begingroup%%----------------------------------->>
\fontsize{9pt}{11pt}\selectfont
\textbf{\small FACT} \ 
If \bC is a triangulated category with coproducts, then the class of compact objects in \bC is the object class of a thick subcategory of \bC.\\
\endgroup %%------------------------------------<<

Notation: Let \bC be a triangulated category with coproducts.  Suppose given an object $(\bX,\bff)$ in $\bFIL(\bC)$ 
$-$then tel$(\bX,\bff)$ is any completion of 
%\begin{tikzcd}%[ sep=small]
%{\coprod\limits_n X_n} \ar{r}{\text{sf}} &{\coprod\limits_n X_n}
%\end{tikzcd}
$\coprod\limits_n X_n \overset{\text{sf}}{\lra} \coprod\limits_n X_n$ 
to an exact triangle (cf. TR$_3$), the $n^{\text{th}}$ component of sf being the arrow
\begin{tikzcd}%[sep=large]
{X_n} \ar{rr}{
\begin{pmatrix}
\id \\
-f_n\\
\end{pmatrix}
}
&&{X_n \coprod\limits_n X_{n+1}}
\end{tikzcd}
.
\vspace{0.25cm}

\begin{proposition} \ %13
Let \bC be a triangulated category with coproducts.  Fix an $(\bX,\bff)$ in $\bFIL(\bC)$ 
$-$then $\forall$ compact \mX, the arrow 
$\colimx \Mor(X,X_n) \ra \Mor(X,\tel(\bX,\bff))$ is an isomorphism.
\end{proposition}

%%----------------------------------------------------------------------------------------------14
[First consider the exact sequence 
$\Mor(X,\coprod\limits_n X_n) \overset{\Phi}{\lra}$ 
$\Mor(X,\tel(\bX,\bff)) \ra $
$\Mor(X,$ $\coprod\limits_n \Sigma X_n) \ra$
$\Mor(X,\coprod\limits_n \Sigma X_n)$ 
(cf. Proposition 4).  
Due to the compactness of \mX, in the commutative diagram
\begin{tikzcd}%[sep=large]
{\bigoplus\limits_n \Mor(X,\Sigma X_n)} \ar{d} \ar{r} &{\bigoplus\limits_n \Mor(X,\Sigma X_n)} \ar{d}\\
{\Mor(X,\coprod\limits_n \Sigma X_n)} \ar{r} &{\Mor(X,\coprod\limits_n \Sigma X_n)}
\end{tikzcd}
, 
the vertical arrows are isomorphisms.  
Because the horizontal arrow on the top is injective, the same holds for the horizontal arrow on the bottom.  
Therefore $\Phi$ is surjective.  
Now write down the commutative diagram
\[
\begin{tikzcd}%[sep=large]
{\bigoplus\limits_n \Mor(X,X_n)} \ar{d} \ar{r}{\phi} &{\bigoplus\limits_n \Mor(X,X_n)} \ar{d}\\
{\Mor(X,\coprod\limits_n  X_n)} \ar{r} 
&{\Mor(X,\coprod\limits_n  X_n)} \ar{r}{\Phi}
&{\Mor(X,\tel(\bX,\bff))} \ar{r}
&{0}
\end{tikzcd}
\]
and observe that $\colimx \Mor(X,X_n)$ can be identified with the cokernel of $\phi$.]\\

\begingroup%%----------------------------------->>
\fontsize{9pt}{11pt}\selectfont
\textbf{\small FACT} \  
Let \bC be a triangulated category with coproducts.  Fix an $(\bX,\bff)$ in$\bFIL(\bC)$ $-$then $\forall$ \mY, there is an exact sequence 
$0 \ra$ 
$\lim^1 \Mor(\Sigma X_n,Y) \ra$ 
$\Mor(\tel(\bX,\bff),Y) \ra$ 
$\lim \Mor(X_n,Y) \ra 0$.\\
\endgroup %%------------------------------------<<

A triangulated category \bC is said to be 
\un{compactly generated}
\index{compactly generated (triangulated category)} 
if it has coproducts and \Ob\bC contains a set $\sU = \{U\}$ of compact objects such that $\Mor(U,X) = 0$ $\forall \ U \in \sU$ 
$\implies$ $X = 0$.

[Note: \ The \un{closure} $\ov{\sU} = \{\ov{U}\}$ of $\sU$ is the set 
$\bigcup\limits_U \{\Sigma^n U: n \geq 0\} \cup \bigcup\limits_U \{\Omega^n U: n \geq 0\}$.]\\

\begingroup%%----------------------------------->>
\fontsize{9pt}{11pt}\selectfont
The stable homotopy category is a compactly generated triangulated category.\\
\endgroup %%------------------------------------<<

\begingroup%%----------------------------------->>
\fontsize{9pt}{11pt}\selectfont
\textbf{\small EXAMPLE}  \ 
Let \mX be a scheme, $\sO_X$ its structure sheaf.  
Denote by $\sO_X\text{-}\bMOD)$ 
\index{$\sO_X\text{-}\bMOD)$}
the category of $\sO_X$-modules and write $\bQ\bC/X$ for the full subcategory whose objects are quasicoherent $-$then 
$\sO_X$-\bMOD and $\bQ\bC/X$ are abelian categories and the inclusion 
$\bQ\bC/X \ra \sO_X\text{-}\bMOD$ is exact.  
In addition, $\sO_X\text{-}\bMOD$ is Grothendieck and has a separator, thus the derived 
category $\bD(\sO_X\text{-}\bMOD)$ exists.  
When \mX is quasicompact (= compact) and separated, $\bQ\bC/X$ is Grothendieck and has a separator, thus in this situation, the derived category $\bD(\bQ\bC/X)$ also exists.  
Moreover, $\bD(\bQ\bC/X)$ is compactly generated, the compact objects being those objects which are isomorphic to perfect complexes 
(Neeman\footnote[2]{\textit{J. Amer. Math. Soc.} \textbf{9} (1996), 205-236.}).\\
\endgroup %%------------------------------------<<

\label{17.6}
\index{Theorem: Brown Representability Theorem}
\index{Brown Representability Theorem}
\textbf{\small BROWN REPRESENTABILITY THEOREM} \quad
Let \bC be a compactly generated triangulated category  $-$then an exact cofunctor $F:\bC \ra \bAB$ is representable iff it converts coproducts into products.

%%----------------------------------------------------------------------------------------------15
[The condition is clearly necessary and the proof of the sufficiency is a variation on the argument used in Proposition 27 of $\S 5$.  Thus setting 
$X_0 = \coprod\limits_{\ov{U}} F\ov{U} \cdot \ov{U}$, one has 
$FX_0 = \coprod\limits_{\ov{U}} (F\ov{U})^{F\ov{U}}$.  
Call $\xi_0$ that element of the product defined by 
$\xi_{0,\ov{U}} = \id_{F\ov{U}}$ $\forall \ \ov{U}$ and let 
$\Xi_0:\Mor(-,X_0) \ra F$ be the natural transformation associated with $\xi_0$ via Yoneda.  
Note that 
$\Xi_{0,\ov{U}}:\Mor(\ov{U},X_0) \ra F\ov{U}$ is surjective $\forall \ \ov{U}$.  
Proceeding inductively, we shall construct an object $(\bX,\bff)$ in $\bFIL(\bC)$ and natural transformations 
$\Xi_n:\Mor(-,X_n) \ra F$ such that $\forall \ n$, the triangle 
\begin{tikzcd}%[sep=large]
{\Mor(-,X_n)} \ar{d} \ar{rd}{\Xi_n}\\
{\Mor(-,X_{n+1})} \ar{r}[swap]{\Xi_{n+1}} &{F}
\end{tikzcd}
commutes.  
To this end, put
$K_n =\ds \coprod\limits_{\ov{U}} (\ker \Xi_{n,\ov{U}}) \cdot \ov{U}$ 
and complete the canonical arrow $K_n \ra X_n$ to an exact triangle
$K_n \ra X_n$ 
$\overset{f_n}{\lra} X_{n+1}$ 
$\ra \Sigma K_n$ (cf. TR$_3$).  
If $\xi_n \in FX_n$ corresponds to $\Xi_n$, then $\xi_n \in \ker(FX_n \ra FK_n)$ and since the sequence 
$FX_{n+1} \ra FX_n \ra FK_n$ is exact, $\exists$ $\xi_{n+1} \in FX_{n+1}: \xi_{n+1} \ra \xi_n$.  
Definition: $\Xi_{n+1} \leftrightarrow \xi_{n+}$, which finishes the induction.  
Abbreviating tel(\bX,\bff) to $X_\omega$, there is a natural transformation $\Xi_\omega:\Mor(-,X_\omega) \ra F$ rendering the triangle 
\begin{tikzcd}%[sep=large]
{\Mor(-,X_n)} \ar{d} \ar{rd}{\Xi_n}\\
{\Mor(-,X_\omega)} \ar{r}[swap]{\Xi_\omega} &{F}
\end{tikzcd}
commutative $\forall \ n$.  Proof: Consider the diagram
\begin{tikzcd}%[sep=large]
{FX_\omega} \ar{r} 
&{F(\coprod\limits_n X_n)} \arrow[d,shift right=0.5,dash] \arrow[d,shift right=-0.5,dash] \ar{r}{F\text{sf}}
&{F(\coprod\limits_n X_n)} \arrow[d,shift right=0.5,dash] \arrow[d,shift right=-0.5,dash]\\
&{\prod\limits_n FX_n} \ar{r}
&{\prod\limits_n FX_n}
\end{tikzcd}
.  
Because $\prod\limits_n \xi_n$ lies in the kernel of $\prod\limits_n FX_n \ra \prod\limits_n FX_n$, exactness gives a 
$\xi_\omega \in FX_\omega:\xi_\omega \ra \prod\limits_n \xi_\omega$, hence 
$\Xi_\omega \leftrightarrow \xi_\omega$
has the stated property.  
The final step is to establish that 
$\Xi_{\omega,X}:\Mor(X,X_\omega) \ra FX$ is bijective $\forall \ X$.  But it is certainly true that 
$\Xi_{\omega,\ov{U}}$ is bijective $\forall \ \ov{U}$ (injectivity follows from the construction of $X_\omega$ 
(cf. Proposition 13)) while $\Mor(\ov{U},X_0) \ra F\ov{U}$ surjective $\implies$
$\Mor(\ov{U},X_\omega) \ra F\ov{U}$ surjective and this turns out to be enough (cf. infra).]\\
%\vspace{0.25cm}

The assumption that $\Mor(U,X) = 0$ $\forall \ U \in \sU$ $\implies$ $X = 0$ has yet to be employed.  
To do so, let $\bC_F$ be the full, isomorphism closed subcategory of \bC whose objects are those \mX such that 
$\Xi_{\omega,\Sigma^n X}:\Mor(\Sigma^n X,X_\omega) \ra F\Sigma^n X$ is bijective $\forall \ n \geq 0$ and 
$\Xi_{\omega,\Omega^n X}:\Mor(\Omega^n X,X_\omega) \ra F\Omega^n X$  is bijective $\forall \ n \geq 0$.  
Obviously, $\bC_F$ contains 0 and $\sU$.

Claim: $\bC_F$ is stable under $\Sigma$ $\&$ $\Omega$.

[To check stability under $\Sigma$, fix an $X \in \Ob\bC_F$ $-$then $\forall \ n \geq 0$, 
$\Mor(\Sigma^n\Sigma X,X_\omega)$ = $\Mor(\Sigma^{n+1}X,X_\omega)$ 
$\approx F\Sigma^{n+1} X =$ $F\Sigma^n\Sigma X$.  On the other hand, the arrow of adjunction 
$X \ra \Omega\Sigma X$ is an isomorphism, thus one sees inductively from the commutative diagram
%%----------------------------------------------------------------------------------------------16
\begin{tikzcd}%[sep=large]
{\Mor(\Omega^n\Sigma X,X_\omega)} \ar{d} \ar{r} 
&{\Mor(\Omega^{n-1}\Sigma X,X_\omega)} \ar{d}\\
{F\Omega^n\Sigma X} \ar{r} &{F\Omega^{n-1} X}
\end{tikzcd}
that $\Mor(\Omega^n\Sigma X,X_\omega) \approx F\Omega^n \Sigma X$.  
Therefore $\Sigma X \in \bC_F$.]

Claim: If $X \overset{u}{\ra} Y \overset{v}{\ra} Z \overset{w}{\ra} \Sigma X$ is an exact triangle with 
$X,Y \in \Ob\bC_F$, then $Z \in \Ob\bC_F$.

[Use the five lemma.]

Claim: $\bC_F$ is closed under the formation of coproducts in \bC.

[E.g.: 
$\Mor(\Sigma^n \coprod\limits_i X_i,X_\omega) \approx$ 
$\Mor(\coprod\limits_i X_i,\Omega^nX_\omega) \approx$ 
$\prod\limits_i\Mor(X_i,\Omega^nX_\omega) \approx$ 
$\prod\limits_i\Mor(\Sigma^n X_i,$ $X_\omega) \approx$ 
$\prod\limits_i F\Sigma^n X_i \approx$ 
$F(\coprod\limits_i \Sigma^n X_i) \approx$ 
$F(\Sigma^n \coprod\limits_i X_i)$.]

In summary, $\bC_F$ is a triangulated subcategory of \bC containing $\sU$ and closed under the formation of coproducts in \bC.  
To conclude that 
$\Xi_{\omega,X}:\Mor(X,X_\omega) \ra FX$ 
is bijective $\forall \ X$, it need only be shown that $\bC_F = \bC$, which is a special case of the following result.\\

\begin{proposition} \ %14
Let \bC be a compactly generated triangulated category.  
Suppose that \bD is a triangulated subcategory of \bC containing $\sU$ and closed under the formation of coproducts in \bC $-$then $\bD = \bC$.
\end{proposition}

[Let $\ov{\bD}$ be the smallest triangulated subcategory of \bC containing $\sU$ and closed under the formation of coproducts in \bC.  
Fix an \mX in \bC $-$then the restriction of $\Mor(-,X)$ to $\ov{\bD}$ is an exact cofunctor.  
Applying what has been proved above about Brown representability to $\Mor(-,X)$ one concludes that there exists an 
$\ov{X}_\omega$ in $\ov{\bD}$ and a natural isomorphism $\Mor(-,\ov{X}_\omega) \ra \Mor(-,X)$ (the minimality of 
$\ov{\bD}$ enters the picture at this point).  
Accordingly, $\exists$ a morphism $\ov{X}_\omega \ra X$ such that 
$\forall \ \ov{X}$ in $\ov{\bD}$, the arrow $\Mor(\ov{X},\ov{X}_\omega) \ra \Mor(\ov{X},X)$ is bijective.  Complete 
$\ov{X}_\omega \ra X$ to an exact triangle 
$\ov{X}_\omega \ra X \ra$ $Y \ra \Sigma \ov{X}_\omega$ in \bC (cf. TR$_3$) $-$then $\forall \ \ov{X} \in \ov{\bD}$, 
$\Mor(\ov{X},Y) = 0$ $\implies$ $\forall \ U \in \sU$, $\Mor(U,Y) = 0$ $\implies$ $Y = 0$.  
Consequently, the morphism $\ov{X}_\omega \ra X$ is an isomorphism 
(cf. p. \pageref{15.15}), 
so 
$X \in \Ob\ov{\bD}$ ($\ov{\bD}$ is isomorphism closed), hence $\ov{\bD} = \bC$ $\implies$ $\bD = \bC$.]\\

\label{15.19}
\label{15.25}
Application: Let \bC be a compactly generated triangulated category.  
Suppose that 
$\Xi:\Mor(-,Y) \ra \Mor(-,Z)$ is a natural transformation such that $\forall$ $\ov{U} \in \ov{\sU}$, $\Xi_{\ov{U}}$ \  is bijective 
$-$then for all \mX in \bC, $\Xi_X:\Mor(X,Y) \ra \Mor(X,Z)$ is bijective.

[Note: \ If $\Xi_f:\Mor(-,Y) \ra \Mor(-,Z)$ is the natural transformation corresponding to $f:Y \ra Z$, then $f$ is an isomorphism whenever $\Xi_{f,\ov{U}}$ is bijective $\forall \ \ov{U} \in \ov{\sU}$.]\\

\label{15.27}
Example: Suppose that \bC be a compactly generated triangulated category.  Let $\Delta:\bI \ra \bC$ be a diagram 
$-$then a weak colimit \mL of $\Delta$ is said to be a 
\un{minimal weak colimit}
\index{minimal weak colimit} 
%%----------------------------------------------------------------------------------------------17
provided that $\forall \ \ov{U} \in \ov{\sU}$, 
$\colimx \Mor(\ov{U},\Delta_i) \approx \Mor(\ov{U},L)$.  If \mL is a minimal weak colimit of $\Delta$ and if \mK is an arbitrary weak colimit of $\Delta$, then there are arrows
$L \overset{\phi}{\ra} K \overset{\psi}{\ra} L$
and $\forall \ \ov{U} \in \ov{\sU}$, $\Xi_{\psi \circx \phi,\ov{U}}:\Mor(\ov{U},L) \ra \Mor(\ov{U},L)$ is bijective, thus by the above, $\psi \circx \phi$ is an isomorphism.  
Corollary: \mL is a direct summand of \mK 
(cf. p. \pageref{15.16}.

[Note: \ \mL $\&$ \mK minimal $\implies$ $L \approx K$.  
Example: $\forall \ (\bX,\bff)$ in $\bFIL(\bC)$, $\telsub(\bX,\bff)$ is a minimal weak colimit of (\bX,\bff) (cf. Proposition 13).]\\

\label{17.16}
\begingroup%%----------------------------------->>
\fontsize{9pt}{11pt}\selectfont
\label{17.28}
\textbf{\small EXAMPLE}  \ 
Suppose that \bC is a compactly generated triangulated category. 
Fix a compact object \mX $-$then for any divisible abelian group \mA, $\Hom(\Mor(X,-),A)$ is an exact cofunctor which converts coproducts into products, thus is representable.\\
\endgroup %%------------------------------------<<

\label{17.53}
\label{17.63}
\index{idempotents split (Example)}
\begingroup%%----------------------------------->>
\fontsize{9pt}{11pt}\selectfont
\textbf{\small EXAMPLE \ (\un{Idempotents Split})} \ 
Suppose that \bC is a compactly generated triangulated category.  Let $e \in \Mor(Y,Y)$ be idempotent $-$then $\exists$ \mX, \mZ and an isomorphism $Y \ra X \amalg Z$ such that the diagram 
\begin{tikzcd}[sep=large]
{Y} \ar{d} \ar{rr}{e} &&{Y} \ar{d}\\
{X \amalg Z} \ar{r} &{X} \ar{r} &{X \amalg Z} 
\end{tikzcd}
commutes.
\vspi
[Using suggestive notation, write $\Mor(-,Y)$ as a direct sum 
$e\Mor(-,Y) \oplus (1 - e)\Mor(-Y)$ of two exact cofunctors which convert coproducts into products and choose \mX,\mZ:
$e\Mor(-,Y) \approx \Mor(-,X)$, 
$(1 - e)\Mor(-,Y) \approx \Mor(-,Z)$.]
\vspi
[Note: \ Defining $r:Y \ra X$ and $i:X \ra Y$ in the obvious way, one has $e = i \circx r$ and $r \circx i = \id_X$.  
Moreover $r:Y \ra X$ is a split coequalizer of $e$, $\id_Y:Y \ra Y$, as can bee seen from the diagram 
\[
\begin{tikzcd}[sep=large]
{Y}     
\ar[shift left=1.5]{rr}{e}
\ar[shift right=1.5]{rr}[swap]{\id_Y}
&&{Y} \ar{rr}{r} \ar{rrd}[swap]{e} \ar[bend right]{ll}[swap]{\id_Y}
&&{X}\ar{d}{i} \ar[bend right]{ll}[swap]{i}\\
&&&&{Y}
\end{tikzcd}
.]
\]
\\
\endgroup %%------------------------------------<<

\index{The Eilenberg Swindle (Example)}
\begingroup%%----------------------------------->>
\fontsize{9pt}{11pt}\selectfont
\textbf{\small EXAMPLE \  (\un{The Eilenberg Swindle})} \ 
Suppose that \bC is a compactly generated triangulated category.  Let \bD be a triangulated subcategory of \bC.  
Assume: \bD is closed under the formation of coproducts in \bC $-$then \bD is thick.
\vspi
[Fix a pair of morphisms $i:X \ra Y$, $r:Y \ra X$ with $r \circx i = \id_X$ and $Y \in \Ob\bD$.  
Put $e = i \circx r$.  
Since $e$ is an idempotent, by the preceding example $Y \approx X \amalg Z$ for some \mZ.  
Write 
$W = X \amalg (Z \amalg X) \amalg (Z \amalg X) \amalg \cdots \approx$ 
$(X \amalg Z) \amalg (X \amalg Z) \amalg \cdots$ 
to get $W \in \Ob\bD$.  
But 
$W \approx $ 
$X \amalg W \approx $
$W \amalg X$ $\implies$ 
$W \amalg X \in \Ob\bD$.  
Because the triangle 
$W \ra W \amalg X \ra X \overset{0}{\ra} \Sigma W$ is exact 
(cf. p. \pageref{15.17} ff.), it follows that $X \in \Ob\bD$.]\\
\endgroup %%------------------------------------<<

\begingroup%%----------------------------------->>
\fontsize{9pt}{11pt}\selectfont
\textbf{\small EXAMPLE}  \ 
Suppose that \bC is a compactly generated triangulated category $-$then \bC has products.  Proof: Given a set of objects 
$X_i$, apply the Brown representability theorem to the exact cofunctor 
$Y \ra \prod\limits_i \Mor(Y,X_i)$.
\vspi
%%----------------------------------------------------------------------------------------------18
[Note: \ The morphism \ 
$t:\ds\coprod\limits_i X_i \ra \prod\limits_i X_i$ of p. 
\pageref{15.18} is an isomorphism iff \ $\forall \ \ov{U} \in \ov{\sU}$ : 
$\#\{i:\Mor(\ov{U},X_i) \neq 0\} < \omega$.  To see this, consider the arrow 
$\Mor(\ov{U},\ds\coprod\limits_i X_i) = \ds\bigoplus\limits_i \Mor(\ov{U},X_i) 
\ra $ $\ds\prod\limits_i \Mor(\ov{U},X_i)  = \Mor(\ov{U},\prod\limits_i X_i)$.]\\
\endgroup %%------------------------------------<<

\begin{proposition} \ %15
Let \bC be a compactly generated triangulated category and let \bD be an arbitrary triangulated category.  
Suppose that 
$F:\bC \ra \bD$ is a triangulated functor which preserves coproducts $-$then $F$ has a right adjoint $G:\bD \ra \bC$.
\end{proposition}

[Given a $Y \in \Ob\bD$, the cofunctor $X \ra \Mor(FX,Y)$ is exact and converts coproducts into products, thus is representable: $\Mor(F-,Y) \approx \Mor(-,GY)$.]\\

\begingroup%%----------------------------------->>
\fontsize{9pt}{11pt}\selectfont
\textbf{\small FACT} \ 
Let \bC be a compactly generated triangulated category and let \bD be an arbitrary triangulated category.  Suppose that 
$F:\bC \ra \bD$ is a triangulated functor which preserves coproducts $-$then its right adjoint $G:\bD \ra \bC$ preserves coproducts iff $\forall \ U \in \sU$, $FU$ is compact.
\vspi
[Necessity: \ 
$\ds\bigoplus\limits_j \Mor(FU,Y_j) \approx$ 
$\ds\bigoplus\limits_j \Mor(U,GY_j) \approx$ 
$\Mor(U,\ds\coprod\limits_j GY_j) \approx$ 
$\Mor(U,G\ds\coprod\limits_j Y_j) \approx$ 
$\Mor(FU,$ $\ds\coprod\limits_j Y_j)$.
\vspi
Sufficiency: \ The natural transformation 
$\Xi:\Mor(-,\ds\coprod\limits_j GY_j) \ra \Mor(-,G\ds\coprod\limits_j Y_j)$ corresponding to the arrow 
$\ds\coprod\limits_j GY_j \ra G\ds\coprod\limits_j Y_j$ has the property that $\Xi_{\ov{U}}$ is bijective 
$\forall \ \ov{U} \in \ov{\sU}$, hence 
$\ds\coprod\limits_j GY_j \approx G\ds\coprod\limits_j Y_j$ 
(cf. p. \pageref{15.19}).]\\
\endgroup %%------------------------------------<<

Notation: $\sU^+$ is the class of objects in \bC that are coproducts of objects in $\ov{\sU}$.

Definition: An object $(\bX,\bff)$ in $\bFIL(\bC)$ is 
\un{completable in $\sU^+$}
\index{completable in $\sU^+$} 
if $X_0 \in \sU^+$ and $\forall \ n \geq 0$, there is an exact triangle 
$X_n \overset{f_n}{\lra} X_{n+1}$ 
$\ra Z_n$ 
$\ra \Sigma X_n$ with $Z_n \in \sU^+$.\\

\begin{proposition} \ %16
Let \bC be a compactly generated triangulated category.  Suppose that $F:\bC \ra \bAB$ is an exact cofunctor which converts coproducts into products $-$then $\exists$ an object $(\bX,\bff)$ in $\bFIL(\bC)$, completable in $\sU^+$, such that 
$\telsub(\bX,\bff)$ represents $F$.
\end{proposition}

[This is implicit in the proof of the Brown representability theorem.  Thus by definition, $X_0 \in \sU^+$.  
Consider the exact triangle 
$K_n \ra X_n$ 
$\overset{f_n}{\ra} X_{n+1} \ra$ 
$\Sigma K_n$.  
Since 
$\Sigma K_n = \Sigma \bigl(\coprod\limits_{\ov{U}} (\ker \Xi_{n,\ov{U}}) \cdot \ov{U}\bigr)$ 
$\approx$ 
$\coprod\limits_{\ov{U}} (\ker \Xi_{n,\ov{U}}) \cdot \Sigma \ov{U}$, 
there is an exact triangle 
$X_n \overset{f_n}{\ra} X_{n+1}$ 
$\ra Z_n \ra$ 
$\Sigma X_n$ 
with $Z_n  \in \sU^+$.]

[Note: \ If $\ov{U} = \Omega^n U$ $(n \geq 1)$, then 
$\Sigma \ov{U} = \Sigma\Omega^nU = \Sigma\Omega(\Omega^{n-1}U) \approx \Omega^{n-1}U \in \ov{\sU}$.]\\

Application: Fix an $X \in \Ob\bC$ $-$then $\exists$ an object $(\bX,\bff)$ in $\bFIL(\bC)$, completable in $\sU^+$, such that $X \approx \telsub(\bX,\bff)$.

[In Proposition 16, take $F = \Mor(-,X)$.]\\

%%----------------------------------------------------------------------------------------------19
Let \bC be a compactly generated triangulated category satisfying the octahedral axiom $-$then one may form 
$\Extx\ov{\sU}$ and $\Extx\sU^+$ 
(cf. p. \pageref{15.20}).  
Example: Using the notation of Proposition 16, $\forall \ n \geq 0$, $X_n \in \Extx\sU^+$.\\

\textbf{\small LEMMA} \ 
Let \bC be a compactly generated triangulated category satisfying the octahedral axiom.  Fix a compact object \mX and suppose that 
$Z^\prime \ra Z \ra Z\pp \ra \Sigma Z^\prime$ is an exact triangle with $Z\pp \in \Extx \sU^+$ $-$then every diagram 
\begin{tikzcd}%[sep=large]
&{X} \ar{d}\\
{Z^\prime} \ar{r}  &{Z}
\end{tikzcd}
can be completed to a commutative diagram 
\begin{tikzcd}%[sep=large]
{X^\prime} \ar{d}  \ar{r} &{X} \ar{d}\\
{Z^\prime} \ar{r}  &{Z}
\end{tikzcd}
in such a way that there is an exact triangle 
$X^\prime \ra X \ra X\pp \ra \Sigma X^\prime$ with $X\pp \in \Extx \ov{\sU}$.

[Argue by induction on the length $l$ of $Z\pp$.

\indent\indent Case 1: \ $l = 1$.  Here $Z\pp \in \sU^+$.  Since \mX is compact, the composition 
$X \ra Z$ $\ra Z\pp$ factors through a finite coproduct $X\pp \subset Z\pp$ and 
\begin{tikzcd}%[sep=large]
&{X} \ar{d}  \ar{r} &{X\pp} \ar{d}\\
{Z^\prime}  \ar{r} &{Z} \ar{r}  &{Z\pp} \ar{r} &{\Sigma Z^\prime}
\end{tikzcd}
extends to a morphism of exact triangles 
\begin{tikzcd}%[sep=large]
{X^\prime} \ar{d} \ar{r} &{X} \ar{d}  \ar{r} &{X\pp} \ar{d} \ar{r} &{\Sigma X^\prime} \ar{d}\\
{Z^\prime}  \ar{r} &{Z} \ar{r}  &{Z\pp} \ar{r} &{\Sigma Z^\prime}
\end{tikzcd}
(cf. Proposition 1).

\indent\indent Case 2: \ $l > 1$.  
By assumption $Z\pp$ occurs in an exact triangle 
$Z_0\pp \ra Z\pp$ 
$\ra Z_1\pp \ra $
$\Sigma Z_0\pp$, where $Z_0\pp, Z_1\pp \in \Extx \sU^+$ and have length $< l$.  
Complete the composite 
$Z \ra Z\pp \ra Z_1\pp$ to an exact triangle 
$Z \ra Z_1\pp$ 
$\ra W \ra \Sigma Z$ (cf. TR$_3$).  
Using the octahedral axiom on 
$Z \ra Z\pp \ra$ $\Sigma Z^\prime  \ra \Sigma Z$,
$Z\pp \ra Z_1\pp \ra$ $\Sigma Z_0\pp \ra \Sigma Z\pp$, construct a factorization 
\begin{tikzcd}%[sep=large]
{Z^\prime} \ar{r} \ar{rd}  &{\ov{Z}} \ar{d}\\
&{Z}
\end{tikzcd}
of $Z^\prime \ra Z$ and exact triangles 
$Z^\prime \ra \ov{Z}$ $\ra Z_0\pp \ra \Sigma Z^\prime$,
$\ov{Z} \ra Z \ra$ $Z_1\pp \ra \Sigma \ov{Z}$.  
Owing to the induction hypothesis, there is a commutative diagram 
\begin{tikzcd}%[sep=large]
{X^\prime} \ar{d}  \ar{r} &{\ov{X}} \ar{d} \ar{r} &{X} \ar{d}\\
{Z^\prime}  \ar{r} &{\ov{Z}} \ar{r}  &{Z}
\end{tikzcd}
and exact triangles 
$X^\prime \ra \ov{X}$ $\ra X_0\pp \ra \Sigma X^\prime$,
$\ov{X} \ra X \ra$ $X_1\pp \ra \Sigma \ov{X}$, 
where $X_0\pp, X_1\pp \in \Extx \ov{\sU}$.
Complete the composite $X^\prime \ra \ov{X} \ra X$ to an exact triangle 
$X^\prime \ra X$ $\ra X\pp$ $\ra \Sigma X^\prime$
(cf. TR$_3$) $-$then the octahedral axiom implies that $X\pp \in \Extx \ov{\sU}$.]\\

\begin{proposition} \ %17
Let \bC be a compactly generated triangulated category satisfying the octahedral axiom 
$-$then every compact object \mX in \bC is a direct summand of an object in $\Extx  \ov{\sU}$.
\end{proposition}

%%----------------------------------------------------------------------------------------------20
[Write $X \approx \tel(\bX,\bff)$ (cf. supra).  Since $\colimx \Mor(X,X_n) \approx \Mor(X,X)$ (cf. Proposition 13), 
$\id_X$ factors through some $X_n$: 
\begin{tikzcd}[sep=small]
{X} \ar{rdd}[swap]{i}  \ar{rr}{\id_X}  &&{X}\\
\\
&{X_n} \ar{ruu}[swap]{r}
\end{tikzcd}
.  
On the other hand, $X_n \in \Extx \sU^+$ and $0 \ra X_n \overset{\id_{X_n}}{\lra}$ $X_n \ra 0$ is exact.  
One may therefore 
apply the lemma to 
\begin{tikzcd}%[sep=large]
&{X} \ar{d}\\
{0} \ar{r} &{X_n}
\end{tikzcd}
and produce a commutative diagram 
\begin{tikzcd}%[sep=large]
{X^\prime} \ar{d} \ar{r} &{X} \ar{d}\\
{0} \ar{r} &{X_n}
\end{tikzcd}
plus an exact triangle 
$X^\prime \ra X$ 
$\ra X\pp \ra \Sigma X^\prime$ with $X\pp \in \Extx \ov{\sU}$.  But  the arrow $X^\prime \ra X$ is the zero morphism, thus 
$X\pp \approx X \amalg \Sigma X^\prime$ 
(cf. p. \pageref{15.21}).]\\

Notation: \ $\cptx\bC$ is the thick subcategory of \bC whose objects are compact.\\

\index{Theorem: Theorem of Neeman-Ravenel}
\index{Theorem of Neeman-Ravenel}
\textbf{\small THEOREM OF NEEMAN-RAVENEL} \quad
Let \bC be a compactly generated triangulated category satisfying the octahedral axiom $-$then the thick subcategory generated by $\sU$ is cpt\bC.

[This is a consequence of Proposition 10 and Proposition 17.]

[Note: \ The thick subcategory generated by $\sU$ is, of course, the intersection of the conglomerate of thick subcategories of \bC containing $\sU$.]\\

The proof of the Neeman-Ravenel theorem depends on the octahedral axiom (by way of Proposition 17) but its use can be eliminated.  Thus, let $\bLambda$ be the thick subcategory generated by $\sU$ and fix a skeleton $\ov{\bLambda}$ of 
$\bLambda$ $-$then $\ov{\bLambda}$  is small (since $\sU$ is a set) and for any \mX in \bC, $\ov{\bLambda}/X$ is the category whose objects are the arrows $K \ra X$ and whose morphisms $(K \ra X) \ra (L \ra X)$ are the commutative triangles 
\begin{tikzcd}[sep=small]
{K} \ar{rdd}  \ar{rr}  &&{L} \ar{ldd}\\
\\
&{X}
\end{tikzcd}
(\mK, \mL in $\ov{\bLambda}$).\\

\textbf{\small LEMMA} \ 
$\forall \ X$, the category $\ov{\bLambda}/X$ is filtered.

[Note: \ The assignment $X \ra \ov{\bLambda}/X$ defines a functor $\bC \ra \bCAT$.]\\

In what follows, $\underset{X}{\colim}$ stands for a colimit calculated over $\ov{\bLambda}/X$.\\

\begin{proposition} \ %18
Let \bC be a compactly generated triangulated category.  
Suppose that $F:\bLambda \ra \bAB$ is an exact functor.  
Given an
$X \in \Ob\bC$, put $\ov{F}X = \underset{X}{\colimx} FK$ $-$then $\ov{F}:\bC \ra \bAB$ is an exact functor which converts coproducts into direct sums.
\end{proposition}

%%----------------------------------------------------------------------------------------------21
[Note: \ $\forall \ K$ in $\bLambda$, $\ov{F}K \approx F K$.]\\

\label{15.40} 
Remark: \ Suppose that $F:\bC \ra \bAB$ is an exact functor which converts coproducts into direct sums $-$then the natural transformation $\ov{\restr{F}{\bLambda}} \ra F$ is a natural isomorphism.  
Proof: The \mX such that the arrows 
$
\begin{cases}
\ \ov{\restr{F}{\bLambda}} \Sigma^n X \ra F \Sigma^n X\\
\ \ov{\restr{F}{\bLambda}} \Omega^n X \ra F \Omega^n X
\end{cases}
$
are isomorphisms $\forall \ n \geq 0$ constitute the object class of a triangulated subcategory of \bC containing $\sU$ and closed under the formation of coproducts in \bC, thus is all of \bC (cf. Proposition 14).\\

\index{Theorem: Theorem of Neeman-Ravenel (bis)}
\index{Theorem of Neeman-Ravenel (bis)}
\textbf{\small THEOREM OF NEEMAN-RAVENEL (bis)} \quad
Let \bC be a compactly generated triangulated category $-$then the thick subcategory generated by $\sU$ is cpt\bC.

[$\forall$ compact \mX, the exact functor $\Mor(X,-)$ converts coproducts into direct sums.  
Therefore, by the above remark, 
$\ov{\restr{\Mor(X,-)}{\bLambda}} \approx \Mor(X,-)$, so $\id_X$ factors through some \mK in $\ov{\bLambda}$:
\begin{tikzcd}[sep=small]
{X} \ar{rdd}[swap]{i}  \ar{rr}{\id_X}  &&{X}\\
\\
&{K} \ar{ruu}[swap]{r}
\end{tikzcd}
.]\\

\begin{proposition} \ %19
Let \bC be a compactly generated triangulated category $-$then $\cptx\bC$ has a small skeleton.\\
\end{proposition}

Let \bC be a compactly generated triangulated category $-$then the additive functor category 
$[(\cptx\bC)^{\OP},\bAB]^+$ is a complete and cocomplete abelian category and has enough projectives 
(cf. p. \pageref{15.22}).  
Call $\bEX[(\cptx \bC)^\OP,\bAB]^+$
\index{$\bEX[(\cptx \bC)^\OP,\bAB]^+$} 
the full subcategory of 
$[(\cpt\bC)^{\OP},\bAB]^+$ whose objects are the exact cofunctors $F:\cptx \bC \ra \bAB$.\\

\begin{proposition} \ %20
Let \bC be a compactly generated triangulated category $-$then all the projective objects of 
$[(\cpt \bC)^\OP,\bAB]^+$ lie in $\bEX[(\cptx \bC)^\OP,\bAB]^+$.
\end{proposition}

[Every projective object of $[(\cpt \bC)^\OP,\bAB]^+$ is a direct summand of a coproduct of representable cofunctors.]\\

\begin{proposition} \ %21
Let \bC be a compactly generated triangulated category $-$then every object in $[(\cptx \bC)^\OP,\bAB]^+$ of finite projective dimension belongs to $\bEX[(\cptx \bC)^\OP,\bAB]^+$.\\
\end{proposition}

Notation: Write $h_X$ for the restriction $\restr{\Mor(-,X)}{\cptx \bC}$ and write $h_f:h_X \ra h_Y$ for the natural transformation induced by the morphism $f:X \ra Y$.\\

\begingroup%%----------------------------------->>
\fontsize{9pt}{11pt}\selectfont
\label{15.30}
\textbf{\small FACT} \ 
Let \bC be a compactly generated triangulated category $-$then the functor 
$h:\bC \ra [(\cpt \bC)^\OP,$ $\bAB]^+$ is exact, conservative, and preserves products $\&$ coproducts.\\
\endgroup %%------------------------------------<<

%%----------------------------------------------------------------------------------------------22
Let \bC be a compactly generated triangulated category $-$then \bC is said to admit 
\un{Adams representability}
\index{Adams representability} 
if the following conditions are satisfied.

\indent\indent (ADR$_1$) \  Every exact cofunctor $F:\cptx \bC \ra \bAB$ is 
\un{representable in the large}
\index{representable in the large} 
i.e., $\exists$ an $X \in \Ob\bC$ and a natural isomorphism 
$h_X \ra F$.

\indent\indent (ADR$_2$) \  Every natural transformation $h_X \ra h_Y$ is induced by a morphism $f:X \ra Y$.\\

\begingroup%%----------------------------------->>
\fontsize{9pt}{11pt}\selectfont
\textbf{\small FACT} \ 
Suppose that \bC admits Adams representability $-$then $\bIND(\cptx\bC)$ is equivalenct to \\
$\bEX[(\cptx\bC)^\OP,\bAB]^+$.\\
\endgroup %%------------------------------------<<

\textbf{\small LEMMA} \ 
Let \bC be a compactly generated triangulated category.  Assume: \bC admits Adams representability $-$then 
$h_X \approx h_Y$ $\implies$ $X \approx Y$, thus an object representing a given exact cofunctor 
$F:\cptx \bC \ra \bAB$ is unique up to isomorphism.\\

\label{17.13}
Suppose that \bC admits Adams representability $-$then $\forall  \ X, \ Y \in \Ob\bC$, there is a surjection 
$\Mor(X,Y) \ra \Nat(h_X,h_Y)$, viz. $f \ra h_f$.  \ 
Definition: $f$ is said to be a 
\un{phantom} \un{map}
\index{phantom map} 
provided that $h_f = 0$.  So, if 
$\Ph(X,Y)$
\index{\Ph(X,Y)} 
is the subgroup of $\Mor(X,Y)$ consisting of phantom maps, then the sequence 
$0 \ra \Ph(X,Y) \ra$ 
$\Mor(X,Y) \ra $ 
$\Nat(h_X,h_Y) \ra 0$ is short exact.

[Note: Let $f \in \Ph(X,Y)$ $-$then for any $\phi:X^\prime \ra X$, $f \circx \phi \in \Ph(X^\prime,Y)$, and for any 
$\psi:Y \ra Y^\prime$, $\psi \circx f \in \Ph(X,Y^\prime)$.  This has the consequence that it makes sense to form the quotient category $\bC/\bPh$, where the set of morphisms from \mX to \mY is $\Mor(X,Y)/\Ph(X,Y)$.]\\

\textbf{\small LEMMA} \ 
Let \bC be a compactly generated triangulated category.  Assume: \bC admits Adams representability $-$then 
$h_X$ is projective iff \mX is a direct summand of a coproduct of compact objects.\\

\begingroup%%----------------------------------->>
\fontsize{9pt}{11pt}\selectfont
\textbf{\small EXAMPLE}  \ 
Consider any exact triangle 
$W \overset{w}{\ra} \ds\coprod\limits_i X_i \overset{t}{\ra} \ds\prod\limits_i X_i \ra \Sigma W$ 
($t$ as on p. \pageref{15.23}) 
$-$then $w$ is a phantom map.\\
\endgroup %%------------------------------------<<

\begingroup%%----------------------------------->>
\fontsize{9pt}{11pt}\selectfont
\textbf{\small FACT} \ 
Suppose that \bC admits Adams representability $-$then $f:X \ra Y$ is a phantom map iff $\forall$ compact \mK and every 
$\phi:K \ra X$, the composite $f \circx \phi$ vanishes.\\
\endgroup %%------------------------------------<<

\begingroup%%----------------------------------->>
\fontsize{9pt}{11pt}\selectfont
\textbf{\small EXAMPLE}  \ 
Given an $X \in \Ob\bC$, complete 
$\ds\coprod\limits_{\ov{\bLambda}/X} K \ra X$ to an exact triangle 
$\ds\coprod\limits_{\ov{\bLambda}/X} K \ra X$ 
$\overset{\theta}{\lra} \ov{X}$ $\ra \ds\coprod\limits_{\ov{\bLambda}/X} \Sigma K$ (cf. TR$_3$) 
$-$then $\Theta$ is a phantom map.  
Moreover, every $f \in \Ph(X,Y)$ factors through $\Theta$.\\
%%----------------------------------------------------------------------------------------------23
Corollary:  \ All phantom maps out of \mX vanish iff $\Theta = 0$.  And, when $\Theta = 0$, \mX is a direct summand of 
$\ds\coprod\limits_{\ov{\bLambda}/X} K$.
\vspi
[Note: \ Therefore $\Theta$ is a ``universal'' phantom map 
(cf. p. \pageref{15.24}).]\\
\endgroup %%------------------------------------<<

\label{15.29}
\begingroup%%----------------------------------->>
\fontsize{9pt}{11pt}\selectfont
\textbf{\small FACT} \ 
Suppose that \bC admits Adams representability $-$then $f:X \ra Y$ is a phantom map iff $\forall$ exact functor 
$F:\bC \ra \bAB$ which convertes coproducts into direct sums, $F f = 0$.\\
\endgroup %%------------------------------------<<

\begin{proposition} \ %22
Let \bC be a compactly generated triangulated category.  
Assume: \bC admits Adams representability.  Let 
$\Delta:\bI \ra \bC$ be a diagram, where \bI is filtered and $\forall \ i \in \Ob\bI$, $\Delta_i$ is compact $-$then $\Delta$ has a miminal weak colimit.
\end{proposition}

[Put $F = \colimx h_{\Delta_i}$ (thus $\forall$ compact \mK, $FK = \colimx \Mor(K,\Delta_i)$).  
Since \bAB is 
Grothendieck, \mF is exact, so by ADR$_1$, $\exists$ an $X \in \Ob\bC$ and a natural isomorphism $h_X \ra F$.  
Claim: \mX is a minimal weak colimit of $\Delta$.  
Indeed, $\forall \ i$, there is a natural transformation 
$\Xi_i:h_{\Delta_i} \ra h_X$ and, by ADR$_2$, $\Xi_i = h_{f_i}$ $(\exists \ f_i:\Delta_i \ra X)$. 
 Moreover, $f_i$ is determined up to an element of $\Ph(\Delta_i,X)$.  
But $\Delta_i$ compact $\implies$ $\Ph(\Delta_i,X) = 0$, hence $f_i$ is unique.  
Consequently, 
$\{\Delta_i \overset{f_i}{\lra} X\}$ is a natural sink.  
If now 
$\{\Delta_i \overset{g_i}{\lra} X\}$ is another natural sink, then 
$\exists$ $\Xi\in \Nat(h_X,h_Y)$: $\forall \ i$, $h_{g_i} = \Xi \circx h_{f_i}$.   
However $\Xi = h_\phi$ for some $\phi:X \ra Y$ 
(cf. ADR$_2$) and this means that $g_i = \phi \circx f_i$.  
Therefore \mX is a weak colimit of $\Delta$.  
Minimality is obvious.]\\

\begingroup%%----------------------------------->>
\fontsize{9pt}{11pt}\selectfont
\textbf{\small EXAMPLE}  \  Suppose that \bC admits Adams representability.  
Fix an $X \in \Ob\bC$ and consider the functor 
$\ov{\bLambda}/X \ra \bC$ that sends $K \ra X$ to \mK.  
Since $\ov{\bLambda}/X$ is filtered, this functor has a minimal weak colimit $L_X$ (cf. Proposition 22).  
There is an arrow $L_X \ra X$ and $\forall \ \ov{U} \in \ov{\sU}$, 
$\Mor(\ov{U},L_X) \approx$ 
$\underset{X}{\colimx} \Mor(\ov{U},K) \approx$ $\Mor(\ov{U},X)$ $\implies$ 
$L_X \approx X$ (cf. p. \pageref{15.25}).\\
\endgroup %%------------------------------------<<

\begingroup%%----------------------------------->>
\fontsize{9pt}{11pt}\selectfont
\textbf{\small FACT} \ 
Let \bC be a compactly generated triangulated category.  
Assume: Every functor from a filtered category \bI to \bC with compact values has a minimal weak colimit $-$then \bC admits Adams representability.\\ 
\endgroup %%------------------------------------<<

\textbf{\small LEMMA} \ 
Let \bC be a compactly generated triangulated category.  
Assume: \bC admits Adams representability $-$then for any $X \in \bC$, there is an exact triangle 
$P \ra Q \ra X \ra \Sigma P$ such that $h_P$ $\&$ $h_Q$ are projective and the sequence 
$0 \ra h_P \ra h_Q \ra h_X \ra 0$ is short exact.

[The functor $\ov{\bLambda}/X \ra \bC$  that sends $K \ra X$ to \mK has a minimal weak colimit, viz. \mX (see the preceeding example).  
It also has a weak colimit \mY constructed via the procedure on 
p. \pageref{15.26}: 
$\coprod\limits_{K \ra L} K \ra \coprod\limits_{\ov{\bLambda}/X} K \ra Y \ra \coprod\limits_{K \ra L} \Sigma K$.
Since \mX is minimal, $\exists$ arrows $\phi:X \ra Y$,
%%----------------------------------------------------------------------------------------------24
$\psi:Y \ra X$, such that $\psi \circx \phi$ is an isomorphism and the triangles 
\begin{tikzcd}%[sep=large]
{\coprod\limits_{\ov{\bLambda}/X} K} \ar{rd} \ar{r} &{X} \ar{d}{\phi}\\
&{Y}
\end{tikzcd}
,
\begin{tikzcd}%[sep=large]
{\coprod\limits_{\ov{\bLambda}/X} K} \ar{rd} \ar{r} &{Y} \ar{d}{\psi}\\
&{X}
\end{tikzcd}
\ 
commute (cf. p. \pageref{15.27} ff.).  
Define \mP by requiring that \ 
$P \ra \coprod\limits_{\ov{\bLambda}/X} K$ 
$\ra X$ 
$\ra$ $\Sigma P$ be exact.  
Using Proposition 1, determine arrows 
$f:P \ra \coprod\limits_{K \ra L} K$,
$g:\coprod\limits_{K \ra L} K \ra P$ such that the diagram
\[
\begin{tikzcd}%[sep=large]
{P}\ar{d}{f}\ar{r}
&{\coprod\limits_{\ov{\bLambda}/X} K}\arrow[d,shift right=0.5,dash] \arrow[d,shift right=-0.5,dash] \ar{r}
&{X}\ar{d}{\phi}\ar{r}
&{\Sigma P}\ar{d}{\Sigma f}\\
{\coprod\limits_{K \ra L} K}\ar{d}{g}\ar{r}
&{\coprod\limits_{\ov{\bLambda}/X} K}\arrow[d,shift right=0.5,dash] \arrow[d,shift right=-0.5,dash]\ar{r}
&{Y}\ar{d}{\psi}\ar{r}
&{\coprod\limits_{K \ra L} \Sigma K}\ar{d}{\Sigma g}\\
{P}\ar{r}
&{\coprod\limits_{\ov{\bLambda}/X} K}\ar{r}
&{X}\ar{r}
&{\Sigma P}
\end{tikzcd}
\]
commutes $-$then $g \circx f$ is an isomorphism (cf. p. \pageref{15.28}), hence $h_P$ is a direct summand of 
$\coprod\limits_{K \ra L} h_K$ which implies that $h_P$ is projective.  And with 
$Q = \coprod\limits_{\ov{\bLambda}/X} K$, the sequence 
$0 \ra h_P \ra$
$h_Q \ra h_X \ra 0$ is short exact.]\\

Remark: The arrow $X \ra \Sigma P$ is a phantom map and if $f:X \ra Y$ is a phantom map, then there is a commutative triangle
\begin{tikzcd}%[sep=large]
{X} \ar{d}[swap]{f} \ar{r} &{\Sigma P} \ar{ld}\\
{Y}
\end{tikzcd}
(cf. p. \pageref{15.29} ff.).\\

Example: $f \in \Ph(X,Y)$ $\&$ $g \in \Ph(Y,Z)$ $\implies$ $g \circx f = 0$.  Proof: $h_P$ projective $\implies$ 
$h_{\Sigma P}$ projective $\implies$ $\Ph(\Sigma P,Z) = 0$.\\

\begin{proposition} \ %23
Let \bC be a compactly generated triangulated category.  Assume: \bC admits Adams representability $-$then 
$\bEX[(\cpt \bC)^\OP,\bAB]^+$ is the full subcategory of $[(\cpt \bC)^\OP,\bAB]^+$ whose objects have projective dimension 
$\leq 1$.
\end{proposition}

[On account of Proposition 21, it need only be shown that every \mF in 
$\bEX[(\cpt \bC)^\OP,$ $\bAB]^+$ has projective 
dimension $\leq 1$.  But by ADR$_1$, $\exists$ $X: h_X \approx F$ and the lemma implies that $h_X$ has a projective resolution of length $\leq 1$.]\\

\begingroup%%----------------------------------->>
\fontsize{9pt}{11pt}\selectfont
\textbf{\small FACT} \ 
Suppose that \bC admits Adams representability $-$then $\forall \ X, Y \in \Ob\bC$, 
$\Ph(\Omega X,Y) \approx \Extx(h_X,h_Y)$.\\
\endgroup %%------------------------------------<<

%%----------------------------------------------------------------------------------------------25
\textbf{\small LEMMA} \ 
Let \bC be a compactly generated triangulated category $-$then every exact cofunctor $F:\cpt \bC \ra \bAB$ of projective dimension $\leq 1$ has a projective resolution $0 \ra H \ra G \ra F \ra 0$, where \mG, \mH are coproducts of representable cofunctors.

[By hypothesis, there is a projective resolution 
$0 \ra F\pp \ra F^\prime \ra F \ra 0$.  Here $F^\prime$ is a coproduct of representable cofunctors, while $F\pp$ is a direct summand of a coproduct of representable cofunctors, say 
$F\pp \amalg \ov{F}\pp \approx \Phi$.  Noting that 
$\coprod\limits_1^\infty \Phi \approx F\pp \amalg \coprod\limits_1^\infty \Phi$, consider 
$0 \ra F\pp \amalg \coprod\limits_1^\infty \Phi \ra F^\prime \amalg \coprod\limits_1^\infty \Phi \ra F \ra 0$.]\\

\begin{proposition} \ %24
Let \bC be a compactly generated triangulated category.  
Assume: Every exact cofunctor $F:\cpt \bC \ra \bAB$ has projective dimension $\leq 1$ $-$then \bC admits Adams representability.
\end{proposition}

[It is a question of checking the validity of ADR$_1$ and ADR$_2$.

\indent\indent Re: \ ADR$_1$.  Fix an exact cofunctor $F:\cpt\bC \ra \bAB$ and resolve it per the lemma: 
$0 \ra H \ra G \ra F \ra 0$.  
Write 
$G = \amalg \ \Mor(-,K)$, $H = \amalg \ \Mor(-,L)$ $-$then the arrow $H \ra G$ gives rise to 
a morphism $\amalg L \ra \amalg K$ which can be completed to an exact triangle 
$\amalg L \ra \amalg K$ 
$\ra X \ra \amalg \Sigma L$ (cf. TR$_3$) and $h_X \approx F$.

\indent\indent Re: \ ADR$_2$.  Fix a natural transformation $\Xi:h_X \ra h_Y$.  Choose projective resolutions 
$0 \ra H_X \ra G_X \ra h_X \ra 0$, 
$0 \ra H_Y \ra G_Y \ra h_Y \ra 0$ per the lemma and lift $\Xi$ to a commutative diagram 
\begin{tikzcd}%[sep=large]
{0} \ar{r} &{H_X} \ar{d}  \ar{r}  &{G_X} \ar{d}  \ar{r} &{h_X} \ar{d}{\Xi} \ar{r} &{0}\\
{0} \ar{r} &{H_Y}\ar{r}  &{G_Y} \ar{r} &{h_Y} \ar{r} &{0}
\end{tikzcd}
.  Write 
$G_X = \amalg \hspace{0.02cm} \Mor(-,K_X)$, 
$H_X = \amalg \hspace{0.02cm} \Mor(-,L_X)$, 
$G_Y = \amalg \hspace{0.02cm} \Mor(-,K_Y)$, 
$H_Y = \amalg \hspace{0.02cm} \Mor(-,L_Y)$, $-$then there is a commutative diagram 
\begin{tikzcd}%[sep=large]
{\amalg L_X} \ar{d} \ar{r} &{\amalg K_X} \ar{d}\\
{\amalg L_Y} \ar{r} &{\amalg K_Y}
\end{tikzcd}
in \bC.  Complete it to a morphism
\begin{tikzcd}%[sep=large]
{\amalg L_X} \ar{d} \ar{r} &{\amalg K_X} \ar{d} \ar{r}  &{X^\prime} \ar{d}{f^\prime} \ar{r} &{\amalg \Sigma L_X} \ar{d}\\
{\amalg L_Y} \ar{r} &{\amalg K_Y} \ar{r}  &{Y^\prime} \ar{r} &{\amalg \Sigma L_Y}
\end{tikzcd}
of exact triangles (cf. TR$_3$ $\&$ TR$_5$).  The rows in the commutative diagram 
\begin{tikzcd}%[sep=large]
{0} \ar{r} &{H_X} \ar{d}  \ar{r}  &{G_X} \ar{d}  \ar{r} &{h_{X^\prime}} \ar{d}{h_{f^\prime}} \ar{r} &{0}\\
{0} \ar{r} &{H_Y}\ar{r}  &{G_Y} \ar{r} &{h_{Y^\prime}} \ar{r} &{0}
\end{tikzcd}
are short exact.  Working with \mX, the composite $\amalg L_X \ra \amalg K_X \ra X$ is a phantom map, hence vanishes, thus 
$\exists$ a commutative triangle
\begin{tikzcd}%[sep=large]
{\amalg K_X} \ar{d}\ar{r} &{X^\prime} \ar{ld}{\phi}\\
{X}
\end{tikzcd}
and $h_\phi:h_{X^\prime} \ra h_X$ is a natural isomorphism, so $\phi:X^\prime \ra X$ is an isomorphism 
($h$ is conservative (cf. p. \pageref{15.30}).
%%----------------------------------------------------------------------------------------------26
Similar considerations apply to \mY.  Since 
\begin{tikzcd}%[sep=large]
{h_{X^\prime}} \ar{d}[swap]{h_{f^\prime}} \ar{r}{h_\phi} &{h_X} \ar{d}{\Xi}\\
{h_{Y^\prime}}  \ar{r}[swap]{h_\psi} &{h_Y}
\end{tikzcd}
commutes, it follows that $\Xi = h_{\psi \circx f^\prime \circx \phi^{-1}}$.]\\

Let \bC be a compactly generated triangulated category $-$then Propositions 23 and 24 tell us that \bC admits Adams representability iff every object in $\bEX[(\cpt \bC)^\OP,\bAB]^+$ has projective 
dimension $\leq 1$ in $[(\cpt \bC)^\OP,\bAB]^+$.  
And this condition can be realized.  Indeed, it suffices that cpt\bC possess a countable skeleton, (cf. infra).
[Note: \ Recall that in any event cpt\bC has a small skeleton (cf. Proposition 19).]\\

\index{Neeman's countability criterion}
\textbf{\small NEEMAN'S COUNTABILITY CRITERION} \quad
Let \bC be a triangulated category with finite coproducts and a countable skeleton $-$then every object of 
$\bEX[\bC^\OP,\bAB]^+$ has projective dimension $\leq 1$ in $[(\bC^\OP,\bAB]^+$.

[Note: \ $\bEX[\bC^\OP,\bAB]^+$ is the full subcategory of $[\bC^\OP,\bAB]^+$ whose objects are the exact cofunctors $F:\bC \ra \bAB$.]\\

\begingroup%%----------------------------------->>
\fontsize{9pt}{11pt}\selectfont
The stable homotopy category is a compactly generated triangulated category and its full subcategory of compact objects has a countable skeleton.  Therefore the stable homotopy category admits Adams representability.\\
\endgroup %%------------------------------------<<

The proof of Neeman's countability criterion requires some preparation.  Call an object of $[\bC^\OP,\bAB]^+$ 
\un{free}
\index{free (object in $[\bC^\OP,\bAB]^+$)} 
if it is a coproduct of representable cofunctors.  
Definition: $\forall$ \mF in $[\bC^\OP,\bAB]^+$, $\#(F)$ is the smallest infinite cardinal $\kappa$ for which there is a free presentation 
$F\pp \ra F^\prime \ra F \ra 0$, where $F^\prime$, $F\pp$ are coproducts of $\leq \kappa$ representable cofunctors.

Observation: If 
$0 \ra F\pp \ra F^\prime \ra F \ra 0$
is a short exact sequence in $[\bC^\OP,\bAB]^+$ and if 
$\#(F\pp) \leq \kappa$, 
$\#(F^\prime) \leq \kappa$, then 
$\#(F) \leq \kappa$.

Let $\kappa$ be an infinite cardinal $-$then \bC is said to satisfy 
\un{condition $\kappa$}
\index{condition $\kappa$} 
if for any \mF in $\bEX[\bC^\OP,\bAB]^+$ and any morphism $\Phi \ra F$, where $\#(\Phi) \leq \kappa$, there is a factorization 
$\Phi \ra \Psi \ra F$ such that $\Psi \ra F$ is a monomorphism and $\Psi$ has a free resolution 
$0 \ra $ 
$\Psi\pp \ra $ 
$\Psi^\prime \ra $ 
$\Psi \ra 0$, 
where $\Psi^\prime$, $\Psi\pp$ are coproducts of $\leq \kappa$ representable cofunctors ($\implies$   
$\#(\Psi) \leq \kappa$).

Observation: Suppose that \bC satisfies condition $\kappa$ $-$then every object \mF of 
$\bEX[\bC^\OP,$ $\bAB]^+$ with 
$\#(F) \leq \kappa$ has a free resolution 
$0 \ra F\pp \ra F^\prime \ra F \ra 0$, 
where $F^\prime$, $F\pp$ are coproducts of $\leq \kappa$ representable cofunctors.  
In particular: The projective dimension of \mF is $\leq 1$.\\

\textbf{\small LEMMA} \ 
Suppose that \bC satisfies condition $\kappa$.  Let $F \ra G$ be a monomorphism of exact cofunctors, where 
$\#(F) \leq \kappa$, 
$\#(G) \leq \kappa$ 
$-$then for any free resolution 
%%----------------------------------------------------------------------------------------------27
$0 \ra F\pp \ra F^\prime \ra F \ra 0$ of \mF, there exists a free resolution 
$0 \ra G\pp \ra G^\prime \ra G \ra 0$ of \mG and a commutative diagram
\begin{tikzcd}%[sep=large]
{0} \ar{r}
&{F\pp} \ar{d} \ar{r}
&{F^\prime} \ar{d} \ar{r}
&{F} \ar{d} \ar{r}
&{0}\\
{0} \ar{r}
&{G\pp} \ar{r}
&{G^\prime}  \ar{r}
&{G} \ar{r}
&{0}
\end{tikzcd}
such that $F\pp \ra G\pp$, $F^\prime \ra G^\prime$ are split monomorphisms.

[Complete $F \ra G$ to a short exact sequence 
$0 \ra F \ra G \ra H \ra 0$.  Since \mF, \mG are exact, so is \mH.  Moreover, 
$\#(F) \leq \kappa$, 
$\#(G) \leq \kappa$, 
$\implies$ 
$\#(H) \leq \kappa$ (cf. supra).  Fix a free resolution 
$0 \ra H\pp \ra H^\prime \ra H \ra 0$, where $H^\prime$, $H\pp$ are coproducts of $\leq \kappa$ representable cofunctors 
and extend
\[
\begin{tikzcd}%[sep=large]
&&&{0} \ar{d}\\
{0} \ar{r}
&{F\pp} \ar{r}
&{F^\prime} \ar{r}
&{F} \ar{d} \ar{r}
&{0}\\
&&&{G} \ar{d}\\
{0} \ar{r}
&{H\pp} \ar{r}
&{H^\prime} \ar{r}
&{H} \ar{d}\ar{r}
&{0}\\
&&&{0}
\end{tikzcd}
\]
in the obvious way: 
$0 \ra F\pp \oplus H\pp \ra F^\prime \oplus H^\prime \ra G \ra 0$.]

[Note: \ Therefore if $F^\prime$ and $F\pp$ are coproducts of $\leq \kappa$ representable cofunctors, then 
$G^\prime = F^\prime \oplus H^\prime$, 
$G\pp = F\pp \oplus H\pp$ are coproducts of $\leq \kappa$ representable cofunctors.]\\

\index{Main Lemma (triangluated categories)}
\textbf{\small MAIN LEMMA} \quad
Let \bC be a countable triangulated category with finite coproducts $-$then \bC satisfies condition $\kappa$ for every $\kappa$, hence Neeman's countability criterion is valid.

[Fix an \mF in $\bEX[\bC^\OP,\bAB]^+$ and a morphism $\Phi \ra F$.

\indent\indent $\#(\Phi) = \omega$.  
There is a free presentation 
$\Phi\pp \ra \Phi^\prime \ra \Phi \ra 0$, where $\Phi^\prime$, $\Phi\pp$ are countable coproducts of representable cofunctors.  
Accordingly, one can assume without loss of generality that $\Phi$ is a  countable coproduct of representable cofunctors (replace 
$\Phi \ra F$ by $\Phi^\prime \ra \Phi \ra F$), say 
$\Phi = \coprod\limits_0^\infty \Mor(-,X_i)$, the morphism 
$\Phi \ra F$ corresponding to a sequence of natural transformations 
$\Mor(-,X_i) \ra F$.  
Put $X_i^0 = X_i$.  
Since \bC is countable, $\forall \ X \in \Ob\bC$, 
$\Mor(X,\coprod\limits_{i=0}^k X_i^0)$ is countable, thus its subset $S_{X,k}$ consisting of the arrows for which the composite 
$\Mor(-,X) \ra \coprod\limits_{i=0}^k \Mor(X,X_i^0) \ra F$ vanishes is countable.  
Enumerate the elements of 
$\bigcup\limits_{X,k} S_{X,k}$.  Supposing that 
$X \ra \coprod\limits_{i=0}^k X_i^0$ is the $l^{th}$ such, define $X_l^1$ by the exact triangle 
$X \ra \coprod\limits_{i=0}^k X_i^0$ $\ra X_l^1$ $\ra \Sigma X$ (cf. TR$_3$).  The natural
%
%%----------------------------------------------------------------------------------------------28
transformation 
$\coprod\limits_{i=0}^k \Mor(-,X_i^0) \ra F$ determines an element 
$x \in F \coprod\limits_{i=0}^k X_i^0$, that, under the arrow 
$F \coprod\limits_{i=0}^k X_i^0 \ra FX$, is sent to 0.  
Since \mF is exact, $\exists$ an element of $FX_l^1$ mapping to $x$.  
This means that 
$\coprod\limits_{i=0}^k \Mor(-,X_i^0) \ra F$ factors as 
$\coprod\limits_{i=0}^k \Mor(-,X_i^0) \ra \Mor(-,X_l^1) \ra F$.  
Iterate the procedure:  From the set $\{X_l^1\}$ one can produce the set $\{X_l^2\}$.  
Continuing, the upshot is a countable filtered category \bI whose objects are the $X_l^k$ and whose morphisms $X_l^k \ra X_{l^\prime}^{k^\prime}$ are the identities and the composites arising from the construction.  
There is a functor 
$\bI \ra [\bC,\bAB]^+$ that sends $X_l^k$ to $\Mor(-,X_l^k)$.  
The natural transformations 
$\Mor(-,X_l^k) \ra F$ constitute a natural sink and the arrow 
$\colimx \Mor(-,X_l^k) \ra F$ is a monomorphism.  
Definition: $\Psi = \colimx \Mor(-,X_l^k)$.  
It is clear that the $\Mor(-,X_i) \ra F$ factor through $\Psi$.  
To show that $\Psi$ has a free resolution
$0 \ra \Psi\pp \ra \Psi^\prime \ra \Psi \ra 0$, where $\Psi^\prime$, $\Psi\pp$ are countable coproducts of representable cofunctors, fix a final functor 
$\nabla:[\N] \ra \bI$ (see below) $-$then 
$\Psi \approx \colimx \Mor(-,\nabla_n)$ and there is a short exact sequence 
$0 \ra \coprod\limits_n \Mor(-,\nabla_n)$ 
$\overset{\text{sf}}{\lra} \coprod\limits_n \Mor(-,\nabla_n)$ $\ra \Psi \ra 0$.  
Here the $n^{th}$ component of sf is the arrow 
\begin{tikzcd}%[sep=large]
{\nabla_n}  \ar{rr}{
\begin{pmatrix}
\id\\
-f_n\\
\end{pmatrix}
}
&&{\nabla_n \amalg \nabla_{n+1}}
\end{tikzcd}
$(f_n:\nabla_n \ra \nabla_{n+1})$.
%{\nabla_n}  \ar{rr}{\left(\overset{\id}{-f_n}\right)} &&{\nabla_n \amalg \nabla_{n+1}}

\indent\indent $\#(\Phi) = \kappa (>\omega)$.  
The induction hypothesis is that \bC satisfies condition $\kappa^\prime$ for all infinite cardinals $\kappa^\prime < \kappa$.  One can assume from the start that $\Phi$ is a coproduct of $\leq \kappa$ representable cofunctors.  
If $\Phi$ is the coproduct of $< \kappa$ representable cofunctors, we are done.  
Suppose, therefore, that
$\Phi = \coprod\limits_{0 \leq \alpha < \kappa} \Mor(-,X_\alpha)$.  
The idea then is to define for each $\alpha \in [\omega,\kappa[$ a subobject $\Psi_\alpha \subset F$ such that 
$\alpha < \beta$ $\implies$ $\Psi_\alpha \subset \Phi_\beta$ and which has a free resolution 
$0 \ra \Psi_\alpha\pp \ra \Psi_\alpha^\prime \ra \Psi_\alpha \ra 0$, where 
$\Psi_\alpha^\prime$, 
$\Psi_\alpha\pp$
are coproducts of $\leq \#(\alpha)$ representable cofunctors.  
Matters will be arranged so as to ensure that 
$\coprod\limits_{i < \alpha} \Mor(-,X_i) \ra F$
factors as
$\coprod\limits_{i < \alpha} \Mor(-,X_i) \ra \Psi_\alpha \ra F$.
In addition, when $\alpha < \beta$, there will be a commutative diagram
\begin{tikzcd}%[sep=large]
{0} \ar{r}
&{\Psi_\alpha\pp} \ar{d} \ar{r}
&{\Psi_\alpha^\prime} \ar{d}\ar{r}
&{\Psi_\alpha} \ar{d} \ar{r}
&{0}\\
{0} \ar{r}
&{\Psi_\beta\pp} \ar{r}
&{\Psi_\beta^\prime} \ar{r}
&{\Psi_\beta} \ar{r}
&{0}
\end{tikzcd}
with
$\Psi_\alpha^\prime \ra \Psi_\beta^\prime$,
$\Psi_\alpha\pp \ra \Psi_\beta\pp$
split monomorphisms, and when 
$\alpha < \beta < \gamma$, 
the composite
\[
\begin{tikzcd}%[sep=large]
{0} \ar{r}
&{\Psi_\alpha\pp} \ar{d} \ar{r}
&{\Psi_\alpha^\prime} \ar{d}\ar{r}
&{\Psi_\alpha} \ar{d} \ar{r}
&{0}\\
{0} \ar{r}
&{\Psi_\beta\pp} \ar{d}\ar{r}
&{\Psi_\beta^\prime} \ar{d}\ar{r}
&{\Psi_\beta} \ar{d}\ar{r}
&{0}\\
{0} \ar{r}
&{\Psi_\gamma\pp} \ar{r}
&{\Psi_\gamma^\prime} \ar{r}
&{\Psi_\gamma} \ar{r}
&{0}
\end{tikzcd}
\]
%%----------------------------------------------------------------------------------------------29
will equal
$
\begin{tikzcd}%[sep=large]
{0} \ar{r}
&{\Psi_\alpha\pp} \ar{d}\ar{r}
&{\Psi_\alpha^\prime} \ar{d}\ar{r}
&{\Psi_\alpha} \ar{d} \ar{r}
&{0}\\
{0} \ar{r}
&{\Psi_\gamma\pp} \ar{r}
&{\Psi_\gamma^\prime} \ar{r}
&{\Psi_\gamma} \ar{r}
&{0}
\end{tikzcd}
.$  
Thus determine $\Psi_\omega$ by applying the above to the arrow 
$\coprod\limits_{i < \omega} \Mor(-,X_i) \ra F$.  
Proceeding, let $\omega < \alpha$, the supposition being that the $\Psi_i$ have been defined $\forall \ i < \alpha$.  
If $\alpha$ is a successor ordinal, say $\alpha = \beta + 1$, set $\kappa^\prime = \#(\Psi_\beta)$ and consider the morphism 
$\Psi_\beta \oplus \Mor(-,X_\beta) \ra F$.  
Appeal to the induction hypothesis to secure a factorization 
$\Psi_\beta \oplus \Mor(-,X_\beta) \ra \Phi_{\beta+1} \ra F$.  $\Psi_\beta$ is obviously a subobject of 
$\Psi_{\beta+1}$ and since 
\bC satisfies condition $\kappa^\prime$, the lemma guarantees that the free resolution 
$0 \ra \Psi_\beta\pp \ra \Psi_\beta^\prime \ra \Psi_\beta \ra 0$ can be extended to a map of free resolutions 
\begin{tikzcd}%[sep=large]
{0} \ar{r}
&{\Psi_\beta\pp} \ar{d}\ar{r}
&{\Psi_\beta^\prime} \ar{d}\ar{r}
&{\Psi_\beta} \ar{d} \ar{r}
&{0}\\
{0} \ar{r}
&{\Psi_{\beta + 1}\pp} \ar{r}
&{\Psi_{\beta + 1}^\prime} \ar{r}
&{\Psi_{\beta + 1}} \ar{r}
&{0}
\end{tikzcd}
with 
$\Psi_\beta^\prime \ra \Psi_{\beta+1}^\prime$, 
$\Psi_\beta\pp \ra \Psi_{\beta+1}\pp$
split monomorphisms and 
$\Psi_{\beta+1}^\prime$, $\Psi_{\beta+1}\pp$ 
(as well as
$\Psi_{\beta+1}^\prime/\Psi_\beta^\prime$, $\Psi_{\beta+1}\pp/\Psi_\beta\pp$) 
a coproduct of $\leq \kappa^\prime$ representable cofunctors. 
If $\alpha$ is a limit ordinal, put 
$\Psi_\alpha = \colimx \Psi_i$, 
$\Psi_\alpha^\prime = \colimx \Psi_i^\prime$, 
$\Psi_\alpha\pp = \colimx \Psi_i\pp$.
That 
$\Psi_{\alpha}^\prime$, $\Psi_{\alpha}\pp$ 
are in fact coproducts of $\leq \#(\alpha)$ representable cofunctors follows upon observing that 
$ \Psi_{\alpha}^\prime = \Psi_{\omega}^\prime  \oplus
\bigl\{\coprod\limits_{\omega \leq i < \alpha} \Psi_{i+1}^\prime/\Psi_{i}^\prime\bigr\}$, 
$ \Psi_{\alpha}\pp = \Psi_{\omega}\pp \oplus 
\bigl\{\coprod\limits_{\omega \leq i < \alpha} \Psi_{i+1}\pp/\Psi_{i}\pp\bigr\}$.  
Conclusion: \bC satisfies condition $\kappa$.]\\

\begingroup%%----------------------------------->>
\fontsize{9pt}{11pt}\selectfont
\textbf{\small LEMMA} \ 
Suppose that \bI is a countable filtered category $-$then $\exists$ a final functor $[\N] \ra \bI$.
\vspi
[One can find a directed set $(J,\leq)$ and a final functor $\bJ \ra \bI$ 
(cf. p. \pageref{15.31}).  
Since \bI is 
countable, so is \bJ (this fact is contained in the passage from \bI to \bJ 
(Cordier-Porter\footnote[2]{\textit{Shape Theory}, Ellis Horwood (1989), 42-44.})).  
Arrange the elements of \bJ in a sequence $j_0, j_1, \ldots$, and take $k_0 = j_0$, $k_n \geq k_{n-1}$, $j_n$ 
$(n \geq 1)$ to get a final functor $[\N] \ra \bJ$.]\\
\endgroup %%------------------------------------<<

\begingroup%%----------------------------------->>
\fontsize{9pt}{11pt}\selectfont
\textbf{\small EXAMPLE}  \ 
Consider $\bD(\bAMOD)$, where \mA is commutative and noetherian $-$then if $\bD(\bAMOD)$ admits Adams representability, 
every flat $A$-module has projective dimension $\leq 1$ 
(Neeman\footnote[2]{\textit{Topology} \textbf{36} (1997), 619-645.}).  
Example: Take $A = \C[x,y]$ $-$then the projective dimension of $\C(x,y)$ is 2, %dmc - verify we really want round brackets
therefore in this case $\bD(\bAMOD)$ does not admit Adams representability.
\vspi
[Note: Recalling the characterization of compact objects in $\bD(\bAMOD)$ mentioned  on
 p. \pageref{15.32}, Neeman's countability criterion implies that $\bD(\bAMOD)$ admits Adams representability provided that \mA is countable,]\\
\endgroup %%------------------------------------<<

Let \bC be a compactly generated triangulated category.  Suppose that \bD is a reflective subcategory of \bC, \mR a reflector for \bD.  Put $T = \iota \circx R$, where $\iota:\bD \ra \bC$ is the inclusion, 
%%----------------------------------------------------------------------------------------------30
and let $(S,D)$ be the associated orthogonal pair (cf. p. \pageref{15.33}) $-$then \mT is said to be a 
\un{localization functor}
\index{localization functor} 
if \mT is a triangulated functor.

[Note: \ The elements of \mS are the 
\un{\mT-equivalences}.
\index{T-equivalences}  
The elements of \mD (i.e., the \mX such that $\epsilon_X:X \ra TX$ is an isomorphism) are the 
\un{\mT-local}
\index{T-local (objects)} 
objects and the elements of ker \mT (i.e., the \mX such that $TX = 0$) are the 
\un{\mT-acyclic}
\index{T-acyclic (objects)} objects.]

Observation: If \mX is \mT-acyclic and if \mY is \mT-local, then 
$
\begin{cases}
\ \Mor(\Sigma^n X,Y) = 0\\
\ \Mor(\Omega^n X,Y) = 0
\end{cases}
(n \geq 0).
$
\\

\begin{proposition} \ %25
Let \bC be a compactly generated triangulated category.  Suppose that \mT is a localization functor $-$then 
$\forall \ X \in \Ob\bC$, $\exists$ an exact triangle 
$X_T \ra X \overset{\epsilon_X}{\ra} TX \ra \Sigma X_T$, where $X_T$ is \mT-acyclic.
\end{proposition}

[Place $X \overset{\epsilon_X}{\lra} TX$ in an exact triangle 
$X_T \ra $
$X \overset{\epsilon_X}{\lra}$ 
$TX \ra $
$\Sigma X_T$
and apply \mT to get an exact triangle 
$TX_T \ra $
$TX \overset{T\epsilon_X}{\lra}$ 
$T^2X \ra $
$\Sigma TX_T$.
Since $T\epsilon_X$ is an isomorphism, $TX_T = 0$.]\\

\begingroup%%----------------------------------->>
\fontsize{9pt}{11pt}\selectfont
The following lemma has been implicitly used in the proof of Proposition 25.\\
\endgroup %%------------------------------------<<

\begingroup%%----------------------------------->>
\fontsize{9pt}{11pt}\selectfont
\textbf{\small LEMMA} \ 
Let \bC be a triangulated category.  Suppose that 
$X \overset{u}{\ra} Y \overset{v}{\ra} Z \overset{w}{\ra} \Sigma X$ is an exact triangle, 
where $v$ is an isomorphism $-$then $X = 0$.
\vspi
[The triangle 
$Y \overset{v}{\ra} Z \ra 0 \ra \Sigma Y$ 
is exact (cf. p. \pageref{15.33a}), thus the triangle 
$0 \ra Y \overset{v}{\ra} Z \ra 0$ is exact (cf. p. \pageref{15.33b}) and 
$
\begin{cases}
\ v \circx u = 0 \ \implies \ u = 0\\
\ w \circx v = 0 \ \implies \ w = 0
\end{cases}
$
(cf. Proposition 3).  Therefore the diagram 
\begin{tikzcd}[sep=large]
{0} \ar{r} \ar{d} 
&{Y} \ar{r}{v} \arrow[d,shift right=0.5,dash] \arrow[d,shift right=-0.5,dash] 
&Z \ar{r} \arrow[d,shift right=0.5,dash] \arrow[d,shift right=-0.5,dash] 
&{0} \ar{d}\\
{X} \ar{r}[swap]{0} 
&{Y} \ar{r}[swap]{v} 
&{Z} \ar{r}[swap]{0} 
&{\Sigma X}
\end{tikzcd}
commutes, so $0 \ra X$ is an isomorphism (cf. p. \pageref{15.33c}).]\\
\endgroup %%------------------------------------<<

\begin{proposition} \ %26
Let \bC be a compactly generated triangulated category.  
Suppose that \mT is a localization functor $-$then the \mT-acyclic objects are the object class of a 
coreflective subcategory of \bC, the coreflector being the functor that sends \mX to $X_T$.
\end{proposition}

[Note: \ There is a natural isomorphism 
$(\Sigma X)_T \ra \Sigma X_T$ and $X \ra Y \ra Z \ra \Sigma X$ exact $\implies$ 
$X_T \ra Y_T \ra Z_T \ra \Sigma X_T$ exact.]\\

\begin{proposition} \ %27
Let \bC be a compactly generated triangulated category.  Suppose that \mT is a localization functor $-$then \mX is \mT-local iff 
$\Mor(Y,X) = 0$ for all \mT-acyclic \mY and \mX is acyclic iff $\Mor(X,Y) = 0$ for all \mT-local \mY.
\end{proposition}

[To see that the condition characterizes the \mT-local objects, take $Y = X_T$.  Thus the arrow $X_T \ra X$ is the zero morphism, so the isomorphism $(X_T)_T \ra X_T$ is the zero morphism, hence $X_T = 0$, which implies that 
$\epsilon_X:X \ra TX$ is an isomorphism.]\\

%%----------------------------------------------------------------------------------------------31
\begingroup%%----------------------------------->>
\fontsize{9pt}{11pt}\selectfont
Using the notation on p. \pageref{15.34}, take for $\sT$ the class of $T$-acyclic objects and take for $\sF$ the class of $T$-local objects 
$-$then Ann$_\tL \sF = \sT$ and Ann$_\tR \sT = \sF$ (cf. Proposition 27), i.e., the pair $(\sT,\sF)$ is a torsion theory on \bC.\\
\endgroup %%------------------------------------<<

\begin{proposition} \ %28
Let \bC be a compactly generated triangulated category.  
Suppose that \mT is a localization functor $-$then the class of $T$-local objects is the object class of a thick subcategory of \bC which is closed under the formation of products in \bC.
\end{proposition}

[Given an exact triangle $X \ra Y \ra Z \ra \Sigma X$, there is a commutative diagram 
\begin{tikzcd}%[sep=large]
{X} \ar{d}{\epsilon_X}\ar{r}
&{Y}\ar{d}{\epsilon_Y}  \ar{r}
&{Z} \ar{d}{\epsilon_Z} \ar{r}
&{\Sigma X}\ar{d}{\Sigma\epsilon_X}\\
{TX} \ar{r}
&{TY} \ar{r}
&{TZ} \ar{r}
&{\Sigma TX}
\end{tikzcd}
of exact triangles, thus if two 
$\epsilon_X$, 
$\epsilon_Y$, 
$\epsilon_Z$ are isomorphisms, so is the third (cf. p. \pageref{15.35}).  
Therefore \bD is a triangulated subcategory of \bC 
(cf. Proposition 7).  
Next, for any pair of morphisms $i:X \ra Y$, $r:Y \ra X$ with $r \circx i = \id_X$, there is a commutative diagram 
\begin{tikzcd}%[sep=large]
{X} \ar{d}{\epsilon_X} \ar{r}{i}
&{Y}\ar{d}{\epsilon_Y} \ar{r}{r}
&{X}\ar{d}{\epsilon_X} \\
{TX} \ar{r}[swap]{T_i}
&{TY} \ar{r}[swap]{T_r}
&{TX}
\end{tikzcd}
.  
Accordingly, $\epsilon_X$ is a retract of $\epsilon_Y$ 
(cf. p. \pageref{15.36}) and if $\epsilon_Y$ is an isomorphism, then the same is true of $\epsilon_X$, hence \bD is thick.]

[Note: \ Analogously, the class of $T$-acyclic objects is the object class of a thick subcategory of \bC which is closed under the formation of coproducts in \bC.]\\

Remark: \ \bD is not necessarily compactly generated.  In fact, there may be no nonzero compact objects in \bD at all.\\

\begingroup%%----------------------------------->>
\fontsize{9pt}{11pt}\selectfont
\textbf{\small EXAMPLE}  \ 
Suppose that \bC is a compactly generated triangulated category.  Let $\sK = \{K\}$ be a set of compact objects.  
Denote by \bK the thick subcategory generated by $\sK$ and denote by \bL the smallest triangulated subcategory of \bC containing $\sK$ and closed under the formation of coproducts in \bC $-$then \bK is a subcategory of \bL (via the Eilenberg swindle) and there is a localization functor $T_{\sK}$ whose acyclic objects are the objects of \bL.  
Moreover, every compact object in \bC which lies in \bL must lie in \bK.
\vspi
[Write $\ov{\sK} = \{\ov{K}\}$ for the set 
$\ds\bigcup\limits_K \{\Sigma^n K: n \geq 0\} \cup \ds\bigcup\limits_K \{\Omega^n K: n \geq 0\}$ 
and let $\sK^+$ be the class of objects in \bC that are coproducts of objects in $\ov{\sK}$ $-$then $\forall \ X \in \Ob\bC$, 
$\exists$ an object (\bX,\bff) in $\bFIL(\bC)$, completable in $\sK^+$ (obvious definition), and an arrow 
$\tel(\bX,\bff) \ra X$ such that $\Mor(Y,\tel(\bX,\bff)) \approx \Mor(Y,X)$ for all \mY in \bL 
(proceed as in the proof of the Brown representability theorem) (cf. Proposition 16)).  Taking 
$X_{\sK} = \tel(\bX,\bff)$, define $T_{\sK}X$ by the exact triangle 
$X_{\sK} \ra X \ra T_{\sK}X \ra \Sigma X_{\sK}$.]
\vspi
[Note: \ The $T_{\sK}$ are the 
\un{compact}
\index{compact (localization functor)} 
localization functors.]\\
\endgroup %%------------------------------------<<

\begingroup%%----------------------------------->>
\fontsize{9pt}{11pt}\selectfont
Let \bC be a compactly generated triangulated category $-$then a localization functor \mT is said to be 
%%----------------------------------------------------------------------------------------------32
\un{smashing}
\index{smashing (localization functor)} 
if it preserves coproducts or, equivalently, if \bD is closed under the formation of coproducts in \bC (recall Proposition 12).

Example: A compact localization functor is smashing.
\vspi
[Note: \ The 
\un{telescope conjecture}
\index{telescope conjecture} 
is said to hold for \bC if every smashing localization functor is compact.  
In the stable homotopy category, the telescope conjecture is false but in the derived category 
$\bD(\bAMOD$), where \mA is commutative and noetherian, the telescope conjecture is true.]\\
\endgroup %%------------------------------------<<

\begingroup%%----------------------------------->>
\fontsize{9pt}{11pt}\selectfont
\textbf{\small FACT} \ 
Suppose that \bC is a compactly generated triangulated category.  Let \mT be a localization functor $-$then \mT is smashing iff \mK compact in \bC $\implies$ $RK$ compact in \bD.\\
\endgroup %%------------------------------------<<

\begingroup%%----------------------------------->>
\fontsize{9pt}{11pt}\selectfont
Application: If \mT is smashing, then \bD is a compactly generated triangulated category.\\
\endgroup %%------------------------------------<<

\begingroup%%----------------------------------->>
\fontsize{9pt}{11pt}\selectfont
\textbf{\small FACT} \ 
Suppose that \bC admits Adams representability.  Let \mT be a localization functor $-$then \bD admits Adams representability provided that \mT is smashing.\\
\endgroup %%------------------------------------<<

Notation: Let \bC be a triangulated category with products.  
Suppose given an object $(\bX,\bff)$ in $\bTOW(\bC)$ $-$then 
$\Sigma \hspace{0.02cm} \mic(\bX,\bff)$ is any completion of 
$\prod\limits_n X_n \overset{\text{sf}}{\lra} \prod\limits_n X_n$ to an exact triangle (cf. TR$_3$), where 
$\pr_n \circx \text{sf} = \pr_n - f_n \circx \pr_{n+1}$.\\

\begingroup%%----------------------------------->>
\fontsize{9pt}{11pt}\selectfont
\textbf{\small EXAMPLE}  \ 
Suppose that \bC is a compactly generated triangulated category.  
Let \mT be a localization functor and let (\bX,\bff) be an object in $\bTOW(\bC)$ such that $\forall \ n$, $X_n$ is $T$-local $-$then mic(\bX,\bff) is $T$-local.\\
\endgroup %%------------------------------------<<

Let \bC be a compactly generated triangulated category.  
Suppose that $F:\bC \ra \bAB$ is an exact functor.  Let $S_F$ be the class of morphisms $X \overset{u}{\ra} Y$ such that 
$\forall \ n \geq 0$, 
$
\begin{cases}
\ F\Sigma^n u\\
\ F\Omega^n u
\end{cases}
$
is an isomorphism $-$then 
(1) $S_F$ admits a calculus of left and right fractions and contains the isomorphisms of \bC; 
(2) $u \in S_F$  $\implies$ $\Sigma u$ $\&$ $\Omega u \in S_F$; 
(3) $f, g \in S_F$ $\implies$ $\exists$ $h \in S_F$ (data as in TR$_5$); 
(4) $u \in S_F$ iff $\exists$ $f, g \in \Mor \bC$ : $u \circx f \in S_F$, $g \circx u \in S_F$.  
Therefore the metacategory 
$S_F^{-1}\bC$ is triangulated and $L_{S_F}:\bC \ra S_F^{-1}\bC$ is a triangulated functor.

[Note: \ In the terminology of p. \pageref{15.37}, $S_F$ is multiplicative.]\\

\begin{proposition} \ %29
Let \bC be a compactly generated triangulated category.  Suppose that $F:\bC \ra \bAB$ is an exact functor which converts coproducts into direct sums.  Assume: The metacategory $S_F^{-1} \bC$ is isomorphic to a category $-$then $S_F^\perp$ is the object class of a reflective subcategory of \bC.
\end{proposition}

[Argue as in the example on p. \pageref{15.38}.  Thus the triangulated functor 
$L_{S_F}:\bC \ra S_F^{-1}\bC$ preserves coproducts, so $\forall \ Y \in \Ob S_F^{-1}\bC$, $\Mor(L_{S_F}-,Y)$ is an 
exact cofunctor $\bC \ra \bAB$
%%----------------------------------------------------------------------------------------------33
which converts coproducts into products, hence by the Brown representability theorem, $\exists$ $Y_{S_F} \in \Ob\bC$ : 
$\Mor(L_{S_F}X,Y) \approx \Mor(X,Y_{S_F})$.]

\label{15.47}
[Note: \ The procedure generates an idempotent triple $\bT_F = (T_F,m,\epsilon)$ in \bC ($T_F:\bC \ra \bC$ is a localization functor, $S_F$ is the class of $T_F$-equivalences, and $O_F = \ker T_F$ (i.e., \mX is $T_F$-acyclic iff 
$\forall \ n \geq 0$, 
$
\begin{cases}
\ F\Sigma^n X = 0\\
\ F\Omega^n X = 0
\end{cases}
$
(cf. p. \pageref{15.39}))).]\\

Maintaining the assumption that \bC is a compactly generated triangulated category, given any $X \in \Ob\bC$, put 
$\kappa_X = \sum\limits_{\ov{U}} \#(\Mor(\ov{U},X))$ and for $\kappa$ an infinite cardinal 
$\geq \kappa_{\sU} \equiv \sum\limits_{U} \kappa_U$, let $\bC_\kappa$ be the full subcategory of \bC whose objects are the \mX such that $\kappa_X \leq \kappa$ $-$then $\bC_\kappa$ is a thick subcategory of \bC which is closed under the formation of coproducts in \bC indexed by sets of cardinality $\leq \kappa$ and $\bC = \bigcup\limits_\kappa \bC_\kappa$.

[Note: \ $\bC_\kappa$ contains $\sU$, hence $\bC_\kappa$ contains $\cptx\bC$ 
(by the theorem of Neeman-Ravenel).]

[Notation: $\sU_\kappa^+$ is the class of objects in \bC that are coproducts of $\leq \kappa$ objects in $\ov{\sU}$.\\

\textbf{\small LEMMA} \ 
Let $\{G_n\}$ be a sequence of abelian groups.  Assume: $\forall \ n$, $\#(G_n) \leq \kappa$, where $\kappa$ is an infinite cardinal $-$then the cardinality of $\bigoplus\limits_n G_n$ is bounded by $\kappa$.

[Note: \ Another triviality is the fact that if $G^\prime \ra G \ra G\pp$ is an exact sequence of abelian groups and if 
$\#(G^\prime) \leq \kappa$, 
$\#(G\pp) \leq \kappa$, where $\kappa$ is an infinite cardinal, then $\#(G) \leq \kappa^2 = \kappa$.]\\

\begin{proposition} \ %30
Let \bC be a compactly generated triangulated category.  
Fix an infinite cardinal $\kappa \geq \kappa_{\sU}$ $-$then 
$X \in \Ob\bC_\kappa$ iff $X \approx \tel(\bX,\bff)$ where $(\bX,\bff)$ is completable in $\sU_\kappa^+$.
\end{proposition}

[The sufficiency is clear (cf. Proposition 13) and the necessity can be established by reworking the proof of Proposition 16 (with $F = \Mor(-,X))$.]

[Note: \ It is a corollary that $\bC_\kappa$ has a small skeleton in $\ov{\bC}_\kappa$.]\\

\textbf{\small LEMMA} \ 
Let \bC be a compactly generated triangulated category.  Suppose that $F:\bC \ra \bAB$ is an exact functor which converts coproducts into direct sums.  Put 
$H = \bigoplus\limits_{n \geq 0} F \circx \Sigma^n \oplus \bigoplus\limits_{n > 0} F \circx \Omega^n$ 
$-$then $H:\bC \ra \bAB$ is an exact functor which converts coproducts into direct sums and a morphism 
$X \overset{u}{\ra} Y$ is in $S_F$ iff $Hu:HX \ra HY$ is an isomorphism.\\

\begin{proposition} \ %31
Let \bC be a compactly generated triangulated category.  Suppose that $F:\bC \ra \bAB$ is an exact functor which converts coproducts into direct sums $-$then $\forall$ infinite cardinal $\kappa \gg \kappa_{\sU}$, $\exists$ an infinite cardinal 
$\delta(\kappa) \ge \kappa$ such that $\forall \ Y$: $\#(HY) \leq \kappa$, $\exists$ $X \in \Ob\bC_{\delta(\kappa)}$ $\&$ 
$X \overset{u}{\ra} Y$ with $Hu:HX \ra HY$ an isomorphism.
\end{proposition}

%%----------------------------------------------------------------------------------------------34
[Bearing in mind that $\cptx\bC$ has a small skeleton $\ov{\cptx\bC}$ (cf. Proposition 19), fix an infinite cardinal 
$\kappa_H > \sup\{\#(H\ov{K}): \ov{K} \in \Ob \ov{\cptx\bC}\}$ and take 
$\kappa = \delta_0(\kappa) > \max\{\kappa_H,\kappa_{\sU}\}$.  
Since 
$HY \approx \underset{Y}{\colimx} HL$ 
(cf. p. \pageref{15.40}) $\forall \ y \in HY$, $\exists$ an object 
$L \ra Y$ in $\ov{\bLambda}/Y$: $y \in \im(HL \ra HY)$.  
Therefore one can choose objects $L_i \ra Y$ in 
$\ov{\bLambda}/Y$ indexed by a set \mI of cardinality $\leq \delta_0(\kappa)$ such that 
$Hu_0:HX_0 \ra HY$ is surjective.  Here $X_0 = \coprod\limits_I L_i$ and $u_0:X_0 \ra Y$ is the coproduct of the $L_i \ra Y$.  
Because the $L_i$ are compact and $\#(I) \leq \delta_0(\kappa)$, $X_0 \in \Ob\bC_{\delta_0(\kappa)}$.  
Embed $X_0 \overset{u_0}{\lra} Y$ in an exact triangle 
$Y^\prime \overset{u^\prime}{\ra} X_0 \overset{u_0}{\lra} Y \ra \Sigma Y^\prime$.  
Claim: $\exists$ an infinite cardinal $\delta_1(\kappa) \geq \delta_0(\kappa)$ for which $\#(HY^\prime) \leq \delta_1(\kappa)$ 
independently of the choices (i.e., the bound is a function only of the initial supposition that $\#(HY) \leq \kappa$).  
To see this, note that 
$\#(H\Sigma^n Y) \leq \kappa$, 
$\#(H\Omega^n Y) \leq \kappa$, and 
$\#(H\Sigma^n X_0) \leq \kappa_H^\kappa$,
$\#(H\Omega^n X_0) \leq \kappa_H^\kappa$
and use the long exact sequence generated by $H$.  
Repeat the process: 
$u_0^\prime:\coprod\limits_{I^\prime} L_{i^\prime}^\prime \ra Y^\prime$ 
$(\#(I^\prime) \leq \delta_1(\kappa))$ and place $u^\prime \circx u_0^\prime$ in an exact triangle 
$Z \ra $ 
\begin{tikzcd}%[sep=large]
{\coprod\limits_{I^\prime} L_{i^\prime}^\prime} \ar{r}{u^\prime \circx u_0^\prime} &{X_0}
\end{tikzcd}
$\ra \Sigma Z$.  
Consider now the diagram 
\begin{tikzcd}%[sep=large]
{X_0} \arrow[d,shift right=0.45,dash] \arrow[d,shift right=-0.45,dash]  \ar{r}
&{\Sigma Z} \ar{r}
&{\coprod\limits_{I^\prime} \Sigma L_{i^\prime}^\prime} \ar{d}{\Sigma u_0^\prime}
\ar{rr}{-\Sigma (u^\prime \circx u_0^\prime)}
&&{\Sigma X_0} \arrow[d,shift right=0.45,dash] \arrow[d,shift right=-0.5,dash]\\
{X_0}  \ar{r}[swap]{u_0}
&{Y} \ar{r} 
&{\Sigma Y^\prime} \ar{rr}[swap]{-\Sigma u^\prime}
&&{\Sigma X_0}
\end{tikzcd}
.  
The rows being in $\Delta$, one can find a filler $u_1:\Sigma Z \ra Y$ (cf. Proposition 2).  
Put $X_1 = \Sigma Z$ (thus $X_1 \in \Ob\bC_{\delta_1(\kappa)}$) 
and let $f_0$ be the arrow $X_0 \ra X_1$.  By construction 
\begin{tikzcd}%[sep=large]
{X_0} \ar{d}[swap]{u_0} \ar{r}{f_0} &{X_{1}}\ar{d}{u_1}\\
Y \arrow[r,shift right=0.45,dash] \arrow[r,shift right=-0.45,dash]  &{Y}
\end{tikzcd}
commutes and $\ker H f_0 = \ker H u_0$.  Continuing, one produces $\forall \ n$ a commutative diagram
\begin{tikzcd}%[sep=large]
{X_n} \ar{d}[swap]{u_n} \ar{r}{f_n} &{X_{n+1}}\ar{d}{u_{n+1}}\\
Y \arrow[r,shift right=0.45,dash] \arrow[r,shift right=-0.45,dash]  &{Y}
\end{tikzcd}
, 
where $\ker H f_n = \ker H u_n$ and $X_n \in \Ob\bC_{\delta_n(\kappa)}$ $(\delta_n(\kappa) \leq \delta_{n+1}(\kappa))$.  
Definition: $X = \tel(\bX,\bff)$ $-$then $X \in \Ob\bC_{\delta(\kappa)}$
$(\delta(\kappa) \geq (\sup\{\delta_n(\kappa)\})^\omega$ (cf. infra)), 
$HX \approx \colimx HX_n$ and there is an arrow $X \overset{u}{\ra} Y$ with $Hu:HX \ra HY$
 an isomorphism (injectivity from the condition on the kernels, surjectivity from the surjectivity of $Hu_0$).]\\

\begingroup%%----------------------------------->>
\fontsize{9pt}{11pt}\selectfont
Thanks to Proposition 13, $\forall$  $\ov{U} \in \ov{\sU}$, 
$\colimx\Mor(\ov{U},X_n) \approx \Mor(\ov{U},\telsub(\bX,\bff))$, 
hence 
$\#(\Mor(\ov{U},\tel(\bX,$ $\bff))) \leq \ds\prod\limits_n \#(\Mor(\ov{U},X_n)) 
\leq \ds\prod\limits_n \delta_n(\kappa) \leq (\sup\{\delta_n(\kappa)\})^\omega$.\\
\endgroup %%------------------------------------<<

\index{Bousfield-Margolis Localization Theorem}
\index{Theorem: Bousfield-Margolis Localization Theorem}
\textbf{\small BOUSFIELD-MARGOLIS LOCALIZATION THEOREM} \quad
Let \bC be a compactly generated triangulated category.  Suppose that $F:\bC \ra \bAB$ is an exact functor which converts coproducts into direct sums $-$then there exists a localization functor $T_F$ such that $S_F^\perp$ is the class of $T_F$-local objects.

%%----------------------------------------------------------------------------------------------35
[In view of Proposition 29, the point is to show that the metacategory $S_F^{-1}$\bC is isomorphic to a category.  
Thus fix $X, Y \in \Ob S_F^{-1}\bC$ $(= \Ob\bC)$ and $\kappa \gg \kappa_{\sU}$ : $X, Y \in \Ob\bC_\kappa$ $\&$ 
$\#(HX) \leq \kappa$, $\#(HY) \leq \kappa$.  
By definition, $\Mor(X,Y)$ is a conglomerate of equivalence classes of pairs 
$(s,f)$ : $X \overset{f}{\ra} Y^\prime \overset{s}{\la} Y$ 
(cf. p. \pageref{15.41}).  
Given such a pair $(s,f)$, consider an exact triangle 
$Z \ra X \amalg Y \ra Y^\prime \ra \Sigma Z$.  Since $HY \approx HY^\prime$, $\#(HZ) \leq \kappa$.  
Using Proposition 31, 
choose $W \in \Ob\bC_{\delta(\kappa)}$ $\&$ $W \overset{u}{\ra} Z$ with $Hu:HW \ra HZ$ an isomorphism.  
There is a diagram 
\begin{tikzcd}%[sep=large]
{W} \ar{d}{u}  \ar{r}
&{X \amalg Y} \arrow[d,shift right=0.5,dash] \arrow[d,shift right=-0.5,dash] \ar{r}{\pi\pp}  
&{Y\pp} \ar{r}
&{\Sigma W}\ar{d}{\Sigma u}\\
{Z} \ar{r}  
&{X \amalg Y}  \ar{r}[swap]{\pi^\prime}
&{Y^\prime} \ar{r} 
&{\Sigma Z}
\end{tikzcd}
and a filler $\phi:Y\pp \ra Y^\prime$ (cf. TR$_5$) which is necessarily in $S_F$.  
Note too that 
$Y\pp \in \Ob\bC_{\delta(\kappa)}$.  Put $g = \pi\pp \circx \ini_X$, $t = \pi\pp \circx \ini_Y$ $-$then 
$\phi \circx g = f$, $\phi \circx t= s$, and $t \in S_F$, so the pair $(s,f)$ is equivalent to the pair $(t,g)$.  But 
$\bC_{\delta(\kappa)}$ has a small skeleton $\ov{\bC}_{\delta(\kappa)}$ 
(cf. Proposition 30) and there is just a set of diagrams of the form
$X \overset {\bar{g}}{\ra} \ov{Y}\pp \overset{\bar{t}}{\la} Y$, 
where $\ov{Y}\pp \ \in \Ob\ov{\bC}_{\delta(\kappa)}$.]\\

\begingroup%%----------------------------------->>
\fontsize{9pt}{11pt}\selectfont
\textbf{\small EXAMPLE}  \  
Take for \bC the stable homotopy category \bHSPEC and fix an $\bX \in \Ob\bC$ $-$then 
$H_{\bX}(\bY) = [\bS^0,\bX \wedge \bY]$ is an exact functor $\bC \ra \bAB$
 which converts coproducts into direct sums and by the Bousfield-Margolis localization theorem, 
$S_{\bX}^\perp$ is the object class of a reflective subcategory of \bC, where 
$S_{\bX}$ is the class of morphisms $\bY^\prime \ra \bY\pp$ such that $\forall \ n \in \Z$, 
$[\bS^n,\bX \wedge \bY^\prime] \approx [\bS^n,\bX \wedge \bY\pp]$.\\
\endgroup %%------------------------------------<<

Given a closed category \bC, the 
\un{dual}
\index{dual (of an object in a closed category)} 
$D X$ of an object \mX is $\hom(X,e)$.

\indent\indent (DU$_1$) \ $\forall \ X, X^\prime \in \Ob\bC$, $\exists$ a natural morphism 
$DX \otimes DX^\prime \ra D(X \otimes X^\prime)$.

[In the pairing 
$\hom(X,Y) \otimes \hom(X^\prime,Y^\prime) \ra \hom(X \otimes X^\prime, Y \otimes Y^\prime)$, specialize and take 
$Y = e$, $Y^\prime = e$.]

\indent\indent (DU$_2$) \ $\forall \ X \in \Ob\bC$, $\exists$ a natural morphism $X \ra D^2 X$.

[
$\Mor(X,D^2X) \approx \Mor(X,\hom(DX,e))$ 
$\approx \Mor(X\otimes DX,e)$ 
$\approx \Mor( DX,\hom(X,e))$ 
$\approx \Mor(DX,DX)$.]\\

\textbf{\small LEMMA} \ 
Suppose that \bC is a closed category $-$then there is an arrow 
$\hom(X,Y) \otimes Z \ra \hom(X,Y \otimes Z)$ natural in \mX, \mY, \mZ.\\

Given a closed category \bC, an object \mX is said to be 
\un{dualizable}
\index{dualizable (object in a closed category)}
if $\forall \ Y \in \Ob\bC$, the arrow $DX \otimes Y \ra \hom(X,Y)$ is an isomorphism.  
Example: $e$ is dualizable.

[Note: \ When \mX is dualizable, $DX \otimes -$ is a right adjoint for $-\otimes X$, hence 
$DX \otimes - \approx \hom(X,-)$.]\\

\begingroup%%----------------------------------->>
\fontsize{9pt}{11pt}\selectfont
\textbf{\small EXAMPLE}  \  
Let \mA be a commutative ring with unit $-$then an object \mX in \bAMOD is dualizable iff \mX is finitely generated and projective.\\
\endgroup %%------------------------------------<<

%%----------------------------------------------------------------------------------------------36
\begingroup%%----------------------------------->>
\fontsize{9pt}{11pt}\selectfont
Let \bC be a closed category $-$then an object \mX in \bC is 
\un{invertible}
\index{invertible (object in a closed category)} 
if there is an object $X^{-1}$ in \bC and an isomorphism 
$X \otimes X^{-1} \ra e$.\\
\endgroup %%------------------------------------<<

\begingroup%%----------------------------------->>
\fontsize{9pt}{11pt}\selectfont
\textbf{\small FACT} \  
Every invertible element \mX in \bC is dualizable and $DX \approx X^{-1}$.\\
\endgroup %%------------------------------------<<

\begin{proposition} \ %32
Suppose that \bC is a closed category.  Assume: \mX is dualizable $-$then $D X$ is dualizable and the morphism 
$X \ra D^2 X$ is an isomorphism.\\
\end{proposition}

Remark: If \bC has coproducts, then $\forall \ Y$, 
$\coprod\limits_i Y \otimes X_i \approx Y \otimes \coprod\limits_i X_i$.  If \bC has products, then $\forall$ dualizable \mX, 
$X \otimes \prod\limits_i Y_i \approx \prod\limits_i X \otimes Y_i$.  
Proof: 
$X \otimes \prod\limits_i Y_i \approx D^2 X \otimes \prod\limits_i Y_i$ 
$\approx \hom(DX,\prod\limits_i Y_i)$ 
$\approx \prod\limits_i \hom(DX,Y_i)$
$\approx \prod\limits_i D^2 X \otimes Y_i$ 
$\approx \prod\limits_i X \otimes Y_i$.\\

\textbf{\small LEMMA} \ 
Suppose that \bC is a closed category $-$then the pairing 
$\hom(X,Y) \otimes \hom(X^\prime,Y^\prime) \ra \hom(X \otimes X^\prime,Y \otimes Y^\prime)$ is an isomorphism if $X$ and $X^\prime$ are dualizable or if \mX $(X^\prime)$ is dualizable and $Y = e$ $(Y^\prime = e)$.\\

\begin{proposition} \ %33
Suppose that \bC is a closed category $-$then \mX, $X^\prime$ dualizable $\implies$ $X \otimes X^\prime$ dualizable.
\end{proposition}

[$\forall \ Y$, 
$D(X \otimes X^\prime) \otimes Y \approx$ 
$DX \otimes DX^\prime \otimes Y \approx$ 
$DX \otimes \hom(X^\prime,Y) \approx$ 
$\hom(X,\hom(X^\prime,Y)) \approx$ 
$\hom(X \otimes X^\prime,Y)$.]\\

\textbf{\small LEMMA} \ 
Suppose that \bC is a closed category $-$then the arrow $\hom(X,Y) \otimes Z \ra \hom(X,Y\otimes Z)$ is an isomorphism if either \mX or \mZ is dualizable.\\

\begin{proposition} \ %34
Suppose that \bC is a closed category $-$then \mX, $X^\prime$ dualizable $\implies$ $\hom(X,X^\prime)$ dualizable.
\end{proposition}

[$\forall \ Y$, 
$D\hom(X,X^\prime) \otimes Y \approx$ 
$\hom(\hom(X,X^\prime),e) \otimes Y \approx$ 
$\hom(DX \otimes X^\prime,e) \otimes Y \approx$ 
$\hom(DX,\hom(X^\prime,e)) \otimes Y \approx$ 
$\hom(DX,DX^\prime)\otimes Y \approx$
$\hom(DX,DX^\prime \otimes Y) \approx$
$\hom(DX,\hom(X^\prime,\\Y)) \approx$ 
$\hom(DX \otimes X^\prime, Y) \approx$
$\hom(\hom(X,X^\prime)Y)$.]\\

\begingroup%%----------------------------------->>
\fontsize{9pt}{11pt}\selectfont
\label{15.45} 
\label{15.50} 
\textbf{\small FACT} \  
Let \bC be a closed category.  Assume: \mX is dualizable $-$then \mX is a retract of $X \otimes DX \otimes X$.\\
\endgroup %%------------------------------------<<

Let \bC be a category with finite coproducts.  Assume: \bC is closed and triangulated $-$then \bC is said to be a 
\un{closed triangulated category}
\index{closed triangulated category} 
(CTC
\index{CTC}
) if there is a natural isomorphism 
$\zeta$, where 
$\zeta_{X,Y}:\Sigma X \otimes Y \ra \Sigma(X \otimes Y)$, subject to the following assumptions.

[Note: \ From the existence of $\zeta$, one derives the existence of a natural isomorphism $\eta$, where 
$\eta_{X,Y}:\Omega \hom(X,Y) \ra \hom(\Sigma X,Y)$.]\\

%%----------------------------------------------------------------------------------------------37
\indent\indent (CTC$_1$) \ The diagram
\begin{tikzcd}%[sep=large]
{\Sigma X \otimes e} \ar{rd}[swap]{R_{\Sigma X }}  \ar{r}{\zeta_{X,e}}  
&{\Sigma (X \otimes e)} \ar{d}{\Sigma R_X}\\
&{\Sigma X}
\end{tikzcd}
commutes.

\indent\indent (CTC$_2$) \ 
The diagram
\[
\begin{tikzcd}%[sep=large]
{(\Sigma X \otimes Y) \otimes Z}   \ar{rr}{\zeta_{X,Y} \otimes \id_Z}  
&&{\Sigma (X \otimes Y) \otimes Z} \ar{rr}{\zeta_{X \otimes Y,Z}}
&&{\Sigma ((X \otimes Y) \otimes Z)}\\
{\Sigma X \otimes (Y \otimes Z)} \ar{u}{A} \ar{rrrr}[swap]{\zeta_{X,Y \otimes Z}}
&&&&{\Sigma (X \otimes (Y \otimes Z))} \ar{u}[swap]{\Sigma A}
\end{tikzcd}
\]
commutes.

\indent\indent (CTC$_3$) \ If 
$X \overset{u}{\ra} Y \overset{v}{\ra}Z \overset{w}{\ra} \Sigma X$ 
is an exact triangle, then $\forall \ W \in \Ob\bC$, the triangle
\begin{tikzcd}%[sep=large]
{X \otimes W} \ar{rr}{u \otimes \id_W}
&&{Y \otimes W} \ar{rr}{v \otimes \id_W}
&&{Z \otimes W} \ar{rrr}{\zeta_{X,W} \circx (w \otimes \id_W)}
&&&{\Sigma (X  \otimes W)}
\end{tikzcd}
is exact.

\indent\indent (CTC$_4$) \ If 
$X \overset{u}{\ra} Y \overset{v}{\ra}Z \overset{w}{\ra} \Sigma X$ 
is an exact triangle, then $\forall \ W \in \Ob\bC$, the triangle
\begin{tikzcd}%[sep=small]
{\Omega \hom(X,W)}  \ar{rrr}{-(w^* \circx \eta_{X,W})}
&&&{\hom(Z,W)} \ar{r}{v^*} \ar{r}{v^*}
&{\hom(Y,W)} \ar{rrr}{\nu_{\hom(X,W)}^{-1} \circx u^*}
&&&{\Sigma \Omega \hom(X,W)} 
\end{tikzcd}
is exact.\\
\indent\indent (CTC$_5$) \ The diagram
\begin{tikzcd}%[sep=large]
{\Sigma e \otimes \Sigma e} \ar{d}[swap]{\Tee} \ar{r}{\approx}
&{\Sigma^2 e} \ar{d}{-1}\\
{\Sigma e \otimes \Sigma e} \ar{r}[swap]{\approx}
&{\Sigma^2 e}
\end{tikzcd}
commutes.

Remarks: 
(1) If 
$\ov{\zeta}_{e,X}$ is the composite $\Sigma \Tee_{X,e} \circx \zeta_{X,e} \circx \Tee_{e,\Sigma X}$, 
then the diagram \\
\begin{tikzcd}%[sep=large]
{e \otimes \Sigma X} \ar{d}[swap]{L_{\Sigma X }}  \ar{r}{\ov{\zeta}_{e,X}}  
&{\Sigma (e \otimes X)} \ar{ld}{\Sigma L_X}\\
{\Sigma X}
\end{tikzcd}
commutes;
(2) The additive functor 
$-\otimes W:\bC \ra \bC$ is a triangulated functor (this is the content of CTC$_3$); 
(3) The additive functor 
$\hom(-,W):\bC \ra \bC^\OP$  is a triangulated functor (this is the content of CTC$_4$); 
(4) If $m,n \in \N$, then the diagram 
\begin{tikzcd}%[sep=large]
{\Sigma^m e \otimes \Sigma^n e} \ar{d}[swap]{\Tee} \ar{r}{\approx}
&{\Sigma^{m+n} e} \ar{d}{(-1)^{mn}}\\
{\Sigma^n e \otimes \Sigma^m e} \ar{r}[swap]{\approx}
&{\Sigma^{m+n} e}
\end{tikzcd}
commutes.

Example: $D: \bC \ra \bC^\OP$ is a triangulated functor.\\

\begingroup%%----------------------------------->>
\fontsize{9pt}{11pt}\selectfont
Since the additive functor $\hom(W,-):\bC \ra \bC$ is a right adjoint for $-\otimes W$, it is necessarily triangulated 
(cf. p. \pageref{15.42}).\\
\endgroup %%------------------------------------<<

Notation: du\bC is the full, isomorphism closed subcategory of \bC whose objects are dualizable.\\

%%----------------------------------------------------------------------------------------------38
\begin{proposition} \ %35
Let \bC be a CTC $-$then $\dux\bC$ is a thick subcategory of \bC.
\end{proposition}

[Observe first that 0 is dualizable.  This said, take any morphism $X \overset{u}{\ra} Y$ in $\dux\bC$ and complete it to an exact triangle 
$X \overset{u}{\ra} $
$Y \overset{v}{\ra} $
$Z \overset{w}{\ra} $
$\Sigma X$
(cf. TR$_3$) $-$then $\forall \ W \in \Ob\bC$, there is a commutative diagram
\[
\begin{tikzcd}%[sep=large]
{\Omega DX \otimes W} \ar{d} \ar{r}
&{DZ \otimes W} \ar{d}\ar{r}
&{DY \otimes W} \ar{d} \ar{r}
&{\Sigma(\Omega DX \otimes W)} \ar{d}\\
{\Omega \hom(X,W)} \ar{r}
&{\hom(Z,W)} \ar{r}
&{\hom(Y,W)} \ar{r}
&{\Sigma\Omega\hom(X,W)}
\end{tikzcd}
,
\]
where, by CTC$_3$ $\&$ CTC$_4$, the rows are exact.  
Specialized to the case $X = X$, $Y = X$, $Z = 0$, $u = \id_X$ 
(cf. TR$_2$), it follows that the arrow 
$\Omega DX \otimes W \ra \Omega \hom(X,W)$ is an isomorphism 
(cf. p. \pageref{15.43}) i.e., that the arrow 
$\hom(\Sigma X,e) \otimes W \ra \hom(\Sigma X,W)$ is an isomorphism, so \mX dualizable $\implies$ $\Sigma X$ dualizable.  Next, \mX dualizable $\implies$ $\Omega X$ dualizable.  
Proof: 
$X \approx \hom(e,X)$ $\implies$ $\Omega X \approx\Omega \hom(e,X)$ $\approx \hom(\Sigma e,X)$ and $e$ dualizable 
$\implies$ $\Sigma e$ dualizable, hence Proposition 34 is applicable.  
Returning to $X \overset{u}{\ra} Y$, one concludes that the arrow 
$DZ \otimes W \ra \hom(Z,W)$ is an isomorphism 
(cf. p. \pageref{15.44}), thus \mZ is dualizable.  
Therefore $\dux\bC$ is a triangulated subcategory of \bC.  
Finally, suppose given a pair of morphisms $i:X \ra Y$, $r:Y \ra X$ with $r \circx i = \id_X$ and Y dualizable 
$-$then $\forall \ W \in \Ob\bC$, there is a commutative diagram 
\begin{tikzcd}%[sep=large]
{DX \otimes W} \ar{d} \ar{r}{Dr}
&{DY \otimes W} \ar{d} \ar{r}{Di}
&{DX \otimes W} \ar{d}\\
{ \hom(X,W)} \ar{r}[swap]{r^*}
&{\hom(Y,W)}\ar{r}[swap]{i^*}
&{\hom(X,W)}
\end{tikzcd}
, which shows that the arrow 
$DX \otimes W \ra \hom(X,W)$ is a retract of the arrow 
$DY \otimes W \ra \hom(Y,W)$.  
But the retract of an isomorphism is an isomorphism and this means that \mX is dualizable.  
Therefore $\dux\bC$ is a thick subcategory of \bC.]\\

\begingroup%%----------------------------------->>
\fontsize{9pt}{11pt}\selectfont
\textbf{\small EXAMPLE}  \ 
Let \bC be a CTC $-$then $e$ dualizable $\implies$ $\Sigma e$ dualizable and 
$D\Sigma e =$ 
$ \hom(\Sigma e,e) \approx$ 
$\Omega \hom(e,e) \approx$ 
$\Omega e$.  
Therefore 
$\Mor(Y,X \otimes \Omega e) \approx$ 
$\Mor(Y,D \Sigma e \otimes X) \approx$ 
$\Mor(Y,\hom(\Sigma e,X)) \approx$ 
$\Mor(Y \otimes \Sigma e,X) \approx$ 
$\Mor(\Sigma Y, X) \approx$ 
$\Mor(Y,\Omega X)$ 
$\implies$ 
$X \otimes \Omega e \approx \Omega X$.  
Consequently, 
$\hom(\Sigma X,Y) \approx$ 
$\hom(X,\hom(\Sigma e,$ $Y)) \approx$ 
$\hom(X,D\Sigma e \otimes Y) \approx$ 
$\hom(X,\Omega e \otimes Y) \approx$ 
$\hom(X,\Omega Y)$.\\
\endgroup %%------------------------------------<<

Suppose that \bC is a CTC $-$then \bC is said to be a 
\un{compactly generated CTC}
\index{compactly generated CTC} 
if \bC is compactly generated and every $U \in \sU$ is dualizable.\\

\begin{proposition} \ %36
Let \bC be a compactly generated CTC $-$then \mX compact $\implies$ \mX dualizable.
\end{proposition}

[The thick subcategory generated by $\sU$ is $\cptx\bC$ (theorem of Neeman-Ravenel).  
On the other hand, $\dux\bC$ is thick 
(cf. Proposition 35) and contains $\sU$.]\\
\vspace{0.25cm}

%%----------------------------------------------------------------------------------------------39
\begingroup%%----------------------------------->>
\fontsize{9pt}{11pt}\selectfont
\textbf{\small FACT} \ 
Suppose that \bC is a compactly generated CTC $-$then \mX is dualizable iff $\forall$ collection $\{X_i\}$ of objects in \bC, 
the arrow 
$\coprod\limits_i \hom(X,X_i) \ra \hom\bigl(X,\coprod\limits_iX_i\bigr)$ is an isomorphism.
\vspi
[Necessity: 
$\ds\coprod\limits_i \hom(X,X_i) \approx$ 
$\ds\coprod\limits_i DX \otimes X_i \approx$ 
$DX \otimes \ds\coprod\limits_i X_i \approx$ 
$\hom(X,\ds\coprod\limits_i X_i)$.
\vspi
Sufficiency: Let \bD be the full, isomorphism closed subcategory of \bC consisting of those \mY for which the arrow 
$DX \otimes Y \ra \hom(X,Y)$ is an isomorphism $-$then \bD is triangulated and closed under the formation of coproducts in \bC.  
Moreover, \bD contains all the dualizable objects, so $\sU \subset \Ob\bD$.  Therefore $\bD = \bC$ (cf Proposition 14).]\\
\endgroup %%------------------------------------<<

\textbf{\small LEMMA} \ 
Let \bC be a CTC with coproducts $-$then \mX compact and \mY dualizable $\implies$ $X \otimes Y$ compact.

[
$\bigoplus\limits_i \Mor(X \otimes Y,Z_i) \approx$ 
$\bigoplus\limits_i \Mor(X,\hom(Y,Z_i)) \approx$ 
$\bigoplus\limits_i \Mor(X, DY \otimes Z_i) \approx$ 
$\Mor(X,\coprod\limits_i DY$ $\otimes Z_i) \approx$ 
$\Mor(X,DY \otimes \coprod\limits_i Z_i) \approx$ 
$\Mor(X \otimes Y,\coprod\limits_i Z_i)$.]\\

\label{15.46}
Application: Let \bC be a compactly generated CTC $-$then \mX compact $\implies$ $DX$ compact.

[\mX is dualizable (cf. Proposition 36), so $DX$ is dualizable (cf. Proposition 32), hence $DX$ is a retract of 
$DX \otimes D^2 X \otimes DX$ (cf. p. \pageref{15.45} or still, is a retract of 
$DX \otimes X \otimes DX$ (cf. Proposition 32) and the lemma implies that 
$DX \otimes X \otimes DX$ is compact.]\\

Suppose that \bC is a compactly generated CTC $-$then \bC is said to be 
\un{unital}
\index{unital (compactly generated CTC)} 
provided that $e$ is compact.\\

\begin{proposition} \ %37
Let \bC be a unital compactly generated CTC $-$then \mX dualizable $\implies$ \mX compact.
\end{proposition}

[By the lemma, $e \otimes X$ is compact.]\\

\label{16.24}
Consequently, in a unital compactly generated CTC, ``compact'' = ``dualizable''.\\

\begingroup%%----------------------------------->>
\fontsize{9pt}{11pt}\selectfont
The stable homotopy category is a unital compactly generated CTC.\\
\endgroup %%------------------------------------<<

\begingroup%%----------------------------------->>
\fontsize{9pt}{11pt}\selectfont
\textbf{\small EXAMPLE}  \ 
Let \mA be a commutative ring with unit $-$then $\bD(\bAMOD)$ is a unital compactly generated CTC 
(B\"okstedt-Neeman\footnote[2]{\textit{Compositio Math.} \textbf{86} (1993), 209-234.}).\\
\endgroup %%------------------------------------<<

Suppose that \bC is a compactly generated CTC $-$then a 
\un{cohomology theory}
\index{cohomology theory (compactly generated CTC)} 
is an exact cofunctor 
$F:\bC \ra \bAB$ which converts coproducts into products and a 
\un{homology theory}
\index{homology theory (compactly generated CTC)}
 is
%%----------------------------------------------------------------------------------------------40
an exact functor 
$F:\bC \ra \bAB$ which converts coproducts into direct sums.  According to the Brown representability theorem, every cohomology theory is representable.  
The situation for homology theories is different.  
Put 
$H_e(X) = \underset{X}{\colimx} \Mor(e,K)$ $(= \ov{\restr{\Mor(e,-)}{\bLambda}} X$) and 
$H_X(Y) = H_e(X \otimes Y)$ $(X, Y \in \Ob\bC)$.  
Proposition 18 guarantees that $H_e$ is a homology theory, thus $H_X$ is also a homology theory 
(cf. CTC$_3$), and there is an arrow 
$H_X(Y) \ra \Mor(e,X \otimes Y)$.

[Note: \ When \bC is unital, $H_X(Y) \approx \Mor(e,X \otimes Y)$.]\\

\textbf{\small LEMMA} \  
The arrow $H_X(Y) \ra \Mor(e,X \otimes Y)$ is an isomorphism if \mX is compact.

[\mX compact $\implies$ \mX dualizable (cf. Proposition 36) $\implies$ 
$\Mor(e,X \otimes Y) \approx$ 
$\Mor(e,D^2 X \otimes Y) \approx$ 
$\Mor(e,D(DX) \otimes Y) \approx$ 
$\Mor(e,\hom(DX,Y)) \approx$ 
$\Mor(DX,Y)$.  
Since $DX$ is compact (cf. p. \pageref{15.46}), $\Mor(DX,-)$ is a homology theory.  
Therefore $\Mor(e,X \otimes -)$ is a homology theory.  
But \mY compact $\implies$ $X \otimes Y$ compact $\implies$ $H_X(Y) \approx \Mor(e, X \otimes Y)$.  
In other words, the arrow $H_X \ra \Mor(e,X \otimes -)$ is an isomorphism for compact \mY, hence for all \mY.]\\

\begingroup%%----------------------------------->>
\fontsize{9pt}{11pt}\selectfont
\label{15.49} %dmc mnft
\textbf{\small FACT} \ 
Suppose that \bC is a compactly generated CTC.  Fix $X \in \Ob\bC$ $-$then $X \otimes Y = 0$ iff $\forall \ Z$, 
$H_X(Y \otimes Z) = 0$.\\
\endgroup %%------------------------------------<<

\begin{proposition} \ %38
Let \bC be a compactly generated CTC.  Assume: \bC admits Adams representability.  Suppose that $F:\bC \ra \bAB$ is a homology theory $-$then $\exists$ an $X \in \Ob\bC$ and a natural isomorphism $H_X \ra F$.  
\end{proposition}

[The composite $F \circx D:\cptx\bC \ra \bAB$ is an exact functor, thus by ADR$_1$, $\exists$ an $X \in \Ob\bC$ and a natural isomorphism $h_X \ra F \circx D$.  
And: $\forall$ compact \mK, 
$H_X(K) \approx$ 
$H_K(X) \approx$ 
$\Mor(e,K \otimes X) \approx$ 
$\Mor(DX,X) \approx$
$h_X(DK) \approx$ 
$FD^2 K \approx$ 
$FK$.]

[Note: \ It follows from ADR$_2$ that 
$\Nat(H_X,H_Y) \approx$ 
$\Mor(X,Y)/\Ph(X,Y)$.  Of course $H_X \approx H_Y$ $\implies$ $X \approx Y$.]\\

\begingroup%%----------------------------------->>
\fontsize{9pt}{11pt}\selectfont
\textbf{\small EXAMPLE}  \ 
Suppose that \bC is a compactly generated CTC which admits Adams representability.  Let $\Delta:\bI \ra \bC$ be a diagram, where \bI is filtered $-$then a weak colimit \mL of $\Delta$ is a minimal weak colimit iff for every homology theory 
$F:\bC \ra \bAB$, the arrow $\colimx F\Delta_i \ra FL$ is an isomorphism.\\
\endgroup %%------------------------------------<<

Suppose that \bC is a compactly generated CTC.  Let \mT be a localization functor $-$then \mT is said to have the 
\un{ideal property}
\index{ideal property} 
(IP)
\index{IP} 
if $TX = 0$ $\implies$ $T(X \otimes Y) = 0$ $\forall \ Y$.\\

\begin{proposition} \ %39
Let \bC be a compactly generated CTC.  Suppose that \mT is a localization functor with the IP $-$then \mX \mT-acyclic and $Y$ 
\mT-local $\implies$ $\hom(X,Y) = 0$.
\end{proposition}

[$\forall \ Z$, 
$\Mor(Z,\hom(X,Y)) \approx$ 
$\Mor(Z \otimes X,Y) \approx$ 
$\Mor(X \otimes Z,Y) \approx$ 
$\Mor(T(X \otimes Z),Y) = 0$.]

%%----------------------------------------------------------------------------------------------41
[Note: \ Conversely, \mX is $T$-local if $\hom(Y,X) = 0$ for all $T$-acyclic $Y$.  In fact, 
$\Mor(Y,X) \approx$ 
$\Mor(e \otimes Y, X) \approx$ 
$\Mor(e,\hom(Y,X)) \approx$ 
$\hom(Y,X) = 0$, so Proposition 27 is applicable.  
Example: \mX $T$-local $\implies$ $\hom(Y,X)$ $T$-local $\forall \ Y$.]\\

\label{15.48} %dmc mnft
\label{17.30} %dmc mnft
\label{17.50} %dmc mnft
\label{17.74} %dmc mnft
\label{17.82} %dmc mnft
Assuming still that \mT is a localization functor with the IP, consider the exact triangle 
$e_T \ra$ 
$e \overset{\epsilon_e}{\ra}$ 
$Te \ra \Sigma_{e_T}$ (cf. Proposition 25) $-$then by CTC$_3$, $\forall \ X \in \Ob\bC$, the triangle 
$e_T \otimes X \ra$ 
\begin{tikzcd}%[sep=large]
{e \otimes X} \ar{rr}{\epsilon_e \otimes \id_X} &&{Te \otimes X}
\end{tikzcd}
$\ra \Sigma(e_T \otimes X)$ 
is exact.  But $T(e_T \otimes X) = 0$, hence 
$TX \approx T(Te \otimes X)$.  
On the other hand, $Te \otimes X$ is $T$-local if \mX is dualizable.  Proof: $Te \otimes X \approx \hom(DX,Te)$ and 
$\forall \ T$-acyclic $Y$, $\hom(Y,\hom(DX,Te))$ $\approx$ $\hom(Y\otimes DX,Te) = 0$ (cf. Proposition 39).\\

\begingroup%%----------------------------------->>
\fontsize{9pt}{11pt}\selectfont
\textbf{\small EXAMPLE}  \ 
Suppose that \bC is a compactly generated CTC.  
Let \mT be a localization functor with the IP $-$then \mT is smashing iff 
$\forall \ X$, the composite 
$T _e \otimes X \ra T(T_e \otimes X) \overset{\approx}{\lra} TX$ is an isomorphism.
\vspi
[By the above, $\sU$ is contained in the class of \mX for which the composite in question is an isomorphism.]\\
\endgroup %%------------------------------------<<

\begingroup%%----------------------------------->>
\fontsize{9pt}{11pt}\selectfont
\textbf{\small FACT} \ 
Suppose that \bC is a compactly generated CTC.  Let \mT be a localization functor with the IP $-$then there is a canonical arrow 
$TX \otimes TY \ra T(X \otimes Y)$.
\vspi
[Working with the exact triangles 
$X \otimes Y_T \ra X \otimes Y \ra $
$X \otimes TY \ra \Sigma (X \otimes Y_T)$, 
$X_T \otimes TY \ra X \otimes TY \ra$  
$TX \otimes TY \ra \Sigma (X_T \otimes TY)$, 
one finds that 
$T(\epsilon_X \otimes \epsilon_Y):T(X \otimes Y) \ra T(TX \otimes TY)$ is an isomorphism.]\\
\endgroup %%------------------------------------<<

\begingroup%%----------------------------------->>
\fontsize{9pt}{11pt}\selectfont
\textbf{\small FACT} \ 
Suppose that \bC is a compactly generated CTC.  Let \mT be a localization functor with the IP $-$then \bD is a CTC.
\vspi
[Define 
$\otimes_T:\bD \times \bD \ra \bD$ by $X \otimes_T Y = R(X \otimes Y)$.  Thus $Re$ serves as the unit and the internal hom functor 
$\hom_T:\bD^\OP \times \bD \ra \bD$ sends $(X,Y)$ to $\hom(X,Y)$ (which is automatically $T$-local).]
\vspi
[Note: \ \mX dualizable in $\bC \implies RX$ dualizable in \bD.]\\
\endgroup %%------------------------------------<<

\begingroup%%----------------------------------->>
\fontsize{9pt}{11pt}\selectfont
\textbf{\small EXAMPLE}  \ 
Suppose that \bC is a compactly generated CTC.  Let \mT be a localization functor with the IP.  
Assume: \mT is smashing $-$then \bD is a compactly generated CTC.  
In addition, \bD is a coreflective subcategory of \bC.
\vspi
[The coreflector $\bC \ra \bD$ is the assignment $ X \ra \hom(Te,X)$.]\\
\endgroup %%------------------------------------<<

Suppose that \bC is a compactly generated CTC $-$then \bC is said to be 
\un{monogenic}
\index{monogenic (compactly generated CTC)} 
if \bC is unital and 
$
\begin{cases}
\ \Mor(\Sigma^n e, X) = 0\\
\ \Mor(\Omega^n e, X) = 0
\end{cases}
$
$\forall  \ n \geq 0$ $\implies$ $X = 0$.\\

\begingroup%%----------------------------------->>
\fontsize{9pt}{11pt}\selectfont
\label{17.43}
The stable homotopy category is monogenic.\\
\endgroup %%------------------------------------<<

\begingroup%%----------------------------------->>
\fontsize{9pt}{11pt}\selectfont
\textbf{\small FACT} \ 
Suppose that \bC is a monogenic compactly generated CTC.  Let \bD be a thick subcategory of \bC $-$then $\forall$ compact \mX, $X \otimes \Ob\bD \subset \Ob\bD$.\\
\endgroup %%------------------------------------<<

%%----------------------------------------------------------------------------------------------42
\label{17.2}
Notation: When \bC is monogenic, write \mS in place of $e$ and $\Sigma^{-1}$ in place of $\Omega$, letting 
$S^{\pm n} = \Sigma^{\pm n} S$ $(\implies S^k \otimes S^l \approx S^{k+l} \ \forall \ k,l \in \Z)$ so 
$\forall \ X$, $\Sigma^{\pm 1} X \approx X \otimes S^{\pm 1}$.

[Note: \ The 
\un{$n^\text{th}$ homotopy group}
\index{$n^\text{th}$ homotopy group (\bC monogenic)} 
$\pi_n(X)$ of \mX $(n \in \Z)$ is $\Mor(S^n,X)$.]\\

\textbf{\small LEMMA} \ 
Let \bC be a monogenic compactly generated CTC $-$then a morphism $f:X \ra Y$ in \bC is an isomorphism iff $\forall \ n$, 
$\pi_n(f):\pi_n(X) \ra \pi_n(Y)$ is bijective.\\

\begingroup%%----------------------------------->>
\fontsize{9pt}{11pt}\selectfont
\textbf{\small EXAMPLE}  \ 
Let \mA be a commutative ring with unit $-$then $\bD(\bAMOD)$ is monogenic.  Here the role of \mS is played by \mA concentrated in degree 0 and $\pi_n(X) = H^{-n}(X)$.\\
\endgroup %%------------------------------------<<

\begin{proposition} \ %40
Let \bC be a monogenic compactly generated CTC.  Suppose that $F:\bC \ra \bAB$ is a homology theory $-$then $T_F$ has the IP (notation per the Bousfield-Margolis localization theorem).
\end{proposition}

[The class of $T_F$-acyclic objects coincides with $O_F$, the class of \mX such that 
$F\Sigma^n X = 0$ $\forall \ n \in \Z$ (cf. p. \pageref{15.47}.  
Therefore the claim is that for all such \mX, 
$F(\Sigma^n(X \otimes Y))$ $(= F(\Sigma^n X \otimes Y)) = 0$ $\forall \ n \in \Z$.  
To see this, note that 
$F(\Sigma^n X \otimes -):\bC \ra \bAB$ is a homology theory with the property that 
$F(\Sigma^n X \otimes S^k) =$ $F(\Sigma^{n+k} X) = 0$ $\forall \ k \in \Z$, thus, as \bC is monogenic, 
$F(\Sigma^n X \otimes -) = 0$.]\\

\begingroup%%----------------------------------->>
\fontsize{9pt}{11pt}\selectfont
\textbf{\small FACT} \ 
Suppose that \bC is a monogenic compactly generated CTC.  Let \mT be a localization functor.  Assume: \mT is smashing 
$-$then \mT has the IP.
\vspi
[Fix an \mX in $\ker T$ and consider the class of \mY : $T(X \otimes \Sigma^n Y) = 0$ $\forall \ n \in \Z$.  
This class is the object class of a triangulated subcategory of \bC containing the $S^n$ and is closed under the formation of coproducts in \bC (\mT being smashing), hence equals \bC (cf. Proposition 14).]\\
\endgroup %%------------------------------------<

Suppose that \bC is a monogenic compactly generated CTC.  
Fix an $X \in \Ob\bC$ $-$then an object $Y$ is said to be 
\un{$X$-acyclic}
\index{X-acyclic (monogenic compactly generated CTC)}
if $X \otimes Y = 0$ and an object \mZ is said to be 
\un{$X$-local}
\index{X-local (monogenic compactly generated CTC)} 
if $\hom(Y,Z) = 0$ for all $X$-acyclic \mY.  
The 
\un{Bousfield class}
\index{Bousfield class} 
$\langle X\rangle$ of \mX is the class of $X$-local objects.

Example: Let \mT be a localization functor.  
Assume: \mT is smashing $-$then $\langle TS\rangle$ is the class of 
$T$-local objects.

[Since \mT has the IP, $TS \otimes Y \approx TY$ (cf. p. \pageref{15.48}), thus \mY is $TS$-acyclic iff \mY is 
$T$-acyclic.]

[Note: \ Another point is that $\forall \ X \in \Ob\bC$,
$\langle TX\rangle = \langle TS\rangle \cap \langle X\rangle$.]\\

\textbf{\small LEMMA} \ 
$\langle X\rangle$ is a thick subcategory of \bC which is closed under the formation of products in \bC.  
And: $\forall \ Y \in \Ob\bC$ 
$\&$ $\forall \ Z \in \langle X\rangle$, $\hom(Y,Z) \in \langle X\rangle$.

%%----------------------------------------------------------------------------------------------43
[Note: \ To interpret $\langle X\rangle$, define a homology theory 
$H_X:\bC \ra \bAB$ by the rule $H_X(Y) = \pi_0(X \otimes Y)$ $-$then \mY is $X$-acyclic iff 
$H_X(Y \otimes Z) = 0$ $\forall \ Z$ (cf. p. \pageref{15.49}).  
Letting $T_X$ be the localization functor attached to $H_X$ by the Bousfield-Margolis localization theorem and taking into account Proposition 40, it follows that \mY is $X$-acyclic iff \mY is $T_X$-acyclic.  Therefore $\langle X\rangle$ is the class of 
$T_X$-local objects.]\\

Write 
$\langle X\rangle \leq \langle Y\rangle$ if $\langle X\rangle \subseteq \langle Y\rangle$ 
calling \mX, \mY 
\un{Bousfield equivalent}
\index{Bousfield equivalent} 
when 
$\langle X\rangle = \langle Y\rangle$.\\

\begin{proposition} \ 
$\langle X\rangle \leq \langle Y\rangle$ iff $Y \otimes Z = 0 \implies X \otimes Z = 0$.
\end{proposition}

[Note: \ Consequently $\langle S\rangle$ is the largest Bousfield class and $\langle 0\rangle$ is the smallest.]\\

Notation: $\langle X\rangle \amalg \langle Y\rangle = \langle X \amalg Y\rangle$ 
and 
$\langle X\rangle \otimes \langle Y\rangle = \langle X \otimes Y\rangle$.

[Note: \ Both operations are welldefined.  
Examples: 
(1) $\langle X\rangle \amalg \langle 0\rangle = \langle X\rangle$, 
$\langle X\rangle \amalg \langle S\rangle = \langle S\rangle$; 
(2) $\langle X\rangle \otimes \langle 0\rangle = \langle 0\rangle$, 
$\langle X\rangle \otimes \langle S\rangle = \langle X\rangle$.]\\

\begingroup%%----------------------------------->>
\fontsize{9pt}{11pt}\selectfont
\textbf{\small FACT} \ 
If $X \ra Y \ra Z \ra \Sigma X$ is an exact triangle, then 
$\langle Y\rangle \leq \langle X\rangle \amalg \langle Z\rangle$.\\
\endgroup %%------------------------------------<<

Maintaining the assumption that \bC is monogenic, let 
$\langle\bC\rangle$ 
be the conglomerate whose elements are the Bousfield classes.  
Denote by $\bDL\langle\bC\rangle$ the subconglomerate of 
$\langle\bC\rangle$ 
consisting of those 
$\langle X\rangle$ 
with 
$\langle X\rangle \otimes \langle X\rangle = \langle X\rangle$ 
and denote by $\bBA\langle\bC\rangle$  the subconglomerate of 
$\langle\bC\rangle$ 
consisiting of those 
$\langle X\rangle$ 
that admit a complement, i.e., for which $\exists$ 
$\langle Y\rangle$ : 
$\langle X\rangle \otimes \langle Y\rangle = \langle0\rangle$ 
and 
$\langle X\rangle \amalg \langle Y\rangle = \langle S\rangle$.

[Note: \ $\bDL\langle\bC\rangle$ is a ``distributive lattice'' and $\bBA\langle\bC\rangle$ is a ``boolean algebra''.]\\

\begingroup%%----------------------------------->>
\fontsize{9pt}{11pt}\selectfont
Complements, if they exist, are unique.  Thus suppose that $\langle X \rangle$ admits two complements 
$\langle Y^\prime\rangle$ and $\langle Y\pp\rangle$ $-$then 
$\langle Y^\prime\rangle = $
$\langle Y^\prime\rangle \otimes \langle S\rangle =$ 
$\langle Y^\prime\rangle \otimes (\langle X\rangle \amalg \langle Y\pp\rangle) =$ 
$(\langle Y^\prime\rangle \otimes \langle X\rangle) \amalg (\langle Y^\prime\rangle \otimes \langle Y\pp\rangle) =$ 
$\langle 0\rangle \amalg (\langle Y^\prime\rangle) \otimes \langle Y\pp\rangle) =$
$\langle Y^\prime\rangle \otimes \langle Y\pp\rangle = \langle Y\pp\rangle$ (by symmetry).\\
\endgroup %%------------------------------------<<

Notation: Given $\langle X\rangle \in \bBA\langle \bC\rangle$, 
let 
$\langle X\rangle^c$ be its complement.\\

\textbf{\small LEMMA} \  
$\bBA\langle \bC\rangle$ is contained in $\bDL\langle \bC\rangle$.

[$\langle X\rangle = \langle X\rangle \otimes (\langle X\rangle \amalg \langle X\rangle^c) =$
$(\langle X\rangle \otimes \langle X\rangle) \amalg (\langle X\rangle \otimes \langle X\rangle^c) =$
$\langle X\rangle \otimes \langle X\rangle$.]\\

\begingroup%%----------------------------------->>
\fontsize{9pt}{11pt}\selectfont
Examples in the stable homotopy category show that the inclusions 
$\bBA\langle\bC\rangle \subset \bDL\langle\bC\rangle \subset \langle\bC\rangle$ 
are strict (Bousfield\footnote[2]{\textit{Comment. Math. Helv.} \textbf{54} (1979), 368-377.}).\\
\endgroup %%------------------------------------<

%%----------------------------------------------------------------------------------------------44
\begingroup%%----------------------------------->>
\fontsize{9pt}{11pt}\selectfont
\textbf{\small EXAMPLE}  \ 
Let \mT be a localization functor $-$then there is an exact triangle 
$S_T \ra $
$S \overset{\epsilon_S}{\ra}$
$TS \ra \Sigma S_T$, 
where $S_T$ is $T$-acyclic (cf. Proposition 25), hence 
$\langle S \rangle = \langle S_T\rangle \amalg \langle TS\rangle$.  
If further \mT is smashing, then 
$\langle S_T\rangle \otimes \langle TS\rangle =$ 
$\langle S_T \otimes  TS\rangle =$ 
$\langle TS_T\rangle = \langle 0\rangle$ 
$\implies$ 
$\langle S_T\rangle^c = \langle TS\rangle$.
\vspi
[Note: \ Take for \bC the stable homotopy category $-$then \mX compact $\implies$ 
$\langle X\rangle \in \bBA\langle\bC\rangle$
and 
$T_Y(\langle Y \rangle = \langle X \rangle^c)$ is smashing (Bousfield (ibid.)).]\\
\endgroup %%------------------------------------<

\begingroup%%----------------------------------->>
\fontsize{9pt}{11pt}\selectfont
\textbf{\small EXAMPLE}  \ 
If \mX is dualizable, then $\langle X\rangle = \langle DX\rangle$.  Indeed, \mX is a retract of 
$X \otimes DX \otimes X$ (cf. p. \pageref{15.50}), thus 
$\langle X\rangle \leq$ 
$\langle X \otimes DX \otimes X\rangle \leq \langle DX\rangle$.  
But $DX$ is dualizable, so 
$\langle DX\rangle \leq$ 
$\langle D^2X\rangle = \langle X\rangle$ (cf. Proposition 32).\\
\endgroup %%------------------------------------<<

\label{17.25}
Suppose that \bC is a monogenic compactly generated CTC $-$then a 
\un{ring object}
\index{ring object (in a monogenic compactly generated CTC)} 
in \bC is an object \mR equipped with a product 
$R \otimes R \ra R$ 
and a unit 
$S \ra R$ 
such that 
\begin{tikzcd}%[sep=large]
{R \otimes R \otimes R} \ar{d}\ar{r} &{R \otimes R} \ar{d}\\
{R \otimes R} \ar{r} &{R}
\end{tikzcd}
and 
\begin{tikzcd}%[sep=large]
{S \otimes R} \ar{rd}\ar{r} &{R \otimes R} \ar{d} &{R \otimes S}\ar{l} \ar{ld}\\
&{R}
\end{tikzcd}
commute.  A ring object \mR is 
\un{commutative}
\index{commutative ring object (in a monogenic compactly generated CTC)}
if 
\begin{tikzcd}[sep=small]
{R \otimes R} \ar{rdd}\ar{rr}{\Tee} &&{R \otimes R} \ar{ldd}\\
\\
&{R}
\end{tikzcd}
commutes.

Example: $\forall \ X \in \Ob\bC$, $\hom(X,X)$ is a ring object, hence $DX \otimes X$ is a ring object if \mX is dualizable.\\

\begingroup%%----------------------------------->>
\fontsize{9pt}{11pt}\selectfont
\textbf{\small EXAMPLE}  \ 
If \mR is a ring object, then 
$\langle R\rangle \otimes \langle R\rangle = \langle R\rangle$ 
(\mR is a retract of $R \otimes R$).\\
\endgroup %%------------------------------------<<

\textbf{\small LEMMA} \ 
If \mR is a ring object, then $\pi_*(R)$ is a graded ring with unit which is graded commutative provided that \mR is commutative.\\

Given a ring object \mR, a (left) 
\un{\mR-module}
\index{\mR-module (in a monogenic compactly generated CTC)}
is an object \mM equipped with an arrow 
$R \otimes M \ra M$ such that 
\begin{tikzcd}%[sep=large]
{R \otimes R \otimes M} \ar{d}\ar{r} &{R \otimes M} \ar{d}\\
{R \otimes M} \ar{r} &{M}
\end{tikzcd}
and 
\begin{tikzcd}%[sep=large]
{S \otimes M} \ar{rd}\ar{r} &{R \otimes M} \ar{d}\\
&{M}
\end{tikzcd}
commute.  
Example: $\forall \ X \in \Ob\bC$, $R \otimes X$ and $\hom(X,R)$ are $R$-modules.  $R$-MOD
\index{\bRMOD}
is the category whose objects are the $R$-modules.

[Note: \ If $f:M \ra N$ is a morphism of $R$-modules and if 
$M \overset{f}{\ra} N \ra C_f \ra \Sigma M$ is exact, then $C_f$ need not admit an $R$-module structure.]\\

\begingroup%%----------------------------------->>
\fontsize{9pt}{11pt}\selectfont
\textbf{\small EXAMPLE}  \ 
If \mR is a ring object and if \mM is an $R$-module, then 
$\langle M\rangle \leq \langle R\rangle$
(\mM is a retract of $R \otimes M$).
\vspi
[Note: \ \mM is necessarily $T_R$-local.]\\
\endgroup %%------------------------------------<

%%----------------------------------------------------------------------------------------------45
\label{17.32} %dmc mnft
\label{17.44} %dmc mnft

\begingroup%%----------------------------------->>
\fontsize{9pt}{11pt}\selectfont
\textbf{\small EXAMPLE}  \ 
Let \mT be a localization functor with the IP $-$then $TS$ is a commutative ring object 
(via $TS \otimes TS \ra T(S \otimes S) = TS$ and $\epsilon_S:S \ra TS$).  
Moreover, every $T$-local object \mX is a $TS$-module 
(via $TS \otimes X = TS \otimes TX \ra T(S \otimes X) = TX = X$).\\
\endgroup %%------------------------------------<

\begingroup%%----------------------------------->>
\fontsize{9pt}{11pt}\selectfont
\textbf{\small EXAMPLE}  \ 
If \mR is a ring object with the property that the product $R \otimes R \ra R$ is an isomorphism, then $T_R$ is smashing.  
Proof: $\forall \ X \in \Ob\bC$, $R \otimes X$ is $T_R$-local and here $T_R X = R \otimes X$ 
(since $R \otimes R \approx R$), thus $T_R$ preserves coproducts.\\
\endgroup %%------------------------------------<

Definitions: 
(1) \ An $R$-module \mM is 
\un{free}
\index{free ($R$-module in a monogenic compactly generated CTC) } 
if it is isomorphic to a coproduct 
$\coprod\limits_i \Sigma^{n_i} R$; 
(2) \ A nonzero ring object \mR is a 
\un{skew field object}
\index{skew field object} 
if every \mM in $R$-MOD is free; 
(3) \ A skew field object \mR is a 
\un{field object}
\index{field object} if \mR is commutative.\\

\begin{proposition} \ %42
Let \bC be a monogenic compactly generated CTC.  Suppose that \mR is a nonzero ring object in \bC.  
Assume: The homogeneous elements of $\pi_*(R)$ are invertible $-$then \mR is a skew field object.
\end{proposition}

[Fix an \mM in \RMOD.  Owing to our assumption, 
$\pi_*(M) = \bigoplus\limits_i \Sigma^{n_i} \pi_*(R)$, where 
$(\Sigma^{n_i} \pi_*(R))_n = \Mor(S^{n - n_i},R) =$ 
$\Mor(S^n,\Sigma^{n_i}R) = \pi_n(\Sigma^{n_i} R)$.  
Thus there is a morphism 
$\coprod\limits_i \Sigma^{n_i} R \ra M$ of \mR-modules inducing an isomorphism 
$\bigoplus\limits_i \pi_{* - n_i}(R) \ra \pi_*(M)$ in homotopy, hence 
$\coprod\limits_i \Sigma^{n_i} R \approx M$.]\\

\begingroup%%----------------------------------->>
\fontsize{9pt}{11pt}\selectfont
In the stable homotopy category, the $n^{th}$ 
Morava K-theory spectrum $\bK(n)$ at the prime $p$ is a skew field object.\\
\endgroup %%------------------------------------<<

\begingroup%%----------------------------------->>
\fontsize{9pt}{11pt}\selectfont
\textbf{\small EXAMPLE}  \ 
Let \mR be a skew field object.  Assume: $\langle R\rangle \in \bBA\langle\bC\rangle$ $-$then 
$\langle R\rangle$ is minimal among all nontrivial Bousfield classes.
\vspi
[Note: \ In the stable homotopy category, the Eilenberg-MacLane spectrum $\bH(\F_p)$ is a field object but 
$\langle\bH(\F_p)\rangle$ is not minimal.]\\
\endgroup %%------------------------------------<<

Suppose that \bC is a monogenic compactly generated CTC.  Given $X \in \Ob\bC$ and 
$f \in \Mor(\Sigma^n X,X)$, let $X/f$ be a completion of 
$\Sigma^n X \overset{f}{\ra} X$ to an exact triangle (cf. TR$_3$) and write $f^{-1}X$ for $\telsub(\bX,\bff)$, where 
$(\bX,\bff)$ is the object in $\bFIL(\bC)$ defined by 
$X \ra \Sigma^{-n} X \ra \Sigma^{-2n} X \ra \cdots$.\\

\textbf{\small LEMMA} \ 
If $f:\Sigma^n X \ra X$ is an isomorphism, then $X \approx f^{-1} X$.\\

\begin{proposition} \ %43
For every $f:\Sigma^n X \ra X$, $\langle X\rangle = \langle X/f\rangle \amalg \langle f^{-1} X\rangle$.
\end{proposition}

%%----------------------------------------------------------------------------------------------46
[To prove that 
$\langle X\rangle \leq \langle X/f\rangle \amalg \langle f^{-1} X\rangle$, 
one must show that $X/f \otimes Z = 0$ $\&$ $f^{-1} X \otimes Z = 0$ $\implies$ 
$X \otimes Z = 0$.  
But
$\Sigma^n X \overset{f}{\ra} X \ra X/f \ra \Sigma(\Sigma^n X)$ exact $\implies$ 
$\Sigma^n X \otimes Z \ra$
$X \otimes Z \ra$ 
$X/f \otimes Z \ra \Sigma(\Sigma^n X \otimes Z)$ exact (cf. CTC$_3$) $\implies$ 
$\Sigma^n X \otimes Z \approx$ 
$X \otimes Z$ (cf. p. \pageref{15.51}) $\implies$ 
$X \otimes Z \approx$ 
$(f \otimes \id_Z)^{-1}(X \otimes Z)$ (by the lemma).  
And: 
$(f \otimes \id_Z)^{-1}(X \otimes Z) = f^{-1}X \otimes Z = 0$.]\\

\label{17.22}
\label{17.23}
\label{17.46} %dmc mnft
\label{17.47} %dmc mnft
\label{17.52}
\label{17.58} %dmc mnft
\label{17.62} %dmc mnft
\label{17.65} %dmc mnft
\label{17.69}
\begingroup%%----------------------------------->>
\fontsize{9pt}{11pt}\selectfont
\textbf{\small FACT} \ 
Suppose that \mX is compact $-$then $f^{-1}X = 0$ iff $\exists$ $k$ such that the composite 
$f \circx \Sigma^n f \circx \cdots \circx \Sigma^{(k-1)n} f$: $\Sigma^{kn}X \overset{f^k}{\lra}X$ vanishes.\\
\endgroup %%------------------------------------<<

\begingroup%%----------------------------------->>
\fontsize{9pt}{11pt}\selectfont
\textbf{\small FACT} \ 
Let \mR be a ring object.  Fix $\alpha \in \pi_n(R)$ and let $\ov{\alpha}$ be the map 
\begin{tikzcd}%[sep=large]
{S^n \otimes R} \ar{r}{\alpha \otimes \id_R}
&{R \otimes R} \ar{r}
&{R}
\end{tikzcd}
$-$then $\alpha$ is nilpotent in $\pi_*(R)$ iff $\ov{\alpha}^{-1}R = 0$.\\
\endgroup %%------------------------------------<<

\begingroup%%----------------------------------->>
\fontsize{9pt}{11pt}\selectfont
\textbf{\small FACT} \ 
Given $f:S \ra X$, write $X_f^{(\infty)}$ for tel(\bX,\bff), where (\bX,\bff) is the object in $\bFIL(\bC)$ defined by 
$S \overset{f}{\lra} $
\begin{tikzcd}%[sep=large]
{X} \ar{r}{f \otimes \id}
&{X \otimes X} \ar{r}{f \otimes \id}
&{\cdots}
\end{tikzcd}
, and let $f^{(\infty)}$ be the arrow 
$S \ra X_f^{(\infty)}$ $-$then $X_f^{(\infty)} = 0$ iff $f^{(\infty)} = 0$.\\
\endgroup %%------------------------------------<<

Let \bC be a triangulated category; let $\bC^{\leq 0},\bC^{\geq 0}$ be full, isomorphism closed subcategories of \bC containing 0 and denote by $\bC^{\leq -1},\bC^{\geq 1}$ the isomorphism closure of 
$\Sigma \bC^{\leq 0}$, 
$\Omega \bC^{\geq 0}$ 
$-$then the pair $(\bC^{\leq 0},\bC^{\geq 0})$ is said to be a 
\un{t-structure}
\index{t-structure} on \bC if the following conditions are satisfied.

\indent\indent(t-st$_1$) \quad $\bC^{\leq -1}$ is a subcategory of $\bC^{\leq 0}$ and 
$\bC^{\geq 1}$ is a subcategory of $\bC^{\geq 0}$.\\
\indent\indent(t-st$_2$) \quad 
$\forall \ X \in \Ob\bC^{\leq 0}$,
$\forall \ Y \in \Ob\bC^{\geq 1}$, $\Mor(X,Y) = 0$.\\
\indent\indent(t-st$_3$) \quad $\forall \ X \in \Ob\bC$, $\exists$ an exact triangle 
$X_0 \ra X \ra X_1 \ra \Sigma X_0$ with $X_0 \in \Ob\bC^{\leq 0}$, $X_1 \in \Ob\bC^{\geq 1}$.

[Note: $\bH(\bC) = \bC^{\leq 0} \cap \bC^{\geq 0}$ is called the 
\un{heart}
\index{heart (of a t-structure)} 
of the t-structure.]

Remark: If $(\bC^{\leq 0},\bC^{\geq 0})$ is a t-structure on \bC, then 
$((\bC^{\geq 0})^\OP,(\bC^{\leq 0})^\OP)$ is a t-structure on $\bC^\OP$.\\

\begingroup%%----------------------------------->>
\fontsize{9pt}{11pt}\selectfont
\textbf{\small EXAMPLE}  \ \ 
Let \bA be an abelian category.  \ 
Given an \mX in $\bC\bX\bA$, $\forall \ n \in \Z$, define the $n^\text{th}$ truncated cochain complexes 
$\tau^{\leq n}X$ $\&$ 
$\tau^{\geq n}X$ of \mX by 
$\cdots \ra X^{n-2} \ra$ 
$X^{n-1} \ra$ 
$\ker d_X^n \ra$ 
$0 \ra \cdots$ 
$\&$ 
$\cdots \ra$ 
$0 \ra $ 
$\coker d_X^{n-1} \ra$ 
$X^{n+1} \ra$ 
$X^{n+2} \ra \cdots$.  
So, the cohomology of $\tau^{\leq n}X$ is trivial in degree $> n$ and 
the cohomology of $\tau^{\geq n}X$ is trivial in degree $< n$ and there is an arrow 
$\tau^{\leq n}X \ra X$ which induces an isomorphism in cohomology in degree $\leq n$ 
and there is an arrow 
$X \ra \tau^{\geq n} X$ which induces an isomorphism in cohomology in degree $\geq n$.    \ 
The functors \ 
$
\begin{cases}
\ \tau^{\leq n}:\bC\bX\bA \ra \bC\bX\bA\\
\ \tau^{\geq n}:\bC\bX\bA \ra \bC\bX\bA
\end{cases}
$
pass through $\bK(\bA)$ to the derived category 
$
\bD(\bA):
\begin{cases}
\ \tau^{\leq n}:\bD(\bA) \ra \bD(\bA)\\
\ \tau^{\geq n}:\bD(\bA) \ra \bD(\bA)
\end{cases}
$
and $\forall \ X$, $\exists$ an exact triangle 
$\tau^{\leq n}X \ra$ 
$X \ra$
$\tau^{\geq {n+1}}X \ra$ 
$\Sigma \tau^{\leq n}X$.
Write 
$\bD^{\leq 0}(\bA)$ for the full subcategory of $\bD(\bA)$ consisting of those \mX such that 
$H^q(X) = 0$ $(q > 0)$ and write
$\bD^{\geq 0}(\bA)$ for the full subcategory of $\bD(\bA)$ consisting of those \mX such that 
%%----------------------------------------------------------------------------------------------47
$H^q(X) = 0$ $(q < 0)$ $-$then the pair $(\bD^{\leq 0}(\bA),\bD^{\geq 0}(\bA))$ is a t-structure on $\bD(\bA)$ and its heart is equivalent to \bA.\\
\endgroup %%------------------------------------<<

Given a t-structure $(\bC^{\leq 0},\bC^{\geq 0})$  on \bC, let 
$
\begin{cases}
\ \bC^{\leq n}\\
\ \bC^{\geq n}
\end{cases}
$
be the isomorphism closure of 
$
\begin{cases}
\ \Omega^n\bC^{\leq 0}\\
\ \Omega^n \bC^{\geq 0}
\end{cases}
(n > 0)
$
and let
$
\begin{cases}
\ \bC^{\leq n}\\
\ \bC^{\geq n}
\end{cases}
$
be the isomorphism closure of 
$
\begin{cases}
\ \Sigma^{\abs{n}}\bC^{\leq 0}\\
\ \Sigma^{\abs{n}} \bC^{\geq 0}
\end{cases}
(n < 0)
$
$-$then $\forall \ n \in \Z$, the pair $(\bC^{\leq n},\bC^{\geq n})$ is a t-structure on \bC.\\

\begin{proposition} \ %44
Suppose that $(\bC^{\leq 0},\bC^{\geq 0})$ is a t-structure on \bC $-$then $\forall \ n \in \Z$, $\bC^{\leq n}$ is a coreflective subcategory of \bC with coreflector $\tau^{\leq n} X \ra X$ and $\bC^{\geq n}$ is a reflective subcategory of \bC with reflector 
$X \ra \tau^{\geq n} X$.
\end{proposition}

[It suffices to construct $\tau^{\leq 0}$.  Thus for any $X \in \Ob\bC$, $\exists$ an exact triangle 
$X_0 \ra$ 
$X \ra$
$X_1 \ra$
$\Sigma X_0$, 
where $X_0 \in \Ob\bC^{\leq 0}$ $\&$ $X_1 \in \Ob\bC^{\geq 1}$ (cf. t-st$_3$), so 
$\forall \ Y \in \Ob\bC^{\leq 0}$, there is an exact sequence 
$\Mor(Y,\Omega X_1) \ra$ 
$\Mor(Y,X_0) \ra$
$\Mor(Y,X) \ra$
$\Mor(Y,X_1)$.  Here $\Mor(Y,X_1) = 0$ (cf. t-st$_2$).  
In addition, 
$\Mor(Y,\Omega X_1) \approx \Mor(\Sigma Y,X_1)$ and 
$\Sigma \bC^{\leq 0} \subset \bC^{\leq -1}$ $\subset$ $\bC^{\leq 0}$ (cf. t-st$_1$) $\implies$ 
$\Mor(\Sigma Y,X_1) = 0$ (cf. t-st$_2$).  Therefore, $\forall \ Y \in \Ob\bC^{\leq 0}$, 
$\Mor(Y,X_0) \approx \Mor(Y,X)$ and we can let $\tau^{\leq 0}X = X_0$.]

[Note: \ Similar reasoning gives $\tau^{\geq 1}X = X_1$.]\\

The functors $\tau^{\leq n}$, $\tau^{\geq n}$ figuring in Proposition 44 are called the 
\un{truncation functors}
\index{truncation functors} 
of the t-structure.

[Note: \ $\forall \ X$, $\exists$ an exact triangle 
$\tau^{\leq n}X \ra$
$X \ra$ 
$\tau^{\geq n+1}X \ra$
$\Sigma \tau^{\leq n}X$
and since 
$\Mor(\Sigma \tau^{\leq n}X,\tau^{\geq n+1}X) = 0$, the arrow 
$\tau^{\geq n+1}X \ra \Sigma \tau^{\leq n}X$ is unique (cf. p. \pageref{15.52}).]\\

\begingroup%%----------------------------------->>
\fontsize{9pt}{11pt}\selectfont
\textbf{\small EXAMPLE}  \ 
Let \bA be an abelian category.  
Working with the t-structure on $\bD(\bA)$ spelled out above, 
$\tD^{\leq n}(\bA)$ is the coreflective subcategory of $\bD(\bA)$ consisting of those \mX such that 
$H^q(X) = 0$ $(q > n)$ and 
$\tD^{\geq n}(\bA)$ is the reflective subcategory of $\bD(\bA)$  consisting of those \mX such that 
$H^q(X) = 0$ $(q < n)$.\\
\endgroup %%------------------------------------<<

\begingroup%%----------------------------------->>
\fontsize{9pt}{11pt}\selectfont
Observations: Let $m,n \in \Z$ $-$then 
(1) $m \leq n$ $\implies$ 
$\tau^{\geq n} \circx$ $\tau^{\geq m}$ 
$\approx$ 
$\tau^{\geq m} \circx$ $\tau^{\geq n}$ 
$\approx$ 
$\tau^{\geq n}$ 
and 
$\tau^{\leq n} \circx \tau^{\leq m}$ 
$\approx$ 
$\tau^{\leq m} \circx \tau^{\leq n}$ 
$\approx$ 
$\tau^{\leq m}$; 
(2) $m > n$ $\implies$ 
$\tau^{\leq n} \circx \tau^{\geq m} = 0$ and 
$\tau^{\geq m} \circx \tau^{\leq n} = 0$.\\
\endgroup %%------------------------------------<<

\begingroup%%----------------------------------->>
\fontsize{9pt}{11pt}\selectfont
\textbf{\small FACT} \ 
If $m \leq n$, then $\forall \ X \in \Ob\bC$, $\exists$ a unique arrow 
$\tau^{\geq m} \tau^{\leq n} X \ra \tau^{\leq n} \tau^{\geq m}X$ such that the diagram 
\begin{tikzcd}%[sep=large]
{\tau^{\leq n}X} \ar{d}  \ar{r}
&{X} \ar{r}
&{\tau^{\geq m}X}\\
{\tau^{\geq m}\tau^{\leq n}X} \ar{rr} &&{\tau^{\leq n}\tau^{\geq m}X} \ar{u}
\end{tikzcd}
commutes.
\vspi
%%----------------------------------------------------------------------------------------------48
[Note: \ The arrow
$\tau^{\geq m} \tau^{\leq n} X \ra \tau^{\leq n} \tau^{\geq m}X$
is an isomorphism provided that \bC satisfies the octahedral axiom.  To see this, consider the exact triangles 
$\tau^{\leq m-1}X \ra$ 
$\tau^{\leq n}X \ra$ 
$\tau^{\geq m}\tau^{\leq n}X \ra$ 
$\Sigma\tau^{\leq m-1}X$, 
$\tau^{\leq n}X \ra$ 
$X \ra$ 
$\tau^{\geq n+1}X \ra$ 
$\Sigma\tau^{\leq n}X$, 
$\tau^{\leq m-1}X \ra$ 
$X \ra$ 
$\tau^{\geq m}X \ra$ 
$\Sigma\tau^{\leq m-1}X$.]\\
\endgroup %%------------------------------------<<

Notation: Write 
$
\begin{cases}
\ \bC^{<n}\\
\ \bC^{>n}
\end{cases}
$
in place of 
$
\begin{cases}
\ \bC^{\leq n-1}\\
\ \bC^{\geq n+1}
\end{cases}
$
and 
$
\begin{cases}
\ \tau^{<n}\\
\ \tau^{>n}
\end{cases}
$
in place of 
$
\begin{cases}
\ \tau^{\leq n-1}\\
\ \tau^{\geq n+1}
\end{cases}
.
$
\\
\vspace{0.2cm}

\textbf{\small LEMMA} \ 
Let $X \in \Ob\bC$ $-$then 
$
X \in 
\begin{cases}
\ \Ob\bC^{\leq n}\\
\ \Ob\bC^{\geq n}
\end{cases}
$
iff 
$
\begin{cases}
\ \tau^{>n} X = 0\\
\ \tau^{<n} X = 0
\end{cases}
$
.\\
\vspace{0.25cm}

\begin{proposition} \ %45
Suppose that $(\bC^{\leq 0},\bC^{\geq 0})$ is a t-structure on \bC.  Let 
$X^\prime \ra X \ra X\pp \ra \Sigma X^\prime$ be an exact triangle $-$then 
$
\begin{cases}
\ X^\prime\\
\ X\pp
\end{cases}
\in \Ob\bC^{\leq 0}
$
$\implies$ $X \in \Ob\bC^{\leq 0}$ $\&$ 
$
\begin{cases}
\ X^\prime\\
\ X\pp
\end{cases}
$
$\in \Ob\bC^{\geq 0}$ $\implies$ $X \in \Ob\bC^{\geq 0}$.\\
\end{proposition}

\begingroup%%----------------------------------->>
\fontsize{9pt}{11pt}\selectfont
Let \bA be an additive category.  Given a class $O \subset \Ob\bA$, the 
$
\begin{cases}
\ \text{\un{left annihilator}}\  \Ann_\text{L} O\\
\ \text{\un{right annihilator}} \ \Ann_\text{R} O
\end{cases}
$
\index{left annihilator}\index{right annihilator}
of $O$ is 
$
\begin{cases}
\ \{Y:\Mor(Y,X) = 0 \ \forall \ X \in O\}\\
\ \{Y:\Mor(X,Y) = 0 \ \forall \ X \in O\}
\end{cases}
$
.\\
\vspace{0.2cm}
\endgroup %%------------------------------------<<

\label{15.34}
\begingroup%%----------------------------------->>
\fontsize{9pt}{11pt}\selectfont
\textbf{\small EXAMPLE}  \ 
Let \bA be an additive category.  Suppose that $\sT$, $\sF$ are subclasses of $\Ob\bA$ $-$then the pair $(\sT,\sF)$ is said to be a \un{torsion theory}
\index{torsion theory} 
on \bA if 
$\Ann_L\sF = \sT$ and 
$\Ann_R\sT = \sF$.  
Example: $\forall$ t-structure $(\bC^{\leq 0},\bC^{\geq 0})$ on \bC, 
$\Ann_\text{L} \bC^{\geq 1} = \bC^{\leq 0}$ and 
$\Ann_\text{R} \bC^{\leq 0} = \bC^{\geq 1}$, i.e., $(\bC^{\leq 0},\bC^{\geq 1})$ is a torsion theory on \bC.\\
\endgroup %%------------------------------------<<

\textbf{\small LEMMA} \ 
Let \bC be a triangulated category satisfying the octahedral axiom.  Suppose that $(\bC^{\leq 0},\bC^{\geq 0})$ is a t-structure on \bC $-$then $\forall \ X \in \Ob\bC$, 
$\tau^{\geq 0} \tau^{\leq 0} X \approx \tau^{\leq 0} \tau^{\geq 0} X$.\\

\index{Theorem of the Heart}
\index{Theorem: Theorem of the Heart}
\textbf{\small THEOREM OF THE HEART} \quad
Let \bC be a triangulated category with finite coproducts satisfying the octahedral axiom.  
Suppose that $(\bC^{\leq 0},\bC^{\geq 0})$ is a t-structure on \bC $-$then its heart $\bH(\bC)$ is an abelian category.

[$\bH(\bC)$ is closed under the formation of finite coproducts in \bC (use the exact triangle 
$X \ra X \amalg Y \ra Y \overset{0}{\ra} \Sigma X$ and quote Proposition 45).  
To prove that $\bH(\bC)$ has kernels and cokernels and that parallel morphisms are isomorphisms, take an arrow 
$f:X \ra Y$ in $\bH(\bC)$ and place it in an exact triangle 
$X \overset{f}{\ra} Y \ra Z \ra \Sigma X$ ($\implies$ $Z \in \Ob\bC^{\leq 0} \cap \Ob\bC^{\geq -1}$ (cf. Proposition 45)).  
For any $W \in \Ob\bH(\bC)$, there are exact sequences 
$\Mor(W,\Omega Y) \ra$ 
$\Mor(W,\Omega Z) \ra$ 
$\Mor(W,X) \ra$ 
$\Mor(W,Y)$, 
$\Mor(\Sigma X,W) \ra$ 
$\Mor(Z,W) \ra$ 
$\Mor(Y,W) \ra$ 
$\Mor(X,W)$.  
Since 
$\Mor(W,\Omega Y) = 0$, 
$\Mor(\Sigma X,W) = 0$ 
and 
$\Mor(W,\Omega Z) \approx \Mor(W,\tau^{\leq 0}\Omega Z)$, 
$\Mor(Z,W) \approx \Mor(\tau^{\geq 0} Z,W)$, 
it follows that 
$\ker f \approx \tau^{\leq 0}\Omega Z$, 
$\coker f \approx \tau^{\geq 0}Z$.
In this connection, note that 
$Z \in \Ob\bC^{\leq 0}$ 
$\implies$ 
$\tau^{\geq 0}Z \approx$ 
$\tau^{\geq 0}\tau^{\leq 0}Z \approx$ 
$\tau^{\leq 0}\tau^{\geq 0}Z$ $\implies$ 
$\coker f \in \Ob\bH(\bC)$
and 
$Z \in \Ob\bC^{\geq -1}$ $\implies$ 
$\Omega Z \in \Ob\bC^{\geq 0}$ $\implies$
%%----------------------------------------------------------------------------------------------49
$\tau^{\leq 0}\Omega Z \approx$ 
$\tau^{\leq 0}\tau^{\geq 0}\Omega Z \approx$ 
$\tau^{\geq 0}\tau^{\leq 0}\Omega Z$ $\implies$ $\ker f \in \Ob\bH(\bC)$.  
Now fix an exact triangle 
$I \ra Y \ra \tau^{\geq 0}Z \ra \Sigma I$ ($\implies I \in \Ob\bC^{\geq 0}$ (cf. Proposition 45)).  
Applying the octahedral axiom to 
$Y \ra Z \ra \Sigma X \ra \Sigma Y$, 
$Z \ra \tau^{\geq 0}Z \ra \Sigma \tau^{< 0}Z \ra \Sigma Z$, 
$Y \ra \tau^{\geq 0}Z \ra \Sigma I \ra \Sigma Y$, 
one gets an exact triangle 
$\Sigma X \ra \Sigma I \ra \Sigma \tau^{< 0}Z \ra \Sigma^2 X$, 
which leads to an exact triangle 
$\tau^{\leq 0}\Omega Z \ra X \ra I \ra \Sigma \tau^{\leq 0}\Omega Z$, thus $I \in \Ob\bC^{\leq 0}$ (cf. Proposition 45) 
and so $I \in \Ob\bH(\bC)$.  
Finally, 
$I \approx \coim f$ (consider $\ker f \ra X \ra I \ra \Sigma \ker f$) and 
$I \approx \im f$ (consider $I \ra Y \ra \coker f \ra \Sigma I$).  
Therefore $\bH(\bC)$ is abelian.]

[Note: \ In general, there is no a priori connection between \bC and the derived category of $\bH(\bC)$.]\\

\begingroup%%----------------------------------->>
\fontsize{9pt}{11pt}\selectfont
\label{17.8}
\textbf{\small EXAMPLE}  \ 
Take for \bC the stable homotopy category and let 
$
\begin{cases}
\ \bC^{\geq 0} = \{\bX: \pi_q(\bX) = 0 \ (q > 0)\}\\
\ \bC^{\leq 0} = \{\bX: \pi_q(\bX) = 0 \ (q < 0)\}
\end{cases}
$
$-$then $(\bC^{\leq 0},\bC^{\geq 0})$ is a t-structure on \bC.  Its heart is equivalent to \bAB 
(cf. p. \pageref{15.53}).
\vspi
[Note: $\tau^{\leq 0}\bX$ is called the
\un{connective cover}
\index{connective cover} of \bX 
(the arrow  $\tau^{\leq 0}\bX \ra \bX$ induces an isomorphism 
$\pi_n(\tau^{\leq 0}\bX) \ra \pi_n(\bX)$ for $n \geq 0$).]\\
\endgroup %%------------------------------------<<

\begingroup%%----------------------------------->>
\fontsize{9pt}{11pt}\selectfont
Let \bC be a triangulated category with finite coproducts satisfying the octahedral axiom.  Suppose that 
$(\bC^{\leq 0},\bC^{\geq 0})$ is a t-structure on \bC $-$then $H^0:\bC \ra \bH(\bC)$ is the functor that sends \mX to 
$\tau^{\geq 0}\tau^{\leq 0}X$ $\approx$ $\tau^{\leq 0}\tau^{\geq 0}X$.\\
\endgroup %%------------------------------------<<

\begingroup%%----------------------------------->>
\fontsize{9pt}{11pt}\selectfont
\textbf{\small FACT} \ 
$H^0$ is an exact functor.
\vspi
[Fix an exact triangle $X \ra Y \ra Z \ra \Sigma X$ and proceed in stages.
\\
\indent\indent (I) \ Assume that $X,Y,Z \in \Ob\bC^{\geq 0}$ $-$then 
$0 \ra H^0(X) \ra H^0(Y) \ra H^0(Z)$ is exact.
\\
\indent\indent (II$^{\geq 0}$) \ 
Assume that $Z \in \Ob\bC^{\geq 0}$ $-$then 
$0 \ra H^0(X) \ra H^0(Y) \ra H^0(Z)$ is exact.
\vspi
[For $\tau^{<0}X \approx \tau^{<0}Y$ and the octahedral axiom furnishes an exact triangle 
$\tau^{\geq 0}X \ra$ 
$\tau^{\geq 0}Y \ra $ 
$Z \ra$ 
$\Sigma\tau^{\geq 0}X$.]
\\
\indent\indent (II$^{\leq 0}$) \ Assume that $X \in \Ob\bC^{\leq 0}$ $-$then
$H^0(X) \ra H^0(Y) \ra H^0(Z) \ra 0$ is exact.
\vspi
Reduce the general case to II$^{\geq 0}$ $\&$ II$^{\leq 0}$.]\\
\endgroup %%------------------------------------<<

\begingroup%%----------------------------------->>
\fontsize{9pt}{11pt}\selectfont
Notation: $H^q:\bC \ra \bH(\bC)$ is the functor that sends \mX to 
$
\begin{cases}
\ H^0(\Sigma^q X) \quad (q > 0)\\
\ H^0(\Omega^q X) \quad (q < 0)
\end{cases}
.
$
\vspace{0.25cm}
\endgroup %%------------------------------------<<

\begingroup%%----------------------------------->>
\fontsize{9pt}{11pt}\selectfont
\textbf{\small FACT} \ 
Assume: The intersections \ 
$\ds\bigcap\limits_n \Ob\bC^{\leq n}$,
$\ds\bigcap\limits_n \Ob\bC^{\geq n}$ 
contain only zero objects $-$then $H^q(X) = 0$ $\forall \ q$ $\implies$ $X = 0$, thus the $H^q$ comprise a conservative system of functors (i.e., $f$ is an isomorphism iff $H^q(f)$ is an isomorphism $\forall \ q$).
\vspi
[Note: \ The objects of $\bC^{\leq n}$ are characterized by the condition that $H^q(X) = 0$ $(q > n)$ and 
the objects of $\bC^{\geq n}$ are characterized by the condition that $H^q(X) = 0$ $(q < n)$.]\\
\endgroup %%------------------------------------<<

%%%%%%%%%%%%%%%%%%%%%%%%%%%%%%%%%%%%%%
%%%%%%%%%%%%%%%%%%%%%%%%%%%%%%%%%%%%%%
%%%%%%%%%%%%%%%%%%%%%%%%%%%%%%%%%%%%%%

\begin{center}
$\S \ 15$
\\[0.5cm]
$\mathcal{REFERENCES}$\\[-.2cm]
\end{center}

\[
\textbf{BOOKS}
\]

\begingroup
\fontsize{9pt}{11pt}\selectfont
\setlength\parindent{0 cm}

[1] \quad Gelfand, S. and Manin, Y., 
\textit{Methods of Homological Algebra}, Springer Verlag (1996).
\\[-.2cm]

[2] \quad Iverson, B., 
\textit{Cohomology of Sheaves}, Springer Verlag (1986).
\\[-.2cm]

[3] \quad Kashiwara, M. and Schapira, P., 
\textit{Sheaves on Manifolds}, Springer Verlag (1990).
\\[-.2cm]

[4] \quad Margolis, H., 
\textit{Spectra and the Steenrod Algebra}, North Holland (1983).
\\[-.2cm]

[5] \quad Weibel, C., 
\textit{An Introduction to Homological Algebra}, Cambridge University Press (1994).
\\[-.2cm]
\endgroup

\[
\textbf{ARTICLES}
\]

\begingroup
\fontsize{9pt}{11pt}\selectfont
\setlength\parindent{0 cm}

[1] \quad Beilinson, A., Bernstein, J., and Deligne, P., Faisceaux Pervers, 
\textit{Ast\'erisque} \textbf{100} (1982), 1-172.
\\[-.2cm]

[2] \quad Christensen, D. and Strickland, N., Phantom Maps and Homology Theories, 
\textit{Topology} \textbf{37} (1998), 

\hspace{0.8cm}339-364.
\\[-.2cm]

[3] \quad Dold, A. and Puppe, D., Duality, Trace, and Transfer, In: 
\textit{Proceedings of the International Conference}

\hspace{0.8cm}\textit{on Geometric Topology}, K. Borsuk and A. Kirkor (ed.), PWN (1980), 81-102.
\\[-.2cm]

[4] \quad Franke, J., Uniqueness Theorems for Certain Triangulated Categories Possessing an Adams Spectral 

\hspace{0.8cm}Sequence, 
\textit{K-Theory Preprint Archives}, \textbf{139} (1996).
\\[-.2cm]

[5] \quad Grivel, P., Cat\'egories D\'eriv\'ees et Foncteurs D\'eriv\'es, In: 
\textit{Algebraic $D$-Modules}, A. Borel (ed.), Aca-

\hspace{0.8cm}demic Press (1987), 1-108.
\\[-.2cm]

[6] \quad Hovey, M., Palmieri, J., and Strickland, N., Axiomatic Stable Homotopy Theory, 
\textit{Memoirs Amer.}

\hspace{0.8cm}\textit{Math. Soc.} \textbf{610} (1997), 1-114.
\\[-.2cm]

[7] \quad Illusie, L., Cat\'egories D\'eriv\'ees et Dualit\'e, 
\textit{Enseign. Math.} \textbf{36} (1990), 369-391.
\\[-.2cm]

[8] \quad Keller, B., Derived Categories and Their Uses, In: 
\textit{Handbook of Algebra}, M. Hazewinkel (ed.), North 

\hspace{0.8cm}Holland (1996), 671-701.
\\[-.2cm]

[9] \quad Neeman, A., The Connection between the K-Theory Localization Theorem of Thomason, 

\hspace{0.8cm}Trobaugh and Yao and the Smashing Subcategories of Bousfield and Ravenel, 
\textit{Ann. Sci. \'Ecole Norm.}

\hspace{0.8cm}\textit{Sup.} \textbf{25} (1992), 547-566.
\\[-.2cm]

[10] \quad Neeman, A., The Grothendieck Duality Theorem via Bousfield's Techniques and Brown Repre-

\hspace{0.95cm}sentability, 
\textit{J. Amer. Math. Soc.} \textbf{9} (1996), 205-236.
\\[-.2cm]

[11] \quad Neeman, A., On a Theorem of Brown and Adams, 
\textit{Topology} \textbf{36} (1997), 619-645.
\\[-.2cm]

[12] \quad Thomason, R., The Classification of Triangulated Subcategories, 
\textit{Compositio Math.} \textbf{105} (1997), 1-27.
\\[-.2cm]

[13] \quad Verdier, J-L., Cat\'egories D\'eriv\'ees, 
\textit{SLN} \textbf{569} (1977), 262-311; see also 
\textit{Ast\'erisque} \textbf{239} (1996), 1-253.

\setlength\parindent{2em}

\endgroup

\chapter{
$\boldsymbol{\S}$\textbf{16}.\quadx  SPECTRA}
\setlength\parindent{2em}
\setcounter{proposition}{0}
\setcounter{chapter}{16}

%%----------------------------------------------------------------------------------------------01
$\text{ }$\\[-1.25cm]

In this $\S$, I shall give a concise exposition of the theory of spectra, concentrating on foundational issues and using model category theoretic methods whenever possible to ease the way.

A prespectrum \bX is said to be 
\un{separated}
\index{separated prespectrum}
 if $\forall \ q$, 
$\sigma_q:X_q \ra \Omega X_{q+1}$ is a \bCG embedding.  
\bSEPPRESPEC 
\index{\bSEPPRESPEC} 
is the full subcategory of \bPRESPEC whose objects are the separated prespectra.

Notation: Given a continuous function $f:X \ra Y$, where \mX $\&$ \mY are compactly generated, write 
$\im f$ for $kf(X)$ (so $f:X \ra Y$ factors as $X \ra \im f \ra Y$ and $\im f \ra Y$ is a \bCG embedding).\\

\begin{proposition} \ 
\bSEPPRESPEC is a reflective subcategory of \bPRESPEC.
\end{proposition}

[We shall construct the reflector $E^\infty$ by transfinite induction.

Claim: \ There is a functor 
$E: \bPRESPEC \ra \bPRESPEC$ and a natural transformation 
$\Xi:\id \ra E$ such that $\forall \ \bX$, $\Xi_{\bX}:\bX \ra E\bX$ is a levelwise surjection, \bX being separated iff $\Xi_{\bX}$ is a levelwise homeomorphism.  
In addition, if $\bff:\bX \ra \bY$ is a morphism of prespectra and if \bY is separated, then \bff factors uniquely through $\Xi_{\bX}$.

[Let 
$(E\bX)_q = \im(X_q \overset {\sigma_q}{\lra} \Omega X_{q+1})$ and determine the arrow 
$(E\bX)_q \ra \Omega (E \bX)_{q+1}$ from the commutative diagram
\begin{tikzcd}%[sep=small]
{X_q} \ar{d}[swap]{\Xi_{\bX,q}} \ar{r}{\sigma_q} &{\Omega X_{q+1}} \ar{d}{\Omega \Xi_{\bX,q+1}}\\
{(E\bX)_q} \ar[dashed]{r} &{\Omega(E\bX)_{q+1}}\\
{\bigcap} &{\bigcap}\\
%{\Omega X_{q+1}} \ar{r}[swap]{\Omega \sigma_{q+1}} &{\Omega\Omega X_{\vspace{0.3cm}}}
{\Omega X_{q+1}} \ar{r}[swap]{\Omega \sigma_{q+1}} &{\Omega\Omega X_{q+2}}
\end{tikzcd}
.  It is clear that \mE is functorial and $\Xi$ is natural.]

Claim: For each ordinal $\alpha$, there is a functor 
$E^\alpha:\bPRESPEC \ra \bPRESPEC$ and for each pair $\alpha \leq \beta$ of ordinals, there is a natural transformation 
$\Xi^{\alpha,\beta}:E^\alpha \ra E^\beta$ such that $\forall \ \bX$, 
$\Xi_{\bX}^{\alpha,\beta}:E^\alpha\bX \ra E^\beta\bX$ is a levelwise surjection, $E^\alpha \bX$ being separated iff 
$\Xi_{\bX}^{\alpha,\alpha+1}:E^\alpha\bX \ra E^{\alpha+1}\bX$ is a levelwise homeomorphism.  In addition, if 
$\bff:\bX \ra \bY$ is a morphism of prespectra and if \bY is separated, then \bff factors uniquely through 
$\Xi_{\bX}^{0,\alpha}$.

[Here, 
$E^0 = \id$, 
$E^1 = E$, 
$\Xi^{0,1} = \Xi$, 
$\Xi^{\alpha,\alpha} = \id$, 
$E^{\alpha + 1} = E \circx E^\alpha$, 
and 
$\Xi^{\alpha,\beta+1} = \Xi \circx \Xi^{\alpha,\beta}$ $(\alpha \leq \beta)$.  
At a limit ordinal $\lambda$, put 
$E^\lambda \bX = \underset{\alpha < \lambda}{\colimx} E^\alpha \bX$ and define
$\Xi_{\bX}^{\alpha,\lambda}:E^\alpha \bX \ra E^\lambda \bX$ in the obvious manner.]

[Note: \ If $E^\alpha \bX$ is separated, then $\forall \ \beta \geq \alpha$, 
$\Xi_{\bX}^{\alpha,\beta}:E^\alpha \bX \ra E^\beta \bX$ is a levelwise homeomorphism.]

%%----------------------------------------------------------------------------------------------02
To finish the proof, it suffices to show that $\forall \ \bX$, $\exists$ an $\alpha_{\bX}$ such that $E^{\alpha_{\bX}}\bX$ is separated.  But for this, one can take $\alpha_{\bX}$ to be any infinite cardinal greater than the cardinality of 
$\bigl(\coprod\limits_q X_q \times X_q \bigr) \amalg \bigl(\coprod\limits_q \tau(X_q)\bigr)$ 
$(\tau(X_q)$ the set of open subsets of $X_q$).]

[Note: \ The arrow of reflection $\bX \ra E^\infty \bX$ is a levelwise surjection.  It is a levelwise homeomorphism iff \bX is separated.]\\

\begingroup%%----------------------------------->>
\fontsize{9pt}{11pt}\selectfont
The existence of the reflector $E^\infty$ can be established by applying the general adjoint functor theorem: 
\bSEPPRESPEC is a priori complete, the inclusion $\bSEPPRESPEC \ra \bPRESPEC$ preserves limits, and the solution set condition is satisfied.  The drawback to this approach is that it provides no information about the behavior of $E^\infty$ with respect to finite limits, a situation that can be partially clarified by using the iterative definition of $E^\infty$ in terms of the 
$E^\alpha$.\\
\endgroup %%------------------------------------<<

\begingroup%%----------------------------------->>
\fontsize{9pt}{11pt}\selectfont
\textbf{\small LEMMA} \  
Suppose that $(I,\leq)$ is a nonempty directed set, regarded as a filtered category \bI.  Let 
$\Delta^\prime, \ \Delta\pp: \bI \ra \dcg$ be diagrams $-$then the arrow 
$\colim_{\bI}(\Delta^\prime \times \Delta\pp) \ra$ 
$\colim_{\bI} \Delta^\prime \times_k \colim_{\bI} \Delta\pp$ is a homeomorphism.
\vspi
[Note: \ The directed colimit in $\dcg_*$ is formed by assigning the evident base point to the corresponding directed colimit in $\dcg$, thus the lemma is valid in $\dcg_*$ as well.]\\
\endgroup %%------------------------------------<<

\begingroup%%----------------------------------->>
\fontsize{9pt}{11pt}\selectfont
\textbf{\small FACT} \  
$E^\infty$ preserves finite products.
\vspi
[Note: \ $E^\infty$ does not preserve equalizers.]\\
\endgroup %%------------------------------------<<

\begingroup%%----------------------------------->>
\fontsize{9pt}{11pt}\selectfont
\textbf{\small LEMMA} \  
Suppose that $(I,\leq)$ is a nonempty directed set, regarded as a filtered category \bI.  Let 
$\Delta: \bI \ra \dcg$ be a diagram such that $\forall$ 
$i \overset{\delta}{\ra} j$, 
$\Delta\delta:\Delta_i \ra \Delta_j$ is an injection $-$then $\colim_{\bI} \Delta$ in 
$\dcg =$ $\colim_{\bI} \Delta$ in \bCG ($= \colim_{\bI} \Delta$ in \bTOP) and $\forall \ i$, the canonical arrow 
$\Delta_i \ra \colim_{\bI} \Delta$ is one-to-one.
\vspi
[Note: \ The set underlying $\colim_{\bI} \Delta$ is therefore the colimit of the underlying diagram in \bSET.]\\
\endgroup %%------------------------------------<<

\begingroup%%----------------------------------->>
\fontsize{9pt}{11pt}\selectfont
\textbf{\small LEMMA} \  
In $\dcg$, directed colimits of diagrams whose arrows are injections commute with finite limits.
\vspi
[Note: \ A finite limit in $\dcg_*$ is formed by assigning the evident base point to the corresponding finite limit in 
$\dcg$, thus the lemma is valid in $\dcg_*$ as well.]\\
\endgroup %%------------------------------------<<

\begingroup%%----------------------------------->>
\fontsize{9pt}{11pt}\selectfont
A prespectrum \bX is said to be 
\un{injective}
\index{injective prespectrum} 
if $\forall \ q$, 
$\sigma_q:X_q \ra \Omega X_{q+1}$ is an injection.  \bINJPRESPEC 
\index{\bINJPRESPEC} 
 is the full subcategory of \bPRESPEC whose objects are the injective prespectra.
\vspi
 [Note: \bSEPPRESPEC is a full subcategory of \bINJPRESPEC.]\\
\endgroup %%------------------------------------<<

\begingroup%%----------------------------------->>
\fontsize{9pt}{11pt}\selectfont
\textbf{\small FACT} \  
The arrow of reflection $\bX \ra E^\infty\bX$ is a levelwise injection iff \bX is injective.
\vspi
[If \bX is injective, then so are the $E^\alpha \bX$.  
Moreover, 
$\Xi_{\bX}^{\alpha,\beta}:E^\alpha \bX \ra E^\beta \bX$ $(\alpha \leq \beta)$ is one-to-one.]
\vspi
%%----------------------------------------------------------------------------------------------03
[Note: \ It therefore follows that the arrow of reflection $\bX \ra E^\infty\bX$ is a levelwise bijection iff \bX is injective.]\\
\endgroup %%------------------------------------<<

\begingroup%%----------------------------------->>
\fontsize{9pt}{11pt}\selectfont
\textbf{\small FACT} \  
The restriction of $E^\infty$ to \bINJPRESPEC preserves finite limits.\\
\endgroup %%------------------------------------<<

\textbf{\small LEMMA} \  
Suppose given a sequence $\{X_n,f_n\}$, where $X_n$ is a \dsp compactly generated space and 
$f_n:X_n \ra X_{n+1}$ is a \bCG embedding $-$then $\forall$ compact Hausdorff space $K$, 
$\colimx X_n^K \approx$ $(\colimx X_n)^K$ (exponential objects in $\dcg$).

[Note: \ There is an analogous assertion in the pointed category.]\\

\begin{proposition} \ %02
\bSPEC is a reflective subcategory of \bSEPPRESPEC.
\end{proposition}

[The reflector sends \bX to $e\bX$, the latter being defined by the rule 
$q \ra \colim\Omega^n X_{n+q}$.]\\

\begingroup%%----------------------------------->>
\fontsize{9pt}{11pt}\selectfont
\textbf{\small LEMMA} \  
Suppose that $(I, \leq)$ is a nonempty directed set, regarded as a filtered category \bI.  
Let 
$\Delta:\bI \ra \dcg$ be a diagram such that $\forall$ 
$i \overset{\delta}{\ra} j$, 
$\Delta\delta: \Delta_i \ra \Delta_j$ is a \bCG embedding $-$then $\forall \ i$, the canonical arrow 
$\Delta_i \ra$ $\colim_{\bI}\Delta$ is a \bCG embedding.
\vspi
[Note: \ Changing the assumption to ``closed embedding'' changes the conclusion to ``closed embedding''.]\\
\endgroup %%------------------------------------<<

\begingroup%%----------------------------------->>
\fontsize{9pt}{11pt}\selectfont
\textbf{\small FACT} \  
The arrow of reflection $\bX \ra e\bX$ is a levelwise \bCG embedding.\\ 
\endgroup %%------------------------------------<<

\begingroup%%----------------------------------->>
\fontsize{9pt}{11pt}\selectfont
\textbf{\small FACT} \  
$e$ preserves finite limits.\\
\endgroup %%------------------------------------<<

\begin{proposition} \ 
\bSPEC is a reflective subcategory of \bPRESPEC.
\end{proposition}

[This is implied by Propositions 1 and 2.]

[Note: \ The composite \quad
$\bPRESPEC \overset{E^\infty}{\lra} \bSEPPRESPEC \overset{e}{\lra} \bSPEC$ \quad is the 
\un{spectrification functor}:
\index{spectrification functor} 
$\bX \ra s\bX$ $(s = e \circx E^\infty)$.]\\

\label{16.6}
Application: \bSPEC is complete and cocomplete.

[Note: \ The colimit of a diagram $\Delta:\bI \ra \bSPEC$ is the spectrification of its colimit in \bPRESPEC.  
Example: The coproduct in \bPRESPEC or \bSPEC is denoted by a wedge.  
If $\{\bX_i\}$ is a set of spectra, then its coproduct in \bPRESPEC is separated, so 
$e\bigl(\bigvee\limits_i \bX_i\bigr)$ is the coproduct $\bigvee\limits_i \bX_i$ of the $\bX_i$ in \bSPEC.]\\

\begingroup%%----------------------------------->>
\fontsize{9pt}{11pt}\selectfont
\textbf{\small FACT} \  
Spectrification preserves finite products and its restriction to \bINJPRESPEC preserves finite limits.\\
\endgroup %%------------------------------------<<

\begingroup%%----------------------------------->>
\fontsize{9pt}{11pt}\selectfont
\textbf{\small EXAMPLE} \  
Let \mX be in $\dcg_*$ $-$then the 
\un{suspension prespectrum}
\index{suspension prespectrum} 
of \mX is the assignment 
$q \ra \Sigma^q X$, where 
$\Sigma^q X \ra$
$\Omega\Sigma\Sigma^q X \approx$ 
$\Omega\Sigma^{q+1} X$ (a \bCG embedding).  Its spectrification is the 
\un{suspension spectrum}
\index{suspension spectrum}
%%----------------------------------------------------------------------------------------------04
of $X$.  Thus, in the notation of p. \pageref{16.1}, the suspension spectrum of \mX is 
$\bQ^\infty X$: $(\bQ^\infty X)_q =$ $\colimx \Omega^n \Sigma^{n+q} X =$ 
$\Omega^\infty\Sigma^\infty\Sigma^q X.$\\
\endgroup %%------------------------------------<<

\begingroup%%----------------------------------->>
\fontsize{9pt}{11pt}\selectfont
\textbf{\small EXAMPLE} \  
Fix $q \geq 0$.  Given an $X$ in $\dcg_*$, let $\bQ_q^\infty X$ be the spectrification of the prespectrum 
$
p \ra 
\begin{cases}
\ \Sigma^{p-q}X \quad (p \geq q)\\
\ * \qquad\quad \ \  (p < q)
\end{cases}
, \ 
$
where
$\Sigma^{p-q}X \ra \Omega\Sigma\Sigma^{p-q}X \approx \Omega\Sigma^{p+1 - q}X$ $(p \geq q)$ (if $p < q$, 
the arrow is the inclusion of the base point).  Viewed as a functor from $\dcg_*$ to $\bSPEC$, $\bQ_p^\infty$ is a left adjoint for the $q^\text{th}$ space functor $\bU_q^\infty:\bSPEC \ra \dcg_*$ that sends $\bX = \{X_q\}$ to $X_q$.  
Special case: 
$\bQ_0^\infty = \bQ^\infty$, 
$\bU_0^\infty = \bU^\infty$.

[Note: \ $\forall \ X$, $q^\prime \leq q\pp$ $\implies$ 
$\bQ_{q^\prime}^\infty X \approx$ 
$\bQ_{q\pp}^\infty \Sigma^{q\pp - q^\prime}X$.]\\ 
\endgroup %%------------------------------------<<

\begingroup%%----------------------------------->>
\fontsize{9pt}{11pt}\selectfont

\textbf{\small FACT} \  
Suppose that \bX is a prespectrum $-$then $s\bX \approx \colimx \bQ_q^\infty X_q$.
\vspi
[For any spectrum \bY, 
$\Mor(\colimx \bQ_q^\infty X_q,\bY) \approx$ 
$\lim \Mor(\bQ_q^\infty X_q, \bY) \approx$ 
$\lim\Mor(X_q,Y_q) \approx$
$\Mor(\bX,\bY) \approx$
$\Mor(s\bX,\bY)$.]\\
\endgroup %%------------------------------------<<

\begingroup%%----------------------------------->>
\fontsize{9pt}{11pt}\selectfont
\label{16.8}
\label{16.13}
\textbf{\small FACT} \  
Let $(\bX,\bff)$ be an object in $\bFIL(\bSPEC)$ (cf. p. \pageref{16.2}).  
Assume $\forall \ n$, 
$\bff_n:\bX_n \ra \bX_{n+1}$ is a levelwise \bCG embedding $-$then $\forall$ pointed compact Hausdorff space $K$, 
$\colimx \Mor(\bQ_q^\infty K,\bX_n) \approx$
$\Mor(\bQ_q^\infty K,\colimx \bX_n)$.
\vspi
[The assumption guarantees that the prespectrum colimit of $(\bX,\bff)$ is a spectrum.  Therefore 
$\colimx $ $\Mor(\bQ_q^\infty K,\bX_n)$ \hspace{0.03cm} $\approx$ \hspace{0.03cm}
$\colimx \Mor(K,\bU_q^\infty \bX_n)$ \hspace{0.03cm} $\approx$ \hspace{0.03cm}
$\Mor(K,\colimx \bU_q^\infty \bX_n)$ \hspace{0.03cm} $\approx$ \hspace{0.03cm}
$\Mor(K,\bU_q^\infty \colimx \bX_n)$ \hspace{0.03cm} $\approx$  \hspace{0.03cm}
$\Mor(\bQ_q^\infty K,\colimx \bX_n)$.]\\
\endgroup %%------------------------------------<<

\begingroup%%----------------------------------->>
\fontsize{9pt}{11pt}\selectfont
\textbf{\small FACT} \  
Let $\{\bX_i\}$ be a set of spectra, $K$ a pointed compact Hausdorff space $-$then every morphism 
$\bff:\bQ_q^\infty K \ra \ds\bigvee\limits_i \bX_i$ factors through a finite subwedge.
\vspi
[Since 
$\Mor \big(\bQ_q^\infty K, \ds\bigvee\limits_i \bX_i \big) \approx$ 
$\Mor \big(K, \bU_q^\infty \big(\ds\bigvee\limits_i \bX_i \big)\big)$, \bff corresponds to an arrow 
$g:K \ra \bU_q^\infty \big(\ds\bigvee\limits_i \bX_i \big)$ $\big(= \big(\ds\bigvee\limits_i \bX_i\big)_q\big)$, 
i.e., to an arrow 
$g:K \ra \colimx \Omega^n  \big(\ds\bigvee\limits_i \bX_i \big)_{n+q}$, 
which factors through 
$\Omega^n  \big(\ds\bigvee\limits_i (\bX_i)_{n+q}\big)$ for some $n$: 
$
\begin{tikzcd}%[sep=small]
{K} \ar{rd}[swap]{g} \ar{r}{g_n} &{\Omega^n  \big(\ds\bigvee\limits_i \big(\bX_i \big)_{n+q}}\big) \ar{d}\\
&{\big(\ds\bigvee\limits_i \bX_i \big)_q}
\end{tikzcd}
.  
$ \ 
The adjoint 
$\ov{g}_n:\Sigma^n K \ra \ds\bigvee\limits_i \big(\bX_i \big)_{n+q}$ factors through a finite subwedge 
$\ds\bigvee\limits_k \big(\bX_{i_k}\big)_{n+q}$, so \bff factors through $\ds\bigvee\limits_k \bX_{i_k}$.]\\
\endgroup %%------------------------------------<<

Notation: Given \bX, \bY in \bPRESPEC, write $\HOM(\bX,\bY)$ for $\Mor(\bX,\bY)$ topologized via the equalizer diagram 
$\Mor(\bX,\bY) \ra \prod\limits_q Y_q^{X_q} \rightrightarrows \prod\limits_q (\Omega Y_{q+1})^{X_q}$.

\begin{proposition} \ 
Spectrification is a continuous functor in the sense that $\forall$ \bX, \bY in \bPRESPEC, the arrow 
$\HOM(\bX,\bY) \ra \HOM(s\bX,s\bY)$ is a continuous function.\\
\end{proposition}

%%----------------------------------------------------------------------------------------------05
\indent\indent ($\bbox$ and $\wedge$) \ Fix a $K$ in $\dcg_*$.  
Given an \bX in \bPRESPEC, let 
$\bX \bbox K$ be the prespectrum 
$q \ra X_q\#_k K$, where 
$X_q\#_k K \ra$ $\Omega(X_{q+1} \#_k K)$ is 
$X_q \#_k K \ra$ 
$\Omega X_{q+1} \#_k  K \ra$ 
$\Omega(X_{q+1} \#_k K)$, and given an \bX in \bSPEC, let $\bX \wedge K$ be the spectrification of 
$\bX \bbox K$.

Examples: 
(1) $\Gamma \bX = \bX \bbox [0,1]$ or $\bX \wedge [0,1]$, the 
\un{cone}
\index{cone (prespectra)} 
of \bX; 
(2) $\Sigma \bX = \bX \bbox \bS^1$ or $\bX \wedge \bS^1$, the 
\un{suspension} 
\index{suspension (prespectra)} 
of \bX.

\indent\indent ($\texttt{HOM}$) \ Fix a $K$ in $\dcg_*$.  Given an \bX in \bPRESPEC, let $\texttt{HOM}(K,\bX)$ be the prespectrum 
$q \ra X_q^K$, where $X_q^K \ra \Omega X_{q+1}^K$ is 
$X_q^K \ra$ 
$(\Omega X_{q+1})^K \approx$ $\Omega X_{q+1}^K$.

[Note: \ If \bX is a spectrum, then $\texttt{HOM}(K, \bX)$ is a spectrum.]

Example: $\forall \ \bX$, $\Omega \bX = \texttt{HOM}(\bS^1,\bX)$ (cf. p. \pageref{16.3}).\\

\begin{proposition} \ %5
For \bX, \bY in \bPRESPEC and \mK in $\dcg_*$, there are natural homeomorphisms 
$\HOM(\bX \bbox K,\bY) \approx \HOM(\bX,\bY)^K \approx \HOM(\bX,\texttt{HOM}(K,\bY))$.
\end{proposition}

[Note: \ Consequently, the functor 
$\bX \bbox -:\dcg_* \ra \bPRESPEC$ has a right adjoint, viz. $\HOM(\bX,-)$, and the functor 
$-\bbox K:\bPRESPEC \ra \bPRESPEC$ has a right adjoint, viz. $\texttt{HOM}(K,-)$.]\\

\begin{proposition} \ %6
For \bX, \bY in \bSPEC and \mK in $\dcg_*$, there are natural homeomorphisms 
$\HOM(\bX \wedge K,\bY) \approx \HOM(\bX,\bY)^K \approx \HOM(\bX,\texttt{HOM}(K,\bY))$.
\end{proposition}

[Note: \ Consequently, the functor 
$\bX \wedge -:\dcg_* \ra \bSPEC$ has a right adjoint, viz. $\HOM(\bX,-)$ and the functor 
$-\wedge K:\bSPEC \ra \bSPEC$ has a right adjoint, viz. $\texttt{HOM}(K,-)$.]\\

Examples: 
(1) 
$\bQ_q^\infty (K \#_k L) \approx (\bQ_q^\infty K) \wedge L$ and 
$\bU_q^\infty \texttt{HOM}(K,\bX) \approx (\bU_q^\infty \bX)^K$;
(2) 
$s(\bX \bbox K) \approx s\bX \wedge K$.

Example: $(\Sigma,\Omega)$ is an adjoint pair.\\

\begingroup%%----------------------------------->>
\fontsize{9pt}{11pt}\selectfont
\textbf{\small EXAMPLE} \  
(1) $\bX \wedge \bS^0 \approx \bX$; 
(2) $\texttt{HOM}(\bS^0,\bX) \approx \bX$; 
(3) $(\bX \wedge K) \wedge L \approx \bX \wedge (K \#_k L)$; 
(4) $\texttt{HOM}(K \#_k L,\bX) \approx \texttt{HOM}(K,\texttt{HOM}(L,\bX))$.\\
\endgroup %%------------------------------------<<

\begingroup%%----------------------------------->>
\fontsize{9pt}{11pt}\selectfont
\textbf{\small FACT} \  
Suppose that \bX is an injective prespectrum $-$then $\forall \ K$, $\bX \bbox K$ is an injective prespectrum.\\
\endgroup %%------------------------------------<<

\begingroup%%----------------------------------->>
\fontsize{9pt}{11pt}\selectfont
\textbf{\small FACT} \  
Suppose that \bX is a separated prespectrum $-$then $\forall$ nonempty compact Hausdorff space $K$, $\bX \bbox K_+$ is a separated prespectrum.\\
\endgroup %%------------------------------------<<

\begingroup%%----------------------------------->>
\fontsize{9pt}{11pt}\selectfont
\textbf{\small LEMMA} \  
Suppose that 
\begin{tikzcd}[sep=large]
{P} \ar{d}[swap]{\xi} \ar{r}{\eta} &{Y} \ar{d}{g}\\
{X} \ar{r}[swap]{f} &{Z}
\end{tikzcd}
is a pullback square in $\dcg$.  Assume: $g$ is a closed embedding $-$then $\xi$ is a closed embedding.\\
\endgroup %%------------------------------------<<

%%----------------------------------------------------------------------------------------------06
\begingroup%%----------------------------------->>
\fontsize{9pt}{11pt}\selectfont
\textbf{\small EXAMPLE} \  
Let $\bff:\bX \ra \bY$ be a morphism of prespectra $-$then the 
\un{mapping cylinder}
\index{mapping cylinder (of a morphism of prespectra)} 
$\bM_\bff$ of \bff is defined by the pushout square 
\begin{tikzcd}[sep=large]
{\bX \bbox \{0\}_+} \ar{d} \ar{r} &{\bY \bbox \{0\}_+} \ar{d}\\
{\bX \bbox I_+} \ar{r} &{\bM_\bff}
\end{tikzcd}
$(I_+ = [0,1] \amalg *$ (cf. p. \pageref{16.4})).  
There is a natural arrow $\bM_\bff \ra \bY \bbox I_+$ and the commutative diagram 
\begin{tikzcd}[sep=large]
{\bX} \ar{d} \ar{r} &{\bM_\bff} \ar{d}\\
{\bY} \ar{r} &{\bY \bbox I_+}
\end{tikzcd}
is a pullback square.  
Definition: \bff is a 
\un{prespectral cofibration}
\index{prespectral cofibration} 
if 
$\bM_\bff \ra \bY \bbox I_+$ has a left inverse.  
Every prespectral cofibration is a levelwise closed embedding.\\
\endgroup %%------------------------------------<<

\begingroup%%----------------------------------->>
\fontsize{9pt}{11pt}\selectfont
\textbf{\small FACT} \  
Let $\bff: \bX \ra \bY$ be a morphism of prespectra.  Assume: 
$
\begin{cases}
\ \bX\\
\ \bY
\end{cases}
$
are injective $-$then $\bM_{\bff}$ is injective.\\
\endgroup %%------------------------------------<<

\begingroup%%----------------------------------->>
\fontsize{9pt}{11pt}\selectfont
\textbf{\small EXAMPLE} \  
Let $\bff:\bX \ra \bY$ be a morphism of spectra $-$then the 
\un{mapping cylinder}
\index{mapping cylinder (of a morphism of spectra)} 
$\bM_\bff$ of \bff is defined by the pushout square 
\begin{tikzcd}%[sep=small]
{\bX \wedge \{0\}_+} \ar{d} \ar{r} &{\bY \wedge \{0\}_+} \ar{d}\\
{\bX \wedge I_+} \ar{r} &{\bM_\bff}
\end{tikzcd}
$(I_+ = [0,1] \amalg *$ (cf. p. \pageref{16.5})).  
There is a natural arrow $\bM_\bff \ra \bY \wedge I_+$ and the commutative diagram 
\begin{tikzcd}%[sep=small]
{\bX} \ar{d} \ar{r} &{\bM_\bff} \ar{d}\\
{\bY} \ar{r} &{\bY \wedge I_+}
\end{tikzcd}
is a pullback square.  
Indeed, the mapping cylinder of \bff in \bSPEC is the spectrification of the mapping cylinder of \bff in \bPRESPEC.  And: 
All data is injective, so 
$s\left(
\begin{tikzcd}%[sep=small]
{\bX} \ar{d} \ar{r} &{\bM_\bff} \ar{d}\\
{\bY} \ar{r} &{\bY \bbox I_+}
\end{tikzcd}
\right)$ is a pullback square in \bSPEC (cf. p. \pageref{16.6}).  
\label{16.9}
Definition: \bff is a 
\un{spectral cofibration}
\index{spectral cofibration}
 if 
$\bM_{\bff} \ra \bY \wedge I_+$ has a left inverse.  
\label{16.12}
\label{16.27}
Every spectral cofibration is a levelwise closed embedding.

[Note: \ The arrow $\bff:\bX \ra \bY$ is a spectral cofibration iff the commutative diagram
\begin{tikzcd}%[sep=small]
{\bX \wedge \{0\}_+} \ar{d} \ar{r} &{}\\
{\bX \wedge I_+} \ar{r} &{}
\end{tikzcd}
\begin{tikzcd}%[sep=small]
{\bY \wedge \{0\}_+} \ar{d}\\
{\bY \wedge I_+}
\end{tikzcd}
\ \ is a weak pushout square or, equivalently, iff $\forall \ \bZ$, \bff has the LLP w.r.t. 
$\texttt{HOM}(I_+,\bZ) \overset{p_0}{\lra} \bZ$.  
Example: Suppose that $L \ra K$ is a pointed cofibration $-$then \ $\forall \ \bX$, \ 
$\bX \wedge L \ra \bX \wedge K$ is a spectral cofibration.]\\
\endgroup %%------------------------------------<<

Notation: For $n \geq 0$, put 
$\bS^n = \bQ^\infty \bS^n$ and for $n > 0$, put 
$\bS^{-n} = \bQ_n^\infty \bS^0$.

[Note: \ $\forall \ n$ $\&$ $\forall \ m \geq 0$, 
$\Sigma^m\bS^n$ $(= \bS^n \wedge \bS^m) \approx \bS^{m+n}$ and 
$\forall \ n \geq 0$ $\&$ $\forall \ m \geq 0$,
$\bS^{-m} \wedge \bS^n \approx$ 
$(\bQ_m^\infty \bS^0) \wedge \bS^n \approx$ 
$\bQ_m^\infty (\bS^0 \#_k \bS^n) \approx$ 
$\bQ_m^\infty \bS^n \approx$ 
$\bS^{n - m}$.]\\

%%----------------------------------------------------------------------------------------------07
\begingroup%%----------------------------------->>
\fontsize{9pt}{11pt}\selectfont
\textbf{\small EXAMPLE} \  
$\forall \ X$, $\bQ_q^\infty X = \bS^{-q} \wedge X$.  So, the arrow of adjunction 
$\id \ra \bU_q^\infty  \circx \bQ_q^\infty$ is given by 
$X \ra (\bS^{-q} \wedge X)_q$ and the arrow of adjunction 
$\bQ_q^\infty \circx \bU_q^\infty \ra \id$ is given by 
$\bS^{-q} \wedge X_q \ra \bX$.\\
\endgroup %%------------------------------------<<

\begin{proposition} \ 
The $q^\text{th}$ space functor $\bU_q^\infty:\bSPEC \ra \dcg_*$ is represented by $\bS^{-q}$.
\end{proposition}

[$\forall \ \bX$, $\Mor(\bS^{-q},\bX) = \Mor(\bQ_q^\infty \bS^0,\bX) \approx 
\Mor(\bS^0,\bU_q^\infty\bX) = \bU_q^\infty\bX$.]\\

A 
\un{homotopy}
\index{homotopy in \bSPEC} 
in \bSPEC is an arrow $\bX \wedge I_+ \ra \bY$.  
Homotopy is an equivalence relation which respects composition, so there is an associated quotient category $\bSPEC/\simeq$: 
$[\bX,\bY]_0 = \Mor(\bX,\bY)/\simeq$, i.e., $[\bX,\bY]_0 = \pi_0(\HOM(\bX,\bY))$.\\

\label{17.42} %dmc mnft
\begingroup%%----------------------------------->>
\fontsize{9pt}{11pt}\selectfont
\label{16.62}
\index{Homotopy Groups of Spectra (example)}
\textbf{\small EXAMPLE (\un{Homotopy Groups of Spectra})} \  
Let \bX be a spectrum $-$then the 
\un{$n^\text{th}$ homotopy} \un{group}
\index{n$^\text{th}$ homotopy group (spectra)} 
$\pi_n(\bX)$ of \bX $(n \in \Z)$ is $[\bS^n,\bX]_0$.  
The $\pi_n(\bX)$ are necessarily abelian.  And: $\forall \ n \geq 0$, 
$\pi_n(\bX) =$ $\pi_n(X_0)$, while $\pi_{-n}(\bX) = \pi_0(X_n)$.  
Therefore \bX is connective iff 
$\pi_n(\bX) = 0$ for $n \leq -1$.  
Example: $\forall \ \bX$ in $\dcg_{*\bc}$, the suspension spectrum $\bQ^\infty X$ of X is connective.  
Proof: $\Sigma X$ is path connected and wellpointed ($\implies$ $\Sigma^2 X$ is simply connected), thus $\forall \ n \geq 1$, 
$\pi_q(\Sigma^{q+n} X) = *$ (by the suspension isomorphism and Hurewicz), so 
$\pi_{-n}(\bQ^\infty X) = \pi_0(\Omega^\infty\Sigma^\infty\Sigma^n X) =$ 
$\colimx \pi_q(\Sigma^{q+n} X) = *$.
\vspi
[Note: \ The 
\un{stable homotopy groups}
\index{stable homotopy groups} 
$\pi_n^s(X)$ $(n \geq 0)$ of \mX are the \ 
$\pi_n(\bQ^\infty X)$ \ $(= \pi_n(\Omega^\infty\Sigma^\infty X))$.  
Example: $\pi_0^s(X) \approx \widetilde{H}_0(X)$.]\\
\endgroup %%------------------------------------<<

\begingroup%%----------------------------------->>
\fontsize{9pt}{11pt}\selectfont
\textbf{\small FACT} \  
Let $(\bX,\bff)$ be an object in $\bFIL(\bSPEC)$ (cf. p. \pageref{16.7}).  Assume: $\forall \ n$, 
$\bff_n:\bX_n \ra \bX_{n+1}$ is a levelwise \bCG embedding $-$then $\forall$ pointed compact Hausdorff space $K$, 
$\colimx[\bQ_q^\infty K,\bX_n]_0 \approx$ 
$[\bQ_q^\infty K,\colimx \bX_n]_0$ (cf. p. \pageref{16.8}). \\
\endgroup %%------------------------------------<<

\begingroup%%----------------------------------->>
\fontsize{9pt}{11pt}\selectfont
\label{16.17}
\textbf{\small EXAMPLE} \  
Imitating the construction in pointed spaces, one can attach to each object $(\bX,\bff)$ in $\bFIL(\bSPEC)$ a spectrum 
$\telsub(\bX,\bff)$, its 
\un{mapping telescope}.
\index{mapping telescope $(\bFIL(\bSPEC))$}   
Thus 
$\telsub(\bX,\bff) = \colimx \telsub_n(\bX,\bff)$ and the arrow 
$\telsub_n(\bX,\bff) \ra \telsub_{n+1}(\bX,\bff)$ is a spectral cofibration (hence is a levelwise closed embedding 
(cf. p. \pageref{16.9})).  Since there are canonical homotopy equivalences 
$\telsub_n(\bX,\bff) \ra \bX_n$, it follows that $\forall$ pointed compact Hausdorff space $K$, 
$\colimx[\bQ_q^\infty K,\bX_n]_0 \approx$ 
$[\bQ_q^\infty K,\telsub(\bX,\bff)]_0$.\\
\endgroup %%------------------------------------<<

\textbf{\small LEMMA} \  
Suppose that $\bff:\bX \ra \bY$ is a homotopy equivalence $-$then $\forall \ q$, $f_q:X_q \ra Y_q$ is a 
homotopy equivalence.

[The $q^\text{th}$ space functor 
$\bU_q^\infty:\bSPEC \ra \dcg_*$ is a \bV-functor $(\bV = \dcg_*)$, hence preserves homotopies.]\\

\begingroup%%----------------------------------->>
\fontsize{9pt}{11pt}\selectfont
\textbf{\small FACT} \  
\bSPEC is a cofibration category if weak equivalence = homotopy equivalence, cofibration = spectral cofibration.  
All objects are cofibrant and fibrant.
\vspi
%%----------------------------------------------------------------------------------------------08
[Note: \ One way to proceed is t0 show that \bSPEC is an \mI-category in the sense of 
Baues\footnote[2]{\textit{Algebraic Homotopy}, Cambridge University Press (1989), 18-27.}.]\\
\endgroup %%------------------------------------<<

A prespectrum \bX is said to satisfy the 
\un{cofibration condition}
\index{cofibration condition (on prespectrum)} 
if $\forall \ q$, the arrow $\Sigma X_q \ra X_{q+1}$ adjoint to $\sigma_q$ is a pointed cofibration.  An \bX which satisifes the cofibration condition is necessarily separated (for then $\sigma_q$ is a closed embedding).  
Example: $\forall$ \bX, $M\bX$ satisfes the cofibration condition (cf. p. \pageref{16.10}).\\

\begingroup%%----------------------------------->>
\fontsize{9pt}{11pt}\selectfont
\textbf{\small EXAMPLE} \  
Equip \bPRESPEC with the model category structure supplied by Proposition 56 in $\S 14$ $-$then every cofibrant \bX satisfies the cofibration condition.
\vspi
[Note: \ The converse is false.  To see this, take any \mX in $\dcg_*$ and consider the prespectrum whose spectrification is $\bQ_q^\infty X$, bearing in mind that the inclusion of a point is always a pointed cofibration.]\\
\endgroup %%------------------------------------<<

A spectrum \bX is said to be 
\un{tame}
\index{tame (spectrum)} 
if it is homotopy equivalent to a spectrum of the form 
$s\bY$, where \bY is a prespectrum satisfying the cofibration condition $(\implies s\bY \approx e\bY)$.\\

\textbf{\small LEMMA} \  
Let $\bff:\bX \ra \bY$ be a morphism of spectra.  Assume: \bff is a levelwise pointed homotopy equivalence $-$then $\forall$ tame spectrum \bZ, $\bff_*:[\bZ,\bX]_0 \ra [\bZ,\bY]_0$ is bijective.\\

Application: A levelwise pointed homotopy equivalence between tame spectra is a homotopy equivalence of spectra.\\

\begingroup%%----------------------------------->>
\fontsize{9pt}{11pt}\selectfont
\label{18.5} %dmc mnft
\textbf{\small FACT} \  
Let $\bff:\bX \ra \bY$ be a morphism of prespectra.  Assume: 
$
\begin{cases}
\ \bX\\
\ \bY
\end{cases}
$
satisfy the cofibration condition and \bff is a levelwise pointed homotopy equivalence $-$then 
$s\bff:s\bX \ra s\bY$ is a homotopy equivalence of spectra.\\
\endgroup %%------------------------------------<<

Equip $\bDelta$-$\bCG_*$ with its singular structure.\\

\label{16.26}
\textbf{\small LEMMA} \  
Let $\bff:\bX \ra \bY$ be a morphism of spectra $-$then \bff is a levelwise fibration iff \bff has the RLP w.r.t. the spectral cofibrations \ 
$\bS^{-q} \wedge [0,1]_+^n \ra \bS^{-q} \wedge I[0,1]_+^n$ $(n \geq 0, q \geq 0)$.\\

\textbf{\small LEMMA} \  
Let $\bff:\bX \ra \bY$ be a morphism of spectra $-$then \bff is a levelwise acyclic fibration iff \bff has the RLP w.r.t. the spectral cofibrations \ 
$\bS^{-q} \wedge \bS_+^{n-1} \ra \bS^{-q} \wedge \bD_+^n$ $(n \geq 0, q \geq 0)$.\\

%%----------------------------------------------------------------------------------------------09
\begingroup%%----------------------------------->>
\fontsize{9pt}{11pt}\selectfont
Since $(\bQ_q^\infty,\bU_q^\infty)$ is an adjoint pair, the lifting problem \ 
\begin{tikzcd}[sep=large]
{\bQ_q^\infty L} \ar{d}[swap]{\bQ_q^\infty f} \ar{r} &{\bX} \ar{d}{\bff}\\
{\bQ_q^\infty K} \ar[dashed]{ru} \ar{r} &{\bY}
\end{tikzcd}
is equivalent to the lifting problem \ 
$
\begin{tikzcd}[sep=large]
{L} \ar{d}[swap]{f} \ar{r} &{\bU_q^\infty \bX} \ar{d}{\bU_q^\infty \bff}\\
{K} \ar[dashed]{ru}  \ar{r} &{\bU_q^\infty \bY}
\end{tikzcd}
.
$\\
\vspace{0.25cm}
\endgroup %%------------------------------------<<

\begin{proposition} \ %8
Equip $\dcg_*$ with its singular structure $-$then \bSPEC is a model category if weak equivalences and fibrations are levelwise, the cofibrations being those morphisms which have the LLP w.r.t. the levelwise acyclic fibrations.
\end{proposition}

[The proof is basically the same as that for the singular structure on \bTOP (cf. p. \pageref{16.11} ff.).  
Thus there are two claims.

Claim: Every morphism $\bff:\bX \ra \bY$ can be written as a composite 
$\bff_\omega \circx \bi_\omega$, where 
$\bi_\omega:\bX \ra \bX_\omega$ is a weak equivalence and has the LLP w.r.t. all fibrations and 
$\bff_\omega:\bX_\omega \ra \bY$ is a fibration.

[In the small object argument, take 
$S_0 = \{\bS^{-q} \wedge [0,1]_+^n \ra \bS^{-q} \wedge I[0,1]_+^n \}$ $(n \geq 0, q\geq 0)$ 
$-$then $\forall \ k$, the arrow $\bX_k \ra \bX_{k+1}$ is a spectral cofibration, 
hence is a levelwise closed embedding (cf. p. \pageref{16.12}).  
Since
$\bQ_q^\infty [0,1]_+^n \approx$ 
$\bS^{-q} \wedge [0,1]_+^n$, it follows that 
$\colim\  \Mor(\bS^{-q} \wedge [0,1]_+^n,\bX_k) \approx$ 
$\Mor(\bS^{-q} \wedge [0,1]_+^n,\bX_\omega)$ $\forall \ n$ (cf. p. \pageref{16.13}), 
so $\bff_\omega$ has the RLP w.r.t. the 
$\bS^{-q} \wedge [0,1]_+^n \ra$
$\bS^{-q} \wedge I[0,1]_+^n$, i.e., is a fibration.  
The assertions regarding $\bi_\omega$ are implicit in its construction.]

Claim:  Every morphism $\bff:\bX \ra \bY$ can be written as a composite 
$\bff_\omega \circx \bi_\omega$, where 
$\bi_\omega:\bX \ra \bX_\omega$ has the LLP w.r.t. levelwise acyclic fibrations and $\bff_\omega$ is both a weak equivalence and a fibration.

[Run the small object argument once again, taking 
$S_0 = \{\bS^{-q} \wedge \bS_+^{n-1} \ra \bS^{-q} \wedge \bD_+^n$ $(n \geq 0, q\geq 0)\}$.]

Combining the claims gives MC-5 and the nontrivial half of MC-4 can be established in the usual way.]

[Note: All objects are fibrant and every cofibration is a spectral cofibration.]\\

\begingroup%%----------------------------------->>
\fontsize{9pt}{11pt}\selectfont
True or false:  The model category structure on \bSPEC is proper.\\
\endgroup %%------------------------------------<<

\bHSPEC is the homotopy category of \bSPEC (cf. p. \pageref{16.14} ff.).  
In this situation, $I\bX = \bX \wedge I_+$ is a cylinder object when \bX is cofibrant while 
$P\bX = \texttt{HOM}(I_+,\bX)$ serves as a path object.  
And: It can be assumed that the ``cofibrant replacement'' $\sL\bX$ is functorial in \bX, so
$\sL:\bSPEC \ra \bSPEC_{\bc}$.

%%----------------------------------------------------------------------------------------------10
[Note: \ Recall too that the inclusion $\bHSPEC_{\bc} \ra \bHSPEC$ is an equivalence of categories 
(cf. $\S 12$, Proposition 13).]

Remark: Suppose that \bX is cofibrant $-$then for any \bY, 
$[\bX,\bY]_0 \approx [\bX,\bY]$ (cf. p. \pageref{16.15}) (all objects are fibrant), thus if $\bY \ra \bZ$ is a weak equivalence, then $[\bX,\bY]_0 \approx [\bX,\bZ]_0$.

Example: Let $(K,k_0)$ be a pointed CW complex $-$then $\bQ_p^\infty K$ is cofibrant.\\

\begingroup%%----------------------------------->>
\fontsize{9pt}{11pt}\selectfont
\textbf{\small FACT} \  
Let $\bff:\bX \ra \bY$ be a morphism of spectra $-$then \bff is a weak equivalence iff $\forall \ n$, 
$\pi_n(\bff):\pi_n(\bX) \ra \pi_n(\bY)$ is an isomorphism.\\
\endgroup %%------------------------------------<<

\begingroup%%----------------------------------->>
\fontsize{9pt}{11pt}\selectfont
\label{16.20}
\textbf{\small LEMMA} \  
$\bHSPEC_\bc$ has coproducts and weak pushouts.
\vspi
[Note: \ The wedge $\ds\bigvee\limits_i \bX_i$ is the coproduct of the $\bX_i$ in $\bHSPEC_\bc$.  
Proof: $\ds\bigvee\limits_i \bX_i$ is cofibrant and for any cofibrant \bY, 
$\ds[\bigvee\limits_i \bX_i, \bY] \approx$ 
$[\ds\bigvee\limits_i \bX_i,\bY]_0 \approx$ 
$\pi_0(\HOM(\ds\bigvee\limits_i \bX_i,\bY)) \approx$ 
$\pi_0(\ds\prod\limits_i \HOM(\bX_i,\bY)) \approx$ 
$\ds\prod\limits_i  \pi_0(\HOM(\bX_i,\bY)) \approx$ 
$\ds\prod\limits_i [\bX_i,\bY]_0 \approx$ 
$\ds\prod\limits_i [\bX_i,\bY]$.]\\
\endgroup %%------------------------------------<<

\index{Theorem: Brown Representability Theorem (Spectra)}
\index{Brown Representability Theorem (Spectra)}
\begingroup%%----------------------------------->>
\fontsize{9pt}{11pt}\selectfont
\textbf{\small BROWN REPRESENTABILITY THEOREM} \quad 
A cofunctor $F:\bHSPEC_{\bc} \ra \bSET$ is representable iff it converts coproducts into products and weak pushouts into weak pullbacks.
\vspi
[In the notation of p. \pageref{16.16}, let $\sU = \{\bS^n:n \in \Z\}$.  If $\bff:\bX \ra \bY$ is a morphism such that 
$\forall \ n$, the arrow 
$[\bS^n,\bX] \ra [\bS^n,\bY]$ is bijective, then \bff is a weak equivalence (cf. supra), thus is a homotopy equivalence 
(cf. $\S 12$, Proposition 10).  Therefore $\sU_1$ holds.  As for $\sU_2$, given an object $(\bX,\bff)$ in 
$\bFIL(\bHSPEC_\bc)$, $\telsub(\bX,\bff)$ is a weak colimit and $\forall \ n$, the arrow 
$\colim[\bS^n,\bX_k] \ra [\bS^n,\telsub(\bX,\bff)]$ is bijective (cf. p. \pageref{16.17}).]\\
\endgroup %%------------------------------------<<

\begingroup%%----------------------------------->>
\fontsize{9pt}{11pt}\selectfont
\textbf{\small EXAMPLE} \  
$\bHSPEC_\bc$ has products.  For if $\{\bX_i\}$ is a set of cofibrant spectra, then the cofunctor 
$\bY \ra \ds\prod\limits_i[\bY,\bX_i]$ satisfies the hypotheses of the Brown representability theorem.\\
\endgroup %%------------------------------------<<

\label{16.47}
\begin{proposition} \ %09
Suppose that $\bA \ra \bY$ is a cofibration and $\bX \ra \bB$ is a fibration $-$then the arrow 
$\HOM(\bY,\bX) \ra \HOM(\bA,\bX) \times_{\HOM(\bA,\bB)} \HOM(\bY,\bB)$ is a Serre fibration which is a weak homotopy equivalence if $\bA \ra \bY$ or $\bX \ra \bB$ is acyclic.\\
\end{proposition}

\begingroup%%----------------------------------->>
\fontsize{9pt}{11pt}\selectfont
Proposition 9 implies (and is implied by) the following equivalent statements (cf. $\S 13$, Propositions 31 and 32).\\
\endgroup %%------------------------------------<<

\begingroup%%----------------------------------->>
\fontsize{9pt}{11pt}\selectfont
\textbf{\small FACT} \  
If $\bA \ra \bY$ is a cofibration in \bSPEC and if $L \ra K$ is a cofibration in $\dcg_*$, then the arrow 
$\bA \wedge K \underset{\bA \wedge L}{\sqcup} \bY \wedge L \ra \bY \wedge K$ is a cofibration in \bSPEC which is acyclic if 
$\bA \ra \bY$ or $L \ra K$ is acyclic.\\
\endgroup %%------------------------------------<<

\begingroup%%----------------------------------->>
\fontsize{9pt}{11pt}\selectfont
\textbf{\small FACT} \  
If $L \ra K$ is a cofibration in $\dcg_*$ and if $\bX \ra \bB$ is a fibration in \bSPEC, then the arrow 
$\texttt{HOM}(K,\bX) \ra \texttt{HOM}(L,\bX) \times_{\texttt{HOM}(L,\bB)} \texttt{HOM}(K,\bB)$ is a fibration in \bSPEC which is acyclic if $L \ra K$ or $\bX \ra \bB$ is acyclic.\\
\endgroup %%------------------------------------<<

%%----------------------------------------------------------------------------------------------11
The 
\un{shift suspension}
\index{shift suspension} 
is the functor $\Lambda:\bSPEC \ra \bSPEC$ defined by 
$(\Lambda \bX)_q = X_{q+1}$ $(q \geq 0)$ and the  
\un{shift desuspension}
\index{shift desuspension}  
is the functor $\Lambda^{-1}:\bSPEC \ra \bSPEC$ defined by 
$
(\Lambda^{-1} \bX)_q =
\begin{cases}
\ X_{q-1} \hspace{0.5cm} (q > 0)\\
\ \Omega X_0 \hspace{0.6cm} (q = 0)
\end{cases}
\hspace{-.2cm}.
$
\\[0.25cm]

\begin{proposition} \ 
The pair $(\Lambda,\Lambda^{-1})$ is an adjoint equivalence of categories.\\ 
\end{proposition}

\begingroup%%----------------------------------->>
\fontsize{9pt}{11pt}\selectfont
\textbf{\small EXAMPLE} \  
$\Lambda^q$ is a left adjoint for $\Lambda^{-q}$ and, by Proposition 10, $\Lambda^{-q}$ is a left adjoint for $\Lambda^{q}$.  
On the other hand, $\bQ^\infty$ is a left adjoint for $\bU^\infty$.  Therefore 
$\Lambda^{-q} \circx \bQ^\infty$ is a left adjoint for $\bU^\infty \circx \Lambda^q$.  But 
$\bU^\infty \circx \Lambda^q = \bU_q^\infty$, thus $\forall \ q \geq 0$, 
$\Lambda^{-q} \circx \bQ^\infty \approx \bQ_q^\infty$.\\
\endgroup %%------------------------------------<<

Remarks: 
(1) $\Lambda$ preserves weak equivalences, so 
$Q \circx \Lambda:\bSPEC \ra \bHSPEC$ sends weak equivalences to isomorphisms and there is a commutative triangle 
\begin{tikzcd}%[sep=small]
{\bSPEC} \ar{d}[swap]{Q} \ar{r}{Q \circx \Lambda} &{\bHSPEC}\\
{\bHSPEC} \ar{ru}[swap]{\bL \Lambda \qquad ,}
\end{tikzcd}
$\bL \Lambda$ the total left derived functor for $\Lambda$; 
(2) $\Lambda^{-1}$ preserves weak equivalences, so 
$Q \circx \Lambda^{-1}:\bSPEC \ra \bHSPEC$  sends weak equivalences to isomorphisms and there is a commutative triangle 
$
\begin{tikzcd}%[sep=small]
{\bSPEC} \ar{d}[swap]{Q} \ar{r}{Q \circx \Lambda^{-1}} &{\bHSPEC}\\
{\bHSPEC} \ar{ru}[swap]{\bR \Lambda^{-1}}
\end{tikzcd}
,
$
$\bR \Lambda^{-1}$ the total right derived functor for $\Lambda^{-1}$.\\
\vspace{0.25cm}

\begin{proposition} \  %11
The pair $(\bL\Lambda,\bR\Lambda^{-1})$ is an adjoint equivalence of categories.
\end{proposition}

[$\Lambda^{-1}$ preserves fibrations and acyclic fibrations (the data is levelwise).  Therefore $\Lambda$ preserves cofibrations and the TDF theorem implies that $(\bL\Lambda,\bR\Lambda^{-1})$ is an adjoint pair.  
Consider now the bijection of adjunction 
$\Xi_{\bX,\bY}:\Mor(\Lambda\bX,\bY) \ra \Mor(\bX,\Lambda^{-1}\bY)$, so \ 
$\Xi_{\bX,\bY}\bff$ is the composition 
%$\bX \ra$
\begin{tikzcd}[sep=small]
{\bX} \ar{r}
&{\Lambda^{-1}\Lambda\bX}  \ar{rr}{\Lambda^{-1}\bff} &&{\Lambda^{-1}\bY.}
\end{tikzcd}
Since the arrow $\bX \ra \Lambda^{-1}\Lambda\bX$ is an isomorphism, $\Xi_{\bX,\bY}\bff$ is a weak equivalence  
iff $\Lambda^{-1} \bff$ is a weak equivalence, 
i.e., iff $\bff$ is a weak equivalence.  
Therefore the pair $(\bL\Lambda,\bR\Lambda^{-1})$ is an adjoint equivalence of categories 
(cf. p. \pageref{16.18}).]\\

$\Lambda^{-1}$ is naturally isomorphic to $\ov{\Omega}$.  Here $(\ov{\Omega}\bX)_q = \Omega X_q$, the arrow of structure $\Omega X_q \ra \Omega \Omega X_{q+1}$ being $\Omega \sigma_q$.  Therefore the difference between 
$\ov{\Omega}$ and $\Omega$ is the twist $\Tee$ 
(cf. p. \pageref{16.19}).  Define a pseudo natural weak equivalence 
$\Xi_{\bX}:\Omega\bX \ra \ov{\Omega}\bX$ by letting $\Xi_{\bX,q}:\Omega X_q \ra \Omega X_q$ 
be the identity for even $q$ and the negative of the identity for odd $q$ (i.e., coordinate reversal).\\

\textbf{\small LEMMA} \  
Let \bC be a category and let $F,G:\bC \ra \bPRESPEC$ be functors.  Suppose given a pseudo natural weak equivalence 
$\Xi:F \ra G$ $-$then in the notation of the conversion principle, there are natural transformations 
$sFX \overset{sr}{\lla}$ 
$sMFX \overset{sm\Xi}{\lra}$ 
$sMGX \overset{sr}{\ra}$ 
$sGX$.

%%----------------------------------------------------------------------------------------------12
[Note: \ $sM\Xi$ is a weak equivalence.  Moreover, the $sr$ are weak equivalences if \mF, \mG factor through 
\bSEPPRESPEC.]\\

Application: $\forall$ \bX in \bSPEC, $\Omega\bX$ is naturally weakly equivalent to $\ov{\Omega}\bX$ or still, is naturally weakly equivalent to $\Lambda^{-1}\bX$.

Example: In \bHSPEC, $\bS^{-n} \approx \Omega^n\bS^0$ $(n \geq 0)$.\\

\begin{proposition} \ 
The total left derived functor $\bL\Sigma$ for $\Sigma$ exists and the total right derived functor $\bR\Omega$ for $\Omega$ exists.  And: $(\bL\Sigma,\bR\Omega)$ is an adjoint pair.
\end{proposition}

[$\Sigma$ preserves cofibrations and $\Omega$ preserves fibrations.  Now quote the TDF theorem.]

[Note: \ Since $\Omega$, $\ov{\Omega}$  preserve weak equivalences, there are commutative triangles 
\begin{tikzcd}%[sep=small]
{\bSPEC} \ar{d}[swap]{Q} \ar{r}{Q \circx \Omega} &{\bHSPEC}\\
{\bHSPEC} \ar{ru}[swap]{\bR \Omega}
\end{tikzcd}
,
\begin{tikzcd}%[sep=small]
{\bSPEC} \ar{d}[swap]{Q} \ar{r}{Q \circx \ov{\Omega}} &{\bHSPEC}\\
{\bHSPEC} \ar{ru}[swap]{\bR\ov{\Omega}}
\end{tikzcd}
and, by the above, natural isomorphisms, 
$\bR\Omega \ra \bR \ov{\Omega}$, 
$\bR \ov{\Omega} \ra \bR \Lambda^{-1}$.]\\

\begingroup%%----------------------------------->>
\fontsize{9pt}{11pt}\selectfont
$\Sigma$ preserves weak equivalences between cofibrant objects.  So, unraveling the definitions, one finds that 
$\bL\Sigma (=L(Q \circx \Sigma))$ ``is'' $\bL(\Sigma \circx \iota \circx \sL)$ 
$(\bL(\Sigma \circx \iota \circx \sL) \circx Q =$ 
$Q \circx \Sigma \circx \iota \circx \sL )$, 
$\iota:\bSPEC_\bc \ra \bSPEC$ the inclusion.  In particular: $\forall \ \bX$, 
$\bL \Sigma \bX = \Sigma \sL \bX$.\\
\endgroup %%------------------------------------<<

\begin{proposition} \ 
The pair $(\bL\Sigma,\bR\Omega)$ is an adjoint equivalence of categories.
\end{proposition}

[According to Proposition 11, the arrows of adjunction 
$\id \overset{\mu}{\lra} \bR\Lambda^{-1}\circx \bL\Lambda$, 
$\bL\Lambda \circx \bR\Lambda^{-1} \overset{\nu}{\lra} \id$ 
are natural isomorphisms and the claim is that the arrows of adjuction
$\id \overset{\mu}{\lra} \bR\Omega \circx \bL \Sigma$, 
$\bL \Sigma \circx \bR\Omega \overset{\nu}{\lra} \id$ 
are natural isomorphisms.  Thus fix a natural isomorphism $\bR\Omega \ra \bR\Lambda^{-1}$ $-$then there exists a unique natural isomorphism $\bL\Lambda \ra \bL\Sigma$ characterized by the commutativity of 
\begin{tikzcd}%[sep=small]
{[\bL\Sigma\bX,\bY]} \ar{d} \ar{r} &{[\bX,\bR\Omega\bY]} \ar{d}\\
{[\bL\Lambda\bX,\bY]}  \ar{r} &{[\bX,\bR\Lambda^{-1}\bY]}
\end{tikzcd}
$\forall$ \bX, \bY.  It remains only to note that the diagrams 
\begin{tikzcd}%[sep=small]
{\id} \ar{d}[swap]{\mu} \ar{r}{\mu} &{\bR\Lambda^{-1} \circx \bL\Lambda}\ar{d}\\
{\bR\Omega \circx \bL \Sigma} \ar{r} &{\bR\Lambda^{-1} \circx \bL\Sigma}
\end{tikzcd}
,
\begin{tikzcd}%[sep=small]
{\bL\Lambda \circx \bR\Omega} \ar{d} \ar{r} &{\bL\Lambda \circx \bR\Lambda^{-1}}\ar{d}{\nu}\\
{\bL \Sigma \circx \bR\Omega} \ar{r}[swap]{\nu} &{\id}
\end{tikzcd}
of natural transformations commute.]\\

Application: \bHSPEC is an additive category and $\bL\Sigma$ is an additive functor.

[Note: \ \bHSPEC has coproducts and products (since $\bHSPEC_{\bc}$ does 
(cf. p. \pageref{16.20}).  
Standard categorical generalities then imply that the arrow 
$\bX \vee \bY \ra \bX \times \bY$ is an isomorphism for all \bX, \bY in \bHSPEC (cf. p. \pageref{16.21}).]\\

%%----------------------------------------------------------------------------------------------13
\label{17.3}
Notation: 
Write
$
\begin{cases}
\ \Sigma\\
\ \Omega
\end{cases}
$
in place of 
$
\begin{cases}
\ \bL \Sigma\\
\ \bR \Omega
\end{cases}
$
and 
$
\begin{cases}
\ \Lambda\\
\ \Lambda^{-1}
\end{cases}
$
in place of 
$
\begin{cases}
\ \bL \Lambda\\
\ \bR \Lambda^{-1}
\end{cases}
.
$
\\
\vspace{0.25cm}

\begin{proposition} \ 
\bHSPEC is a triangulated category satisfying the octahedral axiom.
\end{proposition}

[Working in $\bHSPEC_\bc$ , stipulate that a triangle 
$\bX^\prime \overset{\bu^\prime}{\ra}$ 
$\bY^\prime \overset{\bv^\prime}{\ra}$ 
$\bZ^\prime \overset{\bw^\prime}{\ra}$ 
$\Sigma\bX^\prime$
 is exact if it is isomorphic to a triangle
$\bX \overset{\bff}{\ra}$ 
$\bY \overset{\bj}{\ra}$ 
$\bC_\bff \overset{\bpi}{\ra}$ 
$\Sigma\bX$
for some \bff ($\bC_{\bff} =$ mapping cone of \bff) (obvious definition).  
Since TR$_1$-TR$_5$ are immediate, it will be enough to deal just with the octahedral axiom.  
Suppose given exact triangles 
$\bX \overset{\bu}{\ra}$ 
$\bY \ra$ 
$\bZ^\prime \ra$ 
$\Sigma\bX$,
$\bY \overset{\bv}{\ra}$ 
$\bZ \ra$ 
$\bX^\prime \ra$ 
$\Sigma\bY$,
\begin{tikzcd}%[sep=small]
{\bX} \ar{r}{\bv \circx \bu} &{\bZ}
\end{tikzcd}
$\ra \bY^\prime$ 
$\ra \Sigma \bX$, where without loss of generality, 
$\bZ^\prime = \bC_{\bu}$, 
$\bX^\prime = \bC_{\bv}$, 
$\bY^\prime = \bC_{\bv \circx \bu}$.  
Starting at the prespectrum level, define a pointed continuous function 
$f_n:C_{u_n} \ra C_{v_n \circx u_n}$ by letting $f_n$ be the identity on $\Gamma X_n$ and $v_n$ on $Y_n$ and define a pointed continuous function $g_n:C_{v_n \circx u_n} \ra C_{v_n}$ by letting 
$g_n$ be $\Gamma u_n$ on $\Gamma X_n$ and the identity on 
$Z_n$ $-$then the $f_n$ and the $g_n$ combine to give morphisms of prespectra, 
so applying $s$, $\exists$ morphisms 
$\bff:\bZ^\prime \ra \bY^\prime$ and 
$\bg:\bY^\prime \ra \bX^\prime$ of spectra.  By construction, the composite 
$\bZ^\prime  \overset{\bff}{\lra}$ 
$\bY^\prime \ra $ 
$\Sigma \bX$ 
is the arrow 
$\bZ^\prime \ra \Sigma \bX$ and the composite 
$\bZ \ra $
$\bY^\prime \overset{\bg}{\lra}$
$\bX^\prime$ is the arrow 
$\bZ \ra \bX^\prime$.  
Letting 
$\bh:\bX^\prime \ra \Sigma \bZ^\prime$ be the composite 
$\bX \ra$ 
$\Sigma \bY \ra$ 
$\Sigma \bZ^\prime$, 
one sees that all the commutativity required of the octahedral axiom is present, 
thus the final task is to establish that the triangle 
$\bZ^\prime \overset{\bff}{\ra}$ 
$\bY^\prime \overset{\bg}{\ra}$ 
$\bX^\prime \overset{\bh}{\ra}$ 
$\Sigma\bZ^\prime$ is exact.  
But there is a canonical commutative diagram
\hspace{0.25cm}
\begin{tikzcd}%[sep=small]
{\bZ^\prime} \arrow[d,shift right=0.5,dash] \arrow[d,shift right=-0.5,dash]  \ar{r}{\bff}
&{\ \bY^\prime\ } \arrow[d,shift right=0.5,dash] \arrow[d,shift right=-0.5,dash]  \ar{r}{\bg}
&{\bX^\prime} \ar{d}{\bphi} \ar{r}{\bh}
&{\Sigma\bZ^\prime} \arrow[d,shift right=0.5,dash] \arrow[d,shift right=-0.5,dash] \\
{\bZ^\prime} \ar{r}[swap]{\bff}
&{\ \bY^\prime\ } \ar{r}[swap]{\bj}
&{\bC_{\bff}} \ar{r}[swap]{\bpi}
&{\Sigma\bZ^\prime}
\end{tikzcd}
\hspace{-.2cm}.
And: $\bphi$ is a homotopy equivalence.]\\

Application: An exact triangle \ 
$\bX \overset{\bu}{\ra}$ 
$\bY \overset{\bv}{\ra}$ 
$\bZ \overset{\bw}{\ra}$ 
$\Sigma\bX$ \ in \bHSPEC gives rise to a long exact sequence in homotopy 
$\cdots \ra$
$\pi_{n+1}(\bZ) \ra$ 
$\pi_{n}(\bX) \ra$ 
$\pi_{n}(\bY) \ra$ 
$\pi_{n}(\bZ) \ra$ 
$\pi_{n-1}(\bX) \ra \cdots$.\\

\begingroup%%----------------------------------->>
\fontsize{9pt}{11pt}\selectfont
\label{16.54}
\textbf{\small EXAMPLE} \  If $\bff:\bX \ra \bY$, $\bg:\bY \ra \bZ$ are morphisms in \bHSPEC, then there is an exact triangle 
$\bC_{\bff} \ra$ 
$\bC_{\bg \circx \bff} \ra$ 
$\bC_{\bg} \ra$ 
$\Sigma\bC_{\bff}$.\\
\endgroup %%------------------------------------<<

Remark: \bHSPEC is compactly generated (take $\sU = \{\bS^n:n \in \Z\}$) and admits Adams representability (by Neeman's countability criterion).\\

\begingroup%%----------------------------------->>
\fontsize{9pt}{11pt}\selectfont
\textbf{\small EXAMPLE} \  
The homotopy groups of a compact spectrum are finitely generated.
\vspi
[The thick subcategory of \bHSPEC whose objects are those \bX such that $\pi_q(\bX)$ is finitely generated $\forall \ q$ 
contains the $\bS^n$.]\\
\endgroup %%------------------------------------<<

It is also true that \bHSPEC is a closed category (indeed, a CTC) but the proof requires some preliminary work which is best carried out in a more general context.\\

%%----------------------------------------------------------------------------------------------14
\begingroup%%----------------------------------->>
\fontsize{9pt}{11pt}\selectfont
The main difficulty lies in equipping \bHSPEC with the structure of a closed category (cf. p. \pageref{16.22}).  Granted this, the fact that \bHSPEC is a CTC can be seen as follows.
\vspi
Recall that if $f:X \ra Y$ is a map in the pointed category, then there is a homotopy commutative diagram
\endgroup

\begingroup
\fontsize{9pt}{11pt}\selectfont
\[
\begin{tikzcd}%[sep=small]
&&{\Sigma\Omega X} \ar{d} \ar{r} 
&{\Sigma\Omega Y} \ar{d} \ar{r} 
&{\Sigma E_f} \ar{d} \ar{r} 
&{\Sigma X}\arrow[d,shift right=0.5,dash] \arrow[d,shift right=-0.5,dash]\\
{\Omega Y} \arrow[d,shift right=0.5,dash] \arrow[d,shift right=-0.5,dash] \ar{r}
&{E_f} \ar{d} \ar{r} 
&{X} \ar{d} \ar{r} 
&{Y} \ar{d} \ar{r} 
&{C_f} \ar{r} 
&{\Sigma X}\\
{\Omega Y} \ar{r} 
&{\Omega C_f} \ar{r}
&{\Omega\Sigma X} \ar{r}
&{\Omega\Sigma Y}
\end{tikzcd}
,
\]
\endgroup

\begingroup
\fontsize{9pt}{11pt}\selectfont
a formalism which also holds in the category of prespectra or spectra.  Of course, when viewed in \bHSPEC, the arrows 
$\bE_{\bff} \ra \Omega \bC_{\bff}$, 
$\Sigma\bE_{\bff} \ra \bC_{\bff}$ 
are isomorphisms (cf. Proposition 13).  
Turning to the axioms for a CTC, the only one that is potentially troublesome is CTC$_4$.  
In order not to obscure the issue, we shall proceed informally, omitting all mention of $\sL$ and the underlying total derived functors.  
Thus given \ 
$\bX \overset{\bff}{\ra}$ 
$\bY \overset{\bj}{\ra}$ 
$\bC_{\bff} \overset{\bpi}{\ra} \Sigma \bX$, one has to show that $\forall \ \bZ$, the triangle 
\begin{tikzcd}%[sep=small]
{\Omega\hom(\bX,\bZ)}  \ar{rr}{-(\bpi^*\circx \eta_{\bX,\bZ})} 
&&{{\hom(\bC_{\bff},\bZ)}} \ar{r}{\bj^*} &{}
\end{tikzcd}
\begin{tikzcd}%[sep=small]
{\hom(\bY,\bZ)} \ar{rrr}{\nu^{-1}_{\hom(\bX,\bZ)}\circx\bff^*}
&&&{\Sigma\Omega\hom(\bX,\bZ)}
\end{tikzcd}
is exact.  Consider the commutative diagram
\[
\begin{tikzcd}[sep=large]
{\ \Omega \bY\ } \arrow[d,shift right=0.5,dash] \arrow[d,shift right=-0.5,dash] \ar{r}{\Omega \bj}
&{\ \Omega \bC_{\bff}\ } \arrow[d,shift right=0.5,dash] \arrow[d,shift right=-0.5,dash] \ar{rr}{\Omega \bpi}
&&{\Omega\Sigma \bX} \ar{d}{\mu_{\bX}^{-1}} \ar{rrr}{\nu_{\bY}^{-1} \circx \mu_{\bY}^{-1}\circx \Omega\Sigma \bff}
&&&{\ \Sigma\Omega \bY\ } \arrow[d,shift right=0.5,dash] \arrow[d,shift right=-0.5,dash]\\
{\ \Omega \bY\ } \ar{r}[swap]{\Omega \bj}
&{\ \Omega \bC_{\bff}\ }  \ar{rr}[swap]{\mu_{\bX}^{-1}\circx \Omega\bpi}
&&{\bX} \ar{rrr}[swap]{\nu_{\bY}^{-1}\circx \bff}
&&&{\ \Sigma\Omega \bY\ }
\end{tikzcd}
.
\]
\endgroup

\begingroup
\fontsize{9pt}{11pt}\selectfont
Since the triangle on the bottom is exact (cf. p. \pageref{16.23}), so is the triangle on the top.  
But then, on the basis of the commutative diagram 
\[
\begin{tikzcd}[sep=large]
{\ \Omega \bY\ } \arrow[d,shift right=0.5,dash] \arrow[d,shift right=-0.5,dash] \ar{r}
&{\ \bE_{\bff}\ } \ar{d} \ar{r}
&{\bX} \ar{d}{\mu_{\bX}} \ar{r}{\bff}
&{\bY} \ar{rr}{\nu_{\bY}^{-1}} \ar{d}{\mu_{\bY}}
&&{\Sigma\Omega \bY} \arrow[d,shift right=0.5,dash] \arrow[d,shift right=-0.5,dash]\\
{\ \Omega \bY\ } \ar{r}
&{\ \Omega \bC_{\bff}\ }  \ar{r}
&{\Omega\Sigma \bX} \ar{r}[swap]{\Omega\Sigma \bff}
&{\Omega\Sigma \bY} \ar{rr}[swap]{\nu_{\bY}^{-1}\circx \mu_{\bY}^{-1}}
&&{\Sigma\Omega \bY}
\end{tikzcd}
,
\]
\endgroup

\begingroup
\fontsize{9pt}{11pt}\selectfont
the triangle 
$\Omega \bY \ra \bE_{\bff} \ra$ 
\begin{tikzcd}%[sep=small]
{\bX} \ar{rr}{\nu_{\bY}^{-1}\circx \bff} &&{\Sigma\Omega Y}
\end{tikzcd}
is exact.  In particular: The triangle 
$\Omega\hom(\bX,\bZ) \ra \bE_{\bff^*} \ra$
\begin{tikzcd}%[sep=small]
{\hom(\bY,\bZ)} \ar{rrr}{\nu^{-1}_{\hom(\bX,\bZ)}\circx\bff^*} &&&{\Sigma\Omega\hom(\bX,\bZ)}
\end{tikzcd}
is exact.  However, there is an isomorphism $\bE_{\bff^*} \ra \hom(\bC_{\bff},\bZ)$ and a commutative diagram
\[
\begin{tikzcd}[sep=large]
{\Omega\hom(\bX,\bZ)} \arrow[d,shift right=0.5,dash] \arrow[d,shift right=-0.5,dash] \ar{rr}
&&{\bE_{\bff^*}} \ar{d} \ar{r}
&{\hom(\bY,\bZ)}  
\arrow[d,shift right=0.5,dash] \arrow[d,shift right=-0.5,dash] \ar{rr}{\nu^{-1}_{\hom(\bX,\bZ)}\circx\bff^*}
&&{\ \Sigma\Omega\hom(\bX,\bZ) \ } \arrow[d,shift right=0.5,dash] \arrow[d,shift right=-0.5,dash]\\
{\Omega\hom(\bX,\bZ)} \ar{rr}[swap]{-(\bpi^*\circx \eta_{\bX,\bZ})}
&&{\hom(\bC_{\bff},\bZ)} \ar{r}[swap]{\bj^*}
&{\hom(\bY,\bZ)} \ar{rr}[swap]{\nu^{-1}_{\hom(\bX,\bZ)}\circx\bff^*}
&&{\ \Sigma\Omega\hom(\bX,\bZ) \ }
\end{tikzcd}
,
\]
%%----------------------------------------------------------------------------------------------15
hence the triangle on the bottom is exact, this being the case of the triangle on the top.\\
\endgroup %%------------------------------------<<

\begingroup%%----------------------------------->>
\fontsize{9pt}{11pt}\selectfont
\textbf{\small LEMMA} \  \bHSPEC is a compactly generated CTC.
\vspi
[In general, \bX dualizable $\implies$
$
\begin{cases}
\ \Sigma \bX\\
\ \Omega \bX
\end{cases}
$
dualizable (cf. $\S 15$, Proposition 35).  But trivially the unit $\bS^0$ is dualizable, thus $\forall n > 0$, 
$\bS^n \approx \Sigma^n \bS^0$ $\&$ $\bS^{-n} \approx \Omega^n\bS^0$ are dualizable, i.e., all the elements of 
$\sU = \{\bS^n:n \in \Z\}$ are dualizable.]
\vspi
[Note: \ Observe too that $\forall \ n$, $D\bS^n \approx \bS^{-n}$.]\\
\endgroup %%------------------------------------<<

\label{17.1}
\begingroup%%----------------------------------->>
\fontsize{9pt}{11pt}\selectfont
Remark: \ \bHSPEC is a unital compactly generated CTC (since $\bS^0$ is compact).  \ Accordingly, 
du\bHSPEC = cpt\bHSPEC (cf. p. \pageref{16.24}), the thick subcategory generated by the $\bS^n$ (theorem of Neeman-Ravenel).
\vspi
[Note: \ It is clear that \bHSPEC is actually monogenic.]\\
\endgroup %%------------------------------------<<

\label{17.40}
\begingroup%%----------------------------------->>
\fontsize{9pt}{11pt}\selectfont
\textbf{\small EXAMPLE} 
The compact objects in \bHSPEC are those objects which are isomorphic to a $\bQ_q^\infty K$, where \mK is a pointed finite CW complex.\\
\endgroup %%------------------------------------<<

Notation: Given a real finite dimensional inner product space \mV, let $\bS^V$ denote its one point compactification (base point at $\infty$) and for any \mX in $\dcg_*$, put $\Sigma^V X =  X \#_k \bS^V$, $\Omega^V X = X^{\bS^V}$.

[Note: \  If \mV and \mW are two real finite dimensional inner product spaces such that $V \subset W$, write $W - V$ for the orthogonal complement of \mV in \mW $-$then $\forall$ \mX, 
$\Sigma^{W - V}\Sigma^V X \approx \Sigma^W X$ and 
$\Omega^V\Omega^{W - V} X \approx \Omega^W X$.]\\

A 
\un{universe}
\index{universe} 
is a real inner product space $\sU$ with $\dim \sU = \omega$ equipped with the finite topology.  \bUN is the category whose objects are the universes and whose morphisms are the linear isometries.  
An 
\un{indexing set}
\index{indexing set (in a universe)} 
in a universe $\sU$ is a set $\sA$ of finite dimensional subspaces of $\sU$ such that each finite dimensional subspace $V$ of $\sU$ is contained in some $U \in \sA$.  
The 
\un{standard indexing set}
\index{standard indexing set (in a universe)} 
is the set of all finite dimensional subspaces of $\sU$.  
Example: Take $\sU = \R^\infty$ $-$then $\{\R^q:q \geq 0\}$ is an indexing set in $\R^\infty$.

Let $\sA$ be an indexing set in a universe $\sU$ $-$then a 
\un{$(\sU,\sA)$-prespectrum}
\index{prespectrum! $(\sU,\sA)$-prespectrum} 
\bX is a collection of pointed \dsp compactly generated spaces 
$X_U$ $(U \in \sA)$ and a collection of pointed continuous functions 
\begin{tikzcd}%[sep=small]
{X_V} \ar{r}{\sigma_{V,W}} &{\Omega^{W-V}X_W}
\end{tikzcd}
($V,W \in \sA$ $\&$ $V \subset W$) such that 
\begin{tikzcd}%[sep=small]
{X_V} \ar{r}{\sigma_{V,V}} &{X_V}
\end{tikzcd}
is the identity and for $U \subset V \subset W$ in $\sA$, the diagram
\begin{tikzcd}%[sep=small]
{X_U} \ar{d}[swap]{\sigma_{U,W}} \ar{r}{\sigma_{U,V}} &{\Omega^{V- U}X_V} \ar{d}{\Omega^{V- U}\sigma_{V,W}}\\
{\Omega^{W - U}X_W} \arrow[r,shift right=0.5,dash] \arrow[r,shift right=-0.5,dash]  
&{\Omega^{V- U}\Omega^{W - V}X_W}
\end{tikzcd}
commutes.  
$\bPRESPEC_{\sU,\sA}$ is the category whose objects are the $(\sU,\sA)$-prespectra and whose morphisms 
$\bff:\bX \ra \bY$ are collections of pointed continuous functions $f_U:X_U \ra Y_U$
%%----------------------------------------------------------------------------------------------16
such that the diagram 
\begin{tikzcd}%[sep=small]
{X_V} \ar{d} \ar{r}{f_V} &{Y_V} \ar{d}\\
{\Omega^{W - V}X_W} \ar{r}[swap]{\Omega^{W - V}f_W} &{\Omega^{W - V}Y_W}
\end{tikzcd}
commutes for $V \subset W$ in $\sA$.  A $(\sU,\sA)$-prespectrum \bX is a 
\un{$(\sU,\sA)$-spectrum}
\index{spectrum! $(\sU,\sA)$-spectrum} 
if the $\sigma_{V,W}$ are homeomorphisms.  
$\bSPEC_{\sU,\sA}$ is the full subcategory of $\bPRESPEC_{\sU,\sA}$ with object class the $(\sU,\sA)$-spectra.  
Example: Take $\sU = \R^\infty$, $\sA = \{\R^q:q \geq 0\}$ $-$then 
$\bPRESPEC_{\sU,\sA} =$ $\bPRESPEC$, 
$\bSPEC_{\sU,\sA} =$ $\bSPEC$.

[Note: When $\sA$ is the standard indexing set, write $\bPRESPEC_{\sU}$, $\bSPEC_{\sU}$, is place of 
$\bPRESPEC_{\sU,\sA}$, $\bSPEC_{\sU,\sA}$.]

What has been said earlier can now be said again.  
Thus introduce the notion of a separated $(\sU,\sA)$-prespectrum by requiring that the 
$\sigma_{V,W}:X_V \ra \Omega^{W-V}X_W$ be \bCG embeddings.  
This done, repeat the proof of Proposition 1 to see that 
$\bSEPPRESPEC_{\sU,\sA}$ is a reflective subcategory of $\bPRESPEC_{\sU,\sA}$ with reflector $E^\infty$.  
Next, as in Proposition 2, $\bSPEC_{\sU,\sA}$ is a reflective subcategory of $\bSEPPRESPEC_{\sU,\sA}$ (the reflector sends 
\bX to e\bX, where 
$(e\bX)_V = \underset{W \supset V}{\colimx} \Omega^{W-V}X_W$).  
Conclusion: $\bSPEC_{\sU,\sA}$ is a reflective subcategory of $\bPRESPEC_{\sU,\sA}$ (cf. Proposition 3), hence is complete and cocomplete.

[Note: \ The composite 
$\bPRESPEC_{\sU,\sA} \overset{E^\infty}{\lra}$
$\bSEPPRESPEC_{\sU,\sA} \overset{e}{\lra}$
$\bSPEC_{\sU,\sA}$
is the 
\un{spectrification functor}:
\index{spectrification functor} 
$\bX \ra s\bX$ $(s = e \circx E^\infty)$.]\\

\begingroup%%----------------------------------->>
\fontsize{9pt}{11pt}\selectfont
\textbf{\small EXAMPLE} \  
Fix $U \in \sA$.  Given an $X$ in $\dcg_*$, let $\bQ_U^\infty X$ be the spectrification of the prespectrum 
$V \ra$ 
$
\begin{cases}
\ \Sigma^{V - U} X \hspace{0.35cm} (V \supset U)\\
\ *  \hspace{1.35cm} (V \not\supset U)
\end{cases}
\hspace{-.2cm},
$
where 
$\Sigma^{V - U} X \ra$ 
$\Omega^{W-V}\Sigma^{W-V}\Sigma^{V - U}X \approx$
$\Omega^{W-V}\Sigma^{W-U}X$ ($V,W \in \sA$, $\&$ $U \subset V \subset W$) 
(otherwise, the arrow is the inclusion of the base point).  
Viewed as a functor from $\dcg_*$ to $\bSPEC_{\sU,\sA}$, $\bQ_U^\infty$ 
is a left adjoint for the $U^\text{th}$ space functor 
$\bU_U^\infty:\bSPEC_{\sU,\sA} \ra$ $\dcg_*$ that sends $\bX = \{X_u\}$ to $X_U$.\\
\endgroup %%------------------------------------<<

\begingroup%%----------------------------------->>
\fontsize{9pt}{11pt}\selectfont
\textbf{\small FACT} \  
If \bX is a $(\sU,\sA)$-spectrum and if $\dim V_1 = \dim V_2$ $(V_1, V_2 \in \sA)$, then 
$X_{V_1} \approx X_{V_2}$.
\vspi
[Embed $V_1$ and $V_2$ in a common finite dimensional $W \in \sA$ and observe that 
$X_{V_1} \approx$ 
$\Omega^{W-V_1} X_W \approx$
$\Omega^{W-V_2} X_W \approx$
$X_{V_2}$.]\\
\endgroup %%------------------------------------<<

Notation: Given \bX, \bY in $\bPRESPEC_{\sU,\sA}$, write $\HOM(\bX,\bY)$ for $\Mor(\bX,\bY)$ topologized via the equalizer diagram
$\Mor(\bX,\bY) \ra $ 
$\prod\limits_{V \in \sA} Y_V^{X_V} \rightrightarrows$ 
$\prod\limits_{\substack{V,W \in \sA \\ V \subset W}} (\Omega^{W-V}Y_W)^{X_V}$.

So, just as before, spectrification is a continuous functor (cf. Proposition 4) and there are analogs of Propositions 5 and 6 
($\bbox$ $(\wedge)$ and \texttt{HOM} being defined in the obvious way).

Remark: $\bPRESPEC_{\sU,\sA}$ and $\bSPEC_{\sU,\sA}$ are $\bV$-categories, where 
$\bV = \dcg_*$.  Accordingly, to say that $s$ is continuous simply means that $s$ is a \bV-functor.

%%----------------------------------------------------------------------------------------------17
[Note: \ The interpretation of $\bbox$ $(\wedge)$ and \texttt{HOM} is that 
$\bPRESPEC_{\sU,\sA}$ and $\bSPEC_{\sU,\sA}$ admit a closed $\dcg_*$ action (the topological parallel of closed simplicial action).]\\

\textbf{\small LEMMA} \  
Let $\sA$ and $\sB$ be indexing sets in a universe $\sU$ with $\sA \subset \sB$ $-$then the arrow of restriction 
$i^*:\bPRESPEC_{\sU,\sB} \ra \bPRESPEC_{\sU,\sA}$ has a left adjoint $i_*$ and a right adjoint $i_!$.

[For \bX in $\bPRESPEC_{\sU,\sA}$ and $W$ an element of $\sB$, 
$(i_*\bX)_W$ is the coequalizer of 
$\coprod\limits_{\substack{V\pp \subset V^\prime\\ V^\prime \subset W}} 
\Sigma^{W - V^\prime}\Sigma^{V^\prime - V\pp}X_{V\pp} \rightrightarrows$ 
$\coprod\limits_{\substack{V\\ V \subset W}} 
\Sigma^{W - V}X_V$ 
and $(i_!\bX)_W$ is the equalizer of 
$\prod\limits_{\substack{V\\ W \subset V}} \Omega^{V-W} X_V$ $\rightrightarrows$ 
$\prod\limits_{\substack{V^\prime \subset V\pp\\ W \subset V^\prime}}$
$\Omega^{V^\prime - W}\Omega^{V\pp - V^\prime} X_{V\pp}$ 
$(V, V^\prime, V\pp \in \sA)$.]\\
\vspace{0.25cm}

\begingroup%%----------------------------------->>
\fontsize{9pt}{11pt}\selectfont
The formulas figuring in the lemma can be understood in terms of ``enriched'' Kan extensions.  Thus let $\bI_{\sA}$ be the category whose objects are the elements of $\sA$, with 
$\Mor(V^\prime,V\pp) =$
$
\begin{cases}
\ \bS^{V\pp - V^\prime} \hspace{0.35cm} (V\pp \supset V^\prime)\\
\ * \hspace{1.2cm} \ (V\pp \not\supset V^\prime)
\end{cases}
$
(composition comes from the identification  
$\bS^{V - U} \#_k \bS^{W - V} \approx \bS^{W - U}$) 
$-$then $\bI_{\sA}$ is a small \bV-category and $\bPRESPEC_{\sU,\sA}$ ``is'' 
$\bV[\bI_{\sA},\dcg_*]$ (cf. p. \pageref{16.25}) ($\bV = \dcg_*$).  So, if $\sA \subset \sB$ and 
$i:\bI_{\sA} \ra \bI_{\sB}$ is the inclusion, 
$i_* = \lan$ $\&$ $i_! = \ran$, i.e., 
$i_* \bX = \lan \bX$ (the left Kan extension of \bX along $i$) $\&$ 
$i_! \bX = \ran \bX$ (the right Kan extension of \bX along $i$).\\
\endgroup %%------------------------------------<<

\begin{proposition} \ %15
Let $\sA$ and $\sB$ be indexing sets in a universe $\sU$ with $\sA \subset \sB$ $-$then the arrow of restriction 
$i^*:\bSPEC_{\sU,\sB} \ra \bSPEC_{\sU,\sA}$ is an equivalence of categories.
\end{proposition}

[The functor $s \circx i_*$ is a left adjoint for $i^*$ and the arrows of adjunction 
$\id \overset{\mu}{\ra} i^* \circx (s \circx i_*)$, 
$(s \circx i_*) \circx i^* \overset{\nu}{\ra} \id$ are natural isomorphisms.]\\

Application: Let $\sU$ be a universe $-$then $\forall$ indexing set $\sA$ in $\sU$, $\bSPEC_{\sU,\sA}$ is equivalent to 
$\bSPEC_{\sU}$.\\

\begingroup%%----------------------------------->>
\fontsize{9pt}{11pt}\selectfont
\index{Thom Spectra (example)}
\textbf{\small EXAMPLE \ (\un{Thom Spectra})} \  
If $\sU$ is a universe and if $\bG_n(\sU)$ is the grassmannian of $n$-dimensional subspaces of $\sU$, then $\bG_n(\sU)$ is topologized as the colimit of the $\bG_n(U)$ ($U \subset \sU$ $\&$ $\dim U < \omega$), so every compact subspace of 
$\bG_n(\sU)$ is contained in some $\bG_n(U)$.  Let $K$ be a compact Hausdorff space and suppose that 
$f:K \ra \bG_n(\sU)$ is a continuous function.  Write 
$\sA_f$ for the set of $U: f(K) \subset \bG_n(U)$ $-$then $\sA_f$ is an indexing set in $\sU$.  Give $U \in \sA_f$, call 
$K^{U-f}$ the Thom space of the vector bundle defined by the pullback square 
\begin{tikzcd}[sep=large]
{f^*\gamma_n^\perp} \ar{d} \ar{r} &{\gamma_n^\perp} \ar{d}\\
{K} \ar{r}[swap]{f} &{\bG_n(U)}
\end{tikzcd}
($\gamma_n$ the canonical $n$-plane bundle over $\bG_n(U)$).  
The assignment 
$U \ra K^{U-f}$ defines an object in $\bPRESPEC_{\sU,\sA_f}$.  
Pass to its spectrification $\bSPEC_{\sU,\sA_f}$, thence by the above to an object in $\bSPEC_{\sU}$, say $K^{-f}$.  
In general, an arbitrary $X$ in $\dcg$ can be represented
%%----------------------------------------------------------------------------------------------18
as the colimit of its compact subspaces $K$: $X \approx \colimx K$.  
Accordingly, for $f:X \ra \bG_n(\sU)$ a continuous function, put 
$X^{-f} = \colimx K^{\restr{-f}{K}}$, the 
\un{Thom spectrum}
\index{Thom spectrum} 
of the virtual vector bundle $-f$.  
Example: An $n$-dimensional $U$ determines a map 
$* \overset{U}{\ra} \bG_n(\sU)$ and $*^{-U} \approx \bS^{-U}$.\\
\endgroup %%------------------------------------<<

The $U^\text{th}$ space functor 
$\bU_U^\infty:\bSPEC_{\sU} \ra \dcg_*$ is represented by $\bS^{-U}$, where 
$\bS^{-U} =$ $\bQ_U^\infty \bS^0$ (cf. Proposition 7).  Equipping $\dcg_*$ with its singular structure, if 
$\bff:\bX \ra \bY$ is a morphism of $\sU$-spectra, then \bff is a levelwise fibration iff \bff has the RLP w.r.t. the spectral cofibrations 
$\bS^{-U} \wedge [0,1]_+^n \ra$ 
$\bS^{-U} \wedge I[0,1]_+^n$ and \bff is a levelwise acyclic fibration iff \bff has the RLP w.r.t. the spectral cofibrations 
$\bS^{-U} \wedge \bS_+^{n-1} \ra$ 
$\bS^{-U} \wedge \bD_+^n$ ($n \geq 0$, $U \subset \sU$ $\&$ $\dim U < \omega$) (cf. p. \pageref{16.26}).  
Using this, it follows that $\bSPEC_{\sU}$ is a model category if weak equivalences and fibrations are levelwise, the cofibrations being those morphisms which have the LLP w.r.t. the levelwise acyclic fibrations (cf. Proposition 8) (bear in mind that a spectral cofibration is  necessarily a levelwise closed embedding (cf. p. \pageref{16.27})).  Proposition 9 and its variants go through without change.

[Note: \ $\bHSPEC_{\sU}$
\index{$\bHSPEC_{\sU}$} 
is the homotopy category of $\bSPEC_{\sU}$ 
(cf. p. \pageref{16.28} ff).]

Remark: The functor $\bU_U^\infty$ preserves fibrations and acyclic fibrations, thus the TDF theorem implies that 
$\bL \bQ_U^\infty$ and $\bR \bU_U^\infty$ exist and 
$(\bL \bQ_U^\infty,\bR \bU_U^\infty)$ is an adjoint pair (the requisite assumptions are validated by the generalities on p. \pageref{16.29} ff.).\\

\begingroup%%----------------------------------->>
\fontsize{9pt}{11pt}\selectfont
\label{16.40}
\textbf{\small EXAMPLE} \  
Take $\sU = \R^\infty$ $-$then $i^*:\bSPEC_{\sU} \ra \bSPEC$ preserves fibrations and acyclic fibrations, so the hypotheses of the TDF theorem are satisfied (cf. p. \pageref{16.30} ff.).  Therefore 
$\bL i_*$ and $\bR i^*$ exist and $(\bL i_*,\bR i^*)$ is an adjoint pair.  Dissecting the bijection of adjunction 
$\bXi_{\bX,\bY}:\Mor(i_*\bX,\bY) \ra$ 
$\Mor(\bX,i^*\bY)$, it follows that $\bXi_{\bX,\bY} \bff$ is a weak equivalence iff \bff is a weak equivalence, thus the pair 
$(\bL i_*,\bR i^*)$ is an adjoint equivalence of categories 
(cf. p. \pageref{16.31}).\\
\endgroup %%------------------------------------<<

\label{16.33} %dmc mnft
Let $\sU$, $\sU^\prime$ be universes, $f:\sU \ra \sU^\prime$ a linear isometry $-$then there is a functor 
$f^*:\bPRESPEC_{\sU^\prime} \ra \bPRESPEC_{\sU}$ which assigns to each 
$\bX^\prime$ in $\bPRESPEC_{\sU^\prime}$ the $\sU$-prespectrum 
$f^*\bX^\prime$ specified by 
$(f^*\bX^\prime)_U =$ $\bX^\prime_{f(U)}$, where
$(f^*\bX^\prime)_V \ra$ $\Omega^{W - V}(f^*\bX^\prime)_W$ is the composite 
$\bX^\prime_{f(V)} \ra$ 
$\Omega^{f(W) - f(V)}\bX^\prime_{f(W)} \ra$ 
$\Omega^{W - V}\bX^\prime_{f(W)}$.  
It has a left adjoint 
$f_*:\bPRESPEC_{\sU} \ra \bPRESPEC_{\sU^\prime}$, viz. 
$(f_*\bX)_{U^\prime} =$ 
$\Sigma^{U^\prime - f(U)}X_U$ $(U = f^{-1}(U^\prime))$, where 
$(f_*\bX)_{V^\prime} \ra$ 
$\Omega^{W^\prime - V^\prime}(f_*\bX)_{W^\prime}$ is the composite 
$\Sigma^{V^\prime - f(V)}X_V \ra$
$\Omega^{W^\prime - V^\prime}\Sigma^{W^\prime - V^\prime}\Sigma^{V^\prime - f(V)}X_V \ra$
$\Omega^{W^\prime - V^\prime}\Sigma^{W^\prime - f(W)}\Sigma^{f(W)-f(V)}X_V \ra$ 
$\Omega^{W^\prime - V^\prime}\Sigma^{W^\prime - f(W)}\Sigma^{W - V}X_V \ra$ 
$\Omega^{W^\prime - V^\prime}\Sigma^{W^\prime - f(W)}X_W$ 
($V = f^{-1}(V^\prime)$, $W = f^{-1}(W^\prime)$).  
Since $f^*$ sends $\sU^\prime$-spectra to $\sU$-spectra, there is an induced functor 
$f^*:\bSPEC_{\sU^\prime} \ra \bSPEC_{\sU}$ and a left adjoint for it is $s \circx f_*$, denoted still by $f_*$.\\

\begingroup%%----------------------------------->>
\fontsize{9pt}{11pt}\selectfont
Let $\bI_{\sU}$, $\bI_{\sU^\prime}$ be the small \bV-categories associated with the standard indexing sets in 
$\sU$, $\sU^\prime$  $-$then the 
%%----------------------------------------------------------------------------------------------19
linear isometry $f:\sU \ra \sU^\prime$ determines a continuous functor 
$F_f:\bI_{\sU} \ra \bI_{\sU^\prime}$.  Viewing $\bPRESPEC_{\sU}$  as 
$\bV[\bI_{\sU},\dcg_*]$ and $\bPRESPEC_{\sU^\prime}$  as 
$\bV[\bI_{\sU^\prime},\dcg_*]$, $f^*$ becomes the precomposition with $F_f$ and $f_* =$ lan.\\
\endgroup %%------------------------------------<<

\begingroup%%----------------------------------->>
\fontsize{9pt}{11pt}\selectfont
\textbf{\small EXAMPLE} \  
$f_*(\bX \wedge K) \approx (f_*(\bX) \wedge K$ and 
$f_*(\bQ_U^\infty X) \approx \bQ_{f(U)}^\infty X$.\\
\endgroup %%------------------------------------<<

\begingroup%%----------------------------------->>
\fontsize{9pt}{11pt}\selectfont
\textbf{\small FACT} \  
Let $\sU$, $\sU^\prime$ be universes, $f:\sU \ra \sU^\prime$ a linear isometric isomorphism $-$then the pair 
$(f_*,f^*)$ is an adjoint isomorphism of categories.
\vspi
[Note: \ Here, of course, it is a question of spectra, not prespectra.]\\
\endgroup %%------------------------------------<<

Let $\sU$, $\sU^\prime$ be universes $-$then a
\un{$(\sU,\sU^\prime)$-spectrum}
\index{spectrum! $(\sU,\sU^\prime)$-spectrum} 
$\bX^\prime$ is a collection of 
$\sU^\prime$-spectra $\bX_U^\prime$ indexed by the finite dimensional subspaces \mU of 
$\sU$ and a collection of isomorphisms 
\begin{tikzcd}%[sep=small]
\hspace{-.15cm}{\Sigma^{W - V}\bX_W^\prime} \ar{r}{\rho_{W,V}} &{\bX_V^\prime}
\end{tikzcd}
$(V \subset W)$ such that 
\begin{tikzcd}%[sep=small]
{\bX_V^\prime} \ar{r}{\rho_{V,V}} &{\bX_V^\prime}
\end{tikzcd}
is the identity and for $U \subset$ $V \subset$ $W$, \ the diagram \ \ 
\begin{tikzcd}%[sep=small]
{\Sigma^{V - U}\Sigma^{W - V}\bX_W^\prime}  \ar{d}[swap]{\Sigma^{V - U}\rho_{W,V}} \ar[equals]{r} 
&{\Sigma^{W - U}\bX_W^\prime} \ar{d}{\rho_{W,U}}\\
{\Sigma^{V - U} \bX_V^\prime} \ar{r}[swap]{\rho_{V,U}}
&{\bX_U^\prime}
\end{tikzcd}
\ 
commutes.  \ $\bSPEC(\sU^\prime,\sU)$
\index{$\bSPEC(\sU^\prime,\sU)$} 
is the c ategory whose objects are the $(\sU^\prime,\sU)$-spectra and whose morphisims 
$\bff:\bX^\prime \ra \bY^\prime$ are collections of morphisms of $\sU^\prime$ -spectra 
$\bff_U^\prime:\bX_U^\prime \ra \bY_U^\prime$ such that the diagram
\begin{tikzcd}%[sep=small]
{\Sigma^{W - V}\bX_W^\prime} \ar{d} \ar{rr}{\Sigma^{W - V}\bff_W^\prime} &&{\Sigma^{W - V}\bY_W^\prime}\ar{d}\\
{\bX_V^\prime} \ar{rr}[swap]{\bff_V^\prime} &&{\bY_V^\prime}
\end{tikzcd}
commutes for $V \subset W$.

[Note: \ It makes sense to suspend a $\sU^\prime$-spectrum by a finite dimensional subspace of $\sU$ (this being an instance of smashing with an object in $\dcg_*$).]\\

\begingroup%%----------------------------------->>
\fontsize{9pt}{11pt}\selectfont
\textbf{\small EXAMPLE} \  
Let $\sU$, $\sU^\prime$ be universes, $\bff:\sU \ra \sU^\prime$ a linear isometry.  
Given an $X$ in $\dcg_*$, let 
$\bQ_f^\prime X$ be the object in $\bSPEC(\sU^\prime,\sU)$ defined by 
$(\bQ_f^\prime X)_U =$ $\bQ_{f(U)}^\infty X$, where 
$\Sigma^{W - V}(\bQ_f^\prime X)_W \ra$ $(\bQ_f^\prime X)_V$ is the identification 
$\Sigma^{W - V}\bQ_{f(W)}^\infty X \approx$ 
$\Sigma^{f(W) - f(V)}\bQ_{f(W)}^\infty X \approx$ 
$\bQ_{f(V)}^\infty X$.\\
\endgroup %%------------------------------------<<

Notation: Given $\bX^\prime$, $\bY^\prime$ in $\bSPEC(\sU^\prime,\sU)$, write 
$\HOM(\bX^\prime,\bY^\prime)$ for $\Mor(\bX^\prime,\bY^\prime)$ topologized via the equalizer diagram
$\Mor(\bX^\prime,\bY^\prime) \ra $ 
$\prod\limits_V \HOM(\bX_V^\prime,\bY_V^\prime) \rightrightarrows$ 
$\prod\limits_{\substack{V,W\\ V \subset W}}\HOM(\Sigma^{W-V}\bX_W^\prime,$ $\bY_V^\prime)$.

\indent\indent $(\wedge)$ \ Fix a \mK in $\dcg_*$.  Given an $\bX^\prime$ in $\bSPEC(\sU^\prime,\sU)$, let 
$\bX^\prime \wedge K$ be the $(\sU^\prime,\sU)$-spectrum $U \ra \bX_U^\prime \wedge K$, where 
$\Sigma^{W-V}(\bX_W^\prime \wedge K) \approx$ 
$(\bX_W^\prime \wedge K) \wedge \bS^{W-V} \approx$ 
$(\bX_W^\prime \wedge \bS^{W-V}) \wedge K \approx$ 
$\bX_V^\prime \wedge K$.\\

\begin{proposition} \ %16
For $\bX^\prime$, $\bY^\prime$ in $\bSPEC(\sU^\prime,\sU)$ and \mK in $\dcg_*$, 
there is a natural homeomorphism 
$\HOM(\bX^\prime \wedge K,\bY^\prime) \approx \HOM(\bX^\prime,\bY^\prime)^K$.\\
\end{proposition}

%%----------------------------------------------------------------------------------------------20
\indent\indent $(\texttt{HOM})$ \quad 
Fix an 
$\bX^\prime$ in $\bSPEC(\sU^\prime,\sU)$.  
Given a 
$\bY^\prime$ in $\bSPEC_{\sU^\prime}$, let 
$\texttt{HOM}(\bX^\prime,\bY^\prime)$ be the $\sU$-spectrum 
$U \ra \HOM(\bX_U^\prime,\bY^\prime)$, where 
$\HOM(\bX_V^\prime,\bY^\prime) \approx$ 
$\HOM(\Sigma^{W-V}\bX_W^\prime,\bY^\prime)  \approx$ 
$\HOM(\bX_W^\prime,\Omega^{W-V}\bY^\prime) \approx$ 
$\Omega^{W-V}\HOM(\bX_W^\prime,\bY^\prime)$.

Observation: $\forall \ \bX$ in $\bSPEC_{\sU}$, 
$\Mor(\bX,\texttt{HOM}(\bX^\prime,\bY^\prime)) \approx$ 
$\lim\Mor(X_U,\HOM(\bX_U^\prime,\bY^\prime)) \approx$ 
$\lim\Mor(\bX_U^\prime \wedge X_U,\bY^\prime) \approx$ 
$\Mor(\colimx \bX_U^\prime \wedge X_U,\bY^\prime)$, 
the colimit being taken over the arrows 
$\bX_V^\prime \wedge X_V \approx$ 
$\Sigma^{W-V}\bX_W^\prime \wedge X_V \approx$ 
$\bX_W^\prime \wedge \Sigma^{W-V} X_V \ra$ 
$\bX_W^\prime \wedge X_W$.

Definition: $\bX^\prime \wedge \bX$ is the $\sU^\prime$-spectrum $\colimx \bX_U^\prime \wedge X_U$.\\

\begin{proposition} \ 
For \bX in $\bSPEC_{\sU}$, $\bY^\prime$ in $\bSPEC_{\sU^\prime}$, and 
$\bX^\prime$ in $\bSPEC(\sU^\prime,\sU)$, there is a natural homeomorphism 
$\HOM(\bX^\prime \wedge \bX,\bY^\prime) \approx \HOM(\bX,\texttt{HOM}(\bX^\prime,\bY^\prime))$.\\
\end{proposition}

\begingroup%%----------------------------------->>
\fontsize{9pt}{11pt}\selectfont
\textbf{\small EXAMPLE} \  
(1) $\bX^\prime \wedge \bQ_U^\infty X \approx \bX_U^\prime \wedge X$; 
(2) $(\bX^\prime \wedge \bX) \wedge K \approx 
\bX^\prime \wedge (\bX \wedge K) \approx 
(\bX^\prime \wedge K) \wedge \bX$.\\
\endgroup %%------------------------------------<<

Notation: \  Given a vector bundle $\xi:E \ra B$, $T(\xi)$ is its Thom space.

[Note: \ If $S^\xi$ is the sphere bundle obtained from $\xi$ by fiberwise one point compactification, then 
$T(\xi) = S^\xi/S_\infty$, where $S_\infty$ is the section at infinity.  
Example: If $\un{V}$ is the trivial vector bundle 
$B \times V \ra B$, then $T(\xi \oplus \un{V}) \approx \Sigma^V T(\xi)$.]

Let $\sU$, $\sU^\prime$ be universes.  
Fix an object 
$A \overset{\alpha}{\lra} \sI(\sU,\sU^\prime)$ in $\dcg/\sI(\sU,\sU^\prime)$ 
($\sI(\sU,\sU^\prime)$ topologized as on p. \pageref{16.32}).  
Given finite dimensional 
$U \subset \sU$, 
$U^\prime \subset \sU^\prime$, define $A_{U,U^\prime}$ by the pullback square 
\begin{tikzcd}%[sep=small]
{A_{U,U^\prime}} \ar{d} \ar{rr} &&{\sI(U,U^\prime)} \ar{d}\\
{A} \ar{r}{\alpha} &{\sI(\sU,\sU^\prime)} \ar{r} &{\sI(U,\sU^\prime)}
\end{tikzcd}
(so $A_{U,U^\prime} =$ $\{a \in A:\alpha(a) U \subset U^\prime\}$, which can be empty).  
Write $\xi(\alpha)_{U,U^\prime}$ for the vector bundle over $A_{U,U^\prime}$ with total space 
$\{(a,u^\prime) \in A_{U,U^\prime} \times U^\prime$: $u^\prime \perp \alpha(a)U\}$ and let 
$T \alpha_{U,U^\prime}$ be the associated Thom space (if  $A_{U,U^\prime}$ is empty, then the Thom space is a singleton).  
For each $U$, the assignment 
$U^\prime \ra$ $T \alpha_{U,U^\prime}$ specifies a $\sU^\prime$-prespectrum, call it $\bT^\prime \alpha_U$ 
(the arrow 
$T \alpha_{U,V^\prime} \ra$ $\Omega^{W^\prime - V^\prime}T \alpha_{U,W^\prime}$ 
is the adjoint of the arrow 
$\Sigma^{W^\prime - V^\prime} T \alpha_{U,V^\prime} \ra$ $T \alpha_{U,W^\prime}$ 
induced by the morphism 
$\xi(\alpha)_{U,V^\prime} \oplus (\un{W^\prime} - \un{V^\prime}) \ra$ 
$\xi(\alpha)_{U,W^\prime}$ of vector bundles).  Let 
$\bM^\prime \alpha_U$ be the spectrification of $\bT^\prime \alpha_U$ $-$then there are morphisms 
$\Sigma^{W - V} \bM^\prime \alpha_W \ra$ 
$\bM^\prime \alpha_V$ of $\sU^\prime$-spectra arising from the morphisms 
$\xi(\alpha)_{W,V^\prime} \oplus (\un{W} - \un{V}) \ra$ 
$\xi(\alpha)_{V,V^\prime}$ of vector bundles.\\

\begin{proposition} \ %18
The morphisms $\Sigma^{W-V} \bM^\prime \alpha_W \ra \bM^\prime \alpha_V$ are isomorphisms, thus the collection 
$\bM_\alpha^\prime = \{\bM^\prime\alpha_U\}$ is an object in $\bSPEC(\sU^\prime,\sU)$.
\end{proposition}

[Since all the constructions are natural in $\dcg/\sI(\sU,\sU^\prime)$ and commute with colimits, one can assume that $A$ is compact.  
But then, for $V \subset W$, $\exists$ $V^\prime$ : 
$A_{V,V^\prime} = A_{W,V^\prime} = A$, hence 
$\Sigma^{W- V} \bM^\prime \alpha_W \approx$ 
$\Sigma^{W- V} \bQ_{V^\prime}^\infty T \alpha_{W,V^\prime} \approx$ 
$\bQ_{V^\prime}^\infty \Sigma^{W- V} T \alpha_{W,V^\prime} \approx$ 
$\bQ_{V^\prime}^\infty T \alpha_{V,V^\prime} \approx$ 
$\bM^\prime \alpha_V$.]\\

Example: There is an isomorphism $\bM^\prime\alpha_{\{0\}} \approx \bQ_{\{0\}}^\infty A_+$ natural in $\alpha$.

%%----------------------------------------------------------------------------------------------21
[In fact, $\xi(\alpha)_{\{0\},U^\prime}$ is the trivial vector bundle $A \times U^\prime \ra A$.]\\

\begingroup%%----------------------------------->>
\fontsize{9pt}{11pt}\selectfont
\textbf{\small EXAMPLE} \  
Suppose that $\alpha$ is the constant map at $f \in \sI(\sU,\sU^\prime)$ $-$then 
$\bM^\prime \alpha \approx \bQ_f^\prime A_+$.\\
\endgroup %%------------------------------------<<

Let $\sU$, $\sU^\prime$ be universes.  Fix an 
$A \overset{\alpha}{\lra} \sI(\sU,\sU^\prime)$ in $\dcg/\sI(\sU,\sU^\prime)$.

\indent\indent $(\ltimes)$ \quad Given an \bX in $\bSPEC_{\bU}$, let 
$\alpha \ltimes \bX$ be the $\sU^\prime$-spectrum $\bM^\prime\alpha \wedge \bX$.

\indent\indent $(\sHOM)$ \quad Given an $\bY^\prime$ in $\bSPEC_{\bU^\prime}$, let 
$\sHOM[\alpha,\bY^\prime)$ be the $\sU$-spectrum $\texttt{HOM}(\bM^\prime\alpha,$ $\bY^\prime)$.

Remark: $\alpha \ltimes \bX \approx \colimx \restr{\alpha}{K} \ltimes \bX$ 
and 
$\sHOM[\alpha,\bY^\prime) \approx \lim \sHOM[\restr{\alpha}{K},\bY^\prime)$, where $K$ runs over the compact subspaces of $A$.

[Note: \ 
$\ltimes:\dcg/\sI(\sU,\sU^\prime) \times \bSPEC_{\sU} \ra \bSPEC_{\sU^\prime}$ and 
$\sHOM:(\dcg/\sI(\sU,\sU^\prime))^\OP \times \bSPEC_{\sU^\prime} \ra \bSPEC_{\sU}$
are continuous functors of their respective arguments.  Moreover, 
$\alpha \ltimes \bX$ preserves colimits in $\alpha$ and \bX, while 
$\sHOM[\alpha,\bY^\prime)$ converts colimits in $\alpha$ to limits and preserves limits in $\bY^\prime$.]\\

\begin{proposition} \ %19
For \bX in $\bSPEC_{\sU}$, $\bY^\prime$ in $\bSPEC_{\sU^\prime}$, and $\alpha$ in $\dcg/\sI(\sU,\sU^\prime)$, there is a natural homeomorphism 
$\HOM(\alpha \ltimes \bX,\bY^\prime) \approx$ $\HOM(\bX,\sH\sO\sM[\alpha,\bY^\prime))$ (cf. Proposition 17).\\
\end{proposition}

Example: \ Fix a linear isometry $f:\sU \ra \sU^\prime$, viewed as an object in $* \ra \sI(\sU,\sU^\prime)$ $-$then 
$f \ltimes \bX \approx f_*\bX$ and 
$\sHOM[f,\bY^\prime) \approx$ $f^*\bY^\prime$ (cf. p. \pageref{16.33}).

[E.g.: $\bM^\prime f_U \approx \bQ_{f(U)}^\infty \bS^0$ $\implies$ 
$\sHOM[f,\bY^\prime) \approx$ 
$\HOM(\bQ_{f(U)}^\infty \bS^0,\bY^\prime) \approx$ 
$\bY_{f(U)}^\prime$.]

Examples 
(1) $(\alpha \ltimes \bX) \wedge K \approx \alpha \ltimes (\bX \wedge K)$; 
(2) $\texttt{HOM}(K,\sHOM[\alpha,\bY^\prime)) \approx$ 
$\sHOM[\alpha,\texttt{HOM}(K,$ $\bY^\prime))$.\\

\label{16.41}
Addendum: Let $\bHOM(\bX,\bY^\prime)$ be the set of ordered pairs $(f,\bff)$, where 
$f \in \sI(\sU,\sU^\prime)$ and $\bff:\bX \ra f^*\bY^\prime$ is a morphism of $\sU$-spectra, and let 
$\epsilon:\bHOM(\bX,\bY^\prime) \ra \sI(\sU,\sU^\prime)$ be the projection 
$(f,\bff) \ra f$ $-$then 
Elmendorf\footnote[2]{\textit{J. Pure Appl. Algebra} \textbf{54} (1988), 37-94.} 
has shown that one may equip $\bHOM(\bX,\bY^\prime)$  with the strucure of a \dsp compactly generated space in such a way that $\epsilon$ is continuous (and $\epsilon^{-1}(f) \approx \HOM(\bX,f^*\bY^\prime)$ $\forall \ f$).  
Moreover, there are natural homeomorphisms 
$\HOM(\alpha \ltimes \bX,\bY^\prime) \approx$ 
$\HOM(\alpha,\epsilon) \approx$ 
$\HOM(\bX,\sHOM[\alpha,\bY^\prime)$).

[Note:  $\HOM(\alpha,\epsilon)$ is the set of all continuous functions 
\begin{tikzcd}[sep=small]
{A} \ar{rdd}[swap]{\alpha} \ar{rr} &&{\bHOM(\bX,\bY^\prime)} \ar{ldd}{\epsilon}\\
\\
&{\sI(\sU,\sU^\prime)}
\end{tikzcd}
regarded as a closed subspace of $\bHOM(\bX,\bY^\prime)^A$ (viz., the fiber of 
$\bHOM(\bX,\bY^\prime)^A \overset{\epsilon_*}{\lra} \sI(\sU,\sU^\prime)^A$ over $\alpha$).]\\

%%----------------------------------------------------------------------------------------------22
\begingroup%%----------------------------------->>
\fontsize{9pt}{11pt}\selectfont
\textbf{\small FACT} \  
Suppose given $\alpha:A \ra \sI(\sU,\sU^\prime)$.  Let $B$ be in $\dcg$ and call $\pi$ the projection 
$A \times_k B \ra A$ $-$then 
$(\alpha \circx \pi) \ltimes \bX \approx$ 
$(\alpha \ltimes \bX) \wedge B_+$ and 
$\sHOM(\alpha \circx \pi,\bY^\prime) \approx$ 
$\texttt{HOM}(B_+,\sHOM(\alpha ,\bY^\prime))$.\\
\endgroup %%------------------------------------<<

\label{16.56}
\label{16.59}
\begingroup%%----------------------------------->>
\fontsize{9pt}{11pt}\selectfont
\textbf{\small FACT} \  
Suppose given 
$\alpha:A \ra \sI(\sU,\sU^\prime)$ and 
$\beta:B \ra \sI(\sU^\prime,\sU\pp)$.  Let $\beta \times_c \alpha$ be the composite 
\begin{tikzcd}%[sep=small]
{B \times_k A} \ar{r}{\beta \times_k \alpha} &{\sI(\sU^\prime,\sU\pp) \times_k \sI(\sU,\sU^\prime)}
\end{tikzcd}
$\overset{c}{\lra} \sI(\sU,\sU\pp)$ $-$then 
$(\beta \times_c \alpha) \ltimes \bX \approx$ 
$\beta \ltimes (\alpha \ltimes \bX)$ and 
$\sHOM[\beta \times_c \alpha,\bY\pp) \approx$ 
$\sHOM[\alpha,\sHOM[\beta,\bY\pp) )$.\\
\endgroup %%------------------------------------<<

\begin{proposition} \ %20
Fix an $\alpha$ in $\dcg/\sI(\sU,\sU^\prime)$ $-$then for \bX in $\bSPEC_{\sU}$ and $\bY^\prime$ in 
$\bSPEC_{\sU^\prime}$, a morphism $\bphi:\alpha \ltimes \bX \ra \bY^\prime$ determines and is determined by morphisms 
$\bphi(a):\bX \ra \alpha(a)*\bY^\prime$ $(a \in A)$ such that the functions 
$T\alpha_{U,U^\prime} \#_k X_U \ra \Sigma^{U^\prime - \alpha(a)U} Y_{\alpha(a)U}^\prime \ra$ 
$Y_{U^\prime}^\prime$ are continuous, the first arrow being the assignment 
$(a,u^\prime) \#_k x \ra \phi(a)_U(x) \#_k u^\prime$ 
$(a \in A_{U,U^\prime}, u^\prime \in U^\prime - \alpha(a)U, x \in X_U)$.
\end{proposition}

[Write 
$\bM^\prime \alpha_U =$ 
$\colimx_{U^\prime} \bQ_{U^\prime} ^\infty T\alpha_{U,U^\prime}$ 
to get 
$\Mor(\alpha \ltimes \bX,\bY^\prime)$ \hspace{0.05cm} $\approx$  \hspace{0.05cm}
$\Mor(\colimx_U \ \bM^\prime \alpha_U \wedge X_U, \bY^\prime)$ \hspace{0.05cm} $\approx$ \hspace{0.05cm} 
$\lim_{U} \Mor(\bM^\prime \alpha_U \wedge X_U, \bY^\prime)$ $\approx$  
$\lim_{U} \Mor(\colimx_{U^\prime} \bQ_{U^\prime}^\infty (T\alpha_{U,U^\prime} \#_k X_U),\bY^\prime)$ $\approx$ 
$\lim_U\lim_{U^\prime} \Mor(\bQ_{U^\prime}^\infty (T\alpha_{U,U^\prime} \#_k X_U),\bY^\prime)$ $\approx$
$\lim_U\lim_{U^\prime} \Mor(T\alpha_{U,U^\prime} \#_k X_U,Y^\prime_{U^\prime})$.  
Take now a 
$\bphi:\alpha \ltimes \bX \ra \bY^\prime$ and let $\bphi(a)$ be the adjoint of the composite 
$\alpha(a)_* \bX \ra$ 
$\alpha \ltimes \bX \overset{\bphi}{\lra} \bY^\prime$.  
Projecting from the double limit thus gives rise to continuous functions 
$T\alpha_{U,U^\prime} \#_k X_U \ra$ $Y_{U^\prime}^\prime$ as stated.  
Conversely, a collection of morphisms 
$\bphi(a):\bX \ra$ $\alpha(a)^*\bY^\prime$ $(a \in A)$ 
satisfying the hypotheses define continuous functions compatible with the maps in the double limit, 
hence specify a morphism $\bphi:\alpha \ltimes \bX \ra \bY^\prime$.]\\

Given a universe $\sU$, $\bO(\sU)$ is its orthogonal group, so topologically, 
$\bO(\sU) =$ $\colimx \bO(U)$, where $\bO(U)$ is the orthogonal group of the ambient finite dimensional subspace \mU of $\sU$.\\

\textbf{\small LEMMA} \  
Let $\sU$ be a universe $-$then $\forall$ finite dimensional $U \subset \sU$, the arrow of restriction 
$\bO(\sU) \ra \sI(U,\sU)$ is a Serre fibration.\\

Application: $\sec_{\sI(U,\sU)}(\bO(\sU))$ is not empty.

[$\sI(U,\sU)$ is a CW complex and, being contractible (cf. p. \pageref{16.34}), the identity map 
$\sI(U,\sU) \ra \sI(U,\sU)$ admits a lifting $\sI(U,\sU) \ra \bO(\sU)$ (cf. p. \pageref{16.35}).]\\

\index{Untwisting Lemma}
\textbf{\small UNTWISTING LEMMA} \quad 
Let $\sU$, $\sU^\prime$ be universes.  Fix $U \subset \sU$, $U^\prime \subset \sU^\prime$ such that 
$U \approx U^\prime$ $-$then there is an isomorphism 
$\bM^\prime \alpha_U \approx \bQ_{U^\prime}^\infty A_+$ natural in $\alpha$.

[Choose a linear isometric isomorphism $f:U \ra U^\prime$ and a section 
$s^\prime:\sI(U^\prime,\sU^\prime) \ra$ 
$\bO(\sU^\prime)$.  Put 
$s = s^\prime \circx (f^*)^{-1}$.  Define 
$A_{[U,V^\prime]}$ by the pullback square 
%\begin{tikzcd}%[sep=small]
%{A_{[U,V^\prime]}} \ar{d} \ar{rr} &&{\bO(V^\prime)} \ar{d}\\
%%\\
%{A} \ar{r}{\restr{\alpha}{U}} &{\sI(U,\sU^\prime)} \ar{r}{s} &{\bO(\sU^\prime)} 
%\end{tikzcd}
\begin{tikzcd}[sep=small]
{A_{[U,V^\prime]}} \ar{dd} \ar{rr} &&{}\\
\\
{A} \ar{r}{\restr{\alpha}{U}} &{\sI(U,\sU^\prime)} \ar{r}{s} &{} 
\end{tikzcd}
\begin{tikzcd}[sep=small]
{\bO(V^\prime)} \ar{dd}\\
\\
{\bO(\sU^\prime)} 
\end{tikzcd}
%%----------------------------------------------------------------------------------------------23
if $U^\prime \subset V^\prime$ and let $A_{[U,V^\prime]} = \emptyset$ otherwise (thus 
$A_{[U,V^\prime]} \subset A_{U,V^\prime}$).  
Write 
$\xi(\alpha)_{[U,V^\prime]}$ for the trivial vector bundle 
$A_{[U,V^\prime]} \times (V^\prime - U^\prime)$ and, passing to Thom spaces, let $\bT^\prime \alpha_{[U]}$ be the 
$\sU^\prime$-prespectrum 
$V^\prime \ra T(\xi(\alpha)_{[U,V^\prime]}) \approx$ 
$\Sigma^{V^\prime - U^\prime}A_{[U,V^\prime]^+}$.  Call 
$\bM^\prime \alpha_{[U]}$ the spectrification of $\bT^\prime \alpha_{[U]}$ $-$then there are two claims: 
(1) $\bM^\prime \alpha_{[U]} \approx$ $\bQ_{U^\prime}^\infty A_+$; 
(2) $\bM^\prime \alpha_{[U]} \approx$ $\bM^\prime \alpha_U$.  
For the first, one can assume that $A$ is compact, in which case 
$A_{[U,V^\prime]} = A$ for $V^\prime$ large enough and the claim follows.  
Turning to the second, define a morphism 
$\xi(\alpha)_{[U,V^\prime]} \ra$ $\xi(\alpha)_{U,V^\prime}$ of vector bundles by sending 
$(a,v^\prime)$ to $(a,s(\restr{\alpha(a)}{U})(v^\prime))$.  
These morphisms lead to a morphism 
$\bT^\prime \alpha_{[U]} \ra$
$\bT^\prime \alpha_{U}$
of $\sU^\prime$-prespectra or still, to a morphism 
$\bM^\prime \alpha_{[U]} \ra$ $\bM^\prime \alpha_{U}$ of $\sU^\prime$-spectra.  
But when $A$ is compact and $A_{[U,V^\prime]} = A$, the bundle map is an isomorphism.]\\

\begin{proposition} \ %21
Let $\sU$, $\sU^\prime$ be universes.  Fix $U \subset \sU$, $U^\prime \subset \sU^\prime$ $-$then there is an isomorphism 
$\alpha \ltimes \bQ_U^\infty X \approx \bQ_{U^\prime}^\infty$ $(A_+ \#_k X)$ natural in $\alpha$ and \mX.
\end{proposition}

[For 
$\alpha \ltimes \bQ_U^\infty X = \bM^\prime \alpha \wedge \bQ_U^\infty X$ 
$\approx \bM^\prime \alpha_U \wedge X$ and, by the untwisting lemma, 
$\bM^\prime \alpha_U \wedge X$ 
$\approx \bQ_{U^\prime}^\infty A_+ \wedge X$.]\\

\begingroup%%----------------------------------->>
\fontsize{9pt}{11pt}\selectfont
\textbf{\small EXAMPLE} \  
Fix $U \subset \sU$, $U^\prime \subset \sU^\prime$ such that $U \approx U^\prime$ $-$then the functor
$\bM^\prime -_U:\dcg/\sI(\sU,\sU^\prime) \ra$ $\bSPEC_{\sU^\prime}$ 
has for a right adjoint the functor 
$\bM -_{U^\prime}:\bSPEC_{\sU^\prime} \ra \dcg/\sI(\sU,\sU^\prime)$ 
that sends $\bY^\prime$ to 
$\sI(\sU,\sU^\prime)  \times_k Y_{U^\prime}^\prime \ra \sI(\sU,\sU^\prime)$.
\vspi
[
$\Mor(\bM^\prime \alpha_U,\bY^\prime) \approx$
$\Mor(\bQ_{U^\prime}^\infty A_+,\bY^\prime) \approx$ 
$\Mor(A_+,Y_{U^\prime}^\prime) \approx$ 
$\Mor(A,Y_{U^\prime}^\prime) \approx$
$\Mor(\alpha,\bM \bY_{U^\prime}^\prime)$.]\\
\endgroup %%------------------------------------<<

\begingroup%%----------------------------------->>
\fontsize{9pt}{11pt}\selectfont
\textbf{\small FACT} \  
Suppose that $A$ is a CW complex $-$then the functor 
$\sHOM[\alpha,-)$ preserves weak equivalences.
\vspi
[Let 
$\bff^\prime:\bX^\prime \ra \bY^\prime$ be a weak equivalence of $\sU^\prime$-spectra and consider the induced morphism 
$\sHOM[\alpha, \bX^\prime) \ra$
$\sHOM[\alpha, \bY^\prime)$ of $\sU$-spectra.  
Given $U \subset \sU$, $\exists$ $U^\prime \subset \sU^\prime$ : $U \approx U^\prime$ $\implies$ 
$\sHOM[\alpha, \bX^\prime)_U \approx$  
$(X_{U^\prime}^\prime)^{A_+}$, 
$\sHOM[\alpha, \bY^\prime)_U \approx$  
$(Y_{U^\prime}^\prime)^{A_+}$ 
(cf. Proposition 21).  
Since $A_+$ is a CW complex and 
$X_{U^\prime}^\prime \ra Y_{U^\prime}^\prime$ is a weak homotopy equivalence, 
$(X_{U^\prime}^\prime)^{A_+} \ra (Y_{U^\prime}^\prime)^{A_+}$ is also a weak homotopy equivalence 
(cf. p. \pageref{16.36}).\\
\endgroup %%------------------------------------<<

Rappel: $\dcg/\sI(\sU,\sU^\prime)$ is a model category (singular structure) (cf. p. \pageref{16.37}).\\

\begin{proposition} \ %22
If $\bX \ra \bY$ is a cofibration in $\bSPEC_{\sU}$ and if 
\begin{tikzcd}[sep=small]
{A} \ar{rdd}[swap]{\alpha} \ar{rr} &&{B} \ar{ldd}{\beta}\\
\\
&{\sI(\sU,\sU^\prime)}
\end{tikzcd}
is a cofibration in $\dcg/\sI(\sU,\sU^\prime)$, then the arrow 
$\beta \ltimes \bX \underset{\alpha \ltimes \bX}{\sqcup} \alpha \ltimes \bY \ra \beta \ltimes \bY$ is a cofibration in 
$\bSPEC_{\sU^\prime}$ which is acyclic if $\bX \ra \bY$ or 
\begin{tikzcd}[sep=small]
{A} \ar{rdd}[swap]{\alpha} \ar{rr} &&{B} \ar{ldd}{\beta}\\
\\
&{\sI(\sU,\sU^\prime)}
\end{tikzcd}
is acyclic (cf. $\S 13$, Proposition 31).\\
\end{proposition}

%%----------------------------------------------------------------------------------------------24
\begin{proposition} \ %23
If
\begin{tikzcd}[sep=small]
{A} \ar{rdd}[swap]{\alpha} \ar{rr} &&{B} \ar{ldd}{\beta}\\
\\
&{\sI(\sU,\sU^\prime)}
\end{tikzcd}
is a cofibration in $\dcg/\sI(\sU,\sU^\prime)$ and if $\bY^\prime \ra \bX^\prime$ is a fibration in 
$\bSPEC_{\sU^\prime}$ then the arrow
$\sH\sO\sM[\beta,\bY^\prime) \ra$ $\sH\sO\sM[\alpha,\bY^\prime) \times_{\sH\sO\sM[\beta,\bX^\prime)}$ 
$\sH\sO\sM[\beta,\bX^\prime)$
is a fibration in $\bSPEC_{\sU}$ which is acyclic if 
\begin{tikzcd}[sep=small]
{A} \ar{rdd}[swap]{\alpha} \ar{rr} &&{B} \ar{ldd}{\beta}\\
\\
&{\sI(\sU,\sU^\prime)}
\end{tikzcd}
or 
$\bY^\prime \ra \bX^\prime$ 
is acyclic (cf. $\S 13$, Proposition 32).\\
\end{proposition}

\begingroup%%----------------------------------->>
\fontsize{9pt}{11pt}\selectfont
Propositions 22 and 23 are formally equivalent.  To establish the fibration contention in Proposition 23, use Proposition 21 and convert the lifting problem
\[
\begin{tikzcd}[sep=large]
{\bS^{-U} \wedge [0,1]_+^n} \ar{d} \ar{rr} 
&&{\sHOM[\beta,\bY^\prime)} \ar{d}\\
{\bS^{-U} \wedge I[0,1]_+^n} \ar[dashed]{rru}\ar{rr} 
&&{\sHOM[\alpha,\bY^\prime) \times_{\sHOM[\alpha,\bX^\prime)} \sHOM[\beta,\bX^\prime)}
\end{tikzcd}
\]
in $\bSPEC_{\sU}$ to the lifting problem
\[
\begin{tikzcd}[sep=large]
{[0,1]^n} \ar{d} \ar{rr} 
&&{(Y^\prime_{U^\prime})^B} \ar{d}\\
{I[0,1]^n} \ar[dashed]{rru}\ar{rr} 
&&{(Y^\prime_{U^\prime})^A \times_{(X^\prime_{U^\prime})^A} (X^\prime_{U^\prime})^B}
\end{tikzcd}
\]
in $\bDelta$-\bCG.\\
\endgroup %%------------------------------------<<

\begingroup%%----------------------------------->>
\fontsize{9pt}{11pt}\selectfont
\textbf{\small LEMMA} \  
Let \mA, \mB be cofibrant objects in $\bDelta$-\bCG and suppose that 
\begin{tikzcd}[sep=small]
{A} \ar{rdd}[swap]{\alpha} \ar{rr} &&{B} \ar{ldd}{\beta}\\
\\
&{\sI(\sU,\sU^\prime)}
\end{tikzcd}
is an acyclic cofibration in $\dcg/\sI(\sU,\sU^\prime)$.  Fix a cofibrant object \bX in $\bSPEC_{\sU}$ and consider the cummutative diagram
\begin{tikzcd}[sep=large]
{\bX} \ar{d} \ar{r} &{\sHOM[\alpha,\alpha \ltimes \bX)} \ar{d}\\
{\sHOM[\beta,\beta \ltimes \bX)} \ar{r} &{\sHOM[\alpha,\beta \ltimes \bX)}
\end{tikzcd}
$-$then the arrow of adjunction 
$\bX \ra \sHOM[\alpha,\alpha \ltimes \bX)$
is a weak equivalence iff the arrow of adjunction  
$\bX \ra \sHOM[\beta,\beta \ltimes \bX)$
is a weak equivalence.
\vspi
[Since the arrow $\beta \ltimes \bX \ra *$ is a fibration, it follows from Proposition 23 that 
$\sHOM[\beta,\beta \ltimes \bX) \ra \sHOM[\alpha,\beta \ltimes \bX)$ is an acyclic fibration.  
On the other hand, since the arrow $* \ra \bX$ is a cofibration, it follows from Proposition 22 that the arrow 
$\alpha \ltimes \bX \ra$ $\beta \ltimes \bX$ is an ayclic cofibration.  But from the assumptions, 
$\alpha \ltimes \bX$ and 
$\beta \ltimes \bX$ are cofibrant, thus as fibrancy is automatic, tha arrow 
$\alpha \ltimes \bX \ra$ $\beta \ltimes \bX$ is a homotopy equivalence (cf. $\S 12$, Proposition 10).  Therefore 
$\sHOM[\alpha,\alpha \ltimes \bX) \ra$ 
$\sHOM[\alpha,\beta \ltimes \bX)$ is a homotopy equivalence .]\\
\endgroup %%------------------------------------<<

%%----------------------------------------------------------------------------------------------25

\begingroup%%----------------------------------->>
\fontsize{9pt}{11pt}\selectfont
\textbf{\small EXAMPLE} \  
Let $\sU$, $\sU^\prime$ be universes, $f:\sU \ra \sU^\prime$ a linear isometry $-$then 
$f^*:\bSPEC_{\sU^\prime} \ra \bSPEC_{\sU}$ preserves fibrations and acyclic fibrations so the hypotheses of the TDF theorem are satisfied (cf. p. \pageref{16.38}ff).  
Therefore $\bL f_*$ and $\bR f^*$ exist and $(\bL f_*,\bR f^*)$ is an adjoint pair.  
Claim: $\forall$ cofibrant \bX in $\bSPEC_{\sU}$, the arrow of adjunction $\bX \ra f^*f_*\bX$ is a weak equivalence.  To see this, choose a linear isometric isomorphism $\phi \in \sI(\sU,\sU^\prime)$ and a path $H:[0,1] \ra \sI(\sU,\sU^\prime)$ such that $H \circx i_0 = \phi$ and $H \circx i_1 = f$.  Because
\begin{tikzcd}[sep=small]
{*} \ar{rdd}[swap]{\phi} \ar{rr}{i_0} &&{[0,1]} \ar{ldd}{H}\\
\\
&{\sI(\sU,\sU^\prime)}
\end{tikzcd}
is an acyclic cofibration in $\dcg/\sI(\sU,\sU^\prime)$ with $*, [0,1]$ cofibrant and because the arrow of adjunction 
$\bX \ra \phi^* \phi_* \bX$ is an isomorphism, the lemma implies that the arrow of adjunction 
$\bX \ra \sHOM[H,H \ltimes \bX)$ is a weak equivalence.  
Another application of the lemma to
\begin{tikzcd}[sep=small]
{*} \ar{rddd}[swap]{f} \ar{rr}{i_1} &&{[0,1]} \ar{lddd}{H}\\
\\
\\
&{\sI(\sU,\sU^\prime)}
\end{tikzcd}
then leads to the conclusion that the arrow of adjunction $\bX \ra f^*f_*\bX$ is indeed a weak equivalence.  
Since 
$\bX^\prime \ra \bY^\prime$ is a weak equivalence iff 
$f^*\bX^\prime \ra f^*\bY^\prime$ is a weak equivalence, the pair $(\bL f_*,\bR f^*)$ is an adjoint equivalence of categories (see the note on p. \pageref{16.39} to the TDF theorem).  
Example: $\forall$ universe $\sU$, 
$\bHSPEC_{\sU}$ ``is'' $\bHSPEC$.  
Proof: $\bHSPEC_{\sU}$ ``is'' $\bHSPEC_{\R^\infty}$ which ``is'' $\bHSPEC$ (cf. p. \pageref{16.40}).
\vspi
\label{16.50a}
[Note: \ The functors $\bL f_*:\bHSPEC_{\sU} \ra \bHSPEC_{\sU^\prime}$ obtained from the 
$f \in \sI(\sU,\sU^\prime)$ are naturally isomorphic.  
Thus let $g \in \sI(\sU,\sU^\prime)$ and choose a path 
$H:[0,1] \ra \sI(\sU,\sU^\prime)$ such that $H \circx i_0 = f$ 
and  
$H \circx i_1 = g$ $-$then for cofibrant \bX, there are natural homotopy equivalences 
$f_*\bX \ra H \ltimes \bX \la g_* \bX$ and the natural isomorphism 
$\bL f_* \approx \bL g_*$ is independent of the choice of \mH.  
In effect, if $\sigma, \tau:[0,1] \ra \sI(\sU,\sU^\prime)$ are paths in $\sI(\sU,\sU^\prime)$  such that
$
\begin{cases}
\ \sigma(0) = f\\
\ \sigma(1) = g
\end{cases}
$
,
$
\begin{cases}
\ \tau(0) = f\\
\ \tau(1) = g
\end{cases}
$
and if $\Phi:[0,1]^2 \ra \sI(\sU,\sU^\prime)$ is a homotopy between $\sigma, \tau$ 
through paths from $f$ to $g$, then there is a commutative diagram
\begin{tikzcd}%[sep=small]
{f_*\bX} \ar{d} \ar{r}
&{f_*\bX \wedge I_+}  \ar{d} \arrow[r, dashed, shift right=1]
&{f_*\bX} \ar{d} \arrow[l, shift right=1]\\
{\sigma \ltimes \bX} \ar{r} 
&{\Phi \ltimes\bX}
&{\tau \ltimes \bX} \ar{l}\\
{g_*\bX} \ar{u} \ar{r} 
&{g_*\bX \wedge I_+} \ar{u} \arrow[r, dashed, shift left=1]
&{g_*\bX} \ar{u} \arrow[l, shift left=1]
\end{tikzcd}    
of natural homotopy equivalences, where 
\begin{tikzcd}%[sep=small]
{} \ar[dashed]{r} &{\circx} \ar{r} &{= \id}
\end{tikzcd}
.  Similar remarks apply to the $\bR f^*:\bHSPEC_{\sU^\prime} \ra \bHSPEC_{\sU}$.]\\
\endgroup %%------------------------------------<<

\begingroup%%----------------------------------->>
\fontsize{9pt}{11pt}\selectfont
\textbf{\small FACT} \  
If $\bX \ra \bY$ is a cofibration in $\bSPEC_{\sU}$ and if $\bY^\prime \ra \bX^\prime$ is a fibration  in 
$\bSPEC_{\sU^\prime}$, then the arrow 
$\HOM(\bY,\bY^\prime) \ra \HOM(\bX,\bY^\prime) \times_{\HOM(\bX,\bX^\prime)} \HOM(\bY,\bX^\prime)$ is a fibration in 
$\dcg/\sI(\sU,\sU^\prime)$ which is a weak equivalence if 
$\bX \ra \bY$ or $\bY^\prime \ra \bX^\prime$ is acyclic (the notation is that of the addendum on 
p. \pageref{16.41}).\\
\endgroup %%------------------------------------<<

\begin{proposition} \ %24
Suppose that \mA is a cofibrant object in $\dcg$ $-$then the functor 
$\sH\sO\sM[\alpha,-)$ preserves fibrations and acyclic fibrations (cf. Proposition 23).  Therefore the assumptions of the TDF theorem are met (cf. p. \pageref{16.42} ff.), so $\bL\alpha\ltimes -$ and $\bR\sH\sO\sM[\alpha,-)$ exist and 
$(\bL\alpha\ltimes - ,\bR\sH\sO\sM[\alpha,-))$ is an adjoint pair.\\
\end{proposition}

%%----------------------------------------------------------------------------------------------26
\begingroup%%----------------------------------->>
\fontsize{9pt}{11pt}\selectfont
\textbf{\small FACT} \  
Fix a cofibrant object \mA in $\dcg$ and let $H:IA \ra \sI(\sU,\sU^\prime)$ be a homotopy $-$then $\forall$ cofibrant \bX in 
$\bSPEC_{\sU}$, the arrow $H \circx i_t \ltimes \bX \ra H \ltimes \bX$ is a homotopy equivalence $(t \in \{0,1\})$.
\vspi
[Note: \ Consequently the functors 
$\bL \alpha \ltimes -:\bHSPEC_{\sU} \ra \bHSPEC_{\sU^\prime}$ corresponding to $\alpha:A \ra \sI(\sU,\sU^\prime)$ are naturally isomorphic, as are the functors 
$\bR\sHOM[\alpha,-):\bHSPEC_{\sU^\prime} \ra \bHSPEC_{\sU}$.]\\
\endgroup %%------------------------------------<<

\begingroup%%----------------------------------->>
\fontsize{9pt}{11pt}\selectfont
\textbf{\small FACT} \  
Let \mA, \mB be cofibrant objects in $\dcg$ and suppose that $\phi:A \ra B$ is a homotopy equivalence $-$then $\forall$ 
$\beta:B \ra \sI(\sU,\sU^\prime)$, the arrow $\beta \circx \phi \ltimes \bX \ra \beta \ltimes \bX$ is a homotopy equivalence provided that \bX is cofibrant.
\vspi
[Fix a homotopy inverse $\psi:B \ra A$ for $\phi$, choose $H:IA \ra A$ such that 
$
\begin{cases}
\ H \circx i_0 = \id_A\\
\ H \circx i_1 = \psi \circx \phi
\end{cases}
\&
$
$G:IB \ra B$ such that 
$
\begin{cases}
\ G \circx i_0 = \id_B\\
\ G \circx i_1 = \phi \circx \psi
\end{cases}
, \ 
$
and, keeping in mind the preceding result, use the commutative diagrams 
\[
\begin{tikzcd}%[sep=small]
{\beta \circx \phi \circx \psi \circx \phi \ltimes \bX} \ar{d} \ar{r}{i_1 \ltimes \id} \ar{d}
&{\beta \circx \phi \circx H \ltimes \bX} \ar{ldd}{H \ltimes \id}\\
{\beta \circx \phi \circx \psi \ltimes \bX} \ar{d}\\
{\beta \circx \phi \ltimes \bX} \arrow[r,shift right=0.5,dash] \arrow[r,shift right=-0.5,dash] 
&{\beta \circx \phi \ltimes \bX} \ar{uu}[swap]{i_0 \ltimes \id}
\end{tikzcd}
\hspace{1.0cm}
\begin{tikzcd}%[sep=small]
{\beta \circx \phi \circx \psi \ltimes \bX} \ar{d} \ar{r}{i_1 \ltimes \id} \ar{d}
&{\beta \circx G \ltimes \bX} \ar{ldd}{G \ltimes \id}\\
{\beta \circx \phi \ltimes \bX} \ar{d}\\
{\beta  \ltimes \bX} \arrow[r,shift right=0.5,dash] \arrow[r,shift right=-0.5,dash] 
&{\beta \ltimes \bX} \ar{uu}[swap]{i_0 \ltimes \id}
\end{tikzcd}
\]
to deduce that the arrow 
$\beta \circx \phi \ltimes \bX \ra \beta \ltimes \bX$ is a weak equivalence, hence a homotopy equivalence.]
\vspi
[Note: \ The cofibrancy assumption on \mA, \mB can be dropped.  Thus let $\bY^\prime$ be any 
$\sU^\prime$-spectrum.  
Given $U \subset \sU$, $\exists$ $U^\prime \subset \sU^\prime$: $U \approx U^\prime$ $\implies$ 
$\sHOM[\beta \circx \phi,\bY^\prime)_U \approx$ 
$(Y_{U^\prime}^\prime)^{A_+}$, 
$\sHOM[\beta,\bY^\prime)_U \approx$ $(\bY_{U^\prime}^\prime)^{B_+}$ (cf. Proposition 21).  
Because $\phi:A \ra B$ is a homotopy equivalence, it follows that 
$\sHOM[\beta,\bY^\prime)_U \ra$ 
$\sHOM[\beta\circx\phi,\bY^\prime)_U$ is a homotopy equivalence $\forall \ U$.  
But \bX is cofibrant, so 
$[\bX,-]_0 \approx$ $[\bX,-]$ (cf. p. \pageref{16.43}) (all objects are fibrant).  
Therefore 
$[\bX,\sHOM[\beta,\bY^\prime)]_0 \approx$ 
$[\bX,\sHOM[\beta\circx\phi,\bY^\prime)]_0$ $\implies$ 
$[\beta \ltimes \bX,\bY^\prime]_0 \approx$
$[\beta\circx\phi \ltimes \bX,\bY^\prime]_0$ (cf. Proposition 19).  
And this means that the arrow 
$\beta\circx\phi \ltimes \bX \ra$ $\beta \ltimes \bX$ is a homotopy equivalence ($\bY^\prime$ being arbitrary).  
Variant: The same conclusion obtains if \bX is tame.]\\
\endgroup %%------------------------------------<<

\begingroup%%----------------------------------->>
\fontsize{9pt}{11pt}\selectfont
\textbf{\small EXAMPLE} \  \ \ 
Take \ $\sU = \sU^\prime$ \ $-$then $\forall \ f \ \ \in \ \sI(\sU,\sU)$, \ there is a commutative diagram \quad
$
\begin{tikzcd}[sep=small]
{*} \ar{rddd}[swap]{f} \ar{rr}{f} &&{\sI(\sU,\sU)} \ar{lddd}{\id}\\
\\
\\
&{\sI(\sU,\sU)}
\end{tikzcd}
,
$
thus $\forall$ cofibrant \bX, the arrow $f_*\bX \ra \id \ltimes \bX$ is a homotopy equivalence.
\vspi
[Note: \ The point here is this: $\sI(\sU,\sU)$ is contractible but it is unknown whether it is a cofibrant object in $\dcg$.]\\
\endgroup %%------------------------------------<<

\begingroup%%----------------------------------->>
\fontsize{9pt}{11pt}\selectfont
\textbf{\small FACT} \  
Let $A,B$ be objects in $\dcg$ and suppose that $\phi:A \ra B$ is a closed cofibration $-$then $\forall$ 
$\beta:B \ra \sI(\sU,\sU^\prime)$, the arrow 
$\beta \circx \phi \ltimes \bX \ra \beta \ltimes \bX$ is a spectral cofibration provided that \bX is cofibrant.
\vspi
%%----------------------------------------------------------------------------------------------27
[Given a $\sU^\prime$-spectrum $\bY^\prime$, finding a filler for the diagram 
\begin{tikzcd}%[sep=small]
{\beta \circx \phi \ltimes \bX} \ar{dd} \ar{rr} &&{\beta \ltimes \bX} \ar{ld}\ar{dd}\\
&{\bY^\prime}\\
{(\beta \circx \phi \ltimes \bX) \wedge I_+} \ar{ru}\ar{rr} &&{(\beta \ltimes \bX) \wedge I_+} \ar[dashed]{lu} 
\end{tikzcd}
amounts to finding a filler for the diagram
$
\begin{tikzcd}%[sep=small]
{\bX} \ar{dd} \ar{r} &{\sHOM[\beta,\bY^\prime)} \ar{dd}\\
\\
{\bX \wedge I_+} \ar[dashed]{ruu} \ar{r} &{\sHOM[\beta \circx \phi,\bY^\prime)}
\end{tikzcd}
.  \ 
$
However, the arrow
$\sHOM[\beta,\bY^\prime) \ra$ 
$\sHOM[\beta \circx \phi,\bY^\prime)$ is a levelwise \bCG fibration, therefore is a levelwise Serre fibration, and, as \bX is cofibrant, the arrow $\bX \ra \bX \wedge I_+$ is an acyclic cofibration in our model category structure on $\bSPEC_{\sU}$ 
(cf. p. \pageref{16.44} ff.).]\\
\endgroup %%------------------------------------<<

\begingroup%%----------------------------------->>
\fontsize{9pt}{11pt}\selectfont
\textbf{\small EXAMPLE} \quad
Take \quad $\sU = \sU^\prime$ \quad $-$then \  $\forall \ f \in \sI(\sU,\sU)$, \quadx there is   a \ 
commutative diagram \quad
\begin{tikzcd}%[sep=small]
{*} \ar{rdd}[swap]{f} \ar{rr}{f} &&{\sI(\sU,\sU)} \ar{ldd}{\id}\\
\\
&{\sI(\sU,\sU)}
\end{tikzcd}
, thus $\forall$ cofibrant \bX, the arrow $f_*\bX \ra \id \ltimes \bX$ is a spectral cofibration.
\vspi
[In fact, $\sI(\sU,\sU)$ is $\Delta$-cofibered (cf. p. \pageref{16.45}), so $\forall \ f \in \sI(\sU,\sU)$, 
$\{f\} \ra \sI(\sU,\sU)$ is a closed cofibration (cf. p. \pageref{16.46}).]\\
\endgroup %%------------------------------------<<

Let $\sU$, $\sV$ be univserses.  Put 
$\sA \oplus \sB = \{U \oplus V: U \subset \sU$ $\&$ $\dim U < \omega$, $V \subset \sV$ $\&$ $\dim V < \omega \}$ 
(which is not the standard indexing set in $\sU \oplus \sV$).

\indent\indent $(\un{\wedge})$ \ 
Given \bX in $\bSPEC_{\sU}$ and \bY in $\bSPEC_{\sV}$, the data $\{X_U \#_k Y_V\}$ defines a 
$(\sU \oplus \sV,\sA \oplus \sB)$-prespectrum.  Spectrify and let $\bX \un{\wedge} \bY$ be its image in 
$\bSPEC_{\sU \oplus \sV}$ under the canonical equivalence 
$\bSPEC_{\sU \oplus \sV,\sA \oplus \sB} \ra$ 
$\bSPEC_{\sU \oplus \sV}$
provided by Proposition 15.

Examples: 
(1) 
$\bQ_U^\infty X \un{\wedge} \bQ_V^\infty Y \approx$ 
$\bQ_{U \oplus V}^\infty (X \#_k Y)$ ;
(2) 
$(\bX \wedge K) \un{\wedge} \bY \approx$ 
$(\bX \un{\wedge} \bY) \wedge K \approx$ 
$\bX \un{\wedge} (\bY \wedge K)$.

[Note: \ Take $X = Y = \bS^0$ in $(1)$ to get 
$\bS^{-U} \un{\wedge} \bS^{-V} \approx$ 
$\bS^{-(U \oplus V)}$.]

Remark: It is not literally true that $\un{\wedge}$ is an associative, commutative operation.  Consider, e.g., commutativity.  
If $\Tee:\sU \oplus \sV \ra \sV \oplus \sU$ is the switching map, then $\Tee_*(\bX \un{\wedge} \bY)$ 
is naturally isomorphic to 
$\bY \un{\wedge} \bX$.  
The situation for associativity is analogous (consider the isomomophism 
$\sU \oplus (\sV \oplus \sW) \approx$
$(\sU \oplus \sV) \oplus \sW$
of universes).\\

\begingroup%%----------------------------------->>
\fontsize{9pt}{11pt}\selectfont
Another way to proceed is this.  Write $\bX \un{\bbox} \bY$ for the composite 
\begin{tikzcd}%[sep=small]
{\bI_{\sU} \times \bI_{\sV}} \ar{r}{\bX \times \bY} &{\dcg_* \times \dcg_*}
\end{tikzcd}
\begin{tikzcd}%[sep=small]
{} \ar{r}{\#_k} &{\dcg_*}
\end{tikzcd}
$-$then, relative to the arrow 
$\bI_{\sU} \times \bI_{\bV} \ra \bI_{\sU \oplus \sV}((U,V) \ra U \oplus V)$, 
$\lan \bX \un{\bbox} \bY$ is a $\sU \oplus \sV$-prespectrum, i.e., an object of \ 
$\bV[\bI_{\sU \oplus \sV},\dcg_*]$, \ and its spectrification can be identified with $\bX \un{\wedge} \bY$.  \ 
Therefore 
$\un{\wedge}: \bSPEC_{\sU} \times \bSPEC_{\sV} \ra \bSPEC_{\sU \oplus \sV}$ is a continuous functor.\\
\endgroup %%------------------------------------<<

\begingroup%%----------------------------------->>
\fontsize{9pt}{11pt}\selectfont
\textbf{\small FACT} \  
Suppose given $\alpha:A \ra \sI(\sU,\sU^\prime)$ and $\beta:B \ra \sI(\sV,\sV^\prime)$.  Let 
$\alpha \times_{\oplus} \beta$ be the composite 
\begin{tikzcd}%[sep=small]
{A \times_k B} \ar{r}{\alpha \times_{k} \beta} &{\sI(\sU,\sU^\prime) \times_k \sI(\sV,\sV^\prime)}
\end{tikzcd}
$\overset{\oplus}{\lra} \sI(\sU \oplus \sV,\sU^\prime \oplus \sV^\prime)$ $-$then 
$(\alpha \times_{\oplus} \beta) \ltimes (\bX \un{\wedge} \bY) \approx$ 
$(\alpha \ltimes \bX) \un{\wedge} (\beta \ltimes \bY)$.\\
\endgroup %%------------------------------------<<

%%----------------------------------------------------------------------------------------------28
Given \bY in $\bSPEC_{\sV}$ and \bZ in $\bSPEC_{\sU \oplus \sV}$, let $\bZ^{\bY}$ be the $\sU$-spectrum 
$U \ra \HOM(\bS^{-U} \un{\wedge}$ $\bY,\bZ)$ $-$then there is a natural homeomorphism 
$\HOM(\bX \un{\wedge} \bY,\bZ) \approx \HOM(\bX,\bZ^{\bY})$.

Example: 
$(\bZ^{\bS^{-V}})_U =$ 
$\HOM(\bS^{-U} \un{\wedge} \bS^{-V},\bZ) \approx$ 
$\HOM(\bS^{-(U \oplus V)},\bZ) \approx$ 
$\bZ_{U \oplus V}$.\\

\begin{proposition} \ %25
If $\bA \ra \bX$ is a cofibration in $\bSPEC_{\sU}$ and if $\bB \ra \bY$ is a cofibration in $\bSPEC_{\sV}$, then the arrow 
$\bA \un{\wedge}\bY \underset{\bA \un{\wedge}\bB}{\sqcup} \bX \un{\wedge}\bB \ra \bX \un{\wedge}\bY$ is a cofibration in 
$\bSPEC_{\sU \oplus \sV}$ which is acyclic if $\bA \ra \bX$ or $\bB \ra \bY$ is acyclic.\\
\end{proposition}

\begin{proposition} \ %26
If $\bB \ra \bY$ is a cofibration in $\bSPEC_{\sV}$ and if $\bZ \ra \bC$ is a fibration in $\bSPEC_{\sU \oplus \sV}$, then the arrow 
$\bZ^{\bY} \ra \bZ^{\bB} \times_{\bC^{\bB}} \bC^{\bY}$ is a fibration in $\bSPEC_{\sU}$ which is acyclic if $\bB \ra \bY$ or 
$\bZ \ra \bC$ is acyclic.\\
\end{proposition}

\begingroup%%----------------------------------->>
\fontsize{9pt}{11pt}\selectfont
\quad Propositions 25 and 26 are formally equivalent.  
To establish the fibration contention in Proposition 26, one can assume that $\bB \ra \bY$ has the form 
$\bS^{-V} \wedge L \ra \bS^{-V} \wedge K$, where $L \ra K$ is a cofibration in $\dcg_*$.  
The fact that $\bZ \ra \bC$ is a fibration in $\bSPEC_{\sU \oplus \sV}$ implies that the arrow 
$\texttt{HOM}(K,\bZ) \ra$ 
$\texttt{HOM}(L,\bZ) \times_{\texttt{HOM}(L,\bC)} \texttt{HOM}(K,\bC)$ is a fibration in 
$\bSPEC_{\sU \oplus \sV}$ which is acyclic if 
$L \ra K$ or $\bZ \ra \bC$ is acyclic (cf. p. \pageref{16.47}).  But the functor 
$(-)^{\bS^{-V}}$ preserves fibrations and acyclic fibrations and $\forall \ X$, 
$\hom(X,\bZ)^{\bS^{-V}} \approx \bZ^{\bS^{-V} \wedge X}$, thus the arrow 
$\bZ^{\bS^{-V} \wedge K} \ra \cdots$ is a fibration in $\bSPEC_{\sU}$ which is acyclic if $L \ra K$ or $\bZ \ra \bC$ is acyclic.
\vspi
[Note: \ The functor $\bQ_V^\infty = \bS^{-V} \wedge -$ preserves cofibrations and acyclic cofibrations.]
\vspi
\label{16.51}
Example: 
$
\begin{cases}
\ \bX\\
\ \bY
\end{cases}
$
cofibrant $\implies$ $\bX \un{\wedge}\bY$ cofibrant (cf. Proposition 25).\\
\endgroup %%------------------------------------<<

\begin{proposition} \ %27
Suppose that \bY is a cofibrant object in $\bSPEC_{\sV}$ $-$then the functor $(-)^{\bY}$ preserves fibrations and acyclic fibrations (cf. Proposition 26).  Therefore the assumptions of the TDF theorem are met 
(cf. p. \pageref{16.48} ff.), so 
$\bL(-\un{\wedge}\bY)$ and $\bR(-)^{\bY}$ exist and $(\bL(-\un{\wedge}\bY),\bR(-)^{\bY})$ is an adjoint pair.
\end{proposition}

[Note: \ Since all objects are fibrant, $(-)^{\bY}$ necessarily preserves weak equivalences
 (cf. p. \pageref{16.49}).]\\

\quad If \bC and \bD are model categories, then $\bC \times \bD$ becomes a model category upon imposing the obvious slotwise structure.  In particular: $\bSPEC_{\sU} \times \bSPEC_{\sV}$ is a model category.\\

\begin{proposition} \ %28
The functor $\un{\wedge}: \bSPEC_{\sU} \times \bSPEC_{\sV} \ra \bSPEC_{\sU \oplus \sV}$ 
sends weak equivalences between cofibrant objects to weak equivalences, thus the total left derived functor 
$\bL\un{\wedge}:\bHSPEC_{\sU} \times \bHSPEC_{\sV} \ra \bHSPEC_{\sU \oplus \sV}$ exists (cf. $\S 12$, Proposition 14).  
\end{proposition}

[Suppose that $\bA \ra \bX$ is an acyclic cofibration in $\bSPEC_{\sU}$ and $\bB \ra \bY$ is an acyclic cofibration in 
$\bSPEC_{\sV}$, where $\bA$ $\&$ $\bB$ (hence $\bX$ $\&$ $\bY$) are cofibrant.  Factor the arrow 
%%----------------------------------------------------------------------------------------------29
$\bA \un{\wedge} \bB \ra$ 
$\bX \un{\wedge} \bY$ as the composite 
$\bA \un{\wedge} \bB \ra$ 
$\bX \un{\wedge} \bB \ra$ 
$\bX \un{\wedge} \bY$.  Owing to Proposition 25, 
$\bA \un{\wedge} \bB \ra$ 
$\bX \un{\wedge} \bB \ra$ 
and 
$\bX \un{\wedge} \bB \ra$ 
$\bX \un{\wedge} \bY$ are acyclic cofibrations.  Therefore 
$\bA \un{\wedge} \bB \ra$ 
$\bX \un{\wedge} \bY$ is an acyclic cofibration.  The lemma on p. \pageref{16.50} then implies that 
$\un{\wedge}$ preserves weak equivalences between cofibrant objects.]

[Note: \ $\bL \un{\wedge}(\bX,\bY) = \sL\bX \un{\wedge} \sL\bY$, the value of the total left derived functor of 
$-\un{\wedge} \sL \bY$ at \bX (cf. Proposition 27).]\\

Take in the above $\sU = \sV$ and choose any $f \in \sI(\sU^2,\sU)$ $(\sU^2 = \sU \oplus \sU)$.  
Definition: 
$\bX \wedge_f \bY =$
$f_*(\bX \un{\wedge} \bY)$, 
$\hom_f(\bY,\bZ) = (f^*\bZ)^{\bY}$.  So: 
$\HOM(\bX \wedge_f \bY,\bZ) =$ 
$\HOM(f_*(\bX \un{\wedge}_f \bY),\bZ) \approx$ 
$\HOM(\bX \un{\wedge} \bY,f^*\bZ) \approx$ 
$\HOM(\bX,(f^* \bZ)^{\bY}) =$
$\HOM(\bX,\hom_f(\bY,\bZ))$.

[Note: \ While each of the functors $- \wedge_f \bY$ has a right adjoint 
$\bZ \ra \hom_f(\bY,\bZ)$, $\bSPEC_{\sU}$ is definitely not a symmetric monoidal category under 
$\otimes = \wedge_f$.]\\

\label{17.15}
\label{17.41}
\begingroup%%----------------------------------->>
\fontsize{9pt}{11pt}\selectfont
\textbf{\small EXAMPLE} \  
Write $\bQ^\infty$ in place of $\bQ_{\{0\}}^\infty$ and put $\bS = \bQ^\infty\bS^0$.  Letting 
$i:\sU \ra \sU \oplus \sU$ be the inclusion of $\sU$ onto the first summand, one has 
$i_*(\bX \wedge \bS^0) \approx$ 
$\bX \un{\wedge} \bS$, thus 
$(f \circx i)_*(\bX \wedge \bS^0) \approx$ 
$f_* \circx i_* (\bX \wedge \bS^0) \approx$ 
$f_*(\bX \un{\wedge} \bS) =$ 
$\bX \wedge_f \bS$.  And, when \bX is cofibrant, 
$\bX \wedge \bS^0 \approx$ 
$(f \circx i)_*(\bX \wedge \bS^0)$ in $\bHSPEC_{\sU}$, i.e., 
$\bX \approx \bX \wedge_f \bS$ in $\bHSPEC_{\sU}$.\\
\endgroup %%------------------------------------<<

Definition: $\bX \wedge \bY = \bL f_* (\bL \un{\wedge}(\bX,\bY))$, 
$\hom(\bY,\bZ) = \bR(\bR f^*\bZ)^{\sL \bY}$ ($= (f^*(\bZ)^{\sL\bY}$, all objects being fibrant).

[Note: \ This apparent abuse of notation is justified on the grounds that, up to natural isomorphism, these functors are independent of the choice of $f$ (cf. p. \pageref{16.50a}).  
Terminology: Call $\wedge$ the 
\un{smash product}
\index{smash product (spectra)}.]

Observation: Since $f_*$ sends cofibrant objects to cofibrant objects and $\sL\bX \un{\wedge} \sL \bY$ is cofibrant 
(cf. p. \pageref{16.51}), 
$[\bX \wedge \bY,\bZ]$  $=$ \hspace{0.03cm}
$[\bL f_*(\bL\un{\wedge} (\bX,\bY)),\bZ]$ $\approx$ \hspace{0.03cm}
$[\bL f_*(\sL\bX \un{\wedge}  \sL \bY),\bZ]$ $\approx$ \hspace{0.03cm}
$[f_*(\sL \bX \un{\wedge}  \sL \bY),\bZ]$ $\approx$ 
$\pi_0(\HOM(f_*(\sL \bX \un{\wedge}  \sL \bY),\bZ))$ $\approx$ 
$\pi_0(\HOM(\sL \bX \un{\wedge}  \sL \bY,f^*\bZ))$ $\approx$ 
$\pi_0(\HOM(\sL \bX,$ $(f^* \bZ)^{\sL \bY}))$ $\approx$ 
$[\sL \bX, (f^* \bZ)^{\sL \bY}]$ $\approx$ 
$[\bX, (f^* \bZ)^{\sL \bY}]$ $\approx$ 
$[\bX,\bR(\bR f^*\bZ)^{\sL \bY}] =$ 
$[\bX,\hom(\bY,\bZ)]$.\\

\begingroup%%----------------------------------->>
\fontsize{9pt}{11pt}\selectfont
\textbf{\small FACT} \  
In $\bHSPEC_{\sU}$, $\bX \wedge \bY \approx$ $\bX \wedge \bQ^\infty Y$, hence 
$\bQ^\infty(K \#_k L) \approx$ 
$(\bQ^\infty K) \wedge L \approx$ 
$\bQ^\infty K \wedge \bQ^\infty L$ and 
$\texttt{HOM}(K,\bX) \approx$ 
$\hom(\bQ^\infty K,\bX)$.\\
\endgroup %%------------------------------------<<

\begin{proposition} \ 
$\bHSPEC_{\sU}$ is a monoidal category.
\end{proposition}

[Taking $\otimes = \wedge$ and $e = \bS$ $(= \bQ^\infty \bS^0)$, one has to define natural isomorphisms
$
\begin{cases}
\ R_{\bX}:\bX \wedge \bS \ra \bX\\
\ L_{\bX}:\bS \wedge \bX \ra \bX
\end{cases}
$
and $A_{\bX,\bY,\bZ}$: $\bX \wedge (\bY \wedge \bZ) \ra (\bX \wedge \bY) \wedge \bZ$ satisfying MC$_1$ and MC$_2$ on 
p. \pageref{16.52}.  
The definitions of $R_{\bX}$ and $L_{\bX}$ are clear (cf. supra).  Letting $\Phi$ be the isomorphism 
$(\sU \hthree \oplus \hthree \sU) \hthree \oplus \hthree \sU \ra$
$\sU \hthree \oplus \hthree (\sU \hthree \oplus \hthree \sU)$, 
define $A_{\bX,\bY,\bZ}$ for cofibrant \bX, \bY ,\bZ via the following string of natural isomorphisms in $\bHSPEC_{\bU}$: 
$\bX \wedge (\bY \wedge \bZ) = \bL f_*(\bL \un{\wedge}(\bX,\bY \wedge \bZ)) \approx$ 
%%----------------------------------------------------------------------------------------------30
$\bL f_*(\bX \un{\wedge} \sL(\bY \wedge \bZ))$ \hthree $\approx$ \hthree
$\bL f_*(\bX \un{\wedge} \sL(\bL f_*(\bL\un{\wedge}(\bY,\bZ))))$ \hthree $\approx$ \hthree
$\bL f_*(\bX \un{\wedge} \sL(\bL f_*(\bY \un{\wedge} \bZ)))$ \hthree $\approx$ \hthree
$\bL f_*(\bX \un{\wedge}  f_*(\bY \un{\wedge} \bZ))$ \hthree $\approx$ \hthree
$f_*(\bX \un{\wedge} f_*(\bY \un{\wedge} \bZ))$ \hthree $\approx$ 
$f_* \circx (\id_{\sU} \oplus f)_* \circx \Phi_* ((\bX \un{\wedge} \bY)  \un{\wedge} \bZ)$ $\approx$
$f_* \circx (f \oplus \id_{\sU})_*((\bX \un{\wedge} \bY)  \un{\wedge} \bZ) \approx$
$f_*(f_*(\bX \un{\wedge} \bY)  \un{\wedge} \bZ) \approx$ 
$(\bX \wedge \bY) \wedge \bZ$ (reverse the steps).  
That MC$_1$ and MC$_2$  obtain can then be established by using the contractibility of $\sI(\sU^n,\sU)$.]

[Note: \ $\bHSPEC_{\sU}$ admits an evident compatible symmetry, thus is a symmetric monoidal category 
(cf. p. \pageref{16.53}).  Since each of the functors 
$-\wedge \bY:\bHSPEC_{\sU} \ra \bHSPEC_{\sU}$ has a right adjoint 
$\bZ \ra \hom(\bY,\bZ)$, it follows that $\bHSPEC_{\sU}$ is a closed category.]\\

\begingroup%%----------------------------------->>
\fontsize{9pt}{11pt}\selectfont
\label{16.22}
Therefore \bHSPEC is a closed category.\\
\endgroup %%------------------------------------<<

\label{17.75}
\begingroup%%----------------------------------->>
\fontsize{9pt}{11pt}\selectfont
\textbf{\small EXAMPLE} \  
If $\bff: \bX \ra \bY$, $\bg:\bZ \ra \bW$ are morphisms in \bHSPEC, then there is an exact triangle 
$\bX \wedge \bC_{\bg} \ra$ 
$\bC_{\bff \wedge \bg} \ra$ 
$\bC_{\bff}\wedge \bW \ra$ 
$\Sigma(\bX \wedge \bC_{\bg})$.]
\vspi
[Consider the factorization 
$\bff \wedge \bg =$
$\bff \wedge \id_{\bW} \circx \id_{\bX} \wedge \bg$ and use the result on p. \pageref{16.54}.]\\
\endgroup %%------------------------------------<<

\begingroup%%----------------------------------->>
\fontsize{9pt}{11pt}\selectfont
\textbf{\small FACT} \  
$\bX \wedge \bY$ is connective if $\bX$ $\&$ $\bY$ are connective.\\
\endgroup %%------------------------------------<<

Given a finite dimensional subspace $U$ of $\sU$, put 
$\Sigma^U \bX = \bX \wedge \bS^U$, 
$\Omega^U \bX = \texttt{HOM}(\bS^U,\bX)$ $-$then 
$(\Sigma^U,\Omega^U)$ is an adjoint pair.\\

\begin{proposition} \  %30
The total left derived functor $\bL\Sigma^U$ for $\Sigma^U$ exists and the total right derived functor $\bR\Omega^U$ for 
$\Omega^U$ exists.  And: $(\bL\Sigma^U,\bR\Omega^U)$ is an adjoint pair (cf. Proposition 12).\\
\end{proposition}

\begin{proposition} \ %31
The pair $(\bL\Sigma^U,\bR\Omega^U)$ is an adjoint equivalence of categories (cf. Proposition 13).
\end{proposition}

[Suppose that \bX is cofibrant $-$then in $\bHSPEC_{\sU}$ there are, on the one hand, natural isomorphisms 
$\Sigma^U(\bX \wedge \bS^{-U}) \approx$ 
$f_*(\bX \un{\wedge} \bS^{-U})  \wedge \bS^U \approx$ 
$f_*((\bX \un{\wedge} \bS^{-U})  \wedge \bS^U) \approx$ 
$f_*(\bX \un{\wedge}(\bS^{-U}) \wedge \bS^U)) \approx$ 
$f_*(\bX \un{\wedge} \bQ^\infty \bS^0) \approx \bX$, 
and on the other, natural isomorphisms 
$\Sigma^U \bX \wedge \bS^{-U} \approx$ 
$f_*(\Sigma^U \bX \un{\wedge} \bS^{-U}) \approx$ 
$f_*((\bX \wedge \bS^{U}) \un{\wedge} \bS^{-U}) \approx$ 
$f_*(\bX \un{\wedge}(\bS^{-U} \wedge \bS^U) \approx$ 
$f_*(\bX \un{\wedge}  \bQ^\infty \bS^0) \approx \bX$.  
Therefore $\bL \Sigma^U$ is an equivalence of categories and $- \wedge \bS^{-U} \approx \bR \Omega^U$.]\\

Fix a universe $\sU$ $-$then $S_n$ operates to the left on $\sU^n$ by permutations, hence 
$\forall \ \sigma \in S_n$ there are functors 
$\sigma_*:\bSPEC_{\sU^n} \ra$ $\bSPEC_{\sU^n}$.  Agreeing to write 
$S_n \ltimes -$ for the functor corresponding to the arrow 
$\chi_n:S_n \ra \sI(\sU^n,\sU^n)$, one has 
$S_n \ltimes \bX \approx$ 
$\bigvee\limits_{\sigma \in S_n} \sigma_* \bX$.  
The multiplication and unit of $S_n$ induce natural transformations 
$m_n:S_n \ltimes S_n \times - \ra$ $S_n \ltimes -$ \ $\&$ $\epsilon_n:\id \ra S_n \ltimes -$, so 
$(S_n \ltimes -,m_n,\epsilon_n)$ is a triple in $\bSPEC_{\sU^n}$. 
Its associated category of algebras is called the category of 
\un{$S_n$-spectra}
\index{S$_n$-spectra (category of)} 
(relative to $\sU)$: $S_n$-$\bSPEC_{\sU^n}$.  An $S_n$-spectrum 
%%----------------------------------------------------------------------------------------------31
is therefore a $\sU^n$-spectrum \bX equipped with a morphism 
$\bxi:S_n \ltimes \bX \ra \bX$ satisfying TA$_1$ and TA$_2$ (cf. p. \pageref{16.55} ff.), i.e., equipped with morphisms 
$\bxi_\sigma:\sigma_* \bX \ra \bX$ such that $\bxi_e = \id_{\bX}$ and 
$\bxi_\sigma \circx \sigma_*(\bxi_\tau) = \bxi_{\sigma \tau}$.

\label{16.60}
[Note: \ Given $(\bX,\bxi)$, $(\bY,\bbeta)$ in $S_n\text{-}\bSPEC_{\sU^n}$, write 
$S_n$-$\HOM(\bX,\bY)$ for $\Mor((\bX,\bxi),(\bY,$ $\bbeta))$ topologized via the equalizer diagram 
$\Mor((\bX,\bxi),(\bY,\bbeta)) \ra$ 
$\HOM(\bX,\bY)$ $\rightrightarrows$ 
$\HOM(S_n \ltimes \bX,\bY)$.]

Example: $\forall \ \bX$ in $\bSPEC_{\sU}$, 
$\bX^{(n)} \equiv \bX \un{\wedge}\cdots \un{\wedge} \bX$ (n factors) is an $S_n$-spectrum.

[Note: \ $\forall \ X \in \dcg_*$, $X^{(n)} \equiv X \#_k \cdots \#_k X$ (n factors) and 
$(\bQ^\infty X)^{(n)} \approx$ $\bQ^\infty(X^{(n)})$.]\\

\begingroup%%----------------------------------->>
\fontsize{9pt}{11pt}\selectfont
The functor $S_n \ltimes -$ is a left adjoint, hence preserves colimits.  
Since $\bSPEC_{\sU^n}$ is complete and cocomplete, 
specialization of the following generality allows one to conclude that $S_n$-$\bSPEC_{\sU^n}$ is complete and cocomplete.\\
\endgroup %%------------------------------------<<

\begingroup%%----------------------------------->>
\fontsize{9pt}{11pt}\selectfont
\textbf{\small LEMMA} \  
Suppose that \bC is a complete and cocomplete category.  Let $\bT = (T,m,\epsilon)$ be a triple in \bC.  
Assume: $T$ preserves filtered colimits $-$then $\bT$-\bALG is complete and cocomplete.
\vspi
[A proof can be found in 
Borceux\footnote[2]{\textit{Handbook of Categorical Algebra 2}, Cambridge University Press (1994), 206-211.}.]\\
\endgroup %%------------------------------------<<

\textbf{\small LEMMA} \  
Suppose that $A$ is a right $S_n$-space in $\dcg$.  Let 
$\alpha:A \ra \sI(\sU^n,\sU)$ be $S_n$-equivariant $-$then for every $S_n$-spectrum \bX, there is a coequalizer diagram 
$\alpha \ltimes S_n \ltimes \bX \rightrightarrows$ 
$\alpha \ltimes \bX \ra$ 
$\alpha \ltimes_{S_n} \bX$.

[One of the arrows is 
$\id_\alpha \ltimes \bxi$.  As for the other, 
$\alpha \ltimes S_n \ltimes \bX \approx$ 
$(\alpha \times_c \chi_n) \ltimes \bX$ (cf. p. \pageref{16.56}) and the diagram
\begin{tikzcd}[sep=small]
{A \times S_n} \ar{rdd}[swap]{\alpha \times_c \chi_n} \ar{rr}{\pi} &&{A} \ar{ldd}{\alpha}\\
\\
&{\sI(\sU^n,\sU)}
\end{tikzcd}
commutes ($\pi(a,\sigma) =$ $a \cdot \sigma$).  
Proof: 
$\alpha \times_c \chi_n(a,\sigma) =$ 
$\alpha(a) \circx \chi_n(\sigma)$, 
$\alpha \circx \pi(a,\sigma) =$ 
$\alpha(a \cdot \sigma) = \alpha(a) \cdot \sigma$ and $\forall \ u \in \sU^n$, 
$(\alpha(a) \circx \chi_n(\sigma))(u) =$
$\alpha(a)(\sigma \cdot u) =$
$(\alpha(a) \cdot \sigma) (u)$ (by the very definition of the right action of $S_n$ on $\sI(\sU^n,\sU))$.]\\

Remark: $\alpha \ltimes_{S_n} -$ is a functor from $S_n$-$\bSPEC_{\sU^n}$ to $\bSPEC_{\sU}$.  
On the other hand, 
$\sHOM[\alpha,-)$ is a functor from $\bSPEC_{\sU}$ to $S_n$-$\bSPEC_{\sU^n}$.  And: 
$\HOM(\alpha \ltimes_{S_n} \bX,\bY) \approx$ 
$S_n$-$\HOM(\bX,\sHOM[\alpha,\bY))$.\\

\begingroup%%----------------------------------->>
\fontsize{9pt}{11pt}\selectfont
It is sometimes necessary to consider 
\un{$G$-spectra},
\index{G-spectra} where $G$ is a subgroup of $S_n$ (the objects of $G$-$\bSPEC_{\sU^n}$ are thus the algebras per $G \ltimes -$).  
Given a subgroup $K$ of $G$, there is a forgetful functor 
$G$-$\bSPEC_{\sU^n}$ $\ra$ $K$-$\bSPEC_{\sU^n}$ and, in obvious notation, it has a left adjoint 
$G \ltimes_K -$,  so that 
$G$-$\HOM(G \ltimes_K \bX,\bY) \approx$ $K$-$\HOM(\bX,\bY)$.\\
\endgroup %%------------------------------------<<

%%----------------------------------------------------------------------------------------------32
\begingroup%%----------------------------------->>
\fontsize{9pt}{11pt}\selectfont
\textbf{\small FACT} \  
Let 
$U$: $G$-$\bSPEC_{\sU^n} \ra \bSPEC_{\sU^n}$ be the forgetful functor.  
Call a morphism $\bff:\bX \ra \bY$ of $G$-spectra a 
weak equivalence if $U\bff$ is a weak equivalence, a fibration if $U\bff$ is a fibration, and a cofibration if \bff has the LLP w.r.t. acyclic fibrations $-$then with these choices, $G$-$\bSPEC_{\sU^n}$ is a model category.
\vspi
[Note: \ This is the 
\un{external structure}
\index{external structure ($G$-$\bSPEC_{\sU^n}$)}.  
To define the 
\un{internal structure}
\index{internal structure ($G$-$\bSPEC_{\sU^n}$)}, 
stipulate that $\bff:\bX \ra \bY$ is a weak equivalence or a fibration if for each finite dimensional $G$-stable $U \subset \sU^n$, and for each subgroup $K \subset G$, the induced map of fixed point spaces 
$X_U^K \ra Y_U^K$ is a weak equivalence or a fibration and let the cofibrations be the \bff which have the LLP w.r.t. acyclic fibrations.  
Example: Take $G = S_n$ $-$then $\forall$ cofibrant \bX in $\bSPEC_{\sU}$, $\bX^{(n)}$ is cofibrant in the internal structure on 
$S_n$-$\bSPEC_{\sU^n}$.]\\
\endgroup %%------------------------------------<<

The preceding formalities are the spectral counterpart of a standard topological setup.  
Thus given a right $S_n$-space $A$ in $\dcg$ and a left $S_n$-space $X$ in $\dcg_*$, define 
$A \ltimes_{S_n} X$ by the coequalizer diagram 
$(A \times S_n)_+ \#_k X \rightrightarrows$ 
$A_+ \#_k X \ra$ 
$A \ltimes_{S_n} X$ $((A \times S_n)_+ \approx$ 
$A_+ \#_k S_{n+}$)
$-$then $A \ltimes_{S_n}-$ is a functor from the category of pointed \dsp compactly generated left $S_n$-spaces to the category of pointed \dsp compactly generated spaces.  It has a right adjoint, viz. the functor that sends Y to $Y^{A_+}$ 
($(\sigma \cdot f)(a) = f(a \cdot \sigma)$, with trivial action on the disjoint base point).

Example: Let $\sC$ be a \dsp creation operator, i.e., a functor 
$\sC:\bGamma_{\ini}^\OP \ra \dcg$ such that $\sC_0 = *$ $-$then in the notation of $\S 14$, Proposition 27, the filtration quotient 
$\sC_n[X]/\sC_{n-1}[X]$ is homeomorphic to $\sC_n \ltimes_{S_n} X^{(n)}$.\\

\begingroup%%----------------------------------->>
\fontsize{9pt}{11pt}\selectfont
\textbf{\small FACT} \  
\label{16.61}
\label{16.66}
$\forall \ X$ in $\dcg_*$, $\alpha \ltimes_{S_n} (\bQ^\infty X)^{(n)} \approx \bQ^\infty (A \ltimes_{S_n} X^{(n)})$.\\
\endgroup %%------------------------------------<<

\index{Extended Powers (example)}
\begingroup%%----------------------------------->>
\fontsize{9pt}{11pt}\selectfont
\textbf{\small EXAMPLE \ (\un{Extended Powers})} \  
Take $A = XS_n$ (which is $S_n$-universal) and fix an equivariant arrow 
$XS_n \ra \sI(\sU^n,\sU)$.  Using suggestive notation, the assignment 
$\bX \ra XS_n \ltimes_{S_n} X^{(n)}$ specifies a functor 
$D_n:\bSPEC_{\sU} \ra \bSPEC_{\sU}$ (conventionally, $D_0\bX = \bS$), the $n^\text{th}$ 
\un{extended power}.  Defining 
$D_n:\dcg_* \ra \dcg_*$ in exactly the same way, one has 
$D_n\bQ^\infty X =$ 
$XS_n \ltimes_{S_n} (\bQ^\infty X)^{(n)} \approx$ 
$\bQ^\infty(XS_n \ltimes_{S_n} X^{(n)}) =$ 
$\bQ^\infty(D_n X)$.  
Example: $D_n\bS^0 = B S_{n+}$ ($\implies$ 
$\ds\bigvee\limits_{n \geq 0} D_n\bS^0 =$ 
$B\bM_{\infty +}$, $\bM_\infty$ the permutative category of p. \pageref{16.57}).]
\vspi
[Note: \ Extended powers have many applications in homotopy theory.  For an account, see Bruner
\footnote[2]{\textit{SLN} \textbf{1176} (1986).}
et al..]\\
\endgroup %%------------------------------------<<

Let $\sC$ be a $\Delta$-separated creation operator $-$then $\forall \ X$ in $\dcg_*$, the 
\un{realization}
\index{realization} 
$\sC[X]$ of $\sC$ at \mX is 
$\ds\int^{\bn} \sC_n \times_k X^n$ (cf.  p. \pageref{16.58} (the assumption there that $(X,x_0)$ be wellpointed has been omitted here)), so $\sC[X]$ can be described by the coequalizer diagram 
$\coprod\limits_{\gamma:\bm \ra \bn} \sC_n \times_k X^m$ 
$\overset{u}{\underset{v}{\rightrightarrows}}$
$\coprod\limits_{m \geq 0} \sC_m \times_k X^m \ra$
$\sC[X]$
%%----------------------------------------------------------------------------------------------33
(on the term indexed by $\gamma:\bm \ra \bn$, 
$u$ is the arrow 
$\sC_n \times_k X^m \ra$ $\sC_n \times_k X^n$ 
and $v$ is the arrow 
$\sC_n \times_k X^m \ra$ $\sC_m \times_k X^m$).  
It is this interpretation of $\sC[X]$ that carries over to spectra provided they are unital.

Definition: A 
\un{unital $\sU$-spectrum}
\index{unital $\sU$-spectrum} 
is a pair $(\bX,\be)$, where 
$\be:\bS \ra \bX$ is a morphism of $\sU$-spectra.  Therefore the unital $\sU$-spectra are simply the objects of the category 
$\bS\backslash \bSPEC_{\sU}$.  
Example: $\forall \ X$ in $\dcg_*$, map $\bS^0$ to $X_+$ by 
$
\begin{cases}
\ 0 \ra x_0\\
\ 1 \ra *
\end{cases}
$
$-$then $\bQ^\infty X_+$ is unital.

[Note: \ Morphisms in $\bS\backslash \bSPEC_{\sU}$ are termed unital.]

Let \bX be a unital $\sU$-spectrum.  Viewing $\sU$ as an object in $\dcg_*$ with base point 0, each 
$\gamma:\bm \ra \bn$ in $\bGamma_{\ini}$ induces a linear isometry 
$\gamma:\sU^m \ra \sU^n$ and $\gamma_*\bX^{(m)}$ can be identified with 
$\bX_1 \un{\wedge} \cdots \un{\wedge} \bX_n$, $\bX_j$ being \bX if $\gamma^{-1}(j) \neq \emptyset$ and \bS if 
$\gamma^{-1}(j) = \emptyset$.  
There is an arrow 
$\gamma_*\bX^{(m)} \approx$ 
$\bX_1 \un{\wedge} \cdots \un{\wedge} \bX_n \ra$ $\bX^{(n)}$ which is $\id_{\bX}$ or \be according to whether 
$\bX_j = \bX$ or \bS.

Suppose now that $\phi:\sC \ra \sL$ is a morphism of creation operators, where $\sL$ is the linear isometries operad attached to our universe (recall that $\sL$ extends to a \dsp creation operator (cf. $\S 14$, Proposition 35)) $-$then 
$\forall \ n$, $\phi_n:\sC_n \ra \sL_n$ ($= \sI(\sU^n,\sU)$) is $S_n$-equivariant.  
Given a morphism $\gamma:\bm \ra \bn$ in $\bGamma_{\ini}$, let $\phi_\gamma:\sC_n \ra \sL_m$ be either composite in the commutative diagram 
\begin{tikzcd}%[sep=small]
{\sC_n} \ar{d} \ar[dashed]{rd} \ar{r} &{\sC_m} \ar{d}\\
{\sL_n} \ar{r} &{\sL_m}
\end{tikzcd}
and for \bX in $\bS\backslash \bSPEC_{\sU}$, put 
$\sC_\gamma \ltimes \bX^{(m)} =$ $\phi_\gamma \ltimes \bX^{(m)}$, 
$\sC_m \ltimes \bX^{(m)} =$ $\phi_m \ltimes \bX^{(m)}$
to get an arrow 
$\sC_\gamma \ltimes \bX^{(m)} \ra$ $\sC_m \ltimes \bX^{(m)}$.  
The 
\un{realization}
\index{realization} 
$\sC[\bX]$ of $\sC$ at \bX is then defined by the coequalizer diagram 
$\bigvee\limits_{\gamma:\bm \ra \bn} \sC_\gamma \ltimes \bX^{(m)}$ 
$\overset{u}{\underset{v}{\rightrightarrows}}$
$\bigvee\limits_{m \geq 0} \sC_m \ltimes \bX^{(m)} \ra$
$\sC[\bX]$
(on the term indexed by $\gamma:\bm \ra \bn$, 
$\bu$ is the arrow 
$\sC_\gamma \ltimes \bX^{(m)} \approx$ 
$\sC_n \ltimes \gamma_*\bX^{(m)} \ra$ 
$\sC_n \ltimes \bX^{(n)}$
and $\bv$ is the arrow 
$\sC_\gamma \ltimes \bX^{(m)} \ra$  
$\sC_m \ltimes \bX^{(m)}$).

[Note: \ The isomorphism $\sC_\gamma \ltimes \bX^{(m)} \approx$ 
$\sC_n \ltimes \gamma_*\bX^{(n)}$
is an instance of the ``composition rule'' on p. \pageref{16.59}.  To see this, consider 
$* \overset{\gamma}{\ra}$ 
$\sI(\sU^m,\sU^n)$ and 
$\sC_n \overset{\phi_n}{\lra}$
$\sI(\sU^n,\sU)$:
$\phi_n \times_c \gamma =$ 
$\phi_\gamma$ $\implies$ 
$\phi_\gamma \ltimes \bX^{(m)} \approx$ $\phi_n \ltimes \gamma_* \bX^{(m)}$.]

Remark: $\sC[\bX]$ is unital (since $\bS = \sC_0 \ltimes \bX^{(0)}$ and $\sC[?]$ is functorial.\\

\begin{proposition} \ %32
Let $\sC$ be a \dsp creation operator, augmented over $\sL$ via $\phi:\sC \ra \sL$ $-$then $\forall$ \mX in $\dcg_*$, 
$\sC[\bQ^\infty X_+] \approx \bQ^\infty\sC[X]_+$.
\end{proposition}

[Apply $\bQ^\infty$ to the coequalizer diagram 
$\bigvee\limits_{\gamma:\bm \ra \bn} \sC_{n+} \#_k (X_+)^{(m)} \rightrightarrows$ 
$\bigvee\limits_{m \geq 0} \sC_{m+} \#_k (X_+)^{(m)}$ $\ra$ 
$\sC[X]_+$.]

[Note: \ The isomorphism is natural in \mX.]\\

%%----------------------------------------------------------------------------------------------34
The coequalizer diagram describing $\sC[X]$ can be reduced to 
$\coprod\limits_{n \geq 0} \coprod\limits_{0 \leq i \leq n} \sC_{n+1} \times_k X^n$ 
$\overset{u}{\underset{v}{\rightrightarrows}}$
$\coprod\limits_{n \geq 0}\sC_n \times_{S_n} X^n \ra$
$\sC[X]$ 
and the coequalizer diagram describing $\sC[\bX]$ can be reduced to 
$\bigvee\limits_{n \geq 0}\bigvee\limits_{0 \leq i \leq n} \sC_{\sigma_i} \ltimes \bX^{(n)} 
\overset{u}{\underset{v}{\rightrightarrows}}$
$\bigvee\limits_{n \geq 0} \sC_n \ltimes_{S_n} \bX^{(n)} \ra$
$\sC[\bX]$, the $(n,i)^\text{th}$ term being indexed on $\sigma_i:\bn \ra \bn + \textbf{1}$ $(0 \leq i \leq n)$ 
(notation as in the proof of Proposition 35 in $\S 14$).  
There is also a coequalizer diagram
$\coprod\limits_{m \leq n-1} \coprod\limits_{0 \leq j \leq m} \sC_{m+1} \times_k X^m \overset{u}{\underset{v}{\rightrightarrows}}$
$\coprod\limits_{m \leq n}\sC_m \times_{S_m} X^m \ra$
$\sC_n[X]$ 
(cf. $\S 14$, Proposition 27).  
Here, 
$\sC_0[X] = *$, 
$\sC[X] = \colimx \sC_n[X]$, and the arrows 
$\sC_n[X] \ra \sC_{n+1}[X]$ are closed embeddings.  
Proceeding by analogy, define $\sC_n[\bX]$ by the coequalizer diagram 
$\bigvee\limits_{m \leq n-1}\bigvee\limits_{0 \leq j \leq m} \sC_{\sigma_j} \ltimes \bX^{(m)}$  
$\overset{u}{\underset{v}{\rightrightarrows}}$
$\bigvee\limits_{m \leq n} \sC_m \ltimes_{S_m} \bX^{(m)} \ra$
$\sC_n[\bX]$
$-$then $\sC_0[\bX] = \bS$, 
$\sC[\bX] = \colimx \sC_n[\bX]$, and the arrows 
$\sC_n[\bX] \ra \sC_{n+1}[\bX]$ 
are levelwise closed embeddings if $\be:\bS \ra \bX$ is a levelwise closed embedding.

Recalling that $X_*^{n+1}$ is the subspace of $X^{n+1}$ consisting of those points having at least one coordinate the base point $x_0$, the commutative diagram
%\begin{tikzcd}%[sep=small]
%{\sC_{n+1} \ltimes_{S_{n+1}} X_*^{(n+1)}} \ar{d} \ar{r} &{\sC_n[X]} \ar{d}\\
%{\sC_{n+1} \ltimes_{S_{n+1}} X^{(n+1)}} \ar{r} &{\sC_{n+1}[X]}
%\end{tikzcd}
\begin{tikzcd}%[sep=small]
{\sC_{n+1} \ltimes_{S_{n+1}} X_*^{(n+1)}} \ar{d} \ar{r} &{}\\
{\sC_{n+1} \ltimes_{S_{n+1}} X^{(n+1)}} \ar{r} &{}
\end{tikzcd}
\begin{tikzcd}%[sep=small]
{\sC_n[X]} \ar{d}\\
{\sC_{n+1}[X]}
\end{tikzcd}
is a pushout square.  
To formulate its spectral analog, one first has to define $\bX_*^{n+1}$.  The arrow 
$\bX^{(n)} \un{\wedge} \bS \ra \bX^{(n+1)}$ is a morphism of $S_n$-spectra $(S_n \subset S_{n+1})$, 
hence determines by adjointness a morphism 
$\btheta:S_{n+1} \ltimes_{S_n} (\bX^{(n)} \un{\wedge} \bS) \ra$ 
$\bX^{(n+1)}$ of $S_{n+1}$-spectra.  
Noting that 
$S_{n+1} \ltimes_{S_n} (\bX^{(n)} \un{\wedge} \bS) \approx$ 
$\bigvee\limits_{0 \leq i \leq n} \bX^{(i)} \un{\wedge} \bS \un{\wedge} \bX^{(n-i)}$, the arrows 
$\bX^{(n-1)} \un{\wedge} \bS \un{\wedge} \bS \ra$ 
$\bX^{(n)} \un{\wedge} \bS \subset$ 
$S_{n+1} \ltimes_{S_n} (\bX^{(n)} \un{\wedge} \bS)$,
$\bX^{(n-1)} \un{\wedge} \bS \un{\wedge} \bS \ra$ 
$\bX^{(n-1)} \un{\wedge} \bS \un{\wedge} \bX \subset$ 
$S_{n+1} \ltimes_{S_n} (\bX^{(n)} \un{\wedge} \bS)$ are morphisms of $S_{n-1}$-spectra 
$(S_{n-1} \subset$ $S_n \subset$ $S_{n+1})$, hence determine by adjointness morphisms 
$\bff,\bg: S_{n+1} \ltimes_{S_{n-1}} (\bX^{(n-1)} \un{\wedge} \bS \un{\wedge} \bS) \ra$ 
$S_{n+1} \ltimes_{S_n} (\bX^{(n)} \un{\wedge} \bS)$ of $S_{n+1}$-spectra.  
One then defines 
$\bX_*^{(n+1)}$ by the coequalizer diagram 
$S_{n+1} \ltimes_{S_{n-1}} (\bX^{(n-1)} \un{\wedge} \bS \un{\wedge} \bS) 
\overset{\bff}{\underset{\bg}{\rightrightarrows}}$
$S_{n+1} \ltimes_{S_n} (\bX{(n)} \un{\wedge} \bS) \ra$ 
$\bX_*^{(n+1)}$
(calculated in $S_{n+1}$-$\bSPEC_{\sU^{n+1}}$ (cf. p. \pageref{16.60})).
Since $\btheta$ coequalizes $(\bff,\bg)$, there is a morphism 
$\bX_*^{(n+1)} \ra \bX^{(n+1)}$ of $S_{n+1}$-spectra (which is a levelwise closed embedding if this is the case of 
$\be:\bS \ra \bX$).  Finally, the composites 
$\sC_{n+1} \ltimes (\bX^{(i)} \un{\wedge} \bS \un{\wedge} \bX^{(n-i)}) \approx$ 
$\sC_{\sigma_i} \ltimes \bX^{(n)} \ra$ 
$\sC_n \ltimes \bX^{(n)} \ra \sC_n[\bX]$ give rise to an arrow 
$\sC_{n+1} \ltimes_{S_{n+1}} \bX_*^{(n+1)} \ra$ 
$\sC_n[\bX]$ and the commutative diagram
\begin{tikzcd}%[sep=small]
{\sC_{n+1} \ltimes_{S_{n+1}} \bX_*^{(n+1)}} \ar{d} \ar{r} &{\sC_n[\bX]} \ar{d}\\
{\sC_{n+1} \ltimes_{S_{n+1}} \bX^{(n+1)}} \ar{r} &{\sC_{n+1}[\bX]}
\end{tikzcd}
is a pushout square.

Observation: The forgetful functor 
$\bS \backslash \bSPEC_{\sU} \ra \bSPEC_{\sU}$ has a left adjoint 
$\bX \ra$ 
$\bS \vee \bX$ ($\be:\bX \ra \bS \vee \bX$ is the inclusion of the wedge summand \bS).\\

\begin{proposition} \ %33
Let $\sC$ be a \dsp creation operator, augmented over $\sL$ via $\phi:\sC \ra \sL$ $-$then there is an isomorphism 
$\sC[\bS \vee \bX] \approx \bigvee\limits_{n\geq 0} \sC_n \ltimes_{S_n} \bX^{(n)}$ natural in \bX.
\end{proposition}

%%----------------------------------------------------------------------------------------------35
[In fact, 
$(\bS \vee \bX)^{(n+1)} \approx$ 
$(\bS \vee \bX)_*^{(n+1)}  \vee \bX^{(n+1)}$ as $S_{n+1}$-spectra, thus by induction, 
$\sC_n[\bS \vee \bX] \approx$ 
$\bigvee\limits_{m \leq n} \sC_m \ltimes_{S_m} \bX^{(m)}$ $(n \geq 0)$.]\\

\begingroup%%----------------------------------->>
\fontsize{9pt}{11pt}\selectfont
The spacewise version of Proposition 33 is the relation 
$\sC[X_+] \approx \ds\coprod\limits_{n \geq 0} \sC_n \times_{S_n} X^n$.\\
\endgroup %%------------------------------------<<

\textbf{\small LEMMA} \  
Suppose that $(X,x_0)$ is \dsp and wellpointed $-$then there are unital morphisms 
$\bQ^\infty X_+ \ra \bS \vee \bQ^\infty X$ and 
$\bS \vee \bQ^\infty X \ra \bQ^\infty X_+$ which are unital homotopy equivalences.

[Note: \ A homotopy \bH is unital if $\forall \ t$, $\bH_t$ is unital.]\\

\begin{proposition} \ %34
Let $\sC$ be a \dsp creation operator, augmented over $\sL$ via $\phi:\sC \ra \sL$ 
$-$then $\forall$ \dsp, wellpointed \mX, there is a natural weak equivalence 
$\bQ^\infty\sC[X] \ra \bigvee\limits_{n\geq 1} \bQ^\infty (\sC_n \ltimes_{S_n} X^{(n)})$ 
of $\sU$ -spectra.
\end{proposition}

[$\sC[X]$  is \dsp and wellpointed (cf. $\S 14$, Proposition 27).  
The lemma thus provides a weak equivalence 
$\bS \vee \bQ^\infty \sC[X] \ra$
$\bQ^\infty \sC[X]_+ \approx$ 
$\sC[\bQ^\infty X_+]$ (cf. Proposition 32).  
But 
$\sC[?]: \bS \backslash \bSPEC_{\sU} \ra \bS \backslash \bSPEC_{\sU}$ is a continuous functor, so it's homotopy preserving.  
Accordingly, there is a weak equivalence 
$\sC[\bQ^\infty X_+] \ra$ 
$\sC[\bS \vee \bQ^\infty X] \approx$ 
$\bigvee\limits_{n \geq 0} \sC_n \ltimes_{S_n} (\bQ^\infty X)^{(n)}$ 
(cf. Proposition 33).  And: 
$\bigvee\limits_{n \geq 0} \sC_n \ltimes_{S_n} (\bQ^\infty X)^{(n)} \approx$ 
$\bS \vee \bigvee\limits_{n \geq 1} \sC_n \ltimes_{S_n} (\bQ^\infty X)^{(n)} \approx$ 
$\bS \vee \bigvee\limits_{n \geq 1} \bQ^\infty (\sC_n \ltimes_{S_n} X^{(n)})$  
(cf. p. \pageref{16.61}).  
The weak equivalence in question now follows upon quotienting out by \bS.]\\

Application: 
$\bQ^\infty \sC[X]$ and 
$\bigvee\limits_{n \geq 1} \bQ^\infty(\sC_n[X]/\sC_{n-1}[X])$ 
are isomorphic in $\bHSPEC_{\sU}$.\\

\textbf{\small LEMMA} \  
Let $X$, $Y$, be in $\dcg_{*c}$ and let $f:X \ra Y$ be a pointed continuous function.  Assume: $f$ is a weak homotopy equivalence $-$then 
$\bQ^\infty f: \bQ^\infty X \ra \bQ^\infty Y$ is a weak equivalence.

[Since it suffices to work in \bSPEC, one has only to show that the 
$\pi_n^s(f):\pi_n^s(X) \ra \pi_n^s(Y)$ $(n \geq 0)$ are bijective ($\bQ^\infty X$, $\bQ^\infty Y$ being connective 
(cf. p. \pageref{16.62})).  But 
$\pi_n^s(X) =$ 
$\colimx\pi_{n+q}(\Sigma^q X)$, 
$\pi_n^s(Y) =$ 
$\colimx  \pi_{n+q}(\Sigma^q Y)$,
and 
$\Sigma^q f: \Sigma^q X \ra \Sigma^q Y$ is a weak homotopy equivalence (cf. p. \pageref{16.63}).]\\

\begin{proposition} \ %35
Let
$
\begin{cases}
\ \sC\\
\ \sD
\end{cases}
$
be creation operators, where $\forall \ n$, 
$
\begin{cases}
\ \sC_n\\
\ \sD_n
\end{cases}
$
is a compactly generated Hausdorff space and the action of $\sS_n$ is free.  Suppose given an arrow $\phi:\sC \ra \sD$ such that 
$\forall \ n$, $\phi_n:\sC_n \ra \sD_n$ is a weak homotopy equivalence $-$then $\forall$ \dsep, wellpointed \mX, there is a weak equivalence $\bQ^\infty\sC[X] \ra \bQ^\infty\sD[X]$.
\end{proposition}

%%----------------------------------------------------------------------------------------------36
[$\sC[X]$ and $\sD[X]$ are \dsp and wellpointed (cf. $\S 14$, Proposition 27).  But the hypotheses imply that $\phi$ induces a weak homotopy equivalence $\sC[X] \ra \sD[X]$ (cf. p. \pageref{16.64}).]\\

Application: Let $\sC$ be a creation operator, where $\forall \ n$, $\sC_n$ is a compactly generated Hausdorff space and the action of $S_n$ is free $-$then $\forall$ \dsp, wellpointed \mX, there is a natural weak equivalence 
$\bQ^\infty \sC[X] \ra \bigvee\limits_{n \geq 1} \bQ^\infty (\sC_n \ltimes_{S_n} X^{(n)})$ of $\sU$-spectra.

[The projection $\sC \times \sL \ra \sL$ augments $\sC \times \sL$ over $\sL$.  
On the other hand, $\forall \ n$, the projection $\sC_n \times_k \sL_n \ra \sC_n$ is a weak homotopy equivalence.  Quote Propositions 34 and 35.]

[Note: \ To justify the tacit use of the lemma, it is necessary to observe that 
$(\sC_n \times_k \sL_n) \ltimes_{S_n} X^{(n)}$, 
$\sC_n \ltimes_{S_n} X^{(n)}$ are wellpointed and the arrow 
$(\sC_n \times_k \sL_n) \ltimes_{S_n} X^{(n)} \ra$
$\sC_n \ltimes_{S_n} X^{(n)}$ is a weak homotopy equivalence.]\\

Example:  In $\bHSPEC_{\sU}$, 
$\bQ^\infty\BV^q[X] \approx \bigvee\limits_{n \geq 1} \bQ^\infty(\BV(R(q),n)\ltimes_{S_n} X^{(n)})$.

[Note: \ $\BV^q[X]$ can be replaced by $\Omega^q\Sigma^q X$ 
if \mX is path connected (May's approximation theorem).]

Example: \ In $\bHSPEC_{\sU}$, 
$\bQ^\infty\BV^\infty[X] \approx \bigvee\limits_{n \geq 1} \bQ^\infty(\BV(R(\infty),n)\ltimes_{S_n} X^{(n)})$.

[Note: \ $\BV^\infty[X]$ can be replaced by $\Omega^\infty\Sigma^\infty X$ if \mX is path connected and $\Delta$-cofibered (cf. $\S 14$, Proposition 33) (\mX $\Delta$-cofibered $\implies$ $\Omega^\infty\Sigma^\infty X$ wellpointed 
(cf. p. \pageref{16.65}).]\\

\begingroup%%----------------------------------->>
\fontsize{9pt}{11pt}\selectfont
\textbf{\small EXAMPLE} \  
Take $\sC = \bPER$ $-$then in $\bHSPEC_{\sU}$, 
$\bQ^\infty \bPER [X] \approx$ 
$\ds\bigvee\limits_{n \geq 1} \bQ^\infty(X S_n \ltimes_{S_n} X^{(n)}) \approx$ 
$\ds\bigvee\limits_{n \geq 1} D_n\bQ^\infty X$ (cf. p. \pageref{16.66}).\\
\endgroup %%------------------------------------<<

\textbf{\small LEMMA} \  
Let \bS be a triple in a category \bC and let \bT be a triple in the category \bSALG of \bS algebras $-$then the category 
\bT-(\bSALG) of \bT-algebras in \bSALG is isomorphic to the category $\bT \circx \bSALG$ of $\bT \circx \bS$ algebras in \bC.\\

Let $\sO$ be a reduced operad in $\bDelta$-\bCG, augmented over $\sL$ via $\phi:\sO \ra \sL$ $-$then $\sO$ determines a triple $\bT_{\sO} = (T_{\sO},m,\epsilon)$ in $\bS\backslash \bSPEC_{\sU}$ (cf. $\S 14$, Proposition 36) 
($T_{\sO} \bX = \sO[\bX]$, the realization of $\sO$ at \bX).  
But $\sO$ also determines a triple 
$\ov{\bT}_{\sO} = (\ov{T}_{\sO},\ov{m},\ov{\epsilon})$ in $\bSPEC_{\sU}$, where 
$\ov{T}_{\sO} [\bX] = \bigvee\limits_{n \geq 0} \sO_n \ltimes_{S_n} \bX^{(n)}$.  
To explain the connection between the two, note that $\bS\backslash \bSPEC_{\sU} = \bSALG$, \mS the functor that sends \bX to $\bS \wedge \bX$.  
And, according to Proposition 33, $T_{\sO} \circx S$ ``is'' $\ov{T}_{\sO}$, so by the lemma, the categories 
$\bT_{\sO}$-\bALG, $\ov{\bT}_{\sO}$-\bALG are isomorphic.\\
%%%%%%%%%%%%%%%%%%%%%%%%%%%%%%%%%%%%%%
%%%%%%%%%%%%%%%%%%%%%%%%%%%%%%%%%%%%%%
%%%%%%%%%%%%%%%%%%%%%%%%%%%%%%%%%%%%%%

\begin{center}
$\S \ 16$
\\[0.5cm]
$\mathcal{REFERENCES}$\\[-.5cm]
\end{center}
\vspace{-.7cm}

\[
\textbf{BOOKS}
\]

\begingroup
\fontsize{9pt}{11pt}\selectfont
\setlength\parindent{0 cm}

[1] \quad Adams, J., 
\textit{Stable Homotopy and Generalised Homology}, University of Chicago (1974).
\\[-.2cm]

[2] \quad Cole, M., 
\textit{Complex Oriented $RO(G)$-Graded Equivariant Cohomology Theories and Their Formal}

\hspace{0.8cm}{Group Laws}, Ph.D. Thesis, University of Chicago, Chicago (1996).
\\[-.2cm]

[3] \quad Elmendorf, A., Kriz, I., Mandell, M., and May, J., 
\textit{Rings, Modules, and Algebras in Stable Homotopy}

\hspace{0.8cm}\textit{Theory}, American Mathematical Society (1997).
\\[-.2cm]

[4] \quad Hopkins, M., 
\textit{Course Notes for Elliptic Cohomology}, MIT (1996).
\\[-.2cm]

[5] \quad Lewis, L., 
\textit{The Stable Category and Generalized Thom Spectra}, Ph.D. Thesis, University of 

\hspace{0.8cm}Chicago, Chicago (1978).
\\[-.2cm]

[6] \quad Lewis, L., May, J., and Steinberger, M., 
\textit{Equivariant Stable Homotopy Theory}, Springer Verlag (1986).
\\[-.2cm]

[7] \quad May, J. et al., 
\textit{Equivariant Homotopy and Cohomology Theory}, American Mathematical Society (1996).
\\[-.7cm]
\endgroup

\[
\textbf{ARTICLES}
\]

\begingroup
\fontsize{9pt}{11pt}\selectfont
\setlength\parindent{0 cm}

[1] \quad Block, J. and Lazarev, A., Homotopy Theory and Generalized Duality for Spectral Sheaves, 
\textit{Internat.}

\hspace{0.8cm}\textit{Math. Res. Notices} \textbf{20} (1996), 983-996.
\\[-.2cm]

[2] \quad Carlsson, G., A Survey of Equivariant Stable Homotopy Theory, 
\textit{Topology} \textbf{31} (1992), 1-27.
\\[-.2cm]

[3] \quad Elmendorf, A., Greenlees, J., Kriz, I., and May, J., Commutative Algebra in Stable Homotopy Theory 

\hspace{0.8cm}and a Completion Theorem, 
\textit{Math. Res. Lett.} \textbf{1} (1994), 225-239.
\\[-.2cm]

[4] \quad Elmendorf, A., Kriz, I., Mandell, M., and May, J., Modern Foundations for Stable Homotopy Theory, 

\hspace{0.8cm}In: 
\textit{Handbook of Algebraic Topology}, I. James (ed.), North Holland (1995), 213-253.
\\[-.2cm]

[5] \quad Goerss, P. and Hopkins, M., Notes on Multiplicative Stable Homotopy Theory, 
\textit{unpublished notes}.
\\[-.2cm]

[6] \quad Greenlees, J. and May, J., Equivariant Stable Homotopy Theory, In: 
\textit{Handbook of Algebraic Topology}, 

\hspace{0.8cm}I. James (ed.), North Holland (1995), 277-323.
\\[-.2cm]

[7] \quad Hovey, M., Shipley, B., and Smith, J., Symmetric Spectra, 
\textit{J. Amer. Math. Soc.} \textbf{13} (2000) 149-208.
\\[-.2cm]

[8] \quad May, J., Derived Categories in Algebra and Topology, 
\textit{Rend. Istit. Mat. Univ. Trieste} \textbf{25} (1993), 

\hspace{0.8cm}363-377.
\\[-.2cm]

[9] \quad Puppe, D., On the Stable Homotopy Category, In: 
\textit{Topology and its Applications}, D. Kurepa (ed.), 

\hspace{0.8cm}Budva (1972) 200-212.
\\[-.2cm]

[10] \quad Thomason, R., Symmetric Monoidal Categories Model All Connective Spectra, 
\textit{Theory Appl. Categ.} 

\hspace{0.95cm}\textbf{1} (1995), 78-118.

\setlength\parindent{2em}

\endgroup

\chapter{
$\boldsymbol{\S}$\textbf{17}.\quadx  STABLE HOMOTOPY THEORY}
\setlength\parindent{2em}
\setcounter{proposition}{0}

%%----------------------------------------------------------------------------------------------01
$\text{ }$\\[-1.25cm]

A complete treatment  of stable homotopy theory would require a book of many pages.  
Therefore, to avoid getting bogged down in a welter of detail, 
I shall admit some of the results without proof and keep the calculations to a minimum.  
Despite working within these limitations, it is nevertheless still possible to gain a reasonable understanding of the subject in the ``large''.

Recapitulation: The stable homotopy category \bHSPEC is a triangulated category satisfying the octahedral axiom 
(cf. $\S 16$, Proposition 14).  
Furthermore, \bHSPEC is a monogenic compactly generated CTC 
(cf. p. \pageref{17.1}) and admits Adams representability (by Neeman's countability criterion).

[Note: \bS is the unit in \bHSPEC and $\Sigma^{-1}$ stands for $\Omega$ 
(cf. p. \pageref{17.2}), so 
$\Lambda^{\pm1} \approx \Sigma^{\pm1}$ (recall the convention on 
p. \pageref{17.3}).]\\

\begingroup%%----------------------------------->>
\fontsize{9pt}{11pt}\selectfont
\textbf{\small EXAMPLE \  (\un{Complex K-Theory})} \  
Let $\bU = \colimx \bU(n)$ be the infinite unitary group $-$then \bU is a pointed CW complex and there is a pointed homotopy equivalence 
$\bU \ra \Omega^2 \bU$ (Bott periodicity).  Therefore the prescription 
$X_q = \Omega^k \bU$ ($q \equiv 1 - k \mod 2$ $(0 \leq k \leq 1)$) defines an $\Omega$-prespectrum \bX and by definition 
$\bK\bU = eM\bX$ 
(cf. p. \pageref{17.4}) is the spectrum of complex K-theory.\\
\endgroup %%------------------------------------<<

\begingroup%%----------------------------------->>
\fontsize{9pt}{11pt}\selectfont
\textbf{\small EXAMPLE \ (\un{Real K-Theory})} \  
Let $\bO = \colimx \bO(n)$ be the infinite orthogonal group $-$then \bO is a pointed CW complex and there is a pointed homotopy equivalence 
$\bO \ra \Omega^8 \bO$ (Bott periodicity).  Therefore the prescription 
$X_q = \Omega^k \bO$ ($q \equiv 7 - k \mod 8$ $(0 \leq k \leq 7)$) defines an $\Omega$-prespectrum \bX and by definition 
$\bK\bO = eM\bX$ 
(cf. p. \pageref{17.5}) is the spectrum of real K-theory.\\
\endgroup %%------------------------------------<<

A 
\un{$\Z$-graded cohomology theory $E^*$ on \bSPEC}
\index{Z-graded cohomology theory $E_*$ on \bSPEC} 
is a sequence of exact cofunctors $E^n:\bHSPEC \ra \bAB$ and a sequence of natural isomorphisms 
$\sigma^n:E^{n+1} \circ \Sigma \ra E^n$ such that the $E^n$ convert coproducts into products.  
$\bCT_{\Z}(\bSPEC)$ is the category whose objects are the $\Z$-graded cohomology theories on \bSPEC and whose morphisms 
$\Xi^*:E^* \ra F^*$ are sequences of natural transformations $\Xi^n:E^n \ra F^n$ such that the diagram
\begin{tikzcd}[sep=large]
{E^{n+1} \circ \Sigma} \ar{d}[swap]{\sigma^n} \ar{r}{\Xi^{n+1} \Sigma} &{F^{n+1} \circ \Sigma} \ar{d}{\sigma^n}\\
{E^n} \ar{r}[swap]{\Xi^n} &{F^n}
\end{tikzcd}
commutes $\forall \ n$.\\

Definition: The $\Z$-graded cohomology theory $\bE^*$ on \bSPEC attached to a spectrum \bE is given by 
$\bE^n(\bX) = [\bX,\Sigma^n\bE]$ $(= \pi_{-n}(\hom(\bX,\bE)))$.

[Note: The 
\un{coefficient groups}
\index{coefficient groups (of a $\Z$-graded cohomology theory)} 
of $\bE^*$ are the 
$\bE^n(\bS)$ $(= \pi_{-n}(\bE))$, i.e., $\bE^*(\bS) = \pi_{-*}(\bE)$ $(= \pi_*(\bE)^\OP)$.]\\

%%----------------------------------------------------------------------------------------------02
%
Remark: Owing to the Brown representability theorem 
(cf. p. \pageref{17.6}), every $\Z$-graded cohomology theory on \bSPEC is naturally isomorphic to some $\bE^*$, thus \bHSPEC is the represented equivalent of $\bCT_{\Z}(\bSPEC)$.

[Note: \ Needless to say, $\Mor(\bE^*,\bF^*) \approx [\bE,\bF]$.]\\

\begingroup%%----------------------------------->>
\fontsize{9pt}{11pt}\selectfont
\textbf{\small EXAMPLE}  \ 
Take $\bE = \bS$ $-$then the corresponding $\Z$-graded cohomology theory on \bSPEC is called 
\un{stable cohomotopy}
\index{stable cohomotopy}, 
the coefficient groups being
$
\begin{cases}
\ 0 \quad  \ \ \ (n > 0)\\
\ \Z \quad  \ \ \hspace{0.03cm} (n = 0)\\
\ \pi_{-n}^s \ \ (n < 0)
\end{cases}
.
$
\vspi
[Note: \ As on p. \pageref{17.7}, the $\pi_{-n}^s$ are the stable homotopy groups of spheres.]\\
\endgroup %%------------------------------------<<

\begingroup%%----------------------------------->>
\fontsize{9pt}{11pt}\selectfont
\textbf{\small LEMMA}  \  
If 
$
\begin{cases}
\ \pi_n(\bX) = 0 \quad (n < 0)\\
\ \pi_n(\bY) = 0 \quad (n > 0)
\end{cases}
, \ 
$
then $\pi_0[\bX,\bY] \ra \Hom(\pi_0(\bX),\pi_0(\bY))$ is an isomorphism.\\
\endgroup %%------------------------------------<<
\vspace{0.2cm}

\label{15.53}
\begingroup%%----------------------------------->>
\fontsize{9pt}{11pt}\selectfont
\label{17.29}
\textbf{\small EXAMPLE}  \ 
\bHSPEC carries a t-structure (cf. p. \pageref{17.8}), and the elements of its heart are the 
\un{Eilenberg-MacLane spectra}.
\index{Eilenberg-MacLane spectra}  
An explanation for the terminology is that 
$\pi_0:\bH(\bHSPEC) \ra \bAB$ is an equivalence of categories.  
To see this, consider the functor 
$\bH:\bAB \ra \bH(\bHSPEC)$ that sends $\Z$ to 
$\tau^{\geq 0}\tau^{\leq 0}\bS \approx \tau^{\leq 0}\tau^{\geq 0} \bS$, 
defining $\bH(\pi)$ for an arbitrary abelian group $\pi$ by the exact triangle 
$\ds\bigvee\limits_j \bH(\Z) \ra$ 
$\ds\bigvee\limits_i \bH(\Z) \ra$ 
$\bH(\pi) \ra$ 
$\ds\bigvee\limits_j \Sigma\bH(\bZ)$, where 
$0 \ra$ 
$\ds\bigoplus\limits_j \Z \ra$ 
$\ds\bigoplus\limits_i \Z \ra$ 
$\pi \ra 0$ is a presentation of $\pi$ (the lemma implies that 
$\pi_0:[\ds\bigvee\limits_j \bH(\Z), \ds\bigvee\limits_i \bH(\Z)] \ra$ 
$\Hom\bigl(\ds\bigoplus\limits_j \Z,\ds\bigoplus\limits_i \Z\bigr)$ 
is an isomorphism).  Therefore 
$\pi_0(\bH(\pi)) = \pi$, 
$\pi_n(\bH(\pi)) = 0$ $(n \neq 0)$ and 
$[\bH(\pi^\prime), \bH(\pi\pp)] = \Hom(\pi^\prime,\pi\pp)$.  
Example: 
$[\Sigma^{-1}\bH(\pi^\prime),\bH(\pi\pp)] =$ $\Ext(\pi^\prime,\pi\pp)$ but 
$\Ph(\Sigma^{-1}\bH(\pi^\prime),\bH(\pi\pp)) =$ $\Pur\Ext(\pi^\prime,\pi\pp)$ 
%(Christensen-Strickland\footnote[2]{dmc missing reference}).\\
(Christensen-Strickland\footnote[2]{\textit{Topology} \textbf{37} (1998), 339-364.}).
\vspi
[Note: \ Given $\pi$, $\exists$ an $\Omega$-prespectrum $\bK(\pi)$ such that 
$K(\pi)_q = K(\pi,q)$ (realized as a pointed CW complex with $K(\pi,0) = \pi$ (discrete topology)).  
Since $\pi_n(eM\bK(\pi)) = \colimx \pi_{n+q}(K(\pi)_q) =$
$
\begin{cases}
\ \pi \quad (n = 0)\\
\ 0  \quad \hspace{0.03cm} (n > 0)
\end{cases}
,
$
$eM\bK(\pi)$ ``is'' $\bH(\pi)$ ($M$ the cylinder functor of p. \pageref{17.9}).]\\
\endgroup %%------------------------------------<<

\begingroup%%----------------------------------->>
\fontsize{9pt}{11pt}\selectfont
\textbf{\small EXAMPLE}  \ 
Lin\footnote[3]{\textit{Proc. Amer. Math. Soc.} \textbf{56} (1976), 291-299.}
has shown that $\bS^*(\bH(\F_p)) = 0$, hence $D\bH(\F_p)$ is trivial and $[\bH(\F_p),\bK] = 0$ for all compact \bK.  Therefore the stable cohomotopy $\bS^*(\bH(\pi))$ of $\bH(\pi)$ vanishes if $\pi$ is torsion (but not in general (consider $\pi = \Z)$).
\vspi
[Note: \ $\Ph(\bH(\F_p),\bY)$ is a vector space over $\F_p$ which is nonzero for some \bY.  
Reason: If the contrary held, then $h_{\bH(\F_p)}$ would be projective and since 
$[\bH(\F_p),\bK] = 0$ for all compact \bK, it would follow that $\bH(\F_p) = 0$.]\\
\endgroup %%------------------------------------<<

\begin{proposition} \ %01
The graded abelian group $\bE^*(\bE)$ is a graded ring with unit.
\end{proposition}

%%----------------------------------------------------------------------------------------------03
[Given $\bff \in \bE^n(\bE)$, $\bg \in \bE^m(\bE)$, let $\bff\cdot \bg \in \bE^{n+m}(\bE)$ be the composite 
$\bE \overset{\bg}{\lra} \Sigma^m \bE$
$\overset{\Sigma^m \bff}{\lra} \Sigma^{n+m} \bE$ ($\id_{\bE} \in \bE^0(\bE)$ thus serves as the unit).]

[Note: \ $\forall \ \bX$, $\bE^*(\bX)$ is a graded left $\bE^*(\bE)$-module.]\\

\begingroup%%----------------------------------->>
\fontsize{9pt}{11pt}\selectfont
\textbf{\small EXAMPLE}  \ 
The $\F_p$-algebra $\bH(\F_p)^*(\bH(\F_p))$ is isomorphic to $\sA_p$, the mod $p$ Steenrod algebra.\\
\endgroup %%------------------------------------<<

\label{17.10} %dmc mnft
\begin{proposition} \ 
Fix a spectrum \bE $-$then $\forall \ n$ and $\forall \ \bX$, there is a short exact sequence 
$0 \ra$ 
$\lim^1 \bE^{n+q-1} (\bQ^\infty X_q) \ra$ 
$\bE^n(\bX) \ra$ 
$\lim \bE^{n+q}(\bQ^\infty X_q) \ra 0$.\\
\end{proposition}

\begingroup%%----------------------------------->>
\fontsize{9pt}{11pt}\selectfont
Specialized to the case $n = 0$, the conclusion is that the homomorphism 
$[\bX,\bE] \ra \lim[X_q,E_q]$ is surjective with kernel $\lim^1[\Sigma X_q,E_q]$.
\vspi
[Note: \ This is a recipe for the calculation of morphisms in \bHSPEC by means of morphisms in 
$\bH\dcg_*$.]\\
\endgroup %%------------------------------------<<

A 
\un{$\Z$-graded cohomology theory $E^*$ on $\bCW_*$} 
\index{Z-graded cohomology theory on $\bCW_*$} 
is a sequence of cofunctors $E^n:\bCW_* \ra \bAB$ and a sequence of natural isomorphisms  
$\sigma^n:E^{n+1} \circ \Sigma \ra E^n$ such that the $E^n$ convert coproducts into products and satisfy the following conditions.

\indent\indent (Homotopy) \quad If $f,g:X \ra Y$ are homotopic, then 
$E^n(f) = E^n(g):E^n(Y) \ra E^n(X)$ $\forall \ n$.

\indent\indent (Exactness) \quad If $(X,A,x_0)$ is a pointed \bCW pair, then the sequence 
$E^{n}(X/A) \ra$ 
$E^n(X) \ra$ 
$E^n(A) \ra$ is exact $\forall \ n$.

\indent\indent (Isotropy) \quad If $f:X \ra Y$ is a homotopy equivalence, then 
$E^n(f):E^n(Y) \ra E^n(X)$ is an isomorphism $\forall \ n$.

[Note: \ The homotopy axiom implies that a $\Z$-graded cohomology theory on $\bCW_*$ passes to 
$\bHCW_*$, thus the isotropy axiom is redundant.]

Example: Given a spectrum \bE, the assignment $X \ra \bE^n(\bQ^\infty X)$ defines a $\Z$-graded cohomology theory 
on $\bCW_*$.

$\bCT_{\Z}(\bCW_*)$ is the category whose objects are the $\Z$-graded cohomology theories on $\bCW_*$ and whose morphisms $\Xi^*:E^* \ra F^*$ are sequences of natural transformations $\Xi^n:E^n \ra F^n$ such that the diagram
\begin{tikzcd}%[sep=small]
{E^{n+1} \circ \Sigma} \ar{d}[swap]{\sigma^n} \ar{r}{\Xi^{n+1} \Sigma} &{F^{n+1} \circ \Sigma} \ar{d}{\sigma^n}\\
{E^n} \ar{r}[swap]{\Xi^n} &{F^n}
\end{tikzcd}
commutes $\forall \ n$.\\

Let $E^*$ be a $\Z$-graded cohomology theory on $\bCW_*$ $-$then the 
\un{coefficient groups}
\index{coefficient groups 
(of a $\Z$-graded cohomology theories on $\bCW_*$)} of $E^*$ are the $E^n(\bS^0)$.  
Example: Reduced singular cohomology with coefficients in an abelian group $\pi$ is a $\Z$-graded cohomology theory on 
$\bCW_*$ whose only nontrivial coefficient group is $\pi$ itself.

%%----------------------------------------------------------------------------------------------04	
[Note: \ $E^n(*) = 0$ $\forall \ n$.  Proof: $* \approx */*$, so the composite 
$E^n(*) \ra E^n(*) \ra E^n(*)$ is both the identity map and the zero map.]\\

\begingroup%%----------------------------------->>
\fontsize{9pt}{11pt}\selectfont
\textbf{\small FACT}  \  
Let $\pi$ be an abelian group.  Suppose that $E_1^*$, $E_2^*$ are $\Z$-graded cohomology theories on $\bCW_*$ such that 
$E_1^0(\bS^0) = \pi$, 
$E_2^0(\bS^0) = \pi$ and 
$E_1^n(\bS^0) = 0$,
$E_2^n(\bS^0) = 0$ $(n \neq 0)$ $-$then $E_1^*$, $E_2^*$  are naturally isomorphic.\\
\endgroup %%------------------------------------<<

\begingroup%%----------------------------------->>
\fontsize{9pt}{11pt}\selectfont
\textbf{\small EXAMPLE}  \ 
The $\Z$-graded cohomology theory  on $\bCW_*$ determined by $\bH(\pi)$ is naturally isomorphic to reduced singular cohomology $\widetilde{H}^*(-;\pi)$.\\
\endgroup %%------------------------------------<<

Notation: Let $T:\bCW^2 \ra \bCW^2$ be the functor that sends $(X,A)$ to $(A,\emptyset)$.

[Note: \ The 
\un{lattice}
\index{lattice (of $(X,A)$ in $\bCW^2$} 
of $(X,A)$ is the diagram
\[
\begin{tikzcd}%[sep=small]
&&{(X,\emptyset)} \ar{rd}\\
{(\emptyset,\emptyset)} \ar{r} &{(A,\emptyset)} \ar{ru} \ar{rd} &&{(X,A)} \ar{r}&{(X,X)}\\
&&{(A,A)} \ar{ru}
\end{tikzcd}
.]
\]

A 
\un{$\Z$-graded cohomology theory $H^*$ on $\bCW^2$} 
\index{Z-graded cohomology theory $H^*$ on $\bCW^2$} 
is a sequence of cofunctors $H^n:\bCW^2 \ra \bAB$ and a sequence of natural transformations 
$d^n:H^{n-1} \circ T \ra H^n$ such that the $H^n$ convert coproducts into products and satisfy the following conditions.

\indent\indent (Homotopy) \quad If $f,g:(X,A) \ra (Y,B)$ are homotopic, then 
$H^n(f) = H^n(g):H^n(Y,B) \ra H^n(X,A)$ $\forall \ n$.

\indent\indent (Exactness) \quad If $(X,A)$ is a \bCW pair, then the sequence 
$\cdots \ra$ 
$H^{n-1}(A,\emptyset) \overset{d^n}{\lra}$ 
$H^n(X,A) \ra$ 
$H^n(X,\emptyset) \ra$ 
$H^n(A,\emptyset) \overset{d^{n+1}}{\lra}$ 
$H^{n+1}(X,A) \ra$ 
$\cdots$ is exact.

\indent\indent (Excision) \quad If \mA, \mB are subcomplexes of \mX, then the arrow 
$H^n(A \cup B,B) \ra$ 
$H^n(A,A \cap B)$ is an isomorphism $\forall \ n$.

\indent\indent (Isotropy) \quad If $f:(X,A) \ra (Y,B)$ is a homotopy equivalence, then 
$H^n(f):H^n(Y,B) \ra H^n(X,A)$ is an isomorphism $\forall \ n$.

[Note: \ The homotopy axiom implies that a $\Z$-graded cohomology theory on $\bCW^2$ passes to 
$\bHCW^2$, thus the isotropy axiom is redundant.]

$\bCT_{\Z}(\bCW^2)$ is the category whose objects are the $\Z$-graded cohomology theories on $\bCW^2$ and whose morphisms $\Xi^*:H^* \ra G^*$ are sequences of natural transformations $\Xi^n:H^n \ra G^n$ such that the diagram
\begin{tikzcd}%[sep=small]
{H^{n-1} \circ T} \ar{d}[swap]{d^n} \ar{r}{\Xi^{n-1} T} &{G^{n-1} \circ T} \ar{d}{d^n}\\
{H^n} \ar{r}[swap]{\Xi^n} &{G^n}
\end{tikzcd}
commutes $\forall \ n$.\\
\vspace{0.25cm}

%%----------------------------------------------------------------------------------------------05
\begin{proposition} \ 
$\bCT_{\Z}(\bCW_*)$ and $\bCT_{\Z}(\bCW^2)$ are equivalent categories.
\end{proposition}

[On objects, consider the functor 
$\bCT_{\Z}(\bCW_*) \ra \bCT_{\Z}(\bCW^2)$ that sends $E^*$ to $H^*$, where 
$H^n(X,A) = E^n(X_+/A_+)$, and the functor 
$\bCT_{\Z}(\bCW^2) \ra \bCT_{\Z}(\bCW_*)$ that sends $H^*$ to $E^*$, where
$E^n(X) = H^n(X,\{x_0\})$.]

 [Note: \ Consult Whitehead\footnote[2]{\textit{Elements of Homotopy Theory}, Springer Verlag (1978), 571-600.}
for a verification  down to the last detail.]\\

The definition of a $\Z$-graded homology theory $E_*$ on $\bCW_*$, $\bCW^2$ is dual and, in obvious notation, the categories 
$\bHT_{\Z}(\bCW_*)$, $\bHT_{\Z}(\bCW^2)$ are equivalent (cf. Proposition 3).\\

\begingroup%%----------------------------------->>
\fontsize{9pt}{11pt}\selectfont
\textbf{\small FACT}  \  
Fix a $\Z$-graded cohomology theory $H^*$ on $\bCW^2$.  Let $(X,A)$ be a \bCW pair.  Suppose given a sequence $\{X_q\}$ of subcomplexes of \mX such that 
$A \subset X_0$, $X_q \subset X_{q+1}$, and $X = \colim X_q$ $-$then $\forall \ n$, there is a short exact sequence 
$0 \ra$ 
$\lim^1 H^{n-1}(X_q,A) \ra$ 
$H^n(X,A) \ra$ 
$\lim H^n(X_q,A) \ra 0$.
\vspi
[Note: \ Modulo some additional assumptions on $H^*$, one can establish a variant involving the finite subcomplexes of \mX which contain \mA 
(Huber-Meier\footnote[3]{\textit{Comment. Math. Helv.} \textbf{53} (1978), 239-257; 
see also Yosimura, \textit{Osaka J. Math.} \textbf{25} (1988), 881-890,
and Ohkawa, \textit{Hiroshima Math. J.} \textbf{23} (1993) 1-14.}).]\\
\endgroup %%------------------------------------<<

\begin{proposition} \ 
Let \bE be an $\Omega$-prespectrum $-$then the prescription 
$
E^n(X) = 
\begin{cases}
\ [X,E_n] \hspace{1.0cm} (n \geq 0)\\
\ [X,\Omega^{-n}E_0] \hspace{0.35cm} (n < 0)
\end{cases}
$
specifies a $\Z$-graded cohomology  theory on $\bCW_*$.
\end{proposition}

[Note: \ When \bE is a spectrum, $E^n(X) = \bE^n(\bQ^\infty X)$ 
(cf. p. \pageref{17.10}).]\\

\begin{proposition} \ 
Every $\Z$-graded cohomology theory $E^*$ on $\bCW_*$ is represented by an  $\Omega$-prespectrum \bE.
\end{proposition}

[Let $U:\bAB \ra \bSET$ be the forgetful functor $-$then $\forall \ n$, $U \circ E^n$ is representable 
(cf. p. \pageref{17.11} ff.): $U \circ E^n(X) \approx [X,E_n]$.  And: The $E_n$ $(n \geq 0)$ assemble into an 
$\Omega$-prespectrum.]\\

\begingroup%%----------------------------------->>
\fontsize{9pt}{11pt}\selectfont
The precise connection between $\Omega$-prespectra, spectra, and $\Z$-graded cohomology theories on $\bCW_*$ can be pinned down.  Thus let 
\bWPRESPEC 
\index{\bWPRESPEC} 
be the category whose objects are the prespectra and whose morphisms 
$\bff:\bX \ra \bY$ are sequences of pointed continuous functions $f_q:X_q \ra Y_q$ such that the diagram 
\begin{tikzcd}[sep=large]
{X_q} \ar{d} \ar{r}{f_q} &{Y_q} \ar{d}\\
{\Omega X_{q+1}} \ar{r}[swap]{\Omega f_{q+1}} &{\Omega Y_{q+1}}
\end{tikzcd}
is pointed homotopy commutative $\forall \ q$.  
Denote by 
\bHWPRESPEC 
\index{\bHWPRESPEC} 
the localization of \bWPRESPEC at the class of levelwise weak homotopy equivalences 
(there is no difficulty 
%%----------------------------------------------------------------------------------------------06
in seeing that this procedure leads to a category).  Write
\bHWOMEGAPRESPEC 
\index{\bHWOMEGAPRESPEC}
for the full subcategory of \bHWPRESPEC whose objects are the $\Omega$-prespectra $-$then 
$\Mor(\bX,\bY) = \lim[X_q,Y_q]$, where the limit is taken with respect to the composites 
$[X_{q+1},Y_{q+1}] \ra [\Omega X_{q+1},\Omega Y_{q+1}] \ra [X_q,Y_q]$.\\
\endgroup %%------------------------------------<<

\begingroup%%----------------------------------->>
\fontsize{9pt}{11pt}\selectfont
\textbf{\small FACT}  \  
\bHWOMEGAPRESPEC is the represented equivalent of $\bCT_{\Z}(\bCW_*)$ .\\
\endgroup %%------------------------------------<<

\begingroup%%----------------------------------->>
\fontsize{9pt}{11pt}\selectfont
Let \bHWSPEC 
\index{\bHWSPEC} 
be the full subcategory of \bHWOMEGAPRESPEC whose objects are the spectra.\\
\endgroup %%------------------------------------<<

\begingroup%%----------------------------------->>
\fontsize{9pt}{11pt}\selectfont
\textbf{\small FACT}  \  
The inclusion 
$\bHWSPEC \ra \bHWOMEGAPRESPEC$ is an equivalence of categories.
\vspi
[Consider the functor that on objects sends an $\Omega$-prespectrum \bX to $eM\bX$ (\mM as on 
p. \pageref{17.12}).]
\vspi
[Note: \ If $E^*$ is a $\Z$-graded cohomology theory on $\bCW_*$ which is represented by an $\Omega$-prespectrum
\bE, then $eM\bE$ is a spectrum which also represents $E^*$.]\\
\endgroup %%------------------------------------<<

\begingroup%%----------------------------------->>
\fontsize{9pt}{11pt}\selectfont
Summary: 
$\bHSPEC \leftrightarrow \bCT_{\Z}(\bSPEC)$, 
$\bHWSPEC \leftrightarrow \bCT_{\Z}(\bCW_*)$ 
and there is a functor 
$\bHSPEC \ra \bHWSPEC$ that on morphisms is the arrow $[\bX,\bY] \ra \lim[X_q,Y_q]$.  Accordingly, every  $\Z$-graded cohomology theory on $\bCW_*$ lifts to a $\Z$-graded cohomology theory on \bSPEC and every morphism of 
$\Z$-graded cohomology theories on $\bCW_*$ lifts to a morphism of $\Z$-graded cohomology theories on \bSPEC 
(but not uniquely due to the potential nonvanishing of $\lim^1[\Sigma X_q,Y_q]$ (cf. Proposition 2)).\\
\endgroup %%------------------------------------<<

A \ 
\un{$\Z$-graded homology theory $E_*$ on \bSPEC} \ 
\index{Z-graded homology theory $E_*$ on \bSPEC} 
is a sequence of exact functors \ 
$E_n:\bHSPEC$ $\ra$ $\bAB$ and a sequence of natural isomorphisms 
$\sigma_n:E_n \ra E_{n+1} \circ \Sigma$ such that the $E_n$ convert coproducts into direct sums.  
$\bHT_\Z(\bSPEC)$
\index{$\bHT_\Z(\bSPEC)$} 
is the category whose objects are the $\Z$-graded homology theories on \bSPEC 
and whose morphisms $\Xi_*:E_* \ra F_*$ are the sequences of natural transformations  $\Xi_n:E_n \ra F_n$ such that the diagram
\begin{tikzcd}%[sep=small]
{E_n} \ar{d}[swap]{\sigma_n} \ar{r}{\Xi_n} &{F_n} \ar{d}{\sigma_n}\\
{E_{n+1} \circ \Sigma} \ar{r}[swap]{\Xi_{n+1} \Sigma} &{F_{n+1} \circ \Sigma}
\end{tikzcd}
commutes $\forall \ n$.

Definition: The $\Z$-graded homology theory $\bE_*$ on \bSPEC attached to a spectrum \bE is given by 
$\bE_n(\bX) = \pi_n(\bE \wedge \bX)$.

[Note: The 
\un{coefficient groups}
\index{coefficient groups (of a $\Z$-graded homology theory)} 
of $\bE_*$ are the 
$\bE_n(\bS)$ $(= \pi_n(\bE))$, i.e., $\bE_*(\bS) = \pi_*(\bE)$.]

Remark: \label{17.26}
Because \bHSPEC admits Adams representability, every $\Z$-graded homology theory on \bSPEC is naturally isomorphic to some $\bE_*$ (cf. $\S 15$, Proposition 38), thus \bHSPEC/\bPh 
(cf. p. \pageref{17.13}) is the represented equivalent of $\bHT_\Z(\bSPEC)$.

[Note: \ Here $\Mor(\bE_*,\bF_*) \approx [\bE,\bF]/\bPh(\bE,\bF)$.]\\

%%----------------------------------------------------------------------------------------------07
\begingroup%%----------------------------------->>
\fontsize{9pt}{11pt}\selectfont
\textbf{\small EXAMPLE}  \ 
Take $\bE = \bS$ $-$then the corresponding $\Z$-graded homology theory on \bSPEC is called 
\un{stable homotopy}
\index{stable homotopy}, 
the coefficients groups being 
$
\begin{cases}
\ \pi_n^s \quad (n > 0)\\
\ \Z \ \quad \hspace{0.03cm} (n = 0)\\
\ 0 \  \quad \ (n < 0)
\end{cases}
.
$
\\
\endgroup %%------------------------------------<<

\begingroup%%----------------------------------->>
\fontsize{9pt}{11pt}\selectfont

\textbf{\small EXAMPLE}  \ 
For any two spectra \bE, \bF, the arrow 
$\pi_*(\bE) \otimes \pi_*(\bF) \otimes \Q \ra \pi_*(\bE \wedge \bF) \otimes \Q$ is an isomorphism.
\vspi
[Fix \bE and let \bF vary $-$then the arrow 
$\pi_*(\bE) \otimes \pi_*(-) \otimes \Q \ra \pi_*(\bE \wedge -) \otimes \Q$
is a morphism of $\Z$-graded homology theories on \bSPEC.  But 
$\pi_0^s(\bS) = \Z$ and 
$\pi_n^s(\bS)$ is finite if $n > 0$ (cf. p. \pageref{17.14}), hence 
$\pi_*(\bE) \otimes \pi_*(\bS) \otimes \Q \approx \pi_*(\bE \wedge \bS) \otimes \Q$.]\\
\endgroup %%------------------------------------<<

\begin{proposition} \ %6
Let 
$
\begin{cases}
\ \bE\\
\ \bF
\end{cases}
$
,
$
\begin{cases}
\ \bX\\
\ \bY
\end{cases}
$
be spectra $-$then there is an external product 
$\bE^*(\bX) \otimes \bF^*(\bY) \ra (\bE \wedge \bF)^*(\bX \wedge \bY)$ in cohomology.
\end{proposition}

[Work with the arrow 
$\hom(\bX,\bE) \wedge \hom(\bY,\bF) \ra \hom(\bX \wedge \bY,\bE \wedge \bF)$.]\\

\begin{proposition} \ %7
Let 
$
\begin{cases}
\ \bE\\
\ \bF
\end{cases}
$
,
$
\begin{cases}
\ \bX\\
\ \bY
\end{cases}
$
be spectra $-$then there is an external product 
$\bE_*(\bX) \otimes \bF_*(\bY) \ra (\bE \wedge \bF)_*(\bX \wedge \bY)$ in homology.
\end{proposition}

[Work with the arrow 
$\bE \wedge \bX \wedge \bF \wedge \bY \ra \bE \wedge \bF \wedge \bX \wedge \bY $.]\\

\begin{proposition} \ %8
Let 
$
\begin{cases}
\ \bE\\
\ \bF
\end{cases}
$
,
$
\begin{cases}
\ \bX\\
\ \bY
\end{cases}
$
be spectra $-$then there is an external slant product 
\begin{tikzcd}%[sep=small]
{\bE^*(\bX \wedge \bY) \otimes \bF_*(\bX)}  \ar{r}{/} &{(\bE \wedge \bF)^*(\bY).}
\end{tikzcd}
\end{proposition}

[Use the commutative diagram
\[
\begin{tikzcd}%[sep=small]
{\hom(\bX \wedge \bY),\bE) \wedge \bF \wedge \bX} \ar{d} \ar{r}{/} &{\hom(\bY,\bE \wedge \bF)} \ar{d}\\
{\hom(\bX,\hom(\bY,\bE)) \wedge \bX \wedge \bF} \ar{r} &{\hom(\bY,\bE) \wedge \bF}
\end{tikzcd}
.]
\]
\vspace{0.05cm}

\begin{proposition} \ 
Let 
$
\begin{cases}
\ \bE\\
\ \bF
\end{cases}
$
,
$
\begin{cases}
\ \bX\\
\ \bY
\end{cases}
$
be spectra $-$then there is an external slant product 
\begin{tikzcd}%[sep=small]
{\bE_*(\bX \wedge \bY) \otimes \bF^*(\bX)}  \ar{r}{\backslash} &{(\bE \wedge \bF)_*(\bY).}
\end{tikzcd}
\end{proposition}

[Use the commutative diagram
\[
\begin{tikzcd}%[sep=small]
{\bE \wedge \bX \wedge \bY \wedge \hom(\bX,\bF)} \ar{d} \ar{r}{\backslash} &{\bE \wedge \bF \wedge \bY}\\
{\bE \wedge \hom(\bX,\bF) \wedge \bX \wedge \bY} \ar{ru} 
\end{tikzcd}
.]
\]

%%----------------------------------------------------------------------------------------------08
\begingroup%%----------------------------------->>
\fontsize{9pt}{11pt}\selectfont
The external products are morphisms of graded abelian groups but this is not the case of the slant products.  Explicated: 
$\bE^n(\bX \wedge \bY) \otimes \bF_m(\bX) \overset{/}{\lra} (\bE \wedge \bF)^{n-m}(\bY)$ and 
$\bE_n(\bX \wedge \bY) \otimes \bF^m(\bX) \overset{\backslash}{\lra} (\bE \wedge \bF)_{n-m}(\bY)$, 
thus to get a morphism of graded abelian groups one must give $\bF_*(\bX)$ and $\bF^*(\bX)$ the opposite gradings.\\
\endgroup %%------------------------------------<<

\label{17.61}
A 
\un{ring spectrum}
\index{ring spectrum} 
is a ring object in \bHSPEC.  
Example: \bS is a commutative ring spectrum and every spectrum is an \bS-module.\\

\label{17.20} %dmc mnft
\begingroup%%----------------------------------->>
\fontsize{9pt}{11pt}\selectfont
\textbf{\small EXAMPLE}  \ 
Let \bk be a commutative ring with unit $-$then $\bH(\bk)$ is a commutative ring spectrum and for any \bk-module \mM, 
$\bH(M)$ is an $\bH(\bk)$-module.\\
\endgroup %%------------------------------------<<

\begingroup%%----------------------------------->>
\fontsize{9pt}{11pt}\selectfont
\textbf{\small EXAMPLE}  \ 
McClure\footnote[2]{\textit{SLN} \textbf{1176} (1986), 241-242.} has shown that \bKU is a commutative ring spectrum.  The homotopy $\pi_*(\bKU)$ of \bKU has period 2 and 
$\pi_0(\bKU) = \Z$, 
$\pi_1(\bKU) = 0$.  
In addition, there exists a multiplicatively invertible generator $\bb_{\bU} \in \pi_2(\bKU) \approx \Z$ inducing the homotopy periodicity and as a graded ring, 
$\pi_*(\bKU) \approx \Z[\bb_{\bU},\bb_{\bU}^{-1}]$.
\vspi
[Note: \ $\bKO$ is also a commutative ring spectrum.]\\
\endgroup %%------------------------------------<<

\begingroup%%----------------------------------->>
\fontsize{9pt}{11pt}\selectfont
\textbf{\small EXAMPLE}  \ 
For any \mX in $\bDelta$-$\bCG_*$, $(\Omega X)_+$ ($= \Omega X \amalg *$) is wellpointed, $\bQ^\infty((\Omega X)_+)$ is a ring spectrum, and $\pi_0(\Omega^\infty \Sigma^\infty(\Omega X)_+) \approx \Z[\pi_1(X)]$ (as rings).
\vspi
[To define the product, note that 
$\bQ^\infty(\Omega X)_+ \wedge \bQ^\infty(\Omega X)_+ \approx$ 
$\bQ^\infty((\Omega X)_+ \#_k(\Omega X)_+)$ (cf. p. \pageref{17.15}), which is isomorphic to 
$\bQ^\infty((\Omega X)_+ \times_k \Omega X)_+)$.]\\
\endgroup %%------------------------------------<<

\begingroup%%----------------------------------->>
\fontsize{9pt}{11pt}\selectfont
\textbf{\small FACT}  \  
If \bE is a connective ring spectrum, then 
$\Hom(\pi_0(\bE),\pi_0(\bE)) \approx [\bE,\bH(\pi_0(\bE))]$ and the arrow 
$\bE \ra \bH(\pi_0(\bE))$ realizing the identity $\pi_0(\bE) \ra \pi_0(\bE)$ is a morphism of ring spectra.\\
\endgroup %%------------------------------------<<

\begingroup%%----------------------------------->>
\fontsize{9pt}{11pt}\selectfont
\textbf{\small FACT}  \  
If \bE is a ring spectrum and $\be \ (= \tau^{\leq 0} \bE)$ is its connective cover, then \be admits a unique ring spectrum structure such that the arrow $\be \ra \bE$ is a morphism of ring spectra.\\
\endgroup %%------------------------------------<<

If \bE is a ring spectrum and \bF is an \bE-module, then the products figuring in the preceding propositions can be made ``internal'' through $\bE \wedge \bF \ra \bF$.

Example: Take $\bE = \bF$ and fix an \bX $-$then Proposition 8 furnishes an arrow 
$\bE^*(\bX) \otimes \bE_*(\bX) \overset{/}{\lra}$ 
$(\bE \wedge \bE)^*(\bS) \ra$ 
$\bE^*(\bS) = \pi_{-*}(\bE)$ 
and Proposition 9 furnishes an arrow 
$\bE_*(\bX) \otimes \bE^*(\bX) \overset{\backslash}{\lra}$ 
$(\bE \wedge \bE)_*(\bS) \ra$ 
$\bE_*(\bS) = \pi_{*}(\bE)$.\\

\label{17.24}
\label{17.54}
\begingroup%%----------------------------------->>
\fontsize{9pt}{11pt}\selectfont
\textbf{\small EXAMPLE}  \ 
Let \bE be a ring spectrum $-$then for spectra \bF $\&$ \bX, the 
\un{Hurewicz homomorphism}
\index{Hurewicz homomorphism (ring spectra)}
$\bF_*(\bX) \ra (\bE \wedge \bF)_*(\bX)$ is defined by the arrow 
$\bF_n(\bX) =$ 
$\pi_n(\bF \wedge \bX) \approx$ 
$\pi_n(\bS \wedge \bF \wedge \bX) \ra$ 
$\pi_n(\bE \wedge \bF \wedge \bX) =$
%%----------------------------------------------------------------------------------------------09
$(\bE \wedge \bF)_n(\bX)$ and the 
\un{Boardman homomorphism}
\index{Boardman homomorphism (ring spectra)}  
$\bF^*(\bX) \ra$ 
$(\bE \wedge \bF)^*(\bX)$ is defined by the arrow 
$\bF^n(\bX) =$ 
$[\bX,\Sigma^n \bF] \approx$ 
$[\bX,\Sigma^n(\bS \wedge \bF)] \ra$ 
$[\bX,\Sigma^n(\bE \wedge \bF)] =$ 
$(\bE \wedge \bF)^n(\bX)$.  
Assuming that both \bE and \bF are ring spectra, the commutative diagram 
\begin{tikzcd}[sep=large]
{\bF_n(\bX) \otimes \bF^m(\bX)} \ar{d} \ar{r} &{\pi_{n-m}(\bF)} \ar{d}\\
{(\bE \wedge \bF)_n(\bX) \otimes (\bE \wedge \bF)^m(\bX)} \ar{r} &{\pi_{n-m}(\bE \wedge \bF)}
\end{tikzcd}
serves to relate the two.
\vspi
[Note: \ In particular, there are arrows 
$
\begin{cases}
\ \bS_*(\bX) \ra \bE_*(\bX)\\
\ \bS^*(\bX) \ra \bE^*(\bX)
\end{cases}
$
.]\\
\vspace{0.25cm}
\endgroup %%------------------------------------<<

\label{17.19}
If \bE is a ring spectrum and \bF is an \bE-module, then $\forall \ \bX$, 
$
\begin{cases}
\ \bF^*(\bX) \\
\ \bF_*(\bX)
\end{cases}
$
is a graded 
$
\begin{cases}
\ \bE^*(\bS)\text{-module}\\
\ \bE_*(\bS)\text{-module}
\end{cases}
$
(cf. Propositions 6 and 7).

[Note: \ The structure is on the left.  Observe, however, that
$
\begin{cases}
\ \bE^*(\bX) \\
\ \bE_*(\bX)
\end{cases}
$
is a graded left and right
$
\begin{cases}
 \bE^*(\bS)\text{-module}\\
 \bE_*(\bS)\text{-module}
\end{cases}
$
, in fact, 
$
\begin{cases}
\ \bE^*(\bX) \\
\ \bE_*(\bX)
\end{cases}
$
is a graded
$
\begin{cases}
\ \bE^*(\bS)\text{-bimodule}\\
\ \bE_*(\bS)\text{-bimodule}
\end{cases}
$
.]\\
\vspace{0.5cm}

\begingroup%%----------------------------------->>
\fontsize{9pt}{11pt}\selectfont
In view of the associativity of the operations, the arrows 
$
\begin{cases}
\ \bE^*(\bX) \otimes \bE^*(\bY) \ra \bE^*(\bX \wedge \bY)\\
\ \bE_*(\bX) \otimes \bE_*(\bY) \ra \bE_*(\bX \wedge \bY)
\end{cases}
$
pass to the quotient, thereby giving the arrows
$
\begin{cases}
\ \bE^*(\bX) \otimes_{\bE^*(\bS)} \bE^*(\bY) \ra \bE^*(\bX \wedge \bY)\\
\ \bE_*(\bX) \otimes_{\bE_*(\bS)} \bE_*(\bY) \ra \bE_*(\bX \wedge \bY)
\end{cases}
.
$
\\
\endgroup %%------------------------------------<<

\begin{proposition} \ %10
Suppose that \bE is a ring spectrum.  Let 
$
\begin{cases}
\ \bX\\
\ \bY
\end{cases}
$
be spectra.  Assume: Either $\bE_*(\bX)$, as a graded right $\bE_*(\bS)$-module, is flat or $\bE_*(\bY)$, as a graded left 
$\bE_*(\bS)$-module, is flat $-$then the arrow 
$\bE_*(\bX) \otimes_{\bE_*(\bS)} \bE_*(\bY) \ra \bE_*(\bX \wedge \bY)$ is an isomorphism.
\end{proposition}

[The situation being symmetric, take \bY fixed and $\bE_*(\bY)$ flat $-$then the arrow\\ 
$\bE_*(-) \otimes_{\bE_*(\bS)} \bE_*(\bY)$ $\ra$ 
$\bE_*(-\wedge \bY)$ is a morphism of $\Z$-graded homology theories on 
\bSPEC.  But 
$\bE_*(\bS) \otimes_{\bE_*(\bS)} \bE_*(\bY) \approx$ 
$\bE_*(\bS \wedge \bY)$.]\\

\begingroup%%----------------------------------->>
\fontsize{9pt}{11pt}\selectfont
\textbf{\small FACT}  \  
Let \bE be a ring spectrum, \bF an \bE-module.  
Assume $\pi_*(\bF)$, as a graded left $\pi_*(\bE)$-module, is flat $-$then 
$\forall \ \bX$, the arrow $\bE_*(\bX) \otimes_{\pi_*(\bE)} \pi_*(\bF) \ra \bF_*(\bX)$ is an isomorphism.\\
\endgroup %%------------------------------------<<

\begingroup%%----------------------------------->>
\fontsize{9pt}{11pt}\selectfont
Notation: Given an abelian group $\pi$, put 
$\bH_*(\bX,\pi) = \bH(\pi)_*(\bX)$ and 
$\bH^*(\bX,\pi) = \bH(\pi)^*(\bX)$.\\
\endgroup %%------------------------------------<<

\begingroup%%----------------------------------->>
\fontsize{9pt}{11pt}\selectfont
\textbf{\small EXAMPLE}  \ 
Let \mA be a PID, \mM an $A$-module $-$then $\forall \ \bX$, there is an exact sequence 
$0 \ra$ 
$H_n(\bX,A) \otimes_A M \ra$ 
$\bH_n(\bX;M) \ra$ 
$\Tor^A(\bH_{n-1}(\bX;A),M) \ra 0$.
\vspi
[Since \mA is a PID, the projective dimension of \mM is $\leq 1$, so $\exists$ an exact sequence 
$0 \ra Q \ra P \ra M \ra 0$, where \mP and \mQ are projective, hence flat.  
Applying the above result then gives 
$\bH_*(\bX;A) \otimes_A P \approx \bH_*(\bX;P)$
%%----------------------------------------------------------------------------------------------10
and 
$\bH_*(\bX;A) \otimes_A Q \approx \bH_*(\bX;Q)$.  On the other hand, the exact triangle 
$\bH(Q) \ra$ 
$\bH(P) \ra$ 
$\bH(M) \ra$ 
$\Sigma \bH(Q)$ leads to an exact sequence 
$\bH_n(\bX;Q) \ra$
$\bH_n(\bX;P) \ra$
$\bH_n(\bX;M) \ra$
$\bH_{n-1}(\bX;Q) \ra$
$\bH_{n-1}(\bX;P)$ .]
\vspi
[Note: \ Under the same hypotheses, there is an exact sequence 
$0 \ra$ 
$\Ext_A(\bH_{n-1}(\bX;A),M) \ra$ 
$\bH^n(\bX;M) \ra$ 
$\Hom_A(\bH_n(\bX;A),M) \ra 0$.]\\
\endgroup %%------------------------------------<<

\begingroup%%----------------------------------->>
\fontsize{9pt}{11pt}\selectfont
\textbf{\small FACT}  \  
Suppose that \mA is a PID $-$then $\forall \ \bX$, 
$\bX \wedge \bH(A) \approx \ds\bigvee\limits_n \Sigma^n \bH(G_n)$, where 
$G_n = \bH_n(\bX;A)$.
\vspi
[Here 
$\ds\bigvee\limits_n \Sigma^n \bH(G_n) \approx \ds\prod\limits_n \Sigma^n \bH(G_n)$ 
(cf. p. \pageref{17.16} ff.), thus it suffices to specify arrows 
$\bff_n:\bX \wedge \bH(A) \ra$ 
$\Sigma^n \bH(G_n)$ such that $\pi_n(\bff_n)$ is an isomorphism $\forall \ n$.]\\
\endgroup %%------------------------------------<<

\begingroup%%----------------------------------->>
\fontsize{9pt}{11pt}\selectfont
\textbf{\small EXAMPLE}  \ 
Let \mA be a PID $-$then $\forall \ \bX, \bY$, $\&$ $\forall, i, j$, there is an exact sequence 
$0 \ra$ 
$\bH_i(\bX;A) \otimes_A \bH_j(\bY;A) \ra$ 
$\bH_i(\bX;\bH_j(\bY;A)) \ra$ 
$\Tor^A(\bH_{i-1}(\bX;A),\bH_j(\bY;A)) \ra 0$.  
Now sum over all $(i,j)$: $i + j = k$.  
Setting aside the flanking terms and putting $G_j = \bH_j(\bY;A)$, the middle term assumes the form 
$\ds\bigoplus\limits_{i + j = k} \bH_i(\bX;\bH_j(\bY;A)) =$  \ 
$\ds\bigoplus\limits_{j} \pi_k(\bX \wedge \Sigma^j \bH(G_j)) =$ \ 
$\pi_k(\bX \wedge \ds\bigvee\limits_j \Sigma^j \bH(G_j)) =$ \ 
$\pi_k(\bX \wedge \bY \wedge \bH(A)) =$ \ 
$\bH_k(\bX \wedge \bY;A)$.\\
\endgroup %%------------------------------------<<

In a category \bC with pushouts, one has the notion of an internal cocategory (or a cocategory object) 
(cf. p. \pageref{17.17}), which can be specialized to the notion of an internal cogroupoid (or a cogroupoid object).  
Definition: Let \bk be a commutative ring with unit $-$then a
\un{graded Hopf algebroid over \bk}
\index{graded Hopf algebroid over \bk} 
is a cogroupoid object in the category of graded commutative \bk-algebras with unit.  So, a graded Hopf algebroid over \bk consists of a pair $(A,\Gamma)$ of graded commutative \bk-algebras with unit and morphisms 
$\eta_R:A \ra \Gamma$ (right unit = ``cosource''), 
$\eta_L:A \ra \Gamma$ (left unit = ``cotarget''), 
$\epsilon:\Gamma \ra A$ (augmentation =``coidentity''), 
$\Delta:\Gamma \ra \Gamma \otimes_A \Gamma$ (diagonal =``cocomposition''),
$c:\Gamma \ra \Gamma$ (conjugation =``coinversion'') 
satisfying the dual of the usual category theoretic relations (cf infra).  
Therefore $(A,\Gamma)$ attaches to a graded commutative \bk-algebra \mT with unit a groupoid $\bG_T$, where 
$\Ob\bG_T = \Hom(A,T)$ and 
$\Mor\bG_T = \Hom(\Gamma,T)$.  
Example: $(\bk,\bk)$ is a graded Hopf algebroid over \bk (trivial grading).

[Note: \ When $A = \bk$ and $\eta_L = \eta_R$, $\Gamma$ is a graded commutative Hopf algebra over \bk or still, a cogroup object in the category of graded commutative \bk-algebras with unit.]

Remark: Graded Hopf algebroids over \bk can be organized into a (large) double category 
(Borceux\footnote[2]{\textit{Handbook of Categorical Algebra 1}, Cambridge University Press (1994), 327-328.}).\\

\begingroup%%----------------------------------->>
\fontsize{9pt}{11pt}\selectfont
There is a coequalizer diagram 
\begin{tikzcd}%[sep=small]
{\Gamma \otimes_\bk A \otimes_\bk \Gamma}  
\arrow[rr,swap,"\id_\Gamma \otimes \eta_L" , shift right = 1]
\arrow[rr,"\eta_R \otimes \id_\Gamma", shift left=1]
&&{\Gamma \otimes_\bk \Gamma} 
\end{tikzcd}
$\ra \Gamma \otimes_A \Gamma$ 
and 
\begin{tikzcd}[sep=large]
{A} \ar{d}[swap]{\eta_L} \ar{r}{\eta_R} &{\Gamma} \ar{d}{\ini_R}\\
{\Gamma} \ar{r}[swap]{\ini_L} &{\Gamma \otimes_A \Gamma}
\end{tikzcd}
is a pushout square.
\vspi
%%----------------------------------------------------------------------------------------------11
[Note: \ Tacitly, one uses $\eta_R$ to equip $\Gamma$ with the structure of a graded right $A$-module and $\eta_L$ to equip $\Gamma$ with the structure of a graded left $A$-module.]
\vspi
As for $\eta_R$, $\eta_L$, $\epsilon$, $\Delta$, and $c$, they must have the following properties: 
$\epsilon \circ \eta_R = \id_A = \epsilon \circ \eta_L$, 
$\Delta \circ \eta_R =$ $\ini_L \circ \eta_R$, 
$\Delta \circ \eta_L =$ $\ini_R \circ \eta_L$, 
$(\id_\Gamma \otimes \epsilon) \circ \Delta = \id_\Gamma$, 
$(\epsilon \otimes \id_\Gamma ) \circ \Delta = \id_\Gamma$,  
$(\id_\Gamma \otimes \Delta) \circ \Delta = (\Delta \otimes \id_\Gamma) \circ \Delta$, 
$c \circ \eta_R = \eta_L$, 
$c \circ \eta_L = \eta_R$, 
$(c \otimes \id_\Gamma) \circ \Delta = \eta_R \circ \epsilon$, and 
$(\id_\Gamma \otimes c) \circ \Delta = \eta_L \circ \epsilon$.
\vspi
[Note: \ The formulas relating $c$ to the other arrows are the duals of those on 
p. \pageref{17.18} (the role of $\chi$ in the groupoid object situation is played here by $c$).  
Corollaries: 
(1) $c \circ c = \id_\Gamma$; 
(2) $\epsilon \circ c = \epsilon$.]\\
\endgroup %%------------------------------------<<

\begingroup%%----------------------------------->>
\fontsize{9pt}{11pt}\selectfont
\textbf{\small EXAMPLE}  \ 
The dual of the mod $p$ Steenrod algebra is isomorphic to 
$\bH(\F_p)_*(\bH(\F_p))$, a graded commutative Hopf algebra over $\bF_p$.  One has 
$\bH(\F_2)_*(\bH(\F_2)) \approx$ $\F_2[\xi_1, \xi_2, \ldots]$, where \ $\abs{\xi_k} = 2^k - 1$ and 
$\Delta(\xi_k) =$ $\ds\sum\limits_{i = 0}^k \xi_{k-i}^{2^i} \otimes \xi_i$, \ and for $p > 2$, \ 
$\bH(\F_p)_*(\bH(\F_p)) \approx$ 
$\F_p[\xi_1, \xi_2, \ldots] \otimes_{\F_p} \ds\bigwedge (\tau_0, \tau_1, \ldots)$, where 
$\abs{\xi_k} = 2(p^k - 1)$,
$\abs{\tau_k} = 2p^k - 1$
and \ 
$\Delta(\xi_k) =$ $\ds\sum\limits_{i = 0}^k \xi_{k-i}^{p^i} \otimes \xi_i$, 
$\Delta(\tau_k) =$ $\tau_k \otimes 1 + \ds\sum\limits_{i = 0}^k \xi_{k-i}^{p^i} \otimes \tau_i$.  
The unit and augmentation are isomorphisms in degree 0 and the conjugation $c$ is given recursively by 
$\sum\limits_{i = 0}^k \xi_{k-i}^{p^i} c(\xi_i) = 0$ $(k > 0)$ and 
$\tau_k + \ds\sum\limits_{i = 0}^k \xi_{k-i}^{p^i} c(\tau_i) = 0$ $(k \geq 0)$.
\vspi
[Note: \ In the above, it is understood that $\xi_0 = 1$.]\\
\endgroup %%------------------------------------<<

\begin{proposition} \ %11
Suppose that \bE is a ring spectrum.  Assume: \bE is commutative and $\bE_*(\bE)$, as a graded right 
$\bE_*(\bS)$-module, is flat $-$then the pair 
$(\bE_*(\bS),\bE_*(\bE))$ is a graded Hopf algebroid over $\Z$.
\end{proposition}

[$\bE_*(\bE)$ is a graded commutative $\Z$-algebra with unit.  
Proof: The product is defined by 
$\bE_*(\bE) \otimes \bE_*(\bE) \ra$ 
$(\bE \wedge \bE)_*(\bE \wedge \bE) \ra$
$\bE_*(\bE \wedge \bE) \ra$ 
$\bE_*(\bE)$ and the unit 
$\Z \ra \bE_0(\bE)$ is defined by sending 1 to the arrow 
$\bS = \bS \wedge \bS$ $\ra$ $\bE \wedge \bE$.  This said, let
$
\begin{cases}
\ \eta_R:\bE_*(\bS) \approx \pi_*(\bS \wedge \bE) \ra \pi_*(\bE \wedge \bE) = \bE_*(\bE)\\
\ \eta_L:\bE_*(\bS) \approx \pi_*(\bE \wedge \bS) \ra \pi_*(\bE \wedge \bE) = \bE_*(\bE)
\end{cases}
$
and 
$\epsilon:\bE_*(\bE) =$ 
$\pi_*(\bE \wedge \bE) \ra$ 
$\pi_*(\bE) =$ 
$\bE_*(\bS)$.  
Next, take for $\Delta$ the composite 
$\bE_*(\bE) =$ 
$\pi_*(\bE \wedge \bE) \approx$ 
$\pi_*(\bE \wedge \bS \wedge \bE) \ra$ 
$\pi_*(\bE \wedge \bE \wedge \bE) \ra$ 
$\bE_*(\bE \wedge \bE) \approx$ 
$\bE_*(\bE)  \otimes_{\bE_*(\bS)} \bE_*(\bE)$ (cf. Proposition 10).  
Finally, 
$c:\bE_*(\bE) =$
$\pi_*(\bE \wedge \bE) \ra$ 
$\pi_*(\bE \wedge \bE) \ra$ 
$\bE_*(\bE)$ is induced by the interchange 
$\Tee:\bE \wedge \bE \ra$ $\bE \wedge \bE$.]

[Note: \ Due to the presence of $c$ and the relations 
$
\begin{cases}
\ c \circ \eta_R = \eta_L\\
\ c \circ \eta_L = \eta_R
\end{cases}
,
$
$\bE_*(\bE)$, as a graded right $\bE_*(\bS)$-module, is flat iff $\bE_*(\bE)$, as a graded left $\bE_*(\bS)$-module, is flat 
(the $\bE_*(\bS)$-module structures on $\bE_*(\bE)$ per $\eta_R$ and $\eta_L$ are the same as those introduced on 
 p. \pageref{17.19}).  
 Example: The flatness assumption is met if 
 $\bE \wedge \bE \approx$ $\ds\bigvee\limits_i \Sigma^{n_i} \bE$ (isomorphism of $\bE$-modules) (for then 
 $\pi_*(\bE \wedge \bE) \approx$ $\ds\bigoplus\limits_i  \pi_{* - n_i}(\bE)$, thus is a graded free $\pi_*(\bE)$-module.]\\

\begingroup%%----------------------------------->>
\fontsize{9pt}{11pt}\selectfont
Tied to the definitions are various diagrams and a complete proof of Proposition 11 entails checking
%%----------------------------------------------------------------------------------------------12
that these diagrams commute, which is straighforward if tedious (a discussion can be found in 
Adams\footnote[2]{\textit{SLN} {99} (1969), 56-71.}).\\
\endgroup %%------------------------------------<<

\begingroup%%----------------------------------->>
\fontsize{9pt}{11pt}\selectfont
\textbf{\small EXAMPLE}  \ 
$\bKU_*(\bKU)$ is a graded free $\bKU_*(\bS)$-module 
(Adams-Clarke\footnote[3]{\textit{Illinois J. Math.} \textbf{21} (1977), 826-829.}), 
thus the hypotheses of Proposition 11 are met in this case.
\vspi
[Note: \ The structure of $\bKU_*(\bKU)$ had been worked out by 
Adams-Harris-Switzer\footnote[6]{\textit{Proc. London Math. Soc.} \textbf{23} (1971), 385-408.}.]\\
\endgroup %%------------------------------------<<

Given a graded Hopf algebroid $(A,\Gamma)$ over \bk, a (left) 
\un{$(A,\Gamma)$-comodule}
\index{A,$\Gamma$-comodule} 
is a graded left $A$-module \mM equipped with a morphism $M \ra \Gamma \otimes_A M$ of graded left $A$-modules such that 
\begin{tikzcd}%[sep=small]
{M} \ar{d} \ar{r} &{\Gamma \otimes_A M} \ar{d}\\
{\Gamma \otimes_A M} \ar{r} &{\Gamma \otimes_A \Gamma \otimes_A M}
\end{tikzcd}
\ 
and 
\ 
\begin{tikzcd}%[sep=small]
{M} \ar{rd} \ar{r} &{\Gamma \otimes_A M} \ar{d}\\
&{A \otimes_A M}
\end{tikzcd}
commute.\\
\vspace{0.5cm}

\begin{proposition} \ 
Suppose that \bE is a ring spectrum.  
Assume \bE is commutative and $\bE_*(\bE)$, as a graded right $\bE_*(\bS)$-module, is flat 
$-$then $\forall$ $\bX$, 
$\bE_*(\bX)$ is an $(\bE_*(\bS),\bE_*(\bE))$-comodule.
\end{proposition}

[The arrow 
$\bE_*(\bX) \ra \bE_*(\bE) \otimes_{\bE_*(\bS)} \bE_*(\bX)$ is the composite 
$\bE_*(\bX)  = \pi_*(\bE \wedge \bX) \approx$ 
$\pi_*(\bE \wedge \bS \wedge \bX) =$ 
$\pi_*(\bE \wedge \bE \wedge \bX) =$ 
$\bE_*(\bE \wedge \bX) \approx$ 
$\bE_*(\bE) \otimes_{\bE_*(\bS)} \bE_*(\bX)$ (cf. Proposition 10).]\\

Rappel: A spectrum \bE defines a $\Z$-graded cohomology theory $E^*$ on $\bCW_*$ (cf. Proposition 4) and 
$\forall \ X$ in $\bCW_*$, $E^n(X_+) \approx E^n(X) \oplus E^n(\bS^0) $.

[Note: \ When \bE is a ring spectrum, there is a cup product $\cup$, viz. the composite 
$E^*(X) \otimes E^*(X) \ra$
$E^*(X \#_k X) \ra$ 
$E^*(X)$, where 
$X \ra X \#_k X$ is the reduced diagonal.  Therefore $E^*(X)$ is a graded ring and $E^*(X_+)$ is a graded ring with unit (both are graded and commutative if \bE is commutative).]

Let \bE be a commutative ring spectrum $-$then \bE is said to be 
\un{complex orientable}
\index{complex orientable (commutative ring spectrum)}
if $\exists$ an element $x_{\bE} \in E^2(\bP^\infty(\C))$ with the property that the arrow of restriction 
$E^2(\bP^\infty(\C)) \ra$ 
$E^2(\bP^1(\C)) \approx$ 
$\pi_0(\bE)$ sends $x_{\bE}$ to the unit 
$\bS \ra \bE$ of \bE.  One calls $x_{\bE}$ a 
\un{complex orientation}
\index{complex orientation (commutative ring spectrum)} 
of \bE.

[Note: \ 
$\pi_0(\bE) = [\bS,\bE] \approx [\bS^0,E_0] \approx [\bS^0,\Omega^2 E_2] \approx [\Sigma^2\bS^0, E_2] =$
$E^2(\bS^2)$ and 
$\bS^2 \approx$ 
$\bP^1(\C)$.]

%%----------------------------------------------------------------------------------------------13 %dup
Remark: Identify 
$\pi_0(\bE) = [\bS,\bE]$ with 
$[\bQ_{2n}^\infty\bS^{2n},\bE] \approx [\bS^{2n},E_{2n}]$ and let 
top: $\bP^n(\C) \ra$ 
$\bS^{2n}$ $(=\bP^n(\C)/\bP^{n-1}(\C))$ be the top cell map $-$then the arrow of restriction 
%%----------------------------------------------------------------------------------------------13
$E^{2n}(\bP^\infty(\C)) \ra$
$E^{2n}(\bP^n(\C))$ sends $x_{\bE}^n$ to the image of the unit of \bE under the precomposition arrow 
\begin{tikzcd}%[sep=small]
{[\bS^{2n},E_{2n}]} \ar{r}{\text{top}^*} &{[\bP^n(\C),E_{2n}].}
\end{tikzcd}

[The diagram 
\begin{tikzcd}%[sep=small]
{\bP^n(\C)} \ar{d}[swap]{\text{top}} \ar{r} &{\bP^n(\C) \#_k \cdots \#_k \bP^n(\C)} \\
{\bS^{2n}} \ar{r} &{\bS^2 \#_k \cdots \#_k \bS^2}\ar{u}
\end{tikzcd}
is pointed homotopy commutative.]

Example: Let \mA be a commutative ring with unit $-$then $\bH(A)$ is complex orientable.

[Recall that $H^*(\bP^\infty(\C);A) \approx A[x]$, $\abs{x} = 2$.]\\

\begin{proposition} \ %13
Suppose that \bE is a commutative ring spectrum.  Assume: \bE is complex orientable with complex orientation $x_{\bE}$ $-$then 
$E^*(\bP^\infty(\C)_+) \approx \bE^*(\bS)[[x_{\bE}]]$.
\end{proposition}

[Note: \ $\bE^*(\bS)[[x_{\bE}]]$ is the graded $\bE^*(\bS)$-algebra of formal power series in $x_{\bE}$ 
$(\abs{x_{\bE}} = 2)$.  So: A typical element in $\bE^q(\bS)[[x_{\bE}]]$ has the form 
$\sum\limits_{i=0}^\infty \lambda_i x_{\bE}^i$, where $\lambda_i \in \bE^{q-2i}(\bS)$.]\\

\begin{proposition} \ %14
Suppose that \bE is a commutative ring spectrum.  Assume: \bE is complex orientable with complex orientation $x_{\bE}$ $-$then 
$E^*((\bP^\infty(\C) \times_k \bP^\infty(\C))_+) \approx$ 
$\bE^*(\bS)[[x_{\bE} \otimes 1,1 \otimes x_{\bE}]]$.\\
\end{proposition}

\begingroup%%----------------------------------->>
\fontsize{9pt}{11pt}\selectfont
Cole\footnote[2]{Ph.D. Thesis, University of Chicago, Chicago (1996).} has given a proof of these propositions that does not involve the Atiyah-Hirzebruch spectral sequence.
\vspi
[Note: \ The method is to show from first principles that there are splittings 
$\bE \wedge \bP^n(\bC) = \ds\bigvee\limits_{i=1}^n \Sigma^{2i}\bE$, 
$\texttt{HOM}(\bP^n(\C),\bE) \approx \ds\prod\limits_{i=1}^n  \Omega^{2i}\bE$  in \bE-\bMOD.]\\
\endgroup %%------------------------------------<<

\begingroup%%----------------------------------->>
\fontsize{9pt}{11pt}\selectfont
\textbf{\small EXAMPLE}  \ 
If \bE is complex orientable, then $E_*(\bP^\infty(\bC)_+)$ is a graded free $\bE_*(\bS)$-module and 
$E_*(\bP^\infty(\C)_+) \otimes_{\bE_*(\bS)} E_*(\bP^\infty(\C)_+) \approx$ 
$E_*(\bP^\infty(\C)_+ \#_k (\bP^\infty(\C)_+)$ (cf. Proposition 10).\\
\endgroup %%------------------------------------<<

The standard reference for the theory of formal groups is 
Hazewinkel\footnote[2]{\textit{Formal Groups and Applications}, Academic Press (1978); 
see also Ravenel, \textit{Complex Cobordism and Stable Homotopy Groups of Spheres}, Academic Press (1986), 354-379.}.  
There the reader can look up the proofs but to establish notation, I shall review some of the definitions.

Let \mA be a graded commutative ring with unit.  Consider $A[[x,y]]$, where 
$
\begin{cases}
\ \abs{x} = 2\\
\ \abs{y} = 2
\end{cases}
-
$
then a \un{formal group law}
\index{formal group law} 
(FGL) 
\index{FGL} 
\un{over \mA} is an element 
$F(x,y) \in A[[x,y]]$ of the form $x + y + 
%%----------------------------------------------------------------------------------------------14
\sum\limits_{i,j \geq 1} a_{ij}x^i y^j$, where $a_{ij} \in A_{2 - 2i - 2j}$, such that 
$F(x,F(y,z)) = F(F(x,y),z)$ (associativity) and 
$F(x,y) = F(y,x)$ (commutativity).

[Note: \ In algebra, one does not usually work in the graded setting, the standing assumption being that \mA is a commutative ring with unit (as, e.g., in Hazewinkel).  
Of course, if \mA is a graded commutative ring with unit, then 
$A_{\text{even}} (= \bigoplus\limits_n A_{2n})$ is a commutative ring with unit and every FGL over \mA is a FGL over 
$A_{\text{even}}$.  
Example: 
$F(x,y) = x + y + uxy$ $(u \in A_{-2})$ is a FGL over \mA, hence over $A_{\text{even}}$, while
$F(x,y) = x + y + xy$ is not a FGL over \mA (but is a FGL over $A_{\text{even}}$).]

Notation: Write 
$F(x,y) = x +_F y$, so 
$
\begin{cases}
\ x +_F 0 = x\\
\ 0 +_F y = y
\end{cases}
,
$
$x +_F (y +_F z) = (x +_F y)  +_F z$, and 
$x +_F y = y +_F x$.

Definition: An element 
$\phi(x) = \ds\sum\limits_{i \geq 1} \phi_i x^i \in A[[x]]$ $(\abs{x} = 2)$ is said to be 
\un{homogeneous}
\index{homogeneous (element in $A[[x]]$)} 
if $\phi_i \in A_{2 - 2i}$ $\forall \ i$.\\

\begingroup%%----------------------------------->>
\fontsize{9pt}{11pt}\selectfont
\textbf{\small FACT}  \  
If $F(x,y)$ is a FGL over \mA, then there is a unique homogeneous element $\iota(x) \in A[[x]]$ such that 
$x +_F \iota(x) = 0 = \iota(x) +_F x$.
\vspi
[There exist unique homogeneous elements
$
\begin{cases}
\ \iota_L(x)\\
\ \iota_R(x)
\end{cases}
\in A[[x]]
$ such that 
$
\begin{cases}
\ \iota_L(x) +_F x = 0\\
\ x +_F \iota_R(x) = 0
\end{cases}
$, thus 
$\iota_L(x)  =$ 
$\iota_L(x)  +_F 0 =$ 
$\iota_L(x)  +_F (x +_F \iota_R(x)) =$ 
$(\iota_L(x) +_F x) +_F \iota_R(x) =$
$0 +_F \iota_R(x) =$
$\iota_R(x)$ and one can take 
$\iota(x) =$ 
$\iota_L(x)  =$ 
$\iota_R(x)$.]\\
\endgroup %%------------------------------------<<

\begin{proposition} \ %15
Let 
$m:(\bP^\infty(\C) \times_k \bP^\infty(\C))_+ \ra \bP^\infty(\C)_+$ be the multiplication classifying the tensor product of complex line bundles $-$then $\forall$ complex orientable \bE, 
$F_{\bE} = m^*(x_{\bE})$ is a FGL over $\bE^*(\bS)$.\\
\end{proposition}

Example: The FGL attached to $\bH(\bk)$ by Proposition 15, where \bk is a commutative ring with unit, is the ``additive'' FGL, viz. 
$x + y$.\\

\label{17.27} %dmc mnft
\begingroup%%----------------------------------->>
\fontsize{9pt}{11pt}\selectfont
\textbf{\small EXAMPLE}  \ 
\bKU is complex orientable and the associated FGL is $x + y + \bb_{\bU}xy$ (cf. p. \pageref{17.20}).\\
\endgroup %%------------------------------------<<

Let \mA be a graded commutative ring with unit.  Suppose that \mF, \mG are formal group laws over \mA $-$then a 
\un{homomorphism}
\index{homomorphism (of formal group laws)} 
$\phi:F \ra G$ is a homogeneous element $\phi \in A[[x]]$ such that 
$\phi(x +_F y) = \phi(x) +_G \phi(y)$, i.e., 
$\phi(F(x,y)) = G(\phi(x),\phi(y))$.  
A homomorphism 
$\phi:F \ra G$ is an 
\un{isomorphism}
\index{isomorphism (of formal group laws)} 
if $\phi^\prime(0)$ (the coefficient of $x$) belongs to $A_0^\times$.  
An isomorphism 
$\phi:F \ra G$ is a 
\un{strict isomorphism}
\index{strict isomorphism (of formal group laws)}  
if $\phi^\prime(0) = 1$.

[Note: \ A homomorphism $\phi:F \ra G$ is an isomorphism iff $\exists$ a homomorphism $\psi:G \ra F$ such that 
$\phi(\psi(x)) = x = \psi(\phi(x))$.]\\

%%----------------------------------------------------------------------------------------------15
$\FGL_A$ is the set of formal group laws over \mA and $\bFGL_A$ is the category whose objects are the elements of 
$\FGL_A$ and whose morphisms are the homomorphisms.

[Note: \ If $f:A \ra A^\prime$ is a homomorphism of graded commutative rings with unit, then $f$ induces a functor 
$f_*:\bFGL_A \ra \bFGL_{A^\prime}$ (on objects, 
$f_* F(x,y) = x + y + \sum\limits_{i,j \geq 1} f(a_{i j}) x^i y^j$, and on morphisms 
$f_* \phi(x) = \sum\limits_{i \geq 1} f(\phi_i) x^i$).]\\

\begingroup%%----------------------------------->>
\fontsize{9pt}{11pt}\selectfont
\textbf{\small FACT}  \  
If \bE is complex orientable and if $x_{\bE}^\prime$, $x_{\bE}\pp
$ are two complex orientations of \bE, then the associated formal group laws 
$F_{\bE}^\prime$, $F_{\bE}\pp$ over $\bE^*(\bS)$ are strictly isomorphic.\\
\endgroup %%------------------------------------<<

Let \mA be a graded commutative ring with unit.  Write $\IPS_A$ for the set of homogeneous elements $\phi$ in $A[[x]]$ such that $\phi^\prime(0) = 1$ $-$then $\IPS_A$ is a group under composition, functorially in $A$.

Notation: $B = \Z[b_1, b_2, \ldots]$, where $\abs{b_i} = -2i$.\\

\begin{proposition} \ %16
\mB is a graded Hopf algebra over $\bZ$.
\end{proposition}

[In fact, $\Hom(B,A) \approx \IPS_A$, so \mB is a cogroup object in the category of graded commutative rings 
with unit.]\\

Remark: $\IPS_A$ operates to the left on $\bFGL_A$, viz. $(\phi,F) \ra \phi \cdot F = F_\phi$, where 
$F_\phi(x,y) = \phi(F(\phi^{-1}(x),\phi^{-1}(y)))$.\\

\begingroup%%----------------------------------->>
\fontsize{9pt}{11pt}\selectfont
Let \mA be a graded commutative ring with unit $-$then \mA is said to be 
\un{graded coherent}
\index{graded coherent (graded commutative ring with unit)} 
if each finitely generated graded ideal of \mA is 
finitely presented.  
Example: \mA graded noetherian $\implies$ \mA graded coherent.
\vspi
[Note: \ $\pi_*(\bS)$ is not graded coherent 
(Cohen\footnote[2]{\textit{Comment. Math. Helv.} \textbf{44} (1969), 217-228.}).]
\vspi
Remark: \ Suppose that \mA is graded coherent $-$then a finitely generated graded $A$-module \mM is finitely presented iff it and its finitely generated graded  submodules are finitely presented.\\
\endgroup %%------------------------------------<<

\begingroup%%----------------------------------->>
\fontsize{9pt}{11pt}\selectfont
\textbf{\small EXAMPLE}  \ 
Let \bk be a commutative ring with unit.  Consider $\bk[x_1, x_2, \ldots]$, where $\abs{x_i} = -2i$ $-$then 
$\bk[x_1, x_2, \ldots]$ is not graded noetherian but is graded coherent provided that \bk is noetherian.\\
\endgroup %%------------------------------------<<

\index{Theorem: Lazard's Theorem}
\index{Lazard's Theorem}
\textbf{\small LAZARD'S THEOREM} \quad 
The functor from the category of graded commutative rings with unit to the category of sets which sends \mA to $\FGL_A$ is representable.  Accordingly, there is a graded commutative ring \mL with unit and a FGL $F_L$ over \mL such that $\forall \ A$ and 
$\forall \ F \in \FGL_A$, $\exists$! $f \in \Hom(L,A)$ : $f_*F_L = F$.

%%----------------------------------------------------------------------------------------------16
[Note: \ The structure of \mL can be determined, viz. $L = \Z[x_1, x_2, \ldots]$, where $\abs{x_i} = -2i$, hence \mL is graded coherent (cf. supra).]\\

\begingroup%%----------------------------------->>
\fontsize{9pt}{11pt}\selectfont
The mere existence of \mL is a formality.  
Thus fix indeterminates $t_{i j}$ of degree $2 - 2i - 2j$ and put 
$\mu(x,y) = x + y + \ds\sum\limits_{i,j \geq 1} t_{i j} x^i y^j$.  
Define homogeneous polynomials $p_{i j k}$ in the $t_{i j}$ by writing 
$\mu(x,\mu(y,z))$ 
$-$ 
$\mu(\mu(x,y),z) =  \ds\sum\limits_{i,j,k \geq 1} p_{i j k}x^i y^j z^k$ $-$then 
$L = \Z[t_{i j}: i, j \geq 1]/I$, where \mI is the graded ideal generated by the 
$t_{i j} - t_{j i}$ and the $p_{i j k}$, and $\mu$ induces a FGL $F_L$ over \mL having the universal property in question.
\vspi
Determining the structure of \mL is more difficult and depends in part on the following construction.  Fix indeterminates 
$b_i$ of degree $-2i$ and consider, as above, $B = \Z[b_1, b_2, \ldots]$.  Let 
$\exp x = x + \ds\sum\limits_{i \geq 1} b_i x^{i+1} \in$ $B[[x]]$ $(\abs{x} = 2)$ and let $\log x$ be its inverse (so 
$\exp(\log x) = x = \log(\exp x))$ $-$then 
$F_B(x,y) =$ $\exp(\log x + \log y)$ is a FGL over \mB and the homomorphism $L \ra B$ classifying $F_B$ is injective .\\
\endgroup %%------------------------------------<<

\begingroup%%----------------------------------->>
\fontsize{9pt}{11pt}\selectfont
\textbf{\small FACT}  \  
If $A \ra A^\prime$ is a surjective map of graded commutative rings with unit, then any FGL over $A^\prime$ lifts to a FGL over \mA.\\
\endgroup %%------------------------------------<<

Put 
$LB = L[b_1, b_2, \ldots]$, where $b_i$ is an indeterminate of degree $-2i$ ($\implies$ 
$LB =$ 
$L \otimes_Z \Z[b_1, b_2, \ldots] =$ $L \otimes_Z B$).\\

\begin{proposition} \ %17
The pair $(L,LB)$ is a graded Hopf algebroid over $\Z$.
\end{proposition}

[Let \mA be a graded commutative ring with unit.  
Denoting by $\bG_A$ the groupoid whose objects are the formal group laws over \mA and whose morphisms are the strict isomorphisms, the functor from the category of graded commutative rings with unit to the category of groupoids which sends \mA to $\bG_A^\OP$ is represented by $(L,LB)$.  For Lazard gives 
$\Hom(L,A) \leftrightarrow$ 
$\FGL_A =$
$\Ob\bG_A$ $(= \Ob(\bG_A^\OP)$ and this identifies the objects.  
Turning to the morphisms, suppose that $f \in \Hom(LB,A)$.
Put $F = (\restr{f}{L})_*F_L$ and 
$\phi(x) =$ 
$x + \sum\limits_{i \geq 1} f(b_i) x^{i+1}$ $-$then
$\phi^\OP:G \ra F$ is a strict isomoprhism, where 
$G(x,y) =$ 
$\phi(F(\phi^{-1}(x),\phi^{-1}(y)))$.]

[Note: \ $\eta_L$ is the inclusion $L \ra LB$ but there is no simple explicit formula for $\eta_R$.  However, using definitions only, one can write down explicit formulas for $\epsilon$, $\Delta$, and $c$.]\\

\begingroup%%----------------------------------->>
\fontsize{9pt}{11pt}\selectfont
A groupoid \bG is said to be 
\un{split}
\index{split (groupoid)} 
if there exists a group \mG and a left $G$-set \mY such that \bG is isomorphic to 
$\tran Y$, the translation category of \mY (cf. p. \pageref{17.21}).
\vspi
Example: Take $G =$ IPS$_A$, $Y = $FGL$_A$ $-$then the translation category of $\FGL_A$ is isomorphic to $\bG_A$, i.e., $\bG_A$ is split.\\
\endgroup %%------------------------------------<<

I shall now review the theory of \bMU, referring the reader to 
Adams\footnote[2]{\textit{Stable Homotopy and Generalized Homology}, University of Chicago (1974), 32-93.}
for the details and further information.\\

%%----------------------------------------------------------------------------------------------17
Let $\bG_n(\C^\infty)$ be the grassmannian  of complex $n$-dimensional subpaces of $\C^\infty$, $\gamma_n$ the canonical complex $n$-plane bundle over $\bG_n(\C^\infty)$.  Put $MU(n) = T(\gamma_n)$, the Thom space of $\gamma_n$ $-$then 
$i^*(\gamma_{n+1}) = \gamma_n \oplus \un{\C}(\bG_n(\C^\infty) \overset{i}{\ra} \bG_{n+1}(\C^\infty))$ and 
$T(\gamma_n \oplus \un{\C}) \approx$ 
$\Sigma^2 T(\gamma_n) =$ 
$\Sigma^2 MU(n)$, so there is an arrow 
$\Sigma^2 MU(n) \ra$ 
$MU(n+1)$.  The prescription 
$X_{2n} = MU(n)$, 
$X_{2n+1} =$ $\Sigma MU(n)$ thus defines a separated prespectrum \bX and by definition, $\bMU = e\bX$.\\

\begingroup%%----------------------------------->>
\fontsize{9pt}{11pt}\selectfont
\textbf{\small EXAMPLE}  \ 
\bMU and \bKU are connected by the fact that the arrow 
$\bMU_*(\bX) \otimes_{\bMU_*(\bS)} \bKU(\bS) \ra$ 
$\bKU_*(\bX)$ induced by the Todd genus is an isomorphism of graded  $\bKU_*(\bS)$-modules for all \bX 
(Conner-Floyd\footnote[3]{\textit{The Relation of Cobordism to K-Theories}, Springer Verlag (1966); 
see also Hopkins-Hovey, \textit{Math. Zeit.} \textbf{210} (1992), 181-196.}).\\
\endgroup %%------------------------------------<<

\index{Theorem: MU Theorem}
\index{MU Theorem}
\textbf{\small MU THEOREM} \quad 
\bMU is a commutative ring spectrum with complex orientation $x_{\bMU}$.  And: The map 
$L \ra \bMU^*(\bS)$ classifying $F_{\bMU}$ is an isomorphism of graded commutative rings with unit.

[Note: \ The pair 
$(\bMU_*(\bS), \bMU_*(\bMU))$ satisfies the hypotheses of Proposition 11 ($\bMU_*(\bMU)$ is a graded free $(\bMU_*(\bS)$-module), hence is a graded Hopf algebroid over $\Z$. 
As such, it is isomorphic to $(L,LB)^\OP$ (reversal of gradings).]\\

\label{17.55}
\begingroup%%----------------------------------->>
\fontsize{9pt}{11pt}\selectfont
An arrow 
$\bff:\Sigma^n \bX \ra \bX$ is said to be 
\un{composition nilpotent}
\index{composition nilpotent} 
if $\exists \ k$ such that the composite 
$\bff \circ \Sigma^n \bff \circ \cdots \circ \Sigma^{(k-1)n} \bff: \Sigma^{kn}\bX \overset{\bff^k}{\lra} \bX$ vanishes.
Example: Take \bX compact $-$then \bff is composition nilpotent iff $\bff^{-1}\bX = 0$ (cf. p. \pageref{17.22}).
\vspi
[Note: \ The same terminology is used in the category of graded abelian groups.  
Example: Take \bX compact and let \bE be a ring spectrum $-$then $\bE_*(\bff)$ is composition nilpotent iff 
$\bE \wedge \bff^{-1}\bX = 0$.]
\vspi
An arrow 
$\bff:\bX \ra \bY$ is said to be 
\un{smash nilpotent}
\index{smash nilpotent} 
if $\exists \ k$ such that the $k$-fold smash product 
$\bff^{(k)}:\bX^{(k)} \ra \bY^{(k)}$ vanishes.  
Example: $\bff:\bS \ra \bY$ is smash nilpotent iff $\bY_{\bff}^{\infty} = 0$ (cf. p. \pageref{17.23}).\\
\endgroup %%------------------------------------<<

\label{17.66}
\begingroup%%----------------------------------->>
\fontsize{9pt}{11pt}\selectfont
\index{MU Nilpotence Technology}
\textbf{\small FACT  \ (\un{\bMU Nilpotence Technology})} \  
Let \bE be a ring spectrum and consider the Hurewicz homomorphism 
$\bS_*(\bE) \ra \bMU_*(\bE)$ 
(cf. p. \pageref{17.24} ff.) $-$then the homogeneous elements of its kernel are nilpotent 
(Devinatz-Hopkins-Smith\footnote[2]{\textit{Ann. of Math.} \textbf{128} (1988), 207-241.}).\\
\endgroup %%------------------------------------<<

\begingroup%%----------------------------------->>
\fontsize{9pt}{11pt}\selectfont
Application: If \bX is compact and if $\bff:\Sigma^n \bX \ra \bX$ is an arrow such that 
$\bMU_*(\bff) = 0$, then \bff is composition nilpotent.
\vspi
[$\bMU_*(\bff) = 0$ $\implies$ 
$\bMU \wedge \bff^{-1} \bX = 0$ $\implies$ 
$\exists \ k$: 
$\Sigma^{kn}\bX \overset{\bff^k}{\lra}$ 
$\bX \ra$ 
$\bMU \wedge \bX$ vanishes.  Calling $\ov{\bff}^k \in \pi_{kn}(D\bX \wedge \bX)$ the adjoint of $\bff^k$ and noting that 
$D\bX \wedge \bX$ is a ring spectrum (cf. p. \pageref{17.25}) (\bX compact 
%%----------------------------------------------------------------------------------------------18
$\implies$ \bX dualizable), \bMU nilpotent technology secures a $d$ such that 
\begin{tikzcd}[ sep=small]
{(\bS^{kn})^{(d)}} \ar{rrr}{\ov{\bff}^k \wedge \cdots \wedge \ov{\bff}^k} &&&{(D\bX \wedge \bX)^{(d)}}
\end{tikzcd}
$\ra D\bX \wedge \bX$
is trivial, so 
\begin{tikzcd}[ sep=small]
{\Sigma^{dkn}\bX} \ar{r}{\bff^{dk}} &{\bX}
\end{tikzcd}
is trivial.]
\vspi
[Note: \ The compactness assumption on \bX cannot be dropped 
(Ravenel\footnote[3]{\textit{Amer. J. Math.} \textbf{106} (1984), 351-414 (cf. 400-401).}).]\\
\endgroup %%------------------------------------<<

\begingroup%%----------------------------------->>
\fontsize{9pt}{11pt}\selectfont
A corollary to the foregoing is that every element of positive degree in $\pi_*(\bS)$ is nilpotent.  
Proof: The elements of $\pi_*(\bS)$ $(n > 0)$ are torsion and $\bMU_*(\bS)$ has no torsion.\\
\endgroup %%------------------------------------<<

\label{17.56}
\label{17.67}
\label{17.68}
\begingroup%%----------------------------------->>
\fontsize{9pt}{11pt}\selectfont
Application: If \bX is compact and if $\bff:\bX \ra \bY$ is an arrow such that 
$\id_{\bMU} \wedge \bff = 0$, then \bff is smash nilpotent.
\vspi
[Suppose that 
$\ov{\bff}:\bS \ra D\bX \wedge \bY$ corresponds to \bff under the identifications 
$[\bX,\bY] \approx$ 
$[\bS \wedge \bX,\bY] \approx $ 
$[\bS,\hom(\bX,\bY)] \approx $ 
$[\bS,D\bX \wedge \bY]$ (\bX compact $\implies$ \bX dualizable) $-$then \bff is smash nilpotent iff $\ov{\bff}$ is smash nilpotent  and 
$\id_{\bMU} \wedge \bff = 0$ iff
$\id_{\bMU} \wedge \ov{\bff} = 0$.  
This allows one to reduce to the case when $\bX = \bS$, the assumption becoming that the composite 
$\bS \overset{\bff}{\lra}$ 
$\bY \ra$ 
$\bMU \wedge \bY$ vanishes.  Put 
$\bE\bY = \ds\bigvee\limits_{i \geq 0} \bY^{(i)}$ $(\bY^{(0)} = \bS)$
 and view $\bE\bY$ as a ring spectrum with multiplication given by concatenation .  
\bMU nilpotence technology now implies that the element of $\pi_*(\bE\bY)$ determined by \bff is nilpotent.]\\
\endgroup %%------------------------------------<<

\begingroup%%----------------------------------->>
\fontsize{9pt}{11pt}\selectfont
\textbf{\small FACT}  \  
Suppose that \bE is complex orientable $-$then the set of complex orientations of \bE is in a one-to-one correspondence with the set of morphisms $\bMU \ra \bE$ of ring spectra.
\vspi
[Note: \ If $\bff:\bMU \ra \bE$ corresponds to $x_{\bE}$, then $\bff_* F_{\bMU} = F_{\bE}$.]\\
\endgroup %%------------------------------------<<

Notation: Given $F \in \FGL_A$, define homogeneous elements $[n]_F(x) \in A[[x]]$ by 
$[1]_F(x)$ $=$ $x$, 
$[n]_F(x) = x +_F[n-1]_F(x)$ $(n > 1)$, and for each prime $p$, write 
$[p]_F(x) = v_0x + \cdots + v_1x^p + \cdots + v_n x^{p^n} + \cdots$ 
($\implies$ $v_0 = p$, $v_n \in A_{2(1 - p^n)}$).

Specialized to $A = \bMU^*(\bS)$, $F = F_{\bMU}$, the $v_n$ can and will be construed as elements of $\bMU_*(\bS)$.\\

\index{Theorem: Exact Functor Theorem}
\index{Exact Functor Theorem}
\textbf{\small EXACT FUNCTOR THEOREM} \quad Let \mM be a graded left $\bMU_*(\bS)$-module $-$then 
$\bMU_*(-) \otimes_{\bMU_*(\bS)} M$ is a $\Z$-graded homology theory on \bSPEC if $\forall \ p \in \bPi$, 
the sequence $\{v_n\}$ is 
$M$-regular, i.e., multiplication by $v_0 = p$ on $M$ and by 
$v_n$ on $M/(v_0M + \cdots + v_{n-1}M)$ for $n \geq 1$ is injective.

[Note: \ This result is due to 
Landweber\footnote[2]{\textit{Amer. J. Math.} \textbf{98} (1976), 591-610; 
see also Rudyak, \textit{Math. Notes} \textbf{40} (1986), 562-569.}.]\\

Remark: Since $\bHSPEC/\bPh$ is the represented equivalent of $\bHT_{\Z}(\bSPEC)$ 
(cf. p. \pageref{17.26}), the exact functor theorem implies that $\exists$ a spectrum \bEM such that 
$\bEM_*(\bX)$ $\approx$ 
$\bMU_*(\bX)$ $\otimes_{\bMU_*(\bS)} M$ $\forall \ \bX$ ($\implies$ $\bEM_*(\bS) \approx M$).

%
%%----------------------------------------------------------------------------------------------19
[Note: \ \bEM is unique up to isomorphism (but is not necessarily unique up to unique isomorphism).  To force the latter, it suffices that \mM be countable and concentrated in even degrees 
(Franke\footnote[3]{\textit{Math. Nachr.} \textbf{158} (1992), 43-65.}).]

Remark: Franke (ibid.) has shown that if \mR is a countable graded $\bMU^*(\bS)$-algebra with unit which, when viewed as a graded left $\bMU_*(\bS)$-module, satisfies the hypotheses of the exact functor theorem, then \bER is a ring spectrum (commutative  if \mR is graded commutative).\\

\begingroup%%----------------------------------->>
\fontsize{9pt}{11pt}\selectfont
Suppose given an $F \in \bFGL_A$ $-$then the homomorphism $f:\bMU^*(\bS) \ra A$ classifying \mF serves to equip 
$A^\OP$ with the structure of a graded  left $\bMU_*(\bS)$-module and the $f(v_n)$ are the $v_n \in A$ per $F$.\\
\endgroup %%------------------------------------<<

\begingroup%%----------------------------------->>
\fontsize{9pt}{11pt}\selectfont
\textbf{\small EXAMPLE}  \ 
Take $A = \Q$ (trivial grading) and let $f:\bMU^*(\bS) \ra \Q$ classify the FGL $x + y$ $-$then $\forall \ p \in \bPi$, 
$f(v_0) = p$ is a unit and $f(v_n) = 0$ $(n \geq 1)$.  Therefore the sequence $\{f(v_n)\}$ is $\Q$-regular and the spectrum produced by the exact functor theorem is $\bH(\Q)$.
\vspi
[Note: \ This would not work if $\Q$ were replaced by $\Z$.]\\
\endgroup %%------------------------------------<<

\begingroup%%----------------------------------->>
\fontsize{9pt}{11pt}\selectfont
\textbf{\small EXAMPLE}  \ Take $A = \Z[u, u^{-1}]$ $(\abs{u} = - 2)$ and let 
$f:\bMU_*(\bS) \ra  \Z[u, u^{-1}]$ classify the FGL $x + y + uxy$.  
Here $f(v_0) = p$, $f(v_1) = u^{p-1}$, $f(v_n) = 0$ $(n > 1)$, thus the conditions of the exact functor theorem are met and the representing spectrum is \bKU (cf. p. \pageref{17.27}).\\
\endgroup %%------------------------------------<<

Let \mA be a divisible abelian group $-$then $\Hom([\bS,-,A)$ is an exact cofunctor which converts coproducts into products, thus is representable (cf. p. \pageref{17.28}) (\bS is compact).  
So: $\exists$ a spectrum $\bS[A]$ such that $\forall \ \bX$, 
$[\bX,\bS[A]] \approx$ 
$\Hom(\pi_0(\bX),A)$.  
Definition: The 
\un{$A$-dual}
\index{A-dual}
$\nabla_A \bX$ 
\index{$\nabla_A \bX$} 
of \bX is $\hom(\bX,\bS[A])$.

Observation:  There is a canonical arrow  \ 
$\bX \lra \nabla_A^2 \bX$, \ and \ $\forall \ n$, 
$\bS[A]^n(\bX)$ $\approx$ 
$\Hom(\pi_n(\bX),A)$.\\

\begin{proposition} \ %18
There are no nonzero phantom maps to $\nabla_A \bX$.
\end{proposition}

[Written out, the claim is that $\Ph(\bY,\nabla_A \bX) = 0$ $\forall \ \bY$, i.e., that the kernel of the arrow 
$[\bY,\nabla_A \bX] \ra \Nat(h_{\bY},h_{\nabla_A \bX})$ is trivial.  
But
$h_{\bY} = \underset{\bY}{\colim}\  h_{\bL}$ $\implies$ 
$\Nat(h_{\bY},h_{\nabla_A \bX}) \approx$ 
$\lim\limits_{\bY} \Nat(h_{\bL},h_{\nabla_A \bX}) \approx$ 
$\lim\limits_{\bY}[\bL,\nabla_A \bX]$.  On the other hand, there is an arrow 
$\Hom(\pi_0(\bY \wedge \bX),A) \ra$ 
$\lim\limits_{\bY} \Hom(\pi_0(\bL \wedge \bX),A)$ and a commutative diagram \quad %\\
%\begin{tikzcd}%[sep=small]
%{[\bY,\nabla_A \bX]} \ar{d} \ar{r} &{\Hom(\pi_0(\bY \wedge \bX),A)} \ar{d}\\
%{\lim\limits_{\bY}[\bL,\nabla_A \bX]} \ar{r} &{\lim\limits_{\bY} \Hom(\pi_0(\bL \wedge \bX),A)}
%\end{tikzcd}
\begin{tikzcd}%[sep=small]
{[\bY,\nabla_A \bX]} \ar{d} \ar{r} &{}\\
{\lim\limits_{\bY}[\bL,\nabla_A \bX]} \ar{r} &{}
\end{tikzcd}
\begin{tikzcd}%[sep=small]
{\Hom(\pi_0(\bY \wedge \bX),A)} \ar{d}\\
{\lim\limits_{\bY} \Hom(\pi_0(\bL \wedge \bX),A)}
\end{tikzcd}
.  
%%----------------------------------------------------------------------------------------------20
The horizontal arrows are isomorphisms, as is the vertical arrow on the right (cf. $\S 15$, Proposition 18 and subsequent remark).  Therefore the vertical arrow on the left is an isomorphism, hence $\Ph(\bY,\nabla_A \bX) = 0$.]\\

\begingroup%%----------------------------------->>
\fontsize{9pt}{11pt}\selectfont
\textbf{\small EXAMPLE}  \ 
Take $A = \Q/\Z$ $-$then $\nabla_{\Q/\Z} \bX$ is the 
\underline{Brown-Comenetz}
\footnote[2]{\textit{Amer. J. Math.} \textbf{98} (1976), 1-27.}
\un{dual}
\index{Brown-Comenetz dual} 
of \bX and, thanks to the Pontryagin duality theorem, the canonical arrow 
$\bX \ra \nabla_{\Q/\Z}^2 \bX$ is an isomorphism if the homotopy groups of \bX are finite..  
Example: $\nabla_{\Q/\Z} \bH(\Z/p\Z) \approx \bH(\Z/p\Z)$.
\vspi
[Note: \ In homotopy, the canonical arrow 
$\pi_n(\bX) \ra \pi_n(\nabla_{\Q/\Z}^2 \bX)$ is the inclusion of $\pi_n(\bX)$ into its double dual per $\Q/\Z$ and if $\pi_n(\bX)$ is finitely generated, then 
$\pi_n(\nabla_{\Q/\Z}^2 \bX) = \pro \pi_n(\bX)$, the profinite completion of $\pi_n(\bX)$.]\\
\endgroup %%------------------------------------<<

\begingroup%%----------------------------------->>
\fontsize{9pt}{11pt}\selectfont
\textbf{\small FACT}  \ 
Take $\bC = \bHSPEC$ $-$then $\forall \ \bX$, $h_{\nabla_{\Q/\Z}\bX}$ is an injective object of $[(\cpt \bC)^\OP,\bAB]^+$.

[It follows from the definitions (and Yoneda) that this is true if \bX is compact.  
In general, there are compact objects 
$\bK_i$ and an arrow 
$\nabla_{\Q/\Z}\bX \overset{\bff}{\lra} \ds\prod\limits_i \nabla_{\Q/\Z}\bK_i$ such that $h_{\bff}$ is a monomorphism 
($\Q/\Z$ is an injective coseparator in \bAB).  
Consider now the exact triangle 
$\bY \overset{\bphi}{\lra}$ 
$\nabla_{\Q/\Z}\bX \overset{\bff}{\lra}$ 
$\ds\prod\limits_i \nabla_{\Q/\Z}\bK_i \ra$ 
$\Sigma \bY$.  
Since $\bff \circ \bphi = 0$ (cf. $\S 15$, Proposition 3), 
$h_{\bff} \circ h_{\bphi} = 0$ $\implies$ 
$h_{\bphi} = 0$ $\implies$ 
$\bphi \in \Ph(\bY,\nabla_{\Q/\Z}\bX)$ $\implies$ $\bphi = 0$ (cf. Proposition 18), so 
$\nabla_{\Q/\Z}\bX$ is a retract of $\ds\prod\limits_i \nabla_{\Q/\Z}\bK_i$.]\\
\endgroup %%------------------------------------<<

\label{17.36} %dmc mnft
\begingroup%%----------------------------------->>
\fontsize{9pt}{11pt}\selectfont
\textbf{\small EXAMPLE}  \ 
Define $\bS[\Z]$ by the exact triangle 
$\bS[\Z] \overset{\bu}{\lra}$ 
$\bS[\Q] \overset{\bv}{\lra}$ 
$\bS[\Q/\Z] \overset{\bw}{\lra}$ 
$\Sigma \bS[\Z]$ 
where 
$\bv_*:\pi_0(\bS[\Q]) \ra \pi_0(\bS[\Q/\Z])$ corresponds to the projection $\Q \ra \Q/\Z$ $-$then 
$\pi_0(\bS[\Z]) \approx \Z$ and 
$\bu_*:\pi_0(\bS[\Z] ) \ra \pi_0(\bS[\Q])$ corresponds to the inclusion $\Z \ra \Q$.  
Definition: The 
\un{Anderson dual}
\index{Anderson dual} 
of $\nabla_{\Z}\bX$ of \bX is $\hom(\bX,\bS[\Z] )$.  
There is a canonical arrow $\bX \ra \nabla_{\Z}^2\bX$ which is an isomorphism if the homotopy groups of \bX are finitely generated.  
Examples: 
(1) $\nabla_{\Z} \bH(\Z) \approx \bH(\Z)$; 
(2) $\nabla_{\Z} \bKU \approx \bKU$.\\
\endgroup %%------------------------------------<<

\begingroup%%----------------------------------->>
\fontsize{9pt}{11pt}\selectfont
\textbf{\small FACT}  \  
Suppose that the homotopy groups of \bX are finite $-$then 
$\Sigma \nabla_{\Z} \bX \approx \nabla_{\Q/\Z} \bX$.\\
\endgroup %%------------------------------------<<

Given an abelian group $G$, define the 
\un{Moore spectrum of type $G$}
\index{Moore spectrum of type $G$} 
by the exact triangle 
$\bigvee\limits_j \bS \ra$ 
$\bigvee\limits_i \bS \ra$ 
$\bS(G) \ra$ 
$\bigvee\limits_j \Sigma \bS \ra$, where 
$0 \ra$ 
$\bigoplus\limits_j \Z \ra$ 
$\bigoplus\limits_i \Z \ra$ 
$G \ra 0$ is a presentation of $G$ $-$then $\bS(G)$ is connective and $\pi_0(\bS(G)) = G$.  
Example: $\bS(\Z) = \bS$.\\

\begin{proposition} \ %19
Given a spectrum \bX and an abelian group $G$, there are short 
%%----------------------------------------------------------------------------------------------21
exact sequences 
\[
\begin{cases}
\ 0 \lra \quad \  \pi_n(\bX) \otimes G \ \quad \lra \pi_n(\bX \wedge \bS(G)) \lra \Tor(\pi_{n-1}(\bX),G) \lra 0\\
\ 0 \lra \Ext(G,\pi_{n+1}(\bX)) \lra \ [\Sigma^n \bS(G),\bX]  \ \lra \ \Hom(G,\pi_n(\bX))\  \lra 0
\end{cases}
.
\]
\end{proposition}
\vspace{0.25cm}

\label{17.35}
Application: $\bH(\Z) \wedge \bS(G) \approx \bH(G)$, the Eilenberg-MacLane spectrum attached to \mG 
(cf. p. \pageref{17.29}).\\

\label{17.37}
\label{17.39}
\begingroup%%----------------------------------->>
\fontsize{9pt}{11pt}\selectfont
\textbf{\small EXAMPLE} \quad
Take $G = \Z_P$ $-$then $\bS(\Z_P)$ is a commutative ring spectrum.
\vspi
[Note: \ $\bS(\Q) \approx \bH(\Q)$ (since $\pi_n(\bS) \otimes \Q = 0$ for $n \neq 0$).]\\
\endgroup %%------------------------------------<<

\begingroup%%----------------------------------->>
\fontsize{9pt}{11pt}\selectfont
\textbf{\small EXAMPLE}  \ 
Take $G = \Z/p\Z$, where $p$ is odd $-$then 
$\bS(\Z/p\Z) \wedge \bS(\Z/p\Z) \approx$ 
$\bS(\Z/p\Z) \vee \Sigma \bS(\Z/p\Z)$ and $\bS(\Z/p\Z)$ is a commutative ring spectrum if $p > 3$.
\vspi
[Note: \ When $p = 3$, $\bS(\Z/3\Z)$ admits a commutative multiplication with unit but associativity breaks down.]\\
\endgroup %%------------------------------------<<

\begingroup%%----------------------------------->>
\fontsize{9pt}{11pt}\selectfont
\textbf{\small EXAMPLE}  \ 
Take $G = \Z/2\Z$ $-$then $\bS(\Z/2\Z)$ has no multiplication with unit ($\bS(\Z/2\Z)$ is not a retract of 
$\bS(\Z/2\Z) \wedge \bS(\Z/2\Z)$).
\vspi
[Note: \ $\Hom(\Z/2\Z,\Z/2\Z) = \Z/2\Z$ whereas 
$[\bS(\Z/2\Z),\bS(\Z/2\Z)] = \Z/4\Z$.  Because of this, one cannot construct an additive functor 
$\bAB \overset{F}{\ra} \bHSPEC$ such that $FG = \bS(G)$ (there is no ring homomorphism $\Z/2\Z \ra \Z/4\Z$).]\\
\endgroup %%------------------------------------<<

\begingroup%%----------------------------------->>
\fontsize{9pt}{11pt}\selectfont
\label{17.34}
\textbf{\small EXAMPLE}  \ 
Fix $p \in \bPi$ $-$then 
$\bS(\Z/p^\infty\Z) \approx$ 
$\tel(\bS(\Z/p\Z) \ra$
$\bS(\Z/p^2\Z) \ra \cdots)$ $\implies$ 
$\Sigma^{-1}\bS(\Z/p^\infty\Z) \approx$ 
$\tel(\Sigma^{-1} \bS(\Z/p\Z) \ra$
$\Sigma^{-1}\bS(\Z/p^2\Z) \ra \cdots)$.  But since 
$\bS \ra \bS \ra$ 
$\bS(\Z/p^n\Z) \ra \Sigma \bS$ is exact, 
$\bS(\Z/p^n\Z) \approx$ 
$\Sigma D \bS(\Z/p^n\Z)$, \ so \ 
$\Sigma^{-1}\bS(\Z/p^\infty\Z)$ $\approx$ \ 
$\tel ( D \bS(\Z/p\Z) \ra$ \ 
$D \bS(\Z/p^2\Z) \ra \cdots)$.  \ 
Accordingly, \ $\forall \ \bX$, 
$\hom(\Sigma^{-1}\bS(\Z/p^\infty\Z),$$\bX)$ $\approx$ 
$\mic\hom(D \bS(\Z/p\Z),\bX) \la$ 
$\hom(D \bS(\Z/p^2\Z),\bX)  \la \cdots)$.  
However, $\forall \ n$, $\bS(\Z/p^n\Z)$ is compact, hence dualizable $\implies$ 
$D \bS(\Z/p^n\Z)$ dualizable (cf. $\S 15$, Proposition 32) $\implies$ 
$\hom(D \bS(\Z/p^n\Z),\bX)$ $\approx$ 
$\bS(\Z/p^n\Z) \wedge \bX$.  
Thus, $\forall \ \bX$, 
$\hom(\Sigma^{-1}\bS(\Z/p^\infty\Z), \bX)$ $\approx$ 
$\mic(\bS(\Z/p\Z) \wedge \bX \la$ 
$\bS(\Z/p^2\Z) \wedge \bX \la \cdots)$.
Example: 
$\mic(\bS(\Z/p\Z)  \la$ 
$\bS(\Z/p^2\Z) \la \cdots) \approx$ 
$\bS(\widehat{\Z}_p)$ $\implies$ 
$\Sigma D \bS(\Z/p^\infty\Z)  \approx \bS(\widehat{\Z}_p)$.\\
\endgroup %%------------------------------------<<

Fix a spectrum \bE $-$then a morphism $\bff:\bX \ra \bY$ in \bHSPEC is said to be an 
\un{$\bE_*$-equivalence}
\index{E$_*$-equivalence (\bE a spectrum)} 
if $\bff_*:\bE_*(\bX) \ra \bE_*(\bY)$ is an isomorphism.  
Denoting by $\sS_{\bE}$ the class of $\bE_*$-equivalences, the 
Bousfield-Margolis localization theorem guarantees the existence of a localization functor $T_{\bE}$ such that $S_{\bE}^\perp$ is the class of $\bE_*$-local ($= T_{\bE}$-local) spectra.  
In this connection, recall that \bX is $\bE_*$-local  iff 
$[\bY,\bX] = 0$ for all $\bE_*$-acyclic ($= T_{\bE}$-acyclic) \bY (cf. $\S 15$, Proposition 27) and the class of 
$\bE_*$-local spectra is the object class of a thick 
%%----------------------------------------------------------------------------------------------22
subcategory of \bHSPEC which is closed under the formation of products in \bHSPEC (cf. $\S 15$, Proposition 28).  Let us also bear in mind that $T_{\bE}$ has the IP (cf. $\S 15$, Proposition 40). 

Notation: $\bHSPEC_{\bE}$ is the full subcategory of \bHSPEC whose objects are the $\bE_*$-local spectra, 
$L_{\bE}:\bHSPEC \ra \bHSPEC_{\bE}$ is the associated reflector, and 
$l_{\bE}:\bX \ra L_{\bE}\bX$ is the arrow of localization.

[Note: \ The objects of $\bHSPEC_{\bE}$ are the objects of $\langle \bE \rangle$, the Bousfield class of \bE, and 
$L_{\bE} \approx$ $L_{\bF}$ iff $\langle \bE \rangle = \langle \bF \rangle$.  $\bHSPEC_{\bE}$ is a CTC 
(cf. p. \pageref{17.30}) but need not be compactly generated 
(Strickland\footnote[2]{\textit{No Small Objects}, Preprint.}
).]

Remark: Ohkawa\footnote[3]{\textit{Hiroshima Math. J.} \textbf{19} (1989), 631-639.}
has shown that the conglomerate $\langle \bHSPEC \rangle$ whose elements are the Bousfield classes is codable by a set.\\

\label{17.81}
\textbf{\small LEMMA}  \  
Given spectra \bE and \bF, suppose that 
$\langle \bE \rangle \leq \langle \bF \rangle$ $-$then $\forall \ \bX$, 
$T_{\bE}T_{\bF} \bX \approx$ 
$T_{\bE} \bX \approx$ 
$T_{\bF}T_{\bE}\bX$.\\

\begingroup%%----------------------------------->>
\fontsize{9pt}{11pt}\selectfont
\textbf{\small EXAMPLE}  \ 
Suppose that \bX is connective $-$then $\bX = 0$ iff \bX is $\bH(\Z)_*$-acyclic.
\vspi
[Note: \ $\nabla_{\Q/\Z}\bS$ ($= \bS[\Q/\Z]$) is $\bH(\Z)_*$-acyclic and nonzero (although 
$\nabla_{\Q/\Z}\bS \wedge \nabla_{\Q/\Z}\bS = 0$).]\\
\endgroup %%------------------------------------<<

\begingroup%%----------------------------------->>
\fontsize{9pt}{11pt}\selectfont
Instead of working with $\bE_*$-equivalences, one could work instead with $\bE^*$-equivalences and then define 
$\bE^*$-local spectra in the obvious way.  
Problem: Do the $\bE^*$-local spectra constitute the object class of a reflective subcategory of \bHSPEC?  
While the answer is unknown in general, one does have the following partial result due to 
Bousfield\footnote[2]{\textit{Cohomological Localizations of Spaces and Spectra}, Preprint.}
.\\
\endgroup %%------------------------------------<<

\begingroup%%----------------------------------->>
\fontsize{9pt}{11pt}\selectfont
\index{Theorem: Cohomological Localization Theorem}
\index{Cohomological Localization Theorem}
\textbf{\small COHOMOLOGICAL LOCALIZATION THEOREM} \quad 
Suppose that \bE has the following property: $\forall \ n$
$\Z/p\Z \otimes \pi_n(\bE)$ and 
$\Tor(\Z/p\Z,\pi_n(\bE))$ are finite $\forall \ p \in \bPi$ $-$then there exists an \bF such that the $\bE^*$-equivalences are the same as the $\bF_*$-equivalences, so cohomological localization with respect to \bE exists and is given by homological localization with respect to \bF.
\vspi
[Note: \ When the $\pi_n(\bE)$ are finitely generated, one can take $\bF = \nabla_{\Z} \bE$.]\\
\endgroup %%------------------------------------<<

Given an abelian group \mG, call $\sS(G)$ 
\index{S(G) (\mG an abelian group)} 
the class of abelian groups \mA such that 
$A \otimes G =$ 
$0 = \Tor(A,G)$ (cf. p. \pageref{17.31}).\\

\begin{proposition} \ %20
$\sS(G^\prime) = \sS(G\pp)$ iff 
$\langle\bS(G^\prime)\rangle = \langle\bS(G\pp)\rangle$.
\end{proposition}

%%----------------------------------------------------------------------------------------------23
This result reduces the problem of inventoring the $L_{\bS(G)}$ to when 
$G = \Z_P$ or 
$G = \bigoplus\limits_{p \in P} \Z/p\Z$.\\

\begingroup%%----------------------------------->>
\fontsize{9pt}{11pt}\selectfont
\textbf{\small EXAMPLE}  \ 
$\langle \bS(\Z_P)\rangle  = \langle \bS(\Q)\rangle \vee \ds\bigvee\limits_{p \in P} \langle \bS(\Z/p\Z)\rangle$ 
$(\implies$ 
$\langle \bS\rangle  = \langle \bS(\Q)\rangle \vee \ds\bigvee\limits_{p \in P} \langle \bS(\Z/p\Z)\rangle$).  
And: 
$\langle \bS(\Q) \rangle \wedge \langle \bS(\Z/p\Z) \rangle =$ $\langle 0 \rangle$ $\&$ 
$\langle \bS(\Z/p\Z) \rangle \wedge \langle \bS(\Z/q\Z) \rangle = \langle 0 \rangle$ $(p \neq q)$.\\
\endgroup %%------------------------------------<<

\begin{proposition} \ %21
Let $G = \Z_P$ $-$then 
$L_{\bS(\Z_P)} \bX = \bS(\Z_P) \wedge \bX$ and 
$\pi_*(L_{\bS(\Z_P)} \bX =$ 
$\bZ_P \otimes \pi_*(\bX)$.
\end{proposition}

[$\bS(\Z_P)$ is a commutative ring spectrum with the property that the product 
$\bS(\Z_P) \wedge$ $\bS(\Z_P)$ $\ra \bS(\Z_P)$ is an isomorphism, thus $T_{\bS(\Z_P)}$ is smashing 
(cf. p. \pageref{17.32}) and 
$\bX \approx \bS \wedge \bX \ra$ $\bS(\Z_P) \wedge \bX$ is the arrow of localization.]\\

\label{17.38}
\begingroup%%----------------------------------->>
\fontsize{9pt}{11pt}\selectfont
\textbf{\small FACT}  \  
Suppose that \bX is connective $-$then 
$L_{\bS(\Z_P)} \bX  \approx L_{\bH(\Z_P)} \bX$.
\vspi
[Note: \ Take $P = \bPi$ to see that 
$L_{\bS(\Z)} \bX  \approx L_{\bH(\Z)} \bX$, i.e., $\bX \approx L_{\bH(\Z)} \bX$.]\\
\endgroup %%------------------------------------<<

Write $\bHSPEC_{P}$ for the full subcategory of \bHSPEC whose objects are $P$-local ($= \bS(\Z_P)_*$-local) 
(use the symbol $\bHSPEC_{\Q}$ if $P = \emptyset$) $-$then the objects of $\bHSPEC_{P}$ are those \bX which are $P$-local in homotopy, i.e., $\forall \ n$, $\pi_n(\bX)$ is $P$-local and $\bHSPEC_{P}$ is a monogenic compactly generated CTC.\\

\begingroup%%----------------------------------->>
\fontsize{9pt}{11pt}\selectfont
\textbf{\small FACT}  \  
The category $\bHSPEC_{\Q}$ is equivalent to the category of graded vector spaces over \bQ.
\vspi
[Note: \ The objects of $\bHSPEC_{\Q}$ are the 
\un{rational spectra}
\index{rational spectra}.]\\
\endgroup %%------------------------------------<<

\begin{proposition} \ %22
Let $G = \Z/p\Z$ $-$then 
$L_{\bS(\Z/p\Z)} \bX = \hom(\Sigma^{-1}\bS(\Z/p^\infty\Z),\bX)$ and there is a split short exact sequence 
$0 \ra$ 
$\Ext(\Z/p^\infty\Z,\pi_*(\bX)) \ra$ 
$\pi_*(L_{\bS(\Z/p\Z)} \bX) \ra$ 
$\Hom(\Z/p^\infty\Z,\pi_{*-1}(\bX)) \ra 0$.
\end{proposition}

[Consider the exact triangle 
$\hom\bigl(\bS\bigl(\Z\bigl[\ds\frac{1}{p}\bigr]\bigr),\bX\bigr) \ra$
$\hom(\bS,\bX)$ $\ra$ 
$\hom(\Sigma^{-1}\bS(\Z/p^\infty\Z),\bX)$ $\ra$ 
$\Sigma\hom\bigl(\bS\bigl(\Z\left[\ds\frac{1}{p}\right]\bigr),\bX\bigr)$.  
On the one hand, 
$\hom(\Sigma^{-1}\bS(\Z/p^\infty\Z)\bX)$ is $\bS(\Z/p\Z)_*$-local (for $\bS(\Z/p\Z) = \bS(\Z/p^\infty\Z)$) 
and, on the other, 
$\hom\bigl(\bS\bigl(\Z\left[\ds\frac{1}{p}\right]\bigr),\bX\bigr)$ is $\bS(\Z/p\Z)_*$-acyclic 
(its homotopy groups are uniquely $p$-divisible).  
Therefore 
$\bX \approx$ 
$\hom(\bS,\bX) \ra$ 
$\hom(\Sigma^{-1}\bS(\Z/p^\infty\Z),$ $\bX)$ is the arrow of localization.]

[Note: \ The $\bS(\Z/p\Z)_*$-local spectra are those \bX such that $\forall \ n$, $\pi_n(\bX)$ is $p$-cotorsion.  
Proof: 
$\hom\bigl(\bS\bigl(\Z\left[\ds\frac{1}{p}\right]\bigr),\bX\bigr) = 0$ iff $\forall \ n$, 
$\Hom\bigl(\Z\left[\ds\frac{1}{p}\right],\pi_n(\bX)\bigr) = 0$ $\&$ 
$\Ext\bigl(\Z\left[\ds\frac{1}{p}\right],\pi_n(\bX)\bigr) = 0$.]\\

If the homotopy groups of \bX are finitely generated, put 
$\widehat{\bX}_p = L_{\bS(\Z/p\Z)} \bX$ and call $\widehat{\bX}_p$ the 
\un{$p$-adic completion}
\index{p-adic completion (spectra)} 
of \bX.  
Justification: 
$\forall \ n$, $\pi_n(\widehat{\bX}_p) \approx$ 
$\pi_n(\bX)_p^{\widehat{\  }}$ (cf. p. \pageref{17.33}).  
Example: 
%%----------------------------------------------------------------------------------------------24
$\widehat{\bS}_p = L_{\bS(\Z/p\Z)} \bS =$ 
$\hom(\Sigma^{-1}\bS(\Z/p^\infty\Z),\bS) =$ 
$D \Sigma^{-1}\bS(\Z/p^\infty\Z) =$ 
$\Sigma D\bS(\Z/p^\infty\Z) =$ 
$\bS(\widehat{\Z}_p)$ (cf. p. \pageref{17.34}).\\

\begin{proposition} \ 
The arrow of localization per 
$\bigoplus\limits_{p \in P} \Z/p\Z$ is $\bX \ra \prod\limits_{p \in P} L_{\bS(\Z/p\Z)} \bX$ (cf. $\S 9$, Proposition 22).\\
\end{proposition}

\begingroup%%----------------------------------->>
\fontsize{9pt}{11pt}\selectfont
\textbf{\small FACT}  \  
$\forall \ \bX$, there is an exact triangle 
$\hom(\bS(\Q),\bX) \ra$ 
$\bX \ra$ 
$\ds\prod\limits_p \hom(\Sigma^{-1} \bS(\Z/p^\infty\bZ),\bX) \ra$
$\Sigma \hom(\bS(\Q),\bX)$.\\
\endgroup %%------------------------------------<<

\begingroup%%----------------------------------->>
\fontsize{9pt}{11pt}\selectfont
\textbf{\small FACT}  \  
$\forall \ \bX$, there is an exact triangle 
$\ds\bigvee\limits_p  \bX \wedge \Sigma^{-1} \bS(\Z/p^\infty\bZ) \ra$ 
$\bX \ra$ 
$\bX \wedge \bS(\Q) \ra$ 
$\ds\bigvee\limits_p \Sigma(\bX \wedge \Sigma^{-1} \bS(\Z/p^\infty\bZ))$.\\
\endgroup %%------------------------------------<<

\begin{proposition} \ %24
Let $G$, $K$ be abelian groups such that $\sS(G) = \sS(K)$ $-$then $\forall \ \bX$, 
$\langle \bX \wedge \bS(G) \rangle = \langle \bX \wedge \bS(K) \rangle$.\\
\end{proposition}

\begingroup%%----------------------------------->>
\fontsize{9pt}{11pt}\selectfont
\textbf{\small EXAMPLE}  \ 
Let $G$, $K$ be abelian groups such that $\sS(G) = \sS(K)$ $-$then
$\langle \bH(G) \rangle = \langle \bH(K)\rangle$.  In fact, 
$
\begin{cases}
\ \bH(G) = \bH(\Z) \wedge \bS(G)\\
\ \bH(K) = \bH(\Z) \wedge \bS(K)
\end{cases}
$
(cf. p. \pageref{17.35}).\\
\vspace{0.25cm}
\endgroup %%------------------------------------<<

\label{17.45}
\begingroup%%----------------------------------->>
\fontsize{9pt}{11pt}\selectfont
\textbf{\small FACT}  \  
Suppose that $\bE \wedge \bS(\Q) \neq 0$ $-$then $\forall$ \bX, 
$L_{\bE \wedge \bS(\Q)} \bX \approx L_{\bS(\Q)} \bX$.\\ 
\endgroup %%------------------------------------<<

\textbf{\small LEMMA}  \  
Given a connective spectrum \bE, put 
$\pi \bE = \bigoplus\limits_n \pi_n(\bE)$ $-$then 
$\langle \bH(\pi\bE) \rangle \leq$
$\langle \bE \rangle \leq$
$\langle \bS(\pi\bE) \rangle$.

[$\langle \bH(\pi(\bE)) \rangle \leq \langle \bE \rangle$: 
Since \bE is connective, 
$\sS(\pi\bE) = \sS(\bigoplus\limits_n \bH_n (\bE;\Z))$, so 
$\langle \bH(\pi(\bE) \rangle =$
$\langle \bH(\bigoplus\limits_n \bH_n(\bE; \bZ)) \rangle =$
$\langle \bigvee\limits_n \Sigma^n \bH(\bH_n(\bE;\Z)) \rangle =$
$\langle \bE \wedge \bH(\Z) \rangle$ (cf. p. \pageref{17.36}), which is 
$\leq \langle \bE \rangle$.

$\langle \bE \rangle \leq \langle \bS(\pi\bE) \rangle$: 
Let $G_1$ be the direct sum of the groups in the set 
$\{\Q,\Z/p\Z \ (p \in \bPi\}$ with $\sS(G_1) = \sS(\pi\bE)$ 
and let $G_2$ be the direct sum of what remains $-$then
$\langle \bS(G_1) \rangle \wedge \langle \bS(G_2) \rangle = \langle 0 \rangle$ $\&$ 
$\langle \bS(G_1) \rangle \vee \langle \bS(G_2) \rangle = \langle \bS \rangle$.  And: 
$\bE \wedge \bS(G_2) = 0$, hence 
$\langle \bE \rangle =$ 
$\langle \bE \rangle \wedge \langle \bS \rangle = $
$\langle \bE \rangle \wedge (\langle \bS(G_1) \rangle \vee \langle \bS(G_2) \rangle) = $ 
$(\langle \bE \rangle \wedge \langle \bS(G_1) \rangle) \vee (\langle \bE \rangle \wedge \langle \bS(G_2) \rangle) =$ 
$\langle \bE \rangle \wedge \langle \bS(G_1) \rangle = $
$\langle \bE \rangle \wedge \langle \bS(\pi(\bE) \rangle \leq$ 
$\langle \bS(\pi(\bE) \rangle $.]\\

\begin{proposition} \ %25
Let \bE, \bX be connective $-$then $L_{\bE}\bX \approx L_{\bS(\pi\bE)}\bX$, where $\pi\bE =$ 
$\bigoplus\limits_n \pi_n(\bE)$.
\end{proposition}

[The lemma implies that the arrow of localization 
$\bX \ra L_{\bS(\pi\bE)}\bX$ is an $\bE_*$-equivalence.  But
$L_{\bS(\pi\bE)}\bX = L_{\bH(\pi\bE)}\bX$ (cf. infra) and $L_{\bH(\pi\bE)}\bX$ is $\bE_*$-local (by the lemma).]\\

\begingroup%%----------------------------------->>
\fontsize{9pt}{11pt}\selectfont
\textbf{\small LEMMA}  \  
Let \bE, \bX be spectra and let $G$ be an abelian group $-$then the arrow 
$L_{\bS(G)}L_{\bE}\bX \ra$ $L_{\bE \wedge \bS(G)}\bX$ 
is an isomorphism if $G$ is torsion or if $\bE \wedge \bS(\Q) \neq 0$.
\vspi
%%----------------------------------------------------------------------------------------------25
[Suppose first that $G$ is torsion, say \ 
$G = \ds\bigoplus\limits_{p \in P} \Z/p\Z$ \ (this entails no loss of generality). \  Since 
$L_{\bS(G)}L_{\bE}\bX \ra L_{\bE \wedge \bS(G)}\bX$ is an $(\bE \wedge \bS(G))_*$-equivalence, it suffices to prove that 
$L_{\bS(G)}L_{\bE}\bX$ is $(\bE \wedge \bS(G))_*$-local or still, that \ 
$[\bY,L_{\bS(G)}L_{\bE}\bX] = 0$ \ for all $(\bE \wedge \bS(G))_*$-acyclic \bY.  \ But 
$[\bY,L_{\bS(G)}L_{\bE}\bX] =$ 
$\bigl[\bY,\hom\bigl(\ds\bigvee\limits_{p \in P} \Sigma^{-1}\bS(Z/p^\infty\Z),L_{\bE}\bX\bigr)\bigr] =$
$\bigl[\bY \wedge \ds\bigvee\limits_{p \in P} \Sigma^{-1} \bS(Z/p^\infty\Z),L_{\bE}\bX\bigr]$ and 
$\bY \wedge \ds\bigvee\limits_{p \in P} \Sigma^{-1} \bS(Z/p^\infty\Z)$ is $\bE_*$-acyclic 
$\bigl( \sS \bigl(\ds\bigoplus\limits_{p \in P} Z/p^\infty\Z\bigr) = \sS(G)\bigr)$.  
To discuss the case, viz. when $\bE \wedge \bS(\Q) \neq 0$, one can take $G = \Z_P$.  
Because 
$L_{\bE \wedge \bS(G)} \bX$ is $\bS(G)$-local, it need only be shown that 
$L_{\bE}\bX \ra L_{\bE \wedge \bS(G)}\bX$ is an $\bS(G)_*$-equivalence.  
However 
$\langle \bS(G) \rangle =$ 
$\langle \bS(\Q) \rangle \vee \ds\bigvee\limits_{p \in P} \langle \bS(Z/p\Z) \rangle$, 
which reduces the problem to showing that 
$L_{\bE}\bX \ra L_{\bE \wedge \bS(G)}\bX$ is an $\bS(\Q)_*$-equivalence and an  
$\bS(\Z/p\Z)_*$-equivalence for each $p \in P$.  
Due to our assumption that $\bE \wedge \bS(\Q) \neq 0$ just the second possibility is at issue.  
For this, consider the commutative triangle 
\begin{tikzcd}[sep=large]
{L_{\bE}\bX} \ar{rd} \ar{r} &{L_{\bE \wedge \bS(G)}} \ar{d}\\
&{L_{\bE \wedge \bS(\Z/p\Z)}}
\end{tikzcd}
.  Here, the arrow $L_{\bE \wedge \bS(G)} \ra L_{\bE \wedge \bS(\Z/p\Z)}$ is an $\bS(\Z/p\Z)_*$-equivalence 
$(L_{\bS(\Z/p\Z)}L_{\bE \wedge \bS(G)} \bX \approx$ 
$L_{\bE \wedge \bS(G) \wedge \bS(\Z/p\Z)} \approx$ 
$L_{\bE \wedge \bS(\Z/p\Z)}$), as is the arrow 
$L_{\bE}\bX \ra L_{\bE \wedge \bS(\Z/p\Z)}$ 
$(L_{\bS(\Z/p\Z)}L_{\bE}\bX \approx$ $L_{\bE \wedge \bS(\Z/p\Z)})$.  Therefore the arrow 
$L_{\bE}\bX \ra$ 
$L_{\bE \wedge \bS(G)}\bX$ is an $\bS(\Z/p\Z)_*$-equivalence.]
\vspi
[Note: \ The assumption that $\bE \wedge \bS(\Q) \neq 0$ cannot be dropped.  
Example: $L_{\bS(\Q)} L_{\bS(\Z/p\Z)} \bH(\Z) \neq 0$, yet $L_{\bS(\Z/p\Z) \wedge \bS(\Q)} \bH(\Z) = 0$.]\\
\endgroup %%------------------------------------<<

\begingroup%%----------------------------------->>
\fontsize{9pt}{11pt}\selectfont
To tie up the loose end in the proof of Proposition 25, observe that 
$\bH(\Z) \wedge \bS(\Q) \approx$ 
$\bH(\Q) \neq 0$ (cf. p. \pageref{17.37}).  In addition, since \bX is connective, 
$\bX \approx$ $L_{\bH(\Z)}\bX$ (cf.  p. \pageref{17.38}), hence 
$L_{\bS(\pi\bE)} \bX \approx$ 
$L_{\bS(\pi\bE)} L_{\bH(\Z)} \bX \approx$ 
$L_{\bH(\Z) \wedge \bS(\pi\bE)} \bX$ $\approx$
$L_{\bH(\pi\bE)} \bX$.
(cf. p. \pageref{17.39}).\\
\endgroup %%------------------------------------<<

\textbf{\small LEMMA}  \  
Let \mA be a ring with unit, \mM a left $A$-module $-$then $\sS(A) = \sS(A \oplus M)$.\\

Application: Suppose that \bE is a ring spectrum $-$then 
$\sS(\pi_0(\bE)) = \sS\bigl( \bigoplus\limits_n \pi_n(\bE) \bigr)$.\\

Example: Take $\bE = \bMU$ $-$then 
$\sS(\Z) = \sS\bigl(\bigoplus\limits_n \pi_n(\bMU) \bigr)$, thus for any connective \bX, 
$L_{\bMU} \bX \approx$ 
$L_{\bS(\Z)} \bX \approx$
$L_{\bS} \bX \approx \bX$.

[Note: \ It follows that all compact spectra are $\bMU_*$-local.  Indeed, a compact object in \bHSPEC is isomorphic to a 
$\bQ_q^\infty K$, where \mK is a pointed finite CW complex (cf. p. \pageref{17.40}).  And: 
$\bQ_q^\infty K \approx$ 
$\bS^{-q} \wedge K \approx$ 
$\bS^{-q} \wedge \bQ^\infty K \approx$ (cf. p. \pageref{17.41}).  
But $\bQ^\infty K$ is connective (cf. p. \pageref{17.42}) (\mK is wellpointed).  
Therefore $\bQ^\infty K$ is $\bMU_*$-local, hence 
$\bS^{-q} \wedge \bQ^\infty K$ is too  (cf. p. \pageref{17.43}) $(\bS^{-q}$ is compact and 
\bHSPEC is a monogenic compactly generated CTC).]\\

\begingroup%%----------------------------------->>
\fontsize{9pt}{11pt}\selectfont
\textbf{\small FACT}  \  
Let \bX, \bY be spectra with $\bY^*(\bX) = 0$.  Assume: The homotopy groups of \bY are finite $-$then 
$\pi_*(\bX \wedge \nabla_{\Q/\Z} \bY) = 0$.
\vspi
%%----------------------------------------------------------------------------------------------26
%%----------------------------------------------------------------------------------------------26
%%following s missing from online version
[$\bY \approx \nabla_{\Q/\Z}^2 \bY$ $\implies$ 
$0 = [\bX,\Sigma^n \bY] = [\bX,\Sigma^n \nabla_{\Q/\Z}^2 \bY] =$ 
$\Hom(\pi_n(\bX \wedge \nabla_{\Q/\Z} \bY,\Q/\Z)$.]\\
\endgroup %%------------------------------------<<

\begingroup%%----------------------------------->>
\fontsize{9pt}{11pt}\selectfont
\textbf{\small EXAMPLE}  \ 
The assumptions of the preceding result are met if 
$\bX = \bM\bU$, $\bY = \bS$.  Therefore 
$\nabla_{\Q/\Z} \bS$ is $\bM\bU_*$-acyclic, so $\langle \bM\bU\rangle < \langle \bS \rangle$.\\
\endgroup %%------------------------------------<<

%%----------------------------------------------------------------------------------------------26
One also has a good understanding of homological localization with respect to \bKU.  
Here though, I shall merely provide a summary (proofs can be found in 
Bousfield\footnote[2]{\textit{Topology} \textbf{18} (1979), 257-281; see also 
\textit{J. Pure Appl. Algebra} \textbf{66} (1990), 121-163.}).

[Note: \ There is no need to distinguish between $L_{\bKU}$ and $L_{\bKO}$ since 
$\langle \bKU \rangle = \langle \bKO \rangle$ 
(Meier\footnote[3]{\textit{J. Pure Appl. Algebra} \textbf{14} (1979), 59-71.}).]

\label{17.77}
Put $\bM(p) = \bS(\Z/p\Z)$ $-$then there is a $\bKU_*$-equivalence 
$\bA_p: \Sigma^d \bM(p) \ra \bM(p)$, where $d = 8$ if $p = 2$ $\&$ $d = 2p - 2$ if $p > 2$.  
Using the notation on 
p. \pageref{17.44}, the arrow 
$\bM(p) \ra \bA_p^{-1} \bM(p)$ is a $\bKU_*$-equivalence and 
$\bA_p^{-1} \bM(p)$ is $\bKU_*$-local (
$\implies$ $L_{\bKU}\bM(p)$ $=$ $\bA_p^{-1} \bM(p)$).

[Note: \ Define $\bco\bA_p$ by the exact triangle 
$\Sigma^d \bM(p) \overset{\bA_p}{\lra}$ 
$\bM(p) \ra$ 
$\bco\bA_p \ra$ 
$\Sigma^{d+1}\bM(p)$
$-$then 
$\langle \bKU \rangle = \langle \bigvee\limits_p\bco\bA_p \rangle^c$.]

Remark: $T_{\bKU}$ is smashing and the $\pi_n(L_{\bKU}\bS)$ can be calculated in closed form 
($L_\bKU \bS$ is not connective, e.g., 
$\pi_{-2}(L_\bKU\bS) = \Q/\Z$).

Examples: 
(1)
$L_{\bKU}(\bX \wedge \bM(p)) \approx$ 
$L_{\bKU} \bS \wedge \bX \wedge \bM(p) \approx$ 
$\bX \wedge L_{\bKU} \bS \wedge \bM(p) \approx$ 
$\bX \wedge L_{\bKU} \bM(p) \approx$ 
$\bX \wedge \bA_p^{-1}\bM(p)$; 
(2) 
$L_{\bE \wedge \bM(p)} \bX \approx$ 
$L_{\bM(p)} L_{\bE}\bX$ (cf. p. \pageref{17.45}).\\

\index{Theorem: Bousfield's First KU Theorem}
\index{Bousfield's First KU Theorem}
\textbf{\small BOUSFIELD'S FIRST KU THEOREM} \quad 
Fix an \bX $-$then \bX is $\bKU_*$-local iff $\forall \ p$ $\&$ $\forall \ n$, the arrow 
$[\Sigma^n \bM(p),\bX] \ra [\Sigma^{n+d} \bM(p),\bX]$ induced by $\bA_p$ is bijective or, equivalently, iff $\forall \ p$ $\&$ $\forall \ n$, the arrow 
$\pi_n(\bM(p) \wedge \bX) \ra \pi_{n + d}(\bM(p) \wedge \bX)$ induced by $\bA_p$ is bijective.

[Note: \ Therefore \bX is $\bKU_*$-local iff  
$\pi_*(\bM(p) \wedge \bX) \approx$
$\pi_{*}(\bA_p^{-1}\bM(p) \wedge \bX)$ under the $\bKU_*$-equivalence $\bM(p) \ra \bA_p^{-1}\bM(p)$.]\\

\index{Theorem: Bousfield's Second KU Theorem}
\index{Bousfield's Second KU Theorem}
\textbf{\small BOUSFIELD'S SECOND KU THEOREM} \quad 
Fix an $\bff: \bX \ra \bY$ $-$then $\bff$ is a $\bKU_*$-equivalence iff 
$\bff_*:\pi_*(\bX) \otimes \Q \ra$ $\pi_*(\bY) \otimes \Q$ is bijective and $\forall \ p$, 
$\bff_*:\pi_*(\bA_p^{-1}(\bM(p) \wedge \bX) \ra$
$\pi_*(\bA_p^{-1}(\bM(p) \wedge \bY)$ is bijective.\\

\begingroup%%----------------------------------->>
\fontsize{9pt}{11pt}\selectfont
\textbf{\small FACT}  \  
Let \bku be the connective cover of \bKU $-$then \bku is a ring spectrum (cf. p. \pageref{17.46}) and 
$\bKU \approx$ $\ov{b}_{\bU}^{-1}\bku$ (cf. p. \pageref{17.47}).\\
%dmc note one set of notes has  $\bK\bU \approx$ $\ov{b}_{\bU}^{-1}\bku$ supra
\endgroup %%------------------------------------<<

Fix a prime $p$ $-$then the objects of $\bHSPEC_p$ ($= \bHSPEC_{\{p\}})$ are the 
\un{$p$-local spectra}
\index{p-local spectra} 
and one writes $\bX_p$ in place of $L_{\bS(\Z_p)} \bX$, $\bX_p$ being the 
\un{$p$-localization}
\index{p-localization (spectra)} 
of \bX.  
Example: $\bM(p)$ is $p$-local.

[Note: \ In $\bHSPEC_p$, 
$\bX \wedge_p\bY =(\bX \wedge \bY)_p$ (cf. p. \pageref{17.50}), i.e., 
$\bX \wedge_p \bY = \bX \wedge \bY$ 
($T_{\bS(\Z_p)}$ is smashing ), and $\bS_p$ is the unit.  
Example: $\langle \bS_p \rangle = \langle \bM(p)\rangle \vee \langle \bS(\Q)\rangle$.]\\

%%----------------------------------------------------------------------------------------------27
\label{17.59}
\label{17.72}
\begingroup%%----------------------------------->>
\fontsize{9pt}{11pt}\selectfont
\textbf{\small EXAMPLE}  \ 
Consider $\bKU_p$ $-$then 
Adams\footnote[2]{\textit{SLN} \textbf{99} (1969), 77-98; 
see also Bousfield, \textit{Amer. J. Math.} \textbf{107} (1985), 895-932.}
has shown that there is a splitting 
$\bKU_p \approx \bKU_p(1) \vee$ 
$\Sigma^2 \bKU_p(1) \vee$ 
$\cdots \vee$ 
$\Sigma^{2(p-2)} \bKU_p(1)$ where $\bKU_p(1)$ is a $p$-local spectrum with 
$\pi_*(\bKU_p(1)) \approx$ 
$\Z_p[v_1,v_1^{-1}]$ $(\abs{v_1} = 2(p-1))$.\\
\endgroup %%------------------------------------<<

\begin{proposition} \ %26
Suppose that $\bX_p = 0$ $\forall \ p$ $-$then $\bX = 0$.
\end{proposition}

[$\bX_p  = 0$ $\forall \ p$ $\implies$ 
$\Z_p \otimes \pi_*(\bX) = 0$ $\forall \ p$ $\implies$ 
$\pi_*(\bX) = 0$ (cf. p. \pageref{17.51}). $\implies$ 
$\bX = 0$.]

[Note: \ The converse is trivial.]\\

\label{17.57}
\label{17.83}
\begingroup%%----------------------------------->>
\fontsize{9pt}{11pt}\selectfont
The objects of cpt \bHSPEC$_p$ are the 
\un{$p$-compact spectra}
\index{p-compact spectra}.\\
\endgroup %%------------------------------------<<

\begingroup%%----------------------------------->>
\fontsize{9pt}{11pt}\selectfont
\textbf{\small FACT}  \  
A $p$-local spectrum is $p$-compact iff it is isomorphic to the $p$-localization of a compact spectrum.\\
\endgroup %%------------------------------------<<

\begingroup%%----------------------------------->>
\fontsize{9pt}{11pt}\selectfont
\textbf{\small EXAMPLE}  \ 
Take \mX compact $-$then $\bff:\Sigma^n \bX \ra \bX$ is composition nilpotent iff $\forall \ p$, 
$\bff_p:\Sigma^n\bX_p \ra$ $\bX_p$ is composition nilpotent.
\vspi
[\bff is composition nilpotent iff $\bff^{-1}\bX = 0$ (cf. p. \pageref{17.52}).  
But $\bff^{-1} \bX = 0$ iff $\forall \ p$, $(\bff^{-1}\bX)_p = 0$ (cf. Proposition 26).  And: 
$(\bff^{-1} \bX)_p = \bff_p^{-1} \bX_p$.]\\
\endgroup %%------------------------------------<<

\index{Theorem: BP Theorem}
\index{BP Theorem}
\textbf{\small BP THEOREM} \quad 
Formal group law theory furnishes a canonical idempotent 
$\be_p \in$ 
$[\bMU_p,\bMU_p]$ (the 
\un{Quillen idempotent}
\index{Quillen idempotent}) 
which is a morphism of ring spectra.  Thus, since idempotents split (cf .p. \pageref{17.53}), $\exists$ a commutative ring spectrum \bBP (called the 
\un{Brown-Peterson spectrum at the prime $p$})
\index{Brown-Peterson spectrum at the prime $p$} 
and morphisms 
$\bi:\bBP \ra \bMU_p$, 
$\br:\bMU_p \ra \bBP$ of ring spectra such that 
$\br \circ \bi = \id_{\bBP}$ and $\be_p = \bi \circ \br$.  
\bBP is complex orientable and 
$\bBP^*(\bS) = \bZ_p[v_1, v_2, \ldots]$, where $\abs{v_i} = -2(p^i - 1)$.  And: $\bMU_p$ is isomorphic to a wedge of suspensions of \bBP, hence 
$\langle \bMU_p \rangle = \langle \bBP \rangle$.

[Note: \ The construction is spelled out in 
Adams\footnote[2]{\textit{Stable Homotopy and Generalized Homology}, University of Chicago (1974), 104-116; 
see also Wilson, \textit{CBMS Regional Conference} \textbf{48} (1982), 1-86.}
(a sketch of the underlying ideas is given below).]\\

\begingroup%%----------------------------------->>
\fontsize{9pt}{11pt}\selectfont
Notation: \mA is a commutative $\Z_p$-algebra with unit, $\FGL_A$ is the set of formal group laws over \mA, and 
$\FGL_{A,p}$ is the set of $p$-typical formal group laws over \mA.
\vspi
[Note: \ Initially, it is best to keep the graded picture in the background.]\\
\endgroup %%------------------------------------<<

\begingroup%%----------------------------------->>
\fontsize{9pt}{11pt}\selectfont
\index{Theorem: Cartier's Theorem}
\index{Cartier's Theorem}
\textbf{\small CARTIER'S THEOREM} \quad 
There is an idempotent 
$\epsilon_A:\bFGL_A \ra \bFGL_{A,p}$, functorial in $A$, such that 
$\epsilon_A(\FGL_A) = \FGL_{A,p}$.  Furthermore, there is a natural strict isomorphism 
$F \ra \epsilon_A F$ such that if \mF is $p$-typical, then 
$\epsilon_A F = F$ and $F \ra \epsilon_A F$ is the identity.\\
\endgroup %%------------------------------------<<
%%----------------------------------------------------------------------------------------------28
%% Note duplication removed

\begingroup%%----------------------------------->>
\fontsize{9pt}{11pt}\selectfont
Using this result, one can establish a $p$-typical variant of Lazard's theorem: The functor from the category of commutative $\Z_p$-algebras with unit to the category of sets which sends \mA to FGL$_{A,p}$ is representable.  
Proof: Let $\epsilon_p:L \otimes \Z_p \ra L \otimes \Z_p$ be the homomorphism classifying $\epsilon_{L \otimes \Z_p} F_L$ $-$then $\epsilon_p$ is idempotent, $F_V = \epsilon_{L \otimes \Z_p} F_L$ is defined over $V = \im \epsilon_p$, and $F_V$ is the universal $p$-typical FGL.
\vspi
[Note: \ Structurally, $V = \Z_p[v_1, v_2, \ldots]$, a polynomial algebra on generators $v_i$ od degree $-2(p^i - 1)$.]
\vspi
Remark: To explain the origin of the Quillen idempotent, identify $L \otimes \Z_p$ with $\bMU^*(\bS) \otimes \Z_p$, so 
$F_L \leftrightarrow F_{\bMU}$.  
Let $\phi_p:F_{\bMU} \ra F_V$ be the natural strict isomorphism provided by Cartier, put 
$x_{\bMU_p} = \phi_p(x_{\bMU}) \in \bMU_p^2(\bP^\infty(\C))$ (a complex orientation of $\bMU_p$), and let 
$\be_p:\bMU_p \ra \bMU_p$ be the unique morphism of ring spectra such that 
$\be_p \circ x_{\bMU} = x_{\bMU_p}$ $-$then from the definitions, 
$\be_p \circ \be_p \circ x_{\bMU} = \be_p \circ x_{\bMU}$, hence $\be_p$ is idempotent $\be_p \circ \be_p =$ $\be_p$.
\vspi
[Note: \ \bBP is a commutative ring spectrum with complex orientation $x_{\bBP}$.  
The associated FGL $F_{\bBP}$ is $p$-typical and the map $V \ra \bBP^*(\bS)$ classifying $F_{\bBP}$ is an isomorphism of graded commutative $\Z_p$-algebras with unit.  
Therefore 
$\pi_*(\bMU_p) =$ 
$\pi_*(\bBP) \otimes_{\Z_p} \Z_p[x_1, \ldots, \widehat{x}_{p-1}, x_p, \ldots,\widehat{x}_{p^2-1}, x_{p^2}, \ldots]$.  
Now let $S$ be the set of monomials drawn from 
$\{x_k:k \neq p^i - 1 \ \forall \ i\}$.  Given an $x_I \in S$, write $d_I$ for its degree and call $\bff_I$ the composite 
$\bS^{d_I} \wedge \bBP \ra$ 
$\bMU_p \wedge \bMU_p \ra$ 
$\bMU_p $ $-$then the wedge of the $\bff_I$ defines a morphism 
$\ds\bigvee\limits_{x_I \in S} \Sigma^{d_I} \bBP \ra \bMU_p$ which induces an isomorphism in homotopy.]
\vspi
Rappel: If $F \in \FGL_{A,p}$ and if 
$\phi(x) = \ds\sum\limits_{i \geq 1} \phi_ix^i \in A[[x]]$ with $\phi^\prime(0) = 1$, then the formal group law 
$G(x,y) = \phi(F(\phi^{-1}(x),\phi^{-1}(y)))$ is $p$-typical iff $\phi^{-1}(x)$ has the form 
$x +_F a_1x^p +_F a_2 x^{p^2} +_F \cdots$ $(a_i \in A)$.
\vspi
Set $VT = V[t_1, t_2, \ldots]$, a polynomial algebra on indeterminates $t_i$ $(\abs{t_i} = -2(p^i - 1))$ $-$then the pair 
$(V,VT)$ is a Hopf algebroid over $\Z_p$, i.e., is a cogroupoid object in the category of commutative $\Z_p$-algebras with unit (cf. Proposition 17).  Thus let \mA be a commutative $\Z_p$-algebra with unit.  Denoting by $\bG_{A,p}$ the groupoid whose objects are the $p$-typical formal group laws over $A$ and whose morphisms are the strict isomorphisms, the functor from the category of commutative $\Z_p$-algebras with unit to the category of groupoids which sends \mA to $\bG_{A,p}^\OP$ is represented by $(V,VT)$ .  Indeed, 
$\Hom(V,A) \leftrightarrow \FGL_{A,p} =$ $\Ob\bG_{A,p}$ $(=\Ob\bG_{A,p}^\OP)$ and this identifies the objects.  
Turning to the morphisms, suppose that $f \in \Hom(VT,A)$.  Put $F = (\restr{f}{V})_*F_V$ and let $\phi:F \ra G$ be the morphism 
$\phi^{-1}(x)  = x +_Ff(t_1)x^p +_Ff(t_2)x^{p^2} +_F \cdots$, so $\phi^\OP:G \ra F$ is a strict isomorphism, where 
$G(x,y) = \phi(F(\phi^{-1}(x) ,\phi^{-1}(y)))$ is again $p$-typical.
\vspi
%%----------------------------------------------------------------------------------------------29
[Note: \ $\eta_L$ is the inclusion $V \ra VT$ but there is no simple explicit formula for $\eta_R$.  Incidentally, the groupoid $\bG_{A,p}$ is not split.]
\vspi
To understand the grading on \mV and $VT$, define an action 
$A^\times \times \Ob\bG_{A,p}^\OP \ra$ $\Ob\bG_{A,p}^\OP$ \  by $(u,F) \ra F^u$, where \ 
$F^u(x,y) = uF(u^{-1}x,u^{-1}y)$, \ and define an action 
$A^\times \times \Mor\bG_{A,p}^\OP \ra$ $\Mor\bG_{A,p}^\OP$ \ by 
$(u,\phi^\OP) \ra$ $(\phi^u)^\OP$, where $\phi^u(x) = u\phi(u^{-1}X)$ $-$then this action grades $V$ and $V T$ and one can check that $\abs{v_i} =$ $-2(p^i - 1) =$ $\abs{t_i}$.  
Because the five arrows of structure $\eta_R$, $\eta_L$, $\epsilon$, $\Delta$, $c$ are gradation preserving, it follows that 
$(V,VT)$ is a graded Hopf algebroid over $\Z_p$.
\vspi
[Note: \ Therefore $(V,VT)^\OP$ is but another name for 
$(\bBP_*(\bS)$, $\bBP_*(\bBP))$ and $\bBP_*(\bBP)$ is a graded free $\bBP_*(\bS)$-module.]\\
\endgroup %%------------------------------------<<

\begingroup%%----------------------------------->>
\fontsize{9pt}{11pt}\selectfont
\index{BP Nilpotence Technology}
\textbf{\small FACT \ (\un{$\bB\bP$ Nilpotence Technology})} \  
Let \bE be a $p$-local ring spectrum and consider the Hurewicz homomorphism 
$\bS_*(\bE) \ra \bBP_*(\bE)$ (cf. p. \pageref{17.54}ff) 
$-$then the homogeneous elements of its kernel are nilpotent 
(Devinatz-Hopkins-Smith\footnote[2]{\textit{Ann. of Math.} \textbf{128} (1988), 207-241.}
).\\
\endgroup %%------------------------------------<<

\label{17.64}
\begingroup%%----------------------------------->>
\fontsize{9pt}{11pt}\selectfont
Application:  If \bX is $p$-compact and if $\bff:\Sigma^n \bX \ra \bX$ is an arrow such that 
$\bBP_*(\bff) = 0$, then \bff is composition nilpotent (cf. p. \pageref{17.55}ff).\\
\endgroup %%------------------------------------<<

\begingroup%%----------------------------------->>
\fontsize{9pt}{11pt}\selectfont
Application: If \bX is $p$-compact and \bY is $p$-local and if $\bff:\bX \ra \bY$ is an arrow such that 
$\id_{\bBP} \wedge \bff = 0$, then \bff is smash nilpotent (cf. p. \pageref{17.56}).
\vspi
[Note: \ Write $\bX = \ov{\bX}_p$, where $\ov{\bX}$ is compact (cf. p. \pageref{17.57}) $-$then 
$\hom(\bX,\bY) \approx$
$\hom(\ov{\bX},\bY) \approx$
$D\ov{\bX} \wedge \bY \approx$ 
$D\ov{\bX} \wedge \bS_p \wedge \bY \approx$ 
$\hom(\ov{\bX},\bS_p) \wedge \bY \approx$ 
$\hom(\ov{\bX}_p,\bS_p) \wedge \bY \approx$ 
$\hom(\bX,\bS_p) \wedge \bY$ and $\hom(\bX,\bS_p)$ is the dual of \bX in $\bHSPEC_p$.]\\
\endgroup %%------------------------------------<<

There are two particularly important classes of spectra attached to \bBP, viz. the $\bK(n)$ and the $\bP(n)$ 
$(0 < n < \infty)$ with 
$\pi_*(\bK(n)) = \F_p[v_n,v_n^{-1}]$ and
$\pi_*(\bP(n)) = \F_p[v_n,v_{n+1}, \ldots]$.  
Both are $p$-local ring spectra (commutative if $p > 2$) and \bBP-module spectra but the exact details of their construction need not detain us since all that really counts are the properties possessed by them, which will be listed below.  
Example: $\bP(1) \approx \bBP \wedge \bM(p)$.

[Note: \ The theory has been surveyed by 
W\"urgler\footnote[2]{\textit{SLN} \textbf{1474} (1991), 111-138.}
.]

The role of the $\bP(n)$ is basically technical.  Since $v_n \in \pi_{2(p^n-1)}(\bP(n))$, one can form
$\ov{v}_n:\Sigma^{2(p^n-1)}\bP(n) \ra \bP(n)$ (cf. p. \pageref{17.58}) $-$then there is an exact triangle 
$\Sigma^{2(p^n-1)}\bP(n) \overset{\ov{v}_n}{\lra}$ 
$\bP(n) \ra$ 
$\bP(n+1) \ra$ 
$\Sigma^{2p^n-1}\bP(n)$.  
Moreover, 
$\langle \bK(n) \rangle =$ 
$\langle \ov{v}_n^{-1}\bP(n) \rangle$ and 
$\bH(\F_p) \approx$ 
$\tel(\bP(1) \ra \bP(2) \ra \cdots)$.  
On the other hand, 
$\langle \bBP \rangle =$ 
$\langle \bH(\Q) \rangle \vee \langle \bP(1) \rangle$ and 
$\langle \bP(n) \rangle =$ 
$\langle \bK(n)\rangle \vee \langle \bP(n+1)\rangle$ (cf. $\S 15$, Proposition 43), hence 
$\langle \bBP \rangle =$ 
$\langle \bH(\Q) \rangle \vee \langle \bK(1) \rangle \vee \cdots \vee
%%----------------------------------------------------------------------------------------------30
\langle \bK(n) \rangle \vee \langle \bP(n+1)\rangle$.  
In addition, 
$\langle \bH(\Q) \rangle \wedge \langle \bP(1) \rangle = \langle 0 \rangle$, 
$\langle \bK(i) \rangle \wedge \langle \bP(n+1) \rangle = \langle 0 \rangle$
$(i = 1, \ldots, n)$.

By contrast, $\bK(n)$ (called the 
\un{$n^\text{th}$ Morava K-theory at the prime $p$}
\index{n$^\text{th}$ Morava K-theory at the prime $p$}) 
is a major player.

\indent\indent (Mo$_1$) \quad $\bK(n)$ is a skew field object in \bHSPEC.

[This is because the homogeneous elements of $\pi_*(\bK(n))$ are invertible (cf. $\S 15$, Proposition 42).]

\indent\indent (Mo$_2$) \quad $\forall$ \bX, $\bK(n) \wedge \bX$ is isomorphic to a wedge of suspensions of $\bK(n)$.

[
$\bK(n) \wedge \bX$ is a $\bK(n)$-module, thus the assertion follows from the definition of a skew field object (to accommodate 
$\bK(n) \wedge \bX = 0$, use the empty wedge).]

\indent\indent (Mo$_3$) \quad $\forall$ \bX $\&$ $\forall$ \bY, 
$\bK(n)_*(\bX) \otimes_{\bK(n)_*(\bS)} \bK(n)_*(\bY) \approx$ $\bK(n)_*(\bX \wedge \bY)$.

[This is a special case of Proposition 10.]

\indent\indent (Mo$_4$) \quad 
$\langle \bK(n) \rangle \wedge \langle \bK(m)  \rangle = \langle 0 \rangle$ $(m \neq n)$.

[Suppose that $n < m$ $-$then 
$\langle \bK(m)  \rangle \leq \langle \bP(m)  \rangle \leq \langle \bP(n+1)  \rangle$ and 
$\langle \bK(n)  \rangle \wedge \langle \bP(n+1)  \rangle =$  $\langle 0 \rangle$.]

\indent\indent (Mo$_5$) \quad 
$\langle \bH(\Q) \rangle \wedge \langle \bK(n) \rangle= \langle 0 \rangle$ $\&$
$\langle \bH(\F_p) \rangle \wedge \langle \bK(n) \rangle= \langle 0 \rangle$.

[$\langle \bH(\Q) \rangle \wedge \langle \bP(1) \rangle = \langle 0 \rangle$ and 
$\langle \bK(n) \rangle \leq$ $\langle \bP(n) \rangle$ $\implies$ 
$\langle \bH(\Q) \rangle \wedge \langle \bK(n) \rangle= \langle 0 \rangle$
And: 
$\bH(\F_p) \approx$ 
$\tel(\bP(1) \ra \bP(2) \ra \cdots)$ $\implies$ 
$\langle\bH(\F_p)\rangle \leq \langle \bP(n+1) \rangle$ $\implies$ 
$\langle \bH(\F_p) \rangle \wedge \langle \bK(n) \rangle= \langle 0 \rangle$.]

\indent\indent (Mo$_6$) \quad $\forall$ compact $\bX$, 
$\bK(n)_*(\bX) \approx \bK(n)_*(\bS) \otimes_{\F_p} \bH_*(\bX;\F_p)$ $\forall \ n \gg 0$.

[Apply the Atiyah-Hirzebruch spectral sequence.]

Remarks: 
(1) $\bK(n)$ is complex orientable if $p$ is odd; 
(2) $\bK(1)$ can be identified with $\bKU_p(1) \wedge \bM(p)$ (cf. p. \pageref{17.59}).\\

\label{17.73}
\label{18.36}
\begingroup%%----------------------------------->>
\fontsize{9pt}{11pt}\selectfont
\index{Algebraic K-Theory (example)}
\textbf{\small EXAMPLE \ (\un{Algebraic K-Theory})} \  Suppose that \mA is a ring with unit and let $\bW A$ be the 
$\Omega$-prespectrum attached to \mA by algebraic K-theory (cf. p. \pageref{17.60}).  Consider 
$\bK A = e M \bW A$ $-$then 
Mitchell\footnote[2]{\textit{K-Theory} \textbf{3} (1990), 607-626.}
has shown that $\forall \ p$ $\&$ $\forall \ n \geq 2$, the connective cover of $\bK A$ is $\bK(n)_*$-acyclic.\\
\endgroup %%------------------------------------<<

\begingroup%%----------------------------------->>
\fontsize{9pt}{11pt}\selectfont
\textbf{\small FACT}  \  Let $\bk(n)$ be the connective cover of $\bK(n)$ $-$then $\bk(n)$ is a ring spectrum 
(cf. p. \pageref{17.61}) and 
$\bK(n) \approx \ov{v}_n^{-1}\bk(n)$ (cf. p. \pageref{17.62}).
\vspi
[Note: \ There is an exact triangle 
$\Sigma^{2(p^n-1)}\bk(n) \overset{\ov{v}_n}{\lra}$ 
$\bk(n) \ra$ 
$\bH(\F_p) \ra$ 
$\Sigma^{2p^n-1}\bk(n)$, so by $\S 15$, Proposition 43, 
$\langle \bk(n)\rangle = \langle \bH(\F_p)\rangle \vee \langle \bK(n)\rangle$.]\\
\endgroup %%------------------------------------<<

\label{17.71}
\textbf{\small LEMMA}  \  Any retract of a $\bK(n)$-module is a $\bK(n)$-module.\\

\begingroup%%----------------------------------->>
\fontsize{9pt}{11pt}\selectfont
\textbf{\small EXAMPLE}  \ 
A spectrum \bY is 
\un{indecomposable}
\index{indecomposable (spectrum)} 
if it has no nontrivial direct summands, i.e., 
$\bY \approx \bX \vee \bZ$ $\implies$ $\bX = 0$ or $\bZ = 0$.  
Since idempotents split (cf. p. \pageref{17.63}), \bY 
is indecomposable iff $[\bY,\bY]$ has no nontrivial idempotents.  
Example: $\bK(n)$ is indecomposable.
\vspi
%%----------------------------------------------------------------------------------------------31
[Note: \ One can also prove that \bBP is indecomposable.]\\
\endgroup %%------------------------------------<<

\label{17.70} %dmc mnft
\label{17.84} %dmc mnft
Notation: For uniformity of statement, it is convenient to put 
$\bK(0) = \bH(\Q)$, 
$\bK(\infty) = \bH(\F_p)$.\\

\begingroup%%----------------------------------->>
\fontsize{9pt}{11pt}\selectfont
Hovey\footnote[3]{\textit{Contemp. Math.} \textbf{181} (1995), 230.}
has shown that $\langle\bK(n)\rangle$ is minimal if $n < \infty$ (but this is false if $n = \infty$).\\
\endgroup %%------------------------------------<<

\textbf{\small LEMMA}  \  
Given $\bff:\bX \ra \bY$, suppose that $\bK(n)_*(\bff) = 0$, where $n \in [0,\infty]$ $-$then the composite 
$\bX \overset{\bff}{\lra}$ 
$\bY \approx$ 
$\bS \wedge \bY \ra$
$\bK(n) \wedge \bY$ vanishes.

[For any $\bK(n)$-module \bE, 
$
\begin{cases}
\ \bE^*(\bX) \approx \Hom_{\pi_*(\bK(n))}(\bK(n)_*(\bX),\pi_*(\bE))\\
\ \bE^*(\bY) \approx \Hom_{\pi_*(\bK(n))}(\bK(n)_*(\bY),\pi_*(\bE))
\end{cases}
,
$
hence the induced map 
$\bE^*(\bY) \ra \bE^*(\bX)$ is the zero map.  Now specialize to $\bE = \bK(n) \wedge \bY$.]\\

\begin{proposition} \ %27
If \bX is $p$-compact and \bY is $p$-local and if $\bff:\bX \ra \bY$ is an arrow such that 
$\bK(n)_*(\bff) = 0$ $\forall \ n \in [0,\infty]$, then \bff is smash nilpotent.
\end{proposition}

[It is enough to prove that $\id_{\bBP} \wedge \bff^{(k)} = 0$ $(\exists \ k \gg 0)$ 
(cf. p. \pageref{17.64}), and for this, one can take $\bX = \bS_p$.  So, passing to $\bY_{\bff}^{(\infty)}$ (defined by $\bS_p$ instead of \bS
(cf. p. \pageref{17.65})), it suffices to show that 
$\bBP \wedge \bY_{\bff}^{(\infty)} = 0$.  
But 
$\langle \bBP \rangle =$ 
$\langle \bK(0) \rangle \vee \cdots \vee$ 
$\langle \bK(n) \rangle \vee \langle \bP(n+1)\rangle$
and from our hypotheses and the lemma, 
$\bK(m) \wedge \bY_{\bff}^{(\infty)} = 0$ $(m \leq n)$, thus we are left with proving that  
$\bP(n) \wedge \bY_{\bff}^{(\infty)} = 0$ $(n \gg 0)$, which however is clear since 
$\bH(\F_p) \wedge \bY_{\bff}^{(\infty)} = 0$ and 
$\bH(\F_p) \approx$ 
$\tel(\bP(1) \ra \bP(2) \ra \cdots)$.]\\

Application: If $\bE \neq 0$ is a $p$-local ring spectrum, then for some $n \in [0,\infty]$, 
$\bK(n)_*(\bE) \neq$ $0$.

[Consider the unit $\bS_p \ra \bE$.]\\

\begingroup%%----------------------------------->>
\fontsize{9pt}{11pt}\selectfont
Let \bR be a ring spectrum $-$then \bR is said to 
\un{detect nilpotence}
\index{detect nilpotence (ring spectrum)} 
if for any ring spectrum \bE, the homogeneous elements of the kernel of the Hurewicz homomorphism 
$\bS_*(\bE) \ra \bR_*(\bE)$ are nilpotent.  
Example: \bMU detects nilpotence (cf. p. \pageref{17.66}).\\
\endgroup %%------------------------------------<<

\begingroup%%----------------------------------->>
\fontsize{9pt}{11pt}\selectfont
\textbf{\small LEMMA}  \  
\bR detects nilpotence iff for all compact \bX and any $\bff:\bX \ra \bY$ such that $\id_{\bR} \wedge \bff = 0$, \bff is smash nilpotent.
\vspi
[Necessity: Argue as on p. \pageref{17.67}, with \bMU replaced by \bR.
\vspi
Sufficiency: Given a ring spectrum \bE, fix a homogeneous element $\bff:\bS^n \ra \bE$ in the kernel of the Hurewicz homomorphism $\bS_*(\bE) \ra \bR_*(\bE)$ $-$then $\id_{\bR} \wedge \bff = 0$, so \bff is smash nilpotent, thus nilpotent.]
\vspi
%%----------------------------------------------------------------------------------------------32
Remark: For a compact $\bX$, 
$\bff:\bX \ra \bY$ is smash nilpotent iff 
$\ov{\bff}:\bS \ra D\bX \wedge \bY$ is smash nilpotent (cf. \pageref{17.68}).  This said, the problem of determining the smash nilpotency of $\bff:\bS \ra \bY$ is local, i.e., one has only to check that $\bff_p:\bS_p \ra \bY_p$ is smash nilpotent 
$\forall \ p$.  
Proof:  $\bff:\bS \ra \bY$ is smash nilpotent iff $\bY_{\bff}^{(\infty)} = 0$ (cf. \pageref{17.69}).  But 
$\bY_{\bff}^{(\infty)} = 0$ iff $(\bY_{\bff}^{(\infty)})_p = 0$ $\forall \ p$ (cf. Proposition 26).  And: 
$(\bY_{\bff}^{(\infty)})_p = \bY_{\bff_p}^{(\infty)}$.\\
\endgroup %%------------------------------------<<

\begingroup%%----------------------------------->>
\fontsize{9pt}{11pt}\selectfont
\textbf{\small EXAMPLE}  \ 
A ring spectrum \bR detects nilpotence iff $\forall \ p$ $\&$ $\forall \ n \in [0,\infty]$, $\bK(n)_*(\bR) \neq 0$.
\vspi
[Consider an $\bff:\bS \ra \bY$ such that $\id_{\bR} \wedge \bff = 0$.  
Fixing $p$, one has 
$\bK(n)_*(\bff_p) = 0$ $\forall \ n \in [0,\infty]$ ($\bK(n) \wedge \bR$ is isomorphic to a wedge of suspensions of $\bK(n)$), thus by Proposition 27, $\bff_p$ is smash nilpotent.  Therefore \bR detects nilpotence.]\\
\endgroup %%------------------------------------<<

\begingroup%%----------------------------------->>
\fontsize{9pt}{11pt}\selectfont
\textbf{\small FACT}  \  
Suppose that \bE is a skew field object in \bHSPEC $-$then \bE is isomorphic to a wedge of suspensions of some $\bK(n)$ 
$(\exists \ n \in [0,\infty])$.
\vspi
[$\exists \ p$: $\bE_p \neq 0$ (cf. Proposition 26) $\implies$ 
$\bK(n)_*(\bE) \neq 0$ $(\exists \ n \in [0,\infty])$ (cf. p. \pageref{17.70}).  Since $\bK(n)$ and \bE are both skew field objects, $\bK(n) \wedge \bE \neq 0$ is simultaneously a wedge of suspensions of $\bK(n)$ and a wedge of suspensions of \bE.  
Deduce that \bE is a retract of a wedge of suspensions of $\bK(n)$, hence is a $\bK(n)$-module 
(cf. p. \pageref{17.71}).]\\
\endgroup %%------------------------------------<<

\begingroup%%----------------------------------->>
\fontsize{9pt}{11pt}\selectfont
A skew field object in \bHSPEC is said to be 
\un{prime}
\index{prime (skew field object in \bHSPEC)} 
if it is indecomposable.  
The $\bK(n)$ $(n \in [0,\infty])$ for $p \in \bPi$ are prime and the preceding result implies that, up to isomorphism, they are the only primes in \bHSPEC.\\
\endgroup %%------------------------------------<<

\begingroup%%----------------------------------->>
\fontsize{9pt}{11pt}\selectfont
\textbf{\small EXAMPLE}  \ 
Suppose that $p$ is odd $-$then $\bKU_p(1) \wedge \bM(p)$ is a field object (being isomorphic to 
$\bK(1) \vee \Sigma^2\bK(1) \vee \cdots \Sigma^{2(p-2)}\bK(1)$ (cf. p. \pageref{17.72})) but it is not prime.\\
\endgroup %%------------------------------------<<

\begin{proposition} \ %28
Fix a prime $p$ $-$then $\bH(\F_p)$ is $\bK(n)_*$-acyclic $(n \in [0,\infty[)$.
\end{proposition}

[Trivially, $\bH(\Q) \wedge \bH(\F_p) = 0$.  Proceeding by contradiction, assume that 
$\bK(n) \wedge \bH(\F_p) \neq 0$ for some $n \in [1,\infty[$.  Since $\bH(\F_p)$ is a field object, $\bH(\F_p)$ is isomorphic to a wedge of suspensions of $\bK(n)$ (cf. supra), an impossibility.]\\

\begingroup%%----------------------------------->>
\fontsize{9pt}{11pt}\selectfont
\textbf{\small FACT}  \  
Let \bX be a spectrum with the property that $\exists \ N$: $\pi_n(\bX) = 0$ $(n > N)$ $-$then \bX is 
$\bK(n)_*$-acyclic $(n \in [1,\infty])$.
\vspi
[Using Proposition 28, prove it first under the assumption that $\pi_*(\bX)$ is torsion.  To handle the general case, smash 
$\bS \overset{p}{\ra}$ 
$\bS \ra$
$\bM(p) \ra \Sigma \bS$ with $\bK(n) \wedge \bX$ to see that $\pi_*(\bK(n) \wedge \bX)$ injects into 
$\pi_*(\bK(n) \wedge \bX \wedge \bM(p))$.  But Proposition 19 implies that $\bX \wedge \bM(p)$, like \bX, is ``bounded above'' and 
$\pi_n(\bX \wedge \bM(p))$ is torsion.]
\vspi
[Note: \ In particular, $\bK(n) \wedge \bH(\pi) = 0$ $(n \in [1,\infty[)$, $\pi$ any abelian group.]\\
\endgroup %%------------------------------------<<

\begingroup%%----------------------------------->>
\fontsize{9pt}{11pt}\selectfont
Application: If \bX is a spectrum and $\bx\ (= \tau^{\leq 0} \bX)$ is its connective cover, then the arrow 
$\bx \ra \bX$ is a $\bK(n)_*$-equivalence $(n \in [1,\infty[)$.
\vspi
%%----------------------------------------------------------------------------------------------33
[For $\bK(n) \wedge \bF = 0$, where \bF is defined by the exact triangle 
$\bF \ra \bx \ra \bX \ra \Sigma\bF$.]
\vspi
[Note: \ Let \mA be a ring with unit $-$then $\forall \ p$ $\&$ $\forall \ n \geq 2$, the connective cover of $\bK A$ is 
$\bK(n)_*$-acyclic (cf. \pageref{17.73}), hence so is $\bK A$ itself.]\\
\endgroup %%------------------------------------<<

\begin{proposition} \ %29
If \bX is $p$-compact and if $\bff:\Sigma^d \bX \ra \bX$ is an arrow such that $\bK(n)_*(\bff) = 0$ $\forall \ n \in [0,\infty[$, then 
\bff is composition nilpotent.
\end{proposition}

[This is a consequence of Proposition 27 (one doesn't need the $n = \infty$ case).]\\

\begingroup%%----------------------------------->>
\fontsize{9pt}{11pt}\selectfont
\textbf{\small EXAMPLE}  \ 
If \bX is $p$-compact and if  $\bK(n)_*(\bX) = 0$ ($\forall \ n \in [0,\infty[$), then $\bX = 0$ (in Proposition 29, take 
$\bff = \id_{\bX}$.)
\vspi
[Note: \ Accordingly, if \bX is compact and if 
$\forall \ p$ $\&$ $\forall \ n \in [0,\infty[$, $\bK(n)_*(\bX) = 0$, then $\bX = 0$.  
In fact, 
$\bK(n)_*(\bX) = $
$\pi_*(\bK(n) \wedge \bX) =$ 
$\pi_*(\bK(n) \wedge \bX_p) =$ 
$\bK(n)_*(\bX_p)$ $\implies$ 
$\bX_p = 0$ $\forall \ p$ $\implies$  $\bX = 0$ (cf. Proposition 26).]\\
\endgroup %%------------------------------------<<

Given a prime $p$, write $\bC(0)$ for $\cptx \bHSPEC_p$ and let $\bC(n)$ be the thick subcategory of $\bC(0)$  whose objects are those \bX such that 
$\bK(n-1)_*(\bX) = 0$ $(n \in [1,\infty[)$ (conventionally, the objects of $\bC(\infty)$ are the zero objects) $-$then 
$\bC(n+1) \subset \bC(n)$, i.e., $\bK(n)_*(\bX) = 0$ $\implies$ $\bK(n-1)_*(\bX) = 0$ 
(Ravenel\footnote[2]{\textit{Amer. J. Math.} \textbf{106} (1984), 351-414 (cf. 366-367).}) 
and the containment is strict (Mitchell\footnote[3]{\textit{Topology} \textbf{24} (1985), 227-246; 
see also Palmieri-Sadofsky, \textit{Math. Zeit.} \textbf{215} (1994), 477-490.}).

[Note: \ A $p$-compact \bX is said to have 
\un{type $n$}
\index{type $n$ ($p$-compact spectrum)}
 if 
$n = \min\{m:\bK(m)_*(\bX) \neq 0\}$ ($\bX = 0$ has type $\infty$).  The objects of type $n$ are the objects in $\bC(n)$ which are not in $\bC(n+1)$.  
Examples: 
(1) $\bS_p$ has type 0; 
(2) $\bM(p)$ has type 1; 
(3) $\bco\bA_p$ has type 2.]\\

\textbf{\small LEMMA}  \  
Let \bX be a $p$-compact spectrum, \bE a $p$-local ring spectrum.  Suppose given a $p$-local spectrum \bZ and a morphism 
$\bff:\bX \ra \bE \wedge \bZ$ in $\bHSPEC_p$ such that 
$\bK(n)_*(\bff) =$ 
$0$ $\forall \ n \in [0,\infty]$ $-$then the composite 
\begin{tikzcd}%[sep=small]
{\bX^{(N)}} \ar{r}{\bff^{(N)}} &{(\bE \wedge \bZ)^{(N)}}
\end{tikzcd}
$\approx \bE^{(N)} \wedge \bZ^{(N)}$
$\ra \bE \wedge \bZ^{(N)}$ 
vanishes if $N \gg 0$ (cf. Proposition 27).\\

Application: Let \bX, \bY be $p$-compact spectra.  
Suppose given a $p$-local spectrum \bZ and a morphism 
$\bff:\bX \ra \bZ$ in $\bHSPEC_p$ such that $\bK(n)_*(\bff \wedge \id_{\bY}) = 0$ $\forall \ n \in [0,\infty]$ $-$then 
$\bff^{(N)} \wedge \id_{\bY}: \bX^{(N)} \wedge \bY \ra \bZ^{(N)} \wedge \bY$ vanishes if $N \gg 0$.

[One has 
$[\bX \wedge \bY,\bZ \wedge \bY] \approx [\bX,\hom(\bY,\bZ \wedge \bY)]$.  But \bY is $p$-compact so 
$\hom(\bY,\bZ \wedge \bY) \approx$ 
$\hom(\bY,\bS_p) \wedge \bY \wedge \bZ \approx$ 
$\hom(\bY,\bY) \wedge \bZ$.  Now specialize the lemma to $\bE = \hom(\bY,\bY)$.]\\

\index{Theorem: Thick Subcategory Theorem}
\index{Thick Subcategory Theorem}
\textbf{\small THICK SUBCATEGORY THEOREM} \quad 
The thick subcategories of $\bC(0)$ are the $\bC(n)$.

%%----------------------------------------------------------------------------------------------34
[Fix a thick subcategory of \bC of $\bC(0)$ and let $n_{\bC} = \min\{n:\bC(n) \subset \bC\}$.  
Claim: If $\bX \in \Ob\bC$ has type n, then $\bC(n) \subset \bC$ ($\implies$ $\bC = \bC(n_{\bC})$).    
Define \bF, \bff, by the exact triangle 
$\bF \overset{\bff}{\lra}$ 
$\bS_p \ra$ 
$\hom(\bX,\bX) \ra \Sigma \bF$.  
Because $\bHSPEC_p$ is monogenic ($\implies$ unital), $\hom(\bX,\bS_p)$ is $p$-compact, so 
$\hom(\bX,\bX) \approx$ $\hom(\bX,\bS_p) \wedge \bX$ $\in \Ob\bC$ (\bC being thick 
(cf. p. \pageref{17.74})).  
Putting $\bC_{\bff} = \hom(\bX,\bX)$, one thus concludes that $\bF \wedge \bC_{\bff} \in \Ob\bC$ (here again the assumption that \bC is thick comes in).  
But there is an exact triangle 
$\bF \wedge \bC_{\bff^{(N-1)}} \ra$
$\bC_{\bff^{(N)}} \ra$
$\bC_{\bff} \wedge \bS_p^{(N-1)} \ra$ 
$\Sigma(\bF \wedge \bC_{\bff^{(N-1)}})$ 
(cf. p. \pageref{17.75}), from which inductively, 
$\bC_{\bff^{(N)}} \in \Ob\bC$ $\forall \ N \geq 1$.  
Take a \bY in $\bC(n)$.  
Since $\bK(m)_*(\bff \wedge \id_{\bY}) = 0$ 
$\forall \ m \in [0,\infty]$ $(\bK(m)_*(\bX) \neq 0$ 
$\forall \ m \geq n)$, $\forall \ N \gg 0$, $\bff^{(N)} \wedge \id_{\bY} = 0$ 
(cf. supra).  Working with the exact triangle 
\begin{tikzcd}%[sep=small]
{\bF^{(N)} \wedge \bY} \ar{rr}{\bff^{(N)} \wedge \id_{\bY}}
&&{\bS_p^{(N)} \wedge \bY} \ar{r}
&{\bC_{\bff^{(N)}} \wedge \bY} \ar{r}
&{\Sigma(\bF^{(N)} \wedge \bY)}
\end{tikzcd}
, it then follows that 
$\bC_{\bff^{(N)}} \wedge \bY \approx$ 
$(\bS_p^{(N)} \wedge \bY) \vee \Sigma(\bF^{(N)} \wedge \bY)$ (cf. p. \pageref{17.76}) 
And: 
$\bC_{\bff^{(N)}} \wedge \bY \in \Ob\bC$ 
$\implies$ 
$\bS_p^{(N)} \wedge \bY \in \Ob\bC$
$\implies$ 
$\bY \in \Ob\bC$.]\\

\begingroup%%----------------------------------->>
\fontsize{9pt}{11pt}\selectfont
\textbf{\small EXAMPLE}  \ 
Fix a spectrum \bE and write $\bACY_p(\bE)$ 
\index{$\bACY_p(\bE)$} 
for the class of $p$-compact \bX such that 
$\bE \wedge \bX = 0$ $-$then $\bACY_p(\bE)$ is the object class of a thick subcategory of $\bC(0)$, hence 
$\bACY_p(\bE) = \Ob\bC(n)$ for some $n$.\\
\endgroup %%------------------------------------<<

\begingroup%%----------------------------------->>
\fontsize{9pt}{11pt}\selectfont
\index{Class Invariance Principle}
\textbf{\small FACT \ (\un{Class Invariance Principle})} \  
Let \bX, \bY be $p$-compact.  Suppose that \bX has type $n$ and \bY has type $m$ $-$then 
$\langle \bX \rangle = \langle \bY \rangle$ iff $n = m$.
\vspi
[The necessity is obvious.  
To establish the sufficiency, note that the full, isomorphism closed subcategory of 
$\cptx \bHSPEC_p$ whose objects are the \bZ with $\langle \bZ \rangle \leq \langle \bX \rangle$ is thick.]\\
\endgroup %%------------------------------------<<

Given a prime $p$ and a $p$-compact \bX, an arrow 
$\bff:\Sigma^d\bX \ra \bX$ is said to be a 
\un{$v_n$-map}
\index{v$_n$-map} 
$(n \in [0,\infty])$ if $\bK(n)_*(\bff)$ is an isomorphism and 
$\bK(m)_*(\bff) = 0$ $\forall \ m \neq n$ $(m \in [0,\infty])$ (cf. Proposition 29).  
Example: $\bX \overset{p}{\ra} \bX$ is a $v_0$-map.

[Note: \ For $m \gg 0$, $\bK(m)_*(\bff) = \bH(\F_p)(\bff) \otimes_{\F_p} \id_{\bK(m)_*}$ $\implies$ 
$\bH(\F_p)(\bff) = 0$.]

Example: $\bA_p:\Sigma^d \bM(p) \ra \bM(p)$ is a $v_1$-map ($d = 8$ if $p = 2$ $\&$ $d = 2p - 2$ if $p > 2$ 
(cf. p. \pageref{17.77})).\\

\begin{proposition} \ %30
Let \bX be $p$-compact and fix $n \geq 1$.  Suppose that \bX admits a $v_n$-map $-$then \bX belongs to $\bC(n)$, i.e., 
$\bK(n-1)_*(\bX) = 0$.
\end{proposition}

[Defining \bY by the exact triangle 
$\Sigma^d \bX \overset{\bff}{\lra} \bX \ra \bY \ra \Sigma^{d+1}\bX$, one has 
$\bK(n)_*(\bY) = 0$, thus 
$0 = \bK(n-1)_*(\bY) =$
$\bK(n-1)_*(\bX) \oplus \bK(n-1)_*(\Sigma^{d+1}\bX)$ $\implies$ $\bK(n-1)_*(\bX) = 0$.]\\

I shall omit the proof of the following result as it is quite involved.\\

%%----------------------------------------------------------------------------------------------35
\index{Theorem: Hopkins-Smith Existence Theorem}
\index{Hopkins-Smith Existence Theorem}
\textbf{\small HOPKINS-SMITH\footnote[2]{\textit{Ann. of Math.} \textbf{148} (1998), 1-49; 
see also, \textit{Ravenel, Nilpotence and Periodicity in Stable Homotopy Theory}, Princeton University Press (1992), 53-68.}
EXISTENCE THEOREM} \quad 
Given $n \geq 1$, $\exists$ a $p$-compact \bX of type $n$ which admits a $v_n$-map.

[Note: \ In fact, \bX admits a $v_n$-map 
$\bff:\Sigma^{P^N2(p^n - 1)}\bX \ra \bX$ such that $\bK(n)_*(\bff) = v_n^{p^N}$ $(N \gg 0)$.]\\

\label{17.79}
Remark: A $p$-compact \bX admits a $v_n$-map iff \bX is in $\bC(n)$.  
To see this, call $\bV_n$ the full, isomorphism closed subcategory of 
$\bC(0)$ (= $\cptx \bHSPEC_p$) whose objects are those \bX which admit a $v_n$-map.  
Owing to Proposition 30, $\bC(n) \supset \bV_n$.  
On the other hand, $\bX \overset{0}{\ra} \bX$ is a $v_n$-map if $\bK(n)_*(\bX) = 0$, so 
$\bV_n \supset \bC(n+1)$.  
However, $\bV_n$ is thick (cf. p. \pageref{17.78}), hence by the thick subcategory theorem, either 
$\bV_n = \bC(n)$ or 
$\bV_n = \bC(n+1)$.  
Since the containment $\bC(n+1) \subset \bC(n)$ is proper, the Hopkins-Smith existence theorem eliminates the 
second possibility.

Notation: Write $[\bX,\bX]_*$ for the graded ring with unit defined by 
$[\bX,\bX]_n = [\Sigma^n\bX,\bX]$ (cf. Proposition 1).

[Note: \ An arrow $\bff:\Sigma^n\bX \ra \bX$ is composition nilpotent iff $\bff^k = 0$ for some $k$ or still, is nilpotent when viewed as an element of $[\bX,\bX]_*$.]\\

\begin{proposition} \ %31
Let \bX be $p$-compact and fix $n \geq 1$.  Suppose that 
$\bff:\Sigma^d \bX \ra \bX$, 
$\bg:\Sigma^e \bX \ra \bX$ 
are $v_n$-maps $-$then $\exists \ i, j$ : $\bff^i = \bg^j$.\\
\end{proposition}

\begingroup%%----------------------------------->>
\fontsize{9pt}{11pt}\selectfont
The proof of Propositon 31 rests on the following considerations.
\vspi
Given a $p$-compact \bX in $\bC(n)$ $(n \geq 1)$, put 
$\bR\bX = \hom(\bX,\bS_p) \wedge \bX$ $(\approx \hom(\bX,\bX))$ $-$then $\bR\bX$ is a $p$-compact ring spectrum, 
$\bH(\Q) \wedge \bR\bX = 0$, and $[\bX,\bX]_* \approx \pi_*(\bR\bX)$.
\vspi
Definition: An element $\alpha \in \pi_d(\bR\bX)$ is a 
\un{$v_n$-element}
\index{v$_n$-element} 
provided that its image $\bK(m)_*(\alpha)$ under the Hurewicz homomorphism 
$\bS_*(\bR\bX) \ra \bK(m)_*(\bR\bX)$ is a unit if $m = n$ and vanishes otherwise ($m \in [1,\infty[$).
\vspi
[Note: \ By contrast, if $\bK(m)_*(\alpha) = 0$ $\forall \ m \in [0,\infty[$, then $\alpha$ is nilpotent.]
\vspi
Example: The adjoint $\ov{\bff} \in \pi_d(\bR\bX)$ of a $v_n$-map $\bff \in [\bX,\bX]_d$ is a $v_n$-element (and conversely).
\vspi
Claim: Fix a $v_n$-element $\alpha$ $-$then $\exists \ i$ such that $\bK(n)_*(\alpha^i) = v_n^N$ for some $N$.
\vspi
[The ungraded quotient $\bK(n)_*(\bR\bX)/(v_n - 1)$ is a finite dimensional $\F_P$-algebra, thus its group of units is finite.]
\vspi
Claim: Fix a $v_n$-element $\alpha$ $-$then $\exists \ i$ such that $\alpha^i$ is in the center of $\pi_*(\bR\bX)$.
\vspi
[There is no loss of generality in supposing that $\bK(m)_*(\alpha)$ is in the center of $\bK(m)_*(\bR\bX)$ 
$\forall \ m \in [0,\infty[$.  Letting 
$\ad(\alpha):\Sigma^d \bR\bX \ra \bR\bX$ be the composite
\begin{tikzcd}[sep=small]
{\bS^d \wedge \bR\bX} \ar{rr}{\alpha \wedge \id}
&&{\bR\bX \wedge \bR\bX} \ar{rr}{\id - \Tee}
&&{\bR\bX \wedge \bR\bX}
\end{tikzcd}
%%----------------------------------------------------------------------------------------------36
$\ra \bR\bX$, $\ad(\alpha)_*(\beta) = \alpha\beta - \beta\alpha$ and $\forall \ i$, 
$\ad(\alpha^i)_*(\beta) =$ 
$\ds\sum\limits_j \binom{i}{j}\ad^j(\alpha)_*(\beta)\alpha^{i-j}$.  
Since $p^k\alpha = 0$ for some $k$ and 
$\ad(\alpha) \in [\bR\bX,\bR\bX]_*$ is nilpotent (cf. Proposition 29), one can take $i = p^N$ $(N \gg 0)$ to get that 
$\alpha^i\beta - \beta\alpha^i = 0$ $\forall \ \beta \in \pi_*(\bR\bX)$.]
\vspi
Claim: Fix $v_n$-elements $\alpha$,$\beta$ $-$then $\exists \ i, j$ such that $\alpha^i = \beta^j$.
\vspi
[Assuming, as is permissible, that $\alpha\beta = \beta\alpha$ and $\bK(m)_*(\alpha - \beta) = 0$ 
$\forall \ m \in [0,\infty[$, use the binomial theorem on 
$\alpha^{p^N} = (\beta + (\alpha - \beta))^{p^N}$ $(N \gg 0)$, observing that $\alpha - \beta$ is both torsion and nilpotent.]

The last claim serves to complete the proof of Proposition 31.\\
\endgroup %%------------------------------------<<

\begin{proposition} \ %32
Let \bX, \bY be $p$-compact  and fix $n \geq 1$.  Suppose that 
$\bff:\Sigma^d \bX \ra \bX$, 
$\bg:\Sigma^e \bY \ra \bY$ 
are $v_n$-maps $-$then $\exists \ i, j$ such that $\forall \ \bh \in [\bX,\bY]$ the diagram 
\begin{tikzcd}%[sep=small]
{\Sigma^{\id} \bX} \ar{d}[swap]{\bff^i} \ar{rr}{\Sigma^{\id} \bh = \Sigma^{je} \bh}
&&{\Sigma^{je} \bY} \ar{d}{\bg^j}\\
{\bX} \ar{rr}[swap]{\bh} &&{\bY}
\end{tikzcd}
commutes.\\
\end{proposition}

[Pass to $\hom(\bX,\bS_p ) \wedge \bY$ and apply Proposition 31.]\\

\label{17.78}
To round out the discussion on p. \pageref{17.79}, we shall now verify that $\bV_n$ is thick.  
Obviously, $\bV_n$ contains 0 and is stable under $\Sigma^{\pm 1}$.  Next, let \bX, \bY be objects of $\bV_n$ with $v_n$-maps 
$\bff:\Sigma^d \bX \ra \bX$, 
$\bg:\Sigma^e \bY \ra \bY$.  
Choose $i, j$ per Proposition 32 and put $k = id$ $(= je)$.  Take 
$\bX \overset{\bu}{\lra} \bY$ and complete it to an exact triangle 
$\bX \overset{\bu}{\lra}$ 
$\bY \overset{\bv}{\lra}$ 
$\bZ \overset{\bw}{\lra}$ 
$\Sigma \bX$ $-$then the claim is that $\bZ$ admits a $v_n$-map.  For consider the diagram 
\begin{tikzcd}%[sep=small]
{\Sigma^{\id} \bX}  \ar{d}[swap]{\bff^i} \ar{r}{\Sigma^k \bu} 
&{\Sigma^{je} \bY} \ar{d}{\bg^j} \ar{r}{\Sigma^k \bv}
&{\Sigma^k \bZ}\\
{\bX} \ar{r}[swap]{\bu}
&{\bY} \ar{r}[swap]{\bv}
&{\bZ} 
\end{tikzcd}
.  Since 
\begin{tikzcd}%[sep=small]
{\Sigma^{je} \bY} \ar{r}{\Sigma^k \bv} &{\Sigma^k \bZ}
\end{tikzcd}
is a weak cokernel of 
\begin{tikzcd}%[sep=small]
{\Sigma^{\id} \bX} \ar{r}{\Sigma^k \bu} &{\Sigma^{je} \bY}
\end{tikzcd}
and since 
$\bv \circ \bg^j \circ \Sigma^k \bu = \bv \circ \bu \circ \bff^i =$ $0$, $\exists$ an arrow 
$\bh:\Sigma^k\bZ \ra \bZ$ such that $\bh \circ \Sigma^k \bv = \bv \circ \bg^j$ 
(cf. p. \pageref{17.80} ff.).  
The five lemma gives that $\bK(n)_*(\bh)$ is an isomorphism.  
And: $\forall \ m \neq n$ $(m \in [0,\infty[)$, 
$\bK(m)_*(\bh^2) = 0$.  
Therefore $\bh^2$ is a $v_n$-map, so \bZ is in $\bV_n$, which means that $\bV_n$ is triangulated.  
Finally, if $\bY \in \Ob\bV_n$ and $\bY \approx \bX \vee \bZ$ with $\bi:\bX \ra \bY$, $\br:\bY \ra \bX$ and $\br \circ \bi = \id_{\bX}$, then 
$\bX \in \Ob\bV_n$.  
Thus fix a $v_n$-map $\bg:\Sigma^e\bY \ra \bY$.  
By raising \bg to a sufficiently high power, it can be arranged that the diagram
\begin{tikzcd}%[sep=large]
{\Sigma^e\bY} \ar{d}[swap]{\bg} \ar{rr}{\Sigma^e(\bi \circ \br)} &&{\Sigma^e\bY} \ar{d}{\bg}\\
{\bY} \ar{rr}[swap]{\bi \circ \br} &&{\bY}
\end{tikzcd}
commutes (cf. Proposition 32).  Applying $\bK(n)_*$ to 
\begin{tikzcd}%[sep=small]
{\Sigma^e\bX} \ar{d}{\bff} \ar{r}{\Sigma^e\bi} 
&{\Sigma^e\bY} \ar{d}{\bg} \ar{r}{\Sigma^e\br}
&{\Sigma^e\bX} \ar{d}{\bff} \\
{\bX} \ar{r}[swap]{\bi}
&{\bY} \ar{r}[swap]{\br}
&{\bX}
\end{tikzcd}
, where $\bff = \br \circ \bg \circ \Sigma^e \bi$, and using 
%%----------------------------------------------------------------------------------------------37
the fact that the retract of an isomorphism is an isomorphism, one concludes that \bff is a $v_n$-map.  
Accordingly, $\bV_n$ is thick.\\

\begin{proposition} \ %33
If \bE is $p$-local, then $\forall \ \bX$, 
$L_{\bE} \bX_p \approx L_{\bE} \bX \approx (L_{\bE} \bX)_p$.
\end{proposition}

[Since \bE is $p$-local, $\bE \approx \bE \wedge \bS(\bZ_p)$, hence 
$\langle \bE \rangle \leq \langle \bS(\bZ_p) \rangle$, and the lemma on p. \pageref{17.81} can be quoted.]

[Note: \ In order that \bX be $\bE_*$-local, it is therefore necessary that \bX be $p$-local.]\\

Application: If \bE is $p$-local and if 
$\L_{\bE}\bX \approx$ 
$\bX \wedge L_{\bE}\bS_p$ $\forall \ p$-local \bX, then $T_{\bE}$ is smashing.

[Given an arbitrary \bX, 
$\L_{\bE}\bX \approx$ 
$\L_{\bE}\bX_p \approx$ 
$\bX_p \wedge L_{\bE}\bS_p \approx$ 
$\bX \wedge \bS(\Z_p) \wedge L_{\bE} \bS_p \approx$ 
$\bX \wedge (L_{\bE} \bS_p )_p \approx$ 
$\bX \wedge L_{\bE}(\bS_p)_p \approx$ 
$\bX \wedge L_{\bE} \bS_p \approx$ 
$\bX \wedge L_{\bE} \bS$.]\\

Recall that for any \bE and any compact \bX, 
$L_{\bE}\bX \approx$ 
$\bX \wedge L_{\bE} \bS$ (cf. p. \pageref{17.82})
Corollary: For any $p$-local \bE and for any $p$-compact \bX, 
$L_{\bE}\bX \approx$
$\bX \wedge L_{\bE} \bS_p$.  
Proof: Write $\bX = \ov{\bX}_p $, where $\ov{\bX}$ is compact (cf. p. \pageref{17.83}) $-$then 
$L_{\bE}\bX \approx$ 
$L_{\bE} \ov{\bX}_p \approx$ 
$L_{\bE} \ov{\bX} \approx$ 
$\ov{\bX} \wedge L_{\bE}\bS \approx$ 
$\ov{\bX} \wedge L_{\bE}\bS_p \approx$ 
$\ov{\bX} \wedge \bS(\Z_p) \wedge L_{\bE} \bS_p \approx$ 
$\ov{\bX}_p \wedge  L_{\bE}\bS_p  \approx$ 
$\bX \wedge  L_{\bE} \bS_p$.  
Example: Taking $\bE = \bS(\Z/p\Z)$ $(= \bM(p))$, 
$L_{\bS(\Z/p\Z)} \bX \approx$ 
$\bX \wedge \widehat{\bS}_p$ if \bX is $p$-compact.\\

\begingroup%%----------------------------------->>
\fontsize{9pt}{11pt}\selectfont
\textbf{\small EXAMPLE}  \ 
Let $\bE \neq 0$ be $p$-local and suppose that there exists an $\bE_*$-local object in $\bC(n)$ for some $n < \infty)$.  
Case 1: $\bH(\Q) \wedge \bE \neq 0$ $-$then $L_{\bE}\bX \approx \bX$ $\forall \ p$-compact \bX.  
Case 2: $\bH(\Q) \wedge \bE = 0$ $-$then 
$L_{\bE}\bX \approx \bX \wedge \widehat{\bS}_p$ $\forall \ p$-compact \bX.
\vspi
[The class of all $p$-local \bX which are $\bE_*$-local must contain $\Ob\bC(1)$.  In addition, $\widehat{\bS}_p$ is 
$\bE_*$-local  (consider the exact triangle 
$\widehat{\bS}_p \overset{p}{\ra}$ 
$\widehat{\bS}_p \ra$ 
$\bM(p) \ra$ 
$\Sigma\widehat{\bS}_p$) and if \bF is defined by the exact triangle 
$\bF \ra$ 
$\bS_p \ra$ 
$\widehat{\bS}_p \ra \Sigma \bF$, then $\bF$ is $\bE_*$-local or $\bE_*$-acyclic depending on whether 
$\bH(\Q) \wedge \bE \neq 0$ or 
$\bH(\Q) \wedge \bE = 0$ (\bF is rational).  
Working now with the commutative diagram 
\begin{tikzcd}%[sep=small]
{\bF} \ar{d} \ar{r} 
&{\bS_p} \ar{d} \ar{r}
&{\widehat{\bS}_p} \ar{d} \ar{r}
&{\Sigma\bF} \ar{d}\\
{T_{\bE}\bF} \ar{r} 
&{T_{\bE}\bS_p} \ar{r}
&{T_{\bE}\widehat{\bS}_p} \ar{r}
&{\Sigma T_{\bE}\bF}
\end{tikzcd}
one thus sees that in case 1, $\bS_p$ is $\bE_*$-local ($\implies$ 
$L_{\bE}\bX \approx$
$\bX  \wedge L_{\bE}\bS_p \approx$ 
$\bX \wedge \bS_p \approx \bX$) while in case 2, 
$L_{\bE}\bS_p \approx$ 
$\widehat{\bS}_p$  ($\implies$ 
$L_{\bE}\bX \approx$ 
$\bX \wedge L_{\bE}\bS_p \approx$ 
$\bX \wedge \widehat{\bS}_p$).]\\
\endgroup %%------------------------------------<<

\begingroup%%----------------------------------->>
\fontsize{9pt}{11pt}\selectfont
\textbf{\small EXAMPLE}  \ 
Let $\bE \neq 0$ be a $p$-local ring spectrum with the property that $\bACY_p(\bE) = 0$.  
Case 1: $\bH(\Q) \wedge \bE \neq 0$ $-$then $L_{\bE}\bX \approx \bX$ $\forall \ p$-compact \bX.  
Case 2: $\bH(\Q) \wedge \bE = 0$ $-$then $L_{\bE}\bX \approx \bX \wedge \widehat{\bS}_p$ $\forall \ p$-compact \bX.
\vspi
[In view of the preceding example, one has only to exhibit an $\bE_*$-local object in $\bC(1)$.  
Choose $n \in [0,\infty]$: $\bK(n)_*(\bE) \neq 0$ (cf. p. \pageref{17.84}).  
If 
$\bK(\infty)_*(\bE) = \bH(\F_p)(\bE) \neq 0$, then 
$\langle \bH(\F_p) \rangle \leq \langle \bE \rangle$ and $\bM(p)$ is $\bH(\F_p)_*$-local, hence is $\bE_*$-local.  
So suppose that 
$\bH(\F_p) \wedge \bE = 0$.  
Claim: $\exists$ a sequence $k_1 < k_2 < \cdots$ such that 
$\bE \wedge \bK(k_i) \neq 0$ $(i = 1, 2, \ldots)$.  
Proof: $\forall \ n < \infty$, $\exists$ a $p$-compact ring spectrum $\bX_n$ of type n 
%%----------------------------------------------------------------------------------------------38
and $\bE \wedge \bX_n \neq 0$ 
(by hypothesis) $\implies$ $\bK(m)_*(\bE \wedge \bX_n) = 0$ 
($m < n$ or $m = \infty$) $\implies$ $\bK(m)_*(\bE \wedge \bX_n) \neq 0$ $(\exists \ m \in [n,\infty[)$.  
But
$\langle \bK \rangle \leq \langle \bE \rangle$ and $\bM(p)$ is $\bK_*$-local, where $\bK = \ds\bigvee\limits_i \bK(k_i)$.]\\
\endgroup %%------------------------------------<<

\begingroup%%----------------------------------->>
\fontsize{9pt}{11pt}\selectfont
\textbf{\small FACT}  \  
Let $\bE \neq 0$ be $p$-local.  Assume $\bACY_p(\bE) = 0$ and $T_{\bE}$ is smashing $-$then 
$\langle \bE \rangle = \langle \bS_p \rangle$.
\vspi
[Since $T_{\bE}$ is smashing, 
$\langle \bE \rangle = \langle L_{\bE}\bS \rangle = \langle L_{\bE}\bS_p \rangle$.  
However $L_{\bE}\bS_p \neq 0$ is a $p$-local ring spectrum with the property that $\bACY_p(L_{\bE}\bS_p) = 0$.  
Therefore 
$L_{L_{\bE}\bS_p}\bS_p \approx$ 
$L_{\bE}\bS_p \approx$ 
$\bS_p$ or $\widehat{\bS}_p$.  And: 
$\langle \bS_p \rangle = \langle \widehat{\bS}_p \rangle$ $\implies$ 
$\langle \bE \rangle = \langle \bS_p \rangle$.]\\
\endgroup %%------------------------------------<<

Let $\bX(n)$ be a $p$-compact spectrum of type $n$ $-$then by the class invariance principle, 
$\langle \bX(n) \rangle$ depends only on $n$.  
Write $\bT(n)$ for $\bff^{-1}\bX(n)$, where $\bff:\Sigma^d \bX(n) \ra \bX(n)$ is a $v_n$-map. 
Thanks to Proposition 31, $\bT(n)$ is independent of the choice of \bff.  
Moreover, its Bousfield class
$\langle \bT(n) \rangle$ is independent of the choice of $\bX(n)$ and applying Proposition 43 in $\S 15$ repeatedly, 
one obtains a decomposition 
$\langle \bS_p \rangle =$ 
$\langle \bT(0) \rangle \vee \langle \bT(1) \rangle$ 
$\vee \cdots \vee$ 
$\langle \bT(n) \rangle \vee \langle \bX(n+1) \rangle$ with 
$\langle \bT(i) \rangle \wedge \langle \bX(n+1) \rangle =$  
$\langle 0 \rangle$ 
$(i = 0, 1, \ldots, n)$,
$\langle \bT(n) \rangle \wedge \langle \bT(m)) \rangle =$  
$\langle 0 \rangle$ 
$(m \neq n)$ (here, $\bT(0) = \bH(\Q)$).  
Examples: 
(1) $\langle \bBP \rangle  \wedge \langle \bX(n) \rangle =$ $\langle \bP(n)\rangle$; 
(2) $\langle \bBP \rangle  \wedge \langle \bT(n) \rangle =$ $\langle \bK(n)\rangle$.

Notation: Put $\bT$ $(\leq n)$ $= \bT(0) \vee \bT(1) \vee \cdots \vee \bT(n)$, call $T_n^f$ the corresponding localization functor and let $L_n^f$ be the associated reflector.\\

\begin{proposition} \ %34
$T_n^f$ is smashing, so $\forall \ \bX$, $L_n^f \bX \approx \bX \wedge L_n^f \bS$.
\end{proposition}

[The Bousfield classes of $L_n^f \bS_p$ $(= L_n^f \bS)$ and $\bT (\leq n)$ are one and the same.]\\

\begingroup%%----------------------------------->>
\fontsize{9pt}{11pt}\selectfont
\textbf{\small FACT}  \  
Suppose that \bX is $p$-compact and has type $n$ $-$then 
$L_n^f \bX \approx \bff^{-1}\bX$, $\bff:\Sigma^d\bX \ra \bX$ is a $v_n$-map.\\
\endgroup %%------------------------------------<<

%dmc something troubleing in the below relating to the `latter'

Notation: Put $\bK (\leq n)$ $= \bK(0) \vee \bK(1) \vee \cdots \vee \bK(n)$, call $T_n$ the corresponding localization functor, and let $L_n$ be the associated reflector.

There are similarities between the ``$L_n^f$-theory'' and the ``$L_n$-theory'' (but the proofs for the latter are much more diificult).  Thus, e.g., it turns out that $T_n$ is smashing (cf. Proposition 34).  
Moreover, one can attach to any \bX a tower 
$L_0\bX \la L_1 \bX \la \cdots$ and 
$\bX \approx \mic(L_0\bX \la L_1 \bX \la \cdots)$ 
if \bX is $p$-compact (it is unknown whether the analog of this with $L_n$ replaced by $L_n^f$ is true or not).  On the other hand, $L_n^f$ and $L_n$ are connected by a natural transformation 
$L_n^f \ra L_n$ and $\forall \ \bX$,
$L_n^f \bX \ra L_n \bX$ is a $\bBP_*$-equivalence.

[Note: \ These assertions are detailed in 
Ravenel\footnote[2]{\textit{Nilpotence and Periodicity in Stable Homotopy Theory}, Princeton University Press (1992), 81-98.}
.  They represent the point of departure for the study of the ``chromatic'' aspects of $\bHSPEC$.]\\

\begingroup%%----------------------------------->>
\fontsize{9pt}{11pt}\selectfont
\textbf{\small FACT}  \  
Suppose that \bX is $p$-compact and has type $n$ $-$then 
$L_n^f \bX \approx L_{\bT(n)} \bX$ and 
$L_n\bX \approx L_{\bK(n)} \bX$.\\
\endgroup %%------------------------------------<<
%%%%%%%%%%%%%%%%%%%%%%%%%%%%%%%%%%%%%%
%%%%%%%%%%%%%%%%%%%%%%%%%%%%%%%%%%%%%%
%%%%%%%%%%%%%%%%%%%%%%%%%%%%%%%%%%%%%%

\begin{center}
$\S \ 17$
\\[0.5cm]
$\mathcal{REFERENCES}$\\
\end{center}

\[
\textbf{BOOKS}
\]

\begingroup
\fontsize{9pt}{11pt}\selectfont
\setlength\parindent{0 cm}

[1] \quad Adams, J., 
\textit{Stable Homotopy and Generalized Homology}, University of Chicago (1974).
\\[-.2cm]

[2] \quad Devinatz, E., 
\textit{A Nilpotence Theorem in Stable Homotopy}, Ph.D. Thesis, MIT, Cambridge (1985).
\\[-.2cm]

[3] \quad Hu, S-T., 
\textit{Homology Theory}, Holden-Day (1966).
\\[-.2cm]

[4] \quad Hu, S-T., 
\textit{Cohomology Theory}, Markham (1968).
\\[-.2cm]

[5] \quad Landweber, P. (ed.), 
\textit{Elliptic Curves and Modular Forms in Algebraic Topology}, Springer Verlag (1988).
\\[-.2cm]

[6] \quad Margolis, H., 
\textit{Spectra and the Steenrod Algebra}, North Holland (1983).
\\[-.2cm]

[7] \quad Ravenel, D., 
\textit{Complex Cobordism and Stable Homotopy Groups of Spheres}, Academic Press (1986).
\\[-.2cm]

[8] \quad Ravenel, D., 
\textit{Nilpotence and Periodicity in Stable Homotopy Theory}, Princeton University Press (1992).
\\[-.2cm]

[9] \quad Switzer, R., 
\textit{Algebraic Topology-Homotopy and Homology}, Springer Verlag (1975).
\\[-.2cm]

[10] \quad Vogt, R., 
\textit{Boardman's Stable Homotopy Category}, Aarhus Universitet (1970).
\\[-.2cm]
\endgroup

\[
\textbf{ARTICLES}
\]

\begingroup
\fontsize{9pt}{11pt}\selectfont
\setlength\parindent{0 cm}

[1] \quad Bauer, F., Classifying Spectra for Generalized Homology Theories, 
\textit{Ann. Mat. Pura Appl.} \textbf{164} (1993), 

\hspace{0.8cm}365-399.
\\[-.2cm]

[2] \quad Boardman, J., Stable Operations in Generalized Cohomology, In: 
\textit{Handbook of Algebraic Topology}, I. 

\hspace{0.8cm}James (ed.), North Holland (1995), 585-686.
\\[-.2cm]

[3] \quad Bott, R., The Periodicity Theorem for the Classical Groups and Some of its Applications, 
\textit{Adv. Math.} 

\hspace{0.8cm}\textbf{4} (1970), 353-411.
\\[-.2cm]

[4] \quad Bousfield, A., The Localization of Spectra with respect to Homology, 
\textit{Topology} \textbf{18} (1979), 257-281.
\\[-.2cm]

[5] \quad Bousfield, A., Uniqueness of Infinite Deloopings for K-Theoretic Spaces, 
\textit{Pacific J. Math.} \textbf{129} (1987), 

\hspace{0.8cm}1-31.
\\[-.2cm]

[6] \quad Connell, E., Characteristic Classes, 
\textit{Illinois J. Math.} \textbf{14} (1970), 496-521.
\\[-.2cm]

[7] \quad D\u ad\u arlat, M. and N\u emethi, A., Shape Theory and (Connective) K-Theory, 
\textit{J. Operator Theory} \textbf{23} 

\hspace{0.8cm}(1990), 207-291.
\\[-.2cm]

[8] \quad Devinatz, E., K-Theory and the Generating Hypothesis, 
\textit{Amer. J. Math.} \textbf{112} (1990), 787-804.
\\[-.2cm]

[9] \quad Devinatz, E., Morava's Change of Rings Theorem, 
\textit{Contemp. Math.} \textbf{181} (1995), 83-118.
\\[-.2cm]

[10] \quad Devinatz, E., The Generating Hypothesis Revisited, 
\textit{Fields Institute Communications} \textbf{19} (1998), 

\hspace{0.95cm}73-92.
\\[-.2cm]

[11] \quad Devinatz, E., Hopkins, M., and Smith, J., Nilpotence and Stable Homotopy Theory I, 
\textit{Ann. of Math.} 

\hspace{0.95cm}\textbf{128} (1988),207-241.
\\[-.2cm]

[12] \quad Dold, A., Chern Classes in General Cohomology, 
\textit{Symposia Mathematica} V (1969-70), 385-410.
\\[-.2cm]

[13] \quad Giffen, C., Bott Periodicity and the $Q$-Construction, 
\textit{Contemp. Math.} \textbf{199} (1996), 107-124.
\\[-.2cm]

[14] \quad Hopkins, M., Global Methods in Homotopy Theory, In: 
\textit{Homotopy Theory}, E. Rees and J. Jones 

\hspace{0.95cm}(ed.), Cambridge University Press (1987), 73-96.
\\[-.2cm]

[15] \quad Hopkins, M. and Ravenel, D., Suspension Spectra are Harmonic, 
\textit{Bol. Soc. Mat. Mexicana} \textbf{37} 

\hspace{0.95cm}(1992), 271-279.
\\[-.2cm]

%[16] \quad Hopkins, M. and Smith, J., Nilpotence and Stable Homotopy Theory II, Ann. of Math.\\[-.2cm]
[16] \quad Hopkins, M. and Smith, J., Nilpotence and Stable Homotopy Theory II, 
\textit{Ann. of Math.} \textbf{148}, 1 

\hspace{0.95cm}(1998), 1-49.\\[-.2cm]

%[17] \quad Hovey, M., Sadofsky, H., and Strickland, N., Morava K-Theories and Localization,\\[-.2cm]
[17] \quad Hovey, M.,  and Strickland, N., Morava K-Theories and Localization, 
\textit{Memoirs Amer. Math. Soc.} 

\hspace{0.95cm}(1999).
\\[-.2cm]

[18] \quad Kahn,  D., Kaminker, J., and Schochet, C., Generalized Homology Theories on Compact Metric 

\hspace{0.95cm}Spaces, 
\textit{Michigan Math. J.} \textbf{24} (1977), 203-224.
\\[-.2cm]

[19] \quad Kuhn, N., Morava K-Theories and Infinite Loop Spaces, 
\textit{SLN} \textbf{1370} (1989), 243-257.
\\[-.2cm]

[20] \quad Landweber, P., A Survey of Bordism and Cobordism, 
\textit{Math. Proc. Cambridge Philos. Soc.} \textbf{100} 

\hspace{0.95cm}(1986), 207-223.
\\[-.2cm]

[21] \quad Landweber, P., Ravenel, D., and Stong, R., Periodic Cohomology Theories Defined by Elliptic Curves,

\hspace{0.95cm}\textit{Contemp. Math.} \textbf{181} (1995), 317-337.
\\[-.2cm]

[22] \quad Mahowald, M. and Ravenel, D., Toward a Global Understanding of the Homotopy Groups of Spheres, 

\hspace{0.95cm}\textit{Contemp. Math.} \textbf{58} (1987), 57-74.
\\[-.2cm]

[23] \quad Milnor, J., On Axiomatic Homology Theory, 
\textit{Pacific J. Math.} \textbf{12} (1962), 337-341.
\\[-.2cm]

[24] \quad Ravenel, D., Localization with respect to Certain Periodic Homology Theories, 
\textit{Amer. J. Math.} \textbf{106} 

\hspace{0.95cm}(1984), 351-414.
\\[-.2cm]

[25] \quad Ravenel, D., Life After the Telescope Conjecture, In: 
\textit{Algebraic K-Theory and Algebraic Topology}, P. 

\hspace{0.95cm}Goerss and J.\ Jardine (ed.), Kluwer (1993), 205-222.
\\[-.2cm]

%[26] \quad Ravenel, D., Wilson, W., and Yagita, N., Brown-Peterson Cohomology from Morava K-Theory,\\[-.2cm]
[26] \quad Ravenel, D., Wilson, W., and Yagita, N., Brown-Peterson Cohomology from Morava K-Theory, 

\hspace{0.95cm}\textit{K-Theory} \textbf{15(2)} (1998), 147-199. 
\\[-.2cm]
%\textit{K-Theory} \textbf{324} (1998), 1-52 \\[-.2cm]

[27] \quad Sklyarenko, E., Homology and Cohomology Theories of General Spaces, In: 
\textit{General Topology}, EMS 

\hspace{0.95cm}\textbf{50}, Springer Verlag (1996), 119-246.
\\[-.2cm]

%[28] \quad Strickland, N., Functorial Philosophy for Formal Phenomena,\\[-.2cm]
[28] \quad Strickland, N., Functorial Philosophy for Formal Phenomena, 
\textit{Preprint} (1994). 
\\[-.2cm]

[29] \quad Taylor, J., Banach Algebras and Topology, In: 
\textit{Algebras in Analysis}, J. Williamson (ed.), Academic 

\hspace{0.95cm}Press (1975) 118-186.
\\[-.2cm]

[30] \quad Totaro, B., Torsion Algebraic Cycles and Complex Cobordism, 
\textit{J. Amer. Math. Soc.} \textbf{10} (1997), 

\hspace{0.95cm}467-493.
\\[-.2cm]

[31] \quad Whitehead, G., Generalized Homology Theories, 
\textit{Trans. Amer. Math. Soc.} \textbf{102} (1962), 227-283.

\setlength\parindent{2em}

\endgroup

\chapter{
$\boldsymbol{\S}$\textbf{18}.\quadx  ALGEBRAIC K$-$THEORY}
\setlength\parindent{2em}
\setcounter{proposition}{0}

%%----------------------------------------------------------------------------------------------01
$\text{ }$\\[-1.25cm]

My objective in this $\S$ is to provide an introduction to algebraic K-theory, placing the emphasis on its homotopical underpinnings.

Consider a skeletally small category \bC equipped with two composition closed classes of morphisms termed 
\un{weak equivalences}
\index{weak equivalences} 
(denoted $\overset{\sim}{\ra}$) and 
\un{cofibrations}
\index{cofibrations} 
(denoted $\rightarrowtail$), 
each containing the isomorphisms of \bC $-$then \bC is said to be a 
\un{Waldhausen category}
\index{Waldhausen category} provided that the following axioms are satisfied.

\indent\indent (WC-1) \ \bC has a zero object 0.

\indent\indent (WC-2) \ All the objects of \bC are cofibrant, i.e., $\forall \ X \in \Ob\bC$, the arrow $0 \ra X$ is a cofibration.

\indent\indent (WC-3) \ Every 2-source $X \overset{f}{\la} Z \overset{g}{\ra} Y$, where $f$ is a cofibration, admits a pushout 
$X \overset{\xi}{\ra} P \overset{\eta}{\la} Y$, where $\eta$ is a cofibration.

\indent\indent (WC-4) \ If 
\begin{tikzcd}%[sep=small]
{X} \ar{d} &{Z} \ar{l}[swap]{f} \ar{d} \ar{r}{g} &{Y}  \ar{d}\\
{X^\prime} &{Z^\prime} \ar{l}{f^\prime} &{Y^\prime} \ar{l}{g^\prime}
\end{tikzcd}
%dmc is g^prime going the wrong way?
is a commutative diagram, where 
$
\begin{cases}
\ f\\
\ f^\prime
\end{cases}
$
are cofibrations and the vertical arrows are weak equivalences, then the induced morphism $P \ra P^\prime$ of pushouts is a weak equivalence.

[Note: \ The opposite of a Waldhausen category need not be Waldhausen.]

Remark: \bC has finite coproducts (define $X \ \amalg \ Y$ by the pushout square
\begin{tikzcd}%[sep=small]
{0} \ar{d} \ar{r} &{Y} \ar{d}\\
{X} \ar{r} &{X \amalg Y}
\end{tikzcd}
($\implies$ $\ini_X$ $\&$ $\ini_Y$ are cofibrations)).

[Note: Every cofibration $X \rightarrowtail Y$ has a cokernel $Y/X$, viz. $Y \underset{X}{\sqcup} 0$.]

Example: A finitely cocomplete pointed skeletally small category is Waldhausen if the weak equivalences are the isomorphisms and the cofibrations are the morphisms.\\

\begingroup%%----------------------------------->>
\fontsize{9pt}{11pt}\selectfont
\textbf{\small EXAMPLE} \  
Take for \bC the category whose objects are the pointed finite sets $-$then \bC is a Waldhausen category if weak equivalence = isomorphism, cofibration = pointed injection.\\
\endgroup %%------------------------------------<<

\begingroup%%----------------------------------->>
\fontsize{9pt}{11pt}\selectfont
\textbf{\small EXAMPLE} \  
Take for \bC the category whose objects are the pointed finite simplicial sets $-$then \bC is a Waldhausen category if weak equivalence = weak homotopy equivalence, cofibration = pointed injective simplicial map.\\
\endgroup %%------------------------------------<<

\begingroup%%----------------------------------->>
\fontsize{9pt}{11pt}\selectfont
\textbf{\small EXAMPLE} \  
Let \mA be a ring with unit.  Denote by $\bP(A)$ the full subcategory of \bAMOD whose objects are finitely generated and projective $-$then $\bP(A)$ is a Waldhausen category if weak equivalence = isomorphism, cofibration = split injection.\\
\endgroup %%------------------------------------<<

%%----------------------------------------------------------------------------------------------02
\begingroup%%----------------------------------->>
\fontsize{9pt}{11pt}\selectfont
\textbf{\small EXAMPLE} \  
Let \mA be a ring with unit.  Denote by $\bF(A)$ the full subcategory of \bAMOD whose objects are finitely generated and free $-$then $\bF(A)$ is a Waldhausen category if weak equivalence = isomorphism, cofibration = split injection with free quotient.\\
\endgroup %%------------------------------------<<

\begingroup%%----------------------------------->>
\fontsize{9pt}{11pt}\selectfont
\textbf{\small FACT} \  
The cofibrant objects in a pointed skeletally small cofibration category are the objects of a Waldhausen category 
(cf. $\S 12$, Proposition 3 and p. \pageref{18.1}).\\
\endgroup %%------------------------------------<<

\begin{proposition} \ %01
Any skeleton of a Waldhausen category is a small Waldhausen category.\\
\end{proposition}

There are two other conditions which are sometimes imposed on a Waldhausen category.

\indent\indent (Saturation Axiom) \ Given composable morphisms $f$, $g$ if any two of 
$f, g,  g \circx f$ are weak equivalences, then so is the third.\\
\indent\indent (Extension Axiom) \ Given a \cd 
\begin{tikzcd}%[sep=small]
{X} \ar{d} \ar[r, tail] &{Y} \ar{d} \ar{r} &{Y/X} \ar{d}\\
{X^\prime} \ar[r, tail]  &{Y^\prime} \ar{r} &{Y^\prime/X^\prime}
\end{tikzcd}
, if 
$X \ra X^\prime$ $\&$ $Y/X \ra Y^\prime/X^\prime$ are weak equivalences, then $Y \ra Y^\prime$ is a weak equivalence.\\

\begingroup%%----------------------------------->>
\fontsize{9pt}{11pt}\selectfont
Neither the saturation axiom nor the extension axiom is a consequence of the other axioms.\\
\endgroup %%------------------------------------<<

Observation: If \bC is a Waldhausen category, then its arrow category $\bC(\ra)$ is a Waldhausen category.

[The weak equivalences and cofibrations are levelwise.]

Let \bC be a Waldhausen category $-$then a 
\un{mapping cylinder}
\index{mapping cylinder (Waldhausen category)} 
is a functor $M:\bC(\ra)$ $\ra$ $\bC$ together with natural transformations $i:S \ra M$, $j:T \ra M$, $r:M \ra T$, where 
$S:\bC(\ra) \ra \bC$ is the source functor and 
$T:\bC(\ra) \ra \bC$ is the target functor, all subject to the following assumptions.

[Note: \ Spelled out, \mM assigns to each object $X \overset{f}{\ra} Y$ in $\bC(\ra)$ an object $M_f \in \bC$ and to 
each morphism $(\phi,\psi):f \ra f^\prime$ in $\bC(\ra)$ a morphism 
$M_{\phi,\psi}:M_f \ra M_{f^\prime}$ in \bC.]\\

\indent\indent (MCy$_1$) \ For every object $X \overset{f}{\ra} Y$ in $\bC(\ra)$, the diagrams
\begin{tikzcd}%[sep=small]
{X} \ar{rd}[swap]{f} \ar{r}{i} &{M_f} \ar{d}{r}\\
&{Y}
\end{tikzcd}
,
\begin{tikzcd}%[sep=small]
{M_f} \ar{d}[swap]{r}  &{Y} \ar{l}[swap]{j} \ar[equals]{dl} \\
{Y}
\end{tikzcd}
commute and $i \amalg j:X \amalg Y \ra M_f$ is a cofibration (hence $i$ $\&$$ j$ are cofibrations).

\indent\indent (MCy$_2$) \ For every object \mY in \bC, $M_{0 \ra Y} = Y$ with $r = \id_Y$ and $j = \id_Y$.

\indent\indent (MCy$_3$) \ For every morphism $(\phi,\psi):f \ra f^\prime$ in $\bC(\ra)$, 
$M_{\phi,\psi}:M_f \ra M_{f^\prime}$ is a weak equivalence of morphisms if $\phi, \psi$ are weak equivalences.

%%----------------------------------------------------------------------------------------------03
\indent\indent (MCy$_4$) \ For every morphism $(\phi,\psi):f \ra f^\prime$ in $\bC(\ra)$, 
$M_{\phi,\psi}:M_f \ra M_{f^\prime}$ is a cofibration if $\phi, \psi$ are cofibrations.

\indent\indent (MCy$_5$) \ For every morphism $(\phi,\psi):f \ra f^\prime$ in $\bC(\ra)$, the diagram
\begin{tikzcd}%[sep=small]
{X \amalg Y} \ar{d} \ar{r}{i \amalg j} &{}\\
{X^\prime \amalg Y^\prime} \ar{r}[swap]{i \amalg j}  &{}
\end{tikzcd}
\begin{tikzcd}%[sep=small]
{M_f} \ar{d} \ar{r}{r} &{Y} \ar{d}\\
{M_{f^\prime}} \ar{r}[swap]{r} &{Y^\prime}
\end{tikzcd}
commutes and if $\phi, \psi$ are cofibrations, then the arrow 
$(X^\prime \amalg Y^\prime) \underset{X \amalg Y}{\sqcup} M_f \ra M_{f^\prime}$ is a cofibration.

Example: The 
\un{cone functor}
\index{cone functor (Waldhausen category)} 
$\Gamma:\bC \ra \bC$ sends \mX to $\Gamma X$, where 
$\Gamma X = M_{X \ra 0}$ and the 
\un{suspension functor}
\index{suspension functor (Waldhausen category)} 
$\Sigma:\bC \ra \bC$ sends \mX to $\Sigma X = \Gamma X /X$ (per $X \overset{i}{\ra} \Gamma X$).\\

\begingroup%%----------------------------------->>
\fontsize{9pt}{11pt}\selectfont
\textbf{\small EXAMPLE} \  
The category of pointed finite simplicial sets, where  weak equivalence = weak homotopy equivalence and cofibration = pointed injective simplicial map, has a mapping cylinder.\\
\endgroup %%------------------------------------<<

\index{Mapping Cylinder Axiom (Waldhausen category)}
\indent\indent (Mapping Cylinder Axiom) 
Assume that \bC admits a mapping cylinder $-$then $\forall \ X \overset{f}{\ra} Y \in \Ob\bC(\ra)$, $r:M_f \ra Y$ is a weak equivalence.\\

\begingroup%%----------------------------------->>
\fontsize{9pt}{11pt}\selectfont
\textbf{\small EXAMPLE} \  
The category of pointed finite simplicial sets, where  weak equivalence = isomorphism and cofibration = pointed injective simplicial map, has a mapping cylinder which does not satisfy the mapping cylinder axiom.\\
\endgroup %%------------------------------------<<

In a Waldhausen category, an 
\un{acyclic cofibration}
\index{acyclic cofibration (Waldhausen category)} 
is a morphism 
which is both a weak equivalence and a cofibration.\\

\begin{proposition} \ 
If $X \overset{f}{\la} Z \overset{g}{\ra} Y$ is a 2-source, where $f$ is an acyclic cofibration, then $Y \overset{\eta}{\ra} P$ is an acyclic cofibration.
\end{proposition}

[Bearing in mind WC-3, consider the commutative diagram
\begin{tikzcd}%[sep=small]
{Z} \ar{d}[swap]{f}  &{Z} \arrow[d,shift right=0.5,dash] \arrow[d,shift right=-0.45,dash] 
\ar{l}[swap]{\id_Z}\ar{r}{g} 
&{Y} \arrow[d,shift right=0.5,dash] \arrow[d,shift right=-0.45,dash]\\
{X}  &{Z} \ar{l}{f} \ar{r}[swap]{g} &{Y}
\end{tikzcd}
and apply WC-4.]

[Note: Therefore $0 \ra Y/X$ is an acyclic cofibration if $X \ra Y$ is an acyclic cofibration.]\\

Remark: If \bC satisfies the saturation axiom and the mapping cylinder axiom, then $j$ is an acyclic cofibration and $i$ is an acyclic cofibration provided that $f$ is a weak equivalence.\\

%%----------------------------------------------------------------------------------------------04
Notation: Given a Waldhausen category \bC, $\bw\bC$ is the subcategory of \bC having morphisms the weak equivalences, $\bco\bC$ is the subcategory of \bC having morphisms the cofibrations, and $\bwco\bC$  is the subcategory of \bC having morphisms the acyclic cofibrations.\\

\begin{proposition} \ 
Suppose that \bC is a small Waldhausen category satisfying the saturation axiom and the mapping cylinder axiom $-$then the inclusion $\iota:\textnormal{$\bwco\bC$} \ra \textnormal{$\bw\bC$}$ induces a pointed homotopy equivalence $B\iota:B\bwco\bC \ra B\bw\bC$.
\end{proposition}

[Owing to Quillen's theorem A, it suffices to show that $\iota$ is a strictly initial functor, i.e., that 
$\forall \ Y \in \Ob\bw\bC$, the comma category $\iota/Y$ is contractible.  An object of $\iota/Y$ is a pair $(X,f)$ where 
$f:X \ra Y$ is a weak equivalence.  
Specify a functor $m:\iota/Y \ra \iota/Y$ by sending $(X,f)$ to $(M_f,r)$ $-$then $i$ defines 
a natural transformation $\id_{\iota/Y} \ra m$ and $j$ defines a natural transformation $K_{(Y,\id_Y)} \ra m$.  
Therefore $B\iota/Y$ is contractible (cf. p. \pageref{18.2}).]

[Note: The base point is the 0-cell corresponding to 0.]\\

Let \bC be an additive category $-$then a pair of composable morphisms 
$X \overset{i}{\ra} Y \overset{p}{\ra} Z$ is 
\un{exact}
\index{exact (composable morphisms  in an additive category}
if $i$ is a kernel of $p$ and $p$ is a cokernel of $i$, a
\un{morphism of exact pairs}
\index{morphism of exact pairs (additive category} 
being a triple $(f, g, h)$ such that the diagram 
\begin{tikzcd}%[sep=small]
{X} \ar{d}{f} \ar{r}{i}  &{Y}\ar{d}{g} \ar{r}{p}  &{Z} \ar{d}{h}\\
{X^\prime} \ar{r}[swap]{i^\prime}  &{Y^\prime}  \ar{r}[swap]{p^\prime} &{Z^\prime}
\end{tikzcd}
commutes.

[Note: \ The first component of an exact pair is called an 
\un{inflation}
\index{inflation (first component of an exact pair)} 
(denoted $\rightarrowtail$),
the second component a  
\un{deflation}
\index{deflation (second component of an exact pair)} 
(denoted $\twoheadrightarrow$)
(terminology as in 
Gabriel-Roiter\footnote[2]{\textit{Representations of Finite Dimensional Algebras}, Springer Verlag (1992).}).

Let \bC be a skeletally small additive category $-$then \bC is said to be a 
\un{category with} 
\un{exact sequences}
\index{category with exact sequences} 
(category WES 
\index{WES}
) if there is given an 
isomorphism closed class $\sE$ of exact pairs satisfying the following conditions.

\indent\indent (ES-1) \ The pair $0 \overset{\id_0}{\lra} 0 \overset{\id_0}{\lra} 0$ is in $\sE$.

\indent\indent (ES-2) \ The composition of two inflations is an inflation and the composition of two deflations is a deflation.

\indent\indent (ES-3) \ Every 2-source $X \overset{f}{\la} Z \overset{g}{\ra} Y$, where $f$ is an inflation, admits a pushout 
$X \overset{\xi}{\ra} P \overset{\eta}{\la} Y$, where $\eta$ is an inflation, and every two sink 
$X \overset{f}{\ra} Z \overset{g}{\la} Y$, where $g$ is a deflation, admits a pullback 
$X \overset{\xi}{\la} P \overset{\eta}{\ra} Y$, 
where $\xi$ is a deflation.

[Note: \ The opposite of a category WES is again a category WES.]

A full, additive subcategory \bC of an abelian category \bD is 
\un{closed under extensions}
\index{closed under extensions}
if for every short exact sequence 
$0 \ra X \ra Y \ra Z \ra 0$ in \bD, where $X, Z \in \Ob\bC$, $\exists$ an
%%----------------------------------------------------------------------------------------------05
object in \bC which is isomorphic to \mY.

[Note: \ Such a \bC necessarily has finite coproducts.]

Example: Let \bC be a full, skeletally small additive subcategory of an abelian category \bD.  
Assume: \bC is closed under 
extensions.  
Declare a sequence 
$X \overset{i}{\ra} Y \overset{p}{\ra} Z$ in \bC to be exact iff 
$0 \ra X \overset{i}{\ra} Y \overset{p}{\ra} Z \ra 0$ is short exact in \bD $-$then \bC is a category WES.

[Note: \ This example is prototypical.  Thus suppose that \bC is a category WES $-$then $\exists$ an abelian category 
$\bG\text{-}\bQ$ and an additive functor $\iota:\bC \ra \bG\text{-}\bQ$ which is full and faithful such  that 
$X \overset{i}{\ra} Y \overset{p}{\ra} Z$ is exact iff 
$0 \ra \iota X \overset{\iota i}{\ra} \iota Y \overset{\iota p}{\ra} \iota Z \ra 0$ is short exact.  
And: $\iota \bC$ is closed under extensions.  
Specifically: $\bG\text{-}\bQ$ is the full subcategory of $[\bC^\OP,\bAB]^+$ whose objects are those \mF such that 
$X \overset{i}{\ra} Y \overset{p}{\ra} Z$ exact $\implies$  
$0 \ra FZ \ra FY \ra FX$ exact and $\iota$ is the Yoneda embedding.  
For a proof, consult  
Thomason-Trobaugh\footnote[2]{\textit{The Grothendieck Festschrift}, vol. III Birkh\"auser (1990), 247-435 (cf. 399-406); 
see also Keller, \textit{Manuscripta Math.} \textbf{67} (1990), 379-417 (cf. 408-409).}($\bG\text{-}\bQ$ = Gabriel-Quillen).]\\

\textbf{\small LEMMA} \  
Let \bC be a category WES $-$then $\forall \ X \in \Ob\bC$, $\id_X$ is both an inflation and a deflation.

[Consider the pushout square 
\begin{tikzcd}%[sep=small]
{0} \ar{d}\ar{r}
&{X} \ar{d}{\id_X}\\
{0}\ar{r}
&{X}
\end{tikzcd}
to see that $\id_X$ is an inflation and consider the pullback square 
\begin{tikzcd}%[sep=small]
{X} \ar{d}[swap]{\id_X}\ar{r}
&{0} \ar{d}\\
{X}\ar{r}
&{0}
\end{tikzcd}
to see that $\id_X$ is a deflation.]

[Note: \ Similarly, $0 \ra X$ is an inflation and $X \ra 0$ is a deflation.  Therefore 
$0 \ra X \overset{\id_X}{\ra} X$ and $X \overset{\id_X}{\ra} X \ra 0$ are exact.]\\

Application: \ Every morphism $\phi:X \ra Y$ is both an inflation and a deflation.

[By assumption, $\sE$ is isomorphism closed and there are commutative diagrams \\
\quad
\begin{tikzcd}%[sep=small]
{X} \arrow[d,shift right=0.5,dash] \arrow[d,shift right=-0.5,dash] \ar{r}{\phi}
&{Y} \ar{d}{\phi^{-1}}\ar{r}
&{0} \arrow[d,shift right=0.5,dash] \arrow[d,shift right=-0.5,dash]\\
{X}\ar{r}[swap]{\id_X}
&{X}\ar{r}
&{0}
\end{tikzcd}
, \hspace{0.5cm}
\begin{tikzcd}%[sep=small]
{0} \arrow[d,shift right=0.5,dash] \arrow[d,shift right=-0.5,dash]\ar{r}
&{X} \arrow[d,shift right=0.5,dash] \arrow[d,shift right=-0.5,dash] \ar{r}{\phi}
&{Y} \ar{d}{\phi^{-1}}\\
{0}\ar{r}
&{X}\ar{r}[swap]{\id_X}
&{X}
\end{tikzcd}
.]\\[.5cm]

\begin{proposition} \ %04
A category WES is a Waldhausen category.
\end{proposition}

[Take for the weak equivalences the isomorphisms and take for the cofibrations the inflations.]

%%----------------------------------------------------------------------------------------------06
[Note: \ This interpretation entails a loss of structure.]\\

Remark: Any skeleton of a category WES is a small category WES (cf. Proposition 1).\\

\begingroup%%----------------------------------->>
\fontsize{9pt}{11pt}\selectfont
Let \bC be a category WES.\\
\endgroup %%------------------------------------<<

\begingroup%%----------------------------------->>
\fontsize{9pt}{11pt}\selectfont
\textbf{\small FACT} \  
Consider a pushout square 
\begin{tikzcd}[sep=large]
{Z} \ar{d}[swap]{f} \ar{r}{g}
&{Y} \ar{d}{\eta}\\
{X} \ar{r}[swap]{\xi} &{P}
\end{tikzcd}
, where $f$ is an inflation $-$then 
\begin{tikzcd}[sep=small]
{Z}  \ar{rr}{
\begin{pmatrix}
-f\\
g\\
\end{pmatrix}
}
&&{X \oplus Y} %\ar{rr}{(\xi,\eta)}
\end{tikzcd}
\begin{tikzcd}[sep=small]
{} 
\ar{rr}{(\xi,\eta)}
&&{P}
\end{tikzcd}
is exact.\\
\endgroup %%------------------------------------<<

\begingroup%%----------------------------------->>
\fontsize{9pt}{11pt}\selectfont
\textbf{\small FACT} \  
Consider a pullback square 
\begin{tikzcd}[sep=large]
{P} \ar{d}[swap]{\xi} \ar{r}{\eta}
&{Y} \ar{d}{g}\\
{X} \ar{r}[swap]{f} &{Z}
\end{tikzcd}
, where $g$ is a deflation $-$then 
\begin{tikzcd}[sep=small]
{P}  \ar{rr}{
\begin{pmatrix}
-\eta\\
\xi\\
\end{pmatrix}
}
&&{X \oplus Y} %\ar{rr}{(\xi,\eta)}
\end{tikzcd}
\begin{tikzcd}[sep=small]
{}
\ar{rr}{(f,g)}
&&{Z}
\end{tikzcd}
is exact.\\
\endgroup %%------------------------------------<<

\begingroup%%----------------------------------->>
\fontsize{9pt}{11pt}\selectfont
\textbf{\small FACT} \  
If $f:X \ra Y$ has a cokernel and if $g \circx f$ is an inflation for some morphism $g$, then $f$ is an inflation.\\
\endgroup %%------------------------------------<<

\begingroup%%----------------------------------->>
\fontsize{9pt}{11pt}\selectfont
\textbf{\small FACT} \  
If $f:X \ra Y$ has a kernel and if $f \circx g$ is a deflation for some morphism $g$, then $f$ is a deflation.\\
\endgroup %%------------------------------------<<

\begingroup%%----------------------------------->>
\fontsize{9pt}{11pt}\selectfont
\textbf{\small FACT} \  
$\forall \ X, Y \in \Ob\bC$, 
\begin{tikzcd}%[sep=small]
X \ar{r}{\ini_X} &{X \oplus Y} \ar{r}{\pr_Y} &{Y}
\end{tikzcd}
is exact.\\
\endgroup %%------------------------------------<<

\begingroup%%----------------------------------->>
\fontsize{9pt}{11pt}\selectfont
\textbf{\small EXAMPLE} \  
Let \mA be a ring with unit $-$then $\bP(A)$ and $\bF(A)$ are categories WES.\\
\endgroup %%------------------------------------<<

\begingroup%%----------------------------------->>
\fontsize{9pt}{11pt}\selectfont
\textbf{\small EXAMPLE} \  Let \mX be a scheme, $\sO_X$ its structure sheaf $-$then the category of locally free 
$\sO_X$-modules of finte rank is a category WES.\\
\endgroup %%------------------------------------<<

\begingroup%%----------------------------------->>
\fontsize{9pt}{11pt}\selectfont
\textbf{\small EXAMPLE} \  Let \mX be a topological space $-$then the category of real or complex vector bundles over \mX is a category WES.\\
\endgroup %%------------------------------------<<

Let \bC be a category WES $-$then a pair $(\bA,\iota)$, where \bA is an abelian category and $\iota:\bC \ra \bA$ is an additive functor which is full and faithful, satisfies the 
\un{embedding condition}
\index{embedding condition} provided that 
$X \overset{i}{\ra} Y \overset{p}{\ra} Z$ is exact iff 
$0 \ra \iota X \overset{\iota i}{\lra} \iota Y \overset{\iota p}{\lra} \iota Z \ra 0$ is short exact.  
And: $\iota \bC$ is closed under extensions.  
Example: The pair $(\bG$-$\bQ,\iota)$ satisfies the embedding condition.

%%----------------------------------------------------------------------------------------------07
\indent\indent ($E \Rightarrow D$ Axiom) \ 
Under the assumption that the pair $(\bA,\iota)$ satisfies the embedding condition, an 
$f \in \Mor\bC$ is a deflation whenever $\iota f \in \Mor \bA$ is an epimorphism.\\

\begingroup%%----------------------------------->>
\fontsize{9pt}{11pt}\selectfont
\textbf{\small EXAMPLE} \  Let \mX be a scheme, $\sO_X$ its structure sheaf.  With \bC the category of locally free 
$\sO_X$-modules of finite rank, let \bA be either the abelian category of $\sO_X$-modules or the abelian category of quasicoherent $\sO_X$-modules $-$then in either case, the pair $(\bA,\iota)$ satisfies the embedding condition and the 
$E \Rightarrow D$ axiom.\\
\endgroup %%------------------------------------<<

A 
\un{pseudoabelian category}
\index{pseudoabelian category} 
is an additive category \bC with finite coproducts such that 
every idempotent has a kernel.  
Example: Let \mA be a ring with unit $-$then $\bP(A)$ is pseudoabelian (but this need not be the case of $\bF(A)$).

[Note: \ If \bC is pseudoabelian and if $e:X \ra X$ is an idempotent, then 
$X \approx \ker e \oplus \ker(1 - e)$ and 
$e \leftrightarrow 0 \oplus 1$.]\\

\textbf{\small LEMMA} \  
Let \bC be a category WES.  Assume: \bC is pseudoabelian $-$then $f \in \Mor\bC$ is a deflation if $f$ has a right inverse.\\

\label{18.9}
Remark: Let \bC be a category WES $-$then, while the pair $(\bG$-$\bQ,\iota)$ satisfies the embedding condition, it is not automatic that the $E \Rightarrow D$ axiom holds.  To ensure this, it suffices that retracts be deflations 
(Thomason-Trobaugh (ibid.)) which, by the lemma, will be true if \bC is pseudoabelian.\\

\begingroup%%----------------------------------->>
\fontsize{9pt}{11pt}\selectfont
\textbf{\small EXAMPLE} \  
Let \mX be a topological space $-$then the category of real or complex vector bundles over \mX is pseudoabelian.\\
\endgroup %%------------------------------------<<

Rappel: \ Let \bC be an additive category with finite coproducts $-$then there exists a pseudoabelian category 
$\bC_{\pa}$ and an additive functor $\Phi:\bC \ra \bC_{\pa}$ which is full and faithful such that for any pseudoabelian category  
\bD and any additive functor $F:\bC \ra \bD$, there exists an additive functor $F_{\pa}:\bC_{\pa} \ra \bD$ such that 
$F \approx F_{\pa} \circx \Phi$.  
And: $\bC_{\pa}$ is unique up to equivalence.

[One model for $\bC_{\pa}$ is the category whose objects are the pairs $(X,e)$, where $X \in \Ob\bC$ and 
$e \in \Mor(X,X)$ is idempotent, and whose morphisms $(X,e) \ra (X^\prime,e^\prime)$ are the $f \in \Mor(X,X^\prime)$ 
such that $f = e ^\prime \circx f \circx e$.  Here $\id_{(X,e)} = e$ and 
$(X,e) \oplus (X^\prime,e^\prime) =$ 
$(X \oplus X^\prime,e \oplus e^\prime)$.  As for $\Phi:\bC \ra \bC_{\pa}$, it is defined by 
$\Phi X = (X,\id_X)$ $\&$ $\Phi f = f$.

[Note: \ Every object in $\bC_{\pa}$ is a direct summand of an object in $\Phi \bC$.  Indeed, 
$(X,e) \oplus (X,1 - e) =$
$(X \oplus X,e \oplus (1 - e)) \approx$ 
$(X,\id_X) = \Phi X$.]\\

%%----------------------------------------------------------------------------------------------08
\begingroup%%----------------------------------->>
\fontsize{9pt}{11pt}\selectfont
\textbf{\small FACT} \  
If \bD is a pseudoabelian category and $F:\bC \ra \bD$ is an additive functor which is full and faithful such that every object in 
\bD is a direct summand of an object in $F\bC$, then $F_{\pa}:\bC_{\pa} \ra \bD$ is an equivalence of categories.\\
\endgroup %%------------------------------------<<

\begingroup%%----------------------------------->>
\fontsize{9pt}{11pt}\selectfont
\textbf{\small EXAMPLE} \  
Suppose that \mX is a compact Hausdorff space.  Let \bC be the category of real or complex trivial vector bundles over \mX $-$then $\bC_{\pa}$ is equivalent to the category of real or complex vector bundles over \mX.
\vspi
[Since \mX is compact Hausdorff, $\forall \ E \ra X$ $\exists$ $E^\prime \ra X$ such that 
$E \oplus E^\prime$ is trivial.]\\
\endgroup %%------------------------------------<<

Let \bC, \bD be categories WES.  Assume: \bC is a full, additive subcategory of \bD with the property that a pair 
$X \overset{i}{\ra} Y \overset{p}{\ra} Z$ is exact in \bC iff it is exact in \bD $-$then \bC is said to be 
\un{cofinal}
\index{cofinal (categories WES)} 
in \bD if for every exact pair 
$X \overset{i}{\ra} Y \overset{p}{\ra} Z$ in \bD, where $X, Z \in \Ob\bC$, $\exists$ an object in \bC which is isomorphic to \mY, and $\forall \ X \in \Ob\bD$, $\exists$ $Z \in \Ob\bD$ such that $X \oplus Z$ is isomorphic to an object in \bC.  
Example: Given a ring \mA with unit, $\bF(A)$ is cofinal in $\bP(A)$.\\

\begingroup%%----------------------------------->>
\fontsize{9pt}{11pt}\selectfont
\label{18.10}
\label{18.27}
\textbf{\small EXAMPLE} \  
Let \bC be a category WES.  Viewing \bC as a full, additive subcategory of $\bC_{\pa}$, stipulate that the elements of 
$\sE_{\pa}$ are those pairs which are direct summands of elements of $\sE$ $-$then $\bC_{\pa}$ is a category WES and \bC is cofinal in $\bC_{\pa}$.\\
\vspace{0.25cm}
\endgroup %%------------------------------------<

If 
$
\begin{cases}
\ \bC\\
\ \bD
\end{cases}
$
are Waldhausen categories and if $F:\bC \ra \bD$ is a functor, then \mF is said to be a 
\un{model functor}
\index{model functor} 
provided that $F0 = 0$, \mF sends weak equivalences to weak equivalences and cofibrations to cofibrations, and \mF preserves pushouts along a cofibration, i.e., for any 2-source 
$X \overset{f}{\la} Z \overset{g}{\ra} Y$, where $f$ is a cofibration, the arrow 
$FX \underset{FZ}{\sqcup} FY \ra F(X \underset{Z}{\sqcup} Y)$ is an isomorphism.\\
\vspace{0.25cm}

\begingroup%%----------------------------------->>
\fontsize{9pt}{11pt}\selectfont
\textbf{\small FACT} \  
Let 
$
\begin{cases}
\ \bC\\
\ \bD
\end{cases}
$
be categories WES viewed as Waldhausen categories (cf. Proposition 4) $-$then an additive functor $F:\bC \ra \bD$ is a model functor iff $X \overset{i}{\ra} Y \overset{p}{\ra} Z$ exact $\implies$ 
$FX \overset{Fi}{\ra} FY \overset{Fp}{\ra} FZ$ exact.
\vspi
[Note: \ In this context, a model functor called an 
\un{exact functor}.]
\index{exact functor (model functor)}\\
\endgroup %%------------------------------------<<

\bWALD is the category whose objects are the small Waldhausen categories and whose morphisms are the model functors between them.\\

\begingroup%%----------------------------------->>
\fontsize{9pt}{11pt}\selectfont
\label{18.13}
\textbf{\small EXAMPLE} \  
Let \bC be a small Waldhausen category $-$then the functor category $[[n],\bC]$ is again in \bWALD (the weak equivalences and cofibrations are levelwise) and $\Ob [[n],\bC] = \ner_n\bC$.  
Write $\bw\bC(n)$ for the full subcategory of $[[n],\bC]$ consisting of those functors that take values in $\bw\bC$, i.e., the diagrams of the form 
$X_0 \overset{f_0}{\ra} X_1$ 
$\ra \cdots \ra $ 
\begin{tikzcd}%[sep=small]
{X_{n-1}} \ar{r}{f_{n-1}} &{X_n}
\end{tikzcd}
, where the $f_i$ are weak equivalences (thus $\Ob\bw\bC(n) = \ner_n\bw\bC$
%%----------------------------------------------------------------------------------------------09
and $[[n],\bw\bC]$ is a subcategory of $\bw\bC(n)$).  Since pushouts are levelwise, $\bw\bC(n)$ inherits the structure of a Waldhausen category  from $[[n],\bC]$.
\vspi
[If $X \overset{f}{\la} Z \overset{g}{\ra} Y$ is a 2-source in $\bw\bC(n)$, where $f$ is a cofibration, then there are commutative diagrams
\begin{tikzcd}[sep=large]
{Z_i} \ar{d}[swap]{f_i} \ar{r}{g_i} &{Y_i} \ar{d}{\eta_i}\\
{X_i} \ar{r}[swap]{\xi_i} &{P_i}
\end{tikzcd}
and the claim is that 
$P_0 \ra P_1$ 
$\ra \cdots \ra $ 
$P_{n-1} \ra P_n$ $\in \Ob\bw\bC(n)$.  But this is implied by WC-4.]\\
\endgroup %%------------------------------------<<

Let \bC be a small Waldhausen category.  
Recalling that $[n](\ra)$ is the arrow category of $[n]$ 
(cf. p. \pageref{18.3}), denote by $\bS_n\bC$ the full subcategory of $[[n](\ra),\bC]$ consisting of those functors 
$F:[n](\ra) \ra \bC$ such that $F(i \ra i) = 0$ $(0 \leq i \leq n)$ and for every triple $i \leq j \leq k$ in $[n]$, 
$F(i \ra j) \ra F(i \ra k)$ is a cofibration and the commutative diagram 
\begin{tikzcd}%[sep=small]
{F(i \ra j)}\ar{d} \ar{r} &{F(j \ra j)}\ar{d}\\
{F(i \ra k)} \ar{r}  &{F(j \ra k)}
\end{tikzcd}
is a pushout square $-$then the assignment $[n] \ra \bS_n\bC$ defines an internal category in \bSISET, call it $\bS\bC$.

[Note: \ Each $\alpha:[m] \ra [n]$ in $\Mor\bDelta$ determines a functor 
$\alpha(\ra):[m](\ra) \ra [n](\ra)$ from which a functor $\bS_n\bC \ra \bS_m\bC$, viz. $F \ra F \circx \alpha(\ra)$.]\\

\textbf{\small LEMMA} \  
$\bS_n\bC$ is a small Waldhausen category.

[The weak equivalences are those natural transformations $\Xi:F \ra G$ such that 
$\Xi_{i \ra j}:F(i \ra j) \ra G(i \ra j)$ is a weak equivalence and the cofibrations are those natural transformations 
$\Xi:F \ra G$ such that $\Xi_{i \ra j}:F(i \ra j) \ra G(i \ra j)$ is a cofibration and for every triple $i \leq j \leq k$ in $[n]$, 
the arrow $F(i \ra k) \underset{F(i \ra j)}{\sqcup} G(i \ra j) \ra G(i \ra k)$ is a cofibration.]

[Note: \ $\bS_0\bC \approx 1$, and $\bS_1\bC \approx \bC$.]\\

\begingroup%%----------------------------------->>
\fontsize{9pt}{11pt}\selectfont
\label{18.12a}
Given a \bC in \bWALD, define a simplicial set $W\bC$ by putting $W_n\bC = \Ob\bS_n\bC$.\\
\endgroup %%------------------------------------<<

\begingroup%%----------------------------------->>
\fontsize{9pt}{11pt}\selectfont
\textbf{\small FACT} \  
Suppose that 
$
\begin{cases}
\ \bC\\
\ \bD
\end{cases}
$
are small Waldhausen categories.  Let $F:\bC \ra \bD$ be a model functor $-$then \mF induces a simplicial map 
$WF:W\bC \ra W\bD$.\\
\endgroup %%------------------------------------<<

\begingroup%%----------------------------------->>
\fontsize{9pt}{11pt}\selectfont
\textbf{\small FACT} \  
Suppose that 
$
\begin{cases}
\ \bC\\
\ \bD
\end{cases}
$
are small Waldhausen categories.  Let $F, G:\bC \ra \bD$ be model functors, $\Xi:F \ra G$ a natural isomorphism 
$-$then $\Xi$ induces a simplicial homotopy between $WF$ and $WG$.\\ 
\endgroup %%------------------------------------<<

\begingroup%%----------------------------------->>
\fontsize{9pt}{11pt}\selectfont
\textbf{\small EXAMPLE} \  
Let \bC be a small Waldhausen category.  Denote by $i\bC(\ra)$ the full subcategory of $\bC(\ra)$ whose objects are the 
$X \overset{f}{\ra} Y$ such that $f$ is an isomorphism $-$then there is a model functor $F:\bC 
%%----------------------------------------------------------------------------------------------10
\ra i\bC(\ra)$, viz.  
$FX =$
\begin{tikzcd}%[sep=small]
{X} \ar{r}{id_X} &{X,}
\end{tikzcd}
and a model functor $G:i\bC(\ra): \ra \bC$, viz. $G(X \overset{f}{\ra} Y) = X$.  Obviously, $G \circx F = \id_{\bC}$ and 
$F \circx G \approx \id_{i\bC(\ra)}$, so $\abs{W\bC}$ and $\abs{Wi\bC(\ra)}$ have the same pointed homotopy type.\\
\endgroup %%------------------------------------<<

\begin{proposition} \ %05
Let \bC be a small Waldhausen category $-$then $\bS\bC$ is a simplicial object in \bWALD.
\end{proposition}

[The $d_i$ and the $s_i$ are model functors.]

[Note: \ A model functor $\bC \ra \bD$ induces a model functor $\bS\bC \ra \bS\bD$.  Therefore \bS is a functor from 
\bWALD to \bSIWALD ($=[\bDelta^\OP,\bWALD]$).]\\

\label{18.19}
\label{18.20}
Given a small Waldhausen category \bC, let \ $B\bw\bS\bC = \abs{[n] \ra B\bw\bS_n\bC}$ 
$-$then $B\bw\bS\bC$ is path connected and there is a closed embedding 
$\Sigma B\bw\bC \ra B\bw\bS\bC$.  \ Now iterate the process, i.e., form 
$\bS^{(2)}\bC = \bS\bS\bC$, a bisimplicial object in \bWALD, and in general, 
$\bS^{(q)}\bC = \bS \cdots \bS\bC$, a multisimplicial object in \bWALD.  
Write $\bw\bS^{(q)}\bC$ for the weak equivalences in 
$\bS^{(q)}\bC$.  \ 
If $B\bw\bS^{(q)}\bC$ is its classifying space (see below), 
then $B\bw\bS^{(q)}\bC$ is $(q-1)$-connected $(q > 1)$ and there is a closed embedding 
$\Sigma B\bw\bS^{(q)}\bC \ra B\bw\bS^{(q+1)}\bC$ whose adjoint 
$B\bw\bS^{(q)}\bC \ra \Omega B\bw\bS^{(q+1)}\bC$ is a pointed homotopy equivalence 
(cf. p. \pageref{18.4}).  
The data can be assembled into a separated prespectrum $\bW\bC$, where $(\bW\bC)_0 = B\bw\bC$ and 
$(\bW\bC)_q = B\bw\bS^{(q)}\bC$ $(q \geq 1)$.  
Definition: The spectrum $\bK\bC = e\bW\bC$ is the 
\un{algebraic K-theory}
\index{algebraic K-theory} 
of \bC, its homotopy groups 
$\pi_n(\bK\bC) (\approx \pi_n(\Omega B\bw\bS\bC))$ being the 
\un{algebraic K-groups}
\index{algebraic K-groups} 
$K_n(\bC)$ of \bC.

[Note: \ $\bK\bC$ is connective.  In addition, $\bK\bC$ is tame (since $\bW\bC$ satisfies the cofibration condition).]

Remark: A model functor $F:\bC \ra \bD$ determines a morphism $\bW\bC \ra \bW\bD$ of prespectra , hence a morphism $\bK\bC \ra \bK\bD$ of spectra.  Therefore 
$\bK:\bWALD \ra \bSPEC$ is a functor.

[Note: \ If $B\bw\bS\bC \ra B\bw\bS\bD$ is a weak homotopy equivalence, then $\forall \ q$,
$B\bw\bS^{(q)}\bC \ra B\bw\bS^{(q)}\bD$ is a weak homotopy equivalence or still, a pointed homotopy equivalence, so 
$\bK\bC \ra \bK\bD$ is a homotopy equivalence of spectra (cf. p. \pageref{18.5}).]

Convention: If \bC is an aribitrary Waldhausen category, then \bC is not necessarily small.  However \bC is skeletally small (by definition) and all of the above is applicable to a skeleton $\ov{\bC}$, thus $\bK\bC \equiv \bK\ov{\bC}$ and 
$K_n(\bC) \equiv K_n(\ov{\bC})$.

[Note: \ If \bC is small to begin with, then $B\bw\bS\bC$ and $B\bw\bS\ov{\bC}$ have the same pointed homotopy type, so this is a consistent agreement.]\\

\begingroup%%----------------------------------->>
\fontsize{9pt}{11pt}\selectfont
If $X:(\bDelta \times \cdots \times \bDelta)^\OP \ra \bCG$ is a compactly generated multisimplicial space, then its 
\un{geometric} \un{realization}
\index{geometric realization (compactly generated multisimplicial space)} 
is the coend 
$X \otimes_{\bDelta \times \cdots \times \bDelta}(\Delta^? \times_k \cdots \times_k \Delta^?)$, which is homeomorphic to 
$\abs{\di X}$, the geometric realization of $\di X$ (the diagonal of \mX (cf. p. \pageref{18.6})).\\
\endgroup %%------------------------------------<<

%%----------------------------------------------------------------------------------------------11
\label{18.34} %dmc mnft
\begingroup%%----------------------------------->>
\fontsize{9pt}{11pt}\selectfont
\textbf{\small EXAMPLE} \  
If \bC is an internal category in \bSISET, i.e., a simplicial object in \bCAT, then $\ner \bC$ is a bisimplicial set or still, a functor 
$(\bDelta \times \bDelta)^\OP \ra \bSET$ ($\subset \bCG$) and its geometric realization is the classifying  space $B\bC$ of \bC 
(thus $B\bC \approx \abs{[n] \ra B\bC_n}$).
\vspi
[Note: \ Analogous considerations apply to the multisimplicial objects in \bCAT.]\\
\endgroup %%------------------------------------<<

\begingroup%%----------------------------------->>
\fontsize{9pt}{11pt}\selectfont
\textbf{\small EXAMPLE} \  
If \bC is an internal category in \bCAT , i.e., a double category, then the classifying  space $B\bC$ of \bC is the geometric realization of the bisimplicial set $\ner(\ner \bC)$ (cf. p. \pageref{18.7}).  
Example: Let \bA be a subcategory of $\bB,$ where \bB is small.  
Call $\bA \cdot \bB$ the double category whose objects are those of \textbf{B}, with horizontal morphisms $= \Mor\bB$ and vertical morphisms $= \Mor \bA$, and whose bimorphisms are the commutative squares with horizontal arrows in \bB and vertical arrows in \bA.  View \bB as the double category 
\begin{tikzcd}%[sep=small]
\ {\text{\textbullet}} \arrow[d,shift right=0.5,dash] \arrow[d,shift right=-0.5,dash]\ar{r}
&{\text{\textbullet}} \arrow[d,shift right=0.5,dash] \arrow[d,shift right=-0.5,dash]\\
\ {\text{\textbullet}} \ar{r}&{\text{\textbullet}}
\end{tikzcd}
$-$then the inclusion $\bB \ra \bA \cdot \bB$ induces a homotopy equivalence $B\bB \ra B \bA \cdot \bB$.\\
\endgroup %%------------------------------------<<

\label{18.33}
\begingroup%%----------------------------------->>
\fontsize{9pt}{11pt}\selectfont
\textbf{\small FACT} \  
If \bC is a small Waldhausen category, then there is a pointed homotopy equivalence $\abs{W\bC} \ra B\iso \bS\bC$.\\
\endgroup %%------------------------------------<<

\begingroup%%----------------------------------->>
\fontsize{9pt}{11pt}\selectfont
\textbf{\small EXAMPLE} \  
Let \bC be the Waldhausen category whose objects are the pointed finite sets, where weak equivalence = isomorphism and cofibration = pointed injection $-$then $\bGamma$ is a skeleton of \bC, hence is a small Waldhausen category 
(cf. Proposition 1), 
and a model for $\abs{W\bGamma}$ in the pointed homotopy category is 
$\Omega^\infty\Sigma^\infty\bS^1$.  
Proof: Thanks to the homotopy colimit theorem, $\Omega^\infty\Sigma^\infty\bS^1$ can be identified with 
$\ohc \pow \bS^1$.  But, in the notation of p. \pageref{18.8}, 
$\ohc \pow \bS^1 \approx \ov{\pow}\bS^1 \otimes_{\bGamma} \gamma^\infty$ 
$\approx \abs{\gamma^\infty}_{\bGamma}$ 
$\approx B\abs{M_\infty}$ 
$\approx \abs{W \bGamma}$, where $\abs{M_\infty} = \ds\coprod\limits_{n \geq 0} BS_n$.  Therefore the loop space of 
$B\iso\bS\bGamma$ is pointed homotopy equivalent to 
$\Omega\Omega^\infty\Sigma^\infty\bS^1 \approx \Omega^\infty\Sigma^\infty\bS^0$, so the algebraic K-groups 
$K_*(\bGamma)$ of $\bGamma$ ``are'' the $\pi_*^s$, the stable homotopy groups of spheres.
\vspi
[Note: \ More is true, namely $\bK\bGamma$ and \bS, when viewed as objects in \bHSPEC, are isomorphic 
(Rognes\footnote[2]{\textit{Topology} \textbf{31} (1992), 813-845}).]\\
\endgroup %%------------------------------------<<

\begingroup%%----------------------------------->>
\fontsize{9pt}{11pt}\selectfont
\textbf{\small EXAMPLE} \  
Let \bC be a small category WES, $\bC\bX\bC^\tb$ the category of bounded cochain complexes over \bC.  
Suppose that 
$(\bA,\iota)$ is a pair satisfying the embedding condition and the $E \Rightarrow D$ axiom.  
Equip $\bC\bX\bC^\tb$ with the structure of a small Waldhausen category by stipulating that the weak equivalences are the arrows in $\bC\bX\bC^\tb$  which are quasiisomorphisms in \bA and the cofibrations are the levelwise inflations 
$-$then the exact functor 
$\bC \ra \bC\bX\bC^\tb$ sending \mX to \mX concentrated in degree 0 induces a homotopy equivalence 
$\bK\bC \ra \bK\bC\bX\bC^\tb$ of spectra 
(Thomason-Trobaugh\footnote[3]{\textit{The Grothendieck Festschrift}, vol. III, Birkh\"auser (1990), 247-435 (cf. 278-283).}).
\vspi
%%----------------------------------------------------------------------------------------------12
[Note: \ The definition of weak equivalence is independent of the choice $(\bA,\iota)$.  Recall that when \bC is pseudoabelian one can take for $(\bA,\iota)$ the pair $(\bG$-$\bQ,\iota)$ (cf. p. \pageref{18.9}).]\\
\endgroup %%------------------------------------<<

\begin{proposition} \ %06
Let \bC be a small Waldhausen category $-$then $K_0(\bC)$ is the free abelian group on generators $[X]$ 
$(X \in \Ob\bC)$ subject to the relations 
(i) $[X] = [Y]$ if $\exists$ a weak equivalence $X \ra Y$ and 
(ii) $[Y] = [X] + [Y/X]$ for every sequence 
$X \rightarrowtail Y \ra Y/X$.
\end{proposition}

[Since $K_0(\bC) \approx \pi_1(B\bw\bS\bC)$, $K_0(\bC)$ is the free group on generators $[X]$ $(X \in \Ob\bC)$ 
subject to the relations 
(i) $[X] = [Y]$ if $\exists$ a weak equivalence $X \ra Y$ and 
(ii) $[Y] = [X]\cdot[Y/X]$ for every sequence $X \rightarrowtail Y \ra Y/X$.  Applying the second relation to 
\begin{tikzcd}%[sep=small]
{X} \ar{r}{\ini_X} &{}%&{X \amalg Y} 
\end{tikzcd}
\begin{tikzcd}%[sep=small]
{X \amalg Y} \ar{r}{\pr_Y} &{Y}
\end{tikzcd}
$\&$ 
\begin{tikzcd}%[sep=small]
{Y} \ar{r}{\ini_Y} &{X \amalg Y} \ar{r}{\pr_X} &{X}
\end{tikzcd}
gives 
$[X\amalg Y] = [X] \cdot [Y]$ $\&$ $[X \amalg Y] = [Y] \cdot [X]$, thus $K_0(\bC)$ is abelian and one uses additive notation 
$([0] = 0)$.]

[Note: \ If $X \overset{f}{\la} Z \overset{g}{\ra} Y$ is a 2-source, where $f$ is a cofbration, then 
$[P] = [Y] + [P/Y]$ $=$ $[Y] + [X/Z] =$ $[X] + [Y] - [Z]$.]\\

Example: Suppose that \bC satisfies the mapping cylinder axiom $-$then $\forall \ X \in \Ob\bC$, there is a weak equivalence $\Gamma X \ra 0$, hence $[X] = -[\Sigma X]$.

[Note: \ Under these circumstances, every element of $K_0(\bC)$ is a $[X]$ for some $X \in \Ob\bC$.  
Proof: $[Y] - [Z] = [Y \amalg \Sigma Z]$.]\\

\begingroup%%----------------------------------->>
\fontsize{9pt}{11pt}\selectfont
\textbf{\small EXAMPLE} \  
Let \bC be the category whose objects are the pointed finite CW complexes and whose morphisms are the pointed skeletal maps $-$then \bC is a Waldhausen category if the weak equivalences are the weak homotopy equivalences and the cofibrations are the closed cofibrations which are isomorphic to the inclusion of a subcomplex.  Put 
$A(*) = \Omega B\bw\bS\ov{\bC}$ (the 
\un{algebraic K-theory of a point}
\index{algebraic K-theory of a point}) 
$-$then the 
reduced Euler characteristic $\widetilde{\chi}$ defined by $K \ra \chi(K) - 1$ is an isomorphism from $\pi_0(A(*))$ onto $\Z$.
\vspi
[Note: Dwyer\footnote[2]{\textit{Ann. of Math.} \textbf{111} (1980), 239-251.}
has shown that the homotopy groups of $A(*)$ are finitely generated.  Structurally, in the pointed homotopy category there exists a splitting 
$A(*) \approx \Omega^\infty\Sigma^\infty\bS^0 \times Wh^\text{DIFF}(*)$ 
(Waldhausen\footnote[3]{\textit{Ann. of Math. Studies} \textbf{113} (1987), 392-417.}), 
so 
$\pi_q(A(*)) \approx \pi_q^s \oplus \pi_q(Wh^\text{DIFF}(*))$.  Here $Wh^\text{DIFF}(*)$ is the Whitehead space of a point.  
It has the property that there is a pointed homotopy equivalence 
$\Omega^2 Wh^\text{DIFF}(*) \ra P(*)$, the stable smooth pseudoisotopy space of $*$.  Rationally, it is known that 
$\pi_q(Wh^\text{DIFF}(*)) \otimes \Q = \Q$ if $q \equiv 5 \mod 4$ and is zero otherwise, but the explicit determination of the torsion is difficult and unresolved.]\\
\endgroup %%------------------------------------<<

\begingroup%%----------------------------------->>
\fontsize{9pt}{11pt}\selectfont
\textbf{\small EXAMPLE} \  
Let \bC be a small category WES $-$then \bC has finite coproducts (= finite products), thus \bC can be viewed as a symmetric monoidal category.  Therefore the isomorphism classes of \bC constitute  
%%----------------------------------------------------------------------------------------------13
an abelian monoid, call it \mM.  
Definition: $K^\oplus(\bC) = \ov{M}$, the group completion of \mM.  So $K_0(\bC)$ is a quotient of $K_0^\otimes(\bC)$, the two being the same if every exact pair $X \overset{i}{\ra} Y \overset{p}{\ra} Z$ splits.  (i.e., is isomorphic to 
$X \overset{\ini_X}{\lra} X \oplus Z \overset{\pr_Z}{\lra} Z$).\\
\endgroup %%------------------------------------<<

\begingroup%%----------------------------------->>
\fontsize{9pt}{11pt}\selectfont
\textbf{\small FACT} \  
Let \bC, \bD be small categories WES.  
Assume: \bC is cofinal in \bD $-$then $K_0(\bC)$ is a subgroup of  $K_0(\bD)$.
\vspi
[Observe first that $K^\oplus(\bC)$ is a subgroup of $K^\oplus(\bD)$.  
This said, suppose in addition that \bC is isomorphism closed in \bD.  Given an exact pair 
$X \overset{i}{\ra} Y \overset{p}{\ra} Z$ in \bD, choose $X^\prime, Z^\prime$ in \bD such that 
$X \oplus X^\prime$, $Z \oplus Z^\prime$ are in \bC $-$then 
$X \oplus X^\prime \ra Z^\prime \oplus Y \oplus X^\prime$
$\ra  Z^\prime \oplus Z$ is exact in \bD, hence $Z^\prime \oplus Y \oplus X^\prime \in \Ob\bC$.  
Consequently, in $K^\oplus(\bD)$, 
$[Z^\prime \oplus Y \oplus X^\prime] - [X \oplus X^\prime] - [Z^\prime \otimes Z] = $ 
$[Z^\prime] + [Y] + [X^\prime] - [X] - [X^\prime] - [Z^\prime] - [Z] = $ $[Y] - [X] - [Z]$, 
thus the kernel of 
$K_0^\oplus(\bC) \ra K_0(\bC)$ equals the kernel of  $K_0^\oplus(\bD) \ra K_0(\bD)$, 
which implies that the arrow $K_0(\bC) \ra K_0(\bD)$ is one-to-one.]\\
\endgroup %%------------------------------------<<

\label{18.26}
\begingroup%%----------------------------------->>
\fontsize{9pt}{11pt}\selectfont
\textbf{\small EXAMPLE} \  
Let \bC be a small category WES $-$then \bC is cofinal in $\bC_{\pa}$ (cf. p. \pageref{18.10}), so $K_0(\bC)$ is a subgroup of $K_0(\bC_{\pa})$.
\vspi
[Note: \ Let \mA be a ring with unit $-$then $K_0(\bP(A)) = K_0(A)$ and $\bF(A)$ is cofinal in $\bP(A)$.  The arrow 
$\Z_{\geq 0} \ra \bP(A)$ that sends $n$ to $A^n$ induces a homomorphism $\Z \ra K_0(A)$ of groups (injective iff \mA has the invariant basis property (i.e., $m \neq n$ $\implies$ $A^m \not\approx A^n$)).  Since $\bF(A)_{\pa} = \bP(A)$, it follows that the cyclic group $K_0(\bF(A))$ is a subgroup of $K_0(A)$.]\\
\endgroup %%------------------------------------<<

\begin{proposition} \ %07
Suppose that 
$
\begin{cases}
\ \bC\\
\ \bD
\end{cases}
$
are small Waldhausen categories.  Let $F, G:\bC \ra \bD$ be model functors, $\Xi:F \ra G$ be a natural transformation such that 
$\forall \ X \in \Ob\bC$, $\Xi_X:FX \ra GX$ is a weak equivalence in \bD $-$then $\Xi$ induces a spectral homotopy between 
$\bK F$ and $\bK G$ (cf. p. \pageref{18.11} and $\S 14$, Proposition 12).
\end{proposition}

[Note: \ One starts from the pointed homotopy $B\bw\bS F \simeq B \bw\bS G$.]\\

\begingroup%%----------------------------------->>
\fontsize{9pt}{11pt}\selectfont
\textbf{\small EXAMPLE} \  
Suppose that \bC satisfies the mapping cylinder axiom $-$then $\forall \ X \in \Ob\bC$, there is a weak equivalence 
$\Gamma X \ra 0$.  But $\Gamma:\bC \ra \bC$ is a model functor, hence the induced map 
$B\bw\bS\bC \ra B\bw\bS\bC$ is nullhomotopic.\\
\endgroup %%------------------------------------<<

Let \bC, $\bC^\prime$, $\bC\pp$ be small Waldhausen categories.  Assume $\bC^\prime$ and $\bC\pp$ are subcategories of \bC with the property that the inclusions $\bC^\prime \ra \bC$, $\bC\pp \ra \bC$ are model functors.  
Denote by $\bE(\bC^\prime, \bC, \bC\pp)$ the category whose objects are the pushout squares 
\begin{tikzcd}%[sep=small]
{X^\prime} \ar[d,tail] \ar{r} &{0} \ar{d}\\
{X} \ar{r} &{X\pp}
\end{tikzcd}
in \bC, where $X^\prime \in \Ob\bC^\prime$, $X \in \Ob\bC$, $X\pp \in \Ob\bC\pp$, and whose morphisms are the commutative
%%----------------------------------------------------------------------------------------------14
digarams 
\begin{tikzcd}%[sep=small]
{X^\prime} \ar{d} \ar[r,tail]  &{X} \ar{d} \ar{r} &{X\pp} \ar{d}\\
{Y^\prime} \ar[r,tail] &{Y}  \ar{r} &{Y\pp}
\end{tikzcd}
in \bC, where $X^\prime \ra Y^\prime \in \Mor \bC^\prime$, $X \ra Y \in \Mor \bC$, $X\pp \ra Y\pp \in \Mor \bC\pp$.

[Note: \ When $\bC^\prime = \bC$ and $\bC\pp = \bC$, put $\bE\bC = \bE(\bC,\bC,\bC)$.]\\

\textbf{\small LEMMA} \  
$\bE(\bC^\prime, \bC, \bC\pp)$  is a small Waldhausen category.

[A morphism in $\bE(\bC^\prime, \bC, \bC\pp)$ is a weak equivalence if $X^\prime \ra Y^\prime$ is a weak equivalence in $\bC^\prime$, $X \ra Y$ is a weak equivalence in \bC, $X\pp \ra Y\pp$ is a weak equivalence in $\bC\pp$ and a morphism in $\bE(\bC^\prime, \bC, \bC\pp)$ is a cofibration if 
$X^\prime \ra Y^\prime$ is a cofibration in $\bC^\prime$, 
$Y^\prime \underset{X^\prime}{\sqcup} X \ra Y$ is a cofibration in \bC, 
$X\pp \ra Y\pp$ is a cofibration in $\bC\pp$.

[Note: \ $X \ra Y$ is then a cofibration in \bC (being the composite 
$X^\prime \underset{X^\prime}{\sqcup} X \rightarrowtail Y^\prime \underset{X^\prime}{\sqcup} X \rightarrowtail Y$ 
(cf. $\S 12$, Proposition 4)).]\\

There are model functors 
$s: \bE(\bC^\prime,\bC,\bC\pp) \ra \bC^\prime$,
$t: \bE(\bC^\prime,\bC,\bC\pp) \ra \bC$,
$Q: \bE(\bC^\prime,\bC,\bC\pp) \ra \bC\pp$,
viz. 
$s(X^\prime \rightarrowtail X \ra X\pp) = X^\prime$,
$t(X^\prime \rightarrowtail X \ra X\pp) = X$,
$Q(X^\prime \rightarrowtail X \ra X\pp) = X\pp$.  
In the other direction, there is a model functor 
$I:\bC^\prime \times \bC\pp \ra \bE(\bC^\prime,\bC,\bC\pp)$ which sends 
$(X^\prime,X\pp)$ to $X^\prime \rightarrowtail X^\prime \amalg X\pp \ra X\pp$.  
Agreeing to write $(s,Q)$ for the model functor 
$\bE(\bC^\prime,\bC,\bC\pp) \ra \bC^\prime \times \bC\pp$ defined by s and Q, viz. 
$(s,Q)(X^\prime \rightarrowtail X \ra X\pp) = (X^\prime,X\pp)$, one has 
$(s,Q) \circx I = \id_{\bC^\prime \times \bC\pp}$.\\

\index{Theorem: Relative Additivity Theorem}
\index{Relative Additivity Theorem}
\textbf{\small RELATIVE ADDITIVITY THEOREM} \quad 
The model functor $(s,Q)$ induces a homotopy equivalence 
$\bK(s,Q):\bK\bE(\bC^\prime,\bC,\bC\pp) \ra \bK\bC^\prime \times \bK\bC\pp$ of spectra.\\

\index{Theorem: Absolute Additivity Theorem}
\index{Absolute Additivity Theorem}
\textbf{\small ABSOLUTE ADDITIVITY THEOREM} \quad 
The model functor $(s,Q)$ induces a homotopy equivalence 
$\bK(s,Q):\bK\bE\bC \ra \bK\bC \times \bK\bC$ of spectra.\\

It \ is \ a\  question\  of proving \ that \  $(s,Q)$  \ induces \ a weak \  homotopy \  equivalence 
$B\bw\bS\bE(\bC^\prime,\bC,\bC\pp)$ $\ra$ 
$B\bw\bS\bC^\prime \times_k B\bw\bS\bC\pp$ of classifying spaces.  To this end, 
we shall proceed via a series of lemmas.\\

\index{Homotopy Lemma}
\textbf{\small HOMOTOPY LEMMA} \quad 
Grant the truth of the absolute additivity theorem $-$then 
$B\bw\bS t:B\bw\bS\bE\bC \ra B\bw\bS\bC$ is pointed homotopic to 
$B\bw\bS(s\amalg Q):B\bw\bS\bE\bC \ra B\bw\bS\bC$.

[Note: \ Here $(s \amalg Q) (X ^\prime \rightarrowtail X \ra X\pp) = X^\prime \amalg X\pp$.]\\

\index{Triad Lemma}
\textbf{\small TRIAD LEMMA} \quad 
Grant the truth of the absolute additivity theorem.  Suppose given a small Waldhausen category \bD, model functors 
$G, G^\prime, G\pp: \bD \ra \bC$, and natural transformations $G^\prime \ra G$, $G \ra G\pp$.  
Assume: 
(i) For every object \mX in \bD, the arrow $G^\prime X \ra GX$
%%----------------------------------------------------------------------------------------------15
is a cofibration and the commutative diagram 
\begin{tikzcd}%[sep=small]
{G^\prime X}\ar[d,tail] \ar{r}&{0}\ar{d}\\
{G X}\ar{r}  &{G\pp X}
\end{tikzcd}
is a pushout square;
(ii) For every cofibration $X \rightarrowtail Y$ in \bD, the arrow $G^\prime Y \underset{G^\prime X}{\sqcup} GX \ra GY$ is a cofibration $-$then $B\bw\bS G$ is pointed homotopic to $B\bw\bS(G^\prime \amalg G\pp)$.

[There exists a model functor 
$\bPhi:\bD \ra \bE\bC$ with $G^\prime = s \circx \bPhi$, $G = t \circx \bPhi$, 
$G\pp = Q \circx \bPhi$.  The assertion thus follows from the homotopy lemma by naturality.]\\

\begingroup%%----------------------------------->>
\fontsize{9pt}{11pt}\selectfont
\textbf{\small EXAMPLE} \  
Let \bC be a Waldhausen category whose objects are the pointed finite CW complexes and whose morphisms are the pointed skeletal maps $-$then the arrow 
$B\bw\bS\ov{\bC} \ra B\bw\bS\ov{\bC}$ induced by $\Sigma$ is a pointed homotopy equivalence.
\vspi
[In the triad lemma, take $\ov{\bC} = \ov{\bD}$ and let $G^\prime = \id_{\ov{\bC}^\prime} G = \Gamma$, 
$G\pp = \Sigma$.]
\vspi
[Note: \ The full subcategory $\bC_0$ of $\ov{\bC}$ whose objects are path connected is Waldhausen (WC-3 is a consequence of AD$_1$ (cf. p. \pageref{18.12})).  Since there is a commutative diagram 
\begin{tikzcd}[sep=large]
{B\bw\bS\ov{\bC}}\ar{r}{B\bw\bS\Sigma}  \ar[rd,dashed] &{B\bw\bS\ov{\bC}}\\
{B\bw\bS\bC_0} \ar{u}{B\iota} \ar{r}[swap]{B\bw\bS\Sigma}  &{B\bw\bS\bC_0} \ar{u}[swap]{B\iota}
\end{tikzcd}
, it follows that $B\iota$ is a pointed homotopy equivalence.  
Therefore, the algebraic K-theory of a point can be defined using path connected objects.  
If now $\bC_1$ is the full subcategory of $\bC_0$ whose objects are simply connected, then $\bC_1$ is Waldhausen (WC-3 is implied by the Van Kampen theorem).  
Repeating the argument, one concludes that the algebraic K-theory of a point can be defined by using simply connected objects.  As an aside, observe that $\bC_1$ satisfies the extension axiom (via the Whitehead theorem) but $\ov{\bC}$ does not.]\\ 
\endgroup %%------------------------------------<<

\index{Lemma of Reduction}
\textbf{\small LEMMA OF REDUCTION} \quad 
The absolute additivity theorem implies the relative additivity theorem.

[Since $(s,Q) \circx I = \id_{\bC^\prime \times \bC\pp}$, it suffices to show that $B\bw\bS(I \circx (s,Q))$ is pointed homotopic to the identity.   Accordingly, to apply the triad lemma, define model functors 
$G^\prime, G , G\pp: \bE(\bC^\prime,\bC,\bC\pp) \ra \bE(\bC^\prime,\bC,\bC\pp)$ by 
$G^\prime(X^\prime \rightarrowtail X \ra X\pp) =$
% X^\prime \overset{\id_{X^\prime}}{\rightarrowtail} X^\prime \ra 0$,  
\begin{tikzcd}[sep=small]
X^\prime \arrow[rr,tail,"\id_{X^\prime}"] &&X^\prime \arrow[r] &0,
\end{tikzcd}
$G(X^\prime \rightarrowtail X \ra X\pp) = X^\prime \rightarrowtail X \ra X\pp$, 
$G\pp(X^\prime \rightarrowtail X \ra X\pp) = $
%0 \rightarrowtail X\pp \overset{\id_{X^{\prime\prime}}}{\rightarrowtail}  X\pp$ 
\begin{tikzcd}[sep=small]
0 \arrow[r,tail] &X\pp \arrow[rr,tail,"\id_{X\pp}"] &&X\pp 
\end{tikzcd}
and note that 
$B\bw\bS(I \circx (s,Q)) = B\bw\bS(G^\prime \amalg G\pp)$.]\\

\index{Additivity Lemma}
\textbf{\footnotesize ADDITIVITY LEMMA} \quad 
The simplicial map 
$W(s,Q): W\bE\bC \ra W\bC \times W\bC$ induced by $(s,Q)$ is a weak homotopy equivalence 
(notation as on p. \pageref{18.12a}).\\

The additivity lemma implies the absolute additivity theorem.  To see this, introduce $\bw\bC(n)$ 
(cf. p. \pageref{18.13} ff.) $-$then $\forall \ n$, the arrow 
$W\bE\bw\bC(n) \ra W\bw\bC(n) \times W\bw\bC(n)$ is 
%%----------------------------------------------------------------------------------------------16
a weak homotopy equivalence.  Therefore the diagonal of the bisimplicial map 
$([n] \ra W\bE\bw\bC(n)) \ra ([n] \ra W\bw\bC(n)) \times ([n] \ra W\bw\bC(n))$ is a weak homotopy equivalence 
(cf. p. $\S 13$, Proposition 51) or still, the induced map of geometric realizations is a weak homotopy equivalence.  It remains only to observe that $\Ob\bS_m\bw\bC(n) \approx \ner_n\bw\bS_m\bC$.\\

\textbf{\small LEMMA} \  
The projection $W\bE\bC \overset{p}{\ra} W\bC$ induced by $s$ is a homotopy fibration (cf. infra).\\

This result leads to the additivity lemma.  In fact, $\forall \ n$ $\&$ $\forall \ x \in W_n\bC$, the pullback square 
\begin{tikzcd}[sep=large]
{F_x} \ar{d} \ar{r} &{W\bE\bC} \ar{d}{p} \\
{\Delta [n]} \ar{r}[swap]{\Delta_x} &{W\bC}
\end{tikzcd}
$(F_x = W\bE\bC_x)$ is a homotopy pullback (cf. p. \pageref{18.14}).  
Now take $n = 0$ and recall that $W_0\bC = *$ $-$then $F_0 \ra W\bE\bC \overset{p}{\ra} W\bC$ is a homotopy pullback and 
$F_0$ can be identified with $W\bF_0\bC$, $\bF_0\bC$ being the full subcategory of $\bE\bC$ whose objects are the 
$0 \rightarrowtail X \ra X\pp$ $(\implies X \approx X\pp)$.  But the model functor $\bF_0\bC \ra \bC$ defined by 
\begin{tikzcd}%[sep=small]
{0}\ar[r,tail]
&{X} \ar{d} \ar{r}
&{X\pp}\\
&{X}
\end{tikzcd}
gives rise to a homotopy equivalence $W\bF_0\bC \ra W\bC$ of simplicial sets.  Therefore the sequence 
$W\bC \ra W\bE\bC \overset{p}{\ra} W\bC$ is a homotopy pullback (the arrow $W\bC \ra W\bE\bC$ corresponds to the insertion 
$\bC \ra \bE\bC$ which sends \mX to $0 \rightarrowtail X \overset{\id_X}{\lra} X$).  Consider the diagram
\ 
\begin{tikzcd}%[sep=small]
{W\bC}\arrow[d,shift right=0.5,dash] \arrow[d,shift right=-0.5,dash] \ar{r}
&{W\bC \times W\bC} \ar{d} \ar{r}
&{W\bC}\arrow[d,shift right=0.5,dash] \arrow[d,shift right=-0.5,dash]\\
{W\bC}  \ar{r}
&{W\bE\bC} \ar{r}
&{W\bC}
\end{tikzcd}
\ 
, where the vertical arrow is determined by \mI.  Passing to the geometric realizations, the top and the bottom rows become fibrations up to homotopy (per \bCGH (singular structure) (cf. p. \pageref{18.15})), thus 
$\abs{WI}: \abs{W\bC} \times_k \abs{W\bC} \ra \abs{W\bE\bC}$ is a pointed homotopy equivalence.  Since 
$\abs{W(s,Q)} \circx \abs{WI} = \id_{\abs{W\bC} \times_k \abs{W\bC} }$, it follows that $\abs{W(s,Q)}$ is also a pointed homotopy equivalence, the assertion of the additivity lemma.\\

\begingroup%%----------------------------------->>
\fontsize{9pt}{11pt}\selectfont
Put $X = W\bE\bC$, $B = W\bC$ $-$then to prove the lemma, one must show that for every commutative diagram
\begin{tikzcd}[sep=large]
{X_{b^\prime}} \ar{d}\ar{r}
&{X_b} \ar{d}\ar{r}
&{X} \ar{d}{p}\\
{\Delta[n^\prime]} \ar{r}
&{\Delta[n]} \ar{r}[swap]{\Delta_b}
&{B}
\end{tikzcd}
, the arrow $X_{b^\prime} \ra X_b$ is a weak homotopy equivalence (cf. p. \pageref{18.16}).  
Since any map $[n^\prime] \ra [n]$ can be placed in a commutative triangle 
$
\begin{tikzcd}[sep=small]
&{[0]} \ar{ldd}\ar{rdd}\\
\\
{[n^\prime]} \ar{rr} &&{[n]}
\end{tikzcd}
, 
$
there is no loss of generality in supposing that $n^\prime = 0$, thus our objective may be recast.\\
\endgroup %%------------------------------------<<

%%----------------------------------------------------------------------------------------------17
\textbf{\small LEMMA} \  
\begingroup%%----------------------------------->>
\fontsize{9pt}{11pt}\selectfont
Fix an element $b \in B_n$ and let $v_i:X_{b^\prime} \ra X_b$ be the simplicial map attached to the $i^\text{th}$ vertex operator $\epsilon_i:[0] \ra [n]$ $(0  \leq i \leq n)$ $-$then $v_i$ is a homotopy equivalence.
\vspi
[From the definitions 
$x \in X_m$ $(= W_m\bE\bC) \leftrightarrow F^\prime \rightarrowtail F \ra  F\pp \in \Ob\bE\bS_m\bC$.  
And: An element of 
$(X_b)_m$ consists of an element of $X_m$ plus a map $\alpha:[m] \ra [n]$ such that $F^\prime$ is equal to the composite 
$[m](\ra) \overset{\alpha_*}{\lra} [n](\ra) \overset{b}{\ra} \bC$.  
There is an evident homotopy equivalence 
$W\bC \overset{f}{\ra} X_{b^\prime}$ and $\forall \ i$, $q \circx v_i \circx f = \id_{W\bC}$, where $q:X_b \ra W\bC$ is induced 
by the functor that takes $F^\prime \rightarrowtail F \ra F\pp$ to $F\pp$.  
It will be enough to show that $q$ is a homotopy equivalence and for this it will be enough to show that $\id_{X_b} \simeq v_n \circx f \circx q$.  Let $X_b^*$ be the composite 
$(\bDelta/[1])^\OP \ra \bDelta^\OP \overset{X_b}{\lra} \bSET$ 
and define a natural transformation $H:X_b^* \ra X_b^*$ by 
assigning to 
$\beta:[m] \ra [1]$ the function $H_\beta \in \Mor((X_b)_m,(X_b)_m)$ 
which sends 
$(F^\prime \rightarrowtail F \ra F\pp, \alpha:[m] \ra [n])$ to 
$(\ov{F}{}^\prime \rightarrowtail \ov{F} \ra \ov{F}{}\pp, \ov{\alpha}:[m] \ra [n])$.  
Here $\ov{\alpha}$ is the composite 
\begin{tikzcd}%[sep=small]
{[m]} \ar{r}{(\alpha,\beta)} &{[n] \times [1]}
\end{tikzcd}
$\overset{\gamma}{\ra} [n]$  $(\gamma(j,0) = j$, $\gamma(j,1) = n)$ and $\ov{F}{}^\prime = b \circx \ov{\alpha}_*$.  
Because $\alpha \leq \ov{\alpha}$, $\exists$ a natural transformation $\alpha_* \ra \ov{\alpha}_*$, 
hence $\exists$ a natural transformation $F^\prime \ra \ov{F}{}^\prime$ and $\ov{F}$ is given by the pushout square 
\begin{tikzcd}[sep=large]
{F^\prime}\ar[d,tail] \ar{r} &{ \ov{F}{}^\prime}\ar{d}\\
{F} \ar{r} &{\ov{F}}
\end{tikzcd}
in $\bS_m\bC$ with $\ov{F}{}\pp = \ov{F}/\ov{F}{}^\prime$.  
Needless to say, this procedure involves certain choices and it is necessary to check that they can be made in such a way that 
\mH really is natural.  
Leaving this as an exercise, let us note only that matters can be arranged so that the homotopy starts at the identity 
(viz., if $F^\prime \ra \ov{F}{}^\prime$ is the identity, choose $F \ra \ov{F}$ to be the identity) 
and that the image of $v_n \circx f$ is fixed under the homotopy 
(viz., if $(\ov{F}{}^\prime = 0$, choose $\ov{F} \ra \ov{F}{}\pp$ to be the indentity).]\\
\endgroup %%------------------------------------<<

Rappel: Given a simplicial set \mX, $TX$ is its translate (cf. p. \pageref{18.17}).

[Note: \ $T_0X = X_1$, so there is a simplicial map $\si X_1 \ra TX$.  
On the other hand, the 
$d_0:X_{n+1} \ra X_n$ define a simplicial map $TX \ra X$.]

Example: If \bC is a simplicial object in \bCAT, then $T\bC \leftrightarrow (TM,TO)$, 
where $\bC \leftrightarrow (M,O)$ 
(an internal category in \bSISET) and there is a sequence $\si \bC_1 \ra T\bC \ra \bC$.

[Note: \ This applies to $\bw\bS\bC$, where \bC is a small Waldhausen category.  
Since $\bw\bS_1\bC$ is isomorphic to $\bw\bC$, there is a sequence 
$\si\bw\bC \ra T\bw\bS\bC \ra \bw\bS\bC$ and since $B\bw\bS_0\bC = *$, $BT\bw\bS\bC$ is 
contractible (cf. p. \pageref{18.18}).  
Thus one is led again to the arrow 
$B\bw\bC \ra \Omega B \bw\bS\bC$ whose adjoint $\Sigma B\bw\bC \ra B\bw\bS\bC$ is the closed embedding on 
p. \pageref{18.19}.  
By naturality, \bC can be replaced by $\bS\bC$, which produces another sequence 
$\si \bw\bS\bC \ra T\bw\bS^{(2)}\bC \ra \bw\bS^{(2)}\bC$.  It follows from Proposition 8 below that the sequence 
$B\bw\bS\bC \ra BT \bw\bS^{(2)}\bC$ $\ra B \bw\bS^{(2)}\bC$ of classifying spaces is a fibration up to homotopy 
(per \bCGH (singular structure)).  Therefore the arrow 
$B \bw\bS\bC \ra \Omega B\bw\bS^{(2)}\bC$ is a weak homotopy equivalence or still, a pointed homotopy equivalence.  
Continuing, one sees that 
\label{18.4}
$B\bw\bS^{(q)}\bC \ra \Omega B\bw\bS^{(q+1)}\bC$ is a pointed homotopy equivalence $\forall \ q$ 
(cf. p. \pageref{18.20}).]\\

Let 
$
\begin{cases}
\ \bC\\
\ \bD
\end{cases}
$
be small Waldhausen categories, $F:\bC \ra \bD$ a model functor.  Define

%%----------------------------------------------------------------------------------------------18
$\bS(\bC \overset{F}{\ra} \bD)$ by the pullback square
\begin{tikzcd}%[sep=small]
{\bS(\bC \overset{F}{\ra} \bD)} \ar{d} \ar{r} &{T\bS\bD}\ar{d}\\
{\bS\bC} \ar{r} &{\bS\bD}
\end{tikzcd}
in $[\Delta^\OP,\bCAT]$, so $\forall\  n$, 
\begin{tikzcd}%[sep=small]
{\bS_n(\bC \overset{F}{\ra} \bD)} \ar{d} \ar{r} &{\bS_{n+1}\bD}\ar{d}\\
{\bS_n\bC} \ar{r}  &{\bS_n\bD}
\end{tikzcd}
is a pullback square in \bCAT.

[Note: \ There is a sequence $\si\bD \ra \bS(\bC \overset{F}{\ra} \bD) \ra \bS\bC$.]\\

\textbf{\small LEMMA} \  
$\bS_n(\bC \overset{F}{\ra} \bD)$ is a small Waldhausen category.

[The weak equivalences are given by the pullback square 
\begin{tikzcd}%[sep=small]
{\bw\bS_n(\bC \overset{F}{\ra} \bD)} \ar{d} \ar{r} &{\bw\bS_{n+1}\bD} \ar{d}\\
{\bw\bS_n\bC} \ar{r} &{\bw\bS_n\bD}
\end{tikzcd}
and the cofibrations are given by the pullback square
\begin{tikzcd}%[sep=small]
{\bco\bS_n(\bC \overset{F}{\ra} \bD)} \ar{d} \ar{r} &{\bco\bS_{n+1}\bD} \ar{d}\\
{\bco\bS_n\bC} \ar{r} &{\bco\bS_n\bD}
\end{tikzcd}
.]

[Note: \ $\bS(\bC \overset{F}{\ra} \bD)$  is a simplicial object in \bWALD.]\\

\begingroup%%----------------------------------->>
\fontsize{9pt}{11pt}\selectfont
\textbf{\small EXAMPLE} \  
Taking $\bC = \bD$ and $F = \id_{\bC}$ gives nothing new
\begin{tikzcd}%[sep=small]
{(\bS(\bC} \ar{r}{\id_{\bC}} &{\bC) = T\bS\bC)}
\end{tikzcd}
but there is a variant which is of some interest.  Thus define $\bG\bC$ by the pullback square
\begin{tikzcd}%[sep=small]
{\bG\bC} \ar{d} \ar{r}&{T\bS\bC} \ar{d}\\
{T\bS\bC} \ar{r} &{\bS\bC}
\end{tikzcd}
$-$then $\bG_n\bC$ is a small Waldhausen category and $\bG\bC$ is a simplicial object in \bWALD.  The significance of 
$\bG\bC$ lies in the fact that the arrow 
$B\bw\bG\bC \ra \Omega B\bw\bS\bC$ is a weak homotopy equivalence if \bC is a category WES 
(Gillet-Grayson\footnote[2]{\textit{Illinois J. Math.} \textbf{31} (1987), 574-597; 
see also Gunnarsson et al., \textit{J. Pure Appl. Algebra} \textbf{79} (1992), 255-270.}).\\
\endgroup %%------------------------------------<<

\begin{proposition} \ %08
Let
$
\begin{cases}
\ \bC\\
\ \bD
\end{cases}
$
be small Waldhausen categories, $F:\bC \ra \bD$ a model functor, $-$then the sequence 
$B\bw\bS\bD \ra B\bw\bS^{(2)}(\bC \overset{F}{\ra} \bD) \ra B\bw\bS^{(2)}\bC$ of classifying spaces is a fibration up to homotopy (per \bCGH (singular structure)).
\end{proposition}

[It suffices to verify that $\forall \ n$, the sequence 
$B\bw\bS\bD \ra B\bw\bS\bS_n(\bC \overset{F}{\ra} \bD) \ra B\bw\bS\bS_n\bC$ is a fibration up to homotopy (per \bCGH (singular structure)) (cf. p. \pageref{18.21})
($\pi_0(B\bw\bS\bS_n\bC) = *$ $\forall \ n$).  
Do this by comparing it with the sequence 
$B\bw\bS\bD \ra B\bw\bS\bD \times_k B\bw\bS\bS_n\bC$
$\ra B\bw\bS\bS_n\bC$, using the triad lemma to establish that the arrow 
$B\bw\bS\bD \times_k B\bw\bS\bS_n\bC \ra B\bw\bS\bS_n(\bC \overset{F}{\ra} \bD)$ is a 
``retraction up to homotopy''.]\\

%%----------------------------------------------------------------------------------------------19
\textbf{\small LEMMA} \  
Equip \bCGH with its singular structure.  Suppose given a commutative diagram 
\ 
\begin{tikzcd}%[sep=large]
{A} \ar{d} \ar{r}{g} &{X} \ar{d}{\phi} \ar{r}{f}&{Y}\arrow[d,shift right=0.5,dash] \arrow[d,shift right=-0.5,dash]\\
{A^\prime} \ar{r}[swap]{g^\prime} &{X^\prime} \ar{r}[swap]{f^\prime} &{Y}
\end{tikzcd}
\ 
of pointed compactly generated Hausdorff spaces.  Assume: The rows are fibrations up to homotopy $-$then the square 
\ 
\begin{tikzcd}%[sep=large]
{A} \ar{d} \ar{r}{g} &{X} \ar{d}{\phi}\\
{A^\prime} \ar{r}[swap]{g^\prime} &{X^\prime}
\end{tikzcd}
\ 
is a homotopy pullback.

[The claim is that the arrow $A \ra W_{g^\prime,\phi}$ is a weak homotopy equivalence.  Consider the commutative diagram 
\begin{tikzcd}%[sep=small]
{A^\prime} \ar{d}\ar{r}{g^\prime} 
&{X^\prime} \arrow[d,shift right=0.5,dash] \arrow[d,shift right=-0.5,dash] 
&{X} \arrow[d,shift right=0.5,dash] \arrow[d,shift right=-0.5,dash] \ar{l}[swap]{\phi}\\
{E_{f^\prime}} \ar{r}[swap]{\pi^\prime}
&{X^\prime}
&{X}\ar{l}{\phi}
\end{tikzcd}
.  
By hypothesis, the arrow $A^\prime \ra E_{f^\prime}$, is a weak homotopy equivalence, so the induced map 
$W_{g^\prime,\phi}\ra W_{\pi^\prime,\phi}$ is a weak homotopy equivalence (cf. p. \pageref{18.22}).  
On the other hand, the projection 
$\pi^\prime:E_{f^\prime} \ra X^\prime$ is a pointed \bCG fibration 
(cf. p. \pageref{18.23}), hence is a \bCG fibration 
(cf. p. \pageref{18.24}).  Therefore the arrow 
$E_{f^\prime} \times_{X^\prime} X \ra W_{\pi^\prime,\phi}$
is a homotopy equivalence. (cf. $\S 4$, Proposition 18).  But 
$E_{f^\prime} \times_{X^\prime} X = \{y_0\} \times_Y W_{f^\prime} \times_{X^\prime} X$ 
$= \{y_0\} \times_Y W_f = E_f$ and by hypothesis, the arrow $A \ra E_f$ is a weak homotopy equivalence.]\\

\begin{proposition} \ %09
Let $\bC^\prime$, \bC, $\bC\pp$ be small Waldhausen categories.  Suppose given model functors 
$\bC^\prime \ra \bC$, 
$\bC \ra \bC\pp$ $-$then the square
\begin{tikzcd}%[sep=small]
{B\bw\bS\bC} \ar{d} \ar{r} &{B\bw\bS^{(2)}(\bC^\prime \ra \bC)} \ar{d}\\
{B\bw\bS\bC\pp} \ar{r} &{B\bw\bS^{(2)}(\bC^\prime \ra \bC\pp)}
\end{tikzcd}
is a homotopy pullback (per \bCGH (singular structure)).
\end{proposition}

[Bearing in mind Proposition 8, apply the lemma to the commutative diagram 

\[
\begin{tikzcd}%[sep=small]
{B\bw\bS\bC} \ar{d} \ar{r} 
&{B\bw\bS^{(2)}(\bC^\prime \ra \bC)} \ar{d} \ar{r}
&{B\bw\bS^{(2)}\bC^\prime}\arrow[d,shift right=0.5,dash] \arrow[d,shift right=-0.5,dash]\\
{B\bw\bS\bC\pp} \ar{r} 
&{B\bw\bS^{(2)}(\bC^\prime \ra \bC\pp)} \ar{r}
&{B\bw\bS^{(2)}\bC^\prime}
\end{tikzcd}
.]
\]
\vspace{0.15cm}

\begingroup%%----------------------------------->>
\fontsize{9pt}{11pt}\selectfont
Suppose given a small category \bC carrying the structure of two Waldhausen categories, both having the same subcategory of cofibrations but potentially distinct subcategories of weak equivalences, say $\bv\bC$ and $\bw\bC$, with 
$\bv\bC \subset \bw\bC$ (e.g., $\bv\bC$ might be $\iso\bC$).  Let $\bC^\bw$ be the full subcategory of \bC whose objects are the \mX such that $0 \ra X$ is in $\bw\bC$, put 
$\bv\bC^\bw = \bv\bC \cap \bC^\bw$ $\&$ 
$\bw\bC^\bw = \bw\bC \cap \bC^\bw$, and 
$\bco\bC^\bw = \bco\bC \cap \bC^\bw$ $-$then $\bC^\bw$ is Waldhausen relative to either notion of weak equivalence.\\
\endgroup %%------------------------------------<<

\index{Theorem: Localization Theorem}
\index{Localization Theorem}
\textbf{\small LOCALIZATION THEOREM} \quad 
\begingroup%%----------------------------------->>
\fontsize{9pt}{11pt}\selectfont
Assume that \bC admits a functor $M:\bC(\ra) \ra \bC$ that is
%%----------------------------------------------------------------------------------------------20
a mapping cylinder in the v-structure and the w-structure.  Suppose further that in the $w$-structure, the saturation axiom, the extension axiom, and the mapping cylinder axiom all hold $-$then the square
\begin{tikzcd}[sep=large]
{B\bv\bS\bC^\bw} \ar{d} \ar{r} &{B\bw\bS\bC^\bw} \ar{d} \\
{B\bv\bS\bC}  \ar{r} &{B\bw\bS\bC}
\end{tikzcd}
is a homotopy pullback (per \bCGH (singular structure)).
\vspi
[The proof, which depends on Proposition 9, is detailed in 
Waldhausen\footnote[2]{\textit{SLN} \textbf{1126} (1985), 350-352.}.]
\vspi
[Note: \ $\forall \ n$, $\bw\bS_n\bC^\bw$ has an initial object, thus $B\bw\bS\bC^\bw$ is contractible.]\\
\endgroup %%------------------------------------<<

\begingroup%%----------------------------------->>
\fontsize{9pt}{11pt}\selectfont
Remark: Proposition 3 enters into the proof through the assumption that the w-structure on \bC satisfies the saturation axiom and the mapping cylinder axiom.  
As for the role of the extension axiom, recall that if $X \ra Y$ is an acyclic cofibration, then $0 \ra Y/X$ is an acyclic cofibration (cf. Proposition 2), i.e., $Y/X \in \Ob\bC^\bw$.  
Conversely, if $X \ra Y$ is a cofibration for which $Y/X \in \Ob\bC^\bw$, then the extension axiom implies that $X \ra Y$ is a weak equivalence (consider the commutative diagram 
$
\begin{tikzcd}[sep=large]
{X}
\arrow[d,shift right=0.5,dash] \arrow[d,shift right=-0.5,dash]
\arrow[r,shift right=0.5,dash] \arrow[r,shift right=-0.5,dash]   
&{X}\ar{d} \ar{r} &{0}\ar{d}\\
{X} \ar{r} &{Y} \ar{r} &{Y/X}
\end{tikzcd}
).
$
\vspi
[Note: \ For an interesting application of the localization theorem to the algebraic K-theory of a ring with unit, see 
Weibel-Yao\footnote[3]{\textit{Contemp. Math.} \textbf{126} (1992), 219-230.}
.]\\
\endgroup %%------------------------------------<<

\begin{proposition} \ %10
Let
$
\begin{cases}
\ \bC\\
\ \bD
\end{cases}
$
be small Waldhausen categories, $F:\bC \ra \bD$ a model functor, $-$then there is a long exact sequence 
$\cdots \ra$ 
$\pi_{n+1}(B\bw\bS^{(2)}(\bC \overset{F}{\ra} \bD)) \ra$ 
$\pi_n(B\bw\bS\bC) \ra $
$\pi_n(B\bw\bS\bD) \ra $
$\pi_n(B\bw\bS^{(2)}(\bC \overset{F}{\ra} \bD)) \ra $
$\cdots \ra$
$\pi_{2}(B\bw\bS^{(2)}(\bC \overset{F}{\ra} \bD)) \ra$ 
$\pi_1(B\bw\bS\bC) \ra $
$\pi_1(B\bw\bS\bD) \ra $
$\pi_1(B\bw\bS^{(2)}(\bC \overset{F}{\ra} \bD)) \ra $
$\pi_0(B\bw\bS\bC) \ra $
$\pi_0(B\bw\bS\bD)$ 
in homotopy.
\end{proposition}

[Proposition 9 implies that the square
\begin{tikzcd}%[sep=small]
{B\bw\bS\bC} \ar{d} \ar{r}&{B\bw\bS^{(2)}(\bC \overset{\id_C}{\ra} \bC)} \ar{d}\\
{B\bw\bS\bD} \ar{r} &{B\bw\bS^{(2)}(\bC \overset{F}{\ra} \bD)}
\end{tikzcd}
is a homotopy pullback (per \bCGH (singular structure)), thus the Mayer-Vietoris sequence is applicable 
(cf. p. \pageref{18.25}).  And: $B\bw\bS^{(2)}(\bC \overset{\id_{\bC}}{\lra} \bC)$ is contractible.]\\

\index{Cofinality Principle}
\begingroup%%----------------------------------->>
\fontsize{9pt}{11pt}\selectfont
\textbf{\small COFINALITY PRINCIPLE} \quad 
Let \bC, \bD be small categories WES.  Assume \bC is cofinal in \bD $-$then $K_0(\bC)$ is a subgroup of $K_0(\bD)$ 
(cf. p. \pageref{18.26}) and $\forall \ n \geq 1$, $K_n(\bC) \approx K_n(\bD)$.
\vspi
[Since by defnintion, $K_n(\bC) \approx \pi_{n+1}(B\bw\bS\bC)$ $\&$ $K_n(\bD) \approx \pi_{n+1}(B\bw\bS\bD)$, 
one can invoke Proposition 10 if the higher homotopy groups of $B\bw\bS^{(2)}(\bC \overset{\iota}{\ra} \bD)$ are trivial. 
This is established by showing
%%----------------------------------------------------------------------------------------------21
that $B\bw\bS^{(2)}(\bC \overset{\iota}{\ra} \bD)$  has the same pointed homotopy type as 
$B(K_0(\bD)/K_0(\bC))$, the classifying space of $K_0(\bD)/K_0(\bC)$.]
\vspi
[Note: \ All the particulars can be found in 
Staffeldt\footnote[2]{\textit{K-theory} \textbf{1} (1989), 511-532; 
see also Grayson, \textit{Illinois J. Math.} \textbf{31} (1987), 598-617.}
.]\\
\endgroup %%------------------------------------<<

\begingroup%%----------------------------------->>
\fontsize{9pt}{11pt}\selectfont
\textbf{\small EXAMPLE} \  
Let \bC be a small category WES $-$then \bC is cofinal in $\bC_{\pa}$ (cf. p. \pageref{18.27}), hence $\forall \ n \geq 1$,  $K_n(\bC) \approx K_n(\bC_{\pa})$.
\vspi
[Note: \ Let \mA be a ring with unit $-$then $\bF(A)$ is cofinal in $\bP(A)$, so the higher algebraic K-groups of 
$\bF(A)$ can be identified with the  higher algebraic K-groups of $\bP(A)$.]\\
\endgroup %%------------------------------------<<

Let \bC, \bD be small Waldhausen categories, $F:\bC \ra \bD$ a model functor $-$then \mF is said to have the 
\un{approximation property}
\index{approximation property (functor of small Waldhausen categories)} 
provided that the following conditions are satisfied.

\indent\indent (App$_1$) \ A morphism $f$ in \bC is in $\bw\bC$ if $F f$ is in $\bw\bD$.

\indent\indent (App$_2$) \ Given $X \in \Ob\bC$ and $f \in \Mor(FX,Y)$, there is a $g \in \Mor(X,X^\prime)$ and a weak equivalence $h:FX^\prime \ra Y$ such that $f = h \circx F g$: 
\begin{tikzcd}%[sep=large]
{FX} \ar{d}[swap]{Fg} \ar{r}{f} &{Y}\\
{FX^\prime} \ar{ru}[swap]{h}
\end{tikzcd}
.\\

Remarks: 
(1) Since \mF is a model functor, $F f$ is in $\bw\bD$ if $f$ is in $\bw\bC$; 
(2) When \bC satisfies the mapping cylincer axiom, $\exists$ a commutative triangle 
\begin{tikzcd}%[sep=large]
{X} \ar{rd}[swap]{g} \ar{r}{i} &{M_g} \ar{d}{r}\\
&{X^\prime}
\end{tikzcd}
, where $r$ is a weak equivalence, hence in this case one can assume that the ``$g$'' is a cofibration.\\

%\begin{doublespacing}
\index{Theorem: Approximation Theorem}
\index{Approximation Theorem}
\textbf{\small APPROXIMATION THEOREM} \quad 
Let \bC, \bD be small Waldhausen categories satisfying the saturation axiom, $F:\bC \ra \bD$ a model functor.  Suppose that \bC 
satisfies the mapping cylinder axiom and \mF has the approximation property $-$then 
$B\bw\bS F:B\bw\bS\bC \ra B\bw\bS\bD$ is a pointed homotopy equivalence.

[This result is due to 
Waldhausen\footnote[3]{\textit{SLN} \textbf{1126} (1985), 354-358.}.  
I shall omit the proof (which is long and technical) but by way of simplification, it suffices that 
$B\bw F:B\bw\bC \ra B\bw\bD$ be a pointed homotopy equivalence.  Reason: $\bS_n\bC$ and $\bS_n\bD$ inherit the assumptions made on \bC and \bD, thus $\forall \ n$, 
$B\bw\bS_n F:B\bw\bS_n\bC \ra B\bw\bS_n\bD$ is a pointed homotopy equivalence and so 
$B\bw\bS F:B\bw\bS\bC \ra B\bw\bS\bD$ is a pointed homotopy equivalence (cf. p. \pageref{18.29}).  One then proceeds to the crux, viz. the verification that 
$\bw F:\bw\bC \ra \bw\bD$
is a strictly initial functor, and concludes by appealing to Quillen's theorem A.]\\
%\end{doublespacing}

%%----------------------------------------------------------------------------------------------22
\begingroup%%----------------------------------->>
\fontsize{9pt}{11pt}\selectfont
\textbf{\small EXAMPLE} \  
Let \bC be the Waldhausen category whose objects are the pointed finite CW complexes and whose morphisms are the pointed skeletal maps.  
Let \bD be the category whose objects are the wellpointed spaces with closed base point which have the pointed homotopy type of a pointed finite CW complex and whose morphisms are the pointed continuous functions $-$then \bD satisfies the axioms for a  Waldhausen category if weak equivalence  = weak homotopy equivalence, cofibration = closed cofibration.  
However, while \bC is skeletally small, \bD is definitely not.  
Still, it will be convenient to ignore this detail since the situation can be rectified by the insertion of some additional language.  We claim that the inclusion $\iota:\bC \ra \bD$ has the approximation 
property.   
App$_1$ is, of course trivial.  
To check the validity of App$_2$, fix a \mK in \bC and suppose given a pointed continuous function $f:K \ra X$, where \mX is in \bD.  By definition, $\exists$ an \mL in \bC and pointed continuous functions 
$\phi:X \ra L$, $\psi:L \ra X$ such that $\psi \circx \phi \simeq \id_X$, $\phi \circx \psi \simeq \id_L$.  
Using the skeletal approximation theorem, choose a pointed skeletal $g:K \ra L$ for which $g \simeq \phi \circx f$.  
Display the data in a commutative diagram 
\begin{tikzcd}[sep=large]
{K} \ar{dr}[swap]{g}\ar{r}{i} &{M_g} \ar{d}{r} &{L} \ar{l}[swap]{j} \ar[equals]{dl}\\
&{L}
\end{tikzcd}
and consider the composite 
$M_g \overset{r}{\ra} L \overset{\psi}{\ra} X$.  
Since 
$\psi \circx r \circx i = \psi \circx g$, 
$\psi \circx r \circx j = \psi$, 
the restriction of 
$\psi \circx r$ to $K \vee L$ equals $\psi \circx g \vee \psi$ (identify \mK $\&$ $i(K)$, \mL $\&$ $j(L)$).  
But 
$g \simeq \phi \circx f$ $\implies$ $\psi \circx g \simeq \psi \circx \phi \circx f \simeq f$ $\implies$ 
$\psi \circx g \vee \psi \simeq f \vee \psi$.  
Because $K \vee L \ra M_g$ is a closed cofibration, it follows that $f \vee \psi$ 
admits an extension to $M_g$, call it h:
$
\begin{tikzcd}[sep=large]
{K} \ar{dr}[swap]{f}\ar{r}{i} &{M_g} \ar{d}{h} &{L} \ar{l}[swap]{j}\ar{dl}{\psi}\\
&{X}
\end{tikzcd}
. 
$
From the triangle on the right, one sees that $h$ is a weak homotopy equivalence.  
On the other hand, $f = h \circx i$ and $i$ is skeletal.\\
\endgroup %%------------------------------------<<

\begingroup%%----------------------------------->>
\fontsize{9pt}{11pt}\selectfont
\textbf{\small EXAMPLE} \  
Let \bC be the Waldhausen category whose objects are the pointed finite simplicial sets with weak equivalence  = weak homotopy equivalence, cofibration = pointed injective simplicial map and let \bD be as in the preceding example.  
We claim that the geometric realization $\abs{?}:\bC \ra \bD$ has the approximation property.  App$_1$ is is true by definition.  
Turning to App$_2$, fix an \mX in \bC and suppose given a pointed continuous function $f:\abs{X} \ra Y$, where \mY is in \bD.  
Let us assume for the moment that it is possible to fulfill App$_2$ up to homotopy, i.e., $\exists$ a pointed finite simplicial set 
$X^\prime$, a simplicial map $g:X \ra X^\prime$, and a weak homotopy equivalence $h:\abs{X^\prime} \ra Y$ such that 
$f \simeq h \circx \abs{g}$ $-$then App$_2$ holds on the nose.  Indeed, $\abs{M_g} \approx M_{\abs{g}}$ and there is a commutative diagram 
\begin{tikzcd}[sep=large]
{\abs{X}}\ar{dr}[swap]{\abs{g}} \ar{r}{\abs{i}}
&{M_{\abs{g}}} \ar{d}{\abs{r}} 
&{\abs{X^\prime}} \ar{l}[swap]{\abs{j}} \ar[equals]{dl}\\
&{\abs{X^\prime}}
\end{tikzcd}
.  Obviously, 
$h \circx \abs{r} \circx \abs{i} = h \circx \abs{g}$, 
$h \circx \abs{r} \circx \abs{j} = h$, and 
$h \circx \abs{g} \vee h \simeq f \vee h$, hence $f \vee h$ can be extended to $\approx M_{\abs{g}}$ call it 
$
H:
\begin{tikzcd}[sep=large]
{\abs{X}} \ar{dr}[swap]{f}\ar{r}{\abs{i}} 
&{M_{\abs{g}}} \ar{d}{H} 
&{\abs{X^\prime}} \ar{l}[swap]{\abs{j}} \ar{dl}{h}\\
&{Y}
\end{tikzcd}
.
$
But \mH is a weak homotopy equivalence and $f = H \circx \abs{i}$, as desired.  Proceedin$g$, there exists a pointed CW complex having the pointed homotopy type of \mY and without loss of
%%----------------------------------------------------------------------------------------------23
generality, one can assume that it is the geometric realization of a pointed finite simplicial set \mK (cf. $\S 5$, Proposition 3 and use the barycentric subdivision of the relevant vertex scheme), thus \mY may be replaced by $\abs{K}$.  
Because \mX is finite, the argument employed in the proof of the simplicial approximation theorem produces a simplicial map 
$g:X \ra \text{Ex}^n K$ $(\exists \ n)$  for which $\abs{g} \simeq \abs{e_K^n} \circx f$.  And: 
$\abs{e_K^n}:\abs{K} \ra \abs{\text{Ex}^n K}$ is a pointed homotopy equivalence 
(cf. p. \pageref{18.30}).\\
\endgroup %%------------------------------------<<

\begingroup%%----------------------------------->>
\fontsize{9pt}{11pt}\selectfont
Remark: The above considerations therefore imply that the algebraic K-theory of a point can also be defined in terms of pointed finite simplicial sets.\\
\endgroup %%------------------------------------<<

Let \mA be a ring with unit $-$then it is clear that $K_0(\bP(A)) = K_0(A)$.\\

\index{Consistency  Principle}
\textbf{\small CONSISTENCY  PRINCIPLE} \quad 
There is a pointed homotopy equivalence 
$\Omega B\bw\sS\ov{\bP(A)} \ra K_0(A) \times B\bGL(A)^+$, 
hence $\forall \ n \geq 1$, $K_n(\bP(A)) \approx K_n(A)$.

[Note: Recall that $K_n(A) = \pi_n(B\bGL(A)^+)$ (cf. p. \pageref{18.31} ff.).]\\

This is not obvious and the existing proofs are quite roundabout in that they do not directly invole 
$B\bw\sS\ov{\bP(A)}$.  Instead, one replaces it with $B\bQ\ov{\bP(A)}$, where $\bQ\ov{\bP(A)}$ is the 
``\bQ construction'' on $\ov{\bP(A)}$ (cf. infra), and then introduces yet another artifice, namely the 
``$\bS^{-1}\bS$ construction'' which, in effect, is a bridge between these two very different ways of defining the higher algebraic K-groups of \mA.  For the ``classical'' approach to these matters, consult the seventh chapter of 
Srinivas\footnote[2]{\textit{Algebraic K-Theory}, Birkh\"auser (1991); 
see also Gillet-Grayson, \textit{Illinois J. Math.} \textbf{31} (1987), 574-597 (cf. 591-593).}
(a sophisticated variant has been given by 
Jardine\footnote[3]{\textit{J. Pure Appl. Algebra} \textbf{75} (1991), 103-194; 
see also Thomason, \textit{Comm. Algebra} \textbf{10} (1982), 1589-1668.}).

Example: Form the monoid $\coprod\limits_P B \Aut P$, where \mP runs through the objects in $\ov{\bP(A)}$ $-$then 
in the pointed homotopy category, 
$\Omega B \coprod\limits_P B \Aut P \approx K_0(A) \times B\bGL(A)^+$ 
(cf. p. \pageref{18.32} ff.).\\

\begingroup%%----------------------------------->>
\fontsize{9pt}{11pt}\selectfont
Let \bC be a small category WES $-$then $\bQ\bC$ is the category with the same objects as \bC, a morphism from \mX to \mY in $\bQ\bC$ being an equivalence class of diagrams of the form
$X \twoheadleftarrow A \rightarrowtail Y$, where
$X \twoheadleftarrow A^\prime \rightarrowtail Y$ $\&$ 
$X \twoheadleftarrow A\pp \rightarrowtail Y$ 
are equivalent if $\exists$ an isomorphism $A^\prime \ra A\pp$ rendering 
\ 
\begin{tikzcd}[sep=large]
{X} \arrow[d,shift right=0.5,dash] \arrow[d,shift right=-0.5,dash]
&{A^\prime} \ar{d} \ar[l,two heads]  \ar[r,tail]
&{Y} \arrow[d,shift right=0.5,dash] \arrow[d,shift right=-0.5,dash]\\
{X} 
&{A\pp} \ar[l,two heads]  \ar[r,tail]
&{Y}
\end{tikzcd}
\ 
commutative.  To compose 
$X \twoheadleftarrow A \rightarrowtail Y$ and 
$Y \twoheadleftarrow B \rightarrowtail Z$,
form the pullback $A \times_Y B$ and project to \mX and 
%%----------------------------------------------------------------------------------------------24
Z, i.e., 
\begin{tikzcd}%[sep=small]
{A \times_Y B}\ar[d,two heads]\ar[r,tail]
&{B}\ar[r,tail]\ar[d,two heads]
&{Z}\\
{A} \ar[d,two heads]\ar[r,tail]
&{Y}\\
{X}
\end{tikzcd}
.\\[.5cm]
\endgroup %%------------------------------------<<

\begingroup%%----------------------------------->>
\fontsize{9pt}{11pt}\selectfont
Observation: If \bC, \bD are small categories WES and if $F:\bC \ra \bD$ is an exact functor, then there is an induced functor $\bQ F:\bQ\bC \ra \bQ\bD$.\\
\endgroup %%------------------------------------<<

\index{Proposition W}
\textbf{\small PROPOSITION W} \quad 
\begingroup%%----------------------------------->>
\fontsize{9pt}{11pt}\selectfont
Let \bC be a small category WES $-$then $B\bw\bS\bC$ and $B\bQ\bC$ have the same pointed homotopy type.\\
\endgroup %%------------------------------------<<

\begingroup%%----------------------------------->>
\fontsize{9pt}{11pt}\selectfont
The proof of Proposition W depends on an auxiliary device.
Let $\sd:\bDelta \ra \bDelta$ be the functor that sends $[n]$ to $[2n+1]$ and $\alpha:[m] \ra [n]$ to the arrow 
$[2m+1] \ra [2n+1]$ defined by the prescription 
$0 \ra \alpha(0), \ldots, m \ra \alpha(m)$, 
$m + 1 \ra 2n+1 - \alpha(m)$, $\ldots, 2m+1 \ra$ $2n+1 - \alpha(0)$.\\
\vspi
Given a simplicial space \mX, put $\sd X = X \circx \sd^\OP$, the 
\un{edgewise subdivision}
\index{edgewise subdivision} 
of \mX.  
So, $(\sd X)_n = X_{2n+1}$ and the 
$
\begin{cases}
\ d_i\\
\ s_i
\end{cases}
$
per $\sd X$ are the 
$
\begin{cases}
\ d_i \circx d_{2n+1-i} \ (0 \leq i \leq n, n > 0)\\
\ s_i \circx s_{2n+1-i} \ (0 \leq i \leq n, n \geq 0)
\end{cases}
$
per \mX.\\
\endgroup %%------------------------------------<<
\vspace{0.2cm}

\textbf{\small LEMMA} \  
\begingroup%%----------------------------------->>
\fontsize{9pt}{11pt}\selectfont
Specify a continuous function 
$\theta_n:(\sd X)_n \times \Delta^n \ra X_{2n+1} \times \Delta^{2n+1}$ via the formula 
$\theta_n(x,t_0, \ldots, t_n) = (x, \frac{1}{2}t_0,\ldots, \frac{1}{2}t_n, \frac{1}{2}t_n, \ldots, \frac{1}{2}t_0)$ $-$then 
the $\theta_n$ induce a homeomorphism $\abs{\sd X} \ra \abs{X}$.\\
\endgroup %%------------------------------------<<

\begingroup%%----------------------------------->>
\fontsize{9pt}{11pt}\selectfont
Let \bC be a small category WES $-$then the weak equivalences are isomorphisms (cf. Proposition 4), hence 
$B\bw\bS\bC = B\iso\bS\bC$ and there is a pointed homotopy equivalence 
$\abs{W\bC} \ra B\iso\bS\bC$ 
(cf. p. \pageref{18.33}).  
On the other hand, from the lemma, $\abs{\sd W\bC} \approx \abs{W\bC}$, thus to prove 
Proposition W, it suffices to construct a pointed homotopy equivalence $\abs{\sd W\bC} \ra B\bQ\bC$.  
An element \mF of $(\sd W\bC)_n$ is an element of $W_{2n+1}\bC = \Ob\bS_{2n+1}\bC$.  Writing 
$F_{i,j}$ for $F(i \ra j)$, send \mF to that element of $\nersub_n\bQ\bC$ represented by the diagram
\begin{tikzcd}[sep=small]
&{F_{n-1,n+1}}\ar[lddd,two heads]\ar[rddd,tail]\\
\\
\\
{F_{n,n+1}}&&{F_{n-1,n+2}}
\end{tikzcd}
$\cdots$
\begin{tikzcd}[sep=small]
&{F_{0,2n}}\ar[lddd,two heads]\ar[drdd,tail]\\
\\
\\
{F_{1,2n}}&&{F_{0,2n+1}}
\end{tikzcd}
, i.e., to the string 
$F_{n,n+1} \ra $
$F_{n-1,n+2} \ra $
$\cdots \ra $
$F_{1,2n} \ra $
$F_{0,2n+1}$
in $\nersub_n\bQ\bC$.  
This assignment defines a simplicial map $\sd W\bC \ra \ner \bQ\bC$ and the claim is that its geometric 
realization is a pointed homotopy equivalence.
\vspi
Introduce the double category $i\bQ\bC \equiv \iso \bQ\bC \cdot \bQ\bC$ and recall that there is a pointed homotopy 
equivalence $B\bQ\bC \ra Bi\bQ\bC$ (cf. p. \pageref{18.34}).  
Call $i\bQ_n\bC$ the category whose objects are the functors $[n] \ra \bQ\bC$ and whose morphisms are the natural isomorphisms  
($\implies$ $i\bQ_n\bC = \iso[[n],\bQ\bC]$) $-$then $\forall \ n$, the functor 
$\iso\sd\bS_n\bC \ra i\bQ_n\bC$ is an equivalence of categories.  
Contemplation of the diagram 
\begin{tikzcd}%[sep=small]
{\abs{\sd W\bC}} \ar{d}\ar{r}&{B\bQ\bC} \ar{d}\\
{B\iso\sd\bS\bC} \ar{r} &{Bi\bQ\bC}
\end{tikzcd}
finishes the argument.\\
\vspace{0.25cm}
\endgroup %%------------------------------------<<

%%----------------------------------------------------------------------------------------------25
Let \mA be a ring with unit $-$then by definition, $\bW A$ is the $\Omega$-prespectrum with $q^\text{th}$ space 
$K_0(\Sigma^q A) \times B\bGL(\Sigma^q A)^+$ (cf. p. \pageref{18.35}) and $\bK A = eM\bW A$ 
(cf. p. \pageref{18.36}), thus $\pi_n(\bK A) = K_n(A)$ $(n \geq 0)$.  And: 
$\pi_{-n}(\bK A) = K_0(\Sigma^n A) = (L^nK_0)(A)$ $(n \geq 0)$, the negative algebraic K-groups of \mA in the sense of 
Bass (compare, e.g., Karoubi\footnote[2]{\textit{Ann. Sci. \'Ecole Norm. Sup.} \textbf{4} (1971), 63-95.}
).

[Note: \ The $\pi_{-n}(\bK A)$ vanish if \mA is left  noetherian and every finitely generated left $A$-module has finite projective dimension.]\\

\begingroup%%----------------------------------->>
\fontsize{9pt}{11pt}\selectfont
The consistency principle can be generalized: $\exists$ a morphism of spectra $\bK\bP(A) \ra \bK(A)$ such that the induced map 
$\pi_n(\bK\bP(A)) \ra \pi_n(\bK A)$ is an isomorphism $\forall \ n \geq 0$.\\
\endgroup %%------------------------------------<<

To conclude this $\S$, I shall say a few words about topological K-theory.

[Note: \ A reference is the book of 
Karoubi\footnote[2]{\textit{K-Theory: An Introduction}, Springer Verlag (1978); 
see also N. Wegge-Olsen, \textit{K-Theory and $C^*$-Algebras}, Oxford University Press (1993).}.]

Let \mA be a Banach algebra with unit over \bk, where $\bk = \R$ or $\C$.  Write $\bGL(A)^\text{top}$ for $\bGL(A)$ in its canonical topology $-$then $\bGL(A)^\text{top}$ is a topological group and $\pi_0(\bGL(A)^\text{top})$ is abelian.  
Definition: $\forall \ n > 0$, $K_n^\text{top}(A) = \pi_n(B\bGL(A)^\text{top})$, the $n^\text{th}$ topological K-group of \mA 
(put $K_0^\text{top}(A) = K_0(A)$).\\

\index{Theorem: Bott Periodicity Theorem}
\index{Bott Periodicity Theorem}
\textbf{\small BOTT PERIODICITY THEOREM} \quad 
Let \mA be a Banach algebra with unit over \bk.\\  
$
\indent\indent
\begin{cases}
\ (\bk = \C) \ \forall \ n \geq 0, \ K_n^\text{top}(A) \approx K_{n+2}^\text{top}(A)\\
\ (\bk = \R) \ \forall \ n \geq 0, \ K_n^\text{top}(A) \approx K_{n+8}^\text{top}(A)
\end{cases}
$
.\\[.5cm]

\begingroup%%----------------------------------->>
\fontsize{9pt}{11pt}\selectfont
For instance, one can take for \mA the Banach algebra with unit whose elements are the real or complex valued continuous functions on a compact Hausdorff space \mX.\\
\endgroup %%------------------------------------<<

The identity 
$\bGL(A) \ra \bGL(A)^\text{top}$ induces a map 
$B\bGL(A) \ra B\bGL(A)^\text{top}$, 
from which an arrow 
$B\bGL(A)^+ \ra B\bGL(A)^\text{top}$.  
Passing to homotopy, this gives a homomorphism 
$K_n(A) \ra K_n^\text{top}(A)$ that connects the algebraic K-groups of \mA to the topological K-groups of \mA.

[Note: \ \ The fundamental \ group \ of\ $B\bGL(A)^\text{top}$ is abelian 
($\pi_1(B\bGL(A)^\text{top})$ $\approx$
$\pi_0(\Omega B\bGL(A)^\text{top})$ $\approx$ 
$\pi_0(\bGL(A)^\text{top})$, thus $B\bGL(A)^\text{top}$ is insensitive to the plus construction.]\\

%%----------------------------------------------------------------------------------------------26
\index{Theorem: Theorem of Fischer-Prasolov}
\index{Theorem of Fischer-Prasolov}
\textbf{\small THEOREM OF FISCHER}\footnote[3]{\textit{J. Pure Appl. Algebra} \textbf{69} (1990), 33-50.}
\textbf{-PRASOLOV}\footnote[6]{\textit{Amer. Math. Soc. Transl.} \textbf{154} (1992), 133-137.} 
Let \mA be a commutative Banach algebra over \bk with unit $-$then $\forall \ n \geq 1$, the arrow 
\[
\pi_n(B\bGL(A)^+;\Z/k\Z) \ra \pi_n(B\bGL(A)^\text{top};\Z/k\Z)
\]
is an isomorphism.

[Note: \ The notation is that of p. \pageref{18.37} 
($B\bGL(A)^+$ and $B\bGL(A)^\text{top}$ are H spaces).]\\

Therefore, in the commutative case, the algebraic and topologial K-groups of \mA are indistinguishable if one sticks to finite coefficients.\\

%%%%%%%%%%%%%%%%%%%%%%%%%%%%%%%%%%%%%%
%%%%%%%%%%%%%%%%%%%%%%%%%%%%%%%%%%%%%%
%%%%%%%%%%%%%%%%%%%%%%%%%%%%%%%%%%%%%%

\begin{center}
$\S \ 18$
\\[0.5cm]
$\mathcal{REFERENCES}$\\
\end{center}

\[
\textbf{BOOKS}
\]

\begingroup
\fontsize{9pt}{11pt}\selectfont
\setlength\parindent{0 cm}

[1] \quad Inassaridze, H., \textit{Algebraic K-Theory}, Kluwer (1995).
\\[-.2cm]

[2] \quad Jardine, J., \textit{Generalized \'Etale Cohomology Theories}, Birkh\"auser (1997).
\\[-.2cm]

[3] \quad Lluis-Puebla, E. et al., \textit{Higher Algebraic K-Theory: An Overview}, Springer Verlag (1992).
\\[-.2cm]

[4] \quad Loday, J-L., \textit{Cyclic Homology}, Springer Verlag (1992).
\\[-.2cm]

[5] \quad Paluch, M., \textit{Algebraic and Topological K-Theory}, Ph.D. Thesis, University of Illinois at 

\hspace{0.8cm}Chicago, Chicago (1991).
\\[-.2cm]

[6] \quad Rosenberg, J., \textit{Algebraic K-Theory and its Applications}, Springer Verlag (1994).
\\[-.2cm]

[7] \quad Srinivas, V., \textit{Algebraic K-Theory}, Birkh\"auser (1991).
\\[-.2cm]
\endgroup

\[
\textbf{ARTICLES}
\]

\begingroup
\fontsize{9pt}{11pt}\selectfont
\setlength\parindent{0 cm}

[1] \quad Betley, S. and Pirashvili, T., Stable K-Theory as a Derived Functor, 
\textit{J. Pure Appl. Algebra} \textbf{96} (1994), 

\hspace{0.8cm}245-258
.\\[-.2cm]

[2] \quad B\"okstedt, M., Hsiang, W., and Madsen, I., The Cyclotomic Trace and Algebraic K-Theory of Spaces,

\hspace{0.8cm}\textit{Invent. Math.} \textbf{111} (1993), 465-539.
\\[-.2cm]

[3] \quad Dundas, B. and McCarthy, R., Stable K-Theory and Topological Hochschild Homology, 
\textit{Ann. of Math.} 

\hspace{0.8cm}\textbf{140} (1994), 685-701 and 142 (1995), 425-426.
\\[-.2cm]

[4] \quad Dundas, B. and McCarthy, R., Topological Hochschild Homology of Ring Functors and Exact Cate-

\hspace{0.8cm}gories, 
\textit{J. Pure Appl. Algebra} \textbf{109} (1996), 231-294.
\\[-.2cm]

[5] \quad Dwyer, W. and Mitchell, S., On the K-Theory Spectrum of a Smooth Curve over a Finite Field,

\hspace{0.8cm}\textit{Topology} \textbf{36} (1997), 899-929.
\\[-.2cm]

[6] \quad Fiedorowicz, Z., A Note on the Spectra of Algebraic K-Theory, 
\textit{Topology} \textbf{16} (1977), 417-421.
\\[-.2cm]

[7] \quad Gersten, S., On the Spectrum of Algebraic K-Theory, 
\textit{Bull. Amer. Math. Soc.} \textbf{78} (1972), 216-219.
\\[-.2cm]

[8] \quad Gersten, S., Higher K-Theory of Rings, 
\textit{SLN} \textbf{341} (1973), 1-40.
\\[-.2cm]

[9] \quad Giffen, C., Loop Spaces for the \bQ-Construction, 
\textit{J. Pure Appl. Algebra} \textbf{52} (1988), 1-30.
\\[-.2cm]

[10] \quad Gillet, H., Riemann Roch Theorems for Higher Algebraic K-Theory, 
\textit{Adv. Math.} \textbf{40} (1981), 203-289.
\\[-.2cm]

[11] \quad Grayson, D., Higher Algebraic K-Theory: II, 
\textit{SLN} \textbf{551} (1976), 217-240.
\\[-.2cm]

[12] \quad Grayson, D., On the K-Theory of Fields, 
\textit{Contemp. Math.} \textbf{83} (1989), 31-55.
\\[-.2cm]

[13] \quad Gunnarson, T., Algebraic K-Theory of Spaces as K-Theory of Monads, 
\textit{Preprint Series} \textbf{21} (1981/82), 

\hspace{0.95cm}Aarhus Universitet.
\\[-.2cm]

[14] \quad Higson, N., Algebraic K-Theory of Stable $C^*$-Algebras, 
\textit{Adv. Math.} \textbf{67} (1988), 1-140.
\\[-.2cm]

[15] \quad Hiller, H., $\lambda$-Rings and Algebraic K-Theory, 
\textit{J. Pure Appl. Algebra} \textbf{20} (1981), 241-266.
\\[-.2cm]

[16] \quad Hinich, V. and Schechtman, V., Geometry of a Category of Complexes and Algebraic K-Theory,

\hspace{0.95cm}\textit{Duke Math. J.} \textbf{52} (1985), 399-430.
\\[-.2cm]

[17] \quad Iqusa, K., The Stability Theorem for Smooth Pseudoisotopies, 
\textit{K-Theory} \textbf{2} (1988), 1-355.
\\[-.2cm]

[18] \quad Jardine, J., The Lichtenbaum-Quillen Conjecture for Fields, 
\textit{Canad. Math. Bull.} \textbf{36} (1993), 426-441.
\\[-.2cm]

[19] \quad Jardine, J., The K-Theory of Finite Fields, Revisited, 
\textit{K-Theory} \textbf{7} (1993), 579-595.
\\[-.2cm]

[20] \quad Kahn, B., Bott Elements in Algebraic K-Theory, 
\textit{Topology} \textbf{36} (1997), 963-1006.
\\[-.2cm]

[21] \quad Kassel, C., La K-Th\'eorie Stable, 
\textit{Bull. Soc. Math. France} \textbf{110} (1982), 381-416.
\\[-.2cm]

%[22] \quad Keller, B., On the Cyclic Hommology of Exact Categories,\\[-.2cm]
[22] \quad Keller, B., On the Cyclic Homology of Exact Categories, 
\textit{J. Pure and Applied Alg.} \textbf{136} (1999) p. 1-56.
\\[-.2cm]

[23] \quad Klein, J. and Rognes, J., The Fiber of the Linearization Map $A(*) \ra K(\Z)$, 
\textit{Topology} \textbf{36} (1997), 

\hspace{0.95cm}829-848.
\\[-.2cm]

[24] \quad Landsburg, S., K-Theory and Patching for Categories of Complexes, 
\textit{Duke Math. J.} \textbf{62} (1991), 359-

\hspace{0.95cm}384.
\\[-.2cm]

[25] \quad Loday, J-L., K-Th\'eorie Alg\'ebrique et Repr\'esentations de Groupes, 
\textit{Ann. Sci. \'Ecole Norm. Sup.} 

\hspace{0.95cm}\textbf{9} (1976), 309-377.
\\[-.2cm]

[26] \quad Lydakis, M., Fixed Point Problems, Equivariant Stable Homotopy, and a Trace Map for the Algebraic 

\hspace{0.95cm}K-Theory of a Point, 
\textit{Topology} \textbf{34} (1995), 959-999.
\\[-.2cm]

[27] \quad Madsen, I., Algebraic K-Theory and Traces, In: 
\textit{Current Developments in Mathematics}, 1995, R. 

\hspace{0.95cm}Bott et al. (ed.), International Press (1994), 191-321.
\\[-.2cm]

[28] \quad McCarthy, R., On Fundamental Theorems of Algebraic K-Theory, 
\textit{Topology} \textbf{32} (1993), 325-328.
\\[-.2cm]

[29] \quad McCarthy, R., The Cyclic Homology of an Exact Category, 
\textit{J. Pure Appl. Algebra} \textbf{93} (1994), 251-296.
\\[-.2cm]

[30] \quad Mitchell, S., On the Lichtenbaum-Quillen Conjectures from a Stable Homotopy Theoretic Viewpoint, 

\hspace{0.95cm}In: 
\textit{Algebraic Topology and its Applications}, G. Carlsson et al. (ed.), Springer Verlag (1994), 163-240.
\\[-.2cm]

[31] \quad Mitchell, S., Hypercohomology Spectra and Thomason's Descent Theorem, In: 
\textit{Algebraic K-Theory}, 

\hspace{0.95cm}V. Snaith (ed.), American Mathematical Society (1997), 221-277.
\\[-.2cm]

[32] \quad Nenashev, A., Comparison Theorems for $\lambda$-Operations in Higher Algebraic K-Theory, 
\textit{Ast\'erisque} 

\hspace{0.95cm}\textbf{226} (1994), 335-369.
\\[-.2cm]

[33] \quad Quillen, D., Higher Algebraic K-Theory: I, 
\textit{SLN} \textbf{341} (1973), 85-147.
\\[-.2cm]

[34] \quad Rognes, J., A Spectrum Level Rank Filtration in Algebraic K-Theory, 
\textit{Topology} \textbf{31} (1992), 813-845.
\\[-.2cm]

[35] \quad Staffeldt, R., On Fundamental Theorems of Algebraic K-Theory, 
\textit{K-Theory} \textbf{1} (1989), 511-532.
\\[-.2cm]

[36] \quad Steiner, R., Infinite Loop Structures on the Algebraic K-Theory of Spaces, 
\textit{Math. Proc. Cambridge}

\hspace{0.95cm}\textit{Philos. Soc.} \textbf{90} (1981), 85-111.
\\[-.2cm]

[37] \quad Suslin, A. and Wodzicki, M., Excision in Algebraic K-Theory, 
\textit{Ann. of Math.} \textbf{136} (1992), 51-122.
\\[-.2cm]

[38] \quad Swan, R., Higher Algebraic K-Theory, 
\textit{Proc. Sympos. Pure Math.} \textbf{58} (1995) (Part 1), 247-293.
\\[-.2cm]

[39] \quad Thomason, R., Algebraic K-Theory and \'Etale Cohomology, 
\textit{Ann. of Sci. \'Ecole Norm.  Sup.} \textbf{18}

\hspace{0.95cm}(1985), 437-552 and \textbf{22} (1989), 675-677.
\\[-.2cm]

[40] \quad Thomason, R. and Trobaugh, T., Higher Algebraic K-Theory of Schemes and of Derived Categories, 

\hspace{0.95cm}In: 
\textit{The Grothendieck Festschrift, vol. III}, Birkh\"auser (1990), 247-435.
\\[-.2cm]

[41] \quad Wagoner, J., Delooping Classifying Spaces in Algebraic K-Theory, 
\textit{Topology} \textbf{11} (1972), 349-370.
\\[-.2cm]

[42] \quad Waldhausen, F., Algebraic K-Theory of Generalized Free Products, 
\textit{Ann. of Math.} \textbf{108} (1978), 

\hspace{0.95cm}135-256.
\\[-.2cm]

[43] \quad Waldhausen, F., Algebraic K-Theory of Spaces, 
\textit{SLN} \textbf{1126} (1985), 318-419.
\\[-.2cm]

[44] \quad Weibel, C., A Survey of Products in Algebraic K-Theory, 
\textit{SLN} \textbf{854} (1981), 494-517.
\\[-.2cm]

[45] \quad Yao, D., Higher Algebraic K-Theory of Admissible Abelian Categories and Localization Theorems,

\hspace{0.95cm}\textit{J. Pure Appl. Algebra} \textbf{77} (1992), 263-339.

\setlength\parindent{2em}

\endgroup

\chapter{
$\boldsymbol{\S}$\textbf{19}.\quadx  DIMENSION THEORY}
\setlength\parindent{2em}
\setcounter{proposition}{0}
\setcounter{chapter}{19}

%%----------------------------------------------------------------------------------------------01
$\text{ }$\\[-1.25cm]

Dimension theory enables one to associate with each nonempty normal Hausdorff space \mX a topological invariant 
$\dim X \in \{0, 1, \ldots\} \cup \{\infty\}$ called its topological dimension.  Classically, there are two central theorems, namely:\\
\indent\indent (1) \ The topological dimension of $\R^n$ is exactly $n$, hence as a corollary, $\R^n$, and $\R^m$ are homeomorphic iff $n = m$.\\
\indent\indent (2) \ Every second countable normal Hausdorff space of topological dimension $n$ can be embedded in 
$\R^{2n+1}$.

Although I shall limit the general discussion to what is needed to prove these results, some important applications will be given, e.g., to the ``invariance of domain'' and the ``superposition question''.  On the theoretical side, \u Cech cohomology makes an initial apperance but it does not really come to the fore until $\S 20$.\

Let \mX be a nonempty normal Hausdorff space.  
Consider the following statement.

\indent\indent $(\dim X \leq n)$ \ There exists an integer $n = 0 , 1, \ldots$ such that every finite open covering 
of \mX has a finite open refinement of order $\leq n+1$.

If $\dim X \leq n$ is true for some $n$, then the 
\un{topological dimension}
\index{topological dimension} 
of \mX, 
denoted by $\dim X$, is the smallest value of $n$ for which $\dim X \leq n$.

[Note: \ By convention, $\dim X = -1$ when $X = \emptyset$.  If the statement $\dim X \leq n$ is false for every $n$, then we put $\dim X = \infty$.]

Our primary emphasis will be on spaces of finite topological dimension.  
A simple example of a compact metrizable space of infinite topological dimension is the 
\un{Hilbert cube}\index{Hilbert cube} 
$[0,1]^\omega$.\\

\begingroup%%----------------------------------->>
\fontsize{9pt}{11pt}\selectfont
Why work with finite open coverings?  Answer: The concept of dimension would be very different otherwise.  
Example: Take $X = [0,\Omega[$ $-$then $\dim [0,\Omega[ \ = 0$ (cf. p. \pageref{19.1}).  But the open covering 
$\{[0,\alpha[\ :0 < \alpha < \Omega\}$ has no point finite open refinement, so $[0,\Omega[$ 
would be ``infinite dimensional'' if arbitrary open coverings were allowed.
\vspi
Why work with normal \mX? A priori, this is not necessary since the definition evidently makes sense for any CRH space \mX.  But observe: If $\dim X = 0$, then \mX must be normal.  So, no new spaces of ``dimension zero'' are produced by just formally extending the definition to nonnormal \mX.  Such an agreement would also introduce a degree of pathology.  
Example: The topological dimension of $X = [0,\Omega] \times [0,\omega]$ is zero (cf. p. \pageref{19.2}) but the ``topological dimension'' of $X - \{(\Omega,\omega)\}$, the Tychonoff plank (which is not normal), is one.  The escape from this predicament is to reformulate the definition of dim in such a way that it is naturally applicable to the class of all nonempty CRH spaces.  The topological dimension of the Tychonoff plank then turns out to be zero, as might be expected (cf. p. \pageref{19.3}).\\
\vspi
%%----------------------------------------------------------------------------------------------02
\label{19.23}
Let \mX be a nonempty CRH space.  Consider the following statement.
\\
\indent\indent $(\dim X \leq n)$  There exists an integer $n = 0, 1, \ldots$ such that every finite numerable open covering of \mX has a finite numerable open refinement of order $\leq n+1$.
\vspi
If $\dim X \leq n$ is true for some $n$, then the 
\un{topological dimension}
\index{topological dimension (CRH space)} 
of 
\mX, denoted by $\dim X$, is the smallest value of $n$ for which $\dim X \leq n$.
\vspi
[Note: By convention, $\dim X = -1$ iff $X = \emptyset$.  If the statement $\dim X \leq n$ is false for every $n$, then we put $\dim X = \infty$.]\\
\endgroup %%------------------------------------<<

\label{19.36}
\begingroup%%----------------------------------->>
\fontsize{9pt}{11pt}\selectfont
Since a nonempty CRH space \mX is normal iff every finite open covering of \mX is numerable, this agreement is a consistent extension of dim.  On the other hand, the price to pay for increasing the generality is that more things can go wrong (e.g., every subspace of \mX now has a topological dimension).  
Because of this, my policy will be to concentrate on the normal case and simply indicate as we go along what changes, if any, must be made to accommodate the completely regular situation.  
The omitted details are invariably straightforward.
\vspi
\label{19.13}
[Note: \ By repeating what has been said above verbatim, an arbitratry nonempty topological space \mX acquires a ``topological dimension'' $\dim X$.  One can then show that $\dim X = \dim \crg X$, where $\crg X$ is the complete regularization of \mX (cf. p. \pageref{19.4}).  
Example: $\dim[0,1]/[0,1[ \ = \ 0$.]\\
\endgroup %%------------------------------------<<

\begin{proposition} \ %01
The topological dimension of \mX is equal to the topological dimension of $\beta X$.
\end{proposition}

[$\dim \beta X \leq n$ $\implies$ $\dim X \leq n$ : Let $\sU = \{U\}$ be a finite open covering of \mX.  Since $\sU$ is numerable, one can assume that the \mU are cozero sets.  The collection $\{\beta X - \cl_{\beta X}(X - U)\}$ is then a finite open covering of $\beta X$, thus admits a precise open refinement of order $\leq n+1$ which, when restricted to \mX, is a precise open refinement of $\sU$ of order $\leq n+1$.

$\dim X \leq n$ $\implies$ $\dim \beta X \leq n$ : Let $\sU = \{U\}$ be a finite open covering of $\beta X$.  
Choose a partition of unity $\{\kappa_U\}$ on $\beta X$ subordinate to $\sU$.  
The collection $\{X \cap \kappa_U^{-1}(]0,1])\}$ is a finite open covering of \mX, 
hence has a precise open refinement $\sV = \{V\}$ of order $\leq n+1$.  
Let $\{\kappa_V\}$ be a partition of unity on \mX subordinate to $\sV$ $-$then the collection 
$\{\beta X - \cl_{\beta X}(X - \kappa_V^{-1}(]0,1]))\}$ is a precise open refinement of $\sU$ of order $\leq n + 1$.]\\

\begingroup%%----------------------------------->>
\fontsize{9pt}{11pt}\selectfont
The argument used in Proposition 1 carries over directly to the completely regular situation, so the result holds in that setting too.\\
\endgroup %%------------------------------------<<

A nonempty Hausdorff space is said to be 
\un{zero dimensional}
\index{zero dimensional} 
if it has a basis consisting of clopen sets.  Every zero dimensional space is necessarily completely regular.  The class of zero dimensional spaces is closed under the formation of nonempty products and coproducts.

%%----------------------------------------------------------------------------------------------03
[Note: \ Recall that for nonempty LCH spaces, the notions of zero dimensional and totally disconnected are equivalent.]\\

\begingroup%%----------------------------------->>
\fontsize{9pt}{11pt}\selectfont
A nonempty subspace of the real line is zero dimensional iff it contains no open interval.
\vspi
The Isbell$-$Mr\'owka space, the van Douwen line, the van Douwen space, and the Kunen line are all zero dimensional, But of these, only the Kunen line is normal.\\
\endgroup %%------------------------------------<<

\begingroup%%----------------------------------->>
\fontsize{9pt}{11pt}\selectfont
\textbf{\small FACT} \  
Let \mX be a zero dimensional normal LCH space.  Suppose that \mX is metacompact $-$then \mX is subparacompact.\\
\endgroup %%------------------------------------<<

\begingroup%%----------------------------------->>
\fontsize{9pt}{11pt}\selectfont
Any metric space $(X,d)$ for which $d(x,z) \leq \max(d(x,y),d(y,z))$ is zero dimensional.  Such a metric space is said to be \un{nonarchimedean }.
\index{nonarchimedean  (metric space)}    
They are common fare in algebraic number theory and $p$-adic analysis.  
Example: Suppose that \mX is zero dimensional and second countable $-$then \mX admits a compatible nonarchimedean  metric.  Indeed, let $\sU = \{U_n\}$ be a clopen basis for \mX and put 
$d(x,y) = \max\limits_n\left\{\ds\frac{\abs{\chi_n(x) - \chi_n(y)}}{n}\right\}$, $\chi_n$ the characteristic function on $U_n$.
\vspi
[Note: \ Suppose that \mX is metrizable $-$then 
de Groot\footnote[2]{\textit{Proc. Amer. Math. Soc.} \textbf{7} (1956), 948-953; 
see also Nagata, \textit{Fund. Math.} 55 (1964), 181-194.}
has shown that $\dim X = 0$ iff \mX admits a compatible nonarchimedean  metric.]\\
\endgroup %%------------------------------------<<

\label{19.21}
\begingroup%%----------------------------------->>
\fontsize{9pt}{11pt}\selectfont
\textbf{\small EXAMPLE} \  
Let $\kappa$ be an infinite cardinal $-$then the 
\un{Cantor cube}
\index{Cantor cube} 
$C_\kappa$ is the space 
$\{0,1\}^\kappa$, where $\{0,1\}$ has the discrete topology.  
It is a compact Hausdorff space of weight $\kappa$ and is zero dimensional.  
Of course, the Cantor cube associated with $\kappa = \omega$ is homeomorphic to the usual Cantor set.  Every zero dimensional space \mX of weight $\kappa$ can be embedded in $C_\kappa$, hence has a zero dimensional compactification 
$\zeta X$ of weight $\kappa$.
\vspi
[Let $\sU = \{U_i:i \in I\}$ be a clopen basis for \mX such that $\#(I) = \kappa$.  
Agreeing to denote by $\chi_i$ the characteristic function of $U_i$, call $\chi$ the diagonal of the $\chi_i$ $-$then $\chi:X \ra C_\kappa$ is an embedding.  
The closure  $\zeta X$ of the image of \mX in $C_\kappa$ is a zero dimensional compactification of \mX of weight $\kappa$.  Viewing \mX as a subspace of $\zeta X$, to within topological equivalence $\zeta X$ is the only zero dimensional compactification of \mX with the property: 
For every zero dimensional compact Hausdorff space \mY and for every continuous function 
$f:X \ra Y$ there exists a continuous function $\zeta f: \zeta X \ra Y$ such that $\zeta \restr{f}{X} = f$.]
\vspi
[Note: \ Consider the Cantor cube $C_\omega$.  Since $C_\omega \hookrightarrow \R$, it follows that if \mX is zero dimensional and second countable, then there is an embedding $X \ra \R$.]\\
\endgroup %%------------------------------------<<

Suppose that $\dim X = 0$ $-$then it is clear that \mX is zero dimensional.  To what extent is the converse true?\\

%%----------------------------------------------------------------------------------------------04
\textbf{\small LEMMA} \  
If for every pair $(A,B)$ of disjoint closed subsets of \mX there exists a clopen set $U \subset X$ such that 
$A \subset U \subset X - B$, then $\dim X = 0$.

[Let $\sU = \{U_i: i \in I\}$ be a finite open covering of \mX of cardinality $\#(I) = k$.  
To establish the existence of a finite refinement of $\sU$ by pairwise disjoint clopen sets, we shall argue by induction on $k$.  For $k = 1$, the assertion is trivial.  
Assume that $k > 1$ and that the assertion is true for all open coverings of cardinality $k - 1$.  Enumerate the elements of $\sU$: 
$U_1, \ldots, U_k$ and pass to $\{U_1, \ldots, U_{k-1} \cup U_k\}$, which thus has a precise refinement 
$\{V_1, \ldots, V_{k-1}\}$ by pairwise disjoint clopen sets.  Noting that 
$
\begin{cases}
\ V_{k-1} - U_{k-1}\\
\ V_{k-1} - U_k
\end{cases}
$
are disjoint closed subsets of \mX, choose a clopen set $U \subset X$: 
$
\begin{cases}
\ V_{k-1} - U_{k-1} \subset U\\
\ U \subset (X - V_{k-1}) \cup U_k
\end{cases}
. \ 
$
Consideration of the covering  $\{V_1, \ldots, V_{k-1} - U, V_{k-1} \cap U\}$ then finishes the induction.]\\

\begin{proposition} \ \hspace{-.15cm} % 02
Suppose that \mX is zero dimensional and Lindel\"of $-$then $\dim X = 0$.
\end{proposition}

[Let $(A,B)$ be a pair of disjoint closed subsets of \mX.  
Given $x \in X$, choose a clopen neighborhood $U_x \subset X$ of $x$ such that either 
$A \cap U_x = \emptyset$ or $B \cap U_x = \emptyset$.  Let $\{U_{x_i}\}$ be a countable subcover of 
$\{U_x\}$ $-$then the $U_i = U_{x_i} - \bigcup\limits_{j<i} U_{x_j}$ are pairwise disjoint clopen subsets of \mX and 
$\bigcup\limits_i U_i = X$.  
Put $U = \st(A,\{U_i\})$: \mU is clopen and $A \subset U \subset X - B$.  
The lemma therefore implies that $\dim X = 0$.]\\

\label{19.2}
Take $X = [0,\Omega]$ $-$then \mX is zero dimensional and compact, thus in view of Proposition 2, 
$\dim [0,\Omega]= 0$.  Take next $X = [0,\Omega[$ $-$then $\beta X = [0,\Omega]$, so $\dim [0,\Omega[  \ = 0$ too 
(cf. Proposition 1).\\

\begingroup%%----------------------------------->>
\fontsize{9pt}{11pt}\selectfont 
\textbf{\small LEMMA} \  
Let \mX be a nonempty CRH space $-$then $\dim X = 0$ iff for every pair of disjoint zero sets in \mX there exists a clopen set in \mX containing the one and not the other.\\
\endgroup %%------------------------------------<<

\label{19.1}
\label{19.3}
\begingroup%%----------------------------------->>
\fontsize{9pt}{11pt}\selectfont
Consequently, Proposition 2 is valid as it stands in the completely regular situtation.  
Example: Consider $[0,\Omega] \times [0,\omega]$ and conclude that the topological dimension of the Tychonoff plank is zero.\\
\endgroup %%------------------------------------<<

\begingroup%%----------------------------------->>
\fontsize{9pt}{11pt}\selectfont
\label{19.33}
\textbf{\small LEMMA} \  
Let \mX be a nonempty CRH space $-$then $\dim X = 0$ iff every zero set in \mX is a countable intersection of clopen sets.\\
\endgroup %%------------------------------------<<

\begingroup%%----------------------------------->>
\fontsize{9pt}{11pt}\selectfont
\textbf{\small EXAMPLE} \  
Let $\kappa$ be a cardinal $-$then $\N^\kappa$ is paracompact if $\kappa$ is countable but is neither normal nor submetacompact if $\kappa$ is uncountable.  Claim: $\forall \ \kappa$, $\dim \N^\kappa = 0$.  
For this, it can be assumed that $\kappa$ is uncountable.  
Let $Z(f)$ be a zero set in $\N^\kappa$ $-$then there exists a countable subproduct through which $f$ factors, 
i.e., there exists a continuous $g:\N^\omega \ra \R$ such that $f = g \circx p$, $p:\N^\kappa \ra \N^\omega$ the
%%----------------------------------------------------------------------------------------------05
projection.  
Obviously, $Z(g) = p(Z(f))$.  
Choose a sequence $\{V_n\}$ of clopen sets in $\N^\omega$: 
$Z(g) = \ds\bigcap\limits_n V_n$.  Put $U_n = p^{-1}(V_n)$ $-$then $U_n$ is clopen in $\N^\kappa$ and 
$Z(f) = \ds\bigcap\limits_n U_n$.
\vspi
[Note: \ Suppose that $\kappa$ is uncountable $-$then every open subspace of $\N^\kappa$ has topological dimension zero but this need not be the case of closed subspaces (cf. p. \pageref{19.5}).]\\
\endgroup %%------------------------------------<<

\begingroup%%----------------------------------->>
\fontsize{9pt}{11pt}\selectfont
\textbf{\small FACT} \  
Let \mX be a nonempty CRH space $-$then $\dim X = 0$ iff the real valued continuous functions on \mX with finite range are uniformly dense in $BC(X)$.
\vspi
[There is no loss of generality in assuming that \mX is compact.  If \mX is totally disconnected, use Stone-Weierstrass; If \mX is not totally disconnected, consider the functions constant on some connected subset of \mX that has more than one point.]\\
\endgroup %%------------------------------------<<

It is false that unconditionally:  $X$ zero dimensional $\implies$ $\dim X = 0$, even if \mX is a metric space 
(Roy\footnote[2]{\textit{Trans. Amer. Math. Soc.} \textbf{134} (1968), 117-132; 
see also Kulesza, \textit{Topology Appl.} \textbf{35} (1990), 109-120.}).

[Note: \ The topological dimension of Roy's metric space is equal to 1.  
Does there exist for each $n > 1$ a zero dimensional metric space \mX such that $\dim X = n$?  The answer is unknown.]\\

\begingroup%%----------------------------------->>
\fontsize{9pt}{11pt}\selectfont
\index{Dowker's Example ``M''}
\textbf{\small EXAMPLE \  (\un{Dowker's Example ``M''})} \  
In $[0,1]$, write $x \sim y$ iff $x - y \in \Q$, so $[0,1]/\sim \ = \ds\coprod\limits_\alpha \Q_\alpha$.  
There are $2^\omega$ equivalence classes $\Q_\alpha$.  
Each is a countable dense subset of $[0,1]$.  
Take a subcollection 
$\{\Q_\alpha: \alpha< \Omega\}$, where $\forall \ \alpha < \Omega$ : $\Q_\alpha \neq \Q \cap [0,1]$.  
Put 
$S_\alpha = [0,1] - \ds\bigcup \{\Q_\beta: \alpha \leq \beta < \Omega\}$ and consider the subspace 
$X = \{(\alpha,s):\alpha < \Omega, s \in S_\alpha\}$ of $[0,\Omega[ \ \times [0,1]$ $-$then \mX is zero dimensional and the claim is that \mX is normal, yet $\dim X > 0$.  
To see this, form $X^* = X \cup (\{\Omega\} \times [0,1])$, a subspace of 
$[0,\Omega] \times [0,1]$ which is normal.  
In addition, if \mA and \mB are disjoint closed subsets of \mX, then their closures 
$A^*$ and $B^*$ in $X^*$ are also disjoint.  It follows that \mX is normal.  
If $\dim X = 0$, then there exists a clopen set 
$U \subset X$ such that $[0,\Omega[ \ \times \{0\} \subset U$ and $[0,\Omega[ \ \times \{1\} \subset X - U$.  But 
$U^* \cap (X - U)^* = \emptyset$ $\&$ 
$
\begin{cases}
\ (\Omega,0) \in U^*\\
\ (\Omega,1) \in (X - U)^*
\end{cases}
$ 
, and this contradicts the connectedness of $\{\Omega\} \times [0,1]$.  
Therefore $\dim X > 0$.  One can be precise: 
$\dim X = 1$.  For if $\{U\}$ is a finite open covering of \mX, then $\forall \ t \in [0,1]$, there exists a neighborhood \mO of t and an $\alpha$ such that $X \cap (]\alpha,\Omega[ \ \times O)$ is contained in some \mU, which implies that there exists a finite open covering $\{O\}$ of $[0,1]$ of order $\leq 2$ and an $\alpha$ such that 
$X  \cap  (]\alpha,\Omega[ \ \times O)$ is contained in some \mU.  
Therefore $\dim X \leq 1$.
\vspi
[Note: \ \mX has a zero dimensional compactification $\zeta X$ and the latter has topological dimension zero (cf. Proposition 2).  
So: A compact Hausdorff space of zero topological dimension can have a normal subspace of positive topological dimension.  Another aspect is that while \mX is zero dimensional, $\beta X$ is not.  
In fact, $\dim X = \dim \beta X$ (cf. Proposition 1), which is $> 0$, thus Proposition 2 is applicable.  
Here is a final 
%%----------------------------------------------------------------------------------------------06
remark: By appropriately adjoining to \mX a single point, one can destroy its zero dimensionality or reduce its topological dimension to zero without, in either case, losing normality.]\\
\endgroup %%------------------------------------<<

\begingroup%%----------------------------------->>
\fontsize{9pt}{11pt}\selectfont
Modify the preceding construction, replacing 
$
\begin{cases}
\ [0,1] \text{ by } [0,1]^\omega\\
\ S_\alpha \text{ by } S_\alpha^\omega
\end{cases}
$ 
and conclude that there exists a compact Hausdorff space of zero topological dimension with a normal subspace of infinite topological dimension.\\
\endgroup %%------------------------------------<<

\begingroup%%----------------------------------->>
\fontsize{9pt}{11pt}\selectfont
\textbf{\small FACT} \  
Suppose that $\dim X = 0$ and \mX is paracompact.  Let \mA be a closed subset of \mX; let \mY be a complete metric space 
$-$then every (bounded) continuous function $f:A \ra Y$ has a (bounded) continuous extension $F:X \ra Y$.
\vspi
[For $n = 1, 2, \ldots$, let $\sV_n$ be the covering of \mY by open 1/n balls.  Let $\sA_n = \{A_{i,n}:i \in I_n\}$ be an open partition of $\sA$ that refines $f^{-1}(\sV_n)$.  
Inductively determine an open partition 
$\sU_n = \{U_{i,n}: i \in I_n\}$ of \mX that refines $\sU_{n-1}$ and $\forall \ i \in I_n$: $A \cap U_{i,n} = A_{i,n}$.  
Assign to a given $x \in X$ an index $i(x,n) \in I_n$: $x \in U_{i(x,n),n}$.  
Choose points $y_{i,n} \in f(A_{i,n})$.  
Observe that 
$\{y_{i(x,n),n}\}$ is Cauchy.  Put $F(x) = \lim y_{i(x,n),n}$.]\\
\endgroup %%------------------------------------<<

\begingroup%%----------------------------------->>
\fontsize{9pt}{11pt}\selectfont
Provided that \mY is a separable complete metric space, the preceding result retains its validity if only $\dim X = 0$ and \mX is normal.\\
\endgroup %%------------------------------------<<

\begin{proposition} \ % 03
Suppose that \mX is a nonempty paracompact LCH space $-$then \mX is zero dimensional iff $\dim X = 0$.
\end{proposition}

[Since \mX is paracompact, \mX admits a representation $X = \coprod\limits_i X_i$, where the $X_i$ are nonempty pairwise disjoint open $\sigma$-compact (= Lindel\"of) subspaces of \mX (cf. p. \pageref{19.6}).  
But obviously, \mX is zero dimensional iff each of the $X_i$ is zero dimensional.  Now use Proposition 2.]\\

\begingroup%%----------------------------------->>
\fontsize{9pt}{11pt}\selectfont
Proposition 3 can fail for an arbitrary normal LCH space.  
Consider the space \mX of Dowker's Example ``M''.  
It is not locally compact.  
To get around this, let $p:X \ra [0,\Omega[$ be the projection, form $\beta p:\beta X \ra \beta [0,\Omega[$ 
$= [0,\Omega]$ and put $A = (\beta p)^{-1}([0,\Omega[)$.  
One can check that \mA is normal and zero dimensional.  And: 
$X \subset A \subset \beta X$ $\implies$ $\beta A = \beta X$ $\implies$ $\dim A = \dim X > 0$ (cf. Proposition 1).  
But \mA, being open in $\beta A$, is a LCH space.
\vspi
[Note: \ \mA zero dimensional $\implies$ $A_\infty$ zero dimensional $\implies$ $\dim(A_\infty) = 0$ (cf. Proposition 2).  So: A compact Hausdorff space of zero topological dimension can have an open subspace of positive 
topological dimension.]\\
\endgroup %%------------------------------------<<

\begingroup%%----------------------------------->>
\fontsize{9pt}{11pt}\selectfont
Let \mX be a CRH space.  
Suppose that $\sA$ is a collection of subsets of \mX closed under the formation of finite unions and finite intersections.  A subcollection $\sF \subset \sA$ is said to be an 
\label{2.19} 
\un{$\sA$-filter}
\index{A-filter} if
(i) $\emptyset \notin \sF$, 
(ii) $A \in \sF$ $\&$  $A \subset B \in \sA$ $\implies$ $B \in \sF$, and 
(iii) $\forall \ A, B \in \sF$: $A \cap B \in \sF$.  Example: $\sA =$ all zero sets in \mX or
%%----------------------------------------------------------------------------------------------07
$\sA =$ all clopen sets in \mX, the associated $\sA$-filters then being the zero set filters and the clopen set filters, respectively.\\
\label{2.20} 
\label{2.21}
\indent\indent (Fil$_1$) An $\sA$-filter $\sF$ is said to be an 
\un{$\sA$-ultrafilter}
\index{A-ultrafilter}
if $\sF$ is a maximal $\sA$-filter.  The maximality of $\sF$ is equivalent to the condition: If $B \in \sA$ and if 
$A \cap B \neq \emptyset$ $\forall \ A \in \sF$, then $B \in \sF$.  An $\sA$-ultrafilter $\sF$ is 
\un{prime}
\index{prime ($\sA$-ultrafilter)},  
i.e., if \mA and \mB belong to $\sA$ and if $A \cup B \in \sF$, then $A \in \sF$ or $B \in \sF$.  Every $\sA$-filter is contained in an $\sA$-ultrafilter.
\\
\indent\indent (Fil$_2$) An $\sA$-filter $\sF$ is said to be 
\un{fixed}
\index{fixed ($\sA$-filter)} if $\cap \sF$ is nonempty.
\\
\indent\indent (Fil$_3$) An $\sA$-filter $\sF$ is said to have the 
\un{countable intersection property}
\index{countable intersection property ($\sA$-filter)} 
if for every sequence $\{A_n\} \subset \sF$, $\ds\bigcap\limits_n A_n \neq \emptyset$.
\vspi
[Note: \ The zero sets in \mX are closed under the formation of countable intersections.  Therefore every zero set ultrafilter on \mX with the countable intersection property is closed under the formation of countable intersections.]
\vspi
The following standard characterizations illustrate the terminology.
\\
\indent\indent ($\R$) Let \mX be a CRH space $-$then \mX is $\R$-compact iff every zero set ultrafilter on \mX with the countable intersection property is fixed.
\\
\indent\indent ($\N$) Let \mX be a CRH space.  Suppose that \mX is zero dimensional $-$then \mX is $\N$-compact iff every clopen set ultrafilter on \mX with the countable intersection property is fixed.\\
\endgroup %%------------------------------------<<

\begingroup%%----------------------------------->>
\fontsize{9pt}{11pt}\selectfont
\textbf{\small LEMMA} \  
Let \mX be a nonempty CRH space.  Suppose that $\dim X = 0$ and \mX is $\R$-compact $-$then \mX is $\N$-compact.
\vspi
[Let $\sU$ be a clopen set ultrafilter on \mX with the countable intersection property $-$then the claim is that $\sU$ is fixed.  Choose a zero set ultrafilter $\sZ$ on \mX: $\sZ \supset \sU$.  
Take any sequence $\{Z_n\} \subset \sZ$ and write 
$Z_n = \ds\bigcap\limits_m U_{mn}$, $U_{mn}$ clopen.  Each  $U_{mn}$ meets every element of $\sU$, thus each $U_{mn}$ is in $\sU$.  But $\sU$ has the countable intersection property, so 
$\ds\bigcap\limits_n Z_n = \ds\bigcap\limits_n \ds\bigcap\limits_m U_{mn} \neq \emptyset$.  Therefore $\sZ$ has the countable intersection property, hence is fixed, and this implies that $\sU$ is fixed as well.]\\ 
\endgroup %%------------------------------------<<

\begingroup%%----------------------------------->>
\fontsize{9pt}{11pt}\selectfont
The converse to this lemma is false: There exist $\N$-compact spaces of positive topological dimension.\\
 \endgroup %%------------------------------------<<

\begingroup%%----------------------------------->>
\fontsize{9pt}{11pt}\selectfont
\index{Mysior Space}
\textbf{\small EXAMPLE \ (\un{Mysior Space})} \  
Let \mX be the subspace of $\ell^2$ consisting of all sequences $\{x_n\}$, with $x_n$ rational $-$then \mX is the textbook example of a totally disconnected space that is not zero dimensional (Erd\"os).  
Fix a countable dense subset \mD of \mX.  
For each $S \subset D$ with $\#(\ov{S} \cap \ov{D - S}) = 2^\omega$, choose a point $x_S \in \ov{S} \cap \ov{D - S}$ subject to: $S^\prime \neq S\pp$ $\implies$ $x_{S^\prime} \neq x_{s\pp}$.  
In addition, given $x \in X - D$, let 
$\{s_k(x)\}$ be a sequence in \mD having limit $x$ such that if $x = x_S$ for some $S \subset D$, then both \mS and $D - S$ contain infinitely many terms of $\{s_k(x)\}$.  Topologize \mX as follows: Isolate the points of \mD and take for the basic neighborhoods of $x \in X - D$ the sets $K_k(x) = \{x\} \cup \{s_l(x): l \geq k\}$ $(k = 1, 2, \ldots)$.  
The resulting topology 
$\tau$ on \mX is finer than the metric topology.  
And the space $X_\tau$ thereby produced is a nonnormal zero dimensional LCH  
%%----------------------------------------------------------------------------------------------08
space possessing a basis comprised of countable clopen compact sets.  
To see that $X_\tau$ is $\N$-compact, let $\sU$ be a clopen set ultrafilter on $X_\tau$ with the countable intersection property.  The collection $\{U \in \sU: U \ \text{clopen in \mX}\}$ is a clopen set ultrafilter on \mX with the countable intersection property, hence there exists a point $x_0$ in its intersection 
(\mX is Lindel\"of).  
This $x_0$ is then the intersection of countably many elements of $\sU$, thus $\sU$ is fixed and so $X_\tau$ is 
$\N$-compact.  Still, $\dim X_\tau > 0$.  
Observe first that since \mD is dense in $X_\tau$, the frontier in \mX of any clopen subset of $X_\tau$  has cardinality $< 2^\omega$.  Consider the disjoint zero sets
$
\begin{cases}
\ Z_1 = \{x: \norm{x} \leq 1\\
\ Z_2 = \{x: \norm{x} \geq 2
\end{cases}
. \ 
$ 
Let \mU be a clopen subset of $X_\tau$: $Z_1 \subset U \subset X - Z_2$ $-$then its frontier in \mX necessarily has cardinality 
$2^\omega$.\\
\endgroup %%------------------------------------<<

\begingroup%%----------------------------------->>
\fontsize{9pt}{11pt}\selectfont
\textbf{\small FACT} \  
Let \mX be a nonempty CRH space $-$then \mX is $\N$-compact iff \mX is zero dimensional and there exists a closed embedding $X \ra \ds\prod (\N \times [0,1])$.\\
\endgroup %%------------------------------------<<

\begingroup%%----------------------------------->>
\fontsize{9pt}{11pt}\selectfont
There exist zero dimensional $\R$-compact normal LCH spaces that are not $\N$-compact.  Owing to the lemma, such a space must have positive topological dimension (cf. Proposition  3).\\
\endgroup %%------------------------------------<<

\begingroup%%----------------------------------->>
\fontsize{9pt}{11pt}\selectfont
\index{The Kunen Plane}
\textbf{\small EXAMPLE \ [Assume CH] \ (\un{The Kunen Plane})} \  
The construction of the Kunen line starting from $X = \R$ can be carried out with no change whatsoever starting instead with 
$X = \R^2$, the upshot being that the Kunen plane $X_\Omega$, a space with the same general topological properties as the Kunen line.  So: $X_\Omega$ is a zero dimensional perfectly normal LCH space that is not paracompact but is first countable, hereditarily separable, and collectionwise normal.  
The topology $\tau_\Omega$ on $X_\Omega$ is finer than the usual topology on $\R^2$.  
And, $\forall \ S \subset \R^2$: $\#(\cl_{\R^2}(S) - \cl_\Omega(S)) \leq \omega$.  
It follows from this that if \mA and \mB are disjoint closed subsets of $X_\Omega$, 
then $\#(\ov{A} \cap \ov{B}) \leq \omega$, the bar denoting closure in $\R^2$.
\vspi
Claim: $X_\Omega$ is $\R$-compact.
\vspi
[Let $\sZ_\Omega$ be a zero set ultrafilter on $X_\Omega$ with the countable intersection property.  
Let 
$\sZ \subset \sZ_\Omega$ be the subcollection consisting of the $\R^2$-closed elements of $\sZ_\Omega$.  
Fix a point 
$z_0 \in \cap \hspace{0.03cm} \sZ$ and choose a continuous function $\phi:\R^2 \ra [0,1]$ such that $\phi^{-1}(0) = \{x_0\}$.  The sets
$
\begin{cases}
\ \phi^{-1}([0,1/n])\\
\ \phi^{-1}([1/n,1])
\end{cases}
$ 
are zero sets in $\R^2$, hence are zero sets in $X_\Omega$.  Of course, 
$X_\Omega = \phi^{-1}([0,1/n]) \cup \phi^{-1}[1/n,1])$.  But obviously,
$\phi^{-1}([1/n,1]) \notin \sZ$, thus 
$\phi^{-1}([1/n,1]) \notin \sZ_\Omega$ and so 
$\phi^{-1}([0,1/n]) \in \sZ_\Omega$, $\sZ_\Omega$ being prime.  Consequently, 
$\{z_0\} = \ds\bigcap\limits_n \phi^{-1}([0,1/n])$ $\in \sZ_\Omega$, which means that $\sZ_\Omega$ is fixed.]
\vspi
Claim: $X_\Omega$ is not $\N$-compact.
\vspi
[Let $U \subset X_\Omega$ be clopen $-$then $\#(\ov{U} \cap \ov{X_\Omega - U}) \leq \omega$.  Therefore the plane 
$\R^2$ is not disconnected by $\ov{U} \cap \ov{X_\Omega - U}$, so either $\#(U) \leq \omega$ or 
$\#(X_\Omega - U) \leq \omega$.  Consider the collection $\sU$ of all clopen $U \subset X_\Omega$ for which 
$\#(X_\Omega - U) \leq \omega$ $-$then $\sU$ is a clopen set ultrafilter on $X_\Omega$ with the countable intersection property such that $\cap \hspace{0.03cm} \sU = \emptyset$ (every $x \in X_\Omega$ has a countable clopen neighborhood).]
\vspi
[Note: \ The Kunen line $X_\Omega$ is $\R$-compact (same argument as above) but, in contrast to the Kunen plane, it is also $\N$-compact.  For this, it need only be shown that $\dim X_\Omega = 0$.
\vspi
%%----------------------------------------------------------------------------------------------09
Claim: Let $A \subset X_\Omega$ be countable and closed $-$then there exists a countable open 
$U \subset X_\Omega$ : $A \subset U$ $\&$ $U = \ov{U}$, the bar denoting closure in $\R$.
\vspi
[One can assume that \mA is closed in $\R$.  
Write $A = \ds\bigcap\limits_n O_n = \ds\bigcap\limits_n  \ov{O}_n$, where the $O_n$ are $\R$-open and $\forall \ n$: $O_n \supset O_{n+1}$.  Enumerate $A: \{a_n\}$, and for each $n$ choose a compact countable open $U_n \subset X_\Omega$: $a_n \in U_n$ and $U_n \subset O_n$.  Consider 
$U = \ds\bigcup\limits_n U_n$.]
\vspi
To prove that $\dim X_\Omega = 0$, it suffices to take an arbitrary pair $(A,B)$ of disjoint closed subsets of 
$X_\Omega$ and construct a pair $(U_A, U_B)$ of disjoint clopen subsets of $X_\Omega$: 
$
\begin{cases}
\ A \subset U_A\\
\ B \subset U_B
\end{cases}
. \ 
$
Since $\#(\ov{A} \cap \ov{B}) \leq \omega$, by the claim there exists a countable open $O \subset X_\Omega$: 
$\ov{A} \cap \ov{B} \subset O$ $\&$ $O = \ov{O}$.  Pick disjoint $\R$-open sets $O_A$ and $O_B$: 
$
\begin{cases}
\ A - O \subset \ov{A} - O \subset O_A \subset \R - O\\
\ B - O \subset \ov{B} - O \subset O_B \subset \R - O
\end{cases}
, \ 
$
with $\#((\R - O) - (O_A \cup O_B)) \leq \omega$ (possible because it is a question of $\R$ as opposed to $\R^2$).  
Pass to $\R- (O_A \cup O_B)$ and use the claim once again to choose a countable open $P \subset X_\Omega$: 
$\R- (O_A \cup O_B) \subset$  
$P \subset \R - ((\ov{A} - O) \cup (\ov{B} - O))$ $\&$ $P = \ov{P}$ $-$then 
$
\begin{cases}
\ (O_A \cup P) \cap (\R - O)\\
\ (O_B - P) \cap (\R - O)
\end{cases}
$
are disjoint clopen subsets of $X_\Omega$ containing 
$
\begin{cases}
\ A - O\\
\ B - O
\end{cases}
, \ 
$
respectively.  On the other hand, \mO is a normal subspace of $X_\Omega$ of zero topological dimension (cf. Proposition 2), so we can find disjoint clopen sets $P_A$ and $P_B$ in $X_\Omega$: 
$
\begin{cases}
\ A \cap O \subset P_A \subset O\\
\ B \cap O \subset P_B \subset O
\end{cases}
. \ 
$
Now put 
$
\begin{cases}
\ U_A = ((O_A \cup P) \cap (\R - O)) \cup P_A\\
\ U_B = ((O_B - P) \cap (\R - O)) \cup P_B
\end{cases}
$
.]
\\[.5cm] %? need thisxtra space why is this not behaving well here
\endgroup %%------------------------------------<<

\begingroup%%----------------------------------->>
\fontsize{9pt}{11pt}\selectfont
\textbf{\small EXAMPLE \  (\un{The van Douwen Plane})} \ 
\index{The van Douwen Plane} 
\quad 
The object is to equip $X = \R^2$ with a first countable, separable topology that is finer than the usual topology (hence Hausdorff) and under which $X = \R^2$ is locally compact and normal and zero dimensional and $\R$-compact but not 
$\N$-compact.  
Let $\{U_n\}$ be a countable basis for $\R^2$ with $U_0 = \R^2$.  
Assign to each $x \in \R^2$ the sets 
$O_k(x) = \ds\bigcap\limits_n \{U_n: n \leq k \ \& \ x \in U_n\}$ 
$-$then the collection $\{O_k(x)\}$ is a neighborhood basis at $x$ in $\R^2$.  
Obviously, $x \in O_l(y)$ $\implies$ $O_k(x) \subset O_l(y)$ $(\forall \ k \geq l)$.  Let 
$\{x_\alpha: \alpha < 2^\omega\}$ be an enumeration of $\R^2$ and put $X_\alpha = \{x_\beta: \beta < \alpha\}$ 
$-$then 
$X_c = \R^2$ $(c = 2^\omega)$.  \ 
We shall assume that $X_\omega = \Q^2$.  \ 
Fix an enumeration 
$\{(A_\alpha,B_\alpha)$ : $\alpha < 2^\omega\}$ of the set of all pairs $(A,B)$, where \mA and \mB are countable subsets of $\R^2$ with $\#(\ov{A} \cap \ov{B}) = 2^\omega$, arranging matters in such a way that each pair is listed $2^\omega$ times.  
Here (and below) the bar stands for closure in $\R^2$, while $\cl_c$ will denote the closure operator relative to the upcoming topology $\tau_c$ on $X_c$.  
Define an injection $\Gamma:2^\omega \ra 2^\omega - \omega$ by the prescription
\endgroup

\begingroup%%----------------------------------->>
\fontsize{9pt}{11pt}\selectfont
\vspace{-.25cm}
\[
\Gamma(\gamma) = \min(\{\alpha \in 2^\omega - \omega:A_\gamma \cup B_\gamma \subset X_\alpha, 
x_\alpha \in \ov{A}_\gamma \cap \ov{B}_\gamma\} - \{\Gamma(\beta): \beta < \gamma\}).
\]
\endgroup

\vspace{-.15cm}
\begingroup%%----------------------------------->>
\fontsize{9pt}{11pt}\selectfont
Given $\alpha \in 2^\omega - \omega$, choose a sequence 
$\{s_k(\alpha)\} \subset X_\alpha: \forall \ k, s_k(\alpha) \in O_k(x_\alpha)$, having the property that if 
$\alpha = \Gamma(\gamma)$, then $\{s_k(\alpha)\} \subset \Q^2\cup A_\gamma \cup B_\gamma$ and each of 
$\Q^2$, $A_\gamma$, and $B_\gamma$ contains infinitely many terms of $\{s_k(\alpha)\}$, otherwise
$\{s_k(\alpha)\} \subset \Q^2$.  
Topologize $X = \R^2$ as follows: 
Inductively take for the basic neighborhoods of $x_\alpha$ the sets $K_k(x_\alpha)$, $K_k(x_\alpha)$ being $\{x_\alpha\}$ if 
$\alpha \in \omega$ and 
$\{x_\alpha\} \cup \ds\bigcup\limits_{l \geq k} K_l(x_{s_l(\alpha)})$ if $\alpha \in 2^\omega - \omega$ 
$(k = 1, 2, \ldots)$.  
Needless to say, $\forall \ \alpha$ : $K_k(x_\alpha) \subset O_k(x_\alpha)$, and $\forall \ \alpha,\beta$: 
$x_\alpha \in K_l(x_\beta)$ $\implies$ $K_k(x_\alpha) \subset K_l(x_\beta)$ $(\exists \ k)$.  
Observe too that the 
$K_k(x_\alpha)$ are compact and countable.  
Therefore $X_c$ is a zero dimensional LCH space that is in addition first countable and separable.
\vspi
%%----------------------------------------------------------------------------------------------10
Claim: Let $S, T \subset X_c$.  Suppose that $\ov{S} \cap \ov{T}$ is uncountable $-$ then 
$\cl_c(S) \cap \cl_c(T)$ is uncountable.
\vspi
[There are countable $A, B \subset \R^2$: 
$
\begin{cases}
\ A \subset S \subset \ov{A}\\
\ B \subset T \subset \ov{B}
\end{cases}
. \ 
$
From the definitions, $(A,B) = (A_\alpha,B_\alpha)$ for $2^\omega$ ordinals $\alpha$ and, by construction, 
$x_{\Gamma(\alpha)} \in \cl_c(A_\alpha) \cap \cl_c(B_\alpha)$.  But $\Gamma$ is one-to-one.]
\vspi
[To establish that $X_c$ is normal, it suffices to show that if \mA and \mB are two disjoint closed subsets of $X_c$, 
then there exists a countable open covering 
$\sO = \{O\}$ of $X_c$ such that $\forall \ O \in \sO$ : $\cl_c(O) \cap A = \emptyset$ or 
$\cl_c(O) \cap B = \emptyset$.  \ 
In view of the claim, $\ov{A} \cap \ov{B}$ is countable.  \ 
Let $x \in \ov{A} \cap \ov{B}$ 
$-$then $x \notin A \cup B$, so by regularity there exists an open set $O_x \subset X_c$ containing $x$ : 
$\cl_c(O_x) \cap (A \cup B) = \emptyset$.  \ 
It is equally plain that for any $x \in \R^2 - \ov{A} \cap \ov{B}$  there exists an 
$\R^2$-open set $O_x$ containing 
$x$ : $\ov{O}_x \cap \ov{A} = \emptyset$ or $\ov{O}_x \cap \ov{B} = \emptyset$.  \ 
Select a countable subcollection of $\{O_x: x \in \R^2 - \ov{A} \cap \ov{B}\}$ that covers $\R^2 - \ov{A} \cap \ov{B}$ and combine it with $\{O_x: x \in \ov{A} \cap \ov{B}\}$.
\vspi
Arguing as before, one proves that $X_c$ is $\R$-compact but not $\N$-compact.
\vspi
[Note: \ The van Douwen plane exists in ZFC.  But unlike the Kunen plane, it is not perfect.  Reason: 
$\Q^2 \cup \{x_{\Gamma(\alpha)}: A_\alpha \cup B_\alpha \subset \Q^2\}$ is not a normal subspace of $X_c$.  
However, every closed discrete subspace of $X_c$ is countable, so $X_c$, like the Kunen plane, is collectionwise normal.  
Of course, $X_c$ is not Lindel\"of, thus is not paracompact (being separable), although $X_c$ is countably paracompact.  By the way, if the same procedure is applied to $X = \R$, then the endproduct is a space very different from what was termed the van Douwen line in $\S 1$.]\\
\endgroup %%------------------------------------<<

Is it true that for every normal subspace $Y \subset X$, $\dim Y \leq \dim X$?  In other words, is dim monotonic?  
On closed subspaces, this is certainly the case but, as has been seen above, this is not the case in general.\\

\label{19.5}
\begingroup%%----------------------------------->>
\fontsize{9pt}{11pt}\selectfont
It is false that dim is monotonic on closed subspaces of a nonnormal \mX.  For example, the topological dimension of 
Mysior space is positive but it embeds as a closed subspace of some $\N^\kappa$ and $\dim \N^\kappa = 0$.\\
\endgroup %%------------------------------------<<

\begingroup%%----------------------------------->>
\fontsize{9pt}{11pt}\selectfont
\textbf{\small LEMMA} \  
Let \mX be a nonempty CRH space.  Suppose that \mA is a subspace of \mX which has the EP w.r.t. [0,1] $-$then 
$\dim A \leq \dim X$.\\
\vspace{0.05cm}\
\endgroup %%------------------------------------<<

\begin{proposition} \ %04
Suppose that \mX is hereditarily normal $-$then dim is monotonic iff for every open $U \subset X$: 
$\dim U \leq \dim X$.\\
\end{proposition}

One might conjecture that dim is monotonic if \mX is hereditarily normal.  
This is false: 
Pol-Pol\footnote[2]{\textit{Fund. Math.} \textbf{102} (1979), 137-142.} 
has given an example of a hereditarily normal \mX that has topological dimension zero but which contains for every 
$n = 1, 2, \ldots$ a subspace $X_n: \dim X_n = n$.  Since $\beta X$
%%----------------------------------------------------------------------------------------------11
also has topological dimension zero (cf. Proposition  1), dim is dramatically nonmonotonic even for compact Hausdorff spaces.\\

\begingroup%%----------------------------------->>
\fontsize{9pt}{11pt}\selectfont
Consider the Kunen plane $X_\Omega$ $-$then its one point compactification is hereditarily normal and has topological dimension zero, although $X_\Omega$ appears as an open subspace of positive topological dimension.\\
\endgroup %%------------------------------------<<

\begingroup%%----------------------------------->>
\fontsize{9pt}{11pt}\selectfont
\textbf{\small EXAMPLE} \  
The Isbell$-$Mr\'owka space $\Psi(\N)$ is a nonnormal LCH space.  
While zero dimensional, its ``finer'' topological properties depend on the choice of $\sS$.  \ 
Claim: $\exists$ $\sS$ for which $\dim \Psi(\N) > 0$.  \ 
To this end, replace $\N$ by \ 
$\Q_{[0,1]} \equiv \Q \cap [0,1]$.  \ 
Attach to each $r$, \ $0 < r < 1$, a bijection 
$\iota_r$ : $\{q \in \Q_{[0,1]}: q < r\} \ra \{q \in \Q_{[0,1]}:q > r\}$ such that $q^\prime < q\pp$ iff 
$\iota_r(q^\prime) > \iota_r(q\pp)$.  
Let SEQ be the collection of all sequences $s$ of distinct elements of $\Q_{[0,1]}$ 
satisfying one of the following two conditions: 
(i) $\lim s = 0$ or $\lim s = 1$; 
(ii) $s = t \cup \iota_r(t)$ $(0 < r < 1)$, where $t$ converges to $r$ from the left.  
Because $[0,1]$ is compact, there is a maximal infinite collection $\sS \subset$ SEQ of almost disjoint infinite subsets of 
$\Q_{[0,1]}$.  
Consider the corresponding Isbell-Mr\'owka space $X = \Psi(\Q_{[0,1]})$, i.e., $X = \sS \cup \Q_{[0,1]}$ 
$-$then $\dim \sS = 0$ and $\dim \Q_{[0,1]} = 0$, yet $\dim X > 0$.
To see this, define a continuous function $f:X \ra [0,1]$ by 
$
\begin{cases}
\ f(q)  = q \hspace{0.85cm} (q \in \Q_{[0,1]})\\
\ f(s) = \lim s \hspace{0.35cm} (s \in \sS)
\end{cases}
\hspace{-.2cm}.
$
Verify that there is no clopen subset of \mX containing $f^{-1}(0)$ and missing $f^{-1}(1)$.
\vspi
\label{19.7}
\label{19.12}
[Note: \ Mr\'owka \footnote[3]{\textit{Fund. Math.} \textbf{94} (1977), 83-92.} 
has shown that for certain choices of $\sS$, $\beta(\Psi(\N)) = \Psi(\N)_\infty$, hence $\dim \Psi(\N) = 0$.  At the opposite extreme,  
Terasawa\footnote[6]{\textit{Topology Appl.} \textbf{11} (1980), 93-102.}
proved that for any $n = 1, 2, \ldots$ or $\infty$, it is possible to find an $\sS$ such that the associated $\Psi(\N)$ has topological dimension $n$ but at the same time is expressible as the union of two zero sets, each having topological dimension zero.]\\
\endgroup %%------------------------------------<<

\textbf{\small LEMMA} \  
Let $\sU$ be a finite open covering of \mX $-$then $\sU$ has a finite open refinement of order $\leq n+1$ iff $\sU$ has a finite closed refinement of order $\leq n+1$.

[Suppose that  \ 
$\sU = \{U_1, \ldots, U_k\}$.  \ 
Let \ $\sV = \{V_1, \ldots, V_k\}$ be a precise open refinement of 
$\sU$ of order $\leq n + 1$ $-$then $\sV$ has a precise open refinement 
$\sW = \{W_1, \ldots, W_k\}$ such that $\forall \ i$: 
$\ov{W}_i \subset V_i$.  
And the order of $\ov{\sW}$ is $\leq n + 1$.  \ 
To go the other way, let $\sA = \{A_1, \ldots, A_k\}$ be a precise closed refinement of $\sU$ of order $\leq n + 1$ 
$-$then it will be enough to produce a precise open refinement 
$\sV = \{V_1, \ldots, V_k\}$ of  \ $\sU$ such that $\forall \ i$: $A_i \subset V_i \subset U_i$ and 
$A_{i_1} \cap \cdots \cap A_{i_m} \neq \emptyset$ iff $V_{i_1} \cap \cdots \cap  V_{i_m} \neq \emptyset$.  Here 
$i_1, \ldots, i_m$ are natural numbers, each $\leq k$.  This can be done by a simple iterative procedure.  Denote by $B_1$ the union of all intersections of members of the collection $\{A_1, \ldots, A_k\}$ which are disjoint from $A_1$ and choose an open set 
$
V_1 : \ 
\begin{cases}
\ A_1 \subset V_1\\
\ \ov{V}_1 \subset U_1
\end{cases}
$
$\&$ $B_1 \cap \ov{V}_1 = \emptyset$.  Denote by $B_2$ the union of all 
%%----------------------------------------------------------------------------------------------12
intersections of members of the collection 
$\{\ov{U}_1, A_2, \ldots, A_k\}$ which are disjoint from $A_2$ and choose an open set $V_2$ : 
$
\begin{cases}
\ A_2 \subset V_2\\
\ \ov{V}_2 \subset U_2
\end{cases}
\& \ B_2 \cap \ov{V}_2 = \emptyset.
$
ETC.]\\
\vspace{0.25cm}

\index{Countable Union Lemma}
\textbf{\small COUNTABLE UNION LEMMA} \quad 
Suppose that $X = \bigcup\limits_1^\infty A_j$, where the $A_j$ are closed subspaces of \mX such that $\forall \ j$, 
$\dim A_j \leq n$ $-$then $\dim X \leq n$, hence $\dim X = \sup\dim A_j$.

[Let $\sU = \{U_i\}$ be a finite open covering of \mX.  Put $A_0 = \emptyset$.  Claim: There exists a sequence 
$\sU_0, \sU_1, \ldots$, of finite open coverings $\sU_j = \{U_{i,j}\}$ of \mX such that $U_{i,0} \subset U_i$ but 
\[
\ov{U}_{i,j} \subset U_{i,j-1} \ \& \ \ord(\{A_j \cap \ov{U}_{i,j}\}) \leq n + 1
\]
if $j \geq 1$.  
To prove this, we shall proceed by induction on $j$, setting \ $\sU_0 = \sU$ and then assuming that the $\sU_j$ have been defined for all $j < j_0$, where $j_0 \geq 1$.  
Since  \ $\{A_{j_0} \cap U_{i,j_0-1}\}$ is a finite open covering of 
$A_{j_0}$ and since $\dim A_{j_0} \leq n$, there exist open subsets 
$V_i \subset A_{j_0} \cap U_{i,j_0 - 1}$ of $A_{j_0}$ 
such that $A_{j_0} = \bigcup\limits_i V_i$ and $\ord(\{V_i\}) \leq n + 1$.  Let 
$W_i = (U_{i, j_0 - 1} - A_{j_0}) \cup V_i$ $-$then $\{W_i\}$ is a finite open covering of \mX and 
$\ord(\{A_{j_0} \cap W_i\}) \leq n + 1$.
The induction is completed by choosing the elements $U_{i,j_0}$ of $\sU_{j_0}$ subject to $\ov{U}_{i,j_0} \subset W_i$.  
By construction, the collection 
$\{\bigcap\limits_{j \geq 1} \ov{U}_{i,j}\}$ is a precise closed refinement of $\sU = \{U_i\}$ of order $\leq n + 1$, 
so from the lemma $\dim X \leq n$.]\\

Example: \ $\dim[0,1] = 1$ $\implies$ $\dim \R = 1$.\\

\begingroup%%----------------------------------->>
\fontsize{9pt}{11pt}\selectfont 
\textbf{\small FACT} \  
Suppose that \mX is normal of topological dimension $n \geq 1$ $-$then there exists a sequence of pairwise disjoint closed subspaces $A_j$ of \mX such that $\forall \ j$, $\dim A_j = n$.\\
\endgroup %%------------------------------------<<

\label{19.18}
\begingroup%%----------------------------------->>
\fontsize{9pt}{11pt}\selectfont 
A CRH space \mX is said to be 
\un{strongly paracompact}
\index{strongly paracompact} if every open covering of \mX 
has a star finite open refinement.  Any paracompact LCH space \mX is strongly paracompact (cf. $\S 1$, Proposition 2).  
Also: \mX Lindel\"of $\implies$ \mX strongly paracompact and \mX connected + strongly paracompact $\implies$ 
\mX Lindel\"of.  Not every metric space is strongly paracompact (consider the star space $S(\kappa), \ \kappa > \omega)$.\\
\endgroup %%------------------------------------<<

\label{19.26}
\begingroup%%----------------------------------->>
\fontsize{9pt}{11pt}\selectfont 
\textbf{\small FACT} \  
Suppose that \mX is normal and \mY is a strongly paracompact subspace of \mX $-$then $\dim Y \leq \dim X$.
\vspi
[The assertion is trivial if $\dim X = \infty$, so assume that $\dim X = n$ is finite.  \ 
Let $\{U_i\}$ be a finite open covering of \mY; 
let $O_i$ be an open subset of \mX such that $U_i = Y \cap O_i$ and put $O = \ds\bigcup\limits_i O_i$.  
Assign to each $y \in Y$ a neighborhod $O_y$ of $y$ in \mX \  : $\ov{O}_y \subset O$ $-$then $\{Y \cap O_y\}$ is an open covering of \mY, thus has a star finite open refinement $\sP$.  
Write $\sP = \ds\coprod\limits_j \sP_j$, the equivalence relation corresponding to
%%----------------------------------------------------------------------------------------------13
this partition being $P^\prime \sim P\pp$ iff there exists a finite collection of sets $P_1, \ldots, P_r$ in $\sP$ with 
$P_1 = P^\prime$, $P_r = P\pp$ and $P_1 \cap P_2 \neq \emptyset$, $\ldots$, $P_{r-1} \cap P_r \neq \emptyset$.  
Since 
$\sP$ is star finite, each of the $\sP_j$ is countable.  
Let $Y_j = \ds\bigcup\{\ov{P}:P \in \sP_j\}$, where $\ov{P}$ is the closure of \mP in \mX.  
Being an $F_\sigma$, $Y_j$ is normal and therefore, by the countable union lemma, $\dim Y_j \leq n$.  
But 
$Y_j$ is contained in $O = \ds\bigcup\limits_i O_i$, so there exists an open covering $\{O_{i,j}\}$ of $Y_j$ such that 
$\forall \ i$: $O_{i,j} \subset O_i$ $\&$ $\ord(\{O_{i,j}\}) \leq n + 1$.  Let 
$V_i = Y \cap \ds\bigcup\limits_j (O_{i,j} \bigcap \cup \sP_j)$ $-$then $\{V_i\}$ is a precise open refinement of $\{U_i\}$ of order $\leq n+1$.]\\
\endgroup %%------------------------------------<<

\begingroup%%----------------------------------->>
\fontsize{9pt}{11pt}\selectfont 
The preceding result if false if ``paracompact'' is substituted for ``strongly paracompact''.  
Example: Consider Roy's metric space \mX sitting inside its zero dimensional compactification $\zeta X$.\\
\endgroup %%------------------------------------<<

\label{19.30}
\label{19.32}
\label{19.37}
\begingroup%%----------------------------------->>
\fontsize{9pt}{11pt}\selectfont 
The countable union lemma retains its validity in the completely regular situation provided the $A_j$ are subspaces of \mX which have the EP w.r.t. $[0,1]$.  Proof: The closure of $A_j$ in $\beta X$ is $\beta A_j$, so if 
$Y = \ds\bigcup\limits_1^\infty \beta A_j$, then \mY is normal and therefore, by the countable union lemma, $\dim Y \leq n$, from which $\dim X = \dim \beta X = \dim \beta Y = \dim Y \leq n$.

[Note: \  According to Terasawa (cf. p. \pageref{19.7}), there exists a completely regular \mX of topological dimension $n$ such that $X = X_1 \cup X_2$, where $X_1$ and $X_2$ are zero sets with 
$
\begin{cases}
\ \dim X_1 = 0\\
\ \dim X_2 = 0
\end{cases}
$
.  
Therefore the countable union lemma can fail even when the hypothesis ``closed set'' is strengthened to ``zero set''.]\\
\endgroup %%------------------------------------<<

\begingroup%%----------------------------------->>
\fontsize{9pt}{11pt}\selectfont 
\textbf{\small LEMMA} \  
Let \mX be a nonempty CRH space.  Suppose that \mA is a $\sZ$-embedded subspace of \mX $-$then $\dim A \leq \dim X$.
\vspi
[Assume that $\dim X \leq n$.  Let $\{U_i\}$ be a finite cozero set covering of \mA; let $O_i$ be a cozero set in 
$\beta X$ such that $U_i = A \cap O_i$ and put $O = \ds\bigcup\limits_i O_i$ $-$then \mO is a cozero set in $\beta X$, so by the countable union lemma, $\dim O \leq \dim \beta X = \dim X \leq n$.  Therefore there exists a cozero set covering $\{P_i\}$ of \mO of order $\leq n + 1$ such that $\forall \ i$: $P_i \subset O_i$.  Consider the collection $\{A \cap P_i\}$.]\\
\endgroup %%------------------------------------<<

Recall: Every subspace of a perfectly normal space is perfectly normal.  So: \mX perfectly normal $\implies$ \mX hereditarily normal.  The conjunction  perfectly normal + paracompact is hereditary to all subspaces.  Reason: Every open set is an 
$F_\sigma$ and an $F_\sigma$ in a paracompact space is paracompact.  For example, the class of stratifiable spaces or the class of CW complexes realize this conjunction .

[Note: \ The ordinal space $[0,\Omega]$ is hereditarily normal but not perfectly normal and its product with $[0,1]$ is normal but not hereditarily normal.]\\

\begin{proposition} \ %05
Suppose that \mX is perfectly normal $-$then dim is monotonic.
\end{proposition}

[Apply the countable union lemma to an open subset of \mX and then quote Proposition 4.]\\

%%----------------------------------------------------------------------------------------------14
\label{19.16}
\begingroup%%----------------------------------->>
\fontsize{9pt}{11pt}\selectfont 
Working under CH, the procedure for manufacturing the Kunen line or the Kunen plane is just a specialization to $\R$ or $\R^2$ of a general ``machine'' for refining topologies.  Thus suppose that \mX is a set of cardinality $\Omega$ equipped with a Hausdorff topology $\tau$ which is first countable, hereditarily separable and perfectly normal $-$then a 
\un{Kunen modification}
\index{Kunen modification} 
of $\tau$ is a topology $K\tau$ on \mX finer than $\tau$ which is zero dimensional, locally compact, first countable, hereditarily separable and perfectly normal (but not Lindel\"of) such that each 
$x \in X$ has a countable clopen neighborhood and $\forall \ S \subset X$: $\#(\cl_r(\sS) - \cl_{K\tau}(\sS))$ $\leq \omega$.
\vspi
[Note: \ Any $\tau$ having the stated properties admits a Kunen modification $K\tau$ (cf. p. \pageref{19.8}).]\\
\endgroup %%------------------------------------<<

\begingroup%%----------------------------------->>
\fontsize{9pt}{11pt}\selectfont 
\label{19.35} %dmc mnft
\textbf{\small FACT} \  
[Assume CH]  If $\dim(X,\tau) \geq n$, then $\dim(X,K\tau) \geq n-1$ and if $\dim(X,\tau) \leq n$ then 
$\dim(X,K\tau) \leq n$.\\
\endgroup %%------------------------------------<<

\begin{proposition} \ %06
The statement $\dim X \leq n$ is true iff every neighborhood finite open covering of \mX has a numerable open refinement of order $\leq n+1$.
\end{proposition}

[Let $\sU$ be neighborhood finite open covering of \mX $-$then $\sU$ is numerable, hence has a numerable open refinement that is both neighborhood finite and $\sigma$-discrete, say $\sV = \bigcup\limits_n \sV_n$ (cf. $\S 1$, Proposition 12).  
Choose a partition of unity $\{\kappa_V\}$ on \mX subordinate to $\sV$.
Put 
$\ds f_n = \ds\sum\limits_{V \in \sV_n} \kappa_V$: The collection $\{f_n^{-1}(]0,1])\}$ is a countable cozero set covering of \mX, thus has a countable star finite cozero set refinement $\{O_k\}$ (cf. p. \pageref{19.9}).  
Fix a sequence of integers 
$1 = n_1 < n_2 \cdots$: $O_k \cap O_l = \emptyset$ if $k \leq n_i$ and $l \geq n_{i+2}$ $(i = 1, 2, \ldots)$.  
The subspace 
$\bigcup\limits_{k \leq n_2} O_k$ is a cozero set and so by the countable union lemma its topological dimension is $\leq n$.  
Accordingly, there exists a covering 
$\sW_1 = \{W_1, \ldots, W_{n_1}, W_{n_1 + 1}^\prime, \ldots, W_{n_2}^\prime\}$ of 
$\bigcup\limits_{k \leq n_2} O_k$ by cozero sets of order $\leq n + 1$ such that 
$
\begin{cases}
\ W_k \subset O_k \quad (k \leq n_1)\\
\ W_k^\prime \subset O_k \quad (n_1 < k \leq n_2)
\end{cases}
. \ 
$
Next, there exists a covering 
$\sW_2 = \{W_{n_1 + 1},\ldots, W_{n_2}, W_{n_2 + 1}^\prime, \ldots, W_{n_3}^\prime\}$ of 
$W_{n_1 + 1}^\prime \cup \cdots \cup W_{n_2}^\prime \cup O_{n_2 + 1} \cup \cdots \cup O_{n_3}$
by cozero sets of order $\leq n + 1$ such that 
$
\begin{cases}
\ W_k \subset W_k^\prime \quad (n_1 < k \leq n_2)\\
\ W_k^\prime \subset O_k \quad (n_2 < k \leq n_3)
\end{cases}
. \ 
$
Iterate to get a covering $\sW = \{W_k\}$ of \mX by cozero sets of order $\leq n + 1$ such that $\forall \ k$: 
$W_k \subset O_k$.  
The collection $\bigcup\limits_k \sU \cap W_k$ is a numerable open refinement of $\sU$ of 
order $\leq n + 1$.]\\

Suppose that \mX is paracompact $-$then it follows from Proposition 6 that $\dim X \leq n$ iff every open covering of \mX has an open refinement of order $\leq n + 1$.\\

\begingroup%%----------------------------------->>
\fontsize{9pt}{11pt}\selectfont 
Since cozero sets are $\sZ$-embedded and since dim is monotonic on $\sZ$-embedded subspaces, Proposition 6 goes through without change in the completely regular situation provided one works with numerable open coverings and numerable open refinements.\\
\endgroup %%------------------------------------<<

%%----------------------------------------------------------------------------------------------15
\textbf{\small SUBLEMMA} \quad 
The statement $\dim X \leq n$ is true iff every open covering $\{U_1, \ldots, U_{n+2}\}$ of \mX has a precise open refinement 
$\{V_1, \ldots, V_{n+2}\}$ such that $\bigcap\limits_1^{n+2} V_i = \emptyset$.

[When turned around, the nontrivial assertion is that if $\dim X > n$, 
then there exists an open covering 
$\{U_1, \ldots, U_{n+2}\}$ of \mX, every precise open refinement $\{V_1, \ldots, V_{n+2}\}$  
of which satisfies the condition 
$\bigcap\limits_1^{n+2} V_i \neq \emptyset$.  
But $\dim X > n$ means that there exists an open covering 
$\{O_1, \ldots, O_k\}$ of \mX that has no precise open refinement of order $\leq n + 1$.  
By making at most a finite number of replacements, matters can be arranged so as to ensure that if $\{P_1, \ldots, P_k\}$ is a precise open refinement of 
$\{O_1, \ldots, O_k\}$, then $P_{i_1}\cap \cdots \cap P_{i_m} \neq \emptyset$ whenever 
$O_{i_1}\cap \cdots \cap O_{i_m} \neq \emptyset$.  
Here $i_1, \ldots, i_m$ are natural numbers, each $\leq k$.  We can and will assume that $\bigcap\limits_1^{n+2}O_i \neq \emptyset$.  
Put $U_i = O_i$ $(i \leq n+1)$, 
$U_{n+2} = \bigcap\limits_{n+2}^k O_i$ $-$then $\{U_1, \ldots, U_{n+2}\}$ is an open covering of \mX with the property in question.  
In fact, let $\{V_1, \ldots, V_{n+2}\}$ be an open covering of \mX such that $\forall \ i$: $V_i \subset U_i$.  
The covering $\{V_1, \ldots, V_{n+1}, V_{n+2}\cap O_{n+2}, \ldots V_{n+2} \cap O_k\}$ is a precise open refinement of 
$\{O_1, \ldots, O_k\}$ and 
$\bigcap\limits_1^{n+2} V_i \supset \bigl(\bigcap\limits_1^{n+1} V_i\bigr) \cap (V_{n+2} \cap O_{n+2}) \neq \emptyset$.]\\

\textbf{\small LEMMA} \  
The statement $\dim X \leq n$ is true iff for every collection $\{(A_i,B_i)$:$i = 1, \ldots, n+1\}$ of $n + 1$ pairs of disjoint closed subsets of \mX there exists a collection $\{\phi_i:i = 1, \ldots, n+1\}$ of $n + 1$ continuous functions 
$\phi_i:X \ra [0,1]$ such that 
$
\begin{cases}
\ \restr{\phi_i}{A_i} = 0\\
\ \restr{\phi_i}{B_i} = 1
\end{cases}
$
and $\bigcap\limits_1^{n+1} \phi_i^{-1}(1/2) = \emptyset$.

[Necessity: Put $B_{n+2} = \ds\bigcap\limits_1^{n+1} A_i$ 
$-$then $\ds\bigcap\limits_1^{n+2} B_i = \emptyset$, so there exists an open covering $\{U_1, \ldots, U_{n+2}\}$ of \mX such that $B_i \subset U_i$ and 
$\bigcap\limits_1^{n+2} U_i = \emptyset$.  
Since $A_i \subset U_{n+2}$, we can replace $U_i$ by $U_i - A_i$ and force 
$A_i \subset X - U_i$.  
Fix a precise closed refinement $\{C_1, \ldots, C_{n+2}\}$ of $\{U_1, \ldots, U_{n+2}\}$ with 
$B_i \subset C_i$.  Let $\phi_i:X \ra [0,1]$ be a continuous function such that $\restr{\phi_i}{X - U_i} = 0$ and 
$\restr{\phi_i}{C_i} = 1$.  Obviously, 
$
\begin{cases}
\ \restr{\phi_i}{A_i} = 0\\
\ \restr{\phi_i}{B_i} = 1
\end{cases}
.
$
And finally, 
$\bigcap\limits_1^{n+1} \phi_i^{-1}(1/2) \subset$ 
$\bigcap\limits_1^{n+1}(U_i - C_i) \subset$
$\bigcap\limits_1^{n+2} U_i = \emptyset$.

Sufficiency: Let $\{U_1, \ldots, U_{n+2}\}$ be an open covering of \mX.  
Fix a precise closed refinement $\{C_1, \ldots, C_{n+2}\}$ for it and let
$
\begin{cases}
\ A_i = X - U_i\\
\ B_i = C_i
\end{cases}
(i = 1, \ldots, n+1).\ 
$
The pairs $(A_i,B_i)$ satisfy our hypotheses, so choose the $\phi_i$ as there and then let 
$
\begin{cases}
\ O_i = \{x: \phi_i(x) < 1/2\}\\
\ P_i = \{x: \phi_i(x) > 1/2\}
\end{cases}
. \ 
$
Note that 
$\bigcap\limits_1^{n+1} (X - (O_i \cup P_i)) = $ 
$\bigcap\limits_1^{n+1} \phi_i^{-1}(1/2) = \emptyset$, hence that 
$X = \bigcup\limits_1^{n+1} O_i \cup \bigcup\limits_1^{n+1} P_i$.  Put
%%----------------------------------------------------------------------------------------------16
$V_i = P_i$ $(i \leq n+1)$, $V_{n+2} = U_{n+2} \cap \bigcup\limits_1^{n+1} O_i$ 
$-$then $\{V_1, \ldots, V_{n+2}\}$ is a 
precise open refinement of $\{U_1, \ldots, U_{n+2}\}$ such that $\bigcap\limits_1^{n+2} V_i = \emptyset$.  
The sublemma therefore implies that $\dim X \leq n$.]\\

\begingroup%%----------------------------------->>
\fontsize{9pt}{11pt}\selectfont 
The characterization of $\dim X \leq n$ given by the lemma extends to the completely regular situation so long as it is formulated in terms of disjoint pairs $(A_i,B_i)$ of zero sets.\\
\endgroup %%------------------------------------<<

When the context dictates, we shall abuse the notation and write $\bS^n$ for the frontier of $[0,1]^{n+1}$.\\

\index{Alexandroff's Criterion}
\textbf{\small ALEXANDROFF'S CRITERION} \quad 
The statement $\dim X \leq n$ is true iff every closed subset $A \subset X$ has the EP w.r.t $\bS^n$.

[Necessity: Given $f \in C(A,\bS^n)$: $f = (f_1, \ldots, f_{n+1})$, let 
$
\begin{cases}
\ A_i = \{x: f_i(x) = 0\}\\
\ B_i = \{x:f_i(x) = 1\}
\end{cases}
$
$-$then \mA is the union $\bigcup\limits_i(A_i \cup B_i)$ and the preceding lemma is applicable to the pairs $(A_i,B_i)$.  
The corresponding $\phi_i:X \ra [0,1]$ combine to determine a continuous function $\phi:X \ra [0,1]^{n+1}$, the restriction of which to \mA defines an element $\psi \in C(A,\bS^n)$.  
Put
$H(x,t) = (1 - t)\psi(x) + tf(x)$ $((x,t) \in IA)$ $-$then $H \in C(IA,\bS^n)$, so $\psi$ and $f$ are homotopic.  
On the other hand, $\bS^n$ is a retract of $[0,1]^{n+1}$ punctured at its center $(1/2, \ldots, 1/2)$.  
Since 
$\bigcap\limits_1^{n+1} \phi_i^{-1}(1/2) = \emptyset$, it follows that $\psi$ has an extension $\Psi \in C(X,\bS^n)$.  
But \mA has the HEP w.r.t. $\bS^n$ 
(cf. p. \pageref{19.10}), therefore $f$ has an extension $F \in C(X,\bS^n)$.

Sufficiency: Consider an arbitrary collection $\{(A_i,B_i):i = 1, \ldots, n+1\}$ of $n+1$ pairs of disjoint closed subsets of \mX.  Put $A = \bigcup\limits_i(A_i \cup B_i)$.  Choose $f_i \in C(A,[0,1])$ such that 
$
\begin{cases}
\ \restr{f_i}{A_i} = 0\\
\ \restr{f_i}{B_i} = 1
\end{cases}
$
and then combine the $f_i$ to determine a continuous function $f:A \ra \bS^n$.  
By assumption, $f$ has an extension 
$F \in C(X,\bS^n)$.  
Write $\phi_i$ for the $i^\text{th}$ component of \mF $-$then $\restr{\phi_i}{A} = f_i$ and 
$\bigcap\limits_1^{n+1} \phi_i^{-1}(1/2) = \emptyset$.  
That $\dim X \leq n$ is thus a consequence of the preceding 
lemma.]\\

\begingroup%%----------------------------------->>
\fontsize{9pt}{11pt}\selectfont 
\textbf{\small EXAMPLE} \  
Take for \mX the long ray $L^+$ $-$then $\dim X = 1$.
\vspi
[Since $\dim X > 0$, one need only show that $\dim X \leq 1$.  But real valued continuous functions are constant on 
``tails'', so Alexandroff's criterion is applicable.]\\
\endgroup %%------------------------------------<<

\begingroup%%----------------------------------->>
\fontsize{9pt}{11pt}\selectfont 
\textbf{\small FACT} \  Let \mX be a compact Hausdorff space.  
Suppose that $X = \ds\bigcup\limits_1^\infty A_j$, where the
$A_j$ are closed subspaces of \mX such that $\forall \ i \neq j$: $\dim (A_i \cap A_j) < n$ $-$then each $A_j$ has the EP w.r.t. $\bS^n$.
\vspi
[Recall that if \mX is a connected compact Hausdorff space admitting a disjoint decompostion 
$\ds\bigcup\limits_1^\infty A_j$ by closed subspaces $A_j$, then $A_j = X$ for some $j$.]\\
\endgroup %%------------------------------------<<

%%----------------------------------------------------------------------------------------------17
\begingroup%%----------------------------------->>
\fontsize{9pt}{11pt}\selectfont 
Application: Because the identity map $\bS^n \ra \bS^n$ cannot be extended continuously over $[0,1]^{n+1}$, 
$\R^{n+1}$ cannot be covered by a sequence $\{K_j\}$ of compact sets such that $\forall \ i \neq j$: 
$\dim (K_i \cap K_j) < n$.
\vspi
[Note: \ With more work, one can do better in that ``compact'' can be replaced by ``closed'' 
(cf. p. \pageref{19.11}).]\\
\endgroup %%------------------------------------<<

\begingroup%%----------------------------------->>
\fontsize{9pt}{11pt}\selectfont 
The compactness assumption on \mX in the preceding result is essential.  
Example: Take for \mX a one dimensional connected locally compact subspace of the plane admitting a disjoint decomposition 
$\ds\bigcup\limits_1^\infty A_j$ by nonempty closed proper subspaces $A_j$, fix two indices $i \neq j$, and consider the continuous function $f:A_i \cup A_j \ra \bS^0$ which is 0 on $A_i$ and 1 on $A_j$.\\
\endgroup %%------------------------------------<<

\begingroup%%----------------------------------->>
\fontsize{9pt}{11pt}\selectfont 
Using Alexandroff's criterion, 
Cantwell\footnote[2]{\textit{Proc. Amer. Math. Soc.} \textbf{19} (1968), 821-825.} 
proved that the statement $\dim X \leq n$ is true iff the closed unit ball in $BC(X,\R^{n+1})$ is the convex hull of its extreme points $(n = 1, 2, \ldots)$.
\vspi
[Note: \ Let \mX be a nonempty CRH space $-$then the extreme points of the closed unit ball in $BC(X,\R^{n+1})$  are the functions whose range is a subset of $\bS^n$ and it is always true that the closed unit ball in $BC(X,\R^{n+1})$ is the closed convex hull of its extreme points $(n = 1, 2, \ldots)$, a purely topological assertion.  By contrast, the closed unit ball in 
$BC(X)$ is the closed convex hull of its extreme points iff $\dim X = 0$.]\\
\endgroup %%------------------------------------<<

\begingroup%%----------------------------------->>
\fontsize{9pt}{11pt}\selectfont 
In the completely regular situation, there is only a partial analog to Alexandroff's criterion.
\\
\indent\indent\indent (1)  Suppose that every zero set $A \subset X$ has the EP w.r.t. $\bS^n$ $-$then $\dim X \leq n$.  Proof: 
Since for any pair $(A,B)$ of disjoint zero sets there exists a continuous function $f:X \ra [0,1]$ such that 
$
\begin{cases}
\ \restr{f}{A} = 0\\
\ \restr{f}{B} = 1
\end{cases}
,\ 
$
the argument used in the normal case can be transcribed in the obvious way.
\\
\indent\indent\indent (2) Suppose that $\dim X \leq n$ $-$then every subset $A \subset X$ which has the EP w.r.t. $[0,1]$ has the EP w.r.t. $\bS^n$.  Proof: Since $\dim X = \dim \beta X$, $\beta A$, the closure of \mA in $\beta X$, has the EP w.r.t. $\bS^n$.
\vspi
[Note: \ This need not be true if \mA is a zero set.  
Example: Take, after Terasawa (cf. p. \pageref{19.12}), $X = X_1 \cup X_2$, where $\dim X = 1$ and $X_1$ and $X_2$ are zero sets with  
$
\begin{cases}
\ \dim X_1 = 0\\
\ \dim X_2 = 0
\end{cases}
$
$-$then either $X_1$ or $X_2$ fails to have the EP w.r.t. $[0,1]$ (otherwise $\dim X = \max\{\dim X_1,\dim X_2\}$).  To be specific, assume that it is $X_1$.  Put $A = X_1$ and choose a continuous function $\phi:A \ra [0,1]$ that does not extend to a continuous function $\Phi:X \ra [0,1]$ $-$then $f = (\phi,0)$ is a continuous function $A \ra \bS^1$ that does not extend to a continuous function $F:X \ra \bS^1$.]\\
\endgroup %%------------------------------------<<

Let \mY be a topological space $-$then a map $f \in C(X,Y)$ is said to be 
\un{universal}
\index{universal (map)} 
if $\forall \ g \in C(X,Y)$ $\exists$ $x \in X$: $f(x) = g(x)$.  A universal map is clearly surjective.  Note too that if there is a universal map $X \ra Y$, then every element of $C(Y,Y)$ must have a fixed point.\\

%%----------------------------------------------------------------------------------------------18

\textbf{\small LEMMA} \  
A continuous function $f:X \ra [0,1]^{n+1}$ is universal iff the restriction $f^{-1}(\bS^n) \ra \bS^n$ has no extension 
$F \in C(X,\bS^n)$.

[Necessity: To get a contradiction, suppose that there exists a continuous function $F:X \ra \bS^n$ which agrees with $f$ on $f^{-1}(\bS^n)$ and then postcompose \mF with the antipodal map $\bS^n \ra \bS^n$.

Sufficiency: To get a contradiction, suppose that there exists a continuous function $g:X \ra [0,1]^{n+1}$ such that 
$f(x) \neq g(x)$ for every $x \in X$ and define a continuous function $F:X \ra \bS^n$ by setting $F(x)$ equal to the intersection of $\bS^n$ with the ray containing $f(x)$ which emanates from $g(x)$.]\\

\label{20.4}
It therefore follows that $\dim X \geq n$ iff there exists a universal map $f:X \ra [0,1]^n$.  
Example: $\dim [0,1]^n \geq n$.  Indeed, the Brouwer fixed point theorem says that the identity map 
$[0,1]^n \ra [0,1]^n$ is universal.  
\label{19.15}
Example: $\dim [0,1]^n \geq n$ $\implies$ $\dim \R^n \geq n$.\\

\begingroup%%----------------------------------->>
\fontsize{9pt}{11pt}\selectfont 
The equivalence $\dim X \geq n$ iff there exists a universal map $f:X \ra [0,1]^n$ holds for any completely regular 
\mX.\\
\endgroup %%------------------------------------<<

\label{19.38}
\label{19.43}
\label{19.56} %dmc mnft
\textbf{\small LEMMA} \  
Let \mA be a closed subset of \mX.  Suppose that $\dim B \leq n$ for every closed subset $B \subset X$ which does not meet 
\mA $-$then each $f \in C(A,\bS^n)$ has an extension $F \in C(X,\bS^n)$.

[Choose an open $U \supset A$ and a $\phi \in C(U,\bS^n)$ such that $\restr{\phi}{A} = f$.  
Choose an open 
$V:A \subset V \subset \ov{V} \subset U$ $-$then $\ov{V} - V$ is closed in $X - V$, so Alexandroff's criterion says there exists a $\Phi \in C(X - V,\bS^n)$ : $\restr{\Phi}{\ov{V}- V} = \restr{\phi}{\ov{V}- V}$.  \ 
Consider the function $F \in C(X,\bS^n)$ defined by 
$
F(x) = 
\begin{cases}
\ \phi(x) \quad (x \in \ov{V})\\
\ \Phi(x) \quad (x \in X - V)
\end{cases}
.]
$
\\
\vspace{0.25cm}

\index{Control Lemma}
\textbf{\small CONTROL LEMMA} \quad 
Let \mA be a closed subset of \mX.  Suppose that $\dim A \leq n$ and that $\dim B \leq n$ for every closed subset $B \subset X$ which does not meet \mA $-$then $\dim X \leq n$.

[Fix a closed subset $A_0 \subset X$ and take an $f_0 \in C(A_0,\bS^n)$.  
Claim: $f_0$ has an extension $f \in C(A \cup A_0,\bS^n)$.  Assuming that $A \cap A_0 \neq \emptyset$, in view of Alexandroff's criterion, the restriction $\restr{f_0}{A \cap A_0}$ has an extension $F_0 \in C(A,\bS^n)$.  Define $f \in C(A \cup A_0,\bS^n)$ piecewise: 
$
\begin{cases}
\ \restr{f}{A} = F_0\\
\ \restr{f}{A_0} = f_0
\end{cases}
. \ 
$
Now let \mB be a closed subset of \mX disjoint from $A \cup A_0$.  
By hypothesis, $\dim B \leq n$ so the lemma implies that $f$ has an extension $F \in C(X,\bS^n)$.  
But $\restr{F}{A_0} = f_0$.  
Invoke Alexandroff's criterion to conclude that $\dim X \leq n$.]\\

\label{20.3}
Suppose that $A \subset X$ is closed $-$then the quotient $X/A$ is a normal Hausdorff space and it follows from the control lemma that $\dim X = \max\{\dim A,\dim X/A\}$.

%%----------------------------------------------------------------------------------------------19
[Note: \ If \mA is a closed $G_\delta$, then $X - A$ is an open $F_\sigma$, thus is normal, and 
$\dim X/A = \dim(X - A)$.]\\

\begingroup%%----------------------------------->>
\fontsize{9pt}{11pt}\selectfont 
The position of quotients in the completely regular situation is complicated by the fact that $X/A$ need not be completely regular even under favorable circumstances, e.g., when \mA has the EP w.r.t. $[0,1]$ or \mA is closed.  Still, 
$\dim X/A$ is meaningful (cf. p. \pageref{19.13}) and nothing more than that is really needed.
\vspi
Given a nonempty $A \subset X$, write $*_A$ for the image of \mA under the projection $p: X \ra X/A$.\\
\endgroup %%------------------------------------<<

\label{19.44}
\begingroup%%----------------------------------->>
\fontsize{9pt}{11pt}\selectfont 
\textbf{\small LEMMA} \  
Let \mX be a nonempty CRH space.  Suppose that \mA is a nonempty subspace of \mX $-$then $\dim X/A \leq \dim X$.
\vspi
[Assume that $\dim X \leq n$.  Take a finite cozero set covering $\sU = \{U_1, \ldots, U_k\}$ of $X/A$.  Choose a continuous function $\phi:X/A \ra [0,1]$ such that $\phi^{-1}(]0,1]) = \ds\bigcap\limits_i \{U_i: *_A \in U_i\}$.  Let 
$q = \phi(*_A)$.  Put 
$V_0 = \{x: \phi(x) > q/2\}$, 
$V_i = U_i - \{x: \phi(x) \geq q\}$ $(i > 0)$ $-$then $\sV = \{V_0, \ldots, V_k\}$ is a finite cozero set refinement of $\sU$ and 
$*_A \notin V_i$ $(i > 0)$.  The collection 
$p^{-1}(\sV) = \{p^{-1}(V_0), \ldots, p^{-1}(V_k)\}$ is a finite cozero set covering of \mX, hence has a precise cozero set refinement $\sW = \{W_0, \ldots W_k\}$ of order $\leq n+1$, which in turn has a precise cozero set refinement 
$\sZ = \{Z_0, \ldots, Z_k\}$ of order $\leq n+1$.  Since $Z_i$ and $X - W_i$ are disjoint zero sets, there exists a continuous function $\phi_i:X \ra [0,1]$ with  
$
\begin{cases}
\ \restr{\phi_i}{Z_i} = 1\\
\ \restr{\phi_i}{X - W_i} = 0
\end{cases}
. \ 
$
But $A \subset Z_0$ and $A \cap W_i = \emptyset$ $(i > 0)$.  
Therefore each $\phi_i$ factors through $X/A$ to give a continuous function $\psi_i: X/A \ra [0,1]$.  
The collection $\{\psi_i^{-1}(]0,1])\}$ is a finite cozero set refinement of $\sU$ of order $\leq n + 1$.]\\
\endgroup %%------------------------------------<<

\begingroup%%----------------------------------->>
\fontsize{9pt}{11pt}\selectfont 
\textbf{\small LEMMA} \  
Let \mX be a nonempty CRH space.  Suppose that \mA is a nonempty subspace of \mX which has the EP w.r.t. $[0,1]$ $-$then 
$\dim X = \max\{\dim A,\dim X/A\}$.
\vspi
[The point here is that every finite cozero set covering of \mA is refined by the restriction to \mA of a finite cozero set covering of \mX (cf. $\S 6$, Proposition 4).]\\
\endgroup %%------------------------------------<<

\begingroup%%----------------------------------->>
\fontsize{9pt}{11pt}\selectfont 
The relation $\dim X = \max\{\dim A, \dim X/A\}$ need not hold if \mA is merely $\sZ$-embedded in \mX.  Indeed, 
Pol\footnote[2]{\textit{Fund. Math.} \textbf{102} (1979), 29-43.}
has constructed an example of a completely regular \mX having the following properties:
(i) $\dim X > 0$; 
(ii) $X = X_1 \cup X_2$, where $X_1$ and $X_2$ are zero sets with 
$
\begin{cases}
\ \dim X_1 = 0\\
\ \dim X_2 = 0
\end{cases}
;\ 
$
(iii) 
$
\begin{cases}
\ X_1 = U_1 \cup D\\
\ X_2 = U_2 \cup D
\end{cases}
,\ 
$
where $U_1$ and $U_2$ are cozero sets and \mD is discrete; 
(iv) $U_1 \cup U_2$ is a countable dense subset of \mX.  Consider $A = U_1 \cup U_2$.\\
\endgroup %%------------------------------------<<

\begin{proposition} \ %07
Suppose that $X = Y \cup Z$, where \mY and \mZ are normal $-$then $\dim X \leq \dim Y + \dim Z + 1$.
\end{proposition}

[There is nothing to prove if either $\dim Y = \infty$ or $\dim Z = \infty$, so assume that $\dim Y \leq r$ and 
$\dim Z \leq s$.  Owing to the control lemma, it will be enough to show that 
%%----------------------------------------------------------------------------------------------20
$\dim \ov{Y} \leq r + s + 1$.  Let $\sU = \{U_i\}$ be a finite open covering of $\ov{Y}$.  
Since $\dim Y \leq r$, there exists a collection $\sV = \{V_i\}$ of open subsets of $\ov{Y}$ such that $V_i \subset U_i$, $Y \subset \bigcup\limits_i V_i$, and 
$\ord(\{Y \cap V_i\}) \leq r + 1$.  
Put $D = \ov{Y} - \bigcup\limits_i V_i$.  
Because $\dim D \leq s$, there exists a closed covering $\sA = \{A_i\}$ of \mD of order $\leq s + 1$ such that $A_i \subset U_i$.  Without changing the order, expand 
$\sA$ to a collection $\sW = \{W_i\}$ of open subsets of $\ov{Y}$ such that $A_i \subset W_i \subset U_i$.  
The union 
$\sV \cup \sW$ covers $\ov{Y}$, refines $\sU$, and is of order $\leq r + 1 + s + 1$.]

[Note: \ When \mX is metrizable, there is another way to argue.  
Assume: 
$
\begin{cases}
\ \dim Y = r\\
\ \dim Z = s
\end{cases}
$ $-$then every closed subset of 
$
\begin{cases}
\ Y\\
\ Z
\end{cases}
$
has the EP w.r.t. 
$
\begin{cases}
\ \bS^s\\
\ \bS^r
\end{cases}
$
, thus every closed subset of \mX has the EP w.r.t. $\bS^s*\bS^r = \bS^{s+r+1}$
(cf. p. \pageref{19.14}).]\\

By way of application, suppose that \mX is hereditarily normal and that $X = \bigcup\limits_0^n X_i$, where 
$\forall \ i$: $\dim X_i \leq 0$ $-$then $\dim X \leq n$.

This remark can be used to prove that $\dim \R^n \leq n$, from which $\dim \R^n = n$.  
(cf. p. \pageref{19.15}).  
Thus suppose that $n \geq 1$ and that $0 \leq m \leq n$.  
Denote by $\Q_m^n$ the subspace of $\R^n$ consisting of all points with exactly $m$ rational coordinates $-$then 
$\R^n = \Q_0^n \cup \cdots \cup \Q_n^n$.  
Claim: $\forall \ m$, $\dim \Q_m^n = 0$.  
This is immediate if $m = n$ (cf. Proposition 2), so assume that $m < n$.
For any choice of $m$ distinct natural numbers $i_1, \ldots, i_m$, each $\leq n$, and any choice of $m$ rational numbers 
$r_1, \ldots, r_m$, the space $\prod\limits_{i=1}^n R_i$, where $R_{i_j} = \{r_j\}$ for $j = 1, \ldots, m$ and 
$R_i = \R$ for $i \neq i_j$, is a closed subspace of $\R^n$.  
Therefore
$\Q_m^n \cap \prod\limits_{i=1}^n R_i$ is a closed subspace of $\Q_m^n$.  
On the other hand , 
$\Q_m^n \cap \prod\limits_{i=1}^n R_i$ is homeomorphic to the subspace of $\R^{n-m}$ consisting of all points with irrational coordinates, hence 
$\dim \bigl( \Q_m^n \cap \prod\limits_{i=1}^n R_i\bigr) = 0$ (cf. Proposition 2).  
Since the collection of all sets of the form 
$\Q_m^n \ \cap \  \prod\limits_{i=1}^n R_i$ is a countable closed covering of $\Q_m^n$, 
the countable union lemma implies that 
$\dim \Q_m^n = 0$.\\

\index{Theorem: Fundamental Theorem of Dimension Theory}
\index{Fundamental Theorem of Dimension Theory}
\textbf{\small FUNDAMENTAL THEOREM OF DIMENSION THEORY} \  
The topological dimension of $\R^n$ is exactly $n$.\\

One consequence is the evaluation $\dim[0,1]^n = n$.  Corollary: Take $X = \bS^n$ $-$then $\dim X = n$.  In fact, 
$X = X_1 \cup X_2$, where $X_1$ and $X_2$ are closed and homeomorphic to $[0,1]^n$.

Another consequence is the evaluation
$
\begin{cases}
\ \dim(\Q_0^n \cup \cdots \cup \Q_m^n) = m\\
\ \dim(\Q_m^n \cup \cdots \cup \Q_n^n) = n-m
\end{cases}
$
.\\
\vspace{0.25cm}

%%----------------------------------------------------------------------------------------------21
\begingroup%%----------------------------------->>
\fontsize{9pt}{11pt}\selectfont 
\textbf{\small EXAMPLE}  \ [Assume CH] \  
Take $X = [0,1]^n$ $-$then the topological dimension of \mX in any Kunen modification of its euclidean topology is $n$-1 
(cf. p. \pageref{19.16}).\\
\endgroup %%------------------------------------<<

\begingroup%%----------------------------------->>
\fontsize{9pt}{11pt}\selectfont 
\textbf{\small FACT} \  
Let \mX and \mY be normal.  Let $A \ra X$ be a closed embedding and let $f:A \ra Y$ be a continuous function.  
Assume: $\dim X \leq n$ $\&$ $\dim Y \leq n$ $-$then $\dim (X \sqcup_f Y) \leq n$.
\vspi
[Use the control lemma $(X \sqcup_f Y$ is a normal Hausdorff space (cf. p. \pageref{19.17})).]\\
\endgroup %%------------------------------------<<

\begingroup%%----------------------------------->>
\fontsize{9pt}{11pt}\selectfont 
Application: If \mX is obtained from a normal \mA by attaching $n$-cells, then $\dim X = n$ provided that 
$\dim A \leq n$ and the index set is not empty.
\vspi
[\mX contains an embedded copy of $\bB^n$ which is strongly paracompact, thus a priori, $\dim X \geq n$ 
(cf. p. \pageref{19.18}).]\\
\endgroup %%------------------------------------<<

\label{5.0b}
\label{19.27}
\label{20.11}
\label{19.571}
\index{CW complexes (example)} 
\begingroup%%----------------------------------->>
\fontsize{9pt}{11pt}\selectfont 
\textbf{\small EXAMPLE \  (\un{CW Complexes})} \  
Let \mX be a CW complex $-$then by the countable union lemma, $\dim X = \sup\dim X^{(n)}$ and $\forall \ n$, 
$\dim X^{(n)} \leq n$.  Therefore the combinatorial dimension of \mX is equal to the topological dimension of \mX.\\
\endgroup %%------------------------------------<<

\begingroup%%----------------------------------->>
\fontsize{9pt}{11pt}\selectfont 
\textbf{\small FACT} \  
Suppose that \mX is normal.  Let $\sA = \{A_j:j \in J\}$ be an absolute closure preserving closed covering of \mX such that 
$\forall \ j$, $\dim A_j \leq n$ $-$then $\dim X \leq n$, hence $\dim X = \sup\dim A_j$.
\vspi
[Use Alexandroff's criterion.  Let \mA be a closed subset of \mX, take an $f \in C(A,\bS^n)$, and let $\sF$ be the set of continuous functions \mF that are extensions of $f$ and have domains of the form $A \cup X_I$, where 
$X_I = \bigcup\limits_i A_i$ $(I \subset J)$.  Order $\sF$ by writing $F^\prime \leq F\pp$ iff $F\pp$ is an extension of 
$F^\prime$.  Every chain in $\sF$ has an upper bound, so by Zorn, $\sF$ has a maximal element $F_0$.  But the domain of 
$F_0$ is necessarily all of \mX and $\restr{F_0}{A} = f$.]\\
\endgroup %%------------------------------------<<

\begingroup%%----------------------------------->>
\fontsize{9pt}{11pt}\selectfont 
\textbf{\small EXAMPLE \ (\un{Vertex Schemes})} \ 
\index{Vertex Schemes (example)}  
Let $K = (V,\Sigma)$ be a vertex scheme $-$then one can attach to \mK its combinatorial dimension $\dim K$, as well as the topological dimension of $\abs{K}$ (Whitehead topology) and $\abs{K}_b$ (barycentric topology).  
The claim is that these are all equal.  
Note that in any event, if $\sigma$ is an $n$-simplex of \mK, 
then $\dim \abs{\sigma} = n$, so, $\abs{\sigma}$ being a closed subspace of both $\abs{K}$ and $\abs{K}_b$, 
$
\begin{cases}
\ \dim \abs{K} \geq \dim K\\
\ \dim \abs{K}_b \geq \dim K
\end{cases}
. \ 
$
Regarding the inequalities in the opposite direction, first observe that 
$\{\abs{\sigma} \}$ is an absolute closure preserving closed covering of 
$\abs{K}$, thus in this case the preceding result is immediately applicable.  
Turning to $\abs{K}_b$, $\{\abs{\sigma}\}$ is still closure preserving.  
To exploit this, consider the $n$-skeleton $K^{(n)}$.  
Assertion: $\forall \ n$, 
$\dim |K^{(n)}|_b \leq n$.  
Obviously, $\dim |K^{(0)}|_b = 0$.  
Suppose that $n \geq 1$ and 
$\dim |K^{(n-1)}|_b \leq n-1$.  
Let $\Sigma_n$ be the set of $n$-simplexes of \mK.  
The collection 
$\{\langle \sigma \rangle: \sigma \in \Sigma_n\}$ is an open covering of 
$|K^{(n)}|_b - |K^{(n-1)}|_b$.  
Write $\langle \sigma \rangle = \ds\bigcup\limits_j A_{\sigma j}$, where the 
$A_{\sigma j} \subset \abs{\sigma}$ are compact.  
The collection $\{A_{\sigma j}: \sigma \in \Sigma_n\}$ is discrete.
Let $A_j$ be its union $-$then $\dim A_j \leq n$.  
Finish the induction via the countable union lemma: 
$|K^{(n)}|_b = |K^{(n-1)}|_b \cup \ds\bigcup\limits_j A_j$.
\vspi
%%----------------------------------------------------------------------------------------------22
[Note: \ It is therefore a corollary that the combinatorial dimension of $\abs{K}$ viewed as a CW complex is equal to 
$\dim K$.]\\
\endgroup %%------------------------------------<<

\label{20.12}
\begingroup%%----------------------------------->>
\fontsize{9pt}{11pt}\selectfont 
Let \mX be an $n$-manifold.  Since compact subsets of a nonempty CRH space have the EP w.r.t. $[0,1]$ and since \mX contains a compact subset homeomorphic to $[0,1]^n$, of necessity $\dim X \geq n$, the euclidean dimension of \mX.  
To reverse the inequality $\dim X \geq n$ when \mX is paracompact or, equivalently, metrizable (cf. $\S 1$, Proposition 11), one can assume that \mX is connected.  
But then \mX is second countable (cf. p. \pageref{19.19}), thus admits a covering by a countable collection of closed sets, each of topological dimension $n$, so $\dim X \leq n$.
\vspi
[Note: \ Using the combinatorial principal \ $\Diamond$, 
Fedorchuk\footnote[2]{\textit{Topology Appl.} \textbf{54} (1993), 221-239; 
see also \textit{Math. Sbornik} \textbf{186} (1995), 151-162.}
has constructed a perfectly normal $n$-manifold \mX such that $n < \dim X$.]\\
\endgroup %%------------------------------------<<

\textbf{\small LEMMA} \  
$\R^n$ is homogeneous with repect to countable dense subsets, i.e., if \mA and \mB are two countable dense subsets of 
$\R^n$, then there exists a homeomorphism $f:\R^n \ra \R^n$ such that $f(A) = B$.\\

\begin{proposition} \ %08
Let \mX be a subspace of $\R^n$ $-$then $\dim X = n$ iff \mX has a nonempty interior.
\end{proposition}

[Suppose that the interior of \mX is empty.  
Since $\R^n - X$ is dense in $\R^n$, there exists a countable set 
$A \subset \R^n - X$: $\ov{A} = \R^n$.  
Choose a homeomorphism $f:\R^n \ra \R^n$ such that $f(A) = \Q_n^n$ $-$then 
$f(X) \subset \bigcup\limits_{m<n} \Q_m^n$, which gives $\dim X \leq n - 1$.]\\

It follows from this result that if \mX is a subspace of $[0,1]^n$ or $\bS^n$, then $\dim X = n$ iff \mX has a nonempty interior.\\

\textbf{\small SUBLEMMA} \quad 
Suppose that \mX is  Lindel\"of.  Let $\sO = \{O\}$ be a basis for \mX 
$-$then for every pair $(A,B)$ of disjoint closed subsets of \mX there exists an open set 
$P \subset X$ and a sequence $\{O_j\} \subset \sO$ such that 
$A \subset P \subset \ov{P} \subset X - B$ and $\fr P \subset \bigcup\limits_j \fr O_j$.

[Given $x \in X$, choose a neighborhood $O_x \in \sO$ of $x$ such that either $A \cap \ov{O}_x = \emptyset$ or 
$B \cap \ov{O}_x = \emptyset$.  
Let $\{O_j\}$ be a countable subcover of $\{O_x\}$.  Divide $\{O_j\}$ into two subcollections 
$\{O_j^\prime\}$ and $\{O_j\pp\}$ according to whether $\ov{O}_j$ does or does not meet $A$.  Put 
$
\begin{cases}
\ P_i = {O}_i^\prime - \bigcup\limits_{j < i} \ov{O}_j\pp\\
\ Q_i = {O}_i\pp - \bigcup\limits_{j \leq i} \ov{O}_j^\prime \
\end{cases}
$ $-$then 
$
\begin{cases}
\ P = \bigcup\limits_i P_i\\
\ Q = \bigcup\limits_i Q_i
\end{cases}
$
are disjoint open subsets of \mX and 
$A \subset$  
$P \subset$  
$\ov{P} \subset X - B$ 
with $\fr P \subset X - (P \cup Q)$.  Let $x \in X - (P \cup Q)$.  
Denote by \mS the first element of the
%%----------------------------------------------------------------------------------------------23
sequence $\ov{O_1^\prime},\ov{O_1\pp},\ov{O_2^\prime}, \ov{O_2\pp},\ldots$ that contains $x$.  
If $S = \ov{O_i^\prime}$, 
then $x \notin P_i$ and $x \notin \ov{O_i\pp}$ $(j < i)$, so $x \in \fr O_i^\prime$; 
if $S = \ov{O_i\pp}$, 
then $x \notin Q_i$ and $x \notin \ov{O_j^\prime}$ $(j \leq i)$, so $x \in \fr O_i\pp$.  
Therefore 
$x \in \bigcup\limits_i \fr O_i^\prime \cup \bigcup\limits_i \fr O_i\pp$ or still, $x \in \bigcup\limits_j \fr O_j$.]\\

\textbf{\small LEMMA} \  Suppose that \mX is  Lindel\"of.  Let $\sO = \{O\}$ be a basis for \mX such that $\forall \ O$: 
$\dim\fr O \leq n - 1$ $-$then $\dim X \leq n$.

[Let $\sU = \{U_i\}$ be a finite open covering of \mX; let $\sA = \{A_i\}$ be a precise closed refinement of $\sU$.  Use the sublemma and for each $i$, choose an open set $P_i \subset X$ and a sequence $\{O_{i,j}\} \subset \sO$: 
$A_i \subset P_i \subset \ov{P}_i \subset U_i$ and $\fr P_i \subset \bigcup\limits_j \fr O_{i,j}$.  Put 
$D = \ds\bigcup\limits_i \fr P_i$.  
The countable union lemma implies that 
$\dim D \leq n - 1$, so there exists a collection $\sV = \{V_i\}$ of open subsets of \mX such that 
$\ov{V}_i \subset U_i$, 
$D \subset \bigcup\limits_i V_i$, and $\ord(\{\ov{V}_i\}) \leq n$.  
Write $B_i$ in place of 
$\ov{P}_i - \bigl(\cup \sV \cup \bigcup\limits_{j < i} P_j\bigr)$. 
 Since the $B_i$ are pairwise disjoint, it follows that the collection $\{B_i\} \cup \{\ov{V}_i\}$ is a finite closed refinement of $\sU$ of order $\leq n+1$.]\\

\begin{proposition} \ %09
Let \mU be a nonempty, nondense open subset of $\R^n$ $-$then $\dim\fr U = n-1$.
\end{proposition}

[Suppose that \mU is bounded.  
In this case, \mU has a basis consisting of sets homeomorphic to itself, so if 
$\dim \fr U < n-1$, then by the lemma, $\dim U \leq n - 1$, a contradiction.

Suppose that \mU is not bounded.  
Fix a point $x$ in the interior of the complement of \mU and choose an open ball \mB centered at $x$ which is entirely contained therein.  
The associated inversion $\R^n - \{x\} \ra \R^n - \{x\}$ carries \mU onto a nonempty open set $O \subset B$.  
Obviously, $\fr O - \{x\}$ is homeomorphic to $\fr U$.  
On the other hand, by the above, 
$\dim \fr O = n - 1$.  So, from the control lemma, $\dim \fr U = n - 1$.]\\

\textbf{\small LEMMA} \  
The following conditions are equivalent.

\indent\indent (1) \mX can be disconnected by a closed subset of topological dimension $\leq n$.

\indent\indent (2) \mX contains a nonempty, nondense open subset whose frontier has topological dimension $\leq n$.

\indent\indent (3) $X = A \cup B$, where \mA and \mB are closed proper subsets of \mX such that $\dim(A \cap B) \leq n$.\\

Take $X = \R^n$ $-$then, in view of Proposition 9, $\R^n$ cannot be disconnected by a closed subset of topological dimension $\leq n - 2$.  Ths same is true of $[0,1]^n$ and of $\bS^n$.\\

%%----------------------------------------------------------------------------------------------24
\label{19.11}
\begingroup%%----------------------------------->>
\fontsize{9pt}{11pt}\selectfont 
Let \mX be a LCH space.  Suppose that \mX is connected and locally connected $-$then \mX is said to be 
\un{$n$-solid}
\index{n-solid} 
$(n \geq 1)$ if for every $x \in X$ and for every neighborhood \mU of $x$ there is a connected relatively compact neighborhood \mV of $x$ such that 
$
\begin{cases}
\ \ov{V} \subset U\\
\ \dim \ov{V} \geq n
\end{cases}
$
and $\ov{V}$ cannot be disconnected by a closed subset of topological dimension $\leq n - 2$.  
Examples: $\R^n$, $[0,1]^n$, and $\bS^n$ are $n$-solid.
\vspi
[Note: \ A LCH space \mX that is both connected and locally connected is necessarily 1-solid.  
Specialization of the argument infra then leads to the conclusion that \mX does not admit a disjoint decomposition 
$\ds\bigcup\limits_1^\infty A_j$ by nonempty closed proper subspaces $A_j$.  
If \mX is compact, then the assumption of local connectedness is unnecessary but simple examples show that it is not superfluous in general.]\\
\endgroup %%------------------------------------<<

\begingroup%%----------------------------------->>
\fontsize{9pt}{11pt}\selectfont 
\textbf{\small FACT} \  
Suppose that \mX is $n$-solid and perfectly normal $-$then \mX cannot be covered by a sequence $\{A_j\}$ of nonempty closed proper subsets such that $\forall \ i \neq j$: $\dim (A_i \cap A_j) \leq n - 2$.
\vspi
[Proceed by contradiction, so $X = \ds\bigcup\limits_1^\infty A_j$, 
where the $A_j$ satisfy the conditions set forth above.  
Claim: There exists a sequence $\{x_0, x_1, \ldots\} \subset X$ subject to: 
(1) $x_i \in V_i$, $V_i$ as in the definition of ``$n$-solid''; 
(2) $\forall \ j$: $\ov{V}_i \not\subset A_j$; 
(3) $\ov{V}_i \subset \ov{V}_{i-1}$; 
(4) $\ov{V}_i \cap A_i = \emptyset$.  
Here
$
\begin{cases}
\ V_{-1} = X\\
\ A_0 = \emptyset
\end{cases}
. \ 
$
Granted the claim, $\ds\bigcap\limits_0^\infty \ov{V}_i = \emptyset$, an impossibility.  
The $x_i$ can be constructed inductively.  
Start by fixing an index $j_0$ such that the interior of $A_{j_0}$ is not empty (Baire).  
Choose a point $x_0$ in the frontier of the interior of $A_{j_0}$ and take a neighborhood $V_0$ of $x_0$ as in the definition of ``$n$-solid'' $-$then the pair 
$(x_0,V_0)$ satisfies (1)-(4).  
Given $x_i$ and $V_i$ $(i > 0)$, look at a component \mY of $V_i - A_{i+1}$.  
Show that \mY is not a subset of any $A_j$ and then get $x_{i+1}$ and $V_{i+1}$ by repeating the process used to get $x_0$ and $V_0$.]
\vspi
[Note: \ Proposition 5 is tacitly used at several points.  
When $n = 1$, the assumption of perfect normality plays no role, hence can be dropped.]\\
\endgroup %%------------------------------------<<

\textbf{\small LEMMA} \  
Let \mX be a closed subspace of $\R^n$; let $x \in X$ $-$then $x$ belongs to the frontier of \mX iff 
$x$ has a neighborhood basis 
$\{U\}$ in \mX such that $\forall \ U$: $X = U$ has the EP w.r.t. $\bS^{n-1}$.

[Necessity: Let $x$ be an element of the frontier of \mX.  
Assuming that $x$ is the origin, put $U  = X \cap \epsilon \bB^n$ 
$(\epsilon > 0)$.  
To simplify, take $\epsilon = 1$.  
Fix a point $x_0 \in \bB^n - X$ and write $r_0$ for the radial retraction 
$\bD^n - \{x_0\} \ra \bS^{n-1}$.  
Choose an $f \in C(X - U,\bS^{n-1})$.  
Since $A = (X - U) \cap \bS^{n-1}$ is a closed subset of 
$\bS^{n-1}$, Alexandroff's criterion implies that $\restr{f}{A}$ can be extended to a continuous function 
$g:\bS^{n-1} \ra \bS^{n-1}$.  The function $F:X \ra \bS^{n-1}$ defined by 
$
\begin{cases}
\ \restr{F}{X-U} = f\\
\ \restr{F}{U} = g \circx r_0
\end{cases}
$
is then a continuous extension of $f$ to \mX.

Sufficiency: Let $x$ be an element of the interior of \mX.  
Assuming that $x$ is the origin, fix an $\epsilon > 0$: 
$\epsilon \bD^n \subset X$.  Let \mU be a neighborhood of $x$ in \mX: $U \subset \epsilon\bB^n$ $-$then the claim is that there exists an $f \in C(X - U,\bS^{n-1})$ that has no extension $F \in C(X,\bS^{n-1})$.  
To see this, identify the frontier of 
$\epsilon \bD^n$ with $\bS^{n-1}$ and consider the projection $X - U \ra \bS^{n-1}$
%%----------------------------------------------------------------------------------------------25
determined by $x$ which, if extendible, would lead to a retraction of $\epsilon \bD^n$ onto its frontier.]\\

Let \mX and \mY be closed subspaces of $\R^n$ $-$then the characterization provided by the lemma tells us that any homeomorphism $f:X \ra Y$ necessarily carries the frontier of \mX onto the frontier of \mY.\\

\index{Theorem: Theorem of Invariance of Domain}
\index{Theorem of Invariance of Domain}
\textbf{\small THEOREM OF INVARIANCE OF DOMAIN} \ \  
Let \mU be an open subsest of $\R^n$ $-$then every continuous injective map $U \ra \R^n$ is an open embedding.\\

\begingroup%%----------------------------------->>
\fontsize{9pt}{11pt}\selectfont 
This result does not extend to an infinite dimensional normed linear space \mX.  
Indeed, for such an \mX, there always exists an embedding $f:X \ra X$ that is not open and there always exists a bijective continuous map $f:X \ra X$ that is not a homeomorphism 
(van Mill\footnote[2]{\textit{Proc. Amer. Math. Soc.} \textbf{101} (1987), 173-180.}).\\
\endgroup %%------------------------------------<<

\begingroup%%----------------------------------->>
\fontsize{9pt}{11pt}\selectfont 
\textbf{\small FACT} \  
Let $f:\R^n \ra \R^n$ be continuous and locally one-to-one.  
Assume that $\norm{f(x)} \ra \infty$ as $\norm{x} \ra \infty$ 
$-$then $f(\R^n) = \R^n$.\\
\endgroup %%------------------------------------<<

Let \mX and \mY be $n$-manifolds; let 
$
\begin{cases}
\ U \subset X\\
\ V \subset Y
\end{cases}
$
and suppose that $f:U \ra V$ is a homeomorphism $-$then from the domain invariance of $\R^n$, \mU open in \mX 
$\implies$ \mV open in \mY.  
\label{1.13}
Corollary: Homeomorphic topological manifolds have the same euclidean dimension.\\

\label{19.25}
\begingroup%%----------------------------------->>
\fontsize{9pt}{11pt}\selectfont 
Let \mX be a CRH space.  
Suppose that $\dim X = n$ $(n \geq 1)$ $-$then \mX is said to be a 
\un{Cantor $n$-space}
\index{Cantor $n$-space} 
if \mX cannot be disconnected by a closed subset of topological dimension 
$\le n - 2$.  Since $\dim \emptyset = -1$, a Cantor $n$-space is necessarily connected.  
For example $\R^n$ is a Cantor $n$-space.  
So too are $[0,1]^n$ and $\bS^n$.  
The tubular arrangement
\[
\bigcup\limits_1^\infty \left(\left[-\frac{1}{2n-1},-\frac{1}{2n}\right] \times [-1,1]\right) 
\cup
\bigcup\limits_1^\infty \left(\left[-\frac{1}{2n},-\frac{1}{2n+1}\right] \times \left[-\frac{1}{n},\frac{1}{n}\right]\right) 
\cup
([0,1] \times [-1,1])
\]
is a Cantor 2-space.  It remains connected after removal of the origin but what's left is no longer path connected.\\
\endgroup %%------------------------------------<<

\begingroup%%----------------------------------->>
\fontsize{9pt}{11pt}\selectfont 
\textbf{\small FACT} \  
Suppose that \mX is compact, with $\dim X = n$ $(n \geq 1)$ $-$then \mX contains a Cantor $n$-space, thus \mX has a component of topological dimension $n$.
\vspi
[There exists a closed subset $A \subset X$ and a continuous function $f:A \ra \bS^{n-1}$ that has no continuous extension $F:X \ra \bS^{n-1}$.  Use Zorn and construct a closed subset $B_f \subset X$ such that 
(i) $f$ does not have a continuous extension to $A \cup B_f$ and
(ii) $f$ does have a continuous extension to $A \cup B$ for each closed
%%----------------------------------------------------------------------------------------------26
proper subset \mB of $B_f$.  In view of condition (i), $\dim B_f = n$.  
Claim: $B_f$ is a Cantor $n$-space.  
Assume not and write 
$B_f = B^\prime \cup B\pp$, where $B^\prime$ and $B\pp$ are proper closed subsets of $B_f$ with 
$\dim (B^\prime \cap B\pp) \leq n - 2.$  On account of condition (ii), $f$ has a continuous extension 
$
\begin{cases}
\ f^\prime \text{ to } A \cup B^\prime\\
\ f\pp \text{ to } A \cup B\pp
\end{cases}
. \ 
$
Therefore $f$ has a continuous extension to $A \cup B_f$ (cf. Proposition 15).  Contradiction.]
\vspi
[Note: \ One cannot expect in general that a noncompact \mX will contain a compact Cantor $n$-space.  Reason: For each 
$n \geq 1$, there exists a zero dimensional \mX of topological dimension $n$ 
(consider an ``$n$-dimensional'' variant of Dowker's Example ``M''.]\\
\endgroup %%------------------------------------<<

\begingroup%%----------------------------------->>
\fontsize{9pt}{11pt}\selectfont 
Suppose that \mX is compact and perfectly normal, with $\dim X = n$ $(n \geq 1)$.  
Denote by $C_X$ the union of all Cantor $n$-spaces in \mX $-$then $\dim (X - C_X) \leq \dim X$ but if $n > 1$ equality can obtain even when \mX is metrizable 
(Pol\footnote[2]{\textit{Fund. Math.} \textbf{136} (1990), 127-131.}
).\\
\endgroup %%------------------------------------<<

\begingroup%%----------------------------------->>
\fontsize{9pt}{11pt}\selectfont 
\textbf{\small FACT} \  
Suppose that \mX is a compact connected homogeneous ANR of topological dimension 
$n \geq 1$ $-$then \mX is a Cantor $n$-space.
\vspi
[Note: \ Is such an \mX actually an $n$-manifold?  This is true if $n = 1$ or 2 
(Bing-Borsuk\footnote[3]{\textit{Ann. of Math.} \textbf{81} (1965), 100-111.}) but is a mystery if $n > 2$.  The three dimensional case is related to the Poincar\'e conjecture 
(Jakobsche\footnote[6]{\textit{Fund. Math.} \textbf{106} (1980), 127-134.}).]\\
\endgroup %%------------------------------------<<

\index{Marde\u si\' c Factorization Lemma}
\textbf{\small MARDE\u SI\' C FACTORIZATION LEMMA} \  
Let \mX and \mY be compact Hausdorff spaces $-$then for every $f \in C(X,Y)$ there exists a compact Hausdorff space \mZ with 
$
\begin{cases}
\dim Z \leq \dim X\\
\wt Z \leq \wt Y
\end{cases}
$
and functions 
$
\begin{cases}
\ g \in C(X,Z)\\
\ h \in C(Z,Y)
\end{cases}
$
such that $f = h \circx g$ and $g(X) = Z$.

[Assume that $\dim X = n$ is finite and $\wt Y \geq \omega$.  
Fix a basis $\sV$ for \mY of cardinality $\wt Y$.  Denote by \bV the collection of all finite open coverings of \mY made up of members of $\sV$ and put $\bU_0 = f^{-1}(\bV)$.  
Inductively define a sequence $\bU_1, \bU_2, \ldots$ of collections of finite open coverings of \mX by assigning to each pair 
$
\begin{cases}
\ \sU^\prime\\
\ \sU\pp
\end{cases}
\in \bU_{i - 1}
$ 
a finite open covering $\sU$ of \mX of order $\leq n + 1$ that is a star refinement of both $\sU^\prime$ and $\sU\pp$ and write $\bU_i$ for $\{\sU\}$.  The declaration $x \sim y$ iff 
$y \in [x] \equiv \bigcap\limits_1^\infty \bigcap\{\st(x,\sU):\sU \in \bU_i\}$ is an equivalence relation on \mX and for any open set $U \subset X$ and any $[x] \subset U$, $\exists$ $\sU_x \in \bU_{i_x}$:
\[
[x] \subset \st(x,\sU_x) \subset \bigcup\limits_{\st(x,\sU_x)} [y] \subset \st(\st(x,\sU_x),\sU_x) \subset U.
\]
%%----------------------------------------------------------------------------------------------27
Therefore the union of the equivalence classes that are contained in \mU is open in \mX.  
Give $Z = X/\sim$ the quotient topology.  
Since the projection $g:X \ra Z$ is a closed map, \mZ is a compact Hausdorff space.  
By construction, $f$ is constant on equivalence classes so there is a continuous factorization $f = h \circx g$.  
Assign to each $\sU = \{U\}$ in $\bU_i$ the collection 
$\sU^* = \{U^*\}$, where $U^* = Z - g(X - U)$ $-$then $\sU^*$ is a finite open covering of \mZ of order $\leq n + 1$.  Moreover, every finite open covering $\sP = \{P\}$ of \mZ has a refinement of the form $\sU^*$, hence $\dim Z \leq n$.  
In fact, 
$\forall \ x \in X$ $\exists \ P_x \in \sP$: $[x] \subset g^{-1}(P_x)$.  
Choose $\sU_x \in \bU_{i_x}$: 
$O_x \equiv \st(\st(x,\sU_x),\sU_x) \subset g^{-1}(P_x)$.  
Let $\{O_{x_j}\}$ be a finite subcover of $\{O_x\}$.  Take a 
$\sU \in \bU_i$ that refines the $\sU_{x_j}$ and consider the associated $\sU^*$.  
Finally, the collection 
$\bigcup\limits_1^\infty \bigcup \{\sU^*: \sU \in \bU_i\}$ is a basis for \mZ of cardinality $\leq \wt Y$.]\\

\begin{proposition} \ %10
X has a compactification $\Delta X$ such that 
$
\begin{cases}
\ \dim \Delta X \leq \dim X\\
\ \wt \Delta X \leq \wt X
\end{cases}
. \ 
$
\end{proposition}

[Assume that $\wt X \geq \omega$.  Choose an embedding $X \ra [0,1]^{\wt X}$ and denote by $f$ its extension 
$\beta X \ra [0,1]^{\wt X}$.  Apply the Marde\u si\' c factorization lemma to get a compact Hausdorff space $\Delta X$ and functions 
$
\begin{cases}
\ g \in C(\beta X, \Delta X)\\
\ h \in C(\Delta X,[0,1]^{\wt X})
\end{cases}
$
:
$
\begin{cases}
\ \dim \Delta X \leq \dim \beta X = \dim X\\
\ \wt \Delta X \leq \wt [0,1]^{\wt X} = \wt X
\end{cases}
$
and $f = h \circx g$ \ $(g(\beta X) = \Delta X)$.  Look at $\restr{g}{X}$.]\\

\begingroup%%----------------------------------->>
\fontsize{9pt}{11pt}\selectfont 
Since the normality of \mX was not used in the proof, Proposition 10 is true in the completely regular situation.\\
\endgroup %%------------------------------------<<

\begingroup%%----------------------------------->>
\fontsize{9pt}{11pt}\selectfont 
\textbf{\small FACT} \  
For every integer $n \geq 0$ and for every cardinal $\kappa \geq \omega$, there exists a compact Hausdorff space 
$K(n,\kappa)$: 
$
\begin{cases}
\ \dim K(n,\kappa) \leq n\\
\ \wt K(n,\kappa) \leq \kappa
\end{cases}
$
having the property that if \mX is a nonempty CRH space of topological dimension $\leq n$ and weight $\leq \kappa$, then there is an embedding $X \ra K(n,\kappa)$.
\vspi
[Consider the collection $\{X_i:i \in I\}$ of all subspaces $X_i \subset [0,1]^\kappa$, where $\dim X_i \leq n$.  Let $f$ be the natural map $\ds\coprod\limits_i X_i \ra [0,1]^\kappa$.  Work with $\beta f$.]\\
\endgroup %%------------------------------------<<

\begingroup%%----------------------------------->>
\fontsize{9pt}{11pt}\selectfont 
Does every subspace $X \subset \R^n$ have a dimension preserving compactification that embeds in $\R^n$?  \quad
This is an open question.\\
\endgroup %%------------------------------------<<

A set $S \subset \R^n$ is said to be in 
\un{general position}
\index{general position} 
if every subset $T \subset S$ of cardinality $\leq n+1$ is geometrically independent.\\

\textbf{\small LEMMA} \  
$\R^n$ contains a countable dense set in general position.\\

Suppose that \mX is second countable $-$then there is an embedding \ $X \ra \R^\omega$.  \ 
If $\dim X = n$, then one can say more: There is an embedding $X \ra \R^{2n+1}$.

%%----------------------------------------------------------------------------------------------28
Start with an initial reduction: Take \mX compact (cf. Proposition 10). 
 Fix a compatible metric $d$ on \mX.  Attach to each 
$f \in C(X,\R^{2n+1})$ its ``injectivity deviation''
\index{injectivity deviation}
\[
\dev f = \sup \{\diam f^{-1}(p): p \in \R^{2n+1}\}.
\]
Given $\epsilon > 0$, put $D_\epsilon = \{f: \dev f < \epsilon\}$.  Claim: $\forall \ \epsilon > 0$, $D_\epsilon$ is open and dense in $C(X,\R^{2n+1})$.  Admit this $-$then $\bigcap\limits_1^\infty D_{1/k}$ is dense in $C(X,\R^{2n+1})$ (Baire), thus is nonempty.  
But $\bigcap\limits_1^\infty D_{1/k}$ is the set of embeddings $X \ra \R^{2n+1}$.\\
\indent\indent (1) $D_\epsilon$ is open in $C(X,\R^{2n+1})$.  
Proof: Let $f \in D_\epsilon$.  
Choose $r: \dev f < r < \epsilon$.  Set $A_r = \{(x,y): d(x,y) \geq r\}$.  Call $\delta_f$ the minimum of 
$\frac{1}{2}\norm{f(x) - f(y)}$ on $A_r$ $-$then $\{g : \norm{f - g} < \delta_f\} \subset D_\epsilon$.\\
\indent\indent (2) $D_\epsilon$ is dense in $C(X,\R^{2n+1})$. 
 Proof: Fix $f \in C(X,\R^{2n+1})$.  Given $\delta > 0$, let 
$\sU = \{U_i\}$ be a finite open covering of \mX of order $\leq n + 1$: $\forall \ i$, 
$
\begin{cases}
\ \diam U_i < \epsilon/2\\
\ \diam f(U_i) < \delta/2
\end{cases}
$
and denote by $\{\kappa_i\}$ a partition of unity on \mX subordinate to $\sU$.  
Choose a point $x_i \in U_i$ and then choose a point $p_i \in \R^{2n+1}$ within $\delta /2$ of $f(x_i)$, using the lemma to arrange matters so that in addition $\{p_i\}$ is in general position.  
Put $g = \sum\limits_i \kappa_ip_i$ $-$then 
\[
f(x) - g(x) = \sum\limits_i \kappa_i(x)(f(x_i) - p_i) + \sum\limits_i \kappa_i(x)(f(x) - f(x_i)),
\]
hence $\norm{f - g} < \delta$.  
There remains the verification: $g \in D_\epsilon$.  
For this, it need only be shown that if 
$g(x) = g(y)$, then $\exists \ i$: $x, y \in U_i$.  
Consider the relation $\sum\limits_i(\kappa_i(x) - \kappa_i(y)) p_i = 0$.  
Because the order of $\sU$ is $\leq n + 1$, at most $2n + 2$ of these terms are nonzero.  
However, 
$\sum\limits_i(\kappa_i(x) - \kappa_i(y)) = 0$, from which $\kappa_i(x) - \kappa_i(y) = 0$ $\forall \ i$, $\{p_i\}$ being in general position.  But $\exists \ i$: $\kappa_i(x) > 0$.  
Therefore both $x$ and $y$ belong to $U_i$.\\

%dmc remark - consider adding discussion Whitney and Cohen and Nash
\index{Theorem: Embedding Theorem}
\index{Embedding Theorem}
\label{5.1}
\textbf{\small EMBEDDING THEOREM} \quad 
Every second countable normal Hausdorff space of topological dimension $n$ can be embedded in 
$\R^{2n+1}$.\\

\begingroup%%----------------------------------->>
\fontsize{9pt}{11pt}\selectfont 
\textbf{\small EXAMPLE} \  
The exponent ``$2n + 1$'' is sharp.  Indeed, if $K = (V,\Sigma)$, where $\#(V) = 2n + 3$ and $\Sigma$ is the set of all nonempty subsets of \mV, then $\abs{K^{(n)}}$ cannot be embedded in $\R^{2n}$.
\vspi
[Assuming the contrary, work with the cone $\Gamma |{K^{(n)}}|$ of $|{K^{(n)}}|$ (which would embed in 
$\R^{2n+1}$) and construct a continuous function $f:\bS^{2n+1} \ra \R^{2n+1}$ that does not fuse antipodal points, in violation of the Borsuk-Ulam theorem.]\\
\endgroup %%------------------------------------<<

\begingroup%%----------------------------------->>
\fontsize{9pt}{11pt}\selectfont 
\textbf{\small EXAMPLE} \  
Suppose that \mX and \mY are second countable normal Hausdorff spaces of finite topological dimension $-$then the coarse join $X *_c Y$ is a second countable normal Hausdorff space of finite 
%%----------------------------------------------------------------------------------------------29
topological dimension.  
In fact, there exist positive integers $p$ and $q$ such that \mX embeds in $\bS^p$ and \mY embeds in 
$\bS^q$.  Therefore $X *_c Y$ embeds in $\bS^p *_c \bS^q = \bS^{p+q+1}$.\\
\endgroup %%------------------------------------<<

\label{20.25} %dmc mnft
\begingroup%%----------------------------------->>
\fontsize{9pt}{11pt}\selectfont 
Suppose that \mX is a second countable compact Hausdorff space of topological dimension $n > 1$ $-$then, from the proof of the embedding theorem, the set of embeddings $X \ra \R^{2n+1}$ is dense in $C(X,\R^{2n+1})$.  
What can be said about the set of embeddings $X \ra \R^{2n}$?  
Answer: This set can be empty (cf. supra) or nonempty and nowhere dense (cf. infra) or nonempty and dense.  
As regards the latter point, there is a characterization 
(Krasinkiewicz\footnote[2]{\textit{Fund. Math.} \textbf{133} (1989), 247-253.}
, 
Spiez\footnote[3]{\textit{Fund. Math.} \textbf{134} (1990), 105-115; 
see also, \textit{Fund. Math.} \textbf{135} (1990), 127-145.}): 
The set of embeddings $X \ra \R^{2n}$ is dense in $C(X,\R^{2n})$ iff $\dim (X \times X) < 2n$.  
Examples of spaces satisfying this condition are given in $\S 20$ (cf. p. \pageref{19.20}).
\vspi
[Note: \ It can happen that $\forall \ \epsilon > 0$ 
$\exists \ f \in C(X,\R^{2n})$ with $\dev f < \epsilon$ and yet \mX does not embed in $\R^{2n}$.  
Here is an example when $n = 1$.  
Identify $\R^2$ with the set of $(x,y,z) \in \R^3$: 
$z = 0$.  Put 
$A = \ds\bigcup\limits_1^\infty (1/n) \bS^1$, 
$B = \{(x,0,0): \abs{x} \leq 1\} \cup \{(0,y,0): \abs{y} \leq 1\}$, 
$C = \{(0,0,z): 0 \leq z \leq 1\}$  and set $X = A \cup B \cup C$.  
Given $\epsilon > 0$, select $k: 1/2k < \epsilon$.  Denote by 
$X_k$ the quotient $X/K$, \mK the subset of $A \cup B$ consisting of those points whose distance from the origin is 
$\leq 1/2 k$.  
Let $p$ be the projection $X \ra X_k$, choose an embedding $f_k:X_k \ra \R^2$ and consider $f = f_k \circx p$.  
Nevertheless, \mX cannot be embedded in $\R^2$.]\\
\endgroup %%------------------------------------<<

\begingroup%%----------------------------------->>
\fontsize{9pt}{11pt}\selectfont 
\textbf{\small EXAMPLE} \  
The set of embeddings $[0,1]^n \ra \R^{2n}$ is nonempty and nowhere dense in $C([0,1]^n,\R^{2n})$.
\vspi
[Show that there exists a function $f_0 \in C([0,1]^n,\R^{2n})$ and an $\epsilon_0 > 0$ such that if 
$f \in C([0,1]^n,\R^{2n})$ and if $\norm{f_0 - f} < \epsilon_0$, then $f$ is not one-to-one.]\\
\endgroup %%------------------------------------<<

\begingroup%%----------------------------------->>
\fontsize{9pt}{11pt}\selectfont 
\textbf{\small FACT} \  
Suppose that \mX is a second countable normal Hausdorff space of topological dimension $n$.  Equip the function space 
$C(X,\R^{2n+1})$ with the limitation topology $-$then the set of embeddings $X \ra \R^{2n+1}$ contains a dense 
$G_\delta$ in $C(X,\R^{2n+1})$.\\
\endgroup %%------------------------------------<<

Suppose that $\dim X = n$ $-$then there is a closed embedding $X \ra \R^{2n+1}$ iff \mX is second countable and locally compact.  For $X_\infty$ is second countable  and $\dim X = \dim X_\infty$ (by the control lemma).  
Embed $X_\infty$ in $\R^{2n+1}$.  Add to $\R^{2n+1}$ a point at infinity and remove the point corresponding to 
$X_\infty - X$.  This gives another copy of $\R^{2n+1}$ containing \mX as a closed subset.\\

Put 
$\N_n^{2n+1} = \Q_0^{2n+1} \cup \cdots \cup \Q_n^{2n+1}$, the subspace of $\R^{2n+1}$ consisting of all points 
with at most $n$ rational coordinates $-$then $\dim \N_n^{2n+1} = n$.\\

\textbf{\small LEMMA} \  
Every second countable normal Hausdorff space of topological dimension $n$ can be embedded in $\N_n^{2n+1}$.

%%----------------------------------------------------------------------------------------------30
[The complement $\R^{2n+1} - \N_n^{2n+1}$ has the form $\bigcup\limits_1^\infty H_k$, where $\forall \ k$, 
$H_k$ is a plane of euclidean dimension $n$.  Take \mX compact, let 
$D_{1/k}(H_k) = D_{1/k} \cap \{f:f(X) \cap H_k = \emptyset\}$, and consider $\bigcap\limits_1^\infty D_{1/k}(H_k)$.]\\

Application: Every second countable normal Hausdorff space of topological dimension $n$ can be written as a union of 
$n + 1$ subspaces, each of topological dimension $\leq 0$.

[Note: \ Filippov\footnote[2]{\textit{Soviet Math. Dokl.} \textbf{11} (1970), 687-691.}
has constructed an example of a compact perfectly normal \mX: $\dim X = 1$, which cannot be written as a union 
$X_1 \cup X_2$, where 
$
\begin{cases}
\ \dim X_1 = 0\\
\ \dim X_2 = 0
\end{cases}
$
.]\\

\begingroup%%----------------------------------->>
\fontsize{9pt}{11pt}\selectfont 
When $n = 0$, the space $\N_n^{2n+1}$ becomes the set of irrationals, the latter being homeomorphic to 
$\N^\omega$.  The Cantor cube $C_\omega$ embeds in $\N^\omega$ and, as has been noted on 
p. \pageref{19.21}, 
if 
$
\begin{cases}
\ \dim X = 0\\
\ \wt X \leq \omega
\end{cases}
, \ 
$
then \mX embeds in $C_\omega$.  
There is a higher dimensional counterpart to this in that one can construct a compact subspace 
$M_n^{2n+1} \subset \R^{2n+1}$ of topological dimension $n$ which embeds in $\N_n^{2n+1}$ and has the property that if 
$
\begin{cases}
\ \dim X = n\\
\ \wt X \leq \omega
\end{cases}
, \ 
$
then \mX embeds in $M_n^{2n+1}$.  
In a word: Subdivide $[0,1]^{2n+1}$ into cubes of side length 1/3, retain those that meet the $n$-faces of $[0,1]^{2n+1}$, repeat the process on each element of their union $K_0$ and continue to the limit: 
$M_n^{2n+1} = \ds\bigcap\limits_0^\infty K_i$ 
(Bothe\footnote[3]{\textit{Fund. Math.} \textbf{52} (1963), 209-224; 
see also Bestvina, \textit{Memoirs  Amer. Math. Soc.} 380 (1988), 1-110.}).\\
\endgroup %%------------------------------------<<

\begingroup%%----------------------------------->>
\fontsize{9pt}{11pt}\selectfont 
Denote by $\N_n(\kappa)$ the subspace of $\bS(\kappa)^\omega$ consisting of those points which have at most $n$ nonzero rational coordinates $-$then 
$
\begin{cases}
\ \wt \N_n(\kappa) = \kappa\\
\ \dim \N_n(\kappa) = n\\
\end{cases}
$
.\\
\endgroup %%------------------------------------<<

\begingroup%%----------------------------------->>
\fontsize{9pt}{11pt}\selectfont 
\textbf{\small FACT} \  
Every metrizable space \mX of weight $\leq \kappa$ and of topological dimension $\leq n$ can be embedded in 
$\N_n(\kappa)$.
\vspi
[Note: \ By comparison, recall that every metrizable space \mX of weight $\leq \kappa$ can be embedded in 
$\bS(\kappa)^\omega$ (cf. p. \pageref{19.22}).]\\
\endgroup %%------------------------------------<<

\begingroup%%----------------------------------->>
\fontsize{9pt}{11pt}\selectfont 
Suppose that \mX is metrizable (completely metrizable) of weight $\kappa$.  Equip the function space 
$C(X,$ $\bS(\kappa)^\omega)$ with the limitation topology $-$then 
Pol\footnote[6]{\textit{Topology Appl.} \textbf{39} (1991), 189-204.} 
has shown that the set of embeddings (closed embeddings) $X \ra \bS(\kappa)^\omega$ contains a dense 
$G_\delta$ in $C(X,\bS(\kappa)^\omega)$.\\
\endgroup %%------------------------------------<<

Can one characterize dim by a set of axioms on the class $\sE$, the subspaces of euclidean space?  
The answer is ``yes''.\\

%%----------------------------------------------------------------------------------------------31
Consider a function $d:\sE \ra \{-1,0,1,\ldots\}$ subject to:

\indent\indent (d$_1$) \ (Normalization Axiom) $d(\emptyset) = -1$, $d([0,1]^n) = n$ $(n = 0, 1, \ldots)$.

\indent\indent (d$_2$) \ (Topological Invariance Axiom) \  \ If $X, Y \in \sE$ are homeomorphic, then \quad $d(X) = d(Y)$.

\indent\indent (d$_3$) \ (Monotonicity Axiom) If $X, Y \in \sE$ with $X \subset Y$, then $d(X) \leq d(Y)$.

\indent\indent (d$_4$) \ (Countable Union Axiom) If $X \in \sE$ is the union of a sequence of closed subspaces $X_i$, then 
$d(X) \leq \sup\limits_i d(X_i)$.

\indent\indent (d$_5$) \ (Compactification Axiom) If $X \in \sE$, then there is a compactification $\widetilde{X} \in \sE$ of \mX 
such that $d(X) = d(\widetilde{X})$.

\indent\indent (d$_6$) \ (Decomposition Axiom) If $X \in \sE$ and $d(X) = n$, then there exist $n + 1$ sets $X_i \subset X$ such that $X = \bigcup\limits_0^n X_i$ and $\forall \ i$, $d(X_i) \leq 0$.

Hayashi\footnote[2]{\textit{Topology Appl.} \textbf{37} (1990), 83-92.} 
has shown that these axioms are independent and serve to characterize the topological dimension dim on the class $\sE$.

[Note: \ The key here is the last axiom on the list.  The first five are satisfied by the cohomological dimension 
$\dim_G$ with respect to a nonzero finitely generated abelian group \mG.]\\

While it is not true in general that an arbitrary normal \mX of topological dimension $n$ can be written as a union of $n + 1$ normal subspaces, each of topological dimension $\leq 0$, there is nevertheless a partial substitute in that every neighborhood finite open covering of \mX of order $\leq n + 1$ has an open refinement that can be written as a union of $n + 1$ collections, each of order $\leq 1$.  This is a consequence of the following statement.\\

\index{Decomposition Lemma}
\textbf{\small DECOMPOSITION LEMMA} \  
Let $\sU = \{U_i: i \in I\}$ be a neighborhood finite open covering of \mX of order $\leq n+1$ 
$-$then there exists an open covering $\sV$ of \mX which can be represented as a union of $n+1$ collections 
$\sV_0, \ldots, \sV_n$, where 
$\sV_j = \{V_{i,j}: i \in I\}$ consists of pairwise disjoint open sets such that $\forall \ i: V_{i,j} \subset U_i$.

[There is nothing to prove if $n = 0$.  
Proceeding by induction, assume the validity
 of the assertion for all normal spaces and for all neighborhood finite open coverings of order 
 $< \ n+1$ $(n \geq 1)$.  
 Choose a precise open refinement 
$\sO = \{O_i$ : $i \in I\}$ of \ $\sU = \{U_i$ : $i \in I\}$ : $\forall \ i$, $A_i \equiv \ov{O}_i \subset U_i$.  Put 
$\sF = \{F: F \subset I \ \& \ \#(F) = n + 1\}$.  
Assign to each $F \in \sF$: $U_F = \bigcap\limits_{i \in F} U_i$ and 
$
\begin{cases}
\ O_F = \bigcap\limits_{i \in F} O_i\\
\ A_F = \bigcap\limits_{i \in F} A_i
\end{cases}
. \ 
$
Select a point $i_F \in F$ and let $V_{i,n} = \bigcup \{U_F: i_F = i\}$ 
%%----------------------------------------------------------------------------------------------32
$-$then the order of $\sV_n = \{V_{i,n} : i \in I\}$ is $ \leq 1$ and $\forall \ i$ : $V_{i,n} \subset U_i$.  
The subspace 
$Y = X - \bigcup\limits_F O_F$ is closed, hence normal.  
Since the order of the neighborhood finite open covering 
$\{Y \cap O_i: i \in I\}$ of \mY is $\leq n$, there exists an open covering 
$\sV^\prime$ of \mY which can be represented as a union of $n$ collections 
$\sV_0^\prime, \ldots, \sV_{n-1}^\prime$, where 
$\sV_j^\prime = \{\sV_{i,j}^\prime: i \in I\}$ 
consists of pairwise disjoint open sets such that $\forall \ i$: $\sV_{i,j}^\prime \subset Y \cap O_i$.  
The subspace 
$Z = X - \bigcup\limits_F A_F$ is open ($\{A_F\}$ is neighborhood finite) and is contained in \mY.  
For $j = 0, \ldots, n-1$, let 
$V_{i,j} = Z \bigcap V_{i,j}^\prime$ and $\sV_j = \{V_{i,j}: i \in I\}$.  
Consideration of the union 
$\sV = \bigcup\limits_0^n \sV_j$ completes the induction.]\\

\begin{proposition} \ %11
Suppose that $\dim X \leq n$. Let $\sU = \{U_i: i \in I\}$ be a neighborhood finite open covering of \mX $-$then there exist sequences 
$
\begin{cases}
\ \sV_0, \sV_1, \ldots\\
\ \sW_0, \sW_1, \ldots
\end{cases}
$
of discrete collections of open subsets 
$\sV_j = \{V_{i,j}:i \in I\}$ $\&$ $\sW_j = \{W_{i,j}:i \in I\}$ of \mX such that any  
$n + 1$ of the $\sV_j$ cover \mX and $\forall \ i$: $\ov{V_{i,j}} \subset W_{i,j} \subset U_i$.
\end{proposition}

[Bearing in mind Proposition 6, normality and the decomposition lemma provide us with the $\sV_j$ and $\sW_j$ for 
$j \leq n$.  
Now argue by induction, assuming that $\sV_j$ and $\sW_j$ have been defined for $j \leq m-1$, $m-1$ being 
$\geq n$.  
Assign to each $M \subset \{0,\ldots, m-1\}$ of cardinality $n$ the closed subset 
$A_M = X - \bigcup\limits_{j \in M} \cup \hspace{0.03cm} \sV_j$ 
$-$then the $A_M$ are pairwise disjoint because any $n+1$ of the $\sV_j$ cover 
\mX.  Determine open
$
\begin{cases}
\ V_M\\
\ W_M
\end{cases}
$
: 
$A_M \subset V_M \subset \ov{V}_M \subset W_M$, where $M^\prime \neq M\pp$ $\implies$ 
$\ov{W}_{M^\prime} \cap \ov{W}_{M\pp} = \emptyset$.  
Select a point $j_M \leq m-1$: $j_M \notin M$.  
Note that 
$A_M \subset \cup\hspace{0.03cm} \sV_{j_M}$.  Put
$
\begin{cases}
\ V_{i,m} = \bigcup\limits_M V_M \cap V_{i,j_M}\\
\ W_{i,m} = \bigcup\limits_M W_M \cap W_{i,j_M}
\end{cases}
. \ 
$
The associated collections $\sV_m$ and $\sW_m$ are discrete and open with $\ov{V}_{i,m} \subset W_{i,m} \subset U_i$.  
And since any $n$ of the $\sV_j$ $(j \leq m-1)$ cover 
$X - \bigcup\limits_M A_M$, any $n+1$ of the $\sV_j$ $(j \leq m)$ cover 
\mX.]\\

The Kolmogorov superposition theorem, which resolved Hilbert's 13th problem in the negative, says that for each 
$n \geq 1$ there exist functions $\phi_1, \ldots, \phi_{2n+1}$ in $C([0,1]^n)$ such that every $f \in C([0,1]^n)$ can be represented in the form $f = \sum\limits_i g_i \circx \phi_i$ for certain $g_i \in C(\R)$ (depending on $f$).  
Objective: Isolate the dimension theoretic content of this result.

Suppose that \mX is a second countable compact Hausdorff space.  Let $\phi_i \in C(X)$ $(i = 1, \ldots, k)$ $-$then
the collection $\{\phi_i\}$ is said to be 
\un{basic}
\index{basic (collection of functions)} 
if for every $f \in C(X)$ there exist continuous functions $g_i:\R \ra \R$ such that $f = \sum\limits_i g_i \circx \phi_i$.  
A 
\un{basic embedding}
\index{basic embedding (\mX in $\R^k$ )} 
of \mX in $\R^k$ is an embedding $X \ra \R^k$ corresponding to a basic collection $\{\phi_i\}$.  So, e.g., according to Kolmogorov, $X = [0,1]^n$ can be basically embedded in $\R^{2n+1}$.\\

%%----------------------------------------------------------------------------------------------33
\index{Theorem: Basic Embedding Theorem}
\index{Basic Embedding Theorem}
\textbf{\small BASIC EMBEDDING THEOREM} \  
Every second countable compact Hausdorff space of topological dimension $n$ can be basically embedded in $\R^{2n+1}$.

[Note: \ Sternfeld\footnote[2]{\textit{Israel J. Math.} \textbf{50} (1985), 13-53; 
see also Levin, \textit{Israel J. Math.} \textbf{70} (1990), 205-218.} 
has shown that if $\dim X = n$ $(n > 1)$, then \mX cannot be basically embedded in $\R^{2n}$.  
Example: Let $X = \{(x,0):\abs{x}\leq 1\} \cup \{(0,y):\abs{y} \leq 1\}$ $-$then $\dim X = 1$ and \mX can be basically embedded in $\R^2$.]\\

The proof of the basic embedding theorem is not a general position argument.  It depends instead on Proposition 11 and some elementary functional analysis.\\

There is a simple interpretation of what it means for $\{\phi_i\}$ to be basic in terms of the dual $C(X)^*$ of $C(X)$.  Thus put $Y_i = \phi_i(X)$ and let $Y = \coprod\limits_i Y_i$ 
$-$then the collection $\{\phi_i\}$ determines a bounded linear operator 
$T:C(Y) \ra C(X)$, viz. $T(g_1, \ldots, g_k) = \sum\limits_i g_i \circx \phi_i$ with adjoint 
$T^*:C(X)^* \ra C(Y)^*$, viz. $T^*\mu = \sum\limits_i \mu_i$, \ 
$\mu_i$ the image of $\mu$ under $\phi_i$.  
Note that 
$\norm{T^*\mu} = \sum\limits_i \norm{\mu_i}$.  
Obviously, $\{\phi_i\}$ is basic iff \mT is surjective or still, iff $\exists$ 
$\lambda:0 < \lambda \leq 1$ such that $\forall \ \mu \in C(X)^*$ 
$\exists \ i$ : 
$\norm{\mu_i} \geq \lambda \norm{\mu}$.  
When this occurs, call $\{\phi_i\}$ 
\un{$\lambda$-basic}
\index{lambda-basic, $\lambda$-basic}.

Fix a compatible metric $d$ on \mX.  
Given a finite discrete collection $\sU = \{U\}$ of open subsets of \mX, we shall write $d(\sU)$ for $\sup\{\diam U: U \in \sU\}$ and agree that a function $\phi \in C(X)$ separates $\sU$ if $\forall \ U \neq V$ in 
$\sU$: $\phi(\ov{U}) \cap \phi(\ov{V}) = \emptyset$.\\

\textbf{\small LEMMA} \  
Let $\phi_i \in C(X)$ $(i = 1, \ldots, k)$.  Suppose that 
$\forall \ \epsilon > 0$ and $\forall \ i$, there exists a finite discrete collection 
$\sU_i$ of open subsets of \mX with $d(\sU_i) < \epsilon$ such that $\phi_i$ separates $\sU_i$ and 
\[
\forall \ x \ \in X : \quad \sum\limits_i \ord(x,\sU_i) \geq \left[\frac{k}{2}\right] + 1.
\]
Then $\{\phi_i\}$ is $1/k$-basic.

[The set of $\mu \in C(X)^*$ for which 
$\spt(\mu^+) \cap \hspace{0.05cm} \spt(\mu^-) \hspace{0.05cm} = \hspace{0.05cm} \emptyset$ 
is dense in $C(X)^*$ \ (Hahn \quad plus regularity).  
Therefore take a $\mu \in C(X)^*$ of norm one, assume that $\epsilon = d(\spt(\mu^+)$, $\spt(\mu^-)) > 0$, 
and choose the \ $\sU_i$ accordingly.  
If as usual $\abs{\mu} = \mu^+ + \mu^-$, then $\abs{\mu}$ is a probability measure on \mX and 
$\sum\limits_i \abs{\mu}(\cup \hspace{0.03 cm} \sU_i) \geq [k/2] + 1$, implying that for some $i_0$, 
$\abs{\mu}(\cup \hspace{0.03 cm} \sU_{i_0}) \geq$ 
$(1/k)([k/2] + 1) \geq$ 
$1/2 + 1/2k.$ 
On the other hand, 
$\forall$ $U \in \sU_{i_0}$, $\abs{\mu}(U) = \abs{\mu(U)}$, thus 
$\abs{\mu}(\cup \hspace{0.03 cm} \sU_{i_0}) = \sum\limits_U \abs{\mu(U)}$ and so 
$\norm{\mu_{i_0}} \geq 1/2 + 1/2k - \abs{\mu}(X - \cup \hspace{0.03 cm} \sU_{i_0}) \geq 1/k$.]\\

Let $\sU(p)$ be a finite discrete collection of open subsets of \mX with $d(\sU(p)) < 1/p$ $(p = 1, 2, \ldots)$.  
Claim: There exists a dense set of $\phi \in C(X)$ separating $\sU(p)$ for infinitely 
%%----------------------------------------------------------------------------------------------34
many $p$.  
To see this, let $\Phi_q$ be the set of $\phi \in C(X)$ separating $\sU(p)$ for some $p \geq q$ $(q = 1, 2, \ldots)$ 
$-$then it need only be shown that $\forall \ q$, $\Phi_q$ is open and dense in $C(X)$ (consider 
$\bigcap\limits_1^\infty \Phi_q$ and quote Baire).

\indent\indent (1) \ $\Phi_q$ is open in $C(X)$.  Proof: Let $\phi \in \Phi_q$.  
Choose $p$ per $\phi$.  Let 
$2\epsilon = \inf\{\dis(\phi(\ov{U}),\phi(\ov{V}))$: $U \neq V$ in $\sU(p)\}$.  
Suppose that $\norm{\phi - f} < \epsilon/4$ 
$-$then $U \neq V$ in $\sU(p)$ $\implies$ $\dis(f(\ov{U}),f(\ov{V})) > \epsilon$.

\indent\indent (2) \ $\Phi_q$ is dense in $C(X)$.  Proof: Fix $f \in C(X)$.  Given $\epsilon > 0$, choose $p \geq q$: 
$\osc(\restr{f}{\ov{U}}) < \epsilon/2$ $(U \in \sU(p))$.  Define a continuous function $g: \cup \ \ov{U} \ra \R$ by picking distinct constants $c_U:$ 
$
\begin{cases}
\ \restr{g}{\ov{U}} = c_U\\
\ \norm{\restr{f}{\ov{U}} - \restr{g}{\ov{U}}} < \epsilon
\end{cases}
. \ 
$
Use Tietze and extend $\restr{f}{\cup \ \ov{U}} - g$ to an $h \in C(X)$: $\norm{h} < \epsilon$.  Put 
$\phi = f - h$: $\phi \in \Phi_q$ $\&$ $\norm{f - \phi} < \epsilon$.\\
%\vspace{0.75cm}

To prove the basic embedding theorem, take $k = 2n+1$ $-$then, in view of Proposition 11, there exist finite discrete collections $\sU_i(p)$ $(i = 1, \ldots, k)$ of open subsets of \mX with $d(\sU_i(p)) < 1/p$ $(p = 1, 2, \ldots)$ such that for each $p$ the union of any $n + 1$ of the $\sU_i(p)$ is a covering of \mX, so
\[
\forall \ x \in X \ : \quad \sum\limits_i \ord(x,\sU_i(p)) \geq \left[\frac{k}{2}\right] + 1.
\]
Thanks to the preceeding remarks, it is possible to select integers $p_1 < p_2 < \cdots$ and functions $\phi_i \in C(X)$ 
$(i = 1, \ldots, k)$ having the property that $\phi_i$ separates $\sU_i(p_j)$ $(j = 1, 2, \ldots)$.  
Apply the lemma and conclude that $\{\phi_i\}$ is 1/k-basic $(k = 2n + 1)$.\\

\begingroup%%----------------------------------->>
\fontsize{9pt}{11pt}\selectfont 
When $X = [0,1]^n$, one can explicate, at least to some extent, the analytic structure of the $\phi_i$.  Precisely put: Given rationally independent real numbers $r_1, \ldots, r_n$, there exist increasing continuous functions 
$\psi_1, \ldots, \psi_{2n+1}$ on $[0,1]$ such that the
\[
\phi_i(x_1, \ldots, x_n) = \sum\limits_{j = 1}^n r_j \psi_i(x_j) \quad (1 \leq i \leq 2n+1)
\]
constitute a $1/k$-basic collection $(k = 2n + 1)$.  Moreover, the $g_i$ can be chosen independently of $i$, so 
$\forall \ f \in C([0,1]^n)$ there exists a $g \in C(\R)$: 
\[
f(x_1, \ldots, x_n) = \sum\limits_{i = 1}^{2n + 1} g\left( \sum\limits_{j= 1}^n r_j \psi_i(x_j) \right).
\]
\vspi
[Note: The ``inner functions'' can even be taken in Lip$_1([0,1])$.  
Reason: There exists a homeomorphism 
$\iota:[0,1] \ra [0,1]$ such that $\forall \ i$, $\psi_i \circx \iota \in $Lip$_1([0,1])$.  
Consider, e.g., the inverse to the assignment 
$x \ra C(x + \ds\sum\limits_i(\psi_i(x) - \psi_i(0)))$, where \mC is the reciprocal of 
$1 + \sum\limits_i (\psi_i(1) - \psi_i(0))$.]
\vspi
To avoid trivialities, assume that $n > 1$.  There are then three steps to the proof.
\\
\indent\indent (I) For $p = 1, 2, \ldots$, partition $[0,1]$ into $p$ closed subintervals \mI of length $1/p$ indexed by the natural order and for $1 \leq i \leq k$, 
let $\sI_i(p)$ denote the collection of closed subintervals of $[0,1]$ obtained by removing from 
$[0,1]$ the interior of those \mI whose index is congruent to $i \mod k$.  
Write $\sC_i(p)$ for the set of
%%----------------------------------------------------------------------------------------------35
all products $C_i(p) = I_1(p) \times \cdots \times I_n(p)$: $\forall \ j$, $I_j(p) \in \sI_i(p)$.  
It is clear that $\sC_i(p)$ is a discrete collection of closed $n$-cubes in $[0,1]^n$.  
Furthermore, every $x \in [0,1]^n$ belongs to at least 
$[k/2] + 1 \equiv n + 1$ of the $\cup C_i(p)$.
\\
\indent\indent (II) Let $\Psi$ stand for the set of increasing continuous functions on $[0,1]$, equipped with the uniform norm.  Attach to each $\epsilon > 0$: $0 < \epsilon < 1/2k$, and to each $f \in C([0,1]^n)$: $\norm{f} \neq 0$, the set 
$\Omega_f(\epsilon)$ of all $\{\psi_i\} \in \Psi^k$ for which there exists an $h \in C(\R)$ : $\norm{h} \leq \norm{f}$ $\&$ 
$||{f - \sum\limits_i h \bigl(\sum\limits_j r_j\psi_i\bigr)}|| <$ 
%$\norm{f - \sum\limits_i\left(\sum\limits_j r_j\psi_i\right)} <$ 
$(1 - \epsilon)\norm{f}$.  
Claim: $\Omega_f(\epsilon)$ is open and dense.  
Of course, only the density is at issue.  
And for this, it suffices to fix a nonempty open 
$\Omega \subset \Psi^k$ and show that 
$\Omega \cap \Omega_f(\epsilon) \neq \emptyset$.  
Let $\Psi^k(p)$ be the subset of $\Psi^k$ consisting of the $\{\psi_i\}$ such that 
$\forall \ i$: $\psi_i$ is constant on the elements of $\sI_i(p)$.  
Choose $p \gg 0$: $\Omega \cap \Psi^k(p) \neq \emptyset$ $\&$ 
$\osc(\restr{f}{C_i(p)}) < \epsilon \norm{f}$ $\forall \ C_i(p) \in \sC_i(p)$.  
Fix $\{\psi_i\} \in \Omega \cap \Psi^k(p)$.  Because the 
$r_j$ are rationally independent, there is no loss of generality in supposing that $\phi_i \equiv \ds\sum\limits_j r_j \phi_i$ takes 
different values on different elements of $\sC_i(p)$ and that in addition these values are distinct for distinct $i$.  
We shall now 
construct an $h \in C(\R)$ in terms of the $\phi_i$ and deduce that $\{\psi_i\} \in \Omega_f(\epsilon)$.  Call $M_i$ the value of $f$ at the center of $C_i(p)$.  
Let $h(\phi_i(C_i(p))) = 2 \epsilon M_i$ and extend $h$ continuously to all of $\R$: 
$\norm{h} \leq 2 \epsilon \norm{f}$.  
Using the fact that every $x \in [0,1]^n$ belongs to at least $n + 1$ of the 
$\cup \hspace{0.03 cm} \sC_i(p)$, one has
\endgroup

\begingroup
\fontsize{9pt}{11pt}\selectfont 
\vspace{-.25cm}
\begin{align*}
%$
|{f(x) - \sum\limits_i h(\phi_i(x))}| &\leq (1 - 2(n+1)\epsilon)\abs{f(x)} + 2(n+1) \epsilon^2 \norm{f} + 2n\epsilon\norm{f}\\
&\leq (1 - 2 \epsilon + 2(n + 1)\epsilon^2)\norm{f} < (1 - \epsilon) \norm{f}.
\end{align*}
\endgroup

\begingroup
\fontsize{9pt}{11pt}\selectfont 
\vspace{-.15cm}
Therefore $\{\psi_i\} \in \Omega_f(\epsilon)$.
\\
\indent\indent (III) Let $D = \{f_d\}$ be a countable dense subset of $C([0,1]^n)$, not containing the zero function $-$then 
$\ds\bigcap\limits_1^\infty \Omega_{f_d}(\epsilon)$ is dense in $\Psi^k$ (Baire).  
Fix 
$\{\psi_i\} \in \ds\bigcap\limits_1^\infty \Omega_{f_d}(\epsilon)$.  Let $f \in C([0,1]^n)$: $\norm{f} \neq 0$.  Choose 
$f_d \in D$: $\norm{(1 - \epsilon/4)f - f_d} < (\epsilon/4)\norm{f}$, so
$
\begin{cases}
\ \norm{f_d} \leq \norm{f}\\
\ \norm{f - f_d} < (\epsilon/2) \norm{f}
\end{cases}
$
and choose $h_d \in C(\R)$ : $\norm{h_d} \leq \norm{f_d}$ $\&$ 
$||{f_d - \ds\sum\limits_i h_d\bigl(\ds\sum\limits_j r_j \psi_i\bigr)}|| < (1 - \epsilon) \norm{f_d}$.  
Conclusion: 
$\exists$ $h = \gamma(f) \in C(\R)$ such that $\norm{h} \leq$ $\norm{f}$ $\&$ 
$||{f - \sum\limits_i h\bigl(\ds\sum\limits_j r_j \psi_i\bigr)}|| <$ $(1 - \epsilon/2) \norm{f}$. 
Recursively define a sequence $\chi_0, \chi_1, \ldots$ in $C([0,1]^n)$ by $\chi_0 =f$, 
$\chi_{m+1} =$ $\chi_m - \ds\sum\limits_i h_m\bigl(\ds\sum\limits_j r_j \psi_i\bigr)$, where 
$h_m = \gamma(\chi_m)$ $(\gamma(0) = 0)$.  
The series $\sum\limits_0^\infty h_m$ is uniformly convergent, thus its sum $g$ is continuous and satisfies the relation 
$f = \ds\sum\limits_i g \bigl(\ds\sum\limits_j r_j \psi_i\bigr)$.
\vspi
[Note: Let $C^1([0,1]^n)$ be the set of continuously differentiable functions on $[0,1]^n$ $-$then 
Kaufman\footnote[2]{\textit{Proc. Amer. Math. Soc.} \textbf{46} (1974), 360-362.}
has shown that for $n > 1$, no finite subset of $C^1([0,1]^n)$ can be basic.]\\
\endgroup %%------------------------------------<<

\begingroup%%----------------------------------->>
\fontsize{9pt}{11pt}\selectfont 
\textbf{\small FACT} \  
There exist real valued continuous functions $\phi_i$ $(i = 1, \ldots, 2n+1)$ on $\R^n$ such that $\forall \ f \in BC(\R^n)$ 
$\exists$ $g \in C(\R)$: $f = \ds\sum\limits_i g \circx \phi_i$.
\vspi
%%----------------------------------------------------------------------------------------------36
[Note: \ This result remains true if $\R^n$ is replaced by a noncompact second countable LCH space \mX of topological dimension $n$.]\\
\endgroup %%------------------------------------<<

If \mX and \mY are nonempty normal Hausdorff spaces, what is the relation between $\dim (X \times Y)$ and 
$
\begin{cases}
\ \dim X\\
\ \dim Y
\end{cases}
? \ 
$
An initial difficulty is that $X \times Y$ need not be normal so formally $\dim (X \times Y)$ can be undefined.\\

\begingroup%%----------------------------------->>
\fontsize{9pt}{11pt}\selectfont 
This is not a serious problem.  Reason $X \times Y$ is at least completely regular, therefore in this context 
$\dim(X \times Y)$ is meaninful (cf. p. \pageref{19.23}).\\
\endgroup %%------------------------------------<<

Examples: 
(1) \ Take $X = Y = $ Sorgenfrey line $-$then \mX is perfectly normal and paracompact but $X \times X$ is not normal 
(cf. p. \pageref{19.24}); 
(2) \ Take $X = [0,\Omega[$, $Y = [0,\Omega]$ $-$then \mX is normal and \mY is compact but $X \times Y$ is not normal; 
(3) \ Take $X = $ Michael line, $Y = \PP$ $-$then \mX is paracompact and \mY is metrizable but $X \times Y$ is not normal 
(cf. \pageref{19.24a} ff.);
(4) \ Take $X = $ Rudin's Dowker space, $Y = [0,1]$ $-$then $X \times [0,1]$ is not normal.\\

Here are some conditions on \mX and \mY that ensure that the product $X \times Y$ is normal.
\indent\indent (1) \ Suppose that \mX is perfectly normal (perfectly normal and paracompact) and \mY is metrizable $-$then $X \times Y$ is perfectly normal (perfectly normal and paracompact).

\indent\indent (2) \ Suppose that \mX is normal and countably compact and \mY is metrizable $-$then $X \times Y$ is normal.

\indent\indent (3) \ Suppose that \mX is normal and countably paracompact and \mY is metrizable and $\sigma$-locally compact  $-$then $X \times Y$ is normal.

\label{20.14}
\indent\indent (4) \ Suppose that \mX is paracompact and \mY is paracompact and $\sigma$-locally compact $-$then $X \times Y$ is paracompact.

[Note: \ A CRH space is said to be 
\un{$\sigma$-locally compact}
\index{sigma-locally compact, $\sigma$-locally compact}  
if it can be written as a countable union of closed locally compact subspaces.  
Example: Every CW complex is $\sigma$-locally compact.]\\

If enough pathology is built into \mX and \mY, then it can happen that $\dim X + \dim Y < \dim (X \times Y)$.  
Examples illustrating the point are given below.  Because of this, one looks instead for conditions on \mX and \mY that serve to force $\dim(X \times Y) \leq \dim X + \dim Y$.\\

\index{Theorem: Product Theorem}
\index{Product Theorem}
\textbf{\small PRODUCT THEOREM} \quad 
Suppose that \mX is normal and \mY is paracompact and $\sigma$-locally compact.  
Assume: $X \times Y$ is normal 
$-$then 
$\dim (X \times Y) \leq \dim X + \dim Y$.

[Note: \ Tacitly, $X \neq \emptyset$ $\&$ $Y \neq \emptyset$.]\\

%%----------------------------------------------------------------------------------------------37
\begingroup%%----------------------------------->>
\fontsize{9pt}{11pt}\selectfont 
The inequality in the product theorem can be strict even if \mX and \mY are compact ARs 
(Dranishnikov\footnote[2]{\textit{Soviet Math. Dokl.} \textbf{37} (1988), 769-773.}).\\
\endgroup %%------------------------------------<<

The proof of the product theorem is carried out in stages under the supposition that 
$
\begin{cases}
\ n = \dim X\\
\ m = \dim Y
\end{cases}
< \infty.
$
\\

\begin{proposition} \ 
Suppose that both \mX and \mY are compact $-$then $\dim(X \times Y) \leq \dim X + \dim Y$.
\end{proposition}

[Let $\sW$ be a finite open covering of $X \times Y$.  Choose finite open coverings 
$
\begin{cases}
\ \sU\\
\ \sV
\end{cases}
$
of 
$
\begin{cases}
\ X\\
\ Y
\end{cases}
:\sU \times \sV
$
refines $\sW$.  Attach to 
$
\begin{cases}
\ \sU\\
\ \sV
\end{cases}
$
sequences 
$
\begin{cases}
\ \sO_0, \sO_1, \ldots\\
\ \sP_0, \sP_1, \ldots
\end{cases}
$
of discrete collections of open subsets of 
$
\begin{cases}
\ X\\
\ Y
\end{cases}
$
having the properties delineated in Proposition 11.  
In particular: Each $x \in X$ can fail to belong to at most $n$ of the 
$\cup \hspace{0.03 cm} \sO_k$ and each $y \in Y$ can fail to belong to at most $m$ of the $\cup \hspace{0.03 cm}  \sP_k$.  
The union 
$\sO_0 \times \sP_0 \bigcup \cdots \bigcup \sO_{n+m} \times \sP_{n+m}$ is therefore an open refinement of 
$\sU \times \sV$ of order $\leq n + m + 1$.]\\

\begingroup%%----------------------------------->>
\fontsize{9pt}{11pt}\selectfont 
If \mX and \mY are compact and metrizable and if $f:X \ra Y$ is continuous and surjective, then there exists a Baire class one function $g:Y \ra X$ such that $f \circx g = \id_Y$ 
(Engelking\footnote[3]{\textit{Bull. Acad. Polon. Sci.} \textbf{16} (1968), 277-282.}
Since $g \circx f$ is a function of the first Baire class, its graph is a $G_\delta$ in $X \times X$, 
which implies that the range of $g$, viz. 
$\{x:g(f(x)) = x\}$, is a $G_\delta$ in \mX that intersects each fiber of $f$ in exactly one point.\\
\endgroup %%------------------------------------<<

\begingroup%%----------------------------------->>
\fontsize{9pt}{11pt}\selectfont 
\textbf{\small EXAMPLE} \  
Let $\sK$ be the collection of all nonempty closed subsets of $[0,1] \times [0,1]$ equipped with the Vietoris topology, so $\sK$ 
is compact and metrizable.  
Write $p$ for the vertical projection $-$then the collection $\sC$ of all compact connected subsets of 
$[0,1] \times [0,1]$ that meet both $p^{-1}(0)$ and $p^{-1}(1)$ is a closed subspace of $\sK$, hence is compact.  
Therefore there exists a continuous surjection $\Gamma$ from the Cantor set $C \subset [0,1]$ to $\sC$.  
Because $C \times C$ is 
homeomophic to \mC, one can assume that the fibers of $\Gamma$ have cardinality $2^\omega$.  
If now 
$X = \ds\bigcup \{p^{-1}(t) \cap \Gamma(t): t \in C\}$, then \mX is a compact subspace of $[0,1] \times [0,1]$ and 
$f \equiv \restr{p}{X}:X \ra C$ is surjective.  
From the remark above, there exists a Baire class one function $g:C \ra X$ such that 
$f \circx g = \id_C$.  
Define $\phi:C \ra [0,1]$ by $g(t) = (t,\phi(t))$: $\phi$ is a function of the first Baire class and its graph 
$\gr_\phi$ is a $G_\delta$ in \mX that intersects each fiber of $f$ in exactly one point.  
Consequently, $\gr_\phi$ is completely metrizable, thus is a $G_\delta$ in $C \times [0,1]$.  
Note too that $\gr_\phi$ is totally disconnected and intersects each element of $\sC$ in a set of cardinality $2^\omega$.  
Claim: $\dim \gr_\phi = 1$.  \
In fact, by Proposition 12, $\dim \gr_\phi \leq \dim C + \dim [0,1] =$ $0 + 1 = 1$.
%%----------------------------------------------------------------------------------------------38
To see that $\dim \gr_\phi > 0$, write $q$ for the horizontal projection, put
$
\begin{cases}
\ A = \gr_\phi \cap q^{-1}([0,1/7])\\
\ B = \gr_\phi \cap q^{-1}([6/7,1])
\end{cases}
$
and let \mU be any open subset of $\gr_\phi$: 
$
\begin{cases}
\ A \subset U\\
\ B \cap \ov{U} = \emptyset
\end{cases}
$
$-$then $\#(\fr U) = 2^\omega$.
\vspi
[Note: \ Working instead with $[0,1]^{n+1} = [0,1] \times [0,1]^n$, one can modify the  preceding construction and produce an example of a second countable completely metrizable totally disconnected space of topological dimension $n$.  
Such a space cannot contain a compact Cantor $n$-space (cf. p. \pageref{19.25}).]\\
\endgroup %%------------------------------------<<

\begingroup%%----------------------------------->>
\fontsize{9pt}{11pt}\selectfont 
\textbf{\small FACT} \  
Let \mX and \mY be nonempty CRH spaces.  Suppose that $X \times Y$ is strongly paracompact $-$then 
$\dim(X \times Y) \leq \dim X + \dim Y$.
\vspi
[View $X \times Y$ as a subspace of $\beta X \times \beta Y$ to get 
$\dim (X \times Y) \leq \dim (\beta X \times \beta Y)$ (cf. p. \pageref{19.26}), which is 
$\leq \dim \beta X + \dim \beta Y$ (cf. Proposition 12) or still, $\leq \dim X + \dim Y$ (cf. Proposition 11).]
\vspi
[Note: \ Is it sufficient that $X \times Y$ be paracompact?  The answer is unknown.]\\
\endgroup %%------------------------------------<<

\begingroup%%----------------------------------->>
\fontsize{9pt}{11pt}\selectfont 
Application: Suppose that \mX and \mY are second countable and metrizable 
$-$then $\dim(X \times Y) \leq \dim X + \dim Y$.\\
\endgroup %%------------------------------------<<

\label{19.42}
\begingroup%%----------------------------------->>
\fontsize{9pt}{11pt}\selectfont 
\textbf{\small EXAMPLE} \  
Take for \mX the subspace of $l^2$ consisting of all sequences $\{x_n\}$, with $x_n$ rational $-$then $\dim X = 1$.  
But \mX is homeomorphic to $X \times X$ , so $\dim(X \times X) = 1$, which is $< 2 = \dim X + \dim X$.
\vspi
[Note: \ Given any $n \in \N$, there exists an $X \subset \R^{n+1}$ such that 
$\dim X = \dim (X \times X) = n$ 
(Anderson-Keisler\footnote[2]{\textit{Proc. Amer. Math. Soc.} \textbf{18} (1967), 709-713.}).]\\
\endgroup %%------------------------------------<<

\begingroup%%----------------------------------->>
\fontsize{9pt}{11pt}\selectfont 
\textbf{\small FACT} \  
Let \mX and \mY be nonempty CRH spaces.  Suppose that \mX and \mY are infinite and $X \times Y$ is pseudocompact 
$-$then $\dim(X \times Y) \leq \dim X + \dim Y$.
\vspi
[Glicksberg's theorem says that if \mX and \mY are infinite CRH spaces, then the product $X \times Y$ is pseudocompact iff $\beta (X \times Y) = \beta X \times \beta Y$, the equal sign meaning that the two compactifications of $X \times Y$ are equivalent (and not just homeomorphic).  Recall that the product of two pseudocompact spaces need not be pseudocompact but this will be the case if one of the factors is compactly generated.  
Example: $\dim ([0,\Omega[ \times [0,\Omega]) = 0$.]\\
\endgroup %%------------------------------------<<

\begin{proposition} \ %13
Suppose that \mX is a CW complex and \mY is compact $-$then $\dim(X \times Y) \leq \dim X + \dim Y$.
\end{proposition}

[Argue by induction on $\dim X$.  
There is nothing to prove if $\dim X = 0$.  
If $\dim X > 0$, then, since the combinatorial and topological dimensions of \mX coincide 
(cf. p. \pageref{19.27}), $X = X^{(n)}$.  
Thus one can write 
$X = X^{(n-1)} \cup \bigcup\limits_1^\infty A_j$ where each $A_j$ is closed and expressible 
%%----------------------------------------------------------------------------------------------39
as a disjoint union  
$\bigcup\limits_i K_{i,j}$, $\{K_{i,j}\}$ being a discrete collection of compacta, with $\dim K_{i,j} \leq n$.  
From the induction hypothesis, 
$\dim (X^{(n-1)} \times Y) \leq$ 
$\dim X^{(n-1)} + \dim Y \leq$ 
$n - 1 + m$.  On the other hand, Proposition 12 implies that 
$\dim(K_{i,j} \times Y) \leq$ 
$\dim K_{i,j} + \dim Y \leq$ 
$n + m$, so $\dim (A_j \times Y) \leq n + m$.  Now apply the countable union lemma.]\\

\index{Stacking Lemma}
\textbf{\small STACKING LEMMA} \  
Let \mX and \mY be nonempty CRH spaces.  
Suppose that \mY is compact $-$then for every numerable open covering $\sW$ of 
$X \times Y$, there exists a numerable open covering 
$\sU = \{U_i:i \in I\}$ of \mX and $\forall \ i \in I$, a finite open covering 
$\sV_i = \{V_{i,j}:j \in J_i\}$ of \mY such that the collection $\{U_i \times \sV_i: i \in I\}$ refines $\sW$.

[The assertion is trivial if \mX is paracompact.  
In general, there exists a metric space \mZ, an open covering $\sZ$ of \mZ, 
and a continuous function $f:X \ra Y \ra Z$ such that $f^{-1}(\sZ)$ refines $\sW$ 
(cf p. \pageref{19.28}).  
Define 
$e:C(Y,Z) \times Y \ra Z$ by $e(\phi,y) = \phi(y)$ $-$then $e^{-1}(\sZ)$ is a numerable open covering of 
$C(Y,Z) \times Y$.  
Since $C(Y,Z) \times Y$ is paracompact, one can find a numerable open covering 
$\sO = \{O_i:i \in I\}$ of $C(Y,Z)$ and $\forall \ i \in I$, 
a finite open covering $\sV_i = \{V_{i,j}: j \in J_i\}$ of \mY such that the collection $\{O_i \times \sV_i: i \in I\}$ refines $e^{-1}(\sZ)$.  
Put 
$F(x)(y) = f(x,y)$: $F \in C(X,C(Y,Z))$ $\&$ $f = e \circx (F \times \id_Y)$.  
Consider $\sU = \{U_i:i \in I\}$, where 
$U_i = F^{-1}(O_i)$.]

[Note: \ The complete regularity of \mX plays no role in the proof.]\\

To establish the product theorem, first employ the countable union lemma and make the obvious reductions to the case when \mY is compact.  
This done, let $\sW$ be a finite open covering of $X \times Y$.  
According to the stacking lemma, there exists a neighborhood finite open covering 
$\sU = \{U_i:i \in I\}$ of \mX and for each $i \in I$, a finite open covering 
$\sV_i = \{V_{i,j}: j \in J_i\}$ of \mY such that the collection $\{U_i \times \sV_i: i \in I\}$ refines $\sW$.  
Fix a precise open refinement $\sO = \{O_i: i \in I\}$ of $\sU$ of order $\leq n + 1$ (cf Proposition 6) $-$then $\dim\abs{N(\sO)} \leq n$, 
$N(\sO)$ the nerve of $\sO$.  
Choose an $\sO$-map $f$, i.e., a continuous function $f:X \ra \abs{N(\sO)}$ with the property that 
$\forall \ O_i \in \sO$: $(b_{O_i} \circx f)^{-1} (]0,1]) \subset O_i$ 
(cf. p. \pageref{19.29}).  
Put $F = f \times \id_Y$.  
Since $\dim(\abs{N(\sO)} \times Y) \leq n + m$ (cf. Proposition 13), the open covering 
$\{b_{O_i}^{-1}(]0,1]) \times \sV_i: i \in I\}$ of $\abs{N(\sO)} \times Y$ has an open refinement $\sP$ of order 
$\leq n + m + 1$.  Consider $F^{-1}(\sP)$.\\

\begingroup%%----------------------------------->>
\fontsize{9pt}{11pt}\selectfont 
The product theorem holds if \mX is merely completely regular.  Indeed, once the reductions to the case ``\mY compact'' have been carried out, the argument proceeds as when \mX is normal.  The reductions depend in turn on the countable union lemma which retains its validity in the completely regular situation provided the subspaces in question have the EP w.r.t. $[0,1]$ 
(cf. p. \pageref{19.30}).  Two results are relevant for the transition.\\
\endgroup %%------------------------------------<<

%%----------------------------------------------------------------------------------------------40
\begingroup%%----------------------------------->>
\fontsize{9pt}{11pt}\selectfont 
\textbf{\small LEMMA} \  
Let \mX be a topological space.  Let \mB be a compact subspace of a CRH space \mY $-$then $X \times B$, as a subspace of 
$X \times Y$, has the EP w.r.t. $[0,1]$.

[Recalling that $B \subset Y$ has the EP w.r.t. $[0,1]$ (cf p. \pageref{19.31}), let $\sO$ be a finite numerable open covering of $X \times B$.  Use the stacking lemma and construct a numerable open covering $\sW$ of $X \times Y$ such that 
$\sW \cap (X \times B)$ is a refinement of $\sO$.  Apply $\S 6$, Proposition 4 (the proof of sufficiency does not require a cardinality assumption on $\sW$).]\\
\endgroup %%------------------------------------<<

\begingroup%%----------------------------------->>
\fontsize{9pt}{11pt}\selectfont 
\textbf{\small LEMMA} \  
Let \mX be a topological space.  Let \mB be a closed subspace of a paracompact LCH space \mY $-$then $X \times B$, as a subspace of 
$X \times Y$, has the EP w.r.t. $[0,1]$.
\vspi
[Note: \ Paracompactness  of \mY alone is not enough.  
Example: Take $X = \PP$, $Y =$ Michael line and  $B = \Q$ $-$then $X \times B$, as a subspace of $X \times Y$, does not have the EP w.r.t. $[0,1]$.  One can, however, drop local compactness if some other assumption on \mY is imposed, e.g., stratifiability.]\\
\endgroup %%------------------------------------<<

\begingroup%%----------------------------------->>
\fontsize{9pt}{11pt}\selectfont 
Its utility notwithstanding, there are limitations to the product theorem.  
For example, it is not necessarily applicable if both factors are metrizable.  
However, this possibility (and others) can be readily placed in a general framework.
\vspi
Let \mX and \mY be nonempty CRH spaces $-$then a 
\un{cozero set rectangle}
\index{cozero set rectangle} in 
$X \times Y$ is a set of the form $U \times V$, where
$
\begin{cases}
\ U\\
\ V
\end{cases}
$
is a cozero set in 
$
\begin{cases}
\ X\\
\ Y
\end{cases}
$
.\\[.5cm]
\endgroup %%------------------------------------<<

\begingroup%%----------------------------------->>
\fontsize{9pt}{11pt}\selectfont 
\textbf{\small LEMMA} \  
$X \times Y$ is $\sZ$-embedded in $X \times \beta Y$ iff every cozero set in $X \times Y$ can be written as the union of a collection of cozero set rectangles $U \times V$, where $\{U\}$ is $\sigma$-neighborhood finite.
\vspi
[Use the stacking lemma and the fact that the union of a $\sigma$-neighborhood finite collection of cozero sets is a cozero set.]
\vspi
[Note: \ $X \times Y$ is $\sZ$-embedded in $\beta X \times \beta Y$ iff every cozero set in $X \times Y$ can be written as the union of a countable collection of cozero set rectangles in $U \times V$.]\\
\endgroup %%------------------------------------<<

\begingroup%%----------------------------------->>
\fontsize{9pt}{11pt}\selectfont 
The following conditions are equivalent.
\\
\indent\indent (a) Every cozero set in $X \times Y$ can be written as the union of a collection of cozero set rectangles $U \times V$, where 
$\{U\}$ is $\sigma$-neighborhood finite.
\\
\indent\indent (b) Given any $f \in C(X \times Y)$ and any $\epsilon > 0$, there exists a covering of $X \times Y$ 
by cozero set rectangles $U \times V$ such that $\osc(\restr{f}{U \times V}) < \epsilon$ and $\{U\}$ is 
$\sigma$-neighborhood finite.
\vspi
[(a) $\implies$ (b):  Fix a sequence of open intervals $]a_n,b_n[$, each of length $< \epsilon/2$: 
$\R = \ds\bigcup\limits_1^\infty \ ]a_n,b_n[$ \ $-$then 
$X \times Y = \ds\bigcup\limits_1^\infty f^{-1}(]a_n,b_n[)$.  Write
$f^{-1}(]a_n,b_n[)$ as the union of a collection of cozero set rectangles 
$U_i \times V_i$, where $\{U_i: i \in I_n\}$ is $\sigma$-neighborhood finite.  
Obviously, $\osc(\restr{f}{U_i \times V_i}) < \epsilon$ and 
$\ds\bigcup\limits_1^\infty \{U_i: i \in I_n\}$ is 
$\sigma$-neighborhood finite.
\vspi
%%----------------------------------------------------------------------------------------------
(b) $\implies$ (a):  Take an $f \in C(X \times Y)$.  
Pick a cozero set rectangle covering $\sW_n = \{U \times V\}$ of 
$X \times Y$ such that $\osc(\restr{f}{U \times V}) < 1/n$ and $\{U\}$ is $\sigma$-neighborhood finite.  
Denote by $\sW_n(f)$ the subset of $\sW_n$ consisting of the $U \times V$ that are contained in $X \times Y - Z(f)$ 
$-$then 
$\ds\bigcup\limits_1^\infty \sW_n(f)$ covers $X \times Y - Z(f)$.]\\
\vspace{0.25cm}
\endgroup %%------------------------------------<<

\begingroup%%----------------------------------->>
\fontsize{9pt}{11pt}\selectfont 
Assume: Every open subset of 
$
\begin{cases}
\ X\\
\ Y
\end{cases}
$ 
is $\sZ$-embedded in 
$
\begin{cases}
\ X\\
\ Y
\end{cases}
$
$-$then (a) and (b) above are equivalent to the following conditions.
\\
\indent\indent $(a)_{\sZ}$  Every cozero set in $X \times Y$ can be written as the union of a collection of open rectangles 
$U \times V$, where $\{U\}$ is $\sigma$-neighborhood finite.
\\
\indent\indent $(b)_{\sZ}$  Given any $f \in C(X \times Y)$ and any $\epsilon > 0$, there exists a covering of $X \times Y$ by open rectangles $U \times V$ such that 
$\osc(\restr{f}{U \times V}) < \epsilon$ and $\{U\}$ is $\sigma$-neighborhood finite.
\vspi
[That (a) $\implies$ $(a)_{\sZ}$ is clear, as is $(a)_{\sZ}$ $\implies$ $(b)_{\sZ}$.  To prove that 
$(b)_{\sZ}$ $\implies$ (b), let $f \in C(X \times Y)$ and $\epsilon > 0$ but with 
$\osc(\restr{f}{U \times V}) < \epsilon/2$.  The assumption on 
$
\begin{cases}
\ X\\
\ Y
\end{cases}
$
implies that the interior of 
$
\begin{cases}
\ \ov{U}\\
\ \ov{V}
\end{cases}
$
is a cozero set in 
$
\begin{cases}
\ X\\
\ Y
\end{cases}
$
.  
The corresponding collection of cozero set rectangles thereby produced covers $X \times Y$ and the oscillation of $f$ on any one of them is $< \epsilon$.]\\
\endgroup %%------------------------------------<<

\begingroup%%----------------------------------->>
\fontsize{9pt}{11pt}\selectfont 
In a CRH space, every open subset is $\sZ$-embedded iff every open subset which is the interior of its closure is cozero.  The latter property is evidently a weakening of perfect normality and, e.g., is possessed by an arbitrary product of metrizable spaces (\u S\u cepin\footnote[2]{\textit{Soviet Math. Dokl.} \textbf{17} (1976), 152-155; 
see also Blair-Swardson, \textit{Topology Appl.} \textbf{36} (1990), 73-92.}) but not by $[0,\Omega[$ \ or $\beta \R$.\\
\endgroup %%------------------------------------<<

\begingroup%%----------------------------------->>
\fontsize{9pt}{11pt}\selectfont 
\textbf{\small LEMMA} \  
Suppose that \mX is metrizable and that every open subset of \mY is $\sZ$-embedded in \mY $-$then $X \times Y$ is 
$\sZ$-embedded in $X \times \beta Y$.
\vspi
[It suffices to check $(b)_{\sZ}$, so let $f \in C(X \times Y)$ and $\epsilon > 0$.  \ 
Enumerate $\Q$: $\{q_n\}$ and put 
$I_n = ]q_n - \epsilon/3,q_n + \epsilon/3[$.  \ 
Fix a $\sigma$-neighborhood finite basis $\{U\}$ for \mX. \ 
Let $Y(U,n)$ be the subset of \mY made up of those points which admit a  neighborhood 
\mV: $f(U \times V) \subset I_n$ $-$then $Y(U,n)$ is open in \mY, 
$\osc(\restr{f}{U \times Y(U,n)}) < \epsilon$, and since $\forall \ (x,y) \in X \times Y$ $\exists$ $q_n \in \Q$: 
$\abs{f(x,y) - q_n} < \epsilon/6$, the open rectangles $U \times Y(U,n)$ cover $X \times Y$.]\\
\endgroup %%------------------------------------<<

\begingroup%%----------------------------------->>
\fontsize{9pt}{11pt}\selectfont 
\textbf{\small FACT} \  
Let \mX and \mY be nonempty CRH spaces.  Suppose that $X \times Y$ is $\sZ$-embedded in $X \times \beta Y$ $-$then 
$\dim (X \times Y) \leq \dim X + \dim Y$.
\vspi
[Simply note that 
$\dim (X \times Y) \leq \dim (X \times \beta Y)$ 
(cf. p. \pageref{19.32}), which, by the product theorem, is $\leq \dim X + \dim \beta Y = \dim X + \dim Y$.]\\
\endgroup %%------------------------------------<<

\begingroup%%----------------------------------->>
\fontsize{9pt}{11pt}\selectfont 
Application: Suppose that \mX and \mY are metrizable $-$then $\dim(X \times Y) \leq \dim X + \dim Y$.\\
\endgroup %%------------------------------------<<

%%----------------------------------------------------------------------------------------------42
\begingroup%%----------------------------------->>
\fontsize{9pt}{11pt}\selectfont 
\textbf{\small EXAMPLE} \  
Let \mX and \mY be nonempty M complexes $-$then $X \times_k Y$ is an M complex and 
$\dim(X \times_k Y) \leq \dim X + \dim Y$.
\vspi
[Assume first that \mX is an M$_n$ space and \mY is an M$_m$ space, proceed by induction on $n+m$.]\\
\endgroup %%------------------------------------<<

\begingroup%%----------------------------------->>
\fontsize{9pt}{11pt}\selectfont 
That dim is monotonic on $\sZ$-embedded subspaces is the key to the preceding method.  But one can get away with even less.  In general, a subspace \mA of a topological space \mX is said to be 
\un{weakly $\sZ$-embedded}
\index{weakly $\sZ$-embedded} 
in \mX if for any cozero set \mO in \mA there exists a $\sigma$-neighborhood finite collection $\{O_i:i \in I\}$ of cozero sets $O_i$ in \mA, each of which is the intersection of \mA with a cozero set in \mX, such that $O = \ds\bigcup\limits_i O_i$.\\
\endgroup %%------------------------------------<<

\begingroup%%----------------------------------->>
\fontsize{9pt}{11pt}\selectfont 
\textbf{\small LEMMA} \  
Let \mX be a nonempty CRH space.  Suppose that \mA is a weakly $\sZ$-embedded subspace of \mX $-$then 
$\dim A \leq \dim X$.\\
\endgroup %%------------------------------------<<

\begingroup%%----------------------------------->>
\fontsize{9pt}{11pt}\selectfont 
Let \mX and \mY be nonempty CRH spaces $-$then $X \times Y$ is said to be 
\un{rectangular}
\index{rectangular (CRH spaces)} 
if every cozero set in $X \times Y$ can be written as the union of a 
$\sigma$-neighborhood finite collection of cozero set rectangles $U \times V$.  
If $X \times Y$ is $\sZ$-embedded in $X \times \beta Y$, 
then $X \times Y$ is rectangular (the converse is false).\\
\endgroup %%------------------------------------<<

\begingroup%%----------------------------------->>
\fontsize{9pt}{11pt}\selectfont 
\textbf{\small EXAMPLE} \  
Suppose that \mX and \mY are paracompact Hausdorff spaces satisifying Arhangel'ski\u i's condition 
$-$then $X \times Y$ is rectangular.\\
\endgroup %%------------------------------------<<

\begingroup%%----------------------------------->>
\fontsize{9pt}{11pt}\selectfont 
\textbf{\small FACT} \  
Let \mX and \mY be nonempty CRH spaces. Suppose that $X \times Y$ is rectangular $-$then 
$\dim(X \times Y) \leq \dim X + \dim Y$.
\vspi
[Indeed $X \times Y$ as a subspace of $\beta X \times \beta Y$ is weakly $\sZ$-embedded.]\\
\endgroup %%------------------------------------<<

\begingroup%%----------------------------------->>
\fontsize{9pt}{11pt}\selectfont 
\textbf{\small EXAMPLE} \  
Rectangularity of $X \times Y$ is not a necessary condition for the validity of the relation 
$\dim(X \times Y) \leq \dim X + \dim Y$.
\index{The Sorgenfrey Plane (example)}
\\
\indent\indent (1) \ (\un{The Sorgenfrey Plane}) 
Let \mX be the Sorgenfrey line $-$then \mX is zero dimensional and Lindel\"of, hence $\dim X = 0$ (cf. Proposition 2).  The Sorgenfrey plane $X \times X$ is zero dimensional but not normal and is ``asymmetrical'' in that every line with negative slope is discrete but every line with positive slope is homeomorphic to \mX.  
Moreover, it is not rectangular as may be seen by considering points on or above the line $x + y = 1$.  Still, $\dim (X \times X) = 0$.  As a preliminary, show that if \mO is any open subset of $X \times X$, then there exists a sequence of clopen sets $O_n$ such that 
$O \subset \ds\bigcup\limits_n O_n$ $\subset \ov{O}$ and from this deduce that every cozero set in $X \times X$ is a countable union of clopen sets (cf. p. \pageref{19.33}).
\\
\indent\indent (2) \ 
(\un{The Michael Line $\times$ the Irrationals})
\index{The Michael Line $\times$ the Irrationals (example)} 
Let \mX be the Michael line $-$then \mX is hereditarily paracompact, hence hereditarily normal, so it follows from the control lemma that $\dim X = 0$.  The product $X \times \PP$ is zero dimensional but not normal.  
Nor is it rectangular: Otherwise, 
$\PP$ would be an $F_\sigma$ in $\R$.  However, one an show that $\dim (X \times \PP) = 0$.\\
\endgroup %%------------------------------------<<

%%----------------------------------------------------------------------------------------------43
\begingroup%%----------------------------------->>
\fontsize{9pt}{11pt}\selectfont 
Let \mX and \mY be nonempty CRH spaces  $-$then $X \times Y$ is said to be 
\un{piecewise rectangular}
\index{piecewise rectangular} 
if every cozero set in $X \times Y$ can be written as the union of a 
$\sigma$-neighborhood finite collection $\{W\}$, where each \mW is a clopen subset of some cozero set rectangle $U \times V$.  In this terminlology, 
Pasynkov\footnote[2]{\textit{London Math. Soc. Lecture Notes} \textbf{93} (1985), 227-250.}
proved that if 
$
\begin{cases}
\ \dim X = 0\\
\ \dim Y = 0
\end{cases}
, \ 
$
then $\dim (X \times Y) = 0$ iff $X \times Y$ is piecewise rectangular.
\vspi
[Note: \ For every pair of positive integers $(n,m)$, 
Tsuda\footnote[3]{\textit{Canad. Math. Bull.} \textbf{30} (1987), 49-56.}
has constructed a normal 
$
\begin{cases}
\ X: \dim X = n\\
\ Y: \dim Y = m
\end{cases}
$
for which $X \times Y$ is also normal with $\dim (X \times Y) = n + m$ but such that $X \times Y$ is not 
piecewise rectangular.]\\
\endgroup %%------------------------------------<<

\begingroup%%----------------------------------->>
\fontsize{9pt}{11pt}\selectfont 
\textbf{\small EXAMPLE} \ [Assume CH] \  
There exist nonempty perfectly normal locally compact
$
\begin{cases}
\ X\\
\ Y
\end{cases}
: \ 
$
$X \times Y$ is a perfectly normal LCH space and $\dim X + \dim Y < \dim (X \times Y)$.  For this, use the notation of the example following Proposition 12, letting $\Delta_C$ be the diagonal of \mC in $C^2$, which will then be identified with \mC when convenient.  Transfer the topology on $\gr_\phi$ back to \mC to get a second countable completely metrizable topology $\tau_\phi$ on \mC finer than the euclidean topology $\tau$.
\vspi
Claim: There exists a second countable metrizable topology $\Lambda$ on $C^2$  finer then the euclidean topology 
$\tau^2$ with $\restr{\Lambda}{\Delta_C} = \tau_\phi$ $\&$ 
$\restr{\Lambda}{C^2 - \Delta_C} = \restr{\tau^2}{C^2 - \Delta_C}$ 
such that every element of $\Lambda$ containing a point $(x,x) \in \Delta_C$ also contains the intersection with $C^2$ of two disjoint open disks, tangent to $\Delta_C$ at $(x,x)$.
\vspi
[Fix a countable basis $\{U_i\}$ for $\tau_\phi$.  
Since $\phi$ is Baire one, each $U_i$ is a euclidean $F_\sigma$: 
$U_i = \ds\bigcup\limits_1^\infty A_{ij}$, $A_{ij}$ $\tau$-closed.  
Enumerate the $A_{ij}$: $\{K_n\}$.  Given $r > 0$, let 
$K_n(r)$ be the union of all $B_r \cap C^2$, where $B_r$ is an open disk of radius $r$ tangent to $\Delta_C$ at some point of 
$K_n$.  
Recursively determine a sequence of positive real numbers $r_n$: $r_n > r_{n+1}$ $\&$ $\lim r_n = 0$, subject to 
$K_n \cap K_m = \emptyset$ $\implies$ $K_n(r_n) \cap K_m(r_m) = \emptyset$.  Put 
$O_i = \ds\bigcup \{K_n(2^{-i}r_n):K_n \subset U_i\}$.  
Consider the topology on $C^2$ generated by the $O_i$ and a countable basis for the euclidean topology on $C^2 - \Delta_C$.]
\vspi
Construct Kunen modifications $\tau^\prime$ and $\tau\pp$ of $\tau$ such that 
$\tau^\prime \times \tau\pp$ is a perfectly normal locally compact topology finer than $\Lambda$ whose restriction 
$\tau^\prime \times \restr{\tau\pp}{\Delta_C}$ is a Kunen modification of $\tau_\phi$ 
(cf. p. \pageref{19.34}).  
In so doing, work with an enumeration 
$\{x_\alpha: \alpha < \Omega\}$ of \mC, letting $\{C_\alpha: \alpha < \Omega\}$ be an enumeration of the countable subsets of $C^2$ such that $\forall \ \alpha$: $C_\alpha \subset \{x_\beta: \beta < \alpha\}^2$.  While $\tau^\prime \times \tau\pp$ is not a Kunen modification of $\Lambda$, local compactness is, of course, automatic.  As for perfect normality, the essential preliminary is that $\forall \ S \subset C^2$ $\exists$ $\alpha < \Omega$: 
$\cl_\Lambda(S) \cap \{x_\beta: \beta > \alpha\}^2 =$ 
$\cl_{\tau^\prime \times \tau\pp}(S) \cap \{x_\beta: \beta > \alpha\}^2$.  
This said, let $S \subset C^2$ be 
$\tau^\prime \times \tau\pp$-closed and choose a sequence $\{O_n\}$ of $\Lambda$-open sets: 
$\cl_\Lambda(S) = \ds\bigcap\limits_n O_n =$ 
$\ds\bigcap\limits_n \cl_\Lambda(O_n)$ $-$then $\exists$ $\alpha < \Omega$: 
$\cl_\Lambda(S) \cap \{x_\beta: \beta > \alpha\}^2 =$
$\ds\bigcap\limits_n \cl_{\tau^\prime \times \tau\pp}(O_n) \cap \{x_\beta: \beta > \alpha\}^2$.  On the other hand, for each 
$\beta \leq \alpha$ there are countable collections 
$
\begin{cases}
\ \{P_n^\prime(\beta)\}\\
\ \{P_n\pp(\beta)\}
\end{cases}
$
of $\tau^\prime \times \tau\pp$-open sets: 
 $
\begin{cases}
\ (C \times \{x_\beta\}) \cap (C^2 - S) \subset \bigcup\limits_n P_n^\prime(\beta)\\
\ (\{x_\beta\} \times C) \cap (C^2 - S) \subset \bigcup\limits_n P_n\pp(\beta)
\end{cases}
$
%%----------------------------------------------------------------------------------------------44
$\&$
$
\begin{cases}
\ \cl_{\tau^\prime \times \tau\pp}(P_n^\prime(\beta)) \cap S = \emptyset\\
\ \cl_{\tau^\prime \times \tau\pp}(P_n\pp(\beta)) \cap S = \emptyset
\end{cases}
. \ 
$
Form 
$
\begin{cases}
\ O_n^\prime(\beta) = C^2 - \cl_{\tau^\prime \times \tau\pp}(P_n^\prime(\beta))\\
\ O_n\pp(\beta) = C^2 - \cl_{\tau^\prime \times \tau\pp}(P_n\pp(\beta))
\end{cases}
$
and combine the 
$
\begin{cases}
\ O_n^\prime(\beta)\\
\ O_n\pp(\beta)
\end{cases}
(\beta \leq \alpha)
$
with the $O_n$ to obtain a single countable collection $\{U_n\}$ of $\tau^\prime \times \tau\pp$-open sets: 
$S = \ds\bigcap\limits_n U_n =$ $\ds\bigcap\limits_n \cl_{\tau^\prime \times \tau\pp}(U_n)$.
\vspi
Claim: Let 
$
\begin{cases}
\ X = (C,\tau^\prime)\\
\ Y = (C,\tau\pp)
\end{cases}
$ 
$-$then 
$
\begin{cases}
\ \dim X = 0\\
\ \dim Y = 0
\end{cases}
$ 
(cf. p. \pageref{19.35}) and $\dim (X \times Y) > 0$.
\vspi
[It is enough to show that $\Delta_C \subset (C \times C,\tau^\prime \times \tau\pp)$ has positive topological dimension.  Return to \mC, which thus carries three topologies, namely $\tau$, $\tau_\phi$, and 
$\tau^* \equiv \tau^\prime \times \restr{\tau\pp}{C}$, a Kunen modification of $\tau_\phi$.  Let
$
\begin{cases}
\ A = \phi^{-1}([0,1/7])\\
\ B = \phi^{-1}([6/7,1])
\end{cases}
; \ 
$
let
$
\begin{cases}
\ A^* = \phi^{-1}([0,1/3])\\
\ B^* = \phi^{-1}([2/3,1])
\end{cases}
. \ 
$
To arrive at a contradiction, suppose that $O^*$ is a $\tau^*$-clopen set:
$
\begin{cases}
\ A^* \subset O^*\\
\ B^* \cap O^* = \emptyset
\end{cases}
. \ 
$
If the bar denotes closure in $\tau_\phi$ and if $V = C - \ov{O^*}$, then 
$
\begin{cases}
\ A \subset \ov{V} = \emptyset\\
\ B \subset V
\end{cases}
$
$\&$ $\#(\fr(V) > \omega$.  But $\fr V \subset \ov{O^*} \cap \ov{C - O^*}$ and 
$\#( \ov{O^*} \cap \ov{C - O^*}) \leq \omega$.]
\vspi
[Note: CH is not necessary here.  Examples of this type exist in ZFC 
(Przymusi\'nski\footnote[2]{\textit{Proc. Amer. Math. Soc.} \textbf{76} (1979), 315-321; 
see also Tsuda, \textit{Math. Japon.} \textbf{27} (1982), 177-195.}), the main difference being that the product $X \times Y$ is not perfectly normal but rather is a normal countably paracompact LCH space.]\\
\endgroup %%------------------------------------<<

\label{19.50}
\begingroup%%----------------------------------->>
\fontsize{9pt}{11pt}\selectfont 
One final point: The product theorem holds if \mX is an arbitrary nonempty topological space.  In fact, if $A \subset X$ has the EP w.r.t. $[0,1]$, then its image $\crg A$ in $\crg X$ ``is'' the complete regularization of \mA and as such has the EP w.r.t.$ [0,1]$, so 
$\dim A =$ 
$\dim \crg A \leq$ 
$\dim \crg X = \dim X$ (cf. p. \pageref{19.36}).  
The countable union lemma is therefore applicable provided that the
$A_j \subset X$ have EP w.r.t. $[0,1]$ (cf. p. \pageref{19.37}).  
It is then easy to fall back to the completely regular case since for any LCH space \mY, $\crg(X \times Y) = \crg X \times Y$.\\
\endgroup %%------------------------------------<<

\textbf{\small LEMMA} \  
Suppose that \mX is a compact Hausdorff space.  Let $f,g \in C(X,\bS^n)$ and put 
$D = \{x: f(x) \neq g(x)\}$.  Assume: $\dim D \leq n - 1$ $-$then $f \simeq g$.

[Since $ID$ is an $F_\sigma$ in $IX$, hence is normal, it follows from the product theorem that 
$\dim ID \leq n$.  Set $Y = i_0X \cup I(X - D) \cup i_1X$ and define $h:Y \ra \bS^n$ by 
$
\begin{cases}
\ h(x,0) = f(x)\\
\ h(x,1) = g(x)
\end{cases}
\&
$
$h(x,t) = f(x) = g(x)$ $-$then $h$ is continuous and has a continuous extension $H \in C(IX,\bS^n)$ 
(cf. p. \pageref{19.38}).]\\

\begin{proposition} \ %14
Let $f,g \in C(X,\bS^n)$ and put $D = \{x: f(x) \neq g(x)\}$.  
Assume: $\dim D \leq n-1$ $-$then $f \simeq g$.
\end{proposition}

[The subset of $\beta X$ on which $\beta f \neq \beta g$ can be written as a countable union 
$\bigcup\limits_1^\infty \ov{U_j}$, each $U_j$ being open in $\beta X$.  And: 
$\dim (\ov{U_j} \cap X) \leq n-1$ $\implies$ $\dim \ov{U_j} = \dim \beta(\ov{U_j} \cap X) \leq n - 1$ $\implies$ 
$\dim \bigcup\limits_1^\infty \ov{U_j} \leq n - 1$, thus from the lemma, $\beta f \simeq \beta g$.]\\

%%----------------------------------------------------------------------------------------------45
Application: If $\dim X \leq n - 1$, then $[X,\bS^n] = *$.\\

\begingroup%%----------------------------------->>
\fontsize{9pt}{11pt}\selectfont
\textbf{\small FACT} \  
Suppose that \mX is normal and $\dim X$ is finite $-$then the natural map 
$[\beta X,\bS^n] \ra [X,\bS^n]$ is bijective if $n > 1$ but if $n = 1$ and \mX is connected, there is an exact sequence 
$0 \ra C(X)/BC(X) \ra$ $[\beta X,\bS^1]$ $\ra [X,\bS^1] \ra 0$.
\vspi
[To discuss the second assertion, observe that \mX is connected iff $\beta X$ is connected and form the commutative diagram 
\[
\begin{tikzcd}%[ sep=small]
{0} \ar{r}
&{C(\beta X)/\Z} \ar{d} \ar{r}{\exp} 
&{C(\beta X, \bS^1)} \ar{d} \ar{r}
&{[\beta X,\bS^1]} \ar{d} \ar{r}
&{0}\\
{0} \ar{r}
&{C(X)/\Z} \ar{r}[swap]{\exp}
&{C(X, \bS^1)} \ar{r}
&{[X,\bS^1]} \ar{r}
&{0}
\end{tikzcd}
.
\]
Since the rows are exact and the middle vertical arrow is an isomorphism, the ker-coker lemma gives 
$\ker([\beta X,\bS^1] \ra [X,\bS^1])$ $\approx$ 
$\coker (C(\beta X)/\Z \ra C(X)/\Z)$ $\approx$ $C(X)/BC(X)$.  
As for the need of the connectedness assumption, take $X = \N$: $\dim \N = 0$ $\implies$ 
$[\N,\bS^1] = * = [\beta\N,\bS^1]$.]
\vspi
[Note: \ The exact sequence 
$0 \ra C(X)/BC(X) \ra [\beta X,\bS^1] \ra [X,\bS^1] \ra 0$ translates to 
$0 \ra $
$C(X)/BC(X) \ra $ 
$\cH{}^1(\beta X) \ra $ 
$\cH{}^1(X)$ 
$\ra 0$.  
Because the quotient $C(X)/BC(X)$ is torsion free and divisible when nontrivial, 
it follows that if \mX is not pseudocompact, then 
$\cH{}^1(\beta X) \approx \oplus \Q$ and is in fact uncountable.  
Proof: Let $f:X \ra \R$ be an unbounded continuous function, put $f_r = r\cdot f$ $(r \in \R)$ and consider $f_r + BC(X)$.  
Example: $\cH{}^1(\beta \R) \approx  C(\R)/BC(\R)$.]\\
\endgroup %%------------------------------------<<

\begingroup%%----------------------------------->>
\fontsize{9pt}{11pt}\selectfont
Let \mY be a connected CW space $-$then 
Bartik\footnote[2]{\textit{Quart. J. Math.} \textbf{29} (1978), 77-91; 
see also Calder-Siegel, \textit{Trans. Amer. Math. Soc.} \textbf{235} (1978), 245-270 
and \textit{Proc. Amer. Math. Soc.} \textbf{78} (1980), 288-290.}
has shown that the arrow $[\beta X,Y] \ra [X,Y]$ is bijective for every nonempty CRH space \mX with $\dim X$ finite iff 
$\pi_1(Y)$ is finite and $\forall \ q > 1$, $\pi_q(Y)$ is finitely generated or still, iff $\pi_1(Y)$ is finite and \mY has the homotopy type of a connected CW complex \mK such that $\forall \ n$, $K^{(n)}$ is finite 
(cf. p. \pageref{19.39}).\\
\endgroup %%------------------------------------<<

\begingroup%%----------------------------------->>
\fontsize{9pt}{11pt}\selectfont
Application: Suppose that $\pi$ is a finitely generated abelian group.  Let \mX be a nonempty CRH space of finite topological dimension $-$then $\forall \ n > 1$, $\cH{}^n(\beta X;\pi) \approx \cH{}^n(X;\pi)$.\\
\endgroup %%------------------------------------<<

\begingroup%%----------------------------------->>
\fontsize{9pt}{11pt}\selectfont
\textbf{\small EXAMPLE} \  
Take $X = Y = \bP^\infty(\C)$ $-$then $\dim X = \infty$ and the natural map $[\beta X,X] \ra [X,X]$ is not surjective 
(consider $\id_X$).\\
\endgroup %%------------------------------------<<

\index{Theorem: Dowker Extension Theorem}
\index{Dowker Extension Theorem}
\begingroup%%----------------------------------->>
\fontsize{9pt}{11pt}\selectfont
\textbf{\small DOWKER EXTENSION THEOREM} \  
Let \mX be normal with $\dim X \leq n + 1$ $(n \geq 1)$ and let \mA be a closed subspace of \mX.  Suppose that 
$f \in C(A,\bS^n)$ $-$then $\exists$ $F \in C(X,\bS^n)$: $\restr{F}{A} = f$ iff 
$f^*(\cH{}^n(\bS^n)) \subset i^*(\cH{}^n(X))$, $i:A \ra X$ the inclusion.
\vspi
[The argument splits into two parts.
\\
%%----------------------------------------------------------------------------------------------46
\indent\indent $(n = 1)$\ In this case $[X,A;\bS^1,s_1] \approx \cH{}^1(X,A)$, 
so one can proceed directly (\mA has the HEP w.r.t. $\bS^1$ (cf. p. \pageref{19.40}).)
\\
\indent\indent $(n > 1)$\ To reduce to the compact situation, 
use the fact that the extendability of $f:A \ra \bS^n$ to \mX is equivalent to the extendability of 
$\beta f:\beta A \ra \bS^n$ to $\beta X$ and consider the commutative diagram  
\begin{tikzcd}%[ sep=small]
{\cH{}^n(\beta X)} \ar{d} \ar{r} &{\cH{}^n(\beta A)} \ar{d} 
&{\cH{}^n(\bS^n)} \ar{l} \arrow[d,shift right=0.5,dash] \arrow[d,shift right=-0.5,dash]\\
{\cH{}^n(X)} \ar{r} &{\cH{}^n(A)}  &{\cH{}^n(\bS^n)} \ar{l}
\end{tikzcd}
.]\\
\vspace{0.5cm}
\endgroup %%------------------------------------<<

\index{Theorem: Dowker Classification Theorem}
\index{Dowker Classification Theorem}
\begingroup%%----------------------------------->>
\fontsize{9pt}{11pt}\selectfont
\textbf{\small DOWKER CLASSIFICATION THEOREM} \  
Let \mX be normal with $\dim X \leq n$ $(n \geq 1)$ and let \mA be a closed subspace of \mX.  Fix a generator 
$\iota \in \cH{}^n(\bS^n,s_n;\Z)$ $-$then the assignment $[f] \ra f^*\iota$ defines a bijection 
$[X,A;\bS^n,s_n] \ra \cH{}^n(X,A;\Z)$.
\vspi
[Show that $\forall \ n > 1$, $[\beta X,\beta A;\bS^n,s_n] \approx [X,A;\bS^n,s_n]$.]\\
\endgroup %%------------------------------------<<

\begin{proposition} \ %15
Suppose that $X = A \cup B$, where \mA and \mB are closed.  Let 
$
\begin{cases}
\ f \in C(A,\bS^n)\\
\ g \in C(B,\bS^n)
\end{cases}
$ 
and put $D = \{x \in A \cap B: f(x) \neq g(x)\}$.  Assume: $\dim D \leq n - 1$ $-$then $\exists$
$
\begin{cases}
\ F \in C(X,\bS^n): \restr{F}{A} = f\\
\ G \in C(X,\bS^n): \restr{G}{A} = g
\end{cases}
\& \ F \simeq G.
$ 
\end{proposition}

[Using Proposition 14, fix a homotopy $h:I(A \ \cap \  B) \ra \bS^n$ such that \ 
$
\begin{cases}
\ h(x,0) = f(x) \\
\ h(x,1) = g(x)
\end{cases}
$
$(x \in A \cap B)$.  \ 
Since $A \cap B$ as a subspace of 
$
\begin{cases}
\ A \\
\ B 
\end{cases}
$
has the HEP w.r.t. $\bS^n$, there exist continuous functions 
$
\begin{cases}
\ \phi:IA \ra \bS^n\\
\ \psi:IB \ra \bS^n
\end{cases}
$
with 
$
\begin{cases}
\ \phi(x,0) = f(x) \ (x \in A) \\
\ \psi(x,1) = g(x) \ (x \in B)
\end{cases}
$
and $\restr{\phi}{I(A \cap B)} = h = \restr{\psi}{I(A \cap B)}$.  
Define $H \in C(IX,\bS^n)$ by 
$
\begin{cases}
\ \restr{H}{IA} = \phi \\
\ \restr{H}{IB} = \psi
\end{cases}
$ 
and consider 
$
\begin{cases}
\ F(x) = H(x,0) \\
\ G(x) = H(x,1) 
\end{cases}
$
$(x \in X)$.]
\\

\begingroup%%----------------------------------->>
\fontsize{9pt}{11pt}\selectfont
\textbf{\small FACT} \  
Let \mA be a closed subset of \mX and let $f \in C(A,\bS^n)$.  
Assume $X = \ds\bigcup\limits_1^\infty O_j$, where the $O_j$ are open, 
$\dim\fr O_j \leq n - 1$, and $\forall \ j$,  $f$ has a continuous extension to $A \cup \ov{O}_j$ $-$then 
$\exists$ $F \in C(X,\bS^n)$ : $\restr{F}{A} = f$.\\
\endgroup %%------------------------------------<<

Suppose that $\dim X = n$ is positive.  Let $f:X \ra [0,1]^n$ be universal $-$then the restriction 
$f^{-1}(\bS^{n-1}) \ra \bS^{n-1}$ has no continuous extension to \mX, thus is essential.  Put 
$X_f = X/f^{-1}(\bS^{n-1})$, identify $\bS^n$ with $[0,1]^n/\bS^{n-1}$ and let $F_f:X_f \ra \bS^n$ be the induced map.\\

\textbf{\small LEMMA} \  
$F_f$ is essential, hence $\dim X_f = n$.

[Put $A = f^{-1}(\bS^{n-1})$ $-$then there is a \cd
\[
\begin{tikzcd}%[ sep=small]
{(X,A)} \ar{d}[swap]{f} \ar{r} &{(X/A,*_A)} \ar{d}{F_f}\\
{([0,1]^n,\bS^{n-1})} \ar{r} &{(\bS^n,s_n)}
\end{tikzcd}
.
\]
%\vspace{0.05cm}
%%----------------------------------------------------------------------------------------------47
\indent\indent $(n = 1)$ \ \  To get a contradiction, assume that $F_f$ is inessential.  
Choose $\phi \in C(X_f) :$  
$F_f(x) = \exp(2\pi i \phi(x))$.  \ Since $F_f(x) = 1$ only if $x = *_A$, \  $\phi(x) \in \Z$ only if \ $x = *_A$.  \ 
Normalize and take $\phi(*_A) = 0$.  Let 
$S = f^{-1}(0) \cup p^{-1}(\phi^{-1}(]0,1[))$.  Noting that $f(x) = \phi(p(x))$ mod 1, write 
$S = f^{-1}([0,1/2])$ $\cap$ $p^{-1}(\phi^{-1}([0,1/2]))$ $\cup$ $p^{-1}(\phi^{-1}([1/2,1]))$ 
to see that \mS is closed and write
$S = f^{-1}([0,1/2[)$ $\cap$ $p^{-1}(\phi^{-1}(]-1/4,1/2[))$ $\cup$ $p^{-1}(\phi^{-1}(]1/4,1[))$ 
to see that \mS is open.  
The characteristic function of the complement of \mS is thus a continuous extension to \mX of the restriction 
$f^{-1}(\{0,1\}) \ra \{0,1\}$.

\indent\indent $(n > 1)$ \quad The \cd
\[
\begin{tikzcd}%[ sep=small]
{\cH{}^n(\bS^n,s_n) } \ar{r}
&{\cH{}^n([0,1]^n,\bS^{n-1})} \ar{r}{f^*} 
&{\cH{}^n(X,A)}\\
&{\cH{}^{n-1}(\bS^{n-1})} \ar{u}  \ar{r}
&{\cH{}^{n-1}(A)} \ar{u}
&{\cH{}^{n-1}(X)} \ar{l}
\end{tikzcd}
\]
displays the data (cf. p. \pageref{19.41}).  
In view of the Dowker extension theorem, $f^*$ is not the zero homomorphism.  
Since the arrow $\cH{}^n(\bS^n,s_n) \ra \cH{}^n([0,1]^n,\bS^{n-1})$ is an isomorphism, 
it follows that $F_f$ is essential.]\\

Suppose that $IX$ is normal $-$then by the product theorem, $\dim IX \leq \dim X + 1$.  One can also go the other way: $\dim IX \geq \dim X + 1$.  This is obvious if $\dim X = 0$, so assume that $\dim X = n$ is positive.  
Claim: $\dim IX_f \geq n + 1$.  Indeed, if $\dim IX_f \leq n$, then Alexandroff's criterion would imply that the continuous function $\phi:i_0X_f \cup i_1X_f \ra \bS^n$ defined by 
$
\begin{cases}
\ \phi(x,0) = F_f(x)\\
\ \phi(x,1) = s_n
\end{cases}
(x \in X_f)
$
has a continuous extension to $IX_f$, meaning that $F_f$ is homotopic to a constant map and this contradicts the lemma.  
Now write $X - f^{-1}(\bS^{n-1}) = \bigcup\limits_1^\infty A_j$, where the $A_j$ are closed subspaces of \mX.  Let 
$*_f$ be the image of $f^{-1}(\bS^{n-1})$ in $X_f$ $-$then 
$X_f = \{*_f\} \cup \bigcup\limits_1^\infty A_j$ $\implies$ 
$IX_f = I\{*_f\} \cup \bigcup\limits_1^\infty IA_j$ $\implies$ $\exists \ j$: 
$\dim IX_f = \dim IA_j$ $\implies$ $\dim IX \geq \dim IA_j =$ $\dim IX_f \geq n + 1 = $ $\dim X + 1$.\\

Application: Suppose that $X \times [0,1]^m$ is normal $-$then $\dim(X \times [0,1]^m) = \dim X + m$.\\

\begin{proposition} \ %16
Suppose that \mX is normal and \mY is a CW complex.  Assume: $X \times Y$ is normal $-$then 
$\dim (X \times Y) = \dim X + \dim Y$.
\end{proposition}

[If \mB is a compact subspace of \mY which is homeomorphic to $[0,1]^m$, where $m = \dim Y$, then 
$\dim (X \times B) = \dim X + m$.]

[Note: \ The same conclusion obtains if \mY is a metrizable topological manifold.]\\

\begingroup%%----------------------------------->>
\fontsize{9pt}{11pt}\selectfont
\textbf{\small EXAMPLE} \  
Let \mX and \mY be nonempty CW complexes $-$then $X \times_k Y$ is a CW complex and 
$\dim(X \times_k Y) = \dim(X \times Y)$.\\
\endgroup %%------------------------------------<<

%%----------------------------------------------------------------------------------------------48
\begin{proposition} \ %17
Suppose that \mX is normal with $\dim X = 1$ and \mY is paracompact and $\sigma$-locally compact.   Assume: $X \times Y$ is normal $-$then $\dim (X \times Y) = \dim X + \dim Y$.
\end{proposition}

[Switch the roles of \mX and \mY and reduce to the case when \mX is compact.  Since $\dim Y = 1$, there exist disjoint closed sets 
$
\begin{cases}
\ B^\prime \subset Y\\
\ B\pp \subset Y
\end{cases}
$
such that $\ov{V} - V \neq \emptyset$ for any open $V \subset Y$: $B^\prime \subset V \subset Y - B\pp$.  
Arguing as above, it need only be shown that $\dim (X_f \times Y) \geq n + 1$ $(n > 0)$.  If instead, 
$\dim (X_f \times Y) \leq n$, define a continuous function $\phi:X_f \times (B^\prime \cup B\pp) \ra \bS^n$ by 
$
\begin{cases}
\ \phi(x,y) = F_f(x) \hspace{0.25cm} ((x,y) \in X_f \times B^\prime)\\
\ \phi(x,y) = s_n  \hspace{0.85cm} ((x,y) \in X_f \times B\pp)
\end{cases}
$
and use Alexandroff's criterion to get a continuous extension 
$\Phi:X_f \times Y \ra \bS^n$.  Let $V \subset Y$ be the set of all $y$ with the property that the section 
$
\Phi_y : 
\begin{cases}
\ X_f \ra \bS^n\\
\ x \ra \Phi(x,y)
\end{cases}
$
is essential $-$then $B^\prime \subset V \subset Y - B\pp$ and \mV is clopen, $X_f$ being compact.  Contradiction.]\\

\begingroup%%----------------------------------->>
\fontsize{9pt}{11pt}\selectfont
\textbf{\small EXAMPLE} \  
Take, after Anderson-Keisler (cf. p. \pageref{19.42}), an $X \subset \R^2$: $\dim X = \dim (X \times X) = 1$ $-$then 
$\dim \beta(X \times X) = 1$ but $\dim (\beta X \times \beta X) = \dim \beta X  + \dim \beta X  = 2$ (cf. Proposition 17).\\
\endgroup %%------------------------------------<<

\begingroup%%----------------------------------->>
\fontsize{9pt}{11pt}\selectfont
While there is no reason to suppose that $X_f$ is completely regular if \mX is, nevertheless the lemma and Propositions 16 and 17 are still true in this setting, although some changes in the proofs are necessary 
(Morita\footnote[2]{\textit{Fund. Math.} \textbf{87} (1975), 31-52.}).  
Consider, e.g., Proposition 17.  Having made the reduction and the switch (so \mX is compact and $\dim Y = 1$) choose 
a continuous function $h:Y \ra [0,1]$ such that $\ov{V} - V \neq \emptyset$ for any open $V \subset Y$:
$h^{-1}(0) \subset V \subset Y - h^{-1}(1)$.  Define $H:X_f \times Y \ra [0,1]^{n+1}$ by 
$H(x,y) = (1 - h(y)) F_f(x) + h(y)s_n$.  If $\dim (X \times Y) \leq n$ (where $n \geq 1$), then $\dim (X_f \times Y) \leq n$, 
therefore \mH is not universal.  
Accordingly (cf. p. \pageref{19.43}), $\exists$ $\Phi \in C(X_f \times Y,\bS^n)$: 
$
\begin{cases}
\ \Phi(x,y) = F_f(x) \hspace{0.35cm} (y \in h^{-1}(0))\\
\ \Phi(x,y) = s_n \hspace{0.85cm} (y \in h^{-1}(1))
\end{cases}
$ 
and this suffices.\\
\endgroup %%------------------------------------<<

\begingroup%%----------------------------------->>
\fontsize{9pt}{11pt}\selectfont
\textbf{\small EXAMPLE} \  
Let \mX be an arbitrary nonempty topological space $-$then 
$\dim IX =$  
$\dim \crg IX =$ 
$\dim I\crg X =$ 
$\dim \crg X + 1 = \dim X + 1$.  
This fact can be used to compute $\dim \Gamma X$ and $\dim \Sigma X$, both of which have the value $\dim X + 1$.  
Observe first that the two lemmas on p. \pageref{19.44} hold ``in general''.  Therefore 
$\dim X + 1 = $
$\dim IX = $ 
$\max\{\dim i_1 X,\dim IX/i_1X\} = $ 
$\max \{\dim X, \dim \Gamma X\} = $ 
$\dim \Gamma X$.  And then 
$\dim \Gamma X = \max\{\dim X,\dim\Gamma X/X\} = $
$\max\{\dim X, \Sigma X\} = $
$\dim \Sigma X$.  
Corollary: If $f:X \ra Y$ is a continuous function and if $M_f$ is its mapping cylinder, then 
$\dim M_f = \max\{1 + \dim X,\dim Y\}$.
\vspi
[Note: \ Recall that a cofibered subspace has the EP w.r.t. $\R$, hence w.r.t. $[0,1]$ 
(cf. p. \pageref{19.45}).]\\
\endgroup %%------------------------------------<<

\textbf{\small LEMMA} \  
Let \mX be normal.  Suppose that there exists a sequence 
$\sU_1, \sU_2, \ldots$ of open coverings of \mX such that 
$\sU_{i+1}$ is a refinement of \hspace{0.0cm} $\sU_i$, the collection 
$\{\st(U,\sU_i): U \in \sU_i \ (i = 1, 2, \ldots)\}$ is a basis for \mX, and 
$\forall \ i$: $\ord(\sU_i) \leq n+1$ $-$then $\dim X \leq n$.

%%----------------------------------------------------------------------------------------------49
[Let $\sU = \{U_1, \ldots, U_k\}$ be a finite open covering of \mX.  
Denote by $X_i$ the union of all $U \in \sU_i$ : 
$\st(U,\sU_i)$ is contained in some element of  \hspace{0.0cm} $\sU$.  
Each $X_i$ is open; moreover, $X = \bigcup\limits_i X_i$.  Fix a map 
$f_i^{i+1}:\sU_{i+1} \ra \sU_i$ such that $\forall \ U \in \sU_{i+1}$: $f_i^{i+1} (U) \supset U$.  
Set 
$f_i^i = \id_{\sU_i}$ and for $i < j$, put $f_i^j = f_i^{i+1} \circx \cdots \circx f_{j-1}^j$.  
Introduce
\[
\sU(j) = \{U \in \sU_j: U \cap X_j \neq \emptyset\} \text{ and } 
\sV(j) = \{U \in \sU(j): U \cap \bigl(\bigcup\limits_{i < j} X_i\bigr) = \emptyset\}.
\]
Obviously, $\sV(j) \subset \sU(j) \subset \sU_j$ and $j^\prime \neq \j\pp$ $\implies$ 
$\sV(j^\prime) \cap \sV(j\pp) = \emptyset$.  
Given $U \in \sU(j)$, let $i(U)$ be the smallest integer $i \leq j$: 
$f_i^j(U) \cap X_i \neq \emptyset$, so $f_{i(U)}^j(U) \in \sV(i(U))$.
 Corresponding to any $V \in \sV(i)$ is the open set
\[
V^* = \bigcup\limits_{j \geq i} \bigcup \{U \cap X_j: U \in \sU(j), f_i^j(U) = V \ \& \ i(U) = i\}.
\]
Note that $V^* \subset V$ and $\forall \ U \in \sU(j)$, $U \cap X_j \subset f_{i(U)}^j(U)^*$.  
In addition, $\exists$ 
$U \in \sU_i$: $U \cap V \neq \emptyset$ and $\exists$ $k(V) \leq k$: $V \subset \st(U,\sU_i) \subset U_{k(V)}$,  
hence $V^* \subset U_{k(V)}$.  
The collection $\sV^* = \{V^*: V \in \bigcup\limits_i \sV(i)\}$ is therefore an open refinement of $\sU$.  
The claim then is that 
$\ord(\sV^*) \leq n + 1$.  
To this end, consider a generic nonempty intersection $V_1^* \cap \cdots \cap V_p^*$, 
where $V_1 \in \sV(i_1), \ldots, V_p \in \sV(i_p)$ are distinct elements of $\bigcup\limits_i \sV(i)$.  
Take an $x$ in 
$V_1^* \cap \cdots \cap V_p^*$ and choose $j: x \in X_j - \bigcup\limits_{i < j} X_i$ 
($\implies$ $i_1 \leq j, \ldots, i_p \leq j$).  
From the definitions, there exist \\
$U_1 \in \sU(j_1)$: 
$
\begin{cases}
\ f_{i_1}^{j_1}(U_1) = V_1 \\
\ i(U_1) = i_1
\end{cases}
$
$\& \ x \in U_1 \cap X_{j_1}, \ldots, U_p \ \in \sU(j_p)$: %dmc is this really just U_p
$
\begin{cases}
\ f_{i_p}^{j_p}(U_p) = V_p \\
\ i(U_p) = i_p
\end{cases}
$
$\&$ $x \in U_p \cap X_{j_p}$.  
But $x \in f_j^{j_1} (U_1) \cap \cdots \cap f_j^{j_p}(U_p)$ and since 
$ f_j^{j_1} (U_1) , \ldots, f_j^{j_p}(U_p)$ are all different, the claim is thus seen to follow from the fact that 
$\ord(\sU_j) \leq n + 1$.]\\

Application: Let \mX be normal.  Suppose that \mX admits a development $\{\sU_i\}$ such that $\{\sU_i\}$ is a star sequence and $\forall \ i$: $\ord(\sU_i) \leq n+1$ $-$then $\dim X \leq n$.\\

\index{Pasynkov Factorization Lemma}
\textbf{\small PASYNKOV FACTORIZATION LEMMA} \  \ 
Suppose that \mX is normal and \mY is metrizable $-$then for every $f \in C(X,Y)$ there exists a metrizable space \mZ with 
$
\begin{cases}
\ \dim Z \leq \dim X\\
\ \wt Z \leq \wt Y
\end{cases}
$
and functions 
$
\begin{cases}
\ g \in C(X,Z)\\
\ h \in C(Z,Y)
\end{cases}
$ 
such that $f = h \circx g$ with $h$ uniformly continuous and $g(X) = Z$.

[Assume that $\dim X = n$ is finite and $\wt Y \geq \omega$.  
Fix a sequence $\{\sV_i\}$ of neighborhood finite open coverings of \mY such that $\forall \ i$: $\#(\sV_i) \leq \wt Y$, arranging matters so that the diameter of each $V \in \sV_i$ is 
$< 1/i$.  Inductively construct a star sequence $\{\sU_i\}$ of neighborhood
%%----------------------------------------------------------------------------------------------50
finite open coverings of \mX such that $\forall \ i:$
$
\begin{cases}
\ \ord(\sU_i) \leq n + 1\\
\ \#(\sU_i) \leq \wt Y
\end{cases}
$
and $\sU_i$ is a star refinement of $f^{-1}(\sV_i)$.  
Justification: Quote Proposition 6 and recall $\S 1$, Proposition 13 (the proof of which allows one to say that the cardinality of 
$\sU_i$ remains $\leq \wt Y$).  
Let $\delta$ be a continuous pseudometric on \mX associated with $\{\sU_i\}$ as on p. 
\pageref{19.46}.  
The claim is that one can take for \mZ the metric space $X_\delta$ obtained from \mX by identifying points at a zero distance from one another.  
Granted this, it is clear what $g$ and $h$ have to be.  
Denote by $X(\delta)$ the set \mX equipped with the topology determined by $\delta$.  
Given $U \in \sU_i$, write $U(\delta)$ for its interior in $X(\delta)$ and put 
$\sU_i(\delta) = \{U(\delta):U \in \sU_i\}$ 
$-$then $\{\sU_i(\delta)\}$ is a development for $X(\delta)$ and is a star sequence such that 
$\forall \ i$: $\ord(\sU_i(\delta)) \leq n + 1$.  
The projection $p:X(\delta) \ra Z$ is an open map (every open subset of 
$X(\delta)$ is $p$-saturated), thus $\sW_i \equiv p(\sU_i(\delta))$ is an open covering of \mZ.  
Furthermore, $\{\sW_i\}$ is a development for \mZ and is a star sequence such that $\forall \ i :$ $\ord(\sW_i) \leq n + 1$.  Therefore $\dim Z \leq n$.  
As for the assertion $\wt Z \leq \wt Y$, note that the $\sW_i$ are point finite and the collection 
$\bigcup\limits_1^\infty \{\st(z,\sW_i):z \in Z\}$ is a basis for \mZ.]\\

There are two related results, applicable to pairs $(X,A)$.

\indent\indent (A) \ Suppose that \mX is normal and \mY is metrizable of weight $\leq \kappa$.  Let \mA be a subspace of \mX having the EP w.r.t. $\sB(\kappa)$ $-$then for every $f \in C(A,Y)$ there exists a metrizable space $Z_A$ of weight 
$\leq \kappa$ and functions
$
\begin{cases}
\ g \in C(X,Z_A)\\
\ h_A \in C(g(A),Y)
\end{cases}
$
such that $f = h_A \circx (\restr{g}{A})$ with $h_A$ uniformly continuous and $g(X) = Z_A$.

[Argue as in $\S 6$, Proposition 15 (proof of sufficiency).]

\indent\indent (X/A) \ Suppose that \mX is normal and \mY is metrizable of weight $\leq \kappa$.  Let \mA be a closed subspace of \mX: $\dim(X/A) \leq n$ $-$then for every $f \in C(X,Y)$ there exists a metrizable space \mZ of weight $\leq \kappa$ and functions 
$
\begin{cases}
\ g \in C(X,Z)\\
\ h \in C(Z,Y)
\end{cases}
$
such that $f = h \circx g$ with $h$ uniformly continuous and $g(X) = Z$, $\dim(Z - \ov{g(Z)} \leq n$.

[This is the relative version of the Pasynkov factorization lemma.  The proof is the same as for the absolute case modulo the following remark: Every neighborhood finite open covering $\sU = \{U_i:i \in I\}$ of \mX has a neighborhood finite open refinement $\sO$ such that the order of the collection $\{O,\st(A,\sO):O \in \sO \ \& \ O \cap A = \emptyset\}$ is $\leq n + 1$.  
Proof: Assuming that the $U_i$ are cozero sets, let $\sZ = \{Z_i: i \in I\}$ be a precise zero set refinement of $\sU$ 
(cf. p. \pageref{19.47}).  Define $I_0 = \{i \in I:U_i \cap A \neq \emptyset\}$ and put 
$
\begin{cases}
\ Z_0 = \bigcup \{Z_i:i \in I_0\}\\
\ U_0 = \bigcup \{U_i:i \in I_0\}\
\end{cases}
$
$-$then 
$
\begin{cases}
\ Z_0\\
\ U_0
\end{cases}
$ is a 
$
\begin{cases}
\ \text{zero set}\\
\ \text{cozero set}
\end{cases}
$ 
(cf. p. \pageref{19.48}).  
Choose $\phi \in C(X,[0,1])$: $Z_0 = \phi^{-1}(0)$ $\&$ $X - U_0 = \phi^{-1}(1)$.  
Let $X_0 = \{x: \phi(x) \leq 1/2\}$.  Since \mA is contained in $Z_0$ and $Z_0$ is contained in the interior of $X_0$, the collection $\{U_i - X_0, U_0: i \in I - I_0\}$ is the inverse
%%----------------------------------------------------------------------------------------------51
image of a neighborhood finite cozero set covering of $X/A$ under the projection $X \ra X/A$.  Therefore there exists a neighborhood finite cozero set covering $\{O_i, O_0: i \in I - I_0\}$ of \mX whose order does not exceed $n + 1$ such that 
$O_i \subset U_i - X_0$ $(i \in I - I_0)$ and $A \subset O_0 \subset U_0$.  If 
$\sO = \{O_i: i \in I - I_0\} \cup \{O_0 \cap U_i: i \in I_0\}$, then $O_0 = \st(A,\sO)$ and $\sO$ is a  neighborhood finite cozero set refinement of $\sU$ with the stated property.]\\

\begin{proposition} \ 
Suppose that \mX is normal and \mY is completely metrizable of weight $\leq \kappa$ and locally $n$-connected 
($n$-connected and locally $n$-connected).  Let \mA be a closed subspace of \mX having the EP w.r.t. $\sB(\kappa)$.  
Assume: $\dim X/A \leq n + 1$ $-$then \mA has the NEP (EP) w.r.t. \mY.
\end{proposition}

[Take an $f \in C(A,Y)$ and write $f = h_A \circx (\restr{g}{A})$.  Since $g \in C(X,Z_A)$ and since $\wt Z_A \leq \kappa$, $g$ can in turn be factored: $g = \psi \circx \phi$, where 
$
\begin{cases}
\ \phi \in C(X,Z)\\
\ \psi \in C(Z,Z_A)
\end{cases}
. \ 
$
Here, of course, $\dim (Z - \ov{\phi(A)}) \leq n + 1$.  On the other hand, 
$h_A \circx (\restr{\psi}{\phi(A)})$ is uniformly continuous, hence extends to a continuous function $H_A:\ov{\phi(A)} \ra Y$.  Now apply the results of Dugundji cited on 
p. \pageref{19.49}.]\\

\begin{proposition} \ 
Suppose that $IX$ is normal and \mY is completely metrizable of weight $\leq \kappa$ and locally $n$-connected.  
Let \mA be a closed subspace of \mX having the EP w.r.t. $\sB(\kappa)$.  
Assume: $\dim X/A \leq n$ $-$then \mA has the HEP w.r.t. \mY.
\end{proposition}

[Let $f:i_0X \cup IA \ra Y$ be continuous.  
Since $i_0X \cup IA$, as a subspace of $IX$, has the EP w.r.t. $\sB(\kappa)$ 
(cf. $\S 6$, Proposition 16) and since 
$\dim IX/i_0X$ $\cup$ $IA \leq$ 
$\dim IX/IA \leq$ 
$\dim I(X/A) \leq $ 
$\dim X/A + 1 \leq n + 1$, Proposition 18 implies that there exists a cozero set 
$O \subset IX :$ $O \supset i_0X$ $\cup$ $IA$ and a continuous function $g:O \ra Y$ extending $f$.  
Fix a cozero set $U \subset X$: $IA \subset IU \subset O$.  Choose 
$\phi \in C(X,[0,1])$: 
$
\begin{cases}
\ \restr{\phi}{A} = 1\\
\ \restr{\phi}{X - U} = 0
\end{cases}
. \ 
$
Define $F \in C(IX,Y)$ by $F(x,t) = g(x,\phi(x)t)$: \mF is a continuous extension of $f$.]\\

\begingroup%%----------------------------------->>
\fontsize{9pt}{11pt}\selectfont
The normality of \mX can be dispensed with in Pasynkov's factorization lemma: Everything goes through in the completely regular situation.
\vspi
[Note: \ Pasynkov's factorization lemma is then valid for an arbitrary topological space as may be seen by passing to its complete regularization.]
\vspi
As for Propositions 18 and 19, they too are true if \mX is a nonempty CRH space.  The assumption that \mA is closed was made only to ensure that the quotient $X/A$ is normal.  Therefore it can be dropped.  Likewise the assumption that $IX$ is normal was made only to use the product theorem.  This, however, is of no real consequence, as the product theorem holds for an arbitrary nonempty topological space (cf. p. \pageref{19.50}).
\vspi
%%----------------------------------------------------------------------------------------------52
For another application of these methods, suppose that \mY is completely metrizable of weight 
$\leq \kappa$ and is $n$-connected and locally $n$-connected.  Assume: $\dim X/A \leq n$.  Let 
$
\begin{cases}
\ f:X \ra Y\\
\ g:X \ra Y
\end{cases}
$
be continuous functions such that $\restr{f}{A} \simeq \restr{g}{A}$ $-$then $f \simeq g$.  Corollary: If \mX is $\kappa$-collectionwlse normal with $\dim X \leq n$, then $[X,Y] = *$.\\
\endgroup %%------------------------------------<<

\begingroup%%----------------------------------->>
\fontsize{9pt}{11pt}\selectfont
\textbf{\small FACT} \  
Suppose that \mX is a nonempty metrizable space.  Let \mA be a nonempty closed subspace of \mX: $\dim (X - A) = 0$ 
$-$then there exists a retraction $r:X \ra A$.\\
\endgroup %%------------------------------------<<

A compact connected ANR $Y$ is said to be a 
\un{test space}
\index{test space} 
for dimension $n$ $(n \geq 1)$ provided that the statement $\dim X \leq n$ is true iff every closed subset $A \subset X$ has the EP w.r.t. \mY.  
Example: $\bS^n$ is a test space for dimension $n$ (Alexandroff's criterion).

[Note: \ No compact connected AR \mY can be a test space for dimension n.]\\

\textbf{\small LEMMA} \  Let 
$
\begin{cases}
\ Y^\prime\\
\ Y\pp
\end{cases}
$
be compact connected ANRs of the same homotopy type $-$then $Y^\prime$ is a test space for dimension $n$ iff $Y\pp$ is a test space for dimension $n$.

[If \mX is normal and $A \subset X$ is closed, then \mA has the HEP w.r.t. 
$
\begin{cases}
\ Y^\prime\\
\ Y\pp
\end{cases}
$
(cf. p. \pageref{19.51}).]\\

A finite wedge $\bigvee \bS^n$ of $n$-spheres is a test space for dimension $n$.  
Indeed, $\bigvee \bS^n$ is a compact connected ANR of topological dimension $n$.  
Moreover, $\bigvee \bS^n$ is $(n-1$)-connected (since for $n > 1$, 
$\pi_q(\bigvee \bS^n) = \oplus \pi_q(\bS^n)$ $(q < 2n - 1)$), 
thus Proposition 18 implies that if $\dim X \leq n$ and if 
$A \subset X$ is closed, then \mA has the EP w.r.t. $\bigvee \bS^n$.  
Here it is necessary to recall that \mA has the EP w.r.t. 
$\sB(\omega)$ (cf. p. \pageref{19.52}).  
On the other hand, there is a retraction $r:\bigvee \bS^n \ra \bS^n$ so if 
$A \subset X$ is closed and has the EP w.r.t. $\bigvee \bS^n$ then \mA has the EP w.r.t. $\bS^n$, from which 
$\dim X \leq n$.\\

\index{Theorem: Test Space Theorem}
\index{Test Space Theorem}
\textbf{\small TEST SPACE THEOREM} \  
Let \mY be a compact connected ANR of topological dimension $n$ $(n \geq 1)$ $-$then \mY is a test space for dimension $n$ iff \mY has the homotopy type of a finite wedge of $n$-spheres.

[Only the necessity need be dealt with.  There are two cases: $n = 1$ or $n > 1$.  
If $n = 1$, then $\pi_1(Y) \neq 1$ (otherwise, \mY would be an AR), hence \mY has the homotopy type of a finite wedge of 1-spheres (cf. p. \pageref{19.53}). 
If $n > 1$, then for $q > n$, $H_q(Y) = 0$ (cf. p. \pageref{19.54}) and \mY must be $(n-1)$-connected 
(cf. p. \pageref{19.55} $\&$ p. \pageref{19.56}).  
Accordingly, by Hurewicz, $H_q(Y) = 0$ $(0 < q < n)$ and 
$H_n(Y) = \pi_n(Y)$, a nontrivial finitely generated free abelian group.  
Picking a set of base point preserving maps 
$\bS^n \ra Y$ which generate $\pi_n(Y)$ then leads to a homology equivalence $\bigvee \bS^n \ra Y$ that, by the Whitehead theorem, is a homotopy equivalence.]\\

%%----------------------------------------------------------------------------------------------53
\begingroup%%----------------------------------->>
\fontsize{9pt}{11pt}\selectfont
If \mY is a compact connected ANR which is a test space for dimension $n$, then $\dim Y \geq n$ (look at the proof of the test space theorem).  There are test spaces for dimension $n$ of topological dimension $n + k$ $(k > 0)$.  
Consider e.g., $[0,1]^{n+k} \vee \bS^n$.\\
\endgroup %%------------------------------------<<

\begingroup%%----------------------------------->>
\fontsize{9pt}{11pt}\selectfont
\textbf{\small EXAMPLE} \  
Let $\alpha \in \pi_{n + k}(\bS^n)$ $(k > 0, n \geq 1)$.  Choose a representative $f \in \alpha$ and put 
$Y_\alpha = D^{n + k + 1} \sqcup_f \bS^n$ $-$then $Y_\alpha$ is a compact connected ANR 
(cf. p. \pageref{19.57})
with 
$\dim Y_\alpha = n + k + 1$ 
(cf. p. \pageref{19.571})
and 
Dranishnikov\footnote[2]{\textit{Tsukuba J. Math.} \textbf{14} (1990), 247-262.}
has shown that $Y_\alpha$ is a test space for dimension $n$.
\vspi
[Note: \ The preceding considerations break down if $k = 0$.  
Example: $\bP^2(\R)$ is not a test space for dimension 1.]\\
\endgroup %%------------------------------------<<

%%%%%%%%%%%%%%%%%%%%%%%%%%%%%%%%%%%%%%
%%%%%%%%%%%%%%%%%%%%%%%%%%%%%%%%%%%%%%
%%%%%%%%%%%%%%%%%%%%%%%%%%%%%%%%%%%%%%

\begin{center}
$\S \ 19$
\\[0.5cm]
$\mathcal{REFERENCES}$\\
\end{center}

\[
\textbf{BOOKS}
\]

\begingroup
\fontsize{9pt}{11pt}\selectfont
\setlength\parindent{0 cm}

[1] \quad Aarts, J. and Nishiura, T., \textit{Dimension and Extensions}, North Holland (1993).
\\[-.2cm]

[2] \quad Engelking, R., \textit{Theory of Dimensions: Finite and Infinite}, Heldermann Verlag (1995).
\\[-.2cm]

[3] \quad Hurewicz, W. and Wallman, H., \textit{Dimension Theory}, Princeton University Press (1948).
\\[-.2cm]

[4] \quad Isbell, J., \textit{Uniform Spaces}, American Mathematical Society (1964).
\\[-.2cm]

[5] \quad van Mill, J., \textit{Infinite Dimensional Topology}, North Holland (1989).
\\[-.2cm]

[6] \quad Nagami, K., \textit{Dimension Theory}, Academic Press (1970).
\\[-.2cm]

[7] \quad Nagata, J., \textit{Modern Dimension Theory}, Heldermann Verlag (1983).
\\[-.2cm]

[8] \quad Pears, A., \textit{Dimension Theory of General Spaces}, Cambridge University Press (1975).
\\[-.2cm]
\endgroup

\vspace{-.25cm}
\[
\textbf{ARTICLES}
\]

\begingroup
\fontsize{9pt}{11pt}\selectfont
\setlength\parindent{0 cm}

 [1] \quad Alexandroff, P., The Present Status of the Theory of Dimension, \textit{Amer. Math. Soc. Transl. Ser.} \textbf{1} 
 
 \hspace{0.8cm}(1955), 1-26.
 \\[-.2cm]

[2] \quad Chigogidze, A., Kawamura, K., and Tymchatyn E., Menger Manifolds, In: \textit{Continua}, H. Cook et al. 

\hspace{0.8cm}(ed.), Marcel Dekker (1995), 37-88.
\\[-.2cm]

[3] \quad Dowker, C., Mapping Theorems for Noncompact Spaces, \textit{Amer. J. Math.} \textbf{69} (1947), 200-242.
\\[-.2cm]

[4] \quad Fedorchuk, V., The Fundamentals of Dimension Theory, In: \textit{General Topology}, EMS \textbf{17}, Springer 

\hspace{0.8cm}Verlag (1990), 91-192.
\\[-.2cm]

[5] \quad James, I., On Category in the Sense of Lusternik-Schnirelmann, \textit{Topology} \textbf{17} (1978), 331-348.
\\[-.2cm]

[6] \quad Lelek, A., Dimension Inequalities for Unions and Mappings of Separable Metric Spaces, \textit{Colloq. Math.} 

\hspace{0.8cm}\textbf{23} (1971), 69-91.
\\[-.2cm]

[7] \quad Lon\u car, I., Some Results on the Dimension of Inverse Limit Spaces, \textit{Glasnik Mat.} \textbf{21} (1986), 213-223.
\\[-.2cm]

[8] \quad Luukkainen, J., Embeddings of $n$-Dimensional Locally Compact Metric Spaces to $2n$-Manifolds, \textit{Math.} 

\hspace{0.8cm}\textit{Scand.} \textbf{68} (1991), 193-209.
\\[-.2cm]

[9] \quad Montejano, L., Categorical and Contractible Covers of Polyhedra: Some Topological Invariants Re-

\hspace{0.8cm}lated to the Lusternik-Schnirelmann Category, \textit{Rend. Sem. Fac. Sci. Univ. Cagliari} \textbf{58} (1988), 

\hspace{0.8cm}suppl., 177-264.
\\[-.2cm]

[10] \quad Morita, K., Dimension of General Topological Spaces, In: \textit{Surveys in General Topology}, G. Reed 

\hspace{0.95cm}(ed.), Academic Press (1980), 297-336.
\\[-.2cm]

[11] \quad Nagata, J., A Survey of Dimension Theory, \textit{Proc. Steklov Inst. Math.} \textbf{154} (1984), 201-213 
and \textit{Q.}

\hspace{0.95cm}\textit{A. General Topology} \textbf{8} (1990), 61-77.
\\[-.2cm]

[12] \quad Pasynkov, B., Fedorchuk, V., and Filippov, V., Dimension Theory, \textit{J. Soviet Math.} \textbf{18} (1982), 789-

\hspace{0.95cm}841.
\\[-.2cm]

[13] \quad Pol, R., Questions in Dimension Theory, In: \textit{Open Problems in Topology}, J. van Mill and G. Reed 

\hspace{0.95cm}(ed.), North Holland (1990), 279-291.\\[-.2cm]

[14] \quad Sternfeld, Y., Hilbert's 13th Problem and Dimension, \textit{SLN} \textbf{1376} (1989), 1-49.
\\[-.2cm]

[15] \quad Vitushkin, A., On Representation of Functions by Means of Superpositions and Related Topics, 

\hspace{0.95cm}\textit{Enseign. Math.} \textbf{23} (1977), 255-320.

\setlength\parindent{2em}

\endgroup

\chapter{
$\boldsymbol{\S}$\textbf{20}.\quadx  COHOMOLOGICAL DIMENSION THEORY}
\setlength\parindent{2em}
\setcounter{proposition}{0}

%%----------------------------------------------------------------------------------------------01
$\text{ }$\\[-1.25cm]

Cohomological dimension theory enables one to associate with each nonempty normal Hausdorff space \mX and every nonzero abelian group \mG a topological invariant $\dim_G X \in \{0, 1, \ldots\} \cup \{\infty\}$ called its cohomological dimension with respect to \mG.  
It turns out that $\dim_\Z X = \dim X$ if $\dim X < \infty$ 
(but this can fail if $\dim X = \infty$) and when \mX is  CW complex, $\dim_G X = \dim X$ $\forall \ G \neq 0$.

Let \mG be an abelian group $-$then for any topological pair 
$(X,A)$, $\cH{}^n(X,A;G)$ is the $n^\text{th}$ 
\u Cech cohomology group of $(X,A)$ with coefficients in \mG calculated per numerable open coverings (rather than arbitrary open coverings).

[Note: \ As was shown by Morita\footnote[2]{\textit{Fund. Math.} \textbf{87} (1975), 31-52.}, 
\label{5.0ah}
$[X,A;K(G,n),k_{G,n}] \approx \cH{}^n(X,A;G)$ 
(cf. p. \pageref{20.1}) which, however, need not be true if the usual definition of 
``$\cH{}^n$'' is employed 
(Bredon\footnote[3]{\textit{Proc. Amer. Math. Soc.} \textbf{19} (1968), 396-398.}).  
Bear in mind that when $n = 0$, the agreement is that $K(G,0) = G$ (discrete topology).]\\

\textbf{\small LEMMA} \ 
If \mA is a nonempty subspace of \mX, then $\forall \ n \geq 1$, $\cH{}^n(X,A;G) \approx$ $\cH{}^n(X/A;G)$.\\

\begingroup%%----------------------------------->>
\fontsize{9pt}{11pt}\selectfont
\quad Let \mA be a subspace of \mX $-$then \mA is said to be 
\un{numerably embedded}
\index{numerably embedded} 
in \mX if for every numerable open covering $\sO$ of \mA there exists a numerable open covering $\sU$ of \mX such that 
$\sU \cap A$ is a refinement of $\sO$ (cf. $\S 6$, Proposition 15).  
Example: If \mX is a collectionwise normal Hausdorff space, then every closed subspace \mA of \mX is numerably embedded in \mX 
(cf. p. \pageref{20.2}).\\
\endgroup %%------------------------------------<<

\label{19.41}
\textbf{\small LEMMA} \ 
\begingroup%%----------------------------------->>
\fontsize{9pt}{11pt}\selectfont
Suppose that \mA is numerably embedded in \mX $-$then $\forall \ G$, there is a long exact sequence 
$\cdots \ra \cH{}^{n-1}(A;G) \ra \cH{}^n(X,A;G) \ra \cH{}^n(X;G) \ra \cH{}^n(A;G) \ra \cdots$.\\
\endgroup %%------------------------------------<<

\begingroup%%----------------------------------->>
\fontsize{9pt}{11pt}\selectfont
Remark: If $G = \Z$ (or, more generally, is finitely generated), one can get away with less, viz. it suffices that \mA have the EP w.r.t. $\R$.
\vspi
[Note: \ Working with countable numerable open coverings, an appeal to Proposition 4 in $\S 6$ leads to the definition of the coboundary operator $\cH{}^{n-1}(A) \ra \cH{}^n(X,A)$.]
\vspi
Example: If \mX is a normal Hausdorff space and if $A \subset X$ is closed, then there is a long exact sequence 
$\cdots \ra \cH{}^{n-1}(A) \ra \cH{}^n(X,A) \ra \cH{}^n(X) \ra \cH{}^n(A) \ra \cdots$.\\
\endgroup %%------------------------------------<<

\begingroup%%----------------------------------->>
\fontsize{9pt}{11pt}\selectfont
\textbf{\small FACT} \  
Suppose that \mA is numerably embedded in \mX $-$then $IA$ is numerably embedded in $IX$.\\
\endgroup %%------------------------------------<<

%%----------------------------------------------------------------------------------------------02
\begingroup%%----------------------------------->>
\fontsize{9pt}{11pt}\selectfont
It is known that $\cH{}^*(-;G)$, restricted to the full subcategory of $\bTOP^2$ whose objects are the pairs $(X,A)$, where \mA is closed and numerably embedded in \mX, satisifies the seven axioms of Eilenberg-Steenrod for a cohomology theory (Watanabe\footnote[2]{\textit{Glas. Mat.} \textbf{22} (1987), 187-238; 
see also \textit{SLN} \textbf{1283} (1987), 221-239.}).\\
\endgroup %%------------------------------------<<

\begin{proposition} \ %01
Let \mX be a nonempty normal Hausdorff space.  Assume: $\dim X \leq n$ $-$then $\cH{}^q(X;G) = 0$ $(q > n)$.
\end{proposition}

[This is a consequence of the definitions (cf. $\S 19$, Proposition 6).]\\

\textbf{\small PROPOSITION 1} \textnormal{(\textbf{bis})} \quad
Let \mA be a nonempty closed subspace of \mX.  Assume: $\dim X/A \leq n$ $-$then $\cH{}^q(X,A;G) = 0$ $(q > n)$.\\

\begingroup%%----------------------------------->>
\fontsize{9pt}{11pt}\selectfont
If \mX is a locally contractible paracompact Hausdorff space (e.g., a CW complex or an ANR), then $\forall \ n$, 
$\cH{}^n(X;G) \approx H^n(X;G)$.  In general, though, \u Cech cohomology and singular cohomology can differ even if \mX is compact Hausdorff (Barratt-Milnor\footnote[3]{\textit{Proc. Amer. Math. Soc.} \textbf{13} (1962), 293-297.}).
\vspi
[Note: \ Proposition 1 is a key property of \u Cech cohomology.  It is not shared by singular cohomology.]\\
\endgroup %%------------------------------------<<

Fix an abelian group $G \neq 0$ and let \mX be a nonempty normal Hausdorff space.  
Consider the following statement.

\indent\indent $(\dim_G X \leq n)$ \quad There exists an integer $n = 0, 1, \ldots$ such that $\cH{}^q(X,A;G) = 0$ 
$(q > n)$ for all closed subsets \mA of \mX.

If $\dim_G X \leq n$ is true for some n, then the 
\un{cohomological dimension}
\index{cohomological dimension} 
of \mX with respect to \mG, 
denoted by $\dim_G X$, 
is the smallest value of $n$ for which 
$\dim_G X \leq n$.

[Note: \ By convention, $\dim_G X = -1$ when $X = \emptyset$ or when $G = 0$.  If the statement $\dim_G X \leq n$ is false for every n, then we put $\dim_G X = \infty$.]\\

\begingroup%%----------------------------------->>
\fontsize{9pt}{11pt}\selectfont
\textbf{\small EXAMPLE} \quad 
Let \mX be a metrizable compact Hausdorff space of finite topological dimension, \mK a simply connected CW complex $-$then 
$\dim_{H_q(K)} X \leq q$ $\forall \ q \geq 1$ 
iff 
$\dim_{\pi_q(K)} X \leq q$ $\forall \ q \geq 1$ 
and both are equivalent to every closed subset $A \subset X$ having the EP w.r.t. 
\mK (Dranishnikov\footnote[6]{\textit{Math. Sbornik} \textbf{74} (1993), 47-56; 
see also Dydak, \textit{Trans. Amer. Math. Soc.} \textbf{337} (1993), 219-234.}).  
Example: One can take $K = M(G,n)$ $(n \geq 2)$ (realized as a simply connected CW complex) provided that 
$\dim_G X \leq n$.\\
\endgroup %%------------------------------------<<

\begin{proposition} \ %02
Suppose tht $\dim X \leq n$ $-$then $\dim_G X \leq n$.
\end{proposition}

[In fact, for $A \neq \emptyset$, $\dim X \leq n$ $\implies$ $\dim X/A \leq n$ 
(cf. p. \pageref{20.3}) $\implies$ 
$\cH{}^q(X,A;G) = 0$ $(q > n)$ (cf. Proposition 1 (bis)) $\implies$ $\dim_G X \leq n$ (Proposition 1 covers the case when 
$A = \emptyset$.]\\

%%----------------------------------------------------------------------------------------------03
\begin{proposition} \ %03
Suppose that $\dim X < \infty$ $-$then $\dim_{\Z} X = \dim X$.
\end{proposition}

[In view of Proposition 2, $\dim_{\Z} X \leq \dim X$.  Now argue by contradiction and assume that 
$\dim_{\Z} X \leq n$, $\dim X = n+1$.  Choose a universal map $f:X \ra [0,1]^{n+1}$ 
(cf. p. \pageref{20.4}) $-$then 
on the basis of the Dowker extension theorem, the arrow \ 
$\cH{}^{n+1}([0,1]^{n+1},\bS^n;\Z)$ $\overset{f^*}{\lra} \cH{}^{n+1}(X,f^{-1}(\bS^n);\Z)$ \  is not the zero homomorphism.  
But 
$\dim_{\Z} X \leq n$ $\implies$ $\cH{}^{n+1}(X,f^{-1}(\bS^n);\Z) = 0$.]\\

Application: If the topological dimension of \mX is finite, then $\forall \ G$, $\dim_G X \leq \dim_{\Z} X$.

[Note: \ For any compact Hausdorff space \mX (possibly of infinite topological dimension), one has 
$\dim_G X \leq \dim_{\Z} X$ (immediate from the universal coefficient theorem (cf. infra)).]\\

\begingroup%%----------------------------------->>
\fontsize{9pt}{11pt}\selectfont
\textbf{\small EXAMPLE} \quad 
The validity of the relation $\dim_{\Z} X = \dim X$ depends on the assumption that $\dim X < \infty$.  Indeed, 
Dranishnikov\footnote[2]{\textit{Math. Sbornik} \textbf{63} (1989), 539-545; 
see also Chigogidze, \textit{Inverse Spectra}, North Holland (1996), 100-116.} 
has given an example of a compact metric space \mX such that $\dim X = \infty$, while $\dim_{\Z} X < \infty$.
\vspi
\label{20.18}
[Note: \ According to Watanabe\footnote[3]{\textit{Proc. Amer. Math. Soc.} \textbf{123} (1995), 2883-2885.} 
, $\dim_{\Z} X = \dim X$ if \mX is a compact ANR (no restriction on $\dim X$).]\\
\endgroup %%------------------------------------<<

There is not a great deal that can be said about $\dim_G X$ if \mX is merely normal, so we shall restrict ourselves in what follows to paracompact \mX and begin with a review of \u Cech cohomology in this situation (all open coverings thus being numerable).\\

\index{Mayer-Vietoris Sequence}
\textbf{\small MAYER-VIETORIS SEQUENCE} \quad 
Let \mX be a paracompact Hausdorff space.  Suppose that \mA, \mB are closed subsets of \mX with $X = A \ \cup \ B$ $-$then the sequence 
$\cdots \ra \cH{}^n(X;G)$ 
$\ra \cH{}^n(A;G) \oplus \cH{}^n(B;G)$ 
$\ra \cH{}^n(A \cap B;G)$ 
$\ra \cH{}^{n+1}(X;G) \ra \cdots$ is exact.\\

\index{Bockstein Sequence}
\textbf{\small BOCKSTEIN SEQUENCE} \quad 
Let \mX be a paracompact Hausdorff space, \mA a closed subset.  Suppose that 
$0 \ra G^\prime \ra G \ra G\pp \ra 0$ is a short exact sequence of abelian groups $-$then there is a long exact sequence 
$\cdots \ra \cH{}^n(X,A;G^\prime) \ra$ 
$\cH{}^n(X,A;G) \ra \cH{}^n(X,A;G\pp)$
$\ra \cH{}^{n+1}(X,A;G^\prime) \ra \cdots$.\\

\index{Theorem: Universal Coefficient Theorem}
\index{Universal Coefficient Theorem}
\textbf{\small UNIVERSAL COEFFICIENT THEOREM} \quad 
Let \mX be a compact Hausdorff space, \mA a closed subset $-$then there is a split short exact sequence 
$0 \ra \cH{}^n(X,A;\Z) \otimes G$ 
$\ra \cH{}^n(X,A;G)$ 
$\ra \Tor(\cH{}^{n+1}(X,A;\Z),G) \ra 0$.\\

%%----------------------------------------------------------------------------------------------04
\index{K\"unneth Formula}
\textbf{\small K\"UNNETH FORMULA} \quad 
Let \mX be a paracompact Hausdorff space, \mA a closed subset; Let \mY be a compact Hausdorff space, $B$ a closed subset 
$-$then
$\cH{}^n((X,A)\times (Y,B);G) \approx$ 
$\bigoplus\limits_{q=0}^n \ \cH{}^q(X,A;\cH{}^{n-q}(Y,B;G))$.

[Note: \ The product $X \times Y$ is a paracompact Hausdorff space and, as usual, 
$(X,A) \times (Y,B) =$ 
$(X \times Y,X \times B \cup A \times Y)$.]\\

\begingroup%%----------------------------------->>
\fontsize{9pt}{11pt}\selectfont
Let \mX be a paracompact Hausdorff space of finite topological dimension.  Suppose that \mG is finitely generated 
$-$then Bartik\footnote[2]{\textit{Quart. J. Math.} \textbf{29} (1978), 77-91}
has shown that for every closed subset \mA of \mX, the arrow 
$\cH{}^n(\beta X, \beta A;G) \ra \cH{}^n(X,A;G)$ is surjective for $n = 1$ and bijective for $n > 1$.
\vspi
[Note: \ More is true if \mG is finite:  The arrow 
$\cH{}^n(\beta X, \beta A;G) \ra \cH{}^n(X,A;G)$ is bijective $\forall \ n \geq 0$.]\\ 
\endgroup %%------------------------------------<<

\label{20.9}
\label{20.10}
\begingroup%%----------------------------------->>
\fontsize{9pt}{11pt}\selectfont
\textbf{\small EXAMPLE} \quad 
Let \mX be a paracompact Hausdorff space of finite topological dimension $-$then 
$\dim_G X \leq \dim_G \beta X$ provide that \mG is finitely generated.
\vspi
[This is clear if $\dim_G X \leq 0$, so let $n = \dim_G X$ be positive and choose a closed subset \mA of \mX such that 
$\cH{}^n(X,A;G) \neq 0$.  
By the above, 
$\cH{}^n(\beta X, \beta A;G) \neq 0$, thus $n \leq \dim_G \beta X$.]\\
\endgroup %%------------------------------------<<

Notation: Let \mX be a paracompact Hausdorff space, $A \subset X$ a closed subset.  
Given $e \in \cH{}^n(X;G)$, write $\restr{e}{A}$ for the image of e under the arrow $\cH{}^n(X;G) \ra \cH{}^n(A;G)$.\\

\index{Restriction Principle}
\textbf{\small RESTRICTION PRINCIPLE} \quad 
Let $e$ be an element of $\cH{}^n(X;G)$.  Assume $\restr{e}{A} = 0$ $-$then $\exists$ an open $U \supset A$: 
$\restr{e}{\ov{U}} = 0$.\\

\index{Extension Principle}
\textbf{\small EXTENSION PRINCIPLE} \quad 
Suppose that $\alpha \in \cH{}^n(A;G)$ $-$then $\exists$ an open $U \supset A$ and an $e \in \cH{}^n(\ov{U};G)$: 
$\restr{e}{A} = \alpha$.\\

These two principles date back to 
Wallace\footnote[3]{\textit{Duke Math J.} \textbf{19} (1952), 177-182.}
who used them to establish the following result.\\

\index{Theorem: Relative Homeomorphism Theorem}
\index{Relative Homeomorphism Theorem}
\textbf{\small RELATIVE HOMEOMORPHISM THEOREM} \quad 
Let
$
\begin{cases}
\ X\\
\ Y
\end{cases}
$
be paracompact Hausdorff spaces; let 
$
\begin{cases}
\ A \subset X\\
\ B \subset Y
\end{cases}
$
be closed subsets.  Suppose given a closed map $f:(X,A) \ra (Y,B)$ such that $\restr{f}{X - A}$ 
is a homeomorphism of $X - A$ onto 
$Y - B$ $-$then $f^*:\cH{}^n(Y,B;G) \ra \cH{}^n(X,A;G)$ is an isomorphism.\\

%%----------------------------------------------------------------------------------------------05
Application: Let \mX be a paracompact Hausdorff space; let 
$
\begin{cases}
\ A\\
\ B
\end{cases}
\subset X
$
be closed subsets $-$then the arrow 
$\cH{}^n(A \cup B,A) \ra \cH{}^n(B,A \cap B)$ induced by the inclusion $(B,A \cap B) \ra (A \cup B,A)$ is an isomorphism.\\

\label{20.15}
\label{20.15a}
\begingroup%%----------------------------------->>
\fontsize{9pt}{11pt}\selectfont
It is possible to expand the level of generality and incorporate sheaves (of abelian groups) into the theory.  While this is definitely of interest, I shall limit the discussion to a few elementary observations.
\vspi
Let \mX be a paracompact Hausdorff space.  Given a sheaf $\sF \neq 0$ on \mX, write $\dim_{\sF} X \leq n$ if $\exists$ 
an integer $n = 0, 1, \ldots$ such that $\cH{}^q(X;\restr{\sF}{U}) = 0$ $(q > n)$ for all open subsets \mU of \mX.  
Example: $\dim X \leq n$ $\implies$ $\dim_{\sF} X \leq n$ (cf. Proposition 2).
\vspi
Remark: Let $G \neq 0$ be an abelian group, \bG the constant sheaf on \mX determined by \mG $-$then $\forall$ closed subset $A \subset X$, $\cH{}^n(X,A;G) \approx \cH{}^n(X;\restr{G}{X-A})$ 
(Godement\footnote[2]{\textit{Th\'eorie des Faisceaux}, Hermann (1964), 234-236.}).\\
\endgroup %%------------------------------------<<

\begingroup%%----------------------------------->>
\fontsize{9pt}{11pt}\selectfont
\textbf{\small FACT} \  
Let $\sF \neq 0$ be a sheaf on \mX $-$then $\dim_{\sF} X \leq n$ iff $\sF$ admits a soft resolution 
$0 \ra \sF$ 
$\ra \sS^0$ 
$\ra \sS^1$
$\ra \cdots \ra \sS^n$ of length $n$.\\
\endgroup %%------------------------------------<<

\begingroup%%----------------------------------->>
\fontsize{9pt}{11pt}\selectfont
\textbf{\small LEMMA} \ 
Let $\{\sF_\alpha\}$ be a collection of soft subsheaves of a sheaf $\sF$ which is directed by inclusion.  Assume: 
$\sF = \colim \sF_\alpha$ $-$then $\sF$ is soft.\\
\endgroup %%------------------------------------<<

\begingroup%%----------------------------------->>
\fontsize{9pt}{11pt}\selectfont
\textbf{\small FACT} \  
Let $\{\sF_\alpha\}$ be a collection of subsheaves of a sheaf $\sF$ which is directed by inclusion.  Assume: 
$\sF = \colim \sF_\alpha$ $-$then $\dim_{\sF} X \leq n$ if $\forall \ \alpha$, $\dim_{\sF_\alpha} X \leq n$, hence 
$\dim_{\sF} X \leq \sup\dim_{\sF_\alpha} X.$
\vspi
[Work with the canonical (=Godement) resolution of the $\sF_\alpha$.]\\
\endgroup %%------------------------------------<<

\begingroup%%----------------------------------->>
\fontsize{9pt}{11pt}\selectfont
If $\sF = \sF^\prime \oplus \sF\pp$, then 
$\cH{}^n(X;\sF) = \cH{}^n(X;\sF^\prime) \oplus \cH{}^n(X;\sF\pp)$.  
Therefore, 
$\dim_{\sF^\prime} X \leq n$ $\&$ $\dim_{\sF\pp} X \leq n$ $\implies$ $\dim_{\sF} X \leq n$.
\vspi
Suppose now that $\{\sF_i\}$ is a collection of sheaves indexed by a set \mI.  
Given a finite subset $F \subset I$, put 
$\sF_F = \ds\bigoplus\limits_{i \in F} \sF_i$ $-$then $\sF \equiv \ds\bigoplus\limits_i \sF_i$ $= \colim \sF_F$.  
So, under the assumption that $\dim_{\sF_i} X \leq n$ $\forall \ i$, one has $\dim_{\sF} X \leq n$ as well.\\
\endgroup %%------------------------------------<<

Fix an abelian group \mG and let \mX be a paracompact Hausdorff space $-$then \mX is said to satisfy 
\un{Okuyama's condition at $n$}
\index{Okuyama's condition at $n$} 
if $\forall \ q \geq n$ and each closed subset \mA of \mX, the arrow $\cH{}^q(X;G) \ra \cH{}^q(A;G)$ is surjective.\\

\textbf{\small SUBLEMMA} \ Suppose that \mX satisifes Okuyama's condition at $n$ $-$then every closed subspace \mA of \mX satisfies Okuyama's condition at $n$.

%%----------------------------------------------------------------------------------------------06
[Given a closed subset \mB of \mA, consider the commutative triangle\\
\begin{tikzcd}%[ sep=small]
{\cH{}^q(X;G)} \ar{d} \ar{r} &{\cH{}^q(B;G)}\\
{\cH{}^q(A;G)} \ar{ru}
\end{tikzcd}
.]
\\
\vspace{0.25cm}

\textbf{\small LEMMA} \ 
Suppose that \mX satisifes Okuyama's condition at $n$.  Let 
$
\begin{cases}
\ A\\
\ B
\end{cases}
$
be closed subspaces of \mX and let $e$ be an element of $\cH{}^q(X;G)$ such that 
$
\begin{cases}
\ \restr{e}{A} = 0\\
\ \restr{e}{B} = 0
\end{cases}
$
for some $q > n$ $-$then $\restr{e}{A \cup B} = 0$.

[Consider the Mayer-Vietoris sequence 
$\cdots \ra \cH{}^{q-1}(A;G) \oplus \cH{}^{q-1}(B;G)$
$\overset{i}{\ra} \cH{}^{q-1}(A \cap B;G)$ 
$\ra \cH{}^{q-1}(A \cup B;G)$ 
$\overset{j}{\ra} \cH{}^{q}(A;G) \oplus \cH{}^{q}(B;G)$ 
$\ra \cdots$.  Thanks to the sublemma, $i$ is surjective.  Therefore $j$ is injective.  But 
$j(\restr{e}{A \cup B}) = 0$, so $\restr{e}{A \cup B} = 0$.]\\

\begin{proposition} \ %04
Let \mX be a paracompact Hausdorff space $-$then $\dim_G X \leq n$ iff \mX satisifes Okuyama's condition at $n$.
\end{proposition}

[Necessity: \ Inspect the exact sequence 
$\cdots \ra \cH{}^q(X,A;G)$ 
$\ra \cH{}^q(X;G)$ 
$\ra \cH{}^q(A;G)$ 
$\cH{}^{q+1}(X,A;G) \ra \cdots$.

Sufficiency: Fix $q \geq n$ $-$then since $\cH{}^q(X;G) \ra \cH{}^q(A;G)$ is surjective, 
$\cH{}^{q+1}(X,A;G) \ra \cH{}^{q+1}(X;G)$ is injective, thus it need only be shown that $\cH{}^{q+1}(X;G) = 0$.  
Take an $e \in \cH{}^{q+1}(X;G)$.  
Because $\cH{}^{q+1}(\{x\};G) = 0$, $\exists$ a neighborhood $U_x$ of x such that 
$\restr{e}{\ov{U}_x} = 0$ (restriction principle) 
and by paracompactness, the open covering $\{U_x:x \in X\}$ admits a $\sigma$-discrete closed refinement 
$\sA = \bigcup\limits_k \sA_k$.  
Put $A_k = \cup \sA_k$ and inductively determine a sequence 
$\{U_k\}$ of open sets: $A_k \cup \ov{U}_{k-1} \subset U_k$ $\&$ 
$\restr{e}{\ov{U}_k} = 0$, where $U_0 = \emptyset$.  
Noting that $\restr{e}{A_k} = 0$ $\forall \ k$, proceed as follows.  
First, $\exists$ an open $U_1 \supset A_1$: $\restr{e}{\ov{U}_1} = 0$, 
hence 
$\restr{e}{A_2 \cup \ov{U}_1} = 0$ (apply the preceding lemma). 
Assuming that 
$U_k \supset A_k \cup \ov{U}_{k-1}$ with  $\restr{e}{\ov{U}_k} = 0$ has been constructed, 
one has again $\restr{e}{A_{k+1} \cup {\ov{U}_k}} = 0$, so 
$\exists$ an open set $U_{k+1}$ : $U_{k+1} \supset A_{k+1} \cup \ov{U}_k = 0$ $\&$ $\restr{e}{\ov{U}_{k+1}} = 0$, 
which pushes the induction forward.  
Now let $W_k = \ov{U}_k - U_{k-1}$ : $W_k$ is closed, $\restr{e}{W_k} = 0$, and $X = \bigcup\limits_k W_k$.  
On the other hand, the collections $\{W_1, W_3, \ldots\}$, $\{W_2, W_4, \ldots\}$ are discrete.  
Therefore the restriction of $e$ to their respective unions must vanish, thus from the lemma, $e = 0$.]\\

Notation: Write $K(G,q)$ for an Eilenberg-MacLane space of type $(G,q)$ realized as an ANR in NES(paracompact) 
(cf. p. \pageref{20.5}).\\

\begin{proposition} \ %05
Let \mX be a paracompact Hausdorff space $-$then \mX satisifes Okuyama's condition at $n$ iff every closed subset $A \subset X$ has the EP w.r.t. $K(G,q)$ $\forall$ $q \geq n$.
\end{proposition}

[There are two points: 
(1) $\cH{}^q(X;G) \approx [X,K(G,q)]$, $\cH{}^q(A;G) \approx [A,K(G,q)]$; 
(2) \mA has the HEP w.r.t. $K(G,q)$ 
(cf. p. \pageref{20.6}).]\\

%%----------------------------------------------------------------------------------------------07
Application: Let \mX be a paracompact Hausdorff space $-$then $\dim_G X \leq n$
iff every closed subset $A \subset X$ has the EP w.r.t. $K(G,q)$ $\forall$ $q \geq n$.\\

\begin{proposition} \ %06
The following conditions on  a paracompact Hausdorff space \mX are equivalent.  
$(1)_n$ \ $\forall$ closed $A \subset X$ : $\cH{}^q(X,A;G) = 0$ $(q > n)$;
$(2)_n$ \ $\forall$ closed $A \subset X$ : $\cH{}^{n+1}(X,A;G) = 0$; 
$(3)_n$ \ $\forall$ closed $A \subset X$ : $\cH{}^q(X;G) \twoheadrightarrow \cH{}^n(A;G)$.
\end{proposition}

[Trivially, $(1)_n$ $\implies$ $(2)_n$, $(2)_n$ $\implies$ $(3)_n$.  
And: 
$(4)_n$ $\implies$ $(3)_{n+1}$, $(3)_n \wedge (4)_n$ $\implies$ $(2)_n$, where $(4)_n$ is the condition 
$\cH{}^{n+1}(A;G) = 0$ $\forall$ closed $A \subset X$.  
In addition, $(1)_n = \bigwedge\limits_{q \geq n} (2)_q$.  
Suppose that $(3)_n$ holds $-$then the claim is that $(3)_q \wedge (4)_q$ holds for $q \geq n$, hence that $(1)_n$ holds.  
Here is the pattern for the argument: 
$(3)_n$ $\implies$ $(4)_n$ $\implies$ $(3)_{n+1}$ $\implies$ $(4)_{n+1}$ $\implies$ $\cdots$.  
Therefore one has to show that $(3)_q$ $\implies$ $(4)_q$ $\forall \ q \geq n$.  But $(3)_q$ gives 
$\cH{}^{q+1}(X;G) = 0$ (see the proof of the sufficiency in Proposition 4) and since $(3)_q$ is inherited by \mA, 
$\cH{}^{q+1}(A;G) = 0$ too.]\\

Application: Let \mX be a paracompact Hausdorff space $-$then $\dim_G X \leq n$
iff every closed subset $A \subset X$ has the EP w.r.t. $K(G,n)$.

[Note: \ This result is the cohomological counterpart to Alexandroff's criterion.  If \mX is compact or stratifiable, then one can take $K(G,n)$ as a CW complex 
(cf. p. \pageref{20.7}).]\\

\label{20.16}
\begingroup%%----------------------------------->>
\fontsize{9pt}{11pt}\selectfont
\textbf{\small EXAMPLE} \quad 
Suppose that \mX is an ANR and let $G = \ds\prod\limits_i G_i$ be the direct product of abelian groups $G_i \neq 0$ $-$then 
$\dim_G X = \sup\dim_{G_i} X$.
\vspi
[Since each $G_i$ is a direct summand of \mG, $\dim_G X \geq \dim_{G_i} X$ $\forall \ i$, so if 
$\sup\dim_{G_i} X = \infty$, we are done.  Assume therefore that 
$\sup\dim_{G_i} X = n$.  Consider the product 
$Y = \ds\prod\limits_i K(G_i,n)$ $-$then every closed subset $A \subset X$ has the EP w.r.t. \mY, hence every closed subset 
$A \subset X$ has the EP w.r.t. $\abs{\sin Y}$ 
(cf. p. \pageref{20.8}).  But $\abs{\sin Y}$ is a CW complex and, as such, is an Eilenberg-MacLane space of type $(G,n)$.]\\
\endgroup %%------------------------------------<<

\begin{proposition} \ %07
Let \mX be a nonempty paracompact Hausdorff space $-$then $\dim X = 0$ iff $\dim_G X = 0$ $\forall \ G \neq 0$.
\end{proposition}

[By Proposition 2, $\dim X = 0$ $\implies$ $\dim_G X = 0$.  
Conversely, since \mG (discrete topology) $= K(G,0) \in$ NES(paracompact) contains $\bS^0$ as a retract (\mG being nontrivial), every closed subset $A \subset X$ has the EP w.r.t. $\bS^0$, hence $\dim X \leq 0$ (Alexandroff's criterion).]\\

Examples: $\forall \ G \neq 0$, 
(1) $\dim_G [0,1] = 1$; 
(2) $\dim_G \R = 1$; 
(3) $\dim_G \bS^1 = 1$.\\

\begingroup%%----------------------------------->>
\fontsize{9pt}{11pt}\selectfont
\textbf{\small EXAMPLE} \quad 
Let \mX be a paracompact Hausdorff space of finite topological dimension $-$then 
$\dim_G \beta X \leq \dim_G X$ provided that \mG is finitely generated.
\vspi
%%----------------------------------------------------------------------------------------------08
[It suffices to show that $\dim_G X \leq n$ $\implies$ $\dim_G \beta X \leq n$.  This is trivial if $X = \emptyset$ or 
$G = 0$, so take \mX nonempty and \mG nonzero.  Because $\dim \beta X = \dim X$ (cf. $\S 19$, Proposition 1), from Proposition 7, $\dim_G X = 0$ 
$\implies$ $\dim X = 0$ 
$\implies$ $\dim \beta X = 0$ 
$\implies$ $\dim_G \beta X = 0$.  
Suppose now that $n$ is positive and let \mA be a closed subset of $\beta X$.  
Claim: $\cH{}^{n+1}(A;G) = 0$, which is enough (cf. Proposition 6: $(1)_n \Leftrightarrow (4)_n)$.  
To verify this, fix an $\alpha \in \cH{}^{n+1}(A;G)$ and, using the extension principle, choose an open $U \supset A$ and an 
$e \in \cH{}^{n+1}(\ov{U};G)$: $\restr{e}{A} = \alpha$.  Since 
$\beta(\ov{U} \cap X) = \cl_{\beta X}(\ov{U} \cap X) = \ov{U}$, 
$\cH{}^{n+1}(\ov{U};G) \approx \cH{}^{n+1}(\ov{U} \cap X;G)$ 
(cf. p. \pageref{20.9}).  But 
$\dim_G X \leq n$ $\implies$ $\cH{}^{n+1}(\ov{U}\cap X;G) = 0$, so $e = 0$, thus $\alpha  = 0$.]
\vspi
[Note: \ Consequently, under the stated hypotheses on \mX and \mG, $\dim_G X = \dim_G \beta X$ 
(cf. p. \pageref{20.10}).]\\
\endgroup %%------------------------------------<<

\begingroup%%----------------------------------->>
\fontsize{9pt}{11pt}\selectfont
Remark: If the topological dimension of \mX is infinite, then one can find examples for which 
$\dim_{\Z} X \neq \dim_{\Z} \beta X$ 
(Dranishnikov\footnote[2]{\textit{C.R. Acad. Bulgare Sci.} \textbf{41} (1988), 9-10; 
see also, Dydak-Walsh, \textit{Proc. Amer. Math. Soc.} \textbf{113} (1991), 1155-1162 
and Dydak, \textit{Topology Appl.} \textbf{50} (1993), 1-10.}
).\\
\endgroup %%------------------------------------<<

\begingroup%%----------------------------------->>
\fontsize{9pt}{11pt}\selectfont
\textbf{\small EXAMPLE} \quad
For any nonempty paracompact Hausdorff space \mX, $\dim_{\Z} X = 1$ iff $\dim X = 1$.
\vspi
[If  $\dim_{\Z} X = 1$, then every closed subset $A \subset X$ has the EP w.r.t $\bS^1 = K(\Z,1)$, hence 
$\dim X \leq 1$ (Alexandroff's criterion) and $\dim X = 0$ is untenable (cf. Proposition 7).]\\
\endgroup %%------------------------------------<<

\begin{proposition} \ %08
Let \mX be a paracompact Hausdorff space $-$then for any closed subspace \mA of \mX, $\dim_G A \leq \dim_G X$.\\
\end{proposition}

\begingroup%%----------------------------------->>
\fontsize{9pt}{11pt}\selectfont
\textbf{\small EXAMPLE} \quad 
Let \mX be a paracompact LCH space $-$then $\dim_G X = \sup\dim_G K$, where $K \subset X$ is compact.
\vspi
[Since $\dim_G X \geq \dim_G K$ $\forall \ K$ (cf. Proposition 8), 
$\sup\dim_G K = \infty$ $\implies$ $\dim_G X = \infty$.  So suppose that $\sup\dim_G K = n$.  Write 
$X =\ds \bigcup\limits_i K_i$, where $K_i$ is compact and $\{K_i:i \in I\}$ is neighborhood finite.  Well order I and deduce by transfinite induction that every closed subset $A \subset X$ has the EP w.r.t. $K(G,n)$, hence that $\dim_G X \leq n$.]\\
\endgroup %%------------------------------------<<

\begingroup%%----------------------------------->>
\fontsize{9pt}{11pt}\selectfont
\textbf{\small FACT} \  
Let \mX be a closed subset of $\R^n$ $-$then $\dim X = n - 1$ iff $\dim_G X = n - 1$ $\forall \ G \neq 0$.\\
\endgroup %%------------------------------------<<

\begin{proposition} \ %09
Let \mX be a paracompact Hausdorff space.  Suppose that $X = \bigcup\limits_1^\infty A_j$, where the $A_j$ are closed subspaces of \mX such that $\forall \ j$, $\dim_G A_j \leq n$ $-$then $\dim_G X \leq n$, hence 
$\dim_G X = \sup\dim_G A_j$.
\end{proposition}

[Fix a closed subset $A \subset X$ and a continuous function $f:A \ra K(G,n)$.  
Put $U_0 = A$ and $F_0 = f$ $-$then 
since $\dim_G A_1 \leq n$, $\restr{F_0}{U_0 \cap A_1}$ has a continuous extension $\Phi_0:A_1 \ra K(G,n)$.  
Define a continuous 
function $f_1:U_0 \cup A_1 \ra K(G,n)$ by $\restr{f_1}{U_0} = F_0$ $\&$ 
%%----------------------------------------------------------------------------------------------09
$\restr{f_1}{A_1} = \Phi_0$.  Recalling that $K(G,n) \in $ NES(paracompact), choose an open 
$U_1 \supset U_0\  \cup \  A_1$ and a continuous function 
$F_1:\ov{U}_1 \ra K(G,n)$ such that $\restr{F_1}{U_0 \cup A_1} = f_1$.  
Since 
$\dim_G A_2 \leq n$, $\restr{F_1}{\ov{U}_1 \cap A_2}$ has a continuous extension $\Phi_1:A_2 \ra K(G,n)$.  
Define a continuous function $f_2:\ov{U}_1 \cup A_2 \ra K(G,n)$ by 
$\restr{f_2}{\ov{U}_1} = F_1$ $\&$ $\restr{f_2}{A_2} = \Phi_1$.  \ 
Choose an open 
$U_2 \supset$ $\ov{U}_1$ $\cup$  $A_2 $ 
and a continuous function 
$F_2:\ov{U}_2 \ra K(G,n)$ such that $\restr{F_2}{\ov{U}_1 \cup A_2} = f_2$.  
Continue the process to get a sequence of open sets $U_j$ $(j \geq 1)$: $\ov{U}_j \cup A_{j+1} \subset U_{j+1}$ and a sequence of continuous functions $F_j:\ov{U}_j \ra K(G,n)$ $(j \geq 1)$: $\restr{F_{j+1}}{\ov{U}_j} = F_j$.  
Finally, if $F:X \ra K(G,n)$ is defined by $\restr{F}{\ov{U}_j} = F_j$, then \mF is a continuous extension of $f$ 
($X = \bigcup\limits_j \ov{U}_j$ has the final topology corresponding to the inclusions $\ov{U}_j \ra X$).]\\

Proposition 9 is the analog for cohomological dimension of the countable union lemma for topological dimension but there are instances where the parallel breaks down.  
Here is a case in point.  
Suppose that $X = Y \cup Z$ is metrizable $-$then 
$\dim X \leq \dim Y + \dim Z + 1$ (cf. $\S 19$, Proposition 7).  
The situation for cohomological dimension is more complicated.

\indent\indent (R) \quad For any ring \mR with unit. 
$\dim_R X \leq \dim_R Y + \dim_R Z + 1$.

\indent\indent (G) \quad For any abelian group $G \neq 0$, $\dim_G X \leq \dim_G Y + \dim_G Z + 2$.

[Note: \ These estimates cannot be improved.  
See Dydak\footnote[2]{\textit{Trans. Amer. Math. Soc.} \textbf{348} (1996), 1647-1661.}
for details and references.]\\

\begingroup%%----------------------------------->>
\fontsize{9pt}{11pt}\selectfont
\label{20.13}
\textbf{\small FACT} \  
Suppose that \mX is a paracompact Hausdorff space.  Let $\sA = \{A_j:j \in J\}$ be an absolute closure preserving closed covering of \mX such that $\forall \ j$, $\dim_G A_j \leq n$ $-$then $\dim_G X \leq n$, hence 
$\dim_G X = \sup\dim_G A_j$.\\
\endgroup %%------------------------------------<<

\begingroup%%----------------------------------->>
\fontsize{9pt}{11pt}\selectfont
\textbf{\small LEMMA} \ 
If \mK is a finite CW complex, then $\dim_G K = \dim K$ $\forall \ G \neq 0$.
\vspi
[On general grounds, $\dim_G K \leq \dim K$ (cf. Proposition 2).  Taking $K \neq \emptyset$, let 
$n = \dim_G K > 0$ (cf. Proposition 7), and suppose that $k = \dim K > n$.  
Fix a $k$-cell $e \subset K$ and let $\bS^n$ be an $n$-sphere contained in $e$.  
Since $G \neq 0$, $\exists$ a map $f:\bS^n \ra K(G,n)$ which induces a nontrivial homomorphism 
$\Z = \pi_n(\bS^n) \ra \pi_n(K(G,n)) = G$.  
But $f$ admits a continuous extension $K \ra K(G,n)$.  
Therefore $\pi_n(f)$ is trivial, $\bS^n$ being contractible in \mK.  Contradiction.]\\ 
\endgroup %%------------------------------------<<

\begingroup%%----------------------------------->>
\fontsize{9pt}{11pt}\selectfont
\textbf{\small EXAMPLE} \quad 
Let \mX be a CW complex $-$then the collection $\{K\}$ of finite subcomplexes of 
\mX is an absolute closure preserving closed covering of \mX, thus $\dim X = \sup\dim K$ 
(cf. p. \pageref{20.11}).  On the other hand, $\forall \ G \neq 0$, 
$\dim_G X = \sup\dim_G K$ (cf. supra) and by the lemma, $\dim_G K = \dim K$.  Therefore $\dim_G X = \dim X$.\\
\endgroup %%------------------------------------<<

\begingroup%%----------------------------------->>
\fontsize{9pt}{11pt}\selectfont
Examples: $\forall \ G \neq 0$, 
(1) $\dim_G [0,1]^n = n$; 
(2) $\dim_G \R^n = n$; 
(3) $\dim_G \bS^n = n$.\\
\endgroup %%------------------------------------<<

%%----------------------------------------------------------------------------------------------10
\begingroup%%----------------------------------->>
\fontsize{9pt}{11pt}\selectfont

\textbf{\small EXAMPLE} \quad 
Let \mX be a paracompact $n$-manifold $-$then $\forall \ G \neq 0$, $\dim_G X = n$ 
(cf. p. \pageref{20.12}).\\

\endgroup %%------------------------------------<<

\begin{proposition} \ %10
Let \mX be a paracompact Hausdorff space.  Assume: \mX is hereditarily paracompact $-$then for any subspace \mY of \mX, 
$\dim_G Y \leq \dim_G X$.\\
\end{proposition}

\begin{proposition} \ %11
Let \mX be a paracompact Hausdorff space. Suppose that \mY is a strongly paracompact subspace $-$then 
$\dim_G Y \leq \dim_G X$.\\
\end{proposition}

\begingroup%%----------------------------------->>
\fontsize{9pt}{11pt}\selectfont
\textbf{\small EXAMPLE} \quad 
Suppose that \mX contains an embedded copy of $\R^n$ $-$then $\forall \ G \neq 0$, $\dim_G X \geq n$.\\
\endgroup %%------------------------------------<<

\textbf{\small LEMMA} \ 
Let \mX be a nonempty paracompact Hausdorff space, \mY a nonempty compact Hausdorff space.  Assume: 
$
\begin{cases}
\ \dim X\\
\ \dim Y
\end{cases}
< \infty
$
$-$then $\forall \ G \neq 0$, $\dim_G(X \times Y)$ is the largest integer $n$ such that
$\cH{}^n((A^\prime,A) \times (B^\prime,B);G) \neq 0$ for certain closed sets 
$
\begin{cases}
\ A \subset A^\prime \subset X\\
\ B \subset B^\prime \subset Y
\end{cases}
.
$

[By the product theorem, $\dim (X \times Y) \leq \dim X + \dim Y$, so Proposition 2 implies that 
$\dim_G (X \times Y)$ is finite.  This said, to prove the lemma, it suffices to show that whenever $m > n$ and 
$\cH{}^m((A^\prime,A)  \times (B^\prime,B);G) = 0$ for all closed sets
$
\begin{cases}
\ A \subset A^\prime \subset X\\
\ B \subset B^\prime \subset Y
\end{cases}
,\ 
$
then $\dim_G (X \times Y) \leq n$.  Thus let $W \subset X \times Y$ be closed.  
Fix a continuous function 
$f:W \ra K(G,n)$ $-$then $\exists$ an open $U \supset W$: $f$ is continuously extendable over \mU.  
The open covering 
$\sW = \{U,X \times Y - W\}$ is numerable, hence by the stacking lemma, there exists a neighborhood finite open covering 
$\sU = \{U_i: i \in I\}$ of \mX and $\forall \ i \in I$ a finite open covering 
$\sV_i = \{V_{i,j}:j \in J_i\}$ of \mY such that the collection 
$\{U_i \times \sV_i: i \in I\}$ refines $\sW$.
Choose a neighbhorhood finite open covering $\sO = \{O_\lambda:\lambda \in \Lambda\}$ of \mX of order 
$\leq \dim X + 1$ such that $\{\st(x,\sO):x \in X\}$ is a refinement of $\sU$ (cf. $\S 19$,  Proposition 6).
Given $\xi = (\lambda_1, \ldots, \lambda_p) \in \Lambda^p$, put
$A_\xi = X - \bigcup\limits_{\lambda \in \Lambda} \{O_\lambda:\lambda \neq \lambda_i \ \ (1 \leq i \leq p)\}$ if 
$\bigcap\limits_{i=1}^p O_{\lambda_i} \neq \emptyset$, otherwise put 
$A_\xi = \emptyset$.  The covering 
$\sA = \bigcup\limits_{p=1}^d \sA_p$, where $\sA_p = \{A_\xi:\xi \in \Lambda^p\}$ and $d = \dim X + 1$, is a neighborhood finite closed refinement of $\sU$.  For each $A_\xi \in \sA$, determine 
$U_{i(\xi)} \in \sU$: $A_\xi \subset U_{i(\xi)}$.  
Let 
$\sB_i = \{B_{i,j}:j \in J_i\}$ be a precise closed refinement of $\sV_i$.  
The collection 
$\{A_\xi \times B_{i(\xi),j}:A_\xi \in \sA, j \in J_{i(\xi)}\}$ is therefore a neighborhood finite closed refinement of 
$\sW = \{U, X \times Y - W\}$.  Set 
$A_0 = \bigcup\{A_\xi \times B_{i(\xi),j}: W \cap (A_\xi \times B_{i(\xi),j}) \neq \emptyset\}$.  
Since $W \subset A_0 \subset U$, $\exists$ a continuous function $f_0:A_0 \ra K(G,n)$ such that $\restr{f_0}{W} = f$.  
Now put 
$A_p = A_0 \cup \{A_\xi \times Y: \xi \in \Lambda^p\}$ $(1 \leq p \leq d)$ and assume that $f_0$ has a continuous extension 
$f_{p-1}:A_{p-1} \ra K(G,n)$ for some $p \geq 1$.   
Claim: $f_{p-1}$ can be continuously extended over $A_p$.  To see this, note first that 
$\xi,\ \xi^\prime \in \Lambda^p$ $\&$ $\xi \neq \xi^\prime$ $\implies$ 
$A_\xi \cap A_{\xi^\prime} \in \sA_{p-1}$, so it will be enough to establish that 
$f_{p-1,\xi} \equiv \restr{f_{p-1}}{(A_{p-1} \cap (\A_\xi \times Y))}$ is continuously extendable over $A_\xi \times Y$ for each 
$\xi \in \Lambda^p$.  Write 
$J_{i(\xi)} = \{j:1 \leq j \leq j_{i(\xi)}\}$.  Suppose inductively that
%%----------------------------------------------------------------------------------------------11
$f_{p-1,\xi}$ has been continuously extended over 
$(A_{p-1} \cap (A_\xi \times Y)) \cup \bigcup\limits_{1 \leq j \leq k-1} A_\xi \times B_{i(\xi),j}$ 
for some $k \leq j_{i(\xi)}$.  
Let 
$A^\prime = A_\xi$, 
$A = A_\xi \cap \bigl( \bigcup\limits_{\xi^\prime \in \Lambda^{p-1}} A_{\xi^\prime}\bigr)$, 
$B^\prime = B_{i(\xi),k}$, 
$B =  \bigcup\limits_{1 \leq j \leq k-1} B_{i(\xi),k} \cap B_{i(\xi),j}$ $\cup$ 
$\bigcup$ $\{B_{i(\xi),k} \cap B_{i(\xi^\prime),j}$: 
$W \cap (A_{\xi^\prime} \times B_{i(\xi^\prime),j}) \neq \emptyset\}$ 
$-$then from the exact sequence
$\cdots \ra $ 
$\cH{}^n(A^\prime \times B^\prime;G) \ra$ 
$\cH{}^n(A^\prime \times B \cup A \times B^\prime;G) \ra $ 
$\cH{}^{n+1}((A^\prime,A) \times (B^\prime,B);G) \ra $ 
$\cdots$, it follows that the arrow 
$\cH{}^n(A^\prime \times B^\prime;G) \ra$ $\cH{}^n(A^\prime \times B \cup A \times B^\prime;G)$ is surjective.  
Accordingly, every continuous function 
$A^\prime \times B \cup A \times B^\prime \ra K(G,n)$ can be extended to a continuous function 
$A^\prime \times B^\prime \ra K(G,n)$.  
In particular: $f_{p-1,\xi}$ is continuously extendable over 
$(A_{p-1} \cap (A_\xi \times Y)) \cup \bigcup\limits_{1 \leq j \leq k} A_\xi \times B_{i(\xi),j}$, which completes the induction.  Consequently, $f_{p-1}$ extends to a continuous function 
$f_p:A_p \ra K(G,n)$, hence by induction once again, $f$ extends to a continuous function 
$f_d:X \times Y \ra K(G,n)$.]\\

\begin{proposition} \ %12
Let \mX be a nonempty paracompact Hausdorff space, \mY a nonempty compact Hausdorff space.  Assume: 
$
\begin{cases}
\ \dim X\\
\ \dim Y
\end{cases}
< \infty
$
$-$then $\forall \ G \neq 0$, $\dim_G (X \times Y)$ is the largest integer $n$ such that 
$\cH{}^n((X,A) \times (Y,B);G) \neq 0$ for certain closed sets
$
\begin{cases}
\ A \subset X\\
\ B \subset Y
\end{cases}
. \ 
$
\end{proposition}

[Suppose that $n = \dim_G (X \times Y)$.  Using the lemma, choose closed sets 
$
\begin{cases}
\ A \subset A^\prime \subset X\\
\ B \subset B^\prime \subset Y
\end{cases}
$
such that 
$\cH{}^n((A^\prime,A) \times (B^\prime,B);G) \neq 0$.  
Put
$C = A^\prime \times B^\prime \cup X \times B \cup A \times Y$ 
$-$then 
$(A^\prime \times B^\prime) \cap (X \times B \cup A \times Y) =$ 
$A^\prime \times B \cup A \times B^\prime$, thus by the relative homeomorphism theorem, 
$\cH{}^n(C,X \times B \cup A \times Y;G) \approx$ 
$\cH{}^n(A^\prime \times B^\prime, A^\prime \times B \cup A \times B^\prime;G) \neq 0$.  
Consider the exact sequence 
$\cdots \ra \cH{}^n((X,A) \times (Y,B);G)$ 
$\ra \cH{}^n(C, X \times B \cup A \times Y; G)$ 
$\ra \cH{}^{n+1}(X \times Y,C; G)$
$\ra \cdots$ corresponding to the triple $(X \times Y, C, X \times B \cup A \times Y)$.  Since $n = \dim_G (X \times Y)$, 
$\cH{}^{n+1}(X \times Y, C;G) = 0$, hence $\cH{}^n((X,A) \times (Y,B);G) \neq 0$.]\\

Application: Under the preceding hypotheses on $X \ \& \ Y$, $\dim_G (X \times Y) \leq n$ iff 
$\cH{}^q((X,A) \times (Y,B);G) = 0$ $\forall \ q > n$ and for all closed sets
$
\begin{cases}
\ A \subset X\\
\ B \subset Y
\end{cases}
.
$
\\

\label{20.19}
\label{20.20}
\label{20.21}
\label{20.26}
\begingroup%%----------------------------------->>
\fontsize{9pt}{11pt}\selectfont
\textbf{\small EXAMPLE} \quad 
With \mX $\&$ \mY as in Proposition 12, suppose that $\exists \ k$: 
$\dim_{\cH{}^{k-i}(Y,B;G)} X \leq i$ $\forall \ i \geq 0$ and all closed subsets $B \subset Y$ $-$then
$\dim_G (X \times Y) \leq k$.
\vspi
[It is a question of verifying that 
$\cH{}^l((X,A) \times (Y,B);G) = 0$ $\forall \ l \geq k + 1$.  
But by the K\"unneth formula, 
$\cH{}^l((X,A) \times (Y,B);G) \approx$ 
$\ds\bigoplus\limits_{q= 0}^l \cH{}^q(X,A;\cH{}^{l-q}(Y,B;G)) \approx$ 
$\cH{}^q(X,A;\cH{}^{k-(q-l+k)}(Y,B;G)) = 0$.]\\
\endgroup %%------------------------------------<<

\begingroup%%----------------------------------->>
\fontsize{9pt}{11pt}\selectfont
\textbf{\small EXAMPLE} \quad 
With \mX $\&$ \mY as in Proposition 12, suppose that $\dim_{\cH{}^m(Y,B;G)} X \geq n$ for some closed subset 
$B \subset Y$ $-$then 
$\dim_G(X \times Y) \geq n + m$.
\vspi
[Choose a closed subset $A \subset X$: 
$\cH{}^n(X,A;\cH{}^m(Y,B;G)) \neq 0$ $\implies$ 
$\cH{}^{n+m}((X,A) \times (Y,B);G) \approx$ 
$\ds\bigoplus\limits_{q=0}^{n+m} \cH{}^q(X,A;\cH{}^{n+m-q}(Y,B;G)) \neq 0$, hence 
$\dim_G (X \times Y) \geq n + m$.]\\
\endgroup %%------------------------------------<<

%%----------------------------------------------------------------------------------------------12
\begin{proposition} \ 
Let \mX be a nonempty paracompact Hausdorff space of finite topological dimension $-$then $\forall \ G \neq 0$, 
$\dim_G IX = \dim_G X + 1$.
\end{proposition}

[$\dim_G IX \geq \dim_G X + 1$:\ Choose a closed subset $A \subset X$: 
$\cH{}^n(X,A;G) \neq 0$, where $n = \dim_G X$.  Applying the K\"unneth formula, we have 
$\cH{}^{n+1}((X,A) \times ([0,1],\{0,1\});G) \approx$ 
$\ds\bigoplus\limits_{q=0}^{n+1} \cH{}^q(X,A;\cH{}^{n+1-q}([0,1],\{0,1\};G) \approx$ 
$\cH{}^{n}(X,A;\cH{}^{1}([0,1],\{0,1\};G)) \approx$ 
$\cH{}^n(X,A;G) \neq 0$, 
which implies that $\dim_G IX \geq \dim_G X + 1$.

$\dim_G X + 1 \geq \dim_G IX$: \ Fix $m \geq n + 2$ $(n = \dim_G X)$ and let
$
\begin{cases}
\ A \subset X\\
\ B \subset Y
\end{cases}
$
be closed $(I = [0,1])$.  Utilization of  the K\"unneth formula then gives 
$\cH{}^m((X,A) \times (I,B);G) \approx$
$\cH{}^m(X,A;\cH{}^0(I,B;G))$ $\oplus$ $\cH{}^{m-1}(X,A;\cH{}^1(I,B;G))$.  
Case 1: $B = \emptyset$.  Here, 
$\cH{}^0(I,\emptyset;G) = G$, 
$\cH{}^1(I,\emptyset;G) = 0$, hence
$\cH{}^m((X,A) \times (I,B);G) =0$ 
Case 2: $B \neq \emptyset$.  Here, 
$\cH{}^0(I,B;G) = 0$,  
$\cH{}^1(I,B;G) = \cH{}^1(I,B;\Z)  \otimes G$ (by the universal coefficient theorem), hence 
$\cH{}^m((X,A) \times (I,B);G) \approx$ 
$\cH{}^{m-1}(X,A;\cH{}^1(I,B;\Z) \otimes G) = 0$ (cf. Proposition 18 $(m - 1 \geq n + 1)$).  
Therefore $\dim_G X + 1 \geq \dim_G IX$.]\\

Application: Let \mX be a nonempty paracompact Hausdorff space, \mY a nonempty CW complex.  Assume: 
$
\begin{cases}
\ \dim X \\
\ \dim Y
\end{cases}
< \infty
$
$-$then $\forall \ G \neq 0$, $\dim_G(X \times Y) = \dim_G X + \dim_G Y$ ($= \dim_G X + \dim Y$ 
(cf. p. \pageref{20.13})).

[If \mB is a compact subspace of \mY which is homeomorphic to $[0,1]^m$, where $m = \dim_G Y$, then 
$\dim_G (X \times B) = \dim_B X + m = \dim_G X + \dim_G Y$.]

[Note: \ \mY is paracompact and $\sigma$-locally compact, thus $X \times Y$ is paracompact 
(cf. p. \pageref{20.14}).]\\

\begingroup%%----------------------------------->>
\fontsize{9pt}{11pt}\selectfont
\textbf{\small FACT} \  
Let \mX be a nonempty paracompact Hausdorff space, \mY a nonempty compact Hausdorff space.  
Assume: 
$\dim X < \infty$ $\&$ $\dim Y = 0$ $-$then $\forall \ G \neq 0$, $\dim_G(X \times Y) = \dim_G X$.
\vspi
[It is clear that $\dim_G (X \times Y) \geq \dim_G X$ (cf. Proposition 8).  With $n = \dim_G X$, fix 
$m \geq n + 1$ and let 
$
\begin{cases}
\ A \subset X\\
\ B \subset Y
\end{cases}
$
be closed.  From the K\"unneth formula, 
$\cH{}^m((X,A) \times (Y,B); G) \approx \ds\bigoplus\limits_{q=0}^m \cH{}^q(X,A;\cH{}^{m-q}(Y,B;G))$.  
But $\dim Y = 0$ $\implies$ $\dim_G Y = 0$ (cf. Proposition 2), so 
$\cH{}^{m-q}(Y,B;G) = 0$ if $q \leq m-1$, thus 
$\cH{}^m((X,A) \times (Y,B); G) \approx$ 
$\cH{}^m(X,A;\cH{}^{0}(Y,B;G)) \approx$ 
$\cH{}^m(X,A;\cH{}^{0}(Y,B;\Z)\ \otimes G) = 0$ (cf. Proposition 18).  Therefore 
$\dim_G(X \times Y) \leq n$.]\\
\endgroup %%------------------------------------<<

\begin{proposition} \ %14
Let \mX be a paracompact Hausdorff space.  Suppose that $\{G_\alpha\}$ is a collection of subgroups of an abelian group \mG which is directed by inclusion.  Assume: $G = \colimx G_\alpha$ $-$then 
$\dimx_G X \leq n$ if $\forall \ \alpha$, $\dimx_{G_\alpha} \leq n$, hence 
$\dimx _G X \leq \sup\dimx_{G_\alpha} X$.
\end{proposition}

[This is a special case of the generalities on 
p. \pageref{20.15}.]\\

%%----------------------------------------------------------------------------------------------13
\index{Direct Sum Criterion}
\textbf{\small DIRECT SUM CRITERION} \quad 
Let \mX be a paracompact Hausdorff space $-$then 
$\dim_{\underset{i}{\oplus} G_i} X = \sup\dimx_{G_i} X$.

[Apply Proposition 14 
(cf. p. \pageref{20.15a}).]\\

\begingroup%%----------------------------------->>
\fontsize{9pt}{11pt}\selectfont
\textbf{\small EXAMPLE} \quad Since $\widehat{\Z}_p/\Z_p$ is a vector space over $\Q$, 
$\dimx_{\widehat{\Z}_p/\Z_p} X = \dimx_{\Q} X$.\\
\endgroup %%------------------------------------<<

\begin{proposition} \ %15
Let \mX be a paracompact Hausdorff space.  
Suppose that $0 \ra G^\prime \ra G \ra G\pp \ra 0$ is a short exact sequence of abelian groups $-$then 
$\dimx_G X \leq$ $\max\{\dimx_{G^\prime} X ,\dimx_{G\pp} X\}$, \ 
$\dimx_{G^\prime} X$ \ $\leq$ \ $\max\{\dimx_{G} X ,\dimx_{G\pp} X + 1\}$, 
and \ 
$\dimx_{G\pp} X$ $\leq$
$\max\{\dimx_{G} X ,\dimx_{G^\prime} X - 1\}$.
\end{proposition}

[Use the Bockstein sequence.]\\

\begingroup%%----------------------------------->>
\fontsize{9pt}{11pt}\selectfont
\index{Bockstein's Inequalities}
\textbf{\small EXAMPLE \  (\un{Bockstein's Inequalities})} \  
Let \mX be a paracompact Hausdorff space and fix a prime $p$.
\\
\indent\indent (BO$_1$) \quad $\dimx_{\Z/p\Z} X = \dimx_{\Z/p^n\Z} X$.
\vspi
[From the short exact sequence 
$0 \ra \Z/p^n\Z \ra \Z/p^{n+1}\Z \ra \Z/p\Z \ra 0$, it follows that 
$\dimx_{\Z/p^{n+1}\Z} X \leq \max\{\dimx_{\Z/p^n\Z} X, \dimx_{\Z/p\Z} X\}$
and 
$\dim_{\Z/p\Z} X \leq \max\{\dim_{\Z/p^{n+1}\Z} X, \dim_{\Z/p^n\Z} X - 1\}$.  
Now argue by induction.]
\\
\indent\indent (BO$_2$) \quad $\dimx_{\Z/p^\infty\Z} X \leq \dimx_{\Z/p\Z} X$.
\vspi
[Since $\Z/p^\infty\Z = \colimx \Z/p^n\Z$, Proposition 14 implies that 
$\dim_{\Z/p^\infty\Z} X \leq \sup \dimx_{\Z/p^n\Z} X = \dim_{\Z/p\Z}X$.]
\\
\indent\indent (BO$_3$) \quad $\dim_{\Z/p\Z} X \leq \dimx_{\Z/p^\infty\Z } X + 1$.
\vspi
[Consider the short exact sequence 
$0 \ra \Z/p\Z \ra \Z/p^\infty\Z \overset{p}{\ra} \Z/p^\infty\Z \ra 0$.]
\\
\indent\indent (BO$_4$) \quad $\dimx_{\Z/p\Z} X \leq \dimx_{\Z_p} X$.
\vspi
[Consider the short exact sequence 
$0 \ra \Z_p \overset{p}{\ra} \Z_p \ra \Z/p\Z \ra 0$.]
\\
\indent\indent (BO$_5$) \quad $\dimx_{\Q} X \leq \dimx_{\Z_p} X$.
\vspi
[Consider the short exact sequence 
$0 \ra \Z_p \ra \Q \ra \Z/p^\infty\Z \ra 0$.]
\vspi
[Note: \ In addition, 
$\dimx_{\Z_p} X \leq \max\{\dimx_{\Q} X, \dimx_{\Z/p^\infty\Z} X + 1\}$, 
$\dimx_{\Z/p^\infty\Z} X \leq \max\{\dimx_{\Q} X, \dimx_{\Z_p} X - 1\}$.]
\\[.1cm]
\endgroup %%------------------------------------<<

\begingroup%%----------------------------------->>
\fontsize{9pt}{11pt}\selectfont
Warning: Bockstein's inequalities are used without citation in the sequel.
\\[.1cm]
\endgroup %%------------------------------------<<

\begingroup%%----------------------------------->>
\fontsize{9pt}{11pt}\selectfont
\textbf{\small FACT} \  
Let \mX be a compact Hausdorff space.  Suppose that 
$0 \ra G^\prime \ra G \ra G\pp \ra 0$ is a short exact sequence of abelian groups.  Assume: $G\pp$ is torsion free $-$then $\dimx_G X = \max\{\dimx_{G^\prime} X ,\dim_{G\pp} X\}$.
\\[.1cm]
\endgroup %%------------------------------------<<

\label{20.17}
\begingroup%%----------------------------------->>
\fontsize{9pt}{11pt}\selectfont
\textbf{\small EXAMPLE} \quad 
Let \mX be a compact Hausdorff space $-$then $\dimx_{\Z_p} X = \dimx_{\widehat{\Z}_p} X$.
\vspi
[From the short exact sequence 
$0 \ra \Z_p \ra \widehat{\Z}_p \ra \widehat{\Z}_p/\Z_p \ra 0$, we have 
$\dimx_{\widehat{\Z}_p} X = \max\{\dimx_{\Z_p} X, \dimx_{\widehat{\Z}_p/\Z_p} X\}$.  But 
$\dimx_{\widehat{\Z}_p/\Z_p} X = \dimx_{\Q} X$ and $\dimx_{\Q} X \leq \dimx_{\Z_p} X$.]
\\[.1cm]
\endgroup %%------------------------------------<<

%%----------------------------------------------------------------------------------------------14
\label{20.24}
\begingroup%%----------------------------------->>
\fontsize{9pt}{11pt}\selectfont
A 
\un{Bockstein function} \ 
\index{Bockstein function} 
is a function \mD defined on \ 
$\{\Q\} \cup  \ds\bigcup\limits_p \{\Z_p,\Z/p\Z,\Z/p^\infty\Z\}$ with values in \ $\Z_{\geq 0\}} \cup \{\infty\}$ such that 
$D(\Z/p^\infty\Z) \leq D(\Z/p\Z)$, 
$D(\Z/p\Z) \leq D(\Z/p^\infty\Z) + 1$, 
$D(\Z/p\Z) \leq D(\Z_p)$,
$D(\Q) \leq D(\Z_p)$, 
$D(\Z_p) \leq \max\{D(\Q),D(\Z/p^\infty\Z) + 1\}$, 
$D(\Z/p^\infty\Z) \leq \max\{D(\Q),D(\Z_p) -1\}$, and 
$D$ is $\equiv 0$ if $D(G) = 0$ $(\exists \ G)$ (cf. Proposition 7).

Example: Every nonempty paracompact Hausdorff space \mX gives rise to a Bockstein function $D_X$, viz. 
$D_X(G) = \dim_G X$.
\\[.1cm]

\index{Theorem: Dranishnikov's Realization Theorem}
\index{Dranishnikov's Realization Theorem}
\textbf{\small 
DRANISHNIKOV'S\footnote[2]{\textit{Siberian Math J.} \textbf{29} (1988), 24-29, \textbf{30} (1989), 74-79, 
and \textbf{32} (1991), 145-147; see also Dydak, \textit{Topology Appl.} \textbf{65} (1995), 1-7.} REALIZATION THEOREM} \quad 
Given a Bockstein function $D$, $\exists$ a metrizable compact Hausdorff space \mX such that $D = D_X$ and 
$\dim X = \sup D$.\\

\textbf{\small EXAMPLE} \quad 
The 
\un{fundamental compacta}
\index{fundamental compacta} 
are those metrizable compact Hausdorff spaces which realize the Bockstein functions define by the table below. 

\endgroup %%------------------------------------<<
\[
\begin{array}{|l| c c c c c c c c c c c c c c c}
\hline
D &&\Z_p &&\Z/p\Z &&\Z/p^\infty\Z &&\Q &&\Z_q &&\Z/q\Z &&Z/q^\infty\Z\\
\hline
\Phi(\Q,n) &&n &&1 &&1 &&n &&n &&1 &&1\\
\\
\Phi(\Z_p,n) &&n &&n &&n &&n &&n &&1 &&1\\
\\
\Phi(\Z/p\Z,n) &&n &&n &&n-1 &&1 &&1 &&1 &&1\\
\\
\Phi(\Z/p^\infty\Z,n) &&n &&n-1 &&n-1 &&1 &&1 &&1 &&1\\
\end{array}
\]

\begingroup%%----------------------------------->>
\fontsize{9pt}{11pt}\selectfont
[Note: \ Here $p$, $q$ are primes, $q$ runs over all primes $\neq p$, and $\Phi(G,n)$ is the Bockstein function corresponding to the pair $(G,n)$, where $G = \Q$, $\Z_p$, $\Z/p\Z$, $\Z/p^\infty\Z$.]\\
\endgroup %%------------------------------------<<

Notation: Given an abelian group \mG, $G_\tor$ is its torsion subgroup and 
$G_\tor(p)$ is the $p$-primary component of 
$G_\tor$ (so $G_\tor \approx \bigoplus\limits_p G_\tor(p)$).

[Note: \ Accordingly, for a paracompact Hausdorff space $X$, 
$\dim_{G_\tor} X =$ $\sup\dimx_{G_{\tor(p)}} X$ (direct sum criterion).]

Given an abelian group \mG, its 
\un{Bockstein basis}
\index{Bockstein basis} 
$\sigma(G)$ is the subset of 
$\{\Q\} \cup \bigcup\limits_p \{\Z_p,$ $\Z/p\Z,\Z/p^\infty\Z\}$ defined as follows:

\indent\indent $(\Q)$ \qquad \quad \ $\Q \in \sigma(G)$ iff $G/G_\tor \neq 0$.

\indent\indent $(\Z_p)$ \qquad \quad $\Z_p \in \sigma(G)$ iff $G/G_\tor$ is not divisible by $p$.

\indent\indent $(\Z/p\Z)$ \quad \ \    $\Z/p\Z \in \sigma(G)$ iff $G_\tor(p)$ is not divisible by $p$.

%%----------------------------------------------------------------------------------------------15
\indent\indent $(\Z/p^\infty\Z)$ \quad $\Z/p^\infty\Z \in \sigma(G)$ iff $G_\tor(p) \neq 0$ is divisible by $p$.

Examples: 
(1) $\sigma(\Q) = \{\Q\}$;
(2) $\sigma(\Z_p) = \{\Q,\Z_p\}$;
(3) $\sigma(\Z/p\Z) = \{\Z/p\Z\}$;
(4) $\sigma(\Z/p^\infty\Z) = \{\Z/p^\infty\Z\}$;
(5) $\sigma(\Z) = \{\Q\} \cup \bigcup\limits_p\{\Z_p\}$;
(6) $\sigma(\widehat{\Z}_p) = \{\Q,\Z_p\}$.

Remark: $\forall \ G \neq 0$, $\sigma(G)$ is nonempty.  
Indeed, if $G \neq G_\tor$, then $\Q \in \sigma(G)$ and if $G = G_\tor$, then $\exists \ p$: $G_\tor(p) \neq 0$, so either $\Z/p\Z \in \sigma(G)$ or $\Z/p^\infty\Z \in \sigma(G)$.\\

\textbf{\small LEMMA} \ 
Given an abelian group \mG, $\sigma(G) = \sigma(G/G_\tor) \cup \bigcup\limits_p \sigma(G_\tor(p))$.\\

\begingroup%%----------------------------------->>
\fontsize{9pt}{11pt}\selectfont
\textbf{\small FACT} \  
If $G_\tor(p)$ is not divisible by $p$, then $\exists \ n \geq 1$ : $\Z/p^n\Z$ is a direct summand of \mG.\\
\endgroup %%------------------------------------<<

\begingroup%%----------------------------------->>
\fontsize{9pt}{11pt}\selectfont
\textbf{\small FACT} \  
If $G_\tor(p) \neq 0$ is divisible by $p$, then $G_\tor(p) \approx \oplus \Z/p^\infty\Z$ and $G_\tor(p)$ is a direct summand of \mG.\\
\endgroup %%------------------------------------<<

\begin{proposition} \ %16
Let \mX be a paracompact Hausdorff space.  Suppose that $G \neq 0$ is torsion $-$then 
$\dim_G X = \sup\limits_{H \in \sigma(G)} \dim_H X$.
\end{proposition}

[From what has been said above, one can assume that $G = G(p)$ $(\exists \ p)$.

\indent\indent $(\Z/p\Z)$ \ If $\Z/p\Z \in \sigma(G)$, then 
$\dim_{\Z/p\Z} X = \max\limits_{H \in \sigma(G)} \dim_H X$.  
But $\Z/p^n\Z$ is a direct summand of \mG for some 
$n \geq 1$, thus 
$\dim_G X \geq \dim_{\Z/p^n\Z}X = \dim_{\Z/p\Z} X$.  On the other hand, \mG is a colimit of its finite subgroups.  As these are direct sums of groups of the form $\Z/p^k\Z$, $\dim_G X \leq \dim_{\Z/p\Z}X$ by Proposition 14.

\indent\indent $(\Z/p^\infty\Z)$ \  In this case, \mG is isomorphic to a direct sum of copies of $\Z/p^\infty\Z$ and the direct sum criterion is applicable.]\\

\begin{proposition} \ %17
Let \mX be a paracompact Hausdorff space $-$then for any $G \neq 0$, 
$\dim_G X = \max\{\dim_{G/G_\tor} X, \dim_{G_\tor} X\}$.
\end{proposition}

[The short exact sequence 
$0 \ra G_\tor \ra G \ra G/G_\tor \ra 0$ leads to the inequalities 
$\dim_G X \leq \max\{\dim_{G_\tor} X, \dim_{G/G_\tor} X\}$, 
$\dim_{G/G_\tor} X \leq \max\{\dim_G X, \dim_{G_\tor} X - 1\}$ (cf Proposition 15), thus it suffices to prove that 
$\dim_G X \geq \dim_{G_\tor} X$.  But if $\Z/p\Z \in \sigma(G)$, then $\Z/p^n\Z$ is a direct summand of \mG 
$(\exists \ n \geq 1)$, while if $\Z/p^\infty\Z \in \sigma(G)$, then $\Z/p^\infty\Z$ is a direct summand of \mG.  Therefore $\dim_G X \geq \dim_{G_\tor} X$ (cf Proposition 16).]\\

\begin{proposition} \ %18
Let \mX be a paracompact Hausdorff space $-$then $\dim_{G \otimes K}X \leq \dim_G X$ for any two abelian groups \mG $\&$ \mK.
\end{proposition}

[This is obvious if either \mG or \mK is trivial, so assume $G \neq 0$ $\&$ $K \neq 0$.

\indent\indent (I) $K = \Z^k$ $(k \geq 1)$.  Here $G \otimes \Z^k$ is a direct sum of copies of \mG, thus the direct sum criterion is applicable.

%%----------------------------------------------------------------------------------------------16
\indent\indent (II) $K = \Z/p^k\Z$ $(k \geq 1)$.  
Case 1: $G_\tor(p) = 0$.  Since $G \otimes \Z/p^k\Z = G/p^kG$, the exactness of 
$0 \ra G \overset{p^k}{\ra} G \ra G/p^kG \ra 0$ gives $\dim_{G \otimes K}X \leq \dim_G X$ (cf. Proposition 15).  
Case 2: $G_\tor(p) \neq 0$.  There are two possibilities: 
$\Z/p\Z \in \sigma(G)$ or 
$\Z/p^\infty\Z \in \sigma(G)$.  
If $\Z/p\Z  \in \sigma(G)$, then $\dim_{\Z/p\Z } X \leq \dim_G X$ (cf. Proposition 17).  And: 
$\dim_{G \otimes K}X \leq \dim_{\Z/p\Z } X$ $(G \otimes \Z/p^k\Z$ is $p$-torsion and 
$\Z/p\Z  \in \sigma(G \otimes \Z/p^k\Z)$ (see the proof of Proposition 16)).  
If $\Z/p^\infty\Z \in \sigma(G)$, then $G = G_\tor(p) \oplus H$, where $G \approx \oplus \Z/p^\infty\Z$, so
$G \otimes K = H \otimes K$.    
Because $H_\tor(p) = 0$, it follows that 
$\dim_{G \otimes K}X \leq$ 
$\dim_{H \otimes K}X \leq$  
$\dim_H X \leq \dim_G X$.

\indent\indent (III) \ Taking into account the direct sum criterion, parts I and II cover the case when \mK is finitely generated.  Finally, an arbitrary \mK is a colimit of its finitely generated subgroups, thus this situation can be handled by an appeal to Proposition 14.]\\

\begingroup%%----------------------------------->>
\fontsize{9pt}{11pt}\selectfont
\textbf{\small EXAMPLE} \quad 
If $G \neq G_\tor$, then $\dim_{\Q} X \leq \dim_G X$.
\vspi
[Proposition 18 implies that $\dim_{G \otimes \Q} X \leq \dim_G X$.  But $G \otimes \Q$ contains $\Q$ as a direct summand.]\\
\endgroup %%------------------------------------<<

\begingroup%%----------------------------------->>
\fontsize{9pt}{11pt}\selectfont
\textbf{\small EXAMPLE} \quad 
Suppose that \mX is an ANR $-$then $\dim_{\Z/p\Z} X = \dim_{\widehat{\Z}_p}X$.
\vspi
[Since $\widehat{\Z}_p \otimes \Z/p\Z \approx \widehat{\Z}_p/p\widehat{\Z}_p \approx \F_p$ 
and 
$\Z/p\Z \in \sigma(\F_p)$, one has 
$\dim_{\Z/p\Z} X \leq$ 
$\dim_{\widehat{\Z}_p \otimes \Z/p\Z} X \leq$ 
$\dim_{\widehat{\Z}_p} X$.  
To establish the inequality in the other direction, put 
$G = \ds\prod\limits_1^\infty \Z/p^n\Z$ $-$then 
$\dim_{\Z/p\Z} X = \dim_G X$ 
(cf. p. \pageref{20.16}) and 
$\dim_G X \geq \dim_{\widehat{\Z}_p} X$ ($G/G_\tor$ is not divisible by $p$.).]
\vspi
[Note: \ If \mX is compact, then 
$\dim_{\Z/p\Z} X = \dim_{\Z_p} X$ 
(cf. p. \pageref{20.17}).]\\ 
\endgroup %%------------------------------------<<

\begingroup%%----------------------------------->>
\fontsize{9pt}{11pt}\selectfont
\textbf{\small EXAMPLE} \quad 
Suppose that \mX is an ANR $-$then $\dim_{\Q} X \leq \dim_G X$ $\forall \ G \neq 0$.\\
\endgroup %%------------------------------------<<

\index{Theorem: Bockstein}
\index{Bockstein Theorem}
\textbf{\small BOCKSTEIN THEOREM} \quad 
Let \mX be a compact Hausdorff space $-$then for any $G \neq 0$, 
$\dim_G X = \sup\limits_{H \in \sigma (G)} \dim_H X$.

[One can suppose for this that \mG is torsion free (cf. Propositions 16 and 17), hence that the elements of 
$\sigma(G)$ are $\Q$ and the $\Z_p$: $pG \neq G$.  We then claim that 
$\dim_G X \leq n$ iff $\dim_{\Q} X \leq n$ $\&$ $\dim_{\Z_p} X \leq n$ $\forall \ p$: $pG \neq G$.  
Indeed, for a given closed subset \mA of \mX, by the universal coefficient theorem, 
$\cH{}^{n+1}(X,A;G) = 0$ iff 
$\cH{}^{n+1}(X,A;\Z) \otimes G = 0$ or still, iff 
$\cH{}^{n+1}(X,A;\Z) \otimes \Q = 0$ $\&$  $\cH{}^{n+1}(X,A;\Z)\otimes \Z_p = 0$ $\forall \ p$: $pG \neq G$, i.e., iff 
$\cH{}^{n+1}(X,A;\Q) = 0$ $\&$  $\cH{}^{n+1}(X,A;\Z_p) = 0$ $\forall \ p$ : $pG \neq G$, as claimed.]

[Note: \ The compactness assumption on \mX in the Bockstein theorem can be relaxed to ``paracompact $\&$ 
$\sigma$-locally compact'' 
(Goto\footnote[2]{\textit{Sci. Rep. Tokyo Kyoiku Daigaku Sect. A} \textbf{10} (1969), 17-23.}).  
However, the Bockstein theorem is not true
%%----------------------------------------------------------------------------------------------17
for an arbitrary metrizable \mX, even if \mX has finite topological dimension 
(Dranishnikov-Repov\u s-\u S\u cepin\footnote[3]{\textit{Topology Proc.} \textbf{18} (1993), 57-73.}).]\\

To illustrate the Bockstein theorem, take $G = \Z$.  Since 
$\sigma(\Z) = \{\Q\} \cup \bigcup\limits_p\{\Z_p\}$ and $\dim_{\Q} X \leq \dim_{\Z_p} X$ $\forall \ p$, it follows that 
$\dim_{\Z} X = \dim_{\Z_p} X$ $(\exists \ p)$.

[Note: \ If $\dim X < \infty$, then $\dim X = \dimx_{\Z} X$ (cf. Proposition 3) and either 
$\dim X - 1 \leq \dimx_{\Q} X$ or $\dim X - 1 \leq \dimx_{\Z/p\Z} X$ $(\exists \ p)$.  Thus fix $p$: 
$\dimx_{\Z} X = \dimx_{\Z_p} X$.  There are now two possibilities: 
$\dimx_{\Z_p} X = \dimx_{\Q} X$ 
from which $\dim X - 1 \leq \dimx_{\Q} X$ 
or 
$\dimx_{\Q} X < \dimx_{\Z_p} X$, from which 
$\dimx_{\Z_p} X \leq \max\{\dimx_{\Q} X,\dimx_{\Z/p^\infty\Z} X + 1\}$ 
$= \dim_{\Z/p^\infty\Z} X + 1$ 
$\leq \dim_{\Z/p\Z} X + 1$ 
$\implies$ 
$\dim X - 1 \leq \dimx_{\Z/p\Z} X$.]\\

\label{20.23}
\begingroup%%----------------------------------->>
\fontsize{9pt}{11pt}\selectfont
\textbf{\small EXAMPLE} \quad 
If \mX is a compact ANR, then $\dim X = \dimx_{\Z/p\Z} X$ $(\exists \ p)$.
\vspi
[For $\dimx_{\Z} X = \dim_{\Z_p} X$ $(\exists \ p)$ and, as noted above 
$\dimx_{\Z_p} X = \dim_{\Z/p\Z} X$.  But here $\dimx_{\Z} X = \dim X$ 
(cf. p. \pageref{20.18}).]\\
\endgroup %%------------------------------------<<

\begingroup%%----------------------------------->>
\fontsize{9pt}{11pt}\selectfont
\label{20.22}
\textbf{\small EXAMPLE} \quad 
Let 
$
\begin{cases}
\ X\\
\ Y
\end{cases}
$
be compact Hausdorff spaces.  Assume: $\dim_G X \leq n$ $-$then $\dim_{\cH{}^{i}(Y;G)}X \leq$ $n + 1$ 
$\forall \ i \geq 0$.
\vspi
[Consider the short exact sequence 
$0 \ra \cH{}^i(Y;\Z) \otimes G \ra \cH{}^i(Y;G) \ra \Tor(\cH{}^{i+1}(Y;\Z),G)$ $\ra 0$ coming from the universal coefficient theorem.  By Proposition 18, 
$\dim_{\cH{}^i(Y;\Z) \otimes G} X \leq \dim_G X$ 
$\leq n$, so it suffices to show that 
$\dim_{\Tor(\cH{}^{i+1}(Y;\Z),G)} X \leq n + 1$ (cf. Proposition 15).  Assuming that 
$\Tor(\cH{}^{i+1}(Y;\Z),G) \neq 0$, $\exists \ p$: $G_{\tor}(p) \neq 0$, hence either $\Z/p\Z \in \sigma(G)$ or 
$\Z/p^\infty\Z \in \sigma(G)$.  But
$\dim_{\Z/p\Z} X \leq \dim_G X$ $\&$ $\dim_{\Z/p^\infty\Z} X \leq \dim_G X$ (Bockstein theorem).  And: 
$\dim_{\Z/p^\infty\Z} X \leq \dim_{\Z/p\Z} X$, 
$\dim_{\Z/p\Z} X \leq \dim_{\Z/p^\infty\Z} X + 1 \leq n + 1$.]\\
\endgroup %%------------------------------------<<

\begingroup%%----------------------------------->>
\fontsize{9pt}{11pt}\selectfont
\textbf{\small FACT} \  
Let \mX be a paracompact Hausdorff space $-$then for any $G \neq 0$, 
$\max\{\dimx_G X, \dimx_{\Q} X + 1\} \geq \sup\limits_{H \in \sigma(G)} \dim_H X$.
\vspi
[Take \mG torsion free and consider the case when \ $H$ $=$ $\Z_p$ \ $(pG \neq G)$.  
One has \ 
$\dim_{\Z_p} X \leq$ 
$\max\{\dim_{\widehat{\Z_p}} X, \dim_{\widehat{\Z_p}/\Z_p} X + 1\}$ 
$= \max\{\dim_{\widehat{\Z_p}} X, \dim_{\Q} X + 1\}$.  
Moreover, 
$\dim_G X \leq n$ $\implies$ $\dim_{\widehat{\Z_p}} X \leq n$.]\\
\endgroup %%------------------------------------<<

\begin{proposition} \ %19
Let 
$
\begin{cases}
\ X\\
\ Y
\end{cases}
$
be nonempty compact Hausdorff spaces.  Assume: 
$
\begin{cases}
\ \dim X\\
\ \dim Y
\end{cases}
< \infty
$
$-$then $\dim_G(X \times Y) \leq \dim_G X + \dim_G Y$ if \mG is torsion free.
\end{proposition}

[With $n = \dim_G X$ $\&$ $m = \dim_G Y$, put $k = n \ + \ m$: $\dim_G(X \times Y) \leq k$ if 
$\dim_{\cH{}^{k-i}(Y,B;G)} X \leq i$ $\forall \ i \geq 0$ and all closed subsets $B \subset Y$ 
(cf. p. \pageref{20.19}).  
Case 1: 
%%----------------------------------------------------------------------------------------------18
$i \leq n - 1$.  Since $k - i \geq m + 1$, we have $\cH{}^{k-i}(Y,B;G) = 0$.  
Case 2: $i \geq n$.  By the universal coefficient theorem,  
$\cH{}^{k-i}(Y,B;G) \approx \cH{}^{k-i}(Y,B;\Z) \otimes G$, hence 
$\dim_{\cH{}^{k-i}(Y,B;G)} X$ $\leq$ $\dim_G X \leq i$ (cf. Proposition 18).]

[Note: \  This inequality is also true if $G = \Z/p\Z$.  For 
$\sigma(\cH{}^{k-i}(Y,B;G))$ $\subset$ $\{\Z/p\Z,\Z/p^\infty\Z\}$ and by the Bockstein theorem 
$\dim_{\cH{}^{k-i}(Y,B;\Z/p\Z)} X = \dim_{\Z/p\Z} X$ (because $\dim_{\Z/p^\infty\Z} X \leq \dim_{\Z/p\Z} X$).]\\

\textbf{\small LEMMA} \ 
Let 
$
\begin{cases}
\ X\\
\ Y
\end{cases}
$
be nonempty compact Hausdorff spaces.  Assume: 
$
\begin{cases}
\ \dim X\\
\ \dim Y
\end{cases}
< \infty
$
$-$then $\dim_G(X \times Y) \geq \dim_G X + \dim_G Y$ for any field \mG.

[Let $n = \dim_G X$, $m = \dim_G Y$ and choose closed subsets $A \subset X$, $B \subset Y$ such that 
$\cH{}^{n}(X,A;G) \neq 0$, $\cH{}^{m}(Y,B;G) \neq 0$.  
The universal coefficient theorem then gives 
$\cH{}^{n}(X,A;\cH{}^{m}(Y,B;G)) \approx$ 
$\cH{}^{n}(X,A;\Z) \otimes \cH{}^{m}(Y,B;G)$.  \ 
But \ 
$\cH{}^{m}(Y,B;G) \approx \otimes \hspace{0.05cm} G$, so 
$\cH{}^{n}(X,A;\cH{}^{m}(Y,B;G)) \neq 0$, which means that 
$\dim_{\cH{}^{m}(Y,B;G)} X \geq n$, thus $\dim_G(X \times Y) \geq n + m$ 
(cf. p. \pageref{20.20}).]\\

\begin{proposition} \ 
Let 
$
\begin{cases}
\ X\\
\ Y
\end{cases}
$
be nonempty compact Hausdorff spaces.  Assume: 
$
\begin{cases}
\ \dim X\\
\ \dim Y
\end{cases}
< \infty
$
$-$then $\dim_G(X \times Y) = \dim_G X + \dim_G Y$ for any field \mG.
\end{proposition}

[This is implied by Proposition 19 and the lemma.]\\

\begin{proposition} \ %20
Let 
$
\begin{cases}
\ X\\
\ Y
\end{cases}
$
be nonempty compact Hausdorff spaces.  Assume: 
$
\begin{cases}
\ \dim X\\
\ \dim Y
\end{cases}
< \infty
$
$-$then $\forall \ G \neq 0$,  $\dim_G(X \times Y) \leq \dim_G X + \dim_G Y + 1$.
\end{proposition}

[With $n = \dim_G X$ $\&$ $m = \dim_G Y$, put $k = n + m + 1$: $\dim_G(X \times Y) \leq k$ if 
$\dim_{\cH{}^{k-i}(Y,B;G)} X \leq i$ $\forall \ i \geq 0$ and all closed subsets $B \subset Y$ 
(cf. p. \pageref{20.21}).  
The case $i \leq n$ being trivial, suppose that $i \geq n + 1$.  
Taking $j \geq i$ and $A \subset X$ closed, repeated use of the universal coefficient theorem leads to 
$\cH{}^{j}(X,A;\cH{}^{k-i}(Y,B;G))$
$\approx$ 
$\cH{}^{j}(X,A;\Z)$ $\otimes$ 
$\cH{}^{k-i}(Y,B;G)$ $\otimes$ 
$\Tor(\cH{}^{j+1}(X,A;\Z) ,\cH{}^{k-i}(Y,B;G))$ 
$\approx$ 
$\cH{}^{j}(X,A;\Z)$ $\otimes$ $[\cH{}^{k-i}(Y,B;\Z)$  
$\otimes$ $G$ $\oplus$ 
$\Tor(\cH{}^{k-i+1}(Y,B;\Z),G)]$ $\oplus$ 
$\Tor(\cH{}^{j+1}(X,A;\Z),\cH{}^{k-i}(Y,B;\Z))$ $\otimes$ $G$ $\oplus$ 
$\Tor(\cH{}^{k - i + 1}(Y,B,$ $\Z),G))$ 
$\approx$ 
$[\cH{}^{j}(X,A;\Z)$ $\otimes$ $\cH{}^{k-i}(Y,B;\Z)$ $\otimes$ $G$ $\oplus$
$\Tor(\cH{}^{j+1}(X,A;\Z),\cH{}^{k-i}(Y,B;\Z)\otimes G)]$ 
$\oplus$ 
$[\cH{}^{j}(X,A;\Z)$ $\otimes$ $\Tor(\cH{}^{k-i+1}(Y,B;\Z),G)$ $\oplus$ 
$\Tor(\cH{}^{j+1}(X,A;\Z),\Tor(\cH{}^{k-i+1}(Y,B;\Z),G))]$
$\approx$ 
$\cH{}^{j}(X,A;\cH{}^{k-i}(Y,B;\Z) \otimes G)$ $\oplus$ 
$\cH{}^{j}(X,A;\Tor(\cH{}^{k-i+1}(Y,B;\Z),G))$.  
By Proposition 18, 
$\dim_{\cH{}^{k-i}(Y,B;\Z),\otimes G} X$ $\leq \dimx_G X < i$, so 
$\cH{}^{j}(X,A;\cH{}^{k-i}(Y,B;\Z) \otimes G) = 0$.  On the other hand, 
$\dim_{\Tor(\cH{}^{k-i+1}(Y,B;\Z),G)} X \leq \dim_G X + 1 \leq i$ 
(imitate the argument used in the second example 
(cf. p. \pageref{20.22}), thus 
$\cH{}^{j}(X,A;\Tor(\cH{}^{k-i+1}(Y,B;\Z),G)) = 0$.  Therefore 
$\dim_{\cH{}^{k-i}(Y,B;G)} X \leq i$, as desired.]\\

%%----------------------------------------------------------------------------------------------19
\begingroup%%----------------------------------->>
\fontsize{9pt}{11pt}\selectfont
Let \mX, \mY be nonempty compact Hausdorff spaces of finite topological dimension.\\
\endgroup %%------------------------------------<<

\begingroup%%----------------------------------->>
\fontsize{9pt}{11pt}\selectfont
\textbf{\small FACT} \  \ 
$\dim_{\Z/p^\infty\Z}(X \times Y)$ $=$ 
$\dim_{\Z/p^\infty\Z} X$ $+$ $\dim_{\Z/p^\infty\Z} Y$
 \ if \ 
$\dim_{\Z/p^\infty\Z} X$ $=$ $\dim_{\Z/p\Z} X$ \  or \ 
$\dim_{\Z/p^\infty\Z} Y$ $=$ $\dim_{\Z/p\Z} Y$, \  otherwise \ 
$\dim_{\Z/p^\infty\Z}(X \times Y)$ $=$ $\dim_{\Z/p^\infty\Z} X$ $+$ $\dim_{\Z/p^\infty\Z} Y + \ 1$ 
$= \dim_{\Z/p\Z}(X \times Y) - 1$.
\vspi
[If the second eventuality obtains, then 
$\dim_{\Z/p^\infty\Z} X < \dim_{\Z/p\Z} X$ 
$\&$  
$\dim_{\Z/p^\infty\Z} Y < \dim_{\Z/p\Z} Y$
$\implies$ 
$\dim_{\Z/p\Z} X + \dim_{\Z/p\Z} Y - 1 = \dim_{\Z/p\Z}(X \times Y) - 1$ (cf. Proposition 20) 
$\leq \dim_{\Z/p^\infty\Z}(X \times Y)$ 
$\leq  \dim_{\Z/p^\infty\Z} X + \dim_{\Z/p^\infty\Z} Y + 1$ (cf. Proposition 21) 
$= (\dim_{\Z/p^\infty\Z} X + 1) + (\dim_{\Z/p^\infty\Z} Y + 1) - 1$ 
$= \dim_{\Z/p\Z} X + \dim_{\Z/p\Z} Y - 1$.]\\
\endgroup %%------------------------------------<<

\begingroup%%----------------------------------->>
\fontsize{9pt}{11pt}\selectfont
\textbf{\small FACT} \  
$\dim_{\Z_p}(X \times Y) = \dim_{\Z_p}X + \dim_{\Z_p}Y$ if 
$\dim_{\Z/p^\infty\Z}X = \dim_{\Z_p}X$ 
and 
$\dim_{\Z/p^\infty\Z}Y = \dim_{\Z_p}Y$, otherwise 
$\dim_{\Z_p}(X \times Y) = \max\{\dim_{\Q}(X \times Y), \dim_{\Z/p^\infty\Z}(X \times Y) + 1\}$.
\vspi
[If the first eventuality obtains, then 
$\dim_{\Z_p}X + \dim_{\Z_p}Y \geq \dim_{\Z_p}(X \times Y)$ (cf. Proposition 19) which is 
$\geq \dim_{\Z/p\Z}(X \times Y) = \dim_{\Z/p\Z}X + \dim_{\Z/p\Z}Y$ (cf. Proposition 20) which is 
$\geq \dim_{\Z/p^\infty\Z}X + \dim_{\Z/p^\infty\Z}Y$ 
$= \dim_{\Z_p}X + \dim_{\Z_p}Y$.]\\
\endgroup %%------------------------------------<<

\begingroup%%----------------------------------->>
\fontsize{9pt}{11pt}\selectfont
\textbf{\small EXAMPLE} \quad 
Given $m$, $n$ and $q$ such that $n \leq m < q \leq n + m$, $\exists$ metrizable compact Hausdorff spaces $X_m$, $X_n$ : 
$\dim X_m = m$, $\dim X_n = n$, and $\dim (X_m \times X_n) = q$.
\vspi
[Specify two Bockstein functions $D_m$, $D_n$ by the following table
\endgroup %%------------------------------------<<
\[
\begin{array}{c c c c c c c c c c c c c c c}
\hline
\Z_2 &&\Z/2\Z &&\Z/2^\infty\Z &&\Q &&\Z_p &&\Z/p\Z &&Z/p^\infty\Z\\
\hline
 m &&1 &&1 &&m &&m &&1 &&1\\
\\[-.25cm]
n &&n &&n-1 &&q-m &&q-m &&q-m &&q-m\\
\end{array}
\]
\begingroup%%----------------------------------->>
\fontsize{9pt}{11pt}\selectfont
and consider the metrizable compact Hausdorff spaces produced by the Dranishnikov realization theorem.]\\
\endgroup %%------------------------------------<<

\begin{proposition} \ %22
Let \mX be a nonempty compact Hausdorff space of finite topological dimension.  Assume: 
$\dim X = \dim_{\Q} X$ or $\dim X = \dim_{\Z/p\Z} X$ $(\exists \ p)$ $-$then $\dim X^n = n \cdot \dim X$.
\end{proposition}

[If $\dim X = \dim_G X$, where $G = \Q$ or $\Z/p\Z$ $(\exists \ p)$, then 
$n \cdot \dim X \geq \dim X^n$ (product theorem) $\geq \dim_G X^n$ (cf. Proposition 2) 
$= n \cdot \dim_G X$ (cf. Proposition 20) $= n \cdot \dim X$.]\\

\begingroup%%----------------------------------->>
\fontsize{9pt}{11pt}\selectfont
\textbf{\small EXAMPLE} \quad 
If \mX is a compact ANR of finite topological dimension, then $\dim X^n = n \cdot \dim X$.
[This is because $\dim X = \dim_{\Z/p\Z}$ $(\exists \ p)$ 
(cf. p. \pageref{20.23}).\\
\endgroup %%------------------------------------<<

\begingroup%%----------------------------------->>
\fontsize{9pt}{11pt}\selectfont
\textbf{\small FACT} \  
Let \mX be a nonempty compact Hausdorff space of finite topological dimension.  Assume: 
$\dim X > \dim_G X$ for $G = \Q$ and $G = \Z/p\Z$ $(\forall \ p)$ $-$then $\dim X^n = n \cdot \dim X - (n - 1)$.\\
\endgroup %%------------------------------------<<

%%----------------------------------------------------------------------------------------------20
\label{19.20}
\begingroup%%----------------------------------->>
\fontsize{9pt}{11pt}\selectfont
\textbf{\small EXAMPLE} \quad 
Suppose that \mX realizes the Bockstein function $\Phi(\Z/p^\infty\Z,n)$ 
(cf. p. \pageref{20.24}) $-$then 
$\dim X = n$ and \mX satisfies the assumption of the preceding result.  Therefore 
$\dim (X \times X) = 2n - 1 < 2n$ 
(cf. p. \pageref{20.25}).\\
\endgroup %%------------------------------------<<

\begin{proposition} \ %23
Let 
$
\begin{cases}
\ X\\
\ Y
\end{cases}
$
be nonempty compact Hausdorff spaces.  Assume: 
$
\begin{cases}
\ \dim X\\
\ \dim Y
\end{cases}
< \infty
$
$-$then $\forall \ G,K \neq 0$, $\dimx_{G \otimes K}(X \times Y) \leq \dim_G X + \dim_K Y$.
\end{proposition}

[Take $k = \dim_G X + \dim_K Y$ and show that 
$\dimx_{\cH{}^{k-1}(Y,B;G \otimes K)} X \leq i$ $\forall \ i \geq 0$ and all closed subsets $B \subset Y$ 
(cf. p. \pageref{20.26}).]\\

Application: Under the assumptions of the preceding proposition, $\dim_R(X \times Y ) \leq \dim_R X + \dim_R Y$ for any ring \mR with unit.

[In fact, \mR is a retract of $R \otimes_{\Z} R$, thus is a direct summand, so 
$\dim_R(X \times Y ) \leq \dim_{R \otimes_{\Z} R} (X \times Y)$ 
$\leq \dim_R X + \dim_R Y$.]\\

\begin{proposition} \ %24
Let 
$
\begin{cases}
\ X\\
\ Y
\end{cases}
$
be nonempty compact Hausdorff spaces.  Assume: 
$
\begin{cases}
\ \dim X\\
\ \dim Y
\end{cases}
< \infty
$
$-$then $\forall \ G,K \neq 0$, $\dimx_{\Tor(G,K)}(X \times Y) \leq \dimx_G X + \dimx_K Y + 1$.
\end{proposition}

[Since $\Tor(G,K) = \Tor(G_\tor,K_\tor)$, one can assume that \mG and \mK are torsion (cf. Proposition 17).  Making the obvious reductions, one can assume further that \mG and \mK are $p$-primary (tacitly, $\Tor(G,K) \neq 0$).  
Case 1: $\Tor(G,K)$ is not divisible by $p$.  
In this situation, either \mG or \mK is not divisible by $p$.  And: 
$\dimx_{\Tor(G,K)}(X \times Y) = \dimx_{\Z/p\Z}(X \times Y)$ (Bockstein theorem) 
$\leq \dimx_{\Z/p\Z} X + \dimx_{\Z/p\Z} Y$.  But either 
$\dimx_{\Z/p\Z} X = \dimx_G X$ or 
$\dimx_{\Z/p\Z} Y = \dimx_K Y$ and at worst, 
$\dimx_{\Z/p\Z} X \leq \dimx_G X + 1$ $\&$
$\dimx_{\Z/p\Z} Y \leq \dimx_K Y + 1$, \mG $\&$ \mK being $p$-primary.  
Case2: $\Tor(G,K)$ is divisible by $p$.  Here, 
$\dim_{\Tor(G,K)}(X \times Y) = \dimx_{\Z/p^\infty \Z}(X \times Y)$ (Bockstein theorem) 
$\leq \dim_{\Z/p^\infty \Z} X + \dimx_{\Z/p^\infty \Z} Y + 1$.  
But 
$\dimx_{\Z/p^\infty \Z} X \leq \dimx_G X$ $\&$
$\dimx_{\Z/p^\infty \Z} Y \leq \dimx_K Y$, \mG $\&$ \mK being $p$-primary.]\\

%%%%%%%%%%%%%%%%%%%%%%%%%%%%%%%%%%%%%%
\setcounter{proposition}{0}
%%%%%%%%%%%%%%%%%%%%%%%%%%%%%%%%%%%%%%
%%%%%%%%%%%%%%%%%%%%%%%%%%%%%%%%%%%%%%

\begin{center}
$\S \ 20$
\\[0.5cm]
$\mathcal{REFERENCES}$\\
\end{center}

\[
\textbf{BOOKS}
\]

\begingroup
\fontsize{9pt}{11pt}\selectfont
\setlength\parindent{0 cm}

[1] \quad Alexandroff, P., 
\textit{Introduction \`a la Th\'eorie Homologique de la Dimension et la Topologie Combinatoire}, 

\hspace{0.8cm}\'Editions Mir (1977).
\\[-.2cm]

[2] \quad Chigogidze, A., 
\textit{Inverse Spectra}, North Holland (1996).
\\[-.2cm]

%Fedorchuk V, and Chigogidze A, Absolute Retracts and Infinite-Dimensional Manifolds, Nauka 1992\\
%{
%[3]  \cyr{Fedorcuk, B. i Qigogidze, A., Absoliutnye retrakty i beskonechnomernye mnoobraziia,}\cyr{Nauka} (1992).\\
%}
[3] \quad {\fontencoding{OT2}\selectfont
Fedorchuk V, Chigogidze A, Absolyutnyy retrakty i beskonechnomernaya mnogoobraziya, Nauka 

\hspace{0.8cm}(1992)}
\\[-.2cm]

[4] \quad Nagami, K., 
\textit{Dimension Theory}, Academic Press (1970).
\\[-.2cm]
\endgroup

\[
\textbf{ARTICLES}
\]

\begingroup
\fontsize{9pt}{11pt}\selectfont
\setlength\parindent{0 cm}

[1] \quad Dranishnikov, A., Homological Dimension Theory, 
\textit{Russian Math. Surveys} \textbf{43} (1988), 11-63.
\\[-.2cm]

[2] \quad Dranishnikov, A. and Dydak, J., Extension Dimension and Extension Types, \textit{Proc. Steklov Inst.} 

\hspace{0.8cm}\textit{Math.} \textbf{212} (1996), 55-88.
\\[-.2cm]

[3] \quad Dranishnikov, A., Repov\u s, D., and \u S\u cepin, E., On Approximation and Embedding Problems for Co-

\hspace{0.8cm}homological Dimension, 
\textit{Topology Appl.} \textbf{55} (1994), 67-86.
\\[-.2cm]

[4] \quad Dranishnikov, A. and \u S\u cepin, E., Cell-Like Maps: The Problem of Raising Dimension, 
\textit{Russian Math}

\hspace{0.8cm}\textit{Surveys} \textbf{41} (1986), 59-111.
\\[-.2cm]

[5] \quad Dydak, J., Cohomological Dimension Theory, 
\textit{Handbook of General Topology}, Elseveir, Amsterdam, 

\hspace{0.8cm}to appear.
\\[-.2cm]
%not clear [5] every came out

[6] \quad Dydak, J. and Walsh, J., Complexes that Arise in Cohomological Dimension Theory: A Unified Ap-

\hspace{0.8cm}proach, 
\textit{J. London Math. Soc.} \textbf{48} (1992), 329-347.
\\[-.2cm]

[7] \quad Dydak, J. and Walsh, J., Dimension Theory and the Sullivan Conjecture, In: 
\textit{Topology Hawaii} K. 

\hspace{0.8cm}Dovermann (ed.), World Scientific (1992), 75-89.
\\[-.2cm]

[8] \quad Edwards, D. and Hastings, H., \u Cech Theory: Its Past, Present, and Future, 
\textit{Rocky Mountain J. Math.} 

\hspace{0.8cm}\textbf{10} (1980), 429-468.
\\[-.2cm]

[9] \quad Kuz'minov, V., Homological Dimension Theory, 
\textit{Russian Math. Surveys} \textbf{23} (1968), 1-45.
\\[-.2cm]

[10] \quad Lacher, R., Cell-Like Mappings and their Generalizations, 
\textit{Bull. Amer. Math. Soc.} \textbf{83} (1977), 

\hspace{0.95cm}495-552.
\\[-.2cm]

[11] \quad Mitchell, W. and Repov\u s, D., The Topology of Cell-Like Mappings, 
\textit{Rend. Sem. Fac. Sci. Univ.}

\hspace{0.95cm}\textit{Cagliari (Suppl.)} \textbf{58} (1988), 265-300.
\\[-.2cm]

[12] \quad Sklyarenko, E., General Homology and Cohomology Theories: Current State and Typical Applica-

\hspace{0.95cm}tions, 
\textit{J. Soviet Math.} \textbf{57} (1991), 3009-3066.
\\[-.2cm]

[13] \quad Sklyarenko, E., Homology and Cohomology Theories of General Spaces, In: 
\textit{General Topology}, EMS 

\hspace{0.95cm}\textbf{50}, Springer Verlag (1996), 119-246.
\\[-.2cm]

[14] \quad Walsh, J., Dimension, Cohomological Dimension, and Cell-Like Mappings, 
\textit{SLN} \textbf{870} (1981), 105-118.
\\[-.2cm]

[15] \quad  Williams, R., A Useful Functor and Three Famous Examples in Topology, 
\textit{Trans. Amer. Math. Soc.} 

\hspace{0.95cm}\textbf{106} (1963), 319-329.

\setlength\parindent{2em}

\endgroup

\newpage
\setcounter{page}{1}
\renewcommand{\thepage}{Index-\arabic{page}}
\backmatter
\bibliography{}
\printindex
\end{document}